\documentclass[letterpaper]{EMSMono}
\usepackage{mathrsfs}
\usepackage[all]{xy}
\usepackage[dvips]{graphicx}
\usepackage{makeidx}
\usepackage{manfnt}
\usepackage{eucal}
\usepackage[bookmarks=false]{hyperref}
\makeindex
\def\N{\mathbb{N}}
\def\Z{\mathbb{Z}}
\def\Q{\mathbb{Q}}
\def\R{\mathbb{R}}
\def\C{\mathbb{C}}

\def\G{\mathbb{G}}
\def\P{\mathbb{P}}
\def\A{\mathbb{A}}
\def\D{\mathbb{D}}

\def\O{\mathscr{O}}
\def\m{\mathfrak{m}}
\def\til#1{\widetilde{#1}}
\def\ovl#1{\overline{#1}}
\def\codim{\mathop{\mathrm{codim}}\nolimits}

\def\card{\mathop{\mathrm{card}}\nolimits}
\def\cpt{\mathop{\mathrm{cpt}}\nolimits}

\def\pc{\mathop{\mathrm{pc}}\nolimits}
\def\obj{\mathop{\mathrm{obj}}\nolimits}
\def\End{\mathop{\mathrm{End}}\nolimits}
\def\Isom{\mathop{\mathrm{Isom}}\nolimits}

\def\loc{\mathop{\mathrm{loc}}\nolimits}
\def\sep{\mathop{\mathrm{sep}}\nolimits}
\def\red{\mathop{\mathrm{red}}\nolimits}
\def\ad{\mathop{\mathrm{ad}}\nolimits}
\def\cl{\mathop{\mathrm{cl}}\nolimits}
\def\rig{\mathop{\mathrm{rig}}\nolimits}
\def\sp{\mathop{\mathrm{sp}}\nolimits}

\def\adic{\mathop{\mathrm{adic}}\nolimits}

\def\dist{\mathop{\mathrm{dist}}\nolimits}
\def\int{\mathop{\mathrm{int}}\nolimits}
\def\can{\mathop{\mathrm{can}}\nolimits}
\def\pf{\mathop{\mathrm{FP}}\nolimits}
\def\Sp{\mathop{\mathrm{Sp}}\nolimits}
\def\Spec{\mathop{\mathrm{Spec}}\nolimits}
\def\spec{\mathop{\mathbf{Spec}}\nolimits}
\def\Spf{\mathop{\mathrm{Spf}}\nolimits}
\def\Spm{\mathop{\mathrm{Spm}}\nolimits}

\def\Spz{\mathop{\mathrm{Spz}}\nolimits}
\def\Spa{\mathop{\mathrm{Spa}}\nolimits}
\def\Proj{\mathop{\mathrm{Proj}}\nolimits}
\def\Frac{\mathop{\mathrm{Frac}}\nolimits}
\def\Sym{\mathop{\mathrm{Sym}}\nolimits}
\def\Aut{\mathop{\mathrm{Aut}}\nolimits}
\def\Hom{\mathop{\mathrm{Hom}}\nolimits}
\def\Homcont{\mathop{\mathrm{Hom.cont}}\nolimits}
\def\bHom{\mathop{\mathbf{Hom}}\nolimits}
\def\Der{\mathop{\mathrm{Der}}\nolimits}
\def\Dercont{\mathop{\mathrm{Der.cont}}\nolimits}
\def\Ext{\mathop{\mathrm{Ext}}\nolimits}
\def\Tor{\mathop{\mathrm{Tor}}\nolimits}
\def\ker{\mathop{\mathrm{ker}}\nolimits}
\def\coker{\mathop{\mathrm{coker}}\nolimits}
\def\image{\mathop{\mathrm{image}}\nolimits}
\def\gr{\mathop{\mathrm{gr}}\nolimits}
\def\lHom{\mathop{\mathcal{H}\!\mathit{om}}\nolimits}

\def\Supp{\mathop{\mathrm{Supp}}\nolimits}
\def\H{\mathrm{H}}
\def\id{\mathrm{id}}
\def\an{\mathrm{an}}
\def\et{\mathrm{\acute{e}t}}
\def\Zar{\mathrm{Zar}}
\def\Itor{\textrm{$I$-}\mathrm{tor}}
\def\Jtor{\textrm{$J$-}\mathrm{tor}}
\def\ator{\textrm{$a$-}\mathrm{tor}}
\def\btor{\textrm{$b$-}\mathrm{tor}}
\def\ftor{\textrm{$f$-}\mathrm{tor}}
\def\opp{\mathrm{opp}}
\def\dl{\langle\!\langle}
\def\dr{\rangle\!\rangle}

\def\Sets{\mathbf{Sets}}
\def\Mod{\mathbf{Mod}}
\def\Alg{\mathbf{Alg}}
\def\Ab{\mathbf{Ab}}
\def\ASh{\mathbf{ASh}}

\def\Bir{\mathbf{Bir}}
\def\Coh{\mathbf{Coh}}

\def\QCoh{\mathbf{QCoh}}
\def\AQCoh{\mathbf{AQCoh}}
\def\coh{\mathrm{coh}}

\def\fpa{\mathrm{FPA}}
\def\FPA{\mathbf{FPA}}

\def\qcoh{\mathrm{qcoh}}
\def\aqcoh{\mathrm{aqcoh}}
\def\Sch{\mathbf{Sch}}
\def\As{\mathbf{As}}
\def\CAs{\mathbf{CAs}}
\def\PSch{\mathbf{PSch}}
\def\PAs{\mathbf{PAs}}

\def\Fs{\mathbf{Fs}}
\def\ARf{\mathbf{ARf}}

\def\CFs{\mathbf{CFs}}
\def\FAs{\mathbf{FAs}}
\def\CFAs{\mathbf{CFAs}}
\def\CRf{\mathbf{CRf}}

\def\CZs{\mathbf{CZs}}
\def\Zs{\mathbf{Zs}}
\def\AZs{\mathbf{AZs}}
\def\CRz{\mathbf{CRz}}
\def\Rf{\mathbf{Rf}}

\def\LRsp{\mathbf{LRsp}}
\def\Rsp{\mathbf{Rsp}}
\def\Vsp{\mathbf{Vsp}}
\def\RVsp{\mathbf{RVsp}}
\def\Top{\mathbf{Top}}
\def\STop{\mathbf{STop}}
\def\CSTop{\mathbf{CSTop}}

\def\DLat{\mathbf{DLat}}
\def\Cov{\mathbf{Cov}}
\def\FM{\mathbf{M}}
\def\FN{\mathbf{FN}}
\def\DC{\mathbf{D}}
\def\CC{\mathbf{C}}
\def\KC{\mathbf{K}}
\def\Ht{\mathrm{Htp}}
\def\Isol{\mathrm{Isol}}
\def\cone{\mathrm{cone}}
\def\Ac{\mathbf{Ac}}

\def\Af{\mathbf{Af}}
\def\Noe{\mathbf{Noe}}
\def\RNoe{\mathbf{RigNoe}}
\def\Adh{\mathbf{Adh}}
\def\fin{\mathrm{fin}}

\def\Qis{\mathbf{Qis}}
\def\amp{\mathrm{amp}}
\def\bd{\mathrm{b}}
\def\RD{\mathrm{R}}
\def\LD{\mathrm{L}}
\def\for{\mathrm{for}}
\def\het#1{{#1}^h}
\def\zat#1{{#1}^{\mathrm{Zar}}}
\def\Cont{\mathrm{cont}}
\def\LT{\mathrm{LT}}
\def\ZR#1{\langle #1\rangle}
\def\ZRT{\mathbf{ZR}}
\def\BL{\mathbf{BL}}
\def\MD{\mathbf{Mdf}}
\def\AId{\mathrm{AId}}
\def\Emb{\mathbf{Emb}}

\def\Ouv{\mathrm{Ouv}}
\def\QCOuv{\mathrm{QCOuv}}
\def\CN{\mathrm{CN}}
\def\top{\mathbf{top}}
\def\sob{\mathrm{sob}}
\def\TOPOI{\mathbf{TOPOI}}
\def\LocTOPOI{\mathbf{LocTOPOI}}

\def\pts{\mathrm{pts}}
\def\open{\mathrm{ouv}}
\def\Open{\mathbf{Ouv}}
\def\Cart{\mathbf{Cart}}
\def\Id{\mathrm{Id}}
\def\Aff{\mathcal{A}\mathcal{R}}
\def\Triples{\mathbf{Tri}}
\def\VTriples{\mathbf{VTri}}
\def\AnTriples{\mathbf{AnTri}}
\def\AnAcsp{\mathbf{An}\Acsp}
\def\Acsp{\mathbf{Adsp}}

\def\Bsp{\mathbf{Bsp}}
\def\Berk#1{#1_{\mathrm{B}}}
\newcommand{\varinjLim}{\displaystyle\mathop{\mathrm{Lim}}_{\longrightarrow}}
\newcommand{\varprojLim}{\displaystyle\mathop{\mathrm{Lim}}_{\longleftarrow}}
\def\qed{\hfill$\square$}
\def\longhookrightarrow{\lhook\joinrel\longrightarrow}
\def\longhookleftarrow{\longleftarrow\joinrel\rhook}

\def\dbend{{\manfntsymbol{127}}}
\def\d@nger{\medbreak\begingroup\clubpenalty=10000
  \def\par{\endgraf\endgroup\medbreak} \noindent\hangindent3em\hangafter=-2
  \hbox to0pt{\hskip-\hangindent\dbend\hfill}}
\outer\def\danger{\d@nger}
\def\dd@nger{\medbreak\begingroup\clubpenalty=10000
  \def\par{\endgraf\endgroup\medbreak} \noindent\hangindent4em\hangafter=-2
  \hbox to0pt{\hskip-\hangindent\dbend\kern1pt\dbend\hfill}}
\outer\def\ddanger{\dd@nger}
\renewenvironment{itemize}%
{%
   \begin{list}{\parbox{1em}{$\bullet$}}
   {%
      \setlength{\topsep}{.5em}
      \setlength{\itemindent}{0em}
      \setlength{\leftmargin}{3em}
      \setlength{\rightmargin}{0em}
      \setlength{\labelsep}{.5em}
      \setlength{\labelwidth}{3em}
      \setlength{\itemsep}{0.2em}
      \setlength{\parsep}{0em}
      \setlength{\listparindent}{0em}
   }
}{%
   \end{list}%
}
\theoremstyle{plain}
\newtheorem{thm}{Theorem}[subsection]
\newtheorem{prop}[thm]{Proposition}
\newtheorem{lem}[thm]{Lemma}
\newtheorem{cor}[thm]{Corollary}
\theoremstyle{definition}
\newtheorem{rem}[thm]{Remark}
\newtheorem{dfn}[thm]{Definition}
\newtheorem{ntn}[thm]{Notation}
\newtheorem{exa}[thm]{Example}
\newtheorem{exas}[thm]{Examples}
\newtheorem{exer}{Exercise}[subsection]
\newtheorem{sit}[thm]{Situation}

\newtheorem{const}[thm]{Construction}
\renewcommand{\thechapter}{\ifnum\arabic{chapter}=0
\arabic{chapter}\else\Roman{chapter}\fi}
\renewcommand{\theexer}{{\bf \thechapter}.\arabic{section}.\arabic{exer}}
\renewcommand{\thesubsubsection}{\thesubsection.\,(\alph{subsubsection})}
\renewcommand{\thefigure}{\arabic{figure}}
\begin{document}
\newlength{\centeroffset}
\setlength{\centeroffset}{-0.5\oddsidemargin}
\addtolength{\centeroffset}{0.5\evensidemargin}
%
%
\thispagestyle{empty}
\vspace*{\stretch{1}}
\begin{center}
{\huge {\bf Foundations of Rigid Geometry I}}

\vspace{7ex}
{\Large {\bf ArXiv version}}
\end{center}

\vspace{\stretch{1}}
\noindent\hspace*{\centeroffset}\makebox[0pt][l]{\begin{minipage}{\textwidth}
\flushright
{\Large Kazuhiro Fujiwara}
\\ \ \\
{\small Graduate School of Mathematics\\
Nagoya University\\
Nagoya $464$-$8502$\\
Japan\\
{\tt fujiwara@math.nagoya-u.ac.jp}}
\end{minipage}}

\vspace{7ex}
\noindent\hspace*{\centeroffset}\makebox[0pt][l]{\begin{minipage}{\textwidth}
\flushright
{\Large Fumiharu Kato}
\\ \ \\
{\small Department of Mathematics\\
Tokyo Institute of Technology\\
Tokyo $152$-$8551$\\
Japan\\
{\tt bungen@math.titech.ac.jp}}
\end{minipage}}

\vspace{\stretch{1}}

\clearpage

\thispagestyle{empty}
\vspace*{40ex}
\hfill
{\Large {\it To the memory of Professor Masayoshi Nagata}}

\vspace{3ex}
\hfill
{\Large {\it and Professor Masaki Maruyama}}

\clearpage

%
%

\setcounter{tocdepth}{3}
\tableofcontents
\pagenumbering{arabic}
\setcounter{page}{1}
\addcontentsline{toc}{chapter}{Introduction}
\markboth{Introduction}{Introduction}
\chapter*{Introduction}

In the early stage of its history, {\em rigid geometry} has been first envisaged in an attempt to construct a {\em non-archimedean analytic geometry}, an analogue over non-archimedean valued fields, such as $p$-adic fields, of complex analytic geometry. 
Later, in the course of its development, rigid geometry has acquired several rich structures, much richer than being merely as a `copy' of complex analytic geometry, which endowed the theory with a great potential of applications. 
This theory is nowadays recognized by many mathematicians of various research fields to be an important and independent discipline in arithmetic and algebraic geometry.
This book is the first volume of our prospective book project, which aims to discuss the rich overall structures of rigid geometry, and to explore its applications. 

Before explaining our general perspective of this book project, we first begin with an overview of the past developments of the theory. 

\medskip\noindent
{\bf 0.\ Background.}
After K.\ Hensel introduced $p$-adic numbers by the end of the 19th century, there had emerged the idea of constructing $p$-adic analogues of already existing mathematical theories that were formerly considered only over real or complex number field.
One of such theories was the theory of complex analytic functions, which had by then already matured to be one of the most successful and rich branches of mathematics. 
Complex analysis saw further developments and innovations later on.
Most importantly, from extensive works on complex analytic spaces and analytic sheaves by H.\ Cartan and J.P.\ Serre in the mid-20th century, after the pioneering work by K.\ Oka, arose a new idea that the theory of complex analytic functions should be regarded as part of complex analytic geometry. 
According to this view, it was only natural to expect the notion of $p$-adic analytic geometry, or more generally, non-archimedean analytic geometry.

However, all first attempts had to encounter with essential difficulties, especially in establishing reasonable local-to-global linkage of the notion of analytic functions.
Such a naive approach is, generally speaking, characterized by its inclination to produce a faithful imitation of complex analytic geometry, which can be already seen at the level of point-sets and topology of the putative analytic spaces.
For example, for the `complex plane' over $\C_p$ ($=$ the completion of the algebraic closure $\ovl{\Q}_p$ of $\Q_p$), one takes the naive point-set, that is, $\C_p$ itself, and the topology just induced from the $p$-adic metric.
Starting from the situation like this, it goes on to construct a locally ringed space $X=(X,\O_X)$ by introducing the sheaf $\O_X$ of `holomorphic functions', of which a conventional definition is something like as follows: $\O_X(U)$ for any open subset $U$ is the set of all functions on $U$ that admit the convergent power series expansion at every point. 
But this leads to an extremely cumbersome situation.
Indeed, since the topology of $X$ is totally disconnected, there are too many open subsets, and this causes the patching of the functions to be extremely `wobbly', so much so that one fails to have good control of global behavior of the analytic functions.
For example, if $X$ is the `$p$-adic Riemann sphere' $\C_p\cup\{\infty\}$, one would expect that $\O_X(X)$ consists only of constant functions, which is, however, far from being true in this situation.

Let us call the problem of this kind the {\em Globalization Problem}.\footnote{This problem is, in classical literature, usually referred to as the problem in analytic continuation.}
Although the problem in its essence may be seen, inasmuch as being concerned with patching of analytic functions, as a topological one, as it will turn out, it deeply links with the issue of notion of points.
In the prehistory of rigid geometry, this Globalization Problem has been one, and perhaps the most crucial one, of the obstacles in the quest for a good non-archimedean analytic geometry.\footnote{In his pioneering works \cite{Kras}\cite{Kras10}, M.\ Krasner conducted deep research into the problem and gave a first general recipe to manage a meaningful analytic continuation of non-archimedean analytic functions.}

\medskip\noindent
{\bf 1.\ Tate's rigid analytic geometry.}
The Globalization Problem found its fundamental solution when J.\ Tate introduced his {\em rigid analytic geometry} \cite{Tate1} at a seminar at Harvard University in 1961.
Tate's motivation was to justify his construction of the so-called {\em Tate curves}, a non-archimedean analogue of $1$-dimensional complex tori, constructed by means of an infinite quotient \cite{Tate2}.\footnote{Elliptic curves and elliptic functions over $p$-adic fields have already been studied by \'E.\ Lutz under the suggestion of A.\ Weil, who is inspired by classical works of Eisenstein (cf.\ \cite[p.\ 538]{Weil}).}
Tate's solution to the problem consists of the following items: 
\begin{itemize}
\item `reasonable' and `sufficiently large' class of analytic functions,
\item `correct' notion of analytic coverings.
\end{itemize}
Here, one can find behind this idea an influence of A.\ Grothendieck in at least two ways: First, Tate introduced spaces via local characterization by means of their function rings, as typified in scheme theory; second, he used the machinery of {\em Grothendieck topology} to define analytic coverings. 

Now, let us briefly review Tate's theory. 
First of all, Tate introduced the category $\mathbf{Aff}_K$ of so-called {\em affinoid algebras} over a complete non-archimedean valuation field $K$.
Each affinoid algebra $\mathcal{A}$, which is a $K$-Banach algebra, is considered to be the ring of `reasonable' analytic functions over the `space' $\Sp\mathcal{A}$, called the {\em affinoid}, which is the corresponding object in the dual category $\mathbf{Aff}^{\opp}_K$ of $\mathbf{Aff}_K$.
Moreover, based on the notion of {\em admissible coverings}, he introduced a new `topology', in fact, a Grothendieck topology, on $\Sp\mathcal{A}$, which we call the {\em admissible} topology.
The admissibility imposes, most importantly, a strong finiteness condition on analytic coverings, which establishes the closer ties between local and global behaviors of analytic functions, as well-described by the famous {\em Tate's acyclicity theorem} ({\bf \ref{ch-rigid}}.\ref{thm-tateacyclicitytheoremclassical}).
An important consequence of this nice local-to-global linkage is the good notion of `patching' affinoids, by which Tate was able to solve the Globalization Problem, and thus to construct global analytic spaces.

In summary, Tate overcame the difficulty by `rigidifying' the topology itself by imposing the admissibility condition, which puts strong restriction on patching of local analytic functions.
It is for this reason that this theory is nowadays called {\em rigid} analytic geometry.

Aside from the fact that it gave a beautiful solution to the Globalization Problem, one should also find it remarkable that Tate's rigid anaytic geometry proved it possible to apply Grothendieck's way of constructing geometric objects to the situation of non-archimedean analytic geometry. 
Thus, rather than complex analytic geometry, Tate's rigid analytic geometry resembles scheme theory.
There seemed to be, however, one technical difference between scheme theory and rigid analytic geometry, which was considered to be quite essential at the time when rigid analytic geometry appeared: Rigid analytic geometry had to use Grothendieck topology, not classical point-set topology.

There is yet another aspect of rigid analytic geometry reminiscent of algebraic geometry.
In order to have a better grasp of the abstractly defined analytic spaces, Tate introduced a notion of points.
He defined points of an affinoid $\Sp\mathcal{A}$ to be maximal ideals of the affinoid algebra $\mathcal{A}$; viz., his affinoids are {\em visualized} by the maximal spectra, that is, the set of all maximal ideals of affinoid algebras, just like affine varieties in the classical algebraic geometry are visualized by the maximal spectra of finite type algebras over a field.
Notice that this choice of points is essentially the same as the naive one that we have mentioned before.
This notion of points was, despite its naivety, considered to be natural, especially in view of his construction of Tate curves, and practically good enough as far as being concerned with rigid analytic geometry over a fixed non-archimedean valued field.\footnote{One might be apt to think that Tate's choice of points is an `easygoing' analogue of the spectra of complex commutative Banach algebras, for which the justification, Gelfand-Mazur theorem, is, however, only valid in complex analytic situation, and actually fails in $p$-adic situation (see below).}

\medskip\noindent
{\bf 2.\ Functoriality and topological visualization.}
Tate's rigid analytic geometry has, since its first appearance, proven itself to be useful for many purposes, and been further developed by several authors.
For example, Grauert-Remmert \cite{GrRe1} laid foundations of topological and ring theoretic aspects of affinoid algebras, and R.\ Kiehl \cite{Kieh1}\cite{Kieh3} promoted the theory of coherent sheaves and their cohomologies on rigid analytic spaces.

However, it was widely perceived that rigid analytic geometry still has some essential difficulties.
For example:

\medskip
$\bullet$ {\sl Functoriality of points does not hold}: If $K'/K$ is an extension of complete non-archimedean valuation fields, then one expects to have, for any rigid analytic space $X$ over $K$, a mapping from the points of the base change $X_{K'}$ to the points of $X$, which, however, does not exist in general in Tate's framework.

\medskip
Let us call this problem the {\em Functoriality Problem}.
The problem is linked with the following more fundamental one:

\medskip
$\bullet$ {\sl The analogue of the Gelfand-Mazur theorem does not hold}: The Gelfand-Mazur theorem states that there exist no Banach field extension of $\C$ other than itself. In non-archimedean situation, in contrast, there exist many Banach $K$-fields other than finite extensions of $K$. This would imply that there should be plenty of `valued points' of an affinoid algebra not factoring through the residue field of a maximal ideal; in other words, there should be much more points than those that Tate has introduced.

\medskip
It is clear that, in order to overcome the difficulties of this kind, one has to change the notion of points.
More precisely, the problem lies in what to choose as the spectrum of an affinoid algebra. 
To this, there are at least two solutions:

\begin{itemize}
\item[(I)] Gromov-Berkovich\index{Berkovich, V.G.} style spectrum;
\item[(II)] Stone-Zariski\index{Zariski, O.} style spectrum.
\end{itemize}

The spectrum of the first style, which turns out to be the `smallest' spectrum to solve the Functoriality Problem in the category of Banach algebras, consists of height one valuations, that is, {\em seminorms} (of a certain type) on affinoid algebras.
The resulting point-sets carry a natural topology, the so-called Gelfand topology.
This kind of spectra is adopted by V.G.\ Berkovich\index{Berkovich, V.G.} in his approach to non-archimedean analytic geometry, so-called {\em Berkovich\index{Berkovich, V.G.} analytic geometry} \cite{Berk1}.
A nice point of this approach is that it can deal with, in principle, a wide class of Banach $K$-algebras, including affinoid algebras, and thus solve the Functoriality Problem (in the category of Banach algebras).
Moreover, the spectra of affinoid algebras in this approach are Hausdorff, hence providing intuitively familiar spaces as the underlying topological spaces of the analytic spaces.

However, the Gelfand topology differs from the admissible topology; it is even weaker, in the sense that, as we will see later, the former topology is a {\em quotient} of the latter. 
Therefore, this topology does not solve the Globalization Problem for affinoid algebras compatibly with Tate's solution, and, in order to do analytic geometry, one still has to use the Grothendieck topology just imported from Tate's theory.

It is thus necessary, in order to solve the Globalization Problem (for affinoids) and Functoriality Problem at the same time, to further improve the notion of points and the topology.
In the second style, the Stone-Zariski style, which we will take up in this book, each spectrum has more points by valuations, not only of height one, but of higher height.\footnote{Notice that this height tolerance is necessary even for rigid spaces defined over complete valuation fields of height one.}
It turns out that the topology on the point-set thus obtained coincides with the admissible topology on the corresponding affinoid, thus solving the Globalization Problem without using Grothendieck topology.
Moreover, the spectra have plenty enough points to solve the Functoriality Problem as well.

As we have seen, to sum up, both the Globalization Problem and the Functoriality Problem are closely linked with the more fundamental issue concerned with the notion of points and topology, that is, the problem for the choice of spectra. 
What lies behind all this is the philosophical tenet that every notion of space in {\em commutative} geometry should be accompanied with `visualization' by means of topological spaces, which we call the {\em topological visualization} (Figure \ref{fig-naivevis}).
\begin{figure}[ht]
$$
\xymatrix@-2ex@R-1ex{\fbox{\begin{minipage}{6em}\begin{center}\smallskip {\small  Commutative geometry}\smallskip\end{center}\end{minipage}}\ar@{=>}[r]&\fbox{\begin{minipage}{6em}\begin{center}\smallskip {\small Topological spaces}\smallskip\end{center}\end{minipage}}}
$$
\caption{Topological visualization}
\label{fig-naivevis}
\end{figure}
It can be said, therefore, that the original difficulties in the early non-archimedean analytic geometry in general, Globalization and Functoriality, are rooted in the lack of adequate topological visualizations.
We will discuss more on this topic later.

\medskip\noindent
{\bf 3.\ Raynaud's\index{Raynaud, M.} approach to rigid analytic geometry.}
To adopt the spectra as in the Stone-Zariski style, in which points are described in terms of valuation rings of arbitrary height, it is more or less inevitable to deal with finer structures, somewhat related to integral structures, of affinoid algebras.\footnote{Such a structure, which we call a {\em rigidification}, will be discussed in detail in {\bf \ref{ch-rigid}}, \S\ref{subsub-rigidificationfadicrings}.
In the original Tate's rigid analytic geometry, the rigidifications are canonically determined by classical affinoid algebras themselves, and this fact should come as the reason why Tate's rigid analytic geometry, unlike more general Huber's adic geometry, could work without reference to integral models of affinoid algebras.}
The approach is, then, further divided into the following two branches:

\begin{itemize}
\item[(II-a)] R.\ Huber's\index{Huber, R.} {\em adic spaces}\footnote{Notice that Huber's theory is based on the choice of integral structures of topological rings. We will give, mainly in {\bf \ref{ch-rigid}}, \S\ref{sec-adicspaces}, a reasonably detailed account of Huber's theory.} \cite{Hube1}\cite{Hube2}\cite{Hube3};
\item[(II-b)] M.\ Raynaud's\index{Raynaud, M.} viewpoint via formal geometry\footnote{By {\em formal geometry}, we mean in this book the geometry of formal schemes, developed by A.\ Grothendieck.} as a {\em model geometry} \cite{Rayn2}.
\end{itemize}

The last approach, which we will adopt in this book, fits in the general framework in which a geometry as a whole is a package derived from a so-called model geometry. 
Here is a toy model that exemplifies the framework:
Consider, for example, the category of finite dimensional $\Q_p$-vector spaces.
We observe that this category is equivalent to the quotient category of the category of finitely generated $\Z_p$-modules mod out by the Serre subcategory consisting of $p$-torsion $\Z_p$-modules, since any finite dimensional $\Q_p$-vector space has a $\Z_p$-lattice, that is, a `model' over $\Z_p$.
This suggests that the overall theory of finite dimensional $\Q_p$-vector spaces is derived from the theory of models, in this case, the theory of finitely generated $\Z_p$-modules.

In our context, what Raynaud discovered on rigid analytic geometry consists of the following:
\begin{itemize}
\item {\em Formal geometry}, which has already been established by Grothendieck prior to Tate's work, can be adopted as a model geometry for Tate's rigid analytic geometry.
\item Consequently, {\em the overall theory of rigid analytic geometry arises from Grothendiek's formal geometry} ({\sc Figure \ref{fig-bridge}}), from which one obtains an extremely useful idea that, between formal geometry and Tate's rigid analytic geometry, one can use theorems in one side to prove theorems in the other.
\end{itemize}

\begin{figure}[ht]
$$
\xymatrix@-2ex@R-1ex{\fbox{\begin{minipage}{6em}\begin{center}\smallskip {\small  Formal geometry}\smallskip\end{center}\end{minipage}}\ar@{=>}[r]&\fbox{\begin{minipage}{6em}\begin{center}\smallskip {\small Rigid analytic geometry}\smallskip\end{center}\end{minipage}}}
$$
\caption{Raynaud's approach to rigid geometry}
\label{fig-bridge}
\end{figure}

To make more precise what it means to say formal geometry can be a model geometry for rigid analytic geometry, consider, just as in the toy model as above, the category of rigid analytic spaces over $K$.
Raynaud showed that the category of Tate's rigid analytic spaces (with some finiteness conditions) is equivalent to the quotient category of the category of finite type formal schemes over the valuation ring $V$ of $K$.
Here the `quotient' means inverting all `modifications' (especially, blow-ups) that are `isomorphisms over $K$', the so-called {\em admissible modifications $($blow-ups$)$}.

There are several impacts of Raynaud's discovery; let us mention a few of them.
First, guided by the principle that rigid analytic geometry is derived from formal geometry, one can build the theory afresh, starting from {\em defining} the category of rigid analytic spaces as the quotient category of the category of formal schemes mod out by all admissible modifications.\footnote{The rigid spaces obtained in this way are, more precisely, what we call {\em coherent} ($=$ quasi-compact and quasi-separated) rigid spaces, from which general rigid spaces are constructed by patching.}
Second, Raynaud's theorem says that {\em rigid analytic geometry can be seen as birational geometry of formal schemes}, a novel viewpoint, which attracts one to explore the link with traditional birational geometry.
Third, as already mentioned above, the bridge between formal schemes and rigid analytic spaces, established by Raynaud's viewpoint, gives rise to fruitful interactions between these theories.
Especially useful is the fact that theorems in the rigid analytic side can be deduced, at least in the situation over complete discrete valuation rings, from theorems in the formal geometry side, available in EGA and SGA works by Grothendieck et al, at least in the Noetherian situation.

\medskip\noindent
{\bf 4.\ Rigid geometry of formal schemes.}
We can now describe, along the line of Raynaud's discovery, the basic framework of our rigid geometry that we are to promote in this book project.
Here is what rigid geometry is for us: {\em Rigid geometry is a geometry obtained from a birational geometry of model geometries.}
This being so, the main purpose of this book project is to develop such a theory for formal geometry, thus generalizing Tate's rigid analytic geometry and providing more general analytic geometry.
Thus to each formal scheme $X$ is associated an object of a resulting category, denoted as usual by $X^{\rig}$, which itself should already be regarded as a rigid space.
Then we define general rigid spaces by patching these objects.
Notice that, here, the rigid spaces are introduced as an `absolute' object without reference to a base space.

Among several classes of formal schemes we start with, one of the most important is the class of what we call {\em locally universally rigid-Noetherian formal schemes} ({\bf \ref{ch-formal}}.\ref{dfn-formalsch}).
The rigid spaces obtained from this class of formal schemes are called {\em locally universally Noetherian rigid spaces} ({\bf \ref{ch-rigid}}.\ref{dfn-universallyadhesiverigidspaces}), which cover most of the analytic spaces that appear in contemporary arithmetic geometry.
Notice that the formal schemes of the above kind are not themselves locally Noetherian.
A technical point imposed from the demand of removing Noetherian hypothesis is that one has to treat {\em non-Noetherian adic rings} of fairly general kind, for which classical theories including EGA do not give us enough tools; for example, valuation rings of arbitrary height are necessary for describing points on rigid spaces, and we accordingly need to treat fairly wide class of adic rings over them for describing fibers of finite type morphisms.

Besides, we would like to propose another viewpoint, which classical theory does not offer.
Among what Raynaud's theory suggests, the most inspiring is, we think, the suggestion that rigid geometry should be a birational geometry of formal schemes.
We would like to put this perspective to be one of the core ideas of our theory.
In fact, as we will see soon below, it tells us what should be the most natural notion of points of the rigid spaces, and thus leads to an extremely rich structures concerned with visualizations (that is, spectra), whereby to obtain a satisfactory solution to the above-mentioned Globalization and Functoriality problems.
Let us see this next in the sequel.

\medskip\noindent
{\bf 5.\ Revival of Zariski's\index{Zariski, O.} idea.}
The birational geometric aspect of our rigid geometry is best explained by means of O.\ Zariski's\index{Zariski, O.} classical approach to birational geometry as a model example.
Around 1940's, in his attempt to attack the desingularization problem for algebraic varieties, Zariski\index{Zariski, O.} introduced abstract Riemann spaces for function fields, which we call {\em Zariski-Riemann spaces}, generalizing the classical valuation-theoretic construction of Riemann surfaces by Dedekind-Weber.
This idea has been applied to several other problems in algebraic geometry, including, for example, Nagata's compactification theorem for algebraic varieties. 

Let us briefly overview Zariski's idea.
Let $Y\hookrightarrow X$ be a closed immersion of schemes (with some finiteness conditions), and set $U=X\setminus Y$.
We consider {\em $U$-admissible modifications} of $X$, which are by definition proper birational maps $X'\rightarrow X$ that is an isomorphism over $U$.
This class of morphisms contains the subclass consisting of {\em $U$-admissible blow-ups}, that is, blow-ups along closed subschemes contained in $Y$.
In fact, $U$-admissible blow-ups are cofinal in the set of all $U$-admissible modifications (due to flattening theorem; cf.\ {\bf \ref{ch-rigid}}, \S\ref{subsub-birationalgeomblowups}).
The Zariski-Riemann space, denoted by $\ZR{X}_U$, is the topological space defined as the projective limit taken along the ordered set of all $U$-admissible modifications, or equivalently, $U$-admissible blow-ups, of $X$.
Especially important is the fact that the Zariski-Riemann space $\ZR{X}_U$ is quasi-compact (essentially due to Zariski \cite{Zar4}; cf.\ {\bf \ref{ch-rigid}}.\ref{thm-classicalZRsp1}), a fact that is crucial in proving many theorems, for example, the above-mentioned Nagata's theorem.\footnote{Zariski-Riemann spaces are also used in O.\ Gabber's\index{Gabber, O.} unpublished works in 1980's on algebraic geometry problems. Its first appearance in literature in the context of rigid geometry seems to be in \cite{Fujiw1}.}

As is classically known, points of the Zariski-Riemann space $\ZR{X}_U$ are described in terms of valuation rings.
More precisely, these points are in one-to-one correspondence with the set of all morphisms, up to equivalence by `domination', of the form $\Spec V\rightarrow X$ where $V$ is a valuation ring (possibly of height $0$) that map the generic point to points in $U$ (see {\bf \ref{ch-rigid}}, \S\ref{subsub-pointsZRsp} for details).
Since the spectra of valuation rings are viewed as `long paths' (cf.\ {\sc Figure \ref{fig-specclosedsetval}} in {\bf \ref{ch-pre}}, \S\ref{sec-val}), one can say intuitively that the space $\ZR{X}_U$ is like a `path space' in algebraic geometry ({\sc Figure \ref{fig-tubularnbd}}).
\begin{figure}[ht]
\begin{center}
\includegraphics[width=8em,clip]{sdu2.eps}
\end{center}
\caption{Set-theoretical description of $\ZR{X}_U$}
\label{fig-tubularnbd}
\end{figure}

Now, what we have meant by putting birational geometry into one of the core ideas in our theory is that we apply Zariski's approach to birational geometry to the main body of our rigid geometry. 
Our basic dictionary for doing this, e.g., for rigid geometry over the $p$-adic field, is as follows:
\begin{itemize}
\item $X$ $\leftrightarrow$ formal scheme of finite type over $\Spf\Z_p$; 
\item $Y$ $\leftrightarrow$ the closed fiber, that is, the closed subscheme defined by `$p=0$'.
\end{itemize}
In this, the notion of $U$-admissible blow-ups corresponds precisely to the admissible blow-ups of the formal schemes.

\medskip\noindent
{\bf 6.\ Birational approach to rigid geometry.}
As we have already mentioned above, our approach to rigid geometry, called the {\em birational approach to rigid geometry}, is, so to speak, the combination of Raynaud's\index{Raynaud, M.} algebro-geometric interpretation of rigid analytic geometry, which regards rigid geometry as a birational geometry of formal schemes, and Zariski's classical birational geometry ({\sc Figure \ref{fig-approach}}).
\begin{figure}[ht]
$$
\fbox{\begin{minipage}{11em}\begin{center}\smallskip {Raynaud's viewpoint of rigid geometry} \smallskip\end{center}\end{minipage}}
+
\fbox{\begin{minipage}{11em}\begin{center}\smallskip {Zariski's viewpoint of birational geometry} \smallskip\end{center}\end{minipage}}
$$
\caption{Birational approach to rigid geometry}
\label{fig-approach}
\end{figure}
Most notably, it will turn out that this approach naturally gives rise to the Stone-Zariski style spectrum, which we have already mentioned before. 

A nice point in combining Raynaud's viewpoint and Zariski's viewpoint is that, while the former gives the fundamental recipe for defining rigid spaces, the latter provides them with a `visualization'. 
Let us see this more precisely, and alongside, explain what kind of visualization we mean here to attach to rigid spaces.

As already described earlier, from an adic formal scheme $X$ (of finite ideal type; cf.\ {\bf \ref{ch-formal}}.\ref{dfn-adicformalschemesoffiniteidealtype}), we obtain the associated rigid space $\mathscr{X}=X^{\rig}$.
Then, suggested by what we have seen in the previous section, we define the {\em associated Zariski-Riemann space} $\ZR{\mathscr{X}}$ as the projective limit 
$$
\ZR{\mathscr{X}}=\varprojlim X',
$$
taken in the category of topological spaces, of all admissible blow-ups $X'\rightarrow X$ (Definition {\bf \ref{ch-rigid}}.\ref{dfn-ZRtriple}). 
We adopt this space $\ZR{\mathscr{X}}$ as the topological visualization of the rigid space $\mathscr{X}$.
In fact, this space is exactly what we have expected as the topological visualization in the case of Tate's theory, since it can be shown that the canonical topology (the projective limit topology) of $\ZR{\mathscr{X}}$ actually coincides with the admissible topology.

To explain more about the visualization of rigid spaces, we would like to introduce three kinds of visualizations in general context.
One is the topological visualization, which we have already discussed.
The second one, which we name {\em standard visualization}, is the one that appears in ordinary geometries, as typified by scheme theory; that is, visualization by locally ringed spaces.
Recall that an affine scheme, first defined abstractly as an object of the dual category of the category of all commutative rings, can be visualized by a locally ringed space supported on the prime spectrum of the corresponding commutative ring.
The third visualization, which we call the {\em enriched visualization}, or just {\em visualization} in this book, is given by what we call {\em triples}\footnote{See {\bf \ref{ch-rigid}}, \S\ref{sub-triples} for the generalities of triples.}: this is an object of the form $(X,\O^+_X,\O_X)$ consisting of a topological space $X$ and two sheaves of topological rings together with an injective ring homomorphism $\O^+_X\hookrightarrow\O_X$ that identifies $\O^+_X$ with an open subsheaf of $\O_X$ such that the pairs $X=(X,\O_X)$ and $X^+=(X,\O^+_X)$ are locally ringed spaces; normally speaking, $\O_X$ is regarded as the structure sheaf of $X$, while $\O^+_X$ represents the enriched structure, such as an integral structure (whenever it makes sense) of $\O_X$.

The enriched visualization is typified by rigid spaces.
Indeed, the Zariski-Riemann space $\ZR{\mathscr{X}}$ has two natural structure sheaves, the {\em integral structure sheaf} $\O^{\int}_{\mathscr{X}}$, defined as the inductive limit of the structure sheaves of all admissible blow-ups of $X$, and the {\em rigid structure sheaf} $\O_{\mathscr{X}}$, obtained from $\O^{\int}_{\mathscr{X}}$ by `inverting the ideal of definition'.
What is intended here is that, while the rigid structure sheaf $\O_{\mathscr{X}}$ should, as in Tate's rigid analytic geometry, normally come as the `genuine' structure sheaf of the rigid space $\mathscr{X}$, the integral structure sheaf $\O^{\int}_{\mathscr{X}}$ represents its integral structure.
These data comprise the triple
$$
\ZRT(\mathscr{X})=(\ZR{\mathscr{X}},\O^{\int}_{\mathscr{X}},\O_{\mathscr{X}}),
$$
called the {\em associated Zariski-Riemann triple}, which gives the enriched visualization of the rigid space $\mathscr{X}$.
That the rigid structure sheaf should be {\em the} structure sheaf of $\mathscr{X}$ means that the locally ringed space $(\ZR{\mathscr{X}},\O_{\mathscr{X}})$ visualizes the rigid space in an ordinary sense, that is, in the sense of standard visualization.

Notice that the Zariski-Riemann triple $\ZRT(\mathscr{X})$ for a rigid space $\mathscr{X}$ coincides with Huber's\index{Huber, R.} adic space associated to $\mathscr{X}$; in fact, the notion of Zariski-Riemann triple not only gives an interpretation of adic spaces, but it also gives a foundation for them via formal geometry, which we establish in this book; see {\bf \ref{ch-rigid}}, \S\ref{sub-rigidgeomadicspaces} for more details.

Figure \ref{fig-approach2} illustrates the basic design of our birational approach to rigid geometry, summarizing all what we have discussed so far.
\begin{figure}[ht]
$$
\begin{xy}
(0,0)="FG"*{\fbox{\begin{minipage}{6em}\begin{center}\smallskip {\small  Formal Geometry}\smallskip\end{center}\end{minipage}}},+/r9em/="RG"*{\fbox{\begin{minipage}{6em}\begin{center}\smallskip {\small Rigid Geometry}\smallskip\end{center}\end{minipage}}},+/r9em/="TR"*{\fbox{\begin{minipage}{6em}\begin{center}\smallskip {\small Triples}\smallskip\end{center}\end{minipage}}},+<0em,-12ex>="LRS"*{\fbox{\begin{minipage}{6em}\begin{center}\smallskip {\small Locally ringed spaces}\smallskip\end{center}\end{minipage}}},+<0em,-12ex>="TS"*{\fbox{\begin{minipage}{6em}\begin{center}\smallskip {\small Topological spaces}\smallskip\end{center}\end{minipage}}},"LRS"+<-3.9em,.1ex>*{>},"TS"+<-3.9em,0ex>*{>}
\ar@{=>}^{(\ast1)}"FG"+/r3.7em/;"RG"+/l3.7em/
\ar@{=>}^{(\ast2)}"RG"+/r3.7em/;"TR"+/l3.7em/
\ar@{=>}"TR"+<0em,-3.2ex>;"LRS"+<0em,4.6ex>
\ar@{=>}"LRS"+<0em,-4.6ex>;"TS"+<0em,4.6ex>
\ar@/_.8pc/@{=}_{(\ast3)}"RG"+<2em,-4.6ex>;"LRS"+/l3.9em/
\ar@(d,l)@{=}_{(\ast4)}"RG"+<0em,-4.6ex>;"TS"+/l3.9em/
\end{xy}
$$
\caption{Birational approach to rigid geometry}
\label{fig-approach2}
\end{figure}
The figure shows a `commutative' diagram, in which the arrow $(\ast1)$ is Raynaud's approach to rigid geometry (Figure \ref{fig-bridge}), and the arrow $(\ast2)$ is the enriched visualization by Zariski-Riemann triples, coming from Zariski's viewpoint. The other visualizations are also indicated in the diagram, the standard visualization by $(\ast3)$, and the topological visualization by $(\ast4)$; the right-hand vertical arrows represent the respective forgetful functors.

All these are the outline of what we will discuss in this volume.
Here, before finishing this overview, let us add a few words on the outgrowth of our theory.
Our approach to rigid geometry, in fact, gives rise to a new perspective of rigid geometry itself: {\em Rigid geometry in general is an analysis along a closed subspace in a ringed topos}. 
This idea, which tells us what {\em rigid-geometrical idea} in mathematics should ultimately be, is linked with the idea of {\em tubular neighborhoods} in algebraic geometry, already discussed in \cite{Fujiw1}.
From this viewpoint, Raynaud's choice, for example, of formal schemes as models of rigid spaces can be interpreted as capturing the `tubular neighborhoods' along a closed subspace by means of the formal completion. 
Now that there are several other ways to capture such structures, e.g., henselian schemes etc., there are several other choices for the model geometry of rigid geometry.\footnote{There is, in addition to formal geometry and henselian geometry, the third possibility for the model geometry, by {\em Zariskian schemes}.
We put a general account of the theory of Zariskian schemes and the associated rigid spaces, so-called, {\em rigid Zariskian spaces}, in the appendices {\bf \ref{ch-formal}}, \S\ref{sec-zariskianschemes} and {\bf \ref{ch-rigid}}, \S\ref{sec-rigidzariskiansp}.}
This yields several variants, e.g.\ rigid henselian geometry, rigid Zariskian geometry, etc., all of which are encompassed within our birational approach.\footnote{The reader might notice that this idea is also related to the cdh-topology in the theory of motivic cohomology.}

\medskip\noindent
{\bf 7.\ Relation with other theories.}
In the first three sections {\bf \ref{ch-rigid}}, \S\ref{sec-adicspaces}, {\bf \ref{ch-rigid}}, \S\ref{sec-berkovich}, and {\bf \ref{ch-rigid}}, \S\ref{sec-valspecbanach} of the appendices to Chapter {\bf \ref{ch-rigid}}, we give the comparisons of our theory with other theories related to rigid geometry.
Here we give a digest of the contents of these sections for the reader's convenience.\footnote{A.\ Abbes has recently published another foundational book \cite{Abbes} on rigid geometry, in which, similarly to ours, he developed and generalized Raynaud's approach to rigid geometry.}

\medskip
$\bullet$ {\bf Relation with Tate's rigid analytic geometry.}
Let $V$ be an $a$-adically complete valuation ring of height one, and set $K=\Frac(V)$ (the fractional field), which is a complete non-archimedean valued field with a non-trivial valuation $\|\cdot\|\colon K\rightarrow\R_{\geq 0}$.
In {\bf \ref{ch-rigid}}, \S\ref{subsub-classicalpoints} we will define the notion of {\em classical points} (in the sense of Tate) for rigid spaces of a certain kind including locally of finite type rigid spaces over $\mathscr{S}=(\Spf V)^{\rig}$.
If $\mathscr{X}$ is a rigid space of the latter kind, it will turn out that the classical points of $\mathscr{X}$ are reduced zero dimensional closed subvarieties in $\mathscr{X}$ (cf.\ {\bf \ref{ch-rigid}}.\ref{prop-pointlikestr2}).

We set $\mathscr{X}_0$ to be the set of all classical points of $\mathscr{X}$.
The assignment $\mathscr{X}\mapsto\mathscr{X}_0$ has several nice properties, some of which are put together into the notion of {\em $($continuous$)$ spectral functor} (cf.\ {\bf \ref{ch-rigid}}, \S\ref{sub-pointspoints}).
Among them is an important property that classical points detect quasi-compact open subspaces: for quasi-compact open subspaces $\mathscr{U},\mathscr{V}\subseteq\mathscr{X}$, $\mathscr{U}_0=\mathscr{V}_0$ implies $\mathscr{U}=\mathscr{V}$.
In view of all this, one can introduce on $\mathscr{X}_0$ a Grothendieck topology $\tau_0$ and sheaf of rings $\O_{\mathscr{X}_0}$, which are naturally constructed from the topology and the structure sheaf of $\mathscr{X}$; for example, for a quasi-compact open subspace $\mathscr{U}\subseteq\mathscr{X}$, $\mathscr{U}_0$ is an admissible open subset of $\mathscr{X}_0$, and we have $\O_{\mathscr{X}_0}(\mathscr{U}_0)=\O_{\mathscr{X}}(\ZR{\mathscr{U}})$.
It will turn out that the resulting triple $\mathscr{X}_0=(\mathscr{X}_0,\tau_0,\O_{\mathscr{X}_0})$ is a Tate's rigid analytic variety over $K$, and thus one has the canonical functor
$$
\mathscr{X}\longmapsto\mathscr{X}_0
$$
from the category of locally of finite type rigid spaces over $\mathscr{S}$ to the category of rigid analytic varieties over $K$.

\medskip\noindent
{\bf Theorem} (Theorem {\bf \ref{ch-rigid}}.\ref{thm-rigidanalyticcomparison}, Corollary {\bf \ref{ch-rigid}}.\ref{cor-rigidanalyticcomparison}). {\it The functor $\mathscr{X}\mapsto\mathscr{X}_0$ is a categorical equivalence from the category of quasi-separated locally of finite type rigid spaces over $\mathscr{S}=(\Spf V)^{\rig}$ to the category of quasi-separated Tate analytic varieties over $K$. Moreover, under this functor, affinoids $($resp.\ coherent spaces$)$ correspond to affinoid spaces $($resp.\ coherent analytic spaces$)$.}

\medskip
Notice that the Raynaud's theorem (the existence of formal models) gives the canonical quasi-inverse functor to the above functor.\footnote{To show the theorem, we need Gerritzen-Grauert theorem \cite{GG}, which we assume whenever discussing Tate's rigid analytic geometry. Notice that, when it comes to the rigid geometry over valuation rings, this volume is self-contained only with this exception. We will prove Gerritzen-Grauert theorem without vicious circle in the next volume.}

\medskip
$\bullet$ {\bf Relation with Huber's adic geometry.}
As we have already remarked above, the Zariski-Riemann triple $\ZRT(\mathscr{X})$, at least in the situation as before, is an adic space. 
This is true in much more general situation, for example, in case $\mathscr{X}$ is {\em locally universally Noetherian} ({\bf \ref{ch-rigid}}.\ref{dfn-universallyadhesiverigidspaces}).
In fact, by the enriched visualization, we have the functor
$$
\ZRT\colon \mathscr{X}\longmapsto\ZRT(\mathscr{X})
$$
from the category of locally universally Noetherian rigid spaces to the category of adic spaces (Theorem {\bf \ref{ch-rigid}}.\ref{thm-rigidgeomadicspaces}), which gives rise to a categorical equivalence in most important cases.
In particular, we have:

\medskip\noindent
{\bf Theorem} (Theorem {\bf \ref{ch-rigid}}.\ref{thm-rigidgeomadicspacesff1}){\bf .}
{\it Let $\mathscr{S}$ be a locally universally Noetherian rigid space. 
Then $\ZRT$ gives a categorical equivalence from the category of locally of finite type rigid spaces over $\mathscr{S}$ to the category of adic spaces locally of finite type over $\ZRT(\mathscr{S})$.}

\medskip
$\bullet$ {\bf Relation with Berkovich analytic geometry.}
Let $V$ and $K$ be as before.
We will construct a natural functor
$$
\mathscr{X}\longmapsto\Berk{\mathscr{X}}
$$
from the category of locally quasi-compact\footnote{Note that, if $\mathscr{X}$ is quasi-separated, then $\mathscr{X}$ is locally quasi-compact if and only if $\ZR{\mathscr{X}}$ is {\em taut} in the sen of Huber \cite[5.1.2]{Hube3} (cf.\ {\bf \ref{ch-pre}}.\ref{rem-taut}).} ({\bf \ref{ch-rigid}}.\ref{dfn-locallycompactspacerigid}) and locally of finite type rigid spaces over $\mathscr{S}=(\Spf V)^{\rig}$ to the category of strictly $K$-analytic spaces (in the sense of Berkovich).

\medskip\noindent
{\bf Theorem} (Theorem {\bf \ref{ch-rigid}}.\ref{thm-berkovichcomparison}){\bf .}
{\it The functor $\mathscr{X}\mapsto\Berk{\mathscr{X}}$ gives a categorical equivalence from the category of all locally quasi-compact locallly of finite type rigid spaces over $(\Spf V)^{\rig}$ to the category of all strictly $K$-analytic spaces.
Moreover, $\Berk{\mathscr{X}}$ is Hausdorff $($resp.\ paracompact Hausdorff, resp.\ compact Hausdorff$)$ if and only if $X$ is quasi-separated $($resp.\ paracompact and quasi-separated, resp.\ coherent$)$.}

\medskip
The underlying topological space of $\Berk{\mathscr{X}}$ is what we call the {\em separated quotient} ({\bf \ref{ch-rigid}}, \S\ref{subsub-separation}) of $\ZR{\mathscr{X}}$, denoted by $[\mathscr{X}]$, which comes with the quotient map $\sep_{\mathscr{X}}\colon\ZR{\mathscr{X}}\rightarrow[\mathscr{X}]$ ({\em separation map}).
In particular, the topology of $\Berk{\mathscr{X}}$ is the {\em quotient} topology of the topology of $\ZR{\mathscr{X}}$.

\medskip
Figure \ref{fig-relationother} illustrates the interrelations among those theories we have discussed so far.
In the diagram,
\begin{itemize}
\item the functors $(\ast1)$ $(\ast2)$ are fully faithful; the functor $(\ast3)$, defined on locally quasi-compact rigid analytic spaces, is fully faithful to the category of strictly $K$-analytic spaces;
\item the functor $(\ast4)\colon\mathscr{X}\rightarrow\mathscr{X}_0$, defined on locally of finite type rigid spaces over $(\Spf V)^{\rig}$, is quasi-inverse to $(\ast1)$ restricted on quasi-separated spaces;
\item the functor $(\ast5)$ is given by the enriched visualization, defined on locally universally Noetherian rigid spaces; it is fully faithful in practical situations including those of locally of finite type rigid spaces over a fixed locally universally Noetherian rigid space, and of rigid spaces of type (N) ({\bf \ref{ch-rigid}}.\ref{thm-rigidgeomadicspacesff2});
\item the functor $(\ast6)\colon\mathscr{X}\mapsto\Berk{\mathscr{X}}$, defined on locally quasi-compact locally of finite type rigid spaces over $(\Spf V)^{\rig}$, gives a categorical equivalence to the category of strictly $K$-analytic spaces.
\end{itemize}
\begin{figure}[ht]
$$
\begin{xy}
(0,0)="CR"*{\fbox{\begin{minipage}{8em}\begin{center}\smallskip {\small Tate's rigid analytic varieties}\smallskip\end{center}\end{minipage}}},+<12em,12ex>="RS"*{\fbox{\begin{minipage}{6em}\begin{center}\smallskip {\small Rigid spaces (in our sense)}\smallskip\end{center}\end{minipage}}},+<0em,-12ex>="AS"*{\fbox{\begin{minipage}{6em}\begin{center}\smallskip {\small Adic spaces}\smallskip\end{center}\end{minipage}}},+<0em,-12ex>="BS"*{\fbox{\begin{minipage}{6em}\begin{center}\smallskip {\small Berkovich spaces}\smallskip\end{center}\end{minipage}}}
\ar@/^.4pc/^{(\ast1)}"CR"+<2em,4.4ex>;"RS"+<-3.7em,0ex>
\ar@{.>}@/_.3pc/^{(\ast4)}"RS"+<-3.7em,-1.1ex>;"CR"+<2.7em,4.4ex>
\ar_{(\ast2)}"CR"+<4.7em,0ex>;"AS"+<-3.7em,0ex>
\ar@{.>}@/_.4pc/_{(\ast3)}"CR"+<2em,-4.4ex>;"BS"+<-3.7em,0ex>
\ar@{.>}^{(\ast5)}"RS"+<0em,-4.4ex>;"AS"+<0em,3.2ex>
\ar@{.>}@/^1.6pc/^{(\ast6)}"RS"+<3.7em,0ex>;"BS"+<3.7em,0ex>
\end{xy}
$$
\caption{Relation with other theories}
\label{fig-relationother}
\end{figure}

Finally, we would like to mention that it has recently become known to the experts that some of the non-archimedean spaces that come naturally in contemporary arithmetic geometry may not possibly handled in Berkovich's analytic geometry (e.g.\ \cite[4.4]{Hell}).
This state of affair makes it important to investigate in detail the relationship between Berkovich's analytic geometry and rigid geometry (or adic geometry).
In {\bf \ref{ch-rigid}}, \S\ref{sub-vsbr-NAASBT}, we will study a spectral theory of filtered rings and introduce a new category of spaces, the so-called metrized analytic spaces.
This new notion of spaces generalizes Berkovich's $K$-analytic spaces, and gives a clear picture of the comparison; see {\bf \ref{ch-rigid}}, \S\ref{subsub-vsbr-BerkovichRmetrizedanalyticspace}.
Also, the newly introduced spaces turn out to be equivalent to Kedlaya's reified adic spaces \cite{Kedlaya}, to which our filtered ring approach in this book offers a new perspective. 

\medskip\noindent
{\bf 8.\ Applications.}
We expect that our rigid geometry will have rich applications, not only in arithmetic geometry, but also in various other fields.
A few of them have already been sketched in \cite{FK2}, which include
\begin{itemize}
\item arithmetic moduli spaces (e.g.\ Shimura varieties) and their compactifications,
\item trace formula in characteristic $p>0$ (Deligne's conjecture).
\end{itemize}

In addition to these, since our theory has set out from Zariski's birational geometry, applications to problems in birational geometry, modern or classical, are also expected.
For example, this volume already contains Nagata's compactification theorem for schemes and a proof of it ({\bf \ref{ch-rigid}}, \S\ref{sec-nagataembedding}), as an application of the general idea of our rigid geometry to algebraic geometry.

Some other prospective applications may be to $p$-adic Hodge theory (cf.\ \cite{Scholze1}\cite{Scholze2}) and to rigid cohomology theory for algebraic varieties in positive characteristic.
Here the visualization in our sense of rigid spaces will give concrete pictures for tubes and the dagger construction.
One of such applications in this direction may be named as:
\begin{itemize}
\item $p$-adic Hodge theory vs.\ rigid cohomology.
\end{itemize}
Finally, let us mention that the applications to
\begin{itemize}
\item moduli of Galois representations,
\item mirror symmetry,
\end{itemize}
the second of which has been first envisaged by M.\ Kontsevich, should be among the future challenges.

\medskip\noindent
{\bf 9.\ Contents of this book.}
There are two basic policies in designing the contents of this book, both of which may justify its length.
First, in addition to the role as a front-line exposition presenting new theories and results, we would also like to endow this book an encyclopedic role.
It contains, consequently, as many notions and concepts, hopefully with little omission, that should come about as basic and important ones for present and future use, as possible.

Second, we have paid much effort to make this book as self-contained as possible. 
All results that sit properly inside the main body of our argument are always followed by proofs, except for some minor or not-too-difficult lemmas, some of which are put in the end of each section as an exercise; even in this case, if the result is used in the main text, we give a detailed hint in the end of the book, which, in many cases, almost proves the assertion.
Notice that, because of several laborious requirements on the groundworks, such as removing Noetherian hypothesis, the self-containedness norm should also apply to many of the preliminary parts.

This volume consists of the following three chapters:
\begin{itemize}
\item Chapter {\bf \ref{ch-pre}}. Preliminaries
\item Chapter {\bf \ref{ch-formal}}. Formal geometry
\item Chapter {\bf \ref{ch-rigid}}. Rigid spaces
\end{itemize}
Let us briefly overview the contents of each chapter.
More detailed summary will be given at the beginning of each chapter.

Chapter {\bf \ref{ch-pre}} collects preliminaries, which, however, contain new results. 
The sections {\bf \ref{ch-pre}}, \S\ref{sec-language} to {\bf \ref{ch-pre}}, \S\ref{sec-topologicalringsmodules} give necessary preliminaries on set theory, category theory, general topology, homological algebra, etc.
In the general topology section, we put emphasis on Stone duality between topological spaces and lattices.
In {\bf \ref{ch-pre}}, \S\ref{sec-pairs} and {\bf \ref{ch-pre}}, \S\ref{sec-aadicallycompval}, we will conduct thorough study on topological and algebraic aspects of topological rings and modules. 
This part of the preliminaries will be the bases of the next chapter, the general theory of formal geometry.

Chapter {\bf \ref{ch-formal}} is devoted to formal geometry.
Here it is essential to treat non-Noetherian formal schemes of a certain kind, e.g.\ finite type formal schemes over an $a$-adically complete valuation ring of arbitrary height, for a functoriality reason (as stated in 4 above).
Since this kind of generalities seems missing in the past literature, we provide here a self-contained and fully systematic theory of formal geometry, generalizing many of the theorems in \cite[$\mathbf{III}$]{EGA}.
To this end, we will introduce several new notions of finiteness conditions {\em outside the ideal of definition} and show that they provide feasible and versatile theory of formal schemes.

Chapter {\bf \ref{ch-rigid}} is the main part of this volume, in which we develop rigid geometry, based on the foundational works done in the previous chapters. 
The geometrical theory of rigid spaces that we give in this chapter include
\begin{itemize}
\item cohomology theory of coherent sheaves ({\bf \ref{ch-rigid}}, \S\ref{sec-coherentsheavesrigid}, \S\ref{sec-affinoids}); finiteness ({\bf \ref{ch-rigid}}.\ref{thm-finirigid});
\item local and global studies of morphisms ({\bf \ref{ch-rigid}}, \S\ref{sub-basicmorproprigid});
\item classification of points ({\bf \ref{ch-rigid}}, \S\ref{sec-localring}, \S\ref{sub-classifyingvaluations});
\item GAGA ({\bf \ref{ch-rigid}}, \S\ref{sec-GAGA});
\item relation with other theories ({\bf \ref{ch-rigid}}, \S\ref{sec-adicspaces}, \S\ref{sec-berkovich}, \S\ref{sec-valspecbanach}).
\end{itemize}

There are of course many other important topics that are not dealt with in this volume.
As we will see later, some of them, including several important applications, will be contained in the future volumes.

\medskip\noindent
{\bf 10.\ Use of algebraic spaces.}
In {\bf \ref{ch-formal}}, \S\ref{sec-formalalgsp} we develop a full-fledged theory of formal algebraic spaces.
It is, in fact, one of the characteristic features of our approach to rigid geometry that we allow formal algebraic spaces, not only formal schemes, to be formal models of rigid spaces.
The motivation mainly comes from the applications to algebraic geometry.

In algebraic geometry, while it is often difficult to show that spaces, such as moduli spaces, are represented by schemes, the representability by algebraic spaces is relatively easy to see, thanks to M.\ Artin's formal algebraization theorem \cite{Art1}.
Therefore, taking algebraic spaces into the scope increases the flexibility of the theory.
In order to incorporate algebraic spaces into our rigid geometry, one first needs to discuss formal algebraic spaces, some of which appear as the formal completion of algebraic spaces, and then to proceed to the rigid spaces associated to them.
Now the important fact is the following: Although formal algebraic spaces seem to constitute, via Raynaud's recipe, a new category of rigid spaces that enlarges the already existing category of rigid spaces derived from formal schemes, they actually do not; viz., we do not have to enlarge the category of rigid spaces by this generalization.
The clue for this fact is the following theorem, which we shall prove in the future volume:

\medskip\noindent
{\bf Theorem (Equivalence Theorem).}
{\it Let $X$ be a coherent adic formal algebraic space of finite ideal type. 
Then there exists an admissible blow-up $X'\rightarrow X$ from a formal scheme $X'$.
Therefore, the canonical functor 
$$
\left\{\begin{minipage}{7em}\smallskip{\rm {\small coherent adic formal schemes of finite ideal type}}\smallskip\end{minipage}\right\}_{\bigg/\left\{\begin{minipage}{4em}\smallskip{\rm {\small admissible blow-ups}}\smallskip\end{minipage}\right\}}
\longrightarrow\left\{\begin{minipage}{8em}\smallskip{\rm {\small coherent adic formal algebraic spaces of finite ideal type}}\smallskip\end{minipage}\right\}_{\bigg/\left\{\begin{minipage}{4em}\smallskip{\rm {\small admissible blow-ups}}\smallskip\end{minipage}\right\}}
$$
is a categorical equivalence.}

\medskip
The theorem shows that, up to admissible blow-ups, formal algebraic spaces simply fall into formal schemes, and thus define the associated rigid space $X^{\rig}$ just `as usual'.
As for GAGA, we can generalize the definition of GAGA functor for algebraic spaces (using compactification theorem of Nagata type for algebraic spaces).\footnote{This `analytification of algebraic spaces' is already deeply considered and established by Conrad\index{Conrad, B.}-Temkin\index{Temkin, M.} \cite{CT} over complete non-archimedean valued fields.}

\medskip\noindent
{\bf 11.\ Properness in rigid geometry.}
In rigid geometry, we have the following three natural definitions of properness: A morphism $\varphi\colon\mathscr{X}\rightarrow\mathscr{Y}$ of coherent rigid spaces is proper if either one of the following conditions is satisfied:
\begin{itemize}
\item[(1)] $\varphi$ is universally closed ({\bf \ref{ch-rigid}}.\ref{dfn-seppropmorrigid1}), separated, and of finite type;
\item[(2)] ({\em Raynaud\index{Raynaud, M.} properness}) there exists a proper formal model $f\colon X\rightarrow Y$ of $\varphi$.
\item[(3)] ({\em Kiehl\index{Kiehl, R} properness}) $\varphi$ is separated of finite type, and there exist an affinoid covering $\{\mathscr{U}_i\}_{i\in I}$ and, for each $i\in I$, a pair of finite affinoid coverings $\{\mathscr{V}_{ij}\}_{j\in J_i}$ and $\{\mathscr{V}'_{ij}\}_{j\in J_i}$ of $\varphi^{-1}(\mathscr{U}_i)$ indexed by a common set $J_i$ such that, for any $j\in J_i$, $\mathscr{V}_{ij}\subseteq\mathscr{V}'_{ij}$ and $\mathscr{V}_{ij}$ is relatively compact in $\mathscr{V}'_{ij}$ over $\mathscr{U}_i$ (in the sense of Kiehl).
\end{itemize}

Historically, properness in Tate's rigid geometry has been first defined by Kiehl\index{Kiehl, R} by the condition (3) in his work \cite{Kieh1} on finiteness theorem.
This condition, existence of affinoid enlargements, stems from the general idea by Cartan-Serre and H.\ Grauert for proving finiteness of cohomologies of coherent sheaves.
While the equivalence of (1) and (2) is an easy exercise, the equivalence of (2) and (3), especially (2) $\Rightarrow$ (3), is a very deep theorem. 
L\"utkebohmert's\index{Lutkebohmert@L\"utkebohmert, W} 1990 paper \cite{Lutk1} proves this for rigid spaces of finite type over $(\Spf V)^{\rig}$ where $V$ is a complete discrete valuation ring.
In this book, we temporarily define properness by the condition (1) (and hence equivalently by (2)), and postpone the proof of the equivalence of these three conditions, especially (2) $\Rightarrow$ (3), in the so-called {\em adhesive} ({\bf \ref{ch-rigid}}.\ref{dfn-universallyadhesiverigidspaces}) situation, to the next volume, in which we will show the {\em Enlargement Theorem} by expanding L\"utkebohmert's\index{Lutkebohmert@L\"utkebohmert, W} technique.

\medskip\noindent
{\bf 12.\ Contents of the future volumes.}
The continuation of our project will be given in the future volumes.
The next volume will contain the following chapters:
\begin{itemize}
\item Chapter 3: Formal flattening theorem --- this chapter will also contain several applications of the formal flattening theorem, such as Gerritzen-Grauert theorem;
\item Chapter 4: Enlargement theorem --- this chapter will contain the proof of the equivalence of the three `definitions' of properness;
\item Chapter 5: Equivalence theorem and analytification of algebraic spaces --- this chapter will give the proof of Equivalence Theorem stated above and the definition of the GAGA functor for algebraic spaces.
\end{itemize}

\medskip\noindent
{\bf 13.\ General conventions.}
Chapter numbers are denoted by bold-face roman numerals, while the numerals for section and subsection numbers are written arabic; subsubsections are numbered by alphabets in parentheses; for example, `{\bf I}, \S3.2.(b)' refers to the second subsubsection of the second subsection in \S3 of Chapter {\bf I}.
Cross-references will be given by sequences of numerals, like {\bf I}.3.2.1, which specify the places of the statements in the text.
The chapter numbers are omitted when referring to places in the same chapter.

Almost all sections are equipped with some exercises in the end, which are selected in order to help the reader's understanding of the content.
We insert hints for some of the exercises in the end of this volume.

Let us list some mathematical conventions:
\begin{itemize}
\item We fix once for all a Grothendieck universe $\mathsf{U}$ (\cite[Expos\'e I, 0]{SGA4-1}); cf.\ {\bf \ref{ch-pre}}, \S\ref{subsub-existuniverse}.
\item By a Grothendieck topology (or simply by a topology) on a category $\mathscr{C}$ we always mean a functor $J\colon x\mapsto J(x)$, assigning to each $x\in\obj(\mathscr{C})$ a collection of sieves, as in \cite[III, \S2, Def.\ 1]{MacMoe}.
In many places, however, Grothendieck topologies are introduced by means of basis (covering families) as in \cite[III, \S2, Def.\ 2]{MacMoe} ({\em pr\'etopologie} by the terminology in \cite[Expos\'e II, (1.3)]{SGA4-1}); in this situation, we consider, without explicit mentioning, the Grothendieck topology generated by the basis.
\item A site will always mean a {\em $\mathsf{U}$-site} (cf.\ \cite[Expos\'e II, (3.0.2)]{SGA4-1}), that is, a pair $(\mathscr{C},J)$ consisting of a $\mathsf{U}$-category $\mathscr{C}$ (\cite[Expos\'e I, Def.\ 1.1]{SGA4-1}) and a Grothendieck topology on $\mathscr{C}$.
\item All compact topological spaces are assumed to be Hausdorff; that is, we adopt the Bourbaki convention: 
$$
\textrm{quasi-compact}\ +\ \textrm{Hausdorff}\ =\ \textrm{compact}.
$$
However, we sometimes use such a phrase `compact Hausdorff' just for emphasis. 
Other conventions, in which we do {\em not} follow Bourbaki, are:
\begin{itemize}
\item locally compact spaces are only assumed to be {\em locally} Hausdorff\footnote{Note that, in \cite[Chap.\ I, \S9.7, Def.\ 4]{Bourb4}, locally compact spaces are assumed to be Hausdorff.}; more precisely, a topological space $X$ is said to be locally compact if every point of $X$ has a compact neighborhood contained in a Hausdorff neighborhood; 
\item paracompact spaces are {\em not} assumed to be Hausdorff; see {\bf \ref{ch-pre}}, \S\ref{subsub-paracompactspaces}. 
\end{itemize}
\item Whenever we say $A$ is a ring, we always mean, unless otherwise clearly stated, that $A$ is a commutative ring having the multiplicative unit $1=1_A$. We also assume that any ring homomorphism $f\colon A\rightarrow B$ is unital, that is, maps $1_A$ to $1_B$.
Moreover, 
\begin{itemize}
\item for a ring $A$ we denote by $\Frac(A)$ the total ring of fractions of $A;$
\item for a ring $A$ the Krull dimension of $A$ is denoted by $\dim(A);$
\item when $A$ is a local ring, its unique maximal ideal is denoted by $\m_A$.
\end{itemize}
\item Let $A$ be a ring, and $I\subseteq A$ an ideal. When we say $A$ is $I$-adically complete or complete with respect to the $I$-adic topology, we always mean, unless otherwise clearly stated, that $A$ is {\em Hausdorff complete} with respect to the $I$-adic topology.
\item By an exact functor between derived categories (of any sort) we always mean an exact functor of triangulated categories that preserves the canonical $t$-structures (hence also the canonical cohomology functors), which are clearly specified by the context.
\item We will often write, by abuse of notation, the equality symbol `$=$' for `isomorphic by a canonical morphism'.
\end{itemize}

\medskip\noindent
{\bf Acknowledgements.} 
There are so many people and institutions whom we owe this work very much.
We first of all would like to thank Ofer Gabber for many valuable and inspiring discussions and results, partly obtained by collaboration, which improved very much our original results and initial blue-prints of this project.
We thank Ahmed Abbes, Yves Andr\'e, Matt Baker, Ken-ichi Bannai, Francesco Baldassarri, Vladimir Berkovich, Siegfried Bosch, Bruno Chiarellotto, Gunther Cornelissen, Michel Gros, Quentin Guignard, Bernard Le Stum, Yoichi Mieda, Johannes Nicaise, Frans Oort, Michel Raynaud, Takeshi Saito, Peter Ullrich, and Gerard van der Geer for valuable discussions and friendly encouragement.
We are grateful to Institut des Hautes \'{E}tudes Scientifiques, Research Institute for Mathematical Sciences, Centre de Recerca Matem\`{a}tica, \'Ecole Normale Sup\'erieure, Universit\'{e} de Rennes 1, Universiteit Utrecht, and Nagoya University for nice hospitality.
Finally, let us mention that this work was supported by JSPS KAKENHI Grant Numbers 13440004, 17340002, 21340004, 17740014, 20540015, 23540015, 26400050, and 15H03607.
\setcounter{chapter}{-1}
\renewcommand{\thesection}{\arabic{section}}
\chapter{Preliminaries}\label{ch-pre}
This chapter collects basic notions and results from various fields, which are prepared not only for the rest of this volume, but also for later volumes. 
In spite of its preliminary nature, this chapter contains new notions, results and techniques.
From \S\ref{sec-language} to \S\ref{sec-topologicalringsmodules} are for the background concepts and results in set theory, category theory, general topology, homological algebras, ringed spaces, schemes and algebraic spaces, valuation rings, and topological rings and modules.
One of the new items that come in these sections is the notion of {\em valuative spaces} discussed in \S\ref{sub-valuativespace}.
It will be shown, later in Chapter {\bf \ref{ch-rigid}}, that valuative spaces give a nice topological visualization (cf.\ Introduction) of rigid spaces; viz., the Zariski-Riemann spaces associated to rigid spaces are all valuative.
In this sense, valuative spaces should be regarded as an abstraction, in non-archimedean geometry, of spectral spaces, whereby leading to a spectral geometry that arises from non-archimedean geometries.
We will see that valuative spaces have many rich topological structures, such as separated quotients, overconvergent subsets, tube subsets, etc.

In the sections \S\ref{sec-pairs} and \S\ref{sec-aadicallycompval}, we give a general treatment of adically topologized rings and modules.
The main object of these sections is commutative rings equipped with an adic topology by a finitely generated ideal. 
At first, we develop the theory of these objects in the most general setting, and later, we will consider various kinds of finiteness conditions imposed {\em outside} the ideal of definition.
This `finiteness conditions outside $I$' will form the central part of our discussion.
For example, $I$-adically complete and {\em Noetherian-outside-$I$} rings form a nice class of adic rings that enjoy many of the pleasant properties known to be satisfied by Noetherian topological rings (due to a theorem by  Gabber (Theorem \ref{thm-gabberIHES2008}\index{Gabber, O.})), such as `preservation of adicness', a property similar to Artin-Rees property.
Of particular importance among these topological rings are what we call {\em topologically universally pseudo-adhesive} and {\em topologically universally adhesive} (abbr.\ {\em t.u.\ adhesive}) rings.
All these new notions and techniques will provide the basis for the next chapter, in which we develop the general theory of formal schemes.

In \S\ref{sec-aadicallycompval} we focus on topologically of finite type algebras over $V$, a valuation ring that is $a$-adically complete by a non-zero $a\in\m_V$.
It will be shown, due to another theorem by Gabber\index{Gabber, O.}, that such an algebra is always t.u.\ adhesive (Corollary \ref{cor-convadh}).
Notice that, here, we do not assume that the valuation ring $V$ is of height one.
We will moreover show several ring-theoretically important results, such as {\em Noether normalization theorem} in height one situation (Theorem \ref{thm-noethernormalizationtype(V)}).

In the end of this chapter, we put appendices for further concepts and techniques including Huber's f-adic rings and basics on derived categories. 

\section{Basic Languages}\label{sec-language}
This section gives a short glossary of set theory and category theory.
Like in the modern approaches to algebraic geometry, we postulate Grothendieck's axiom {\bf (UA)} in \cite[Expos\'e I, 0]{SGA4-1} on existence of Grothendieck universe\index{Grothendieck!universe@--- universe} and fix one universe once for all. 
Some of the related technical notions, such as $\mathsf{U}$-small sets and $\mathsf{U}$-categories, are briefly reviewed in \S\ref{sub-sets} and \S\ref{sub-categories}.
In \S\ref{sub-limdef} we discuss limits and colimits, especially, filtered (cofiltered) and essentially small limits (colimits).
The final subsection gives an overview of general categorical frameworks for several stabilities of properties of arrows, which are mainly taken from \cite[I]{Knu} and further developed. 

\subsection{Sets and ordered sets}\label{sub-sets}
\subsubsection{Sets}\label{subsub-existuniverse}
In this book, we entirely work in the ZFC-set theory (cf.\ \cite[Chap.\ 1, \S3]{Drake}) with the language of {\em classes}\index{class} (cf.\ \cite[Chap.\ 1, \S5.3]{Drake}), assuming Grothendieck's axiom {\bf (UA)} for existence of {\em Grothendieck universes}\index{Grothendieck, A.}\index{Grothendieck!universe@--- universe}\index{universe|see{Grothendieck universe}} (\cite[Expos\'e I, 0]{SGA4-1}).
We fix once for all a Grothendieck universe $\mathsf{U}$ containing at least one infinite ordinal.
As in \cite[Expos\'e I, 1.0]{SGA4-1}, a set $x$ is said to be {\em $\mathsf{U}$-small}\index{small!Usmall@$\mathsf{U}$-{---}}\index{set!Usmall set@$\mathsf{U}$-small ---}, if it is isomorphic to a member of $\mathsf{U}$.

\subsubsection{Ordered sets and order types}\label{subsub-orderings}
Recall that an {\em ordering}\index{ordering@order(ing)} on a (not necessarily $\mathsf{U}$-small) set $X$ is a relation $\leq$ on $X$ satisfying the following conditions:
\begin{itemize}
\item[{\bf (O1)}] $x\leq x$ for any $x\in X$;
\item[{\bf (O2)}] for $x,y\in X$, $x\leq y$ and $y\leq x$ imply $x=y$;
\item[{\bf (O3)}] for $x,y,z\in X$, $x\leq y$ and $y\leq z$ imply $x\leq z$.
\end{itemize}
As usual, we write $x<y$ if $x\leq y$ and $x\neq y$.

A {\em $($partially$)$ ordered set}\index{ordered!set@--- set}\index{set!ordered set@ordered ---}  {\em $($poset$)$}\index{partially ordered set|see{ordered set}} is a pair $(X,\leq)$ consisting of a set and an ordering.
It is clear that an ordered set $(X,\leq)$ is $\mathsf{U}$-small (resp.\ a member of $\mathsf{U}$) if and only if so is the set $X$.

For an ordered set $X=(X,\leq)$, we denote by $X^{\opp}=(X,\leq^{\opp})$ the ordered set having the same underlying set $X$ with all the inequalities reversed, that is, for $x,y\in X$, $x\leq^{\opp}y$ if and only if $y\leq x$.

Let $(X,\leq)$ and $(Y,\leq)$ be ordered sets.
A map $f\colon X\rightarrow Y$ is said to be {\em ordered}\index{ordered!map@--- map} if, for any $x,y\in X$ with $x\leq y$, we have $f(x)\leq f(y)$.

An ordered set $(X,\leq)$ is said to be {\em totally ordered}\index{ordered!totally@totally ---} if, for any $x,y\in X$, either one of $x<y$, $x=y$, or $x>y$ holds.
For example, any ordinal (cf.\ \cite[Chap.\ 2, \S2]{Drake}) is a totally ordered set by the membership relation $\in$.

Consider the set ${\rm \mathsf{U}\textrm{-}Ord}$ of all isomorphism classes of totally ordered sets belonging to $\mathsf{U}$.
For a totally ordered set $(X,\leq)$ in $\mathsf{U}$, the unique element $\rho\in{\rm \mathsf{U}\textrm{-}Ord}$ such that $(X,\leq)\in\rho$ is called the {\em order type}\index{order type} of $(X,\leq)$.
Any finite ordinal $n=\{0,1,\ldots,n-1\}$ defines a unique order type, denoted again by $n$.
Order types of this form are said to be {\em finite}\index{order type!finite@finite ---}.

\subsection{Categories}\label{sub-categories}
\subsubsection{Conventions}\label{subsub-categoryconv}
In this book, categories\index{category} are always considered within set theory; our standard reference to category theory is \cite{Mac}.
For a category $\mathscr{C}$ we denote by $\obj(\mathscr{C})$ the class of objects of $\mathscr{C}$, and for each pair $(x,y)$ of objects of $\mathscr{C}$, we denote by $\Hom_{\mathscr{C}}(x,y)$ the class of arrows from $x$ to $y$.
A a category $\mathscr{C}$ is called a {\em $\mathsf{U}$-category}\index{category!Ucategory@$\mathsf{U}$-{---}} if $\Hom_{\mathscr{C}}(x,y)$ is $\mathsf{U}$-small for any $x,y\in\obj(\mathscr{C})$ (\cite[Expos\'e I, Def.\ 1.1]{SGA4-1}).
Almost all categories in this book are $\mathsf{U}$-categories; moreover, they most of the time satisfy the following conditions (\cite[Expos\'e I, 1.1.2]{SGA4-1}):
\begin{itemize}
\item[{\bf (C1)}] the class of objects $\obj(\mathscr{C})$ is a {\em subset} of $\mathsf{U}$;
\item[{\bf (C2)}] for any $x,y\in\obj(\mathscr{C})$ the set $\Hom_{\mathscr{C}}(x,y)$ is a {\em member} of $\mathsf{U}$.
\end{itemize}

Let $\mathscr{C}$ be a $\mathsf{U}$-category.
For each pair $(x,y)$ of objects, we denote by $\Isom_{\mathscr{C}}(x,y)$ the subset of $\Hom_{\mathscr{C}}(x,y)$ consisting of all isomorphisms.
Also, for an object $x$, we write $\End_{\mathscr{C}}(x)=\Hom_{\mathscr{C}}(x,x)$ and $\Aut_{\mathscr{C}}(x)=\Isom_{\mathscr{C}}(x,x)$.

For a category $\mathscr{C}$, we denote by $\mathscr{C}^{\opp}$ the {\em opposite category}\index{category!opposite@opposite ---} of $\mathscr{C}$ (\cite[Chap.\ II, \S2]{Mac}), that is, the category such that 
\begin{itemize}
\item $\obj(\mathscr{C}^{\opp})=\obj(\mathscr{C})$;
\item $\Hom_{\mathscr{C}^{\opp}}(x,y)=\Hom_{\mathscr{C}}(y,x)$ for any $x,y\in\obj(\mathscr{C})$.
\end{itemize}

By a functor\index{functor} $\mathscr{C}\rightarrow\mathscr{D}$ we always mean a {\em covariant} functor, unless otherwise clearly stated.
Contravariant functors will be written as a covariant functor from the opposite category of the domain category.

\subsubsection{Frequently used categories}\label{subsub-categorynot}
The following categories are frequently used in this book:
\begin{itemize}
\item $\Sets=$ the category of all sets in $\mathsf{U}$;
\item $\Top=$ the category of all topological spaces in $\mathsf{U}$;
\item $\Ab=$ the category of all abelian groups in $\mathsf{U}$.
\end{itemize}
For a ring (commutative with unit) $A$,
\begin{itemize}
\item $\Mod_A=$ the category of all $A$-modules in $\mathsf{U}$;
\item $\Alg_A=$ the category of all $A$-algebras in $\mathsf{U}$.
\end{itemize}
These are $\mathsf{U}$-categories satisfying the conditions {\bf (C1)} and {\bf (C2)} in \ref{subsub-categoryconv}.

\subsubsection{Functor category}\label{subsub-functorcat}\index{functor!category@--- category}\index{category!functor@functor ---}
Let $\mathscr{C}$ and $\mathscr{D}$ be categories.
We denote by $\mathscr{C}^{\mathscr{D}}$ the category of functors from $\mathscr{D}$ to $\mathscr{C}$ and natural transformations. 
Notice that (cf.\ \cite[Expos\'e I, 1.1.1]{SGA4-1}):
\begin{itemize}
\item[(1)] if both $\mathscr{C}$ and $\mathscr{D}$ are members of $\mathsf{U}$ (resp.\ $\mathsf{U}$-small), then so is $\mathscr{C}^{\mathscr{D}}$;
\item[(2)] if $\mathscr{D}$ is $\mathsf{U}$-small and $\mathscr{C}$ is a $\mathsf{U}$-category, then $\mathscr{C}^{\mathscr{D}}$ is a $\mathsf{U}$-category.
\end{itemize}

\subsubsection{Groupoids and discrete categories}\label{subsub-groupoinds}
A category $\mathscr{C}$ is said to be a {\em groupoid}\index{groupoid} if all arrows are isomorphisms.
A category $\mathscr{C}$ is said to be {\em discrete}\index{category!discrete@discrete ---} if the set $\Hom_{\mathscr{C}}(x,y)$ is empty unless $x=y$ and, for each object $x$, we have $\End_{\mathscr{C}}(x)=\{\id_x\}$.
Notice that a discrete category $\mathscr{C}$ is completely determined by its class of objects.

\subsubsection{Category associated to an ordered set}\label{subsub-catposet}\index{category!associated to an ordered set@--- associated to an ordered set}
To an ordered set\index{ordered!set@--- set}\index{set!ordered set@ordered ---} $X=(X,\leq)$ is canonically associated a category, denoted again by $X$, whose objects are the elements of $X$, and for $x,y\in X$,
$$
\Hom_X(x,y)=
\begin{cases}
\{(x,y)\}&\textrm{if $x\leq y$,}\\
\emptyset&\textrm{otherwise.}
\end{cases}
$$
If the set $X$ is a member of $\mathsf{U}$ (resp.\ $\mathsf{U}$-small), then so is the associated category.
Notice that the category associated to the inverse ordered set $X^{\opp}$ (\S\ref{subsub-orderings}) is the opposite category of the category associated to $X$ (\S\ref{subsub-categoryconv}).

\subsection{Limits}\label{sub-limdef}
\subsubsection{Definition and universal property}\label{subsub-limdefuniv}
Let $\mathscr{C}$ and $\mathscr{D}$ be categories.
The {\em diagonal functor}\index{functor!diagonal functor@diagonal ---}
$$
\Delta\colon\mathscr{C}\longrightarrow\mathscr{C}^{\mathscr{D}}
$$
is the functor defined as follows: $\Delta(x)\colon\mathscr{D}\rightarrow\mathscr{C}$ for $x\in\obj(\mathscr{C})$ is the constant functor given by $\Delta(x)(y)=x$ for any object $y$ of $\mathscr{D}$ and $\Delta(x)(f)=\id_x$ for any arrow $f$ of $\mathscr{D}$.
We denote the right (resp.\ left) adjoint to $\Delta$, if it exists, by $\varprojlim$ (resp.\ $\varinjlim$).
For a functor $F\colon\mathscr{D}\rightarrow\mathscr{C}$, the object $\varprojlim F$ (resp.\ $\varinjlim F$) is called the {\em limit}\index{limit} (resp.\ {\em colimit}\index{limit!colimit@co---}\index{colimit}) of $F$.

To describe the mapping universality of the limit $\varprojlim F$, consider a morphism $\Delta(x)\rightarrow F$ of functors, which amounts to the same as a collection of arrows $x\rightarrow F(y)$ in $\mathscr{C}$ for all $y\in\obj(\mathscr{D})$ such that for any arrow $f\colon y\rightarrow z$ in $\mathscr{D}$ the resulting triangle
$$
\xymatrix@R-4ex@C-1ex{&F(y)\ar[dd]^{F(f)}\\ x\ar[ur]\ar[dr]\\ &F(z)}
$$
is commutative.
Replacing $x$ by $\varprojlim F$, and considering the adjunction morphism $\Delta(\varprojlim F)\rightarrow F$, one gets the compatible collection of arrows $\varprojlim F\rightarrow F(y)$ ($y\in\obj(\mathscr{D})$).
The limit $\varprojlim F$ is then characterized up to isomorphism by the following universal mapping property: Whenever a collection of arrows $x\rightarrow F(y)$ as above is given, there exists a unique arrow $x\rightarrow\varprojlim F$ such that the diagram 
$$
\xymatrix@R-4ex@C-1ex{&\varprojlim F\ar[dd]\\ x\ar[ur]\ar[dr]\\ &F(y)}
$$
is commutative for any $y\in\obj(\mathscr{D})$.

The mapping universality of the colimit $\varinjlim F$ can be described similarly; the details are left to the reader.

\subsubsection{Limits over ordered sets}\label{subsub-limdefsets}
We will most frequently deal with limits and colimits with the index category being an ordered set\index{ordered!set@--- set}\index{set!ordered set@ordered ---} (\S\ref{subsub-catposet}).
If $I$ is an ordered set\index{ordered!set@--- set}\index{set!ordered set@ordered ---}, then the functor $F\colon I\rightarrow\mathscr{C}$ as above amounts to the same as what is usually called the {\em inductive system}\index{system!inductive system@inductive ---}\index{system!direct system@direct ---|see{inductive system}} (synonym: {\em direct system}) $\{X_i,f_{ij}\}$ of objects and arrows in $\mathscr{C}$.
Similarly, a functor of the form $G\colon I^{\opp}\rightarrow\mathscr{C}$ corresponds to what is called a {\em projective system}\index{system!projective system@projective ---}\index{system!inverse system@inverse ---|see{projective system}} (synonym: {\em inverse system}).
The corresponding limits, inductive and projective\index{limit!inductive limit@inductive ---}\index{limit!projective limit@projective ---}, are respectively written in the usual way as
$$
\varinjlim_{i\in I}X_i\qquad\textrm{and}\qquad\varprojlim_{i\in I}X_i.
$$

\subsubsection{Final and cofinal functors}\label{subsub-finalcofinal}
A category $I$ is said to be {\em filtered}\index{category!filtered category@filtered ---} if it is non-empty and satisfies the following conditions:
\begin{itemize}
\item[{\rm (a)}] for any $x,x'\in\obj(I)$, there exist arrows $x\rightarrow y$ and $x'\rightarrow y$;
\item[{\rm (b)}] given two arrows $f,f'\colon x\rightarrow y$, there exists $g\colon y\rightarrow z$ such that $g\circ f=g\circ f'$.
\end{itemize}

An ordered set\index{ordered!set@--- set}\index{set!ordered set@ordered ---} $I$ is said to be {\em directed}\index{set!directed set@directed ---}\index{directed set} if any finite subset has an upper bound. 
If $I$ is a directed set, then it is, viewed as a category as in \S\ref{subsub-catposet}, a filtered category.

A functor $L\colon J\rightarrow I$ between filtered categories is said to be {\em final}\index{functor!final functor@final ---} if
\begin{itemize}
\item[{\bf (F)}] for any $i\in\obj(I)$, there exists $j\in\obj(J)$ and an arrow $i\rightarrow L(j)$. 
\end{itemize}

\begin{prop}[{\cite[Chap.\ IX, \S3, Theorem 1]{Mac}}]\label{prop-final}
Let $L\colon J\rightarrow I$ be a final functor between filtered categories, and $F\colon I\rightarrow\mathscr{C}$ a functor. 
If $\varinjlim F\circ L$ exists, then so does $\varinjlim F$, and the canonical morphism 
$$
\varinjlim (F\circ L)\longrightarrow\varinjlim F
$$
is an isomorphism.\hfill$\square$
\end{prop}

Thus, when taking the colimit, one can replace the index category $I$ by a category $J$ that admits a final functor $J\rightarrow I$. 
If such a category $J$ can be chosen as a directed set, the index category $I$ is said to be {\em essentially small}\index{small!essentially small@essentially ---}.
In this case, the limit `along category' can be replaced by a limit `along set', which is, needless to say, easier to handle. 

Similarly, by duality, one has the notions of {\em cofiltered} categories\index{category!cofiltered category@cofiltered ---} and {\em cofinal} functors\index{functor!cofinal functor@cofinal ---} $L\colon J\rightarrow I$.
For example, if $I$ is a directed set, the category $I^{\opp}$ is cofiltered.
A cofiltered category $I$ is said to be {\em essentially small}\index{small!essentially small@essentially ---} if it admits a cofinal functor $J^{\opp}\rightarrow I$, where $J$ is a directed set.
The dual statement of \ref{prop-final} holds alike:  
\begin{prop}\label{prop-cofinal}
Let $L\colon J\rightarrow I$ be a cofinal functor between cofiltered categories, and $F\colon I\rightarrow\mathscr{C}$ a functor. 
If $\varprojlim F\circ L$ exists, then so does $\varprojlim F$, and the canonical morphism 
$$
\varprojlim F\longrightarrow\varprojlim (F\circ L)
$$
is an isomorphism.\hfill$\square$
\end{prop}

\subsection{Several stabilities for properties of arrows}\label{sub-propertyarrow}
\subsubsection{Base-change stability}\label{subsub-basechangestable}
Let $\mathscr{C}$ be a category.
We consider a subcategory $\mathscr{D}$ of $\mathscr{C}$ containing all isomorphisms of $\mathscr{C}$ (hence, in particular, $\obj(\mathscr{D})=\obj(\mathscr{C})$).
Typically, such a subcategory comes about in the following way: Let $P$ be a property of arrows in the category $\mathscr{C}$ such that
\begin{itemize}
\item[{\bf (I)}] any isomorphism satisfies $P$; 
\item[{\bf (C)}] if $a\colon x\rightarrow y$ and $b\colon y\rightarrow z$ satisfy $P$, then the composition $b\circ a$ satisfies $P$.
\end{itemize}
Then the subcategory consisting of all arrows satisfying $P$, denoted by $\mathscr{D}_P$, is a subcategory of the type as above.
\begin{prop}\label{prop-basechangestable}
Let $\mathscr{D}$ and $\mathscr{E}$ be subcategories of $\mathscr{C}$ containing all isomorphisms.
We assume that the following condition holds$:$
\begin{itemize}
\item if $a\colon x\rightarrow y$ belongs to $\mathscr{E}$, then for any $b\colon z\rightarrow y$, the fiber product $x\times_yz$ is representable in $\mathscr{C}$, and the arrow $x\times_yz\rightarrow z$ belongs to $\mathscr{E}$.
\end{itemize}
Consider the following conditions$:$
\begin{itemize}
\item[{\boldmath $(\mathbf{B}_1)$}] Suppose $\xymatrix@R-3ex@C-4ex{x\ar[rr]^a\ar[dr]&&x'\ar[dl]\\ &z}$ and $\xymatrix@R-3ex@C-4ex{y\ar[rr]^b\ar[dr]&&y'\ar[dl]\\ &z}$ are commutative diagrams in $\mathscr{C}$ such that {\rm (i)} the arrows $a$ and $b$ are in $\mathscr{D}$, and {\rm (ii)} either the arrows $x\rightarrow z$ and $x'\rightarrow z$ or the arrows $x\rightarrow z$ and $y'\rightarrow z$ belong to $\mathscr{E}$. Then the induced arrow
$$
a\times_zb\colon x\times_zy\longrightarrow x'\times_zy'
$$
belongs to $\mathscr{D}$.
\item[{\boldmath $(\mathbf{B}_2)$}] Suppose a commutative diagram $\xymatrix@R-3ex@C-4ex{x\ar[rr]^a\ar[dr]&&y\ar[dl]\\ &z}$ in $\mathscr{C}$ and an arrow $z'\rightarrow z$ in $\mathscr{C}$ are given such that the arrow $a$ belongs to $\mathscr{D}$.
Suppose, moreover, that either one of the following holds$:$  {\rm (i)} the arrows $x\rightarrow z$ and $y\rightarrow z$ belong to $\mathscr{E}$, or {\rm (ii)} the arrow $z'\rightarrow z$ belongs to $\mathscr{E}$.
Then the induced arrow
$$
a_{z'}\colon x\times_zz'\longrightarrow y\times_zz'
$$
belongs to $\mathscr{D}$.
\item[{\boldmath $(\mathbf{B}_3)$}] Suppose a diagram $x\stackrel{a}{\rightarrow}y\leftarrow y'$
in $\mathscr{C}$ is given such that {\rm (i)} the arrow $a$ belongs to $\mathscr{D}$, and {\rm (ii)} either the arrow $x\rightarrow y$ or the arrow $y'\rightarrow y$ belongs to $\mathscr{E}$. 
Then the induced arrow
$$
a_{y'}\colon x\times_yy'\longrightarrow y'
$$
belongs to $\mathscr{D}$.
\end{itemize}
Then we have the implications
$$
\textrm{{\boldmath $(\mathbf{B}_1)$}}\quad\Longleftrightarrow\quad\textrm{{\boldmath $(\mathbf{B}_2)$}}\quad\Longrightarrow\quad\textrm{{\boldmath $(\mathbf{B}_3)$}}.
$$
If we assume, moreover, the following condition, then the implication {\boldmath $(\mathbf{B}_3)$} $\Rightarrow$ {\boldmath $(\mathbf{B}_2)$} also holds$:$
\begin{itemize}
\item for arrows $a\colon x\rightarrow y$ and $b\colon y\rightarrow z$ in $\mathscr{C}$, if $b\circ a$ and $b$ belong to $\mathscr{E}$, then $a$ belongs to $\mathscr{E};$
\end{itemize}
\end{prop}

\begin{proof}
First let us show {\boldmath $(\mathbf{B}_1)$} $\Rightarrow$ {\boldmath $(\mathbf{B}_2)$}.
In the situation as in {\boldmath $(\mathbf{B}_2)$}, if the arrows $x\rightarrow z$ and $y\rightarrow z$ belong to $\mathscr{E}$, apply {\boldmath $(\mathbf{B}_1)$} with $y=y'$ replaced by $z'$ and with $x'$ replaced by $y$.
If the arrow $z'\rightarrow z$ belongs to $\mathscr{E}$, apply {\boldmath $(\mathbf{B}_1)$} with $x=x'$ replaced by $z'$ and with $y$ and $y'$ by $x$ and $y$, respectively.
Conversely, to show {\boldmath $(\mathbf{B}_2)$} $\Rightarrow$ {\boldmath $(\mathbf{B}_1)$}, we use the fact that the arrow $a\times_zb$ coincides with the composites
$$
x\times_zy\longrightarrow x'\times_zy\longrightarrow x'\times_zy'\qquad\textrm{and}\qquad x\times_zy\longrightarrow x\times_zy'\longrightarrow x'\times_zy'
$$
(cf.\ \cite[($\mathbf{0}$, 1.3.9)]{EGAInew}).
In either case where the arrows $x\rightarrow z$ and $x'\rightarrow z$ or the arrows $x\rightarrow z$ and $y'\rightarrow z$ belong to $\mathscr{E}$, one can apply {\boldmath $(\mathbf{B}_2)$} to show that the arrow $a\times_zb$ belongs to $\mathscr{D}$.

The implication {\boldmath $(\mathbf{B}_2)$} $\Rightarrow$ {\boldmath $(\mathbf{B}_3)$} follows easily, since {\boldmath $(\mathbf{B}_3)$} is the special case of {\boldmath $(\mathbf{B}_2)$} with $z=y$.
To show {\boldmath $(\mathbf{B}_3)$} $\Rightarrow$ {\boldmath $(\mathbf{B}_2)$}, we use the fact that the arrow $a_{z'}$ is isomorphic to 
$$
x\times_y(y\times_zz')\longrightarrow y\times_zz'
$$
(cf.\ \cite[($\mathbf{0}$, 1.3.4)]{EGAInew}).
If the arrows $x\rightarrow z$ and $y\rightarrow z$ belong to $\mathscr{E}$, then one can apply {\boldmath $(\mathbf{B}_3)$} with $y'$ replaced by $y\times_zz'$, since $x\rightarrow y$ belongs to $\mathscr{E}$.
If the arrow $z'\rightarrow z$ belongs to $\mathscr{E}$, one can apply {\boldmath $(\mathbf{B}_3)$} with $y'$ replaced by $y\times_zz'$, since $y\times_zz'\rightarrow y$ belongs to $\mathscr{E}$.
\end{proof}

\begin{dfn}\label{dfn-basechangestability}{\rm 
Suppose $\mathscr{C}$ has all fiber products.
Then the subcategory $\mathscr{D}$, or the property $P$ with $\mathscr{D}=\mathscr{D}_P$, is said to be {\em base-change stable}\index{stable!base-change stable@base-change ---} in $\mathscr{C}$ if the conditions in \ref{prop-basechangestable} with $\mathscr{E}=\mathscr{C}$ are satisfied.}
\end{dfn}

\subsubsection{Topology associated to base-change stable subcategory}\label{subsub-topologygenconst}
We insert here a brief review of the Grothendieck topology\index{topology!Grothendieck topology@Grothendieck ---}\index{Grothendieck!Grothendieck topology@--- topology}\footnote{See General conventions in the introduction for our convention for Grothendieck topologies.} associated to a base-change stable\index{stable!base-change stable@base-change ---} subcategory (cf.\ \cite[I.1]{Knu}).
Let $\mathscr{C}$ be a category with all fiber products.
\begin{dfn}\label{dfn-UEEF}{\rm 
{\rm (1)} A family $\{u_{\alpha}\rightarrow u\}_{\alpha\in L}$ of arrows in $\mathscr{C}$ is said to be {\em universally effectively epimorphic}\index{universally effectively epimorphic} if, for any arrow $w\rightarrow u$ in $\mathscr{C}$ and for any object $v\in\obj(\mathscr{C})$, the induced sequence of maps
$$
\xymatrix@C-3ex{\Hom_{\mathscr{C}}(w,v)\ar[r]&\prod_{\alpha\in L}\Hom_{\mathscr{C}}(w\times_uu_{\alpha},v)\ar@<.5ex>[r]\ar@<-.5ex>[r]&\prod_{\alpha,\beta\in L}\Hom_{\mathscr{C}}(w\times_uu_{\alpha\beta},v),}\leqno{(\ast)}
$$
where $u_{\alpha\beta}=u_{\alpha}\times_uu_{\beta}$, is exact.

{\rm (2)} An arrow $a\colon v\rightarrow u$ is called a {\em universally effectively epimorphism}\index{universally effectively epimorphic!universally effectively epimorphism@--- --- epimorphism} if $\{a\}$ is a universally effectively epimorphic family.}
\end{dfn}

\begin{dfn}\label{dfn-topologygenconst}{\rm 
Let $\mathscr{D}$ be a base-change stable subcategory of $\mathscr{C}$ (\ref{dfn-basechangestability}).
The {\em topology associated to $\mathscr{D}$}\index{topology!topology associated to a subcategory@--- associated to a subcategory} (or to $P$ when $\mathscr{D}=\mathscr{D}_P$) is the Grothendieck topology on $\mathscr{C}$, denoted by $J_{\mathscr{D}}$ (or by $J_P$), such that the covering families are given by all universally effectively epimorphic families in $\mathscr{C}$ consisting of arrows in $\mathscr{D}$.}
\end{dfn}

\begin{prop}\label{prop-topologygenconst}
The topology $J=J_{\mathscr{D}}$ satisfies the following condition$:$
\begin{itemize}
\item[{\boldmath $(\mathbf{A}_0)$}] any representable presheaf on the site $\mathscr{C}_J=(\mathscr{C},J)$ is a sheaf.
\end{itemize}
In other words, the topology $J$ is coarser than the canonical topology on $\mathscr{C}$ {\rm (\cite[Expos\'e II, (2.5)]{SGA4-1}\cite[p.\ 126]{MacMoe})}.\hfill$\square$
\end{prop}

\noindent
{\bf Convention.} {\sl Throughout this book, whenever a site $\mathscr{C}_J=(\mathscr{C},J)$ satisfies the condition {\boldmath $(\mathbf{A}_0)$}, the sheaf on $\mathscr{C}_J$ associated to an object $x\in\obj(\mathscr{C})$ is simply denoted by the same symbol $x$.}

\medskip
There are several other useful conditions for the topology $J=J_{\mathscr{D}}$, some of which come from properties of the subcategory $\mathscr{D}$.
It will be convenient for us to put some of them into axioms.
\begin{dfn}\label{dfn-disjointunion}{\rm 
Let $\mathscr{C}$ be a category with all fiber products.

{\rm (1)} An object $\emptyset\in\obj(\mathscr{C})$ is said to be {\em strictly initial}\index{strictly initial} if it is an initial object and for any $x\in\obj(\mathscr{C})$ the set $\Hom_{\mathscr{C}}(x,\emptyset)$ is empty unless $x$ is isomorphic to $\emptyset$.
(It can be shown that, if a strictly initial object exists, then any initial object is strictly initial.)

{\rm (2)} Let $\{x_{\alpha}\}_{\alpha\in L}$ be a family of objects in $\mathscr{C}$.
The coproduct $\coprod_{\alpha}x_{\alpha}$ is called the {\em disjoint sum}\index{disjoint sum} if each $x_{\alpha}\times_xx_{\beta}$ for $(\alpha,\beta)\in L\times L$ is a strictly initial object.}
\end{dfn}

Now the axioms for the base-change stable subcategory $\mathscr{D}$ that we will use in the sequel are listed as follows:
\begin{itemize}
\item[{\boldmath $(\mathbf{S}_1)$}] For a family $\{x_{\alpha}\rightarrow y\}_{\alpha\in L}$ of arrows in $\mathscr{C}$ such that the disjoint sum $x=\coprod_{\alpha\in L}x_{\alpha}$ exists, the induced arrow $x\rightarrow y$ belongs to $\mathscr{D}$ if and only if each $x_{\alpha}\rightarrow x$ belongs to $\mathscr{D}$.
\item[{\boldmath $(\mathbf{S}_2)$}] An arrow in $\mathscr{D}$ is a universally effectively epimorphism\index{universally effectively epimorphic!universally effectively epimorphism@--- --- epimorphism} in $\mathscr{C}$ if and only if it is an epimorphism in $\mathscr{C}$.
\item[{\boldmath $(\mathbf{S}_3)$}] For a commutative diagram
$\xymatrix@R-2ex@C-3.5ex{x\ar[rr]^a\ar[dr]_c&&y\ar[dl]^b\\ &z}$
in $\mathscr{C}$ such that $c$ belongs to $\mathscr{D}$,
\begin{itemize}
\item[{\boldmath $(\mathbf{S}_3{\rm (a)})$\hspace{1pt}}] if $a$ is a covering of $J_{\mathscr{D}}$, then $b$ belongs to $\mathscr{D}$;
\item[{\boldmath $(\mathbf{S}_3{\rm (b)})$}] if $b$ belongs to $\mathscr{D}$, then so does $a$.
\end{itemize}
\end{itemize}

\subsubsection{Stability and effective descent}\label{subsub-stabilitydescent}
Let $\mathscr{C}$ be a category with all fiber products, and $J$ a Grothendieck topology on $\mathscr{C}$ satisfying the condition {\boldmath $(\mathbf{A}_0)$} in \ref{prop-topologygenconst}.

\begin{dfn}\label{dfn-catequivrelstableobj}{\rm 
A subclass $\mathscr{S}$ of objects in $\mathscr{C}$ is said to be {\em stable $($under $J)$}\index{stable!stable under a topology@--- under a topology} if, for any covering family $\{u_{\alpha}\rightarrow u\}_{\alpha\in L}$, the object $u$ belongs to $\mathscr{S}$ if and only if all $u_{\alpha}$ for $\alpha\in L$ belong to $\mathscr{S}$.}
\end{dfn}

\begin{dfn}\label{dfn-catequivrelstablearrow}{\rm
Let $\mathscr{D}\subseteq\mathscr{C}$ be a subcategory of $\mathscr{C}$ containing all isomorphisms in $\mathscr{C}$.

{\rm (1)} The subcategory $\mathscr{D}$ (or $P$ when $\mathscr{D}=\mathscr{D}_P$) is said to be {\em local on the target $($under $J)$}\index{local!local on the target@--- on the target} if, for any $a\colon x\rightarrow y$ in $\mathscr{C}$ and for any covering family $\{y_{\alpha}\rightarrow y\}_{\alpha\in L}$, the arrow $a$ belongs to $\mathscr{D}$ if the base change $a_{\alpha}\colon y_{\alpha}\times_yx\rightarrow y_{\alpha}$ belongs to $\mathscr{D}$ for each $\alpha\in L$.
If it is, moreover, base-change stable, then $\mathscr{D}$ (or $P$ when $\mathscr{D}=\mathscr{D}_P$) is said to be {\em stable $($under $J)$}\index{stable!stable under a topology@--- under a topology}.

{\rm (2)} The subcategory $\mathscr{D}$ (or $P$ when $\mathscr{D}=\mathscr{D}_P$) is said to be {\em local on the domain $($under $J)$}\index{local!local on the domain@--- on the domain} if it is stable under $J$ and, for any $a\colon x\rightarrow y$ in $\mathscr{C}$ and for any covering family $\{x_{\alpha}\rightarrow x\}_{\alpha\in L}$, the arrow $a$ belongs to $\mathscr{D}$ if and only if the composition $x_{\alpha}\rightarrow x\rightarrow y$ belongs to $\mathscr{D}$ for each $\alpha\in L$.

{\rm (3)} The subcategory $\mathscr{D}$ is called an {\em effective descent class $($with respect to $J)$}\index{descent!effective descent class@effective --- class}\index{effective!effective descent class@--- descent class} if it is stable under $J$ and the following condition is satisfied: Let $\{u_{\alpha}\rightarrow u\}_{\alpha\in L}$ be a covering family, and $\mathscr{F}$ a set-valued sheaf on $\mathscr{C}_J$ together with a map $\mathscr{F}\rightarrow u$ of sheaves; if each sheaf fiber product $u_{\alpha}\times_u\mathscr{F}$ is represented by an object $w_{\alpha}\in\obj(\mathscr{C})$, and if all arrows $w_{\alpha}\rightarrow u_{\alpha}$ for $\alpha\in L$ belong to $\mathscr{D}$, then $\mathscr{F}$ is representable (consequently, if $w$ represents $\mathscr{F}$, the resulting arrow $w\rightarrow u$ belongs to $\mathscr{D}$).}
\end{dfn}

\begin{dfn}\label{dfn-localconstruction}{\rm
Let $x\in\obj(\mathscr{C})$.

$(1)$ A {\em local construction}\index{local!local construction@--- construction} $\Phi$ on $x$ consists of the following data:
\begin{itemize}
\item[{\rm (a)}] a cofinal set $\mathscr{V}=\{\{x_{\lambda\alpha}\rightarrow x\}_{\lambda\in\Lambda_{\alpha}}\}_{\alpha\in L}$ of covering families of $x$;
\item[{\rm (b)}] an arrow $\Phi_{\lambda\alpha}\rightarrow x_{\lambda\alpha}$ in $\mathscr{C}$ for each $x_{\lambda\alpha}\rightarrow x$ in a covering family in $\mathscr{V}$;
\item[{\rm (c)}] an arrow $\Phi_{\lambda\alpha}\rightarrow\Phi_{\mu\beta}$ in $\mathscr{C}$ for any arrow $x_{\lambda\alpha}\rightarrow x_{\mu\beta}$ over $x$ such that 
$$
\xymatrix@-1ex{\Phi_{\lambda\alpha}\ar[r]\ar[d]&\Phi_{\mu\beta}\ar[d]\\ x_{\lambda\alpha}\ar[r]&x_{\mu\beta}}
$$
is Cartesian in $\mathscr{C}$.
\end{itemize}

$(2)$ A local construction $\Phi$ on $x$ is said to be {\em effective}\index{effective} if there exist an arrow $y\rightarrow x$ and isomorphisms $x_{\lambda\alpha}\times_xy\cong\Phi_{\lambda\alpha}$ for all $\alpha\in L$ and $\lambda\in\Lambda_{\alpha}$ that are compatible with the arrows $\Phi_{\lambda\alpha}\rightarrow\Phi_{\mu\beta}$ as in (c) above.}
\end{dfn}

\begin{prop}[{\cite[(1.12)]{Knu}}]\label{prop-localconstruction}
Let $\Phi$ be a local construction on $x\in\obj(\mathscr{C})$, and $\mathscr{D}$ a stable subcategory of $\mathscr{C}$ that is an effective descent class\index{descent!effective descent class@effective --- class}\index{effective!effective descent class@--- descent class}.
Suppose that, for any $x_{\lambda\alpha}\rightarrow x$ in a covering family in the cofinal set $\mathscr{V}$ of $\Phi$ as above, the arrow $\Phi_{\lambda\alpha}\rightarrow x_{\lambda\alpha}$ belongs to $\mathscr{D}$.
Then $\Phi$ is effective, and the arrow $y\rightarrow x$ as in {\rm \ref{dfn-localconstruction} (2)} belongs to $\mathscr{D}$. \hfill$\square$
\end{prop}

\subsubsection{Categorical equivalence relations}\label{subsub-catequivrel}
Finally, let us include some generalities of categorical equivalence relations (cf.\ \cite[I.5]{Knu}).

Let $\mathscr{C}_J=(\mathscr{C},J)$ be a site with the underlying category $\mathscr{C}$ with a final object and all fiber products.
\begin{dfn}\label{dfn-tauequivrel}{\rm 
A {\em $J$-equivalence relation}\index{equivalence relation} in $\mathscr{C}$ is a diagram
$$
\xymatrix{R\ar@<.5ex>[r]^(.5){p_1}\ar@<-.5ex>[r]_(.5){p_2}&Y}
$$
consisting of arrows in $\mathscr{C}$ such that
\begin{itemize}
\item[(a)] for any object $Z$ of $\mathscr{C}$, the induced map
$$
\Hom_{\mathscr{C}}(Z,R)\longrightarrow\Hom_{\mathscr{C}}(Z,Y)\times\Hom_{\mathscr{C}}(Z,Y)
$$
is injective and defines an equivalence relation on the set $\Hom_{\mathscr{C}}(Z,Y)$;
\item[(b)] the maps $p_1,p_2$ are covering maps with respect to the topology $J$.
\end{itemize}}
\end{dfn}

Let $\mathscr{C}$ be a category with a final object and all fiber products, and $\mathscr{B}$ a base-change stable\index{stable!base-change stable@base-change ---} subcategory of $\mathscr{C}$ enjoying the conditions {\boldmath $(\mathbf{S}_1)$}, {\boldmath $(\mathbf{S}_3{\rm (a)})$}, and {\boldmath $(\mathbf{S}_3{\rm (b)})$} in \S\ref{subsub-topologygenconst}.
We consider the associated topology\index{topology!topology associated to a subcategory@--- associated to a subcategory} $J=J_{\mathscr{B}}$ (\ref{dfn-topologygenconst}).

As usual, a map $\mathscr{F}\rightarrow\mathscr{G}$ of sheaves on the site $\mathscr{C}_J$ is said to be {\em representable}\index{representable!representable map of sheaves@--- map of sheaves} if, for any map $z\rightarrow\mathscr{G}$ of sheaves from a representable sheaf\index{representable!representable (pre)sheaf@--- (pre)sheaf}, the sheaf fiber product $z\times_{\mathscr{G}}\mathscr{F}$ is representable.
Let $\mathscr{D}$ be a base-change stable subcategory of $\mathscr{C}$.
We say that a map $\mathscr{F}\rightarrow\mathscr{G}$ of sheaves is represented by a map in $\mathscr{D}$ if it is representable and, for any map $z\rightarrow\mathscr{G}$ of sheaves from a representable sheaf, the morphism $w\rightarrow z$ in $\mathscr{C}$ representing the base-change $z\times_{\mathscr{G}}\mathscr{F}\rightarrow z$ belongs to $\mathscr{D}$.
If, moreover, $\mathscr{D}=\mathscr{B}$ and all morphisms $w\rightarrow z$ as above are covering maps, then we say $\mathscr{F}\rightarrow\mathscr{G}$ is represented by a covering map.

\begin{prop}[{\cite[I.5.5]{Knu}}]\label{prop-fundamentalquotient0}
Let $p_1,p_2\colon R\rightarrow Y$ be a $J$-equivalence relation in the site $\mathscr{C}_J=(\mathscr{C},J)$, and $q\colon Y\rightarrow\mathscr{F}$ the cokernel of $p_1,p_2$ in the category of sheaves on $\mathscr{C}_J$.
Then for any map of sheaves $\alpha\colon T\rightarrow\mathscr{F}$ from a representable sheaf, there exists a covering map $\pi\colon U\rightarrow T$ in $\mathscr{C}_J$ sitting in a commutative diagram
$$
\xymatrix{U\times_TU\ar@<.5ex>[r]^(.65){\mathrm{pr}_1}\ar@<-.5ex>[r]_(.65){\mathrm{pr}_2}\ar[d]_{\gamma}&U\ar[r]^{\pi}\ar[d]^{\beta}&T\ar[d]^{\alpha}\\ R\ar@<.5ex>[r]^(.5){p_1}\ar@<-.5ex>[r]_(.5){p_2}&Y\ar[r]_{q}&\mathscr{F}\rlap{$;$}}
$$
here, by the commutativity of the left-hand square, we mean $\beta\circ\mathrm{pr}_1=p_1\circ\gamma$ and $\beta\circ\mathrm{pr}_2=p_2\circ\gamma$. \hfill$\square$
\end{prop}

\begin{prop}\label{prop-fundamentalquotient}
Let $p_1,p_2\colon R\rightarrow Y$ be a $J$-equivalence relation in the site $\mathscr{C}_J=(\mathscr{C},J)$, and suppose that the induced arrow $R\rightarrow Y\times Y$ belongs to an effective descent class $\mathscr{D}$ $($cf.\ {\rm \ref{dfn-catequivrelstablearrow} (3))}.
Let $q\colon Y\rightarrow\mathscr{F}$ be the cokernel of $(p_1,p_2)$ in the category of sheaves on $\mathscr{C}_J$.
Then$:$
\begin{itemize}
\item[{\rm (1)}] the map $q$ is represented by a covering map$;$
\item[{\rm (2)}] the canonical map $R\rightarrow Y\times_{\mathscr{F}}Y$ of sheaves is an isomorphism$;$
\item[{\rm (3)}] the diagonal map $\Delta_{\mathscr{F}}\colon\mathscr{F}\rightarrow\mathscr{F}\times\mathscr{F}$ is representable by an arrow in $\mathscr{D}$.
\end{itemize}
\end{prop}

\begin{proof}
The assertion (1) follows from \cite[I.5.9]{Knu}, while (2) follows from \cite[I.5.4]{Knu}.
Finally, (3) follows easily from \cite[I.5.10]{Knu}.
\end{proof}

\addcontentsline{toc}{subsection}{Exercises}
\subsection*{Exercises}
\begin{exer}\label{exer-doublecolimits}
Let $\mathscr{C}$ be a category with all essentially small\index{small!essentially small@essentially ---} filtered colimits, and $I$ a directed set.
We consider an inductive system $\{J_i,h_{ij}\}_{i\in I}$ of directed sets.
Suppose we are given the following data:
\begin{itemize}
\item[{\rm (a)}] for any $i\in I$ an inductive system $\{X_{i,\alpha}\}_{\alpha\in J_i}$ consisting of objects and arrows in $\mathscr{C}$ indexed by $J_i$;
\item[{\rm (b)}] for any $i\leq j$ an arrow $X_{i,\alpha}\rightarrow X_{j,h_{ij}(\alpha)}$ in $\mathscr{C}$ such that, whenever $\alpha\leq\beta$ in $J_i$, the diagram
$$
\xymatrix@-1ex{X_{i,\alpha'}\ar[r]&X_{j,h_{ij}(\alpha')}\\ X_{i,\alpha}\ar[u]\ar[r]&X_{j,h_{ij}(\alpha)}\ar[u]}
$$
commutes.
\end{itemize}

(1) There exists for any $i\leq j$ a canonical arrow
$$
\varinjlim_{\alpha\in J_i}X_{i,\alpha}\longrightarrow\varinjlim_{\beta\in J_j}X_{j,\beta},
$$
by which one can consider the double inductive limit $\varinjlim_{i\in I}\varinjlim_{\alpha\in J_i}X_{i,\alpha}$.

(2) Set $\Lambda=\{(i,\alpha)\,|\,i\in I,\ \alpha\in J_i\}$, and consider the ordering on $\Lambda$ defined as follows: $(i,\alpha)\leq (j,\beta)$ if and only if $i\leq j$ and $h_{ij}(\alpha)\leq\beta$. Then $\Lambda$ is a directed set, and $\{X_{i,\alpha}\}_{(i,\alpha)\in\Lambda}$ is an inductive system indexed by $\Lambda$.

(3) There exists a canonical isomorphism
$$
\varinjlim_{i\in I}\varinjlim_{\alpha\in J_i}X_{i,\alpha}\stackrel{\sim}{\longrightarrow}\varinjlim_{(i,\alpha)\in\Lambda}X_{i,\alpha}.
$$
\end{exer}

\begin{exer}\label{exer-finalcountable}
Let $I$ be a non-empty directed set that admits a final and at most countable subset.
Show that there exists an ordered final map\index{map!final map@final ---} $L\colon\N\rightarrow I$.
\end{exer}

\begin{exer}\label{exer-equivalencerelationinvolution}{\rm 
Let $\mathscr{C}_J=(\mathscr{C},J)$ be a site, where $\mathscr{C}$ has a final object and all fiber products, and  
$$
\xymatrix{R\ar@<.5ex>[r]^(.5){p_1}\ar@<-.5ex>[r]_(.5){p_2}&Y}\leqno{(\ast)}
$$
a diagram of arrows in $\mathscr{C}$ such that for any $Z\in\obj(\mathscr{C})$ the induced map 
$$
\Hom_{\mathscr{C}}(Z,R)\longrightarrow\Hom_{\mathscr{C}}(Z,Y)\times\Hom_{\mathscr{C}}(Z,Y)
$$
is injective and that $p_1$ and $p_2$ are covering maps with respect to $J$.
Define $T$ by the Cartesian diagram
$$
\xymatrix{T\ar[r]^{p_{12}}\ar[d]_{p_{23}}&R\ar[d]^{p_2}\\ R\ar[r]_{p_1}&Y\rlap{.}}
$$
Show that the diagram $(\ast)$ gives a $J$-equivalence relation if and only if the following conditions are satisfied:
\begin{itemize}
\item[{\rm (a)}] there exists an arrow $\iota\colon Y\rightarrow R$ such that $p_1\circ\iota=p_2\circ\iota=\id_Y$;
\item[{\rm (b)}] there exists an automorphism $\sigma\colon R\stackrel{\sim}{\rightarrow}R$ such that $p_1\circ\sigma=p_2$ and $p_2\circ\sigma=p_1$ (in particular, we have $\sigma^2=\id_R$);
\item[{\rm (c)}] there exists an arrow $p_{13}\colon T\rightarrow R$ such that $p_1\circ p_{13}=p_1\circ p_{12}$ and $p_2\circ p_{13}=p_2\circ p_{23}$.
\end{itemize}}
\end{exer}

\section{General topology}\label{sec-gentop}
In this section, we discuss generality of the topological spaces that appear typically (if not generally) in algebraic geometry.
One of such typical `algebraic topological spaces' is the so-called {\em $($locally$)$ coherent\footnote{`Coherent' means `quasi-compact and quasi-separated'.} topological spaces}, introduced in \S\ref{sub-coherenttopologicalspace}.
For example, the underlying topological spaces of schemes are always locally coherent.
Coherent and sober topological spaces can be regarded as an algebraic analogue of compact Hausdorff spaces, as indicated in the following important fact: The filtered projective limit of coherent sober topological spaces by quasi-compact transition maps is again a coherent sober topological space.
This fact, which we will prove in \S\ref{subsub-projlimcoherentspaces} by means of Stone duality between topological spaces and lattices, contains the famous theorem by Zariski\index{Zariski, O.}, proved in his 1944 paper \cite{Zar4}, which asserts that his `generalized Riemann space' is quasi-compact.

The topological spaces that appear in rigid geometry as a topological visualization (cf.\ Introduction) are locally coherent topological spaces of a special kind, called {\em valuative spaces}, which we will introduce in \S\ref{sub-valuativespace}.
Already at the general topology level, valuative spaces have several interesting special features, such as overconvergent subsets, separated quotients, etc., all of which will play significant roles in rigid geometry.
Moreover, by means of a certain class of valuative spaces, the so-called reflexive (\ref{dfn-reflexivevaluativespaces}) locally strongly compact (\ref{dfn-locallycompactspace}) valuative spaces, one can develop an interesting variant of Stone duality, given in \S\ref{sub-valuationslocallyhausdorff}, which will be useful in understanding the relationship between our rigid geometry and Berkovich analytic geometry.

This section ends with a brief exposition of topos theory with special emphasis on the so-called {\em coherent topoi} and their limits.\footnote{Due to lack of space, we limit ourselves to a brief overview of the theory of locales and localic topoi. For these subjects, the reader is advised to see \cite{Johnstone2} and \cite{Johnstone}.}

{\sl In the sequel, we will follow the following conventions$:$
\begin{itemize}
\item compact spaces are always assumed to be Hausdorff$;$ that is, being compact is, by definition, equivalent to being quasi-compact and Hausdorff$;$ 
\item locally compact spaces are always assumed to be locally Hausdorff$;$ 
\item paracompact spaces are, however, {\em not} assumed to be Hausdorff $($cf.\ {\rm \S\ref{subsub-paracompactspaces})}. 
\end{itemize}}

\subsection{Some general prerequisites}\label{sub-gentopgenpre}
\subsubsection{Generization and specialization}\label{subsub-genspetopsp}
Let $X$ be a topological space.
We say that $y\in X$ is a {\em generization}\index{generization} of $x\in X$ or, equivalently, that $x$ is a {\em specialization}\index{specialization} of $y$, if $x$ belongs to the closure $\ovl{\{y\}}$ of $\{y\}$ in $X$ or, equivalently, if $y$ is contained in any open subset containing $x$.
Let us denote by $G_x$ the set of all generizations of $x$, that is, the intersection of all open neighborhoods of $x$:
$$
G_x=\bigcap_{x\in U}U.
$$

Clearly, if $z$ is a generization of $y$ and $y$ is a generization of $x$, then $z$ is a generization of $x$.
If $X$ is a $\mathrm{T}_0$-space\footnote{A topological space $X$ is said to be a {\em $\mathrm{T}_0$-space} (or {\em Kolmogorov space}) if, for any pair of two distinct points $x\neq y$ of $X$, there exists an open neighborhood of one of them that does not contain the other.}, the set $G_x$ is an ordered set\index{ordered!set@--- set}\index{set!ordered set@ordered ---} by the following ordering: for $y,z\in G_x$, $y\leq z$ if and only if $z$ is a generization of $y$.
A point $x$ said to be {\em maximal}\index{point!maximal point@maximal ---} (resp.\ {\em minimal}\index{point!minimal point@minimal ---}) if there exists no other generization (resp.\ specialization) of $x$ than $x$ itself.
For a point $x$ to be minimal it is necessary and sufficient that the singleton set $\{x\}$ is a closed subset of $X$, and hence, minimal points are also said to be {\em closed points}\index{point!closed point@closed ---}.
Notice that $x$ is the unique closed point of the subspace $G_x\subseteq X$ (that is, the topological subspace endowed with the subspace topology).

\subsubsection{Sober spaces}\label{subsub-sober}
Recall that a topological space $X$ is said to be {\em irreducible}\index{space@space (topological)!irreducible topological space@irreducible ---}\index{irreducible!topological space@--- topological space} if it is non-empty and cannot be the union of two closed subsets distinct from $X$.
If $X$ is not irreducible, it is said to be {\em reducible}\index{space@space (topological)!reducible topological space@reducible ---}.
Let us list some of the basic facts on irreducible spaces and irreducible subsets (cf.\ \cite[Chap.\ II, \S4.1]{Bourb1}):
\begin{itemize}
\item a non-empty topological space $X$ is irreducible if and only if any non-empty open subset is dense in $X$, and in this case, any non-empty open subset is again irreducible (and hence is connected);
\item a subset $Y$ of a topological space is irreducible if and only if its closure $\ovl{Y}$ is irreducible;
\item the image of an irreducible subset by a continuous mapping is again irreducible. 
\end{itemize}

Let $X$ be a topological space, $Z\subseteq X$ an irreducible closed subset, and $U\subseteq X$ an open subset such that $Z\cap U\neq\emptyset$.
Then $Z\cap U$ is an irreducible closed subset of $U$, and the closure $\ovl{Z\cap U}$ of $Z\cap U$ in $X$ coincides with $Z$.
In particular, if $Z_1$ and $Z_2$ are irreducible closed subsets of $X$, and if $U$ is an open subset such that $Z_i\cap U\neq\emptyset$ for $i=1,2$, then $Z_1\cap U=Z_2\cap U$ implies $Z_1=Z_2$.

For a topological space $X$ and a point $x\in X$, the subset $\ovl{\{x\}}$ consisting of all specializations of $x$ is an irreducible closed subset.
For a closed subset $Z\subseteq X$, a point $x\in Z$ is called a {\em generic point}\index{point!generic point@generic ---} of $Z$ if $Z=\ovl{\{x\}}$.
A topological space $X$ is said to be {\em sober}\index{space@space (topological)!sober topological space@sober ---} if it is a $\mathrm{T}_0$-space and any irreducible closed subset has a (unique) generic point.\footnote{The uniqueness of the generic point follows from $\mathrm{T}_0$-ness.}
Notice that any $\mathrm{T}_2$-space ($=$ Hausdorff space) is sober, since irreducible $\mathrm{T}_2$-space is precisely a singleton set (with the unique topology). 
Thus we have 
$$
\mathrm{T}_2\ \Longrightarrow\ \textrm{sober}\ \Longrightarrow\ \mathrm{T}_0,
$$
while $\mathrm{T}_1$-ness\footnote{A topological space $X$ is said to be a {\em $\mathrm{T}_1$-space} (or {\em Kuratowski space}) if for any pair of two distinct points $x\neq y$ of $X$ there exist an open neighborhood $U$ of $x$ and an open neighborhood $V$ of $y$ such that $x\not\in V$ and $y\not\in U$.} and soberness are not comparable (cf.\ Exercise \ref{exer-soberT1counterexample}).

The proofs of the following propositions are straightforward and left to the reader:
\begin{prop}\label{prop-propsoberness}
Every locally closed\footnote{A subset of a topological space $X$ is said to be {\em locally closed} if it is the intersection of an open subset and a closed subset.} subspace of a sober space is sober. \hfill$\square$
\end{prop}

\begin{prop}\label{prop-corsoberness}
If a topological space $X$ admits an open covering $X=\bigcup_{\alpha\in L}U_{\alpha}$ such that each $U_{\alpha}$ is sober, then $X$ is sober. \hfill$\square$
\end{prop}

We denote by $\STop$ the full subcategory of the category $\Top$ of topological spaces consisting of sober spaces.
\begin{prop}\label{prop-associatedsoberspace}
The inclusion functor $i\colon\STop\hookrightarrow\Top$ admits the left-adjoint functor
$$
\cdot^{\sob}\colon\Top\longrightarrow\STop,\qquad X\longmapsto X^{\sob}.\eqno{\square}
$$
\end{prop}

Briefly, $X^{\sob}$ is the set of all irreducible closed subsets of $X$ endowed with the topology with respect to which open subsets are of the form $\til{U}$, each associated to an open subset $U\subseteq X$, consisting of irreducible closed subsets of $X$ that intersect $U$.
For more detail (and the proof of \ref{prop-associatedsoberspace}), see \cite[$\mathbf{0}$, \S2.9]{EGAInew}.

\subsubsection{Completely regular spaces}\label{subsub-crspaces}
Recall that a {\em completely regular}\index{space@space (topological)!completely regular topological space@completely regular ---} space is a $\mathrm{T}_1$-space $X$ that enjoys the following property: for any point $x\in X$ and a closed subset $Z\subseteq X$ not containing $x$, there exists a continuous function $f\colon X\rightarrow [0,1]\subset\R$ such that $f(x)=0$ and $f(C)=\{1\}$.
By Urysohn's lemma, {\em normal} (that is, $\mathrm{T}_1$ and $\mathrm{T}_4$\footnote{A topological space $X$ is said to be a {\em $\mathrm{T}_4$-space} if for any disjoint pair $F_1,F_2$ of closed subsets of $X$ there exists a disjoint pair $U_1,U_2$ of open neighborhoods of $F_1,F_2$, respectively.}) topological spaces and locally compact Hausdorff spaces are completely regular (consider the one-point compactification for the latter).

\subsubsection{Quasi-compact spaces and quasi-separated spaces}\label{subsub-quasicompactmappings}
\begin{dfn}\label{dfn-quasicompactness}{\rm 
(1) A topological space $X$ is said to be {\em quasi-compact}\index{space@space (topological)!quasi-compact topological space@quasi-compact ---}\index{quasi-compact!quasi-compact topological space@--- (topological) space} if for any open covering $\{U_{\alpha}\}_{\alpha\in L}$ of $X$ there exists a finite subset $L'\subseteq L$ of the indices such that $\{U_{\alpha}\}_{\alpha\in L'}$ already covers $X$.

(2) A continuous map $f\colon X\rightarrow Y$ of topological spaces is said to be {\em quasi-compact}\index{map@map (continuous)!quasi-compact map@quasi-compact ---}\index{quasi-compact!quasi-compact map@--- map} if for any quasi-compact open subset $V$ of $Y$ its pull-back $f^{-1}(V)$ is a quasi-compact open subset of $X$.}
\end{dfn}

The set of all quasi-compact open subsets on a topological space is closed under finite union.
The composition of two quasi-compact maps is again quasi-compact.
If $X$ is a quasi-compact space, and if $Z\subseteq X$ is a closed subset, then $Z$ is quasi-compact, and the inclusion map $i\colon Z\hookrightarrow X$ is quasi-compact.

\begin{prop}\label{prop-quasicompactness}
Let $X$ be a topological space, and $Y\subseteq X$ a subspace.
Suppose that $X$ admits an open basis consisting of quasi-compact open subsets.
Then for any quasi-compact open subset $V$ of $Y$, there exists a quasi-compact open subset $U$ of $X$ such that $V=Y\cap U$. \hfill$\square$
\end{prop}

\begin{cor}\label{cor-quasicompactness}
In the situation as in $\ref{prop-quasicompactness}$, let $g\colon Z\rightarrow Y$ be a continuous map.
If $j\circ g\colon Z\rightarrow X$ is quasi-compact, then so is $g\colon Z\rightarrow Y$.\hfill$\square$
\end{cor}

\begin{dfn}\label{dfn-retrocompact}{\rm 
Let $X$ be a topological space. A subset $Z\subseteq X$ is said to be {\em retrocompact}\index{retrocompact} if for any quasi-compact open subset $U\subset X$ the intersection $Z\cap U$ is quasi-compact.}
\end{dfn}

That is, $Z$ is retrocompact if and only if the inclusion map $Z\hookrightarrow X$, where $Z$ is endowed with the subspace topology, is quasi-compact (\ref{dfn-quasicompactness} (2)).
The union of finitely many retrocompact subsets is again retrocompact.
The intersection of finitely many retrocompact open subsets is again a retrocompact open subset.

\begin{dfn}\label{dfn-quasiseparatedness}{\rm 
A topological space $X$ is said to be {\em quasi-separated}\index{space@space (topological)!quasi-separated topological space@quasi-separated ---}\index{quasi-separated!quasi-separated topological space@--- (topological) space} if, for any two quasi-compact open subsets $U,V\subseteq X$, the intersection $U\cap V$ is again quasi-compact.}
\end{dfn}

In other words, $X$ is quasi-separated if and only if any quasi-compact open subset of $X$ is retrocompact\index{retrocompact} {\rm (\ref{dfn-retrocompact}).
The set of all quasi-compact open subsets in a quasi-separated space is closed under finite intersection.
Any open subset of a quasi-separated space is again quasi-separated.

\subsection{Coherent spaces}\label{sub-coherenttopologicalspace}
\index{space@space (topological)!coherent topological space@coherent ---|(}\index{coherent!coherent topological space@--- (topological) space|(}
\subsubsection{Definition and first properties}\label{subsub-dfncohtop}
\begin{dfn}\label{dfn-quasicompact1}{\rm 
A topological space $X$ is said to be {\em coherent}\index{space@space (topological)!coherent topological space@coherent ---}\index{coherent!coherent topological space@--- (topological) space} if the following conditions are satisfied$:$
\begin{itemize}
\item[{\rm (a)}] $X$ has an open basis consisting of quasi-compact open subsets$;$
\item[{\rm (b)}] $X$ is quasi-compact\index{space@space (topological)!quasi-compact topological space@quasi-compact ---}\index{quasi-compact!quasi-compact topological space@--- (topological) space} and quasi-separated\index{space@space (topological)!quasi-separated topological space@quasi-separated ---}\index{quasi-separated!quasi-separated topological space@--- (topological) space}.
\end{itemize}
(See \S\ref{subsub-coherenttopos} for the topos-theoretic interpretation.)}
\end{dfn}

\begin{exas}\label{exas-coherenttopsp}{\rm 
(1) Empty set and singleton sets (endowed with the unique topology) are coherent; more generally, any finite space is coherent.
Any continuous mapping from a finite space to an arbitrary coherent space is quasi-compact.

(2) The underlying topological space of a scheme admits an open basis consisting of quasi-compact open subsets.
Hence, the underlying topological space of a quasi-compact and quasi-separated scheme (e.g.\ Noetherian schemes) is coherent (due to \cite[$\mathbf{IV}$, (1.2.7)]{EGA}); in particular, the underlying topological space of an affine scheme\index{affine!affine scheme@--- scheme} is coherent.}
\end{exas}

For the reason mentioned in (2), quasi-compact and quasi-separated schemes are often called {\em coherent schemes}\index{scheme!coherent scheme@coherent ---}\index{coherent!coherent scheme@--- scheme}.

\begin{prop}\label{prop-coherentclosedqcopen}
Let $X$ be a coherent topological space.
Then any quasi-compact locally closed subset $Z\subseteq X$, endowed with the subspace topology, is coherent, and is retrocompact\index{retrocompact} in $X$. 
\end{prop}

\begin{proof}
By \ref{prop-quasicompactness}, any quasi-compact locally closed subset of $X$ is the intersection of a closed subset and a quasi-compact open subset.
Hence it suffices to check the assertion in the cases where $Z$ is closed, and where $Z$ is quasi-compact open.
Both cases are easy to verify.
\end{proof}

\begin{rem}\label{rem-quasicompact}{\rm 
(1) Coherent and sober spaces are also called {\em spectral spaces}\index{space@space (topological)!spectral topological space@spectral ---|see{coherent space}} by some authors.

(2) Hochster \cite{Hochster} has shown that any coherent sober space is {\em homeomorphic} to the prime spectrum $\Spec A$ of a commutative ring $A$.}
\end{rem}

One of the most remarkable features of coherent spaces is that, as we will see soon later (\ref{thm-projlimcohsch1} (1)), the small filtered projective limit of a projective system consisting of coherent sober topological spaces with quasi-compact transition maps is again coherent and sober.
In connection with this, the reader is invited, before proceeding to the next paragraph, to try out Exercise \ref{exer-projectivelimitaffine}, which deals with an analogous topic on filtered projective limits of schemes.

\subsubsection{Stone's representation theorem}\label{subsub-strcohtopsp}
For the reader's convenience, we include basic facts on the relationship between distributive lattices and coherent sober spaces.
Our basic reference to the theory of distributive lattices\index{lattice}\index{lattice!distributive lattice@distributive ---} is \cite{Johnstone}.

Let $A$ be a distributive lattice (we use the binary symbols $\vee$, $\wedge$, and $\leq$ as in \cite{Johnstone} and always assume that $A$ has $0$ and $1$).
We view $A$ as a category with respect to the partial order of $A$ in the manner mentioned in \S\ref{subsub-catposet}.
Then, as a category, $A$ is stable under finite limits and finite colimits, and finite disjoint sum is universally disjoint (\cite[Expos\'e II, 4.5]{SGA4-1}).
A lattice homomorphism $A\rightarrow B$ is, in the categorical language, a functor that commutes with finite limits and colimits. 
Let $\DLat$ be the category of distributive lattices with lattice homomorphisms.

\begin{exa}\label{exa-latticeexafoherent}{\rm 
For a topological space $X$ we denote by $\Ouv(X)$ the set of all open subsets of $X$.
The set $\Ouv(X)$ forms a distributive lattice with the structure $(0,1,\vee,\wedge,\leq)=(\emptyset,X,\cup,\cap,\subseteq)$.
If $X$ is coherent, then the set $\QCOuv(X)$ of all quasi-compact open subsets of $X$ forms a distributive sublattice of $\Ouv(X)$.}
\end{exa}

There are notions of {\em ideals}\index{ideal!ideal of lattices@(of lattices)} and {\em filters}\index{filter (of lattices)}, which are dual to each other:
\begin{dfn}\label{dfn-idealfilter}{\rm 
Let $A$ be a distributive lattice.

(1) A subset $I\subseteq A$ is called an {\em ideal} if (a) $0\in I$, (b) $a,b\in I$ implies $a\vee b\in I$, and (c) $a\in I$ and $b\leq a$ imply $b\in I$. 

(2) A subset $F\subseteq A$ is called a {\em filter} if (a) $1\in F$, (b) $a,b\in F$ implies $a\wedge b\in F$, and (c) $a\in F$ and $b\geq a$ imply $b\in F$.}
\end{dfn}

Accordingly, by an ideal (resp.\ a filter) of a topological space $X$, we mean an ideal (resp.\ a filter) of the complete distributive lattice $\Ouv(X)$.

For $A\in\obj(\DLat)$ the lattice $\Id(A)$ of all ideals of $A$ forms a complete Heyting algebra, which admits the embedding $A\hookrightarrow\Id(A)$ of lattices that maps $a$ to the principal ideal $\{x\,|\,x\leq a\}$ generated by $a$. 
The image of $A$ in $\Id(A)$ consists exactly of the sets of finite elements (\cite[Chap.\ II, \S3.1]{Johnstone}).

\begin{dfn}\label{dfn-latticeprimeideal}{\rm 
An ideal $P$ is said to be {\em prime}\index{ideal!prime ideal of lattices@prime --- (of lattices)} if (a) $1\not\in P$, and (b) $a\wedge b\in P$ implies $a\in P$ or $b \in P$.}
\end{dfn}

In other words, an ideal $P$ is prime if and only if $A\setminus P$ is a filter. 
A filter of the form $A\setminus P $ (where $P$ is a prime ideal) is called a {\em prime filter}\index{filter (of lattices)!prime filter of lattices@prime ---}. 
Here is another description of prime ideals: 
Let $\mathbf{2}=\{0,1\}$ be the Boolean lattice\index{lattice!Boolean lattice@Boolean ---} consisting of two elements.
Then any prime ideal corresponds to the kernel of a surjective lattice homomorphism $A\rightarrow\mathbf{2}$.

Any non-zero distributive lattice has at least one prime ideal (one needs, similarly to the case of commutative rings, the axiom of choice to prove this).

Let $\Spec A$ denote the set of prime ideals of $A$. 
For any ideal $I\subset A$ we define $V(I)$ to be the set of prime ideals containing $I$. 
The Zariski topology on $\Spec A$ is the topology of which $\{V(I)\}_{I\in\Id(A)}$ is the system of closed sets. 
We set $D(I)=\Spec A\setminus V(I)$.

\begin{thm}[Stone's representation theorem; cf.\ {\cite[Chap.\ II, \S3.4]{Johnstone}}]\label{thm-stonerepresentationthm}\index{Stone's representation theorem}\ 

{\rm (1)} For any distributive lattice $A$, the topological space $\Spec A$ is coherent\index{space@space (topological)!coherent topological space@coherent ---}\index{coherent!coherent topological space@--- (topological) space} and sober\index{space@space (topological)!sober topological space@sober ---}, and we have the functor
$$
\Spec\colon\DLat^{\opp}\longrightarrow\CSTop,\leqno{(\ast)}
$$
where $\CSTop$ denotes the category of coherent sober topological spaces and quasi-compact maps\index{map@map (continuous)!quasi-compact map@quasi-compact ---}\index{quasi-compact!quasi-compact map@--- map}.
The distributive lattice $\Ouv(\Spec A)$ consisting of all open subsets of $\Spec A$ is isomorphic to $\Id(A)$ by $I\mapsto D(I)$.
Moreover, the lattice $A$ is identified, by means of the embedding $A\hookrightarrow\Id(A)$ and the isomorphism $\Id(A)\stackrel{\sim}{\rightarrow}\Ouv(\Spec A)$, with the distributive lattice of all quasi-compact open subsets of $\Spec A$.

{\rm (2)} The functor $(\ast)$ is an equivalence of categories.
The quasi-inverse functor is given by $X\mapsto\QCOuv(X)$ {\rm (\ref{exa-latticeexafoherent})}. \hfill$\square$
\end{thm}

\subsubsection{Projective limit of coherent sober spaces}\label{subsub-projlimcoherentspaces}
The category $\Top$ is closed under small projective limits\index{limit!projective limit@projective ---}.
For a projective system $\{X_i,p_{ij}\}$ of topological spaces indexed by an ordered set\index{ordered!set@--- set}\index{set!ordered set@ordered ---} $I$ (cf.\ \cite[Chap.\ I, \S4.4]{Bourb4}), the projective limit $X=\varprojlim_{i\in I}X_i$ has, as its underlying set, the set-theoretic projective limit of the underlying sets of $X_i$'s, endowed with the coarsest topology for which all projection maps $p_i\colon X\rightarrow X_i$ are continuous, or equivalently, the topology generated by the subsets of the form $p^{-1}_i(U)$ for $i\in I$ and open subsets $U\subseteq X_i$.

\begin{prop}\label{prop-projlimcohtopspqcptopen}
Let $\{X_i,p_{ij}\colon X_j\rightarrow X_i\}$ be a projective system of topological spaces, indexed by a directed set\index{set!directed set@directed ---}\index{directed set} $(I,\leq )$.
Suppose that the topology of each $X_i$ $(i\in I)$ is generated by quasi-compact open subsets, and that each transition map $p_{ij}$ $(i\leq j)$ is quasi-compact\index{map@map (continuous)!quasi-compact map@quasi-compact ---}\index{quasi-compact!quasi-compact map@--- map} {\rm (\ref{dfn-quasicompactness} (2))}.
Set $X=\varprojlim_{i\in I}X_i$, and let $U\subseteq X$ be a quasi-compact open subset.
Then there exists $i\in I$ and a quasi-compact open subset $U_i\subseteq X_i$ such that $p^{-1}_i(U_i)=U$, where $p_i\colon X\rightarrow X_i$ is the canonical projection map.
\end{prop}

\begin{proof}
By the definition of the topology of $X$ and by the quasi-compactness of $U$, there exist finitely many $i_1,\ldots,i_n\in I$ and, for each $k=1,\ldots,n$, a quasi-compact open subset $U_{i_k}\subseteq X_{i_k}$, such that $U=\bigcup^n_{k=1}p^{-1}_{i_k}(U_{i_k})$ (cf.\ \cite[Chap.\ I, \S4.4, Prop.\ 9]{Bourb4}).
Take $i\in I$ such that $i\geq i_1,\ldots,i_n$, and set $U_i=\bigcup^n_{k=1}p^{-1}_{i_ki}(U_{i_k})$.
Then $U_i$ is quasi-compact, and $U=p^{-1}_i(U_i)$, as desired.
\end{proof}

\begin{thm}\label{thm-projlimcohsch1}
Let $(I,\leq)$ be a directed set\index{set!directed set@directed ---}\index{directed set}, and $\{X_i,p_{ij}\}_{i\in I}$ a projective system\index{system!projective system@projective ---} of coherent sober\index{space@space (topological)!sober topological space@sober ---} spaces\index{space@space (topological)!coherent topological space@coherent ---}\index{coherent!coherent topological space@--- (topological) space} with quasi-compact\index{map@map (continuous)!quasi-compact map@quasi-compact ---}\index{quasi-compact!quasi-compact map@--- map} transition maps $p_{ij}$ for $i\geq j$.
Set $X=\varprojlim_{i\in I}X_i$.

{\rm (1)} The topological space $X$ is coherent and sober, and the canonical projection maps $p_i\colon X\rightarrow X_i$ for $i\in I$ are quasi-compact.

{\rm (2)} If, moreover, each $X_i$ $(i\in I)$ is non-empty, then $X$ is non-empty.
\end{thm} 

By (1) one finds that, in particular, the category $\CSTop$ is closed under small filtered projective limits.
Notice that the first assertion implies that the projective limit space $X$ is, in particular, quasi-compact. 

\begin{proof}
Set $A_i=\QCOuv(X_i)$ for each $i\in I$. 
The transition maps $p_{ij}$ for $i\leq j$ induce the lattice homomorphisms $u_{ij}\colon A_i\rightarrow A_j$, which form an inductive system $\{A_i\}_{i\in I}$ of lattices. 
Let $A=\varinjlim_{i\in I}A_i$ be the inductive limit. 
By \ref{thm-stonerepresentationthm} the topological space $X'=\Spec A$ is the projective limit of $\{X_i,p_{ij}\}$ in category $\CSTop$. 
If all $X_i$'s are non-empty, each distributive lattice $A_i$ is non-zero, that is, $0\neq 1$ in $A_i$.
This implies that $0\neq 1$ also in $A$, hence $X'\neq\emptyset$.

Hence, to prove the assertions, it remains to show that $X$ and $X'$ are homeomorphic.
By the mapping universality of projective limits, we have the canonical map $g\colon X'=\Spec A\rightarrow X$.
We are first going to see that $g$ is bijective.
To see this, notice that each point $x$ of $X$ is uniquely represented by a sequence $\{x_i\in X_i\}_{i\in I}$ of points such that $p_{ij}(x_j)=x_i$ for $i\leq j$, which is, furthermore, uniquely interpreted into a system $\{h_i\colon A_i\rightarrow\mathbf{2}\}_{i\in I}$ of surjective lattice homomorphisms such that $h_j\circ u_{ij}=h_i$ for $i\leq j$ ($h_i$ is defined as follows: it maps $U\in A_i=\QCOuv(X_i)$ to $1\in\mathbf{2}$ if $x_i\in U$, and to $0$ otherwise).
Giving the last data is then equivalent to giving a surjective lattice homomophism $h\colon A\rightarrow\mathbf{2}$, whose kernel corresponds to a point $x'\in X'=\Spec A$.
It is easy to see that $x\mapsto x'$ thus constructed gives the inverse mapping of $g\colon X'\rightarrow X$, and hence $g$ is a continuous bijection.
Moreover, quasi-compact open subsets of $X'$ are of the form $p^{\prime -1}_i(U)$ by a quasi-compact open subset $U$ of some $X_i$ (\ref{prop-projlimcohtopspqcptopen}).
Since $g$ is bijective, we have $g(p^{\prime -1}_i(U))=p^{-1}_i(U)$, and hence $g$ maps any quasi-compact open subset of $X'$ to an open subset of $X$.
Since quasi-compact open subsets of $X'$ form an open basis of the topology on $X'$, we conclude that $g\colon X'\rightarrow X$ is a homeomorphism.
\end{proof}

\begin{rem}\label{rem-projlimcohsch1}{\rm 
The above theorem can be more properly formulated in topos theory, where the projective system $\{X_i\}_{i\in I}$ is translated into the projective system of corresponding coherent topoi $\{\top(X_i)\}_{i\in I}$.
Under this interpretation, the theorem follows from \cite[Exp.\ VI, (8.3.13)]{SGA4-2} and Deligne's theorem \cite[Exp.\ VI, (9.0)]{SGA4-2}; see (\ref{cor-fiberedtopos21}) below.}
\end{rem}

Theorem \ref{thm-projlimcohsch1} is not only important in its own light, but has many useful consequences:
\begin{cor}\label{cor-projlimcohsch11}
Let $\{X_i,p_{ij}\}_{i\in I}$ be a filtered projective system of coherent sober topological spaces and quasi-compact maps indexed by a directed set $I$.
Let $i\in I$ be an index, and $U,V\subseteq X_i$ open subsets, where $U$ is quasi-compact.
Then the following conditions are equivalent$:$
\begin{itemize}
\item[{\rm (a)}] $p^{-1}_i(U)\subseteq p^{-1}_i(V)$, where $p_i\colon X=\varprojlim_{j\in I}X_j\rightarrow X_i$ is the canonical mapping$;$
\item[{\rm (b)}] there exists an index $j$ with $i\leq j$ such that $p^{-1}_{ij}(U)\subseteq p^{-1}_{ij}(V)$.
\end{itemize}
\end{cor}

\begin{proof}
The implication (b) $\Rightarrow$ (a) is clear.
To show (a) $\Rightarrow$ (b), set $Z_j=p^{-1}_{ij}(U)\setminus p^{-1}_{ij}(V)$ for each $j\in I$ with $i\leq j$.
We need to show that $Z_j$ is empty for some $j$.
Suppose that all $Z_j$ are non-empty.
By \ref{prop-coherentclosedqcopen} and \ref{prop-propsoberness}, each $Z_j$ is coherent and sober.
Applying \ref{thm-projlimcohsch1} (2), we deduce that the projective limit $Z=\varprojlim_{j\geq i}Z_j$ is non-empty.
On the other hand, since $p^{-1}_i(U)=\varprojlim_{j\geq i}p^{-1}_{ij}(U)$ and $p^{-1}_i(V)=\varprojlim_{j\geq i}p^{-1}_{ij}(V)$, we have $Z=p^{-1}_i(U)\setminus p^{-1}_i(V)$.
But, then, $Z\neq\emptyset$ contradicts the assumption $p^{-1}_i(U)\subseteq p^{-1}_i(V)$, whence the result.
\end{proof}

\begin{thm}\label{thm-projlimcohspacepres}
Let $\{X_i,p_{ij}\}_{i\in I}$ and $\{Y_i,q_{ij}\}_{i\in I}$ be two filtered projective systems of coherent sober topological spaces with quasi-compact transition maps indexed by a common directed set $I$, and $\{f_i\}_{i\in I}$ a map of the systems, that is, a collection of continuous maps $f_i\colon X_i\rightarrow Y_i$ such that the diagram
$$
\xymatrix{X_j\ar[r]^{f_j}\ar[d]_{p_{ij}}&Y_j\ar[d]^{q_{ij}}\\ X_i\ar[r]_{f_i}&Y_i}
$$
commutes for any $i\leq j$.
Let
$$
f=\varprojlim_{i\in I}f_i\colon X=\varprojlim_{i\in I}X_i\longrightarrow Y=\varprojlim_{i\in I}Y_i
$$
be the induced map.

{\rm (1)} If each $f_i$ for $i\in I$ is quasi-compact, then $f$ is quasi-compact.

{\rm (2)} If each $f_i$ for $i\in I$ is quasi-compact and surjective, then $f$ is surjective.

{\rm (3)} If each $f_i$ for $i\in I$ is quasi-compact and closed, then $f$ is closed.
\end{thm}

Before proving the theorem, let us present a useful corollary of it:
\begin{cor}\label{cor-projlimcohspacelattice2}
Let $\{X_i,p_{ij}\}_{i\in I}$ be a projective system of coherent sober spaces with quasi-compact transition maps indexed by a directed set $I$.
Suppose that all transition maps $p_{ij}$ are surjective $($resp.\ closed$)$. 
Then for any $j\in I$ the canonical projection $p_j\colon X=\varprojlim X_i \rightarrow X_j$ is surjective $($resp.\ closed$)$. 
\end{cor}

\begin{proof}
For each $j\in I$, replace $I$ by the cofinal subset $\{i\in I\,|\,i\geq j\}$ and apply \ref{thm-projlimcohspacepres} (2) and (3) to the constant projective system $\{Y_i,q_{ij}\}$ with $Y_i=X_j$.
\end{proof}

The rest of this subsection is devoted to the proof of \ref{thm-projlimcohspacepres}.
We denote by $p_i\colon X\rightarrow X_i$ and $q_i\colon Y\rightarrow Y_i$ the canonical projections for $i\in I$.

\begin{proof}[Proof of Theorem {\rm \ref{thm-projlimcohspacepres} (1)}]
We need to show that, for any quasi-compact open subset $V$ of $Y$, its pull-back $f^{-1}(V)$ is quasi-compact.
In view of \ref{prop-projlimcohtopspqcptopen} one can choose $i\in I$ and a quasi-compact open subset $V_i$ of $Y_i$ such that $V=q^{-1}_i(V_i)$.
Then $f^{-1}(V)=p^{-1}_i(f^{-1}_i(V_i))$, which is quasi-compact by \ref{thm-projlimcohsch1} (1), and this shows (1) of \ref{thm-projlimcohspacepres}.
\end{proof}

To show the other assertions of the theorem, we need a few lemmas:
\begin{lem}\label{lem-projlimcohspacelattice001}
Let $X$ be a coherent sober topological space, and $F$ a filter\index{filter (of lattices)} of the distributive lattice $\Ouv(X)$ {\rm (\ref{dfn-idealfilter} (2))} generated by quasi-compact open subsets.
Then the subspace $\bigcap_{U\in F}U$ is coherent and sober, and the inclusion map $\bigcap_{U\in F}U\hookrightarrow U$ for any $U\in F$ $($in particular, $\bigcap_{U\in F}U\hookrightarrow X)$ is quasi-compact.
\end{lem}

\begin{proof}
Consider the subfilter $F'=F\cap\QCOuv(X)$ of $F$. 
Since topology on $X$ is generated by quasi-compact open subsets, the inclusion map $F'\hookrightarrow F$ is cofinal.
Since $\bigcap_{U\in F}U=\varprojlim_{U\in F'}U$, the lemma follows from \ref{thm-projlimcohsch1} (1).
\end{proof}

\begin{cor}\label{cor-projlimcohspacelattice1}
Let $X$ be a coherent sober space, and $x\in X$.
Then the subset $G_x$ of all generizations\index{generization} of $x$ with the subspace topology is coherent and sober, and the inclusion map $G_x\hookrightarrow X$ is quasi-compact. 
\end{cor}

\begin{proof}
Since $X$ has an open basis consisting of quasi-compact open subsets, and any quasi-compact open subset is coherent (\ref{prop-coherentclosedqcopen}) and sober, the subset $G_x$ coincides with the intersection of all coherent open neighborhoods of $x$ (cf.\ \S\ref{subsub-genspetopsp}), and hence the assertion follows from \ref{lem-projlimcohspacelattice001}.
\end{proof}

\begin{lem}\label{lem-pprojlimcohspacelattice1}
Let $f\colon X\rightarrow Y$ be a quasi-compact map between coherent sober topological spaces. 

{\rm (1)} Let $C$ be a closed subset of $Y$. 
Then $C$ and $f^{-1}(C)$ are coherent and sober. 
The map $f^{-1}(C)\rightarrow C$ induced by $f$ is quasi-compact. 

{\rm (2)} Let $F$ be a filter of $\Ouv(Y)$ generated by quasi-compact open sets. 
Then $\bigcap_{U\in F}U$ and $f^{-1}(\bigcap_{U\in F}U)$ are coherent and sober, and the map $f^{-1}(\bigcap_{U\in F}U)\rightarrow\bigcap_{U\in F}U$ induced by $f$ is quasi-compact. 
\end{lem}

In particular, it follows from (2) and \ref{lem-projlimcohspacelattice001} that the map $f^{-1}(\bigcap_{U\in F}U)\rightarrow X$ is quasi-compact.

\begin{proof}
The first assertion follows immediately from \ref{prop-coherentclosedqcopen}, and the second from \ref{lem-projlimcohspacelattice001}.
\end{proof}

\begin{lem}\label{lem-projlimcohspacelattice1}
{\rm (1)} Let $f\colon X\rightarrow Y$ be a quasi-compact map between coherent sober topological spaces.  
Then for any $y\in Y$ the set $f^{-1}(y)$ with the subspace topology is coherent and sober, and the inclusion $f^{-1}(y)\hookrightarrow X$ is quasi-compact. 

{\rm (2)} Consider a commutative diagram of coherent sober spaces and quasi-compact maps
$$
\xymatrix{X'\ar[r]^{f'}\ar[d]&Y'\ar[d]\\ X\ar[r]_f&Y\rlap{.}}
$$
Let $y'\in Y'$, and $y\in Y$ the image of $y'$ by the map $Y'\rightarrow Y$.
Then the natural map ${f'}^{-1}(y')\rightarrow f^{-1}(y)$ is quasi-compact. 
\end{lem}

\begin{proof}
(1) By \ref{cor-projlimcohspacelattice1} the subset $G_y$ of all generizations of $y$ is coherent and sober, and the inclusion map $G_y\hookrightarrow Y$ is quasi-compact.
Moreover, by \ref{lem-pprojlimcohspacelattice1}, $f^{-1}(G_y)$ is coherent sober, and the inclusion map $f^{-1}(G_y)\hookrightarrow X$ is quasi-compact.
Since $y$ is the closed point of $G_y$, $f^{-1}(y)$ is a closed subset of $f^{-1}(G_y)$.
Hence by \ref{prop-coherentclosedqcopen}, $f^{-1}(y)$ is coherent and sober, and the inclusion map $f^{-1}(y)\hookrightarrow f^{-1}(G_y)$ is quasi-compact.
Then the inclusion $f^{-1}(y)\hookrightarrow X$, the composition of $f^{-1}(y)\hookrightarrow f^{-1}(G_y)$ and $f^{-1}(G_y)\hookrightarrow X$, is quasi-compact.

(2) It suffices to show the following: If (i) $g\colon X'\rightarrow X$ is a quasi-compact map between coherent sober spaces, (ii) $Z'\subseteq X'$ and $Z\subseteq X$ are subspaces such that $g(Z')\subseteq Z$, and (iii) the inclusions $Z'\hookrightarrow X'$ and $Z\hookrightarrow X$ are quasi-compact, then the map $g\colon Z'\rightarrow Z$ is quasi-compact (indeed, $Z'=f^{\prime -1}(y')$ and $Z=f^{-1}(y)$ give our original situation).
This follows from \ref{cor-quasicompactness}.
\end{proof}

Now, we can prove the second part of \ref{thm-projlimcohspacepres}:
\begin{proof}[Proof of Theorem {\rm \ref{thm-projlimcohspacepres} (2)}]
Take $y\in Y$, and set $y_i=q_i(y)$ for each $i\in I$. 
By our assumptions and \ref{lem-projlimcohspacelattice1}, the set $Z_i=f_i^{-1}(y_i)$ is coherent, sober, and non-empty, and the map $Z_j\rightarrow Z_i$ induced by $p_{ij}$ is quasi-compact. 
Then it follows from \ref{thm-projlimcohsch1} (2) that 
$$
f^{-1}(y)=\varprojlim_{i\in I}Z_i\neq\emptyset,
$$
which verifies the claim.
\end{proof}

For the proof of the third part, we need the following lemma; the proof is straightforward and left to the reader:
\begin{lem}\label{lem-projlimclosedmapslem}
Let $\{X_i,p_{ij}\}_{i\in I}$ be a projective system of topological spaces indexed by a directed set $I$. 
We set $X=\varprojlim_{i\in I}X_i$ and denote the projection $X\rightarrow X_i$ by $p_i$. 

{\rm (1)} Suppose we are given a projective system $\{C_i,{p_{ij}}|_{C_i}\}_{i\in I}$ consisting of subspaces $C_i\subseteq X_i$ for $i\in I$.
Then we have
$$
\varprojlim_{i\in I}C_i=\bigcap_{i\in I}p_i^{-1}(C_i).
$$

{\rm (2)} For any subset $C\subseteq X$ we have
$$
\ovl{C}=\varprojlim_{i\in I}\ovl{p_i(C)}=\bigcap_{i\in I}p_i^{-1}(\ovl{p_i(C)}).\eqno{\square}
$$
\end{lem}

\begin{proof}[Proof of Theorem {\rm \ref{thm-projlimcohspacepres} (3)}]
We need to show that for any closed subset $C\subseteq X$ the image $f(C)$ is closed.
Set $C_i=\ovl{p_i(C)}$ and $D_i=f_i(C_i)$.
These are closed subsets of $X_i$ and $Y_i$, respectively.
By \ref{lem-projlimclosedmapslem} we have $C=\varprojlim_{i\in I}C_i$.  
On the other hand, $\varprojlim_{i\in I}D_i=\bigcap_{i\in I}q^{-1}_i(D_i)$ is a closed subset of $Y$.
Hence it suffices to show the equality 
$$
f(C)=\varprojlim_{i\in I}D_i.
$$
It is clear that the left-hand side is contained in the right-hand side.
In order to show the converse, take $y\in\varprojlim D_i=\bigcap_{i\in I}q^{-1}_i(D_i)$.
For each $i$, the set $Z_i=f^{-1}_i(q_i(y))\cap C_i$ is non-empty, and $\{Z_i\}_{i\in I}$ forms a projective system of coherent sober spaces with quasi-compact transitions maps by \ref{lem-projlimcohspacelattice1} and \ref{lem-pprojlimcohspacelattice1} (1).
Then it follows that $\varprojlim Z_i\neq\emptyset$, and any point $x$ in this projective limit, regarded as a point of $C$, is mapped to $y$ by the construction.
\end{proof}
\index{coherent!coherent topological space@--- (topological) space|)}\index{space@space (topological)!coherent topological space@coherent ---|)}

\begin{small}
Here we include a theorem, known as Steenrod's theorem (cf.\ \cite{WY}), informed to the authors from O.\ Gabber\index{Gabber, O.}, as well as the proof suggested by him, on quasi-compactness of the projective limit spaces, although we will not use it in our later discussion.
Later in \S\ref{subsub-quasicompacttopos}, we will indicate a topos-theoretic proof.
\begin{thm}\label{thm-gabbertheoremprojectivelimits}
Let $(I,\leq)$ be a directed set, and $\{X_i,p_{ij}\}$ a projective system indexed by $I$ consisting of quasi-compact sober topological spaces.
$($Here we do not assume that the transition maps $p_{ij}$ are quasi-compact.$)$

{\rm (1)} The projective limit $X=\varprojlim_{i\in I}X_i$ is quasi-compact and sober.

{\rm (2)} If each $X_i$ $(i\in I)$ is non-empty, then so is the limit $X$.
\end{thm}

\begin{proof}
First we show (2). 
Let $\Phi$ be the set of all projective systems consisting of non-empty closed subsets of $X_i$'s; that is, each element of $\Phi$ is a collection $\{Y_i\}_{i\in I}$ of non-empty closed subsets $Y_i\subseteq X_i$ for $i\in I$ such that $p_{ij}(Y_j)\subseteq Y_i$ for any $i\leq j$.
We define a partial ordering $\leq$ on the set $\Phi$ as follows: for $y=\{Y_i\}$ and $y'=\{Y'_i\}$, 
$$
y\leq y'\quad\Longleftrightarrow\quad Y'_i\subseteq Y_i\ \textrm{for any $i\in I$}.
$$
One deduces from quasi-compactness of each $X_i$ that any totally ordered subset of $\Phi$ has an upper bound, and thus that there exists a maximal element $z=\{Z_i\}$ by Zorn's lemma.

\medskip
{\sc Claim 1.} {\it Let $i\in I$, and consider a non-empty closed subset $W\subseteq Z_i$.
Suppose there exists a cofinal subset $J$ of $I_{\geq i}=\{j\in I\,|\,j\geq i\}$ such that $p^{-1}_{ij}(W)\cap Z_j\neq\emptyset$ for any $j\in J$.
Then we have $W=Z_i$.}

\medskip
To show this, set 
$$
Z_k=\bigcap_{j\in J\cap I_{\geq k}}\ovl{p_{kj}(p^{-1}_{ij}(W)\cap Z_j)},
$$
for any $k\in I$, and consider the collection $C_{i,W,J}=\{Z'_k\}_{k\in I}$.
Note that for any cofinal subset $J'$ of $J$ we have $C_{i,W,J'}=C_{i,W,J}$; in particular, for $k\leq l$, $Z'_l$ coincides with the $l$-th component of $C_{i,W,J\cap I_{\geq k}}$, and thus we have $p_{kl}(Z'_l)\subseteq Z'_k$. 
Moreover, by our assumption and the quasi-compactness of each $X_k$, the closed subsets $Z'_k$ are non-empty for any $k\in I$.
Hence we have $C_{i,W,J}\in\Phi$, and thus $C_{i,W,J}=z$ by the minimality of $z$.
Since the $i$-th component of $C_{i,W,J}$ is a closed subset of $Y$, we conclude $Y=Z_i$, as desired.

\medskip
{\sc Claim 2.} {\it Each $Z_i$ for $i\in I$ is irreducible.}

\medskip
Indeed, if $Z_i=W_1\cup W_2$ where $W_1,W_2$ are closed subsets, then at least one of the subsets $J_t=\{j\in I_{\geq i}\,|\,p^{-1}_{ij}(W_t)\cap Z_j\neq\emptyset\}$ $(t=1,2)$ is cofinal in $I_{\geq i}$; by {\sc Claim 1} we have $Z_i=W_1$ or $Z_i=W_2$.

One can similarly show that $\ovl{p_{ij}(Z_j)}=Z_i$ for $i\leq j$.
Since $X_i$ is sober, $Z_i$ admits a unique generic point $\eta_i$, and then the collection of points $\{\eta_i\}_{i\in I}$ determines a point in the limit $X=\varprojlim_{i\in I}X_i$.
Hence $X$ is non-empty, as desired.

Next, we show (1).
Let $\{Z_{\alpha}\}_{\alpha\in L}$ be a collection of closed subsets of $X$ such that $\bigcap_{\alpha\in L}Z_{\alpha}=\emptyset$.
We want to show that there exists a finite subset $L'\subseteq L$ such that $\bigcap_{\alpha\in L'}Z_{\alpha}=\emptyset$.
We may assume that each $Z_{\alpha}$ is of the form $p^{-1}_i(W)$ by a closed subset $W\subseteq X_i$ for some $i$ (where $p_i\colon X\rightarrow X_i$ is the canonical projection).
We choose by the axiom of choice a function $\alpha\mapsto (i(\alpha),W_{i(\alpha)})$ such that $Z_{\alpha}=p^{-1}_{i(\alpha)}(W_{i(\alpha)})$.
Now, for any $i\in I$, we set 
$$
L_i=\{\alpha\in L\,|\,i(\alpha)\leq i\},\qquad Z_i=\bigcap_{\alpha\in L_i}p^{-1}_{i(\alpha)\, i}(W_{i(\alpha)}).
$$
Then we have $p_{ij}(Z_j)\subseteq Z_i$ for any $i\leq j$, and hence $\{Z_i\}_{i\in I}$ forms a projective system such that $\varprojlim_{i\in I}Z_i=\bigcap_{\alpha\in L}Z_{\alpha}=\emptyset$.
Since each $Z_i$ is quasi-compact and sober (\ref{prop-propsoberness}), there exists $i\in I$ such that $Z_i=\emptyset$ (due to (2) proved above).
Since $X_i$ is quasi-compact, there exists a finite subset $L'\subseteq L_i$ such that $\bigcap_{\alpha\in L'_i}p^{-1}_{i(\alpha)\, i}(W_{i(\alpha)})=\emptyset$, and hence $\bigcap_{\alpha\in L'}Z_{\alpha}=\emptyset$, as desired.
\end{proof}
\end{small}

\subsubsection{Locally coherent spaces}\label{subsub-locallycoherentspaces}
\index{coherent!coherent topological space@--- (topological) space!locally coherent topological space@locally --- ---|(}\index{space@space (topological)!coherent topological space@coherent ---!locally coherent topological space@locally --- ---|(}
\begin{dfn}\label{dfn-locallycoherent}{\rm 
A topological space $X$ is {\em locally coherent} if $X$ admits an open covering by coherent subspaces.}
\end{dfn}

Since any coherent space admits an open basis by quasi-compact open subsets, and since quasi-compact open subsets of a coherent space is again coherent (\ref{prop-coherentclosedqcopen}), a locally coherent topological space admits an open basis consisting of coherent open subsets.
By this we deduce that, for a locally coherent space $X$ and a point $x\in X$, the set $G_x$ of all generizations of $x$ coincides with the intersection of all coherent open neighborhoods of $x$ (cf.\ \S\ref{subsub-genspetopsp}); in particular, if $X$ is sober, then $G_x$ is coherent by \ref{cor-projlimcohspacelattice1}.

Notice that by \ref{exas-coherenttopsp} (2) the underlying topological space of any scheme is locally coherent.

\begin{prop}\label{prop-locallycoherentbasicproperties}
Let $X$ be a topological space.
Then the following conditions are equivalent$:$
\begin{itemize}
\item[{\rm (a)}] $X$ is locally coherent$;$
\item[{\rm (b)}] every open subset of $X$ is locally coherent$;$
\item[{\rm (c)}] $X$ admits an open covering $X=\bigcup_{\alpha\in L}U_{\alpha}$ by locally coherent spaces$;$
\item[{\rm (d)}] the topology on $X$ is generated by coherent open subsets. \hfill$\square$
\end{itemize}
\end{prop}

The proof is straightforward and left to the reader.
The following proposition follows immediately from \ref{prop-coherentclosedqcopen}:
\begin{prop}\label{prop-locallycoherentbasicproperties2}
Any locally closed subspace of a locally coherent space is again locally coherent. \hfill$\square$
\end{prop}

In particular, any open subset of a coherent space is locally coherent.

\danger{For a locally coherent space $X$ to be coherent it is necessary and sufficient that $X$ is quasi-compact and quasi-separated.
But notice that a quasi-compact locally coherent space is not necessarily coherent.
Indeed, if $Y$ is a coherent space, and $U\subseteq Y$ is an open subset that is not quasi-compact, then $X=Y\amalg_UY$ (the gluing of two copies of $Y$ along $U$) is quasi-compact and locally coherent, but is not coherent.}

As the Stone duality (\S\ref{subsub-strcohtopsp}) indicates, quasi-compact maps give a good notion for morphisms between coherent topological spaces.
As for locally coherent spaces, similarly, we should have the following notion of maps:
\index{map@map (continuous)!quasi-compact map@quasi-compact ---!locally quasi-compact map@locally --- ---|(}\index{quasi-compact!quasi-compact map@--- map!locally quasi-compact map@locally --- ---|(}
\begin{dfn}\label{dfn-locallyquasicompactmaps}{\rm 
A continuous map $f\colon X\rightarrow Y$ between locally coherent spaces is said to be {\em locally quasi-compact} if, for any pair $(U,V)$ consisting of coherent open subsets $U\subseteq X$ and $V\subseteq Y$ with $f(U)\subseteq V$, the map $f|_U\colon U\rightarrow V$ is quasi-compact (\ref{dfn-quasicompactness} (2)).
}
\end{dfn}

The following facts are easy to see:
\begin{itemize}
\item A continuous map $f\colon X\rightarrow Y$ between coherent spaces is locally quasi-compact if and only if it is quasi-compact.
\item A coutinuous map $f\colon X\rightarrow Y$ between locally coherent spaces is locally quasi-compact if and only if there exist an open covering $\{V_{\alpha}\}_{\alpha\in L}$ of $Y$ by coherent open subsets and, for each $\alpha\in L$, an open covering $\{U_{\alpha,\lambda}\}_{\lambda\in\Lambda_{\alpha}}$ of $f^{-1}(V_{\alpha})$ by coherent open subsets, such that each $f|_{U_{\alpha,\lambda}}\colon U_{\alpha,\lambda}\rightarrow V_{\alpha}$ is quasi-compact.
\end{itemize}
\index{map@map (continuous)!quasi-compact map@quasi-compact ---!locally quasi-compact map@locally --- ---|)}\index{quasi-compact!quasi-compact map@--- map!locally quasi-compact map@locally --- ---|)}

\begin{prop}\label{prop-retrocompactpullback}
Let $f\colon X\rightarrow Y$ be a locally quasi-compact map between locally coherent spaces.
If $V\subseteq Y$ is a retrocompact\index{retrocompact} open subset of $Y$, then $f^{-1}(V)$ is retrocompact in $X$.
\end{prop}

\begin{proof}
It suffices to show that, for any coherent open subset $U\subseteq X$, $f^{-1}(V)\cap U$ is quasi-compact; here we use the fact that the topology of $X$ is generated by coherent open subsets.
Since $f(U)$ is quasi-compact, there exists a coherent open subset $W\subseteq Y$ such that $f(U)\subset W$.
Since $f|_U\colon U\rightarrow W$ is quasi-compact, $f^{-1}(V)\cap U=(f|_U)^{-1}(V\cap W)$ is quasi-compact.
\end{proof}

\begin{thm}\label{thm-ZRpoints61x}
Let $X$ be a locally coherent sober space, $Y$ a locally coherent space, and $f\colon X\rightarrow Y$ a quasi-compact continuous map.  Then we have
$$
\ovl{f(X)}=\bigcup_{y\in f(X)}\ovl{\{y\}}.
$$
In other words, the closure $\ovl{f(X)}$ is the set of all specializations of points of $f(X)$.
\end{thm}

\begin{proof}
The inclusion $\ovl{f(X)}\supseteq\bigcup_{y\in f(X)}\ovl{\{y\}}$ is clear. 
Let us show the converse. 
Considering an open covering of $Y$ by quasi-compact open subsets, one reduces to the case where $Y$ is quasi-compact.
In this case, as the map $f$ is quasi-compact, $X$ is also quasi-compact.
Hence there exists a finite open covering $X=\bigcup^n_{i=1}V_i$ by coherent open subsets. 
Replacing $X$ by the disjoint union $\coprod^n_{i=1}V_i$, we may assume that $X$ is coherent and sober.

For any $y\in\ovl{f(X)}$, let $\CN_y$ be the set of coherent open neighborhoods of $y$. 
We view $\CN_y$ as a directed set by the reversed inclusion order. 
For any $U\in\CN_y$, since $f$ is quasi-compact, the set $f^{-1}(U)$ is non-empty and coherent. 
Moreover, for $U \subseteq U'$, the inclusion $f^{-1}(U)\hookrightarrow f^{-1}(U')$ is quasi-compact. 
By \ref{thm-projlimcohsch1} (2) we have
$$
\varprojlim_{U\in\CN_y} f^{-1}(U)=f^{-1}(\bigcap_{U\in\CN_y}U)\neq\emptyset. 
$$
On the other hand, the set $\bigcap_{U\in\CN_y}U$ is the set $G_y$ of all generizations of $y$. 
This means that there is a generization of $y$ in $f(X)$, and hence the claim follows. 
\end{proof}

\begin{cor}\label{cor-ZRpoints61a}
Let $X$ be a locally coherent sober space, and $U\subseteq X$ a retrocompact\index{retrocompact} {\rm (\ref{dfn-retrocompact})} open subset of $X$.
Then we have
$$
\ovl{U}=\bigcup_{x\in U}\ovl{\{x\}}.
$$
In other words, the closure $\ovl{U}$ is the set of all specializations\index{specialization} of points in $U$. \hfill$\square$
\end{cor}

\begin{cor}\label{cor-propZRpoints61x001prop}
In the situation as in {\rm (\ref{thm-ZRpoints61x})}, the following conditions are equivalent$:$
\begin{itemize}
\item[{\rm (a)}] $f$ is a closed map$;$ 
\item[{\rm (b)}] $f(\ovl{\{x\}})=\ovl{\{f(x)\}}$ for any $x\in X$.
\end{itemize}
\end{cor}

\begin{proof}
(a) $\Rightarrow$ (b) is clear.
Suppose (b) holds, and let $C$ be a closed subset of $X$.
We have $\ovl{f(C)}=\bigcup_{x\in C}\ovl{\{f(x)\}}=\bigcup_{x\in C}f(\ovl{\{x\}})$ by \ref{thm-ZRpoints61x}.
Since $\ovl{\{x\}}\subseteq C$, one has $\ovl{f(C)}\subseteq f(C)$.
\end{proof}
\index{space@space (topological)!coherent topological space@coherent ---!locally coherent topological space@locally --- ---|)}\index{coherent!coherent topological space@--- (topological) space!locally coherent topological space@locally --- ---|)}

\subsection{Valuative spaces}\label{sub-valuativespace}
\index{valuative!valuative topological space@--- (topological) space|(}\index{space@space (topological)!valuative topological space@valuative ---|(}
\subsubsection{Valuative spaces}\label{subsub-valspdef}
\begin{dfn}\label{dfn-valuativespace}{\rm 
A topological space $X$ is said to be {\em valuative} if the following conditions are satisfied$:$
\begin{itemize}
\item[{\rm (a)}] $X$ is locally coherent {\rm (\ref{dfn-locallycoherent})}\index{coherent!coherent topological space@--- (topological) space!locally coherent topological space@locally --- ---}\index{space@space (topological)!coherent topological space@coherent ---!locally coherent topological space@locally --- ---} and sober {\rm (\S\ref{subsub-sober})}\index{space@space (topological)!sober topological space@sober ---}$;$
\item[{\rm (b)}] for any point $x\in X$ the ordered set $G_x$ of all generizations of $x$ $(${\rm \S\ref{subsub-sober}}$)$ is totally ordered\index{ordered!totally@totally ---}$;$
\item[{\rm (c)}] every point $x\in X$ has a maximal generization $\til{x}\in G_x$.
\end{itemize}}
\end{dfn}

Note that, in view of (b), the maximal generization of $x$ as in (c) is uniquely determined.
The order type\index{order type} (\S\ref{subsub-orderings}) of the totally ordered set $G_x$ is called the {\em height}\index{height!height of a point@--- (of a point)} of the point $x$.
\begin{rem}\label{rem-valuativespace}{\rm 
(1) Notice that, under the axiom of choice, the condition (c) follows automatically from (a).
Indeed, to show this, we may assume that $X$ is coherent and sober, hence is of the form $X=\Spec A$ for a distributive lattice\index{lattice}\index{lattice!distributive lattice@distributive ---} $A=\QCOuv(X)$ by Stone duality (\S\ref{subsub-strcohtopsp}).
In this situation, points of $G_x$ correspond to prime filters\index{filter (of lattices)!prime filter of lattices@prime ---} $F_y=\{U\in A\,|\,y\in U\}$ $(y\in G_x)$.
Any totally ordered subset $I$ of $G_x$ admits an upper bound $\bigcup_{F\in I}F$ in $G_x$; here, notice that $y\in X$ is a generization of $z\in X$ if and only if $F_y\supseteq F_z$.
Hence Zorn's lemma implies that there exists a maximal element in $G_x$.

(2) It follows that any locally coherent and sober subspace of a valuative space is again valuative.
In particular, by \ref{prop-locallycoherentbasicproperties2} and \ref{prop-propsoberness}, any locally closed subspace of a valuative space is a valuative space.}
\end{rem}

The following proposition follows easily from \ref{prop-locallycoherentbasicproperties}:
\begin{prop}\label{proop-valuativespacedef}
Let $X$ be a topological space. 
Then the following conditions are equivalent$:$
\begin{itemize}
\item[{\rm (a)}] $X$ is valuative$;$
\item[{\rm (b)}] every open subset of $X$ is valuative$;$
\item[{\rm (c)}] $X$ admits an open covering by valuative spaces.\hfill$\square$
\end{itemize}
\end{prop}

An important example of valuative spaces, which motivates the terminology `valuative', is the underlying topological spaces of the spectrum $\Spec V$ of a valuation ring $V$; cf.\ \ref{prop-spectrumval}.
In this case, the height of the closed point of $\Spec V$ is nothing but the height of the valuation ring $V$ (\ref{dfn-height1} (1)).

\subsubsection{Closures and tubes}\label{subsub-tubesgen}
\index{tube|(}
\begin{dfn}\label{dfn-tubes1a}{\rm 
Let $X$ be a valuative space. 
A {\em tube closed subset}\index{tube!tube closed subset@--- closed subset} of $X$ is a closed subset of the form $\ovl{U}$ by a retrocompact\index{retrocompact} {\rm (\ref{dfn-retrocompact})} open subset $U$.
The complement of a tube closed subset is called a {\em tube open subset}\index{tube!tube open subset@--- open subset}.
Tube closed and tube open subsets are collectively called {\em tube subsets}\index{tube!tube subset@--- subset}.}
\end{dfn}

The following proposition is an immediate consequence of \ref{cor-ZRpoints61a}:
\begin{prop}\label{prop-tubes11112a}
Let $X$ be a valuative space, and $U\subseteq X$ a retrocompact\index{retrocompact} open subset.
Then for $x\in X$ to belong to the tube closed subset\index{tube!tube closed subset@--- closed subset} $\ovl{U}$ it is necessary and sufficient that the maximal generization $\til{x}$ of $x$ belongs to $U$. \hfill$\square$
\end{prop}

\begin{cor}\label{cor-tubes11111a}
Let $X$ be a valuative space, and $C=(X\setminus U)^{\circ}$ $($where $U\subseteq X$ is retrocompact open$)$ a tube open subset\index{tube!tube open subset@--- open subset}.
Then for $x\in X$ to belong to $C$ it is necessary and sufficient that the maximal generization $\til{x}$ of $x$ does not belong to $U$. \hfill$\square$
\end{cor}

Now, recall the following general-topology fact: For two open subsets $U_1,U_2\subseteq X$, we have
$$
\ovl{U_1\cup U_2}=\ovl{U}_1\cup\ovl{U}_2,
$$
but not in general for the intersections; in our situation, however, we have:
\begin{prop}\label{prop-tubes121a} 
Let $X$ be a valuative space, and $\{U_{\alpha}\}_{\alpha\in L}$ a family of retrocompact\index{retrocompact} open subsets of $X$. 
Then we have
$$
\ovl{\bigcap_{\alpha\in L}U_{\alpha}}=\bigcap_{\alpha\in L}\ovl{U}_{\alpha}.
$$
\end{prop}

\begin{proof}
The inclusion $\ovl{\bigcap_{\alpha\in L}U_{\alpha}}\subseteq\bigcap_{\alpha\in L}\ovl{U}_{\alpha}$ is clear.  
To show the converse, take any point $y\in\bigcap_{\alpha\in L}\ovl{U}_{\alpha}$. 
Then by \ref{prop-tubes11112a} the maximal generization $\til{y}$ of $y$ belongs to $U_{\alpha}$ for any $\alpha\in L$. 
Hence we have $y\in\ovl{\{\til{y}\}}\subseteq\ovl{\bigcap_{\alpha\in L}U_{\alpha}}$, as desired.
\end{proof}

\begin{cor}\label{cor-proptubes10a}
{\rm (1)} Any finite union of tube closed $($resp.\ tube open$)$ subsets\index{tube!tube closed subset@--- closed subset}\index{tube!tube open subset@--- open subset} is a tube closed $($resp.\ tube open$)$ subset.

{\rm (2)} Any finite intersection of tube closed $($resp.\ tube open$)$ subsets\index{tube!tube closed subset@--- closed subset}\index{tube!tube open subset@--- open subset} is a tube closed $($resp.\ tube open$)$ subset.
\end{cor}

\begin{proof}
In view of \ref{prop-tubes121a}, the closure operator $\ovl{\,\cdot\,}$ commutes with finite intersection and finite union of retrocompact open subsets.
Hence, to show the corollary, it suffices to show that the finite intersection and the finite union of retrocompact\index{retrocompact} open subsets are retrocompact, which is clear.
\end{proof}
\index{tube|)}

\subsubsection{Separated quotients and separation maps}\label{subsub-separationgen}
Let $X$ be a valuative space, and let $[X]$ denote the subset of $X$ consisting of all maximal points\index{point!maximal point@maximal ---} of $X$. 
We have the canonical retraction map 
$$
\sep_{X}\colon X\longrightarrow[X],\qquad x\longmapsto\textrm{the maximal generization\index{generization!maximal generization@maximal ---} of $x$},
$$
which we call the {\em separation map}\index{map@map (continuous)!separation map@separation ---}\index{separation map}.
The separation map $\sep_X$ is clearly surjective.
We endow $[X]$ with the quotient topology induced from the topology on $X$.
Then the space $[X]$ is a $\mathrm{T}_1$-space (Exercise \ref{exer-valuativespaceT1}).
The topological space $[X]$ thus obtained is called the {\em separated quotient}\index{separated quotient} (or {\em $\mathrm{T}_1$-quotient}) of $X$.

\danger{Notice that, with this topology, the inclusion map $[X]\hookrightarrow X$ is not continuous in general.}

\begin{prop}[Universality of separated quotients]\label{prop-separation3}
Let $X$ be a valuative space.
Suppose we are given a continuous mapping $\varphi\colon X \rightarrow T$ to a $\mathrm{T}_1$-space $T$.
Then there exists a unique continuous map $\psi\colon[X]\rightarrow T$ such that the resulting diagram
$$
\xymatrix@-1ex{
X \ar[r]^{\varphi}\ar[d]_{\sep_{X}}&T\\
[X]\ar[ur]_{\psi}}
$$
commutes.
\end{prop}

\begin{proof}
Let $x\in X$, and $y$ the maximal generization of $x$.
We need to show $\varphi(x)=\varphi(y)$.
Since $T$ is $\mathrm{T}_1$, the subset $\varphi^{-1}(\varphi(y))$ is closed and hence contains $\ovl{\{y\}}$.
But this means $x\in\varphi^{-1}(\varphi(y))$, that is, $\varphi(x)=\varphi(y)$.
\end{proof}

\begin{cor}[Functoriality]\label{cor-functorialitysep1a}
Any continuous map $ f\colon X \rightarrow  Y $ between valuative spaces induces a unique continuous map $[f]\colon [X]\rightarrow [Y]$ such that diagram
$$
\xymatrix@-1ex{X\ar[r]^f\ar[d]_{\sep_X}&Y\ar[d]^{\sep_Y}\\ [X]\ar[r]_{[f]}&[Y]}
$$
commutes. $($Hence, in particular, the formation $X\mapsto [X]$ is functorial.$)$ \hfill$\square$
\end{cor}

\subsubsection{Overconvergent sets}\label{subsub-overconvsets}
\begin{dfn}\label{dfn-separation2a}{\rm 
Let $X$ be a valuative space.
A closed (resp.\ an open) subset $S$ of $X$ is said to be {\em overconvergent}\index{overconvergent!overconvergent subset@--- subset} if for any $x\in S$ any generization (resp.\ specialization) of $x$ belongs to $S$.}
\end{dfn}

For example, if $x$ is a maximal point of $X$, then $\ovl{\{x\}}$ is an overconvergent closed subset.
Notice that, if $S\subseteq X$ is overconvergent, then so is the complement $X\setminus S$.
Notice also that, if $S$ is an overconvergent closed or open subset and $x\in S$, then both $G_x$ (the set of all generizations of $x$) and $\ovl{\{x\}}$ (the set of all specializations of $x$) are contained in $S$.

The following propositions are easy to see:
\begin{prop}\label{prop-overconvergentbycovering}
Let $X$ be a valuative space, and $X=\bigcup_{\alpha\in L}U_{\alpha}$ an open covering. 
Then a subset $S\subseteq X$ is an overconvergent closed $($resp.\ open$)$ subset if and only if $S\cap U_{\alpha}$ is an overconvergent closed $($resp.\ open$)$ subset of $U_{\alpha}$ for every $\alpha\in L$. \hfill$\square$
\end{prop}

\begin{prop}\label{prop-separation21a}
Let $X$ be a valuative space, and $S$ a closed or an open subset of $X$.
Then $S$ is overconvergent if and only if $S=\sep^{-1}_{X}(\sep_{X}(S))$.
Hence, in particular, $T\mapsto\sep^{-1}_X(T)$ gives a bijection between the set of all open $($resp.\ closed$)$ subsets of $[X]$ and the set of all overconvergent open $($resp.\ closed$)$ subsets of $X$. \hfill$\square$
\end{prop}

\begin{cor}\label{cor-separation211a}
{\rm (1)} Any finite intersection of overconvergent open subsets is an overconvergent open subset.
The union of arbitrarily many overconvergent open subsets is an overconvergent open subset.

{\rm (2)} Any finite union of overconvergent closed subsets is an overconvergent closed subset.
The intersection of arbitrarily many overconvergent closed subsets is an overconvergent closed subset. \hfill$\square$
\end{cor}

By \ref{prop-tubes11112a} and \ref{cor-tubes11111a}, we have:
\begin{prop}\label{prop-separation22a}
Any tube closed $($resp.\ tube open$)$ subset\index{tube!tube closed subset@--- closed subset}\index{tube!tube open subset@--- open subset} of a valuative space $X$ is overconvergent closed $($resp.\ overconvergent open$)$ in $X$. \hfill$\square$
\end{prop}

\begin{dfn}\label{dfn-tubes111a}{\rm 
Let $X$ be a valuative space, and consider the separated quotient $[X]$ of $X$.
A subset $T$ of $[X]$ is said to be a {\em tube closed} $($resp.\ {\em tube open}$)$ {\em subset}\index{tube!tube closed subset@--- closed subset}\index{tube!tube open subset@--- open subset} if $\sep_{X}^{-1}(T)$ is a tube closed $($resp.\ tube open$)$ subset of $X$ $(\ref{dfn-tubes1a})$.}
\end{dfn}

Thus, by \ref{prop-separation21a}, $S\mapsto\sep_X(S)$ gives a canonical order preserving bijection between the sets of all tube closed (resp.\ tube open) subsets of $X$ and of tube closed (resp.\ tube open) subsets of $[X]$.

\begin{prop}\label{prop-tubes2a} 
Let $ X$ be a coherent valuative space. 

{\rm (1)} For any overconvergent closed set $F$, the set of all tube open subset\index{tube!tube open subset@--- open subset}s containing $F$ form a fundamental system of neighborhoods of $F$. 

{\rm (2)} For overconvergent closed subsets $F_1,F_2$ with $F_1\cap F_2=\emptyset$, there exist tube open subsets $U_1$ and $U_2$ such that $F_i\subseteq U_i$ $(i=1,2)$ and $U_1\cap U_2=\emptyset$.
\end{prop}

\begin{proof}
(1) Take a quasi-compact open neighborhood $U$ of $F$. 
We want to find a tube open subset $C$ such that $F\subseteq C\subseteq U$. 
Take a quasi-compact open set $V$ such that  $X\setminus U\subseteq V\subseteq X\setminus F$. 
Since $F$ is overconvergent, $\overline{V}\subseteq X\setminus F$ by \ref{prop-separation22a}, and $U$ and $\overline{V}$ cover $X$. 
This implies $F\subseteq X\setminus\overline{V}\subseteq U$. 

(2) Since $X=(X\setminus F_1)\cup(X\setminus F_2)$, there exist quasi-compact open subsets $U'_i$ ($i=1,2$) such that $U'_1\cup U'_2=X$ and $U'_i\cap F_i=\emptyset$.
As $F_i$ is overconvergent, we have $\ovl{U'}_i\cap F_i=\emptyset$ for $i=1,2$ by \ref{cor-ZRpoints61a}.
Set $U_i=X \setminus\ovl{U'}_i$ for $i=1,2$, which are tube open subsets.
We have $F_i\subseteq U_i$ $(i=1,2)$ and $U_1\cap U_2=\emptyset$, as desired.
\end{proof}

\begin{cor}\label{cor-tubes2acor}
Let $X$ be a coherent valuative space.
\begin{itemize}
\item[{\rm (1)}] The separated quotient $[X]$ is a normal topological space\index{space@space (topological)!normal topological space@normal ---}. 
\item[{\rm (2)}] The space $[X]$ is compact $($and hence is Hausdorff$)$. 
\end{itemize}
\end{cor}

\begin{proof} (1) follows from \ref{prop-tubes2a} (2), for $[X]$ is a $\mathrm{T}_1$-space (Exercise \ref{exer-valuativespaceT1}). 
(2) follows from (1) and the quasi-compactness of $X$. 
\end{proof}

\begin{cor}\label{cor-tubes211a}
Let $X$ be a coherent valuative space. 
Then the set of all tube open subsets\index{tube!tube open subset@--- open subset} of $[X]$ forms an open basis of the topological space $[X]$.
\end{cor}

\begin{proof}
This follows from \ref{prop-tubes2a} (1) and the fact that, as $[X]$ is $\mathrm{T}_1$, any singleton set $\{x\}$ for $x\in[X]$ is a closed subset.
\end{proof}

Recall here that a continuous map $f\colon X\rightarrow Y$ between topological spaces is said to be {\em proper} if for any topological space $Z$ the induced mapping $f\times\id_Z\colon X\times Z\rightarrow Y\times Z$ is closed (cf.\ \cite[Chap.\ I, \S10.1, Def.\ 1]{Bourb4}).

\begin{cor}\label{cor-tubes211acor}
Let $X$ be a coherent valuative space. 
Then $\sep_{X}$ is a proper map of topological spaces.
\end{cor}

\begin{proof}
By \cite[Chap.\ I, \S10.2, Cor.\ 2]{Bourb4} every continuous mapping from a quasi-compact space to a Hausdorff space is proper\index{map@map (continuous)!proper map@proper ---}.
Hence the corollary follows from \ref{cor-tubes2acor} (2).
\end{proof}

\subsubsection{Valuative maps}\label{subsub-valuativemaps}
\begin{dfn}\label{dfn-valuativemaps}{\rm 
A continuous map $f\colon X\rightarrow Y$ between valuative spaces is said to be {\em valuative}\index{valuative!valuative map@--- map} if $f([X])\subset [Y]$, that is, $f$ maps maximal points to maximal points.}
\end{dfn}

For example, open immersions between valuative spaces are valuative.
Clearly, if $f\colon X\rightarrow Y$ is valuative, then $[f](x)=f(x)$ for any $x\in [X]$ ($\subseteq X$). 

\begin{prop}\label{prop-valuativemapstubes2}
Let $f\colon X\rightarrow Y$ be a valuative map between valuative spaces.

{\rm (1)} For any subset $C\subseteq Y$ that is stable under generization, $f^{-1}(C)$ is stable under generization in $X$, and the diagram
$$
\xymatrix{f^{-1}(C)\cap[X]\ar[r]\ar[d]&[X]\ar[d]^{[f]}\\ C\cap[Y]\ar[r]&[Y]}
$$
is Cartesian. 

{\rm (2)} We have $f^{-1}(\sep^{-1}_Y(C\cap[Y]))=\sep^{-1}_X(f^{-1}(C)\cap[X])$.
\end{prop}

\begin{proof}
(1) is obvious by the definition of valuative maps.
To show (2), first observe that the set $F=f^{-1}(\sep^{-1}_Y(C\cap[Y]))$ is stable under specialization and generization in $X$ and that $F\cap[X]=[f]^{-1}(C\cap[Y])$.
By (1) we have $F\cap[X]=f^{-1}(C)\cap[X]$, that is, $F$ coincides with the set of all specializations of points in $f^{-1}(C)\cap[X]$, whence the desired equality.
\end{proof}

\begin{cor}\label{cor-valuativemapstubes3}
Let $f\colon X\rightarrow Y$ be a valuative map between valuative spaces.
Then, for any overconvergent closed $($resp.\ overconvergent open$)$ subset $C\subseteq Y$, $f^{-1}(C)$ is overconvergent closed $($resp.\ overconvergent open$)$ in $X$. \hfill$\square$
\end{cor}

\begin{cor}\label{cor-valuativemapstubes2}
Let $f\colon X\rightarrow Y$ be a locally quasi-compact\index{map@map (continuous)!quasi-compact map@quasi-compact ---!locally quasi-compact map@locally --- ---}\index{quasi-compact!quasi-compact map@--- map!locally quasi-compact map@locally --- ---} {\rm (\ref{dfn-locallyquasicompactmaps})} valuative map between valuative spaces.

{\rm (1)} For any retrocompact\index{retrocompact} open subset $U$ of $Y$, we have
$$
f^{-1}(\ovl{U})=\ovl{f^{-1}(U)}.
$$

{\rm (2)} For any tube closed\index{tube!tube closed subset@--- closed subset} $($resp.\ tube open$)$\index{tube!tube open subset@--- open subset} {\rm (\ref{dfn-tubes1a})} subset $S$ of $Y$, $f^{-1}(S)$ is a tube closed $($resp.\ tube open$)$ subset of $X$.
\end{cor}

\begin{proof}
(1) Since $U$ and $f^{-1}(U)$ are retrocompact by \ref{prop-retrocompactpullback}, it follows from \ref{cor-ZRpoints61a} that $\ovl{U}=\sep^{-1}_Y([U])$ and $\ovl{f^{-1}(U)}=\sep^{-1}_X([f^{-1}(U)])$.
Then the desired equality follows from \ref{prop-valuativemapstubes2} (2).
(2) follows immediately from (1).
\end{proof}

\subsubsection{Structure of separated quotients}\label{subsub-strsepquot}
Let $X$ be a valuative space, and $U\subseteq X$ an open subset.
Then $U$ is again a valuative space (\ref{proop-valuativespacedef}), and $U\hookrightarrow X$ induces a continuous injection $[U]\hookrightarrow [X]$ with the image $\sep_X(U)$ (\ref{cor-functorialitysep1a}).

\danger{The continuous bijection $[U]\rightarrow\sep_X(U)$, where $\sep_X(U)$ is endowed with the subspace topology from $[X]$, may not be a homeomorphism; it is homeomorphic if, for example, $U$ is overconvergent, or $X$ and $U$ are coherent; cf.\ Exercise \ref{exer-locallystronglycompact2}.}

\danger{If $U\subseteq X$ is overconvergent open, then $[U]=\sep_X(U)$ is open in $[X]$.
But, in general, $\sep_X(U)$ for an open subset $U\subset X$ may not be open in $[X]$.
Note that, if $X$ is coherent and $U$ is quasi-compact, then $[U]$ is compact in the Hausdorff space $[X]$ (\ref{cor-tubes2acor} (2)), and hence is closed.}

\begin{prop}[{Continuity of $[\,\cdot\,]$}]\label{prop-structuresepquot2}
Let $\mathscr{C}$ be a $\mathsf{U}$-small\index{small!Usmall@$\mathsf{U}$-{---}} category, and $F\colon\mathscr{C}\rightarrow\Top$ a functor.
Suppose that $F(X)$ is a valuative space for any object $X$ in $\mathscr{C}$ and that $F(f)$ is an open immersion for any morphism $f$ in $\mathscr{C}$. 
\begin{itemize}
\item[{\rm (1)}] The colimit\index{limit!colimit@co---}\index{colimit} $X=\varinjlim_{\mathscr{C}}F$ is representable by a valuative space.
\item[{\rm (2)}] $[X]\cong\varinjlim_{\mathscr{C}}([\,\cdot\,]\circ F)$. 
\end{itemize}
\end{prop}

\begin{proof}
(1) is clear.
To show (2), observe first that the limit $Q=\varinjlim_{\mathscr{C}}[F(Z)]$ is a $\mathrm{T}_1$-space (easy to see).
The maps $F(x)\rightarrow X$ for $x\in\obj(\mathscr{C}$) induce the maps $[F(x)]\rightarrow[X]$, and hence the canonical map $\alpha\colon Q\rightarrow[X]$.
On the other hand, the composition $F(x)\rightarrow[F(x)]\rightarrow Q$ for $x\in\obj(\mathscr{C})$ give rise to $X\rightarrow Q$.
As $Q$ is a $\mathrm{T}_1$-space, we have $\beta\colon[X]\rightarrow Q$ by the universality of separated quotients (\ref{prop-separation3}).
It is easy to see that $\alpha$ and $\beta$ are inverse to each other.
\end{proof}

Let $X$ be a valuative space, and $\{U_{\alpha}\}_{\alpha\in L}$ an open covering of $X$.
Consider the coequalizer sequence
$$
\xymatrix@-1ex{
R \ar@<.5ex>[r]\ar@<-.5ex>[r]&\coprod_{\alpha\in L} U_{\alpha}\ar[r]&X,}
$$
where $R=\coprod_{\alpha,\beta\in L}U_{\alpha}\cap U_{\beta}$, and the induced sequence
$$
\xymatrix@-1ex{
[R] \ar@<.5ex>[r]\ar@<-.5ex>[r]&\coprod_{\alpha\in L}[U_{\alpha}]\ar[r]&[X]}; 
$$
notice that the functor $X\mapsto [X]$ commutes with disjoint union and with finite intersection of open subsets.
One sees easily that, in fact, $[R]$ defines an equivalence relation on $\coprod_{\alpha\in L}[U_{\alpha}]$ (cf.\ Exercise \ref{exer-equivalencerelationinvolution}).
By \ref{prop-structuresepquot2} we immediately have:
\begin{cor}\label{cor-propstructuresepquot}
The topological space $[X]$ is the quotient of $\coprod_{\alpha\in L}[U_{\alpha}]$ by the equivalence relation $[R]$. \hfill$\square$
\end{cor}

\begin{prop}\label{prop_structuresepquot3}
Let $f\colon X\rightarrow Y$ be a quasi-compact valuative map of valuative spaces. 
Then the induced map $[f]$ is proper\index{map@map (continuous)!proper map@proper ---} as a map of topological spaces.
\end{prop}

\begin{proof} 
First we assume $Y$ is coherent. 
Then $[X]$ is a quasi-compact space.
Since $[Y]$ is Hausdorff (\ref{cor-tubes2acor} (2)), $[f]$ is proper by \cite[Chap.\ I, \S10.2, Cor.\ 2]{Bourb4}, and the claim is shown in this case.  

In general, we take an open covering $\{U_{\alpha}\}_{\alpha\in L}$ of $Y$ by coherent spaces.
We claim that $[f]$ is a closed map.
Let $F\subseteq[X]$ be a closed subset.
By \ref{prop-overconvergentbycovering}, to show that $[f](F)$ is closed in $[Y]$, it suffices to show that $[f](F)\cap[U_{\alpha}]$ is closed in $[U_{\alpha}]$ for each $\alpha\in L$.
To show this, first observe that, since $f$ is valuative, we have the equality $[f]^{-1}([U_{\alpha}])=[f^{-1}(U_{\alpha})]$; similarly, one has $[f](F)\cap[U_{\alpha}]=[f](F\cap[f^{-1}(U_{\alpha})])$.
Since $[f]\vert_{[f^{-1}(U_{\alpha})]}$ is the separated quotient of $f\vert_{f^{-1}(U_{\alpha})}\colon f^{-1}(U_{\alpha})\rightarrow U_{\alpha}$, it follows from the coherent case that $[f](F)\cap[U_{\alpha}]$ is closed in $[U_{\alpha}]$ for any $\alpha\in L$, as desired.

To finish the proof, in view of \cite[Chap.\ I, \S10.2, Thm.\ 2]{Bourb4}, it suffices to show that $[f]^{-1}(y)$ is quasi-compact for any $y\in [Y]$.
But this is reduced to showing the properness of $[f]\vert_{[f^{-1}(U_{\alpha})]}$ with $y\in [U_{\alpha}]$, since $\{y\}$ is closed both in $[Y]$ and $[U_{\alpha}]$, and hence follows again from the coherent case. 
\end{proof}

\subsubsection{Overconvergent interior}\label{subsub-gentopoverconvint}
\begin{dfn}\label{dfn-interior1a}{\rm 
Let $X$ be a valuative space, and $F$ a subset of $X$.
The maximal overconvergent open subset in $X$ contained in $F$, of which the existence is guaranteed by \ref{cor-separation211a} (1), is called the {\em overconvergent interior}\index{overconvergent!overconvergent interior@--- interior} of $F$ in $X$ and denoted by $\int_{X}(F)$.}
\end{dfn}

\begin{prop}\label{prop-interior2a} 
Let $X$ be a valuative space, and $U\subseteq X$ a quasi-compact open subset.
Suppose that the following condition is satisfied$:$
\begin{itemize}
\item[$(\ast)$] there exists a coherent open subset $V\subseteq X$ such that $\ovl{U}\subseteq V$. 
\end{itemize}
Let $y\in X$ be a maximal point.
Then $y\in\int_{X}(U)$ if and only if $\ovl{\{y\}}\subseteq U$. 
\end{prop}

Notice that the condition $(\ast)$ is automatic if $X$ is quasi-separated and locally strongly compact\index{valuative!valuative topological space@--- (topological) space!locally strongly compact valuative topological space@locally strongly compact --- ---} (\ref{dfn-locallycompactspace}); cf.\ \ref{prop-locstringcompsp2a} below.

\begin{proof}
Suppose $\ovl{\{y\}}\subseteq U$. 
Since $\ovl{\{y\}}$ is overconvergent, there exists by \ref{prop-tubes2a} (1) an open subset $W\subseteq V$, overconvergent in $V$, such that $\ovl{\{y\}}\subseteq W\subseteq U$.
Since $W\subseteq\ovl{U}\subseteq V$, $W$ is also overconvergent in $X$.
Hence we have $y\in\int_{X}(U)$, as desired.
The converse is clear. 
\end{proof}

\begin{cor}\label{cor-interior3a} 
In the situation as in {\rm \ref{prop-interior2a}} we have
$$
\int_{X}(U)=\sep^{-1}_{X}([U]\setminus\sep_{X}(\partial U)),
$$
where $\partial U=\ovl{U}\setminus U$. 
\end{cor}

\begin{proof}
A point $y\in X$ lies in $\sep^{-1}_{X}([U]\setminus\sep_{X}(\partial U))$ if and only if the maximal generization $\til{y}$ of $y$ satisfies $\til{y}\in U$ and $\ovl{\{\til{y}\}}\cap\partial U=\emptyset$. 
By \ref{prop-interior2a} this is equivalent to $y\in\int_{X}(U)$.
\end{proof}

\begin{cor}\label{cor-interior4a}
Let $X$ be a quasi-separated valuative space, $U$ a quasi-compact open subset of $X$, and $y\in X$ a maximal point.
Assume $\ovl{U}$ is quasi-compact. 
Then $y\in\int_{X}(U)$ if and only if $\ovl{\{y\}}\subseteq U$. 
\end{cor}

\begin{proof}
Since $X$ is quasi-separated, any quasi-compact open subset of $X$ is coherent.
Since $\ovl{U}$ is quasi-compact, there exists a coherent open subset $V$ such that $\ovl{U}\subseteq V$.
Hence the assertion follows from \ref{prop-interior2a}.
\end{proof}

\subsection{Reflexive valuative spaces}\label{sub-reflexivevaluativespaces}
\index{valuative!valuative topological space@--- (topological) space!reflexive valuative topological space@reflexive --- ---|(}
\subsubsection{Reflexive valuative spaces}\label{subsub-reflexivevaluativespaces}
\begin{dfn}\label{dfn-reflexivevaluativespaces}{\rm 
A valuative space (\ref{dfn-valuativespace}) $X$ is said to be {\em reflexive} if the following condition holds$:$ for two coherent open subsets $U\subseteq V$ of $X$, $[U]=[V]$ implies $U=V$.}
\end{dfn}

\begin{exa}\label{exa-reflexivevaluativespaces1}{\rm 
Let $V$ be an $a$-adically complete valuation ring\index{valuation!valuation ring@--- ring!a-adically complete valuation ring@$a$-adically complete --- ---} $(a\in\m_V\setminus\{0\})$.
Then the formal spectrum\index{formal spectrum} $\Spf V$ (cf.\ {\bf \ref{ch-formal}}, \S\ref{subsub-formalnotformalspec}) of $V$ is a valuative space with the unique maximal point $\mathfrak{p}_V=\sqrt{(a)}$, the associated height one prime\index{associated height one prime} (cf.\ \ref{dfn-maxspe2}).
It is reflexive if and only if $V$ is of height\index{height!height of a valuation ring@--- (of a valuation (ring))} one (cf.\ \ref{dfn-height1} (1)).}
\end{exa}

Notice that the reflexiveness is a local property on $X$. Indeed, we have:
\begin{prop}\label{prop-reflexivevaluativespaces1}
A valuative space $X$ is reflexive if and only if, for any pair of open subsets $U\subseteq V$ of $X$ such that the inclusion $U\hookrightarrow V$ is quasi-compact, $[U]=[V]$ implies $U=V$.
\end{prop}

\begin{proof}
The `if' part is clear.
Suppose $X$ is reflexive, and let $U\hookrightarrow V$ be as above.
For any coherent open subset $W$ of $V$, $U\cap W$ is coherent.
Since $[U\cap W]=[U]\cap[W]=[V]\cap[W]=[W]$, we have $U\cap W=W$.
Since $V$ has an open basis consisting of coherent open subsets, this implies that $U=V$.
\end{proof}

\begin{prop}\label{prop-lemlocallystronglycompactconverse}
Let $X$ be a reflexive valuative space.
Then any retrocompact\index{retrocompact} {\rm (\ref{dfn-retrocompact})} open subset $U\subseteq X$ is regular, that is, $(\ovl{U})^{\circ}=U$. 
\end{prop}

\begin{proof}
Considering an open covering of $X$ by coherent open subsets, one can easily reduce to the case where $X$ and $U$ are coherent.
Since $U\subseteq(\ovl{U})^{\circ}$ trivially holds, we only need to show $U\supseteq(\ovl{U})^{\circ}$.
By \ref{cor-ZRpoints61a}, we have $\sep^{-1}_X([U])=\ovl{U}$.
In particular, we have $[U]=[(\ovl{U})^{\circ}]$.
For any quasi-compact open subset $V\subseteq(\ovl{U})^{\circ}$, $U\cap V$ is coherent, and we have $[U\cap V]=[U]\cap[V]=[(\ovl{U})^{\circ}]\cap[V]=[V]$.
Hence $U\cap V=V$.
Since this holds for any quasi-compact open subsets of $(\ovl{U})^{\circ}$, we have $U=(\ovl{U})^{\circ}$, as desired.
\end{proof}

\begin{prop}\label{prop-lemlocallyquasicompactmaps}
Let $X,Y$ be coherent valuative spaces, and $f,g\colon X\rightarrow Y$ two valuative quasi-compact maps.
Suppose $X$ is reflexive.
Then $[f]=[g]$ implies $f=g$.
\end{prop}

\begin{proof}
In view of the Stone duality (\ref{thm-stonerepresentationthm}), it suffices to show that, for any quasi-compact open subset $V\subseteq Y$, we have $f^{-1}(V)=g^{-1}(V)$.
By \ref{cor-valuativemapstubes2} (1), we have $\ovl{f^{-1}(V)}=f^{-1}(\ovl{V})$ and $\ovl{g^{-1}(V)}=g^{-1}(\ovl{V})$.
Since we have $\ovl{V}=\sep^{-1}_Y(V)$, and since $\sep_Y\circ f=\sep_Y\circ g$ by the assumption, we have $\ovl{f^{-1}(V)}=\ovl{g^{-1}(V)}$.
Then by \ref{prop-lemlocallystronglycompactconverse}
we have $f^{-1}(V)=g^{-1}(V)$, as desired.
\end{proof}

\subsubsection{Reflexivization}\label{subsub-reflexivization}
\index{reflexivization|(}
Let us denote by $\Vsp$ the category of valuative spaces and valuative and locally quasi-compact maps, and by $\RVsp$ the full subcategory of $\Vsp$ consisting of reflexive valuative spaces.
\begin{thm}\label{thm-reflexivization}
The canonical inclusion $I\colon\RVsp\hookrightarrow\Vsp$ admits a left adjoint functor $(\,\cdot\,)^{\mathrm{ref}}\colon\Vsp\rightarrow\RVsp$, which has the following properties$:$
\begin{itemize}
\item[{\rm (a)}] $(\,\cdot\,)^{\mathrm{ref}}$ preserves open immersions$;$
\item[{\rm (b)}] the adjunction map $X^{\mathrm{ref}}\rightarrow X$ induces a homeomorphism $[X^{\mathrm{ref}}]\cong[X];$
\item[{\rm (c)}] the adjunction map $X^{\mathrm{ref}}\rightarrow X$ is quasi-compact$;$
\item[{\rm (d)}] if $X$ is quasi-separated, then so is $X^{\mathrm{ref}}$.
\end{itemize}
\end{thm}

The rest of this subsection will be devoted to showing this theorem.

\subsubsection{Coherent case}\label{subsub-reflexivizationcoherent}
First, we construct $X^{\mathrm{ref}}$ for a coherent valuative space $X$.
Let $A=\QCOuv(X)$ be the distributive lattice\index{lattice!distributive lattice@distributive ---} of quasi-compact open subsets of $X$.
We have $X\cong\Spec A$ by Stone duality (\ref{thm-stonerepresentationthm}).
Consider the map 
$$
\varphi\colon A\longrightarrow\mathbf{2}^{[X]}
$$
to the power-set lattice of $[X]$ given by $U\mapsto\sep_X(U)=[U]$.
This is a lattice homomorphism, for we have $\sep_X(U)=U\cap [X]$ set-theoretically.
Hence the image $B$ of this mapping is a distributive lattice, giving a coherent sober space $\Spec B$, which we denote by $X^{\mathrm{ref}}$.
It is then easy to see that the surjective homomorphism $\varphi\colon A\rightarrow B$ gives rise to a quasi-compact map $i_X\colon X^{\mathrm{ref}}\rightarrow X$, which induces a homeomorphism from $X^{\mathrm{ref}}$ to its image with the subspace topology.
We thereby identify $X^{\mathrm{ref}}$ with its image in $X$.

\begin{prop}\label{prop-lemeflexivizationcoherent3}
$X^{\mathrm{ref}}$ is a valuative space, and the map $i_X\colon X^{\mathrm{ref}}\hookrightarrow X$ is valuative, inducing a homeomorphism $[X^{\mathrm{ref}}]\approx [X]$.
\end{prop}

To show this, we need the following lemma: 
\begin{lem}\label{lem-Fujiwara1.11}
Let $S$ be a compact space, and $D$ a distributive sublattice of $\mathbf{2}^S$ consisting of closed subspaces of $S$.
Suppose $D$ satisfies the following conditions$:$
\begin{itemize}
\item[{\rm (a)}] for any $x\in S$, there exists $C\in D$ such that $x\in C;$
\item[{\rm (b)}] for any $x,y\in S$ with $x\neq y$, there exists $C\in D$ such that $x\in C$ and $y\not\in C$.
\end{itemize}
For any $x\in S$, set $F_x=\{C\in D\,|\,x\in C\}$ $($non-empty due to {\rm (a)}$)$.

{\rm (1)} Any prime filter\index{filter (of lattices)!prime filter of lattices@prime ---} {\rm (cf.\ \ref{subsub-strcohtopsp})} is contained in a filter of the form $F_x$ for some $x\in S$, and, for any $x\in S$, $F_x$ is a maximal filter {\rm (\ref{dfn-idealfilter} (2))} of $D$.

{\rm (2)} The map $x\mapsto F_x$ gives a bijection between $S$ to the set of all maximal filters of $D$.
\end{lem}

\begin{proof}
(1) It is easy to see that $F_x$ for $x\in S$ is a prime filter.
Let $F\subseteq D$ be a prime filter.
For any finitely many elements $C_1,\ldots,C_r\in F$, we have $C_1\cap\cdots\cap C_r\in F$, and since $\emptyset\not\in F$, we have $C_1\cap\cdots\cap C_r\neq\emptyset$.
This implies, since $S$ is compact, that the intersection of all $C$'s in $F$ is non-empty, containing a point $x\in S$.
This shows $F\subseteq F_x$.

To show that $F_x$ is maximal, suppose $F_x$ is contained in a prime filter $F$.
By what we have seen above, there exists $y\in S$ such that $F\subseteq F_y$, hence $F_x\subseteq F_y$.
But the assumption (b) implies $\bigcap_{C\in F_x}C=\{x\}$, from which we deduce $x=y$, and thus $F=F_x$.
Hence $F_x$ is a maximal filter.

(2) By (1), the map $x\mapsto F_x$ maps $S$ surjectively onto the set of all maximal filters.
Since $\bigcap_{C\in F_x}C=\{x\}$, it is injective, too.
\end{proof}

\begin{proof}[Proof of Proposition {\rm \ref{prop-lemeflexivizationcoherent3}}]
Let us first confirm that the distributive lattice $B$ satisfies the hypotheses of \ref{lem-Fujiwara1.11}.
The condition (a) clearly holds.
For $x,y\in[X]$ ($x\neq y$), since $[X]$ is Hausdorff (\ref{cor-tubes2acor} (2)), $\sep^{-1}_X(\{y\})$ is closed, and hence there exists a coherent open neighborhood $U$ of $x$ such that $y\not\in [U]$, whence verifying (b).

To show that $X^{\mathrm{ref}}$ is valuative, since the topology of $X^{\mathrm{ref}}$ is the subspace topology from $X$, which is valuative, we only need to show that any $x\in X^{\mathrm{ref}}$ has a maximal generization (as in \ref{rem-valuativespace}, this is automatic, but allows a direct proof as follows).
By \ref{lem-Fujiwara1.11}, we know that $X^{\mathrm{ref}}=\Spec B$ contains all points in $[X]$, in which any point has a maximal generization. 
This means that $X^{\mathrm{ref}}$ is valuative, that the map $i_X$ is valuative, and that $[i_X]$ gives a continuous bijection from $[X^{\mathrm{ref}}]$ to $[X]$.
Since both $[X^{\mathrm{ref}}]$ and $[X]$ are compact (Hausdorff)  (\ref{cor-tubes2acor} (2)), $[i_X]$ is a homeomorphism.
\end{proof}

\begin{prop}\label{prop-lemeflexivizationcoherent4}
For any quasi-compact valuative map $f\colon Z\rightarrow X$ from a reflexive coherent valuative space $Z$, there exists uniquely a valuative map $h\colon Z\rightarrow X^{\mathrm{ref}}$ such that $f=i_X\circ h$.
\end{prop}

\begin{proof}
For any subset $S\subseteq [X^{\mathrm{ref}}]$, set $H(S)=(\sep^{-1}_Z([f]^{-1}(S)))^{\circ}$, where $(\cdot)^{\circ}$ denotes the interior kernel.
For $S=[U]$ with $U\in A$, since $f$ is valuative, we have $H([U])=(f^{-1}(\sep^{-1}_X([U])))^{\circ}=(f^{-1}(\ovl{U}))^{\circ}$ (due to \ref{cor-ZRpoints61a}), which is equal to $f^{-1}(U)$ by \ref{cor-valuativemapstubes2} (1) and \ref{prop-lemlocallystronglycompactconverse}.
Since $f$ is quasi-compact, the map $H$ gives a lattice homomorphism from $B=\QCOuv(X^{\mathrm{ref}})$ to $\QCOuv(Z)$, which defines a quasi-compact map $h\colon Z\rightarrow X^{\mathrm{ref}}$.
Since $f$ is valuative, so is $h$.
Since the composition $H\circ\varphi\colon A=\QCOuv(X)\rightarrow\QCOuv(Z)$ clearly coincides with the homomorphism corresponding to $f\colon Z\rightarrow X$, we have $f=i_X\circ h$.
Finally, the uniqueness of $h$ follows from \ref{prop-lemlocallyquasicompactmaps} and the fact that $i_X$ is a homeomorphism.
\end{proof}

By the last proposition, we see that $X\mapsto X^{\mathrm{ref}}$ for $X$ coherent is functorial, giving the left adjoint functor of the inclusion functor from the category of coherent reflexive valuative spaces with quasi-compact maps to the category of coherent valuative spaces with quasi-compact maps.

\subsubsection{General case}\label{subsub-reflexivizationgeneral}
Let $X$ be a valuative space, and consider the functor 
$$
F_X\colon \RVsp^{\opp}\longrightarrow\Sets
$$
defined as follows: for any reflexive valuative space $Z$, we set $F_X(Z)=\Hom_{\mathscr{C}}(Z,X)$, the set of all locally quasi-compact valuative maps from $Z$ to $X$.
Let us consider the category $\mathscr{D}$ of quasi-separated reflexive valuative space with quasi-compact valuative maps, and the category $\mathbf{CRVsp}$ of coherent reflexive valuative space with quasi-compact valuative maps.
The discussion in the previous paragraph shows that the functor $F_X|_{\mathbf{CRVsp}^{\opp}}$ is representable for any coherent valuative space $X$.

\begin{lem}\label{lem-reflexivizationgeneral1}
Suppose that the functor $F_X|_{\mathbf{CRVsp}^{\opp}}$ is representable by a coherent reflexive valuative space $X'$.
Then $F_X$ itself is representable by $X'$.
\end{lem}

\begin{proof}
Let us first show that $F_X|_{\mathscr{D}^{\opp}}$ is representable by $X'$.
For any quasi-separated reflexive valuative space $Z$, write $Z$ as a filtered inductive limit $\varinjlim_{i\in I}Z_i$ of coherent open subspaces.
Then giving $f\colon Z\rightarrow X$ is equivalent to giving the collection $\{f_i=f|_{Z_i}\colon Z_i\rightarrow X\}_{i\in I}$ of maps satisfying $f_j|_{Z_i}=f_i$ for $i\leq j$.
By the assumption, each $f_i$ factors through $X'$ by $f'_i\colon Z_i\rightarrow X'$, and by the functoriality, we have $f'_j|_{Z_i}=f'_i$ for $i\leq j$.
We thus obtain $f'\colon Z\rightarrow X'$, which shows the representability of $F_X|_{\mathscr{D}^{\opp}}$.

In general, given a reflexive valuative space $Z$, we take an open covering $\{Z_{\alpha}\}_{\alpha\in L}$ of $Z$ by coherent open subspaces.
Then giving $f\colon Z\rightarrow X$ is equivalent to giving the collection of maps $\{f_{\alpha}\colon Z_{\alpha}\rightarrow X\}_{\alpha\in L}$ in such a way that $f_{\alpha}|_{Z_{\alpha}\cap Z_{\beta}}=f_{\beta}|_{Z_{\alpha}\cap Z_{\beta}}$.
By the assumption, each $f_{\alpha}$ is canonically factored through $X'$ by $f'_{\alpha}\colon Z_{\alpha}\rightarrow X'$.
Since $Z_{\alpha}\cap Z_{\beta}$ are quasi-separated, what we have already shown above implies $f'_{\alpha}|_{Z_{\alpha}\cap Z_{\beta}}=f'_{\beta}|_{Z_{\alpha}\cap Z_{\beta}}$.
Hence we obtain $f'\colon Z\rightarrow X'$ by patching $f'_{\alpha}$'s, which verifies the representability of $F_X$ by $X'$.
\end{proof}

By the lemma, we have already shown that the functor $F_X$ is representable for any coherent valuative space $X$, and whence the existence of $(\,\cdot\,)^{\mathrm{ref}}$ as a functor $\mathbf{CVsp}\rightarrow\mathbf{CRVsp}$ of coherent valuative spaces with quasi-compact maps.

The following lemma shows that the reflexivizaion functor $(\,\cdot\,)^{\mathrm{ref}}$, if it exists, commutes with open immersions.
\begin{lem}\label{lem-reflexivizationgeneral2}
Suppose that $F_X$ for a valuative space $X$ is representable by $(X^{\mathrm{ref}},i_X)$.
Then for any open subset $U\subseteq X$, $(i^{-1}_X(U),i_X|_{i^{-1}_X(U)})$ represents $F_U$.
\end{lem}

\begin{proof}
It is enough to invoke that any open subspace of a reflexive valuative space is again reflexive, which is trivial.
\end{proof}

Now we are going to construct $X^{\mathrm{ref}}$ for any valuative space.
Take an open covering $\{X_{\alpha}\}_{\alpha\in L}$ of $X$ by coherent open subspaces.
For each $\alpha\in L$, we have $(X^{\mathrm{ref}}_{\alpha},i_{X_{\alpha}})$.
By \ref{lem-reflexivizationgeneral2}, $i^{-1}_{X_{\alpha}}(X_{\alpha}\cap X_{\beta})$ gives the reflexivization of $X_{\alpha}\cap X_{\beta}$, hence being equal to $i^{-1}_{X_{\beta}}(X_{\alpha}\cap X_{\beta})$, which one can consistently denote by $(X_{\alpha}\cap X_{\beta})^{\mathrm{ref}}$.
One can then patch $X^{\mathrm{ref}}_{\alpha}$'s along $(X_{\alpha}\cap X_{\beta})^{\mathrm{ref}}$'s to obtain a valuative space $X'$, which is easily seen to be reflexive, and a map $i_X\colon X'\rightarrow X$.

In order to see that $X'$ gives the reflexivization of $X$, it remains to show that $X'$ enjoys the desired functoriality (which, at the same time, also confirms that the formation of $X'$ is independent, up to canonical isomorphism, of all choices we have made).

For a given $f\colon Z\rightarrow X$, where $Z$ is reflexive, define $Z_{\alpha}=f^{-1}(X_{\alpha})$ and $f_{\alpha}=f|_{Z_{\alpha}}$ for any $\alpha\in L$.
For each $\alpha\in L$, we have the canonical factorization 
$$
Z_{\alpha}\stackrel{f^{\mathrm{ref}}_{\alpha}}{\longrightarrow}X^{\mathrm{ref}}_{\alpha}\stackrel{i_{X_{\alpha}}}{\longrightarrow}X_{\alpha}
$$
of the map $f_{\alpha}$.
By the functoriality, $f^{\mathrm{ref}}_{\alpha}|_{Z_{\alpha}\cap Z_{\beta}}$ is equal to $(f_{\alpha}|_{Z_{\alpha}\cap Z_{\beta}})^{\mathrm{ref}}$, and hence is equal to $f^{\mathrm{ref}}_{\beta}|_{Z_{\alpha}\cap Z_{\beta}}$.
Hence, by patching, we obtain a map $f'\colon Z\rightarrow X'$ such that $i_X\circ f'=f$.
The map $f$ is valuative and locally quasi-compact, since it is obtained by patching local maps having these properties.
Since the representability holds locally, and the above construction is functorial, it follows that $(X',i_X)$ represents the functor $F_X$, and that $X^{\mathrm{ref}}=X'$ gives the reflexivization of $X$.

To conclude the proof of Theorem \ref{thm-reflexivization}, it remains to verify (a) $\sim$ (d).
The condition (a) follows immediately from \ref{lem-reflexivizationgeneral2}.
For (b), notice that the functor $[\,\cdot\,]$ has the continuity (\ref{prop-structuresepquot2}), and, by our construction of $X^{\mathrm{ref}}$, the functor $(\,\cdot\,)^{\mathrm{ref}}$ has the similar continuity.
Hence checking (b) is reduced to that in the case where $X$ is coherent, which has already been shown in \ref{prop-lemeflexivizationcoherent3}.
The condition (c) follows from (a) and the fact that $X^{\mathrm{ref}}$ is coherent if $X$ is coherent.
Finally, (d) follows from (a) and (c).
\hfill$\square$

\begin{cor}\label{cor-proplemlocallyquasicompactmapsgen}
Let $X,Y$ be valuative spaces, and $f,g\colon X\rightarrow Y$ two valuative locally quasi-compact maps.
Suppose $X$ is reflexive.
Then $[f]=[g]$ implies $f=g$. \hfill$\square$
\end{cor}
\index{reflexivization|)}\index{valuative!valuative topological space@--- (topological) space!reflexive valuative topological space@reflexive --- ---|)}

\subsection{Locally strongly compact valuative spaces}\label{sub-locallystronglycompact}
\index{valuative!valuative topological space@--- (topological) space!locally strongly compact valuative topological space@locally strongly compact --- ---|(}
\subsubsection{Locally strongly compact valuative spaces}\label{subsub-locallystronglycompact}
\begin{dfn}\label{dfn-locallycompactspace}{\rm 
A valuative space $X$ is said to be {\em locally strongly compact} if for any $x\in X$ there exists a pair $(W_x,V_x)$ consisting of an overconvergent open subset\index{overconvergent!overconvergent subset@--- subset} $W_x$ and a coherent open subset $V_x$ such that $x\in W_x\subseteq V_x$.}
\end{dfn}

If $X$ is coherent, then one can take $W_x=V_x=X$ for any $x\in X$.
Hence: 
\begin{prop}\label{prop-locallycompactspacecoherent}
Any coherent valuative space is locally strongly compact. \hfill$\square$
\end{prop}

\begin{prop}\label{prop-quasicompactnessclosuresingletonset}
Let $X$ be a locally strongly compact valuative space.
Then for any $x\in X$ the closure $\ovl{\{x\}}$ of $\{x\}$ in $X$ is quasi-compact.
\end{prop}

\begin{proof}
Take $(W_x,V_x)$ as in \ref{dfn-locallycompactspace}.
Since $W_x$ is overconvergent, $\ovl{\{x\}}$ is contained in $W_x$, and hence in $V_x$.
Thus $\ovl{\{x\}}$ is a closed subset of the quasi-compact space $V_x$, and hence is quasi-compact.
\end{proof}

\begin{prop}\label{prop-locallycompactnessformerdef}
Let $X$ be a quasi-separated valuative space.
Then $X$ is locally strongly compact if and only if the following condition is satisfied$:$ For any $x\in X$ there exists a pair $(U_x,V_x)$ of coherent open neighborhoods of $\ovl{\{x\}}$ such that $\ovl{U_x}\subseteq V_x$.
\end{prop}

\begin{proof}
Suppose $X$ is locally strongly compact, and let $x\in X$.
Take $(W_x,V_x)$ as in \ref{dfn-locallycompactspace}.
Since $\ovl{\{x\}}$ is quasi-compact, one can take a quasi-compact open neighborhood $U_x$ of $\ovl{\{x\}}$ inside $W_x$; since $X$ is quasi-separated, $U_x$ is coherent.
By \ref{cor-ZRpoints61a} we have $\ovl{U_x}\subseteq W_x\subseteq V_x$.

Conversely, suppose $X$ satisfies the condition as in the assertion.
For $x\in X$, take $(U_x,V_x)$ as above.
We may assume without loss of generality that $x$ is a maximal point. 
Then, by \ref{cor-interior4a}, the pair $(W_x,V_x)$ with $W_x=\int_X(U_x)$ gives a pair as in \ref{dfn-locallycompactspace}.
\end{proof}

\begin{prop}\label{prop-locstringcompsp2a}
Let $X$ be a quasi-separated valuative space. 
Then the following conditions are equivalent$:$
\begin{itemize}
\item[{\rm (a)}] $X$ is locally strongly compact$;$
\item[{\rm (b)}] the closure $\ovl{U}$ of any quasi-compact open subset $U\subseteq X$ is quasi-compact$;$ 
\item[{\rm (c)}] there exists an open covering $\{U_{\alpha}\}_{\alpha\in L}$ of $X$ such that $U_{\alpha}$ and $\ovl{U}_{\alpha}$ are quasi-compact for any $\alpha\in L$.
\end{itemize}
\end{prop}

\begin{proof}
First we prove (a) $\Rightarrow$ (c). 
By \ref{prop-locallycompactnessformerdef}, $X$ has an open covering $\{U_{\alpha}\}_{\alpha\in L}$ such that each $U_{\alpha}$ is coherent and its closure $\ovl{U}_{\alpha}$ is contained in a coherent open subset $V_{\alpha}$.
Since $V_{\alpha}$ is quasi-compact, each $\ovl{U}_{\alpha}$ is quasi-compact.

Next we show (c) $\Rightarrow$ (b). 
Since $U$ is quasi-compact, there exist finitely many $\alpha_1,\ldots,\alpha_n\in L$ such that $U=\bigcup_{j=1}^nU\cap U_{\alpha_j}$.
By \ref{prop-tubes121a} we have $\ovl{U}=\bigcup_{j=1}^n\ovl{U}\cap\ovl{U}_{\alpha_j}$.
Since each $\ovl{U}\cap\ovl{U}_{\alpha_j}$ is quasi-compact, one deduces that $\ovl{U}$ is quasi-compact.

Finally, let us show (b) $\Rightarrow$ (a).
For any $x\in X$, take a quasi-compact open neighborhood $W$ of $x$.
We have $\ovl{\{x\}}\subset\ovl{W}$.
Since $\ovl{W}$ is quasi-compact, so is $\ovl{\{x\}}$.
Then one has a quasi-compact (hence coherent) open subset $U$ that contains $\ovl{\{x\}}$.
Since $\ovl{U}$ is again quasi-compact, it is further contained in a quasi-compact open subset $V$.
In view of \ref{prop-locallycompactnessformerdef}, this shows (a).
\end{proof}

\begin{rem}\label{rem-taut}{\rm 
Proposition \ref{prop-locstringcompsp2a} shows that, if $X$ is a quasi-separated valuative space, then $X$ is locally strongly compact if and only if it is {\em taut} in the sense of \cite[5.1.2]{Hube3}.}
\end{rem}

\subsubsection{Characteristic properties}\label{subsub-locallystronglycompactcharacterization}
\begin{thm}\label{thm-locallystronglycompacttheorem}
Let $X$ be a locally strongly compact valuative space.
Then $X$ satisfies the following properties$:$
\begin{itemize}
\item[{\rm (a)}] the separated quotient $[X]$ is locally compact $($hence locally Hausdorff$)$ space$;$
\item[{\rm (b)}] the separation map $\sep_X\colon X\rightarrow[X]$ is proper$;$
\item[{\rm (c)}] each $x\in[X]$ admits an open neighborhood $\mathscr{Z}_x\subseteq[X]$ such that $\sep^{-1}_X(\mathscr{Z}_x)$ is quasi-separated.
\end{itemize}
Conversely, if a valuative space $X$ satisfies {\rm (b)} and {\rm (c)}, then $X$ is locally strongly compact.
\end{thm}

\begin{proof}
For any $x\in[X]$, take $(W_x,V_x)$ as in \ref{dfn-locallycompactspace}.
Since $W_x$ is overconvergent also in $V_x$, we have $\sep_X(W_x)=\sep_{V_x}(W_x)$, which is open both in $[V_x]$ and $[X]$.
Now, by \ref{cor-tubes2acor} (2), we know that $[V_x]$ is compact. 
Hence $\mathscr{W}_x=\sep_X(W_x)=\sep_{V_x}(W_x)$ is a Hausdorff open neighborhood of $x$ both in $[V_x]$ and $[X]$, and $[V_x]$ gives a compact neighborhood of $x$ in $[X]$. 
This shows that $\mathscr{W}_x$ is a locally compact Hausdorff space, whence the property (a).
To show (b), it suffices to show that $W_x=\sep^{-1}_X(\mathscr{W}_x)\rightarrow\mathscr{W}_x$ is proper for any $x\in[X]$.
Since $W_x=\sep^{-1}_{V_x}(\mathscr{W}_x)$, this follows from the properness of $V_x\rightarrow[V_x]$ (\ref{cor-tubes211acor}).
One verifies the property (c) with $\mathscr{Z}_x=\mathscr{W}_x$; indeed, $W_x=\sep^{-1}_{V_x}(\mathscr{W}_x)$ is an open subspace of the coherent space $V_x$, and hence is quasi-separated.

Now we show the converse. 
Suppose $X$ is a valuative space that satisfies (b) and (c).
For any point $x\in X$, we need to find a pair $(W_x,V_x)$ as in \ref{dfn-locallycompactspace}.
To this end, we may assume that $x$ is a maximal point, that is, $x\in[X]$.
Take an open neighborhood $\mathscr{Z}_x$ of $x$ in $[X]$ such that $Z_x=\sep^{-1}_X(\mathscr{Z}_x)$ is quasi-separated.
It follows from (b) that the closure $\ovl{\{x\}}$ of $\{x\}$ in $X$, being equal to the pull-back $\sep^{-1}_X(\{x\})$, is quasi-compact.
Since $\ovl{\{x\}}$ is contained in $Z_x$, we have a quasi-compact open subset $U_x\subseteq Z_x$ that contains $\ovl{\{x\}}$.
Since $Z_x$ is quasi-separated, in view of \ref{cor-ZRpoints61a}, the closure $\ovl{U}$ of $U$ in $Z_x$ is equal to $\sep^{-1}_X(\sep_X(U_x))$.
This shows that $[U_x]$ is closed in $[X]$, and hence $\ovl{U}_x$ is quasi-compact by the properness of $\sep_X$.
Thus there exists a quasi-compact open neighborhood $V_x$ of $\ovl{U}_x$ in $Z_x$.
By \ref{cor-interior4a}, we have $x\in\int_{Z_x}(U_x)$.
Set $W_x=\int_{Z_x}(U_x)$.
Since $W_x$ is overconvergent in $Z_x$ and $Z_x$ is overconvergent in $X$, $W_x$ is overconvergent in $X$.
Hence the pair $(W_x,V_x)$ is a desired one.
\end{proof}

\begin{cor}\label{cor-locallystronglycompacttheorem}
Let $X$ be a valuative space. Then the following conditions are equivalent to each other$:$
\begin{itemize}
\item[{\rm (a)}] $X$ is locally strongly compact$;$
\item[{\rm (b)}] any overconvergent open subset\index{overconvergent!overconvergent subset@--- subset} $W\subseteq X$ is locally strongly compact$;$
\item[{\rm (c)}] there exists an open covering $\{W_{\alpha}\}_{\alpha\in L}$ of $X$ consisting of overconvergent open subsets such that each $W_{\alpha}$ $(\alpha\in L)$ is locally strongly compact. \hfill$\square$
\end{itemize}
\end{cor}

\begin{cor}\label{cor-thmlocstringcompsp1a}
The separated quotient $[X]$ of a locally strongly compact valuative space $X$ is locally compact $($hence locally Hausdorff$)$.
It is, moreover, Hausdorff, if $X$ is quasi-separated.
\end{cor}

\begin{proof}
The first assertion is already proven in \ref{thm-locallystronglycompacttheorem}.
Suppose $X$ is quasi-separated.
Take two distinct points $x_1\neq x_2$ in $[X]$ and, for each $i=1,2$, a pair $(W_{x_i},V_{x_i})$ as in \ref{dfn-locallycompactspace}. 
Since $U=V_{x_1}\cup V_{x_2}$ is quasi-compact open and $X$ is quasi-separated, $U$ is coherent open.
In view of \ref{prop-tubes2a} (2), take for each $i=1,2$ an open neighbourhoods $W'_i$ of $y_i$ that is contained and overconvergent in $U$, such that $W'_1\cap W'_2=\emptyset$.
Then $W_i\cap W'_i$ for $i=1,2$ are overconvergent in $X$ and separate $\ovl{\{y_1\}}$ and $\ovl{\{y_2\}}$.
\end{proof}

\begin{cor}\label{cor-locstringcompsp1cora}
Let $X$ be a quasi-separated locally strongly compact valuative space. 
Then $[X] $ is a completely regular\index{space@space (topological)!completely regular topological space@completely regular ---} topological space {\rm (cf.\ \S\ref{subsub-crspaces})}.\hfill$\square$
\end{cor}

\begin{prop}\label{prop-locallystronglycompactconverse}
Let $X$ be a reflexive valuative space\index{valuative!valuative topological space@--- (topological) space!reflexive valuative topological space@reflexive --- ---} {\rm (\ref{dfn-reflexivevaluativespaces})}. 
Suppose there exists a family $\{U_{\alpha}\}_{\alpha\in L}$ of open subsets $X$ such that$:$
\begin{itemize}
\item[{\rm (a)}] each $U_{\alpha}$ is locally strongly compact, and retrocompact\index{retrocompact} {\rm (\ref{dfn-retrocompact})} in $X;$
\item[{\rm (b)}] $[X]=\bigcup_{\alpha\in L}[U_{\alpha}]^{\circ}$, where $(\cdot)^{\circ}$ denotes the interior kernel.
\end{itemize}
Then $X$ is locally strongly compact.
\end{prop}

\begin{proof}
Set $\mathscr{U}_{\alpha}=[U_{\alpha}]^{\circ}$ and $W_{\alpha}=\sep^{-1}_X(\mathscr{U}_{\alpha})$ for $\alpha\in L$.
We have $W_{\alpha}\subseteq(\sep^{-1}_X([U_{\alpha}]))^{\circ}=(\ovl{U_{\alpha}})^{\circ}$ (\ref{cor-ZRpoints61a}).
By \ref{prop-lemlocallystronglycompactconverse} we have $W_{\alpha}\subseteq U_{\alpha}$ for any $\alpha\in L$.
Thus we have the following diagram
$$
\xymatrix{W_{\alpha}\ar@{^{(}->}[r]\ar[d]_{\sep_X|_{W_{\alpha}}}&U_{\alpha}\ar@{^{(}->}[r]\ar[d]_{\sep_{U_{\alpha}}}&X\ar[d]^{\sep_X}\\ \mathscr{U}_{\alpha}\ar@{^{(}->}[r]&[U_{\alpha}]\ar@{^{(}->}[r]&[X].}
$$
By Exercise \ref{exer-locallystronglycompact2}, the topology of $[U_{\alpha}]$ coincides with the subspace topology from $[X]$.
Since $W_{\alpha}$ is overconvergent in $U_{\alpha}$, it is locally strongly compact due to \ref{cor-locallystronglycompacttheorem}.
Since $\{W_{\alpha}\}_{\alpha\in L}$ covers $X$, again by \ref{cor-locallystronglycompacttheorem}, we deduce that $X$ is locally strongly compact.
\end{proof}

\begin{cor}\label{cor-locallystronglycompactconversecor1}
Let $X$ be a reflexive and quasi-separated valuative space having a family $\{U_{\alpha}\}_{\alpha\in L}$ of coherent open subsets such that $[X]=\bigcup_{\alpha\in L}[U_{\alpha}]^{\circ}$.
Then $X$ is locally strongly compact. \hfill$\square$
\end{cor}

\subsubsection{Paracompact spaces}\label{subsub-paracompactspaces}
\index{space@space (topological)!paracompact topological space@paracompact ---|(}
\begin{dfn}\label{dfn-paracompactgentop}{\rm 
Let $X$ be a topological space.

(1) An open covering $\{U_{\alpha}\}_{\alpha\in L}$ of $X$ is said to be {\em locally finite}\index{locally finite covering@locally finite (covering)} if any $x\in X$ has an open neighborhood $V$ such that $V\cap U_{\alpha}\neq\emptyset$ for at most finitely many indices $\alpha\in L$.

(2) The space $X$ is said to be {\em paracompact} if any open covering admits a refinement by a locally finite covering.}
\end{dfn}

\danger{Notice that our definition of `paracompact' differs from that in \cite[Chap.\ I, \S9.10, Def.\ 6]{Bourb4} in that we do not assume Hausdorffness.}

The following lemma is easy to see, and the proof is left to the reader:
\begin{lem}\label{lem-paracompact1}
Let $X$ be a locally coherent space.

{\rm (1)} The space $X$ is paracompact if there is a locally finite covering by quasi-compact open subsets.

{\rm (2)} Let $\{U_{\alpha}\}_{\alpha\in L}$ be a locally finite covering of $X$ consisting of quasi-compact open subsets.
Then for any $\alpha\in L$ the set of all indices $\beta\in L$ such that $U_{\beta}$ intersects $U_{\alpha}$ is finite. \hfill$\square$
\end{lem}

\begin{prop}\label{prop-paracompact1}
Let $X$ be a paracompact quasi-separated valuative space. 
Then $X$ is locally strongly compact.
\end{prop}

\begin{proof}
By \ref{prop-locstringcompsp2a} it suffices to show that for any quasi-compact open subset $U$ its closure $\ovl{U}$ is again quasi-compact.
Let $\{V_{\alpha}\}_{\alpha\in L}$ be a covering of $\ovl{U}$, that is, $\ovl{U}\subseteq\bigcup_{\alpha\in L}V_{\alpha}$.
Together with $X\setminus\ovl{U}$ this gives an open covering of $X$, and hence, there exists a locally finite refinement $\{W_{\lambda}\}_{\lambda\in\Lambda}$ of $\{V_{\alpha}\}_{\alpha\in L}$ by quasi-compact open subsets.
Let $\Lambda'$ be the subset of $\Lambda$ consisting of all $\lambda$ such that $\ovl{U}\cap W_{\lambda}\neq\emptyset$.
Since $U\cap V_{\lambda}\neq\emptyset$ for any $\lambda\in\Lambda'$ and $U$ is quasi-compact, the set $\Lambda'$ is actually finite, and thus $\bigcup_{\lambda\in\Lambda'}V_{\lambda}$ is quasi-compact.
Since $\ovl{U}\subseteq\bigcup_{\lambda\in\Lambda'}V_{\lambda}$, $\ovl{U}$ is quasi-compact, as desired.
\end{proof}

The following proposition can be shown by an argument similar to that in \cite[Chap.\ 1, \S9.10, Theorem 5]{Bourb4}:
\begin{prop}\label{prop-paracompactsigmacompact}
Let $X$ be a paracompact quasi-separated locally strongly compact valuative space.
Then $X$ is decomposed into a disjoint sum $X=\coprod_{\lambda\in\Lambda}X_{\lambda}$ such that for each $\lambda\in\Lambda$ there exists an increasing sequence $U_{\lambda,0}\subseteq U_{\lambda,1}\subseteq\cdots$ of quasi-compact open subsets satisfying$:$
\begin{itemize}
\item[{\rm (a)}] $\ovl{U}_{\lambda,n}\subseteq U_{\lambda,n+1}$ for any $n\geq 0;$
\item[{\rm (b)}] $X_{\lambda}=\bigcup_{n\geq 0}U_{\lambda,n}$.\hfill$\square$
\end{itemize}
\end{prop}

\begin{prop}\label{prop-paracompactseparatedquotient}
Let $X$ be a quasi-separated locally strongly compact valuative space.
Then $X$ is paracompact if and only if $[X]$ is paracompact $($in the sense as in {\rm \cite[Chap.\ I, \S9.10, Def.\ 6]{Bourb4}}$;$ notice that, due to {\rm \ref{cor-thmlocstringcompsp1a}}, $[X]$ is locally compact and Hausdorff$)$.
\end{prop}

\begin{proof}
Suppose $X$ is paracompact.
We may assume in view of \ref{prop-paracompactsigmacompact} that there exists an increasing sequence $U_0\subseteq U_1\subseteq\cdots$ of quasi-compact open subsets satisfying the conditions similar to (a) and (b) in \ref{prop-paracompactsigmacompact}.
For each $n\geq 0$ let $V_n$ be the interior kernel of $[U_n]$ in $[X]$.
Since $\ovl{U}_n$ is overconvergent (\ref{cor-ZRpoints61a}), we have $\sep^{-1}_X(\ovl{V}_n)\subseteq\ovl{U}_n$ and hence $\sep^{-1}_X(\ovl{V}_n)\subseteq\int_X(U_{n+1})$ (cf.\ \ref{cor-interior4a}).
Thus we have $\ovl{V}_n\subseteq V_{n+1}$ for any $n\geq 0$, and hence $[X]$ is paracompact due to \cite[Chap.\ 1, \S9.10, Theorem 5]{Bourb4}.

Conversely, suppose $[X]$ is paracompact, and take three open coverings $X=\bigcup_{\alpha\in L}U_{\alpha}=\bigcup_{\alpha\in L}V_{\alpha}=\bigcup_{\alpha\in L}W_{\alpha}$ of $X$ by quasi-compact open subsets such that for any $\alpha\in L$ we have $\ovl{U}_{\alpha}\subseteq V_{\alpha}$ and $\ovl{V}_{\alpha}\subseteq W_{\alpha}$ (possible due to \ref{prop-locstringcompsp2a} (b)).
By \ref{cor-interior4a} we have $[X]=\bigcup_{\alpha\in L}[V_{\alpha}]^{\circ}=\bigcup_{\alpha\in L}[W_{\alpha}]^{\circ}$ (where $(\cdot)^{\circ}$ denotes the interior kernel).
Take a locally finite open coverings $[X]=\bigcup_{i\in I}\til{V}_i$ and $[X]=\bigcup_{j\in J}\til{W}_j$ that refine $\{[V_{\alpha}]^{\circ}\}_{\alpha\in L}$ and $\{[W_{\alpha}]^{\circ}\}_{\alpha\in L}$, respectively.
Since for each $i\in I$ there exists an $\alpha\in L$ such that $\sep^{-1}_X(\til{V}_i)\subseteq\ovl{V}_{\alpha}$, we deduce that the closure of $\sep^{-1}_X(\til{V}_i)$ is quasi-compact.
Moreover, since $[X]=\bigcup_{j\in J}\til{W}_j$ is a locally finite covering, $\sep^{-1}_X(\til{V}_i)\cap\sep^{-1}_X(\til{W}_j)$ is non-empty for only finitely many $j\in J$, say, $j_{i,1},\ldots,j_{i,m_i}$.
One can take a finite collection $\{U_{i_1},\ldots,U_{i_{n_i}}\}$ of quasi-compact open subsets contained in $\bigcup^{m_i}_{k=1}\sep^{-1}_X(\til{W}_{j_{i,k}})$ such that $\sep^{-1}_X(\til{V}_i)\subseteq\bigcup^{n_i}_{k=1}U_{i_k}$.
Then $\{U_{i_k}\}_{i\in I,1\leq k\leq n_i}$ is a locally finite covering of $X$ by quasi-compact open subsets.
\end{proof}
\index{valuative!valuative topological space@--- (topological) space!locally strongly compact valuative topological space@locally strongly compact --- ---|)}\index{space@space (topological)!paracompact topological space@paracompact ---|)}\index{space@space (topological)!valuative topological space@valuative ---|)}\index{valuative!valuative topological space@--- (topological) space|)}

\subsection{Valuations of locally Hausdorff spaces}\label{sub-valuationslocallyhausdorff}
\index{valuation!valuation of a locally Hausdorff space@--- of a locally Hausdorff space|(}
In this subsection, we give a description of the category of reflexive\index{valuative!valuative topological space@--- (topological) space!reflexive valuative topological space@reflexive --- ---} (\ref{dfn-reflexivevaluativespaces}) locally strongly compact\index{valuative!valuative topological space@--- (topological) space!locally strongly compact valuative topological space@locally strongly compact --- ---} (\ref{dfn-locallycompactspace}) valuative spaces with valuative\index{valuative!valuative map@--- map} (\ref{dfn-valuativemaps}) and locally quasi-compact\index{map@map (continuous)!quasi-compact map@quasi-compact ---!locally quasi-compact map@locally --- ---}\index{quasi-compact!quasi-compact map@--- map!locally quasi-compact map@locally --- ---} (\ref{dfn-locallyquasicompactmaps}) maps.
We will show, in \ref{thm-RLSCVsp1} below, that this category is completely described, via separated quotients\index{separated quotient} (\S\ref{subsub-separationgen}), by locally compact (and locally Hausdorff) spaces with an extra structure, which we call a {\em valuation}.
The description we give here can be regarded as a variant of Stone duality (cf.\ \S\ref{subsub-strcohtopsp}).

\subsubsection{Nets and coverings}\label{sub-nets}
\begin{dfn}[cf.\ {\rm \cite{Berk2}\cite{Berk4}}]\label{dfn-netsfortopologicalspaces}
Let $X$ be a topological space.

{\rm (1)} A collection $\tau$ of subsets of $X$ is called a {\em quasi-net}\index{net!quasinet@quasi-{---}} on $X$ if for each $x\in X$ there exists finitely many $U_1,\ldots,U_n\in\tau$ such that $x\in U_1\cap\cdots\cap U_n$ and that $U_1\cup\cdots\cup U_n$ is a (possibly not open) neighborhood of $x$;

{\rm (2)} A quasi-net $\tau$ is called a {\em net}\index{net} on $X$ if for any $U,U'\in\tau$ the subset $\tau|_{U\cap U'}=\{U''\in\tau\,|\,U''\subseteq U\cap U'\}$ of $\tau$ is a quasi-net on the subspace $U\cap U'$.
\end{dfn}

Notice that, if $X$ is Hausdorff, then one can omit the condition `$x\in U_1\cap\cdots\cap U_n$' in (1).
Let us recall some of the first properties of quasi-nets and nets (cf.\ \cite[\S1.1]{Berk2}).
Let $\tau$ be a quasi-net on a topological space $X$:
\begin{itemize}
\item a subset $W\subseteq X$ is open if and only if for any $U\in\tau$ the intersection $U\cap W$ is open in $U$;
\item if any element of $\tau$ is compact, then $X$ is Hausdorff if and only if for any $U,V\in\tau$ the intersection $U\cap V$ is again compact.
\end{itemize}
Suppose, moreover, that $X$ is locally Hausdorff and $\tau$ is a net consisting of compact subsets, then:
\begin{itemize}
\item for any $U,V\in\tau$, the intersection $U\cap V$ is locally closed in $U$ and $V$.
\end{itemize}

\begin{prop}\label{prop-coveringandnets1topspace}
Let $X$ be a locally strongly compact\index{space@space (topological)!locally strongly compact topological space@locally strongly compact ---} {\rm (\ref{dfn-locallycompactspace})} valuative space\index{valuative!valuative topological space@--- (topological) space}\index{space@space (topological)!valuative topological space@valuative ---} {\rm (\ref{dfn-valuativespace})}, and $X=\bigcup_{\alpha\in L}U_{\alpha}$ an open covering by quasi-compact open subsets.
Then the collection $\tau=\{[U_{\alpha}]\}_{\alpha\in L}$ gives a quasi-net\index{net!quasinet@quasi-{---}} on the separated quotient\index{separated quotient} $[X]$ {\rm (\S\ref{subsub-separationgen})}.
It is a net\index{net} if and only if for any $\alpha,\beta\in L$ the collection $\{U_{\gamma}\,|U_{\gamma}\subseteq U_{\alpha}\cap U_{\beta}\}$ is a covering of $U_{\alpha}\cap U_{\beta}$.
\end{prop}

Notice that, due to \ref{cor-thmlocstringcompsp1a}, the separated quotient $[X]$ is a locally compact, hence locally Hausdorff, space.

\begin{proof}
Let $x\in[X]$ be a maximal point of $X$, and consider the closure $\ovl{\{x\}}$ in $X$ of the singleton set $\{x\}$.
Since $\ovl{\{x\}}$ is quasi-compact due to \ref{prop-quasicompactnessclosuresingletonset}, there exists an overconvergent open neighborhood $W$ of $\ovl{\{x\}}$ such that $\ovl{W}$ is quasi-compact.
Replacing $W$ by a smaller overconvergent open neighborhood of $\ovl{\{x\}}$, if necessary, one can take $\alpha_1,\ldots,\alpha_n\in L$ such that $\ovl{W}\subseteq U_{\alpha_1}\cup\cdots\cup U_{\alpha_n}$ and $\ovl{\{x\}}\cap U_{\alpha_i}\neq\emptyset$ for each $i=1,\ldots,n$.
Then $[U_{\alpha_1}]\cup\cdots\cup[U_{\alpha_n}]$ is a neighborhood of $x$ in $[X]$, since it contains $[W]$, which is open in $[X]$.
Moreover, since $\ovl{\{x\}}\cap U_{\alpha_i}\neq\emptyset$ for each $i=1,\ldots,n$, we have $x\in[U_{\alpha_1}]\cap\cdots\cap[U_{\alpha_n}]$, showing that $\tau=\{[U_{\alpha}]\}_{\alpha\in L}$ is a quasi-net on $[X]$.
The other assertion is clear.
\end{proof}

\subsubsection{Valuations of compact spaces}\label{subsub-valuationscompactHausdorff}
If $X$ is a coherent reflexive valuative space, then $X\cong X^{\mathrm{ref}}$ is completely determined by the distributive lattice $\QCOuv(X)$ of coherent open subsets, which is isomorphic to the distributive lattice $v=\{[U]\,|\,U\in\QCOuv(X)\}$ consisting of compact subsets of $[X]$ (cf.\ \S\ref{subsub-reflexivizationcoherent}).
This motivates the following definition:
\begin{dfn}\label{dfn-valuationscompactHausdorff}
Let $S$ be a compact (in particular, Hausdorff) space.

(1) A distributive sublattice $v$ of $\mathbf{2}^S$ is called a {\em valuation}\index{valuation!valuation of a compact space@--- of a compact space} of $S$ if there exists a continuous map $\pi\colon X\rightarrow S$ from a coherent reflexive valuative space $X$ such that the following conditions are satisfied:
\begin{itemize}
\item[{\rm (a)}] the map $[X]\rightarrow S$ induced from $\pi$ (cf.\ \ref{prop-separation3}) is a homeomorphism;
\item[{\rm (b)}] $v$ coincides, through the identification $[X]\cong S$, with the lattice $\{[U]\,|\,U\subseteq X\in\QCOuv(X)\}$.
\end{itemize}

(2) If $v$ is a valuation of $S$, then the pair $(S,v)$ is called a {\em valued compact space}\index{space@space (topological)!valued compact space@valued compact ---}.
\end{dfn}

Notice that, if $(S,v)$ is a valued compact space, then the coherent reflexive valuative space $X$ is uniquely determined up to canonical homeomorphisms and isomorphic to $\Spec v$ by Stone duality (cf.\ \ref{thm-stonerepresentationthm}).
Notice also that, if $v$ is a valuation of $S$, then every member $T\in v$ of $v$ is a compact (and hence closed) subset of $S$.
\begin{exa}\label{exa-valuationscompactHausdorff}
The distributive lattice $\mathbf{2}^S$ for a singleton set $S=\{\ast\}$ with the obvious compact topology is a valuation of $S$, for one has the obvious continuous map $\Spf V\rightarrow S$ from the formal spectrum of an $a$-adically complete valuation ring\index{valuation!valuation ring@--- ring!a-adically complete valuation ring@$a$-adically complete --- ---} $(a\in\m_V\setminus\{0\})$ of height one (cf.\ \ref{exa-reflexivevaluativespaces1}).
\end{exa}

\begin{prop}\label{prop-Fujiwara1.10}
Let $S$ be a compact $($hence Hausdorff$)$ space, and $v$ a distributive lattice consisting of compact subspaces of $S$.
Then $v$ is a valuation of $S$ if and only if the following conditions are satisfied$:$
\begin{itemize}
\item[{\rm (a)}] $\Spec v$ is a valuative space$;$
\item[{\rm (b)}] for any $x\in S$ and any neighborhood $U$ of $x$, there exists $T\in v$ such that $x\in T\subseteq U$.
\end{itemize}
Moreover, in this situation, there exists a canonical continuous map $\pi\colon\Spec v\rightarrow S$ that induces $[\Spec v]\cong S$.
\end{prop}

\begin{proof}
The `only if' part is clear.
For any $p\in\Spec v$, consider the prime filter $F_p$ corresponding to $p$, that is, the complement of $p$ in $v$.
As in the proof of \ref{lem-Fujiwara1.11}, points in $\bigcap_{T\in F_p}T$ canonically correspond to maximal filters $F'$ that contain $F_p$.
The condition (a) implies that such $F'$ is unique, and hence we have $\bigcap_{T\in F_p}T=\{\pi(p)\}$ for some $\pi(p)\in S$.
This gives a map
$$
\pi\colon\Spec v\longrightarrow S,\qquad p\longmapsto\pi(p),
$$
which is surjective due to \ref{lem-Fujiwara1.11}.
For each $T\in v$, set $U_T=\{p\in\Spec v\,|\,T\in F_p\}$.
Then $U_T$ is open in $\Spec v$, whose image under $\pi$ is contained in $T$.
Since $\{U_T\,|\,T\in v\}$ forms an open basis of $\Spec v$, the map $\pi$ is continuous by the condition (b).
Moreover, $\pi$ restricted on $[\Spec v]$ gives the inverse to the map $x\mapsto F_x$ as in \ref{lem-Fujiwara1.11} (2), and thus induces the continuous bijection $[\Spec v]\rightarrow S$ between compact Hausdorff spaces.
Hence $v$ is gives a valuation of $S$, and the proposition is proved.
\end{proof}

\begin{cor}\label{cor-Fujiwara1.12}
Let $S$ be a compact space, and $v$ and $v'$ distributive lattices consisting of compact subsets of $S$.
Suppose that $v\subseteq v'$, that $v$ is a valuation of $S$, and that $\Spec v'$ is valuative.
then $v'$ is a valuation of $S$. \hfill$\square$
\end{cor}

Let $S$ be a compact space, $v$ a valuation of $S$, and $T=[U]\in v$, where $U\subseteq X=\Spec v$ is a quasi-compact open subset.
Then the following statements are easy to see:
\begin{itemize}
\item the collection 
$$
v\vert_T=\{T'\in v\,|\,T'\subseteq T\}
$$
is a net on $T$ (due to \ref{prop-coveringandnets1topspace}), and is a valuation of $T$; in particular, $v$ is a net on $S$;
\item the homomorphism
$$
v\longrightarrow v\vert_T,\qquad T'\longmapsto T'\cap T
$$
gives rise to the open immersion $\Spec v\vert_T\hookrightarrow\Spec v$ corresponding to the inclusion map $U\hookrightarrow X$.
\end{itemize}

\begin{dfn}\label{dfn-valuationscompactHausdorffmorph}
A {\em morphism} $f\colon (S,v)\rightarrow (S',v')$ between valued compact spaces is a continuous mapping $f\colon S\rightarrow S'$ such that the following conditions are satisfied:
\begin{itemize}
\item[(a)] $f$ induces $f^{\ast}\colon v'\rightarrow v$ by $T\mapsto f^{-1}(T)$;
\item[(b)] $[\Spec f^{\ast}]=f$.
\end{itemize}
\end{dfn}

Note that a morphism of valued compact spaces $f\colon (S,v)\rightarrow (S',v')$ induces a valuative map $\Spec f^{\ast}\colon X=\Spec v\rightarrow X'=\Spec v'$ of coherent reflexive valuative spaces such that the diagram
$$
\xymatrix{X\ar[r]\ar[d]&X'\ar[d]\\ S\ar[r]_f&S'}
$$
commutes.
Now the following proposition is easy to see, essentially due to Stone duality (\S\ref{subsub-strcohtopsp}): 
\begin{prop}\label{prop-valuationscompactHausdorff1}
The functor $(S,v)\mapsto\Spec v$ gives a categorical equivalence from the category of valued compact spaces $($with the above-defined morphisms$)$ to the category of coherent reflexive valuative spaces with valuative quasi-compact maps.
The quasi-inverse to this functor is given by $X\mapsto ([X],v)$, where $v=\{[U]\,|\,U\in\QCOuv(X)\}$. \hfill$\square$
\end{prop}

\subsubsection{Valuations of locally Hausdorff spaces}\label{subsub-valuationslocallycompactlocallyHausdorff}
\begin{dfn}\label{dfn-valuationlocallyHausdorff}
Let $X$ be a locally Hausdorff space.

{\rm (1)} A {\em pre-valuation}\index{prevaluation@pre-valuation} $v=(\tau(v),\{v_S\}_{S\in\tau(v)})$ of $X$ consists of
\begin{itemize}
\item[{\rm (a)}] a net\index{net} (\ref{dfn-netsfortopologicalspaces} (2)) $\tau(v)$ of compact subspaces of $X$,
\item[{\rm (b)}] a valuation (\ref{dfn-valuationscompactHausdorff} (1)) $v_S$ of $S$ for each $S\in\tau(v)$,
\end{itemize}
satisfying the following condition:
\begin{itemize}
\item for any $S,S'\in\tau(v)$ with $S\subseteq S'$, we have $v_S=\{T\in v_{S'}\,|\,T\subseteq S\}$.
\end{itemize}

{\rm (2)} A pre-valuation $v=(\tau(v),\{v_S\}_{S\in\tau(v)})$ of $X$ is said to be {\em saturated}\index{prevaluation@pre-valuation!saturated prevaluation@saturated ---} if $\bigcup_{S\in\tau(v)}v_S=\tau(v)$ holds.

{\rm (3)} A pre-valuation $v=(\tau(v),\{v_S\}_{S\in\tau(v)})$ is called a {\em valuation} if it is saturated and, for any finitely many elements $S_1,\ldots,S_n\in\tau(v)$, $\bigcup_{i=1}^nS_i$ again belongs to $\tau(v)$ whenever it is Hausdorff as a subspace of $X$.

{\rm (4)} If $v$ is a valuation of $X$, then the pair $(X,v)$ is called a {\em valued locally Hausdorff space}\index{space@space (topological)!valued locally Hausdorff space@valued locally Hausdorff ---}.
\end{dfn}

Notice that:
\begin{itemize}
\item a locally Hausdorff space that admits a pre-valuation is locally compact;
\item if a pre-valuation $v$ of $X$ is saturated, then $v$ is determined from the net $\tau(v)$; indeed, for any $S\in\tau(v)$, one has $v_S=\{T\in\tau(v)\,|\,T\subseteq S\}$;
\item if $X$ is compact, then the notion of valuations defined here coincides with the previous one (\ref{dfn-valuationscompactHausdorff} (1)) (whence the abuse of terminology); indeed, if $v$ is a valuation of $X$ in the sense of \ref{dfn-valuationscompactHausdorff} (1), then $(v,\{v\vert_S\}_{S\in v})$ gives a valuation in the sense of \ref{dfn-valuationlocallyHausdorff}; conversely, given a valuation $v$ of $X$ in the latter sense, since $X$ is a finite union of elements of $\tau(v)$, $X$ itself belongs to the net $\tau(v)$.
\end{itemize}

\begin{rem}\label{rem-valuationlocallyHausdorff}
Notice that, if $(X,v)$ is a valued locally Hausdorff space, the net $\tau(v)$ does not, in general, give a distributive lattice, for it may fail to contain $X$ (the unit element). In other words, $\tau(v)$ is a distributive lattice if and only if $X$ is compact, in which case the valued locally Hausdorff space $(X,v)$ is a valued compact space.
\end{rem}

\begin{dfn}\label{dfn-valuationlocallyHausdorffmorph}
A {\em morphism} $f\colon (X,v)\rightarrow (X',v')$ between valued locally Hausdorff spaces is a continuous mapping $f\colon X\rightarrow X'$ such that, for any $S\in\tau(v)$ and $S'\in\tau(v')$ such that $f(S)\subseteq S'$, the map $f|_S\colon S\rightarrow S'$ with the valuations $v_S$ and $v_{S'}$ of $S$ and $S'$, respectively, is a morphism of valued compact spaces (\ref{dfn-valuationscompactHausdorffmorph}).
\end{dfn}

\subsubsection{Saturation and associated valuations}\label{subsub-extvaluationslocallycompactlocallyHausdorff}
\begin{dfn}\label{dfn-prevaluationsordering}
Let $X$ be a locally Hausdorff space $X$, and $v_1,v_2$ pre-valuations of $X$.
We say that $v_2$ is an {\em extension} of $v_1$ or that $v_1$ is a {\em restriction} of $v_2$, denoted as
$$
v_1\leq v_2,
$$
if $\tau(v_1)\subseteq\tau(v_2)$ and $(v_1)_S=(v_2)_S$ for any $S\in\tau(v_1)$.
\end{dfn}

\begin{prop}\label{prop-extensionprevaluation1}
Any pre-valuation of a compact space extends uniquely to a valuation.
\end{prop}

\begin{proof}
Let $X$ be a compact space, and $v=(\tau(v),\{v_S\}_{S\in\tau(v)})$ a pre-valuation of $X$.
Let $\til{v}$ be the distributive sublattice of $\mathbf{2}^X$ generated by all $v_S$'s.
More precisely, a subset $T\subseteq X$ belongs to $\til{v}$ if it is a finite union of finite intersections of elements in $\bigcup_{S\in\tau(v)}v_S$.
Then $\til{v}$ is a valuation of $X$ such that $\til{v}\vert_S=v_S$ for any $S\in\tau(v)$; indeed, one can glue the valuative spaces $\Spec v_S$ by open immersions (cf.\ \S\ref{subsub-valuationscompactHausdorff}).
It is then clear that $\til{v}=(\til{v},\{\til{v}\vert_T\}_{T\in\til{v}})$ extends $v$.
The uniqueness is clear.
\end{proof}

\begin{prop}\label{prop-saturationprevaluation}
Let $X$ be a locally Hausdorff space, and $v=(\tau(v),\{v_S\}_{S\in\tau(v)})$ a pre-valuation of $X$.
Set
$$
\tau(v^{\mathrm{sat}})=\bigcup_{S\in\tau(v)}v_S
$$
and for any $T\in\tau(v^{\mathrm{sat}})$
$$
v^{\mathrm{sat}}_T=\{T'\in\tau(v^{\mathrm{sat}})\,|\,T'\subseteq T\}.
$$
Then $v^{\mathrm{sat}}=(\tau(v^{\mathrm{sat}}),\{v^{\mathrm{sat}}_S\}_{S\in\tau(v^{\mathrm{sat}})})$ gives the smallest saturated extension of $v$ $($called the {\em saturation} of $v)$.
\end{prop}

\begin{proof}
It is clear that $\tau(v^{\mathrm{sat}})$ is a quasi-net\index{net!quasinet@quasi-{---}}.
To show that it is a net, take $S_1,S_2\in\tau(v^{\mathrm{sat}})$ and set $S=S_1\cap S_2$.
We want to show that $\tau(v^{\mathrm{sat}})\vert_S$ is a quasi-net.
Since $X$ is locally compact, $S$ is locally closed both in $S_1$ and $S_2$ and also in $X$.
This allows one to discuss locally on $X$ and to assume that $S$ is compact. 
Similarly to the proof of \ref{prop-extensionprevaluation1}, consider the distributive sublattice $v'_S$ of $\mathbf{2}^S$ generated by all $T'\in v_T$ for $T\in\tau(v)$ and $T\subseteq S$.
Then $v'_S$ gives a valuation of $S$, and hence is a net.
Since each element of $v'_S$ is a finite union of finite intersections of $T$'s as above, one verifies the conditions as in \ref{dfn-netsfortopologicalspaces} for any $x\in S$, showing that $\tau(v^{\mathrm{sat}})\vert_S$ is a quasi-net, as desired.
The rest of the proof is obvious.
\end{proof}

Let $X$ be a locally Hausdorff space, $v=(\tau(v),\{v_S\}_{S\in\tau(v)})$ a pre-valuation of $X$, and $C\subseteq X$ a subset.
Set
$$
v_C=(\{S\in\tau(v)\,|\,S\subseteq C\},\{v_S\}_{S\in\tau(v),\ S\subseteq C}).
$$
If $v_C$ defines a pre-valuation of $C$, we call it the {\em restriction} of $v$ to $C$.
The following assertions are easy to see; cf.\ Exercise \ref{exer-resrictionprevaluation}:
\begin{prop}\label{prop-associatedvaluationtopologicalspaces1}
Let $X$ be a locally Hausdorff space, $v=(\tau(v),\{v_S\}_{S\in\tau(v)})$ a pre-valuation of $X$.
Define $\tau(v^{\mathrm{val}})$ to be the set of all compact subsets $S\subseteq X$ such that the restriction, denoted consistently by $v^{\mathrm{sat}}_S$, of $v^{\mathrm{sat}}$ to $S$ exists.
For any $S\in\tau(v^{\mathrm{val}})$, set $v^{\mathrm{val}}_S$ to be the unique valuation of $S$ that extends $v^{\mathrm{sat}}_S$ {\rm (\ref{prop-extensionprevaluation1})}.
Then the data
$$
v^{\mathrm{val}}=(\tau(v^{\mathrm{val}}),\{v^{\mathrm{val}}_S\}_{S\in\tau(v^{\mathrm{val}})})
$$
is the minimal valuation that extends $v$ $($called the valuation associated to $v)$. \hfill$\square$
\end{prop}

\begin{thm}\label{thm-associatedvaluationtopologicalspaces2}
Let $X$ be a locally Hausdorff space, $v=(\tau(v),\{v_S\}_{S\in\tau(v)})$ a pre-valuation of $X$.
Define 
$$
\Spec v=\varinjlim_{S\in\tau(v)}\Spec v_S.
$$
Then $\Spec v$ is a reflexive valuative space with a canonical homeomorphism $[\Spec v]\cong X$.
Moreover, $\tau(v^{\mathrm{val}})$ coincides with the set of all subsets of $X$ that are images under the separation map $\sep$ of coherent open subsets of $\Spec v$. \hfill$\square$
\end{thm}

\begin{cor}\label{cor-associatedvaluationtopologicalspaces1}
Any pre-valuation $v$ of a locally Hausdorff space $X$ has an extension to a valuation.
Moreover, such an extension is unique. \hfill$\square$
\end{cor}

\subsubsection{Reflexive locally strongly compact valuative spaces}\label{subsub-RLSCVsp}
For a valued locally Hausdorff space $(X,v)$, we have constructed the reflexive valuative space $\mathscr{X}=\Spec v$ with the canonical map $\mathscr{X}\rightarrow X$ inducing a homeomorphism $[\mathscr{X}]\approx X$ (\ref{thm-associatedvaluationtopologicalspaces2}).
\begin{prop}\label{prop-RLSCVsp1}
The valuative space $\mathscr{X}=\Spec v$ is locally strongly compact.
Moreover, if $X$ is Hausdorff $($resp.\ compact$)$, then $\mathscr{X}$ is quasi-separated $($resp.\ coherent$)$.
\end{prop}

\begin{proof}
Since $\tau(v)=\{[U]\,|\,U\subseteq\mathscr{X}\ \textrm{coherent}\}$, for any open subset $W\subseteq[X]$ the restriction $v|_W$ of $v$ to $W$ exists, giving $\Spec v|_W\approx\sep^{-1}_{\mathscr{X}}(W)$ (cf.\ Exercise \ref{exer-resrictionprevaluation}).
Hence, to show that $\mathscr{X}$ is locally strongly compact, we may assume that $X$ is Hausdorff.

Since $\tau(v)$ is a net, there exists an open covering $\{U_{\alpha}\}_{\alpha\in L}$ consisting of coherent open subsets of $\mathscr{X}$ such that $[X]=\bigcup_{\alpha\in L}[U_{\alpha}]^{\circ}$, where $(\cdot)^{\circ}$ denotes the interior kernel.
Moreover, since $v$ is a valuation, we have $[U_{\alpha}]\in\tau(v)$ for any $\alpha\in L$.
For any $\alpha,\beta\in L$, the restriction $\tau(v)|_{[U_{\alpha}\cap U_{\beta}]}$ of $\tau(v)$ to $[U_{\alpha}\cap U_{\beta}]=[U_{\alpha}]\cap[U_{\beta}]$ is a quasi-net.
Hence there exists a family $\{V_{\lambda}\}_{\lambda\in\Lambda}$ of coherent open subsets of $U_{\alpha}\cap U_{\beta}$ such that $[U_{\alpha}\cap U_{\beta}]=\bigcup_{\lambda\in\Lambda}[V_{\lambda}]^{\circ}$.
Then $U_{\alpha}\cap U_{\beta}$ is, since it is reflexive and quasi-separated valuative space, locally strongly compact due to \ref{cor-locallystronglycompactconversecor1}.
Now by Exercise \ref{exer-locallystronglycompact2}, the topology on $[U_{\alpha}\cap U_{\beta}]$ coincides with the subspace topology from $[U_{\alpha}]$.
Since $[U_{\alpha}]$ and $[U_{\beta}]$ are compact, and since $X=[\mathscr{X}]$ is Hausdorff, we deduce that $[U_{\alpha}\cap U_{\beta}]=[U_{\alpha}]\cap[U_{\beta}]$ is compact, and hence, by \ref{thm-locallystronglycompacttheorem} (b), $U_{\alpha}\cap U_{\beta}$ is quasi-compact.
This shows that $X$ is quasi-separated, and hence is locally strongly compact by \ref{cor-locallystronglycompactconversecor1}.

Finally, if $[X]$ is compact, then by \ref{thm-locallystronglycompacttheorem} (b), $X$ is quasi-compact, hence coherent.
\end{proof}

\begin{thm}\label{thm-RLSCVsp1}
By $(X,v)\mapsto\Spec v$ we have a categorical equivalence from the category of valued locally Hausdorff spaces $($with the morphisms defined as in {\rm \ref{dfn-valuationlocallyHausdorffmorph}}$)$ to the category of reflexive locally strongly compact valuative spaces with valuative locally quasi-compact maps. The quasi-inverse to this functor is given by $\mathscr{X}\mapsto ([\mathscr{X}],v=(\tau(v),\{v_S\}_{S\in\tau(v)}))$, where $\tau(v)=\{[U]\,|\,U\subseteq\mathscr{X}\ \textrm{coherent}\}$ and, for any $[U]\in\tau(v)$, $v_{[U]}=\{[V]\subseteq[U]\,|\,V\in\QCOuv(U)\}$.
Moreover, $X$ in $(X,v)$ is Hausdorff $($resp.\ compact$)$ if and only if $\Spec v$ is quasi-separated $($resp.\ coherent$)$.
\end{thm}

\begin{proof}
We first show that any morphism $f\colon (X,v)\rightarrow (X',v')$ of valued locally Hausdorff spaces induces a valuative locally quasi-compact map $\Spec f\colon\Spec v\rightarrow\Spec v'$.
Since we already know from \ref{thm-associatedvaluationtopologicalspaces2} that $\tau(v)$ (resp.\ $\tau(v')$) is the set of all $[U]$'s by coherent open subsets $U$ of $\mathscr{X}=\Spec v$ (resp.\ $\mathscr{X}'=\Spec v'$), $f$ induces, for any $[U]\in\tau(v)$ and $[U']\in\tau(v')$ with $f([U])\subseteq[U']$, a valuative quasi-compact map $\Spec f|_{[U]}\colon U=\Spec v_{[U]}\rightarrow U'=\Spec v'_{[U']}$ such that $[\Spec f|_{[U]}]=f|_{[U]}$.
Hence, in view of \ref{cor-proplemlocallyquasicompactmapsgen}, we have the desired map $\Spec f\colon\Spec v\rightarrow\Spec v'$ by $\Spec f=\varinjlim_{S\in\tau(v)}\Spec f|_S$, which is valuative and locally quasi-compact.
In view of \ref{prop-RLSCVsp1}, we thereby have the desired functor by $(X,v)\mapsto\Spec v$ between the prescribed categories.

Now, consider, for a reflexive locally strongly compact valuative space $\mathscr{X}$, the data $([\mathscr{X}],v=(\tau(v),\{v_S\}_{S\in\tau(v)}))$ as above.
By \ref{thm-locallystronglycompacttheorem}, we know that $[\mathscr{X}]$ is a locally Hausdorff space.
Moreover, due to \ref{prop-valuationscompactHausdorff1}, each $v_{[U]}$ gives the valuation of $[U]$ that corresponds to $U$.
Hence $v$ gives a pre-valuation of $[\mathscr{X}]$, and is easily shown to be a valuation.
It is clear, then, that $\mathscr{X}\mapsto ([\mathscr{X}],v=(\tau(v),\{v_S\}_{S\in\tau(v)}))$ gives a functor.
By \ref{prop-valuationscompactHausdorff1} and a straightforward patching argument, this gives a quasi-inverse functor to the above functor.

Finally, the last statement follows from \ref{prop-RLSCVsp1}, \ref{cor-thmlocstringcompsp1a}, and \ref{cor-tubes2acor} (2).
\end{proof}
\index{valuation!valuation of a locally Hausdorff space@--- of a locally Hausdorff space|)}

\subsection{Some generalities on topoi}\label{sub-generalpretopostheory}
In this book, a {\em site}\index{site} will always mean a $\mathsf{U}$-site (\cite[Expos\'e II, (3.0.2)]{SGA4-1}).\footnote{We sometimes use such commonly-used expressions as `large' sites and `small' sites, which, as usual, have nothing to do with the set-theoretic size, such as $\mathsf{U}$-small.}
{\em Topos} in this book is always a Grothendieck topos; a {\em Grothendieck topos}\index{Grothendieck!topos@--- topos}\index{topos|see{Grothendieck topos}}, or more precisely, {\em Grothendieck} {\em $\mathsf{U}$-topos} is the category of $\mathsf{U}$-sheaves of sets (\cite[Expos\'e II, (2.1)]{SGA4-1}) on a site (\cite[Expos\'e IV, (1.1)]{SGA4-1}).
Our general reference to topos theory is \cite[Expos\'e IV]{SGA4-1} and \cite{MacMoe}.
We denote by $\TOPOI$ the $2$-category of topoi, that is, the $2$-category of which the objects are topoi, $1$-morphisms are morphisms of topoi, and for two $1$-morphisms $f,g\colon X\rightarrow Y$ the $2$-morphisms are the natural transformations from $f^{\ast}$ to $g^{\ast}$
(see \cite[Chap.\ I, \S1]{Hakim} for the generalities of $2$-categories).

\subsubsection{Spacial topoi}\label{subsub-toposasstopsp}\index{topos!spacial topos@spacial ---|(}
For a topological space $X$, we denote by $\top(X)$ the associated topos, that is, the category of $\mathsf{U}$-sheaves of sets over $X$\footnote{Notice that, since we only deal with the topological spaces in $\mathsf{U}$ (cf.\ \ref{subsub-categorynot}), the topos $\top(X)$ is actually a $\mathsf{U}$-topos.}; we call $\top(X)$ the {\em sheaf topos}\index{topos!sheaf topos@sheaf ---} of $X$.
The category $\Sets$ of sets in $\mathsf{U}$ (\S\ref{subsub-categorynot}) is a topos, as it is identical with the sheaf topos of the topological space consisting of one point (endowed with the unique topology).

The sheaf topos construction gives rise to a $2$-functor 
$$
\top\colon\Top\longrightarrow\TOPOI.
$$
A topos that belongs to the essential image of $\top$ is said to be {\em spacial}; that is, a topos $E$ is spacial if it is equivalent to the sheaf topos of a topological space.
Notice that the canonical map $X\rightarrow X^{\sob}$ (\S\ref{subsub-sober}) induces an equivalence of topoi
$$
\top(X)\stackrel{\sim}{\longrightarrow}\top(X^{\sob}).
$$
Thus a spacial topos is always equivalent to the sheaf topos of a sober topological space.

\begin{thm}[{\cite[Expos\'e IV, (4.2.3)]{SGA4-1}}]\label{thm-pointstopos1}
The $2$-functor 
$$
\top\colon\STop\longrightarrow\TOPOI
$$
is $2$-faithful$;$ that is, for any sober topological spaces\index{space@space (topological)!sober topological space@sober ---} $X$ and $Y$, the functor
$$
\Hom_{\STop}(X,Y)\longrightarrow\bHom_{\TOPOI}(\top(X),\top(Y))
$$
is an equivalence of categories, where the left-hand side is considered to be a discrete category $(${\rm \S\ref{subsub-groupoinds}}$)$. \hfill$\square$
\end{thm}
\index{topos!spacial topos@spacial ---|)}

\subsubsection{Points}\label{subsub-topoipopints}
Let us recall the notion of points of topoi.
Let $E$ be a topos. 
A {\em point}\index{point!point of a topos@--- (of a topos)} of $E$ is a morphism of topoi
$$
\xi\colon\Sets\longrightarrow E
$$
or, equivalently, a functor (so-called {\em fiber functor}\index{functor!fiber functor@fiber ---})
$$
\xi^{\ast}\colon E\rightarrow\Sets
$$
that commutes with finite projective limits and arbitrary inductive limits.
A topos $E$ is said to have {\em enough points}\index{enough points} if, for any arrow $u$ in $E$, $u$ is a monomorphism (resp.\ an epimorphism) if and only if $\xi^{\ast}(u)$ is injective (resp.\ a surjective) for any point $\xi$.

For a $\mathsf{U}$-small topos $E$ we denote respectively by 
$\pts(E)$ and $\open(E)$ the set of all isomorphism classes of points of $E$ and the set of all isomorphism classes of subobjects of a fixed final object of $E$.
For any $U\in\open(E)$ we set 
$$
|U|=\{\xi\in\pts(E)\,|\,\xi^{\ast}(U)\neq\emptyset\}.
$$
Then the collection $\{|U|\,|\,U\in\open(E)\}$ is stable under finite intersections and under arbitrary unions and hence gives a topology on the set $\pts(E)$.
We denote the resulting topological space by $\sp(E)$.
The topological space $\sp(E)$ is sober; if $E=\top(X)$ by a topological space $X$, then
$$
X^{\sob}\cong\sp(\top(X))
$$
(\cite[Expos\'e IV, (7.1.6)]{SGA4-1}).

\subsubsection{Localic topoi}\label{subsub-localictopoi}\index{topos!localic topos@localic ---|(}
A spacial topos\index{topos!spacial topos@spacial ---} $E$ satisfies the following condition (cf.\ \cite[\S5.3]{Johnstone2}):
\begin{itemize}
\item[{\bf (SG)}] $E$ is generated by subobjects of a fixed final object.
\end{itemize}
A (Grothendieck) topos that satisfies the condition {\bf (SG)} is said to be {\em localic}.
Notice that localic topoi are not always spacial.
\begin{thm}[{\rm cf.\ \cite[Theorem 7.25]{Johnstone2}}]\label{thm-localicspacial}
A localic topos $E$ is spacial if and only if $E$ has enough points.
\end{thm}

\begin{proof}
Here we give a sketch of the proof of the ``if'' part; see [loc.\ cit.] for the details.
Since $E$ is enough points, any open set $U$ of $X=\sp(E)$ corresponds bijectively to an element $U_E$ of $\open(E)$.
For any object $u$ of $E$, $U\mapsto\Hom_E(U_E,u)$ gives a sheaf on $X$.
The construction gives a functor $E\rightarrow\top(X)$.
Since $U_E$'s for any open $U$ of $X$ generate $E$, it follows that the functor gives a categorical equivalence.
\end{proof}

Let us finally remark that the notion of localic topoi is closely related with the theory of {\em locales}; cf.\ \cite[Chap.\ II]{Johnstone}.
\index{topos!localic topos@localic ---|)}

\subsubsection{Coherent topoi}\label{subsub-coherenttopos}
\index{coherent!coherent topos@--- topos|(}\index{topos!coherent topos@coherent ---|(}
In this subsection all topoi are considered with the canonical topology (\cite[Expos\'e II, (2.5)]{SGA4-1}), whenever considered as a site.
\begin{dfn}\label{dfn-coherenttopos1}{\rm 
(1) An object $X$ of a site (or a topos; cf.\ \cite[Expos\'e VI, (1.2)]{SGA4-2}) is said to be {\em quasi-compact}\index{quasi-compact!quasi-compact object of a site@--- object (of a site)} if any covering family $\{X_i\rightarrow X\}_{i\in I}$ has a finite covering subfamily.

(2) An object $X$ of a topos $E$ is said to be {\em quasi-separated}\index{quasi-separated!quasi-separated object of a topos@--- object (of a topos)} if, for any arrows $S\rightarrow X$ and $T\rightarrow X$ in $E$ with $S$ and $T$ quasi-compact, the fiber product $S\times_XT$ is quasi-compact.

(3) An object $X$ of a topos $E$ is said to be {\em coherent}\index{coherent!coherent object of a topos@--- object (of a topos)} if it is quasi-compact and quasi-separated.}
\end{dfn}

\begin{dfn}\label{dfn-coherenttopos12}{\rm 
Let $E$ be a topos.

(1) An arrow $f\colon X\rightarrow Y$ in $E$ is said to be {\em quasi-compact}\index{quasi-compact!quasi-compact arrow of a topos@--- arrow (of a topos)} if, for any arrow $Y'\rightarrow Y$ in $E$ with $Y'$ quasi-compact, the fiber product $X\times_YY'$ is quasi-compact.

(2) An arrow $f\colon X\rightarrow Y$ in $E$ is said to be {\em quasi-separated}\index{quasi-separated!quasi-separated arrow of a topos@--- arrow (of a topos)} if the diagonal arrow $X\rightarrow X\times_YX$ is quasi-compact.

(3) An arrow $f\colon X\rightarrow Y$ in $E$ is said to be {\em coherent}\index{coherent!coherent arrow of a topos@--- arrow (of a topos)} if it is quasi-compact and quasi-separated.}
\end{dfn}

We define basic notions for topoi based on finiteness conditions of objects and morphisms listed above:
\begin{dfn}\label{dfn-coherenttopos2}{\rm 
A topos $E$ is said to be {\em quasi-separated}\index{quasi-separated!quasi-separated topos@--- topos}\index{topos!quasi-separated topos@quasi-separated ---} (resp.\ {\em coherent}\index{coherent!coherent topos@--- topos}\index{topos!coherent topos@coherent ---}) if it satisfies the following conditions:
\begin{itemize}
\item[{\rm (a)}] there exists a generating full subcategory consisting of coherent objects; 
\item[{\rm (b)}] every object $X$ of $E$ is quasi-separated over the final object, that is, the diagonal morphism $X\rightarrow X\times X$ is quasi-compact.
\item[{\rm (c)}] the final object of $E$ is quasi-separated (resp.\ coherent).
\end{itemize}}
\end{dfn}

Let $E$ be a topos, and consider the following conditions: $E$ admits a generating full subcategory $C$ consisting of quasi-compact objects that is
\begin{itemize}
\item[(a)] stable by fiber products;
\item[(b)] stable by fiber products and by products of two objects;
\item[(c)] stable by any finite projective limits.
\end{itemize}
Then:
\begin{equation*}
\begin{split}
\textrm{(b)}&\Longleftrightarrow\textrm{$E$ is quasi-separated;}\\
\textrm{(c)}&\Longleftrightarrow\textrm{$E$ is coherent.}
\end{split}
\end{equation*}

\begin{dfn}\label{dfn-coherenttopos2loc}{\rm 
We say that the topos $E$ is {\em locally coherent}\index{topos!coherent topos@coherent ---!locally coherent topos@locally --- ---}\index{coherent!coherent topos@--- topos!locally coherent topos@locally --- ---} if it admits a generating full subcategory $C$ consisting of quasi-compact objects satisfying (a).}
\end{dfn}

Note that, in this situation, any object of $C$ is coherent (\cite[Expos\'e VI, (2.1)]{SGA4-2}).

\begin{dfn}\label{dfn-coherenttopos3}{\rm 
Let $f\colon E'\rightarrow E$ be a morphism of topoi.
We say $f$ is {\em quasi-compact}\index{quasi-compact!quasi-compact morphism of topoi@--- morphism (of topoi)}\index{morphism of topoi@morphism (of topoi)!quasi-compact morphism of topoi@quasi-compact ---} (resp.\ {\em quasi-separated}\index{quasi-separated!quasi-separated morphism of topoi@--- morphism (of topoi)}\index{morphism of topoi@morphism (of topoi)!quasi-separated morphism of topoi@quasi-separated ---}) if, for any quasi-compact (resp.\ quasi-separated) object $X$ of $E$, $f^{\ast}(X)$ is quasi-compact (resp.\ quasi-separated).
If $f$ is quasi-compact and quasi-separated, we say $f$ is {\em coherent}\index{coherent!coherent morphism of topoi@--- morphism (of topoi)}\index{morphism of topoi@morphism (of topoi)!coherent morphism of topoi@coherent ---}.}
\end{dfn}

\begin{prop}\label{prop-coherenttopos1}
Let $X$ be a topological space that admits an open basis consisting of quasi-compact open subsets, and consider the associated topos $\top(X)$.
Then $\top(X)$ is quasi-separated $($resp.\ locally coherent, resp.\ coherent$)$ if and only if $X$ is quasi-separated $($resp.\ locally coherent, resp.\ coherent$)$.\hfill$\square$
\end{prop}

Here we insert the important result by Deligne\index{Deligne, P.}:
\begin{thm}[{\cite[Expos\'e VI, (9.0)]{SGA4-2}}]\label{thm-deligne}
A locally coherent topos $E$ has enough points\index{enough points}.\hfill$\square$
\end{thm}

By \ref{thm-localicspacial} and \ref{thm-deligne}, we have:
\begin{cor}\label{cor-localicspacial2}
A coherent localic\index{topos!localic topos@localic ---} topos is spacial.\index{topos!spacial topos@spacial ---} \hfill$\square$
\end{cor}
\index{coherent!coherent topos@--- topos|)}\index{topos!coherent topos@coherent ---|)}

\subsubsection{Fibered topoi and projective limits}\label{subsub-fiberedtopos}
\begin{dfn}[{\cite[Expos\'e VI, (7.1.1)]{SGA4-2}}]\label{dfn-fiberedtopos1}{\rm 
Let $I$ be a category.
A fibered category (considered with a cleavage)
$$
p\colon E\longrightarrow I
$$
over $I$ is said to be a {\em fibered topos}\index{topos!fibered topos@fibered ---} or an {\em $I$-topos}\index{topos!Itopos@$I$-{---}} if every fiber $F_i$ $(i\in\obj(I))$ is a topos and for any arrow $f\colon i\rightarrow j$ of $I$ there exists a morphism of topoi $E_f\colon E_i\rightarrow E_j$ such that the pull-back morphism $f^{\ast}\colon E_j\rightarrow E_i$ by $f$ coincides with $E^{\ast}_f$.}
\end{dfn}

Hence, by the cleavage-construction of fibered categories (\cite[Expos\'e VI, \S8]{SGA1}), giving a fibered topos is equivalent to giving a functor
$$
E_{\bullet}\colon I\longrightarrow\TOPOI.
$$
We denote by $\TOPOI/I$ the $2$-category of $I$-topoi.

Recall that for a fibered category $p\colon E\rightarrow I$ a morphism $\varphi\colon x\rightarrow y$ over $f\colon i\rightarrow j$ is said to be {\em Cartesian} if, for any object $z$ with $p(z)=i$, the canonical map 
$$
\Hom_i(z,x)\longrightarrow\Hom_f(z,y),\qquad \alpha\mapsto\varphi\circ\alpha,
$$
is a bijection, where the left-hand set consists of arrows $z\rightarrow x$ over $\id_i$, and the right-hand set consists of arrows $z\rightarrow y$ over $f$.
An $I$-functor $E\rightarrow F$ of fibered categories over $I$ is said to be {\em Cartesian} if it maps Cartesian arrows to Cartesian arrows.
The category of Cartesian functors from $E$ to $F$ is denoted by 
$$
\Cart_{/I}(E,F).
$$
\begin{dfn}[{\cite[Expos\'e VI, (8.1.1)]{SGA4-2}}]\label{dfn-fiberedtopos2}{\rm 
Let $p\colon E\rightarrow I$ be a fibered topos.
A couple $(F,m)$ consisting of a topos $F$ and a cartesian morphism $m\colon F\times I\rightarrow E$ of topoi over $I$ is called a {\em projective limit of the fibered topos $F$} if for any topos $D$ the functor
$$
\bHom_{\TOPOI}(D,F)\longrightarrow\Cart_{\TOPOI/I}(D\times I,E)
$$
obtained by the composition of $\bHom_{\TOPOI}(D,F)\rightarrow\bHom_{\TOPOI/I}(D\times I,F\times I)$ followed by $m$ is an equivalence of categories.}
\end{dfn}

The projective limit is determined up to natural equivalences. 
We denote it by 
$$
\varprojLim E=\varprojLim{}_IE.
$$

In case $I$ is cofiltered and is essentially small (cf.\ \S\ref{subsub-finalcofinal}), we have the following more down-to-earth description of the projective limit. 
Let $p\colon E\rightarrow I$ be an $I$-topos, and $S$ the set of all cartesian arrows in $E$.
We denote by $\varinjLim{}_{I^{\opp}}E$ the category obtained from $E$ by inverting all arrows in $S$ (see \cite[Expos\'e VI, (6.2)]{SGA4-2} for details of the construction).
We endow with $\varinjLim{}_{I^{\opp}}E$ the weakest topology so that the canonical functor $E_i\rightarrow\varinjLim{}_{I^{\opp}}E$ for each $i\in\obj(I)$ is continuous.
Then $\varinjLim{}_{I^{\opp}}E$ becomes a site, and we have the natural equivalence
$$
\varprojLim{}_IE\cong(\varinjLim{}_{I^{\opp}}E)^{\sim}
$$
(\cite[Expos\'e VI, (8.2.3)]{SGA4-2}).

What we have described here can be furthermore boiled down to the following description, which starts from a simple observation: Assuming $\varprojLim{}_IE$ exists as a topos, there is, for each $i\in\obj(I)$, a projection $p_i\colon\varprojLim{}_IE\rightarrow E_i$ of topoi.
Then, for any object $\mathscr{F}$ of $\varprojLim{}_IE$, $\mathscr{F}_i=p_{i\ast}\mathscr{F}$ satisfies $(E_f)_{\ast}\mathscr{F}_i=\mathscr{F}_j$ for any arrow $f\colon i\rightarrow j$ of $I$.
Thus it is natural to define $\varprojLim{}_IE$ to be the category whose objects are collections $\{\mathscr{F}_i\}_{i\in\obj(I)}$ of objects $\mathscr{F}_i$ of $E_i$ that satisfu the compatibility $(E_f)_{\ast}\mathscr{F}_i=\mathscr{F}_j$ for any $f\colon i\rightarrow j$.
The verification of the equivalence with the above construction is left to the reader.

\begin{thm}[{\cite[Expos\'e VI, (8.3.13)]{SGA4-2}}]\label{thm-fiberedtopos10}
Let $p\colon E\rightarrow I$ be an $I$-topos, where $I$ is cofiltered and essentially small.
Suppose that$:$
\begin{itemize}
\item each fiber $E_i$ $(i\in\obj(I))$ is a coherent topos\index{coherent!coherent topos@--- topos}\index{topos!coherent topos@coherent ---} {\rm (\ref{dfn-coherenttopos2})}$;$
\item for any $f\colon i\rightarrow j$ the morphism $E_f\colon E_i\rightarrow E_j$ of topoi is coherent\index{coherent!coherent morphism of topoi@--- morphism (of topoi)}\index{morphism of topoi@morphism (of topoi)!coherent morphism of topoi@coherent ---} {\rm (\ref{dfn-coherenttopos3})}.
\end{itemize}
Then the projective limit $\varprojLim{}_IE$ is a coherent topos, and for each $i\in\obj(I)$ the canonical projection $\varprojLim{}_IE\rightarrow E_i$ is coherent.
Moreover, if we denote by $E_{\coh}$ the full subcategory of $E$ consisting of objects that are coherent in their fibers, then $E_{\coh}$ is a fibered category over $I$, and $\varinjLim{}_{I^{\opp}}E_{\coh}$ $($defined similarly as above$)$ is canonically equivalent to the category $(\varprojLim{}_IE)_{\coh}$, the full subcategory of $\varprojLim{}_IE$ consisting of coherent objects.\hfill$\square$
\end{thm}

\subsubsection{Projective limit of spacial topoi}\label{subsub-projectivelimitspacialtopoi}
\index{topos!spacial topos@spacial ---|(}
Let $I$ be a cofiltered and essentially small category, and consider a functor $X_{\bullet}\colon I\rightarrow\STop$, which we denote by $i\mapsto X_i$.
Then one can consider the projective limit $X_{\infty}=\varprojlim_IX_i$ in the category $\STop$ (in fact, it is easy to see that the projective limit of sober spaces taken in the category $\Top$ is sober).
We are interested in comparing the topos $\top(X_{\infty})$ and the topos-theoretic projective limit $\varprojLim{}_I\top(X_i)$; note that by the functor $X_{\bullet}$ we have the $I$-topos $\top(X_{\bullet})\rightarrow I$ with the fiber over $i$ being $\top(X_i)$.
Notice that these two topoi may not be equivalent in general, for the topos $\varprojLim{}_I\top(X_i)$ may not be spacial in general.

Let us first discuss the topos theoretic limit.
\begin{thm}\label{thm-localiclimitlocalic}
Let $I$ be a cofiltered and essentially small category, and consider a functor $X_{\bullet}\colon I\rightarrow\STop$. 

{\rm (1)} The topos $\varprojLim{}_I\top(X_i)$ is localic\index{topos!localic topos@localic ---}. 

{\rm (2)} The topological space $\sp(\varprojLim{}_I\top(X_i) ) $ is homeomorphic to $X_{\infty}=\varprojlim_IX_i$. 
\end{thm}

\begin{proof}
Set $E= \varprojLim{}_I\top(X_i)$. 
As for (1), by the definition of the Grothendieck topology of $E$, $E$ is generated by all $p_i^{\ast}(U_i)$ for all $i\in\obj(I)$ and subobjects $U_i$ of a final object of $\top (X_i)$, where $p_i\colon E\rightarrow\top(X_i)$ is the canonical projection. Thus $E$ is localic.  
For (2), consider the canonical map 
$$
F\colon\sp(E)\longrightarrow\varprojlim_{i\in\obj(I)}\sp(\top(X_i)), 
$$
where each $\sp(\top(X_i))$ is, since $X_i$ is sober, identified with $X_i$ and the right hand side is given the projective limit topology. 
We claim that $F$ is bijective: Indeed, since each category of the points of $\top(X_i)$ is discrete, any elements in the right hand side lifts to a projective system of points of topoi $\top(X_i)$, and hence to a point of $E$. 
Thus $F$ is surjective, and the injectivity is shown in a similar way. 
With the same notations as above, subobjects of the form $p_i^{\ast}(U_i)$ of the final object of $E$ generates the topology of $\sp(E)$, thus the map $F$ gives a homeomorphism with respect to this topology.
\end{proof}

\begin{rem}\label{rem-localiclimitlocalic}{\rm 
Notice that the above proof of (1) shows that the projective limit of localic topoi is again localic.}
\end{rem}

\begin{cor}\label{cor-fiberedtopos213} 
For a projective system of sober spaces $X_{\bullet}\colon I\rightarrow\STop$, the limit $\varprojLim{}_I\top(X_i)$ is spacial if and only if it has enough points. 
If it is spacial, the limit is equivalent to $\top(\varprojlim{}_I X_i)$. \hfill$\square$
\end{cor} 

An important case is provided by the following: Suppose, moreover, that 
\begin{itemize}
\item each $X_i$ is coherent and sober, and each transition map $X_i\rightarrow X_j$ is quasi-compact. 
\end{itemize}
Then, by \ref{thm-fiberedtopos10}, $\varprojLim{}_I\top(X_i)$ is coherent; applying \ref{cor-localicspacial2}, \ref{thm-fiberedtopos10}, and \ref{cor-fiberedtopos213}, we conclude that $\varprojLim{}_I\top(X_i)$ is canonically equivalent to $\top(\varprojlim{}_I X_i)$, and $\varprojlim{}_I X_i$ is coherent. 
Thus we have given a topos-theoretic proof of the first half of the following statement, which recasts \ref{thm-projlimcohsch1}:
\begin{cor}\label{cor-fiberedtopos21}
Let $\{X_i\}_{i\in I}$ be a projective system of topological spaces indexed by a directed set such that each $X_i$ is coherent sober and each transition map $X_i\rightarrow X_j$ is quasi-compact.
Set $X_{\infty}=\varprojlim_{i\in I}X_i$.
Then $X_{\infty}$ is coherent and sober.
If, moreover, each $X_i$ is non-empty, then $X_{\infty}$ is non-empty.
\end{cor}

\begin{proof}
The first assertion has been already shown above.
For the second, we use \ref{prop-coherentprojectivesystemindlimits1} shown later independently:
$$
\Gamma (X, \Z) = \varinjlim _I \Gamma (X_i , \Z).
$$
Notice that, if each $X_i$ is non-empty, then the inductive limit is clearly non-zero, whence the non-emptiness of $X$ follows.
\end{proof}

\begin{rem}\label{rem-variantringed}{\rm 
There is the `ringed'-version of the above argument, that is, we have the notion of {\em ringed $I$-topoi} and their projective limits as a ringed topos.
The argument is quite similar to as above.
For the details, see \cite[Expos\'e VI, (8.6)]{SGA4-2}.}
\end{rem}
\index{topos!spacial topos@spacial ---|)}

\subsubsection{Quasi-compact topoi and projective limits}\label{subsub-quasicompacttopos}
\index{topos!quasi-compact topos@quasi-compact ---|(}
Finally, let us include a few facts on quasi-compact topoi and their projective limits.

\begin{dfn}\label{dfn-quasicompacttopos1}{\rm 
A topos $E$ is said to be {\em quasi-compact} if it admits a quasi-compact final object.}
\end{dfn}

\begin{rem}\label{rem-quasicompacttopos1}{\rm 
In literature (e.g.\ \cite{Johnstone3}), quasi-compact topoi are called {\em compact topoi}. 
In this book, however, we prefer to use the terminology {\em quasi-compact}, which is consistent with the standard usage in algebraic geometry.}
\end{rem}

Consider the $2$-category $\LocTOPOI$ of localic topoi\index{topos!localic topos@localic ---} (\S\ref{subsub-localictopoi}) and the inclusion functor $\LocTOPOI\hookrightarrow\TOPOI$.
By collecting objects generated by subobjects of a fixed final objects for each topos $E$, we have a right $2$-adjoint functor $E\mapsto E^{\loc}$ from $\TOPOI$ to $\LocTOPOI$.
This implies that the formation of projective limits in $\TOPOI$ and $\LocTOPOI$ are compatible, and that, in particular, projective limits of localic topoi are again localic.
Notice that a topos $E$ is quasi-compact if and only if so is $E^{\loc}$.

\begin{thm}
Let $p\colon E\rightarrow I$ be an $I$-topos, where $I$ is cofiltered and essentially small.
Suppose that each fiber $E_i$ $(i\in\obj(I))$ is a quasi-compact topos. 

{\rm (1)} The projective limit $\varprojLim{}_IE$ is a quasi-compact topos. 

{\rm (2)} If each $E_i$ is non-empty, $\varprojLim{}_IE$ is non-empty. 
\end{thm}

\begin{proof} When each $E_i $ for $i\in\obj(I)$ is localic, the claim is proven in \cite[Theorem 2.3 \& Cor.\ 2.4]{WY}.
In general, by what we have remarked above, one can reduce to this case.
\end{proof}

It follows from the adjunction between $\top\colon\STop\rightarrow\LocTOPOI$ and $\sp\colon\LocTOPOI\rightarrow\STop$ that $\sp$ preserves projective limits.
Hence, as a corollary of theorem and the following easy lemma, we obtain a topos-theoretic proof of \ref{thm-gabbertheoremprojectivelimits}:
\begin{lem}
Let $E$ be a quasi-compact localic topos.

{\rm (1)} If $E$ is non-empty, $\sp(E)$ is non-empty.  

{\rm (2)} The topological space $\sp(E)$ is quasi-compact. \hfill$\square$
\end{lem}
\index{topos!quasi-compact topos@quasi-compact ---|)}

\addcontentsline{toc}{subsection}{Exercises}
\subsection*{Exercises}
\begin{exer}\label{exer-soberT1counterexample}
Let $X$ be an infinite set, and topologize it in such a way that a subset $Y\subseteq X$ is closed if and only if either it is $X$ itself, or is a finite subset.
Then show that the topological space $X$ is $\mathrm{T}_1$, but is not sober.
\end{exer}

\begin{exer}\label{exer-quasiseparatedtopsp}{\rm 
Show that a topological space $X$ is quasi-separated if and only if the diagonal mapping $\Delta\colon X\rightarrow X\times X$ is quasi-compact.}
\end{exer}

\begin{exer}\label{exer-projectivelimitaffine}
Let $\{X_i,p_{ij}\}_{i\in I}$ be a projective system of topological spaces indexed by a directed set $I$.
We assume:
\begin{itemize}
\item[(a)] each $X_i$ is the underlying topological space of a coherent ($=$ quasi-compact and quasi-separated) scheme;
\item[(b)] each transition map $p_{ij}\colon X_j\rightarrow X_i$ for $i\leq j$ is the underlying continuous mapping of an affine morphism\index{affine!affine morphism@--- morphism} of schemes.
\end{itemize}
Then show that the projective limit $X=\varprojlim_{i\in I}X_i$ (taken in $\Top$) is a coherent sober topological space.
\end{exer}

\begin{exer}\label{exer-latticespec1}{\rm 
Let $A$ be a distributive lattice. 

(1) A {\em covering} of $\alpha\in A$ is a finite subset $C$ of $A$ such that $\bigvee_{\beta\in C}\beta=\alpha$. 
Show that $A$ with this notion of coverings is a site.  The associated topos is denoted by $\spec A$.

(2) Show that the topos $\spec A$ is coherent.

(3) Show that the topos associated to $\Spec A$ is canonically equivalent to $\spec A$.}
\end{exer}

\begin{exer}\label{exer-latticespec2}{\rm 
Show that any coherent sober space is homeomorphic to a projective limit of finite sober spaces.}
\end{exer}

\begin{exer}\label{exer-Booleanlattices}{\rm

(1) Let $A$ be a Boolean distributive lattice.
Then show that $\Spec A$ is a profinite set and is homeomorphic to the spectrum of $A$ regarded as a commutative ring.

(2) Show that $A\mapsto\Spec A$ gives a categorical equivalence between the category of Boolean distributive lattices and the opposite category of totally disconnected Hausdorff spaces.}
\end{exer}

\begin{exer}\label{exer-generizationmaximal}{\rm 
Let $X$ be a coherent sober space.

(1) Let $\{x_i\}_{i\in I}$ be a system of points of $X$ indexed by a directed set $I$ such that, if $i\leq j$, $x_j$ is a generization of $x_i$.
Show that there exists a point $x\in X$ that is a generization of all $x_i$'s.

(2) Show that for any $x\in X$ the set $G_x$ of all generizations of $x$ contains a maximal generization of $x$.}
\end{exer}

\begin{exer}\label{exer-coherentprojlimopencovering}{\rm 
Let $\{X_i,p_{ij}\}_{i\in I}$ be a filtered projective system of coherent sober topological spaces and quasi-compact maps indexed by a directed set $I$, $X=\varprojlim_{i\in I}X_i$, and $i\in I$ an index.
Let $U_1,\ldots,U_n$ be finitely many quasi-compact open subsets of $X_i$ such that $X=\bigcup^n_{k=1}p^{-1}_i(U_k)$, where $p_i\colon X\rightarrow X_i$ is the projection map.
Then show that there exists $j\in I$ with $i\leq j$ such that $X_j=\bigcup^n_{k=1}p^{-1}_{ij}(U_k)$.}
\end{exer}

\begin{exer}\label{exer-limitcoherentsoberspaceconnected}{\rm 
Let $\{X_i,p_{ij}\}_{i\in I}$ be a projective system of coherent sober spaces with quasi-compact surjective transition maps indexed by a directed set $I$.
Show that $X=\varprojlim_{i\in I}X_i$ is connected if and only if all $X_i$ for $i\in I$ are connected.}
\end{exer}

\begin{exer}\label{exer-valuativespaceT1}{\rm Show that for any valuative space $X$ the separated quotient $[X]$ is a $\mathrm{T}_1$-space.}
\end{exer}

\begin{exer}\label{exer-locallystronglycompact1}{\rm Let $X$ be a locally strongly compact valuative space.

(1) Show that any overconvergent open subset $U\subseteq X$ is locally strongly compact.

(2) Show that for any quasi-compact open immersion $U\hookrightarrow X$ and any locally strongly compact open subset $V\subseteq X$, $U\cap V$ is locally strongly compact.
In particular, an open subset $U\subseteq X$ is locally strongly compact if the inclusion $U\hookrightarrow X$ is quasi-compact.}
\end{exer}

\begin{exer}\label{exer-locallystronglycompact2}{\rm (1) Let $X$ be a locally compact (hence locally Hausdorff) space, $Y$ a locally Hausdorff space, and $f\colon X\rightarrow Y$ a continuous injective map.
Then show that $f$ induces a homeomorphism onto its image endowed with the subspace topology induced from the topology of $Y$.

(2) Suppose, moreover, that $Y$ is locallly compact.
Then show that $f(X)$ is locally closed in $Y$.

(3) Let $X$ be a locally strongly compact valuative space, and $U\subseteq X$ a locally strongly compact open subset.
Then the injective map $[U]\rightarrow[X]$ maps $[U]$ homeomorphically onto a locally closed subspace of $[X]$.}
\end{exer}

\begin{exer}[Structure Theorem]\label{exer-locallystronglycompact3}{\rm Let $X$ be a locally strongly compact valuative space, and $U\subseteq X$ a locally strongly compact open subset.

(1) The inclusion map $j\colon U\rightarrow X$ is quasi-compact if and only if $[j]\colon[U]\hookrightarrow[X]$ is a closed immersion.

(2) In general, there exists an overconvergent open subset $Z\subseteq X$ containing $U$ such that the inclusion map $U\hookrightarrow Z$ is quasi-compact.}
\end{exer}

\begin{exer}\label{exer-locallystronglycompact4}{\rm Show that any finite intersection of locally strongly compact open subsets of a locally strongly compact valuative space is again locally strongly compact.}
\end{exer}

\begin{exer}\label{exer-locallystronglycompact5}{\rm Let $X$ be a locally strongly compact valuative space such that $[X]$ is Hausdorff. 
Then show that $X$ is quasi-separated.
In particular, $X$ is coherent if and only if $[X]$ is compact.}
\end{exer}

\begin{exer}\label{exer-valuationcharacterization}{\rm 
Show that a distributive sublattice $v$ of $\mathbf{2}^S$ for a compact space $S$ gives a valuation of $S$ if and only if the following conditions are satisfied:
\begin{itemize}
\item[(a)] $\emptyset,S\in v$ and all element of $v$ are compact;
\item[(b)] for any $x\in S$ the family $\{T\in v\,|\,x\in T^{\circ}\}$ forms a fundamental system of neighborhoods of $x$;
\item[(c)] for any prime filter (cf.\ \S\ref{subsub-strcohtopsp}) $F\subseteq v$, there exists a unique maximal prime filter containing $F$;
\item[(d)] for two distinct maximal prime filters $F,F\subseteq v$, we have $\bigcap_{T\in F}T\neq\bigcap_{T\in F'}T$;
\item[(e)] for any prime filter $F\subseteq v$ and a maximal filter $\til{F}$ containing $F$, the set of all prime filters between $F$ and $\til{F}$ is totally ordered with respect to the inclusion order.
\end{itemize}}
\end{exer}

\begin{exer}\label{exer-prevaluationsaturatedintersection}{\rm 
Let $v=(\tau(v),\{v_S\}_{S\in\tau(v)})$ be a pre-valuation of a locally Hausdorff space $X$, and $S_1,\ldots,S_n\in\tau(v)$.
Show that, if $\bigcap_{i=1}^nS_i$ is Hausdorff, then $\bigcap_{i=1}^nS_i\in\tau(v)$.}
\end{exer}

\begin{exer}\label{exer-resrictionprevaluation}{\rm 
Let $X$ be a locally Hausdorff space, $v=(\tau(v),\{v_S\}_{S\in\tau(v)})$ a pre-valuation of $X$.

(1) Show that, for any finitely many $S_1,\ldots,S_n\in\tau(v)$, the restriction of $v$ to the intersection $C=\bigcap_{i=1}^nS_i$ exists.

(2) Show that, for any finitely many $S_1,\ldots,S_n\in\tau(v)$, the union $C=\bigcup_{i=1}^nS_i$ is locally compact, and the restriction of $v$ to $C$ exists.

(3) Suppose $v$ is saturated. Then show that, for any open subset $U\subseteq X$, the restriction of $v$ to $U$ exists.}
\end{exer}


\section{Homological algebra}\label{sec-homologicalalgebrawithlimits}
In this section, we discuss two topics on homological algebra.
The first topic is on inductive and projective limits of sheaves and their cohomologies (\S\ref{sub-directlimits} and \S\ref{sub-projectivelimits}).
Given a projective system of topological spaces and an inductive system of sheaves on them, one has the inductive limit sheaf on the projective limit space.
In this situation, we will give a general recipe to calculate the cohomologies of the inductive limit sheaf (\S\ref{subsub-sheavestopologicalspacelimits} and \S\ref{subsub-cohsheavestopologicalspacelimits}).
This subsection also discuss Noetherness of inductive limit rings (\S\ref{subsub-directlimitsnoetherianness}). 
As for projective limit sheaves, we discuss the so-called {\em Mittag-Leffler condition}\index{Mittag-Leffler condition} and some of its consequences.

The second topic, discussed in \S\ref{sub-cohringsmodules}, is on {\em coherent} rings and modules.
Here, a ring $A$ is said to be coherent if every finitely generated ideal is finitely presented, or equivalently, the full subcategory of the category of $A$-modules consisting of finitely presented modules is an abelian subcategory (\ref{prop-cohringsmodules1}). 

\medskip\noindent
{\bf Convention.} {\sl Throughout this book, whenever we say $A$ is a ring, we always mean that $A$ is a commutative ring having the multiplicative unit $1=1_A$, unless otherwise clearly stated$;$ we also assume that any ring homomorphism $f\colon A\rightarrow B$ is unitary, that is, maps $1$ to $1$.
Here are other conventions$:$
\begin{itemize}
\item for a ring $A$ we denote by $\Frac(A)$ the total ring of fractions of $A;$
\item for a ring $A$ the Krull dimension of $A$ is denoted by $\dim(A);$
\item when $A$ is a local ring, we denote by $\m_A$ the unique maximal ideal of $A;$
\item a ring homomorphism $f\colon A\rightarrow B$ between local rings is said to be {\em local} if $f(\m_A)\subseteq\m_B$.
\end{itemize}}

\subsection{Inductive limits}\label{sub-directlimits}\index{limit!inductive limit@inductive ---|(}
\subsubsection{Preliminaries}\label{subsub-directlimitspre}
First, we collect basic known facts on inductive limits of rings and modules, which we quote, without proofs, mainly from \cite{Bourb3}.

\begin{prop}[{\cite[Chap.\ I, \S10.3, Prop.\ 3]{Bourb3}}]\label{prop-directlimits1}
Let $\{A_i,\phi_{ij}\}_{i\in I}$ be an inductive system of rings indexed by a directed set $I$, and set $A=\varinjlim_{i\in I}A_i$.

{\rm (1)} If $A_i\neq 0$ for each $i\in I$, then $A\neq 0$.

{\rm (2)} If each $A_i$ is an integral domain, then so is $A$.

{\rm (3)} If each $A_i$ is a field, then so is $A$. \hfill$\square$
\end{prop}

Notice that (1) is the direct consequence of our convention; since each transition map $\phi_{ij}$ maps $1_{A_i}$ to $1_{A_j}$, we immediately conclude that the inductive limit $A$ has the element $1_A$ not equal to $0$.
Notice also that (1) is the basis for the following well-known fact (somewhat similar to \ref{thm-projlimcohsch1}): if $\{X_i,f_{ij}\}_{i\in I}$ is a projective system of non-empty affine schemes\index{affine!affine scheme@--- scheme} indexed by a directed set, then the projective limit $X=\varprojlim_{i\in I}X_i$ in the category of schemes exists and is affine and non-empty (cf.\ Exercise \ref{exer-projectivelimitaffine}).

\begin{prop}[{\cite[Chap.\ II, \S3, Exer.\ 16]{Bourb1}}]\label{prop-directlimits11}
Let $\{A_i,\phi_{ij}\}_{i\in I}$ be an inductive system of local rings and local homomorphisms indexed by a directed set $I$.

{\rm (1)} The inductive limit $A=\varinjlim_{i\in I}A_i$ is a local ring with the maximal ideal $\m_A=\varinjlim_{i\in I}\m_{A_i}$.

{\rm (2)} Let $k_i=A_i/\m_{A_i}$ be the residue field of $A_i$ for each $i\in I$.
Then $k=\varinjlim_{i\in I}k_i$ gives the residue field of $A$.

{\rm (3)} If, moreover, $\m_{A_j}=\m_{A_i}A_j$ for $j\geq i$, then $\m_A=\m_{A_i}A$ for any $i\in I$. \hfill$\square$
\end{prop}

Let $\textrm{{\boldmath $A$}}=\{A_i,\phi_{ij}\}_{i\in I}$ be an inductive system of rings.
By an {\em inductive system of {\boldmath $A$}-modules} we mean an inductive system $\textrm{{\boldmath $M$}}=\{M,f_{ij}\}_{i\in I}$ such that each $M_i$ is an $A_i$-module and that each $f_{ij}\colon M_i\rightarrow M_j$ is compatible with $\phi_{ij}\colon A_i\rightarrow A_j$, that is, for $x\in M_i$ and $a\in A_i$ we have $f_{ij}(ax)=\phi_{ij}(a)f_{ij}(x)$.
For another inductive system {\boldmath $M'$} of {\boldmath $A$}-modules, we can define the notion of homomorphism $\textrm{{\boldmath $M$}}\rightarrow\textrm{{\boldmath $M'$}}$ in an obvious manner.

\begin{prop}[{\cite[Chap.\ II, \S6.2, Prop.\ 3]{Bourb3}}]\label{prop-directlimits2}
Let $\textrm{{\boldmath $M$}}=\{M,f_{ij}\}_{i\in I}$, $\textrm{{\boldmath $M'$}}=\{M',f'_{ij}\}_{i\in I}$, and $\textrm{{\boldmath $M''$}}=\{M'',f''_{ij}\}_{i\in I}$ be filtered inductive systems of {\boldmath $A$}-modules indexed by a directed set $I$.
Consider an exact sequence
$$
\textrm{{\boldmath $M'$}}\longrightarrow\textrm{{\boldmath $M$}}\longrightarrow\textrm{{\boldmath $M''$}}
$$
$($that is, $M'_i\rightarrow M_i\rightarrow M''_i$ is exact for every $i)$.
Then the induced sequence
$$
\varinjlim_{i\in I}M'_i\longrightarrow\varinjlim_{i\in I}M_i\longrightarrow\varinjlim_{i\in I}M''_i
$$ is exact.\hfill$\square$
\end{prop}

The following proposition is a formal implication of the fact that the inductive limit functor has the right adjoint (cf.\ \S\ref{subsub-limdefuniv}):
\begin{prop}[{cf.\ \cite[Chap.\ II, \S6.4, Prop.7]{Bourb3}}]\label{prop-directlimits3}
Let $\textrm{{\boldmath $M$}}=\{M,f_{ij}\}_{i\in I}$, $\textrm{{\boldmath $M'$}}=\{M',f'_{ij}\}_{i\in I}$, and $\textrm{{\boldmath $M''$}}=\{M'',f''_{ij}\}_{i\in I}$ be inductive systems of {\boldmath $A$}-modules, and suppose we  are given homomorphisms
$$
\textrm{{\boldmath $M'$}}\longleftarrow\textrm{{\boldmath $M$}}\longrightarrow\textrm{{\boldmath $M''$}}
$$
of {\boldmath $A$}-modules.
Then the canonical map
$$
\varinjlim_{i\in I}M'_i\otimes_{M_i}M''_i\longrightarrow(\varinjlim_{i\in I}M'_i)\otimes_{(\varinjlim_{i\in I}M_i)}(\varinjlim_{i\in I}M''_i)
$$
is an isomorphism.\hfill$\square$
\end{prop}

\subsubsection{Inductive limits and Noetherness}\label{subsub-directlimitsnoetherianness}
Let us include here a useful technique invented by M.\ Nagata\index{Nagata, M.} (cf.\ \cite[(43.10)]{Nagata1}) to show that an inductive limit ring in a certain situation is Noetherian:
\begin{prop}\label{prop-nagataquintessence}
Let $\{R_{\alpha}\}_{\alpha\in L}$ be a filtered inductive system of rings, and $R=\varinjlim_{\alpha\in L}R_{\alpha}$.
Suppose that the following conditions are satisfied$:$
\begin{itemize}
\item[{\rm (a)}] $R_{\alpha}$ is a Noetherian local ring for any $\alpha\in L;$
\item[{\rm (b)}] for any $\alpha\leq\beta$ the transition map $R_{\alpha}\rightarrow R_{\beta}$ is local and flat$;$
\item[{\rm (c)}] for any $\alpha\leq\beta$ we have $\m_{\alpha}R_{\beta}=\m_{\beta}$ $($where $\m_{\alpha}$ denotes the maximal ideal of $R_{\alpha})$.
\end{itemize}
Then $R$ is Noetherian.
\end{prop}

\begin{proof}
By (c) the maximal ideal $\m$ of $R$ coincides with $\m_{\alpha}R$ for any $\alpha\in L$; in particular, $\m$ is a finitely generated ideal of $R$.
Hence the $\m$-adic completion $\widehat{R}$ of $R$ exists (cf.\ \ref{prop-Iadiccompletioncomplete1} below); moreover, since $\m\widehat{R}$ is finitely generated, $\widehat{R}$ is Noetherian (cf.\ \cite[(29.4)]{Matsu}).
Since $R_{\alpha}/\m^{k+1}_{\alpha}\rightarrow R/\m^{k+1}$ for any $\alpha\in L$ and $k\geq 0$ is faithfully flat, by local criterion of flatness\index{local criterion of flatness} (cf.\ \S\ref{subsub-localcriterionofflatness} below), we deduce that $\widehat{R}_{\alpha}\rightarrow\widehat{R}$ is faithfully flat for any $\alpha\in L$.

Now, to prove the proposition, it suffices to show that any increasing sequence
$$
J_0\subseteq J_1\subseteq J_2\subseteq\cdots
$$
of finitely generated ideals of $R$ is stationary.
Since $\widehat{R}$ is Noetherian, the induced sequence
$$
J_0\widehat{R}\subseteq J_1\widehat{R}\subseteq J_2\widehat{R}\subseteq\cdots
$$
is stationary, that is, there exists $N\geq 0$ such that $J_n\widehat{R}=J_m\widehat{R}$ whenever $n,m\geq N$.
We want to show the equality $J_n=J_m$ for $n,m\geq N$.
Hence we only have to show the following: if $J,J'$ are finitely generated ideals of $R$ such that $J\subseteq J'$ and $J\widehat{R}=J'\widehat{R}$, then $J=J'$.
Take $\alpha\in L$ sufficiently large so that we have ideals $J_{\alpha},J'_{\alpha}\subseteq R_{\alpha}$ with $J_{\alpha}R=J$ and $J'_{\alpha}R=J'$.
Since $J_{\alpha}\widehat{R}=J\widehat{R}=J'\widehat{R}=J'_{\alpha}\widehat{R}$, we have $J_{\alpha}=J_{\beta}$ (since $\widehat{R}$ is faithfully flat over $R_{\alpha}$), and thus $J=J'$, as desired.
\end{proof}

\begin{cor}\label{cor-propnagataquintessencevar}
Let $\{R_{\alpha}\}_{\alpha\in L}$ be a filtered inductive system of rings, and $R=\varinjlim_{\alpha\in L}R_{\alpha}$.
Suppose that the following conditions are satisfied$:$
\begin{itemize}
\item[{\rm (a)}] $R_{\alpha}$ is a Noetherian local ring for any $\alpha\in L;$
\item[{\rm (b)}] for any $\alpha\leq\beta$ the transition map $R_{\alpha}\rightarrow R_{\beta}$ is local and formally smooth $($that is, $R_{\beta}$ is $\m_{R_{\beta}}$-smooth over $R_{\alpha}$ in the terminology as in {\rm \cite[\S28]{Matsu}}$);$
\item[{\rm (c)}] the set of integers $\{\dim(R_{\alpha})\}_{\alpha\in L}$ is bounded.
\end{itemize}
Then $R$ is Noetherian.
\end{cor}

\begin{proof}
Replacing $\{R_{\alpha}\}_{\alpha\in L}$ by a cofinal subsystem if necessary, we may assume that the numbers $\dim(R_{\alpha})$ are all equal.
Then for any $\alpha\leq\beta$ the closed fiber of $\Spec R_{\beta}\rightarrow \Spec R_{\alpha}$ is of dimension $0$.
Since $\Spec R_{\beta}\rightarrow\Spec R_{\alpha}$ is formally smooth, it is flat \cite[$\mathbf{IV}$, (19.7.1)]{EGA}, and hence we deduce that $R_{\beta}/\m_{\alpha}R_{\beta}$ are fields, that is, $\m_{\alpha}R_{\beta}=\m_{\beta}$.
Now the assertion follows from \ref{prop-nagataquintessence}.
\end{proof}

\subsubsection{Inductive limit of sheaves}\label{subsub-directlimitssheaves}
Let $X$ be a topological space, and consider a filtered inductive system $\{\mathscr{F}_i,\phi_{ij}\}_{i\in I}$ of sheaves (of sets, abelian groups, rings, etc.)\ on $X$ indexed by a directed set $I$.
The inductive limit sheaf $\varinjlim_{i\in I}\mathscr{F}_i$ in the category of sheaves (of sets, abelian groups, rings, etc.)\ is described as follows (cf.\ e.g.\ \cite[II, 1.11]{Gode}): Define a presheaf $\mathscr{F}$ by 
$$
\mathscr{F}(U)=\varinjlim_{i\in I}\mathscr{F}_i(U)
$$
for any open subset $U\subseteq X$; then the desired sheaf is the sheafification of the presheaf $\mathscr{F}$.
By the construction we have
$$
\varinjlim_{i\in I}\mathscr{F}_{i,x}=(\varinjlim_{i\in I}\mathscr{F}_i)_x
$$
for any point $x\in X$.
From this and \ref{prop-directlimits11} (1) we deduce:
\begin{prop}\label{prop-directlimits111}
Let $X$ be a topological space, and $\{\mathscr{F}_i,\phi_{ij}\}_{i\in I}$ a filtered inductive system consisting of sheaf of local rings $($that is, every stalk $\mathscr{F}_{i,x}$ is a local ring$)$ and local homomorphisms $($that is, $\phi_{ij,x}$ is a local homomorphism$)$.
Then the inductive limit $\mathscr{F}=\varinjlim_{i\in I}\mathscr{F}_i$ is a sheaf of local rings.\hfill$\square$
\end{prop}

\begin{prop}\label{prop-projlimcohsch21}
Let $X$ be a topological space, and $\{\mathscr{F}_i,\varphi_{ij}\}_{i\in I}$ a filtered inductive system of sheaves on $X$ indexed by a directed set $I$.
Consider the canonical map
$$
\Phi\colon\varinjlim_{i\in I}\Gamma(X,\mathscr{F}_i)\longrightarrow\Gamma(X,\varinjlim_{i\in I}\mathscr{F}_i).
$$

{\rm (1)} If $X$ is quasi-compact\index{space@space (topological)!quasi-compact topological space@quasi-compact ---}\index{quasi-compact!quasi-compact topological space@--- (topological) space} {\rm (\ref{dfn-quasicompactness} (1))}, then $\Phi$ is injective.

{\rm (2)} If $X$ is coherent\index{space@space (topological)!coherent topological space@coherent ---}\index{coherent!coherent topological space@--- (topological) space} {\rm (\ref{dfn-quasicompact1})}, then $\Phi$ is bijective.
\end{prop}

This proposition can be seen as a special case of a topos-theoretic result \cite[Expos\'e  VI, Th\'eor\`eme 1.23]{SGA4-2}. 
We give here a proof for the reader's convenience.
The ringed space version will be given in \ref{prop-projlimcohsch22}.

\begin{proof}
Set $\mathscr{F}=\varinjlim_{i\in I}\mathscr{F}_i$.

(1) Let $\{s_i\}_{i\in I}$ and $\{t_i\}_{i\in I}$ be inductive systems of sections $s_i,t_i\in\mathscr{F}_i(X)$ whose images by $\Phi$ in $\mathscr{F}(X)$ coincide.
(Here, if necessary, we replace $I$ by a cofinal subset.)
We need to show that $s_k$ and $t_k$ coincide for $k$ sufficiently large.
Since $\mathscr{F}_x=\varinjlim_{i\in I}\mathscr{F}_{i,x}$ for any $x\in X$, there exists $j$ (depending on $x$) such that $s_{k,x}=t_{k,x}$ for any $k\geq j$.
There exists an open neighborhood $U$ of $x$ such that $s_j|_U=t_j|_U$.
Since $X$ is quasi-compact, there exists finite open covering $X=\bigcup^n_{\alpha=1}U_{\alpha}$ and indices $j_{\alpha}$ such that $s_k$ and $t_k$ coincide on $U_{\alpha}$ for $k\geq j_{\alpha}$.
Taking $j$ to be the maximum of $\{j_1,\ldots,j_n\}$, we deduce that $s_k$ and $t_k$ coincide over $X$ for $k\geq j$, as desired.

(2) We only need to show that $\Phi$ is surjective.
Take $s\in\mathscr{F}(X)$.
The germ $s_x$ at $x\in X$ can be written as an inductive system $\{s_{x,i}\}_{i\in I}$ in $\{\mathscr{F}_{i,x}\}_{i\in I}$.
(Here, again, we replace $I$ by a cofinal subset if necessary.)
Take an index $i\in I$ and a section $t_i\in\mathscr{F}_i(U)$ over a quasi-compact open neighborhood $U$ of $x$ such that $s_{x,i}=t_{i,x}$.
Then for any $j\geq i$ one sets $t_j=\varphi_{ij}(U)(t_i)$ so that one gets an inductive system $\{t_j\}_{j\geq i}$ in $\{\mathscr{F}_j(U)\}_{j\geq i}$.
Replacing $U$ by a smaller quasi-compact open neighborhood if necessarily, we may assume that the system $\{t_j\}_{j\geq i}$ is mapped by $\Phi(U)$ to $s|_U$.

Thus, replacing $I$ by a cofinal subset, we get a finite open covering $X=\bigcup^n_{\alpha=1}U_{\alpha}$ consisting of quasi-compact open subsets and for each $\alpha$ an inductive system $\{t_{\alpha,i}\}_{i\in I}$ sitting in $\{\mathscr{F}_i(U_{\alpha})\}_{i\in I}$ that is mapped by $\Phi(U_{\alpha})$ to $s|_{U_{\alpha}}$.
Since $X$ is coherent, each $U_{\alpha\beta}=U_{\alpha}\cap U_{\beta}$ is quasi-compact, and hence by (1) there exists $j_{\alpha\beta}\in I$ such that $t_{\alpha,k}$ and $t_{\beta,k}$ coincide over $U_{\alpha\beta}$ for $k\geq j_{\alpha\beta}$.
Taking $j$ to be the maximum of all $j_{\alpha\beta}$, the local sections $t_{\alpha,k}$ glue together to a section $t_k$ on $X$ for each $k\geq j$.
Then the inductive system $\{t_k\}_{k\geq j}$ thus obtained is mapped by $\Phi$ to $s$.
\end{proof}

\begin{cor}\label{cor-projlimcohsch21}
Let $f\colon X\rightarrow Y$ be a continuous mapping between topological spaces.
Suppose that the following conditions are satisfied$:$
\begin{itemize}
\item[{\rm (a)}] $X$ is quasi-separated\index{space@space (topological)!quasi-separated topological space@quasi-separated ---}\index{quasi-separated!quasi-separated topological space@--- (topological) space} {\rm (\ref{dfn-quasiseparatedness})} and has an open basis consisting of quasi-compact open subsets$;$
\item[{\rm (b)}] $Y$ has an open basis consisting of quasi-compact open subsets$;$
\item[{\rm (c)}] $f$ is quasi-compact\index{map@map (continuous)!quasi-compact map@quasi-compact ---}\index{quasi-compact!quasi-compact map@--- map} {\rm (\ref{dfn-quasicompactness} (2))}.
\end{itemize}
Then for any filtered inductive system $\{\mathscr{F}_i,\varphi_{ij}\}_{i\in I}$ of sheaves on $X$ indexed by a directed set $I$, the canonical morphism 
$$
\Phi\colon\varinjlim_{i\in I}f_{\ast}\mathscr{F}_i\longrightarrow f_{\ast}(\varinjlim_{i\in I}\mathscr{F}_i)
$$
is an isomorphism.
\end{cor}

In other words, the direct image functor $f_{\ast}$ commutes with arbitrary small filtered inductive limits.

\begin{proof}
The sheaf $\varinjlim_{i\in I}f_{\ast}\mathscr{F}_i$ is the sheafification of $V\mapsto\varinjlim_{i\in I}\Gamma(f^{-1}(V),\mathscr{F}_i)$; here, by the assumption, it is enough to consider only quasi-compact $V$'s.
Since $f$ is quasi-compact, and since $X$ is quasi-separated, $f^{-1}(V)$ is coherent (\ref{prop-propsoberness}).
Hence we may apply \ref{prop-projlimcohsch21} (2) to deduce that the map $\varinjlim_{i\in I}\Gamma(f^{-1}(V),\mathscr{F}_i)\rightarrow\Gamma(f^{-1}(V),\varinjlim_{i\in I}\mathscr{F}_i)$ is bijective.
This means, in particular, that the map $\Phi$ is stalkwise bijective, since the formation of taking stalks commutes with the inductive limit $\varinjlim_{i\in I}$, and hence that $\Phi$ is an isomorphism.
\end{proof}

\subsubsection{Sheaves on limit spaces}\label{subsub-sheavestopologicalspacelimits}
In this paragraph, we consider 
\begin{itemize}
\item a filtered projective system of topological spaces $\{X_i,p_{ij}\colon X_j\rightarrow X_i\}_{i\in I}$ indexed by a directed set $I$
\end{itemize}
such that
\begin{itemize}
\item[{\rm (a)}] for any $i\in I$ the topological space $X_i$ is coherent\index{space@space (topological)!coherent topological space@coherent ---}\index{coherent!coherent topological space@--- (topological) space} {\rm (\ref{dfn-quasicompact1})} and sober\index{space@space (topological)!sober topological space@sober ---} {\rm (\S\ref{subsub-sober})}$;$
\item[{\rm (b)}] for any $i\leq j$ the transition map $p_{ij}\colon X_j\rightarrow X_i$ is quasi-compact {\rm (\ref{dfn-quasicompactness} (2))}.
\end{itemize}
Note that by \ref{thm-projlimcohsch1} (1) the limit space $X=\varprojlim_{i\in I}X_i$ is coherent and sober, and that the canonical projection maps $p_i\colon X\rightarrow X_i$ for $i\in I$ are quasi-compact.
Suppose, moreover, that we are given the following data:
\begin{itemize}
\item for each $i\in I$ a sheaf $\mathscr{F}_i$ (of sets, abelian groups, etc.)\ on $X_i$;
\item for each pair $(i,j)$ of indices in $I$ with $i\leq j$, a morphism $\varphi_{ij}\colon p^{-1}_{ij}\mathscr{F}_i\rightarrow\mathscr{F}_j$ of sheaves such that $\varphi_{ik}=\varphi_{jk}\circ p^{-1}_{jk}\varphi_{ij}$ holds whenever $i\leq j\leq k$.
\end{itemize}
Then one has the inductive system $\{p^{-1}_i\mathscr{F}_i\}_{i\in I}$ of sheaves on $X$ indexed by $I$, and thus the sheaf 
$$
\mathscr{F}=\varinjlim_{i\in I}p^{-1}_i\mathscr{F}_i
$$
on $X$.

\begin{prop}\label{prop-coherentprojectivesystemindlimits1}
The canonical map
$$
\varinjlim_{i\in I}\Gamma(X_i,\mathscr{F}_i)\longrightarrow\Gamma(X,\mathscr{F})
$$
is an isomorphism.
\end{prop}

\begin{proof}
We already know by \ref{prop-projlimcohsch21} that $\Gamma(X,\mathscr{F})\cong\varinjlim_{i\in I}\Gamma(X,p^{-1}_i\mathscr{F}_i)$ holds.
We want to show that the map $\varinjlim_{i\in I}\Gamma(X_i,\mathscr{F}_i)\rightarrow\varinjlim_{i\in I}\Gamma(X,p^{-1}_i\mathscr{F}_i)$ is bijective.

\medskip
{\sc Step 1.} We first claim that the canonical map
$$
\varinjlim_{i\in I}\Gamma(X_i,\mathscr{F}_i)\longrightarrow\varinjlim_{i\in I}\varinjlim_U\Gamma(U,\mathscr{F}_i),
$$
where $U$ in the right-hand side runs through open subsets of $X_i$ such that $p_i(X)\subseteq U$, is an isomorphism.
By Exercise \ref{exer-doublecolimits} the double inductive limit in the right-hand side is canonically isomorphic to the inductive limit taken over the directed set 
$$
\Lambda=\{(i,U)\,|\,i\in I,\ p_i(X)\subseteq U\subseteq X_i\},
$$
where $(i,U)\leq (j,V)$ if and only if $i\leq j$ and $p^{-1}_{ij}(U)\supseteq V$.
Hence the desired result follows if one shows that the subset of $\Lambda$ consisting of elements of the form $(i,X_i)$ is cofinal.

To see this, take any $(i,U)\in\Lambda$.
Since $p_i(X)$ is quasi-compact and $X_i$ is coherent, we may assume that $U$ is quasi-compact.
The condition $p_i(X)\subseteq U$ is equivalent to $p^{-1}_i(U)=X$.
Then by \ref{cor-projlimcohsch11} there exists $j\in I$ with $i\leq j$ such that $p^{-1}_{ij}(U)=X_j$.
Hence we have $(i,U)\leq (j,X_j)$, as desired.

\medskip
{\sc Step 2.} Take $\{s_i\}_{i\in I}\in\varinjlim_{i\in I}\Gamma(X,p^{-1}_i\mathscr{F}_i)$.
By \ref{prop-projlimcohtopspqcptopen} there exist a finite open covering $X=\bigcup_{\alpha\in L}U_{\alpha}$ by quasi-compact open subsets and an index $i\in I$ such that
\begin{itemize}
\item each $U_{\alpha}$ is of the form $U_{\alpha}=p^{-1}_{\alpha}(U_{\alpha i})$ for a quasi-compact open subset $U_{\alpha i}$ of $X_i$;
\item $s_i|_{U_{\alpha}}$ for each $\alpha$ lies in $\varinjlim_V\Gamma(V,\mathscr{F}_i)$, where $V$ runs through all open subsets of $X_i$ containing $p_i(U_{\alpha})$.
\end{itemize}
Then $\{s_j|_{U_{\alpha}}\}_{j\geq i}$ defines a section in $\varinjlim_{j\geq i}\varinjlim_V\Gamma(V,\mathscr{F}_j)$, where $U$ runs through all open subsets of $U_{\alpha j}=p^{-1}_{ij}(U_{\alpha i})$ containing $p_j(U_{\alpha})$.
This defines by {\sc Step 1} a unique section in $\varinjlim_{j\geq i}\Gamma(U_{\alpha j},\mathscr{F}_j)$.
Since the spaces of the form $U_{\alpha j}\cap U_{\beta j}$ ($\alpha,\beta\in L$, $j\geq i$) are all coherent, these sections glue to a unique section in $\varinjlim_{i\in I}\Gamma(X_i,\mathscr{F}_i)$, as desired.
\end{proof}

\begin{cor}\label{cor-coherentprojectivesystemindlimits01}
Let $\{X_i,p_{ij}\}_{i\in I}$ and $p_i\colon X\rightarrow X_i$ for $i\in I$ be as above, and $Z$ a topological space.
Suppose we are given a system of continuous mappings $\{g_i\colon X_i\rightarrow Z\}_{i\in I}$ such that we have $g_j=g_i\circ p_{ij}$ whenever $i\leq j$.
Then for any sheaf $\mathscr{G}$ on $Z$ the canonical map
$$
\varinjlim_{i\in I}\Gamma(X_i,g^{-1}_i\mathscr{G})\longrightarrow\Gamma(X,g^{-1}\mathscr{G}),
$$
where $g=\varprojlim_{i\in I}g_i$, is an isomorphism.
\end{cor}

\begin{proof}
Apply \ref{prop-coherentprojectivesystemindlimits1} to the situation where $\mathscr{F}_i=g^{-1}_i\mathscr{G}$ for $i\in I$.
\end{proof}

\begin{cor}\label{cor-coherentprojectivesystemindlimits1}
Let $\{X_i,p_{ij}\}_{i\in I}$ and $p_i\colon X\rightarrow X_i$ for $i\in I$ be as above, and $\mathscr{F}$ a sheaf on $X$.
Then the canonical morphism 
$$
\varinjlim_{i\in I}p^{-1}_ip_{i\ast}\mathscr{F}\longrightarrow\mathscr{F}
$$
is an isomorphism.
\end{cor}

\begin{proof}
By a similar reasoning as in the proof of \ref{cor-projlimcohsch21}, it suffices to show that the morphism in question induces an isomorphism between the set of sections over each quasi-compact open subset $U\subseteq X$.
By \ref{prop-projlimcohtopspqcptopen} there exist an index $i_0\in I$ and quasi-compact open subset $U_0$ of $X_{i_0}$ such that $p^{-1}_{i_0}(U_0)=U$.
Since each $p^{-1}_{ii_0}(U_0)$ for $i\geq i_0$ is a coherent sober space, and since $U=\varprojlim_{i\geq i_0}p^{-1}_{ii_0}(U_0)$, we may assume that $U=X$ without loss of generality. 
(Here we replace the index set $I$ by the cofinal subset $\{i\in I\,|\, i\geq i_0\}$.)
Then we have
$$
\Gamma(X,\varinjlim_{i\in I}p^{-1}_ip_{i\ast}\mathscr{F})\cong\varinjlim_{i\in I}\Gamma(X_i,p_{i\ast}\mathscr{F})=\varinjlim_{i\in I}\Gamma(X,\mathscr{F})=\Gamma(X,\mathscr{F}),
$$
where the first isomorphism is due to \ref{prop-coherentprojectivesystemindlimits1}.
\end{proof}

Next, in addition to the data fixed in the beginning of this paragraph, we consider the following data:
\begin{itemize}
\item another filtered projective system of topological spaces $\{Y_i,q_{ij}\}_{i\in I}$, indexed by the same directed set $I$, that satisfies the conditions similar to (a) and (b) in the beginning of this paragraph;
\item a map $\{f_i\}_{i\in I}$ of projective systems from $\{X_i,p_{ij}\}_{i\in I}$ to $\{Y_i,q_{ij}\}_{i\in I}$, that is, a collection of continuous maps $f_i\colon X_i\rightarrow Y_i$ such that $q_{ij}\circ f_j=f_i\circ p_{ij}$ for any $i\leq j$.
\end{itemize}
We set $Y=\varprojlim_{i\in I}Y_i$ and denote the canonical projection by $q_i\colon Y\rightarrow Y_i$ for each $i\in I$.
Moreover, we have the continuous map
$$
f=\varprojlim_{i\in I}f_i\colon X\longrightarrow Y.
$$
We assume that 
\begin{itemize}
\item[{\rm (c)}] for any $i\in I$ the map $f_i$ is quasi-compact.
\end{itemize}
Note that the map $f$ is quasi-compact due to \ref{thm-projlimcohspacepres} (1).

\begin{cor}\label{cor-coherentprojectivesystemindlimits2}
The canonical morphism of sheaves
$$
\varinjlim_{i\in I}q^{-1}_if_{i\ast}\mathscr{F}_i\longrightarrow f_{\ast}\mathscr{F}
$$
is an isomorphism.
\end{cor}

\begin{proof}
By a similar reasoning as in the proof of \ref{cor-projlimcohsch21}, we may restrict ourselves to showing that the map
$$
\varinjlim_{i\in I}\Gamma(V,q^{-1}_if_{i\ast}\mathscr{F}_i)\longrightarrow\Gamma(f^{-1}(V),\mathscr{F})
$$
is an isomorphism for any quasi-compact open subset $V$ of $Y$.
Sinilarly to the proof of \ref{cor-coherentprojectivesystemindlimits1}, we may assume that $V=Y$ without loss of generality.
In this situation, one can replace the left-hand side by the double inductive limit
$$
\varinjlim_{i\in I}\varinjlim_U\Gamma(f^{-1}_i(U),\mathscr{F}_i),
$$
where $U$ runs through open subsets of $Y_i$ containing $q_i(Y)$.
By an argument similar to that in the proof of \ref{prop-coherentprojectivesystemindlimits1}, one sees that this limit is isomorphic to $\varinjlim_{i\in I}\Gamma(X_i,\mathscr{F}_i)$.
Then the desired result follows from \ref{prop-coherentprojectivesystemindlimits1}.
\end{proof}

\begin{cor}\label{cor-coherentprojectivesystemindlimits21}
Let $\{X_i,p_{ij}\}_{i\in I}$ and $p_i\colon X\rightarrow X_i$ for $i\in I$ be as above, and $Z$ a coherent sober space.
Suppose we are given a system of quasi-compact maps $\{g_i\colon X_i\rightarrow Z\}_{i\in I}$ such that $g_j=g_i\circ p_{ij}$ whenever $i\leq j$.
Then for any sheaf $\mathscr{G}$ on $Z$ the canonical morphism
$$
\varinjlim_{i\in I}g_{i\ast}g^{-1}_i\mathscr{G}\longrightarrow g_{\ast}g^{-1}\mathscr{G},
$$
where $g=\varprojlim_{i\in I}g_i$, is an isomorphism.
\end{cor}

\begin{proof}
Apply \ref{cor-coherentprojectivesystemindlimits2} to the situation where $\{Y_i,q_{ij}\}_{i\in I}$ is the constant system ($Y_i=Z$) and $\mathscr{F}_i=g^{-1}_i\mathscr{G}$.
\end{proof}

\begin{cor}\label{cor-coherentprojectivesystemindlimits3}
Let $\{X_i,p_{ij}\}_{i\in I}$, $\{Y_i,q_{ij}\}_{i\in I}$, and $f_i$'s be as above, and $\mathscr{F}$ a sheaf on $X$.
Then the canonical morphism
$$
\varinjlim_{i\in I}q^{-1}_if_{i\ast}p_{i\ast}\mathscr{F}\longrightarrow f_{\ast}\mathscr{F}
$$
is an isomorphism.
\end{cor}

\begin{proof}
By \ref{cor-coherentprojectivesystemindlimits1} we have $\varinjlim_{i\in I}q^{-1}_if_{i\ast}p_{i\ast}\mathscr{F}=\varinjlim_{i\in I}q^{-1}_iq_{i\ast}f_{\ast}\mathscr{F}\cong f_{\ast}\mathscr{F}$.
\end{proof}

\subsubsection{Canonical flasque resolution}\label{subsub-canflasqueres}
Let us recall the {\em canonical flasque resolution}\index{resolution!canonical flasque resolution@canonical flasque ---} for abelian sheaves; we only recall its basic properties, and refer to \cite[II.4.3]{Gode} for the construction.
Let $X$ be a topological space.
We denote by $\ASh_X$ the abelian category of sheaves of abelian groups on $X$ and by $\CC^+(\ASh_X)$ the abelian category of complexes bounded below consisting of objects and arrows in $\ASh_X$ (cf.\ \S\ref{subsub-complexcategorydef}). 
The canonical flasque resolution is a functor 
$$
\mathscr{C}^{\bullet}(X,\,\cdot\,)\colon\ASh_X\longrightarrow\CC^+(\ASh_X)
$$
such that for any abelian sheaf $\mathscr{F}$ on $X$:
\begin{itemize}
\item[{\rm (a)}] $\mathscr{C}^q(X,\mathscr{F})=0$ for $q<0$; 
\item[{\rm (b)}] $\mathscr{C}^q(X,\mathscr{F})$ for $q\in\Z$ is flasque; 
\item[{\rm (c)}] $\mathscr{C}^{\bullet}(X,\mathscr{F})$ is equipped with the augmentation $\mathscr{F}\rightarrow\mathscr{C}^{\bullet}(X,\mathscr{F})$ such that the sequence 
$$
0\longrightarrow\mathscr{F}\longrightarrow\mathscr{C}^0(X,\mathscr{F})\longrightarrow\mathscr{C}^1(X,\mathscr{F})\longrightarrow\mathscr{C}^2(X,\mathscr{F})\longrightarrow\cdots
$$
is exact.
\end{itemize}
Moreover, the functor $\mathscr{C}^{\bullet}(X,\,\cdot\,)$ is exact.

The importance of the canonical flasque resolution lies in its canonicity.
For instance, given an ordered set\index{ordered!set@--- set}\index{set!ordered set@ordered ---} $I$, one can consider the category $\ASh^I_X$ of inductive systems of abelian sheaves on $X$; the canonicity of the canonical flasque resolution allows one to construct the exact functor 
$$
\mathscr{C}^{\bullet}(X,\,\cdot\,)\colon\ASh^I_X\longrightarrow\CC^+(\ASh_X)^I,
$$
which maps each inductive system $\{\mathscr{F}_i,\varphi_{ij}\}_{i\in I}$ of abelian sheaves to the inductive system of complexes $\{\mathscr{C}^{\bullet}(X,\mathscr{F}_i),\varphi^{\bullet}_{ij}\}$ consisting of canonical flasque resolutions of $\mathscr{F}_i$'s.

\subsubsection{Inductive limit and cohomology}\label{subsub-directlimitandcohomology}
\begin{prop}\label{prop-directlimcohcoh11}
Let $X$ be a coherent topological space\index{space@space (topological)!coherent topological space@coherent ---}\index{coherent!coherent topological space@--- (topological) space} {\rm (\ref{dfn-quasicompact1})}, and $\{\mathscr{F}_i,\varphi_{ij}\}_{i\in I}$ a filtered inductive system of sheaves of abelian groups on $X$ indexed by a directed set $I$. 
Then the canonical map 
$$
\Phi^q\colon\varinjlim_{i\in I}\H^q(X,\mathscr{F}_i)\longrightarrow\H^q(X,\varinjlim_{i\in I}\mathscr{F}_i)
$$
is bijective for any $q\geq 0$.
\end{prop}

\begin{proof}
First note that the case $q=0$ has already been proved in \ref{prop-projlimcohsch21} (2).
In order to show the general case, let us take the filtered inductive system of complexes $\{\mathscr{C}^{\bullet}(X,\mathscr{F}_i),\varphi^{\bullet}_{ij}\}$ as above.
By \ref{prop-directlimits2} 
$$
0\longrightarrow\varinjlim_{i\in I}\mathscr{F}_i\longrightarrow\varinjlim_{i\in I}\mathscr{C}^{\bullet}(X,\mathscr{F}_i)
$$
is exact.
This gives a quasi-flasque resolution of $\varinjlim_{i\in I}\mathscr{F}_i$ (cf.\ Exercise \ref{exer-injlimcohcoh}), and hence the desired result follows from \ref{prop-projlimcohsch21} (2) and exactness of inductive limits.
\end{proof}

\begin{cor}\label{cor-directlimcohcoh11}
Let $f\colon X\rightarrow Y$ be a continuous mapping between topological spaces.
Suppose that the following conditions are satisfied$:$
\begin{itemize}
\item[{\rm (a)}] $X$ is quasi-separated\index{space@space (topological)!quasi-separated topological space@quasi-separated ---}\index{quasi-separated!quasi-separated topological space@--- (topological) space} {\rm (\ref{dfn-quasiseparatedness})} and has an open basis consisting of quasi-compact open subsets$;$
\item[{\rm (b)}] $Y$ has an open basis consisting of quasi-compact open subsets$;$
\item[{\rm (c)}] $f$ is quasi-compact\index{map@map (continuous)!quasi-compact map@quasi-compact ---}\index{quasi-compact!quasi-compact map@--- map} {\rm (\ref{dfn-quasicompactness} (2))}.
\end{itemize}
Then for any filtered inductive system $\{\mathscr{F}_i,\varphi_{ij}\}_{i\in I}$ of sheaves of abelilan groups on $X$ indexed by a directed set, the canonical morphism 
$$
\Phi^q\colon\varinjlim_{i\in I}\RD^qf_{\ast}\mathscr{F}_i\longrightarrow \RD^qf_{\ast}(\varinjlim_{i\in I}\mathscr{F}_i)
$$
is an isomorphism. \hfill$\square$
\end{cor}

The proof is similar to that of \ref{cor-projlimcohsch21}, where we use \ref{prop-directlimcohcoh11} instead of \ref{prop-projlimcohsch21}.

\subsubsection{Cohomology of sheaves on limit spaces}\label{subsub-cohsheavestopologicalspacelimits}
Let us now return to the situation as in \S\ref{subsub-sheavestopologicalspacelimits}.
\begin{lem}\label{lem-cohinjlim}
Suppose that each $\mathscr{F}_i$ is a flasque sheaf on $X_i$.
Then the sheaf $\mathscr{F}=\varinjlim_{i\in I}p^{-1}_i\mathscr{F}_i$ is quasi-flasque.\hfill$\square$
\end{lem}

This follows easily from Exercise \ref{exer-injlimcohcoh} and that each $p^{-1}_i\mathscr{F}_i$ is flasque.
In the following statements, all sheaves are supposed to be sheaves of abelian groups, and morphisms of sheaves are morphisms of abelian sheaves.
\begin{prop}\label{prop-cohcoherentprojectivesystemindlimits1}
The canonical map
$$
\varinjlim_{i\in I}\H^q(X_i,\mathscr{F}_i)\longrightarrow\H^q(X,\mathscr{F})
$$
is an isomorphism for $q\geq 0$. \hfill$\square$
\end{prop}

This follows from \ref{prop-coherentprojectivesystemindlimits1} and \ref{lem-cohinjlim} as in the proof of \ref{prop-directlimcohcoh11}.
By this and \ref{cor-coherentprojectivesystemindlimits1} we have:

\begin{cor}\label{cor-cohcoherentprojectivesystemindlimits11}
Let $\{X_i,p_{ij}\}_{i\in I}$ and $p_i\colon X\rightarrow X_i$ for $i\in I$ be as in the beginning of {\rm \S\ref{subsub-sheavestopologicalspacelimits}}, and $\mathscr{F}$ a sheaf on $X$.
Then the canonical morphism 
$$
\varinjlim_{i\in I}\H^q(X_i,p_{i\ast}\mathscr{F})\longrightarrow\H^q(X,\mathscr{F})
$$
is an isomorphism for $q\geq 0$. \hfill$\square$
\end{cor}

\begin{cor}\label{cor-cohcoherentprojectivesystemindlimits01}
Let $\{X_i,p_{ij}\}_{i\in I}$ and $p_i\colon X\rightarrow X_i$ for $i\in I$ be as in the beginning of {\rm \S\ref{subsub-sheavestopologicalspacelimits}}, and $Z$ a topological space.
Suppose we are given a system of continuous maps $\{g_i\colon X_i\rightarrow Z\}_{i\in I}$ such that, whenever $i\leq j$, we have $g_j=g_i\circ p_{ij}$.
Then for any sheaf $\mathscr{G}$ on $Z$ the canonical map
$$
\varinjlim_{i\in I}\H^q(X_i,g^{-1}_i\mathscr{G})\longrightarrow\H^q(X,g^{-1}\mathscr{G}),
$$
where $g=\varprojlim_{i\in I}g_i$, is an isomorphism for $q\geq 0$.\hfill$\square$
\end{cor}

This follows from \ref{prop-cohcoherentprojectivesystemindlimits1} applied to the situation where $\mathscr{F}_i=g^{-1}_i\mathscr{G}$ for $i\in I$.
Similarly to \ref{cor-coherentprojectivesystemindlimits2} and its subsequent results, one can show the following statements:
\begin{cor}\label{cor-cohcoherentprojectivesystemindlimits2}
Let $\{X_i,p_{ij}\}_{i\in I}$, $\{Y_i,q_{ij}\}_{i\in I}$, and $f_i$'s be as in {\rm \S\ref{subsub-sheavestopologicalspacelimits}}.
Then the canonical morphism of sheaves
$$
\varinjlim_{i\in I}q^{-1}_i\RD^qf_{i\ast}\mathscr{F}_i\longrightarrow\RD^qf_{\ast}\mathscr{F}
$$
is an isomorphism for $q\geq 0$. \hfill$\square$
\end{cor}

\begin{cor}\label{cor-cohcoherentprojectivesystemindlimits21}
Let $\{X_i,p_{ij}\}_{i\in I}$ be as in {\rm \S\ref{subsub-sheavestopologicalspacelimits}}, and $Z$ a coherent sober space. 
Suppose we are given a system of quasi-compact maps $\{g_i\colon X_i\rightarrow Z\}_{i\in I}$ such that $g_j=g_i\circ p_{ij}$ whenever $i\leq j$.
Then for any sheaf $\mathscr{G}$ on $Z$ the canonical morphism
$$
\varinjlim_{i\in I}\RD^qg_{i\ast}(g^{-1}_i\mathscr{G})\longrightarrow\RD^qg_{\ast}(g^{-1}\mathscr{G}),
$$
where $g=\varprojlim_{i\in I}g_i$, is an isomorphism for $q\geq 0$. \hfill$\square$
\end{cor}

\begin{cor}\label{cor-cohcoherentprojectivesystemindlimits3}
Let $\{X_i,p_{ij}\}_{i\in I}$, $\{Y_i,q_{ij}\}_{i\in I}$, and $f_i$'s be as in {\rm \S\ref{subsub-sheavestopologicalspacelimits}}, and $\mathscr{F}$ a sheaf on $X$.
Then the canonical morphism
$$
\varinjlim_{i\in I}q^{-1}_i\RD^qf_{i\ast}(p_{i\ast}\mathscr{F})\longrightarrow\RD^qf_{\ast}\mathscr{F}
$$
is an isomorphism for $q\geq 0$. \hfill$\square$
\end{cor}
\index{limit!inductive limit@inductive ---|)}

\subsection{Projective limits}\label{sub-projectivelimits}\index{limit!projective limit@projective ---|(}
\subsubsection{The Mittag-Leffler condition}\label{subsub-ML}
Let $\mathscr{C}$ be either the category of sets or an abelian category that has all small products (that is, for any family $\{A_i\}_{i\in I}$ of objects indexed by a small set, the product $\prod_{i\in I}A_i$ is representable).
Then the category $\mathscr{C}$ is {\em small complete}\index{small!small complete@--- complete}, that is, the limit $\varprojlim F$ for a functor $F\colon\mathscr{D}\rightarrow\mathscr{C}$ exists whenever $\mathscr{D}$ is essentially small\index{small!essentially small@essentially ---} (\cite[Chap.\ V, \S2, Cor.\ 2]{Mac}).

Let $\textrm{{\boldmath $A$}}=\{A_i,f_{ij}\colon A_j\rightarrow A_i\}$ be a projective system\index{system!projective system@projective ---} of objects in $\mathscr{C}$ indexed by a directed set $I$, and set $A=\varprojlim_{i\in I}A_i$.
\begin{dfn}\label{dfn-ML1}{\rm 
(1) The projective system {\boldmath $A$} is said to be {\em strict}\index{strict projective system@strict (projective system)}\index{system!projective system@projective ---!strict projective system@strict --- ---} if all transition maps $f_{ij}$ for $i\leq j$ are epimorphic.

(2) Suppose $I=\N$, the ordered set of all natural numbers. 
Then the projective system {\boldmath $A$} is said to be {\em essentially constant}\index{essentially constant}, if there exists $N$ such that for $N\leq i\leq j$ the transition maps $f_{ij}$ are isomorphisms.}
\end{dfn}

For each $i\in I$ we set
$$
A'_i=\inf_{i\leq j}f_{ij}(A_j),
$$
and call it the {\em universal image}\index{universal image} in $A_i$; notice that here the infimum is nothing but the projective limit of $\{f_{ij}(A_j)\}_{j\geq i}$ and hence is a subobject of $A_i$.
Clearly, we have $f_{ij}(A'_j)\subseteq A'_i$ for any $i\leq j$, and $f_i(A)\subseteq A'_i$ for any $i\in I$.
Thus we have the projective system $\textrm{{\boldmath $A'$}}=\{A'_i,f_{ij}|_{A'_j}\}$ such that $A=\varprojlim_{i\in I}A'_i$.
Notice that, if {\boldmath $A$} is strict, then we have $\textrm{{\boldmath $A'$}}=\textrm{{\boldmath $A$}}$.

The {\em Mittag-Leffler condition}\index{Mittag-Leffler condition} for a projective system $\textrm{{\boldmath $A$}}=\{A_i,f_{ij}\}$ is:
\begin{itemize}
\item[{\bf (ML)}] for any $i\in I$ there exists $j\geq i$ such that $f_{ik}(A_k)=f_{ij}(A_j)$ for any $k\geq j$.
\end{itemize}
Clearly, any strict projective system satisfies {\bf (ML)}. 
Conversely, if {\boldmath $A$} satisfies {\bf (ML)}, then the induced projective system {\boldmath $A'$} of universal images is strict.
The condition {\bf (ML)} is closely related to the non-emptiness of projective limits; by Exercise \ref{exer-finalcountable} and \cite[Chap.\ III, \S7.4, Prop.\ 5]{Bourb2} we have:
\begin{prop}\label{prop-ML2}
Let $\textrm{{\boldmath $X$}}=\{X_i,f_{ij}\}$ be a projective system of sets indexed by a directed set $I$, and $X=\varprojlim_{i\in I}X_i$.
Suppose that $I$ has a cofinal and at most countable subset and that the system {\boldmath $X$} satisfies $\mbox{{\bf (ML)}}$. 
Then for any $i\in I$ the canonical projection $f_i\colon X\rightarrow X_i$ maps $X$ surjectively onto the universal image\index{universal image} $X'_i$.
$($Hence, in particular, if $\textrm{{\boldmath $X$}}=\{X_i,f_{ij}\}$ is strict, then all maps $f_i$ for $i\in I$ are surjective.$)$ \hfill$\square$
\end{prop}

\begin{prop}[{\cite[$\mathbf{0}_{\mathbf{III}}$, (13.2.1)]{EGA}}]\label{prop-ML3}
Let 
$$
0\longrightarrow\textrm{{\boldmath $M'$}}\longrightarrow\textrm{{\boldmath $M$}}\longrightarrow\textrm{{\boldmath $M''$}}\longrightarrow 0,\leqno{(\ast)}
$$
where $\textrm{{\boldmath $M$}}=\{M_i\}$ etc., be an exact sequence in $\Ab^{I^{\opp}}$ consisting of projective systems of abelian groups indexed by a common directed set $I$.
\begin{itemize}
\item[$(1)$] If {\boldmath $M$} satisfies $\mbox{{\bf (ML)}}$, then so does {\boldmath $M''$}.
\item[$(2)$] If {\boldmath $M'$} and {\boldmath $M''$} satisfy $\mbox{{\bf (ML)}}$, then so does {\boldmath $M$}. \hfill$\square$
\end{itemize}
\end{prop}

Since the functor $\varprojlim_{i\in I}$ has a left-adjoint, we have:
\begin{prop}\label{prop-projlimleftexact}
The functor $\varprojlim_{i\in I}\colon\Ab^{I^{\opp}}\rightarrow\Ab$ is left-exact.\hfill$\square$
\end{prop}

Hence one can consider the right derived functors $\varprojlim^{(q)}_{i\in I}$ for $q\geq 0$.
\begin{lem}\label{lem-ML3}
Let $I$ be a directed set that contains a cofinal and at most countable subset.
Consider an exact sequence
$$
\cdots\longrightarrow\textrm{{\boldmath $M$}}^{q-1}\stackrel{d^{q-1}}{\longrightarrow}\textrm{{\boldmath $M$}}^q\stackrel{d^q}{\longrightarrow}\textrm{{\boldmath $M$}}^{q+1}\longrightarrow\cdots
$$
of projective system of abelian groups indexed by $I$, where $\textrm{{\boldmath $M$}}^q=\{M^q_i,f^q_{ij}\}$.
Suppose that for any $q$ the projective system $\textrm{{\boldmath $M$}}^q$ satisfies $\mbox{{\bf (ML)}}$.
Then the induced sequence
$$
\cdots\longrightarrow\varprojlim_{i\in I}M^{q-1}_i\stackrel{d^{q-1}_{\infty}}{\longrightarrow}\varprojlim_{i\in I}M^q_i\stackrel{d^q_{\infty}}{\longrightarrow}\varprojlim_{i\in I}M^{q+1}_i\longrightarrow\cdots
$$
is exact.
\end{lem}

\begin{proof}
It is easy to see that $d^q_{\infty}\circ d^{q-1}_{\infty}=0$.
Let $\{x_i\}\in\varprojlim_{i\in I}M^q_i$ be such that $d^q_{\infty}(\{x_i\})=0$.
We consider the following diagram with exact rows for $i\leq j$:
$$
\xymatrix{M^{q-2}_j\ar[d]_{f^{q-2}_{ij}}\ar[r]^{d^{q-2}_j}&M^{q-1}_j\ar[d]_{f^{q-1}_{ij}}\ar[r]^{d^{q-1}_j}&M^q_j\ar[d]^{f^q_{ij}}\ar[r]^{d^q_j}&M^{q+1}_j\ar[d]^{f^{q+1}_{ij}}\\ M^{q-2}_i\ar[r]_{d^{q-2}_i}&M^{q-1}_i\ar[r]_{d^{q-1}_i}&M^q_i\ar[r]_{d^q_i}&M^{q+1}_i\rlap{.}}
$$
For each $i$ we have $d^q_i(x_i)=0$. 
Hence the subset $S_i=(d^{q-1}_i)^{-1}(x_i)$ of $M^{q-1}_i$ is non-empty, and $\{S_i,f^{q-1}_{ij}|_{S_j}\}$ is a projective system of sets.

\medskip
{\sc Claim.} {\it The projective system $\{S_i,f^{q-1}_{ij}|_{S_j}\}$ satisfies} {\bf (ML)}.

\medskip
Let $T_i$ be the image of $d^{q-2}_i$. 
The projective system $\{T_i,f^{q-1}_{ij}|_{T_j}\}$ satisfies {\bf (ML)} (\ref{prop-ML3} (1)).
Hence for any $i\in I$ there exists $j\geq i$ such that $f^{q-1}_{ij}(T_j)=f^{q-1}_{ik}(T_k)$ for any $j\leq k$.
We want to show that $f^{q-1}_{ij}(S_j)=f^{q-1}_{ik}(S_k)$ for any $j\leq k$.
Fix $y_k\in S_k$, and set $y_j=f^{q-1}_{jk}(y_k)$ and $y_i=f^{q-1}_{ik}(y_k)$.
Take any $z_i\in f^{q-1}_{ij}(S_j)$, and put $z_i=f^{q-1}_{ij}(z_j)$.
We have $z_j-y_j\in T_j$, and hence there exists $w_k\in T_k$ such that $z_i-y_i=f^{q-1}_{ik}(w_k)$. 
Hence $z_i=f^{q-1}_{ik}(w_k)+y_i=f^{q-1}_{ik}(w_k+y_k)\in f^{q-1}_{ik}(S_k)$, which shows
$f^{q-1}_{ij}(S_j)\subseteq f^{q-1}_{ik}(S_k)$.
As the other inclusion is clear, we have $f^{q-1}_{ij}(S_j)=f^{q-1}_{ik}(S_k)$, and thus the claim follows.

Now, by \ref{prop-ML2} and the fact that the universal images\index{universal image} of the system $\{S_i\}$ are non-empty, there exists $\{y_i\}\in\varprojlim_{i\in I}M^{q-1}_i$ such that $d^{q-1}_{\infty}(\{y_i\})=\{x_i\}$.
\end{proof}

\begin{cor}\label{cor-ML3cornew}
The functor $\varprojlim_{i\in I}$ maps any acyclic complex in $\CC^+(\Ab^{I^{\opp}})$ consisting of objects satisfying $\mbox{{\bf (ML)}}$ to an acyclic complex in $\CC^+(\Ab)$.\hfill$\square$
\end{cor}

\subsubsection{Canonical strict resolution}\label{subsub-canstrictres}
Let $\textrm{{\boldmath $M$}}=\{M_i,f_{ij}\}_{i\in\N}$ be a projective system of abelian groups indexed by $\N$.
We are going to construct a short exact sequence 
$$
0\longrightarrow\textrm{{\boldmath $M$}}\longrightarrow\textrm{{\boldmath $J$}}^0\longrightarrow\textrm{{\boldmath $J$}}^1\longrightarrow 0
$$
in $\Ab^{\N^{\opp}}$ such that $\textrm{{\boldmath $J$}}^0$ and $\textrm{{\boldmath $J$}}^1$ are strict.

The system $\textrm{{\boldmath $J$}}^0=\{J^0_i\}_{i\in\N}$ is constructed as follows.
For $i\in\N$ we set $J^0_i=M_0\oplus\cdots\oplus M_i$ and define the transition map $J^0_j\rightarrow J^0_i$ for $i\leq j$ by
$$
M_0\oplus\cdots\oplus M_j\ni(x_0,\ldots,x_j)\longmapsto (x_0,\ldots,x_{i_1},\sum^j_{k=i}f_{ki}(x_k))\in M_0\oplus\cdots\oplus M_i.
$$
We have the obvious inclusion $\textrm{{\boldmath $M$}}\hookrightarrow\textrm{{\boldmath $J$}}^0$.
The system $\textrm{{\boldmath $J$}}^1=\{J^1_i\}_{i\in\N}$ is defined to be the cokernel of this map.
Explicitly, it is given by $J^1_i=J^0_{i-1}$ for $i>0$ and $J^1_0=0$; the transition maps are the canonical projections.
Clearly, with the systems $\textrm{{\boldmath $J$}}^0$ and $\textrm{{\boldmath $J$}}^1$ thus constructed, we have the desired exact sequence of projective systems as above.

Notice that, by the construction, the formation of the resolution $0\rightarrow\textrm{{\boldmath $M$}}\rightarrow\textrm{{\boldmath $J$}}^{\bullet}$ is functorial and hence defines a functor 
$$
\Ab^{\N^{\opp}}\longrightarrow\CC^{[0,1]}(\Ab^{\N^{\opp}}),
$$
which is, moreover, exact.

We call this resolution $0\rightarrow\textrm{{\boldmath $M$}}\rightarrow\textrm{{\boldmath $J$}}^{\bullet}$ the {\em canonical strict resolution}\index{resolution!canonical strict resolution@canonical strict ---} of {\boldmath $M$}.

In view of \ref{lem-ML3}, we can use the canonical strict resolution to compute the right derived functors $\varprojlim^{(q)}_{i\in I}$ for $q\geq 0$ attached to $\varprojlim_{i\in I}$.

\begin{prop}\label{prop-ML6}
Let $\textrm{{\boldmath $M$}}=\{M_i,f_{ij}\}_{i\in\N}$ be a projective system of abelian groups indexed by a directed set $I$ that has a cofinal and at most countable subset.
\begin{itemize}
\item[$(1)$] We have $\varprojlim^{(q)}_{i\in I}M_i=0$ for $q\geq 2$.
\item[$(2)$] If {\boldmath $M$} satisfies $\mbox{{\bf (ML)}}$, then we have $\varprojlim^{(1)}_{i\in I}M_i=0$.
\end{itemize}
\end{prop}

\begin{proof}
By Exercise \ref{exer-finalcountable} we may assume that $I=\N$.
Since the canonical strict resolution is of length $1$, (1) follows immediately.
If {\boldmath $M$} satisfies {\bf (ML)}, then by \ref{lem-ML3} the induced sequence of projective limits is exact, whence (2).
\end{proof}

\begin{cor}[{\cite[$\mathbf{0}_{\mathbf{III}}$, (13.2.2)]{EGA}}]\label{cor-ML4}
Consider the exact sequence of projective systems of abelian groups $(\ast)$ in $\ref{prop-ML3}$ with $I$ a directed set that has a cofinal and at most countable subset.
If {\boldmath $M'$} satisfies $\mbox{{\bf (ML)}}$, then the induced sequence
$$
0\longrightarrow\varprojlim_{i\in I}M'_i\longrightarrow\varprojlim_{i\in I}M_i\longrightarrow\varprojlim_{i\in I}M''_i\longrightarrow 0
$$
of abelian groups is exact. \hfill$\square$
\end{cor}

\subsubsection{Projective limit of sheaves}\label{subsub-projlimsheaves}
Let $X$ be a topological space, and consider a projective system $\{\mathscr{F}_i,p_{ij}\}_{i\in I}$ of sheaves (of sets, abelian groups, rings, etc.)\ on $X$ indexed by an ordered set $I$.
The projective limit sheaf $\varprojlim_{i\in I}\mathscr{F}_i$ in the category of sheaves (of sets, abelian groups, rings, etc.)\ is described as follows: For any open subset $U\subseteq X$ set
$$
\mathscr{F}(U)=\varprojlim_{i\in I}\mathscr{F}_i(U).
$$
This defines a presheaf $\mathscr{F}$, which is easily seen to be a sheaf.
The desired projective limit $\varprojlim_{i\in I}\mathscr{F}_i$ is given by the sheaf $\mathscr{F}$ with the canonical maps $p_i\colon\mathscr{F}\rightarrow\mathscr{F}_i$ for $i\in I$.
Notice that, unlike the case of inductive limits, the canonical map 
$$
(\varprojlim_{i\in I}\mathscr{F}_i)_x\longrightarrow\varprojlim_{i\in I}\mathscr{F}_{i,x}
$$
for $x\in X$ is not in general surjective nor injective.
\begin{prop}\label{prop-projlimsheafleftexact0}
$(1)$ For a topological space $X$ and an ordered set $I$, we have the canonical isomorphism 
$$
\varprojlim_{i\in I}\circ\Gamma_X\cong\Gamma_X\circ\varprojlim_{i\in I}
$$
of functors from $\ASh^{I^{\opp}}_X$ to $\Ab$, where $\Gamma_X$ is the global section functor $\Gamma_X(\mathscr{F})=\Gamma(X,\mathscr{F})$.

$(2)$ For a continuous map $f\colon X\rightarrow Y$ between topological spaces and an ordered set $I$, we have the canonical isomorphism
$$
\varprojlim_{i\in I}\circ f_{\ast}\cong f_{\ast}\circ\varprojlim_{i\in I}
$$
of functors from $\ASh^{I^{\opp}}_X$ to $\ASh_Y$. \hfill$\square$
\end{prop}

Finally, notice that the projective limit functor $\varprojlim_{i\in I}\colon\ASh^{I^{\opp}}_X\rightarrow\ASh_X$ is, as one can show similarly to \ref{prop-projlimleftexact}, left-exact.

\subsubsection{Canonical s-flasque resolution}\label{subsub-canstrictflasqueres}
Let $X$ be a topological space.
\begin{dfn}\label{dfn-strictflasque}{\rm 
A projective system $\{\mathscr{F}_i,p_{ij}\}_{i\in I}$ of abelian sheaves on $X$ (that is, an object of $\ASh^{I^{\opp}}_X$) is said to be {\em s-flasque}\index{s-flasque} if the following conditions are satisfied:
\begin{itemize}
\item[(a)] each $\mathscr{F}_i$ for $i\in I$ is flasque;
\item[(b)] for each $i<j$ the map $p_{ji}\colon\mathscr{F}_j\rightarrow\mathscr{F}_i$ is surjective and has the flasque kernel.
\end{itemize}}
\end{dfn}

Let $\{\mathscr{F}_i,p_{ij}\}_{i\in\N}$ be a projective system of abelian sheaves on $X$ indexed by $\N$.
We are going to construct an s-flasque resolution 
$$
0\longrightarrow\{\mathscr{F}_i\}_{i\in\N}\longrightarrow\{\mathscr{G}^{\bullet}_i\}_{i\in\N}.
$$

First we construct a resolution 
$$
0\longrightarrow\{\mathscr{F}_i\}_{i\in\N}\longrightarrow\{\mathscr{J}^0_i\}_{i\in\N}\longrightarrow\{\mathscr{J}_i^1\}_{i\in\N}\longrightarrow 0
$$
just as in the construction in \S\ref{subsub-canstrictres}; that is, we set 
\begin{equation*}
\begin{split}
\mathscr{J}^0_i&=\mathscr{F}_0\oplus\cdots\oplus\mathscr{F}_i\\
\mathscr{J}^1_i&=
\begin{cases}
0&\textrm{if $i=0$}\\
\mathscr{F}_0\oplus\cdots\oplus\mathscr{F}_{i-1}&\textrm{if $i>0$,}
\end{cases}
\end{split}
\end{equation*}
for each $i$ and define the maps similarly to that for canonical strict resolutions (cf.\ \S\ref{subsub-canstrictres}).
Next we take the canonical flasque resolutions of sheaves in $\{\mathscr{J}^0_i\}_{i\in\N}$ and in $\{\mathscr{J}^1_i\}_{i\in\N}$.
Thus we get the double complex $\{\mathscr{C}^{\bullet}(X,\mathscr{J}^{\bullet}_i)\}_{i\in\N}$ of projective systems.
The desired resolution $\{\mathscr{G}^{\bullet}_i\}_{i\in\N}$ is the single complex associated to this double complex; thus we have
$$
\{\mathscr{G}^q_i\}_{i\in\N}=
\begin{cases}
\{\mathscr{C}^0(X,\mathscr{J}^0_i)\}_{i\in\N}&\textrm{if $q=0$}\\
\{\mathscr{C}^q(X,\mathscr{J}^0_i)\oplus\mathscr{C}^{q-1}(X,\mathscr{J}^1_i)\}_{i\in\N}&\textrm{if $q>0$.}
\end{cases}
$$

We want to show that this resolution is s-flasque.
The condition (a) in \ref{dfn-strictflasque} is clear.
The condition (b) in \ref{dfn-strictflasque} follows from the fact that the formation of canonical flasque resolution of sheaves is an exact functor (\S\ref{subsub-canflasqueres}).
Indeed, if 
$$
0\longrightarrow\mathscr{K}\longrightarrow\mathscr{J}^k_j\longrightarrow\mathscr{J}^k_i\longrightarrow 0
$$
is exact for $k=0,1$ and $i<j$, then we have the induced exact sequnce
$$
0\longrightarrow\mathscr{C}^q(X,\mathscr{K})\longrightarrow\mathscr{C}^q(X,\mathscr{J}^k_j)\longrightarrow\mathscr{C}^q(X,\mathscr{J}^k_i)\longrightarrow 0.
$$
The condition (b) follows from this.

By the construction one sees that the formation of the above resolution defines a functor 
$$
\ASh^{\N^{\opp}}_X\longrightarrow\CC^+(\ASh^{\N^{\opp}}_X),
$$
which is, moreover, exact.
We call the resolution $0\rightarrow\{\mathscr{F}_i\}_{i\in\N}\rightarrow\{\mathscr{G}^{\bullet}_i\}_{i\in\N}$ the {\em canonical s-flasque resolution}\index{resolution!canonical s-flasque resolution@canonical s-flasque ---}.

\begin{prop}\label{prop-canstrictflasqueres}
Consider the exact sequence 
$$
\cdots\longrightarrow\{\mathscr{F}^{p-1}_i\}_{i\in\N}\longrightarrow\{\mathscr{F}^p_i\}_{i\in\N}\longrightarrow\{\mathscr{F}^{p+1}_i\}_{i\in\N}\longrightarrow\cdots
$$
of projective systems of abelian sheaves on $X$ indexed by $\N$.
Suppose that for any $p$ the projective system $\{\mathscr{F}^p_i\}_{i\in\N}$ is s-flasque.
Then the induced sequence
$$
\cdots\longrightarrow\mathscr{F}^{p-1}(X)\longrightarrow\mathscr{F}^p(X)\longrightarrow\mathscr{F}^{p+1}(X)\longrightarrow\cdots
$$
of abelian groups is exact, where $\mathscr{F}^p=\varprojlim_{i\in\N}\mathscr{F}^p_i$.
\end{prop}

\begin{proof}
The only non-trivial part of the proof is to show that any section $s\in\mathscr{F}^p(X)$ that is mapped to $0$ in $\mathscr{F}^{p+1}(X)$ lies in the image of $\mathscr{F}^{p-1}(X)\rightarrow\mathscr{F}^p(X)$.
Let $s_i$ be the image of $s$ by the projection $\mathscr{F}^p(X)\rightarrow\mathscr{F}^p_i(X)$ for each $i\in\N$.
Since $\mathscr{F}^{p-1}(X)=\varprojlim_{i\in\N}\mathscr{F}^{p-1}_i(X)$, we want to construct the compatible system $\{t_i\}_{i\in\N}$ of liftings $t_i\in\mathscr{F}^{p-1}_i(X)$ of $s_i$.
By \cite[Chap.\ II, 3.1.3]{Gode} one has a lifting $t_0\in\mathscr{F}^{p-1}_0(X)$ of $s_0$.
Suppose one has already constructed compatible liftings up to $t_{k-1}\in\mathscr{F}^{p-1}_{k-1}(X)$.
Take any lifting $t'_k$ of $s_k$ (possible due to \cite[Chap.\ II, 3.1.3]{Gode}).
The image $t'_{k-1}$ of $t'_k$ in $\mathscr{F}^{p-1}_{k-1}(X)$ has the same image in $\mathscr{F}^p_{k-1}(X)$ as $t_{k-1}$, and hence there exists an element $u_{k-1}\in\mathscr{F}^{p-2}_{k-1}(X)$ that is mapped to $t_{k-1}-t'_{k-1}$.
Since $\{\mathscr{F}^{p-2}_i\}_{i\in\N}$ is s-flasque, the projection $\mathscr{F}^{p-2}_k(X)\rightarrow\mathscr{F}^{p-2}_{k-1}(X)$ is surjective, and hence we have an element $u_k\in\mathscr{F}^{p-2}_k(X)$ that lifts $u_{k-1}$.
Let $v_k$ be the image of $u_k$ in $\mathscr{F}^{p-1}_k(X)$, and set $t_k=t'_k+v_k$.
Then $t_k$ is a lifting of both $t_{k-1}$ and $s_k$.
By induction we get the desired system of liftings.
\end{proof}

By the proposition we see that the canonical s-flasque resolution can be used to compute the right derived functors of the left-exact functor $\Gamma_X\circ\varprojlim_{i\in\N}$ ($\cong\varprojlim_{i\in\N}\circ\Gamma_X$) (cf.\ \ref{prop-projlimsheafleftexact0} (1)), that is, if 
$$
0\longrightarrow\{\mathscr{F}_i\}_{i\in\N}\longrightarrow\{\mathscr{G}^{\bullet}_i\}_{i\in\N}
$$
is the canonical s-flasque resolution, then we have
$$
\RD^q(\Gamma_X\circ\varprojlim_{i\in\N})(\{\mathscr{F}_i\}_{i\in\N})=\H^q(\Gamma(X,\mathscr{G}^{\bullet})),
$$
where $\mathscr{G}^{\bullet}=\varprojlim_{i\in\N}\mathscr{G}^{\bullet}_i$.

\begin{prop}\label{prop-strictflasque22}
{\rm (1)} If $\{\mathscr{G}_i\}_{i\in\N}$ is an s-flasque projective system of abelian sheaves, then $\mathscr{G}=\varprojlim_{i\in\N}\mathscr{G}_i$ is flasque.

{\rm (2)} Suppose 
$$
0\longrightarrow\{\mathscr{F}_i\}_{i\in\N}\longrightarrow\{\mathscr{G}^{\bullet}_i\}_{i\in\N}
$$
is an s-flasque resolution, and set $\mathscr{F}=\varprojlim_{i\in\N}\mathscr{F}_i$ and $\mathscr{G}^{\bullet}=\varprojlim_{i\in\N}\mathscr{G}^{\bullet}_i$.
Then the induced sequence 
$$
0\longrightarrow\mathscr{F}\longrightarrow\mathscr{G}^{\bullet}
$$
gives a flasque resolution of $\mathscr{F}$.
\end{prop}

\begin{proof}
(1) can be shown by an easy diagram similar to that in the proof of \ref{prop-canstrictflasqueres}.
By an argument similar to that in the proof of \ref{prop-canstrictflasqueres}, one can see that for any open subset $U$ the induced sequence
$$
0\longrightarrow\mathscr{F}(U)\longrightarrow\mathscr{G}^{\bullet}(U)
$$
is exact.
Taking stalks at every point, one concludes the exactness of $0\rightarrow\mathscr{F}\rightarrow\mathscr{G}^{\bullet}$, as desired.
\end{proof}

\subsubsection{Projective limit and cohomology}\label{subsub-projectivelimitandcohomology}
Let $X$ be a topological space, and $I$ a directed set containing a countable cofinal set.
We consider the following conditions for a projective system $\{\mathscr{F}_i\}_{i\in I}$ of sheaves of abelian groups on $X$ indexed by $I$:
\begin{itemize}
\item[{\bf (E1)}] there exists an open basis $\mathfrak{B}$ of $X$ such that for any $U\in\mathfrak{B}$ and $q\geq 0$ the projective system $\{\H^q(U,\mathscr{F}_i)\}_{i\in I}$ satisfies $\mbox{{\bf (ML)}};$\index{Mittag-Leffler condition}
\item[{\bf (E2)}] for any $x\in X$ and any $q>0$ we have 
$$
\varinjlim_{x\in U\in\mathfrak{B}}\varprojlim_{i\in I}\H^q(U,\mathscr{F}_i)=0,
$$
where $U$ varies over the set of open neighborhoods of $x$ in $\mathfrak{B}$. 
\end{itemize}

\begin{rem}\label{rem-coheffaceable}{\rm 
(1) In practice, these conditions are considered in the situation where $X$ has an open basis $\mathfrak{B}$ such that $\H^q(U,\mathscr{F}_i)=0$ for any $U\in\mathfrak{B}$, $i\in I$, and $q>0$.
In this case, the condition {\bf (E1)} for $q>0$ and the condition {\bf (E2)} are trivially satisfied, and hence one only need to check that the projective system $\{\H^0(U,\mathscr{F}_i)\}_{i\in I}$ satisfies $\mbox{{\bf (ML)}}$.

(2) If $X$ is a scheme, and if $\{\mathscr{F}_i\}_{i\in I}$ is a projective system consisting of quasi-coherent sheaves, then one can choose as $\mathfrak{B}$ the set of all affine open subsets\index{affine!affine open subspace@--- open subset (subspace)} of $X$ (due to \cite[$\mathbf{III}$, (1.3.1)]{EGA}).
In this case, for example, if $\{\mathscr{F}_i\}_{i\in I}$ is strict, then $\{\mathscr{F}_i\}_{i\in I}$ satisfies the conditions {\bf (E1)} and {\bf (E2)}; indeed, for $i\leq j$ and $U\in\mathfrak{B}$, since the kernel of the surjective map $\mathscr{F}_j\rightarrow\mathscr{F}_i$ is quasi-coherent, we deduce that $\H^0(U,\mathscr{F}_j)\rightarrow\H^0(U,\mathscr{F}_i)$ is surjective, and hence $\{\H^0(U,\mathscr{F}_i)\}_{i\in I}$ is strict.}
\end{rem}

\begin{prop}\label{prop-ML5ML5}
Let $X$ be a topological space, and $I$ a directed set $I$ containing a countable cofinal subset.
Let $\{\mathscr{F}_i\}_{i\in I}$ be a projective system of sheaves of abelian groups on $X$ indexed by $I$, and set $\mathscr{F}=\varprojlim_{i\in I}\mathscr{F}_i$.
Suppose $\{\mathscr{F}_i\}_{i\in I}$ satisfies the conditions {\bf (E1)} and {\bf (E2)}.
\begin{itemize}
\item[$(1)$] For any $q>0$ we have
$$
{\textstyle \varprojlim^{(q)}_{i\in I}\mathscr{F}_i=0.}
\leqno{(\ast)_q}
$$
\item[$(2)$] For any $q\geq 0$ there exist canonical isomorphisms
$$
\H^q(X,\mathscr{F})\cong\RD^q(\Gamma_X\circ\varprojlim_{i\in I})(\{\mathscr{F}_i\}_{i\in I})\cong\RD^q(\varprojlim_{i\in I}\circ\Gamma_X)(\{\mathscr{F}_i\}_{i\in I}).
$$
\end{itemize}
\end{prop}

\begin{rem}\label{rem-ML6*}{\rm 
Combined with \ref{prop-canstrictflasqueres}, the last statement (in the case $I=\N$) shows that one can use s-flasque\index{s-flasque} resolutions (\S\ref{subsub-canstrictflasqueres}) to compute the cohomologies $\H^q(X,\mathscr{F})$.
Indeed, if 
$$
0\longrightarrow\{\mathscr{F}_i\}_{i\in\N}\longrightarrow\{\mathscr{G}^{\bullet}_i\}_{i\in\N}
$$
is the canonical s-flasque resolution, we have 
$$
\H^q(X,\mathscr{F})\cong\H^q(\Gamma(X,\mathscr{G}^{\bullet})),
$$
where $\mathscr{G}^{\bullet}=\varprojlim_{i\in\N}\mathscr{G}^{\bullet}_i$.}
\end{rem}

\begin{proof}[Proof of Proposition {\rm \ref{prop-ML5ML5}}]
We assume without loss of generality that $I=\N$ (Exercise \ref{exer-finalcountable} and Exercise \ref{exer-cofinalML}).
Let $U$ be an open subset of $X$.
Since $\Gamma_U\circ\varprojlim_{i\in\N}\cong\varprojlim_{i\in\N}\circ\Gamma_U$, we have the following two spectral sequences that converge to the same infinity terms:
\begin{equation*}
\begin{split}
&{\textstyle {}^{{\it I}}\!E^{p,q}_2(U)=\varprojlim^{(p)}_{i\in\N}\H^q(U,\mathscr{F}_i)\Rightarrow{}^{{\it I}}\!E^{p+q}_{\infty}(U)=\RD^{p+q}(\varprojlim_{i\in\N}\circ\Gamma_U)(\{\mathscr{F}_i\}_{i\in\N}),}\\
&{\textstyle {}^{{\it II}}\!E^{p,q}_2(U)=\H^p(U,\varprojlim^{(q)}_{i\in\N}\mathscr{F}_i)\Rightarrow{}^{{\it II}}\!E^{p+q}_{\infty}(U)=\RD^{p+q}(\Gamma_U\circ\varprojlim_{i\in\N})(\{\mathscr{F}_i\}_{i\in\N})}
\end{split}
\end{equation*}
with ${}^{{\it I}}\!E^{p+q}_{\infty}(U)={}^{{\it II}}\!E^{p+q}_{\infty}(U)$.
For any point $x\in X$ we set 
$$
{}^{{\it I}}\!E_x=\varinjlim_{x\in U\in\mathfrak{B}}{}^{{\it I}}\!E(U),
$$
where $U$ runs through all open neighborhoods of $x$ in $\mathfrak{B}$.
Since the inductive limit functor is exact, this defines a spectral sequence that converges to the inductive limits of the corresponding $\infty$-terms.
The similar notation will be used also for ${}^{{\it II}}\!E(U)$.

From the conditions {\bf (E1)} and {\bf (E2)} we have
\begin{itemize}
\item[{\bf (E1)}] $\Rightarrow$ ${}^{{\it I}}\!E^{p,q}_{2,x}=0$ for $p>0$; 
\item[{\bf (E2)}] $\Rightarrow$ ${}^{{\it I}}\!E^{p,q}_{2,x}=0$ for $p=0$ and $q>0$.
\end{itemize}
Hence the only non-zero $E_2$-term of the spectral sequence ${}^{{\it I}}\!E_x$ is ${}^{{\it I}}\!E^{0,0}_{2,x}$.
This shows the following:
$$
{}^{{\it I}}\!E^{p+q}_{\infty,x}={}^{{\it II}}\!E^{p+q}_{\infty,x}=
\begin{cases}
\varinjlim_{x\in U}\varprojlim_{i\in\N}\Gamma(U,\mathscr{F}_i)&\textrm{if $p+q=0$}\\
0&\textrm{otherwise.}
\end{cases}
\leqno{(\dagger)}
$$

Now we claim that ${}^{{\it II}}\!E^{p,0}_{2,x}=0$ for $p>0$.
Take the canonical s-flasque resolution $0\rightarrow\{\mathscr{F}_i\}_{i\in\N}\rightarrow\{\mathscr{G}^{\bullet}_i\}_{i\in\N}$.
By \ref{prop-canstrictflasqueres} we have 
$$
0={}^{{\it II}}\!E^p_{\infty,x}=\varinjlim_{x\in U}\H^p(\Gamma(U,\mathscr{G}^{\bullet})),
$$
where $\mathscr{G}^{\bullet}=\varprojlim_{i\in\N}\mathscr{G}^{\bullet}_i$.
But then, since $0\rightarrow\mathscr{F}\rightarrow\mathscr{G}^{\bullet}$ gives a flasque resolution (\ref{prop-strictflasque22} (2)), we deduce that 
$$
{}^{{\it II}}\!E^{p,0}_{2,x}=\varinjlim_{x\in U}\H^p(U,\mathscr{F})=\varinjlim_{x\in U}\H^p(\Gamma(U,\mathscr{G}^{\bullet}))=0,
$$
as desired.

Now we start proving (1).
This is shown by induction with respect to $q$.
Since 
$$
0={}^{{\it II}}\!E^{-2,2}_{2,x}\stackrel{d}{\longrightarrow}{}^{{\it II}}\!E^{0,1}_{2,x}\stackrel{d}{\longrightarrow}{}^{{\it II}}\!E^{2,0}_{2,x}=0,
$$
one has ${}^{{\it II}}\!E^{0,1}_{2,x}={}^{{\it II}}\!E^{0,1}_{\infty,x}$.
But ${}^{{\it II}}\!E^{0,1}_{\infty,x}$ is a subquotient of ${}^{{\it II}}\!E^1_{\infty,x}$, which is zero by $(\dagger)$, and hence we have $(\varprojlim^{(1)}_{i\in\N}\mathscr{F}_i)_x={}^{{\it II}}\!E^{0,1}_{2,x}=0$.
Since this is true for any point $x\in X$, we get the desired vanishing $(\ast)_1$.
By induction we assume that $(\ast)_k$ is true for $k=1,\ldots,q-1$.
Then for $k\leq q-1$ we have ${}^{{\it II}}\!E^{p,k}_{2,x}=0$.
Since 
$$
0={}^{{\it II}}\!E^{-2,q+1}_{2,x}\stackrel{d}{\longrightarrow}{}^{{\it II}}\!E^{0,q}_{2,x}\stackrel{d}{\longrightarrow}{}^{{\it II}}\!E^{2,q-1}_{2,x}=0,
$$
one has ${}^{{\it II}}\!E^{0,q}_{2,x}={}^{{\it II}}\!E^{0,q}_{\infty,x}$.
But by the similar reasoning as above, one sees the last term is zero for any $x\in X$, and hence we have $(\ast)_q$, as desired.

Now we look at the global spectral sequence ${}^{{\it II}}\!E(X)$.
Since (1) holds, this spectral sequence degenerates at $E_2$-terms.
By this and (1), we get (2).
\end{proof}

\begin{cor}[{cf.\ \cite[$\mathbf{0}_{\mathbf{III}}$, (13.3.1)]{EGA}}]\label{cor-ML5}
Let $X$ be a topological space, $I$ a directed set $I$ containing a countable cofinal subset, and $\{\mathscr{F}_i\}_{i\in I}$ a projective system of abelian sheaves on $X$ indexed by $I$.
Suppose $\{\mathscr{F}_i\}_{i\in I}$ satisfies the conditions {\bf (E1)} and {\bf (E2)}.
Then for any $q>0$ the canonical morphism 
$$
\H^q(X,\mathscr{F})\longrightarrow\varprojlim_{i\in I}\H^q(X,\mathscr{F}_i)\leqno{(\ast\ast)_q}
$$
is surjective.
If, moreover, the projective system $\{\H^{q-1}(X,\mathscr{F}_i)\}_{i\in I}$ satisfies $\mbox{{\bf (ML)}}$, then $(\ast\ast)_q$ is bijective.
\end{cor}

\begin{rem}\label{rem-ML6ML6}{\rm 
(1) Notice that the isomorphy of $(\ast\ast)_0$, not mentioned in the above statement, is always true due to \ref{prop-projlimsheafleftexact0} (1).

(2) In \cite[$\mathbf{0}_{\mathbf{III}}$, (13.3.1)]{EGA} it is assumed furthermore that the projective system $\{\mathscr{F}_i\}_{i\in I}$ is strict.
But this assumption is not necessary, since one can always replace $\{\mathscr{F}_i\}_{i\in I}$ by a complex of strict systems.}
\end{rem}

\begin{proof}[Proof of Corollary {\rm \ref{cor-ML5}}]
We use the spectral sequence $E={}^I\!E(X)$ as in the proof of \ref{prop-ML5ML5}.
By \ref{prop-ML5ML5} (2) this spectral sequence converges to $E^{p+q}_{\infty}=\H^{p+q}(X,\mathscr{F})$.
By \ref{prop-ML6} (1) one sees that this spectral sequence degenerates at $E_2$-terms.
Hence we have the exact sequence
$$
0\longrightarrow E^{1,q-1}_2\longrightarrow\H^q(X,\mathscr{F})\longrightarrow E^{0,q}_2\longrightarrow 0,
$$
which shows, in particular, the surjectivity of $(\ast\ast)_q$.
If $\{\H^{q-1}(X,\mathscr{F}_i)\}_{i\in\N}$ satisfies {\bf (ML)}, then by \ref{prop-ML6} (2) we have $E^{1,q-1}_2=0$, whence the isomorphism $\H^q(X,\varprojlim_k\mathscr{F}_k)\cong E^{0,q}_2$.
\end{proof}

\begin{prop}\label{prop-ML7}
Consider the situation as in $\ref{prop-ML5ML5}$ with $I=\N$, and let $q\geq 0$ be an integer.
For each $i\in\N$ set $\mathscr{N}_i=\ker(\mathscr{F}\rightarrow\mathscr{F}_i)$, and let $F^i$ be the image of the map $\H^q(X,\mathscr{N}_i)\rightarrow\H^q(X,\mathscr{F})$.
Then the following conditions are equivalent$:$
\begin{itemize}
\item[{\rm (a)}] the map $(\ast\ast)_q$ in $\ref{cor-ML5}$ is injective$;$ 
\item[{\rm (b)}] the map $(\ast\ast)_q$ in $\ref{cor-ML5}$ is bijective$;$ 
\item[{\rm (c)}] $\bigcap_{i\in\N}F^i=0;$ in other words, $\H^q(X,\mathscr{F})$ is separated with respect to the filtration $F^{\bullet}=\{F^i\}_{i\in\N}$.
\end{itemize}
\end{prop}

\begin{proof}
The equivalence of (a) and (b) follows from \ref{cor-ML5}.
Since the map $\H^q(X,\mathscr{F})/F^i\rightarrow\H^q(X,\mathscr{F}_i)$ is injective, we have the injective homomorphism 
$$
\varprojlim_{i\in\N}\H^q(X,\mathscr{F})/F^i\longhookrightarrow\varprojlim_{i\in\N}\H^q(X,\mathscr{F}_i).\leqno{(\ddagger)}
$$
The canonical map $(\ast\ast)_q$ clearly factors through $\varprojlim_{k\in\N}\H^q(X,\mathscr{F})/F^i$ by the canonical map 
$$
\H^q(X,\mathscr{F})\longrightarrow\varprojlim_{i\in\N}\H^q(X,\mathscr{F})/F^i\leqno{(\ddagger\ddagger)}
$$
followed by the map $(\ddagger)$, and hence $(\ddagger)$ is bijective.
Therefore, $(\ast\ast)_q$ is injective if and only if $(\ddagger\ddagger)$ is injective, and the last condition is equivalent to (c).
\end{proof}

\begin{prop}\label{prop-ML5rel}
Let $f\colon X\rightarrow Y$ be a continuous mapping between topological spaces, $I$ a directed set containing a countable cofinal subset, and $\{\mathscr{F}_i\}_{i\in I}$ a projective system of abelian sheaves on $X$ indexed by $I$.
Set $\mathscr{F}=\varprojlim_{i\in I}\mathscr{F}_i$.
Suppose that$:$
\begin{itemize}
\item $\{\mathscr{F}_i\}_{i\in I}$ satisfies the conditions {\bf (E1)} and {\bf (E2)}$;$
\item the projective systems $\{\RD^qf_{\ast}\mathscr{F}_i\}_{i\in I}$ $(q\geq 0)$ of abelian sheaves on $Y$ satisfies the conditions {\bf (E1)} and {\bf (E2)}.
\end{itemize}
Then the canonical morphism 
$$
\RD^qf_{\ast}\mathscr{F}\longrightarrow\varprojlim_{i\in I}\RD^qf_{\ast}\mathscr{F}_i
$$
is an isomorphism.
\end{prop}

By what we have seen in \ref{rem-coheffaceable}, the hypotheses in the proposition are satisfied in the following situation: $f\colon X\rightarrow Y$ is a coherent ($=$ quasi-compact and quasi-separated) morphism of schemes, and $\{\mathscr{F}_i\}_{i\in I}$ is a strict projective system consisting of quasi-coherent sheaves on $X$ such that the induced system $\{\RD^qf_{\ast}\mathscr{F}_i\}_{i\in I}$ is also strict for $q\geq 0$.
Indeed, by the well-known fact in the theory of schemes (cf.\ \ref{prop-cohqcoh-0} below), the sheaves $\RD^qf_{\ast}\mathscr{F}_i$ for $q\geq 0$ and $i\in I$ are quasi-coherent.

\begin{proof}[Proof of Proposition {\rm \ref{prop-ML5rel}}]
We consider the following two spectral sequences:
\begin{equation*}
\begin{split}
&{\textstyle {}^{{\it I}}\!\mathscr{E}^{p,q}_2=\varprojlim^{(p)}_{i\in\N}\RD^qf_{\ast}\mathscr{F}_i\Rightarrow{}^{{\it I}}\!\mathscr{E}^{p+q}_{\infty}=\RD^{p+q}(\varprojlim_{i\in\N}\circ f_{\ast})(\{\mathscr{F}_i\}_{i\in\N}),}\\
&{\textstyle {}^{{\it II}}\!\mathscr{E}^{p,q}_2=\RD^pf_{\ast}\varprojlim^{(q)}_{i\in\N}\mathscr{F}_i\Rightarrow{}^{{\it II}}\!\mathscr{E}^{p+q}_{\infty}=\RD^{p+q}(f_{\ast}\circ\varprojlim_{i\in\N})(\{\mathscr{F}_i\}_{i\in\N})}\rlap{.}
\end{split}
\end{equation*}
By \ref{prop-projlimsheafleftexact0} (2) we have ${}^{{\it I}}\!\mathscr{E}^{p+q}_{\infty}={}^{{\it II}}\!\mathscr{E}^{p+q}_{\infty}$.
By \ref{prop-ML5ML5} (1) we have ${}^{{\it II}}\!\mathscr{E}^{p+q}_{\infty}\cong\RD^{p+q}f_{\ast}\mathscr{F}$ and ${}^{{\it I}}\!\mathscr{E}^{p+q}_{\infty}\cong\varprojlim_{i\in I}\RD^{p+q}f_{\ast}\mathscr{F}_i$.
\end{proof}
\index{limit!projective limit@projective ---|)}

\subsection{Coherent rings and modules}\label{sub-cohringsmodules}
\index{coherent!coherent module@--- module|(}\index{coherent!coherent ring@--- ring|(}
\begin{dfn}[{cf.\ \cite[Chap.\ I, \S 3, Exercise 11]{Bourb1}}]\label{dfn-cohringsmodules1}{\rm 

(1) Let $A$ be a ring. A finitely generated $A$-module $M$ is said to be {\em coherent} if every finitely generated $A$-submodule of $M$ is finitely presented.
We denote by $\Coh_A$ the full subcategory of $\Mod_A$ consisting of coherent $A$-modules.

(2) A ring $A$ is said to be {\em coherent} if it is coherent as an $A$-module or, what amounts to the same, every finitely generated ideal of $A$ is finitely presented.}
\end{dfn}

For example, any Noetherian ring is a coherent ring, and any finitely generated module over a Noetherian ring is a coherent module.
It is easy to see that, if $A$ is a coherent ring and $S\subseteq A$ is a multiplicative subset, then $S^{-1}A$ is again a coherent ring.
In Exercise \ref{exer-coherencedirectlimit} we will see that the integral closure $\til{\Z}$ of $\Z$ in the ring of algebraic numbers $\ovl{\Q}$ is a coherent ring, which is, however, not Noetherian (Exercise \ref{exer-algintegernonNoe}).

Let us first insert here an easy proposition; the second part follows promptly from \cite[Chap.\ I, \S3.6, Prop.\ 11]{Bourb1}:
\begin{prop}\label{prop-cohringsmodules1a}
$(1)$ Let $A$ and $B$ be rings, $M$ a finitely generated $A$-module, and $N$ a finitely generated $B$-module.
Then $M\times N$ is a coherent $A\times B$-module if and only if $M$ and $N$ are coherent over $A$ and $B$, respectively.
In particular, $A\times B$ is coherent if and only if $A$ and $B$ are coherent.

$(2)$ Let $A$ be a ring, and $B$ a faithfully flat $A$-algebra.
Then a finitely generated $A$-module $M$ is coherent if and only if $M\otimes_AB$ is a coherent $B$-module.
In particular, $A$ is coherent if $B$ is coherent. \hfill$\square$
\end{prop}

\begin{prop}\label{prop-cohringsmodules1}
The following conditions for a ring $A$ are equivalent$:$
\begin{itemize}
\item[{\rm (a)}] $A$ is coherent$;$ 
\item[{\rm (b)}] any finitely presented $A$-module is coherent$;$ 
\item[{\rm (c)}] the category of finitely presented $A$-modules is an abelian subcategory of the abelian category $\Mod_A$ of $A$-modules$;$ 
\item[{\rm (d)}] if $f\colon M\rightarrow N$ is a homomorphism between finitely presented $A$-modules, $\ker(f)$ is finitely presented$;$ 
\item[{\rm (e)}] if the $A$-modules $M$ and $N$ sitting in the exact sequence
$$
0\longrightarrow L\longrightarrow M\longrightarrow N\longrightarrow 0\leqno{(\ast)}
$$
of $A$-modules are finitely presented, the so is $L;$ 
\item[{\rm (f)}] if two of $L$, $M$, and $N$ in the exact sequence $(\ast)$ are finitely presented, so is the rest$;$ 
\item[{\rm (g)}] for any exact sequence of $A$-modules
$$
M_1\longrightarrow M_2\longrightarrow M_3\longrightarrow M_4\longrightarrow M_5, 
$$
if $M_1$ is finitely generated, and if $M_2$, $M_4$, and $M_5$ are finitely presented, then $M_3$ is finitely presented$;$ 
\item[{\rm (h)}] for any complex 
$$
M^{\bullet}=(\cdots\longrightarrow M^{k-1}\longrightarrow M^k\longrightarrow M^{k+1}\longrightarrow\cdots)
$$
of $A$-modules, if every $M^k$ is finitely presented, then $\H^q(M^{\bullet})$ is finitely presented for any $q$. \hfill$\square$
\end{itemize}
\end{prop}

The proof is easy with the aid of the following elementary lemma:
\begin{lem}\label{lem-pf1}
Let $A$ be a ring, and 
$$
0\longrightarrow L\longrightarrow M\longrightarrow N\longrightarrow 0\leqno{(\ast)}
$$
an exact sequence of $A$-modules.

$(1)$ If $L$ is a finitely generated and $M$ is a finitely presented, then $N$ is finitely presented.

$(2)$ If $L$ and $N$ are finitely presented, then $M$ is finitely presented. \hfill$\square$
\end{lem}

\begin{cor}\label{cor-cohringsmodules11}
Let $A$ be a coherent ring.
Then an $A$-module $M$ is coherent if and only if it is finitely presented. \hfill$\square$
\end{cor}

If $A$ is a coherent ring, then any coherent $A$-module has a free resolution; hence we have:
\begin{cor}\label{cor-cohringsmodules12}
Let $A$ be a coherent ring.
Then for any coherent $A$-modules $M$ and $N$, $\Ext^p_A(M,N)$ and $\Tor^A_p(M,N)$ are coherent for $p\geq 0$. \hfill$\square$
\end{cor}

The following notion will be particularly important later:
\begin{dfn}\label{dfn-universallycoherent}{\rm A ring $A$ is said to be {\em universally coherent}\index{coherent!coherent ring@--- ring!universally coherent ring@universally --- ---} if every finitely presented $A$-algebra is coherent.}
\end{dfn}

Clearly, if $A$ is universally coherent, then any finitely presented $A$-algebra is again universally coherent.
Note also that, if $A$ is universally coherent, then any localization of $A$ is again universally coherent.
\index{coherent!coherent ring@--- ring|)}\index{coherent!coherent module@--- module|)}

\addcontentsline{toc}{subsection}{Exercises}
\subsection*{Exercises}
\begin{exer}\label{exer-injlimcohcoh}{\rm 
Let $X$ be a coherent topological space.
A sheaf $\mathscr{F}$ of sets on $X$ is said to be {\it quasi-flasque} (\cite{Kempf}) if, for any quasi-compact open subset $U\subseteq X$, the restriction map $\mathscr{F}(X)\rightarrow\mathscr{F}(U)$ is surjective.

(1) Let $0\rightarrow\mathscr{F}'\rightarrow\mathscr{F}\rightarrow\mathscr{F}''\rightarrow 0$ be an exact sequence of abelian sheaves on $X$, and suppose $\mathscr{F}'$ is quasi-flasque.
Then show that 
$$
0\longrightarrow\Gamma(X,\mathscr{F}')\longrightarrow\Gamma(X,\mathscr{F})\longrightarrow\Gamma(X,\mathscr{F}'')\longrightarrow 0
$$
is exact.

(2) Show that, for a quasi-flasque sheaf of abelian groups $\mathscr{F}$, $\H^q(X,\mathscr{F})=0$ for $q>0$.

(3) Let $\{\mathscr{G}_i,\varphi_{ij}\}_{i\in I}$ be a filtered inductive system of quasi-flasque sheaves of sets indexed by a directed set $I$.
Show that the inductive limit $\mathscr{G}=\varinjlim_{i\in I}\mathscr{G}_i$ is quasi-flasque.}
\end{exer}

\begin{exer}\label{exer-cofinalML}{\rm 
Let $\{A_i,f_{ij}\}_{i\in I}$ be a projective system of sets indexed by a directed set $I$, and $J\rightarrow I$ a cofinal ordered map.
Show that, if $\{A_i,f_{ij}\}_{i\in I}$ satisfies the condition {\bf (ML)}, then so does $\{A_i,f_{ij}\}_{i\in J}$.}
\end{exer}

\begin{exer}\label{exer-coherencedirectlimit}{\rm 
Let $\{A_i,\phi_{ij}\}_{i\in I}$ be an inductive system of rings indexed by a directed set.
Assume that each $A_i$ is coherent (resp.\ universally coherent) and that each transition map is flat.
Then show that $A=\varinjlim_{i\in I}A_i$ is coherent (resp.\ universally coherent).}
\end{exer}

\begin{exer}\label{exer-algintegernonNoe}{\rm 
Show that the subring of all algebraic integers in $\ovl{\Q}$ is coherent, but not Noetherian.}
\end{exer}

\begin{exer}\label{exer-adhesivecounterexas}{\rm 
Let $k$ be a field, and consider the polynomial ring $k[x_1,x_2,\ldots]$ of countably many indeterminacies.
Set $J=(x_1x_2,x_1x_3,\ldots, x_1x_n,\ldots)$.
Then show that the ring $A=k[x_1,x_2,\ldots]/J$ is not coherent.}
\end{exer}


\section{Ringed spaces}\label{sec-ringedspaces}\index{ringed space|(}
In most of ordinary commutative geometries, spaces are `visualized' by means of locally ringed spaces (as referred to as {\em standard visualization} in Introduction).
Homological algebra of $\O_X$-modules, especially that of quasi-coherent and coherent sheaves, is an important tool for analyzing the spaces.
Coherent sheaves are particularly useful if the structure sheaf $\O_X$ itself is coherent.
So the coherency of the structure sheaf is one of the fundamental conditions for (locally) ringed spaces.
In \S\ref{sub-ringedspgenpre} we discuss ringed spaces satisfying this condition in general, which we call {\em cohesive ringed spaces}.

From \S\ref{sub-projlimringedspsheaf} we begin the study of module sheaves and cohomologies of them in the context of filtered projective limits.
Notice that most of the results in these sections are, in fact, rehashes of what are already done in \cite[Expos\'e VI]{SGA4-2} in the topos-theoretic language.

\subsection{Generalities}\label{sub-ringedspgenpre}
\subsubsection{Ringed spaces and locally ringed spaces}\label{subsub-ringedsplocalringedsp}
A {\em ringed space}\index{ringed space} is a couple $(X,\O_X)$ consisting of a topological space $X$ and a sheaf of rings $\O_X$ on $X$, called the {\em structure sheaf}.
Given two ringed spaces $X=(X,\O_X)$ and $Y=(Y,\O_Y)$, a {\em morphism}\index{morphism of ringed spaces@morphism (of ringed spaces)} of ringed spaces from $X$ to $Y$ is a couple $(f,\varphi)$ consisting of a continuous map $f\colon X\rightarrow Y$ and a morphism $\varphi\colon f^{-1}\O_Y\rightarrow\O_X$ of sheaves of rings on $X$ or, equivalently by adjunction, a morphism $\O_Y\rightarrow f_{\ast}\O_X$ of sheaves of rings on $Y$.
We denote by $\Rsp$ the category of ringed spaces.

A ringed space $(X,\O_X)$ is said to be a {\em locally ringed space}\index{ringed space!locally ringed space@locally --- ---} if for any point $x\in X$ the ring $\O_{X,x}$ is a local ring.
In this case we denote by $\m_{X,x}$ the maximal ideal of the local ring $\O_{X,x}$.
A morphism $(f,\varphi)\colon (X,\O_X)\rightarrow(Y,\O_Y)$ of ringed spaces, where $X=(X,\O_X)$ and $Y=(Y,\O_Y)$ are locally ringed spaces, is said to be {\em local}\index{morphism of ringed spaces@morphism (of ringed spaces)!local morphism of ringed spaces@local ---} if for any $x\in X$ the map $\varphi_x\colon\O_{Y,f(x)}\rightarrow\O_{X,x}$ is a local homomorphism (that is, $\varphi_x(\m_{Y,f(x)})\subseteq\m_{X,x}$).
We denote by $\LRsp$ the category of locally ringed spaces and local morphisms.

\medskip\noindent
{\bf Convention.} {\sl In the sequel, whenever we deal with locally ringed spaces, simply by a morphism we mean a local morphism, unless otherwise clearly stated.}

\medskip
Let $(X,\O_X)$ be a ringed space.
An {\em open $($ringed$)$ subspace}\index{subspace!open subspace of a ringed space@open --- (of a ringed space)} of $(X,\O_X)$ is a ringed space of the form $(U,\O_X|_U)$, where $U$ is an open subset of $X$.
An {\em open immersion}\index{immersion!open immersion of ringed spaces@open --- (of ringed spaces)} is a morphism of ringed spaces $(Y,\O_Y)\rightarrow(X,\O_X)$ that factors through the canonical map from an open subspace $(U,\O_U)$ of $(X,\O_X)$ by an isomorphism $(Y,\O_Y)\stackrel{\sim}{\rightarrow}(U,\O_U)$.

A ringed space $(X,\O_X)$ is said to be {\em reduced}\index{ringed space!reduced ringed space@reduced --- ---} if for any point $x\in X$ the ring $\O_{X,x}$ has no non-zero nilpotent element; in other words, if $\mathscr{N}_X$ is the subsheaf of $\O_X$ consisting of nilpotent sections, we have $\mathscr{N}_X=0$.

\subsubsection{Generization map}\label{subsub-generizationmap}
Let $X$ be a topological space, and $\mathscr{F}$ a sheaf (of sets) on $X$.
If $y\in G_x$ is a generization of a point $x\in X$, then any open neighborhood of $x$ is an open neighborhood of $y$.
Hence we have a canonical map
$$
\mathscr{F}_x\longrightarrow\mathscr{F}_y
$$
between the stalks, which we call the {\em generization map}\index{generization!generization map@--- map}.
In particular, if $X=(X,\O_X)$ is a ringed space, then the generization map 
$$
\O_{X,x}\longrightarrow\O_{X,y}
$$
is a ring homomorphism.
If, moreover, $X$ is a locally ringed space, then the generization map induces a local homomorphism
$$
(\O_{X,x})_{\mathfrak{q}_y}\longrightarrow\O_{X,y},\leqno{(\ast)}
$$
where $\mathfrak{q}_y$ is the pull-back of the maximal ideal $\m_{X,y}$ at $y$.

\begin{exas}\label{exas-generizationmapsch}{\rm 
(1) If $X$ is a scheme, then $(\ast)$ is an isomorphism.

(2) If $X$ is a locally Noetherian formal scheme (\cite[$\mathbf{I}$, (10.4.2)]{EGA}) or, more generally, {\em locally universally rigid-Noetherian formal schemes}\index{formal scheme!universally rigid-Noetherian formal scheme@universally rigid-Noetherian ---!locally universally rigid-Noetherian formal scheme@locally --- ---} (defined later in {\bf \ref{ch-formal}}.\ref{dfn-formalsch} below), then $(\ast)$ is faithfully flat (Exercise \ref{exer-generizationmapsch}). }
\end{exas}

\subsubsection{Sheaves of modules}\label{subsub-ringedspsheavesmodules}
For a ringed space $X=(X,\O_X)$ we denote by $\Mod_X$ the category of $\O_X$-modules.
This is an abelian category with tensor products and internal Hom's.
Any morphism $(f,\varphi)\colon(X,\O_X)\rightarrow(Y,\O_Y)$ of ringed spaces induces the two functors, adjoint to each other:
$$
\xymatrix{\Mod_X\ar@<.5ex>[r]^(.5){f_{\ast}}&\Mod_Y\rlap{;}\ar@<.5ex>[l]^(.5){f^{\ast}}}
$$
here $f^{\ast}\mathscr{G}$ for an $\O_Y$-module $\mathscr{G}$ is defined as $f^{\ast}\mathscr{G}=f^{-1}\mathscr{G}\otimes_{f^{-1}\O_Y,\varphi}\O_X$.

\begin{dfn}\label{dfn-idealpullback}{\rm 
Let $(f,\varphi)\colon(X,\O_X)\rightarrow(Y,\O_Y)$ be a morphism of ringed spaces, and $\mathscr{J}$ an ideal sheaf of $\O_Y$.
The {\em ideal pull-back}\index{ideal pull-back} of $\mathscr{J}$ is the ideal sheaf of $\O_X$ generated by the image of $f^{-1}\mathscr{J}$ by the map $\varphi$.}
\end{dfn}

The ideal pull-back is denoted by $(f^{-1}\mathscr{J})\O_X$ or more simply by $\mathscr{J}\O_X$.

\begin{dfn}\label{dfn-ftypeoxmodule}{\rm 
(1) Let $\mathscr{F}$ be an $\O_X$-module, and $n\in\Z$ a non-negative integer.
We say that $\mathscr{F}$ is {\em of $($finite$)$ $n$-presentation}\index{module (OX)@module ($\O_X$-{---})!module of finite n-presentation@--- of (finite) $n$-presentation} or {\em $n$-presented} if for any point $x\in X$ there exists an open neighborhood $U$ of $x$ on which there exists an exact sequence of $\O_U$-modules
$$
\mathscr{E}^{n}\longrightarrow\cdots\longrightarrow\mathscr{E}^1\longrightarrow\mathscr{E}^0\longrightarrow\mathscr{F}|_U\longrightarrow 0,
$$
where each $\mathscr{E}^i$ $(0\leq i\leq n)$ is a free $\O_U$-module of finite rank. 

(2) An $\O_X$-module $\mathscr{F}$ is said to be {\em of finite type}\index{module (OX)@module ($\O_X$-{---})!module of finite type@--- of finite type} if it is of $0$-presentation.
If it is of $1$-presentation, we say that $\mathscr{F}$ is {\em of finite presentation}\index{module (OX)@module ($\O_X$-{---})!module of finite presentation@--- of finite presentation} or {\em finitely presented}.
Finally, if $\mathscr{F}$ is $n$-presented for any non-negative integer $n$, we say that $\mathscr{F}$ is {\em of $\infty$-presentation}\index{module (OX)@module ($\O_X$-{---})!module of infty presentation@--- of $\infty$-presentation} or {\em $\infty$-presented}.}
\end{dfn}

Note that, if $(f,\varphi)\colon(X,\O_X)\rightarrow(Y,\O_Y)$ is a morphism of ringed spaces, then the functor $f^{\ast}\colon\Mod_Y\rightarrow\Mod_X$ maps $\O_Y$-modules of finite type (resp.\ of finite presentation) to $\O_X$-modules of finite type (resp.\ of finite presentation).

\begin{dfn}\label{dfn-qcohcohringedspace}{\rm 
(1) An $\O_X$-module $\mathscr{F}$ is said to be {\em quasi-coherent}\index{quasi-coherent!quasi-coherent OX module@--- sheaf (on a ringed space)} if for any $x\in X$ there exists an open neighborhood $U$ of $x$ on which we have an exact sequence of $\O_X|_U$-modules of the following form:
$$
\O^{\oplus J}_X|_U\longrightarrow\O^{\oplus I}_X|_U\longrightarrow\mathscr{F}|_U\longrightarrow 0,
$$
where $\O^{\oplus I}_X$ for a set $I$ (and $\O^{\oplus J}_X$ similarly) denotes the direct sum of copies of $\O_X$ indexed by $I$.

(2) An $\O_X$-module $\mathscr{F}$ is said to be {\em coherent}\index{coherent!coherent sheaf@--- sheaf (on a ringed space)} if it satisfies the following conditions:
\begin{itemize}
\item[{\rm (a)}] $\mathscr{F}$ is of finite type;
\item[{\rm (b)}] the kernel of arbitrary morphism $\O^{\oplus n}_X|_U\longrightarrow\mathscr{F}|_U$, where $U\subseteq X$ is an open subset and $n\geq 0$, is of finite type.
\end{itemize}}
\end{dfn}

Finitely presented $\O_X$-modules are quasi-coherent, but finite type $\O_X$-modules are not in general quasi-coherent.
Coherent $\O_X$-modules are finitely presented and hence quasi-coherent, but finitely presented $\O_X$-modules are not necessarily coherent.
We denote by $\QCoh_X$ (resp.\ $\Coh_X$) the full subcategory of $\Mod_X$ consisting of quasi-coherent (resp.\ coherent) $\O_X$-modules.
These are abelian {\em thick}\index{subcategory!thick subcategory@thick ---} (\S\ref{sub-subcategoryderived}) subcategories of $\Mod_X$ (cf.\ \cite[(1.4.7)]{EGAInew}).
Note that, if $(f,\varphi)\colon(X,\O_X)\rightarrow(Y,\O_Y)$ is a morphism of ringed spaces, then the functor $f^{\ast}\colon\Mod_Y\rightarrow\Mod_X$ maps $\QCoh_Y$ to $\QCoh_X$.
Notice also that, if $\mathscr{F}$ and $\mathscr{G}$ are coherent $\O_X$-modules, then so are $\mathscr{F}\otimes_{\O_X}\mathscr{G}$ and $\lHom_{\O_X}(\mathscr{F},\mathscr{G})$.

\begin{dfn}\label{dfn-flatnessringedsp}{\rm Let $f=(f,\varphi)\colon(X,\O_X)\rightarrow(Y,\O_Y)$ be a morphism of ringed spaces.

(1) An $\O_X$-module $\mathscr{F}$ is said to be {\em $f$-flat $($or $Y$-flat$)$}\index{flatness} at a point $x\in X$ if $\mathscr{F}_x$ is flat as a module over $\O_{Y,f(x)}$.
If $\mathscr{F}$ is $f$-flat at all points of $X$, we simply say that $\mathscr{F}$ is $f$-flat (or $Y$-flat).
In particular, if $(X,\O_X)=(Y,\O_Y)$ and $f=\id$, then we say $\mathscr{F}$ is {\em flat}. 

(2) If $\O_X$ is $f$-flat, that is, if the morphism $\varphi_x\colon\O_{Y,f(x)}\rightarrow\O_{X,x}$ is flat for any point $x\in X$, then the $f=(f,\varphi)$ is said to be {\em flat}\index{morphism of ringed spaces@morphism (of ringed spaces)!flat morphism of ringed spaces@flat ---}.}
\end{dfn}

If $f$ is flat, then the functor $f^{\ast}\colon\Mod_Y\rightarrow\Mod_X$ is exact.

\begin{prop}\label{prop-projlimcohsch22}
Let $(X,\O_X)$ be a ringed space, and $\mathscr{F}$ an $\O_X$-module.
Let $\{\mathscr{G}_i\}_{i\in I}$ be a filtered inductive system of $\O_X$-modules indexed by a directed set $I$.
Consider the natural map
$$
\Phi\colon\varinjlim_{i\in I}\Hom_{\O_X}(\mathscr{F},\mathscr{G}_i)\longrightarrow\Hom_{\O_X}(\mathscr{F},\varinjlim_{i\in I}\mathscr{G}_i).
$$

{\rm (1)} If $X$ is quasi-compact and $\mathscr{F}$ is of finite type, then $\Phi$ is injective.

{\rm (2)} If $X$ is coherent {\rm (\ref{dfn-quasicompact1})} and $\mathscr{F}$ is finitely presented, then $\Phi$ is bijective. \hfill$\square$
\end{prop}

Similarly to \ref{prop-projlimcohsch21}, this proposition can be seen as a special case of \cite[Expos\'e  VI, Th\'eor\`em 1.23]{SGA4-2}.
One can prove it in a similar way to the proof of \ref{prop-projlimcohsch21}, and the checking is left to the reader.

\subsubsection{Cohesive ringed spaces}\label{subsub-cohesiveringedspaces}
For a ringed space $X$ the structure sheaf $\O_X$ is always quasi-coherent of finite type, but not necessarily coherent.
\begin{dfn}\label{dfn-cohesive}{\rm 
A ringed space $X=(X,\O_X)$ is said to be {\em cohesive}\index{ringed space!cohesive ringed space@cohesive --- ---}\index{cohesive!cohesive ringed space@--- (ringed space)} if the structure sheaf $\O_X$ is coherent as an $\O_X$-module.}
\end{dfn}

Then one readily sees the following (cf.\ \ref{prop-cohringsmodules1}): 
\begin{prop}\label{prop-cohesive}
Let $X=(X,\O_X)$ be a cohesive ringed space\index{ringed space!cohesive ringed space@cohesive --- ---}\index{cohesive!cohesive ringed space@--- (ringed space)}.
Then an $\O_X$-module $\mathscr{F}$ is coherent\index{coherent!coherent sheaf@--- sheaf (on a ringed space)} if and only if it is finitely presented.
$($In this case, moreover, $\mathscr{F}$ admits an $\infty$-presentation.$)$ \hfill$\square$
\end{prop}

\begin{cor}\label{cor-cohesive}
Let $(f,\varphi)\colon(X,\O_X)\rightarrow(Y,\O_Y)$ be a morphism of ringed spaces, and suppose $(X,\O_X)$ is cohesive\index{ringed space!cohesive ringed space@cohesive --- ---}\index{cohesive!cohesive ringed space@--- (ringed space)}.
Then the functor $f^{\ast}$ maps $\Coh_Y$ to $\Coh_X$. \hfill$\square$
\end{cor}

\subsubsection{Filtered projective limit of ringed spaces}\label{subsub-projlimringedspexist}
We will need to consider filtered projective limits\index{limit!projective limit@projective ---} of (locally) ringed spaces.
In fact, as one can easily check, such a limit always exists in the category of (locally) ringed spaces:
\begin{prop}\label{prop-limLRS}
Let $\{X_i=(X_i,\O_{X_i})\}_{i\in I}$ be a projective system of ringed $($resp.\ locally ringed$)$ spaces indexed by a directed set\index{set!directed set@directed ---}\index{directed set} $I$.
Then the projective limit $X=\varprojlim_{i\in I}X_i$ exists in the category of ringed $($resp.\ locally ringed$)$ spaces.
Moreover, the underlying topological space of $X$ is isomorphic to the projective limit of underlying topological spaces of $X_i$'s, and $\O_{X,x}$ for any $x\in X$ is canonically isomorphic to the inductive limit $\varinjlim_{i\in I}\O_{X_i,p_i(x)}$, where $p_i\colon X\rightarrow X_i$ for each $i\in I$ is the canonical projection. \hfill$\square$
\end{prop}

That is to say, the topological space $X=\varprojlim_{i\in I}X_i$ coupled with the inductive limit sheaf $\O_X=\varinjlim_{i\in I}p^{-1}_i\O_{X_i}$, which is a sheaf of local rings if $\{X_i\}_{i\in I}$ is a projective system of locally ringed spaces (\ref{prop-directlimits111}), gives the desired projective limit.

\danger{Notice that, in the local-ringed case, it is due to presence of the filtering that the underlying topological space of the limit coincides with the limit of underlying topological spaces.
In fact, this is not in general the case for projective limits that are not filtered.
For instance, fiber products of schemes are taken in the category of locally ringed spaces (cf.\ \cite[(3.2.1)]{EGAInew}), and, as is well-known, their underlying topological spaces do not necessarily coincide with the fiber products of the underlying topological spaces.}

\begin{cor}\label{cor-lemlimLRS}
Let $\textrm{{\boldmath $X$}}=\{X_i=(X_i,\O_{X_i})\}_{i\in I}$ and $\textrm{{\boldmath $U$}}=\{U_i=(U_i,\O_{U_i})\}_{i\in I}$ be two projective systems of ringed $($resp.\ locally ringed$)$ spaces indexed by a directed set $I$, and $\iota\colon\textrm{{\boldmath $U$}}\rightarrow\textrm{{\boldmath $X$}}$ a morphism of projective systems $\iota=\{\iota_i\}$ consisting of open immersions $\iota_i\colon U_i\hookrightarrow X_i$.
Suppose we have $U_i\times_{X_i}X_j\stackrel{\sim}{\rightarrow}U_j$ for each pair $(i,j)$ with $i\leq j$.
Then the induced map
$$
\varprojlim_{i\in I}\iota_i\colon\varprojlim_{i\in I}U_i\longrightarrow\varprojlim_{i\in I}X_i
$$
is an open immersion, and the image coincides with $p^{-1}_i(U_i)$ for any $i$, where  $p_i\colon\varprojlim_{j\in I}X_j\rightarrow X_i$ is the projection map.\hfill$\square$
\end{cor}
\index{ringed space|)}

\subsection{Sheaves on limit spaces}\label{sub-projlimringedspsheaf}
\subsubsection{Finitely presented sheaves on limit spaces}\label{subsub-projlimringedspsheaf}
Let us consider
\begin{itemize}
\item a filtered projective system of ringed spaces $\{X_i=(X_i,\O_{X_i}),p_{ij}\}_{i\in I}$ indexed by a directed set $I$.
\end{itemize}
We set $X=\varprojlim_{i\in I}X_i$, and denote by $p_i\colon X\rightarrow X_i$ the canonical projection for each $i\in I$.
As we saw in \S\ref{subsub-projlimringedspexist}, the ringed space $X$ is supported on the projective limit of the underlying topological spaces of $X_i$'s and has $\O_X=\varinjlim_{i\in I}p^{-1}_i\O_{X_i}$ as its structure sheaf.
We assume:
\begin{itemize}
\item[{\rm (a)}] for any $i\in I$ the underlying topological space of $X_i$ is coherent\index{space@space (topological)!coherent topological space@coherent ---}\index{coherent!coherent topological space@--- (topological) space} {\rm (\ref{dfn-quasicompact1})} and sober\index{space@space (topological)!sober topological space@sober ---} {\rm (\S\ref{subsub-sober})}$;$
\item[{\rm (b)}] for any $i\leq j$ the underlying continuous mapping of the transition map $p_{ij}\colon X_j\rightarrow X_i$ is quasi-compact {\rm (\ref{dfn-quasicompactness} (2))}.
\end{itemize}
By \ref{thm-projlimcohsch1} (1) the underlying topological space of the limit $X$ is coherent and sober, and each canonical projection map $p_i\colon X\rightarrow X_i$ for $i\in I$ is quasi-compact.

\begin{thm}\label{thm-injlimmodpf}
{\rm (1)} For any finitely presented $\O_X$-module $\mathscr{F}$ there exist an index $i\in I$ and a finitely presented $\O_{X_i}$-module $\mathscr{F}_i$ such that $\mathscr{F}\cong p^{\ast}_i\mathscr{F}_i$.

{\rm (2)} For any morphism $\varphi\colon\mathscr{F}\rightarrow\mathscr{G}$ of finitely presented $\O_X$-modules, there exist an index $i\in I$ and a morphism $\varphi_i\colon\mathscr{F}_i\rightarrow\mathscr{G}_i$ of finitely presented $\O_{X_i}$-modules such that $\varphi\cong p^{\ast}_i\varphi_i$.
Moreover, if $\varphi$ is an isomorphism $($resp.\ epimorphism$)$, then one can take $\varphi_i$ to be an isomorphism $($resp.\ epimorphism$)$.
\end{thm}

The theorem is a consequence of the following theorem:
\begin{thm}\label{thm-projlimcohsch31}
Let $0\in I$ be an index, and $\mathscr{F}_0$ and $\mathscr{G}_0$ two $\O_{X_0}$-modules.
Suppose $\mathscr{F}_0$ is finitely presented.
Then the canonical map 
$$
\varinjlim_{i\geq 0}\Hom_{\O_{X_i}}(p^{\ast}_{0i}\mathscr{F}_0,p^{\ast}_{0i}\mathscr{G}_0)\longrightarrow\Hom_{\O_X}(p^{\ast}_0\mathscr{F}_0,p^{\ast}_0\mathscr{G}_0)
$$
is bijective.
\end{thm}

Notice that the last map is obtained by the inductive limit of the maps
$$
\Hom_{\O_{X_i}}(p^{\ast}_{0i}\mathscr{F}_0,p^{\ast}_{0i}\mathscr{G}_0)\longrightarrow\Hom_{\O_X}(p^{\ast}_0\mathscr{F}_0,p^{\ast}_0\mathscr{G}_0),
$$
which are simply induced by the canonical projections $p_i\colon X\rightarrow X_i$ for $i\geq 0$.

\begin{rem}\label{rem-injlimmodpf}{\rm 
In $2$-categorical language the theorems \ref{thm-injlimmodpf} and \ref{thm-projlimcohsch31} mean that the category of finitely presented $\O_X$-modules is equivalent to the `inductive limit category'
$$
\varinjLim\,\{\textrm{finitely presented $\O_{X_i}$-modules}\}
$$
in the sense of \cite[Expos\'e VI, (6.3)]{SGA4-2}.
Note that there exists an essentially unique fibered category over $I$ such that each fiber over $i\in I$ is the category of $\O_{X_i}$-modules of finite presentation and that Cartesian morphisms are given by pull-backs; see \cite[Expos\'e VI, \S 8]{SGA1} for the construction.}
\end{rem}

Now, let us show that the first theorem follows from the second:
\begin{proof}[Proof of Theorem {\rm \ref{thm-projlimcohsch31}} $\Rightarrow$ Theorem {\rm \ref{thm-injlimmodpf}}]
To show (1), we first consider the case where $\mathscr{F}$ has a finite presentation 
$$
\O^{\oplus p}_X\longrightarrow\O^{\oplus q}_X\longrightarrow\mathscr{F}\longrightarrow 0
$$
over $X$.
By \ref{thm-projlimcohsch31} the morphism $\O^{\oplus p}_X\rightarrow\O^{\oplus q}_X$ is the pull-back of a morphism $\O^{\oplus p}_{X_i}\rightarrow\O^{\oplus q}_{X_i}$ for some $i\geq 0$.
Let $\mathscr{F}_i$ be its cokernel. 
Then we have $\mathscr{F}\cong p^{\ast}_i\mathscr{F}_i$.

In general, we take a finite open covering $X=\bigcup_{k=1}^nU_k$ by quasi-compact open subsets such that $\mathscr{F}|_{U_k}$ for each $k=1,\ldots,n$ has a finite presentation.
By \ref{prop-projlimcohtopspqcptopen} there exists an index $i_k\in I$ for each $k=1,\ldots,n$ and a quasi-compact open subset $V_k\subseteq X_{i_k}$ such that $p^{-1}_{i_k}(V_k)=U_k$.
Taking an upper bound $i$ of $\{i_1,\ldots,i_n\}$ and replacing each $V_k$ by $p^{-1}_{i_ki}(V_k)$, we may assume that $V_k$ is a quasi-compact open subset of $X_i$ for any $k=1,\ldots,n$.
We have $\bigcup^n_{k=1}V_k\subseteq X_i$ and $X=\bigcup^n_{k=1}p^{-1}_i(V_k)$.
Hence, replacing $i$ by a larger index, we may assume that $X_i=\bigcup_{k=1}^nV_k$ (Exercise \ref{exer-coherentprojlimopencovering}).

By what we have seem above, replacing $i$ by a larger index if necessary, there exists a finitely presented $\O_{V_k}$-module $\mathscr{F}_k$ such that $p^{\ast}_i\mathscr{F}_k\cong\mathscr{F}|_{U_k}$ for any $k=1,\ldots,n$.
Since $p^{\ast}_i\mathscr{F}_k$'s patch together on $X$, there exists due to \ref{thm-projlimcohsch31} an index $j\geq i$ such that $p^{\ast}_{ij}\mathscr{F}_k$'s patch to a finitely presented $\O_{X_j}$-module $\mathscr{F}_j$.
Since we have $p^{\ast}_j\mathscr{F}_j\cong\mathscr{F}$, we have shown (1).

(2) follows from (1) and \ref{thm-projlimcohsch31}; to show that we can take $\varphi_i$ to be an epimorphism, we observe the following: If $\mathscr{H}_i$ is finitely presented $\O_X$-module such that $p^{\ast}_i\mathscr{H}_i=0$, then there exists an index $j\in I$ with $i\leq j$ such that $p^{\ast}_{ij}\mathscr{H}_i=0$. 
This follows from \ref{thm-projlimcohsch31}.
\end{proof}

We henceforth devote ourselves to the proof of \ref{thm-projlimcohsch31}.
\begin{lem}\label{lem-projlimcohsheaf1}
Let $Z=(Z,\O_Z)$ be a ringed space, and $\{g_i\colon X_i\rightarrow Z\}_{i\in I}$ a collection of morphisms of ringed spaces such that $g_j=g_i\circ p_{ij}$ whenever $i\leq j$.
Then for any $\O_Z$-module $\mathscr{G}$ we have a canonical isomorphism
$$
\varinjlim_{i\geq 0}p^{-1}_ig^{\ast}_i\mathscr{G}\stackrel{\sim}{\longrightarrow}g^{\ast}\mathscr{G}
$$
of $\O_X$-modules, where $g=\varprojlim_{i\in I}g_i$.
\end{lem}

Notice that the left-hand side of the last isomorphism can be regarded as an $\O_X$-module via the equality $\O_X=\varinjlim_{i\geq 0}p^{-1}_i\O_{X_i}$.

\begin{proof}
There exists a canonical map $p^{-1}_ig^{\ast}_i\mathscr{G}\rightarrow p^{\ast}_ig^{\ast}_i\mathscr{G}=g^{\ast}\mathscr{G}$, by which one has the morphism $\varinjlim_{i\geq 0}p^{-1}_ig^{\ast}_i\mathscr{G}\rightarrow g^{\ast}\mathscr{G}$.
In order to show that this is an isomorphism, we only need to check stalkwise.
Let $x\in X$.
We have 
$$
(\varinjlim_{i\geq 0}p^{-1}_ig^{\ast}_i\mathscr{G})_x=\varinjlim_{i\geq 0}(g^{\ast}_i\mathscr{G})_{p_i(x)}=\varinjlim_{i\geq 0}\mathscr{G}_{g(x)}\otimes_{\O_{Y,g(x)}}\O_{X_i,p_i(x)}.
$$
By \ref{prop-directlimits3} we have
\begin{equation*}
\begin{split}
\varinjlim_{i\geq 0}\mathscr{G}_{g(x)}\otimes_{\O_{Y,g(x)}}\O_{X_i,p_i(x)}&=\mathscr{G}_{g(x)}\otimes_{\O_{Y,g(x)}}\varinjlim_{i\geq 0}\O_{X_i,p_i(x)}\\
&=\mathscr{G}_{g(x)}\otimes_{\O_{Y,g(x)}}\O_{X,x},
\end{split}
\end{equation*}
whence the result.
(Notice that here we do not use the conditions (a) and (b) in the beginning of this paragraph.)
\end{proof}

\begin{lem}\label{lem-projlimcohsch31}
Suppose in the situation as in {\rm \ref{lem-projlimcohsheaf1}} that the underlying topological space of $Z$ is coherent and sober and that the underlying continuous mapping of each $g_i\colon X_i\rightarrow Z$ is quasi-compact.
Then the canonical morphism 
$$
\varinjlim_{i\in I}g_{i\ast}g^{\ast}_i\mathscr{G}\longrightarrow g_{\ast}g^{\ast}\mathscr{G}
$$
is an isomorphism of $\O_X$-modules.
\end{lem}

\begin{proof}
By \ref{lem-projlimcohsheaf1} one can apply \ref{cor-coherentprojectivesystemindlimits2} to the situation where $Y_i=Z$ for all $i\in I$ and $\mathscr{F}_i=f^{\ast}_i\mathscr{G}$ (hence $\mathscr{F}=f^{\ast}\mathscr{G}$).
(Notice that here we need to use the conditions (a) and (b).)
\end{proof}

\begin{proof}[Proof of Theorem {\rm \ref{thm-projlimcohsch31}}]
We have the equalities (up to canonical isomorphisms) 
$$
\Hom_{\O_X}(p^{\ast}_0\mathscr{F}_0,p^{\ast}_0\mathscr{G}_0)=\Hom_{\O_{X_0}}(\mathscr{F}_0,p_{0\ast}p^{\ast}_0\mathscr{G}_0)=\Hom_{\O_{X_0}}(\mathscr{F}_0,\varinjlim_{i\geq 0}p_{0i\ast}p^{\ast}_{0i}\mathscr{G}_0),
$$
where the last equality is due to \ref{lem-projlimcohsch31}.
Now by \ref{prop-projlimcohsch22} (2) we have
\begin{equation*}
\begin{split}
\Hom_{\O_{X_0}}(\mathscr{F}_0,\varinjlim_{i\geq 0}p_{0i\ast}p^{\ast}_{0i}\mathscr{G}_0)&=\varinjlim_{i\geq 0}\Hom_{\O_{X_0}}(\mathscr{F}_0,p_{0i\ast}p^{\ast}_{0i}\mathscr{G}_0)\\
&=\varinjlim_{i\geq 0}\Hom_{\O_{X_i}}(p^{\ast}_{0i}\mathscr{F}_0,p^{\ast}_{0i}\mathscr{G}_0),
\end{split}
\end{equation*}
as desired.
\end{proof}

\subsubsection{Limits and direct images}\label{subsub-limitsanddirectimages}
Next, in addition to the data fixed in the beginning of previous paragraph, we fix
\begin{itemize}
\item another filtered projective system of ringed spaces $\{Y_i,q_{ij}\}_{i\in I}$ indexed by the same directed set $I$ that satisfies the conditions similar to (a) and (b) as in the beginning of \S\ref{subsub-projlimringedspsheaf};
\item a map $\{f_i\}_{i\in I}$ of projective systems from $\{X_i,p_{ij}\}_{i\in I}$ to $\{Y_i,q_{ij}\}_{i\in I}$, that is, a collection of morphisms $f_i\colon X_i\rightarrow Y_i$ such that $q_{ij}\circ f_j=f_i\circ p_{ij}$ whenever $i\leq j$.
\end{itemize}
We set $Y=\varprojlim_{i\in I}Y_i$, and denote the canonical projection by $q_i\colon Y\rightarrow Y_i$ for each $i\in I$.
We have by passage to the projective limits the continuous map
$$
f=\varprojlim_{i\in I}f_i\colon X\longrightarrow Y.
$$
We moreover assume that 
\begin{itemize}
\item[{\rm (c)}] for any $i\in I$ the underlying continuous mapping of $f_i$ is quasi-compact.
\end{itemize}
Note that, then, the underlying continuous map of $f$ is quasi-compact due to \ref{thm-projlimcohspacepres} (1).

\begin{prop}\label{prop-projlimcohsch52}
For any $\O_X$-module $\mathscr{F}$ the canonical morphism 
$$
\varinjlim_{i\in I}q_i^{\ast}f_{i\ast}p_{i\ast}\mathscr{F}\longrightarrow f_{\ast}\mathscr{F}
$$
induced from the canonical morphisms $q_i^{\ast}f_{i\ast}p_{i\ast}\mathscr{F}\rightarrow f_{\ast}\mathscr{F}$ for $i\in I$, defined by adjunction from $f_{i\ast}p_{i\ast}\cong q_{i\ast}f_{\ast}$, is an isomorphism. \hfill$\square$
\end{prop}

The proposition follows from \ref{cor-coherentprojectivesystemindlimits3} and the following lemma, which can be viewed as a corollary of \ref{lem-projlimcohsheaf1}:
\begin{lem}\label{lem-corprojlimcohsheaf1}
Suppose we are given a system $\{\mathscr{F}_i,\varphi_{ij}\}$ consisting of 
\begin{itemize}
\item[$\mathrm{(a)}$] an $\O_{X_i}$-module $\mathscr{F}_i$ for each $i\in I$, 
\item[$\mathrm{(b)}$] a morphism $\varphi_{ji}\colon p^{\ast}_{ij}\mathscr{F}_i\rightarrow\mathscr{F}_j$ of $\O_{X_j}$-modules for each $j\geq i$, 
\end{itemize}
such that $\varphi_{kj}\circ p^{\ast}_{kj}\varphi_{ji}=\varphi_{ki}$ for $k\geq j\geq i$.
Then the canonical morphism 
$$
\varinjlim_{i\in I}p^{-1}_i\mathscr{F}_i\longrightarrow\varinjlim_{i\in I}p^{\ast}_i\mathscr{F}_i
$$
is an isomorphism of $\O_X$-modules.
\end{lem}

Here we regard $\{p^{-1}_i\mathscr{F}_i\}$ and $\{p^{\ast}_i\mathscr{F}_i\}$ as filtered inductive systems in an obvious manner.
Notice that we do not use in the following proof the conditions (a) and (b) in the beginning of \S\ref{subsub-projlimringedspsheaf} for $\{X_i,p_{ij}\}_{i\in I}$.
\begin{proof}
By \ref{lem-projlimcohsheaf1} one has
$$
\varinjlim_{i\in I}p^{\ast}_i\mathscr{F}_i=\varinjlim_{i\in I}\varinjlim_{j\geq i}p^{-1}_jp^{\ast}_{ji}\mathscr{F}_i
$$
up to canonical isomorphisms, where the last double inductive limit can be seen as a single inductive limit taken along the directed set $J=\{(i,j)\,|\, j\geq i\}$ where $(i,j)\leq (i',j')$ if and only if $i\leq i'$ and $j\leq j'$ (cf.\ Exercise \ref{exer-doublecolimits}).
Since the diagonal subset $\{(i,i)\}$ is evidently cofinal, we have 
$$
\varinjlim_{i\in I}\varinjlim_{j\geq i}p^{-1}_jp^{\ast}_{ji}\mathscr{F}_i=\varinjlim_{i\in I}p^{-1}_i\mathscr{F}_i,
$$
as desired.
\end{proof}

Applying \ref{prop-projlimcohsch52} to the situation $X_i=Y_i$ and $f_i=\id_{X_i}$ ($i\in I$), we have:
\begin{cor}\label{cor-projlimcohsch51}
For any $\O_X$-module $\mathscr{F}$ the canonical morphism 
$$
\varinjlim_{i\in I}p^{\ast}_ip_{i\ast}\mathscr{F}\longrightarrow\mathscr{F}
$$
is an isomorphism. \hfill$\square$
\end{cor}

\subsection{Cohomologies of sheaves on ringed spaces}\label{sub-cohomologiesringedspaces}
\subsubsection{Derived category formalism}\label{subsub-ringedspacederivedcatformalism}
Let $X=(X,\O_X)$ be a ringed space.
We denote by $\DC^{\ast}(X)$ (where $\ast=$``\ \ '', $+$, $-$, $\bd$) the derived category\index{category!derived category@derived ---}\index{derived category!derived category of OX modules@--- of $\O_X$-modules} associated to the abelian category $\Mod_X$ of $\O_X$-modules (cf.\ \S\ref{sub-derivedcategory}).
Inside $\Mod_X$ are thick\index{subcategory!thick subcategory@thick ---} (\S\ref{sub-subcategoryderived}) abelian full subcategories $\QCoh_X$ and $\Coh_X$ of quasi-coherent sheaves and coherent sheaves, respectively.
We denote by $\DC^{\ast}_{\qcoh}(X)$ (resp.\ $\DC^{\ast}_{\coh}(X)$) the full subcategory of $\DC^{\ast}(X)$ consisting of objects $F$ such that the cohomology sheaves $\mathcal{H}^k(F)$ are quasi-coherent (resp.\ coherent) for all $k\in\Z$ (where $\mathcal{H}^0$ is the canonical cohomology functor\index{functor!cohomology functor@cohomology ---!canonical cohomology functor@canonical --- ---} on $\DC^{\ast}(X)$ (cf.\ \ref{prop-homotopypre6})).
The full subcategories $\DC^{\ast}_{\qcoh}(X)$ and $\DC^{\ast}_{\coh}(X)$ are triangulated subcategories of $\DC^{\ast}(X)$ with the induced cohomology functor $\mathcal{H}^0$ and the induced $t$-structure\index{t-structure@$t$-structure} (cf.\ \ref{prop-subcategoryderived1}).

\subsubsection{Calculation of cohomologies}\label{subsub-ringedspacecalculationcoh}
In this and next paragraphs we make a few small, but perhaps at least need-to-know, remarks on cohomology groups of $\O_X$-modules.

Let $X=(X,\O_X)$ be a ringed space, and consider the commutative diagram of functors
$$
\xymatrix@C-5ex@R-1ex{\Mod_X\ar[rr]^u\ar[dr]_{\Gamma_X}&&\ASh_X\ar[dl]^{\Gamma_X}\\ &\Ab\rlap{,}}
$$
where $\Gamma_X$'s are the global section functors and $u\colon\Mod_X\rightarrow\Ab$ is the forgetful functor.
In this situation, one finds that there are at least two ways for obtaining the cohomology groups $\H^q(X,\mathscr{F})$ of an $\O_X$-module $\mathscr{F}$; one is by $\RD^q\Gamma_X(\mathscr{F})$, that is, by applying directly the right derived functors of $\Gamma_X\colon\Mod_X\rightarrow\Ab$, and the other by $\RD^q\Gamma_X(u(\mathscr{F}))$, calculated from the right derived functors of $\Gamma_X\colon\ASh_X\rightarrow\Ab$ applied to the underlying abelian sheaf of $\mathscr{F}$.
These two ways of calculating the cohomology groups, indeed, lead to the same result, but by a non-trivial reason.

The functor $u$, being clearly exact, induces the exact functor $\DC^+(u)\colon\DC^+(X)\rightarrow\DC^+(\ASh_X)$ (\ref{prop-derivedcategory71}), defined simply by term-by-term application of $u$ to complexes of $\O_X$-modules.
Then the issue lies in comparison of the exact functors $\RD^+\Gamma_X$ and $\RD^+\Gamma_X\circ\DC^+(u)$:
\begin{prop}\label{prop-ringedspacecalculationcoh1}
There exists a canonical isomorphism
$$
\RD^+\Gamma_X\cong\RD^+\Gamma_X\circ\DC^+(u)
$$
of exact functors $\DC^+(X)\rightarrow\DC^+(\Ab)$.
\end{prop}

\begin{proof}
The key-point is that the functor $u$ maps injective objects of $\Mod_X$, that is, injective $\O_X$-modules, to flasque (but not necessarily injective) abelian sheaves.
We have by \cite[C.D., \S2, Prop.\ 3.1]{SGA4.5} the canonical isomorphism 
$$
\RD^+(\Gamma_X\circ u)\cong\RD^+\Gamma_X\circ\RD^+u.
$$
Now, since $u$ is exact, the right derived functor $\RD^+u$ coincides with the induced functor $\DC^+(u)$, and thus we get the desired isomorphism of functors.
\end{proof}

A similar remark can be made also for higher direct images.
Let $f\colon X=(X,\O_X)\rightarrow Y=(Y,\O_Y)$ be a morphism of ringed spaces, and consider the commutative diagram
$$
\xymatrix{\Mod_X\ar[r]^u\ar[d]_{f_{\ast}}&\ASh_X\ar[d]^{f_{\ast}}\\ \Mod_Y\ar[r]_u&\ASh_Y\rlap{.}}
$$
Then by \cite[C.D., \S2, Prop.\ 3.1]{SGA4.5} and the fact that $f_{\ast}$ maps flasque sheaves to flasque sheaves, we have:
\begin{prop}\label{prop-ringedspacecalculationcoh2}
There exists a canonical isomorphism 
$$
\DC^+(u)\circ\RD^+f_{\ast}\cong\RD^+f_{\ast}\circ\DC^+(u)
$$
of exact functors $\DC^+(X)\rightarrow\DC^+(\ASh_Y)$. \hfill$\square$
\end{prop}

That is to say, the underlying abelian sheaves of the higher direct images $\RD^qf_{\ast}\mathscr{F}$ of an $\O_X$-module $\mathscr{F}$ taken in $\DC^+(Y)$ coincide up to isomorphisms with the higher direct images of $\mathscr{F}$ regarded as an abelian sheaf.

\subsubsection{Module structures on cohomologies}\label{subsub-ringedspacemodulestrcoh}
Let $X=(X,\O_X)$ be a ringed space.
Once given a ring homomorphism $A\rightarrow\Gamma_X(\O_X)$ from a ring $A$, one can put the canonical $A$-module structure to the global section $\Gamma_X(\mathscr{F})$ of an arbitrary $\O_X$-module $\mathscr{F}$, and thus obtain $\Gamma_X\colon\Mod_X\rightarrow\Mod_A$ and the associated right derived functor
$$
\RD^+\Gamma_X\colon\DC^+(X)\longrightarrow\DC^+(\Mod_A).\leqno{(\ast)}
$$

One has, on the other hand, the commutative diagram
$$
\xymatrix{\Mod_X\ar[d]_{\Gamma_X}\ar[dr]^{\Gamma_X}\\ \Mod_A\ar[r]_u&\Ab}
$$
consisting of global section functors and the forgetful functor; by \cite[C.D., \S2, Prop.\ 3.1]{SGA4.5} we have a canonical isomorphism
$$
\RD^+\Gamma_X\cong\DC^+(u)\circ\RD^+\Gamma_X
$$
of exact functors $\DC^+(X)\rightarrow\DC^+(\Ab)$.
This amounts to saying that the cohomology groups $\H^q(X,\mathscr{F})$ of an $\O_X$-module $\mathscr{F}$, defined originally as abelian groups, carry the canonical $A$-module structure induced from the above-fixed ring homomorphism $A\rightarrow\Gamma_X(\O_X)$.

Let $f\colon X=(X,\O_X)\rightarrow Y=(Y,\O_Y)$ be a morphism of ringed spaces, and suppose we are given a ring homomorphism $A\rightarrow\Gamma_Y(\O_Y)$ from a ring $A$; the homomorphism $\Gamma_Y(\O_Y)\rightarrow\Gamma_X(\O_X)$ induces a ring homomorphism $A\rightarrow\Gamma_X(\O_X)$.
As above, the global section functors yield $\Gamma_X\colon\Mod_X\rightarrow\Mod_A$ and $\Gamma_Y\colon\Mod_Y\rightarrow\Mod_A$, sitting in the following commutative diagram
$$
\xymatrix{\Mod_X\ar[d]_{f_{\ast}}\ar[dr]^{\Gamma_X}\\ \Mod_Y\ar[r]_{\Gamma_Y}&\Mod_A\rlap{.}}
$$
Just similarly to \ref{prop-ringedspacecalculationcoh2} one has:
\begin{prop}\label{prop-ringedspacecalculationcoh3}
There exists a canonical isomorphism
$$
\RD^+\Gamma_X\cong\RD^+\Gamma_Y\circ\RD^+f_{\ast}
$$ of exact functors $\DC^+(X)\rightarrow\DC^+(\Mod_A)$.\hfill$\square$
\end{prop}

In other words, the two ways of computing cohomology groups $\H^q(X,\mathscr{F})$ of an $\O_X$-module $\mathscr{F}$, together with the $A$-module structure, one being done directly on $X$, and the other via $Y$ by way of the higher direct images of $f$, lead to the same result.

\subsection{Cohomologies of module sheaves on limit spaces}\label{sub-cohomologiessheaveslimitspaces}
In this subsection we consider the data fixed in the beginning of \S\ref{subsub-projlimringedspsheaf} and \S\ref{subsub-limitsanddirectimages} and, furthermore, 
\begin{itemize}
\item for each $i\in I$ an $\O_{X_i}$-module $\mathscr{F}_i$;
\item for each $i\leq j$ a morphism $\varphi_{ij}\colon p^{\ast}_{ij}\mathscr{F}_i\rightarrow\mathscr{F}_j$ of $\O_{X_j}$-modules such that $\varphi_{ik}=\varphi_{jk}\circ p^{\ast}_{jk}\varphi_{ij}$ whenever $i\leq j\leq k$.
\end{itemize}
Then one has the inductive system $\{p^{\ast}_i\mathscr{F}_i\}_{i\in I}$ of $\O_X$-modules indexed by $I$ and the $\O_X$-module
$$
\mathscr{F}=\varinjlim_{i\in I}p^{\ast}_i\mathscr{F}_i.
$$
Moreover, we consider
\begin{itemize}
\item a ring $A$; 
\item a collection of ring homomorphisms $A\rightarrow\Gamma(X_i,\O_{X_i})$ for $i\in I$ such that the following diagram commutes whenever $i\leq j$:
$$
\xymatrix@R-4ex@C-2ex{&\Gamma(X_j,\O_{X_j})\\ A\ar[ur]\ar[dr]\\ &\Gamma(X_i,\O_{X_i})\rlap{.}\ar[uu]}
$$
\end{itemize}

By \ref{prop-cohcoherentprojectivesystemindlimits1}, \ref{lem-corprojlimcohsheaf1}, and what we have seen in \S\ref{subsub-ringedspacemodulestrcoh} (canonicity of $A$-module structure on the cohomologies) we have:
\begin{prop}\label{prop-ringedspacecohcoherentprojectivesystemindlimits1}
Then the canonical map
$$
\varinjlim_{i\in I}\H^q(X_i,\mathscr{F}_i)\longrightarrow\H^q(X,\mathscr{F})
$$
is an isomorphism of $A$-modules for $q\geq 0$. \hfill$\square$
\end{prop}

\begin{cor}\label{cor-ringedspacecohcoherentprojectivesystemindlimits11}
Let $\mathscr{H}$ be an $\O_X$-module.
Then the canonical morphism 
$$
\varinjlim_{i\in I}\H^q(X_i,p_{i\ast}\mathscr{H})\longrightarrow\H^q(X,\mathscr{H})
$$
is an isomorphism of $A$-modules for $q\geq 0$. \hfill$\square$
\end{cor}

\begin{cor}\label{cor-ringedspacecohcoherentprojectivesystemindlimits01}
Let $Z$ be a ringed space with the coherent\index{space@space (topological)!coherent topological space@coherent ---}\index{coherent!coherent topological space@--- (topological) space} underlying topological space, and $\{g_i\colon X_i\rightarrow Z\}_{i\in I}$ a system of morphisms such that $g_j=g_i\circ p_{ij}$ whenever $i\leq j$.
Suppose that the underlying continuous mapping of each $g_i$ is quasi-compact.
Then for any $\O_Z$-module $\mathscr{G}$ the canonical map
$$
\varinjlim_{i\in I}\H^q(X_i,g^{\ast}_i\mathscr{G})\longrightarrow\H^q(X,g^{\ast}\mathscr{G}),
$$
where $g=\varprojlim_{i\in I}g_i$, is an isomorphism of $A$-modules for $q\geq 0$. \hfill$\square$
\end{cor}

\begin{cor}\label{cor-ringedspacecohcoherentprojectivesystemindlimits2}
Let $\{X_i,p_{ij}\}_{i\in I}$, $\{Y_i,q_{ij}\}_{i\in I}$, and $f_i$'s be as in {\rm \S\ref{sub-projlimringedspsheaf}} $($with the conditions {\rm (a)} and {\rm (b)} for $\{X_i,p_{ij}\}_{i\in I}$ and $\{Y_i,q_{ij}\}_{i\in I}$ and the condition {\rm (c)} for $f_i$'s$)$, and $\mathscr{F}_i$ and $\mathscr{F}$ as above.
The canonical morphism of $\O_Y$-modules
$$
\varinjlim_{i\in I}q^{-1}_i\RD^qf_{i\ast}\mathscr{F}_i\longrightarrow\RD^qf_{\ast}\mathscr{F}
$$
is an isomorphism for $q\geq 0$. \hfill$\square$
\end{cor}

This follows from \ref{cor-cohcoherentprojectivesystemindlimits2} and \ref{lem-corprojlimcohsheaf1}.
One can, moreover, show the following results in a similar manner:
\begin{cor}\label{cor-ringedspacecohcoherentprojectivesystemindlimits21}
Let $Z$ be a ringed space with the coherent underlying topological space, and $\{g_i\colon X_i\rightarrow Z\}_{i\in I}$ a system of morphisms such that $g_j=g_i\circ p_{ij}$ whenever $i\leq j$.
Suppose that the underlying continuous mapping of each $g_i$ is quasi-compact.
Then for any $\O_Z$-module $\mathscr{G}$ the canonical morphism of $\O_Z$-modules
$$
\varinjlim_{i\in I}\RD^qg_{i\ast}(g^{\ast}_i\mathscr{G})\longrightarrow\RD^qg_{\ast}(g^{\ast}\mathscr{G}),
$$
where $g=\varprojlim_{i\in I}g_i$, is an isomorphism for $q\geq 0$. \hfill$\square$
\end{cor}

\begin{cor}\label{cor-ringedspacecohcoherentprojectivesystemindlimits3}
Let $\{X_i,p_{ij}\}_{i\in I}$, $\{Y_i,q_{ij}\}_{i\in I}$, and $f_i$'s be as in {\rm \ref{cor-ringedspacecohcoherentprojectivesystemindlimits2}}, and $\mathscr{H}$ an $\O_X$-module.
Then the canonical morphism of $\O_Y$-modules
$$
\varinjlim_{i\in I}q^{\ast}_i\RD^qf_{i\ast}(p_{i\ast}\mathscr{H})\longrightarrow\RD^qf_{\ast}\mathscr{H}
$$
is an isomorphism for $q\geq 0$. \hfill$\square$
\end{cor}

\addcontentsline{toc}{subsection}{Exercises}
\subsection*{Exercises}
\begin{exer}\label{exer-generizationmapsch}{\rm 
Let $X$ be a locally Noetherian formal scheme, $x\in X$ a point, and $y\in G_x$ a generization of $x$.
Then show that the map 
$$
(\O_{X,x})_{\mathfrak{q}_y}\longrightarrow\O_{X,y}
$$
induced from the generization map (\S\ref{subsub-generizationmap}), where $\mathfrak{q}_y$ is the pull-back of the maximal ideal $\m_{X,y}$ of $\O_{X,y}$, is faithfully flat\index{flatness!faithfully-flatness@faithfully-{---}}.}
\end{exer}

\begin{exer}\label{exer-limitflatnesspreserve}{\rm 
(1) Let $\{X_i,p_{ij}\}_{i\in I}$ be a filtered projective system of ringed spaces indexed by a directed set $I$, and $X=\varprojlim_{i\in I}X_i$.
Suppose that for any $i\leq j$ the transition map $p_{ji}\colon X_j\rightarrow X_i$ is flat.
Then show that for each $i\in I$ the canonical projection $p_i\colon X\rightarrow X_i$ is flat.

(2) Let $\{X_i,p_{ij}\}_{i\in I}$ and $\{Y_i,q_{ij}\}_{i\in I}$ be two filtered projective system of ringed spaces indexed by a directed set, and $\{f_i\}_{i\in I}$ a projective system of morphisms $f_i\colon X_i\rightarrow Y_i$ of ringed spaces.
Let $f\colon X\rightarrow Y$ be the limits, and $p_i\colon X\rightarrow X_i$ and $q_i\colon Y\rightarrow Y_i$ the canonical maps for each $i\in I$.
Suppose that for any $i\in I$ the map $f_i$ is flat.
Then show that the map $f$ is flat.}
\end{exer}

\begin{exer}\label{exer-limitcohesivepreserve}{\rm 
Let $\{X_i,p_{ij}\}_{i\in I}$ be a filtered projective system of ringed spaces indexed by a directed set $I$, and $X=\varprojlim_{i\in I}X_i$.
Suppose that the underlying topological space of each $X_i$ is coherent sober, and the underlying continuous map of each $p_{ij}$ is quasi-compact.
Suppose, moreover, that each $X_i=(X_i,\O_{X_i})$ is cohesive and that each $p_{ij}$ is flat.
Then show that $X=(X,\O_X)$ is cohesive.}
\end{exer}


\section{Schemes and algebraic spaces}\label{sec-schemes}
In \cite{EGA} coherent sheaves are always discussed under the assumption that the schemes in the context are locally Noetherian.
But it is true that many non-Noetherian schemes may have the coherent structure sheaves, to which, therefore, one can apply many of the results on locally Noetherian schemes.
In \S\ref{sub-schemes} and \S\ref{sub-algebraicspaces}, we will introduce the notion of {\em universally cohesive} schemes and algebraic spaces; a scheme is said to be universally cohesive if any scheme locally of finite presentation over it has the coherent structure sheaf.
Noetherian schemes are of course universally cohesive.
Some of the non-trivial examples of universally cohesive schemes will appear in \S\ref{subsub-coherencyadhesive}.
\danger{Unlike the case of Noetherian schemes, quasi-coherent sheaves of finite type on universally cohesive schemes are not necessarily coherent.}

In \S\ref{sub-formalismderived} we discuss some fundamental topics of cohomology calculus in the derived categories.
We discuss, for example, the comparison of two derived categories $\DC^{\ast}(\Coh_X)$, the derived category of the category of coherent sheaves, and $\DC^{\ast}_{\coh}(X)$, the full subcategory of the derived category of $\O_X$-modules consisting of objects having coherent cohomologies.

This section ends with a collection of known facts on cohomologies of quasi-coherent sheaves (\S\ref{sub-cohprerequi-0}) and on generalities of algebraic spaces (\S\ref{sub-morebasicsalgspaces}).

\subsection{Schemes}\label{sub-schemes}
\index{scheme|(}
\subsubsection{Schemes}\label{subsub-schemesbasics}
We denote by $\Sch$ the category of schemes, and by $\Sch_S$ (where $S$ is a scheme) the category of $S$-schemes.
Notice that $\Sch$ is a full subcategory of $\LRsp$, the category of locally ringed spaces.

For an affine scheme\index{affine!affine scheme@--- scheme} $X=\Spec A$ and an $A$-module $M$, we denote, as usual, by $\til{M}$ the associated quasi-coherent $\O_X$-module.\index{quasi-coherent!quasi-coherent OX module on schemes@--- sheaf (on a scheme)}
As is well-known (\cite[\S1.4]{EGAInew}), $M\mapsto\til{M}$ gives an exact equivalence of abelian categories
$$
\til{\,\cdot\,}\colon\Mod_A\stackrel{\sim}{\longrightarrow}\QCoh_X
$$
preserving the tensor products and internal hom's.

\subsubsection{Universally cohesive schemes}\label{subsub-schemesunivcoh}
\index{scheme!universally cohesive scheme@universally cohesive ---|(}
\begin{dfn}\label{dfn-universallycohesive}{\rm 
A scheme $X$ is said to be {\em universally cohesive}\index{cohesive!universally cohesive@universally --- (schemes)} if any $X$-scheme locally of finite presentation is cohesive (\ref{dfn-cohesive})\index{cohesive!cohesive ringed space@--- (ringed space)} as a ringed space.}
\end{dfn}

\begin{prop}\label{prop-cohschemes2}
Let $A$ be a ring, and $X=\Spec A$.
Then $X$ is universally cohesive if and only if $A$ is universally coherent {\rm (\ref{dfn-universallycoherent})}\index{coherent!coherent ring@--- ring!universally coherent ring@universally --- ---}.
\end{prop}

\begin{proof}
Suppose $X$ is universally cohesive, and let $B$ be a finitely presented $A$-algebra.
Suppose we are given an exact sequence 
$$
0\longrightarrow K\longrightarrow B^{\oplus m}\longrightarrow B.
$$
We need to show that $K$ is a finitely generated $B$-module.
By \cite[(1.3.11)]{EGAInew} this gives the exact sequence 
$$
0\longrightarrow \til{K}\longrightarrow\O_Y^{\oplus m}\longrightarrow\O_Y,
$$
where $Y=\Spec B$.
Since $\O_Y$ is coherent, $\til{K}$ is of finite type.
Then by \cite[(1.4.3)]{EGAInew} we deduce that $K$ is finitely generated.

Conversely, suppose $A$ is universally coherent, and let $Y\rightarrow X$ be an $X$-scheme locally of finite presentation.
We need to show that for any open subset $U\subseteq Y$ and any exact sequence of the form 
$$
0\longrightarrow\mathscr{K}\longrightarrow\O_U^{\oplus m}\longrightarrow\O_U,
$$
the quasi-coherent sheaf $\mathscr{K}$ on $U$ is of finite type.
To this end, we may assume that $U$ is affine $U=\Spec B$, where $B$ is an $A$-algebra of finite presentation and hence is universally coherent.
Set $\mathscr{K}=\til{K}$, where $K$ is a (uniquely determined up to isomorphism) $B$-module.
Then by \cite[(1.3.11)]{EGAInew} we see that $K$ is finitely generated, and by \cite[(1.4.3)]{EGAInew} we conclude that $\mathscr{K}$ is of finite type, as desired.
\end{proof}

\begin{prop}\label{prop-univcohesiveaffine1}
Consider an affine scheme\index{affine!affine scheme@--- scheme} $X=\Spec A$, and suppose $A$ is universally coherent.
Then the functor $\til{\,\cdot\,}$ induces an exact equivalence of abelian categories
$$
\til{\,\cdot\,}\colon\Coh_A\stackrel{\sim}{\longrightarrow}\Coh_X.
$$
\end{prop}

\begin{proof}
By \cite[(1.4.3)]{EGAInew} the functor $\til{\,\cdot\,}$ gives a categorical equivalence from the category of finitely presented $A$-modules to the category of finitely presented $\O_X$-modules, which are precisely coherent $\O_X$-modules due to \ref{prop-cohesive}.
\end{proof}
\index{scheme!universally cohesive scheme@universally cohesive ---|)}\index{scheme|)}

\subsection{Algebraic spaces}\label{sub-algebraicspaces}
\index{algebraic space|(}
\subsubsection{Conventions}\label{subsub-algebraicspacesconv}
Our basic reference to the generality of algebraic spaces is Knutson's\index{Knutson, D.} book \cite{Knu}.
Accordingly, we adopt the following convention:

\medskip\noindent
{\bf Convention.} {\sl In this book, all algebraic spaces are assumed to be quasi-separated.}

\medskip
This assumption will be particularly useful when we compare algebraic spaces with schemes, using several effective \'etale descent (\cite[II.3]{Knu}) and several local constructions such as open complement of closed subspaces (\cite[II.5]{Knu}).
Note that, with this assumption, the two definitions of algebraic spaces, the one as a sheaf on the large \'etale site of affine schemes and the other as the quotient of schemes by an \'etale equivalence relation, coincide (\cite[II.1]{Knu}).
In particular, if $X$ is a quasi-separated scheme, we have:

\medskip
$\bullet$ By \'etale descent of quasi-coherent sheaves, the category of quasi-coherent sheaves on $X$ regarded as an algebraic space (with respect to \'etale topology) and the category of quasi-coherent sheaves on $X$ (as usual with respect to Zariski topology) are canonically equivalent.

\medskip
$\bullet$ The cohomologies of quasi-coherent sheaves computed on the scheme $X$ by means of Zariski topology coincide up to canonical isomorphism with those computed from sheaves on $X$ regarded as an algebraic space with respect to \'etale topology\index{topology!etale topology@\'etale ---} (cf.\ \cite[Exp.\ VII, Prop.\ 4.3]{SGA4-2}).

\medskip
By these facts, as long as we are concerned with quasi-coherent sheaves and their cohomologies, we do not have to distinguish the two standpoints for a quasi-separated scheme $X$, the one on which $X$ is regarded as a scheme (Zariski topologized) and the other on which $X$ is regarded as an algebraic space (considered with \'etale topology).

\medskip\noindent
{\bf Convention.} {\sl When we say `$X$ is an algebraic space,' we always mean either {\rm (a)} $X$ is a scheme $($not necessarily quasi-separated$)$ or {\rm (b)} $X$ is an algebraic space $($necessarily quasi-separated$);$ when we like to exclude the non-quasi-separated schemes, we say `$X$ is an algebraic space or a quasi-separated scheme.'}

\medskip
Accordingly, in the case (a), all quasi-coherent sheaves on $X$ and their cohomologies are considered with respect to Zariski topology (unless otherwise clearly stated) and, in the case (b), they are considered with respect to \'etale topology\index{topology!etale topology@\'etale ---}.

In some situation where algebraic spaces and schemes are mixed, e.g.\ when taking fiber products of algebraic spaces with schemes, one has to (and we do) assume that the schemes are quasi-separated, even if not explicitly stated.

In the sequel we denote by $\As$ the category of algebraic spaces and by $\As_S$ (where $S$ is an algebraic space) the category of $S$-algebraic spaces.

\subsubsection{Basic notions}\label{subsub-algebraicspacesbasics}
A morphism $f\colon X\rightarrow Y$ of algebraic spaces is said to be {\em coherent}\index{coherent!coherent morphism of algebraic spaces@--- morphism (of algebraic spaces)} if it is quasi-compact and quasi-separated.
By saying that an algebraic space $X$ is {\em coherent}\index{coherent!coherent algebraic space@--- algebraic space}\index{algebraic space!coherent algebraic space@coherent ---} we always understand that $X$ is coherent over $\Spec\Z$.

Let $X$ be an algebraic space, and $\mathscr{B}$ a quasi-coherent $\O_X$-algebra sheaf.
Then as in \cite[II.5]{Knu} we have the algebraic space $\Spec\mathscr{B}$ affine over $X$; if $X=\Spec A$ and $\Gamma(X,\mathscr{B})=B$, then it is the associated algebraic space to the affine scheme $\Spec B$.
The existence of $\Spec\mathscr{B}$ in general follows from effective \'etale descent of affine maps (cf.\ \cite[Expos\'e VIII, 2]{SGA1}) and the well-known fact on local constructions (\ref{prop-localconstruction}).

Similarly, for a quasi-coherent $\O_X$-module $\mathscr{E}$ of finite type, one has the algebraic space $\P(\mathscr{E})$ projective over $X$ together with the invertible sheaf $\O_{\P(\mathscr{E})}(1)$.
The construction in the case where $X$ is scheme is due to \cite[{$\mathbf{II}$}, (4.1.1)]{EGA}; the general case follows from \ref{prop-localconstruction} and effective \'etale descent of schemes with relatively ample invertible sheaves (cf.\ \cite[Expos\'e VIII, 7.8]{SGA1}).

\subsubsection{Universally cohesive algebraic spaces}\label{subsub-algebraicspacesunivcoh}
\index{algebraic space!universally cohesive algebraic space@universally cohesive ---|(}
\begin{dfn}\label{dfn-universallycohesive2}{\rm 
An algebraic space $X$ is said to be {\em universally cohesive} if for any algebraic space $Y$ locally of finite presentation over $X$, $\O_Y$ is a coherent $\O_Y$-module.}
\end{dfn}

This definition is consistent with \ref{dfn-universallycohesive} above, because of the following `fppf descent of cohesiveness'.
\begin{prop}\label{prop-flatdescentcoh}
Let $f\colon Y\rightarrow X$ be a faithfully flat and finitely presented morphism of schemes.
Then $X$ is universally cohesive if and only if so is $Y$.
\end{prop}

\begin{proof}
The `only if' part is obvious.
To show the `if' part, it suffices to show that $X$ is cohesive if so is $Y$, which follows easily from \cite[Expos\'e VIII, 1.10]{SGA1}.
\end{proof}
\index{algebraic space!universally cohesive algebraic space@universally cohesive ---|)}
\index{algebraic space|)}

\subsection{Derived category calculus}\label{sub-formalismderived}
\subsubsection{Quasi-coherent sheaves on affine schemes}\label{subsub-formalismderived}
The following proposition is easy to verify (cf.\ \ref{prop-derivedcategory71}):
\begin{prop}\label{prop-formalismderived1}
Let $A$ be a ring, and set $X=\Spec A$.

{\rm (1)} The exact equivalence $\Mod_A\stackrel{\sim}{\rightarrow}\QCoh_X$ by $M\mapsto\til{M}$ induces an exact $($cf.\ {\rm \S\ref{sub-triangulatedcategory}}$)$ equivalence
$$
\DC^{\ast}(\Mod_A)\stackrel{\sim}{\longrightarrow}\DC^{\ast}(\QCoh_X)
$$
of triangulated categories.

{\rm (2)} If $A$ is universally coherent, then the exact equivalence $\Coh_A\stackrel{\sim}{\rightarrow}\Coh_X$ by $M\mapsto\til{M}$ induces an exact equivalence
$$
\DC^{\ast}(\Coh_A)\stackrel{\sim}{\longrightarrow}\DC^{\ast}(\Coh_X)
$$
of triangulated categories. \hfill$\square$
\end{prop}

Let us denote the composite functor
$$
\DC^{\ast}(\Mod_A)\stackrel{\sim}{\longrightarrow}\DC^{\ast}(\QCoh_X)\stackrel{\delta^{\ast}}{\longrightarrow}\DC^{\ast}_{\qcoh}(X)
$$
(cf.\ \S\ref{sub-subcategoryderived}) by 
$$
M\longmapsto M_X.
$$

\begin{prop}\label{prop-formalismderived2}
Let $A$ be a universally coherent ring, and set $X=\Spec A$.
Let $\delta^{\bd}$ be the canonical exact functor
$$
\delta^{\bd}\colon\DC^{\bd}(\Coh_X)\longrightarrow\DC^{\bd}_{\coh}(X).
$$
Then the functor $\delta^{\bd}$ is a categorical equivalence.
\end{prop}

\begin{proof}
By \ref{prop-subcategoryderived4} it suffices to show that $\delta^{\bd}$ is fully faithful.
Since we have $\Hom_{\DC^{\bd}(X)}(K,L)=\mathcal{H}^0(\RD\Hom_{\O_X}(K,L))$ etc., it suffices to show the following: for any objects $K,L$ of $\DC^{\bd}(\Coh_X)$ we have
$$
\RD\Hom_{\Coh_X}(K,L)=\RD\Hom_{\O_X}(\delta^{\bd}(K),\delta^{\bd}(L))\leqno{(\ast)}
$$
up to isomorphism in $\DC^{+}(\Ab)$.
By \ref{prop-formalismderived1} there exists an object $M$ of $\DC^{\bd}(\Coh_A)$ such that $M_X=K$ (up to isomorphism).
By \ref{cor-cohringsmodules11} and \cite[Chap.\ III, (2.4.1) (b)]{Verd1} we have a finite free resolution $F\in\obj(\DC^-(\Coh_A))$ of $M$.
Then, by \ref{thm-vanishcohaff-0} (1), the both sides of $(\ast)$ is isomorphic to $\Hom_{\O_X}(F_X,L)$ in $\DC^{+}(\Ab)$, whence the claim.
\end{proof}

\begin{cor}\label{cor-formalismderived21}
Let $A$ be a universally coherent ring, and $X=\Spec A$.
Then the functor 
$$
\DC^{\bd}(\Coh_A)\longrightarrow\DC^{\bd}_{\coh}(X),\qquad M\longmapsto M_X,
$$
is an exact equivalence of the triangulated category. \hfill$\square$
\end{cor}

In other words, any object $M$ of $\DC^{\bd}(X)$ whose cohomologies are all coherent can be represented (in the sense as in \ref{dfn-derivedcategory5}) by a complex consisting of coherent sheaves and hence by a complex consisting of finitely presented $A$-modules.

Finally, we include here the related result quoted from \cite{SGA6}:
\begin{prop}[{\cite[Expos\'e II, Cor.\ 2.2.2.1]{SGA6}}]\label{prop-formalderived3}
Let $X$ be a Noetherian scheme. 
Then the canonical exact functor
$$
\delta^{\bd}\colon\DC^{\bd}(\Coh_X)\longrightarrow\DC^{\bd}_{\coh}(X)
$$
is an equivalence of triangulated categories. \hfill$\square$
\end{prop}

\subsubsection{Permanence of coherency}\label{sub-calculusderived}
\begin{prop}\label{prop-schpairadh3}
Let $X$ be a universally cohesive algebraic space.

{\rm (1)} For $F,G\in\obj(\DC^-_{\coh}(X))$, $F\otimes^{\LD}_{\O_X}G$ belongs to $\DC^-_{\coh}(X)$.

{\rm (2)} For $F\in\obj(\DC^-_{\coh}(X))$ and $G\in\obj(\DC^+_{\coh}(X))$, $\RD\lHom_{\O_X}(F,G)$ belongs to $\DC^+_{\coh}(X)$.
\end{prop}

\begin{proof}
Replacing $X$ by an \'etale or Zariski covering, one can reduce to the case where $X$ is affine $X=\Spec A$. 
By \ref{prop-cohschemes2} the ring $A$ is universally coherent.
(In the following, we may work with either \'etale topology or Zariski topology without essential difference.)

(1) By suitable shifts we may assume that $\mathcal{H}^k(F)=\mathcal{H}^k(G)=0$ for $k>0$.
Let $n$ be a negative integer, and consider the natural morphism
$$
f^n\colon F\otimes^{\LD}_{\O_X}G\longrightarrow\tau^{\geq n}F\otimes^{\LD}_{\O_X}\tau^{\geq n}G.
$$
Then $\mathcal{H}^k(f^n)$ is an isomorphism for $k\geq n$.
Hence, to show that $F\otimes^{\LD}_{\O_X}G$ has coherent cohomologies, we may assume that $F$ and $G$ belong to $\DC^{\bd}_{\coh}(X)$.
By \ref{cor-formalismderived21} $F$ and $G$ are represented by bounded complexes of finitely presented (that is, coherent) $A$-modules $M^{\bullet}$ and $N^{\bullet}$, respectively.
Then the assertion follows from \ref{cor-cohringsmodules12} and an easy homological algebra.

(2) Similarly, we may assume that $F$ and $G$ are in $\DC^{\bd}_{\coh}(X)$ so that they are represented by bounded complexes of coherent $A$-modules.
Then apply \ref{cor-cohringsmodules12} in a similar way.
\end{proof}

\begin{prop}\label{prop-schpairadh31}
Let $f\colon X\rightarrow Y$ be a morphism of universally cohesive algebraic spaces.
Then the functor $\LD f^{\ast}$ maps $\DC^-_{\coh}(Y)$ to $\DC^-_{\coh}(X)$.
\end{prop}

\begin{proof}
Recall that the functor $\LD f^{\ast}\colon\DC^-(Y)\rightarrow\DC^-(X)$ is given by the following composition:
$$
\DC^-(Y)\longrightarrow\DC^-(\Mod_{(X,f^{-1}\O_Y)})\longrightarrow\DC^-(X)
$$
mapping $M\mapsto f^{-1}M\mapsto f^{-1}M\otimes^{\LD}_{f^{-1}\O_Y}\O_X$, where the first functor is the one obtained from the exact functor $f^{-1}\colon\Mod_Y\rightarrow\Mod_{(X,f^{-1}\O_Y)}$.
Since the first functor preserves the canonical $t$-structures, we have
$$
f^{-1}\tau^{\geq n}M=\tau^{\geq n}f^{-1}M
$$
for $n\in\Z$.
By this and the similar reasoning as in the proof of \ref{prop-schpairadh3}, we may assume that $M$ lies in $\DC^{\bd}_{\coh}(Y)$.
We may also assume, as in the proof of \ref{prop-schpairadh3}, that $Y$ is affine.
In this situation, $M$ has a finite free resolution, and hence $f^{-1}M$ has a free $f^{-1}\O_Y$-resolution.
Hence, by a standard homological algebra we deduce that $\LD f^{\ast}M$ is coherent.
\end{proof}

\subsection{Cohomology of quasi-coherent sheaves}\label{sub-cohprerequi-0}
\subsubsection{Cohomologies on affine schemes}\label{subsub-cohprerequi-0-affine}
\index{affine!affine scheme@--- scheme!cohomology of affine schemes@cohomology of --- ---s|(}
We first include the following well-known facts:
\begin{thm}\label{thm-qcohaff-0}
Let $A$ be a ring, and set $X=\Spec A$.

{\rm (1)} {\rm (\cite[$\mathbf{I}$, (1.4.1)]{EGA})} The following conditions for an $\O_X$-module $\mathscr{F}$ are equivalent to each other$:$
\begin{itemize}
\item[{\rm (a)}] $\mathscr{F}$ is quasi-coherent\index{quasi-coherent!quasi-coherent OX module on schemes@--- sheaf (on a scheme)}$;$ 
\item[{\rm (b)}] there exists an $A$-module $M$ such that $\mathscr{F}\cong\til{M};$ 
\item[{\rm (c)}] there exist a finite open covering $X=\bigcup_iU_i$ by open sets of the form $U_i=D(f_i)$ with $f_i\in A$ and for each $i$ an $A_{f_i}$-module $M_i$ such that $\mathscr{F}|_{U_i}\cong\til{M}_i$.
\end{itemize}

{\rm (2)} {\rm (\cite[$\mathbf{I}$, (1.3.7)]{EGA})} For any $A$-module $M$ we have $\Gamma(X,\til{M})\cong M$.

{\rm (3)} {\rm (\cite[(1.4.3)]{EGAInew})} For any $A$-module $M$ the quasi-coherent $\O_X$-module $\til{M}$ is of finite type $($resp.\ of finite presentation$)$ if and only if $M$ is finitely generated $($resp.\ finitely presented$)$ over $A$. \hfill$\square$
\end{thm}

\begin{thm}[{\cite[$\mathbf{III}$, (1.3.1), (1.3.2)]{EGA}, \cite[II.4.8]{Knu}}]\label{thm-vanishcohaff-0}
{\rm (1)} Let $A$ be a ring, and $\mathscr{F}$ a quasi-coherent sheaf on $X=\Spec A$.
Then for $q>0$ we have:
$$
\H^q(X,\mathscr{F})=0.
$$

{\rm (2)} Let $f\colon X\rightarrow Y$ be an affine morphism between algebraic spaces, and $\mathscr{F}$ a quasi-coherent sheaf on $X$.
Then for $q>0$ we have:
$$
\RD^qf_{\ast}\mathscr{F}=0.\eqno{\square}
$$
\end{thm}

Note that, as mentioned in \S\ref{subsub-algebraicspacesconv}, the vanishing as in (1) is also true for the cohomology calculated in terms of \'etale topology\index{topology!etale topology@\'etale ---}.

Let $X$ be an algebraic space, and $\mathscr{F}^{\bullet}$ a complex of $\O_X$-modules such that $\mathscr{F}^q=0$ for $q\ll 0$.
Then $F=Q^+(\mathrm{h}^+(\mathscr{F}^{\bullet}))$ is clearly an object of $\DC^+(X)$ (see \S\ref{subsub-homopotypre} and \S\ref{subsub-derivedcategorydef} for the notation).
We write 
$$
\RD^+\Gamma_X(\mathscr{F}^{\bullet})=\RD^+\Gamma_X(F),
$$
which is an object of $\DC^+(\Ab)$.
If $X=\Spec A$, then $\RD^+\Gamma_X(\mathscr{F}^{\bullet})$ is canonically regarded as an object of $\DC^+(\Mod_A)$ (\S\ref{subsub-ringedspacemodulestrcoh}).
If $f\colon X\rightarrow Y$ is a morphism of algebraic spaces, then we write
$$
\RD^+f_{\ast}\mathscr{F}^{\bullet}=\RD^+f_{\ast}F,
$$
which is an object of $\DC^+(Y)$.

\begin{thm}\label{thm-vanishcohaffder-0}
{\rm (1)} Let $A$ be a ring, and $\mathscr{F}^{\bullet}$ a complex of quasi-coherent sheaves on $X=\Spec A$ such that $\mathscr{F}^q=0$ for $q\ll 0$ $($resp.\ $|q|\gg 0)$.
Then $\H^q(X,\mathscr{F}^{\bullet})=0$ for $q\ll 0$ $($resp.\ $|q|\gg 0)$, and the object $\RD^+\Gamma_X(\mathscr{F}^{\bullet})$ of $\DC^+(\Mod_A)$ is represented {\rm (\ref{dfn-derivedcategory5})} by the complex $\Gamma_X(\mathscr{F}^{\bullet})$.

{\rm (2)} Let $f\colon X\rightarrow Y$ be an affine morphism between algebraic spaces, and $\mathscr{F}^{\bullet}$ a complex of quasi-coherent sheaves on $X$ such that $\mathscr{F}^q=0$ for $q\ll 0$ $($resp.\ $|q|\gg 0)$.
Then $\RD^q f_{\ast}\mathscr{F}^{\bullet}=0$ for $q\ll 0$ $($resp.\ $|q|\gg 0)$, and the object $\RD^+ f_{\ast}\mathscr{F}^{\bullet}$ of $\DC^+(Y)$ is represented by the complex $f_{\ast}\mathscr{F}^{\bullet}$.
\end{thm}

\begin{proof}
We only present the proof of (1), for (2) can be shown similarly.
First we deal with the bounded case.
By a suitable shift we may assume that the complex $\mathscr{F}^{\bullet}$ is of the form
$$
\mathscr{F}^{\bullet}=(\cdots\rightarrow 0\rightarrow\mathscr{F}^0\rightarrow\mathscr{F}^1\rightarrow\cdots\rightarrow\mathscr{F}^{l-1}\rightarrow0\rightarrow\cdots).
$$
Consider the distinguished triangle 
$$
\sigma^{\geq 1}\mathscr{F}^{\bullet}\longrightarrow\mathscr{F}^{\bullet}\longrightarrow\mathscr{F}^0\stackrel{+1}{\longrightarrow}
$$
in $\KC^{\bd}(\QCoh_X)$, where $\sigma^{\geq n}$ denotes the stupid truncation (\S\ref{subsub-complexcategorytruncations}).
Taking the cohomology exact sequence, we find by induction with respect to the length of $\mathscr{F}^{\bullet}$ that we may assume $\mathscr{F}$ is a single sheaf; but then the theorem in this case is nothing but \ref{thm-vanishcohaff-0}.

The general case can be reduced to the bounded case similarly by the stupid truncation; if $\mathscr{F}^{\bullet}=(\cdots\rightarrow 0\rightarrow\mathscr{F}^0\rightarrow\mathscr{F}^1\rightarrow\cdots)$, then to detect $\H^q(X,\mathscr{F}^{\bullet})$ for a fixed $q$, we may replace $\mathscr{F}^{\bullet}$ by $\sigma^{\leq q+1}\mathscr{F}^{\bullet}$.
\end{proof}
\index{affine!affine scheme@--- scheme!cohomology of affine schemes@cohomology of --- ---s|)}

\subsubsection{Some finiteness results}\label{subsub-cohprerequi-0-finiteness}
The following is a corollary of \ref{thm-vanishcohaffder-0} (2):
\begin{cor}\label{cor-calculusderived2-0}
Let $X$ be a universally cohesive algebraic space, and $i\colon Y\hookrightarrow X$ a closed immersion of finite presentation $($hence $Y$ is again universally cohesive$)$.
Then $\RD^+i_{\ast}$ maps $\DC^{\bd}_{\coh}(Y)$ to $\DC^{\bd}_{\coh}(X)$.
\end{cor}

The proof is easy: First reduce to the affine situation $Y=\Spec B\hookrightarrow X=\Spec A$, where $A\rightarrow B$ is a surjective homomorphism with the finitely generated kernel; then apply \ref{thm-vanishcohaffder-0} (2) with the aid of the following easy lemma:
\begin{lem}\label{lem-pf2}
Let $A\rightarrow B$ be a surjective ring homomorphism, and $M$ a $B$-module.

{\rm (1)} If $M$ is finitely presented as an $A$-module, then so is as a $B$-module.

{\rm (2)} Suppose that the kernel of the map $A\rightarrow B$ is a finitely generated ideal of $A$.
Then, if $M$ is finitely presented over $B$, it is finitely presented over $A$. \hfill$\square$
\end{lem}

\begin{prop}\label{prop-cohqcoh-0}
Let $f\colon X\rightarrow Y$ be a coherent morphism\index{coherent!coherent morphism of algebraic spaces@--- morphism (of algebraic spaces)} between algebraic spaces {\rm (\S\ref{subsub-algebraicspacesbasics})}.
Then $\RD^+ f_{\ast}$ maps an object of $\DC^+_{\qcoh}(X)$  to an object of $\DC^{+}_{\qcoh}(Y)$.
\end{prop}

\begin{proof}
Here we present the proof for the case where $X$ and $Y$ are schemes.
(This is actually enough, due to the standard technique (cohomological \'etale descent; cf.\ \cite[Expos\'e $\textrm{V}^{\mathrm{bis}}$, (4.2.1)]{SGA4-2}).)
Let $F$ be an object of $\DC^+_{\qcoh}(X)$.
We need to show that $\RD^qf_{\ast}F$ for any $q\geq 0$ is a quasi-coherent sheaf on $Y$.
But, to show it for a fixed $q$, we may replace $F$ by $\tau^{\leq q+1}F$, and thus we may assume that $F$ belongs to $\DC^{\bd}_{\qcoh}(X)$.
By shift we may assume $\mathcal{H}^q(F)=0$ for $q<0$.
Then by the distinguished triangle
$$
F\longrightarrow\tau^{\geq 1}F\longrightarrow G\stackrel{+1}{\longrightarrow},
$$
we may assume by induction that $F$ is concentrated in degree $0$ (cf.\ \ref{dfn-derivedcategory6} (2)) and hence that $F$ is represented by a single quasi-coherent sheaf $\mathscr{F}$ on $X$ (\ref{prop-derivedcategory61}).
What to prove is that $\RD^qf_{\ast}\mathscr{F}$ is quasi-coherent for every $q$, which is shown in \cite[$\mathbf{III}$, (1.4.10)]{EGA} (cf.\ \cite[$\mathbf{IV}$, (1.7.21)]{EGA}).
\end{proof}

\begin{cor}[cf.\ {\cite[$\mathbf{III}$, (1.4.11)]{EGA}}]\label{cor-cohqcoh1-0}
In the situation as in $\ref{prop-cohqcoh-0}$, let $i\colon V\rightarrow Y$ be an \'etale morphism from an affine scheme $V$.
Denote by $j\colon U=X\times_YV\rightarrow X$ the induced morphism.
Let $\mathscr{F}^{\bullet}$ be a complex of quasi-coherent $\O_X$-modules such that $\mathscr{F}^q=0$ for $q\ll 0$.
Then for any $q$ we have:
$$
\H^q(U,j^{\ast}\mathscr{F}^{\bullet})\cong\Gamma(V,i^{\ast}\RD^qf_{\ast}\mathscr{F}^{\bullet}).\eqno{\square}
$$
\end{cor}

\subsubsection{Cohomologies on projective spaces}\label{subsub-cohprerequi-0-projective}
We will also need some results on (cohomologies of) quasi-coherent sheaves\index{quasi-coherent!quasi-coherent OX module on schemes@--- sheaf (on a scheme)} on projective spaces. 
We fix the following notation: Let $A$ be a ring, $r$ a positive integer, and $S=A[T_0,\ldots,T_r]$ the polynomial ring over $A$.
Let $X=\P^r_A=\Proj S$. 
For $n>0$ we set $\textrm{{\boldmath $T$}}^n=\{T_0^n,\ldots,T_r^n\}$, and let $(\textrm{{\boldmath $T$}}^n)$ be the ideal of $S$ generated by the set $\textrm{{\boldmath $T$}}^n$.
For $0\leq n\leq m$ we have the map $\varphi_{nm}\colon S/(\textrm{{\boldmath $T$}}^n)\rightarrow S/(\textrm{{\boldmath $T$}}^m)$ defined by the multiplication by $(T_0\cdots T_r)^{m-n}$, and thus we get an inductive system $\{S/(\textrm{{\boldmath $T$}}^n),\varphi_{nm}\}$ indexed by the set of non-negative integers.
For any $(r+1)$-tuple $(p_0,\ldots,p_r)$ of positive integers and an integer $n$ such that $n\geq\sup_{0\leq i\leq r}p_i$, we denote by $\xi^{(n)}_{p_0\ldots p_r}$ the modulo $(\textrm{{\boldmath $T$}}^n)$ class of the monomial $T^{n-p_0}_0\cdots T^{n-p_r}_r$.
Clearly, the sequence $\{\xi^{(n)}_{p_0\ldots p_r}\}_n$ defines the unique element of the inductive limit $\varinjlim_nS/(\textrm{{\boldmath $T$}}^n)$, which we denote by $\xi_{p_0\ldots p_r}$.

\begin{thm}[{\cite[$\mathbf{III}$, (2.1.12)]{EGA}}]\label{thm-cohprojsp-0}
Set 
$$
\H^q(X,\O_X(\ast))=\bigoplus_{n\in\Z}\H^q(X,\O_X(n))
$$
for any $q\geq 0$, and regard it as a graded $A$-module.

{\rm (1)} We have $\H^q(X,\O_X(\ast))=0$ for $q\neq 0,r$.

{\rm (2)} There exists a canonical isomorphism $S\cong\H^0(X,\O_X(\ast))$ of graded $A$-modules, where $S$ is regarded as a graded $A$-module in the standard way.

{\rm (3)} The graded $A$-module $\H^r(X,\O_X(\ast))$ is canonically isomorphic to the inductive limit $\varinjlim_nS/(\textrm{{\boldmath $T$}}^n)$, which is free with the basis $\{\xi_{p_0\ldots p_r}\}_{p_0,\ldots,p_r>0}$ and equipped with the grading such that the degree of $\xi_{p_0\ldots p_r}$ is $-(p_0+\cdots+p_r)$. \hfill$\square$
\end{thm}

\begin{cor}[{\cite[$\mathbf{III}$, (2.1.13)]{EGA}}]\label{cor-cohprojsp-0}
In the situation as above the cohomology group $\H^q(X,\O_X(n))$ is a free $A$-module of finite type$;$ it is non-zero if and only if either one of the following conditions holds$:$
\begin{itemize}
\item[{\rm (a)}] $q=0$ and $n\geq 0;$ 
\item[{\rm (b)}] $q=r$ and $n\leq -(r+1)$. \hfill$\square$
\end{itemize}
\end{cor}

\begin{prop}[{\cite[$\mathbf{II}$, (2.7.9)]{EGA}}]\label{prop-genglobalsec-0}
Let $A$ be a ring, and $\mathscr{F}$ be a quasi-coherent sheaf of finite type on $X=\P^r_A$. 
Then there exists an integer $N$ such that for any $n\geq N$ the sheaf $\mathscr{F}(n)$ is generated by global sections$;$ more precisely, there exists a surjective morphism $\O_X^{\oplus k}\rightarrow\mathscr{F}(n)$, where $k$ is a positive integer depending on $n$. \hfill$\square$
\end{prop}

\subsubsection{Ample and very ample sheaves}\label{subsub-cohprerequi-0-ample}
Let us briefly recall the definitions of ample and very ample sheaves (cf.\ \cite[$\mathbf{II}$, (4.4.2), (4.6.11)]{EGA}).
Let $f\colon X\rightarrow Y$ be a morphism of finite type between coherent algebraic spaces, and $\mathscr{L}$ an invertible sheaf on $X$.
\begin{itemize}
\item We say that $\mathscr{L}$ is {\em very ample relative to $f$}\index{ample invertible sheaf@ample (invertible sheaf)!very ample invertible sheaf@very ---} if there exist a quasi-coherent $\O_X$-module $\mathscr{E}$ of finite type and a factorization 
$$
\xymatrix{X\ar[r]_{i}\ar@/^1pc/[rr]^f&\P(\mathscr{E})\ar[r]_{\pi}&Y}
$$
of $f$ by an immersion $i$ such that $\mathscr{L}$ is isomorphic to $i^{\ast}\O_{\P(\mathscr{E})}(1)$.
\item We say that $\mathscr{L}$ is {\em ample relative to $f$}\index{ample invertible sheaf@ample (invertible sheaf)} (or {\em $f$-ample}) if for any quasi-compact open subset $V$ of $Y$, there exists a positive integer $k$ such that $\mathscr{L}^{\otimes k}|_{f^{-1}(V)}$ is very ample relative to $f_V\colon f^{-1}(V)\rightarrow V$.
\end{itemize}

\subsection{More basics on algebraic spaces}\label{sub-morebasicsalgspaces}
Let us include some useful facts on coherent algebraic spaces, which will be used in our later discussion.
\subsubsection{The stratification by subschemes}\label{subsub-amazingstratification}
The following theorem, which shows that coherent algebraic spaces are `tangible', is very useful to reduce many arguments concerning with algebraic spaces to scheme cases, and thus has a lot of important applications:
\begin{thm}[{\cite[Premi\`ere partie, (5.7.6)]{RG}}]\label{thm-RGamazing}
Let $X$ be a coherent algebraic space.
Then there exists a finite sequence $Z_1,\ldots,Z_r$ of locally closed subspaces of $X$ that fulfills the following properties$:$
\begin{itemize}
\item[{\rm (a)}] $Z_i$ for each $i=1,\ldots,r$ is reduced and quasi-compact$;$
\item[{\rm (b)}] $Z_i$'s are disjoint to each other and cover $X$, that is, $X=\bigcup_{i=1}^rZ_i;$
\item[{\rm (c)}] $Y_i=\bigcup_{j\geq i}Z_j$ for each $i=1,\ldots,r$ is an open subspace of $X;$
\item[{\rm (d)}] for each $i=1,\ldots,r$ there exists a separated and quasi-compact elementary \'etale neighborhood $Y'_i$ of $Z_i$ in $Y_i$ such that the image of $Y'_i$ in $X$ coincides with $Y_i$. \hfill$\square$
\end{itemize}
\end{thm}

Here, for an algebraic space $Y$ and closed subspace $Z$ of $Y$, an {\em elementary \'etale neighborhood}\index{elementary etale neighborhood@elementary \'etale neighborhood} of $Z$ in $Y$ is an \'etale map $u\colon Y'\rightarrow Y$ from a {\em scheme} $Y'$ such that the induced map $Y'\times_YZ\rightarrow Z$ is an isomorphism.
In particular, if an elementary \'etale neighborhood of $Z$ in $Y$ exists, then $Z$ is a scheme, since it is a closed subscheme of $Y'$.

\subsubsection{Affineness criterion}\label{subsub-KnudsonSerrecriterion}
\index{affine!affineness criterion@---ness criterion}
\begin{thm}[Affineness criterion; cf.\ {\cite[III.2.3]{Knu}}]\label{thm-KnudsonSerrecriterion}
Let $X$ be a coherent algebraic space.
Then the following conditions are equivalent$:$
\begin{itemize}
\item[{\rm (a)}] the global section functor 
$$
\Gamma_X\colon\QCoh_X\longrightarrow\Ab
$$
from the category of quasi-coherent sheaves on $X$ is exact and faithful, that is, for any quasi-coherent sheaf $\mathscr{F}$ on $X$, $\Gamma_X(\mathscr{F})\neq 0$ whenever $\mathscr{F}\neq 0;$
\item[{\rm (b)}] $X$ is an affine scheme. \hfill$\square$
\end{itemize}
\end{thm}

\begin{rem}{\rm 
In \cite[III.2.3]{Knu} the theorem is proved under the slightly stronger hypothesis that $X$ is quasi-compact and separated.
One can modify the proof therein to show the above version of the theorem, changing `quasi-compact and separated' for the map $\gamma\colon X\rightarrow\Spec A$ into `coherent'.
In view of \ref{prop-cohqcoh-0} this is enough for the rest of the proof.}
\end{rem}

The same remark shows that one can similarly drop `separated' in \cite[III.2.5]{Knu} as follows:
\begin{thm}[The Serre criterion; cf.\ {\cite[III.2.5]{Knu}}]\label{thm-KnudsonSerrecriterion2}
Let $X$ be a Noetherian algebraic spaces.
Then the following conditions are equivalent$:$
\begin{itemize}
\item[{\rm (a)}] the global section functor 
$$
\Gamma_X\colon\Coh_X\longrightarrow\Ab
$$
from the category of coherent sheaves on $X$ is exact$;$
\item[{\rm (b)}] $X$ is an affine scheme. \hfill$\square$
\end{itemize}
\end{thm}

\begin{cor}\label{cor-KnudsonSerrecriterion}
Let $X$ be a coherent algebraic space, and $X_0$ the closed subspace of $X$ defined by a nilpotent quasi-coherent ideal $\mathscr{J}$.
If $X_0$ is a scheme, then $X$ is a scheme.
\end{cor}

\begin{proof}
By induction with respect to $s\geq 1$ with $\mathscr{J}^s=0$, we reduce to the case $\mathscr{J}^2=0$.
We may also assume that $X_0$ is affine.
For any quasi-coherent sheaf $\mathscr{F}$ on $X$ there exists an exact sequence
$$
0\longrightarrow\mathscr{J}\mathscr{F}\longrightarrow\mathscr{F}\longrightarrow\mathscr{F}/\mathscr{J}\mathscr{F}\longrightarrow 0,
$$
where the first and the third sheaves can be regarded as quasi-coherent sheaves on $X_0$.
Since $X_0$ is affine, we have $\H^1(X,\mathscr{J}\mathscr{F})=0$ and $\H^1(X,\mathscr{F}/\mathscr{J}\mathscr{F})=0$ (\ref{thm-vanishcohaff-0} (1)).
Hence $\H^1(X,\mathscr{F})=0$ for any quasi-coherent sheaf $\mathscr{F}$ on $X$.
In particular, the functor $\Gamma_X$ is exact.
To show that $\Gamma_X$ is faithful, suppose $\Gamma(X,\mathscr{F})=0$.
Then $\Gamma(X,\mathscr{J}\mathscr{F})=\Gamma(X,\mathscr{F}/\mathscr{J}\mathscr{F})=0$ by the exactness.
Hence $\mathscr{J}\mathscr{F}=0$ and $\mathscr{F}/\mathscr{J}\mathscr{F}=0$, which means $\mathscr{F}=0$.
Now by \ref{thm-KnudsonSerrecriterion} we deduce that $X$ is an affine scheme.
\end{proof}

\subsubsection{Limit theorem}\label{subsub-grusonraynaudlimit}
\begin{thm}[Raynaud {\cite[Premi\`ere partie, (5.7.8)]{RG}}]\label{thm-grusonraynaudlimit}
\index{Raynaud, M.}
Let $X$ be a coherent algebraic space.
Then any quasi-coherent sheaf $\mathscr{F}$ on $X$ is the inductive limit $\varinjlim_{\alpha\in I}\mathscr{F}_{\alpha}$ of quasi-coherent subsheaves $\mathscr{F}_{\alpha}\subseteq\mathscr{F}$ of finite type. \hfill$\square$
\end{thm}

This theorem has been proved by Knutson\index{Knutson, D.} \cite[III.1.1]{Knu} in case $X$ is a Noetherian locally separated algebraic space.
The proof for the general case uses \ref{thm-RGamazing}.
By this and by an argument similarly to the proof of \cite[$\mathbf{I}$, (9.4.3)]{EGA}, one has:
\begin{cor}[Extension theorem]\label{cor-grusonraynaudlimit}
Let $X$ be a coherent algebraic space, and $U$ a quasi-compact open subspace of $X$.
Then for any quasi-coherent sheaf $\mathscr{F}$ of finite type on $U$, there exists a quasi-coherent sheaf $\mathscr{G}$ of finite type on $X$ such that $\mathscr{G}|_U=\mathscr{F}$. \hfill$\square$
\end{cor}

The first author of this book has proved in \cite{Fujiw3} the following absolute affine limit theorem for algebraic spaces, which generalizes \cite[C.9]{TT}:
\begin{thm}[Affine limit theorem]\label{thm-absoluteaffinelimittheoremalgebraicspace}
Let $S$ be a coherent algebraic space, and $f\colon X\rightarrow S$ a coherent morphism of algebraic spaces.

{\rm (1)} There exists a projective system $\{X_i\}_{i\in I}$ of $S$-schemes indexed by a category $I$ such that$:$
\begin{itemize}
\item[{\rm (a)}] each $X_i$ is coherent and finitely presented over $S;$
\item[{\rm (b)}] for any arrows $i\rightarrow j$ in $I$ the transition map $X_j\rightarrow X_i$ is affine$;$
\item[{\rm (c)}] $X\cong\varprojlim_{i\in I}X_i$.
\end{itemize}

{\rm (2)} If $X$ is a scheme, then the algebraic spaces $X_i$'s can be taken to be schemes.

{\rm (3)} If $S$ is Noetherian, then the index category $I$ can be replaced by a directed set in such a way that each transition map $X_j\rightarrow X_i$ is scheme-theoretically dominant.\hfill$\square$
\end{thm}

By this theorem and \ref{cor-KnudsonSerrecriterion} we have:
\begin{cor}\label{cor-KnudsonSerrecriterionabsolute}
Let $X$ be a quasi-separated algebraic space.
If $X_{\red}$ is a scheme $($resp.\ an affine scheme$)$, then so is $X$.\hfill$\square$
\end{cor}

\addcontentsline{toc}{subsection}{Exercises}
\subsection*{Exercises}
\begin{exer}\label{exer-lemqcptcomplementsupp}{\rm 
Let $X$ be a coherent scheme, and $U\subseteq X$ an open subset.
Then show that $U$ is quasi-compact if and only if $X\setminus U$ is the support of a closed subscheme of $X$ of finite presentation.}
\end{exer}

\begin{exer}\label{exer-closedimmersionnilpotentthickening}{\rm 
Let $X$ be a scheme, $\mathscr{I}\subseteq\O_X$ a nilpotent quasi-coherent ideal of finite type, and $Z$ the closed subscheme of $X$ defined by $\mathscr{I}$.
Let $f\colon Y\rightarrow X$ be a morphism of schemes such that $f_Z\colon Y\times_XZ\rightarrow Z$ is a closed immersion.
Then show that $f$ is a closed immersion.}
\end{exer}

\begin{exer}\label{exer-derivedformalism3}{\rm 
Let $A$ be a universally coherent ring, and $X$ a projective finitely presented $A$-scheme.
Then show that the canonical exact functor
$$
\delta^{\bd}\colon\DC^{\bd}(\Coh_X)\longrightarrow\DC^{\bd}_{\coh}(X)
$$
is a categorical equivalence.}
\end{exer}


\section{Valuation rings}\label{sec-val}
In this section, we give a brief overview of the theory of valuation rings.
Our basic references for valuation rings are \cite[Chap.\ VI]{ZSII} and \cite[Chap.\ VI]{Bourb1}. 
We also refer to \cite{Vaqu} as a short but concise survey.
Since almost all what we need to know about valuation rings are already well-documented in these references, we most of the time limit ourselves to be sketchy, and omit many of the proofs.

In \S\ref{sub-valpre} we will discuss two prerequisites, totally ordered commutative groups and invertible ideals.
The basic definitions and first properties of valuations and valuation rings will be given in \S\ref{sub-valringsval}, which also includes the definitions of height (also called `rank' in some literatures) and rational rank. 
As described in \S\ref{sub-schematic}, the spectrum of a valuation ring is a {\em path-like} object; this might suggest that valuation rings are algebro-geomtric analogue of paths, perhaps more precisely `long paths', reflecting the fact that the valuations may possible be of large height.
This subsection also gives a detailed description of valuation rings of finite height, especially of height one, together with the concept of the so-called {\em non-archimedean norms}.
In \S\ref{sub-composition} we discuss composition and decomposition of valuation rings, which are one of the most characteristic features of valuation rings and their spectra.
In the next two subsections \S\ref{sub-estimate} and \S\ref{sub-valexample}, we recall some techniques for studying structures of valuation rings, which enable us to give a rough classification of them, which we do in the end of \S\ref{sub-valexample}.

The valuation rings discussed in \S\ref{sub-adicsepval} are in fact the most important for our purpose, that is, valuation rings equipped with a separated adic topology defined by a principal ideal. 
Such valuation rings have many significant features, which will be of fundamental importance in our later studies of formal and rigid geometries.

\index{valuation!valuation ring@--- ring|(}
\subsection{Prerequisites}\label{sub-valpre}
\subsubsection{Totally ordered commutative group}\label{subsub-ordab}
\index{ordered!totally@totally ---!totally ordered commutative group@--- --- commutative group|(}
A {\em totally ordered commutative group} is an abelian group together with an ordering $(\Gamma,\geq)$ (cf.\ \S\ref{subsub-orderings}) such that the following conditions are satisfied:
\begin{itemize}
\item[{\rm (a)}] if $a\geq b$ for $a,b\in\Gamma$, then $a+c\geq b+c$ for any $c\in\Gamma$; 
\item[{\rm (b)}] for any $a,b\in\Gamma$ either $a\geq b$ or $b\geq a$ holds; that is, $(\Gamma,\geq)$ is a totally ordered set (cf.\ \S\ref{subsub-orderings}).
\end{itemize}

An {\em ordered homomorphism}\index{ordered!ordered homomorphism@--- homomorphism} of totally ordered commutative groups is a group homomorphism that is an ordered map (cf.\ \S\ref{subsub-orderings}).

Let $\Gamma$ be a totally ordered commutative group. 
An element $a\in\Gamma$ is said to be positive (resp.\ negative) if $a>0$ (resp.\ $a<0$).
A non-empty subset $\Delta$ of $\Gamma$ is said to be a {\em segment}\index{segment} if for any element $a\in\Delta$ any $b\in\Gamma$ with $-a\leq b\leq a$ or $a\leq b\leq-a$ belongs to $\Delta$.
If a subgroup $\Delta$ of $\Gamma$ is a segment, then it is called an {\em isolated subgroup}\index{isolated subgroup}.
Notice that, unlike in \cite[\S 3]{Vaqu}, we allow $\Gamma$ itself to be an isolated subgroup.

The basic role of isolated subgroups is explained in the following `Homomorphism Theorem' for totally ordered commutative groups:
\begin{prop}[{\cite[Chap.\ VI, \S4.2, Prop.\ 3]{Bourb1}}]\label{prop-isolated1}
Let $\Gamma$ be a totally ordered commutative group.

{\rm (1)} The kernel of an ordered homomorphism of $\Gamma$ to an ordered group is an isolated subgroup of $\Gamma$.

{\rm (2)} Conversely, for an isolated subgroup $\Delta\subseteq\Gamma$ the quotient $\Gamma/\Delta$ is again a totally ordered commutative group by the induced ordering, and the canonical map $\Gamma\rightarrow\Gamma/\Delta$ is an ordered homomorphism. \hfill$\square$
\end{prop}

Moreover, it is easy to see that the canonical map $\Gamma\rightarrow\Gamma/\Delta$ in (2) induces the bijection between the set of all isolated subgroups of $\Gamma$ containing $\Delta$ and the set of all isolated subgroups of $\Gamma/\Delta$.

For a totally ordered commutative group $\Gamma$ we denote by $\Isol(\Gamma)$ the set of all proper isolated subgroups of $\Gamma$.
Then $\Isol(\Gamma)$ together with the inclusion order is a totally ordered set; indeed, if there were two isolated subgroups $\Delta$ and $\Delta'$ such that none of $\Delta\subseteq\Delta'$ and $\Delta\supseteq\Delta'$ holds, then there would be positive elements $a\in\Delta\setminus\Delta'$ and $a'\in\Delta'\setminus\Delta$; if, for example, $a\geq a'$, then $a'$ must belong to $\Delta$, for $\Delta$ is isolated, which is absurd.
The order type\index{order type} (cf.\ \S\ref{subsub-orderings}) of $\Isol(\Gamma)$ is called the {\em height}\index{height!height of totally ordered commutative group@--- (of a totally ordered commutative group)} of $\Gamma$ and is denoted by
$$
\mathrm{ht}(\Gamma).
$$
If it is finite and equal to $n$ (the order type of the totally ordered set $\{0,1,\ldots,n-1\}$ with the obvious ordering), then we say that $\Gamma$ is of finite height $n$. 
Otherwise, $\Gamma$ is said to be of infinite height.

Note that the height $0$ totally ordered commutative group is the trivial group $\{0\}$.
As for height one groups, we have the following characterization:
\begin{prop}[{\cite[Chap.\ VI, \S4.5, Prop.\ 8]{Bourb1}}]\label{prop-isolated2}
The following conditions for a totally ordered commutative group $\Gamma$ are equivalent$:$
\begin{itemize}
\item[{\rm (a)}] $\mathrm{ht}(\Gamma)=1;$ 
\item[{\rm (b)}] for any $a,b\in\Gamma$ with $a>0$ and $b\geq 0$ there exists an integer $n\geq 0$ such that $b\leq na;$ 
\item[{\rm (c)}] $\Gamma$ is ordered isomorphic to a non-zero subgroup of $\R$, the additive group of real numbers endowed with the usual order. \hfill$\square$
\end{itemize}
\end{prop}

To discuss totally ordered commutative groups of higher height, the following construction will be useful:
\begin{exa}\label{exa-height}
{\rm Let $h$ be a positive integer, and $\Gamma_i$ totally ordered commutative groups for $i=1,\ldots,h$.
Consider the direct sum
$$
\Gamma=\Gamma_1\oplus\cdots\oplus\Gamma_h
$$
endowed with the so-called {\em lexicographical order}\index{ordering@order(ing)!lexicographical ordering@lexicographical ---}: for two elements $a=(a_1,\ldots,a_h)$ and $b=(b_1,\ldots,b_h)$ of $\Gamma$, 
$$
a\leq b\quad\Longleftrightarrow\quad\left\{
\begin{minipage}{20em}
{\small the first entries $a_i$ and $b_i$ in $a$ and $b$ from the left that are different from each other satisfy $a_i\leq b_i$.}
\end{minipage}\right.
$$
For $j=1,\ldots,h$ the subgroup of $\Gamma$ of the form $\Gamma_j\oplus\cdots\oplus\Gamma_h$ ($=0\oplus\cdots\oplus 0\oplus\Gamma_j\oplus\cdots\oplus\Gamma_h$) is an isolated subgroup, and the quotient of $\Gamma$ by this subgroup is ordered isomorphic to $\Gamma_1\oplus\cdots\oplus\Gamma_{j-1}$.
Hence by induction one sees $\mathrm{ht}(\Gamma)=\sum^h_{i=1}\mathrm{ht}(\Gamma_i)$.}
\end{exa}

In particular, if all $\Gamma_i$ are of height one, the resulting $\Gamma$ as above gives a totally ordered commutative group of height $h$.
While it is not true that any totally ordered commutative group of finite height is of this form, we have the following useful fact:
\begin{prop}[{\cite[Chap.\ II, Prop.\ 2.10]{Abhya1}}]\label{prop-isolated3}
Let $\Gamma$ be a totally ordered commutative group of height $n<\infty$. Then $\Gamma$ is ordered isomorphic to a subgroup of $\R^n=\R\oplus\cdots\oplus\R$. Suppose, moreover, that $\Gamma$ satisfies the following condition$:$ for any $a\in\Gamma$ and any non-zero integer $m$ there exists $b\in\Gamma$ such that $a=mb$. Then there exists a subgroup $\Gamma_i$ of $\R$ such that $\Gamma$ is order isomorphic to $\Gamma_1\oplus\cdots\oplus\Gamma_n$ with the lexicographical order. \hfill$\square$
\end{prop}

The {\em rational rank}\index{rational rank!rational rank of a totally ordered commutative group@--- (of a totally ordered commutative group)} of a totally ordered commutative group $\Gamma$ is the dimension of the $\Q$-vector space $\Gamma\otimes_{\Z}\Q$ and is denoted by
$$
\mathrm{rat\textrm{-}rank}(\Gamma).
$$

\begin{prop}\label{prop-rationalrank}
Let $\Gamma$ be a totally ordered commutative group, and suppose the rational rank of $\Gamma$ is finite. 
Then the height of $\Gamma$ is finite, and following inequality holds$:$
$$
\mathrm{ht}(\Gamma)\leq\mathrm{rat\textrm{-}rank}(\Gamma).\eqno{\square}
$$
\end{prop}

For the proof, see \cite[Chap.\ VI, \S10.2]{Bourb1} or \cite[Chap.\ VI, \S10, Note]{ZSII}.
\begin{exa}\label{exa-rationalrank1}
{\rm Consider the totally ordered commutative group $\Gamma$ as in \ref{exa-height}, where each $\Gamma_i$ is of height one. Then the rational rank of $\Gamma$ is finite if and only if the rational rank of each $\Gamma_i$ ($i=1,\ldots,h$) is finite. In this case we have $\mathrm{rat\textrm{-}rank}(\Gamma)=\sum^h_{i=1}\mathrm{rat\textrm{-}rank}(\Gamma_i)\geq h=\mathrm{ht}(\Gamma)$.}
\end{exa}
\index{ordered!totally@totally ---!totally ordered commutative group@--- --- commutative group|)}

\subsubsection{Invertible ideals}\label{subsub-invertible}
\index{ideal!invertible ideal@invertible ---|(}
Let $A$ be a ring, and $F=\Frac(A)$ the total ring of fractions of $A$. 
An $A$-submodule $I\subseteq F$ of $F$ is said to be {\em non-degenerate} if $F\cdot I=F$.
Hence, in particular, an ideal $I$ of $A$ is non-degenerate if and only if $I$ contains a non-zero-divisor of $A$.
It can be shown that for a non-degenerate $A$-submodule $I\subseteq F$ the following conditions are equivalent (cf.\ \cite[Chap.\ II, \S5.6, Theorem 4]{Bourb1}):
\begin{itemize}
\item[{\rm (a)}] there exists an $A$-submodule $J$ of $F$ such that $I\cdot J=A$; 
\item[{\rm (b)}] $I$ is projective; 
\item[{\rm (c)}] $I$ is finitely generated, and for any maximal ideal $\m\subseteq A$ the $A_{\m}$-module $I_{\m}$ is principally generated.
\end{itemize}

If these conditions are satisfied, we say that $I$ is an {\em invertible fractional ideal}.
In case $I=(a)$ is a principal fractional ideal, then $I$ is invertible if and only if $a$ is a non-zero-divisor of $F$.
If an invertible fractional ideal $I$ is an ideal of $A$, we say that $I$ is an {\em invertible ideal} of $A$. 

\begin{prop}\label{prop-invertible1}
Let $I\subseteq F$ be an invertible fractional ideal. 
Then the $A$-submodule $J$ as in the condition {\rm (a)} above is unique and is given by
$$
J=(A:I)\ (=\{x\in F\,|\,xI\subseteq A\}).\eqno{\square}
$$
\end{prop}

In particular, if $I=aA$ for $a\in F$ is invertible, $J=(A:I)=(a^{-1})$.
The set of all invertible fractional ideals of $F$ forms a group by multiplication:
\begin{lem}\label{lem-invertible3}
Let $A$ be a ring, and $I,J$ ideals of $A$. Then both $I$ and $J$ are invertible if and only if $IJ$ is invertible.
\end{lem}

\begin{proof}
Suppose $K=IJ$ is invertible.
Since $IF\subseteq F=IJF\subseteq IF$, we have $IF=F$, that is, $I$ is non-degenerate.
Take the fractional ideal $L$ such that $KL=A$.
Then $JL$ gives the inverse of $I$, and hence $I$ is invertible.
One sees similarly that $J$ is also invertible.
The converse is clear.
\end{proof}

The case where $A$ is a local domain will be of particular importance. 
In this situation, by virtue of the condition {\rm (c)} above, any invertible ideal is principal. 
\begin{prop}\label{prop-invertible2}
If $A$ is a local domain, then there exists a canonical bijection between the set of all invertible fractional ideals of $F=\Frac(A)$ and the set $F^{\times}/A^{\times};$ the bijection is established by $(a)\mapsto [a]=(a\ \mathrm{mod}\ A^{\times})$. Moreover, if we order the former set by inclusion and in the latter set by $[x]\leq[y]$ $\Leftrightarrow$ $x=zy$ for some $z\in A$, then this bijection gives an isomorphism of ordered groups. \hfill$\square$
\end{prop}
\index{ideal!invertible ideal@invertible ---|)}

\subsection{Valuation rings and valuations}\label{sub-valringsval}
\subsubsection{Valuation rings}\label{subsub-valdfn}
Let $B$ be a local ring, and $A$ a subring of $B$ that is again a local ring. We say that {\em $B$ dominates $A$}\index{dominate}, written $A\preceq B$, if $\m_A\subseteq\m_B$ or, equivalently, $\m_A=A\cap\m_B$.
For a field $K$ the relation $\preceq$ gives an ordering on the set of all local subrings of $K$. 

\begin{dfn}\label{dfn-val}
{\rm Let $V$ be an integral domain, and $K$ a field containing $V$ as a subring. Then $V$ is said to be a {\em valuation ring for $K$}\index{valuation!valuation ring@--- ring} if it satisfies one of the following equivalent conditions:
\begin{itemize}
\item[{\rm (a)}] $V$ is a maximal with respect to $\preceq$ in the set of all local subrings of $K$; 
\item[{\rm (b)}] for any $x\in K\setminus\{0\}$ either $x$ or $x^{-1}$ belongs to $V$; 
\item[{\rm (c)}] we have $\Frac(V)=K$, and the set of all ideals of $V$ with the ordering by the inclusion order is totally ordered; 
\item[{\rm (d)}] we have $\Frac(V)=K$, and the set of all principal ideals of $V$ with the ordering by the inclusion order is totally ordered.
\end{itemize}}
\end{dfn}

For the equivalence of these conditions, see the references mentioned at the beginning of this section.
When we just say $V$ is a valuation ring, we always mean that $V$ is a valuation ring for its field of fractions.
By the property (b) we have the following easy but useful fact: any subring of a field $K$ that contain at least one valuation ring for $K$ is again a valuation ring for $K$.

\danger{Notice that, according to our definition of valuation rings, we allow fields to be valuation rings.
This case is usually ruled out from the notion of valuation rings, but in some places this convention will be useful for the sake of formality.}

As it will turn out in \S\ref{sub-composition}, this convention is consistent with that we allowed in \S\ref{subsub-ordab} a totally ordered commutative group $\Gamma$ itself to be an isolated subgroup of itself.

\begin{prop}\label{prop-val1}
{\rm (1)} Any valuation ring is integrally closed.

{\rm (2)} Any finitely generated ideal of a valuation ring is principal.
\end{prop}

\begin{proof}
By \cite[Chap.\ V, \S2.1, Theorem 1]{Bourb1} we can find a prime ideal $\mathfrak{p}$ of the integral closure $\til{V}$ of $V$ in $K$ that lies over $\m_V$, that is, $V\cap\mathfrak{p}=\m_V$. 
If $\til{V}\neq V$, then $\til{V}_{\mathfrak{p}}$ would be a local subring of $K$ strictly larger than and dominating $V$ and thus leads to a contradiction.
Thus we conclude $\til{V}=V$, thereby (1). 
(2) follows easily from the property \ref{dfn-val} {\rm (c)}.
\end{proof}

\begin{prop}\label{prop-lemvalaval}
Let $V$ be a valuation ring, and $I\subsetneq V$ a finitely generated ideal not equal to $V$ itself.
Then there exists the minimal prime ideal $\mathfrak{p}\subseteq V$ among the prime ideals containing $I;$ more explicitly, $\mathfrak{p}=\sqrt{I}$.
\end{prop}

\begin{proof}
By \ref{prop-val1} (2) we have $I=(a)$ for $a\in\m_V$.
It suffices to show that the ideal $\sqrt{(a)}$ is a prime ideal.
Suppose $bc\in\sqrt{(a)}$ and $b\not\in\sqrt{(a)}$.
This implies that there exists $n\geq 0$ such that $(bc)^n=ad$ for some $d\in V$ and that $a/b^m\in\m_V$ for any $m\geq 0$.
Then $c^{2n}=a\cdot(a/b^{2n})\cdot d^2\in (a)$, and hence $c\in\sqrt{(a)}$.
\end{proof}

\subsubsection{Valuations}\label{subsub-valuation}
\index{valuation|(}
Let $V$ be a valuation ring, and $K=\Frac(V)$.
Since $V$ is a local domain, the ordered group of all invertible fractional ideals is isomorphic to $K^{\times}/V^{\times}$ (\ref{prop-invertible2}).
By \ref{dfn-val} {\rm (d)} this is a totally ordered commutative group.
We set $\Gamma_V=K^{\times}/V^{\times}$ and write the group operation additively.
Consider the mapping
$$
{\textstyle v\colon K\longrightarrow\Gamma_V\cup\{\infty\}},\quad v(x)=
\begin{cases}
[x]\ (=x\ \mathrm{mod}\ V^{\times})&\textrm{if}\ x\neq 0,\\
\infty&\textrm{if}\ x=0.
\end{cases}
$$
As it is easily verified, the map $v$ enjoys the following properties:
\begin{itemize}
\item[{\rm (a)}] $v(xy)=v(x)+v(y)$ for $x,y\in K$; 
\item[{\rm (b)}] $v(x+y)\geq\inf\{v(x),v(y)\}$ for $x,y\in K$; 
\item[{\rm (c)}] $v(1)=0$ and $v(0)=\infty$.
\end{itemize}

\begin{dfn}\label{dfn-valuation1}
{\rm Let $A$ be a ring, and $\Gamma$ a totally ordered commutative group. 
A mapping $v\colon A\rightarrow\Gamma\cup\{\infty\}$ is called a {\em valuation} on $A$ with values in $\Gamma$ if it satisfies the conditions (a), (b), and (c) as above (with $K$ replaced by $A$).
In this situation, we call the totally ordered group $\Gamma$ the {\em value target group} of the valuation $v$.}
\end{dfn}

Thus any valuation ring $V$ induces the canonical valuation $v\colon K\rightarrow\Gamma_V$ on its fractional field; we call this the {\em valuation associated to $V$}\index{valuation!associated valuation@associated ---} and call the totally ordered commutative group $\Gamma_V$ the {\em value group}\index{value group} of $V$.
Notice that in this situation we have $V=\{x\in K\,|\,v(x)\geq 0\}$.
\begin{prop}\label{prop-valuation1}
Let $v\colon K\rightarrow\Gamma\cup\{\infty\}$ be a valuation on a field with values in a totally ordered commutative group $\Gamma$. Then $V=\{x\in K\,|\,v(x)\geq 0\}$ is a valuation ring for $K$, and $\{x\in K\,|\,v(x)>0\}$ is the maximal ideal of $V$. \hfill$\square$
\end{prop}

Moreover, the group $\Gamma$ contains the isomorphic copy of $\Gamma_V$, which coincides with the image of $v$.
Note that, if $K=V$, the corresponding valuation $v$ maps all elements in $K^{\times}$ to $0$. 
Such a valuation is called the {\em trivial valuation}\index{valuation!trivial valuation@trivial ---}. 
\index{valuation|)}

\subsubsection{Height and rational rank of valuation rings}\label{subsub-height}
\begin{dfn}\label{dfn-height1}
{\rm Let $V$ be a valuation ring, $v$ the associated valuation, and $\Gamma_V$ the corresponding value group. 

(1) The height of the value group $\Gamma_V$ (\S\ref{subsub-ordab}) is called the {\em height}\index{height!height of a valuation ring@--- (of a valuation (ring))} of $V$ (or of $v$) and denoted by 
$$
\mathrm{ht}(V)=\mathrm{ht}(v)=\mathrm{ht}(\Gamma_V).
$$

(2) The rational rank of the value group $\Gamma_V$ (\S\ref{subsub-ordab}) is called the {\em rational rank}\index{rational rank!rational rank of a valuation ring@ --- (of a valuation (ring))} of $V$ (or of $v$) and denoted by
$$
\mathrm{rat\textrm{-}rank}(V)=\mathrm{rat\textrm{-}rank}(v)=\mathrm{rat\textrm{-}rank}(\Gamma_V).
$$}
\end{dfn}

\begin{rem}
{\rm We prefer to use the terminology `height' following \cite{Bourb1}, while in \cite{ZSII} it is called `rank'.}
\end{rem}

\begin{prop}[cf.\ {\ref{prop-rationalrank}}]\label{prop-rationalrank1}
Let $V$ be a valuation ring, and suppose the rational rank of $V$ is finite. 
Then the height of $V$ is finite, and the following inequality holds$:$
$$
\mathrm{ht}(V)\leq\mathrm{rat\textrm{-}rank}(V).\eqno{\square}
$$
\end{prop}

\begin{prop}\label{prop-height2}
Let $V$ be a valuation ring with $K$ and $\Gamma_V$ as above.
Then there exist canonical order-preserving bijections among the following sets$:$ 
\begin{itemize}
\item[{\rm (a)}] the set of all prime ideals of $V$ with the inclusion order$;$ 
\item[{\rm (b)}] the set of all subrings $(\neq V)$ lying in between $V$ and $K$ $($which are automatically valuation rings$)$ with the reversed inclusion order$;$ 
\item[{\rm (c)}] the set of all proper isolated subgroups\index{isolated subgroup} of $\Gamma_V$ with the reversed inclusion order. \hfill$\square$
\end{itemize}
\end{prop}

The bijections are established as follows:

\medskip\noindent
$\bullet$ {\sl From the set {\rm (a)} to the set {\rm (b)}}: to a prime ideal $\mathfrak{p}\subseteq V$ we associate the local ring $V_{\mathfrak{p}}$; note that by \ref{dfn-val} (b) any local ring lying in between $V$ and $K$ is a valuation ring.

\medskip\noindent
$\bullet$ {\sl From the set {\rm (b)} to the set {\rm (c)}}: for $W$ as in (b) we consider the subgroup $v(W^{\times})$ of $\Gamma_V$, where $v\colon K\rightarrow\Gamma_V\cup\{\infty\}$ is the valuation associated to $V$; then it is easy to see that $v(W^{\times})$ is a segment.

\medskip
See \cite[Chap.\ VI, Theorem 15]{ZSII} and \cite[Chap.\ VI, \S4.3, Prop.\ 4]{Bourb1} for more details.
From this proposition it follows that the cardinality of $\mathrm{ht}(V)$ coincides with the Krull dimension of $V$.
In particular, a valuation ring of height $0$ is nothing but a field.

\subsection{Spectrum of valuation rings}\label{sub-schematic}
\subsubsection{General description}\label{subsub-schematicgeneral}
From \ref{dfn-val} (c) it follows in particular that prime ideals of a valuation ring $V$ are totally ordered\index{ordered!totally@totally ---} with respect to the inclusion order\index{ordering@order(ing)!inclusion ordering@inclusion ---}. 
Hence the spectrum $\Spec V$ can be understood as a `path' with two extremities $(0)$ and $\m_V$; all the other points sit in between these points in such a way that, if $\mathfrak{p}\subseteq\mathfrak{q}$, then $\mathfrak{p}$ sits in between $(0)$ and $\mathfrak{q}$; see Figure \ref{fig-specclosedsetval}.
The `length', so to speak, of the path is (the cardinality of) the height of $V$.
\begin{figure}[htb]
\begin{center}
\setlength{\unitlength}{1pt}
\begin{picture}(130,35)(0,-10)
\put(5,5){\line(1,0){120}}
\put(5,5){\circle*{4}}
\put(65,5){\circle*{4}}
\put(125,5){\circle*{4}}
\put(0,-5){$\scriptstyle{(0)}$}
\put(63,-5){$\scriptstyle{\mathfrak{p}}$}
\put(119,-5){$\scriptstyle{\mathfrak{m}_V}$}
\put(67,8){$\overbrace{\hspace{56pt}}^{V(\mathfrak{p})}$}
\end{picture}
\end{center}
\caption{$\Spec V$ and a closed set}\label{fig-specclosedsetval}
\end{figure}

By \ref{prop-height2} such a linear pattern is reflected in the set of all subrings of $K$ containing $V$ and also in the set of all isolated subgroups of the value group $\Gamma_V$.

The Zariski topology on $\Spec V$ can also be understood intuitively by this picture. 
For two ideals $I,J$ of $V$ with $I\subseteq J$, define the subsets $[I,J]_V$, $(I,J)_V$, $[I,J)_V$, and $(I,J]_V$ as follows:
\begin{equation*}
\begin{split}
[I,J]_V&=\{\mathfrak{p}\in\Spec V\,|\,I\subseteq\mathfrak{p}\subseteq J\},\\ 
(I,J)_V&=\{\mathfrak{p}\in\Spec V\,|\,I\subsetneq\mathfrak{p}\subsetneq J\},\\ 
[I,J)_V&=\{\mathfrak{p}\in\Spec V\,|\,I\subseteq\mathfrak{p}\subsetneq J\},\\ 
(I,J]_V&=\{\mathfrak{p}\in\Spec V\,|\,I\subsetneq\mathfrak{p}\subseteq J\}.
\end{split}
\end{equation*}
Then for any ideal $I\subseteq V$ the set $[I,\m_V]_V$ is exactly the closed set $V(I)$, and any closed sets are of this form.
Hence open sets are exactly the subsets of the form $[(0),I)_V$.
For any point $\mathfrak{p}\in\Spec V$ the Zariski closure of $\{\mathfrak{p}\}$ is given by $[\mathfrak{p},\m_V]_V$; in other words, a point is the specialization\index{specialization} of points sitting on the left and is the generization\index{generization} of points sitting on the right, when $\Spec V$ is oriented as in Figure \ref{fig-specclosedsetval}.

\begin{prop}\label{prop-spectrumval}
The underlying topological space of $\Spec V$, where $V$ is a valuation ring, is a valuative space\index{valuative!valuative topological space@--- (topological) space}\index{space@space (topological)!valuative topological space@valuative ---} {\rm (\ref{dfn-valuativespace})}.
\end{prop}

\begin{proof}
First notice that the underlying topological spaces of affine schemes are coherent and sober.
For any point $\mathfrak{p}\in\Spec V$ the set $G_{\mathfrak{p}}$ of all generizations of $\mathfrak{p}$ is totally ordered, as described above.
\end{proof}

\subsubsection{Valuation rings of finite height}\label{subsub-heightone}
\index{valuation!valuation ring@--- ring!valuation ring of finite height@--- --- of finite height|(}
The valuation ring $V$ is of finite height if and only if $\Spec V$ consists of finitely many points (Figure \ref{fig-valspecfinheight}). 
\begin{figure}[htb]
\begin{center}
\setlength{\unitlength}{1pt}
\begin{picture}(130,35)(0,-10)
\put(5,5){\line(1,0){45}}
\put(100,5){\line(1,0){25}}
\put(5,5){\circle*{4}}
\put(25,5){\circle*{4}}
\put(45, 5){\circle*{4}}
\multiput(55,5)(6,0){8}{\circle*{1}}
\put(105,5){\circle*{4}}
\put(125,5){\circle*{4}}
\put(0,-5){$\scriptstyle{(0)}$}
\put(119,-5){$\scriptstyle{\mathfrak{m}_V}$}
\end{picture}
\end{center}
\caption{Spectrum of a valuation ring of finite height}\label{fig-valspecfinheight}
\end{figure}
In this case, one is able to speak about the adjacent points, that is to say, the `next' generization or specialization of an arbitrary given point of $\Spec V$.
This feature allows us, aided by composition and decomposition of valuation rings (explained later in \S\ref{sub-composition}), to carry out inductive arguments with respect to height, reducing many situations to the height one case.
It is therefore important to study the case $\mathrm{ht}(V)=1$, the case where $\Spec V$ consists only of $(0)$, the open point, and $\m_V$, the closed point:

\begin{prop}[cf.\ {\ref{prop-isolated2}}]\label{prop-heightone1}
Let $V$ be a valuation ring of non-zero height. 
Then the following conditions are equivalent$:$
\begin{itemize}
\item[{\rm (a)}] the height of $V$ is $1;$ 
\item[{\rm (b)}] for any $x\in\m_V\setminus\{0\}$ and $y\in V\setminus\{0\}$ there exist an integer $n\geq 0$ and an element $z\in V$ such that $yz=x^n;$
\item[{\rm (c)}] the value group $\Gamma_V$ is isomorphic to a non-zero subgroup of the ordered additive group $\R$ of real numbers.\hfill$\square$
\end{itemize}
\end{prop}
\index{valuation!valuation ring@--- ring!valuation ring of finite height@--- --- of finite height|)}

\subsubsection{Non-archimedean norms}\label{subsub-nonarchnorms}
\index{norm!non-archimedean norm@non-archimedean ---|(}
By \ref{prop-heightone1} (3) the associated valuation to a height one valuation ring $V$ is of the form $v\colon K\rightarrow\R\cup\{\infty\}$, where $K=\Frac(V)$.
In this situation, the function $|\cdot|\colon K\rightarrow\R_{\geq 0}$ defined by $|x|=e^{-v(x)}$ (where $e>1$ is a fixed real number) is more familiar object, the associated (non-archimedean) {\em norm}. 
In literature the pair $(K,|\cdot|)$ is called a {\em non-archimedean valued field}, which is, in our language, equivalent to the fractional field of a valuation ring {\em of height one}.
For more details of norms we refer to some of the first chapters of \cite{BGR}.

\begin{prop}[{\cite[Chap.\ VI, Theorem 16]{ZSII}}]\label{prop-heightone2}
Let $V$ be a valuation ring of non-zero height.
Then $V$ is Noetherian if and only if the value group $\Gamma_V$ is isomorphic to $\Z$ $($and hence, in particular, $V$ is of height one$)$.
$($In this case $V$ is called a {\em discrete valuation ring}\index{valuation!valuation ring@--- ring!discrete valuation ring@discrete --- ---} $($acronym$:$ {\em DVR}$).)$\hfill$\square$
\end{prop}
\index{norm!non-archimedean norm@non-archimedean ---|)}

\subsection{Composition and decomposition of valuation rings}\label{sub-composition}
Let $V$ be a valuation ring.
As mentioned in \S\ref{subsub-schematicgeneral}, $\Spec V$ looks like a path with the extremities $(0)$ and $\m_V$. 
Any prime ideal $\mathfrak{p}\subseteq V$ divides $\Spec V$ into two segments $[(0),\mathfrak{p}]$ and $[\mathfrak{p},\m_V]$.
The former subsegment corresponds to $\Spec V_{\mathfrak{p}}$ and the latter to $\Spec V/\mathfrak{p}$ (Figure \ref{fig-subvidisionvalspec}).
\begin{figure}[htb]
\begin{center}
\setlength{\unitlength}{1pt}
\begin{picture}(130,45)(0,-10)
\put(5,5){\line(1,0){120}}
\put(5,5){\circle*{4}}
\put(65,5){\circle*{4}}
\put(125,5){\circle*{4}}
\put(0,-5){$\scriptstyle{(0)}$}
\put(63,-5){$\scriptstyle{\mathfrak{p}}$}
\put(119,-5){$\scriptstyle{\mathfrak{m}_V}$}
\put(7,8){$\overbrace{\hspace{56pt}}^{V_{\mathfrak{p}}}$}
\put(67,8){$\overbrace{\hspace{56pt}}^{V/\mathfrak{p}}$}
\end{picture}
\end{center}
\caption{Subdivision of $\Spec V$ into two spectra of valuation rings}\label{fig-subvidisionvalspec}
\end{figure}

\begin{prop}\label{prop-composition1}
{\rm (1)} The ring $V_{\mathfrak{p}}$ is a valuation ring for $K=\Frac(V)$, and $V/\mathfrak{p}$ is a valuation ring for the residue field at $\mathfrak{p}$. 

{\rm (2)} The value group of $V/\mathfrak{p}$ is the isolated subgroup $\Delta$ corresponding to $\mathfrak{p}$ by the correspondence in $\ref{prop-height2}$, and the value group of $V_{\mathfrak{p}}$ is $\Gamma_V/\Delta$.\hfill$\square$
\end{prop}

The proof is easy; see \cite[Chap.\ VI, \S10]{ZSII}\cite[Prop.\ 4.1]{Vaqu}.
In this situation, one recovers $V$ by the formula
$$
V=\{x\in V_{\mathfrak{p}}\,|\,(x\ \mathrm{mod}\ \mathfrak{p}V_{\mathfrak{p}})\in V/\mathfrak{p}\}.
$$
Note that we have $\mathfrak{p}V_{\mathfrak{p}}=\mathfrak{p}$ (shown easily by \ref{dfn-val} (b)).
\begin{prop}\label{prop-composition2}
If we are given a valuation ring $\til{V}$ and a valuation ring $W$ for the residue field of $\til{V}$, the subring $V$ of $\til{V}$ consisting of elements $x$ such that $(x$ mod $\m_{\til{V}})\in W$ is a valuation ring.
Moreover, $\mathfrak{p}=\m_{\til{V}}$ is the prime ideal of $V$ such that $\til{V}\cong V_{\mathfrak{p}}$ and $W\cong V/\mathfrak{p}$. \hfill$\square$
\end{prop}

The proof is easy and left to the reader.
The valuation ring $V$ in this situation is said to be the {\em composite}\index{valuation!valuation ring@--- ring!composite of valuation rings@composite of --- ---s} of the valuation rings $\til{V}$ and $W$.
Schematically, the composition is like gluing two segments $\Spec V_{\mathfrak{p}}$ and $\Spec V/\mathfrak{p}$ at their end points to make a new segment $\Spec V$.
The following proposition is clear:
\begin{prop}\label{prop-composition3}
If a valuation ring $V$ is the composite of $\til{V}$ and $W$, then we have $\mathrm{ht}(V)=\mathrm{ht}(\til{V})+\mathrm{ht}(W)$ and $\mathrm{rat\textrm{-}rank}(V)=\mathrm{rat\textrm{-}rank}(\til{V})+\mathrm{rat\textrm{-}rank}(W)$.\hfill$\square$
\end{prop}

Consistently, we have the following:
\begin{prop}[{\cite[Chap.\ VI, Theorem 17]{ZSII}}]\label{prop-compositionvalvalgroup}
If a valuation ring $V$ is the composite of two valuation rings $\til{V}$ and $W$ as in $\ref{prop-composition2}$, then there exists a canonical exact sequence
$$
0\longrightarrow\Gamma_W\longrightarrow\Gamma_V\longrightarrow\Gamma_{\til{V}}\longrightarrow 0
$$
of totally ordered commutative groups.
\end{prop}

\begin{proof}
We may suppose $\til{V}=V_{\mathfrak{p}}$ and $W=V/\mathfrak{p}$, where $\mathfrak{p}$ is a prime ideal of $V$.
Since $\Gamma_V=K^{\times}/V^{\times}$ and $\Gamma_{\til{V}}=K^{\times}/V^{\times}_{\mathfrak{p}}$ (where $K=\Frac(V)$), we have the canonical surjection $\Gamma_V\rightarrow\Gamma_{\til{V}}$.
Hence it is enough to show that the kernel $V^{\times}_{\mathfrak{p}}/V^{\times}$ is isomorphic to $k^{\times}/W^{\times}$, where $k$ is the residue field of $V_{\mathfrak{p}}$.
But this is clear, since $V=\{x\in V_{\mathfrak{p}}\,|\,(x\ \mathrm{mod}\ \mathfrak{p}V_{\mathfrak{p}})\in W\}$.
\end{proof}

\begin{prop}\label{prop-schematic1}
Let $V$ and $W$ be valuation rings, and $f\colon V\rightarrow W$ be a ring homomorphism.

{\rm (1)} The induced map $\Spec f\colon\Spec W\rightarrow\Spec V$ maps $\Spec W$ surjectively onto the set $[\ker(f),f^{-1}(\m_W)]$ preserving the ordering by inclusion.

{\rm (2)} The map $f$ is local if and only if $\Spec f$ maps $\Spec W$ onto the set $[\ker(f),\m_V]$.

{\rm (3)} The map $f$ is injective if and only if $\Spec f$ maps $\Spec W$ onto the set $[(0),f^{-1}(\m_W)]$.

{\rm (4)} The map $f$ is injective and local $($that is, $W$ dominates $V)$ if and only if the map $\Spec f$ is surjective.
\end{prop}

The proof uses the following easy fact: an inclusion $V\hookrightarrow W$ of valuation rings is local if and only if $W\cap\Frac(V)=V$.

\begin{proof}
First notice that (2) and (3) are immediate consequences of (1).
To show (1), consider the valuation ring $V'=V/\ker(f)$ and the prime ideal $\mathfrak{q}=(f^{-1}(\m_W))V'$ of $V'$.
Set $V''=V'_{\mathfrak{q}}$.
Then the morphism $\Spec V''\rightarrow \Spec V$ maps $\Spec V''$ surjectively onto the subset $[\ker(f),f^{-1}(\m_W)]$ of $\Spec V$; since $V\rightarrow W$ induces a local injective homomorphism $V''\rightarrow W$, (1) follows from (4)

Thus it suffices to show (4).
The `if' part is easy.
Suppose that $V\hookrightarrow W$ be a local injective homomorphism and that $\mathfrak{p}$ a prime ideal of $V$.
We need to find a prime ideal $\mathfrak{q}$ of $W$ such that $\mathfrak{q}\cap V=\mathfrak{p}$.
Consider the localization $V'=V_{\mathfrak{p}}$, and set $W'=W\otimes_VV_{\mathfrak{p}}$.
Since $W\subseteq W'\subseteq\Frac(W)$, $W'$ is a valuation ring for $\Frac(W)$.
Then it suffices to show that $W'$ dominates $V'$.
Since $W$ dominates $V$, we have $W\cap\Frac(V)=V$.
By this we calculate $W'\cap\Frac(V')=W'\cap\Frac(V)=V'$, and thus we conclude that $W'$ dominates $V'$.
\end{proof}

\subsection{Center of a valuation and height estimates for Noetherian domains}\label{sub-estimate}
Let $R$ be an integral domain, and $K=\Frac(R)$. 
Let $v\colon K\rightarrow\Gamma\cup\{\infty\}$ be a valuation on $K$ (\ref{dfn-valuation1}), and $R_v=\{x\in K\,|\,v(x)\geq 0\}$ and $\m_v=\{x\in K\,|\,v(x)>0\}$ the corresponding valuation ring and the maximal ideal, respectively.
Suppose $R\subseteq R_v$, that is, the valuation $v$ takes non-negative values on $R$. 
The prime ideal $\mathfrak{p}=\m_v\cap R$ is said to be the {\em center}\index{valuation!center of a valuation@center of a ---} of the valuation $v$ in $R$.
Note that the valuation ring $R_v$ dominates the local ring $R_{\mathfrak{p}}$.

\begin{prop}\label{prop-center1}
Given an integral domain $R$ and a prime ideal $\mathfrak{p}$, we can always find a valuation $v$ on $K=\Frac(R)$ whose center is $\mathfrak{p}$.\hfill$\square$
\end{prop}

The proof uses Zorn's lemma; see \cite[Chap.\ VI, \S1.2]{Bourb1}.

Let the situation be as above, and set $k=R_{\mathfrak{p}}/\mathfrak{p}R_{\mathfrak{p}}$ and $k_v=R_v/\m_v$. 
Clearly, $k$ is a subfield of $k_v$.
Set 
$$
\mathrm{tr.deg}_kv=\mathrm{tr.deg}_kk_v.
$$

\begin{thm}[Abhyankar]\label{thm-estimate}
Let $R$ be a {\em Noetherian} local domain, $K=\Frac(R)$ the fractional field, and $k=R/\m_R$ the residue field. Suppose $v$ is a valuation on $K$ that dominates $R$. 

{\rm (1)} We have the inequality
$$
\mathrm{rat\textrm{-}rank}(v)+\mathrm{tr.deg}_kv\leq\dim(R),
$$
where $\dim(R)$ denotes the Krull dimension of $R$.

{\rm (2)} If the equality $\mathrm{rat\textrm{-}rank}(v)+\mathrm{tr.deg}_kv=\dim(R)$ holds, then the value group $\Gamma$ of $v$ is isomorphic as a group to $\Z^d$ $($for some $d\geq 0)$, and $k_v$ is a finitely generated extension of $k$.

{\rm (3)} If moreover the stronger equality $\mathrm{ht}(v)+\mathrm{tr.deg}_kv=\dim(R)$ holds, then $\Gamma$ is isomorphic as an ordered group to $\Z^d$ $($with $d=\mathrm{ht}(R))$ equipped with the lexicographical order\index{ordering@order(ing)!lexicographical ordering@lexicographical ---} $($cf.\ $\ref{exa-height})$. \hfill$\square$
\end{thm}

For the proof we refer to \cite{Abhya2} and \cite{Spiva1}; see also \cite[\S 9]{Vaqu}.
We will use the following relative version of the inequality; see \cite[\S 5]{Vaqu} for the proof:
\begin{prop}\label{prop-extensionvaluationestimate}
Let $K\subseteq K'$ be an extension of fields, $v'$ a valuation on $K'$, and $V'\subseteq K'$ the valuation ring of $v'$.
Consider the restriction $v=v|_K$ and the associated valuation ring $V=V'\cap K$.
Let $k'=V'/\m_{V'}$ and $k=V/\m_V$ be the respective residue fields.
Then we have the inequality
$$
\dim_{\Q}(\Gamma_{V'}/\Gamma_V)\otimes\Q+\mathrm{tr.deg}_kk'\leq\mathrm{tr.deg}_KK'.
$$
If we have the equality and if $K'$ is a finitely generated field over $K$, then the group $\Gamma_{V'}/\Gamma_V$ is a finitely generated $\Z$-module, and $k'$ is a finitely generated field over $k$. \hfill$\square$
\end{prop}

\subsection{Examples of valuation rings}\label{sub-valexample}
Theorem \ref{thm-estimate} is powerful for classifying possible valuations centered in a given Noetherian ring. 
Here we exhibit the well-known classification of valuations on an  algebraic function field of dimension $\leq 2$.
The classification in this case was already known to Zariski\index{Zariski, O.} \cite{Zar1}.
For the details and proofs, we refer to \cite[Chap.\ VI, \S15]{ZSII}, \cite{Abhya1}, and \cite{Spiva1}.

\subsubsection{Divisorial valuations}\label{subsub-divval}
Let $R$ be a Noetherian local domain, and $v$ a valuation on $K=\Frac(R)$ centered in $R$. 
The valuation $v$ is said to be {\em divisorial}\index{valuation!divisorial valuation@divisorial ---} if $\mathrm{ht}(v)=1$ and $\mathrm{tr.deg}_kv=\dim(R)-1$.
By \ref{thm-estimate} (3) we know that the value group $\Gamma$ for divisorial $v$ is isomorphic to $\Z$ with the usual ordering.

Divisorial valuations appear mostly in geometric context as follows: Let $X\rightarrow\Spec R$ be a birational morphism, and $D$ a height one regular point of $X$ in the closed fiber. 
Then $\O_{X,D}$ carries canonically the discrete valuation $v$, whose valuation ring coincides with $\O_{X,D}$ and dominates $R$. 
The valuation $v$ thus obtained is clearly a divisorial valuation on $K$ and, in fact, all divisorial valuations are constructed in this way.

Now, let $R$ be a Noetherian local regular domain with the fractional field $K$, and $v$ a non-trivial valuation on $K$ dominating $R$.
Other notations are as in \S\ref{sub-estimate}.

\subsubsection{The case $\dim(R)=1$}\label{subsub-casekrull1}
Since $v$ is not trivial, the height of $v$ must be positive and hence is $1$. 
As $\mathrm{tr.deg}_kv=0$, the valuation $v$ is divisorial and hence is discrete; $R$ is a discrete valuation ring, and $v$ is the associated valuation.
Notice that this valuation is characterized by the formula 
$$
v(f)=\max\{n\geq 0\,|\,f\in\m^n_R\}
$$
for any non-zero $f\in R$.

\subsubsection{The case $\dim(R)=2$}\label{subsub-casekrull2}
Since $v$ is not trivial, we have $\mathrm{tr.deg}_kv\leq 2$.

(1) (Divisorial case) Suppose $\mathrm{tr.deg}_kv=1$. Since the valuation $v$ is non-trivial, we have $\mathrm{rat\textrm{-}rank}(v)=1$.
Hence $v$ is divisorial. 

(2) (Subject-to-divisorial case) If $\mathrm{ht}(v)=\mathrm{rat\textrm{-}rank}(v)=2$, we may apply \ref{thm-estimate} (3) and deduce that the value group $\Gamma$ is isomorphic to $\Z^{\oplus 2}$ with lexicographical order\index{ordering@order(ing)!lexicographical ordering@lexicographical ---}. 
Let $\mathfrak{p}\subseteq R_v$ be the prime ideal lying in between $(0)$ and $\m_{R_v}$, and $\Delta$ the corresponding isolated subgroup of $\Gamma$. 
Then the valuation ring $R_{v,\mathfrak{p}}$, having the value group $\Gamma/\Delta$, is of height one.
Let $v'$ be the valuation associated to $R_{v,\mathfrak{p}}$.
Then we have $\mathrm{tr.deg}_kv'=1$, and hence $v'$ is divisorial. 
The valuation $v$ is the composite of $v'$ and a valuation $v''$ of the residue field $k'$ of $R_{v,\mathfrak{p}}$.
Since $\Delta\cong\Z$ is the value group of $v''$, $v''$ is again divisorial.
As a result, the valuation $v$ in this case is the composite of two divisorial valuations.

Geometrically, such a valuation comes through the following picture. 
Let $X\rightarrow\Spec R$ and $D$ as before such that $\O_{X,D}$ carries the valuation $v'$.
Let $E\rightarrow\ovl{D}$ be a birational morphism to the closed subscheme $\ovl{D}$, the closure of the point $D$ in $X$, and $x$ a regular closed point of $E$.
As $x$ is a height one point of $E$, we may consider the natural discrete valuation $v''$ on $\O_{E,x}$ and the composite $v$ with $v'$.

In fact, it is a general feature for valuations having lexicographically ordered group $\Z^d$ as value groups to have such an inductive structure; they are obtained by successive composition of divisorial valuations. 

(3) (Irrational case) If $\mathrm{rat\textrm{-}rank}(v)=2$ and $\mathrm{ht}(v)=1$, the value group $\Gamma$ is a subgroup of $\R$ of the form $\Z+\Z\tau$, where $\tau$ is an irrational number.

(4) (Limit case) If $\mathrm{rat\textrm{-}rank}(v)=\mathrm{ht}(v)=1$ and $\mathrm{tr.deg}_kv=0$, then the value group is a subgroup of $\Q$.

For the actual construction of the valuations of type (3) and (4), we refer to \cite[\S6]{Zar1}.
\begin{rem}\label{rem-valclassificationberkovich}{\rm 
We will find in {\bf \ref{ch-rigid}}.\ref{subsub-classifyingvaluationsdisk} below a similar list of valuations when classifying points of the unit-disk in rigid geometry; see {\bf \ref{ch-rigid}}.\ref{prop-classifyingvaluationsdisk}.}
\end{rem}

\subsection{$a$-adically separated valuation rings}\label{sub-adicsepval}
\index{valuation!valuation ring@--- ring!a-adically separated valuation ring@$a$-adically separated --- ---|(}
\begin{prop}\label{prop-associatedsepval}
Let $V$ be a valuation ring, and $a\in\m_V\setminus\{0\}$.
Then the ideal $J=\bigcap_{n\geq 1}(a^n)$ is a prime ideal of $V$.
\end{prop}

\begin{proof}
Since $J\neq V$, it suffices to show that $bc\in J$ and $c\not\in J$ imply $b\in J$.
The assumption says that $bc/a^n\in V$ for any $n\geq 1$ and that $a^m/c\in\m_V$ for some $m\geq 1$.
If $b\not\in J$, then $a^n/b\in\m_V$ for some $n\geq 1$, and hence $a^{n+m}/bc\in\m_V$, which is absurd.
\end{proof}

By the proposition we deduce that the associated $a$-adically separated ring $V/\bigcap_{n\geq 1}(a^n)$ is again a valuation ring (\ref{prop-composition1} (1)).

\begin{prop}\label{prop-sep}
Let $V$ be a valuation ring of arbitrary height, and $a\in\m_V\setminus\{0\}$.
Then the following conditions are equivalent$:$
\begin{itemize}
\item[{\rm (a)}] $V$ is $a$-adically separated, that is, $\bigcap_{n\geq 1}(a^n)=0;$ 
\item[{\rm (b)}] $V[\frac{1}{a}]$ is a field $($hence is the fractional field of $V)$.
\end{itemize}
\end{prop}

\begin{proof}
First let us show the implication (a) $\Rightarrow$ (b).
Let $K$ be the fractional field of $V$, and take $x\in K\setminus V$.
We have $x^{-1}\in\m_V$.
Suppose $x^{-1}$ does not divide $a^n$ for any $n\geq 1$. 
Then the all $a^n$ divide $x^{-1}$, which would imply the absurd $x^{-1}=0$.
Hence there exists $n\geq 1$ such that $(x^{-1})\supseteq (a^n)$, thereby $xa^n\in V$.

Next, we show (b) $\Rightarrow$ (a).
Set $J=\bigcap_{n\geq 1}(a^n)$. 
Suppose there exists a non-zero $b\in J$.
The element $b^{-1}$ can be written as $c/a^m$ for some $c\in V$ and $m\geq 0$.
Since $b\in J$, there exists $d\in V$ such that $b=a^{m+1}d$. 
Consequently, we have $a^m=bc=a^{m+1}cd$, which implies $acd=1$. 
But this contradicts $a\in\m_V$.
\end{proof}

An $a$-adically separated valuation ring has a remarkable property, which will be of particular importance in our later discussion:
\begin{prop}\label{prop-maxspe}
Suppose $V$ is an $a$-adically separated valuation ring for $a\in\m_V\setminus\{0\}$.
Then $V$ has a unique height one prime ideal$;$ explicitly, it is $\mathfrak{p}=\sqrt{(a)}$. 
\end{prop}

\begin{proof}
By \ref{prop-lemvalaval} $\mathfrak{p}=\sqrt{(a)}$ is a prime ideal.
To see it is the unique height one prime, it suffices to show that it is actually the `minimum' among all the non-zero prime ideals of $V$.
But, to show it, it suffices to show the following: for any non-zero $b\in V$ there exists $n\geq 1$ such that $a^n\in(b)$. 
Otherwise, we have $a^n/b\not\in V$ and hence $b\in(a^n)$ for any $n$. 
But this means, since $V$ is $a$-adically separated, that $b=0$, which is absurd. 
\end{proof}

By the proposition one can depict the spectrum of an $a$-adically separated valuation ring as in Figure \ref{fig-specvalcomplete}.
As the figure shows, if $V$ is $a$-adically separated, then the generic point $(0)$ has the `adjacent' specialization $\mathfrak{p}=\sqrt{(a)}$.
In particular, any $a$-adically separated valuation ring $V$ is the composite of the height one valuation ring $V_{\mathfrak{p}}$ and the valuation ring $V/\mathfrak{p}$.
\begin{figure}[htb]
\begin{center}
\setlength{\unitlength}{1pt}
\begin{picture}(130,40)(0,-5)
\put(5,5){\line(1,0){20}}
\put(25,4){\line(1,0){100}}
\put(25,6){\line(1,0){100}}
\put(5,5){\circle*{4}}
\put(25,5){\circle*{4}}
\put(125,5){\circle*{4}}
\put(-1,12){$\scriptstyle{(0)}$}
\put(14,12){$\scriptstyle{\sqrt{(a)}}$}
\put(119,12){$\scriptstyle{\mathfrak{m}_V}$}
\put(50,12){{\scriptsize other points}}
\end{picture}
\end{center}
\caption{$\Spec V$ for an $a$-adically separated valuation ring $V$}\label{fig-specvalcomplete}
\end{figure}

\begin{dfn}\label{dfn-maxspe2}{\rm 
Let $V$ be an $a$-adically separated valuation ring for $a\in\m_V\setminus\{0\}$.
We call the prime ideal $\mathfrak{p}=\sqrt{(a)}$ the {\em associated height one prime}\index{associated height one prime} of $V$ (or of the pair $(V,a)$). }
\end{dfn}

\begin{prop}\label{prop-maxspe3}
Let $V$ be an $a$-adically separated valuation ring for $a\in\m_V\setminus\{0\}$.
Then every non-zero prime ideal of $V$ is open with respect to $a$-adic topology.
\end{prop}

\begin{proof}
By \ref{prop-maxspe} every non-zero prime ideal $\mathfrak{q}$ of $V$ contains the associated height one prime $\mathfrak{p}=\sqrt{(a)}$, which contains $a$.
\end{proof}

\begin{prop}\label{prop-maxspe4}
Let $f\colon V\rightarrow W$ be a ring homomorphism between valuation rings, and $a\in\m_V\setminus\ker(f)$.
Suppose that $V$ is $a$-adically separated and that $W$ is $f(a)$-adically separated.
Then the map $f$ is injective, and the map $\Spec f\colon\Spec W\rightarrow\Spec V$ maps the set of all $a$-adically open prime ideals surjectively onto the set $[\sqrt{(a)},f^{-1}(\m_W)]$.
If, moreover, $f$ is local $($that is, $W$ dominates $V)$, then it maps the set of all non-zero prime ideals surjectively onto the set of all non-zero prime ideals.
\end{prop}

\begin{proof}
The map $f$ induces $V[\frac{1}{a}]\rightarrow W[\frac{1}{f(a)}]$; by the assumption this is the non-zero map between the fractional fields (\ref{prop-sep}) and hence is injective.
It follows therefore that $f$ is injective.
If we show $\sqrt{(f(a))}\cap V=\sqrt{(a)}$, then the other assertions follow from \ref{prop-schematic1} (3) and \ref{prop-schematic1} (4).
It suffices to show $\sqrt{(f(a))}\cap V\subseteq\sqrt{(a)}$.
Let $b\in\sqrt{(f(a))}\cap V$ and suppose $b\neq 0$; there exists $n\geq 1$ such that $b^n=ac$ with $c\in W$.
Since $c=b^n/a\in\Frac(V)$, $c$ or $c^{-1}$ lies in $V$.
If $c\in V$, then we have $b\in\sqrt{(a)}$.
Suppose $c^{-1}\in V$. 
Then $a=b^nc^{-1}\in\sqrt{(a)}$; since $\sqrt{(a)}$ is a prime ideal, $b^n$ or $c^{-1}$ lies in $\sqrt{(a)}$.
If $c^{-1}\in\sqrt{(a)}$, we have $c^{-1}\in\sqrt{(f(a))}\subseteq\m_W$, which is absurd.
Hence $b^n\in\sqrt{(a)}$, that is, $b\in\sqrt{(a)}$.
\end{proof}
\index{valuation!valuation ring@--- ring!a-adically separated valuation ring@$a$-adically separated --- ---|)}
\index{valuation!valuation ring@--- ring|)}

\addcontentsline{toc}{subsection}{Exercises}
\subsection*{Exercises}
\begin{exer}\label{exer-valuationdrill}{\rm 
Let $A$ be a ring, and $v\colon A\rightarrow\Gamma\cup\{\infty\}$ a valuation on $A$ with values in a totally ordered commutative group $\Gamma$.
Show that for $x,y\in A$ with $v(x)\neq v(y)$ we have $v(x+y)=\inf\{v(x),v(y)\}$.}
\end{exer}

\begin{exer}\label{exer-valtorsionflat}{\rm 
Let $V$ be a valuation ring, and $K=\Frac(A)$. 

(1) A $V$-module $M$ is flat if and only if it is torsion free. 

(2) If $Z\subset\P^1_K$ is finite over $K$, the closure $Z'$ in $\P^1_V$ is finite flat over $V$.

(3) For any finitely generated flat $V$-algebra $A$ such that $A\otimes_VK$ is finite over $K$, $A$ is quasi-finite and finitely presented over $V$.}
\end{exer}

\begin{exer}\label{exer-valflatatorfree}{\rm 
Let $V$ be an $a$-adically separated valuation ring for $a\in\m_V\setminus\{0\}$.
Show that a $V$-module $M$ is flat if and only if $M$ is $a$-torsion free, that is, for any $x\in M\setminus\{0\}$ and any $n\geq 0$ we have $a^nx\neq 0$.}
\end{exer}

\begin{exer}\label{exer-valuationdrill2}{\rm 
Let $V$ be a valuation ring such that $0<\mathrm{ht}(V)<+\infty$.

$(1)$ Show that there exists $a\in\m_V\setminus\{0\}$ such that $V$ is $a$-adically separated.

$(2)$ Show that for any non-maximal prime ideal $\mathfrak{p}\subseteq V$, there exists $b\in V$ such that $\mathfrak{p}=\bigcap_{n\geq 0}b^nV$.}
\end{exer}

\begin{exer}\label{exer-valuationdrill2ht1}{\rm 
Show that, if $V$ is a valuation ring of height one, then $V$ is $a$-adically separated for any $a\in\m_V\setminus\{0\}$.}
\end{exer}


\section{Topological rings and modules}\label{sec-topologicalringsmodules}
In this section we give a quick but reasonably detailed overview of the theory of {\em linearly topologized} rings and modules. 
First in \S\ref{sub-topologicalringsmodules} we recall the general definitions and treatments of linearly topologized rings and modules.
In \ref{subsub-completionfiltration}, we discuss {\em Hausdorff completion}\index{completion!Hausdorff completion@Hausdorff ---}.

From \S\ref{sub-adictopologyringsmodules} onward, we will be mainly interested in the so-called {\em adic topologies}.
One of the most important topics in this subsection is the notion of $I$-adic completions, which is at first defined by a ring-theoretic mapping universality.
The existence of the $I$-adic completion is a delicate matter, especially when treating non-Noetherian rings.
In fact, the $I$-adic completion exists under a mild condition, but in general, as we will see below, it may fail to exist. 

After briefly discussing henselian and Zariskian rings in \S\ref{sub-henselianzariskiantopologicalring}, we proceed in \S\ref{sub-adictop} to one of the most delicate issues, the preservation of adicness on passage to subspace topologies.
Classically, this is guaranteed in Noetherian situation by the well-known Artin-Rees lemma\index{Artin-Rees lemma}.
As it will be necessary for us to treat non-Noetherian rings, we consider analogous conditions, slightly weaker than the Artin-Rees condition, and prove many pratically useful results.

\subsection{Topology defined by a filtration}\label{sub-topologicalringsmodules}
\subsubsection{Filtrations}\label{subsub-topfromfil}
\index{filtration by submodules@filtration (by submodules)|(}
Let $A$ be a ring, and $M$ an $A$-module.
We consider a {\em descending filtration by $A$-submodules} $F^{\bullet}=\{F^{\lambda}\}_{\lambda\in\Lambda}$, that is, a collection of $A$-submodules $F^{\lambda}\subseteq M$ indexed by a directed set\index{set!directed set@directed ---}\index{directed set} $\Lambda$ such that for $\lambda,\mu\in\Lambda$
$$
\lambda\geq\mu\ \Longrightarrow\ F^{\lambda}\subseteq F^{\mu}.
$$

A descending filtration $F^{\bullet}=\{F^{\lambda}\}_{\lambda\in\Lambda}$ by $A$-submodules on an $A$-module $M$ is said to be {\em separated}\index{filtration by submodules@filtration (by submodules)!separated filtration by submodules@separated ---} (resp.\ {\em exhaustive}\index{filtration by submodules@filtration (by submodules)!exhaustive filtration by submodules@exhaustive ---}) if $\bigcap_{\lambda\in\Lambda}F^{\lambda}=\{0\}$ (resp.\ $\bigcup_{\lambda\in\Lambda}F^{\lambda}=M$).

Let $f\colon M\rightarrow L$ be a morphism of $A$-modules, and suppose that $M$ is endowed with a descending filtration by $A$-submodules $F^{\bullet}=\{F^{\lambda}\}_{\lambda\in\Lambda}$ indexed by a directed set $\Lambda$.
Then one has the descending filtration by $A$-submodules on $L$ given by $f(F^{\bullet})=\{f(F^{\lambda})\}_{\lambda\in\Lambda}$, which we call the {\em induced filtration}.
If $L$ is of the form $L=M/N$ by an $A$-submodule $N\subseteq M$ and $f$ is the canonical projection, then we have $f(F^{\bullet})=\{(N+F^{\lambda})/N\}_{\lambda\in\Lambda}$.

Similarly, for a morphism $g\colon N\rightarrow M$ of $A$-modules where $M$ is endowed with a descending filtration by $A$-submodules as above, one has the descending filtration by $A$-submodules on $N$ given by $g^{-1}(F^{\bullet})=\{g^{-1}(F^{\lambda})\}_{\lambda\in\Lambda}$, which we also call the {\em induced filtration}.
If $N$ is an $A$-submodule of $M$ and $g$ is the canonical inclusion, then we have $g^{-1}(F^{\bullet})=\{N\cap F^{\lambda}\}_{\lambda\in\Lambda}$.
\index{filtration by submodules@filtration (by submodules)|)}

\subsubsection{Topology defined by a filtration}\label{subsub-topologicalringsmodules}
\begin{prop}[{cf.\ \cite[Chap.\ III, \S1.2]{Bourb4}}]\label{prop-topologyfromfiltrationexist}
Let $A$ be a ring, and $M$ an $A$-module endowed with a descending filtration by $A$-submodules $F^{\bullet}=\{F^{\lambda}\}_{\lambda\in\Lambda}$ indexed by a directed set $\Lambda$.
Then there exists a unique topology on $M$ satisfying the following conditions$:$
\begin{itemize}
\item[{\rm (a)}] it is compatible with the additive group structure, that is, the addition $M\times M\rightarrow M$ is continuous$;$
\item[{\rm (b)}] $\{F^{\lambda}\}_{\lambda\in\Lambda}$ gives a fundamental system of open neighborhoods of $0\in M$.
\end{itemize}
Moreover, for any $a\in A$ the selfmap $x\mapsto ax$ of $M$ is a continuous endomorphism with respect to this topology.\hfill$\square$
\end{prop}

The topology on $M$ characterized as in \ref{prop-topologyfromfiltrationexist} is called the {\em topology defined by the filtration $F^{\bullet}$}\index{filtration by submodules@filtration (by submodules)!topology defined by a filtration by submodules@topology defined by a ---}\index{topology!topology defined by a filtration@--- defined by a filtration}.
It is explicitly described as follows: A subset $U\subseteq M$ is open if and only if for any $x\in U$ there exists $\lambda\in\Lambda$ such that $x+F^{\lambda}\subseteq U$.
In particular, an $A$-submodule $N\subseteq M$ is open if and only if it contains some $F^{\lambda}$.
In more formal terms, the topology thus obtained is the one defined by the uniform structure (cf.\ \cite[Chap.\ II, \S1.2]{Bourb4}) having $\{\til{F}^{\lambda}\}_{\lambda\in\Lambda}$ where $\til{F}^{\lambda}=\{(x,y)\in M\times M\,|\,x-y\in F^{\lambda}\}$ for $\lambda\in\Lambda$ as its fundamental system of entourages (cf.\ \cite[Chap.\ III, \S3.1]{Bourb4}); notice that, since $M$ is a commutative group by addition, the left and right uniformities coincide with each other.

The following facts are easy to verify, and we leave the proofs to the reader:
\begin{prop}\label{prop-topologyfromfiltrationproperty1}
Let $M$ be an $A$-module endowed with a descending filtration by $A$-submodules $F^{\bullet}=\{F^{\lambda}\}_{\lambda\in\Lambda}$ indexed by a directed set $\Lambda$.

{\rm (1)} Let $f\colon M\rightarrow L$ be a morphism of $A$-modules.
Then $f$ is continuous with respect to the topology on $M$ defined by $F^{\bullet}$ and the topology on $L$ defined by $f(F^{\bullet})$ $(${\rm \S\ref{subsub-topfromfil}}$)$. 
If, moreover, $f$ is surjective, then the topology on $L$ coincides with the quotient topology, that is, the strongest topology on $L$ such that the map $f$ is continuous.

{\rm (2)} Let $g\colon N\rightarrow M$ be a morphism of $A$-modules.
Then $g$ is continuous with respect to the topology on $N$ defined by $g^{-1}(F^{\bullet})$ $(${\rm \S\ref{subsub-topfromfil}}$)$ and the topology on $M$ defined by $F^{\bullet}$. 
If, moreover, $g$ is injective, then the topology on $N$ coincides with the subspace topology, that is, the weakest topology on $N$ such that the map $g$ is continuous. \hfill$\square$
\end{prop}

\begin{prop}\label{prop-topologyfromfiltrationproperty2}
Let $f\colon M\rightarrow N$ be a morphism of $A$-modules, and consider descending filtrations $\{F^{\lambda}\}_{\lambda\in\Lambda}$ and $\{G^{\sigma}\}_{\sigma\in\Sigma}$ by $A$-submodules on $M$ and $N$, respectively.
We consider the topologies on $M$ and $N$ defined by these filtrations.

{\rm (1)} The map $f$ is continuous if and only if for any $\sigma\in\Sigma$ there exists $\lambda\in\Lambda$ such that $f(F^{\lambda})\subseteq G^{\sigma}$.

{\rm (2)} The map $f$ is an open map if and only if for any $\lambda\in\Lambda$ there exists $\sigma\in\Sigma$ such that $G^{\sigma}\subseteq f(F^{\lambda})$. \hfill$\square$
\end{prop}

\begin{cor}\label{cor-topologyfromfiltrationproperty2}
Let $M$ be an $A$-module, and consider two descending filtrations by $A$-submodules $\{F^{\lambda}\}_{\lambda\in\Lambda}$ and $\{G^{\sigma}\}_{\sigma\in\Sigma}$.
Then these filtrations define the same topology on $M$ if and only if any $F^{\lambda}$ contains some $G^{\sigma}$ and any $G^{\sigma}$ contains some $F^{\lambda}$. \hfill$\square$
\end{cor}

\begin{prop}\label{prop-topologyfromfiltrationproperty3}
Let $M$ be an $A$-module endowed with a descending filtration by $A$-submodules $F^{\bullet}=\{F^{\lambda}\}_{\lambda\in\Lambda}$ indexed by a directed set $\Lambda$.
Then the topology on $M$ defined by $F^{\bullet}$ is Hausdorff if and only if the filtration $F^{\bullet}$ is separated\index{filtration by submodules@filtration (by submodules)!separated filtration by submodules@separated ---} $(${\rm \S\ref{subsub-topfromfil}}$)$. \hfill$\square$
\end{prop}

Let $M$ be a module over a ring $A$ endowed with a descending filtration of $A$-submodules $F^{\bullet}=\{F^{\lambda}\}_{\lambda\in\Lambda}$.
We leave the proof of the following (easy) proposition to the reader as an exercise (Exercise \ref{exer-filtrationtopclosure}):
\begin{prop}\label{prop-colsuresubmodule1}
Let $N\subseteq M$ be an $A$-submodule.
Then the closure $\ovl{N}$ of $N$ in $M$ with respect to the topology defined by the filtration $F^{\bullet}$ is given by
$$
\ovl{N}=\bigcap_{\lambda\in\Lambda}(N+F^{\lambda}).\eqno{\square}
$$
\end{prop}

\subsubsection{Hausdorff completion}\label{subsub-completionfiltration}
\index{completion!Hausdorff completion@Hausdorff ---|(}
Let us first briefly recall some generalities on {\em uniform spaces} (\cite[Chap.\ II]{Bourb4}).
A uniform space $X$ is said to be {\em Hausdorff complete}\index{complete!Hausdorff complete@Hausdorff ---} if the topology is Hausdorff and any Cauchy filter on $X$ is a convergent filter. 
It is known (\cite[Chap.\ II, \S3.7]{Bourb4}) that any uniform space has the {\em Hausdorff completion} $\widehat{X}$ and the canonical uniformly continuous mapping $i_X\colon X\rightarrow\widehat{X}$ in such a way that the pair $(\widehat{X},i_X)$ is uniquely characterized up to canonical isomorphisms by a universal mapping property with respect to uniformly continuous mappings to Hausdorff complete uniform spaces.
As a set, $\widehat{X}$ is the set of all minimal Cauchy filters on $X$.
In particular, $X$ is Hausdorff complete if and only if the canonical map $i_X$ is an isomorphism of uniform spaces.

Let $A$ be a ring, and $M$ an $A$-module equipped with a descending filtration by $A$-submodules $F^{\bullet}=\{F^{\lambda}\}_{\lambda\in\Lambda}$ indexed by a directed set $\Lambda$.
Since the topology on $M$ defined by the filtration $F^{\bullet}$ is a uniform topology, one can consider the Hausdorff completion, which we denote by 
$$
M^{\wedge}_{F^{\bullet}}.
$$
The Hausdorff completion $M^{\wedge}_{F^{\bullet}}$ is canonically a commutative group (\cite[Chap.\ III, \S3.5, Theorem 2]{Bourb4}); moreover, since any continuous additive endomorphisms of $M$, which are automatically uniformly continuous, lifts uniquely to that of $M^{\wedge}_{F^{\bullet}}$, one sees that $M^{\wedge}_{F^{\bullet}}$ is an $A$-module and that the canonical map $i_M\colon M\rightarrow M^{\wedge}_{F^{\bullet}}$ is an $A$-module homomorphism. 
Notice that the notion of Hausdorff completeness depends only on the {\em topologies} induced from uniform structures, and hence the Hausdorff completions $M^{\wedge}_{F^{\bullet}}$ depend up to isomorphisms only on the topologies defined by filtrations.

By \cite[Chap.\ III, \S7.3, Cor.\ 2]{Bourb4} we know that the Hausdorff completion $M^{\wedge}_{F^{\bullet}}$ is canonically identified (as an $A$-module) with the filtered projective limit\index{limit!projective limit@projective ---} 
$$
\varprojlim_{\lambda\in\Lambda}M/F^{\lambda}
$$
and the canonical map $i_M$ with the map induced from the canonical projections $M\rightarrow M/F^{\lambda}$ for $\lambda\in\Lambda$; the uniform structure on the projective limit $\varprojlim_{\lambda\in\Lambda}M/F^{\lambda}$ is the one induced from the descending filtration $\widehat{F}^{\bullet}=\{\widehat{F}^{\lambda}\}_{\lambda\in\Lambda}$ (called: the {\em induced filtration}\index{filtration by submodules@filtration (by submodules)!induced filtration by submodules@induced ---}) given by 
$$
\widehat{F}^{\lambda}=\ker(M^{\wedge}_{F^{\bullet}}\rightarrow M/F^{\lambda})
$$
for each $\lambda\in\Lambda$, where $M^{\wedge}_{F^{\bullet}}=\varprojlim_{\lambda\in\Lambda}M/F^{\lambda}\rightarrow M/F^{\lambda}$ is the canonical projection map.
In other words, it is the uniform structure by which the induced topology is the projective limit topology, where each $M/F^{\lambda}$ is considered with the discrete topology; more briefly, the topology on $M^{\wedge}_{F^{\bullet}}$ is the weakest one such that all projection maps $M^{\wedge}_{F^{\bullet}}\rightarrow M/F^{\lambda}$ is continuous.
(Note that, since the composition $M\stackrel{i_M}{\rightarrow}M^{\wedge}_{F^{\bullet}}\rightarrow M/F^{\lambda}$ is surjective, the canonical projection $M^{\wedge}_{F^{\bullet}}\rightarrow M/F^{\lambda}$ is surjective.)

Notice that, since $M^{\wedge}_{F^{\bullet}}/\widehat{F}^{\lambda}\cong M/F^{\lambda}$, we have $i_M(M)+\widehat{F}^{\lambda}=M^{\wedge}_{F^{\bullet}}$ for any $\lambda\in\Lambda$.
In particular, in view of \ref{prop-colsuresubmodule1}, the image $i_M(M)$ is dense in $M^{\wedge}_{F^{\bullet}}$.
More generally, we have $\widehat{F}^{\lambda}/\widehat{F}^{\mu}\cong F^{\lambda}/F^{\mu}$ for $\lambda\leq\mu$, and hence $i_M(F^{\lambda})+\widehat{F}^{\mu}=\widehat{F}^{\lambda}$; in particular, $\widehat{F}^{\lambda}$ coincides with the closure of $i_M(F^{\lambda})$ in $M^{\wedge}_{F^{\bullet}}$:
$$
\widehat{F}^{\lambda}=\ovl{i_M(F^{\lambda})}.
$$

\begin{prop}\label{prop-topologyfromfiltration31}
For any $\lambda\in\Lambda$ the submodule $\widehat{F}^{\lambda}$ is the closure of $i_M(F^{\lambda})$ in $M^{\wedge}_{F^{\bullet}}$ and coincides up to canonical isomorphisms with the Hausdorff completion of $F^{\lambda}$ with respect to the filtration induced from $F^{\bullet}$.
\end{prop}

\begin{proof}
The first statement has been already shown above.
For $\mu\geq\lambda$ we have the exact sequence
$$
0\longrightarrow F^{\lambda}/F^{\mu}\longrightarrow M/F^{\mu}\longrightarrow M/F^{\lambda}\longrightarrow 0.
$$
Since the subset $\{\mu\in\Lambda\,|\,\mu\geq\lambda\}$ of $\Lambda$ is cofinal, applying the projective limits $\varprojlim_{\mu}$, we have the exact sequence
$$
0\longrightarrow\varprojlim_{\mu}F^{\lambda}/F^{\mu}\longrightarrow M^{\wedge}_{F^{\bullet}}\longrightarrow M/F^{\lambda}\longrightarrow 0;
$$
notice that the passage to projective limits is a left-exact functor and that the canonical projection $M^{\wedge}_{F^{\bullet}}\rightarrow M/F^{\lambda}$ is surjective.
This shows that $\widehat{F}^{\lambda}$ coincides with $\varprojlim_{\mu}F^{\lambda}/F^{\mu}$, which is nothing but the Hausdorff completion of $F^{\lambda}$ with respect to the induced filtration.
\end{proof}

\begin{prop}\label{prop-topologyfromfiltration5}
Let $M$ be an $A$-module endowed with a descending filtration by $A$-submodules $F^{\bullet}=\{F^{\lambda}\}_{\lambda\in\Lambda}$ indexed by a directed set $\Lambda$, $M^{\wedge}_{F^{\bullet}}$ the associated Hausdorff completion, and $\widehat{F}^{\bullet}=\{\widehat{F}^{\lambda}\}_{\lambda\in\Lambda}$ the induced filtration on $M^{\wedge}_{F^{\bullet}}$.

{\rm (1)} For any $\lambda\in\Lambda$ we have $F^{\lambda}=i^{-1}_M(\widehat{F}^{\lambda})$, that is, the induced filtration $i^{-1}(\widehat{F}^{\bullet})$ coincides with the original one $F^{\bullet}$.

{\rm (2)} The $A$-module $M^{\wedge}_{F^{\bullet}}$ is Hausdorff complete with respect to the topology defined by the induced filtration $\widehat{F}^{\bullet}$.
\end{prop}

\begin{proof}
Since $M^{\wedge}_{F^{\bullet}}/\widehat{F}^{\lambda}=M/F^{\lambda}$ for $\lambda\in\Lambda$, we have 
$$
\varprojlim_{\lambda\in\Lambda}M^{\wedge}_{F^{\bullet}}/\widehat{F}^{\lambda}=\varprojlim_{\lambda\in\Lambda}M/F^{\lambda}=M^{\wedge}_{F^{\bullet}},
$$
which shows (2).
(1) follows easily from the commutative diagram 
$$
\xymatrix@-2ex{M\ar[r]^i\ar[d]&M^{\wedge}_{F^{\bullet}}\ar[d]\\ M/F^{\lambda}\ar@{=}[r]&M^{\wedge}_{F^{\bullet}}/\widehat{F}^{\lambda}\rlap{,}}
$$
for any $\lambda\in\Lambda$, where the vertical arrows are the canonical projections.
\end{proof}

One can easily show the following mapping universality by \cite[Chap.\ II, \S3.7, Theorem 3]{Bourb4}, aided by the fact that the product of Hausdorff complete uniform spaces is Hausdorff complete (\cite[Chap.\ II, \S3.5, Prop.\ 10]{Bourb4}) and the fact that the addition $(x,y)\mapsto x+y$, the inversion $x\mapsto -x$, and the scalar multiplication $x\mapsto ax$ (for $a\in A$) are uniformly continuous:
\begin{prop}[Mapping universality of Hausdorff completions]\label{prop-uniformgroupuniv1}
Let $M$ be an $A$-module equipped with a descending filtration $F^{\bullet}=\{F^{\lambda}\}_{\lambda\in\Lambda}$ indexed by a directed set $\Lambda$, and consider the canonical map $i_M\colon M\rightarrow M^{\wedge}_{F^{\bullet}}$.
Let $N$ be another $A$-module,  Hausdorff complete with respect to a descending filtration by $A$-submodules.
Then for any continuous $A$-module homomorphism $M\rightarrow N$, there exists a unique continuous $A$-module homomorphism $M^{\wedge}_{F^{\bullet}}\rightarrow N$ such that the resulting diagram
$$
\xymatrix@-1ex{M^{\wedge}_{F^{\bullet}}\ar@{-->}[dr]\\ M\ar[u]^{i_M}\ar[r]&N}
$$
commutes. \hfill$\square$
\end{prop}

\begin{rem}\label{rem-completionfiltration2}{\rm 
(1) An $A$-module $M$ linearly topologized by a descending filtration $F^{\bullet}=\{F^{\lambda}\}_{\lambda\in\Lambda}$ by $A$-submodules is separated if and only if the canonical map $i_M\colon M\rightarrow M^{\wedge}_{F^{\bullet}}$ is injective.

(2) If a ring $A$ is linearly topologized by a descending filtration $F^{\bullet}=\{F^{\lambda}\}_{\lambda\in\Lambda}$ by ideals, then the Hausdorff completion $A^{\wedge}_{F^{\bullet}}$ is canonically a ring, and the canonical map $i_A\colon A\rightarrow A^{\wedge}_{F^{\bullet}}$ is a ring homomorphism.}
\end{rem}
\index{completion!Hausdorff completion@Hausdorff ---|)}

\subsubsection{Hausdorff completion and exact sequence}\label{subsub-completionexactsequence}
Let us consider an exact sequence
$$
0\longrightarrow N\stackrel{g}{\longrightarrow}M\stackrel{f}{\longrightarrow}L\longrightarrow 0
$$
of $A$-modules and a descending filtration $F^{\bullet}=\{F^{\lambda}\}_{\lambda\in\Lambda}$ by $A$-submodules of $M$.
We consider the induced filtrations\index{filtration by submodules@filtration (by submodules)!induced filtration by submodules@induced ---} on $N$ and $L$ (cf.\ \S\ref{subsub-topfromfil}); the one on $N$ is given by $g^{-1}(F^{\bullet})=\{N\cap F^{\lambda}\}_{\lambda\in\Lambda}$ (where we regard $N$ as an $A$-submodule of $M$), and the one on $L$ is $f(F^{\bullet})=\{(N+F^{\lambda})/N\}_{\lambda\in\Lambda}$ (where $L$ is identified with $M/N$).
For each $\lambda\in\Lambda$ we have the induced exact sequence
$$
0\longrightarrow N/N\cap F^{\lambda}\longrightarrow M/F^{\lambda}\longrightarrow M/(N+F^{\lambda})\longrightarrow 0.
$$

Now we suppose:
\begin{itemize}
\item the directed set $\Lambda$ has a cofinal and at most countable subset.
\end{itemize}
Then by \ref{cor-ML4} the induced sequence of Hausdorff completions
$$
0\longrightarrow N^{\wedge}_{g^{-1}(F^{\bullet})}\longrightarrow M^{\wedge}_{F^{\bullet}}\longrightarrow L^{\wedge}_{f(F^{\bullet})}\longrightarrow 0
$$
is exact.

To proceed, let us set for brevity $G^{\lambda}=g^{-1}(F^{\lambda})=N\cap F^{\lambda}$ and $E^{\lambda}=f(F^{\lambda})$ for $\lambda\in\Lambda$.
Let $\widehat{G}^{\bullet}$ (resp.\ $\widehat{F}^{\bullet}$, resp.\ $\widehat{E}^{\bullet}$) be the induced filtration on the Hausdorff completion $N^{\wedge}_{G^{\bullet}}$ (resp.\ $M^{\wedge}_{F^{\bullet}}$, resp.\ $L^{\wedge}_{E^{\bullet}}$) (cf.\ \S\ref{subsub-completionfiltration}).

\begin{prop}\label{prop-completionexactsequence1}
The filtration $\widehat{G}^{\bullet}$ $($resp.\ $\widehat{E}^{\bullet})$ coincides with the one induced from $\widehat{F}^{\bullet}$ by the map $N^{\wedge}_{G^{\bullet}}\rightarrow M^{\wedge}_{F^{\bullet}}$ $($resp.\ $M^{\wedge}_{F^{\bullet}}\rightarrow L^{\wedge}_{E^{\bullet}})$ as in {\rm \S\ref{subsub-topfromfil}}.
\end{prop}

\begin{proof}
For each $\lambda\in\Lambda$ we have the commutative diagram with exact rows 
$$
\xymatrix{0\ar[r]&N^{\wedge}_{G^{\bullet}}\ar[d]\ar[r]&M^{\wedge}_{F^{\bullet}}\ar[d]\ar[r]&L^{\wedge}_{E^{\bullet}}\ar[d]\ar[r]&0\\ 0\ar[r]&N/G^{\lambda}\ar[r]&M/F^{\lambda}\ar[r]&L/E^{\lambda}\ar[r]&0\rlap{.}}
$$
The assertion on the filtration $\widehat{G}^{\lambda}$ follows easily.
One can also show the other assertion easily by diagram chasing, using the fact that the left-hand vertical arrow is surjective (\ref{prop-ML2}).
\end{proof}

\begin{prop}\label{prop-completionexactsequence2}
{\rm (1)} The image of $N^{\wedge}_{G^{\bullet}}\rightarrow M^{\wedge}_{F^{\bullet}}$ coincides with the closure of the image of $N$ by the canonical map $i_M\colon M\rightarrow M^{\wedge}_{F^{\bullet}}$.

{\rm (2)} The closure $\ovl{N}$ of $N$ in $M$ coincides with the pull-back of $N^{\wedge}_{G^{\bullet}}$ by the canonical map $i_M\colon M\rightarrow M^{\wedge}_{F^{\bullet}}$.
\end{prop}

\begin{proof}
(1) Consider the $A$-submodule $i_M(N)+\widehat{F}^{\lambda}$ for each $\lambda\in\Lambda$.
Since $i_M(N)+\widehat{F}^{\lambda}$ is the pull-back of the image of $N/N\cap F^{\lambda}$ by the canonical projection $M^{\wedge}_{F^{\bullet}}\rightarrow M/F^{\lambda}$, we have the exact sequence
$$
0\longrightarrow i_M(N)+\widehat{F}^{\lambda}\longrightarrow M^{\wedge}_{F^{\bullet}}\longrightarrow L/f(F^{\lambda})\longrightarrow 0.
$$
Applying the projective limits $\varprojlim_{\lambda\in\Lambda}$ (\ref{prop-projlimleftexact}), we have the exact sequence 
$$
0\longrightarrow\bigcap_{\lambda\in\Lambda}(i_M(N)+\widehat{F}^{\lambda})\longrightarrow M^{\wedge}_{F^{\bullet}}\longrightarrow L^{\wedge}_{E^{\bullet}}\longrightarrow 0.
$$
Hence the assertion follows from \ref{prop-colsuresubmodule1}.

(2) By \ref{prop-topologyfromfiltration5} (1) one sees easily that the equality 
$$
i^{-1}_M(i_M(N)+\widehat{F}^{\lambda})=N+F^{\lambda}
$$
holds for any $\lambda\in\Lambda$.
Then the assertion follows from (1) and \ref{prop-colsuresubmodule1}.
\end{proof}

\subsubsection{Completeness of sub and quotient modules}\label{subsub-completenesspermanence}
Consider as before the exact sequence 
$$
0\longrightarrow N\stackrel{g}{\longrightarrow}M\stackrel{f}{\longrightarrow}L\longrightarrow 0
$$
of $A$-modules, and a descending filtration $F^{\bullet}=\{F^{\lambda}\}_{\lambda\in\Lambda}$ by $A$-submodules of $M$, which yields the induced filtrations $G^{\bullet}=N\cap F^{\bullet}$ and $E^{\bullet}=f(F^{\bullet})$ on $N$ and $L$, respectively.
Furthermore, we continue with the assumption that the directed set $\Lambda$ contains a cofinal and at most countable subset.

\begin{prop}\label{prop-completenesspermanence1}
Suppose $M$ is Hausdorff complete\index{complete!Hausdorff complete@Hausdorff ---} with respect to the topology defined by the filtration $F^{\bullet}$.
Then the following conditions are equivalent$:$
\begin{itemize}
\item[{\rm (a)}] $N$ is closed in $M$ with respect to the topology defined by $F^{\bullet};$
\item[{\rm (b)}] $N$ is Hausdorff complete with respect to the topology defined by $G^{\bullet};$
\item[{\rm (c)}] $L$ is Hausdorff complete with respect to the topology defined by $E^{\bullet}$.
\end{itemize}
\end{prop}

\begin{proof}
Since $M$ is Hausdorff complete, we have the exact sequence
$$
0\longrightarrow N^{\wedge}_{G^{\bullet}}\longrightarrow M\longrightarrow L^{\wedge}_{E^{\bullet}}\longrightarrow 0.
$$
Suppose $N$ is closed in $M$.
Since the induced filtration $\widehat{F}^{\bullet}$ defined as in \S\ref{subsub-completionfiltration} coincides with $F^{\bullet}$, we deduce by \ref{prop-completionexactsequence2} (1) that $N=N^{\wedge}_{G^{\bullet}}$, which shows that $N$ is Hausdorff complete with respect to the topology defined by $G^{\bullet}$.
Then by the above exact sequence we conclude that $L=L^{\wedge}_{E^{\bullet}}$, that is, $L$ is Hausdorff complete with respect to the topology defined by $E^{\bullet}$.
Finally, if $L$ is complete with respect to the topology defined by $E^{\bullet}$, then again by the above exact sequence we know that $N=N^{\wedge}_{G^{\bullet}}$, which is closed in $M$ by \ref{prop-completionexactsequence2} (1).
\end{proof}

\subsection{Adic topology}\label{sub-adictopologyringsmodules}
\subsubsection{Adic filtration and adic topology}\label{subsub-adicfiltrationtopology}
\index{topology!adic topology@adic ---|(}\index{adic!adic topology@--- topology|(}
Let $A$ be a ring, $I\subseteq A$ an ideal, and $M$ an $A$-module.
One has the descending filtration on $M$ by $A$-submodules given by
$$
I^{\bullet}M=\{I^nM\}_{n\geq 0}.
$$
We call this filtration the {\em $I$-adic filtration}\index{filtration by submodules@filtration (by submodules)!I-adic filtration@$I$-adic ---} on $M$.
The $I$-adic filtration $\{I^n\}_{n\geq 0}$ in the case $M=A$ is simply denoted by $I^{\bullet}$.

Let $A$ be a ring equipped with the $I$-adic filtration by an ideal $I\subseteq A$.
Then an ideal $J\subseteq A$ is called an {\em ideal of definition}\index{ideal of definition} if there exist positive integers $m,n>0$ such that 
$$
I^m\subseteq J^n\subseteq I.
$$

The following proposition promptly follows from \ref{cor-topologyfromfiltrationproperty2}:
\begin{prop}\label{prop-adicfiltrationtopology2}
Let $A$ be a ring equipped with the $I$-adic filtration by an ideal $I\subseteq A$.

{\rm (1)} Let $J\subseteq A$ be an ideal of definition of $A$.
Then for any $A$-module $M$ the $I$-adic filtration $I^{\bullet}M$ and the $J$-adic filtration $J^{\bullet}M$ define the same topology.

{\rm (2)} Let $J\subseteq A$ be an ideal. If the filtrations $J^{\bullet}$ and $I^{\bullet}$ define the same topology on $A$, then $J$ is an ideal of definition of $A$. \hfill$\square$
\end{prop}

Let $A$ be a ring, and $M$ a linearly topologized $A$-module.
If the topology on $M$ is the same as the one defined by the $I$-adic filtration for an ideal $I$ of $A$, we say that the topology is an {\em adic topology}\index{topology!adic topology@adic ---}\index{adic!adic topology@--- topology}; if we like to spell out the ideal $I$, we say it is the {\em $I$-adic topology}.
Explicitly, a topology on $M$ defined by a descending filtration $F^{\bullet}=\{F^{\lambda}\}_{\lambda\in\Lambda}$ by $A$-submodules indexed by a directed set $\Lambda$ is $I$-adic if and only if for any $\lambda\in\Lambda$ there exists $n\geq 0$ such that $I^nM\subseteq F^{\lambda}$ and for any $n\geq 0$ there exists $\lambda\in\Lambda$ such that $F^{\lambda}\subseteq I^nM$ (\ref{cor-topologyfromfiltrationproperty2}).

\begin{rem}\label{rem-adicfiltrationtopology3}{\rm 
According to EGA terminology (cf.\ \cite[$\mathbf{0}_{\mathbf{I}}$, \S7]{EGA}), 
$$
\textrm{preadic}\ +\ \textrm{separated and complete}\ =\ \textrm{adic}.
$$
The terminology `preadic' is, however, not commonly used nowadays. 
Also in this book we avoid this terminology and suppose that `adic' does not imply `separated and complete', except for {\em adic ring}\index{adic!adic ring@--- ring}, which is already customarily considered to be separated and complete with respect to an adic topology; cf.\ {\bf \ref{ch-formal}}.\ref{dfn-admissibleringsadicrings} (2).}
\end{rem}

Let $A$ and $B$ be rings with adic topologies.
A ring homomorphism $f\colon A\rightarrow B$ is said to be {\em adic} if for some ideal of definition $I$ of $A$ the ideal $IB=f(I)B$ is an ideal of definition of $B$.
It is, in fact, easy to see that the condition for $f$ to be adic is equivalent to that $IB$ is an ideal of definition of $B$ for {\em any} ideal of definition $I$ of $A$.
Note that, whereas `adic' implies `continuous', the converse is not true. 
For example, any ring homomorphism $A\rightarrow B$ is continuous, if $A$ is equipped with the $0$-adic topology ($=$ discrete topology), but is not adic unless $(0)$ is an ideal of definition of $B$.

Let $f\colon A\rightarrow B$ and $g\colon B\rightarrow C$ be ring homomorphisms between adically topologized rings. 
Then one can easily see the following:
\begin{itemize}
\item if $f$ and $g$ are adic, then so is the composition $g\circ f$;
\item if $g\circ f$ and $f$ are adic, then so is $g$.
\end{itemize}

Let $A$ be a ring, $I\subseteq A$ an ideal, and $M$ an $A$-module considered with the $I$-adic topology. 
We say:
\begin{itemize}
\item $M$ is {\em $I$-adically separated}\index{separated!I-adically separated@$I$-adically ---} if $\bigcap_{n\geq 0}I^nM=\{0\}$ or, equivalently, the canonical map $i_M\colon M\rightarrow\varprojlim_{n\geq 0}M/I^nM$ is injective;
\item $M$ is {\em $I$-adically complete}\index{complete!I-adically complete@$I$-adically ---} if it is Hausdorff complete\index{complete!Hausdorff complete@Hausdorff ---} with respect to the $I$-adic topology, that is, the canonical map $i_M\colon M\rightarrow\varprojlim_{n\geq 0}M/I^nM$ is an isomorphism.
\end{itemize}
Notice that {\sl what we mean by `complete' with respect to $I$-adic topology is `{\em Hausdorff} complete'$;$ that is, in our convention `$I$-adically complete' implies `$I$-adically separated'.}

\begin{lem}\label{lem-bt1vis}
Let $A$ be a ring, $I\subseteq A$ an ideal, and $M$ an $A$-module considered with the $I$-adic topology.
Suppose there exists an integer $n\geq 1$ such that $I^nM=0$.
Then $M$ is $I$-adically complete. \hfill$\square$
\end{lem}

The following proposition will be frequently used later:
\begin{prop}[{\cite[Theorem 8.4]{Matsu}}]\label{prop-complpair1}
Let $A$ be a ring, $I\subseteq A$ an ideal, and $M$ an $A$-module.
Suppose that $A$ is $I$-adically complete and that $M$ is $I$-adically separated.
Let $\{x_1,\ldots,x_r\}\subseteq M$. 
If $\{(x_1\ \mathrm{mod}\ IM),\ldots,(x_r\ \mathrm{mod}\ IM)\}$ generates $M/IM$ over $A/I$, then $\{x_1,\ldots,x_r\}$ generates $M$ over $A$.
In particular, $M$ is finitely generated if and only if $M/IM$ is finitely generated over $A/I$. \hfill$\square$
\end{prop}

\begin{prop}\label{prop-exertopologicallynilpotentadic}
Let $A$ be a ring topologized by a descending filtration $F^{\bullet}=\{F^{\lambda}\}_{\lambda\in\Lambda}$ by ideals, and suppose it is Hausdorff complete.
Let $I\subseteq A$ be a {\em topologically nilpotent}\index{topologically nilpotent} ideal, that is, for any $\lambda\in\Lambda$ there exists $n\geq 0$ such that $I^n\subseteq F^{\lambda}$.
Then $A$ is $I$-adically complete if either one of the following conditions is satisfied$:$
\begin{itemize}
\item[{\rm (a)}] $I^n$ is closed in $A$ for any $n\geq 0;$
\item[{\rm (b)}] $I$ is finitely generated.
\end{itemize}
\end{prop}

\begin{proof}
(a) We have $I^n=\bigcap_{\lambda\in\Lambda}(I^n+F^{\lambda})$ for any $n\geq 0$.
Take for any $\lambda\in\Lambda$ an integer $N_{\lambda}\geq 0$ such that $I^n\subseteq F^{\lambda}$ for $n\geq N_{\lambda}$.
Then we calculate 
$$
\varprojlim_{n\geq 0}A/I^n=\varprojlim_{n\geq 0}\varprojlim_{\lambda\in\Lambda}A/(I^n+F^{\lambda})=\varprojlim_{n,\lambda}A/(I^n+F^{\lambda})=\varprojlim_{\lambda}\varprojlim_{n\geq N_{\lambda}}A/F^{\lambda}=A
$$
(up to canonical isomorphisms), where the second equality is due to (dual of) Exercise \ref{exer-doublecolimits}, and the third equality is obtained by replacing the index set by a cofinal one.

(b) Clearly, $A$ is $I$-adically separated.
Let $\{x_n\}_{n\geq 0}$ be a sequence in $A$ such that for $n\leq m$ we have $x_n\equiv x_m$ mod $I^n$.
There exists $x\in A$ such that $x\equiv x_n$ mod $F^{\lambda}$ for sufficiently large $\lambda$ (depending on $n$).
We want to show that $x$ is the limit of $\{x_n\}$ with respect to the $I$-adic topology.
Set $I=(a_1,\ldots,a_r)$, and let $k\geq 1$ be an arbitrary fixed positive integer.
Let $\xi_1,\ldots,\xi_s$ be the monomials in $a_i$'s of degree $k$.
Then one finds inductively $s$ elements $y_{l,j}\in A$ ($l\geq 0$, $j=1,\ldots,s$) such that 
\begin{itemize}
\item $x_{k+l}-x_k=\sum^s_{j=1}y_{l,j}\xi_j$;
\item $y_{l,j}\equiv y_{l',j}$ mod $I^l$ for $l'\geq l$ and $j=1,\ldots,s$.
\end{itemize}
Since $y_{l,j}$ for each $j$ converges in $A$ to an element $y_j$, we have $x=x_k+\sum^s_{j=1}y_j\xi_j$, which belongs to $I^k$, as desired.
\end{proof}

\danger{In connection to the above proposition, we remark here that there is a flaw in \cite[$\mathbf{0}_{\mathbf{I}}$, (7.2.4)]{EGA}, which has been repaired in the Springer version \cite[$\mathbf{0}$, (7.2.4)]{EGAInew}.}
\index{adic!adic topology@--- topology|)}\index{topology!adic topology@adic ---|)}

\subsubsection{$I$-adic completion}\label{subsub-Iadiccompletioncomplete}
\index{completion!I-adic completion@$I$-adic ---|(}
Our definition of $I$-adic completions is given by universal mapping property and, a priori, presented independently from the notion of Hausdorff completions:
\begin{dfn}\label{dfn-Iadiccompletiondelicate1}{\rm 
Let $A$ be a ring, $I\subseteq A$ an ideal, and $M$ an $A$-module.
An {\em $I$-adic completion}\index{completion!I-adic completion@$I$-adic ---} of $M$ is an $A$-module $\widehat{M}$ together with an $A$-module homomorphism $i_M\colon M\rightarrow\widehat{M}$ such that the following conditions are satisfied:
\begin{itemize}
\item[{\rm (a)}] $\widehat{M}$ is $I$-adically complete;
\item[{\rm (b)}] for any $I$-adically complete $A$-module $N$ and any $A$-module homomorphism $f\colon M\rightarrow N$, there exists uniquely an $A$-module homomorphism $\widehat{f}\colon\widehat{M}\rightarrow N$ such that the resulting diagram 
$$
\xymatrix@-1ex{\widehat{M}\ar@{-->}[dr]^{\widehat{f}}\\ M\ar[u]^{i_M}\ar[r]_f&N}
$$
commutes.
\end{itemize}}
\end{dfn}

It is clear by definition that an $I$-adic completion of $M$ is, if it exists, unique up to canonical isomorphisms.
\ddanger{As we will see in \ref{exa-Iadiccompletionifexists} below, the $I$-adic completion thus defined may fail to exist in general.
It exists, however, under a mild condition; see \ref{prop-Iadiccompletioncomplete2} and \ref{prop-Iadiccompletioncomplete1} below.}

Here the subtlety in the existence can be more clearly illustrated as follows.
Consider the Hausdorff completion of $M$ with respect to the $I$-adic topology, which we denote in what follows simply by\footnote{The notation $M^{\wedge}_{I^{\bullet}}$ is the abbreviation of $M^{\wedge}_{I^{\bullet}M}$, the Hausdorff completion of $M$ with respect to the $I$-adic topology on $M$.}
$$
M^{\wedge}_{I^{\bullet}}=\varprojlim_{n\geq 0}M/I^nM.
$$
Then, as we are soon going to see, it turns out that:
\begin{itemize}
\item if $M^{\wedge}_{I^{\bullet}}$ is $I$-adically complete, then this together with the canonical map $M\rightarrow M^{\wedge}_{I^{\bullet}}$ actually gives an $I$-adic completion of $M$;
\item conversely, if the $I$-adic completion $\widehat{M}$ of $M$ exists, then it is isomorphic to $M^{\wedge}_{I^{\bullet}}$.
\end{itemize}
Hence the clue for the existence of the $I$-adic completion lies in whether or not the Hausdorff completion with respect to the $I$-adic topology is $I$-adically complete or not.
The above dangerous-bend says that, actually, there exists such an example in which the Hausdorff completion with respect to the $I$-adic topology is {\em not} $I$-adically complete.

Suppose that the $I$-adic completion $i_M\colon M\rightarrow\widehat{M}$ of $M$ exists, and consider the canonical projections 
$$
p_n\colon M\longrightarrow M/I^nM,\quad q_n\colon\widehat{M}\longrightarrow\widehat{M}/I^n\widehat{M},
$$
for $n\geq 0$.
The $A$-module homomorphism $i_M\colon M\rightarrow\widehat{M}$ induces for each $n\geq 0$ a map $i_n\colon M/I^nM\rightarrow\widehat{M}/I^n\widehat{M}$ such that the diagram
$$
\xymatrix{M\ar[r]^{i_M}\ar[d]_{p_n}&\widehat{M}\ar[d]^{q_n}\\ M/I^nM\ar[r]_{i_n}&\widehat{M}/I^n\widehat{M}}
$$
commutes.
Notice that, since $M/I^nM$ is $I$-adically complete by \ref{lem-bt1vis}, we have the unique maps $\widehat{p}_n\colon\widehat{M}\rightarrow M/I^nM$ such that $\widehat{p}_n\circ i_M=p_n$ for $n\geq 0$, which then induce by passage to the projective limits the morphism
$$
p\colon\widehat{M}\longrightarrow M^{\wedge}_{I^{\bullet}}=\varprojlim_{n\geq 0}M/I^nM
$$
such that $p\circ i_M=p_n$.

\begin{prop}\label{prop-Iadiccompletionifexists}
The morphism $p$ is an isomorphism.
\end{prop}

Since $\widehat{M}$ is $I$-adically complete, the canonical morphism $\widehat{M}\rightarrow\varprojlim_{n\geq 0}\widehat{M}/I^n\widehat{M}$ is an isomorphism.
Hence the proposition follows from the following lemma:
\begin{lem}\label{lem-Iadiccompletiondelicate111}
For each $n\geq 0$ the map $i_n\colon M/I^nM\rightarrow\widehat{M}/I^n\widehat{M}$ is an isomorphism.
\end{lem}

\begin{proof}
The $A$-module homomorphism $\widehat{p}_n$ induces $j_n\colon\widehat{M}/I^n\widehat{M}\rightarrow M/I^nM$ such that $j_n\circ q_n=\widehat{p}_n$.
We want to show that the maps $i_n$ and $j_n$ are inverse to each other.

Consider first the composition $j_n\circ i_n$.
Since $j_n\circ i_n\circ\widehat{p}_n\circ i_M=j_n\circ i_n\circ p_n=j_n\circ q_n\circ i_M=\widehat{p}_n\circ i_M$, we have $j_n\circ i_n\circ\widehat{p}_n=\widehat{p}_n$ by the mapping universality of $I$-adic completions.
But since $\widehat{p}_n$ is surjective (as $p_n$ is surjective), we deduce $j_n\circ i_n=\id_{M/I^nM}$.

Next we discuss $i_n\circ j_n$.
We first notice that $\widehat{M}/I^n\widehat{M}$ is $I$-adically complete (by \ref{lem-bt1vis}), and hence the canonical projection $q_n$ is the unique map that makes the above diagram commute.
We calculate $i_n\circ j_n\circ q_n\circ i_M=i_n\circ\widehat{p}_n\circ i_M=i_n\circ p_n=q_n\circ i_M$ and hence have $i_n\circ j_n\circ q_n=q_n$.
Since $q_n$ is surjective, we deduce $i_n\circ j_n=\id_{\widehat{M}/I^n\widehat{M}}$, as desired.
\end{proof}

Let $A$ be a ring, $I\subseteq A$ an ideal, and $M$ an $A$-module.
Consider the Hausdorff completion $M^{\wedge}_{I^{\bullet}}$ of $M$ with respect to the $I$-adic topology.
Similarly, for each $k\geq 1$ we set $\widehat{F}^k=\varprojlim_{n\geq 0}I^kM/I^{k+n}M$, which is the Hausdorff completion of $I^kM$ with respect to the $I$-adic topology.
As in \ref{prop-topologyfromfiltration31}, $\{\widehat{F}^k\}_{k\geq 1}$ gives the descending filtration by $A$-submodules on $M^{\wedge}_{I^{\bullet}}$ induced from the $I$-adic filtration on $M$.
The following statement, which is a corollary of \ref{prop-Iadiccompletionifexists} (and \ref{lem-Iadiccompletiondelicate111}), gives an existence criterion of the $I$-adic completion of $M$:
\begin{cor}\label{cor-Iadiccompletionifexists1}
The following conditions are equivalent$:$
\begin{itemize}
\item[{\rm (a)}] an $I$-adic completion $i_M\colon M\rightarrow\widehat{M}$ exists$;$
\item[{\rm (b)}] the $A$-module $M^{\wedge}_{I^{\bullet}}$ is $I$-adically complete$;$
\item[{\rm (c)}] $I^kM^{\wedge}_{I^{\bullet}}$ is closed in $M^{\wedge}_{I^{\bullet}}$ with respect to the topology defined by the filtration $\{\widehat{F}^n\}_{n\geq 1}$ for any $k\geq 1$.
\item[{\rm (d)}] $\widehat{F}^k=I^kM^{\wedge}_{I^{\bullet}}$ for any $k\geq 1$.
\end{itemize}
\end{cor}

\begin{proof}
The implication (a) $\Rightarrow$ (b) follows from \ref{prop-Iadiccompletionifexists}.
Conversely, if (b) holds, then one can check that $M^{\wedge}_{I^{\bullet}}$ together with the canonical morphism $M\rightarrow M^{\wedge}_{I^{\bullet}}$ gives an $I$-adic completion by the mapping universality of projective limits, whence the implication (b) $\Rightarrow$ (a).

Next we show the equivalence (c) $\Leftrightarrow$ (d).
The implication (d) $\Rightarrow$ (c) follows immediately.
By \ref{prop-topologyfromfiltration31} $\widehat{F}^k$ is the closure of the image of $I^kM$ by the map $M\rightarrow M^{\wedge}_{I^{\bullet}}$; since $I^kM^{\wedge}_{I^{\bullet}}$ contains the image, and since $I^kM^{\wedge}_{I^{\bullet}}\subseteq\widehat{F}^k$, $\widehat{F}^k$ is the closure of $I^kM^{\wedge}_{I^{\bullet}}$.
Hence we have the other implication (c) $\Rightarrow$ (d).

Suppose (a) holds. 
Since $\widehat{F}^k$ is the kernel of the map $M^{\wedge}_{I^{\bullet}}\rightarrow M/I^kM$, we have $\widehat{F}^k=I^kM^{\wedge}_{I^{\bullet}}$ for any $k\geq 1$ by \ref{prop-Iadiccompletionifexists} (and \ref{lem-Iadiccompletiondelicate111}), whence the implication (a) $\Rightarrow$ (d).
Finally, if (d) holds, then $M^{\wedge}_{I^{\bullet}}=\varprojlim_{k\geq 1}M^{\wedge}_{I^{\bullet}}/\widehat{F}^k=\varprojlim_{k\geq 1}M^{\wedge}_{I^{\bullet}}/I^kM^{\wedge}_{I^{\bullet}}$, which shows that $M^{\wedge}_{I^{\bullet}}$ is $I$-adically complete, whence (d) $\Rightarrow$ (b).
\end{proof}

\begin{exa}\label{exa-Iadiccompletionifexists}{\rm 
Let $A=k[x_1,x_2,x_3,\ldots]$ be the polynomial ring of countably many indeterminacies over a field $k$, and $I=(x_1,x_2,x_3,\ldots)$ the ideal of $A$ generated by all the indeterminacies. 
We claim that the $I$-adic completion of $A$ does not exist.
Indeed, if it exists, then $B=A^{\wedge}_{I^{\bullet}}=\varprojlim_{n\geq 0}A/I^n$ would be $IB$-adically complete and hence $IB=\widehat{F}^1=\ker(B\rightarrow A/I)$; but one can show that the infinite sum $x_1+x^2_2+x^3_3+\cdots=\sum^{\infty}_{i=1}x^i_i$ does not converge with respect to the $IB$-adic topology as follows: Consider the partial sums $s_k=\sum^k_{i=1}x^i_i$ for $k\geq 1$ and suppose $\{s_k\}_{k\geq 1}$ has the limit $s\in B$ with respect to the $IB$-adic topology; since $s-s_k$ lies in $\widehat{F}^{k+1}=\ker(B\rightarrow A/I^{k+1})$ for each $k$, $s$ is also the limit of $\{s_k\}_{k\geq 1}$ with respect to the topology defined by the filtration $\{\widehat{F}^k\}$; then one sees $s\in\widehat{F}^1=IB$, which is absurd.}
\end{exa}

\subsubsection{Criterion for adicness}\label{subsub-criterionadicness}
It is in general difficult to determine whether a given topology defined by a filtration is adic or not. 
Let us state some criteria for adicness.
The following proposition is a rehash of some contents in \cite[$\mathbf{0}$, \S7.2]{EGAInew} and \cite[Chap.\ III, \S2.11]{Bourb1}:
\begin{prop}\label{prop-criterionadicness1}
Let $A$ be a ring endowed with a descending filtration by ideals $\{F^{(n)}\}_{n\geq 0}$ indexed by non-negative integers such that $F^{(0)}=A$.
We set $I=F^{(1)}$.
Suppose that the following conditions are satisfied$:$
\begin{itemize}
\item[{\rm (a)}] $A$ is Hausdorff complete with respect to the topology defined by the filtration $\{F^{(n)}\}_{n\geq 0};$ in other words, the canonical map 
$$
i\colon A\longrightarrow\varprojlim_{n\geq 0}A/F^{(n)}
$$
is an isomorphism$;$
\item[{\rm (b)}] for any $n>0$ the induced filtration on $A/F^{(n)}$ $($cf.\ {\rm \S\ref{subsub-topfromfil}}$)$ is $I$-adic$;$ in other words, we have
$$
I^m(A/F^{(n)})\ (=(F^{(1)}/F^{(n)})^m)\ =F^{(m)}/F^{(n)}
$$
for any $m$ and $n$ with $0\leq m\leq n$.
\item[{\rm (c)}] $F^{(1)}/F^{(2)}$ is finite generated as an ideal of $A/F^{(2)}$.
\end{itemize}
Then we have 
$$
F^{(n)}=I^n
$$
for any $n\geq 0$, and thus the filtration $\{F^{(n)}\}_{n\geq 0}$ is $I$-adic.
In particular, $A$ is $I$-adically complete, and $I$ is finitely generated. \hfill$\square$
\end{prop}

As a first implication of the proposition, we have the following useful:
\begin{prop}\label{prop-exeramazingfactonalgebra}
Let $B$ be a ring, and $J\subseteq B$ a finitely generated ideal.
We suppose that $B$ is $J$-adically complete, that is, $B=\varprojlim_{k\geq 0}B_k$, where $B_k=B/J^{k+1}$ for $k\geq 0$.
Let $\{A_k\}_{k\geq 0}$ be a projective system of $B$-algebras such that for $k\leq l$ the transition map $A_l\rightarrow A_k$ is surjective with the kernel equal to $J^{k+1}A_l$.
Let $A=\varprojlim_{k\geq 0}A_k$, and consider the descending filtration $\{F^{(n)}\}_{n\geq 0}$ by ideals given by $F^{(0)}=A$ and $F^{(n)}=\ker(A\rightarrow A_{n-1})$ for $n\geq 1$. Set $I=F^{(1)}$.
Then the ring $A$ is $I$-adically complete, and we have $I=JA$.
\end{prop}

To show this, we first need the following elementary fact:
\begin{lem}\label{lem-Iadiccompletioncomplete21}
Let $A$ be a ring, and $I\subseteq A$ an ideal, and suppose $A$ is $I$-adically complete.
Then we have $1+I\subseteq A^{\times}$. 
In particular, $I$ is contained in the Jacobson radical of $A$.
\end{lem}

\begin{proof}
The inverse of $1-a$ for $a\in I$ is given by $\sum^{\infty}_{n=0}a^n$, which belongs to $A$, for $A$ is $I$-adically complete.
\end{proof}

\begin{proof}[Proof of Proposition {\rm \ref{prop-exeramazingfactonalgebra}}]
First notice that each projection map $A\rightarrow A_{n-1}$ for $n\geq 1$ is surjective (\ref{prop-ML2}).
Notice also that we have $J^mA\subseteq F^{(m)}$ for any $m\geq 1$.
We have
$$
A_{n-1}/IA_{n-1}=A/I=A_0=A_{n-1}/JA_{n-1},
$$
which shows that $I=JA+F^{(n)}$.
This implies that $IA_{n-1}=JA_{n-1}$ and hence that
$$
I^mA_{n-1}=J^mA_{n-1}=F^{(m)}A_{n-1}=F^{(m)}/F^{(n)}
$$
for $0\leq m\leq n$.
Hence we can apply \ref{prop-criterionadicness1} to conclude that $I^n=F^{(n)}$ for $n\geq 1$ and that $A$ is $I$-adically complete.
Moreover, in view of \ref{lem-Iadiccompletioncomplete21} the equality $JA+I^2=I$ implies $I=JA$ by Nakayama's lemma.
\end{proof}

\begin{prop}\label{prop-criterionadicness2}
Let $A$ be a ring, $I\subseteq A$ a finitely generated ideal, and $M$ an $A$-module endowed with a descending filtration by $A$-submodules of the form $\{F^{(n)}\}_{n\geq 0}$ with $F^{(0)}=M$ such that $I^nM\subseteq F^{(n)}$ for any $n\geq 0$.
Suppose that the following conditions are satisfied$:$
\begin{itemize}
\item[{\rm (a)}] $A$ is $I$-adically complete, and $M$ is Hausdorff complete with respect to the topology defined by the filtration $F^{\bullet};$ in other words, the canonical map
$$
i\colon M\longrightarrow\varprojlim_{n\geq 0}M/F^{(n)}
$$
is an isomorphism$;$
\item[{\rm (b)}] for any $n>0$ the induced filtration on $M/F^{(n)}$ $($cf.\ {\rm \S\ref{subsub-topfromfil}}$)$ is $I$-adic$;$ in other words, we have
$$
I^m(M/F^{(n)})=F^{(m)}/F^{(n)}
$$
for any $0\leq m\leq n$.
\end{itemize}
Then we have
$$
F^{(n)}=I^nM
$$
for any $n\geq 0$, and thus the filtration $\{F^{(n)}\}_{n\geq 0}$ is $I$-adic.
In particular, $M$ is $I$-adically complete.
If moreover
\begin{itemize}
\item[{\rm (c)}] $M/F^{(1)}$ is finitely generated as an $A$-module
\end{itemize}
is satisfied, then $M$ is finitely generated as an $A$-module.
\end{prop}

\begin{proof}
Consider $B=A\oplus M$ and regard it as an $A$-algebra by the multiplication $(a,x)\cdot (a',x')=(aa',ax'+a'x)$ for $a,a'\in A$ and $x,x'\in M$; this is a ring with the square-zero ideal $M\subseteq B$.
Consider a descending filtration $\{J^{(n)}\}_{n\geq 0}$ by ideals given by $J^{(n)}=I^n\oplus F^{(n)}$ for each $n\geq 0$, and set $J=J^{(1)}$.
By the condition (a) the ring $B$ is Hausdorff complete with respect to the topology defined by this filtration, since, clearly, Hausdorff completion commutes with direct sums. 
Moreover, we have $I^nB\subseteq J^{(n)}$ for $n\geq 0$ by the assumption.
By the condition (b) with $m=1$ one has $I(M/F^{(n)})=F^{(1)}/F^{(n)}$ and hence $I^m(M/F^{(n)})=I^{m-1}(F^{(1)}/F^{(n)})$ for any $m\geq 1$.
By this one calculates
\begin{equation*}
\begin{split}
J^m(B/J^{(n)})=I^m/I^n\oplus I^{m-1}&(F^{(1)}/F^{(n)})=I^m(A/I^n\oplus M/F^{(n)})\\
&=I^m\oplus F^{(m)}/I^n\oplus F^{(n)}=J^{(m)}/J^{(n)}.
\end{split}
\end{equation*}
(Notice that the first line shows that $J^m(B/J^{(n)})=I^m(B/J^{(n)})$.)
Moreover, $J^{(1)}/J^{(2)}=I(A/I^2\oplus M/F^{(2)})=I(B/J^{(2)})$ is finitely generated.
Hence it follows from \ref{prop-criterionadicness1} that $J^{(n)}=J^n$ for $n\geq 1$, that $B$ is $J$-adically complete, and that $J\subseteq B$ is finitely generated.
Now, setting $B_k=B/J^{k+1}$ for $k\geq 0$, we see that the kernel of the surjective map $B_l\rightarrow B_k$ for $k\geq l$ is given by $J^{k+1}/J^{l+1}=I^{k+1}B_l$ and hence that by \ref{prop-exeramazingfactonalgebra} we have $J=IB$, and hence $J^n=I^nB$ for $n\geq 1$.
In particular, $I^nM$ is closed in $M$, and hence we have $F^{(n)}=I^nM$, as desired.
If $M/F^{(1)}=M/IM$ is finitely generated, $M$ is finitely generated due to \ref{prop-complpair1}.
\end{proof}

\subsubsection{Existence of $I$-adic completions}\label{subsub-Iadiccompletionsexistence}
\index{completion!I-adic completion@$I$-adic ---!existence of I-adic completions@existence of --- ---s|(}
As indicated in \ref{exa-Iadiccompletionifexists}, it is highly non-trivial whether or not the $I$-adic completion exists for a given adically topologized module. 
The following propositions shows that the existence is true if the ideal $I$ is finitely generated:
\begin{prop}\label{prop-Iadiccompletioncomplete2}
Let $A$ be a ring, and $I\subseteq A$ a {\em finitely generated} ideal.
Then the Hausdorff completion $A^{\wedge}_{I^{\bullet}}=\varprojlim_{n\geq 0}A/I^n$ of $A$ with respect to the $I$-adic topology is $I$-adic complete.
In particular, the $I$-adic completion of $A$ exists.
\end{prop}

\begin{prop}\label{prop-Iadiccompletioncomplete1}
Let $A$ be a ring, $I\subseteq A$ a finitely generated ideal, and $M$ an $A$-module.
Then the Hausdorff completion $M^{\wedge}_{I^{\bullet}}$ of $M$ with respect to the $I$-adic topology is $I$-adically complete.
In particular, the $I$-adic completion $\widehat{M}$ of $M$ exists.
If, moreover, $M/IM$ is finitely generated, then $\widehat{M}$ is finitely generated over the $I$-adic completion $\widehat{A}$ of $A$.
\end{prop}

\begin{proof}[Proof of Proposition {\rm \ref{prop-Iadiccompletioncomplete2}}]
Consider the Hausdorff completion $B=A^{\wedge}_{I^{\bullet}}$ of $A$ with respect to the $I$-adic topology, and set $F^{(k)}=\varprojlim_{n\geq 0}I^k/I^{k+n}$ for $k\geq 1$.
By \ref{cor-Iadiccompletionifexists1} it suffices to show the equality $F^{(k)}=I^kB$ for each $k\geq 0$.
To this end, we need to check $F^{(1)}=IB$ and to verify the conditions (a), (b), and (c) in \ref{prop-criterionadicness1} (with $A$ replaced by $B$).
The condition (a) is clear.
Since one has $F^{(m)}/F^{(n)}=I^m/I^n$ for $0\leq m\leq n$, the condition (b) is also verified.
Finally, by the assumption, $F^{(1)}/F^{(2)}=I/I^2$ is finitely generated as an ideal of $B/F^{(2)}=A/I^2$, and thus the condition (c) holds.
Thus we have shown that $B$ is $J$-adically complete, where $J=F^{(1)}$.
Moreover, by \ref{prop-complpair1} we deduce that $J$ is finitely generated, since $J/J^2=F^{(1)}/F^{(2)}=I/I^2$ is finitely generated as an ideal of $B/F^{(2)}=A/I^2$.

Next we show that the equality $J=IB$ holds.
Since the image of $I$ in $B$ is dense in $J$, we know that the finitely generated ideal $IB$ is dense in $J$.
In particular, we have $IB+J^2=J$.
Then the desired equality follows, in view of \ref{lem-Iadiccompletioncomplete21}, from Nakayama's lemma.
\end{proof}

\begin{proof}[Proof of Proposition {\rm \ref{prop-Iadiccompletioncomplete1}}]
By \ref{prop-Iadiccompletioncomplete2} we already know that the $I$-adic completion $\widehat{A}$ of $A$ exists.
Let $N=M^{\wedge}_{I^{\bullet}}$ be the Hausdorff completion of $M$ with respect to the $I$-adic topology.
We consider the filtration $\{F^{(k)}\}_{k\geq 1}$ on $N$ given by $F^{(k)}=\varprojlim_{n\geq 0}I^kM/I^{k+n}M$ for $k\geq 1$.
In view of \ref{cor-Iadiccompletionifexists1} we only need to show that $N$ is $I$-adically complete and, to this end, to check the conditions in \ref{prop-criterionadicness2} with $M$ replaced by $N$, $A$ by $\widehat{A}$, and $I$ by $I\widehat{A}$.
The condition (a) is clear.
Since $N/F^{(n)}=M/I^nM$ for $n\geq 0$ and $F^{(m)}/F^{(n)}=I^mM/I^nM$ for $0\leq m\leq n$, we also verify the condition (b).
\end{proof}
\index{completion!I-adic completion@$I$-adic ---!existence of I-adic completions@existence of --- ---s|)}
\index{completion!I-adic completion@$I$-adic ---|)}

\subsection{Henselian rings and Zariskian rings}\label{sub-henselianzariskiantopologicalring}
\subsubsection{Henselian rings}\label{subsub-henselianpairs}
\index{henselian!I-adically henselian@$I$-adically ---|(}
Recall that (cf.\ \cite{Cox1}\cite{Greco1}\cite{KPR}) a ring $A$ endowed with an $I$-adic topology ($I\subseteq A$) is said to be {\em henselian with respect to $I$} or {\em $I$-adically henselian}\index{henselian!I-adically henselian@$I$-adically ---} if it satisfies either one of the following equivalent conditions: 
\begin{itemize}
\item[(a)] for any \'etale morphism $X\rightarrow\Spec A$, any section $\sigma_0$ of the morphism $X_0\ (=X\times_{\Spec A}\Spec A_0)\rightarrow\Spec A_0$, where $A_0=A/I$, lifts to a section $\sigma$ of $X\rightarrow\Spec A$:
$$
\xymatrix{X_0\ar@{^{(}->}[r]\ar[d]&X\ar[d]\\ \Spec A_0\ar@{^{(}->}[r]\ar@/^1pc/[u]^{\sigma_0}&\Spec A;\ar@/^1pc/@{-->}[u]^{\sigma}}
$$
\item[(b)] the ideal $I$ is contained in the Jacobson radical of $A$, and for any monic polynomial $F(T)\in A[T]$ such that $F(0)\equiv 0$ mod $I$ and $F'(0)$ is invertible in $A/I$, there exists $a\in I$ such that $F(a)=0$ (that is, Hensel's lemma holds for monic polynomials).
\end{itemize}
(It is easy to see that the property `henselian' is actually a topological one.)
Similarly to completion, there is the notion of {\em henselization}\index{henselization}
$$
i_A\colon A\longrightarrow\het{A},
$$
characterized up to isomorphisms by a universal mapping property similar to that for $I$-adic completions (\ref{dfn-Iadiccompletiondelicate1}).
The henselization $\het{A}$ always exists, even if the ideal $I$ is not finitely generated.
Here is a rough sketch of the construction: The pair $(\het{A},I\het{A})$ is the inductive limit of all pairs $(B,IB)$ with $B$ an \'etale A-algebra such that $B/IB\cong A/I$; see \cite{Greco1} and \cite[\S2.8]{KPR} for other (equivalent) constructions.

\begin{prop}\label{prop-relpair31} 
The property `henselian' is preserved by filtered inductive limits.
More precisely, if $\{A_{\lambda}\}_{\lambda\in\Lambda}$ is an inductive system of rings with adic topologies with adic transition maps {\rm (\S\ref{subsub-adicfiltrationtopology})} and if all $A_{\lambda}$ are henselian, then the inductive limit $A=\varinjlim_{\lambda\in\Lambda}A_{\lambda}$ with the induced adic topology is henselian. \hfill$\square$
\end{prop}

The proof is easy and left to the reader (cf.\ \cite[3.6]{Greco1}).
Notice that the analogous statement for `adically complete' is not true.
\index{henselian!I-adically henselian@$I$-adically ---|)}

\subsubsection{Zariskian rings}\label{subsub-zariskianpairs}
\index{Zariskian!I-adically Zariskian@$I$-adically ---|(}
The following proposition is easy to see, and the proof is left to the reader:
\begin{prop}\label{prop-zariskipair1}
The following conditions for a ring $A$ and an ideal $I\subseteq A$ are equivalent to each other$:$
\begin{itemize}
\item[{\rm (a)}] for any $a\in I$ the element $1+a$ is invertible in $A$, that is, $1+I\subseteq A^{\times};$ 
\item[{\rm (b)}] an element $a\in A$ is invertible if and only if $a\ \mathrm{mod}\ I$ is invertible in $A/I;$ 
\item[{\rm (c)}] $I$ is contained in the Jacobson radical of $A$. \hfill$\square$
\end{itemize}
\end{prop}

A ring $A$ endowed with the $I$-adic topology defined by an ideal $I\subseteq A$ is said to be {\em Zariskian with respect to $I$} or {\em $I$-adically Zariskian}\index{Zariskian!I-adically Zariskian@$I$-adically ---} if it satisfies the equivalent conditions in \ref{prop-zariskipair1}.
We have already seen in \ref{lem-Iadiccompletioncomplete21} that any $I$-adically complete ring is $I$-adically Zariskian.
Notice that, due to the condition (c) in \ref{prop-zariskipair1}, `Zariskian' is a topological property.

For an arbitrary ring $A$ and an ideal $I\subseteq A$ one can construct the {\em associated Zariskian ring $\zat{A}$}\index{Zariskian!associated Zariskian@associated ---} simply by setting $\zat{A}=S^{-1}A$, where $S=1+I$ is the multiplicative subset consisting of all elements of the form $1+a$ for $a\in I$.
It is clear that this construction gives a unique solution to the universal mapping property similar to the universal mapping properties of completion and of henselization.

\begin{prop}\label{prop-relpair32}
Let $A$ be a ring, and $I\subseteq A$ an ideal.
Then the following conditions are equivalent$:$
\begin{itemize}
\item[{\rm (a)}] $A$ is $I$-adically Zariskian$;$ 
\item[{\rm (b)}] every maximal ideal of $A$ is open with respect to the $I$-adic topology, and for any $I$-adically open prime ideal $\mathfrak{p}\subseteq A$, the localization $A_{\mathfrak{p}}$ is $IA_{\mathfrak{p}}$-adically Zariskian$;$ 
\item[{\rm (c)}] for any maximal ideal $\mathfrak{m}\subseteq A$, $A_{\mathfrak{m}}$ is $IA_{\mathfrak{m}}$-adically Zariskian.
\end{itemize}
\end{prop}

\begin{proof}
If $A$ is $I$-adically Zariskian, then any maximal ideal $\mathfrak{m}$ contains $I$ and hence is open.
Let $\mathfrak{p}$ be an open prime ideal.
Since $\mathfrak{p}$ contains $I^n$ for some $n>0$, it contains $I$.
We need to show that any element of $A$ of the form $f+a$ where $f\not\in\mathfrak{p}$ and $a\in I$ is invertible in $A_{\mathfrak{p}}$.
Suppose it is not. 
Then $f+a\in\mathfrak{p}A_{\mathfrak{p}}$.
Since $f\not\in\mathfrak{p}A_{\mathfrak{p}}$, we have $a\not\in\mathfrak{p}A_{\mathfrak{p}}$.
But this is absurd, since $I\subseteq\mathfrak{p}$.
Hence we have the implication (a) $\Rightarrow$ (b).
The implication (b) $\Rightarrow$ (c) is trivial, and (c) $\Rightarrow$ (a) is easy to see.
\end{proof}

\begin{rem}{\rm 
The term `Zariskian' is coined from the already widespread term `Zariski ring' (\cite[Chap.\ VIII, \S4]{ZSII}).
The required condition $1+I\subseteq A^{\times}$ is, if $A$ is Noetherian, equivalent to several other conditions as in \cite[Chap.\ VIII, Theorem 9]{ZSII}.
When $A$ is not Noetherian, however, it is not necessarily equivalent to all of them; in fact, a slight consideration over the proof of \cite[Chap.\ VIII, Theorem 9]{ZSII} leads one to the question of validity of Artin-Rees lemma\index{Artin-Rees lemma}, which will be at the center of our later observation in \S\ref{sub-adictop}; cf.\ \ref{prop-zariskipair2} below.}
\end{rem}

\subsubsection{Interrelation of the conditions}\label{subsub-interrelationconditions}
\index{henselian!I-adically henselian@$I$-adically ---|(}
\begin{prop}\label{prop-relpair1}
Let $A$ be a ring endowed with an adic topology defined by an ideal $I\subseteq A$.

{\rm (1)} The following implications hold$:$
$$
\textrm{`complete'}\ \Longrightarrow\ \textrm{`henselian'}\ \Longrightarrow\ \textrm{`Zariskian'}.
$$

{\rm (2)} There exists a unique adic homomorphism {\rm (}cf.\ {\rm \S\ref{subsub-adicfiltrationtopology})} $\zat{A}\rightarrow\het{A}$ such that the diagram $\xymatrix@R-3.5ex@C-5.5ex{\zat{A}\ar[rr]&&\het{A}\\ &A\ar[ur]\ar[ul]}$ commutes$;$ if the $I$-adic completion $\widehat{A}$ of $A$ exists, then there exists a unique adic homomorphism $\het{A}\rightarrow\widehat{A}$ such that the diagram $\xymatrix@R-3.5ex@C-5.5ex{\het{A}\ar[rr]&&\widehat{A}\\ &A\ar[ur]\ar[ul]}$ commutes.
\end{prop}

\begin{proof}
By Hensel's lemma (cf.\ \cite[Chap.\ III, \S4.3, Theorem 1]{Bourb1}) we know that $I$-adically complete rings are $I$-adically henselian. 
It is clear by the definition (cf.\ \S\ref{subsub-henselianpairs}) that $I$-adically henselian rings are $I$-adically Zariskian.
This proves (1).
Then (2) follows immediately from the universal mapping properties.
\end{proof}

In particular, it follows that, for example, the $I\het{A}$-adic completion of $\het{A}$ coincides up to canonical isomorphism with $\widehat{A}$, and the $I\zat{A}$-adic henselization of $\zat{A}$ with $\het{A}$, etc.

In the rest of this subsection we collect some basic facts on the canonical maps $\zat{A}\rightarrow\het{A}$ and $\zat{A}\rightarrow\widehat{A}$, etc.\ concerning with flatness.

\begin{lem}\label{lem-propzariskianfaithfulltflat}
Let $A$ and $B$ be rings with adic topologies, $I\subseteq A$ an ideal of definition of $A$, and $f\colon A\rightarrow B$ an adic homomorphism.
Suppose that $A$ is $I$-adically Zariskian and that $B$ is $A$-flat.
Then the following conditions are equivalent$:$
\begin{itemize}
\item[{\rm (a)}] the induced morphism $A/I\rightarrow B/IB$ is faithfully flat\index{flatness!faithfully-flatness@faithfully-{---}}$;$
\item[{\rm (b)}] the morphism $f\colon A\rightarrow B$ is faithfully flat.
\end{itemize}
\end{lem}

\begin{proof}
The implication (b) $\Rightarrow$ (a) is clear.
To show the converse, let $N$ be a finitely generated $A$-module such that $N\otimes_AB=0$.
We have $N\otimes_AB/IB=0$, and hence $N/IN=N\otimes_AA/I=0$ due to (a).
Since $A$ is $I$-adically Zariskian, we deduce $N=0$ by Nakayama's lemma.
\end{proof}

\begin{prop}\label{prop-relpair311}
Let $A$ be a ring with an adic topology, and $\het{A}$ the henselization.
Then the map $A\rightarrow\het{A}$ is flat.
\end{prop}

\begin{proof}
As $\het{A}$ is isomorphic to the inductive limit of a system of rings \'etale over $A$, $\het{A}$ is flat over $A$.
\end{proof}

\danger{Note that the analogous statement for `completion' may fail to hold; the canonical map $A\rightarrow\widehat{A}$ may not be flat in general (cf.\ Exercise \ref{exer-counterexaflatness1}).}

\begin{prop}\label{prop-relpair21}
Let $A$ be a ring endowed with an adic topology defined by an ideal $I\subseteq A$.

{\rm (1)} The canonical map $\zat{A}\rightarrow\het{A}$ {\rm (\ref{prop-relpair1} (2))} is faithfully flat.

{\rm (2)} If the $I$-adic completion $\widehat{A}$ exists and if the canonical map $A\rightarrow\widehat{A}$ is flat, then the canonical map $\zat{A}\rightarrow\widehat{A}$ {\rm (\ref{prop-relpair1} (2))} is faithfully flat.
\end{prop}

\begin{proof}
(1) It follows from \ref{prop-relpair311} that the map $\zat{A}\rightarrow\het{A}$ is flat.
By the construction of the henselization we know that $\zat{A}/I\zat{A}\cong\het{A}/I\het{A}$.
Hence the assertion follows from \ref{lem-propzariskianfaithfulltflat}.

(2) Since $\zat{A}/I\zat{A}\cong\widehat{A}/\widehat{I}$, again the assertion follows immediately from \ref{lem-propzariskianfaithfulltflat}.
\end{proof}
\index{henselian!I-adically henselian@$I$-adically ---|)}
\index{Zariskian!I-adically Zariskian@$I$-adically ---|)}

\subsection{Preservation of adicness}\label{sub-adictop}
\subsubsection{General observation}\label{subsub-ARgeneral}
Let $A$ be a ring, and consider an exact sequence of $A$-modules:
$$
0\longrightarrow N\stackrel{g}{\longrightarrow}M\stackrel{f}{\longrightarrow}L\longrightarrow 0.\leqno{(\ast)}
$$
Let $I\subseteq A$ be an ideal of $A$, and consider the $I$-adic filtration $I^{\bullet}M=\{I^nM\}_{n\geq 0}$ on $M$.
This filtration induces as in \S\ref{subsub-topfromfil} the filtration $\{G^{(n)}=N\cap I^nM\}_{n\geq 0}$ on $N$ and the filtration $\{E^{(n)}=(N+I^nM)/N\}_{n\geq 0}$ on $L$.
It will be important in many places in the sequel to know whether or not these induced topologies are the $I$-adic ones.
In order to prepare for such situations, let us here give a general observation.

To this end, for a while, let us consider the following situation.
Let $A$ be a ring, $I\subseteq A$ an ideal, and $M$ an $A$-module equipped with a descending filtration by $A$-submodules of the form $F^{\bullet}=\{F^{(n)}\}_{n\in\Z}$ indexed by the directed set of all integers.
We consider the following conditions for the filtration:
\begin{itemize}
\item[{\bf (F1)}] $IF^{(m)}\subseteq F^{(m+1)}$ for any $m\in\Z$;
\item[{\bf (F2)}] there exist $p\geq 0$ and $q\in\Z$ such that $I^pM\subseteq F^{(q)}$.
\end{itemize}
Notice that the second condition is satisfied if, for example, there exists $q\in\Z$ such that $F^{(q)}=M$.

\begin{lem}\label{lem-generalizedAR1}
The topology on $M$ is $I$-adic if the filtration $F^{\bullet}$ satisfies {\bf (F1)}, {\bf (F2)}, and the following condition$:$
\begin{itemize}
\item[{\bf (F3)}] there exists $c\geq 0$ and $d\in\Z$ such that $F^{(n+d)}\subseteq I^nM$ for any $n\geq c$.
\end{itemize}
\end{lem}

\begin{proof}
First notice that the conditions {\bf (F1)} and {\bf (F2)} imply that $I^{m+p}M\subseteq F^{(m+q)}$ for any $m\geq 0$.
Set $e=\max\{c,p\}$.
Then we have the inclusions
$$
F^{(n+d)}\subseteq I^nM\subseteq F^{(n+q-p)},\quad I^{n+d+p}M\subseteq F^{(n+d+q)}\subseteq I^{n+q}M
$$
for any $n\geq e$.
Hence the assertion follows from \ref{cor-topologyfromfiltrationproperty2}.
\end{proof}

\begin{dfn}\label{dfn-Igood}{\rm 
Let $A$ be a ring, $I\subseteq A$ an ideal, and $M$ an $A$-module equipped with a descending filtration by $A$-submodules $F^{\bullet}=\{F^{(n)}\}_{n\in\Z}$.
We say that the filtration $F^{\bullet}$ is {\em $I$-good}\index{I-good@$I$-good (filtration)}\index{filtration by submodules@filtration (by submodules)!I-good filtration by submodules@$I$-good ---} if it satisfies the condition {\bf (F1)} and the following condition$:$
\begin{itemize}
\item[{\bf (F4)}] there exists an integer $N$ such that $IF^{(n)}=F^{(n+1)}$ for $n\geq N$.
\end{itemize}}
\end{dfn}

\begin{lem}\label{lem-generalizedAR2}
If $F^{\bullet}$ is $I$-good, then it satisfies {\bf (F3)}.
In particular, if $F^{\bullet}$ is $I$-good and satisfies {\bf (F2)}, then the topology on $M$ defined by the filtration $F^{\bullet}$ is $I$-adic.
\end{lem}

\begin{proof}
If $F^{\bullet}$ is $I$-good, we deduce by induction that $I^nF^{(d)}=F^{(n+d)}$ for any $d\geq N$ and $n\geq 0$.
Then we have $F^{(n+d)}=I^nF^{(d)}\subseteq I^nM$.
The last assertion follows from \ref{lem-generalizedAR1}.
\end{proof}

\subsubsection{$I$-adicness of quotient topologies}\label{subsub-ARIgoodquot}
Now we turn back to the exact sequence $(\ast)$ in \S\ref{subsub-ARgeneral}.
First we study the filtration $E^{\bullet}$ on the quotient module $L$ induced from the $I$-adic filtration $I^{\bullet}M$ on $M$.
In fact, one can readily show that the topology on $L$ is always $I$-adic, since it is completely elementary to check that $(N+I^nM)/N=I^n(M/N)$ for $n\geq 0$:
\begin{lem}\label{lem-ARIgoodquotfiltration1}
The induced filtration $E^{\bullet}$ on the quotient module $L=M/N$ coincides with the $I$-adic filtration $I^{\bullet}L$. \hfill$\square$
\end{lem}

\begin{prop}\label{prop-qconsistency1-1p}
Let $A$ be a ring, $I\subseteq A$ an ideal, and $f\colon M\rightarrow L$ a surjective morphism of $A$-modules.

{\rm (1)} The induced map $f^{\wedge}_{I^{\bullet}}\colon M^{\wedge}_{I^{\bullet}}\rightarrow L^{\wedge}_{I^{\bullet}}$ between the Hausdorff completions with respect to the $I$-adic topologies is again surjective.

{\rm (2)} If $N=\ker(f)$, then we have
$$
\ker(f^{\wedge}_{I^{\bullet}})=\varprojlim_{n\geq 0}N/(N\cap I^nM),
$$
which coincides with the closure $($with respect to the topology defined by the filtration $\{\widehat{F}^{(n)}=\ker(M^{\wedge}_{I^{\bullet}}\rightarrow M/I^nM)\}_{n\geq 0})$ of the image of $N$ in $M^{\wedge}_{I^{\bullet}}$ by the canonical map $M\rightarrow M^{\wedge}_{I^{\bullet}}$.
\end{prop}

\begin{proof}
Since the topology on $L$ by the induced filtration $E^{\bullet}$ is $I$-adic, we have a canonical isomorphism $L^{\wedge}_{E^{\bullet}}\cong L^{\wedge}_{I^{\bullet}}$ between the Hausdorff completion with respect to the topology defined by the filtration $E^{\bullet}$ and the Hausdorff completion with respect to the $I$-adic topology.
Hence we get as in \S\ref{subsub-completionexactsequence} the exact sequence
$$
0\longrightarrow N^{\wedge}_{G^{\bullet}}\longrightarrow M^{\wedge}_{I^{\bullet}}\longrightarrow L^{\wedge}_{I^{\bullet}}\longrightarrow 0.
$$
In particular, we have (1).
By \ref{prop-completionexactsequence2} (1) we also have (2).
\end{proof}

\begin{cor}\label{cor-qconsistency1111}
Let $A$ be an $I$-adically complete ring by an ideal $I\subseteq A$, and $M$ an $I$-adically complete $A$-module.
Let $N\subseteq M$ be an $A$-submodule of $M$.
Then the quotient $M/N$ is $I$-adically complete if and only if $N$ is closed in $M$ with respect to the $I$-adic topology. \hfill$\square$
\end{cor}

\begin{cor}\label{cor-qconsistency111}
Let $A$ be an $I$-adically complete ring by an ideal $I\subseteq A$.
Then any finitely generated $I$-adically separated $A$-module is $I$-adically complete.
\end{cor}

\begin{proof}
Let $M$ be a finitely generated $A$-module, and write $M\cong A^{\oplus m}/K$ for some $m>0$ and an $A$-submodule $K\subseteq A^{\oplus m}$.
By \ref{prop-qconsistency1-1p} we have $M^{\wedge}_{I^{\bullet}}\cong A^{\oplus m}/\ovl{K}$, where $\ovl{K}$ is the closure of $K$ in $A^{\oplus m}$ with respect to the $I$-adic topology; notice that $A^{\oplus m}$ is $I$-adically complete.
Hence the canonical map $M\rightarrow M^{\wedge}_{I^{\bullet}}$ is surjective.
If, furthermore, $M$ is $I$-adically separated, it is injective.
\end{proof}

Since any $A$-submodule of an $I$-adically separated $A$-module is $I$-adically separated, we have immediately the following:
\begin{cor}\label{cor-qconsistency113}
Let $A$ be an $I$-adically complete ring by an ideal $I\subseteq A$, and $M$ an $I$-adically separated $A$-module.
Then any finitely generated $A$-submodule $N\subseteq M$ is $I$-adically complete. \hfill$\square$
\end{cor}

\danger{Notice that the submodule $N$ in \ref{cor-qconsistency113} is not necessarily closed in $M$, whereas it is $I$-adically complete; in fact, the subspace topology on $N$ induced from the $I$-adic topology on $M$ may fail to coincide with the $I$-adic topology on $N$. Consequently, the quotient $M/N$ may not be $I$-adically complete.}

\subsubsection{$I$-adicness of subspace topologies}\label{subsub-ARIgoodnesssub}
Let $A$ be a ring, $I\subseteq A$ an ideal, $M$ an $A$-module, and $N\subseteq M$ an $A$-submodule of $M$.
We consider the $I$-adic topology on $M$.
It is usually a very delicate problem to determine whether the subspace topology on $N$ induced from the $I$-adic topology on $M$ or, what amounts to the same, the topology defined by the induced filtration $G^{\bullet}=\{N\cap I^nM\}_{n\geq 0}$ (cf.\ \ref{prop-topologyfromfiltrationproperty1} (2)) is $I$-adic or not.

Let us start with the following lemma, which follows immediately from \ref{cor-topologyfromfiltrationproperty2} and the obvious inclusions $I^mN\subseteq N\cap I^mM$ $(m\geq 0)$:
\begin{lem}\label{lem-ARIgoodnesssubfiltration1}
The topology on $N$ defined by the induced filtration $G^{\bullet}=\{N\cap I^nM\}_{n\geq 0}$ is $I$-adic if and only if the following condition is satisfied$:$
\begin{itemize}
\item[$(\ast)$] for any $n\geq 0$ there exists $m\geq 0$ such that $N\cap I^mM\subseteq I^nN$. \hfill$\square$
\end{itemize}
\end{lem}

\begin{prop}\label{prop-qconsistency1-2p-1}
Let $A$ be a ring, and $I\subseteq A$ an ideal.
Consider an exact sequence
$$
0\longrightarrow N\stackrel{g}{\longrightarrow}M\stackrel{f}{\longrightarrow}L\longrightarrow 0
$$
of $A$-modules.
If $N\subseteq M$ satisfies the condition $(\ast)$ in {\rm \ref{lem-ARIgoodnesssubfiltration1}}, then the sequence
$$
0\longrightarrow N^{\wedge}_{I^{\bullet}}\stackrel{g^{\wedge}_{I^{\bullet}}}{\longrightarrow}M^{\wedge}_{I^{\bullet}}\stackrel{f^{\wedge}_{I^{\bullet}}}{\longrightarrow}L^{\wedge}_{I^{\bullet}}\longrightarrow 0
$$
consisting of Hausdorff completions with respect to the $I$-adic topologies is exact.
\end{prop}

\begin{proof}
We have already obtained, in the proof of \ref{prop-qconsistency1-1p}, the exact sequence
$$
0\longrightarrow N^{\wedge}_{G^{\bullet}}\longrightarrow M^{\wedge}_{I^{\bullet}}\longrightarrow L^{\wedge}_{I^{\bullet}}\longrightarrow 0.
$$
Since the topology defined by the filtration $G^{\bullet}$ on $N$ is $I$-adic by \ref{lem-ARIgoodnesssubfiltration1}, we have $N^{\wedge}_{G^{\bullet}}\cong N^{\wedge}_{I^{\bullet}}$.
\end{proof}

It is easy to see that, in the situation as in \ref{lem-ARIgoodnesssubfiltration1}, the condition $(\ast)$ is satisfied if $N$ is an open submodule of $M$.
Indeed, if $I^sM\subseteq N$ $(s\geq 0)$, then for any $n\geq 0$ one has $N\cap I^{n+s}M=I^{n+s}M\subseteq I^nN$.
However, for a general $N$ the condition $(\ast)$ in \ref{lem-ARIgoodnesssubfiltration1} is not necessarily satisfied.
In this connection, it is useful to axiomatize some practical conditions on $I$-adicness of subspace topologies: Let $A$ be a ring, and $I\subseteq A$ an ideal;
\begin{itemize}
\item[{\bf (AP)}] any $A$-submodule $N\subseteq M$ of a finitely generated $A$-module $M$ satisfies the condition $(\ast)$ in \ref{lem-ARIgoodnesssubfiltration1};
\item[{\bf (APf)}] any {\em finitely generated} $A$-submodule $N\subseteq M$ of a finitely generated $A$-module $M$ satisfies the condition $(\ast)$ in \ref{lem-ARIgoodnesssubfiltration1}.
\end{itemize}

\begin{prop}\label{prop-qconsistency1-2p-2}
Let $A$ be a ring, and $I\subseteq A$ a finitely generated ideal.
Suppose that $A$ with the $I$-adic topology satisfies the condition {\bf (APf)}.
Then if 
$$
N\stackrel{g}{\longrightarrow}M\stackrel{f}{\longrightarrow}L
$$
is an exact sequence of finitely generated $A$-modules, the induced sequence of the $I$-adic completions\index{completion!I-adic completion@$I$-adic ---}
$$
\widehat{N}\stackrel{\widehat{g}}{\longrightarrow}\widehat{M}\stackrel{\widehat{f}}{\longrightarrow}\widehat{L}
$$
is exact.
Hence the $I$-adic completion functor $M\mapsto\widehat{M}$ is exact on the full subcategory of $\Mod_A$ consisting of finitely generated $A$-modules.
\end{prop}

\begin{proof}
Let $N_1$ (resp.\ $L_1$) be the image of the map $g\colon N\rightarrow M$ (resp.\ $f\colon M\rightarrow L$).
Then \ref{prop-qconsistency1-2p-1} implies that 
$$
0\longrightarrow\widehat{N}_1\longrightarrow\widehat{M}\longrightarrow\widehat{L}_1\longrightarrow 0
$$
is exact.
We need to show:
\begin{itemize}
\item the map $\widehat{N}\rightarrow\widehat{N}_1$ is surjective; 
\item the map $\widehat{L}_1\rightarrow\widehat{L}$ is injective.
\end{itemize}
The former assertion follows promptly from \ref{prop-qconsistency1-1p} (1).
For the latter we apply \ref{prop-qconsistency1-2p-1} to the exact sequence $0\rightarrow L_1\rightarrow L\rightarrow L/L_1\rightarrow 0$, where $L_1\subseteq L$ satisfies $(\ast)$ in \ref{lem-ARIgoodnesssubfiltration1}.
\end{proof}

\begin{prop}\label{prop-ARconseq1-1}
Let $A$ be a ring, $I\subseteq A$ an ideal, $M$ an $A$-module, and $N\subseteq M$ an $A$-submodule of $M$.
Suppose that the following conditions are satisfied$:$
\begin{itemize}
\item[{\rm (i)}] $M$ is $I$-adically complete$;$
\item[{\rm (ii)}] $N\subseteq M$ satisfies $(\ast)$ in {\rm \ref{lem-ARIgoodnesssubfiltration1}}.
\end{itemize}
Then the following conditions are equivalent$:$
\begin{itemize}
\item[{\rm (a)}] $N$ is closed in $M$ with respect to the $I$-adic topology$;$
\item[{\rm (b)}] $N$ is $I$-adically complete$;$
\item[{\rm (c)}] $M/N$ is $I$-adically complete.\hfill$\square$
\end{itemize}
\end{prop}

This follows immediately from \ref{prop-completenesspermanence1} in view of the fact that the filtration on $M/N$ induced from the $I$-adic filtration on $M$ is $I$-adic (due to \ref{lem-ARIgoodquotfiltration1}) and the induced topology on $N$ is $I$-adic (due to the assumption).
Notice that the assumption (ii) is automatic in either of the following cases:
\begin{itemize}
\item $M$ is finitely generated and $A$ satisfies {\bf (AP)};
\item $M$ and $N$ are finitely generated and $A$ satisfies {\bf (APf)};
\item $N$ is open in $M$ with respect to the $I$-adic topology.
\end{itemize}

\begin{prop}\label{prop-corpropARconseq1-1-1}
Let $A$ be a ring, and $I\subseteq A$ a finitely generated ideal.
Consider a finitely generated $A$-module $M$ and an $A$-submodule $N\subseteq M$.
Suppose that either one of the following conditions is satisfied$:$
\begin{itemize}
\item $A$ with the $I$-adic topology and the $I$-adic completion $\widehat{A}$ with the $I\widehat{A}$-adic topology satisfy the condition {\bf (AP)}$;$
\item $N$ is open in $M$ with respect to the $I$-adic topology.
\end{itemize}
Then the $I$-adic completion $\widehat{N}$ coincides with the closure of the image of $N$ in $\widehat{M}$. 
Moreover, we have the exact sequence
$$
0\longrightarrow\widehat{N}\longrightarrow\widehat{M}\longrightarrow\widehat{M/N}\longrightarrow 0.
$$
\end{prop}

\begin{proof}
By the assumption the $I$-adic completion $\widehat{N}$ coincides with the Hausdorff completion of $N$ with respect to the induced filtration $\{N\cap I^nM\}_{n\geq 0}$.
Hence the canonical map $\widehat{N}\rightarrow\widehat{M}$ is injective, and by \ref{prop-completionexactsequence2} the image coincides with the closure of the image of $N$ in $\widehat{M}$.
We obtain the desired exact sequence due to the observation in \S\ref{subsub-completionexactsequence} and \ref{lem-ARIgoodquotfiltration1}.
\end{proof}

Finally, let us mention here that any adically topologized {\em Noetherian} ring $A$ satisfies the condition {\bf (AP)} (which amounts to the same as {\bf (APf)}):
\begin{prop}\label{prop-APforNoetherian}
Let $A$ be a Noetherian ring endowed with the $I$-adic topology defined by an ideal $I\subseteq A$.
Then $A$ satisfies {\bf (AP)}.\hfill$\square$
\end{prop}

The proposition is classical, verified by Artin-Rees lemma, which we will mention below (\ref{prop-ARclassic}).

\subsubsection{Useful consequences of the conditions}\label{subsub-ARconsequence}
\begin{prop}\label{prop-ARconseq0}
Let $A$ be a ring, and $I\subseteq A$ a finitely generated ideal.
Suppose that the ring $A$ together with the $I$-adic topology satisfy the condition {\bf (APf)} in {\rm \S\ref{subsub-ARIgoodnesssub}}.
Then for any {\em finitely presented} $A$-module $M$ the canonical morphism
$$
M\otimes_A\widehat{A}\longrightarrow\widehat{M}
$$
is an isomorphism.
\end{prop}

\begin{proof}
Take a finite presentation $A^{\oplus q}\rightarrow A^{\oplus p}\rightarrow M\rightarrow 0$, and consider the commutative diagram with exact rows
$$
\xymatrix{A^{\oplus q}\otimes_A\widehat{A}\ar[d]\ar[r]&A^{\oplus p}\otimes_A\widehat{A}\ar[d]\ar[r]&M\otimes_A\widehat{A}\ar[d]\ar[r]&0\\ \widehat{A}^{\oplus q}\ar[r]&\widehat{A}^{\oplus p}\ar[r]&\widehat{M}\ar[r]&0\rlap{.}}
$$
Here the exactness of the second row is a consequence of {\bf (APf)} via \ref{prop-qconsistency1-2p-2}.
Since the first two vertical arrows are clearly isomorphisms, the rest is also an isomorphism.
\end{proof}

\begin{prop}\label{prop-zariskipair2}
Let $A$ be a ring with the $I$-adic topology defined by an ideal $I\subseteq A$ satisfying the condition {\bf (APf)} in {\rm \S\ref{subsub-ARIgoodnesssub}}.
Then the following conditions are equivalent$:$
\begin{itemize}
\item[{\rm (a)}] $A$ is $I$-adically Zariskian$;$
\item[{\rm (b)}] any finitely generated $A$-module is $I$-adically separated$;$
\item[{\rm (c)}] any $A$-submodule of a finitely generated $A$-module $M$ is closed in $M$ with respect to the $I$-adic topology$;$
\item[{\rm (d)}] any maximal ideal of $A$ is closed with respect to the $I$-adic topology.
\end{itemize}
\end{prop}

\begin{proof}
The proposition can be shown similarly to \cite[Chap.\ III, \S3.3, Prop.\ 6]{Bourb1}; we include the proof for the reader's convenience.

First we prove (a) $\Rightarrow$ (b).
Let $M$ be a finitely generated $A$-module, and take $x\in\bigcap_{n\geq 0}I^nM$. 
We want to show that $x=0$.
Consider the submodule $Ax$ of $M$.
By virtue of {\bf (APf)} the topology on $Ax$ induced from $M$ is the $I$-adic topology, which is, however, the coarsest topology on $Ax$, for $Ax$ is contained in all neighborhoods of $M$.
Hence we have $Ix=Ax$, which means $(1-a)x=0$ for some $a\in I$.
Since $1+I\subseteq A^{\times}$, we have $x=0$, as desired.

Next we show (b) $\Rightarrow$ (c).
Let $N$ be an $A$-submodule of $M$. 
Then $M/N$ is $I$-adically separated, whence (c) (cf.\ Exercise \ref{exer-filtrationtopclosure}).

As the implication (c) $\Rightarrow$ (d) is clear, it only remains to show (d) $\Rightarrow$ (a).
Let $\mathfrak{m}$ be a maximal ideal of $A$. 
We need to show that $\mathfrak{m}$ contains $I$.
Consider the field $A/\mathfrak{m}$, which is $I$-adically separated, that is, $I(A/\mathfrak{m})\neq A/\mathfrak{m}$.
Hence we have $I(A/\mathfrak{m})=0$, that is, $I\subseteq\mathfrak{m}$.
\end{proof}

\begin{cor}\label{cor-qconsistency11}
Let $A$ be a ring, and $I\subseteq A$ an ideal.
Suppose $A$ is $I$-adically complete and satisfies the condition {\bf (APf)}.
Then any finitely generated $A$-module $M$ is $I$-adically complete.
\end{cor}

\begin{proof}
Since $A$ is $I$-adically Zariskian, by \ref{prop-zariskipair2} any finitely generated $A$-module is $I$-adically separated.
Then the corollary follows from \ref{cor-qconsistency111}.
\end{proof}

By this, \ref{prop-completenesspermanence1}, and \ref{prop-ARconseq1-1}, we have:
\begin{cor}\label{cor-propARconseq1-2}
Let $A$ be an $I$-adically complete ring by an ideal $I\subseteq A$, and $M$ a finitely generated $A$-module.
Suppose $A$ with the $I$-adic topology satisfies {\bf (APf)}.
Then any $A$-submodule $N\subseteq M$ is closed in $M$ with respect to the $I$-adic topology.
If, moreover, $A$ satisfies {\bf (AP)}, then $N$ is $I$-adically complete. \hfill$\square$
\end{cor}

Here is another useful corollary of \ref{prop-zariskipair2}:
\begin{cor}\label{cor-ARconseq2}
Let $A$ be a ring, and $I\subseteq A$ an ideal.
Suppose $A$ is $I$-adically Zariskian and satisfies {\bf (AP)}.
Then the following conditions for a finitely generated $A$-module $M$ are equivalent$:$
\begin{itemize}
\item[{\rm (a)}] $M$ is finitely presented over $A;$
\item[{\rm (b)}] $M/I^nM$ is finitely presented over $A/I^n$ for any $n\geq 1$.
\end{itemize}
\end{cor}

\begin{proof}
The implication (a) $\Rightarrow$ (b) is obvious.
Suppose (b) holds, and write $M\cong A^{\oplus m}/K$.
We need to show that $K$ is finitely generated.
We have for any $n>0$ the exact sequence
$$
0\longrightarrow K/K\cap I^nA^{\oplus m}\longrightarrow (A/I^n)^{\oplus m}\longrightarrow M/I^nM\longrightarrow 0.
$$
Take an $n>0$ such that $K\cap I^nA^{\oplus m}\subseteq IK$ (here we use the condition {\bf (AP)}).
Since $K/K\cap I^nA^{\oplus m}$ is finitely generated (cf.\ \cite[Chap.\ I, \S2.8, Lemma 9]{Bourb1}), we deduce that $K/IK$ is finitely generated.

Take a finitely generated $A$-submodule $K'\subseteq K$ that is mapped surjectively onto $K/IK$.
Then we have $K/K'=I(K/K')$.
On the other hand, since $K/K'$ is a finitely generated $A$-submodule of $A^{\oplus m}/K'$, it is $I$-adically separated by \ref{prop-zariskipair2}.
Hence we have $K/K'=0$, that is, $K=K'$ is finitely generated.
\end{proof}

Finally, let us summarize some of the important consequences of the conditions {\bf (AP)} and {\bf (APf)}.
Let $A$ be a ring, and $I\subseteq A$ an ideal.
Suppose $A$ is $I$-adically complete.
\begin{itemize}
\item The condition {\bf (APf)} for $A$ implies:
\begin{itemize}
\item[{\rm (a)}] any finitely generated $A$-module is $I$-adically complete;
\item[{\rm (b)}] any $A$-submodule $N$ of a finitely generated $A$-module $M$ is closed in $M$.
\end{itemize}
\item If, moreover, $A$ satisfies the condition {\bf (AP)}, then:
\begin{itemize}
\item[{\rm (c)}] any $A$-submodule $N$ of a finitely generated $A$-module $M$ is $I$-adically complete.
\end{itemize}
\end{itemize}

\subsection{Rees algebra and $I$-goodness}\label{sub-Reescone}
Let $A$ be a ring, and $I\subseteq A$ an ideal.
In this situation, the associated {\em Rees algebra}\index{algebra!Rees algebra@Rees ---}\index{Rees algebra} is the graded algebra
$$
R(A,I)=\bigoplus_{n\geq 0}I^n
$$
(where $I^0=A$) over $A$.
Clearly, if the ideal $I$ is finitely generated, then $R(A,I)$ is an $A$-algebra of finite type.
Let $M$ be an $A$-module equipped with a descending filtration $F^{\bullet}=\{F^{(n)}\}_{n\in\Z}$ by $A$-submodules.
Suppose that the condition {\bf (F1)} in \S\ref{subsub-ARgeneral} is satisfied, that is, 
$$
I^qF^{(n)}\subseteq F^{(n+q)}
$$
holds for any $q\geq 0$ and $n\in\Z$.
Then for any $k\in\Z$ 
$$
N_{\geq k}=\bigoplus_{n\geq k}F^{(n)}
$$
is a graded $R(A,I)$-module.

\begin{prop}\label{prop-Reescone1}
Suppose in the above situation that the ideal $I$ is finitely generated.
If $N_{\geq k}$ is finitely generated as an $R(A,I)$-module for some $k\in\Z$, then the filtration $F^{\bullet}$ is $I$-good\index{I-good@$I$-good (filtration)}\index{filtration by submodules@filtration (by submodules)!I-good filtration by submodules@$I$-good ---} {\rm (\ref{dfn-Igood})}.
\end{prop}

\begin{proof}
By shift of indices we may assume that $k=0$; write $N=N_{\geq 0}$.
Let $\{l_1,\ldots,l_s\}$ be a set of generators of $N$ consisting of homogenous elements, and set $d_i=\deg(l_i)$ for $1\leq i\leq s$.
Let $c=\max\{d_1,\ldots,d_s\}$, and let $a_1,\ldots,a_r$ generate $I$.
Suppose $n>c$.
Any $x\in F^{(n)}$ is written as
$$
x=\sum^s_{i=1}f_i(a_1,\ldots,a_r)l_i,
$$
where $f_i$ for each $i$ is a homogeneous polynomial of degree $n-d_i$.
If $h(a_1,\ldots,a_r)$ is a monomial of degree $n-d_i$, divide $h=h_1h_2$ with $\deg(h_1)=n-c$ and $\deg(h_2)=c-d_i$.
Since $c-d_i\geq e$, we have $h_2(a)l_i\in F^{(c)}$. 
Hence we have $h(a)l_i\in I^{n-c}F^{(c)}$ and thus have the inclusion $F^{(n)}\subseteq I^{n-c}F^{(c)}$, whence the equality $F^{(n)}=I^{n-c}F^{(c)}$.
By this one verifies the condition {\bf (F4)} in \S\ref{subsub-ARgeneral}.
\end{proof}

It is sometimes useful to consider the so-called {\em conormal cone}\index{conormal!conormal cone@--- cone} associated to $R(A,I)$:
$$
\gr^{\bullet}_I(A)=\bigoplus_{n\geq 0}I^n/I^{n+1}.
$$
If $F^{\bullet}$ is a descending filtration as above, then 
$$
\gr^{\bullet}_F(M)=\bigoplus_{n\in\Z}F^{(n)}/F^{(n+1)}
$$
is a graded $\gr^{\bullet}_I(A)$-module.

\begin{prop}[{cf.\ \cite[Chap.\ III, \S3.1, Prop.\ 3]{Bourb1}}]\label{prop-complpair2}
Suppose that the ideal $I$ is finitely generated, that $A$ is $I$-adically complete, and that the filtration $F^{\bullet}$ is separated and exhaustive $($cf.\ {\rm \S\ref{subsub-topfromfil}}$)$.
Then the following conditions are equivalent$:$
\begin{itemize}
\item[{\rm (a)}] $\gr^{\bullet}_F(M)$ is a finitely generated $\gr^{\bullet}_I(A)$-module$;$
\item[{\rm (b)}] $M$ and $F^{(n)}$ for all $n\in\Z$ are finitely generated as $A$-modules, and the filtration $F^{\bullet}$ is $I$-good\index{I-good@$I$-good (filtration)}\index{filtration by submodules@filtration (by submodules)!I-good filtration by submodules@$I$-good ---}.
\end{itemize}
\end{prop}

\begin{proof}
We want to prove (a) $\Rightarrow$ (b) (the converse (b) $\Rightarrow$ (a) is easy (cf.\ \cite[Chap.\ III, \S3.1, Prop.\ 3]{Bourb1})).
Suppose $\gr^{\bullet}_F(M)$ is finitely generated over $\gr^{\bullet}_I(A)$.
Since $\gr^{\bullet}_I(A)=\bigoplus_{m\geq 0}I^m/I^{m+1}$ is positively graded, any homogeneous element $x$ of $\gr^{\bullet}_F(M)$ is written as a linear combination of the generators of degree less than or equal to the degree of $x$. 
Hence $F^{(m)}=M$ for $m$ sufficiently small, as $F^{\bullet}$ is exhaustive.
Renumbering the indices if necessary, we may assume $F^{(0)}=M$.
Since $I$ is finitely generated, $\bigoplus_{0\leq m\leq k}I^m/I^{m+1}$ is a finitely generated $A$-module, and hence $\bigoplus_{0\leq m\leq k}F^{(m)}/F^{(m+1)}$ is also finitely generated over $A$ for any $k$.

Hence, to prove the assertion, it is enough to show that there exists an integer $N\geq 0$ such that $F^{(n+1)}$ for $n\geq N$ is finitely generated $A$-module and that $IF^{(n)}=F^{(n+1)}$.

We take $N$ such that $\gr^{\bullet}_F(M)$ is generated by homogeneous elements of degree less than $N$.
If $n\geq N$, then $\gr^n_F(M)=F^{(n)}/F^{(n+1)}$ is the image of $I^{n-N}F^{(N)}$.
Consider $F^{(n+1)}$ with the filtration $F'$ induced from $F^{\bullet}$.
Then $\gr^{\bullet}_{F'}(F^{(n+1)})=\bigoplus_{m\geq n+1}F^{(m)}/F^{(m+1)}$ is generated by $IF^{(n)}/F^{(n+2)}$ as a $\gr^{\bullet}_I(A)$-module.
Since the $A$-modules $F^{(n)}/F^{(n+1)}$ and $I$ are finitely generated, we can find a finite set $\{u_j\}_{j\in J}$ of elements of $IF^{(n)}$ that generates $IF^{(n)}/F^{(n+2)}$.
Then by \cite[Chap.\ III, \S2.9, Prop.\ 12]{Bourb1} $\{u_j\}_{j\in J}$ generates $F^{(n+1)}$.
Hence $F^{(n+1)}$ is finitely generated and, moreover, we deduce that $F^{(n+1)}\subseteq IF^{(n)}$, as desired.
\end{proof}

\addcontentsline{toc}{subsection}{Exercises}
\subsection*{Exercises}
\begin{exer}\label{exer-filtrationtopclosure}{\rm 
Let $M$ be a module over a ring $A$, and consider the topology on $M$ defined by a descending filtration $F^{\bullet}=\{F^{\lambda}\}_{\lambda\in\Lambda}$ of $A$-submodules of $M$.
Then show that for an $A$-submodule $N\subseteq M$ the closure $\ovl{N}$ of $N$ is given by
$$
\ovl{N}=\bigcap_{\lambda\in\Lambda}(N+F^{\lambda}).
$$}
\end{exer}

\begin{exer}\label{exer-principalclosed}{\rm 
Let $A$ be a ring, $\{I^{(\lambda)}\}_{\lambda\in\Lambda}$ a descending filtration by ideals, and $g\in A$ an element. 
Suppose that:
\begin{itemize}
\item $A$ is complete with respect to the topology defined by $\{I^{(\lambda)}\}_{\lambda\in\Lambda}$;
\item $(g$ mod $I^{(\lambda)})$ is a non-zero-divisor in $A/I^{(\lambda)}$ for any $\lambda\in\Lambda$.
\end{itemize}
Then show that the principal ideal $(g)\subseteq A$ is closed.}
\end{exer}

\begin{exer}\label{exer-lemcriterionadicness3}{\rm 
Let $A$ be a ring, $I\subseteq A$ an ideal, and $B$ a faithfully flat $A$-algebra.
Let $\{J^{(\lambda)}\}_{\lambda\in\Lambda}$ be a descending filtration by ideals of $A$.
The topology on $B$ defined by the filtration $\{J^{(\lambda)}B\}_{\lambda\in\Lambda}$ is $IB$-adic if and only if the topology on $A$ defined by the filtration $\{J^{(\lambda)}\}_{\lambda\in\Lambda}$ is $I$-adic.}
\end{exer}

\begin{exer}\label{exer-Fujiwaralemma2.5}{\rm 
Let $A$ be a ring, $I=(a_1,\ldots,a_r)\subseteq A$ a finitely generated ideal, and $M$ an $A$-module.
Show that, if $M$ is $a_i$-adically complete for each $i=1,\ldots,r$, then $M/\bigcap_{n\geq 0}I^nM$ is $I$-adically complete.
In particular, if $M$ is $I$-adically separated, then $M$ is $I$-adically complete if and only if it is $a_i$-adically complete for each $i=1,\ldots,r$.}
\end{exer}

\begin{exer}\label{exer-filtrationcompletion}{\rm 
Let $M$ be a module over a ring $A$, and $N\subseteq M$ an $A$-submodule.
Consider a descending filtration $F^{\bullet}=\{F^{(n)}\}_{n\geq 0}$ by $A$-submodules of $N$, and topologize $M$ and $N$ by this filtration.
Let $N^{\wedge}_{F^{\bullet}}$ and $M^{\wedge}_{F^{\bullet}}$ be the Hausdorff completions of $N$ and $M$, respectively, with respect to this topology.

(1) The $A$-module $N^{\wedge}_{F^{\bullet}}$ is canonically an $A$-subgroup of $M^{\wedge}_{F^{\bullet}}$, and the canonical map $M/N\rightarrow M^{\wedge}_{F^{\bullet}}/N^{\wedge}_{F^{\bullet}}$ is an isomorphism.

(2) If the filtration $F^{\bullet}$ is separated, then $N^{\wedge}_{F^{\bullet}}\cap M=N$.}
\end{exer}

\begin{exer}\label{exer-topologicaltensorproducts}
Let $A$ be a ring endowed with a descending filtration $\{I^{(\lambda)}\}_{\lambda\in\Lambda}$ by ideals.
Let $M$ (resp.\ $N$) be an $A$-module endowed with a descending filtration $F^{\bullet}=\{F^{\alpha}\}_{\alpha\in\Sigma}$ (resp.\ $G^{\bullet}=\{G^{\beta}\}_{\beta\in T}$) by $A$-submodules.
We suppose that for any $\alpha\in\Sigma$ (resp.\ $\beta\in T$) there exists $\lambda\in\Lambda$ such that $I^{(\lambda)}N\subseteq F^{\alpha}$ (resp.\ $I^{(\lambda)}M\subseteq G^{\beta}$).
Consider the tensor product $M\otimes_AN$ together with the descending filtration $\{H^{\alpha,\beta}\}_{(\alpha,\beta)\in\Sigma\times T}$ by $A$-submodules given by
$$
H^{\alpha,\beta}=\image(F^{\alpha}\otimes_AN\rightarrow M\otimes_AN)+\image(M\otimes_AG^{\beta}\rightarrow M\otimes_AN)
$$
for any $(\alpha,\beta)\in\Sigma\times T$; here we regard $\Sigma\times T$ as a directed set by the ordering 
$$
(\alpha,\beta)\leq(\alpha',\beta')\quad\Longleftrightarrow\quad\alpha\leq\alpha',\ \beta\leq\beta'.
$$
The Hausdorff completion of $M\otimes_AN$ with respect to the topology defined by this filtration is denoted by
$$
M\widehat{\otimes}_AN,
$$
and is called the {\em complete tensor product}\index{complete!complete tensor product@--- tensor product}.
Show that there exist canonical isomorphisms
$$
M\widehat{\otimes}_AN\stackrel{\sim}{\longrightarrow}M^{\wedge}_{F^{\bullet}}\widehat{\otimes}_AN^{\wedge}_{G^{\bullet}}\stackrel{\sim}{\longrightarrow}M^{\wedge}_{F^{\bullet}}\widehat{\otimes}_{A^{\wedge}_{I^{\bullet}}}N^{\wedge}_{G^{\bullet}}.
$$
\end{exer}

\begin{exer}\label{exer-counterexaflatness1}{\rm 
Let $V$ be a valuation ring of height $2$, and $a\in V$ an element such that $\bigcap_{n\geq 0}(a^n)$ is the prime ideal of height $1$ (cf.\ \ref{prop-associatedsepval}).
Then show that the $a$-adic completion $\widehat{V}$ is not flat over $V$.}
\end{exer}

\begin{exer}\label{exer-APverifyinglemma}{\rm 
Let $A$ be a ring, and $I\subseteq A$ an ideal.
Show that the condition {\bf (AP)} (resp.\ {\bf (APf)}) in \S\ref{subsub-ARIgoodnesssub} is equivalent to the following one: for any finitely generated $A$-module $M$ and any $A$-submodule (resp.\ finitely generated $A$-submodule) $N\subseteq M$ such that $I^nN=0$ for some $n\geq 0$, there exists $m\geq 0$ such that $N\cap I^mM=0$.}
\end{exer}

\begin{exer}\label{exer-auxiliaryexeroncompletion}{\rm 
Let $A$ be a ring, and $I\subseteq A$ a finitely generated ideal.
Suppose that the $I$-adic completion $\widehat{A}$ with the $I\widehat{A}$-adic topology satisfies the condition {\bf (AP)} (resp.\ {\bf (APf)}).
Then show that for any finitely generated (resp.\ finitely presented) $A$-module $M$ the canonical map $M\otimes_A\widehat{A}\rightarrow\widehat{M}$ is an isomorphism.}
\end{exer}

\begin{exer}\label{exer-Fujiwaralemma2.6}{\rm 
Let $T$ be an abelian group topologized by a descending filtration $\{F^{(n)}\}_{n\in\Z}$ by subgroups. 

(1) Suppose that the filtration $F^{\bullet}$ is separated. 
For a fixed real number $0<c<1$, define a function $d\colon T\times T\rightarrow\R$ by
$$
d(x,y)=\begin{cases}0&(x=y),\\ c^a&(x\neq y,\ a=\min\{n\,|\,x-y\in F^{(n)}\}).\end{cases}
$$
Show that the function $d$ gives a metric on $T$.

(2) Show that $T$ is complete with respect to the filtration $F^{\bullet}$ if and only if the metric space $(T,d)$ is complete. 

(3) (Baire's category theorem\index{Baires category theorem@Baire's category theorem}) Suppose that $T$ is complete. Show that, if $U_1,U_2,\ldots$ is a sequence of dense open subsets of $T$, the intersection $\bigcap_{i=1}^{\infty}U_i$ is non-empty and dense in $T$ (cf.\ \cite{Bourb5}, IX, \S5.3, Theorem 1).

(4) Suppose $T$ is complete. Show that, if $T_1,T_2,\ldots$ is a sequence of closed subgroups such that $T=\bigcup^{\infty}_{i=1}T_i$, then $T_i$ for some $i$ is an open subgroup.}
\end{exer}


\section{Pairs}\label{sec-pairs}
Simply by a {\em pair}, we mean throughout this book a couple $(A,I)$ consisting of a commutative ring $A$ and an ideal $I\subseteq A$.
Given a pair $(A,I)$, one has the $I$-adic topology\index{topology!adic topology@adic ---}\index{adic!adic topology@--- topology} (\S\ref{subsub-adicfiltrationtopology}) on $A$-modules.
In fact, according to our definition of morphisms between pairs, the notion of pairs is equivalent to the notion of adically topologized rings, and hence, the theory of pairs, which we develop in this section, can be seen as a continuation of what we have done in the previous section, being more focused on adic topology, especially in the context of homological algebra.

After briefly discussing generalities on pairs in \S\ref{sub-pairspairs}, we discuss in \S\ref{sub-btcondition} the so-called {\em bounded torsion condition} {\bf (BT)} on pairs, which will turn out to be closely related with preservation of adicness (e.g.\ the condition {\bf (AP)}) discussed in \S\ref{subsub-ARIgoodnesssub}.
This subsection contains (without proof) a significant result by Gabber\index{Gabber, O.} (\ref{thm-gabberIHES2008}), which says that any $I$-adically complete Noetherian-outside-$I$ ring satisfies both {\bf (BT)} and {\bf (AP)}.
In \S\ref{sub-powerseries}, we recall the so-called {\em restricted formal power series rings}, an important object in the theory of formal schemes.

In \S\ref{sub-adhesive}, we will introduce new notions, {\em $($pseudo-$)$ adhesive} and {\em universally $($pseudo-$)$ adhesive} pairs.
It will turn out that adhesive pairs satisfy the condition {\bf (AP)}, hence having many useful properties as topological rings.
Moreover, interestingly enough, universally adhesivness guarantees some ring-theoretic properties of the underlying ring.
For example, if $(A,I)$ is universally adhesive and $A$ is $I$-torsion free, then the ring $A$ is universally coherent (\ref{dfn-universallycoherent}), as we will see in \S\ref{subsub-coherencyadhesive}.
This fact might imply that the notion of universally adhesive pairs provides a good generalization of Noetherian rings.

The subsection \S\ref{sub-valuative} discusses a different matter, the so-called {\em $I$-valuative rings}, which give, as we shall see later, the prototype for the local rings attached to the visualization of rigid spaces.

In the final subsection \S\ref{sub-paircomplex}, we collect several useful results on homological algebra in interplay with $I$-adic topology and related filtrations, which provide technical bases for our later calculations of cohomologies on formal schemes.

\subsection{Pairs}\label{sub-pairspairs}
\subsubsection{Generalities}\label{subsub-pairspairsgen}
By a {\em pair}\index{pair} we mean a datum $(A,I)$ consisting of a commutative ring $A$ and an ideal $I$ of $A$.
When the ideal $I$ is principal, say $I=(a)$, then we often write $(A,a)$ in place of $(A,I)$.
A {\em morphism of pairs}\index{morphism of pairs@morphism (of pairs)} $f\colon(A,I)\rightarrow(B,J)$ is a ring homomorphism $f\colon A\rightarrow B$ satisfying the following condition: there exists an integer $n\geq 1$ such that $I^n\subseteq f^{-1}(J)$ or, equivalently, $I^nB\subseteq J$ holds.

Given a pair $(A,I)$, one can consider the $I$-adic topology (\S\ref{subsub-adicfiltrationtopology}) on $A$.
If $(B,J)$ is another pair, a ring homomorphism $f\colon A\rightarrow B$ gives rise to a morphism of pairs $(A,I)\rightarrow(B,J)$ if and only if the map $f$ is continuous with respect to the $I$-adic topology on $A$ and the $J$-adic topology on $B$ (cf.\ \ref{prop-topologyfromfiltrationproperty2} (1)).
For example, the pairs $(A,I)$ and $(A,I^n)$ for $n\geq 1$ are isomorphic to each other.

Let $(A,I)$ be a pair, and $J\subseteq A$ an ideal. Then $J$ is said to be an {\em ideal of definition}\index{pair!ideal of definition of a pair@ideal of definition a ---}\index{ideal of definition}, if the identity map $\id_A$ gives an isomorphism of pairs between $(A,I)$ and $(A,J)$ or, equivalently, there exist positive integers $m,n$ such that $I^m\subseteq J^n\subseteq I$.
A morphism $(A,I)\rightarrow (B,J)$ of pairs is said to be {\em adic}\index{morphism of pairs@morphism (of pairs)!adic morphism of pairs@adic ---} if $IB$ is an ideal of definition of $(B,J)$; this definition is consistent with the definition of adic maps already given in \S\ref{subsub-adicfiltrationtopology}.
Notice that any isomorphism of pairs is adic and that the composition of two adic maps is adic.

A pair $(A,I)$ is said to be {\em complete}\index{pair!complete pair@complete ---} (resp.\ {\em henselian}\index{pair!henselian pair@henselian ---}\index{henselian!henselian pair@--- pair}, resp.\ {\em Zariskian}\index{pair!Zariskian pair@Zariskian ---}\index{Zariskian!Zariskian pair@--- pair}) if the ring $A$ is $I$-adically complete (resp.\ henselian, resp.\ Zariskian).
Given an arbitrary pair $(A,I)$, one can construct the henselian pair $(\het{A},I\het{A})$ (cf.\ \S\ref{subsub-henselianpairs}), called the {\em henselization}\index{henselization!henselization of a pair@--- (of a pair)} of $(A,I)$, and the Zariskian pair $(\zat{A},I\zat{A})$ (cf.\ \S\ref{subsub-zariskianpairs}), called the {\em associated Zariskian pair}\index{Zariskian!associated Zariskian@associated ---!associated Zariskian of a pair@--- --- (of a pair)} of $(A,I)$.

\subsubsection{Pairs of finite ideal type}\label{subsub-pairsoffiniteidealtype}
A pair $(A,I)$ is said to be {\em of finite ideal type}\index{pair!pair of finite ideal type@--- of finite ideal type} if there exists a finitely generated ideal of definition.
If $(A,I)$ is a pair of finite ideal type, then replacing it by an isomorphic one, we may assume that $I$ is finitely generated.
Notice that, if a pair $(A,I)$ is of finite ideal type, then one can construct the complete pair $(\widehat{A},I\widehat{A})$, the so-called {\em completion}\index{completion!completion of a pair@--- (of a pair)} (\ref{prop-Iadiccompletioncomplete1}), which is again of finite ideal type.

The following proposition says that the property `of finite ideal type' is local with respect to the flat topology:
\begin{prop}\label{prop-finiteidealtype2}
Let $(A,I)$ be a pair, and $B$ a faithfully flat $A$-algebra.
If $(B,IB)$ is of finite ideal type, then so is $(A,I)$.
\end{prop}

\begin{proof}
Since $\Spec B\setminus V(IB)$ is quasi-compact, $\Spec A\setminus V(I)$ is quasi-compact.
Hence there exists a finitely generated ideal $J\subseteq A$ such that $V(J)=V(I)$.
Since $\sqrt{J}=\sqrt{I}$, and since $J$ is finitely generated, there exists $n>0$ such that $J^n\subseteq I$.
On the other hand, since $(B,IB)$ is of finite ideal type, there exists a finitely generated subideal $K\subseteq IB$ such that $I^mB\subseteq K$ for some $m>0$.
Since $V(JB)=V(K)$, there exists $l>0$ such that $K^l\subseteq JB$.
Hence we have $I^{lm}B\subseteq JB$ and thus $I^{lm}\subseteq J$.
\end{proof}

\begin{dfn}\label{dfn-adm}
{\rm Let $(A,I)$ be a pair of finite ideal type. An ideal $J\subseteq A$ is said to be {\em $I$-admissible}\index{admissible!I-admissible@$I$-{---}} or, more simply, {\em admissible} if $J$ is finitely generated and there exists an integer $n\geq 1$ such that $I^n\subseteq J$.}
\end{dfn}

In other words, $I$-admissible ideals are precisely the finitely generated open ideals with respect to the $I$-adic topology.

\subsubsection{Torsions and saturation}\label{subsub-pairstorsionsandsaturation}
Let $A$ be a ring, $a\in A$ an element, and $M$ an $A$-module.
An element $x\in M$ is said to be an {\em $a$-torsion} if there exists an integer $n>0$ such that $a^nx=0$.
The subset of all $a$-torsion elements in $M$, denoted by $M_{\ator}$, is an $A$-submodule, called the {\em $a$-torsion part} of $M$.

Let $I\subseteq A$ be an ideal.
An element $x\in M$ is called an {\em $I$-torsion}\index{torsion!Itorsion@$I$-{---}} if it is $a$-torsion for all $a\in I$.
The $A$-submodule of all $I$-torsion elements, the so-called {\em $I$-torsion part}, is denoted by $M_{\Itor}$\index{torsion!Itorsion@$I$-{---}!Itorsion part@--- part}, which is the the intersection of all $M_{\ator}$ for $a\in I$.
It is easy to see that, if $I$ is generated by $a_1,\ldots,a_r\in A$, $M_{\Itor}$ coincides with the intersection of all $a_i$-torsion parts for $i=1,\ldots,r$.
We say that $M$ is {\em $I$-torsion free}\index{torsion!Itorsion@$I$-{---}!Itorsion free@--- free} (resp.\ {\em $I$-torsion}\index{torsion!Itorsion@$I$-{---}}) if $M_{\Itor}=0$ (resp.\ $M_{\Itor}=M$).

Let $N\subseteq M$ be an $A$-submodule.
The so-called {\em $I$-saturation of $N$ in $M$}\index{saturation!Isaturation@$I$-{---}} is the $A$-submodule 
$$
\til{N}=\{x\in M\,|\,\textrm{for any $a\in I$ there exists $n\geq 0$ such that $a^nx\in N$}\},
$$
or, equivalently, the inverse image of $(M/N)_{\Itor}$ by the canonical map $M\rightarrow M/N$.
We say that $N$ is {\em $I$-saturated in $M$}\index{saturated!Isaturated@$I$-{---}} if $N=\til{N}$.
It is clear that the $I$-saturation $\til{N}$ is the smallest $I$-saturated $A$-submodule containing $N$.
Note that these notions are topological, that is, stable by change of ideals of definition.
Note also that the $I$-saturation of $\{0\}$ in $M$ is nothing but $M_{\Itor}$ and hence that the $I$-torsion part $M_{\Itor}$ is always $I$-saturated.

If $(A,I)$ is of finite ideal type, then these notions are described in terms of schemes as follows.
Set $X=\Spec A$ and $U=X\setminus V(I)$, and let $j\colon U\hookrightarrow X$ be the open immersion. 
Let $\mathscr{M}$ (resp.\ $\mathscr{N}$) be the quasi-coherent sheaf on $X$ associated to $M$ (resp.\ $N$). 
Then $M_{\Itor}$ corresponds to the kernel of the canonical morphism $\mathscr{M}\rightarrow j_{\ast}j^{\ast}\mathscr{M}$, and $\til{N}$ to the pull-back of $j_{\ast}j^{\ast}\mathscr{N}$ by the same morphism.
Note that, since $(A,I)$ is of finite ideal type, the morphism $j$ is quasi-compact and hence that these sheaves are quasi-coherent (cf.\ \cite[$\mathbf{I}$, (9.2.2)]{EGA}).

\begin{dfn}\label{dfn-bt1}{\rm
Let $(A,I)$ be a pair, and $M$ an $I$-torsion $A$-module.
We say that $M$ is {\em bounded $I$-torsion}\index{torsion!bounded torsion@bounded ---} if there exists an integer $n\geq 1$ such that $I^nM=0$.}
\end{dfn}

Clearly, if $I$ is a finitely generated ideal, then any finitely generated $I$-torsion module is bounded.
We have already stated in \ref{lem-bt1vis} that any bounded $I$-torsion $A$-module is $I$-adically complete; the following is a slightly enhanced version:
\begin{lem}\label{lem-bt1}
Let $(A,I)$ be a pair, and $M$ a bounded $I$-torsion $A$-module.
Then $M$ is $I$-adically complete.
Moreover, the canonical map $M\rightarrow M\otimes_AA^{\wedge}_{I^{\bullet}}$ is an isomorphism, where $A^{\wedge}_{I^{\bullet}}$ is the Hausdorff completion of $A$ with respect to the $I$-adic topology.
\end{lem}

\begin{proof}
The first part of the lemma is clear, since $I^nM=0$ for $n\gg 1$.
Since $M$ is automatically an $A^{\wedge}_{I^{\bullet}}$-module, and any $A$-module homomorphism from $M$ to an $A^{\wedge}_{I^{\bullet}}$-module is automatically an $A^{\wedge}_{I^{\bullet}}$-module homomorphism, we have $M\cong M\otimes_AA^{\wedge}_{I^{\bullet}}$ by the universality of tensor products.
\end{proof}

\begin{dfn}\label{dfn-outsideI}{\rm Let $(A,I)$ be a pair.

(1) We say $A$ is {\em locally Noetherian outside $I$}\index{Noetherian!Noetherian outside I@--- outside $I$!locally Noetherian outside $I$@locally --- ---} if the scheme $\Spec A\setminus V(I)$ is locally Noetherian; if, moreover, $(A,I)$ is of finite ideal type, then $A$ is said to be {\em Noetherian outside $I$}\index{Noetherian!Noetherian outside I@--- outside $I$}.

(2) An $A$-module $M$ is said to be {\em finitely generated}\index{finitely generated!finitely generated outside I@--- outside $I$} (resp.\ {\em finitely presented}\index{finitely presented!finitely presented outside I@--- outside $I$}) {\em outside $I$} if the quasi-coherent sheaf $\til{M}|_U$ on $U=\Spec A\setminus V(I)$ (cf.\ \S\ref{subsub-schemesbasics} for the notation) is of finite type (resp.\ of finite presentation).}
\end{dfn}

The following lemma will be useful later:
\begin{lem}\label{lem-finitelypresentedoutside}
Let $(A,I)$ be a pair of finite ideal type, and $M$ an $A$-module.

{\rm (1)} If $M$ is finitely generated outside $I$, then there exists a finitely generated submodule $N\subseteq M$ such that $M/N$ is $I$-torsion.

{\rm (2)} Suppose $M$ is finitely generated as an $A$-module.
If $M$ is finitely presented outside $I$, then there exist a finitely presented $A$-module $N$ and a surjective morphism $N\rightarrow M$ whose kernel is $I$-torsion.
\end{lem}

\begin{proof}
Set $X=\Spec A$, $U=X\setminus V(I)$, and $\mathscr{F}=\til{M}$.

(1) Since $\mathscr{F}|_U$ is a quasi-coherent sheaf of finite type, by \cite[$\mathbf{I}$, (9.4.7) \& $\mathbf{IV}$, (1.7.7)]{EGA} we have a quasi-coherent subsheaf $\mathscr{G}\subseteq\mathscr{F}$ of finite type such that $\mathscr{G}|_U=\mathscr{F}|_U$.
If $N\subseteq M$ is the finitely generated $A$-submodule such that $\mathscr{G}=\til{N}$, then $M/N$ is $I$-torsion, since $\mathscr{G}|_U=\mathscr{F}|_U$.

(2) Take a surjective morphism $A^n\rightarrow M$, and let $K$ be the kernel.
Then $K$ is finitely generated outside $I$ (cf.\ \cite[Chap.\ I, \S2.8, Lemma 9]{Bourb1}).
Hence by (1) there exists a finitely generated $A$-submodule $K_0\subseteq K$ such that $K/K_0$ is $I$-torsion.
Set $N=A^n/K_0$, which is a finitely presented $A$-module, and consider the surjective map $N\rightarrow M$.
Its kernel is $K/K_0$, which is $I$-torsion.
\end{proof}

\subsection{Bounded torsion condition and preservation of adicness}\label{sub-btcondition}
\subsubsection{Bounded torsion condition}\label{subsub-BTfirstprop}
\index{torsion!bounded torsion@bounded ---!bounded torsion condition@--- --- condition|(}\index{bounded torsion condition|(}
The following conditions for a pair $(A,I)$, the so-called {\em bounded torsion conditions}, will be important in what follows:
\begin{itemize}
\item[{\bf (BT)}] for any finitely generated $A$-module $M$, $M_{\Itor}$ is bounded $I$-torsion\index{torsion!bounded torsion@bounded ---}; 
\item[{\bf (UBT)}] for any $A$-algebra $B$ of finite type, $(B,IB)$ satisfies the condition {\bf (BT)}.
\end{itemize}

The following proposition is easy to see:
\begin{prop}\label{prop-bt31}
Let $(A,I)$ be a pair satisfying the condition {\bf (BT)} $($resp.\ {\bf (UBT)}$)$.
If $B$ is a finite $A$-algebra $($resp.\ an $A$-algebra of finite type$)$, then $(B,IB)$ satisfies {\bf (BT)} $($resp.\ {\bf (UBT)}$)$. \hfill$\square$
\end{prop}

\begin{cor}\label{cor-bt32}
For a pair $(A,I)$ to satisfy {\bf (UBT)} it is necessary and sufficient that the pairs $(A[X_1,\ldots,X_n],IA[X_1,\ldots,X_n])$ by the polynomial rings satisfy {\bf (BT)} for any $n\geq 0$. \hfill$\square$
\end{cor}

\begin{prop}\label{prop-bt33}
{\rm (1)} Let $(A,I)$ and $(B,J)$ be pairs.
Then $(A,I)$ and $(B,J)$ satisfy {\bf (BT)} $($resp.\ {\bf (UBT)}$)$ if and only if $(A\times B,I\times J)$ satisfies {\bf (BT)} $($resp.\ {\bf (UBT)}$)$.

{\rm (2)} Let $(A,I)$ be a pair, and $B$ a faithfully flat $A$-algebra.
If $(B,IB)$ satisfies {\bf (BT)} $($resp.\ {\bf (UBT)}$)$, then $(A,I)$ satisfies {\bf (BT)} $($resp.\ {\bf (UBT)}$)$.
\end{prop}

\begin{proof}
(1) The `if' part follows from \ref{prop-bt31}.
To show the converse, notice that any finitely generated $A\times B$-module $M$ is a product of finitely generated mouldes $M=M_A\oplus M_B$, and we clearly have $M_{\textrm{$I\times J$-}\mathrm{tor}}=(M_A)_{\Itor}\times(M_B)_{\Jtor}$.
Hence if $(A,I)$ and $(B,J)$ satisfy {\bf (BT)}, then {\bf (BT)} for $A\times B$ follows immediately.
A similar argument works for finitely generated algebras, and hence {\bf (UBT)} for $A\times B$ follows from that for $A$ and $B$.

(2) Clearly, it suffices to check the case of {\bf (BT)}.
Let $M$ be a finitely generated $A$-module, and consider $M\otimes_AB$.
We have $I^n(M\otimes_AB)_{\Itor}=0$ for some $n>0$.
Since $B$ is flat over $A$, we have $(M\otimes_AB)_{\Itor}\subseteq M_{\Itor}\otimes_AB$, but since the right-hand module is clearly $I$-torsion, we have $(M\otimes_AB)_{\Itor}=M_{\Itor}\otimes_AB$.
Since $B$ is faithfully flat over $A$, we have $I^nM_{\Itor}=0$.
\end{proof}

\begin{prop}\label{prop-BTZariskian}
Let $(A,I)$ be a pair of finite ideal type.
Then $(A,I)$ satisfies {\bf (BT)} if and only if so does the associated Zariskian pair $(\zat{A},I\zat{A})$.
\end{prop}

\begin{proof}
We may assume that $I$ is finitely generated.
Suppose $(A,I)$ satisfies {\bf (BT)}.
For a finitely generated $\zat{A}$-module $N$ one can find a finitely generated $A$-module $M$ such that $M\otimes_A\zat{A}\cong N$.
For any element $x/(1+a)\in N_{\Itor}$ ($x\in M$, $a\in I$), since $I$ is finitely generated, one can find $b\in I$ such that $(1+b)x\in M_{\Itor}$.
Then one can find $n\geq 0$, independently of the element $(1+b)x$, such that $I^n\cdot(1+b)x=0$, thereby $I^n\cdot x/(1+a)=0$.

Conversely, suppose $(\zat{A},I\zat{A})$ satisfies {\bf (BT)}, and for a finitely generated $A$-module $M$, consider $M\otimes_A\zat{A}$.
Take any $x\in M_{\Itor}$.
Since $I$ is finitely generated, there exists $a\in I$ such that $(1-a)x$ is annihilated by $I^n$ (where $n$ is independent of $x$).
For any $b\in I$ we have $(1-a)b^nx=0$ and thus $b^nx=ab^nx=a^2b^nx=\cdots=0$, which shows that $I^n\cdot x=0$.
\end{proof}

\begin{prop}\label{prop-btred}
The following conditions for a pair $(A,I)$ are equivalent$:$
\begin{itemize}
\item[{\rm (a)}] $(A,I)$ satisfies {\bf (BT)} $($resp.\ {\bf (UBT)}$);$
\item[{\rm (b)}] $A_{\Itor}$ is bounded $I$-torsion, and $(A/A_{\Itor},I(A/A_{\Itor}))$ satisfies {\bf (BT)} $($resp.\ {\bf (UBT)}$)$. 
\end{itemize}
\end{prop}

\begin{proof}
Since $(A/A_{\Itor})[X]\cong A[X]/A[X]_{\Itor}$, it suffices to show the assertion for {\bf (BT)}.
(a) $\Rightarrow$ (b) is obvious by \ref{prop-bt31}.
To show the converse, let $M$ be a finitely generated $A$-module, and set $M'=M\otimes_A(A/A_{\Itor})=M/A_{\Itor}M$.
Then it is easy to verify that we have the exact sequence
$$
0\longrightarrow A_{\Itor}M\longrightarrow M_{\Itor}\longrightarrow M'_{\Itor}\longrightarrow 0.
$$
Since $A_{\Itor}M$ and $M'_{\Itor}$ are bounded $I$-torsion, so is $M_{\Itor}$.
\end{proof}

\begin{prop}\label{prop-btarfflat2}
Let $(A,I)$ be a pair of finite ideal type, and suppose that the following conditions are satisfied$:$
\begin{itemize}
\item[{\rm (a)}] $(\widehat{A},I\widehat{A})$ satisfies {\bf (BT)}$;$ 
\item[{\rm (b)}] $A\rightarrow\widehat{A}$ is flat. 
\end{itemize}
Then $(A,I)$ satisfies {\bf (BT)}.
\end{prop}

\begin{proof}
We may assume that $I$ finitely generated.
Since $A\rightarrow\widehat{A}$ is flat, $\zat{A}\rightarrow\widehat{A}$ is faithfully flat (\ref{prop-relpair21} (2)). 
By \ref{prop-bt33} (2) the pair $(\zat{A},I\zat{A})$ satisfies {\bf (BT)}, and hence the assertion follows from \ref{prop-BTZariskian}.
\end{proof}

\subsubsection{Preservation of adicness}\label{subsub-APfirstprop}
We say that a pair $(A,I)$ satisfies the condition {\bf (AP)} if $A$ endowed with the $I$-adic topology satisfies the condition {\bf (AP)} in \S\ref{subsub-ARIgoodnesssub}.
We will also consider the following condition for pairs $(A,I)$: 
\begin{itemize}
\item[{\bf (UAP)}] for any $A$-algebra $B$ of finite type, the induced pair $(B,IB)$ satisfies the condition {\bf (AP)}.
\end{itemize}

\begin{prop}\label{prop-bt31AP}
Let $(A,I)$ be a pair satisfying the condition {\bf (AP)} $($resp.\ {\bf (UAP)}$)$.
If $B$ is a finite $A$-algebra $($resp.\ an $A$-algebra of finite type$)$, then $(B,IB)$ satisfies {\bf (AP)} $($resp.\ {\bf (UAP)}$)$. \hfill$\square$
\end{prop}

\begin{cor}\label{cor-bt32AP}
For a pair $(A,I)$ to satisfy {\bf (UAP)} it is necessary and sufficient that the pairs $(A[X_1,\ldots,X_n],IA[X_1,\ldots,X_n])$ by the polynomial rings for $n\geq 0$ satisfy {\bf (AP)}. \hfill$\square$
\end{cor}

One can show the following proposition by an argument similar to that in \ref{prop-bt33}:
\begin{prop}\label{prop-bt33AP}
{\rm (1)} Let $(A,I)$ and $(B,J)$ be pairs.
Then $(A,I)$ and $(B,J)$ satisfy {\bf (AP)} $($resp.\ {\bf (UAP)}$)$ if and only if $(A\times B,I\times J)$ satisfies {\bf (AP)} $($resp.\ {\bf (UAP)}$)$.

{\rm (2)} Let $(A,I)$ be a pair, and $B$ a faithfully flat $A$-algebra.
If $(B,IB)$ satisfies {\bf (AP)} $($resp.\ {\bf (UAP)}$)$, then $(A,I)$ satisfies {\bf (AP)} $($resp.\ {\bf (UAP)}$)$. \hfill$\square$
\end{prop}

\begin{prop}\label{prop-btarfflat2AP}
Let $(A,I)$ be a pair of finite ideal type, and suppose that the following conditions are satisfied$:$
\begin{itemize}
\item[{\rm (a)}] $(\widehat{A},I\widehat{A})$ satisfies {\bf (AP)}$;$ 
\item[{\rm (b)}] $A\rightarrow\widehat{A}$ is flat$;$ 
\item[{\rm (c)}] $A$ is Noetherian outside $I$\index{Noetherian!Noetherian outside I@--- outside $I$} {\rm (\ref{dfn-outsideI} (1))}.
\end{itemize}
Then $(A,I)$ satisfies {\bf (AP)}.
\end{prop}

\begin{proof}
We may assume that $I$ is finitely generated.
Since $A\rightarrow\widehat{A}$ is flat, $\zat{A}\rightarrow\widehat{A}$ is faithfully flat (\ref{prop-relpair21} (2)). 
By \ref{prop-bt33AP} (2) the pair $(\zat{A},I\zat{A})$ satisfies {\bf (AP)}.
On the other hand, since $I$ is finitely generated, the open set $\Spec A\setminus V(I)$ is the union of finitely many affine open subsets of the form $D(f_i)$ ($i=1,\ldots,r$), where $f_1,\ldots,f_r$ generate $I$.
Since each $A_{f_i}$ is Noetherian, the pair $(A_{f_i},IA_{f_i})$ clearly satisfies {\bf (AP)} (cf.\ \ref{prop-ARclassic}).
Let $B=\zat{A}\times A_{f_1}\times\cdots\times A_{f_r}$. 
Then the pair $(B,IB)$ satisfies {\bf (AP)} (\ref{prop-bt33AP} (1)), and $A\rightarrow B$ is faithfully flat.
Then again by \ref{prop-bt33AP} (2) the pair $(A,I)$ satisfies {\bf (AP)}.
\end{proof}

For a pair $(A,I)$ the condition {\bf (AP)} can be refined as follows.
In classical commutative algebra the condition $(\ast)$ in {\rm \ref{lem-ARIgoodnesssubfiltration1}} is usually verified by means of the $I$-goodness (cf.\ \ref{dfn-Igood}) of the induced filtration $G^{\bullet}=\{N\cap I^nM\}_{n\geq 0}$.
For an $A$-module $M$ and an $A$-submodule $N\subseteq M$, the $I$-goodness in question is equivalently paraphrased into the following form:
\begin{itemize}
\item[$(\ast\ast)$] {\it there exists a non-negative integer $c$ such that for any $n>c$ we have}
$$
{\textstyle N\cap I^nM=I^{n-c}(N\cap I^cM).}
$$
\end{itemize}
Similarly to the conditions {\bf (AP)} and {\bf (APf)}, accordingly, we can consider the following conditions for a pair $(A,I)$:
\begin{itemize}
\item[{\bf (AR)}] the condition $(\ast\ast)$ holds for any finitely generated $A$-module $M$ and any $A$-submodule $N\subseteq M$; 
\item[{\bf (ARf)}] the condition $(\ast\ast)$ holds for any finitely generated $A$-module $M$ and any {\em finitely generated} $A$-submodule $N\subseteq M$.
\end{itemize}
Obviously, these conditions, together with {\bf (AP)} and {\bf (APf)}, sit in the following diagram of implications:
$$
\xymatrix@R-3ex@C-4ex{{\bf (AR)}\ar@{=>}[rr]\ar@{=>}[dr]&&{\bf (AP)}\ar@{=>}[dr]\\ &{\bf (ARf)}\ar@{=>}[rr]&&{\bf (APf)}}
$$

It is well-known that the condition {\bf (AR)} is always satisfied if $A$ is Noetherian.
This is exactly what classical Artin-Rees lemma\index{Artin-Rees lemma} asserts:
\begin{prop}[Classical Artin-Rees lemma; cf.\ e.g.\ {\cite[Theorem 8.5]{Matsu}}]\label{prop-ARclassic}
Any pair $(A,I)$ with $A$ Noetherian satisfies the condition {\bf (AR)}. \hfill$\square$
\end{prop}

\danger{The conditions {\bf (AR)} and {\bf (ARf)} depend on the choice of the ideal of definition $I\subseteq A$; hence, when considering these conditions for adically topologized ring, one has to specify an ideal of definition.}

\subsubsection{The properties {\bf (BT)} and {\bf (AP)}}\label{subsub-BTAR}
\begin{prop}\label{prop-APimpliesBT}
The condition {\bf (AP)} implies {\bf (BT)}.
\end{prop}

\begin{proof}
Suppose a pair $(A,I)$ satisfies {\bf (AP)}, and let $M$ be a finitely generated $A$-module.
Define $F_k$ for any $k\geq 1$ to be the $A$-submodule of $M$ consisting of elements annihilated by $I^k$.
Clearly, we have $M_{\Itor}=\bigcup_{k\geq 1}F_k$.
By Exercise \ref{exer-APverifyinglemma} there exists $m>0$ such that $F_1\cap I^mM=0$.
Then for $k\geq m+1$ we have $I^{k-1}F_k\subseteq F_1\cap I^mM=0$ and hence $F_m=F_{m+1}=F_{m+2}=\cdots$.
Thus we deduce $M_{\Itor}=F_m$, which is of bounded $I$-torsion.
\end{proof}

\begin{prop}\label{prop-ARbt1}
Let $(A,I)$ be a pair with $I=(a)$ principal.
Suppose $(A,I)$ satisfies the condition {\bf (BT)}.
Then it satisfies {\bf (AR)}.
\end{prop}

The main part of the proof relies on the following lemma:
\begin{lem}\label{lem-ARprincipalgeneralized}
Let $(A,I)$ be a pair with $I=(a)$ principal, $M$ an $A$-module, and $N\subseteq M$ an $A$-submodule. 
Suppose $a^n(M/N)_{\ator}=0$ for some $n\geq 0$.
Then for any $m\geq 0$ we have
$$
N\cap a^{n+m}M=a^m(N\cap a^nM),
$$
that is, the condition $(\ast\ast)$ in {\rm \S\ref{subsub-ARIgoodnesssub}} is satisfied.
\end{lem}

\begin{proof}
Set $L=M/N$, and denote the canonical projection $M\rightarrow L$ by $\pi$. 
The inclusion $a^m(N\cap a^nM)\subseteq N\cap a^{n+m}M$ is trivial. 
To show the converse inclusion, take $x=a^my\in N\cap a^{n+m}M$ ($y\in a^nM$).
Since $a^m\pi(y)=\pi(x)=0$ in $L$, we have $\pi(y)\in a^nL\cap L_{\ator}$.
But since $a^nL\cap L_{\ator}=a^nL_{\ator}=0$, we have $y\in N$, that is, $x=a^my\in a^m(N\cap a^nM)$, as desired.
\end{proof}

\begin{proof}[Proof of Proposition {\rm \ref{prop-ARbt1}}]
Let $M$ be a finitely generated $A$-module, and $N\subseteq M$ an $A$-submodule.
The condition {\bf (BT)} implies $a^n(M/N)_{\ator}=0$ for some $n\geq 0$.
Hence by \ref{lem-ARprincipalgeneralized} the condition $(\ast\ast)$ in {\rm \S\ref{subsub-ARIgoodnesssub}} is satisfied.
\end{proof}

\begin{prop}\label{prop-APandBTnonprincipal}
Let $A$ be a ring, and $I=(a_1,\ldots,a_r)\subseteq A$ a finitely generated ideal.
Suppose that for each $i=1,\ldots,r$ the pair $(A,a_i)$ satisfies {\bf (BT)}.
Then the pair $(A,I)$ satisfies {\bf (BT)} and {\bf (AP)}.
\end{prop}

\begin{proof}
For any finitely generated $A$-module $M$ the $a_i$-torsion part of $M$ is bounded for each $i=1,\ldots,r$.
As the $I$-torsion part of $M$ is the intersection of the $a_i$-torsion parts for $1\leq i\leq r$, it is bounded $I$-torsion.
Thus $(A,I)$ satisfies {\bf (BT)}.

To show that $(A,I)$ satisfies {\bf (AP)}, we apply induction with respect to $r$.
If $r=1$, the assertion follows from \ref{prop-ARbt1}.
Set $a=a_1$ and $J=(a_2,\ldots,a_r)$, and let $N\subseteq M$ be an $A$-submodule such that $I^nN=0$.
We want to show that $N\cap I^mM=0$ for some $m\geq 0$ (cf.\ Exercise \ref{exer-APverifyinglemma}).
Take $s\geq 0$ so that $a^s$ annihilates the $a$-torsion part of $M$.
Since $N$ is $a$-torsion, we have $N\cap a^sM=0$.
Hence one can regard $N$ as a submodule of $\ovl{M}=M/a^sM$.
Note that the $I$-adic topology on $\ovl{M}$ coincides with the $J$-adic topology.
By induction we know that $N\cap I^t\ovl{M}=0$ for some $t\geq 0$ (cf.\ Exercise \ref{exer-APverifyinglemma}).
Then for $m\geq\max\{s,t\}$ one has $N\cap I^mM=0$, as desired.
\end{proof}

\begin{prop}\label{prop-APandBToutsidenoetherian}
Let $(A,I)$ be a pair of finite ideal type, and suppose that $A$ is Noetherian outside $I$\index{Noetherian!Noetherian outside I@--- outside $I$} {\rm (\ref{dfn-outsideI} (1))}.
If $(A,I)$ satisfies {\bf (BT)}, then for any finitely generated subideal $J\subseteq I$, $(A,J)$ satisfies {\bf (BT)} and {\bf (AP)}. 
\end{prop}

\begin{proof}
In view of \ref{prop-APandBTnonprincipal} we only have to show that for any {\em principal} subideal $J=(a)\subseteq I$, $(A,a)$ satisfies {\bf (BT)}.
Let $M$ be a finitely generated $A$-module.
Since $\Spec A\setminus V(I)$ is Noetherian, there exists a finitely generated submodule $N$ of $M_{\ator}$ such that $M_{\ator}/N$ is supported on $V(I)$ (\ref{lem-finitelypresentedoutside} (1)).
Now the $a$-torsion part of $M/N$ is $M_{\ator}/N$, which is nothing but the $I$-torsion part of $M/N$.
Hence $M_{\ator}/N$ is bounded $a$-torsion.
Since $N$ is finitely generated and hence is bounded $a$-torsion, we deduce that $M_{\ator}$ is bounded $a$-torsion.
\end{proof}

\begin{cor}\label{cor-APandBToutsidenoetherian}
Let $(A,I)$ be a pair of finite ideal type.
If $A$ is Noetherian outside $I$\index{Noetherian!Noetherian outside I@--- outside $I$} and $(A,I)$ satisfies {\bf (BT)}, then it also satisfies {\bf (AP)}. \hfill$\square$
\end{cor}

\begin{prop}\label{prop-btarf1}
Let $(A,I)$ be a pair of finite ideal type satisfying {\bf (BT)}, and suppose $A$ is Noetherian outside $I$\index{Noetherian!Noetherian outside I@--- outside $I$}.

{\rm (1)} For any finitely generated $A$-module $M$ the canonical map $M\otimes_A\widehat{A}\rightarrow\widehat{M}$ is an isomorphism.

{\rm (2)} The canonical map $A\rightarrow\widehat{A}$ is flat.
\end{prop}

\begin{proof}
(1) By \ref{lem-finitelypresentedoutside} (2) there exists a surjective morphism $N\rightarrow M$ of $A$-modules such that $N$ is finitely presented and the kernel $K$ is $I$-torsion.
By {\bf (BT)} the kernel $K$ is bounded $I$-torsion.
Consider the exact sequence $0\rightarrow K\rightarrow N\rightarrow M\rightarrow 0$, which yields the following commutative diagram with exact rows:
$$
\xymatrix{0\ar[r]&\widehat{K}\ar[r]&\widehat{N}\ar[r]&\widehat{M}\ar[r]&0\\ &K\otimes_A\widehat{A}\ar[u]\ar[r]&N\otimes_A\widehat{A}\ar[u]\ar[r]&M\otimes_A\widehat{A}\ar[u]\ar[r]&0\rlap{;}}
$$
Notice that the exactness of the first row is due to \ref{prop-qconsistency1-2p-2} (here we use \ref{cor-APandBToutsidenoetherian}).
The first vertical arrow is an isomorphism due to \ref{lem-bt1}, and the second one is an isomorphism due to \ref{prop-ARconseq0}.
Hence the last vertical arrow is an isomorphism, which is nothing but we wanted to show.

(2) Let $\mathfrak{a}\subseteq A$ be a finitely generated ideal.
We want to show that the map $\mathfrak{a}\otimes_A\widehat{A}\rightarrow\widehat{A}$ is injective.
By (1) it follows that $\mathfrak{a}\otimes_A\widehat{A}=\widehat{\mathfrak{a}}$.
By \ref{prop-qconsistency1-2p-2}, on the other hand, $\widehat{\mathfrak{a}}$ is an ideal of $\widehat{A}$, as desired.
\end{proof}

\subsubsection{Bounded torsion condition for complete pairs}\label{subsub-BTARcomp}
For the proof of the following significant theorem by Gabber, we refer to \cite{FGK}:
\begin{thm}[O.\ Gabber; {\rm \cite[Theorem 5.1.2]{FGK}}]\label{thm-gabberIHES2008}\index{Gabber, O.}
Let $(A,I)$ be a complete pair\index{pair!complete pair@complete ---} of finite ideal type, and suppose that $A$ is Noetherian outside $I$\index{Noetherian!Noetherian outside I@--- outside $I$!locally Noetherian outside $I$@locally --- ---}.
Then $(A,I)$ satisfies {\bf (BT)} and {\bf (AP)}. \hfill$\square$
\end{thm}

In particular, by \ref{cor-qconsistency11} and \ref{cor-propARconseq1-2} we have:
\begin{cor}\label{cor-gabberIHES2008}
Let $(A,I)$ be a complete pair\index{pair!complete pair@complete ---} of finite ideal type, and suppose that $A$ is Noetherian outside $I$\index{Noetherian!Noetherian outside I@--- outside $I$!locally Noetherian outside $I$@locally --- ---}.
Then any finitely generated $A$-module is $I$-adically complete.
Moreover, if $M$ is a finitely generated $A$-module, then any $A$-submodule $N\subseteq M$ is closed in $M$ with respect to the $I$-adic topology. \hfill$\square$
\end{cor}
\index{bounded torsion condition|)}\index{torsion!bounded torsion@bounded ---!bounded torsion condition@--- --- condition|)}

\subsection{Pairs and flatness}\label{sub-gluingofflatness}
\subsubsection{Gluing of flatness}\label{subsub-gluingofflatness}
\index{gluing of flatness|(}
\begin{thm}[Gluing of flatness (I)]\label{thm-newversiongluingofflatness}
Let $(A,I)$ be a pair of finite ideal type, and $M$ an $A$-module.
Then $M$ is flat over $A$ if and only if the following conditions hold$:$
\begin{itemize}
\item[{\rm (a)}] $\Tor^A_q(M,N)=0$ for any $q\geq 1$ and any $A$-module $N$ supported in $V(I)\subseteq\Spec A;$
\item[{\rm (b)}] $\til{M}$ is flat over $\Spec A\setminus V(I)$.
\end{itemize}
\end{thm}

\begin{proof}
We may assume that the ideal $I\subseteq A$ is finitely generated; let $a_1,\ldots,a_n$ generate $I$ $(n\geq 0)$.
For $k=0,\ldots,n$ we set $I_k=(a_1,\ldots,a_k)$.
Consider the condition: 
\begin{itemize}
\item[$(\ast)_k$] $\Tor^A_q(M,N)=0$ for any $q\geq 1$ and any $A$-module $N$ supported in $V(I_k)\subseteq\Spec A$.
\end{itemize}
The condition $(\ast)_n$ is nothing but the assumption (a), and what to show is $(\ast)_0$ (cf.\ \cite[Chap.\ I, \S4, Prop.\ 1]{Bourb1}).
It suffices therefore to show $(\ast)_k$ by descending induction with respect to $k$.

Suppose $(\ast)_{k+1}$ is true, and let $N$ be an $A$-module supported in $V(I_k)$.
Set $a=a_{k+1}$, and consider the exact sequence
$$
0\longrightarrow N_{\ator}\longrightarrow N\longrightarrow N/N_{\ator}\longrightarrow 0.
$$
Since $N_{\ator}$ is supported in $V(I_{k+1})$, by induction we have $\Tor^A_q(M,N_{\ator})=0$ for $q\geq 1$.
Hence it suffices to show $\Tor^A_q(M,N/N_{\ator})=0$ for $q\geq 1$, and thus we may assume that $N$ is $a$-torsion free.
Now, consider the exact sequence
$$
0\longrightarrow N\longrightarrow N[{\textstyle \frac{1}{a}}]\longrightarrow C\longrightarrow 0.
$$
Here the $A$-module $C$ is supported in $V(I_{k+1})$.
Since $N[{\textstyle \frac{1}{a}}]$ is flat over $A[\frac{1}{a}]$ (by the assumption (b)), we have $\Tor^A_q(M,N)=0$ for $q\geq 1$, as desired.
\end{proof}

\begin{cor}\label{cor-glfl}
Let $(A,a)$ be a pair with $a\in A$ being a non-zero divisor, and $M$ an $A$-module.
Then the following conditions are equivalent$:$
\begin{itemize}
\item[{\rm (a)}] $M$ is $A$-flat$;$
\item[{\rm (b)}] $M[\frac{1}{a}]$ is $A[\frac{1}{a}]$-flat, $M/aM$ is $(A/aA)$-flat, and $M$ is $a$-torsion free.
\end{itemize}
\end{cor}

\begin{proof}
(a) $\Rightarrow$ (b) is obvious. 
To show the converse, in view of \ref{thm-newversiongluingofflatness} it suffices to show that $\Tor^A_q(M,N)$ vanishes for $q\geq 1$ and for any $A$-module $N$ supported in $V((a))$. 
Since the functor $\Tor^A_q(M,\textrm{--})$ commutes with inductive limits, we may assume that $N$ is killed by some $a^n$. 
But if so, we may further reduce the situation where $N$ is killed even by $a$ by the inductive argument using the filtration $\{a^mN\}_{m\geq 0}$.

Since $a$ is a non-zero-divisor and $M$ is $a$-torsion free, we have $\Tor^A_q(M,A/aA)=0$ for $q\geq 1$.
This implies the last isomorphy of the following:
$$
M\otimes^{\LD}_AN\cong M\otimes^{\LD}_A(A/aA)\otimes^{\LD}_{A/aA}N\cong (M/aM)\otimes^{\LD}_{A/aA}N.
$$
Now, since $M/aM$ is $(A/aA)$-flat, we deduce that $\Tor^A_q(M,N)=0$ for $q\geq 1$, as desired.
\end{proof}

\begin{prop}[Gluing of flatness (II)]\label{prop-glueingofflatnessrelative}
Let $(A,I)$ be a pair of finite ideal type, $B$ an $A$-algebra, and $M$ a $B$-module.
Suppose that the following conditions are satisfied$:$
\begin{itemize}
\item[{\rm (a)}] $B$ and $M$ are flat over $A;$
\item[{\rm (b)}] $M/IM$ is flat over $B/IB$, and $\til{M}$ is flat over $\Spec B\setminus V(IB)$.
\end{itemize}
Then $M$ is $B$-flat.
\end{prop}

\begin{proof}
In case the ideal $I$ is principal, the assertion follows immediately from \cite[5.2.1]{GabRam}.
In general, in order to apply induction with respect to the number of generators, set $I=(a_1,\ldots,a_s,b)$ and $J=(a_1,\ldots,a_s)$.
Since $\ovl{A}=A/bA$, $B/bB$, and $M/bM$ together with the ideal $\ovl{I}=J\ovl{A}$ satisfy the conditions, we deduce by induction that $M/bM$ is $B/bB$-flat.
Considering next the situation by $A$, $B$, $M$, and $(b)\subseteq A$, we conclude that $M$ is $B$-flat, as desired.
\end{proof}
\index{gluing of flatness|)}

\subsubsection{Local criterion of flatness}\label{subsub-localcriterionofflatness}
\index{local criterion of flatness|(}
\begin{prop}[Local criterion of flatness]\label{prop-flatmorformal1}
Let $(A,I)$ be a pair, and $M$ an $A$-module.
We suppose that the following conditions are satisfied$:$
\begin{itemize}
\item[{\rm (i)}] for any finitely generated ideal $\mathfrak{a}$ of $A$, the topology on $\mathfrak{a}$ induced from that of $A$ coincides with the $I$-adic topology$;$ that is, for any $n\geq 0$ there exists $k\geq 0$ such that 
$$
I^k\cap\mathfrak{a}\subseteq I^n\mathfrak{a};
$$
\item[{\rm (ii)}] $M$ is {\em idealwise separated}\index{separated!idealwise separated@idealwise ---} for $I;$ that is, for any finitely generated ideal $\mathfrak{a}$ of $A$, the $A$-module $\mathfrak{a}\otimes_AM$ is $I$-adically separated\index{separated!I-adically separated@$I$-adically ---}.
\end{itemize}
We write $A_k=A/I^{k+1}$ and $M_k=M/I^{k+1}M$ for any $k\geq 0$.
Then the following conditions are equivalent$:$
\begin{itemize}
\item[{\rm (a)}] $M$ is $A$-flat$;$
\item[{\rm (b)}] for any $A_0$-module $N$ we have $\Tor^A_1(N,M)=0;$
\item[{\rm (c)}] $M_0$ is $A_0$-flat, and we have $\Tor^A_1(A_0,M)=0;$
\item[{\rm (d)}] $M_k$ is $A_k$-flat for any $k\geq 0$.
\end{itemize}
\end{prop}

The main idea of the proof is borrowed from \cite[Chap.\ III, \S5.3]{Bourb1}. 
\begin{proof}
The implications (a) $\Rightarrow$ (b) $\Rightarrow$ (c) $\Rightarrow$ (d) are shown in \cite[Chap.\ III, \S5.2, Theorem 1]{Bourb1}, and here we omit the proofs.
We are going to show (d) $\Rightarrow$ (a) by verifying that for any finitely generated ideal $\mathfrak{a}$ of $A$, the map $\iota\colon\mathfrak{a}\otimes_AM\rightarrow M$ is injective.

Take $x\in\ker(\iota)$.
As we have $\bigcap_nI^n(\mathfrak{a}\otimes_AM)=0$ by the hypothesis (ii), it suffices to show that $x\in I^n(\mathfrak{a}\otimes_AM)$ for any $n$.
By (i) there exists $k$ such that $I^{k+1}\cap\mathfrak{a}\subseteq I^n\mathfrak{a}$.
Hence it is enough to show that $x$ belongs to the image of the canonical map $(I^{k+1}\cap\mathfrak{a})\otimes_AM\rightarrow\mathfrak{a}\otimes_AM$.
To this end, we consider the following commutative diagram with an exact row:
$$
\xymatrix{(I^{k+1}\cap\mathfrak{a})\otimes_AM\ar[r]&\mathfrak{a}\otimes_AM\ar[d]\ar[r]^(.36){(\ast)}&(\mathfrak{a}/\mathfrak{a}\cap I^{k+1})\otimes_AM\ar[d]\ar[r]&0\\ &M\ar[r]&A_k\otimes_AM\rlap{.}}
$$
Here by the assumption, the first vertical arrow maps $x$ to $0$.
Hence it is mapped to $0$ in $A_k\otimes_AM$.
On the other hand, the second vertical arrow coincides with the morphism $(\mathfrak{a}/\mathfrak{a}\cap I^{k+1})\otimes_{A_k}M_k\rightarrow M_k$, which is injective by our assumption that $M_k$ is $A_k$-flat.
Hence $x$ belongs to the kernel of the map $(\ast)$, which is nothing but the image of the first horizontal arrow in the upper row.
But this is what we wanted to prove.
\end{proof}

\begin{cor}\label{cor-flatmorformal11}
Let $(A,I)$ be a pair, and $M$ an $A$-module.
We suppose that the conditions {\rm (i)} and {\rm (ii)} in $\ref{prop-flatmorformal1}$ are satisfied and that $(A,I)$ is a Zariskian pair.
Then the following conditions are equivalent$:$
\begin{itemize}
\item[{\rm (a)}] $M$ is faithfully flat over $A;$
\item[{\rm (b)}] $M_k$ is faithfully flat over $A_k$ for any $k\geq 0$.
\end{itemize}
\end{cor}

\begin{proof}
The implication (a) $\Rightarrow$ (b) is clear.
Conversely, by \ref{prop-flatmorformal1} we know that $M$ is $A$-flat.
To show that $M$ is faithfully flat, it suffices to show that for any maximal ideal $\mathfrak{m}$ of $A$ we have $M\otimes_A(A/\mathfrak{m})\neq 0$.
Since $(A,I)$ is Zariskian, we have $I\subseteq\mathfrak{m}$.
Hence $M\otimes_A(A/\mathfrak{m})=M_0\otimes_{A_0}(A_0/\mathfrak{m}_0)$, where $\mathfrak{m}_0=\mathfrak{m}/I$.
Since $M_0$ is faithfully flat over $A_0$, the last module is non-zero, as desired.
\end{proof}

\begin{rem}\label{rem-localcriterionflatness}{\rm 
The conditions (i) and (ii) in \ref{prop-flatmorformal1} are trivially satisfied if $I$ is a nilpotent ideal.
In this case, moreover, one can show easily that the conditions (a) $\sim$ (d) are equivalent to the following one (cf.\ \cite[Expos\'e IV, Prop.\ 5.1]{SGA1}):
\begin{itemize}
\item[(e)] $M_0$ is $A_0$-flat, and the canonical surjective morphism
$$
\gr^0_I(M)\otimes_{A_0}\gr^{\bullet}_I(A)\longrightarrow\gr^{\bullet}_I(M)
$$
is an isomorphism.
\end{itemize}}
\end{rem}

In particular, we have (cf.\ \cite[Expos\'e IV, Cor.\ 5.9]{SGA1}):
\begin{prop}\label{prop-localcriterionflatuseful}
Let $(A,I)$ be a pair where $I$ is nilpotent, $B$ a flat $A$-algebra, and $M$ a $B$-module.
Then the following conditions are equivalent$:$
\begin{itemize}
\item[{\rm (a)}] $M$ is flat over $B;$
\item[{\rm (b)}] $M$ is flat over $A$, and $M_0=M/IM$ is flat over $B_0=B/IB$.
\end{itemize}
\end{prop}

\begin{proof}
The implication (a) $\Rightarrow$ (b) is trivial.
Suppose (b) holds.
In order to apply \ref{prop-flatmorformal1}, we check the condition (e) in \ref{rem-localcriterionflatness}.
For $n\geq 0$ the $n$-th graded piece of $\gr^0_{IB}(M)\otimes_{B_0}\gr^{\bullet}_{IB}(B)$ is $M\otimes_B(I^nB/I^{n+1}B)$.
Since $B$ is $A$-flat, $I^nB/I^{n+1}B=(I^n/I^{n+1})\otimes_AB$.
Since $M$ is $A$-flat, we have $M\otimes_B(I^nB/I^{n+1}B)=M\otimes_A(I^n/I^{n+1})=I^nM/I^{n+1}M$, which is the $n$-th graded piece of $\gr^{\bullet}_{IB}(M)$.
\end{proof}

Finally, we give a useful sufficient condition to verify (i) and (ii) in \ref{prop-flatmorformal1}:
\begin{prop}\label{prop-flatmorformal12cor}
Let $(A,I)$ be a pair, $B$ an $A$-algebra, and $M$ a finitely generated $B$-module.
Suppose that $(A,I)$ and $(B,IB)$ satisfy {\bf (APf)} and that $B$ is $IB$-adically Zariskian.
Then the conditions {\rm (i)} and {\rm (ii)} in {\rm \ref{prop-flatmorformal1}} are satisfied, and therefore the conditions {\rm (a)} $\sim$ {\rm (d)} in {\rm \ref{prop-flatmorformal1}} are all equivalent.
\end{prop}

\begin{proof}
Clearly, the condition (i) in \ref{prop-flatmorformal1} is satisfied.
To verify (ii), let $\mathfrak{a}\subseteq A$ be a finitely generated ideal, and consider $N=\mathfrak{a}\otimes_AM$.
As $N$ is finitely generated over $B$, it is $I$-adically separated due to \ref{prop-zariskipair2}.
\end{proof}

\begin{cor}\label{cor-flatformal2}
Let $(A,I)$ and $(B,IB)$ be as in {\rm \ref{prop-flatmorformal12cor}}.
Then $A\rightarrow B$ is flat if and only if $A_k=A/I^{k+1}\rightarrow B_k=B/I^{k+1}B$ is flat for any $k\geq 0$. \hfill$\square$
\end{cor}
\index{local criterion of flatness|)}

\subsubsection{Formal fpqc descent of `Noetherian outside $I$'}
The following proposition, which gives formal fpqc patching principle for the property `Noetherian outside $I$', will be of fundamental importance in our later discussion:
\begin{prop}[O.\ Gabber; {\rm \cite[Prop.\ 5.2.1]{FGK}}]\label{prop-fpqcdescentrigidNoetherian}\index{Gabber, O.}
Let $(A,I)\rightarrow (B,IB)$ be an adic morphism between complete pairs of finite ideal type such that for any $k\geq 0$ the induced map $A/I^{k+1}\rightarrow B/I^{k+1}B$ is faithfully flat.
Suppose $B$ is Noetherian outside $IB$ {\rm (\ref{dfn-outsideI} (1))}\index{Noetherian!Noetherian outside I@--- outside $I$}.
Then $A$ is Noetherian outside $I$.
Moreover, the map $A\rightarrow B$ is faithfully flat.
\end{prop}

\begin{proof}
To prove the first assertion, we want to show that any ideal $J\subseteq A$ is finitely generated outside $I$ (that is, the associated quasi-coherent ideal $\til{J}$ on $\Spec A$ is of finite type over $\Spec A\setminus V(I)$).
Considering an approximation of $J$ by finitely generated subideals of $J$, one has a finitely generated subideal $J_0\subseteq J$ such that $J_0B$ and $JB$ coincides outside $I$.
This means that $JB/J_0B$ is $I$-torsion; since the $I$-torsion part of $B/J_0B$ is bounded (due to \ref{thm-gabberIHES2008}), there exists $n>0$ such that $I^nJB\subseteq J_0B$.
Moreover, by the assumption we have $I^nJA/I^{k+1}\subseteq J_0A/I^{k+1}$ for any $k\geq 0$, that is,
$$
\bigcap_{k\geq 0}(I^nJ+I^{k+1})\subseteq\bigcap_{k\geq 0}(J_0+I^{k+1})=\ovl{J_0}\leqno{(\ast)}
$$
(the closure of $J_0$ in $A$).

We want to show the inclusion $I^nJ\subseteq J_0$, for this implies that $J$ and $J_0$ coincides with each other outside $I$.
By $(\ast)$ it suffices to show that $J_0$ is closed in $A$.
As $J_0$ is $I$-adically complete (\ref{cor-qconsistency111}), it is enough to show that the subspace topology on $J_0$ induced from the $I$-adic topology on $A$ is $I$-adic. 

To show this, we first use the condition {\bf (AP)} for $B$ to deduce that for any $i>0$ there exists $m=m(i)>i$ such that
$$
J_0B\cap I^{m(i)}B\subseteq I^iJ_0B.
$$
Again by the assumption we have
$$
J_0\cap I^{m(n)}\subseteq\ovl{I^nJ_0}.\leqno{(\ast\ast)}
$$
We want to show that the left-hand side is actually contained in $I^nJ_0$.
Suppose $x$ is in the left-hand side.
By $(\ast\ast)$ the element $x$ is decomposed as
$$
x=z_1+x_1,\qquad(z_1\in I^nJ_0,\ x_1\in I^{m(n+1)}).
$$
As $x_1$ also lies in $J_0$, we again apply $(\ast\ast)$ to decompose $x_1$ into the sum of $z_2\in I^{n+1}J_0$ and $x_2\in I^{m(n+2)}$.
One can repeat this procedure to get sequences $\{z_k\}$ and $\{x_k\}$ such that
$$
x_k=z_{k+1}+x_{k+1},\qquad(z_{k+1}\in I^{n+k}J_0,\ x_{k+1}\in I^{m(n+k+1)}).
$$
Hence $x$ is equal to the infinite series $\sum_{k\geq 1}z_k$, which converges in the $I$-adically complete $I^nJ_0$ (\ref{cor-qconsistency111}).
This means that $x\in I^nJ_0$, which shows the first assertion of the proposition.
The other assertion follows from \ref{prop-flatmorformal12cor}, \ref{thm-gabberIHES2008}, and \ref{cor-flatmorformal11}.
\end{proof}

\subsection{Restricted formal power series ring}\label{sub-powerseries}
\index{restricted formal power series|(}
Let $(A,I)$ be a complete pair of finite ideal type, and $M$ an $I$-adically complete $A$-module.
We denote by
$$
M\dl X_1,\ldots,X_n\dr
$$
the set of all formal power series $f=\sum a_{i_1,\ldots,i_n}X^{i_1}_1\cdots X^{i_n}_n$ with all coefficients $a_{i_1,\ldots,i_n}$ in $M$ such that for any $m\geq 1$ there exists $N\geq 1$ such that $a_{i_1,\ldots,i_n}\in I^mM$ whenever $i_1+\cdots+i_n>N$.
This is an $I$-adically complete $A$-module.
In particular, $A\dl X_1,\ldots,X_n\dr$ is an $I$-adically complete $A$-algebra, the so-called {\em restricted formal power series ring}\index{restricted formal power series!restricted formal power series ring@--- ring} (\cite[$\mathbf{0}_{\mathbf{I}}$, (7.5.1)]{EGA}), which is isomorphic to the $I$-adic completion of the polynomial ring $A[X_1,\ldots,X_n]$ and yields the complete pair $(A\dl X_1,\ldots,X_n\dr,IA\dl X_1,\ldots,X_n\dr)$ of finite ideal type.
Clearly, we have
$$
M\widehat{\otimes}_AA\dl X_1,\ldots,X_n\dr\cong M\dl X_1,\ldots,X_n\dr,
$$
and if $Y_1,\ldots,Y_m$ is another set of indeterminacies, then
$$
M\dl X_1,\ldots,X_n\dr\widehat{\otimes}_AA\dl Y_1,\ldots,Y_m\dr\cong M\dl X_1,\ldots,X_n,Y_1,\ldots,Y_m\dr.
$$

\begin{dfn}\label{dfn-topfinigen}{\rm 
Let $(A,I)$ be a complete pair of finite ideal type.
An $I$-adically complete $A$-algebra $B$ is said to be a {\em topologically finitely generated}\index{finitely generated!topologically finitely generated@topologically ---} $A$-algebra or an $A$-algebra {\em topologically of finite type} if $B$ is isomorphic to an $A$-algebra of the form $A\dl X_1,\ldots,X_n\dr/\mathfrak{a}$. 
If, moreover, $\mathfrak{a}$ is finitely generated, we say that $B$ is a {\em topologically finitely presented}\index{finitely presented!topologically finitely presented@topologically ---} $A$-algebra or an $A$-algebra {\em topologically of finite presentation}.}
\end{dfn}

\danger{In this book, as indicated in the above definition, any topologically finitely generated/presented algebras are assumed to be {\em complete} (and hence the ideal $\mathfrak{a}$ as above has to be closed due to \ref{cor-qconsistency1111}.)
As we will see soon below (\ref{rem-topadhesive}), this hypothesis is, in practice, not restrictive.}

\begin{prop}\label{prop-formalnot3}
Let $(A,I)$ be a complete pair with a finitely generated ideal $I\subseteq A$, and $B$ an $I$-adically complete $A$-algebra.
Then the following conditions are equivalent$:$
\begin{itemize}
\item[{\rm (a)}] $B$ is topologically finitely generated over $A;$
\item[{\rm (b)}] $B/IB$ is an $(A/I)$-algebra of finite type.
\end{itemize}
\end{prop}

\begin{proof}
The implication (a) $\Rightarrow$ (b) is clear. 
Suppose (b) holds, and take $c_1,\ldots,c_n\in B$ whose images in $B/IB$ generates $B/IB$ as an $(A/I)$-algebra.
Since $B$ is $IB$-adically complete, there exists a morphism $A'=A\dl X_1,\ldots,X_n\dr\rightarrow B$ that maps $X_i$ to $c_i$ for $i=1,\ldots,n$.
By \ref{prop-complpair1} this map is surjective.
\end{proof}

\begin{dfn}\label{dfn-topologicallyuniversallynoetherian}{\rm 
Let $(A,I)$ be a pair of finite ideal type.
We say that $A$ is {\em topologically universally Noetherian outside $I$}\index{Noetherian!Noetherian outside I@--- outside $I$!topologically universally Noetherian outside I@topologically universally --- ---} if it is Noetherian outside $I$ (\ref{dfn-outsideI} (1)) and for any $n\geq 0$, $\widehat{A}\dl X_1,\ldots,X_n\dr$ (the $I$-adic completion of $A[X_1,\ldots,X_n]$) is Noetherian outside $I\widehat{A}\dl X_1,\ldots,X_n\dr$.}
\end{dfn}

It follows from the definition that, if $A$ is topologically universally Noetherian outside $I$, then any topologically finitely generated $\widehat{A}$-algebra $B$ is topologically universally Noetherian outside $IB$. 
\begin{prop}\label{rem-topadhesive}
Let $(A,I)$ be a complete pair of finite ideal type, and suppose that $A$ is topologically universally Noetherian outside $I$\index{Noetherian!Noetherian outside I@--- outside $I$!topologically universally Noetherian outside I@topologically universally --- ---}. 
Then for any $n\geq 0$ any ideal $\mathfrak{a}\subseteq A\dl X_1,\ldots,X_n\dr$ is closed $($hence any ring of the form $A\dl X_1,\ldots,X_n\dr/\mathfrak{a}$ is topologically finitely generated over $A$ in our sense$)$.
\end{prop}

\begin{proof}
Let $B=A\dl X_1,\ldots,X_n\dr$.
Since $B$ is Noetherian outside $IB$, it satisfies {\bf (AP)} (\ref{thm-gabberIHES2008}).
Then by \ref{cor-propARconseq1-2} any ideal of $A\dl X_1,\ldots,X_n\dr$ is closed.
\end{proof}

\begin{prop}\label{prop-formalnot31}
Let $(A,I)$ be a complete pair of finite ideal type, and suppose that $A$ is topologically universally Noetherian outside $I$\index{Noetherian!Noetherian outside I@--- outside $I$!topologically universally Noetherian outside I@topologically universally --- ---}.
Let $B$ be an $I$-adically complete $A$-algebra.
Then if $B/I^nB$ is an $A/I^n$-algebra of finite presentation for any $n\geq 1$, $B$ is topologically finitely presented over $A$.
\end{prop}

\begin{proof}
By \ref{prop-formalnot3} we already know that $B$ is isomorphic to an $A$-algebra of the form $A\dl X_1,\ldots,X_n\dr/\mathfrak{a}$.
It then follows from \ref{cor-ARconseq2} that $B$ is finitely presented as an $A\dl X_1,\ldots,X_n\dr$-module.
\end{proof}

\begin{prop}\label{prop-powerbtar}
Let $(A,I)$ be a complete pair of finite ideal type, and suppose that $A$ is Noetherian outside $I$\index{Noetherian!Noetherian outside I@--- outside $I$}.

{\rm (1)} Any finitely generated $A$-module is $I$-adically complete, and for any finitely generated $A$-odule $M$ the canonical map 
$$
M\otimes_AA\dl X_1,\ldots,X_n\dr\longrightarrow M\dl X_1,\ldots,X_n\dr
$$
is an isomorphism.

{\rm (2)} The map $A\rightarrow A\dl X_1,\ldots,X_n\dr$ is flat.
\end{prop}

\begin{proof}
First notice that by \ref{thm-gabberIHES2008} $(A,I)$ satisfies {\bf (AP)} and {\bf (BT)}.
To show (1), we first notice that the completeness of $M$ is already proved in \ref{cor-gabberIHES2008}.
To show the other statement in (1), we first assume that $M$ is finitely presented.
Take a presentation $A^{\oplus q}\rightarrow A^{\oplus p}\rightarrow M\rightarrow 0$, and let $K$ be the image of $A^{\oplus q}\rightarrow A^{\oplus p}$.
Then the subspace topology on $K$ coincides with the $I$-adic topology (by the property {\bf (AP)}), and $K$ is complete.
Hence we have the exact sequence:
$$
0\longrightarrow K\dl X_1,\ldots,X_n\dr\longrightarrow A\dl X_1,\ldots,X_n\dr^{\oplus p}\longrightarrow M\dl X_1,\ldots,X_n\dr\longrightarrow 0.
$$
Thus we get an exact sequence
$$
A\dl X_1,\ldots,X_n\dr^{\oplus q}\longrightarrow A\dl X_1,\ldots,X_n\dr^{\oplus p}\longrightarrow M\dl X_1,\ldots,X_n\dr\longrightarrow 0.
$$
On the other hand, since $A^{\oplus p}\otimes_AA\dl X_1,\ldots,X_n\dr\cong A\dl X_1,\ldots,X_n\dr^{\oplus p}$, one can show the desired isomorphism by an argument similar to that in the proof of \ref{prop-ARconseq0}.

In general, one first observe that, since $M$ is finitely presented outside $I$, we have a surjective map $N\rightarrow M$ from a finitely presented $A$-module such that the kernel is $I$-torsion (\ref{lem-finitelypresentedoutside} (2)).
We can apply the argument as in the proof of \ref{prop-btarf1} (1), once we know that the assertion is true for bounded $I$-torsion modules (here we use the property {\bf (BT)}).
But the assertion in this case is easy, for both $M\otimes_AA\dl X_1,\ldots,X_n\dr$ and $M\dl X_1,\ldots,X_n\dr$ are isomorphic to $M[X_1,\ldots,X_n]$.

The second assertion can be shown by an argument similar to that in the proof of \ref{prop-btarf1} (2).
\end{proof}

\begin{prop}\label{prop-btarfflat1}
Let $(A,I)$ be a complete pair of finite ideal type, and consider the restricted power series ring $A\dl X\dr$ on one variable.
Suppose that $A$ is Noetherian outside $I$.
Then $A\dl X\dr$ is flat over $A[X]$.
\end{prop}

\begin{proof}
First we notice that both $A\dl X\dr$ and $A[X]$ are flat over $A$ (\ref{prop-powerbtar} (2)).
Since $A\dl X\dr/IA\dl X\dr\cong(A/I)[X]$, in view of \ref{prop-glueingofflatnessrelative} we only need to show the flatness over the points outside $I$.
To this end, take $x\in U=\Spec A\setminus V(I)$, and set $R=\O_{U,x}$, which is a Noetherian local ring.
Then it suffices to show that the map 
$$
R[X]\longrightarrow A\dl X\dr\otimes_AR\leqno{(\ast)}
$$
is flat.

\medskip
{\sc Claim.} If $(\ast)$ is flat outside the maximal ideal $\m_R$ of $R$, then $(\ast)$ is flat.

\medskip
Indeed, again applying \ref{prop-glueingofflatnessrelative} with $A$ replaced by $R$ and the ideal $I$ by $\m_R$, we find that we only have to show that the induced map $(\ast)\otimes_Rk$ is flat, where $k$ is the residue field of $R$.
Let $\mathfrak{p}\subseteq A$ be the prime ideal corresponding to $x$.
By \ref{prop-powerbtar} (1) we have $A\dl X\dr\otimes_A(A/\mathfrak{p})\cong (A/\mathfrak{p})\dl X\dr$.
Hence $A\dl X\dr\otimes_Ak$ is regarded as a subring of $k[\![X]\!]$.
Since $k[\![X]\!]$ is torsion free as a $k[X]$-module, $A\dl X\dr\otimes_Ak$ is also torsion free.
But since $k[X]$ is PID, this means that $A\dl X\dr\otimes_Ak$ is flat over $k[X]$, as desired.

Now we want to show the flatness of $(\ast)$ by induction with respect to $\dim(R)$.
If $\dim(R)=0$, then $\Spec R\setminus V(\m_R)$ is empty, and thus the premise of {\sc Claim} is automatically fulfilled.
Hence the desired flatness follows from the claim.

If $\dim(R)>0$, by induction with respect to $\dim(R)$ the map $(\ast)$ with $R$ replaced by $R_{\mathfrak{p}}$ for any non-maximal prime ideal $\mathfrak{p}$ (that is, the local ring at the point $y\in U$ corresponding to $\mathfrak{p}$) is flat.
In particular, $(\ast)$ is flat outside $\m_R$.
Hence again by {\sc Claim} we deduce that $(\ast)$ is flat.
\end{proof}

\begin{thm}\label{thm-btarf1}
Let $(A,I)$ be a complete pair of finite ideal type.
Suppose that $A$ is topologically universally Noetherian outside $I$\index{Noetherian!Noetherian outside I@--- outside $I$!topologically universally Noetherian outside I@topologically universally --- ---}.
Then for any $n\geq 0$ the pair $(A\dl X_1,\ldots,X_n\dr,IA\dl X_1,\ldots,X_n\dr)$ satisfies {\bf (UBT)} and {\bf (UAP)} {\rm (\S\ref{subsub-APfirstprop})}.
\end{thm}

\begin{proof}
By \ref{cor-bt32}, \ref{cor-bt32AP}, and \ref{prop-APandBToutsidenoetherian} we only need to check that a pair of the form 
$$
(A\dl X_1,\ldots,X_r\dr[Y_1,\ldots,Y_s],IA\dl X_1,\ldots,X_r\dr[Y_1,\ldots,Y_s])
$$
satisfies {\bf (BT)}.
We claim that the map 
$$
A\dl X_1,\ldots,X_r\dr[Y_1,\ldots,Y_s]\longrightarrow A\dl X_1,\ldots,X_r,Y_1,\ldots,Y_s\dr
$$
is flat.
The case $s=1$ follows from \ref{prop-btarfflat1}.
In general, this follows by induction with respect to $s$ from the following factorization:
\begin{equation*}
\begin{aligned}
A\dl X_1,\ldots,X_r\dr[Y_1,\ldots,Y_s]\longrightarrow A\dl X_1,\ldots,X_r,Y_1,\ldots,Y_{s-1}\dr[Y_s]\\
\longrightarrow A\dl X_1,\ldots,X_r,Y_1,\ldots,Y_s\dr\rlap{.}
\end{aligned}
\end{equation*}
Then we apply \ref{thm-gabberIHES2008} and \ref{prop-btarfflat2} to deduce the desired result.
\end{proof}
\index{restricted formal power series|)}

\subsection{Adhesive pairs}\label{sub-adhesive}
\subsubsection{Adhesive pairs and universally adhesive pairs}\label{subsub-adhesivedfn}
\index{pair!adhesive pair@adhesive ---|(}\index{adhesive!adhesive pair@--- pair|(}
\index{pair!adhesive pair@adhesive ---!pseudo adhesive pair@pseudo-{---} ---|(}\index{adhesive!adhesive pair@--- pair!pseudo adhesive pair@pseudo-{---} ---|(}
\index{pseudo-adhesive|(}
\begin{dfn}\label{dfn-adhesive}{\rm 
A pair $(A,I)$ of finite ideal type is said to be {\em pseudo-adhesive}\index{adhesive!adhesive pair@--- pair!pseudo adhesive pair@pseudo-{---} ---} if the following conditions are satisfied:
\begin{itemize}
\item[{\rm (a)}] $A$ is Noetherian outside $I$ {\rm (\ref{dfn-outsideI})}\index{Noetherian!Noetherian outside I@--- outside $I$}$;$
\item[{\rm (b)}] $(A,I)$ satisfies the condition {\bf (BT)} in \S\ref{subsub-BTfirstprop}\index{torsion!bounded torsion@bounded ---!bounded torsion condition@--- --- condition}\index{bounded torsion condition}.
\end{itemize}
A pseudo-adhesive pair $(A,I)$ is said to be {\em adhesive} if it satisfies the following stronger condition than (b)\index{adhesive!adhesive pair@--- pair}:
\begin{itemize}
\item[{\rm (c)}] the $I$-torsion part $M_{\Itor}$ of any finitely generated $A$-module $M$ is finitely generated.
\end{itemize}}
\end{dfn}

It is clear that the pseudo-adhesiveness and the adhesiveness depend only on the topology on $A$, and not on the ideal $I$ itself.
In the sequel we often say that the ring $A$ is {\em $I$-adically adhesive}\index{adhesive!Iadically adhesive@$I$-adically ---} (resp.\ {\em $I$-adically pseudo-adhesive})\index{adhesive!Iadically adhesive@$I$-adically ---!Iadically pseudo adhesive@--- pseudo-{---}} to mean that the pair $(A,I)$ is adhesive (resp.\ pseudo-adhesive).
Notice that the theorem of Gabber (\ref{thm-gabberIHES2008}\index{Gabber, O.}) implies:
\begin{prop}\label{prop-gabberIHES2008pseudoadhesive}
A complete pair $(A,I)$ of finite ideal type is pseudo-adhesive if and only if $A$ is Noetherian outside $I$. \hfill$\square$
\end{prop}

\begin{prop}\label{prop-adhesive}
Let $(A,I)$ be a pair of finite ideal type.
Then the following conditions are equivalent$:$
\begin{itemize}
\item[{\rm (a)}] the pair $(A,I)$ is adhesive$;$
\item[{\rm (b)}] for any finitely generated $A$-module $M$, $M/M_{\Itor}$ is finitely presented$;$
\item[{\rm (c)}] for any finitely generated $A$-module $M$ and any $A$-submodule $N$ of $M$, the $I$-saturation\index{saturation!Isaturation@$I$-{---}} $\til{N}$ of $N$ in $M$ is finitely generated.
\end{itemize}
\end{prop}

\begin{proof}
We may assume that the ideal $I\subseteq A$ is finitely generated.
First we show the equivalence of (b) and (c). 
Suppose (b) holds, and let $M$ and $N$ be as in (c). 
The $I$-saturation $\til{N}$ (\S\ref{subsub-pairstorsionsandsaturation}) sits in the following exact sequence:
$$
0\longrightarrow\til{N}\longrightarrow M\longrightarrow(M/N)/(M/N)_{\Itor}\longrightarrow 0.
$$
By \cite[Chap.\ I, \S 2.8, Lemma 9]{Bourb1} we deduce $\til{N}$ is finitely generated, whence (c).
The converse ((c) $\Rightarrow$ (b)) follows from the fact that, if $\phi\colon F\rightarrow M/M_{\Itor}$ is a surjective morphism from a finitely generated $A$-module, then $\ker(\phi)$ is saturated.

Next let us show the equivalence of (a) and (c).
Suppose (c) holds, and let $M$ be a finitely generated $A$-module. 
Since $M_{\Itor}$ is $I$-saturated, it is finitely generated.
To show that $\Spec A\setminus V(I)$ is a Noetherian scheme, it suffices to show that $A[\frac{1}{a}]$ is Noetherian for any $a\in I$.
Let $J$ be an ideal of $A[\frac{1}{a}]$, and $J'$ the pull-back of $J$ by $A\rightarrow A[\frac{1}{a}]$. 
Then $J'$ is easily seen to be $I$-saturated, and we have $J'A[\frac{1}{a}]=J$. 
Since $J'$ is finitely generated, so is $J$.
Conversely, suppose (a) holds.
Let $N$ be a submodule of a finitely generated $A$-module $M$.
By \ref{lem-finitelypresentedoutside} (1) we can find a finitely generated submodule $N'$ of $N$ such that $N/N'$ is $I$-torsion.
Since we have $\til{N'}=\til{N}$, we may replace $N$ by $N'$ and hence may suppose $N$ is finitely generated.
Then the exact sequence
$$
0\longrightarrow N\longrightarrow\til{N}\longrightarrow (M/N)_{\Itor}\longrightarrow 0
$$
gives that $\til{N}$ is finitely generated.
\end{proof}

\begin{dfn}\label{dfn-universallyadhesive}{\rm 
A pair $(A,I)$ is said to be {\em universally adhesive}\index{pair!adhesive pair@adhesive ---!universally adhesive pair@universally --- ---}\index{adhesive!universally adhesive@universally ---} (resp.\ {\em universally pseudo-adhesive})\index{pair!adhesive pair@adhesive ---!universally pseudo adhesive pair@universally pseudo-{---} ---}\index{adhesive!universally pseudo adhesive@universally pseudo-{---}}\index{pseudo-adhesive!universally pseudo adhesive@universally ---} if for any $n\geq 0$ the pair $(A[X_1,\ldots,X_n],IA[X_1,\ldots,X_n])$ is adhesive (resp.\ pseudo-adhesive).}
\end{dfn}

In this situation we also say that the ring $A$ is {\em $I$-adically universally adhesive} (resp.\ {\em $I$-adically universally pseudo-adhesive}).
Notice that for a pseudo-adhesive pair $(A,I)$ to be universally pseudo-adhesive it is necessary and sufficient that $(A,I)$ satisfies {\bf (UBT)} (cf.\ \ref{cor-bt32}).

\begin{prop}\label{prop-adhesivenessidealchange}
{\rm (1)} Let $(A,I)$ be an adhesive pair $($resp.\ universally adhesive$)$. Then for any finitely generated subideal $J\subseteq I$ the pair $(A,J)$ is adhesive $($resp.\ universally adhesive$)$.

{\rm (2)} Let $I=(a_1,\ldots,a_n)\subseteq A$ be a finitely generated ideal, and suppose that $(A,a_i)$ is adhesive $($resp.\ universally adhesive$)$ for any $i=1,\ldots,n$.
Then $(A,I)$ is adhesive $($resp.\ universally adhesive$)$.
\end{prop}

\begin{proof}
Clearly, it is enough to prove the assertions only in the `adhesive' case.

(1) Any finitely generated $J$-torsion free $A$-module is $I$-torsion free, which verifies the condition (b) of \ref{prop-adhesive}.

(2) By induction with respect to $n$ we may reduce to the case $I=(a,b)$; suppose $(A,a)$ and $(A,b)$ are adhesive.
It is clear that $A$ is Noetherian outside $I$, for it is Noetherian outside $(a)$ and $(b)$.
For a finitely generated $A$-module $M$, the $a$-torsion part $M_{\ator}$ is finitely generated.
Hence its $b$-torsion part $(M_{\ator})_{\btor}$, which is nothing but $M_{\Itor}$, is finitely generated, whence verifying the condition (c) in \ref{dfn-adhesive}.
\end{proof}

\begin{prop}\label{prop-adhesive0}
{\rm (1)} If $(A,I)$ and $(B,J)$ are adhesive $($resp.\ pseudo-adhesive, resp.\ universally adhesive, resp.\ universally pseudo-adhesive$)$ pairs, then $(A\times B,I\times J)$ is adhesive $($resp.\ pseudo-adhesive, resp.\ universally adhesive, resp.\ universally pseudo-adhesive$)$.

{\rm (2)} Let $(A,I)$ be a pair, and $B$ a faithfully flat $A$-algebra.
If $(B,IB)$ is adhesive $($resp.\ pseudo-adhesive, resp.\ universally adhesive, resp.\ universally pseudo-adhesive$)$, then so is $(A,I)$.
\end{prop}

\begin{proof}
(1) This can be shown by an argument similar to that in the proof of \ref{prop-bt33} (1) in view of the fact that $A\times B$ is Noetherian outside $I\times J$ if and only if $A$ is Noetherian outside $I$ and $B$ is Noetherian outside $J$.

(2) It suffices to verify the assertion only in `adhesive' and `pseudo-adhesive' cases.
By \ref{prop-finiteidealtype2} we see that $(A,I)$ is of finite ideal type.
Since $\Spec B\setminus V(IB)\rightarrow\Spec A\setminus V(I)$ is faithfully flat, it follows that $A$ is Noetherian outside $I$.
If $(B,IB)$ is pseudo-adhesive, then it follows from \ref{prop-bt33} (2) that $(A,I)$ is pseudo-adhesive.
Suppose $(B,IB)$ is adhesive, and let $M$ be an $I$-torsion free finitely generated $A$-module.
Take a surjection $A^{\oplus m}\rightarrow M$, and let $K$ be the kernel.
Since $B$ is $A$-flat, $M\otimes_AB$ is $IB$-torsion free, and we have the exact sequence
$$
0\longrightarrow K\otimes_AB\longrightarrow B^{\oplus m}\longrightarrow M\otimes_AB\longrightarrow 0.
$$
By the assumption $K\otimes_AB$ is a finitely generated $B$-module (here we used \cite[Chap.\ I, \S 2.8, Lemma 9]{Bourb1}).
By \cite[Chap.\ I, \S3.1, Prop.\ 2]{Bourb1} we deduce that $K$ is finitely generated and hence that $(A,I)$ is adhesive.
\end{proof}

\begin{prop}\label{prop-adhesive1}
Let $(A,I)$ be an adhesive $($resp.\ a pseudo-adhesive$)$ pair.

{\rm (1)} For any multiplicative subset $S\subseteq A$ the induced pair $(S^{-1}A,IS^{-1}A)$ is adhesive $($resp.\ pseudo-adhesive$)$.

{\rm (2)} For any quasi-finite $A$-algebra $B$ the induced pair $(B,IB)$ is adhesive $($resp.\ pseudo-adhesive$)$.
\end{prop}

\begin{proof}
(1) As $\Spec S^{-1}A\setminus V(S^{-1}I)$ is clearly Noetherian, it suffices to check that for any finitely generated $S^{-1}A$-module $M$, its $S^{-1}I$-torsion part is finitely generated (resp.\ bounded).
Let $x_1,\ldots,x_n$ be a generator of $M$.
Set $M'=Ax_1+\cdots+Ax_n$; we have $M'\otimes_AS^{-1}A=M$.
As one can verify easily, $M_{S^{-1}\Itor}$ coincides with $M'_{\Itor}\otimes_AS^{-1}A$.
Since $M'_{\Itor}$ is finitely generated (resp.\ bounded $I$-torsion), $M_{S^{-1}\Itor}$ is finitely generated over $S^{-1}A$ (resp.\ bounded $S^{-1}I$-torsion).

(2) The assertion is clear if $B$ is finite over $A$.
In general, we apply Zariski's Main Theorem (\cite[$\mathbf{IV}$, (18.12.13)]{EGA}) to reduce to this case, using (1) and \ref{prop-adhesive0} as follows: $\Spec B\rightarrow\Spec A$ is the composition of an open immersion followed by a finite morphism.
Hence there exists a finite open covering $\Spec A=\bigcup_{i\in I}\Spec A_i$ such that for each $i\in I$
\begin{itemize}
\item $A_i$ is of the form $S^{-1}_iA$ by a multiplicative subset $S_i\subseteq A$;
\item there exists a finite $A$-algebra $B'$ such that $B_i=B\otimes_AA_i$ is isomorphic to an $A$-algebra of the form $T^{-1}_iB'$ by a multiplicative subset $T_i\subseteq B'$.
\end{itemize}
By (1) each $(A_i,IA_i)$ is adhesive (resp.\ pseudo-adhesive), and hence $(B_i,IB_i)$ is adhesive (resp.\ pseudo-adhesive).
Now by \ref{prop-adhesive0} (1) and (2) applied to $\coprod\Spec B_i\rightarrow\Spec B$, it follows that $(B,IB)$ is adhesive (resp.\ pseudo-adhesive).
\end{proof}

\begin{prop}\label{prop-uadhesive1}
Let $(A,I)$ be a universally adhesive $($resp.\ universally pseudo-adhesive$)$ pair.

{\rm (1)} For any multiplicative subset $S\subseteq A$ the induced pair $(S^{-1}A,S^{-1}I)$ is universally adhesive $($resp.\ universally pseudo-adhesive$)$.

{\rm (2)} For any $A$-algebra $B$ of finite type the induced pair $(B,IB)$ is universally adhesive $($resp.\ universally pseudo-adhesive$)$.
\end{prop}

The second assertion says that universally-adhesiveness (resp.\ universally-pseudo-adhesiveness) is stable under finite type extensions.
\begin{proof}
(1) follows easily from \ref{prop-adhesive1} (1).
To show (2), notice first that any polynomial ring over $A$ is universally adhesive (resp.\ universally pseudo-adhesive) and apply \ref{prop-adhesive1} (2).
\end{proof}

\begin{prop}\label{prop-adhesivered}
Let $(A,I)$ be a pair of finite ideal type.
The following two conditions are equivalent$:$
\begin{itemize}
\item[{\rm (a)}] $(A,I)$ is adhesive $($resp.\ pseudo-adhesive, resp.\ universally adhesive, resp.\ universally pseudo-adhesive$);$
\item[{\rm (b)}] $A_{\Itor}$ is finitely generated over $A$ $($resp.\ bounded $I$-torsion, resp.\ finitely generated over $A$, resp.\ bounded $I$-torsion$)$ and $(A/A_{\Itor},I(A/A_{\Itor}))$ is adhesive $($resp.\ pseudo-adhesive, resp.\ universally adhesive, resp.\ universally pseudo-adhesive$)$. 
\end{itemize}
\end{prop}

\begin{proof}
Since $(A/A_{\Itor})[X]\cong A[X]/A[X]_{\Itor}$, it suffices to show the assertion in `adhesive' and `pseudo-adhesive' cases.
As (a) $\Rightarrow$ (b) is obvious by \ref{prop-adhesive1} (2), we show the converse.
First we claim that the scheme $\Spec A\setminus V(I)$ is Noetherian.
For any $f\in I$ the canonical map $A\rightarrow A_f$ factors through $A/A_{\Itor}$, whence $(A/A_{\Itor})_f\cong A_f$.
This means that $A_f$ for any $f\in I$ is Noetherian and hence that $\Spec A\setminus V(I)$ is Noetherian.

Let $M$ be a finitely generated $A$-module, and set $M'=M\otimes_A(A/A_{\Itor})=M/A_{\Itor}M$.
We have the exact sequence
$$
0\longrightarrow A_{\Itor}M\longrightarrow M_{\Itor}\longrightarrow M'_{\Itor}\longrightarrow 0.
$$
Since $A_{\Itor}M$ and $M'_{\Itor}$ are finitely generated over $A$ in `adhesive' case or are bounded $I$-torsion in `pseudo-adhesive' case, so is $M_{\Itor}$.
\end{proof}

\begin{prop}\label{prop-adhesive4}
Let $(A,I)$ be a pair of finite ideal type.
Consider an inductive system $\{B_i,f_{ij}\colon B_i\rightarrow B_j\}$ of $A$-algebras indexed by a directed set $J$, and set $B=\varinjlim_{i\in I}B_i$.
Suppose the following conditions are satisfied$:$
\begin{itemize}
\item[{\rm (a)}] $(B_i,IB_i)$ is adhesive $($resp.\ pseudo-adhesive, resp.\ universally adhesive, resp.\ universally pseudo-adhesive$)$ for each $i;$
\item[{\rm (b)}] for each pair of indices $i,j$ such that $i\leq j$ the morphism $f_{ij}$ is flat$;$
\item[{\rm (c)}] $\Spec B\setminus V(IB)$ is Noetherian.
\end{itemize}
Then the pair $(B,IB)$ is adhesive $($resp.\ pseudo-adhesive, resp.\ universally adhesive, resp.\ universally pseudo-adhesive$)$.
\end{prop}

\begin{proof}
It suffices to check the proposition for `adhesive' and `pseudo-adhesive' cases.
Let $M$ be a finitely generated $B$-module, and consider an exact sequence of the form 
$$
0\longrightarrow N\longrightarrow B^{\oplus m}\longrightarrow M\longrightarrow 0.
$$
Let $\til{N}$ be the $I$-saturation of $N$ in $B^{\oplus m}$.
Then we have $\til{N}/N\cong M_{\Itor}$.
Since $N$ is finitely generated outside $I$, one can find a finitely generated submodule $N'\subseteq N$ such that its $I$-saturation in $B^{\oplus m}$ coincides with $\til{N}$ (\ref{lem-finitelypresentedoutside} (1)).

Take an index $i\in J$ and a finitely generated $B_i$-submodule $N'_i$ of $B^{\oplus m}_i$ such that $N'_i\otimes_{B_i}B=N'$ (note that by (b) the map $B_i\rightarrow B$ is flat).
Take the $I$-saturation $\til{N'_i}$ of $N_i$ in $B^{\oplus m}_i$. 
Since $\til{M'_i}=B^{\oplus m}_i/\til{N'_i}$ is an $I$-torsion free $B_i$-module isomorphic outside $I$ to $M'_i=B^{\oplus m}_i/N'_i$, it follows from the flatness of $B_i\rightarrow B$ that $\til{M'_i}\otimes_{B_i}B$ is an $I$-torsion free $B$-module isomorphic outside $I$ to $M'=M'_i\otimes_{B_i}B\cong B^{\oplus m}/N'$ and hence also isomorphic outside $I$ to $M$.
This implies that $\til{M'_i}\otimes_{B_i}B$ coincides with $\til{M}=B^{\oplus m}/\til{N}$ and hence that $\til{N'_i}\otimes_{B_i}B=\til{N}$.
Now the assertion in each case follows from the surjection 
$$
(\til{N'_i}/N'_i)\otimes_{B_i}B\cong\til{N}/N'\longrightarrow\til{N}/N\cong M_{\Itor}
$$
and that $\til{N'_i}/N'_i\cong(M'_i)_{\Itor}$.
\end{proof}

One can prove the following proposition by an argument similar to that in \ref{prop-btarfflat2AP}, using \ref{prop-adhesive0}:
\begin{prop}\label{prop-adhesive3}
Let $(A,I)$ be a pair of finite ideal type.
Suppose the following conditions are satisfied$:$
\begin{itemize}
\item[{\rm (a)}] $(\widehat{A},I\widehat{A})$ is adhesive $($resp.\ pseudo-adhesive$)$, where $\widehat{A}$ is the $I$-adic completion of $A;$
\item[{\rm (b)}] $A\rightarrow\widehat{A}$ is flat$;$ 
\item[{\rm (c)}] $A$ is Noetherian outside $I$\index{Noetherian!Noetherian outside I@--- outside $I$}.
\end{itemize}
Then the pair $(A,I)$ is adhesive $($resp.\ pseudo-adhesive$)$.\hfill$\square$
\end{prop}

\begin{rem}\label{rem-compl}{\rm 
We do not know whether the converse of \ref{prop-adhesive3} holds or not, that is: if $(A,I)$ is adhesive (resp.\ pseudo-adhesive), then is $(\widehat{A},I\widehat{A})$ adhesive (resp.\ pseudo-adhesive), too?}
\end{rem}

\subsubsection{Some examples}\label{subsub-adhesiveexam}
Here we collect a few examples of adhesive pairs, which will be of particular importance in our later arguments.

\begin{exa}[Adhesive pairs of type (N)]\label{exa-adhesivetypeN}{\rm 
Needless to say, any Noetherian ring $A$ is $I$-adically universally adhesive for any ideal $I\subseteq A$.
Notice that a ring $A$ is $1$-adically adhesive if and only if $A$ is an Noetherian ring.}
\end{exa}

\begin{exa}[Adhesive pairs of type (V)]\label{exa-adhesivetypeV}{\rm 
One of the most interesting examples of adhesive pairs is given by a pair $(V,a)$ consisting of a valuation ring\index{valuation!valuation ring@--- ring!a-adically separated valuation ring@$a$-adically separated --- ---} (of arbitrary height) and an element $a\in\m_V\setminus\{0\}$ such that $V$ is $a$-adically separated (cf.\ \S\ref{sub-adicsepval}).
}
\end{exa}

Let us prove the last-mentioned fact:
\begin{prop}\label{prop-exaadhesiveval}
Let $V$ be a valuation ring, and $a\in\m_V\setminus\{0\}$.
Then the following conditions are equivalent$:$
\begin{itemize}
\item[{\rm (a)}] $V$ is $a$-adically adhesive$;$
\item[{\rm (b)}] $V$ is $a$-adically pseudo-adhesive$;$
\item[{\rm (c)}] $V$ is $a$-adically separated$;$
\item[{\rm (d)}] $V[\frac{1}{a}]$ is a field $(=\Frac(V))$.
\end{itemize}
\end{prop}

\begin{proof}
The implication (a) $\Rightarrow$ (b) is clear.
Let us show (b) $\Rightarrow$ (c). 
Suppose that $V$ is $a$-adically pseudo-adhesive, and consider the ideal $J=\bigcap_{n\geq 1}(a^n)$.
Take a finitely generated subideal $J_0\subseteq J$ such that $J/J_0$ is $a$-torsion.
As $J/J_0$ is contained in $V/J_0$, $J/J_0$ is bounded $a$-torsion.
This means that there exists $n\geq 0$ such that $a^nJ\subseteq J_0$; but since $a^nJ=J$, it follows that $J$ itself is finitely generated and hence is principal.
Then we easily see that $J=(0)$.

We have already shown (c) $\Rightarrow$ (d) in \ref{prop-sep}.
To show (d) $\Rightarrow$ (a), first notice that the condition (d) implies that any non-zero element $b\in V\setminus\{0\}$ divides some power $a^n$ ($n\geq 1$) of $a$.
Hence for any $V$-module $M$, $M$ is torsion free if and only if it is $a$-torsion free (that is, $M_{\ator}=0$).
In particular, any $a$-torsion free finitely generated $V$-module $M$ is $V$-flat and hence is a free $V$-module. 
\end{proof}

We will see later in \S\ref{sub-convpreadh} that, if $V$ as above is, moreover, $a$-adically complete, then for any $n\geq 0$ the pair $(V\dl X_1,\ldots,X_n\dr,a)$ is universally adhesive (\ref{cor-convadh}).

\subsubsection{Preservation of adicness}\label{subsub-adhesiveAR}
By \ref{prop-APandBToutsidenoetherian} and \ref{prop-ARbt1} we have:
\begin{prop}\label{prop-AR}
Let $(A,I)$ be a pseudo-adhesive pair.
Then $(A,I)$ satisfies the condition {\bf (AP)}.
If, moreover, $I$ is principal, then $(A,I)$ satisfies {\bf (AR)}.\hfill$\square$
\end{prop}

Hence by what we have already seen before, a pseudo-adhesive pair $(A,I)$ enjoys the following properties: 
\begin{itemize}
\item the functor $M\mapsto\widehat{M}$ by $I$-adic completion on the full subcategory of $\Mod_A$ consisting of finitely generated $A$-modules is exact (\ref{prop-qconsistency1-2p-2});
\item $M\otimes_A\widehat{A}\cong\widehat{M}$ for any finitely generated $A$-module $M$ (\ref{prop-btarf1} (1));
\item the $I$-adic completion map $A\rightarrow\widehat{A}$ is flat (\ref{prop-btarf1} (1)).
\end{itemize}
If, moreover, $(A,I)$ is complete, then (\S\ref{subsub-ARconsequence}):
\begin{itemize}
\item any finitely generated $A$-module is $I$-adically complete;
\item any $A$-submodule $N$ of a finitely generated $A$-module $M$ is closed in $M$ and $I$-adically complete.
\end{itemize}

\subsubsection{Topologically universally adhesive pairs}\label{subsub-adhesivetua}
\begin{dfn}\label{dfn-topadhesive}{\rm 
We say that a pair $(A,I)$ is {\em topologically universally adhesive}\index{adhesive!universally adhesive@universally ---!topologically universally adhesive@topologically --- --- (t.u.\ adhesive)}\index{pair!adhesive pair@adhesive ---!topologically universally adhesive pair@topologically universally --- ---} (resp.\ {\em topologically universally pseudo-adhesive}\index{pair!adhesive pair@adhesive ---!topologically universally pseudo adhesive pair@topologically universally pseudo-{---} ---}\index{adhesive!universally pseudo adhesive@universally pseudo-{---}!topologically universally pseudo adhesive@topologically --- --- (t.u.\ pseudo-adhesive)}\index{pseudo-adhesive!universally pseudo adhesive@universally ---!topologically universally pseudo adhesive@topologically --- --- (t.u.\ pseudo-adhesive)}) or that $A$ is {\em $I$-adically topologically universally adhesive} (resp.\ {\em $I$-adically topologically universally pseudo-adhesive}) if $(A,I)$ is universally adhesive (resp.\ universally pseudo-adhesive) and for any $n\geq 0$ the $I$-adic completion of $(A[X_1,\ldots,X_n],IA[X_1,\ldots,X_n])$ is again universally adhesive (resp.\ universally pseudo-adhesive).}
\end{dfn}

The lengthy name `topologically universally adhesive' (resp.\ `topologically universally pseudo-adhesive') will be often shortened into {\em `t.u.\ adhesive'} (resp.\ {\em `t.u.\ pseudo-adhesive'}) for brevity.\footnote{As indicated in the definition, the property `$I$-adically t.u.\ adhesive' does not imply `$I$-adically complete'; but later in {\bf \ref{ch-formal}}, \S\ref{subsub-tuarings} we define for the sake of terminological brevity the notion of {\em t.u.\ adhesive rings}, which are complete by definition. Compare to the the terminology {\em adic rings}; in \cite[$\mathbf{0}_{\mathbf{I}}$, \S7.1]{EGA} adic rings are complete by definition, whereas `$I$-adic' does not imply completeness; see \ref{rem-adicfiltrationtopology3} and the warning after {\bf \ref{ch-formal}}.\ref{dfn-admissibleringsadicrings}.}
Similarly to (pseudo-) adhesiveness and universally (pseudo-) adhesiveness, these notions depend only on the topology on $A$.
Notice that, if $(A,I)$ is t.u.\ pseudo-adhesive, then $A$ is topologically universally Noetherian outside $I$\index{Noetherian!Noetherian outside I@--- outside $I$!topologically universally Noetherian outside I@topologically universally --- ---} (\ref{dfn-topologicallyuniversallynoetherian}).
Moreover, by \ref{thm-btarf1} we have:
\begin{prop}\label{prop-tupaautomatic}
A complete pair $(A,I)$ of finite ideal type is t.u.\ pseudo-adhesive if and only if it is topologically universally Noetherian outside $I$\index{Noetherian!Noetherian outside I@--- outside $I$!topologically universally Noetherian outside I@topologically universally --- ---}. \hfill$\square$
\end{prop}

\begin{prop}\label{prop-topadhesive001}
Let $(A,I)$ be t.u.\ adhesive $($resp.\ t.u.\ pseudo-adhesive$)$. 

{\rm (1)} The completion\index{completion!completion of a pair@--- (of a pair)} $(\widehat{A},I\widehat{A})$ is t.u.\ adhesive $($resp.\ t.u.\ pseudo-adhesive$)$.

{\rm (2)} For any $A$-algebra $B$ of finite type the induced pair $(B,IB)$ is t.u.\ adhesive $($resp.\ t.u.\ pseudo-adhesive$)$.

{\rm (3)} For any $\widehat{A}$-algebra $B$ topologically of finite type {\rm (\ref{dfn-topfinigen})} the induced pair $(B,IB)$ is t.u.\ adhesive $($resp.\ t.u.\ pseudo-adhesive$)$.
\end{prop}

\begin{proof}
(1) is clear.
(2) follows easily from \ref{prop-uadhesive1} (2).
In the situation as in (3), the $I$-adic completion of a polynomial ring over $B$ is finite over an $A$-algebra of the form $\widehat{A}\dl Y_1,\ldots,Y_m\dr$ and hence is universally adhesive (resp.\ universally pseudo-adhesive) by (2).
\end{proof}

\begin{thm}\label{thm-tuaunivadhe}
Let $(A,I)$ be a complete pair of finite ideal type.
Then the following conditions are equivalent$:$
\begin{itemize}
\item[{\rm (a)}] $(A,I)$ is t.u.\ adhesive $($resp.\ t.u.\ pseudo-adhesive$);$ 
\item[{\rm (b)}] $(A\dl X_1,\ldots,X_n\dr,IA\dl X_1,\ldots,X_n\dr)$ is adhesive $($resp.\ pseudo-adhesive$)$ for any $n\geq 0$.
\end{itemize}
\end{thm}

\begin{proof}
The proof is done in a similar way to that of \ref{thm-btarf1}.
Only (b) $\Rightarrow$ (a) calls for the proof.
We need to show that a pair of the form 
$$
(A\dl X_1,\ldots,X_r\dr[Y_1,\ldots,Y_s],IA\dl X_1,\ldots,X_r\dr[Y_1,\ldots,Y_s])
$$
is adhesive (resp.\ pseudo-adhesive).
Since $A\dl X_1,\ldots,X_n\dr$ is Noetherian outside $I$, one can show as in the proof of \ref{thm-btarf1} that the map 
$$
A\dl X_1,\ldots,X_r\dr[Y_1,\ldots,Y_s]\longrightarrow A\dl X_1,\ldots,X_r,Y_1,\ldots,Y_s\dr
$$
is flat.
Then we apply \ref{prop-adhesive3} to deduce the desired result.
\end{proof}

\subsubsection{Adhesiveness and coherence}\label{subsub-coherencyadhesive}
\begin{prop}\label{prop-cohringsmodules6}
Let $(A,I)$ be an adhesive pair, and suppose $A$ is $I$-torsion free.
Then the ring $A$ is coherent\index{coherent!coherent ring@--- ring} {\rm (\ref{dfn-cohringsmodules1} (2))}.
\end{prop}

\begin{proof}
We verify the condition (e) in \ref{prop-cohringsmodules1}.
By \ref{lem-pf1} (1) and by easy homological algebra we may assume $M=A^{\oplus m}$ (for some $m\geq 0$) without loss of generality. 
Then $L$ is an $I$-torsion free finitely generated $A$-module and hence is finitely presented. 
\end{proof}

\begin{dfn}\label{dfn-cohringsmodules5}{\rm 
Let $(A,I)$ be a t.u.\ pseudo-adhesive pair.
Then the ring $A$ is said to be {\em topologically universally coherent with respect to $I$}\index{coherent!topologically universally coherent@topologically universally ---} if any topologically finitely presented $\widehat{A}$-algebra is universally coherent (\ref{dfn-universallycoherent}), where $\widehat{A}$ denotes the $I$-adic completion\index{completion!I-adic completion@$I$-adic ---} of $A$.}
\end{dfn}

\begin{prop}\label{prop-cohringsmodules52}
Let $A$ be a topologically universally coherent ring with respect to an ideal $I\subseteq A$.

{\rm (1)} Any finitely presented $A$-algebra $B$ is topologically universally coherent with respect to $IB$. 

{\rm (2)} Any topologically finitely presented $\widehat{A}$-algebra $B$ is topologically universally coherent with respect to $IB$.
\end{prop}

\begin{proof}
To show (1), write $B=A[Y_1,\ldots,Y_m]/\mathfrak{a}$ by a finitely generated $\mathfrak{a}$.
Then $(B,IB)$ is t.u.\ pseudo-adhesive (\ref {prop-topadhesive001}).
Let $R=A[Y_1,\ldots,Y_m,X_1,\ldots,X_n]$.
Then we have $R/\mathfrak{a}R=B[X_1,\ldots,X_n]$.
Since $\widehat{R}$ satisfies {\bf (AP)}, the ideal $\mathfrak{a}\widehat{R}$ is closed in $\widehat{R}$ due to \ref{cor-propARconseq1-2}, and we have $\widehat{B}\dl X_1,\ldots,X_n\dr=\widehat{R}/\mathfrak{a}\widehat{R}$. 
Hence $\widehat{B}\dl X_1,\ldots,X_n\dr$ is a topologically finitely presented $\widehat{A}$-algebra and hence is universally coherent.
By this we deduce that $B$ is topologically universally coherent with respect to $IB$.
(2) can be shown in a similar manner.
\end{proof}

\begin{prop}\label{prop-cohringsmodules51a}
Let $A$ be $I$-adically t.u.\ adhesive and topologically universally coherent with respect to $I\subseteq A$.
Then $A$ is universally coherent\index{coherent!coherent ring@--- ring!universally coherent ring@universally --- ---}.
\end{prop}

\begin{proof}
Since $(A,I)$ is t.u.\ adhesive, the canonical map $\zat{A}\rightarrow\widehat{A}$ is faithfully flat (\ref{prop-btarf1}, \ref{prop-relpair21} (2)).
We may assume that $I$ is finitely generated; set $I=(f_1,\ldots,f_r)$.
Then $\Spec A\setminus V(I)$ is the union of the open subsets $\Spec A_{f_i}$ $(i=1,\ldots,r)$.
Set $B=\zat{A}\times A_{f_1}\times\cdots\times A_{f_r}$, which is a faithfully flat algebra over $A$.
Since each $A_{f_i}$ is a Noetherian ring, it is clearly a coherent ring.
By the assumption the ring $\widehat{A}$ is coherent.
Then by \ref{prop-cohringsmodules1a} (1) and (2) we deduce that $A$ is coherent.
Since, as we saw above, `topologically universally coherent' is closed under finitely presented extension, we have the assertion.
\end{proof}

\begin{thm}\label{thm-pf02a}
Let $(A,I)$ be a pair.

{\rm (1)} Suppose $A$ is $I$-adically universally adhesive and $I$-torsion free.
Then $A$ is universally coherent\index{coherent!coherent ring@--- ring!universally coherent ring@universally --- ---}.

{\rm (2)} Suppose $A$ is $I$-adically complete, $I$-adically t.u.\ adhesive, and $I$-torsion free.
Then $A$ is topologically universally coherent with respect to $I$\index{coherent!topologically universally coherent@topologically universally ---}.
\end{thm}

Hence in (1) any finitely presented $A$-algebra $B$ is again universally coherent.
Likewise, in (2) any finitely presented algebra $B$ over a topologically finitely presented $A$-algebra is again topologically universally coherent.
Note that in both cases $B$ may not be $I$-torsion free.

\begin{proof}
(1) Let $B$ be a finitely presented $A$-algebra. 
We are going to check the condition (e) in \ref{prop-cohringsmodules1}.
Take a surjective map $B'\rightarrow B$ with the finitely generated kernel, where $B'$ is a polynomial ring over $A$.
Notice that $B'$ is $IB'$-torsion free.
Since for a $B$-module to be finitely presented over $B$ is equivalent to being finitely presented over $B'$ (cf.\ \ref{lem-pf2}), we may assume $B=B'$.
But then the assertion of this case is nothing but \ref{prop-cohringsmodules6}.

(2) Let $B$ be a topologically finitely presented $A$-algebra, and take a surjective map of the form $B'=A\dl X_1,\ldots,X_n\dr\rightarrow B$ with the finitely generated kernel.
Note that $B'$ is $IB'$-torsion free.
Then the rest of the proof goes similarly to that of (1) as above.
\end{proof}
\index{pseudo-adhesive|)}
\index{adhesive!adhesive pair@--- pair!pseudo adhesive pair@pseudo-{---} ---|)}\index{pair!adhesive pair@adhesive ---!pseudo adhesive pair@pseudo-{---} ---|)}
\index{adhesive!adhesive pair@--- pair|)}\index{pair!adhesive pair@adhesive ---|)}

\subsection{Scheme-theoretic pairs}\label{sub-schpair}
Let us mention briefly that the notion of pairs has the obvious interpretation into the language of schemes.
A (scheme-theoretic) {\em pair}\index{pair!scheme-theoretic pair@(scheme-theoretic) ---} is a couple $(X,Y)$ consisting of a scheme $X$ and a closed subscheme $Y$ of $X$.
More generally, one can consider a pair $(X,Y)$ consisting of an algebraic space $X$ and a closed subspace $Y$.
By a morphism of scheme-theoretic pairs $f\colon (X,Y)\rightarrow (Z,W)$, we mean a morphism $f\colon X\rightarrow Z$ of schemes (or algebraic spaces) such that $\mathscr{I}^n_W\O_X\subseteq\mathscr{I}_Y$ for some $n\geq 1$, where $\mathscr{I}_Y$ (resp.\ $\mathscr{I}_W$) is the defining ideal of $Y$ (resp.\ $W$) in $X$ (resp.\ $Z$).

An appropriate scheme-theoretic counterpart of complete pairs is provided by adic formal schemes:
\begin{dfn}\label{dfn-formalschemepair}
{\rm A {\em formal pair} is a couple $(\mathfrak{X},Y)$ consisting of an adic formal scheme\index{formal scheme!adic formal scheme@adic ---}\index{adic!adic formal scheme@--- formal scheme} $\mathfrak{X}$ and its closed subscheme $Y$ defined by an ideal of definition.}
\end{dfn}
(See {\bf \ref{ch-formal}}.\ref{dfn-formalschemesadicformalschemes} for the definition of {\em adic} formal schemes.)
For the practical use, however, a more handy definition of formal pairs may be the following one: A formal pair is a pair $(\mathfrak{X},\mathscr{I})$ consisting of an adic formal scheme $\mathfrak{X}$ and an ideal of definition $\mathscr{I}$.
If, moreover, we are only interested in properties that does not depend on particular choices of ideals of definition, which is most frequently the case in the sequel, we even do not have to spell out $\mathscr{I}$, and just consider an adic formal scheme $\mathfrak{X}$ itself as a scheme-theoretic counterpart of complete pairs.

The above definition of formal pairs fits in with the so-called {\em formal completion}\index{completion!formal completion@formal ---} (\cite[$\mathbf{I}$, \S10.8]{EGA}): For a scheme $X$ and a closed subscheme $Y\subseteq X$ of finite presentation, the formal completion 
$$
\widehat{X}|_Y
$$
of $X$ along $Y$ is an adic formal scheme (cf.\ \ref{prop-Iadiccompletioncomplete2}); the associated pair $(\widehat{X}|_Y,Y)$ is a formal pair.

As for henselian pairs the correct substitution is given by {\em henselian schemes}\index{scheme!henselian scheme@henselian ---}\index{henselian!henselian scheme@--- scheme} (cf.\ \cite[Chap.\ 7]{KPR}).
\begin{dfn}\label{dfn-henselschemepair}
{\rm A {\em henselian pair}\index{pair!henselian pair@henselian ---}\index{henselian!henselian pair@--- pair} is a couple $(\mathscr{X},Y)$ consisting of a henslian scheme $\mathscr{X}$ and a closed subscheme $Y$ defined by an ideal of definition.}
\end{dfn}

Similarly to the case of formal pairs, it might be more practically useful to think henselian schemes themselves as the scheme-theoretic counterpart of henselian pairs.

To consider the scheme-theoretic counterpart of Zariskian pairs, we need the notion of {\em Zarikian schemes}: For a scheme-theoretic pair $(X,Y)$ we denote by $Z=\zat{X}|_Y$ the locally ringed space supported on the underlying topological space of $Y$ and the structure sheaf $\O_Z=i^{\ast}\O_X$, where $i\colon Y\rightarrow X$ is the closed immersion.
The sheaf $\O_Z$ comes with the topology induced from the $\mathscr{I}$-adic topology on $\O_X$, where $\mathscr{I}$ is the defining ideal of $Y$; moreover, $\mathscr{I}$ gives rise to a quasi-coherent sheaf of ideals of $\O_Z$, denoted again by $\mathscr{I}$.
Note that we have $1+\mathscr{I}\subseteq\O^{\times}_Z$.
A quasi-coherent ideal $\mathscr{J}\subseteq\O_Z$ is said to be an {\em ideal of definition}\index{ideal of definition} if, locally, there exists positive integers $n,m\in\Z$ such that $\mathscr{I}^n\subseteq\mathscr{J}^m\subseteq\mathscr{I}$.
\begin{dfn}\label{dfn-zariskianscheme1s}{\rm 
A {\em Zariskian scheme}\index{Zariskian!Zariskian scheme@--- scheme}\index{scheme!Zariskian scheme@Zariskian ---} is a topologically locally ringed space $X=(X,\O_X)$ that is locally isomorphic (as a topologically locally ringed space) to a Zariskian scheme associated to a pair.}
\end{dfn}

In {\bf \ref{ch-formal}}, \S\ref{sec-zariskianschemes} we give more generalities on Zariskian schemes.

\begin{dfn}\label{dfn-schemepair}{\rm 
A {\em Zariskian pair}\index{Zariskian!Zariskian pair@--- pair}\index{pair!Zariskian pair@Zariskian ---} is a pair $(X,Y)$ consisting of a Zariskian scheme $X$ and a closed subscheme $Y$ defined by an ideal of definition.}
\end{dfn}

For any pair $(X,Y)$ we have the associated Zariskian\index{Zariskian!associated Zariskian@associated ---}
$$
\zat{X}|_Y=(\zat{X}|_Y,Y)
$$
and the henselization\index{henselization!henselization of a pair@--- (of a pair)}
$$
\het{X}|_Y=(\het{X}|_Y,Y).
$$

By \ref{prop-adhesive1} (1) and \ref{prop-uadhesive1} (1) one can obviously think about (pseudo-) adhesiveness and universally (pseudo-) adhesiveness for pairs of schemes:
\begin{dfn}\label{dfn-schpairadh}{\rm 
A scheme-theoretic pair $(X,Y)$ is said to be {\em adhesive}\index{pair!adhesive pair@adhesive ---}\index{adhesive!adhesive pair@--- pair} (resp.\ {\em pseudo-adhesive}, resp.\ {\em universally adhesive}\index{pair!adhesive pair@adhesive ---!universally adhesive pair@universally --- ---}\index{adhesive!universally adhesive@universally ---}, resp.\ {\em universally pseudo-adhesive}) if there exists an affine open covering\index{affine!affine open covering@--- open covering} $\{U_i=\Spec A_i\}_{i\in I}$ of $X$ such that for each $i$ the pair $(A_i,I_i)$ is adhesive (resp.\ pseudo-adhesive, resp.\ universally adhesive, resp.\ universally pseudo-adhesive), where $I_i$ is the ideal defining the closed subscheme $Y\cap U_i\hookrightarrow U_i$.}
\end{dfn}

The following proposition follows immediately from \ref{prop-uadhesive1}:
\begin{prop}\label{prop-schpairadh2}
Let $(X,Y)$ be a universally adhesive $($resp.\ universally pseudo-adhesive$)$ pair of schemes, and $X'\rightarrow X$ an $X$-scheme locally of finite type.
Then $(X',X'\times_XY)$ is universally adhesive $($resp.\ universally pseudo-adhesive$)$. \hfill$\square$
\end{prop}

As the notion of (universally-) adhesiveness is local with respect to \'etale topology\index{topology!etale topology@\'etale ---} (\ref{prop-adhesive0} (2), \ref{prop-adhesive1} (2)), (universally-, pseudo-) adhesiveness for pairs $(X,Y)$ of algebraic spaces (where $Y$ is a closed subspace of $X$) can be defined in the obvious way (see \S\ref{subsub-algebraicspacesconv} to confirm what we mean precisely by algebraic spaces); we leave the details to the reader.

By \ref{prop-finiteidealtype2}, \ref{prop-adhesive0}, \ref{prop-adhesive1} (1), and \ref{prop-uadhesive1} (1) one sees easily the following result (and its obvious analogue for pairs of algebraic spaces):
\begin{prop}\label{prop-schpairadh}
A pair $(X,Y)$ of schemes is adhesive $($resp.\ pseudo-adhesive, resp.\ universally adhesive, resp.\ universally pseudo-adhesive$)$ if and only if for any \'etale neighborhood $U=\Spec A$ of $X$ the induced pair $(A,I)$ is adhesive $($resp.\ pseudo-adhesive, resp.\ universally adhesive, resp.\ universally pseudo-adhesive$)$, where $I$ is the ideal of $A$ corresponding to the closed subscheme $Y\times_XU$ of $U$. \hfill$\square$
\end{prop}

Finally, let us include a useful fact on coherency of structure sheaves, which follows immediate from \ref{thm-pf02a} (1) and \ref{prop-cohschemes2}:
\begin{prop}\label{prop-thm-cor-cohschemes21}
Let $(X,Y)$ be a universally adhesive pair of algebraic spaces such that $\O_X$ is $\mathscr{I}$-torsion free, where $\mathscr{I}=\mathscr{I}_Y$ is the defining ideal of the closed subspace $Y$.
Then $X$ is universally cohesive {\rm (\ref{dfn-universallycohesive})}\index{cohesive!universally cohesive@universally --- (schemes)}. \hfill$\square$
\end{prop}

\subsection{$I$-valuative rings}\label{sub-valuative}
\subsubsection{$I$-valuative rings}\label{subsub-Ivalnotion}
\index{valuative!Ivaluative ring@$I$-{---} ring|(}
Let $(A,I)$ be a pair of finite ideal type. 
Recall that an ideal $J$ of $A$ is said to be $I$-admissible\index{admissible!I-admissible@$I$-{---}} if it is finitely generated and contains a power of $I$ (cf.\ \ref{dfn-adm}).

\begin{dfn}\label{dfn-valuative1}{\rm 
Let $(A,I)$ be a pair of finite ideal type.
The ring $A$ is said to be {\em $I$-valuative}\index{valuative!Ivaluative ring@$I$-{---} ring} if any $I$-admissible ideal is invertible\index{ideal!invertible ideal@invertible ---} $($cf.\ \S\ref{subsub-invertible}).}
\end{dfn}

If one replaces $I$ by a finitely generated ideal of definition (without loss of generality), the above definition requires that also the ideal $I$ itself is invertible. 
If $I$ can be taken to be a principal ideal $I=(a)$, we will often say that $A$ is {\em $a$-valuative}.

The notion of $I$-valuative rings becomes particularly simple in the local ring case: a local ring $A$ is $I$-valuative if and only if $I$ is a principal ideal $I=(a)$ generated by a non-zero-divisor $a\in A$ and every $I$-admissible ideal is principal. 
For example, valuation rings are $I$-valuative for any non-zero finitely generated ideal $I$.
The following propositions are useful to reduce many situations to the local ring case:
\begin{prop}\label{prop-valuative11}
Let $(A,I)$ be a pair of finite ideal type.
If $A$ is $I$-valuative, then for any multiplicative subset $S\subseteq A$ the localization $B=S^{-1}A$ is $IB$-valuative.
\end{prop}

\begin{proof}
Let $J$ be a finitely generated ideal of $B$.
Then one can take a finitely generated ideal $J'$ of $A$ such that $J'B=J$.
Suppose $I^nB\subseteq J$ for some $n>0$.
Then one can replace $J'$ by $J'+I^n$, and hence we may assume that $J'$ is $I$-admissible.
Since $J'$ is invertible, $J=J'B$ is invertible.
\end{proof}

\begin{prop}\label{prop-valuative12}
Let $(A,I)$ be a pair of finite ideal type.
The following conditions are equivalent$:$
\begin{itemize}
\item[{\rm (a)}] $A$ is $I$-valuative$;$
\item[{\rm (b)}] $A_{\mathfrak{p}}$ is $IA_{\mathfrak{p}}$-valuative for any prime ideal $\mathfrak{p}$ of $A;$
\item[{\rm (c)}] $A_{\m}$ is $IA_{\m}$-valuative for any maximal ideal $\m$ of $A$.
\end{itemize}
\end{prop}

\begin{proof}
The implication (a) $\Rightarrow$ (b) follows from \ref{prop-valuative11}, and (b) $\Rightarrow$ (c) is trivial.
(c) $\Rightarrow$ (a) follows easily from \cite[Chap.\ II, \S 5.6, Theorem 4]{Bourb1}.
\end{proof}

By \ref{prop-valuative12} $I$-valuativeness has an obvious translation into the language of schemes:
\begin{dfn}\label{dfn-valuative2}
{\rm Let $(S,T)$ be a scheme-theoretic pair (\S\ref{sub-schpair}) such that the defining ideal $\mathscr{I}$ of $T$ is of finite type. We say that $S$ is {\em $\mathscr{I}$-valuative} (synonymously, $S$ is {\em $T$-valuative} or $(S,T)$ is {\em valuative}) if for any $x\in S$ the local ring $\O_{S,x}$ is $\mathscr{I}_x$-valuative.}
\end{dfn}

\begin{prop}\label{prop-valuative3}
Let $S$ be a scheme, and $\mathscr{I}$ a quasi-coherent ideal of finite type. 
Suppose that $S$ is $\mathscr{I}$-valuative, and set $U=S\setminus V(\mathscr{I})$. 

{\rm (1)} The scheme $S$ is integral if and only if $U$ is integral.

{\rm (2)} The scheme $S$ is integrally closed if and only if $U$ is integrally closed.
\end{prop}

\begin{proof}
(1) The `only if' part is clear. Suppose $U$ is integral. Let $j\colon U\hookrightarrow S$ be the open immersion.
Since $\mathscr{I}$ is invertible, $\O_S\rightarrow j_{\ast}\O_U$ is injective.
Let $\xi\in U$ be the generic point, and $i\colon\{\xi\}\hookrightarrow\O_S$ the inclusion.
Since $U$ is integral, $\O_U\rightarrow i_{\ast}\O_{U,\xi}$ is injective.
It follows that $\O_S\rightarrow i_{\ast}\O_{S,\xi}$ is injective, which shows that $S$ is integral.

(2) Only the `if' part calls for the proof.
We may assume that $S$ is affine $S=\Spec A$ and that $\mathscr{I}$ comes from a principal ideal $I=(a)$ with $a$ being a non-zero-divisor.
We may further assume that $A$ is a local ring.
Suppose $A[\frac{1}{a}]$ is integrally closed (in its total ring of fractions). 
Let $f(x)=x^n+\alpha_1x^{n-1}+\cdots+\alpha_{n-1}x+\alpha_n$ be a monic polynomial in $A[x]$, and suppose an element $x=c/b$ where $b$ is a non-zero-divisor satisfies $f(x)=0$.
Since $A[\frac{1}{a}]$ is integrally closed, one can take $c$ and $b$ such that $b$ is invertible in $A[\frac{1}{a}]$ and hence that $(b)$ is $I$-admissible. 
Since $(b,c)$ is also $I$-admissible and since $A$ is a local $I$-valuative ring, $(b,c)$ is a principal ideal $(d)$ generated by a non-zero-divisor $d$.
Replacing $b$ and $c$ respectively by $b/d$ and $c/d$, we may assume that $(b,c)=(1)$.
To show that $b$ is invertible, we suppose $b\in\m_A$. 
Since $(b,c)=(1)$, $c$ must be a unit. 
But then by the equation $b^nf(c/b)=0$, $c^n$ is divisible by $b$, which is absurd.
\end{proof}

\begin{prop}\label{prop-valuative5}
Let $S$ be an $\mathscr{I}$-valuative scheme, and $T$ an integral subscheme of $S$. 
Suppose that $T$ is not contained in $V(\mathscr{I})$. 
Then the scheme $T$ is $\mathscr{I}\O_T$-valuative.
\end{prop}

\begin{proof}
By \ref{prop-valuative12} we reduce to the case $S=\Spec A$ with $A$ local, $T=\Spec A/\mathfrak{p}$ where $\mathfrak{p}\subseteq A$ is a prime ideal, and $I=(a)$ with $a\not\in\mathfrak{p}$.
Then the assertion is straightforward.
\end{proof}

\begin{prop}\label{prop-valuative4}
Let $A$ be an $a$-valuative local ring with $a\in\m_A$, and suppose that $A[\frac{1}{a}]$ is a field. 
Then $A$ is an $a$-adically separated valuation ring. 
Conversely, if $A$ is a valuation ring, separated with respect to the adic topology defined by a finitely generated ideal $I$, then $\varinjlim_{n\geq 1}\Hom(I^n,A)$ is a field.
\end{prop}

Note that $\varinjlim_{n\geq 1}\Hom(I^n,A)\cong A[\frac{1}{a}]$ by a choice of a generator $a\in I$.
Hence the second part of the proposition is already shown in \ref{prop-sep}.

\begin{proof}
By \ref{prop-valuative3} (1) $A$ is an integral domain, and $A[\frac{1}{a}]$ is the fractional field of $A$.
Let $b/a^n\in A[\frac{1}{a}]$.
Since $(b,a^n)$ is invertible, we set $(b,a^n)=(d)$ and set $b=b'd$ and $a^n=c'd$. 
If $c'$ is invertible in $A$, then $b/a^n=b'/c'\in A$. 
If not, since $(b',c')=(1)$, $b'$ is a unit and hence $(b/a^n)^{-1}\in A$.
We thus verify the condition {\rm (b)} in \ref{dfn-val} and have shown that $A$ is a valuation ring.
To show it is $a$-adically separated, suppose $J=\bigcap_{n\geq 1}(a^n)\neq 0$, and take $f\in J\setminus\{0\}$. 
Since $f$ is invertible in $A[\frac{1}{a}]$, we can write $f^{-1}=c/a^n$ for $c\in A$, that is, $a^n$ is divisible by $f$.
But since $f$ is divisible by all $a^m$ ($m\geq 1$) at the same time, this is absurd.
\end{proof}

\subsubsection{Structure theorem}\label{subsub-valuativecomp}
The most important feature of $I$-valuative local rings is that they are `composite' of local rings and valuation rings\index{valuation!valuation ring@--- ring}:
\begin{thm}\label{thm-valuative}
{\rm (1)} Let $A$ be an $I$-valuative local ring, where $I\subseteq A$ is a non-zero proper finitely generated ideal. 
Set $J=\bigcap_{n\geq 1}I^n$.
Then$:$
\begin{itemize}
\item[{\rm (a)}] $B={\textstyle \varinjlim_{n\geq 1}}\Hom(I^n,A)$ is a local ring, and $V=A/J$ is an $\ovl{a}$-adically separated\index{valuation!valuation ring@--- ring!a-adically separated valuation ring@$a$-adically separated --- ---} valuation ring $($where $IV=(\ovl{a}))$ for the residue field $K$ of $B;$ 
\item[{\rm (b)}] $A=\{f\in B\,|\, (f\ \mathrm{mod}\ \m_B)\in V\};$
\item[{\rm (c)}] $J=\m_B$.
\end{itemize}

{\rm (2)} Conversely, for a local ring $B$ and an $\ovl{a}$-adically separated valuation ring $V$ for the residue field $K$ of $B$ with $\ovl{a}\neq 0$, the subring $A$ of $B$ defined by the equality in {\rm (b)} above is an $I$-valuative local ring for any finitely generated ideal $I$ such that $IV=(\ovl{a});$ we have, moreover, $B={\textstyle \varinjlim_{n\geq 1}}\Hom(I^n,A)$.
\end{thm}

Note that in (1), since $I\subseteq\m_A$, we have $I\neq I^2$ and hence $\ovl{a}\neq 0$.
Needless to say, composition/decomposition of valuation rings is the basic example of the theorem; recall that a valuation ring $V$ is $I$-valuative for any finitely generated ideal $I\subseteq V$ and that $J=\bigcap_{n\geq 1}I^n$ is a prime ideal (\ref{prop-associatedsepval}); the theorem in this particular situation is equivalent to \ref{prop-composition1} (1) and \ref{prop-composition2}.

\begin{proof}
(1) Let $I=(a)$; we have $B=A[\frac{1}{a}]$. 
Let $x$ be a closed point of $\Spec B$, and $S$ the Zariski closure of $\{x\}$ in $\Spec A$. 
If $\m$ denotes the maximal ideal of $B$ corresponding to $x$, we have $S=\Spec A/\mathfrak{p}$ with $\mathfrak{p}=A\cap\m$. 
By \ref{prop-valuative5} the ring $A/\mathfrak{p}$ is $I(A/\mathfrak{p})$-valuative.
Since $(A/\mathfrak{p})[\frac{1}{a}]=B/\m$ is a field, $A/\mathfrak{p}$ is an $\ovl{a}$-adically separated valuation ring due to \ref{prop-valuative4}, where $\ovl{a}$ is the image of $a$ in $A/\mathfrak{p}$.

\medskip
{\sc Claim 1.} {\it $V=A/\mathfrak{p}$ or, equivalently, $\mathfrak{p}=J$.}

\medskip
Indeed, since $A/\mathfrak{p}$ is $\ovl{a}$-adically separated, we have $(J$ mod $\mathfrak{p})=0$ and thus $J\subseteq\mathfrak{p}$.
For $f\in\mathfrak{p}$ consider the admissible ideal $(a^n,f)$ of $A$ for each $n\geq 1$, and set $(d)=(a^n,f)$ with $d\in A$. 
Since $f$ vanishes at $x$, $d$ must not be a unit. 
We need to show $f\in (a^n)=I^n$, and for this it is enough to show $d\in (a^n)$.
Since $a^n\in (d)$, there exists $b\in A$ such that $a^n=bd$.
On the other hand, since $f\in\mathfrak{p}$, $(b$ mod $\mathfrak{p})$ is a unit in $A/\mathfrak{p}$.
Since $A$ is a local ring, this implies that $b$ is a unit of $A$ and hence that $d\in (a^n)$, as desired.

\medskip
{\sc Claim 2.} {\it $B$ is a local ring.}

\medskip
We want show that any $f\in B\setminus\m$ is invertible in $B$.
Since $a$ is invertible in $B$, we may assume $f\in A$.
Let $\ovl{f}$ be the image of $f$ by the canonical map $A\rightarrow A/\mathfrak{p}=V$.
Since $V$ is an $\ovl{a}$-adically separated valuation ring, we have $\ovl{a}^n\in(\ovl{f})$ for some $n\geq 1$.
Let $H=(a^n,f)$, which is an admissible ideal of $A$ such that $HV=(\ovl{f})$.
There exists $h\in A$ such that $H=(h)$, and hence $f=gh$ for some $g\in A$.
Since the image of $g$ in $V$ is a unit, $g$ is a unit of $A$, and hence $a^n\in (f)$.
But this means that $fB=B$ and hence that $f$ is a unit of $B$, as desired.

Now since $V=A/(A\cap\m_B)$ ($\m_B=\m$), one has $V[\frac{1}{\ovl{a}}]=K$, where $K=B/\m_B$; that is, the residue field of $B$ is the fractional field of $V$.
Therefore, (a) is proved.

Next, let us show (b).
The inclusion $A\subseteq\{f\in B\,|\,(f\ \mathrm{mod}\ \m_B)\in V\}$ is clear.
Let $f$ be an element in the right-hand side, and set $f=g/a^n$ with $g\in A$.
There exists $h\in A$ such that $(h$ mod $J)$=$(f$ mod $\m_B)$, that is, $h-g/a^n\in\m_B$, and hence that $a^nh-g\in J$.
In particular, one can find $s\in A$ such that $a^nh-g=a^ns$, which gives $f=h-s\in A$, as desired.
By (b), in particular, we know that $\m_B\subseteq A$.
Hence we have $J=\mathfrak{p}=A\cap\m_B=\m_B$, whence (c).

(2) Let $A=\{f\in B\,|\,(f\ \mathrm{mod}\ \m_B)\in V\}$ and take finitely generated $I$ such that $IV=(\ovl{a})$.
Set $J=\m_B\cap A$. We have $A/J=V$.
First we claim that $A$ is a local ring. 
Set $\mathfrak{q}=\{x\in A\,|\,(x\ \mathrm{mod}\ J)\in\m_V\}$, and take $x\in A\setminus\mathfrak{q}$.
Then $\ovl{x}=x$ mod $J$ is a unit in $V$, which implies $x\in B^{\times}$. 
But since $(x^{-1}$ mod $\m_B)$ belongs to $V$, we deduce $x^{-1}\in A$.
Hence $A$ is a local ring, and $\m_A=\mathfrak{q}$.

Take $a\in I$ such that $(a$ mod $J)=\ovl{a}$.
We claim that $A[\frac{1}{a}]=B$.
The inclusion $A[\frac{1}{a}]\subseteq B$ is clear.
Take any $f\in B$, and set $\ovl{f}=(f$ mod $\m_B)$.
By \ref{prop-sep} there exists $n\geq 0$ such that $\ovl{f}\ovl{a}^n\in V$ and hence that $fa^n\in A$.
Hence $A[\frac{1}{a}]\supseteq B$, as desired.

Next we claim that $\bigcap_{n\geq 1}(a^n)=J$.
We have $\bigcap_{n\geq 1}(a^n)\subseteq J$, as $V$ is $\ovl{a}$-adically separated.
Suppose $x\in A\setminus\bigcap_{n\geq 1}(a^n)$, that is, $x\not\in(a^n)$ for some $n\geq 1$.
Since $x/a^n\not\in A$, we have $\ovl{a}^n/\ovl{x}\in\m_V$.
There exists $y\in\m_A$ such that $\ovl{y}=\ovl{a}^n/\ovl{x}$, that is, $\ovl{x}\ovl{y}=\ovl{a}^n$.
Since $\ovl{a}\neq 0$, $x\not\in J$. 
Hence $\bigcap_{n\geq 1}(a^n)\supseteq J$, thereby the claim.
In particular, since $(a)\subseteq I$ and $V$ is $\ovl{a}$-adically separated (which implies $\bigcap_{n\geq 1}I^n\subseteq J$), we have $\bigcap_{n\geq 1}I^n=J$.

Now, since any $I$-admissible ideals of $A$ contains $J$, they are in inclusion preserving one to one correspondence with finitely generated (hence principal) ideals of $V$. 
Let $H\supseteq I^n$ be an $I$-admissible ideal, and $(\ovl{b})$ the corresponding ideal of $V$, where $b\in H$.
To show $H=(b)$, it suffices to show that $(b)$ contains $J=\bigcap_{n\geq 1}(a^n)$.
Since $(\ovl{b})\supseteq (\ovl{a}^n)$, there exists $c\in A$ such that $a^n-bc\in J$.
Hence there exists $e\in 1+J$ such that $bc=a^ne$.
As $e$ is a unit in $A$, we deduce that $(a^n)\subseteq (b)$ and hence that $J\subseteq (b)$, as desired.
Hence any $I$-admissible ideals of $A$ are principal and are invertible, as $a$ is a non-zero-divisor.
(In particular, we have $I=(a)$.)
Therefore, $A$ is $I$-valuative, and now, all the assertions are proved.
\end{proof}

\subsubsection{Patching method}\label{subsub-patching}
\index{patching method|(}
Let $A$ be an $I$-valuative local ring with $I=(a)\subseteq\m_A$, and $B$ and $V$ the corresponding local ring and the valuation ring, respectively, determined as in \ref{thm-valuative}. 
The ring $V$ is a valuation ring for the residue field $K=B/\m_B$ of $B$ and is $\ovl{a}$-adically separated where $\ovl{a}=(a$ mod $J)$.
In particular, we have $K=V[\frac{1}{\ovl{a}}]$ (\ref{prop-sep}).

Let $X$ be an object over $A$, e.g.\ scheme, algebra, module, etc. 
Then, by base change, it induces the objects $X_B$ over $B$ and $X_V$ over $V$ with an isomorphism $\phi\colon X_B\otimes K\stackrel{\sim}{\rightarrow}X_V\otimes K$.
In many occasions it is essential to regard $X$ as obtained by `patching' of $X_B$ and $X_V$ along $\phi$. 
While it will turn out in many situations that the functor $X\mapsto (X_B,X_V,\phi)$ is essentially surjective, that is to say, any triple $(X_B,X_V,\phi)$ can be anyway patched together to an object over $A$, we like to have more precise picture of the patching; for example, we would like to ask for an oppositely-oriented functor, so to speak, the `patching functor', from the category of the triples as above to the category of $A$-objects.

Let us formulate the situation more precisely.
For a ring $R$ we denote by $\Mod_R$ the category of $R$-modules. 
By base change, we have the commutative (more precisely, 2-commutative) diagram of categories
$$
\xymatrix{\Mod_A\ar[r]\ar[d]&\Mod_B\ar[d]\\ \Mod_V\ar[r]&\Mod_K\rlap{,}}
$$
which gives rise to a functor
$$
\beta\colon\Mod_A\longrightarrow\Mod_B\times_{\Mod_K}\Mod_V,
$$
where the right-hand side denotes the $2$-fiber product of the categories, that is, the category of triples $(L,M,\iota)$ consisting of
\begin{itemize}
\item a $B$-module $L$,
\item a $V$-module $M$,
\item an isomorphism $\iota\colon L\otimes_BK\stackrel{\sim}{\rightarrow} M\otimes_VK$ of $K$-modules,
\end{itemize}
and a morphism $(L_1,M_1,\iota_1)\longrightarrow (L_2,M_2,\iota_2)$ is a pair $(f,g)$ consisting of a $B$-morphism $f\colon L_1\rightarrow L_2$ and a $V$-morphism $g\colon M_1\rightarrow M_2$ such that the following diagram commutes:
$$
\xymatrix{L_1\otimes_BK\ar[d]_{\iota_1}\ar[r]^{f\otimes_B\mathrm{id}_K}&L_2\otimes_BK\ar[d]^{\iota_2}\\ M_1\otimes_VK\ar[r]_{g\otimes_V\mathrm{id}_K}&M_2\otimes_VK\rlap{.}}
$$
We say that $(f,g)$ is injective (resp.\ surjective) if both $f$ and $g$ are injective (resp.\ surjective).
The functor $\beta$ is then given by 
$$
\beta(N)=(N\otimes_AB,N\otimes_AV,\mathrm{can})
$$
for any $A$-module $N$, where $\mathrm{can}$ is the canonical isomorphism.

Now we define the {\em patching functor}\index{functor!patching functor@patching ---} 
$$
\alpha\colon\Mod_B\times_{\Mod_K}\Mod_V\longrightarrow\Mod_A
$$
by
$$
\alpha(L,M,\iota)=\{(x,y)\in L\times M\,|\,\iota(x\ \mathrm{mod}\ \m_BL)=y\otimes 1\ \mathrm{in}\ M\otimes_VK\}.
$$
Thus we have two functors between these categories
$$
\xymatrix{\Mod_A\ar@<.5ex>[r]^(.33){\beta}&\Mod_B\times_{\Mod_K}\Mod_V\rlap{.}\ar@<.5ex>[l]^(.67){\alpha}}
$$

These functors can be similarly defined between the corresponding categories of algebras:
$$
\xymatrix{\Alg_A\ar@<.5ex>[r]^(.33){\beta}&\Alg_B\times_{\Alg_K}\Alg_V\rlap{.}\ar@<.5ex>[l]^(.67){\alpha}}
$$

\begin{thm}\label{thm-patching}
{\rm (1)} The functor $\beta$ is the left-adjoint to the functor $\alpha$, and the adjunction morphism $\beta\circ\alpha\rightarrow\mathrm{id}$ is a natural equivalence. 
In particular, $\beta$ is essentially surjective, and $\alpha$ is fully faithful.

{\rm (2)} The essential image of $\alpha$ consists of $A$-modules $N$ such that $JN$ is $a$-torsion free.

Moreover, the similar assertions hold for algebras.
\end{thm}

\begin{proof}
Clearly, it suffices to show the theorem in the case of modules.

(1) First we are to show that the adjunction map $\beta\circ\alpha(L,M,\iota)\rightarrow (L,M,\iota)$ is an isomorphism.
This map is defined by the pair of morphisms $N\otimes_AB\rightarrow L$ and $N\otimes_AV\rightarrow M$, where $N=\alpha(L,M,\iota)\subseteq L\times M$, induced respectively by the first and second projections.
It is easy to see that these maps are injective and that $N\otimes_AV\rightarrow M$ is surjective.
To show that $N\otimes_AB\rightarrow L$ is surjective, take $x\in L$ and set $\ovl{x}=x$ mod $\m_BL$. 
Since $K=V[\frac{1}{\ovl{a}}]$ (\ref{prop-sep}), we can find $n\geq 0$ such that $\iota(\ovl{a^nx})$ is of the form $y\otimes 1$. 
Then one has the element $(a^nx,y)\otimes a^{-n}$ in $N\otimes_AB$, which maps to $x$, as desired.
Now, by means of the natural transformation $\beta\circ\alpha\rightarrow\mathrm{id}$ thus obtained, we have the canonical map
$$
\Hom_{\Mod_A}(N,\alpha(L,M,\iota))\longrightarrow\Hom_{\Mod_B\times_{\Mod_K}\Mod_V}(\beta(N),(L,M,\iota)),
$$
which is bijective due to the presence of the other adjunction map $N\rightarrow\alpha\circ\beta(N)$ defined in an obvious way.

(2) Let us first show that the adjunction map $N\rightarrow\alpha\circ\beta(N)$ is surjective.
Every element of $\alpha\circ\beta(N)$ is of the form $(x\otimes a^{-n},y\otimes 1)\in N_B\times N_V$ such that $\ovl{x}=\ovl{a}^n\ovl{y}$, where $\ovl{\,\cdot\,}$ denotes the mod-$JN$ class.
Since $x-a^ny\in JN$, there exists $z\in N$ such that $x=a^nz$ and that $\ovl{z}=\ovl{y}$.
Hence we have $(x\otimes a^{-n},y\otimes 1)=(z\otimes 1, z\otimes 1)$, which is the image of $z$ by the adjunction map, thereby the claim.
Note that the kernel of this map is $JN\cap N_{\ator}=(JN)_{\ator}$.
The assertion follows from these facts, combined with that the other adjunction map $\beta\circ\alpha(L,M,\iota)\rightarrow (L,M,\iota)$ is an isomorphism.
\end{proof}

\begin{prop}\label{prop-patching1}
Let $N$ be an $A$-module. Then $N$ lies in the essential image of the functor $\alpha$ if $\Tor^B _1(N\otimes_AB, K)=0$ or, sufficiently, if $N$ is flat outside $I$.
\end{prop}

\begin{proof}
The exact sequence
$$
0\longrightarrow J\longrightarrow A\longrightarrow V\longrightarrow 0
$$
gives rise to the exact sequence
$$
0\longrightarrow\Tor^A_1(N,V)\longrightarrow J\otimes_AN\longrightarrow N\longrightarrow N/JN\longrightarrow 0.
$$
Since $\Tor^A_1(N,V)\otimes K\cong\Tor^A_1(N,K)\cong\Tor^B_1(N\otimes_AB,K)=0$, $\Tor^A_1(N,V)$ is an $a$-torsion module.
But since $J\otimes_AN$ is $a$-torsion free (as $x\mapsto ax$ is bijective on $J$), we have $\Tor^A_1(N,V)=0$.
Then it follows that $J\otimes_AN\cong JN$ and that $JN$ is $a$-torsion free.
\end{proof}

\begin{prop}\label{prop-patching2}
{\rm (1)} Let $(L,M,\iota)$ be an object of $\Mod_B\times_{\Mod_K}\Mod_V$, and set $N=\alpha(L,M,\iota)$.
Suppose $L$ is flat over $B$. 
Then the following conditions are equivalent$:$
\begin{itemize}
\item[{\rm (a)}] $N$ is a finitely generated $A$-module$;$
\item[{\rm (b)}] $L$ is a finitely generated $B$-module, and $M$ is a finitely generated $V$-module.
\end{itemize}

{\rm (2)} Let $(P,Q,\iota)$ be an object of $\Alg_B\times_{\Alg_K}\Alg_V$, and set $R=\alpha(P,Q,\iota)$.
Suppose $P$ is flat over $B$. 
Then the following conditions are equivalent$:$
\begin{itemize}
\item[{\rm (a)}] $R$ is a finitely generated $($resp.\ finitely presented$)$ $A$-algebra$;$
\item[{\rm (b)}] $P$ is a finitely generated $($resp.\ finitely presented$)$ $B$-algebra, and $Q$ is a finitely generated $($resp.\ finitely presented$)$ $V$-algebra.
\end{itemize}
\end{prop}

\begin{proof}
Both in (1) and (2) the direction (a) $\Rightarrow$ (b) is easy. 
We want to show the converse. 
The reasonings for proving (b) $\Rightarrow$ (a) in (1) and (2) (for `finitely generated') are quite similar, and we confine ourselves to present it only in (2). 

Suppose $P$ (resp.\ $Q$) is a finitely generated algebra over $B$ (resp.\ $V$), and take the generators $x_1,\ldots,x_s$ (resp.\ $\ovl{y}_1,\ldots,\ovl{y}_t$).
Since $P=R\otimes_AB$, multiplying a power of $a$ if necessary, we may assume that all $x_i$ belong to $R$. 
Similarly, since $Q=R\otimes_AV=R/JR$, we can take $y_i\in R$ that lifts $\ovl{y}_i$ for $i=1,\ldots,t$.
Then we want to show that $\{x_1,\ldots,x_s,y_1,\ldots,y_t\}$ generates $R$.
To this end, we consider the morphism $f\colon C=A[X_1,\ldots,X_s,Y_1,\ldots,Y_t]\rightarrow R$ given by $f(X_i)=x_i$ and $f(Y_j)=y_j$.

\medskip
{\sc Claim 1.} {\it The induced map $f\colon JC\rightarrow JR$ is surjective.}

\medskip
It is clear that, if the claim is settled, we can prove that $R$ is finitely generated; indeed, $Q$ is finitely generated, and $Q=R/JR$.

To prove the claim, we first note that the flatness assumption implies that $\Tor^B_1(R\otimes_AB,K)=0$.
The same reasoning as in the proof of \ref{prop-patching1} shows that $\Tor^A_1(R,V)=0$ and $J\otimes_AR=JR$. 
Similarly, we have $J\otimes_AC=JC$.
But since $\mathrm{id}_J\otimes f\colon J\otimes_AC\rightarrow J\otimes_AR$ is surjective, the claim follows.

Next, suppose $P$ (resp.\ $Q$) is a finitely presented algebra over $B$ (resp.\ $V$).
To show that $R$ is finitely presented, it suffices to show that for any surjection $A[X]=A[X_1,\ldots,X_s]\rightarrow R$ from a polynomial ring its kernel $H$ is finitely generated.
Since $H\otimes_AB$ is the kernel of the induced map $B[X_1,\ldots,X_s]\rightarrow P=R\otimes_AB$ and since $\Tor^B_i(P,K)=0$ for all $i\geq 1$, we have $\Tor^B_1(H\otimes_AB,K)=0$.
Hence $H$ lies in the essential image of $\alpha$ by \ref{prop-patching1}.
In particular, $H=\alpha\circ\beta(H)$.

\medskip
{\sc Claim 2.} {\it $\beta(H)\rightarrow\beta(A[X])$ is injective.}

\medskip
Indeed, since $\Tor^A_1(R,V)=0$ as above, we see that $H\otimes_AV\rightarrow A[X]\otimes_AV$ is injective. 
Since the other map $H\otimes_AB\rightarrow B[X]$ is clearly injective, the claim is settled.

The claim implies that $H\otimes_AV$ is the kernel of $V[X]\rightarrow Q$.
By the assumption $H\otimes_AB$ (resp.\ $H\otimes_AV$) is a finitely generated $B[X]$ (resp.\ $V[X]$)-module. 
Hence by (1) $H=\alpha\circ\beta(H)$ is a finitely generated $A$-module.
\end{proof}

\begin{prop}\label{prop-patching3}
Let $(P,Q,\iota)$ be an object of $\Alg_B\times_{\Alg_K}\Alg_V$, and set $R=\alpha(P,Q,\iota)$.
Then the following conditions are equivalent$:$
\begin{itemize}
\item[{\rm (a)}] $R$ is flat $($resp.\ smooth, resp.\ \'etale$)$ over $A;$
\item[{\rm (b)}] $P$ is flat $($resp.\ smooth, resp.\ \'etale$)$ over $B$, and $Q$ is flat $($resp.\ smooth, resp.\ \'etale$)$ over $V$.
\end{itemize}
\end{prop}

\begin{proof}
(a) $\Rightarrow$ (b) is clear.
We need to show the converse.
We have the diagram
$$
\xymatrix{\Spec Q\ar[d]\ar@{^{(}->}(7,0);(15,0)&\Spec R\ar[d]&\Spec P\ar[d]\ar@{_{(}->}(35,0);(27,0)\\ \Spec V\ar@{^{(}->}(7,-13.5);(15,-13.5)&\Spec A&\Spec B\ar@{_{(}->}(35,-13.5);(27,-13.5)\rlap{,}}
$$
where the squares are Cartesian (due to \ref{thm-patching} (1)) and the right-hand (resp.\ left-hand) inclusions are open (resp.\ closed) immersions.
Then in view of \ref{prop-patching2} (2) and \cite[$\mathbf{IV}$, (17.5.1), (17.6.1)]{EGA} it is enough to show that, if $P$ is $B$-flat and $Q$ is $V$-flat, then $R$ is $A$-flat.
Since $P=R[\frac{1}{a}]$ and $Q=R/JR$, it suffices to show that $R$ is $a$-torsion free; indeed, since $Q$ is $V$-flat, $R/aR$ is flat over $V/aV=A/aA$, and hence we can apply \ref{cor-glfl}.

Let $x\in R$ and suppose $a^nx=0$ for some $n>0$.
Then in $V$ we have $\ovl{a^nx}=0$; since $\ovl{a}\neq 0$ and $V$ is an integral domain, we have $x\in JR$.
But by \ref{thm-patching} (2) we know that $JR$ is $a$-torsion free and hence that $x=0$, as desired.
\end{proof}

\begin{prop}\label{prop-patching4}
Let $A$ be an $I$-valuative local ring by a non-zero $I=(a)\neq A$, and $B$, $J$, $V$, and $K$ as in {\rm \ref{thm-valuative}}.
Then the following conditions are equivalent$:$
\begin{itemize}
\item[{\rm (a)}] $A$ is $I$-adically henselian\index{henselian!I-adically henselian@$I$-adically ---}$;$
\item[{\rm (b)}] $V$ is $IV$-adically henselian, and $B$ is $\mathfrak{m}_B$-adically henselian. 
\end{itemize}
\end{prop}

\begin{proof}
First we prove (a) $\Rightarrow$ (b).
It is obvious that $V$ is henselian, for it is a quotient of $A$.
To see that $B$ is henselian, we verify the following: for any $B$-algebra $P$ \'etale over $B$ such that $P\otimes_BF=P/\mathfrak{m}_BP=K$, there exists a homomorphism $P\rightarrow B$ such that the composition $B\rightarrow P\rightarrow B$ is the identity map.

If $P$ is as above, then $(P,V,\iota)$ (for some $\iota$) is an object of $\Alg_B\times_{\Alg_K}\Alg_V$, and hence we have an $A$-algebra $R=\alpha(P,V,\iota)$.
By \ref{prop-patching3} $R$ is an \'etale $A$-algebra, and we moreover know that $R/IR=V/IV=A/I$.
Hence there exists a homomorphism $R\rightarrow A$ such that the composition $A\rightarrow R\rightarrow A$ is the identity map.
Since $R[\frac{1}{a}]=P$, we have the desired morphism $P\rightarrow B$ by thebase change.

Next, we show (b) $\Rightarrow$ (a).
Let $R$ be an \'etale $A$-algebra such that $R/IR=A/I$.
Then by the base change (that is, by the functor $\beta$) we have the \'etale $B$-algebra $P=R[\frac{1}{a}]$ and the \'etale $V$-algebra $Q=R/JR$.
Since $V$ is henselian with respect to $IV$-adic topology, the obvious equalities $Q/IQ=R/IR=V/IV$ implies a morphism $Q\rightarrow V$; similarly, as $J=\m_B$, we have a morphism $P/JP=Q[\frac{1}{\ovl{a}}]\rightarrow V[\frac{1}{\ovl{a}}]=K$ (cf.\ \ref{prop-sep}), which gives rise to a map $P\rightarrow B$.
Hence we have homomorphisms $P\rightarrow B$ and $Q\rightarrow V$ that give respectively the right-inverses of $B\rightarrow P$ and $V\rightarrow Q$.
These maps form a morphism $(P,Q,\mathrm{can})\rightarrow (B,V,\mathrm{can})$ in $\Alg_B\times_{\Alg_K}\Alg_V$, giving a right-inverse to the structural map $(B,V,\mathrm{can})\rightarrow (P,Q,\mathrm{can})$.
Then by \ref{prop-patching1} the functor $\alpha$ gives rise to a morphism $R\rightarrow A$ such that the composition $A\rightarrow R\rightarrow A$ is the identity map, whence the claim.
\end{proof}
\index{patching method|)}
\index{valuative!Ivaluative ring@$I$-{---} ring|)}

\subsection{Pairs and complexes}\label{sub-paircomplex}
\subsubsection{Set-up}\label{subsub-paircomplexgen}
Let us first fix notations that we will use throughout this subsection.
Let $(A,I)$ be a pair, and 
$$
K^{\bullet}\colon\cdots\longrightarrow K^{q-1}\stackrel{d^{q-1}}{\longrightarrow}K^q\stackrel{d^q}{\longrightarrow}K^{q+1}\stackrel{d^{q+1}}{\longrightarrow}\cdots
$$
a complex\index{complex} of $A$-modules.
We consider the $I$-adic filtration $\{I^nK^{\bullet}\}_{n\geq 0}$ on the complex $K^{\bullet}$ and define the complex $K^{\bullet}_k$ by the exact sequence
$$
0\longrightarrow I^{k+1}K^{\bullet}\longrightarrow K^{\bullet}\longrightarrow K^{\bullet}_k\longrightarrow 0\leqno{(\ast)}
$$
of complexes for each $k\geq 0$.
This induces the cohomology long exact sequence 
$$
\cdots\longrightarrow\H^q(I^{k+1}K^{\bullet})\stackrel{\iota^q_k}{\longrightarrow}\H^q(K^{\bullet})\stackrel{\sigma^q_k}{\longrightarrow}\H^q(K^{\bullet}_k)\stackrel{\delta^q_k}{\longrightarrow}\H^{q+1}(I^{k+1}K^{\bullet})\longrightarrow\cdots.
$$

\begin{dfn}\label{dfn-inducedfiltrationpair}{\rm 
For $q\in\Z$ we define the filtrations 
$$
\{F^n\H^q(I^{k+1}K^{\bullet})\}_{n\geq 0}\quad\textrm{and}\quad\{F^n\H^q(K^{\bullet}_k)\}_{n\geq 0}
$$
on the cohomologies $\H^q(I^{k+1}K^{\bullet})$ and $\H^q(K^{\bullet}_k)$, respectively, for each $k\geq 0$ by$:$
\begin{equation*}
\begin{split}
&F^n\H^q(I^kK^{\bullet})=\image(\H^q(I^nK^{\bullet}\cap I^kK^{\bullet})\rightarrow\H^q(I^kK^{\bullet}));\\
&F^n\H^q(K^{\bullet}_k)=\image(\H^q(K^{\bullet}/I^nK^{\bullet}\cap I^{k+1}K^{\bullet})\rightarrow\H^q(K^{\bullet}_k)).
\end{split}
\end{equation*}
We call these filtrations the {\em induced filtrations on the cohomologies}\index{filtration by submodules@filtration (by submodules)!induced filtration by submodules@induced ---!induced filtration by submodules on the cohomologies@--- --- on the cohomology}.}
\end{dfn}

\subsubsection{Results in case $I$ is finitely generated}\label{subsub-paircomplexgenresult}
\begin{lem}\label{lem-steinfactadqformalsch2exer}
Let $(A,I)$ be a complete pair of finite ideal type, and $q\in\Z$.
Suppose that the following conditions are satisfied$:$
\begin{itemize}
\item[{\rm (a)}] $K^q$ is $I$-adically separated\index{separated!I-adically separated@$I$-adically ---}, and $K^{q-1}$ is $I$-adically complete$;$
\item[{\rm (b)}] $\H^q(K^{\bullet})/I\H^q(K^{\bullet})$ is finitely generated as an $A$-module.
\end{itemize}
Then the $A$-module $\H^q(K^{\bullet})$ is finitely generated.
More precisely, if $\beta_1,\ldots,\beta_m$ are elements of $\H^q(K^{\bullet})$ that generate $\H^q(K^{\bullet})/I\H^q(K^{\bullet})$, then $\beta_1,\ldots,\beta_m$ generate $\H^q(K^{\bullet})$ as an $A$-module.
\end{lem}

\begin{proof}
We may assume that $I$ is finitely generated: $I=(a_1,\ldots,a_r)$.
Consider the projection $p\colon\ker(d^q)\rightarrow\H^q(K^{\bullet})$ and for each $j=1,\ldots,m$ an element $\alpha_j\in\ker(d^q)$ such that $p(\alpha_j)=\beta_j$.
Set $M=A^{\oplus m}$, and let $\Phi\colon M\rightarrow\ker(d^q)$ be the map given by the elements $\alpha_1,\ldots,\alpha_m$.
Set $\varphi=p\circ\Phi$.
We have $\H^q(K^{\bullet})=\varphi(M)+I\H^q(K^{\bullet})$.
Hence, in particular, 
$$
\ker(d^q)=\Phi(M)+\image(d^{q-1})+I\ker(d^q).
$$
For $x=x_0\in\ker(d^q)$ we have $y_1\in M$, $z_1\in K^{q-1}$, and $x_1^{(1)},\ldots,x_1^{(r)}\in\ker(d^q)$ such that
$$
x=\Phi(y_1)+d^{q-1}(z_1)+\sum_{i=1}^ra_ix_1^{(i)}.
$$
We do the same for each $x_1^{(i)}$ to get a similar equalities $x_1^{(i)}=\Phi(y^{(i)}_2)+d^{q-1}(z^{(i)}_2)+\sum^r_{j=1}a_jx^{(ij)}_2$.
Then we have
$$
x=\Phi(y_2)+d^{q-1}(z_2)+\sum^r_{i,j=1}a_ia_jx^{(ij)}_2, 
$$
where $y_2=y_1+\sum^r_{i=1}a_iy^{(i)}_2$ and $z_2=z_1+\sum^r_{i=1}a_iz^{(i)}_2$.

We repeat this to obtain the sequences $\{y_n\}_{n\geq 1}$, where $y_n$ is a polynomial of $a_1,\ldots,a_r$ of degree $n-1$ with the coefficients in $M$, and $\{z_n\}$, where $z_n$ is a polynomial of $a_1,\ldots,a_r$ of degree $n-1$ wi the coefficients in $K^{q-1}$ such that, for any $n\geq 1$, $x-(\Phi(y_n)+d^{q-1}(z_n))$ is a homogeneous polynomial of $a_1,\ldots,a_r$ of degree $n$ with the coefficients in $\ker(d^q)$, hence belonging to $I^n\ker(d^q)$.
Set $y=\lim y_n\in M$ and $z=\lim z_n\in K^{q-1}$.
Since $K^q$, and hence $\ker(d^q)$ also, is $I$-adically separated, we have
$$
x=\Phi(y)+d^{q-1}(z).
$$
Thus we have shown that $\ker(d^q)=\Phi(M)+\image(d^{q-1})$, that is, $\varphi\colon M\rightarrow\H^q(K^{\bullet})$ is surjective.
\end{proof}

\begin{rem}\label{rem-steinfactadqformalsch2exer}{\rm 
If we further assume in \ref{lem-steinfactadqformalsch2exer} that $I$ is principal $I=(a)$ and that $\H^q(K^{\bullet})_{\ator}$ is bounded $a$-torsion, then the consequence of \ref{lem-steinfactadqformalsch2exer} follows by a simpler argument: by \ref{prop-complpair1} it suffices to show that $\H^q(K^{\bullet})$ is $a$-adically separated, which is verified by an easy argument (cf.\ Exercise \ref{exer-aadictopcompleteseparatedclosed}).}
\end{rem}

\begin{lem}\label{lem-steinfactadqformalsch2}
Let $(A,I)$ be a complete pair of finite ideal type, and $q\in\Z$ an integer.
Suppose that the following conditions are satisfied$:$
\begin{itemize}
\item[{\rm (a)}] $K^q$ is $I$-adically separated, and $K^{q-1}$ is $I$-adically complete$;$
\item[{\rm (b)}] the topology on $\H^q(K^{\bullet})$ defined by $\{F^n\H^q(K^{\bullet})\}_{n\geq 0}$ is $I$-adic$;$
\item[{\rm (c)}] for any $k\geq 0$ the image of the map $\H^q(K^{\bullet})\rightarrow\H^q(K^{\bullet}_k)$ is finitely generated as an $A$-module.
\end{itemize}
Then $\H^q(K^{\bullet})$ is finitely generated as an $A$-module.
\end{lem}

\begin{proof}
By (c) the $A$-module $\H^q(K^{\bullet})/F^{k+1}\H^q(K^{\bullet})$ is finitely generated for any $k\geq -1$.
By (b) there exists an integer $s=s(k)\geq 0$ such that $F^{k+s+1}\H^q(K^{\bullet})\subseteq I^{k+1}\H^q(K^{\bullet})$ for any $k\geq 0$.
Hence $\H^q(K^{\bullet})/I\H^q(K^{\bullet})$ is a quotient of the finitely generated $A$-module $\H^q(K^{\bullet})/F^{s(0)+1}\H^q(K^{\bullet})$, and hence is finitely generated.
Now, applying \ref{lem-steinfactadqformalsch2exer}, we deduce the assertion.
\end{proof}

\begin{lem}\label{lem-steinfactadqformalsch1}
Suppose that the pair $(A,I)$ is pseudo-adhesive\index{adhesive!adhesive pair@--- pair!pseudo adhesive pair@pseudo-{---} ---} {\rm (\ref{dfn-adhesive})}.

{\rm (1)} We have the equality
$$
(\delta^q_k)^{-1}(F^n\H^{q+1}(I^{k+1}K^{\bullet}))=F^n\H^q(K^{\bullet}_k)
$$
for any $k,n\geq 0$ and $q\in\Z$.

{\rm (2)} For a fixed $q\in\Z$ suppose that $\H^{q+1}(I^{k+1}K^{\bullet})$ is finitely generated and that the topology on it defined by the induced filtration $\{F^n\H^{q+1}(I^{k+1}K^{\bullet})\}_{n\geq 0}$ is $I$-adic for any $k\geq 0$.
Then the projective system\index{system!projective system@projective ---} $\{\H^q(K^{\bullet}_k)\}_{k\geq 0}$ satisfies {\bf (ML)}\index{Mittag-Leffler condition}.
\end{lem}

\begin{proof}
To show (1), we may assume $n\geq k$.
Consider the commutative diagram 
$$
\xymatrix{\H^q(K^{\bullet})\ar@{=}[d]\ar[r]&\H^q(K^{\bullet}_n)\ar[d]_{\alpha}\ar[r]^(.4){\delta^q_n}&\H^{q+1}(I^{n+1}K^{\bullet})\ar[d]^{\beta}\ar[r]&\H^{q+1}(K^{\bullet})\ar@{=}[d]\\ \H^q(K^{\bullet})\ar[r]&\H^q(K^{\bullet}_k)\ar[r]_(.4){\delta^q_k}&\H^{q+1}(I^{k+1}K^{\bullet})\ar[r]&\H^{q+1}(K^{\bullet})}
$$
with the exact rows.
What to prove is the equality $(\delta^q_k)^{-1}(\image(\beta))=\image(\alpha)$, which follows from an easy diagram chasing.

To show (2), first notice that we have
$$
(\delta^q_k)^{-1}(F^n\H^{q+1}(I^{k+1}K^{\bullet})\cap\H^{q+1}(I^{k+1}K^{\bullet})_{\Itor})=F^n\H^q(K^{\bullet}_k)
$$
due to (1).
Since $(A,I)$ is pseudo-adhesive, $\H^{q+1}(I^{k+1}K^{\bullet})_{\Itor}$ is bounded $I$-torsion.
Since the filtration $\{F^n\H^{q+1}(I^{k+1}K^{\bullet})\}_{n\geq 0}$ defines the $I$-adic topology and since $(A,I)$ is pseudo-adhesive, we have in view of \ref{prop-AR}
$$
F^n\H^{q+1}(I^{k+1}K^{\bullet})\cap\H^{q+1}(I^{k+1}K^{\bullet})_{\Itor}=0
$$
for sufficiently large $n$.
This implies that the filtration $\{F^n\H^q(K^{\bullet}_k)\}_{n\geq 0}$ is stationary.
In other words, for any $k\geq -1$ there exists $l\geq k$ such that for any $m\geq l$ the maps $\H^q(K^{\bullet}_l)\rightarrow\H^q(K^{\bullet}_k)$ and $\H^q(K^{\bullet}_m)\rightarrow\H^q(K^{\bullet}_k)$ have the same image, which is nothing but the condition {\bf (ML)}.
\end{proof}

\subsubsection{Results in case $I$ is principal}\label{subsub-paircomplexprinresult}
When the ideal of definition $I$ is principal $I=(a)$, one can prove stronger results.

\begin{lem}\label{lem-ullrichlemma2}
Consider the situation as in {\rm \S\ref{subsub-paircomplexgen}}, and suppose that $I$ is principal $I=(a)$.
Let $q\in\Z$ be an integer.

{\rm (1)} Suppose $a^nK^{q+1}_{\ator}=0$ for some $n\geq 0$.
Then for any $k\geq n$ we have
$$
a^{k+1}\H^q(K^{\bullet})\subseteq F^{k+1}\H^q(K^{\bullet})\subseteq a^{k-n+1}\H^q(K^{\bullet}).
$$

{\rm (2)} Suppose, moreover, that $A$ and $K^{q-1}$ is $a$-adically complete and $K^q$ is $a$-adically separated.
Let $\beta_1,\ldots,\beta_m\in\H^q(K^{\bullet})$ be elements such that for some $k\geq n$ the elements $\sigma^q_k(\beta_1),\ldots,\sigma^q_k(\beta_m)$ generate $\sigma^q_k(\H^q(K^{\bullet}))$ as an $A$-module.
Then $\H^q(K^{\bullet})$ is generated by the elements $\beta_1,\ldots,\beta_m$ over $A$.
In particular, the $A$-module $\H^q(K^{\bullet})$ is finitely generated, and the topology on $\H^q(K^{\bullet})$ defined by the induced filtration $\{F^l\H^q(K^{\bullet})\}_{l\geq 0}$ is $I$-adic.
\end{lem}

\begin{proof}
Applying \ref{lem-ARprincipalgeneralized} to 
$$
0\longrightarrow\ker(d^q)\longrightarrow K^q\longrightarrow\image(d^q)\longrightarrow 0,
$$
one obtains the equality $\ker(d^q)\cap a^{k+1}K^q=a^{k-n+1}(\ker(d^q)\cap a^nK^q)$.
By this one has the inclusions $a^{k+1}\ker(d^q)\subseteq\ker(d^q)\cap a^{k+1}K^q\subseteq a^{k-n+1}\ker(d^q)$, whence (1).
This means that the topology on $\H^q(K^{\bullet})$ defined by the filtration $\{F^l\H^q(K^{\bullet})\}_{l\geq 0}$ coincides with the $a$-adic topology, which shows the last assertion of (2).
Then the first assertion of (2) follows from \ref{lem-steinfactadqformalsch2exer}.
\end{proof}

\begin{rem}\label{rem-steinfactadqformalsch2exer2}{\rm 
Similarly to Remark \ref{rem-steinfactadqformalsch2exer}, the proof of \ref{lem-ullrichlemma2} can be much simplified if we further assume, for example, that the $a$-torsion part $K^q_{\ator}$ is bounded $a$-torsion (which we can assume in most of the later applications), since in this situation one can show that the cohomology $\H^q(K^{\bullet})$ is $a$-adically separated.}
\end{rem}

\begin{lem}\label{lem-finiformalspecial2}
Consider the situation as in {\rm \S\ref{subsub-paircomplexgen}}, where the ideal $I$ is assumed to be principal $I=(a)$, and suppose, moreover, that $A$ is $a$-adically complete t.u.\ adhesive {\rm (\ref{dfn-topadhesive})}\index{adhesive!universally adhesive@universally ---!topologically universally adhesive@topologically --- --- (t.u.\ adhesive)} and $a$-torsion free $($hence, in particular, $A$ is topologically universally coherent with respect to $I$ {\rm (\ref{dfn-cohringsmodules5})}\index{coherent!topologically universally coherent@topologically universally ---} by $\ref{thm-pf02a})$.
Let $q\in\Z$, and suppose there exists $n\geq 0$ such that $a^nK^{q+1}_{\ator}=0$.

{\rm (1)} If $K^{q-1}$ is $a$-adically complete, $K^q$ is $a$-adically separated, $\H^q(K^{\bullet}_k)$ is finitely generated as an $A$-module, and $\H^{q+1}(a^{k+1}K^{\bullet})$ is a coherent $A$-module for any $k\geq n$, then $\H^q(K^{\bullet})$ is a finitely generated $A$-module.

{\rm (2)} If $\H^q(K^{\bullet})$ is finitely generated as an $A$-module and $\H^q(K^{\bullet}_k)$ is a coherent $A$-module for any $k\geq n$, then $\H^q(K^{\bullet})$ is a coherent $A$-module.
\end{lem}

\begin{proof}
(1) Since $\delta^q_k(\H^q(K^{\bullet}_k))$ is finitely presented, $\ker(\delta^q_k)=\image(\sigma^q_k)$ is finitely generated.
Then the assertion follows from \ref{lem-ullrichlemma2} (2).

(2) Set $T=\H^q(K^{\bullet})_{\ator}$.
Since $(A,a)$ is adhesive, $\H^q(K^{\bullet})/T$ is finitely presented and hence is coherent due to \ref{prop-cohringsmodules1}.
In view of \ref{lem-pf1} (2) it suffices to show that $T$ is finitely presented.
Since $T$ is finitely generated, there exists $m\geq 0$ such that $a^{m+1}T=T\cap a^{m+1}\H^q(K^{\bullet})=0$.
Let $k=m+n$.
Then $\ker(\sigma^q_k)=\image(\iota^q_k)\subseteq a^{m+1}\H^q(K^{\bullet})$ by \ref{lem-ullrichlemma2} (1).
Hence the map $T\rightarrow\H^q(K^{\bullet}_k)$ by $\sigma^q_k$ is injective, and thus $T$ can be regarded as a finitely generated $A$-submodule of the coherent $A$-module $\H^q(K^{\bullet}_k)$.
Hence $T$ is finitely presented, as desired.
\end{proof}

The following lemma will be used later:
\begin{lem}\label{lem-ullrichlemma3}
Consider the situation as in {\rm \S\ref{subsub-paircomplexgen}}, where the ideal $I$ is assumed to be principal $I=(a)$ and $A$ is $a$-adically complete.
Let $B$ be an $a$-adically complete flat $A$-algebra such that the pair $(B,IB)$ satisfies {\bf (APf)}.
Let $L^{\bullet}$ be a complex of $B$-modules, and $K^{\bullet}\rightarrow L^{\bullet}$ a morphism of complexes of $A$-modules.
By the commutative diagram
$$
\xymatrix@-1ex{0\ar[r]&a^{k+1}K^{\bullet}\ar[d]\ar[r]&K^{\bullet}\ar[d]\ar[r]&K^{\bullet}_k\ar[d]\ar[r]&0\\ 0\ar[r]&a^{k+1}L^{\bullet}\ar[r]&L^{\bullet}\ar[r]&L^{\bullet}_k\ar[r]&0}
$$
of complexes with exact rows, one has the commutative diagram
$$
\xymatrix@C-2ex@R-2ex{\cdots\ar[r]&\H^q(a^{k+1}K^{\bullet})_B\ar[d]_{\varphi^q_k}\ar[r]^(.55){\iota^q_{k,B}}&\H^q(K^{\bullet})_B\ar[d]_{\varphi^q}\ar[r]^{\sigma^q_{k,B}}&\H^q(K^{\bullet}_k)_B\ar[d]_{\Phi^q_k}\ar[r]^(.39){\delta^q_{k,B}}&\H^{q+1}(a^{k+1}K^{\bullet})_B\ar[d]^{\varphi^{q+1}_k}\ar[r]&\cdots\\ \cdots\ar[r]&\H^q(a^{k+1}L^{\bullet})\ar[r]_(.55){\eta^q_k}&\H^q(L^{\bullet})\ar[r]_{\tau^q_k}&\H^q(L^{\bullet}_k)\ar[r]_(.39){\epsilon^q_k}&\H^{q+1}(a^{k+1}L^{\bullet})\ar[r]&\cdots}
$$
of $B$-modules with exact rows $($due to the flatness of $A\rightarrow B);$ here the exact sequence in the first row is the one obtained by base change by $A\rightarrow B$ of the cohomology exact sequence as in {\rm \S\ref{subsub-paircomplexgen}}.
Let $q\in\Z$ and suppose that $\H^q(K^{\bullet})$ is finitely generated as an $A$-module.

{\rm (1)} If $a^nK^{q+1}_{\ator}=0$ for some $n\geq 0$ and $\Phi^q_k$ is injective for any $k\geq n$, then $\varphi^q$ is injective.

{\rm (2)} Suppose that $K^{q-1}$ is $a$-adically complete, that $L^q$ is $a$-adically separated, and that $a^nL^{q+1}_{\ator}=0$ for some $n\geq 0$.
If for any $k\geq n$ the map $\Phi^q_k$ is surjective and the map $\varphi^{q+1}_k$ is injective, then $\varphi^q$ is surjective.
\end{lem}

\begin{proof}
(1) For any $k\geq n$ we have $\ker(\varphi^q)\subseteq\ker(\tau^q_k\circ\varphi^q)=\ker(\Phi^q_k\circ\sigma^q_{k,B})=\ker(\sigma^q_{k,B})=\image(\iota^q_{k,B})=\image(\iota^q_k)\otimes_AB\subseteq a^{k-n+1}(H^q(K^{\bullet})\otimes_AB)$, where the last inclusion is due to \ref{lem-ullrichlemma2} (1).
By \ref{prop-zariskipair2} we know that $H^q(K^{\bullet})\otimes_AB$ is $a$-adically separated.
Hence $\ker(\varphi^q)=0$.

(2) By an easy diagram chasing one deduces that $\tau^q_k(\varphi^q(\H^q(K^{\bullet})\otimes_AB))=\tau^q_k(\H^q(L^{\bullet}))$ for $k\geq n$.
Let $\alpha_1,\ldots,\alpha_m$ be elements of $\H^q(K^{\bullet})$ that generate $\H^q(K^{\bullet})$ as an $A$-module, and set $\beta_j=\varphi^q(\alpha_j\otimes 1_B)$ for each $j=1,\ldots,m$.
We know that the elements $\tau^q_k(\beta_1),\ldots,\tau^q_k(\beta_m)$ generate $\tau^q_k(\H^q(L^{\bullet}))$ as a $B$-module.
By \ref{lem-ullrichlemma2} (2) we deduce that $H^q(L^{\bullet})$ is generated by $\beta_1,\ldots,\beta_m$, that is, $\varphi^q$ is surjective.
\end{proof}

\addcontentsline{toc}{subsection}{Exercises}
\subsection*{Exercises}
\begin{exer}\label{exer-zriskiancovering}{\rm 
Let $(A,I)$ be a Zariskian pair, and $f_1,\ldots,f_r\in A$ finitely many elements such that $(f_1,\ldots,f_r)=A$.
Then show that 
$$
\coprod^r_{i=1}\Spec\zat{(A_{f_i})}\longrightarrow\Spec A
$$
gives a flat covering of $\Spec A$.}
\end{exer}

\begin{exer}\label{exer-aadictopcompleteseparatedclosed}{\rm 
Let $(A,I)$ be a complete pair with $I=(a)$ principal, and $f\colon M\rightarrow N$ a morphism of $A$-modules.
Suppose:
\begin{itemize}
\item[(a)] $M$ is $a$-adically complete, and $N$ is $a$-adically separated;
\item[(b)] $(\coker(f))_{\ator}$ is bounded $a$-torsion.
\end{itemize}
Then show that $\coker(f)$ is $a$-adically separated.}
\end{exer}

\begin{exer}\label{exer-closuresaturated}{\rm 
Let $(A,I)$ be a pair, $\mathfrak{a}$ an ideal of $A$, and $f\in A$ an element of $A$ such that the ideal $(f)$ is $I$-adically open.
Show that, if $\mathfrak{a}$ is $f$-saturated, then so is the closure $\ovl{\mathfrak{a}}$ in $A$ with respect to the $I$-adic topology.}
\end{exer}

\begin{exer}\label{exer-adhesivecounterexamples}{\rm 
Let $W$ be a discrete valuation ring, and $V$ a valuation ring for the residue field of $W$ with $0<\mathrm{ht}(V)<\infty$.
Let $\til{V}$ be the composite of $W$ and $V$.
Show that there exists $a\in\m_V$ such that $\til{V}$ is Noetherian outside $I=(a)$ and $a$-torsion free, but is not $a$-adically separated.
In particular, $\til{V}$ is not $a$-adically adhesive.}
\end{exer}

\begin{exer}{\rm 
(1) Show that adhesiveness $($resp.\ pseudo-adhesiveness$)$ is stable under \'etale extensions of rings.

(2) If $(A,I)$ is adhesive $($resp.\ pseudo-adhesive, resp.\ universally adhesive, resp.\ universally pseudo-adhesive$)$ and $\Spec\het{A}\setminus V(\het{I})$ is Noetherian, then show that $(\het{A},\het{I})$ is adhesive $($resp.\ pseudo-adhesive, resp.\ universally adhesive, resp.\ universally pseudo-adhesive$)$.}
\end{exer}

\begin{exer}\label{exer-intersecfg}{\rm 
Let $(A,I)$ be an adhesive pair, $M$ an $I$-torsion free $A$-module, and $N,P\subseteq M$ finitely generated $A$-submodules.
Show that $N\cap P$ is finitely generated.}
\end{exer}

\begin{exer}{\rm 
Let $A$ be a ring, and $I\subseteq A$ a finitely generated ideal.

(1) The following conditions are equivalent (cf.\ \ref{prop-adhesive}):
\begin{itemize}
\item[(a)] for any finitely generated $A$-module $M$ that is finitely presented outside $I$, $M_{\Itor}$ is finitely generated;
\item[(b)] for any finitely generated $A$-module $M$ that is finitely presented outside $I$, $M/M_{\Itor}$ is finitely presented;
\item[(c)] for any finitely generated $A$-module $M$ and any finitely generated outside $I$ $A$-submodule $N\subseteq M$, the $I$-saturation $\til{N}$ is finitely generated.
\end{itemize}

If one of (and hence all) these conditions are satisfied, then we say that the pair $(A,I)$ is {\em pre-adhesive}.
If, moreover, any polynomial rings over $A$ together with the ideal induced from $I$ is pre-adhesive, we call $(A,I)$ {\em universally pre-adhesive}.
Clearly, a pair $(A,I)$ is adhesive (resp.\ universally adhesive) if and only if it is pre-adhesive (resp,\ universally pre-adhesive) and $A$ is Noetherian outside $I$.

(2) Show that the statements analogous to \ref{prop-adhesive1} and \ref{prop-uadhesive1} hold; in particular, universally pre-adhesiveness is stable under finite type extensions.

(3) Consider an inductive system $\{B_i,f_{ij}\colon B_i\rightarrow B_j\}$ of $A$-algebras indexed by a directed set $J$, and set $B=\varinjlim_{i\in I}B_i$.
Suppose:
\begin{itemize}
\item[{\rm (i)}] $(B_i,IB_i)$ is pre-adhesive $($resp.\ universally pre-adhesive$)$ for each $i;$
\item[{\rm (ii)}] for each pair of indices $i,j$ such that $i\leq j$ the morphism $f_{ij}$ is flat.
\end{itemize}
Then show that $(B,IB)$ is pre-adhesive $($resp.\ universally pre-adhesive$)$ (cf.\ \ref{prop-adhesive4}).

(4) Show that, if $(A,I)$ is pre-adhesive (resp.\ universally pre-adhesive), then so is the henselization $(\het{A},\het{I})$.}
\end{exer}

\begin{exer}\label{exer-Ivaluativeringsmaps}{\rm 
Let $A$ $($resp.\ $A')$ be $a$-valuative $($resp.\ $a'$-valuative$)$ local ring, where $a\in\m_A$ $($resp.\ $a'\in\m_{A'})$, and $h\colon A\rightarrow A'$ a local homomorphism that is adic with respect to the $a$-adic and $a'$-adic topologies.
Let $J=\bigcap_{n\geq 1}a^nA$, $B=A[\frac{1}{a}]$, $V=A/J$, and $K=\Frac(V)$ be as in \ref{thm-valuative}, and similarly, $J'=\bigcap_{n\geq 1}a^{\prime n}A'$, $B'=A'[\frac{1}{a'}]$, $V'=A'/J'$, and $K'=\Frac(V')$.

{\rm (1)} Show that the map $h$ induces a local injection $V\hookrightarrow V'$ and, moreover, we have $V=K\cap V'$ in $K'$.

{\rm (2)} Show that the map $h$ induces a local homomorphism $g=h[\frac{1}{a}]\colon B\rightarrow B'$ and, moreover, we have $g^{-1}(A')=A$.}
\end{exer}


\section{Topological algebras of type (V)}\label{sec-aadicallycompval}
\index{algebra!topological algebra of type (V)@topological --- of type (V)|(}
By a {\em topological algebra of type} (V), we mean a topologically finitely generated algebra over an $a$-adically complete valuation ring.
The terminology echos our later introducing notion, the {\em rigid spaces of type} (V) ({\bf \ref{ch-rigid}}, \S\ref{subsub-examplesV}), which will be defined consistently as rigid spaces locally of finite type over an $a$-adically complete valuation ring\index{valuation!valuation ring@--- ring!a-adically complete valuation ring@$a$-adically complete --- ---}.
Like that finitely generated algebras over a field are the `coordinate rings' and play a dominating role in classical algebraic geometry, topological algebras of type (V), especially those over an $a$-adically complete valuation ring of height one, are the central object in the classical rigid geometry of Tate\index{Tate, J.} and Raynaud\index{Raynaud, M.}.
More precisely, if $A$ is a topologically finitely generated algebra over an $a$-adically complete valuation ring $V$, then the algebra $A[\frac{1}{a}]$ is an {\em affinoid algebra}\index{algebra!affinoid algebra@affinoid ---} in Tate's rigid analytic geometry (\cite{Tate1}).
Many of the known properties of affinoid algebras, ring-theoretic (e.g.\ Noetherianness) or Banach algebra-theoretic, are, as we will see in \S\ref{subsub-convpreadh}, actually the consequences of the fact that $V$ is $a$-adically topologically universally adhesive (\ref{dfn-topadhesive}).

In \S\ref{sub-aadicallycompletevalrings} we give a generality of $a$-adically complete valuation rings and $a$-adic completion of valuation rings.
The main theorem (\ref{thm-compval2006ver1}) describes in detail the structure of the $a$-adic completion of a valuation ring, which will be at the basis of the discussions in the subsequent subsections.

In \S\ref{subsub-convpreadhnormalization} we will discuss the Noether normalization theorem for topologically finitely generated algebra over an $a$-adically complete valuation ring of height one, which is one of the `special techniques' valid only in height one situation.

The final subsection \S\ref{sub-classicalaffinoidalgebras} surveys several important and basic properties of affinoid algebras.

\medskip\noindent
{\bf Convention.} {\sl In this book, affinoid algebras obtained as $A[\frac{1}{a}]$ from a topologically finitely generated algebra $A$ over an $a$-adically complete valuation ring are called {\em classical} affinoid algebras, in order to distinguish them from more general algebras that arise from affinoids\index{affinoid} {\rm ({\bf \ref{ch-rigid}}, \S\ref{sec-affinoids})} in our approach to rigid geometry.}

\subsection{$a$-adic completion of valuation rings}\label{sub-aadicallycompletevalrings}
\subsubsection{Fundamental structure theorem}\label{subsub-aadicallycompletevalrings}
Throughout this section we fix a valuation ring\index{valuation!valuation ring@--- ring} $V$ and a non-zero element $a\in\mathfrak{m}_V\setminus\{0\}$.
(Notice that the valuation ring $V$ is therefore of non-zero height, that is, not a field.)
We consider the $a$-adic topology\index{topology!adic topology@adic ---}\index{adic!adic topology@--- topology} on $V$ and denote by 
$$
\widehat{V}=\varprojlim_{n\geq 0}V/a^nV
$$
the $a$-adic completion\index{completion!I-adic completion@$I$-adic ---} of $V$ (cf.\ \ref{prop-Iadiccompletioncomplete2}).
Furthermore, we set $K=\Frac(V)$.
\begin{thm}\label{thm-compval2006ver1}
Suppose $V$ is $a$-adically separated\index{valuation!valuation ring@--- ring!a-adically separated valuation ring@$a$-adically separated --- ---}.

{\rm (1)} The $a$-adic completion $\widehat{V}$ is a valuation ring of non-zero height, and $\Frac(\widehat{V})$ is canonically isomorphic to the Hausdorff completion\index{completion!Hausdorff completion@Hausdorff ---} $\widehat{K}$ of $K$ with respect to the filtration $\{a^nV\}_{n\geq 0}$.

{\rm (2)} The canonical map $V\hookrightarrow\widehat{V}$ is a local homomorphism, that is, $\mathfrak{m}_{\widehat{V}}\cap V=\mathfrak{m}_V$.
Moreover, the canonical map $V/\m_V\rightarrow\widehat{V}/\m_{\widehat{V}}$ is an isomorphism.

{\rm (3)} Let $\mathfrak{p}$ be a non-zero prime ideal of $V$.
Then the localization $V_{\mathfrak{p}}$ is $a$-adically separated, and $\widehat{V}$ is isomorphic to the composite of $\widehat{V}_{\mathfrak{p}}$, the $a$-adic completion of $V_{\mathfrak{p}}$, and $V/\mathfrak{p}$.

{\rm (4)} For any prime ideal $\mathfrak{p}$ of $V$, $\mathfrak{p}\widehat{V}$ is a prime ideal of $\widehat{V}$, and we have $\mathfrak{p}\widehat{V}\cap V=\mathfrak{p}$.

{\rm (5)} The underlying continuous mapping of the canonical morphism $\Spec\widehat{V}\rightarrow\Spec V$ is a homeomorphism$;$ in particular, we have $\mathrm{ht}(\widehat{V})=\mathrm{ht}(V)$.
Moreover, the canonical map $\Gamma_V=K^{\times}/V^{\times}\rightarrow\Gamma_{\widehat{V}}={\widehat{K}}^{\times}/{\widehat{V}}^{\times}$ between the value groups is an isomorphism of totally ordered commutative groups\index{ordered!totally@totally ---!totally ordered commutative group@--- --- commutative group}.
\end{thm}

Let us make the following small remark before proceeding to the proof: In the situation as in the theorem, we have the commutative square 
$$
\begin{xy}
(0,0)="Vhat"*{\widehat{V}},+<5em,0ex>="Khat"*{\widehat{K}},+<0em,-10ex>="K"*{K},+<-5em,0ex>="V"*{V}
\ar@{^{(}->}"Vhat"+<1em,-.5ex>;"Khat"+<-.8em,-.5ex>
\ar@{^{(}->}"V"+<1em,-.3ex>;"K"+<-.8em,-.3ex>
\ar@{^{(}->}"V"+<.1em,2.3ex>;"Vhat"+<.1em,-2ex>
\ar@{^{(}->}"K"+<.1em,2.3ex>;"Khat"+<.1em,-2ex>
\end{xy}
$$
consisting of injective ring homomorphisms. 
The injectivity of the vertical arrows follows from that the filtration $\{a^nV\}_{n\geq 0}$ is separated, and the injectivity of the upper horizontal arrow follows from Exercise \ref{exer-filtrationcompletion}.
It follows, moreover, from the same exercise that the above square is Cartesian, that is, $V=\widehat{V}\cap K$, and that $K/V\cong\widehat{K}/\widehat{V}$.
In particular, we have $\widehat{V}\neq\widehat{K}$.

\subsubsection{Proof of Theorem \ref{thm-compval2006ver1}}\label{subsub-aadicallycompletevalringsproof}
We first show the assertions (1) and (2) of the theorem.
To show that the completion $\widehat{K}$ of $K$ is a field, take $x\in\widehat{K}$.
The element $x$ of $\widehat{K}$ is the limit of a Cauchy sequence $\{x_n\}_{n\geq 0}$ in $K$.
If $x\neq 0$, replacing $\{x_n\}_{n\geq 0}$ by a cofinal subsequence, we may assume that there exists $N\geq 1$ such that $x_n\not\in a^NV$ for any $n\geq 0$.
Let us consider the sequence $\{x^{-1}_n\}_{n\geq 0}$, which we want to prove to be a Cauchy sequence.
For any $k\geq N$ there exists $y\in V$ such that $x_n-x_m=(x^{-1}_m-x^{-1}_n)x_nx_m=a^{k+2N}y$ for $n$ and $m$ large enough.
Since $x_n,x_m\not\in a^NV$, there exists $z_n,z_m\in V$ such that $a^N=x_nz_n=x_mz_m$.
Then we have $x^{-1}_m-x^{-1}_n=a^kz_nz_my\in a^kV$.
Hence $\{x^{-1}_n\}_{n\geq 0}$ is a Cauchy sequence that converges to the multiplicative inverse of $x$, and thus we showed that $\widehat{K}$ is a field.

Next, consider the maximal ideal $\m_V$ of $V$.
Since $a^nV\subseteq\m_V$ for each $n\geq 0$, one can consider the completion $\widehat{\m}_V$ of $\m_V$ by the filtration $\{a^nV\}_{n\geq 0}$.
By Exercise \ref{exer-filtrationcompletion} $\widehat{\m}_V$ is a subgroup of $\widehat{V}$; moreover, by the construction of completions $\widehat{\m}_V$ is canonically a $\widehat{V}$-module and hence is an ideal of $\widehat{V}$.
By the same exercise we have $\widehat{V}/\widehat{\m}_V\cong V/\m_V$ and $\widehat{\m}_V\cap V=\m_V$.
By the first equality we deduce that $\widehat{\m}_V$ is a maximal ideal of $\widehat{V}$, and the second one shows that this maximal ideal dominates the maximal ideal of $V$.
Moreover, one can show that $\widehat{\m}_V$ is the unique maximal ideal of $\widehat{V}$ as follows.
Consider an element $1+x$ with $x\in\widehat{\m}_V$.
Since $x$ is the limit of a Cauchy sequence $\{x_n\}_{n\geq 0}$ with each $x_n\in\m_V$, and since $1+x_n$ is invertible in $V$, we deduce that $1+x$ is invertible in $\widehat{V}$.
But this shows that $\widehat{\m}_V$ is contained in the Jacobson radical of $\widehat{V}$ (cf.\ \ref{prop-zariskipair1}) and hence that it is the unique maximal ideal.

Thus we have shown that $\widehat{V}$ is a local ring with the maximal ideal $\m_{\widehat{V}}=\widehat{\m}_V$ and that the inclusion map $V\hookrightarrow\widehat{V}$ is a local homomorphism.
To show that $\widehat{V}$ is a valuation ring for $\widehat{K}$, take a non-zero $x\in\widehat{K}$, for which we are going to check the condition {\rm (b)} in \ref{dfn-val}.
Since $\widehat{K}/\widehat{V}\cong K/V$, one can write $x=y+z$, where $y\in\widehat{V}$ and $z\in (K\setminus V)\cup\{0\}$.
If $z=0$, then $x\in\widehat{V}$.
If not, there exists $w\in\m_V$ such that $z=w^{-1}$.
Then $x^{-1}=w/(1+yw)$.
Since $w\in\m_V$ and since the inclusion map $V\hookrightarrow\widehat{V}$ is local, we have $yw\in\m_{\widehat{V}}$, and hence $1+yw$ is invertible in $\widehat{V}$.
By this we conclude $x^{-1}\in\widehat{V}$, and thus we have shown that $\widehat{V}$ is a valuation ring for $\widehat{K}$.
Moreover, since $\widehat{V}\neq\widehat{K}$, $\widehat{V}$ is of non-zero height.
Thus we have shown all assertions in (1) and (2) of the theorem.

To continue the proof, we need to prepare some auxiliary results:
\begin{lem}\label{lem-compval2006ver1}
Let $V$ be an $a$-adically separated valuation ring, and $S\subseteq V$ a multiplicative subset such that $S\cap\sqrt{(a)}=\emptyset$.
Then for any $n\geq 0$ we have
$$
a^{n+1}V_S\subseteq a^nV\subseteq a^nV_S.
$$
In particular, the $a$-adic topology on the ring $V_S$ coincides with the one given by the filtration $\{a^nV\}_{n\geq 0}$.
\end{lem}

\begin{proof}
The inclusion $a^nV\subseteq a^nV_S$ is obvious.
Let $a^{n+1}b/c\in a^{n+1}V_S$ where $c\in S$.
Since $c\not\in (a)$, there exists $d\in V$ such that $a=cd$.
Hence $a^{n+1}b/c=a^nbd\in a^nV$, thereby the other inclusion.
\end{proof}

\begin{cor}\label{cor-compval2006ver11}
Let $V$ be an $a$-adically separated valuation ring, and $\mathfrak{p}$ a non-zero prime ideal of $V$.
Then $V_{\mathfrak{p}}$ is $a$-adically separated, and we have $\Frac(\widehat{V})=\Frac(\widehat{V}_{\mathfrak{p}})$, where $\widehat{V}_{\mathfrak{p}}$ denotes the $a$-adic completion of $V_{\mathfrak{p}}$.
\end{cor}

\begin{proof}
In view of \ref{prop-maxspe} we know that $\mathfrak{p}$ contains $\sqrt{(a)}$; applying \ref{lem-compval2006ver1}, we have $\bigcap_{n\geq 0}a^nV_{\mathfrak{p}}=\bigcap_{n\geq 0}a^nV=0$, which shows that $V_{\mathfrak{p}}$ is $a$-adically separated.
Moreover, by \ref{lem-compval2006ver1} and Exercise \ref{exer-filtrationcompletion} we have $\widehat{V}\subseteq\widehat{V}_{\mathfrak{p}}\subseteq\widehat{K}=\Frac(\widehat{V})$, whence the other assertion.
\end{proof}

\begin{lem}\label{lem-compval2006ver2}
For any $x\in\widehat{K}$ there exists $u\in\widehat{V}^{\times}$ such that $ux\in K$.
Moreover, if $x\in\widehat{V}$, then $ux\in V$.
\end{lem}

\begin{proof}
We may assume $x\neq 0$.
The element $x$ is the limit of a Cauchy sequence $\{x_n\}_{n\geq 0}$ consisting of elements in $K$; since $x\neq 0$, we may assume that $x_n\neq 0$ for any $n\geq 0$.
Since the sequence $\{x_n/x\}_{n\geq 0}$ converges to $1$, one finds a sufficiently large $n$ such that $u=x_n/x$ belongs to the open neighborhood $1+a\widehat{V}$ of $1$.
Since $1+a\widehat{V}\subseteq 1+\m_{\widehat{V}}$, we have $u\in\widehat{V}^{\times}$, and thus the first claim is shown.
If $x\in\widehat{V}$, then $ux\in\widehat{V}\cap K=V$.
\end{proof}

Now we proceed to the proof of the assertions (3), (4), and (5) of \ref{thm-compval2006ver1}.
Let $\mathfrak{p}$ be a non-zero prime ideal of $V$.
By \ref{cor-compval2006ver11} we already know that $V_{\mathfrak{p}}$ is $a$-adically separated.
By (2) we know that the residue field of $\widehat{V}_{\mathfrak{p}}$ is equal to that of $V_{\mathfrak{p}}$, and hence one can take the composite valuation ring $\til{V}$ of $\widehat{V}_{\mathfrak{p}}$ and $V/\mathfrak{p}$.
By \ref{cor-compval2006ver11} $\til{V}$ is a valuation ring for $\Frac(\widehat{V})=\widehat{K}$.
Moreover, by the construction we have the canonical inclusion map $\widehat{V}\hookrightarrow\til{V}$, which is local.
Hence, by the characterization {\rm (a)} in \ref{dfn-val} we have $\widehat{V}=\til{V}$, whence (3).

To show (4), let us consider the completion $\widehat{\mathfrak{p}}$ of $\mathfrak{p}$ with respect to the filtration $\{a^nV\}_{n\geq 1}$ (notice that $a^nV$ is contained in $\mathfrak{p}$ for any $n\geq 1$ (\ref{prop-maxspe})).
By Exercise \ref{exer-filtrationcompletion} we know that $\widehat{\mathfrak{p}}$ is an ideal of $\widehat{V}$; since $\widehat{V}/\widehat{\mathfrak{p}}\cong V/\mathfrak{p}$, it is a prime ideal; we have, moreover, $\widehat{\mathfrak{p}}\cap V=\mathfrak{p}$.
Now, by the construction we have the canonical inclusion $\mathfrak{p}\widehat{V}\subseteq\widehat{\mathfrak{p}}$.
Moreover, since $\mathfrak{p}\widehat{V}$ contains $a\widehat{V}$, we have $V/\mathfrak{p}\cong\widehat{V}/\mathfrak{p}\widehat{V}$.
Hence we have $\mathfrak{p}\widehat{V}=\widehat{\mathfrak{p}}$, which shows (4).

Since $\widehat{V}$ is $a$-adically separated, any non-zero prime ideal contains $a$.
Hence, to show the first assertion of (5), it suffices to show that the map $\Spec\widehat{V}/a\widehat{V}\rightarrow\Spec V/aV$ is a homeomorphism, which is clear, since $\widehat{V}/a\widehat{V}\cong V/aV$.
This shows the equality $\mathrm{ht}(\widehat{V})=\mathrm{ht}(V)$ between the heights.

Finally, Let us show $\Gamma_V\cong\Gamma_{\widehat{V}}$.
Since $\widehat{V}\cap K=V$ and $\m_{\widehat{V}}\cap V=\m_V$, it is easy to see that the map $K^{\times}/V^{\times}\rightarrow\widehat{K}^{\times}/\widehat{V}^{\times}$ is injective.
The surjectivity follows from \ref{lem-compval2006ver2}. \hfill$\square$

\subsubsection{Corollaries}\label{subsub-aadicallycompletevalringscor}
Theorem \ref{thm-compval2006ver1} has several useful corollaries:
\begin{cor}\label{cor-compval2006ver2}
Let $V$ be a valuation ring, and $a\in\m_V\setminus\{0\}$.
Then the $a$-adic completion $\widehat{V}=\varprojlim_{n\geq 0}V/a^nV$ is a valuation ring of non-zero height.
\end{cor}

\begin{proof}
By \ref{prop-associatedsepval} $J=\bigcap_{n\geq 0}a^nV$ is a prime ideal, and hence $V/J$ is a valuation ring (\ref{prop-composition1} (1)).
Since $\widehat{V}$ is isomorphic to the $a$-adic completion of $V/J$ and $a\not\in J$, the assertion follows from \ref{thm-compval2006ver1} (1).
\end{proof}

\begin{rem}\label{rem-compval2006ver2}{\rm 
Here is another proof of \ref{cor-compval2006ver2}.
Since we have a canonical bijection between the set of all $a$-admissible\index{admissible!I-admissible@$I$-{---}} ideals (\ref{dfn-adm}) of $\widehat{V}=\varprojlim_{n\geq 0}V/a^nV$ and the set of $a$-admissible ideals of $V$, we see that $\widehat{V}$ is $a$-valuative\index{valuative!Ivaluative ring@$I$-{---} ring} (\ref{dfn-valuative1}).
Hence by \ref{thm-valuative} $\widehat{V}=\widehat{V}/\bigcap_{n\geq 0}a^n\widehat{V}$ is a valuation ring.}
\end{rem}

\begin{cor}\label{cor-compval2006ver21}
Let $V$ be an $a$-adically separated valuation ring with $a\in\m_V\setminus\{0\}$, and $S\subseteq V$ a multiplicative subset such that $S\cap\sqrt{(a)}=\emptyset$.
Then the canonical map
$$
(\widehat{V})_S\longrightarrow\widehat{V}_S
$$
is an isomorphism, where the left-hand ring is the ring of fractions of $\widehat{V}$ with respect to the image of $S$ by the canonical map $V\hookrightarrow\widehat{V}$, and the right-hand ring is the $a$-adic completion of $V_S$.
\end{cor}

\begin{proof}
By \ref{lem-compval2006ver1} and Exercise \ref{exer-filtrationcompletion} the canonical map $\widehat{V}\rightarrow\widehat{V}_S$ is injective, and hence the map in question is injective.
To show it is surjective, we write each element $x$ of $\widehat{V}_S$ as a power series of $a$
$$
x=\sum_{n\geq 0}\frac{b_n}{s_n}a^n,
$$
where $b_n\in V$ and $s_n\in S$ for $n\geq 0$.
Since $s_n\not\in (a)$, there exists $t_n\in V$ such that $a=s_nt_n$.
Hence we have
$$
x=\frac{b_0}{s_0}+\sum_{n\geq 1}b_nt_na^{n-1},
$$
which belongs to $(\widehat{V})_S$.
\end{proof}

\begin{cor}\label{cor-compval2006ver221}
If $V$ is an $a$-adically complete valuation ring\index{valuation!valuation ring@--- ring!a-adically complete valuation ring@$a$-adically complete --- ---} $(a\in\m_V\setminus\{0\})$ and $S\subseteq V$ is a multiplicative subset such that $S\cap\sqrt{(a)}=\emptyset$, then $V_S$ is $a$-adically complete. \hfill$\square$
\end{cor}

\begin{cor}\label{cor-compval2006ver22}
Let $V$ be an $a$-adically separated valuation ring $(a\in\m_V\setminus\{0\})$, and $\mathfrak{p}$ a non-zero prime ideal of $V$.
Then we have $(\widehat{V})_{\mathfrak{p}\widehat{V}}\cong\widehat{V}_{\mathfrak{p}}$.
\end{cor}

\begin{proof}
In view of \ref{cor-compval2006ver21} it suffices to show that the localization of $\widehat{V}$ by $S=V\setminus\mathfrak{p}$ is equal to the localization by $\widehat{V}\setminus\mathfrak{p}\widehat{V}$.
By \ref{lem-compval2006ver2}, for any element $x\in\widehat{V}\setminus\mathfrak{p}\widehat{V}$ there exists $u\in\widehat{V}$ such that $ux\in V\cap(\widehat{V}\setminus\mathfrak{p}\widehat{V})$.
Since $\mathfrak{p}\widehat{V}\cap V=\mathfrak{p}$ (\ref{thm-compval2006ver1} (4)), we have $ux\in S$.
The claim follows immediately from this.
\end{proof}

\begin{cor}\label{cor-compval2006ver23}
Let $V$ be an $a$-adically separated valuation ring $(a\in\m_V\setminus\{0\})$.
Let $\mathfrak{p}=\sqrt{aV}$ $($resp.\ $\widehat{\mathfrak{p}}=\sqrt{a\widehat{V}})$ be the associated height one prime\index{associated height one prime} of $V$ $($resp.\ $\widehat{V})$ $($cf.\ $\ref{dfn-maxspe2})$.
Then we have $(\widehat{V})_{\widehat{\mathfrak{p}}}\cong\widehat{V}_{\mathfrak{p}}$.
\end{cor}

\begin{proof}
In view of \ref{cor-compval2006ver22} it suffices to show that $\widehat{\mathfrak{p}}=\mathfrak{p}\widehat{V}$.
Since $\mathfrak{p}\widehat{V}\cap V=\mathfrak{p}$ (\ref{thm-compval2006ver1} (4)) and since $\Spec\widehat{V}\rightarrow\Spec V$ is a homeomorphism (\ref{thm-compval2006ver1} (5)), both $\mathfrak{p}\widehat{V}$ and $\widehat{\mathfrak{p}}$ are height one primes of $\widehat{V}$.
Hence we have $\widehat{\mathfrak{p}}=\mathfrak{p}\widehat{V}$.
\end{proof}

\subsection{Topologically finitely generated {\boldmath $V$-algebras}}\label{sub-convpreadh}
\index{finitely generated!topologically finitely generated@topologically ---}
\subsubsection{Adhesiveness}\label{subsub-convpreadh}
We consider in the sequel of this section an {\em $a$-adically complete}\index{complete!I-adically complete@$I$-adically ---} valuation ring $V$ with $a\in\m_V\setminus\{0\}$.
As defined in \ref{dfn-topfinigen}, a topologically finitely generated $V$-algebra is an $a$-adically complete $V$-algebra of the form 
$$
A=V\dl X_1,\ldots,X_n\dr/\mathfrak{a},
$$
where $\mathfrak{a}$ is a closed ideal of the restricted formal power series ring\index{restricted formal power series} $V\dl X_1,\ldots,X_n\dr$.

We have already seen in \ref{prop-exaadhesiveval} that the an $a$-adically complete valuation ring\index{valuation!valuation ring@--- ring!a-adically complete valuation ring@$a$-adically complete --- ---} is $a$-adically adhesive\index{adhesive!adhesive pair@--- pair}\index{adhesive!Iadically adhesive@$I$-adically ---}.
The aim of this paragraph is to show a much stronger result, that $(V,a)$ is topologically universally adhesive\index{adhesive!universally adhesive@universally ---!topologically universally adhesive@topologically --- --- (t.u.\ adhesive)}\index{pair!adhesive pair@adhesive ---!topologically universally adhesive pair@topologically universally --- ---} (\ref{cor-convadh}).
Due to \ref{thm-tuaunivadhe} it suffices to show the following:
\begin{thm}\label{thm-cpsadhesive}
The pair $(V\dl X_1,\ldots,X_n\dr,a)$ is adhesive.
$($Hence any topologically finitely generated algebra $($cf.\ {\rm \ref{dfn-topfinigen})} over $V$ is $a$-adically adhesive.$)$
\end{thm}

The proof of the theorem will be done in two steps.
The first step deals with the special case where $V$ is of finite height, and the second step discusses the general case.
The argument of the first step is almost the same as \cite[Lemma 1.1.2]{Fujiw1}, of which the proof is only valid when the height of $V$ is finite. 
Our proof here fills in this gap and is valid for arbitrary height.

Before the proof we include here a useful result by Raynaud\index{Raynaud, M.} and Gruson\index{Gruson, L.} for the reader's convenience:
\begin{thm}[{\cite[Premi\`ere partie, Th\'eor\`eme (3.4.6)]{RG}}]\label{thm-RG}
Let $f\colon X\rightarrow S$ be a locally of finite presentation morphism of schemes, and $\mathscr{M}$ an $\O_X$-module of finite type. 
Suppose that the set of associated primes $\mathrm{Ass}(S)$ of $S$ to $\O_S$, that is, the set of point $s\in S$ such that $\mathfrak{m}_{S,s}$ is up to radical the annihilator of an element of $\O_{S,s}$, is a locally finite set.
Then the set $U$ of all points of $X$ where $\mathscr{M}$ is $S$-flat is open, and $\mathscr{M}|_U$ is a finitely presented $\O_U$-module. \hfill$\square$
\end{thm}

\begin{cor}[{\cite[Premi\`ere partie, Corollaire (3.4.7)]{RG}}]\label{cor-RG}
Any flat of finite type algebra over an integral domain is finitely presented. \hfill$\square$
\end{cor}

\begin{prop}\label{prop-convpreadhloc1}
Let $A$ be an $a$-adically complete $V$-algebra, and $S\subseteq V$ a multiplicative subset such that $S\cap\sqrt{aV}=\emptyset$.
Then the ring $A_S$ is $a$-adically complete, and hence we have the equalities
$$
A_S=A\otimes_VV_S=A\widehat{\otimes}_VV_S
$$
up to canonical isomorphisms.
\end{prop}

\begin{proof}
By a similar reasoning as in the proof of \ref{lem-compval2006ver1}, for any $s\in S$ there exists $t\in V$ such that $st=a$.
Hence for $n\geq 1$ the subset $a^nA_S$ of $A_S$ is contained in the image of $a^{n-1}A$ by the canonical map $A\rightarrow A_S$.
This shows that $A_S$ is $a$-adically separated, and hence that the canonical map $A_S\rightarrow\widehat{A}_S$ to the $a$-adic completion is injective.
To show that this map is bijective, observe that any element $x$ of $\widehat{A}_S$ is written as a power series $x=\sum_{n\geq 0}\frac{b_n}{s_n}a^n$ as in the proof of \ref{cor-compval2006ver21}, where $b_n\in A$ and $s_n\in S$ for each $n\geq 0$.
Then by the similar reasoning one finds that $x$ lies in $A_S$.
\end{proof}

\begin{proof}[Proof of Theorem {\rm \ref{thm-cpsadhesive}}]
Set $A=V\dl X_1,\ldots,X_n\dr$, and let $M$ be an $a$-torsion free finitely generated $A$-module.
We need to show that $M$ is finitely presented (cf.\ \ref{prop-adhesive} (b)).
Take a surjection $A^{\oplus N}\rightarrow M$ and consider the exact sequence
$$
0\longrightarrow L\longrightarrow A^{\oplus N}\longrightarrow M\longrightarrow 0.\leqno{(\ast)}
$$
We want to prove is that $L$ is finitely generated. This will be shown in two steps.

\medskip
{\sc Step 1.} Suppose that the height of $V$ is finite. 
Due to \ref{prop-complpair1} it is enough to show that $L/aL$ is a finitely generated $A/aA$-module or, equivalently, that $M/aM$ is a finitely presentated $A/aA$-module \cite[Chap.\ I, \S 2.8, Lemma 9]{Bourb1}.
Since $M$ has no $a$-torsion, it is $V$-flat and hence $M/aM$ is $(V/aV)$-flat.
Now, $A/aA$ is the polynomial ring $A/aA\cong (V/aV)[X_1,\ldots,X_n]$, and $(V/aV)_{\mathrm{red}}$ is a valuation ring (possibly a field), since $\mathfrak{p}=\sqrt{(a)}$ is a prime. 
Hence we can apply \ref{thm-RG} (here we use the hypothesis that the height is finite). 
Since $M/aM$ is a finitely generated flat $(A/aA)$-module, it follows that $M/aM$ is of finite presentation as desired, and the proof in this case is done.

\medskip
{\sc Step 2.} In general, consider be the associated height one prime $\mathfrak{p}=\sqrt{(a)}$ (\ref{dfn-maxspe2}), and set $V'=V_{\mathfrak{p}}$.
By \ref{cor-compval2006ver23} $V'$ is an $a$-adically complete valuation ring of height one.
Moreover, by \ref{prop-convpreadhloc1} we know that $A\otimes_VV'$ is isomorphic to $V'\dl X_1,\ldots,X_n\dr$.
Hence one can apply the argument of {\sc Step 1} to conclude that $M\otimes_VV'$ is an $A\otimes_VV'$-module of finite presentation and hence that $L\otimes_VV'$ is finitely generated over $A\otimes_VV'$.
Take $x_1,\ldots,x_d\in L$ that generate $L\otimes_VV'$.
We look at the exact sequence induced from $(\ast)$
$$
0\longrightarrow L/\mathfrak{p}L\longrightarrow (A/\mathfrak{p}A)^{\oplus N}\longrightarrow M/\mathfrak{p}M\longrightarrow 0;\leqno{(\ast\ast)}
$$
note that this is exact, since $L$ is $a$-saturated (and hence $M$ is $V$-flat).
Since $(\ast\ast)$ is an exact sequence of modules over a polynomial ring $A/\mathfrak{p}A\cong (V/\mathfrak{p})[X_1,\ldots,X_n]$, and since $M/\mathfrak{p}M$ is flat over $V/\mathfrak{p}$, it follows from \ref{cor-RG} that $L/\mathfrak{p}L$ is finitely generated over $(A/\mathfrak{p}A)$.
Hence one can take $y_1,\ldots,y_e\in L$ that generates $L/\mathfrak{p}L$.

Now we claim that $x_1,\ldots,x_d,y_1,\ldots,y_e$ generates $L$ as an $A$-module.
Take any $z\in L$. 
There exist $\alpha_1,\ldots,\alpha_e\in A$ such that 
$$
z-(\alpha_1y_1+\cdots+\alpha_ey_e)\in\mathfrak{p}L.
$$
Set $\gamma y=z-(\alpha_1y_1+\cdots+\alpha_ey_e)$ ($\gamma\in\mathfrak{p}$, $y\in L$).
We can find $\beta_1,\ldots,\beta_d\in A\otimes_VV'$ such that $y=\beta_1x_1+\cdots+\beta_dx_d$. 
Since $\mathfrak{p}V'\subseteq V$, we have $\gamma\beta_i\in A$ for $i=1,\ldots,d$, and hence $z=\alpha_1y_1+\cdots+\alpha_ey_e+(\gamma\beta_1)x_1+\cdots+(\gamma\beta_d)x_d$ gives an expression of $z$ by $A$-linear combination.
\end{proof}

\begin{rem}\label{rem-trick}{\rm 
It deserves to be pointed out that the above proof has the following basic feature, which seems applicable to many other situations.
It is divided into two parts, and this division arises from taking the associated height one prime $\mathfrak{p}=\sqrt{(a)}$ of the valuation ring $V$; that is to say, the division of the discussion precisely corresponds to the decomposition of the valuation ring $V$ into $V_{\mathfrak{p}}$ and $V/\mathfrak{p}$.
Unlike the first part, which discusses things over the height one generalization $V_{\mathfrak{p}}$, the second part deals with the situation over $V/\mathfrak{p}$ in which, since $a=0$, the issues usually come down to classical algebraic geometry or classical commutative algebra.
In this sense the first part is usually the essential part of the proof and much harder than the second part.
It is, however, manageable in principle and can be done with the aid of the finiteness of height; moreover, we have several special techniques valid in height one situation.
It is by this reason that it is important to develop specialized techniques for dealing with topological algebras of type (V) over height one base; in \S\ref{subsub-convpreadhnormalization} and the appendix \S\ref{sec-furthertech}, we will explain some of these techniques. 
For example, by using some of the results in \S\ref{sec-furthertech}, we can give an alternative argument for {\sc Step 1} of the above proof that does not use the result of Raynaud\index{Raynaud, M.} and Gruson\index{Gruson, L.}; see \ref{prop-cpsadhesiveht1}.}
\end{rem}

Let us also mention that there is a shorter proof of \ref{thm-cpsadhesive} due to O.\ Gabber\index{Gabber, O.}.
In this proof one only needs to argue as in {\sc Step 1} of the above proof; to show that $M/aM$ is finitely presented (even in case $V$ is of infinite height), we first notice that $M\otimes_V(V/aV)_{\mathrm{red}}$ is finite presented by \ref{cor-RG} and apply the following lemma, which complements results by Raynaud-Gruson:
\begin{lem}\label{lem-propcomplementsRaynaudGruson}
Let $A$ be a ring, $B=A[X_1\ldots X_n]$ a polynomial ring, and $M$ a finitely generated $B$-module. 
Suppose $M$ is $A$-flat. 
Then $M$ is finitely presented over $B$ if and only if $M\otimes_AA_{\red}$ is finitely presented over $B\otimes_AA_{\red}$.
\end{lem}

\begin{proof}
Take a finitely presented $B$-module $N$ and a surjective map $\varphi\colon N\rightarrow M$ such that $N\otimes_AA_{\red}\cong M\otimes_AA_{\red}$; to obtain these, write $M\cong B^{\oplus n}/K$ and take a finitely generated $B$-submodule $K_0\subseteq K$ such that $K_0\otimes_AA_{\red}=K\otimes_AA_{\red}$; then $N=B^{\oplus n}/K_0\rightarrow M$ gives the desired map.
By \cite[Premi\`ere partie, Th\'eor\`eme (3.4.1)]{RG} $M_{\mathfrak{p}}$ for any prime ideal $\mathfrak{p}\subseteq B$ is finitely presented over $B_{\mathfrak{p}}$.
Hence if $L=\ker(\varphi)$, $L_{\mathfrak{p}}$ is finitely generated over $B$.
Since $M$ is flat over $A$, we have the exact sequence
$$
0\longrightarrow L_{\mathfrak{p}}\otimes_AA_{\red}\longrightarrow N_{\mathfrak{p}}\otimes_AA_{\red}\longrightarrow M_{\mathfrak{p}}\otimes_AA_{\red}\longrightarrow 0, 
$$
whence $L_{\mathfrak{p}}\otimes_AA_{\red}=0$.
But since $L_{\mathfrak{p}}$ is finitely generated, one applies Nakayama's lemma to deduce $L_{\mathfrak{p}}=0$.
Since this holds for any prime ideal $\mathfrak{p}$ of $B$, we conclude that $\varphi$ is an isomorphism.
\end{proof}

Now by \ref{thm-tuaunivadhe} and \ref{thm-cpsadhesive} we have:
\begin{cor}[O.\ Gabber]\label{cor-convadh}
\index{Gabber, O.}
Let $V$ be an $a$-adically complete valuation ring\index{valuation!valuation ring@--- ring!a-adically complete valuation ring@$a$-adically complete --- ---} of arbitrary height\index{height!height of a valuation ring@--- (of a valuation (ring))}.
Then $V$ is $a$-adically topologically universally adhesive. 
$($Hence any topologically finitely generated $V$-algebra is $a$-adically topologically universally adhesive.$)$ \hfill$\square$
\end{cor}

From this and \ref{thm-pf02a} we deduce:
\begin{cor}\label{cor-cohschemes21ver2}
Let $V$ be an $a$-adically complete valuation ring of arbitrary height.
Then any algebraic space locally of finite presentation over $V$ is universally cohesive. \hfill$\square$
\end{cor}

Moreover, since any $V$-flat algebra is $a$-torsion free, we have:
\begin{cor}\label{cor-convadh2}
Let $V$ be an $a$-adically complete valuation ring of arbitrary height.
Then any $V$-flat topologically finitely generated algebra is topologically finitely presented. \hfill$\square$
\end{cor}

\subsubsection{Noether normalization}\label{subsub-convpreadhnormalization}
\index{Noether normalization|(}
\begin{thm}[Noether normalization theorem for topologically finitely generated $V$-algebras]\label{thm-noethernormalizationtype(V)}
Let $V$ be an $a$-adically complete $($for $a\in\m_V\setminus\{0\})$ valuation ring {\em of height one}, and $A$ a $V$-flat topologically finitely generated $V$-algebra.
Then there exists a finite injective $V$-morphism
$$
V\dl X_1,\ldots,X_d\dr\longhookrightarrow A
$$
with the $V$-flat cokernel.
\end{thm}

Before proceeding to the proof the theorem, let us introduce a useful notion. 
Consider the restricted formal power series ring $B=V\dl Y_1,\ldots,Y_n\dr$ (where $V$ is not necessarily of height one).
Each element $f\in B$ can be written as a power series
$$
f=\sum_{\nu_1,\ldots,\nu_n\geq 0}b_{\nu_1,\ldots,\nu_n}Y^{\nu_1}_1\cdots Y^{\nu_n}_n,
$$
where $b_{\nu_1,\ldots,\nu_n}\in V$ for any $\nu_1,\ldots,\nu_n\geq 0$.
The {\em content ideal}\index{ideal!content ideal@content ---}\index{content ideal} $\Cont(f)$ of $f$ is the ideal of $V$ generated by all the coefficients $b_{\nu_1,\ldots,\nu_n}$.
This ideal is actually finitely generated (Exercise \ref{exer-contentidealfg}) and hence is principal.
Notice that, if $\Cont(f)=(b)$ and $f\neq 0$, then $f$ is divisible by $b$ and $\Cont(f/b)=V$.

\begin{lem}\label{lem-noethernormalizationtype(V)}
Suppose $V$ is of height one, and let $\mathfrak{a}\subseteq V\dl Y_1,\ldots,Y_n\dr$ be an $a$-saturated ideal.
If $\mathfrak{a}\neq 0$, then there exists a non-zero $f\in\mathfrak{a}$ such that $\Cont(f)=V$.
\end{lem}

\begin{proof}
Take $f\in\mathfrak{a}\setminus\{0\}$, and consider the content ideal $\Cont(f)=(b)$.
Since $V$ is of height one, there exists $m\geq 0$ and $c\in V$ such that $bc=a^m$ holds (\ref{prop-heightone1}).
Since $\mathfrak{a}$ is $a$-saturated, the element $cf/a^n$ belongs to $\mathfrak{a}$ and $\Cont(cf/a^n)=V$.
\end{proof}

\begin{proof}[Proof of Theorem {\rm \ref{thm-noethernormalizationtype(V)}}]
Let $k=V/\m_V$ be the residue field of $V$, and consider $A_0=A\otimes_Vk$, which is a finite type algebra over $k$.
By classical Noether normalization theorem (\cite[Chap.\ I, \S14]{Nagata1}) we can find elements $\ovl{x}_1,\ldots,\ovl{x}_d\in A_0$ algebraically independent over $k$ such that the map $k[\ovl{x}_1,\ldots,\ovl{x}_d]\hookrightarrow A_0$ is finite.
Take $x_1,\ldots,x_d\in A$ such that $\ovl{x}_i=(x_i$ mod $\m_V)$ for $i=1,\ldots,d$, and consider the $V$-subalgebra $A'\subseteq A$ topologically generated by them; $A'$ is the image of the unique morphism $\varphi\colon V\dl X_1,\ldots,X_d\dr\rightarrow A$ mapping each $X_i$ to $x_i$ ($i=1,\ldots,d$) (cf.\ \ref{lem-weaklygenerated} in the appendix).
Clearly, we have $A'_0=A'\otimes_Vk\cong k[\ovl{x}_1,\ldots,\ovl{x}_d]\subseteq A_0$.

We first claim that the map $\varphi$ is injective.
Indeed, if not, $\ker(\varphi)$ is a non-zero $a$-saturated ideal (since $A$ and $A'$ are $a$-torsion free by our hypothesis). 
Then by \ref{lem-noethernormalizationtype(V)} one has an element $f\in\ker(\varphi)$ such that $\ovl{f}=(f$ mod $\m_V)\neq 0$; but this would imply that the map $k[\ovl{x}_1,\ldots,\ovl{x}_d]\rightarrow A_0$ has non-zero kernel, contradicting that $\ovl{x}_1,\ldots,\ovl{x}_d$ are algebraically independent over $k$.
Now, since $\varphi_0=\varphi\otimes_Vk\colon A'_0\hookrightarrow A_0$ is finite and since $A'$ and $A$ are $a$-adically complete, we readily deduce that $\varphi$ is finite by \ref{prop-complpair1}.

Finally, let us show that the cokernel $A/A'$ is $a$-torsion free.
Since $V$ is of height one, one has $\m_V=\sqrt{(a)}$.
By this one verifies easily that for an $V$-module $M$ to be $a$-torsion free it is necessary and sufficient that $\Tor^V_1(M,k)=0$.
Hence, to verify the claim, it suffices to invoke that the map $A'_0=A'\otimes_Vk\hookrightarrow A_0=A\otimes_Vk$ is injective.
\end{proof}

\begin{cor}\label{cor-noethernormalizationtype(V)}
Let $V$ be as in {\rm \ref{thm-noethernormalizationtype(V)}}, and $A$ a $V$-flat quotient of $V\dl X_1,\ldots,X_n\dr$ such that the image of the closed immersion $\Spec A[\frac{1}{a}]\hookrightarrow\Spec V\dl X_1,\ldots,X_n\dr[\frac{1}{a}]$ is a finite set of closed points.
Then $A$ is finite over $V$.
\end{cor}

\begin{proof}
We have a finite injective $V$-morphism $V\dl Y_1,\ldots,Y_d\dr\hookrightarrow A$.
Since $\Spec A[\frac{1}{a}]$ consists of a single point, we have $d=0$.
\end{proof}
\index{Noether normalization|)}

\subsection{Classical affinoid algebras}\label{sub-classicalaffinoidalgebras}
\index{algebra!affinoid algebra@affinoid ---!classical affinoid algebra@classical --- ---|(}
\subsubsection{Tate algebra and classical affinoid algebras}\label{subsub-classicalaffinoidalgebras}
As in the previous subsections $V$ denotes an $a$-adically complete valuation ring with $a\in\m_V\setminus\{0\}$, and $K=\Frac(V)$ the field of fractions of $V$.
\begin{dfn}\label{dfn-classicalaffinoidalgebras}{\rm 
A {\em classical affinoid algebra}\index{algebra!affinoid algebra@affinoid ---!classical affinoid algebra@classical --- ---} is the $K$-algebra of the form
$$
\mathcal{A}=A\otimes_VK\ ({\textstyle =A[\frac{1}{a}]})
$$
(cf.\ \ref{prop-sep}), where $A$ is a topologically finitely generated $V$-algebra.}
\end{dfn}

Let $\mathfrak{p}=\sqrt{aV}$ be the associated height one prime of $V$ (\ref{dfn-maxspe2}), and $V'=V_{\mathfrak{p}}$ the associated height one localization of $V$.
Since $A\otimes_VK=(A\otimes_VV')\otimes_{V'}K$ and since $A\otimes_VV'$ is $a$-adically complete (hence is a topologically finitely generated $V'$-algebra) due to \ref{prop-convpreadhloc1}, we may always assume, whenever discussing classical affinoid algebras, that $V$ is of height one.
This reduction will be very helpful in developing generalities of classical affinoid algebras, for then one can use special techniques valid only in height one case, such as Noether normalization (\ref{thm-noethernormalizationtype(V)}).

The first general result on classical affinoid algebras is:
\begin{prop}\label{prop-classicalaffinoidalgnoe}
Any classical affinoid algebra is a Noetherian ring.
\end{prop}

\begin{proof}
Let $\mathcal{A}=A\otimes_VK$ be a classical affinoid algebra.
By \ref{thm-cpsadhesive} we know that $A$ is $a$-adically adhesive and hence is Noetherian outside $(a)$.
\end{proof}

Let us mention a special kind of affinoid algebras; this is the case where $A$ as above is the ring $V\dl X_1,\ldots,X_n\dr$ of restricted formal power series.
In this case the corresponding classical affinoid algebra $A\otimes_VK$ is the one usually called the {\em Tate algebra}\index{Tate, J.}\index{algebra!Tate algebra@Tate ---} and is denoted by $K\dl X_1,\ldots,X_n\dr$; explicitly, 
$$
K\dl X_1,\ldots,X_n\dr=\left\{
\begin{array}{r|l}
\sum_{\nu_1,\ldots,\nu_n\geq 0}a_{\nu_1,\ldots,\nu_n}X^{\nu_1}_1\cdots X^{\nu_n}_n&
|a_{\nu_1,\ldots,\nu_n}|\rightarrow 0\ \textrm{as}\\
\in K[\![X_1,\ldots,X_n]\!]&\nu_1+\cdots+\nu_n\rightarrow\infty
\end{array}
\right\}.
$$
Here the function $|\cdot|\colon K\rightarrow\R_{\geq 0}$ denotes a non-archimedean valuation associated to the height one localization $V'$ (cf.\ \S\ref{subsub-nonarchnorms}).
That $A$ is $a$-adically complete is then interpreted as that the algebra $K\dl X_1,\ldots,X_n\dr$ is a $K$-Banach algebra\index{algebra!Banach algebra@Banach ---} with respect to the norm
$$
\|\sum_{\nu_1,\ldots,\nu_n\geq 0}a_{\nu_1,\ldots,\nu_n}X^{\nu_1}_1\cdots X^{\nu_n}_n\|=\sup_{\nu_1,\ldots,\nu_n\geq 0}|a_{\nu_1,\ldots,\nu_n}|,
$$
called the {\em Gauss norm}\index{norm!Gauss norm@Gauss ---}.
Indeed, we have:
\begin{lem}\label{lem-tatealgtopcomp}
The topology on $K\dl X_1,\ldots,X_n\dr$ given by the Gauss norm $\|\cdot\|$ is equivalent to the $a$-adic topology, that is, the topology given by the filtration $\{a^nV\dl X_1,\ldots,X_n\dr\}_{n\geq 0}$.
\end{lem}

\begin{proof}
By \ref{lem-compval2006ver1} we may replace $V$ by the height one localization and thus may assume that $V$ is of height one; let $A=V\dl X_1,\ldots,X_n\dr$ and $B=A[\frac{1}{a}]=K\dl X_1,\ldots,X_n\dr$.
Set $\alpha=|a|$.
Then we have $0<\alpha<1$, since $a\in\m_V\setminus\{0\}$.
Since $a^nA\subseteq\{f\in B\,|\,\|f\|\leq\alpha^n\}$, it suffices to show that for any $n\geq 0$ there exists $\beta>0$ such that $\{f\in B\,|\,\|f\|\leq\beta\}\subseteq a^nA$.
If $g\in B$ satisfies $\|g\|\leq\alpha^n$, then every coefficient of $g$ has to be divisible by $a^n$.
Hence we have $\{f\in B\,|\,\|f\|\leq\alpha^n\}=a^nA$.
\end{proof}

\begin{prop}\label{prop-tatealgidealclosed}
Any ideal $\mathfrak{a}$ of $K\dl X_1,\ldots,X_n\dr$ is closed.
\end{prop}

\begin{proof}
Since $K\dl X_1,\ldots,X_n\dr$ is Noetherian, any ideal $\mathfrak{a}$ is finitely generated, and hence one can choose generators of $\mathfrak{a}$ from $V\dl X_1,\ldots,X_n\dr$.
Thus we have a finitely generated ideal $\til{\mathfrak{a}}$ of $V\dl X_1,\ldots,X_n\dr$ such that $\mathfrak{a}=\til{\mathfrak{a}}K\dl X_1,\ldots,X_n\dr$; replacing $\til{\mathfrak{a}}$ by its $a$-saturation, which is again finitely generated due to the adhesiveness of $V\dl X_1,\ldots,X_n\dr$, we may assume that $\til{\mathfrak{a}}$ is an $a$-saturated ideal.
Since $V\dl X_1,\ldots,X_n\dr$ is $a$-adically adhesive, we know by \ref{prop-AR} and \ref{cor-propARconseq1-2} that $\til{\mathfrak{a}}$ is closed with respect to the $a$-adic topology, that is, we have $\bigcap_{n\geq 0}(\til{\mathfrak{a}}+a^nV\dl X_1,\ldots,X_n\dr)=\til{\mathfrak{a}}$ (cf.\ Exercise \ref{exer-filtrationtopclosure}).
Now, to show the assertion, in view of \ref{lem-tatealgtopcomp} we only need to check the inclusion $\bigcap_{n\geq 0}(\mathfrak{a}+a^nV\dl X_1,\ldots,X_n\dr)\subseteq\mathfrak{a}$ (since the other inclusion is obvious).
Since $\til{\mathfrak{a}}$ is $a$-saturated, it is easy to see that 
$$
V\dl X_1,\ldots,X_n\dr\cap\bigcap_{n\geq 0}(\mathfrak{a}+a^nV\dl X_1,\ldots,X_n\dr)
=\bigcap_{n\geq 0}(\til{\mathfrak{a}}+a^nV\dl X_1,\ldots,X_n\dr)=\til{\mathfrak{a}}.
$$
Now for any $f\in\bigcap_{n\geq 0}(\mathfrak{a}+a^nV\dl X_1,\ldots,X_n\dr)$ there exists $n\geq 0$ such that $a^nf\in V\dl X_1,\ldots,X_n\dr$, and hence $a^nf\in\til{\mathfrak{a}}$; by this we conclude $f\in\mathfrak{a}$.
\end{proof}

By the proposition we deduce that any $K$-algebra $\mathcal{A}$ of the form
$$
\mathcal{A}=K\dl X_1,\ldots,X_n\dr/\mathfrak{a},
$$
where $\mathfrak{a}$ is an ideal of $K\dl X_1,\ldots,X_n\dr$, is a $K$-Banach algebra by the induced norm
$$
\|f\|=\inf_{F\mapsto f}\|F\|
$$
for $f\in\mathcal{A}$, where $F$ varies among all elements $F\in K\dl X_1,\ldots,X_n\dr$ that are mapped to $f$ by the canonical surjection $K\dl X_1,\ldots,X_n\dr\rightarrow\mathcal{A}$.

Moreover, any $K$-algebra $\mathcal{A}$ of this form is a classical affinoid algebra; indeed, as in the proof of \ref{prop-tatealgidealclosed}, one can take a finitely generated ideal $\til{\mathfrak{a}}$ of $V\dl X_1,\ldots,X_n\dr$ such that, if we set $A=V\dl X_1,\ldots,X_n\dr/\til{\mathfrak{a}}$ (which is topologically finitely generated over $V$), we have $\mathcal{A}=A\otimes_VK=A[\frac{1}{a}]$.
Replacing $\til{\mathfrak{a}}$ by its $a$-saturation (as in the proof of \ref{prop-tatealgidealclosed}), we may moreover find such an $A$ to be $a$-torsion free, hence flat over $V$ (cf.\ Exercise \ref{exer-valflatatorfree}).

\begin{dfn}\label{dfn-formalmodelclassical}{\rm 
Let $V$ be an $a$-adically complete valuation ring for $a\in\m_V\setminus\{0\}$ $($not necessarily of height one$)$, and $\mathcal{A}$ a classical affinoid algebra over $K=\Frac(V)$.

{\rm (1)} A {\em formal model}\index{formal model!formal model of a classical affinoid algebra@--- (of a classical affinoid algebra)} of $\mathcal{A}$ over $V$ is a topologically finitely generated $V$-algebra $A$ such that $A\otimes_VK\cong\mathcal{A}$.

{\rm (2)} A formal model $A$ of $\mathcal{A}$ is said to be {\em distinguished}\index{formal model!formal model of a classical affinoid algebra@--- (of a classical affinoid algebra)!distinguished formal model of a classical affinoid algebra@distinguished --- ---} if $A$ is $a$-torsion free.}
\end{dfn}

The term `distinguished' is coined in order to maintain consistency with our later terminology (cf.\ {\bf \ref{ch-rigid}}, \S\ref{subsub-cohrigspaceformalmodels}).
Distinguished formal models may also be called {\em flat formal models}, since $a$-torsion freeness is equivalent to $V$-flatness (Exercise \ref{exer-valflatatorfree}).
By what we have seen above, any classical affinoid algebra has a distinguished formal model.
Note that due to the adhesiveness of $V\dl X_1,\ldots,X_n\dr$ any distinguished formal model is topologically finitely presented (cf.\ proof of \ref{prop-tatealgidealclosed}).
\index{algebra!affinoid algebra@affinoid ---!classical affinoid algebra@classical --- ---|)}

\subsubsection{Ring-theoretic properties}\label{subsub-classicalaffinoidalgringprop}
Finally, let us discuss some ring-theoretic properties of classical affinoid algebras. 
We assume, without loss of generality, that the $a$-adically complete valuation ring $V$ is of height one.

First, by the Noether normalization theorem for topologically finitely presented $V$-algebras (\ref{thm-noethernormalizationtype(V)}), we get the similar theorem for classical affinoid algebras:
\begin{thm}[Noether normalization theorem for classical affinoid algebras]\label{thm-northernormaclassaff}
\index{Noether normalization}
For any classical affinoid algebra $\mathcal{A}$ over $K$ there exists an injective finite map of the form 
$$
K\dl T_1,\ldots,T_d\dr\longhookrightarrow\mathcal{A}.\eqno{\square}
$$
\end{thm}

\begin{cor}\label{cor-Tatealgebraprime}
For any maximal ideal $\mathfrak{m}$ of the Tate algebra\index{algebra!Tate algebra@Tate ---}\index{Tate, J.} $K\dl X_1,\ldots,X_n\dr$, the residue field $K\dl X_1,\ldots,X_n\dr/\mathfrak{m}$ is a finite extension of $K$. \hfill$\square$
\end{cor}

\begin{cor}[Weak Hilbert Nullstellensatz]\label{cor-dimensioncal2}
\index{Hilbert Nullstellensatz@(Weak) Hilbert Nullstellensatz}
Suppose $K$ is algebraically closed.
Then maximal ideals of $K\dl X_1,\ldots,X_n\dr$ are precisely the ideals of the form $(X_1-a_1,\ldots,X_n-a_n)$ with $a_1,\ldots,a_n\in V$.
\end{cor}

\begin{proof}
{\it Proof.} It is clear that the ideals of the form $(X_1-a_1,\ldots,X_n-a_n)$ $(a_1,\ldots,a_n\in V)$ are maximal. Let $\m\subseteq K\dl X_1,\ldots,X_n\dr$ be any maximal ideal. Since $K$ is algebraically closed, $K\dl X_1,\ldots,X_n\dr/\m$ is isomorphic to $K$ as topological rings. For $i=1,\ldots,n$, let $a_i\in K$ be the image of $X_i$. Since $X_i$ is power-bounded (cf.\ \S\ref{subsub-fadicringsgeneralities} in the appendix), so is $a_i$, and hence $a_i\in V$.
Since $\m$ contains $(X_1-a_1,\ldots,X_n-a_n)$, which is maximal, we have $\m=(X_1-a_1,\ldots,X_n-a_n)$.
\end{proof}

\begin{prop}\label{prop-tatealgebradim}
For any closed point $z$ of $\Spec K\dl X_1,\ldots,X_n\dr$ we have
$$
\dim_z(\Spec K\dl X_1,\ldots,X_n\dr)=n.
$$
In particular, the Krull dimension of the ring $K\dl X_1,\ldots,X_n\dr$ is equal to $n$.
\end{prop}

\begin{proof}
Consider the canonical map $f\colon\Spec K\dl X_1,\ldots,X_n\dr\rightarrow\Spec K[X_1,\ldots,X_n]$, which is flat since $V[X_1,\ldots,X_n]\hookrightarrow V\dl X_1,\ldots,X_n\dr$ is flat due to \ref{prop-btarf1} (2).
By Exercise \ref{exer-BGRlemma2pp261} (1) and \cite[$\mathbf{IV}$, (6.1.2)]{EGA} we have 
$$
\dim_z(\Spec K\dl X_1,\ldots,X_n\dr)=\dim_{f(z)}(\Spec K[X_1,\ldots,X_n])=n,
$$
as desired.
The last part of the proposition follows from \cite[$\mathbf{IV}$, (5.1.4)]{EGA}.
\end{proof}

\begin{prop}\label{prop-classicalaffringjacobson}
Classical affinoid algebras are Jacobson\index{Jacobson ring}.
\end{prop}

For the proof we need the following lemma:
\begin{lem}\label{lem-classicalaffringjacobson}
Let $F\in V\dl X_1,\ldots,X_n\dr$ be such that $\Cont(F)=V$ or, what amounts to the same, $\|F\|=1$, where $\|\cdot\|$ is the Gauss norm.
Then the following conditions are equivalent$:$
\begin{itemize}
\item[{\rm (a)}] $F$ is invertible in $K\dl X_1,\ldots,X_n\dr;$
\item[{\rm (b)}] $F$ is invertible in $V\dl X_1,\ldots,X_n\dr;$
\item[{\rm (c)}] the constant term $F_0$ is invertible in $V$, and $\Cont(F-F_0)\neq V$.
\end{itemize}
\end{lem}

\begin{proof}
Since $\|\cdot\|$ is a norm, the inverse $F^{-1}$ of $F$ in $K\dl X_1,\ldots,X_n\dr$ satisfies $\|F^{-1}\|=1$, that is, $F^{-1}\in V\dl X_1,\ldots,X_n\dr$, whence the equivalence of (a) and (b).
Suppose $F$ is invertible in $V\dl X_1,\ldots,X_n\dr$.
Then the image of $F$ by the residue map $V\dl X_1,\ldots,X_n\dr\rightarrow V\dl X_1,\ldots,X_n\dr/\mathfrak{m}_VV\dl X_1,\ldots,X_n\dr=k[X_1,\ldots,X_n]$ modulo the maximal ideal $\mathfrak{m}_V=\sqrt{(a)}$ of $V$, where $k$ is the residue field of $V$, is invertible, and this shows the implication (b) $\Rightarrow$ (c).
If, conversely, the conditions in (c) is satisfied, then we may assume that $F_0=1$, and then $F$ belongs to the subset of the form $1+a^rV\dl X_1,\ldots,X_n\dr$ (where $(a^r)=\Cont(F-F_0)$).
Since $V\dl X_1,\ldots,X_n\dr$ is $a$-adically complete and hence is $a$-adically Zariskian, we deduce that $F$ is invertible.
\end{proof}

\begin{proof}[Proof of Proposition {\rm \ref{prop-classicalaffringjacobson}}]
It is enough to show the proposition for Tate algebras $\mathcal{A}=K\dl X_1,\ldots,X_n\dr$ (\cite[$\mathbf{IV}$, (10.4.6)]{EGA}).
The proof is done by induction with respect to $n$.
If $n=0$, then $\mathcal{A}=K$ is a field, which is obviously Jacobson.
In general, what to prove is the following: for any prime ideal $\mathfrak{p}\subseteq\mathcal{A}$ the intersection of all maximal ideals containing $\mathfrak{p}$ coincides with $\mathfrak{p}$.
Suppose $\mathfrak{p}\neq 0$, and consider $\mathcal{B}=\mathcal{A}/\mathfrak{p}$.
Since $\mathfrak{p}$ has the positive height, the dimension of the ring $\mathcal{B}$ is strictly less than $n$, and hence by \ref{thm-northernormaclassaff} there exists an injective finite map $K\dl Y_1,\ldots,Y_d\dr\hookrightarrow\mathcal{B}$ with $d<n$.
By induction $K\dl Y_1,\ldots,Y_d\dr$ is Jacobson, and hence so is $\mathcal{B}$ by \cite[$\mathbf{IV}$, (10.4.6)]{EGA}; in other words, $\mathfrak{p}$ is the intersection of all maximal ideals that contain $\mathfrak{p}$.

It remains to show that the intersection of all maximal ideals of $\mathcal{A}$ is $(0)$.
Suppose $F\neq 0$ belongs to the intersection of all maximal ideals of $\mathcal{A}$.
We may assume that $F\in V\dl X_1,\ldots,X_n\dr$ and that $\Cont(F)=V$.
Let $F_0$ be the constant term of $F$.
Then $F-F_0$ belongs to the maximal ideal $(X_1,\ldots,X_n)$ of $\mathcal{A}$, and hence so does $F_0=F+(F_0-F)$.
Hence we must have $F_0=0$.
But then $1+F$ is not invertible, since $\Cont(F)=V$ (\ref{lem-classicalaffringjacobson}).
Hence there exists a maximal ideal $\mathfrak{m}$ of $\mathcal{A}$ such that $1+F\in\mathfrak{m}$.
But again this is absurd, since we would have $1=(1+F)-F\in\mathfrak{m}$.
\end{proof}

Finally, let us include without proofs a few more ring-theoretic facts, which will be needed in our later discussion:
\begin{prop}[{Exercise \ref{exer-BGRlemma2pp261} \& \cite[(5.2.6/1)]{BGR}}]\label{prop-exerTateUFD}
Tate algebras are regular and factorial. \hfill$\square$
\end{prop}

\begin{thm}[R.\ Kiehl {\cite{Kieh2}}]\label{thm-propclassicalaffinoidjapanese}
Any classical affinoid algebra is excellent. \hfill$\square$
\end{thm}
\index{algebra!topological algebra of type (V)@topological --- of type (V)|)}

\addcontentsline{toc}{subsection}{Exercises}
\subsection*{Exercises}
\begin{exer}\label{exer-valuapproxfiniteheight0}{\rm 
Let $V\hookrightarrow V'$ be an inclusion of $a$-adically separated valuation rings, where $a\in V$.
Show that the induced map $\widehat{V}\rightarrow\widehat{V'}$ between the $a$-adic completions is injective.}
\end{exer}

\begin{exer}\label{exer-valuapproxfiniteheight}{\rm 
Let $V$ be a valuation ring for $K=\Frac(V)$.

{\rm (1)} Show that $V$ is the filtered inductive limit\index{limit!inductive limit@inductive ---} of subrings $V=\varinjlim_{\lambda\in\Lambda}V_{\lambda}$ with each $V_{\lambda}\subseteq V$ being a valuation ring of finite height.

{\rm (2)} If, moreover, $V$ is $a$-adically separated (resp.\ complete) for $a\in\m_V\setminus\{0\}$, then one can find the inductive system of subrings $\{V_{\lambda}\}_{\lambda\in\Lambda}$ as above consisting of $a_{\lambda}$-adically separated (resp.\ complete) valuation rings of finite height.}
\end{exer}

\begin{exer}\label{exer-contentidealfg}{\rm 
Let $V$ be an $a$-adically complete valuation ring, and consider the restricted formal power series ring $V\dl X_1,\ldots,X_n\dr$.
Let $f\in V\dl X_1,\ldots,X_n\dr$.
Then show that the content ideal $\Cont(f)$ (\S\ref{subsub-convpreadhnormalization}) is finitely generated.}
\end{exer}

\begin{exer}\label{exer-BGRlemma2pp261}{\rm 
Consider the Tate algebra $\mathcal{A}=K\dl X_1,\ldots,X_n\dr$ and the canonical inclusion $A_0=K[X_1,\ldots,X_n]\hookrightarrow\mathcal{A}$.

(1) Show that for any maximal ideal $\m\subseteq\mathcal{A}$ its contraction $\m_0=\m\cap A_0$ is a maximal ideal of $A_0$ such that $\mathcal{A}/\m^{k+1}\cong A_0/\m^{k+1}_0$ and $\m^{k+1}_0\mathcal{A}=\m^{k+1}$ for any $k\geq 0$.

(2) Show that $\mathcal{A}=K\dl X_1,\ldots,X_n\dr$ is a regular ring.}
\end{exer}


\setcounter{section}{0}
\renewcommand{\thesection}{\Alph{section}}
\renewcommand{\theexer}{{\bf \thechapter}.\Alph{section}.\arabic{exer}}
\section{Appendix: Further techniques for topologically of finite type algebras}\label{sec-furthertech}
\subsection{Nagata's idealization trick}\label{sub-nagata1}
The so-called {\em Nagata's\index{Nagata, M.} idealization trick}\footnote{Nagata's idealization is the method to represent a given quasi-coherent sheaf on an affine scheme as the conormal sheaf of a closed immersion, viz.\ the zero section of the affine cone of the sheaf, thus regarding it as an ideal of the first infinitesimal neighborhood.}, stated in his book \cite[pp.\ 2]{Nagata1}, reduces many situations involving finitely generated modules to the special case of finitely generated ideals.
Let us first state its general principle. 

A ring homomorphism $f\colon S\rightarrow R$ is said to be a {\em thickening of order $\leq 1$}\index{thickening} if it is surjective and $K=\ker(f)$ satisfies $K^2=(0)$.
Let $R$ be a ring, and $M$ a finitely generated $R$-module.
We set $S=\Sym_RM/(M^2)$; $S$ is the direct sum $R\oplus M$ as an $R$-module with the multiplication by $(r+m)(r'+m')=rr'+rm'+r'm$. 
Then $S$ is a finite $R$-algebra containing $M$ as an ideal. 
Moreover, the projection map $S=R\oplus M\rightarrow R$ is a thickening of order $\leq 1$. 
Hence we get the diagram 
$$
\xymatrix{\Spec R\ar@<.5ex>[r]^(.5){\sigma}&\Spec S\rlap{,}\ar@<.5ex>[l]^(.5){\tau}}
$$
where $\tau$ is finite and $\sigma$ is the `zero section' to $\tau$.
Now the quasi-coherent sheaf $\til{M}$ of finite type on $\Spec R$ coincides with the pull-back of the quasi-coherent ideal $\til{M}$ of finite type on $\Spec S$; indeed, the two homomorphisms between $R$ and $S$ give rise to the isomorphism between $M$ and $M\otimes_SR$ as $R$-modules. 

By this observation the following proposition is now clear:
\begin{prop}[Nagata's idealization trick]\label{prop-nagata1}
\index{Nagata's idealization trick}
Let $A$ be a ring, and consider a property $P(R,M)$ involving an $A$-algebra $R$ of finite type and finitely generated $R$-module $M$. 
Suppose that for any thickening $S\rightarrow R$ of order $\leq 1$ between $A$-algebras of finite type, we have the implication $P(S,M)\Rightarrow P(R,M\otimes_SR)$.
Then the following conditions are equivalent$:$
\begin{itemize}
\item[{\rm (a)}] $P(R,I)$ holds for any $R$ and a finitely generated ideal $I\subseteq R;$
\item[{\rm (b)}] $P(R,M)$ holds for any $R$ and $M$. \hfill$\square$
\end{itemize}
\end{prop}

\subsection{Standard basis and division algorithm}\label{sub-stbasis}
The notion of standard (or Gr\"obner) basis, useful both in theories and calculations in algebraic geometry, is also useful to deal with the algebras $V\dl X_1,\ldots,X_n\dr$ over an $a$-adically complete valuation ring $V$ {\em of height one}.
As H.\ Hironaka\index{Hironaka, H.} has first envisaged in complex analytic situation, the resulting division theorem gives a generalization of Weierstrass preparation-division theorem; also in our situation it also gives an analogous division theorem, which will be of broad use as a theoretical and computational device.

\subsubsection{Setting}\label{subsub-stbasissit}
We consider $\N^n$ as an additive monoid in the standard way, and we equip it with a {\it term ordering}\index{ordering@order(ing)!term ordering@term ---} (see, e.g., \cite[1.4]{AL}), for example, the lexicographical order\index{ordering@order(ing)!lexicographical ordering@lexicographical ---} (cf.\ \ref{exa-height}).
For $\nu=(\nu_1,\ldots,\nu_n)\in\N^n$ we write $X^{\nu}=X^{\nu_1}_1\cdots X^{\nu_n}_n$.

Let $R$ be a ring and $I=(a_1,\ldots,a_r)\subseteq R$ a finitely generated ideal, and suppose that $R$ is $I$-adically complete.
For a restricted power series 
$$
f=\sum_{\nu\in\N^n}a_{\nu}X^{\nu}\in R\dl X_1,\ldots,X_n\dr,\leqno{(\ast)}
$$
the {\it content ideal}\index{ideal!content ideal@content ---}\index{content ideal} of $f$, denoted by $\Cont(f)$, is the ideal of $R$ generated by all the coefficients $a_{\nu}$ of $f$. 
It is easy to see that $\Cont(f)$ is actually finitely generated (cf.\ Exercise \ref{exer-contentidealfg}). 
We say that $f$ is {\it primitive}\index{primitive} if $\Cont(f)=R$.

If $R=V$ is a valuation ring, complete with respect to the $a$-adic topology ($I=(a)$), then the content ideal $\Cont(f)$ is always principal, and we can set
$$
\nu(f)=
\begin{cases}
\nu(f_0\ \mathrm{mod}\ \mathfrak{m}_V),\ \textrm{where $f_0=f/\alpha$ with $\Cont(f)=(\alpha)$,}&\textrm{if}\ f\neq 0,\\
-\infty&\textrm{otherwise.}
\end{cases}
$$
Here the {\it leading degree}\index{leading degree} $\nu(g)$ of a polynomial $g\in A[X_1,\ldots,X_n]$ over any ring $A$ is the maximal (with respect to the term order) multidegree among those appear in non-zero terms in $g$.
Note that the definition of $\nu(f)$ for $f\neq 0$ is also given by the following:
$$
\nu(f)=\sup\{\nu:a_{\nu}\ \textrm{generates}\ \Cont(f)\}
$$
for $f$ presented as $(\ast)$.
Define, furthermore,
$$
\LT(f)=a_{\nu(f)}X^{\nu(f)},
$$
and call it the {\it leading term}\index{leading term} of $f$.

\subsubsection{Division algorithm}\label{subsub-stbasis}
As in the previous paragraph, we consider an $I$-adically complete ring $R$, where $I=(a_1,\ldots,a_r)\subseteq R$ is finitely generated.

\begin{thm}[Division Lemma]\label{thm-division}\index{division lemma}
Let $\{g_{\lambda}\}_{\lambda\in\Lambda}$ be a finite collection of elements in $R\dl X_1,\ldots,X_n\dr$ such that, for each $\lambda\in\Lambda$, $g_{0\lambda}=(g_{\lambda}\ \mathrm{mod}\ IR\dl X_1,\ldots,X_n\dr)$ is a monic polynomial in $(R/I)[X_1,\ldots,X_n]$.
Set $\nu_{\lambda}=\nu(g_{0\lambda})$ for any $\lambda\in\Lambda$, and $M=\bigcup_{\lambda\in\Lambda}(\nu_{\lambda}+\N^n)$.
Then, for any $f\in R\dl X_1,\ldots,X_n\dr$, there exists $g\in R\dl X_1,\ldots,X_n\dr$ with no exponents in $M$ such that $f-g$ belongs to $\sum_{\lambda\in\Lambda}g_{\lambda}R\dl X_1,\ldots,X_n\dr$.
\end{thm}

\begin{proof}
By performing the division algorithm by $\{g_{0\lambda}\}_{\lambda\in\Lambda}$ in the polynomial ring $(R/I)[X_1,\ldots,X_n]$, there exists an expression
$$
f=g'+\sum_{\lambda\in\Lambda}\alpha_{\lambda}g_{\lambda}+\sum^r_{j=1}a_jf_j
$$
in $R\dl X_1,\ldots,X_n\dr$, where $g'$ has no exponents in $M$.
Repeating this, by induction, one deduces that, for any $m\geq 1$,
$$
f=g'_m+\sum_{\lambda\in\Lambda}\alpha_{\lambda,m}g_{\lambda}+\sum_{h\in P_m}h(a_1,\ldots,a_r)f_h,
$$
where the last sum is taken over all monomials $h$ of $r$ variables of degree $m$, such that:
\begin{itemize}
\item $g'_m$ has no exponents in $M$;
\item $g'_{m+1}-g'_m$ and $\alpha_{\lambda,m+1}-\alpha_{\lambda,m}$ belong to $I^mR\dl X_1,\ldots,X_n\dr$.
\end{itemize}
Then, with $g=\lim_{m\rightarrow\infty}g'_m$ and $\alpha_{\lambda}=\lim_{m\rightarrow\infty}\alpha_{\lambda,m}$, we have
$$
f=g+\sum_{\lambda\in\Lambda}\alpha_{\lambda}g_{\lambda},
$$
where $g$ has no exponents in $M$.
\end{proof}

\subsubsection{Standard bases}\label{subsub-standardbasis}
Set $R=V\dl X_1,\ldots,X_n\dr$, where $V$ is an $a$-adically complete valuation ring $(a\in\m_V\setminus\{0\})$.
For a non-zero ideal $I\subseteq R$ we denote by $\LT(I)$ the ideal of $R$ generated by the leading terms of all non-zero elements in $I$.

\begin{dfn}\label{dfn-stbasis}{\rm
Let $I$ be a non-zero ideal of $R$, and $g_1,\ldots,g_d\in I$ non-zero primitive elements.
The set $\{g_1,\ldots,g_d\}$ is said to be a {\em standard basis}\index{standard basis} for $I$ if 
$$
\LT(I)=(\LT(g_1),\ldots,\LT(g_d)).
$$}
\end{dfn}

\begin{cor}\label{cor-division}
If $\{g_1,\ldots,g_d\}$ is a standard basis of a non-zero ideal $I$ of $A$, we have $I=(g_1,\ldots,g_d)$.
\end{cor}

\begin{proof}
We only need to check the inclusion $I\subseteq(g_1,\ldots,g_d)$.
Let $f\in I$, and take $q_1,\ldots,q_d\in R$ such that $h=f-\sum^d_{i=1}q_ig_i$ has no exponent in $M=\bigcup^d_{i=1}(\nu(g_i)+\N^n)$ as in \ref{thm-division}.
Suppose $h\neq 0$.
As $h\in I$, we have $\LT(h)\in\LT(I)$.
But since the leading degree $\nu(h)$ does not belong to $M$, this is absurd.
Hence $h=0$ and thus we have $f\in(g_1,\ldots,g_d)$.
\end{proof}

Since $R$ is not necessarily Noetherian, it is not always the case that an ideal $I\subseteq R$ has a standard basis\index{standard basis}; however, we have:
\begin{prop}\label{prop-stbasis}
Let $I\subseteq R$ be a non-zero and $a$-saturated {\rm (\S\ref{subsub-pairstorsionsandsaturation})} ideal.
Then $I$ has a standard basis $($and hence is finitely generated$)$.
\end{prop}

\begin{proof}
Define the subset $L\subset\N^n$ by
$$
L=\{\nu(f)\,|\,f\in I\setminus\{0\}\}.
$$
Then $L$ is an ideal of the monoid $\N^n$, that is, for any $\nu\in L$ and $\mu\in\N^n$ we have $\nu+\mu\in L$.
By Dickson's lemma (cf.\ e.g.\ \cite[Excecise 1.4.12]{AL}) $L$ is finitely generated, that is, there exist $\nu_1,\ldots,\nu_r\in L$ such that 
$$
L=\bigcup^r_{i=1}(\nu_i+\N^n).
$$
Take $g_1,\ldots,g_r\in I$ such that $\nu(g_i)=\nu_i$ for $i=1,\ldots,r$. 
Since $I$ is $a$-saturated, we may assume that the coefficient of $\LT(g_i)$ is $1$ for each $i$.
Then for $f\in I$ there exists $i$ and a monomial $h$ such that $\LT(f)=h\LT(g_i)$.
Hence we have $\LT(I)=(\LT(g_1),\ldots,\LT(g_r))$.
\end{proof}

Finally, let us mention one application of the Nagata's\index{Nagata, M.} trick and the division algorithm.
We give an elementary proof, without referring to the result by Raynaud\index{Raynaud, M.} and Gruson\index{Gruson, L.} (\ref{thm-RG}), of \ref{thm-cpsadhesive} in the case $V$ is of height one:

\begin{prop}\label{prop-cpsadhesiveht1}
Let $V$ be an $a$-adically complete valuation ring of height one. 
Then any topologically finitely generated $V$-algebra is $a$-adically adhesive.
\end{prop}

\begin{proof}
By \ref{prop-adhesive1} (2), It suffices to show that $V\dl X_1,\ldots,X_n\dr$ is $a$-adically adhesive.
We apply Nagata's trick as in \ref{prop-nagata1}.
For any ring of the form $R=V\dl X_1,\ldots,X_n\dr/J$, for some $n\geq 0$ and some ideal $J$, and any finitely generated $R$-module $M$, let $P(R,M)$ be the property that for any $R$-submodule $N\subseteq M$ its $a$-saturation $\til{N}$ is finitely generated over $R$.
By \ref{prop-nagata1} we only have to prove $P(R,I)$ for any finitely generated ideal $I$ of $R$; moreover, one reduces to $P(R,R)$, that is, that for any ideal $I$ of $R$ its $a$-saturation $\til{I}$ is again finitely generated.
We may further assume that $R=V\dl X_1,\ldots,X_n\dr$.
But the assertion in this case is proved in \ref{prop-stbasis}.
\end{proof}

\addcontentsline{toc}{subsection}{Exercises}
\subsection*{Exercises}\label{sub-exercisefurthertech}
\begin{exer}\label{exer-adhesivenesscheckingnagata}{\rm 
Show that a complete pair of finite ideal type $(A,I)$ is adhesive if the following condition holds: if $B$ is an $I$-adically complete finite $A$-algebra, any $IB$-saturated ideal $J\subseteq B$ is finitely generated.}
\end{exer}

\begin{exer}\label{exer-exercisefurthertech2}{\rm 
Let $V$ be an $a$-adically complete valuation ring, and $V'$ the associated height one valuation ring (\S\ref{sub-adicsepval}).  

(1) Let $A$ be a topologically finitely generated $V$-algebra such that $A[\frac{1}{a}]$ is finite over $V[\frac{1}{a}]$. 
Show that $A$ is finite type over $V$ and that $(A/A_{\ator})\otimes_VV'$ is finite over $V'$. 

(2) Show that, conversely, a finite type $V$-algebra $A$ is $a$-adically complete if and only if $(A/A_{\ator})\otimes_VV'$ is finite over $V'$. 

(3) Show that any finite type $V$-algebra $A$ such that $A[\frac{1}{a}]$ is finite over $V[\frac{1}{a}]$ is canonically decomposed as $A=A'\times A''$, where $A'$ is $a$-adically complete and $A''\otimes_V(V/aV)=0$.}
\end{exer}

\begin{exer}[Weierstrass preparation theorem]\label{exer-preparationthm}
\index{Weierstrass preparation theorem@Weierstrass preparation theorem}{\rm 
Let $V$ be an $a$-adically complete valuation ring of height one, and $f\in V\dl X_1,\ldots,X_n\dr$ ($n\geq 1$) a primitive (that is, $\Cont(f)=V$) element.
Consider the lexicographical order (\ref{exa-height}) as the term order for the exponents of monomials, and set $\nu(f)=(\nu_1(f),\ldots,\nu_n(f))$.
Suppose that, if we write $f=\sum^{\infty}_{n=0}f_n(X_2,\ldots,X_n)X^n_1$, then $f_{\nu_1(f)}$ is a unit in $V\dl X_2,\ldots,X_n\dr$.
Then show that there exists a unique monic polynomial $g\in V\dl X_2,\ldots,X_n\dr[X_1]$ of degree $\nu_1(f)$ in $X_1$ and a unique unit element $u\in V\dl X_1,\ldots,X_n\dr^{\times}$ such that $g=u\cdot f$.}
\end{exer}


\section{Appendix: f-adic rings}\label{sec-fadicrings}
In \cite{Hube1} R.\ Huber\index{Huber, R.} introduced the notion of `f-adic rings' in an attempt to give a broad generalization of the notion of (classical) affinoid algebras and thus to develop a new geometry that contains the classical rigid analytic geometry as a special case.
This appendix gives a brief survey of generalities on f-adic rings and thus prepares for our later discussion on Huber's adic spaces in {\bf \ref{ch-rigid}}, \S\ref{sec-adicspaces}.

\subsection{f-adic rings}\label{sub-fadicrings}
\index{adic!adic ring@--- ring!fadic ring@f-{---} ---|(}\index{fadic ring@f-adic ring|(}
\subsubsection{Extension of adic topologies}\label{subsub-fadicringsgeneralities1}
Let $A$ be a ring, $B\subseteq A$ a subring, and $I\subseteq B$ a finitely generated ideal of $B$.
The $I$-adic filtration $\{I^n\}_{n\geq 1}$ (cf.\ \S\ref{subsub-topologicalringsmodules}) gives rise to a linear topology on $A$ regarded as a $B$-module.
\begin{prop}\label{prop-lemfadicbasicfact11}
The underlying ring structure makes the topological $B$-module $A$ with the topology defined by the filtration $\{I^n\}_{n\geq 1}$ into a topological ring if and only if 
$$
A=\bigcup_{n\geq 0}[B:I^n]
$$
holds.
Moreover, we have $\Spec A\setminus V(IA)=\Spec B\setminus V(I)$ in this situation. 
\end{prop}

\begin{proof}
Suppose $A$ is a topological ring, and take $x\in A$.
Since the selfmap $y\mapsto xy$ on $A$ is continuous, there exists $n\geq 0$ such that $I^nx\subseteq B$, that is, $x\in[B:I^n]$.
The `if' part is easy to see.
Finally, in this situation, the inclusion $B\hookrightarrow A$ has $I$-torsion cokernel, whence the equality $\Spec A\setminus V(IA)=\Spec B\setminus V(I)$.
\end{proof}

Aside from the trivial case $B=A$, in which the topology in question is nothing but the $I$-adic topology itself, there are plenty of examples of the situations as in \ref{prop-lemfadicbasicfact11}.
For example, if $B$ is $a$-torsion free for an element $a\in B$, then $A=B[\frac{1}{a}]$ with the topology defined by the $a$-adic filtration on $B$ is a topological ring.
This, needless to say, gives grounds for justifying the topology on classical affinoid algebras, already discussed in \S\ref{sub-classicalaffinoidalgebras} (cf.\ \ref{lem-tatealgtopcomp}).

\subsubsection{f-adic rings}\label{subsub-fadicringsgeneralities}
\begin{dfn}\label{dfn-fadic}{\rm 
An {\em f-adic} ring\index{adic!adic ring@--- ring!fadic ring@f-{---} ---}\index{fadic ring@f-adic ring} is a topological ring $A$ that admits an open subring $A_0\subseteq A$ such that the induced topology on $A_0$ is an adic topology\index{topology!adic topology@adic ---}\index{adic!adic topology@--- topology} by a finitely generated ideal $I_0$ of $A_0$.}
\end{dfn}

In this situation, the subring $A_0$ is called a {\em ring of definition}\index{ring of definition}, and the ideal $I_0$ is called an {\em ideal of definition} of $A$.
If, for example, $I\subseteq B\subseteq A$ are as in \ref{prop-lemfadicbasicfact11}, then $A$ is an f-adic ring if and only if $A=\bigcup_{n\geq 0}[B:I^n]$; moreover, in this situation, $B$ is a ring of definition, and $I$ is an ideal of definition of $A$.
It can be shown, moreover, that:
\begin{itemize}
\item a subring $B\subseteq A$ of an f-adic ring $A$ is a ring of definition of $A$ if and only if it is open and bounded (\cite[Prop.\ 1 (ii)]{Hube1});
\item every ring of definition of an f-adic ring has at least one finitely generated ideal of definition (\cite[Prop.\ 1 (iii)]{Hube1}).
\end{itemize}
Here, a subset $S$ of a topological ring $A$ is said to be {\em bounded} if for any neighborhood $U$ of $0$ in $A$ there exists a neighborhood $V$ of $0$ in $A$ such that $V\cdot S\subseteq U$.

In general, for a topological ring $A$
\begin{itemize}
\item an element $a\in A$ is said to be {\em power-bounded}\index{power-bounded} if the subset $\{a^n\}_{n\geq 0}$ is bounded;
\item an element $a\in A$ is said to be {\em topologically nilpotent}\index{topologically nilpotent} if for any open neighborhood $V$ of $0$ there exists $N\geq 0$ such that $a^n\in V$ whenever $n\geq N$.
\end{itemize}
We denote by $A^o$ (resp.\ $N(A)$) the subset of $A$ consisting of all power-bounded (resp.\ topologically nilpotent) elements.
The following propositions are easy:
\begin{prop}\label{prop-lempowerboundedfadic2}
Let $A$ be an f-adic ring.

{\rm (1)} Let $A_0\subseteq A$ be a ring of definition.
Then $f\in A$ is power-bounded\index{power-bounded} if and only if $A_0[f]\subseteq A$ is a ring of definition of $A$.

{\rm (2)} The set of power-bounded elements $A^o$ coincides with the union of all rings of definition of $A$. \hfill$\square$
\end{prop}

\begin{prop}\label{prop-lempowerboundedfadic}
Let $A$ be an f-adic ring, $A_0$ a ring of definition of $A$, and $I_0$ an ideal of definition of $A_0$.

{\rm (1)} The subset $A^o$ is an open subring of $A$, and $N(A)$ is an open ideal of $A^o$.

{\rm (2)} Any element of $I_0$ is topologically nilpotent. 
Conversely, for any topologically nilpotent element $a$ of $A$ there exists $n>0$ such that $a^n\in I_0$. \hfill$\square$
\end{prop}

Let $A$ be an f-adic ring, $A_0$ a ring of definition, and $I_0\subseteq A_0$ an ideal of definition.
Then every open subring $A'_0\subseteq A_0$ is a ring of definition of $A$.
Moreover, due to \ref{prop-lempowerboundedfadic} (2) there exists $n\geq 0$ such that $I^n_0\subseteq A'_0$, and hence the topology on $A'_0$ is $I^n_0$-adic.

Let $\varphi\colon A\rightarrow B$ be a continuous ring homomorphism between f-adic rings, and $A_0\subseteq A$ and $B_0\subseteq B$ respective rings of definition.
Then $A'_0=A_0\cap\varphi^{-1}(B_0)$ is a ring of definition of $A$, since it is open and bounded.
Thus, whenever we are given a continuous ring homomorphism as above, we can always take rings of definition $A_0$ and $B_0$ in such a way that $\varphi(A_0)\subseteq B_0$.

If the restriction $\varphi|_{A_0}\colon A_0\rightarrow B_0$ is an adic homomorphism ({\bf \ref{ch-pre}}, \S\ref{subsub-adicfiltrationtopology}), we say that the morphism $\varphi\colon A\rightarrow B$ is {\em adic}\index{adic!adic morphism@--- morphism}.
Notice that this notion does not depend on the choice of the rings of definition; indeed, for another choice $A'_0\subseteq A$ and $B'_0\subseteq B$ with $\varphi(A'_0)\subseteq B'_0$, taking the intersections $A_0\cap A'_0$ and $B_0\cap B'_0$, which are again rings of definitions, we may assume that $A'_0\subseteq A_0$ and $B'_0\subseteq B_0$ and thus reduce to showing that $A_0\rightarrow B_0$ is adic if and only if $A'_0\rightarrow B'_0$, which is easy to see.

\subsubsection{Extremal f-adic rings}\label{subsub-extremalfadicrings}
\index{fadic ring@f-adic ring!extremal fadic ring@extremal ---|(}
The notion of f-adic rings contains, as a special case, usual linearly topologized rings with adic topology defined by finitely generated ideals.
Also, as remarked at the end of \S\ref{subsub-fadicringsgeneralities1}, classical affinoid algebras are another example of f-adic rings. 
These special cases indicates that there are two interesting classes of f-adic rings as follows.

One of them is comprised of f-adic rings that themselves are bounded (hence `bounded f-adic rings' should be the logical name); if $A$ is such an f-adic ring, then $A$ itself is a ring of definition, and hence f-adic rings of this type are nothing but rings with adic topology defined by finitely generated ideals.

The other interesting type of f-adic rings is the following:
\begin{itemize}
\item an f-adic ring $A$ is said to be {\em extremal}\index{fadic ring@f-adic ring!extremal fadic ring@extremal ---} if it has an ideal of definition $I$ (of a ring of definition) such that $IA=A$.
\end{itemize}
(It can be shown easily that this notion does not depend on the choice of $I$.)
If $A$ is extremal f-adic and $A_0\subseteq A$ is a ring of definition with an ideal of definition $I_0\subseteq A_0$, then it follows from \ref{prop-lemfadicbasicfact11} that $\Spec A=\Spec A_0\setminus V(I_0)$.
Conversely, if $B$ is a ring endowed with an adic topology defined by a finitely generated ideal $I\subseteq B$ such that $\Spec B\setminus V(I)$ is affine, then $A$ with $\Spec A=\Spec B\setminus V(I)$ is an extremal f-adic ring by the topology defined by the filtration $\{I^n\}_{n\geq 1}$.

The f-adic rings of the type mentioned at the end of \S\ref{subsub-fadicringsgeneralities1}, and hence all classical affinoid algebras, are the further special type of extremal f-adic rings: 
\begin{itemize}
\item an extremal f-adic ring is called a {\em Tate ring}\index{Tate ring} if it has a principal ideal of definition; 
\end{itemize}
equivalently, Tate rings are f-adic rings having topologically nilpotent units (\cite[\S1, Definition (ii)]{Hube1}).

\begin{prop}\label{prop-typeRfadicrings}
An f-adic ring $A$ is extremal f-adic if and only if it has a ring of definition $A_0$ and finitely generated topologically nilpotent ideal $J_0\subseteq A_0$ such that $J_0A=A$.
In this situation, moreover, the ideal $J_0$ is open in $A_0$ and the topology on $A_0$ is $J_0$-adic .
\end{prop}

\begin{proof}
The `only if' part is trivial.
To show the converse, replacing $J_0$ by an ideal of the form $J^n_0$, we may assume that $J_0$ is contained in an ideal of definition $I_0$.
Then obviously we have $I_0A=A$.
Since $V(I_0)\subseteq V(J_0)$, we have the chain of morphisms
$$
\Spec A=\Spec A\setminus V(J_0A)\longrightarrow\Spec A_0\setminus V(J_0)\hookrightarrow\Spec A\setminus V(I_0)=\Spec A
$$
of which the composition is equal to the identity map of $\Spec A$.
Hence we deduce that $V(I_0)=V(J_0)$ and that there exists $n\geq 0$ such that $I^n_0\subseteq J_0$.
\end{proof}

\begin{cor}[cf.\ {\rm \cite[Prop.\ 1.10]{Hube1}}]\label{cor-typeRfadicrings}
Let $A$ and $B$ be f-adic rings, and $\varphi\colon A\rightarrow B$ a continuous ring homomorphism.
Suppose $A$ is an extremal f-adic $($resp.\ Tate$)$ ring\index{Tate ring}.
Then $B$ is also an extremal f-adic $($resp.\ Tate$)$ ring, and the map $\varphi$ is adic.
\end{cor}

\begin{proof}
Take rings of definition $A_0\subseteq A$ and $B_0\subseteq B$ such that $\varphi(A_0)\subseteq B_0$.
Let $I_0\subseteq A_0$ be an ideal of definition, and set $J_0=I_0B_0$.
Then $J_0$ is a finitely generated topologically nilpotent ideal such that $J_0B=B$, and hence $B_0$ is extremal f-adic due to \ref{prop-typeRfadicrings}; moreover, $J_0$ is an ideal of definition of $B_0$.
\end{proof}
\index{fadic ring@f-adic ring!extremal fadic ring@extremal ---|)}

\subsubsection{Complete f-adic rings}\label{subsub-completeaffinoidrings}
Let $A$ be an f-adic ring, $B$ a ring of definition\index{ring of definition}, and $I\subseteq B$ a finitely generated ideal of definition of $B$.
By \ref{prop-lemfadicbasicfact11} $J=\bigcap_{n\geq 1}I^n$ is an ideal of $A$, and $A$ is separated if and only if $J=0$.
In general, the associated separated f-adic ring is the one given by $A/J$, which is again an f-adic ring with an ideal of definition $B/J$.

The f-adic ring $A$ is said to be {\em complete}\index{fadic ring@f-adic ring!complete fadic ring@complete ---} if it is separated and complete.
Notice that, if $A$ is a complete f-adic ring, then any ring of definition $A_0\subseteq A$ is $I_0$-adically complete for a finitely generated ideal of definition $I_0\subseteq A_0$, since $A_0$ is an open subring of $A$.

Let $A$, $B$, and $I$ be as above.
Then $\widehat{B}$, the $I$-adic completion of $B$ (cf.\ {\bf \ref{ch-pre}}.\ref{prop-Iadiccompletioncomplete2}), can be seen as a subring of the completion $\widehat{A}$ of $A$ (cf.\ \S\ref{subsub-completionfiltration}).

\begin{prop}\label{prop-completefadicbasic1}
The canonical map
$$
\widehat{B}\otimes_BA\longrightarrow\widehat{A}
$$
is an isomorphism of rings.
\end{prop}

\begin{proof}
The inverse map is constructed as follows. 
Any element $x$ of $\widehat{A}$ is represented by a Cauchy sequence $\{x_i\}_{i=0}^{\infty}$ by elements of $A$.
We may assume that for any $i,j$ we have $x_i-x_j\in B$.
Then $\{x_i-x_0\}_{i=0}^{\infty}$ is a Cauchy sequence in $B$, and hence define a unique element $y\in\widehat{B}$.
Then the inverse mapping in question is given by $x\mapsto y+x_0\in\widehat{B}\otimes_BA$.
For more details of the proof, see \cite[1.6]{Hube1}.
\end{proof}

By \ref{prop-Iadiccompletioncomplete2} the topology on $\widehat{B}$ is the $I\widehat{B}$-adic topology, and hence the completion $\widehat{A}$ is again an f-adic ring having $\widehat{B}$ and $I\widehat{B}$ as a ring of definition and an ideal of definition, respectively.
It is clear that, if $A$ is extremal f-adic\index{fadic ring@f-adic ring!extremal fadic ring@extremal ---} (resp.\ Tate\index{Tate ring}) then so is $\widehat{A}$.

Let $A$ be a complete f-adic ring, and $n\geq 0$ an integer.
We denote by $A\dl X_1,\ldots,X_n\dr$ the restricted formal power series ring\index{restricted formal power series!restricted formal power series ring@--- ring} with coefficients in $A$, that is, the completion of $A[X_1,\ldots,X_n]$; if $A_0\subseteq A$ is a ring of definition, then one has
$$
A\dl X_1,\ldots,X_n\dr=A_0\dl X_1,\ldots,X_n\dr\otimes_{A_0}A,
$$
where $A_0\dl X_1,\ldots,X_n\dr$ is the one defined as in \S\ref{sub-powerseries}.
Notice that $A\dl X_1,\ldots,X_n\dr$ is again an f-adic ring having $A_0\dl X_1,\ldots,X_n\dr$ as a ring of definition.
The following lemma is easy to see:
\begin{lem}\label{lem-weaklygenerated}
Let $A\rightarrow B$ be an adic map between complete f-adic rings, and $f_1,\ldots,f_n\in B^o$ power-bounded elements of $B$.
Then there exists an adic $A$-algebra homomorphism $A\dl X_1,\ldots,X_n\dr\rightarrow B$ that maps each $T_i$ to $f_i$ $(i=1,\ldots,n)$. \hfill$\square$
\end{lem}

Let $A\rightarrow B$ be an adic map between complete f-adic rings, and $f_1,\ldots,f_n\in B$.
Suppose that the $A$-algebra homomorphism $A[X_1,\ldots,X_n]\rightarrow B$ that maps each $X_i$ to $f_i$ ($i=1,\ldots,n$) extends to an adic map $A\dl X_1,\ldots,X_n\dr\rightarrow B$.
In this situation, we say that the image of the last map is the subring {\em weakly generated} by $f_1,\ldots,f_n$ over $A$.
If, moreover, the induced subspace topology on the image coincides with the quotient topology induced from the topology on $A\dl X_1,\ldots,X_n\dr$, we say that the image is {\em generated} by $f_1,\ldots,f_n$ over $A$.

\subsubsection{Banach f-adic rings and classical affinoid algebras}\label{subsub-banachfadicrings}
\index{fadic ring@f-adic ring!Banach fadic ring@Banach ---|(}
Let $V$ be an $a$-adically complete valuation ring\index{valuation!valuation ring@--- ring!a-adically complete valuation ring@$a$-adically complete --- ---} of height one, where $a$ is a non-zero element of $\m_V$; in view of \ref{prop-maxspe} we have $\m_V=\sqrt{(a)}$ in this situation.
The fractional field of $V$ is denoted by $K=\Frac(V)$.
Then as in \S\ref{subsub-nonarchnorms} we have a valuation\index{valuation} 
$$
|\cdot|\colon K\longrightarrow\R_{\geq 0}
$$
such that $|a|<1$; here the valuation is written multiplicatively.

A {\rm $K$-Banach algebra}\index{algebra!Banach algebra@Banach ---!KBanach algebra@$K$-{---} ---} is a pair $(A,\|\cdot\|)$ consisting of a $K$-algebra and a non-archimedean ring norm\index{norm!non-archimedean norm@non-archimedean ---} $\|\cdot\|$ such that for any $x\in A$ and $u\in K$ we have $\|ux\|\leq  |u|\|x\|$ and that $A$ is complete with respect to the norm $\|\cdot\|$.

\begin{prop}\label{prop-banachfadicrings1}
Let $A=(A,\|\cdot\|)$ be a $K$-Banach algebra, and set
$$
A_0=\{x\in A\,|\,\|x\|\leq 1\}.
$$
Then $A_0$ is an open subring, and the induced topology on $A_0$ is $a$-adic.
In particular, $A$ is a complete Tate ring. \hfill$\square$
\end{prop}

The proof is easy and is left to the reader.
An f-adic ring $A$ thus obtained is called a {\em Banach f-adic ring}.

Now let us consider a topologically finitely generated algebra $A$ over $V$ (\ref{dfn-topfinigen}) and the related classical affinoid algebra\index{algebra!affinoid algebra@affinoid ---!classical affinoid algebra@classical --- ---} $\mathcal{A}=A[\frac{1}{a}]$ ({\bf \ref{ch-pre}}, \S\ref{subsub-classicalaffinoidalgebras}).

In case $A=V\dl T_1,\ldots,T_n\dr$ the associated classical affinoid algebra $A[\frac{1}{a}]$ is the Tate algebra $K\dl T_1,\ldots,T_n\dr$ equipped with  the Gauss norm (\S\ref{subsub-classicalaffinoidalgebras})\index{norm!Gauss norm@Gauss ---}
$$
\|\sum_{\nu_1,\ldots,\nu_n\geq 0}a_{\nu_1,\ldots,\nu_n}T^{\nu_1}_1\cdots T^{\nu_n}_n\|=\sup_{\nu_1,\ldots,\nu_n\geq 0}|a_{\nu_1,\ldots,\nu_n}|.
$$
Obviously, we have
$$
K\dl T_1,\ldots,T_n\dr^o=V\dl T_1,\ldots,T_n\dr.
$$

In general, a classical affinoid algebra $\mathcal{A}$, written $\mathcal{A}=K\dl T_1,\ldots,T_n\dr/\mathfrak{a}$, is a $K$-Banach algebra with respect to the norm induced from the Gauss norm, the so-called {\em residue norm}), defined as follows: For $f\in\mathcal{A}$ we set
$$
\|f\|_{\mathrm{res}}=\sup_{F\mapsto f}\|F\|,
$$
where $F$ varies in the set of all elements of $K\dl T_1,\ldots,T_n\dr$ that are mapped to $f$ by the quotient map.

\begin{lem}\label{lem-banachfadicringstate2}
The topology induced by the residue norm $\|\cdot\|_{\mathrm{res}}$ on $\mathcal{A}$ coincides with the $a$-adic topology.
In particular, the topology on $\mathcal{A}$ by the residue norm does not depend on the choice of the presentation of $\mathcal{A}$ as a quotient of a Tate algebra.
\end{lem}

\begin{proof}
As in \S\ref{subsub-ARIgoodquot} the $a$-adic topology on $\mathcal{A}$ is the one induced from the filtration $\{\mathfrak{a}+a^nV\dl T_1,\ldots,T_n\dr\}_{n\geq 0}$ by the subgroup consisting of modulo $\mathfrak{a}$ residue classes.
Since $\mathfrak{a}+a^nV\dl T_1,\ldots,T_n\dr=\{f\in A\,|\,\|f\|_{\mathrm{res}}\leq |a^n|\}$ for each $n\geq 0$, we have the lemma.
\end{proof}
\index{fadic ring@f-adic ring!Banach fadic ring@Banach ---|)}

\subsection{Modules over f-adic rings}\label{sub-analyticfadicrings}
\subsubsection{Topological modules}\label{subsub-topologicalmodulesfadicrings}
Let $A$ be an f-adic ring with a ring of definition $A_0\subseteq A$ and a finitely generated ideal of definition $I_0\subseteq A_0$.
Let $M$ be an $A$-module.
To topologize $M$, let $M_0\subseteq M$ be an $A_0$-submodule, and consider the $I_0$-adic filtration $\{I^n_0M_0\}_{n\geq 1}$ (cf.\ \S\ref{subsub-topologicalringsmodules}).
\begin{prop}\label{prop-topologicalmodulesfadicrings1}
The topological group $M$ endowed with the topology given by the filtration $\{I^n_0M_0\}_{n\geq 1}$ is a topological $A$-module if and only if 
$$
M=\bigcup_{n\geq 0}[M_0:I^n_0].\eqno{\square}
$$
\end{prop}

The proof is similar to that of \ref{prop-lemfadicbasicfact11}; one only has to verify that the scalar multiplication $A\times M\rightarrow M$ is continuous if and only if the above equality holds.

In general, a subset $W\subseteq M$ of a topological $A$-module $M$ is said to be {\em bounded} if for any open neighborhood $U\subseteq M$ of $0$ in $M$ there exists an open neighborhood $V\subseteq A$ of $0$ in $A$ such that $V\cdot W\subseteq U$.
\begin{prop}\label{prop-topologicalmodulesfadicrings2}
Let $M$ be a topological $A$-module, and $M_0\subseteq M$ an $A_0$-submodule.
Then the subspace topology on $M_0$ induced from that of $M$ coincides with the $I_0$-adic topology if and only if $M_0$ is open and bounded in $M$.
\end{prop}

\begin{proof}
If the subspace topology on $M_0$ is $I_0$-adic, then $M_0$ is clearly open in $M$.
Since the map $A\times M\rightarrow M$ is continuous, for any open neighborhood $U\subseteq M$ of $0$ there exists $n\geq 1$ such that $I^nM_0\subseteq U$, which shows that $M_0$ is bounded.
Conversely, if $M_0$ is open and bounded, then so is $I^n_0M_0$ for any $n\geq 1$.
Hence $\{I^n_0M_0\}_{n\geq 1}$ gives a fundamental system of open neighborhoods of $M$.
\end{proof}

\subsubsection{Open mapping theorem}\label{subsub-openmappingtheorem}
\begin{thm}[Open mapping theorem]\label{thm-openmappingtheoremmodules}
Let $A$ be a complete f-adic ring, and $A_0\subseteq A$ a ring of definition.
Let $M$ and $N$ be complete topological $A$-modules, and $M_0\subseteq M$ and $N_0\subseteq N$ open and bounded $A_0$-submodules.
Then for a continuous $A$-linear homomorphism $\varphi\colon M\rightarrow N$ to be topological isomorphism is equivalent to that $\varphi$ is bijective.
\end{thm}

In the following proof of the theorem we regard $M$ and $N$ as metric spaces in the way as in Exercise \ref{exer-Fujiwaralemma2.6} defined by the decreasing filtration 
$$
F^{(n)}=\begin{cases}I^n_0M_0&(n\geq 0)\\ [M_0:I^{-n}_0]&(n<0)\end{cases}
$$
and similarly for $N$, where $I_0\subseteq A_0$ is a finitely generated ideal of definition.
Note that $M$ and $N$ are complete as metric spaces.

\begin{proof}[Proof of Theorem {\rm \ref{thm-openmappingtheoremmodules}}]
Suppose that $\varphi$ is bijective.
Due to the continuity of $\varphi$ we may assume $\varphi(M_0)\subseteq N_0$.
By identifying $M$ with the image $\varphi(M)$, we may assume that $M=N$ as an $A$-module and that $\varphi$ is the identity map.

We need to show that $M_0$ is open in $M$ with respect to the topology defined by $N_0$.
To this end, let us regard $M$ as a topological group by the topology induced from that of $N_0$.
The resulting $M$ is a complete metric space.
Since we have
$$
M=\bigcup_{n\geq 1}[M_0:I^n_0]\subseteq\bigcup_{n\geq 1}\ovl{[M_0:I^n_0]},
$$
by Baire's category theorem (Exercise \ref{exer-Fujiwaralemma2.6} (3)), there exist integers $m,n\geq 1$ such that $I^m_0N_0\subseteq\ovl{[M_0:I^n_0]}$.
It follows from this that
$$
I^{m+n}_0N_0\subseteq M_0+I^{n+k}_0N_0
$$
for any $k\geq 0$.
For $k=m+1$ and $c=m+n$ we have
$$
I^c_0N_0\subseteq M_0+I_0\cdot I^c_0N_0,
$$
from which it follows that for any $x\in I^c_0N_0$ and $l\geq 0$ there exists an element $y_l\in M_0$ such that $x-y_l\in I^l_0\cdot I^c_0N_0$ and $y_{l+1}-y_l\in I^{l+1}_0M_0$.
Since $M_0$ is $I_0$-adically complete, the sequence $\{y_l\}_{l\geq 0}$ converges to an element $y\in M_0$.
On the other hand, the same sequence converges by $I_0N_0$-adic topology to the element $x$, equal to $y$.
Thus we have $I^c_0N_0\subseteq M_0$, and $M_0$ is open with respect to the topology defined by $N_0$, as desired.
\end{proof}

It follows from \ref{prop-topologicalmodulesfadicrings1} that for a continuous morphism $\varphi\colon A\rightarrow B$ of $f$-adic rings, $B$ can be regarded as a topological $A$-module if $\varphi$ is adic.
Hence we have:
\begin{cor}[Open mapping theorem for complete f-adic rings]\label{cor-thmopenmappingtheorem}
Let $\varphi\colon A\rightarrow B$ be an adic homomorphism between complete f-adic rings.
Then $\varphi$ is a topological isomorphism if and only if it is bijective. \hfill$\square$
\end{cor}

By \ref{cor-typeRfadicrings} we have:
\begin{cor}\label{cor-openmappingtheorem}
Let $\varphi\colon A\rightarrow B$ be a continuous homomorphism between two complete f-adic rings.
Suppose $A$ is extremal f-adic.
Then $\varphi$ is a topological isomorphism if and only if it is bijective. \hfill$\square$
\end{cor}
\index{adic!adic ring@--- ring!fadic ring@f-{---} ---|)}\index{fadic ring@f-adic ring|)}


\section{Appendix: Addendum on derived category}\label{sec-derivedcategory}
We refer to Verdier's expositions \cite[C.D.]{SGA4.5} and \cite{Verd1} as our basic references for generalities of derived categories.
In this appendix we collect some (but not all) of the facts and materials on derived categorical calculi used in our discussion.
The proofs of almost all assertions below can be found, in addition to the above-mentioned references, also in \cite{BBD} and \cite{KS}.

\subsection{Prerequisites on triangulated categories}\label{sub-triangulatedcategory}
\index{category!triangulated category@triangulated ---|(}
In the sequel we write triangulated categories as
$$
\mathscr{D}=(\mathscr{D},T,\mathscr{T}),
$$
where $\mathscr{D}$ is the underlying additive category, $T$ the automorphism of $\mathscr{D}$ (the {\em shift operator}\index{shift functor@shift functor}), and $\mathscr{T}$ is the collection of all distinguished triangles.
An additive functor $F\colon\mathscr{D}=(\mathscr{D},T,\mathscr{T})\rightarrow\mathscr{D}'=(\mathscr{D}',T',\mathscr{T}')$ between triangulated categories is said to be {\em exact}\index{functor!exact functor@exact --- (of triangulated categories)} if there exists an isomorphism of functors $p\colon F\circ T\stackrel{\sim}{\rightarrow}T'\circ F$ such that for any distinguished triangle $(X,Y,Z,u,v,w)$ of $\mathscr{D}$ 
$$
F(X)\stackrel{F(u)}{\longrightarrow}F(Y)\stackrel{F(v)}{\longrightarrow}F(Z)\stackrel{p(X)\circ F(w)}{\longrightarrow}T'(F(X))
$$
is a distinguished triangle of $\mathscr{D}'$.

Let $\mathscr{D}=(\mathscr{D},T,\mathscr{T})$ be a triangulated category, and $\mathscr{A}$ an abelian category.
A functor $\mathcal{H}\colon\mathscr{D}\rightarrow\mathscr{A}$ is called a {\em cohomology functor}\index{functor!cohomology functor@cohomology ---} if for any distinguished triangle $(X,Y,Z,u,v,w)$ of $\mathscr{D}$ the sequence
$$
\mathcal{H}(X)\stackrel{\mathcal{H}(u)}{\longrightarrow}\mathcal{H}(Y)\stackrel{\mathcal{H}(v)}{\longrightarrow}\mathcal{H}(Z)\stackrel{\mathcal{H}(w)}{\longrightarrow}\mathcal{H}(T(X))
$$
is exact.
In this situation, we have the following long exact sequence in $\mathscr{A}$:
$$
\cdots\longrightarrow\mathcal{H}^k(X)\stackrel{\mathcal{H}^k(u)}{\longrightarrow}\mathcal{H}^k(Y)\stackrel{\mathcal{H}^k(v)}{\longrightarrow}\mathcal{H}^k(Z)\stackrel{\mathcal{H}^k(w)}{\longrightarrow}\mathcal{H}^{k+1}(X)\longrightarrow\cdots,
$$
where $\mathcal{H}^k=\mathcal{H}\circ T^k$.
For example, for any object $X$ of $\mathscr{D}$ the functors
$$
\Hom_{\mathscr{D}}(X,\,\cdot\,)\colon\mathscr{D}\longrightarrow\Ab, \quad\Hom_{\mathscr{D}}(\,\cdot\,,X)\colon\mathscr{D}^{\opp}\longrightarrow\Ab
$$
are cohomology functors (\cite[Chap.\ II, (1.2.1)]{Verd1}).

For a triangulated category $\mathscr{D}$ a pair $(\mathscr{D}^{\leq 0},\mathscr{D}^{\geq 0})$ of full subcategories of $\mathscr{D}$ is called a {\em $t$-structure}\index{t-structure@$t$-structure} on $\mathscr{D}$ if the following conditions are satisfied (where we write $\mathscr{D}^{\leq n}=T^{-n}(\mathscr{D}^{\leq 0})$ and $\mathscr{D}^{\geq n}=T^{-n}(\mathscr{D}^{\geq 0})$):
\begin{itemize}
\item[{\rm (a)}] $\mathscr{D}^{\leq -1}\subseteq\mathscr{D}^{\leq 0}$ and $\mathscr{D}^{\geq 1}\subseteq\mathscr{D}^{\geq 0}$; 
\item[{\rm (b)}] if $X\in\obj(\mathscr{D}^{\leq 0})$ and $Y\in\obj(\mathscr{D}^{\geq 1})$, then $\Hom_{\mathscr{D}}(X,Y)=0$; 
\item[{\rm (c)}] for any $X\in\obj(\mathscr{D})$ there exists a distinguished triangle of the form 
$$
X_0\longrightarrow X\longrightarrow X_1\stackrel{+1}{\longrightarrow},
$$
where $X_0\in\obj(\mathscr{D}^{\leq 0})$ and $X_1\in\obj(\mathscr{D}^{\geq 1})$.
\end{itemize}

If $(\mathscr{D}^{\leq 0},\mathscr{D}^{\geq 0})$ is a $t$-structure, we set $\mathscr{A}=\mathscr{D}^{\leq 0}\cap\mathscr{D}^{\geq 0}$ and call it the {\em core}\index{core} of the $t$-structure.
The core $\mathscr{A}$ is an abelian category.

In this book. triangulated categories, such as derived categories, are almost always equipped with the `canonical' cohomology functor and the `canonical' $t$-structure (cf.\ \S\ref{subsub-derivedcategorycanonical}), and accordingly, exact functors between them preserve these structures in the following sense: Let $\mathscr{D}=(\mathscr{D},T,\mathscr{T})$ and $\mathscr{D}'=(\mathscr{D}',T',\mathscr{T}')$ be triangulated categories having respective cohomology functors $\mathcal{H}\colon\mathscr{D}\rightarrow\mathscr{A}$ and $\mathcal{H}'\colon\mathscr{D}'\rightarrow\mathscr{A}'$ and respective $t$-structures $(\mathscr{D}^{\leq 0},\mathscr{D}^{\geq 0})$ and $(\mathscr{D}^{\prime\leq 0},\mathscr{D}^{\prime\geq 0})$; let $F\colon\mathscr{D}\rightarrow\mathscr{D}'$ be an exact functor; then,
\begin{itemize}
\item $F$ is said to {\em preserve the cohomology functors} with respect to an exact functor $q\colon\mathscr{A}\rightarrow\mathscr{A}'$ of abelian categories if there exists an isomorphism $q\circ\mathcal{H}\stackrel{\sim}{\rightarrow}\mathcal{H}'\circ F$ of the functors; 
\item  $F$ is said to {\em preserve the $t$-structures} if $F$ maps $\mathscr{D}^{\leq 0}$ $($resp.\ $\mathscr{D}^{\geq 0})$ to $\mathscr{D}^{\prime\leq 0}$ $($resp.\ $\mathscr{D}^{\prime\geq 0})$.
\end{itemize}

\medskip\noindent
{\bf Convention.} {\sl In the sequel, when discussing triangulated categories endowed with the `canonical' cohomology functors and `canonical' $t$-structures that are clear by the context, all exact functors are assumed to preserve the cohomology functors and $t$-structures, unless otherwise clearly stated.}
\index{category!triangulated category@triangulated ---|)}

\subsection{The category of complexes}\label{sub-complexcategory}
\subsubsection{Definitions}\label{subsub-complexcategorydef}
Let $\mathscr{A}$ be an additive category.
\begin{itemize}
\item A {\em complex}\index{complex} with entries in $\mathscr{A}$ is a data $F^{\bullet}=\{F^n,d^n_F\}_{n\in\Z}$ consisting of objects $F^n$ of $\mathscr{A}$ and arrows $d^n_F\colon F^n\rightarrow F^{n+1}$ such that $d^{n+1}_F\circ d^n_F=0$ for $n\in\Z$.
\item A {\em morphism} of complexes\index{complex!morphism of complexes@morphism of ---es} $f^{\bullet}\colon F^{\bullet}\rightarrow G^{\bullet}$ is a collection of arrows $f^{\bullet}=\{f^n\colon F^n\rightarrow G^n\}_{n\in\Z}$ in $\mathscr{A}$ such that $d^n_G\circ f^n=f^{n+1}\circ d^n_F$ for $n\in\Z$.
\end{itemize}

We denote by $\CC(\mathscr{A})$ the additive category of all complexes in $\mathscr{A}$; $\CC(\mathscr{A})$ is an additive category.
If $\mathscr{A}$ is an abelian category, then $\CC(\mathscr{A})$ is an abelian category.

A complex $F^{\bullet}$ is said to be {\em bounded}\index{complex!bounded complex@bounded ---} (resp.\ {\em bounded below}, resp.\ {\em bounded above}) if $F^n=0$ for $|n|\gg 0$ (resp.\ $n\ll 0$, resp.\ $n\gg 0$).
The full subcategory of $\CC(\mathscr{A})$ consisting of bounded complexes (resp.\ complexes bounded below, resp.\ complexes bounded above) is denotes by $\CC^{\bd}(\mathscr{A})$ (resp.\ $\CC^+(\mathscr{A})$, resp.\ $\CC^-(\mathscr{A})$).

There exists a canonical functor $\mathscr{A}\rightarrow\CC^{\bd}(\mathscr{A})$ that maps an object $F$ of $\mathscr{A}$ to the complex $F^{\bullet}$, where $F^n=0$ unless $n=0$ and $F^0=F$.
In this way, $\mathscr{A}$ can be regarded as a full subcategory of $\CC^{\bd}(\mathscr{A})$, and thus we get the diagram of inclusions of categories:
$$
\xymatrix@R-4ex@C-3ex{
&&\CC^+(\mathscr{A})\ \ \ar@{^{(}->}[dr]\\
\mathscr{A}\ \ \ar@{^{(}->}[r]&\CC^{\bd}(\mathscr{A})\ \ \ar@{^{(}->}[ur]\ar@{^{(}->}[dr]&&\ \ \CC(\mathscr{A})\rlap{.}\\
&&\CC^-(\mathscr{A})\ \ \ar@{^{(}->}[ur]}
$$

The category $\CC^{\ast}(\mathscr{A})$ ($\ast=$``\ \ '', $\bd$, $+$, $-$) is canonically endowed with: 
\begin{itemize}
\item the {\em shift functor} $[k]\colon\CC^{\ast}(\mathscr{A})\rightarrow\CC^{\ast}(\mathscr{A})$ for $k\in\Z$.
\end{itemize}
If $\mathscr{A}$ is abelian, $\CC^{\ast}(\mathscr{A})$ has, moreover, the following structures:
\begin{itemize}
\item the {\em cohomology functor} $\mathcal{H}^0\colon\CC^{\ast}(\mathscr{A})\rightarrow\mathscr{A}$; 
\item the {\em truncations} $\tau^{\leq n}$ and $\tau^{\geq n}\colon\CC^{\ast}(\mathscr{A})\rightarrow\CC^{\ast}(\mathscr{A})$.
\end{itemize}

\subsubsection{Shifts}\label{subsub-complexcategoryshifts}
\index{shift functor@shift functor|(}
For any object $F^{\bullet}$ of $\CC(\mathscr{A})$ the shift (by $k\in\Z$) $F[k]^{\bullet}$ is the complex defined as
$$
F[k]^n=F^{n+k}\qquad\textrm{and}\qquad d^n_{F[k]}=(-1)^kd^{n+k}_F
$$
for any $n\in\Z$.
Any morphism $f^{\bullet}\colon F^{\bullet}\rightarrow G^{\bullet}$ in $\CC(\mathscr{A})$ canonically induces $f[k]^{\bullet}\colon F[k]^{\bullet}\longrightarrow G[k]^{\bullet}$ given by 
$$
f[k]^n=f^{n+k}
$$
for any $n\in\Z$.
\index{shift functor@shift functor|)}

\subsubsection{Cohomology functor}\label{subsub-complexcategorycohfunc}
\index{functor!cohomology functor@cohomology ---|(}
Suppose $\mathscr{A}$ is an abelian category.
For an object $F^{\bullet}$ of $\CC(\mathscr{A})$ we set
$$
\mathcal{H}^0(F^{\bullet})=\ker(d^0_F)/\image(d^{-1}_F).
$$
For an arrow $f^{\bullet}\colon F^{\bullet}\rightarrow G^{\bullet}$ in $\CC(\mathscr{A})$ the definition $\mathcal{H}^0(f^{\bullet})\colon\mathcal{H}^0(F^{\bullet})\rightarrow\mathcal{H}^0(G^{\bullet})$ is obvious.
For $k\in\Z$ we set $\mathcal{H}^k(F^{\bullet})=H^0(F[k]^{\bullet})$.
\index{functor!cohomology functor@cohomology ---|)}

\subsubsection{Truncations}\label{subsub-complexcategorytruncations}
\index{truncation|(}
The {\em stupid truncations}\index{truncation!stupid truncation@stupid ---} are the operators $\sigma^{\leq n}$ and $\sigma^{\geq n}$ ($n\in\Z$) defined as follows: For an object $F^{\bullet}$ in $\CC(\mathscr{A})$
\begin{equation*}
\begin{split}
\sigma^{\leq n}F^{\bullet}&=(\cdots\longrightarrow F^{n-1}\longrightarrow F^n\longrightarrow 0\longrightarrow 0\longrightarrow\cdots),\\
\sigma^{\geq n}F^{\bullet}&=(\cdots\longrightarrow 0\longrightarrow 0\longrightarrow F^n\longrightarrow F^{n+1}\longrightarrow\cdots).
\end{split}
\end{equation*}
Any arrow $f^{\bullet}\colon F^{\bullet}\rightarrow G^{\bullet}$ in $\CC(\mathscr{A})$ induces $\sigma^{\leq n}f^{\bullet}\colon \sigma^{\leq n}F^{\bullet}\rightarrow\sigma^{\leq n}G^{\bullet}$ (and similarly $\sigma^{\geq n}f^{\bullet}$) in the obvious way.

The stupid truncation is a handy device to manipulate complexes, but since they change cohomologies of the complexes, they do not conform with other structures.
More conceptually important is the following one.
Suppose $\mathscr{A}$ is an abelian category.
For an object $F^{\bullet}$ of $\CC(\mathscr{A})$
\begin{equation*}
\begin{split}
\tau^{\leq n}F^{\bullet}&=(\cdots\longrightarrow F^{n-2}\longrightarrow F^{n-1}\longrightarrow\ker(d^n_F)\longrightarrow 0\longrightarrow\cdots),\\
\tau^{\geq n}F^{\bullet}&=(\cdots\longrightarrow 0\longrightarrow\coker(d^{n-1}_F)\longrightarrow F^{n+1}\longrightarrow F^{n+2}\longrightarrow\cdots).
\end{split}
\end{equation*}
For an arrow $f^{\bullet}\colon F^{\bullet}\rightarrow G^{\bullet}$ in $\CC(\mathscr{A})$ we have $\tau^{\leq n}f^{\bullet}\colon \tau^{\leq n}F^{\bullet}\rightarrow\tau^{\leq n}G^{\bullet}$ (and similarly $\tau^{\geq n}f^{\bullet}$) given in the obvious way.
There exists obvious arrows
$$
\tau^{\leq n}F^{\bullet}\longrightarrow F^{\bullet}\qquad\textrm{and}\qquad F^{\bullet}\longrightarrow\tau^{\geq n}F^{\bullet}
$$
in $\CC(\mathscr{A})$.
(In the sequel we often write $\tau^{<n}=\tau^{\leq n-1}$ and $\tau^{>n}=\tau^{\geq n+1}$ and also the stupid truncations similarly.)

\begin{prop}\label{prop-complexcategory2}
Let $F^{\bullet}$ be an object of $\CC(\mathscr{A})$.
Then we have
$$
\mathcal{H}^k(\tau^{\leq n}F^{\bullet})=
\begin{cases}
\mathcal{H}^k(F^{\bullet})&(k\leq n),\\
0&(k>n),
\end{cases}
\qquad
\mathcal{H}^k(\tau^{\geq n}F^{\bullet})=
\begin{cases}
0&(k<n),\\
\mathcal{H}^k(F^{\bullet})&(k\geq n),
\end{cases}
$$
for $k,n\in\Z$. \hfill$\square$
\end{prop}
\index{truncation|)}

\subsection{The triangulated category $\KC(\mathscr{A})$}\label{sub-homopotypre}
\subsubsection{Homotopies}\label{subsub-homopotypre}
We say that a morphism $f\colon F^{\bullet}\rightarrow G^{\bullet}$ in $\CC(\mathscr{A})$ is said to be {\em homotopic to zero}\index{homotopic!homotopic to zero@--- to zero} if there exist a collection of arrows $\{s^n\colon F^n\rightarrow G^{n-1}\}_{n\in\Z}$ such that
$$
f^n=s^{n+1}\circ d^n_F+d^{n-1}_G\circ s^n
$$
for $n\in\Z$.
Two arrows $f^{\bullet},g^{\bullet}\colon F^{\bullet}\rightarrow G^{\bullet}$ are said to be {\em homotopic}\index{homotopic} to each other if $f-g$ is homotopic to zero.

We denote by $\Ht(F^{\bullet},G^{\bullet})$ the set of all homotopic-to-zero arrows $F^{\bullet}\rightarrow G^{\bullet}$.
This is a subgroup of $\Hom_{\CC(\mathscr{A})}(F^{\bullet},G^{\bullet})$; moreover, the composition map 
$$
\Hom_{\CC(\mathscr{A})}(F^{\bullet},G^{\bullet})\times\Hom_{\CC(\mathscr{A})}(G^{\bullet},H^{\bullet})\longrightarrow\Hom_{\CC(\mathscr{A})}(F^{\bullet},H^{\bullet})
$$
maps $\Ht(F^{\bullet},G^{\bullet})\times\Hom_{\CC(\mathscr{A})}(G^{\bullet},H^{\bullet})$ and $\Hom_{\CC(\mathscr{A})}(F^{\bullet},G^{\bullet})\times\Ht(G^{\bullet},H^{\bullet})$ to $\Ht(F^{\bullet},H^{\bullet})$.
Hence we can define the following category.
\begin{dfn}\label{dfn-homotopypre2}{\rm 
We define the category $\KC(\mathscr{A})$ as follows:
\begin{itemize}
\item The objects of $\KC(\mathscr{A})$ are the same as $\CC(\mathscr{A})$:
$$
\obj(\KC(\mathscr{A}))=\obj(\CC(\mathscr{A})); 
$$
\item For two objects $F^{\bullet}$ and $G^{\bullet}$ we set
$$
\Hom_{\KC(\mathscr{A})}(F^{\bullet},G^{\bullet})=\Hom_{\CC(\mathscr{A})}(F^{\bullet},G^{\bullet})/\Ht(F^{\bullet},G^{\bullet}).
$$
\end{itemize}}
\end{dfn}

The categories $\KC^{\ast}(\mathscr{A})$ for $\ast=+$, $-$, $\bd$ are defined similarly; they are canonically regarded as full subcategories of $\KC(\mathscr{A})$.
By definition we easily see that $\KC^{\ast}(\mathscr{A})$ ($\ast=$``\ \ '', $\bd$, $+$, $-$) is canonically an additive category.
We denote the canonical functor $\CC^{\ast}(\mathscr{A})\rightarrow\KC^{\ast}(\mathscr{A})$ by 
$$
\mathrm{h}^{\ast}\colon\CC^{\ast}(\mathscr{A})\longrightarrow\KC^{\ast}(\mathscr{A})
$$
for $\ast=$``\ \ '', $\bd$, $+$, $-$.
The canonical functor $\mathscr{A}\rightarrow\KC(\mathscr{A})$, the composition of the inclusion $\mathscr{A}\hookrightarrow\CC(\mathscr{A})$ followed by $\mathrm{h}$, is fully faithful.
Hence we have the following diagram of categories consisting of fully faithful arrows:
$$
\xymatrix@R-4ex@C-3ex{
&&\KC^+(\mathscr{A})\ \ \ar@{^{(}->}[dr]\\
\mathscr{A}\ \ \ar@{^{(}->}[r]&\KC^{\bd}(\mathscr{A})\ \ \ar@{^{(}->}[ur]\ar@{^{(}->}[dr]&&\ \ \KC(\mathscr{A})\rlap{.}\\
&&\KC^-(\mathscr{A})\ \ \ar@{^{(}->}[ur]}
$$

\begin{prop}\label{prop-homotopyore3}
Let $\ast=$``\ \ '', $\bd$, $+$, $-$.

{\rm (1)} The shift operator $[k]$ on $\CC(\mathscr{A})$ maps homotopic arrows to homotopic arrows. 
Consequently, there exists a unique self-functor $[k]$ on $\KC^{\ast}(\mathscr{A})$ such that $\mathrm{h}^{\ast}\circ [k]=[k]\circ\mathrm{h}^{\ast}$.

{\rm (2)} Suppose $\mathscr{A}$ is abelian.
If $f^{\bullet}\colon F^{\bullet}\rightarrow G^{\bullet}$ is a morphism in $\CC(\mathscr{A})$ homotopic to zero, then $\mathcal{H}^{0}(f^{\bullet})$ is a zero map.
Consequently, there exists a unique functor $\mathcal{H}^{0}\colon\KC^{\ast}(\mathscr{A})\rightarrow\mathscr{A}$ such that $\mathcal{H}^0\circ\mathrm{h}^{\ast}=\mathcal{H}^0$.

{\rm (3)} Suppose $\mathscr{A}$ is abelian.
Let $f^{\bullet}\colon F^{\bullet}\rightarrow G^{\bullet}$ be a morphism in $\CC(\mathscr{A})$, and $n$ an integer.
Then if $f^{\bullet}$ is homotopic to zero, so are $\tau^{\leq n}f^{\bullet}$ and $\tau^{\geq n}g^{\bullet}$.
Consequently, there exists uniquely self-functors $\tau^{\leq n},\tau^{\geq n}$ on $\KC^{\ast}(\mathscr{A})$ such that $\tau^{\leq n}\circ\mathrm{h}^{\ast}=\mathrm{h}^{\ast}\circ\tau^{\leq n}$ and $\tau^{\geq n}\circ\mathrm{h}^{\ast}=\mathrm{h}^{\ast}\circ\tau^{\geq n}$. \hfill$\square$
\end{prop}

\subsubsection{Mapping cones}\label{subsub-homopotypremappingcone}
\index{mapping cone|(}
For a morphism $f\colon F^{\bullet}\rightarrow G^{\bullet}$ in $\CC(\mathscr{A})$ the {\em mapping cone} of $f$, denoted by $\cone(f)^{\bullet}$, is the object of $\CC(\mathscr{A})$ defined as follows:
$$
\cone(f)^n=F[1]^n\oplus G^n,\qquad d^n_{\cone(f)}=
\begin{bmatrix}
d^n_{F[1]}&0\\
f^{n+1}&d^n_G
\end{bmatrix}.
$$
This complex admits the canonical morphisms
$$
q^{\bullet}_f\colon G^{\bullet}\longrightarrow\cone(f)^{\bullet}, \qquad
p^{\bullet}_f\colon\cone(f)^{\bullet}\longrightarrow F[1]^{\bullet}
$$
given respectively by the collection of canonical inclusions $q^{\bullet}_f=\{q^n_f\}$ and by the collection of canonical projections $p^{\bullet}_f=\{p^n_f\}$.
Note that, if $f^{\bullet}$ is an arrow in $\CC^{\ast}(\mathscr{A})$ then $\cone(f)^{\bullet}$ belongs to $\CC^{\ast}(\mathscr{A})$.

\begin{prop}\label{prop-homotopypre41}
Let $f^{\bullet},g^{\bullet}\colon F^{\bullet}\rightarrow G^{\bullet}$ be two morphisms in $\CC(\mathscr{A})$ homotopic to each other. 
Then there exists an isomorphism 
$$
r^{\bullet}\colon\cone(f)^{\bullet}\stackrel{\sim}{\longrightarrow}\cone(g)^{\bullet}
$$
in $\CC(\mathscr{A})$ such that the following diagram commutes$:$
$$
\xymatrix@C-2ex@R-4ex{
&\cone(f)^{\bullet}\ar[dd]^{r^{\bullet}}\ar[dr]^{p^{\bullet}_f}\\
G^{\bullet}\ar[ur]^{q^{\bullet}_f}\ar[dr]_{q^{\bullet}_g}&&F[1]^{\bullet}\rlap{.}\\
&\cone(g)^{\bullet}\ar[ur]_{p^{\bullet}_g}}
$$
\end{prop}

\begin{proof}
If $\{s^n\colon F^n\rightarrow G^{n-1}\}$ gives the homotopy from $f^{\bullet}$ to $g^{\bullet}$, that is, $f^n-g^n=s^{n+1}\circ d^n_F+d^{n-1}_G\circ s^n$ for $n\in\Z$, then the map $F[1]^n\oplus G^n\ni (u,v)\mapsto (u,v+s^{n+1}(u))\in F[1]^n\oplus G^n$ induces a morphism $r^{\bullet}\colon\cone(f)^{\bullet}\rightarrow\cone(g)^{\bullet}$.
Clearly, it has the inverse map and hence is an isomorphism of complexes.
It is then straightforward to check that the map $r^{\bullet}$ thus constructed have the desired property.
\end{proof}

By \ref{prop-homotopypre41} the mapping cone $\cone(f)^{\bullet}$ together with the canonical arrows $q^{\bullet}_f$ and $p^{\bullet}_f$ can be defined for an arrow $f^{\bullet}$ in $\KC(\mathscr{A})$.
Thus for any morphism $f^{\bullet}\colon F^{\bullet}\rightarrow G^{\bullet}$ in $\KC(\mathscr{A})$ we have a triangle
$$
F^{\bullet}\stackrel{f^{\bullet}}{\longrightarrow}G^{\bullet}\stackrel{q^{\bullet}_f}{\longrightarrow}\cone(f)^{\bullet}\stackrel{q^{\bullet}_f}{\longrightarrow}F[1]^{\bullet}.\leqno{(\ast)}
$$
\begin{prop}[{\cite[Chap.\ II, (1.3.2)]{Verd1}}]\label{prop-homotopypre5}
Let $\ast=$``\ \ '', $\bd$, $+$, $-$.
Let $\mathscr{T}^{\ast}$ be the family of triangles in $\KC^{\ast}(\mathscr{A})$ that are isomorphic to a triangle of the form $(\ast)$.
Then $\KC^{\ast}(\mathscr{A})=(\KC^{\ast}(\mathscr{A}),[1],\mathscr{T}^{\ast})$ is a triangulated category. \hfill$\square$
\end{prop}
\index{mapping cone|)}

\subsection{The derived category $\DC(\mathscr{A})$}\label{sub-derivedcategory}
\index{category!derived category@derived ---|(}\index{derived category|(}
\subsubsection{Definition and first properties}\label{subsub-derivedcategorydef}
Let $\mathscr{A}$ be an abelian category, and consider the category $\KC^{\ast}(\mathscr{A})$ where $\ast=$``\ \ '', $\bd$, $+$, $-$.
Recall that a complex $F^{\bullet}\in\obj(\CC(\mathscr{A}))$ is said to be {\em acyclic}\index{complex!acyclic complex@acyclic ---} if $\mathcal{H}^k(F^{\bullet})=0$ for any $k\in\Z$.
By \ref{prop-homotopyore3} (2) one can also define the acyclicity of objects in $\KC(\mathscr{A})$.
We denote by $\Ac^{\ast}(\mathscr{A})$ the full subcategory of $\KC^{\ast}(\mathscr{A})$ consisting of acyclic objects.
This is a saturated full subcategory (\cite[Chap.\ II, (2.1.5)]{Verd1}) of $\KC^{\ast}(\mathscr{A})$.
\begin{dfn}\label{dfn-derivedcategory2}{\rm 
A morphism $f^{\bullet}\colon F^{\bullet}\rightarrow G^{\bullet}$ in $\KC(\mathscr{A})$ is called a {\em quasi-isomorphism}\index{quasi-isomorphism} if there exists a distinguished triangle of the form 
$$
F^{\bullet}\stackrel{f^{\bullet}}{\longrightarrow}G^{\bullet}\longrightarrow H^{\bullet}\stackrel{+1}{\longrightarrow},
$$
where $H^{\bullet}$ is acyclic.}
\end{dfn}

In other words, a quasi-isomorphism is a morphism $f^{\bullet}\colon F^{\bullet}\rightarrow G^{\bullet}$ such that for any $k\in\Z$ the induced morphism $\mathcal{H}^k(F^{\bullet})\rightarrow\mathcal{H}^k(G^{\bullet})$ is an isomorphism in $\mathscr{A}$.
We denote by $\Qis^{\ast}(\mathscr{A})$ the set of all quasi-isomorphisms in $\KC^{\ast}(\mathscr{A})$.
This is a multiplicative system of the triangulated category $\KC^{\ast}(\mathscr{A})$ compatible with the triangulation and is the one corresponding to $\Ac^{\ast}(\mathscr{A})$ by the correspondence as in \cite[Chap.\ II, (2.1.8)]{Verd1}.

\begin{dfn}\label{dfn-derivedcategory3}{\rm 
We set
$$
\DC^{\ast}(\mathscr{A})=\KC^{\ast}(\mathscr{A})/\Ac^{\ast}(\mathscr{A})
$$
for $\ast=$``\ \ '', $\bd$, $+$, $-$ and call them the derived categories of $\mathscr{A}$.}
\end{dfn}

We denote by 
$$
Q^{\ast}\colon\KC^{\ast}(\mathscr{A})\longrightarrow\DC^{\ast}(\mathscr{A})
$$
the quotient functor.
The category $\DC^{\ast}(\mathscr{A})$ is a triangulated category by:
\begin{itemize}
\item the shift operator $[k]$, which simply comes from the shift operator $[k]$ of $\KC^{\ast}(\mathscr{A})$; 
\item the set of all distinguished triangles $\mathscr{T}^{\ast}$, that is, the image of the set of all distinguished triangles of $\KC^{\ast}(\mathscr{A})$ by the quotient functor $Q$.
\end{itemize}

\begin{prop}\label{prop-derivedcategory32}
The canonical functors form a diagram of categories
$$
\xymatrix@R-4ex@C-3ex{
&&\DC^+(\mathscr{A})\ \ \ar[dr]\\
\mathscr{A}\ar[r]&\DC^{\bd}(\mathscr{A})\ \ \ar[ur]\ar[dr]&&\ \ \DC(\mathscr{A})\\
&&\DC^-(\mathscr{A})\ \ \ar[ur]}
$$
consisting of fully faithful functors. \hfill$\square$
\end{prop}

\subsubsection{Canonical cohomology functor and canonical $t$-structure}\label{subsub-derivedcategorycanonical}
\index{functor!cohomology functor@cohomology ---!canonical cohomology functor@canonical --- ---|(}\index{t-structure@$t$-structure!canonical t-structure@canonical ---|(}
The triangulated category $\DC^{\ast}(\mathscr{A})$ possesses the following canonical structures:
\begin{itemize}
\item the {\em canonical cohomology functor} $\mathcal{H}^0\colon\DC^{\ast}(\mathscr{A})\rightarrow\mathscr{A}$; 
\item the {\em canonical $t$-structure} $(\DC^{\ast}(\mathscr{A})^{\leq 0},\DC^{\ast}(\mathscr{A})^{\geq 0})$.
\end{itemize}

\begin{prop}\label{prop-homotopypre6}
Let $\ast=$``\ \ '', $\bd$, $+$, $-$.
The functor $\mathcal{H}^0\colon\KC^{\ast}(\mathscr{A})\rightarrow\mathscr{A}$ as in {\rm \ref{prop-homotopyore3} (2)} gives rise to a cohomology functor functor 
$$
\mathcal{H}^0\colon\DC^{\ast}(\mathscr{A})\longrightarrow\mathscr{A}.\eqno{\square}
$$
\end{prop}

\begin{prop}\label{prop-homotopypre7}
Let $\ast=$``\ \ '', $\bd$, $+$, $-$.

{\rm (1)} Consider the full subcategory $\DC^{\ast}(\mathscr{A})^{\leq 0}$ $($resp.\ $\DC^{\ast}(\mathscr{A})^{\geq 0})$ of $\DC^{\ast}(\mathscr{A})$ consisting of objects $F$ such that $\mathcal{H}^k(F)=0$ for $k>0$ $($resp.\ $k<0)$.
Then $(\DC^{\ast}(\mathscr{A})^{\leq 0},\DC^{\ast}(\mathscr{A})^{\geq 0})$ gives a $t$-structure on $\DC^{\ast}(\mathscr{A})$.

{\rm (2)} Let $\tau^{\leq n},\tau^{\geq n}\colon\KC^{\ast}(\mathscr{A})\rightarrow\KC^{\ast}(\mathscr{A})$ be as in {\rm \ref{prop-homotopyore3}} $(3)$.
Then $\tau^{\leq n}$ and $\tau^{\geq n}$ maps $\Ac^{\ast}(\mathscr{A})$ to itself.
Consequently, they induce self-functors 
$$
\tau^{\leq n},\tau^{\geq n}\colon\DC^{\ast}(\mathscr{A})\rightarrow\DC^{\ast}(\mathscr{A}),
$$
respectively.

{\rm (3)} The functor $\tau^{\leq n}$ $($resp.\ $\tau^{\geq n})$ gives a right-adjoint $($resp.\ left-adjoint$)$ to the inclusion functor $\DC^{\ast}(\mathscr{A})^{\leq n}\rightarrow\DC^{\ast}(\mathscr{A})$ $($resp.\ $\DC^{\ast}(\mathscr{A})^{\geq n}\rightarrow\DC^{\ast}(\mathscr{A}))$. \hfill$\square$
\end{prop}

\begin{prop}[{cf.\ \cite[Chap.\ III, 1.2.7]{Verd1}}]\label{prop-derivedcategory71}
Let $\mathscr{A}$ and $\mathscr{B}$ be abelian categories, and $F\colon\mathscr{A}\stackrel{\sim}{\rightarrow}\mathscr{B}$ an exact functor.
Then $F$ induces canonically the commutative diagram
$$
\xymatrix{\CC^{\ast}(\mathscr{A})\ar[d]_{\mathrm{h}^{\ast}}\ar[r]^{\CC^{\ast}(F)}&\CC^{\ast}(\mathscr{B})\ar[d]^{\mathrm{h}^{\ast}}\\ \KC^{\ast}(\mathscr{A})\ar[d]_{Q^{\ast}}\ar[r]^{\KC^{\ast}(F)}&\KC^{\ast}(\mathscr{B})\ar[d]^{Q^{\ast}}\\ \DC^{\ast}(\mathscr{A})\ar[r]_{\DC^{\ast}(F)}&\DC^{\ast}(\mathscr{B})\rlap{.}}
$$
Moreover, $\DC^{\ast}(F)$ is exact and preserves the canonical cohomology functors $($with respect to $F)$ and the canonical $t$-structures $($cf.\ {\rm \S\ref{sub-triangulatedcategory})}. \hfill$\square$
\end{prop}
\index{t-structure@$t$-structure!canonical t-structure@canonical ---|)}\index{functor!cohomology functor@cohomology ---!canonical cohomology functor@canonical --- ---|)}

\subsubsection{Representation by complexes and amplitude}\label{subsub-derivedcategoryrep}
\begin{prop}\label{prop-derivedcategory4}
Let $\ast=$``\ \ '', $\bd$, $+$, $-$.
The canonical functor $\mathscr{A}\rightarrow\DC^{\ast}(\mathscr{A})$ gives an exact categorical equivalence between $\mathscr{A}$ and the core\index{core} of the $t$-structure, that is, the full subcategory of $\DC^{\ast}(\mathscr{A})$ consisting of objects $F$ such that $\mathcal{H}^k(F)=0$ unless $k=0$. \hfill$\square$
\end{prop}

\begin{dfn}\label{dfn-derivedcategory5}{\rm 
Let $\ast=$``\ \ '', $\bd$, $+$, $-$.

(1) Let $F$ be an object of $\DC^{\ast}(\mathscr{A})$.
We say that a complex $F^{\bullet}$ in $\CC^{\ast}(\mathscr{A})$ {\em represents} $F$ if $F$ and $Q^{\ast}\circ\mathrm{h}^{\ast}(F^{\bullet})$ are isomorphic in $\DC^{\ast}(\mathscr{A})$.

(2) Let $f\colon F\rightarrow G$ be an arrow in $\DC^{\ast}(\mathscr{A})$.
We say that a morphism $f^{\bullet}\colon F^{\bullet}\rightarrow G^{\bullet}$ of complexes {\em represents} $f$ if $f$ and $Q^{\ast}\circ\mathrm{h}^{\ast}(f^{\bullet})$ are isomorphic in $\DC^{\ast}(\mathscr{A})$.}
\end{dfn}

\begin{dfn}\label{dfn-derivedcategory6}{\rm 
(1) Let $F$ be an object of $\DC^{\bd}(\mathscr{A})$. Then the {\em amplitude}\index{amplitude} of $F$, denoted by $\amp(F)$, is the following number (if it exists):
$$
\amp(F)=\sup\{k\,|\,\mathcal{H}^k(F)\neq 0\}-\inf\{k\,|\,\mathcal{H}^k(F)\neq 0\}.
$$

(2) An object $F$ of $\DC(\mathscr{A})$ is said to be {\em concentrated in degree $n$} if $\mathcal{H}^k(F)=0$ unless $k=n$.}
\end{dfn}

\begin{prop}\label{prop-derivedcategory61}
Let $F$ be an object of $\DC(\mathscr{A})$.
Then $F$ is concentrated in degree $0$ if and only if it is represented by a complex contained in the image of $\mathscr{A}\hookrightarrow\CC(\mathscr{A})$. \hfill$\square$
\end{prop}
\index{derived category|)}\index{category!derived category@derived ---|)}

\subsection{Subcategories of $\DC(\mathscr{A})$}\label{sub-subcategoryderived}
Let $\mathscr{A}$ be an abelian category, and $\mathscr{B}$ an abelian full subcategory of $\mathscr{A}$.
We say $\mathscr{B}$ is {\em thick}\index{subcategory!thick subcategory@thick ---} in $\mathscr{A}$ if for any exact sequence of the form
$$
X_0\longrightarrow X_1\longrightarrow X_2\longrightarrow X_3\longrightarrow X_4
$$
in $\mathscr{A}$ with $X_0,X_1,X_3,X_4\in\obj(\mathscr{B})$, we have $X_2\in\obj(\mathscr{B})$.

Let $\mathscr{A}$ be an abelian category, and $\DC^{\ast}(\mathscr{A})$ the associated derived category where $\ast=$``\ \ '', $\bd$, $+$, $-$.
We consider the canonical cohomology functor $\mathcal{H}^0$ (\ref{prop-homotopypre6}) and the canonical $t$-structure (\ref{prop-homotopypre7}).
For an abelian full subcategory $\mathscr{B}$ of $\mathscr{A}$ we denote by
$$
\DC^{\ast}_{\mathscr{B}}(\mathscr{A})
$$
the full subcategory of $\DC^{\ast}(\mathscr{A})$ consisting of objects $F$ such that $\mathcal{H}^k(F)\in\obj(\mathscr{B})$ for any $k\in\Z$.

\begin{prop}\label{prop-subcategoryderived1}
Suppose $\mathscr{B}$ is thick in $\mathscr{A}$.
Then $\DC^{\ast}_{\mathscr{B}}(\mathscr{A})$ together with the shift operator $[k]|_{\DC^{\ast}_{\mathscr{B}}(\mathscr{A})}$ and the set of distinguished triangles of $\DC^{\ast}(\mathscr{A})$
$$
F\longrightarrow G\longrightarrow H\stackrel{+1}{\longrightarrow}
$$
such that $F,G,H\in\obj(\DC^{\ast}_{\mathscr{B}}(\mathscr{A}))$ is a triangulated category.
Moreover$:$
\begin{itemize}
\item[{\rm (1)}] the composition 
$$
\mathcal{H}^0\colon\DC^{\ast}_{\mathscr{B}}(\mathscr{A})\longhookrightarrow\DC^{\ast}(\mathscr{A})\stackrel{\mathcal{H}^0}{\longrightarrow}\mathscr{A}
$$
is a cohomology functor $($again called the {\em canonical cohomology functor}\index{functor!cohomology functor@cohomology ---!canonical cohomology functor@canonical --- ---}$);$
\item[{\rm (2)}] if we set
$$
\DC^{\ast}_{\mathscr{B}}(\mathscr{A})^{\leq 0}=\DC^{\ast}(\mathscr{A})^{\leq 0}\cap\DC^{\ast}_{\mathscr{B}}(\mathscr{A}),\quad
\DC^{\ast}_{\mathscr{B}}(\mathscr{A})^{\geq 0}=\DC^{\ast}(\mathscr{A})^{\geq 0}\cap\DC^{\ast}_{\mathscr{B}}(\mathscr{A}),
$$
then $(\DC^{\ast}_{\mathscr{B}}(\mathscr{A})^{\leq 0},\DC^{\ast}_{\mathscr{B}}(\mathscr{A})^{\geq 0})$ gives a $t$-structure on $\DC^{\ast}_{\mathscr{B}}(\mathscr{A})$ $($again called the {\em canonical $t$-structure}\index{t-structure@$t$-structure!canonical t-structure@canonical ---}$)$. \hfill$\square$
\end{itemize}
\end{prop}

In the sequel we fix an abelian category $\mathscr{A}$ and a thick abelian full subcategory $\mathscr{B}$ of $\mathscr{A}$.
By construction of the categories $\DC^{\ast}(\mathscr{A})$ and $\DC^{\ast}(\mathscr{B})$ we have a natural functor $\delta^{\ast}\colon\DC^{\ast}(\mathscr{B})\rightarrow\DC^{\ast}(\mathscr{A})$ such that the following diagram commutes:
$$
\xymatrix{\CC^{\ast}(\mathscr{B})\ar[d]_{Q^{\ast}\circ\mathrm{h}^{\ast}}\ar@{^{(}->}[r]&\CC^{\ast}(\mathscr{A})\ar[d]^{Q^{\ast}\circ\mathrm{h}^{\ast}}\\ \DC^{\ast}(\mathscr{B})\ar[r]&\DC^{\ast}(\mathscr{A})\rlap{.}}
$$
Clearly, this functor maps $\DC^{\ast}(\mathscr{B})$ to $\DC^{\ast}_{\mathscr{B}}(\mathscr{A})$.
Thus we get the functor 
$$
\delta^{\ast}\colon\DC^{\ast}(\mathscr{B})\longrightarrow\DC^{\ast}_{\mathscr{B}}(\mathscr{A}), 
$$
which is clearly exact.

\begin{prop}[{cf.\ \cite[1.7.11, 1.7.12]{KS}}]\label{prop-subcategoryderived3}
{\rm (1)} Suppose that the following condition is satisfied$:$
\begin{itemize}
\item for any monomorphism $f\colon G\rightarrow F$ of $\mathscr{A}$ such that $G\in\obj(\mathscr{B})$ there exists a morphism $g\colon F\rightarrow H$ with $H\in\obj(\mathscr{B})$ such that the composition $g\circ f\colon G\rightarrow H$ is a monomorphism.
\end{itemize}
Then the functors $\delta^+$ and $\delta^{\bd}$ are equivalences.

{\rm (2)} Suppose the following condition is satisfied$:$
\begin{itemize}
\item for any epimorphism $f\colon F\rightarrow G$ of $\mathscr{A}$ such that $G\in\obj(\mathscr{B})$ there exists a morphism $g\colon H\rightarrow F$ with $H\in\obj(\mathscr{B})$ such that the composition $f\circ g\colon H\rightarrow G$ is an epimorphism.
\end{itemize}
Then the functors $\delta^-$ and $\delta^{\bd}$ are equivalences.\hfill$\square$
\end{prop}

\begin{cor}\label{cor-subcategoryderived31}
{\rm (1)} If $\mathscr{B}$ has enough $\mathscr{A}$-injectives $($that is, for any object $F$ of $\mathscr{A}$ there exists a monomorphism $F\rightarrow G$, where $G\in\obj(\mathscr{B})$ and is an injective object in $\mathscr{A})$, then $\delta^+$ and $\delta^{\bd}$ are equivalences.

{\rm (2)} If $\mathscr{B}$ has enough $\mathscr{A}$-projectives $($that is, for any object $F$ of $\mathscr{A}$ there exists an epimorphism $G\rightarrow F$, where $G\in\obj(\mathscr{B})$ and is a projective object in $\mathscr{A})$, then $\delta^-$ and $\delta^{\bd}$ are equivalences.\hfill$\square$
\end{cor}

\begin{prop}\label{prop-subcategoryderived4}
If the functor $\delta^{\bd}$ is fully faithful, then it is an equivalence.
\end{prop}

\begin{proof}
Let $F\in\obj(\DC^{\bd}_{\mathscr{B}}(\mathscr{A}))$.
We are going to show that $F$ belongs to the essential image of $\delta^{\bd}$ by induction with respect to the amplitude $\amp(F)$.
If $\amp(F)=0$, we may assume by suitable shifts (which does not change the nature of our assertion) that $F$ is concentrated in degree $0$.
Then the claim follows from \ref{prop-derivedcategory61}.
In general, consider the distinguished triangle
$$
\tau^{\geq n+1}F[-1]\longrightarrow\tau^{\leq n}F\longrightarrow F\stackrel{+1}{\longrightarrow}
$$
induced from $\tau^{\leq n}F\rightarrow F\rightarrow\tau^{\geq n+1}F\stackrel{+1}{\rightarrow}$.
By induction we may assume that $\tau^{\geq n+1}F[-1]=\delta^{\bd}(L)$ and $\tau^{\leq n}F=\delta^{\bd}(M)$ for $L,M\in\obj(\DC^{\bd}(\mathscr{B}))$.
Since $\delta^{\bd}$ is fully faithful, there exists a unique arrow $L\rightarrow M$ in $\DC^{\bd}(\mathscr{B})$ mapped to the arrow $\tau^{\geq n+1}F[-1]\rightarrow\tau^{\leq n}F$.
Take a distinguished triangle
$$
L\longrightarrow M\longrightarrow N\stackrel{+1}{\longrightarrow}
$$
of $\mathscr{B}$.
Then since $\delta^{\bd}$ is exact, $F$ is isomorphic to $\delta^{\bd}(N)$.
\end{proof}


\setcounter{section}{0}
\renewcommand{\thesection}{\arabic{section}}
\renewcommand{\theexer}{{\bf \thechapter}.\arabic{section}.\arabic{exer}}
\chapter{Formal geometry}\label{ch-formal}
This chapter is devoted to formal geometry, the geometry of formal schemes.
As we pointed out in Introduction, it is for us essential to treat non-Noetherian formal schemes, e.g.\ finite type formal schemes over an $a$-adically complete valuation ring of arbitrary height, the importance of which stems from the requirement for the functoriality of taking fibers of finite type morphisms in rigid geometry.
Our aim of this chapter is to give sufficiently general theory of formal schemes, including GFGA theorems, which does not limit itself to the locally Noetherian situation.

The section \S\ref{sec-basicformsch} collects basic notions in formal geometry.
Our central objects in this and following sections are the so-called {\em adic formal schemes of finite ideal type}, that is, formal schemes that are locally isomorphic to the formal spectrum $\Spf A$ by an adic ring having a finitely generated ideal of definition.
In \S\ref{sec-adequateformalschemes} we will introduce the so-called {\em universally rigid-Noetherian} and {\em universally adhesive} formal schemes, corresponding respectively to topologically universally pseudo-adhesive (via Gabber's theorem) and topologically universally adhesive rings ({\bf \ref{ch-pre}}, \S\ref{sub-adhesive}), which sit in the following hierarchy of classes of formal schemes:
$$
\bigg\{\begin{minipage}{6.6em}\begin{center}{\small univ.\ adhesive formal schemes}\end{center}\end{minipage}\bigg\}\ \subseteq\bigg\{\begin{minipage}{9.2em}\begin{center}{\small univ.\ rigid-Noetherian formal schemes}\end{center}\end{minipage}\bigg\}\ \subseteq\bigg\{\begin{minipage}{9.4em}\begin{center}{\small adic formal schemes of finite ideal type}\end{center}\end{minipage}\bigg\}.
$$
For example, any finite type formal schemes over an $a$-adically complete valuation ring, including what has been called {\em admissible formal schemes} in Tate-Raynaud's classical rigid analytic geometry, are all universally adhesive.

The sections \S\ref{sec-adicallyqcoh} to \S\ref{sec-pairdifferential} are devoted largely to basic aspects of formal geometry.
Among them, what we do in \S\ref{sec-adicallyqcoh} is worth noting; we will give in this section a systematic treatment of what we call {\em adically quasi-coherent sheaves}, which seems missing in the past literature, even in \cite[($\mathbf{0}$, \S10)]{EGAInew}.
An adically quasi-coherent sheaf on an adic formal scheme $X$ with an ideal of definition $\mathscr{I}$ is an $\O_X$-module $\mathscr{F}$ that is (1) complete $\mathscr{F}\cong\varprojlim_{k\geq 0}\mathscr{F}/\mathscr{I}^{k+1}\mathscr{F}$, and (2) $\mathscr{F}_k=\mathscr{F}/\mathscr{I}^{k+1}\mathscr{F}$ for every $k\geq 0$ is a quasi-coherent sheaf on the scheme $X_k=(X,\O_X/\mathscr{I}^{k+1})$.
It will be shown (Theorem \ref{thm-adicqcoh1}) that, if $X=\Spf A$ is affine, then an $\O_X$-module $\mathscr{F}$ is adically quasi-coherent if and only if it is isomorphic to the $\O_X$-module given by the `$\Delta$-construction', that is, $M^{\Delta}$ for a complete $A$-module $M$ (cf.\ \cite[$\mathbf{I}$, (10.10.1)]{EGA}); moreover, we have $M=\Gamma(X,\mathscr{F})$ in this situation.
The functor $M\mapsto M^{\Delta}$ from the category of finitely generated $A$-modules is exact if $X=\Spf A$ is universally rigid-Noetherian (Theorem \ref{thm-adicqcohpre1}), and this property improves significantly the homological algebra of adically quasi-coherent sheaves of finite type on locally universally rigid-Noetherian formal schemes.
This section, moreover, contains the following fundamental result (Corollary \ref{cor-extension2}): {\em Any coherent $(=$ quasi-compact and quasi-separated$)$ formal scheme of finite type admits an ideal of definition of finite type.}

In \S\ref{sec-formalalgsp} we develop the theory of formal algebraic spaces, more precisely, {\em quasi-separated adic formal algebraic spaces of finite ideal type}.
It seems that, among the past literature, we have only \cite[Chap.\ 5]{Knu} as the general and systematic reference for formal algebraic spaces.
But the class of formal algebraic spaces covered in this reference consists only of separated and Noetherian formal algebraic spaces, which is evidently not general enough for our purpose. 
It is therefore necessary to implement the whole theory from the very beginning, and this is what we do in this section.
Notice that, as we have explained in Introduction, there are important reasons for us to deal with formal algebraic spaces, not only formal schemes, in our framework of rigid geometry.

The rest of the chapter (\S\ref{sec-cohomologycoherent} to \S\ref{sec-finiformal}) gives the main body of formal geometry, which consists, roughly speaking, of generalizations of the contents of \cite[$\mathbf{III}$]{EGA}, including GFGA theorems.
The results in these sections generalize the classical results in the following two ways:
\begin{itemize}
\item we replace the Noetherian hypothesis with a weaker one;
\item the theorems (finiteness theorem, GFGA theorems) are formulated and proved in terms of derived category language. 
\end{itemize}
As for the first, more precisely, the GFGA theorems will be proved in the universally adhesive situation and furthermore generalized to the universally rigid-Noetherian situation in the appendix \S\ref{sec-weakcoherency}.

\section{Formal schemes}\label{sec-basicformsch}
In this section we overview the fundamental concepts in the theory of formal schemes.
This section is, therefore, mostly a rehash of the already well-documented accounts such as \cite[$\mathbf{I}$, \S10]{EGA} and \cite[\S10]{EGAInew}.
Our main objects to discuss are so-called {\em adic formal schemes of finite ideal type}\index{formal scheme!adic formal scheme@adic ---!adic formal scheme of finite ideal type@--- --- of finite ideal type}\index{adic!adic formal scheme@--- formal scheme!adic formal scheme of finite ideal type@--- --- of finite ideal type} (\ref{dfn-formalschemesadicformalschemes}, \ref{dfn-adicformalschemesoffiniteidealtype}), that is, adic formal schemes that admit, Zariski locally, ideals of definition of finite type.

In \S\ref{sub-formalnot} we collect some basic definitions and properties of formal schemes and ideals of definition.
As  discussed in \S\ref{sub-fiberproductsformalschemes}, the category of formal schemes has fiber products\index{fiber product!fiber product of formal schemes@--- (of formal schemes)}\index{formal scheme!fiber product of formal schemes@fiber product of ---s}, which extends the notion of fiber products of schemes.

The above-mentioned references define the notion of adic morphisms only for morphisms between locally Noetherian formal schemes (cf.\ \cite[\S10.12]{EGAInew}).
In \S\ref{sub-formalnotadicmor} we extend the definition to morphisms between general adic formal schemes of finite ideal type.

In \S\ref{sub-formalcompletionschemes} we discuss formal completions\index{completion!formal completion@formal ---} of schemes. 
Then, after briefly discussing several categories of formal schemes in \S\ref{sub-categoryofformalschemes}, we finish this section by introducing locally of finite type morphisms between formal schemes in \S\ref{sub-finitypeformal}.

\subsection{Formal schemes and ideals of definition}\label{sub-formalnot}
\subsubsection{Admissible rings}\label{subsub-formalnotadmrings}
We basically refer to \cite[$\mathbf{0}_{\mathbf{I}}$, \S7]{EGA} and \cite[$\mathbf{I}$, \S10]{EGA} for most of the fundamental notions in the theory of formal schemes. 
Here we recall some of them.

Let $A$ be a ring endowed with the topology defined by a descending filtration $F^{\bullet}=\{F^{\lambda}\}_{\lambda\in\Lambda}$ by ideals (cf.\ {\bf \ref{ch-pre}}, \S\ref{subsub-topfromfil}).
\begin{dfn}\label{dfn-admissibleringidealofdefinition1}{\rm 
An ideal $I\subseteq A$ is said to be an {\em ideal of definition}\index{ideal of definition} of the topological ring $A$ if the following conditions are satisfied:
\begin{itemize}
\item[{\rm (a)}] $I$ is {\em open}$;$ that is, there exists $\lambda\in\Lambda$ such that $F^{\lambda}\subseteq I$;
\item[{\rm (b)}] $I$ is {\em topologically nilpotent}\index{topologically nilpotent}$;$ that is, for any $\mu\in\Lambda$ there exists $n\geq 0$ such that $I^n\subseteq F^{\mu}$.
\end{itemize}}
\end{dfn}

Clearly, any open ideal contained in an ideal of definition is again an ideal of definition.
Hence, if $A$ admits at least one ideal of definition, it has a fundamental system of open neighborhoods of $0$ consisting of ideals of definition, called a {\em fundamental system of ideals of definition}\index{ideal of definition!fundamental system of ideals of definition@fundamental system of ---s ---}.

\begin{lem}\label{lem-admissibleringidealofdefinition1}
Let $A$ be a ring endowed with the topology defined by a descending filtration $F^{\bullet}=\{F^{\lambda}\}_{\lambda\in\Lambda}$ by ideals, and $I\subseteq A$ an ideal.
Then the following conditions are equivalent$:$
\begin{itemize}
\item[{\rm (a)}] the topology on $A$ is $I$-adic {\rm ({\bf \ref{ch-pre}}, \S\ref{subsub-adicfiltrationtopology})}\index{topology!adic topology@adic ---}\index{adic!adic topology@--- topology}$;$
\item[{\rm (b)}] $I$ is an ideal of definition, and $I^n$ is open for any $n\geq 0;$
\item[{\rm (c)}] $I^n$ is an ideal of definition for any $n\geq 1;$
\item[{\rm (d)}] $\{I^n\}_{n\geq 0}$ is a fundamental system of open neighborhoods of $0$.
\end{itemize}
Moreover, if these conditions are fulfilled, then for any ideal of definition $J\subseteq A$ the topology on $A$ is $J$-adic.
\end{lem}

\begin{proof}
The equivalence of (a) and (d) follows from the definition of adic topologies. 
The condition (a) is equivalent to that the following two conditions are satisfied:
\begin{itemize}
\item[(i)] for any $\lambda\in\Lambda$ there exists $n\geq 0$ such that $I^n\subseteq F^{\lambda}$;
\item[(ii)] for any $n\geq 0$ there exists $\mu\in\Lambda$ such that $F^{\mu}\subseteq I^n$.
\end{itemize}
The first condition says exactly that $I$ is topologically nilpotent, and the second one says that all $I^n$ are open, thereby the equivalence of (a) and (b).
The equivalence of (b) and (c) is clear.
If the topology on $A$ is $I$-adic and $J$ is an ideal of definition, then there exists $n\geq 0$ such that $I^n\subseteq J$; by this we have $I^{nm}\subseteq J^m$, that is, $J^m$ for any $m\geq 0$ is open.
\end{proof}

\begin{dfn}\label{dfn-admissibleringsadicrings}{\rm 
Let $A$ be a ring endowed with the topology defined by a descending filtration $F^{\bullet}=\{F^{\lambda}\}_{\lambda\in\Lambda}$ by ideals.

(1) We say that the topological ring $A$ is an {\em admissible ring}\index{admissible!admissible topological ring@--- (topological) ring} if the following conditions are satisfied:
\begin{itemize}
\item[{\rm (a)}] $A$ admits an ideal of definition;
\item[{\rm (b)}] $A$ is Hausdorff complete\index{complete!Hausdorff complete@Hausdorff ---} with respect to the topology defined by the filtration $F^{\bullet}=\{F^{\lambda}\}_{\lambda\in\Lambda}$ (cf.\ {\bf \ref{ch-pre}}, \S\ref{subsub-completionfiltration}).
\end{itemize}

(2) An admissible ring $A$ is said to be an {\em adic ring}\index{adic!adic ring@--- ring} if the topology on $A$ is $I$-adic for some ideal $I\subseteq A$.}
\end{dfn}

Notice that giving an adic ring $A$ amounts to the same as giving an isomorphism class of complete pairs $(A,I)$ ({\bf \ref{ch-pre}}, \S\ref{subsub-pairspairsgen}).

\danger{In {\bf \ref{ch-pre}}, \S\ref{subsub-adicfiltrationtopology} we used the terminology `adic topology' even when the topologies in question are not necessarily Hausdorff complete\index{complete!Hausdorff complete@Hausdorff ---}.
However, {\em adic rings} are always required to be Hausdorff complete; in other words:
$$
\text{adic ring}\ =\ \text{a ring with {\em Hausdorff complete} adic topology},
$$
which the reader should always keep in mind; cf.\ {\bf \ref{ch-pre}}.\ref{rem-adicfiltrationtopology3}.}

\begin{exa}\label{exa-admissibleringformalcompletion}{\rm 
Let $A$ be a ring, and $I\subseteq A$ an ideal.
We consider the $I$-adic topology\index{topology!adic topology@adic ---}\index{adic!adic topology@--- topology} on $A$ ({\bf \ref{ch-pre}}, \S\ref{subsub-adicfiltrationtopology}).
Then the Hausdorff completion $A^{\wedge}_{I^{\bullet}}$ of $A$ with respect to the $I$-adic topology is an admissible ring, and the closure $J$ of the image of $I$ in $A^{\wedge}_{I^{\bullet}}$ (cf.\ {\bf \ref{ch-pre}}, \S\ref{subsub-completionfiltration}) gives an ideal of definition (\cite[$\mathbf{0}_{\mathbf{I}}$, (7.2.2)]{EGA}).}
\end{exa}

\begin{dfn}\label{dfn-admissibleringsmorphisms}{\rm 
Let $A$ and $B$ be admissible rings.

(1) A {\em morphism}\index{morphism of admissible rings@morphism (of admissible/adic rings)} of admissible rings $f\colon A\rightarrow B$ is a continuous ring homomorphism.

(2) Suppose $A$ and $B$ are adic rings.
Then a morphism $f\colon A\rightarrow B$ of admissible rings is said to be {\em adic}\index{adic!adic morphism@--- morphism}\index{morphism of admissible rings@morphism (of admissible/adic rings)!adic morphism of admissible rings@adic ---} if there exists an ideal of definition $I$ of $A$ such that $IB$ is an ideal of definition of $B$ (cf.\ {\bf \ref{ch-pre}}, \S\ref{subsub-adicfiltrationtopology}).}
\end{dfn}

Notice that the continuity in (1) is equivalent to that for any ideal of definition $J$ of $B$ there exists an ideal of definition $I$ of $A$ such that $IB$ is contained in $J$ (cf.\ {\bf \ref{ch-pre}}.\ref{prop-topologyfromfiltrationproperty2} (1)).
It is easy to see that in the situation as in (2) a morphism $f\colon A\rightarrow B$ is adic if and only if $IB$ is an ideal of definition of $B$ for {\em any} ideal of definition $I$ of $A$.

\begin{dfn}\label{dfn-admissibleringoffiniteidealtype}{\rm 
An admissible ring $A$ is said to be {\em of finite ideal type}\index{admissible!admissible topological ring@--- (topological) ring!admissible topological ring of finite ideal type@--- --- of finite ideal type}\index{adic!adic ring@--- ring!adic ring of finite ideal type@--- --- of finite ideal type} if it has a fundamental system of ideals of definition consisting of finitely generated ideals.}
\end{dfn}

If the topology on $A$ is adic, then in view of \ref{lem-admissibleringidealofdefinition1} the condition is equivalent to that $A$ has at least one finitely generated ideal of definition.
\begin{exa}\label{exa-admissibleringformalcompletion2}{\rm 
Suppose in the situation as in \ref{exa-admissibleringformalcompletion} that the ideal $I$ is finitely generated.
Then $\widehat{A}=A^{\wedge}_{I^{\bullet}}$ is the $I$-adic completion\index{completion!I-adic completion@$I$-adic ---} of $A$ ({\bf \ref{ch-pre}}.\ref{prop-Iadiccompletioncomplete2}) and hence is an adic ring of finite ideal type. 
Notice that in this situation the closure $J$ of the image of $I$ in $\widehat{A}$ coincides with $I\widehat{A}$.}
\end{exa}

Let $A$ be an admissible ring, $\{F^{\lambda}\}_{\lambda\in\Lambda}$ a fundamental system of ideals of definition, and $S\subseteq A$ a multiplicative subset.
Consider the ring of fractions $A_S$ endowed with the topology defined by $\{F^{\lambda}A_S\}_{\lambda\in\Lambda}$.
Let $A_{\{S\}}$ denote the Hausdorff completion\index{completion!Hausdorff completion@Hausdorff ---}: 
$$
A_{\{S\}}=\varprojlim_{\lambda\in\Lambda}A_S/F^{\lambda}A_S.
$$
It has the induced filtration $\{\widehat{F}^{\lambda}_S\}_{\lambda\in\Lambda}$ defined as in {\bf \ref{ch-pre}}, \S\ref{subsub-completionfiltration}; each $\widehat{F}^{\lambda}_S$ is the closure of the image of $F^{\lambda}A_S$ by the canonical map $A_S\rightarrow A_{\{S\}}$ ({\bf \ref{ch-pre}}.\ref{prop-topologyfromfiltration31}).

\begin{prop}[{\cite[$\mathbf{0}_{\mathbf{I}}$, (7.6.11)]{EGA}}]\label{prop-completelocalizationadmissble}
{\rm (1)} The topological ring $A_{\{S\}}$ is an admissible ring, and $\{\widehat{F}^{\lambda}_S\}_{\lambda\in\Lambda}$ gives a fundamental system of ideals of definition.

{\rm (2)} Suppose $A$ is adic of finite ideal type, and let $I\subseteq A$ be a finitely generated ideal of definition.
Then $A_{\{S\}}$ is again adic of finite ideal type, and $IA_{\{S\}}$ is an ideal of definition. \hfill$\square$
\end{prop}

Notice that (2) follows from (1) and {\bf \ref{ch-pre}}.\ref{prop-criterionadicness1}; indeed, if one defines $J^{(n)}\subseteq A_{\{S\}}$ ($n\geq 1$) to be the closure of the image of $I^nA_S$ in $A_{\{S\}}$, then {\bf \ref{ch-pre}}.\ref{prop-criterionadicness1} implies that $J^{(n)}=I^nA_{\{S\}}$ for each $n\geq 0$.

The admissible rings $A_{\{S\}}$ in the case $S=\{f^n\,|\,n\geq 0\}$, denoted simply by $A_{\{f\}}$, will be frequently used.
If $A$ is an adic ring\index{adic!adic ring@--- ring} of finite ideal type, then the ring $A_{\{f\}}$ allows a more explicit description as follows.
Since $A_f=A[f^{-1}]$, we have a morphism
$$
A\dl T\dr\longrightarrow A_{\{f\}}
$$
that sends $X$ to $f^{-1}$, where $A\dl T\dr$ denotes the restricted formal power series ring\index{restricted formal power series!restricted formal power series ring@--- ring} ({\bf \ref{ch-pre}}, \S\ref{sub-powerseries}) with the variable $T$.

\begin{lem}\label{lem-formalnot0}
Let $A$ be an adic ring of finite ideal type. Then the above-defined morphism induces an isomorphism
$$
A\dl T\dr/(fT-1)\stackrel{\sim}{\longrightarrow}A_{\{f\}}
$$
of adit rings. 
\end{lem}

\begin{proof}
See \cite{SBosch}, 7.1, Remark 10. 
\end{proof}

\subsubsection{Formal spectrum}\label{subsub-formalnotformalspec}
\index{formal spectrum|(}
Let $A$ be an admissible ring.
The {\em formal spectrum}\index{formal spectrum} $\Spf A$ is a topologically locally ringed space with the underlying set consisting of all open prime ideals of $A$.
Notice that a prime ideal $\mathfrak{p}\subseteq A$ is open if and only if it contains at least one (hence all) ideals of definition (\cite[$\mathbf{0}_{\mathbf{I}}$, (7.1.5)]{EGA}).
Hence, $\Spf A$ as a point-set is nothing but the closed subset $V(I)$ of $\Spec A$ defined by an ideal of definition $I$ and, in fact, the topology (the Zariski topology) of $\Spf A$ is the subspace topology induced from that of $\Spec A$.
Moreover, $X=\Spf A$ is endowed with the sheaf of topological rings (considered with the pseudo-discrete topology (\cite[$\mathbf{0}_{\mathbf{I}}$, (3.8.1)]{EGA}))
$$
\O_X=\varprojlim_{I}\til{A/I}|_X,
$$
where $I$ runs through all ideals of definition of $A$ (with the reversed inclusion order). 
Here the projective limit is filtered, and any fundamental system of ideals of definition is cofinal in the set of all ideals of definition.
In particular, if $A$ is adic, the collection $\{I^{k+1}\}_{k\geq 0}$ by an ideal of definition $I$ is cofinal, and thus the above sheaf coincides with the projective limit of $\til{A_k}|_X$ $(k\geq 0)$, where $A_k=A/I^{k+1}$ for $k\geq 0$.

\begin{dfn}\label{dfn-affineformalschemedefinition}{\rm 
(1) A topologically locally ringed space isomorphic to $\Spf A$ for an admissible ring $A$ is called an {\em affine formal scheme}\index{formal scheme!affine formal scheme@affine ---}\index{affine!affine formal scheme@--- formal scheme}. 

(2) A {\em morphism}\index{morphism of formal schemes@morphism (of formal schemes)} $f\colon X\rightarrow Y$ between two affine formal schemes is a morphism of topologically locally ringed spaces.}
\end{dfn}

Like in the case of schemes, the functor 
$$
A\longmapsto\Spf A
$$
gives rise to a categorical equivalence between the opposite category of the category of admissible rings and the category of affine formal schemes (\cite[$\mathbf{I}$, \S10.2]{EGA}); in particular, from $X=\Spf A$ we recover $A=\Gamma(X,\O_X)$ as a topological ring (\cite[$\mathbf{I}$, (10.1.3)]{EGA}).

Let $A$ be an admissible ring, and consider the formal spectrum $X=\Spf A$.
As $X$ is a subspace of the topological space $\Spec A$, the subsets of the form $\mathfrak{D}(f)=D(f)\cap X$ for $f\in A$ give an open basis of $X$; moreover, since $X$ is closed in $\Spec A$, such an open subset is quasi-compact.
The space $\mathfrak{D}(f)$, equipped with the topologically locally ringed structure as an open subspace of $X$ (cf.\ {\bf \ref{ch-pre}}, \S\ref{subsub-ringedsplocalringedsp}), is an affine formal scheme isomorphic to $\Spf A_{\{f\}}$.

\begin{lem}\label{lem-formalnot1}
The following conditions for finitely many elements $f_1,\ldots,f_r\in A$ are equivalent to each other$:$
\begin{itemize}
\item[{\rm (a)}] the open sets $\mathfrak{D}(f_i)$ $(i=1,\ldots,r)$ cover $X;$
\item[{\rm (b)}] for any ideal of definition $I$ of $A$ the open sets $D(\ovl{f}_i)$ where $\ovl{f}_i=(f_i\ \mathrm{mod}\ I)$ $(i=1,\ldots,r)$ cover $\Spec A/I;$
\item[{\rm (c)}] the open sets $D(f_i)$ $(i=1,\ldots,r)$ cover $\Spec A;$
\item[{\rm (d)}] the ideal generated by $f_1,\ldots,f_r$ coincides with $A$.
\end{itemize}
\end{lem}

\begin{proof}
The following implications are clear: (a) $\Leftrightarrow$ (b) $\Leftarrow$ (c) $\Leftrightarrow$ (d); the equivalence of (a) and (b) is due to the fact that $V(I)=\Spec A/I$ and $\Spf A$ are homeomorphic to each other.
If (b) holds, then the ideal $(f_1,\ldots,f_r)$ contains an element in $1+I$. 
Then we deduce (d) by the fact $1+I\subset A^{\times}$ (that is, $A$ is $I$-adically Zariskian; cf.\ {\bf \ref{ch-pre}}.\ref{lem-Iadiccompletioncomplete21}).
\end{proof}

\begin{prop}\label{prop-lemformalnot111}
An open subset $U$ of $X=\Spf A$ is quasi-compact if and only if $U$ is of the form $U=X\setminus V(\mathfrak{a})$ $($in $\Spec A)$ by a finitely generated ideal $\mathfrak{a}\subseteq A$.
\end{prop}

\begin{proof}
The `if' part follows from \cite[(1.1.4)]{EGAInew} and the fact that $X$ is a closed subset of $\Spec A$.
If $U$ is quasi-compact, then $U=\bigcup^r_{i=1}\mathfrak{D}(f_i)$ for $f_1,\ldots,f_r\in A$.
Set $V=\bigcup^r_{i=1}D(f_i)$.
Then $V$ is a quasi-compact open set of $\Spec A$ and is equal to $\Spec A\setminus V(\mathfrak{a})$, where $\mathfrak{a}=(f_1,\ldots,f_r)$.
Hence $U=V\cap X=X\setminus V(\mathfrak{a})$.
\end{proof}
\index{formal spectrum|)}

\subsubsection{Formal schemes}\label{subsub-formalnotformalsch}
\begin{dfn}\label{dfn-formalschemesformalschemes}{\rm 
{\rm (1)} A {\em formal scheme}\index{formal scheme} is a topologically locally ringed space that is locally isomorphic to an affine formal scheme\index{affine!affine formal scheme@--- formal scheme}.

{\rm (2)} Let $X$ and $Y$ be formal schemes. A {\em morphism}\index{morphism of formal schemes@morphism (of formal schemes)} $f\colon X\rightarrow Y$ of formal schemes is a morphism of topologically locally ringed spaces.}
\end{dfn}

An {\em open formal subscheme}\index{formal subscheme!open formal subscheme@open ---} of a formal scheme $X$ is a formal scheme of the form $(U,\O_X|_U)$ by an open subset $U$ of the underlying topological space of $X$.
An {\em open immersion}\index{immersion!open immersion of formal schemes@open --- (of formal schemes)} of formal schemes is defined in a similar way as in the case of ringed spaces; cf.\ {\bf \ref{ch-pre}}, \S\ref{subsub-ringedsplocalringedsp}.
An open formal subscheme $U\subseteq X$ is said to be {\em affine}\index{affine!affine open subspace@--- open subset (subspace)} if it is an affine formal scheme.
Thus any formal scheme $X$ allows an open covering $X=\bigcup_{\alpha\in L}U_{\alpha}$ consisting of affine open formal subschemes; an open covering of this form is called an {\em affine $($open$)$ covering}\index{affine!affine open covering@--- open covering}.

\begin{dfn}\label{dfn-formalschemesadicformalschemes}{\rm 
A formal scheme $X$ is said to be {\em adic}\index{formal scheme!adic formal scheme@adic ---}\index{adic!adic formal scheme@--- formal scheme} if it allows an affine open covering $X=\bigcup_{\alpha\in L}U_{\alpha}$ such that each $U_{\alpha}$ is isomorphic to $\Spf A_{\alpha}$ by an adic ring $A_{\alpha}$ (\ref{dfn-admissibleringsadicrings} (2)).}
\end{dfn}

\begin{rem}\label{rem-clear}{\rm 
Any scheme\index{scheme} can be regarded as an adic formal scheme. 
Indeed, any ring is an adic ring by the ideal $(0)$ (hence a $0$-adic ring).
Hence, schemes are naturally regarded as $0$-adic formal schemes.
Thus the category of all formal schemes contain as a full subcategory the category of schemes.
Most importantly, the category of all formal schemes has a final object $\Spec\Z$.}
\end{rem}

\begin{dfn}\label{dfn-adicformalschemesoffiniteidealtype}{\rm 
An adic formal scheme $X$ is said to be {\em of finite ideal type}\index{formal scheme!adic formal scheme@adic ---!adic formal scheme of finite ideal type@--- --- of finite ideal type}\index{adic!adic formal scheme@--- formal scheme!adic formal scheme of finite ideal type@--- --- of finite ideal type} if there exists an affine open covering $X=\bigcup_{\alpha\in L}U_{\alpha}$ such that each $U_{\alpha}$ is isomorphic to $\Spf A_{\alpha}$ by an adic ring $A_{\alpha}$ of finite ideal type (\ref{dfn-admissibleringoffiniteidealtype}).}
\end{dfn}

The following proposition follows from \ref{prop-completelocalizationadmissble} (2):
\begin{prop}\label{prop-adicformalschemesoffiniteidealtype}
Let $X$ be an adic formal scheme of finite ideal type.
Then there exists open basis of the topology on $X$ consisting of affine open subschemes of the form $\Spf A$ by adic rings $A$ of finite ideal type. \hfill$\square$
\end{prop}

We will see later in \ref{cor-affineadicformalschemebyadicring} that, if $X=\Spf A$ is an affine adic formal scheme of finite ideal type, then $A$ itself is an adic ring of finite ideal type.

\subsubsection{Ideals of definition}\label{subsub-formalschemesidealsofdefinition}
Let $A$ be an admissible ring.
For any open ideal $J\subseteq A$ one defines the sheaf $J^{\Delta}$ on $X=\Spf A$ by 
$$
J^{\Delta}=\varprojlim_{I\subseteq J}\til{J/I}
$$
(cf.\ {\bf \ref{ch-pre}}, \S\ref{subsub-schemesbasics}), where $I$ runs through all ideals of definition contained in $J$.
Since $\til{J/I}$ is an ideal of $\til{A/I}=\O_{\Spec A/I}$, $J^{\Delta}$ is an ideal sheaf of $\O_X$ (cf.\ {\bf \ref{ch-pre}}.\ref{prop-projlimleftexact}).
Notice that, in view of \cite[$\mathbf{I}$, (1.3.7)]{EGA} and the definition of projective limit sheaves (cf.\ {\bf \ref{ch-pre}}, \S\ref{subsub-projlimsheaves}), we have
$$
\Gamma(X,J^{\Delta})=\varprojlim_{I\subseteq J}\Gamma(X,\til{J/I})=\varprojlim_{I\subseteq J}J/I=J,
$$
where the last equality follows from the fact that $J$ is open (and hence is closed) (cf. Exercise \ref{exer-filtrationcompletion}).

\begin{dfn}\label{dfn-idealofdefinitionaffineformalschemes}{\rm 
(1) Let $A$ be an admissible ring, and $X=\Spf A$.
An ideal sheaf $\mathscr{I}$ of $\O_X$ is said to be an {\em ideal of definition}\index{ideal of definition} if for any $x\in X$ there exists an open neighborhood $U$ of $x$ of the form $U=\mathfrak{D}(f)$ $(f\in A)$ such that $\mathscr{I}|_U=J^{\Delta}$ for some ideal of definition $J$ of the admissible ring $A_{\{f\}}$.

(2) Let $X$ be a formal scheme, and $\mathscr{I}$ an ideal sheaf of $\O_X$.
We say that $\mathscr{I}$ is an {\em ideal of definition}\index{ideal of definition} of $X$ if for any $x\in X$ there exists an affine open neighborhood $U\cong\Spf A$ of $x$ such that $\mathscr{I}|_U$ is isomorphic to an ideal of definition of $\Spf A$ as in (1).
An ideal of definition $\mathscr{I}$ is said to be {\em of finite type}\index{ideal of definition!ideal of definition of finite type@--- of finite type} if it is of finite type as an $\O_X$-module.}
\end{dfn}

\begin{prop}[{\cite[$\mathbf{I}$, (10.3.5)]{EGA}}]\label{prop-grothendieckEGAnew10.3.5}
Any ideal of definition\index{ideal of definition} of $\Spf A$ is of the form $I^{\Delta}$ for a uniquely determined ideal of definition $I$ of $A$. \hfill$\square$
\end{prop}

Thus, if $\mathscr{I}$ is an ideal of definition of a formal scheme $X$ and $X=\bigcup_{\alpha\in L}U_{\alpha}$ is an affine open covering with $U_{\alpha}\cong\Spf A_{\alpha}$ for each $\alpha\in L$, then $\mathscr{I}|_{U_{\alpha}}\cong I^{\Delta}_{\alpha}$ for $\alpha\in L$, where $I_{\alpha}\subseteq A_{\alpha}$ is an ideal of definition of the admissible ring $A_{\alpha}$.
In particular, the locally ringed space $(X,\O_X/\mathscr{I})$ is a scheme, which has $\{(U_{\alpha},\O_{U_{\alpha}}/\mathscr{I}|_{U_{\alpha}})=\Spec A/I_{\alpha}\}_{\alpha\in L}$ as an affine open covering.

In this situation, a collection $\{\mathscr{I}^{(\lambda)}\}_{\lambda\in\Lambda}$ of ideals of definition of $X$ indexed by a directed set is called a {\em fundamental system of ideals of definition}\index{ideal of definition!fundamental system of ideals of definition@fundamental system of ---s ---} of $X$ if for each $\alpha\in L$ the collection of ideals $\{I^{(\lambda)}_{\alpha}\}_{\lambda\in\Lambda}$ given by $\mathscr{I}^{(\lambda)}|_{U_{\alpha}}\cong (I^{(\lambda)}_{\alpha})^{\Delta}$ for each $\lambda\in\Lambda$ is a fundamental system of ideals of definition of $A_{\alpha}$ (cf.\ \cite[$\mathbf{I}$, (10.3.7), (10.5.1)]{EGA}).

\begin{prop}\label{prop-adicformalschemeidealofdefinitionfinitetype4}
Let $A$ be an adic ring of finite ideal type {\rm (\ref{dfn-admissibleringoffiniteidealtype})}, and $I\subseteq A$ a finitely generated ideal of definition.
Then we have
$$
I^{\Delta}=I\O_X.
$$
\end{prop}

\begin{proof}
For any $f\in A$ the ring $A_{\{f\}}$ coincides with the $IA_f$-adic completion of $A_f$ ({\bf \ref{ch-pre}}.\ref{prop-Iadiccompletioncomplete2}).
Then by {\bf \ref{ch-pre}}.\ref{cor-Iadiccompletionifexists1} we see that the closure of the image of $IA_f$ in $A_{\{f\}}$ coincides with $IA_{\{f\}}$.
The assertion follows from this.
\end{proof}

\begin{cor}\label{cor-adicformalschemeidealofdefinitionfinitetype11}
In the situation as in {\rm \ref{prop-adicformalschemeidealofdefinitionfinitetype4}}, the associated ideal of definition $I^{\Delta}$ of $X=\Spf A$ is an ideal sheaf of $\O_X$ of finite type and, moreover, we have 
$$
(I^{\Delta})^n=(I^n)^{\Delta}
$$
for any $n\geq 1$ $($cf.\ {\rm \cite[$\mathbf{I}$, (10.3.6)]{EGA}}$)$.
In particular, $\{\mathscr{I}^n\}_{n\geq 1}$ with $\mathscr{I}=I^{\Delta}$ gives a fundamental system of ideals of definition of $X=\Spf A$. \hfill$\square$
\end{cor}

\begin{prop}\label{prop-adicformalschemeidealofdefinitionfinitetype22}
Let $A$ be an adic ring of finite ideal type, and $\mathscr{I}$ an ideal of definition of finite type on $X=\Spf A$.
Then there exists a finitely generated ideal of definition $I\subseteq A$ such that $\mathscr{I}=I^{\Delta}$.
\end{prop}

The following lemma is useful not only for proving this proposition, but also for our later arguments:
\begin{lem}\label{lem-vanishingcohomologyadicallyuseful}
Let $X$ be a formal scheme, and $\{\mathscr{I}^{(\lambda)}\}_{\lambda\in\Lambda}$ a fundamental system of ideals of definition\index{ideal of definition!fundamental system of ideals of definition@fundamental system of ---s ---}, and suppose that the directed set $\Lambda$ has a final and at most countable subset.
Let $\{\mathscr{F}_{\lambda}\}_{\lambda\in\Lambda}$ be a projective system\index{system!projective system@projective ---} of $\O_X$-modules such that the following conditions are satisfied$:$
\begin{itemize}
\item[{\rm (a)}] for any $\lambda\in\Lambda$ we have $\mathscr{I}^{(\lambda)}\mathscr{F}_{\lambda}=0$, and $\mathscr{F}_{\lambda}$ is a quasi-coherent sheaf on the scheme $X_{\lambda}=(X,\O_X/\mathscr{I}^{(\lambda)});$
\item[{\rm (b)}] the projective system $\{\mathscr{F}_{\lambda}\}_{\lambda\in\Lambda}$ is strict\index{strict projective system@strict (projective system)}\index{system!projective system@projective ---!strict projective system@strict --- ---}, that is, all transitions maps $\mathscr{F}_{\mu}\rightarrow\mathscr{F}_{\lambda}$ for $\lambda\leq\mu$ are surjective.
\end{itemize}

{\rm (1)} For $q\geq 1$ we have
$$
{\textstyle \varprojlim^{(q)}_{\lambda\in\Lambda}\mathscr{F}_{\lambda}=0}.
$$

{\rm (2)} If, moreover, $X$ is affine, then for $q\geq 1$ we have 
$$
\H^q(X,\mathscr{F})=0.
$$
\end{lem}

\begin{proof}
In view of {\bf \ref{ch-pre}}.\ref{rem-coheffaceable} (2) the projective system $\{\mathscr{F}_{\lambda}\}_{\lambda\in\Lambda}$ satisfies the conditions {\bf (E1)} and {\bf (E2)} in {\bf \ref{ch-pre}}, \S\ref{subsub-projectivelimitandcohomology}.
Hence (1) follows from {\bf \ref{ch-pre}}.\ref{prop-ML5ML5} (1).
To show (2), we first notice that for $\lambda\leq\mu$ the surjection $\mathscr{F}_{\mu}\rightarrow\mathscr{F}_{\lambda}$ of quasi-coherent sheaves on the affine scheme $X_{\mu}=\Spec A/I^{(\mu)}$ (where $(I^{(\mu)})^{\Delta}=\mathscr{I}^{(\mu)}$) induces the surjective map $\Gamma(X,\mathscr{F}_{\mu})\rightarrow\Gamma(X,\mathscr{F}_{\lambda})$; that is, the projective system $\{\Gamma(X,\mathscr{F}_{\lambda})\}_{\lambda\in\Lambda}$ is strict.
Hence the assertion follows from {\bf \ref{ch-pre}}.\ref{cor-ML5} and {\bf \ref{ch-pre}}.\ref{thm-vanishcohaff-0} (1).
\end{proof}

\begin{proof}[Proof of Proposition {\rm \ref{prop-adicformalschemeidealofdefinitionfinitetype22}}]
Take the unique $I\subseteq A$ such that $I^{\Delta}=\mathscr{I}$ (\ref{prop-grothendieckEGAnew10.3.5}).
We want to show that $I$ is finitely generated.
Take a finitely generated ideal of definition $J\subseteq A$ such that $J\subseteq I$, and set $\mathscr{J}=J^{\Delta}$.
Consider the exact sequence
$$
0\longrightarrow\til{J/J^n}\longrightarrow\til{I/J^n}\longrightarrow\til{I/J}\longrightarrow 0
$$
for any $n>0$.
By \ref{lem-vanishingcohomologyadicallyuseful} (1) we deduce that $\til{I/J}\cong I^{\Delta}/J^{\Delta}$, which is a quasi-coherent sheaf on the scheme $\Spec A/J$ of finite type.
Hence $I/J$ is finitely generated, and thus $I$ is finitely generated, as desired.
\end{proof}

\begin{cor}\label{cor-adicformalschemeidealofdefinitionfinitetype22}
Let $X$ be an adic formal scheme of finite ideal type, and $\mathscr{I}$ an ideal of definition of finite type.
Then $\{\mathscr{I}^n\}_{n\geq 1}$ gives a fundamental system of ideals of definition of $X$.
\end{cor}

\begin{proof}
Obviously, we may assume that $X$ is affine of the form $X=\Spf A$, where $A$ is an adic ring of finite ideal type.
Then by \ref{prop-adicformalschemeidealofdefinitionfinitetype22} $\mathscr{I}$ is of the form $I^{\Delta}$ by a finitely generated ideal of definition $I\subseteq A$.
The assertion in this case has already been shown in \ref{cor-adicformalschemeidealofdefinitionfinitetype11}.
\end{proof}

\subsubsection{Noetherian formal schemes}\label{subsub-noetherianformalschemes}
\begin{dfn}[cf.\ {\cite[$\mathbf{I}$, (10.4.2)]{EGA}}]\label{dfn-northeroanformalschemes}{\rm 
A formal scheme $X$ is said to be {\em locally Noetherian}\index{formal scheme!Noetherian formal scheme@Noetherian ---!locally Noetherian formal scheme@locally --- ---} if it has an affine open covering $X=\bigcup_{\alpha\in L}U_{\alpha}$ such that each $U_{\alpha}$ is isomorphic to $\Spf A_{\alpha}$ by a Noetherian adic ring.
A locally Noetherian formal scheme is said to be {\em Noetherian}\index{formal scheme!Noetherian formal scheme@Noetherian ---} if the underlying topological space is quasi-compact\index{space@space (topological)!quasi-compact topological space@quasi-compact ---}\index{quasi-compact!quasi-compact topological space@--- (topological) space}.}
\end{dfn}

Thus any locally Noetherian formal scheme is, by definition, an {\em adic} formal scheme {\em of finite ideal type}.
One of the most remarkable properties of locally Noetherian formal schemes is that they always have a global ideal of definition:
\begin{prop}[{\cite[$\mathbf{I}$, (10.5.4)]{EGA}}]\label{prop-northeroanformalschemes}
Any locally Noetherian formal scheme $X$ has a unique ideal of definition $\mathscr{I}$ of finite type such that the induced scheme $(X,\O_X/\mathscr{I})$ is reduced.\hfill$\square$
\end{prop}

\subsection{Fiber products}\label{sub-fiberproductsformalschemes}
\subsubsection{Complete tensor product of admissible rings}\label{subsub-completetensorproducts}
Consider the diagram of admissible rings 
$$
B\stackrel{f}{\longleftarrow}A\stackrel{g}{\longrightarrow}C, 
$$
and let $\{I^{(\lambda)}\}_{\lambda\in\Lambda}$ (resp.\ $\{J^{(\alpha)}\}_{\alpha\in\Sigma}$, resp.\ $\{K^{(\beta)}\}_{\beta\in T}$) be a fundamental system of ideals of definition of $A$ (resp.\ $B$, resp.\ $C$).
Since $f$ and $g$ are continuous, for any $\alpha\in\Sigma$ and $\beta\in T$ there exists $\lambda\in\Lambda$ such that $I^{(\lambda)}B\subseteq J^{(\alpha)}$ and $I^{(\lambda)}C\subseteq K^{(\beta)}$.
We consider the complete tensor product\index{complete!complete tensor product@--- tensor product} $B\widehat{\otimes}_AC$ sitting in the diagram
$$
\xymatrix@-1ex{B\ar[r]&B\widehat{\otimes}_AC\\ A\ar[u]^f\ar[r]_g&C\rlap{;}\ar[u]}
$$
the ring $B\widehat{\otimes}_AC$ is the Hausdorff completion\index{completion!Hausdorff completion@Hausdorff ---} of the tensor product $B\otimes_AC$ with respect to the topology defined by the filtration $\{H^{\alpha,\beta}\}_{(\alpha,\beta)\in\Sigma\times T}$, where
$$
H^{\alpha,\beta}=f\otimes g(J^{(\alpha)}\otimes_AC)+f\otimes g(B\otimes_AK^{(\beta)})
$$
for $(\alpha,\beta)\in\Sigma\times T$ (cf.\ Exercise \ref{exer-topologicaltensorproducts}).
Let $\widehat{H}^{\alpha,\beta}$ for each $(\alpha,\beta)\in\Sigma\times T$ be the closure of the image of $H^{\alpha,\beta}$ in $B\widehat{\otimes}_AC$; then $B\widehat{\otimes}_AC$ is Hausdorff complete with respect to the topology defined by the filtration $\{\widehat{H}^{\alpha,\beta}\}_{(\alpha,\beta)\in\Sigma\times T}$ (cf.\ {\bf \ref{ch-pre}}, \S\ref{subsub-completionfiltration}, {\bf \ref{ch-pre}}.\ref{prop-topologyfromfiltration5} (2)).
More explicitly, 
$$
B\widehat{\otimes}_AC=\varprojlim_{(\alpha,\beta)\in\Sigma\times T}B/J^{(\alpha)}\otimes_AC/K^{(\beta)}=\varprojlim_{(\lambda,\alpha,\beta)\in L}B/J^{(\alpha)}\otimes_{A/I^{(\lambda)}}C/K^{(\beta)},
$$
where $\Sigma\times T$ is considered with the ordering defined by
$$
(\alpha,\beta)\leq(\alpha',\beta')\quad\Longleftrightarrow\quad\alpha\leq\alpha'\ \textrm{and}\ \beta\leq\beta',
$$
and $L$ is the directed set given by
$$
L=\{(\lambda,\alpha,\beta)\in\Lambda\times\Sigma\times T\,|\,I^{(\lambda)}B\subseteq J^{(\alpha)},\ I^{(\lambda)}C\subseteq K^{(\beta)}\}
$$
considered with the similar ordering; notice that the map $L\rightarrow\Sigma\times T$ given by the canonical projection is cofinal.

\begin{prop}[{\cite[$\mathbf{0}_{\mathbf{I}}$, (7.7.7)]{EGA}}]\label{prop-completetensorproductadmissiblerings}
The topological ring $B\widehat{\otimes}_AC$ is an admissible ring, and $\{\widehat{H}^{\alpha,\beta}\}_{(\alpha,\beta)\in\Sigma\times T}$ gives a fundamental system of ideals of definition. \hfill$\square$
\end{prop}

\begin{lem}\label{lem-fiverprodformaladicadic}
If $B$ and $C$ are adic of finite ideal type, then so is $B\widehat{\otimes}_AC$.
If $J$ $($resp.\ $K)$ is a finitely generated ideal of definition of $B$ $($resp.\ $C)$ and if we set $H^{m,n}=f\otimes g(J^m\otimes_AC)+f\otimes g(B\otimes_AK^n)$ for $m,n\geq 0$, then $H=H^{1,1}(B\widehat{\otimes}_AC)$ gives a finitely generated ideal of definition of $B\widehat{\otimes}_AC$.
\end{lem}

\begin{proof}
Clearly, the filtration $\{H^{n,n}\}_{n\geq 0}$ gives a fundamental system of neighborhoods of $0$ for the ring $B\otimes_AC$, for the diagonal map $\N\hookrightarrow\N^2$ is cofinal.
We calculate for $k\geq 0$:
\begin{equation*}
\begin{split}
(H^{1,1})^{2k}&=[f\otimes g(J\otimes_AC)+f\otimes g(B\otimes_AK)]^{2k}\\
&\subseteq{\textstyle \sum^{2k}_{i=0}f\otimes g(J^{2k-i}\otimes_AC)\cdot f\otimes g(B\otimes_AK^i)}\\
&\subseteq f\otimes g(J^k\otimes_AC)+f\otimes g(B\otimes_AK^k)\\
&=H^{k,k}.
\end{split}
\end{equation*}
On the other hand, we clearly have $H^{2k,2k}\subseteq (H^{1,1})^{2k}$.
Hence $H^{2k,2k}\subseteq (H^{1,1})^{2k}\subseteq H^{k,k}$ holds for any $k\geq 0$, and thus the topology on $B\otimes_AC$ given by $\{H^{mn}\}_{m,n\geq 0}$ is $H^{1,1}$-adic.
Now since $H^{1,1}$ is a finitely generated ideal of definition, the Hausdorff completion $B\widehat{\otimes}_AC$ is actually the $H^{1,1}$-adic completion ({\bf \ref{ch-pre}}.\ref{prop-Iadiccompletioncomplete2}), and hence $H=H^{1,1}(B\widehat{\otimes}_AC)$ is an ideal of definition ({\bf \ref{ch-pre}}.\ref{cor-Iadiccompletionifexists1}).
\end{proof}

\begin{cor}\label{lem-prodformal1}
In the situation as in {\rm \ref{lem-fiverprodformaladicadic}}, suppose furthermore that $A$ is adic.
Let $I$ be an ideal of definition of $A$ such that $IB\subseteq J$ and $IC\subseteq K$.
Then 
$$
B\widehat{\otimes}_AC=\varprojlim_kB_k\otimes_{A_k}C_k, 
$$
where $A_k=A/I^{k+1}$, $B_k=B/J^{k+1}$, and $C_k=C/K^{k+1}$ for $k\geq 0$.
Moreover, $B\widehat{\otimes}_AC$ is an adic ring with a finitely generated ideal of definition $H$ generated by the images of $J\otimes_AC\rightarrow B\widehat{\otimes}_AC$ and $B\otimes_AK\rightarrow B\widehat{\otimes}_AC$, and we have $B\widehat{\otimes}_AC/H^{k+1}\cong B_k\otimes_{A_k}C_k$ for any $k\geq 0$. \hfill$\square$
\end{cor}

\subsubsection{Fiber products of formal schemes}\label{subsub-fiberproductsformalschemesexistence}
\index{fiber product!fiber product of formal schemes@--- (of formal schemes)|(}\index{formal scheme!fiber product of formal schemes@fiber product of ---s|(}
\begin{thm}\label{thm-fiberproductsformalexists}
The category of formal schemes has fiber products. \hfill$\square$
\end{thm}

As in \cite[$\mathbf{I}$, \S10.7]{EGA}, the general construction of fiber products reduces to the case of affine formal schemes, and for a diagram
$$
\Spf B\longrightarrow\Spf A\longleftarrow\Spf C
$$
of affine formal schemes the fiber product is given by $\Spf B\widehat{\otimes}_AC$.
The following statements are corollaries of \ref{lem-fiverprodformaladicadic} and \ref{lem-prodformal1}:
\begin{cor}\label{cor-fiverprodformaladicadic}
Let $Y\rightarrow X\leftarrow Z$ be a diagram of formal schemes, where $Y$ and $Z$ are adic formal schemes of finite ideal type\index{formal scheme!adic formal scheme@adic ---!adic formal scheme of finite ideal type@--- --- of finite ideal type}\index{adic!adic formal scheme@--- formal scheme!adic formal scheme of finite ideal type@--- --- of finite ideal type}.
Then the fiber product $Y\times_XZ$ is an adic formal scheme of finite ideal type.\hfill$\square$
\end{cor}

\begin{cor}\label{cor-fiberprodschemescheme}
Let $Y\rightarrow X\leftarrow Z$ be a diagram of formal schemes, where $Y$ and $Z$ are schemes.
Then the fiber product $Y\times_XZ$ is a scheme.\hfill$\square$
\end{cor}

The last corollary shows, in particular, that fiber products of schemes taken in the category of formal schemes coincide with the ones taken in the category of schemes.
\index{formal scheme!fiber product of formal schemes@fiber product of ---s|)}\index{fiber product!fiber product of formal schemes@--- (of formal schemes)|)}

\subsubsection{Fiber products and open immersions}\label{subsub-fiberproductsopenimmersions}
\begin{prop}\label{prop-sepqsepformal1}
Let $f\colon X\rightarrow Y$ be a morphism of formal schemes, and $V\hookrightarrow Y$ an open immersion\index{immersion!open immersion of formal schemes@open --- (of formal schemes)}.
Then the Cartesian diagram
$$
\xymatrix{V\times_YX\ \ar@{^{(}->}[r]\ar[d]_g&X\ar[d]^f\\ V\, \ar@{^{(}->}[r]&Y}
$$
of formal schemes remains Cartesian on underlying topological spaces.
In particular, we have 
$$
\image(g)=\image(f)\cap V.
$$
\end{prop}

\begin{proof}
By the construction of fiber products, one reduces to the affine situation $X=\Spf A$, $Y=\Spf B$, and $V=\Spf B_{\{h\}}$ for some $h\in B$.
The assertion in this case is easy to see.
\end{proof}

\begin{prop}\label{prop-openimmform1}
{\rm (1)} Let $f\colon X\rightarrow Y$ and $g\colon Y\rightarrow Z$ be morphisms of formal schemes.
If $f$ and $g$ are open immersions, then so is $g\circ f$.
If $g\circ f$ and $g$ are open immersions, then so is $f$.

{\rm (2)} If $S$ is a formal scheme and if $f\colon X\rightarrow X'$ and $g\colon Y\rightarrow Y'$ are two $S$-open immerions of formal schemes, then $f\times_Sg\colon X\times_SY\rightarrow X'\times_SY'$ is an open immersion.

{\rm (3)} If $S$ is a formal scheme and if $f\colon X\rightarrow Y$ is an $S$-open immersion  between formal schemes, then for any morphism $S'\rightarrow S$ of formal schemes, the induced morphism $f_{S'}\colon X\times_SS'\rightarrow Y\times_SS'$ is an open immersion.
\end{prop}

\begin{proof}
(1) is clear.
By {\bf \ref{ch-pre}}.\ref{prop-basechangestable} the assertions (2) and (3) follow from the special case of (3) with $Y=S$, which is already shown in \ref{prop-sepqsepformal1}.
\end{proof}

\subsection{Adic morphisms}\label{sub-formalnotadicmor}
\subsubsection{Adic morphisms}\label{subsub-formalnotadicmor}
\index{morphism of formal schemes@morphism (of formal schemes)!adic morphism of formal schemes@adic ---|(}
\begin{dfn}\label{dfn-adicmor}{\rm
A morphism $f\colon X\rightarrow Y$ of adic formal schemes of finite ideal type is said to be {\em adic}\index{adic!adic morphism@--- morphism} if there exists an open covering $\{V_{\alpha}\}_{\alpha\in L}$ of $Y$ and for each $\alpha\in L$ an ideal of definition $\mathscr{I}_{\alpha}$ of $V_{\alpha}$ of finite type such that for each $\alpha\in L$ the ideal pull-back\index{ideal pull-back} $\mathscr{I}_{\alpha}\O_{f^{-1}(V_{\alpha})}$ {\rm ({\bf \ref{ch-pre}}.\ref{dfn-idealpullback})} is an ideal of definition of the open formal subscheme\index{formal subscheme!open formal subscheme@open ---} $f^{-1}(V_{\alpha})\subseteq X$.}
\end{dfn}

By an argument similar to that in \cite[$\mathbf{I}$, (10,12,1)]{EGA}, in this situation, $
\mathscr{I}_{\alpha}\O_{f^{-1}(V_{\alpha})}$ for any ideal of definition of finite type $\mathscr{I}_{\alpha}$ on $V_{\alpha}$ is an ideal of definition of $f^{-1}(V_{\alpha})$.
Hence, in particular, if $Y$ itself has an ideal of definition $\mathscr{I}$ of finite type, then $\mathscr{I}\O_X$ is an ideal of definition of $X$.
(By this, in particular, we see that our definition of adic morphisms agrees with the one in \cite[$\mathbf{I}$, (10,12,1)]{EGA} in the locally Noetherian case.)

Note that the property `adic' is local on the target\index{local!local on the target@--- on the target} under Zariski topology {\rm (cf.\ {\bf \ref{ch-pre}}.\ref{dfn-catequivrelstablearrow} (1))}.
\begin{prop}\label{prop-adicmorphismbetweenaffineformalschemes}
Let $\varphi\colon A\rightarrow B$ be a morphism of adic rings of finite ideal type, and $f\colon Y=\Spf B\rightarrow X=\Spf A$ the induced morphism.
Then the morphism $f$ is adic if and only if the morphism $\varphi$ is adic {\rm (\ref{dfn-admissibleringsmorphisms} (2))}.
\end{prop}

\begin{proof}
Suppose $\varphi$ is adic, and let $I\subseteq A$ be a finitely generated ideal of definition.
Then $\mathscr{I}=I\O_X$ is an ideal of definition of $X$ of finite type (\ref{prop-adicformalschemeidealofdefinitionfinitetype4}, \ref{cor-adicformalschemeidealofdefinitionfinitetype11}).
Since $IB$ is an ideal of definition of $B$, $\mathscr{I}\O_Y=I\O_Y$ is an ideal of definition, which shows that $f$ is adic.
Conversely, if $f$ is adic, then $\mathscr{I}\O_Y$ is an ideal of definition.
Since $\mathscr{I}\O_Y=(IB)^{\Delta}$, we see that $IB$ is an ideal of definition (cf.\ \cite[$\mathbf{I}$, (10.3.5)]{EGA}).
\end{proof}

\begin{prop}\label{prop-corpropformal12}
Let $X$, $Y$, and $Z$ be adic formal schemes of finite ideal type, and $Y\rightarrow X$ and $Z\rightarrow X$ morphisms of formal schemes.
Suppose $Y\rightarrow X$ is adic.
Then the fiber product $Y\times_XZ$ is an adic formal scheme of finite ideal type, and the morphism $Y\times_XZ\rightarrow Z$ is adic.
\end{prop}

This proposition follows from the following lemma:
\begin{lem}\label{lem-prodformal11}
Consider the situation as in {\rm \ref{lem-prodformal1}}, and suppose that the morphism $A\rightarrow B$ is adic.
Then the adic ring $B\widehat{\otimes}_AC$ has $H=(B\otimes_AK)B\widehat{\otimes}_AC$ as an ideal of definition, and the map $C\rightarrow B\widehat{\otimes}_AC$ is adic.
\end{lem}

\begin{proof}
We may assume $J=IB$.
Then the ideal $H$ of $B\widehat{\otimes}_AC$ as in \ref{lem-prodformal1} coincides with $(B\otimes_AK)B\widehat{\otimes}_AC$, which is clearly the one generated by the image of $K$ under the map $C\rightarrow B\widehat{\otimes}_AC$.
\end{proof}

By \ref{lem-prodformal11} and \ref{lem-prodformal1} we have:
\begin{cor}\label{cor-prodformal21}
Let $X\rightarrow Y$ and $Z\rightarrow Y$ be morphisms of adic formal schemes of finite ideal type.
Suppose that $X\rightarrow Y$ is adic and that there exist ideals of definition $\mathscr{I}$ and $\mathscr{K}$ of finite type of $Y$ and $Z$, respectively, such that $\mathscr{I}\O_Z\subseteq\mathscr{K}$.
Set $Y_k=(Y,\O_Y/\mathscr{I}^{k+1})$, $X_k=(X,\O_X/\mathscr{I}^{k+1}\O_X)$, and $Z_k=(Z,\O_Z/\mathscr{K}^{k+1})$ for $k\geq 0$, and $W=X\times_YZ$.
Then for any $k\geq 0$ the scheme $W_k=(W,\O_W/\mathscr{K}^{k+1}\O_W)$ is isomorphic to the fiber product $X_k\times_{Y_k}Z_k$ of schemes. \hfill$\square$
\end{cor}

\begin{prop}\label{prop-adicmor1}
{\rm (1)} Let $f\colon X\rightarrow Y$ and $g\colon Y\rightarrow Z$ be morphisms of adic formal schemes of finite ideal type.
If $f$ and $g$ are adic, then so is the composition $g\circ f$.
If $g\circ f$ and $g$ are adic, then so is $f$.

{\rm (2)} Let $S$ be an adic formal scheme of finite ideal type, and $f\colon X\rightarrow X'$ and $g\colon Y\rightarrow Y'$ two adic $S$-morphisms of adic formal schemes of finite ideal type over $S$.
Then $f\times_Sg\colon X\times_SY\rightarrow X'\times_SY'$ is adic.

{\rm (3)} Let $S$ be an adic formal scheme of finite ideal type, and $f\colon X\rightarrow Y$ an adic $S$-morphism between adic formal schemes of finite ideal type over $S$. 
Then for any morphism $S'\rightarrow S$ of adic formal schemes of finite ideal type, $f_{S'}\colon X\times_SS'\rightarrow Y\times_SS'$ is adic.
\end{prop}

\begin{proof}
(1) is easy to see.
Applying {\bf \ref{ch-pre}}.\ref{prop-basechangestable} to the category of adic formal schemes (with morphisms being not necessarily adic), we deduce that (2) and (3) follows from the special case of (3) with $S=Y$, which is already shown in \ref{prop-corpropformal12}.
\end{proof}

\danger{Notice that in the statements (2) and (3) of \ref{prop-adicmor1} the formal schemes $X$, $Y$, $X'$, $Y'$, and $S'$ are adic formal scheme over $S$, but not necessarily adic over $S$, that is, the structural maps such as $X\rightarrow S$ are not necessarily adic.}

In the sequel we employ the following convention, which remains in force throughout this book, in order to distinguish adicness of morphisms from that of formal schemes alone: 

\medskip\noindent
{\bf Convention.} {\sl For adic formal schemes $X$ and $S$ of finite ideal type,}
\begin{itemize}
\item {\sl by `$X$ over $S$' we only mean that $X$ is simply considered with a morphism $X\rightarrow S$, {\em not necessarily adic}, of formal schemes}; 
\item {\sl if, however, we say `$X$ adic over $S$', we mean that the structural map $X\rightarrow S$ is {\em adic}.}
\end{itemize}
\index{morphism of formal schemes@morphism (of formal schemes)!adic morphism of formal schemes@adic ---|)}

\subsubsection{Adicness of diagonal maps}\label{subsub-adicnessdiagonalmaps}
\begin{prop}\label{adicmorphismdiagonal}
Let $f\colon X\rightarrow Y$ be a morphism of formal schemes, and suppose that $X$ is adic of finite ideal type.
Then the diagonal map\index{diagonal map}
$$
\Delta_X\colon X\longrightarrow X\times_YX
$$
is an adic morphism.
\end{prop}

Notice that due to \ref{cor-fiverprodformaladicadic} the formal scheme $X\times_YX$ is adic of finite ideal type; notice also that here we do not assume that the map $f$ nor $Y$ is adic.
\begin{proof}
We may assume that $X$ and $Y$ are affine: $X=\Spf B\rightarrow Y=\Spf A$.
This gives the situation as in \ref{lem-fiverprodformaladicadic}, where $B=C$.
Using the notation therein, we know that $B\widehat{\otimes}_AB$ has the ideal $H$ generated by $f\otimes g(J\otimes_AB)+f\otimes g(B\otimes_AJ)$ as an ideal of definition.
The diagonal map $\Delta_X$ as above is induced from the codiagonal map $B\widehat{\otimes}_AB\rightarrow B$, by which we clearly have $HB=J$, thereby the assertion.
\end{proof}

\subsection{Formal completion}\label{sub-formalcompletionschemes}
\subsubsection{Formal schemes as inductive limits of schemes}\label{subsub-formalindlimsch}
\index{limit!inductive limit@inductive ---}
We have already seen in \S\ref{subsub-formalschemesidealsofdefinition} that for any formal scheme $X$ and an ideal of definition $\mathscr{I}$, the locally ringed space $(X,\O_X/\mathscr{I})$ is a scheme.
\begin{prop}[{\cite[$\mathbf{I}$, (10.6.2)]{EGA}}]\label{prop-formalindlimschfact}
Let $X$ be a formal scheme, and suppose there exists a fundamental system of ideals of definition\index{ideal of definition!fundamental system of ideals of definition@fundamental system of ---s ---} $\{\mathscr{I}^{(\lambda)}\}_{\lambda\in\Lambda}$, and set $X_{\lambda}=(X,\O_X/\mathscr{I}^{(\lambda)})$ for $\lambda\in\Lambda$.
Then $\{X_{\lambda}\}_{\lambda\in\Lambda}$ with the canonical closed immersions forms an inductive system\index{system!inductive system@inductive ---} of schemes, and we have
$$
X=\varinjlim_{\lambda\in\Lambda}X_{\lambda}
$$
in the category of formal schemes. \hfill$\square$
\end{prop}

\begin{prop}\label{prop-formalindlimsch}
Let $X$ be a topological space, and consider the projective system $\{\O_i,u_{ij}\}_{i\in I}$ of sheaves of rings on $X$ indexed by a directed set $I\neq\emptyset$ that admits an at most countable final subset.
Fix $0\in I$ and for any $i\geq 0$ let $\mathscr{I}_i$ be the kernel of $u_{0i}\colon\O_i\rightarrow\O_0$.
Suppose that the following conditions are satisfied$:$
\begin{itemize}
\item[{\rm (a)}] for any $i\geq 0$ the ringed space $(X,\O_i)$ is a scheme$;$
\item[{\rm (b)}] for any $i\geq 0$ and any $x\in X$ there exists an open neighborhood $U_i$ of $x$ in $X$ such that the sheaf $\mathscr{I}_i$ restricted on $U_i$ is nilpotent$;$
\item[{\rm (c)}] for any $i\geq j\geq 0$ the morphism $u_{ji}$ is surjective.
\end{itemize}
Then the ringed space $(X,\O)$ equipped with the projective limit sheaf $\O=\varprojlim_{i\in I}\O_i$ is a formal scheme$;$ the canonical morphism $u_i\colon\O\rightarrow\O_i$ for $i\geq 0$ is surjective$;$ if one puts $\mathscr{I}^{(i)}$ to be the kernel of $u_i$, then $\{\mathscr{I}^{(i)}\}_{i\geq 0}$ forms a fundamental system of ideals of definition, and $\mathscr{I}^{(0)}$ coincides with the projective limit $\varprojlim_{i\geq 0}\mathscr{I}_i$.
\end{prop}

\begin{proof}
We may assume $I=\N$ without loss of generality (Exercise \ref{exer-finalcountable}). Then the proposition in this case is nothing but \cite[$\mathbf{I}$, (10.6.3)]{EGA}.
\end{proof}

It is often delicate to judge whether a given formal scheme, given as an inductive limit of schemes as above, is adic or not; the situation is somewhat similar to judging adicness of a given topological ring ({\bf \ref{ch-pre}}, \S\ref{subsub-criterionadicness}).
But similarly to {\bf \ref{ch-pre}}.\ref{prop-criterionadicness1} one has a criterion for the property `adic of finite ideal type', as follows:
\begin{prop}\label{prop-formalindlimschadic}
In the situation as in {\rm \ref{prop-formalindlimsch}} with $I=\N$ we assume that for $i\geq j$ the kernel of the morphism $u_{ji}$ coincides with $\mathscr{I}^{j+1}_i$ and that $\mathscr{I}_1/\mathscr{I}^2_1$ is of finite type over $\O_0=\O_1/\mathscr{I}_1$.
Then the formal scheme $X$ is adic$;$ if one puts $\mathscr{I}^{(k)}$ to be the kernel of the map $\O_X\rightarrow\O_k$ and $\mathscr{I}=\mathscr{I}^{(0)}$, then $\mathscr{I}^{(k)}=\mathscr{I}^{k+1}$ and $\mathscr{I}/\mathscr{I}^2$ is isomorphic to $\mathscr{I}_1$.
In particular, the ideal $\mathscr{I}$ is an ideal of definition of the adic formal scheme $X$. \hfill$\square$
\end{prop}

We refer to \cite[$\mathbf{I}$, (10.6.4)]{EGA} for the proof.
Notice that, in this situation, the ideal of definition $\mathscr{I}$ as above is of finite type, and hence the formal scheme $X$ thus obtained is an adic formal scheme of finite ideal type (Exercise \ref{exer-adicformalschemesidealsofdefinitionfinitude}).

\subsubsection{Formal completion of schemes}\label{subsub-formalcompletionsch}
\index{completion!formal completion@formal ---|(}
Let $X$ be a scheme, and $Y\subseteq X$ a closed subscheme.
We denote by $\mathscr{I}_Y$ the quasi-coherent ideal sheaf of $\O_X$ that defines $Y$.
For each $k\geq 0$ we set $X_k=(Y,\O_X/\mathscr{I}^{k+1}_Y)$.
Then the inductive limit 
$$
\widehat{X}|_Y=\varinjlim_{k\geq 0}X_k
$$
is represented by a formal scheme (\ref{prop-formalindlimsch}), called the {\em formal completion of $X$ along $Y$} (cf.\ {\bf \ref{ch-pre}}, \S\ref{sub-schpair} and \cite[$\mathbf{I}$, \S10.8]{EGA}).
The formal scheme $\widehat{X}|_Y$ is sometimes denoted simply by $\widehat{X}$, when $Y$ is clear by the context.

In the affine situation $X=\Spec A$ and $Y=\Spec A/I$, where $I$ is an ideal of $A$, the formal completion $\widehat{X}|_Y$ is nothing but the affine formal scheme $\Spf A^{\wedge}_{I^{\bullet}}$, where $A^{\wedge}_{I^{\bullet}}$ is the Hausdorff completion of $A$ with respect to the topology defined by $I^{\bullet}=\{I^n\}_{n\geq 0}$, which is an admissible ring (\ref{exa-admissibleringformalcompletion}).

The following proposition is clear; cf.\ \ref{exa-admissibleringformalcompletion2}:
\begin{prop}\label{prop-formalcompletion1}
Suppose that the closed subscheme $Y$ is of finite presentation.
Then the formal scheme $\widehat{X}=\widehat{X}|_Y$ is adic of finite ideal type and has an ideal of definition $\mathscr{I}=\mathscr{I}_Y\O_{\widehat{X}}$ of finite type. \hfill$\square$
\end{prop}

Here we include a basic fact on the formal completion of morphisms of schemes:
\begin{prop}[{\cite[(10.9.9)]{EGAInew}}]\label{prop-formalcomplmor}
Let $Y$ be a scheme, $Z\subseteq Y$ a closed subscheme, and $f\colon X\rightarrow Y$ a morphism of schemes. 
Let $W=f^{-1}(Z)$, and set $\widehat{X}=\widehat{X}|_W$ and $\widehat{Y}=\widehat{Y}|_Z$.
Then the commutative diagram
$$
\xymatrix{X\ar[d]_f&\widehat{X}\ar[l]_j\ar[d]^{\widehat{f}}\\ Y&\widehat{Y}\ar[l]^i}
$$
is Cartesian in the category of formal schemes. \hfill$\square$
\end{prop}

The following proposition follows from \ref{cor-prodformal21}:
\begin{prop}\label{prop-prodformal22}
Let $X\rightarrow Y$ and $Z\rightarrow Y$ be morphisms of schemes, and $W\subseteq Y$ a closed subscheme of $Y$ of finite presentation.
Let $\widehat{X}$ $($resp.\ $\widehat{Y}$, resp.\ $\widehat{Z})$ be the formal completion along the $W$ $($resp.\ $W_X$, resp.\ $W_Z)$.
Then there exists a canonical isomorphism 
$$
\widehat{X}\times_{\widehat{Y}}\widehat{Z}\stackrel{\sim}{\longrightarrow}\widehat{X\times_YZ},
$$
where $\widehat{X\times_YZ}$ denotes the formal completion of $X\times_YZ$ along $W_X\times_YW_Z$. \hfill$\square$
\end{prop}

\subsubsection{Formal completion of quasi-coherent sheaves}\label{subsub-formalcompletionqcohsh}
Let $X$ be a scheme, $Y\subseteq X$ a closed subscheme with the defining ideal sheaf $\mathscr{I}_Y$, and $\mathscr{F}$ a quasi-coherent\index{quasi-coherent!quasi-coherent OX module on schemes@--- sheaf (on a scheme)} $\O_X$-module.
Then the {\em formal completion of $\mathscr{F}$ along $Y$}, denoted by $\widehat{\mathscr{F}}|_Y$ or by $\widehat{\mathscr{F}}$, is given as follows (\cite[(10.8.2)]{EGAInew}):
$$
\widehat{\mathscr{F}}|_Y=\varprojlim_n\mathscr{F}\otimes_{\O_X}(\O_X/\mathscr{I}_Y^n)=\varprojlim_n\mathscr{F}/\mathscr{I}_Y^n\mathscr{F}.
$$
This is an $\O_{\widehat{X}|_Y}$-module; we will see later in \S\ref{sec-adicallyqcoh} that $\widehat{\mathscr{F}}|_Y$ thus obtained, in case $Y$ is a closed subscheme of finite presentation, is an example of what we will later call {\em adically quasi-coherent sheaves} (\ref{prop-adicqchohexa}).
Notice that in the affine situation $X=\Spec A$, $Y=\Spec A/I$, and $\mathscr{F}=\til{M}$, we have
$$
\Gamma(\widehat{X}|_Y,\mathscr{F})=\Gamma(Y,\varprojlim_n\mathscr{F}/\mathscr{I}_Y^n\mathscr{F})=\varprojlim_n M/I^nM=M^{\wedge}_{I^{\bullet}}.
$$

\begin{prop}\label{prop-adqformalprop}
Let $(X,Y)$ be a pseudo-adhesive pair\index{pair!adhesive pair@adhesive ---!pseudo adhesive pair@pseudo-{---} ---}\index{adhesive!adhesive pair@--- pair!pseudo adhesive pair@pseudo-{---} ---}\index{pseudo-adhesive} of schemes {\rm ({\rm {\bf \ref{ch-pre}}.\ref{dfn-schpairadh}})}, and set $\widehat{X}=\widehat{X}|_Y$.

{\rm (1)} For any quasi-coherent sheaf $\mathscr{F}$ on $X$ of finite type, the canonical morphism 
$$
i^{\ast}\mathscr{F}\longrightarrow\widehat{\mathscr{F}}=\widehat{\mathscr{F}}|_Y,
$$
where $i\colon\widehat{X}\rightarrow X$ is the canonical morphism, is an isomorphism.

{\rm (2)} The map $i\colon\widehat{X}\rightarrow X$ of locally ringed spaces is flat, that is, for any $x\in\widehat{X}$ the ring $\O_{\widehat{X},x}$ is flat over $\O_{X,i(x)}$ $($cf.\ {\rm {\bf \ref{ch-pre}}.\ref{dfn-flatnessringedsp} (2))}.
\end{prop}

\begin{proof}
Since the question is local, we may assume that $X=\Spec A$, where we are given a pseudo-adhesive pair $(A,I)$ ({\bf \ref{ch-pre}}.\ref{prop-schpairadh}), and hence that $\widehat{X}=\Spf\widehat{A}$.
Denote the canonical map $A\rightarrow\widehat{A}$ by $j$.

(1) Take the finitely generated $A$-module $M$ such that $\mathscr{F}=\til{M}$.
By Exercise \ref{exer-adqformalprop1} and the equality $\widehat{A_f}=\widehat{A}_{\{j(f)\}}$ (\cite[$\mathbf{0}_{\mathbf{I}}$, (7.6.2)]{EGA}), what to prove is the following:
for any $f\in A$ the canonical morphism
$$
M_f\otimes_{A_f}\widehat{A_f}\longrightarrow\widehat{M_f}
$$
is an isomorphism.
Since $M_f$ is finitely generated over $A_f$ and since $A_f=A[T]/(fT-1)$ is $IA_f$-adically pseudo-adhesive ({\bf \ref{ch-pre}}.\ref{prop-adhesive1}), this follows from {\bf \ref{ch-pre}}.\ref{prop-btarf1} (1).

(2) By {\bf \ref{ch-pre}}.\ref{prop-btarf1} (2) we see that the canonical map $j_f\colon A_f\rightarrow\widehat{A_f}$ is flat.
By Exercise \ref{exer-adqformalprop1} and the equality $\widehat{A_f}=\widehat{A}_{\{j(f)\}}$, we see that the map between stalks 
$$
\O_{X,i(x)}\longrightarrow\O_{\widehat{X},x}
$$
is the filtered inductive limit of the maps of the form $j_f$ as above.
Since the flatness is preserved by filtered inductive limits, the assertion follows.
\end{proof}
\index{completion!formal completion@formal ---|)}

\subsection{Categories of formal schemes}\label{sub-categoryofformalschemes}
\subsubsection{Notation}\label{subsub-categoryformalschnotation}
First of all, we set
\begin{itemize}
\item $\Fs=$ the category of all formal schemes.
\end{itemize}
This category has fiber products\index{fiber product!fiber product of formal schemes@--- (of formal schemes)}\index{formal scheme!fiber product of formal schemes@fiber product of ---s} and the final object $\Spec\Z$.
As mentioned before (\ref{rem-clear}), the category $\Fs$ contains the category $\Sch$ of schemes as the full subcategory consisting of $0$-adic formal schemes.
Moreover, it has a strict initial\index{strictly initial} object (the empty scheme $\emptyset=\Spec 0$) and disjoint sums (cf.\ {\bf \ref{ch-pre}}.\ref{dfn-disjointunion} and Exercise \ref{exer-disjointsumformalsch}).

The category $\Fs$ has the following full subcategories:
\begin{itemize}
\item $\Af\Fs=$ the category of {\em affine} formal schemes;
\item $\Ac\Fs=$ the category of {\em adic} formal schemes {\em of finite ideal type};
\item $\Af\Ac\Fs=$ the category of {\em affine adic} formal schemes {\em of finite ideal type}.
\end{itemize}
Notice that, here, the condition `adic' is only put on {\em objects} and not on {\em morphisms}. We set
\begin{itemize}
\item $\Ac\Fs^{\ast}=$ the category of adic formal schemes of finite ideal type {\em with adic morphisms};
\item $\Af\Ac\Fs^{\ast}=$ the category of affine adic formal schemes of finite ideal type {\em with adic morphisms},
\end{itemize}
and will follow the principle that the superscript `$\ast$' always means that the morphisms between objects under consideration are only adic morphisms.

The categories $\Af\Fs$, $\Ac\Fs$ and $\Af\Ac\Fs$ have fiber products (\ref{cor-fiverprodformaladicadic}), strictly initial objects, disjoint sums, and final objects.
The categories $\Ac\Fs^{\ast}$ and $\Af\Ac\Fs^{\ast}$ have fiber products (\ref{prop-corpropformal12}), strictly initial objects, and disjoint sums.

For any formal scheme $S$ we write $\Fs_S$, $\Ac\Fs_S$, and $\Ac\Fs^{\ast}_S$ etc.\ the respective categories of objects over $S$.
Here, even in case $S$ is adic of finite ideal type, objects of the category $\Ac\Fs^{\ast}_S$ are {\em not necessarily adic over $S$} (that is, the structural maps $X\rightarrow S$ are not assumed to be adic; see the general convention at the end of \S\ref{subsub-formalnotadicmor}).
To specify the categories of adic formal schemes of finite ideal type that are adic over a fixed adic formal scheme $S$ of finite ideal type, we set
\begin{itemize}
\item $\Ac\Fs^{\ast}_{/S}=$ the category of adic formal schemes of finite ideal type {\em adic over $S$};
\item $\Af\Ac\Fs^{\ast}_{/S}=$ the category of affine adic formal schemes of finite ideal type {\em adic over $S$}.
\end{itemize}
In these categories all arrows are automatically adic due to \ref{prop-adicmor1} (1), and hence $\Ac\Fs^{\ast}_{/S}$ and $\Af\Ac\Fs^{\ast}_{/S}$ are full subcategories of $\Fs_S$, $\Ac\Fs_S$, and $\Ac\Fs^{\ast}_S$ etc.

The categories $\Fs_S$, $\Af\Fs_S$, $\Ac\Fs_S$, $\Af\Ac\Fs_S$, $\Ac\Fs^{\ast}_{/S}$, and $\Af\Ac\Fs^{\ast}_{/S}$ have fiber products, strictly initial objects, and disjoint sums; $\Fs_S$ and $\Ac\Fs^{\ast}_{/S}$ have final objects; $\Af\Fs_S$ (resp.\ $\Ac\Fs_S$, resp.\ $\Af\Ac\Fs_S$) has final objects if $S$ is affine (resp.\ adic of finite ideal type, resp.\ affine adic of finite ideal type).

\subsubsection{Properties of morphisms in $\Fs$}\label{subsub-genpropertymorphism}
Consider a property $P$ for morphisms of formal schemes such that 
\begin{itemize}
\item[{\bf (I)}] any isomorphism satisfies $P$;
\item[{\bf (C)}] if $f\colon X\rightarrow Y$ and $g\colon Y\rightarrow Z$ are arrows satisfying $P$, then the composition $g\circ f$ satisfies $P$.
\end{itemize}
As in {\bf \ref{ch-pre}}, \S\ref{subsub-basechangestable} we consider the subcategory $\mathscr{D}=\mathscr{D}_P$ of $\Fs$ consisting of morphisms of formal schemes satisfying $P$.
Then, as we have seen in {\bf \ref{ch-pre}}.\ref{prop-basechangestable}, the following conditions are equivalent, and when they are fulfilled, we say that the property $P$ is {\em base-change stable}\index{stable!base-change stable@base-change ---}:
\begin{itemize}
\item[{\boldmath $(\mathbf{B}_1)$}] for any $Z$-morphisms $f\colon X\rightarrow Y$ and $g\colon X'\rightarrow Y'$ satisfying $P$ of formal schemes over a formal scheme $Z$, the induced morphism
$$
f\times_Zg\colon X\times_ZY\longrightarrow X'\times_ZY'
$$
satisfies $P;$
\item[{\boldmath $(\mathbf{B}_2)$}] for any $Z$-morphism $f\colon X\rightarrow Y$ satisfying $P$ of formal schemes over a formal scheme $Z$ and for any morphism $Z'\rightarrow Z$ of formal schemes, the induced morphism
$$
f_{Z'}\colon X\times_ZZ'\longrightarrow Y\times_ZZ'
$$
satisfies $P;$
\item[{\boldmath $(\mathbf{B}_3)$}] for any morphism $f\colon X\rightarrow Y$ satisfying $P$ and any morphism $Y'\rightarrow Y$ of formal schemes, the induced arrow
$$
f_{Y'}\colon X\times_YY'\longrightarrow Y'
$$
satisfies $P$.
\end{itemize}

\subsubsection{Properties of morphisms in $\Ac\Fs$}\label{subsub-genpropertymorphismadic}
Let $P$ be a property of morphisms of {\em adic} formal schemes {\em of finite ideal type} that satisfies the conditions {\bf (I)} and {\bf (C)} in \S\ref{subsub-genpropertymorphism}.
We assume that the property $P$ implies that the morphisms in question are adic, that is to say, that the corresponding subcategory $\mathscr{D}_P$ is contained in $\Ac\Fs^{\ast}$.
Then the conditions in {\bf \ref{ch-pre}}.\ref{prop-basechangestable}, which are equivalent to each other, are written in the following way:
\begin{itemize}
\item[{\boldmath $(\mathbf{B}_1)$}] for any adic $Z$-morphisms $f\colon X\rightarrow Y$ and $g\colon X'\rightarrow Y'$  satisfying $P$ of adic formal schemes of finite ideal type over an adic formal scheme of finite ideal type $Z$, the induced morphism
$$
f\times_Zg\colon X\times_ZY\longrightarrow X'\times_ZY'
$$
satisfies $P$ (notice that $f\times_Zg$ is automatically adic);
\item[{\boldmath $(\mathbf{B}_2)$}] for any adic $Z$-morphism $f\colon X\rightarrow Y$ satisfying $P$ of adic formal schemes of finite ideal type over an adic formal scheme of finite ideal type $Z$ and for any morphism $Z'\rightarrow Z$ of adic formal schemes of finite ideal type, the induced morphism
$$
f_{Z'}\colon X\times_ZZ'\longrightarrow Y\times_ZZ'
$$
satisfies $P;$
\item[{\boldmath $(\mathbf{B}_3)$}] for any adic morphism $f\colon X\rightarrow Y$ satisfying $P$ and any morphism $Y'\rightarrow Y$ of adic formal schemes of finite ideal type, the induced arrow
$$
f_{Y'}\colon X\times_YY'\longrightarrow Y'
$$
satisfies $P$.
\end{itemize}

\danger{Similarly to the warning at the end of \S\ref{subsub-formalnotadicmor}, the morphisms without being spelled out to be adic are not assumed to be adic.
For example, $X$, $Y$, $X'$, and $Y'$ in {\boldmath $(\mathbf{B}_1)$} are not assumed to be adic over $Z$.}

In the category $\Ac\Fs$, we can, moreover, consider the following conditions:
\begin{itemize}
\item[{\bf (R)}] for any adic morphism $f\colon X\rightarrow Y$ where $Y$ has an ideal of definition $\mathscr{I}$ of finite type, the following conditions are equivalent:
\begin{itemize}
\item[{\rm (a)}] $f$ satisfies $P;$
\item[{\rm (b)}] $f_k$ satisfies $P$ for any $k\geq 0;$
\end{itemize}
here, for any integer $k\geq 0$ we set $X_k=(X,\O_X/\mathscr{I}^{k+1}\O_X)$ and $Y_k=(Y,\O_Y/\mathscr{I}^{k+1})$ and set $f_k\colon X_k\rightarrow Y_k$ to be the induced map of schemes$;$
\item[{\bf (F)}] for any morphism $f\colon X\rightarrow Y$ of schemes satisfying $P$ and any closed subscheme $Z\subseteq Y$ of finite presentation, the induced map between formal completions\index{completion!formal completion@formal ---} $\widehat{f}|_Z\colon\widehat{X}|_{f^{-1}(Z)}\rightarrow\widehat{Y}|_Z$ satisfies $P$.
\end{itemize}

\begin{prop}\label{prop-genpropertymorphismadic1}
Let $P$ be a property that satisfies {\bf (R)} and is local on the target\index{local!local on the target@--- on the target} under Zariski topology {\rm (cf.\ {\bf \ref{ch-pre}}.\ref{dfn-catequivrelstablearrow} (1))}. Then $P$ is base-change stable\index{stable!base-change stable@base-change ---} in $\Ac\Fs$ if and only if the property $P$ restricted on morphisms of schemes is base-change stable in $\Sch$.
\end{prop}

\begin{proof}
The `only if' part is trivial.
Suppose the property $P$ restricted on morphisms of schemes is base-change stable in $\Sch$.
Let $f\colon X\rightarrow Y$ be an adic morphism of adic formal schemes of finite ideal type that satisfies $P$, and $Y'\rightarrow Y$ a (not necessarily adic) morphism of adic formal schemes of finite ideal type.
Let $Y=\bigcup_{\alpha\in L}U_{\alpha}$ be an open covering of $Y$.
Then $Y'$ is covered by open subsets $Y'\times_YU_{\alpha}=Y'_{\alpha}$.
Since $P$ is local on the target under Zariski topology, it suffices to show the property $P$ for the base change $f'_{\alpha}\colon X'_{\alpha}\rightarrow Y'_{\alpha}$, where $X'_{\alpha}=X\times_YY'_{\alpha}$.
Hence we may assume that $Y$ has an ideal of definition $\mathscr{I}$ of finite type.
Similarly, since $Y'$ is covered by an open subsets that have ideals of definition of finite type, we may further assume that $Y'$ has an ideal of definition $\mathscr{J}$ of finite type.
Since we may furthermore assume that $Y$ and $Y'$ are affine, we may assume $\mathscr{I}\O_{Y'}\subseteq\mathscr{J}$.

Now let $f_k\colon X_k\rightarrow Y_k$ for any integer $k\geq 0$ be the induced morphism of schemes as above, and set $Y'_k=(Y',\O_{Y'}/\mathscr{J}^{k+1})$.
We have the induced morphism $Y'_k\rightarrow Y_k$ for any $k\geq 0$.
By \ref{cor-prodformal21} the similarly defined $(X\times_YY')_k$ is isomorphic to $X_k\times_{Y_k}Y'_k$, and $(f\times_YY')_k=f_k\times_{Y_k}Y'_k$ for any $k\geq 0$.
By the assumption the map $f_k$ satisfies $P$, and hence $(f\times_YY')_k=f_k\times_{Y_k}Y'_k$ satisfies $P$.
But this implies $f\times_YY'$ satisfies $P$, as desired.
\end{proof}

\begin{prop}\label{prop-genpropertymorphismadic2}
Let $P$ be a property that satisfies {\bf (R)} and is base-change stable\index{stable!base-change stable@base-change ---} restricted on the morphisms of schemes.
Then $P$ satisfies {\bf (F)}.
\end{prop}

\begin{proof}
Suppose that $f\colon X\rightarrow Y$ is a morphism of schemes that satisfies $P$ and that $Z$ is a closed subscheme of $Y$ of finite presentation.
For any $k\geq 0$ the scheme $(\widehat{X}|_Z)_k$ (defined as above) is the closed subscheme of $X$ defined by the ideal $\mathscr{I}^{k+1}$, where $\mathscr{I}$ is the defining ideal of $Z$ in $X$.
By the base-change stability for morphisms of schemes we know that each $(\widehat{f}|_Z)_k$ satisfies $P$.
Then {\bf (R)} implies that $\widehat{f}|_Z$ satisfies $P$.
\end{proof}

\subsubsection{Adicalization}\label{subsub-adicalization}
\index{adicalization|(}
Let $P$ be a property of morphisms of schemes that satisfies {\bf (I)} and {\bf (C)} in \S\ref{subsub-genpropertymorphism} and is stable under Zariski topology\index{stable!stable under a topology@--- under a topology} ({\bf \ref{ch-pre}}.\ref{dfn-catequivrelstablearrow} (1)).
Let $f\colon X\rightarrow Y$ be a morphism of adic formal schemes of finite ideal type.
Then we say that $f$ satisfies {\em adically $P$} if the following conditions are satisfied: 
\begin{itemize}
\item[{\rm (a)}] $f$ is adic;
\item[{\rm (b)}] there exists an open covering $Y=\bigcup_{\alpha\in L}V_{\alpha}$ and for each $\alpha$ an ideal of definition of finite type $\mathscr{I}_{\alpha}$ on $V_{\alpha}$ such that for any $\alpha\in L$ and $k\geq 0$ the induced morphism of schemes
$$
U_{\alpha,k}=(U_{\alpha},\O_{U_{\alpha}}/\mathscr{I}^{k+1}_{\alpha}\O_{U_{\alpha}})\longrightarrow V_{\alpha,k}=(V_{\alpha},\O_{V_{\alpha}}/\mathscr{I}^{k+1}_{\alpha})
$$
satisfies $P$, where $U_{\alpha}=f^{-1}(V_{\alpha})$.
\end{itemize}

\begin{prop}\label{prop-adicalization1}
Let $f\colon X\rightarrow Y$ be an adic morphism of adic formal schemes of finite ideal type.
Then $f$ satisfies adically $P$ if and only if for any open subspace $V\subseteq Y$ and any ideal of definition of finite type $\mathscr{I}$ on $V$, the induced morphism
$$
U_0=(U,\O_U/\mathscr{I}\O_U)\longrightarrow V_0=(V,\O_V/\mathscr{I})
$$
of schemes satisfies $P$, where $U=f^{-1}(V)$.
\end{prop}

\begin{proof}
The `if' part is clear.
Suppose $f$ is adically $P$.
Then we have an open covering $Y=\bigcup_{\alpha\in L}V_{\alpha}$ as above.
Take for any $\alpha\in L$ a sufficiently large $k\geq 0$ such that $\mathscr{I}^{k+1}_{\alpha}\subseteq\mathscr{I}$ on $V\cap V_{\alpha}$.
Then on $V\cap V_{\alpha}$ the morphism $U_0\rightarrow V_0$ in question is obtained by base change of $U_{\alpha,k}\rightarrow V_{\alpha,k}$; since $P$ is stable under Zariski topology, we deduce that $U_0\rightarrow V_0$ is $P$.
\end{proof}

\begin{prop}\label{prop-adicalization2}
The property `adically $P$' satisfies {\bf (I)}, {\bf (C)} $($in the category $\Ac\Fs)$, {\bf (R)}, and {\bf (F)}.
In particular, it is base-change stable, if the property $P$ is base-change stable.
\end{prop}

\begin{proof}
It is clear that {\bf (I)} is satisfied.
By \ref{prop-adicalization1} we see that also {\bf (C)} and {\bf (R)} are satisfied.
The other assertion follows from \ref{prop-genpropertymorphismadic1} and \ref{prop-genpropertymorphismadic2}.
\end{proof}

In some cases, it is of little problem to drop `adically' from `adically $P$'; for example, `adically of finite type' is, due to \ref{prop-cortopfintype11} below, the same as what we appropriately call `of finite type'.
A similar example is `adically affine', which will turn out to be just `affine (and adic)' (\S\ref{subsub-fundamentalpropertiesaffineadicbasechange}). But in some other cases, it is important to distinguish `adically $P$' from `$P$', for example: 
\begin{exas}\label{exas-adicalizationexamples}{\rm 
(1) $P=$ `flat': {\em adically flat} morphisms\index{morphism of formal schemes@morphism (of formal schemes)!adically flat morphism of formal schemes@adically flat ---} will be of essential importance, since `flat' in formal geometry in general is not such a reasonable notion.
Adically flat morphisms are discussed in more detail in \S\ref{subsub-adicallyflat} below.

(2) $P=$ `quasi-affine': {\em adically quasi-affine} morphisms\index{morphism of formal schemes@morphism (of formal schemes)!adically quasi-affine morphism of formal schemes@adically quasi-affine ---} will be of technical importance in \S\ref{subsub-formalalgebraicspacesdef}, for it satisfies adically flat effective descent.}
\end{exas}
\index{adicalization|)}

\subsection{Quasi-compact and quasi-separated morphisms}\label{sub-qcompqsepformal}
\subsubsection{Quasi-compact morphisms and some preliminary facts on diagonal morphisms}\label{subsub-qcompqsepformalpre}
\begin{dfn}\label{dfn-qcompformal}{\rm 
A morphism $f\colon X\rightarrow Y$ of formal schemes is said to be {\em quasi-compact}\index{morphism of formal schemes@morphism (of formal schemes)!quasi-compact morphism of formal schemes@quasi-compact ---}\index{quasi-compact!quasi-compact morphism of formal schemes@--- morphism (of formal schemes)} if, for any quasi-compact open subset $U$ of $Y$, $f^{-1}(U)$ is quasi-compact.}
\end{dfn}

Here, observe that the condition `quasi-compact' is a topological one (cf.\ {\bf \ref{ch-pre}}.\ref{dfn-quasicompactness} (2)).
In particular, if $f\colon X\rightarrow Y$ is an {\em adic} morphism between adic formal schemes of finite ideal type and $Y$ has an ideal of definition $\mathscr{I}$ of finite type, then $f$ is quasi-compact if and only if the induced morphism of schemes $f_0\colon X_0=(X,\O_X/\mathscr{I}\O_X)\rightarrow Y_0=(Y,\O_Y/\mathscr{I})$ is quasi-compact.

In order to define quasi-separated morphisms, we need to prove some preparatory results:
\begin{lem}\label{lem-sepqsepformal2}
Let $f\colon X\rightarrow Y$ be a morphism of formal schemes, and consider the diagonal map\index{diagonal map} $\Delta_X\colon X\rightarrow X\times_YX$.
Then for any open immersion $U\hookrightarrow X$ the diagram
$$
\xymatrix{U\ar[d]\ar[r]^(.4){\Delta_U}&U\times_YU\ar[d]\\ X\ar[r]_(.4){\Delta_X}&X\times_YX}
$$
is Cartesian in the category of formal schemes.
\end{lem}

\begin{proof}
Let $W=X\times_{(X\times_YX)}(U\times_YU)$.
Since $U\times_YU\rightarrow X\times_YX$ is an open immersion (\ref{prop-openimmform1} (2)), so is $W\rightarrow X$.
We have a morphism $\varphi\colon U\rightarrow W$ such that the resulting diagram
$$
\xymatrix@-1ex{
U\ar@/^/[drr]^{\Delta_U}\ar[dr]^{\varphi}\ar@/_/[ddr]\\
&W\ar[r]\ar[d]&U\times_YU\ar[d]\\
&X\ar[r]_(.4){\Delta_X}&X\times_YX}
$$
commutes. 
Then by \ref{prop-openimmform1} (1) the morphism $\varphi$ is an open immersion.
Consider the composition
$$
W\longrightarrow U\times_YU\stackrel{\mathrm{pr}_i}{\longrightarrow}U,
$$
for $i=1,2$, and denote them by $\psi_i$.
If we denote by $\alpha$ the open immersion $W\hookrightarrow X$, then we have $\alpha\circ\varphi\circ\psi_i=\alpha$.
Then by \ref{prop-openimmform1} (1) the morphism $\psi_i$ is an open immersion.
But since we have $\psi_i\circ\varphi=\id_U$, $\varphi$ must be an isomorphism.
\end{proof}

\begin{cor}\label{cor-sepqsepformal21}
Let $f\colon X\rightarrow Y$ be a morphism of formal schemes.
Then the diagonal map $\Delta_X\colon X\rightarrow X\times_YX$ maps the underlying topological space of $X$ homeomorphically onto its image $\Delta_X(X)$ endowed with the subspace topology induced from the topology on $X\times_YX$.
\end{cor}

\begin{proof}
Since $\mathrm{pr}_1\circ\Delta_X=\id_X$, the diagonal map $\Delta_X$ is clearly injective.
Hence it suffices to show that for any open subset $U\subseteq X$ the image $\Delta_X(U)$ is open in $\Delta_X(X)$.
By \ref{lem-sepqsepformal2} and \ref{prop-sepqsepformal1}, identifying $U\times_YU$ with the open subset of $X\times_YX$ by the open immersion $U\times_YU\rightarrow X\times_YX$, we have
$$
\Delta_X(U)=\Delta_U(U)=\Delta_X(X)\cap U\times_YU,
$$
whence the result.
\end{proof}

\begin{cor}\label{cor-qsepformal001}
The following conditions for a morphism $f\colon X\rightarrow Y$ of formal schemes are equivalent$:$
\begin{itemize}
\item[{\rm (a)}] the diagonal morphism $\Delta\colon X\rightarrow X\times_YX$ is quasi-compact$;$
\item[{\rm (b)}] the inclusion $\Delta(X)\hookrightarrow X\times_YX$ of the underlying topological spaces is quasi-compact $(${\rm {\bf \ref{ch-pre}}.\ref{dfn-quasicompactness} (2)}$)$. \hfill$\square$
\end{itemize}
\end{cor}

\subsubsection{Quasi-separated morphisms and coherent morphisms}\label{subsub-qcompqsepformaldef}
\index{morphism of formal schemes@morphism (of formal schemes)!quasi-compact morphism of formal schemes@quasi-compact ---|(}
\index{quasi-compact!quasi-compact morphism of formal schemes@--- morphism (of formal schemes)|(}
\index{morphism of formal schemes@morphism (of formal schemes)!quasi-separated morphism of formal schemes@quasi-separated ---|(}
\index{quasi-separated!quasi-separated morphism of formal schemes@--- morphism (of formal schemes)|(}
\begin{dfn}\label{dfn-qsepformal}{\rm 
A morphism $f\colon X\rightarrow Y$ of formal schemes is said to be {\em quasi-separated} if it satisfies one of the $($equivalent$)$ conditions in $\ref{cor-qsepformal001}$.
A formal scheme $X$ is said to be {\em quasi-separated}\index{formal scheme!quasi-separated formal scheme@quasi-separated ---} if it is quasi-separated over $\Spec\Z$.}
\end{dfn}

\begin{dfn}\label{dfn-cohformalmorphism}{\rm 
A quasi-compact and quasi-separated morphism $($resp.\ formal scheme$)$ is said to be {\em coherent}\index{morphism of formal schemes@morphism (of formal schemes)!coherent morphism of formal schemes@coherent ---}\index{coherent!coherent morphism of formal schemes@--- morphism (of formal schemes)}\index{formal scheme!coherent formal scheme@coherent ---}.}
\end{dfn}

For example, any affine formal scheme $X=\Spf A$ is coherent.
Indeed, it is quasi-compact, as it is a closed subset of the affine scheme $\Spec A$; since the diagonal map $\Delta\colon X\rightarrow X\times_{\Z}X$ comes from the surjective map $A\widehat{\otimes}_{\Z}A\rightarrow A$, one can easily show that it is quasi-compact (by an argument similar to that in \ref{adicmorphismdiagonal}, one can show that this map satisfies the assumption in Exercise \ref{exer-adicmappingformalschemesfiberproduct}).

\begin{prop}\label{prop-qsepmorformal3}
{\rm (1)} The composition of two quasi-compact $($resp.\ quasi-separated, resp.\ coherent$)$ morphisms of formal schemes is quasi-compact $($resp.\ quasi-separated, resp.\ coherent$)$.

{\rm (2)} If $f\colon X\rightarrow X'$ and $g\colon Y\rightarrow Y'$ are two quasi-compact $($resp.\ quasi-separated, resp.\ coherent$)$ $S$-morphisms of formal schemes, then $f\times_Sg\colon X\times_SY\rightarrow X'\times_SY'$ is quasi-compact $($resp.\ quasi-separated, resp.\ coherent$)$.

{\rm (3)} If $f\colon X\rightarrow Y$ is a quasi-compact $($resp.\ quasi-separated, resp.\ coherent$)$ $S$-morphism of formal schemes and $S'\rightarrow S$ is a morphism, then $f_{S'}\colon X\times_SS'\rightarrow Y\times_SS'$ is quasi-compact $($resp.\ quasi-separated, resp.\ coherent$)$.

{\rm (4)} If the composition $g\circ f$ of two morphisms of formal schemes is quasi-compact and $f$ is surjective, then $g$ is quasi-compact$;$ if $g\circ f$ is quasi-compact and $g$ is quasi-separated, then $f$ is quasi-compact$;$ if $g\circ f$ is quasi-separated, then $f$ is quasi-separated$;$ if, moreover, $f$ is quasi-compact and surjective, then $g$ is quasi-separated.
\end{prop}

\begin{proof}
First let us prove the assertions (1), (2), and (3) for quasi-compactness.
(1) is clear.
As we have seen in \S\ref{subsub-genpropertymorphism}, (2) and (3) follow from the special case of (3) with $S=Y$.
Let $f\colon X\rightarrow Y$ be quasi-compact, and $g\colon Y'\rightarrow Y$ any morphism.
We set $f'\colon X'\rightarrow Y'$ to be the base change; we want to show that $f'$ is quasi-compact.
Let $U'\subseteq Y'$ be a quasi-compact open subset, and $s'\in U'$ a point.
Take an affine neighborhood $T$ in $Y$ of the point $g(s')$ and an affine neighborhood $W$ in $U'\cap g^{-1}(T)$ of $s'$.
Since it suffices to show that $f^{-1}(W)$ is quasi-compact, we may assume $Y$ and $Y'$ are affine; moreover, it is enough only to show that the formal scheme $X'=X\times_YY'$ is quasi-compact.
Since $X$ is quasi-compact, it is covered by finitely many affine open sets $V_1,\ldots,V_n$.
Then $X'$ is covered by the affine open subsets $V_i\times_YY'$ (\ref{prop-openimmform1} (3)), and thus $X'$ is quasi-compact, as desired.

Now we proceed to show the assertions (1), (2), and (3) for quasi-separatedness. 
Let $f\colon X\rightarrow Y$ and $g\colon Y\rightarrow Z$ be quasi-separated morphisms of formal schemes.
We have the following commutative diagram with the Cartesian square:
$$
\xymatrix{X\ar[r]_(.4){\Delta_f}\ar@/^1pc/[rr]^{\Delta_{g\circ f}}\ar@/_/[dr]_f&X\times_YX\ar[r]\ar[d]_q&X\times_ZX\ar[d]^{f\times_Zf}\\ &Y\ar[r]_(.43){\Delta_g}&Y\times_ZY.}\leqno{(\ast)}
$$
Since $\Delta_g$ is quasi-compact, the upper arrow in the square is quasi-compact.
Since $\Delta_f$ is quasi-compact, $\Delta_{g\circ f}$ is quasi-compact, which implies that $g\circ f$ is quasi-separated.
Thus (1) for quasi-separatedness is proved.

Similarly to the quasi-compactness case, the assertions (2) and (3) follow from the special case of (3) with $S=Y$.
Let $f\colon X\rightarrow Y$ be quasi-separated, $g\colon Y'\rightarrow Y$ a morphism, and $f'\colon X'\rightarrow Y'$ the base change.
We have the following commutative diagram:
$$
\xymatrix{X'\times_{Y'}X'\ar@{=}[r]&(X\times_YX)\times_YY'\ar[r]&X\times_YX\\ X'\ar[u]^{\Delta_{f'}}\ar@{=}[r]&X\times_YY'\ar[u]_{\Delta_f\times_YY'}\ar[r]&X.\ar[u]_{\Delta_f}}\leqno{(\ast\ast)}
$$
Since $\Delta_f$ is quasi-compact, $\Delta_{f'}=\Delta_f\times_YY'$ is quasi-compact by (3) for quasi-compactness, which we have already shown; that is, $f'$ is quasi-separated, as desired.

Finally let us prove (4).
The first assertion is easy to see; the second one can be shown similarly to \cite[(6.1.5) (v)]{EGAInew}.
Let $f\colon X\rightarrow Y$ and $g\colon Y\rightarrow Z$ be such that $g\circ f$ is quasi-separated.
The morphism $f$ coincides with the composition
$$
X\stackrel{\Gamma_f}{\longrightarrow}X\times_ZY\stackrel{p_2}{\longrightarrow}Y,
$$
where $\Gamma_f$ is the graph of $f$.
The projection $p_2$ coincides with $(g\circ f)\times_Z\id_Y$ and hence is quasi-separated.
The diagonal morphism of the graph $\Gamma_f$ is isomorphic to $\id_X$, and hence $\Gamma_f$ is quasi-separated.
Thus $f=p_2\circ\Gamma_f$ is quasi-separated.
Suppose that $g\circ f$ is quasi-separated and that $f$ is quasi-compact and surjective.
In the diagram $(\ast)$ above we know that the maps $\Delta_{g\circ f}$ and $f\times_Zf$ are quasi-compact.
Hence the composition $\Delta_g\circ q\circ\Delta_f=\Delta_g\circ f$ is quasi-compact.
Since $f$ is surjective, we deduce that $\Delta_g$ is quasi-compact, that is, $g$ is quasi-separated.
\end{proof}

\begin{prop}\label{prop-qsepmorformal0}
Let $f\colon X\rightarrow Y$ be a morphism of formal schemes, and $Y=\bigcup_{\alpha\in L}V_{\alpha}$ an open covering of $Y$.
Then $f$ is quasi-compact $($resp.\ quasi-separated, resp.\ coherent$)$ if and only if the induced map $f_{\alpha}\colon X_{\alpha}=X\times_YV_{\alpha}\rightarrow V_{\alpha}$ is quasi-compact $($resp.\ quasi-separated, resp.\ coherent$)$ for any $\alpha\in L$.
\end{prop}

\begin{proof}
The `only if' part follows from \ref{prop-qsepmorformal3} (3).
Suppose $f_{\alpha}$ is quasi-compact for any $\alpha\in L$, and let $V\subset Y$ be a quasi-compact open subset of $Y$.
Then $V$ is covered by finitely many quasi-compact open subsets each of which is contained in some $V_{\alpha}$.
Thus $f^{-1}(V)$ is covered by finitely many quasi-compact open subsets and hence is quasi-compact.

Next, suppose $f_{\alpha}$ is quasi-separated for any $\alpha\in L$.
Since the diagonal map 
$$
\Delta_{X\times_YV_{\alpha}}\colon X\times_YV_{\alpha}\longrightarrow(X\times_YV_{\alpha})\times_{V_{\alpha}}(X\times_YV_{\alpha})\cong(X\times_YX)\times_YV_{\alpha}
$$
coincides with the base change of $\Delta_X$ by the open immersion $(X\times_YX)\times_YV_{\alpha}\hookrightarrow X\times_YX$ (\ref{lem-sepqsepformal2}), the assertion follows from the assertion for the quasi-compactness.
\end{proof}

\begin{prop}\label{prop-qsepformal1}
Let $f\colon X\rightarrow Y$ be an adic morphism of adic formal schemes of finite ideal type, and suppose $Y$ has an ideal of definition $\mathscr{I}$ of finite type.
For any $k\geq 0$ set $X_k=(X,\O_X/\mathscr{I}^{k+1}\O_X)$ and $Y_k=(Y,\O_Y/\mathscr{I}^{k+1})$, and let $f_k\colon X_k\rightarrow Y_k$ be the induced morphism of schemes.
Then the following conditions are equivalent$:$
\begin{itemize}
\item[{\rm (a)}] $f$ is quasi-compact $($resp.\ quasi-separated, resp.\ coherent$);$
\item[{\rm (b)}] $f_k$ is quasi-compact $($resp.\ quasi-separated, resp.\ coherent$)$ for any $k\geq 0;$
\item[{\rm (c)}] $f_0$ is quasi-compact $($resp.\ quasi-separated, resp.\ coherent$)$.
\end{itemize}
\end{prop}

\begin{proof}
The assertion for quasi-compact morphisms is clear.
As for quasi-separated morphisms, the assertion follows easily from \ref{cor-prodformal21}.
\end{proof}

By \ref{prop-genpropertymorphismadic2} we have the following:
\begin{cor}\label{cor-qsepformal1}
Let $f\colon X\rightarrow Y$ be a morphism of schemes, and $Z$ a closed subscheme of $Y$ of finite presentation.
If $f$ is quasi-compact $($resp.\ quasi-separated, resp.\ coherent$)$, then the formal completion {\rm (\S\ref{subsub-formalcompletionsch})} $\widehat{f}\colon\widehat{X}|_{f^{-1}(Z)}\rightarrow\widehat{Y}|_Z$ is quasi-compact $($resp.\ quasi-separated, resp.\ coherent$)$.
\end{cor}

\subsubsection{Notation}\label{subsub-notationcategorycoherentformaschemes}
In the sequel we write
\begin{itemize}
\item $\CFs_S=$ the category of formal schemes coherent over $S$.
\end{itemize}
We employ the similar rule as in \S\ref{subsub-categoryformalschnotation} for the notations of categories of coherent formal schemes; for example:
\begin{itemize}
\item $\Ac\CFs_S=$ the category of adic formal schemes of finite ideal type coherent over $S$;
\item $\Ac\CFs^{\ast}_S=$ the category of adic formal schemes of finite ideal type coherent over $S$ with adic $S$-morphisms;
\item $\Ac\CFs^{\ast}_{/S}=$ the category of adic formal schemes of finite ideal type coherent and adic over $S$;
\end{itemize}
When $S=\Spec\Z$, we denote these categories without reference to $S$; e.g.\ $\CFs=\CFs_{\Spec\Z}$.
\index{quasi-separated!quasi-separated morphism of formal schemes@--- morphism (of formal schemes)|)}
\index{morphism of formal schemes@morphism (of formal schemes)!quasi-separated morphism of formal schemes@quasi-separated ---|)}
\index{quasi-compact!quasi-compact morphism of formal schemes@--- morphism (of formal schemes)|)}
\index{morphism of formal schemes@morphism (of formal schemes)!quasi-compact morphism of formal schemes@quasi-compact ---|)}

\subsection{Morphisms of finite type}\label{sub-finitypeformal}
\begin{dfn}\label{dfn-topfintype}{\rm 
A morphism $f\colon X\rightarrow Y$ of adic formal schemes of finite ideal type is said to be {\em locally of finite type}\index{morphism of formal schemes@morphism (of formal schemes)!morphism of formal schemes locally of finite type@--- locally of finite type} if the following conditions are satisfied:
\begin{itemize}
\item[{\rm (a)}] the morphism $f$ is adic {\rm (\ref{dfn-adicmor})}; 
\item[{\rm (b)}] there exist an affine open covering $\{V_{\alpha}\}_{\alpha\in L}$ of $Y$ with $V_{\alpha}=\Spf B_{\alpha}$, where each $B_{\alpha}$ is an adic ring of finite ideal type, and for each $\alpha\in L$ an affine open covering $\{U_{\alpha,\lambda}\}_{\lambda\in\Lambda_{\alpha}}$ of $f^{-1}(V_{\alpha})$ with $U_{\alpha,\lambda}=\Spf A_{\alpha,\lambda}$, where each $A_{\alpha,\lambda}$ is an adic ring topologically finitely generated\index{finitely generated!topologically finitely generated@topologically ---} over $B_{\alpha}$ {\rm ({\rm {\bf \ref{ch-pre}}.\ref{dfn-topfinigen}})}.
\end{itemize}
The morphism $f$ is said to be {\em of finite type}\index{morphism of formal schemes@morphism (of formal schemes)!morphism of formal schemes of finite type@--- of finite type} if it is locally of finite type and quasi-compact (\ref{dfn-qcompformal}).}
\end{dfn}

\begin{prop}\label{prop-topfintype2}
{\rm (1)} An open immersion\index{immersion!open immersion of formal schemes@open --- (of formal schemes)} of adic formal schemes of finite ideal type is locally of finite type.

{\rm (2)} The composition of two morphisms locally of finite type $($resp.\ of finite type$)$ is again locally of finite type $($resp.\ of finite type$)$.
If the composition $g\circ f$ of morphisms $f\colon X\rightarrow Y$ and $g\colon Y\rightarrow Z$ of formal schemes is locally of finite type and $g$ is adic, then $f$ is locally of finite type.

{\rm (3)} Let $S$ be an adic formal scheme of finite ideal type, and $f\colon X\rightarrow X'$ and $g\colon Y\rightarrow Y'$ two adic $S$-morphisms of adic formal schemes of finite ideal type.
Suppose that $f$ and $g$ are locally of finite type $($resp.\ of finite type$)$.
Then $f\times_Sg\colon X\times_SY\rightarrow X'\times_SY'$ is locally of finite type $($resp.\ of finite type$)$.

{\rm (4)} Let $S$ be an adic formal scheme of finite ideal type, and $f\colon X\rightarrow Y$ an adic $S$-morphism between adic formal schemes of finite ideal type. 
Suppose $f$ is locally of finite type $($resp.\ of finite type$)$.
Then for any morphism $S'\rightarrow S$ of adic formal schemes of finite ideal type, $f_{S'}\colon X\times_SS'\rightarrow Y\times_SS'$ is locally of finite type $($resp.\ of finite type$)$.
\end{prop}

\begin{proof} The assertions (1) and (2) are easy to see. 
As we have seen in \S\ref{subsub-genpropertymorphismadic}, the assertions (3) and (4) follow from the special case of (4) with $S=Y$.
To show (4) in this case, we may assume that all formal schemes are affine of the following form: $f\colon X=\Spf B\dl X_1,\ldots,X_n\dr\rightarrow Y=\Spf B$ and $Y'=S'=\Spf R\rightarrow Y$.
Let $I$ be a finitely generated ideal of definition of $B$, and $J$ a finitely generated ideal of definition of $R$ such that $IR\subseteq J$.
Then what to prove is that the adic ring $B\dl X_1,\ldots,X_n\dr\widehat{\otimes}_BR$ is topologically of finitely generated over $R$.
But it is easy to see that the admissible ring in question is isomorphic to $R\dl X_1,\ldots,X_n\dr$, since it is the $J$-adic completion of $B[X_1,\ldots,X_n]\otimes_BR=R[X_1,\ldots,X_n]$.
\end{proof}

If $X$ and $Y$ are locally Noetherian, it can be shown that our definition of `of finite type' agrees with that of \cite[$\mathbf{I}$, (10.13.1)]{EGA} due to \cite[$\mathbf{I}$, (10.13.4)]{EGA}; that is, `of finite type' means that the induced morphism $(X,\O_X/(f^{-1}\mathscr{J})\O_X)\rightarrow (Y,\O_Y/\mathscr{J})$ (where $\mathscr{J}$ is an ideal of definition of $Y$) is of finite type.
In fact, we have, more generally, the following:
\begin{prop}\label{prop-cortopfintype11}
Let $f\colon X\rightarrow Y$ be an adic morphism of adic formal schemes of finite ideal type.
Suppose that $Y$ has an ideal of definition $\mathscr{I}$ of finite type.
Set $X_k=(X,\O_X/\mathscr{I}^{k+1}\O_X)$ and $Y_k=(Y,\O_Y/\mathscr{I}^{k+1})$ for $k\geq 0$, and denote by $f_k\colon X_k\rightarrow Y_k$ the induced morphism of schemes.
Then the following conditions are equivalent$:$
\begin{itemize}
\item[{\rm (a)}] $f$ is locally of finite type $($resp.\ of finite type$);$
\item[{\rm (b)}] $f_k$ is locally of finite type $($resp.\ of finite type$)$ for $k\geq 0;$
\item[{\rm (c)}] $f_0$ is locally of finite type $($resp.\ of finite type$)$.
\end{itemize}
\end{prop}

The proposition follows immediately from the following lemma:
\begin{lem}\label{lem-proptopfintype1}
Let $A$ be an adic ring with a finitely generated ideal of definition $I$, and $B$ an $IB$-adically complete $A$-algebra.
Then the following conditions are equivalent$:$
\begin{itemize}
\item[{\rm (a)}] the morphism $\Spf B\rightarrow\Spf A$ is of finite type$;$
\item[{\rm (b)}] the morphism $\Spec B/IB\rightarrow\Spec A/I$ is of finite type$;$
\item[{\rm (c)}] $B$ is topologically finitely generated over $A$.
\end{itemize}
\end{lem}

\begin{proof}
The implication (c) $\Rightarrow$ (b) is clear; the converse follows from \cite[$\mathbf{I}$, (6.3.3)]{EGA} and {\bf \ref{ch-pre}}.\ref{prop-formalnot3}.
The implication (c) $\Rightarrow$ (a) is obvious. 
Since for any topologically finitely generated $A$-algebra $C$, $C/IC$ is an $(A/I)$-algebra of finite type, (a) $\Rightarrow$ (b) follows.
\end{proof}

\begin{cor}\label{cor-topfintype4}
Let $f\colon X\rightarrow Y$ be a morphism of schemes, and $Z$ a closed subscheme of $Y$ of finite presentation.
If $f$ is locally of finite type $($resp.\ of finite type$)$, then the formal completion $\widehat{f}\colon\widehat{X}|_{f^{-1}(Z)}\rightarrow\widehat{Y}|_Z$ is locally of finite type $($resp.\ of finite type$)$. \hfill$\square$
\end{cor}

\addcontentsline{toc}{subsection}{Exercises}
\subsection*{Exercises}
\begin{exer}\label{exer-adicmappingformalschemesfiberproduct}{\rm 
Let $\varphi\colon A\rightarrow B$ be a continuous homomorphism between admissible rings.
Suppose there exists an ideal of definition $I\subseteq A$ such that $IB$ is an ideal of definition of $B$ (e.g.\ $A$ and $B$ are adic rings of finite ideal type and $\varphi$ is an adic morphism).
Then show that the square
$$
\xymatrix{\Spf B\,\ar@{^{(}->}[r]\ar[d]&\Spec B\ar[d]\\ \Spf A\,\ar@{^{(}->}[r]&\Spec A}
$$
is Cartesian in the category of topological spaces.}
\end{exer}

\begin{exer}\label{exer-disjointsumformalsch}{\rm 
Let $\{X_{\lambda}\}_{\lambda\in\Lambda}$ be a collection of formal schemes.
Show that the functor $\Fs\rightarrow\Sets$ by $Y\mapsto\prod_{\lambda\in\Lambda}\Hom_{\Fs}(X_{\lambda},Y)$ is representable; in other words, the disjoint sum\index{disjoint sum} $\coprod_{\lambda\in\Lambda}X_{\lambda}$ exists in the category $\Fs$.
Show, moreover, that if each $X_{\lambda}$ is adic, then $\coprod_{\lambda\in\Lambda}X_{\lambda}$ is adic.}
\end{exer}

\begin{exer}\label{exer-idealofdefinitionintersection}{\rm 
Let $X$ be a formal scheme.
Show that, if $\mathscr{I},\mathscr{I}'\subseteq\O_X$ are ideals of definition of $X$, then $\mathscr{I}+\mathscr{I}'$ are $\mathscr{I}\cap\mathscr{I}'$ are ideals of definition.}
\end{exer}

\begin{exer}\label{exer-adicformalschemesidealsofdefinitionfinitude}{\rm 
Let $X$ be an adic formal scheme, and $\mathscr{I}$ an ideal of definition of $X$.
Suppose that $\mathscr{I}/\mathscr{I}^2$ is a quasi-coherent ideal sheaf of finite type on the scheme $(X,\mathscr{O}_X/\mathscr{I})$. 
Then show that $\mathscr{I}$ is of finite type.}
\end{exer}

\begin{exer}\label{exer-disjointsumproduct}{\rm 
Let $\{X_{\lambda}\}_{\lambda\in\Lambda}$ and $\{Y_{\mu}\}_{\mu\in M}$ be collections of formal schemes over a fixed formal scheme $S$, and $X=\coprod_{\lambda\in\Lambda}X_{\lambda}$ and $Y=\coprod_{\mu\in M}Y_{\mu}$ the respective disjoint sums\index{disjoint sum}.
Show that the fiber product $X\times_SY$ is canonically isomorphic to the disjoint sum of $\{X_{\lambda}\times_SY_{\mu}\}_{(\lambda,\mu)\in\Lambda\times M}$.}
\end{exer}

\begin{exer}\label{exer-adqformalprop1}{\rm 
Let $(A,I)$ be a pair of finite ideal type, and consider the canonical maps $j\colon A\rightarrow\widehat{A}$ and $i\colon\Spf\widehat{A}\rightarrow\Spec A$.
Show that for any $f\in\widehat{A}$ there exists $g\in A$ such that $\mathfrak{D}(f)=\mathfrak{D}(j(g))=i^{-1}D(g)$.}
\end{exer}

\begin{exer}\label{exer-formalindlimsch}{\rm 
Let $S$ be an adic formal scheme with an ideal of definition $\mathscr{I}$ of finite type, and consider the category $\Ac\Fs^{\ast}_{/S}$ of adic formal schemes of finite ideal type adic over $S$.
For any object $X$ of $\Ac\Fs^{\ast}_{/S}$ and any integer $k\geq 0$, we set $X_k=(X,\O_X/\mathscr{I}^{k+1}\O_X)=X\times_SS_k$ and define the functor $\Ac\Fs^{\ast}_{/S}\rightarrow\Sch_{S_k}$ by $X\mapsto X_k$.
Varying $k$, we get the canonical functor
$$
\Ac\Fs^{\ast}_{/S}\longrightarrow\varprojLim\Sch_{S_k},
$$
where the latter category is the $2$-categorical limit of the categories $\Sch_{S_k}$ for $k\geq 0$.
Show that this functor is fully faithful, that is, for any formal schemes $X,Y$ adic over $S$ the canonical map
$$
\Hom_{\Ac\Fs^{\ast}_{/S}}(X,Y)\longrightarrow\varprojlim_k\Hom_{S_k}(X_k,Y_k)
$$
is a bijection.}
\end{exer}

\begin{exer}\label{exer-adicformalschemesinductivelimitadicmaps}{\rm 
Let $Y$ be an adic formal scheme of finite ideal type, and $\mathscr{I}$ an ideal of definition of finite type of $Y$.
Let $Y_k=(Y,\O_Y/\mathscr{I}^{k+1})$ for $k\geq 0$.
Suppose we are given an inductive system of schemes $\{X_k\}_{k\geq 0}$ over $Y$ such that
\begin{itemize}
\item for $k\leq l$ the map $X_k\rightarrow X_l$ is a closed immersion whose underlying continuous mapping is a homeomorphism;
\item for $k\leq l$ the kernel of $\O_{X_l}\rightarrow\O_{X_k}$ coincides with $\mathscr{I}^{k+1}\O_{X_l}$.
\end{itemize}
Then show that $X=\varinjlim_{k\geq 0}X_k$ is an adic formal schemes of finite ideal type and that the morphism $X\rightarrow Y$ is adic, that is, $\mathscr{I}\O_X$ is an ideal of definition of finite type of $X$.}
\end{exer}

\begin{exer}\label{exer-affinecovaluedpoint}{\rm 
Let $X$ be a formal scheme, $A$ an admissible ring, and $\{\mathscr{I}^{(\lambda)}\}_{\lambda\in\Lambda}$ a fundamental system of ideals of definition on $X$.
Consider on the ring $\Gamma(X,\O_X)$ the topology induced from the filtration $\{\Gamma(X,\mathscr{I}^{(\lambda)})\}_{\lambda\in\Lambda}$ by ideals.
Show that there exists a canonical bijection between the set of all continuous homomorphisms $A\rightarrow\Gamma(X,\O_X)$ and the set of all morphisms $X\rightarrow\Spf A$ of formal schemes (cf.\ \cite[$\mathbf{I}$, (2.2.4)]{EGA}).}
\end{exer}


\section{Universally rigid-Noetherian formal schemes}\label{sec-adequateformalschemes}
In this section, we introduce two new classes of adic formal schemes, the so-called {\em $($locally$)$ universally rigid-Noetherian formal schemes}\index{formal scheme!universally rigid-Noetherian formal scheme@universally rigid-Noetherian ---}\index{formal scheme!universally rigid-Noetherian formal scheme@universally rigid-Noetherian ---!locally universally rigid-Noetherian formal scheme@locally --- ---} and {\em $($locally $)$ universally adhesive formal schemes}\index{formal scheme!universally adhesive formal scheme@universally adhesive ---}\index{formal scheme!universally adhesive formal scheme@universally adhesive ---!locally universally adhesive formal scheme@locally --- ---}\index{adhesive!universally adhesive@universally ---!universally adhesive formal scheme@--- --- formal scheme}.
The former is based on the ring-theoretic notion `topologically universally Noetherian outside $I$'\index{Noetherian!Noetherian outside I@--- outside $I$!topologically universally Noetherian outside I@topologically universally --- ---} ({\bf \ref{ch-pre}}.\ref{dfn-topologicallyuniversallynoetherian}), and the latter on `topologically universally adhesive'\index{adhesive!universally adhesive@universally ---!topologically universally adhesive@topologically --- --- (t.u.\ adhesive)} ({\bf \ref{ch-pre}}.\ref{dfn-topadhesive}).
Hence the notion of (locally) universally adhesive formal schemes is more restrictive than that of (locally) universally rigid-Noetherian formal schemes.
As we will see later in subsequent sections, these kinds of formal schemes enjoy many of the nice properties that locally Noetherian formal schemes possess and thus provide a good generalization of the notion of locally Noetherian formal schemes, which is often too restrictive for developing general rigid geometry.

First in \S\ref{subsub-tuarings} we introduce the ring-theoretic notions, t.u.\ ($=$ topologically universally) rigid-Noetherian rings and t.u.\ adhesive rings, before giving the definition in \S\ref{subsub-adequateformalschadmissible} of (locally) universally rigid-Noetherian formal schemes and universally adhesive formal schemes\index{formal scheme!universally rigid-Noetherian formal scheme@universally rigid-Noetherian ---}\index{formal scheme!universally rigid-Noetherian formal scheme@universally rigid-Noetherian ---!locally universally rigid-Noetherian formal scheme@locally --- ---}.
The notion of locally universally rigid-Noetherian formal schemes allows one to define locally of finite presentation morphisms, which we discuss in \S\ref{sub-finipresformal}.

It turns out that the category of (locally) universally adhesive formal schemes contains the category of the so-called {\em admissible formal schemes}\index{formal scheme!admissible formal scheme@admissible ---}\index{admissible!admissible formal scheme@--- formal scheme}, which play a central role in the classical rigid geometry.
We briefly recall the definition of admissible formal schemes in \S\ref{subsub-admissibleformalschemes}, and overview the interrelations of the several notions of formal schemes so far obtained in \S\ref{subsub-catout}.

\subsection{Universally rigid-Noetherian and universally adhesive formal schemes}\label{sub-formalsch}
\subsubsection{T.u.\ rigid-Noetherian rings and t.u.\ adhesive\ rings}\label{subsub-tuarings}
\index{t.u. rigid-Noetherian ring@t.u.\ rigid-Noetherian ring|(}
\index{t.u.a. ring@t.u.\ adhesive ring|(}
\begin{dfn}\label{dfn-tuaringadmissible}{\rm 
Let $A$ be an adic ring of finite ideal type, and $I\subseteq A$ an ideal of definition.

{\rm (1)} The topological ring $A$ is called a {\em rigid-Noetherian ring}\index{rigid-Noetherian ring@rigid-Noetherian ring} if it is Noetherian outside $I$\index{Noetherian!Noetherian outside I@--- outside $I$} {\rm ({\bf \ref{ch-pre}}.\ref{dfn-outsideI} (1))}; if it is, furthermore, topologically universally Noetherian outside $I$\index{Noetherian!Noetherian outside I@--- outside $I$!topologically universally Noetherian outside I@topologically universally --- ---} ({\bf \ref{ch-pre}}.\ref{dfn-topologicallyuniversallynoetherian}), then it is called a {\em topologically universally} (abbreviated as {\em t.u.}) {\em rigid-Noetherian} ring\index{t.u. rigid-Noetherian ring@t.u.\ rigid-Noetherian ring}.

{\rm (2)} The ring $A$ is called a {\em t.u.\ adhesive ring}\index{t.u.a. ring@t.u.\ adhesive ring} if it is $I$-adically topologically universally adhesive\index{adhesive!universally adhesive@universally ---!topologically universally adhesive@topologically --- --- (t.u.\ adhesive)} ({\bf \ref{ch-pre}}.\ref{dfn-topadhesive}).}
\end{dfn}

It follows immediately from the definition that, if $A$ is t.u.\ adhesive (resp.\ t.u.\ rigid-Noetherian), then any topologically finitely generated\index{finitely generated!topologically finitely generated@topologically ---} $A$-algebra is again t.u.\ adhesive (resp.\ t.u.\ rigid-Noetherian).
Since adhesiveness implies Noetherian outside $I$, it follows that `t.u.\ adhesive' implies `t.u.\ rigid-Noetherian':
$$
\textrm{t.u.\ adhesive}\quad\Longrightarrow\quad\textrm{t.u.\ rigid-Noetherian}\quad\Longrightarrow\quad\textrm{rigid-Noetherian}.
$$
If $A$ is a t.u.\ adhesive ring and $I\subseteq A$ is an ideal of definition, then the pair\index{pair} $(A,I)$ is a complete t.u.\ adhesive pair\index{adhesive!universally adhesive@universally ---!topologically universally adhesive@topologically --- --- (t.u.\ adhesive)}\index{pair!adhesive pair@adhesive ---!topologically universally adhesive pair@topologically universally --- ---}.
Note that, due to {\bf \ref{ch-pre}}.\ref{thm-gabberIHES2008}, rigid-Noetherian rings are pseudo-adhesive\index{pair!adhesive pair@adhesive ---!pseudo adhesive pair@pseudo-{---} ---}\index{adhesive!adhesive pair@--- pair!pseudo adhesive pair@pseudo-{---} ---}\index{pseudo-adhesive} ({\bf \ref{ch-pre}}.\ref{dfn-adhesive}).
By \ref{lem-formalnot0} and {\bf \ref{ch-pre}}.\ref{prop-btarfflat1} we readily see: 
\begin{prop}\label{prop-rigidnoetherianflatness}
Let $A$ be a rigid-Noetherian ring.\index{rigid-Noetherian ring@rigid-Noetherian ring}
Then for any $f\in A$ the $I$-adic completion map $A_f\rightarrow A_{\{f\}}$ is flat.
In particular, the map $\Spf A\rightarrow\Spec A$ of locally ringed spaces is flat $($cf.\ {\rm {\bf \ref{ch-pre}}.\ref{dfn-flatnessringedsp} (2)}$)$.
\end{prop}

\begin{rem}\label{rem-turigidnoetherianbasicproperties}{\rm 
By {\bf \ref{ch-pre}}.\ref{thm-btarf1} if $A$ is t.u.\ rigid-Noetherian ring, then the pair $(A,I)$ is a complete t.u.\ pseudo-adhesive pair\index{pair!adhesive pair@adhesive ---!topologically universally pseudo adhesive pair@topologically universally pseudo-{---} ---}\index{adhesive!universally pseudo adhesive@universally pseudo-{---}!topologically universally pseudo adhesive@topologically --- --- (t.u.\ pseudo-adhesive)}\index{pseudo-adhesive!universally pseudo adhesive@universally ---!topologically universally pseudo adhesive@topologically --- --- (t.u.\ pseudo-adhesive)} ({\bf \ref{ch-pre}}.\ref{dfn-topadhesive}).
Hence, in particular, a t.u.\ rigid-Noetherian ring $A$ enjoys the following property: 
\begin{itemize}
\item for any $n\geq 0$ and $m\geq 0$ the $A$-algebra of the form 
$$
A\dl X_1,\ldots,X_n\dr[Y_1,\ldots,Y_m],
$$
together with the ideal of definition $IA\dl X_1,\ldots,X_n\dr[Y_1,\ldots,Y_m]$, satisfies the conditions {\bf (BT)} in {\bf \ref{ch-pre}}, \S\ref{subsub-BTfirstprop} and {\bf (AP)} in {\bf \ref{ch-pre}}, \S\ref{subsub-ARIgoodnesssub}.
\end{itemize}
Then by {\bf \ref{ch-pre}}.\ref{prop-btarf1} we have:
\begin{itemize}
\item if $A$ is t.u.\ rigid-Noetherian, then any finitely generated $A$-module $M$ is $I$-adically complete and, moreover, any $A$-submodule $N\subseteq M$ is closed in $M$ and $I$-adically complete ({\bf \ref{ch-pre}}.\ref{cor-propARconseq1-2});
\item if $B$ is an $A$-algebra of finite type, then the map $B\rightarrow\widehat{B}$ by $I$-adic completion is flat.
\end{itemize}}
\end{rem}

\begin{exas}\label{exas-tuaringsadmissible}{\rm 
(1) Any Noetherian adic ring is t.u.\ adhesive. 
Indeed, if $A$ is a Noetherian adic ring with an ideal of definition $I\subseteq A$, then it is $I$-adically universally adhesive\index{adhesive!universally adhesive@universally ---} ({\bf \ref{ch-pre}}.\ref{exa-adhesivetypeN}); if $B$ is an $A$-algebra of finite type (hence Noetherian), then the $I$-adic completion of $B$ is again Noetherian and hence is $I$-adically universally adhesive.

(2) Let $V$ be an $a$-adically complete valuation ring (of arbitrary height)\index{valuation!valuation ring@--- ring!a-adically complete valuation ring@$a$-adically complete --- ---} where $a\in\m_V\setminus\{0\}$.
By Gabber's\index{Gabber, O.} theorem ({\bf \ref{ch-pre}}.\ref{cor-convadh}) the topological ring $V$ is t.u.\ adhesive.
Hence any topologically finitely generated $V$-algebra, that is, a topological algebra over type (V)\index{algebra!topological algebra of type (V)@topological --- of type (V)} (cf.\ {\bf \ref{ch-pre}}, \S\ref{sec-aadicallycompval}), is t.u.\ adhesive.}
\end{exas}

Using \ref{prop-rigidnoetherianflatness}, one can prove the following proposition similarly to \ref{prop-adqformalprop}:
\begin{prop}\label{prop-adqformalpropaffine}
Let $A$ be a rigid-Noetherian ring, and $I\subseteq A$ a finitely generated ideal of definition.
Then for $X=\Spec A$, $Y=\Spec A/I$, and $\widehat{X}=\Spf A$, the assertions {\rm (1)} and {\rm (2)} in {\rm \ref{prop-adqformalprop}} are true. \hfill$\square$
\end{prop}

Let us finally mention that the properties `t.u.\ rigid-Noetherian' and `t.u.\ adhesive' have the formal fpqc patching principle:
\begin{prop}\label{prop-formalfpqcpatchingtua}
Let $A\rightarrow B$ be an adic morphism between adic rings of finite ideal type, and $I\subseteq A$ a finitely generated ideal of definition of $A$.
Suppose that for any $k\geq 0$ the induced map $A/I^{k+1}\rightarrow B/I^{k+1}B$ is faithfully flat.
If $B$ is t.u.\ rigid-Noetherian $($resp.\ t.u.\ adhesive$)$, then so is $A$.
Moreover, in this situation, the map $A\dl X_1,\ldots,X_n\dr\rightarrow B\dl X_1,\ldots,X_n\dr$ is faithfully flat for any $n\geq 0$.
\end{prop}

\begin{proof}
By {\bf \ref{ch-pre}}.\ref{prop-fpqcdescentrigidNoetherian} $A\dl X_1,\ldots,X_n\dr$ is Noetherian outside $IA\dl X_1,\ldots,X_n\dr$ and the map $A\dl X_1,\ldots,X_n\dr\rightarrow B\dl X_1,\ldots,X_n\dr$ is faithfully flat.
By {\bf \ref{ch-pre}}.\ref{prop-adhesive0} (2) if $B\dl X_1,\ldots,X_n\dr$ is universally adhesive\index{pair!adhesive pair@adhesive ---!universally adhesive pair@universally --- ---}\index{adhesive!universally adhesive@universally ---}, then so is $A\dl X_1,\ldots,X_n\dr$.
\end{proof}
\index{t.u.a. ring@t.u.\ adhesive ring|)}
\index{t.u. rigid-Noetherian ring@t.u.\ rigid-Noetherian ring|)}

\subsubsection{Universally adhesive and universally rigid-Noetherian formal schemes}\label{subsub-adequateformalschadmissible}
\index{formal scheme!universally adhesive formal scheme@universally adhesive ---|(}\index{adhesive!universally adhesive@universally ---!universally adhesive formal scheme@--- --- formal scheme|(}\index{formal scheme!universally rigid-Noetherian formal scheme@universally rigid-Noetherian ---|(}
\index{formal scheme!universally rigid-Noetherian formal scheme@universally rigid-Noetherian ---!locally universally rigid-Noetherian formal scheme@locally --- ---|(}
\index{formal scheme!universally adhesive formal scheme@universally adhesive ---!locally universally adhesive formal scheme@locally --- ---|(}
\begin{dfn}\label{dfn-formalsch}{\rm 
A formal scheme $X$ is said to be {\em locally universally rigid-Noetherian} (resp.\ {\em locally universally adhesive}) if there exists an affine open covering $X=\bigcup_{\alpha\in L}U_{\alpha}$ such that each $U_{\alpha}$ is isomorphic to $\Spf A_{\alpha}$ with $A_{\alpha}$ t.u.\ rigid-Noetherian (resp.\ t.u.\ adhesive).
If $X$ is, moreover, quasi-compact, we say that $X$ is {\em universally rigid-Noetherian} (resp.\ {universally adhesive}).}
\end{dfn}

Locally universally rigid-Noetherian (resp.\ locally universally adhesive) formal schemes are, by definition, adic\index{formal scheme!adic formal scheme@adic ---}\index{adic!adic formal scheme@--- formal scheme} of finite ideal type\index{formal scheme!adic formal scheme@adic ---!adic formal scheme of finite ideal type@--- --- of finite ideal type}\index{adic!adic formal scheme@--- formal scheme!adic formal scheme of finite ideal type@--- --- of finite ideal type} (\ref{dfn-formalschemesadicformalschemes}, \ref{dfn-adicformalschemesoffiniteidealtype}).
Notice also that locally universally adhesive formal schemes are locally universally rigid-Noetherian.
Since the properties `t.u.\ adhesive' and `t.u.\ right-Noetherian' are closed under topologically finitely generated extension, we have (cf.\ \ref{prop-adicformalschemesoffiniteidealtype}):
\begin{prop}\label{prop-adequateformalscheme1}
Let $X$ be an locally universally rigid-Noetherian $($resp.\ locally universally adhesive$)$ formal scheme, and $X'\rightarrow X$ a locally of finite type morphism of adic formal schemes of finite ideal type.
Then $X'$ is locally universally rigid-Noetherian $($resp.\ locally universally adhesive$)$. \hfill$\square$
\end{prop}

\begin{prop}\label{prop-tuaadeq2}
An affine formal scheme $\Spf A$, where $A$ is an adic ring of finite ideal type, is universally rigid-Noetherian $($resp.\ universally adhesive$)$ if and only if the topological ring $A$ is t.u.\ rigid-Noetherian\index{t.u. rigid-Noetherian ring@t.u.\ rigid-Noetherian ring} {\rm (\ref{dfn-tuaringadmissible} (1))} $($resp.\ t.u.\ adhesive\index{t.u.a. ring@t.u.\ adhesive ring} {\rm (\ref{dfn-tuaringadmissible} (2))}$)$.
\end{prop}

\begin{proof}
The `if' part is clear.
To show the converse, take a finite affine covering $\coprod_{\alpha\in L}\Spf A_{\alpha}\rightarrow\Spf A$ such that each $A_{\alpha}$ is t.u.\ rigid-Noetherian (resp.\ t.u.\ adhesive).
Applying \ref{prop-formalfpqcpatchingtua} to the map $A\rightarrow B=\prod_{\alpha\in L}A_{\alpha}$ (cf.\ {\bf \ref{ch-pre}}.\ref{prop-adhesive0} (1)), we deduce that $A$ is t.u.\ rigid-Noetherian (resp.\ t.u.\ adhesive), as desired.
\end{proof}

\begin{cor}\label{prop-tuaadeq3}
Let $A$ be a t.u.\ rigid-Noetherian $($resp.\ t.u.\ adhesive$)$ ring, and $X$ a formal scheme locally of finite type over $\Spf A$.
Then $X$ is locally universally rigid-Noetherian $($resp.\ locally universally adhesive$)$. \hfill$\square$
\end{cor}

\begin{prop}\label{prop-fiberprodadq}
Let $X\stackrel{f}{\rightarrow}Z\stackrel{g}{\leftarrow}Y$ be a diagram of adic formal schemes.
Suppose that $X$ is locally universally rigid-Noetherian $($resp.\ locally universally adhesive$)$ and that $g$ is locally of finite type.
Then the fiber product $X\times_ZY$ in the category of formal schemes is locally universally rigid-Noetherian $($resp.\ locally universally adhesive$)$, and the projection $X\times_ZY\rightarrow X$ is adic.
\end{prop}

\begin{proof}
The first assertion follows from the fact that $X\times_ZY$ is locally of finite type over $X$ (\ref{prop-topfintype2} (4)).
The other assertion follows from \ref{prop-adicmor1} (3).
\end{proof}

\begin{exas}\label{exas-adqformalscheme}{\rm 
(1) Any locally Noetherian formal schemes\index{formal scheme!Noetherian formal scheme@Noetherian ---!locally Noetherian formal scheme@locally --- ---} are locally universally adhesive; this follows from what we have seen in \ref{exas-tuaringsadmissible} (1).

(2) Let $V$ be as in \ref{exas-tuaringsadmissible} (2), and consider the affine formal scheme $\Spf V$.
By \ref{prop-tuaadeq2} we deduce that $\Spf V$ is an universally adhesive formal scheme.
Hence by \ref{prop-tuaadeq3} any formal scheme locally of finite type over $\Spf V$ is universally adhesive.}
\end{exas}
\index{formal scheme!universally adhesive formal scheme@universally adhesive ---!locally universally adhesive formal scheme@locally --- ---|)}
\index{formal scheme!universally rigid-Noetherian formal scheme@universally rigid-Noetherian ---!locally universally rigid-Noetherian formal scheme@locally --- ---|)}\index{formal scheme!universally rigid-Noetherian formal scheme@universally rigid-Noetherian ---|)}\index{adhesive!universally adhesive@universally ---!universally adhesive formal scheme@--- --- formal scheme|)}\index{formal scheme!universally adhesive formal scheme@universally adhesive ---|)}

\subsubsection{Categories of universally rigid-Noetherian formal schemes}\label{subsub-ntnadequatecategory}
We often use the following notations for categories:
\begin{itemize}
\item $\RNoe\Fs=$ the category of locally universally rigid-Noetherian formal schemes;
\item $\RNoe\Fs^{\ast}=$ the category of locally universally rigid-Noetherian formal schemes with adic morphisms;
\item $\RNoe\Fs^{\ast}_{/S}=$ the category of locally universally rigid-Noetherian formal schemes adic over $S$.
\end{itemize}

The similar categories of locally universally adhesive formal schemes are likewise denoted by $\Adh\Fs$, $\Adh\Fs^{\ast}$, and $\Adh\Fs^{\ast}_{/S}$, respectively.

We also define $\RNoe\Fs_S$, $\RNoe\Fs^{\ast}_S$, $\Adh\Fs_S$, $\Adh\Fs^{\ast}_S$ etc.\ in the usual way for any formal scheme $S$; note that in the categories $\RNoe\Fs^{\ast}_S$ and $\Adh\Fs^{\ast}_S$ all arrows are adic, but neither the base $S$ nor the structural map $X\rightarrow S$ are assumed to be adic.
Note also that in $\RNoe\Fs^{\ast}_{/S}$ and $\Adh\Fs^{\ast}_{/S}$ the base formal scheme $S$, necessarily adic of finite type, is not required to be locally universally adhesive nor locally universally rigid-Noetherian.

We also have the `coherent version':
\begin{itemize}
\item $\RNoe\CFs=$ the category of coherent universally rigid-Noetherian formal schemes;
\item $\RNoe\CFs^{\ast}=$ the category of coherent universally rigid-Noetherian formal schemes with adic morphisms;
\item $\RNoe\CFs^{\ast}_{/S}=$ the category of coherent universally rigid-Noetherian formal schemes adic over $S$;
\end{itemize}
and `affine version' of the above categories:
\begin{itemize}
\item $\Af\RNoe\Fs=$ the category of affine universally rigid-Noetherian formal schemes;
\item $\Af\RNoe\Fs^{\ast}=$ the category of affine universally rigid-Noetherian formal schemes with adic morphisms;
\item $\Af\RNoe\Fs^{\ast}_{/S}=$ the category of affine universally rigid-Noetherian formal schemes adic over $S$.
\end{itemize}

Similar rule of notation is also put on coherent and affine universally adhesive formal schemes; for example, $\Adh\CFs^{\ast}$ denotes the category of coherent universally adhesive formal schemes with adic morphisms, etc.

\subsection{Morphisms of finite presentation}\label{sub-finipresformal}
\begin{dfn}\label{dfn-topfinpres}{\rm 
A morphism $f\colon X\rightarrow Y$ of locally universally rigid-Noetherian\index{formal scheme!universally rigid-Noetherian formal scheme@universally rigid-Noetherian ---!locally universally rigid-Noetherian formal scheme@locally --- ---} formal schemes is said to be {\em locally of finite presentation}\index{morphism of formal schemes@morphism (of formal schemes)!morphism of formal schemes locally of finite presentation@--- locally of finite presentation} if the following conditions are satisfied:
\begin{itemize}
\item[{\rm (a)}] the morphism $f$ is adic {\rm (\ref{dfn-adicmor})}; 
\item[{\rm (b)}] there exist an affine open covering $\{V_{\alpha}\}$ of $Y$ with $V_{\alpha}=\Spf B_{\alpha}$, where each $B_{\alpha}$ is an adic ring of finite ideal type, and for each $\alpha$ an affine open covering $\{U_{\alpha\beta}\}_{\beta}$ of $f^{-1}(V_{\alpha})$ with $U_{\alpha\beta}=\Spf A_{\alpha\beta}$, where each $A_{\alpha\beta}$ is an adic ring topologically finitely presented\index{finitely presented!topologically finitely presented@topologically ---} over $B_{\alpha}$ {\rm ({\rm {\bf \ref{ch-pre}}.\ref{dfn-topfinigen}})}.
\end{itemize}
The morphism $f$ is said to be {\em of finite presentation}\index{morphism of formal schemes@morphism (of formal schemes)!morphism of formal schemes of finite presentation@--- of finite presentation} if it is locally of finite presentation and quasi-compact\index{space@space (topological)!quasi-compact topological space@quasi-compact ---}\index{quasi-compact!quasi-compact topological space@--- (topological) space} {\rm (\ref{dfn-qcompformal})}.}
\end{dfn}

\begin{prop}\label{prop-topfinpres2}
{\rm (1)} An open immersion\index{immersion!open immersion of formal schemes@open --- (of formal schemes)} between locally universally rigid-Noetherian formal schemes is locally of finite presentation.

{\rm (2)} The composition of two morphisms locally of finite presentation $($resp.\ of finite presentation$)$ is again locally of finite presentation $($resp.\ of finite presentation$)$.

{\rm (3)} Let $X$, $X'$, $Y$, and $Y'$ be locally universally rigid-Noetherian formal schemes adic over an adic formal scheme $S$ of finite ideal type, and $f\colon X\rightarrow X'$ and $g\colon Y\rightarrow Y'$ two $S$-morphisms.
Suppose $X'\times_SY'$ is locally universally rigid-Noetherian.
Then if $f$ and $g$ are locally of finite presentation $($resp.\ of finite presentation$)$, so is $f\times_Sg\colon X\times_SY\rightarrow X'\times_SY'$.

{\rm (4)} If $S$ is an adic formal scheme of finite ideal type and if $f\colon X\rightarrow Y$ is an $S$-morphism locally of finite presentation $($resp.\ of finite presentation$)$ between locally universally rigid-Noetherian formal schemes adic over $S$, then for any morphism $S'\rightarrow S$ of adic formal schemes of finite ideal type such that $Y\times_SS'$ is universally rigid Noetherian, $f_{S'}\colon X\times_SS'\rightarrow Y\times_SS'$ is locally of finite presentation $($resp.\ of finite presentation$)$.
\end{prop}

\begin{proof}
For a t.u.\ rigid-Noetherian ring $A$ and $f\in A$, the ring $A_{\{f\}}\cong A\dl T\dr/(fT-1)$ (\ref{lem-formalnot0}) is topologically of finitely presentation over $A$, whence (1).
The assertion (2) is easy to see. 

To show (3), we may work in the affine situation.
Let $(R,I)$ be a complete pair, and consider the morphisms of t.u.\ rigid-Noetherian rings, all adic over $R$:
$$
A\longrightarrow A\dl X_1,\ldots,X_n\dr/\mathfrak{a},\qquad
B\longrightarrow B \dl Y_1,\ldots,Y_m\dr/\mathfrak{b},
$$
where $\mathfrak{a}$ and $\mathfrak{b}$ are finitely generated, hence closed (\ref{rem-turigidnoetherianbasicproperties}), ideals.
We want to show is that the closure of the ideal generated by the image of
$$
\mathfrak{c}=\mathfrak{a}\otimes B \dl Y_1,\ldots,Y_m\dr+A\dl X_1,\ldots,X_n\dr\otimes\mathfrak{b}
$$
in the ring 
$$
A\dl X_1,\ldots,X_n\dr\widehat{\otimes}_RB \dl Y_1,\ldots,Y_m\dr=(A\widehat{\otimes}_RB)\dl X_1,\ldots,X_n,Y_1,\ldots,Y_m\dr
$$
is finitely generated.
But since $A\widehat{\otimes}_RB$ is assumed to be t.u.\ rigid-Noetherian, the ideal generated by the image of $\mathfrak{c}$ is already closed (\ref{rem-turigidnoetherianbasicproperties}), whence (3).

The assertion (4) also follows from the affine local observation as follows.
Let $S=\Spf R$ where $R$ is an adic ring of finite ideal type, and $Y=\Spf B$ where $B$ is t.u.\ rigid-Noetherian and adic over $R$.
Let $X=\Spf A$ with $A=B\dl X_1,\ldots,X_n\dr/\mathfrak{a}$ where $\mathfrak{a}$ is finitely generated ideal, and $S'=\Spf R'$ where $R'$ is an adic ring of finite ideal type that is adic over $R$.
What to prove is that the morphism 
$$
\Spf A\widehat{\otimes}_RR'\longrightarrow\Spf B\widehat{\otimes}_RR'
$$
is finitely presented if $B\widehat{\otimes}_RR'$ is t.u.\ rigid-Noetherian.
Set $B'=B\widehat{\otimes}_RR'$.
We have $A\widehat{\otimes}_RR'=B'\dl X_1,\ldots,X_n\dr/\mathfrak{b}$, where $\mathfrak{b}$ is the closure of $\mathfrak{a}B'\dl X_1,\ldots,X_n\dr$.
Since the ring $A\widehat{\otimes}_RR'$ is t.u.\ rigid-Noetherian, the ideal $\mathfrak{a}B'\dl X_1,\ldots,X_n\dr$ is already closed, and hence $\mathfrak{b}=\mathfrak{a}B'\dl X_1,\ldots,X_n\dr$, which is finitely generated. 
\end{proof}

\begin{prop}\label{prop-topfinpres1}
Let $A$ be a t.u.\ rigid-Noetherian ring, $I\subseteq A$ a finitely generated ideal of definition, and $B$ a topologically finitely generated $A$-algebra.
Then the following conditions are equivalent$:$
\begin{itemize}
\item[{\rm (a)}] the morphism $\Spf B\rightarrow\Spf A$ is of finite presentation$;$
\item[{\rm (b)}] the morphism $\Spec B/I^{k+1}B\rightarrow\Spec A/I^{k+1}$ is of finite presentation for any $k\geq 0;$
\item[{\rm (c)}] $B$ is topologically finitely presented over $A$.
\end{itemize}
\end{prop}

\begin{proof}
The implications (c) $\Rightarrow$ (a) $\Rightarrow$ (b) are obvious, while (b) $\Rightarrow$ (c) follows from {\bf \ref{ch-pre}}.\ref{prop-formalnot31}.
\end{proof}

\begin{cor}\label{cor-topfinpres11}
Let $f\colon X\rightarrow Y$ be an adic morphism of locally universally rigid-Noetherian formal schemes.
Suppose that $Y$ has an ideal of definition $\mathscr{I}$ of finite type.
Set $X_k=(X,\O_X/\mathscr{I}^{k+1}\O_X)$ and $Y_k=(Y,\O_Y/\mathscr{I}^{k+1})$ for $k\geq 0$, and denote by $f_k\colon X_k\rightarrow Y_k$ the induced morphism of schemes.
Then the following conditions are equivalent$:$
\begin{itemize}
\item[{\rm (a)}] $f$ is locally of finite presentation $($resp.\ of finite presentation$);$
\item[{\rm (b)}] $f_k$ is locally of finite presentation $($resp.\ of finite presentation$)$ for any $k\geq 0$. \hfill$\square$
\end{itemize}
\end{cor}

\subsection{Relation with other notions}\label{sub-interrelationwithadmissibleformalschemes}
\subsubsection{Admissible formal schemes}\label{subsub-admissibleformalschemes}
\begin{prop}\label{prop-fintypetorfreefinpres}
Let $Y$ be an locally universally adhesive formal scheme, and $f\colon X\rightarrow Y$ a morphism locally of finite type $($resp.\ of finite type$)$.
Suppose that the following condition is satisfied$:$
\begin{itemize}
\item[$(\ast)$] there exist an open covering $X=\bigcup_{\alpha\in L}U_{\alpha}$ and for each $\alpha\in L$ an ideal of definition $\mathscr{I}_{\alpha}$ of finite type on $U_{\alpha}$ such that $\O_{U_{\alpha}}$ is $\mathscr{I}_{\alpha}$-torsion free.
\end{itemize}
Then $f$ is locally of finite presentation $($resp.\ of finite presentation$)$.
\end{prop}

\begin{proof}
The assertion follows from the following argument.
Let $B$ be a t.u.\ adhesive ring with a finitely generated ideal of definition $I\subseteq B$, and consider $A=B\dl X_1,\ldots,X_n\dr/\mathfrak{a}$, where $\mathfrak{a}$ is an ideal.
Suppose $A$ is $I$-torsion free.
Since $B$ is t.u.\ adhesive, the pair $(B\dl X_1,\ldots,X_n\dr,IB\dl X_1,\ldots,X_n\dr)$ is adhesive.
Since $A$ is $I$-torsion free, we deduce from {\bf \ref{ch-pre}}.\ref{prop-adhesive} that $A$ is finitely presented as a module over $B\dl X_1,\ldots,X_n\dr$ or, equivalently, that $\mathfrak{a}$ is finitely generated.
\end{proof}

\begin{cor}\label{cor-fintypetorfreefinpresvaluation}
Let $V$ be an $a$-adically complete valuation ring\index{valuation!valuation ring@--- ring!a-adically complete valuation ring@$a$-adically complete --- ---} $($of arbitrary height$)$ where $a\in\m_V\setminus\{0\}$, and $X$ a formal scheme locally of finite type over $\Spf V$.
Suppose that the structure sheaf $\O_X$ is $a$-torsion free.
Then $X$ is locally of finite presentation over $\Spf V$. \hfill$\square$
\end{cor}

Notice that any $V$-algebra $A$ is $a$-torsion free if and only if it is flat over $V$ (Exercise \ref{exer-valflatatorfree}).
In classical rigid geometry the so-called {\em admissible formal schemes}\index{formal scheme!admissible formal scheme@admissible ---}\index{admissible!admissible formal scheme@--- formal scheme} play a central role; cf.\ \cite[\S1]{BL1}:
\begin{dfn}\label{dfn-admissibleformalschemes}{\rm 
An {\em admissible formal scheme}\index{formal scheme!admissible formal scheme@admissible ---}\index{admissible!admissible formal scheme@--- formal scheme} is a formal scheme $X$ locally of finite type over $\Spf V$, where $V$ is an $a$-adically complete valuation ring $(a\in\m_V\setminus\{0\})$ of height one, such that the structure sheaf $\O_X$ is $a$-torsion free.}
\end{dfn}

It follows from \ref{cor-fintypetorfreefinpresvaluation} that any admissible formal scheme is locally of finite presentation over $\Spf V$.

\begin{rem}\label{rem-adequateformalschemestypeVadmissible}{\rm
Let $V$ be an $a$-adically complete valuation ring (where $a\in\m_V\setminus\{0\}$) of arbitrary height, and $X$ a formal scheme locally of finite type over $S=\Spf V$ such that $\O_X$ is $a$-torsion free.
Let $\mathfrak{p}=\sqrt{aV}$ be the associated height one prime\index{associated height one prime} of $V$ ({\bf \ref{ch-pre}}.\ref{dfn-maxspe2}), and set $V'=V_{\mathfrak{p}}$, which is an $a$-adically complete valuation ring of height one ({\bf \ref{ch-pre}}.\ref{cor-compval2006ver23}).
Then in view of {\bf \ref{ch-pre}}.\ref{prop-convpreadhloc1} we see that the base change $X\times_SS'$ (where $S'=\Spf V'$) is an admissible formal scheme.}
\end{rem}

\subsubsection{Interrelations between the classes}\label{subsub-catout}
Finally, let us overview the interrelations between the classes of formal schemes so far introduced.
Let $\mathbf{Noe}\Fs^{\ast}$ be the category of locally Noetherian formal schemes with adic morphisms, which is, in fact, a full subcategory of $\Adh\Fs^{\ast}$ (\ref{exas-adqformalscheme} (1)).
Hence one can draw the following diagram:
$$
\begin{xy}
(0,0)="NF"*{\mathbf{Noe}\Fs^{\ast}},+/r7em/="AdF"*{\Adh\Fs^{\ast}},+/r7em/="AcF"*{\Ac\Fs^{\ast}},+/r6em/="F"*{\Fs}, "NF"+<0em,-10ex>="FD"*{\Fs^{\ast}_{\fin/\mathrm{DVR}}},+/r7em/="FV"*{\Fs^{\ast}_{\fin/\mathrm{Val}}.}
\ar@{^{(}->}"NF"+<2.4em,-.3ex>;"AdF"+<-2.3em,-.3ex>
\ar@{^{(}->}"AdF"+<2.4em,-.3ex>;"AcF"+<-1.9em,-.3ex>
\ar@{^{(}->}"AcF"+<2.1em,-.3ex>;"F"+<-.9em,-.3ex>
\ar@{^{(}->}"FD"+<0em,2.6ex>;"NF"+<0em,-2.1ex>
\ar@{^{(}->}"FV"+<0em,2.6ex>;"AdF"+<0em,-2.1ex>
\ar@{^{(}->}"FD"+<2.6em,-.3ex>;"FV"+<-2.3em,-.3ex>
\end{xy}\leqno{(\ast)}
$$
Here $\Fs^{\ast}_{\fin/\mathrm{DVR}}$ denotes the full subcategory of $\mathbf{Noe}\Fs^{\ast}$ consisting of formal schemes locally of finite type over a complete discrete valuation ring, and $\Fs^{\ast}_{\fin/\mathrm{Val}}$ is the full subcategory of $\Adh\Fs^{\ast}$ consisting of formal schemes locally of finite type over an $a$-adically complete valuation ring of arbitrary height (cf.\ \ref{exas-adqformalscheme} (2)).
Notice that the category $\RNoe\Fs^{\ast}$ of locally universally rigid-Noetherian formal schemes lies in between $\Noe\Fs^{\ast}$ and $\Adh\Fs^{\ast}$.
Moreover:
\begin{itemize}
\item the first two inclusions in the first row are fully faithful, whereas the last one is only faithful; 
\item only $\Fs$ in $(\ast)$ has a final object (cf.\ \ref{rem-clear}); 
\end{itemize}



\section{Adically quasi-coherent sheaves}\label{sec-adicallyqcoh}
\index{adically quasi-coherent (a.q.c.) sheaf|(}
From now on, we mainly deal with {\em adic} formal schemes {\em of finite ideal type}\index{formal scheme!adic formal scheme@adic ---!adic formal scheme of finite ideal type@--- --- of finite ideal type}\index{adic!adic formal scheme@--- formal scheme!adic formal scheme of finite ideal type@--- --- of finite ideal type}.
On such formal schemes, a reasonable class of $\O_X$-module sheaves is provided by so-called {\em adically quasi-coherent $($a.q.c.$)$ sheaves}\index{quasi-coherent!adically quasi-coherent OX module@adically --- (a.q.c.) sheaf}\index{adically quasi-coherent (a.q.c.) sheaf}, which we are going to discuss in this section.
They are complete (\ref{dfn-completesheaf} (2)) $\O_X$-modules such that the truncated pieces (that is, the sheaves obtained by modulo ideals of definition) are quasi-coherent sheaves on the induced schemes.
In \S\ref{sub-adicallyqcohdef} we give the general definition of adically quasi-coherent sheaves and discuss some of the general properties.

In the affine situation, as we will discuss in \S\ref{sub-adicqcohaff}, adically quasi-coherent sheaves are obtained by `$\Delta$-construction' (\cite[$\mathbf{I}$, (10.10.1)]{EGA}), that is, the sheaves on $X=\Spf A$ of the form $M^{\Delta}$ by an $A$-module $M$.
In fact, if $A$ is an adic ring of finite ideal type\index{adic!adic ring@--- ring!adic ring of finite ideal type@--- --- of finite ideal type} (\ref{dfn-admissibleringoffiniteidealtype}), then the $\Delta$-construction gives rise to a categorical equivalence between the category of complete $A$-modules and the category of adically quasi-coherent sheaves on $X=\Spf A$ (\ref{thm-adicqcoh1} (2)).

After discussing adically quasi-coherent sheaves as projective limits of quasi-coherent sheaves on schemes in \S\ref{sub-adicallyqcohinductivelimit}, we discuss the case where the formal schemes are locally universally rigid-Noetherian in \S\ref{sub-adicqcohpf}.
The notion of adically quasi-coherent sheaves in this particular situation reveals itself to be much more tractable.
For example, the above-mentioned categorical equivalence restricted on finitely generated modules and a.q.c.\ sheaves of finite type is an {\em exact} equivalence if $X=\Spf A$ is universally rigid-Noetherian (\ref{thm-adicqcohpre1}).

In the final subsection \S\ref{sub-admidealextlem}, we discuss the so-called {\em admissible ideals}\index{admissible!admissible ideal@--- ideal}, which are, so to speak, the sheaf version of $I$-admissible ideals\index{admissible!I-admissible@$I$-{---}} ({\bf \ref{ch-pre}}.\ref{dfn-adm}).
This class of ideals contains ideals of definition and has a lot of nice properties. 
One of them is the `extension property', which is most significantly shown in \ref{prop-extension1} below.
In pursuing the extension properties of admissible ideals, we will obtain in \S\ref{sub-admidealextlem} a technically remarkable result saying that any coherent\index{formal scheme!coherent formal scheme@coherent ---} (\ref{dfn-cohformalmorphism}) adic formal scheme of finite ideal type always admits an ideal of definition of finite type (\ref{cor-extension2}), a generalization of the known fact that any Noetherian formal scheme has a coherent ideal of definition (cf.\ \ref{prop-northeroanformalschemes}).

\subsection{Complete sheaves and adically quasi-coherent sheaves}\label{sub-adicallyqcohdef}
\subsubsection{Hausdorff completion of $\O_X$-modules}\label{subsub-adicallyqcohdefcompl}
Let $X$ be an adic formal scheme of finite ideal type, and suppose for the time being that $X$ has an ideal of definition $\mathscr{I}$ of finite type.
In this situation, $\{\mathscr{I}^{k+1}\}_{k\geq 0}$ gives a fundamental system of ideals of definition of $X$ (\ref{cor-adicformalschemeidealofdefinitionfinitetype22}).
For $k\geq 0$, we set $X_k=(X,\O_X/\mathscr{I}^{k+1})$.
This is a scheme having the same underlying topological space as $X$, and $\{X_k\}_{k\geq 0}$ forms a filtered inductive system of schemes together with the closed immersions $X_k\hookrightarrow X_l$ as the transition maps for all $k\leq l$.
Moreover, the formal scheme $X$ coincides with the inductive limit of this system $X=\varinjlim_kX_k$ (\ref{prop-formalindlimschfact}), or equivalently, the structure sheaf $\O_X$ is the projective limit (cf.\ {\bf \ref{ch-pre}}, \S\ref{subsub-projlimsheaves})
$$
\O_X=\varprojlim_k\O_{X_k}
$$
as a sheaf of topological rings on the topological space $X$.

For an $\O_X$-module $\mathscr{F}$ we set 
$$
\mathscr{F}_k=\mathscr{F}/\mathscr{I}^{k+1}\mathscr{F}
$$
for $k\geq 0$.
This is an $\O_{X_k}$-module on the scheme $X_k$, and the resulting system $\{\mathscr{F}_k\}_{k\geq 0}$ together with the obvious transition morphisms is a filtered projective system of abelian sheaves on the topological space $X$, which admits the compatible action of the projective system of rings $\{\O_{X_k}\}_{k\geq 0}$; that is, each $\mathscr{F}_k$ is an $\O_{X_k}$-module and the transition maps $\mathscr{F}_l\rightarrow\mathscr{F}_k$ ($k\leq l$) are compatible with $\O_{X_l}\rightarrow\O_{X_k}$.
Hence the projective limit
$$
\widehat{\mathscr{F}}=\varprojlim_k\mathscr{F}_k
$$
has the canonical $\O_X$-module structure, and the canonical morphism
$$
i_{\mathscr{F}}\colon\mathscr{F}\longrightarrow\widehat{\mathscr{F}}
$$
is a morphism of $\O_X$-modules.
Notice that the definition of $\widehat{\mathscr{F}}$ does not depend on the choice of the ideal of definition $\mathscr{I}$; indeed, the projective limit as above coincides with the projective limit $\varprojlim_{\mathscr{J}}(\mathscr{F}/\mathscr{J}\mathscr{F})$, where $\mathscr{J}$ runs in the filtered set of all ideals of definition on $X$ ordered by the reversed inclusion order\index{ordering@order(ing)!inclusion ordering@inclusion ---}, and for a fixed ideal of definition $\mathscr{I}$ of finite type, the collection $\{\mathscr{I}^{k+1}\}_{k\geq 0}$ gives a cofinal subset.

By the last-mentioned fact, for a given $\O_X$-module $\mathscr{F}$ we can define $\widehat{\mathscr{F}}$ up to canonical isomorphism even in case $X$ does not have an ideal of definition; indeed, since $X$ has locally an ideal of definition of finite type, one can define $\widehat{\mathscr{F}}$ locally on $X$, and they glue with each other by canonical isomorphisms.
The canonical morphism $\mathscr{F}\rightarrow\widehat{\mathscr{F}}$ can be likewise defined.

\begin{dfn}\label{dfn-completesheaf}{\rm 
Let $X$ be an adic formal scheme of finite ideal type, and $\mathscr{F}$ an $\O_X$-module.

(1) The above-defined $\O_X$-module $\widehat{\mathscr{F}}$ together with the canonical morphism $i_{\mathscr{F}}\colon\mathscr{F}\rightarrow\widehat{\mathscr{F}}$ is called the {\em completion}\index{completion} of $\mathscr{F}$.

(2) We say that $\mathscr{F}$ is {\em complete}\index{complete} if $i_{\mathscr{F}}\colon\mathscr{F}\rightarrow\widehat{\mathscr{F}}$ is an isomorphism.}
\end{dfn}

\subsubsection{Adically quasi-coherent (a.q.c.) sheaves}\label{subsub-adicallyqcohdefdef}
\index{quasi-coherent!adically quasi-coherent OX module@adically --- (a.q.c.) sheaf|(}
\begin{lem}\label{lem-propadicallyqcohdef1}
Let $X$ be an adic formal scheme of finite ideal type, and $\mathscr{I}$ an ideal of definition of finite type of $X$.
Let $\mathscr{F}$ be an $\O_X$-module.
Then the following conditions are equivalent$:$
\begin{itemize}
\item[{\rm (a)}] for any $k\geq 0$ the sheaf $\mathscr{F}_k=\mathscr{F}/\mathscr{I}^{k+1}\mathscr{F}$ is a quasi-coherent sheaf on the scheme $X_k=(X,\O_X/\mathscr{I}^{k+1});$
\item[{\rm (b)}] for any ideal of definition $\mathscr{J}$ of $X$ the sheaf $\mathscr{F}/\mathscr{J}\mathscr{F}$ is a quasi-coherent sheaf on the scheme $(X,\O_X/\mathscr{J})$.
\end{itemize}
\end{lem}

\begin{proof}
The implication (b) $\Rightarrow$ (a) is trivial.
Suppose (a) holds.
Since the conditions are local on $X$, we may assume that $X$ is affine (hence, in particular, the underlying topological space is quasi-compact).
For any ideal of definition $\mathscr{J}$ we can find a positive integer $k$ such that $\mathscr{I}^{k+1}\subseteq \mathscr{J}$.
Consider the closed immersion $i\colon (X,\O_X/\mathscr{J})\hookrightarrow X_k=(X,\O_X/\mathscr{I}^{k+1})$.
The sheaf $\mathscr{F}/\mathscr{J}\mathscr{F}$ on $(X,\O_X/\mathscr{J})$ coincides with $i^{\ast}\mathscr{F}_k$, which is quasi-coherent.
\end{proof}

\begin{dfn}\label{dfn-adicqcoh}{\rm 
(1) Let $X$ be an adic formal scheme of finite ideal type, and $\mathscr{F}$ an $\O_X$-module.
We say that $\mathscr{F}$ is {\em adically quasi-coherent} (acronym: a.q.c.) if the following conditions are satisfied:
\begin{itemize}
\item[{\rm (a)}] $\mathscr{F}$ is complete; 
\item[{\rm (b)}] for any open subset $U\subseteq X$ considered as an open formal subscheme\index{formal subscheme!open formal subscheme@open ---} (\S\ref{subsub-formalnotformalsch}) and for any ideal of definition $\mathscr{I}$ of finite type of $U$, the sheaf $(\mathscr{F}|_U)/\mathscr{I}(\mathscr{F}|_U)$ is a quasi-coherent sheaf\index{quasi-coherent!quasi-coherent OX module on schemes@--- sheaf (on a scheme)} on the scheme $(U,\O_U/\mathscr{I})$.
\end{itemize}

(2) An adically quasi-coherent sheaf $\mathscr{F}$ on $X$ is said to be {\em of finite type}\index{adically quasi-coherent (a.q.c.) sheaf!adically quasi-coherent sheaf of finite type@--- of finite type} if it is of finite type as an $\O_X$-module.

(3) A {\em morphism}\index{adically quasi-coherent (a.q.c.) sheaf!morphism of adically quasi-coherent sheaves@morphism of --es} between adically quasi-coherent sheaves is a morphism of $\O_X$-modules.}
\end{dfn}

By \ref{lem-propadicallyqcohdef1} the condition (b) is equivalent to the following one:
\begin{itemize}
\item[${\rm (b)}'$] there exist an open covering $X=\bigcup_{\alpha\in L}U_{\alpha}$ and for each $\alpha\in L$ an ideal of definition $\mathscr{I}_{\alpha}$ of finite type of $U_{\alpha}$ such that for any $\alpha\in L$ and $k\geq 0$ the sheaf $(\mathscr{F}|_{U_{\alpha}})/\mathscr{I}^{k+1}_{\alpha}(\mathscr{F}|_{U_{\alpha}})$ is a quasi-coherent sheaf on the scheme $(U_{\alpha},\O_{U_{\alpha}}/\mathscr{I}^{k+1}_{\alpha})$.
\end{itemize}

If $X$ itself has an ideal of definition $\mathscr{I}$ of finite type, then again by \ref{lem-propadicallyqcohdef1} the last condition is equivalent to that, with the notation as in \ref{lem-propadicallyqcohdef1}, $\mathscr{F}_k$ is quasi-coherent on $X_k$ for any $k\geq 0$.

Notice that any morphism of adically quasi-coherent sheaves on $X$ is continuous in the following sense: If $X$ has an ideal $\mathscr{I}$ of definition of finite type, such a morphism $f\colon\mathscr{F}\rightarrow\mathscr{G}$ induces the morphism $f_k\colon\mathscr{F}_k\rightarrow\mathscr{G}_k$ of quasi-coherent sheaves on $X_k$ for $k\geq 0$ and coincides with the projective limit of $\{f_k\}$.

We denote by $\AQCoh_X$ the category of adically quasi-coherent sheaves on an adic formal scheme $X$ of finite ideal type.

\begin{prop}\label{prop-adicallyquasicohstr}
Let $X$ be an adic formal scheme of finite ideal type. 

{\rm (1)} The structure sheaf $\O_X$ is adically quasi-coherent of finite type.

{\rm (2)} Any ideal of definition $\mathscr{I}$ of $X$ is an adically quasi-coherent ideal of $\O_X$.

{\rm (3)} If $\mathscr{I}$ is an ideal of definition and $\mathscr{F}$ is an adically quasi-coherent sheaf, then $\mathscr{I}\mathscr{F}$ is again an adically quasi-coherent sheaf.
\end{prop}

\begin{proof}
The assertion (1) is clear.
The assertion (2) follows from (3) applied to $\mathscr{F}=\O_X$.
To show (3), since the question is local on $X$, we may assume that there exists an ideal of definition $\mathscr{J}$ of finite type such that $\mathscr{I}^n\subseteq\mathscr{J}\subseteq\mathscr{I}$.
We first look at the exact sequence
$$
0\longrightarrow\mathscr{I}\mathscr{F}/\mathscr{J}^{k+1}\mathscr{F}\longrightarrow\mathscr{F}/\mathscr{J}^{k+1}\mathscr{F}\longrightarrow\mathscr{F}/\mathscr{I}\mathscr{F}\longrightarrow 0
$$
of $\O_X$-modules for $k\geq 1$.
Since $\mathscr{F}/\mathscr{J}^{k+1}\mathscr{F}$ and $\mathscr{F}/\mathscr{I}\mathscr{F}$ are quasi-coherent sheaves on the scheme $X_k=(X,\O_X/\mathscr{J}^{k+1})$ (\ref{lem-propadicallyqcohdef1}), we deduce that $\mathscr{I}\mathscr{F}/\mathscr{J}^{k+1}\mathscr{F}$ is quasi-coherent.
Then by \ref{lem-vanishingcohomologyadicallyuseful} (1) we see that 
$$
0\longrightarrow\varprojlim_{k\geq 1}\mathscr{I}\mathscr{F}/\mathscr{J}^{k+1}\mathscr{F}\longrightarrow\mathscr{F}\longrightarrow\mathscr{F}/\mathscr{I}\mathscr{F}\longrightarrow 0
$$
is exact and hence that $\mathscr{I}\mathscr{F}=\varprojlim_{k\geq 1}\mathscr{I}\mathscr{F}/\mathscr{J}^{k+1}\mathscr{F}$.
This shows that $\mathscr{I}\mathscr{F}$ is complete, for we have $\mathscr{J}^{k+n-1}\mathscr{I}\mathscr{F}\subseteq\mathscr{J}^{k+1}\mathscr{F}\subseteq\mathscr{J}^k\mathscr{I}\mathscr{F}$ for any $k\geq 1$.
\end{proof}

The following proposition is clear:
\begin{prop}\label{prop-adicqchohexa}
Let $X$ be a scheme, and $Y\subseteq X$ a closed subscheme of finite presentation.
Consider the formal completion\index{completion!formal completion@formal ---} $\widehat{X}|_Y$ $($cf.\ {\rm \ref{prop-formalcompletion1}}$)$.
If $\mathscr{F}$ is a quasi-coherent sheaf on $X$, then its formal completion $\widehat{\mathscr{F}}|_Y$ {\rm (\S\ref{subsub-formalcompletionqcohsh})} is adically quasi-coherent. \hfill$\square$
\end{prop}
\index{quasi-coherent!adically quasi-coherent OX module@adically --- (a.q.c.) sheaf|)}

\subsection{A.q.c.\ sheaves on affine formal schemes}\label{sub-adicqcohaff}
\subsubsection{The $\Delta$-sheaves}\label{subsub-thedeltasheaves}
 We first recall the definition of the sheaf of $\O_X$-modules $M^{\Delta}$ on an adic formal scheme $X=\Spf A$ associated to an $A$-module $M$ (\cite[$\mathbf{I}$, (10.10.1)]{EGA}). 
Let $I\subseteq A$ be an ideal of definition, and set $X_k=\Spec A/I^{k+1}$ for $k\geq 0$; we have $X=\varinjlim_kX_k$ (\ref{prop-formalindlimschfact}).
For an $A$-module $M$, $M_k=M/I^{k+1}M$ defines the quasi-coherent sheaf $\til{M}_k$ on the scheme $X_k$.
Then 
$$
M^{\Delta}=\varprojlim_{k\geq 0}\til{M}_k,
$$
which is a sheaf of $\O_X$-modules on $X$.
Here is an alternative definition of $M^{\Delta}$, equivalent to the above one up to canonical isomorphism: 
the sheaf $M^{\Delta}$ is the formal completion $\widehat{\til{M}}|_{X_0}$ (\S\ref{subsub-formalcompletionqcohsh}) of the quasi-coherent sheaf $\til{M}$ on $\Spec A$ along the closed subscheme $X_0=\Spec A/I$, that is, the sheaf $\varprojlim_k\til{M}/I^{k+1}\til{M}$ restricted on the topological space $X$.

Notice that, for any open ideal $J\subseteq A$, the sheaf $J^{\Delta}$ thus constructed coincides with the one given in \S\ref{subsub-formalschemesidealsofdefinition}; indeed, if $I\subseteq A$ is an ideal of definition contained in $J$, then, since $I^{k+1}J\subseteq I^{k+1}\subseteq I^kJ$ for any $k\geq 0$, we have $\varprojlim_{k\geq 0}\til{J/I^{k+1}}\cong\varprojlim_{k\geq 0}\til{J/I^{k+1}J}$.

The construction of $M^{\Delta}$ induces an additive functor 
$$
\cdot^{\Delta}\colon\Mod_A\longrightarrow\Mod_X,\qquad M\longmapsto M^{\Delta}
$$
from the category of $A$-modules to the category of $\O_X$-modules.

\begin{prop}\label{prop-adicqcohaff1g}
We have 
$$
\Gamma(X,M^{\Delta})=M^{\wedge}_{I^{\bullet}},
$$
where $M^{\wedge}_{I^{\bullet}}$ denotes the Hausdorff completion\index{complete!Hausdorff complete@Hausdorff ---} of $M$ with respect to the $I$-adic topology {\rm ({\bf \ref{ch-pre}}, \S\ref{subsub-completionfiltration})}$;$ more generally, for an affine open set $U=\Spf A_{\{f\}}$ with $f\in A$, we have 
$$
\Gamma(U,M^{\Delta})=(M_f)^{\wedge}_{I^{\bullet}},
$$
where $M_f=M\otimes_AA_f$.
\end{prop}

\begin{proof}
Indeed, one calculates 
$$
\Gamma(U,M^{\Delta})=\Gamma(U,\varprojlim_{k\geq 0}\til{M}_k)=\varprojlim_{k\geq 0}\Gamma(U,\til{M}_k)=\varprojlim_{k\geq 0}M_f/I^{k+1}M_f=(M_f)^{\wedge}_{I^{\bullet}}.
$$
\end{proof}

\subsubsection{Adically quasi-coherent $\Delta$-sheaves}\label{subsub-affinecaseadicallyqcohfinitype}
\begin{prop}\label{prop-adicqcohaff1f}
Let $A$ be an adic ring, and $I\subseteq A$ a finitely generated ideal.
Set $X=\Spf A$ and $\mathscr{I}=I^{\Delta}$.
Let $M$ an $A$-module, and consider the sheaf $M^{\Delta}$ on $X$.
We have 
$$
M^{\Delta}/\mathscr{I}^{k+1}M^{\Delta}\cong\til{M}_k,
$$
where $M_k=M/I^{k+1}M$ for $k\geq 0$, and the sheaf $M^{\Delta}$ is complete {\rm (\ref{dfn-completesheaf} (2))}.
In particular, the sheaf $M^{\Delta}$ is an adically quasi-coherent sheaf on $X$.
\end{prop}

To show the proposition, we first show:
\begin{lem}[{cf.\ {\bf \ref{ch-pre}}.\ref{prop-corpropARconseq1-1-1}; see also \ref{prop-lemdeltasheafadicallyquasicoherent1-3} below}]\label{lem-deltasheafadicallyquasicoherent1-1}
Let $A$ be an adic ring with a finitely generated ideal of definition $I$, $M$ an $A$-module, and $N\subseteq M$ an $A$-submodule open with respect to the $I$-adic topology on $M$.
Then the sequence 
$$
0\longrightarrow N^{\Delta}\longrightarrow M^{\Delta}\longrightarrow(M/N)^{\Delta}\longrightarrow 0
$$
by canonical maps is exact.
\end{lem}

\begin{proof}
Take $n\geq 0$ such that $I^{n+1}M\subseteq N$, and consider the exact sequence
$$
0\longrightarrow\til{N/I^{k+1}M}\longrightarrow\til{M_k}\longrightarrow\til{M/N}\longrightarrow 0
$$
of quasi-coherent sheaves on $X_k$ for $k\geq n$.
Since $I^{k+1}N\subseteq I^{k+1}M\subseteq I^{k-n}N$, we have $\varprojlim_k\til{N/I^{k+1}M}\cong N^{\Delta}$.
Moreover, since $N$ is open in $M$, we have $\til{M/N}=(M/N)^{\Delta}$.
Hence, taking projective limits along $k$, we get the desired exact sequence due to \ref{lem-vanishingcohomologyadicallyuseful} (1).
\end{proof}

\begin{lem}\label{lem-deltasheafadicallyquasicoherent1-2}
Let $A$ be an adic ring with a finitely generated ideal of definition $I\subseteq A$, and $M$ an $A$-module.
Set $\mathscr{I}=I^{\Delta}$.
Then for any $k\geq 0$ we have
$$
(I^{k+1}M)^{\Delta}=\mathscr{I}^{k+1}M^{\Delta}=I^{k+1}M^{\Delta}
$$
as sheaves on $X=\Spf A$.
\end{lem}

\begin{proof}
We first remark that the sheaf $\mathscr{I}^{k+1}M^{\Delta}$ is the associated sheaf of the presheaf given by $U\mapsto\Gamma(U,\mathscr{I}^{k+1})\cdot\Gamma(U,M^{\Delta})$ (\cite[$\mathbf{0}_{\mathbf{I}}$, (4.1.6)]{EGA}).
By \ref{cor-adicformalschemeidealofdefinitionfinitetype11} we have $\Gamma(U,\mathscr{I}^{k+1})\cdot\Gamma(U,M^{\Delta})=I^{k+1}\widehat{M}_f$ for $U=\Spf A_{\{f\}}$, where $\widehat{M}_f$ is the $I$-adic completion\index{completion!I-adic completion@$I$-adic ---} of $M_f$ ({\bf \ref{ch-pre}}.\ref{prop-Iadiccompletioncomplete1}). 
Since $I^{k+1}\widehat{M}_f$ is a closed subset of $\widehat{M}_f$ ({\bf \ref{ch-pre}}.\ref{cor-Iadiccompletionifexists1}), it is $I$-adically complete\index{complete!I-adically complete@$I$-adically ---} ({\bf \ref{ch-pre}}.\ref{prop-ARconseq1-1}); it then coincides with the $I$-adic completion\index{completion!I-adic completion@$I$-adic ---} of $I^{k+1}M_f$ and hence with $\Gamma(U,(I^{k+1}M)^{\Delta})$, whence the result.
\end{proof}

\begin{proof}[Proof of Proposition {\rm \ref{prop-adicqcohaff1f}}]
By \ref{lem-deltasheafadicallyquasicoherent1-1} and \ref{lem-deltasheafadicallyquasicoherent1-2} we have $M^{\Delta}/\mathscr{I}^{k+1}M^{\Delta}=\til{M}_k$ for $k\geq 0$.
Since 
$$
\varprojlim_{k\geq 0}M^{\Delta}/\mathscr{I}^{k+1}M^{\Delta}=\varprojlim_{k\geq 0}\til{M}_k=M^{\Delta},
$$
the sheaf $M^{\Delta}$ is complete, as desired.
\end{proof}

\begin{prop}\label{prop-lemtopqcoh111}
Let $A$ be an adic ring of finite ideal type, and $M\rightarrow N$ a surjective homomorphism of $A$-modules.
Then the induced morphism $M^{\Delta}\rightarrow N^{\Delta}$ of sheaves on $X=\Spf A$ is surjective.
\end{prop}

\begin{proof}
Let $I\subseteq A$ be a finitely generated ideal of definition.
Since for $k\geq 0$ the induced map $M_k=M/I^{k+1}M\rightarrow N_k=N/I^{k+1}N$ is surjective, $\til{M}_k\rightarrow\til{N}_k$ is surjective.
By \ref{lem-deltasheafadicallyquasicoherent1-2} we have $\til{M}_k=(M^{\Delta})_k=M^{\Delta}/I^{k+1}M^{\Delta}$, etc.
Let $\mathscr{K}_k$ be the kernel of the surjection $(M^{\Delta})_k\rightarrow(N^{\Delta})_k$.
Consider the commutative diagram with exact rows and columns 
$$
\xymatrix@R-1ex@C-2ex{&&0\ar[d]&0\ar[d]\\ &&I^{k+1}(M^{\Delta})_l\ar[d]\ar[r]^{(\ast)}&I^{k+1}(N^{\Delta})_l\ar[d]\\ 0\ar[r]&\mathscr{K}_l\ar[d]\ar[r]&(M^{\Delta})_l\ar[d]\ar[r]&(N^{\Delta})_l\ar[d]\ar[r]&0\\ 0\ar[r]&\mathscr{K}_k\ar[r]&(M^{\Delta})_k\ar[d]\ar[r]&(N^{\Delta})_k\ar[d]\ar[r]&0\\ &&0&0}
$$
for $k\leq l$.
Since the map $(\ast)$ is surjective, we see by snake lemma that the map $\mathscr{K}_l\rightarrow\mathscr{K}_k$ is surjective, that is, the projective system $\{\mathscr{K}_k\}_{k\geq 0}$ of quasi-coherent sheaves is strict\index{strict projective system@strict (projective system)}\index{system!projective system@projective ---!strict projective system@strict --- ---}.
Hence by \ref{lem-vanishingcohomologyadicallyuseful} (1) we see that $M^{\Delta}=\varprojlim_{k\geq 0}(M^{\Delta})_k\rightarrow N^{\Delta}=\varprojlim_{k\geq 0}(N^{\Delta})_k$ is surjective, as desired.
\end{proof}

\begin{cor}\label{cor-topqcoh111}
Let $A$ be an adic ring of finite ideal type, and $M$ a finitely generated $A$-module.
Then $M^{\Delta}$ is an adically quasi-coherent sheaf of finite type.
\end{cor}

\begin{proof}
Take a surjection $A^{\oplus n}\rightarrow M$.
Then we have the surjection $\O^{\oplus n}_X=(A^{\oplus n})^{\Delta}\rightarrow M^{\Delta}$.
\end{proof}

\begin{prop}\label{prop-topqcoh1}
Let $A$ be an adic ring with a finitely generated ideal of definition $I\subseteq A$, and set $X=\Spf A$. 
Let $\mathscr{F}$ $($resp.\ $\mathscr{B})$ be an $\O_X$-module $($resp.\ an $\O_X$-algebra$)$.
Then the following conditions are equivalent$:$
\begin{itemize}
\item[{\rm (a)}] there exist an affine open covering $\{U_{\alpha}\}_{\alpha\in L}$ of $X$ with $U_{\alpha}=\Spf A_{\{f_{\alpha}\}}$ and for each $\alpha\in L$ an $IA_{\{f_{\alpha}\}}$-adically complete $A_{\{f_{\alpha}\}}$-module $M_{\alpha}$ $($resp.\ $A_{\{f_{\alpha}\}}$-algebra $B_{\alpha})$ such that $\mathscr{F}|_{U_{\alpha}}\cong M^{\Delta}_{\alpha}$ $($resp.\ $\mathscr{B}|_{U_{\alpha}}\cong B^{\Delta}_{\alpha});$
\item[{\rm (b)}] there exists an $I$-adically complete $A$-module $M$ $($resp.\ $A$-algebra $B)$ such that $\mathscr{F}\cong M^{\Delta}$ $($resp.\ $\mathscr{B}\cong B^{\Delta});$
\item[{\rm (c)}] $\mathscr{F}$ $($resp.\ $\mathscr{B})$ is an adically quasi-coherent sheaf $($resp.\ an adically quasi-coherent $\O_X$-algebra$)$ on $X$.
\end{itemize}
Moreover, if $\mathscr{F}\cong M^{\Delta}$ as in {\rm (b)}, $\mathscr{F}$ is of finite type if and only if $M=\Gamma(X,\mathscr{F})$ is finitely generated.
\end{prop}

\begin{proof}
First we show the assertion for the module sheaf $\mathscr{F}$.
The implication (b) $\Rightarrow$ (a) is trivial, and (a) $\Rightarrow$ (c) follows from \ref{prop-adicqcohaff1f}.
Let us show (c) $\Rightarrow$ (b).
Set $\mathscr{I}=I^{\Delta}$ and $\mathscr{F}_k=\mathscr{F}/\mathscr{I}^{k+1}\mathscr{F}$ for $k\geq 0$; each $\mathscr{F}_k$ is a quasi-coherent sheaf on the scheme $X_k=\Spec A/I^{k+1}$.
Let $M_k=\Gamma(X_k,\mathscr{F}_k)$.
These $A$-modules constitute a strict projective system\index{strict projective system@strict (projective system)}\index{system!projective system@projective ---!strict projective system@strict --- ---} $\{M_k\}_{k\geq 0}$.
Set $M=\varprojlim_{k\geq 0}M_k$, and let $F^{(n)}$ be the kernel of the surjective map $M\rightarrow M_{n-1}$ for each $n\geq 1$.
Then by {\bf \ref{ch-pre}}.\ref{prop-criterionadicness2} we have $F^{(n)}=I^nM$ for $n\geq 1$, and hence $M$ is $I$-adically complete.
Moreover, by \ref{lem-deltasheafadicallyquasicoherent1-2} we have
$$
\mathscr{F}=\varprojlim_{k\geq 0}\mathscr{F}_k=\varprojlim_{k\geq 0}\til{M}_k=\varprojlim_{k\geq 0}M^{\Delta}/\mathscr{I}^{k+1}M^{\Delta}=M^{\Delta}.
$$
The last assertion follows from \ref{cor-topqcoh111} and the following observation: if $\mathscr{F}$ is of finite type, then $M_0=M/F^{(1)}$ is finitely generated, and hence by the last part of {\bf \ref{ch-pre}}.\ref{prop-criterionadicness2} $M$ is finitely generated.

Also for the algebra sheaf $\mathscr{B}$ the implications (b) $\Rightarrow$ (a) and (a) $\Rightarrow$ (c) are clear.
To show (c) $\Rightarrow$ (b), we set $\mathscr{B}_k=\mathscr{B}/\mathscr{I}^{k+1}\mathscr{B}$ and $B_k=\Gamma(X_k,\mathscr{B}_k)$ for each $k\geq 0$ and define $B=\varprojlim_{k\geq 0}B_k$.
We know that the map $B_l\rightarrow B_k$ for $k\leq l$ is surjective with the kernel equal to $I^{k+1}B_l$. 
Hence by {\bf \ref{ch-pre}}.\ref{prop-exeramazingfactonalgebra} we see that $B$ is $IB$-adically complete and that $B/I^{k+1}=B_k$ for $k\geq 0$.
The rest of the argument is similar to the module case.
\end{proof}

\begin{thm}\label{thm-adicqcoh1}
Let $X=\Spf A$, where $A$ is an adic ring of finite ideal type, and consider the functor 
$$
M\longmapsto M^{\Delta}\leqno{(\ast)}
$$
as in {\rm \S\ref{subsub-thedeltasheaves}}.

{\rm (1)} If $M$ is an $I$-adically complete $A$-module, then we have 
$$
\Gamma(X,M^{\Delta})=M.
$$

{\rm (2)} The functor $(\ast)$ gives a categorical equivalence between the category of $I$-adically complete $A$-modules $($resp.\  $I$-adically complete finitely generated $A$-modules, resp.\ $I$-adically complete $A$-algebras$)$ and the category of adically quasi-coherent sheaves $($resp.\ adically quasi-coherent sheaves of finite type, resp.\ adically quasi-coherent $\O_X$-algebras$)$ on $X$.
The quasi-inverse functor is given by $\mathscr{F}\mapsto\Gamma(X,\mathscr{F})$. \hfill$\square$
\end{thm}

The assertion (1) follows from \ref{prop-adicqcohaff1g}.
The other assertion follows from \ref{prop-topqcoh1} and the following lemma:
\begin{lem}\label{lem-corcoradicqcoh1}
In the situation as in {\rm \ref{thm-adicqcoh1}}, let $M$ and $N$ be $I$-adically complete $A$-modules.
Then the canonical map
$$
\Hom_A(M,N)\stackrel{\cdot^{\Delta}}{\longrightarrow}\Hom_{\O_X}(M^{\Delta},N^{\Delta})\leqno{(\ast\ast)}
$$
is an isomorphism.
The analogous statement with $M,N$ replaced by $I$-adically complete $A$-algebras and the morphisms replaced by algebra homomorphisms is also true.
\end{lem}

\begin{proof}
Given a morphism $\varphi\colon M^{\Delta}\rightarrow N^{\Delta}$, we have $\Gamma_X(\varphi)\colon M\rightarrow N$ and thus $\Gamma_X(\varphi)^{\Delta}\colon M^{\Delta}\rightarrow N^{\Delta}$.
We claim that $\varphi=\Gamma_X(\varphi)^{\Delta}$.
To show this, we first show $\Gamma_X(\varphi_k)=\Gamma_X(\varphi)\otimes_AA_k$ for $k\geq 0$, where $\varphi_k\colon\til{M}_k=M^{\Delta}/\mathscr{I}^{k+1}M^{\Delta}\rightarrow\til{N}_k=N^{\Delta}/\mathscr{I}^{k+1}N^{\Delta}$ (where $\mathscr{I}=I^{\Delta}$) is the induced map (cf.\ \ref{lem-deltasheafadicallyquasicoherent1-2}) and $A_k=A/I^{k+1}$.
By \ref{lem-vanishingcohomologyadicallyuseful} (2) we have the commutative diagram with exact rows
$$
\xymatrix{0\ar[r]&\Gamma(X,\mathscr{I}^{k+1}M^{\Delta})\ar[d]\ar[r]&\Gamma(X,M^{\Delta})\ar[d]_{\Gamma(\varphi)}\ar[r]&\Gamma(X,\til{M}_k)\ar[d]^{\Gamma(\varphi_k)}\ar[r]&0\\ 0\ar[r]&\Gamma(X,\mathscr{I}^{k+1}N^{\Delta})\ar[r]&\Gamma(X,N^{\Delta})\ar[r]&\Gamma(X,\til{N}_k)\ar[r]&0\rlap{;}}
$$
here we used the fact that $\mathscr{I}^{k+1}M^{\Delta}$ and $\mathscr{I}^{k+1}N^{\Delta}$ are adically quasi-coherent sheaves on $X$ due to \ref{prop-adicallyquasicohstr} (3) and hence have vanishing higher cohomologies due to \ref{lem-vanishingcohomologyadicallyuseful} (2).
By \ref{lem-deltasheafadicallyquasicoherent1-2} we have $\Gamma(X,\mathscr{I}^{k+1}M^{\Delta})=I^{k+1}M$ and $\Gamma(X,\mathscr{I}^{k+1}N^{\Delta})=I^{k+1}N$ (here we used the assumption that $M$ and $N$ are finitely generated) and hence $\Gamma_X(\varphi_k)=\Gamma_X(\varphi)\otimes_AA_k$, as desired.
Then, 
$$
\Gamma_X(\varphi)^{\Delta}=\varprojlim_{k\geq 0}\til{\Gamma_X(\varphi)\otimes_AA_k}=\varprojlim_{k\geq 0}\til{\Gamma_X(\varphi_k)}=\varprojlim_{k\geq 0}\varphi_k=\varphi,
$$
which shows our claim.
In particular, we have shown that the map $(\ast\ast)$ in question is surjective.
Since for any homomorphism $f\colon M\rightarrow N$ of $A$-modules we clearly have $\Gamma_X(f^{\Delta})=f$ (by {\bf \ref{ch-pre}}.\ref{prop-projlimsheafleftexact0}), $(\ast\ast)$ is injective.
This concludes the proof of the first statement.
The other assertion (for morphisms between $I$-adically complete $A$-algebras) is shown similarly.
\end{proof}

\subsection{A.q.c.\ algebras of finite type}\label{sub-adicallyqcohalgebraoffinitetype}
\begin{dfn}\label{dfn-adicallyqcohalgebraoffinitetype}{\rm 
Let $X$ be an adic formal scheme of finite ideal type.
We say that an adically quasi-coherent $\O_X$-algebra $\mathscr{B}$ is {\em of finite type}\index{adically quasi-coherent (a.q.c.) algebra!adically quasi-coherent (a.q.c.) algebra of finite type@--- of finite type} if the following condition is satisfied: there exists an open covering $X=\bigcup_{\alpha\in L}U_{\alpha}$ and an ideal of definition of finite type $\mathscr{I}_{\alpha}$ on each $U_{\alpha}$ such that $\mathscr{B}/\mathscr{I}_{\alpha}\mathscr{B}$ for each $\alpha$ is a quasi-coherent algebra of finite type on the scheme $(U_{\alpha},\O_{U_{\alpha}}/\mathscr{I}_{\alpha})$.}
\end{dfn}

By {\bf \ref{ch-pre}}.\ref{prop-formalnot3} (applied to the case where $I$ is nilpotent) we readily see that the required condition does not depend on the choice of the open covering and the ideal of definitions. 
In particular, we have the following:
\begin{prop}\label{prop-adicallyqcohalgebraoffinitetype1}
Let $X$ be an adic formal scheme of finite ideal type, and $\mathscr{B}$ an adically quasi-coherent $\O_X$-algebra of finite type.
Then for any open subspace $U\subseteq X$ that admits an ideal of definition of finite type $\mathscr{I}$, $\mathscr{B}/\mathscr{I}\mathscr{B}$ is a quasi-coherent algebra of finite type on the scheme $(U,\O_U/\mathscr{I})$. \hfill$\square$
\end{prop}

\begin{prop}\label{prop-adicallyqcohalgebraoffinitetype2}
Let $A$ be an adic ring of finite ideal type, and set $X=\Spf A$.
Let $\mathscr{B}$ be an adically quasi-coherent $\O_X$-algebra.
Then $\mathscr{B}$ is of finite type if and only if $B=\Gamma(X,\mathscr{B})$ is topologically finitely generated\index{finitely generated!topologically finitely generated@topologically ---} $A$-algebra.
\end{prop}

\begin{proof}
In view of \ref{thm-adicqcoh1} (2) $B$ is an $I$-adically complete $A$-algebra where $I\subseteq A$ is a finitely generated ideal of definition.
Since $\til{B/IB}=\mathscr{B}/\mathscr{I}\mathscr{B}$ on $\Spec A/I$, the result follows from {\bf \ref{ch-pre}}.\ref{prop-formalnot3}.
\end{proof}

\subsection{A.q.c.\ sheaves as projective limits}\label{sub-adicallyqcohinductivelimit}
\begin{prop}\label{prop-adicallyqcohinductivelimit}
Let $X$ be an adic formal scheme of finite ideal type, and $\mathscr{I}$ an ideal of definition of finite type.
Set $X_k=(X,\O_X/\mathscr{I}^{k+1})$.
Suppose we have a projective system $\{\mathscr{F}_k,\varphi_{ij}\}_{i\in\N}$ of $\O_X$-modules such that$:$
\begin{itemize}
\item[{\rm (a)}] for each $k\geq 0$ we have $\mathscr{I}^{k+1}\mathscr{F}_k=0$, and the sheaf $\mathscr{F}_k$ is a quasi-coherent sheaf on the scheme $X_k;$
\item[{\rm (b)}] for any $i\leq j$ the morphism $\varphi_{ij}\colon\mathscr{F}_j\rightarrow\mathscr{F}_i$ is a surjective map with the kernel equal to $\mathscr{I}^{i+1}\mathscr{F}_j;$ 
\end{itemize}
Then the projective limit $\mathscr{F}=\varprojlim_k\mathscr{F}_k$ is an adically quasi-coherent sheaf on $X$ such that $\mathscr{F}/\mathscr{I}^{k+1}\mathscr{F}\cong\mathscr{F}_k$ for each $k\geq 0$.
Moreover, if $\mathscr{F}_0$ is of finite type, then $\mathscr{F}$ is of finite type.
\end{prop}

\begin{proof}
We may assume that $X$ is affine $X=\Spf A$ with $\mathscr{I}=I^{\Delta}$, where $I$ is a finitely generated ideal of definition of $A$.
We have $X_k=\Spec A_k$ where $A_k=A/I^{k+1}$ for any $k\geq 0$.
Take for each $k$ the $A_k$-module $M_k$ such that $\mathscr{F}_k=\til{M}_k$.
Then we have the projective system $\{M_k,f_{ij}\}$ of $A$-modules such that for each $i\leq j$ the transition map $f_{ij}\colon M_j\rightarrow M_i$ is surjective with the kernel $I^{i+1}M_j$.
Set $M=\varprojlim_kM_k$ and denote by $f_i\colon M\rightarrow M_i$ the projection map for each $i$.
Set $F^{(n)}=\ker(f_{n-1})$ for $n\geq 1$.
Then by {\bf \ref{ch-pre}}.\ref{prop-criterionadicness2} $M$ is an $I$-adically complete finitely generated $A$-module, and $F^{(n)}=I^nM$ holds for $n\geq 1$.
If $\mathscr{F}_0$ is of finite type, then $M_0=M/F^{(1)}$ is finitely generated, and hence $M$ is finitely generated.

Now consider the open subset $U=\Spf A_{\{g\}}$ of $X$ for any $g\in A$; we have
$$
\Gamma(U,\mathscr{F})=\varprojlim_k\Gamma(U,\mathscr{F}_k)=\varprojlim_k(M_k\otimes_AA_g)=M\widehat{\otimes}_AA_{\{g\}},
$$
which coincides with $\Gamma(U,M^{\Delta})$, as calculated in \ref{prop-adicqcohaff1g}.
Hence we have $\mathscr{F}=M^{\Delta}$, which is an adically quasi-coherent sheaf (\ref{prop-adicqcohaff1f}).
Moreover, by \ref{lem-deltasheafadicallyquasicoherent1-2} we have
$$
\mathscr{F}/\mathscr{I}^{k+1}\mathscr{F}=M^{\Delta}/\mathscr{I}^{k+1}M^{\Delta}=\til{M}_k=\mathscr{F}_k,
$$
as desired.
\end{proof}

One can show the following analogous statement for adically quasi-coherent algebras; the proof is similar (cf.\ the proof of \ref{prop-topqcoh1}):
\begin{prop}\label{prop-adicallyqcohinductivelimitalgebra}
Let $X$ be an adic formal scheme of finite ideal type, and $\mathscr{I}$ an ideal of definition of finite type.
Set $X_k=(X,\O_X/\mathscr{I}^{k+1})$.
Suppose we have a projective system $\{\mathscr{B}_k,\varphi_{ij}\}_{i\in\N}$ of $\O_X$-algebras such that$:$
\begin{itemize}
\item[{\rm (a)}] for each $k\geq 0$ we have $\mathscr{I}^{k+1}\mathscr{B}_k=0$, and the sheaf $\mathscr{B}_k$ is a quasi-coherent algebra on the scheme $X_k;$
\item[{\rm (b)}] for any $i\leq j$ the morphism $\varphi_{ij}\colon\mathscr{B}_j\rightarrow\mathscr{B}_i$ is a surjective map with the kernel equal to $\mathscr{I}^{i+1}\mathscr{B}_j$. 
\end{itemize}
Then the projective limit $\mathscr{B}=\varprojlim_k\mathscr{B}_k$ is an adically quasi-coherent $\O_X$-algebra on $X$ such that $\mathscr{B}/\mathscr{I}^{k+1}\mathscr{B}\cong\mathscr{B}_k$ for each $k\geq 0$.
Moreover, if $\mathscr{B}_0$ is of finite type, then $\mathscr{B}$ is an adically quasi-coherent $\O_X$-algebra of finite type. \hfill$\square$
\end{prop}

\subsection{A.q.c.\ sheaves on locally universally rigid-Noetherian formal schemes}\label{sub-adicqcohpf}
\index{formal scheme!universally rigid-Noetherian formal scheme@universally rigid-Noetherian ---!locally universally rigid-Noetherian formal scheme@locally --- ---}
\subsubsection{$\Delta$-sheaves on affine universally rigid-Noetherian formal schemes}\label{subsub-deltasheavesonadequate}
\begin{prop}\label{prop-adequatedeltaconstruction}
Let $A$ be a rigid-Noetherian ring\index{rigid-Noetherian ring@rigid-Noetherian ring} {\rm (\ref{dfn-tuaringadmissible} (1))}, and set $X=\Spf A$ and $Y=\Spec A$.
We consider the canonical morphism $i\colon X\rightarrow Y$ of locally ringed spaces\index{morphism of ringed spaces@morphism (of ringed spaces)!local morphism of ringed spaces@local ---}.
Then for any finitely generated $A$-module $M$ we have
$$
M^{\Delta}\cong i^{\ast}\til{M}.
$$
\end{prop}

\begin{proof}
By \ref{prop-adqformalpropaffine} we have 
$$
i^{\ast}\til{M}\cong\widehat{\til{M}}=\varprojlim_{k\geq 0}\til{M}/\mathscr{I}^{k+1}\til{M}=M^{\Delta},
$$
where $\mathscr{I}=I^{\Delta}$ and $I\subseteq A$ is a finitely generated ideal of definition.
\end{proof}

\begin{prop}\label{prop-coradequatedeltaconstruction}
In the situation as in {\rm \ref{prop-adequatedeltaconstruction}}, the functor $M\mapsto M^{\Delta}$ gives an {\em exact} equivalence between the category of finitely generated $A$-modules to the category of adically quasi-coherent sheaves of finite type on $X$.
\end{prop}

\begin{proof}
In view of \ref{thm-adicqcoh1} (2) and the fact that any finitely generated $A$-module is $I$-adically complete (\ref{rem-turigidnoetherianbasicproperties}), only the exactness of the functor is in question.
But this follows from \ref{prop-adequatedeltaconstruction} and \ref{prop-adqformalpropaffine}.
\end{proof}

\begin{prop}[{cf.\ {\bf \ref{ch-pre}}.\ref{prop-corpropARconseq1-1-1}; see also \ref{lem-deltasheafadicallyquasicoherent1-1}}]\label{prop-lemdeltasheafadicallyquasicoherent1-3}
Let $A$ be a rigid-Noetherian ring\index{rigid-Noetherian ring@rigid-Noetherian ring} with a finitely generated ideal of definition $I\subseteq A$, and $M$ a finitely generated $A$-module.
Let $N\subseteq M$ be an $A$-submodule.
Then the sequence 
$$
0\longrightarrow N^{\Delta}\longrightarrow M^{\Delta}\longrightarrow(M/N)^{\Delta}\longrightarrow 0
$$
by canonical maps is exact.
\end{prop}

\begin{proof}
Let $f\in A$, and consider the affine open subset $U=\Spf A_{\{f\}}$ of $X$.
Since $A_f$ is $IA_f$-adically pseudo-adhesive\index{pseudo-adhesive} ({\bf \ref{ch-pre}}.\ref{prop-adhesive1} (1)), it satisfies the condition {\bf (AP)} in {\bf \ref{ch-pre}}, \S\ref{subsub-ARIgoodnesssub}.
Hence by {\bf \ref{ch-pre}}.\ref{prop-corpropARconseq1-1-1} we have the exact sequence
$$
0\longrightarrow N\widehat{\otimes}_AA_{\{f\}}\longrightarrow M\widehat{\otimes}_AA_{\{f\}}\longrightarrow (M/N)\widehat{\otimes}_AA_{\{f\}}\longrightarrow 0
$$
and the exact sequence
$$
0\longrightarrow\Gamma(U,N^{\Delta})\longrightarrow\Gamma(U,M^{\Delta})\longrightarrow\Gamma(U,(M/N)^{\Delta})\longrightarrow 0
$$
by \ref{prop-adicqcohaff1g}.
Since this is valid for any open subsets of the form $U=\Spf A_{\{f\}}$, we have the desired exact sequence.
\end{proof}

\subsubsection{Adically quasi-coherent sheaves of finite presentation}\label{subsub-adicqcohpf}
\begin{prop}\label{prop-adicqcoh3bis}
Let $A$ be a rigid-Noetherian ring\index{rigid-Noetherian ring@rigid-Noetherian ring} {\rm (\ref{dfn-tuaringadmissible} (1))}, and set $X=\Spf A$.
Then any finitely presented $\O_X$-module\index{module (OX)@module ($\O_X$-{---})!module of finite presentation@--- of finite presentation} $\mathscr{F}$ is an adically quasi-coherent sheaf.
Moreover, $M=\Gamma(X,\mathscr{F})$ is a finitely presented $A$-module, and $\mathscr{F}=M^{\Delta}$.
\end{prop}

\begin{proof}
The sheaf $\mathscr{F}$ is locally isomorphic to the cokernel of a map of the form $\O^{\oplus q}_X\rightarrow\O^{\oplus p}_X$.
By \ref{prop-coradequatedeltaconstruction} this means that $\mathscr{F}$ can be locally written as a $\Delta$-sheaf of a finitely generated module.
Hence $\mathscr{F}$ is adically quasi-coherent of finite type.
Moreover, we know that each $\mathscr{F}_k=\mathscr{F}/I^{k+1}\mathscr{F}$ (where $I\subseteq A$ is a finitely generated ideal of definition) is finitely presented on the scheme $\Spec A/I^{k+1}$.
This means that $M=\Gamma(X,\mathscr{F})$ is finitely presented due to {\bf \ref{ch-pre}}.\ref{cor-ARconseq2}, since $M_k=M/I^{k+1}M$ is finitely presented by \cite[(1.4.3)]{EGAInew}.
\end{proof}

\begin{cor}\label{cor-propadicqcoh3bis1}
Any finitely presented $\O_X$-module\index{module (OX)@module ($\O_X$-{---})!module of finite presentation@--- of finite presentation} on a locally universally rigid-Noetherian formal scheme is adically quasi-coherent of finite type. \hfill$\square$
\end{cor}

We can call such a sheaf $\mathscr{F}$, as usual, a {\em quasi-coherent sheaf of finite presentation} (even without `adically'), because it is obviously quasi-coherent in the usual sense (cf.\ {\bf \ref{ch-pre}}.\ref{dfn-qcohcohringedspace} (1)).

\begin{thm}\label{thm-adicqcohpre1}
Let $X=\Spf A$, where $A$ is a rigid-Noetherian ring.
Then the functor 
$$
M\longmapsto M^{\Delta}\leqno{(\ast)}
$$
gives an exact categorical equivalence between the category of finitely generated $($resp.\ finitely presented$)$ $A$-modules and the category of adically quasi-coherent sheaves of finite type $($resp.\ of finite presentation$)$ on $X$.
The functor $\mathscr{F}\mapsto\Gamma(X,\mathscr{F})$ gives the quasi-inverse to $(\ast)$.
\end{thm}

\begin{proof}
The assertion for `finitely generated' is already shown in \ref{prop-coradequatedeltaconstruction}.
The other part follows from \ref{prop-coradequatedeltaconstruction} and \ref{prop-adicqcoh3bis}.
\end{proof}

\begin{prop}\label{prop-adicqcoh7}
Let $A$ be a rigid-Noetherian ring\index{rigid-Noetherian ring@rigid-Noetherian ring}, and $M$ and $N$ finitely generated $A$-modules. 
Set $X=\Spf A$.

{\rm (1)} We have a canonical isomorphism
$$
(M\otimes_AN)^{\Delta}\cong M^{\Delta}\otimes_{\O_X}N^{\Delta}.
$$

{\rm (2)} If $M$ is finitely presented, then we have a canonical isomorphism
$$
\Hom_A(M,N)^{\Delta}\cong\lHom_{\O_X}(M^{\Delta},N^{\Delta}).
$$
\end{prop}

\begin{proof}
Set $Y=\Spec A$, and let $i\colon X\rightarrow Y$ be the canonical morphism of locally ringed spaces.

(1) By \ref{prop-adequatedeltaconstruction}, \cite[$\mathbf{I}$, (1.3.12) (i)]{EGA}, and \cite[$\mathbf{0}_{\mathbf{I}}$, (4.3.3.1)]{EGA} we have
$$
(M\otimes_AN)^{\Delta}\cong i^{\ast}\til{M\otimes_AN}\cong i^{\ast}(\til{M}\otimes_{\O_Y}\til{N})\cong i^{\ast}\til{M}\otimes_{\O_X}i^{\ast}\til{N}\cong M^{\Delta}\otimes_{\O_X}N^{\Delta},
$$
as claimed.

(2) By \ref{prop-adequatedeltaconstruction}, \cite[$\mathbf{I}$, (1.3.12) (ii)]{EGA}, and \cite[$\mathbf{0}_{\mathbf{I}}$, (6.7.6.1)]{EGA} we have
$$
\begin{aligned}
\Hom_A(M,N)^{\Delta}\cong i^{\ast}\til{\Hom_A(M,N)}\cong i^{\ast}\lHom_{\O_Y}(\til{M},\til{N})\cong\lHom_{\O_X}(i^{\ast}\til{M},i^{\ast}\til{N})\\ \cong\lHom_{\O_X}(M^{\Delta},N^{\Delta}),
\end{aligned}
$$
as claimed; here we used the fact that $\til{M}$ is of finite presentation (\ref{thm-adicqcohpre1}).
\end{proof}

\begin{cor}\label{cor-adicqcoh711}
Let $X=\Spf A$ where $A$ is a rigid-Noetherian ring.\index{rigid-Noetherian ring@rigid-Noetherian ring}
Then the functor $(\ast)$ in {\rm \ref{thm-adicqcohpre1}} defined on the category of finitely presented $A$-modules preserves tensor products and internal $\Hom$'s. \hfill$\square$
\end{cor}

\begin{cor}\label{cor-adicqcoh71}
Let $X$ be a locally universally rigid-Noetherian formal scheme, and $\mathscr{F},\mathscr{G}$ adically quasi-coherent sheaves of finite type on $X$.
Then $\mathscr{F}\otimes_{\O_X}\mathscr{G}$ is an adically quasi-coherent sheaf of finite type.
If, moreover, $\mathscr{F}$ is of finite presentation, then $\lHom_{\O_X}(\mathscr{F},\mathscr{G})$ is an adically quasi-coherent sheaf of finite type on $X$.
Moreover, if both $\mathscr{F}$ and $\mathscr{G}$ are finitely presentated, then $\mathscr{F}\otimes_{\O_X}\mathscr{G}$ and $\lHom_{\O_X}(\mathscr{F},\mathscr{G})$ are of finite presentation. \hfill$\square$
\end{cor}

Finally, let us mention that $\mathscr{I}$-torsion free adically quasi-coherent sheaves of finite type on {\em locally universally adhesive}\index{formal scheme!universally adhesive formal scheme@universally adhesive ---!locally universally adhesive formal scheme@locally --- ---} formal schemes are automatically finitely presented:
\begin{prop}\label{prop-exerItorsionfreeaqcsheavesoffinitetype}
Let $X$ be a locally universally adhesive formal scheme\index{formal scheme!universally adhesive formal scheme@universally adhesive ---!locally universally adhesive formal scheme@locally --- ---}.
Suppose that $X$ has an ideal of definition of finite type $\mathscr{I}\subseteq\O_X$.
Then any $\mathscr{I}$-torsion free adically quasi-coherent sheaf of finite type is finitely presented.
\end{prop}

\begin{proof}
We may assume $X=\Spf A$ where $A$ is a t.u.\ adhesive ring\index{t.u.a. ring@t.u.\ adhesive ring} with a finitely generated ideal of definition $I\subseteq A$.
Any $\mathscr{I}$-torsion free adically quasi-coherent sheaf of finite type $\mathscr{F}$ is of the form $\mathscr{F}=M^{\Delta}$, where $M$ is an $I$-torsion free finitely generated $A$-module.
Since $A$ is $I$-adically adhesive, $M$ is finitely presented (cf.\ {\bf \ref{ch-pre}}.\ref{prop-adhesive}).
\end{proof}

\subsubsection{Adically quasi-coherent algebras of finite presentation}\label{subsub-adicallyqcohalgebrasoffinitepres}
\begin{dfn}\label{dfn-adicallyqcohalgebraoffinitepres}{\rm 
Let $X$ be a locally universally rigid-Noetherian formal scheme\index{formal scheme!universally rigid-Noetherian formal scheme@universally rigid-Noetherian ---!locally universally rigid-Noetherian formal scheme@locally --- ---}.
We say that an adically quasi-coherent $\O_X$-algebra $\mathscr{B}$ is {\em of finite presentation}\index{adically quasi-coherent (a.q.c.) algebra!adically quasi-coherent (a.q.c.) algebra of finite presentation@--- of finite presentation} if the following condition is satisfied: there exists an open covering $X=\bigcup_{\alpha\in L}U_{\alpha}$ and an ideal of definition of finite type $\mathscr{I}_{\alpha}$ on each $U_{\alpha}$ such that,\ for each $\alpha$ and $k\geq 0$ the sheaf $\mathscr{B}/\mathscr{I}^{k+1}_{\alpha}\mathscr{B}$ is a quasi-coherent algebra of finite presentation on the scheme $(U_{\alpha},\O_{U_{\alpha}}/\mathscr{I}^{k+1}_{\alpha})$.}
\end{dfn}

One can show the following propositions similarly to \ref{prop-adicallyqcohalgebraoffinitetype1} and \ref{prop-adicallyqcohalgebraoffinitetype2}, using {\bf \ref{ch-pre}}.\ref{prop-formalnot31} instead of {\bf \ref{ch-pre}}.\ref{prop-formalnot3}:
\begin{prop}\label{prop-adicallyqcohalgebraoffinitepres1}
Let $X$ be a locally universally rigid-Noetherian formal scheme\index{formal scheme!universally rigid-Noetherian formal scheme@universally rigid-Noetherian ---!locally universally rigid-Noetherian formal scheme@locally --- ---}, and $\mathscr{B}$ an adically quasi-coherent $\O_X$-algebra of finite presentation.
Then for any open subspace $U\subseteq X$ that admits an ideal of definition of finite type $\mathscr{I}$, the sheaf $\mathscr{B}/\mathscr{I}\mathscr{B}$ is a quasi-coherent algebra of finite presentation on the scheme $(U,\O_U/\mathscr{I})$. \hfill$\square$
\end{prop}

\begin{prop}\label{prop-adicallyqcohalgebraoffinitepres2}
Let $A$ be a t.u.\ rigid-Noetherian ring\index{t.u. rigid-Noetherian ring@t.u.\ rigid-Noetherian ring}, and set $X=\Spf A$.
Let $\mathscr{B}$ be an adically quasi-coherent $\O_X$-algebra.
Then $\mathscr{B}$ is of finite presentation if and only if $B=\Gamma(X,\mathscr{B})$ is a topologically finitely presented\index{finitely presented!topologically finitely presented@topologically ---} $A$-algebra. \hfill$\square$
\end{prop}

\subsection{Complete pull-back of a.q.c.\ sheaves}\label{sub-completepullbackaqcsheaves}
\begin{dfn}\label{dfn-completepullback}{\rm 
Let $f\colon X\rightarrow Y$ be an adic morphism of adic formal schemes of finite ideal type, and $\mathscr{F}$ an $\O_Y$-module.
We define the $\O_X$-module $\widehat{f^{\ast}}\mathscr{F}$ by
$$
\widehat{f^{\ast}}\mathscr{F}=\widehat{f^{\ast}\mathscr{F}},
$$
where the last sheaf is the completion\index{completion} of the $\O_X$-module $f^{\ast}\mathscr{F}$ (\ref{dfn-completesheaf} (1)); we call this sheaf the {\em complete pull-back}\index{complete!complete pull-back@--- pull-back} of $\mathscr{F}$.}
\end{dfn}

\begin{prop}\label{prop-completepullbackaqcsheaves0}
Let $f\colon X\rightarrow Y$ be an adic morphism of locally universally rigid-Noetherian formal schemes, and $\mathscr{F}$ an adically quasi-coherent sheaf of finite type $($resp.\ of finite presentation$)$ on $Y$.
Then we have
$$
\widehat{f^{\ast}}\mathscr{F}=f^{\ast}\mathscr{F},
$$
which is an adically quasi-coherent sheaf of finite type $($resp.\ of finite presentation$)$ on $X$.
\end{prop}

\begin{prop}\label{prop-completepullbackaqcsheaves1}
Let $f\colon X\rightarrow Y$ be an adic morphism of adic formal schemes of finite ideal type, and $\mathscr{F}$ an adically quasi-coherent sheaf $($resp.\ adically quasi-coherent sheaf of finite type, resp.\ adically quasi-coherent $\O_Y$-algebra$)$ on $Y$.
Suppose that $f$ is adically flat\index{morphism of formal schemes@morphism (of formal schemes)!adically flat morphism of formal schemes@adically flat ---} {\rm (\ref{exas-adicalizationexamples} (1))}.
Then the complete pull-back $\widehat{f^{\ast}}\mathscr{F}$ is an adically quasi-coherent sheaf $($resp.\ adically quasi-coherent sheaf of finite type, resp.\ adically quasi-coherent $\O_X$-algebra$)$ on $X$.
\end{prop}

To show the propositions, since the question is local, we may work in the affine situation.
Then the propositions follow from the following lemma:
\begin{lem}\label{lem-completepullbackaqcsheaves1lem}
Let $B\rightarrow A$ be an adic homomorphism of adic rings of finite ideal type, and $I\subseteq B$ a finitely generated ideal of definition of $B$.

{\rm (1)} Suppose $A$ and $B$ are t.u.\ rigid-Noetherian rings\index{t.u. rigid-Noetherian ring@t.u.\ rigid-Noetherian ring}, and let $M$ be a finitely generated $($resp.\ finitely presented$)$ $B$-module.
Then we have $f^{\ast}M^{\Delta}\cong\widehat{f^{\ast}}M^{\Delta}\cong(M\otimes_BA)^{\Delta}$. 

{\rm (2)} Suppose that for any $k\geq 0$ the induced map $B/I^{k+1}\rightarrow A/I^{k+1}A$ is flat, and let $M$ be a $B$-module $($resp.\ $B$-algebra$)$.
Then we have $\widehat{f^{\ast}}M^{\Delta}\cong(M\otimes_BA)^{\Delta}$. 
\end{lem}

\begin{proof}
The assertion (1) follows easily, for in this case $M^{\Delta}$ is the pull-back of $\til{M}$ on $\Spec B$ by the map $\Spf B\rightarrow\Spec B$ (\ref{prop-adequatedeltaconstruction}), and thus the pull-back $f^{\ast}M^{\Delta}$ is the pull-back of $\til{M\otimes_BA}$ by the map $\Spf A\rightarrow\Spec A$, which is already complete.

To show the assertion (2), we look at the canonical morphism
$$
f^{\ast}\mathscr{I}^{k+1}M^{\Delta}\longrightarrow\mathscr{I}^{k+1}f^{\ast}M^{\Delta}\leqno{(\ast)}
$$
for $k\geq 0$, where $\mathscr{I}=I^{\Delta}$.
We first want to show that this is an isomorphism by checking that the map between stalks at each point is bijective.
To this end, it suffices to show that for any $h\in A$ the map
$$
(I^{k+1}M)\widehat{\otimes}_BA_{\{h\}}\longrightarrow I^{k+1}(M\widehat{\otimes}_BA_{\{h\}})
$$
is an isomorphism, which amounts to showing that the descending filtration $\{J^{(n)}\}_{n\geq 1}$ on $M\widehat{\otimes}_BA_{\{h\}}$ defined by 
$$
J^{(n)}=\image((I^nM)\widehat{\otimes}_BA_{\{h\}}\longrightarrow M\widehat{\otimes}_BA_{\{h\}})
$$
for $n\geq 1$ is the $I$-adic filtration.
Since
$$
\varprojlim_{n\geq 1}M\widehat{\otimes}_BA_{\{h\}}/J^{(n)}=\varprojlim_{n\geq 1}(M/I^nM)\otimes_{B/I^n}A/I^nA[{\textstyle \frac{1}{h}}]=M\widehat{\otimes}_BA_{\{h\}},
$$
it suffices to show that the condition (b) of {\bf \ref{ch-pre}}.\ref{prop-criterionadicness2} (resp.\ {\bf \ref{ch-pre}}.\ref{prop-criterionadicness1}; cf.\ {\bf \ref{ch-pre}}.\ref{prop-exeramazingfactonalgebra}) holds.
But since $B/I^n\rightarrow A/I^nA$ is flat, the map
$$
(I^mM/I^nM)\otimes_{B/I^n}A/I^nA[{\textstyle \frac{1}{h}}]\longrightarrow I^m(M/I^nM)\otimes_{B/I^n}A/I^nA[{\textstyle \frac{1}{h}}]
$$
(which is, a priori, surjective) is an isomorphism for $0\leq m\leq n$.
Thus we have shown that the morphism $(\ast)$ is an isomorphism for $k\geq 0$.

Now we calculate
$$
\begin{aligned}
\widehat{f^{\ast}}M^{\Delta}=\varprojlim_{k\geq 0}f^{\ast}M^{\Delta}/\mathscr{I}^{k+1}f^{\ast}M^{\Delta}=\varprojlim_{k\geq 0}f^{\ast}(M^{\Delta}/\mathscr{I}^{k+1}M^{\Delta})\\
=\varprojlim_{k\geq 0}f^{\ast}\til{M}_k=\varprojlim_{k\geq 0}\til{M_k\otimes_{B_k}A_k},
\end{aligned}
$$
where $M_k=M/I^{k+1}M$, which shows that $\widehat{f^{\ast}}M^{\Delta}=(M\otimes_BA)^{\Delta}$, as desired.
\end{proof}

The following is a corollary of the proof of \ref{lem-completepullbackaqcsheaves1lem} (2):
\begin{cor}\label{cor-completepullbackaqcsheaves0}
In the situation as in {\rm \ref{prop-completepullbackaqcsheaves1}}, suppose that $Y$ admits an ideal of definition of finite type $\mathscr{I}$.
Let $f_k\colon X_k=(X,\O_X/\mathscr{I}^{k+1}\O_X)\rightarrow Y_k=(Y,\O_Y/\mathscr{I}^{k+1})$ be the induced morphism of schemes for $k\geq 0$.
Then we have $f^{\ast}\mathscr{I}^{k+1}\mathscr{F}=\mathscr{I}^{k+1}f^{\ast}\mathscr{F}$ for $k\geq 0$ and 
$$
\widehat{f^{\ast}}\mathscr{F}=\varprojlim_{k\geq 0}f^{\ast}_k\mathscr{F}_k,
$$
where $\mathscr{F}_k=\mathscr{F}/\mathscr{I}^{k+1}\mathscr{F}$ for $k\geq 0$. \hfill$\square$
\end{cor}
\index{adically quasi-coherent (a.q.c.) sheaf|)}

\subsection{Admissible ideals}\label{sub-admidealextlem}
\subsubsection{Pull-back of quasi-coherent sheaves on closed subschemes}\label{subsub-pullbackofquasicoherentsheaves}
\begin{prop}\label{prop-openaqcsubmoduleoffinitetype}
Let $X$ be an adic formal scheme, $\mathscr{I}\subseteq\O_X$ an ideal of definition of finite type, and $\mathscr{F}$ an adically quasi-coherent sheaf\index{adically quasi-coherent (a.q.c.) sheaf} on $X$.
Let $\mathscr{G}$ be an adically quasi-coherent subsheaf of $\mathscr{F}$ such that $\mathscr{I}\mathscr{F}\subseteq\mathscr{G}\subseteq\mathscr{F}$.
Then $\mathscr{G}/\mathscr{I}\mathscr{F}$ is a quasi-coherent sheaf on the scheme $X_0=(X,\O_X/\mathscr{I})$.
\end{prop}

\begin{proof}
We may assume that $X$ is affine $X=\Spf A$, where $A$ is an adic ring of finite ideal type, and that $\mathscr{I}=I^{\Delta}$ where $I\subseteq A$ is a finitely generated ideal of definition.
Set $N=\Gamma(X,\mathscr{G})\subseteq M=\Gamma(X,\mathscr{F})$.
In view of \ref{lem-deltasheafadicallyquasicoherent1-2} one has $IM\subseteq N\subseteq M$.
We see by \ref{lem-deltasheafadicallyquasicoherent1-1} that $\mathscr{F}/\mathscr{G}=(M/N)^{\Delta}$.
Since $\mathscr{I}(\mathscr{F}/\mathscr{G})=0$, we deduce that $\mathscr{F}/\mathscr{G}$ is a quasi-coherent sheaf on the scheme $X_0$.
Now we look at the exact sequence
$$
0\longrightarrow\mathscr{G}/\mathscr{I}\mathscr{F}\longrightarrow\mathscr{F}/\mathscr{I}\mathscr{F}\longrightarrow\mathscr{F}/\mathscr{G}\longrightarrow 0.
$$
Since $\mathscr{F}/\mathscr{I}\mathscr{F}$ and $\mathscr{F}/\mathscr{G}$ are quasi-coherent on $X_0$, we conclude that $\mathscr{G}/\mathscr{I}\mathscr{F}$ is quasi-coherent on $X_0$.
\end{proof}

\begin{prop}\label{prop-exeradicallyqcohsheavesbypullback}
Let $X$ be an adic formal scheme of finite ideal type with an ideal of definition $\mathscr{I}$ of finite type, and $\mathscr{F}$ an adically quasi-coherent sheaf $($resp.\ of finite type$)$ on $X$.
Consider the quasi-coherent sheaf $\mathscr{F}_0=\mathscr{F}/\mathscr{I}\mathscr{F}$ on the scheme $X_0=(X,\O_X/\mathscr{I})$ and a quasi-coherent subsheaf $\mathscr{G}_0\subseteq\mathscr{F}_0$ $($resp.\ of finite type$)$.
Then the inverse image $\mathscr{G}\subseteq\mathscr{F}$ of $\mathscr{G}_0$ by the canonical map $\mathscr{F}\rightarrow\mathscr{F}_0$ is an adically quasi-coherent sheaf $($resp.\ of finite type$)$\index{adically quasi-coherent (a.q.c.) sheaf!adically quasi-coherent sheaf of finite type@--- of finite type} on $X$.
\end{prop}

\begin{proof}
We may assume that $X$ is affine $X=\Spf A$ without loss of generality. 
Let $I\subseteq A$ be the finitely generated ideal such that $\mathscr{I}=I^{\Delta}$, and $M$ the $I$-adically complete $A$-module such that $\mathscr{F}=M^{\Delta}$.
The quasi-coherent subsheaf $\mathscr{G}_0\subseteq\mathscr{F}_0$ corresponds to a submodule $N_0\subseteq M_0=M/IM$.
Let $N$ be the pull-back of $N_0$ by the canonical map $M\rightarrow M_0$.
Since $IM\subseteq N$, we have the exact sequence
$$
0\longrightarrow N^{\Delta}\longrightarrow M^{\Delta}\longrightarrow(M/N)^{\Delta}\longrightarrow 0
$$
(\ref{lem-deltasheafadicallyquasicoherent1-1}).
But since $(M/N)^{\Delta}=\til{M_0/N_0}=\mathscr{F}_0/\mathscr{G}_0$ and since 
$$
0\longrightarrow\mathscr{G}\longrightarrow\mathscr{F}\longrightarrow\mathscr{F}_0/\mathscr{G}_0\longrightarrow 0
$$
is exact, we have $\mathscr{G}=N^{\Delta}$.
Notice that, since $N$ is open in $M$ and $M$ is $I$-adically complete, $N$ is $I$-adically complete.
If $\mathscr{F}$ and $\mathscr{G}_0$ are of finite type, then $M$ and $N_0$ are finitely generated.
Since $N/N\cap I^2M=N/I^2M$ is finitely generated due to the exact sequence
$$
0\longrightarrow IM/I^2M\longrightarrow N/I^2M\longrightarrow N/IM\ (=N/N\cap IM=N_0)\longrightarrow 0, 
$$
we deduce that $N/IN$ is finitely generated (for $I^2M\subseteq IN$).
By {\bf \ref{ch-pre}}.\ref{prop-complpair1} $N$ is finitely generated, and hence $\mathscr{G}$ is of finite type.
\end{proof}

\begin{cor}\label{cor-openaqcsubmoduleoffinitetype}
Let $X$ be an adic formal scheme, $\mathscr{I}\subseteq\O_X$ an ideal of definition of finite type, and $\mathscr{F}$ an adically quasi-coherent sheaf $($resp.\ of finite type$)$ on $X$.
Consider the following sets$:$
\begin{itemize}
\item[{\rm (a)}] the set of all adically quasi-coherent subsheaves $\mathscr{G}\subseteq\mathscr{F}$ $($resp.\ of finite type$)$ that contains $\mathscr{I}\mathscr{F};$
\item[{\rm (b)}] the set of all quasi-coherent subsheaves of $\mathscr{F}_0=\mathscr{F}/\mathscr{I}\mathscr{F}$ $($resp.\ of finite type$)$ on the scheme $X_0=(X,\O_X/\mathscr{I})$.
\end{itemize}
Then the map $\mathscr{G}\mapsto\mathscr{G}/\mathscr{I}\mathscr{F}$ from the set {\rm (a)} to the set {\rm (b)} constitutes a bijection between these sets. \hfill$\square$
\end{cor}

\subsubsection{Admissible ideals}\label{subsub-admissibleideals}
\index{admissible!admissible ideal@--- ideal|(}
\begin{dfn}\label{dfn-admissibleideal}{\rm 
Let $X$ be an adic formal scheme of finite ideal type, and $\mathscr{J}$ an ideal sheaf of $\O_X$.
Then $\mathscr{J}$ is said to be an {\em admissible ideal} if it satisfies the following conditions:
\begin{itemize}
\item[(a)] (finiteness) $\mathscr{J}$ is an adically quasi-coherent ideal of finite type;
\item[(b)] (openness) $\mathscr{J}$ contains locally an ideal of definition.
\end{itemize}
We denote by $\AId_X$ the set of all admissible ideals of $X$.}
\end{dfn}

By \ref{lem-deltasheafadicallyquasicoherent1-1}, \ref{lem-deltasheafadicallyquasicoherent1-2}, and \ref{thm-adicqcoh1} we have:
\begin{prop}\label{prop-admissibleideallocalpresentation}
Let $A$ be an adic ring with a finitely generated ideal of definition $I\subseteq A$, and set $X=\Spf A$.
Then for any $I$-admissible\index{admissible!I-admissible@$I$-{---}} ideal $J\subseteq A$ {\rm ({\bf \ref{ch-pre}}.\ref{dfn-adm})}, $J^{\Delta}$ is an admissible ideal of $\O_X$.
Moreover, any admissible ideal of $\O_X$ is of this form by a uniquely determined $I$-admissible ideal $J$. \hfill$\square$
\end{prop}

\begin{prop}\label{prop-admissibleideal1x}
Let $f\colon Y\rightarrow X$ be an adic morphism of adic formal schemes of finite ideal type, and $\mathscr{J}$ an admissible ideal of $\O_X$.
Then the ideal pull-back\index{ideal pull-back} $\mathscr{J}\O_Y=(f^{-1}\mathscr{J})\O_Y$ {\rm ({\bf \ref{ch-pre}}.\ref{dfn-idealpullback})} is an admissible ideal of $\O_Y$.
\end{prop}

By this we have the mapping $\AId_X\rightarrow\AId_Y$ by $\mathscr{J}\mapsto\mathscr{J}\O_Y$.
\begin{proof}
We may work in the affine situation $X=\Spf A$ and $Y=\Spf B$ where $A$ and $B$ are adic rings of finite ideal type.
Let $I\subseteq A$ be a finitely generated ideal of definition of $A$.
Set $\mathscr{J}=J^{\Delta}$, where $J$ is an $I$-admissible ideal of $A$ (\ref{prop-admissibleideallocalpresentation}).
Then $JB$ is an $IB$-admissible ideal of $B$.
We need to show that $\mathscr{J}\O_Y=(JB)^{\Delta}$.
This follows from \ref{lem-deltasheafadicallyquasicoherent1-1} in view of the fact that $\O_Y=B^{\Delta}$ and $\O_Y/\mathscr{J}\O_Y=\til{B/JB}$.
\end{proof}

\begin{prop}[{cf.\ {\rm \ref{prop-admissibleideal2}} below}]\label{prop-admissibleidealquotientclosedsubscheme}
Let $X$ be an adic formal scheme $X$ of finite ideal type.
For an admissible ideal $\mathscr{J}\subseteq\O_X$, let $Y$ be the support of the sheaf $\O_X/\mathscr{J}$.
Then $Y$ is a closed subset of the underlying topological space of $X$, and the locally ringed space $(Y,\O_X/\mathscr{J})$ is a closed subscheme of $X$.
\end{prop}

\begin{proof}
We may assume that $X$ has an ideal of definition of finite type $\mathscr{I}$ such that $\mathscr{I}\subseteq\mathscr{J}$.
Consider the scheme $X_0=(X,\O_X/\mathscr{I})$.
The ideal sheaf $\mathscr{J}/\mathscr{I}$ on $X_0$ defines the closed subscheme of $X_0$ as in the proposition.
\end{proof}

\begin{cor}\label{cor-extension11}
Let $X$ be an adic formal scheme, and $\mathscr{J}\subseteq\O_X$ an admissible ideal.
Let $Y$ be the closed subscheme corresponding to $\mathscr{J}$ $($as in {\rm \ref{prop-admissibleidealquotientclosedsubscheme})}, and $i\colon Y\hookrightarrow X$ the canonical morphism.

{\rm (1)} For any quasi-coherent ideal $\mathscr{K}$ of finite type of $\O_Y$, the inverse image of $i_{\ast}\mathscr{K}$ by the map $\O_X\rightarrow i_{\ast}\O_Y$ is an admissible ideal of $\O_X$.

{\rm (2)} The map $\mathscr{K}\mapsto i^{-1}\mathscr{K}\O_Y$ gives a bijection from the set of all admissible ideals of $\O_X$ containing $\mathscr{J}$ to the set of all quasi-coherent ideals of finite type of $\O_Y$.
The inverse mapping is given by the inverse image by $\O_X\rightarrow i_{\ast}\O_Y$.
\end{cor}

\begin{proof}
We may assume without loss of generality that $X$ has an ideal of definition $\mathscr{I}$ such that $\mathscr{I}\subseteq\mathscr{J}$.
We regard quasi-coherent sheaves on $Y$ as quasi-coherent sheaves on the scheme $X_0=(X,\O_X/\mathscr{I})$ and apply \ref{cor-openaqcsubmoduleoffinitetype}.
\end{proof}

\begin{prop}\label{prop-exerpropadmissibleideal1}
Let $X$ be an adic formal scheme.
If $\mathscr{J}$ and $\mathscr{J}'$ are admissible ideals on $X$, then $\mathscr{J}\mathscr{J}'$ and $\mathscr{J}+\mathscr{J}'$ are admissible ideals.
\end{prop}

\begin{proof}
We may assume that $X$ admits an ideal of definition of finite type $\mathscr{I}$ contained in $\mathscr{J}$ and $\mathscr{J}'$.
Consider the closed subscheme $X_1=(X,\O_X/\mathscr{I}^2)$.
One sees easily that $\mathscr{J}\mathscr{J}'$ coincides with the pull-back of $\mathscr{J}\O_{X_1}\cdot\mathscr{J}'\O_{X_1}$.
Similarly, $\mathscr{J}+\mathscr{J}'$ is the pull-back of $\mathscr{J}\O_{X_1}+\mathscr{J}'\O_{X_1}$.
\end{proof}

\begin{prop}\label{prop-exeridealofdefinitioncompletepullback}
Let $f\colon X\rightarrow Y$ be an adically flat\index{morphism of formal schemes@morphism (of formal schemes)!adically flat morphism of formal schemes@adically flat ---} {\rm (\ref{exas-adicalizationexamples} (1))} morphism of adic formal schemes of finite ideal type, and $\mathscr{J}\subseteq\O_Y$ an admissible ideal of $Y$.
Then we have
$$
\widehat{f^{\ast}}\mathscr{J}=\mathscr{J}\O_X.
$$
In particular, $\widehat{f^{\ast}}\mathscr{J}$ is an admissible ideal of $X$ $($cf.\ {\rm \ref{prop-completepullbackaqcsheaves1}}$)$.
\end{prop}

\begin{proof}
We may assume that $Y$ has an ideal of definition $\mathscr{I}$ of finite type.
Let $f_k\colon X_k=(X,\O_X/\mathscr{I}^{k+1}\O_X)\rightarrow Y_k=(Y,\O_Y/\mathscr{I}^{k+1})$ be the induced morphism of schemes for $k\geq 0$.
By \ref{cor-completepullbackaqcsheaves0} we have
$$
\widehat{f^{\ast}}\mathscr{J}=\varprojlim_{k\geq 1}f^{\ast}_k\mathscr{J}/\mathscr{I}^{k+1}\mathscr{J}=\varprojlim_{k\geq 1}\mathscr{J}\O_{X_k}=\mathscr{J}\O_X,
$$
as desired.
\end{proof}
\index{admissible!admissible ideal@--- ideal|)}

\subsubsection{Extension of admissible ideals}\label{subsub-admissibleidealsextension}
\index{admissible!admissible ideal@--- ideal!extension of admissible ideals@extension of --- ---s|(}
Let us introduce a relation $\sim$ on the set $\AId_X$ as follows: $\mathscr{J}\sim\mathscr{J}'$ for $\mathscr{J},\mathscr{J}'\in\AId_X$ if and only if there exist $m,n>0$ such that $\mathscr{J}^m\subseteq\mathscr{J}^{\prime n}\subseteq\mathscr{J}$.
This gives an equivalence relation compatible with the semi-group structure on $\AId_X$.
Notice that the set of all ideals of definition\index{ideal of definition} of finite type on $X$ (if they exist) forms one single equivalence class. 

\begin{prop}\label{prop-extension1}
Let $X$ be a {\em coherent} {\rm (\ref{dfn-cohformalmorphism})}\index{formal scheme!coherent formal scheme@coherent ---} adic formal scheme, and $X=\bigcup_{\alpha\in L}U_{\alpha}$ a finite covering by quasi-compact open subsets.
Suppose that for each $\alpha\in L$ an admissible ideal $\mathscr{I}_{\alpha}$ on $U_{\alpha}$ is given such that on $U_{\alpha\beta}=U_{\alpha}\cap U_{\beta}$ we have $\mathscr{I}_{\alpha}|_{U_{\alpha\beta}}\sim\mathscr{I}_{\beta}|_{U_{\alpha\beta}}$.
Then there exists an admissible ideal $\mathscr{I}$ on $X$ such that we have $\mathscr{I}|_{U_{\alpha}}\sim\mathscr{I}_{\alpha}$ on each $U_{\alpha}$.
\end{prop}

\begin{proof}
By an easy inductive argument, it suffices to show the following: let $X=U_1\cup U_2$ with $U_1$ and $U_2$ quasi-compact, and $\mathscr{I}_1$ and $\mathscr{I}_2$ admissible ideals on $U_1$ and $U_2$, respectively, such that on $U_{12}=U_1\cap U_2$ (which is quasi-compact, since $X$ is quasi-separated (\ref{dfn-qsepformal})) there exists a positive integer $m$ such that $\mathscr{I}^m_1\subseteq\mathscr{I}_2\subseteq\mathscr{I}_1$; then there exists an admissible ideal $\mathscr{I}$ on $X$ such that $\mathscr{I}^m_1\subseteq\mathscr{I}\subseteq\mathscr{I}_1$ on $U_1$ and $\mathscr{I}=\mathscr{I}_2$ on $U_2$.

Let $V$ be the closed subscheme of $U_1$ defined by $\mathscr{I}^m_1$ (\ref{prop-admissibleidealquotientclosedsubscheme}).
Similarly, let $W$ be the closed subscheme of $U_{12}$ defined by $\mathscr{I}^m_1$ (restricted on $U_{12}$).
Then since $U_{12}$ is quasi-compact, $W$ is a quasi-compact open subscheme of $V$.

Consider the quasi-coherent ideal $\mathscr{I}_2/\mathscr{I}^m_1$ of finite type on $W$.
By \cite[$\mathbf{I}$, (9.4.2) (ii)]{EGA} one can extend this ideal to a quasi-coherent ideal $\ovl{\mathscr{I}}_{21}$ on $V$.
Replacing $\ovl{\mathscr{I}}_{21}$ by $\ovl{\mathscr{I}}_{21}\cap(\mathscr{I}_1/\mathscr{I}^m_1)$ if necessary, we may assume that the extension $\ovl{\mathscr{I}}_{21}$ is contained in $\mathscr{I}_1/\mathscr{I}^m_1$.
Moreover, by \cite[$\mathbf{I}$, (9.4.9) \& $\mathbf{IV}$, (1.7.7)]{EGA} there exists a subideal of $\ovl{\mathscr{I}}_{21}$ of finite type that coincides with $\mathscr{I}_2/\mathscr{I}^m_1$ on $W$ (since $\mathscr{I}_2/\mathscr{I}^m_1$ is of finite type).
Hence we may further assume that $\ovl{\mathscr{I}}_{21}$ is quasi-coherent of finite type on $V$.

Let $\mathscr{I}_{21}$ be the ideal on $U_1$ given by taking the inverse image of $\ovl{\mathscr{I}}_{21}$ by the map $\O_{U_1}\rightarrow\O_V$.
By \ref{cor-extension11} (1) $\mathscr{I}_{21}$ is an admissible ideal of $\O_{U_1}$.
Since $\mathscr{I}_{21}/\mathscr{I}^m_1$ coincides with $\mathscr{I}_2/\mathscr{I}^m_1$ on $W$, we have $\mathscr{I}_{21}=\mathscr{I}_2$ on $U_{12}$ by \ref{cor-extension11} (2).
Moreover, we have $\mathscr{I}_{21}\subseteq\mathscr{I}_1$.
Hence we get, by gluing, the admissible ideal $\mathscr{I}$ with the desired properties.
\end{proof}

\begin{cor}\label{cor-extension2}
Let $X$ be a coherent adic formal scheme of finite ideal type.
Then $X$ has an ideal of definition\index{ideal of definition} $\mathscr{I}$ of finite type.
\end{cor}

\begin{proof}
Take a finite open covering $\{U_{\alpha}\}_{\alpha\in L}$ such that on each $U_{\alpha}$ there exists an ideal of definition\index{ideal of definition} $\mathscr{I}_{\alpha}$ of finite type.
Then apply \ref{prop-extension1}.
\end{proof}

\begin{cor}\label{cor-affineadicformalschemebyadicring}
Let $A$ be an admissible ring\index{admissible!admissible topological ring@--- (topological) ring} {\rm (\ref{dfn-admissibleringsadicrings})} such that $X=\Spf A$ is an affine adic formal scheme of finite ideal type.
Then the ring $A$ is an adic ring of finite ideal type.
\end{cor}

\begin{proof}
There exists an ideal of definition of finite type $\mathscr{I}$ on $X$ (\ref{cor-extension2}).
By \ref{prop-grothendieckEGAnew10.3.5} we have the unique finitely generated ideal of definition $I^{(n)}\subseteq A$ such that $\mathscr{I}^n=(I^{(n)})^{\Delta}$ for each $n>0$.
Since $\{\mathscr{I}^n\}_{n>0}$ gives a fundamental system of ideals of definition on $X$ (\ref{cor-adicformalschemeidealofdefinitionfinitetype22}), $\{I^{(n)}\}_{n>0}$ is a fundamental system of ideals of definition of the admissible ring $A$.
Set $I=I^{(1)}$.
Since each $X_k=(X,\O_X/\mathscr{I}^{k+1})$ ($k\geq 0$) is affine, we have for $m\leq n$
$$
I^{(m)}/I^{(n)}=\Gamma(X,\mathscr{I}^m\O_X/\mathscr{I}^n)=I^m(A/I^{(n)});
$$
moreover, $I^{(1)}/I^{(2)}$ is obviously a finitely generated ideal of $A/I^{(2)}$.
Hence by {\bf \ref{ch-pre}}.\ref{prop-criterionadicness1} we have $I^{(n)}=I^n$ for $n>0$, and hence $A$ is adic of finite ideal type.
\end{proof}

By this and \ref{prop-tuaadeq2} we deduce:
\begin{cor}\label{cor-affineadicformalschemebyadicringtua}
Let $A$ be an admissible ring\index{admissible!admissible topological ring@--- (topological) ring} such that $X=\Spf A$ is universally adhesive $($resp.\ universally rigid-Noetherian$)$ {\rm (\ref{dfn-formalsch})}\index{formal scheme!universally adhesive formal scheme@universally adhesive ---}\index{adhesive!universally adhesive@universally ---!universally adhesive formal scheme@--- --- formal scheme}\index{formal scheme!universally rigid-Noetherian formal scheme@universally rigid-Noetherian ---}.
Then the ring $A$ is t.u.\ adhesive $($resp.\ t.u.\ rigid-Noetherian$)$\index{t.u.a. ring@t.u.\ adhesive ring}\index{t.u. rigid-Noetherian ring@t.u.\ rigid-Noetherian ring} {\rm (\ref{dfn-tuaringadmissible})}. \hfill$\square$
\end{cor}

\begin{prop}[Extension lemma]\label{prop-extension3}
Let $X$ be a coherent adic formal scheme, and $U$ a quasi-compact open subset of $X$.
Then for any admissible ideal $\mathscr{J}$ on $U$ there exists an admissible ideal $\til{\mathscr{J}}$ on $X$ such that $\til{\mathscr{J}}|_{U}=\mathscr{J}$.
$($That is, the restriction map $\AId_X\rightarrow\AId_U$ is surjective.$)$
\end{prop}

\begin{proof}
Let $\mathscr{I}$ be an ideal of definition of $X$ (\ref{cor-extension2}).
We may assume that $\mathscr{I}\subseteq\mathscr{J}$ on $U$, since $U$ is quasi-compact.
Let $V$ (resp.\ $W$) be the closed subscheme of $X$ (resp.\ $U$) defined by $\mathscr{I}$.
Then $W$ is a quasi-compact open subset of $V$.
We apply \cite[$\mathbf{I}$, (9.4.2) (ii)]{EGA} and \cite[$\mathbf{I}$, (9.4.9) \& $\mathbf{IV}$, (1.7.7)]{EGA} to get a quasi-coherent extension $\ovl{\mathscr{J}}$ of finite type on $X$ and an admissible ideal $\til{\mathscr{J}}$ on $X$ (\ref{cor-extension11} (1)).
By \ref{cor-extension11} (2) one sees that this $\til{\mathscr{J}}$ is a desired extension.
\end{proof}
\index{admissible!admissible ideal@--- ideal!extension of admissible ideals@extension of --- ---s|)}

\addcontentsline{toc}{subsection}{Exercises}
\subsection*{Exercises}
\begin{exer}\label{exer-adicallyqcohalgebrafinitetypecriterion}{\rm 
Let $X$ be an adic formal scheme, $\mathscr{I}$ an ideal of definition of finite type of $X$, and $\mathscr{B}$ an adically quasi-coherent $\O_X$-algebra.
Then show that $\mathscr{B}$ is an adically quasi-coherent $\O_X$-module sheaf of finite type if and only if $\mathscr{B}/\mathscr{I}\mathscr{B}$ is a quasi-coherent sheaf of finite type over the scheme $X_0=(X,\O_X/\mathscr{I})$.}
\end{exer}

\begin{exer}\label{exer-adicquotientbynonadic1}{\rm 
Let $X$ be a locally universally rigid-Noetherian formal scheme, and $\varphi\colon\mathscr{F}\rightarrow\mathscr{G}$ a morphism of adically quasi-coherent sheaves of finite type on $X$.
Then show that $\ker(\varphi)$ and $\coker(\varphi)$ are adically quasi-coherent sheaves on $X$.
Show, moreover, that if $X$ is coherent and $\mathscr{I}$ is an ideal of definition of finite type on $X$, then the $\mathscr{I}$-torsion parts of $\ker(\varphi)$ and $\coker(\varphi)$ are bounded $\mathscr{I}$-torsion.}
\end{exer}

\begin{exer}[Extension of adically quasi-coherent sheaves]\label{exer-extopenadiccoh1}{\rm 
Let $X$ be a coherent adic formal scheme with an ideal of definition $\mathscr{I}$ of finite type, and $U$ a quasi-compact open subset of $X$.
Let $\mathscr{F}$ be an adically quasi-coherent sheaf of finite type, and $\mathscr{G}$ an adically quasi-coherent subsheaf of $\mathscr{F}|_U$ of finite type such that $\mathscr{I}\mathscr{F}|_U\subseteq\mathscr{G}$.
Then there exists an adically quasi-coherent subsheaf $\mathscr{G}'$ of $\mathscr{F}$ of finite type such that $\mathscr{G}'|_U=\mathscr{G}$.}
\end{exer}

\begin{exer}\label{exer-extopenadiccoh1approx}{\rm 
Let $X$ be a coherent universally rigid-Noetherian formal scheme with an ideal of definition $\mathscr{I}$ of finite type, $\mathscr{F}$ an adically quasi-coherent sheaf of finite type on $X$, and $\mathscr{G}\subseteq\mathscr{F}$ an adically quasi-coherent subsheaf.
Then show that $\mathscr{G}$ is an inductive limit $\varinjlim_{\lambda\in\Lambda}\mathscr{H}_{\lambda}$ of adically quasi-coherent subsheaves of finite type such that $\mathscr{G}/\mathscr{H}_{\lambda}$ for all $\lambda$ are annihilated by $\mathscr{I}^n$ for some $n>0$.}
\end{exer}

\begin{exer}[cf.\ \ref{prop-lemdeltasheafadicallyquasicoherent1-3}]\label{exer-lemdeltasheafadicallyquasicoherent1-4general}{\rm 
Let $A$ be a rigid-Noetherian ring with a finitely generated ideal of definition $I$, and $N\subseteq M$ an inclusion of $A$-modules, and suppose that $M$ is contained as an $A$-submodule in a finitely generated $A$-module.
Then show the sequence 
$$
0\longrightarrow N^{\Delta}\longrightarrow M^{\Delta}\longrightarrow(M/N)^{\Delta}\longrightarrow 0
$$
by canonical morphisms is exact.}
\end{exer}

\begin{exer}\label{exer-lemdeltasheafadicallyquasicoherent1-application}{\rm 
Let $X$ be a locally universally rigid-Noetherian formal scheme, and consider an exact sequence of $\O_X$-modules
$$
0\longrightarrow\mathscr{F}\longrightarrow\mathscr{G}\longrightarrow\mathscr{H}\longrightarrow 0.
$$

(1) If $\mathscr{G}$ and $\mathscr{H}$ are adically quasi-coherent of finite type, then $\mathscr{F}$ is adically quasi-coherent.

(2) If $\mathscr{F}$ is adically quasi-coherent and $\mathscr{G}$ is adically quasi-coherent of finite type, then $\mathscr{H}$ is adically quasi-coherent of finite type.}
\end{exer}

\begin{exer}\label{exer-adicallyqcohlocallnoetherian}{\rm 
Let $X$ be a locally Noetherian formal scheme\index{formal scheme!Noetherian formal scheme@Noetherian ---!locally Noetherian formal scheme@locally --- ---} (\ref{dfn-northeroanformalschemes}).

(1) Show that the sheaf $\O_X$ is coherent\index{coherent!coherent sheaf@--- sheaf (on a ringed space)}.

(2) Show that an $\O_X$-module sheaf $\mathscr{F}$ is coherent if and only if it is adically quasi-coherent of finite type. 

(3) Show that any adically quasi-coherent ideal sheaf $\mathscr{J}\subseteq\O_X$ is coherent.}
\end{exer}

\begin{exer}\label{exer-propadmissibleideal2}
{\rm Let $X$ be a locally universally adhesive formal scheme, and $\mathscr{J},\mathscr{J}'\subseteq\O_X$ admissible ideals. Suppose that $X$ has an ideal of definition $\mathscr{I}$ and that $\O_X$ is $\mathscr{I}$-torsion free. Then show that $\mathscr{J}\cap\mathscr{J}'$ is admissible.}
\end{exer}


\section{Several properties of morphisms}\label{sec-severalcond}
\subsection{Affine morphisms}\label{sub-affinemorphism}
\index{morphism of formal schemes@morphism (of formal schemes)!affine morphism of formal schemes@affine ---|(}
\subsubsection{Definition of affine morphisms}\label{subsub-affinemorphismdef}
First of all, let us recall the definition of affine morphisms of formal schemes (\cite[(10.16.1)]{EGAInew}):
\begin{dfn}\label{dfn-affmorform}{\rm 
A morphism $f\colon X\rightarrow Y$ of formal schemes is said to be {\em affine}\index{affine!affine morphism@--- morphism} if there exists an affine open covering $Y=\bigcup_{\alpha\in L}V_{\alpha}$ of $Y$ such that for each $\alpha\in L$ the open formal subscheme $f^{-1}(V_{\alpha})$ in $X$ is an affine formal scheme.}
\end{dfn}

For example, any morphism between affine formal schemes is affine.
\begin{prop}\label{prop-affinemorestablezar}
Let $f\colon X\rightarrow Y$ be a morphism of formal schemes, and $Y=\bigcup_{\alpha\in L}V_{\alpha}$ an open covering of $Y$.
Then $f$ is affine if and only if the induced map $f_{\alpha}\colon X_{\alpha}=X\times_YV_{\alpha}\rightarrow V_{\alpha}$ is affine for any $\alpha\in L$.
\end{prop}

\begin{proof}
The `if' part is trivial and the `only if' part follows from the following observation: Let $U=\Spf B$ be an affine open subset of $Y$ such that $f^{-1}(U)$ is affine $f^{-1}(U)=\Spf A$. Then $V_{\alpha}\cap U$ for each $\alpha\in L$ is covered by affine open subsets of the form $W=\Spf B_{\{g\}}$, giving $f^{-1}(W)=\Spf A\widehat{\otimes}_BB_{\{g\}}$. 
\end{proof}

\subsubsection{Affine adic morphisms and adically quasi-coherent sheaves}\label{subsub-affinemorphismaqch}
\begin{prop}\label{prop-affinemorphism2}
Let $f\colon X\rightarrow Y$ be an affine adic {\rm (\ref{dfn-adicmor})} morphism\index{morphism of formal schemes@morphism (of formal schemes)!adic morphism of formal schemes@adic ---} between adic formal schemes of finite ideal type.
Then for any adically quasi-coherent sheaf $\mathscr{F}$\index{adically quasi-coherent (a.q.c.) sheaf} on $X$, $f_{\ast}\mathscr{F}$ is an adically quasi-coherent sheaf on $Y$.
\end{prop}

Since the question is local on $Y$, we may assume that $X$ and $Y$ are affine $X=\Spf A$ and $Y=\Spf B$, where $A$ and $B$ are adic rings of finite ideal type, and that the map $B\rightarrow A$ is adic\index{adic!adic morphism@--- morphism}\index{morphism of admissible rings@morphism (of admissible/adic rings)!adic morphism of admissible rings@adic ---}; we have $\mathscr{F}=M^{\Delta}$ by an $A$-module $M$ (\ref{prop-topqcoh1}). 
Then the proposition follows from the following lemma and \ref{prop-adicqcohaff1f}:
\begin{lem}\label{lem-exeradicallydeltasheavesunited}
In the situation as above, we have
$$
f_{\ast}M^{\Delta}=M_{[B]}^{\Delta},
$$
where $M_{[B]}$ denotes the module $M$ regarded as a $B$-module.
\end{lem}

\begin{proof}
By \ref{prop-adicqcohaff1f} and {\bf \ref{ch-pre}}.\ref{prop-projlimsheafleftexact0} (2) we calculate
$$
f_{\ast}M^{\Delta}=f_{\ast}\varprojlim_{k\geq 0}\til{M/I^{k+1}M}=\varprojlim_{k\geq 0}f_{\ast}\til{M/I^{k+1}M}=\varprojlim_{k\geq 0}\til{M_{[B]}/I^{k+1}M_{[B]}},
$$
where the last projective limit sheaf is equal to $M_{[B]}^{\Delta}$, as desired.
\end{proof}

\subsubsection{Formal spectra of a.q.c.\ algebras}\label{subsub-effectivedescaffineadic}
As a corollary of \ref{prop-affinemorphism2} we have:
\begin{cor}\label{cor-effectivedescaffineadic}
Let $f\colon X\rightarrow Y$ be an affine adic morphism of adic formal schemes of finite ideal type.
Then $f_{\ast}\O_X$ is an adically quasi-coherent $\O_Y$-algebra\index{adically quasi-coherent (a.q.c.) algebra}.
Moreover, we have$:$
\begin{itemize}
\item[{\rm (1)}] $f_{\ast}\O_X$ is of finite type\index{adically quasi-coherent (a.q.c.) algebra!adically quasi-coherent (a.q.c.) algebra of finite type@--- of finite type} {\rm (\ref{dfn-adicallyqcohalgebraoffinitetype})} if and only if $f\colon X\rightarrow Y$ is locally of finite type\index{morphism of formal schemes@morphism (of formal schemes)!morphism of formal schemes locally of finite type@--- locally of finite type} {\rm (\ref{dfn-topfintype})}$;$
\item[{\rm (2)}] suppose $Y$ is locally universally rigid-Noetherian {\rm (\ref{dfn-formalsch})}\index{formal scheme!universally rigid-Noetherian formal scheme@universally rigid-Noetherian ---!locally universally rigid-Noetherian formal scheme@locally --- ---}$;$ then $f_{\ast}\O_X$ is of finite presentation\index{adically quasi-coherent (a.q.c.) algebra!adically quasi-coherent (a.q.c.) algebra of finite presentation@--- of finite presentation} {\rm (\ref{dfn-adicallyqcohalgebraoffinitepres})} if and only if $f\colon X\rightarrow Y$ is locally of finite presentation\index{morphism of formal schemes@morphism (of formal schemes)!morphism of formal schemes locally of finite presentation@--- locally of finite presentation} {\rm (\ref{dfn-topfinpres})}.
\end{itemize}
\end{cor}

\begin{proof}
We may suppose that $X$ and $Y$ are affine; if $f\colon X=\Spf A\rightarrow Y=\Spf B$ comes from an adic morphism\index{adic!adic morphism@--- morphism}\index{morphism of admissible rings@morphism (of admissible/adic rings)!adic morphism of admissible rings@adic ---} $B\rightarrow A$ of adic rings of finite ideal type, then we have $f_{\ast}\O_X=A_{[B]}^{\Delta}$.
Then apply \ref{prop-adicallyqcohalgebraoffinitetype2} and \ref{prop-adicallyqcohalgebraoffinitepres2}.
\end{proof}

\begin{prop}\label{prop-effectivedescaffineadic1}
Let $f\colon X\rightarrow Y$ and $f'\colon X'\rightarrow Y$ be affine adic morphisms of adic formal schemes of finite ideal type.
Then the induced map
$$
\Hom_Y(X,X')\longrightarrow\Hom_{\Alg_Y}(f'_{\ast}\O_{X'},f_{\ast}\O_X)
$$
is bijective.
\end{prop}

\begin{proof}
We may work in the affine situation $Y=\Spf B$, $X=\Spf A$, and $X'=\Spf A'$, and $f$ and $f'$ come respectively from adic morphisms $B\rightarrow A$ and $B\rightarrow A'$ of adic rings of finite ideal type.
Then by \cite[$\mathbf{I}$, (10.1.3)]{EGA} the set $\Hom_Y(X,X')$ is identified with the set of continuous homomorphisms $A'\rightarrow A$ of topological rings over $B$.
But since $A$ and $A'$ are the adic over $B$, the last set is simply the set $\Hom_B(A',A)$ of $B$-algebra homomorphisms.
By \ref{lem-corcoradicqcoh1} this is further isomorphic to the set $\Hom_{\Alg_Y}(A_{[B]}^{\prime\Delta},A_{[B]}^{\Delta})$.
Now, since $A^{\Delta}=f_{\ast}\O_X$ etc., we obtain the desired result.
\end{proof}

\begin{prop}\label{prop-affinemorlocalconst1}
Let $Y$ be an adic formal scheme of finite ideal type, and $\mathscr{A}$ an adically quasi-coherent $\O_Y$-algebra\index{adically quasi-coherent (a.q.c.) algebra}.
Then there exists uniquely up to $Y$-isomorphisms an affine morphism $f\colon X\rightarrow Y$ such that $f_{\ast}\O_X\cong\mathscr{A}$.
\end{prop}

\begin{proof}
By \ref{thm-adicqcoh1}, for any affine open subsets $V=\Spf B_V$ of $Y$ ($B_V$ is an adic ring of finite ideal type due to \ref{cor-affineadicformalschemebyadicring}), $\mathscr{A}(V)=A_V$ is an $I_V$-adically complete $B_V$-algebra, where $I_V$ is an ideal of definition of $B_V$, and $A^{\Delta}_V=\mathscr{A}|_V$ holds.
Take an affine open covering $Y=\bigcup_{\alpha\in L}V_{\alpha}$ of $Y$, and set $X_{\alpha}=\Spf A_{V_{\alpha}}$ for each $\alpha\in L$.
By \ref{prop-effectivedescaffineadic1} these formal schemes glue to an adic formal scheme $X=\bigcup_{\alpha\in L}X_{\alpha}$ adic over $Y$.
The map $f\colon X\rightarrow Y$ thus obtained is affine and adic, and $f_{\ast}\O_X\cong\mathscr{A}$.
Notice that the above construction does not depend on the choice of the affine covering $\{V_{\alpha}\}_{\alpha\in L}$ (due to \ref{prop-effectivedescaffineadic1}).
\end{proof}

Let $Y$ be an adic formal scheme of finite ideal type, and consider the functor 
$$
\Af\Ac\Fs^{\ast,\opp}_{/Y}\longrightarrow\AQCoh\Alg_Y,\leqno{(\ast)}
$$
where $\AQCoh\Alg_Y$ is the category of adically quasi-coherent $\O_Y$-algebras (with $\O_Y$-algebra homomorphisms), which maps each affine adic morphism $f\colon X\rightarrow Y$ to $f_{\ast}\O_X$.
By \ref{prop-effectivedescaffineadic1} and \ref{prop-affinemorlocalconst1} we have:
\begin{thm}\label{thm-effectivedescaffineadic}
The functor $(\ast)$ is a categorical equivalence. \hfill$\square$
\end{thm}

\begin{dfn}\label{dfn-affinespectrumformalrel}{\rm 
Let $X$ be an adic formal scheme of finite ideal type, and $\mathscr{A}$ an adically quasi-coherent $\O_X$-algebra.
Then by \ref{prop-affinemorlocalconst1} there exists uniquely up to isomorphism over $X$ an adic formal scheme $Z$ and an affine adic map $f\colon Z\rightarrow X$ such that $f_{\ast}\O_Z\cong\mathscr{A}$.
We denote this formal scheme $Z$ by $\Spf\mathscr{A}$ and call it the {\em formal spectrum}\index{formal spectrum} of $\mathscr{A}$.}
\end{dfn}

Similarly to the scheme case, the formation of the formal spectra is an example of `effective local construction'\index{local!local construction@--- construction} (with respect to Zariski topology) in {\bf \ref{ch-pre}}.\ref{dfn-localconstruction}.

\subsubsection{Basic properties of affine adic morphisms}\label{subsub-fundamentalpropertiesaffineadicbasechange}
\begin{prop}\label{prop-affmorform1}
Let $f\colon X\rightarrow Y$ be an adic morphism of adic formal schemes of finite ideal type.
Then $f$ is an affine map if and only if for any affine open subset $V=\Spf B$ of $Y$ the open formal subscheme $f^{-1}(V)$ of $X$ is affine.
\end{prop}

\begin{proof}
The `if' part is trivial.
Suppose $f\colon X\rightarrow Y$ is affine adic, and let $\mathscr{A}=f_{\ast}\O_X$.
By \ref{thm-effectivedescaffineadic} we have $X\cong\Spf\mathscr{A}$ over $Y$.
Let $V=\Spf B$ be an affine open subset of $Y$.
By \ref{cor-affineadicformalschemebyadicring} $B$ is an adic ring of finite ideal type, and so is $\mathscr{A}(V)$ due to \ref{thm-adicqcoh1}.
Then by the construction of $\Spf\mathscr{A}$ we find that $X\times_YV\cong\Spf\mathscr{A}|_V=\Spf\mathscr{A}(V)$.
\end{proof}

\begin{prop}\label{prop-affmorform2}
{\rm (1)} The composition of two affine adic morphisms of adic formal schemes of finite ideal type is again affine adic.

{\rm (2)} For any affine adic $S$-morphisms $f\colon X\rightarrow Y$ and $g\colon X'\rightarrow Y'$ of adic formal schemes of finite ideal type over an adic formal scheme $S$ of finite ideal type, the induced morphism $f\times_Sg\colon X\times_SY\rightarrow X'\times_SY'$ is affine adic.

{\rm (3)} For any affine adic $S$-morphism $f\colon X\rightarrow Y$ of adic formal schemes of finite ideal type over an adic formal scheme $S$ of finite ideal type and for any morphism $S'\rightarrow S$ of adic formal schemes of finite ideal type, the induced morphism $f_{S'}\colon X\times_SS'\rightarrow Y\times_SS'$ is affine adic.
\end{prop}

\begin{proof}
(1) follows immediately from \ref{prop-affmorform1}.
By \ref{prop-genpropertymorphismadic1} the assertions (2) and (3) follow from \ref{prop-affinemorphism1} below combined with \ref{prop-affinemorestablezar} and \cite[$\mathbf{II}$, (1.6.2)]{EGA}.
\end{proof}

\begin{prop}\label{prop-affinemorphism1}
Let $f\colon X\rightarrow Y$ be an adic morphism of adic formal schemes of finite ideal type, and suppose $Y$ has an ideal of definition $\mathscr{I}$ of finite type.
For any integer $k\geq 0$ we denote by $f_k\colon X_k=(X,\O_X/\mathscr{I}^{k+1}\O_X)\rightarrow Y_k=(Y,\O_Y/\mathscr{I}^{k+1})$ the induced morphism of schemes.
Then the following conditions are equivalent$:$
\begin{itemize}
\item[{\rm (a)}] $f$ is affine$;$
\item[{\rm (b)}] $f_k$ is affine\index{affine!affine morphism@--- morphism} for any $k\geq 0;$
\item[{\rm (c)}] $f_0$ is affine.
\end{itemize}
\end{prop}

\begin{proof}
By \ref{prop-affinemorestablezar} we may assume that $Y$ is affine $Y=\Spf B$ with the ideal $I$ of definition of $B$.
Suppose $f$ is affine.
Then by \ref{prop-affmorform1} $X$ is an affine formal scheme $X=\Spf A$ and, in this situation, we have $X_k=\Spec A_k$, where $A_k=A/I^{k+1}A$. 
Hence the implication (a) $\Rightarrow$ (b) holds.
(b) $\Rightarrow$ (c) is trivial.
Suppose (c) holds.
By \cite[(2.3.5)]{EGAInew} we deduce that $f_k$ is affine for any $k\geq 0$.
Set $X_k=\Spec A_k$ for each $k$.
Then $\{A_k\}_{k\geq 0}$ forms a projective system of rings such that for $k\leq l$ the transition map $A_l\rightarrow A_k$ is surjective with the kernel equal to $I^{k+1}A_l$.
Then by {\bf \ref{ch-pre}}.\ref{prop-exeramazingfactonalgebra} we see that $A$ is $IA$-adically complete and that $A/I^{k+1}=A_k$ for any $k\geq 0$.
Therefore, we have $X=\Spf A$.
\end{proof}

In view of \ref{prop-genpropertymorphismadic2} we have the following result due to \ref{prop-affinemorphism1} and \cite[$\mathbf{II}$, (1.6.2)]{EGA}:
\begin{cor}\label{cor-affinemorphism11}
Let $f\colon X\rightarrow Y$ be a morphism of schemes, and $Z$ a closed subscheme of $Y$ of finite presentation.
If $f$ is affine, then the formal completion $\widehat{f}\colon\widehat{X}|_{f^{-1}(Z)}\rightarrow\widehat{Y}|_Z$ is affine.\hfill$\square$
\end{cor}
\index{morphism of formal schemes@morphism (of formal schemes)!affine morphism of formal schemes@affine ---|)}

\subsection{Finite morphisms}\label{sub-finitemorformal}
\index{morphism of formal schemes@morphism (of formal schemes)!finite morphism of formal schemes@finite ---|(}
\begin{prop}\label{prop-finitemorform1}
The following conditions for an adic morphism\index{morphism of formal schemes@morphism (of formal schemes)!adic morphism of formal schemes@adic ---} $f\colon X\rightarrow Y$ of adic formal schemes of finite ideal type are equivalent to each other$:$
\begin{itemize}
\item[{\rm (a)}] there exists an affine open covering $\{V_{\alpha}=\Spf A_{\alpha}\}_{\alpha\in L}$ of $Y$, where each $A_{\alpha}$ is an adic ring with a finitely generated ideal of definition $I_{\alpha}\subseteq A_{\alpha}$, such that for each $\alpha\in L$ the induced morphism $X\times_YV_{\alpha,0}\rightarrow V_{\alpha,0}$ $($where $V_{\alpha,0}=\Spec A_{\alpha}/I_{\alpha})$ of schemes is finite$;$
\item[{\rm (a)\rlap{${}'$}}] for any affine open set $V=\Spf A$ of $Y$, where $A$ is an adic ring with a finitely generated ideal of definition $I\subseteq A$, the induced morphism $X\times_YV_0\rightarrow V_0$ $($where $V_0=\Spec A/I)$ of schemes is finite$;$
\item[{\rm (b)}] there exists an affine open covering $\{V_{\alpha}=\Spf A_{\alpha}\}_{\alpha\in L}$ of $Y$, where each $A_{\alpha}$ is an adic ring with a finitely generated ideal of definition $I_{\alpha}\subseteq A_{\alpha}$, such that for each $\alpha\in L$ the inverse image $f^{-1}(V_{\alpha})$ is affine of the form $f^{-1}(V_{\alpha})=\Spf B_{\alpha}$, where $B_{\alpha}$ is finitely generated as an $A_{\alpha}$-module$;$
\item[{\rm (b)\rlap{${}'$}}] for any affine open set $V=\Spf A$ of $Y$, where $A$ is an adic ring with a finitely generated ideal of definition $I\subseteq A$, $f^{-1}(V)$ is affine of the form $f^{-1}(V)=\Spf B$, where $B$ is finitely generated as an $A$-module.
\end{itemize}
\end{prop}

\begin{proof}
The implications (b) $\Rightarrow$ (a) and (b)\rlap{${}'$} $\Rightarrow$ (a)\rlap{${}'$} are immediate.
Let us show (a) $\Rightarrow$ (b) and (a)\rlap{${}'$} $\Rightarrow$ (b)\rlap{${}'$}.
Suppose we are in the situation as in (a)\rlap{${}'$}.
We set $V_k=\Spf A/I^{k+1}$ and $U_k=X\times_YV_k$ for $k\geq 0$.
By \cite[$\mathbf{I}$, (5.1.9)]{EGA} each $X_k$ is an affine scheme $X_k=\Spec B_k$ and for $k\leq l$ we have $B_k=B_l/I^{k+1}B_l$.
Then by {\bf \ref{ch-pre}}.\ref{prop-exeramazingfactonalgebra} we see that $B=\varprojlim_{k\geq 0}B_k$ is an $IB$-adically complete $A$-algebra such that $B/I^{k+1}B=B_k$ for $k\geq 0$.
Since $B_0$ is finitely generated as an $A$-module, by {\bf \ref{ch-pre}}.\ref{prop-complpair1} we find that $B$ is finitely generated as an $A$-module.
Moreover, since $f\colon X\rightarrow Y$ is adic, we have $f^{-1}(V)=\varinjlim_{k\geq 0}U_k=\Spf B$ and thus arrive at the situation as in (b)\rlap{${}'$}.
The other implication (a) $\Rightarrow$ (b) can be verified similarly.

It remains to show the equivalence of (a) and (a)\rlap{${}'$}.
As the implication (a)\rlap{${}'$} $\Rightarrow$ (a) is clear, we want to show the converse.
First notice that in view of \cite[$\mathbf{I}$, (5.1.9)]{EGA} each ideal of definition $I_{\alpha}$ of $A_{\alpha}$ can be replaced by a power $I^n_{\alpha}$ for any $n\geq 1$.
Notice also that the covering $\{V_{\alpha}=\Spf A_{\alpha}\}_{\alpha\in L}$ can be replaced by a refinement by replacing each $V_{\alpha}=\Spf A_{\alpha}$ by a finite affine covering given by affine open subschemes of the form $\Spf(A_{\alpha})_{\{g\}}$.
Hence the desired implication follows from the fact that finite morphisms of schemes is stable under Zariski topology\index{stable!stable under a topology@--- under a topology} ({\bf \ref{ch-pre}}.\ref{dfn-catequivrelstablearrow} (1)).
\end{proof}

\begin{dfn}\label{dfn-finitemorform1}{\rm 
Let $X$ and $Y$ be adic formal schemes of finite ideal type.
Then a morphism $f\colon X\rightarrow Y$ is said to be {\em finite} if it is adic\index{morphism of formal schemes@morphism (of formal schemes)!adic morphism of formal schemes@adic ---} and satisfies the conditions in {\rm \ref{prop-finitemorform1}}.}
\end{dfn}

\begin{prop}\label{prop-finitemorphism1}
Let $f\colon X\rightarrow Y$ be an adic morphism of adic formal schemes of finite ideal type, and suppose $Y$ has a ideal of definition $\mathscr{I}$ of finite type.
For any integer $k\geq 0$ we denote by $f_k\colon X_k=(X,\O_X/\mathscr{I}^{k+1}\O_X)\rightarrow Y_k=(Y,\O_Y/\mathscr{I}^{k+1})$ the induced morphism of schemes.
Then the following conditions are equivalent$:$
\begin{itemize}
\item[{\rm (a)}] $f$ is finite$;$
\item[{\rm (b)}] $f_k$ is finite for any $k\geq 0;$
\item[{\rm (c)}] $f_0$ is finite.
\end{itemize}
\end{prop}

\begin{proof}
We may assume that $Y$ is affine of the form $Y=\Spf A$, where $A$ is an adic ring with the finitely generated ideal of definition $I$ such that $\mathscr{I}=I^{\Delta}$.
Then the equivalence of (b) and (c) follows from \cite[$\mathbf{I}$, (5.1.9)]{EGA} and {\bf \ref{ch-pre}}.\ref{prop-complpair1}.
The equivalence of (a) and (c) follows from the definition.
\end{proof}

\begin{prop}\label{prop-finitemorform2}
{\rm (1)} A finite morphism is affine\index{morphism of formal schemes@morphism (of formal schemes)!affine morphism of formal schemes@affine ---} {\rm (\ref{dfn-affmorform})}.

{\rm (2)} The composition of two finite morphisms is again finite.

{\rm (3)} For any finite $S$-morphisms $f\colon X\rightarrow Y$ and $g\colon X'\rightarrow Y'$ of adic formal schemes of finite ideal type over an adic formal scheme $S$ of finite ideal type, the induced morphism $f\times_Sg\colon X\times_SY\rightarrow X'\times_SY'$ is finite.

{\rm (4)} For any finite $S$-morphism $f\colon X\rightarrow Y$ of adic formal schemes of finite ideal type over an adic formal scheme $S$ of finite ideal type and for any morphism $S'\rightarrow S$ of adic formal schemes of finite ideal type, the induced morphism $f_{S'}\colon X\times_SS'\rightarrow Y\times_SS'$ is finite.
\end{prop}

\begin{proof}
(1) and (2) are clear.
In view of \ref{prop-genpropertymorphismadic1} the assertions (3) and (4) follow from \ref{prop-finitemorphism1}, \ref{prop-finitemorform1}, and \cite[$\mathbf{II}$, (6.1.5)]{EGA}.
\end{proof}

In view of \ref{prop-genpropertymorphismadic2} we deduce the following due to \ref{prop-finitemorphism1} and \cite[$\mathbf{II}$, (6.1.5)]{EGA}:
\begin{cor}\label{cor-finitemorphism11}
Let $f\colon X\rightarrow Y$ be a morphism of schemes, and $Z$ a closed subscheme of $Y$ of finite presentation.
Then if $f$ is finite, the formal completion $\widehat{f}\colon\widehat{X}|_{f^{-1}(Z)}\rightarrow\widehat{Y}|_Z$ is finite. \hfill$\square$
\end{cor}

Since any finite morphism $f\colon X\rightarrow Y$ is affine, it comes from a formal spectrum\index{formal spectrum} $X\cong\Spf\mathscr{A}\rightarrow Y$ (\ref{dfn-affinespectrumformalrel}).
By \ref{thm-adicqcoh1} (2) we have (cf.\ Exercise \ref{exer-adicallyqcohalgebrafinitetypecriterion}):
\begin{prop}\label{prop-finitemorphismlocalconstruction}
Let $Y$ be an adic formal scheme of finite ideal type, and $\mathscr{A}$ an adically quasi-coherent $\O_Y$-algebra\index{adically quasi-coherent (a.q.c.) algebra}.
Then the map $\Spf\mathscr{A}\rightarrow Y$ is finite if and only if $\mathscr{A}$ is an adically quasi-coherent $\O_Y$-module of finite type. \hfill$\square$
\end{prop}
\index{morphism of formal schemes@morphism (of formal schemes)!finite morphism of formal schemes@finite ---|)}

\subsection{Closed immersions}\label{sub-closedimmformal}
\index{immersion!closed immersion of formal schemes@closed --- (of formal schemes)|(}
\subsubsection{A preliminary result}\label{subsub-closedimmformalpre}
\begin{prop}\label{prop-closedimmformal1}
Let $X$ be an adic formal scheme of finite ideal type, $\mathscr{F}$ an adically quasi-coherent sheaf\index{quasi-coherent!adically quasi-coherent OX module@adically --- (a.q.c.) sheaf} on $X$, and $\mathscr{K}\subseteq\mathscr{F}$ an $\O_X$-submodule.
Suppose that the quotient $\mathscr{G}=\mathscr{F}/\mathscr{K}$ is an adically quasi-coherent sheaf.

{\rm (1)} Suppose $X$ has an ideal of definition $\mathscr{I}$ of finite type. 
Then the morphism 
$$
\mathscr{K}\longrightarrow\varprojlim_{k\geq 0}\mathscr{K}/\mathscr{K}\cap\mathscr{I}^{k+1}\mathscr{F}
$$
is an isomorphism.

{\rm (2)} For any affine open set $U=\Spf A$ of $X$, where $A$ is an adic ring of finite ideal type, the sequence
$$
0\longrightarrow\Gamma(U,\mathscr{K})\longrightarrow\Gamma(U,\mathscr{F})\longrightarrow\Gamma(U,\mathscr{G})\longrightarrow 0
$$
is exact.
\end{prop}

\begin{proof}
Consider the exact sequence
$$
0\longrightarrow\mathscr{K}\longrightarrow\mathscr{F}\longrightarrow\mathscr{G}\longrightarrow 0\leqno{(\ast)}
$$
of $\O_X$-modules.
Taking $\otimes_{\O_X}\O_X/\mathscr{I}^{k+1}$, we get the exact sequence
$$
0\longrightarrow\mathscr{K}/\mathscr{K}\cap\mathscr{I}^{k+1}\mathscr{F}\longrightarrow\mathscr{F}_k\longrightarrow\mathscr{G}_k\longrightarrow 0\leqno{(\ast)_k}
$$
(where $\mathscr{F}_k=\mathscr{F}/\mathscr{I}^{k+1}\mathscr{F}$, etc.) for any $k\geq 0$.
The exact sequences $(\ast)_k$ induce, by passage to the projective limits, the exact sequence
$$
0\longrightarrow\varprojlim_{k\geq 0}\mathscr{K}/\mathscr{K}\cap\mathscr{I}^{k+1}\mathscr{F}\longrightarrow\varprojlim_{k\geq 0}\mathscr{F}_k\longrightarrow\varprojlim_{k\geq 0}\mathscr{G}_k.
$$
Since $\mathscr{F}$ and $\mathscr{G}$ are adically quasi-coherent sheaves, we have $\mathscr{F}\cong\varprojlim_{k\geq 0}\mathscr{F}_k$ and $\mathscr{G}\cong\varprojlim_{k\geq 0}\mathscr{G}_k$.
Hence, comparing the last exact sequence with $(\ast)$, we get (1).
To show (2), first notice that, due to the exact sequence $(\ast)_k$, $\mathscr{K}/\mathscr{K}\cap\mathscr{I}^{k+1}\mathscr{F}$ is a quasi-coherent sheaf on the scheme $X_k=(X,\O_X/\mathscr{I}^{k+1})$; then the assertion follows from \ref{lem-vanishingcohomologyadicallyuseful} (2).
\end{proof}

\begin{cor}\label{cor-closedimmformal11}
Let $A$ be an adic ring of finite ideal type, and $I\subseteq A$ a finitely generated ideal of definition.
Let $\mathscr{F}$ be an adically quasi-coherent sheaf\index{quasi-coherent!adically quasi-coherent OX module@adically --- (a.q.c.) sheaf} of finite type on $X=\Spf A$, and $\mathscr{K}\subseteq\mathscr{F}$ an $\O_X$-submodule such that $\mathscr{G}=\mathscr{F}/\mathscr{K}$ is adically quasi-coherent.
Then $\Gamma(X,\mathscr{K})$ is closed in $\Gamma(X,\mathscr{F})$ with respect to the $I$-adic topology.
\end{cor}

\begin{proof}
Since $\mathscr{F}$ and $\mathscr{G}=\mathscr{F}/\mathscr{K}$ are adically quasi-coherent of finite type, $\Gamma(X,\mathscr{F})$ and $\Gamma(X,\mathscr{G})$ are $I$-adically complete finitely generated $A$-modules (\ref{thm-adicqcoh1} (2)).
Then the result follows from \ref{prop-closedimmformal1} (2) and {\bf \ref{ch-pre}}.\ref{cor-qconsistency1111}.
\end{proof}

\subsubsection{Definitions and first properties}\label{subsub-closedimmformaldef}
\begin{dfn}\label{dfn-closedimmform1}{\rm 
Let $X$ be an adic formal scheme of finite ideal type.
A {\em closed formal subscheme}\index{formal subscheme!closed formal subscheme@closed ---} of $X$ is a formal scheme of the form $(Y,(\O_X/\mathscr{K})|_Y)$, where $\mathscr{K}$ is an ideal of $\O_X$ such that $\O_X/\mathscr{K}$ is adically quasi-coherent\index{quasi-coherent!adically quasi-coherent OX module@adically --- (a.q.c.) sheaf} and $Y$ is the support of the sheaf $\O_X/\mathscr{K}$.}
\end{dfn}

Note that the subset $Y$ is closed in $X$ due to \cite[$\mathbf{0}_{\mathbf{I}}$, (5.2.2)]{EGA}.
Moreover, one can show that the ringed space $(Y,(\O_X/\mathscr{K})|_Y)$ is actually an adic formal scheme of finite ideal type as follows.
Let $U=\Spf A$ be an affine open set of $X$, where $A$ is an adic ring with a finitely generated ideal of definition $I\subseteq A$, and set $K=\Gamma(U,\mathscr{K})$. 
Then by \ref{cor-closedimmformal11} $B=A/K$ is an $I$-adically complete algebra over $A$ and hence is an adic ring.
By \ref{prop-closedimmformal1} (2) and \ref{thm-adicqcoh1} (2) we have $(\O_X/\mathscr{K})|_U=B^{\Delta}$.
Hence $Y$ is covered by affine formal scheme of the form $\Spf B$.

\danger{Note that we do not assume in \ref{dfn-closedimmform1} that the ideal $\mathscr{K}$ itself is adically quasi-coherent.}

\begin{dfn}\label{dfn-closedimmform2}{\rm 
A morphism $i\colon Z\rightarrow X$ of adic formal schemes of finite ideal type is said to be a {\em closed immersion} if it admits a factorization $Z\stackrel{\sim}{\rightarrow}Y\hookrightarrow X$ by an isomorphism onto a closed formal subscheme $Y$ of $X$ followed by the canonical morphism.}
\end{dfn}

\begin{prop}\label{prop-closedimmformal3}
A closed immersion is finite\index{morphism of formal schemes@morphism (of formal schemes)!finite morphism of formal schemes@finite ---} {\rm (\ref{dfn-finitemorform1})}. \hfill$\square$
\end{prop}

We will see in \ref{rem-closedimmformal7} below that the notions of closed formal subschemes and closed immersions thus defined coincide with those in \cite[$\mathbf{I}$, \S10.14]{EGA} in the locally Noetherian\index{formal scheme!Noetherian formal scheme@Noetherian ---!locally Noetherian formal scheme@locally --- ---} case.

\begin{prop}\label{prop-closedimmformal4}
Let $i\colon Z\rightarrow X$ be an adic morphism of adic formal schemes of finite ideal type, and suppose $X$ has an ideal of definition $\mathscr{I}$ of finite type.
For any integer $k\geq 0$ we denote by $i_k\colon Z_k=(Z,\O_Z/\mathscr{I}^{k+1}\O_Z)\rightarrow X_k=(X,\O_X/\mathscr{I}^{k+1})$ the induced morphism of schemes.
Then the following conditions are equivalent$:$
\begin{itemize}
\item[{\rm (a)}] $i$ is a closed immersion$;$
\item[{\rm (b)}] $i_k$ is a closed immersion for any $k\geq 0;$
\item[{\rm (c)}] $i_0$ is a closed immersion.
\end{itemize}
\end{prop}

\begin{proof}
Let us first show the implication (a) $\Rightarrow$ (b). 
Consider the exact sequence
$$
0\longrightarrow\mathscr{K}\longrightarrow\O_X\longrightarrow i_{\ast}\O_Z\longrightarrow 0,
$$
which induces for any $k\geq 0$ the exact sequence
$$
0\longrightarrow\mathscr{K}_k\longrightarrow\O_{X_k}\longrightarrow i_{k\ast}\O_{Z_k}\longrightarrow 0,
$$
where $\mathscr{K}_k=\mathscr{K}/\mathscr{K}\cap\mathscr{I}^{k+1}$.
Since $i_{\ast}\O_Z$ is adically quasi-coherent of finite type (\ref{prop-finitemorphismlocalconstruction}), $i_{k\ast}\O_{Z_k}$ is quasi-coherent on $X_k$.
Hence $\mathscr{K}_k$ is a quasi-coherent sheaf on $X$, and $i_k$ is the closed immersion of schemes corresponding to the quasi-coherent ideal $\mathscr{K}_k$.

Next we show the converse. 
By \ref{prop-affinemorphism1} we already know that $i$ is affine.
By \ref{prop-affinemorphism2} we deduce that $i_{\ast}\O_Z$ is an adically quasi-coherent $\O_X$-algebra.
To show that the map $\O_X\rightarrow i_{\ast}\O_Z$ is surjective, let $\mathscr{K}_k$ be the kernel of $\O_{X_k}\rightarrow i_{k\ast}\O_{Z_k}$ for any $k\geq 0$.
It is easy to see that $\mathscr{K}_k\rightarrow\mathscr{K}_l$ is surjective for $k\geq l$ and hence that the projective system $\{\mathscr{K}_k\}_{k\geq 0}$ is strict.
Hence we have the exact sequence
$$
0\longrightarrow\varprojlim_k\mathscr{K}_k\longrightarrow\O_X\longrightarrow\varprojlim_ki_{k\ast}\O_{Z_k}\longrightarrow 0
$$
({\bf \ref{ch-pre}}.\ref{prop-ML5ML5} (1) and {\bf \ref{ch-pre}}.\ref{rem-coheffaceable} (2)).
On the other hand, we have $\varprojlim_ki_{k\ast}\O_{Z_k}=i_{\ast}\O_Z$ due to {\bf \ref{ch-pre}}.\ref{prop-projlimsheafleftexact0} (2) and \ref{prop-formalindlimschfact}.
Hence the map $\O_X\rightarrow i_{\ast}\O_Z$ is surjective, as desired.
Now since $i(Z)$ clearly coincides with the support of $\O_X/\mathscr{K}$ (which is equal to the support of $\O_{X_k}/\mathscr{K}_k$), we deduce that $i$ is a closed immersion.
Thus we have shown the equivalence of (a) and (b).
The equivalence of (b) and (c) follows from Exercise \ref{exer-closedimmersionnilpotentthickening}.
\end{proof}

\begin{cor}\label{cor-closedimmformal42}
Let $A\rightarrow B$ be an adic map between adic rings of finite ideal type, and consider the morphism $f\colon Y=\Spf B\rightarrow X=\Spf A$ of formal schemes.
Then the following conditions are equivalent$:$
\begin{itemize}
\item[{\rm (a)}] $f$ is a closed immersion$;$
\item[{\rm (b)}] $A\rightarrow B$ is surjective.
\end{itemize}
\end{cor}

\begin{proof}
Let us first show (a) $\Rightarrow$ (b).
Suppose $f$ is a closed immersion.
Then we have the surjective map $\O_X\rightarrow f_{\ast}\O_Y$.
Since $f$ is finite (\ref{prop-closedimmformal3}), $f_{\ast}\O_Y$ is an adically quasi-coherent sheaf of finite type on $X$ by \ref{prop-finitemorphismlocalconstruction}.
Since $A^{\Delta}=\O_X$ and $B^{\Delta}=\O_Y$, the surjectivity of $A\rightarrow B$ follows from \ref{prop-closedimmformal1} (2).

Next we show (b) $\Rightarrow$ (a).
Let $I$ be a finitely generated ideal of definition of $A$.
For any $k\geq 0$ the map $A/I^{k+1}\rightarrow B/I^{k+1}B$ is surjective.
Hence by \ref{prop-closedimmformal4} the morphism $f$ is a closed immersion.
\end{proof}

\begin{cor}\label{cor-closedimmformal43}
Let $i\colon Z\rightarrow X$ be a morphism of adic formal schemes of finite ideal type.
Then the following conditions are equivalent$:$
\begin{itemize}
\item[{\rm (a)}] $i$ is a closed immersion$;$
\item[{\rm (b)}] $i$ is adic\index{morphism of formal schemes@morphism (of formal schemes)!adic morphism of formal schemes@adic ---} and affine\index{morphism of formal schemes@morphism (of formal schemes)!affine morphism of formal schemes@affine ---}, and the morphism $\O_X\rightarrow i_{\ast}\O_Z$ is surjective.
\end{itemize}
\end{cor}

\begin{proof}
The implication (a) $\Rightarrow$ (b) follows from \ref{prop-closedimmformal3}.
Let us show the converse. 
By \ref{prop-affinemorphism2} the sheaf $i_{\ast}\O_Z$ is adically quasi-coherent on $X$.
Let $\mathscr{K}$ be the kernel of $\O_X\rightarrow i_{\ast}\O_Z$.
For any affine open subset $U=\Spf A$ of $X$ (where $A$ is an adic ring of finite ideal type), $i^{-1}(U)$ is affine $i^{-1}(U)=\Spf B$ (where $B$ is an adic ring), and by \ref{prop-closedimmformal1} (2) the induced map $A\rightarrow B$ is surjective.
Hence by \ref{cor-closedimmformal42} the base change $f^{-1}(U)\rightarrow U$ is a closed immersion.
As the underlying morphism of $i$ is, therefore, injective, it suffices to show that the set $i(Z)$ coincides with the support of $\O_X/\mathscr{K}=i_{\ast}\O_Y$.
But this is clear, since they coincide on each affine open subsets.
\end{proof}

\begin{cor}\label{cor-closedimmformalinsert1}
Let $f\colon Y\rightarrow X$ be a morphism of adic formal schemes of finite ideal type, and $X=\bigcup_{\alpha\in L} V_{\alpha}$ an open covering of $X$.
Then $f$ is a closed immersion if and only if for any $\alpha\in L$ the base change $f^{-1}(V_{\alpha})\rightarrow V_{\alpha}$ is a closed immersion. \hfill$\square$
\end{cor}

\begin{prop}\label{prop-closedimmformal5}
{\rm (1)} If $f\colon Z\rightarrow Y$ and $g\colon Y\rightarrow X$ are closed immersions, then so is the composition $g\circ f$.

{\rm (2)} For any $S$-closed immersions $f\colon X\rightarrow Y$ and $g\colon X'\rightarrow Y'$ of adic formal schemes of finite ideal type over an adic formal scheme $S$ of finite ideal type, the induced morphism $f\times_Sg\colon X\times_SY\rightarrow X'\times_SY'$ is a closed immersion.

{\rm (3)} For any $S$-closed immersion $f\colon X\rightarrow Y$ of adic formal schemes of finite ideal type over an adic formal scheme $S$ of finite ideal type and for any morphism $S'\rightarrow S$ of adic formal schemes of finite ideal type, the induced morphism $f_{S'}\colon X\times_SS'\rightarrow Y\times_SS'$ is a closed immersion.
\end{prop}

\begin{proof}
(1) follows easily from \ref{prop-adicmor1} (1), \ref{prop-affmorform2} (1), and \ref{cor-closedimmformal43}.
By \ref{prop-genpropertymorphismadic1} the assertions (2) and (3) follow from \ref{prop-closedimmformal4}, \ref{cor-closedimmformalinsert1}, and \cite[(4.3.6)]{EGAInew}.
\end{proof}

\begin{cor}\label{cor-closedimmformal41}
Let $f\colon X\rightarrow Y$ be a morphism of schemes, and $Z$ a closed subscheme of $Y$ of finite presentation.
If $f$ is a closed immersion, then the formal completion $\widehat{f}\colon\widehat{X}|_{f^{-1}(Z)}\rightarrow\widehat{Y}|_Z$ is a closed immersion.
\end{cor}

\begin{proof}
Since closed immersions are closed under composition, we may apply \ref{prop-genpropertymorphismadic2}, and the corollary follows from \ref{prop-closedimmformal4} and \cite[(4.3.6)]{EGAInew}.
\end{proof}

\subsubsection{Universally rigid-Noetherian case}\label{subsub-closedimmformaladq}
\begin{prop}\label{prop-closedimmformal6}
Let $X$ be a locally universally rigid-Noetherian formal scheme\index{formal scheme!universally rigid-Noetherian formal scheme@universally rigid-Noetherian ---!locally universally rigid-Noetherian formal scheme@locally --- ---} {\rm (\ref{dfn-formalsch})}, and $\mathscr{K}\subseteq\O_X$ an adically quasi-coherent ideal $($resp.\ of finite type$)$.
Consider the locally ringed space $Y=(Y,\O_X/\mathscr{K})$ where $Y$ is the support of the sheaf $\O_X/\mathscr{K}$.
Then $Y$ is a closed formal subscheme of $X$, and the canonical morphism $i\colon Y\hookrightarrow X$ is a closed immersion $($resp.\ of finite presentation$)$.
Moreover, any closed immersion $i\colon Z\hookrightarrow X$ $($resp.\ of finite presentation$)$ can be obtained in this way up to isomorphism by a uniquely determined adically quasi-coherent ideal $\mathscr{K}$ $($resp.\ of finite type$)$.
\end{prop}

\begin{proof}
Let $U=\Spf A$ be an affine open set of $X$, where $A$ is a t.u.\ rigid-Noetherian ring\index{t.u. rigid-Noetherian ring@t.u.\ rigid-Noetherian ring} (\ref{dfn-tuaringadmissible} (1)) with a finitely generated ideal of definition $I\subseteq A$.
Set $K=\Gamma(U,\mathscr{K})$.
Then $K$ is an ideal (resp.\ a finitely generated ideal) of $A$ such that $K^{\Delta}=\mathscr{K}|_U$.
The ring $B=A/K$ is a finitely generated $A$-module and hence is $I$-adically complete (\ref{rem-turigidnoetherianbasicproperties}).
By \ref{prop-lemdeltasheafadicallyquasicoherent1-3} we have $B^{\Delta}=(\O_X/\mathscr{K})|_U$, and hence $\O_X/\mathscr{K}$ is adically quasi-coherent.
Thus the first half of the proposition has been proved.

Conversely, if $i\colon Z\rightarrow X$ is a closed immersion (resp.\ of finite presentation), then for any affine open $U=\Spf A$ of $X$, we have $i^{-1}(U)=\Spf B$ where $A\rightarrow B$ is surjective (\ref{cor-closedimmformal42}).
Let $K$ be the kernel, which is a finitely generated ideal if $i$ is of finite presentation (\ref{prop-topfinpres1}).
By \ref{prop-lemdeltasheafadicallyquasicoherent1-3} we see that $K^{\Delta}$ is an adically quasi-coherent sheaf (resp.\ of finite type) isomorphic to the kernel of $\O_X\rightarrow i_{\ast}\O_Z$ restricted on $U$.
Hence the kernel of $\O_X\rightarrow i_{\ast}\O_Z$ is an adically quasi-coherent ideal (resp.\ of finite type) of $\O_X$, as desired.
\end{proof}

\begin{rem}\label{rem-closedimmformal7}{\rm 
Note that by the proposition we see that our notion of `closed immersion' coincides with that of \cite[$\mathbf{I}$, \S10.14]{EGA} in the locally Noetherian\index{formal scheme!Noetherian formal scheme@Noetherian ---!locally Noetherian formal scheme@locally --- ---} case; indeed, if $X$ is locally Noetherian, then $X$ is locally universally adhesive, and any adically quasi-coherent ideal is a coherent ideal (cf.\ Exercise \ref{exer-adicallyqcohlocallnoetherian}).}
\end{rem}

\subsubsection{Closed immersions and admissible ideals}\label{subsub-closedimmadmissibleideals}
As usual, for an adic formal scheme $X$ with an ideal of definition $\mathscr{I}$ of finite type, a {\em closed subscheme} of $X$ means a closed subscheme of $X_k=(X,\O_X/\mathscr{I}^{k+1})$ for some $k\geq 0$.
Since the quotient $\O_X/\mathscr{J}$ by an admissible ideal $\mathscr{J}\subseteq\O_X$ is a quasi-coherent sheaf, we have the following:
\begin{prop}\label{prop-admissibleideal2}
Let $X$ be an adic formal scheme of finite ideal type.
For any admissible ideal $\mathscr{J}$, let $Y$ be the closed formal subscheme of $X$ corresponding to $\mathscr{J}$.
Then $Y$ is a scheme. \hfill$\square$
\end{prop}

\begin{prop}\label{prop-coradmmissibleideal1xx}
Let $i\colon Y\hookrightarrow X$ be a closed immersion of finite presentation between locally universally rigid-Noetherian formal schemes, and $\mathscr{K}$ an admissible ideal of $\O_Y$.
Let $\mathscr{J}$ be the pull-back of $i_{\ast}\mathscr{K}$ by the map $\O_X\rightarrow i_{\ast}\O_Y$.
Then $\mathscr{J}$ is an admissible ideal of $\O_X$.
\end{prop}

\begin{proof}
We may work in the affine situation $X=\Spf B$ and $Y=\Spf A$, where $A,B$ are t.u.\ rigid-Noetherian rings\index{t.u. rigid-Noetherian ring@t.u.\ rigid-Noetherian ring} and $\varphi\colon B\rightarrow A$ is a surjective adic homomorphism with the finitely generated kernel.
Let $K\subseteq A$ be a finitely generated open ideal such that $\mathscr{K}=K^{\Delta}$, and $J=\varphi^{-1}(K)$, which is a finitely generated open ideal of $B$.
By \ref{prop-affinemorphism2} we know that $i_{\ast}\mathscr{K}$ and $i_{\ast}\O_Y$ are adically quasi-coherent sheaves on $X$ of finite type.
Hence by the exactness of the functor $\cdot^{\Delta}$ (\ref{thm-adicqcohpre1}) we have $\mathscr{J}=J^{\Delta}$.
\end{proof}
\index{immersion!closed immersion of formal schemes@closed --- (of formal schemes)|)}

\subsection{Immersions}\label{sub-immformal}
\index{immersion!immersion of formal schemes@--- (of formal schemes)|(}
\begin{dfn}\label{dfn-immersionformal}{\rm 
A morphism $f\colon Y\rightarrow X$ of adic formal schemes of finite ideal type is said to be an {\em immersion} if it is a composition $f=j\circ i$ of a closed immersion\index{immersion!closed immersion of formal schemes@closed --- (of formal schemes)}  $i$ followed by an open immersion\index{immersion!open immersion of formal schemes@open --- (of formal schemes)} $j$.}
\end{dfn}

Notice that immersions are adic\index{morphism of formal schemes@morphism (of formal schemes)!adic morphism of formal schemes@adic ---} morphisms.
\begin{prop}\label{prop-immersionimmformal4}
Let $f\colon Y\rightarrow X$ be an adic morphism of adic formal schemes of finite ideal type, and suppose $X$ has an ideal of definition $\mathscr{I}$ of finite type.
For any $k\geq 0$ we denote by $f_k\colon Y_k=(Y,\O_Y/\mathscr{I}^{k+1}\O_Y)\rightarrow X_k=(X,\O_X/\mathscr{I}^{k+1})$ the induced morphism of schemes.
Then the following conditions are equivalent$:$
\begin{itemize}
\item[{\rm (a)}] $f$ is an immersion$;$
\item[{\rm (b)}] $f_k$ is an immersion for any $k\geq 0$.
\end{itemize}
\end{prop}

\begin{proof}
Let $f=j\circ i$ be an immersion, where $i\colon Y\hookrightarrow U$ is a closed immersion and $j\colon U\hookrightarrow X$ is an open immersion.
Let $U_k$ for $k\geq 0$ be the scheme defined similarly to $X_k$.
Then, clearly, $j_k\colon U_k\rightarrow X_k$ is an open immersion.
Due to \ref{prop-closedimmformal4} one finds that $i_k\colon Y_k\rightarrow U_k$ is a closed immersion, and thus $f_k=j_k\circ i_k$ is an immersion.

Conversely, suppose $f_k$ is an immersion for any $k\geq 0$.
Take an open subset $U$ of $X$ that contains the subset $f_0(Y)$ as a closed subset.
For any $k\geq 0$ we have the morphism $f_k\colon Y_k\rightarrow U_k$ of schemes, which is a closed immersion by \cite[$\mathbf{I}$, (4.2.2) (ii)]{EGA}.
Hence $f\colon Y\rightarrow U$ is a closed immersion due to \ref{prop-closedimmformal4}; this implies that $f$ is an immersion, as desired.
\end{proof}

\begin{prop}\label{prop-immersionformal1}
Let $f\colon Y\rightarrow X$ be a morphism of adic formal schemes, and $\{V_{\alpha}\}_{\alpha\in L}$ an open covering of $X$.
Then $f$ is an immersion if and only if for any $\alpha\in L$ the base change $f^{-1}(V_{\alpha})\rightarrow V_{\alpha}$ is an immersion.
\end{prop}

\begin{proof}
The `only if' part is clear.
Let us show the other part.
Take an open subset $U_{\alpha}$ of $V_{\alpha}$ for each $\alpha\in L$ such that the immersion $f^{-1}(V_{\alpha})\rightarrow V_{\alpha}$ factors by the closed immersion $f^{-1}(V_{\alpha})\hookrightarrow U_{\alpha}$.
Set $U=\bigcup_{\alpha\in L}U_{\alpha}$.
Then by \ref{cor-closedimmformalinsert1} the morphism $Y\rightarrow U$ is a closed immersion.
\end{proof}

\begin{prop}\label{prop-exerimmersionbasechangeunderlyingtopspace}
Let $f\colon X\rightarrow Y$ be a morphism of formal schemes where $Y$ is adic of finite ideal type, and $i\colon Z\hookrightarrow Y$ an immersion.
Then $X\times_YZ$ has the underlying topological space homeomorphic to $f^{-1}(i(Z))$.
\end{prop}

\begin{proof}
The assertion is shown in \ref{prop-sepqsepformal1} in the case where $i$ is an open immersion.
Hence it suffices to show the claim in case $i$ is a closed immersion.
To show this, again by \ref{prop-sepqsepformal1} we may assume that $X$ and $Y$ have ideals of definition.
Then the claim follows from \ref{prop-closedimmformal4} and \cite[$\mathbf{I}$, (4.4.1)]{EGA}.
\end{proof}

\begin{lem}\label{lem-immersionformal5lem}
Let $Z\hookrightarrow Y$ be an open immersion, and $Y\hookrightarrow X$ a closed immersion\index{immersion!closed immersion of formal schemes@closed --- (of formal schemes)}, where $X,Y,Z$ are adic formal schemes of finite ideal type.
Then the composition $Z\hookrightarrow X$ is an immersion.
\end{lem}

\begin{proof}
Take an open subset $U\subseteq X$ such that $Z=Y\cap U$.
Then $Z\hookrightarrow U$ is a closed immersion.
\end{proof}

\begin{prop}\label{prop-immersionformal5}
{\rm (1)} If $f\colon Z\rightarrow Y$ and $g\colon Y\rightarrow X$ are immersions, then so is the composition $g\circ f$.

{\rm (2)} For any $S$-immersions $f\colon X\rightarrow Y$ and $g\colon X'\rightarrow Y'$ of adic formal schemes of finite ideal type over an adic formal scheme $S$ of finite ideal type, the induced morphism $f\times_Sg\colon X\times_SY\rightarrow X'\times_SY'$ is an immersion.

{\rm (3)} For any $S$-immersion $f\colon X\rightarrow Y$ of adic formal schemes of finite ideal type over an adic formal scheme $S$ of finite ideal type and for any morphism $S'\rightarrow S$ of adic formal schemes of finite ideal type, the induced morphism $f_{S'}\colon X\times_SS'\rightarrow Y\times_SS'$ is an immersion.
\end{prop}

\begin{proof}
(1) follows from \ref{lem-immersionformal5lem} combined with \ref{prop-openimmform1} (1) and \ref{prop-closedimmformal5} (1).
In view of \ref{prop-genpropertymorphismadic1} the assertions (2) and (3) follow from \ref{prop-immersionimmformal4}, \ref{prop-immersionformal1}, and \cite[(4.3.6)]{EGAInew}.
\end{proof}

\begin{cor}\label{cor-immersionimmformal41}
Let $f\colon X\rightarrow Y$ be a morphism of schemes, and $Z$ a closed subscheme of $Y$ of finite presentation.
If $f$ is an immersion, then the formal completion $\widehat{f}\colon\widehat{X}|_{f^{-1}(Z)}\rightarrow\widehat{Y}|_Z$ is an immersion.
\end{cor}

\begin{proof}
Since immersions are closed under composition, we may apply \ref{prop-genpropertymorphismadic2}, and the corollary follows from \ref{prop-immersionimmformal4} and \cite[(4.3.6)]{EGAInew}.
\end{proof}
\index{immersion!immersion of formal schemes@--- (of formal schemes)|)}

\subsection{Surjective, closed and universally closed morphisms}\label{sub-surjectiveclosedunivclosed}
\subsubsection{Surjective morphisms}\label{subsub-surjectivemorformal}
\index{morphism of formal schemes@morphism (of formal schemes)!surjective morphism of formal schemes@surjective ---|(}
\begin{dfn}\label{dfn-surjectivemorphismsformalschemes}{\rm 
A morphism $f\colon X\rightarrow Y$ of formal schemes is said to be {\em surjective} if the underlying continuous mapping is surjective.}
\end{dfn}

If there exist ideals of definition $\mathscr{I}$ and $\mathscr{J}$ of $X$ and $Y$, respectively, such that $\mathscr{J}\O_X\subseteq\mathscr{I}$, then $f$ is surjective if and only if the induced morphism $(X,\O_X/\mathscr{I})\rightarrow(Y,\O_Y/\mathscr{J})$ of schemes is surjective.
The following proposition follows immediately from \ref{prop-sepqsepformal1}:
\begin{prop}\label{prop-surjectionstablezar}
Let $f\colon Y\rightarrow X$ be a morphism of formal schemes, and $\{V_{\alpha}\}_{\alpha\in L}$ an open covering of $X$.
Then $f$ is a surjective if and only if for any $\alpha\in L$ the base change $f^{-1}(V_{\alpha})\rightarrow V_{\alpha}$ is surjective. \hfill$\square$
\end{prop}

\begin{prop}\label{prop-surjectivemorformal}
{\rm (1)} The composition of two surjective morphisms of formal schemes is surjective.

{\rm (2)} If $S$ is a formal scheme and if $f\colon X\rightarrow X'$ and $g\colon Y\rightarrow Y'$ are surjective morphisms of $S$-formal schemes, then $f\times_Sg\colon X\times_SY\rightarrow X'\times_SY'$ is surjective.

{\rm (3)} If $S$ is a formal scheme and if $f\colon X\rightarrow Y$ is a surjective morphism of $S$-formal schemes, then for any $S'\rightarrow S$ the induced morphism $f_{S'}\colon X\times_SS'\rightarrow Y\times_SS'$ is surjective.
\end{prop}

The proof uses the following lemma, which can be shown similarly to the scheme case (cf.\ \cite[(3.6.2)]{EGAInew}):
\begin{lem}\label{lem-exerlemsurjectivemorformal}
Let $f\colon X\rightarrow Y$ be a surjective morphism, $K$ a field, and $\Spec K\rightarrow Y$ a map of formal schemes.
Then $X\times_Y\Spec K$ is non-empty. \hfill$\square$
\end{lem}

\begin{proof}[Proof of Proposition {\rm \ref{prop-surjectivemorformal}}]
(1) is clear.
As we saw in \S\ref{subsub-genpropertymorphism}, (2) and (3) follow from the special case of (3) with $S=Y$.
Let $f\colon X\rightarrow Y$ be a surjective map of formal schemes, and $Y'\rightarrow Y$ a map of formal schemes.
We want to show that $f'\colon X'=X\times_YY'\rightarrow Y'$ is surjective.
By \ref{prop-surjectionstablezar} we may assume that $Y$ and $Y'$ are affine; set $Y=\Spf B$ and $Y'=\Spf B'$.
Let $\mathfrak{q}\in Y'$ be an open prime ideal of $B'$.
Let $K=\Frac(B'/\mathfrak{q})$, and $\Spec K\rightarrow Y'$ the natural map.
Then by \ref{lem-exerlemsurjectivemorformal} $X\times_Y\Spec K=X'\times_{Y'}\Spec K$ is non-empty, and any point of $X'\times_{Y'}\Spec K$ is mapped to a point of $X'$ that is mapped to $\mathfrak{q}$ by $f'$.
\end{proof}
\index{morphism of formal schemes@morphism (of formal schemes)!surjective morphism of formal schemes@surjective ---|)}

\subsubsection{Closed and universally closed morphisms}\label{subsub-closeduniversallyclosedmorphisms}
\index{morphism of formal schemes@morphism (of formal schemes)!closed morphism of formal schemes@closed ---|(}
\index{morphism of formal schemes@morphism (of formal schemes)!closed morphism of formal schemes@closed ---!universally closed morphism of formal schemes@universally --- ---|(}
\begin{dfn}\label{dfn-closedunivclosed}{\rm 
Let $f\colon X\rightarrow Y$ be a morphism of formal schemes.

(1) We say that $f$ is {\em closed} if the underlying continuous map of $f$ is closed.

(2) We say that $f$ is {\em universally closed} if for any morphism $S\rightarrow Y$ of formal schemes the induced morphism $f_S\colon X\times_YS\rightarrow S$ is closed.}
\end{dfn}

\begin{prop}\label{prop-closedstablezar}
Let $f\colon Y\rightarrow X$ be a morphism of formal schemes, and $\{V_{\alpha}\}_{\alpha\in L}$ an open covering of $X$.
Then $f$ is universally closed if and only if for any $\alpha\in L$ the base change $f^{-1}(V_{\alpha})\rightarrow V_{\alpha}$ is universally closed.
\end{prop}

\begin{proof}
The `only if' part is clear.
To show the `if' part, by \ref{prop-openimmform1} (3) it suffices to show that the map $f\colon Y\rightarrow X$ is closed.
But this follows easily from \ref{prop-sepqsepformal1}.
\end{proof}

\begin{prop}\label{prop-closedadicenough}
Let $f\colon X\rightarrow Y$ be a morphism of formal schemes, and suppose $Y$ is a scheme.
Then $f$ is universally closed if and only if for any morphism $S\rightarrow Y$ of {\em schemes} the induced morphism $f_S\colon X\times_YS\rightarrow S$ is closed.
\end{prop}

In particular, our notion of `universally closed' restricted to morphisms of schemes coincides with the usual one in scheme theory.
\begin{proof}
The `only if' part is trivial.
Let us show the converse.
Let $S\rightarrow Y$ be a morphism of formal schemes.
By \ref{prop-closedstablezar} we may assume that $S$ is affine $S=\Spf R$ where $R$ is an admissible ring.
Let $J\subseteq R$ be an ideal of definition, and consider $S_0=\Spec R_0$ with $R_0=R/J$.
Then $S_0$ is a scheme, and hence $f_{S_0}\colon X\times_YS_0\rightarrow S_0$ is closed.
To show that $f_S\colon X\times_YS\rightarrow S$ is closed, it suffices to show that $X\times_YS$ and $X\times_YS_0$ have the same underlying topological space.
To show this, we may assume that $X$ and $Y$ are affine $X=\Spf A$ and $Y=\Spec B$.
Let $I\subseteq A$ be an ideal of definition.
Consider the ideal $H=\image(I\otimes R)+\image(A\otimes J)$ of $A\otimes_BR$, and set $(X\times_YS)_0=\Spec A\otimes_BR/H$.
Then $X\times_YS$ and $(X\times_YS)_0$ has the same underlying topological space.
On the other hand, we have $X\times_YS_0=\Spf A\widehat{\otimes}_BR/J$.
Observe that the kernel of $A\widehat{\otimes}_BR\rightarrow A\widehat{\otimes}_BR/J$ is contained in the kernel of $A\widehat{\otimes}_BR\rightarrow A\otimes_BR/H$.
Hence we have the chain of closed immersions $(X\times_YS)_0\hookrightarrow X\times_YS_0\hookrightarrow X\times_YS$, from which the claim follows.
\end{proof}

\begin{prop}\label{prop-univclosedstaility}
{\rm (1)} The composition of two closed $($resp.\ universally closed$)$ morphisms of formal schemes is closed $($resp.\ universally closed$)$.

{\rm (2)} If $f\colon X\rightarrow X'$ and $g\colon Y\rightarrow Y'$ are two universally closed $S$-morphisms of formal schemes, then $f\times_Sg\colon X\times_SY\rightarrow X'\times_SY'$ is universally closed.

{\rm (3)} If $f\colon X\rightarrow Y$ is a universally closed $S$-morphism of formal schemes and $S'\rightarrow S$ is a morphism, then $f_{S'}\colon X\times_SS'\rightarrow Y\times_SS'$ is universally closed.
\end{prop}

\begin{proof}
(1) is clear.
The assertions (2) and (3) follow from the special case of (3) with $S=Y$ (\S\ref{subsub-genpropertymorphism}), which follows immediately from the definition.
\end{proof}

\begin{prop}\label{prop-univclosedadicred}
Let $f\colon X\rightarrow Y$ be an adic morphism of adic formal schemes, and suppose $Y$ has an ideal of definition $\mathscr{I}$.
Set $X_k=(X,\O_X/\mathscr{I}^{k+1}\O_X)$ and $Y_k=(Y,\O_Y/\mathscr{I}^{k+1})$, and denote by $f_k\colon X_k\rightarrow Y_k$ the induced morphism of schemes.
Then the following conditions are equivalent$:$
\begin{itemize}
\item[{\rm (a)}] $f$ is closed $($resp.\ universally closed$);$
\item[{\rm (b)}] $f_k$ is closed $($resp.\ universally closed$)$ for any $k\geq 0;$
\item[{\rm (c)}] $f_0$ is closed $($resp.\ universally closed$)$.
\end{itemize}
\end{prop}

\begin{proof}
The assertion for `closed' is clear.
For the other, use \ref{prop-closedadicenough}.
\end{proof}

\begin{cor}\label{cor-univclosedadicred0}
Let $f\colon X\rightarrow Y$ be an immersion between adic formal schemes of finite ideal type.
Then $f$ is a closed immersion if and only if it is closed.
\end{cor}

\begin{proof}
The ``only if'' part is clear.
Suppose $f$ is closed.
To show that $f$ is a closed immersion, we may work locally on $Y$, and thus we may assume that $Y$ has an ideal of definition $\mathscr{I}$ of finite type.
With the notation as in \ref{prop-univclosedadicred}, the morphism $f_0$ of schemes is closed, and is an immersion (due to \ref{prop-immersionimmformal4}).
Hence $f_0$ is a closed immersion (e.g.\ \cite{EGA}, {\bf I}, (4.2.2)).
Then it follows from \ref{prop-closedimmformal4} that $f$ is a closed immersion.
\end{proof}

\begin{cor}\label{cor-univclosedadicred}
Let $f\colon X\rightarrow Y$ be a morphism of schemes, and $Z$ a closed subscheme of $Y$ of finite presentation.
If $f$ is universally closed, then the formal completion $\widehat{f}\colon\widehat{X}|_{f^{-1}(Z)}\rightarrow\widehat{Y}|_Z$ is universally closed. \hfill$\square$
\end{cor}
\index{morphism of formal schemes@morphism (of formal schemes)!closed morphism of formal schemes@closed ---!universally closed morphism of formal schemes@universally --- ---|)}
\index{morphism of formal schemes@morphism (of formal schemes)!closed morphism of formal schemes@closed ---|)}

\subsection{Separated morphisms}\label{sub-sepqsepformal}
\index{morphism of formal schemes@morphism (of formal schemes)!separated morphism of formal schemes@separated ---|(}
\subsubsection{Definition and fundamental properties}\label{subsub-sepqsepformal}
\begin{prop}\label{prop-diagonalimmersionadic}
Let $f\colon X\rightarrow Y$ be a morphism of adic formal schemes, and suppose $X$ is of finite ideal type.
Then the diagonal map\index{diagonal map} $\Delta_X\colon X\rightarrow X\times_YX$ is an immersion.
\end{prop}

Notice that we have already shown in \ref{adicmorphismdiagonal} that the diagonal map $\Delta_X$ is adic\index{morphism of formal schemes@morphism (of formal schemes)!adic morphism of formal schemes@adic ---}.
\begin{proof}
Let $V\subset Y$ be an open subset, and $U=f^{-1}(V)$.
Then the base change of $\Delta_X$ by the open immersion $U\times_VU\rightarrow X\times_YX$ (\ref{prop-openimmform1} (2)) coincides with $\Delta_U\colon U\rightarrow U\times_VU$.
Since the image of $\Delta_X$ is contained in the union of open subsets of the form $U\times_VU$, by \ref{prop-immersionformal1} we may assume that $Y$ is affine $Y=\Spf B$.
Let $X=\bigcup_{\alpha\in L}U_{\alpha}$ be an affine open covering.
By \ref{prop-openimmform1} (2) the canonical morphism $U_{\alpha}\times_YU_{\beta}\rightarrow X\times_YX$ is an open immersion.
As the image of $\Delta_X$ is contained in the open subset $\bigcup_{\alpha\in L}U_{\alpha}\times_YU_{\alpha}$, it suffices to show that the map $X\rightarrow\bigcup_{\alpha\in L}U_{\alpha}\times_YU_{\alpha}$ is a closed immersion.
By \cite[$\mathbf{0}$, (1.4.8)]{EGAInew} we have $\Delta^{-1}_X(U_{\alpha}\times_YU_{\alpha})=U_{\alpha}$.
Hence by \ref{cor-closedimmformalinsert1} it suffices to show that for each $\alpha\in L$ the diagonal map $\Delta_{U_{\alpha}}\colon U_{\alpha}\rightarrow U_{\alpha}\times_YU_{\alpha}$ is a closed immersion.
Set $U_{\alpha}=\Spf A_{\alpha}$ for $\alpha\in L$, where $A_{\alpha}$ is an adic ring of finite ideal type.
Then the diagonal map corresponds to $A_{\alpha}\widehat{\otimes}_BA_{\alpha}\rightarrow A_{\alpha}$, which is clearly surjective.
Now by \ref{cor-closedimmformal42} the proof is completed.
\end{proof}

\begin{dfn}[{\cite[$\mathbf{I}$, \S10.15]{EGA}}]\label{dfn-sepmorformal}{\rm
A morphism $f\colon X\rightarrow Y$ of formal schemes is said to be {\em separated} if the image $\Delta_X(X)$ of the diagonal map\index{diagonal map} $\Delta_X\colon X\rightarrow X\times_YX$ is a closed subset of $X\times_YX$. 
A formal scheme $X$ is said to be {\em separated}\index{formal scheme!separated formal scheme@separated ---} if it is separated over $\Spec\Z$.}
\end{dfn}

\begin{prop}\label{prop-qsepformal2}
Separated morphisms are quasi-separated\index{morphism of formal schemes@morphism (of formal schemes)!quasi-separated morphism of formal schemes@quasi-separated ---}\index{quasi-separated!quasi-separated morphism of formal schemes@--- morphism (of formal schemes)}. \hfill$\square$
\end{prop}

By \ref{prop-diagonalimmersionadic} we have:
\begin{prop}\label{cor-sepmorformal12}
Let $f\colon X\rightarrow Y$ be a morphism of adic formal schemes, where $X$ is supposed to be of finite ideal type.
Then $f$ is separated if and only if the diagonal map $\Delta_X\colon X\rightarrow X\times_YX$ is a closed immersion. \hfill$\square$
\end{prop}

\begin{prop}\label{prop-sepmorformal24}
{\rm (1)} The composition of two separated morphisms of adic formal schemes of finite ideal type is separated.

{\rm (2)} If the composition $g\circ f$ of two morphisms of adic formal schemes of finite ideal type is separated, then $f$ is separated.

{\rm (3)} An open immersion is separated.
\end{prop}

\begin{proof}
We use the diagram $(\ast)$ in the proof of \ref{prop-qsepmorformal3}.
Suppose $f\colon X\rightarrow Y$ and $g\colon Y\rightarrow Z$ is separated.
Then $\Delta_f$ and $\Delta_g$ are closed immersions, and hence $X\times_YX\rightarrow X\times_ZX$ is a closed immersion (\ref{prop-closedimmformal5} (3)).
This implies that $\Delta_{g\circ f}$ is a closed immersion (\ref{prop-closedimmformal5} (1)), that is, $g\circ f$ is separated, whence (1).

Suppose $g\circ f$ is separated.
Since $X\times_YX\rightarrow X\times_ZX$ is an immersion (\ref{prop-diagonalimmersionadic} and \ref{prop-immersionformal5} (3)), $\Delta_f(X)$ coincides with the pull-back of $\Delta_{g\circ f}(X)$.
Since $\Delta_{g\circ f}(X)$ is closed in $X\times_ZX$, $\Delta_f(X)$ is closed in $X\times_YX$, whence (2).

Finally, if $j\colon U\hookrightarrow X$ is an open immersion, then the diagonal map is an isomorphism, which is clearly a closed immersion, whence (3).
\end{proof}

\begin{prop}\label{prop-sepmorformal0}
Let $f\colon Y\rightarrow X$ be a morphism of adic formal schemes of finite ideal type, and $\{V_{\alpha}\}_{\alpha\in L}$ an open covering of $X$.
Then $f$ is separated if and only if for any $\alpha\in L$ the base change $f^{-1}(V_{\alpha})\rightarrow V_{\alpha}$ is separated.
\end{prop}

\begin{proof}
Set $U_{\alpha}=f^{-1}(V_{\alpha})$.
Then the image $\Delta_X(X)$ of the diagonal map $\Delta_X\colon X\rightarrow X\times_YX$ is contained in the open subset $\bigcup_{\alpha\in L}U_{\alpha}\times_{V_{\alpha}}U_{\alpha}$.
By \cite[$\mathbf{0}$, (1.4.8)]{EGAInew} we have $\Delta^{-1}_X(U_{\alpha}\times_{V_{\alpha}}U_{\alpha})=U_{\alpha}$.
Then the assertion follows from \ref{prop-sepmorformal24} (1), (2), and (3) and \ref{cor-closedimmformalinsert1}.
\end{proof}

\begin{cor}\label{cor-sepmorformal01}
An affine morphism\index{morphism of formal schemes@morphism (of formal schemes)!affine morphism of formal schemes@affine ---} between adic formal schemes of finite ideal type is separated.
\end{cor}

\begin{proof}
In view of \ref{prop-sepmorformal0} it suffices to show that a morphism $X=\Spf A\rightarrow Y=\Spf B$, where $A$ and $B$ are adic rings of finite ideal type, is separated.
Since $B\widehat{\otimes}_AB\rightarrow B$ is clearly surjective, the assertion follows from \ref{adicmorphismdiagonal} and \ref{cor-closedimmformal42}.
\end{proof}

\begin{prop}\label{prop-sepmorformal2}
Let $S$ be an adic formal scheme $S$ of finite ideal type.

{\rm (1)} For separated $S$-morphisms $f\colon X\rightarrow Y$ and $g\colon X'\rightarrow Y'$ of adic formal schemes of finite ideal type over $S$, the morphism $f\times_Sg\colon X\times_SY\rightarrow X'\times_SY'$ is separated.

{\rm (2)} For a separated $S$-morphism $f\colon X\rightarrow Y$ of adic formal schemes of finite ideal type over $S$ and a morphism $S'\rightarrow S$ of adic formal schemes of finite ideal type, the induced morphism $f_{S'}\colon X\times_SS'\rightarrow Y\times_SS'$ is separated.
\end{prop}

\begin{proof}
As we have seen in \S\ref{subsub-genpropertymorphismadic}, (1) and (2) follow from the special case of (2) with $S=Y$.
We use the diagram $(\ast\ast)$ in the proof of \ref{prop-qsepmorformal3}.
By \ref{prop-closedimmformal5} (3), if $\Delta_f$ is a closed immersion, then so is $\Delta_{f'}$.
\end{proof}

\begin{prop}\label{prop-sepmorformal1}
Let $f\colon X\rightarrow Y$ be an adic morphism of adic formal schemes of finite ideal type, and suppose $Y$ has an ideal of definition $\mathscr{I}$ of finite type.
Set $X_k=(X,\O_X/\mathscr{I}^{k+1}\O_X)$ and $Y_k=(Y,\O_Y/\mathscr{I}^{k+1})$, and denote by $f_k\colon X_k\rightarrow Y_k$ the induced morphism of schemes.
Then the following conditions are equivalent$:$
\begin{itemize}
\item[{\rm (a)}] $f$ is separated$;$
\item[{\rm (b)}] $f_k$ is separated for any $k\geq 0;$
\item[{\rm (c)}] $f_0$ is separated.
\end{itemize}
\end{prop}

\begin{proof}
This follows from \ref{cor-prodformal21} and the definition of separatedness.
\end{proof}

By \ref{prop-genpropertymorphismadic2}, \ref{prop-sepmorformal1}, and \cite[(5.3.1)]{EGAInew} we have:
\begin{cor}\label{cor-sepmorformal1}
Let $f\colon X\rightarrow Y$ be a morphism of schemes, and $Z$ a closed subscheme of $Y$ of finite presentation.
If $f$ is separated, then the formal completion $\widehat{f}\colon\widehat{X}|_{f^{-1}(Z)}\rightarrow\widehat{Y}|_Z$ is separated. \hfill$\square$
\end{cor}

\begin{prop}\label{prop-separatedformalaffine}
Let $B$ be an adic ring of finite ideal type, and $f\colon X\rightarrow Y=\Spf B$ a morphism between adic formal schemes of finite ideal type.
Let $X=\bigcup_{\alpha\in L}U_{\alpha}$ be an affine open covering of $X$ where $U_{\alpha}=\Spf A_{\alpha}$ by an adic ring $A_{\alpha}$ of finite ideal type for each $\alpha\in L$.
Then $f$ is separated if and only if the following conditions are satisfied for any $\alpha,\beta\in L:$
\begin{itemize}
\item[{\rm (a)}] the intersection $U_{\alpha}\cap U_{\beta}$ is affine$;$ $U_{\alpha}\cap U_{\beta}=\Spf A_{\alpha\beta}$ by an adic ring $A_{\alpha\beta}$ of finite ideal type$;$
\item[{\rm (b)}] the induced map $A_{\alpha}\widehat{\otimes}_BA_{\beta}\rightarrow A_{\alpha\beta}$ is surjective.
\end{itemize}
\end{prop}

\begin{proof}
First notice that the diagram
$$
\xymatrix{X\ar[r]^(.4){{\Delta}_X}&X\times_YX\\ U_{\lambda}\cap U_{\mu}\ar@{^{(}->}[u]\ar[r]&U_{\lambda}\times_YU_{\mu}\ar@{^{(}->}[u]}
$$
is Cartesian (\ref{prop-sepqsepformal1} and \ref{cor-sepqsepformal21}).
Suppose $f$ is separated, that is, $\Delta_X$ is a closed immersion (\ref{cor-sepmorformal12}).
Then by \ref{prop-closedimmformal5} (3) the map $U_{\alpha}\cap U_{\beta}\rightarrow U_{\alpha}\times_YU_{\beta}$ is a closed immersion.
Since $U_{\alpha}\times_YU_{\beta}=\Spf A_{\alpha}\widehat{\otimes}_BA_{\beta}$ is affine, $U_{\alpha}\cap U_{\beta}$ is affine (\ref{prop-closedimmformal3} (1)) and the map $A_{\alpha}\widehat{\otimes}_BA_{\beta}\rightarrow\Gamma(U_{\alpha}\cap U_{\beta},\O_X)$ is surjective (\ref{cor-closedimmformal42}).

Conversely, if any $U_{\alpha}\cap U_{\beta}=\Spf A_{\alpha\beta}$ is affine and $A_{\alpha}\widehat{\otimes}_BA_{\beta}\rightarrow A_{\alpha\beta}$ is surjective, then $U_{\alpha}\cap U_{\beta}\rightarrow U_{\alpha}\times_YU_{\beta}$ is a closed immersion (\ref{cor-closedimmformal42}).
Since the open subsets of $X\times_YX$ of the form $U_{\alpha}\times_YU_{\beta}$ covers $X\times_YX$, $\Delta_X$ is a closed immersion due to \ref{cor-closedimmformalinsert1}.
\end{proof}

\subsubsection{Separatedness and properties of morphisms}\label{subsub-seppropertymor}
\begin{prop}\label{prop-closedimmform111}
Let $f\colon X\rightarrow Y$ be a morphism of adic formal schemes of finite ideal type, and $g\colon Y\rightarrow Z$ a separated morphism of adic formal schemes of finite ideal type.
Then the graph $\Gamma_f\colon X\rightarrow X\times_ZY$ is a closed immersion.
If, moreover, $X,Y,Z$ are locally universally rigid-Noetherian\index{formal scheme!universally rigid-Noetherian formal scheme@universally rigid-Noetherian ---!locally universally rigid-Noetherian formal scheme@locally --- ---} and $g$ is of finite type, then $\Gamma_f$ is of finite presentation.
\end{prop}

\begin{proof}
Since the diagram
$$
\xymatrix{X\ar[d]_f\ar[r]^(.4){\Gamma_f}&X\times_ZY\ar[d]^{f\times_ZY}\\ Y\ar[r]_(.4){\Delta_g}&Y\times_ZY}
$$
is Cartesian, the first assertion follows from \ref{cor-sepmorformal12}.
The last assertion follows from \ref{prop-closedimmformal6}, \ref{prop-topfinpres2} (4), and the following lemma.
\end{proof}

\begin{lem}\label{lem-sepmorformal4}
Let $A$ be an adic ring of finite ideal type, and $B$ an $A$-algebra topologically of finite type.
Then the kernel of the morphism $B\widehat{\otimes}_AB\rightarrow B$ $(x\otimes y\mapsto xy)$ is generated by $\{1\otimes s-s\otimes 1\,|\,s\in S\}$, where $S\subseteq B$ is a topological generator of $B$ over $A$.
\end{lem}

\begin{proof}
It suffices to show the assertion in the case where $B$ is the restricted power series ring $B=A\dl X_1,\ldots,X_n\dr$.
In this situation, what to show is the following: the kernel of the morphism 
$$
A\dl X_1,\ldots,X_n,Y_1,\ldots,Y_n\dr\longrightarrow A\dl Z_1,\ldots,Z_n\dr 
$$
sending $X_i$ and $Y_i$ to $Z_i$ for $i=1,\ldots,n$ is the ideal $\mathfrak{a}$ generated by the elements $X_i-Y_i$ ($1\leq i\leq n$).
Since the similar assertion is known to hold for polynomial rings, this is equivalent to that the ideal $\mathfrak{a}$ is closed.
We show this by induction with respect to $n$.
The case $n=0$ is trivial.
Consider $A'=A\dl X_n,Y_n\dr/(X_n-Y_n)$, which is $I$-adically complete (where $I\subseteq A$ is an ideal of definition) due to Exercise \ref{exer-principalclosed}; in particular, we have $A'\cong A\dl Z_n\dr$.
The morphism in question factors through $A'\dl X_1,\ldots,X_{n-1},Y_1,\ldots,Y_{n-1}\dr\rightarrow A'\dl Z_1,\ldots,Z_{n-1}\dr\cong A\dl Z_1,\ldots,Z_n\dr$.
By induction we have $A'\dl X_1,\ldots,X_{n-1},Y_1,\ldots,Y_{n-1}\dr/(X_i-Y_i\,|\,i=1,\ldots,n-1)\cong A\dl Z_1,\ldots,Z_n\dr$.
Now again by Exercise \ref{exer-principalclosed} we have $A\dl X_1,\ldots,X_n,Y_1,\ldots,Y_n\dr/(X_n-Y_n)\cong A'\dl X_1,\ldots,X_{n-1},Y_1,\ldots,Y_{n-1}\dr$, and hence we have $A\dl X_1,\ldots,X_n,Y_1,\ldots,Y_n\dr/\mathfrak{a}\cong A\dl Z_1,\ldots,Z_n\dr$, as desired.
\end{proof}

\begin{prop}\label{prop-propertyP3}
Let $P$ be a property of arrows in the category $\Ac\Fs$ satisfying {\bf (I)}, {\bf (C)} in {\rm {\bf \ref{ch-pre}}, \S\ref{subsub-genpropertymorphism}} and the mutually equivalent conditions {\boldmath $(\mathbf{B}_i)$} for $i=1,2,3$ $($with $\mathscr{D}_Q=\Ac\Fs)$ in {\rm {\bf \ref{ch-pre}}, \S\ref{subsub-genpropertymorphismadic}}.
Suppose that any closed immersion satisfies $P$.
Then if $f\colon X\rightarrow Y$ and $g\colon Y\rightarrow Z$ are morphisms of adic formal schemes of finite ideal type such that $g\circ f$ satisfies $P$ and $g$ is separated, $f$ satisfies $P$.
\end{prop}

\begin{proof}
The morphism $f$ coincides with the following composition:
$$
X\stackrel{\Gamma_f}{\longrightarrow}X\times_ZY\stackrel{(g\circ f)\times_Z\id_Y}{\longrightarrow}Y.
$$
Since $g$ is separated, the first arrow is a closed immersion (\ref{prop-closedimmform111}), whence the claim.
\end{proof}

\begin{cor}\label{cor-propertyP3}
Let $f\colon X\rightarrow Y$ and $g\colon Y\rightarrow Z$ be morphisms of adic formal schemes of finite ideal type, and suppose $g$ is separated.
Then if $g\circ f$ satisfies one of the following conditions, then so does $f:$ 
\begin{itemize}
\item[{\rm (a)}] adic\index{morphism of formal schemes@morphism (of formal schemes)!adic morphism of formal schemes@adic ---},
\item[{\rm (b)}] locally of finite type\index{morphism of formal schemes@morphism (of formal schemes)!morphism of formal schemes locally of finite type@--- locally of finite type} $($resp.\ of finite type\index{morphism of formal schemes@morphism (of formal schemes)!morphism of formal schemes of finite type@--- of finite type}$)$,
\item[{\rm (c)}] quasi-compact\index{morphism of formal schemes@morphism (of formal schemes)!quasi-compact morphism of formal schemes@quasi-compact ---}\index{quasi-compact!quasi-compact morphism of formal schemes@--- morphism (of formal schemes)} $($resp.\ quasi-separated\index{morphism of formal schemes@morphism (of formal schemes)!quasi-separated morphism of formal schemes@quasi-separated ---}\index{quasi-separated!quasi-separated morphism of formal schemes@--- morphism (of formal schemes)}, resp.\ coherent\index{morphism of formal schemes@morphism (of formal schemes)!coherent morphism of formal schemes@coherent ---}\index{coherent!coherent morphism of formal schemes@--- morphism (of formal schemes)}$)$,
\item[{\rm (d)}] affine\index{morphism of formal schemes@morphism (of formal schemes)!affine morphism of formal schemes@affine ---} adic\index{morphism of formal schemes@morphism (of formal schemes)!adic morphism of formal schemes@adic ---},
\item[{\rm (e)}] finite\index{morphism of formal schemes@morphism (of formal schemes)!finite morphism of formal schemes@finite ---},
\item[{\rm (f)}] closed immersion\index{immersion!closed immersion of formal schemes@closed --- (of formal schemes)} $($resp.\ immersion\index{immersion!immersion of formal schemes@--- (of formal schemes)}$)$,
\item[{\rm (g)}] universally closed\index{morphism of formal schemes@morphism (of formal schemes)!closed morphism of formal schemes@closed ---!universally closed morphism of formal schemes@universally --- ---}.\hfill$\square$
\end{itemize}
\end{cor}
\index{morphism of formal schemes@morphism (of formal schemes)!separated morphism of formal schemes@separated ---|)}

\subsection{Proper morphisms}\label{sub-propermorformal}
\index{morphism of formal schemes@morphism (of formal schemes)!proper morphism of formal schemes@proper ---|(}
\begin{dfn}\label{dfn-propermorformal}{\rm 
A morphism $f\colon X\rightarrow Y$ of adic formal schemes of finite ideal type is said to be {\em proper} if it is separated of finite type\index{morphism of formal schemes@morphism (of formal schemes)!morphism of formal schemes of finite type@--- of finite type} and is universally closed\index{morphism of formal schemes@morphism (of formal schemes)!closed morphism of formal schemes@closed ---!universally closed morphism of formal schemes@universally --- ---}.}
\end{dfn}

By \ref{prop-closedadicenough} a proper morphism of schemes $f\colon X\rightarrow Y$ is also proper as a morphism of ($0$-adic) formal schemes.
The following proposition follows immediately from \ref{prop-sepmorformal0}, the definition of `of finite type' (\ref{dfn-topfintype}), and \ref{prop-closedstablezar}:
\begin{prop}\label{prop-properstablezar}
Let $f\colon Y\rightarrow X$ be a morphism of adic formal schemes of finite ideal type, and $\{V_{\alpha}\}_{\alpha\in L}$ an open covering of $X$.
Then $f$ is proper if and only if for any $\alpha\in L$ the base change $f^{-1}(V_{\alpha})\rightarrow V_{\alpha}$ is proper.\hfill$\square$
\end{prop}

By \ref{prop-sepmorformal1}, \ref{prop-cortopfintype11}, and \ref{prop-univclosedadicred} we have:
\begin{prop}\label{prop-propermorformal0}
Let $f\colon X\rightarrow Y$ be an adic morphism of adic formal schemes of finite ideal type, and suppose $Y$ has an ideal of definition $\mathscr{I}$ of finite type.
Set $X_k=(X,\O_X/\mathscr{I}^{k+1}\O_X)$ and $Y_k=(Y,\O_Y/\mathscr{I}^{k+1})$ for $k\geq 0$. 
Let $f_k\colon X_k\rightarrow Y_k$ be the induced morphism of scheme.
Then the following conditions are equivalent$:$
\begin{itemize}
\item[{\rm (a)}] $f$ is proper$;$
\item[{\rm (b)}] $f_k$ is proper for any $k\geq 0;$
\item[{\rm (c)}] $f_0$ is proper. \hfill$\square$
\end{itemize}
\end{prop}

\begin{prop}\label{prop-propermorformal10}
A finite morphism\index{morphism of formal schemes@morphism (of formal schemes)!finite morphism of formal schemes@finite ---} between adic formal schemes of finite ideal type is proper.
\end{prop}

\begin{proof}
To check if an adic morphism $f\colon X\rightarrow Y$ of adic formal schemes is proper, by \ref{prop-properstablezar} we may assume that $Y$ has an ideal of definition of finite type.
Then the assertion follows immediately from \ref{prop-finitemorphism1} and \ref{prop-propermorformal0}.
\end{proof}

\begin{prop}\label{prop-propermorformal2}
{\rm (1)} If $f\colon Z\rightarrow Y$ and $g\colon Y\rightarrow X$ are proper, then so is the composition $g\circ f$.

{\rm (2)} For any proper $S$-morphisms $f\colon X\rightarrow Y$ and $g\colon X'\rightarrow Y'$ of adic formal schemes of finite ideal type over an adic formal scheme $S$ of finite ideal type, the induced morphism $f\times_Sg\colon X\times_SY\rightarrow X'\times_SY'$ is proper.

{\rm (3)} For any proper $S$-morphism $f\colon X\rightarrow Y$ of adic formal schemes of finite ideal type over an adic formal scheme $S$ of finite ideal type and for any morphism $S'\rightarrow S$ of adic formal schemes of finite ideal type, the induced morphism $f_{S'}\colon X\times_SS'\rightarrow Y\times_SS'$ is proper.

{\rm (4)} Suppose the composition $g\circ f$ of two adic morphisms of adic formal schemes of finite ideal type is proper.
If $g$ is separated, $f$ is proper.
If $g$ is separated of finite type and $f$ is surjective, then $g$ is proper.
\end{prop}

\begin{proof}
(1), (2), and (3) follow from \ref{prop-sepmorformal24}, \ref{prop-sepmorformal2}, \ref{prop-topfintype2}, and \ref{prop-univclosedstaility}.
The first assertion of (4) follows from \ref{prop-propertyP3} and \ref{prop-propermorformal10}.
For the other assertion, it suffices to show that $g$ is closed.
But this is easy to see, since $f$ is surjective.
\end{proof}

In view of \ref{prop-genpropertymorphismadic2} we have the following due to \ref{prop-propermorformal0} and \cite[$\mathbf{II}$, (5.4.2)]{EGA}:
\begin{cor}\label{cor-propermorformal00}
Let $f\colon X\rightarrow Y$ be a morphism of schemes, and $Z$ a closed subscheme of $Y$ of finite presentation.
If $f$ is proper, then the formal completion $\widehat{f}\colon\widehat{X}|_{f^{-1}(Z)}\rightarrow\widehat{Y}_Z$ is proper.
\end{cor}
\index{morphism of formal schemes@morphism (of formal schemes)!proper morphism of formal schemes@proper ---|)}

\subsection{Flat and faithfully flat morphisms}\label{sub-flatformal}
\index{morphism of formal schemes@morphism (of formal schemes)!flat morphism of formal schemes@flat ---|(}
\index{morphism of formal schemes@morphism (of formal schemes)!flat morphism of formal schemes@flat ---!faithfully flat morphism of formal schemes@faithfully --- ---|(}
Let $f\colon X\rightarrow Y$ be a morphism of formal schemes, and $\mathscr{F}$ an $\O_X$-module.
As usual (cf.\ {\bf \ref{ch-pre}}.\ref{dfn-flatnessringedsp}), we say that $\mathscr{F}$ is $f$-flat (or $Y$-flat) at $x\in X$ if $\mathscr{F}_x$ is flat\index{flatness} as a module over $\O_{Y,f(x)}$; likewise, by saying $\mathscr{F}$ is $f$-flat we mean that $\mathscr{F}$ is flat at every point of $X$.
Similarly, a morphism of formal schemes $f\colon X\rightarrow Y$ is said to be flat if it is flat as a morphism of locally ringed spaces\index{morphism of ringed spaces@morphism (of ringed spaces)!flat morphism of ringed spaces@flat ---}, that is, for any $x\in X$ the induced morphism $\O_{Y,f(x)}\rightarrow\O_{X,x}$ is flat.

\subsubsection{First properties of flatness}\label{subsub-flatformalflat}
\begin{prop}\label{prop-flatformal3}
Let $A\rightarrow B$ be an adic map of t.u.\ rigid-Noetherian rings\index{t.u. rigid-Noetherian ring@t.u.\ rigid-Noetherian ring}, and $I\subseteq A$ a finitely generated ideal of definition.
Consider the associated adic morphism\index{morphism of formal schemes@morphism (of formal schemes)!adic morphism of formal schemes@adic ---} $f\colon X=\Spf B\rightarrow Y=\Spf A$. 
Let $M$ be a finitely generated $A$-module, and $\mathscr{F}=M^{\Delta}$ the associated adically quasi-coherent sheaf on $X$ of finite type\index{adically quasi-coherent (a.q.c.) sheaf!adically quasi-coherent sheaf of finite type@--- of finite type}.
For $k\geq 0$ we set $X_k=\Spec B_k$ and $Y_k=\Spec A_k$, where $B_k=B/I^{k+1}B$ and $A_k=A/I^{k+1}$, and denote by $f_k\colon X_k\rightarrow Y_k$ the induced morphism of schemes$;$ we likewise set $\mathscr{F}_k=\mathscr{F}/I^{k+1}\mathscr{F}$ for $k\geq 0$.
Then the following conditions are equivalent$:$
\begin{itemize}
\item[{\rm (a)}] $\mathscr{F}$ is $f$-flat$;$
\item[{\rm (b)}] $\mathscr{F}_k$ is $f_k$-flat for any $k\geq 0;$
\item[{\rm (c)}] $M$ is $A$-flat.
\end{itemize}
\end{prop}

\begin{proof}
The implication (a) $\Rightarrow$ (b) is easy; indeed, $f_k$ is the base change of $f$ by the closed immersion $Y_k\rightarrow Y$ defined by $I^{k+1}$.
If (b) holds, then $M_k=M/I^{k+1}M$ is $A_k$-flat for any $k\geq 0$.
Then (c) follows by {\bf \ref{ch-pre}}.\ref{prop-flatmorformal12cor}.

Let us show (c) $\Rightarrow$ (a).
Let $M$ be $A$-flat.
Then $M_g$ for any $g\in A$ is $A_g$-flat (where we denote the image of $g$ in $B$ again by $g$).
But then by {\bf \ref{ch-pre}}.\ref{prop-flatmorformal12cor} the completed localization $M_{\{g\}}$ is flat over $A_{\{g\}}$, since $A_{\{g\}}$ is t.u.\ rigid-Noetherian.
Similarly, whenever there exists a canonical morphism $B_g\rightarrow B_h$, $M_{\{h\}}$ is flat over $A_{\{g\}}$; indeed, we have $M_{\{h\}}=M_{\{g\}}\widehat{\otimes}_{B_{\{g\}}}B_{\{h\}}=M_{\{g\}}\otimes_{B_{\{g\}}}B_{\{h\}}$; on the other hand, since $B_g\rightarrow B_h$ is flat, $B_{\{g\}}\rightarrow B_{\{h\}}$ is flat by {\bf \ref{ch-pre}}.\ref{prop-flatmorformal12cor}; thus $M_{\{h\}}$ is flat over $A_{\{g\}}$.

Now let $y=\mathfrak{q}$ be a point of $Y$ ($\mathfrak{q}$ is an open prime ideal of $A$), and set $x=f(y)=\mathfrak{p}$.
For any $g\not\in\mathfrak{p}$, $\varinjlim_{h\not\in\mathfrak{q}}M_{\{h\}}=\mathscr{F}_y$ is flat over $A_{\{g\}}$.
Varying $g$ and taking the inductive limit, we deduce that $\mathscr{F}_y$ is flat over $\O_{X,x}=\varinjlim_{g\not\in\mathfrak{p}}A_{\{g\}}$, whence the assertion (a).
\end{proof}

\begin{cor}\label{cor-flatformal3bis}
Let $f\colon X\rightarrow Y$ be an adic morphism between locally universally rigid-Noetherian formal schemes\index{formal scheme!universally rigid-Noetherian formal scheme@universally rigid-Noetherian ---!locally universally rigid-Noetherian formal scheme@locally --- ---}, and $\mathscr{F}$ an adically quasi-coherent sheaf on $X$ of finite type\index{adically quasi-coherent (a.q.c.) sheaf!adically quasi-coherent sheaf of finite type@--- of finite type}.
Suppose $Y$ has an ideal of definition $\mathscr{I}$ of finite type, and set $X_k=(X,\O_X/\mathscr{I}^{k+1}\O_X)$, $Y_k=(Y,\O_Y/\mathscr{I}^{k+1})$, and $\mathscr{F}_k=\mathscr{F}/\mathscr{I}^{k+1}\mathscr{F}$ for any $k\geq 0;$ let $f_k\colon X_k\rightarrow Y_k$ be the induced morphism of schemes.
Then the following conditions are equivalent$:$
\begin{itemize}
\item[{\rm (a)}] $\mathscr{F}$ is flat over $Y;$
\item[{\rm (b)}] $\mathscr{F}_k$ is flat over $Y_k$ for any $k\geq 0$. \hfill$\square$
\end{itemize}
\end{cor}

\begin{cor}\label{cor-flatformal21}
Let $f\colon X\rightarrow Y$ be adic morphism between locally universally rigid-Noetherian formal schemes\index{formal scheme!universally rigid-Noetherian formal scheme@universally rigid-Noetherian ---!locally universally rigid-Noetherian formal scheme@locally --- ---}, and $\mathscr{F}$ an adically quasi-coherent sheaf of finite type on $X$\index{adically quasi-coherent (a.q.c.) sheaf!adically quasi-coherent sheaf of finite type@--- of finite type}.
Then the following conditions are equivalent$:$
\begin{itemize}
\item[{\rm (a)}] $\mathscr{F}$ is $f$-flat$;$
\item[{\rm (b)}] the functor $\mathscr{G}\mapsto\mathscr{F}\otimes_{f^{-1}\O_Y}f^{-1}\mathscr{G}$ from the category of adically quasi-coherent sheaves on $Y$ of finite type to the category of adically quasi-coherent sheaves on $X$ of finite type $($cf.\ {\rm \ref{cor-adicqcoh71}}$)$ is exact. \hfill$\square$
\end{itemize}
\end{cor}

\begin{cor}\label{cor-flatformal}
{\rm (1)} Let $f\colon X\rightarrow Y$ and $g\colon X'\rightarrow Y'$ be flat $S$-morphisms of locally universally rigid-Noetherian formal schemes that are adic over an adic formal scheme $S$ of finite ideal type.
Suppose $X\times_SY$ and $X'\times_SY'$ are locally universally rigid-Noetherian.
Then the induced morphism $f\times_Sg\colon X\times_SY\rightarrow X'\times_SY'$ is flat.

{\rm (2)} Let $f\colon X\rightarrow Y$ be a flat $S$-morphism of locally universally rigid-Noetherian formal schemes that are adic over an adic formal scheme $S$ of finite ideal type, and $S'\rightarrow S$ an adic morphism of adic formal schems of finite ideal type.
Suppose that $X\times_SS'$ and $Y\times_SS'$ are locally universally rigid-Noetherian.
Then the induced morphism $f_{S'}\colon X\times_SS'\rightarrow Y\times_SS'$ is flat.
\end{cor}

\begin{proof}
Since the the statement is local on $S$, we may assume that $S$ has an ideal of definition of finite type.
Then the corollary follows from \ref{prop-flatformal3}, \ref{cor-prodformal21}, and \cite[$\mathbf{IV}$, (2.1.7)]{EGA}.
\end{proof}

\subsubsection{Faithfully flat morphisms}\label{subsub-flatformalfflat}
\begin{dfn}\label{dfn-ffformal1}{\rm 
Let $f\colon X\rightarrow Y$ be an adic morphism between locally universally rigid-Noetherian formal schemes\index{formal scheme!universally rigid-Noetherian formal scheme@universally rigid-Noetherian ---!locally universally rigid-Noetherian formal scheme@locally --- ---}, and $\mathscr{F}$ an adically quasi-coherent sheaf on $X$ of finite type\index{adically quasi-coherent (a.q.c.) sheaf!adically quasi-coherent sheaf of finite type@--- of finite type}.

(1) We say $\mathscr{F}$ is {\em faithfully flat}\index{flatness!faithfully-flatness@faithfully-{---}} over $Y$ if the functor $\mathscr{G}\mapsto\mathscr{F}\otimes_{f^{-1}\O_Y}f^{-1}\mathscr{G}$ from the category of adically quasi-coherent sheaves on $Y$ of finite type to the category of adically quasi-coherent sheaves on $X$ of finite type is exact and faithful (cf.\ \ref{cor-flatformal21}).

(2) If $\O_X$ is faithfully flat over $Y$, then we say the morphism $f$ is {\em faithfully flat} or $X$ is {\em faithfully flat} over $Y$.}
\end{dfn}

The condition in (1) is equivalent to that $\mathscr{F}$ is $f$-flat and for an adically quasi-coherent sheaf $\mathscr{G}$ of finite type on $Y$, $\mathscr{F}\otimes_{f^{-1}\O_Y}f^{-1}\mathscr{G}=0$ implies $\mathscr{G}=0$.

\begin{prop}\label{prop-ffsupport}
Let $f\colon X\rightarrow Y$ be an adic morphism between locally universally rigid-Noetherian formal schemes, and $\mathscr{F}$ an adically quasi-coherent sheaf on $X$ of finite type.
Then $\mathscr{F}$ is faithfully flat over $Y$ if and only if $f$-flat and $f(\Supp(\mathscr{F}))=Y$.
\end{prop}

To show the proposition, we first need to show:
\begin{lem}\label{lem-exerflatmodulesupport}
Let $f\colon X\rightarrow Y$ be an adic morphism of locally universally rigid-Noetherian formal schemes, and suppose that $Y$ has an ideal of definition of finite type $\mathscr{I}$.
Let $f_k\colon X_k=(X,\O_X/\mathscr{I}^{k+1}\O_X)\rightarrow Y_k=(Y,\O_Y/\mathscr{I}^{k+1})$ be the induced morphism of schemes for $k\geq 0$.
Then for any adically quasi-coherent sheaf $\mathscr{F}$ of finite type on $X$ and $k\geq 0$, $f(\Supp(\mathscr{F}))=f_k(\Supp(\mathscr{F}_k))$ holds, where $\mathscr{F}_k=\mathscr{F}/\mathscr{I}^{k+1}\mathscr{F}$.
\end{lem}

\begin{proof}
We first show the claim in the affine case $X=\Spf A$ and $Y=\Spf B$, where $A,B$ are t.u.\ rigid-Noetherian rings, and $I\subseteq B$ is a finitely generated ideal of definition.
In this case, $\mathscr{F}=M^{\Delta}$ by a finitely generated $A$-module $M$.
For an open prime ideal $\mathfrak{q}$ of $B$, $\mathfrak{q}\in f(\Supp(\mathscr{F}))$ if and only if $M_{\{g\}}=M\otimes_AA_{\{g\}}\neq 0$ for some $g\in B\setminus\mathfrak{q}$.
By Nakayama's lemma this is equivalent to that $M_g/I^{k+1}M_g\neq 0$, whence the claim in this case.
The general case reduces to the affine case, since $y\in Y$ lies in $f(\Supp(\mathscr{F}))$ if and only if there exists $x\in\Supp(\mathscr{F})$ such that $f(x)=y$.
\end{proof}

\begin{proof}[Proof of Proposition {\rm \ref{prop-ffsupport}}]
We may assume that $Y$ has an ideal of definition of finite type $\mathscr{I}$.
Let $f_k\colon X_k=(X,\O_X/\mathscr{I}^{k+1}\O_X)\rightarrow Y_k=(Y,\O_Y/\mathscr{I}^{k+1})$ be the induced morphism of schemes for $k\geq 0$, and set $\mathscr{F}_k=\mathscr{F}/\mathscr{I}^{k+1}\mathscr{F}$.
By \ref{lem-exerflatmodulesupport} we have $f(\Supp(\mathscr{F}))=f_k(\Supp(\mathscr{F}_k))$ for any $k\geq 0$.
Suppose $\mathscr{F}$ is faithfully flat over $Y$.
Then $\mathscr{F}$ is $f$-flat due to \ref{cor-flatformal21}.
To show $f(\Supp(\mathscr{F}))=Y$, it suffices to show that $f(\Supp(\mathscr{F}_k))=Y$.
By \ref{cor-flatformal3bis} the quasi-coherent $\O_{X_k}$-module $\mathscr{F}_k$ is $f_k$-flat.
If $f(\Supp(\mathscr{F}_k))\neq Y$ for some $k$, there exists a point $y\in Y$ such that $\mathscr{F}_k\otimes_{\O_{Y_k}}k(y)=0$, where $k(y)$ is the residue field of $y$ (\cite[$\mathbf{IV}$, (2.2.1)]{EGA}).
This implies, in particular, that $\mathscr{F}\otimes_{\O_Y}k(y)=0$, which contradicts the assumption, since $k(y)$ can be regarded as an adically quasi-coherent sheaf on $Y$ of finite type in an obvious way.

Conversely, suppose $\mathscr{F}$ is $f$-flat and $f(\Supp(\mathscr{F}))=Y$.
By \ref{cor-flatformal3bis} we know that $\mathscr{F}_k$ is $f_k$-flat for any $k\geq 0$.
Suppose $\mathscr{F}\otimes_{\O_Y}\mathscr{G}=0$ for an adically quasi-coherent sheaf $\mathscr{G}$ on $Y$.
Consider for $k\geq 0$ the equalities $(\mathscr{F}\otimes_{\O_Y}\mathscr{G})_k=\mathscr{F}_k\otimes_{\O_{Y_k}}\mathscr{G}_k=0$ (where $\mathscr{G}_k=\mathscr{G}/\mathscr{I}^{k+1}\mathscr{G}$).
Since $f_k(\Supp(\mathscr{F}_k))=Y_k$, $\mathscr{F}_k$ is faithfully flat over $Y_k$ (\cite[$\mathbf{IV}$, (2.2.6)]{EGA}), and hence we have $\mathscr{G}_k=0$.
Since $\mathscr{G}=\varprojlim_{k\geq 0}\mathscr{G}_k$, we have $\mathscr{G}=0$, as desired.
\end{proof}

\begin{cor}\label{cor-ffformal3}
Let $f\colon X\rightarrow Y$ be an adic morphism between locally universally rigid-Noetherian formal schemes.
Then $f$ is faithfully flat if and only if $f$ is flat and is surjective. \hfill$\square$
\end{cor}

\begin{cor}\label{cor-ffformal33}
Let $X$ be an locally universally rigid-Noetherian formal scheme, and $\ovl{X}\rightarrow X$ a Zariski covering, that is, $\ovl{X}=\coprod_{\alpha\in L}U_{\alpha}\rightarrow X$ by an open covering $\{U_{\alpha}\}_{\alpha\in L}$ of $X$.
Then the map $\ovl{X}\rightarrow X$ is faithfully flat. \hfill$\square$
\end{cor}

\begin{cor}[Local Criterion of Flatness]\label{cor-ffformal2}
Let $f\colon X\rightarrow Y$ be an adic morphism between locally universally rigid-Noetherian formal schemes, and $\mathscr{F}$ an adically quasi-coherent sheaf on $X$ of finite type.
Suppose $Y$ has an ideal of definition $\mathscr{I}$ of finite type, and set $X_k=(X,\O_X/\mathscr{I}^{k+1}\O_X)$, $Y_k=(Y,\O_Y/\mathscr{I}^{k+1})$, and $\mathscr{F}_k=\mathscr{F}/\mathscr{I}^{k+1}\mathscr{F}$ for any $k\geq 0;$ let $f_k\colon X_k\rightarrow Y_k$ be the induced morphism of schemes.
Then the following conditions are equivalent$:$
\begin{itemize}
\item[{\rm (a)}] $\mathscr{F}$ is faithfully flat over $Y;$
\item[{\rm (b)}] $\mathscr{F}_k$ is faithfully flat over $Y_k$ for any $k\geq 0$.
\end{itemize}
If, moreover, $X$ and $Y$ are affine $X=\Spf B$ and $Y=\Spf A$, and if $\mathscr{F}=M^{\Delta}$ by a finitely generated $B$-module $M$, then the conditions are equivalent to$:$
\begin{itemize}
\item[{\rm (c)}] $M$ is faithfully flat over $A$. 
\end{itemize}
\end{cor}

\begin{proof}
The equivalence of (a) and (b) follows immediately from \ref{cor-flatformal3bis}, \ref{prop-ffsupport} and \cite[$\mathbf{II}$, (2.2.6)]{EGA}.
The last assertion follows from {\bf \ref{ch-pre}}.\ref{cor-flatmorformal11}.
\end{proof}

\begin{cor}\label{cor-ffformal4}
{\rm (1)} Let $f\colon X\rightarrow Y$ and $g\colon X'\rightarrow Y'$ be faithfully flat $S$-morphisms of locally universally rigid-Noetherian formal schemes that are adic over an adic formal scheme $S$ of finite ideal type.
Suppose $X\times_SY$ and $X'\times_SY'$ are locally universally rigid-Noetherian.
Then the induced morphism $f\times_Sg\colon X\times_SY\rightarrow X'\times_SY'$ is faithfully flat.

{\rm (2)} Let $f\colon X\rightarrow Y$ be a faithfully flat $S$-morphism of locally universally rigid-Noetherian formal schemes that are adic over an adic formal scheme $S$ of finite ideal type, and $S'\rightarrow S$ an adic morphism of adic formal schemes of finite ideal type.
Suppose that $X\times_SS'$ and $Y\times_SS'$ are locally universally rigid-Noetherian.
Then the induced morphism $f_{S'}\colon X\times_SS'\rightarrow Y\times_SS'$ is faithfully flat.

{\rm (3)} Let $f\colon X\rightarrow Y$ and $g\colon Y\rightarrow Z$ two adic morphisms of locally universally rigid-Noetherian formal schemes, and suppose $f$ is faithfully flat.
Then for $g$ to be flat $($resp.\ faithfully flat$)$ it is necessary and sufficient that $g\circ f$ is flat $($resp.\ faithfully flat$)$.
\end{cor}

\begin{proof}
We may assume that $S$ has an ideal of definition.
Then the corollary follows from \ref{cor-ffformal2}, \ref{cor-prodformal21}, and \cite[$\mathbf{IV}$, (2.1.7)]{EGA}.
\end{proof}

\subsubsection{Adically flat morphisms}\label{subsub-adicallyflat}
\index{morphism of formal schemes@morphism (of formal schemes)!adically flat morphism of formal schemes@adically flat ---|(}
\index{morphism of formal schemes@morphism (of formal schemes)!adically faithfully flat morphism of formal schemes@adically faithfully flat ---|(}
As is implicit in the above observations, we need to assume, whenever discussing flatness, that the formal schemes under consideration are locally universally rigid-Noetherian.
In more general situation, the weaker notion `adically flatness' defined as follows is more reasonable.
\begin{dfn}\label{dfn-adicallyflat}{\rm 
(1) Let $f\colon X\rightarrow Y$ be an adic morphism\index{morphism of formal schemes@morphism (of formal schemes)!adic morphism of formal schemes@adic ---} of adic formal schemes of finite ideal type, and $x\in X$ a point.
Then $f$ is said to be {\em adically flat} at $x$ if the following condition is satisfied: 
there exists an affine open neighborhood $V$ of $y=f(x)$ and an ideal of definition $\mathscr{I}$ of finite type of $V$ such that for any $k\geq 0$ the induced morphism of schemes $U_k\rightarrow V_k$ is flat at $x$, where $U=f^{-1}(V)$, $V_k$ is the closed subscheme of $V$ defined by $\mathscr{I}^{k+1}$, and $U_k=V_k\times_VU$.

(2) An adic morphism $f\colon X\rightarrow Y$ between adic formal schemes is said to be {\em adically flat} if $f$ is adically flat at all points of $X$ (cf.\ \ref{exas-adicalizationexamples}).
If $f$ is, moreover, surjective, we say that $f$ is {\em adically faithfully flat}.}
\end{dfn}

It is clear that in (1) the definition does not depend on the choice of $\mathscr{I}$.
Hence we readily have:
\begin{prop}\label{prop-adicallyflatformal00}
Let $f\colon X\rightarrow Y$ be an adic morphism of adic formal schemes of finite ideal type.
Suppose that $Y$ has an ideal of definition $\mathscr{I}$ of finite type.
Set $X_k=(X,\O_X/\mathscr{I}^{k+1}\O_X)$ and $Y_k=(Y,\O_Y/\mathscr{I}^{k+1})$ for $k\geq 0$, and denote by $f_k\colon X_k\rightarrow Y_k$ the induced morphism of schemes.
Then the following conditions are equivalent$:$
\begin{itemize}
\item[{\rm (a)}] $f$ is adically flat$;$
\item[{\rm (b)}] $f_k$ is flat for $k\geq 0$. \hfill$\square$
\end{itemize}
\end{prop}

\begin{cor}\label{cor-adicallyflatformal002}
Let $f\colon X\rightarrow Y$ be a morphism of schemes, and $Z$ a closed subscheme of $Y$ of finite presentation.
If $f$ is flat, then the formal completion $\widehat{f}\colon\widehat{X}|_{f^{-1}(Z)}\rightarrow\widehat{Y}|_Z$ is adically flat. \hfill$\square$
\end{cor}

In particular, by \ref{cor-flatformal3bis} and \ref{cor-ffformal3} we have:
\begin{prop}\label{prop-adicallyflatformal01}
An adically flat morphism between locally universally rigid-Noetherian formal schemes\index{formal scheme!universally rigid-Noetherian formal scheme@universally rigid-Noetherian ---!locally universally rigid-Noetherian formal scheme@locally --- ---} is flat. 
In particular, an adic morphism between locally universally rigid-Noetherian formal schemes is faithfully flat if and only if it is adically faithfully flat. \hfill$\square$
\end{prop}

\begin{prop}\label{prop-adicallyflatformal1}
{\rm (1)} The composition of two adically flat morphisms is adically flat.
Let $f\colon X\rightarrow Y$ and $g\colon Y\rightarrow Z$ be adic morphisms of adic formal schemes of finite ideal type.
Suppose $f$ is adically faithfully flat.
Then if $g\circ f$ is adically flat, so is $g$.

{\rm (2)} For any adically flat $S$-morphisms $f\colon X\rightarrow Y$ and $g\colon X'\rightarrow Y'$ of adic formal schemes of finite ideal type over an adic formal scheme $S$ of finite ideal type, the induced morphism $f\times_Sg\colon X\times_SY\rightarrow X'\times_SY'$ is adically flat.

{\rm (3)} For any adically flat $S$-morphism $f\colon X\rightarrow Y$ of adic formal schemes of finite ideal type over an adic formal scheme $S$ of finite ideal type and for any morphism $S'\rightarrow S$ of adic formal schemes of finite ideal type, the induced morphism $f_{S'}\colon X\times_SS'\rightarrow Y\times_SS'$ is adically flat. \hfill$\square$
\end{prop}

The proof is easy, and is left to the reader.
By {\bf \ref{ch-pre}}.\ref{prop-localcriterionflatuseful} one easily deduces the following useful fact:
\begin{prop}\label{prop-adicallyflatlocalcriterian}
Let $f\colon X\rightarrow Y$ and $g\colon Y\rightarrow Z$ be adic morphisms\index{morphism of formal schemes@morphism (of formal schemes)!adic morphism of formal schemes@adic ---} of adic formal schemes of finite ideal type.
Suppose that $Z$ has an ideal of definition $\mathscr{I}$ of finite type, and set $X_k=(X,\O_X/\mathscr{I}^{k+1}\O_X)$, $Y_k=(Y,\O_Y/\mathscr{I}^{k+1}\O_Y)$, and $Z_k=(Z,\O_Z/\mathscr{I}^{k+1})$ for $k\geq 0$.
Suppose, moreover, that $g$ is adically flat.
Then the following conditions are equivalent$:$
\begin{itemize}
\item[{\rm (a)}] $f$ is adically flat$;$
\item[{\rm (b)}] $g\circ f$ is adically flat, and $f_0\colon X_0\rightarrow Y_0$ is flat. \hfill$\square$
\end{itemize}
\end{prop}

\begin{prop}\label{prop-adicnessadicallyflatdescent}
Suppose that we have the commutative diagram
$$
\xymatrix{Z\ar[d]_f\ar[dr]^{g\circ f}\\ Y\ar[r]_g&X}
$$
consisting of adic formal schemes of finite ideal type and that $f$ is adically faithfully flat and quasi-compact.
Then $g$ is adic if and only if $g\circ f$ is adic\index{morphism of formal schemes@morphism (of formal schemes)!adic morphism of formal schemes@adic ---}.
\end{prop}

\begin{proof}
The `only if' part is clear due to \ref{prop-adicmor1} (1).
To show the `if' part, we may assume that $X$ and $Y$ are affine $X=\Spf A$ and $Y=\Spf B$, where $A$ and $B$ are an adic rings of finite ideal type.
Since $f$ is quasi-compact, one can replace $Z$ by the disjoint union of affine open subspaces, and thus we may assume that $Z$ is affine $Z=\Spf R$ where $R$ is an adic ring of finite ideal type.
Hence we have the following diagram
$$
\xymatrix{&R\\ A\ar[ur]\ar[r]&B\ar[u]\rlap{,}}
$$
where we assume that the map $A\rightarrow R$ is adic.
We need to show that $A\rightarrow B$ is adic.
Take a finitely generated ideal of definition $I\subseteq A$ (resp.\ $J\subseteq B$) of $A$ (resp.\ $B$).
Replacing $I$ by a suitable power, we may assume that $J^mR\subseteq IR\subseteq J^nR$ for some $m,n\geq 1$.
Since, by the assumption, the maps $B/J^{k+1}\rightarrow R/J^{k+1}R$ are faithfully flat for all $k\geq 0$, we have
$$
J^m\subseteq IB+J^{k+1}\subseteq J^n
$$
for any sufficiently large $k$.
This shows that the closure $\ovl{IB}$ of $IB$ in the $J$-adic ring $B$ is an ideal of definition of $B$.
Since $\ovl{IB}/J^{k+1}\cong(IB+J^{k+1})/J^{k+1}$, the ideal $\ovl{IB}$ is finitely generated by {\bf \ref{ch-pre}}.\ref{prop-complpair1}.
Since $IB$ is dense in $\ovl{IB}$ and $\ovl{IB}$ is open, we have $\ovl{IB}=IB+\ovl{IB}^2$.
By Nakayama's lemma we have $IB=\ovl{IB}$, and hence $IB$ is an ideal of definition of $B$.
This shows that $A\rightarrow B$ is adic, as desired.
\end{proof}
\index{morphism of formal schemes@morphism (of formal schemes)!adically faithfully flat morphism of formal schemes@adically faithfully flat ---|)}
\index{morphism of formal schemes@morphism (of formal schemes)!adically flat morphism of formal schemes@adically flat ---|)}
\index{morphism of formal schemes@morphism (of formal schemes)!flat morphism of formal schemes@flat ---!faithfully flat morphism of formal schemes@faithfully --- ---|)}
\index{morphism of formal schemes@morphism (of formal schemes)!flat morphism of formal schemes@flat ---|)}

\addcontentsline{toc}{subsection}{Exercises}
\subsection*{Exercises}
\begin{exer}\label{exer-affinemorphism2idealofdefinition}{\rm 
Let $f\colon X\rightarrow Y$ be an affine adic morphism of adic formal schemes of finite ideal type, and $\mathscr{F}$ (resp.\ $\mathscr{R}$) an adically quasi-coherent sheaf (resp.\ adically quasi-coherent $\O_X$-algebra). 
Suppose that there exists an ideal of definition $\mathscr{I}$ of $Y$ of finite type.
Then show that $f_{\ast}\mathscr{I}^{k+1}\mathscr{F}=\mathscr{I}^{k+1}f_{\ast}\mathscr{F}$ (resp.\ $f_{\ast}\mathscr{I}^{k+1}\mathscr{R}=\mathscr{I}^{k+1}f_{\ast}\mathscr{R}$) holds for $k\geq 0$.}
\end{exer}

\begin{exer}\label{exer-closedimmformal61}
{\rm Let $X=\Spf A$ be an affine universally rigid-Noetherian formal scheme, and $K\subseteq A$ an ideal.
Show that $\mathscr{K}=K^{\Delta}$ is an ideal of $\O_X$ and that $\O_X/\mathscr{K}$ is adically quasi-coherent.
In particular, $(Y,\O_X/\mathscr{K})$ where $Y$ is the support of $\O_X/\mathscr{K}$ is a closed formal subscheme of $X$.}
\end{exer}

\begin{exer}\label{exer-adicnessadicallyflatdescent}{\rm 
Let $A\rightarrow B$ be an adic map of adic rings of finite ideal type, and $I\subseteq A$ a finitely generated ideal of definition.
Suppose that $\Spf B\rightarrow\Spf A$ is adically faithfully flat.
Let $F^{\bullet}=\{F^n\}_{n\in\Z}$ be a descending filtration of $A$, separated\index{filtration by submodules@filtration (by submodules)!separated filtration by submodules@separated ---} and exhaustive\index{filtration by submodules@filtration (by submodules)!exhaustive filtration by submodules@exhaustive ---} (cf.\ {\bf \ref{ch-pre}}, \S\ref{subsub-topfromfil}), consisting of finitely generated ideals such that for any $q\geq 0$ and $n\in\Z$ we have $I^qF^n\subseteq F^{n+q}$ (cf.\ {\bf \ref{ch-pre}}, \S\ref{subsub-ARgeneral}).
Show that the following conditions are equivalent$:$
\begin{itemize}
\item[{\rm (a)}] the filtration $F^{\bullet}$ is $I$-good\index{I-good@$I$-good (filtration)}\index{filtration by submodules@filtration (by submodules)!I-good filtration by submodules@$I$-good ---} $({\bf \ref{ch-pre}}.\ref{dfn-Igood});$
\item[{\rm (b)}] the induced filtration $F^{\bullet}B=\{F^nB\}_{n\in\Z}$ on $B$ is $IB$-good.
\end{itemize}}
\end{exer}


\section{Differential calculus on formal schemes}\label{sec-pairdifferential}
This section aims at establishing some basics on differential calculi on formal schemes.
We define the related notions of morphisms of formal schemes, neat, \'etale, and smooth morphisms, and discuss fundamental properties of them.
In order to do this, we first develop in \S\ref{sub-pairdifferentialgen} the `continuous' version of the theory of derivations and differentials and discuss completions of the differential modules.
The theory of continuous derivations and differentials has already been discussed in \cite[$\mathbf{0}_{\mathbf{IV}}$, \S20]{EGA}, and most of our arguments here will be, therefore, brief rehashes of what have been done in this reference.

In \S\ref{sub-differentialinvformalsch} we define for a morphism $X\rightarrow Y$ of adic formal schemes of finite ideal type the sheaf $\Omega^1_{X/Y}$ of $1$-differentials as an adically quasi-coherent sheaf. 
Based on this, we then proceed to discuss \'etale and smooth morphisms of formal schemes in \S\ref{sub-etaleformal}.
Here the reader should be warned at the fact that, according to our definitions, smooth or even \'etale morphisms are only adically flat, but not necessarily flat.
Needless to say, this defect comes from the failure of local criterion of flatness (cf.\ {\bf \ref{ch-pre}}.\ref{prop-flatmorformal1}).
Hence, in particular, they are flat if the formal schemes under consideration are locally universally rigid-Noetherian\index{formal scheme!universally rigid-Noetherian formal scheme@universally rigid-Noetherian ---!locally universally rigid-Noetherian formal scheme@locally --- ---} (\ref{dfn-formalsch}).

\subsection{Differential calculi for topological rings}\label{sub-pairdifferentialgen}
Let us first recall some generalities of differential calculus for linearly topologized rings, which has already been developed in \cite[$\mathbf{0}_{\mathbf{IV}}$, \S20]{EGA}.
\subsubsection{Continuous derivations}\label{subsub-continuousderivations}
Let $A\rightarrow B$ be a homomorphism of rings.
Suppose that $B$ is endowed with a descending filtration by ideals\index{filtration by submodules@filtration (by submodules)} ({\bf \ref{ch-pre}}, \S\ref{subsub-topfromfil}) $J^{\bullet}=\{J^{\lambda}\}_{\lambda\in\Lambda}$, and consider the topology\index{filtration by submodules@filtration (by submodules)!topology defined by a filtration by submodules@topology defined by a ---}\index{topology!topology defined by a filtration@--- defined by a filtration} on $B$ defined by $J^{\bullet}$ ({\bf \ref{ch-pre}}, \S\ref{subsub-topologicalringsmodules}).
In this situation, one can similarly consider the topology on any $B$-module $M$ defined by the descending filtration by $B$-submodules $J^{\bullet}M=\{J^{\lambda}M\}_{\lambda\in\Lambda}$.

We denote by
$$
\Dercont_A(B,M)
$$
the set of all {\em continuous} $A$-derivations\index{derivation!continuous derivation@continuous ---} of $B$ with values in $M$, that is, continuous additive mappings $\delta\colon B\rightarrow M$ such that $\delta(xy)=x\delta(y)+y\delta(x)$ for $x,y\in B$ and that $\delta(a)=0$ for $a\in A$.
Notice that the continuity of $\delta$ is equivalent to that for any $\lambda\in\Lambda$ there exists $\mu\in\Lambda$ such that $\delta(J^{\mu})\subseteq J^{\lambda}M$ holds.
It is then easy to verify that the set $\Dercont_A(B,M)$ has the canonical $B$-module structure in such a way that it is a $B$-submodule of the $B$-module $\Der_A(B,M)$ of all $A$-derivations of $B$ with values in $M$.

\subsubsection{Differentials and canonical topology}\label{subsub-differentialscantop}
Consider the differential module\index{differential!module of differentials@module of ---s} $\Omega^1_{B/A}$ (without regarding the topologies) together with the $A$-derivation 
$$
d\colon B\longrightarrow\Omega^1_{B/A}.
$$
In order that a $B$-linear morphism $\varphi\colon\Omega^1_{B/A}\rightarrow M$ is such that the composition $\varphi\circ d$ is continuous, it is necessary and sufficient that for any $\lambda\in\Lambda$ there exist $\mu\in\Lambda$ such that $d(J^{\mu})\subseteq\varphi^{-1}(J^{\lambda}M)$.
Hence it is natural to consider the topology on the $B$-module $\Omega^1_{B/A}$ defined by the descending filtration by $B$-submodules $\langle d(J^{\bullet})\rangle_B=\{\langle d(J^{\lambda})\rangle_B\}_{\lambda\in\Lambda}$, where $\langle d(J^{\lambda})\rangle_B$ denotes the $B$-submodule of $\Omega^1_{B/A}$ generated by the image $d(J^{\lambda})$.
We call this topology on $\Omega^1_{B/A}$ the {\em canonical topology}\index{topology!canonical topology on differentials@canonical --- (on differentials)}.
Notice that the $A$-derivation $d\colon B\rightarrow\Omega^1_{B/A}$ is continuous.

\begin{prop}\label{prop-genpairdifferential1}
The map $\varphi\mapsto\varphi\circ d$ induces an isomorphism 
$$
\Homcont_B(\Omega^1_{B/A},M)\stackrel{\sim}{\longrightarrow}\Dercont_A(B,M)
$$
of $B$-modules, where the left-hand side is the $B$-module consisting of {\em continuous} $B$-linear maps and $\Omega^1_{B/A}$ is endowed with the canonical topology. \hfill$\square$
\end{prop}

\begin{prop}\label{prop-genpairdifferential2}
The canonical topology on $\Omega^1_{B/A}$ is coarser than the linear topology defined by $J^{\bullet}\Omega^1_{B/A}$.
These topologies coincide with each other if for any $\lambda\in\Lambda$ there exists $\mu\in\Lambda$ such that $J^{\mu}\subseteq(J^{\lambda})^2$.
\end{prop}

Notice that the hypothesis of the second assertion is satisfied if the topology on $B$ is adic\index{topology!adic topology@adic ---}\index{adic!adic topology@--- topology} ({\bf \ref{ch-pre}}, \S\ref{subsub-adicfiltrationtopology}).
\begin{proof}
For $a\in J^{\lambda}$ and $x\in B$ we have $ad(x)=d(ax)-xd(a)$, which belongs to $\langle d(J^{\lambda})\rangle_B$.
As $\Omega^1_{B/A}$ is generated by elements of the form $d(x)$ $(x\in B)$, this shows the inclusion $J^{\lambda}\Omega^1_{B/A}\subseteq\langle d(J^{\lambda})\rangle_B$, whence the first assertion.
As for the second, observe that $\langle d((J^{\lambda})^2)\rangle_B\subseteq J^{\lambda}\Omega^1_{B/A}$, which shows that each $J^{\lambda}\Omega^1_{B/A}$ is open with respect to the topology defined by the filtration $\{\langle d(J^{\lambda})\rangle_B\}_{\lambda\in\Lambda}$.
\end{proof}

\subsubsection{Completion and differentials}\label{subsub-diffcomplete}
Now we assume that $B$ is endowed with the adic topology\index{topology!adic topology@adic ---}\index{adic!adic topology@--- topology} by a finitely generated ideal $J\subseteq B$ ({\bf \ref{ch-pre}}, \S\ref{subsub-adicfiltrationtopology}).
In this situation, by \ref{prop-genpairdifferential2} the canonical topology on the differential module $\Omega^1_{B/A}$ coincides with the $J$-adic topology.
We denote the $J$-adic completion\index{completion!I-adic completion@$I$-adic ---} ({\bf \ref{ch-pre}}.\ref{dfn-Iadiccompletiondelicate1}, {\bf \ref{ch-pre}}.\ref{prop-Iadiccompletioncomplete1}) of $\Omega^1_{B/A}$ by $\widehat{\Omega}^1_{B/A}$, called the {\em complete differential module}\index{differential!complete differential module@complete --- module} of $B$ relative to $A$.
By the universality of $J$-adic completions we deduce from \ref{prop-genpairdifferential1} the following:
\begin{prop}\label{prop-genpairdifferential1compl}
For any $J$-adically complete $B$-module $M$ the canonical map $\varphi\mapsto\varphi\circ d$ gives rise to an isomorphism 
$$
\Homcont_B(\widehat{\Omega}^1_{B/A},M)\stackrel{\sim}{\longrightarrow}\Dercont_A(B,M)
$$
of $B$-modules. \hfill$\square$
\end{prop}

\begin{cor}\label{cor-genpairdifferential1compl}
Let $A,B$ be rings with adic topologies by finitely generated ideals, and $A\rightarrow B$ a continuous homomorphism.
We have the canonical isomorphisms
$$
\widehat{\Omega}^1_{B/A}\stackrel{\sim}{\longrightarrow}\widehat{\Omega}^1_{\widehat{B}/A}\stackrel{\sim}{\longrightarrow}\widehat{\Omega}^1_{\widehat{B}/\widehat{A}}.
$$
\end{cor}

\begin{proof}
In view of \ref{prop-genpairdifferential1compl} it suffices to show that the natural maps
$$
\Dercont_{\widehat{A}}(\widehat{B},M)\longhookrightarrow\Dercont_A(\widehat{B},M)\longrightarrow\Dercont_A(B,M)
$$
for any $J$-adically complete $B$-module $M$ are bijective.
The bijectivity of the first arrow follows immediately from the fact that $A\rightarrow B$ is continuous.
Since $\delta(J^{k+1})\subseteq J^kM$ for $k\geq 0$, any $A$-derivation $\delta\colon B\rightarrow M$ can be uniquely extended to an $A$-derivation from $\widehat{B}$, whence the bijectivity of the second arrow.
\end{proof}

Let $A$ (resp.\ $B$, resp.\ $C$) be a ring with the adic topology defined by a finitely generated ideal $I\subseteq A$ (resp.\ $J\subseteq B$, resp.\ $K\subseteq C$), and $A\rightarrow B$ and $A\rightarrow C$ continuous homomorphisms such that $IB\subseteq J$ and $IC\subseteq K$.
As we have seen in \S\ref{subsub-completetensorproducts}, $B\otimes_AC$ has the topology defined by $H^{\bullet\bullet}=\{H^{mn}\}_{m,n\geq 0}$, where
$$
H^{m,n}=\image(J^m\otimes_AC\rightarrow B\otimes_AC)+\image(B\otimes_AK^n\rightarrow B\otimes_AC)
$$
for $m,n\geq 0$. 
It is shown in (the proof of) \ref{lem-fiverprodformaladicadic} that this topology is actually $H$-adic, where $H=H^{1,1}$. 
We consider the $B\otimes_AC$-module $\Omega^1_{B/A}\otimes_AC$.
This is endowed with the topology defined similarly to that of $B\otimes_AC$, that is, the topology defined by $\{{H'}^{m,n}\}_{m,n\geq 0}$, where 
$$
{H'}^{mn}=\image(J^m\Omega^1_{B/A}\otimes_AC\rightarrow\Omega^1_{B/A}\otimes_AC)+\image(\Omega^1_{B/A}\otimes_AK^n\rightarrow\Omega^1_{B/A}\otimes_AC)
$$
for $m,n\geq 0$.
Clearly, this topology coincides with the one by $\{H^{m,n}(\Omega^1_{B/A}\otimes_AC)\}_{m,n\geq 0}$, since we have ${H'}^{m,n}=H^{m,n}(\Omega^1_{B/A}\otimes_AC)$ for $m,n\geq 0$.
Hence we now deduce that the topology on $\Omega^1_{B/A}\otimes_AC$ coincides with the $H$-adic topology.
Hence in view of \cite[$\mathbf{0}_{\mathbf{IV}}$, (20.5.5)]{EGA} we have:
\begin{prop}\label{prop-adicdifferentialcompatibility2}
We have the canonical isomorphism
$$
\Omega^1_{B/A}\widehat{\otimes}_AC\stackrel{\sim}{\longrightarrow}\widehat{\Omega}^1_{B\widehat{\otimes}_AC/C}
$$
of topological $B\widehat{\otimes}_AC$-modules. \hfill$\square$
\end{prop}

\begin{cor}\label{cor-adicdifferentialcompatibility2}
In the situation as in {\rm \ref{prop-adicdifferentialcompatibility2}}, we have
$$
\widehat{\Omega}^1_{B\widehat{\otimes}_AC/C}/H^{k+1}\widehat{\Omega}^1_{B\widehat{\otimes}_AC/C}\cong\Omega^1_{B_k/A_k}\otimes_{A_k}C_k
$$
for $k\geq 0$, where $A_k=A/I^{k+1}$, $B_k=B/J^{k+1}$, and $C_k=C/K^{k+1}$.
\end{cor}

\begin{proof}
By \ref{prop-adicdifferentialcompatibility2} we see that $H^{k+1}\cdot(\Omega^1_{B/A}\widehat{\otimes}_AC)$ is closed in $\Omega^1_{B/A}\widehat{\otimes}_AC$, and hence 
$$
\Omega^1_{B/A}\widehat{\otimes}_AC/H^{k+1}\cdot(\Omega^1_{B/A}\widehat{\otimes}_AC)\cong\Omega^1_{B/A}\otimes_AC/H^{k+1}\cdot(\Omega^1_{B/A}\otimes_AC),
$$
from which the desired equality follows.
\end{proof}

\begin{cor}\label{cor-adicdifferentialcompatibility1}
Let $A$ $($resp.\ $B)$ be a ring with the adic topology defined by a finitely generated ideal $I\subseteq A$ $($resp.\ $J\subseteq B)$, and $A\rightarrow B$ a continuous homomorphism such that $IB\subseteq J$.
Set $A_k=A/I^{k+1}$ and $B_k=B/J^{k+1}$ for $k\geq 0$.
Then we have the following equalities up to canonical isomorphisms for any $k\geq 0:$
$$
\Omega^1_{B/A}/J^{k+1}\Omega^1_{B/A}=\Omega^1_{B/A}\otimes_BB_k=\Omega^1_{B_k/A_k}.
$$
In particular, the complete differential module $\widehat{\Omega}^1_{B/A}$ is canonically isomorphic to the projective limit $\varprojlim_{k\geq 0}\Omega^1_{B_k/A_k}$.
\end{cor}

\begin{proof}
We apply \ref{cor-adicdifferentialcompatibility2} with $C=A$.
Since, in this case, we have $H=J$, we have the following equalities up to canonical isomorphisms
$$
\widehat{\Omega}^1_{\widehat{B}/A}/J^{k+1}\widehat{\Omega}^1_{\widehat{B}/A}=\widehat{\Omega}^1_{B/A}/J^{k+1}\widehat{\Omega}^1_{B/A}=\Omega^1_{B/A}/J^{k+1}\Omega^1_{B/A}=\Omega^1_{B_k/A_k},
$$
where the first equality is due to \ref{cor-genpairdifferential1compl}.
\end{proof}

\begin{prop}\label{prop-adicdifferentialcompatibility3}
Let $A\rightarrow B$ be a continuous homomorphism between adic rings of finite ideal type\index{adic!adic ring@--- ring!adic ring of finite ideal type@--- --- of finite ideal type} {\rm (\ref{dfn-admissibleringsadicrings}, \ref{dfn-admissibleringoffiniteidealtype})}, and $S\subseteq B$ a multiplicative subset.
Then we have the canonical isomorphism
$$
\widehat{\Omega}^1_{B/A}\widehat{\otimes}_BB\{S^{-1}\}\stackrel{\sim}{\longrightarrow}\widehat{\Omega}^1_{B\{S^{-1}\}/A}
$$
of topological $B\{S^{-1}\}$-modules $($cf.\ {\rm \S\ref{subsub-formalnotadmrings}} for the definition of $B\{S^{-1}\})$.
\end{prop}

\begin{proof}
Let $J\subseteq B$ be a finitely generated ideal of definition.
By \cite[$\mathbf{0}_{\mathbf{IV}}$, (20.5.9)]{EGA} we have $\Omega^1_{B/A}\otimes_BS^{-1}B\cong\Omega^1_{S^{-1}B/A}$.
Since the topology on both sides are $J$-adic, we obtain the isomorphism $\Omega^1_{B/A}\widehat{\otimes}_BS^{-1}B\cong\widehat{\Omega}^1_{S^{-1}B/A}$ between the $J$-adic completions.
The left-hand side is clearly isomorphic to $\widehat{\Omega}^1_{B/A}\widehat{\otimes}_BB\{S^{-1}\}$, while the right-hand side is isomorphic to $\widehat{\Omega}^1_{B\{S^{-1}\}/A}$ due to \ref{cor-genpairdifferential1compl}.
\end{proof}

\begin{prop}\label{prop-adicdifferentialcompatibility31}
Let $A$ and $B$ be adic rings of finite ideal type, $S\subseteq A$ a multiplicative subset, and $A\{S^{-1}\}\rightarrow B$ a continuous homomorphism.
Then we have the canonical isomorphism
$$
\widehat{\Omega}^1_{B/A}\stackrel{\sim}{\longrightarrow}\widehat{\Omega}^1_{B/A\{S^{-1}\}}.
$$
of topological $B$-modules.
\end{prop}

\begin{proof}
By \cite[$\mathbf{0}_{\mathbf{IV}}$, (20.7.17)]{EGA} (cf.\ {\bf \ref{ch-pre}}.\ref{prop-qconsistency1-1p}) the map in question is surjective and the image of $\widehat{\Omega}^1_{A\{S^{-1}\}/A}\widehat{\otimes}_{A\{S^{-1}\}}B$ is dense in the kernel.
But this is $0$, since $\widehat{\Omega}^1_{A\{S^{-1}\}/A}=\widehat{\Omega}^1_{S^{-1}A/A}=0$, where the first equality is due to \ref{cor-genpairdifferential1compl}.
\end{proof}

\subsubsection{Differentials and finiteness conditions}\label{subsub-difffiniteness}
By \ref{cor-genpairdifferential1compl} applied to $B=A\dl X_1,\ldots,X_n\dr$ (the completion of $A[X_1,\ldots,X_n]$), we have:
\begin{prop}\label{prop-difffiniteness1}
Let $A$ be an adic ring of finite ideal type, and consider $B=A\dl X_1,\ldots,X_n\dr$.
Then the complete differential module $\widehat{\Omega}^1_{B/A}$ is a free $B$-module of rank $n$ with the bases $dX_1,\ldots,dX_n$. \hfill$\square$
\end{prop}

\begin{cor}\label{cor-difffiniteness1}
Let $A$ be an adic ring, and $B$ a topologically finitely generated\index{finitely generated!topologically finitely generated@topologically ---} algebra over $A$ $($cf.\ {\rm {\bf \ref{ch-pre}}.\ref{dfn-topfinigen}}$)$.
Then $\widehat{\Omega}^1_{B/A}$ is a finitely generated $B$-module.
\end{cor}

\begin{proof}
Write $B=R/\mathfrak{a}$ where $R=A\dl X_1,\ldots,X_n\dr$ and $\mathfrak{a}\subseteq R$ is a closed ideal.
Then by \cite[$\mathbf{0}_{\mathbf{IV}}$, (20.7.8)]{EGA} $\Omega^1_{B/A}$ is, as a topological $B$-module, the quotient of $\Omega^1_{R/A}\otimes_RB$ by the image of $\mathfrak{a}$.
Hence $\widehat{\Omega}^1_{B/A}$ is the quotient of $\Omega^1_{R/A}\widehat{\otimes}_RB=\Omega^1_{R/A}\otimes_RB=\bigoplus^n_{i=1}B(dX_i\otimes 1)$ (\ref{prop-difffiniteness1}) by the closure of the image of $\mathfrak{a}$ (cf.\ {\bf \ref{ch-pre}}.\ref{prop-qconsistency1-1p}).
\end{proof}

\begin{thm}[Fundamental exact sequences]\label{thm-fundamentalexactseqdiff}
Let $A$ be a t.u.\ rigid-Noetherian ring\index{t.u. rigid-Noetherian ring@t.u.\ rigid-Noetherian ring} {\rm (\ref{dfn-tuaringadmissible} (1))}.

{\rm (1)} Let $B\rightarrow C$ be an $A$-algebra homomorphism between topologically finitely generated $A$-algebras.
Then we have the canonical exact sequence of $C$-modules
$$
\widehat{\Omega}^1_{B/A}\otimes_BC\longrightarrow\widehat{\Omega}^1_{C/A}\longrightarrow\widehat{\Omega}^1_{C/B}\longrightarrow 0.
$$

{\rm (2)} Let $B$ be a topologically finitely generated $A$-algebra, and $\mathfrak{a}$ a finitely generated ideal of $B$.
Set $C=B/\mathfrak{a}$.
Then we have the canonical exact sequence of $C$-modules
$$
\mathfrak{a}/\mathfrak{a}^2\longrightarrow\widehat{\Omega}^1_{B/A}\otimes_BC\longrightarrow\widehat{\Omega}^1_{C/A}\longrightarrow 0.
$$
\end{thm}

Notice that in (1) and (2), since $C$ is a t.u.\ rigid-Noetherian ring, and $\widehat{\Omega}^1_{B/A}$ is a finitely generated $B$-module, we have $\widehat{\Omega}^1_{B/A}\otimes_BC=\widehat{\Omega}^1_{B/A}\widehat{\otimes}_BC$ (\ref{rem-turigidnoetherianbasicproperties}).
\begin{proof}
(1) Let $I\subseteq A$ be a finitely generated ideal of definition, and for any $k\geq 0$ set $A_k=A/I^{k+1}$, $B_k=B/I^{k+1}B$, and $C_k=C/I^{k+1}C$; $B_k$ and $C_k$ are finite type $A_k$-algebras.
By \cite[$\mathbf{0}_{\mathbf{IV}}$, (20.5.7)]{EGA} we have the exact sequence
$$
\Omega^1_{B_k/A_k}\otimes_{B_k}C_k\longrightarrow\Omega^1_{C_k/A_k}\longrightarrow\Omega^1_{C_k/B_k}\longrightarrow 0.
$$
By {\bf \ref{ch-pre}}.\ref{cor-adicdifferentialcompatibility1} the projective systems $\{\Omega^1_{B_k/A_k}\otimes_{B_k}C_k\}_{k\geq 0}$, $\{\Omega^1_{C_k/A_k}\}_{k\geq 0}$, and $\{\Omega^1_{C_k/B_k}\}_{k\geq 0}$ are strict.
Hence, due to {\bf \ref{ch-pre}}.\ref{cor-adicdifferentialcompatibility1} and {\bf \ref{ch-pre}}.\ref{lem-ML3}, we have the desired exact sequence by applying $\varprojlim_{k\geq 0}$.

(2) Define $B_k$ and $C_k$ ($k\geq 0$) similarly as above.
The kernel of $B_k\rightarrow C_k$ is $\mathfrak{a}_k=\mathfrak{a}B_k=\mathfrak{a}/\mathfrak{a}\cap I^{k+1}$.
Hence we have the exact sequence (\cite[$\mathbf{0}_{\mathbf{IV}}$, (20.5.12)]{EGA})
$$
\mathfrak{a}_k/\mathfrak{a}^2_k\longrightarrow\widehat{\Omega}^1_{B_k/A_k}\otimes_{B_k}C_k\longrightarrow\widehat{\Omega}^1_{C_k/A_k}\longrightarrow 0.
$$
Since $\alpha$ is finitely generated, 
$$
0\longrightarrow\mathfrak{a}^2\longrightarrow\mathfrak{a}\longrightarrow\mathfrak{a}/\mathfrak{a}^2\longrightarrow 0
$$
is an exact sequence of $I$-adically complete modules.
Since $A$ is t.u.\ rigid-Noetherian\index{t.u. rigid-Noetherian ring@t.u.\ rigid-Noetherian ring}, we have $\mathfrak{a}=\varprojlim_{k\geq 0}\mathfrak{a}/I^{k+1}\mathfrak{a}=\varprojlim_{k\geq 0}\mathfrak{a}_k$ etc., and hence we have $\mathfrak{a}/\mathfrak{a}^2=\varprojlim_{k\geq 0}\mathfrak{a}_k/\mathfrak{a}^2_k$.
Then the desired exact sequence can be obtained by an argument similar to that in the proof of (1).
\end{proof}

\begin{cor}\label{cor-difffiniteness2}
Let $A$ be a t.u.\ rigid-Noetherian ring, and $B$ a topologically finitely presented algebra over $A$.
Then $\widehat{\Omega}^1_{B/A}$ is a finitely presented $B$-module.
More precisely, if $B=A\dl X_1,\ldots,X_n\dr/\mathfrak{a}$ where $\mathfrak{a}=(F_1,\ldots,F_m)$ is a finitely generated ideal, then 
$$
\widehat{\Omega}^1_{B/A}\cong{\textstyle \bigoplus^n_{i=1}BdX_i/\sum^m_{j=1}BdF_j}. \eqno{\square}
$$
\end{cor}

\subsection{Differential invariants on formal schemes}\label{sub-differentialinvformalsch}
\subsubsection{The sheaf of differentials}\label{subsub-differentialinvformalsch}
\begin{thm}\label{thm-differentialinvariants}
Let $f\colon X\rightarrow Y$ be a morphism between adic formal schemes of finite type.
Then there exists uniquely up to isomorphism an adically quasi-coherent sheaf\index{quasi-coherent!adically quasi-coherent OX module@adically --- (a.q.c.) sheaf} $\Omega^1_{X/Y}$ of finite type on $X$ such that for any affine open subset $V=\Spf A$ of $Y$ and any affine open subset $U=\Spf B$ of $f^{-1}(V)$ $($where $A$ and $B$ are adic rings of finite ideal type$;$ cf.\ {\rm \ref{cor-affineadicformalschemebyadicring}}$)$, we have $\Gamma(U,\Omega^1_{X/Y})\cong\widehat{\Omega}^1_{B/A}$ or, equivalently $($cf.\ {\rm \ref{thm-adicqcoh1}}$)$, $\Omega^1_{X/Y}|_U\cong(\widehat{\Omega}^1_{B/A})^{\Delta}$.
If $f$ is locally of finite type\index{morphism of formal schemes@morphism (of formal schemes)!morphism of formal schemes locally of finite type@--- locally of finite type}, then $\Omega^1_{X/Y}$ is adically quasi-coherent of finite type.
\end{thm}

\begin{proof}
The uniqueness is clear, and the existence follows from \ref{prop-adicdifferentialcompatibility3} and \ref{prop-adicdifferentialcompatibility31}.
The last assertion follows from \ref{cor-difffiniteness1} and \ref{cor-topqcoh111}.
\end{proof}

We call the adically quasi-coherent sheaf $\Omega^1_{X/Y}$ thus obtained the {\em $($complete$)$ sheaf of differentials over $X$ relative to $Y$}\index{differential!complete sheaf of differentials@(complete) sheaf of ---s}; this sheaf is equipped with a canonical $\O_Y$-module morphism 
$$
d\colon\O_X\longrightarrow\Omega^1_{X/Y},
$$
the so-called {\em canonical derivation}\index{derivation!canonical derivation@canonical ---}.

\begin{exa}\label{exa-differentialinvformalschfintype1}{\rm 
Let $Y$ be an adic formal scheme, and $\mathscr{E}$ a locally free sheaf on $Y$ of finite type.
Consider the vector bundle\index{vector bundle} $X=\widehat{\mathbf{V}}(\mathscr{E})$ associated to $\mathscr{E}$ (Exercise \ref{exer-vectorbundleformal}).
Then due to \ref{prop-difffiniteness1} $\Omega^1_{X/Y}$ is a locally free sheaf on $X$ of finite type.}
\end{exa}

By \ref{prop-adicdifferentialcompatibility2} and \ref{lem-completepullbackaqcsheaves1lem} we readily see: 
\begin{prop}\label{prop-differentialinvariants1}
Consider a Cartesian diagram
$$
\xymatrix{X'\ar[r]^g\ar[d]&X\ar[d]\\ Y'\ar[r]&Y}
$$
of adic formal schemes of finite ideal type, and suppose that either one of the following conditions are satisfied$:$
\begin{itemize}
\item[{\rm (a)}] $X$ and $X'$ are locally universally rigid-Noetherian\index{formal scheme!universally rigid-Noetherian formal scheme@universally rigid-Noetherian ---!locally universally rigid-Noetherian formal scheme@locally --- ---}, and $X\rightarrow Y$ is locally of finite type$;$
\item[{\rm (b)}] $Y'\rightarrow Y$ is adically flat\index{morphism of formal schemes@morphism (of formal schemes)!adically flat morphism of formal schemes@adically flat ---}.
\end{itemize}
Then we have a canonical isomorphism
$$
\widehat{g}^{\ast}\Omega^1_{X/Y}\stackrel{\sim}{\longrightarrow}\Omega^1_{X'/Y'}. \eqno{\square}
$$
\end{prop}

The following proposition follows immediately from \ref{cor-adicdifferentialcompatibility1}:
\begin{prop}\label{prop-differentialinvariants2}
Let $f\colon X\rightarrow Y$ be a morphism between adic formal schemes of finite ideal type, and suppose that $Y$ has an ideal of definition $\mathscr{I}$ of finite type.
Set $X_k=(X,\O_X/\mathscr{I}^{k+1}\O_X)$ and $Y_k=(Y,\O_Y/\mathscr{I}^{k+1})$ for any $k\geq 0$.
Then we have natural isomorphisms
$$
\Omega^1_{X/Y}/\mathscr{J}^{k+1}\Omega^1_{X/Y}\cong\Omega^1_{X/Y}\otimes_{\O_X}\O_{X_k}\cong\Omega^1_{X_k/Y_k}
$$
for $k\geq 0$, where $\Omega^1_{X_k/Y_k}$ denotes the usual differential module for the map of schemes $f_k\colon X_k\rightarrow Y_k$. \hfill$\square$
\end{prop}

\subsubsection{Differentials on universally rigid-Noetherian formal schemes}\label{subsub-differentialinvformalschfintype}
By \ref{cor-difffiniteness2} and \ref{thm-adicqcohpre1}, we have:
\begin{prop}\label{prop-differentialinvformalschfintype2}
Let $Y$ be a locally universally rigid-Noetherian formal scheme\index{formal scheme!universally rigid-Noetherian formal scheme@universally rigid-Noetherian ---!locally universally rigid-Noetherian formal scheme@locally --- ---} {\rm (\ref{dfn-formalsch})}, and $X$ an $Y$-formal scheme locally of finite presentation.
Then $\Omega^1_{X/Y}$ is a finitely presented $\O_X$-module. \hfill$\square$
\end{prop}

By \ref{thm-fundamentalexactseqdiff} we have:
\begin{thm}[Fundamental exact sequences]\label{thm-fundamentalexactseqdiffsch}
Let $Z$ be a locally universally rigid-Noetherian formal scheme\index{formal scheme!universally rigid-Noetherian formal scheme@universally rigid-Noetherian ---!locally universally rigid-Noetherian formal scheme@locally --- ---}.

{\rm (1)} Let $f\colon X\rightarrow Y$ and $Y\rightarrow Z$ be locally of finite type morphisms of adic formal schemes of finite ideal type.
Then we have the natural exact sequence
$$
\widehat{f}^{\ast}\Omega^1_{Y/Z}\longrightarrow\Omega^1_{X/Z}\longrightarrow\Omega^1_{X/Y}\longrightarrow 0
$$
of adically quasi-coherent sheaves on $X$ of finite type, where $\widehat{f^{\ast}}$ denotes the complete pull-back\index{complete!complete pull-back@--- pull-back} {\rm (\ref{dfn-completepullback}; cf.\ \ref{prop-completepullbackaqcsheaves0})}.

{\rm (2)} Let $Y$ be an adic formal scheme locally of finite type over $Z$, and $i\colon X\hookrightarrow Y$ be an immersion of finite presentation.
Let $\mathscr{N}_{X/Y}$ be the conormal sheaf\index{conormal!conormal sheaf@--- sheaf} of $X$ in $Y$ {\rm (Exercise \ref{exer-conormalsheaf})}.
Then we have the canonical exact sequence
$$
\mathscr{N}_{X/Y}\longrightarrow\widehat{i}^{\ast}\Omega^1_{Y/Z}\longrightarrow\Omega^1_{X/Z}\longrightarrow 0
$$
of adically quasi-coherent sheaves on $X$ of finite type. \hfill$\square$
\end{thm}

\subsection{\'Etale and smooth morphisms}\label{sub-etaleformal}
\subsubsection{Neat morphisms}\label{subsub-neat}
\index{morphism of formal schemes@morphism (of formal schemes)!neat morphism of formal schemes@neat ---|(}
\begin{dfn}\label{dfn-adicallyfinitepresented}{\rm 
(1) Let $f\colon X\rightarrow Y$ be an adic morphism between adic formal schemes of finite ideal type.
Then $f$ is said to be {\em adically locally of finite presentation}\index{morphism of formal schemes@morphism (of formal schemes)!adically locally of finite presentation morphism of formal schemes@adically locally of finite presentation ---} if the following condition is satisfied: for any point $x\in X$ there exists an open neighborhood $U$ of $x$ in $X$ and an open neighborhood $V$ of $f(x)$ in $Y$ that has an ideal of definition $\mathscr{I}$ of finite type such that
\begin{itemize}
\item[{\rm (a)}] $f(U)\subseteq V;$
\item[{\rm (b)}] if we set $U_k=(U,\O_U/\mathscr{I}^{k+1}\O_U)$ and $V_k=(V,\O_V/\mathscr{I}^{k+1})$ for $k\geq 0$, then the induced map $U_k\rightarrow V_k$ of schemes is of finite presentation for any $k\geq 0$.
\end{itemize}

(2) We say that $f$ is {\em adically of finite presentation}\index{morphism of formal schemes@morphism (of formal schemes)!adically of finite presentation morphism of formal schemes@adically of finite presentation ---} if it is adically locally of finite presentation and quasi-compact\index{morphism of formal schemes@morphism (of formal schemes)!quasi-compact morphism of formal schemes@quasi-compact ---}\index{quasi-compact!quasi-compact morphism of formal schemes@--- morphism (of formal schemes)} (\ref{dfn-qcompformal}).}
\end{dfn}

Notice that these are the properties of the form `adically $P$' (\S\ref{subsub-adicalization}) where $P=$ `locally of finite presentation' or `of finite presentation'.
By \ref{prop-cortopfintype11} we know that a morphism $f\colon X\rightarrow Y$ of adically finite presentation (resp.\ adically locally of finite presentation) is of finite type (resp.\ locally of finite type).
Moreover, if $Y$ is locally universally rigid-Noetherian\index{formal scheme!universally rigid-Noetherian formal scheme@universally rigid-Noetherian ---!locally universally rigid-Noetherian formal scheme@locally --- ---}, then $f$ is of adically finite presentation (resp.\ adically locally of finite presentation) if and only if it is of finite presentation (resp.\ locally of finite presentation) (\ref{cor-topfinpres11}).

\begin{prop}\label{prop-neat1}
Let $f\colon X\rightarrow Y$ be a morphism adically locally of finite presentation between adic formal schemes of finite ideal type, and suppose that $Y$ has an ideal of definition $\mathscr{I}$ of finite type.
Set $X_k=(X,\O_X/\mathscr{I}^{k+1}\O_X)$ and $Y_k=(Y,\O_Y/\mathscr{I}^{k+1})$ for $k\geq 0$, and let $f_k\colon X_k\rightarrow Y_k$ be the induced map of schemes.
Then the following conditions are equivalent$:$
\begin{itemize}
\item[{\rm (a)}] $\Omega^1_{X/Y}=0;$
\item[{\rm (b)}] $\Omega^1_{X_k/Y_k}=0$ for any $k\geq 0;$
\item[{\rm (c)}] $\Omega^1_{X_0/Y_0}=0;$
\item[{\rm (d)}] the diagonal map $\Delta_X\colon X\rightarrow X\times_YX$ is an open immersion.
\end{itemize}
\end{prop}

\begin{proof}
The implications (a) $\Rightarrow$ (b) $\Rightarrow$ (c) are clear due to \ref{prop-differentialinvariants2}.
To show (c) $\Rightarrow$ (a), we may work locally and may assume that $X$ and $Y$ are affine $X=\Spf B$ and $Y=\Spf A$, where $A$ and $B$ are adic rings of finite ideal type and $I\subseteq A$ is a finitely generated ideal of definition such that $\mathscr{I}=I^{\Delta}$.
Then (c) implies that $\widehat{\Omega}^1_{B/A}/I\widehat{\Omega}^1_{B/A}=0$.
Since $\widehat{\Omega}^1_{B/A}$ is a finitely generated $B$-algebra (\ref{cor-difffiniteness1}), we have $\widehat{\Omega}^1_{B/A}=0$ by Nakayama's lemma, whence (a).
The condition (d) is equivalent to that for any $k\geq 0$ the diagonal map $\Delta_{X_k}\colon X_k\rightarrow X_k\times_{Y_k}X_k$ is an open immersion (cf.\ \ref{cor-prodformal21}).
Hence by \cite[$\mathbf{IV}$, (17.4.2)]{EGA} we have the equivalence (b) $\Leftrightarrow$ (d).
\end{proof}

\begin{dfn}\label{dfn-neat}{\rm 
(1) An adic morphism $f\colon X\rightarrow Y$ of adic formal schemes of finite ideal type is said to be {\em neat} or {\em unramified}\index{morphism of formal schemes@morphism (of formal schemes)!unramified morphism of formal schemes@unramified ---|see{neat morphism}} if it is adically locally of finite presentation and $\Omega^1_{X/Y}=0$.

(2) Let $f\colon X\rightarrow Y$ be an adic morphism of adic formal schemes of finite ideal type, and $x\in X$ a point.
We say that $f$ is {\em neat at $x$} or {\em unramified at $x$} if there exists an open neighborhood $U$ of $x$ in $X$ such that the map $U\rightarrow Y$ is neat.}
\end{dfn}

\begin{prop}\label{prop-neatmorformal0}
Let $f\colon X\rightarrow Y$ be an adic morphism of adic formal schemes of finite ideal type, and suppose $Y$ has an ideal of definition $\mathscr{I}$ of finite type.
Set $X_k=(X,\O_X/\mathscr{I}^{k+1}\O_X)$ and $Y_k=(Y,\O_Y/\mathscr{I}^{k+1})$ for $k\geq 0$.
Let $f_k\colon X_k\rightarrow Y_k$ be the induced morphism of schemes.
Then the following conditions are equivalent$:$
\begin{itemize}
\item[{\rm (a)}] $f$ is neat$;$
\item[{\rm (b)}] $f_k$ is neat for any $k\geq 0.$
\end{itemize}
If we assume, moreover, that $f$ is adically locally of finite presentation, then the conditions are equivalent to 
\begin{itemize}
\item[{\rm (c)}] $f_0$ is neat. \hfill$\square$
\end{itemize}
\end{prop}

This follows immediately from \ref{prop-neat1}.
By \ref{prop-immersionimmformal4} and \cite[$\mathbf{IV}$, (17.3.3) (i)]{EGA} we have the following:
\begin{prop}\label{prop-neatmorformal10}
An immersion\index{immersion!immersion of formal schemes@--- (of formal schemes)} is neat if and only if it is adically locally of finite presentation. \hfill$\square$
\end{prop}

\begin{prop}\label{prop-neatermorformal2}
{\rm (1)} If $f\colon Z\rightarrow Y$ and $g\colon Y\rightarrow X$ are neat, then so is the composition $g\circ f$.

{\rm (2)} For any neat $S$-morphisms $f\colon X\rightarrow Y$ and $g\colon X'\rightarrow Y'$ of adic formal schemes of finite ideal type over an adic formal scheme $S$ of finite ideal type, the induced morphism $f\times_Sg\colon X\times_SY\rightarrow X'\times_SY'$ is neat.

{\rm (3)} For any neat $S$-morphism $f\colon X\rightarrow Y$ of adic formal schemes of finite ideal type over an adic formal scheme $S$ of finite ideal type, and for any morphism $S'\rightarrow S$ of adic formal schemes of finite ideal type, the induced morphism $f_{S'}\colon X\times_SS'\rightarrow Y\times_SS'$ is neat.

{\rm (4)} Suppose that the composition $g\circ f$ of two adic morphisms of adic formal schemes of finite ideal type is neat and that $g$ is adically locally of finite presentation. 
Then $f$ is neat.
\end{prop}

\begin{proof}
(1) follows from \cite[$\mathbf{IV}$, (17.3.3) (ii)]{EGA}.
By \ref{prop-genpropertymorphismadic1} and \ref{prop-neatmorformal0} the assertions (2) and (3) follow from \cite[$\mathbf{IV}$, (17.3.3) (iii), (iv)]{EGA}.
Finally, (4) follows from \ref{prop-neatmorformal0} and \cite[$\mathbf{IV}$, (17.3.3) (v)]{EGA}.
\end{proof}

In view of \ref{prop-genpropertymorphismadic2} we have:
\begin{cor}\label{cor-neatmorformal00}
Let $f\colon X\rightarrow Y$ be a morphism of schemes, and $Z$ a closed subscheme of $Y$ of finite presentation.
If $f$ is neat, then the formal completion $\widehat{f}\colon\widehat{X}|_{f^{-1}(Z)}\rightarrow\widehat{Y}_Z$ is neat. \hfill$\square$
\end{cor}

The following two propositions follow from \ref{thm-fundamentalexactseqdiffsch} (1):
\begin{prop}\label{prop-diffmoduleetaleglue1}
Consider morphisms $f\colon X\rightarrow Y$ and $g\colon Y\rightarrow Z$ of locally universally rigid-Noetherian formal schemes.
Suppose that $f$ is locally of finite type and that $g$ is neat.
Then we have the canonical isomorphism
$$
\Omega^1_{X/Z}\stackrel{\sim}{\longrightarrow}\Omega^1_{X/Y}
$$
of adically quasi-coherent sheaves on $X$. \hfill$\square$
\end{prop}

\begin{prop}\label{prop-diffmoduleneatsurjadq}
Let $Z$ be a locally universally rigid-Noetherian formal scheme, and $f\colon X\rightarrow Y$ a $Z$-morphism between adic formal schemes that are locally of finite type over $Z$.
If $f$ neat, then the canonical map 
$$
\widehat{f}^{\ast}\Omega^1_{Y/Z}\longrightarrow\Omega^1_{X/Z}
$$
is surjective.
The converse holds if $f$ is locally of finite presentation. \hfill$\square$
\end{prop}
\index{morphism of formal schemes@morphism (of formal schemes)!neat morphism of formal schemes@neat ---|)}

\subsubsection{\'Etale morphisms}\label{subsub-etaleformal}
\index{morphism of formal schemes@morphism (of formal schemes)!etale morphism of formal schemes@\'etale ---|(}
\begin{dfn}\label{dfn-etaleformaletale}{\rm 
(1) An adic morphism $f\colon X\rightarrow Y$ between adic formal schemes of finite ideal type is said to be {\em \'etale} if it is neat\index{morphism of formal schemes@morphism (of formal schemes)!neat morphism of formal schemes@neat ---} and adically flat\index{morphism of formal schemes@morphism (of formal schemes)!adically flat morphism of formal schemes@adically flat ---} (\ref{dfn-adicallyflat}).

(2) Let $f\colon X\rightarrow Y$ be an adic morphism of adic formal schemes, and $x\in X$ a point.
We say that $f$ is {\em \'etale at $x$} if there exists an open neighborhood $U$ of $x$ in $X$ such that the map $U\rightarrow Y$ is \'etale.}
\end{dfn}

\'Etale morphisms are adically locally of finitely presentation\index{morphism of formal schemes@morphism (of formal schemes)!adically locally of finite presentation morphism of formal schemes@adically locally of finite presentation ---} {\rm (\ref{dfn-adicallyfinitepresented} (1))} and hence are locally of finite type.
Moreover, if $f\colon X\rightarrow Y$ is \'etale and $Y$ is locally universally rigid-Noetherian\index{formal scheme!universally rigid-Noetherian formal scheme@universally rigid-Noetherian ---!locally universally rigid-Noetherian formal scheme@locally --- ---}, then $f$ is locally of finite presentation (due to \ref{cor-topfinpres11}) and flat (due to \ref{prop-adicallyflatformal01}).

\begin{prop}\label{prop-etaleformal00}
Let $f\colon X\rightarrow Y$ be an adic morphism of adic formal schemes of finite ideal type.
Suppose that $Y$ has an ideal of definition $\mathscr{I}$ of finite type.
Set $X_k=(X,\O_X/\mathscr{I}^{k+1}\O_X)$ and $Y_k=(Y,\O_Y/\mathscr{I}^{k+1})$ for $k\geq 0$, and denote by $f_k\colon X_k\rightarrow Y_k$ the induced morphism of schemes.
Then the following conditions are equivalent$:$
\begin{itemize}
\item[{\rm (a)}] $f$ is \'etale$;$
\item[{\rm (b)}] $f_k$ is \'etale for $k\geq 0$. \hfill$\square$
\end{itemize}
\end{prop}

This follows from \ref{prop-neatmorformal0} and \ref{prop-adicallyflatformal00}.
\begin{cor}\label{cor-etaleformal33}
Let $f\colon X\rightarrow Y$ be a morphism of schemes, and $Z$ a closed subscheme of $Y$ of finite presentation.
If $f$ is \'etale, then the formal completion $\widehat{f}\colon\widehat{X}|_{f^{-1}(Z)}\rightarrow\widehat{Y}|_Z$ is \'etale. \hfill$\square$
\end{cor}

\begin{prop}\label{prop-etaleformal1}
{\rm (1)} An open immersion of adic formal schemes of finite type is \'etale.

{\rm (2)} The composition of two \'etale morphisms is \'etale.
Let $f\colon X\rightarrow Y$ and $g\colon Y\rightarrow Z$ be adically locally of finite presentation morphisms.
Suppose $f$ is adically faithfully flat\index{morphism of formal schemes@morphism (of formal schemes)!adically faithfully flat morphism of formal schemes@adically faithfully flat ---}.
Then if $g\circ f$ is \'etale, so is $g$.

{\rm (3)} For any \'etale $S$-morphisms $f\colon X\rightarrow Y$ and $g\colon X'\rightarrow Y'$ of adic formal schemes of finite ideal type over an adic formal scheme $S$ of finite ideal type, the induced morphism $f\times_Sg\colon X\times_SY\rightarrow X'\times_SY'$ is \'etale.

{\rm (4)} For any \'etale $S$-morphism $f\colon X\rightarrow Y$ of adic formal schemes of finite ideal type over an adic formal scheme $S$ of finite ideal type and for any morphism $S'\rightarrow S$ of adic formal schemes of finite ideal type, the induced morphism $f_{S'}\colon X\times_SS'\rightarrow Y\times_SS'$ is \'etale.
\end{prop}

\begin{proof}
(1) and the first part of (2) are clear.
The second part of (2) follows from \ref{prop-etaleformal00} and \cite[$\mathbf{IV}$, (17.7.7)]{EGA}.
By \ref{prop-genpropertymorphismadic1} the assertions (3) and (4) follow from \cite[$\mathbf{IV}$, (17.3.3)]{EGA}.
\end{proof}

\begin{prop}\label{prop-diffmoduleetaleglue2}
Consider morphisms $f\colon X\rightarrow Y$ and $g\colon Y\rightarrow Z$ of adic formal schemes of finite ideal type.
Suppose that $f$ is \'etale and that $g$ is locally of finite type.
Then we have the canonicall isomorphism
$$
\widehat{f}^{\ast}\Omega^1_{Y/Z}\stackrel{\sim}{\longrightarrow}\Omega^1_{X/Z}
$$
of adically quasi-coherent sheaves on $X$.
\end{prop}

\begin{proof}
As the question is local, we may assume that $Y$ and $Z$ have ideals of definition $\mathscr{J}$ and $\mathscr{I}$, respectively, such that $\mathscr{I}\O_Y\subseteq\mathscr{J}$.
Set $X_k=(X,\O_X/\mathscr{J}^{k+1}\O_X)$, $Y_k=(Y,\O_Y/\mathscr{J}^{k+1})$, and $Z_k=(Z,\O_Z/\mathscr{I}^{k+1})$ for $k\geq 0$.
By \ref{prop-differentialinvariants2} (cf.\ \ref{cor-adicdifferentialcompatibility1}) and \ref{cor-completepullbackaqcsheaves0} the differential module $\widehat{f^{\ast}}\Omega^1_{Y/Z}$ (resp.\ $\Omega^1_{X/Z}$) coincides with the projective limit $\varprojlim_{k\geq 0}f^{\ast}_k\Omega^1_{Y_k/Z_k}$ (resp.\ $\varprojlim_{k\geq 0}\Omega^1_{X_k/Z_k}$).
As $f_k\colon X_k\rightarrow Y_k$ is \'etale (\ref{prop-etaleformal00}), we have $f^{\ast}_k\Omega^1_{Y_k/Z_k}\cong\Omega^1_{X_k/Z_k}$ (\cite[$\mathbf{IV}$, (17.2.4)]{EGA}), and hence we get the desired isomorphism.
\end{proof}
\index{morphism of formal schemes@morphism (of formal schemes)!etale morphism of formal schemes@\'etale ---|)}

\subsubsection{Smooth morphisms}\label{subsub-smoothformal}
\index{morphism of formal schemes@morphism (of formal schemes)!smooth morphism of formal schemes@smooth ---|(}
We first remark that the following proposition holds due to \ref{prop-difffiniteness1}:
\begin{prop}\label{prop-affinedifffree}
Let $Y$ be an adic formal scheme of finite ideal type, and consider the affine $n$-space\index{affine!affine space in formal geometry@--- space (in formal geometry)} $\widehat{\A}^n_Y=\Spec\Z[X_1,\ldots,X_n]\times_{\Spec\Z}Y$ $($fiber product taken in the category $\Fs$ of all formal schemes$)$ over $Y$ {\rm (Exercise \ref{exer-affinespaceformal} (1))}.
Then the differential module $\Omega^1_{\widehat{\A}^n_Y/Y}$ is a free $\O_{\widehat{\A}^n_Y}$-module with the basis $dX_1,\ldots,dX_n$. \hfill$\square$
\end{prop}

\begin{dfn}\label{dfn-smoothformalsmooth}{\rm 
(1) An adic morphism $f\colon X\rightarrow Y$ between adic formal schemes of finite ideal type is said to be {\em smooth} if for any $x\in X$ there exists an open neighborhood $U$ of $x$ in $X$ and a commutative diagram
$$
\xymatrix@R-2ex@C-3.5ex{U\ar[rr]\ar[dr]&&\widehat{\A}^n_Y\ar[dl]\\ &Y}
$$
(for some $n$ that depends on $x$), where the horizontal arrow is \'etale\index{morphism of formal schemes@morphism (of formal schemes)!etale morphism of formal schemes@\'etale ---}.

(2) Let $f\colon X\rightarrow Y$ be an adic morphism of adic formal schemes of finite ideal type, and $x\in X$ a point.
We say that $f$ is {\em smooth at $x$} if there exists an open neighborhood $U$ of $x$ in $X$ such that the map $U\rightarrow Y$ is smooth.}
\end{dfn}

Here are the immediate consequences of the definition: \'etale morphisms are smooth, and smooth morphisms are adically locally of finite presentation\index{morphism of formal schemes@morphism (of formal schemes)!adically locally of finite presentation morphism of formal schemes@adically locally of finite presentation ---} (hence is locally of finite type) and is adically flat\index{morphism of formal schemes@morphism (of formal schemes)!adically flat morphism of formal schemes@adically flat ---}; the affine $n$-space $\widehat{\A}^n_Y$ over an adic formal scheme $Y$ is smooth over $Y$.

The following proposition is a direct consequence of the definition and \ref{prop-diffmoduleetaleglue2}:
\begin{prop}\label{prop-smoothformal1}
If $f\colon X\rightarrow Y$ is smooth, then $\Omega^1_{X/Y}$ is a locally free $\O_X$-module of finite rank.\hfill$\square$
\end{prop}

\begin{prop}\label{prop-smoothformal00}
Let $f\colon X\rightarrow Y$ be an adic morphism of adic formal schemes of finite ideal type.
Suppose that $Y$ has an ideal of definition $\mathscr{I}$ of finite type.
Set $X_k=(X,\O_X/\mathscr{I}^{k+1}\O_X)$ and $Y_k=(Y,\O_Y/\mathscr{I}^{k+1})$ for $k\geq 0$, and denote by $f_k\colon X_k\rightarrow Y_k$ the induced morphism of schemes.
Then the following conditions are equivalent$:$
\begin{itemize}
\item[{\rm (a)}] $f$ is smooth$;$
\item[{\rm (b)}] $f_k$ is smooth for $k\geq 0$.
\end{itemize}
\end{prop}

\begin{proof}
The implication (a) $\Rightarrow$ (b) follows from \ref{prop-etaleformal00} and \cite[$\mathbf{IV}$, (17.11.4)]{EGA}.
To show the converse, we may work locally; we may assume that $X_0$ admits an \'etale $Y_0$-map $X_0\rightarrow\A^n_{Y_0}$ to the affine $n$-space over the scheme $Y_0$.
We may furthermore assume that $X$ and $Y$ are affine; set $X=\Spf B$ and $Y=\Spf A$, where $A$ and $B$ are adic rings of finite ideal type, and let $I\subseteq A$ be the finitely generated ideal of definition with $\mathscr{I}=I^{\Delta}$.
Set $A_k=A/I^{k+1}$ and $B_k=B/I^{k+1}B$ for $k\geq 0$.
We already have an \'etale map $A_0[X_1,\ldots,X_n]\rightarrow B_0$ over $A_0$.
Take $g_1,\ldots,g_n\in B$ whose images in $B_0$ are the images of $X_1,\ldots,X_n$, respectively; since $B$ is $I$-adically complete, one has the $A$-algebra homomorphism $A\dl X_1,\ldots,X_n\dr\rightarrow B$, mapping $X_i$ to $g_i$ for $i=1,\ldots,n$, which extends the map $A_0[X_1,\ldots,X_n]\rightarrow B_0$.
Since this is an $A$-algebra morphism, it is an adic map.
Thus we have the $Y$-morphism $g\colon X\rightarrow\widehat{\A}^n_Y$ such that $g_0$ is \'etale.
Now since $f$ is adically locally of finite presentation, we deduce that $g$ is adically locally of finite presentation; by \ref{prop-neatmorformal0} we know that $g$ is neat.
Now by \ref{prop-adicallyflatlocalcriterian} we deduce that $g$ is adically flat, and hence is \'etale.
\end{proof}

\begin{cor}\label{cor-smoothformal33}
Let $f\colon X\rightarrow Y$ be a morphism of schemes, and $Z$ a closed subscheme of $Y$ of finite presentation.
If $f$ is smooth, then the formal completion $\widehat{f}\colon\widehat{X}|_{f^{-1}(Z)}\rightarrow\widehat{Y}|_Z$ is smooth. \hfill$\square$
\end{cor}

\begin{prop}\label{prop-smoothformal2}
{\rm (1)} The composition of two smooth morphisms is smooth.
Let $f\colon X\rightarrow Y$ and $g\colon Y\rightarrow Z$ be morphisms adically locally of finite presentation.
Suppose $f$ is adically faithfully flat.
Then if $g\circ f$ is smooth, so is $g$.

{\rm (2)} For any smooth $S$-morphisms $f\colon X\rightarrow Y$ and $g\colon X'\rightarrow Y'$ of adic formal schemes of finite ideal type over an adic formal scheme $S$ of finite ideal type, the induced morphism $f\times_Sg\colon X\times_SY\rightarrow X'\times_SY'$ is smooth.

{\rm (3)} For any smooth $S$-morphism $f\colon X\rightarrow Y$ of adic formal schemes of finite ideal type over an adic formal scheme $S$ of finite ideal type and for any morphism $S'\rightarrow S$ of adic formal schemes of finite ideal type, the induced morphism $f_{S'}\colon X\times_SS'\rightarrow Y\times_SS'$ is smooth.
\end{prop}

\begin{proof}
The first part of (1) follows from \ref{prop-smoothformal00} and \cite[$\mathbf{IV}$, (17.3.3) (ii)]{EGA}.
The second part of (1) follows from \ref{prop-smoothformal00} and \cite[$\mathbf{IV}$, (17.7.7)]{EGA}.
By \ref{prop-genpropertymorphismadic1} the assertions (2) and (3) follow from \cite[$\mathbf{IV}$, (17.3.3)]{EGA}.
\end{proof}

\begin{prop}\label{prop-smoothformalexactseq}
Let $Z$ be an adic formal scheme of finite ideal type, and $Y\rightarrow Z$ a locally of finite type morphism of adic formal schemes of finite ideal type.
Let $f\colon X\rightarrow Y$ be a smooth morphism.
Then we have the canonical exact sequence
$$
0\longrightarrow\widehat{f}^{\ast}\Omega^1_{Y/Z}\longrightarrow\Omega^1_{X/Z}\longrightarrow\Omega^1_{X/Y}\longrightarrow 0\leqno{(\ast)}
$$
of adically quasi-coherent sheaves on $X$, where $\Omega^1_{X/Y}$ is locally free of finite type$;$ in particular, this exact sequence splits locally.
\end{prop}

\begin{proof}
As the question is local, we may assume that $Z$ has an ideal of definition of finite type $\mathscr{I}$.
Set $X_k=(X,\O_X/\mathscr{I}^{k+1}\O_X)$, $Y_k=(Y,\O_Y/\mathscr{I}^{k+1}\O_Y)$, and $Z_k=(Z,\O_Z/\mathscr{I}^{k+1})$ for $k\geq 0$.
Since $f_k\colon X_k\rightarrow Y_k$ is smooth (\ref{prop-smoothformal00}), the sequence
$$
0\longrightarrow f^{\ast}_k\Omega^1_{Y_k/Z_k}\longrightarrow\Omega^1_{X_k/Z_k}\longrightarrow\Omega^1_{X_k/Y_k}\longrightarrow 0\leqno{(\ast)_k}
$$
is exact for any $k\geq 0$ (\cite[$\mathbf{IV}$, (17.2.3) (ii)]{EGA}).
By \ref{prop-differentialinvariants2} and \ref{cor-completepullbackaqcsheaves0} the desired sequence $(\ast)$ is the one obtained from $(\ast)_k$ by passage to the projective limits $\varprojlim_{k\geq 0}$.
Since $\{f^{\ast}_k\Omega^1_{Y_k/Z_k}\}_{k\geq 0}$ is strict, the exactness of $(\ast)$ follows from {\bf \ref{ch-pre}}.\ref{prop-ML5ML5} (1) (cf.\ {\bf \ref{ch-pre}}.\ref{rem-coheffaceable} (2)).
The last assertion follows from \ref{prop-smoothformal1}.
\end{proof}
\index{morphism of formal schemes@morphism (of formal schemes)!smooth morphism of formal schemes@smooth ---|)}

\addcontentsline{toc}{subsection}{Exercises}
\subsection*{Exercises}
\begin{exer}\label{exer-affinespaceformal}{\rm 
Let $S$ be a formal scheme.

(1) Define $\widehat{\A}^n_S=\A^n_{\Z}\times_{\Spec\Z}S$ (fiber product taken in $\Fs$), and call it the {\em affine $n$-space over $S$}\index{affine!affine space in formal geometry@--- space (in formal geometry)}.
Show that $\widehat{\A}^n_S\rightarrow S$ is affine and that, if $S$ is adic, then $\widehat{\A}^n_S\rightarrow S$ is affine and adic.

(2) Define $\widehat{\P}^n_S=\P^n_{\Z}\times_{\Spec\Z}S$ (fiber product taken in $\Fs$), and call it the {\em projective $n$-space over $S$}\index{projective!projective space in formal geometry@--- space (in formal geometry)}.
Show that, if $S$ is adic, then $\widehat{\P}^n_S\rightarrow S$ is proper.}
\end{exer}

\begin{exer}\label{exer-vectorbundleformal}{\rm
Let $X$ be an adic formal scheme of finite ideal type.

(1) Show that any locally free $\O_X$-module $\mathscr{E}$ of finite type is an adically quasi-coherent sheaf.

(2) Show that the completion (\ref{dfn-completesheaf} (1)) of the symmetric algebra $\Sym_{\O_X}(\mathscr{E})$ is an adically quasi-coherent $\O_X$-algebra of finite type.
(The corresponding $X$-affine formal scheme (\ref{dfn-affinespectrumformalrel}) is denoted by $\widehat{\mathbf{V}}(\mathscr{E})$ and is called the {\em vector bundle over $X$ associated to $\mathscr{E}$}\index{vector bundle}.)}
\end{exer}

\begin{exer}[Conormal sheaf]\label{exer-conormalsheaf}{\rm 
Let $f\colon Y\hookrightarrow X$ be an immersion between adic formal schemes of finite ideal type.
Take a decomposition $f=j\circ i$ where $i\colon Y\hookrightarrow U$ is a closed immersion and $j\colon U\hookrightarrow X$ is an open immersion.
Let $\mathscr{K}$ be the defining ideal of the closed formal subscheme $Y$ of $U$, that is, the kernel of the map $\O_U\rightarrow i_{\ast}\O_Y$.
Then the $\O_Y$-module $i^{\ast}\mathscr{K}$ is called the {\em conormal sheaf}\index{conormal!conormal sheaf@--- sheaf} of $Y$ in $X$, and is denoted by $\mathscr{N}_{Y/X}$.
Show that, if $X$ is locally universally rigid-Noetherian and if $f$ is of finite presentation, then $\mathscr{N}_{Y/X}$ is an adically quasi-coherent sheaf on $Y$ of finite type.}
\end{exer}

\begin{exer}\label{exer-propermorformalschemex1}{\rm 
Let $A$ be an adic ring of finite ideal type, $X$ and $Y$ $A$-formal schemes adically of finite presentation, and $f\colon X\rightarrow Y$ a morphism adically of finite presentation.
Show that $f$ is proper if and only if for any integer $N\geq 0$ the induced morphism $f\times_Y\id_{\widehat{\A}^N_Y}\colon X\times_Y\widehat{\A}^N_Y\rightarrow\widehat{\A}^N_Y$ is closed.}
\end{exer}


\section{Formal algebraic spaces}\label{sec-formalalgsp}
The aim of this section is to lay the foundations of the theory of formal algebraic spaces.
Here, what we mean by {\em formal algebraic spaces} are more precisely spelled out as {\em quasi-separated adic formal algebraic spaces of finite ideal type}, generalizing {\em Noetherian} formal algebraic spaces already discussed in the last chapter of \cite{Knu}.

Notice that an important and non-trivial thing to do first, when introducing formal algebraic spaces as an \'etale sheaf, is to choose the domain category.
Since we only need, for the sake of our application to rigid geometry, to consider {\em adic} situation, we can obviously restrict to the category of adic formal schemes, but with arbitrary morphisms, not only adic ones, in order for several universalities and functorialities to make reasonable sense.
Notice that, for instance, by considering such a wide category, one can keep having final objects and, moreover, treating schemes at the same time (by this, in particular, the notion of formal algebraic spaces can include that of algebraic spaces).
In view of all this, we adopt as the base category the category $\Ac\Fs_S$ of adic formal schemes of finite ideal type over $S$ with arbitrary morphisms.

In \S\ref{sub-etaledescent} we establish the preliminaries on descent theory in formal geometry; as our interest is not only in defining {\em $($locally$)$ universally rigid-Noetherian} formal algebraic spaces but more general ones, flat descent is not the reasonable notion to discuss; instead, one should consider {\em adically flat descent}. 
In the first part, we are to give the general framework of adically flat descent for adically quasi-coherent sheaves and prove several basic theorems.

Next in \S\ref{sub-etaletopology} we discuss \'etale topology and \'etale sites on adic formal schemes of finite ideal type.
Based on these, the notion of adically quasi-coherent sheaves will be extended on the \'etale site will be given.

The main part of this section is \S\ref{sub-formalalgebraicspaces}, where we give the definition of formal algebraic spaces and discuss basic geometries of them.
Several properties for morphisms between formal algebraic spaces are defined and discussed in \S\ref{sub-propertiesmorformalalgsp}.
Then this section ends up with the subsection \S\ref{sub-adequatealralgsp}, in which we will focus on locally universally rigid-Noetherian formal algebraic spaces and locally universally adhesive formal algebraic spaces.

\subsection{Adically flat descent}\label{sub-etaledescent}
\index{descent!adically flat descent@adically flat ---|(}
\subsubsection{Basic assertions}\label{subsub-etaledescent}
As usual, a diagram of the form 
$$
\xymatrix{M_0\ar[r]^f&M_1\ar@<.5ex>[r]^{g_1}\ar@<-.5ex>[r]_{g_2}&M_2}
$$
consisting of sets (or sheaves of sets, abelian groups, etc.) is said to be {\em exact} if the map $f$ coincides with the difference kernel of $g_1$ and $g_2$, that is, $f$ is injective, $g_1\circ f=g_2\circ f$, and the image of $f$ coincides with the locus of coincidence of the maps $g_1$ and $g_2$.
In case of (sheaves of) abelian groups, this amounts to the same as that 
$$
0\longrightarrow M_0\stackrel{f}{\longrightarrow}M_1\stackrel{g_1-g_2}{\longrightarrow}M_2
$$ 
is exact in the usual sense.

\begin{sit}\label{sit-etaledescent}{\rm 
Let $q\colon Y\rightarrow X$ be a quasi-compact\index{morphism of formal schemes@morphism (of formal schemes)!quasi-compact morphism of formal schemes@quasi-compact ---}\index{quasi-compact!quasi-compact morphism of formal schemes@--- morphism (of formal schemes)} adically faithfully flat\index{morphism of formal schemes@morphism (of formal schemes)!adically faithfully flat morphism of formal schemes@adically faithfully flat ---} (\ref{dfn-adicallyflat} (2)) morphism of adic formal schemes of finite ideal type\index{formal scheme!adic formal scheme@adic ---!adic formal scheme of finite ideal type@--- --- of finite ideal type}\index{adic!adic formal scheme@--- formal scheme!adic formal scheme of finite ideal type@--- --- of finite ideal type}.
We consider the diagram
$$
\xymatrix{X&Y\ar[l]_(.4){q}&Y\times_XY\ar@<-.5ex>[l]_(.55){p_1}\ar@<.5ex>[l]^(.55){p_2},}
$$
where $p_1$ and $p_2$ are the projections, and we set $p=q\circ p_1=q\circ p_2$.
If $X$ has an ideal of definition $\mathscr{I}$ of finite type, then for any non-negative integer $k\geq 0$ we have the induced diagram
$$
\xymatrix{X_k&Y_k\ar[l]_(.4){q_k}&Y_k\times_{X_k}Y_k\ar@<-.5ex>[l]_(.6){p_{1,k}}\ar@<.5ex>[l]^(.6){p_{2,k}}}
$$
of schemes, where $X_k=(X,\O_X/\mathscr{I}^{k+1})$ and $Y_k=(Y,\O_Y/\mathscr{I}^{k+1}\O_Y)$; notice that the map $q_k$ is a quasi-compact faithfully flat morphism between schemes (\ref{prop-adicallyflatformal00}).}
\end{sit}

\begin{prop}\label{prop-etaledescent1}
Consider the situation as in {\rm \ref{sit-etaledescent}}, and let $\mathscr{F}$ be an adically quasi-coherent sheaf on $X$ $($resp.\ an adically quasi-coherent $\O_X$-algebra$)$.
Then the induced diagram
$$
\xymatrix{\Gamma(X,\mathscr{F})\ar[r]&\Gamma(Y,\widehat{q^{\ast}}\mathscr{F})\ar@<.5ex>[r]\ar@<-.5ex>[r]&\Gamma(Y\times_XY,\widehat{p^{\ast}}\mathscr{F})}
$$
is exact.
\end{prop}

Recall that, since $q$ and $p$ are adically flat, the complete pull-backs\index{complete!complete pull-back@--- pull-back} $\widehat{q^{\ast}}\mathscr{F}$ and $\widehat{p^{\ast}}\mathscr{F}$ are adically quasi-coherent sheaves (resp.\ adically quasi-coherent $\O_X$-algebras) (\ref{prop-completepullbackaqcsheaves1}).

\begin{prop}\label{prop-etaledescent2}
Consider the situation as in {\rm \ref{sit-etaledescent}}, and let $\mathscr{F}$ and $\mathscr{G}$ be either adically quasi-coherent sheaves on $X$, or adically quasi-coherent $\O_X$-algebras.
Then the induced diagram 
$$
\xymatrix{\Hom_{\O_X}(\mathscr{F},\mathscr{G})\ar[r]&\Hom_{\O_Y}(\widehat{q^{\ast}}\mathscr{F},\widehat{q^{\ast}}\mathscr{G})\ar@<.5ex>[r]\ar@<-.5ex>[r]&\Hom_{\O_{Y\times_XY}}(\widehat{p^{\ast}}\mathscr{F},\widehat{p^{\ast}}\mathscr{G})}
$$
is exact.
\end{prop}

The proofs of these propositions use the following lemma:
\begin{lem}\label{lem-etaledescent1}
Let $A\rightarrow B$ be an adic homomorphism of adic rings of finite ideal type such that the induced map $q\colon Y=\Spf B\rightarrow X=\Spf A$ is adically faithfully flat.
Let $M$ be an $I$-adically complete $A$-module, where $I\subseteq A$ is a finitely generated ideal of definition.
Then 
$$
\xymatrix{M\ar[r]&M\widehat{\otimes}_AB\ar@<.5ex>[r]\ar@<-.5ex>[r]&M\widehat{\otimes}_A(B\widehat{\otimes}_AB)}\leqno{(\ast)}
$$
is exact.
\end{lem}

\begin{proof}
Let $A_k=A/I^{k+1}$, $B_k=B/I^{k+1}B$, and $M_k=M/I^{k+1}M$ for $k\geq 0$.
Since $q_k\colon Y_k=\Spec B_k\rightarrow X_k=\Spec A_k$ is faithfully flat (\ref{prop-adicallyflatformal00}), the sequence
$$
\xymatrix{M_k\ar[r]&M_k\otimes_{A_k}B_k\ar@<.5ex>[r]\ar@<-.5ex>[r]&M_k\otimes_{A_k}(B_k\otimes_{A_k}B_k)}\leqno{(\ast)_k}
$$
is exact for any $k\geq 0$.
Since $(\ast)$ is the projective limit $\varprojlim_{k\geq 0}(\ast)_k$, the desired exactness follows from the left-exactness of the functor $\varprojlim_{k\geq 0}$ ({\bf \ref{ch-pre}}.\ref{prop-projlimleftexact}).
\end{proof}

\begin{proof}[Proof of Propositions {\rm \ref{prop-etaledescent1}} and {\rm \ref{prop-etaledescent2}}]
Since $\widehat{q^{\ast}}\O_X=\O_Y$ and $\widehat{p^{\ast}}\O_X=\O_{Y\times_XY}$, \ref{prop-etaledescent1} follows from \ref{prop-etaledescent2}.
Therefore, it suffices to show \ref{prop-etaledescent2}.

By a standard argument (cf.\ e.g.\ the proof of \cite[Chap.\ 6, Prop.\ 1]{BLR3}) one can easily reduce to the case where $X$ and $Y$ are affine $X=\Spf A$ and $Y=\Spf B$, where $A$ and $B$ are adic rings of finite ideal type.
Let $I\subseteq A$ be a finitely generated ideal of definition.
Take $I$-adically complete $A$-modules $M$ and $N$ corresponding to $\mathscr{F}$ and $\mathscr{G}$, that is, $\mathscr{F}=M^{\Delta}$ and $\mathscr{G}=N^{\Delta}$ (\ref{thm-adicqcoh1}).
Then by \ref{lem-corcoradicqcoh1} and \ref{lem-completepullbackaqcsheaves1lem} (2) what we need to show is the exactness of
\begin{equation*}
\begin{split}
\Hom_A(M,N)&\xymatrix{{}\ar[r]&{}}\Hom_B(M\widehat{\otimes}_AB,N\widehat{\otimes}_AB)\\ 
&\xymatrix{{}\ar@<.5ex>[r]\ar@<-.5ex>[r]&{}}\Hom_{B\widehat{\otimes}_AB}(M\widehat{\otimes}_AB\widehat{\otimes}_AB,N\widehat{\otimes}_AB\widehat{\otimes}_AB)\rlap{.}
\end{split}
\end{equation*}
This follows, by an easy diagram chasing, from the exactness of 
$$
\xymatrix{N\ar[r]&N\widehat{\otimes}_AB\ar@<.5ex>[r]\ar@<-.5ex>[r]&N\widehat{\otimes}_A(B\widehat{\otimes}_AB),}
$$
which is a consequence of \ref{lem-etaledescent1}.
\end{proof}

\subsubsection{Descent of morphisms}\label{subsub-etaledescentmorphism}
As defined in \S\ref{subsub-categoryformalschnotation}, for an adic formal scheme $X$ of finite ideal type we denote by $\Ac\Fs_X$ the category of adic formal schemes of finite ideal type over $X$ and (not necessarily adic) morphisms over $X$.
\begin{prop}\label{prop-morphismdescent}
Consider the situation as in {\rm \ref{sit-etaledescent}}, and let $Z$ and $W$ be adic formal schemes of finite ideal type over $X$.
Set $Z'=Y\times_XZ$ and $Z''=(Y\times_XY)\times_XZ$, and similarly for $W$.
Then the following sequence is exact$:$
$$
\xymatrix@C-2ex{\Hom_{\Ac\Fs_X}(Z,W)\ar[r]&\Hom_{\Ac\Fs_Y}(Z',W')\ar@<.5ex>[r]\ar@<-.5ex>[r]&\Hom_{\Ac\Fs_{Y\times_XY}}(Z'',W'').}
$$
\end{prop}

\begin{proof}
Obviously, we may assume that $X$ is affine $X=\Spf A$, where $A$ is an adic ring of finite ideal type, with a finitely generated ideal of definition $I\subseteq A$; this implies that $Y$ and $Y\times_XY$ are quasi-compact and have ideals of definition $I\O_Y$ and $I\O_{Y\times_XY}$ of finite type, respectively.

By base change with respect to $Z\rightarrow X$, we may assume $X=Z=\Spf A$.
Let $W=\bigcup_{\alpha\in L}W_{\alpha}$ be an open covering, and set $W'_{\alpha}=Y\times_XW_{\alpha}$ and $W''_{\alpha}=(Y\times_XY)\times_XW_{\alpha}$ for each $\alpha\in L$.
Suppose that the assertion is true with $W$, $W'$, and $W''$ replaced by $W_{\alpha}$, $W'_{\alpha}$, and $W''_{\alpha}$, respectively, for any $\alpha\in L$.
Let $f'\in\Hom_Y(Z',W')$ be a morphism with the same image in $\Hom_{Y\times_XY}(Z'',W'')$.
Set $Z'_{\alpha}=Z'\times_{W'}W'_{\alpha}$ and $Z''_{\alpha}=Z''\times_{W''}W''_{\alpha}$ for each $\alpha\in L$.
By \cite[Expos\'e VIII, (4.4)]{SGA1} there exists for each $\alpha$ an open subset $Z_{\alpha}$ of $Z$ such that $Y\times_XZ_{\alpha}=Z'_{\alpha}$; indeed, since the two maps from $Z''_{\alpha}$ to $Z'_{\alpha}$ are adic, one can reduce it modulo powers of $\mathscr{J}$ to find the inductive system of schemes supported on the open subset of $Z$, which gives the desired open formal subscheme $Z_{\alpha}$.
Moreover, since $Z'\rightarrow Z$ is adically faithfully flat, $\{Z_{\alpha}\}_{\alpha\in L}$ gives an open covering of $Z$.
By our assumption we can find the unique morphism $f_{\alpha}\colon Z_{\alpha}\rightarrow W_{\alpha}$ such that $f_{\alpha}\times_XY=f'|_{Z'_{\alpha}}$ for each $\alpha\in L$.
By the uniqueness we have the unique map $f\colon Z\rightarrow W$ such that $f\times_XY=f'$, as desired.

\medskip
Thus we may assume that both $Z$ and $W$ are affine; let $\mathscr{K}$ be an ideal of definition of finite type of $W$.
We may assume that $\mathscr{I}\O_Z\subseteq\mathscr{J}$ and $\mathscr{I}\O_W\subseteq\mathscr{K}$.
For $k\geq 0$, let $\Hom^{(k)}_X(Z,W)$ (resp.\ $\Hom^{(k)}_Y(Z',W')$, resp.\ $\Hom^{(k)}_{Y\times_XY}(Z'',W'')$) be the subset of $\Hom_X(Z,W)$ (resp.\ $\Hom_Y(Z',W')$, resp.\ $\Hom_{Y\times_XY}(Z'',W'')$) consisting of morphisms such that $\mathscr{K}^{k+1}\O_Z\subseteq\mathscr{J}$ (resp.\ $\mathscr{K}^{k+1}\O_{Z'}\subseteq\mathscr{J}\O_{W'}$, resp.\ $\mathscr{K}^{k+1}\O_{Z''}\subseteq\mathscr{J}\O_{W''}$).
Since $\Hom_X(Z,W)$ is the union $\Hom^{(k)}_X(Z,W)$ for $k\geq 0$, it suffices to show that the sequence
$$
\xymatrix@C-2ex{\Hom^{(k)}_X(Z,W)\ar[r]&\Hom^{(k)}_Y(Z',W')\ar@<.5ex>[r]\ar@<-.5ex>[r]&\Hom^{(k)}_{Y\times_XY}(Z'',W'')}
$$
is exact for each $k\geq 0$; replacing $\mathscr{K}$ by its powers, we may restrict to the case $k=0$.
Set $Z_k=(Z,\O_Z/\mathscr{J}^{k+1})$ and $W_k=(W,\O_W/\mathscr{K}^{k+1})$ for $k\geq 0$.
By \ref{prop-formalindlimschfact} the canonical map
$$
\Hom^{(0)}_X(Z,W)\longrightarrow\varprojlim_{k\geq 0}\Hom_{X_k}(Z_k,W_k)
$$
is bijective, and similarly for the other hom sets.
Hence the desired exactness follows from the scheme case \cite[Expos\'e VIII, \S5]{SGA1} and the left-exactness of $\varprojlim_{k\geq 0}$ ({\bf \ref{ch-pre}}.\ref{prop-projlimleftexact}).
\end{proof}

\begin{cor}\label{cor-morphismdescent1}
In the situation as in {\rm \ref{sit-etaledescent}}, the functor $\Ac\Fs_X\rightarrow\Ac\Fs_Y$ given by $Z\mapsto Z\times_XY$ is faithful. \hfill$\square$
\end{cor}

\subsubsection{Descent of properties of morphisms}\label{subsub-adicallyflatdescentpropertymorphisms}
We continue with working in the situation as in \ref{sit-etaledescent}.
Let $Z,W$ be adic formal schemes of finite ideal type over $X$, and $f\colon Z\rightarrow W$ a morphism over $X$.
We consider the $Y$-morphism $f\times_XY\colon Z\times_XY\rightarrow W\times_XY$ obtained by base change.
By \ref{prop-adicnessadicallyflatdescent} we have the following `descent of adicness' by adically faithfully flat quasi-compact morphisms:
\begin{prop}\label{prop-morphismdescentadic}
The morphism $f$ is adic\index{morphism of formal schemes@morphism (of formal schemes)!adic morphism of formal schemes@adic ---} if and only if $f\times_XY$ is adic. \hfill$\square$
\end{prop}

Let $P$ be a property of morphisms of schemes that satisfies {\bf (I)} and {\bf (C)} in \S\ref{subsub-genpropertymorphism} and is stable under Zariski topology\index{stable!stable under a topology@--- under a topology} ({\bf \ref{ch-pre}}.\ref{dfn-catequivrelstablearrow} (1)).
Then as in \S\ref{subsub-adicalization}, one can consider the property `adically $P$' for morphisms of adic formal schemes of finite ideal type.
\begin{prop}\label{prop-adicallyPadicallyflatdescent}
Suppose that the property $P$ of morphisms of schemes descends by faithfully flat quasi-compact morphisms of schemes.
Then the property `adically $P$' descends by adically faithfully flat quasi-compact morphisms, that is, the morphism $f$ is adically $P$ if and only if $f\times_XY$ is adically $P$.
\end{prop}

\begin{proof}
By \ref{prop-morphismdescentadic} we may assume that $f$ is adic.
We may work locally on the target $W$ and thus may assume that $X$ has an ideal of definition of finite type $\mathscr{I}$.
Let $f_k\colon Z_k=(Z,\O_Z/\mathscr{I}^{k+1}\O_Z)\rightarrow W_k=(W,\O_W/\mathscr{I}^{k+1}\O_W)$ be the induced morphism of schemes for $k\geq 0$; we use similar notation for other formal schemes and morphisms.
Since $(f\times_XY)_k=f_k\times_{X_k}Y_k$ and since $Y_k\rightarrow X_k$ is faithfully flat and quasi-compact, $f_k$ is $P$ if and only if $(f\times_XY)_k$ is $P$.
\end{proof}

\begin{cor}\label{cor-morphismdescent2}
The morphism $f$ is a closed immersion\index{immersion!closed immersion of formal schemes@closed --- (of formal schemes)} $($resp.\ an open immersion\index{immersion!open immersion of formal schemes@open --- (of formal schemes)}, resp.\ a quasi-compact immersion\index{immersion!immersion of formal schemes@--- (of formal schemes)}$)$ if and only if so is $f\times_XY$.
\end{cor}

\begin{proof}
By \ref{prop-closedimmformal4}, \ref{prop-qsepformal1}, and \ref{prop-immersionimmformal4} these properties are of the form `adically $P$'.
Hence the assertion follows from \ref{prop-adicallyPadicallyflatdescent} and \cite[Expos\'e VIII, 5.5]{SGA1}.
\end{proof}

\begin{cor}\label{cor-morphismdescent3}
The morphism $f$ is affine\index{morphism of formal schemes@morphism (of formal schemes)!affine morphism of formal schemes@affine ---} adic $($resp.\ finite\index{morphism of formal schemes@morphism (of formal schemes)!finite morphism of formal schemes@finite ---}$)$ if and only if so is $f\times_XY$.
\end{cor}

\begin{proof}
Due to \ref{prop-affinemorphism1} and \ref{prop-finitemorphism1} we may apply \ref{prop-adicallyPadicallyflatdescent} in view of \cite[Expos\'e VIII, 5.6, 5.7]{SGA1}.
\end{proof}

\subsubsection{Effective descent}\label{subsub-effectiveetaledescent}
\index{descent!adically flat descent@adically flat ---!effective adically flat descent@effective --- ---|(}
\index{effective!effective adically flat descent@--- adically flat descent|(}
Consider the situation as in \ref{sit-etaledescent}, and let 
$$
p_{23},p_{31},p_{12}\colon\xymatrix{Y\times_XY\times_XY\ar@<.9ex>[r]\ar[r]\ar@<-.9ex>[r]&Y\times_XY}
$$
be the projections, where $p_{ij}$ projects the $i$-th and $j$-th components of $Y\times_XY\times_XY$ onto $Y\times_XY$.
For an adically quasi-coherent sheaf $\mathscr{F}$ on $Y$, a {\em descent datum}\index{descent!descent datum@--- datum} on $\mathscr{F}$ is an isomorphism 
$$
\varphi\colon\widehat{p^{\ast}_1}\mathscr{F}\stackrel{\sim}{\longrightarrow}\widehat{p^{\ast}_2}\mathscr{F}
$$
that satisfies the following $1$-cocycle condition:
$$
\widehat{p^{\ast}_{31}}\varphi=\widehat{p^{\ast}_{23}}\varphi\circ\widehat{p^{\ast}_{12}}\varphi.
$$
Pairs $(\mathscr{F},\varphi)$ consisting of adically quasi-coherent sheaves on $Y$ and descent data form a category in an obvious manner (cf.\ \cite[6.1]{BLR3}).
Notice that for any adically quasi-coherent sheaf $\mathscr{G}$ on $X$, $\widehat{q^{\ast}}\mathscr{G}$ is caonically equipped with the standard descent datum.

\begin{prop}[Effective descent for adically quasi-coherent sheaves]\label{prop-etaledescent3}
The functor $\mathscr{G}\mapsto\widehat{q^{\ast}}\mathscr{G}$ induces a categorical equivalence between the category of adically quasi-coherent sheaves $($resp.\ adically quasi-coherent sheaves of finite type, resp.\ adically quasi-coherent $\O_X$-algebras, resp.\ adically quasi-coherent $\O_X$-algebras of finite type$)$ on $X$ and the category of pairs $(\mathscr{F},\varphi)$ consisting of adically quasi-coherent sheaves $($resp.\ adically quasi-coherent sheaves of finite type, resp.\ adically quasi-coherent $\O_Y$-algebras, resp.\ adically quasi-coherent $\O_Y$-algebras of finite type$)$ on $Y$ and descent data.
\end{prop}

\begin{proof}
By \ref{prop-etaledescent2} we only need to show that the functor is essentially surjective.
By a reduction process similar to that in the proof of \ref{prop-etaledescent2}, we may assume that $X$ has an ideal of definition $\mathscr{I}$ of finite type.
Given $(\mathscr{F},\varphi)$, we set $\mathscr{F}_k=\mathscr{F}/\mathscr{I}^{k+1}\mathscr{F}$ for $k\geq 0$.
Since for each $k$ the descent datum $\varphi$ induces a descent datum $\varphi_k$ on the quasi-coherent sheaf $\mathscr{F}_k$ on the scheme $Y_k$, we have a quasi-coherent sheaf $\mathscr{G}_k$ on $X_k$, uniquely up to isomorphism, such that $\mathscr{F}_k\cong q^{\ast}_k\mathscr{G}_k$ (\cite[Expos\'e VIII, Th\'eor\`eme 1.1]{SGA1}).
By the uniqueness the sheaves $\mathscr{G}_k$ form a projective system $\{\mathscr{G}_k\}_{k\geq 0}$ of abelian sheaves on $X$, and for any $i\leq j$ the transition map $\mathscr{G}_j\rightarrow\mathscr{G}_i$ is surjective with the kernel equal to $\mathscr{I}^{i+1}\mathscr{G}_j$.
Hence by \ref{prop-adicallyqcohinductivelimit} (or \ref{prop-adicallyqcohinductivelimitalgebra}) its limit $\mathscr{G}=\varprojlim_{k\geq 0}\mathscr{G}_k$ is an adically quasi-coherent sheaf on $X$ (resp.\ adically quasi-coherent sheaf on $X$ of finite type, resp.\ an adically quasi-coherent $\O_X$-algebra, resp.\ an adically quasi-coherent $\O_X$-algebra of finite type) such that $\mathscr{G}/\mathscr{I}^{k+1}\mathscr{G}\cong\mathscr{G}_k$.
The last property implies $\widehat{q^{\ast}}\mathscr{G}\cong\mathscr{F}$ by \ref{cor-completepullbackaqcsheaves0}, as desired.
\end{proof}

Descent data are defined also for adic formal schemes.
In the above situation, for an adic formal scheme $W$ adic over $Y$ a {\em descent datum}\index{descent!descent datum@--- datum} on $W$ is an isomorphism
$$
\varphi\colon p^{\ast}_1W\stackrel{\sim}{\longrightarrow}p^{\ast}_2W,
$$
where $p^{\ast}_iW=W\times_{Y,p_i}(Y\times_XY)$, satisfying the similar $1$-cocycle condition.

Let $P$ be a property of morphisms of schemes that satisfies {\bf (I)} and {\bf (C)} in \S\ref{subsub-genpropertymorphism} and is stable under the Zariski topology\index{stable!stable under a topology@--- under a topology} ({\bf \ref{ch-pre}}.\ref{dfn-catequivrelstablearrow} (1)).
\begin{prop}\label{prop-effectiveadicallyflatdescentadicationzation}
Suppose that the subcategory of $\Sch$ consisting of all morphisms satisfying $P$ is an effective descent class\index{descent!effective descent class@effective --- class}\index{effective!effective descent class@--- descent class} {\rm ({\bf \ref{ch-pre}}.\ref{dfn-catequivrelstablearrow})} with respect to the fpqc topology.
Then the subcategory of $\Ac\Fs$ consisting of all adically $P$ morphisms {\rm (\S\ref{subsub-adicalization})} is an effective descent class with respect to adically flat descent.
In other words, in the situation as in {\rm \ref{sit-etaledescent}}, the functor $Z\mapsto q^{\ast}Z=Z\times_XY$ gives rise to a categorical equivalence between the following categories$:$
\begin{itemize}
\item the category of adic formal schemes $Z$ of finite ideal type over $X$ such that $Z\rightarrow X$ is adically $P;$
\item the category of adic formal schemes $W$ of finite ideal type over $Y$ such that $W\rightarrow Y$ is adically $P$ together with descent data.
\end{itemize}
\end{prop}

\begin{proof}
By \ref{cor-morphismdescent1} it suffices to show that the functor is essentially surjective.
We may assume that $X$ has an ideal of definition of finite type $\mathscr{I}$; we work with the notation (such as $X_k$, $Y_k$, etc) as above.
An adic $Y$-formal scheme $W$ as above with a descent datum gives rise to a $Y_k$-scheme $W_k$ with a descent datum for each $k\geq 0$.
Hence for each $k\geq 0$ we have an $X_k$-scheme $Z_k$, uniquely up to canonical isomorphisms, such that $Z_k\rightarrow X_k$ satisfies $P$.
By the uniqueness of $Z_k$'s we have an inductive system $\{Z_k\}_{k\geq 0}$ of schemes over $X$ that satisfies the conditions in Exercise \ref{exer-adicformalschemesinductivelimitadicmaps}, whence the assertion.
\end{proof}

We know that the following properties for morphisms of schemes satisfy effective descent with respect to the fpqc topology:
\begin{itemize}
\item quasi-compact immersion (resp.\ open immersion, resp.\ closed immersion);
\item affine;
\item quasi-affine.
\end{itemize}
\begin{cor}\label{cor-effectiveadicallyflatdescentadicationzation}
The following properties for morphisms of adic formal schemes of finite ideal type satisfy effective descent with respect to the adically flat descent$:$
\begin{itemize}
\item quasi-compact immersion $($resp.\ open immersion, resp.\ closed immersion$);$
\item affine adic$;$
\item adically quasi-affine\index{morphism of formal schemes@morphism (of formal schemes)!adically quasi-affine morphism of formal schemes@adically quasi-affine ---} {\rm (\ref{exas-adicalizationexamples} (2))}. \hfill$\square$
\end{itemize}
\end{cor}
\index{effective!effective adically flat descent@--- adically flat descent|)}
\index{descent!adically flat descent@adically flat ---!effective adically flat descent@effective --- ---|)}

\subsubsection{Adically flat descent and finiteness conditions}\label{subsub-etaledescentfincond}
\begin{prop}\label{prop-etaledescentfincond1}
Consider the situation as in {\rm \ref{sit-etaledescent}}, and suppose that $X$ and $Y$ are locally universally rigid-Noetherian\index{formal scheme!universally rigid-Noetherian formal scheme@universally rigid-Noetherian ---!locally universally rigid-Noetherian formal scheme@locally --- ---} {\rm (\ref{dfn-formalsch})}.

{\rm (1)} Let $\mathscr{F}$ be an adically quasi-coherent sheaf of finite type on $X$.
Then $\mathscr{F}$ is of finite presentation $($resp.\ locally free of finite type$)$ if and only if $\widehat{q^{\ast}}\mathscr{F}$ $(=q^{\ast}\mathscr{F})$ is of finite presentation $($resp.\ locally free of finite type$)$ over $Y$.

{\rm (2)} Let $\mathscr{B}$ be an adically quasi-coherent $\O_X$-algebra.
Then $\mathscr{B}$ is of finite presentation\index{adically quasi-coherent (a.q.c.) algebra!adically quasi-coherent (a.q.c.) algebra of finite presentation@--- of finite presentation} {\rm (\ref{dfn-adicallyqcohalgebraoffinitepres})} if and only if $\widehat{q^{\ast}}\mathscr{B}$ is an adically quasi-coherent $\O_Y$-algebra of finite presentation.
\end{prop}

\begin{proof}
Since the question is local, we may assume that $X$ and $Y$ are affine $X=\Spf A$ and $Y=\Spf B$.
By \ref{prop-adicallyflatformal01} the morphism $q$ is faithfully flat.
Hence by \ref{cor-ffformal2} we deduce that $B$ is faithfully flat over $A$.
Then the assertions follow from \ref{thm-adicqcohpre1}, \ref{prop-adicallyqcohalgebraoffinitepres2}, and \ref{lem-completepullbackaqcsheaves1lem} (cf.\ {\bf \ref{ch-pre}}.\ref{cor-ARconseq2}, {\bf \ref{ch-pre}}.\ref{prop-formalnot31}).
\end{proof}

\begin{prop}\label{prop-etaledescentfincond2}
Consider the situation as in {\rm \ref{sit-etaledescent}}, and let $f\colon W\rightarrow Z$ be a morphism of adic formal schemes of finite ideal type that are adic over $X$.

{\rm (1)} The morphism $f$ is surjective if and only if so is $q^{\ast}f$. 
If $q^{\ast}f$ is injective $($resp.\ bijective$)$, then so is $f$.

{\rm (2)} The morphism $f$ is quasi-compact $($resp.\ of finite type$)$ if and only if so is $q^{\ast}f$.
If $X$ and $Y$ are locally universally rigid-Noetherian, then $f$ is of finite presentation if and only if so is $q^{\ast}f$.
\end{prop}

\begin{proof}
We may assume that $X$ has an ideal of definition.
Then (1) follows from \cite[Expos\'e VIII, 3.1, 3.2]{SGA1}.
(2) follows from \cite[Expos\'e VIII, 3.3, 3.6]{SGA1} combined with \ref{prop-qsepformal1}, \ref{prop-cortopfintype11} and \ref{cor-topfinpres11}.
\end{proof}
\index{descent!adically flat descent@adically flat ---|)}

\subsection{\'Etale topology on adic formal schemes}\label{sub-etaletopology}
\index{topology!etale topology@\'etale ---|(}
\subsubsection{\'Etale sites}\label{subsub-etaleformalsites}
\index{site!etale site@\'etale ---|(}
We consider the category $\Ac\Fs_S$ of adic formal schemes of finite ideal type over a fixed adic formal scheme $S$ of finite ideal type.
Let $\mathscr{E}$ be the subcategory of $\Ac\Fs_S$ consisting of all \'etale morphisms; $\mathscr{E}$ is base-change stable due to \ref{prop-etaleformal1}.

\begin{prop}\label{prop-etaledescentmorphism}
A collection $\{U_{\alpha}\rightarrow U\}_{\alpha\in L}$ of \'etale morphisms in $\Ac\Fs_S$ is universally effectively epimorphic {\rm ({\bf \ref{ch-pre}}.\ref{dfn-UEEF})} if and only if it is a surjective family, that is, the union of the images of $U_{\alpha}\rightarrow U$ coincides with $U$.
\end{prop}

\begin{proof}
The `only if' part is easy to see and is left to the reader.
We prove the `if' part.
Let $W\rightarrow U$ be a morphism in $\Ac\Fs_S$, $V$ an object in $\Ac\Fs_S$, and consider the sequence 
\begin{equation*}
\begin{split}
\Hom_S(W,V)&\xymatrix{{}\ar[r]&{}}{\textstyle \prod_{\alpha\in L}\Hom_S(W\times_UU_{\alpha},V)}\\ 
&\xymatrix{{}\ar@<.5ex>[r]\ar@<-.5ex>[r]&{}}{\textstyle \prod_{\alpha,\beta\in L}\Hom_S(W\times_UU_{\alpha\beta},V)}\rlap{,}
\end{split}
\end{equation*}
where $U_{\alpha\beta}=U_{\alpha}\times_UU_{\beta}$.
We need to show that this sequence is exact.
By the usual Zariski descent\index{descent!Zariski descent@Zariski ---} we may assume that $S$ and $U$ have ideals of definition of finite type.
By an argument similar to that in the proof of \ref{prop-morphismdescent}, we may assume $W$ and $V$ are affine and have ideals of definition of finite ideal type.
Define the filtration on the hom sets by the subsets of the form `$\Hom^{(k)}$' as in the proof of \ref{prop-morphismdescent}, and reduce to the case $k=0$.
Since $\{U_{\alpha,k}\rightarrow U_k\}_{\alpha\in L}$ for each $k$ is an \'etale covering family of schemes, we have the exact sequence
\begin{equation*}
\begin{split}
\Hom_{S_k}(W_k,V_k)&\xymatrix{{}\ar[r]&{}}{\textstyle \prod_{\alpha\in L}\Hom_{S_k}(W_k\times_{U_k}U_{\alpha,k},V_k)}\\ 
&\xymatrix{{}\ar@<.5ex>[r]\ar@<-.5ex>[r]&{}}{\textstyle \prod_{\alpha,\beta\in L}\Hom_{S_k}(W_k\times_{U_k}U_{\alpha\beta,k},V_k)}\rlap{.}
\end{split}
\end{equation*}
We have, for example, $\Hom^{(0)}_S(W,V)\cong\varprojlim_k\Hom_{S_k}(W_k,V_k)$.
We apply $\varprojlim_k$ to the above exact sequences; since $\varprojlim_k$ is left-exact and commutes with products, one gets the desired exact sequence as above.
\end{proof}

\begin{prop}\label{prop-etaletopaxioms}
The subcategory $\mathscr{E}$ of $\Ac\Fs_S$ satisfies the conditions {\boldmath $(\mathbf{S}_1)$}, {\boldmath $(\mathbf{S}_2)$}, {\boldmath $(\mathbf{S}_3{\rm (a)})$}, and {\boldmath $(\mathbf{S}_3{\rm (b)})$} in {\rm {\bf \ref{ch-pre}}, \S\ref{subsub-topologygenconst}}.
\end{prop}

\begin{proof}
It is easy to see that $\mathscr{E}$ satisfies {\boldmath $(\mathbf{S}_1)$} and {\boldmath $(\mathbf{S}_2)$}.
The condition {\boldmath $(\mathbf{S}_3{\rm (a)})$} follows from the second assertion of \ref{prop-etaleformal1} (2).
To show the rest, let $f\colon X\rightarrow Y$ and $g\colon Y\rightarrow Z$ be adic $S$-morphisms of adic formal schemes, and suppose that $g\circ f$ is \'etale.
To show that $f$ is \'etale, we may assume that $Z$ has an ideal of definition of finite type.
Then by \ref{prop-etaleformal00} we can easily reduce to the scheme situation, where the assertion is well-known (\cite[$\mathbf{IV}$, (17.3.5)]{EGA}).
\end{proof}

\begin{dfn}[Large \'etale site]\label{dfn-largeetalesites}{\rm 
The topology on the category $\Ac\Fs_S$ associated to the base-change stable subcategory $\mathscr{E}$ consisting of \'etale morphisms is called the {\em \'etale topology} on $\Ac\Fs_S$, and the resulting site, denoted by $\Ac\Fs_{S,\et}$, is called the {\em large \'etale site}\index{site!etale site@\'etale ---!large etale site@large --- ---} over $S$.}
\end{dfn}

Similarly, one can define the large \'etale site $\Ac\CFs_{S,\et}$ (resp.\ $\Af\Ac\Fs_{S,\et}$) with the underlying category $\Ac\CFs_S$ (resp.\ $\Af\Ac\Fs_S$), the category of coherent (resp.\ affine) adic formal schemes of finite ideal type over $S$.
One can also define similarly the site $\Ac\Fs^{\ast}_{S,\et}$ with the underlying category $\Ac\Fs^{\ast}_S$.

As we saw in {\bf \ref{ch-pre}}.\ref{prop-topologygenconst}, the \'etale topology enjoys the property {\boldmath $(\mathbf{A}_0)$}:
\begin{prop}\label{prop-largeetalesitesheaf}
Any representable presheaf\index{representable!representable (pre)sheaf@--- (pre)sheaf} on $\Ac\Fs_{S,\et}$ is a sheaf. \hfill$\square$
\end{prop}

\begin{prop}\label{prop-etaleeffectivedescentclass}
In the category $\Ac\Fs_S$ of adic formal schemes of finite ideal type over $S$, the following properties of morphisms constitute effective descent classes\index{descent!effective descent class@effective --- class}\index{effective!effective descent class@--- descent class} {\rm (cf.\ {\bf \ref{ch-pre}}.\ref{dfn-catequivrelstablearrow})} under the \'etale topology$:$
\begin{itemize}
\item[{\rm (a)}] affine\index{morphism of formal schemes@morphism (of formal schemes)!affine morphism of formal schemes@affine ---} adic\index{morphism of formal schemes@morphism (of formal schemes)!adic morphism of formal schemes@adic ---}$;$
\item[{\rm (b)}] finite\index{morphism of formal schemes@morphism (of formal schemes)!finite morphism of formal schemes@finite ---}$;$
\item[{\rm (c)}] open immersion$;$
\item[{\rm (d)}] closed immersion\index{immersion!closed immersion of formal schemes@closed --- (of formal schemes)}$;$
\item[{\rm (e)}] immersion\index{immersion!immersion of formal schemes@--- (of formal schemes)}.
\end{itemize}
\end{prop}

\begin{proof}
The stability follows from \ref{cor-morphismdescent2} and \ref{cor-morphismdescent3}.
As for `effective descent class' property, we can reduce to the scheme situation (cf.\ \ref{cor-effectiveadicallyflatdescentadicationzation}), where the assertion is classically known (for the fppf topology); here we use \ref{prop-formalindlimschfact}.
\end{proof}

\begin{prop}\label{prop-etaleeffectivedescentclass2}
In the category $\Ac\Fs^{\ast}_S$ of adic formal schemes over $S$ with adic morphisms, the following properties of morphisms are local on the domain {\rm (cf.\ {\bf \ref{ch-pre}}.\ref{dfn-catequivrelstablearrow})} under the \'etale topology$:$
\begin{itemize}
\item[{\rm (a)}] locally of finite type\index{morphism of formal schemes@morphism (of formal schemes)!morphism of formal schemes locally of finite type@--- locally of finite type}$;$
\item[{\rm (b)}] adically flat\index{morphism of formal schemes@morphism (of formal schemes)!adically flat morphism of formal schemes@adically flat ---}$;$
\item[{\rm (c)}] smooth\index{morphism of formal schemes@morphism (of formal schemes)!smooth morphism of formal schemes@smooth ---}$;$
\item[{\rm (d)}] \'etale\index{morphism of formal schemes@morphism (of formal schemes)!etale morphism of formal schemes@\'etale ---}.
\end{itemize}
\end{prop}

\begin{proof}
By a similar idea as in the proof of \ref{prop-etaleeffectivedescentclass}, one can reduce the situation to the case of schemes, where the assertions are known to be true.
\end{proof}

\begin{cor}\label{cor-etaleeffectivedescentclass2}
The properties for morphisms as in {\rm \ref{prop-etaleeffectivedescentclass2}} and, moreover, 
\begin{itemize}
\item[{\rm (e)}] locally of finite presentation\index{morphism of formal schemes@morphism (of formal schemes)!morphism of formal schemes locally of finite presentation@--- locally of finite presentation}
\end{itemize}
are local on the domain under the \'etale topology in the category $\RNoe\Fs_S$ of locally universally rigid-Noetherian formal schemes. \hfill$\square$
\end{cor}

\begin{dfn}[Small \'etale site]\label{dfn-smalletalesite}{\rm 
For an adic formal scheme $X$ of finite ideal type we denote by $X_{\et}$ the category of \'etale formal schemes over $X$. 
We consider the class $\mathscr{E}_X$ of \'etale morphisms in this category (which obviously satisfies the statements analogous to \ref{prop-etaledescentmorphism} and \ref{prop-etaletopaxioms}), and topologize $X_{\et}$ with the associated topology.
We call the resulting site, denoted by $X_{\et}$, the {\em small \'etale site}\index{site!etale site@\'etale ---!small etale site@small --- ---} of $X$.}
\end{dfn}

We denote by $X^{\sim}_{\et}$ the associated \'etale topos\index{topos!etale topos@\'etale ---}.
As usual, any morphism $f\colon X\rightarrow Y$ of adic formal schemes of finite ideal type induces the inverse image functor $f^{-1}\colon Y_{\et}\rightarrow X_{\et}$ in view of \ref{prop-etaleformal1} (4), which induces a morphism 
$$
f^{\sim}_{\et}=(f_{\ast},f^{-1})\colon X^{\sim}_{\et}\longrightarrow Y^{\sim}_{\et}
$$
of the \'etale topoi.

For an adic formal scheme $X$ of finite ideal type we denote the associated Zariski topos\index{topos!Zariski topos@Zariski ---} by $X^{\sim}_{\Zar}$.
We have the morphism of topoi
$$
\varepsilon=(\varepsilon_{\ast},\varepsilon^{-1})\colon X^{\sim}_{\et}\longrightarrow X^{\sim}_{\Zar}
$$
defined in an obvious manner.

\begin{prop}\label{prop-etaletopologyadicqcoh}
Let $X$ be an adic formal scheme of finite ideal type, and $\mathscr{F}$ an adically quasi-coherent sheaf on $X$ $($resp.\ an adically quasi-coherent $\O_X$-algebra$)$.
Define a presheaf $\mathscr{F}_{\et}$ on $X_{\et}$ as follows$:$ for any object $q\colon Y\rightarrow X$ of $X_{\et}$ we set 
$$
\mathscr{F}_{\et}(Y)=\Gamma(Y,\widehat{q^{\ast}}\mathscr{F}).
$$
Then the presheaf $\mathscr{F}_{\et}$ is a sheaf. \hfill$\square$
\end{prop}

This follows immediately from \ref{prop-etaledescent1}.
Applying the proposition to the structure sheaf $\O_X$, we get the {\em structure sheaf} $\O^{\et}_X$ of the site $X_{\et}$.
Any adic morphism $f\colon X\rightarrow Y$ between adic formal schemes induces a morphism 
$$
f^{\sim}_{\et}=(f_{\ast},f^{\ast})\colon(X^{\sim}_{\et},\O^{\et}_X)\longrightarrow(Y^{\sim}_{\et},\O^{\et}_Y)
$$
of ringed topoi\index{topos!ringed topos@ringed ---}.
Moreover, for any adic formal scheme $X$ of finite ideal type we have a morphism of ringed topoi
$$
\varepsilon=(\varepsilon_{\ast},\varepsilon^{\ast})\colon(X^{\sim}_{\et},\O^{\et}_X)\longrightarrow(X^{\sim}_{\Zar},\O_X).
$$
\index{site!etale site@\'etale ---|)}

\subsubsection{Adically quasi-coherent sheaves on the \'etale site}\label{subsub-etaleadicqcoh}
\index{adically quasi-coherent (a.q.c.) sheaf|(}
Let $X$ be an adic formal scheme of finite ideal type, and $\mathscr{I}$ an ideal of definition of finite type.
As usual, we set $X_k=(X,\O_X/\mathscr{I}^{k+1})$ for $k\geq 0$.
Let $i_k\colon X_k\hookrightarrow X$ be the closed immersion.
Then we have the induced morphism $i_k=(i_{k\ast},i^{\ast}_k)\colon X^{\sim}_{k,\et}\rightarrow X^{\sim}_{\et}$ of ringed topoi.
For any $\O^{\et}_X$-module $\mathscr{F}$ and $k\geq 0$ we define
$$
\mathscr{F}_k=i^{\ast}_k\mathscr{F}.
$$
Since we have the canonical morphism $X^{\sim}_{k,\et}\rightarrow X^{\sim}_{l,\et}$ of ringed topoi for $k\leq l$, we obtain a projective system $\{i_{k\ast}\mathscr{F}_k\}_{k\geq 0}$ of $\O^{\et}_X$-modules.
We set 
$$
\widehat{\mathscr{F}}=\varprojlim_ki_{k\ast}\mathscr{F}_k,
$$
which is again an $\O^{\et}_X$-module.
As in \S\ref{subsub-adicallyqcohdefcompl} one sees that the definition of $\widehat{\mathscr{F}}$ does not depend on the choice of the ideal of definition $\mathscr{I}$, and thus one can define $\widehat{\mathscr{F}}$ for any $\O^{\et}_X$-module even in the case $X$ does not have an ideal of definition.
The $\O^{\et}_X$-module $\widehat{\mathscr{F}}$ thus obtained is called the {\em completion}\index{completion} of $\mathscr{F}$.
The completion comes with the canonical morphism $\mathscr{F}\rightarrow\widehat{\mathscr{F}}$; if this morphism is an isomorphism, we say that $\mathscr{F}$ is {\em complete}\index{complete}.

\begin{dfn}\label{dfn-adicallyqcohetale}{\rm
(1) Let $X$ be an adic formal scheme of finite ideal type, and $\mathscr{F}$ an $\O^{\et}_X$-module.
We say that $\mathscr{F}$ is an {\em adically quasi-coherent sheaf} ({\em on $X_{\et}$}) if the following conditions are satisfied:
\begin{itemize}
\item[{\rm (a)}] $\mathscr{F}$ is complete;
\item[{\rm (b)}] for any object $q\colon U\rightarrow X$ of $X_{\et}$ and any ideal of definition $\mathscr{I}$ of finite type of $U$, the sheaf $\widehat{q^{\ast}}\mathscr{F}/\mathscr{I}\widehat{q^{\ast}}\mathscr{F}$ is a quasi-coherent sheaf with respect to \'etale topology on the scheme $(U,\O_U/\mathscr{I})$.
\end{itemize}

(2) An adically quasi-coherent sheaf $\mathscr{F}$ on $X_{\et}$ is said to be {\em of finite type} if it is of finite type as an $\O_{X_{\et}}$-module.

(3) A morphism between adically quasi-coherent sheaves is a morphism of $\O^{\et}_X$-modules.}
\end{dfn}

By \ref{lem-propadicallyqcohdef1} and effective \'etale descent of quasi-coherent sheaves on schemes, the condition (b) is equivalent to the following one:
\begin{itemize}
\item[${\rm (b)}'$] there exist a covering $\{q_{\alpha}\colon U_{\alpha}\rightarrow X\}_{\alpha\in L}$ of $X_{\et}$ and for each $\alpha\in L$ an ideal of definition $\mathscr{I}_{\alpha}$ of finite type of $U_{\alpha}$ such that for any $\alpha\in L$ and $k\geq 0$ the sheaf $\widehat{q^{\ast}_{\alpha}}\mathscr{F}/\mathscr{I}^{k+1}_{\alpha}\widehat{q^{\ast}_{\alpha}}\mathscr{F}$ is a quasi-coherent sheaf with respect to \'etale topology on the scheme $(U_{\alpha},\O_{U_{\alpha}}/\mathscr{I}^{k+1}_{\alpha})$.
\end{itemize}
If, in particular, $X$ itself has an ideal of definition $\mathscr{I}$ of finite type, then the last condition is equivalent to that (notation being as usual) $\mathscr{F}_k$ is quasi-coherent on $X_{k,\et}$ for any $k\geq 0$.

\begin{lem}\label{lem-adicallyqcohetale}
Let $X$ be an adic formal scheme of finite ideal type.
Then the \'etale sheaf $\mathscr{F}_{\et}$ defined as in {\rm \ref{prop-etaletopologyadicqcoh}} is an adically quasi-coherent sheaf on $X_{\et}$.
\end{lem}

\begin{proof}
It follows easily from the definition of $\mathscr{F}_{\et}$ and the definition of projective limits of sheaves in {\bf \ref{ch-pre}}, \S\ref{subsub-projlimsheaves} that the sheaf $\mathscr{F}_{\et}$ is complete.
By \ref{cor-completepullbackaqcsheaves0} one sees that the other condition is satisfied.
\end{proof}

Thus we have the functor
$$
\mathscr{F}\longmapsto\widehat{\varepsilon^{\ast}}\mathscr{F}=\mathscr{F}_{\et}
$$
from the category of adically quasi-coherent sheaves over $X$ (resp.\ adically quasi-coherent $\O_X$-algebras) to the category of adically quasi-coherent sheaves on $X_{\et}$ (resp.\ adically quasi-coherent $\O^{\et}_X$-algebras).
It is clear that the direct image functor $\varepsilon_{\ast}$ maps adically quasi-coherent sheaves on $X_{\et}$ (resp.\ adically quasi-coherent $\O^{\et}_X$-algebras) to adically quasi-coherent sheaves over $X$ (resp.\ adically quasi-coherent $\O_X$-algebras).
\begin{prop}\label{prop-adicallyqcohetaleequiv}
The functors $\varepsilon_{\ast}$ and $\widehat{\varepsilon^{\ast}}$ give a categorical equivalence between the category of adically quasi-coherent sheaves of finite type over $X$ $($resp.\ adically quasi-coherent $\O_X$-algebras$)$ and the category of adically quasi-coherent sheaves of finite type on $X_{\et}$ $($resp.\ adically quasi-coherent $\O^{\et}_X$-algebras$)$.
\end{prop}

\begin{proof}
We already know by \ref{prop-etaledescent2} that the functor $\widehat{\varepsilon^{\ast}}$ is fully faithful.
Moreover, it is easy to see that, for an adically quasi-coherent sheaf $\mathscr{F}$ of finite type on $X_{\et}$, $\varepsilon_{\ast}\mathscr{F}$ is an adically quasi-coherent sheaf of finite type on $X$.
Hence it suffices to show that for any adically quasi-coherent sheaf $\mathscr{F}$ on $X_{\et}$ (resp.\ adically quasi-coherent $\O^{\et}_X$-algebra) we have a canonical isomorphism $\mathscr{F}\cong\widehat{\varepsilon^{\ast}}\varepsilon_{\ast}\mathscr{F}$.
To this end, since the question is local on $X$, one may assume that $X$ has an ideal of definition $\mathscr{I}$ of finite type.
Moreover, since $\mathscr{F}$ and $\widehat{\varepsilon^{\ast}}\varepsilon_{\ast}\mathscr{F}$ are complete, we only need to show that $\mathscr{F}_k$ and $(\widehat{\varepsilon^{\ast}}\varepsilon_{\ast}\mathscr{F})_k$ are canonically isomorphic for any $k\geq 0$, where $\mathscr{F}_k$ etc.\ are defined as above.
We use the notation as above and denote by $\varepsilon_k=(\varepsilon_{k\ast},\varepsilon^{\ast}_k)\colon X^{\sim}_{k,\et}\rightarrow X^{\sim}_{k,\Zar}$ the morphism of ringed topoi for any $k\geq 0$.
Clearly, we have the commutative diagram of ringed topoi
$$
\xymatrix{X^{\sim}_{k,\et}\ar[d]_{\varepsilon_k}\ar[r]^{i_k}&X^{\sim}_{\et}\ar[d]^{\varepsilon}\\ X^{\sim}_{k,\Zar}\ar[r]_{i_k}&X^{\sim}_{\Zar}\rlap{.}}
$$
Hence we have $(\widehat{\varepsilon^{\ast}}\varepsilon_{\ast}\mathscr{F})_k=i^{\ast}_k\widehat{\varepsilon^{\ast}}\varepsilon_{\ast}\mathscr{F}=\varepsilon^{\ast}_ki^{\ast}_k\varepsilon_{\ast}\mathscr{F}$.
Now by the definition of $\varepsilon_{\ast}$ we see that $i^{\ast}_k\varepsilon_{\ast}\mathscr{F}=\varepsilon_{k\ast}i^{\ast}_k\mathscr{F}$ and hence that $(\widehat{\varepsilon^{\ast}}\varepsilon_{\ast}\mathscr{F})_k=\varepsilon^{\ast}_k\varepsilon_{k\ast}i^{\ast}_k\mathscr{F}$.
By the theory of \'etale descent of quasi-coherent sheaves on schemes (cf.\ \cite[Expos\'e VII, 4]{SGA4-2}), we have $\mathscr{F}_k=i^{\ast}_k\mathscr{F}=\varepsilon^{\ast}_k\varepsilon_{k\ast}i^{\ast}_k\mathscr{F}$, as desired.
\end{proof}

By \ref{prop-etaledescent3} and \ref{prop-etaledescentfincond1} we immediately have:
\begin{prop}\label{prop-adicallyqcohetaleequiv2}
By the functors $\varepsilon_{\ast}$ and $\widehat{\varepsilon^{\ast}}$ adically quasi-coherent sheaves of finite type on $X_{\et}$ $($resp.\ adically quasi-coherent $\O^{\et}_X$-algebras of finite type$)$ correspond to adically quasi-coherent sheaves of finite type on $X$ $($resp.\ adically quasi-coherent $\O_X$-algebras of finite type$)$.
If $X$ is locally universally rigid-Noetherian, then adically quasi-coherent sheaves on $X_{\et}$ of finite presentation $($resp.\ adically quasi-coherent $\O^{\et}_X$-algebras of finite presentation$)$ correspond to adically quasi-coherent sheaves on $X$ of finite presentation $($resp.\ adically quasi-coherent $\O_X$-algebras of finite presentation$)$. \hfill$\square$
\end{prop}
\index{adically quasi-coherent (a.q.c.) sheaf|)}
\index{topology!etale topology@\'etale ---|)}

\subsection{Formal algebraic spaces}\label{sub-formalalgebraicspaces}
\subsubsection{Formal algebraic spaces}\label{subsub-formalalgebraicspacesdef}
\index{formal algebraic space|(}
We consider the large \'etale site $\Ac\Fs_{S,\et}$ consisting of adic formal schemes of finite ideal type over a fixed adic formal scheme $S$ of finite ideal type (cf.\ \ref{dfn-largeetalesites}).
As we have seen in \ref{prop-largeetalesitesheaf} any representable presheaf on the site $\Ac\Fs_{S,\et}$ is a sheaf.
As in {\bf \ref{ch-pre}}.\ref{subsub-catequivrel} we say that a map $\mathscr{F}\rightarrow\mathscr{G}$ of sheaves on $\Ac\Fs_{S,\et}$ is representable\index{representable!representable map of sheaves@--- map of sheaves} if for any object $Z$ of $\Ac\Fs_S$ (regarded as a sheaf on $\Ac\Fs_{S,\et}$) and a map $Z\rightarrow\mathscr{G}$ of sheaves, the fiber product $Z\times_{\mathscr{G}}\mathscr{F}$ is represented by an object of $\Ac\Fs_S$. 
For a base-change stable property $P$ for morphisms between adic formal schemes of finite ideal type, we say that a map $\mathscr{F}\rightarrow\mathscr{G}$ of sheaves is `representable and $P$' if it is representable and for any $Z\rightarrow\mathscr{G}$ as above the induced morphism $Z\times_{\mathscr{G}}\mathscr{F}\rightarrow Z$ satisfies $P$.
Especially, we say that $\mathscr{F}\rightarrow\mathscr{G}$ is an open immersion if it is representable and an open immersion. 

Let $\mathscr{F}$ be a sheaf on $\Ac\Fs_{S,\et}$, and suppose that there exist an object $Y$ of $\Ac\Fs_S$ and a representable \'etale and surjective morphism $q\colon Y\rightarrow\mathscr{F}$.
We say that a collection of morphisms of sheaves $\{\mathscr{F}_{\alpha}\rightarrow\mathscr{F}\}_{\alpha\in L}$ is a Zariski covering of $\mathscr{F}$ if each $\mathscr{F}_{\alpha}\rightarrow\mathscr{F}$ is an open immersion and if $\coprod_{\alpha\in L}\mathscr{F}_{\alpha}\rightarrow\mathscr{F}$ is an epimorphism. 
Since open immersions satisfy effective descent (\ref{prop-etaleeffectivedescentclass}), the last condition is equivalent to that for each $\alpha\in L$ the morphism $\mathscr{F}_{\alpha}\times_{\mathscr{F}}Y\rightarrow Y$ is an open immersion and that $\coprod_{\alpha\in L}\mathscr{F}_{\alpha}\times_{\mathscr{F}}Y\rightarrow Y$ is surjective.
Notice that, in this situation, each $\mathscr{F}_{\alpha}\times_{\mathscr{F}}Y$ is representable, and the map $\mathscr{F}_{\alpha}\times_{\mathscr{F}}Y\rightarrow\mathscr{F}_{\alpha}$ is representable and \'etale surjective.

\begin{dfn}\label{dfn-formalalgebraicspaces1}{\rm 
A {\em formal algebraic space}\footnote{The formal algebraic spaces defined here should be called, more precisely, {\em quasi-separated adic} formal algebraic spaces {\em of finite ideal type}.}\index{algebraic space!formal algebraic space@formal ---} over $S$ is a set-valued sheaf $X$ on the site $\Ac\Fs_{S,\et}$ satisfying the following conditions:
\begin{itemize}
\item[{\rm (a)}] the diagonal morphism $\Delta_X\colon X\rightarrow X\times_SX$ is representable and quasi-compact;
\item[{\rm (b)}] there exist an adic formal scheme $Y$ of finite ideal type over $S$ and a representable \'etale surjective morphism $q\colon Y\rightarrow X$.
\end{itemize}}
\end{dfn}

The morphism $q\colon Y\rightarrow X$ as in (b) will be simply called a {\em representable \'etale covering}\index{representable!representable etale covering@--- \'etale covering} of $X$.
It is obvious that, if $X$ is a formal algebraic space, then any Zariski open subsheaf (that is, a subsheaf $Y\subseteq X$ such that the inclusion map is an open immersion) is again a formal algebraic space, called an {\em open subspace}\index{subspace!open subspace of a formal algebraic space@open --- (of a formal algebraic space)}.

We denote by $\Ac\FAs_S$ the category of formal algebraic spaces over $S$; here morphisms between such formal algebraic spaces are morphisms of sheaves on the site $\Ac\Fs_{S,\et}$.
It is clear by definition that any quasi-separated adic formal scheme of finite ideal type over $S$ is canonically a formal algebraic space over $S$.

Let $X$ be a formal algebraic space over $S$, and $q\colon Y\rightarrow X$ a representable \'etale covering map.
Then $R=Y\times_XY$ is representable, and the projections $p_1,p_2\colon R\rightarrow Y$ are \'etale surjective:
$$
\xymatrix{Y\times_XY=R\ar@<.5ex>[r]^(.7){p_1}\ar@<-.5ex>[r]_(.7){p_2}&Y\ar[r]^(.5){q}&X\rlap{.}}
$$

\begin{prop}[{cf.\ \cite[II.1.3 (a)]{Knu}}]\label{prop-etaleequivalencequot}
{\rm (1)} The morphisms $p_1,p_2\colon R\rightarrow Y$ define an \'etale equivalence relation in $\Ac\Fs_{S,\et}$.

{\rm (2)} The morphism $q\colon Y\rightarrow X$ is the cokernel of the morphisms $p_1,p_2\colon R\rightarrow Y$ in the category of sheaves on $\Ac\Fs_{S,\et}$. \hfill$\square$
\end{prop}

Here an {\em \'etale equivalence relation}\index{equivalence relation!etale equivalence relation@\'etale ---} in $\Ac\Fs_{S,\et}$ means a $\tau$-equivalence relation ({\bf \ref{ch-pre}}.\ref{dfn-tauequivrel}) with $\tau$ equal to the \'etale topology.
Notice that, since 
$$
\xymatrix{R\ar[r]^(.37){(p_1,p_2)}\ar[d]&Y\times_SY\ar[d]^{q\times q}\\ X\ar[r]_(.4){\Delta_X}&X\times_SX}
$$
is Cartesian, it follows that $(p_1,p_2)\colon R\rightarrow Y\times_SY$ is quasi-compact.

\begin{prop}[{cf.\ \cite[II.1.4]{Knu}}]\label{prop-formalalgspmorphism}
Let $X$ and $X'$ be formal algebraic spaces over $S$, and $f\colon X\rightarrow X'$ a morphism over $S$.
Then there exist representable \'etale covering morphisms $q\colon Y\rightarrow X$ and $q'\colon Y'\rightarrow X'$ and a commutative diagram
$$
\xymatrix{Y\times_XY\ar@<.5ex>[r]^(.65){p_1}\ar@<-.5ex>[r]_(.65){p_2}\ar[d]_h&Y\ar[r]^q\ar[d]_g&X\ar[d]^f\\ Y'\times_{X'}Y'\ar@<.5ex>[r]^(.65){p'_1}\ar@<-.5ex>[r]_(.65){p'_2}&Y'\ar[r]_{q'}&X'}
$$
$($the commutativity of the left-hand square means $g\circ p_i=p'_i\circ h$ for $i=1,2)$. \hfill$\square$
\end{prop}

\begin{dfn}\label{dfn-formalalgebraicspaceqcpt}{\rm 
A formal algebraic space $X$ is said to be {\em quasi-compact}\index{formal algebraic space!quasi-compact formal algebraic space@quasi-compact ---}\index{quasi-compact!quasi-compact formal algebraic space@--- (formal algebraic) space} $($or {\em coherent}\index{formal algebraic space!coherent formal algebraic space@coherent ---}$)$ if it has a representable \'etale covering $Y\rightarrow X$ with $Y$ quasi-compact.}
\end{dfn}

We denote by $\Ac\CFAs_S$ the full subcategory of $\Ac\FAs_S$ consisting of quasi-compact formal algebraic spaces over $S$.

\subsubsection{Formal algebraic spaces by quotients}\label{subsub-formalalgebraicspacesbyquotients}
Let 
$$
\xymatrix{R\ar@<.5ex>[r]^(.5){p_1}\ar@<-.5ex>[r]_(.5){p_2}&Y}\leqno{(\ast)}
$$
be an \'etale equivalence relation\index{equivalence relation!etale equivalence relation@\'etale ---} (cf.\ {\bf \ref{ch-pre}}.\ref{dfn-tauequivrel}) by adic formal schemes of finite ideal type over $S$ such that the following condition is satisfied:
\begin{itemize}
\item $Y$ is separated and the induced map 
$$
R\stackrel{(p_1,p_2)}{\longrightarrow}Y\times_SY
$$
is quasi-compact.
\end{itemize}

\begin{thm}\label{thm-etaleequivquote}
In the situation as above, let $q\colon Y\rightarrow X$ be the categorical cokernel of $(\ast)$ in the category of sheaves on $\Ac\Fs_{S,\et}$.
Then $X$ is a formal algebraic space over $S$, and $q\colon Y\rightarrow X$ is a representable \'etale covering\index{representable!representable etale covering@--- \'etale covering}.
\end{thm}

In order to show the theorem, we first notice that we may work locally on $Y$ (cf.\ \cite[I.5.7]{Knu}) and hence may assume that $S$ is coherent and that $Y$ is coherent over $S$.

\begin{prop}\label{prop-existenceidealofdefinitiongabber}
In the situation as above, suppose that $S$ is coherent and that $Y$ is coherent over $S$ {\rm (\ref{dfn-cohformalmorphism})}\index{formal scheme!coherent formal scheme@coherent ---}.
Then there exists an ideal of definition $\mathscr{I}\subseteq\O_Y$ of finite type such that 
$$
p_1^{-1}\mathscr{I}\O_R=p_2^{-1}\mathscr{I}\O_R.
$$
\end{prop}

To show the proposition we need:
\begin{lem}\label{lem-schemetheoreticimage}
Consider a Cartesian diagram
$$
\xymatrix{Z\ar[d]_f&V\ar[l]_p\ar[d]^g\\ X&U\ar[l]^q}
$$
of coherent adic formal schemes of finite ideal type and adic morphisms.
Suppose$:$
\begin{itemize}
\item[{\rm (a)}] $p$ and $q$ are \'etale, and $f$ and $g$ are coherent$;$
\item[{\rm (b)}] $Z$ and $V$ are schemes $($that is, $0$-adic formal schemes$).$
\end{itemize}
Then $\mathscr{K}=\ker(\O_X\rightarrow f_{\ast}\O_Z)$ is an ideal of definition of $X$, and we have 
$$
q^{-1}\mathscr{K}\O_U=\ker(\O_U\rightarrow g_{\ast}\O_V).
$$
\end{lem}

\begin{proof}
Let $\mathscr{I}$ be an ideal of definition of finite type on $X$ (\ref{cor-extension2}).
Since $f$ is adic and $Z$ is a quasi-compact scheme, there exists $n\geq 1$ such that $\mathscr{I}^n\O_Z=0$.
Thus, replacing $\mathscr{I}$ by $\mathscr{I}^n$, we may assume that $\mathscr{I}\subseteq\mathscr{K}$.

Let us show that $\mathscr{K}$ is an ideal of definition; to this end, we want to show that $\mathscr{K}$ is adically quasi-coherent.
Since $\mathscr{K}$ contains $\mathscr{I}$, it suffices to show that for any $k\geq 0$ the induced sheaf $\mathscr{K}_k=\mathscr{K}/\mathscr{I}^{k+1}$ is a quasi-coherent sheaf on the scheme $X_k=(X,\O_X/\mathscr{I}^{k+1})$.
But this follows from that the morphism $f_k\colon Z\rightarrow X_k$ induced from $f$ is coherent and that 
$$
\mathscr{K}_k=\ker(\O_{X_k}\rightarrow f_{k\ast}\O_Z).
$$

To show the other assertion, set $U_k=(U,\O_U/\mathscr{I}^{k+1}\O_U)$ for $k\geq 0$, and let $q_k\colon U_k\rightarrow X_k$ and $g_k\colon V\rightarrow U_k$ be the morphisms induced from $q$ and $g$, respectively.
Since both $q^{-1}\mathscr{K}\O_U$ and $\ker(\O_U\rightarrow g_{\ast}\O_V)$ contain $\mathscr{I}\O_U$ (and hence are complete), it suffices to check the equality
$$
q^{-1}_k\mathscr{K}_k\O_{U_k}=\ker(\O_{U_k}\rightarrow g_{k\ast}\O_V)
$$
for any $k\geq 0$.
But this is clear since $p$ and $q_k$ are \'etale morphisms of schemes and $V\cong Z\times_{X_k}U_k$.
\end{proof}

\begin{proof}[Proof of Proposition {\rm \ref{prop-existenceidealofdefinitiongabber}}]
We first claim the following:

\medskip
{\sc Claim.} {\it It suffices to find an ideal of definition $\mathscr{I}\subseteq\O_Y$, not necessarily of finite type, such that $p_1^{-1}\mathscr{I}\O_R=p_2^{-1}\mathscr{I}\O_R$.}

\medskip
Set $Y_k=(Y,\O_Y/\mathscr{I}^{k+1})$ and $R_k=(R,\O_R/\mathscr{I}^{k+1}\O_R)$) for $k\geq 0$.
We have for each $k\geq 0$ the induced morphisms
$$
\xymatrix{R_k\ar@<.5ex>[r]^(.5){p_{1,k}}\ar@<-.5ex>[r]_(.5){p_{2,k}}&Y_k}
$$
of schemes.
By \ref{prop-etaleformal00} and Exercise \ref{exer-equivalencerelationinvolution} one sees easily that this diagram gives an \'etale equivalence relation of schemes. 
Let $X_k$ be the resulting coherent algebraic space with the quotient map $q_k\colon Y_k\rightarrow X_k$.
Set $\O^{\et}_X=\varprojlim_{k\geq 0}\O^{\et}_{X_k}$, and let $\mathscr{I}^{(k)}_X$ be the kernel of $\O^{\et}_X\rightarrow\O^{\et}_{X_k}$ for each $k\geq 0$.
By \'etale descent the sheaf $\mathscr{I}^{(l)}_X/\mathscr{I}^{(k)}_X$ for any $k\geq l\geq 0$ is a quasi-coherent sheaf on $X_k$ such that $q^{\ast}_k\mathscr{I}^{(l)}_X/\mathscr{I}^{(k)}_X=\mathscr{I}^l/\mathscr{I}^k$.
By {\bf \ref{ch-pre}}.\ref{thm-grusonraynaudlimit} $\mathscr{I}^{(0)}_X/\mathscr{I}^{(1)}_X$ is the filtered inductive limit of quasi-coherent ideals of finite type on $X_1$:
$$
\mathscr{I}^{(0)}_X/\mathscr{I}^{(1)}_X=\varinjlim_{\lambda\in\Lambda}\ovl{\mathscr{I}}_{\lambda}.
$$
For each $\lambda\in\Lambda$, let $\mathscr{I}_{\lambda}$ be the pull-back of $\ovl{\mathscr{I}}_{\lambda}$ by the canonical morphism $\O^{\et}_X\rightarrow\O^{\et}_{X_1}$.
Then $\mathscr{I}_{\lambda}\O^{\et}_Y$ is an adically quasi-coherent ideal of $\O^{\et}_Y$ (\ref{prop-exeradicallyqcohsheavesbypullback}).
Take an ideal of definition $\mathscr{I}'\subseteq\O_Y$ of finite type (\ref{cor-extension2}) contained in $\mathscr{I}$.
Then there exists $\lambda\in\Lambda$ such that $\mathscr{I}'\O^{\et}_Y\subseteq\mathscr{I}_{\lambda}\O^{\et}_Y$.
For a sufficiently large $k\gg 0$ we have the inclusions:
$$
\mathscr{I}^{(k+1)}_X\O^{\et}_Y\subseteq\mathscr{I}'\O^{\et}_Y\subseteq\mathscr{I}_{\lambda}\O^{\et}_Y\subseteq\mathscr{I}\O^{\et}_Y.
$$
Since $\mathscr{I}_{\lambda}\O^{\et}_Y/\mathscr{I}^{(k+1)}_X\O^{\et}_Y$ ($=\ovl{\mathscr{I}}_{\lambda}\O^{\et}_{Y_k}$) is of finite type, $\mathscr{I}_{\lambda}\O^{\et}_Y/\mathscr{I}'\O^{\et}_Y$ is of finite type; but since $\mathscr{I}'\O^{\et}_Y$ is of finite type, we deduce that $\mathscr{I}_{\lambda}\O^{\et}_Y$ is of finite type.
Now the ideal sheaf $\mathscr{I}''=\varepsilon_{\ast}\mathscr{I}_{\lambda}\O^{\et}_Y$ with respect to the Zariski topology is an ideal of definition of finite type on $Y$ such that $p_1^{-1}\mathscr{I}''\O_R=p_2^{-1}\mathscr{I}''\O_R$, whence the claim.

\medskip
To proceed, we consider the diagram (cf.\ Exercise \ref{exer-equivalencerelationinvolution}): 
$$
\xymatrix{Y&R\ar[l]_{p_1}\ar[d]_{p_2}&T\ar@<.5ex>[l]^(.5){p_{13}}\ar@<-.5ex>[l]_(.5){p_{12}}\ar[d]^{p_{23}}\\ &Y&R\ar@<.5ex>[l]^(.5){p_2}\ar@<-.5ex>[l]_(.5){p_1}}
$$
consisting of \'etale covering morphisms, where $T=R\times_{p_1,Y,p_2}R$; here we have $p_1\circ p_{23}=p_2\circ p_{12}$ and $p_2\circ p_{23}=p_2\circ p_{13}$.

Now we are going to construct an ideal of definition $\mathscr{I}$ on $Y$ satisfying $p_1^{-1}\mathscr{I}\O_R=p_2^{-1}\mathscr{I}\O_R$.
Let $\mathscr{J}\subseteq\O_Y$ be an ideal of definition (\ref{cor-extension2}), and define
$$
\mathscr{I}=\ker(\O_Y\rightarrow p_{2\ast}\O_R/p^{-1}_1\mathscr{J}\O_R).
$$
Set $R_0=(R,\O_R/p^{-1}_1\mathscr{J}\O_R)$ and $T_0=(T,\O_T/p^{-1}_{12}p^{-1}_1\mathscr{J}\O_T)$; notice the equality $p^{-1}_{12}p^{-1}_1\mathscr{J}\O_T=p^{-1}_{13}p^{-1}_1\mathscr{J}\O_T$.
We have the double Cartesian diagram
$$
\xymatrix{R_0\ar[d]_{p_2}&T_0\ar@<.5ex>[l]^(.5){p_{13}}\ar@<-.5ex>[l]_(.5){p_{12}}\ar[d]^{p_{23}}\\ Y&R\ar@<.5ex>[l]^(.5){p_2}\ar@<-.5ex>[l]_(.5){p_1}}
$$
where all horizontal arrows are \'etale.
Since $\mathscr{I}=\ker(\O_Y\rightarrow p_{2\ast}\O_{R_0})$, we can apply \ref{lem-schemetheoreticimage}.
It follows that $\mathscr{I}$ is an ideal of definition on $Y$ and that 
$$
p^{-1}_1\mathscr{I}\O_R=\ker(\O_R\rightarrow p_{23\ast}\O_{T_0})=p^{-1}_2\mathscr{I}\O_R,
$$
which finishes the proof of the proposition.
\end{proof}

\begin{proof}[Proof of Theorem {\rm \ref{thm-etaleequivquote}}]
We use the notation as in \ref{prop-existenceidealofdefinitiongabber}.
Set $Y_k=(Y,\O_Y/\mathscr{I}^{k+1})$ and $R_k=(R,\O_R/\mathscr{I}^{k+1}\O_R)$ for $k\geq 0$.
Suppose $S$ has an ideal of definition of finite type $\mathscr{J}\subseteq\O_S$.
As we may assume that $\mathscr{J}\O_Y\subseteq\mathscr{I}$, the schemes $Y_k$ and $R_k$ are considered to be $S_k$-schemes, where $S_k=(S,\O_S/\mathscr{J}^{k+1})$, for $k\geq 0$.
We have the induced diagram 
$$
\xymatrix{R_k\ar@<.5ex>[r]^(.5){p_{1,k}}\ar@<-.5ex>[r]_(.5){p_{2,k}}&Y_k}\leqno{(\ast)_k}
$$
for each $k\geq 0$.
It follows from \ref{prop-etaleformal00} and Exercise \ref{exer-equivalencerelationinvolution} that $(\ast)_k$ gives an \'etale equivalence relation by $S_k$-schemes such that the induced map 
$$
R_k\stackrel{(p_{1,k},p_{2,k})}{\longrightarrow}Y_k\times_{S_k}Y_k
$$
is quasi-compact.
Then by \cite[II.1.3 (b)]{Knu} the categorical cokernel $q_k\colon Y_k\rightarrow X_k$ of $(\ast)_k$ is a representable \'etale covering of an algebraic space $X_k$ over $S_k$.

It follows from \cite[I.5.12]{Knu} that the map $(p_{1,k},p_{2,k})\colon R_k\rightarrow Y_k\times_{S_k}Y_k$ is quasi-affine for any $k\geq 0$ and hence that the map $R\stackrel{(p_1,p_2)}{\rightarrow}Y\times_SY$ is adically quasi-affine\index{morphism of formal schemes@morphism (of formal schemes)!adically quasi-affine morphism of formal schemes@adically quasi-affine ---} (\ref{exas-adicalizationexamples} (2)).
Then the assertion is the formal consequence of {\bf \ref{ch-pre}}.\ref{prop-fundamentalquotient} in view of the fact that adically quasi-affine morphisms form an effective descent class (\ref{prop-effectiveadicallyflatdescentadicationzation}).
\end{proof}

\begin{rem}\label{rem-thmetaleequivquote}{\rm 
(1) Note that in the situation as in \ref{thm-etaleequivquote} the diagonal morphism $\Delta\colon X\rightarrow X\times_SX$ is representable and adically quasi-affine.

(2) Notice also that in the situation as in the proof of \ref{thm-etaleequivquote} we have the natural isomorphism 
$$
X\cong\varinjlim_{k\geq 0}X_k
$$
of sheaves on $\Ac\Fs_{S,\et}$.}
\end{rem}
\index{formal algebraic space|)}

\subsubsection{Fiber products}\label{subsub-fiberproductsformalalgsp}
\index{fiber product!fiber product of formal algebraic spaces@--- (of formal algebraic spaces)|(}
\index{formal algebraic space!fiber product of formal algebraic spaces@fiber product of ---s|(}
\begin{prop}\label{prop-fiberprodformalalgsp}
Let $X$, $Y$, and $Z$ be formal algebraic spaces over $S$, and $X\rightarrow Z\leftarrow Y$ morphisms over $S$.
Then the sheaf fiber product $X\times_ZY$ is a formal algebraic space over $S$, and the diagram
$$
\xymatrix{X\times_ZY\ar[d]\ar[r]&X\ar[d]\\ Y\ar[r]&Z}
$$
is Cartesian in $\Ac\FAs_S$.
\end{prop}

\begin{proof}
We first construct the fiber products locally and then globalize by patching; notice that, since $X$, $Y$, and $Z$ are quasi-separated, it follows that the diagonal map for the sheaf $X\times_ZY$ is automatically representable and quasi-compact (indeed, $\Delta_{X\times_ZY}\cong\Delta_X\times_{\Delta_Z}\Delta_Y$).
Hence, in view of \ref{thm-etaleequivquote} the result follows by an argument similar to the proof of \cite[II.1.5]{Knu}.
\end{proof}
\index{formal algebraic space!fiber product of formal algebraic spaces@fiber product of ---s|)}
\index{fiber product!fiber product of formal algebraic spaces@--- (of formal algebraic spaces)|)}

\subsubsection{\'Etale topology on formal algebraic spaces}\label{subsub-etaletopformalalgsp}
\begin{dfn}\label{dfn-etaletopformalalgsp}{\rm 
We say that an $S$-morphism $f\colon X\rightarrow Y$ of formal algebraic spaces is {\em \'etale}\index{morphism of formal algebraic spaces@morphism (of formal algebraic spaces)!etale morphism of formal algebraic spaces@\'etale ---} if there exist a representable \'etale covering $V\rightarrow Y$ and a representable \'etale covering $U\rightarrow X\times_YV$ such that the resulting map $U\rightarrow V$ of adic formal schemes of finite ideal type is \'etale.}
\end{dfn}

\begin{prop}\label{prop-etaleformalalgsp1}
{\rm (1)} The composition of two \'etale morphisms is \'etale.

{\rm (2)} For any \'etale $Z$-morphisms $f\colon X\rightarrow Y$ and $g\colon X'\rightarrow Y'$ of formal algebraic spaces over a formal algebraic space $Z$, the induced morphism $f\times_Zg\colon X\times_ZY\rightarrow X'\times_ZY'$ is \'etale.

{\rm (3)} For any \'etale $Z$-morphism $f\colon X\rightarrow Y$ of formal algebraic spaces over a formal algebraic space $Z$ and for any morphism $Z'\rightarrow Z$ of formal algebraic spaces, the induced morphism $f_{Z'}\colon X\times_ZZ'\rightarrow Y\times_ZZ'$ is \'etale.
\end{prop}

\begin{proof}
(1) is easy to see.
To show (2) and (3), by {\bf \ref{ch-pre}}.\ref{prop-basechangestable} it suffices to show the special case of (3) with $Z=Y$, which is straightforward.
\end{proof}

The following proposition is straightforward (use \ref{prop-etaletopaxioms}):
\begin{prop}\label{prop-etaletopaxiomsforalgsp}
The base-change stable subcategory $\mathscr{E}$ of $\Ac\FAs_S$ consisting of all \'etale morphisms satisfies the conditions {\boldmath $(\mathbf{S}_1)$}, {\boldmath $(\mathbf{S}_2)$}, {\boldmath $(\mathbf{S}_3{\rm (a)})$}, and {\boldmath $(\mathbf{S}_3{\rm (b)})$} in {\rm {\bf \ref{ch-pre}}, \S\ref{subsub-topologygenconst}}. \hfill$\square$
\end{prop}

\begin{dfn}[Large \'etale site]\label{dfn-largeetalesitesforalgsp}{\rm 
The topology on the category $\Ac\FAs_S$ associated to the subcategory $\mathscr{E}$ consisting of \'etale maps is called the {\em \'etale topology} on $\Ac\FAs_S$, and the resulting site, denoted by $\Ac\FAs_{S,\et}$, is called the {\em large \'etale site}\index{site!etale site@\'etale ---!large etale site@large --- ---} over $S$.}
\end{dfn}

Notice that for any formal algebraic space $X$ a representable \'etale covering $Y\rightarrow X$ is a covering map with respect to the \'etale topology.
In particular, in the site $\Ac\FAs_{S,\et}$, covering families consisting of morphisms from {\em adic formal schemes} (of finite ideal type) are cofinal in the set of all covering families.
The following proposition follows from {\bf \ref{ch-pre}}.\ref{prop-topologygenconst}.
\begin{prop}\label{prop-largeetalesitesheafforalgsp}
Any representable presheaf\index{representable!representable (pre)sheaf@--- (pre)sheaf} on $\Ac\FAs_{S,\et}$ is a sheaf. \hfill$\square$
\end{prop}

\begin{dfn}[Small \'etale site]\label{dfn-smalletalesiteforalgsp}{\rm 
For a formal algebraic space $X$ we denote by $X_{\et}$ the category of \'etale formal algebraic spaces over $X$. 
We consider the class $\mathscr{E}$ of \'etale maps in this category (which obviously satisfies the conditions as in \ref{prop-etaletopaxiomsforalgsp}) and topologize $X_{\et}$ with the associated topology.
We call the resulting site, denoted by $X_{\et}$, the {\em small \'etale site}\index{site!etale site@\'etale ---!small etale site@small --- ---} over $X$.}
\end{dfn}

We denote by $X^{\sim}_{\et}$ the associated \'etale topos.
Any $S$-morphism $f\colon X\rightarrow Y$ of formal algebraic spaces induces the inverse image functor $f^{-1}\colon Y_{\et}\rightarrow X_{\et}$, which gives rise to a morphism 
$$
f^{\sim}_{\et}=(f_{\ast},f^{-1})\colon X^{\sim}_{\et}\longrightarrow Y^{\sim}_{\et}
$$
of the \'etale topoi.

For a formal algebraic space $X$ we define the {\em structure sheaf} of $X$ as follows: $\O_X(U)=\Gamma(U,\O_U)$ for any \'etale map $U\rightarrow X$ from an adic formal scheme $U$; this defines a sheaf due to \ref{prop-etaledescent1}.

\subsubsection{Ideal of definition and adic morphisms}\label{subsub-idealofdefinitionforalgsp}
Let $X$ be a formal algebraic space, and $q\colon Y\rightarrow X$ a representable \'etale covering\index{representable!representable etale covering@--- \'etale covering}.
For any $\O_X$-module $\mathscr{F}$ one can define the complete pull-back\index{complete!complete pull-back@--- pull-back} $\widehat{q^{\ast}}\mathscr{F}$ by an obvious manner (cf.\ \S\ref{sub-completepullbackaqcsheaves}).

\begin{dfn}\label{dfn-idealofdefinitionforalgsp}{\rm 
Let $X$ be a formal algebraic space.
An {\em ideal of definition of finite type}\index{ideal of definition} of $X$ is an ideal sheaf $\mathscr{I}\subseteq\O_X$ of finite type such that the following condition is satisfied: for any representable \'etale covering $q\colon Y\rightarrow X$ we have $\widehat{f^{\ast}}\mathscr{I}=\mathscr{I}\O_Y$ (cf.\ \ref{prop-exeridealofdefinitioncompletepullback}) and is an ideal of definition of finite type on $Y$.}
\end{dfn}

By Exercise \ref{exer-idealofdefinitioneffectivedescent} the condition is equivalent to the following: for any representable \'etale covering $q\colon Y\rightarrow X$ we have $\widehat{f^{\ast}}\mathscr{I}=\mathscr{I}\O_Y$, which is an adically quasi-coherent sheaf of finite type on $Y$ and, for at least one such $q\colon Y\rightarrow X$, $\mathscr{I}\O_Y$ is an ideal of definition of finite type.
By \ref{prop-existenceidealofdefinitiongabber} we deduce:
\begin{prop}\label{prop-idealofdefinitionforalgsp}
There exists a Zariski covering $\{X_{\alpha}\rightarrow X\}_{\alpha\in L}$ such that each $X_{\alpha}$ has an ideal of definition of finite type. \hfill$\square$
\end{prop}

\begin{dfn}\label{dfn-adicmapforalgsp}{\rm 
A morphism $f\colon X\rightarrow Y$ of formal algebraic spaces is said to be {\em adic}\index{morphism of formal algebraic spaces@morphism (of formal algebraic spaces)!adic morphism of formal algebraic spaces@adic ---} if there exists a Zariski covering $\{Y_{\alpha}\rightarrow Y\}_{\alpha\in L}$ of $Y$ and for each $\alpha\in L$ an ideal of definition of finite type $\mathscr{I}_{\alpha}$ on $Y_{\alpha}$ such that $\mathscr{I}_{\alpha}\O_{X_{\alpha}}$ is an ideal of definition of $X_{\alpha}=X\times_YY_{\alpha}$.}
\end{dfn}

\begin{prop}\label{prop-dfnadicmorphismformalalgebraicspaces}
Let $f\colon X\rightarrow Y$ be a morphism of formal algebraic spaces.
Then the following conditions are equivalent$:$
\begin{itemize}
\item[{\rm (a)}] $f$ is adic$;$
\item[{\rm (b)}] for any representable \'etale covering $V\rightarrow Y$ and any representable \'etale covering $U\rightarrow X\times_YV$, the morphism $U\rightarrow V$ is adic$;$
\item[{\rm (c)}] there exists a representable \'etale covering $V\rightarrow Y$ and a representable \'etale covering $U\rightarrow X\times_YV$ such that the morphism $U\rightarrow V$ is adic. \hfill$\square$
\end{itemize}
\end{prop}

This follows easily from \ref{prop-adicnessadicallyflatdescent}.
\begin{prop}\label{prop-adicmorforalgsp1}
{\rm (1)} Let $f\colon X\rightarrow Y$ and $g\colon Y\rightarrow Z$ be morphisms of formal algebraic spaces.
If $f$ and $g$ are adic, then so is the composition $g\circ f$.
If $g\circ f$ and $g$ are adic, then so is $f$.

{\rm (2)} Let $Z$ be a formal algebraic space, and $f\colon X\rightarrow X'$ and $g\colon Y\rightarrow Y'$ two adic $Z$-morphisms of formal algebraic spaces.
Then $f\times_Zg\colon X\times_ZY\rightarrow X'\times_ZY'$ is adic.

{\rm (3)} Let $Z$ be a formal algebraic space, and $f\colon X\rightarrow Y$ an adic $Z$-morphism between formal algebraic spaces. Then for any morphism $Z'\rightarrow Z$ of formal algebraic spaces the induces morphism $f_{Z'}\colon X\times_ZZ'\rightarrow Y\times_ZZ'$ is adic.
\end{prop}

\begin{proof}
All assertions except for the second one in (1) follow easily from \ref{prop-adicmor1}.
To show the second assertion of (1), take a representable \'etale covering $U\rightarrow Z$, and let $f'\colon U\times_ZX\rightarrow U\times_ZY$ and $g'\colon U\times_ZY\rightarrow U$ be the induced morphisms of adic formal schemes of finite ideal type.
Suppose $g\circ f$ and $g$ are adic.
Then $g'\circ f'$ and $g'$ are adic.
For $f$ to be adic, it is enough by \ref{prop-adicnessadicallyflatdescent} that $f'$ is adic, which follows from \ref{prop-adicmor1} (1).
\end{proof}

We denote by $\Ac\FAs^{\ast}_S$ (resp.\ $\Ac\FAs^{\ast}_{/S}$) the category of formal algebraic spaces over $S$ with adic morphisms (resp.\ the category of formal algebraic spaces adic over $S$).

Suppose that a formal algebraic space $X$ has an ideal of definition $\mathscr{I}$, and take a representable \'etale covering $q\colon Y\rightarrow X$; replacing $Y$ by a disjoint union of affine subsets of $Y$, we may assume that $Y$ is separated.
Then $Y$ is an adic formal scheme on which the ideal $\mathscr{I}\O_Y$ is an ideal of definition.
Set $R=Y\times_XY$ and consider the resulting diagram
$$
\xymatrix{X&Y\ar[l]_(.45){q}&R.\ar@<-.5ex>[l]_(.45){p_1}\ar@<.5ex>[l]^(.45){p_2}}
$$
Then we notice that we are in a similar situation as in \S\ref{subsub-formalalgebraicspacesbyquotients}.
In particular, we have 
$$
X\cong\varinjlim_{k\geq 0}X_k,
$$
where $X_k$ ($k\geq 0$) is the algebraic space defined as in \S\ref{subsub-formalalgebraicspacesbyquotients}.

\begin{prop}\label{prop-localcriterionrepresentability}
Let $X$ be a coherent formal algebraic space having an ideal of definition $\mathscr{I}$ of finite type, and $X_k$ for $k\geq 0$ the algebraic spaces defined as above.
Then the following conditions are equivalent$:$
\begin{itemize}
\item[{\rm (a)}] $X$ is represented by an adic formal scheme of finite ideal type$;$
\item[{\rm (b)}] $X_k$ is represented by a scheme for any $k\geq 0;$
\item[{\rm (c)}] $X_0$ is represented by a scheme.
\end{itemize}
\end{prop}

\begin{proof}
(a) $\Rightarrow$ (b) is trivial.
(b) $\Rightarrow$ (c) follows from {\bf \ref{ch-pre}}.\ref{cor-KnudsonSerrecriterion}.
Suppose (b) holds.
As we have seen above, $X$ is the inductive limit of $\{X_k\}_{k\geq 0}$ taken in the category of formal algebraic spaces; but by \ref{prop-formalindlimschadic} it is an adic formal scheme.
\end{proof}

\subsubsection{Formal completion of algebraic spaces}\label{subsub-formalcompletionalgsp}
\index{completion!formal completion@formal ---|(}
Let $X$ be an algebraic space over a scheme $S$, and $Z$ a closed subspace of $X$ of finite presentation with the defining ideal sheaf $\mathscr{I}$ (cf.\ \cite[II.5]{Knu}).
Take a representable \'etale covering $q\colon Y\rightarrow X$, which yields the cokernel sequence:
$$
\xymatrix{X&Y\ar[l]_(.4){q}&Y\times_XY\ar@<-.5ex>[l]_(.55){p_1}\ar@<.5ex>[l]^(.55){p_2},}
$$
where $p_1$ and $p_2$ are the projections; we set $R=Y\times_XY$.
Consider the formal completions $\widehat{Y}|_Z$ and $\widehat{R}|_Z$ along the closed subschemes $q^{-1}(Z)$ and $p_i^{-1}q^{-1}(Z)$, respectively.
We get a diagram of adic formal schemes of finite ideal type
$$
\xymatrix{\widehat{Y}|_Z&\widehat{R}|_Z\ar@<-.5ex>[l]_(.45){\widehat{p}_1}\ar@<.5ex>[l]^(.45){\widehat{p}_2}}\leqno{(\ast)}
$$
consisting of \'etale surjective morphisms.
\begin{prop}\label{prop-formalcompletionalgsp1}
The diagram $(\ast)$ gives an \'etale equivalence relation in $\Ac\Fs_S$, and if we put $\widehat{q}\colon\widehat{Y}|_Z\rightarrow\widehat{X}|_Z$ to be the sheaf quotient, then $\widehat{X}|_Z$ is a formal algebraic space over $S$ and $\widehat{q}$ gives a representable \'etale covering.
\end{prop}

\begin{proof}
The first assertion is clear, and the other assertions follow from \ref{thm-etaleequivquote}.
\end{proof}
\index{completion!formal completion@formal ---|)}

\subsubsection{Adically quasi-coherent sheaves on formal algebraic spaces}\label{subsub-etaleadicqcohonformalalgebraicspaces}
\index{adically quasi-coherent (a.q.c.) sheaf|(}
Let $X$ be a formal algebraic space, and $\mathscr{I}$ an ideal of definition of finite type.
As in \S\ref{subsub-idealofdefinitionforalgsp} we have algebraic spaces $X_k$ for $k\geq 0$, which we often denote loosely by $X_k=(X,\O_X/\mathscr{I}^{k+1})$.
For any $\O_X$-module $\mathscr{F}$ we define $\mathscr{F}_k=\mathscr{F}/\mathscr{I}^{k+1}\mathscr{F}$ for $k\geq 0$, regarded as an $\O_{X_k}$-module.
We thus obtain a projective system $\{\mathscr{F}_k\}_{k\geq 0}$ of $\O_X$-modules.
Define an $\O_X$-module $\widehat{\mathscr{F}}$ by $\widehat{\mathscr{F}}=\varprojlim_k\mathscr{F}_k$.
As in \S\ref{subsub-adicallyqcohdefcompl} one sees that the definition of $\widehat{\mathscr{F}}$ does not depend on the choice of the ideal of definition $\mathscr{I}$, and thus one can define $\widehat{\mathscr{F}}$ for any $\O_X$-module even in the case $X$ does not have an ideal of definition.
The $\O_X$-module $\widehat{\mathscr{F}}$ thus obtained is called the {\em completion}\index{completion} of $\mathscr{F}$.
As before, the completion comes with the canonical morphism $\mathscr{F}\rightarrow\widehat{\mathscr{F}}$, and if this morphism is an isomorphism, we say that $\mathscr{F}$ is {\em complete}\index{complete}.

\begin{dfn}\label{dfn-adicallyqcohetaleonformalalgebraicspace}{\rm
(1) We say that an $\O_X$-module $\mathscr{F}$ is an {\em adically quasi-coherent sheaf} if the following conditions are satisfied:
\begin{itemize}
\item[{\rm (a)}] $\mathscr{F}$ is complete;
\item[{\rm (b)}] for any subspace $U\subseteq X$ and any ideal of definition $\mathscr{I}$ of finite type of $U$, the sheaf $(\mathscr{F}|_U)/\mathscr{I}(\mathscr{F}|_U)$ is a quasi-coherent sheaf on the algebraic space $(U,\O_U/\mathscr{I})$.
\end{itemize}

(2) An adically quasi-coherent sheaf $\mathscr{F}$ on $X$ is said to be {\em of finite type} if it is of finite type as an $\O_X$-module.

(3) A morphism between adically quasi-coherent sheaves is a morphism of $\O_X$-modules.}
\end{dfn}

To check the condition (b), it is enough to check the following:
\begin{itemize}
\item[${\rm (b)}'$] there exist a Zariski covering $\{X_{\alpha}\rightarrow X\}_{\alpha\in L}$ of $X$ and for each $\alpha\in L$ an ideal of definition $\mathscr{I}_{\alpha}$ of finite type of $X_{\alpha}$ such that for any $\alpha\in L$ and $k\geq 0$ the sheaf $(\mathscr{F}|_{X_{\alpha}})/\mathscr{I}^{k+1}_{\alpha}(\mathscr{F}|_{X_{\alpha}})$ is a quasi-coherent sheaf on the algebraic space $(X_{\alpha},\O_{X_{\alpha}}/\mathscr{I}^{k+1}_{\alpha})$.
\end{itemize}

If $X$ itself has an ideal of definition $\mathscr{I}$ of finite type, then the last condition is equivalent to that $\mathscr{F}_k$ (defined as above) is quasi-coherent on $X_k$ for any $k\geq 0$.

We denote, as usual, the category of adically quasi-coherent sheaves on a formal algebraic space $X$ by $\AQCoh_X$ (or $\AQCoh^{\et}_X$).
By \ref{prop-adicallyqcohetaleequiv}, if $X$ is a quasi-separated adic formal scheme of finite ideal type, then the notion of adically quasi-coherent sheaves on $X$ as a formal algebraic space coincides with that on $X$ as an adic formal scheme.
By \ref{cor-completepullbackaqcsheaves0} and \ref{prop-etaledescent3} we have:
\begin{prop}\label{prop-adicallyqcohsheavesformalalgebraicspaces}
Let $X$ be a formal algebraic space, and $\mathscr{F}$ a complete $\O_X$-module.
Then $\mathscr{F}$ is an adically quasi-coherent sheaf $($resp.\ adically quasi-coherent sheaf of finite type, resp.\ adically quasi-coherent $\O_X$-algebra, resp.\ adically quasi-coherent $\O_X$-algebra of finite type$)$ if and only if for any representable \'etale covering $q\colon Y\rightarrow X$, $\widehat{q^{\ast}}\mathscr{F}$ is an adically quasi-coherent sheaf $($resp.\ adically quasi-coherent sheaf of finite type, resp.\ adically quasi-coherent $\O_Y$-algebra, resp.\ adically quasi-coherent $\O_Y$-algebra of finite type$)$ on $Y$. \hfill$\square$
\end{prop}
\index{adically quasi-coherent (a.q.c.) sheaf|)}

\subsection{Several properties of morphisms}\label{sub-propertiesmorformalalgsp}
We begin with the following definition, consistent with the one given in the beginning of \S\ref{subsub-formalalgebraicspacesdef}:
\begin{dfn}\label{dfn-propertiesmorformalalgsp1}{\rm 
Let $P$ be one of the following properties:
\begin{itemize}
\item[(a)] affine adic\index{morphism of formal algebraic spaces@morphism (of formal algebraic spaces)!affine morphism of formal algebraic spaces@affine ---};
\item[(b)] finite\index{morphism of formal algebraic spaces@morphism (of formal algebraic spaces)!finite morphism of formal algebraic spaces@finite ---};
\item[(c)] open immersion\index{immersion!open immersion of formal algebraic spaces@open --- (of formal algebraic spaces)};
\item[(d)] closed immersion\index{immersion!closed immersion of formal algebraic spaces@closed --- (of formal algebraic spaces)};
\item[(e)] immersion\index{immersion!immersion of formal algebraic spaces@--- (of formal algebraic spaces)}.
\end{itemize}
We say that the a morphism $f\colon X\rightarrow Y$ of formal algebraic spaces is $P$ if it is adic and for any morphism $V\rightarrow Y$ where $V$ is an adic formal scheme, the fiber product $X\times_YV$ is represented by an adic formal scheme and the resulting map $f_V\colon X\times_YV\rightarrow V$ of adic formal schemes has $P$.}
\end{dfn}

\begin{prop}\label{prop-Pmorforalgsp1}
Let $P$ be one of the properties for morphisms of formal algebraic spaces listed in {\rm \ref{dfn-propertiesmorformalalgsp1}}.

{\rm (1)} Let $f\colon X\rightarrow Y$ and $g\colon Y\rightarrow Z$ be morphisms of formal algebraic spaces.
If $f$ and $g$ satisfy $P$, then so is the composition $g\circ f$.

{\rm (2)} Let $Z$ be a formal algebraic space, and $f\colon X\rightarrow X'$ and $g\colon Y\rightarrow Y'$ two $Z$-morphisms of formal algebraic spaces satisfying $P$.
Then $f\times_Zg\colon X\times_ZY\rightarrow X'\times_ZY'$ satisfies $P$.

{\rm (3)} Let $Z$ be a formal algebraic space, and $f\colon X\rightarrow Y$ a $Z$-morphism between formal algebraic spaces satisfying $P$. Then for any morphism $Z'\rightarrow Z$ of formal algebraic spaces the induced morphism $f_{Z'}\colon X\times_ZZ'\rightarrow Y\times_ZZ'$ satisfies $P$. \hfill$\square$
\end{prop}

The proof is straightforward, and we omit it.
The following proposition follows from \ref{prop-etaleeffectivedescentclass}:
\begin{prop}\label{cor-Pmorforalgspetstable1}
Let $P$ be one of the properties for morphisms of formal algebraic spaces listed in {\rm \ref{dfn-propertiesmorformalalgsp1}}.
Let $f\colon X\rightarrow Y$ be a morphism of formal algebraic spaces, and $V\rightarrow Y$ a representable \'etale covering of $Y$.
Then $f$ satisfies $P$ if and only if $X\times_YV$ is a representable and $f_V\colon X\times_YV\rightarrow V$ satisfies $P$.\hfill$\square$
\end{prop}

\begin{prop}\label{prop-Pmorphismtruncation1}
Let $P$ be one of the properties for morphisms of formal algebraic spaces listed in {\rm \ref{dfn-propertiesmorformalalgsp1}}.
Let $f\colon X\rightarrow Y$ be an adic morphism of formal algebraic spaces, and suppose $Y$ has an ideal of definition $\mathscr{I}$ of finite type.
For any integer $k\geq 0$ we denote by $f_k\colon X_k\rightarrow Y_k$ the induced morphism of algebraic spaces defined as in {\rm \S\ref{subsub-idealofdefinitionforalgsp}}.
Then the following conditions are equivalent$:$
\begin{itemize}
\item[{\rm (a)}] $f$ satisfies $P;$
\item[{\rm (b)}] $f_k$ is satisfies $P$ for any $k\geq 0$.
\end{itemize}
If $P=$ `affine adic', `finite', or `closed immersion', then the conditions are furthermore equivalent to 
\begin{itemize}
\item[{\rm (c)}] $f_0$ is satisfies $P$. \hfill$\square$
\end{itemize}
\end{prop}

The proof is straightforward and is omitted here.
\begin{dfn}\label{dfn-propertiesmorformalalgsp2}{\rm 
Let $P$ be one of the following properties:
\begin{itemize}
\item[(a)] locally of finite type\index{morphism of formal algebraic spaces@morphism (of formal algebraic spaces)!morphism of formal algebraic spaces locally of finite type@--- locally of finite type};
\item[(b)] adically flat\index{morphism of formal algebraic spaces@morphism (of formal algebraic spaces)!adically flat morphism of formal algebraic spaces@adically flat ---};
\item[(c)] smooth\index{morphism of formal algebraic spaces@morphism (of formal algebraic spaces)!adically flat morphism of formal algebraic spaces@smooth ---};
\item[(d)] \'etale\index{morphism of formal algebraic spaces@morphism (of formal algebraic spaces)!etale morphism of formal algebraic spaces@\'etale ---}. 
\end{itemize}
We say that an $S$-morphism $f\colon X\rightarrow Y$ of formal algebraic spaces has $P$ if it is adic and there exist a representable \'etale covering $V\rightarrow Y$ and a representable \'etale covering $U\rightarrow X\times_YV$ such that the resulting morphism $U\rightarrow V$ of adic formal schemes satisfies $P$.}
\end{dfn}

Notice that by \ref{prop-dfnadicmorphismformalalgebraicspaces} the definition of \'etale morphisms as in (d) is consistent with the one given in \ref{dfn-etaletopformalalgsp}.
The following proposition is clear:
\begin{prop}\label{prop-Pmorphismtruncation2}
Let $P$ be one of the properties for morphisms of formal algebraic spaces listed in {\rm \ref{dfn-propertiesmorformalalgsp2}}.
Let $f\colon X\rightarrow Y$ be an adic morphism of formal algebraic spaces, and suppose $Y$ has an ideal of definition $\mathscr{I}$ of finite type.
For any integer $k\geq 0$ we denote by $f_k\colon X_k\rightarrow Y_k$ the induced morphism of algebraic spaces defined as in {\rm \S\ref{subsub-idealofdefinitionforalgsp}}.
Then the following conditions are equivalent$:$
\begin{itemize}
\item[{\rm (a)}] $f$ satisfies $P;$
\item[{\rm (b)}] $f_k$ is satisfies $P$ for any $k\geq 0$.
\end{itemize}
If $P=$ `locally of finite type', then the conditions are furthermore equivalent to 
\begin{itemize}
\item[{\rm (c)}] $f_0$ is satisfies $P$. \hfill$\square$
\end{itemize}
\end{prop}

\begin{prop}\label{cor-Pmorforalgspetstable2}
Let $P$ be one of the properties for morphisms of formal algebraic spaces listed in {\rm \ref{dfn-propertiesmorformalalgsp2}}.
Let $f\colon X\rightarrow Y$ be an adic morphism of formal algebraic spaces, and $V\rightarrow Y$ a representable \'etale covering of $Y$.
Then $f$ is satisfies $P$ if and only if the base change $f_V\colon X\times_YV\rightarrow V$ satisfies $P$. \hfill$\square$
\end{prop}

This follows from \ref{prop-etaleeffectivedescentclass2}.
The following proposition is easy to see:
\begin{prop}\label{prop-Pmorforalgsp2}
Let $P$ be one of the properties for morphisms of formal algebraic spaces listed in {\rm \ref{dfn-propertiesmorformalalgsp2}}.

{\rm (1)} Let $f\colon X\rightarrow Y$ and $g\colon Y\rightarrow Z$ be morphisms of formal algebraic spaces.
If $f$ and $g$ satisfy $P$, then so is the composition $g\circ f$.

{\rm (2)} Let $Z$ be a formal algebraic space, and $f\colon X\rightarrow X'$ and $g\colon Y\rightarrow Y'$ two $Z$-morphisms of formal algebraic spaces satisfying $P$.
Then $f\times_Zg\colon X\times_ZY\rightarrow X'\times_ZY'$ satisfies $P$.

{\rm (3)} Let $Z$ be a formal algebraic space, and $f\colon X\rightarrow Y$ a $Z$-morphism between formal algebraic spaces satisfying $P$. Then for any morphism $Z'\rightarrow Z$ of formal algebraic spaces the induced morphism $f_{Z'}\colon X\times_ZZ'\rightarrow Y\times_ZZ'$ satisfies $P$. \hfill$\square$
\end{prop}

\begin{prop}\label{prop-formalalgebraicspaceqcpt}
The following conditions for a morphism $f\colon X\rightarrow Y$ of formal algebraic spaces are equivalent$:$
\begin{itemize}
\item[{\rm (a)}] for any \'etale morphism $V\rightarrow Y$ from a quasi-compact adic formal scheme of finite ideal type, the formal algebraic space $X\times_YV$ is quasi-compact {\rm (\ref{dfn-formalalgebraicspaceqcpt})}$;$
\item[{\rm (b)}] for any morphism $V\rightarrow Y$ from a quasi-compact adic formal scheme of finite ideal type, the formal algebraic space $X\times_YV$ is quasi-compact.
\end{itemize}
\end{prop}

\begin{proof}
The implication (b) $\Rightarrow$ (a) is trivial.
Suppose (a) holds.
We take a representable \'etale covering $Y'\rightarrow Y$, where $Y'$ is quasi-compact, and set $V'=Y'\times_YV$, which is an adic formal scheme \'etale and surjective over $V$.
Since $V'$ is quasi-compact, there exists a quasi-compact open subset $Y''$ of $Y'$ containing the image of $V'$.
Since $X\times_YY''$ is quasi-compact, we deduce that $X\times_YV'$ is quasi-compact (\ref{prop-qsepmorformal3} (3)).
Thus to show the assertion, it suffices to show the following: let $f\colon X\rightarrow Y$ be a morphism of formal algebraic spaces where $Y$ is a quasi-compact adic formal scheme, and $Y'\rightarrow Y$ an \'etale surjective map; if $X\times_YY'$ is quasi-compact, then so is $X$.
This is easy to see.
\end{proof}

\begin{dfn}\label{dfn-formalalgebraicspaceqcptmor}{\rm 
A morphism $f\colon X\rightarrow Y$ of formal algebraic spaces is said to be {\em quasi-compact}\index{morphism of formal algebraic spaces@morphism (of formal algebraic spaces)!quasi-compact morphism of formal algebraic spaces@quasi-compact ---} if the equivalent conditions in \ref{prop-formalalgebraicspaceqcpt} are satisfied.}
\end{dfn}

Clearly, quasi-compact morphisms are closed under composition and base-change stable.

\begin{dfn}\label{dfn-formalalgebraicspacefintypemor}{\rm 
A morphism $f\colon X\rightarrow Y$ of formal algebraic spaces is said to be {\em of finite type}\index{morphism of formal algebraic spaces@morphism (of formal algebraic spaces)!morphism of formal algebraic spaces of finite type@--- of finite type} if it is locally of finite type and quasi-compact.}
\end{dfn}

Obviously, morphisms of finite type are closed under composition and base-change stable; moreover, the conditions (a), (b), and (c) in \ref{prop-Pmorphismtruncation2} with $P=$ `of finite type' are equivalent.

\begin{dfn}\label{dfn-formalalgegraicspaces2}{\rm 
A morphism $f\colon X\rightarrow Y$ of formal algebraic spaces is said to be {\em locally separated}\index{morphism of formal algebraic spaces@morphism (of formal algebraic spaces)!locally separated morphism of formal algebraic spaces@locally separated ---} (resp.\ {\em separated}\index{morphism of formal algebraic spaces@morphism (of formal algebraic spaces)!separated morphism of formal algebraic spaces@separated ---}) if the diagonal map $\Delta_X\colon X\rightarrow X\times_YX$ is a quasi-compact immersion (resp.\ closed immersion).}
\end{dfn}

\begin{prop}\label{prop-locsepmorphismtruncation2}
Let $f\colon X\rightarrow Y$ be an adic morphism of formal algebraic spaces, and suppose $Y$ has an ideal of definition $\mathscr{I}$.
For any integer $k\geq 0$ we denote by $f_k\colon X_k\rightarrow Y_k$ the induced morphism of algebraic spaces defined as in {\rm \S\ref{subsub-idealofdefinitionforalgsp}}.
Then the following conditions are equivalent$:$
\begin{itemize}
\item[{\rm (a)}] $f$ is locally separated $($resp.\ separated$);$
\item[{\rm (b)}] $f_k$ is locally separated $($resp.\ separated$)$ for any $k\geq 0$.
\end{itemize}
Moreover, $f$ is separated if and only if $f_0$ is separated. \hfill$\square$
\end{prop}

\begin{dfn}\label{dfn-propermapforalgsp}{\rm 
A morphism $f\colon X\rightarrow Y$ of formal algebraic spaces is said to be {\em proper}\index{morphism of formal algebraic spaces@morphism (of formal algebraic spaces)!proper morphism of formal algebraic spaces@proper ---} if it is adic and there exists a Zariski covering $\{Y_{\alpha}\rightarrow Y\}_{\alpha\in L}$ of $Y$ and an ideal of definition $\mathscr{I}_{\alpha}$ of finite type on each $Y_{\alpha}$ such that for each $\alpha\in L$ the induced morphism $X_{\alpha,0}=X\times_YY_{\alpha}\rightarrow Y_{\alpha,0}$ of algebraic spaces $($cf.\ {\rm \S\ref{subsub-idealofdefinitionforalgsp})} is proper.}
\end{dfn}

It is easy to see that the condition does not depend on the choice of the Zariski covering $\{Y_{\alpha}\rightarrow Y\}_{\alpha\in L}$ and the ideal of definitions $\{\mathscr{I}_{\alpha}\}_{\alpha\in L}$.
In particular, if $Y$ has an ideal of definition $\mathscr{I}$, then $f$ is proper if and only if it is adic and the induced maps $f_k\colon X_k\rightarrow Y_k$ for any $k\geq 0$ is proper.
Notice that proper morphisms are separated of finite type.
The following proposition is easy to see:
\begin{prop}\label{prop-propermorforalgsp1}
{\rm (1)} Let $f\colon X\rightarrow Y$ and $g\colon Y\rightarrow Z$ be morphisms of formal algebraic spaces.
If $f$ and $g$ are separated $($resp.\ proper$)$, then so is $g\circ f$.

{\rm (2)} Let $Z$ be a formal algebraic space, and $f\colon X\rightarrow X'$ and $g\colon Y\rightarrow Y'$ two separated $($resp.\ proper$)$ $Z$-morphisms of formal algebraic spaces.
Then $f\times_Zg\colon X\times_ZY\rightarrow X'\times_ZY'$ is separated $($resp.\ proper$)$.

{\rm (3)} Let $Z$ be a formal algebraic space, and $f\colon X\rightarrow Y$ a separated $($resp.\ proper$)$ $Z$-morphism of formal algebraic spaces. Then for any map $Z'\rightarrow Z$ of formal algebraic spaces the induced morphism $f_{Z'}\colon X\times_ZZ'\rightarrow Y\times_ZZ'$ is separated $($resp.\ proper$)$. \hfill$\square$
\end{prop}

\subsection{Universally adhesive and universally rigid-Noetherian formal algebraic spaces}\label{sub-adequatealralgsp}
\index{formal algebraic space!universally adhesive formal algebraic space@universally adhesive ---|(}
\index{adhesive!universally adhesive@universally ---!universally adhesive formal algebraic space@--- --- formal algebraic space|(}
\index{formal algebraic space!universally rigid-Noetherian formal algebraic space@universally rigid-Noetherian ---|(}
\begin{dfn}\label{dfn-adequateformalalgsp}{\rm 
A formal algebraic space $X$ is said to be {\em locally universally adhesive} (resp.\ {\em locally universally rigid-Noetherian}) if it has an \'etale covering $Y\rightarrow X$ by a locally universally adhesive (resp.\ locally universally rigid-Noetherian) formal scheme\index{formal scheme!universally adhesive formal scheme@universally adhesive ---}\index{adhesive!universally adhesive@universally ---!universally adhesive formal scheme@--- --- formal scheme}\index{formal scheme!universally rigid-Noetherian formal scheme@universally rigid-Noetherian ---!locally universally rigid-Noetherian formal scheme@locally --- ---} $Y$ (\ref{dfn-formalsch}).
If $Y$ can be taken to be quasi-compact, we say that $X$ is {\em universally adhesive} (resp.\ {\em universally rigid-Noetherian}).}
\end{dfn}

If a formal algebraic space $X$ is locally universally adhesive (resp.\ locally universally rigid-Noetherian), then any locally of finite type formal algebraic space over $X$ is locally universally adhesive (resp.\ locally universally rigid-Noetherian).
Notice that this definition is consistent with the definition of corresponding notions for formal schemes due to \ref{prop-formalfpqcpatchingtua}.

\begin{prop}\label{prop-adqfiberprodformalalgsp}
Let $X\rightarrow Z\leftarrow Y$ be a diagram of formal algebraic spaces, where $X$ is locally universally adhesive $($resp.\ locally universally rigid-Noetherian$)$ and the map $Y\rightarrow Z$ is locally of finite type.
Then the fiber product $X\times_ZY$ is locally universally adhesive $($resp.\ locally universally rigid-Noetherian$)$. \hfill$\square$
\end{prop}

\begin{dfn}\label{dfn-pfformalalgsp}{\rm 
A morphism $f\colon X\rightarrow Y$ of locally universally rigid-Noetherian formal algebraic spaces is said to be {\em locally of finite presentation}\index{morphism of formal algebraic spaces@morphism (of formal algebraic spaces)!morphism of formal algebraic spaces locally of finite presentation@--- locally of finite presentation} if it is adic and there exists a representable \'etale covering $V\rightarrow Y$ and a representable \'etale covering $U\rightarrow X\times_YV$ such that the resulting map $U\rightarrow V$ (between locally universally rigid-Noetherian formal schemes) is locally of finite presentation.
If moreover $f$ is quasi-compact, then $f$ is said to be {\em of finite presentation}\index{morphism of formal algebraic spaces@morphism (of formal algebraic spaces)!morphism of formal algebraic spaces of finite presentation@--- of finite presentation}.}
\end{dfn}

The following propositions are straightforward:
\begin{prop}\label{prop-pfmorphismtruncation2}
Let $f\colon X\rightarrow Y$ be an adic morphism of locally universally rigid-Noetherian formal algebraic spaces, and suppose $Y$ has an ideal of definition $\mathscr{I}$ of finite type.
For any integer $k\geq 0$ we denote by $f_k\colon X_k\rightarrow Y_k$ the induced morphism of algebraic spaces defined as in {\rm \S\ref{subsub-idealofdefinitionforalgsp}}.
Then the following conditions are equivalent$:$
\begin{itemize}
\item[{\rm (a)}] $f$ is locally of finite presentation $($resp.\ of finite presentation$);$
\item[{\rm (b)}] $f_k$ is locally of finite presentation $($resp.\ of finite presentation$)$ for any $k\geq 0$. \hfill$\square$
\end{itemize}
\end{prop}

\begin{prop}\label{cor-pfmorforalgspetstable2}
Let $f\colon X\rightarrow Y$ be an adic morphism of locally universally rigid-Noetherian formal algebraic spaces, and $V\rightarrow Y$ a representable \'etale covering of $Y$.
Then $f$ is locally of finite presentation $($resp.\ of finite presentation$)$ if and only if so is the base change $f_V\colon X\times_YV\rightarrow V$. \hfill$\square$
\end{prop}

\begin{prop}\label{prop-pfmorforalgsp2}
{\rm (1)} Let $f\colon X\rightarrow Y$ and $g\colon Y\rightarrow Z$ be morphisms of locally universally rigid-Noetherian formal algebraic spaces.
If $f$ and $g$ is locally of finite presentation $($resp.\ of finite presentation$)$, then so is $g\circ f$.

{\rm (2)} Let $Z$ be a formal algebraic space, and $f\colon X\rightarrow X'$ and $g\colon Y\rightarrow Y'$ two $Z$-morphisms of locally universally rigid-Noetherian formal algebraic spaces that are locally of finite presentation $($resp.\ of finite presentation$)$.
Suppose $X'\times_ZY'$ is locally universally rigid-Noetherian.
Then $f\times_Zg\colon X\times_ZY\rightarrow X'\times_ZY'$ is locally of finite presentation $($resp.\ of finite presentation$)$.

{\rm (3)} Let $Z$ be a formal algebraic space, and $f\colon X\rightarrow Y$ a $Z$-morphism between locally universally rigid-Noetherian formal algebraic spaces that is locally of finite presentation $($resp.\ of finite presentation$)$. Then for any morphism $Z'\rightarrow Z$ of formal algebraic spaces such that $Y\times_ZZ'$ is locally universally rigid-Noetherian the induced morphism $f_{Z'}\colon X\times_ZZ'\rightarrow Y\times_ZZ'$ is locally of finite presentation $($resp.\ of finite presentation$)$. \hfill$\square$
\end{prop}

By \ref{prop-adicallyqcohetaleequiv2} one can deal with {\em adically quasi-coherent sheaves of finite presentation} over locally universally rigid-Noetherian formal algebraic spaces; here the details are omitted and left to the reader.
\index{formal algebraic space!universally rigid-Noetherian formal algebraic space@universally rigid-Noetherian ---|)}
\index{adhesive!universally adhesive@universally ---!universally adhesive formal algebraic space@--- --- formal algebraic space|)}
\index{formal algebraic space!universally adhesive formal algebraic space@universally adhesive ---|)}

\addcontentsline{toc}{subsection}{Exercises}
\subsection*{Exercises}
\begin{exer}\label{exer-coralgspformalalgsp}
{\rm Show that any algebraic space\index{algebraic space} $X$ over a scheme $S$ is canonically regarded as a formal algebraic space over $S$.
(See {\bf \ref{ch-pre}}, \S\ref{subsub-algebraicspacesconv} for our convention for algebraic spaces.)}
\end{exer}

\begin{exer}\label{exer-propetaleequivquote2}
{\rm Let $X$ be a (resp.\ locally separated, resp.\ separated) formal algebraic space over an adic formal scheme $S$ of finite ideal type, $Z$ and $W$ adic formal schemes over $S$, and $Z\rightarrow X$ and $W\rightarrow X$ maps of sheaves on $\Ac\Fs_{S,\et}$.
Show that the sheaf fiber product $Z\times_XW$ is representable and that the map $Z\times_XW\rightarrow Z\times_SW$ is quasi-compact $($resp.\ a quasi-compact immersion, resp.\ a closed immersion$)$.}
\end{exer}

\begin{exer}\label{exer-idealofdefinitioneffectivedescent}{\rm 
Let $f\colon X\rightarrow Y$ be an adically faithfully flat morphism between adic formal schemes of finite ideal type, and $\mathscr{J}$ an adically quasi-coherent ideal sheaf on $Y$.
Suppose that $\widehat{f^{\ast}}\mathscr{J}$ is an ideal of definition of $X$.
Then show that $\mathscr{J}$ is an ideal of definition of $Y$.}
\end{exer}

\begin{exer}\label{exer-formalalgebraicspaceszariskicoveringqcpt}{\rm 
Show that for any formal algebraic space $X$ there exists a Zariski covering (that is, a covering by open subspaces) $\{X_{\alpha}\rightarrow X\}_{\alpha\in L}$ such that each $X_{\alpha}$ is a quasi-compact formal algebraic space.
Show, moreover, that in this situation if each $X_{\alpha}$ is a formal scheme, then $X$ is a formal scheme.}
\end{exer}

\begin{exer}\label{exer-pointsformalalgebraicspaces}{\rm 
Let $X$ be a formal algebraic space.

(1) (Cf.\ \cite[II.6]{Knu}.) A {\em point}\index{point!point of a formal algebraic space@--- (of a formal algebraic space)} of $X$ is a morphism of the form $\Spec k\rightarrow X$ (where $k$ is a field) that is a sheaf monomomorphism.
We denote by $|X|$ the set of all isomorphism classes of points of $X$.
Suppose that $X$ has an ideal of definition of finite type $\mathscr{I}$, and let $X_0=(X,\O_X/\mathscr{I})$ be the associated algebraic space.
Then show that $|X_0|=|X|$.
In particular, $|X|\neq\emptyset$ whenever $X$ is non-empty.

(2) The {\em underlying topological space} of $X$ is the set $|X|$ endowed with the topology as follows: a subset of $|X|$ is closed if it is of the form $|Y|$ for a closed subspace $Y$ of $X$.
Show that there exists a canonical one to one correspondence between the set of all open subspaces of $X$ and the set of all open subsets of $|X|$.}
\end{exer}

\begin{exer}\label{exer-pointsformalalgebraicspaces2}{\rm 
Let $X$ be a non-empty formal algebraic space.
Then show that there exists a dense open subspace $Y$ of $X$ that is a formal scheme.
(Note that, according to our convention, all formal algebraic spaces are quasi-separated.)}
\end{exer}


\section{Cohomology theory}\label{sec-cohomologycoherent}
In this section we collect fundamental facts on cohomologies of adically quasi-coherent sheaves on formal algebraic spaces.
After giving some general facts on cohomologies of adically quasi-coherent sheaves, we will discuss in \S\ref{sub-coherentadqfor} coherent sheaves on locally universally adhesive formal algebraic spaces of a certain kind.
The last subsection \S\ref{sub-calculusderivedformal} collects some derived categorical calculi that are needed in the later arguments. 

Here we would like to mention again that when we say `$X$ is a formal algebraic space', we always mean either one of the following (cf.\ {\bf \ref{ch-pre}}, \S\ref{subsub-algebraicspacesconv}):
\begin{itemize}
\item $X$ is an adic formal scheme of finite ideal type (but not necessarily quasi-separated); in this case, unless otherwise clearly stated, all adically quasi-coherent sheaves on $X$ and their cohomologies are considered with respect to the Zariski topology;
\item $X$ is a quasi-separated adic formal algebraic space of finite ideal type; in this case, all adically quasi-coherent sheaves on $X$ and their cohomologies are considered with respect to the \'etale topology (\S\ref{subsub-etaletopformalalgsp}).
\end{itemize}

\subsection{Cohomologies of adically quasi-coherent sheaves}\label{sub-cohomologyaqcoh}
\begin{thm}\label{thm-formalacyclicity}
{\rm (1)} Let $A$ be an adic ring of finite ideal type\index{adic!adic ring@--- ring!adic ring of finite ideal type@--- --- of finite ideal type} {\rm (\ref{dfn-admissibleringsadicrings} (2), \ref{dfn-admissibleringoffiniteidealtype})}, and set $X=\Spf A$.
Then for any adically quasi-coherent sheaf\index{adically quasi-coherent (a.q.c.) sheaf} $\mathscr{F}$ {\rm (\ref{dfn-adicqcoh})} we have 
$$
\H^q(X,\mathscr{F})=0
$$
for $q\geq 1$.

{\rm (2)} Let $f\colon X\rightarrow Y$ be an affine adic morphism between adic formal schemes of finite ideal type, and $\mathscr{F}$ an adically quasi-coherent sheaf on $X$. 
Then 
$$
\RD^qf_{\ast}\mathscr{F}=0
$$
for $q\geq 1$.
\end{thm}

\begin{thm}\label{thm-formalacyclicityder}
{\rm (1)} Let $A$ be an adic ring of finite ideal type, and set $X=\Spf A$.
Then for any complex $\mathscr{F}^{\bullet}$ of adically quasi-coherent sheaves on $X$ such that $\mathscr{F}^q=0$ for $q\ll 0$ $($resp.\ $|q|\gg 0)$, we have $\H^q(X,\mathscr{F}^{\bullet})=0$ for $q\ll 0$ $($resp.\ $|q|\gg 0)$, and the object $\RD^+\Gamma_X(\mathscr{F}^{\bullet})$ of $\DC^+(A)$ is represented by the complex $\Gamma_X(\mathscr{F}^{\bullet})$.

{\rm (2)} Let $f\colon X\rightarrow Y$ be an affine adic morphism between formal schemes of finite ideal type.
Let $\mathscr{F}^{\bullet}$ be a complex of adically quasi-coherent sheaves on $X$ such that $\mathscr{F}^q=0$ for $q\ll 0$ $($resp.\ $|q|\gg 0)$.
Then $\RD^q f_{\ast}\mathscr{F}^{\bullet}=0$ for $q\ll 0$ $($resp.\ $|q|\gg 0)$, and the object $\RD^+ f_{\ast}\mathscr{F}^{\bullet}$ of $\DC^+(Y)$ is represented by the complex $f_{\ast}\mathscr{F}^{\bullet}$.
\end{thm}

Similarly to the proof of {\bf \ref{ch-pre}}.\ref{thm-vanishcohaffder-0} one can deduce \ref{thm-formalacyclicityder} from \ref{thm-formalacyclicity}.
Hence it suffices to show \ref{thm-formalacyclicity}.
\begin{proof}[Proof of Theorem {\rm \ref{thm-formalacyclicity}}]
As (1) is a special case of \ref{lem-vanishingcohomologyadicallyuseful} (2), we only need to show (2).
Since $\RD^qf_{\ast}\mathscr{F}$ is the sheaf on $Y$ associated to the presheaf given by $U\mapsto \H^q(f^{-1}(U),\mathscr{F})$, we may assume $Y$ is affine of the form $Y=\Spf B$, where $B$ is an adic ring that has a finitely generated ideal of definition $I\subseteq B$.
Then it suffices to check that $\H^q(X,\mathscr{F})$ vanishes for $q\geq 0$.
In this situation, $X$ is also of the form $X=\Spf A$, where $A$ is an adic ring having $IA$ as a finitely generated ideal of definition.
Then the desired vanishing follows from (1).
\end{proof}

\begin{cor}\label{cor-formalacyclicity4}
Let $f\colon X\rightarrow Y$ be an adic and separated morphism of adic formal schemes of finite ideal type, and suppose $Y$ is quasi-compact.
Then there exists an integer $r>0$ such that for any adically quasi-coherent sheaf $\mathscr{F}$ on $X$ and any $q\geq r$ we have $\RD^q f_{\ast}\mathscr{F}=0$.
If, moreover, $Y$ is affine, then one can take as $r$ $($the minimum of$)$ the number of affine open sets that cover $X$. $($Hence such an $f$ always has finite cohomological dimension.$)$
\end{cor}

\begin{proof}
We may assume that $Y$ is affine of the form $Y=\Spf B$, where $B$ is an adic ring of finite type.
The formal scheme $X$ can be covered by finitely many affine open subsets; by \ref{cor-sepmorformal12} all intersections of the members of the covering are affine.
Hence in view of \ref{thm-formalacyclicity} the assertion follows from the \v{C}ech calculation of the cohomology using the Leray covering.
\end{proof}

\subsection{Coherent sheaves}\label{sub-coherentadqfor}
\index{coherent!coherent sheaf on a formal algebraic space@--- sheaf (on a formal algebraic space)|(}
\begin{dfn}\label{dfn-cohformalschemes1}{\rm 
A formal algebraic space $X$ is said to be {\em universally cohesive}\index{formal algebraic space!universally cohesive formal algebraic space@universally cohesive ---}\index{cohesive!universally cohesive@universally --- (formal algebraic space)} if it is locally universally rigid-Noetherian\index{formal scheme!universally rigid-Noetherian formal scheme@universally rigid-Noetherian ---!locally universally rigid-Noetherian formal scheme@locally --- ---} (\ref{dfn-adequateformalalgsp}) and for any locally of finite presentation morphism\index{morphism of formal algebraic spaces@morphism (of formal algebraic spaces)!morphism of formal algebraic spaces locally of finite presentation@--- locally of finite presentation} $Y\rightarrow X$ from a formal scheme $Y$, $\O_Y$ is a coherent $\O_Y$-module.}
\end{dfn}

\begin{prop}\label{prop-cohformalschemes2}
Let $A$ be a t.u.\ rigid-Noetherian ring\index{t.u. rigid-Noetherian ring@t.u.\ rigid-Noetherian ring} {\rm (\ref{dfn-tuaringadmissible} (1))}, and $I\subseteq A$ a finitely generated ideal of definition.
Set $X=\Spf A$.
Then $X$ is universally cohesive if $A$ is topologically universally coherent\index{coherent!topologically universally coherent@topologically universally ---} {\rm ({\bf \ref{ch-pre}}.\ref{dfn-cohringsmodules5})} with respect to $I$.
\end{prop}

\begin{proof}
Suppose $A$ is topologically universally coherent, and let $Y\rightarrow X$ be an $X$-formal scheme locally of finite presentation.
In order to show that $\O_Y$ is coherent, we may assume that $Y$ is affine $Y=\Spf B$.
By the assumption $B$ is a topologically finitely presented $B$-algebra and, in view of {\bf \ref{ch-pre}}.\ref{prop-cohringsmodules52} (2), topologically universally coherent with respect to the ideal $IB$.
Hence it suffices only to show that $\O_X$ is coherent (since for $Y$, the situation is the same).
Let $U$ be an open subset of $X$, and suppose we are given an exact sequence
$$
0\longrightarrow\mathscr{K}\longrightarrow\O_U^{\oplus m}\longrightarrow\O_U.
$$
We are to show that $\mathscr{K}$ is of finite type.
We may further assume that $U$ is an affine open subset, and thus suppose $U=X$ without loss of generality.
By \ref{lem-corcoradicqcoh1} the last morphism comes from a map $A^{\oplus m}\rightarrow A$.
Let $K$ be the kernel, which is a finitely generated $A$-module.
Then by \ref{thm-adicqcohpre1}
$$
0\longrightarrow K^{\Delta}\longrightarrow\O_X^{\oplus m}\longrightarrow\O_X
$$
is exact, that is, we have $\mathscr{K}\cong K^{\Delta}$, which is an adically quasi-coherent sheaf of finite type.
\end{proof}

\begin{prop}\label{prop-adicqcoh6}
Let $X$ be a locally universally adhesive formal algebraic space satisfying the following condition$:$ there exists a covering $($Zariski or \'etale$)$ $\coprod_{\alpha\in L}U_{\alpha}\rightarrow X$, where each $U_{\alpha}$ is affine with an ideal of definition $\mathscr{I}_{\alpha}$ of finite type, such that $\O_{U_{\alpha}}$ is $\mathscr{I}_{\alpha}$-torsion free.
Then $X$ is universally cohesive. \hfill$\square$
\end{prop}

This follows from \ref{prop-cohformalschemes2} and {\bf \ref{ch-pre}}.\ref{thm-pf02a} (2).
Hence any formal scheme locally of finite presentation over $X$ as above, even if it may have $\mathscr{I}$-torsions, is universally cohesive.
For example, if $V$ is an $a$-adically complete valuation ring (of arbitrary height), then the $a$-adic formal scheme $\Spf V$, and hence any formal algebraic space locally of finite presentation over $V$, is universally cohesive ({\bf \ref{ch-pre}}.\ref{cor-convadh}).

\begin{prop}\label{prop-adicqcoh55}
Let $X$ be a universally cohesive formal algebraic space, and suppose $X$ has an ideal of definition $\mathscr{I}$ of finite type.
Set $X_k=(X,\O_X/\mathscr{I}^{k+1})$ for $k\geq 0$ $($note that $X_k$ is a universally cohesive scheme$)$.
Let $\mathscr{F}$ be an $\O_X$-module, and set $\mathscr{F}_k=\mathscr{F}/\mathscr{I}^{k+1}\mathscr{F}$ for $k\geq 0$.
Then the following conditions are equivalent$:$
\begin{itemize}
\item[{\rm (a)}] $\mathscr{F}$ is coherent$;$
\item[{\rm (b)}] $\mathscr{F}_k$ is coherent on $X_k$ for any $k\geq 0$ and $\mathscr{F}=\varprojlim_k\mathscr{F}_k$.
\end{itemize}
\end{prop}

\begin{proof}
Notice that in either situation (a) or (b) the sheaf $\mathscr{F}$ is adically quasi-coherent of finite type (cf.\ \ref{prop-adicallyqcohinductivelimit}).
We may work in the affine situation; set $X=\Spf A$ where $A$ is a t.u.\ rigid-Noetherian ring with a finitely generated ideal of definition $I\subseteq A$.
In this case, we have $\mathscr{F}=M^{\Delta}$ for a uniquely determined finitely generated $A$-module $M$.
Then the assertion follows from {\bf \ref{ch-pre}}.\ref{cor-ARconseq2}.
\end{proof}
\index{coherent!coherent sheaf on a formal algebraic space@--- sheaf (on a formal algebraic space)|)}

\subsection{Calculi in derived categories}\label{sub-calculusderivedformal}
Let $X$ be a formal algebraic space.
We denote by $\DC^{\ast}(X)$ $(\ast=$``\ \ '', $+$, $-$, $\bd)$ the derived category\index{derived category} associated to the abelian category $\Mod_X$ of $\O_X$-modules.
In case $X$ is universally cohesive, one can consider the full subcategory $\DC^{\ast}_{\coh}(X)$ of $\DC^{\ast}(X)$ consisting of objects having coherent cohomologies in all degrees.
This is a triangulated category equipped with the canonical cohomology functor\index{functor!cohomology functor@cohomology ---!canonical cohomology functor@canonical --- ---} $\mathcal{H}^0$ and the canonical $t$-structure\index{t-structure@$t$-structure!canonical t-structure@canonical ---} (cf.\ {\bf \ref{ch-pre}}, \S\ref{sub-formalismderived}).
\begin{prop}[{cf.\ {\rm {\bf \ref{ch-pre}}.\ref{prop-derivedcategory71}}}]\label{prop-formalismformalderived1}
Let $A$ be a t.u.\ rigid-Noetherian ring with a finitely generated ideal of definition $I\subseteq A$, and suppose that $A$ is topologically universally coherent with respect to $I$.
Set $X=\Spf A$.
Then the exact equivalence $\Coh_A\stackrel{\sim}{\rightarrow}\Coh_X$ by $M\mapsto M^{\Delta}$ {\rm (\ref{thm-adicqcohpre1})} induces an exact equivalence
$$
\DC^{\ast}(\Coh_A)\stackrel{\sim}{\longrightarrow}\DC^{\ast}(\Coh_X),
$$
where $\Coh_A$ denotes the category of finitely presented $A$-modules $($which is an abelian category due to {\rm {\bf \ref{ch-pre}}.\ref{prop-cohringsmodules1}}$)$. \hfill$\square$
\end{prop}

In the following we denote the composition
$$
\DC^{\ast}(A)\stackrel{\sim}{\longrightarrow}\DC^{\ast}(\Coh_X)\stackrel{\delta^{\ast}}{\longrightarrow}\DC^{\ast}_{\coh}(X)
$$
(cf.\ {\bf \ref{ch-pre}}, \S\ref{sub-subcategoryderived}) by 
$$
M\longmapsto M_X.
$$

\begin{prop}\label{prop-formalismformalderived2}
Let $A$ and $I$ be as in {\rm \ref{prop-formalismformalderived1}}, and set $X=\Spf A$.
Then the canonical exact functor
$$
\delta^{\bd}\colon\DC^{\bd}(\Coh_X)\rightarrow\DC^{\bd}_{\coh}(X)
$$
$($cf.\ {\rm {\bf \ref{ch-pre}}, \S\ref{sub-subcategoryderived})} is a categorical equivalence. \hfill$\square$
\end{prop}

\begin{cor}\label{cor-formalismformalderived21}
In the situation as in {\rm \ref{prop-formalismformalderived1}} the canonical functor 
$$
\DC^{\bd}(\Coh_A)\longrightarrow\DC^{\bd}_{\coh}(X),\qquad M\longmapsto M_X,
$$
is an exact equivalence. 
In other words, any object $M$ of $\DC^{\bd}_{\coh}(X)$ can be represented {\rm ({\bf \ref{ch-pre}}.\ref{dfn-derivedcategory5})} by a complex consisting of coherent sheaves and hence by a complex consisting of finitely presented $A$-modules. \hfill$\square$
\end{cor}

All these can be shown similarly to {\bf \ref{ch-pre}}.\ref{prop-formalismderived2} and {\bf \ref{ch-pre}}.\ref{cor-formalismderived21}.
By an argument similar to {\bf \ref{ch-pre}}.\ref{prop-schpairadh3} one has:
\begin{prop}\label{prop-schformalpairadh3}
Let $X$ be a universally cohesive formal algebraic space.

{\rm (1)} For $F,G\in\obj(\DC^-_{\coh}(X))$ the object $F\otimes^{\LD}_{\O_X}G$ belongs to $\DC^-_{\coh}(X)$.

{\rm (2)} For $F\in\obj(\DC^-_{\coh}(X))$ and $G\in\obj(\DC^+_{\coh}(X))$ the object $\RD\lHom_{\O_X}(F,G)$ belongs to $\DC^+_{\coh}(X)$. \hfill$\square$
\end{prop}

\begin{prop}\label{prop-schformalpairadh31}
Let $f\colon X\rightarrow Y$ be a morphism of universally cohesive formal algebraic spaces.
Then the functor $\LD f^{\ast}$ maps $\DC^-_{\coh}(Y)$ to $\DC^-_{\coh}(X)$. \hfill$\square$
\end{prop}

The proof is similar to that of {\bf \ref{ch-pre}}.\ref{prop-schpairadh31}.
Using \ref{thm-formalacyclicityder} (2) one can show the following by an argument similar to that in {\bf \ref{ch-pre}}.\ref{cor-calculusderived2-0}:
\begin{prop}\label{prop-calculusformalderived2}
Let $X$ be a universally cohesive formal algebraic space, and $i\colon Y\hookrightarrow X$ a closed immersion\index{immersion!closed immersion of formal algebraic spaces@closed --- (of formal algebraic spaces)} of finite presentation $($hence $Y$ is universally cohesive$)$.
Then $\RD^+i_{\ast}$ maps $\DC^{\bd}_{\coh}(Y)$ to $\DC^{\bd}_{\coh}(X)$. \hfill$\square$
\end{prop}

\addcontentsline{toc}{subsection}{Exercises}
\subsection*{Exercises}
\begin{exer}\label{exer-universallycohesivekernelaqc}{\rm 
Let $X$ be a universally cohesive (hence locally universally rigid-Noetherian) formal algebraic space, and $\varphi\colon\mathscr{F}\rightarrow\mathscr{G}$ a morphism of adically quasi-coherent sheaves of finite type on $X$.
Show that, if $\mathscr{G}$ is coherent, then $\ker(\varphi)$ is an adically quasi-coherent sheaf of finite type.}
\end{exer}


\section{Finiteness theorem for proper algebraic spaces}\label{sec-cohomology}
In this section we announce and prove the finiteness theorem for cohomologies of coherent sheaves on universally cohesive ({\bf \ref{ch-pre}}.\ref{dfn-universallycohesive2}) algebraic spaces (not necessarily locally Noetherian), which generalizes the classically known finiteness theorem in, e.g., \cite{EGA}, {\bf III}. This theorem provides an important preliminary for (generalized) GFGA theorems, which will be discussed in the next two sections. We put this entirely scheme-theoretic (not formal-scheme-theoretic) section in this chapter, not only because it gives a preliminary for the GFGA, but also some of the techniques for the proof are common to the proofs of the GFGA theorems.
Most notably, a {\it variation} of {\it Grothendieck's d\'evissage} is discussed, which we call the {\it carving method}. 
As in case of d\'evissage, it consists of a reduction to particular cases by an induction with respect to sequences of closed subspaces, that is, {\it induction with respect to a stratification}. 
Aside from its technical merits (e.g., removing Noetherian hypothesis), the method is best suited for treating algebraic spaces, since an inductive argument with respect to a stratification is a basic and fundamental tool to study these spaces.

\medskip
In the first subsection, \S\ref{sub-finitudes}, we announce the finiteness theorem (\ref{thm-fini}) and some related statements for algebraic spaces. The rest of this section is mainly devoted to the proof.  
By the carving method, the finiteness theorem is reduced to the projective case, and the theorem in this particular situation is verified by a generalization of Serre's theorem \cite{EGA}, $\mathbf{III}$, (2.2.1), which will be announced and proved in \S\ref{sub-projfinitudes}.
In \S\ref{sub-proofmethod}, we formulate the carving method in a derived-categorical setting. This will be contained in Proposition \ref{prop-proofmethod} below, and the proof of this proposition is based on the carving lemma\index{carving lemma} (\ref{lem-proofmethod}). The proof of the finiteness theorem will be finished in \S\ref{sub-pffinitudes}.
In the final subsection \S\ref{sub-appliIgood} we give  an application of the theorem to $I$-goodness (cf.\ {\bf \ref{ch-pre}}.\ref{subsub-ARgeneral}) of the induced filtrations of the cohomology groups of coherent sheaves, which we will need in our later discussion.

\subsection{Finiteness theorem: Announcement}\label{sub-finitudes}
First, let us collect some known finiteness results that are already general enough for our purpose:
\begin{prop}[cf.\ {\cite[$\mathbf{III}$, (1.4.12)]{EGA}}]\label{prop-cohqcoh2}
Let $f\colon X\rightarrow S$ be a quasi-compact and separated morphism of algebraic spaces, and suppose that $S$ is quasi-compact.
Then there exists an integer $r>0$ such that for any quasi-coherent sheaf $\mathscr{F}$ on $X$ and any $q\geq r$, we have $\RD^qf_{\ast}\mathscr{F}=0$.
If, moreover, $S$ is affine and $X$ is a scheme, one can take as $r$ $($the minimum of$)$ the number of affine open subsets that cover $X$.
\end{prop}

This follows readily from {\bf \ref{ch-pre}}.\ref{prop-cohqcoh-0} and {\bf \ref{ch-pre}}.\ref{cor-cohqcoh1-0} (cf.\ \cite[II.2.7, II.3.12]{Knu}).
The proposition says, in other words, that the morphism $f$ as above has {\em finite cohomological dimension}.

\begin{prop}\label{prop-fini1}
Let $f\colon X\rightarrow Y$ be a separated and quasi-compact morphism of algebraic spaces with $X$ quasi-compact, and $\mathscr{L}$ an $f$-ample invertible sheaf on $X$.
Let $\mathscr{F}$ be a quasi-coherent sheaf on $X$ of finite type. 
Then there exists an integer $N$ such that for any $n\geq N$ the canonical morphism $f^{\ast}f_{\ast}\mathscr{F}(n)\rightarrow\mathscr{F}(n)$ is surjective.
\end{prop}

\begin{proof}
One can assume that $Y$ is affine, and hence that $X$ is a scheme.
Then the assertion follows from \cite[$\mathbf{II}$, (4.6.8)]{EGA}.
\end{proof}

Now we announce our finiteness theorem:
\begin{thm}[Finiteness theorem]\label{thm-fini}
Let $Y$ be a universally cohesive\index{algebraic space!universally cohesive algebraic space@universally cohesive ---} {\rm ({\bf \ref{ch-pre}}.\ref{dfn-universallycohesive2})} quasi-compact algebraic space\index{algebraic space}, and $f\colon X\rightarrow Y$ a proper morphism of finite presentation of algebraic spaces $($hence $X$ is also universally cohesive$)$.
Then the functor $\RD f_{\ast}$ maps an object of $\DC^{\bd}_{\coh}(X)$ to an object of $\DC^{\bd}_{\coh}(Y)$.
\end{thm}

Note that the premise `universally cohesive' for $Y$ is fulfilled, in particular, when $Y$ is Noetherian or is locally of finite presentation over an $a$-adically complete valuation ring $V$ of arbitrary height ({\bf \ref{ch-pre}}.\ref{cor-cohschemes21ver2}).
Before the proof of the theorem we include some useful corollaries:
\begin{cor}\label{cor-fini}
Consider the situation as in {\rm \ref{thm-fini}}.
Then the functor $\RD f_{\ast}$ maps $\DC^{\ast}_{\coh}(X)$ to $\DC^{\ast}_{\coh}(Y)$ for $\ast=$``\ \ '', $+$, $-$, $\bd$.
\end{cor}

Here, note that, since $f$ has a finite cohomological dimension (\ref{prop-cohqcoh2}), one can extend the domain of $\RD^+ f_{\ast}$ to the whole $\DC_{\coh}(X)$ (cf.\ \cite[C.D.\ Chap.\ 2, \S2, $\mathrm{n}^{\mathrm{o}}\ 2$, Cor.\ 2]{SGA4.5}).

\begin{proof}
Let $K$ be an object of $\DC^+_{\coh}(X)$, and consider $\RD f_{\ast}K$.
To detect the cohomology $\RD^qf_{\ast}K$, we can always find a sufficiently large $n$ such that $\RD^qf_{\ast}K=\RD^qf_{\ast}(\tau^{\leq n}K)$; therefore, $\RD f_{\ast}$ maps $\DC^+_{\coh}(X)$ to $\DC^+_{\coh}(Y)$.
The others are easy to see (use $\tau^{\geq m}$ instead of $\tau^{\leq n}$).
\end{proof}

\begin{cor}\label{cor-fini2}
Let $B$ be a universally coherent\index{coherent!coherent ring@--- ring!universally coherent ring@universally --- ---} ring $({\bf \ref{ch-pre}}.\ref{dfn-cohringsmodules5})$, and $f\colon X\rightarrow Y=\Spec B$ a proper morphism of finite presentation between algebraic spaces.
Then if $M\in\DC^-_{\coh}(X)$ and $N\in\DC^+_{\coh}(X)$, then $\RD^q\Hom_{\O_X}(M,N)$ $($with the natural $B$-module structure$;$ cf.\ {\rm {\bf \ref{ch-pre}}, \S\ref{subsub-ringedspacemodulestrcoh}}$)$ is a coherent $B$-module for any $q$.
\end{cor}

\begin{proof}
By {\bf \ref{ch-pre}}.\ref{prop-schpairadh3} (2) we know that $\RD\lHom_{\O_X}(M,N)$ belongs to $\DC^+_{\coh}(X)$.
Since $\RD\Hom_{\O_X}(M,N)=\RD\Gamma_X\circ\RD\lHom_{\O_X}(M,N)$, the assertion follows from \ref{cor-fini}.
\end{proof}

\subsection{Generalized Serre's theorem}\label{sub-projfinitudes}
\subsubsection{Announcement}\label{subsub-announcegenSerre}
The following theorem together with \ref{prop-fini1} generalizes \cite[$\mathbf{III}$, (2.2.1)]{EGA}:
\begin{thm}\label{thm-fini2}
Let $Y$ be a universally cohesive\index{algebraic space!universally cohesive algebraic space@universally cohesive ---} quasi-compact algebraic space, $f\colon X\rightarrow Y$ a proper morphism of finite presentation, and $\mathscr{F}$ a coherent sheaf on $X$.
Suppose there exists an $f$-ample invertible sheaf $\mathscr{L}$ on $X$.

{\rm (1)} The sheaf $\RD^qf_{\ast}\mathscr{F}$ for any $q$ is coherent on $Y$.

{\rm (2)} There exists an integer $N$ such that $\RD^qf_{\ast}\mathscr{F}(n)=0$ for $n\geq N$ and $q>0$.
\end{thm}

The rest of this subsection is devoted to the proof of this theorem. 
\subsubsection{Reduction process}\label{subsub-reductionprocessgenSerre}
We may assume that $Y$ is affine $Y=\Spec B$, where $B$ is universally coherent ({\bf \ref{ch-pre}}.\ref{prop-cohschemes2}).
We may further assume that $\mathscr{L}$ is very ample relative to $f$.
Since $f$ is proper and quasi-projective, it is in fact projective (\cite[II.7.8]{Knu}), and hence there exists an closed $Y$-immersion $i\colon X\hookrightarrow\P^r_B$ for some $r>0$ such that $\mathscr{L}\cong i^{\ast}\O_{\P}(1)$.
Let $g$ be the structure morphism $\P^r_B\rightarrow Y$.
As $i$ is affine, we have $\RD f_{\ast}\mathscr{F}(n)\cong\RD g_{\ast}(i_{\ast}\mathscr{F}(n))$; also, by projection formula (cf.\ \cite[$\mathbf{0}_{\mathbf{III}}$, (12.2.3)]{EGA}) we have $i_{\ast}(\mathscr{F}(n))\cong (i_{\ast}\mathscr{F})(n)$.
Thus we may assume $X=\P^r_B$ without loss of generality.
Therefore, to show \ref{thm-fini2}, it suffices to prove:
\begin{prop}\label{prop-fini2}
Let $B$ be a universally coherent ring, $X=\P^r_B$, and $\mathscr{F}$ a coherent $\O_X$-module.

{\rm (1)} For any $q$ the $B$-module $\H^q(X,\mathscr{F})$ is finitely presented.

{\rm (2)} There exists an integer $N$ such that $\H^q(X,\mathscr{F}(n))=0$ for $n\geq N$ and $q>0$.
\end{prop}

\subsubsection{Proof of Proposition \ref{prop-fini2}}\label{subsub-proofgenSerre}
\begin{lem}\label{lem-inftyrepsplitbundles}
Let $B$ be a universally coherent ring, and $X=\P^r_B$.
Then for any coherent $\O_X$-module $\mathscr{F}$ there exists an exact sequence
$$
\cdots\longrightarrow\mathscr{E}_1\longrightarrow\mathscr{E}_0\longrightarrow\mathscr{F}\longrightarrow 0
$$
by sheaves $\mathscr{E}_i$ $(i\geq 0$) of the form $\O_X(n_i)^{\oplus k_i}$.
\end{lem}

\begin{proof}
The proof is done by successive application of {\bf \ref{ch-pre}}.\ref{prop-genglobalsec-0}; by {\bf \ref{ch-pre}}.\ref{prop-genglobalsec-0} there exists a surjective morphism $\O_X(n_0)^{\oplus k_0}\rightarrow\mathscr{F}$ for some integers $n_0$ and $k_0>0$.
Let $\mathscr{K}_0$ be the kernel.
Since $\mathscr{K}_0$ is again coherent, we have similarly a surjective map $\O_X(n_1)^{\oplus k_1}\rightarrow\mathscr{K}_0$.
One can repeat this procedure.
\end{proof}

\begin{lem}\label{lem-pf7}
Let $B$ be a coherent ring, and $X$ an algebraic space over $B$. 
Let 
$$
\mathscr{E}^{\bullet}=(\cdots\rightarrow 0\rightarrow\mathscr{E}^s\rightarrow\mathscr{E}^{s+1}\rightarrow\cdots\rightarrow\mathscr{E}^r\rightarrow0\rightarrow\cdots)
$$
be a bounded complex of $\O_X$-modules such that for all $q$ and $k$ the cohomology group $\H^q(X,\mathscr{E}^k)$ is a finitely presented $B$-module.
Then $\H^q(X,\mathscr{E}^{\bullet})$ for any $q$ is a finitely presented $B$-module.
\end{lem}

\begin{proof}
The proof is done by induction with respect to the length of $\mathscr{E}^{\bullet}$.
By a suitable shift we may assume that the complex $\mathscr{E}^{\bullet}$ is of the form
$$
\mathscr{E}^{\bullet}=(\cdots\rightarrow 0\rightarrow\mathscr{E}^0\rightarrow\mathscr{E}^1\rightarrow\cdots\rightarrow\mathscr{E}^{l-1}\rightarrow0\rightarrow\cdots).
$$
Set $\mathscr{E}^{\prime\bullet}=\sigma^{\geq 1}\mathscr{E}^{\bullet}$ (the stupid truncation), that is, 
$$
\mathscr{E}^{\prime\bullet}=(\cdots\rightarrow 0\rightarrow 0\rightarrow\mathscr{E}^1\rightarrow\cdots\rightarrow\mathscr{E}^{l-1}\rightarrow0\rightarrow\cdots).
$$
Then we have the distinguished triangle 
$$
\mathscr{E}^{\prime\bullet}\longrightarrow\mathscr{E}^{\bullet}\longrightarrow\mathscr{E}^0\stackrel{+1}{\longrightarrow}
$$
in $\KC^{\bd}(\Mod_X)$.
By the cohomology exact sequence
$$
\H^{q-1}(X,\mathscr{E}^0)\longrightarrow\H^q(X,\mathscr{E}^{\prime\bullet})\longrightarrow\H^q(X,\mathscr{E}^{\bullet})\longrightarrow\H^q(X,\mathscr{E}^0)\longrightarrow\H^{q+1}(X,\mathscr{E}^{\prime\bullet}),
$$
and by {\bf \ref{ch-pre}}.\ref{prop-cohringsmodules1}, we deduce that $\H^q(X,\mathscr{E}^{\bullet})$ is finitely presented, as desired.
\end{proof}

By the proposition and {\bf \ref{ch-pre}}.\ref{cor-cohqcoh1-0} we have:
\begin{cor}\label{cor-pf8}
Let $Y$ be a universally cohesive algebraic space, and $f\colon X\rightarrow Y$ a morphism of algebraic spaces.
Let $\mathscr{E}^{\bullet}$ be a bounded complex of $\O_X$-modules such that $\RD^qf_{\ast}\mathscr{E}^k$ are coherent $\O_Y$-modules for all $q$ and $k$.
Then $\RD^qf_{\ast}\mathscr{E}^{\bullet}$ are coherent $\O_Y$-modules for all $q$. \hfill$\square$
\end{cor}

\begin{proof}[Proof of Proposition {\rm \ref{prop-fini2}}]
We take a resolution of $\mathscr{F}$ as in \ref{lem-inftyrepsplitbundles}.
By \ref{prop-cohqcoh2} one already knows that the cohomology group $\H^q(X,\mathscr{F})$ vanishes for $q>r$, since $\P^r_B$ is separated over $B$ and covered by $r+1$ affine open subsets.
Hence it is enough to calculate the cohomology groups for $0\leq q\leq r$. 
To this end, one is allowed to replace $\mathscr{F}$ by the bounded complex 
$$
\cdots\longrightarrow 0\longrightarrow\mathscr{E}_{r+1}\longrightarrow\cdots\longrightarrow\mathscr{E}_1\longrightarrow\mathscr{E}_0\longrightarrow0\longrightarrow\cdots.
$$
By {\bf \ref{ch-pre}}.\ref{cor-cohprojsp-0} and \ref{lem-pf7} we deduce that $\H^q(X,\mathscr{F})$ is a finitely presented $B$-module for $0\leq q\leq r$.
Combined with the fact that $\H^q(X,\mathscr{F})$ for $q>r$ are trivially finitely presented, we get the assertion (1).

Next, let us show (2).
For $q>r$ we already know that $\H^q(X,\mathscr{F}(n))=0$ for any $n$.
Hence one can prove the assertion by descending induction with respect to $q$.
Consider the surjection $\mathscr{E}=\mathscr{E}_0\rightarrow\mathscr{F}$ as above, and let $\mathscr{K}$ be its kernel (finitely presented due to {\bf \ref{ch-pre}}.\ref{prop-cohringsmodules1}).
By induction, there exists $N$ such that $\H^{q+1}(X,\mathscr{K}(n))=0$ for $n\geq N$.
Moreover, one can take $N$ sufficiently large so that for $n\geq N$ the sheaf $\mathscr{E}(n)$ is of the form $\O_X(m)^{\oplus k}$ by a positive integer $m$.
Hence $\H^q(X,\mathscr{E}(n))=0$ and $\H^{q+1}(X,\mathscr{E}(n))=0$ for $q>0$ ({\bf \ref{ch-pre}}.\ref{cor-cohprojsp-0}).
Then the vanishing of $\H^q(X,\mathscr{F}(n))$ follows from the cohomology exact sequence.
\end{proof}

\subsection{The carving method}\label{sub-proofmethod}
\index{carving method|(}
\subsubsection{The main assertion}\label{subsub-carvingannounce}
Let $B$ be a universally coherent ring\index{coherent!coherent ring@--- ring!universally coherent ring@universally --- ---}, and consider the category $\PAs_B$, a full subcategory of the category $\As_B$ of $B$-algebraic spaces consisting of proper and finitely presented $B$-algebraic spaces.
Note that every algebraic space in $\PAs_B$ is universally cohesive\index{algebraic space!universally cohesive algebraic space@universally cohesive ---}.

\begin{prop}[Carving method]\label{prop-proofmethod}
Suppose we are given for each object $f\colon X\rightarrow\Spec B$ of $\PAs_B$ a full subcategory 
$$
\DC_f\subseteq\DC^{\bd}_{\coh}(X)
$$
satisfying the following conditions$:$
\begin{itemize}
\item[{\bf (C0)}] the zero object $0$ belongs to $\DC_f$, and $\DC_f$ is stable under isomorphisms in $\DC^{\bd}_{\coh}(X);$
\item[{\bf (C1)}] let $K\rightarrow L\rightarrow M\stackrel{+1}{\rightarrow}$ be a distinguished triangle in $\DC^{\bd}_{\coh}(X);$
if two of $K$, $L$, and $M$ are in $\DC_f$, then so is the rest$;$
\item[{\bf (C2)}] if $f$ is projective, then $\DC_f=\DC^{\bd}_{\coh}(X);$
\item[{\bf (C3)}] consider a morphism $\pi$ in $\PAs_B$, which amounts to the same as a commutative diagram of the form
$$
\xymatrix@C-5ex{\til{X}\ar[rr]^{\pi}\ar[dr]_g&&X\ar[dl]^f\\ &\Spec B\rlap{,}}
$$
and suppose $g$ and $\pi$ are projective$;$ 
then $\RD\pi_{\ast}$ maps $\DC_g=\DC^{\bd}_{\coh}(\til{X})$ to $\DC_f;$
\item[$\mathbf{(C4)}$] consider a closed immersion $\iota$ in $\PAs_B$, that is, a commutative diagram
$$
\xymatrix@C-5ex{Z\ar[rr]^{\iota}\ar[dr]_g&&X\ar[dl]^f\\ &\Spec B\rlap{,}}
$$
where $\iota$ is a closed immersion$;$
then if $\DC_g=\DC^{\bd}_{\coh}(Z)$, the functor $\RD\iota_{\ast}$ maps $\DC_g$ to $\DC_f$.
\end{itemize}
Then we have $\DC_f=\DC^{\bd}_{\coh}(X)$ for any object $f\colon X\rightarrow\Spec B$ of $\PAs_B$.
\end{prop}

\subsubsection{Preparation for the proof and carving lemma}\label{subsub-carvinglemma}
The rest of this subsection is devoted to showing the proposition.
Let us first establish the following fact, which will be tacitly used in the sequel:

\medskip
{\sc Claim 0.} {\it If $M\in\obj(\DC_f)$, then $M[n]\in\obj(\DC_f)$ for any integer $n$.}

\medskip
Indeed, we have the distinguished triangle 
$$
M\longrightarrow 0\longrightarrow M[1]\stackrel{+1}{\longrightarrow},
$$
whence the claim due to {\bf (C1)}.

\medskip
{\sc Claim 1.} {\it Consider the diagram as in $\mathbf{(C4)}$ where $\iota$ is a closed immersion defined by a nilpotent quasi-coherent ideal $\mathscr{J}\subset\O_X$.
Then if $\DC_g=\DC^{\bd}_{\coh}(Z)$, we have $\DC_f=\DC^{\bd}_{\coh}(X)$.}

\medskip
Since $\iota$ is finitely presented, there exists an integer $k\geq 0$ such that $\mathscr{J}^{k+1}=0$.
By induction with respect to $k$ we may assume $\mathscr{J}^2=0$ without loss of generality.
Let $M\in\obj(\DC^{\bd}_{\coh}(X))$.
Taking the distinguished triangle
$$
\tau^{\leq n}M\longrightarrow M\longrightarrow\tau^{\geq n+1}M\stackrel{+1}{\longrightarrow}
$$
and {\sc Claim 0} into account, by an inductive application of {\bf (C1)} we may assume that $M$ is concentrated in degree $0$, that is, $\mathcal{H}^q(M)=0$ unless $q=0$.

Consider the canonical morphism $M\rightarrow\RD\iota_{\ast}(\tau^{\geq 0}\LD\iota^{\ast}M)$, and embed it into a distinguished triangle
$$
M\longrightarrow\RD\iota_{\ast}(\tau^{\geq 0}\LD\iota^{\ast}M)\longrightarrow N\stackrel{+1}{\longrightarrow}.\leqno{(\ast)}
$$
By {\bf \ref{ch-pre}}.\ref{prop-schpairadh31} we have $\tau^{\geq 0}\LD\iota^{\ast}M\in\obj(\DC^{\bd}_{\coh}(Z))=\obj(\DC_g)$.
By {\bf (C4)} the middle term of $(\ast)$ belongs to $\DC_f$.
Hence by {\bf (C1)} it suffices to show that $N$ belongs to $\DC_f$.
To see this, we shall show that the canonical morphism $N\rightarrow\RD\iota_{\ast}(\tau^{\geq 0}\LD\iota^{\ast}N)$ is an isomorphism.
Consider the cohomology exact sequence of $(\ast)$; since $\iota$ is affine, by {\bf \ref{ch-pre}}.\ref{thm-vanishcohaff-0} (2) we have 
$$
\mathcal{H}^q(\RD\iota_{\ast}(\tau^{\geq 0}\LD\iota^{\ast}M))=
\begin{cases}
\iota_{\ast}\iota^{\ast}\mathcal{H}^0(M)&\textrm{if}\ q=0,\\
0&\textrm{otherwise.}
\end{cases}
$$
Hence the initial portion of the cohomology exact sequence gives
$$
0\longrightarrow\mathcal{H}^{-1}(N)\longrightarrow\mathcal{H}^0(M)\longrightarrow\iota_{\ast}\iota^{\ast}\mathcal{H}^0(M)\longrightarrow\mathcal{H}^0(N)\longrightarrow 0.
$$
Since $\iota_{\ast}\iota^{\ast}\mathcal{H}^0(M)\cong\mathcal{H}^0(M)/\mathscr{J}\mathcal{H}^0(M)$, we deduce that $\mathcal{H}^{-1}(N)=\mathscr{J}\mathcal{H}^0(M)$ and $\mathcal{H}^0(N)=0$.
Other parts of the cohomology exact sequence imply that $\mathcal{H}^q(\RD\iota_{\ast}(\tau^{\geq 0}\LD\iota^{\ast}M))=\mathcal{H}^q(N)=0$ for $q>0$.
Hence $N$ is concentrated in degree $-1$, and $\mathcal{H}^{-1}(N)=\mathscr{J}\mathcal{H}^0(M)$, which is a coherent sheaf on $X$ (for $\iota$ is of finite presentation).
Since $\mathscr{J}^2=0$, we deduce that the canonical morphism $N\rightarrow\RD\iota_{\ast}(\tau^{\geq 0}\LD\iota^{\ast}N)$ is an isomorphism.

Then again by {\bf \ref{ch-pre}}.\ref{prop-schpairadh31} we have $\tau^{\geq 0}\LD\iota^{\ast}N\in\obj(\DC^{\bd}_{\coh}(Z))=\obj(\DC_g)$ and deduce by {\bf (C4)} and {\bf (C0)} that $N$ belongs to $\DC_f$, as desired.

\medskip
{\sc Claim 2.} {\it Consider the diagram as in $\mathbf{(C4)}$ where $\iota$ is a closed immersion.
Suppose $\DC_g=\DC^{\bd}_{\coh}(Z)$.
Then for any $M\in\obj(\DC^{\bd}_{\coh}(X))$ such that $\mathcal{H}^q(M)$ is supported on $Z$ for any $q$, we have $M\in\obj(\DC_f)$.}

\medskip
By a similar reduction process as above, we may assume that $M$ is concentrated in degree $0$.
If $\mathscr{J}$ is the ideal defining $Z$ in $X$, there exists a positive integer $n$ such that $\mathscr{J}^n\mathcal{H}^0(M)=0$.
By {\sc Claim 1} we may further assume $n=1$, that is, $\mathscr{J}\mathcal{H}^0(M)=0$.
Then the canonical morphism $M\rightarrow\RD\iota_{\ast}(\tau^{\geq 0}\LD\iota^{\ast}M)$ is an isomorphism.
Since $\tau^{\geq 0}\LD\iota^{\ast}M$ belongs to $\DC^{\bd}_{\coh}(Z)$ ({\bf \ref{ch-pre}}.\ref{prop-schpairadh31}), we have $M\in\obj(\DC_f)$ by {\bf (C4)} and {\bf (C0)}.

\medskip
To proceed, we need the following lemma:
\begin{lem}[Carving lemma]\label{lem-proofmethod}\index{carving lemma}
Let $B$ be any ring, and $f\colon X\rightarrow Y=\Spec B$ an algebraic space separated and of finite presentation over $B$.
Then there exists a sequence of closed subspaces 
$$
X=X_0\supsetneq X_1\supsetneq\cdots\supsetneq X_r\supsetneq X_{r+1}=\emptyset\leqno{(\ast)}
$$
such that the following conditions are satisfied$:$
\begin{itemize}
\item[{\rm (a)}] for any $0\leq i\leq r$ the closed immersion $\iota_{i+1}\colon X_{i+1}\hookrightarrow X_i$ is of finite presentation$;$
\item[{\rm (b)}] for any $0\leq i\leq r$ the complement $X_i\setminus X_{i+1}$ $($cf.\ {\rm \cite[II.5.12]{Knu}}$)$ is quasi-projective$;$
\item[{\rm (c)}] for any $0\leq i\leq r$ there exists a projective morphism $\pi_i\colon\til{X}_i\rightarrow X_i$ of finite presentation such that $\pi_i$ is an isomorphism over the open subspace $X_i\setminus X_{i+1}$ and that the composition $\til{X}_i\rightarrow Y$ is quasi-projective.
\end{itemize}
\end{lem}

\begin{proof}
Take a subring $B'$ of $B$ of finite type over $\Z$ and an algebraic space $X'$ over $Y'=\Spec B'$ of finite presentation such that there exists a Cartesian diagram as follows:
$$
\xymatrix{X\ar[r]\ar[d]_f&X'\ar[d]^{f'}\\ Y\ar[r]&Y'\rlap{.}}
$$
This is possible by the standard argument well-known for schemes (\cite[$\mathbf{IV}$, \S 8]{EGA}).
We leave the details to the reader as an exercise (Exercise \ref{exer-desalgsp}).
The algebraic space $X'$ is then a Noetherian algebraic space, and hence any descending sequence of closed subspaces terminates (cf.\ \cite[II.5.18]{Knu}).

Suppose $X'=X'_0$ is not quasi-projective over $Y'$.
Take a non-empty open subspace $U'_1$ of $X'$ quasi-projective over $Y'$, and let $X'_1$ be the complement.
If $X'_1$ is not quasi-projective, then we repeat this procedure to take $X'_2$ with non-empty quasi-projective complement.
Since $X'$ is Noetherian, this procedure stops after finitely many steps.
Thus we get the sequence of closed subspaces
$$
X'=X'_0\supsetneq X'_1\supsetneq\cdots\supsetneq X'_r\supsetneq X'_{r+1}=\emptyset.
$$
Now we apply Chow's lemma (stated below) to each $X'_i$ and $U'_i=X'_i\setminus X'_{i+1}$ to get the $U'_i$-admissible blow-up $\til{X'_i}\rightarrow X'_i$ such that the composition $\til{X'_i}\rightarrow Y'$ is quasi-projective.
Pulling back all these data by the map $Y\rightarrow Y'$ onto $Y$, we get the desired sequence of subspaces of $X$.
\end{proof}

\begin{thm}[Chow's lemma {\cite[Premi\`ere partie, (5.7.14)]{RG}}]\label{thm-chow}\index{Chow's lemma}\index{Chow, W-L.}
Let $Y$ be a coherent algebraic space, $X\rightarrow Y$ a separated $Y$-algebraic space of finite type, and $U$ an open subspace quasi-projective over $Y$.
Then there exists a $U$-admissible blow-up $\pi\colon\til{X}\rightarrow X$ such that $\til{X}$ is quasi-projective over $Y$. \hfill$\square$
\end{thm}

Here a blow-up $\til{X}\rightarrow X$ is said to be {\em $U$-admissible}\index{admissible!U-admissible blow-up@$U$-{---} blow-up}\index{blow-up!U-admissible blow-up@$U$-admisible ---} if it is isomorphic to the blow-up\index{blow-up} along a closed subspace of $X$ of finite presentation disjoint from $U$ (cf.\ {\bf \ref{ch-rigid}}.\ref{dfn-Uadmissibleblowups}).

\subsubsection{Proof of Proposition \ref{prop-proofmethod}}
Let $f\colon X\rightarrow\Spec B$ be an object of $\PAs_B$, and apply \ref{lem-proofmethod} to $f$ to get a sequence of closed subspaces $(\ast)$.
Since $f$ is proper, the morphisms $f_i\colon X_i\rightarrow\Spec B$ and $\til{f}_i=f_i\circ\pi_i\colon\til{X}_i\rightarrow\Spec B$ are proper for $0\leq i\leq r$.
In particular, $\til{f}_i$ are projective, and hence $\til{X}_i$ are schemes.
Note also that $X_r=\til{X}_r$, which is therefore a projective scheme over $B$.

We prove the assertion by induction with respect to $r$.
If $r=0$, then $X=X_0$ is projective, and the result follows from {\bf (C2)}.
In general, we may assume by induction that $\DC_{f_1}=\DC^{\bd}_{\coh}(X_1)$.
In what follows, for the sake of brevity, we write $Z=X_1$, $g=f_1$, $\til{X}=\til{X}_0$, $\pi=\pi_0$, and $\iota=\iota_1$.
Thus we are in the situation depicted as follows:
$$
\xymatrix@C-5ex{&&\til{X}\ar[d]^{\pi}\\ Z\ar[rr]^{\iota}\ar[dr]_g&&X\ar[dl]^f\\ &\Spec B\rlap{.}}
$$
Here $\iota$ is a closed immersion, $\pi$ and $f\circ\pi$ are projective, and $\pi$ is identity over $X\setminus\iota(Z)$; moreover, we already know that $\DC_g=\DC^{\bd}_{\coh}(Z)$.

Let $M\in\obj(\DC^{\bd}_{\coh}(X))$. 
We are going to show that $M$ belongs to $\DC_f$.
Similarly to the proof of {\sc Claim 1} we may assume that $M$ is concentrated in degree $0$.
Embed the canonical map $M\rightarrow\RD\pi_{\ast}(\tau^{\geq 0}\LD\pi^{\ast}M)$ in the distinguished triangle
$$
M\longrightarrow\RD\pi_{\ast}(\tau^{\geq 0}\LD\pi^{\ast}M)\longrightarrow N\stackrel{+1}\longrightarrow.\leqno{(\ast)}
$$
By {\bf \ref{ch-pre}}.\ref{prop-schpairadh31}, {\bf (C2)}, and {\bf (C3)} the middle term belongs to $\DC_f$.
Hence by {\bf (C1)} it suffices to show that $N$ belongs to $\DC_f$.
But by {\sc Claim 2} it is enough to show that all cohomologies of $N$ are supported on $Z$.

Since $\tau^{\geq 1}\RD\pi_{\ast}(\tau^{\geq 0}\LD\pi^{\ast}M)\cong\tau^{\geq 1}N$ and since $\pi$ is isomorphic over the complement of $Z$ in $X$, $\tau^{\geq 1}N$ has the cohomologies supported on $Z$.
Moreover, looking at the initial portion of the cohomology exact sequence of $(\ast)$
$$
0\longrightarrow\mathcal{H}^{-1}(N)\longrightarrow\mathcal{H}^0(M)\longrightarrow\mathcal{H}^0(\RD\pi_{\ast}(\tau^{\geq 0}\LD\pi^{\ast}M))\longrightarrow\mathcal{H}^0(N)\longrightarrow0
$$
and taking the equality $\mathcal{H}^0(\RD\pi_{\ast}(\tau^{\geq 0}\LD\pi^{\ast}M))=\pi_{\ast}\pi^{\ast}\mathcal{H}^0(M)$ into account, we deduce that $\mathcal{H}^{-1}(N)$ and $\mathcal{H}^0(N)$ are supported on $Z$.
Hence all cohomologies of $N$ are supported on $Z$, as desired. \hfill$\square$
\index{carving method|)}

\subsection{Proof of Theorem \ref{thm-fini}}\label{sub-pffinitudes}
\subsubsection{Reduction process}\label{subsub-pffinitudesreduction}
We assume without loss of generality that $Y$ is affine $Y=\Spec B$ where $B$ is universally coherent\index{coherent!coherent ring@--- ring!universally coherent ring@universally --- ---}.
For any object $(f\colon X\rightarrow\Spec B)\in\obj(\PAs_B)$ we define 
$$
\DC_f\subseteq\DC^{\bd}_{\coh}(X)
$$
to be the full subcategory consisting of all objects $M\in\obj(\DC^{\bd}_{\coh}(X))$ such that $\RD f_{\ast}M\in\obj(\DC^{\bd}_{\coh}(\Spec B))$.
We want to show that $\DC_f=\DC^{\bd}_{\coh}(X)$ for any $f\in\obj(\PAs_B)$, and hence it suffices to check the conditions in \ref{prop-proofmethod}.

We postpone the checking of {\bf (C2)}, which we suppose is verified for the time being.
The conditions {\bf (C0)} and {\bf (C1)} are obviously satisfied.
As for {\bf (C3)}, since we already know {\bf (C2)}, $\RD\pi_{\ast}$ maps any object $M$ of $\DC^{\bd}_{\coh}(\til{X})$ to an object of $\DC^{\bd}_{\coh}(X)$.
Then since $\RD f_{\ast}(\RD\pi_{\ast}M)=\RD g_{\ast} M$, we have $\RD\pi_{\ast}M\in\DC_f$.
Thus {\bf (C3)} is satisfied.
Finally, in the diagram in {\bf (C4)}, since $\iota$ is of finite presentation, $\RD\iota_{\ast}$ maps an object of $\DC^{\bd}_{\coh}(Z)$ to an object of $\DC^{\bd}_{\coh}(X)$ ({\bf \ref{ch-pre}}.\ref{cor-calculusderived2-0}).
Hence {\bf (C4)} can be verified similarly.

\subsubsection{End of the proof}\label{subsub-finitudesproof}
Now we finish the proof of \ref{thm-fini} by checking the condition {\bf (C2)}.
By a similar reduction argument as in \S\ref{subsub-reductionprocessgenSerre}, we may assume $X=\P^r_B$.
Indeed, if $\iota\colon X\hookrightarrow\P^r_B$ is the $B$-closed immersion of finite presentation, then $\RD\iota_{\ast}$ maps an object $M$ of $\DC^{\bd}_{\coh}(X)$ to $\DC^{\bd}_{\coh}(\P^r_B)$ ({\bf \ref{ch-pre}}.\ref{cor-calculusderived2-0}); if it is proved that the theorem is true for $\P^r_B\rightarrow\Spec B$, then it is also true for $X\rightarrow\Spec B$ by composition formula of right derived functors.

Now if $X=\P^r_B$, then by induction with respect to $\amp(M)$ ({\bf \ref{ch-pre}}.\ref{dfn-derivedcategory6}) using the distinguished triangles of the form 
$$
\tau^{\leq n}M\longrightarrow M\longrightarrow\tau^{\geq n}M\stackrel{+1}{\longrightarrow}
$$
and \ref{lem-pf7}, we may assume $\amp(M)=0$; by a suitable shift we may assume that $M$ is concentrated in degree $0$.
In this case $M$ is represented by a single sheaf.
But the assertion in this case is nothing but \ref{thm-fini2} (1), thereby the theorem. \hfill$\square$

\subsection{Application to $I$-goodness}\label{sub-appliIgood}
\begin{prop}\label{prop-fini122}
Let $f\colon X\rightarrow Y$ be as in $\ref{thm-fini}$, and $\mathscr{S}=\bigoplus_{k\geq 0}\mathscr{S}_k$ a quasi-coherent positively graded $\O_Y$-algebra of finite presentation.
Set $\mathscr{S}'=f^{\ast}\mathscr{S}$.
Let $\mathscr{M}=\bigoplus_{k\in\Z}\mathscr{M}_k$ be a quasi-coherent graded $\mathscr{S}'$-module of finite presentation.
Then $\RD^qf_{\ast}\mathscr{M}=\bigoplus_{k\in\Z}\RD^qf_{\ast}\mathscr{M}_k$ is a quasi-coherent graded $\mathscr{S}$-module of finite presentation for any $q\geq 0$.
Moreover, there exists a positive integer $r$ such that for any $q\geq r$ we have $\RD^qf_{\ast}\mathscr{M}=0$. \hfill$\square$
\end{prop}

The first half of the proposition can be shown by the same reasoning as in \cite[$\mathbf{III}$, (3.3.1)]{EGA} and hence as in \cite[$\mathbf{III}$, (2.4.1)]{EGA}, where we use \ref{thm-fini} instead of \cite[$\mathbf{III}$, (2.2.1)]{EGA}.
Since the proof is quite similar to them, we are not going to repeat it.
The last half follows from the proof of \cite[$\mathbf{III}$, (3.3.1)]{EGA} and \ref{prop-cohqcoh2}.

\begin{prop}\label{prop-fini4}
Let $(B,I)$ be a universally adhesive pair\index{pair!adhesive pair@adhesive ---!universally adhesive pair@universally --- ---}\index{adhesive!universally adhesive@universally ---} {\rm ({\bf \ref{ch-pre}}.\ref{dfn-universallyadhesive})}, and suppose $B$ is universally coherent.
Let $f\colon X\rightarrow Y=\Spec B$ be a proper morphism of finite presentation between algebraic spaces, and $\mathscr{F}$ a coherent $\O_X$-module.
Then for any $q\geq 0$ the filtration $F^{\bullet}=\{F^n\}_{n\geq 0}$ on the $B$-module $\H^q(X,\mathscr{F})$ by $F^n=\image(\H^q(X,I^n\mathscr{F})\rightarrow\H^q(X,\mathscr{F}))$ is $I$-good\index{I-good@$I$-good (filtration)}\index{filtration by submodules@filtration (by submodules)!I-good filtration by submodules@$I$-good ---} {\rm ({\bf \ref{ch-pre}}.\ref{dfn-Igood})}.
\end{prop}

\begin{proof}
It is clear that $I^mF^n\subseteq F^{n+m}$.
Consider the $I$-torsion part $\mathscr{F}_{\Itor}$, and take $m\geq 0$ such that $I^n\mathscr{F}=0$ for any $n\geq m$.
Set $\mathscr{F}_1=\mathscr{F}/\mathscr{F}_{\Itor}$.
We first claim that $I^n\mathscr{F}\cong I^n\mathscr{F}_1$ for $n\gg 0$.
To see this, we may assume that $X$ is affine $X=\Spec A$ where $A$ is a finitely presented $B$-algebra; notice that $A$ is $I$-adically universally adhesive and universally coherent.
In this situation we have $\mathscr{F}=\til{N}$ for a coherent $A$-module $N$.
The canonical morphism $I^nN\rightarrow I^n(N/N_{\Itor})$ is clearly surjective, and its kernel is $I^nN\cap N_{\Itor}$.
By {\bf \ref{ch-pre}}.\ref{prop-AR} the filtration $\{I^nN\cap N_{\Itor}\}_{n\geq 0}$ on $N_{\Itor}$ is equivalent to the $I$-adic one, and we conclude $I^nN\cap N_{\Itor}=0$ for $n\gg 0$, that is, $I^nN\cong I^n(N/N_{\Itor})$, which is what we have claimed.
The claim shows that, in order to check the $I$-goodness of $F^{\bullet}$, we may replace $\mathscr{F}$ by $\mathscr{F}_1$ and hence may assume that $\mathscr{F}$ is $I$-torsion free.

Now, consider the Rees algebra\index{Rees algebra} $S=R(B,I)=\bigoplus_{n\geq 0}I^n$ (cf.\ {\bf \ref{ch-pre}}, \S\ref{sub-Reescone}) and the graded $S$-module $M=\bigoplus_{n\geq 0}\H^q(X,I^n\mathscr{F})$.
We claim that $M$ is finitely generated $S$-module; if this is shown, then $\bigoplus_{n\geq 0}F^n$ is, as a quotient of $M$, a finitely generated $S$-module, and hence the desired result follows from {\bf \ref{ch-pre}}.\ref{prop-Reescone1}.

Set $S_1=S/S_{\Itor}$.
Then $S_1$ is finitely presented $B$-algebra; indeed, it is clearly of finite type, and hence we have a surjective map of the form $B[X_1,\ldots,X_n]\rightarrow S_1$; since $B[X_1,\ldots,X_n]$ is adhesive and $S_1$ is $I$-torsion free, the kernel is finitely generated.
Set $\mathscr{S}=\til{S}$ (resp.\ $\mathscr{S}_1=\til{S}_1$), and $\mathscr{S}'=f^{\ast}\mathscr{S}$ (resp.\ $\mathscr{S}'_1=f^{\ast}\mathscr{S}_1$).
Then $\mathscr{S}_1$ is a quasi-coherent positively graded $\O_Y$-algebra of finite presentation.
Let $\mathscr{M}=\bigoplus_{n\geq 0}I^n\mathscr{F}$, which is a quasi-coherent graded $\mathscr{S}'$-algebra of finite type.
We have $\H^q(X,\mathscr{M})=M$.
Since $\mathscr{F}$ is $I$-torsion free, $\mathscr{M}$ is $I$-torsion free and hence carries the canonical $\mathscr{S}'_1$-module structure; moreover, by the adhesiveness of $\mathscr{S}'_1$ one sees that $\mathscr{M}$ is finitely presented as an $\mathscr{S}'_1$-module.
Hence we deduce that the cohomology $\H^q(X,\mathscr{M})=M$ is a finitely presented $S_1$-module (\ref{prop-fini122}); in particular, it is finitely generated as an $S$-module, as desired.
\end{proof}

\addcontentsline{toc}{subsection}{Exercises}
\subsection*{Exercises}
\begin{exer}\label{exer-fini1}{\rm 
Let $B$ be a ring, $f\colon X\rightarrow Y=\Spec B$ a quasi-compact separated morphism of algebraic spaces, and $\mathscr{L}$ an $f$-ample invertible sheaf on $X$.
Suppose $\O_X$ is coherent on $X$.
Then show that the canonical functor 
$$
\DC^{\ast}(\Coh_X)\longrightarrow\DC^{\ast}_{\coh}(\QCoh_X)
$$
is an equivalence of triangulated categories for $\ast=-$, $\bd$.}
\end{exer}


\begin{exer}\label{exer-desalgsp}{\rm 
Let $A$ be a ring, and $X$ a locally separated (resp.\ separated) $A$-algebraic space of finite presentation. 
Then there exists a subring $A_0$ of $A$ of finite type over $\Z$ and a locally separated (resp.\ separated) $A_0$-algebraic space $X_0$ of finite presentation such that $X_0\otimes_{A_0}A\cong X$.}
\end{exer}


\section{GFGA comparison theorem}\label{sec-GFGAcom}
\index{GFGA!GFGA comparison theorem@--- comparison theorem|(}
In this and the next sections, we discuss GFGA ($=$ g\'eom\'etrie formelle et g\'eom\'etrie alg\'ebrique) theorems\index{GFGA}, which generalize the classical GFGA theorems in \cite[$\mathbf{III}$]{EGA}.
The generalization will be done in the following two directions: 
First, with the application to rigid geometry in mind, we drop the Noetherian hypothesis and replace it by weaker ones like universally adhesive, etc.
Second, our argument will be done entirely in the derived categorical language. 

In this section, we state and prove the GFGA comparison theorem.
The theorem will be announced in the first subsection \S\ref{sub-GFGAcomann}.
In \S\ref{sub-GFGAcomcl} we give the classical version, the comparison of the cohomologies without derived categorical setting but regarding topologies, along with the formalism of \cite[$\mathbf{III}$, \S4.1]{EGA}.
Notice that even this classical version avoids Noetherian hypotheses, and thus gives a generalization of the corresponding theorem in \cite[$\mathbf{III}$, \S4.1]{EGA}.
The subsection \S\ref{sub-proofGFGAcom} is devoted to the proof of the theorem, in which we again use the carving method\index{carving method} developed in \S\ref{sub-proofmethod}.
The final subsection \S\ref{sub-GFGAcomext} presents, as a corollary, a subsidiary result on comparison of Ext modules, which will be referred to in the next section.

\subsection{Announcement of the theorem}\label{sub-GFGAcomann}
\subsubsection{Formal completion functor}\label{subsub-GFGAcomannconst}
Let $(X,W)$ be a pseudo-adhesive pair\index{pseudo-adhesive}\index{adhesive!adhesive pair@--- pair!pseudo adhesive pair@pseudo-{---} ---}\index{pair!adhesive pair@adhesive ---!pseudo adhesive pair@pseudo-{---} ---} of algebraic spaces\footnote{Confer {\bf \ref{ch-pre}}, \S\ref{subsub-algebraicspacesconv} and the beginning of \S\ref{sec-cohomologycoherent} for our general conventions for algebraic spaces and formal algebraic spaces.} ({\bf \ref{ch-pre}}.\ref{dfn-schpairadh}) such that $X$ is universally cohesive\index{cohesive!universally cohesive@universally --- (schemes)} ({\bf \ref{ch-pre}}.\ref{dfn-universallycohesive}). 
Note that such a situation is attained if $(X,W)$ is universally adhesive and $\O_X$ is $\mathscr{I}_W$-torsion free or, more generally, if the pair $(X,W)$ is of finite presentation over another such pair ({\bf \ref{ch-pre}}.\ref{prop-thm-cor-cohschemes21}).
For example, one can a pair $(X,W)$ where $X$ is finitely presented over an $a$-adically complete valuation ring\index{valuation!valuation ring@--- ring!a-adically complete valuation ring@$a$-adically complete --- ---} of arbitrary height and $W$ is defined by the ideal $(a)$.

Set $\widehat{X}=\widehat{X}|_W$, the formal completion\index{completion!formal completion@formal ---} of $X$ along $W$ (cf.\ \S\ref{subsub-formalcompletionalgsp}), and let $j\colon\widehat{X}\rightarrow X$ be the canonical morphism of formal algebraic spaces.
By $\mathscr{I}$-adic completion (where $\mathscr{I}=\mathscr{I}_W$, the defining ideal of $W$ in $X$), we have the functor 
$$
\Mod_X\longrightarrow\Mod_{\widehat{X}},\qquad
\mathscr{F}\longmapsto\widehat{\mathscr{F}}=\widehat{\mathscr{F}}|_W.
$$
But by \ref{prop-adqformalprop} (1) this functor restricted on coherent $\O_X$-modules is canonically equivalent to the functor 
$$
\mathscr{F}\longmapsto j^{\ast}\mathscr{F}.
$$
Guided by this, instead of the completion functor, we consider the functor
$$
\for\colon\Mod_X\longrightarrow\Mod_{\widehat{X}},\qquad
\mathscr{F}\longmapsto\mathscr{F}^{\for}=j^{\ast}\mathscr{F},
$$
which is exact, since $j$ is flat, and hence induces an exact functor of triangulated categories
$$
\DC^{\ast}(X)\longrightarrow\DC^{\ast}(\widehat{X}),
$$
for $\ast=$``\ \ '', $+$, $-$, $\bd$. 
We write the functor thus obtained as
$$
M\longmapsto M^{\for}.
$$
Note that, if $M\in\obj(\DC^{\ast}_{\coh}(X))$, then we have 
$$
\mathcal{H}^q(M^{\for})=j^{\ast}\mathcal{H}^q(M)=\widehat{\mathcal{H}^q(M)}.
$$

\begin{rem}\label{rem-GFGAcomfunctordef}{\rm 
If, moreover, the formal algebraic space $\widehat{X}$ is universally cohesive, then we have the exact functor
$$
\DC^{\ast}_{\coh}(X)\longrightarrow\DC^{\ast}_{\coh}(\widehat{X}), \qquad M\longmapsto M^{\for}
$$
(for $\ast=$``\ \ '', $+$, $-$, $\bd$), the so-called {\em comparison functor}\index{functor!comparison functor@comparison ---}.
Note that under our hypothesis in the beginning of this subsection, $\widehat{X}$ is universally cohesive if $(X,W)$ is finitely presented over a complete t.u.\ adhesive pair\index{pair!adhesive pair@adhesive ---!topologically universally adhesive pair@topologically universally --- ---} $(A,I)$ such that $A$ is $I$-torsion free (\ref{prop-cohformalschemes2}).
For example, if $X$ is finitely presented over an $a$-adically complete valuation ring $V$ of arbitrary height and $W$ is the closed subspace defined by $(a)$, then $\widehat{X}$ is universally cohesive.}
\end{rem}

\subsubsection{The statement}\label{subsub-GFGAcomann}
In the sequel we work in the following situation:
\begin{sit}\label{sit-GFGA}{\rm 
Let $(Y,Z)$ be a universally pseudo-adhesive pair of algebraic spaces with $Y$ universally cohesive, and $f\colon X\rightarrow Y$ a proper and finitely presented morphism.
Set $W=f^{-1}(Z)$.
Set $\widehat{X}=\widehat{X}|_W$ and $\widehat{Y}=\widehat{Y}|_Z$, and consider the commutative diagram
$$
\xymatrix{X\ar[d]_f&\widehat{X}\ar[l]_j\ar[d]^{\widehat{f}}\\ Y&\widehat{Y}\ar[l]^i\rlap{,}}
$$
which is Cartesian in the category of formal algebraic spaces (cf.\ \ref{prop-formalcomplmor}).
For any $k\geq 0$ we denote by $X_k$ (resp.\ $Y_k$) the closed subspace of $X$ (resp.\ $Y$) defined by the ideal $\mathscr{I}^{k+1}\O_X$ (resp.\ $\mathscr{I}^{k+1}$), where $\mathscr{I}$ is the defining ideal of $Z$ on $Y$.}
\end{sit}

By \ref{thm-fini} and \ref{cor-fini} we already know that the functor $\RD f_{\ast}$ maps $\DC^{\ast}_{\coh}(X)$ to $\DC^{\ast}_{\coh}(Y)$ for $\ast=$``\ \ '', $+$, $-$, $\bd$.
Let $\DC^{\ast}_{\aqcoh}(\widehat{X})$ be the full subcategory of $\DC^{\ast}(\widehat{X})$ consisting of objects having adically quasi-coherent cohomologies in all degrees ($\DC^{\ast}_{\aqcoh}(\widehat{X})$ may not be a triangulated category).
We have the diagram:
$$
\xymatrix{\DC^{\ast}_{\coh}(X)\ar[d]_{\RD f_{\ast}}\ar[r]^{\for}&\DC^{\ast}_{\aqcoh}(\widehat{X})\ar[d]^{\RD\widehat{f}_{\ast}}\\ \DC^{\ast}_{\coh}(Y)\ar[r]_{\for}&\DC^{\ast}(\widehat{Y})}\leqno{(\ast)}
$$
for $\ast=+$, $\bd$.
Note that we can extend the domain of the functor $\RD\widehat{f}_{\ast}$ to the whole $\DC(\widehat{X})$ due to the fact that $\RD\widehat{f}_{\ast}$ has finite cohomology dimension (\ref{cor-formalacyclicity4}); cf.\ \cite[C.D.\ Chap.\ 2, \S2, $\mathrm{n}^{\mathrm{o}}\ 2$, Cor.\ 2]{SGA4.5}.
Hence one can consider the diagram $(\ast)$ also for $\ast=$``\ \ '', $-$.

We are going to construct a natural transformation ({\em comparsion map})
$$
\rho=\rho_f\colon \for\circ\RD f_{\ast}\longrightarrow\RD\widehat{f}_{\ast}\circ\for.\leqno{(\dagger)}
$$
As is well-known, there exists a canonical natural transformation
$$
i^{-1}\circ\RD f_{\ast}\longrightarrow\RD\widehat{f}_{\ast}\circ j^{-1}.
$$
Indeed, for any object $M$ of $\DC(X)$ we represent $M$ by a complex $\mathscr{J}^{\bullet}$ consisting of injective $\O_X$-modules.
Then we have the following chain of canonical morphisms:
$$
i^{-1}\RD f_{\ast}M\stackrel{\sim}{\longrightarrow}i^{-1}f_{\ast}\mathscr{J}^{\bullet}\longrightarrow \widehat{f}_{\ast}j^{-1}\mathscr{J}^{\bullet}\stackrel{\sim}{\longleftarrow}\RD\widehat{f}_{\ast}j^{-1}M,
$$
where the first and the last arrows are quasi-isomorphisms; notice that the last quasi-isomorphy follows from that $j^{-1}\mathscr{J}^{\bullet}$ gives a flasque resolution of $j^{-1}M$.
In view of the fact that the maps $i$ and $j$ are flat as maps of locally ringed spaces, one can extend the above morphism to
$$
(\RD f_{\ast}M)^{\for}\longrightarrow\RD\widehat{f}_{\ast}M^{\for}.
$$
Moreover, since the formation of this morphism is canonical, we get the desired natural transformation $\rho$.
Therefore, we finally get the diagram as follows:
$$
\xymatrix@-2ex{
\DC^{\ast}_{\coh}(X)\ar[rr]^{\for}\ar[dd]_{\RD f_{\ast}}&&\DC^{\ast}_{\aqcoh}(\widehat{X})\ar[dd]^{\RD\widehat{f}_{\ast}}\\
&&\ \\
\DC^{\ast}_{\coh}(Y)\ar[rr]_{\for}&\ar@{=>}[ur]^{\rho}&\DC^{\ast}(\widehat{Y})\rlap{.}}\leqno{(\ast\ast)}
$$

Now our main theorem of this section is stated as follows:
\begin{thm}[GFGA comparison theorem]\label{thm-GFGAcom}
The diagram $(\ast\ast)$ is $2$-commutative, that is, $\rho$ gives a natural equivalence for $\ast=$``\ \ '', $+$, $-$, $\bd$.
\end{thm}

\subsection{The classical comparison theorem}\label{sub-GFGAcomcl}
Before proceeding to the proof of \ref{thm-GFGAcom}, let us mention a classical version of the comparison theorem (the generalized classical comparison theorem) along with the formalism of \cite[$\mathbf{III}$, \S4.1]{EGA}.

We continue with working in the situation as in \ref{sit-GFGA}.
Let $\mathscr{F}$ be a coherent $\O_X$-module, and consider the completion $\widehat{\mathscr{F}}=\widehat{\mathscr{F}}|_W$.
For $k\geq 0$ we set $\mathscr{F}_k=\mathscr{F}/\mathscr{I}^{k+1}\mathscr{F}$.
Then we have $\widehat{\mathscr{F}}=j^{-1}\varprojlim_k\mathscr{F}_k$.
In this situation we have the commutative diagram of cohomologies
$$
\xymatrix@C-10ex{\widehat{\RD^qf_{\ast}\mathscr{F}}\ar[rr]^{\rho_q}\ar[dr]_{\varphi_q}&&\RD^q\widehat{f}_{\ast}\widehat{\mathscr{F}}\ar[dl]^{\psi_q}\\ &\varprojlim_k\RD^qf_{\ast}\mathscr{F}_k}\leqno{(\ast)}
$$
for $q\geq 0$, constructed as follows: 

\medskip
{\sl Construction of $\varphi^q$}. The canonical map $\RD^qf_{\ast}\mathscr{F}\rightarrow\RD^qf_{\ast}\mathscr{F}_k$ factors through 
$$
(\RD^qf_{\ast}\mathscr{F})\otimes_{\O_Y}(\O_Y/\mathscr{I}^{k+1}\O_Y)\longrightarrow\RD^qf_{\ast}\mathscr{F}_k,
$$
from which we obtain $\varphi^q$ by passage to the projective limits $\varprojlim_k$.

\medskip
{\sl Construction of $\psi^q$}. We look at the commutative diagram
$$
\xymatrix@R-4ex{&\widehat{X}\ar[dd]^j\\ X_k\ar[ur]^{h_k}\ar[dr]_{\iota_k}\\ &X\rlap{,}}
$$
where $h_k$ and $\iota_k$ are the closed immersions.
Since we have the canonical isomorphism $\mathscr{F}_k\cong(\iota_k)_{\ast}(\iota_k)^{\ast}\mathscr{F}_k$, we have $\H^q(X_k,(\iota_k)^{\ast}\mathscr{F}_k)\cong\H^q(X,\mathscr{F}_k)$.
On the other hand, by \ref{prop-adqformalprop} (1) we know that $\widehat{\mathscr{F}}=j^{\ast}\mathscr{F}$.
Hence we have the map $\H^q(\widehat{X},\widehat{\mathscr{F}})\rightarrow\H^q(X,\mathscr{F}_k)$ and hence the map
$$
\H^q(\widehat{X},\widehat{\mathscr{F}})\longrightarrow\varprojlim_k\H^q(X,\mathscr{F}_k).
$$
Do the same for all schemes of the form $f^{-1}(V)$ with $V$ \'etale over $Y$ and apply {\bf \ref{ch-pre}}.\ref{cor-cohqcoh1-0} to obtain the desired map $\psi^q$.

\medskip
{\sl Construction of $\rho^q$}. This is essentially done in \S\ref{subsub-GFGAcomann}.
By \ref{thm-fini} we already know that $\RD^qf_{\ast}\mathscr{F}$ is a coherent $\O_Y$-module.
By virtue of \ref{prop-adqformalprop} (1) we have $\widehat{\RD^qf_{\ast}\mathscr{F}}\cong i^{\ast}\RD^qf_{\ast}\mathscr{F}$, and hence we get $\rho^q$ by the composition
$$
\widehat{\RD^qf_{\ast}\mathscr{F}}\cong i^{\ast}\RD^qf_{\ast}\mathscr{F}\longrightarrow\RD^q\widehat{f}_{\ast}(j^{\ast}\mathscr{F})\cong\RD^q\widehat{f}_{\ast}\widehat{\mathscr{F}}.
$$

\medskip
The diagram $(\ast)$ thus obtained is commutative; to see this, we replace $Y$ by open subsets to reduce to the commutativity of 
$$
\xymatrix@C-10ex{\widehat{\H^q(X,\mathscr{F})}\ar[rr]^{\rho^q}\ar[dr]_{\varphi^q}&&\H^q(\widehat{X},\widehat{\mathscr{F}})\ar[dl]^{\psi^q}\\ &\varprojlim_k\H^q(X,\mathscr{F}_k)\rlap{.}}\leqno{(\ast\ast)}
$$
But this follows from the commutativity of 
$$
\xymatrix@C-10ex{\H^q(X,\mathscr{F})\ar[rr]^{\rho^q}\ar[dr]_{\varphi^q}&&\H^q(\widehat{X},j^{\ast}\mathscr{F})\ar[dl]^{\psi^q}\\ &\H^q(X_k,\iota^{\ast}_k\mathscr{F})\rlap{,}}
$$
which is obvious by the fact $\iota_k=j\circ h_k$.

Now the classical (but generalized to our situation) comparison theorem is stated as follows:
\begin{thm}\label{thm-GFGAcomcl}
In the above situation, suppose, moreover, that $(Y,Z)$ is universally adhesive.
Then the morphisms $\rho^q$, $\varphi^q$ and $\psi^q$ in the diagram $(\ast)$ are bicontinuous isomorphisms, where the topology on $\widehat{\RD^qf_{\ast}\mathscr{F}}$ is the $\mathscr{I}$-adic topology, and the topologies on $\RD^q\widehat{f}_{\ast}\widehat{\mathscr{F}}$ and $\varprojlim_k\RD^qf_{\ast}\mathscr{F}_k$ are the ones given by the filtration by the kernels of the canonical maps to $\RD^qf_{\ast}\mathscr{F}_k$ for $k\geq 0$.
\end{thm}

\begin{proof}[Proof of Theorem {\rm \ref{thm-GFGAcom}} $\Rightarrow$ Theorem {\rm \ref{thm-GFGAcomcl}}]
We may assume that $Y$ is an affine scheme $Y=\Spec B$, and let $I$ be the ideal of $B$ corresponding to $\mathscr{I}$.
It suffices to show that the maps $\rho^q$, $\varphi^q$ and $\psi^q$ in the diagram $(\ast\ast)$ are bicontinuous isomorphisms.

By \ref{thm-GFGAcom} applied to the sheaf $\mathscr{F}$ (regarded as an object of $\DC^{\bd}_{\coh}(X)$), the map 
$$
\rho^q\colon\widehat{\H^q(X,\mathscr{F})}\longrightarrow\H^q(\widehat{X},\widehat{\mathscr{F}})\leqno{(\dagger)}
$$
is an isomorphism.
Since $\widehat{\mathscr{F}}$ is the projective limit of $\{\mathscr{F}_k\}$, by {\bf \ref{ch-pre}}.\ref{cor-ML5} the map $\psi^q$ is surjective.
On the other hand, the filtration $F^{\bullet}=\{F^n\}_{n\geq 0}$ on the $B$-module $\H^q(X,\mathscr{F})$ given by 
$$
F^n=\image(\H^q(X,I^n\mathscr{F})\rightarrow\H^q(X,\mathscr{F}))=\ker(\H^q(X,\mathscr{F})\rightarrow\H^q(X,\mathscr{F}_{n-1}))
$$
is $I$-good by \ref{prop-fini4}.
Hence the map $\varphi^q$ is injective; indeed, $\widehat{\H^q(X,\mathscr{F})}$ is topologically isomorphic to the projective limit of $\H^q(X,\mathscr{F})/F^{k+1}$ that are mapped injectively into $\H^q(X,\mathscr{F}_k)$. 
Hence the injectivity of $\varphi^q$ follows from the left-exactness of $\varprojlim_k$.

Therefore, since $(\dagger)$ is an isomorphism, we deduce that $\varphi^q$ and $\psi^q$ are isomorphisms.
But, then, once $\varphi^q$ is known to be isomorphic, it is automatically topologically isomorphic, since the filtration $\{F^{\bullet}\}$ as above is $I$-good (and hence both $\widehat{\H^q(X,\mathscr{F})}$ and $\varprojlim_k\H^q(X,\mathscr{F}_k)$ are complete with respect to the same filtration).
Also, since $\psi^q$ is now known to be an isomorphism, it is automatically a bicontinuous isomorphism due to the definition of the topologies.
Therefore, we deduce the other assertions of the theorem.
\end{proof}

\begin{rem}\label{rem-GFGAcompfmethod}{\rm
Notice that our proof of \ref{thm-GFGAcomcl} differs in structure from the proof in \cite[$\mathbf{III}$, (4.1.5)]{EGA}.
In [loc.\ cit.] the proof was given by showing the morphisms $\varphi_q$ and $\psi_q$ in the diagram $(\ast)$ are isomorphisms; the isomorphy of the former was established by the $I$-goodness of the filtration $F^{\bullet}$ as above, and the latter by the condition {\bf (ML)}\index{Mittag-Leffler condition} (cf.\ {\bf \ref{ch-pre}}, \S\ref{subsub-ML}) for the projective system of the cohomologies.
But here we first prove that $\rho_q$ is an isomorphism (disregarding the topologies); then it turns out, as we saw above, that the isomorphy of the rests is almost automatic.}
\end{rem}

\subsection{Proof of Theorem \ref{thm-GFGAcom}}\label{sub-proofGFGAcom}
\subsubsection{Reduction process}\label{subsub-GFGAcompf1}
By a similar argument as in \ref{cor-fini} it suffices to show the theorem only in the case $\ast=\bd$.
We may assume that $Y$ is affine $Y=\Spec B$, and let $I\subseteq B$ the ideal corresponding to $\mathscr{I}$.
The pair $(B,I)$ is universally adhesive, and the ring $B$ is universally coherent.
We apply \ref{prop-proofmethod} to the following situation: For each object $f\colon X\rightarrow\Spec B$ of $\PAs_B$, we define 
$$
\DC_f\subseteq\DC^{\bd}_{\coh}(X)
$$
to be the full subcategory consisting of objects $M$ such that the comparison morphism 
$$
\rho_f(M)\colon(\RD f_{\ast}M)^{\for}\longrightarrow\RD\widehat{f}_{\ast}M^{\for}
$$
is an isomorphism in $\DC(\widehat{Y})$.

We claim that, provided that the condition {\bf (C2)} of \ref{prop-proofmethod} is verified, the other conditions are fulfilled, whence finishing the proof of \ref{thm-GFGAcom}.
The condition {\bf (C0)} is trivially satisfied, while for {\bf (C1)}, consider the morphism of exact triangles
$$
\xymatrix{(\RD f_{\ast}K)^{\for}\ar[r]\ar[d]_{\cong}&(\RD f_{\ast}L)^{\for}\ar[r]\ar[d]&(\RD f_{\ast}M)^{\for}\ar[r]^(.65){+1}\ar[d]^{\cong}&\ \\ \RD\widehat{f}_{\ast}K^{\for}\ar[r]&\RD\widehat{f}_{\ast}L^{\for}\ar[r]&\RD\widehat{f}_{\ast}M^{\for}\ar[r]^(.65){+1}&,}
$$
by which the condition {\bf (C1)} is immediately verified.

Next, let us check {\bf (C3)}; Let $M$ be an object of $\DC^{\bd}_{\coh}(\til{X})=\DC_g$.
Then $\RD\pi_{\ast}M$ belongs to $\DC^{\bd}_{\coh}(X)$ (by \ref{thm-fini}).
Since $\pi$ is projective, by our assumption that {\bf (C2)} is true, we have $(\RD\pi_{\ast}M)^{\for}\cong\RD\widehat{\pi}_{\ast}M^{\for}$.
Then
$$
(\RD f_{\ast}\RD\pi_{\ast}M)^{\for}=(\RD g_{\ast} M)^{\for}\cong\RD\widehat{g}_{\ast}M^{\for}
=\RD\widehat{f}_{\ast}\RD\widehat{\pi}_{\ast}M^{\for}\cong\RD\widehat{f}_{\ast}(\RD\pi_{\ast}M)^{\for},
$$
whence verifying {\bf (C3)}.

Finally, if $\iota$ is a closed immersion as in the diagram of {\bf (C4)}, the comparison map $(\RD\iota_{\ast}M)^{\for}\rightarrow\RD\widehat{\iota}_{\ast}M^{\for}$ is shown to be a quasi-isomorphism as follows.
By induction with respect to $\amp(M)$, taking shifts and distinguished triangles of the form 
$$
\tau^{\leq n}M\longrightarrow M\longrightarrow\tau^{\geq n+1}M\stackrel{+1}{\longrightarrow}
$$
into account (and using {\bf (C1)} verified above), we may assume that $M$ is concentrated in degree $0$.
Then $M$ is represented by a coherent sheaf $\mathscr{F}$ on $X$.
By {\bf \ref{ch-pre}}.\ref{thm-vanishcohaff-0} (2) and \ref{thm-formalacyclicity} (2) the $q$-th cohomologies of the both sides are $\widehat{\iota_{\ast}\mathcal{H}^q(\mathscr{F})}$ (the completion of $\iota_{\ast}\mathcal{H}^q(\mathscr{F})$) and $\widehat{\iota}_{\ast}\mathcal{H}^q(\mathscr{F}^{\for})$, which are isomorphic to each other, for the completion functor is exact for coherent sheaves.
Hence the condition {\bf (C4)} is verified.

\subsubsection{Projective case}\label{subsub-GFGAcompf2}
Thus we may restrict to the projective case; clearly, by the similar reduction process as in \ref{subsub-reductionprocessgenSerre}, using {\bf \ref{ch-pre}}.\ref{thm-vanishcohaff-0} (2) and \ref{thm-formalacyclicity} (2) in a similar way as above, we may assume $X=\P^r_B$.
Let $M$ be an object of $\DC^{\bd}_{\coh}(X)$.
As we saw in checking {\bf (C4)} above, we may assume that $M$ is concentrated in degree $0$ and hence is represented by a coherent sheaf $\mathscr{F}$ on $X$.

By \ref{lem-inftyrepsplitbundles} we have an exact sequence 
$$
\cdots\longrightarrow\mathscr{E}_1\longrightarrow\mathscr{E}_0\longrightarrow\mathscr{F}\longrightarrow 0
$$
by the sheaves $\mathscr{E}_i$ of the form $\O_X(n_i)^{\oplus k_i}$.
For computing $q$-cohomologies for a fixed $q$, we can truncate the complex $\{\mathscr{E}_{\bullet}\}$ into a bounded complex
$$
\cdots\longrightarrow 0\longrightarrow\mathscr{E}_{q+1}\longrightarrow\cdots\longrightarrow\mathscr{E}_1\longrightarrow\mathscr{E}_0\longrightarrow0\longrightarrow\cdots,
$$
and $\mathscr{F}$ can be replaced by this complex.
Hence it suffices to show the theorem in the case $\mathscr{F}$ is a line bundle $\O_X(n)$.

Set $\mathscr{F}=\O_X(n)$, and let $X_k$ and $\mathscr{F}_k=\mathscr{F}/\mathscr{I}^{k+1}\mathscr{F}$ be defined as in \ref{sit-GFGA}.
The cohomology $H=\H^q(X,\mathscr{F})$ (resp.\ $H_k=\H^q(X_k,\mathscr{F}_k)$) is a free $B$-module of finite rank (resp.\ a free $B_k$-module of finite rank); obviously $H_k=H/I^{k+1}H$ due to the explicit description of the cohomology group {\bf \ref{ch-pre}}.\ref{thm-cohprojsp-0}.
Hence we have
$$
\widehat{\H^q(X,\mathscr{F})}=\varprojlim_k\H^q(X_k,\mathscr{F}_k).
$$
Now since the diagram $(\ast\ast)$ in \S\ref{sub-GFGAcomcl} is commutative, it suffices to show that
$$
\H^q(\widehat{X},\varprojlim_k\mathscr{F}_k)\longrightarrow\varprojlim_k\H^q(X_k,\mathscr{F}_k)
$$
is an isomorphism.
But since the projective system $\{\H^q(X_k,\mathscr{F}_k)\}_{k\geq 0}$ is strict for any $q$, the desired isomorphy follows from {\bf \ref{ch-pre}}.\ref{cor-ML5}.
Therefore, the theorem in the projective case is proved and hence, by virtue of \ref{prop-proofmethod}, the proof of Theorem \ref{thm-GFGAcom} is finished. \hfill$\square$
\index{GFGA!GFGA comparison theorem@--- comparison theorem|)}

\subsection{Comparison of Ext modules}\label{sub-GFGAcomext}
\begin{prop}\label{prop-GFGAcom1}
Let $(B,I)$ be a universally adhesive pair with $B$ universally coherent, and $f\colon X\rightarrow Y=\Spec B$ a proper morphism of finite presentation of algebraic spaces.
Then for $M\in\obj(\DC^-_{\coh}(X))$ and $N\in\obj(\DC^+_{\coh}(X))$ we have the canonical isomorphism
$$
\RD\Hom_{\O_X}(M,N)\cong\RD\Hom_{\O_{\widehat{X}}}(M^{\for},N^{\for})
$$
in $\DC^+(B)$ $($cf.\ {\rm {\bf \ref{ch-pre}}, \S\ref{subsub-ringedspacemodulestrcoh}} for the $B$-module structures on both sides$)$.
\end{prop}

\begin{proof}
First notice that by {\bf \ref{ch-pre}}.\ref{prop-schpairadh3} (2) we have $\RD\lHom_{\O_X}(M,N)\in\DC^+_{\coh}(X)$.
Hence by \ref{thm-GFGAcom} we have 
$$
(\RD f_{\ast}\RD\lHom_{\O_X}(M,N))^{\for}\cong\RD\widehat{f}_{\ast}(\RD\lHom_{\O_X}(M,N))^{\for}.
$$
On the other hand, since $j\colon\widehat{X}\rightarrow X$ is flat,
$$
(\RD\lHom_{\O_X}(M,N))^{\for}\cong\RD\lHom_{\O_{\widehat{X}}}(M^{\for},N^{\for})
$$
(\cite[$\mathbf{0}_{\mathbf{III}}$, (12.3.5)]{EGA}).
The result follows from this.
\end{proof}

\begin{cor}\label{cor-GFGAcom11}
Let $(B,I)$ be a complete universally adhesive pair with $B$ universally coherent, and $X$ a proper $B$-algebraic space of finite presentation.
Then the comparison functor\index{functor!comparison functor@comparison ---}
$$
\DC^{\bd}_{\coh}(X)\stackrel{\for}{\longrightarrow}\DC^{\bd}(\widehat{X})
$$
is fully faithful.
\end{cor}

\begin{proof}
What to prove is that the canonical map
$$
\Hom_{\DC(X)}(M,N)\longrightarrow\Hom_{\DC(\widehat{X})}(M^{\for},N^{\for})
$$
is bijective for $M,N\in\DC^{\bd}_{\coh}(X)$.
First notice that the left-hand side is isomorphic to $\H^0(\RD\Hom_{\O_X}(M,N))$, which is a finitely presented $B$-module by \ref{cor-fini}.
Since $B$ is complete, it is complete by {\bf \ref{ch-pre}}.\ref{prop-btarf1} (1).
Then by \ref{prop-GFGAcom1} we have $\H^0(\RD\Hom_{\O_X}(M,N))\cong\H^0(\RD\Hom_{\O_{\widehat{X}}}(M^{\for},N^{\for}))$ (as abelian groups); but the latter module is nothing but $\Hom_{\DC(\widehat{X})}(M^{\for},N^{\for})$.
\end{proof}



\section{GFGA existence theorem}\label{sec-GFGAexist}
\index{GFGA!GFGA existence theorem@--- existence theorem|(}
\subsection{Announcement of the theorem.}\label{sub-GFGAexistann}
\begin{sit}\label{sit-GFGAex}{\rm 
Let $B$ be a t.u.\ adhesive ring\index{t.u.a. ring@t.u.\ adhesive ring} (\ref{dfn-tuaringadmissible} (2)) with a finitely generated ideal of definition $I\subseteq B$, and suppose that $B$ is topologically universally coherent\index{coherent!topologically universally coherent@topologically universally ---} {\rm ({\bf \ref{ch-pre}}.\ref{dfn-cohringsmodules5})} with respect to $I$.
Let $f\colon X\rightarrow Y=\Spec B$ be a proper morphism of algebraic spaces of finite presentation.
We will use the notation as in \ref{sit-GFGA}.
Note that the algebraic space $X$ is universally cohesive\index{algebraic space!universally cohesive algebraic space@universally cohesive ---} ({\bf \ref{ch-pre}}.\ref{dfn-universallycohesive2}), and the formal algebraic space $\widehat{X}$ is universally adhesive and universally cohesive\index{formal algebraic space!universally cohesive formal algebraic space@universally cohesive ---}\index{cohesive!universally cohesive@universally --- (formal algebraic space)} (\ref{prop-cohformalschemes2}).}
\end{sit}

An important special case of the above situation is as follows: $B$ and $I$ are as above, and $B$ is $I$-torsion free or, more generally, $B$ is topologically finitely presented over a t.u.\ adhesive ring of this kind ({\bf \ref{ch-pre}}.\ref{thm-pf02a} (2)).
For example, $B$ can be a topologically finitely presented $V$-algebra with $I=aB$, where $V$ is an $a$-adically complete valuation ring\index{valuation!valuation ring@--- ring!a-adically complete valuation ring@$a$-adically complete --- ---} of arbitrary height.

In this section we are going to prove the following theorem:
\begin{thm}[GFGA existence theorem]\label{thm-GFGAexist}
In the situation as in {\rm \ref{sit-GFGAex}}, the comparison functor\index{functor!comparison functor@comparison ---}
$$
\DC^{\bd}_{\coh}(X)\stackrel{\for}{\longrightarrow}\DC^{\bd}_{\coh}(\widehat{X})\leqno{(\ast)}
$$
is an exact equivalence of triangulated categories.
\end{thm}

Notice that we have already shown in \ref{cor-GFGAcom11} that the functor $(\ast)$ is fully faithful.
Hence what to prove here is that the functor is essentially surjective, in other words, every object of $\DC^{\bd}_{\coh}(\widehat{X})$ is {\em algebraizable}\index{algebraizable}.

\subsection{Proof of Theorem \ref{thm-GFGAexist}}\label{sub-pfexist}
\subsubsection{Modification of the carving method}\label{subsub-pfmodification}
\index{carving method|(}
Theorem \ref{thm-GFGAexist} will be proved similarly to \ref{thm-fini} and \ref{thm-GFGAcom} by means of the carving method introduced in \ref{sub-proofmethod} with a slight modification indicated as follows:

We consider the category $\PAs_B$ $($full subcategory of the cateogory $\As_B$ of $B$-algebraic spaces$)$ consisting of proper and finitely presented $B$-algebraic spaces.
For any object $f\colon X\rightarrow\Spec B$ we denote by $\widehat{f}\colon\widehat{X}\rightarrow\Spf B$ its $I$-adic completion.
\begin{prop}[Carving method (formal version)]\label{prop-proofmethodvar}
Suppose for each object $f\colon X\rightarrow\Spec B$ of $\PAs_B$ we are given a full subcategory 
$$
\DC_f\subseteq\DC^{\bd}_{\coh}(\widehat{X})
$$
such that the following conditions are satisfied$:$
\begin{itemize}
\item[{\bf (C0)}] the zero object $0$ belongs to $\DC_f$, and $\DC_f$ is stable under isomorphisms in $\DC^{\bd}_{\coh}(\widehat{X});$
\item[{\bf (C1)}] let $K\rightarrow L\rightarrow M\stackrel{+1}{\rightarrow}$ be a exact triangle in $\DC^{\bd}_{\coh}(\widehat{X});$
if two of $K$, $L$, and $M$ are in $\DC_f$, then so is the rest$;$
\item[{\bf (C2)}] if $f$ is projective, then $\DC_f=\DC^{\bd}_{\coh}(\widehat{X});$
\item[{\bf (C3)}] consider a morphism $\pi$ in $\PAs_B$, which amounts to the same as a commutative diagram of the form
$$
\xymatrix@C-5ex{\til{X}\ar[rr]^{\pi}\ar[dr]_g&&X\ar[dl]^f\\ &\Spec B\rlap{,}}
$$
and suppose $g$ and $\pi$ are projective$;$
then $\RD\widehat{\pi}_{\ast}$ maps $\DC_g=\DC^{\bd}_{\coh}(\widehat{\til{X}})$ to $\DC_f;$
\item[{\bf (C4)}] consider a closed immersion $\iota$ in $\PAs_B$, that is, a commutative diagram
$$
\xymatrix@C-5ex{Z\ar[rr]^{\iota}\ar[dr]_g&&X\ar[dl]^f\\ &\Spec B\rlap{,}}
$$
where $\iota$ is a closed immersion$;$
then if $\DC_g=\DC^{\bd}_{\coh}(\widehat{Z})$, the functor $\RD\widehat{\iota}_{\ast}$ maps $\DC_g$ to $\DC_f$.
\end{itemize}
Then we have $\DC_f=\DC^{\bd}_{\coh}(\widehat{X})$ for any object $f\colon X\rightarrow\Spec B$ of $\PAs_B$. \hfill$\square$
\end{prop}

The proof can be done parallel to that of \ref{prop-proofmethod}; we only indicate the points of changes and leave the checking of details to the reader:

\medskip
On {\sc Claim 1.} Since $\iota$ is (automatically) finitely presented, the ideal $\mathscr{J}$ is of finite type, and hence $\mathscr{J}\O_{\widehat{X}}$ (closed by {\bf \ref{ch-pre}}.\ref{cor-qconsistency11} and {\bf \ref{ch-pre}}.\ref{cor-propARconseq1-2}) is the defining ideal of the induced closed immersion $\widehat{\iota}\colon\widehat{Z}\hookrightarrow\widehat{X}$ (cf.\ \ref{cor-closedimmformal41}).
In particular, it is a square zero ideal.
Then all the rest of the proof can be done parallel; we use \ref{prop-schformalpairadh31} instead of {\bf \ref{ch-pre}}.\ref{prop-schpairadh31}, and \ref{thm-formalacyclicity} (2) instead of {\bf \ref{ch-pre}}.\ref{thm-vanishcohaff-0} (2).

\medskip
On {\sc Claim 2.} As above, $\mathscr{J}\O_{\widehat{X}}$ is the defining ideal of $\widehat{Z}$.
Using \ref{prop-schformalpairadh31} instead of {\bf \ref{ch-pre}}.\ref{prop-schpairadh31}, one can prove the claim in our version similarly; here we need to show the following lemma, which we can show by an argument similar to that in the proof of \cite[$\mathbf{III}$, (5.3.4)]{EGA}, using our already proven comparison theorem \ref{thm-GFGAcomcl}:
\begin{lem}\label{lem-EGAIII534}
Let $f\colon Z\rightarrow X$ be a proper morphism in the category $\PAs_B$, and $\mathscr{J}\subseteq\O_X$ a coherent ideal.
Set $U=X\setminus V(\mathscr{J})$, and suppose that $f\colon f^{-1}(U)\rightarrow U$ is an isomorphism.
Then for any coherent $\O_{\widehat{X}}$-module $\mathscr{F}$ there exists an integer $n>0$ such that the kernel and the cokernel of the map $\mathscr{F}\rightarrow \widehat{f}_{\ast}\widehat{f}^{\ast}\mathscr{F}$ is annihilated by $\mathscr{J}^n$. \hfill$\square$
\end{lem}

All the rest of the proof goes just parallel, involving the carving lemma (\ref{lem-proofmethod}) in a similar manner as before.
\index{carving method|)}

\subsubsection{Reduction process}\label{subsub-pfexistred}
Thus it suffices to show that the conditions in \ref{prop-proofmethodvar} are satisfied, when we put $\DC_f$ to be the full subcategory of $\DC^{\bd}_{\coh}(\widehat{X})$ consisting of objects $M$ that are algebraizable\index{algebraizable}, that is, there exists $M_0\in\obj(\DC^{\bd}_{\coh}(X))$ such that $M\cong M^{\for}_0$.

The condition {\bf (C0)} is trivially satisfied.
The following proposition verifies the condition {\bf (C1)}:
\begin{prop}\label{prop-GFGAexist2}
Consider the situation as in {\rm \ref{sit-GFGAex}}, and let
$$
K\longrightarrow L\longrightarrow M\stackrel{+1}{\longrightarrow}
$$
be an exact triangle in $\DC^{\bd}_{\coh}(\widehat{X})$.
If two of $K$, $L$, and $M$ are algebraizable, then so is the rest.
\end{prop}

\begin{proof}
Suppose $K$ and $L$ are algebraizable, and take $K_0$ and $L_0$ in $\DC^{\bd}_{\coh}(X)$ such that $K\cong K^{\for}_0$ and $L\cong L^{\for}_0$.
Since $\Hom_{\DC^{\bd}_{\coh}(X)}(K_0,L_0)\cong\Hom_{\DC^{\bd}_{\coh}(\widehat{X})}(K,L)$ (\ref{cor-GFGAcom11}), we have a map $K_0\rightarrow L_0$ that induces $K\rightarrow L$.
Embed $K_0\rightarrow L_0$ into a exact triangle 
$$
K_0\longrightarrow L_0\longrightarrow M_0\stackrel{+1}{\longrightarrow}
$$
in $\DC^{\bd}_{\coh}(X)$.
Then since the functor $\for$ is exact, it induces the exact triangle
$$
K\longrightarrow L\longrightarrow M^{\for}_0\stackrel{+1}{\longrightarrow}\rlap{,}
$$
whence $M^{\for}_0\cong M$ (\cite[Chap.\ II, 1.2.4]{Verd1}), that is, $M$ is algebraizable.
The other cases can be reduced to the above case by shift.
\end{proof}

Assuming {\bf (C2)}, let us verify {\bf (C3)}.
Since we already know $\DC_g=\DC^{\bd}_{\coh}(\widehat{\til{X}})$, we may begin with an object of the form $M^{\for}$ with $M\in\DC^{\bd}_{\coh}(\til{X})$.
Then by \ref{thm-GFGAcom} we have $\RD\widehat{\pi}_{\ast}M^{\for}\cong(\RD\pi_{\ast}M)^{\for}$, which is obviously an object of $\DC^{\bd}_{\coh}(\widehat{X})$ (since $\RD\pi_{\ast}M$ belongs to $\DC^{\bd}_{\coh}(X)$ by \ref{thm-fini}; cf.\ \ref{rem-GFGAcomfunctordef}) and hence belongs to $\DC_f$.

The condition {\bf (C4)} can be checked similarly.

\subsubsection{Projective case}\label{subsub-GFGAexistproofproj}
Now let us prove \ref{thm-GFGAexist} in case $f$ is a projective morphism (hence $X$ is a scheme); by what we have seen above, this finishes the proof of \ref{thm-GFGAexist}.
Clearly, by a similar reduction process as in \ref{subsub-GFGAcompf2} (using \ref{thm-formalacyclicity} (2)), we may assume $X=\P^r_B$.
To proceed to the main routine of the proof, we need to show some preparatory results.
\begin{prop}\label{prop-GFGAexistproj1}
Consider the situation as in {\rm \ref{sit-GFGAex}}, and suppose there exists an $f$-ample line bundle $\mathscr{L}$ on $X$.
Then for any coherent $\O_{\widehat{X}}$-module $\mathscr{F}$ there exists an integer $N\geq 0$ such that$:$
\begin{itemize}
\item[{\rm (a)}] for $n\geq N$ we have $\RD^q\widehat{f}_{\ast}\mathscr{F}(n)=0$ for any $q>0$, and $\widehat{f}_{\ast}\mathscr{F}(n)$ is a coherent sheaf on $\widehat{Y};$
\item[{\rm (b)}] for $n\geq N$ the canonical morphism $\widehat{f}^{\ast}\widehat{f}_{\ast}\mathscr{F}(n)\rightarrow\mathscr{F}(n)$ is surjective.
\end{itemize}
\end{prop}

For the proof we need:
\begin{lem}\label{lem-GFGAexistproj2}
In the situation as in {\rm \ref{prop-GFGAexistproj1}}, suppose $B$ is $I$-torsion free.
Let $\mathscr{F}$ be a coherent $\O_{\widehat{X}}$-module.
There exists an integer $N\geq 0$ such that 
$$
\H^q(X_0,\mathscr{I}^k\mathscr{F}(n)/\mathscr{I}^{k+1}\mathscr{F}(n))=0
$$
for any $n\geq N$, $k\geq 0$, and $q>0$ and that $\bigoplus_{k\geq 0}\H^0(X_0,\mathscr{I}^k\mathscr{F}(n)/\mathscr{I}^{k+1}\mathscr{F}(n))$ for any $n\geq N$ is a graded module of finite presentation over the graded ring $\gr^{\bullet}_I(B)=\bigoplus_{k\geq 0}I^k/I^{k+1}$.
\end{lem}

\begin{proof}
Set $R=\gr^{\bullet}_I(B)=\bigoplus_{k\geq 0}I^k/I^{k+1}$.
We first claim that $R$ is a finitely presented $B$-algebra.
Indeed, if $R(B,I)$ denotes the Rees algebra\index{Rees algebra} ({\bf \ref{ch-pre}}, \S\ref{sub-Reescone}), we have $R=R(B,I)/IR(B,I)$.
Since $R(B,I)$ is $I$-torsion free and of finite type over $B$, it is of finite presentation over $B$.
Since $I$ is finitely generated, $R$ is of finite presentation over $B$, as desired.

Now set $X'=X\otimes_BR$, and let $j\colon X'\rightarrow X$ be the canonical morphism.
Since the ring $R$ is killed by $I$, the formal completion $\widehat{X}'$ coincides with $X'$.
Hence the sheaf $\widehat{j}^{\ast}\mathscr{F}(n)=\bigoplus_{k\geq 0}\mathscr{I}^k\mathscr{F}(n)/\mathscr{I}^{k+1}\mathscr{F}(n)$ is a coherent sheaf on $X'$, and thus we may apply \ref{thm-fini2}.
Set $M^q=\bigoplus_{k\geq 0}\H^q(X_0,\mathscr{I}^k\mathscr{F}(n)/\mathscr{I}^{k+1}\mathscr{F}(n))$.
Then $M^q=\H^q(X',\widehat{j}^{\ast}\mathscr{F}(n))$, and it follows from \ref{thm-fini2} that there exists $N$ such that for $n\geq N$ we have $M^q=0$ for $q>0$ and $M^0$ is a finitely presented $R$-module.
This is exactly what we wanted to show.
\end{proof}

\begin{proof}[Proof of Proposition {\rm \ref{prop-GFGAexistproj1}}]
First, let us claim that it is enough to prove the proposition in case $\mathscr{F}$ is $I$-torsion free.
Indeed, consider the exact sequence
$$
0\longrightarrow\mathscr{F}_{\Itor}(n)\longrightarrow\mathscr{F}(n)\longrightarrow(\mathscr{F}/\mathscr{F}_{\Itor})(n)\longrightarrow 0\leqno{(\ast)}
$$
for any $n$.
Since $\widehat{X}$ is universally adhesive, the $I$-torsion free $(\mathscr{F}/\mathscr{F}_{\Itor})(n)$ is of finite presentation, hence is coherent.
Hence $\mathscr{F}_{\Itor}$ is a coherent sheaf on $\widehat{X}$.
Since $\widehat{X}$ is quasi-compact, one can take $k\geq 0$ such that $I^{k+1}\mathscr{F}_{\Itor}=0$.
Hence if we put $X_k=(X,\O_X/I^{k+1}\O_X)$ and $B_k=B/I^{k+1}$, then $\mathscr{F}_{\Itor}$ is isomorphic to the direct image of a coherent sheaf $\mathscr{G}$ on $X_k$ (here we used {\bf \ref{ch-pre}}.\ref{lem-pf2} (1)); since $B_k$ is a finitely presented $B$-algebra, we may apply \ref{thm-fini2} and \ref{prop-fini1} to $\mathscr{G}$ and $X_k\rightarrow\Spec B_k$ to deduce that the properties (a) and (b) are valid for $\mathscr{F}_{\Itor}$.
Hence, considering the cohomology exact sequence for $(\ast)$ for a sufficiently large $n$, we deduce easily that (a) is valid also for $\mathscr{F}$; as for (b), by the commutative diagram with exact rows
$$
\xymatrix{0\ar[r]&\mathscr{F}_{\Itor}(n)\ar[r]&\mathscr{F}(n)\ar[r]&(\mathscr{F}/\mathscr{F}_{\Itor})(n)\ar[r]&0\\ &\widehat{f}^{\ast}\widehat{f}_{\ast}\mathscr{F}_{\Itor}(n)\ar[u]\ar[r]&\widehat{f}^{\ast}\widehat{f}_{\ast}\mathscr{F}(n)\ar[u]\ar[r]&\widehat{f}^{\ast}\widehat{f}_{\ast}(\mathscr{F}/\mathscr{F}_{\Itor})(n)\ar[u]\ar[r]&0\rlap{,}}
$$
the result follows from snake lemma.

Hence we may assume $\mathscr{F}$ is $I$-torsion free.
But then, in this case, one can further assume that $B$ is $I$-torsion free.
Indeed, set $B'=B/B_{\Itor}$; since $B_{\Itor}$ is finitely generated and $B$ satisfies {\bf (APf)}, $B'$ is $I$-adically complete ({\bf \ref{ch-pre}}.\ref{cor-qconsistency11})).
Since $B'$ is finitely presented and topologically finitely presented $B$-algebra, $B'$ is topologically universally coherent ({\bf \ref{ch-pre}}.\ref{prop-cohringsmodules52} (2)).
Consider the Cartesian diagram
$$
\xymatrix{\widehat{X}'\ \ar@{^{(}->}[r]\ar[d]&\widehat{X}\ar[d]\\ \Spf B'\,\ar@{^{(}->}[r]&\Spf B\rlap{,}}
$$
where the two horizontal arrows are closed immersions (\ref{cor-closedimmformal42}, \ref{prop-closedimmformal5} (2)).
Since $\mathscr{F}$ is $I$-torsion free, there exists a coherent sheaf $\mathscr{G}$ on $\widehat{X}'$ such that $\mathscr{F}$ coincides up to isomorphism with the direct image of $\mathscr{G}$.
Hence if the proposition is proved for $\mathscr{G}$ (and for $\widehat{X}'\rightarrow\Spf B'$), it is also true for $\mathscr{F}$; here we use \ref{thm-formalacyclicityder} (2) and the fact that any finitely presented $B'$-module is finitely presented as $B$-module.

Hence we may assume that $B$ is $I$-torsion free and thus may apply \ref{lem-GFGAexistproj2}.
Take $N\geq 0$ as in \ref{lem-GFGAexistproj2}, and fix $n\geq N$.
Starting from the result of \ref{lem-GFGAexistproj2}, one can show recursively that 
$$
\H^q(X_m,\mathscr{I}^k\mathscr{F}(n)/\mathscr{I}^{k+m+1}\mathscr{F}(n))=0
$$
for $m\geq 0$ and $q>0$.
By this one sees that the map
$$
\H^0(X_{m+1},\mathscr{I}^k\mathscr{F}(n)/\mathscr{I}^{k+m+2}\mathscr{F}(n))\longrightarrow\H^0(X_m,\mathscr{I}^k\mathscr{F}(n)/\mathscr{I}^{k+m+1}\mathscr{F}(n))
$$
is surjective.
Since $\H^q(\widehat{X},\mathscr{I}^k\mathscr{F}(n))=\varprojlim_{m\geq 0}\H^q(X_m,\mathscr{I}^k\mathscr{F}(n)/\mathscr{I}^{k+m+1}\mathscr{F}(n))$, we deduce by {\bf \ref{ch-pre}}.\ref{cor-ML5} that 
$$
\H^q(\widehat{X},\mathscr{I}^k\mathscr{F}(n))=0\leqno{(\ast)}
$$
for $q>0$ and $k\geq 0$, which already proves the first half of (a).

On the other hand, by $(\ast)$, \ref{lem-GFGAexistproj2} and {\bf \ref{ch-pre}}.\ref{prop-complpair2} we see that the filtration $\{\widehat{f}_{\ast}\mathscr{I}^k\mathscr{F}(n)\}_{k\geq 0}$ on $\widehat{f}_{\ast}\mathscr{I}^k\mathscr{F}(n)$ is $I$-good.
This means
$$
{\textstyle \widehat{f}_{\ast}\mathscr{F}(n)=\varprojlim_k\widehat{f}_{\ast}(\mathscr{F}(n)/\mathscr{I}^k\mathscr{F}(n))=\varprojlim_k\widehat{f}_{\ast}\mathscr{F}(n)/I^k\widehat{f}_{\ast}\mathscr{F}(n)},
$$
where the first equality follows from left exactness of $\widehat{f}_{\ast}$ and $(\ast)$.
Since $I^k\widehat{f}_{\ast}\mathscr{F}(n)\subseteq\widehat{f}\mathscr{I}^k\mathscr{F}(n)$, we have the exact sequence
\begin{equation*}
\begin{split}
0\longrightarrow\widehat{f}_{\ast}\mathscr{I}^k\mathscr{F}(n)/I^k\widehat{f}_{\ast}\mathscr{F}(n)\longrightarrow&\widehat{f}_{\ast}\mathscr{F}(n)/I^k\widehat{f}_{\ast}\mathscr{F}(n)\\ &\longrightarrow\widehat{f}_{\ast}\mathscr{F}(n)/\widehat{f}_{\ast}\mathscr{I}^k\mathscr{F}(n)\longrightarrow 0,
\end{split}
\end{equation*}
where the last term, equal to $\widehat{f}_{\ast}(\mathscr{F}(n)/\mathscr{I}^k\mathscr{F}(n))$, is known to be coherent due to \ref{thm-fini}.
Take a sufficiently large $l>0$ such that $\widehat{f}_{\ast}\mathscr{I}^l\mathscr{F}(n)\subseteq I^k\widehat{f}_{\ast}\mathscr{F}(n)$.
We have the exact sequence
\begin{equation*}
\begin{split}
0\longrightarrow I^k(\widehat{f}_{\ast}\mathscr{F}(n)/\widehat{f}_{\ast}\mathscr{I}^l\mathscr{F}(n))\longrightarrow&\widehat{f}_{\ast}\mathscr{I}^k\mathscr{F}(n)/\widehat{f}_{\ast}\mathscr{I}^l\mathscr{F}(n)\\ &\longrightarrow\widehat{f}_{\ast}\mathscr{I}^k\mathscr{F}(n)/I^k\widehat{f}_{\ast}\mathscr{F}(n)\longrightarrow 0,
\end{split}
\end{equation*}
where, by a similar reasoning as above, the first two terms are coherent.
Hence the last term is coherent, and thus we deduce that $\widehat{f}_{\ast}\mathscr{F}(n)/I^k\widehat{f}_{\ast}\mathscr{F}(n)$ is coherent for any $k$.
Hence by \ref{prop-adicqcoh55}, $\widehat{f}_{\ast}\mathscr{F}(n)$ is coherent, whence the rest of (a). 

To show (b), we take an integer $N'\geq 0$ such that in view of {\bf \ref{ch-pre}}.\ref{prop-genglobalsec-0} the sheaf $\mathscr{F}/\mathscr{I}\mathscr{F}(n)$ on $X_0$ is generated by global sections for $n\geq N'$.
If $n\geq\max(N,N')$, then (b) holds, for the map $\H^0(\widehat{X},\mathscr{F}(n))\rightarrow\H^0(\widehat{X},\mathscr{F}/\mathscr{I}\mathscr{F}(n))$ is surjective.
\end{proof}

\begin{cor}\label{cor-GFGAexistproj1}
In the above situation with $X=\P^r_B$, let $\mathscr{F}$ be a coherent sheaf on $\widehat{X}$.
Then we have the exact sequence of the following form$:$
$$
\cdots\longrightarrow\mathscr{E}_m\longrightarrow\mathscr{E}_{m-1}\longrightarrow\cdots\longrightarrow\mathscr{E}_1\longrightarrow\mathscr{E}_0\longrightarrow\mathscr{F}\longrightarrow 0, 
$$
where each $\mathscr{E}_i$ is a free $\O_{\widehat{X}}$-module of the form $\O_{\widehat{X}}(n_i)^{\oplus k_i}$.
\end{cor}

\begin{proof}
By \ref{prop-GFGAexistproj1} the sheaf $\mathscr{F}(n)$ for a large $n$ is generated by global sections, that is, the morphism $\widehat{f}^{\ast}\widehat{f}_{\ast}\mathscr{F}(n)\rightarrow\mathscr{F}(n)$ is surjective.
Now $\widehat{f}_{\ast}\mathscr{F}(n)$ is a coherent sheaf on $\widehat{Y}=\Spf B$, and hence is of the form $M^{\Delta}$ by a $B$-module $M$ of finite presentation (\ref{thm-adicqcohpre1}). 
There exists $n_0>0$ such that we have a surjective map $B^{n_0}\rightarrow M$.
Pulling it back onto $\widehat{X}$ and compose with the surjection $\widehat{f}^{\ast}\widehat{f}_{\ast}\mathscr{F}(n)\rightarrow\mathscr{F}(n)$, we get a surjection $\O_{\widehat{X}}^{\oplus k_0}\rightarrow\mathscr{F}(n)$, whence $\mathscr{E}_0\rightarrow\mathscr{F}\rightarrow 0$ as above.
The kernel of this map is again coherent, and hence one can repeat the procedure.
\end{proof}

Now we are going to finish the proof of the theorem in case $X=\P^r_B$, whence the proof of \ref{thm-GFGAexist}.
Let $M\in\obj(\DC^{\bd}_{\coh}(X))$.
By induction with respect to $\amp(M)$, taking shifts and distinguished triangles of the form 
$$
\tau^{\leq n}M\longrightarrow M\longrightarrow\tau^{\geq n+1}M\stackrel{+1}{\longrightarrow}
$$
into account (and using \ref{prop-GFGAexist2}), we may assume that $M$ is concentrated in degree $0$ and hence is represented by a coherent sheaf $\mathscr{F}$.
Then by \ref{cor-GFGAexistproj1} and \ref{prop-GFGAexist2} (since \ref{prop-GFGAexist2} implies that for any morphism $\mathscr{F}\rightarrow\mathscr{G}$ of algebraizable coherent sheaves the kernel and the cokernel are algebraizable) it suffices to show that the sheaves on $\widehat{X}$ of the form $\O_{\widehat{X}}(n)$ are algebraizable, but this is trivial. \hfill$\square$
\index{GFGA!GFGA existence theorem@--- existence theorem|)}

\subsection{Applications}\label{sub-GFGAexistappli}
\begin{prop}\label{prop-GFGAexistapp1}
Let $(B,I)$ be a complete t.u.\ adhesive pair with $B$ topologically universally coherent with respect to $I$.
Then the functor by $I$-adic formal completion 
$$
\PSch_B\longrightarrow\Fs_B,\qquad X\longmapsto\widehat{X},
$$
from the category of proper $B$-schemes of finite presentation to the category of formal $B$-schemes, is fully faithful.
\end{prop}

\begin{proof}
We need to show that for two objects $X$ and $Y$ of $\PSch_B$ the map
$$
\Hom_{\PSch_B}(X,Y)\longrightarrow\Hom_{\Fs_B}(\widehat{X},\widehat{Y})
$$
is bijective.
By considering the graphs one sees easily that this map is injective; indeed, any graph in $X\times_BY$ is a closed subspace of finite presentation, and hence one can apply \ref{cor-GFGAcom11}.
To show the surjectivity, take a morphism $\widehat{X}\rightarrow\widehat{Y}$ in the right-hand side.
This amounts to the same as taking the graph $\Gamma$; the morphism $\Gamma\rightarrow\widehat{X}\times_B\widehat{Y}$ is a closed immersion of finite presentation due to \ref{prop-closedimmform111}.
Let $\mathscr{J}$ be the ideal defining $\Gamma$.
Then this is a coherent ideal on $\widehat{X}\times_B\widehat{Y}$.
By \ref{prop-prodformal22} and \ref{thm-GFGAexist} we have a coherent ideal $\mathscr{J}_0$ on $X\times_BY$ (uniquely up to isomorphism) such that $\mathscr{J}^{\for}_0=\mathscr{J}$.
The ideal $\mathscr{J}_0$ defines a closed subspace $\Gamma_0$ of $X\times_BY$ such that $\widehat{\Gamma}_0\cong\Gamma$.
The projection $\Gamma_0\rightarrow X$ is isomorphism, as its completion is an isomorphism (Exercise \ref{exer-formalcompletionfaithflness}).
Hence there exists a composite map $X\rightarrow Y$ followed by the other projection, whose formal completion coincides with the morphism $\widehat{X}\rightarrow\widehat{Y}$ we started with.
\end{proof} 

\begin{prop}\label{prop-GFGAexistapp2}
Let $(B,I)$ be a complete t.u.\ adhesive pair with $B$ topologically universally coherent with respect to $I$, and $X$ a proper $B$-formal algebraic space of finite presentation.
Suppose there exists an invertible sheaf $\mathscr{L}$ on $X$ such that $\mathscr{L}_0$ is ample on $X_0$, where $X_k=(X,\O_X/I^{k+1}\O_X)$ and $\mathscr{L}_k=\mathscr{L}/I^{k+1}\mathscr{L}$ for $k\geq 0$.
Then $(X,\mathscr{L})$ is algebraizable\index{algebraizable}, that is, there exists a proper $B$-scheme $Y$ of finite presentation and an ample invertible sheaf $\mathscr{M}$ on $Y$ such that $\widehat{Y}\cong X$ and $\widehat{\mathscr{M}}\cong\mathscr{L}$.
$($In particular, the scheme $Y$ is projective.$)$ \hfill$\square$
\end{prop}

This can be shown in a similar way to the proof of \cite[$\mathbf{III}$, (5.4.5)]{EGA}; in the proof one can use \ref{prop-GFGAexistproj1} instead of \cite[$\mathbf{III}$, (5.2.3)]{EGA}.
We are not going to repeat it here and leave checking details to the reader.

\addcontentsline{toc}{subsection}{Exercises}
\subsection*{Exercises}
\begin{exer}[cf.\ {\cite[$\mathbf{III}$, (4.6.8)]{EGA}}]\label{exer-formalcompletionfaithflness}{\rm
Consider the commutative diagram
$$
\xymatrix@C-4ex@R-2ex{X\ar[rr]^f\ar[dr]_g&&Y\ar[dl]^h\\ &Z}
$$
of schemes, where $g$ and $h$ are proper of finite presentation.
Let $W\subseteq Z$ be a closed subscheme of finite presentation.
Suppose that the pair $(Z,W)$ is universally adhesive.
Set $\widehat{Z}=\widehat{Z}|_W$, $\widehat{Y}=\widehat{Y}|_{h^{-1}(W)}$, and $\widehat{X}=\widehat{X}|_{g^{-1}(W)}$, and consider the resulting diagram
$$
\xymatrix@C-4ex@R-2ex{\widehat{X}\ar[rr]^{\widehat{f}}\ar[dr]_{\widehat{g}}&&\widehat{Y}\ar[dl]^{\widehat{h}}\\ &\widehat{Z}}
$$
of formal schemes.
Then show that $\widehat{f}$ is an isomorphism (resp.\ a closed immersion) if and only if there exists an open neighborhood $U\subseteq Z$ of $W$ such that the morphism $g^{-1}(U)\rightarrow h^{-1}(U)$ induced by $f$ is an isomorphism (resp.\ a closed immersion).}
\end{exer}


\section{Finiteness theorem and Stein factorization}\label{sec-finiformal}
The first part \S\ref{sub-finiformalann} of this section is devoted to the finiteness theorem for proper morphisms between universally adhesive formal schemes.
The proof consists of three steps.
First, we prove the theorem for a formal scheme with {\em invertible} ideal of definition.
The theorem in this particular case has been essentially done in \cite{Ull}; we contain the full proof of this case for the reader's convenience. 
The second step deals with the case where the map is a so-called {\em admissible blow-up}\index{blow-up!admissible blow-up@admissible ---} (cf.\ {\bf \ref{ch-rigid}}, \S\ref{sec-blowup}).
The final step combines these results to show the general case.

The second part \S\ref{sub-steinfactadqformalsch} of this section discusses the  Stein factorizations for universally adhesive formal schemes.
This part refers to the general theory of the Stein factorization for schemes, contained in the appendix \S\ref{sec-steinfactorization} of this chapter.

\subsection{Finiteness theorem for proper morphisms}\label{sub-finiformalann}
\begin{thm}\label{thm-finiformal}
Let $f\colon X\rightarrow Y$ be a proper morphism of finite presentation of universally adhesive formal schemes\index{formal scheme!universally adhesive formal scheme@universally adhesive ---}\index{adhesive!universally adhesive@universally ---!universally adhesive formal scheme@--- --- formal scheme}, and suppose $Y$ is universally cohesive\index{formal algebraic space!universally cohesive formal algebraic space@universally cohesive ---}\index{cohesive!universally cohesive@universally --- (formal algebraic space)}.
Then the functor $\RD f_{\ast}$ maps $\DC^{\ast}_{\coh}(X)$ to $\DC^{\ast}_{\coh}(Y)$ for $\ast=$``\ \ '', $+$, $-$, $\bd$.
\end{thm}

Note that by \ref{cor-formalacyclicity4} the functor $\RD f_{\ast}$ is defined on $\DC_{\coh}(X)$ $($cf.\ \cite[C.D.\ Chap.\ 2, \S2, $\mathrm{n}^{\mathrm{o}}\ 2$, Cor.\ 2]{SGA4.5}$))$.
Clearly, it suffices to show the theorem in case $\ast=\bd$.
But then, by a standard reduction process as in \S\ref{subsub-finitudesproof} (by induction with respect to amplitudes\index{amplitude}), one reduces the theorem to the following one.
\begin{thm}\label{thm-finiformal2}
Let $f\colon X\rightarrow Y$ be a proper morphism of finite presentation of universally adhesive formal schemes, and suppose $Y$ is universally cohesive.
Then for any coherent $\O_X$-module $\mathscr{F}$ the sheaf $\RD^qf_{\ast}\mathscr{F}$ is coherent for any $q\geq 0$.
\end{thm}

\subsection{Proof of Theorem \ref{thm-finiformal}}
\subsubsection{Invertible ideal case}\label{subsub-finiformalproofspecial}
First, we are to show the theorem in the case where $Y$ has an {\em invertible} ideal of definition $\mathscr{I}$ and $\O_Y$ is $\mathscr{I}$-torsion free.
Since one may work locally on $Y$, we can assume, without loss of generality, that $Y$ is affine $Y=\Spf A$ and that $A$ has a principal ideal $I=(a)$ of definition and $A$ is $a$-torsion free.
Thus what we are to show in this paragraph is the following:
\begin{prop}\label{prop-finiformalspecial}
Let $(A,I)$ with $I=(a)$ be a complete t.u.\ adhesive pair such that $A$ is $a$-torsion free, and $f\colon X\rightarrow Y=\Spf A$ a proper morphism of finite presentation.
$($In this situation, $X$ and $Y$ are automatically universally adhesive, and $Y$ is universally cohesive.$)$
Then for any coherent $\O_X$-module $\mathscr{F}$, $\RD^qf_{\ast}\mathscr{F}$ is a coherent $\O_Y$-module for any $q\geq 0$.
\end{prop}

The theorem in this case is proved by P.\ Ullrich\index{Ullrich, P.} \cite{Ull}; the following argument is essentially the same as the one therein.
\begin{lem}\label{lem-finiformalspecial}
In order to show {\rm \ref{prop-finiformalspecial}} it suffices to prove the following statements for $q\geq 0:$
\begin{itemize}
\item[{\rm (a)}] the $A$-module $\H^q(X,\mathscr{F})$ is coherent $($cf.\ {\rm {\bf \ref{ch-pre}}.\ref{dfn-cohringsmodules1} (1)}$);$
\item[{\rm (b)}] for any $g\in A$ the canonical map 
$$
\H^q(X,\mathscr{F})\otimes_AA_{\{g\}}\longrightarrow\H^q(f^{-1}(V),\mathscr{F}),
$$
where $V=\Spf A_{\{g\}}\hookrightarrow Y=\Spf A$, is an isomorphism.
\end{itemize}
\end{lem}

Notice that the statement (b) is necessary, for we do not know, a priori, that the sheaf $\RD^qf_{\ast}\mathscr{F}$ is adically quasi-coherent of finite type.
\begin{proof}
Set $M=\H^q(X,\mathscr{F})$.
Then by (a) $M$ is $a$-adically complete ({\bf \ref{ch-pre}}.\ref{prop-btarf1} (1)), and $M^{\Delta}$ is a coherent $\O_Y$-module.
Since $\RD^qf_{\ast}\mathscr{F}$ is the sheafification of the presheaf given by $V\mapsto\H^q(f^{-1}(V),\mathscr{F})$, there exists a canonical map $\H^q(X,\mathscr{F})\rightarrow\Gamma(X,\RD^qf_{\ast}\mathscr{F})$.
By what we have seen in \S\ref{sub-adicqcohaff} we have a canonical morphism $M^{\Delta}\rightarrow\RD^qf_{\ast}\mathscr{F}$ of $\O_Y$-modules.
By (b) we see that the induced map between stalks at each point of $Y$ is an isomorphism, and thus we have $M^{\Delta}\cong\RD^qf_{\ast}\mathscr{F}$, whence the assertion.
\end{proof}

\begin{proof}[Proof of Proposition {\rm \ref{prop-finiformalspecial}}]
We are going to check (a) and (b) in \ref{lem-finiformalspecial}.
We fix a finite open covering $\mathscr{U}=\{U_{\lambda}\}_{\lambda\in\Lambda}$ of $X$ by affine subsets, and consider the \v{C}ech complex $K^{\bullet}=C^{\bullet}(\mathscr{U},\mathscr{F})$.
For any $q\geq 0$ we have $\H^q(K^{\bullet})=\H^q(X,\mathscr{F})$.
If $V=\Spf A_{\{g\}}$ as in (b) above, set $\mathscr{U}|_V=\{U_{\lambda}\cap f^{-1}(V)\}_{\lambda\in\Lambda}$.
Since $U_{\lambda}\rightarrow Y$ is an affine map for each $\lambda\in\Lambda$, $\mathscr{U}|_V$ is a finite open covering of $f^{-1}(V)$ by affine open subsets.
Thus the \v{C}ech complex $L^{\bullet}=C^{\bullet}(\mathscr{U}|_V,\mathscr{F}|_{f^{-1}(V)})$ gives rise to the cohomologies $\H^q(L^{\bullet})=\H^q(f^{-1}(V),\mathscr{F})$ for $q\geq 0$.
By \ref{cor-formalacyclicity4} there exists $q_0\geq 0$ such that for $q>q_0$ we have $\H^q(K^{\bullet})=\H^q(X,\mathscr{F})=0$ and $\H^q(L^{\bullet})=\H^q(f^{-1}(V),\mathscr{F})=0$.
Thus the conditions (a) and (b) are trivially satisfied for $q>q_0$.
Hence one can check the conditions by descending induction with respect to $q$.
In the following, we assume that (a) and (b) are true with $q$ replaced by $q+1$.

Let us verify (a).
Before doing it, notice that each member $K^q$ of the complex $K^{\bullet}$ is $a$-adically complete, for it is finitely generated over a topologically finitely presented $A$-algebra ({\bf \ref{ch-pre}}.\ref{prop-btarf1}).
By the adhesiveness the $a$-torsion part $K^q_{\ator}$ is again finitely generated over a topologically finitely presented $A$-algebra.
Hence there exists $n\geq 0$ such that $a^nK^q_{\ator}=0$.
Moreover, we have $C^{\bullet}(\mathscr{U},a^{k+1}\mathscr{F})=a^{k+1}K^{\bullet}$ and 
$$
0\longrightarrow C^{\bullet}(\mathscr{U},a^{k+1}\mathscr{F})\longrightarrow C^{\bullet}(\mathscr{U},\mathscr{F})\longrightarrow C^{\bullet}(\mathscr{U},\mathscr{F}_k)\longrightarrow 0
$$
is exact for any $k\geq 0$, where $\mathscr{F}_k=\mathscr{F}/a^{k+1}\mathscr{F}$, since each $U_{\lambda}$ is affine (cf.\ \ref{thm-formalacyclicity}).
Hence for $k\geq 0$ we have $C^{\bullet}(\mathscr{U},\mathscr{F}_k)=K^{\bullet}_k$ $(=K^{\bullet}/a^{k+1}K^{\bullet})$ (cf.\ {\bf \ref{ch-pre}}, \S\ref{subsub-paircomplexgen}).

Since $\mathscr{F}$ and $a^{k+1}\mathscr{F}$ are coherent, we know that $\H^{q+1}(K^{\bullet})$ and $\H^{q+1}(a^{k+1}K^{\bullet})$ are coherent $A$-modules by induction.
Since the sheaf $\mathscr{F}_k$ is coherent on the scheme $X_k=(X,\O_X/a^{k+1}\O_X)$ (\ref{prop-adicqcoh55}) and the map $f_k\colon X_k\rightarrow Y_k=\Spec A_k$ (where $A_k=A/a^{k+1}A$) is proper (\ref{prop-propermorformal0}), we know that the cohomology group $\H^q(K^{\bullet}_k)$ is a coherent $A_k$-module by \ref{thm-fini} and hence is also coherent as an $A$-module.
Now all the hypotheses of {\bf \ref{ch-pre}}.\ref{lem-finiformalspecial2} are satisfied, and thus we conclude that $H^q(K^{\bullet})$ is coherent.

Next we verify (b).
With the notation as above, what to show is that the map 
$$
\H^q(K^{\bullet})\otimes_AB\longrightarrow\H^q(L^{\bullet}),
$$
where $B=A_{\{g\}}$, is an isomorphism.
Notice that, since $(A,a)$ is t.u.\ adhesive, $B$ is flat over $A$.
The checking is done by applying {\bf \ref{ch-pre}}.\ref{lem-ullrichlemma3}.
To this end, what to show is that the maps $\H^{q+1}(a^{k+1}K^{\bullet})\otimes_AB\rightarrow\H^{q+1}(a^{k+1}L^{\bullet})$ and $\H^q(K^{\bullet}_k)\otimes_AB\rightarrow\H^q(L^{\bullet}_k)$ are isomorphisms for any $k\geq n$, where $n\geq 0$ is an integer such that $a^nK^{q+1}_{\ator}=0$ and $a^nL^{q+1}_{\ator}=0$.
The first map coincides with the canonical map $\H^{q+1}(X,a^{k+1}\mathscr{F})\otimes_AB\rightarrow\H^{q+1}(f^{-1}(V),a^{k+1}\mathscr{F})$, which is an isomorphism by induction.
The second map coincides with the canonical map $\H^q(X,\mathscr{F}_k)\otimes_{A_k}B_k\rightarrow\H^q(f^{-1}(V),\mathscr{F}_k)$ where $B_k=B/a^{k+1}B=(A_k)_g$, which is an isomorphism, since $\RD^qf_{k\ast}\mathscr{F}_k$ is coherent  by \ref{thm-fini}.
Hence the condition (b) is verified, and thus the proof of \ref{prop-finiformalspecial} is completed.
\end{proof}

\subsubsection{Blow-up case}\label{subsub-finiformalproofspecial2}
We suppose that $Y$ is affine $Y=\Spf A$ where $A$ is a t.u.\ adhesive ring with a finitely generated ideal of definition $I\subseteq A$.
Consider the blow-up\index{blow-up} $g\colon\Proj\bigoplus_{k\geq 0}I^k\rightarrow\Spec A$, which is proper, since $I$ is finitely generated.
Moreover, it is of finite presentation, since $A$ is universally adhesive and the structure sheaf of $\Proj\bigoplus_{k\geq 0}I^k$ is $I$-torsion free.
Consider the case $f\colon X=\widehat{\Proj\bigoplus_{k\geq 0}I^k}\rightarrow Y$, that is, $X$ is the $I$-adic completion of $\Proj\bigoplus_{k\geq 0}I^k$.
In this case, by \ref{thm-GFGAexist} there exists a coherent sheaf $\mathscr{G}$ such that $\widehat{\mathscr{G}}=\mathscr{F}$.
Due to \ref{thm-fini} $\RD^qg_{\ast}\mathscr{G}$ is coherent on $\Spec A$.
Then by \ref{thm-GFGAcomcl} we have $\RD^qf_{\ast}\mathscr{F}=\widehat{\RD^qg_{\ast}\mathscr{G}}$, which is coherent by \ref{prop-adicqcoh55}.
Hence the theorem is true in this case.

\begin{rem}\label{rem-remarkontheblowup}{\rm 
The map $f\colon X=\widehat{\Proj\bigoplus_{k\geq 0}I^k}\rightarrow Y$ as above is an example of so-called {\em admissible blow-ups}\index{blow-up!admissible blow-up@admissible ---}, which will be discussed thoroughly in {\bf \ref{ch-rigid}}, \S\ref{sec-blowup}.}
\end{rem}

\subsubsection{General case}\label{subsub-finiformalproofgeneral}
We may assume that $Y$ is affine $Y=\Spf A$ where $A$ is a t.u.\ adhesive ring with a finitely generated ideal of definition $I$ such that $\mathscr{I}=I^{\Delta}$.
Let $\pi_Y\colon Y'\rightarrow Y$ be the admissible blow-up along $I$ as in \S\ref{subsub-finiformalproofspecial2}, and $X'\hookrightarrow X\times_YY'$ be the closed formal subscheme defined by dividing out $\mathscr{I}$-torsions (the so-called {\em strict transform}\index{strict transform} of $X$).
Since $X\times_YY'$ is again universally adhesive, $X'\hookrightarrow X\times_YY'$ is of finite presentation, and hence the induced map $f'\colon X'\rightarrow Y'$ is finitely presented and proper.
Denote the map $X'\rightarrow X$ by $\pi_X$; it follows that this map is also an admissible blow-up (cf.\ {\bf \ref{ch-rigid}}.\ref{prop-blowups15}).
Thus we have the following commutative square:
$$
\xymatrix{X'\ar[d]_{\pi_X}\ar[r]^{f'}&Y'\ar[d]^{\pi_Y}\\ X\ar[r]_f&Y\rlap{.}}
$$
Note that, since the ideal of definition $\mathscr{I}\O_{Y'}$ is invertible, $f'\colon X'\rightarrow Y'$ satisfies the hypothesis of \S\ref{subsub-finiformalproofspecial}, and hence the theorem is true for $f'$.

Let $M\in\obj(\DC^{\bd}_{\coh}(X))$.
By induction with respect to the amplitudes (as in \S\ref{subsub-finitudesproof}), we may assume that $M$ is concentrated in degree $0$.
We want to show that $\RD f_{\ast}M$ belongs to $\DC^{\bd}_{\coh}(Y)$ (by \ref{cor-formalacyclicity4} we already know that $\RD f_{\ast}M$ belongs to $\DC^{\bd}(Y)$).
Consider the object $\RD\pi_{X\ast}(\tau^{\geq 0}\LD\pi^{\ast}_XM)$ and the distinguished triangle 
$$
M\longrightarrow\RD\pi_{X\ast}(\tau^{\geq 0}\LD\pi^{\ast}_XM)\longrightarrow N\stackrel{+1}{\longrightarrow}.
$$
Notice that by \ref{prop-schformalpairadh31}, \S\ref{subsub-finiformalproofspecial2} and by the fact that $\pi_X$ is an admissible blow-up, the object in the middle, and hence $N$ also, is an object of $\DC^{\bd}_{\coh}(X)$.

Now since $\pi_X$ is also an admissible blow-up along $I$, all the cohomologies $\mathcal{H}^q(N)$ for $q\in\Z$ are $I$-torsion sheaves, and hence there exists $k\geq 0$ such that $I^{k+1}\mathcal{H}^q(N)=0$ for any $q$.
Consider the map of schemes $f_k\colon X_k=(X,\O_X/I^{k+1}\O_X)\rightarrow Y_k=\Spf A_k$ where $A_k=A/I^{k+1}$, and the commutative diagram
$$
\xymatrix{X\ar[r]^f&Y\\ X_k\ar[u]^{\iota_X}\ar[r]_{f_k}&Y_k\ar[u]_{\iota_Y}\rlap{,}}
$$
where the vertical maps are closed immersions.
Since $I^{k+1}\mathcal{H}^q(N)=0$ for any $q$, the canonical morphism $N\rightarrow\RD\iota_{X\ast}(\tau^{\geq 0}\LD\iota^{\ast}N)$ is an isomorphism in $\DC^{\bd}_{\coh}(X)$.

Now, in order to show the theorem, in view of the distinguished triangle
$$
\RD f_{\ast}M\longrightarrow\RD f_{\ast}\RD\pi_{X\ast}(\tau^{\geq 0}\LD\pi^{\ast}_XM)\longrightarrow \RD f_{\ast}N\stackrel{+1}{\longrightarrow},
$$
we need to show that the second and the third terms belong to $\DC^{\bd}_{\coh}(Y)$.
As for the second term, since $\RD f_{\ast}\RD\pi_{X\ast}(\tau^{\geq 0}\LD\pi^{\ast}_XM)\cong\RD \pi_{Y\ast}\RD f'_{\ast}(\tau^{\geq 0}\LD\pi^{\ast}_XM)$ and $\tau^{\geq 0}\LD\pi^{\ast}_XM$ belongs to $\DC^{\bd}_{\coh}(X')$, we deduce from what we have proved in \S\ref{subsub-finiformalproofspecial} and \S\ref{subsub-finiformalproofspecial2} that $\RD f_{\ast}\RD\pi_{X\ast}(\tau^{\geq 0}\LD\pi^{\ast}_XM)$ belongs to $\DC^{\bd}_{\coh}(Y)$.
As for the third, since 
$$
\RD f_{\ast}N\cong\RD f_{\ast}\RD\iota_{X\ast}(\tau^{\geq 0}\LD\iota^{\ast}N)\cong\RD\iota_{k\ast}\RD f_{k\ast}(\tau^{\geq 0}\LD\iota^{\ast}N),
$$
and since $\tau^{\geq 0}\LD\iota^{\ast}N$ belongs to $\DC^{\bd}_{\coh}(X_k)$, we deduce the desired assertion from \ref{thm-fini} and \ref{prop-calculusformalderived2}.

Thus the proof of \ref{thm-finiformal} is completed. \hfill$\square$

\subsection{Stein factorization}\label{sub-steinfactadqformalsch}
\index{Stein factorization|(}
\subsubsection{Announcement of the theorem}\label{subsub-steinfactadqformalsch}
In this subsection we show:
\begin{thm}[Stein factorization]\label{thm-steinfactadqformalsch}
Let $A$ be a t.u.\ adhesive\index{t.u.a. ring@t.u.\ adhesive ring} $($resp.\ t.u.\ rigid-Noetherian$)$\index{t.u. rigid-Noetherian ring@t.u.\ rigid-Noetherian ring} ring with an invertible ideal of definition of the form $I=(a)$, and $f\colon X\rightarrow Y=\Spf A$ a morphism of finite presentation.
Let $f_k\colon X_k=(X,\O_X/I^{k+1}\O_X)\rightarrow Y_k=\Spec A/I^{k+1}$ be the induced morphism of schemes for $k\geq 0$.
Suppose that $Y$ is universally cohesive, that $f_0$ is pseudo-affine\index{pseudo-affine}\index{affine!pseudo-affine@pseudo-{---}} $(${\rm \ref{dfn-pseudoaffine} (2)}$)$, and that $\O_X$ is $I$-torsion free.
Then the ring $B=\Gamma(X,\O_X)$ is a t.u.\ adhesive $($resp.\ t.u.\ rigid-Noetherian$)$ ring with an ideal of definition $IB$, and in the factorization 
$$
\xymatrix{X\ar[r]_{\pi}\ar@/^1pc/[rr]^f&Z\ar[r]_g&Y,}
$$
where $Z=\Spf B$, the map $\pi$ is proper.
\end{thm}

Notice that the assumption `$Y$ is universally cohesive' is automatic in the t.u.\ adhesive case due to {\bf \ref{ch-pre}}.\ref{thm-pf02a} (2).

The rest of this subsection is devoted to the proof of the theorem.
To this end, we need the following proposition, from which the first assertion of \ref{thm-steinfactadqformalsch} follows:
\begin{prop}\label{prop-steinfactadqformalsch11}
Let $f\colon X\rightarrow Y$ be a morphism of finite presentation between coherent universally rigid-Noetherian formal schemes\index{formal scheme!universally rigid-Noetherian formal scheme@universally rigid-Noetherian ---}, and $\mathscr{I}$ an invertible ideal of definition of $Y$ of finite type.
Let $f_k\colon X_k=(X,\O_X/\mathscr{I}^k\O_X)\rightarrow Y_k=(Y,\O_Y/\mathscr{I}^k)$ be the induced morphism of schemes for each $k\geq 0$.
Suppose that $Y$ is universally cohesive and that $f_0$ is pseudo-affine.
Then $f_{\ast}\O_X$ is an adically quasi-coherent $\O_Y$-algebra of finite type\index{adically quasi-coherent (a.q.c.) algebra!adically quasi-coherent (a.q.c.) algebra of finite type@--- of finite type}.
\end{prop}

To show the proposition, since the question is local with respect to $Y$, one can assume that $Y$ is affine $Y=\Spf A$ where $A$ is a t.u.\ rigid-Noetherian ring with an invertible principal ideal of definition $I=(a)$.
Let $I$ be the ideal of definition of $Y$ such that $\mathscr{I}=I^{\Delta}$.
\begin{prop}\label{prop-steinfactadqformalsch1}
In the situation as in {\rm \ref{thm-steinfactadqformalsch}}, let $\mathscr{F}$ be a coherent sheaf on $X$, and set $\mathscr{F}_k=\mathscr{F}/I^{k+1}\mathscr{F}$ for $k\geq 0$.
Then for any $q\geq 1$ we have$:$
\begin{itemize}
\item[{\rm (a)}] the filtration $\{F^n\H^q(X,\mathscr{F})\}_{n\geq 0}$ on the cohomology $\H^q(X,\mathscr{F})$ defined by 
$$
F^n\H^q(X,\mathscr{F})=\image(\H^q(X,I^n\mathscr{F})\rightarrow\H^q(X,\mathscr{F}))
$$
for $n\geq 0$ gives the $I$-adic topology$;$
\item[{\rm (b)}] the canonical map 
$$
\H^q(X,\mathscr{F})\longrightarrow\varprojlim_{k\geq 0}\H^q(X,\mathscr{F}_k)
$$
is an isomorphism$;$
\item[{\rm (c)}] $\H^q(X,\mathscr{F})$ is coherent as an $A$-module.
\end{itemize}
\end{prop}

\begin{proof}
The proof is done by descending induction with respect to $q\geq 1$ (cf.\ \ref{cor-formalacyclicity4}).
We henceforth assume that the conclusions (a), (b), and (c) are true for any $\mathscr{F}$ with $q$ replaced by $p$ with $p\geq q+1$.

First notice that, since $f_0$ is psuedo-affine, $f_k$ is pseudo-affine for any $k\geq 0$ (\ref{prop-cohomologicalcriterionpsuedoaff}), and that each $Y_k=\Spec A/I^{k+1}$ is universally cohesive.
Notice also that, for any $I$-torsion coherent sheaf $\mathscr{G}$ on $X$, there exists $m\geq 0$ such that $I^m\H^q(X,\mathscr{G})=0$ for all $q\geq 0$, since they are $A/I^m$-modules for some $m$.

To proceed, we fix a finite open covering $\mathscr{U}=\{U_{\alpha}\}_{\alpha\in L}$ of $X$ by affine subsets, and consider the \v{C}ech complex $K^{\bullet}=C^{\bullet}(\mathscr{U},\mathscr{F})$; notice that, since $f_0$ is pseudo-affine, $X$ is separated (cf.\ \ref{prop-sepmorformal1}).
For any $q\geq 0$ we have $\H^q(K^{\bullet})=\H^q(X,\mathscr{F})$.
Notice also that for each $p\geq 0$ the $A$-module $K^p$ is $a$-adically complete and we have $a^nK^p_{\ator}=0$ for a sufficiently large $n$.
We consider the exact sequence 
$$
0\longrightarrow a^{k+1}K^{\bullet}\longrightarrow K^{\bullet}\longrightarrow K^{\bullet}_k\longrightarrow 0
$$
of complexes for each $k\geq 0$ (cf.\ {\bf \ref{ch-pre}}, \S\ref{subsub-paircomplexgen}) and the induced filtration given by $\{F^n\H^q(K^{\bullet})\}_{n\geq 0}$ on each cohomology of $K^{\bullet}$.
Clearly, the this filtration coincides with the one in (a).

Consider the following portion of the cohomology exact sequence
$$
\cdots\longrightarrow\H^q(K^{\bullet})\longrightarrow\H^q(K^{\bullet}_k)\longrightarrow\H^{q+1}(a^{k+1}K^{\bullet})\longrightarrow\H^{q+1}(K^{\bullet})\longrightarrow\cdots.
$$
By the equality $\H^q(K^{\bullet}_k)=\H^q(X,\mathscr{F}_k)$ where $\mathscr{F}_k$ is an $a$-torsion coherent sheaf, and by the induction hypothesis (applied to the cohomologies $\H^{q+1}(a^{k+1}K^{\bullet})=\H^{q+1}(X,I^{k+1}\mathscr{F})$ and $\H^{q+1}(K^{\bullet})=\H^{q+1}(X,\mathscr{F})$), we know that the last three of the above exact sequence are coherent $A$-modules.
Hence the image of $\H^q(K^{\bullet})\rightarrow\H^q(K^{\bullet}_k)$ is finitely generated.
By {\bf \ref{ch-pre}}.\ref{lem-ullrichlemma2} (2) we deduce that the cohomology $\H^q(K^{\bullet})=\H^q(X,\mathscr{F})$ is finitely generated as an $A$-module and, moreover, the induced filtration defines the $I$-adic topology.
Thus we obtain (a).

Next we show (b).
The above argument applied to each coherent sheaf $I^{k+1}\mathscr{F}$ shows that the hypotheses of {\bf \ref{ch-pre}}.\ref{lem-steinfactadqformalsch1} (2) (where $q$ is replaced by $q-1$) are satisfied, and hence we deduce that the projective system $\{\H^{q-1}(X,\mathscr{F}_k)\}_{k\geq 0}$ satisfies {\bf (ML)}\index{Mittag-Leffler condition}.
Hence by {\bf \ref{ch-pre}}.\ref{cor-ML5} the map in question is an isomorphism.

Finally, let us show (c).
By (a) and (b), which we have already proved, $\H^q(X,\mathscr{F})$ is $I$-adically complete and finitely generated (as we have seen above).
Since $\H^q(X,\mathscr{F}_k)$ for each $k\geq 0$ is coherent, it follows from (a) that the $A$-module $\H^q(X,\mathscr{F})/I^{k+1}\H^q(X,\mathscr{F})$ is coherent for sufficiently large $k$, and hence we deduce that $\H^q(X,\mathscr{F})$ is coherent (cf.\ \ref{prop-adicqcoh55}).
\end{proof}

\subsubsection{Proof of Proposition \ref{prop-steinfactadqformalsch11}}\label{subsub-steinfactadqformalschpf21}
We may assume that $Y$ is affine $Y=\Spf A$ where $A$ is a t.u.\ rigid-Noetherian ring and that $\mathscr{I}$ comes from an invertible ideal of definition $I=(a)$ of $A$ generated by a non-zero-divisor $a\in A$.
Moreover, we may assume that $\O_X$ is $I$-torsion free, since $I$-torsion part is a coherent ideal of $\O_X$.
We look at the exact sequence
$$
0\longrightarrow\Gamma(X,I^{k+1}\O_X)\longrightarrow\Gamma(X,\O_X)\longrightarrow\Gamma(X,\O_{X_k})\longrightarrow\H^1(X,I^{k+1}\O_X),
$$
and set
$$
B_k=\Gamma(X,\O_X)/\Gamma(X,I^{k+1}\O_X)=\Gamma(X,\O_X)/I^{k+1}\Gamma(X,\O_X)\subseteq\Gamma(X,\O_{X_k})
$$
for $k\geq 0$.
(Notice that the equality $\Gamma(X,I^{k+1}\O_X)=I^{k+1}\Gamma(X,\O_X)$ follows from our assumption that $I=(a)$ is invertible.)
Since $H^1(X,I^{k+1}\O_X)$ is a coherent $A$-module (\ref{prop-steinfactadqformalsch1}), there exists $m\geq 0$ such that $I$-torsions of $H^1(X,I^{k+1}\O_X)$ are killed by $I^m$.
Then by an easy diagram chasing we deduce that $B_k$ coincides with the image of the map $\Gamma(X,\O_{X_{k+m}})\rightarrow\Gamma(X,\O_{X_k})$.
In particular, $B_k$ is an $A_k$-algebra of finite type (since $f_{k+m}$ is pseudo-affine).

Since $B=\Gamma(X,\O_X)=\varprojlim_{k\geq 0}\Gamma(X,\O_{X_k})$ and the $I$-adic completion $\widehat{B}$ is given by $\varprojlim_{k\geq 0}B_k$, we deduce that there exists an injective morphism $\widehat{B}\hookrightarrow B$ (cf.\ {\bf \ref{ch-pre}}.\ref{prop-projlimleftexact}); since this morphism factorizes the identity map $\id_B$, it is also surjective.
Hence we have $B\cong\widehat{B}$.
Since $B_0$ is finitely generated as an $A_0$-module, one has a map $A\dl X_1,\ldots,X_d\dr\rightarrow B$ that induces the surjective map $A_0[X_1,\ldots,X_d]\rightarrow B_0$.
By \cite[Chap.\ III, \S2.11, Prop.\ 14]{Bourb1} the ring $B$ is topologically finitely generated over $A$.

To show that $f_{\ast}\O_X$ is adically quasi-coherent of finite type, it remains to show that the canonical map 
$$
\Gamma(X,\O_X)\otimes_AA_{\{g\}}\rightarrow\Gamma(f^{-1}(V),\O_X),
$$
where $V=\Spf A_{\{g\}}$ with $g\in A$, is an isomorphism.
This is shown by an argument similar to that in the proof of (b) part in the proof of \ref{prop-finiformalspecial}, which uses {\bf \ref{ch-pre}}.\ref{lem-ullrichlemma3}.
Thus the proposition is proved. \hfill$\square$

\subsubsection{Proof of Theorem $\ref{thm-steinfactadqformalsch}$}\label{subsub-steinfactadqformalschpf2}
It remains to show that the map $\pi\colon X\rightarrow Z=\Spf B$ is proper.
Let $X_k\stackrel{\pi_k}{\rightarrow}Z_k=\Spec A'_k\stackrel{g_k}{\rightarrow}Y_k=\Spec A_k$ be the Stein factorization of $f_k$ for each $k\geq 0$.
Consider the exact sequence
$$
0\longrightarrow I^{k+1}\Gamma(X,\O_X)\longrightarrow\Gamma(X,\O_X)\longrightarrow\Gamma(X,\O_{X_k})\longrightarrow\H^1(X,I^{k+1}\O_X)
$$
as before and the image of $\Gamma(X,\O_X)/I^{k+1}\Gamma(X,\O_X)\hookrightarrow\Gamma(X,\O_{X_k})$.
Since the $I$-torsion part of the coherent $A$-module $\H^1(X,I^{k+1}\O_X)$ is killed by $I^m$ for some $m\geq 0$, the image of $\Gamma(X,\O_X)\rightarrow\Gamma(X,\O_{X_k})$ for sufficiently large $k$ coincides with the image of $\Gamma(X,\O_{X_{k+m}})\rightarrow\Gamma(X,\O_{X_k})$.
Hence we have the commutative diagram
$$
\xymatrix{X_k\ar[d]\ar[r]^{\pi_k}&Z_k\ar[d]\ar[r]&\Spec B_k\ar[dl]^{\gamma}\\ X_{k+m}\ar[r]_{\pi_{k+m}}&Z_{k+m}}
$$
with $B_k=\Gamma(X,\O_X)/I^{k+1}\Gamma(X,\O_X)$, where the left vertical map is a closed immersion.
Since $\gamma$ is a closed immersion, and since the map and $\pi_{k+m}$ is proper, we deduce that $X_k\rightarrow\Spec B_k$ is proper.
Then by \ref{prop-propermorformal0} the morphism $\pi\colon X\rightarrow Z=\Spf B$ is proper, and thus the proof of \ref{thm-steinfactadqformalsch} is finished. \hfill$\square$
\index{Stein factorization|)}



\setcounter{section}{0}
\renewcommand{\thesection}{\Alph{section}}
\renewcommand{\theexer}{{\bf \thechapter}.\Alph{section}.\arabic{exer}}
\section{Appendix: Stein factorization for schemes}\label{sec-steinfactorization}
\index{Stein factorization|(}
\subsection{Pseudo-affine morphisms of schemes}\label{sub-pseudoaffinemorphism}
\subsubsection{Definition and the first properties}\label{subsub-pseudoaffinemorphism}
Let $f\colon X\rightarrow Z$ be a locally of finite type morphism of schemes.
The subset $B_f$ of $Z$ defined by
$$
B_f=\{z\in Z\,|\,\dim_zX\times_Zk(z)\geq 1\}
$$
is called the {\em center}\index{center} of the morphism $f$.
By \cite[$\mathbf{IV}$, (13.1.5)]{EGA} the center $B_f$ of $f$ is a closed subset of $Z$ if $f$ is a closed map.

\begin{dfn}\label{dfn-pseudoaffine}{\rm 
Let $Y$ be a coherent scheme, and $f\colon X\rightarrow Y$ a morphism of finite presentation between schemes.

(1) A {\em pre-Stein factorization}\index{Stein factorization!pre-Stein factorization@pre-{---}} is a commutative diagram
$$
\xymatrix{X\ar[r]_{\pi}\ar@/^1pc/[rr]^f&Z\ar[r]_g&Y,}
$$
such that the following conditions are satisfied:
\begin{itemize}
\item[{\rm (a)}] $\pi$ is proper and $g$ is affine of finite presentation (and hence $\pi$ is of finite presentation);
\item[{\rm (b)}] the center $B_{\pi}$ of $\pi$ is finite over $Y$.
\end{itemize}

(2) If $f$ has a pre-Stein factorization, $f$ is said to be {\em pseudo-affine}\index{pseudo-affine}\index{affine!pseudo-affine@pseudo-{---}}.}
\end{dfn}

We say that two pre-Stein factorizations $f=g\circ\pi=g'\circ\pi'$, where $X\stackrel{\pi}{\rightarrow}Z\stackrel{g}{\rightarrow}Y$ and $X\stackrel{\pi'}{\rightarrow}Z'\stackrel{g'}{\rightarrow}Y$, are isomorphic if there exists an isomorphism $h\colon Z\stackrel{\sim}{\rightarrow}Z'$ of schemes such that the following diagram commutes:
$$
\xymatrix@R-4ex@C-1ex{&Z\ar[dd]_h^{\cong}\ar[dr]^g\\ X\ar[ur]^{\pi}\ar[dr]_{\pi'}&& Y\\ &Z'\rlap{.}\ar[ur]_{g'}}
$$

\begin{dfn}\label{dfn-steinfactorization}{\rm 
A pre-Stein factorization $f=g\circ\pi$ of $f\colon X\rightarrow Y$ as in \ref{dfn-pseudoaffine} (1) is called a {\em Stein factorization} if the following condition is satisfied:
\begin{itemize}
\item[{\rm (c)}] the canonical morphism $\O_Z\rightarrow\pi_{\ast}\O_X$ is an isomorphism.
\end{itemize}}
\end{dfn}

The Stein factorization of $f$ is, if it exists, unique up to isomorphism.
\begin{prop}\label{prop-steinfactunivcoh}
Let $Y$ be a coherent scheme, and $f\colon X\rightarrow Y$ a morphism of finite presentation.
Suppose that $Y$ is universally cohesive\index{scheme!universally cohesive scheme@universally cohesive ---}\index{cohesive!universally cohesive@universally --- (schemes)} $({\bf \ref{ch-pre}}.\ref{dfn-universallycohesive})$ $($e.g.\ Noetherian$)$ and that $f$ is pseudo-affine\index{pseudo-affine}\index{affine!pseudo-affine@pseudo-{---}}.
Then $f$ has a Stein factorization.
\end{prop}

\begin{proof}
Let $X\stackrel{\pi}{\rightarrow}Z\stackrel{g}{\rightarrow}Y$ be a pre-Stein factorization of $f$.
Since $g$ is of finite presentation, $Z$ is universally cohesive.
By \ref{thm-fini} $\pi_{\ast}\O_X$ is a coherent $\O_Z$-module, and hence $Z'=\Spec\pi_{\ast}\O_X\rightarrow Z$ is affine of finite presentation.
Then the factorization $X\rightarrow Z'\rightarrow Y$ gives a Stein factorization of $f$.
\end{proof}

\begin{prop}\label{prop-pseudoaffineproperties}
{\rm (1)} Proper $($resp.\ affine$)$ morphism of finite presentation to a coherent scheme are pseudo-affine.

{\rm (2)} Let $f\colon X\rightarrow Y$ be pseudo-affine, and $Y'\rightarrow Y$ a morphism of coherent schemes.
Then the induced map $f_{Y'}\colon X\times_YY'\rightarrow Y'$ is pseudo-affine. \hfill$\square$
\end{prop}

\begin{prop}\label{prop-steinfactunivcoh2}
Let $Y$ be a coherent and universally cohesive scheme, and $f\colon X\rightarrow Y$ a morphism of finite presentation.
Suppose $f$ is pseudo-affine.
Then for any coherent sheaf $\mathscr{F}$ on $X$, $\RD^qf_{\ast}\mathscr{F}$ for $q\geq 1$ is a coherent sheaf on $Y$, and $f_{\ast}\mathscr{F}$ is a coherent module over a finitely presented $\O_Y$-algebra.
\end{prop}

\begin{proof}
Let $X\stackrel{\pi}{\rightarrow}Z\stackrel{g}{\rightarrow}Y$ be the Stein factorization.
Since $g$ is affine, we have $\RD^qf_{\ast}\mathscr{F}=g_{\ast}\RD^q\pi_{\ast}\mathscr{F}$ for $q\geq 0$.
Since $\pi$ is proper of finite presentation, $\RD^q\pi_{\ast}\mathscr{F}$ is a coherent $\O_Z$-module for $q\geq 0$ (\ref{thm-fini}).
If $q>0$ then the support of the sheaf $\RD^q\pi_{\ast}\mathscr{F}$ is contained in the center $B_{\pi}$, which is finite over $Y$, and hence $\RD^qf_{\ast}\mathscr{F}=g_{\ast}\RD^q\pi_{\ast}\mathscr{F}$ is a coherent $\O_Y$-module.
The sheaf $f_{\ast}\mathscr{F}=g_{\ast}\pi_{\ast}\mathscr{F}$ is coherent over $g_{\ast}\O_Z$, whence the result.
\end{proof}

\subsubsection{Pseudo-affineness and compactifications}
In this paragraph we need the notion of {\em $U$-admissible blow-ups}\index{admissible!U-admissible blow-up@$U$-{---} blow-up}\index{blow-up!U-admissible blow-up@$U$-admisible ---} (already appeared in \S\ref{subsub-carvinglemma}):
\begin{dfn}\label{dfn-Uadmissiblebupsadhoc}{\rm 
Let $X$ be a coherent scheme, and $U\subseteq X$ a quasi-compact open subscheme of $X$.
A {\em $U$-admissible blow-up} is a blow-up $Y\rightarrow X$ of $X$ along a quasi-coherent ideal $\mathscr{J}$ of $\O_X$ of finite type such that the closed subscheme $V(\mathscr{J})$ of $X$ corresponding to $\mathscr{J}$ is disjoint from $U$.}
\end{dfn}

Later in {\bf \ref{ch-rigid}}, \S\ref{subsub-birationalgeomblowups} we will discuss $U$-admissible blow-ups in detail.
\begin{prop}\label{prop-pseudoaffinecomp}
Let $Y$ be a coherent scheme, and $f\colon X\rightarrow Y$ a morphism of finite presentation.
Then the following conditions are equivalent$:$
\begin{itemize}
\item[{\rm (a)}] the morphism $f$ is pseudo-affine\index{pseudo-affine}\index{affine!pseudo-affine@pseudo-{---}}$;$
\item[{\rm (b)}] there exist an open immersion $X\hookrightarrow\ovl{X}$ over $Y$ into a proper $Y$-scheme $\ovl{X}$ and an effective Cartier divisor $D$ of $\ovl{X}$ such that the support of $D$ is $\ovl{X}\setminus X$, that $D$ is semiample over $Y$, and that the normal sheaf $\O_{\ovl{X}}(D)|_D$ on $D$ is ample over $Y;$
\item[{\rm (c)}] for any open immersion $X\rightarrow X'$ over $Y$ into a proper $Y$-scheme $X'$, there exist an $X$-admissible blow-up $X''\rightarrow X'$, and an effective Cartier divisor $D$ on the closure $\til{X}$ of $X$ in $X''$ such that the support of $D$ is $\til{X}\setminus X$, that $D$ is semiample over $Y$, and that the normal sheaf $\O_{\til{X}}(D)|_D$ on $D$ is ample over $Y$.
\end{itemize}
\end{prop}

The proof uses a generalized version of Nagata's embedding theorem\index{Nagata's embedding theorem}\index{Nagata, M.} ({\bf \ref{ch-rigid}}.\ref{thm-nagataembedding}).
\begin{proof}
First, let us show (a) $\Rightarrow$ (b).
Let 
$$
\xymatrix{X\ar[r]_{\pi}\ar@/^1pc/[rr]^f&Z\ar[r]_g&Y}
$$
be a pre-Stein factorization of $f$.
Take a projective compactification $\ovl{Z}$ of $Z$, which is projective of finite presentation over $Y$, and let $\Delta$ be an effective ample Cartier divisor whose support is $\ovl{Z}\setminus Z$.
One can extend the morphism $f$ to a proper morphism $\ovl{X}\rightarrow\ovl{Z}$ in such a way that $\ovl{X}\times_{\ovl{Z}}Z\cong X$; indeed, by Nagata's embedding theorem one has a proper map $\ovl{X}\rightarrow\ovl{Z}$ such that $\ovl{X}$ contains $X$ as a Zariski open subset; then replace $\ovl{X}$ by the scheme-theoretic closure of $X$ in $\ovl{X}$.
We thus arrive at the following commutative diagram
$$
\xymatrix@R-2ex@C-3ex{X\ \ar@{^{(}->}[rr]\ar[d]_{\pi}&&\ovl{X}\ar[d]^{\ovl{\pi}}\\ Z\ \ar@{^{(}->}[rr]\ar[dr]_g&&\ovl{Z}\ar[dl]^{\ovl{g}}\\ &Y,}
$$
where the horizontal arrows are open immersions and the square is Cartesian.
Now, consider the pull-back $D$ of the divisor $\Delta$ to $\ovl{X}$.
Clearly, $D$ is supported on $\ovl{X}\setminus X$ and is $Y$-semiample (since $\O_{\ovl{Z}}(\Delta)$ is generated by global sections).
Since the map $\ovl{\pi}$ is an isomorphism around $\ovl{X}\setminus X$, $\O(D)|_D$ on $D$ is $Y$-ample.

Next we show (b) $\Rightarrow$ (c).
Let $X\hookrightarrow X'$ be an open immersion over $Y$ into a scheme $X'$ proper over $Y$, and take $X\hookrightarrow\ovl{X}$ as in (b).
Then there exists a diagram of proper $Y$-schemes
$$
\xymatrix@-3ex{&X''\ar[dl]\ar[dr]\\ X'&&\ovl{X}}
$$
consisting of $X$-admissible blow-ups (cf.\ {\bf \ref{ch-rigid}}.\ref{prop-correspondencediagramfurther} (2)); let $\til{X}$ be the closure of $X$ in $X''$.
Let $\mathscr{I}$ be the blow-up center of $X''\rightarrow\ovl{X}$.
Then $\mathscr{I}\O_{\til{X}}$ is invertible and $\ovl{X}$-ample.
Moreover, the support of the corresponding divisor $E$ lies in $\til{X}\setminus X$.
Let $D$ be the effective Cartier Divisor of $\ovl{X}$ as in (b), and take the pull-back $\til{D}$ on $\til{X}$.
Then by \cite[$\mathbf{II}$, (4.6.13) (ii)]{EGA} there exists an integer $n>0$ such that the divisor $E+n\til{D}$ satisfies the conditions as in (c).

The converse implication (c) $\Rightarrow$ (b) follows from Nagata's embedding theorem.

It remains to show (b) $\Rightarrow$ (a).
Consider the situation as in (b).
Since $D$ is $Y$-semiample, replacing $D$ by a multiple of itself if necessary, we may assume that $D$ induces a morphism $\ovl{\pi}\colon\ovl{X}\rightarrow P$ over $Y$, where $P$ is projective of finite presentation over $Y$, such that there exists a $Y$-ample divisor $\Delta$ of which the pull-back to $\ovl{X}$ coincides with $D$.
Let $Z$ be the open complement of $\Delta$ in $P$, which is affine and of finite presentation over $Y$.
Since the support of $D$ is $\ovl{X}\setminus X$, we have $\ovl{X}\times_PZ=X$, and thus we get the proper morphism $\pi\colon X\rightarrow Z$.
The center of $\ovl{\pi}$ is a closed subset of $P$.
Since $\O_{\ovl{X}}(D)|_D$ is $Y$-ample, $\ovl{\pi}$ is finite over the boundary $P\setminus Z$.
Hence the center $B_{\ovl{\pi}}$ lies in $Z$ and is finite over $Y$, since it is affine and proper over $Y$ (cf.\ {\bf \ref{ch-rigid}}.\ref{prop-appclassicalZR}).
Thus we have the pre-Stein factorization $X\stackrel{\pi}{\rightarrow}Z\stackrel{g}{\rightarrow}Y$ of $f$ and, therefore, the proposition is proved.
\end{proof}
\index{Stein factorization|)}

\subsection{Cohomological criterion}\label{sub-cohomologicalcriterionpsuedoaff}
For the proof of the following theorem, see \cite[4.6]{Lutk1} (cf.\ \cite{GH}):
\begin{thm}[Cohomological criterion]\label{thm-cohomologicalcriterionpsuedoaff}
Let $Y$ be a Noetherian scheme, and $f\colon X\rightarrow Y$ a separated morphism of finite type.
Then the following conditions are equivalent$:$
\begin{itemize}
\item[{\rm (a)}] the morphism $f$ is pseudo-affine$;$
\item[{\rm (b)}] $\RD^1f_{\ast}\mathscr{F}$ is coherent on $Y$ for any coherent sheaf $\mathscr{F}$ on $X;$
\item[{\rm (c)}] $\RD^qf_{\ast}\mathscr{F}$ are coherent for $q\geq 1$ for any coherent sheaf $\mathscr{F}$ on $X$. \hfill$\square$
\end{itemize}
\end{thm}

\begin{prop}\label{prop-cohomologicalcriterionpsuedoaff}
Let $Y=\Spec A$ be an affine scheme, and $Y_0=\Spec A_0$ a closed subscheme of $Y$ defined by a finitely generated nilpotent ideal $I$ of $A$.
Let $f\colon X\rightarrow Y$ be a morphism of finite presentation, and $f_0\colon X_0=X\times_YY_0\rightarrow Y_0$ the induced map.
Then $f$ is pseudo-affine if and only if $f_0$ is pseudo-affine.
\end{prop}

\begin{proof}
The `only if' part follows from \ref{prop-pseudoaffineproperties} (2).
We are to show the `if' part.
Suppose $f_0$ is pseudo-affine, and let $X_0\stackrel{\pi_0}{\rightarrow}Z_0\stackrel{g_0}{\rightarrow}Y_0$ be a pre-Stein factorization of $f_0$.
Let us take a filtered inductive system $\{A_{\alpha},I_{\alpha}\}_{\alpha\in L}$ indexed by a directed set consisting of Noetherian rings and nilpotent ideals such that $A=\varinjlim_{\alpha\in L}A_{\alpha}$ and $I=\varinjlim_{\alpha\in L}I_{\alpha}$.
Set $Y_{\alpha}=\Spec A_{\alpha}$ and $Y_{0\alpha}=\Spec A_{0\alpha}$ where $A_{0\alpha}=A_{\alpha}/I_{\alpha}$ for each $\alpha\in L$.
Replacing $L$ by a cofinal subset if necessary, we may assume that for each $\alpha\in L$ there exists a $Y_{\alpha}$-scheme $f_{\alpha}\colon X_{\alpha}\rightarrow Y_{\alpha}$ of finite type such that $X\cong X_{\alpha}\times_{Y_{\alpha}}Y$.
Let $f_{0\alpha}\colon X_{0\alpha}=X_{\alpha}\times_{Y_{\alpha}}Y_{0\alpha}\rightarrow Y_{0\alpha}$ be the induced map for each $\alpha\in L$.
Again replacing $L$ by a cofinal subset if necessary, we may assume that for each $\alpha\in L$ there exists an affine $Y_{0\alpha}$-scheme $g_{0\alpha}\colon Z_{0\alpha}\rightarrow Y_{0\alpha}$ of finite type such that $Z_0\cong Z_{0\alpha}\times_{Y_{0\alpha}}Y_0$.
We may also assume that each $f_{0\alpha}$ factors through $Z_{0\alpha}$ by a morphism $\pi_{0\alpha}\colon X_{0\alpha}\rightarrow Z_{0\alpha}$.
By \cite[$\mathbf{IV}$, (8.10.5)]{EGA} we may assume that the morphisms $\pi_{0\alpha}$ are proper and that the center $B_{\pi_{0\alpha}}$ is finite over $Y_{0\alpha}$ (since $B_{\pi_0}=B_{\pi_{0\alpha}}\times_{Y_{0\alpha}}Y_0$).
Hence, in this way, we arrive at the situation where the map $f$ is realized as the filtered projective limit of the maps $f_{\alpha}$ between Noetherian schemes such that all $f_{0\alpha}$ are pseudo-affine.

By \ref{thm-cohomologicalcriterionpsuedoaff} we know that for any coherent sheaf $\mathscr{F}$ on $X_{0\alpha}$ the sheaf $\RD^qf_{0\alpha\ast}\mathscr{F}$ for $q\geq 1$ is coherent on $Y_{0\alpha}$.
To show the analogous statement for each $f_{\alpha}$, let $\mathscr{F}$ be a coherent sheaf on $X_{\alpha}$, and set $\mathscr{F}_0=\mathscr{F}/I_{\alpha}\mathscr{F}$; by induction with respect to $d\geq 0$ such that $I^d_{\alpha}=0$ we may assume $I^2_{\alpha}=0$.
Consider the cohomology exact sequence induced from 
$$
0\longrightarrow I_{\alpha}\mathscr{F}\longrightarrow\mathscr{F}\longrightarrow\mathscr{F}_0\longrightarrow 0,
$$
where the first and the third sheaves are coherent sheaves on $Y_{0\alpha}$.
By \ref{prop-cohqcoh2} and in view of the fact that quasi-coherent sheaves of finite type are coherent (\cite[$\mathbf{I}$, (1.5.1)]{EGA}), one sees that the cohomologies $\RD^qf_{\alpha\ast}\mathscr{F}$ are coherent on $Y_{\alpha}$ for $q\geq 1$.
Hence again by \ref{thm-cohomologicalcriterionpsuedoaff} we deduce that $f_{\alpha}$ is pseudo-affine.
By \ref{prop-pseudoaffineproperties} (2) it follows that $f$ is pseudo-affine, as desired.
\end{proof}

\section{Appendix: Zariskian schemes}\label{sec-zariskianschemes}
\index{Zariskian!Zariskian scheme@--- scheme|(}
\index{scheme!Zariskian scheme@Zariskian ---|(}
\subsection{Zariskian schemes}\label{sub-zariskianschemes}
\subsubsection{Zariskian rings and Zariskian spectra}\label{subsub-zariskianrings}
The notion of Zariskian schemes has already appeared in {\bf \ref{ch-pre}}, \S\ref{sub-schpair}.
In this subsection we give a more systematic and precise account on Zariskian schemes.

Let $A$ be a commutative ring endowed with the $I$-adic topology\index{topology!adic topology@adic ---}\index{adic!adic topology@--- topology} ({\bf \ref{ch-pre}}, \S\ref{subsub-adicfiltrationtopology}) by an ideal $I\subset A$.
The topological ring $A$ is called a {\em Zariskian ring}\index{Zariskian!Zariskian ring@--- ring} if $1+I\subseteq A^{\times}$ (cf.\ {\bf \ref{ch-pre}}, \S\ref{subsub-zariskianpairs}).
The definition does not depend on the choice of the ideal of definition $I$; indeed, if $J$ is another ideal of definition, that is, $I^m\subseteq J^n\subseteq I$ for some $m,n>0$, then since $1+J^n\subseteq A^{\times}$, $J^n$ is contained in the Jacobson radical of $A$, and hence $J$ is contained in the Jacobson radical, that is to say, $1+J\subseteq A^{\times}$ (cf.\ {\bf \ref{ch-pre}}.\ref{prop-zariskipair1}).
For any ring $A$ with $I$-adic topology the {\em associated Zariskian ring}\index{Zariskian!associated Zariskian@associated ---}, denoted by $\zat{A}$, is defined to be $S^{-1}A$ where $S$ is the multiplicative set given by $S=1+I$ (cf.\ {\bf \ref{ch-pre}}, \S\ref{subsub-zariskianpairs}).

\begin{dfn}\label{dfn-zariskianspectrum}{\rm
Let $A$ be a Zariskian ring, and $I$ an ideal of definition.
Then the {\em Zariskian spectrum}\index{Zariskian!Zariskian spectrum@--- spectrum} of $A$, denoted by $X=\Spz A$, is the topologically locally ringed space defined as follows:
\begin{itemize}
\item it is, as a set, the subset $V(I)$ of $\Spec A$;
\item the topology is the induced topology as a subset of $\Spec A$;
\item the structure sheaf is the one given by $i^{-1}\O_{\Spec A}$ considered with the $I$-adic topology, where $i\colon V(I)\hookrightarrow\Spec A$ is the inclusion map.
\end{itemize}}
\end{dfn}

The topology on $X=\Spz A$ has the open basis $\{\mathfrak{D}(f)\}_{f\in A}$ consisting of quasi-compact open subsets, where we set $\mathfrak{D}(f)=D(f)\cap X$ for $f\in A$.
The open subset $\mathfrak{D}(f)$, considered as a topologically locally ringed space with the induced structure sheaf, is isomorphic to the Zariskian spectrum $\Spz\zat{A}_f$, where $\zat{A}_f$ is the associated Zariskian ring of $A_f$ with respect to the $IA_f$-adic topology.
Notice that the underlying topological space of $\Spz A$ is coherent ({\bf \ref{ch-pre}}.\ref{dfn-quasicompact1}).

The proofs of the following lemmas are entirely similar to those of \ref{lem-formalnot1} and \ref{prop-lemformalnot111}, respectively:
\begin{lem}\label{lem-zariskiannot1}
For finitely many elements $f_1,\ldots,f_r\in A$ the following conditions are equivalent$:$
\begin{itemize}
\item[{\rm (a)}] the open sets $\mathfrak{D}(f_i)$ $(i=1,\ldots,r)$ cover $X=\Spz A$, that is, $X=\bigcup^r_{i=1}\mathfrak{D}(f_i);$
\item[{\rm (b)}] for any ideal of definition $I$ of $A$ the open sets $D(\ovl{f}_i)$ $(i=1,\ldots,r)$, where $\ovl{f}_i=(f_i\ \mathrm{mod}\ I)$, cover $\Spec A/I$, that is, $\Spec A/I=\bigcup^r_{i=1}D(\ovl{f}_i);$
\item[{\rm (c)}] the open sets $D(f_i)$ $(i=1,\ldots,r)$ cover $\Spec A$, that is, $\Spec A=\bigcup^r_{i=1}D(f_i);$
\item[{\rm (d)}] the ideal generated by $f_1,\ldots,f_r$ coincides with $A$. \hfill$\square$
\end{itemize}
\end{lem}

\begin{lem}\label{lem-zariskiannot111}
An open subset $U$ of $X=\Spz A$ is quasi-compact if and only if $U$ is of the form $U=X\setminus V(J)$ in $\Spec A$ by a finitely generated ideal $J\subseteq A$. \hfill$\square$
\end{lem}

Let $A$ be a ring, and $I\subset A$ an ideal.
As above, let $i\colon V(I)\hookrightarrow\Spec A$ be the inclusion map.
For an $A$-module $M$ we define the sheaf $M^{\lozenge}$ on the topological space $V(I)$ by
$$
M^{\lozenge}=i^{-1}\til{M}.
$$

\begin{prop}\label{prop-zariskianglobalsection}
For any $f\in A$ we have
$$
\Gamma(V(I)\cap D(f),M^{\lozenge})=M\otimes_AS^{-1}A_f,
$$
where $S=1+IA_f$.
\end{prop}

\begin{proof}
Since $D(f)=\Spec A_f$ and $V(I)\cap D(f)=V(IA_f)$, we may assume without loss of generality that $A=A_f$.
First we show that, if $\mathscr{M}$ denotes the presheaf pull-back of $\til{M}$ by $i$, then $\Gamma(V(I),\mathscr{M})=M\otimes_AS^{-1}A_f$.
The left-hand side is the inductive limit of $\Gamma(U,\til{M})$, where $U$ runs over the set of all open subsets containing $V(I)$.
If an open subset $U=\Spec A\setminus V(J)$ contains $V(I)$, then $V(J)\cap V(I)=\emptyset$, that is, $J+I=A$.
This implies that there exists $a\in I$ such that $1+a\in J$.
Since $V(I)\subseteq D(1+a)\subseteq U$, we deduce that the open subsets of the form $D(1+a)$ with $a\in I$ form a cofinal subset in the set of all open subsets containing $V(I)$.
Hence the module in question is the inductive limit of the modules of the form $M\otimes_AA_{(1+a)}$ with $a\in I$, which is nothing but $M\otimes_A(1+I)^{-1}A$, as desired.

It remains to show that $\Gamma(V(I),\mathscr{M})=\Gamma(V(I),M^{\lozenge})$.
This can be shown by the standard argument as in \cite[$\mathbf{I}$, (1.3.7)]{EGA}; the detail is left to the reader.
\end{proof}

\begin{cor}\label{cor-zariskianglobalsection}
If $A$ is a Zariskian ring, then $M^{\lozenge}$ is an $\O_X$-module, where $X=\Spz A$, and we have $\Gamma(\mathfrak{D}(f),M^{\lozenge})=M\otimes_A\zat{A}_f$ for any $f\in A$. \hfill$\square$
\end{cor}

\begin{prop}\label{prop-zariskianlocalring}
Let $A$ be a Zariskian ring, and set $X=\Spz A$.
Let $M$ be an $A$-module.
Then for any point $x\in X$, that is, an open prime ideal $x=\mathfrak{p}$, we have $M^{\lozenge}_x=M\otimes_AA_{\mathfrak{p}}$.
\end{prop}

To show this, we need the following obvious lemma:
\begin{lem}\label{lem-zariskianlocalring}
Let $\{A_{\lambda}\}_{\lambda\in\Lambda}$ be a filtered inductive system of rings with adic topology such that for any $\lambda\leq\mu$ the map $A_{\lambda}\rightarrow A_{\mu}$ is adic, that is, for any ideal of definition $I$ of $A_{\lambda}$, $IA_{\mu}$ is an ideal of definition of $A_{\mu}$.
Then we have the canonical isomorphism ${\textstyle \zat{(\varinjlim_{\lambda\in\Lambda}A_{\lambda})}\stackrel{\sim}{\longrightarrow}\varinjlim_{\lambda\in\Lambda}\zat{A}_{\lambda}}$ of topological rings. \hfill$\square$
\end{lem}

\begin{proof}[Proof of Proposition {\rm \ref{prop-zariskianlocalring}}]
By \ref{lem-zariskianlocalring} we have $M^{\lozenge}_x=\varinjlim_{x\in\mathfrak{D}(f)}\Gamma(\mathfrak{D}(f),M^{\lozenge})=M\otimes_A\varinjlim_{x\in\mathfrak{D}(f)}\zat{A}_f=M\otimes_A\zat{A}_{\mathfrak{p}}$.
But due to {\bf \ref{ch-pre}}.\ref{prop-relpair32} we have $\zat{A}_{\mathfrak{p}}=A_{\mathfrak{p}}$.
\end{proof}

\subsubsection{Zariskian schemes}\label{subsub-zariskianschemes}
\begin{dfn}\label{dfn-zariskianschemes}{\rm 
(1) A topologically locally ringed space $(X,\O_X)$ is called an {\em affine Zariskian scheme}\index{scheme!Zariskian scheme@Zariskian ---!affine Zariskian scheme@affine --- ---}\index{Zariskian!Zariskian scheme@--- scheme!affine Zariskian scheme@affine --- ---} if it is isomorphic to $\Spz A$ for a Zariskian ring $A$.

(2) A topologically locally ringed space $(X,\O_X)$ is called a {\em Zariskian scheme} if it has an open covering $X=\bigcup_{\alpha\in L}U_{\alpha}$ by affine Zariskian schemes.}
\end{dfn}

Morphisms of Zariskian schemes\index{morphism of Zariskian schemes@morphism (of Zariskian schemes)} are the morphisms of topologically locally ringed spaces.
A Zariskian scheme $X$ is said to be {\em quasi-compact} (resp.\ {\em quasi-separated}, resp.\ {\em coherent}) if the underlying topological space of $X$ is quasi-compact\index{space@space (topological)!quasi-compact topological space@quasi-compact ---}\index{quasi-compact!quasi-compact topological space@--- (topological) space} (resp.\ quasi-separated\index{space@space (topological)!quasi-separated topological space@quasi-separated ---}\index{quasi-separated!quasi-separated topological space@--- (topological) space}, resp.\ coherent\index{space@space (topological)!coherent topological space@coherent ---}\index{coherent!coherent topological space@--- (topological) space}) ({\bf \ref{ch-pre}}.\ref{dfn-quasicompactness}, {\bf \ref{ch-pre}}.\ref{dfn-quasiseparatedness}, {\bf \ref{ch-pre}}.\ref{dfn-quasicompact1}).
Notice that any Zariskian scheme is locally coherent\index{coherent!coherent topological space@--- (topological) space!locally coherent topological space@locally --- ---}\index{space@space (topological)!coherent topological space@coherent ---!locally coherent topological space@locally --- ---} and sober\index{space@space (topological)!sober topological space@sober ---} ({\bf \ref{ch-pre}}.\ref{dfn-locallycoherent}, {\bf \ref{ch-pre}}, \S\ref{subsub-sober}).
We denote by $\Zs$ (resp.\ $\CZs$, resp.\ $\AZs$) the category of Zariskian (resp.\ coherent Zariskian, resp.\ affine Zariskian) schemes.
These are full subcategories of the category of topologically locally ringed spaces.
Notice that schemes are canonically regarded as Zariskian schemes considered with the $0$-adic topology.
The following theorem is straightforward in view of \ref{cor-zariskianglobalsection} and \ref{prop-zariskianlocalring}:
\begin{thm}\label{thm-affinezariskianschemes}
The functor
$$
A\longmapsto\Spz A
$$
gives rise to a categorical equivalence between the opposite category of the category of Zariskian rings and the category $\AZs$ of affine Zariskian schemes. \hfill$\square$
\end{thm}

For a Zariskian scheme $X$ and an open subset $U$ of the underlying topological space of $X$, $U$ has the induced structure as a Zariskian scheme in such a way that the inclusion map extends to a morphism $j\colon U\rightarrow X$ of Zariskian schemes.
The Zariskian schemes obtained in this way is called an {\em open Zariskian subscheme} of $X$.
Similarly, one can define the notion of {\em open immersions} of Zariskian schemes.

\begin{dfn}\label{dfn-zariskianschemesqcpt}{\rm 
A morphism $f\colon X\rightarrow Y$ of Zariskian schemes is said to be {\em quasi-compact}\index{morphism of Zariskian schemes@morphism (of Zariskian schemes)!quasi-compact morphism of Zariskian schemes@quasi-compact ---}\index{quasi-compact!quasi-compact morphism of Zariskian schemes@--- morphism (of Zariskian schemes)} if the map of underlying topological spaces is quasi-compact ({\bf \ref{ch-pre}}.\ref{dfn-quasicompactness} (2)).}
\end{dfn}

We have the so-called {\em Zariskian completion}\index{completion!Zariskian completion@Zariskian ---}: let $X$ be a scheme, and $Y\subseteq X$ a closed subscheme; then $(Y,i^{-1}\O_X)$, where $i\colon Y\rightarrow X$ is the closed immersion, is a Zariskian scheme, which we denote by 
$$
\zat{X}|_Y.
$$

\subsection{Fiber products}\label{sub-zariskianfiberprod}
\begin{thm}\label{thm-fiberproductszariskianexists}
The category of Zariskian schemes has fiber products. \hfill$\square$
\end{thm}

As usual, the construction is reduced to the affine case.
Let 
$$
\Spz B\longrightarrow\Spz A\longleftarrow\Spz C\leqno{(\ast)}
$$
be a diagram of affine Zariskian schemes, where $A$ (resp.\ $B$, resp.\ $C$) is a Zariskian ring with an ideals of definition $I$ (resp.\ $J$, resp.\ $K$).
The fiber product of the diagram $(\ast)$ is given by a diagram of {\em affine} Zariskian schemes corresponding to
$$
\xymatrix@R-3ex@C-6ex{&\zat{\ \ (B\otimes_AC)}\\ B\ar[ur]&&C\ar[ul]\\ &A.\ar[ul]\ar[ur]}
$$
Here the associated Zariskian $\zat{(B\otimes_AC)}$ is taken with respect to the ideal $H$ of $B\otimes_AC$ given by 
$$
H=\image(J\otimes_AC\rightarrow B\otimes_AC)+\image(B\otimes_AK\rightarrow B\otimes_AC).
$$

\subsection{Ideals of definition and adic morphisms}\label{sub-zariskianadic}
\begin{dfn}\label{dfn-zariskianschemesideal}{\rm 
Let $X$ be a Zariskian scheme.
An {\em ideal of definition}\index{ideal of definition} of $X$ is a quasi-coherent ideal $\mathscr{I}$ of $\O_X$ such that for an affine open covering $X=\bigcup_{\alpha\in L}U_{\alpha}$ with $U_{\alpha}\cong\Spz A_{\alpha}$ by a Zariskian ring $A_{\alpha}$ for each $\alpha\in L$, $\mathscr{I}|_{U_{\alpha}}$ is isomorphic to the ideal of the form $I^{\lozenge}_{\alpha}$ by an ideal of definition $I_{\alpha}$ of $A_{\alpha}$.}
\end{dfn}

Clearly, ideals of definition always exist locally.
If $\mathscr{I}$ is an ideal of definition of a Zariskian scheme $X$ and $U=\Spz A$ is an affine open subset of $X$ where $A$ is a Zariskian ring, then $\mathscr{I}|_U\cong I^{\lozenge}$ for an ideal of definition $I$ of $A$.

\begin{dfn}\label{dfn-zariskianschemesadic}{\rm
A morphism $f\colon X\rightarrow Y$ of Zariskian schemes is said to be {\em adic}\index{morphism of Zariskian schemes@morphism (of Zariskian schemes)!adic morphism of Zariskian schemes@adic ---} if for any open subset $V$ of $Y$ having an ideal of definition $\mathscr{I}$ the ideal $\mathscr{I}\O_X|_{f^{-1}(V)}$ is an ideal of definition of the open Zariskian subscheme $f^{-1}(V)$.}
\end{dfn}

We denote by $\Zs^{\ast}$ (resp.\ $\CZs^{\ast}$, resp.\ $\AZs^{\ast}$) the category of Zariskian (resp.\ coherent Zariskian, resp.\ affine Zariskian) schemes and adic morphisms.

\begin{prop}\label{prop-adicmor1zat}
{\rm (1)} Let $f\colon X\rightarrow Y$ and $g\colon Y\rightarrow Z$ be morphisms of Zariskian schemes.
If $f$ and $g$ are adic, then so is the composition $g\circ f$.
If $g\circ f$ and $g$ are adic, then so is $f$.

{\rm (2)} Let $S$ be a Zariskian scheme, and $f\colon X\rightarrow X'$ and $g\colon Y\rightarrow Y'$ two adic $S$-morphisms of Zariskian schemes $($not necessarily adic over $S)$.
Then $f\times_Sg\colon X\times_SY\rightarrow X'\times_SY'$ is adic.

{\rm (3)} Let $S$ be a Zariskian scheme, and $f\colon X\rightarrow Y$ an adic $S$-morphism between Zariskian schemes $($not necessarily adic$)$ over $S$. Then for any $($not necessarily adic$)$ morphism $S'\rightarrow S$ of Zariskian schemes, the induced morphism $f_{S'}\colon X\times_SS'\rightarrow Y\times_SS'$ is adic. \hfill$\square$
\end{prop}

The proof is similar to that of \ref{prop-adicmor1}.

\subsection{Morphism of finite type and morphism of finite presentation}\label{sub-zariskiantypefin}
\begin{dfn}\label{dfn-topfintypezat}{\rm 
A morphism $f\colon X\rightarrow Y$ of Zariskian schemes is said to be {\em locally of finite type}\index{morphism of Zariskian schemes@morphism (of Zariskian schemes)!morphism of Zariskian schemes locally of finite type@--- locally of finite type} (resp.\ {\em locally of finite presentation}\index{morphism of Zariskian schemes@morphism (of Zariskian schemes)!morphism of Zariskian schemes locally of finite presentation@--- locally of finite presentation}) if the following conditions are satisfied:
\begin{itemize}
\item[{\rm (a)}] the morphism $f$ is adic $({\rm \ref{dfn-zariskianschemesadic}})$; 
\item[{\rm (b)}] there exist an affine open covering $\{V_i\}$ of $Y$ with $V_i=\Spz B_i$, where $B_i$ is a Zariskian ring with an ideal of definition $J_i$, and for each $i$ an affine open covering $\{U_{ij}\}_j$ of $f^{-1}(V_i)$ with $U_{ij}=\Spz A_{ij}$, where $A_{ij}$ is a Zariskian ring with the ideal of definition $J_iA_{ij}$ due to (a), such that each $A_{ij}$ is isomorphic to the associated Zariskian of a finitely generated (resp.\ finitely presented) algebra over $B_i$.
\end{itemize}
The morphism $f$ is said to be {\em of finite type}\index{morphism of Zariskian schemes@morphism (of Zariskian schemes)!morphism of Zariskian schemes of finite type@--- of finite type} (resp.\ {\em of finite presentation}\index{morphism of Zariskian schemes@morphism (of Zariskian schemes)!morphism of Zariskian schemes of finite presentation@--- of finite presentation}) if it is locally of finite type and $f$ is quasi-compact.}
\end{dfn}

Notice that, if $B$ is the associated Zariskian $\zat{(A[X_1,\ldots,X_n]/\mathfrak{a})}$ of a finitely presented $A$-algebra (where $\mathfrak{a}$ is a finitely generated ideal), then we have 
$$
B=\zat{A[X_1,\ldots,X_n]}/\mathfrak{a}\zat{A[X_1,\ldots,X_n]},
$$
since, in general, quotient rings (by not necessarily closed ideals) of a Zariskian ring are always Zariskian with respect to the induced ideal of definition.
In particular, for any affine Zariskian scheme $\Spz A$ and any $f\in A$ the open immersion $\Spz\zat{A}_f\hookrightarrow\Spz A$ is of finite presentation.
This fact justifies our definition of `(locally) of finite presentation' as above.
By an argument similar to the proof of \ref{prop-topfintype2}, we have the following:
\begin{prop}\label{prop-topfintype2xzat}
{\rm (1)} An open immersion is locally of finite presentation.

{\rm (2)} The composition of two morphisms locally of finite type $($resp.\ of finite type, resp.\ locally of finite presentation, resp.\ of finite presentation$)$ is again locally of finite type $($resp.\ of finite type, resp.\ locally of finite presentation, resp.\ of finite presentation$)$.
If the composition $g\circ f$ of morphisms $f\colon X\rightarrow Y$ and $g\colon Y\rightarrow Z$ of formal schemes is locally of finite type and $g$ is adic, then $f$ is locally of finite type.
If $g\circ f$ is locally of finite presentation and $g$ is locally of finite type, then $f$ is locally of finite presentation.

{\rm (3)} Let $S$ be a Zariskian scheme, and $f\colon X\rightarrow X'$ and $g\colon Y\rightarrow Y'$ two adic $S$-morphisms of Zariskian schemes.
Suppose $f$ and $g$ are locally of finite type $($resp.\ of finite type, resp.\ locally of finite presentation, resp.\ of finite presentation$)$.
Then $f\times_Sg\colon X\times_SY\rightarrow X'\times_SY'$ is locally of finite type $($resp.\ of finite type, resp.\ locally of finite presentation, resp.\ of finite presentation$)$.

{\rm (4)} Let $S$ be a Zariskian scheme, and $f\colon X\rightarrow Y$ an adic $S$-morphism between Zariskian schemes. Suppose $f$ is locally of finite type $($resp.\ of finite type, resp.\ locally of finite presentation, resp.\ of finite presentation$)$.
Then for any morphism $S'\rightarrow S$ of Zariskian schemes the induced morphism $f_{S'}\colon X\times_SS'\rightarrow Y\times_SS'$ is locally of finite type $($resp.\ of finite type, resp.\ locally of finite presentation, resp.\ of finite presentation$)$. \hfill$\square$
\end{prop}
\index{scheme!Zariskian scheme@Zariskian ---|)}
\index{Zariskian!Zariskian scheme@--- scheme|)}

\section{Appendix: FP-approximated sheaves and GFGA theorems}\label{sec-weakcoherency}
\subsection{Finiteness up to bounded torsion}\label{sub-weaklycoherentmodules}
\subsubsection{Weak isomorphisms}\label{subsub-weakisomorphisms}
\index{weak!weak isomorphism@--- isomorphism|(}
Let $A$ be a ring, and $I\subseteq A$ a finitely generated ideal.
\begin{dfn}\label{dfn-weakisomorphisms}{\rm 
A homomorphism $f\colon M\rightarrow N$ of $A$-modules is said to be {\em weakly injective}\index{weakly!weakly injective@--- injective} (resp.\ {\em weakly surjective}\index{weakly!weakly surjective@--- surjective}) if $\ker(f)$ (resp.\ $\coker(f)$) is bounded $I$-torsion\index{torsion!bounded torsion@bounded ---}, that is, there exists a positive integer $n>0$ such that $I^n\ker(f)=0$ (resp.\ $I^n\coker(f)=0$) (cf.\ {\bf \ref{ch-pre}}.\ref{dfn-bt1}).}
\end{dfn}

The following properties are easy to verify: let $f\colon L\rightarrow M$ and $g\colon M\rightarrow N$ be $A$-linear homomorphisms of $A$-modules; 
\begin{itemize}
\item if $f$ and $g$ are weakly injective (resp.\ weakly surjective), then so is the composition $g\circ f$;
\item if $g\circ f$ is weakly injective, then so is $f$;
\item if $g\circ f$ is weakly surjective, then so is $g$;
\item if $g\circ f$ is weakly injective and $f$ is weakly surjective, then $g$ is weakly injective;
\item if $g\circ f$ is weakly surjective and $g$ is weakly injective, then $f$ is weakly surjective.
\end{itemize}

\begin{dfn}\label{dfn-weakisomorphic}{\rm 
A homomorphism $f\colon M\rightarrow N$ of $A$-modules is called a {\em weak isomorphism} if it is weakly injective and weakly surjective.}
\end{dfn}

By what we have seen above, we easily deduce:
\begin{itemize}
\item for two homomorphisms $f\colon L\rightarrow M$ and $g\colon M\rightarrow N$ of $A$-modules such that $g\circ f$ is a weak isomorphism, if one of $f$ and $g$ is a weakly isomorphism, then so is the other.
\end{itemize}

It is straightforward to see that the family of all weak isomorphisms forms a multiplicative system (cf.\ \cite[1.6.1]{KS}) in the category $\Mod_A$ of all $A$-modules.
We say that two $A$-modules $M$ and $N$ are {\em weakly isomorphic}\index{weakly!weakly isomorphic@--- isomorphic} if they are isomorphic when regarded as objects in the localized category $\Mod_A/\{\textrm{weak isomorphisms}\}$.
\index{weak!weak isomorphism@--- isomorphism|)}

\subsubsection{Weakly finitely presented modules}\label{subsub-weakfinitelygenerated}
\index{weakly!weakly finitely generated modules@--- finitely generated (modules)|(}
\index{weakly!weakly finitely presented modules@--- finitely presented (modules)|(}
\begin{dfn}\label{dfn-weaklyfinitelygenerated}{\rm 
An $A$-module $M$ is said to be {\em weakly finitely generated} (resp.\ {\em weakly finitely presented}) if it is weakly isomorphic to a finitely generated (resp.\ finitely presented) $A$-module.}
\end{dfn}

\begin{prop}\label{prop-weaklyfinitelygenerated}
{\rm (1)} An $A$-module $M$ is weakly finitely generated if and only if it has a finitely generated $A$-submodule $N\subseteq M$ such that $M/N$ is bounded $I$-torsion.

{\rm (2)} An $A$-module $M$ is weakly finitely presented if and only if there exists a weak isomorphism $N\rightarrow M$ from a finitely presented $A$-module.
\end{prop}

\begin{proof}
To prove (1), it suffices to show the following:
\begin{itemize}
\item[(a)] let $g\colon N_1\rightarrow N_2$ be a weak isomorphism, and $L_2\subseteq N_2$ a finitely generated $A$-submodule such that $N_2/L_2$ is bounded $I$-torsion; then there exists a finitely generated $A$-submodule $L_1\subseteq N_1$ such that $N_1/L_1$ is bounded $I$-torsion;
\item[(b)] if $h\colon L\rightarrow N$ be a weak isomorphism with $L$ finitely generated, then $N/h(L)$ is bounded $I$-torsion.
\end{itemize}
The property (b) is clear.
To show (a), let $x_1,\ldots,x_r$ be a set of generators of $L_2$.
There exists $n>0$ such that $ax_j$ belongs of the image of $g$ ($j=1,\ldots,r$) for any $a\in I^n$.
Let $a_1,\ldots,a_s\in I^n$ generate $I^n$, and take $y_{ij}\in N_1$ such that $g(y_{ij})=a_ix_j$ for $i=1,\ldots,s$ and $j=1,\ldots,r$.
Let $L_1$ be the $A$-submodule of $N_1$ generated by all $y_{ij}$'s.
Then since the composition $L_1\rightarrow N_1\rightarrow N_2$ is weakly injective, the map $L_1\rightarrow L_2$ is weakly injective; since $L_1\rightarrow L_2$ is obviously weakly surjective, we deduce that it is a weak isomorphism.
Hence the inclusion $L_1\hookrightarrow N_1$ is a weak isomorphism, as desired.
To show (2), one only needs:
\begin{itemize}
\item[(c)] let $g\colon N_1\rightarrow N_2$ be a weak isomorphism, and $h_2\colon L_2\rightarrow N_2$ a weak isomorphism from a finitely presented $A$-module $L_2$; then there exists a finitely presented $A$-module $L_1$ and weak isomorphisms $h_1\colon L_1\rightarrow N_1$ and $L_1\rightarrow L_2$.
\end{itemize}
The checking is straightforward and is left to the reader.
\end{proof}
\index{weakly!weakly finitely presented modules@--- finitely presented (modules)|)}
\index{weakly!weakly finitely generated modules@--- finitely generated (modules)|)}

\subsection{Global approximation by finitely presented sheaves}\label{sub-weaklycoherentsheaves}
\subsubsection{FP-approximation of sheaves on schemes}\label{subsub-weaklycoherentsheaves}
Let $X$ be a coherent scheme, $Z\hookrightarrow X$ a closed subscheme of finite presentation, and $\mathscr{I}=\mathscr{I}_Z$ the defining ideal of $Z$.
The notions `weakly injective'\index{weakly!weakly injective@--- injective}, `weakly surjective'\index{weakly!weakly surjective@--- surjective}, and `weak isomorphism'\index{weak!weak isomorphism@--- isomorphism} have obvious analogues for morphisms of $\O_X$-modules; for example, a morphism $\varphi\colon\mathscr{F}\rightarrow\mathscr{G}$ of $\O_X$-modules is a weak isomorphism\index{weak!weak isomorphism@--- isomorphism} if and only if there exists a positive integer $n>0$ such that $\mathscr{I}^n\ker(\varphi)=0$ and $\mathscr{I}^n\coker(\varphi)=0$, that is, $\ker(\varphi)$ and $\coker(\varphi)$ are bounded $\mathscr{I}$-torsion.
It is straightforward to see that the family of all weak isomorphisms (resp.\ between quasi-coherent $\O_X$-modules) is a multiplicative system in the category $\Mod_X$ (resp.\ $\QCoh_X$) of all $\O_X$-modules (resp.\ quasi-coherent $\O_X$-modules).

In this paragraph we assume that:
\begin{itemize}
\item $X$ is {\em Noetherian outside $Z$}, that is, the open subscheme $X\setminus Z$ is Noetherian;
\item $(X,Z)$ satisfies {\bf (UBT)}, that is, for any finite type map $f\colon\Spec A\rightarrow X$ the induced pair $(A,I)$ (where $\til{I}=f^{\ast}\mathscr{I}$) satisfies {\bf (BT)} in {\rm {\bf \ref{ch-pre}}, \S\ref{subsub-BTfirstprop}}.
\end{itemize}
\begin{dfn}\label{dfn-FPthickenings}{\rm 
Let $\mathscr{F}$ be a quasi-coherent sheaf on $X$ of finite type.

(1) By an {\em FP-approximation}\index{FPapproximation@FP-approximation} of $\mathscr{F}$ we mean a weak isomorphism $\varphi\colon\mathscr{G}\rightarrow\mathscr{F}$ from a finitely presented $\O_X$-module.

(2) An FP-approximation $\varphi\colon\mathscr{G}\rightarrow\mathscr{F}$ of $\mathscr{F}$ is called an {\em FP-thickening}\index{FPthickening@FP-thickening} if it is surjective.}
\end{dfn}

If $\mathscr{F}$ admits an FP-approximation, then we sometimes say that $\mathscr{F}$ is {\em FP-approximated}\index{FPapproximated@FP-approximated}.
Clearly, a quasi-coherent sheaf that admits an FP-thickening is of finite type.
Conversely, we have:
\begin{prop}\label{prop-labelFPthickeningexist}
In the situation as above, any quasi-coherent sheaf of finite type on $X$ admits an FP-thickening.
\end{prop}

The proposition is easy to see in case $X$ is affine; this follows from what we have seen in the end of the previous paragraph.
To show the proposition in general, we first introduce the category of FP-thickenings by the following notion of morphisms: given two FP-thickenings\index{FPthickening@FP-thickening} $\mathscr{G}_1\rightarrow\mathscr{F}$ and $\mathscr{G}_2\rightarrow\mathscr{F}$ of $\mathscr{F}$, a morphism from the former to the latter is a surjective morphism $\mathscr{G}_1\rightarrow\mathscr{G}_2$ that makes the triangle
$$
\xymatrix@R-4ex{\mathscr{G}_1\ar[dd]\ar[dr]\\ &\mathscr{F}\\ \mathscr{G}_2\ar[ur]}
$$
commutative.

\begin{lem}\label{lem-categoryFPthickeningfiltered}
Let $\mathscr{F}$ be a quasi-coherent sheaf of finite type on $X$.
Then the category of all FP-thickenings\index{FPthickening@FP-thickening} of $\mathscr{F}$ is a filtered\index{category!cofiltered category@cofiltered ---} category $($cf.\ {\rm {\bf \ref{ch-pre}},\ \S\ref{subsub-finalcofinal}}$)$.
\end{lem}

\begin{proof}
Let $\mathscr{G}_1\rightarrow\mathscr{F}$ and $\mathscr{G}_2\rightarrow\mathscr{F}$ be FP-thickenings. 
We need to construct another FP-thickening $\mathscr{H}\rightarrow\mathscr{F}$ that is dominated by $\mathscr{G}_1$ and $\mathscr{G}_2$.
Consider first the fiber product $\mathscr{K}$ of the maps $\mathscr{G}_1\rightarrow\mathscr{F}$ and $\mathscr{G}_2\rightarrow\mathscr{F}$.
It is easy to see that $\mathscr{K}$ is a quasi-coherent subsheaf of the direct product $\mathscr{G}_1\oplus\mathscr{G}_2$, and that the canonical morphisms $\mathscr{G}_i\rightarrow\mathscr{G}_1\oplus\mathscr{G}_2/\mathscr{K}$ $(i=1,2)$ are surjective.
Write $\mathscr{K}$ as an inductive limit $\varinjlim_{\lambda\in\Lambda}\mathscr{K}_{\lambda}$ of quasi-coherent subsheaves of finite type (\cite[$\mathbf{I}$, (9.4.9) \& $\mathbf{IV}$, (1.7.7)]{EGA}).
We need to show that for some $\lambda\in\Lambda$ the morphisms $\mathscr{G}_i\rightarrow\mathscr{G}\oplus\mathscr{H}/\mathscr{K}_{\lambda}$ $(i=1,2)$ are surjective.
To check this, we may assume that $X$ is affine $X=\Spec A$.
Set $\mathscr{F}=\til{M}$ and $\mathscr{G}_i=\til{N}_i$ for $i=1,2$; $M$ is a finitely generated $A$-module, and $N_1$ and $N_2$ are finitely presented $A$-modules.

Let $\{y_{ij}\}_{j\in J_i}$ be a finite set of generators of $N_i$ $(i=1,2)$.
Denote by $\ovl{y}_{ij}$ the image of $y_{ij}$ in $M$.
Take a lift $x_{1j}$ (resp.\ $x_{2j}$) of $\ovl{y}_{1j}$ (resp.\ $\ovl{y}_{2j}$) in $N_2$ (resp.\ $N_1$).
For some $\lambda\in\Lambda$, $K_{\lambda}$ (where $K^{\Delta}_{\lambda}=\mathscr{K}_{\lambda}$) contains the elements $(y_{1j},x_{1j})$ $(j\in J_1)$ and $(x_{2j},y_{2j})$ $(j\in J_2)$.
Then the maps $N_i\rightarrow N_1\oplus N_2/K_{\lambda}$ are surjective.
\end{proof}

\begin{proof}[Proof of Proposition {\rm \ref{prop-labelFPthickeningexist}}]
As we have already seen above, the assertion is true if $X$ is affine.
Consider a finite open covering $X=\bigcup^r_{i=1}U_i$ such that the assertion is true on each $U_i$ (e.g.\ finite affine covering).
By induction with respect to $r$, we may work in the situation where $r=2$.
Let $\mathscr{F}$ be a quasi-coherent sheaf of finite type on $X$, and $\mathscr{G}_i\rightarrow\mathscr{F}|_{U_i}$ $(i=1,2)$ an FP-thickening on $U_i$.
Take, on $U_1\cap U_2$, an FP-thickening $\mathscr{H}\rightarrow\mathscr{F}|_{U_1\cap U_2}$ dominates by $\mathscr{G}_i|_{U_1\cap U_2}$ for $i=1,2$ (\ref{lem-categoryFPthickeningfiltered}).
Consider $\mathscr{K}_i=\ker(\mathscr{G}_i|_{U_1\cap U_2}\rightarrow\mathscr{H})$ $(i=1,2)$, which is a quasi-coherent sheaf of finite type, and is bounded $\mathscr{I}$-torsion.
Hence there exists a quasi-coherent sheaf $\til{\mathscr{K}}_i$ ($\mathscr{I}$-torsion subsheaf of $\mathscr{G}_i$) of finite type on $U_i$ that extends $\mathscr{K}_i$ (for $i=1,2$).
Now the quotient sheaves $\mathscr{G}_1/\til{\mathscr{K}}_1$ and $\mathscr{G}_2/\til{\mathscr{K}}_2$ patch together to a finitely presented sheaf $\mathscr{H}$, which gives an FP-thickening of $\mathscr{F}$.
\end{proof}

\begin{cor}\label{cor-FPapproxfgenough}
A quasi-coherent sheaf $\mathscr{F}$ is FP-approximated if and only if there exists a quasi-coherent subsheaf $\mathscr{G}\subseteq\mathscr{F}$ of finite type such that $\mathscr{F}/\mathscr{G}$ is bounded $\mathscr{I}$-torsion.
\end{cor}

\begin{proof}
The `only if' part is clear; one only has to take the image of an FP-approximation.
Suppose $\mathscr{F}$ has a subsheaf $\mathscr{G}$ as above.
By \ref{lem-categoryFPthickeningfiltered} we have an FP-thickening $\mathscr{H}\rightarrow\mathscr{G}$, where $\mathscr{H}$ is finitely presented.
Then the composition $\mathscr{H}\rightarrow\mathscr{F}$ gives an FP-approximation.
\end{proof}

\begin{thm}\label{thm-FPapproxthickabelian2}
The full subcategory $\FPA_{(X,Z)}$ of the category $\Mod_X$ of $\O_X$-modules consisting of FP-approximated sheaves is a thick\index{subcategory!thick subcategory@thick ---} {\rm ({\bf \ref{ch-pre}}, \S\ref{sub-subcategoryderived})} abelian full subcategory.
\end{thm}

As a first step of the proof, we show:
\begin{lem}\label{lem-FPapproxthickabelian1}
Let $\varphi\colon\mathscr{F}\rightarrow\mathscr{G}$ be a morphism between FP-approximated sheaves on $X$.
Then $\ker(\varphi)$ and $\coker(\varphi)$ are FP-approximated.
\end{lem}

\begin{proof}
First notice that $\ker(\varphi)$ and $\coker(\varphi)$ are quasi-coherent sheaves on $X$.
To show that $\ker(\varphi)$ is FP-approximated, first take an FP-approximation $\alpha\colon\mathscr{F}'\rightarrow\mathscr{F}$ of $\mathscr{F}$, and consider the composition $\varphi\circ\alpha\colon\mathscr{F}'\rightarrow\mathscr{G}$.
Since the image of $\varphi\circ\alpha$ is finitely presented outside $Z$, $\ker(\varphi\circ\alpha)$ is finitely generated outside $Z$.
By \cite[$\mathbf{I}$, (9.4.7) \& $\mathbf{IV}$, (1.7.7)]{EGA} we have a quasi-coherent subsheaf $\mathscr{K}$ of $\ker(\varphi\circ\alpha)$ of finite type such that $\ker(\varphi\circ\alpha)/\mathscr{K}$ is $\mathscr{I}$-torsion; since this is a subsheaf of the finite type $\mathscr{F}'/\mathscr{K}$, the $\mathscr{I}$-torsion is bounded.
Set $\mathscr{G}'=\mathscr{F}'/\mathscr{K}$, which is a finitely presented sheaf sitting the commutative diagram
$$
\xymatrix{\mathscr{F}\ar[r]^{\varphi}&\mathscr{G}\\ \mathscr{F}'\ar[u]^{\alpha}\ar[r]_{\varphi'}&\mathscr{G}'\ar[u]_{\beta}\rlap{,}}
$$
where $\alpha$ has bounded $\mathscr{I}$-torsion kernel and cokernel, and $\beta$ has $\mathscr{I}$-torsion kernel.
By snake lemma the induced morphism $\ker(\varphi')\rightarrow\ker(\varphi)$ has bounded $\mathscr{I}$-torsion cokernel.
Now by \ref{cor-FPapproxfgenough} we deduce that $\ker(\varphi)$ is FP-approximated.

Next, let us show that $\coker(\varphi)$ is FP-approximated.
To this end, take an FP-approximation $\beta\colon\mathscr{G}'\rightarrow\mathscr{G}$ of $\mathscr{G}$; since the composition $\mathscr{G}'\rightarrow\coker(\varphi)$ obviously has bounded $\mathscr{I}$-torsion cokernel, we immediately deduce the desired result again due to \ref{cor-FPapproxfgenough}.
\end{proof}

\begin{proof}[Proof of Theorem {\rm \ref{thm-FPapproxthickabelian2}}]
In view of \ref{lem-FPapproxthickabelian1} the essential point to show is the following: let 
$$
0\longrightarrow\mathscr{F}\longrightarrow\mathscr{G}\longrightarrow\mathscr{H}\longrightarrow 0
$$
an exact sequence of $\O_X$-modules, where $\mathscr{F}$ and $\mathscr{H}$ are FP-approximated; then $\mathscr{G}$ is FP-approximated.
Notice first that $\mathscr{G}$ is quasi-coherent (\cite[(1.4.7)]{EGAInew}).
In view of \ref{cor-FPapproxfgenough} we are going to construct a quasi-coherent subsheaf $\mathscr{G}'\subseteq\mathscr{G}$ of finite type such that $\mathscr{G}/\mathscr{G}'$ is bounded $\mathscr{I}$-torsion.

Take a quasi-coherent subsheaf $\mathscr{H}'$ (resp.\ $\mathscr{F}'$) of $\mathscr{H}$ (resp.\ $\mathscr{F}$) of finite type such that $\mathscr{H}/\mathscr{H}'$ (resp.\ $\mathscr{F}/\mathscr{F}'$) is bounded $\mathscr{I}$-torsion.
Let $\{\mathscr{G}_{\lambda}\}_{\lambda\in\Lambda}$ be directed set of quasi-coherent subsheaves of $\mathscr{G}$ of finite type such that $\mathscr{G}=\varinjlim_{\lambda\in\Lambda}\mathscr{G}_{\lambda}$ (\cite[$\mathbf{I}$, (9.4.9) \& $\mathbf{IV}$, (1.7.7)]{EGA}).
Then there exists $\lambda\in\Lambda$ such that the image of $\mathscr{G}_{\lambda}$ in $\mathscr{H}$ contains $\mathscr{H}'$ and that the preimage $\mathscr{G}_{\lambda}\cap\mathscr{F}$ contains $\mathscr{F}'$.
Now we deduce by snake lemma that $\mathscr{G}/\mathscr{G}_{\lambda}$ is bounded $\mathscr{I}$-torsion, as desired.
\end{proof}

\subsubsection{FP-approximation of sheaves on formal schemes}\label{subsub-weaklycoherentsheavesformal}
Similarly to the scheme case, one can define the notion of weak isomorphisms on formal schemes: Let $X$ be a coherent adic formal scheme of finite ideal type\index{formal scheme!adic formal scheme@adic ---!adic formal scheme of finite ideal type@--- --- of finite ideal type}\index{adic!adic formal scheme@--- formal scheme!adic formal scheme of finite ideal type@--- --- of finite ideal type} (\ref{dfn-adicformalschemesoffiniteidealtype}).
Then a morphism of $\O_X$-modules $\varphi\colon\mathscr{F}\rightarrow\mathscr{G}$ is a weak isomorphism if and only if there exists a positive integer $n>0$ such that $\mathscr{I}^n\ker(\varphi)=0$ and $\mathscr{I}^n\coker(\varphi)=0$, where $\mathscr{I}$ is an ideal of definition of finite type of $X$.

Now let $X$ be a coherent universally rigid-Noetherian formal scheme\index{formal scheme!universally rigid-Noetherian formal scheme@universally rigid-Noetherian ---} {\rm (\ref{dfn-formalsch})}.
By \ref{prop-tuaadeq2}, for any affine open $\Spf A\subseteq X$ the ring $A$ is a t.u.\ rigid-Noetherian ring\index{t.u. rigid-Noetherian ring@t.u.\ rigid-Noetherian ring} {\rm (\ref{dfn-tuaringadmissible} (1))}; in particular, it satisfies (together with a finitely generated ideal of definition) the condition {\bf (UBT)}.

Similarly to the scheme case, we define:
\begin{dfn}\label{dfn-FPthickeningsformal}{\rm 
Let $X$ be a coherent universally rigid-Noetherian formal scheme, and $\mathscr{F}$ an $\O_X$-module.

(1) By an {\em FP-approximation}\index{FPapproximation@FP-approximation} of $\mathscr{F}$ we mean a weak isomorphism $\varphi\colon\mathscr{G}\rightarrow\mathscr{F}$ from a finitely presented $\O_X$-module $\mathscr{G}$ such that the sheaves $\ker(\varphi)$ and $\coker(\varphi)$ are adically quasi-coherent\index{quasi-coherent!adically quasi-coherent OX module@adically --- (a.q.c.) sheaf} (\ref{dfn-adicqcoh}).

(2) An FP-approximation $\varphi\colon\mathscr{G}\rightarrow\mathscr{F}$ of $\mathscr{F}$ is called an {\em FP-thickening}\index{FPthickening@FP-thickening} if it is surjective.}
\end{dfn}

The terminology `FP-approximated'\index{FPapproximated@FP-approximated} will be used similarly. 
Notice that the sheaf $\mathscr{F}$ is not assumed to be adically quasi-coherent.
The following proposition is easy to see:
\begin{prop}\label{prop-FPAcompletion}
Let $X$ be a coherent scheme, $Z\hookrightarrow X$ a closed subscheme of finite presentation such that $X$ is Noetherian outside $Z$ and that $(X,Z)$ satisfies {\bf (UBT)}.
Consider the formal completion $\widehat{X}=\widehat{X}|_Z$, which we assume to be universally rigid-Noetherian.
Let $\mathscr{F}$ be a quasi-coherent sheaf on $X$ that admits an FP-approximation $($resp.\ FP-thickening$)$.
Then the sheaf $\mathscr{F}^{\for}$ on $\widehat{X}$ obtained by the functor $\for\colon\Mod_X\rightarrow\Mod_{\widehat{X}}$ defined in {\rm \ref{subsub-GFGAcomannconst}} admits an FP-approximation $($resp.\ FP-thickening$)$. \hfill$\square$
\end{prop}

If $\mathscr{F}$ is adically quasi-coherent of finite type and if $X$ is affine, then, similarly to the scheme case, $\mathscr{F}$ admits an FP-thickening.
More generally we have:
\begin{prop}\label{prop-labelFPthickeningexistformal}
Let $X$ be a coherent rigid-Noetherian formal scheme.
Then any adically quasi-coherent sheaf of finite type admits an FP-thickening. \hfill$\square$
\end{prop}

The proposition can be shown by an argument similar to the scheme case by means of:
\begin{lem}\label{lem-categoryFPthickeningfilteredformal}
Let $\mathscr{F}$ be a quasi-coherent sheaf of finite type on $X$.
Then the category of all FP-thickenings\index{FPthickening@FP-thickening} of $\mathscr{F}$ is a filtered\index{category!cofiltered category@cofiltered ---} category. \hfill$\square$
\end{lem}

The proof of the lemma is quite similar to that of \ref{lem-categoryFPthickeningfiltered} (use Exercise \ref{exer-lemdeltasheafadicallyquasicoherent1-application} and Exercise \ref{exer-extopenadiccoh1approx}).

\subsection{Finiteness theorem and GFGA theorems}\label{sub-weakGFGAtheorems}
The proofs of the results in this subsection will be given in \cite{FK}.
\subsubsection{Finiteness theorem for FP-approximated sheaves}\label{sub-finitenessweaklycoherentmodules}
Let $Y$ be a quasi-compact scheme, and $Z\hookrightarrow Y$ a closed subspace of finite presentation.
We suppose that $Y$ is Noetherian outside $Z$ and that the pair $(Y,Z)$ satisfies {\bf (UBT)} (cf.\ \S\ref{subsub-weaklycoherentsheaves}).
Consider a morphism $f\colon X\rightarrow Y$ of finite type; $X$ is Noetherian outside $W=f^{-1}(Z)$ and the pair $(X,W)$ satisfies the condition {\bf (UBT)}.
Let us denote by $\DC^{\ast}_{\fpa}(X,W)$ (resp.\ $\DC^{\ast}_{\fpa}(Y,Z)$) for $\ast=$``\ \ '', $+$, $-$, $\bd$ the full subcategory of the derived category $\DC^{\ast}(X)$ (resp.\ $\DC^{\ast}(Y)$) consisting of objects with FP-approximated cohomologies; by \ref{thm-FPapproxthickabelian2} it is a triangulated full subcategory.

\begin{thm}\label{thm-finitenessweaklycoherentmodules}
In the situation as above, we assume that $f$ is proper.
Then $\RD f_{\ast}$ maps $\DC^{\ast}_{\fpa}(X,W)$ to $\DC^{\ast}_{\fpa}(Y,Z)$ for $\ast=$``\ \ '', $+$, $-$, $\bd$. \hfill$\square$
\end{thm}

\subsubsection{GFGA comparison theorem in rigid-Noetherian situation}\label{sub-GFGAcomparisonweaklycoherentmodules}
\index{GFGA!GFGA comparison theorem@--- comparison theorem}
Similarly to the previous paragraph, let $Y$ be a quasi-compact algebraic space, and $Z\hookrightarrow Y$ a closed subspace of finite presentation, and suppose that $Y$ is Noetherian outside $Z$ and that the pair $(Y,Z)$ satisfies {\bf (UBT)}.
For a finite type morphism $f\colon X\rightarrow Y$ we set $W=f^{-1}(Z)$.
We consider the formal completions $\widehat{X}=\widehat{X}|_W$ and $\widehat{Y}=\widehat{Y}|_Z$ and the functors $\for\colon\Mod_X\rightarrow\Mod_{\widehat{X}}$ (and similarly on $Y$) defined by $\mathscr{F}\mapsto j^{\ast}\mathscr{F}$, where $j\colon\widehat{X}\rightarrow X$ is the canonical morphism (cf.\ \S\ref{subsub-GFGAcomannconst}).
If $f$ is proper, we can construct the following diagram in the similar way as in \S\ref{subsub-GFGAcomann} for $\ast=$``\ \ '', $+$, $-$, $\bd$:
$$
\xymatrix@-2ex{
\DC^{\ast}_{\fpa}(X,W)\ar[rr]^{\for}\ar[dd]_{\RD f_{\ast}}&&\DC^{\ast}(\widehat{X})\ar[dd]^{\RD\widehat{f}_{\ast}}\\
&&\ \\
\DC^{\ast}_{\fpa}(Y,Z)\ar[rr]_{\for}&\ar@{=>}[ur]^{\rho}&\DC^{\ast}(\widehat{Y})\rlap{.}}\leqno{(\ast)}
$$

\begin{thm}[GFGA comparison theorem]\label{thm-GFGAcomweakcoherent}
The above diagram $(\ast)$ is $2$-commutative, that is, the natural transformation $\rho$ gives a natural equivalence for $\ast=$``\ \ '', $+$, $-$, $\bd$.\hfill$\square$
\end{thm}

The theorem essentially follows from the following special case:
\begin{thm}\label{thm-GFGAcomweakcoherentspecial}
Let $A$ be a ring, and $I\subseteq A$ a finitely generated ideal. 
Suppose that the pair $(A,I)$ satisfies {\bf (UBT)} and that the ring $A$ is Noetherian outside $I$.
Let $f\colon X\rightarrow\Spec A$ be a proper morphism.
Then for any FP-approximated sheaf $\mathscr{F}$ on $X$ the canonical map
$$
(\RD^qf_{\ast}\mathscr{F})^{\for}\longrightarrow\RD^q\widehat{f}_{\ast}\widehat{\mathscr{F}}
$$
is an isomorphism for any $q\geq 0$, where $\widehat{\,\cdot\,}$ denotes the $I$-adic completion. \hfill$\square$
\end{thm}

We will need the following implication from the theorem:
\begin{cor}\label{cor-GFGAcomweakcoherentHom}
Let $X$ be as in {\rm \ref{thm-GFGAcomweakcoherentspecial}}, and $\mathscr{F}$ and $\mathscr{G}$ two $\O_X$-modules.
Suppose that $\mathscr{F}$ is finitely presented and that $\mathscr{G}$ is FP-approximated.
Then the canonical map
$$
\Hom_{\O_X}(\mathscr{F},\mathscr{G})\longrightarrow\Hom_{\O_{\widehat{X}}}(\mathscr{F}^{\for},\mathscr{G}^{\for})
$$
is an isomorphism. \hfill$\square$
\end{cor}

\subsubsection{GFGA existence theorem in rigid-Noetherian situation}\label{sub-GFGAexistenceweaklycoherentmodules}
\index{GFGA!GFGA existence theorem@--- existence theorem}
Let $A$ be a t.u.\ rigid-Noetherian ring\index{t.u. rigid-Noetherian ring@t.u.\ rigid-Noetherian ring} (\ref{dfn-tuaringadmissible}), and $f\colon X\rightarrow\Spec A$ a proper map of schemes.
By the $I$-adic completion (where $I\subseteq A$ is a finitely generated ideal of definition) we have the universally rigid-Noetherian formal scheme\index{formal scheme!universally rigid-Noetherian formal scheme@universally rigid-Noetherian ---} (\ref{dfn-formalsch}) $\widehat{f}\colon\widehat{X}\rightarrow\Spf A$ over $A$.

As in \S\ref{subsub-pfexistred} we say that an $\O_{\widehat{X}}$-module $\mathscr{F}$ is algebraizable\index{algebraizable} if there exists an $\O_X$-module $\mathscr{G}$ such that $\mathscr{F}\cong\mathscr{G}^{\for}$.

\begin{thm}\label{thm-GFGAexaweakcoherentspecial}
Suppose $f$ is projective.
Then any finitely presented $\O_{\widehat{X}}$-module $\mathscr{F}$ is algebraizable\index{algebraizable}. \hfill$\square$
\end{thm}

\addcontentsline{toc}{subsection}{Exercises}
\subsection*{Exercises}
\begin{exer}\label{exer-weakisomcomplete}{\rm 
Let $A$ be a ring, and $I\subseteq A$ a finitely generated ideal.
We assume that the ring $A$ with the $I$-adic topology satisfies the condition {\bf (AP)} in {\bf \ref{ch-pre}}, \S\ref{subsub-ARIgoodnesssub}.
Let $f\colon N\rightarrow M$ be a weak isomorphism of $A$-modules, where $N$ is assumed to be finitely generated.
Then show that, if either one of $N$ and $M$ is $I$-adically complete, then so is the rest.}
\end{exer}

\begin{exer}\label{exer-FPapproximationfunctoriality}{\rm 
Let $X$ be either:
\begin{itemize}
\item a coherent scheme considered with a closed subscheme $Z\hookrightarrow X$ of finite presentation such that $X\setminus Z$ is Noetherian and $(X,Z)$ satisfies {\bf (UBT)};
\item a coherent universally rigid-Noetherian formal shceme.
\end{itemize}
Let $\varphi\colon \mathscr{F}\rightarrow\mathscr{G}$ a morphism of sheaves on $X$ that admits an FP-approximation (resp.\ FP-thickening).
Then show that there exist a commutative diagram
$$
\xymatrix{\mathscr{F}\ar[r]^{\varphi}&\mathscr{G}\\ \mathscr{F}'\ar[u]^{\alpha}\ar[r]_{\varphi'}&\mathscr{G}'\ar[u]_{\beta}\rlap{,}}
$$
where $\alpha$ and $\beta$ are FP-approximations (resp.\ FP-thickenings).}
\end{exer}
\setcounter{section}{0}
\renewcommand{\thesection}{\arabic{section}}
\renewcommand{\theexer}{{\bf \thechapter}.\arabic{section}.\arabic{exer}}
\chapter{Rigid spaces}\label{ch-rigid}
This chapter is the main part of this volume, where we define rigid spaces and develop the geometry of them.
In the first section \S\ref{sec-blowup} we discuss generalities of admissible blow-ups.
We then give the definition of rigid spaces in \S\ref{sec-cohrigidispace}, according to Raynaud's\index{Raynaud, M.} viewpoint:
We first define the category of {\em coherent} rigid spaces as the quotient of the category of coherent adic formal schemes of finite ideal type mod out by admissible blow-ups.
Thus, any coherent rigid space $\mathscr{X}$ is, by definition, of the form 
$$
X^{\rig}
$$
by a coherent adic formal scheme $X$ of finite ideal type, and $X$ in this situation is called a {\em formal model} of $\mathscr{X}$.
We then define general (not necessarily coherent) rigid spaces by `patching'.
Corresponding to universally adhesive and universally rigid-Noetherian formal schemes ({\bf \ref{ch-formal}}, \S\ref{sub-formalsch}), we have respectively {\em universally adhesive} and {\em universally Noetherian} rigid spaces:
$$
\bigg\{\begin{minipage}{6.4em}\begin{center}{\small univ.\ adhesive rigid spaces}\end{center}\end{minipage}\bigg\}\ \subseteq\bigg\{\begin{minipage}{7em}\begin{center}{\small univ.\ Noetherian rigid spaces}\end{center}\end{minipage}\bigg\}\ \subseteq\bigg\{\begin{minipage}{5em}\begin{center}{\small rigid spaces}\end{center}\end{minipage}\bigg\}.
$$
`Classical' rigid spaces (called {\em rigid spaces of type} (V) in this book), that is, locally of finite type rigid spaces over $(\Spf V)^{\rig}$ where $V$ is an $a$-adically complete valuation ring, are an example of universally adhesive rigid spaces.

In \S\ref{sec-embodying} we introduce the {\em visualization}, the {\em Zariski-Riemann triple} 
$$
\ZRT(\mathscr{X})=(\ZR{\mathscr{X}},\O^{\int}_{\mathscr{X}},\O_{\mathscr{X}})
$$
of a rigid space $\mathscr{X}$, which `visualizes' the rigid space in the sense that the space $\mathscr{X}$, introduced at the first place by an abstract categorical argument, is interpreted into a concrete topological space with the `doubly-ringed structure', which we call a {\em triple}.
The triples thus obtained are, in important cases, Huber's adic spaces (cf.\ \S\ref{sec-adicspaces} in the Appendix).
We also discuss points of Zariski-Riemann spaces by means of {\em rigid points}, which are, similarly to the situation of Zariski's classical birational geometry mentioned in Introduction, described in terms of valuations.
Having thus the notion of visualization of rigid spaces, one is then able to consider several `topological properties' of rigid spaces, some of which we introduce and develop in \S\ref{sec-topologicalproperties}.
We will show, most importantly, that Zariski-Riemann spaces are {\em valuative} (Corollary \ref{cor-valuativerigidZRsp}), hence admitting the so-called {\em separated quotient}
$$
\sep_{\mathscr{X}}\colon\ZR{\mathscr{X}}\longrightarrow[\mathscr{X}].
$$
Roughly speaking, when one regards the Zariski-Riemann space $\ZR{\mathscr{X}}$ as the `space of arbitrary valuations', the separated quotient $[\mathscr{X}]$ is the subset (endowed with, however, the quotient topology by $\sep_{\mathscr{X}}$) of $\ZR{\mathscr{X}}$ consisting of valuations of height one, and the map $\sep_{\mathscr{X}}$ is given by `maximal generization'.
Notice that, therefore, the space $[\mathscr{X}]$ can be regarded as the `space of seminorms', and hence is a space of the same kind as those appearing in Berkovich\index{Berkovich, V.G.} analytic geometry; in fact, as we will explain in \S\ref{sub-vsbr-Berkovichanalyticspaces}, the separated quotient $[\mathscr{X}]$ coincides with the underlying topological space of the Berkovich analytic space associated to $\mathscr{X}$, at least when $\mathscr{X}$ is locally of finite type over $(\Spf V)^{\rig}$, where $V$ is an $a$-adically complete valuation ring of height one (called {\em rigid space of type {\rm ($\mathrm{V_{\R}}$)}} in this book).
Moreover, several related classical features in Tate's rigid analytic geometry or in Berkovich analytic geometry, such as overconvergent structure, can also be `visualized' entirely by means of usual point-set topology techniques. 
It is perhaps one of the most powerful aspects of our visualization that, in this way, many useful concepts in rigid geometry can be simply boiled down to (often elementary) general topology.

From \S\ref{sec-coherentsheavesrigid} we go for `analytic geometry' of our rigid spaces.
After discussing coherent sheaves on universally Noetherian rigid spaces in \S\ref{sec-coherentsheavesrigid}, we then proceed to the theory of {\em affinoids} in \S\ref{sec-affinoids}.
In this book, affinoids are defined as rigid spaces of the form $(\Spf A)^{\rig}$, that is, coherent rigid spaces having an affine formal model.
Among them, especially important are universally Noetherian affinoids of this form where $\Spec A\setminus V(I)$ is affine.
Affinoids of this type are called {\em Stein affinoids}, the name coined from the fact that these affinoids enjoy Theorem A and Theorem B for coherent sheaves and thus can be viewed as an analogue of Stein domains in complex analysis.
The last-mentioned fact is based on the comparison theorem (Theorem \ref{thm-comparisonaffinoid}) for affinoids, which roughly asserts that the cohomology of coherent sheaves on a universally Noetherian affinoid $(\Spf A)^{\rig}$ is isomorphic to the cohomology of the corresponding coherent sheaves on the Noetherian scheme $\Spec A\setminus V(I)$.

In \S\ref{sub-basicmorproprigid} we collect basic properties of morphisms between rigid spaces, such as finite morphisms, immersions, separated morphisms, etc. 
In \S\ref{sec-localring} we develop some useful tools to investigate points of (the visualizations of) rigid spaces, and apply them to the study of {\em classical points} on rigid spaces of type {\rm ($\mathrm{V_{\R}}$)} or on rigid spaces locally having Noetherian formal models (called: {\em rigid spaces of type} (N)).
The notion of classical points will play an important role in establishing the bridge between our rigid geometry and Tate's rigid analytic geometry, which will be explained in the appendix \S\ref{sec-berkovich}.
Another useful application of our study on classical points is Noetherness theorem (Theorem \ref{thm-noetherness}), which asserts that the local ring at each point of a rigid space of type (V) or of type (N) is Noetherian.

In \S\ref{sec-GAGA} we discuss GAGA.
Our GAGA functor is a functor
$$
X\longmapsto X^{\an}
$$
from the category of quasi-separated finite type schemes over $U=\Spec A\setminus V(I)$, where $A$ is an adic ring with a finitely generated ideal of definition $I\subseteq A$, to the category of rigid spaces of finite type over $\mathscr{S}=(\Spf A)^{\rig}$.
The construction of GAGA functor relies on a generalization of Nagata's embedding theorem, which we prove in the appendix (Theorem \ref{thm-nagataembedding}), and thus will be carried out via formal geometry.
Hence the GAGA theorems (comparison and existence) are the implications from the GFGA theorems that we have already established in Chapter {\bf \ref{ch-formal}}; notice that, like GFGA theorems, our GAGA theorems will be presented in the derived category language. 

We briefly give in \S\ref{sec-dimension} the dimension theory in rigid geometry, and in the final section \S\ref{sec-pointsV1}, we discuss maximal modulus principle. 

\section{Admissible blow-ups}\label{sec-blowup}
In this section we discuss the generalities and related topics of the so-called {\em admissible blow-ups} of formal schemes.
As rigid spaces (more precisely, {\em coherent} rigid spaces) are defined roughly as `limits' of admissible blow-ups of formal schemes, the notion of admissible blow-ups is, so to speak, a cornerstone of the bridge from formal geometry to rigid geometry.

After discussing basic properties of admissible blow-ups, we define in \S\ref{sub-stricttransform} the so-called {\em strict transforms} of admissible blow-ups and give a collection of basic properties of them.
In the final subsection \S\ref{sub-categoryadmblow-up} we will see that admissible blow-ups of a fixed coherent formal scheme constitute a cofiltered category that admits a small cofinal set.
This fact will be used later when we define coherent rigid spaces.

\subsection{Admissible blow-ups}\label{sub-blowups}
\index{blow-up!admissible blow-up@admissible ---|(}
\subsubsection{Admissible blow-ups}\label{subsub-blowups}
Let $X$ be an adic formal scheme of finite ideal type\index{formal scheme!adic formal scheme@adic ---!adic formal scheme of finite ideal type@--- --- of finite ideal type}\index{adic!adic formal scheme@--- formal scheme!adic formal scheme of finite ideal type@--- --- of finite ideal type} ({\bf \ref{ch-formal}}.\ref{dfn-adicformalschemesoffiniteidealtype}).
We assume for a while that there exists an ideal of definition of finite type $\mathscr{I}\subseteq\O_X$.
For $k\geq 0$ we set $X_k=(X,\O_X/\mathscr{I}^{k+1})$, which is a closed subscheme of $X$.
Let $\mathscr{J}\subseteq\O_X$ be an admissible ideal\index{admissible!admissible ideal@--- ideal} ({\bf \ref{ch-formal}}.\ref{dfn-admissibleideal}).
Consider for each $k\geq 0$ the projective $X_k$-scheme 
$$
{\textstyle X'_k=\Proj(\bigoplus_{n\geq 0}\mathscr{J}^n\otimes\O_{X_k})\longrightarrow X_k.}
$$
Whenever $k\leq l$ we have obvious closed immersions $X'_k\hookrightarrow X'_l$, compatible with $X_k\hookrightarrow X_l$, and thus get an inductive system of schemes $\{X'_k\}_{k\geq 0}$.
It is easy to see that this inductive system satisfies the conditions in {\bf \ref{ch-formal}}.\ref{prop-formalindlimschadic}, and hence the inductive limit $X'=\varinjlim_{k\geq 0}X'_k$ is an adic formal scheme of finite ideal type endowed with the structural adic map
$$
X'=\varinjlim_{k\geq 0}\Proj\bigg(\bigoplus_{n\geq 0}\mathscr{J}^n\otimes\O_{X_k}\bigg)\stackrel{\pi}{\longrightarrow} X,\leqno{(\ast)}
$$
which is proper\index{morphism of formal schemes@morphism (of formal schemes)!proper morphism of formal schemes@proper ---} due to {\bf \ref{ch-formal}}.\ref{prop-propermorformal0}.
Notice that the $X$-isomorphism class of the map $(\ast)$ does not depend on the choice of the ideal of definition $\mathscr{I}$.
Hence the construction of $\pi\colon X'\rightarrow X$ can be done without a globally defined ideal of definition of finite type on $X$ as follows: We consider an open covering $X=\bigcup_{\alpha\in L}X_{\alpha}$ where each $X_{\alpha}$ has an ideal of definition $\mathscr{I}_{\alpha}$ of finite type.
Then, for an admissible ideal $\mathscr{J}$ on $X$, define $\pi\colon X'\rightarrow X$ to be the gluing of $\pi_{\alpha}\colon X'_{\alpha}\rightarrow X_{\alpha}$ constructed as above by means of the admissible ideal $\mathscr{J}|_{X_{\alpha}}$.

\begin{dfn}\label{dfn-blowups}{\rm 
Let $X$ be an adic formal scheme of finite ideal type, and $\mathscr{J}\subseteq\O_X$ an admissible ideal.
An adic morphism $\pi\colon X'\rightarrow X$ of adic formal schemes of finite ideal type is said to be an {\em admissible $($formal$)$ blow-up along $\mathscr{J}$} if it is locally isomorphic to an morphism of the form $(\ast)$.}
\end{dfn}

The admissible blow-ups are uniquely determined up to canonical isomorphisms by admissible ideals.
Notice that, if $X$ is quasi-compact (resp.\ quasi-separated, resp.\ coherent), then so is $X'$.
In the sequel, when we want to spell out the blow-up center $\mathscr{J}$, we write
$$
\pi_{\mathscr{J}}\colon X_{\mathscr{J}}\longrightarrow X.
$$
The following proposition follows immediately from the fact that admissible blow-ups are of finite type:
\begin{prop}\label{prop-admissibleblowupuniversallyrigidnoetherian}
Let $X$ be an adic formal scheme of finite ideal type, and $\mathscr{J}\subseteq\O_X$ an admissible ideal.
If $X$ is locally universally rigid-Noetherian $($resp.\ locally universally adhesive$)$, then so is the admissible blow-up $X_{\mathscr{J}}$ along $\mathscr{J}$. \hfill$\square$
\end{prop}

\subsubsection{Explicit local description}\label{subsub-blowupsdescription}
The formation of admissible blow-ups is an effective local construction with respect to Zariski topology (cf.\ {\bf \ref{ch-pre}}.\ref{dfn-localconstruction}), and hence most of the properties of admissible blow-ups can be verified by reducing to the affine situation. 

Let $A$ be an adic ring of finite ideal type\index{admissible!admissible topological ring@--- (topological) ring!admissible topological ring of finite ideal type@--- --- of finite ideal type}\index{adic!adic ring@--- ring!adic ring of finite ideal type@--- --- of finite ideal type} ({\bf \ref{ch-formal}}.\ref{dfn-admissibleringoffiniteidealtype}), $I\subseteq A$ a finitely generated ideal of definition, and $J$ is an $I$-admissible ideal of $A$.
Then the admissible blow-up of $X=\Spf A$ along $J^{\Delta}$ is the formal completion\index{completion!formal completion@formal ---} (cf.\ {\bf \ref{ch-formal}}, \S\ref{subsub-formalcompletionsch}) of the usual blow-up of the affine scheme $\Spec A$ 
$$
\Proj R(A,J)\longrightarrow\Spec A
$$
(\cite[$\mathbf{II}$, (8.1.3)]{EGA}), where $R(A,J)=\bigoplus_{n\geq 0}J^n$ is the Rees algebra\index{Rees algebra} (cf.\ {\bf \ref{ch-pre}}, \S\ref{sub-Reescone}).

Suppose $f_0,\ldots,f_r\in J$ generate $J$, and consider the exact sequence
$$
0\longrightarrow K\longrightarrow A[X_0,\ldots,X_r]\stackrel{\varphi}{\longrightarrow}R(A,J)\longrightarrow 0,
$$
where $\varphi$ maps each $X_i$ to $f_i$ in degree $1$.
Let $\mathfrak{a}\subseteq A[X_0,\ldots,X_r]$ be the ideal generated by all elements of the form 
$$
f_iX_j-f_jX_i
$$
for $0\leq i,j\leq r$.
By \cite[$\mathbf{II}$, (2.9.2) (i)]{EGA} the map $\varphi$ gives rise to an $A$-closed immersion 
$$
\Proj R(A,J)\longhookrightarrow\P^r_A, 
$$
which induces, by {\bf \ref{ch-formal}}.\ref{cor-closedimmformal41}, the closed immersion
$$
\widehat{\Proj R(A,J)}\longhookrightarrow\widehat{\P}^r_A,
$$
into the formal projective $r$-space $\widehat{\P}^r_A$ over $\Spf A$ (cf.\ Execise \ref{exer-affinespaceformal} (2)).
As $\widehat{\P}^r_A$ is covered by the affine open subsets $\{U_i=\Spf A\dl \frac{X_0}{X_i},\ldots,\frac{X_r}{X_i}\dr\}_{i=0,\ldots,r}$, we want to describe explicitly the closed immersion 
$$
\widehat{\Proj R(A,J)}\times_{\widehat{\P}^r_A}U_i\longhookrightarrow U_i\leqno{(\dagger)}
$$
for each $0\leq i\leq r$.
By {\bf \ref{ch-formal}}.\ref{prop-prodformal22} the morphism $(\dagger)$ is the formal completion of the closed immersion of schemes induced by the morphism of rings
$$
{\textstyle A[\frac{X_0}{X_i},\ldots,\frac{X_r}{X_i}]\longrightarrow R(A,J)_{(f_i)}=A[\frac{X_0}{X_i},\ldots,\frac{X_r}{X_i}]/\mathfrak{a}_i^{f_i\textrm{-}\mathrm{sat}},}
$$
where $\mathfrak{a}_i=(f_i\frac{X_0}{X_i}-f_0,\ldots,f_i\frac{X_r}{X_i}-f_r)$, and $\mathfrak{a}_i^{f_i\textrm{-}\mathrm{sat}}$ is the $f_i$-saturation of $\mathfrak{a}$; notice that in the ring $R(A,J)_{(f_i)}$ the ideal $JR(A,J)_{(f_i)}=(f_i)$ is invertible (cf.\ \cite[$\mathbf{II}$, (8.1.11)]{EGA}).
Hence the morphism $(\dagger)$ is isomorphic to
$$
{\textstyle 
\Spf B\longhookrightarrow\Spf A\dl\frac{X_0}{X_i},\ldots,\frac{X_r}{X_i}\dr,}
$$
where
$$
{\textstyle 
B=A\dl\frac{X_0}{X_i},\ldots,\frac{X_r}{X_i}\dr/\ovl{\mathfrak{a}_i^{f_i\textrm{-}\mathrm{sat}}A\dl\frac{X_0}{X_i},\ldots,\frac{X_r}{X_i}\dr}}.\leqno{(\ddagger)}
$$
(Recall {\bf \ref{ch-pre}}, \S\ref{sub-powerseries} for the definition of restricted formal power series rings\index{restricted formal power series}.)

\begin{prop}\label{prop-explicitdescription}
Suppose $A$ is a t.u.\ rigid-Noetherian ring\index{t.u. rigid-Noetherian ring@t.u.\ rigid-Noetherian} {\rm ({\bf \ref{ch-formal}}.\ref{dfn-tuaringadmissible} (1))}, and consider the ideal 
$$
{\textstyle \mathfrak{b}_i=(f_i\frac{X_0}{X_i}-f_0,\ldots,f_i\frac{X_r}{X_i}-f_r)\subseteq A\dl\frac{X_0}{X_i},\ldots,\frac{X_r}{X_i}\dr.}
$$

{\rm (1)} We have $\mathfrak{b}_i^{f_i\textrm{-}\mathrm{sat}}=\mathfrak{a}_i^{f_i\textrm{-}\mathrm{sat}}A\dl\frac{X_0}{X_i},\ldots,\frac{X_r}{X_i}\dr$ and $B=A\dl\frac{X_0}{X_i},\ldots,\frac{X_r}{X_i}\dr/\mathfrak{b}_i^{f_i\textrm{-}\mathrm{sat}}$.

{\rm (2)} If $A$ is an $I$-torsion free t.u.\ adhesive ring\index{t.u.a. ring@t.u.\ adhesive ring} {\rm ({\bf \ref{ch-formal}}.\ref{dfn-tuaringadmissible} (2))}, then $\mathfrak{b}^{f_i\textrm{-}\mathrm{sat}}$ is finitely generated and the ring $B$ is $IB$-torsion free.
\end{prop}

\begin{proof}
(1) First notice that the following inclusions hold:
$$
{\textstyle \mathfrak{b}_i\subseteq\mathfrak{a}_i^{f_i\textrm{-}\mathrm{sat}}A\dl\frac{X_0}{X_i},\ldots,\frac{X_r}{X_i}\dr\subseteq\mathfrak{b}_i^{f_i\textrm{-}\mathrm{sat}}}.
$$
Hence, to show the equality $\mathfrak{b}_i^{f_i\textrm{-}\mathrm{sat}}=\mathfrak{a}_i^{f_i\textrm{-}\mathrm{sat}}A\dl\frac{X_0}{X_i},\ldots,\frac{X_r}{X_i}\dr$, it suffices to show that $\mathfrak{a}_i^{f_i\textrm{-}\mathrm{sat}}A\dl\frac{X_0}{X_i},\ldots,\frac{X_r}{X_i}\dr$ is $f_i$-saturated. 
Since $A[\frac{X_0}{X_i},\ldots,\frac{X_r}{X_i}]$ is $I$-adically universally pseudo-adhesive, the map
$$
{\textstyle A[\frac{X_0}{X_i},\ldots,\frac{X_r}{X_i}]\longrightarrow A\dl\frac{X_0}{X_i},\ldots,\frac{X_r}{X_i}\dr}
$$
is flat ({\bf \ref{ch-pre}}.\ref{prop-btarf1}), and hence so is its base change 
$${\textstyle 
A[\frac{X_0}{X_i},\ldots,\frac{X_r}{X_i}]/\mathfrak{a}_i^{f_i\textrm{-}\mathrm{sat}}\longrightarrow
A\dl\frac{X_0}{X_i},\ldots,\frac{X_r}{X_i}\dr/\mathfrak{a}_i^{f_i\textrm{-}\mathrm{sat}}A\dl\frac{X_0}{X_i},\ldots,\frac{X_r}{X_i}\dr.}
$$
Since the left-hand side is $f_i$-torsion free, so is the right-hand side.
Hence $\mathfrak{a}_i^{f_i\textrm{-}\mathrm{sat}}A\dl\frac{X_0}{X_i},\ldots,\frac{X_r}{X_i}\dr$ is $f_i$-saturated, as desired.
Now by {\bf \ref{ch-pre}}.\ref{cor-propARconseq1-2} we know that the ideal $\mathfrak{b}_i^{f_i\textrm{-}\mathrm{sat}}$ is closed, whence the second equality.

(2) If $A$ is $I$-torsion free, then so is the rings $R(A,J)$ and $R(A,J)_{(f_i)}$.
This implies that $\mathfrak{a}_i^{f_i\textrm{-}\mathrm{sat}}$ is $I$-saturated.
Since $A[\frac{X_0}{X_i},\ldots,\frac{X_r}{X_i}]$ is $I$-adically adhesive, $\mathfrak{a}_i^{f_i\textrm{-}\mathrm{sat}}$ is finitely generated ({\bf \ref{ch-pre}}.\ref{prop-adhesive}).
Since the completion map
$$
{\textstyle A[\frac{X_0}{X_i},\ldots,\frac{X_r}{X_i}]/\mathfrak{a}^{f_i\textrm{-}\mathrm{sat}}\longrightarrow B=A\dl\frac{X_0}{X_i},\ldots,\frac{X_r}{X_i}\dr/\mathfrak{b}_i^{f_i\textrm{-}\mathrm{sat}}}
$$
is flat and the ring $A[\frac{X_0}{X_i},\ldots,\frac{X_r}{X_i}]/\mathfrak{a}^{f_i\textrm{-}\mathrm{sat}}$ is $I$-torsion free, we conclude that $B$ is $IB$-torsion free.
\end{proof}

It follows from (2) that, if $X$ is a locally universally adhesive formal scheme\index{formal scheme!universally adhesive formal scheme@universally adhesive ---!locally universally adhesive formal scheme@locally --- ---}\index{adhesive!universally adhesive@universally ---!universally adhesive formal scheme@--- --- formal scheme} {\rm ({\bf \ref{ch-formal}}.\ref{dfn-formalsch})} such that $\O_X$ is $\mathscr{I}$-torsion free, where $\mathscr{I}$ is an ideal of definition of $X$, then any admissible blow-up of $X$ is finitely presented.

\subsubsection{Universal mapping property}\label{subsub-blowupsuniversal}
\begin{prop}\label{prop-blowups1}
Let $X$ be an adic formal scheme of finite ideal type, $\mathscr{J}\subseteq\O_X$ an admissible ideal, and $\pi\colon X'\rightarrow X$ the admissible blow-up of $X$ along $\mathscr{J}$.

{\rm (1)} The morphism $\pi$ is proper. 

{\rm (2)} The ideal $\mathscr{J}\O_{X'}=(\pi^{-1}\mathscr{J})\O_{X'}$ is invertible.
In particular, $\O_{X'}$ is $\mathscr{J}$-torsion free.

{\rm (3)} $(${\em Universal mapping property}$)$ Given an adic morphism $\theta\colon Z\rightarrow X$ of adic formal schemes of finite ideal type such that $(\theta^{-1}\mathscr{J})\O_Z$ is invertible, there exists uniquely a morphism $Z\rightarrow X'$ such that the resulting triangle
$$
\xymatrix@C-4ex{Z\ar@{-->}[rr]\ar[dr]_{\theta}&&X'\ar[dl]^{\pi}\\ &X}
$$
commutes.
\end{prop}

\begin{proof}
(1) is clear by definition. 
To show (2), we may assume that $X$ is affine $X=\Spf A$, where $A$ is an adic ring of finite ideal type.
What to show is that the ring $B$ as in $(\ddagger)$ in \S\ref{subsub-blowupsdescription} is $f_i$-torsion free.
We will show this in \ref{lem-completiontorsionfreeadmissible} below.

To show (3), we may reduce to the affine situation $X=\Spf A$ and $Z=\Spf B$; let $J\subseteq A$ be an admissible ideal such that $J^{\Delta}=\mathscr{J}$.
We may moreover assume that $JB$ is a principal ideal generated by a non-zero-divisor $a\in B$.
By the universality of blow-ups of schemes we have the unique morphism from $\Spec B$ to the blow-up of $\Spec A$ along $J$, from which the desired morphism $Z\rightarrow X'$ is obtained by completion.
The uniqueness is easy to see.
\end{proof}

\begin{lem}\label{lem-completiontorsionfreeadmissible}
Let $A$ be a ring, $I\subseteq A$ a finitely generated ideal, and $f\in A$ an element.
Suppose that the ideal $J=(f)$ is $I$-admissible\index{admissible!I-admissible@$I$-{---}} {\rm ({\bf \ref{ch-pre}}.\ref{dfn-adm})} and that $A$ is $f$-torsion free.
Then the $I$-adic completion $\widehat{A}$ is $f$-torsion free.
\end{lem}

\begin{proof}
We may assume $I\subseteq J$.
We want to show that the map $\widehat{A}\rightarrow\widehat{A}$ given by $x\mapsto f^n x$ is injective for any $n\geq 1$.
Since the image $J^n$ of $A\rightarrow A$ by $x\mapsto f^nx$ is an open ideal of $A$, the subspace topology on $J^n$ induced from the $I$-adic topology on $A$ is the $I$-adic topology.
Hence by {\bf \ref{ch-pre}}.\ref{prop-projlimleftexact} the injectivity of $\widehat{A}\rightarrow\widehat{A}$ follows from the injectivity of $A\rightarrow A$ by $x\mapsto f^n x$.
\end{proof}

As a corollary of \ref{lem-completiontorsionfreeadmissible} we have:
\begin{cor}\label{cor-blowups155}
Let $X$ be an adic formal scheme with an ideal of definition of finite type $\mathscr{I}$, and $X'\rightarrow X$ an admissible blow-up.
If $\O_X$ is $\mathscr{I}$-torsion free, then $\O_{X'}$ is $\mathscr{I}\O_{X'}$-torsion free.
In particular, if $\mathscr{I}$ is invertible on $X$, then $\mathscr{I}\O_{X'}$ is invertible on $X'$. \hfill$\square$
\end{cor}

\begin{proof}
We may assume $X$ is affine $X=\Spf A$; $X'\rightarrow\Spf A$ is the admissible blow-up along an admissible ideal $J\subseteq A$.
Let $I\subseteq A$ be a finitely generated ideal of definition such that $I^{\Delta}=\mathscr{I}$; we may assume $I\subseteq J$.
Using the notation as in \S\ref{subsub-blowupsdescription}, we have that $B$ is $f$-torsion free, where $JB=(f)$.
Since $\widehat{B}$ is $f$-torsion free by \ref{lem-completiontorsionfreeadmissible}, it is $I\widehat{B}$-torsion free, for $I\widehat{B}\subseteq J\widehat{B}$.
\end{proof}

\subsubsection{Some basic properties}\label{subsub-basicpropertyadmblowups}
The following proposition follows immediately from the universality (\ref{prop-blowups1} (3)) and the local description (\S\ref{subsub-blowupsdescription}) of admissible blow-ups:
\begin{prop}\label{prop-admissibleblowupmorphism}
Let $X\rightarrow Y$ be an adic morphism $($resp.\ of finite type$)$ between adic formal schemes of finite ideal type, and $Y'\rightarrow Y$ an admissible blow-up.
Then there exists an admissible blow-up $X'\rightarrow X$ and an adic morphism $X'\rightarrow Y'$ $($resp.\ of finite type$)$ such that the diagram
$$
\xymatrix{X'\ar[r]\ar[d]&Y'\ar[d]\\ X\ar[r]&Y}
$$
commutes.
Explicitly, if $\mathscr{J}\subseteq\O_Y$ is an admissible ideal that gives the admissible blow-up $Y'\rightarrow Y$, then one can take $X'\rightarrow X$ to be the admissible blow-up along the admissible ideal $\mathscr{J}\O_X$ $($cf.\ {\rm {\bf \ref{ch-formal}}.\ref{prop-admissibleideal1x}}$)$. \hfill$\square$
\end{prop}

Since the formation of admissible blow-ups is an effective local construction, we have:
\begin{prop}\label{prop-blowups41}
Let $X$ be an adic formal scheme of finite ideal type, and $U\subseteq X$ an open subset.
Let $\mathscr{J}$ be an admissible ideal on $X$, and $X'\rightarrow X$ and $U'\rightarrow U$ the admissible blow-ups along $\mathscr{J}$ and $\mathscr{J}|_U$, respectively.
Then we have $U'\cong X'\times_XU;$ in particular, the induced morphism $U'\rightarrow X'$ is an open immersion. \hfill$\square$
\end{prop}

\begin{prop}[Extension of admissible blow-up]\label{prop-blowups4111}
Let $X$ be a coherent adic formal scheme of finite ideal type, $U\subseteq X$ a quasi-compact open subset, and $U'\rightarrow U$ an admissible blow-up.
Then there exists an admissible blow-up $X'\rightarrow X$ that admits an open immersion $U'\hookrightarrow X'$ such that the resulting square
$$
\xymatrix{U'\ \ar@{^{(}->}[r]\ar[d]&X'\ar[d]\\ U\ \ar@{^{(}->}[r]&X}
$$
is Cartesian.
\end{prop}

\begin{proof}
Let $\mathscr{J}$ be an admissible ideal on $U$ that gives the admissible blow-up $U'\rightarrow U$.
By {\bf \ref{ch-formal}}.\ref{prop-extension3} there exists an admissible ideal $\til{\mathscr{J}}$ on $X$ that extends $\mathscr{J}$.
By \ref{prop-blowups41} the admissible blow-up $X'\rightarrow X$ along $\til{\mathscr{J}}$ has the desired property.
\end{proof}

\begin{prop}\label{prop-blowups3}
The composition of two admissible blow-ups between coherent adic formal schemes of finite ideal type is again an admissible blow-ups.
\end{prop}

Notice that the analogous fact is known in the scheme case (\cite[Premi\`ere partie, (5.1.4)]{RG}, which we present below in (\ref{prop-blowup3a})).
Before the proof of the proposition we show the following lemma:
\begin{lem}\label{lem-blowups3prooflem}
Let $Y$ be a quasi-compact scheme, and $\mathscr{I}\subseteq\O_Y$ a quasi-coherent ideal of finite type.
Let $X=\widehat{Y}$ be the formal completion of $Y$ along $\mathscr{I}$.
Then for any admissible ideal\index{admissible!admissible ideal@--- ideal} $\mathscr{J}\subseteq\O_X$ there exists a quasi-coherent ideal $\mathscr{K}\subseteq\O_Y$ of finite type such that the following conditions are satisfied$:$
\begin{itemize}
\item[{\rm (a)}] there exists $n\geq 1$ such that $\mathscr{I}^n\subseteq\mathscr{K};$
\item[{\rm (b)}] $\mathscr{K}\O_X=\mathscr{J}$.
\end{itemize}
\end{lem}

\begin{proof}
Take $n\geq 1$ such that $\mathscr{I}^n\O_X\subseteq\mathscr{J}$.
The sheaf $\mathscr{J}/\mathscr{I}^n$ is a quasi-coherent sheaf of finite type on the scheme $(X,\O_X/\mathscr{I}^n\O_X)=(Y,\O_{Y}/\mathscr{I}^n)$.
Let $\mathscr{K}$ be the pull-back of $\mathscr{J}/\mathscr{I}^n$ by the canonical projection $\O_Y\rightarrow\O_Y/\mathscr{I}^n$; this is a quasi-coherent ideal of finite type such that $\mathscr{I}^n\subseteq\mathscr{K}$ and $\mathscr{K}/\mathscr{I}^n=\mathscr{J}/\mathscr{I}^n$.
By {\bf \ref{ch-formal}}.\ref{cor-extension11} (2) we have $\mathscr{K}\O_{X}=\mathscr{J}$.
\end{proof}

\begin{proof}[Proof of Proposition {\rm \ref{prop-blowups3}}]
Consider the sequence of morphisms of coherent adic formal schemes of finite ideal type
$$
X''\stackrel{\pi'}{\longrightarrow}X'\stackrel{\pi}{\longrightarrow}X,
$$
where $\pi$ is an admissible blow-up along $\mathscr{J}\subseteq\O_X$, and $\pi'$ is an admissible blow-up along $\mathscr{J}'\subseteq\O_{X'}$.
We are going to show that the composition $\pi\circ\pi'$ is an admissible blow-up.
Let $\mathscr{I}\subseteq\O_X$ be an ideal of definition of finite type ({\bf \ref{ch-formal}}.\ref{cor-extension2}), assumed $\mathscr{I}\subseteq\mathscr{J}$ without loss of generality.
We want to construct an adically quasi-coherent ideal $\mathscr{J}''\subseteq\O_X$ of finite type such that
\begin{itemize}
\item[{\rm (a)}] $\mathscr{I}^{k+1}\subseteq\mathscr{J}''$ for some $k\geq 0$;
\item[{\rm (b)}] $\mathscr{J}''\O_{X'}=\mathscr{J}'\mathscr{J}^n\O_{X'}$ for some $n\geq 1$.
\end{itemize}
If this is done, then by the universality of admissible blow-ups (\ref{prop-blowups1} (3)), $\pi\circ\pi'$ is easily seen to be isomorphic to the admissible blow-up of $X$ along $\mathscr{J}\mathscr{J}''$.

\medskip
{\sc Step 1.} We first construct an adically quasi-coherent ideal $\mathscr{J}''$ that satisfies the properties (a) and (b), but is not necessarily of finite type.
We first deal with the affine case $X=\Spf A$, where $A$ is an adic ring with the finitely generated ideal of definition $I\subseteq A$ such that $\mathscr{I}=I^{\Delta}=I\O_X$.
Let $J\subseteq A$ be the $I$-admissible ideal such that $\mathscr{J}=J^{\Delta}$; we may assume $I\subseteq J$ without loss of generality.
Let $q\colon Y'\rightarrow Y=\Spec A$ be the blow-up along the ideal $J$.
Then $\pi\colon X'\rightarrow X$ is the $I$-adic formal completion of $q$, and hence by \ref{lem-blowups3prooflem} we have a quasi-coherent ideal $\mathscr{K}'\subseteq\O_{Y'}$ of finite type such that $I^n\O_{Y'}\subseteq\mathscr{K}'$ for some $n\geq 1$ and $\mathscr{K}\O_{X'}=\mathscr{J}'$.
Let $q'\colon Y''\rightarrow Y'$ be the blow-up along $\mathscr{K}'$.
Then $\pi'\colon X''\rightarrow X'$ is the $I$-adic formal completion of $q'$.

Since $J\O_{Y'}$ is an ample invertible ideal, by \cite[$\mathbf{II}$, (4.6.8)]{EGA}, there exists $N_0\geq 1$ such that $\mathscr{K}'\otimes_{\O_{Y'}}J^n\O_{Y'}=\mathscr{K}'J^n\O_{Y'}$ is generated by global sections for any $n\geq N_0$.
Moreover, since $q_{\ast}\O_{Y'}/\O_Y$ is $J$-torsion, there exists $N_1$ ($\geq N_0$) such that for any $n\geq N_1$ the ideal sheaf $\mathscr{K}'J^n\O_{Y'}$ is generated by global sections from $\O_Y$.
Hence if we define $\mathscr{K}''$ to be the kernel of
$$
\O_Y\longrightarrow q_{\ast}\O_{Y'}/q_{\ast}\mathscr{K}'J^n\O_{Y'}
$$
for $n\geq N_1$, it satisfies
\begin{itemize}
\item[$\mathrm{(a)'}$] $\mathscr{I}^{k+1}\O_Y\subseteq\mathscr{K}''$ for some $k\geq 0$;
\item[$\mathrm{(b)'}$] $\mathscr{K}''\O_{Y'}=\mathscr{K}'J^n\O_{Y'}$.
\end{itemize}
Then $\mathscr{J}''=\mathscr{K}''\O_X$ is an adically quasi-coherent ideal satisfying the conditions (a) and (b) above.

\medskip
{\sc Step 2.} Let $B$ be another adic ring that admits an open immersion $U=\Spf B\hookrightarrow X=\Spf A$; notice that open immersions are adically flat\index{morphism of formal schemes@morphism (of formal schemes)!adically flat morphism of formal schemes@adically flat ---} ({\bf \ref{ch-formal}}.\ref{dfn-adicallyflat}).
Set $V=\Spec B$, and let $q\colon V'\rightarrow V$ be the blow-up along $JB$, which admits a canonical map $V'\rightarrow Y'$, and $q'\colon V''\rightarrow V'$ the blow-up along $\mathscr{K}'\O_{V'}$. 
We thereby have the commutative diagram
$$
\xymatrix{V''\ar[r]^{q'}\ar[d]&V'\ar[r]^q\ar[d]&V\ar[d]\\ Y''\ar[r]_{q'}&Y'\ar[r]_q&Y\rlap{,}}
$$
which induces by the $I$-adic completion the commutative diagram
$$
\xymatrix{U''\ar[r]^{\pi'}\ar[d]&U'\ar[r]^{\pi}\ar[d]&U\ar[d]\\ X''\ar[r]_{\pi'}&X'\ar[r]_{\pi}&X}
$$
with all vertical arrows being open immersions.
Similarly to {\sc Step 1}, we define for a sufficiently large $n$ the quasi-coherent ideal $\mathscr{M}''\subseteq\O_V$ to be the kernel of
$$
\O_V\longrightarrow q_{\ast}\O_{V'}/q_{\ast}\mathscr{K}'J^n\O_{V'}.
$$
We claim that $\mathscr{J}''=\mathscr{K}''\O_X$ as defined in {\sc Step 1} (by the same $n$) enjoys the property that $\mathscr{J}''\O_U=\mathscr{M}''\O_U$.
To show this, since both sides are complete sheaves, it is enough to check that $\mathscr{J}''\O_U/I^{l+1}\O_U=\mathscr{M}''\O_U/I^{l+1}\O_U$ for any $l\geq k$, where $k$ is chosen such that $\mathscr{I}^{k+1}\O_Y\subseteq\mathscr{K}''$ and $\mathscr{I}^{k+1}\O_U\subseteq\mathscr{M}''$.
Since $\mathscr{J}''/I^{l+1}\O_X=\mathscr{K}''/I^{l+1}\O_Y$ is the kernel of $\O_Y/I^{l+1}\O_Y\rightarrow q_{\ast}\O_{Y'}/q_{\ast}\mathscr{K}'J^n\O_{Y'}$ and since $U_l=(U,\O_V/I^{l+1}\O_V)\hookrightarrow X_l=(X,\O_Y/I^{l+1}\O_Y)$ is flat, $\mathscr{J}''\O_{U_l}$ coincides with the kernel of $\O_V/I^{l+1}\O_V\rightarrow q_{\ast}\O_{V'}/q_{\ast}\mathscr{K}'J^n\O_{V'}$, that is, we have $\mathscr{J}''\O_{U_l}=\mathscr{M}''/\O_{U_l}$, as desired.

\medskip
{\sc Step 3.} Now we discuss the general case.
Take a finite affine open covering $X=\bigcup_{\alpha\in L}X_{\alpha}$, and consider for each $\alpha,\beta\in L$ a covering of $X_{\alpha\beta}=X_{\alpha}\cap X_{\beta}$ by finitely many affine open subsets.
If one takes $n$ and $k$ to be sufficiently large, one has for each $\alpha\in L$ an adically quasi-coherent ideal $\mathscr{J}''_{\alpha}$ such that $\mathscr{I}^{k+1}\O_{X_{\alpha}}\subseteq\mathscr{J}''_{\alpha}$ and $\mathscr{J}''\O_{X'_{\alpha}}=\mathscr{J}'\mathscr{J}^n\O_{X'_{\alpha}}$, where $X'_{\alpha}$ is the open subset of $X'$ that is the admissible blow-up of $X_{\alpha}$ along $\mathscr{J}\O_{\alpha}$ (cf.\ \ref{prop-blowups41}).
By what we have seen in {\sc Step 2}, replacing $n$ by a larger one if necessary, we have $\mathscr{J}''_{\alpha}\O_{X_{\alpha\beta}}=\mathscr{J}''_{\beta}\O_{X_{\alpha\beta}}$ for any $\alpha,\beta\in L$.
Hence the ideals $\mathscr{J}''_{\alpha}$ glue to an adically quasi-coherent ideal $\mathscr{J}''$ of $\O_X$, which obviously satisfies the conditions (a) and (b).

\medskip
{\sc Step 4.} In view of \cite[$\mathbf{I}$, (9.4.9) \& $\mathbf{IV}$, (1.7.7)]{EGA} we know that the quasi-coherent ideal $\mathscr{J}''/\mathscr{I}^{k+1}$ on the scheme $X_k=(X,\O_X/\mathscr{I}^{k+1})$ is the inductive limit $\varinjlim_{\lambda\in\Lambda}\ovl{\mathscr{J}''}_{\lambda}$ of quasi-coherent ideals of finite type.
Since $\mathscr{J}'\mathscr{J}^n\O_{X'}$ is finite type, there exists $\lambda\in\Lambda$ such that the pull-back $\mathscr{J}''_{\lambda}$ of $\ovl{\mathscr{J}''}_{\lambda}$ by the canonical projection $\O_X\rightarrow\O_{X_k}$ satisfies the conditions (b); it also satisfies (a) by the construction. 
Hence, replacing $\mathscr{J}''$ by this $\mathscr{J}''_{\lambda}$, we finally get a desired admissible ideal, and thus the proof of the proposition is finished.
\end{proof}

\subsection{Strict transform}\label{sub-stricttransform}
\index{strict transform|(}
\begin{dfn}\label{dfn-stricttransform}{\rm 
Let $X$ be an adic formal scheme of finite ideal type, $\mathscr{J}$ an admissible ideal, and $\pi\colon X'\rightarrow X$ the admissible blow-up of $X$ along $\mathscr{J}$.
For an $\O_X$-module $\mathscr{F}$, the {\em strict transform} of $\mathscr{F}$ by $\pi$ is the $\O_{X'}$-module given by 
$$
\pi'\mathscr{F}=\textrm{the completion of}\ \pi^{\ast}\mathscr{F}/(\pi^{\ast}\mathscr{F})_{\textrm{$\mathscr{J}$-}\mathrm{tor}}
$$
(cf.\ {\bf \ref{ch-formal}}.\ref{dfn-completesheaf} for the definition of the completion of sheaves).}
\end{dfn}

\begin{prop}\label{prop-stricttransformaqctoaqc}
If $\mathscr{F}$ is an adically quasi-coherent sheaf\index{adically quasi-coherent (a.q.c.) sheaf} {\rm ({\bf \ref{ch-formal}}.\ref{dfn-adicqcoh})} $($resp.\ adically quasi-coherent $\O_X$-algebra\index{adically quasi-coherent (a.q.c.) algebra}$)$ on $X$, then the strict transform $\pi'\mathscr{F}$ is an adically quasi-coherent sheaf $($resp.\ adically quasi-coherent $\O_X$-algebra\index{adically quasi-coherent (a.q.c.) algebra}$)$ on $X'$.
If, moreover, $\mathscr{F}$ is of finite type, then so is $\pi'\mathscr{F}$.
\end{prop}

To show the proposition, we need to prepare a few lemmas.
The following lemma is a generalization of \ref{lem-completiontorsionfreeadmissible}, and the proof is similar:
\begin{lem}\label{lem-completiontorsionfree}
Let $A$ be a ring, $I\subseteq A$ a finitely generated ideal, and $f\in A$ an element such that the ideal $(f)$ is open with respect to the $I$-adic topology.
If $M$ is an $f$-torsion free $A$-module, then its $I$-adic completion $\widehat{M}$ is $f$-torsion free. \hfill$\square$
\end{lem}

\begin{cor}\label{cor-completiontorsionfree}
Let $A$ be a ring, $I\subseteq A$ a finitely generated ideal, and $f\in A$ an element such that the ideal $(f)$ is open with respect to the $I$-adic topology.
Then for an $A$-module $M$ we have
$$
\ovl{\widehat{M}_{\ftor}}=\ovl{M_{\ftor}},
$$
where the left-hand side denotes the closure in $\widehat{M}$ of its $f$-torsion part, and the right-hand side is the closure of the image of $M_{\ftor}$ in $\widehat{M}$.
\end{cor}

\begin{proof}
The exact sequence $0\rightarrow M_{\ftor}\rightarrow M\rightarrow M/M_{\ftor}\rightarrow 0$ induces the exact sequence
$$
0\longrightarrow\ovl{M_{\ftor}}\longrightarrow \widehat{M}\longrightarrow\widehat{M/M_{\ftor}}\longrightarrow 0.
$$
The last module is $f$-torsion free by \ref{lem-completiontorsionfree}, and hence $\widehat{M}_{\ftor}\subseteq\ovl{M_{\ftor}}$.
On the other hand, since the image of $M_{\ftor}$ in $\widehat{M}$ is contained in $\widehat{M}_{\ftor}$, we have $\ovl{\widehat{M}_{\ftor}}\supseteq\ovl{M_{\ftor}}$, whence the claim.
\end{proof}

\begin{cor}\label{cor-comppletiontorsionfree2}
Let $A$ be a ring, $I\subseteq A$ a finitely generated ideal, and $B$ an $A$-algebra.
Consider the respective $I$-adic completions $\widehat{A}$ and $\widehat{B}$ of $A$ and $B$.
Let $f\in B$ be an element such that the ideal $(f)$ is open with respect to the $IB$-adic topology.
Then for any finitely generated $A$-module $M$ the $I$-adic completion of $\widehat{M}\otimes_{\widehat{A}}\widehat{B}/(\widehat{M}\otimes_{\widehat{A}}\widehat{B})_{\ftor}$ coincides up to canonical isomorphism with the $I$-adic completion of $M\otimes_AB/(M\otimes_AB)_{\ftor}$.
\end{cor}

\begin{proof}
It follows from \ref{cor-completiontorsionfree} that the closure of the image of $(\widehat{M}\otimes_{\widehat{A}}\widehat{B})_{\ftor}$ in $M\widehat{\otimes}_AB$ coincides with the closure of $(M\widehat{\otimes}_AB)_{\ftor}$, which further coincides, again by \ref{cor-completiontorsionfree}, with the closure of the image of $(M\otimes_AB)_{\ftor}$.
\end{proof}

\begin{lem}\label{lem-comppletiontorsionfree2}
Let $A$ be an adic ring with a finitely generated ideal of definition $I\subseteq A$, and $J\subseteq A$ an admissible ideal.
Set $X=\Spf A$ and $Y=\Spec A$, and consider the commutative diagram
$$
\xymatrix{Y'\ar[d]_p&X'\ar[l]_i\ar[d]^{\pi}\\ Y&X\rlap{,}\ar[l]}
$$
where $p$ is the blow-up along $J$, and $\pi$ is the admissible blow-up along $J$,  the formal completion of $p$.
Then for any $A$-module $($resp.\ $A$-algebra$)$ $M$ we have
$$
\pi'M^{\Delta}=\widehat{p'\til{M}},
$$
where $p'\til{M}$ denotes the scheme-theoretic strict transform of $\til{M}$.
In particular, $\pi'M^{\Delta}$ is an adically quasi-coherent sheaf $($resp.\ adically quasi-coherent $\O_{X'}$-algebra$)$.
If, moreover, $M$ is finitely generated, then $\pi'M^{\Delta}$ is of finite type.
\end{lem}

\begin{proof}
Consider any affine open subset $U=\Spec B$ of $Y'$ such that $JB=(f)$ by a non-zero-divisor $f\in B$.
Then $\pi'M^{\Delta}$ is the sheafification of the presheaf 
$$
\widehat{U}\longmapsto\textrm{the completion of}\ \widehat{M}\otimes_A\widehat{B}/(\widehat{M}\otimes_A\widehat{B})_{\ftor}
$$
and $\widehat{p'\til{M}}$ is the sheafification of
$$
\widehat{U}\longmapsto\textrm{the completion of}\ M\otimes_AB/(M\otimes_AB)_{\ftor}.
$$
These two sheaves coincide with each other due to \ref{cor-comppletiontorsionfree2}.
Since the $f$-torsion part of the quasi-coherent sheaf $p^{\ast}\til{M}$ is quasi-coherent, $\widehat{p'\til{M}}$ is a completion of a quasi-coherent sheaf, and hence is adically quasi-coherent.
\end{proof}

\begin{proof}[Proof of Proposition {\rm \ref{prop-stricttransformaqctoaqc}}]
We may assume that $X$ is affine.
Then the admissible blow-up $\pi\colon X'\rightarrow X$ is the formal completion of a blow-up of an affine scheme as in \ref{lem-comppletiontorsionfree2}, hence the assertion follows from {\bf \ref{ch-formal}}.\ref{thm-adicqcoh1} and \ref{lem-comppletiontorsionfree2}.
\end{proof}

\begin{prop}\label{prop-stricttransformaqctoaqcrigidNoetherian}
Let $X$ be a locally universally rigid-Noetherian formal scheme\index{formal scheme!universally rigid-Noetherian formal scheme@universally rigid-Noetherian ---!locally universally rigid-Noetherian formal scheme@locally --- ---} {\rm ({\bf \ref{ch-formal}}.\ref{dfn-formalsch})}, $\mathscr{J}$ an admissible ideal, and $\pi\colon X'\rightarrow X$ the admissible blow-up along $\mathscr{J}$.
Then for any adically quasi-coherent sheaf of finite type $\mathscr{F}$ we have
$$
\pi'\mathscr{F}=\pi^{\ast}\mathscr{F}/(\pi^{\ast}\mathscr{F})_{\textrm{$\mathscr{J}$-}\mathrm{tor}}.
$$
\end{prop}

\begin{proof}
We may assume that $X$ is affine of the form $X=\Spf A$, where $A$ is a t.u.\ rigid-Noetherian ring\index{t.u. rigid-Noetherian ring@t.u.\ rigid-Noetherian ring} ({\bf \ref{ch-formal}}.\ref{dfn-tuaringadmissible}); we use the notation as in \ref{lem-comppletiontorsionfree2}.
Since $p'\til{M}$ is a quasi-coherent sheaf of finite type on $Y'$, its completion coincides with $i^{\ast}p'\til{M}$ due to {\bf \ref{ch-formal}}.\ref{prop-adequatedeltaconstruction}.
Then the assertion follows from the fact that the morphism $i$ is flat.
\end{proof}

We will indicate below in Exercise \ref{exer-secondproofpropstricttransformaqctoaqcrigidNoetherian} another proof of \ref{prop-stricttransformaqctoaqcrigidNoetherian}.
It can be shown that, in the situation as in \ref{prop-stricttransformaqctoaqcrigidNoetherian}, the strict transform of $\mathscr{J}$ coincides with the ideal pull-back $\mathscr{J}\O_{X'}=(\pi^{-1}\mathscr{J})\O_{X'}$ (cf.\ Exercise \ref{exer-admissibleidealstricttransform}).

Let $f\colon Y\rightarrow X$ be an adic morphism of adic formal schemes of finite ideal type, and $\pi\colon X'\rightarrow X$ the admissible blow-up of $X$ along an admissible ideal $\mathscr{J}\subseteq\O_X$.
Consider the base change $Y\times_XX'\rightarrow X'$.
By \ref{prop-stricttransformaqctoaqc} the completion of $\O_{Y\times_XX'}/\O_{Y\times_XX',\textrm{$\mathscr{J}$-}\mathrm{tor}}$ is an adically quasi-coherent $\O_{Y\times_XX'}$-algebra, which is a quotient of the structure sheaf $\O_{Y\times_XX'}$ (due to {\bf \ref{ch-formal}}.\ref{lem-vanishingcohomologyadicallyuseful} (1)).
Hence one can consider the closed formal subspace $Y'\hookrightarrow Y\times_XX'$ defined by the ideal $\mathscr{K}=\O_{Y\times_XX',\textrm{$\mathscr{J}$-}\mathrm{tor}}$ (cf.\ {\bf \ref{ch-formal}}.\ref{dfn-closedimmform1}).
\begin{dfn}\label{dfn-stricttransform2}{\rm 
In the above situation, the composite map 
$$
Y'\longhookrightarrow Y\times_XX'\longrightarrow X'
$$
is called the {\em strict transform} of $f$ by the admissible blow-up $\pi$.}
\end{dfn}

\begin{prop}\label{prop-blowups15}
Let $\pi\colon X'\rightarrow X$ be an admissible blow-up along an admissible ideal $\mathscr{J}$ of an adic formal scheme $X$ of finite ideal type, and $f\colon Y\rightarrow X$ an adic morphism of adic formal schemes of finite ideal type. 
Let $f'\colon Y'\rightarrow X'$ be the strict transform of $f$ by $\pi$.
Then the morphism $Y'\rightarrow Y$ is the admissible blow-up along the admissible ideal $\mathscr{J}\O_Y$ $($cf.\ {\rm {\bf \ref{ch-formal}}.\ref{prop-admissibleideal1x})}.
\end{prop}

\begin{proof}
Let $Y''\rightarrow Y$ be the admissible blow-up along $\mathscr{J}\O_Y$.
By the universality of admissible blow-ups applied to $X'\rightarrow X$, we get an $X$-morphism $Y''\rightarrow X'$.
Hence we get a morphism $Y''\rightarrow Y\times_XX'$, which induces, in view of \ref{prop-blowups1} (2), a $Y$-morphism $Y''\rightarrow Y'$.
On the other hand, since $\mathscr{J}\O_{Y'}$ is invertible, the universality of admissible blow-ups applied to $Y''\rightarrow Y$ implies that there exists a unique $Y$-morphism $Y'\rightarrow Y''$, which is easily seen to be inverse to the morphism $Y''\rightarrow Y'$.
\end{proof}

\begin{prop}\label{prop-blowups15x}
Let $i\colon Y\hookrightarrow X$ be a closed immersion\index{immersion!closed immersion of formal schemes@closed --- (of formal schemes)} of finite presentation of locally universally rigid-Noetherian formal schemes\index{formal scheme!universally rigid-Noetherian formal scheme@universally rigid-Noetherian ---!locally universally rigid-Noetherian formal scheme@locally --- ---}, and $Y'\rightarrow Y$ an admissible blow-up of $Y$.
Then there exists an admissible blow-up $X'\rightarrow X$ and a closed immersion $Y'\hookrightarrow X'$ that coincides with the strict transform of $i$.
\end{prop}

\begin{proof}
Let $\mathscr{K}$ be the admissible ideal of $\O_Y$ that gives $Y'\rightarrow Y$.
Let $\mathscr{J}$ be the pull-back of $\mathscr{K}$ by the map $\O_X\rightarrow i_{\ast}\O_Y$, which is an admissible ideal ({\bf \ref{ch-formal}}.\ref{prop-coradmmissibleideal1xx}).
Since $\mathscr{K}=\mathscr{J}\O_Y$, the assertion follows from \ref{prop-blowups15}.
\end{proof}

\begin{prop}\label{prop-blowups1552}
Let $X$ be a locally universally rigid-Noetherian formal scheme with an ideal of definition $\mathscr{I}$ such that $\O_X$ is $\mathscr{I}$-torsion free, and $\mathscr{F}$ an $\mathscr{I}$-torsion free adically quasi-coherent $\O_X$-module of finite type\index{adically quasi-coherent (a.q.c.) sheaf!adically quasi-coherent sheaf of finite type@--- of finite type}.
Let $\mathscr{J}$ be an admissible ideal of $X$, and $\pi\colon X'\rightarrow X$ the admissible blow-up along $\mathscr{J}$.
Then we have
$$
(\pi^{\ast}\mathscr{F})_{\textrm{$\mathscr{J}$-}\mathrm{tor}}=(\pi^{\ast}\mathscr{F})_{\textrm{$\mathscr{I}$-}\mathrm{tor}}.
$$
\end{prop}

\begin{proof}
We may assume that $X$ is affine $X=\Spf A$ with $\mathscr{I}=I^{\Delta}$, where $A$ is an $I$-torsion free t.u.\ rigid-Noetherian ring.
Let $\mathscr{J}=J^{\Delta}$ with $J=(f_0,\ldots,f_r)$, and $B$ the ring as in \ref{prop-explicitdescription} (1).
Let $M$ be an $I$-torsion free finitely generated $A$-module such that $\mathscr{F}=M^{\Delta}$.
What we need to show is the equality $(M\otimes_AB)_{\textrm{$f_i$-}\mathrm{tor}}=(M\otimes_AB)_{\textrm{$I$-}\mathrm{tor}}$; note that, since $M$ is finitely generated, we have $M\otimes_AB=M\widehat{\otimes}_AB$ ({\bf \ref{ch-pre}}.\ref{prop-btarf1} (1)).
Since $J$ is admissible, the inclusion $(M\otimes_AB)_{\textrm{$f_i$-}\mathrm{tor}}\subseteq (M\otimes_AB)_{\textrm{$I$-}\mathrm{tor}}$ is clear.
The other inclusion is shown by an argument similar to that in the proof of \ref{prop-explicitdescription}.
Let $P\in R(A,J)_{(f_i)}\otimes_AM$ be an $I$-torsion element, and take $N>0$ large enough so that $f^N_iP\in M$.
Since $M$ is $I$-torsion free, we have $f^N_iP=0$, which implies that $P$ is an $f_i$-torsion.
Hence we have $(R(A,J)_{(f_i)}\otimes_AM)_{\textrm{$f_i$-}\mathrm{tor}}=(R(A,J)_{(f_i)}\otimes_AM)_{\textrm{$I$-}\mathrm{tor}}$.
Since $R(A,J)_{(f_i)}\rightarrow B$ is flat ({\bf \ref{ch-pre}}.\ref{prop-btarf1} (2)), we have the desired equality by base change.
\end{proof}

\begin{cor}\label{cor-blowups1552x}
Let $X$ be a locally universally rigid-Noetherian formal scheme, and $\mathscr{I}$ an ideal of definition of finite type of $X$.
Suppose $\O_X$ is $\mathscr{I}$-torsion free.
Let $\mathscr{F}$ be an $\mathscr{I}$-torsion free adically quasi-coherent sheaf of finite type\index{adically quasi-coherent (a.q.c.) sheaf!adically quasi-coherent sheaf of finite type@--- of finite type} on $X$.
Then for an admissible blow-up $\pi\colon X'\rightarrow X$ the strict transform $\pi'\mathscr{F}$ is an $\mathscr{I}\O_{X'}$-torsion free adically quasi-coherent sheaf of finite type on $X'$. \hfill$\square$
\end{cor}

Notice that, if $X$ is locally universally adhesive, then $\pi'\mathscr{F}$ is an adically quasi-coherent sheaf of finite presentation on $X'$ (cf.\ {\bf \ref{ch-formal}}.\ref{prop-exerItorsionfreeaqcsheavesoffinitetype}).
\index{strict transform|)}

\subsection{The cofiltered category of admissible blow-ups}\label{sub-categoryadmblow-up}
Let $X$ be a coherent adic formal scheme of finite ideal type.
We define $\BL_X$ to be the category of all admissible blow-ups of $X$; more precisely:
\begin{itemize}
\item[$\bullet$] objects of $\BL_X$ are the admissible blow-ups $\pi\colon X'\rightarrow X$;
\item[$\bullet$] an arrow $\pi'\rightarrow\pi$ between two objects $\pi\colon X'\rightarrow X$ and $\pi'\colon X''\rightarrow X$ is a morphism $X''\rightarrow X'$ over $X$:
$$
\xymatrix@C-4.5ex{X''\ar[dr]_{\pi'}\ar[rr]&&X'\ar[dl]^(.55){\,\pi}\\ &X\rlap{.}}
$$
\end{itemize}

\begin{prop}\label{prop-catblowups1}
{\rm (1)} The category $\BL_X$ is cofiltered\index{category!cofiltered category@cofiltered ---} $($cf.\ {\rm {\bf \ref{ch-pre}},\ \S\ref{subsub-finalcofinal}}$)$, and $\id_X$ gives the final object.

{\rm (2)} Define the ordering\index{ordering@order(ing)} on the set $\AId_X$ $(=$ the set of all admissible ideals of $X)$ as follows$:$ $\mathscr{J}\preceq\mathscr{J}'$ if and only if there exists an admissible ideal $\mathscr{J}''$ such that $\mathscr{J}=\mathscr{J}'\mathscr{J}''$ $($cf.\ {\rm {\bf \ref{ch-formal}}.\ref{prop-exerpropadmissibleideal1}}$)$.
Then $\AId^{\opp}_X$ is a directed set\index{set!directed set@directed ---}\index{directed set}, and the functor 
$$
\AId_X\longrightarrow\BL_X
$$
$($where $\AId_X$ is regarded as a category$;$ cf.\ {\rm {\bf \ref{ch-pre}},\ \S\ref{subsub-catposet}}$)$ that maps $\mathscr{J}$ to the admissble blow-up along $\mathscr{J}$ is cofinal\index{functor!cofinal functor@cofinal ---}.
\end{prop}

\begin{proof}
To show (1), we need to check the following:
\begin{itemize}
\item[(a)] for two admissible blow-ups $X'\rightarrow X$ and $X''\rightarrow X$, there exists an admissible blow-up $X'''\rightarrow X$ and $X$-morphisms $X'''\rightarrow X'$ and $X'''\rightarrow X''$; 
\item[(b)] for two admissible blow-ups $X'\rightarrow X$ and $X''\rightarrow X$ and two $X$-morphisms $f_0,f_1\colon X''\rightarrow X'$, there exists an admissible blow-up $X'''\rightarrow X$ with an $X$-morphism $g\colon X'''\rightarrow X''$ such that $f_0\circ g=f_1\circ g$.
\end{itemize}
Let $X'\rightarrow X$ and $X''\rightarrow X$ be admissible blow-ups along $\mathscr{J}$ and $\mathscr{J}'$, respectively.
Let $X'''\rightarrow X$ be the admissible blow-up along $\mathscr{J}\mathscr{J}'$ (cf.\ {\bf \ref{ch-formal}}.\ref{prop-exerpropadmissibleideal1}).
Then by \ref{prop-blowups1} (3) one has the morphisms $X'''\rightarrow X'$ and $X'''\rightarrow X''$ as in (a).
If $f_0$ and $f_1$ are as in (b), we have $f_0\circ g=f_1\circ g$ by the uniqueness in \ref{prop-blowups1} (3), whence (b). 
We have thus shown the assertion (1).
(2) is clear.
\end{proof}

\begin{cor}\label{cor-catblowups101}
The category $\BL_X$ is cofiltered\index{category!cofiltered category@cofiltered ---} and essentially small\index{small!essentially small@essentially ---} $($cf.\ {\rm {\bf \ref{ch-pre}},\ \S\ref{subsub-finalcofinal}}$)$. \hfill$\square$
\end{cor}
\index{blow-up!admissible blow-up@admissible ---|)}

\addcontentsline{toc}{subsection}{Exercises}
\subsection*{Exercises}
\begin{exer}\label{exer-propblowups2}
{\rm Let $X$ be an adic formal scheme of finite ideal type, and $\mathscr{J},\mathscr{J'}\subseteq\O_X$ admissible ideals.
Let $\pi\colon X'\rightarrow X$ be the admissible blow-up along $\mathscr{J}$, and $\pi'\colon X''\rightarrow X'$ the admissible blow-up along $(\pi^{-1}\mathscr{J})\O_{X'}$.
Then show that the composition $\pi\circ\pi'\colon X''\rightarrow X$ coincides up to canonical isomorphism with the admissible blow-up of $X$ along $\mathscr{J}\cdot\mathscr{J}'$.}
\end{exer}

\begin{exer}\label{exer-propblowups4a}
{\rm Let $\pi\colon X'\rightarrow X$ be an admissible blow-up along an admissible ideal $\mathscr{J}$ on $X$, and $f\colon Y\rightarrow X$ an adic morphism of adic formal schemes of finite ideal type.
Then there exist admissible blow-ups $Z\rightarrow Y'=X'\times_XY$ and $Z\rightarrow Y$ such that the resulting diagram commutes:
$$
\xymatrix@C-4ex@R-1.5ex{&Z\ar[dl]\ar[dr]\\ Y'\ar[rr]&&Y}
$$}
\end{exer}

\begin{exer}\label{exer-Jtorsionpartquoteaqc}{\rm 
Let $A$ be a t.u.\ rigid-Noetherian ring\index{t.u. rigid-Noetherian ring@t.u.\ rigid-Noetherian ring} {\rm ({\bf \ref{ch-formal}}.\ref{dfn-tuaringadmissible} (1))}, $I\subseteq A$ a finitely generated ideal of definition, and $J\subseteq A$ an $I$-admissible ideal.
Let $M$ be a finitely generated $A$-module.
Set $\mathscr{F}=M^{\Delta}$ and $\mathscr{J}=J^{\Delta}$.
Then show that $\mathscr{F}_{\textrm{$\mathscr{J}$-}\mathrm{tor}}=(M_{\Jtor})^{\Delta}$.}
\end{exer}

\begin{exer}\label{exer-Jtorsionpartquoteaqc2}{\rm 
Let $X$ be a locally universally rigid-Noetherian formal scheme\index{formal scheme!universally rigid-Noetherian formal scheme@universally rigid-Noetherian ---!locally universally rigid-Noetherian formal scheme@locally --- ---}, $\mathscr{F}$ an adically quasi-coherent sheaf of finite type on $X$, and $\mathscr{J}\subseteq\O_X$ an admissible ideal.
Then show that $\mathscr{F}/\mathscr{F}_{\textrm{$\mathscr{J}$-}\mathrm{tor}}$ is an adically quasi-coherent sheaf of finite type.}
\end{exer}

\begin{exer}\label{exer-secondproofpropstricttransformaqctoaqcrigidNoetherian}{\rm 
Prove \ref{prop-stricttransformaqctoaqcrigidNoetherian} by using Exercise \ref{exer-Jtorsionpartquoteaqc2}.}
\end{exer}

\begin{exer}\label{exer-admissibleidealstricttransform}{\rm
Let $A\rightarrow B$ be an adic morphism of t.u.\ rigid-Noetherian rings, $I\subseteq A$ a finitely generated ideal of definition of $A$, and $J\subseteq A$ an $I$-admissible ideal.
Suppose $B$ is $J$-torsion free.
Set $X=\Spf A$, $Y=\Spf B$, and $\mathscr{J}=J^{\Delta}$.
Let $\mathscr{K}$ be the kernel of $\mathscr{J}\otimes_{\O_X}\O_Y\rightarrow\mathscr{J}\O_Y$, and $K$ the kernel of $J\otimes_AB\rightarrow JB$.

(1) Show that $\mathscr{K}=K^{\Delta}$.

(2) Show that $\mathscr{K}=(\mathscr{J}\otimes_{\O_X}\O_Y)_{\textrm{$\mathscr{J}$-}\mathrm{tor}}$.}
\end{exer}


\section{Rigid spaces}\label{sec-cohrigidispace}
\index{rigid space|(}
In this section, we give the definition of rigid spaces and discuss their first properties of them.
We first define in \S\ref{sub-cohrigidspace} what we call {\em coherent rigid spaces} as the objects in the quotient category of coherent adic formal schemes of finite ideal type mod out by all admissible blow-ups. 
The coherent rigid spaces are, therefore, always of the form `$X^{\rig}$', induced from a coherent adic of finite ideal type formal scheme $X$, which we call a {\em formal model} of the coherent rigid space.
There are of course plenty of formal models attached to a single coherent rigid space, but they are always connected by admissible blow-ups and blow-downs.

A coherent rigid space thus defined comes with a natural topology, the so-called {\em admissible topology}, which will be discussed in \S\ref{sub-admissibletop}.
Roughly speaking, this is the most naturally induced topology from the Zariski topology on the formal schemes.
Since we are, at first, only able to speak about `coherent' rigid spaces, open subsets (or equivalently, open immersions) are temporarily restricted only to quasi-compact ones.
Endowed with this topology, the category of coherent rigid spaces gives rise to the so-called {\em coherent admissible site}.
Here, it turns out that the adjective `coherent' is justified by the fact that the topos associated to the coherent admissible site is in fact coherent in the sense of \cite[Expos\'e VI]{SGA4-2} (cf.\ {\bf \ref{ch-pre}}.\ref{dfn-coherenttopos2}).
It also turns out that the terminologies `coherent' rigid space and `coherent' open immersion (as defined in \ref{dfn-cohrigidspaceopenimm1}) can be justified in the sense that they are exactly the coherent objects in the topos, the so-called {\em admissible topos}. 
Since this justification is most fluently done in terms of {\em visualization}\index{visualization}, we postpone it to the next section.

General rigid spaces are defined to be sheaves on the coherent admissible site that satisfy certain `local representability' property.
Here the already defined coherent rigid spaces are in fact rigid spaces in this generalized sense, since on the coherent admissible site the representable presheaves are sheaves (\ref{prop-admissiblesite3}).
This fact is closely linked with one of the most important aspects of our {\em birational} approach to rigid geometry (cf.\ Introduction): The (coherent) rigid spaces, which are at first defined formally as objects in the above-mentioned quotient category, admit `patching', modeled on `birational' ($=$ up to admissible blow-ups) patching of formal schemes. 

The definition of general rigid spaces given in \S\ref{subsub-generalrigidspace} allows one to enhance the admissible topology into a slightly stronger and more consistent one, also called the admissible topology.
It would be helpful, especially for the reader who is familiar to Tate's rigid analytic geometry, to remark that the admissible topology in the former sense is the one somewhat similar to the so-called `weak topology', and the latter one to the `strong topology' in Tate's rigid analytic geometry (cf.\ \cite{BGR}); see \S\ref{sub-rigidanalyticgeometry} for more precise comparison of the topologies.

In \S\ref{sub-rigidmorphismfintype} we will define the notion of (locally) of finite type morphisms.
The final subsection \S\ref{sub-rigidspacefiberprod} is devoted to fiber products of rigid spaces.

\subsection{Coherent rigid spaces and their formal models}\label{sub-cohrigidspace}
\subsubsection{Coherent rigid spaces}\label{subsub-cohrigidspace}
\index{rigid space!coherent rigid space@coherent ---|(}
\begin{dfn}\label{dfn-cohrigidspace1}{\rm 
We define the category $\CRf$ as follows:
\begin{itemize}
\item[$\bullet$] the objects of $\CRf$ are the same as those of $\Ac\CFs^{\ast}$, the category of coherent adic formal schemes of finite ideal type and adic morphisms:
$$
\obj(\CRf)=\obj(\Ac\CFs^{\ast});
$$
for an object $X$ of $\Ac\CFs^{\ast}$ we denote by $X^{\rig}$ the same object regarded as an object of $\CRf$;
\item[$\bullet$] for $X,Y\in\obj(\Ac\CFs^{\ast})$ we define
$$
\Hom_{\CRf}(X^{\rig},Y^{\rig})=\varinjlim\Hom_{\Ac\CFs^{\ast}}(\,\cdot\,,Y),
$$
where $\Hom_{\Ac\CFs^{\ast}}(\,\cdot\,,Y)$ is the functor $\BL^{\opp}_X\rightarrow\mathbf{Sets}$ that maps $\pi\colon X'\rightarrow X$ to the set $\Hom_{\Ac\CFs^{\ast}}(X',Y)$ (cf.\ {\bf \ref{ch-pre}}, \S\ref{subsub-limdefuniv}).
\end{itemize}}
\end{dfn}

We sometimes describe the quotient functor as
$$
Q\colon\Ac\CFs^{\ast}\longrightarrow\CRf, \qquad X\longmapsto Q(X)=X^{\rig}.
$$
Notice that by \ref{prop-catblowups1} the inductive limit in the above definition can be replaced by a filtered inductive limit along the directed set $\AId^{\opp}_X$ (cf.\ {\bf \ref{ch-pre}}.\ref{prop-final}).

In the category $\CRf$ the composition law
$$
\Hom_{\CRf}(X^{\rig},Y^{\rig})\times\Hom_{\CRf}(Y^{\rig},Z^{\rig})\longrightarrow\Hom_{\CRf}(X^{\rig},Z^{\rig})
$$
is described as follows.
A morphism $\varphi$ in $\Hom_{\CRf}(X^{\rig},Y^{\rig})$ is given by a diagram of the form
$$
\xymatrix@-3ex{
&X'\ar[dr]\ar[dl]\\
X&&Y\rlap{,}}
$$
where $X'\rightarrow X$ is an admissible blow-up\index{blow-up!admissible blow-up@admissible ---} (\ref{dfn-blowups}).
Similarly, a morphism $\psi$ in $\Hom_{\CRf}(Y^{\rig},Z^{\rig})$ is given by $Y\leftarrow Y'\rightarrow Z$, where the first arrow is an admissible blow-up.
By \ref{prop-admissibleblowupmorphism} we have an admissible blow-up $X''\rightarrow X'$ and an adic morphism $X''\rightarrow Y'$ such that the square in the following diagram commutes:
$$
\xymatrix@-3ex{
&&X''\ar[dl]\ar[dr]\\
&X'\ar[dl]\ar[dr]&&Y'\ar[dl]\ar[dr]\\
X&&Y&&Z\rlap{.}}
$$
By \ref{prop-blowups3} one sees that the composition $X''\rightarrow X$ is an admissible blow-up, and hence the diagram $X\leftarrow X''\rightarrow Z$ gives an element in $\Hom_{\CRf}(X^{\rig},Z^{\rig})$.
It can be seen that the resulting element in $\Hom_{\CRf}(X^{\rig},Z^{\rig})$ does not depend on the choice of $X'$ and $Y'$ and thus gives the desired composition $\psi\circ\varphi$.

\begin{dfn}\label{dfn-cohrigidspace}{\rm 
(1) Objects of $\CRf$ are called {\em coherent rigid $($formal$)$ spaces}.

(2) For an object $X$ of $\Ac\CFs^{\ast}$ the coherent rigid space $X^{\rig}$ is called the {\em associated $($coherent$)$ rigid space}.
Similarly, for a morphism $f\colon X\rightarrow Y$ of $\Ac\CFs^{\ast}$, the {\em associated morphism} of rigid spaces is denoted by $f^{\rig}\colon X^{\rig}\rightarrow Y^{\rig}$.}
\end{dfn}

\begin{rem}\label{rem-futurepossiblevar}{\rm 
The adjective `formal' in the parenthesis in (1) and also the letter `{\bf f}' in the notation $\CRf$ are put for specifying that the above-defined rigid spaces come from {\em formal} schemes, indicating future variants including, for example, rigid {\em henselian} spaces and rigid {\em Zariskian} spaces (cf.\ \S\ref{sec-rigidzariskiansp} in the appendix), which are the similarly defined spaces associated respectively to henselian schemes and Zariskian schemes.}
\end{rem}

By an easy but deft use of \ref{prop-blowups1} (3), one can easily show the following:
\begin{prop}\label{prop-cohrigidspace3}
Let $X$ and $Y$ be coherent adic formal schemes of finite ideal type, and consider the rigid spaces $X^{\rig}$ and $Y^{\rig}$.
Then there exists an isomorphism $X^{\rig}\stackrel{\sim}{\rightarrow}Y^{\rig}$ in $\CRf$ if and only if there exists a diagram $X\leftarrow Z\rightarrow Y$ consisting of admissible blow-ups. \hfill$\square$
\end{prop}

\begin{cor}\label{cor-cohrigidspace31}
Let $f\colon X\rightarrow Y$ be a morphism in $\Ac\CFs^{\ast}$.
Then $f^{\rig}\colon X^{\rig}\rightarrow Y^{\rig}$ is an isomorphism if and only there exists a commutative diagram 
$$
\xymatrix@-3ex{
&Z\ar[dr]\ar[dl]\\
X\ar[rr]_{f}&&Y\rlap{,}}
$$
where both $Z\rightarrow X$ and $Z\rightarrow Y$ are admissible blow-ups. \hfill$\square$
\end{cor}

By \ref{cor-cohrigidspace31} and Exercise \ref{exer-propblowups4a} we have:
\begin{cor}\label{cor-cohrigidspace32}
Consider the diagram 
$$
\xymatrix{
X\ar[d]_f&X\times_ZY\ar[d]^h\ar[l]\\
Z&Y\rlap{.}\ar[l]}
$$
If $f^{\rig}$ is an isomorphism in $\CRf$, then so is $h^{\rig}$. \hfill$\square$
\end{cor}
\index{rigid space!coherent rigid space@coherent ---|)}

\subsubsection{Formal models}\label{subsub-cohrigspaceformalmodels}
\begin{dfn}\label{dfn-formalmodelcat}{\rm 
(1) Let $\mathscr{X}$ be a coherent rigid space.
A {\em formal model}\index{formal model!formal model of a coherent rigid space@--- (of a coherent rigid space)} of $\mathscr{X}$ is a couple $(X,\phi)$ consisting of $X\in\obj(\Ac\CFs^{\ast})$ and an isomorphism $\phi\colon X^{\rig}\stackrel{\sim}{\rightarrow}\mathscr{X}$ in $\CRf$.
They constitute the category $\FM_{\mathscr{X}}$, the category of formal models of $\mathscr{X}$, in which an arrow $(X,\phi)\rightarrow(X',\phi')$ is defined to be a morphism $f\colon X\rightarrow X'$ in $\Ac\CFs^{\ast}$ such that $\phi'\circ f^{\rig}=\phi$.

(2) Let $\varphi\colon\mathscr{X}\rightarrow\mathscr{Y}$ be a morphism of coherent rigid spaces.
A {\em formal model} of $\varphi$ is a triple $(f,\phi,\psi)$ consisting of a morphism $f\colon X\rightarrow Y$ of $\Ac\CFs^{\ast}$ and isomorphisms $\phi\colon X^{\rig}\stackrel{\sim}{\rightarrow}\mathscr{X}$ and $\psi\colon Y^{\rig}\stackrel{\sim}{\rightarrow}\mathscr{Y}$ such that the resulting square
$$
\xymatrix{X^{\rig}\ar[d]_{\phi}\ar[r]^{f^{\rig}}&Y^{\rig}\ar[d]^{\psi}\\ \mathscr{X}\ar[r]_{\varphi}&\mathscr{Y}}
$$
commutes.
They constitute the category $\FM_{\varphi}$, the category of formal models of $\varphi$, in which an arrow $(f\colon X\rightarrow Y,\phi,\psi)\rightarrow(f'\colon X'\rightarrow Y',\phi',\psi')$ is defined to be a couple of morphisms $(u,v)$ consisting of $u\colon X\rightarrow X'$ and $v\colon Y\rightarrow Y'$ such that $v\circ f=f'\circ u$ and that the following diagram is commutative:
$$
\xymatrix@-2ex{
&X^{\rig}\ar@/_1pc/[ddl]_{\phi}\ar[d]^{u^{\rig}}\ar[r]^{f^{\rig}}&Y^{\rig}\ar@/^1pc/[ddr]^{\psi}\ar[d]_{v^{\rig}}\\
&X^{\prime\rig}\ar[dl]_{\phi'}\ar[r]_{f^{\prime\rig}}&Y^{\prime\rig}\ar[dr]^{\psi'}\\
\mathscr{X}\ar[rrr]_{\varphi}&&&\mathscr{Y}\rlap{.}}
$$}
\end{dfn}

It is readily seen that these categories are cofiltered (due to \ref{cor-catblowups101}).
\begin{dfn}\label{dfn-cohrigidspacedist}{\rm 
(1) Let $\mathscr{X}$ be a coherent rigid space.
A formal model $(X,\phi)$ of $\mathscr{X}$ is said to be {\em distinguished}\index{formal model!formal model of a coherent rigid space@--- (of a coherent rigid space)!distinguished formal model of a coherent rigid space@distinguished --- ---} if $\O_X$ is $\mathscr{I}$-torsion free for some $($hence any$)$ ideal of definition $\mathscr{I}$ of $X$.

(2) Let $\varphi\colon\mathscr{X}\rightarrow\mathscr{Y}$ be a morphism of coherent rigid spaces.
A formal model $(f\colon X\rightarrow Y,\phi,\psi)$ of $\varphi$ is said to be {\em distinguished} if $X$ and $Y$ are distinguished formal models of $\mathscr{X}$ and $\mathscr{Y}$, respectively.}
\end{dfn}

\begin{prop}[{cf.\ {\rm \ref{cor-blowups155}}}]\label{prop-cohrigidspacedist100}
Let $(X,\phi)$ be a distinguished formal module of a coherent rigid space $\mathscr{X}$, and $\pi\colon X'\rightarrow X$ an admissible blow-up.
Then $(X',\phi\circ\pi^{\rig})$ is a distinguished formal model of $\mathscr{X}$.
If, moreover, $X$ has an invertible ideal of definition $\mathscr{I}$, then $\mathscr{I}\O_{X'}$ is invertible. \hfill$\square$
\end{prop}

We denote by $\FM^{\dist}_{\mathscr{X}}$ (resp.\ $\FM^{\dist}_{\varphi}$) the full subcategory of $\FM_{\mathscr{X}}$ (resp.\ $\FM_{\varphi}$) consisting of distinguished formal models.
\begin{prop}\label{prop-cohrigidspacedist1}
The categories $\FM^{\dist}_{\mathscr{X}}$ and $\FM^{\dist}_{\varphi}$ are cofiltered, and the inclusions $\FM^{\dist}_{\mathscr{X}}\hookrightarrow\FM_{\mathscr{X}}$ and $\FM^{\dist}_{\varphi}\hookrightarrow\FM_{\varphi}$ are cofinal.
$($Hence, in particular, any object and morphism of $\CRf$ have distinguished formal models.$)$
\end{prop}

\begin{proof}
The proposition follows from the following observation: let $X$ be an object of $\Ac\CFs^{\ast}$, and $\mathscr{I}$ an ideal of definition of finite type; then the structure sheaf of the admissible blow-up $X'$ of $X$ along $\mathscr{I}$ is $\mathscr{I}\O_{X'}$-torsion free by \ref{prop-blowups1} (2).
Notice the following: for a morphism $\varphi\colon\mathscr{X}\rightarrow\mathscr{Y}$ of coherent rigid spaces with a formal model $f\colon X\rightarrow Y$, let $Y'\rightarrow Y$ and $X'\rightarrow X$ be the admissible blow-ups along an ideal of definition $\mathscr{I}$ of finite type on $Y$ and $\mathscr{I}\O_X$, respectively; then by \ref{prop-blowups1} (3) we have the unique morphism $f'\colon X'\rightarrow Y'$ (cf.\ \ref{prop-admissibleblowupmorphism}), which gives a distinguished formal model of $\varphi$.
\end{proof}

\begin{dfn}\label{dfn-emptyrigidspace}{\rm 
A coherent rigid space $\mathscr{X}$ is said to be {\em empty}\index{rigid space!empty rigid space@empty ---} if it has an empty formal model.}
\end{dfn}

An empty rigid space will be denoted by $\emptyset$.
For instance, if $\mathscr{X}$ has a formal model that is a {\em scheme} ($=$ $0$-adic formal scheme), then $\mathscr{X}$ is empty.
Later in \ref{cor-ZRnonempty} we will see that a coherent rigid space $\mathscr{X}$ is non-empty if and only if it has a non-empty distinguished formal model.
The following proposition follows easily from the existence of formal models for morphisms of coherent rigid spaces:
\begin{prop}\label{prop-strictinitial}
The empty rigid space $\emptyset$ is a strict initial object {\rm (\cite[Expos\'e II, 4.5]{SGA4-1})} of the category $\CRf$.
\end{prop}

\subsubsection{Comma category $\CRf_{\mathscr{S}}$}\label{subsub-cohrigidspacecomma}
\begin{dfn}\label{dfn-cohrigidspacecomma}{\rm 
Let $\mathscr{S}$ be a coherent rigid space.
We define the category $\CRf_{\mathscr{S}}$ as follows:
\begin{itemize}
\item objects of $\CRf_{\mathscr{S}}$ are morphisms $\mathscr{X}\rightarrow\mathscr{S}$ in $\CRf$ with the target $\mathscr{S}$; 
\item an arrow from $\mathscr{X}\rightarrow\mathscr{S}$ to $\mathscr{Y}\rightarrow\mathscr{S}$ is a morphism $\mathscr{X}\rightarrow\mathscr{Y}$ in $\CRf$ such that the following diagram commutes:
$$
\xymatrix@C-4ex{
\mathscr{X}\ar[dr]\ar[rr]&&\mathscr{Y}\ar[dl]\\
&\mathscr{S}\rlap{.}}
$$
\end{itemize}}
\end{dfn}

Let $S$ be a coherent adic formal scheme of finite ideal type, and set $\mathscr{S}=S^{\rig}$.
Let $\Ac\CFs^{\ast}_{/S}$ be the category of coherent adic formal schemes adic over $S$ (cf.\ {\bf \ref{ch-formal}}, \S\ref{subsub-ntnadequatecategory}).
Then we have the canonical functor
$$
\Ac\CFs^{\ast}_{/S}\longrightarrow\CRf_{\mathscr{S}},\qquad X/S\longmapsto X^{\rig}/S^{\rig}.
$$
On the other hand, one can mimic the definition of $\CRf$ to define the quotient category $\Ac\CFs^{\ast}_{/S}/_{{\textstyle \sim}}$, that is, a category consisting of objects of $\Ac\CFs^{\ast}_{/S}$ with the set of arrows defined by
$$
\Hom_{\Ac\CFs^{\ast}_{/S}/_{{\scriptstyle \sim}}}(X,Y)=\varinjlim_{X'\rightarrow X}\Hom_S(X',Y),
$$
where the right-hand inductive limit is taken over all admissible blow-ups $X'\rightarrow X$.
Then the above functor obviously factors through the canonical functor 
$$
\Ac\CFs^{\ast}_{/S}/_{{\textstyle \sim}}\longrightarrow\CRf_{\mathscr{S}}.\leqno{(\ast)}
$$
\begin{prop}\label{prop-cohrigidspacecomma}
The functor $(\ast)$ is a categorical equivalence.
\end{prop}

\begin{proof}
It suffices to show the following: let $Y\rightarrow S$ be a morphism of coherent adic formal schemes of finite ideal type, and set $\mathscr{Y}=Y^{\rig}$; let $\mathscr{X}\rightarrow\mathscr{Y}$ be a morphism of coherent rigid spaces; then there exists a formal model of the form $X\rightarrow Y$ and, moreover, such a formal model is unique up to admissible blow-ups of $X$.
To see this, take an arbitrary formal model $X'\rightarrow Y'$ of $\mathscr{X}\rightarrow\mathscr{Y}$.
Replacing $Y'$ and $X'$ by admissible blow-ups, we may assume that there exists an admissible blow-up $Y'\rightarrow Y$ (cf.\ \ref{prop-admissibleblowupmorphism} and \ref{prop-catblowups1}).
Hence $X=X'\rightarrow Y$ gives the desired formal model.
The uniqueness is clear.
\end{proof}

\subsubsection{Coherent universally Noetherian and universally adhesive rigid spaces}\label{subsub-coherentuniversallyadhesiverigidspaces}
\index{rigid space!universally adhesive rigid space@universally adhesive ---!coherent universally adhesive rigid space@coherent --- ---|(}
\index{rigid space!universally Noetherian rigid space@universally Noetherian ---!coherent universally Noetherian rigid space@coherent --- ---|(}
\begin{dfn}\label{dfn-coherentuniversallyadhesiverigidspaces}{\rm 
A coherent rigid space $\mathscr{X}$ is said to be {\em universally Noetherian} (resp.\ {\em universally adhesive}) if it has a formal model $(X,\phi)$ by a coherent universally rigid-Noetherian\index{formal scheme!universally rigid-Noetherian formal scheme@universally rigid-Noetherian ---} (resp.\ coherent universally adhesive\index{formal scheme!universally adhesive formal scheme@universally adhesive ---}\index{adhesive!universally adhesive@universally ---!universally adhesive formal scheme@--- --- formal scheme}) formal scheme $X$ ({\bf \ref{ch-formal}}.\ref{dfn-formalsch}).}
\end{dfn}

By \ref{prop-admissibleblowupuniversallyrigidnoetherian} we know that, once $\mathscr{X}$ has a universally rigid-Noetherian (resp.\ universally adhesive) formal model, then any admissible blow-up of it is again universally rigid-Noetherian (resp.\ universally adhesive).
Hence, as in the proof of \ref{prop-cohrigidspacedist1}, one can show that a coherent universally Noetherian (resp.\ universally adhesive) rigid space has a {\em distinguished} universally rigid-Noetherian (resp.\ universally adhesive) formal model.
Notice also that, similarly to the case treated in \ref{prop-cohrigidspacecomma}, one can equivalently define universally Noetherian (resp.\ universally adhesive) rigid spaces as objects in the quotient category of the form $\RNoe\CFs^{\ast}/\sim$ (resp.\ $\Adh\CFs^{\ast}/\sim$) (cf.\ {\bf \ref{ch-formal}}.\ref{subsub-ntnadequatecategory} for the notation) constructed similarly from the category of universally rigid-Noetherian (resp.\ universally adhesive) formal schemes.
\index{rigid space!universally adhesive rigid space@universally adhesive ---!coherent universally adhesive rigid space@coherent --- ---|)}
\index{rigid space!universally Noetherian rigid space@universally Noetherian ---!coherent universally Noetherian rigid space@coherent --- ---|)}

\subsection{Admissible topology and general rigid spaces}\label{sub-admissibletop}
\index{admissible!admissible topology@--- topology|(}
\index{topology!admissible topology@admissible ---|(}
\subsubsection{Coherent admissible sites}\label{subsub-admissibletop}
\index{admissible!admissible site@--- site!coherent admissible site@coherent --- ---|(}
\index{site!admissible site@admissible ---!coherent admissible site@coherent --- ---|(}
\begin{prop}\label{prop-cohrigidspaceopenimm1}
The following conditions for a morphism $\mathscr{U}\rightarrow\mathscr{X}$ of coherent rigid spaces are equivalent$:$
\begin{itemize}
\item[{\rm (a)}] there exists a formal model\index{formal model!formal model of a coherent rigid space@--- (of a coherent rigid space)} $(j,\phi,\psi)$ of $\mathscr{U}\rightarrow\mathscr{X}$ such that the morphism $j\colon U\rightarrow X$ in $\Ac\CFs^{\ast}$ is an open immersion\index{immersion!open immersion of formal schemes@open --- (of formal schemes)}$;$
\item[{\rm (b)}] there exists a {\em distinguished} formal model\index{formal model!formal model of a coherent rigid space@--- (of a coherent rigid space)!distinguished formal model of a coherent rigid space@distinguished --- ---} $(j,\phi,\psi)$ of $\mathscr{U}\rightarrow\mathscr{X}$ such that the morphism $j\colon U\rightarrow X$ in $\Ac\CFs^{\ast}$ is an open immersion.
\end{itemize}
\end{prop}

\begin{proof}
We only need to show (a) $\Rightarrow$ (b).
Let $j\colon U\rightarrow X$ be an open immersion of coherent adic formal schemes of finite ideal type, and $\mathscr{I}$ an ideal of definition of finite type.
Let $X'\rightarrow X$ be the admissible blow-up along $\mathscr{I}$.
Then, as we saw in the proof of \ref{prop-cohrigidspacedist1}, $X'$ gives a distinguished formal model of $\mathscr{X}$.
Consider the induced open immersion $j_{X'}\colon U\times_XX'\rightarrow X'$.
By \ref{prop-blowups41} $U\times_XX'\rightarrow U$ coincides up to isomorphism with the admissible blow-up along $\mathscr{I}|_U$.
Hence $U\times_XX'$ is a distinguished formal model of $\mathscr{U}$, and the open immersion $j_{X'}$ gives a distinguished formal model of $\mathscr{U}\rightarrow\mathscr{X}$.
\end{proof}

\begin{dfn}\label{dfn-cohrigidspaceopenimm1}{\rm 
A morphism $\iota\colon\mathscr{U}\rightarrow\mathscr{X}$ of coherent rigid spaces is called a {\em coherent open immersion}\index{immersion!open immersion of rigid spaces@open --- (of rigid spaces)!coherent open immersion of rigid spaces@coherent --- ---} if it fulfills the equivalent conditions in \ref{prop-cohrigidspaceopenimm1}.}
\end{dfn} 

We will justify `coherent' in this terminology later in \ref{prop-consistencytermopenimm2}.
\begin{prop}\label{prop-openimmpullback1}
{\rm (1)} Let $\mathscr{U}\hookrightarrow\mathscr{V}$ and $\mathscr{V}\hookrightarrow\mathscr{X}$ be two coherent open immersions of coherent rigid spaces.
Then the composition $\mathscr{U}\rightarrow\mathscr{X}$ is a coherent open immersion.

{\rm (2)} Let $\mathscr{Y}\rightarrow\mathscr{X}$ be a morphism of coherent rigid spaces, and $\mathscr{U}\hookrightarrow\mathscr{X}$ a coherent open immersion.
Then the fiber product\index{fiber product!fiber product of rigid spaces@--- (of rigid spaces)}\index{rigid space!fiber product of rigid spaces@fiber product of ---s} $\mathscr{U}\times_{\mathscr{X}}\mathscr{Y}$ is representable in the category $\CRf$ and the morphism $\mathscr{U}\times_{\mathscr{X}}\mathscr{Y}\rightarrow\mathscr{Y}$ is a coherent open immersion. $($The general fiber products will be discussed later in {\rm \ref{prop-cohrigidspacefiberprod1}}.$)$
\end{prop}

\begin{proof}
(1) Take open immersions $U\hookrightarrow V$ and $V'\hookrightarrow X$ that provides respective formal models for $\mathscr{U}\hookrightarrow\mathscr{V}$ and $\mathscr{V}\hookrightarrow\mathscr{X}$.
Take admissible blow-ups $V''\rightarrow V$ and $V''\rightarrow V'$ (cf.\ \ref{prop-cohrigidspace3}).
By \ref{prop-blowups4111} one can take an admmissible blow-up $X'\rightarrow X$ such that $V''=V'\times_XX'$ and that $V''\rightarrow X'$ is an open immersion.
By {\bf \ref{ch-formal}}.\ref{prop-openimmform1} (3) $U'=U\times_VV''\rightarrow V''$ is an open immersion and gives a formal model of $\mathscr{U}\hookrightarrow\mathscr{V}$ due to \ref{prop-blowups41}.
The composition $U'\rightarrow X'$, which is an open immersion due to {\bf \ref{ch-formal}}.\ref{prop-openimmform1} (1), gives a formal model of $\mathscr{U}\rightarrow\mathscr{X}$.

(2) Let $U\stackrel{j}{\hookrightarrow} X\stackrel{f}{\leftarrow} Y$ be a diagram in $\Ac\CFs^{\ast}$, where $j$ is an open immersion, that gives rise to $\mathscr{U}\hookrightarrow\mathscr{X}\leftarrow\mathscr{Y}$ by passage to the associated coherent rigid spaces.
Take the fiber product $U\times_XY$ in the category $\Ac\CFs^{\ast}$.
Then one can easily check that the desired fiber product $\mathscr{U}\times_{\mathscr{X}}\mathscr{Y}$ is given by $(U\times_XY)^{\rig}$ and hence that the morphism $\mathscr{U}\times_{\mathscr{X}}\mathscr{Y}\rightarrow\mathscr{Y}$ is a coherent open immersion.
\end{proof}

\begin{prop}\label{prop-admissiblecovering}
Let $\mathscr{X}$ be a coherent rigid space, and $\{\mathscr{U}_{\alpha}\hookrightarrow\mathscr{X}\}_{\alpha\in L}$ a finite family of coherent open immersions of coherent rigid spaces.
Then there exists a $($distinguished$)$ formal model $X$ of $\mathscr{X}$ and a finite family $\{U_{\alpha}\hookrightarrow X\}_{\alpha\in L}$ of open immersions of coherent adic formal schemes that induces by passage to the associated coherent rigid spaces the given family $\{\mathscr{U}_{\alpha}\hookrightarrow\mathscr{X}\}_{\alpha\in L}$.
\end{prop}

\begin{proof}
By induction with respect to the cardinality of $L$, one reduce to the situation where $L=\{0,1\}$.
Take a formal model $U_0\hookrightarrow X$ (resp.\ $U_1\hookrightarrow X'$) of $\mathscr{U}_0\hookrightarrow\mathscr{X}$ (resp.\ $\mathscr{U}_1\hookrightarrow\mathscr{X}$).
There exist admissible blow-ups $X''\rightarrow X$ and $X''\rightarrow X'$ (cf.\ \ref{prop-cohrigidspace3}); moreover, these admissible blow-ups can be taken so that $X''$ is a distinguished formal model of $\mathscr{X}$ (\ref{prop-cohrigidspacedist1}).
Then $U'_0=U_0\times_XX''\rightarrow X''$ and $U'_1=U_1\times_{X'}X''\rightarrow X''$, which are open immersions by {\bf \ref{ch-formal}}.\ref{prop-openimmform1} (3)), give the desired formal models by \ref{prop-blowups41}.
\end{proof}

\begin{dfn}\label{dfn-admissiblecovering11}{\rm 
Let $\mathscr{X}$ be a coherent rigid space, and $\{\mathscr{U}_{\alpha}\hookrightarrow\mathscr{X}\}_{\alpha\in L}$ a finite family of coherent open immersions of coherent rigid spaces.
This family is said to be a {\em covering} of $\mathscr{X}$ if the following condition is satisfied:
there exists a formal model $X$ of $\mathscr{X}$ and a finite Zariski open covering $\{U_{\alpha}\hookrightarrow X\}_{\alpha\in L}$ that induces $\{\mathscr{U}_{\alpha}\hookrightarrow\mathscr{X}\}_{\alpha\in L}$ by passage to the associated coherent rigid spaces.}
\end{dfn}

By an argument similar to that in the proof of \ref{prop-admissiblecovering} and by {\bf \ref{ch-formal}}.\ref{prop-surjectivemorformal} (3), one easily sees the following: if $\{U_{\alpha}\hookrightarrow X\}_{\alpha\in L}$ is a finite Zariski covering as above, then for any admissible blow-up $X'$ one has the finite Zariski covering $\{U'_{\alpha}=U_{\alpha}\times_XX'\hookrightarrow X'\}_{\alpha\in L}$ that also induces $\{\mathscr{U}_{\alpha}\hookrightarrow\mathscr{X}\}_{\alpha\in L}$.
This shows, in particular, that the formal model $X$ in \ref{dfn-admissiblecovering11} can be taken to be distinguished.

\begin{dfn}[Coherent small admissible site for coherent rigid space]\label{dfn-admissiblecovering12}{\rm 
Let $\mathscr{X}$ be a coherent rigid space.
Define the site $\mathscr{X}_{\ad}$ as follows:
\begin{itemize}
\item[$\bullet$] objects of the category $\mathscr{X}_{\ad}$ are coherent open immersions $\mathscr{U}\hookrightarrow\mathscr{X}$ between coherent rigid spaces; 
\item[$\bullet$] arrows from $\mathscr{U}\hookrightarrow\mathscr{X}$ to $\mathscr{V}\hookrightarrow\mathscr{X}$ are morphisms of coherent rigid spaces over $\mathscr{X}$; 
\item[$\bullet$] for any object $\mathscr{U}\hookrightarrow\mathscr{X}$ the set of coverings $\Cov(\mathscr{U})$ consists of finite families of morphisms $\{\mathscr{U}_{\alpha}\hookrightarrow\mathscr{U}\}_{\alpha\in L}$ that gives a covering of the coherent rigid space $\mathscr{U}$ in the sense as in \ref{dfn-admissiblecovering11}.
\end{itemize}
The last-mentioned notion of coverings defines, by \ref{prop-openimmpullback1} and {\bf \ref{ch-formal}}.\ref{prop-surjectivemorformal}, a pretopology on the category $\mathscr{X}_{\ad}$.
The site $\mathscr{X}_{\ad}$ thus obtained is said to be the {\em coherent small admissible site} associated to the coherent rigid space $\mathscr{X}$.
We denote by $\mathscr{X}^{\sim}_{\ad}$ the topos induced from the site $\mathscr{X}_{\ad}$, called the {\em admissible topos}\index{topos!admissible topos@admissible ---} associated to $\mathscr{X}$.}
\end{dfn}

\begin{dfn}[Large admissible site of coherent rigid spaces]\label{dfn-admissiblecovering13}{\rm 
We endow the category $\CRf$ with the following topology: for any object $\mathscr{X}$ of $\CRf$ the set of coverings $\Cov(\mathscr{X})$ consists of finite families of coherent open immersions $\{\mathscr{U}_{\alpha}\hookrightarrow\mathscr{X}\}_{\alpha\in L}$ that gives a covering of $\mathscr{X}$ in the sense of \ref{dfn-admissiblecovering11}.
We denote the site thus obtained by $\CRf_{\ad}$, and the associated topos by $\CRf^{\sim}_{\ad}$.

For a coherent rigid space $\mathscr{S}$ the large admissible site $\CRf_{\mathscr{S},\ad}$, defined on the comma category $\CRf_{\mathscr{S}}$, and its associated topos $\CRf^{\sim}_{\mathscr{S},\ad}$ are defined similarly.}
\end{dfn}

\subsubsection{Properties of coherent admissible sites}\label{subsub-admissibletopproperty}
The following proposition is clear by the definition of the admissible topology:
\begin{prop}\label{prop-admissiblesitecoherent0}
Any object of the site $\mathscr{X}_{\ad}$ $($resp.\ $\CRf_{\ad}$, resp.\ $\CRf_{\mathscr{S},\ad})$ is quasi-compact as an object of $\mathscr{X}^{\sim}_{\ad}$ $($resp.\ $\CRf^{\sim}_{\ad}$, resp.\ $\CRf^{\sim}_{\mathscr{S},\ad})$.
In particular, the topos $\mathscr{X}^{\sim}_{\ad}$ $($resp.\ $\CRf^{\sim}_{\ad}$, resp.\ $\CRf^{\sim}_{\mathscr{S},\ad})$ has a generating full subcategory consisting of quasi-compact objects. \hfill$\square$
\end{prop}

\begin{prop}\label{prop-admissiblesitecoherent}
Let $\mathscr{X}$ be a coherent rigid space.
Then the topos $\mathscr{X}^{\sim}_{\ad}$ is coherent\index{coherent!coherent topos@--- topos}\index{topos!coherent topos@coherent ---} $({\bf \ref{ch-pre}}.\ref{dfn-coherenttopos2})$.
\end{prop}

\begin{proof}
By \ref{prop-admissiblesitecoherent0} and \ref{prop-openimmpullback1} (2) any object of $\mathscr{X}_{\ad}$ is coherent.
Since the final object $\mathscr{X}$ is a coherent objects, \ref{prop-openimmpullback1} (2) also implies that the category $\mathscr{X}_{\ad}$ is stable under finite projective limits.
\end{proof}

\begin{prop}\label{prop-admissiblesite3}
On the site $\CRf_{\ad}$ any representable presheaf is a sheaf.
\end{prop}

To prove this, and for later purpose, here we introduce the notion of {\it patching} of coherent rigid spaces:
\begin{prop}[(Birational) patching of coherent rigid spaces\index{rigid space!coherent rigid space@coherent ---!birational patching of coherent rigid spaces@birational patching of --- ---s}]\label{prop-cohrigidspacepatching1}

{\rm (1)} Let $\mathscr{X}$ be a coherent rigid space, and $\mathscr{U}_0\hookrightarrow\mathscr{X}$ and $\mathscr{U}_1\hookrightarrow\mathscr{X}$ coherent open immersions.
Then the Cartesian diagram 
$$
\xymatrix@-2ex{
&\mathscr{U}_0\ar@{^{(}->}[dr]\\
\mathscr{U}_{01}\ar@{^{(}->}[ur]\ar@{^{(}->}[dr]&&\mathscr{X}\rlap{,}\\
&\mathscr{U}_1\ar@{^{(}->}[ur]}
$$
where $\mathscr{U}_{01}=\mathscr{U}\times_{\mathscr{X}}\mathscr{U}_1$, is co-Cartesian in $\CRf$.

{\rm (2)} Consider the diagram 
$$
\mathscr{X}\stackrel{\alpha}{\longhookleftarrow}\mathscr{U}\stackrel{\beta}{\longhookrightarrow}\mathscr{Y}\leqno{(\ast)}
$$
in $\CRf$ where $\alpha$ and $\beta$ are coherent open immersions.
Then there exists a co-Cartesian and Cartesian square in $\CRf$
$$
\xymatrix@-2ex{
&\mathscr{X}\ar@{^{(}->}[dr]^{\varphi}\\
\mathscr{U}\ar@{^{(}->}[ur]^{\alpha}\ar@{^{(}->}[dr]_{\beta}&&\mathscr{Z}\rlap{,}\\
&\mathscr{Y}\ar@{^{(}->}[ur]_{\psi}}
$$
where $\varphi$ and $\psi$ are coherent open immersions.
\end{prop}

\begin{proof}
(1)  Let $\mathscr{W}$ be a coherent rigid space, and $\gamma_i\colon\mathscr{U}_i\rightarrow\mathscr{W}$ ($i=0,1$) morphisms of coherent rigid spaces.
Suppose $f_0=f_1$ on $\mathscr{U}_{01}$.
Take formal models $g_0\colon U_0\rightarrow W_0$ and $g_1\colon U_1\rightarrow W_1$ of $\gamma_0$ and $\gamma_1$, respectively.
One can replace $W_0$ and $W_1$ by a common admissible blow-up of them (due to \ref{prop-cohrigidspace3}), and $U_0$ and $U_1$ by the strict transforms, so that we may suppose $W_0=W_1$, which we denote by $W$.
We can moreover replace each $U_i$ ($i=0,1$) by an admissible blow-up, which is a coherent open subset of a formal model $X_i$ $(i=0,1)$ of $\mathscr{X}$; replacing $X_0$ and $X_1$ by a common admissible blow-ups, we may assume $X_0=X_1$, which we denote by $X$.
By the assumption, there exists an admissible blow-up $U'_{01}$ of $U_{01}=U_0\cap U_1$ on which $g_0$ and $g_1$ are equal. 
By \ref{prop-blowups4111}, there exists an admissible blow-up $X'\rightarrow X$ that extends the admissible blow-up $U'_{01}\rightarrow U_{01}$ such that $U'_{01}\cong U_{01}\times_XX'$.
Hence, replacing $X$ with $X'$ and $U_i$ ($i=0,1$) with the pull-back by the map $X'\rightarrow X$, we may assume that $g_0=g_1$ on $U_{01}$.
Then we have the unique morphism $g\colon X\rightarrow W$ such that $g|_{U_i}=g_i$ for $i=0,1$, and $\gamma=g^{\rig}\colon\mathscr{X}\rightarrow\mathscr{W}$ such that $\gamma|_{\mathscr{U}_i}=\gamma_i$ for $i=0,1$.
The uniqueness of $\gamma$ is straightforward.

(2) Take formal models $U\hookrightarrow X$ and $U'\hookrightarrow Y$ of $\alpha$ and $\beta$, respectively, which are open immersions.
Since $U^{\rig}\cong U^{\prime\rig}$, there exists admissible blow-ups $U''\rightarrow U$ and $U''\rightarrow U'$ (\ref{prop-cohrigidspace3}).
By \ref{prop-blowups4111} there exists an admissible blow-up $X'\rightarrow X$ of $X$ such that $U''=U\times_XX'$.
Hence one can replace the formal model $U\hookrightarrow X$ by $U''\hookrightarrow X'$.
Doing the same for $U'\hookrightarrow Y$, one sees that we may assume $U=U'$, that is, we may start with a diagram $X\hookleftarrow U\hookrightarrow Y$ of open immersions that induces the diagram $(\ast)$ by passage to the associated coherent rigid spaces.
Now consider the push-out $Z=X\amalg_UY$, the patching of $X$ and $Y$ along $U$, in $\Ac\CFs^{\ast}$, and set $\mathscr{Z}=Z^{\rig}$.
Thus we get the desired commutative diagram, which is Cartesian, for we have $U=X\times_ZY$.
The diagram thus obtained is co-Cartesian as well, due to (1).
\end{proof}

In the situation as in \ref{prop-cohrigidspacepatching1}, the coherent rigid space $\mathscr{Z}$ is denoted by $\mathscr{X}\amalg_{\mathscr{U}}\mathscr{Y}$ and called the coherent rigid space obtained by {\em patching of $\mathscr{X}$ and $\mathscr{Y}$ along $\mathscr{U}$}.
In case $\mathscr{U}$ is empty (\ref{dfn-emptyrigidspace}), we write $\mathscr{X}\amalg\mathscr{Y}$ (disjoint sum) for $\mathscr{X}\amalg_{\mathscr{U}}\mathscr{Y}$.
Since the square diagram in \ref{prop-cohrigidspacepatching1} is Cartesian, we have:
\begin{cor}\label{cor-disjointsumrigid}
The disjoint sum $\mathscr{X}\amalg\mathscr{Y}$ is universally disjoint {\rm (\cite[Expos\'e II, 4.5]{SGA4-1})}. \hfill$\square$
\end{cor}

\begin{cor}\label{cor-cohrigidspacepatching11}
In the category $\CRf$ any finite colimit consisting of coherent open immersions is representable. \hfill$\square$
\end{cor}

To prove \ref{prop-admissiblesite3}, we still need a few more lemmas:
\begin{lem}\label{lem-admissiblesite2}
Let $j\colon\mathscr{U}\hookrightarrow\mathscr{X}$ be a coherent open immersion of coherent rigid spaces.
If the singleton set $\{j\colon\mathscr{U}\hookrightarrow\mathscr{X}\}$ is a covering $($in the sense of $\ref{dfn-admissiblecovering11})$, then $j$ is an isomorphism. \hfill$\square$
\end{lem}

This is clear from the definition of coverings (\ref{dfn-admissiblecovering11}).
\begin{lem}\label{lem-cohrigidspacepatchingx}
Let $\varphi\colon\mathscr{X}\hookrightarrow\mathscr{Z}$ and $\psi\colon\mathscr{Y}\hookrightarrow\mathscr{Z}$ be coherent open immersions of coherent rigid spaces, and set $\mathscr{U}=\mathscr{X}\times_{\mathscr{Z}}\mathscr{Y}$.
Then the canonical morphism $\mathscr{X}\amalg_{\mathscr{U}}\mathscr{Y}\rightarrow\mathscr{Z}$ $($cf.\ {\rm \ref{prop-cohrigidspacepatching1}}$)$ is a coherent open immersion.
\end{lem}

\begin{proof}
Take a formal model $Z$ of $\mathscr{Z}$ and quasi-compact open immersions $X\hookrightarrow Z$ and $Y\hookrightarrow Z$ that induce $\varphi$ and $\psi$, respectively, by the passage to the associated rigid spaces. 
Set $U=X\times_ZY$.
Then we have $\mathscr{U}=U^{\rig}$.
As in the proof of \ref{prop-cohrigidspacepatching1}, the morphism $\mathscr{X}\amalg_{\mathscr{U}}\mathscr{Y}\rightarrow\mathscr{Z}$ of coherent rigid spaces is represented by the open immersion $X\amalg_UY\hookrightarrow Z$, whence the result.
\end{proof}

By \ref{lem-cohrigidspacepatchingx} and \ref{lem-admissiblesite2} we have:
\begin{lem}\label{lem-admissiblesite3}
Let $\mathscr{U}_0\hookrightarrow\mathscr{X}$ and $\mathscr{U}_1\hookrightarrow\mathscr{X}$ be coherent open immersions of coherent rigid spaces. 
Set $\mathscr{U}_{01}=\mathscr{U}_0\times_{\mathscr{X}}\mathscr{U}_1$.
Suppose that the induced morphism $\mathscr{U}_0\amalg\mathscr{U}_1\rightarrow\mathscr{X}$ is a covering of the site $\CRf_{\ad}$.
Then $\mathscr{U}_0\amalg_{\mathscr{U}_{01}}\mathscr{U}_1\rightarrow\mathscr{X}$ is an isomorphism. \hfill$\square$
\end{lem}

\begin{proof}[Proof of Proposition {\rm \ref{prop-admissiblesite3}}]
Let $\mathscr{X}$ and $\mathscr{U}$ be coherent rigid spaces, and suppose a covering $\coprod_{\alpha\in L}\mathscr{U}_{\alpha}\rightarrow\mathscr{U}$ is given.
We need to show that the sequence
$$\xymatrix@-1ex{
\Hom_{\CRf}(\mathscr{U},\mathscr{X})\ar[r]&\prod_{\alpha\in L}\Hom_{\CRf}(\mathscr{U}_{\alpha},\mathscr{X})\ar@<.5ex>[r]\ar@<-.5ex>[r]&\prod_{\alpha,\beta\in L}\Hom_{\CRf}(\mathscr{U}_{\alpha\beta},\mathscr{X})}
$$
is exact, where $\mathscr{U}_{\alpha\beta}=\mathscr{U}_{\alpha}\times_{\mathscr{U}}\mathscr{U}_{\beta}$ (cf.\ \ref{prop-openimmpullback1} (2)).
By induction with respect to the cardinality of the index set $I$, it suffices to show the asserion in the case $I=\{0,1\}$.
But this case follows promptly from \ref{lem-admissiblesite3} and the fact that the square diagram as in \ref{prop-cohrigidspacepatching1} (1) is co-Cartesian.
\end{proof}
\index{site!admissible site@admissible ---!coherent admissible site@coherent --- ---|)}
\index{admissible!admissible site@--- site!coherent admissible site@coherent --- ---|)}
\index{topology!admissible topology@admissible ---|)}
\index{admissible!admissible topology@--- topology|)}

\subsubsection{General rigid space}\label{subsub-generalrigidspace}
\index{rigid space!general rigid space@general ---|(}
\begin{dfn}\label{dfn-admissiblesite31}{\rm 
(1) A sheaf $\mathscr{F}$ of sets on $\CRf_{\ad}$ is said to be a {\em stretch of coherent rigid spaces} if there exists an inductive system $\{\mathscr{U}_i\}_{i\in J}$ of coherent rigid spaces indexed by a directed set such that:
\begin{itemize}
\item[(a)] for any $i,j\in J$ with $i\leq j$ the transitions map $\mathscr{U}_i\rightarrow\mathscr{U}_j$ is a coherent open immersion;
\item[(b)] for any coherent rigid space $\mathscr{X}$ we have
$$
\mathscr{F}(\mathscr{X})=\varinjlim_{i\in J}\Hom(\mathscr{X},\mathscr{U}_i)
$$
(that is, $\mathscr{F}$ is the inductive limit of the sheaves represented by $\mathscr{U}_i$).
\end{itemize}
In this situation, we also say that $\mathscr{F}$ is {\em represented by the stretch of the coherent rigid spaces} $\{\mathscr{U}_i\}_{i\in J}$.

(2) A morphism $\mathscr{F}\rightarrow\mathscr{G}$ of sheaves of sets on $\CRf_{\ad}$ is said to be (represented by) a {\em stretch of coherent open immersions} if for any representable sheaf $\mathscr{X}$ and any morphism $\mathscr{X}\rightarrow\mathscr{G}$ of sheaves, $\mathscr{F}\times_{\mathscr{G}}\mathscr{X}$ is represented by a stretch of coherent rigid spaces $\{\mathscr{U}_i\}_{i\in J}$ and the map $\mathscr{F}\times_{\mathscr{G}}\mathscr{X}\rightarrow\mathscr{X}$ coincides with the inductive limit of $\{\mathscr{U}_i\hookrightarrow\mathscr{X}\}_{i\in J}$.}
\end{dfn}

\begin{dfn}[General rigid space]\label{dfn-generalrigidspace1}{\rm 
A sheaf $\mathscr{F}$ of sets on $\CRf_{\ad}$ is called a {\em $($general$)$ rigid space} if it satisfies the following conditions:
\begin{itemize}
\item[(a)] there exists a surjective morphism of sheaves
$$
\mathscr{Y}=\coprod_{\alpha\in L}\mathscr{Y}_{\alpha}\longrightarrow\mathscr{F}
$$
(with the cardinality of the index set $L$ being arbitrary) where each $\mathscr{Y}_{\alpha}$ is represented by a coherent rigid space; 
\item[(b)] for any $\alpha,\beta\in L$ the projection 
$$
\mathscr{Y}_{\alpha}\times_{\mathscr{F}}\mathscr{Y}_{\beta}\longrightarrow\mathscr{Y}_{\alpha}
$$
is a stretch of coherent open immersions $(\ref{dfn-admissiblesite31}$ $(2))$.
\end{itemize}
A {\em morphism} $\varphi\colon\mathscr{F}\rightarrow\mathscr{G}$ of rigid spaces is, by definition, a morphism of sheaves.}
\end{dfn}

We denote by $\Rf$ the category of rigid spaces.
\begin{prop}\label{prop-generalrigidspace1}
Any sheaf representable by a coherent rigid space is a rigid space.
More generally, any stretch of coherent rigid spaces is a rigid space. \hfill$\square$
\end{prop}

The proof is straightforward. 
Obviously, we have the fully faithful functor
$$
\CRf\longhookrightarrow\Rf
$$
that maps a coherent rigid space $\mathscr{X}$ to the sheaf represented by $\mathscr{X}$.

\begin{dfn}[Open immersion]\label{dfn-generalrigidspace4}{\rm 
Let $\varphi\colon\mathscr{F}\rightarrow\mathscr{G}$ be a morphism of rigid spaces.
We say that $\varphi$ is an {\em open immersion}\index{immersion!open immersion of rigid spaces@open --- (of rigid spaces)} if it is a stretch of coherent open immersions $(\ref{dfn-admissiblesite31}$ $(2))$.
A $\mathscr{G}$-isomorphism class of open immersions is called an {\em open rigid subspace}\index{rigid subspace!open rigid subspace@open ---} of $\mathscr{G}$.}
\end{dfn}

By the definition of open immersions and \ref{prop-openimmpullback1} (2), we immediately deduce the following:
\begin{prop}\label{prop-generalrigidspace4x}
Let $\mathscr{U}\rightarrow\mathscr{X}$ be a coherent open immersion between coherent rigid spaces.
Then, regarded as a morphism between rigid spaces, it is an open immersion. \hfill$\square$
\end{prop}

Later in \ref{prop-consistencytermopenimm2} we will show that, conversely, any open immersion between coherent rigid spaces are coherent open immersions.
By \ref{prop-openimmpullback1} and \ref{lem-cohrigidspacepatchingx}, we have the following proposition.
\begin{prop}\label{prop-generalrigidspace4xx}
{\rm (1)} Let $\mathscr{U}\hookrightarrow\mathscr{V}$ and $\mathscr{V}\hookrightarrow\mathscr{F}$ be two open immersions of rigid spaces.
Then the composition $\mathscr{U}\rightarrow\mathscr{F}$ is an open immersion.

{\rm (2)} Let $\mathscr{G}\rightarrow\mathscr{F}$ be a morphism of rigid spaces, and $\mathscr{U}\hookrightarrow\mathscr{F}$ an open immersion.
Then the fiber product\index{fiber product!fiber product of rigid spaces@--- (of rigid spaces)}\index{rigid space!fiber product of rigid spaces@fiber product of ---s} $\mathscr{U}\times_{\mathscr{F}}\mathscr{G}$ is representable in the category $\Rf$ and the morphism $\mathscr{U}\times_{\mathscr{F}}\mathscr{G}\rightarrow\mathscr{G}$ is an open immersion. \hfill$\square$
\end{prop}

\subsubsection{Universally Noetherian and universally adhesive rigid spaces}\label{subsub-universallyadhesiverigidspaces}
\index{rigid space!universally adhesive rigid space@universally adhesive ---|(}
\index{rigid space!universally Noetherian rigid space@universally Noetherian ---|(}
\index{rigid space!universally adhesive rigid space@universally adhesive ---!locally universally adhesive rigid space@locally --- ---|(}
\index{rigid space!universally Noetherian rigid space@universally Noetherian ---!locally universally Noetherian rigid space@locally --- ---|(}
\begin{dfn}\label{dfn-universallyadhesiverigidspaces}{\rm 
A (general) rigid space $\mathscr{F}$ is called a {\em locally universally Noetherian} (resp.\ {\em locally universally adhesive}) rigid space if it has a covering as in \ref{dfn-generalrigidspace1} (a) such that each $\mathscr{Y}_{\alpha}$ $(\alpha\in L)$ is a coherent universally Noetherian (resp.\ universally adhesive) rigid space\index{rigid space!universally adhesive rigid space@universally adhesive ---!coherent universally adhesive rigid space@coherent --- ---}\index{rigid space!universally Noetherian rigid space@universally Noetherian ---!coherent universally Noetherian rigid space@coherent --- ---} (\ref{dfn-coherentuniversallyadhesiverigidspaces}).
If the covering can be chosen such that the index set $L$ is finite, then $\mathscr{F}$ is said to be a {\em universally Noetherian} (resp.\ {\em universally adhesive}) rigid space.}
\end{dfn}

According to our later terminology in \ref{dfn-generalrigidspace2} below, a locally universally Noetherian (resp.\ locally universally adhesive) rigid space is universally Noetherian (resp.\ universally adhesive) if it is quasi-compact.\index{rigid space!quasi-compact rigid space@quasi-compact ---}

Using some of the results in \S\ref{sec-embodying} below, one can show without vicious circle the following: {\it if $\mathscr{F}$ is represented by a coherent rigid space, then $\mathscr{F}$ is locally universally Noetherian $($resp.\ locally universally adhesive$)$ if and only if it is universally Noetherian $($resp.\ universally adhesive$)$ in the sense as in {\rm \ref{dfn-coherentuniversallyadhesiverigidspaces}}.}
Indeed, if $\mathscr{F}$, assumed to be representable by a coherent rigid space, is locally universally Noetherian (resp.\ locally universally adhesive) in the sense of \ref{dfn-universallyadhesiverigidspaces}, then by \ref{prop-generalrigidspace2} we may assume that the covering $\{\mathscr{Y}_{\alpha}\}_{\alpha\in L}$ is finite.
By \ref{prop-consistencytermopenimm2} we see that each map $\mathscr{Y}_{\alpha}\rightarrow\mathscr{F}$ is represented by a coherent open immersion.
Let $Y_{\alpha}$ be a universally rigid-Noetherian (resp.\ universally adhesive) formal model of $\mathscr{Y}_{\alpha}$ for each $\alpha$, and $X$ a formal model of the coherent rigid space that represents $\mathscr{F}$.
By \ref{prop-zariskiriemanntoptop}, replacing $X$ by an admissible blow-up if necessary, we may assume that $X$ has a Zariski covering $\{U_{\alpha}\}_{\alpha\in L}$ consisting of quasi-compact open subsets $U_{\alpha}$ together with an admissible blow-up $U_{\alpha}\rightarrow Y_{\alpha}$ for each $\alpha\in L$ (here we used \ref{prop-thmadmissiblesite200}).
Since each $U_{\alpha}$ is universally rigid-Noetherian (resp.\ universally adhesive), we deduce that $X$ is universally rigid-Noetherian (resp.\ universally adhesive).
\index{rigid space!universally Noetherian rigid space@universally Noetherian ---!locally universally Noetherian rigid space@locally --- ---|)}
\index{rigid space!universally adhesive rigid space@universally adhesive ---!locally universally adhesive rigid space@locally --- ---|)}
\index{rigid space!universally Noetherian rigid space@universally Noetherian ---|)}
\index{rigid space!universally adhesive rigid space@universally adhesive ---|)}
\index{rigid space!general rigid space@general ---|)}

\subsubsection{Admissible sites}\label{subsub-admissibletopcont}
\index{admissible!admissible site@--- site|(}
\index{site!admissible site@admissible ---|(}
\begin{dfn}[Small admissible site]\label{dfn-admissiblesite3gensmall}{\rm 
For a rigid space $\mathscr{F}$ we denote by $\mathscr{F}_{\ad}$ the site defined as follows. 
As a category, it is the category of all open immersions $\mathscr{U}\rightarrow\mathscr{F}$ and morphisms over $\mathscr{F}$.
For an object $\mathscr{U}\rightarrow\mathscr{F}$ the collection of coverings $\Cov(\mathscr{U})$ consists of familes $\{\mathscr{U}_{\alpha}\rightarrow\mathscr{U}\}_{\alpha\in L}$ (indexed by an arbitrary set) of morphisms such that $\coprod\mathscr{U}_{\alpha}\rightarrow\mathscr{U}$ is an epimorphism of sheaves on $\CRf_{\ad}$.
Note that this defines a pretopology due to \ref{prop-openimmpullback1}.
We denote by $\mathscr{F}^{\sim}_{\ad}$ the topos associated to the site $\mathscr{F}_{\ad}$.}
\end{dfn}

\danger{Notice that, if $\mathscr{X}$ is a coherent rigid space and $\mathscr{F}$ is the general rigid space represented by $\mathscr{X}$ $($that is, the image of $\mathscr{X}$ by the functor $\CRf\hookrightarrow\Rf)$, then the site $\mathscr{X}_{\ad}$ defined in \ref{dfn-admissiblecovering12} and the site $\mathscr{F}_{\ad}$ defined in \ref{dfn-admissiblesite3gensmall} are different.
We will see in \ref{thm-generalrigidspace31-32}, however, there exists a canonical morphism of sites $\mathscr{X}_{\ad}\rightarrow\mathscr{F}_{\ad}$ inducing an equivalence of topoi $\mathscr{X}^{\sim}_{\ad}\stackrel{\sim}{\rightarrow}\mathscr{F}^{\sim}_{\ad}$.}

\begin{dfn}[Large admissible site]\label{dfn-admissiblesite3genlarge}{\rm 
We endow the category $\Rf$ with the following topology: for any object $\mathscr{F}$ of $\Rf$ the collection of coverings $\Cov(\mathscr{F})$ consists of families of open immersions $\{\mathscr{U}_{\alpha}\hookrightarrow\mathscr{F}\}_{\alpha\in L}$ (indexed by an arbitrary set) such that the map of sheaves $\coprod\mathscr{U}_{\alpha}\rightarrow\mathscr{F}$ on $\CRf_{\ad}$ is an epimorphism.
Thanks to \ref{prop-generalrigidspace4xx}, this gives a pretopology on the category $\Rf$.
We denote this site by $\Rf_{\ad}$ and the associated to topos by $\Rf^{\sim}_{\ad}$.

For a rigid space $\mathscr{G}$ the large admissible site $\Rf_{\mathscr{G},\ad}$, defined on the obvious comma category $\Rf_{\mathscr{G}}$, and the associated topos $\Rf^{\sim}_{\mathscr{G},\ad}$ are defined similarly.}
\end{dfn}
\index{site!admissible site@admissible ---|)}
\index{admissible!admissible site@--- site|)}
\index{rigid space|)}

\subsection{Morphism of finite type}\label{sub-rigidmorphismfintype}
\begin{dfn}\label{dfn-cohrigidspacefintype}{\rm 
(1) A morphism $\varphi\colon\mathscr{X}\rightarrow\mathscr{Y}$ of coherent rigid spaces\index{rigid space!coherent rigid space@coherent ---} is said to be {\em of finite type} if it has a formal model\index{formal model!formal model of a coherent rigid space@--- (of a coherent rigid space)} $f\colon X\rightarrow Y$ of finite type\index{morphism of formal schemes@morphism (of formal schemes)!morphism of formal schemes of finite type@--- of finite type} ({\bf \ref{ch-formal}}.\ref{dfn-topfintype}).

(2) Let $\varphi\colon\mathscr{X}\rightarrow\mathscr{Y}$ be a morphism of rigid spaces where $\mathscr{Y}$ is a coherent rigid space. Then $\varphi$ is said to be {\em locally of finite type} if there exists a covering family $\{\mathscr{U}_{\alpha}\hookrightarrow\mathscr{X}\}_{\alpha\in L}$ of the site $\mathscr{X}_{\ad}$, where each $\mathscr{U}_{\alpha}$ is a coherent rigid space, such that for any $\alpha\in L$ the morphism $\mathscr{U}_{\alpha}\rightarrow\mathscr{Y}$ is of finite type in the sense of (1).
If, moreover, such a covering family can be taken to be finite, we say that $\varphi$ is {\em of finite type}.

(3) For a morphism of general rigid spaces $\varphi\colon\mathscr{X}\rightarrow\mathscr{Y}$, we say $\varphi$ is {\em locally of finite type}\index{morphism of rigid spaces@morphism (of rigid spaces)!morphism of rigid spaces locally of finite type@--- locally of finite type} (resp.\ {\em of finite type}\index{morphism of rigid spaces@morphism (of rigid spaces)!morphism of rigid spaces of finite type@--- of finite type}) if there exists a covering family $\{\mathscr{V}_{\alpha}\hookrightarrow\mathscr{Y}\}_{\alpha\in L}$ of the site $\mathscr{Y}_{\ad}$, where each $\mathscr{V}_{\alpha}$ is a coherent rigid space, such that for any $\alpha\in L$ the induced morphism $\mathscr{X}\times_{\mathscr{Y}}\mathscr{V}_{\alpha}\rightarrow\mathscr{V}_{\alpha}$ is locally of finite type (resp.\ of finite type) in the sense of (2).
(Note that $\mathscr{X}\times_{\mathscr{Y}}\mathscr{V}_{\alpha}$ is represented by a rigid space.)}
\end{dfn}

One can show that these definitions are consistent with each other; cf.\ \ref{prop-generalrigidspace2}.

\begin{prop}\label{prop-cohrigidspacefintype0a}
{\rm (1)} An open immersion\index{immersion!open immersion of rigid spaces@open --- (of rigid spaces)} is locally of finite type.

{\rm (2)} The composition of two morphisms locally of finite type $($resp.\ of finite type$)$ is again locally of finite type $($resp.\ of finite type$)$.
If the composition $\psi\circ\varphi$ of morphisms $\varphi\colon\mathscr{X}\rightarrow\mathscr{Y}$ and $\psi\colon\mathscr{Y}\rightarrow\mathscr{Z}$ is locally of finite type, then $\varphi$ is locally of finite type.
\end{prop}

\begin{proof}
(1) is clear.
We show (2). 
Let $\varphi\colon\mathscr{X}\rightarrow\mathscr{Y}$ and $\psi\colon\mathscr{Y}\rightarrow\mathscr{Z}$ be morphisms of finite type.
We may assume without loss of generality that $\mathscr{X}$, $\mathscr{Y}$, and $\mathscr{Z}$ are coherent rigid spaces.
Take respective formal models $f\colon X\rightarrow Y$ and $g\colon Y'\rightarrow Z$.
Since $Y^{\rig}\cong Y^{\prime\rig}$, by \ref{prop-cohrigidspace3} we have the following diagram
$$
\xymatrix@-2ex{
&&Y''\ar[dl]_{\pi}\ar[dr]^{\pi'}\\
X\ar[r]&Y&&Y'\ar[r]&Z\rlap{,}}
$$
where the two oblique arrows are admissible blow-ups.
Consider $f'\colon X'=X\times_YY''\rightarrow Y''$.
By \ref{cor-cohrigidspace32} the morphism $f'$ gives another formal model of $\varphi$, and the composition $g\circ\pi'\circ f'$, which is of finite type ({\bf \ref{ch-formal}}.\ref{prop-topfintype2} (2), (4)), gives a formal model of $\psi\circ\varphi$.

To show the other assertion, take a formal model $h\colon X''\rightarrow Z'$ of finite type of the composition $\psi\circ\varphi$.
Since there exist admissible blow-ups $Z''\rightarrow Z$ and $Z''\rightarrow Z'$, one can assume by the base change to $Z''$ that $Z'=Z$.
Replacing $X''$ by a suitable admissible blow-up, we may assume that there exists a map $\pi''\colon X''\rightarrow X'$ such that $h=g\circ\pi'\circ f'\circ\pi''$.
Since $h$ is of finite type, we deduce that $f'\circ\pi''$, which is a formal model of $\varphi$, is of finite type ({\bf \ref{ch-formal}}.\ref{prop-topfintype2} (2)).
\end{proof}

\begin{prop}\label{prop-finitetypeuniversallyadhesive}
Let $\varphi\colon\mathscr{X}\rightarrow\mathscr{Y}$ be a locally of finite type morphism between rigid spaces.
If $\mathscr{Y}$ is locally universally Noetherian $($resp.\ locally universally adhesive$)$, then so is $\mathscr{X}$.
\end{prop}

\begin{proof}
We may assume without loss of generality that $\mathscr{X}$ and $\mathscr{Y}$ are coherent and that $\varphi$ has a formal model $f\colon X\rightarrow Y$ of finite type; here, since admissible blow-ups are of finite type, we may furthermore assume that $Y$ is universally rigid-Noetherian (resp.\ universally adhesive).
Then $X$ is universally rigid-Noetherian (resp.\ universally adhesive) ({\bf \ref{ch-formal}}.\ref{prop-adequateformalscheme1}).
\end{proof}

\begin{prop}\label{prop-cohrigidspacefintype1}
Let $\varphi\colon\mathscr{X}\rightarrow\mathscr{Y}$ be a finite type morphism of coherent rigid spaces. 
Then there exists a distinguished formal model\index{formal model!formal model of a coherent rigid space@--- (of a coherent rigid space)!distinguished formal model of a coherent rigid space@distinguished --- ---} $f\colon X\rightarrow Y$ of $\varphi$ of finite type.
Moreover, such formal models are cofinal in the category $\FM_{\varphi}$.
If $\mathscr{Y}$ is universally adhesive, then such formal models are always of finite presentation.
\end{prop}

\begin{proof}
The first and the second assertions follow from \ref{prop-admissibleblowupmorphism} applied to the admissible blow-up along an ideal of definition of finite type (cf.\ \ref{prop-cohrigidspacedist1}).
Suppose $\mathscr{Y}$ is universally adhesive (and hence so is $\mathscr{X}$ due to \ref{prop-finitetypeuniversallyadhesive}).
If there exists a distinguished formal model as above of finite type consisting of universally adhesive formal schemes, it is automatically of finite presentation due to {\bf \ref{ch-formal}}.\ref{prop-fintypetorfreefinpres}.
\end{proof}

\subsection{Fiber products of rigid spaces}\label{sub-rigidspacefiberprod}
\index{fiber product!fiber product of rigid spaces@--- (of rigid spaces)|(}
\index{rigid space!fiber product of rigid spaces@fiber product of ---s|(}
\begin{prop}\label{prop-cohrigidspacefiberprod1}
For any diagram in $\Rf$ of the form $\mathscr{X}\stackrel{\varphi}{\rightarrow}\mathscr{S}\stackrel{\psi}{\leftarrow}\mathscr{Y}$ the fiber product $\mathscr{X}\times_{\mathscr{S}}\mathscr{Y}$ is representable in $\Rf$.
\end{prop}

\begin{proof}
We first deal with the case where the rigid spaces $\mathscr{X}$, $\mathscr{S}$, and $\mathscr{Y}$ are coherent.
Take a diagram $X\stackrel{f}{\rightarrow}S\stackrel{g}{\leftarrow}Y$ in $\Ac\CFs^{\ast}$ such that $f^{\rig}=\varphi$ and $g^{\rig}=\psi$.
Then by {\bf \ref{ch-formal}}.\ref{cor-fiverprodformaladicadic} there exists a fiber product $X\times_SY$ in $\Ac\CFs^{\ast}$.
The desired coherent rigid space $\mathscr{X}\times_{\mathscr{S}}\mathscr{Y}$ is then defined to be $(X\times_SY)^{\rig}$ together with the projections induced from the projections $X\times_SY\rightarrow X$ and $X\times_SY\rightarrow Y$.
By \ref{cor-cohrigidspace32} one sees easily that this gives a well-defined coherent rigid space, and $\mathscr{X}\times_{\mathscr{S}}\mathscr{Y}$ thus obtained indeed gives a fiber product of the diagram in question.

In general, first construct the fiber product $\mathscr{X}\times_{\mathscr{S}}\mathscr{Y}$ in the category of sheaves on the site $\CRf_{\ad}$.
Then it is easy to see by the standard argument, with the aid of \ref{prop-openimmpullback1} (2), that, by what we have just seen for coherent rigid spaces, the last sheaf is a rigid space.
\end{proof}

\begin{cor}\label{cor-cohrigidspacefiberprod1}
The quotient functor $Q\colon\Ac\CFs^{\ast}\rightarrow\CRf$ preserves fiber products. \hfill$\square$
\end{cor}

\begin{prop}\label{prop-cohrigidspacefiberprod12}
{\rm (1)} If $\varphi\colon\mathscr{X}\rightarrow\mathscr{X}'$ and $\psi\colon\mathscr{Y}\rightarrow\mathscr{Y}'$ are two open immersions $($resp.\ morphisms locally of finite type, resp.\ morphisms of finite type$)$ over a rigid space $\mathscr{S}$, then the induced morphism $\varphi\times_{\mathscr{S}}\psi\colon\mathscr{X}\times_{\mathscr{S}}\mathscr{Y}\rightarrow\mathscr{X}'\times_{\mathscr{S}}\mathscr{Y}'$ is an open immersion $($resp.\ a morphism locally of finite type, resp.\ a morphism of finite type$)$.

{\rm (2)} If $\varphi\colon\mathscr{X}\rightarrow\mathscr{Y}$ is an open immersion $($resp.\ a morphism locally of finite type, resp.\ a morphism of finite type$)$ over a rigid space $\mathscr{S}$ and $\mathscr{S}'\rightarrow\mathscr{S}$ is a morphism of rigid spaces, then the induced morphism $\varphi_{\mathscr{S}'}\colon\mathscr{X}\times_{\mathscr{S}}\mathscr{S}'\rightarrow\mathscr{Y}\times_{\mathscr{S}}\mathscr{S}'$ is an open immersion $($resp.\ a morphism locally of finite type, resp.\ a morphism of finite type$)$.
\end{prop}

\begin{proof}
By {\bf \ref{ch-pre}}.\ref{prop-basechangestable} we only need to show the following: for a diagram $\mathscr{X}\stackrel{\varphi}{\rightarrow}\mathscr{Y}\stackrel{\psi}{\leftarrow}\mathscr{Z}$ of rigid spaces where $\varphi$ is an open immersion $($resp.\ a morphism locally of finite type, resp.\ a morphism of finite type$)$, the induced morphism $\varphi_{\mathscr{Z}}\colon \mathscr{X}\times_{\mathscr{Y}}\mathscr{Z}\rightarrow\mathscr{Z}$ is an open immersion $($resp.\ a morphism locally of finite type, resp.\ a morphism of finite type$)$.
We may assume that $\mathscr{X}$, $\mathscr{Z}$, and $\mathscr{Y}$ are coherent rigid spaces and that $\varphi$ is a coherent open immersion (resp.\ a morphism of finite type, resp.\ a morphism of finite type).
Then the claim follows from {\bf \ref{ch-formal}}.\ref{prop-openimmform1} (2) and {\bf \ref{ch-formal}}.\ref{prop-topfintype2} (3).
\end{proof}
\index{rigid space!fiber product of rigid spaces@fiber product of ---s|)}
\index{fiber product!fiber product of rigid spaces@--- (of rigid spaces)|)}

\subsection{Examples of rigid spaces}\label{sub-examplesrigidspaces}
\subsubsection{Rigid spaces of type (V)}\label{subsub-examplesV}
\index{rigid space!rigid space of typeV@--- of type (V)|(}
\begin{dfn}\label{dfn-typeV}{\rm 
A {\em rigid space of type} (V) is a rigid space $\mathscr{X}$ locally of finite type over $(\Spf V)^{\rig}$ where $V$ is an $a$-adically complete valuation ring\index{valuation!valuation ring@--- ring!a-adically complete valuation ring@$a$-adically complete --- ---} with $a\in\m_V\setminus\{0\}$ (cf.\ {\bf \ref{ch-pre}}, \S\ref{sub-aadicallycompletevalrings}).
If, in particular, the valuation ring $V$ is of height one, then $\mathscr{X}$ is called a rigid space {\em of type {\rm ($\mathrm{V_{\R}}$)}}\index{rigid space!rigid space of typeV1@--- of type ($\mathrm{V_{\R}}$)}.}
\end{dfn}

Since the topological ring $V$ is t.u.\ adhesive \index{t.u.a. ring@t.u.\ adhesive ring} ({\bf \ref{ch-pre}}.\ref{cor-convadh}), the coherent rigid space $(\Spf V)^{\rig}$ is universally adhesive, and hence by \ref{prop-finitetypeuniversallyadhesive} any rigid space of type (V) is locally universally adhesive\index{rigid space!universally adhesive rigid space@universally adhesive ---!locally universally adhesive rigid space@locally --- ---} (\ref{dfn-universallyadhesiverigidspaces}).
\index{rigid space!rigid space of typeV@--- of type (V)|)}

\subsubsection{Rigid spaces of type (N)}\label{subsub-examplesN}
\index{rigid space!rigid space of typeN@--- of type (N)|(}
\begin{dfn}\label{dfn-typeN}{\rm 
A {\em rigid space of type} (N) is a rigid space $\mathscr{X}$ that admits an open covering $\{\mathscr{U}_{\alpha}\}_{\alpha\in L}$ consisting of coherent rigid spaces having Noetherian formal models.}
\end{dfn}

Since any Noetherian adic ring is t.u.\ adhesive, it follows that any rigid space of type (N) is locally universally adhesive\index{rigid space!universally adhesive rigid space@universally adhesive ---!locally universally adhesive rigid space@locally --- ---} (\ref{dfn-universallyadhesiverigidspaces}).
It can be shown, with the aid of \ref{prop-zariskiriemanntoptop} below, that, if a coherent rigid space $\mathscr{X}$ is of type (N), then $\mathscr{X}$ has a Noetherian formal model (Exercise \ref{exer-proptypeN}).
\index{rigid space!rigid space of typeN@--- of type (N)|)}

\subsubsection{Unit disk over a rigid space}\label{subsub-unitdisk}
\index{unit disk|(}
Let $\mathscr{X}$ be a rigid space.
The so-called {\em $($closed$)$ unit disk} over $\mathscr{X}$, denoted by 
$$
\D^n_{\mathscr{X}},
$$
is the rigid space defined as follows.
When $\mathscr{X}$ is of the form $\mathscr{X}=(\Spf A)^{\rig}$ where $A$ is an adic ring of finite ideal type, then $\D^n_{\mathscr{X}}=(\Spf A\dl X_1,\ldots,X_n\dr)^{\rig}$ (cf.\ {\bf \ref{ch-pre}}, \S\ref{sub-powerseries}) or, what amounts to the same, the coherent rigid space associated to the formal affine $n$-space\index{affine!affine space in formal geometry@--- space (in formal geometry)} $\widehat{\A}^n_A$ over $A$ (Exercise \ref{exer-affinespaceformal} (1)).
Since any rigid space has an admissible covering by open subspaces of the above form, one can define $\D^n_{\mathscr{X}}$ for general $\mathscr{X}$ by patching.
\index{unit disk|)}

\subsubsection{Projective space over a rigid space}\label{subsub-examplesproj}
\index{projective!projective space in rigid geometry@--- space (in rigid geometry)}
The projective space
$$
\P^{n,\an}_{\mathscr{X}}
$$
over a rigid space $\mathscr{X}$ is constructed as follows.
When $\mathscr{X}$ is of the form $\mathscr{X}=(\Spf A)^{\rig}$, then consider the formal projective space $\widehat{\P}^n_X=\widehat{\P}^n_A$ (Exercise \ref{exer-affinespaceformal} (2)); then set $\P^{n,\an}_{\mathscr{X}}=(\widehat{\P}^n_X)^{\rig}$.
For a general rigid space $\mathscr{X}$ one can define $\P^{n,\an}_{\mathscr{X}}$ by gluing.

\addcontentsline{toc}{subsection}{Exercises}
\subsection*{Exercises}
\begin{exer}\label{exer-propgeneralrigidspace0}
{\rm Show that in the category $\Rf$ of rigid spaces any colimit consisting of open immersions is representable.}
\end{exer}

\begin{exer}\label{exer-adequateformaltorigid}{\rm 
Consider the category $\Ac\Fs^{\ast}$ of all adic formal schemes of finite ideal type with adic morphisms.
Then define the natural functor 
$$
\Ac\Fs^{\ast}\longrightarrow\Rf
$$
extending $Q\colon\Ac\CFs^{\ast}\rightarrow\CRf$ ($X\mapsto X^{\rig}$).}
\end{exer}


\section{Visualization}\label{sec-embodying}
\index{visualization|(}
While the introduction of rigid spaces in the previous section is, so to speak, one of the two most fundamental starting points of rigid geometry, the other one is the so-called {\em visualization}, which we will discuss in this section.
Our definition of rigid spaces as `generic fibers' of formal schemes, which we have given in the previous section, is based on the creed that rigid geometry is like a `birational geometry' of formal schemes. 
It being so, one can say that the visualization of rigid spaces is the way to enhance the birational geometric aspect of the rigid geometry.
It does this job by adopting Zariski's old idea of birational geometry, by means of the so-called {\em Zariski-Riemann spaces}, which we will introduce in \S\ref{sub-ZRdef}.
There, the Zariski-Riemann spaces are first constructed for coherent rigid spaces.
Coherent rigid spaces have formal models by definition, and the Zariski-Riemann space associated to a given coherent rigid space is defined to be the projective limit of all formal models.
As the admissible blow-ups of an arbitrary formal model comprise a cofinal part of the totality of all formal models of the coherent rigid space, the Zariski-Riemann space is equivalently defined as the filtered projective limit of admissible blow-ups of a fixed formal model.

One of the most significant topological features of the Zariski-Riemann spaces comes promptly after the definition in Theorem \ref{thm-ZRcompact}, which asserts that the Zariski-Riemann space associated to a coherent rigid space is in fact a coherent sober topological space ({\bf \ref{ch-pre}}.\ref{dfn-quasicompact1}).
As stated in the introduction of {\bf \ref{ch-pre}}, \S\ref{sec-gentop}, the coherence or, especially, the quasi-compactness of the Zariski-Riemann spaces plays a very important role in our theory of rigid geometry, similarly to the fact that the quasi-compactness of Zariski's generalized Riemann spaces, proved by Zariski in 1944, played one of the most important roles in his works of resolution of singularities of algebraic varieties.
It is also worth while remarking that introducing the Zariski-Riemann spaces as above allows one to have more lucid picture of what we have called the `birational patching' in the previous section (cf.\ \ref{prop-cohrigidspacepatching1}).
Indeed, the `birational patching' gives the usual topological patching of the Zariski-Riemann spaces along quasi-compact open subsets.
By this, in particular, the general construction of Zariski-Riemann spaces associated to a general rigid space can be given quite naturally.

In the next subsection \S\ref{sub-ZRstrsheaf}, we discuss the structure sheaves on the Zariski-Riemann spaces.
As the Zariski-Riemann spaces are locally given by the projective limits of formal schemes, they have the natural sheaves of rings, simply given by the inductive limits of the structure sheaves of the formal schemes.
This sheaf, called the {\em integral structure sheaf} and denoted by $\O^{\int}_{\mathscr{X}}$, is needless to say an important object to consider.
But, in view of the central tenet of rigid geometry saying that rigid geometry is the geometry of `generic fibers' of formal schemes, one has to `invert' the ideal of definition in the sheaf $\O^{\int}_{\mathscr{X}}$ to obtain the `correct' structure sheaf of the rigid spaces.
This `inversion' of ideals of definition is possible by the fact that the stalk at each point of the integral structure sheaf $\O^{\int}_{\mathscr{X}}$ is a valuative ring ({\bf \ref{ch-pre}}.\ref{dfn-valuative1}) with respect to the ideal of definition (\ref{prop-ZRstrsheaf21}).
The resulting sheaf of rings, denoted by $\O_{\mathscr{X}}$, is the one we take up as the structure sheaf of the rigid space $\mathscr{X}$, called the {\em rigid structure sheaf} or just {\em structure sheaf} of $\mathscr{X}$.
Since the integral structure sheaf is, from the viewpoint of rigid geometry, regarded as the {\em canonical} formal model of the structure sheaf, it retains the importance as a natural object associated to the rigid space.
For this reason, when considering visualization of rigid spaces, one should consider the triple $\ZRT(\mathscr{X})=(\ZR{\mathscr{X}},\O^{\int}_{\mathscr{X}},\O_{\mathscr{X}})$, the so-called {\em Zariski-Riemann triple}\index{Zariski-Riemann triple}\index{triple!Zariski Riemann triple@Zariski-Riemann ---}, consisting of the Zariski-Riemann space $\ZR{\mathscr{X}}$, the integral structure sheaf, and the rigid structure sheaf, rather than the locally ringed space $(\ZR{\mathscr{X}},\O_{\mathscr{X}})$ alone.

In \S\ref{sub-ZRpoints} we will study the points of the Zariski-Riemann spaces.
There, we will see that the notion of points coincides with the one usually called {\em rigid points}\index{point!rigid point@rigid ---}\index{rigid point} in classical rigid geometry.
This is the place where valuation rings come to play in rigid geometry; remember that, as stated in the introduction of {\bf \ref{ch-pre}}, \S\ref{sec-val}, valuation rings are the most natural `value rings' at points in rigid spaces, and thus take the role of fields in algebraic geometry.

The subsection \S\ref{sub-comparisontopoi} is devoted to the comparison between admissible topology and the topology (in the usual sense) of the Zariski-Riemann spaces.
In the final subsection \S\ref{sub-consistencyterm} we consider quasi-compactness and quasi-separatedness and show that these notions for rigid spaces and those for corresponding Zariski-Riemann spaces coincide.

\subsection{Zariski-Riemann spaces}\label{sub-ZRdef}
\index{Zariski-Riemann space|(}
\subsubsection{Construction in coherent case}\label{subsub-ZRdef}
Let $\mathscr{X}$ be a coherent rigid space\index{rigid space!coherent rigid space@coherent ---}, and take a formal model\index{formal model!formal model of a coherent rigid space@--- (of a coherent rigid space)} $(X,\phi)$ of $\mathscr{X}$.
Consider the functor $S_X\colon\BL_X\rightarrow\LRsp$ from the category of admissible blow-ups\index{blow-up!admissible blow-up@admissible ---} of $X$ (\S\ref{sub-categoryadmblow-up}) to the category of locally ringed spaces\index{ringed space!locally ringed space@locally --- ---} that maps each admissible blow-up $X'\rightarrow X$ to the underlying locally ringed space of $X'$.
Consider the corresponding limit
$$
\varprojlim S_X\leqno{(\ast)}
$$
(cf.\ {\bf \ref{ch-pre}}, \S\ref{sub-limdef} for the notation).
Since the functor $\AId_X\rightarrow\BL_X$ (cf.\ \ref{prop-catblowups1} (2)) is cofinal, one can replace the above limit by the projective limit ({\bf \ref{ch-pre}}.\ref{prop-cofinal})
$$
\varprojlim_{\mathscr{J}\in\AId_X}X_{\mathscr{J}};
$$
note that such a limit exists in the category $\LRsp$ due to {\bf \ref{ch-pre}}.\ref{prop-limLRS}.
Notice also that the limit $(\ast)$ does not depend on the choice of the formal model $X$; indeed, for example, if we change $X$ by an admissible blow-up $X'$ of $X$, then we have a canonical functor $\BL_{X'}\rightarrow\BL_X$ (simply by composition; cf.\ \ref{prop-blowups3}), which is clearly cofinal, and hence the limits taken along $\BL_X$ and $\BL_{X'}$ coincide up to canonical isomorphisms in the category $\LRsp$.
Hence the limit $(\ast)$ only depends on $\mathscr{X}$ up to canonical isomorphisms.

\begin{dfn}\label{dfn-zariskiriemanndfn}{\rm 
The underlying topological space of the limit $(\ast)$ is called the {\em Zariski-Riemann space associated to $\mathscr{X}$} and denoted by $\ZR{\mathscr{X}}$.}
\end{dfn}

Let $\mathscr{X}$ be a coherent rigid space, and $\ZR{\mathscr{X}}$ the associated Zariski-Riemann space.
Let $X$ be any formal model of $\mathscr{X}$.
Since $\ZR{\mathscr{X}}$ is the projective limit of all admissible blow-ups of $X$, we have the projection map $\ZR{\mathscr{X}}\rightarrow X$ of locally ringed spaces, called the {\em specialization map}\index{specialization map} and denoted by
$$
\sp_X\colon\ZR{\mathscr{X}}\longrightarrow X.
$$

Clearly, we have $\ZR{\emptyset}=\emptyset$ (cf.\ \ref{dfn-emptyrigidspace} for the definition of empty rigid space $\emptyset$).
More generally, for a {\em scheme} $X$ (that is, a $0$-adic formal scheme; cf.\ {\bf \ref{ch-formal}}.\ref{rem-clear}) we have $\ZR{X^{\rig}}=\emptyset$.

\subsubsection{Functoriality}\label{subsub-ZRdeffunct}
Let $\varphi\colon\mathscr{X}\rightarrow\mathscr{Y}$ be a morphism in $\CRf$.
Then one can take a formal model $f\colon X\rightarrow Y$ of $\varphi$, and by \ref{prop-admissibleblowupmorphism} we get a canonical continuous map
$$
\ZR{\varphi}\colon\ZR{\mathscr{X}}\longrightarrow\ZR{\mathscr{Y}}
$$
of topological spaces.
Thus the formation of the associated Zariski-Riemann spaces gives rise to a functor
$$
\ZR{\,\cdot\,}\colon\CRf\longrightarrow\Top,\leqno{(\dagger)}
$$
where $\Top$ denotes the category of topological spaces.

\subsubsection{Topological feature}\label{subsub-ZRdeftopological}
\begin{thm}\label{thm-ZRcompact}
Let $\mathscr{X}$ be a coherent rigid space, and $\ZR{\mathscr{X}}$ the associated Zariski-Riemann space.
Then$:$
\begin{itemize}
\item[{\rm (1)}] The topological space $\ZR{\mathscr{X}}$ is coherent\index{space@space (topological)!coherent topological space@coherent ---}\index{coherent!coherent topological space@--- (topological) space} {\rm ({\bf \ref{ch-pre}}.\ref{dfn-coherenttopos2})} and sober {\rm ({\bf \ref{ch-pre}}, \S\ref{subsub-sober})}\index{space@space (topological)!sober topological space@sober ---}.
\item[{\rm (2)}] For any formal model $X$ of $\mathscr{X}$, the specialization map $\sp_X$ is quasi-compact\index{map@map (continuous)!quasi-compact map@quasi-compact ---}\index{quasi-compact!quasi-compact map@--- map} {\rm ({\bf \ref{ch-pre}}.\ref{dfn-quasicompactness}\ (2))} and closed.
\end{itemize}
\end{thm}

\begin{proof}
Since the topological space $\ZR{\mathscr{X}}$ is, by definition, the projective limit of coherent sober topological spaces with all the transition maps quasi-compact, (1) follows from {\bf \ref{ch-pre}}.\ref{thm-projlimcohsch1} (1).
Since admissible blow-ups are proper, it is a closed map.
Hence (2) follows from {\bf \ref{ch-pre}}.\ref{cor-projlimcohspacelattice2}.
\end{proof}

By the definition of the projective limit topology, the collection
$$
\left\{\sp^{-1}_X(U)\,\bigg|\,
\begin{minipage}{14em}
{\small $X$ is a formal model of $\mathscr{X}$ and $U$ is a quasi-compact open subset of $X$}
\end{minipage}
\right\}
$$
gives an open basis of $\ZR{\mathscr{X}}$ consisting of quasi-compact open subsets.
The following proposition shows that this collection is actually the set of all quasi-compact open subsets of $\ZR{\mathscr{X}}$:
\begin{prop}\label{prop-zariskiriemanntoptop}
{\rm (1)} For any quasi-compact open subset $\mathfrak{U}$ of $\ZR{\mathscr{X}}$ there exists a formal model $X$ and a quasi-compact open subset $U\subseteq X$ such that $\mathfrak{U}=\sp^{-1}_X(U)$.
Moreover, the formal model $X$ here can be taken to be distinguished.

{\rm (2)} Let $X$ be a coherent adic formal scheme of finite ideal type, and $U\subseteq X$ a quasi-compact open subset.
Set $\mathscr{X}=X^{\rig}$ and $\mathscr{U}=U^{\rig}$.
Then the induced map $\ZR{\mathscr{U}}\rightarrow\ZR{\mathscr{X}}$ maps $\ZR{\mathscr{U}}$ homeomorphically onto the quasi-compact open subset $\sp^{-1}_X(U)$ of $\ZR{\mathscr{X}}$.
\end{prop}

The proposition says that the coherent small admissible site $\mathscr{X}_{\ad}$ is isomorphic to the projective limit of the Zariski sites $X_{\Zar}$ of the formal models.

\begin{proof}
Obviously, it is enough to consider the case where $\ZR{\mathscr{X}}$ is non-empty.

(1) The first assertion follows from {\bf \ref{ch-pre}}.\ref{prop-projlimcohtopspqcptopen}. 
Once such a model $X$ and a quasi-compact open subset $U$ are found, we can find a similar quasi-compact open subset in any admissible blow-up $X'$ of $X$; indeed, we just take $U'=X'\times_XU$ so that $\mathfrak{U}=\sp^{-1}_{X'}(U')$.
By \ref{prop-cohrigidspacedist1}, in particular, we can assume that $X'$ is distinguished.

(2) Consider the ordered sets $\AId_X$ and $\AId_U$ as in \ref{prop-catblowups1}.
The map $\AId_X\rightarrow\AId_U$ given simply by restriction is an ordered map and is surjective by {\bf \ref{ch-formal}}.\ref{prop-extension3}.
Hence the map $\ZR{\mathscr{U}}\rightarrow \ZR{\mathscr{X}}$ in question is projective limit of the system of morphisms $\{U_{\mathscr{J}|_U}\rightarrow X_{\mathscr{J}}\}_{\mathscr{J}\in\AId_X}$.
By \ref{prop-blowups41} these morphisms are open immersions, and we have $U_{\mathscr{J}|_U}\cong X_{\mathscr{J}}\times_XU$.
Since the fiber product by an open immersion of (formal) schemes and the fiber product of the underlying sets coincide ({\bf \ref{ch-formal}}.\ref{prop-sepqsepformal1}), the desired result follows from {\bf \ref{ch-pre}}.\ref{cor-lemlimLRS}. 
 \end{proof}

\begin{prop}\label{prop-zariskiriemanntop2}
The functor $\ZR{\,\cdot\,}\colon\CRf\rightarrow\Top$ maps coherent open immersions\index{immersion!open immersion of rigid spaces@open --- (of rigid spaces)!coherent open immersion of rigid spaces@coherent --- ---} to open immersions and preserves finite colimits consisting of coherent open immersions.
\end{prop}

\begin{proof}
The first assertion follows from \ref{prop-zariskiriemanntoptop} (2).
To show the other assertion, it suffices to show that the functor $\ZR{\,\cdot\,}$ preserves cofiber products by coherent open immersions.
Consider the patching diagram as in \ref{prop-cohrigidspacepatching1}, and take, as in the proof there, the patching $Z=X\amalg_UY$ of coherent adic formal schemes of finite ideal type by open immersions that represents the patching diagram.
Since the restriction maps $\AId_Z\rightarrow\AId_X$, $\AId_Z\rightarrow\AId_Y$, and $\AId_Z\rightarrow\AId_U$ are surjective ({\bf \ref{ch-formal}}.\ref{prop-extension3}), we have
$$
\varprojlim Z_{\mathscr{J}}=\varprojlim X_{\mathscr{J}}\coprod_{\varprojlim U_{\mathscr{J}}}\varprojlim Y_{\mathscr{J}},
$$
that is, $\ZR{\mathscr{X}\amalg_{\mathscr{U}}\mathscr{Y}}=\ZR{\mathscr{X}}\amalg_{\ZR{\mathscr{U}}}\ZR{\mathscr{Y}}$.
\end{proof}

\subsubsection{Non-emptiness}\label{subsub-ZRdefnonempty}
\begin{prop}\label{prop-ZRpoints4}
Let $\mathscr{X}$ be a coherent rigid space, and $X$ a distinguished formal model\index{formal model!formal model of a coherent rigid space@--- (of a coherent rigid space)!distinguished formal model of a coherent rigid space@distinguished --- ---} of $\mathscr{X}$.
Then the specialization map\index{specialization map} $\sp_X\colon\ZR{\mathscr{X}}\rightarrow X$ is surjective.
\end{prop}

\begin{proof}
We may assume in view of \ref{prop-zariskiriemanntoptop} that $X$ is affine $X=\Spf A$.
By the assumption the adic ring $A$ is $I$-torsion free, where $I\subseteq A$ is an ideal of definition.
Note that, in this situation, $Y\setminus V(I)$ is dense in $Y$.
Set $Y=\Spec A$, and consider the canonical map $X=\Spf A\rightarrow Y=\Spec A$.
Take any admissible ideal $J\subseteq A$, and let $Y'\rightarrow Y$ be the blow-up along $J$.
The formal completion $X'\rightarrow X$ is an admissible blow-up of $X$, and any admissible blow-up is of this form.
Since $Y'\rightarrow Y$ is proper and identity outside $V(I)$, it is surjective.
Hence $X'\rightarrow X$ is surjective.
Now the assertion follows from {\bf \ref{ch-pre}}.\ref{thm-projlimcohspacepres} (2).
\end{proof}

\begin{cor}\label{cor-ZRnonempty}
The following conditions for a coherent rigid space $\mathscr{X}$ are equivalent$:$
\begin{itemize}
\item[{\rm (a)}] the associated Zariski-Riemann space $\ZR{\mathscr{X}}$ is non-empty$;$
\item[{\rm (b)}] there exists a non-empty distinguished formal model of $\mathscr{X};$
\item[{\rm (c)}] there exists a non-empty formal model $X$ of $\mathscr{X}$ with an ideal of definition of finite type $\mathscr{I}$ such that $X$ is not $\mathscr{I}$-torsion$;$
\item[{\rm (d)}] $\mathscr{X}$ is not an empty rigid space $(\ref{dfn-emptyrigidspace})$.
\end{itemize}
\end{cor}

\begin{proof}
The implication (a) $\Rightarrow$ (b) is clear.
(b) $\Rightarrow$ (a) follows from \ref{prop-ZRpoints4}.
(b) $\Rightarrow$ (c) is trivial.
Let us show (c) $\Rightarrow$ (b).
Let $X'\rightarrow X$ be the admissible blow-up along $\mathscr{I}$.
Since $X$ is not $\mathscr{I}$-torsion, $X'$ is non-empty.
Thus $X'$ gives a non-empty distinguished formal model of $\mathscr{X}$.
As (d) $\Rightarrow$ (b) is trivial, it only remains to show (b) $\Rightarrow$ (d).
Suppose (b) holds, and let $X$ be a formal model of $\mathscr{X}$.
Then there exists a non-empty distinguished formal model $X'$ that dominates $X$ (\ref{prop-cohrigidspacedist100} and \ref{prop-cohrigidspacedist1}).
Hence $X$ is non-empty, thereby (b) $\Rightarrow$ (d).
\end{proof}

\subsubsection{General construction}\label{subsub-ZRdefgen}
\begin{dfn}[Zariski-Riemann space for general rigid spaces]\label{dfn-generalrigidspace3}{\rm 
(1) Let $\mathscr{F}$ be a sheaf represented by a stretch of coherent rigid spaces $\{\mathscr{U}_i\}_{i\in J}$ {\rm (\ref{dfn-admissiblesite31} (1))}.
Then the induced inductive system $\{\ZR{\mathscr{U}_i}\}_{i\in J}$ is a filtered inductive system of topological spaces such that each transition map is an open immersion (\ref{prop-zariskiriemanntoptop} (2)).
We set
$$
\ZR{\mathscr{F}}=\varinjlim_{i\in J}\ZR{\mathscr{U}_i},
$$
which is a sober topological space\index{space@space (topological)!sober topological space@sober ---} ({\bf \ref{ch-pre}}.\ref{prop-corsoberness}).
It is clear that $\ZR{\mathscr{F}}$ does not depend, up to canonical isomorphisms, on the choice of an inductive system $\{\mathscr{U}_i\}_{i\in J}$ representing $\mathscr{F}$.
The space $\ZR{\mathscr{F}}$ is called the {\em associated Zariski-Riemann space} of $\mathscr{F}$.
For an open immersion $\varphi\colon\mathscr{F}\hookrightarrow\mathscr{G}$ where $\mathscr{F}$ and $\mathscr{G}$ are rigid space represented by a stretch of coherent rigid spaces, one defines
$$
\ZR{\varphi}\colon\ZR{\mathscr{F}}\longrightarrow\ZR{\mathscr{G}}
$$
in a similar way.

(2) Let $\mathscr{F}$ be a rigid space, and take $\mathscr{Y}=\coprod_{\alpha\in L}\mathscr{Y}_{\alpha}\rightarrow\mathscr{F}$ as in \ref{dfn-generalrigidspace1}.
Set
$$
\ZR{\mathscr{Y}}=\coprod_{\alpha\in L}\ZR{\mathscr{Y}_{\alpha}}.
$$
For each $\alpha,\beta\in L$ the projections $\mathrm{pr}\colon\mathscr{Y}_{\alpha}\times_{\mathscr{F}}\mathscr{Y}_{\beta}\rightarrow\mathscr{Y}_{\alpha}$ are open immersions, and hence one defines $\ZR{\mathscr{Y}_{\alpha}\times_{\mathscr{F}}\mathscr{Y}_{\beta}}$ and $\ZR{\mathrm{pr}}$ as in (1).
Set 
$$
\ZR{\mathscr{Y}\times_{\mathscr{F}}\mathscr{Y}}=\coprod_{\alpha,\beta\in L}\ZR{\mathscr{Y}_{\alpha}\times_{\mathscr{F}}\mathscr{Y}_{\beta}},
$$
and $\ZR{\mathrm{pr}}$ similarly.
The maps $\ZR{\mathrm{pr}}$ are local isomorphisms of topological spaces.
Define 
$$
\ZR{\mathscr{F}}=\ZR{\mathscr{Y}}\amalg_{\ZR{\mathscr{Y}\times_{\mathscr{F}}\mathscr{Y}}}\ZR{\mathscr{Y}}
$$
as a topological space.
It is easy to see that $\ZR{\mathscr{F}}$ does not depend on the choice of the presentation $\mathscr{Y}=\coprod_{\alpha\in L}\mathscr{Y}_i\rightarrow\mathscr{F}$.
We call $\ZR{\mathscr{F}}$ the {\em associated Zariski-Riemann space} of $\mathscr{F}$.}
\end{dfn}

Note that (1) and (2) are consistent with each other due to \ref{prop-zariskiriemanntop2}.
The formation $\mathscr{F}\mapsto\ZR{\mathscr{F}}$ is functorial; in other words, for a morphism $\varphi\colon\mathscr{F}\rightarrow\mathscr{G}$ of rigid spaces we have the canonically induced map $\ZR{\varphi}\colon\ZR{\mathscr{F}}\rightarrow\ZR{\mathscr{G}}$ of topological spaces. 
Thus the functor $(\dagger)$ in \S\ref{subsub-ZRdeffunct} extends to 
$$
\ZR{\,\cdot\,}\colon\Rf\longrightarrow\Top.\leqno{(\dagger)}
$$

By the construction of the Zariski-Riemann spaces and \ref{thm-ZRcompact} (1), we have:
\begin{prop}\label{prop-zariskiriemanntopz}
The topological space $\ZR{\mathscr{F}}$ for a rigid space $\mathscr{F}$ is locally coherent {\rm ({\bf \ref{ch-pre}}.\ref{dfn-locallycoherent})}\index{coherent!coherent topological space@--- (topological) space!locally coherent topological space@locally --- ---}\index{space@space (topological)!coherent topological space@coherent ---!locally coherent topological space@locally --- ---} and sober {\rm ({\bf \ref{ch-pre}}, \S\ref{subsub-sober})}\index{space@space (topological)!sober topological space@sober ---}. \hfill$\square$
\end{prop}

\begin{prop}\label{prop-zariskiriemanntop2gen}
The functor $\ZR{\,\cdot\,}\colon\Rf\rightarrow\Top$ maps open immersions\index{immersion!open immersion of rigid spaces@open --- (of rigid spaces)} to open immersions\index{immersion!open immersion of ringed spaces@open --- (of ringed spaces)} and preserves colimits consisting of open immersions.
\end{prop}

\begin{proof}
By \ref{prop-zariskiriemanntop2} one sees that the functor $\ZR{\,\cdot\,}\colon\CRf\rightarrow\Top$ preserves base change of coherent open immersions by coherent open immersions.
Hence, to show that $\ZR{\,\cdot\,}\colon\Rf\rightarrow\Top$ preserves open immersions, it suffices to check that, if $\varphi\colon\mathscr{F}\rightarrow\mathscr{G}$ is an open immersion where $\mathscr{F}$ and $\mathscr{G}$ are rigid space represented by a stretch of coherent rigid spaces, $\ZR{\varphi}\colon\ZR{\mathscr{F}}\longrightarrow\ZR{\mathscr{G}}$ is an open immersion of topological space; but this is clear from the definition.
In order to check that $\ZR{\,\cdot\,}$ preserves colimits by open immersions, it is enough to show that it preserves finite colimits by open immersion, since $\ZR{\,\cdot\,}$ obviously preserves filtered inductive limits by open immersions. 
But this easily reduces to the case of finite colimits of coherent rigid spaces, which we have already dealt with in \ref{prop-zariskiriemanntop2}.
\end{proof}

\subsubsection{Connectedness}\label{subsub-connectedrigidspaces}
\begin{dfn}\label{dfn-connectedrigidpsaces}{\rm 
We say that a rigid space $\mathscr{X}$ is {\em connected}\index{rigid space!connected rigid space@connected ---} if the associated Zariski-Riemann space $\ZR{\mathscr{X}}$ is connected.}
\end{dfn}

\begin{prop}\label{prop-connectedcoherentrigidspaces}
A coherent rigid space $\mathscr{X}$ is connected if and only if any distinguished formal model is connected.
\end{prop}

\begin{proof}
The ``only if'' part follows from \ref{prop-ZRpoints4}.
By \ref{cor-blowups155} $\ZR{\mathscr{X}}$ is the projective limit of a projective system consisting of distinguished formal models.
Hence the ``if'' part follows from Exercise \ref{exer-limitcoherentsoberspaceconnected}.
\end{proof}

\subsubsection{Notation}\label{subsub-notationcohomology}
Let us use the following simplified notation for cohomologies in the sequel.
\begin{ntn}\label{ntn-cohomologyrigidsp1}{\rm 
(1) Let $\mathscr{X}$ be a rigid space, and $\ZR{\mathscr{X}}$ the associated Zariski-Riemann space.
For an abelian sheaf $\mathscr{F}$ on $\ZR{\mathscr{X}}$ we write
$$
\H^q(\mathscr{X},\mathscr{F})=\H^q(\ZR{\mathscr{X}},\mathscr{F})
$$
for any $q\geq 0$.
When $q=0$, we often denote it by $\Gamma(\mathscr{X},\mathscr{F})$ or by $\Gamma_{\mathscr{X}}(\mathscr{F})$.
Note that $\H^q(\mathscr{X},\textrm{--})$ is the $q$-th right derived functor of $\Gamma_{\mathscr{X}}$.

(2) Let $\varphi\colon\mathscr{X}\rightarrow\mathscr{Y}$ be a morphism of rigid spaces.
Then for any abelian sheaf $\mathscr{F}$ on $\ZR{\mathscr{X}}$ we write
$$
\RD^q\varphi_{\ast}\mathscr{F}=\RD^q\ZR{\varphi}_{\ast}\mathscr{F}
$$
for any $q\geq 0$.
When $q=0$, we often denote it by $\varphi_{\ast}\mathscr{F}$.
Note that $\RD^q\varphi_{\ast}$ is the $q$-th right derived functor of $\varphi_{\ast}$.}
\end{ntn}
\index{Zariski-Riemann space|)}

\subsection{Structure sheaves and local rings}\label{sub-ZRstrsheaf}
\subsubsection{Integral structure sheaf}\label{subsub-ZRstrsheafint}
\index{structure sheaf!integral structure sheaf@integral ---|(}
Let $\mathscr{X}$ be a coherent rigid space, and $X$ a formal model\index{formal model!formal model of a coherent rigid space@--- (of a coherent rigid space)} of $\mathscr{X}$.
Recall that the Zariski-Riemann space\index{Zariski-Riemann space} $\ZR{\mathscr{X}}$ has been defined as the underlying topological space of the limit $(\ast)$ in \S\ref{subsub-ZRdef} taken in the category of locally ringed spaces.
Hence it has the canonical sheaf of rings, denoted by $\O^{\int}_{\mathscr{X}}$, such that
$$
(\ZR{\mathscr{X}},\O^{\int}_{\mathscr{X}})=\varprojlim S_X.
$$
For any point $x\in\ZR{\mathscr{X}}$ we have
$$
\O^{\int}_{\mathscr{X},x}=\varinjlim_{\mathscr{J}\in\AId_X}\O_{X_{\mathscr{J}},\sp_{X_{\mathscr{J}}}(x)},
$$
where the inductive limit is taken in the category of local rings $($with local homomorphisms$)$ 
$($cf.\ {\rm {\bf \ref{ch-pre}}.\ref{prop-limLRS}}$)$.
If $\mathscr{U}\hookrightarrow\mathscr{X}$ is a coherent open immersion of coherent rigid spaces, then by \ref{prop-zariskiriemanntoptop} and {\bf \ref{ch-pre}}.\ref{cor-lemlimLRS} we have $\O^{\int}_{\mathscr{X}}|_{\ZR{\mathscr{U}}}=\O^{\int}_{\mathscr{U}}$.
Hence one can extend the definition of $\O^{\int}_{\mathscr{X}}$ to general rigid spaces $\mathscr{X}$ as follows:
\begin{dfn}[Integral structure sheaf]\label{dfn-ZRstrsheaf1}{\rm 
Let $\mathscr{X}$ be a rigid space, and $\ZR{\mathscr{X}}$ the associated Zariski-Riemann space. 
The {\em integral structure sheaf}, denoted by $\O^{\int}_{\mathscr{X}}$, is the sheaf of local rings on the topological space $\ZR{\mathscr{X}}$ such that for any open immersion $j\colon \mathscr{U}\hookrightarrow\mathscr{X}$ from a coherent rigid space $\mathscr{U}$ we have $j^{-1}\O^{\int}_{\mathscr{X}}=\O^{\int}_{\mathscr{U}}$, where $\O^{\int}_{\mathscr{U}}$ is the one defined as above.}
\end{dfn}

The following proposition is easy to see.
\begin{prop}\label{prop-zariskiriemanntopringed}
{\rm (1)} If $\mathscr{U}\hookrightarrow\mathscr{X}$ is an open immersion of rigid spaces, the induced morphism $(\ZR{\mathscr{U}},\O^{\int}_{\mathscr{U}})\rightarrow(\ZR{\mathscr{X}},\O^{\int}_{\mathscr{X}})$ is an open immersion of locally ringed spaces.

{\rm (2)} Let $\mathscr{X}$ be a coherent rigid space, and $\mathfrak{U}$ a quasi-compact open subset of $\ZR{\mathscr{X}}$. Then there exists an open immersion $\mathscr{U}\hookrightarrow\mathscr{X}$ from a coherent rigid space such that the open immersion $(\mathfrak{U},\O^{\int}_{\mathscr{X}}|_{\mathfrak{U}})\hookrightarrow (\ZR{\mathscr{X}},\O^{\int}_{\mathscr{X}})$ is isomorphic to the induced open immersion $(\ZR{\mathscr{U}},\O^{\int}_{\mathscr{U}})\rightarrow(\ZR{\mathscr{X}},\O^{\int}_{\mathscr{X}})$. \hfill$\square$
\end{prop}

\begin{dfn}[Ideal of definition]\label{dfn-ZRstrsheaf2}{\rm 
(1) Let $\mathscr{X}$ be a coherent rigid space, and $\ZR{\mathscr{X}}$ the associated Zariski-Riemann space.
An ideal $\mathscr{I}\subseteq\O^{\int}_{\mathscr{X}}$ is said to be an {\em ideal of definition of finite type}\index{ideal of definition} if there exists a formal model $X$ and an ideal of definition of finite type $\mathscr{I}_X$ such that $\mathscr{I}=(\sp^{-1}_X\mathscr{I}_X)\O^{\int}_{\mathscr{X}}$.

(2) Let $\mathscr{X}$ be a rigid space, and $\ZR{\mathscr{X}}$ the associated Zariski-Riemann space.
An ideal $\mathscr{I}\subseteq\O^{\int}_{\mathscr{X}}$ is said to be an {\em ideal of definition of finite type} if for any open immersion $j\colon\mathscr{U}\hookrightarrow\mathscr{X}$ from a coherent rigid space, $\mathscr{I}|_{\ZR{\mathscr{U}}}$ is an ideal of definition of finite type in the above sense.}
\end{dfn}

The above two definitions (1) and (2) of ideals of definition coincide if $\mathscr{X}$ is a coherent rigid space, that is, if $\mathscr{I}$ is an ideal of definition of finite type in the sense of (2), then there exists a formal model $X$ of $\mathscr{X}$ and an ideal of definition of finite type $\mathscr{I}_X$ of $X$ such that $\mathscr{I}=(\sp^{-1}_X\mathscr{I}_X)\O^{\int}_{\mathscr{X}}$.
Recall that by {\bf \ref{ch-formal}}.\ref{cor-extension2} any coherent rigid space has an ideal of definition of finite type.

\begin{prop}\label{prop-ZRstrsheaf213}
Let $\mathscr{X}$ be a coherent rigid space, and $\mathscr{I}$ an ideal of definition of finite type.
Then there exists a distinguished formal model\index{formal model!formal model of a coherent rigid space@--- (of a coherent rigid space)!distinguished formal model of a coherent rigid space@distinguished --- ---} $X$ of $\mathscr{X}$ and an {\em invertible} ideal of definition $\mathscr{I}_X$ such that $\mathscr{I}=(\sp^{-1}_X\mathscr{I}_X)\O^{\int}_{\mathscr{X}};$ moreover, such formal models are cofinal in the category of all formal models of $\mathscr{X}$.
\end{prop}

\begin{proof}
Let $X$ be a formal model of $\mathscr{X}$, and $\mathscr{I}_X$ an ideal of definition of finite type on $X$ such that $\mathscr{I}=(\sp_X^{-1}\mathscr{I}_X)\O^{\int}_{\mathscr{X}}$.
Replacing $X$ by the admissible blow-up along $\mathscr{I}_X$, we can make $\mathscr{I}_X$ invertible, and hence $X$ is distinguished.
\end{proof}

\begin{cor}\label{cor-ZRstrsheaf2131}
Let $\mathscr{X}$ be a rigid space, and $\mathscr{I}$ an ideal of definition of finite type.
Then $\mathscr{I}$ is an invertible ideal of $\O^{\int}_{\mathscr{X}}$, and hence $\O^{\int}_{\mathscr{X}}$ is $\mathscr{I}$-torsion free.
\end{cor}

\begin{proof}
We may assume that $\mathscr{X}$ is coherent.
By \ref{prop-ZRstrsheaf213} there exists a formal model $X$ with an invertible ideal of definition $\mathscr{I}_X$ such that $\mathscr{I}=(\sp_X^{-1}\mathscr{I}_X)\O^{\int}_{\mathscr{X}}$.
Let $\pi\colon X'\rightarrow X$ be an admissible blow-up.
Then $\mathscr{I}_{X'}=\mathscr{I}_X\O_{X'}$ is an invertible ideal of definition (\ref{cor-blowups155}).
Since $\varinjlim_{X'}\sp^{-1}_{X'}\mathscr{I}_{X'}$ is an ideal of $\O^{\int}_{\mathscr{X}}$, we have 
$$
\mathscr{I}=\varinjlim_{X'}\sp^{-1}_{X'}\mathscr{I}_{X'}.
$$
By this we easily deduce that $\O^{\int}_{\mathscr{X}}$ is $\mathscr{I}$-torsion free and that $\mathscr{I}$ is invertible.
\end{proof}
\index{structure sheaf!integral structure sheaf@integral ---|)}

\subsubsection{Rigid structure sheaf}\label{subsub-ZRstrsheafrig}
\index{structure sheaf!rigid structure sheaf@(rigid) ---|(}
\begin{prop}\label{prop-ZRstrsheaf21}
Let $\mathscr{X}$ be a coherent rigid space, $\ZR{\mathscr{X}}$ the associated Zariski-Riemann space, and $x\in\ZR{\mathscr{X}}$ a point.
Let $\mathscr{I}$ be an ideal of definition of finite type, and set $I_x=\mathscr{I}_x$.
Then the stalk $A_x=\O^{\int}_{\mathscr{X},x}$ of the integral structure sheaf is $I_x$-valuative\index{valuative!Ivaluative ring@$I$-{---} ring} $(${\rm {\bf \ref{ch-pre}}.\ref{dfn-valuative1}}$)$ and $I_x$-adically henselian\index{henselian!I-adically henselian@$I$-adically ---} {\rm ({\bf \ref{ch-pre}}, \S\ref{subsub-henselianpairs})}.
\end{prop}

\begin{proof}
Let $\mathscr{X}=X^{\rig}$ and $\mathscr{I}=\mathscr{I}_X\O^{\int}_{\mathscr{X}}$, where $\mathscr{I}_X$ is an ideal of definition of $X$ of finite type.
Let $\til{I}=\mathscr{I}_{X,\sp_X(x)}$.
Then $A_x$ is the filtered inductive limit of $\til{I}$-adically complete local rings and hence is henselian with respect to $I_x$ ({\bf \ref{ch-pre}}.\ref{prop-relpair31}).

Let $J$ be an $I_x$-admissible ideal ({\bf \ref{ch-pre}}.\ref{dfn-adm}) of $A_x$.
We need to show that $J$ is an invertible ideal of $A_x$.
Replacing $X$ by an admissible blow-up if necessary, we may assume that there exists a finitely generated ideal $\til{J}\subseteq\O_{X,\sp_X(x)}$ such that $\til{J}A_x=J$.
Now, since we may work locally around $x\in\ZR{\mathscr{X}}$, we can replace $\mathscr{X}$ by any coherent open subspace $\mathscr{U}\subseteq\mathscr{X}$ such that $x\in\ZR{\mathscr{U}}$.
In view of \ref{prop-zariskiriemanntopringed} (2) this allows us to replace $X$ by any quasi-compact open neighborhood $U\subseteq X$ of $\sp_X(x)$.
We may therefore assume that $X$ is affine $X=\Spf A$ and, moreover, that $A$ has a finitely generated ideal of definition $I\subseteq A$ and an $I$-admissible ideal $J_X\subseteq A$ such that $I^{\Delta}=\mathscr{I}_X$ and $J_X\O_{X,\sp_X(x)}=\til{J}$ (hence $J_XA_x=J$).
Now we consider the admissible blow-up $X'\rightarrow X$ along $J_X$.
Then $J_X\O_{X'}$ is invertible.
Since $\O^{\int}_{\mathscr{X}}$ is $\mathscr{I}_X$-torsion free due to \ref{cor-ZRstrsheaf2131}, it is $J_X$-torsion free, and hence $J=J_X\O^{\int}_{\mathscr{X},x}$ is an invertible ideal of $A_x=\O^{\int}_{\mathscr{X},x}$, as desired.
\end{proof}

\begin{cor}\label{cor-ZRstrsheaf211}
Let $\mathscr{X}$ be a coherent rigid space, $\ZR{\mathscr{X}}$ the associated Zariski-Riemann space, and $x\in\ZR{\mathscr{X}}$ a point.
Then for any ideal of definition $\mathscr{I}$ of finite type of $\ZR{\mathscr{X}}$, we have $\mathscr{I}_x=(a)$ by a non-zero-divisor $a\in\O^{\int}_{\mathscr{X},x}$. \hfill$\square$
\end{cor}

By {\bf \ref{ch-pre}}.\ref{thm-valuative} and {\bf \ref{ch-pre}}.\ref{prop-patching4} we have:
\begin{cor}\label{cor-ZRstrsheaf212}
Let $\mathscr{X}$ be a coherent rigid space, $\ZR{\mathscr{X}}$ the associated Zariski-Riemann space, and $x\in\ZR{\mathscr{X}}$ a point.
Let $\mathscr{I}$ be an ideal of definition of finite type of $\ZR{\mathscr{X}}$, and set $(A_x,I_x)=(\O^{\int}_{\mathscr{X},x},\mathscr{I}_x)$.
Set $I_x=(a)$ $($cf.\ {\rm \ref{cor-ZRstrsheaf211}}$)$ and $J_x=\bigcap_{n\geq 1}I^n_x$.
Then$:$
\begin{itemize}
\item[{\rm (a)}] $B_x={\displaystyle \varinjlim_{n\geq 1}}\Hom(I^n_x,A_x)=A_x[\frac{1}{a}]$ is a local ring, and $V_x=A_x/J_x$ is a valuation ring, $\ovl{a}$-adically separated\index{valuation!valuation ring@--- ring!a-adically separated valuation ring@$a$-adically separated --- ---}, for the residue field $K_x$ of $B_x$, where $I_xV_x=(\ovl{a})$.
\item[{\rm (b)}] $A_x=\{f\in B_x\,|\, (f\ \mathrm{mod}\ \m_{B_x})\in V_x\}$.
\item[{\rm (c)}] $J_x=\m_{B_x}$.
\end{itemize}
Moreover, $B_x$ is a henselian local ring, and $V_x$ is henselian with respect to $\ovl{a}$-adic topology. \hfill$\square$
\end{cor}

\begin{dfn}[Rigid structure sheaf]\label{dfn-ZRstrsheaf3}{\rm 
Let $\mathscr{X}$ be a rigid space. 
We define a sheaf of rings $\O_{\mathscr{X}}$ on its associated Zariski-Riemann space $\ZR{\mathscr{X}}$ as follows.

(1) If $\mathscr{X}$ is a coherent rigid space, then 
$$
\O_{\mathscr{X}}=\varinjlim_{n>0}\lHom_{\O^{\int}_{\mathscr{X}}}(\mathscr{I}^n,\O^{\int}_{\mathscr{X}}),
$$
where $\mathscr{I}$ is an ideal of definition of finite type (cf.\ Deligne's formula \cite[(6.9.17)]{EGAInew}).
(Exercise \ref{exer-propZRstrsheaf2} verifies that this does not depend on the choice of $\mathscr{I}$.)

(2) In general, we construct $\O_{\mathscr{X}}$ by patching; that is, $\O_{\mathscr{X}}$ is the sheaf such that for any open immersion $\mathscr{U}\hookrightarrow\mathscr{X}$ from a coherent rigid space, $\O_{\mathscr{X}}|_{\ZR{\mathscr{U}}}$ is the one defined as in (1).}
\end{dfn}

We call the sheaf $\O_{\mathscr{X}}$ the (rigid) structure sheaf of $\ZR{\mathscr{X}}$.
By the local description given in \ref{cor-ZRstrsheaf212} we have the following:
\begin{prop}\label{prop-ZRstrsheaf4}
The sheaf $\O_{\mathscr{X}}$ is a sheaf of henselian local rings, that is, for any $x\in\ZR{\mathscr{X}}$ the stalk $\O_{\mathscr{X},x}$ is a henselian local ring. \hfill$\square$
\end{prop}
\index{structure sheaf!rigid structure sheaf@(rigid) ---|)}

\subsubsection{Zariski-Riemann triple}\label{subsub-ZRstrsheaftriple}
\index{Zariski-Riemann triple|(}
We have so far obtained two sheaves $\O^{\int}_{\mathscr{X}}$ and $\O_{\mathscr{X}}$, by which the space $\ZR{\mathscr{X}}$ is endowed with locally ringed structures in two ways.

\begin{dfn}\label{dfn-ZRtriple}{\rm 
Let $\mathscr{X}$ be a rigid space. 
Then the triple
$$
\ZRT(\mathscr{X})=(\ZR{\mathscr{X}},\O^{\int}_{\mathscr{X}},\O_{\mathscr{X}})
$$
is called the {\em Zariski-Riemann triple} associated to the rigid space $\mathscr{X}$.}
\end{dfn}

See \S\ref{sub-triples} for a general theory of triples.

\medskip\noindent
{\bf Convention.} {\sl By the {\em structure sheaf} of a rigid space $\mathscr{X}$ we always mean the rigid structure sheaf $\O_{\mathscr{X}}$ on $\ZR{\mathscr{X}}$, unless otherwise clearly stated.}
\index{Zariski-Riemann triple|)}

\subsubsection{Reducedness}\label{subsub-reducedrigidspaces}
\begin{dfn}\label{dfn-reducedrigidspaces}{\rm 
We say that a rigid space $\mathscr{X}$ is {\em reduced}\index{rigid space!reduced rigid space@reduced ---} if the ringed space $(\ZR{\mathscr{X}},\O_{\mathscr{X}})$ is reduced in the sense as in {\bf \ref{ch-pre}}, \S\ref{sub-ringedspgenpre}.}
\end{dfn}

It will be shown in \ref{cor-reducedrigidspaces} that a coherent universally Noetherian rigid space $\mathscr{X}$ is reduced if and only if any distinguished formal model of $\mathscr{X}$ is reduced.

\subsubsection{Description of the local rings}\label{subsub-fibersoverrigptslocalrings}
\begin{ntn}\label{ntn-ZRpoints}{\rm 
Let $\mathscr{X}$ be a rigid space, and $x\in\ZR{\mathscr{X}}$.
We will often use the following notations (cf.\ \ref{cor-ZRstrsheaf212}):
\begin{itemize}
\item $A_x=\O^{\int}_{\mathscr{X},x}$.
\end{itemize}
One can always find a coherent open neighborhood of $x$ of the form $\ZR{(\Spf A)^{\rig}}$ such that the adic ring $A$ has a principal invertible ideal of definition $I=(a)$.
Considering the invertible ideal
\begin{itemize}
\item $I_x=IA_x=aA_x$
\end{itemize}
of $A_x$, one has
\begin{itemize}
\item $B_x=A_x[\frac{1}{a}]=\O_{\mathscr{X},x}$, $K_x=B_x/\m_{B_x}$;
\item $J_x=\m_{B_x}=\bigcap_{n\geq 1}I^n_x$,
\item $V_x=A_x/J_x$, $k_x=V_x/\m_{V_x}$.
\end{itemize}
Let $\ovl{a}$ be the image of $a$ in $V_x$.
Then $V_x$ is $\ovl{a}$-adically separated and henselian valuation ring with the residue field $k_x$ such that $\Frac(V_x)=K_x$, and we have $A_x=\{f\in B_x\,|\, (f\ \mathrm{mod}\ \m_{B_x})\in V_x\}$.}
\end{ntn}

Notice that the objects $A_x$, $B_x$, $K_x$, $J_x$, $V_x$, and $k_x$ do not depend on the choice of $I=(a)\subseteq A$.
When we want to describe these local data, the following notion will be often useful:
\begin{dfn}\label{dfn-formalneighborhoods}{\rm 
Let $\mathscr{X}$ be a rigid space, and $x\in\ZR{\mathscr{X}}$ a point.

(1) A {\em formal neighborhood}\index{formal neighborhood} of $x$ in $\mathscr{X}$ is a pair $(U,\iota)$ consisting of a coherent adic formal scheme of finite ideal type $U$ and an open immersion 
$$
\iota\colon\mathscr{U}=U^{\rig}\longhookrightarrow\mathscr{X}
$$
such that $x\in\ZR{\iota}(\ZR{\mathscr{U}})$.

(2) By a morphism between formal neighborhoods $(U,\iota)$ and $(U',\iota')$ we mean an adic morphism $h\colon U\rightarrow U'$ such that the resulting diagram
$$
\xymatrix@R-3ex{U^{\rig}\ar@{^{(}->}[dr]^{\iota}\ar[dd]_{h^{\rig}}\\ &\mathscr{X}\\ {U'}^{\rig}\ar@{^{(}->}[ur]_{\iota'}}
$$
commutes.}
\end{dfn}

We denote by $\FN_{\mathscr{X},x}$ or simply by $\FN_x$ the category of formal neighborhoods of $x$ in $\mathscr{X}$.
Formal neighborhoods are usually considered up to isomorphism in this category.

A formal neighborhood $U=(U,\iota)$ is said to be {\em affine}\index{formal neighborhood!affine formal neighborhood@affine ---} if the formal scheme $U$ is affine.
Notice here that the category $\FN_{\mathscr{X},x}$ is canonically cofiltered, and affine formal neighborhoods give a cofinal family.

Suppose that $\mathscr{X}$ is coherent, and let $X$ be a formal model of $\mathscr{X}$.
Then any quasi-compact Zariski open neighborhood $U\subseteq X$ of $\sp_X(x)$ gives rise to a formal neighborhood $(U,i^{\rig})$ of $x$, where $i\colon U\hookrightarrow X$ is the inclusion map.
By \ref{prop-zariskiriemanntoptop}, considering all quasi-compact Zariski (affine) open neighborhoods of $\sp_{X'}(x)$ in any admissible blow-up $X'$ of $X$, one has the a system of formal neighborhoods that gives a cofinal system of quasi-compact open neighborhoods of $x$ in $\ZR{\mathscr{X}}$.

Now suppose we are in the situation as in \ref{ntn-ZRpoints} with $\mathscr{X}=X^{\rig}$ where $X=\Spf A$.
We fix a system of affine formal neighborhoods considered with embedded formal models
$$
\{(\mathscr{U}_{\alpha}=(\Spf A_{\alpha})^{\rig},j_{\alpha}\colon U_{\alpha}=\Spf A_{\alpha}\hookrightarrow X_{\alpha})\}_{\alpha\in L}
$$
indexed by a directed set $L$; each $U_{\alpha}=\Spf A_{\alpha}$ is an affine open subset of an admissible blow-up $X_{\alpha}$ of $X=\Spf A$ that contains the point $\sp_{\alpha}(x)$, where we denote the specialization map $\sp_{X_{\alpha}}\colon\ZR{\mathscr{X}}\rightarrow X_{\alpha}$ simply by $\sp_{\alpha}$.
(Here we may assume that the directed set $L$ has the minimum, say $0\in L$, and $A_0=A$.)
For $\alpha\leq\beta$ we assume that there exists a commutative diagram
$$
\xymatrix{U_{\beta}\,\ar@{^{(}->}[r]^{j_{\beta}}\ar[d]&X_{\beta}\ar[d]^{\pi_{\beta\alpha}}\\ U_{\alpha}\,\ar@{^{(}->}[r]_{j_{\alpha}}&X_{\alpha}\rlap{,}}
$$
where $\pi_{\beta\alpha}$ is an admissible blow-up, such that for $\alpha\leq\beta\leq\gamma$ we have $\pi_{\gamma\alpha}=\pi_{\beta\alpha}\circ\pi_{\gamma\beta}$.
In this situation we have the inductive system of rings $\{A_{\alpha}\}_{\alpha\in L}$.
Notice that by \ref{cor-blowups155} each $A_{\alpha}$ has the invertible ideal of definition $IA_{\alpha}=(a)$.

Taking $L$ to be the set of all pairs $(X',U')$ consisting of admissible blow-ups $X'$ of $X_0$ and an affine neighborhood $U'\subseteq X'$ of $\sp_{X'}(x)$ with the ordering defined as above, we can construct a system as above such that the image of the forgetful map $(\mathscr{U}_{\alpha},j_{\alpha}\colon U_{\alpha}\hookrightarrow X_{\alpha})\mapsto(\mathscr{U}_{\alpha},j^{\rig}_{\alpha})$ is cofinal in $\FN_{\mathscr{X},x}$.
\begin{prop}\label{prop-descriptionlocalrings}
{\rm (1)} We have the following canonical identifications$:$
$$
A_x=\O^{\int}_{\mathscr{X},x}=\varinjlim_{\alpha\in L}A_{\alpha},\qquad B_x=\O_{\mathscr{X},x}=\varinjlim_{\alpha\in L}{\textstyle A_{\alpha}[\frac{1}{a}]}.
$$

{\rm (2)} Let $J_{\alpha}$ for each $\alpha\in L$ be the ideal of $A_{\alpha}$ that is the pull-back of $J_x$ by the map $A_{\alpha}\rightarrow A_x=\O^{\int}_{\mathscr{X},x}$.
Then we have $J_x=\varinjlim_{\alpha\in L}J_{\alpha}$ and 
$$
V_x=\varinjlim_{\alpha}A_{\alpha}/J_{\alpha}.
$$
\end{prop}

\begin{proof}
By the definition of the sheaf $\O^{\int}_{\mathscr{X}}$ (\S\ref{subsub-ZRstrsheafint}) we have
$$
A_x=\varinjlim_{\alpha\in L}\Gamma(\mathscr{U}_{\alpha},\O^{\int}_{\mathscr{X}})=\varinjlim_{\alpha\in L}\varinjlim_{X'}\Gamma(X',\O_{X'}),
$$
where $X'$ in the second limit runs through the totality of all admissible blow-ups of $\Spf A_{\alpha}$, and then one deduces by a standard argument (cf.\ {\bf \ref{ch-pre}}.\ref{prop-final}), using the extension of admissible blow-ups (\ref{prop-blowups4111}), that the last limit is canonically isomorphic to $\varinjlim_{\alpha\in L}A_{\alpha}$.
Then the second equality follows immediately.
(2) follows from the first equality of (1) and the exactness of filtered inductive limits (cf.\ {\bf \ref{ch-pre}}.\ref{prop-directlimits2}).
\end{proof}

\subsubsection{Generization maps}\label{subsub-fibersoverrigptsbehavior}
Let $\mathscr{X}$ be a rigid space, and $x,x'\in\ZR{\mathscr{X}}$, and suppose $x'$ is a generization\index{generization} of $x$ (cf.\ {\bf \ref{ch-pre}}, \S\ref{subsub-genspetopsp}).
We have the generization map\index{generization!generization map@--- map} ({\bf \ref{ch-pre}}, \S\ref{subsub-generizationmap})
$$
A_x\longrightarrow A_{x'}.\leqno{(\ast)}
$$

\begin{prop}\label{prop-generizationmapsrigidspaces}
The generization map $(\ast)$ maps $J_x$ to $J_{x'}$, and induces a {\em local} homomorphism $B_x\rightarrow B_{x'}$.
In particular, it induces an injective morphism $K_x\hookrightarrow K_{x'}$, which maps $V_x$ injectively into $V_{x'}$.
\end{prop}

\begin{proof}
Since any open neighborhood of $x$ in $\ZR{\mathscr{X}}$ contains $x'$, we may work in the situation as in \ref{ntn-ZRpoints} where $\mathscr{X}=X^{\rig}$ with $X=\Spf A$.
Then it is clear that the generization map $A_x\rightarrow A_{x'}$, which is obviously $a$-adic, maps $J_x$ to $J_{x'}$ and induces $B_x\rightarrow B_{x'}$ and $V_x\rightarrow V_{x'}$.
Since $J_x=\m_{B_x}$ and $J_{x'}=\m_{B_{x'}}$, the map $B_x\rightarrow B_{x'}$ is a local homomorphism and hence induces the injective morphism $K_x\hookrightarrow K_{x'}$ between the residue fields.
Since $K_x=\Frac(V_x)$ and $K_{x'}=\Frac(V_{x'})$, the morphism $V_x\rightarrow V_{x'}$ is injective.
\end{proof}

\begin{thm}\label{thm-fibersoverrigptsbehavior}
Suppose $\mathscr{X}$ is locally universally Noetherian\index{rigid space!universally Noetherian rigid space@universally Noetherian ---!locally universally Noetherian rigid space@locally --- ---} {\rm (\ref{dfn-universallyadhesiverigidspaces})}. 

{\rm (1)} Let $\mathfrak{p}'$ be the prime ideal of $A_x=\O^{\int}_{\mathscr{X},x}$ that is the pull-back of the maximal ideal of $A_{x'}=\O^{\int}_{\mathscr{X},x'}$ by the generization map $A_x\rightarrow A_{x'}$.
Then the induced $a$-adic map $(A_x)_{\mathfrak{p}'}\rightarrow A_{x'}$ is faithfully flat, and $(A_x)_{\mathfrak{p}'}/(a)\rightarrow A_{x'}/(a)$ is an isomorphism.

{\rm (2)} Set $\mathfrak{q}'=\mathfrak{p}'/J_x$.
Then the $a$-adic map $(V_x)_{\mathfrak{q}'}\rightarrow V_{x'}$ induced from the generization map $V_x\rightarrow V_{x'}$ {\rm (\ref{prop-generizationmapsrigidspaces})} is faithfully flat, and the induced map $(V_x)_{\mathfrak{q}'}/(a)\rightarrow V_{x'}/(a)$ is an isomorphism$;$ in particular, we have the natural isomorphism
$$
(V_x)_{\mathfrak{q}'}\!\!\!\!\widehat{\ }\stackrel{\sim}{\longrightarrow}\widehat{V}_{x'}
$$
between the $a$-adic completions.

{\rm (3)} The generization map $B_x\rightarrow B_{x'}$ {\rm (\ref{prop-generizationmapsrigidspaces})} is faithfully flat.
\end{thm}

\begin{proof}
Take a cofinal system of affine formal neighborhoods $\{\mathscr{U}_{\alpha}=(\Spf A_{\alpha})^{\rig}\}_{\alpha\in L}$ of $x$ as in \S\ref{subsub-fibersoverrigptslocalrings} and the filtered inductive system $\{J_{\alpha}\}_{\alpha\in L}$ of ideals as in \ref{prop-descriptionlocalrings} (2).
Since each affinoid neighborhood $\ZR{\mathscr{U}_{\alpha}}$ of $x$ contains $x'$, replacing the index set $L$ by a larger directed set, one can form analogous cofinal system of affine neighborhoods $\{\mathscr{U}'_{\alpha}=(\Spf A'_{\alpha})^{\rig}\}_{\alpha\in L}$ of $x'$ and a filtered inductive system $\{J'_{\alpha}\}_{\alpha\in L}$ of ideals such that $U'_{\alpha}=\Spf A'_{\alpha}$ is an affine open subset of $U_{\alpha}=\Spf A_{\alpha}$ for each $\alpha\in L$.
Indeed, for each $\alpha\in L$ one can form a cofinal system of open neighborhoods of $x'$ of the form $\{\mathscr{U}'_{\alpha,\lambda}=(\Spf A'_{\alpha,\lambda})^{\rig}\}_{\lambda\in\Lambda_{\alpha}}$ in $\mathscr{U}_{\alpha}$; hence one can replace $L$ by the set $\{(\alpha,\lambda)\,|\,\lambda\in\Lambda_{\alpha}\}$, by defining $\mathscr{U}_{\alpha,\lambda}=\mathscr{U}_{\alpha}$.

Let $\mathfrak{p}'_{\alpha}$ be the open prime ideal of $A_{\alpha}$ that is the image of $x'$ by the specialization map $\sp_{\alpha}\colon\ZR{\mathscr{U}_{\alpha}}\rightarrow U_{\alpha}=\Spf A_{\alpha}$.
We can assume that each $A'_{\alpha}$ is a complete localization of $A_{\alpha}$ of the form $A'_{\alpha}=(A_{\alpha})_{\{f_{\alpha}\}}$, where $f_{\alpha}\in A_{\alpha}$ is an element not contained in $\mathfrak{p}'_{\alpha}$.
The sequence of open prime ideals $\{\mathfrak{p}'_{\alpha}\}_{\alpha}$ gives by passage to the inductive limit the open prime ideal $\mathfrak{p}'$ of $\O^{\int}_{\mathscr{X},x}=\varinjlim_{\alpha\in L}A_{\alpha}$ corresponding to the generization $x'$.

(1) Since $A_{\alpha}$ is a t.u.\ rigid-Noetherian\index{t.u. rigid-Noetherian ring@t.u.\ rigid-Noetherian ring} ring, the canonical map $A_{\alpha}\rightarrow (A_{\alpha})_{\{f_{\alpha}\}}$ is flat ({\bf \ref{ch-pre}}.\ref{prop-btarf1}).
Hence by \ref{prop-descriptionlocalrings} (1) it follows that the morphism $A_x=\varinjlim_{\alpha\in L}A_{\alpha}\rightarrow\varinjlim_{\alpha\in L}A'_{\alpha}=A_{x'}$ is flat.
To show that $(A_x)_{\mathfrak{p}'}\rightarrow A_{x'}$ is faithfully flat, it suffices to show that $(A_x)_{\mathfrak{p}'}/(a)\rightarrow A_{x'}/(a)$ is an isomorphism ({\bf \ref{ch-pre}}.\ref{lem-propzariskianfaithfulltflat}).
By the exactness of the functor by filtered inductive limits, we know that the ring $(A_x)_{\mathfrak{p}'}/(a)$ is isomorphic to the inductive limit
$$
\varinjlim_{\alpha\in L}(A_{\alpha}/(a))_{\ovl{\mathfrak{p}}'_{\alpha}},
$$
where $\ovl{\mathfrak{p}}'_{\alpha}=\mathfrak{p}'_{\alpha}/(a)$.
Since this limit is obviously isomorphic to $\varinjlim_{\alpha\in L}A'_{\alpha}/(a)$, we have the desired assertion.

(2) By \ref{prop-generizationmapsrigidspaces} the map $V_x\rightarrow V_{x'}$ is injective, and hence is flat.
The pull-back of the maximal ideal $\mathfrak{m}_{V_{x'}}$ is the open prime ideal $\mathfrak{q}'$.
To show that $(V_x)_{\mathfrak{q}'}\rightarrow V_{x'}$ is faithfully flat, it suffices to show that $(V_x)_{\mathfrak{q}'}/(a)\rightarrow V_{x'}/(a)$ is an isomorphism ({\bf \ref{ch-pre}}.\ref{lem-propzariskianfaithfulltflat}), which follows from the fact that $(A_x)_{\mathfrak{p}'}/(a)\rightarrow A_{x'}/(a)$ is an isomorphism, as shown in (1).

(3) follows immediately from (1) due to the equalities $B_x=A_x[\frac{1}{a}]$ and $B_{x'}=A_{x'}[\frac{1}{a}]$.
\end{proof}



\subsection{Points on Zariski-Riemann spaces}\label{sub-ZRpoints}
\subsubsection{Rigid points}\label{subsub-rigpoints}
\index{point!rigid point@rigid ---|(}\index{rigid point|(}

\begin{dfn}\label{dfn-ZRpoints2}{\rm 
(1) A {\em rigid point} of a rigid space $\mathscr{X}$ is a morphism of rigid spaces of the form
$$
\alpha\colon\mathscr{T}=(\Spf V)^{\rig}\longrightarrow\mathscr{X},
$$
where $V$ is an $a$-adically complete valuation ring\index{valuation!valuation ring@--- ring!a-adically complete valuation ring@$a$-adically complete --- ---} $(a\in \m_V\setminus\{0\})$.

(2) Let $X$ be a coherent adic formal scheme of finite ideal type.
A {\em rigid point} of $X$ is an adic morphism\index{morphism of formal schemes@morphism (of formal schemes)!adic morphism of formal schemes@adic ---} of the form
$$
\alpha\colon\Spf V\longrightarrow X,
$$
where $V$ is an $a$-adically complete valuation ring ($a\in V\setminus\{0\}$).}
\end{dfn}

Two rigid points $\alpha\colon\mathscr{T}=(\Spf V)^{\rig}\rightarrow\mathscr{X}$ and $\beta\colon\mathscr{S}=(\Spf W)^{\rig}\rightarrow\mathscr{X}$ are said to be isomorphic if there exists an isomorphism of rigid spaces $\mathscr{T}\stackrel{\sim}{\rightarrow}\mathscr{S}$ over $\mathscr{X}$; the isomorphy of rigid points of $X$ is defined similarly.

\begin{dfn}\label{dfn-ZRpoints1}{\rm 
Let $\mathscr{X}$ be a rigid space, and $Y$ an adic formal scheme of finite ideal type ({\bf \ref{ch-formal}}.\ref{dfn-adicformalschemesoffiniteidealtype}).
Let $\alpha\colon Y\rightarrow(\ZR{\mathscr{X}},\O^{\int}_{\mathscr{X}})$ be a morphism of locally ringed spaces.
We say that the map $\alpha$ is {\em adic} if for any open immersion $\mathscr{U}\hookrightarrow\mathscr{X}$ and an ideal of definition of finite type $\mathscr{I}$ of $\ZR{\mathscr{U}}$, the ideal $(\alpha^{-1}\mathscr{I})\O_{\alpha^{-1}(\ZR{\mathscr{U}})}$ is an ideal of definition of the open formal subscheme $\alpha^{-1}(\ZR{\mathscr{U}})$ of $Y$.}
\end{dfn}

\begin{prop}\label{prop-ZRpointx1}
There exists a canonical bijection between the set of all isomorphism classes of rigid points of a rigid space $\mathscr{X}$ and the set of all isomorphism classes of adic morphisms $\alpha\colon\Spf V\rightarrow(\ZR{\mathscr{X}},\O^{\int}_{\mathscr{X}})$ of locally ringed spaces, where $V$ is an $a$-adically complete valuation ring, where two such morphisms $\alpha\colon\Spf V\rightarrow\ZR{\mathscr{X}}$ and $\beta\colon\Spf W\rightarrow\ZR{\mathscr{X}}$ are said to be isomorphic if there exists an isomorphism of formal schemes $\Spf V\stackrel{\sim}{\rightarrow}\Spf W$ such that the resulting triangle diagram commutes.
\end{prop}

The bijection maps a rigid point $\alpha\colon\mathscr{T}=(\Spf V)^{\rig}\rightarrow\mathscr{X}$ to the map $\ZR{\alpha}$; note that we have $\ZR{(\Spf V)^{\rig}}=\Spf V$ (Exercise \ref{exer-exavaluationrigidspace}).

\begin{proof}
We give the inverse map.
Let $\alpha\colon\Spf V\rightarrow(\ZR{\mathscr{X}},\O^{\int}_{\mathscr{X}})$ be an adic morphism, and $x$ the image of the closed point.
Take a coherent open subspace $\mathscr{U}\subseteq\mathscr{X}$ such that $x\in\ZR{\mathscr{U}}$.
We have $\alpha(\Spf V)\subseteq\ZR{\mathscr{U}}$.
Let $U$ be a formal model of $\mathscr{U}$, and consider the specialization map $\sp_U\colon\ZR{\mathscr{U}}\rightarrow U$.
Then the composition $\sp_U\circ\alpha$ gives a rigid point of $U$.
Taking $\cdot^{\rig}$, we get a rigid point $\mathscr{T}=(\Spf V)^{\rig}\rightarrow\mathscr{U}\hookrightarrow\mathscr{X}$.
\end{proof}

In the sequel, by a rigid point of a rigid space $\mathscr{X}$ we sometimes mean a morphism of locally ringed spaces as in {\rm \ref{prop-ZRpointx1}}.

\begin{prop}\label{prop-ZRpoint32}
Let $\mathscr{X}$ be a rigid space.

{\rm (1)} For $x\in\ZR{\mathscr{X}}$ there exists a rigid point of the form 
$$
\alpha_x\colon\Spf\widehat{V}_x\longrightarrow\ZR{\mathscr{X}}, \qquad \m_{\widehat{V}_x}\longmapsto x,
$$
where $\widehat{V}_x$ is the $a$-adic completion of $V_x$ $($in the notation as in {\rm \ref{ntn-ZRpoints})}, such that the induced map of stalks at $x$ is the canonical map $A_x=\O^{\int}_{\mathscr{X},x}\rightarrow\widehat{V}_x$.

{\rm (2)} Let $\alpha\colon\Spf V\rightarrow\ZR{\mathscr{X}}$ be a rigid point such that $\alpha(\m_V)=x$.
Then there exists uniquely an injective homomorphism $\widehat{V}_x\hookrightarrow V$ such that $V$ dominates $\widehat{V}_x$ and the diagram
$$
\xymatrix@C-3ex@R-4ex{
\Spf V\ar[dd]\ar[dr]^{\alpha}\\
&\ZR{\mathscr{X}}\\
\Spf\widehat{V}_x\ar[ur]_{\alpha_x}}
$$
commutes.
\end{prop}

\begin{proof}
We may assume $\mathscr{X}$ is coherent.

(1) Take the $\ovl{a}$-adic completion $\widehat{V}_x$ of the $\ovl{a}$-adically separated valuation ring $V_x$.
The map $A_x=\O^{\int}_{\mathscr{X},x}\rightarrow\widehat{V}_x$ induces the inductive system of homomorphisms $\{\O_{X',\sp_{X'}(x)}\rightarrow\widehat{V}_x\}$, where $X'$ runs through all admissible blow-ups $X'\rightarrow X$ of $X$, and hence the projective system of adic morphisms $\{\Spf\widehat{V}_x\rightarrow X'\}$.
The desired map $\alpha_x$ is the projective limit of this system of morphisms.

(2) The morphism $\alpha$ gives an $a$-adic homomorphism $\O^{\int}_{\mathscr{X},x}\rightarrow V$.
Since $V$ is $a$-adically complete, there exists a unique factorizing map $\widehat{V}_x\rightarrow V$.
Since the last map is $a$-adic, we deduce by {\bf \ref{ch-pre}}.\ref{prop-maxspe4} that the map $\widehat{V}_x\rightarrow V$ is injective.
Since $\alpha(\m_V)=\alpha_x(\m_{\widehat{V}_x})=x$, $V$ dominates $\widehat{V}_x$, as desired.
\end{proof}

\begin{dfn}\label{dfn-ZRpoints32}{\rm 
Let $x\in\ZR{\mathscr{X}}$ be a point.
Then the rigid point $\alpha_x\colon\Spf\widehat{V}_x\rightarrow\ZR{\mathscr{X}}$ as in \ref{prop-ZRpoint32} (1) is called the {\em associated rigid point} of $x$.}
\end{dfn}

In the sequel of this subsection, $\mathscr{X}$ denotes a coherent rigid space.
Let $\alpha\colon\Spf V\rightarrow\ZR{\mathscr{X}}$ be a rigid point of $\mathscr{X}$.
It defines the point $\alpha(\m_V)$ in $\ZR{\mathscr{X}}$.
Thus we get a mapping
$$
\left\{
\begin{minipage}{9em}
{\small isomorphism classes of rigid points of $\ZR{\mathscr{X}}$}
\end{minipage}
\right\}\longrightarrow\ZR{\mathscr{X}},\qquad \alpha\longmapsto\alpha(\m_V).\leqno{(\ast)}
$$
If $X$ is a formal model of $\mathscr{X}$, then we have another mapping 
$$
\left\{
\begin{minipage}{9em}
{\small isomorphism classes of rigid points of $\ZR{\mathscr{X}}$}
\end{minipage}
\right\}\longrightarrow\left\{
\begin{minipage}{8.3em}
{\small isomorphism classes of rigid points of $X$}
\end{minipage}
\right\},\qquad\alpha\longmapsto\sp_X\circ\alpha.\leqno{(\ast\ast)}
$$
Note that $\sp_X\circ\alpha$ is a morphism of formal schemes, for it is a morphism of locally ringed spaces, and is adic (hence is continuous).

\begin{prop}\label{prop-ZRpoints3}
{\rm (1)} The mapping $(\ast)$ is surjective.
Define an equivalence relation $\approx$ on the left-hand set generated by the relation $\sim$ given as follows$:$ for rigid points $\alpha\colon\Spf V\rightarrow\ZR{\mathscr{X}}$ and $\beta\colon\Spf W\rightarrow\ZR{\mathscr{X}}$, $\alpha\sim\beta$ if there exists an injective map $f\colon V\hookrightarrow W$ such that $\alpha\circ \Spf f=\beta$ and $W$ dominates $V$.
Then $(\ast)$ induces a bijection 
$$
\left\{
\begin{minipage}{8.3em}
{\small {\rm isomorphism classes of rigid points of $\ZR{\mathscr{X}}$}}
\end{minipage}
\right\}\big/_{{\textstyle \approx}}\longrightarrow\ZR{\mathscr{X}}.\leqno{(\ast)'}
$$

{\rm (2)} The mapping $(\ast\ast)$ is bijective.
\end{prop}

Note that in (1) the map $\Spf f\colon\Spf W\rightarrow\Spf V$ is automatically adic and that $V=\Frac(V)\cap W$ in $\Frac(W)$ (cf.\ {\rm (a)} in {\bf \ref{ch-pre}}.\ref{dfn-val}).
Notice also that, if $\alpha\sim\beta$, then $\alpha(\m_V)=\alpha(\m_W)$; thus the map $(\ast)'$ as above is well-defined.

For the proof of \ref{prop-ZRpoints3} we need the following lemma:
\begin{lem}\label{lem-ZRpoints3}
Let $\pi\colon X'\rightarrow X$ be an admissible blow-up, and $\alpha\colon\Spf V\rightarrow X$ a rigid point.
Then there exists a unique rigid point $\alpha'\colon\Spf V\rightarrow X'$ such that the resulting diagram commutes$:$
$$
\xymatrix@-1ex{
&X'\ar[d]^{\pi}\\
\Spf V\ar@{-->}[ur]^{\alpha'}\ar[r]_{\alpha}&X\rlap{.}}
$$
\end{lem}

\begin{proof}
Let $\mathscr{J}$ be a admissible ideal of $X$ that gives the admissible blow-up $\pi$.
Since $\mathscr{J}$ is of finite type, the ideal $\mathscr{J}V$ of $V$ is invertible ({\bf \ref{ch-pre}}.\ref{prop-val1} (2)).
Hence there exists a unique factoring map $\alpha'$ as above.
By {\bf \ref{ch-formal}}.\ref{prop-adicmor1} (1) the morphism $\alpha'$ is adic.
\end{proof}

\begin{proof}[Proof of Proposition {\rm \ref{prop-ZRpoints3}}]
(1) Let $x\in\ZR{\mathscr{X}}$ be a point.
We use the notation as in \ref{ntn-ZRpoints}.
By \ref{prop-ZRpoint32} (1) we have a mapping $x\mapsto\alpha_x$ from the set $\ZR{\mathscr{X}}$ to the set of all isomorphism classes of rigid points of $\mathscr{X}$.
Clearly, for $x\in\ZR{\mathscr{X}}$ we have $\alpha_x(\m_{\widehat{V}_x})=x$.
Hence $(\ast)$ is surjective.
To show the other half of (1), it suffices to check the following: for a rigid point $\alpha\colon\Spf V\rightarrow\ZR{\mathscr{X}}$ such that $\alpha(\m_V)=x$, we have $\alpha_x\sim\alpha$.
But this follows from \ref{prop-ZRpoint32} (2).

(2) Let $\alpha\colon\Spf V\rightarrow X$ be a rigid point.
Since $\ZR{\mathscr{X}}$ is the projective limit of all admissible blow-ups of $X$, by \ref{lem-ZRpoints3} we have a unique adic map 
$$
\ZR{\alpha}\colon\Spf V\longrightarrow\ZR{\mathscr{X}}
$$
such that $\sp_X\circ\ZR{\alpha}=\alpha$.
By this we get a map from the set of all isomorphism classes of rigid points of $X$ to the set of all isomorphism classes of rigid points of $\ZR{\mathscr{X}}$.
By the relation $\sp_X\circ\ZR{\alpha}=\alpha$ and the uniqueness of the lifting, one sees that this map gives the inverse to $(\ast\ast)$.
\end{proof}
\index{point!rigid point@rigid ---|)}\index{rigid point|)}

\subsubsection{Seminorms associated to points}\label{subsub-seminormsofpoints}
\index{seminorm|(}
Let $\mathscr{X}$ be a coherent rigid space, and $\mathscr{I}\subseteq\O^{\int}_{\mathscr{X}}$ an ideal of definition of finite type (\ref{dfn-ZRstrsheaf2}).
We have the valuation ring $V_x$ at $x$ for the field $K_x$ as in \ref{ntn-ZRpoints}.
Since $V_x$ is $a$-adically separated (where $\mathscr{I}A_x=(a)$), we have the associated height one prime\index{associated height one prime} $\mathfrak{p}=\sqrt{(a)}$ ({\bf \ref{ch-pre}}.\ref{dfn-maxspe2}), and hence we have the corresponding height one valuation on the field $K_x$ with the valuation ring $V_{\mathfrak{p}}$.
We can then define, choosing a once for all fixed real number $0<c<1$, the corresponding (non-archimedean) norm\index{norm!non-archimedean norm@non-archimedean ---} ({\bf \ref{ch-pre}}, \S\ref{subsub-nonarchnorms}), denoted by
$$
\|\cdot\|_{x,\mathscr{I},c}\colon K_x\longrightarrow\R_{\geq 0}, 
$$
uniquely determined by $\|a\|_{x,\mathscr{I},c}=c$.
It is clear that another choice of $\mathscr{I}$ and $c$ only leads to an equivalent norm. 

For $f\in\Gamma(\mathscr{X},\O_{\mathscr{X}})$ and $x\in\ZR{\mathscr{X}}$ we denote by $f(x)$ the image of $f$ by the composite map
$$
\Gamma(\mathscr{X},\O_{\mathscr{X}})\longrightarrow B_x\ (=\O_{\mathscr{X},x})\longrightarrow K_x\ (=B_x/\m_{B_x}),
$$
and write
$$
\|f(x)\|_{\mathscr{I},c}=\|f(x)\|_{x,\mathscr{I},c}.
$$
This construction yields a mapping (denoted also by $\|\cdot\|_{x,\mathscr{I},c}$)
$$
\|\cdot\|_{x,\mathscr{I},c}\colon\Gamma(\mathscr{X},\O_{\mathscr{X}})\longrightarrow\R_{\geq 0},\qquad f\longmapsto\|f(x)\|_{\mathscr{I},c},
$$
which is a {\em multiplicative seminorm}\index{seminorm!multiplicative seminorm@multiplicative ---} on the ring $\Gamma(\mathscr{X},\O_{\mathscr{X}})$, that is,
\begin{itemize}
\item $\|0\|_{x,\mathscr{I},c}=0$; $\|-f\|_{x,\mathscr{I},c}=\|f\|_{x,\mathscr{I},c}$;
\item $\|f+g\|_{x,\mathscr{I},c}\leq\max\{\|f\|_{x,\mathscr{I},c},\|g\|_{x,\mathscr{I},c}\}$;
\item $\|f\cdot g\|_{x,\mathscr{I},c}=\|f\|_{x,\mathscr{I},c}\cdot\|g\|_{x,\mathscr{I},c}$;
\item $\|1\|_{x,\mathscr{I},c}=1$.
\end{itemize}

\begin{prop}\label{prop-lemnormsatpoints}
Let $x,x'\in\ZR{\mathscr{X}}$, and suppose that $x'$ is a generization of $x$.
Then we have the equality of the seminorms on $\Gamma(\mathscr{X},\O_{\mathscr{X}}):$
$$
\|\cdot\|_{x,\mathscr{I},c}=\|\cdot\|_{x',\mathscr{I},c}.
$$
\end{prop}

\begin{proof}
In view of \ref{prop-generizationmapsrigidspaces} we have the generization map $K_x\hookrightarrow K_{x'}$ that maps $V_x$ into $V_{x'}$.
Consider the associated height one prime $\mathfrak{p}'=\sqrt{aV_{x'}}$ of $V_{x'}$; in view of {\bf \ref{ch-pre}}.\ref{prop-maxspe4} the prime $\mathfrak{p}=\mathfrak{p}'\cap V_x$ is equal to $\sqrt{aV_x}$ and hence is the associated height one prime of $V_x$.
In particular, we have the local injective morphism $V_{x,\mathfrak{p}}\hookrightarrow V_{x',\mathfrak{p}'}$ and hence the injective mapping $K^{\times}_x/V^{\times}_{x,\mathfrak{p}}\hookrightarrow K^{\times}_{x'}/V^{\times}_{x',\mathfrak{p}'}$ between the value groups.
Since this is an ordered mapping between the totally ordered sets isomorphic to $\R$ mapping the class of $a$ to the class of $a$, it is an ordered isomorphism.
Then the desired equality follows from this.
\end{proof}

\subsubsection{Spectral seminorms}\label{subsub-spectralseminorms}
\index{seminorm!spectral seminorm@spectral ---|(}
Let $X$ and $\mathscr{I}$ as in \S\ref{subsub-seminormsofpoints}.
For $f\in\Gamma(\mathscr{X},\O_{\mathscr{X}})$ we set
$$
\|f\|_{\Sp,\mathscr{I},c}=\sup_{x\in\ZR{\mathscr{X}}}\|f(x)\|_{\mathscr{I},c},
$$
which defines the map
$$
\|\cdot\|_{\Sp,\mathscr{I},c}\colon\Gamma(\mathscr{X},\O_{\mathscr{X}})\longrightarrow\R_{\geq 0}.
$$
It is easy to see that $\|\cdot\|_{\Sp,\mathscr{I},c}$ gives a seminorm on the ring $\Gamma(\mathscr{X},\O_{\mathscr{X}})$, that is, 
\begin{itemize}
\item $\|0\|_{\Sp,\mathscr{I},c}=0$; $\|-f\|_{\Sp,\mathscr{I},c}=\|f\|_{\Sp,\mathscr{I},c}$;
\item $\|f+g\|_{\Sp,\mathscr{I},c}\leq\max\{\|f\|_{\Sp,\mathscr{I},c},\|g\|_{\Sp,\mathscr{I},c}\}$;
\item $\|f\cdot g\|_{\Sp,\mathscr{I},c}\leq \|f\|_{\Sp,\mathscr{I},c}\cdot\|g\|_{\Sp,\mathscr{I},c}$;
\item $\|1\|_{\Sp,\mathscr{I},c}=1$.
\end{itemize}
We called this the {\em spectral seminorm} on $\mathscr{X}$.
Clearly, changing $\mathscr{I}$ and $c$ only leads to an equivalent seminorm .
\index{seminorm!spectral seminorm@spectral ---|)}
\index{seminorm|)}

\subsection{Comparison of topologies}\label{sub-comparisontopoi}
\begin{prop}\label{prop-thmadmissiblesite200}
Let $\{j_{\alpha}\colon\mathscr{U}_{\alpha}\hookrightarrow\mathscr{F}\}_{\alpha\in L}$ be a collection of open immersions of rigid spaces.
Then the following conditions are equivalent$:$
\begin{itemize}
\item[{\rm (a)}] $\{j_{\alpha}\colon\mathscr{U}_{\alpha}\hookrightarrow\mathscr{F}\}_{\alpha\in L}$ is a covering family in the site $\mathscr{F}_{\ad}$ or, equivalently, in $\Rf_{\ad}$ $($that is, $\coprod_{\alpha\in L}\mathscr{U}_{\alpha}\rightarrow\mathscr{F}$ is an epimorphism of sheaves on $\CRf_{\ad});$
\item[{\rm (b)}] $\{\ZR{j_{\alpha}}\colon\ZR{\mathscr{U}_{\alpha}}\hookrightarrow\ZR{\mathscr{F}}\}_{\alpha\in L}$ is a covering of the topological space $\ZR{\mathscr{F}}$ $($that is, $\ZR{\mathscr{F}}=\bigcup_{\alpha\in L}\ZR{\iota_{\alpha}}(\ZR{\mathscr{U}_{\alpha}}))$.
\end{itemize}
\end{prop}

For the proof we need:
\begin{lem}\label{lem-admissiblesite2zr}
Let $j\colon\mathscr{U}\hookrightarrow\mathscr{X}$ be a coherent open immersion of coherent rigid spaces.
Then the following conditions are equivalent.
\begin{itemize}
\item[{\rm (a)}] $j$ is an isomorphism.
\item[{\rm (b)}] The singleton set $\{j\colon\mathscr{U}\hookrightarrow\mathscr{X}\}$ is a covering in the coherent small site $\mathscr{X}_{\ad}$.
\item[{\rm (c)}] The singleton set $\{j\colon\mathscr{U}\hookrightarrow\mathscr{X}\}$ is a covering in the large admissible site $\Rf_{\ad}$.
\item[{\rm (d)}] $\ZR{j}$ is an isomorphism.
\end{itemize}
\end{lem}

\begin{proof}
The equivalence of (a) and (b) is obvious.
The implication (a) $\Rightarrow$ (c) is clear.
Suppose (c) holds and $\ZR{j}$ is not an isomorphism.
Since $\ZR{j}\colon\ZR{\mathscr{U}}\hookrightarrow\ZR{\mathscr{X}}$ is an open immersion of locally ringed spaces, this means that $\ZR{j}$ is not surjective, and hence there exists $x\in\ZR{\mathscr{X}}$ not lying in the image of $\ZR{j}$.
Consider the associated rigid point $\alpha_x\colon(\Spf \widehat{V_x})^{\rig}\rightarrow\mathscr{X}$.
By definition of the coverings in the site $\Rf_{\ad}$ (\ref{dfn-admissiblesite3genlarge}), this should come from a map $(\Spf \widehat{V_x})^{\rig}\rightarrow\mathscr{U}$, which contradicts what we have supposed, since this implies that $x$ lies in the image of $\ZR{j}$.
Finally, we show (d) $\Rightarrow$ (a).
Suppose $\ZR{j}$ is an isomorphism.
In view of \ref{prop-cohrigidspaceopenimm1} there exists a distinguished formal model $U\hookrightarrow X$ of $j$.
It suffices to show that $U=X$.
Suppose $X\setminus U$ is non-empty, and take a point $x\in X\setminus U$.
By \ref{prop-ZRpoints4} there exists a rigid point $\alpha\colon \Spf V\rightarrow X$ such that $\alpha(\m_V)=x$.
In view of \ref{prop-ZRpoints3} (2) there exists a unique lift $\til{\alpha}\colon\Spf V\rightarrow\ZR{\mathscr{X}}$ of $\alpha$.
By \ref{prop-zariskiriemanntoptop} the point $\alpha(\m_V)$ lies outside $\ZR{\mathscr{U}}$ in $\ZR{\mathscr{X}}$, which is absurd.
\end{proof}

\begin{proof}[Proof of Proposition {\rm \ref{prop-thmadmissiblesite200}}]
(a) $\Rightarrow$ (b) is due to \ref{prop-zariskiriemanntop2gen}.
Let us prove (b) $\Rightarrow$ (a).
If $\{j_{\alpha}\colon\mathscr{U}_{\alpha}\hookrightarrow\mathscr{F}\}_{\alpha\in L}$ is not a covering, then there exists a coherent open rigid subspace $\mathscr{V}$ of $\mathscr{F}$ such that $\{j_{\alpha}\colon\mathscr{U}_{\alpha}\times_{\mathscr{F}}\mathscr{V}\hookrightarrow\mathscr{V}\}_{\alpha\in L}$ is not a covering.
Hence, to show (a), we may assume that $\mathscr{F}$ is coherent.
In this case, since $\ZR{\mathscr{F}}$ is quasi-compact, it is covered by only finitely many of $\ZR{\mathscr{U}_{\alpha}}$.
Hence we can assume that $L$ is a finite set.
Moreover, since $\ZR{\mathscr{F}}$ has an open basis consisting of quasi-compact open subsets, we may assume that each $\mathscr{U}_{\alpha}$ is a coherent rigid space.
Let $\mathscr{U}$ be the rigid space defined by the quotient of sheaves
$$
\xymatrix@-1ex{
\coprod_{\alpha,\beta\in L}\mathscr{U}_{\alpha}\times_{\mathscr{F}}\mathscr{U}_{\beta}\ar@<.5ex>[r]\ar@<-.5ex>[r]&\coprod_{\alpha\in L}\mathscr{U}_{\alpha}\ar[r]&\mathscr{U}.}
$$
Due to \ref{cor-cohrigidspacepatching11} and the exactness of $\CRf\rightarrow\Rf^{\sim}_{\ad}$, $\mathscr{U}$ is represented by a coherent rigid space, and the canonical morphism $j\colon\mathscr{U}\rightarrow\mathscr{F}$ is a coherent open immersion (\ref{lem-cohrigidspacepatchingx}).
Then the assertion follows from \ref{prop-zariskiriemanntop2} and \ref{lem-admissiblesite2zr}.
\end{proof}

\begin{cor}\label{cor-admissiblesite200c1}
Let $\mathscr{X}$ be a coherent rigid space, and $\{\mathscr{U}_{\alpha}\hookrightarrow\mathscr{X}\}_{\alpha\in L}$ a covering family in the small admissible site $\mathscr{X}_{\ad}$ $(\ref{dfn-admissiblesite3gensmall})$, where $\mathscr{X}$ is regarded as a general rigid space.
Then the covering family $\{\mathscr{U}_{\alpha}\hookrightarrow\mathscr{X}\}_{\alpha\in L}$ can be refined by a covering family in the coherent small admissible site {\rm (\ref{dfn-admissiblecovering11})}.
In particular, there exists a finite subset $L'\subseteq L$ such that $\{\mathscr{U}_{\alpha}\hookrightarrow\mathscr{X}\}_{\alpha\in L'}$ gives a covering in the site $\mathscr{X}_{\ad}$.
\end{cor}

\begin{proof}
Considering a covering of each $\mathscr{U}_{\alpha}$ by coherent rigid spaces, we may assume that each $\mathscr{U}_{\alpha}$ is coherent.
Then by \ref{prop-thmadmissiblesite200} and by the fact that $\ZR{\mathscr{X}}$ is quasi-compact due to \ref{thm-ZRcompact} (1), we have a finite subset $L'\subseteq L$ such that $\{\mathscr{U}_{\alpha}\hookrightarrow\mathscr{X}\}_{\alpha\in L'}$ already gives a covering.
By \ref{lem-admissiblesite2zr} we conclude that this gives a covering in the coherent small admissible site, as desired.
\end{proof}

Let $\mathscr{F}$ be a rigid space represented by a coherent rigid space $\mathscr{X}$.
Let $\mathscr{F}_{\ad}$ denote the small admissible site as in \ref{dfn-admissiblesite3gensmall}, and $\mathscr{X}_{\ad}$ the coherent small admissible site as in \ref{dfn-admissiblecovering11}.
There exists a canonical comparison morphism of sites 
$$
\mathscr{X}_{\ad}\longrightarrow\mathscr{F}_{\ad}\leqno{(\ast)}
$$ 
obtained as follows.
Consider the obvious functor of categories $\mathscr{X}_{\ad}\rightarrow\mathscr{F}_{\ad}$ that maps $\mathscr{U}\hookrightarrow\mathscr{X}$ to the associated morphism of the representable functors.
To see that the functor in question is a morphism of sites, it suffices to show that it maps coverings to coverings, since it is obvious that coherent open rigid subspace of $\mathscr{X}$ and coherent open immersions generates the site $\mathscr{F}_{\ad}$.
But this follows from \ref{prop-thmadmissiblesite200}, since a covering in the coherent site $\mathscr{X}_{\ad}$ induces an open covering by passage to the associated Zariski-Riemann spaces (\ref{prop-zariskiriemanntop2}).
Since by \ref{cor-admissiblesite200c1} the topology on $\mathscr{F}_{\ad}$ is generated by the topology on $\mathscr{X}_{\ad}$, we have:
\begin{thm}\label{thm-generalrigidspace31-32}
The morphism of sites $(\ast)$ gives rise to an equivalence of topoi 
$$
\mathscr{X}^{\sim}_{\ad}\stackrel{\sim}{\longrightarrow}\mathscr{F}^{\sim}_{\ad}. \eqno{\square}
$$
\end{thm}

Let $\mathscr{F}$ be a rigid space.
Then for any open immersion $\mathscr{U}\hookrightarrow\mathscr{F}$ the induced map $\ZR{\mathscr{U}}\rightarrow\ZR{\mathscr{F}}$ is an open immersion of topological spaces (\ref{prop-zariskiriemanntop2gen}).
Hence one has a natural functor 
$$
\mathscr{F}_{\ad}\longrightarrow\Ouv(\ZR{\mathscr{F}}),\leqno{(\ast\ast)}
$$
where the right-hand category is the category of open subsets of $\ZR{\mathscr{F}}$.
We claim that there exists a canonical morphism of sites 
$$
\Ouv(\ZR{\mathscr{F}})\longrightarrow\mathscr{F}_{\ad}\leqno{(\dagger)}
$$
underlain by the functor $(\ast\ast)$.
By \ref{prop-zariskiriemanntop2gen} the above functor maps all covering families to open coverings of $\ZR{\mathscr{F}}$.
By \ref{prop-zariskiriemanntopz} the topological space $\ZR{\mathscr{F}}$ has a generator consisting of quasi-compact open subsets coming from coherent open rigid subspaces of $\mathscr{F}$.
Hence we deduce that $(\ast\ast)$ induces a morphism of sites as above.

\begin{thm}\label{thm-generalrigidspace31-31}
Let $\mathscr{F}$ be a rigid space, and $\ZR{\mathscr{F}}$ the associated Zariski-Riemann space.
Then the morphism $(\dagger)$ of sites induces an equivalence of topoi
$$
\top(\ZR{\mathscr{F}})\stackrel{\sim}{\longrightarrow}\mathscr{F}^{\sim}_{\ad},
$$
where $\top(\cdot)$ is the functor giving the associated topos.
\end{thm}

\begin{proof}
We may assume that $\mathscr{F}$ is represented by a coherent rigid space.
Then the assertion follows from \ref{prop-zariskiriemanntoptop} and \ref{prop-thmadmissiblesite200}.
\end{proof}

\begin{thm}\label{prop-admissiblesite3genconsistencyx}
The canonical functor $\CRf\hookrightarrow\Rf$ gives rise to a morphism of sites $\CRf_{\ad}\rightarrow\Rf_{\ad};$ the site $\Rf_{\ad}$ is generated by the objects in the image of $\CRf\hookrightarrow\Rf$.
In particular, the topos $\Rf^{\sim}_{\ad}$ $($and, similarly, $\Rf^{\sim}_{\mathscr{G},\ad}$ for a rigid space $\mathscr{G})$ is generated by quasi-compact objects.
\end{thm}

\begin{proof}
It is clear that the functor $\CRf\hookrightarrow\Rf$ gives rise to a morphism of sites $\CRf_{\ad}\rightarrow\Rf_{\ad}$.
By \ref{cor-admissiblesite200c1} we see that the site $\Rf_{\ad}$ is generated by the objects in the image of $\CRf\hookrightarrow\Rf$.
Finally, by the second assertion of \ref{cor-admissiblesite200c1} the image of $\CRf\hookrightarrow\Rf$ consists of quasi-compact objects.
\end{proof}

\subsection{Finiteness conditions and consistency of terminologies}\label{sub-consistencyterm}
\subsubsection{Finiteness conditions}\label{subsub-consistencyterm}
\begin{dfn}\label{dfn-generalrigidspace2}{\rm 
(1) A rigid space $\mathscr{F}$ is said to be {\em quasi-compact}\index{rigid space!quasi-compact rigid space@quasi-compact ---} if it is quasi-compact as an object of the topos $\CRf^{\sim}_{\ad}$ (cf.\ {\bf \ref{ch-pre}}.\ref{dfn-coherenttopos1} (1)).

(2) A rigid space $\mathscr{F}$ is said to be {\em quasi-separated}\index{rigid space!quasi-separated rigid space@quasi-separated ---} if it is quasi-separated as an object of the topos $\CRf^{\sim}_{\ad}$, that is, the diagonal map $\mathscr{F}\rightarrow\mathscr{F}\times\mathscr{F}$ of sheaves is quasi-compact (cf.\ {\bf \ref{ch-pre}}.\ref{dfn-coherenttopos1} (2)).}
\end{dfn}

Notice that, if a rigid space $\mathscr{F}$ is quasi-compact, then one can take the covering $\mathscr{Y}=\coprod_{\alpha\in L}\mathscr{Y}_{\alpha}$ as in \ref{dfn-generalrigidspace1} (a) with the index set $L$ finite.
\begin{prop}\label{prop-generalrigidspace2}
A rigid space $\mathscr{F}$ is quasi-compact and quasi-separated if and only if it is represented by a coherent rigid space.
\end{prop}

This proposition allows us to say consistently that a rigid space $\mathscr{F}$ is coherent\index{rigid space!coherent rigid space@coherent ---} if it is quasi-compact and quasi-separated. 
\begin{proof}
Suppose $\mathscr{F}$ is quasi-compact and quasi-separated.
Take a surjective map $\mathscr{Y}=\coprod_{\alpha\in L}\mathscr{Y}_{\alpha}\rightarrow\mathscr{F}$ as in \ref{dfn-generalrigidspace1} where $L$ is finite.
Since $\mathscr{F}$ is quasi-separated, each $\mathscr{Y}_{\alpha}\times_{\mathscr{F}}\mathscr{Y}_{\beta}$ is quasi-compact.
Since it is covered by finitely many stretches of coherent rigid spaces, it is covered by finite stretches of coherent rigid spaces (\ref{dfn-admissiblesite31} (1)) and hence is coherent.
Then one can refine the covering $\mathscr{Y}=\coprod_{\alpha\in L}\mathscr{Y}_{\alpha}\rightarrow\mathscr{F}$ so that each projection $\mathscr{Y}_{\alpha}\times_{\mathscr{F}}\mathscr{Y}_{\beta}\rightarrow\mathscr{Y}_{\alpha}$ is a coherent open immersion.
Hence $\mathscr{F}$ is isomorphic to the coherent rigid space $\mathscr{Y}\amalg_{\mathscr{Y}\times_{\mathscr{F}}\mathscr{Y}}\mathscr{Y}$.
Note that, since the the canonical functor from a site to the associated topos is exact, it preserves finite cofiber products.
\end{proof}

\begin{prop}\label{prop-generalrigidspace31-1}
A rigid space $\mathscr{F}$ is quasi-separated if and only if it is a stretch of coherent rigid spaces {\rm (\ref{dfn-admissiblesite31} (1))}.
\end{prop}

\begin{proof}
The `if' part is clear.
Suppose $\mathscr{F}$ is quasi-separated, and take $\mathscr{Y}=\coprod_{\alpha\in L}\mathscr{Y}_{\alpha}\rightarrow\mathscr{F}$ as in \ref{dfn-generalrigidspace1}.
Let $\mathfrak{L}$ be the set of all finite subsets of $L$ considered with the inclusion order.
Set $\mathscr{Y}_{L'}=\coprod_{\alpha\in L'}\mathscr{Y}_{\alpha}$ for any $L'\in\mathfrak{L}$, and let $\mathscr{X}_{L'}$ be the image of $\mathscr{Y}_{L'}$.
Then $\mathscr{X}_{L'}$ is quasi-compact and quasi-separated and hence is represented by a coherent rigid space by \ref{prop-generalrigidspace2}.
Since $\mathscr{F}=\varinjlim_{L'\in\mathfrak{L}}\mathscr{X}_{L'}$, we have done.
\end{proof}

\subsubsection{Consistency of open immersions}\label{subsub-consistencytermopenimm}
The following proposition shows that the terminologies `open immersion'\index{immersion!open immersion of rigid spaces@open --- (of rigid spaces)} and `coherent open immersion'\index{immersion!open immersion of rigid spaces@open --- (of rigid spaces)!coherent open immersion of rigid spaces@coherent --- ---} are consistent.

\begin{prop}\label{prop-consistencytermopenimm2}
Let $\mathscr{U}$ and $\mathscr{X}$ be coherent rigid spaces, and $\iota\colon\mathscr{U}\hookrightarrow\mathscr{X}$ an open immersion.
Then $\iota$ is a coherent open immersion.
Moreover, it is coherent {\rm ({\bf \ref{ch-pre}}.\ref{dfn-coherenttopos12}\ (3))} as an arrow in the topos $\Rf^{\sim}_{\ad}$.
\end{prop}

\begin{proof}
One can write $\mathscr{U}$ as the union of increasing sequence of coherent open rigid subspaces $\{\mathscr{U}_i\}$ such that each $\mathscr{U}_i\hookrightarrow\mathscr{X}$ is a coherent open immersion.
Since $\mathscr{U}$ is quasi-compact as an object of $\Rf^{\sim}_{\ad}$ (\ref{cor-admissiblesite200c1}), we find some $i$ such that $\mathscr{U}=\mathscr{U}_i$, and hence $\iota$ is a coherent open immersion.
The last assertion follows from \ref{prop-openimmpullback1} and \ref{prop-generalrigidspace2}.
\end{proof}

\subsubsection{Rigid space as quotient}\label{subsub-consistencytermrigidquot}
The following proposition follows immediately from the definition of rigid spaces.
\begin{prop}\label{prop-rigidspacequotient}
Let $\mathscr{Y}$ be a rigid space, and $\mathscr{R}$ an equivalence relation\index{equivalence relation} in the sheaf $\mathscr{Y}\times\mathscr{Y}$ on the site $\CRf_{\ad}$ with the projection maps
$$
\xymatrix@-1ex{\mathscr{R}\ar@<.5ex>[r]^{q_1}\ar@<-.5ex>[r]_{q_2}&\mathscr{Y}}
$$
such that there exists a covering family $\{\mathscr{V}_{\alpha}\rightarrow\mathscr{R}\}$ such that $q_i\colon \mathscr{V}_{\alpha}\rightarrow\mathscr{Y}$ for each $i=1,2$ and any $\alpha$ is an open immersion.
Then the quotient $\mathscr{F}$ of $\mathscr{Y}$ by $\mathscr{R}$ is a rigid space. \hfill$\square$
\end{prop}

\subsubsection{Consistency of finiteness conditions}\label{subsub-consistencytermfincond}
\begin{prop}\label{prop-thmadmissiblesite3genalgebraic}
Let $\mathscr{F}$ be a rigid space.
Then the topos $\mathscr{F}^{\sim}_{\ad}$ is algebraic {\rm ({\bf \ref{ch-pre}}.\ref{dfn-coherenttopos2})}.
Moreover, the following conditions are equivalent$:$
\begin{itemize}
\item[{\rm (a)}] $\mathscr{F}$ is quasi-separated $($resp.\ coherent$);$
\item[{\rm (b)}] $\mathscr{F}^{\sim}_{\ad}$ is quasi-separated $($resp.\ coherent$)$.
\end{itemize}
\end{prop}

\begin{proof}
Let $C$ be the set of all objects in $\mathscr{F}^{\sim}_{\ad}$ represented by coherent rigid spaces.
Then $C$ generates the topos $\mathscr{F}^{\sim}_{\ad}$ (by \ref{cor-admissiblesite200c1}).
Clearly, $C$ is stable under fiber products and equalizers. 
Hence $\mathscr{F}^{\sim}_{\ad}$ is algebraic.
Suppose $\mathscr{F}$ is quasi-separated.
Then $C$ is further stable under products (over $\mathscr{F}$).
Hence $\mathscr{F}^{\sim}_{\ad}$ is quasi-separated.
If $\mathscr{F}$ is coherent, then $C$ is stable under all finite colimits, and hence $\mathscr{F}^{\sim}_{\ad}$ is coherent.
Thus we have (a) $\Rightarrow$ (b).
If, conversely, the topos $\mathscr{F}^{\sim}_{\ad}$ is quasi-separated, then $\mathscr{F}$ is a colimit of coherent rigid spaces with all transition maps being open immersions.
Hence by a similar reasoning as in \ref{prop-generalrigidspace31-1} one sees that $\mathscr{F}$ is a stretch of coherent rigid spaces and hence is quasi-separated.
If $\mathscr{F}^{\sim}_{\ad}$ is coherent, then such a colimit is taken to be finite, and hence $\mathscr{F}$ is coherent.
\end{proof}

\begin{prop}\label{prop-thmgeneralrigidspace31-2}
Let $\mathscr{F}$ be a rigid space, and $\ZR{\mathscr{F}}$ the associated Zariski-Riemann space.
Then the following conditions are equivalent$:$
\begin{itemize}
\item[{\rm (a)}] $\mathscr{F}$ is quasi-compact $($resp.\ quasi-separated, resp.\ coherent$);$
\item[{\rm (b)}] $\ZR{\mathscr{F}}$ is quasi-compact $($resp.\ quasi-separated, resp.\ coherent$)$ {\rm (}cf.\ {\rm {\bf \ref{ch-pre}}.\ref{dfn-quasicompactness}, {\bf \ref{ch-pre}}.\ref{dfn-quasiseparatedness}, {\bf \ref{ch-pre}}.\ref{dfn-quasicompact1})}.
\end{itemize}
\end{prop}

\begin{proof}
The implication (a) $\Rightarrow$ (b) is clear by the construction of $\ZR{\mathscr{F}}$ and \ref{thm-ZRcompact} (1).
If $\ZR{\mathscr{F}}$ is quasi-compact, then $\mathscr{F}$ is quasi-compact by \ref{prop-thmadmissiblesite200}.
If $\ZR{\mathscr{F}}$ is quasi-separated (resp.\ coherent), then the topos $\ZR{\mathscr{F}}^{\sim}$ is quasi-separated (resp.\ coherent), and hence so is $\mathscr{F}^{\sim}_{\ad}$ by \ref{thm-generalrigidspace31-31}.
But then by \ref{prop-thmadmissiblesite3genalgebraic} we deduce that $\mathscr{F}$ is quasi-separated (resp.\ coherent).
\end{proof}

\subsubsection{Rigid spaces associated to adic formal schemes}\label{subsub-extensionofrigfuntor}
Let $X$ be an adic formal scheme of finite ideal type, which is not necessarily coherent.
We are going to construct the associated rigid space $X^{\rig}$.

(1) Suppose $X$ is coherent. In this case, we just take $X^{\rig}$ as in \S\ref{subsub-cohrigidspace}; the resulting rigid space $X^{\rig}$ is, therefore, coherent.

(2) Suppose $X$ is quasi-separated. In this case, since $X$ is locally coherent, there exists an increasing family $\{U_i\}_{i\in I}$ (where $I$ is a directed set) of coherent open subsets such that $X=\bigcup U_i$.
Then consider $U^{\rig}_i$ for each $i$, and notice that the induced morphism $U^{\rig}_i\rightarrow U^{\rig}_j$ is a coherent open immersion for any $i\leq j$.
Thus we get the desired rigid space $X^{\rig}$ by the union of $U^{\rig}_i$; notice that the rigid space $X$ in this case is a stretch of coherent rigid spaces and hence is quasi-separated (\ref{prop-generalrigidspace31-1}).

(3) In general, let $X=\bigcup_{\alpha\in L}U_{\alpha}$ be an open covering by coherent open subsets.
For each $\alpha\in L$, consider $U^{\rig}_{\alpha}$ as in (1).
Since each intersection $U_{\alpha}\cap U_{\beta}$ is quasi-separated, one can consider the rigid space $(U_{\alpha}\cap U_{\beta})^{\rig}$ as in (2).
Since $(U_{\alpha}\cap U_{\beta})^{\rig}\rightarrow U^{\rig}_{\alpha}$ is obviously an open immersion, one can define $X^{\rig}$ by patching all $U^{\rig}_{\alpha}$ along $(U_{\alpha}\cap U_{\beta})^{\rig}$.
\index{visualization|)}

\addcontentsline{toc}{subsection}{Exercises}
\subsection*{Exercises}
\begin{exer}[Deligne's\index{Deligne, P.} formula]\label{exer-deligne}{\rm 
Let $(X,Z)$ be a universally pseudo-adhesive pair\index{pair!adhesive pair@adhesive ---!pseudo adhesive pair@pseudo-{---} ---}\index{adhesive!adhesive pair@--- pair!pseudo adhesive pair@pseudo-{---} ---} of schemes such that $X$ is coherent\index{scheme!coherent scheme@coherent ---}\index{coherent!coherent scheme@--- scheme}.
Denote the defining ideal of $Z$ by $\mathscr{I}_Z$.
Set $U=X\setminus Z$, and let $j\colon U\hookrightarrow X$ be the canonical open immersion.
Then for any finitely presented $\O_X$-module $\mathscr{F}$ and a quasi-coherent $\O_X$-module $\mathscr{G}$, the canonical map 
$$
\varinjlim_{n\geq 0}\Hom_{\O_X}(\mathscr{I}^n_Z\mathscr{F},\mathscr{G})\longrightarrow\Hom_{\O_U}(\mathscr{F}|_U,\mathscr{G}|_U)
$$
is bijective.
In particular, we have $\varinjlim_{n\geq 0}\Hom(\mathscr{I}^n_Z,\mathscr{G})\stackrel{\sim}{\rightarrow}\Gamma(U,\mathscr{G})$.}
\end{exer}

\begin{exer}\label{exer-propZRstrsheaf2}
{\rm Let $\mathscr{X}$ be a quasi-compact rigid space, and $\mathscr{I},\mathscr{I}'\subseteq\O^{\int}_{\mathscr{X}}$ ideals of definition of finite type.
Then show that there exist integers $n,m>0$ such that $\mathscr{I}^m\subseteq\mathscr{I}^{\prime n}\subseteq\mathscr{I}$.}
\end{exer}

\begin{exer}\label{exer-exavaluationrigidspace}{\rm 
Let $T=\Spf V$ where $V$ is an $a$-adically complete valuation ring with $a\in\m_V\setminus\{0\}$, and consider the associated coherent rigid space $\mathscr{T}=T^{\rig}$.
Show that the topological space $\ZR{\mathscr{T}}$ is homeomorphic to the underlying topological space of $T$ and, moreover, that$:$ 
\begin{itemize}
\item[(a)] the integral structure sheaf $\O^{\int}_{\mathscr{T}}$ is isomorphic to $\O_T$; 
\item[(b)] the rigid structure sheaf $\O_{\mathscr{T}}$ is isomorphic to $\O_T[\frac{1}{a}]$.
\end{itemize}}
\end{exer}

\begin{exer}\label{exer-proptypeN}
Show that the following conditions for a coherent rigid space $\mathscr{X}$ are equivalent$:$
\begin{itemize}
\item[{\rm (a)}] $\mathscr{X}$ has a Noetherian formal model$;$
\item[{\rm (b)}] there exists a cofinal family $\{X_{\lambda}\}_{\lambda\in\Lambda}$ of formal models of $\mathscr{X}$ consisting of Noetherian formal schemes\index{formal scheme!Noetherian formal scheme@Noetherian ---!locally Noetherian formal scheme@locally --- ---}$;$
\item[{\rm (c)}] any quasi-compact open subspace $\mathscr{U}$ of $\mathscr{X}$ has a Noetherian formal model$;$
\item[{\rm (d)}] there exists an open covering $\{\mathscr{U}_{\alpha}\}_{\alpha\in L}$ of the coherent small admissible site $\mathscr{X}_{\ad}$ such that each $\mathscr{U}_{\alpha}$ has a Noetherian formal model.
\end{itemize}
\end{exer}


\section{Topological properties}\label{sec-topologicalproperties}
In this section we study several topological aspects of rigid spaces.
Most of the topological properties, which capture several geometric structures of rigid spaces, are defined and discussed in terms of the associated Zariski-Riemann spaces and are often easily grasped simply by point-set topology.
Many of the statements in this section will be, therefore, recasts of general topology statements, already discussed in {\bf \ref{ch-pre}}, \S\ref{sec-gentop}.

The main aim of \S\ref{sub-generization} is to show that, for any rigid space $\mathscr{X}$, the associated Zariski-Riemann space $\ZR{\mathscr{X}}$, which we have already shown to be locally coherent and sober (\ref{prop-zariskiriemanntopz}), is a valuative space\index{valuative!valuative topological space@--- (topological) space}\index{space@space (topological)!valuative topological space@valuative ---} ({\bf \ref{ch-pre}}.\ref{dfn-valuativespace}).
The proof of this fundamental result requires a careful study of generizations and specializations of points of adic formal schemes of finite ideal type, especially their behavior on passage to the admissible blow-ups. 
It turns out that, for a given point of the associated Zariski-Riemann space, the corresponding rigid point carries all generizations of the point.

In \S\ref{sub-separation} we study the separated quotient of the associated Zariski-Riemann spaces (cf.\ {\bf \ref{ch-pre}}, \S\ref{subsub-separationgen}).
It will turn out that the separated quotient is, as a set, the subset of the Zariski-Riemann space consisting of the points of height $1$.
The nice aspect of the separated quotients lies mainly in their topological feature; for example, the separated quotient of the Zariski-Riemann space associated to a coherent rigid space is a compact\footnote{Recall that, as mentioned in the introduction, all compact topological spaces are assume to be Hausdorff.} space. 

\subsection{Generization and specialization}\label{sub-generization}
\begin{prop}\label{prop-ZRpoints51}
Let $\mathscr{X}=X^{\rig}$ be a coherent rigid space, and $x,y\in\ZR{\mathscr{X}}$.
Then $y$ is a generization\index{generization} {\rm ({\bf \ref{ch-pre}}, \S\ref{subsub-genspetopsp})} of $x$ if and only if for any admissible blow-up\index{blow-up!admissible blow-up@admissible ---} $X'\rightarrow X$, $\sp_{X'}(y)$ is a generization of $\sp_{X'}(x)$ in $X'$.
\end{prop}

To show this, we need the following preparatory lemma:
\begin{lem}\label{lem-ZRpoints50}
Let $X$ be a coherent adic formal scheme of finite ideal type, and $\mathscr{X}=X^{\rig}$.
Let $F\subseteq\ZR{\mathscr{X}}$ be a subset.

{\rm (1)} For any $X'\in\obj(\BL_X)$ we have
$$
\sp_{X'}(\ovl{F})=\ovl{\sp_{X'}(F)},
$$
where $\ovl{F}$ is the closure of $F$ in $\ZR{\mathscr{X}}$, and $\ovl{\sp_{X'}(F)}$ is the closure of $\sp_{X'}(F)$ in $X'$.

{\rm (2)} We have
$$
\ovl{F}=\bigcap_{X'\in\obj(\BL_X)}\sp^{-1}_{X'}(\ovl{\sp_{X'}(F)})=\varprojlim_{X'\in\obj(\BL_X)}\ovl{\sp_{X'}(F)}.
$$
\end{lem}

\begin{proof}
Since each $\sp_{X'}$ is a closed map (\ref{thm-ZRcompact} (2)), we have (1).
Then (2) follows from {\bf \ref{ch-pre}}.\ref{lem-projlimclosedmapslem} (2).
\end{proof}

\begin{proof}[Proof of Proposition {\rm \ref{prop-ZRpoints51}}]
Suppose $y$ is a generization of $x$, and set $C=\ovl{\{y\}}$.
By \ref{lem-ZRpoints50} we have $\sp_{X'}(C)=\ovl{\{\sp_{X'}(y)\}}$.
Hence $\sp_{X'}(x)$ is a specialization of $\sp_{X'}(y)$.
Conversely, suppose $\sp_{X'}(y)$ is a generization of $\sp_{X'}(x)$ for any $X'$, and set $C=\ovl{\{y\}}$.
By \ref{lem-ZRpoints50} (2) we have
$$
C=\varprojlim_{X'\rightarrow X}\ovl{\{\sp_{X'}(y)\}}.
$$
Hence $x\in C$, as desired.
\end{proof}

\begin{prop}\label{prop-ZRpoints52}
Let $\mathscr{X}$ be a coherent rigid space, and $x\in\ZR{\mathscr{X}}$.
Let $\alpha_x\colon\Spf\widehat{V}_x\rightarrow\ZR{\mathscr{X}}$ be the associated rigid point\index{point!rigid point@rigid ---}\index{rigid point} of $x$ {\rm (\ref{dfn-ZRpoints32})}.
Then $\alpha_x$ maps the set $\Spf\widehat{V}_x$ bijectively onto the set of all generizations of $x$.
\end{prop}

\begin{proof}
Let us first show that any point in the image of $\alpha_x$ is a generization of $x$.
Let $y=\alpha_x(\mathfrak{p})$ for $\mathfrak{p}\in\Spf\widehat{V}_x$, and $X$ a formal model of $\mathscr{X}$.
By \ref{prop-ZRpoints51} it suffices to show that for any admissible blow-up $X'\rightarrow X$ the point $\sp_{X'}(y)$ is a generization of $\sp_{X'}(x)$.
But this is clear, since $\mathfrak{p}$ is a generization (in $\Spf\widehat{V}_x$) of $\m_{\widehat{V}_x}$.

Conversely, suppose $y$ is a generization of $x$, and we want to show that there exists a unique $\mathfrak{p}\in\Spf\widehat{V}_x$ such that $y=\alpha_x(\mathfrak{p})$.
For any admissible blow-up $X'\rightarrow X$, find the prime ideal $\mathfrak{p}_{X'}$ of $\O_{X',\sp_{X'}(x)}$ that corresponds to the generization $\sp_{X'}(y)$.
The system $\{\mathfrak{p}_{X'}\}$ defines an ideal (easily seen to be prime) of the inductive limit $A_x=\O^{\int}_{\ZR{\mathscr{X}}}=\varinjlim_{X'\rightarrow X}\O_{X',\sp_{X'}(x)}$; denote this by $\mathfrak{p}$.
Since each $\mathfrak{p}_{X'}$ is open, $\mathfrak{p}$ is an open ideal with respect to the $I_x$-adic topology.
Moreover, it is straightforward to see that $\mathfrak{p}$ contains $J_x$ (in the notation as in \ref{ntn-ZRpoints}); indeed, as one can show easily, the ideal $\mathfrak{p}_{X'}A_x$ contains $J_x$ for any $X'\rightarrow X$.
Hence $\mathfrak{p}$ determines a prime ideal of $V_x$, which we again denote by $\mathfrak{p}$.
Then by {\bf \ref{ch-pre}}.\ref{thm-compval2006ver1} (5) we find that $\alpha_x$ maps $\Spf\widehat{V}_x$ bijectively onto the set of all generizations of $x$.
\end{proof}

By \ref{prop-ZRpoints52}, \ref{prop-ZRpoint32} (2) and {\bf \ref{ch-pre}}.\ref{prop-maxspe4} we immediately have:
\begin{cor}\label{cor-ZRpoints531}
Let $\mathscr{X}$ be a coherent rigid space, and $\alpha\colon\Spf V\rightarrow\ZR{\mathscr{X}}$ a rigid point.
Then $\alpha$ maps the set $\Spf V$ surjectively onto the set of all generizations of $\alpha(\m_V)$. \hfill$\square$
\end{cor}

As in {\bf \ref{ch-pre}}, \S\ref{subsub-genspetopsp}, for any point $x\in\ZR{\mathscr{X}}$ we denote the set of all generizations of $x$ in $\ZR{\mathscr{X}}$ by $G_x$, and consider the ordering on $G_x$ as follows$:$ $y\leq z$ if and only if $z$ is a generization of $y$.

\begin{cor}\label{cor-ZRpoints53}
The set $G_x$ equipped with the above ordering is a totally ordered set.
\end{cor}

\begin{proof}
This follows immediately from the fact that $G_x$ is the bijective image by an order preserving map of $\Spf\widehat{V}_x$ with the ordering by inclusion, and that $\Spf\widehat{V}_x$ is totally ordered (since $\widehat{V}_x$ is a valuation ring\index{valuation!valuation ring@--- ring}; cf.\ {\bf \ref{ch-pre}}.\ref{dfn-val} (c)).
\end{proof}

\begin{dfn}\label{dfn-ZRpoints6}{\rm 
Let $\mathscr{X}$ be a rigid space, and $x\in\ZR{\mathscr{X}}$ a point of the associated Zariski-Riemann space.

{\rm (1)} The {\em height}\index{height!height of a point@--- (of a point)} of the point $x$ is the height of the associated valuation ring $V_x$ (defined as in \ref{ntn-ZRpoints}).
Notice that this definition coincides with the one already given in {\bf \ref{ch-pre}}, \S\ref{sub-valuativespace}.

{\rm (2)} A generization $y$ of the point $x$ is called the {\em maximal generization}\index{generization!maximal generization@maximal ---} if it is the maximal element in the ordered set $G_x$ as in \ref{cor-ZRpoints53}.}
\end{dfn}

The maximal generization of $x$ is, if it exists, uniquely determined, since $G_x$ is totally ordered (\ref{cor-ZRpoints53}).

\begin{prop}\label{prop-ZRpoints6}
Let $\mathscr{X}$ be a coherent rigid space, and $x\in\ZR{\mathscr{X}}$ a point of the associated Zariski-Riemann space.
Then the maximal generization\index{generization!maximal generization@maximal ---} of $x$ exists$;$ moreover, the following conditions for a point $y\in\ZR{\mathscr{X}}$ are equivalent$:$
\begin{itemize}
\item[{\rm (a)}] $y$ is the maximal generization of $x;$
\item[{\rm (b)}] $y$ is a generization of $x$ and is of height one$;$
\item[{\rm (c)}] there exists a rigid point $\alpha\colon\Spf V\rightarrow\ZR{\mathscr{X}}$ such that $\alpha(\m_V)=x$ and $\alpha(\mathfrak{p}_V)=y$, where $\mathfrak{p}_V=\sqrt{(a)}$ is the associated height one prime $({\bf \ref{ch-pre}}.\ref{dfn-maxspe2});$
\item[{\rm (d)}] for any rigid point $\alpha\colon\Spf V\rightarrow\ZR{\mathscr{X}}$ such that $\alpha(\m_V)=x$, we have $\alpha(\mathfrak{p}_V)=y$, where $\mathfrak{p}_V=\sqrt{(a)}$.
\end{itemize}
\end{prop}

\begin{proof}
Take the associated rigid point $\alpha_x\colon\Spf\widehat{V}_x\rightarrow\ZR{\mathscr{X}}$.
By \ref{prop-ZRpoints52} and {\bf \ref{ch-pre}}.\ref{prop-maxspe3} the image $\til{x}$ of the associated height one prime $\mathfrak{p}_{\widehat{V}_x}=\sqrt{(a)}$ is the maximal generization of $x$, whence the existence.
The equivalence of the above conditions follows immediately from \ref{prop-ZRpoint32} and {\bf \ref{ch-pre}}.\ref{prop-maxspe4}.
\end{proof}

\begin{cor}\label{cor-valuativerigidZRsp}
Let $\mathscr{X}$ be a rigid space, and $\ZR{\mathscr{X}}$ the associated Zariski-Riemann space.
Then the topological space $\ZR{\mathscr{X}}$ is valuative\index{space@space (topological)!valuative topological space@valuative ---}\index{valuative!valuative topological space@--- (topological) space} {\rm ({\bf \ref{ch-pre}}.\ref{dfn-valuativespace})}.
Moreover, for a morphism $\varphi\colon\mathscr{X}\rightarrow\mathscr{Y}$ of rigid spaces the induced map $\ZR{\varphi}\colon\ZR{\mathscr{X}}\rightarrow\ZR{\mathscr{Y}}$ between the associated Zariski-Riemann spaces is valuative\index{valuative!valuative map@--- map} {\rm ({\bf \ref{ch-pre}}.\ref{dfn-valuativemaps})}.
\end{cor}

\begin{proof}
By \ref{prop-zariskiriemanntopz} the topological space $\ZR{\mathscr{X}}$ is locally coherent {\rm ({\bf \ref{ch-pre}}.\ref{dfn-locallycoherent})}\index{coherent!coherent topological space@--- (topological) space!locally coherent topological space@locally --- ---}\index{space@space (topological)!coherent topological space@coherent ---!locally coherent topological space@locally --- ---} and sober {\rm ({\bf \ref{ch-pre}}, \S\ref{subsub-sober})}\index{space@space (topological)!sober topological space@sober ---} and is valuative due to \ref{cor-ZRpoints53} and \ref{prop-ZRpoints6}.
The last assertion follows from \ref{prop-ZRpoints6}.
\end{proof}

\subsection{Tubes}\label{sub-tubes}
\subsubsection{Tubes}\label{subsub-tubesclosure2}
\index{tube|(}
The following statements are the recasts of the already proven {\bf \ref{ch-pre}}.\ref{cor-ZRpoints61a}, {\bf \ref{ch-pre}}.\ref{prop-tubes11112a}, and {\bf \ref{ch-pre}}.\ref{prop-tubes121a}:
\begin{thm}\label{thm-ZRpoints61}
Let $\mathscr{X}$ be a rigid space, and $\mathfrak{U}\hookrightarrow\ZR{\mathscr{X}}$ a retrocompact\index{retrocompact} {\rm ({\bf \ref{ch-pre}}.\ref{dfn-retrocompact}}) open subset of $\ZR{\mathscr{X}}$.
Then we have
$$
\ovl{\mathfrak{U}}=\bigcup_{x\in\mathfrak{U}}\ovl{\{x\}}.
$$
In other words, the closure $\ovl{\mathfrak{U}}$ is the set of all specializations\index{specialization} of points of $\mathfrak{U}$.\hfill$\square$
\end{thm}

\begin{cor}\label{cor-tubes11112}
Let $\mathscr{X}$ be a rigid space, and $\mathfrak{U}\hookrightarrow\ZR{\mathscr{X}}$ a retrocompact open subset of $\ZR{\mathscr{X}}$.
For $x\in\ZR{\mathscr{X}}$ to belong to $\ovl{\mathfrak{U}}$ it is necessary and sufficient that the maximal generization\index{generization!maximal generization@maximal ---} $\til{x}$ of $x$ belongs to $\mathfrak{U}$.\hfill$\square$
\end{cor}

\begin{prop}\label{prop-tubes121}
Let $\mathscr{X}$ be a rigid space, and $\{\mathfrak{U}_{\alpha}\}_{\alpha\in L}$ a family of retrocompact open sets of $\ZR{\mathscr{X}}$. 
Then we have
$$
\ovl{\bigcap_{\alpha\in L}\mathfrak{U}_{\alpha}}=\bigcap_{\alpha\in L}\ovl{\mathfrak{U}_{\alpha}}.\eqno{\square}
$$
\end{prop}

\begin{dfn}[{\rm cf.\ {\bf \ref{ch-pre}}.\ref{dfn-tubes1a}}]\label{dfn-tubes1}{\rm 
Let $\mathscr{X}$ be a rigid space. 
A {\em tube closed subset}\index{tube!tube closed subset@--- closed subset} of $\mathscr{X}$ is a closed subset of $\ZR{\mathscr{X}}$ of the form $\ovl{\mathfrak{U}}$ by a retrocompact\index{retrocompact} open subset $\mathfrak{U}\subseteq\ZR{\mathscr{X}}$.
The complement of a tube closed subset is called a {\em tube open subset}\index{tube!tube open subset@--- open subset}.
Tube closed and tube open subsets are collectively called {\em tube subsets}\index{tube!tube subset@--- subset}.}
\end{dfn}

Let $X$ be a coherent adic formal scheme of finite ideal type, and $Y$ a closed subscheme of $X$ of finite presentation defined by an admissible ideal\index{admissible!admissible ideal@--- ideal} (cf.\ {\bf \ref{ch-formal}}.\ref{prop-admissibleidealquotientclosedsubscheme}).
We set
$$
C_{Y|X}=\sp^{-1}_X(Y)^{\circ}, 
$$
where $(\,\cdot\,)^{\circ}$ means the topological interior kernel.
The subset $C_{Y|X}\subseteq\ZR{X^{\rig}}$ is called the {\em tube of $Y$ in $X$}.
Notice that $C_{Y|X}$ depends only on the topological structure of $Y$. 
There is a canonical projection map $C_{Y|X}\rightarrow Y$ induced from the specialization map $\sp_X$, which is clearly continuous. 
\begin{prop}\label{prop-tubes91}
Let $X$ be a coherent adic formal scheme of finite ideal type, and $Y$ a closed subscheme of $X$ of finite presentation defined by an admissible ideal\index{admissible!admissible ideal@--- ideal}.
Then the subset $C_{Y|X}$ is a tube open subset\index{tube!tube open subset@--- open subset} of $X^{\rig}$.
Conversely, any tube open subset of the coherent rigid space $\mathscr{X}$ is of this form for some $X$ and $Y$.
\end{prop}

\begin{proof}
Let $U=X\setminus Y$, which is a coherent open subset of $X$, and $\mathfrak{U}=\sp^{-1}_X(U)$.
Then $C_{Y|X}$ is the complement of $\ovl{\mathfrak{U}}$, whence the first assertion.
Conversely, let $\mathfrak{U}$ be a quasi-compact open subset of $\ZR{\mathscr{X}}$.
We want to show that its complement $C=\ZR{\mathscr{X}}\setminus\ovl{\mathfrak{U}}$ is of the form $C_{Y|X}$.
Take a formal model $X$ and a quasi-compact open subset $U\subseteq X$ such that $\sp^{-1}_X(U)=\mathfrak{U}$.
Since $U$ is quasi-compact, there exists a closed subscheme $Y$ of $X_0$ of finite presentation with the underlying topological space $X\setminus U$ (Exercise \ref{exer-lemqcptcomplementsupp}). 
Then we have $C=C_{Y|X}$.
\end{proof}

The following statements are transcriptions of {\bf \ref{ch-pre}}.\ref{cor-tubes11111a} and {\bf \ref{ch-pre}}.\ref{cor-proptubes10a}:
\begin{prop}\label{prop-cortubes11111}
Let $\mathscr{X}$ be rigid space, and $C=(\ZR{\mathscr{X}}\setminus\mathfrak{U})^{\circ}$ $($where $\mathfrak{U}\subseteq\ZR{\mathscr{X}}$ is retrocompact open$)$ a tube open subset.
Then a point $x\in\ZR{\mathscr{X}}$ belongs to $C$ if and only if the maximal generization $\til{x}$ does not belong to $\mathfrak{U}$. 
Especially, every height one point of $\ZR{\mathscr{X}}\setminus\mathfrak{U}$ lies in $C$. \hfill$\square$
\end{prop}

\begin{prop}\label{prop-tubes10}
{\rm (1)} Any finite union of tube closed $($resp.\ tube open$)$ subsets\index{tube!tube closed subset@--- closed subset}\index{tube!tube open subset@--- open subset} is a tube closed $($resp.\ tube open$)$ subset.

{\rm (2)} Any finite intersection of tube closed $($resp.\ tube open$)$ subsets\index{tube!tube closed subset@--- closed subset}\index{tube!tube open subset@--- open subset} is a tube closed $($resp.\ tube open$)$ subset.\hfill$\square$
\end{prop}
\index{tube|)}

\subsubsection{Explicit description}\label{subsub-explicitdesc}
\begin{prop}[Closure formula]\label{prop-closureformula0402}
Let $X$ be a coherent adic formal scheme of finite ideal type, $Y\subseteq X$ a closed subscheme of $X$ defined by an admissible ideal $\mathscr{J}_Y$, and $U = X \setminus Y$ the open complement. 
Let $\mathscr{X}=X^{\rig}$ be the associated coherent rigid space, and set $\mathfrak{U}=\sp^{-1}_X(U)$, where $\sp_X\colon\mathscr{X}\rightarrow X$ is the specialization map.
Let $\mathscr{I}$ be an ideal of definition of $X$, and $\pi_n\colon X_n\rightarrow X$ the admissible blow-up along $\mathscr{J}_n=\mathscr{I}+\mathscr{J}^n_Y$.
Then we have
$$
\ovl{\mathfrak{U}}=\bigcap_{n\geq 1}\sp^{-1}_{X_n}(V_n),
$$
where $V_n$ is the maximal open subset of $X_n$ such that $\mathscr{J}_n\O_{X_n}|_{V_n}=\mathscr{J}_Y^n\O_{X_n}|_{V_n}$ for $n\geq 1$.
\end{prop}

\begin{proof}
Let $D$ be the right-hand side of the formula. 
We first show that for $x\in\ZR{\mathscr{X}}$ to belong to $D$ it is necessary and sufficient that the maximal generization $\til{x}$ of $x$ belongs to $\mathfrak{U}$. 
Let $x\in\ZR{\mathscr{X}}$, and take $V=V_x$ (as in \ref{ntn-ZRpoints}).
Let $I$ and $J$ be the pull-back ideals of $V$ of $\mathscr{I}$ and $\mathscr{J}_Y$, respectively, which are principal.
By the definition of $D$ a point $x$ lies in $D$ if and only if $I$ is divisible by $J^n$ for any $n\geq 1$. But the last condition is equivalent to $JV'=V'$, where $V'$ is the height 1 valuation ring associated to $V$, and hence is also equivalent to $\til{x}\in\mathfrak{U}$.
By \ref{cor-tubes11112} we deduce that $x\in D$ if and only if $x\in\ovl{\mathfrak{U}}$ and hence that $D=\ovl{\mathfrak{U}}$.  
\end{proof}

\begin{prop}[Open interior formula]\label{prop-tubes11}
Let $X$ be a coherent adic formal scheme of finite ideal type, and $Y\subseteq X$ a closed subscheme of $X$ defined by an admissible ideal $\mathscr{J}_Y$, and consider the associated coherent rigid space $\mathscr{X}=X^{\rig}$.
Let $\mathscr{I}$ be an ideal of definition of $X$.
For $n\geq 1$ we set $\mathscr{J}_n=\mathscr{I}+\mathscr{J}^n_Y$, and let $\pi_n\colon X_n\rightarrow X$ be the admissible blow-up along $\mathscr{J}_n$.
Then we have
$$
C_{Y|X}=\sp^{-1}_X(Y)^{\circ}=\varinjlim_{n\geq 1}\sp^{-1}_{X_n}(U_n),\leqno{(\ast)}
$$
where $U_n$ for each $n\geq 1$ is the open subset of $X_n$ that is maximal among open subsets $U$ of $X_n$ such that $\mathscr{J}_n\O_{X_n}|_U=\mathscr{I}\O_{X_n}|_U$.
\end{prop}

Note that the open subsets $U_n$ are quasi-compact.
To prove the proposition, we need:
\begin{lem}\label{lem-tubes11}
Let $X$ and $Y$ be as in {\rm \ref{prop-tubes11}}, and denote the right-hand side of $(\ast)$ by $C'$.

{\rm (1)} We have $C'\subseteq\sp^{-1}_X(Y)$.

{\rm (2)} For $x\in\ZR{\mathscr{X}}$ $($where $\mathscr{X}=X^{\rig})$ to belong to $C'$ it is necessary and sufficient that the maximal generization $\til{x}$ of $x$ belongs to $\sp^{-1}_X(Y)$. 
Especially, every height one point of $\sp^{-1}_X(Y)$ belongs to $C'$.
\end{lem}

\begin{proof}
(1) Let $x\in C'$, and take $V=V_x$ (as in \ref{ntn-ZRpoints}).
Let $I$ and $J$ be the pull-back ideals of $V$ of $\mathscr{I}$ and $\mathscr{J}_Y$, respectively, which are principal.
By the definition of $C'$ there exists $n\geq 1$ such that $J^n$ is divisible by $I$.
Hence $\sp_X$ maps $x$ to a point in $Y$.

(2) If $x\in C'$, then there exists $n\geq 1$ such that $J^n$ is divisible by $I$ as before.
This implies, a fortiori, that any point of $\Spf\widehat{V}$ is mapped to a point of $Y$ by the composition $\sp_X\circ\alpha_x$; hence so is, in particular, the height one point of $\Spf\widehat{V}$.
Conversely, if $\til{x}$ is mapped by $\sp_X$ to a point in $Y$, then $J^n$ is divisible by $I$ for some $n$; but this means that $x\in C'$.
\end{proof}

\begin{proof}[Proof of Proposition {\rm \ref{prop-tubes11}}]
By \ref{lem-tubes11} (1) it suffices to show the inclusion $C_{Y|X}\subseteq C'$.
Let $x\in C_{Y|X}$.
Suppose that the maximal generization $\til{x}$ belongs to $\mathscr{U}=\sp^{-1}_X(X\setminus Y)$.
Then $x\in\ovl{\mathscr{U}}$.
On the other hand, since $\mathscr{U}$ and $\sp^{-1}_X(Y)$ are disjoint, we deduce $\ovl{\mathscr{U}}$ and $C_{Y|X}$ are disjoint, which contradicts $x\in\ovl{\mathscr{U}}$.
Hence $\til{x}$ does not belong to $\mathscr{U}$, and thus $\til{x}\in\sp^{-1}_X(Y)$.
Since $\til{x}$ is of height one, we deduce by \ref{lem-tubes11} (2) that $x\in C'$, as desired.
\end{proof}

\begin{prop}[Norm description of tube subsets]\label{prop-tubesbynorms}
Let $X=\Spf A$ where $A$ is an adic ring of finite ideal type with a finitely generated ideal of definition $I\subseteq A$, and set $\mathscr{X}=X^{\rig}$.
Let $Y\subseteq X$ be a closed subscheme defined by an admissible ideal $J=(f_1,\ldots,f_n)\subseteq A$.
Then we have 
$$
C_{Y|X}=\{x\in\ZR{\mathscr{X}}\,|\,\|f_i(x)\|_{\mathscr{I},c}<1\ \textrm{for}\ i=1,\ldots,n\},
$$
where $\mathscr{I}=I\O^{\int}_{\mathscr{X}};$ here we consider $f_i$'s as elements of $\Gamma(\mathscr{X},\O_{\mathscr{X}})$ by the map $A\rightarrow\Gamma(\mathscr{X},\O_{\mathscr{X}})$.
\end{prop}

\begin{proof}
In view of \ref{prop-cortubes11111} and \ref{prop-lemnormsatpoints} it suffices to show that
$$
\sp_X^{-1}(Y)\cap[\mathscr{X}]=\{x\in[\mathscr{X}]\,|\,\|f_i(x)\|_{\mathscr{I},c}<1\ \textrm{for}\ i=1,\ldots,n\}
$$
holds.
For any height one point $x\in[\mathscr{X}]$, consider the associated rigid point (\ref{dfn-ZRpoints32})\index{point!rigid point@rigid ---}\index{rigid point} $\alpha\colon\Spf\widehat{V}_x\rightarrow\ZR{\mathscr{X}}$.
The composition $\sp_X\circ\alpha\colon\Spf\widehat{V}_x\rightarrow X=\Spf A$ is an adic morphism of formal schemes, hence inducing an adic morphism $A\rightarrow\widehat{V}_x$.
Now for $i=1,\ldots,n$ the inequality $\|f_i(x)\|_{\mathscr{I},c}<1$ holds if and only if the image of $f_i$ in $V_x$ lies in the maximal ideal $\m_{V_{x}}$, which is equivalent to, in view of {\bf \ref{ch-pre}}.\ref{thm-compval2006ver1} (2), that the image of $f_i$ in $\widehat{V}_x$ lies in its maximal ideal $\m_{\widehat{V}_x}$.
Since $\sp_X(x)$, as an open prime ideal of $A$, is the inverse image of $\m_{\widehat{V}_x}$ by the map $A\rightarrow\widehat{V}_x$, this is equivalent to $\sp_X(x)\in Y$.
\end{proof}

\subsection{Separation map and overconvergent sets}\label{sub-separation}
\subsubsection{Separation map}\label{subsub-separation}
Let $\mathscr{X}$ be a rigid space, and consider the associated Zariski-Riemann space $\ZR{\mathscr{X}}$.
Consider the separation map\index{separation map} ({\bf \ref{ch-pre}}, \S\ref{subsub-separationgen})
$$
\sep_{\mathscr{X}}\colon\ZR{\mathscr{X}}\longrightarrow[\mathscr{X}],\qquad
x\longmapsto\textrm{the maximal generization of $x$;\index{generization!maximal generization@maximal ---}}
$$
here $[\mathscr{X}]$ denotes the separated quotient consisting of height one points of $\ZR{\mathscr{X}}$ (\ref{prop-ZRpoints6}).

\begin{prop}[Functoriality; cf.\ {\bf \ref{ch-pre}}.\ref{cor-functorialitysep1a}]\label{prop-corfunctorialitysep1}
Any morphism $f\colon\mathscr{X}\rightarrow\mathscr{Y}$ of rigid spaces induces a unique continuous map $[f]\colon[\mathscr{X}]\rightarrow[\mathscr{Y}]$ that makes the following diagram commutative$:$
$$
\xymatrix{\ZR{\mathscr{X}}\ar[r]^{\ZR{f}}\ar[d]_{\sep_{\mathscr{X}}}&\ZR{\mathscr{Y}}\ar[d]^{\sep_{\mathscr{Y}}}\\ [\mathscr{X}]\ar[r]_{[f]}&[\mathscr{Y}]\rlap{$;$}}
$$
in other words, the formation of the separated quotients is functorial. \hfill$\square$
\end{prop}

\subsubsection{Overconvergent sets and tube subsets}\label{subsub-separationover}
\begin{dfn}[cf.\ {\bf \ref{ch-pre}}.\ref{dfn-separation2a}]\label{dfn-separation2}{\rm 
Let $\mathscr{X}$ be a rigid space, and $S$ a closed (resp.\ an open) subset of $\ZR{\mathscr{X}}$.
Then $S$ is said to be {\em overconvergent}\index{overconvergent!overconvergent subset@--- subset} if for any $x\in S$ any generization (resp.\ specialization) of $x$ belongs to $S$.}
\end{dfn}

Obviously, if $S$ is a overconvergent set, then its complement $\ZR{\mathscr{X}}\setminus S$ is also overconvergent.
The following statements are the special cases of {\bf \ref{ch-pre}}.\ref{prop-separation21a}, {\bf \ref{ch-pre}}.\ref{cor-separation211a}, and {\bf \ref{ch-pre}}.\ref{prop-separation22a}:
\begin{prop}\label{prop-separation21}
Let $\mathscr{X}$ be a rigid space, and $S$ a closed or an open subset of $\ZR{\mathscr{X}}$.
Then $S$ is overconvergent if and only if $S=\sep^{-1}_{\mathscr{X}}(\sep_{\mathscr{X}}(S))$.
$($Hence there exists a canonical bijection between the set of all overconvergent open $($resp.\ closed$)$ subsets of $\ZR{\mathscr{X}}$ and the set of all open $($resp.\ closed$)$ subsets of $[\mathscr{X}]$.$)$\hfill$\square$
\end{prop}

\begin{cor}\label{cor-separation211}
{\rm (1)} Any finite intersection of overconvergent open subsets is an overconvergent open subset.
The union of arbitrarily many overconvergent open subsets is an overconvergent open subset.

{\rm (2)} Any finite union of overconvergent closed subsets is an overconvergent closed subset.
The intersection of arbitrarily many overconvergent closed subsets is an overconvergent closed subset. \hfill$\square$
\end{cor}

\begin{prop}\label{prop-separation22}
Any tube closed $($resp.\ tube open$)$ subset\index{tube!tube closed subset@--- closed subset}\index{tube!tube open subset@--- open subset} of a rigid space is an overconvergent closed $($resp.\ overconvergent open$)$ subset. \hfill$\square$
\end{prop}

\begin{dfn}[cf.\ {\bf \ref{ch-pre}}.\ref{dfn-tubes111a}]\label{dfn-tubes111}{\rm 
Let $\mathscr{X}$ be a rigid space, and consider the separated quotient $[\mathscr{X}]$ of $\ZR{\mathscr{X}}$.
A closed (resp.\ open) subset $C$ of $[\mathscr{X}]$ is said to be a {\em tube closed} (resp.\ {\em tube open}) {\em subset}\index{tube!tube closed subset@--- closed subset}\index{tube!tube open subset@--- open subset} if $\sep_{\mathscr{X}}^{-1}(C)$ is a tube closed (resp.\ tube open) subset of $\mathscr{X}$.}
\end{dfn}

By \ref{prop-separation22} there exists a canonical order preserving bijection between the set of all tube closed (resp.\ tube open) subsets of $\mathscr{X}$ and that of tube closed (resp.\ tube open) subsets of $[\mathscr{X}]$.
Due to {\bf \ref{ch-pre}}.\ref{prop-tubes2a} we have:
\begin{prop}\label{prop-tubes2} 
Let $\mathscr X$ be a coherent rigid space. 

{\rm (1)} For any overconvergent closed set $F$ the set of all tube open subsets containing $F$ form a fundamental system of neighborhoods of $F$. 

{\rm (2)} For overconvergent closed subsets $F_1$ and  $F_2$ of $\ZR{\mathscr{X}}$ such that $F_1\cap F_2=\emptyset$, there exist tube open subsets $\mathfrak{U}_1$ and $\mathfrak{U}_2$ such that $F_i\subseteq\mathfrak{U}_i$ for $i=1,2$ and that $\mathfrak{U}_1\cap\mathfrak{U}_2=\emptyset$. \hfill$\square$
\end{prop}

\begin{cor}\label{cor-compacthausdorffsepquot}
Let $\mathscr{X}$ be a coherent rigid space.

{\rm (1)} The separated quotient $[\mathscr{X}]$ is a normal topological space\index{space@space (topological)!normal topological space@normal ---}. 

{\rm (2)} The separated quotient $[\mathscr{X}]$ is a compact $($hence Hausdorff$)$ space. 
\end{cor}

\begin{proof}
(1) follows from proposition \ref{prop-tubes2} (1), since $[\mathscr X]$ satisfied $\mathrm{T}_1$-axiom by definition. (2) follows from (1) and the quasi-compactness of $\ZR{\mathscr X}$. 
\end{proof}

\begin{cor}\label{cor-tubes211}
Let $\mathscr X$ be a coherent rigid space. 
Then the set of all tube open subsets of $[\mathscr{X}]$ forms an open basis of the topological space $[\mathscr{X}]$.
\end{cor}

\begin{proof}
This follows from \ref{prop-tubes2} (2) and the fact that, since $[\mathscr{X}]$ satisfies $\mathrm{T}_1$-axiom, any singleton set $\{x\}$ for $x\in[\mathscr{X}]$ is a closed subset.
\end{proof}

\subsubsection{Overconvergent interior}\label{subsub-rigidoverconvinterior}
\begin{dfn}[{cf.\ {\bf \ref{ch-pre}}.\ref{dfn-interior1a}}]\label{dfn-interior1}{\rm 
Let $\mathscr{X}$ be a rigid space, and $F$ a subset of $\ZR{\mathscr{X}}$.
We denote by $\int_{\mathscr{X}}(F)$ the maximal overconvergent open subset in $F$ and call it the {\em overconvergent interior}\index{overconvergent!overconvergent interior@--- interior} of $F$ in $\mathscr{X}$.}
\end{dfn}

The existence of $\int_{\mathscr{X}}(F)$ follows from \ref{cor-separation211}.
The following statements are already proved in {\bf \ref{ch-pre}}.\ref{prop-interior2a}, {\bf \ref{ch-pre}}.\ref{cor-interior3a}, and {\bf \ref{ch-pre}}.\ref{cor-interior4a}:
\begin{prop}\label{prop-interior2} 
Let $\mathscr{X}$ be a rigid space, and $\mathscr{U}\subseteq\mathscr{X}$ a quasi-compact open subspace.
Suppose that the following condition is satisfied$:$
\begin{itemize}
\item[$(\ast)$] there is a coherent open subspace $\mathscr{V}$ of $\mathscr{X}$ such that $\ovl{\mathscr{U}}\subset\mathscr{V}$. 
\end{itemize}
Then for a height one point $y$ to belong to $\int_{\mathscr{X}}(\mathscr{U})$ it is necessary and sufficient that $\ovl{\{y\}}\subset\mathscr{U}$. \hfill$\square$
\end{prop}

The condition $(\ast)$ is automatic if $\mathscr{X}$ is quasi-separated and locally quasi-compact\index{rigid space!locally quasi-compact@locally quasi-compact ---} (\ref{dfn-locallycompactspacerigid}); see \ref{prop-locstringcompsp2} below.

\begin{cor}\label{cor-interior3} 
In the situation as in {\rm \ref{prop-interior2}} the overconvergent interior $\int_{\mathscr{X}}(\mathscr{U})$ is described as
$$
\int_{\mathscr{X}}(\mathscr{U})=\sep^{-1}_{\mathscr{X}}([\mathscr{U}]\setminus\sep_{\mathscr{X}}(\partial\mathscr{U}),
$$
where $\partial\mathscr{U}=\ovl{\mathscr{U}}\setminus\mathscr{U}$. \hfill$\square$
\end{cor}

\begin{cor}\label{cor-interior4}
Let $\mathscr{X}$ be a quasi-separated rigid space, and $\mathscr{U}$ a quasi-compact open subspace of $\mathscr{X}$. 
Suppose $\ovl{\mathscr{U}}$ is quasi-compact. 
Then for a height one point $y$ to belong to $\int_{\mathscr{X}}(\mathscr{U})$ it is necessary and sufficient that $\ovl{\{y\}}\subset\mathscr{U}$. \hfill$\square$
\end{cor}

\subsection{Locally quasi-compact and paracompact rigid spaces}\label{sub-locallyquasicompact}
\subsubsection{Locally quasi-compact rigid spaces}\label{subsub-locallyquasicompact}
\begin{dfn}\label{dfn-locallycompactspacerigid}{\rm 
A rigid space $\mathscr{X}$ is said to be {\em locally quasi-compact}\index{rigid space!locally quasi-compact@locally quasi-compact ---} if the associated Zariski-Riemann space $\ZR{\mathscr{X}}$ is locally strongly compact ({\bf \ref{ch-pre}}.\ref{dfn-locallycompactspace}).}
\end{dfn}

By {\bf \ref{ch-pre}}.\ref{thm-locallystronglycompacttheorem}, {\bf \ref{ch-pre}}.\ref{prop-locstringcompsp2a}, {\bf \ref{ch-pre}}.\ref{cor-thmlocstringcompsp1a}, and Exercise \ref{exer-locallystronglycompact5} we have the following statements:
\begin{thm}\label{thm-proplocstringcompsp1}
Let $\mathscr{X}$ be a locally quasi-compact rigid space. 

{\rm (1)} The separated quotient $[\mathscr{X}]$ is a locally compact $($in particular, locally Hausdorff$)$ space, and the separation map $\sep_{\mathscr{X}}\colon\ZR{\mathscr{X}}\rightarrow[\mathscr{X}]$ is proper.

{\rm (2)} The rigid space $\mathscr{X}$ is quasi-compact $($resp.\ quasi-separated$)$ if and only if the separated quotient $[\mathscr{X}]$ is quasi-compact $($resp.\ Hausdorff$)$. \hfill$\square$
\end{thm}

\begin{prop}\label{prop-locstringcompsp2}
Let $\mathscr{X}$ be a quasi-separated rigid space. 
Then the following conditions are equivalent$:$
\begin{itemize}
\item[{\rm (a)}] $\mathscr{X}$ is locally quasi-compact$;$
\item[{\rm (b)}] for any quasi-compact open set $\mathscr{U}$ the closure $\ovl{\mathscr{U}}$ is quasi-compact$;$
\item[{\rm (c)}] there is an open covering $\{\mathscr{U}_{\alpha}\}_{\alpha\in L}$ of $\mathscr{X}$ such that $\mathscr{U}_{\alpha}$ and $\ovl{\mathscr{U}_{\alpha}}$ are quasi-compact for any $\alpha\in L$. \hfill$\square$
\end{itemize}
\end{prop}

\subsubsection{Paracompact rigid spaces}\label{subsub-rigidparacompactspaces}
\begin{dfn}[cf.\ {\bf \ref{ch-pre}}.\ref{dfn-paracompactgentop}]\label{dfn-paracompactrigid}{\rm 
A rigid space $\mathscr{X}$ is said to be {\em paracompact}\index{rigid space!paracompact rigid space@paracompact ---} if the topological space $\ZR{\mathscr{X}}$ is paracompact.}
\end{dfn}

For the proofs of the following statements, see {\bf \ref{ch-pre}}.\ref{lem-paracompact1} and {\bf \ref{ch-pre}}.\ref{prop-paracompact1}, respectively:
\begin{lem}\label{lem-paracompactrigid1}
Let $\mathscr{X}$ be a rigid space.

{\rm (1)} The rigid space $\mathscr{X}$ is paracompact if there is a locally finite covering by quasi-compact open sets.

{\rm (2)} If $\{\mathscr{U}_{\alpha}\}_{\alpha\in L}$ is a locally finite covering such that $\mathscr{U}_{\alpha}$ is quasi-compact for any $\alpha\in L$, then for each $\alpha\in L$ there are finitely many $\beta\in L$ such that $\mathscr{U}_{\beta}$ intersect $\mathscr{U}_{\alpha}$. \hfill$\square$
\end{lem}

\begin{prop}\label{prop-paracompact1r}
Let $\mathscr{X}$ be a paracompact quasi-separated rigid space. 
Then $\mathscr{X}$ is locally quasi-compact. \hfill$\square$
\end{prop}

\addcontentsline{toc}{subsection}{Exercises}
\subsection*{Exercises}
\begin{exer}[General closure formula]\label{exer-generalclosure}{\rm 
Let $X$ be a coherent adic formal scheme of finite ideal type, $\mathscr{J}$ and $\mathscr{K}$ admissible ideals, and $U$ the complement of $V(\mathscr{J})$ in $X$.
Let $X_n$ for $n\geq 1$ be the admissible blow-up along the ideal $\mathscr{K}+\mathscr{J}^n$, $V_n$ the maximal open subset of $X_n$ where $\mathscr{J}^n$ generates the pull-back of $\mathscr{K}+\mathscr{J}^n$.
Set $\mathfrak{U}=\sp^{-1}_X(U)$ and $\mathfrak{V}_n=\sp^{-1}_{X_n}(V_n)$ for $n\geq 1$.

(1) Show that $V_n$ is quasi-compact and that the closure $\ovl{\mathfrak{U}}$ is given by the formula
$$
\ovl{\mathfrak{U}}=\bigcap_{n\geq 1}\ovl{\mathfrak{V}}_n.
$$

(2) Suppose, moreover, that $\mathscr{K}$ is an ideal of definition of $X$.
Then show that the closure $C_n$ of $U$ inside $X_n$ is contained in $V_n$ and hence that the morphism $C_{n+1}\rightarrow C_n$ is proper.

(3) Suppose $\mathscr{K}$ is an ideal of definition.
Then show that
$$
\ovl{\mathfrak{U}}=\bigcap_{n\geq 1}\mathfrak{V}_n.
$$

(4) Show that for any quasi-compact open subset $\mathfrak{U}$ of a coherent rigid space $\mathscr{X}$ the closure $\ovl{\mathfrak{U}}$ admits a countable system of neighborhoods.}
\end{exer}

\begin{exer}\label{exer-tubesubsetscoherentmaps}{\rm 
Let $\varphi\colon\mathscr{X}\rightarrow\mathscr{Y}$ be a morphism between rigid spaces, and $T\subseteq\ZR{\mathscr{Y}}$ a tube closed (resp.\ tube open) subset.
Show that $\ZR{\varphi}^{-1}(T)$ is a tube closed (resp.\ tube open) subset of $\ZR{\mathscr{X}}$.}
\end{exer}

\begin{exer}\label{exer-spectralnorm}{\rm 
Let $\mathscr{X}$ be a coherent rigid space, and $\mathscr{I}\subseteq\O^{\int}_{\mathscr{X}}$ an ideal of definition of finite type.
Let $0<c<1$ be a real number.

(1) Show that for any $f\in\Gamma(\mathscr{X},\O_{\mathscr{X}})$ the map $\|f(\cdot)\|_{\mathscr{I},c}\colon\ZR{\mathscr{X}}\rightarrow\R_{\geq 0}$ factors through the separation map $\sep_{\mathscr{X}}\colon\ZR{\mathscr{X}}\rightarrow[\mathscr{X}]$.

(2) Show that the map $\|f(\cdot)\|_{\mathscr{I},c}\colon\ZR{\mathscr{X}}\rightarrow\R_{\geq 0}$ is continuous.}
\end{exer}


\section{Coherent sheaves}\label{sec-coherentsheavesrigid}
The aim of this section is to discuss coherent sheaves on locally universally Noetherian rigid spaces (\ref{dfn-universallyadhesiverigidspaces})\index{rigid space!universally Noetherian rigid space@universally Noetherian ---!locally universally Noetherian rigid space@locally --- ---} and their formal models.
Here, by a coherent sheaf on a rigid space $\mathscr{X}$, we mean a coherent module sheaf over the rigid structure sheaf $\O_{\mathscr{X}}$ (cf.\ \ref{dfn-ZRstrsheaf3}) on the the associated Zariski-Riemann space $\ZR{\mathscr{X}}$.
In order to lay reasonable foundations for the study of coherent sheaves on rigid spaces, it is of course necessary first to establish that the rigid structure sheaf $\O_{\mathscr{X}}$ is coherent as a module over itself, that is, the locally ringed space $(\ZR{\mathscr{X}},\O_{\mathscr{X}})$ is cohesive ({\bf \ref{ch-pre}}.\ref{dfn-cohesive}) (the analogue of Oka's theorem in complex analysis).
To show this fundamental result, we first study in \S\ref{sub-latticemodels} formal models of an $\O_{\mathscr{X}}$-module, that is, when $\mathscr{X}$ is coherent, $\O_X$-modules that give rise to the given $\O_{\mathscr{X}}$-module by passage to the functor `rig' given by `inverting' the ideal of definition.
Then the coherency of the rigid structure sheaf is obtained from the existence theorem (weak form) of finitely presented formal models for finitely presented $\O_{\mathscr{X}}$-modules, which we will discuss in \S\ref{sub-existlattice}.

In \S\ref{sub-tateacyclic} a stronger result on the existence of finitely presented formal models will be stated and proved.
This result asserts, roughly speaking, that the category of coherent sheaves on a coherent universally Noetherian rigid space $\mathscr{X}$ is equivalent to the quotient category of finitely presented sheaves on a fixed formal model $X$ of $\mathscr{X}$ mod out by the so-called {\em weak isomorphisms} (cf.\ {\bf \ref{ch-formal}}, \S\ref{subsub-weaklycoherentsheaves})\index{weak!weak isomorphism@--- isomorphism}, that is, isomorphisms up to torsion by ideals of definition.

\subsection{Formal models of sheaves}\label{sub-latticemodels}
\subsubsection{The `rig' functor for $\O_X$-modules}\label{subsub-latticemodelsrig}
\begin{ntn}\label{ntn-modintmodrigid}{\rm 
For a rigid space $\mathscr{X}$ we denote by 
$$
\Mod^{\int}_{\mathscr{X}}\quad\textrm{and}\quad\Mod_{\mathscr{X}}
$$
the categoriy of $\O^{\int}_{\mathscr{X}}$-modules and the category of $\O_{\mathscr{X}}$-modules, respectively.}
\end{ntn}

Let $\mathscr{X}$ be a coherent rigid space, and $X$ a formal model\index{formal model!formal model of a coherent rigid space@--- (of a coherent rigid space)} of $\mathscr{X}$.
For any $\O_X$-module $\mathscr{F}$ we denote by $\mathscr{F}^{\rig}$ the $\O_{\mathscr{X}}$-module defined by 
$$
\mathscr{F}^{\rig}=\sp^{-1}_X\mathscr{F}\otimes_{\sp^{-1}_X\O_X}\O_{\mathscr{X}}, 
$$
where $\sp_X\colon\ZR{\mathscr{X}}\rightarrow X$ is the specialization map defined in \S\ref{subsub-ZRdef}.
For a morphism $\varphi\colon\mathscr{F}\rightarrow\mathscr{G}$ of $\O_X$-modules we denote by $\varphi^{\rig}$ the induced morphism 
$$
\varphi^{\rig}\colon\mathscr{F}^{\rig}\longrightarrow\mathscr{G}^{\rig}
$$
of $\O_{\mathscr{X}}$-modules.
Thus we get the functor
$$
\cdot^{\rig}\colon\Mod_X\longrightarrow\Mod_{\mathscr{X}}.
$$

\begin{prop}\label{prop-latticemodels1}
Let $X$ be an a coherent adic formal scheme with an ideal of definition $\mathscr{I}_X$ of finite type, and $\mathscr{F}_X$ an $\O_X$-module.
For any admissible blow-up $\pi\colon X'\rightarrow X$, set $\mathscr{I}_{X'}=\mathscr{I}_X\O_{X'}$ and $\mathscr{F}_{X'}=\pi^{\ast}\mathscr{F}_X$.
Then one has
$$
\mathscr{F}^{\rig}_X=\varinjlim_{\pi\colon X'\rightarrow X}\varinjlim_{n\geq 0}\sp^{-1}_{X'}\lHom_{\O_{X'}}(\mathscr{I}^n_{X'},\mathscr{F}_{X'}),
$$
where the first inductive limit\index{limit!inductive limit@inductive ---} is taken along the cofiltered category $\BL_X$ of all admissible blow-ups\index{blow-up!admissible blow-up@admissible ---} of $X$ $(\S\ref{sub-categoryadmblow-up})$.
Moreover, $\mathscr{F}_{X'}$ in the right-hand side can be replaced by the strict transform\index{strict transform} $\pi'\mathscr{F}$ {\rm (\ref{dfn-stricttransform})}.
\end{prop}

To show this, we first need to prove the following lemma:
\begin{lem}\label{lem-cohsheavesZR01lem}
Let $\mathscr{X}$ be a coherent rigid space, and $\mathscr{F}$ an $\O^{\int}_{\mathscr{X}}$-module.
Then we have a canonical isomorphism
$$
\mathscr{F}\otimes_{\O^{\int}_{\mathscr{X}}}\O_{\mathscr{X}}\stackrel{\sim}{\longrightarrow}\varinjlim_{n>0}\lHom_{\O^{\int}_{\mathscr{X}}}(\mathscr{I}^n,\mathscr{F})
$$
of $\O_{\mathscr{X}}$-modules where $\mathscr{I}$ is an ideal of definition of $\ZR{\mathscr{X}}$ of finite type {\rm (\ref{dfn-ZRstrsheaf2})}.
\end{lem}

\begin{proof}
The morphism is constructed as follows.
By {\bf \ref{ch-pre}}.\ref{prop-directlimits3} the left-hand side is equal to
$$
\mathscr{F}\otimes_{\O^{\int}_{\mathscr{X}}}\varinjlim_{n>0}\lHom_{\O^{\int}_{\mathscr{X}}}(\mathscr{I}^n,\O^{\int}_{\mathscr{X}})=\varinjlim_{n>0}\mathscr{F}\otimes_{\O^{\int}_{\mathscr{X}}}\lHom_{\O^{\int}_{\mathscr{X}}}(\mathscr{I}^n,\O^{\int}_{\mathscr{X}}).
$$
On the other hand, we have a canonical morphism (cf.\ \cite[$\mathbf{0}_{\mathbf{I}}$, (5.4.2)]{EGA})
$$
\varinjlim_{n>0}\mathscr{F}\otimes_{\O^{\int}_{\mathscr{X}}}\lHom_{\O^{\int}_{\mathscr{X}}}(\mathscr{I}^n,\O^{\int}_{\mathscr{X}})\longrightarrow\varinjlim_{n>0}\lHom_{\O^{\int}_{\mathscr{X}}}(\mathscr{I}^n,\mathscr{F}).
$$
The desired morphism is the composition of these two morphisms.
To show that it is an isomorphism, we may check it stalkwise at any point $x\in\ZR{\mathscr{X}}$.
By \ref{cor-ZRstrsheaf211} we may put $\mathscr{I}_x=(a)$ by a non-zero divisor $a\in A_x=\O^{\int}_{\mathscr{X},x}$.
Then the map between stalks is described as 
$$
\mathscr{F}_x\otimes_{A_x}A_x{\textstyle [\frac{1}{a}]}=\mathscr{F}_x{\textstyle [\frac{1}{a}]}\longrightarrow\varinjlim_{n>0}\lHom_{\O^{\int}_{\mathscr{X}}}(\mathscr{I}^n,\mathscr{F})_x=\varinjlim_{n>0}\lHom_{A_x}(a^nA_x,\mathscr{F}_x),
$$
where the last equality is due to \cite[(4.1.1)]{Groth1} and the fact that $\mathscr{I}$ is invertible (\ref{cor-ZRstrsheaf2131}).
Now the last module is easily seen to be isomorphic to $\mathscr{F}_x[\frac{1}{a}]$, whence the lemma.
\end{proof}

\begin{proof}[Proof of Proposition {\rm \ref{prop-latticemodels1}}]
Consider the canonical map
\begin{equation*}
\begin{split}
\varinjlim_{\pi\colon X'\rightarrow X}\varinjlim_{n\geq 0}\sp^{-1}_{X'}\lHom_{\O_{X'}}(\mathscr{I}^n_{X'},\mathscr{F}_{X'})&\longrightarrow\\ \varinjlim_{n>0}\lHom_{\O^{\int}_{\mathscr{X}}}&(\mathscr{I}^n,\sp^{-1}_X\mathscr{F}_X\otimes_{\sp^{-1}_X\O_X}\O^{\int}_{\mathscr{X}})\cong\mathscr{F}^{\rig}_X,
\end{split}
\end{equation*}
where the last isomorphism is due to \ref{lem-cohsheavesZR01lem}.
To show that this is an isomorphism, one may check stalkwise.
Replacing $X$ by the admissible blow-up along $\mathscr{I}_X$, one can assume that $\mathscr{I}_X$ is invertible.
Then for any admissible blow-up $\pi\colon X'\rightarrow X$ the ideal $\mathscr{I}_{X'}$ is again invertible (\ref{cor-blowups155}).
Let $x\in\ZR{\mathscr{X}}$, and set $\mathscr{I}_{X,\sp(x)}=(a)$.
Then $\O_{\mathscr{X}}=\O^{\int}_{\mathscr{X},x}[\frac{1}{a}]$, and hence $\mathscr{F}^{\rig}_{X,x}=\sp^{-1}_X\mathscr{F}\otimes\O^{\int}_{\mathscr{X},x}[\frac{1}{a}]$.
On the other hand, by \cite[(4.1.1)]{Groth1} we have
\begin{equation*}
\begin{split}
\varinjlim_{\pi\colon X'\rightarrow X}\varinjlim_{n\geq 0}(\sp^{-1}_{X'}\lHom_{\O_{X'}}(\mathscr{I}^n_{X'},\mathscr{F}_{X'}))_x&=
\varinjlim_{\pi\colon X'\rightarrow X}\varinjlim_{n\geq 0}\Hom((a^n),\mathscr{F}_{X',\sp_{X'}(x)})\\
&=\varinjlim_{\pi\colon X'\rightarrow X}\mathscr{F}_{X',\sp_{X'}(x)}\otimes\O_{X',\sp_{X'}(x)}[{\textstyle \frac{1}{a}}],
\end{split}
\end{equation*}
which coincides with $\mathscr{F}^{\rig}_x$ by {\bf \ref{ch-pre}}.\ref{prop-directlimits3}.
The second assertion is clear.
\end{proof}

\begin{cor}\label{cor-latticemodels20}
Let $X$ be a coherent adic formal scheme of finite ideal type, and set $\mathscr{X}=X^{\rig}$.
Then the functor $\cdot^{\rig}\colon\Mod_X\rightarrow\Mod_{\mathscr{X}}$ is right-exact.
If, moreover, $X$ is universally rigid-Noetherian\index{formal scheme!universally rigid-Noetherian formal scheme@universally rigid-Noetherian ---} {\rm ({\bf \ref{ch-formal}}.\ref{dfn-formalsch})}, then it is exact
\end{cor}

\begin{proof}
The first part is clear.
To show the rest, as the proof of \ref{prop-latticemodels1} indicates, it suffices to show the following: let $\pi\colon X'\rightarrow X$ be an admissible blow-up of a coherent universally rigid-Noetherian formal schemes\index{formal scheme!universally rigid-Noetherian formal scheme@universally rigid-Noetherian ---} {\rm ({\bf \ref{ch-formal}}.\ref{dfn-formalsch})}, where $X'$ is assumed to have an invertible ideal of definition, and $x\in X'$; then the map $\O_{X,\pi(x)}\rightarrow\O_{X',x}[\frac{1}{a}]$ (where $a$ is a local generator of an ideal of definition of $X'$) is flat.
To this end, we may assume $X$ is affine $X=\Spf A$, where $A$ is a t.u.\ rigid-Noetherian ring\index{t.u. rigid-Noetherian ring@t.u.\ rigid-Noetherian} {\rm ({\bf \ref{ch-formal}}.\ref{dfn-tuaringadmissible} (1))}; since Zariski localization is flat on universally rigid-Noetherian formal schemes ({\bf \ref{ch-pre}}.\ref{prop-btarf1}), we may concentrate on showing that, in the notation as in \ref{prop-explicitdescription}, the map $\Spec B\setminus V(IB)\rightarrow\Spec A$ is flat.
The last map can decompose into $\Spec B\setminus V(IB)\rightarrow\Spec A\setminus V(I)\hookrightarrow\Spec A$, and $B$ is the $I$-adic completion of an affine patch of a blow-up along an admissible ideal.
This means due to {\bf \ref{ch-pre}}.\ref{prop-btarf1} that the morphism $\Spec B\setminus V(IB)\rightarrow\Spec A\setminus V(I)$ is flat.
\end{proof}

\begin{cor}\label{cor-latticemodels21}
Let $X$ be a coherent adic formal scheme of finite ideal type.
The following conditions for an $\O_X$-module $\mathscr{F}_X$ of finite type are equivalent$:$
\begin{itemize}
\item[{\rm (a)}] $\mathscr{F}^{\rig}_X=0;$
\item[{\rm (b)}] there exists an admissible blow-up $\pi\colon X'\rightarrow X$ such that $\pi^{\ast}\mathscr{F}_X$ is an $\mathscr{I}_{X'}$-torsion module, where $\mathscr{I}_{X'}$ is an ideal of definition of finite type of $X';$
\item[{\rm (c)}] there exists an admissible blow-up $\pi\colon X'\rightarrow X$ such that the strict transform $\pi'\mathscr{F}_X$ is zero. \hfill$\square$
\end{itemize}
\end{cor}

\begin{cor}\label{cor-existlattice0}
Let $X$ be a coherent universally rigid-Noetherian formal scheme\index{formal scheme!universally rigid-Noetherian formal scheme@universally rigid-Noetherian ---}, and $\mathscr{I}_X$ an ideal of definition of finite type of $X$.
Let $\varphi_X\colon\mathscr{F}_X\rightarrow\mathscr{G}_X$ be a morphism of $\O_X$-modules of finite type.
Then the following condition is equivalent$:$
\begin{itemize}
\item[{\rm (a)}] $\varphi^{\rig}_X=0;$
\item[{\rm (b)}] there exists an admissible blow-up $\pi\colon X'\rightarrow X$ such that the induced morphism between strict transforms $\varphi_{X'}\colon\mathscr{F}_{X'}\rightarrow\mathscr{G}_{X'}$ is a zero map.
\end{itemize}
\end{cor}

\begin{proof}
The implication (b) $\Rightarrow$ (a) is clear.
The converse follows from \ref{cor-latticemodels21} applied to the image of $\varphi_X$.
\end{proof}

\subsubsection{Formal models and lattice models}\label{subsub-latticemodelsdef}
\begin{dfn}\label{dfn-formalmodelcohsheaf1}{\rm 
Let $\mathscr{X}$ be a coherent {\em universally Noetherian} rigid space (\ref{dfn-universallyadhesiverigidspaces}).\index{rigid space!universally Noetherian rigid space@universally Noetherian ---}

(1) Let $\mathscr{F}$ be an $\O_{\mathscr{X}}$-module.
A {\em formal model}\index{formal model!formal model of an OXmodule@--- (of an $\O_{\mathscr{X}}$-module)} of $\mathscr{F}$ is a datum $((X,\phi),(\mathscr{F}_X,\varphi))$ consisting of
\begin{itemize}
\item a formal model $(X,\phi)$ of $\mathscr{X}$ (\ref{dfn-formalmodelcat} (1)),
\item an adically quasi-coherent\index{quasi-coherent!adically quasi-coherent OX module@adically --- (a.q.c.) sheaf} $\O_X$-module $\mathscr{F}_X$ of finite type  ({\bf \ref{ch-formal}}.\ref{dfn-adicqcoh}),
\item and an isomorphism $\varphi\colon\mathscr{F}^{\rig}_X\stackrel{\sim}{\rightarrow}\mathscr{F}$.
\end{itemize}
A formal model $((X,\phi),(\mathscr{F}_X,\varphi))$ is called a {\em lattice model}\index{lattice model} if $(X,\phi)$ is distinguished (\ref{dfn-cohrigidspacedist} (1)) and $\mathscr{F}_X$ is of finite type and $\mathscr{I}_X$-torsion free for some (hence all) ideal of definition $\mathscr{I}_X$ of finite type.

(2) Let $\Phi\colon\mathscr{F}\rightarrow\mathscr{G}$ be a morphism of $\O_{\mathscr{X}}$-modules.
A {\em formal model} of $\Phi$ is a datum $((X,\phi),(\Phi_X\colon\mathscr{F}_X\rightarrow\mathscr{G}_X,\varphi,\psi))$ consisting of
\begin{itemize}
\item a formal model $(X,\phi)$ of $\mathscr{X}$, 
\item a morphism $\Phi_X\colon\mathscr{F}_X\rightarrow\mathscr{G}_X$ of $\O_X$-modules, where $\mathscr{F}_X$ and $\mathscr{G}_X$ are adically quasi-coherent of finite type,
\item and isomorphisms $\varphi\colon\mathscr{F}^{\rig}_X\stackrel{\sim}{\rightarrow}\mathscr{F}$ and 
$\psi\colon\mathscr{G}^{\rig}_X\stackrel{\sim}{\rightarrow}\mathscr{G}$
\end{itemize}
such that the following square is commutative:
$$
\xymatrix{\mathscr{F}^{\rig}_X\ar[d]_{\varphi}\ar[r]^{\Phi^{\rig}_X}&\mathscr{G}^{\rig}_X\ar[d]^{\psi}\\ \mathscr{F}\ar[r]_{\Phi}&\mathscr{G}\rlap{.}}
$$
If both $((X,\phi),(\mathscr{F}_X,\varphi))$ and $((X,\phi),(\mathscr{G}_X,\psi))$ are lattice models, then we say that the formal model $((X,\phi),(\Phi_X\colon\mathscr{F}_X\rightarrow\mathscr{G}_X,\varphi,\psi))$ is a {\em lattice model}.}
\end{dfn}

Lattice models are particularly important when $\mathscr{X}$ is a coherent universally adhesive rigid space\index{rigid space!universally adhesive rigid space@universally adhesive ---}: if $X$ is universally adhesive and if $((X,\phi),(\mathscr{F}_X,\varphi))$ is a lattice model of $\mathscr{F}$, then $\mathscr{F}_X$ is a {\em coherent} $\O_X$-module ({\bf \ref{ch-formal}}.\ref{prop-adicqcoh6}).
\begin{prop}\label{prop-formalmodelcohsheaf101}
If $((X,\phi),(\mathscr{F}_X,\varphi))$ is a formal model of $\mathscr{F}$where $(X,\phi)$ is a distinguished formal model of $\mathscr{X}$, then $((X,\phi),(\mathscr{F}'_X,\varphi))$ with $\mathscr{F}'_X=\mathscr{F}_X/\mathscr{F}_{X,\textrm{$\mathscr{I}$-}\mathrm{tor}}$ gives a lattice model of $\mathscr{F}$, where $\mathscr{I}$ is an ideal of definition of $X$.
\end{prop}

\begin{proof}
By Exercise \ref{exer-Jtorsionpartquoteaqc2} $\mathscr{F}'_X$ is adically quasi-coherent of finite type.
By \ref{cor-latticemodels20} and \ref{cor-latticemodels21} we deduce that $\mathscr{F}'_X$ gives a formal model.
\end{proof}

The following proposition is clear:
\begin{prop}\label{prop-formalmodelcohsheaf1}
Let $\mathscr{X}$ be a coherent universally Noetherian rigid space.\index{rigid space!universally Noetherian rigid space@universally Noetherian ---}

{\rm (1)} Let $\mathscr{F}$ be an $\O_{\mathscr{X}}$-module, and $((X,\phi),(\mathscr{F}_X,\varphi))$ a formal model of $\mathscr{F}$.
Then for any admissible blow-up $\pi\colon X'\rightarrow X$, $((X',\phi\circ\pi^{\rig}),(\pi^{\ast}\mathscr{F}_X,\varphi))$ is a formal model of $\mathscr{F}$.
If $((X,\phi),(\mathscr{F}_X,\varphi))$ is a lattice model, then $((X',\phi\circ\pi^{\rig}),(\pi'\mathscr{F}_X,\varphi))$ is again a lattice model, where $\pi'\mathscr{F}$ denotes the strict transform $(\ref{dfn-stricttransform})$.

{\rm (2)} Let $\mathscr{F},\mathscr{G}$ be $\O_{\mathscr{X}}$-modules, $\Phi\colon\mathscr{F}\rightarrow\mathscr{G}$ a morphism $\O_{\mathscr{X}}$-modules, and $((X,\phi),(\Phi_X,\varphi,\psi))$ a formal model of $\Phi$.
Then for any admissible blow-up $\pi\colon X'\rightarrow X$ the pair $((X',\phi\circ\pi^{\rig}),(\pi^{\ast}\Phi_X,\varphi,\psi))$ is a formal model of $\Phi$.
If $((X,\phi),(\Phi_X,\varphi,\psi))$ is a lattice model, then $((X',\phi\circ\pi^{\rig}),(\pi'\Phi_X,\varphi,\psi))$ is again a lattice model, where $\pi'\Phi_X\colon\pi'\mathscr{F}_X\rightarrow\pi'\mathscr{G}_X$ is the induced morphism between strict transforms. \hfill$\square$
\end{prop}

\begin{prop}\label{prop-cohsheavesZR203}
Let $\mathscr{X}$ be a coherent universally rigid-Noetherian rigid space, and $X$ a distinguished formal model of $\mathscr{X}$.
Let $\mathscr{F}$ be a finitely presented $\O_{\mathscr{X}}$-module, and suppose we have two finitely presented formal models $\mathscr{F}_X$ and $\mathscr{G}_X$ of $\mathscr{F}$ on $X$.
Then, replacing $X$ be an admissible blow-up and $\mathscr{F}_X$ and $\mathscr{G}_X$ by their strict transforms if necessary, one can find positive integers $n,m,l$ such that, up to isomorphisms, the following inclusions hold$:$
$$
\mathscr{I}^n_X\mathscr{F}_X\subseteq\mathscr{I}^m_X\mathscr{G}_X\subseteq\mathscr{I}^l_X\mathscr{F}_X,
$$
where $\mathscr{I}_X$ is an invertible ideal of definition of $X$ of finite type.
\end{prop}

\begin{proof}
First notice that, since the $\mathscr{I}_X$-torsion part of $\mathscr{F}_X$ and $\mathscr{G}_X$ are locally bounded, there exists $l>0$ such that $\mathscr{I}^l_X\mathscr{F}_X$ and $\mathscr{I}^l_X\mathscr{G}_X$ are $\mathscr{I}_X$-torsion free.
Notice also that, since $\mathscr{I}_X$ is invertible, $\mathscr{I}^l_X\mathscr{F}_X$ and $\mathscr{I}^l_X\mathscr{G}_X$ are finitely presented.

Let $x\in X$.
The assumption implies that, replacing $X$ by an admissible blow-up and $\mathscr{F}_X$ and $\mathscr{G}_X$ by their strict transforms, we have $\mathscr{F}_{X,x}[\frac{1}{a}]\cong\mathscr{G}_{X,x}[\frac{1}{a}]$, where $\mathscr{I}_{X,x}=(a)$.
Hence we have
$$
a^n\mathscr{F}_{X,x}\subseteq a^m\mathscr{G}_{X,x}\subseteq a^l\mathscr{F}_{X,x}
$$
for some $n,m>0$,
Hence by \cite[(4.1.1)]{Groth1} there exists an open neighborhood $U$ of $x$ in $X$ such that 
$$
(\mathscr{I}_X|_U)^n\mathscr{F}_X|_U\subseteq(\mathscr{I}_X|_U)^m\mathscr{G}_X|_U\subseteq(\mathscr{I}_X|_U)^l\mathscr{F}_X|_U.
$$
Using the fact that $X$ is quasi-compact, one gets the desired inclusions for sufficiently large $n$ and $m$.
\end{proof}

\subsubsection{Weak isomorphisms}\label{subsub-latticemodelsweakisom}
\begin{dfn}[cf.\ {\bf \ref{ch-formal}}, \S\ref{subsub-weaklycoherentsheaves}]\label{dfn-cohsheavesZR200-1}{\rm 
Let $X$ be a coherent adic formal scheme, and $\mathscr{I}_X$ an ideal of definition of finite type.

(1) A morphism $\varphi\colon\mathscr{F}\rightarrow\mathscr{G}$ of $\O_X$-modules is said to be a {\em weak isomorphism}\index{weak isomorphism} if there exists an integer $s\geq 0$ such that $\mathscr{I}^s_X\ker(\varphi)=0$ and $\mathscr{I}^s_X\coker(\varphi)=0$.

(2) We say that two $\O_X$-modules $\mathscr{F}$ and $\mathscr{G}$ are {\em weakly isomorphic} if they are connected by a chain of weak isomorphisms:
$$
\mathscr{F}\longleftarrow\mathscr{H}_1\longrightarrow\mathscr{H}_2\longleftarrow\cdots\longleftarrow\mathscr{H}_r\longrightarrow\mathscr{G}.
$$}
\end{dfn}

\begin{prop}\label{prop-latticemodels3}
Let $\varphi_X\colon\mathscr{F}_X\rightarrow\mathscr{G}_X$ be a morphism between adically quasi-coherent sheaves of finite type on a coherent universally rigid-Noetherian formal scheme $X$.
Then the following conditions are equivalent$:$
\begin{itemize}
\item[{\rm (a)}] there exists an admissible blow-up $\pi\colon X'\rightarrow X$ such that $\pi^{\ast}\varphi_X$ is a weak isomorphism$;$
\item[{\rm (b)}] $\varphi^{\rig}_X$ is an isomorphism.
\end{itemize}
\end{prop}

\begin{proof}
By \ref{cor-latticemodels20} and \ref{cor-latticemodels21} we readily see that (a) implies (b).
Suppose $\varphi^{\rig}_X$ is an isomorphism.
By \ref{cor-latticemodels21} we have an admissible blow-up $\pi\colon X'\rightarrow X$ such that the kernel and the cokernel of $\pi^{\ast}\varphi_X$ are $\mathscr{I}_{X'}$-torsion, where $\mathscr{I}_{X'}=\mathscr{I}\O_{X'}$.
By \ref{exer-adicquotientbynonadic1} these are bounded $\mathscr{I}_{X'}$-torsion.
\end{proof}

\begin{cor}\label{cor-latticemodels31x}
Let $X$ be a coherent universally rigid-Noetherian formal scheme, and $\varphi_X\colon\mathscr{F}_X\rightarrow\mathscr{G}_X$ and $\psi_X\colon\mathscr{G}_X\rightarrow\mathscr{H}_X$ weak isomorphisms between adically quasi-coherent sheaves of finite type on $X$.
Then there exists an admissible blow-up $\pi\colon X'\rightarrow X$ such that $\pi^{\ast}(\psi_X\circ\varphi_X)$ and $\pi'(\psi_X\circ\varphi_X)$ are weak isomorphisms. \hfill$\square$
\end{cor}

\begin{cor}\label{cor-proplatticemodels4}
Let $\mathscr{X}$ be a coherent universally Noetherian rigid space, and $X$ a distinguished formal model of $\mathscr{X}$.
Let $\mathscr{F}$ be a finitely presented $\O_{\mathscr{X}}$-module, and $\mathscr{F}_X$ a formal model of $\mathscr{F}$ on $X$.
Then any adically quasi-coherent sheaf $\mathscr{G}_X$ of finite type on $X$ weakly isomorphic to $\mathscr{F}_X$ is a formal model of $\mathscr{F}$.
In particular, for any ideal of definition $\mathscr{I}_X$ of finite type of $X$ and any positive integer $n$, $\mathscr{I}^n_X\mathscr{F}_X$ is a formal model of $\mathscr{F}$. \hfill$\square$
\end{cor}

The following proposition is a consequence of {\bf \ref{ch-formal}}.\ref{prop-labelFPthickeningexistformal}:
\begin{prop}\label{prop-corlabelFPthickeningexist}
Let $X$ be a coherent universally rigid-Noetherian formal scheme, and $\mathscr{X}=X^{\rig}$.
If a finite type $\O_{\mathscr{X}}$-module $\mathscr{F}$ has a formal model on $X$, then it has a formal model given by a finitely presented $\O_X$-module. \hfill$\square$
\end{prop}

\subsection{Existence of finitely presented formal models (weak version)}\label{sub-existlattice}
The main goal of this subsection is to prove the following theorem:
\begin{thm}[Existence of finitely presented formal models (weak form)]\label{thm-existlattice}
Let $\mathscr{X}$ be a coherent universally Noetherian\index{rigid space!universally Noetherian rigid space@universally Noetherian ---} rigid space, and $\varphi\colon\mathscr{F}\rightarrow\mathscr{G}$ a morphism of finitely presented $\O_{\mathscr{X}}$-modules.
Then there exists a formal model $X$ of $\mathscr{X}$ and a formal model $\varphi_X\colon\mathscr{F}_X\rightarrow\mathscr{G}_X$ of $\varphi$ by finitely presented $\O_X$-modules.
\end{thm}

As a special case where $\mathscr{F}=\mathscr{G}$ and $\varphi=\id_{\mathscr{F}}$, we have:
\begin{cor}\label{cor-existlatticecor1}
Let $\mathscr{X}$ be a coherent universally Noetherian rigid space, and $\mathscr{F}$ a finitely presented $\O_{\mathscr{X}}$-module.
Then there exists a finitely presented formal model $\mathscr{F}_X$ of $\mathscr{F}$ on a formal model $X$ of $\mathscr{X}$. \hfill$\square$
\end{cor}

By \ref{prop-cohrigidspacedist1} and \ref{prop-formalmodelcohsheaf101} we deduce:
\begin{cor}\label{cor-existlatticecor2}
Let $\mathscr{X}$ and $\varphi$ be as in {\rm \ref{thm-existlattice}}.
Then there exists a distinguished formal model $X$ of $\mathscr{X}$ and a lattice model\index{lattice model} $\varphi'_X\colon\mathscr{F}'_X\rightarrow\mathscr{G}'_X$ of $\varphi$ on $X$. \hfill$\square$
\end{cor}

Notice that, if $X$ is universally adhesive\index{formal scheme!universally adhesive formal scheme@universally adhesive ---}\index{adhesive!universally adhesive@universally ---!universally adhesive formal scheme@--- --- formal scheme} (hence $\mathscr{X}$ is universally adhesive\index{rigid space!universally adhesive rigid space@universally adhesive ---}), then lattice models are automatically finitely presented.
For the proof of the theorem, we first need to show the following:
\begin{lem}\label{lem-existlattice1}
Let $X$ be a coherent universally rigid-Noetherian formal scheme with an ideal of definition $\mathscr{I}_X$ of finite type, and $\mathscr{F}_X$ and $\mathscr{G}_X$ finitely presented $\O_X$-modules.
Suppose there exists a morphism $\varphi\colon\mathscr{F}^{\rig}_X\rightarrow\mathscr{G}^{\rig}_X$ of $\O_{\mathscr{X}}$-modules, where $\mathscr{X}=X^{\rig}$.
Then, replacing $X$ by an admissible blow-up and $\mathscr{F}_X$ and $\mathscr{G}_X$ by their strict transforms, we can find a morphism $\til{\varphi}\colon\mathscr{I}^s_X\mathscr{F}_X\rightarrow\mathscr{G}_X$ of $\O_X$-modules, where $s$ is a positive integer, such that $\til{\varphi}^{\rig}=\varphi$.
\end{lem}

\begin{proof}
{\sc Step 1.} 
Replacing $X$ by the admissible blow-up along $\mathscr{I}_X$, we may assume that the ideal of definition $\mathscr{I}_X$ is invertible.
Let $x\in\ZR{\mathscr{X}}$.
By \ref{prop-latticemodels1} (and using the notation there) we have $\mathscr{G}^{\rig}_x=\varinjlim_{X'\rightarrow X}\mathscr{G}_{X',\sp_{X'}(x)}[\frac{1}{a}]$, where $\mathscr{G}_{X'}$ is the pull-back (total transform) of $\mathscr{G}_X$ by an admissible blow-up $X'\rightarrow X$, and $a$ is the generator of $\mathscr{I}_{X,\sp_X(x)}$.
On the other hand, by the definition of $\mathscr{F}^{\rig}$ we can write $\mathscr{F}^{\rig}_x=\mathscr{F}_{X,\sp_X(x)}\otimes_{\O_{X,\sp_X(x)}}\O^{\int}_{\mathscr{X},x}[\frac{1}{a}]$.
Hence the stalk map $\varphi_x$ gives
$$
\varphi_x\colon {\textstyle \mathscr{F}_{X,\sp_X(x)}\otimes_{\O_{X,\sp_X(x)}}\O^{\int}_{\mathscr{X},x}[\frac{1}{a}]}\longrightarrow\varinjlim_{X'\rightarrow X}{\textstyle \mathscr{G}_{X',\sp_{X'}(x)}[\frac{1}{a}]}.
$$
Let $\psi_x\colon\mathscr{F}_{X,\sp_X(x)}\rightarrow\varinjlim\mathscr{G}_{X',\sp_{X'}(x)}[\frac{1}{a}]$ be the composition of $\varphi_x$ preceded by the canonical map $\mathscr{F}_{X,\sp_X(x)}\rightarrow\mathscr{F}_{X,\sp_X(x)}\otimes_{\O_{X,\sp_X(x)}}\O^{\int}_{\mathscr{X},x}[\frac{1}{a}]$; $\psi_x$ is a morphism of $\O_{X,\sp_X(x)}$-modules.
Since $\mathscr{F}_{X,\sp_X(x)}$ is a finitely presented $\O_{X,\sp_X(x)}$-module, there exist an admissible blow-up $X'\rightarrow X$ and a map $\phi_x\colon\mathscr{F}_{X,\sp_X(x)}\rightarrow\mathscr{G}_{X',\sp_{X'}(x)}[\frac{1}{a}]$ such that $\psi_x$ is equal to the composition of $\phi_x$ followed by the canonical map.
Hence we get a map $\phi'_x\colon\mathscr{F}_{X',\sp_{X'}(x)}\rightarrow\mathscr{G}_{X',\sp_{X'}(x)}[\frac{1}{a}]$ by base change.
Since $\mathscr{F}_{X',\sp_{X'}(x)}$ is finitely generated, one can find an integer $s\geq 0$ and a map $\til{\varphi}_x\colon a^s\mathscr{F}_{X',\sp_{X'}(x)}\rightarrow\mathscr{G}_{X',\sp_{X'}(x)}$ that gives the restriction of $\phi'_x$.
By the construction we have
$$
\til{\varphi}_x\otimes_{\O_{X',\sp_{X'}(x)}}\O_{\mathscr{X},x}=\varphi_x.
$$

{\sc Step 2.}
The previous step implies that for any $x\in\ZR{\mathscr{X}}$ one can replace $X$ by a suitable admissible blow-up of $X$ in such a way that there exists a map $\til{\varphi}_x\colon a^{s_x}\mathscr{F}_{X,\sp_{X}(x)}\rightarrow\mathscr{G}_{X,\sp_{X}(x)}$ of $\O_{X,\sp_X(x)}$-modules that induces the stalk map $\varphi_x$ by base change.
Here $s_x$ is a non-negative integer, which depends on the point $x$.
By \cite[(4.1.1)]{Groth1} (here we use the fact that $\mathscr{F}_X$ is finitely presented) one can take an open neighborhood $U$ of $\sp_X(x)$ in $X$ and a morphism $\til{\varphi}_U\colon\mathscr{I}^s\mathscr{F}_X|_U\rightarrow\mathscr{G}_X|_U$ such that $\til{\varphi}_{U,x}=\til{\varphi}_x$.
Moreover, one can take $U$ sufficiently small such that $\til{\varphi}_{U}^{\rig}=\sp^{-1}_X\til{\varphi}_U\otimes_{\O^{\int}_{\mathscr{X}}|_{\mathfrak{U}}}\O_{\mathscr{X}}|_{\mathfrak{U}}$ coincides with $\varphi|_{\mathfrak{U}}$, where $\mathfrak{U}=\sp^{-1}_X(U)$.

\medskip
{\sc Step 3.}
Thus we have an open covering $\ZR{\mathscr{X}}=\bigcup^n_{\alpha=1}\mathfrak{U}_{\alpha}$ such that we have for each $\alpha$:
\begin{itemize}
\item an admissible blow-up $X_{\alpha}\rightarrow X$ and a quasi-compact open subset $U_{\alpha}$ of $X_{\alpha}$ such that $\mathfrak{U}_{\alpha}=\sp^{-1}_{X_{\alpha}}(U_{\alpha})$; 
\item a positive integer $s_{\alpha}$; 
\item and a morphism $\til{\varphi}_{\alpha}\colon \mathscr{I}^{s_{\alpha}}_X\mathscr{F}_{X_{\alpha}}|_{U_{\alpha}}\rightarrow\mathscr{G}_{X_{\alpha}}|_{U_{\alpha}}$ such that $\til{\varphi}_{\alpha}^{\rig}=\varphi|_{\mathfrak{U}_{\alpha}}$.
\end{itemize}
Since $\ZR{\mathscr{X}}$ is quasi-compact, the open covering as above can be taken to be finite, and hence the positive integer $s=s_{\alpha}$ can be chosen independently of $\alpha$.
Moreover, by successive use of \ref{prop-blowups4111} there exists an admissible blow-up $X'$ that dominates the all $X_{\alpha}$ on which there exists an open covering $X'=\bigcup^n_{\alpha=1}U'_{\alpha}$ consisting of quasi-compact open subsets such that there exists for each $\alpha$ a morphism 
$$
\til{\varphi}_{\alpha}\colon\mathscr{I}_X^s\mathscr{F}_{X'}|_{U'_{\alpha}}\rightarrow \mathscr{G}_{X'}|_{U'_{\alpha}}
$$
such that $\til{\varphi}^{\rig}_{\alpha}=\varphi|_{\mathfrak{U}_{\alpha}}$.
Now by \ref{cor-existlattice0} there exists a further admissible blow-up $\pi\colon X''\rightarrow X'$ of $X'$ on which the strict transforms of $\til{\varphi}_{\alpha}$ glue together to a morphism $\til{\varphi}\colon\mathscr{I}^s_{X''}\pi'\mathscr{F}_{X'}\rightarrow\pi'\mathscr{G}_{X'}$, as desired.
\end{proof}

\begin{proof}[Proof of Theorem {\rm \ref{thm-existlattice}}]
First we claim that any finitely presented $\O_{\mathscr{X}}$-module $\mathscr{F}$ has a finitely presented formal model $\mathscr{F}_X$.
Once this is shown, the assertion of the theorem immediately follows from \ref{lem-existlattice1} since, replacing $X$ by a further admissible blow-up, one can assume that $\mathscr{F}_X$ is $\mathscr{I}_X$-torsion free (where $\mathscr{I}_X$ is an ideal of definition of finite type) and that $X$ has an invertible ideal of definition.
But for this, by \ref{prop-corlabelFPthickeningexist} we only have to show that $\mathscr{F}$ has a formal model.

Let $\mathscr{F}$ be a finitely presented $\O_{\mathscr{X}}$-module.
Then there exist a finite open covering $\ZR{\mathscr{X}}=\bigcup^r_{\alpha=1}\mathfrak{U}_{\alpha}$ consisting of quasi-compact open subsets and, for each $\alpha$, an exact sequence
$$
\O^{\oplus p}_{\mathscr{X}}|_{\mathfrak{U}_{\alpha}}\stackrel{\varphi_{\lambda}}{\longrightarrow}\O^{\oplus q}_{\mathscr{X}}|_{\mathfrak{U}_{\alpha}}\longrightarrow\mathscr{F}|_{\mathfrak{U}_{\alpha}}\longrightarrow 0.
$$
(Here the numbers $p$ and $q$ depend on $\alpha$.)
By \ref{prop-zariskiriemanntoptop}, \ref{prop-cohrigidspacedist1}, and \ref{prop-ZRpoints4} there exists a distinguished formal model $X$ of $\mathscr{X}$ and an open covering $X=\bigcup^r_{\alpha=1}U_{\alpha}$ such that $\mathfrak{U}_{\alpha}=\sp^{-1}_X(U_{\alpha})$ for each $\alpha$.
Moreover, each $\mathscr{U}_{\alpha}$ can be identified with $\ZR{\mathscr{U}_{\alpha}}$, where $\mathscr{U}_{\alpha}=U^{\rig}_{\alpha}$.
Replacing $X$ by the admissible blow-up along $\mathscr{I}_X$, we may assume that $\mathscr{I}_X$ is invertible.

By \ref{lem-existlattice1}, replacing $X$ by an admissible blow-up, we have for each $\alpha=1,\ldots,r$ a morphism 
$$
\til{\varphi}_{\alpha}\colon\mathscr{I}^s_X\O^{\oplus p}_X|_{U_{\alpha}}\longrightarrow\O^{\oplus q}_X|_{U_{\alpha}}
$$
such that $\til{\varphi}_{\alpha}^{\rig}=\varphi_{\alpha}$.
Since the number of the indices $\alpha$ is finite, one can choose such an admissible blow-up independently of $\alpha$.
Set $\til{\mathscr{F}}_{\alpha}=\coker(\til{\varphi}_{\alpha})$.
Then $\til{\mathscr{F}}_{\alpha}$ is a finitely presented $\O_X|_{U_{\alpha}}$-module such that $\til{\mathscr{F}}_{\alpha}^{\rig}\cong\mathscr{F}|_{\mathfrak{U}_{\alpha}}$.

Since $X$ is coherent, $U_{\alpha\beta}=U_{\alpha}\cap U_{\beta}$ are quasi-compact.
By \ref{prop-cohsheavesZR203} we deduce that, replacing $X$ by an admissible blow-up if necessary, there exist, for each pair $\alpha,\beta$ of indices, positive integers $m,l$ such that
$$
\mathscr{I}^n_X\til{\mathscr{F}}_{\alpha}\subseteq\mathscr{I}^m_X\til{\mathscr{F}}_{\beta}\subseteq\mathscr{I}^l_X\til{\mathscr{F}}_{\alpha}
$$
on $U_{\alpha\beta}$ (here we used \ref{prop-blowups4111} to extend the possible admissible blow-up of $U_{\alpha\beta}$ to an admissible blow-up of $X$).

We want to show that there exists a finitely presented $\O_X$-module $\til{\mathscr{F}}_X$ such that for any $\alpha=1,\ldots,r$ we have
$$
\mathscr{I}^s_X\til{\mathscr{F}}_{\alpha}\subseteq\til{\mathscr{F}}_X\subseteq\mathscr{I}^t_X\til{\mathscr{F}}_{\alpha}
$$
on $U_{\alpha}$ for some $s,t>0$.
This is shown by induction with respect to $n$, and thus we reduce to the case $n=2$, that is, $X=U_{\alpha}\cup U_{\beta}$.
But then, the issue is to extend the sheaf $\mathscr{I}^m_X\til{\mathscr{F}}_{\beta}$ to $X$.
It follows from Exercise \ref{exer-extopenadiccoh1} that one can do this to get the desired extension as an adically quasi-coherent sheaf $\mathscr{F}_X$ of finite type.
By the construction we have $\til{\mathscr{F}}^{\rig}\cong\mathscr{F}_X$.
\end{proof}

\begin{cor}\label{cor-existlattice3}
Let $\mathscr{X}$ be a locally universally Noetherian rigid space.\index{rigid space!universally Noetherian rigid space@universally Noetherian ---!locally universally Noetherian rigid space@locally --- ---}
Then $\O_{\mathscr{X}}$ is a coherent $\O_{\mathscr{X}}$-module.
In other words, the locally ringed space $(\ZR{\mathscr{X}},\O_{\mathscr{X}})$ is cohesive\index{cohesive!cohesive ringed space@--- (ringed space)} $({\bf \ref{ch-pre}}.\ref{dfn-cohesive})$.
\end{cor}

\begin{proof}
Consider, for any open subspace $\mathscr{U}\subseteq\mathscr{X}$, a morphism of $\O_{\mathscr{U}}$-modules of the form $\varphi\colon\O^{\oplus p}_{\mathscr{U}}\rightarrow\O_{\mathscr{U}}$. We need to show that $\ker(\varphi)$ is of finite type. As we may work locally on $\mathscr{U}$, we may assume that $\mathscr{U}$ is coherent, which, by a further formal reduction, we may set $\mathscr{X}=\mathscr{U}$.
By \ref{lem-existlattice1} we have a formal model $X$ of $\mathscr{X}$ having an invertible ideal of definition $\mathscr{I}_X$ and a morphism $\varphi_X\colon\O_X^{\oplus p}\rightarrow\O_X$ that gives rise to $\varphi$ by $\cdot^{\rig}$.
Then, as the question is local on $\mathscr{X}$, we may further assume that $X$ is affine $X=\Spf A$.
The map $\varphi_X$ is given by a $A$-linear map $A^{\oplus p}\rightarrow A$.
Let $K$ be the kernel.
Since $A$ is t.u.\ rigid-Noetherian, we have an injective map $K'\hookrightarrow K$ from a finitely generated $A$-module such that $K/K'$ is $I$-torsion (where $I\subseteq A$ is an ideal of definition).
Since $K/K'\subseteq A^{\oplus p}/K'$, we know that $K/K'$ is bounded $I$-torsion.
If we set $\mathscr{K}'=(K')^{\Delta}$ and $\mathscr{K}=K^{\Delta}$, we have a weak isomorphism $\mathscr{K}'\hookrightarrow\mathscr{K}$.
Now by \ref{cor-latticemodels20} the kernel of $\varphi$ is isomorphic to $\mathscr{K}^{\rig}$ and hence to $\mathscr{K}^{\prime\rig}$, which is of finite type.
\end{proof}

\subsection{Existence of finitely presented formal models (strong version)}\label{sub-tateacyclic}
As usual, we denote by $\Coh_{\mathscr{X}}$ the category of coherent $\O_{\mathscr{X}}$-modules.
The goal of this subsection is to show the following theorem:
\begin{thm}\label{thm-tateacyclic20}
Let $\mathscr{X}$ be a coherent universally Noetherian rigid space\index{rigid space!universally Noetherian rigid space@universally Noetherian ---} {\rm (\ref{dfn-universallyadhesiverigidspaces})}, and $X$ a universally rigid-Noetherian formal model\index{formal model!formal model of a coherent rigid space@--- (of a coherent rigid space)} of $\mathscr{X}$.
Then the functor 
$$
\Mod^{\pf}_X\stackrel{\cdot^{\rig}}{\longrightarrow}\Coh_{\mathscr{X}},
$$
where $\Mod^{\pf}_X$ denotes the full subcategory of $\Mod_X$ consisting of finitely presented $\O_X$-modules, induces an equivalence of thecategories
$$
\Mod^{\pf}_X/\{\textrm{weak isomorphisms}\}\stackrel{\sim}{\longrightarrow}\Coh_{\mathscr{X}}.
$$
\end{thm}

As a corollary, we get a stronger result for existence of lattice models:
\begin{cor}[Existence of finitely presented formal models (strong form)]\label{cor-existlatticemodelstrong}
Let $\mathscr{X}$ be a coherent universally Noetherian rigid space, $X$ a universally rigid-Noetherian formal model of $\mathscr{X}$, and $\varphi\colon\mathscr{F}\rightarrow\mathscr{G}$ a morphism of finitely presented $\O_{\mathscr{X}}$-modules.
Then there exists a formal model $\varphi_X\colon\mathscr{F}_X\rightarrow\mathscr{G}_X$ of $\varphi$ on $X$ by finitely presented $\O_X$-modules. \hfill$\square$
\end{cor}

Similarly to \S\ref{sub-existlattice} we immediately have the following two corollaries:
\begin{cor}\label{cor-existlatticemodelstrongcor1}
Let $\mathscr{X}$ be a coherent universally Noetherian rigid space, $X$ a universally rigid-Noetherian formal model of $\mathscr{X}$, and $\mathscr{F}$ a finitely presented $\O_{\mathscr{X}}$-module.
Then there exists a finitely presented formal model $\mathscr{F}_X$ on $X$. \hfill$\square$
\end{cor}

\begin{cor}\label{cor-existlatticemodelstrongcor2}
Let $\mathscr{X}$, $X$, and $\varphi$ be as in {\rm \ref{cor-existlatticemodelstrong}}.
Then there exists a lattice model $\varphi'_X\colon\mathscr{F}'_X\rightarrow\mathscr{G}'_X$ of $\varphi$ on $X$. \hfill$\square$
\end{cor}

In order to show the theorem, we need a few preparatory results:
\begin{prop}\label{prop-tateacyclic1}
Let $X$ be a coherent universally adhesive formal scheme\index{formal scheme!universally adhesive formal scheme@universally adhesive ---}\index{adhesive!universally adhesive@universally ---!universally adhesive formal scheme@--- --- formal scheme} that is $\mathscr{I}_X$-torsion free $($resp.\ coherent universally rigid-Noetherian formal scheme\index{formal scheme!universally rigid-Noetherian formal scheme@universally rigid-Noetherian ---}$)$, where $\mathscr{I}_X$ is an ideal of definition of finite type on $X$.
Let $\pi\colon X'\rightarrow X$ be an admissible blow-up\index{blow-up!admissible blow-up@admissible ---}.

{\rm (1)} For any finitely presented $\O_X$-module $\mathscr{G}$, $\RD^q\pi_{\ast}\pi^{\ast}\mathscr{G}$ is finitely presented $($resp.\ FP-approximated\index{FPapproximated@FP-approximated} {\rm ({\bf \ref{ch-formal}}.\ref{dfn-FPthickeningsformal})}$);$ if $q>0$, it is bounded $\mathscr{I}_X$-torsion.
Moreover, the canonical morphism $\mathscr{G}\rightarrow\pi_{\ast}\pi^{\ast}\mathscr{G}$ is a weak isomorphism.

{\rm (2)} For any finitely presented $\O_{X'}$-module $\mathscr{F}$, $\RD^q\pi_{\ast}\mathscr{F}$ is finitely presented $($resp.\ FP-approximated\index{FPapproximated@FP-approximated}$);$ if $q>0$, it is bounded $\mathscr{I}_X$-torsion.
Moreover, the canonical morphism $\pi^{\ast}\pi_{\ast}\mathscr{F}\rightarrow\mathscr{F}$ is a weak isomorphism.
\end{prop}

\begin{proof}
We may assume that $X$ is affine $X=\Spf A$ (note that the property ``FP-approximated'' is local with respect to Zariski topology) and that $\mathscr{I}_X=I^{\Delta}$, where $A$ is a t.u.\ adhesive ring\index{t.u.a. ring@t.u.a.\ ring} $($resp.\ t.u.\ rigid-Noetherian ring\index{t.u. rigid-Noetherian ring@t.u.\ rigid-Noetherian ring}$)$, and $I\subseteq A$ is a finitely generated ideal of definition.
In this situation, the admissible blow-up $\pi$ is the formal completion of the usual blow-up $p\colon Y'\rightarrow Y=\Spec A$ along an admissible ideal\index{admissible!admissible ideal@--- ideal} $J$ of $A$.

(1) Since $X$ is affine, there exists a finitely presented $\O_Y$-module $\mathscr{H}$ such that $\mathscr{H}^{\for}=\mathscr{G}$, where $\cdot^{\for}$ denotes the `virtual' formal completion functor defined in \S\ref{subsub-GFGAcomannconst}.
By {\bf \ref{ch-formal}}.\ref{prop-adqformalprop} we deduce that $\pi^{\ast}\mathscr{G}=(p^{\ast}\mathscr{H})^{\for}$.
Now we apply comparison theorem ({\bf \ref{ch-formal}}.\ref{thm-GFGAcom}) (resp.\ ({\bf \ref{ch-formal}}.\ref{thm-GFGAcomweakcoherentspecial})) to deduce that $\RD^q\pi_{\ast}\pi^{\ast}\mathscr{G}=(\RD^qp_{\ast}p^{\ast}\mathscr{H})^{\for}$ for $q\geq 0$.
By finiteness theorem ({\bf \ref{ch-formal}}.\ref{thm-fini}) (resp.\ ({\bf \ref{ch-formal}}.\ref{thm-finitenessweaklycoherentmodules})) we deduce that $\RD^qp_{\ast}p^{\ast}\mathscr{H}$ is coherent (resp.\ FP-approximated) and hence that $\RD^q\pi_{\ast}\pi^{\ast}\mathscr{G}$ is coherent (resp.\ FP-approximated by {\bf \ref{ch-formal}}.\ref{prop-FPAcompletion}).
Moreover, when $q>0$, $\RD^qp_{\ast}p^{\ast}\mathscr{H}$ is $I$-torsion, and hence so is $\RD^q\pi_{\ast}\pi^{\ast}\mathscr{G}$.
It is now clear that $\mathscr{G}\rightarrow\pi_{\ast}\pi^{\ast}\mathscr{G}$ is a weak isomorphism, for the kernel and cokernel of the morphism $\mathscr{H}\rightarrow p_{\ast}p^{\ast}\mathscr{H}$ are bounded $I$-torsion.

(2) By existence theorem ({\bf \ref{ch-formal}}.\ref{thm-GFGAexist}) (resp.\ ({\bf \ref{ch-formal}}.\ref{thm-GFGAexaweakcoherentspecial})) one can take a coherent $\O_{Y'}$-module $\mathscr{H}$ such that $\mathscr{H}^{\for}=\mathscr{F}$.
Then the assertion follows by a similar argument as above.
\end{proof}

\begin{cor}\label{cor-tateacyclic}
Let $X$ be a coherent universally adhesive formal scheme\index{formal scheme!universally adhesive formal scheme@universally adhesive ---}\index{adhesive!universally adhesive@universally ---!universally adhesive formal scheme@--- --- formal scheme} that is $\mathscr{I}_X$-torsion free $($resp.\ coherent universally rigid-Noetherian formal scheme\index{formal scheme!universally rigid-Noetherian formal scheme@universally rigid-Noetherian ---}$)$, where $\mathscr{I}_X$ is an ideal of definition of finite type on $X$.
Let $\pi\colon X'\rightarrow X$ be an admissible blow-up.

{\rm (1)} The sheaf $\pi_{\ast}\O_{X'}$ is finitely presented $($resp.\ FP-approximated$)$, and the canonical morphism $\O_X\rightarrow\pi_{\ast}\O_{X'}$ is a weak isomorphism.

{\rm (2)} If $q>0$, $\RD^q\pi_{\ast}\O_{X'}$ is bounded $\mathscr{I}_X$-torsion. \hfill$\square$
\end{cor}

\begin{proof}[Proof of Theorem {\rm \ref{thm-tateacyclic20}}]
First we prove that the functor in question is essentially surjective.
Let $\mathscr{F}$ be an object of $\Coh_{\mathscr{X}}$.
By \ref{thm-existlattice} and \ref{prop-formalmodelcohsheaf1} (1) there exists an admissible blow-up $\pi\colon X'\rightarrow X$ and a finitely presented formal model $\mathscr{F}_{X'}$ of $\mathscr{F}$ on $X'$.
By \ref{prop-tateacyclic1} (2) $\mathscr{F}_X=\pi_{\ast}\mathscr{F}_{X'}$ is finitely presented (resp.\ FP-approximated, which can be replaced by a finitely presented formal model by an FP-approximation ({\bf \ref{ch-formal}}.\ref{dfn-FPthickeningsformal} (1))).
Hence the essential surjectivity is proven.
By a similar argument one can also show that for any morphism $\varphi$ in $\Coh_{\mathscr{X}}$ there exists a morphism $\varphi_X$ of $\Mod^{\pf}_X$ such that $\varphi^{\rig}_X=\varphi$
(resp.\ here we need to use the fact that the category of FP-approximations is filtered; cf.\ Exercise \ref{exer-FPapproximationfunctoriality}).
Hence the functor in question is full and essentially surjective.
The faithfullness follows from \ref{prop-latticemodels3} and \ref{prop-tateacyclic1} (2).
\end{proof}

\subsection{Integral models}\label{sub-integralomodels}
\begin{ntn}\label{ntn-cohsheavesZR}{\rm
Let $\mathscr{X}$ be a coherent rigid space, and $X$ a formal model of $\mathscr{X}$.
For any $\O_X$-module $\mathscr{F}$ we denote by $\sp^{\ast}_X\mathscr{F}$ the $\O^{\int}_{\mathscr{X}}$-module defined by 
$$
\sp^{-1}_X\mathscr{F}\otimes_{\sp^{-1}_X\O_X}\O^{\int}_{\mathscr{X}},
$$
that is, the pull-back by the morphism of locally ringed spaces $(\ZR{\mathscr{X}},\O^{\int}_{\mathscr{X}})\rightarrow(X,\O_X)$ whose underlying morphism of topological spaces is $\sp_X$.}
\end{ntn}

\begin{thm}\label{thm-cohsheavesZR1}
Let $\mathscr{X}$ be a coherent rigid space.

{\rm (1)} For any finitely presented $\O^{\int}_{\mathscr{X}}$-module $\mathscr{F}$ there exist a distinguished formal model $X$ of $\mathscr{X}$ and a finitely presented $\O_X$-module $\mathscr{F}_X$ such that $\mathscr{F}\cong\sp^{\ast}_X\mathscr{F}_X$.

{\rm (2)} Let $X$ be a distinguished formal model of $\mathscr{X}$, and $\mathscr{F}_X$ and $\mathscr{G}_X$ $\O_X$-modules.
Suppose $\mathscr{F}_X$ is finitely presented.
Let $\varphi_X,\psi_X\colon\mathscr{F}_X\rightarrow\mathscr{G}_X$ be two morphisms of $\O_X$-modules such that $\sp^{\ast}_X\varphi_X=\sp^{\ast}_X\psi_X$.
Then there exists an admissible blow-up $\pi\colon X'\rightarrow X$ such that $\pi^{\ast}\varphi_X=\pi^{\ast}\psi_X$.

{\rm (3)} Let $\mathscr{F}$ and $\mathscr{G}$ be finitely presented $\O^{\int}_{\mathscr{X}}$-modules, and $\varphi\colon\mathscr{F}\rightarrow\mathscr{G}$ a morphism of $\O^{\int}_{\mathscr{X}}$-modules.
Then there exist a distinguished formal model $X$ of $\mathscr{X}$, finitely presented $\O_X$-modules $\mathscr{F}_X$ and $\mathscr{G}_X$, and a morphism $\varphi_X\colon\mathscr{F}_X\rightarrow\mathscr{G}_X$ of $\O_X$-modules such that $\sp^{\ast}_X\varphi_X=\varphi$.
\end{thm}

Note that the formal model $X'$ in (2) is again distinguished (\ref{prop-cohrigidspacedist100}).
Here, for the proof of the theorem, and also for the later purpose, let us set up a useful situation.\begin{sit}\label{sit-cohsheavescoh}{\rm
Let $\mathscr{X}$ be a coherent rigid space, and $X$ a {\em distinguished} formal model. 
For notational convenience, we display the set $\AId_X$ of all admissible ideals of $\O_X$ in the form $\AId_X=\{\mathscr{J}_{\alpha}\,|\,\alpha\in L\}$, where $L$ is an index set.
The set $L$ is considered with the ordering defined by 
$$
\alpha\leq\beta\quad\Longleftrightarrow\quad\textrm{there exists $\gamma\in L$ such that}\ \mathscr{J}_{\beta}=\mathscr{J}_{\alpha}\mathscr{J}_{\gamma}.
$$
Note that the ordered set $L$ is isomorphic to $\AId^{\opp}_X$ and hence is a directed set\index{set!directed set@directed ---}\index{directed set} (\ref{prop-catblowups1} (2)).
The set $\AId_X$ has the trivial ideal $\O_X$ to which the corresponding element in $L$ is denoted by $0$, that is, $\mathscr{J}_0=\O_X$.
We write $X_{\alpha}=X_{\mathscr{J}_{\alpha}}$ (the admissible blow-up along the admissible ideal $\mathscr{J}_{\alpha}$) for each $\alpha$, and the corresponding morphism $X_{\alpha}\rightarrow X_0=X$ is denoted by $\pi_{0\alpha}$.
By \ref{prop-catblowups1} (2) the functor 
$$
L^{\opp}\longrightarrow\BL_X,\qquad\alpha\mapsto(\pi_{0\alpha}\colon X_{\alpha}\rightarrow X)
$$
is cofinal.
Note that by \ref{prop-cohrigidspacedist100} each $X_{\alpha}$ is again a distinguished formal model of $\mathscr{X}$.
For $\alpha,\beta\in L$ such that $\alpha\leq\beta$ we have the unique commutative diagram
$$
\xymatrix@-4ex{
X_{\alpha}\ar[ddr]_{\pi_{0\alpha}}&&X_{\beta}\ar[ddl]^{\pi_{0\beta}}\ar[ll]_{\pi_{\alpha\beta}}\\ \\ &X}
$$
by the universality of admissible blow-ups (\ref{prop-blowups1} (3)).
Again by the universality one see easily that the morphism $\pi_{\alpha\beta}$ is also an admissible blow-up; in fact, it is the admissible blow-up along $(\pi^{-1}_{0\alpha}\mathscr{J}_{\gamma})\O_{X_{\alpha}}$, where $\mathscr{J}_{\beta}=\mathscr{J}_{\alpha}\mathscr{J}_{\gamma}$.
Moreover, for $\alpha\leq\beta\leq\gamma$ we have $\pi_{\alpha\beta}\circ\pi_{\beta\gamma}=\pi_{\alpha\gamma}$.
Thus we get a filtered projective system $\{X_{\alpha},\pi_{\alpha\beta}\}$ consisting of admissible blow-ups among distinguished formal models of $\mathscr{X}$ indexed by the directed set $L$.
By \ref{prop-catblowups1} (2) we have
$$
\varprojlim_{\alpha\in L}X_{\alpha}=(\ZR{\mathscr{X}},\O^{\int}_{\mathscr{X}}).
$$
For brevity we write $\sp_{\alpha}\colon\ZR{\mathscr{X}}\rightarrow X_{\alpha}$ in place of $\sp_{X_{\alpha}}$.}
\end{sit}

\begin{proof}[Proof of Theorem {\rm \ref{thm-cohsheavesZR1}}]
In the situation as in \ref{sit-cohsheavescoh}, we may assume that $\ZR{\mathscr{X}}$ is the projective limit along the directed set $L$.
Then theorem follows from the formal results {\bf \ref{ch-pre}}.\ref{thm-injlimmodpf} and {\bf \ref{ch-pre}}.\ref{thm-projlimcohsch31}. 
\end{proof}

For a (not necessarily coherent) rigid space $\mathscr{X}$ we have the functor 
$$
\Mod^{\int}_{\mathscr{X}}\longrightarrow\Mod_{\mathscr{X}},\qquad
\mathscr{F}\longmapsto\mathscr{F}\otimes_{\O^{\int}_{\mathscr{X}}}\O_{\mathscr{X}}, 
$$
where the last $\O_{\mathscr{X}}$-module $\mathscr{F}\otimes_{\O^{\int}_{\mathscr{X}}}\O_{\mathscr{X}}$ is, by a slight abuse of notation, often denoted by $\mathscr{F}^{\rig}$ (resp.\ $\varphi^{\rig}$) (and similarly for morphisms in $\Mod^{\int}_{\mathscr{X}}$).

\begin{dfn}\label{dfn-integralmodelcohsheaf1}{\rm 
Let $\mathscr{X}$ be a coherent universally Noetherian rigid space.

(1) Let $\mathscr{F}$ be a coherent $\O_{\mathscr{X}}$-module.
An {\em integral model}\index{integral model} of $\mathscr{F}$ is a pair of the form $(\til{\mathscr{F}},\varphi)$, where $\til{\mathscr{F}}$ is a finitely presented $\O^{\int}_{\mathscr{X}}$-module and $\varphi$ is an isomorphism $\varphi\colon\til{\mathscr{F}}^{\rig}\stackrel{\sim}{\rightarrow}\mathscr{F}$.

(2) Let $\Phi\colon\mathscr{F}\rightarrow\mathscr{G}$ be a morphism of coherent $\O_{\mathscr{X}}$-modules.
An {\em integral model} of $\Phi$ is a datum of the form $(\til{\Phi}\colon\til{\mathscr{F}}\rightarrow\til{\mathscr{G}},\varphi,\psi)$ consisting of 
\begin{itemize}
\item a morphism $\til{\Phi}\colon\til{\mathscr{F}}\rightarrow\til{\mathscr{G}}$ of $\O^{\int}_{\mathscr{X}}$-modules, where $\til{\mathscr{F}}$ and $\til{\mathscr{G}}$ are finitely presented, 
\item and isomorphisms $\varphi\colon\til{\mathscr{F}}^{\rig}\stackrel{\sim}{\rightarrow}\mathscr{F}$ and 
$\psi\colon\til{\mathscr{G}}^{\rig}\stackrel{\sim}{\rightarrow}\mathscr{G}$
\end{itemize}
such that the following square is commutative:
$$
\xymatrix{\til{\mathscr{F}}^{\rig}\ar[d]_{\varphi}\ar[r]^{\til{\Phi}^{\rig}}&\til{\mathscr{G}}^{\rig}\ar[d]^{\psi}\\ \mathscr{F}\ar[r]_{\Phi}&\mathscr{G}\rlap{.}}
$$}
\end{dfn}

\begin{dfn}\label{dfn-cohsheavesZR200}{\rm 
Let $\mathscr{X}$ be a coherent rigid space, and $\mathscr{I}$ an ideal of definition of $\ZR{\mathscr{X}}$ of finite type.
A morphism $\varphi\colon\mathscr{F}\rightarrow\mathscr{G}$ of $\O^{\int}_{\mathscr{X}}$-modules is said to be a {\em weak isomorphism}\index{weak isomorphism}\index{weak!weak isomorphism@--- isomorphism} if there exists an integer $s\geq 0$ such that $\mathscr{I}^s\ker\varphi=0$ and $\mathscr{I}^s\coker\varphi=0$.}
\end{dfn}

\begin{thm}\label{thm-corcohsheavesZR2101}
Let $\mathscr{X}$ be a coherent universally Noetherian rigid space.
Then the functor 
$$
\Mod^{\int,\pf}_{\mathscr{X}}\stackrel{\cdot^{\rig}}{\longrightarrow}\Coh_{\mathscr{X}}, 
$$
where $\Mod^{\int,\pf}_{\mathscr{X}}$ denotes the category of finitely presented $\O^{\int}_{\mathscr{X}}$-modules, induces a categorical equivalence 
$$
\Mod^{\int,\pf}_{\mathscr{X}}\big/\{\textrm{weak isomorphisms}\}\stackrel{\sim}{\longrightarrow}\Coh_{\mathscr{X}},
$$
where the left-hand category is the localized category by the set of all weak isomorphisms.
Moreover, for a locally free $\O_{\mathscr{X}}$-module $\mathscr{F}$ one can find a locally free $\O^{\int}_{\mathscr{X}}$-module $\til{\mathscr{F}}$ of the same rank such that $\til{\mathscr{F}}^{\rig}\cong\mathscr{F}$. \hfill$\square$
\end{thm}

The theorem says that any morphism $\Phi$ of coherent $\O_{\mathscr{X}}$-modules has an integral model, uniquely determined up to weak isomorphism.
The proof can be done by combination of \ref{thm-existlattice} and \ref{thm-cohsheavesZR1} and by Exercise \ref{exer-corexistlatticelocfree}.
The details are left to the reader.
Notice that, since $\O_{\mathscr{X}}$ is coherent, any finitely presented $\O_{\mathscr{X}}$-module is coherent.

\addcontentsline{toc}{subsection}{Exercises}
\subsection*{Exercises}
\begin{exer}\label{exer-corlatticemodels41}
{\rm Under the notation as in \ref{prop-cohsheavesZR203}, suppose that $X$ is universally adhesive and that there exists an $\mathscr{I}_X$-torsion free $\O_X$-module $\mathscr{M}$ that contains $\mathscr{F}_X$ and $\mathscr{G}_X$ as $\O_X$-submodules.
Then show that both $\mathscr{F}_X+\mathscr{G}_X$ and $\mathscr{F}_X\cap\mathscr{G}_X$ are lattice models of $\mathscr{F}$.}
\end{exer}

\begin{exer}\label{exer-corexistlatticelocfree}
{\rm Let $\mathscr{X}$ be a coherent universally Noetherian rigid space, and $\mathscr{F}$ a locally free $\O_{\mathscr{X}}$-module of finite rank.
Then show that there exists a lattice model $(X,\mathscr{F}_X)$ of $\mathscr{F}$ such that $\mathscr{F}_X$ is a locally free $\O_X$-module of the same rank.}
\end{exer}


\section{Affinoids}\label{sec-affinoids}
\index{affinoid|(}
In this section we discuss the so-called {\em affinoids}, which are, by definition, coherent rigid spaces having an affine formal model. 
Our approach to the basic geometry of affinoids arises from the following viewpoint: Since an affinoid $\mathscr{X}$ has a formal model of the form $X=\Spf A$, which is the $I$-adic completion of the scheme $\Spec A$, the geometries of the rigid space $\mathscr{X}$ should reflect the geometries of scheme $\Spec A\setminus V(I)$, the complement of the closed subset defined by an ideal of definition.
In \S\ref{sub-morbetweenaffinoid} we will see this in the context of morphisms between affinoids; roughly speaking, if $A$ and $B$ are t.u.\ adhesive\index{t.u.a. ring@t.u.\ adhesive ring}, then morphisms $\mathscr{X}=(\Spf A)^{\rig}\rightarrow\mathscr{Y}=(\Spf B)^{\rig}$ between affinoids always come from morphisms of affine formal models of the form $\Spf A'\rightarrow\Spf B$, where $A'$ is finite and isomorphic outside $I$ over $A$.
In \S\ref{sub-affinoidscohsheaf} we will see that, in the context of coherent sheaves, the category of coherent sheaves on an affinoid $\mathscr{X}=(\Spf A)^{\rig}$, where $A$ is a t.u.\ rigid-Noetherian ring\index{t.u. rigid-Noetherian ring@t.u.\ rigid-Noetherian ring}, is equivalent to the category of coherent sheaves on the Noetherian scheme $\Spec A\setminus V(I)$.

The subsection \S\ref{sub-comparisonaffinoid} is devoted to the calculation of the cohomologies of coherent sheaves on affinoids, where we find that the cohomologies on an affinoid $\mathscr{X}=(\Spf A)^{\rig}$ can be calculated by means of the cohomologies of the corresponding coherent sheaf on the scheme $\Spec A\setminus V(I)$.
This fundamental result for the cohomology calculus on locally universally Noetherian rigid spaces leads us to the notion of {\em Stein affinoids}, which is an honest analogue of the classical Stein domains in complex analytic geometry or of affine schemes in algebraic geometry.
In fact, any locally universally Noetherian\index{rigid space!universally Noetherian rigid space@universally Noetherian ---} rigid spaces can have an open covering consisting of Stein affinoids and, moreover, such coverings, called {\em Stein affinoid coverings}, are cofinal in the set of all coverings.
Due to the basic results Theorem A and Theorem B proved in \S\ref{sub-steinaffinoid}, Stein affinoid coverings are Leray coverings for computing the cohomologies of coherent sheaves, and this fact, as in the classical complex analytic geometry, algebraic geometry, etc., provides the foundations for the cohomology calculus on rigid spaces. 
Note that, in the classical situation, every affinoid is a Stein affinoid.

In the final subsection \S\ref{sub-associatedschemes}, we focus on the comparison between universally Noetherian affinoids and their {\em associated schemes}, that is, the Noetherian schemes of the form `$\Spec A\setminus V(I)$' as above.
This will give us a useful bridge between local geometries of rigid spaces and of schemes.

\subsection{Affinoids and affinoid coverings}\label{sub-affinoidsintro}
\subsubsection{Definition and basic properties}\label{subsub-affinoidsintro}
\begin{dfn}\label{dfn-cohomologyrigidsp1}{\rm 
A coherent rigid space $\mathscr{X}$ is called an {\em affinoid} if there exists a formal model $(X,\phi)$ of $\mathscr{X}$ where $X$ is affine.}
\end{dfn}

We denote by $\ARf$ the full subcategory of $\CRf$ consisting of affinoids.
\begin{dfn}\label{dfn-cohomologyrigidsp12}{\rm 
Let $\mathscr{X}$ be a rigid space.

(1) An {\em affinoid open subspace}\index{affinoid!affinoid open subspace@--- open (subspace)} of $\mathscr{X}$ is an isomorphism class over $\mathscr{X}$ of objects $\mathscr{U}\hookrightarrow\mathscr{X}$ in the small site $\mathscr{X}_{\ad}$ (\ref{dfn-admissiblesite3gensmall}) such that $\mathscr{U}$ is an affinoid.
For a point $x\in\ZR{\mathscr{X}}$ an {\em affinoid neighborhood}\index{affinoid!affinoid neighborhood@--- neighborhood} of $x$ is an affinoid open subspace $\mathscr{U}\hookrightarrow\mathscr{X}$ such that the image of $\ZR{\mathscr{U}}$ contains $x$.

(2) An {\em affinoid covering}\index{affinoid!affinoid covering@--- covering} of $\mathscr{X}$ is a covering 
$$
\coprod_{\alpha\in L}\mathscr{U}_{\alpha}\longrightarrow\mathscr{X}
$$
of the site $\mathscr{X}_{\ad}$ such that each $\mathscr{U}_{\alpha}$ is an affinoid.}
\end{dfn}

As the following proposition shows, any rigid space has an open basis consisting of affinoid open subspaces:
\begin{prop}\label{prop-cohomologyrigidsp1}
Let $\mathscr{X}$ be a rigid space, $x\in\ZR{\mathscr{X}}$, and $\mathfrak{V}$ an open neighborhood of $x$ in $\ZR{\mathscr{X}}$.
Then there exists an affinoid neighborhood $\mathscr{U}\hookrightarrow\mathscr{X}$ of $x$ such that the image of $\ZR{\mathscr{U}}$ is contained in $\mathfrak{V}$.
\end{prop}

\begin{proof}
We may assume $\mathscr{X}$ is coherent.
Moreover, since the topology on $\ZR{\mathscr{X}}$ is generated by quasi-compact open subsets, we may assume that $\mathfrak{V}$ is a quasi-compact open neighborhood of $x$.
Take a formal model $X$ of $\mathscr{X}$ that admits a quasi-compact open subset $V\subseteq X$ such that $\sp^{-1}_X(V)=\mathfrak{V}$, and an affine open neighborhood $U\subseteq V$ of $\sp_X(x)$.
Then $\mathscr{U}=U^{\rig}$ admits an open immersion $\mathscr{U}\hookrightarrow\mathscr{X}$ enjoying the desired properties.
\end{proof}

An affinoid $\mathscr{X}$ is said to be {\em distinguished} if it is of the form $\mathscr{X}=(\Spf A)^{\rig}$ by an $I$-torsion free $A$, where $I\subseteq A$ is an ideal of definition.
If $\mathscr{X}=(\Spf A)^{\rig}$ is an affinoid, then considering the admissible blow-up along a finitely generated ideal of definition $I\subseteq A$, one finds that $\mathscr{X}$ is covered by distinguished affinoids.
Hence any rigid space has an open basis consisting of distinguished affinoid open subspaces:
\begin{cor}\label{cor-propaffinoidcov1}
Let $\mathscr{X}$ be a rigid space.
Then any covering of the $($coherent$)$ small admissible site $\mathscr{X}_{\ad}$ is refined by an affinoid covering consisting of distinguished affinoids. \hfill$\square$
\end{cor}

If a rigid space $\mathscr{X}$ is locally universally Noetherian\index{rigid space!universally Noetherian rigid space@universally Noetherian ---!locally universally Noetherian rigid space@locally --- ---} (resp.\ locally universally adhesive\index{rigid space!universally adhesive rigid space@universally adhesive ---!locally universally adhesive rigid space@locally --- ---}), then, clearly, $\mathscr{X}$ has an open basis consisting of affinoid open subspaces by affinoids of the form $(\Spf A)^{\rig}$, where $A$ is a t.u.\ rigid-Noetherian\index{t.u. rigid-Noetherian ring@t.u.\ rigid-Noetherian ring} (resp.\ t.u.\ adhesive)\index{t.u.a. ring@t.u.\ adhesive ring} ring {\rm ({\bf \ref{ch-formal}}.\ref{dfn-tuaringadmissible})}.
In this situation, in view of the following proposition, the difference between affinoids and distinguished affnoids is not important:
\begin{prop}\label{prop-cohomologyrigidsp0}
Let $\mathscr{X}=(\Spf A)^{\rig}$ be an affinoid, where $A$ is a t.u.\ rigid-Noetherian $($resp.\ t.u.\ adhesive$)$ ring, and $I\subseteq A$ an ideal of definition.
Then $A'=A/A_{\Itor}$ is again a t.u.\ rigid-Noetherian $($resp.\ t.u.\ adhesive$)$ ring, and we have $\mathscr{X}=(\Spf A')^{\rig}$.
In particular, $\mathscr{X}$ is a distinguished affinoid.
\end{prop}

\begin{proof}
By {\bf \ref{ch-pre}}.\ref{cor-propARconseq1-2} $A_{\Itor}$ is closed in $A$ with respect to the $I$-adic topology.
Hence $A'=A/A_{\Itor}$ is $I$-adically complete and hence is t.u.\ rigid-Noetherian (resp.\ t.u.\ adhesive).
Let $Y\rightarrow\Spf A$ be the admissible blow-up along $I$, and $Y'\rightarrow\Spf A'$ the admissible blow-up along $IA'$.
The latter map coincides with the strict transform\index{strict transform} of the former by the closed immersion $\Spf A'\hookrightarrow\Spf A$.
Since $I\O_Y$ is invertible, it follows that $Y'\cong Y$, that is, $\Spf A$ and $\Spf A'$ are dominated by a common admissible blow-up.
\end{proof}

\noindent
{\bf Convention.} {\sl In the sequel, whenever we discuss an affinoid $\mathscr{X}$ that we know by the context to be universally Noetherian $($resp.\ universally adhesive$)$, it is always supposed, unless otherwise clearly stated, to be of the form $\mathscr{X}=(\Spf A)^{\rig}$ by a t.u.\ rigid-Noetherian\index{t.u. rigid-Noetherian ring@t.u.\ rigid-Noetherian ring} $($resp.\ t.u.\ adhesive\index{t.u.a. ring@t.u.\ adhesive ring}$)$ ring $A$.}

\medskip
The above convention will not be necessary, if one can show the following statement (for which we do not know the proof): If $\mathscr{X}$ is universally Noetherian and $\mathscr{X}=(\Spf A)^{\rig}$, then $A$ is t.u.\ rigid-Noetherian.

\subsubsection{Affinoid subdomains}\label{subsub-affinoidsubdom}
\begin{dfn}\label{dfn-affinoidsubdomain1}{\rm 
An affinoid open subspace of an affinoid $\mathscr{X}$ is called an {\em affinoid subdomain}\index{affinoid!affinoid subdomain@--- subdomain} of $\mathscr{X}$.}
\end{dfn}

\begin{exas}[Examples of affinoid subdomains; cf.\ {\cite[(7.2.3/2), (7.2.3/5)]{BGR}}]\label{exas-affinoidsubdomain}{\rm Let $\mathscr{X}=(\Spf A)^{\rig}$ be an affinoid, and $I=(a)\subseteq A$ an ideal of definition.

(1) (Weierstrass subdomain\index{affinoid!affinoid subdomain@--- subdomain!weierstrass subdomain@Weierstrass subdomain}\index{Weierstrass subdomain@Weierstrass subdomain}) 
Let $f_1,\ldots,f_m\in\Gamma(\Spec A\setminus V(I),\O_{\Spec A})=A[\frac{1}{a}]$.
We can find a positive integer $k$ such that $a^kf_i\in A$ for $i=1,\ldots,m$.
Set $J=(a^k,a^kf_1,\ldots,a^kf_m)$.
Then $\mathscr{J}=J^{\Delta}$ is an admissible ideal of $X=\Spf A$.
Consider the admissible blow-up $\pi\colon X'\rightarrow X$ along $\mathscr{J}$, and let $U$ be the affine open part of $X'$ where $\mathscr{J}\O_U$ is generated by $a$:
$$
U=\Spf A\dl{\textstyle \frac{a^kf_1}{a^k},\ldots,\frac{a^kf_m}{a^k}}\dr/\Itor.
$$
The associated rigid space $U^{\rig}$ is an affinoid subdomain of $\mathscr{X}$, denoted by 
$$
\mathscr{X}(f)=\mathscr{X}(f_1,\ldots,f_m)
$$
and called a {\em Weierstrass subdomain}.

(2) (Laurent subdomain\index{affinoid!affinoid subdomain@--- subdomain!laurent subdomain@Laurent subdomain}\index{Laurent subdomain}) In the above situation, we moreover take $g_1,\ldots,g_n\in A[\frac{1}{a}]$ and $l\geq 1$ such that $a^lg_j\in A$ for $j=1,\ldots,n$.
Set $\mathscr{K}_j=(a^l,a^lg_j)^{\Delta}$, which are admissible ideals on $X$.
Define $U_j$ $(j=1,\ldots,n)$ inductively as follows: $U_0=U$; let $\pi_j\colon X'_j\rightarrow U_{j-1}$ be the admissible blow-up along the ideal $\mathscr{K}_j\O_{U_{j-1}}$, and $U_j$ $(j=1,\ldots,n)$ the affine open part of $X'_j$ on which $a^lg_j$ generates the ideal $\mathscr{K}_j\O_{X'_j}$.
Then the rigid space $U^{\rig}_n$ is an affinoid subdomain of $\mathscr{X}$, denoted by 
$$
\mathscr{X}(f,g^{-1})=\mathscr{X}(f_1,\ldots,f_m,g^{-1}_1,\ldots,g^{-1}_n)
$$
and called a {\em Laurent subdomain}.

(3) (Rational subdomain\index{affinoid!affinoid subdomain@--- subdomain!rational subdomain@rational subdomain}\index{rational subdomain}) Take elements $f_1,\ldots,f_m,g\in A[\frac{1}{a}]$ such that the ideal $(f_1,\ldots,f_m,g)$ of the ring $A[\frac{1}{a}]$ is the unit ideal.
As before, fix $k\geq 1$ such that $a^kf_1,\ldots,a^kf_m,a^kg\in A$.
Then the ideal $J=(a^kf_1,\ldots,a^kf_m,a^kg)$ is an open ideal, and hence $\mathscr{J}=J^{\Delta}$ is an admissible ideal on $X=\Spf A$.
Take the admissible blow-up $X'\rightarrow X$ along $\mathscr{J}$, and let $U$ be the affine part of $X'$ where the ideal $\mathscr{J}\O_{X'}$ is generated by $a^kg$:
$$
U=\Spf A\dl{\textstyle \frac{a^kf_1}{a^kg},\ldots,\frac{a^kf_m}{a^kg}}\dr/\Itor.
$$
The associated rigid space $U^{\rig}$ is an affinoid subdomain of $\mathscr{X}$, denoted by 
$$
\mathscr{X}({\textstyle \frac{f}{g}})=\mathscr{X}({\textstyle \frac{f_1}{g},\ldots,\frac{f_n}{g}})
$$
and called a {\em rational subdomain}.}
\end{exas}

\subsection{Morphisms between affinoids}\label{sub-morbetweenaffinoid}
\begin{dfn}\label{dfn-morbetweenaffinoid1}{\rm 
Let $A\rightarrow A'$ be an adic morphism of adic rings of finite ideal type, and $I\subseteq A$ an ideal of definition.
We say that the map $A\rightarrow A'$ (or $\Spf A'\rightarrow\Spf A$) is an {\em isomorphism outside $I$}\index{isomorphic outside I@isomorphic outside $I$} if the morphism $\Spec A'\setminus V(IA')\rightarrow\Spec A\setminus V(I)$ of schemes is an isomorphism.}
\end{dfn}

Let $A$ be an adic ring of finite ideal type with a finitely generated ideal of definition $I\subseteq A$, and $\pi\colon X'\rightarrow X=\Spf A$ an admissible blow-up.
We have a ring homomorphism $A\rightarrow A'=\Gamma(X',\O_{X'})$ (cf.\ Exercise \ref{exer-affinecovaluedpoint}).
We consider the $IA'$-adic topology on $A'$.

\begin{prop}\label{prop-affinoidblowupglobalsection}
Suppose $A$ is t.u.\ rigid-Noetherian\index{t.u. rigid-Noetherian ring@t.u.\ rigid-Noetherian ring}.
Then$:$
\begin{itemize}
\item[{\rm (1)}] $A\rightarrow A'$ is a weak isomorphism\index{weak!weak isomorphism@--- isomorphism} {\rm ({\bf \ref{ch-formal}}.\ref{dfn-weakisomorphic})}$;$
\item[{\rm (2)}] $A'$ endowed with the $IA'$-adic topology is rigid-Noetherian\index{rigid-Noetherian ring@rigid-Noetherian ring}  {\rm ({\bf \ref{ch-formal}}.\ref{dfn-tuaringadmissible} (1))}$;$
\item[{\rm (3)}] the morphism $X'\rightarrow\Spf A'$ is an admissible blow-up$;$ in particular, we have $(\Spf A)^{\rig}\cong(\Spf A')^{\rig}$.
\end{itemize}
Moreover, if $A$ is t.u.\ adhesive\index{t.u.a. ring@t.u.\ adhesive ring} {\rm ({\bf \ref{ch-formal}}.\ref{dfn-tuaringadmissible} (2))} and $I$-torsion free, then $A'$ is finite over $A$ $($and hence is again t.u.\ adhesive$)$.
\end{prop}

\begin{proof}
(1) It is clear that $A$ and $A'$ (due to GFGA comparison) are isomorphic outside $I$ (that is, $\Spec A'\setminus V(IA')\cong\Spec A\setminus V(I)$); in particular, $A'$ is Noetherian outside $I$\index{Noetherian!Noetherian outside I@--- outside $I$}.
Hence the kernel of $A\rightarrow A'$ is $I$-torsion; since $A_{\Itor}$ is bounded $I$-torsion, the kernel is also bounded $I$-torsion.
Let $J=(f_0,\ldots,f_r)\subseteq A$ be an admissible ideal that gives the admissible blow-up $\pi$, and $X'=\bigcup^r_{i=0}U_i$ the affine covering as in \S\ref{subsub-blowupsdescription}.
It is clear that the cokernel of each map $A\rightarrow\Gamma(U_i,\O_{U_i})$ is bounded $I$-torsion, and hence the cokernel of $A\rightarrow A'$ is also bounded $I$-torsion.

(2) Since it is already shown that $A'$ is Noetherian outside $IA'$, it suffices to show that $A'$ is $IA'$-adic complete, which follows from Exercise \ref{exer-weakisomcomplete}.

(3) By the universal mapping property (\ref{prop-blowups1} (3)) the morphism $X'\rightarrow\Spf A'$ coincides with the admissible blow-up along the admissible ideal $JA'$.
If $A$ is t.u.\ adhesive and $I$-torsion free, then by \ref{cor-tateacyclic} (1) the sheaf $\pi_{\ast}\O_{X'}$ is a coherent $\O_X$-module, and hence $A'=\Gamma(X,\pi_{\ast}\O_{X'})$ is a coherent $A$-module ({\bf \ref{ch-formal}}.\ref{thm-adicqcohpre1}).
\end{proof}

\begin{dfn}\label{dfn-strictweakisomorphism}{\rm 
Let $f\colon A\rightarrow A'$ be an adic morphism between adic rings of finite ideal type.
We say that $f$ is a {\em strict weak isomorphism}\index{weak!weak isomorphism@--- isomorphism!strict weak isomorphism@strict --- ---} if it is a weak isomorphism and induces an isomorphism $(\Spf A')^{\rig}\stackrel{\sim}{\rightarrow}(\Spf A)^{\rig}$ of rigid spaces.}
\end{dfn}

It follows from \ref{prop-cohomologyrigidsp0} that, if $A$ is a t.u.\ rigid-Noetherian ring and $I\subseteq A$ is an ideal of definition, then $A\rightarrow A'=A/A_{\Itor}$ is a strict weak isomorphism\index{weak!weak isomorphism@--- isomorphism!strict weak isomorphism@strict --- ---}.
The above proposition shows that, if $A$ is t.u.\ rigid-Noetherian, then any admissible blow-up $X'\rightarrow\Spf A$ induces the strict weak isomorphism\index{weak!weak isomorphism@--- isomorphism!strict weak isomorphism@strict --- ---} $A\rightarrow A'=\Gamma(X',\O_{X'})$.
\begin{prop}\label{prop-lemmorbetweenaffnoid1}
Let $p\colon X'=\Spf A'\rightarrow X=\Spf A$ be a morphism of affine adic formal schemes of finite ideal type, and suppose that $p$ is finite\index{morphism of formal schemes@morphism (of formal schemes)!finite morphism of formal schemes@finite ---} $({\bf \ref{ch-formal}}.\ref{dfn-finitemorform1})$ and isomorphic outside $I$.
Then $A\rightarrow A'$ is a strict weak isomorphism\index{weak!weak isomorphism@--- isomorphism!strict weak isomorphism@strict --- ---}.
In particular, a finite weak isomorphism between adic rings is a strict weak isomorphism.
\end{prop}

\begin{proof}
Let $q\colon Y'=\Spec A'\rightarrow Y=\Spec A$ be such that $\widehat{q}=p$.
Then by {\bf \ref{ch-formal}}.\ref{prop-finitemorform1} the morphism $q$ is finite.
Set $U=Y\setminus V(I)$, which is a quasi-compact open subset of $Y$.
By \cite[Premi\`ere partie, (5.7.12)]{RG} (included below in \ref{prop-birationalgeom02}) there exists an $q^{-1}(U)$-admissible blow-up\index{admissible!U-admissible blow-up@$U$-{---} blow-up}\index{blow-up!U-admissible blow-up@$U$-admisible ---} $W\rightarrow Y'$ such that the composition $W\rightarrow Y'\rightarrow Y$ is an $U$-admissible blow-up.\index{admissible!U-admissible blow-up@$U$-{---} blow-up}\index{blow-up!U-admissible blow-up@$U$-admisible ---}
Taking the $I$-adic completions, we get the admissible blow-ups $Z=\widehat{W}\rightarrow\widehat{Y'}=X'$ and $Z\rightarrow\widehat{Y}=X$.
\end{proof}

\begin{prop}\label{prop-lemmorbetweenaffinoid2}
Consider a morphism between affinoids
$$
\varphi\colon (\Spf A)^{\rig}\rightarrow(\Spf B)^{\rig},
$$
where $A$ and $B$ are t.u.\ rigid-Noetherian\index{t.u. rigid-Noetherian ring@t.u.\ rigid-Noetherian ring} rings.

{\rm (1)} There exists a diagram
$$
\xymatrix@C-2ex@R-2ex{&A'\\ A\ar[ur]&&B\ar[ul]}
$$
consisting of adic morphisms of adic rings such that$:$
\begin{itemize}
\item[{\rm (a)}] $A'$ is a rigid-Noetherian\index{rigid-Noetherian ring@rigid-Noetherian ring} ring$;$
\item[{\rm (b)}] $A\rightarrow A'$ is a strict weak isomorphism\index{weak!weak isomorphism@--- isomorphism!strict weak isomorphism@strict --- ---}$;$
\item[{\rm (c)}] the induced morphism $(\Spf A)^{\rig}\stackrel{\sim}{\leftarrow}(\Spf A')^{\rig}\rightarrow(\Spf B)^{\rig}$ coincides with the morphism $\varphi$.
\end{itemize}
If, moreover, $\varphi$ is an isomorphism, then the above diagram can be taken such that $B\rightarrow A'$ is also a strict weak isomorphism.

{\rm (2)} If $A$ is t.u.\ adhesive\index{t.u.a. ring@t.u.\ adhesive ring}, and is $I$-torsion free, where $I\subseteq A$ is an ideal of definition, then the above diagram can be taken such that $A\rightarrow A'$ is a finite weak isomorphism $($and hence $A'$ is again t.u.\ adhesive\index{t.u.a. ring@t.u.\ adhesive ring}$)$.
If, moreover, $\varphi$ is an isomorphism and $B$ is a $J$-torsion free t.u.\ adhesive\index{t.u.a. ring@t.u.\ adhesive ring} ring $($where $J\subseteq B$ is an ideal of definition$)$, then the diagram can be taken such that $B\rightarrow A'$ is also a finite weak isomorphism.
\end{prop}

\begin{proof}
The morphism $\varphi$ comes from a diagram 
$$
\xymatrix@-2ex{
&X'\ar[dr]\ar[dl]\\
\Spf A&&\Spf B\rlap{,}}
$$
where the left-hand arrow is an admissible blow-up.
Set $A'=\Gamma(X',\O_{X'})$.
Then, as we have seen in \ref{prop-affinoidblowupglobalsection}, the ring $A'$ with the $IA'$-adic topology is a rigid-Noetherian ring, and the map $A\rightarrow A'$ is a strict weak isomorphism\index{weak!weak isomorphism@--- isomorphism!strict weak isomorphism@strict --- ---}, since both $\Spf A$ and $\Spf A'$ are dominated by $X'$.
Moreover, we have the canonical map $B\rightarrow A'$ (cf.\ Exercise \ref{exer-affinecovaluedpoint}).
If $\varphi$ is an isomorphism, then one can take a diagram as above such that both $X'\rightarrow\Spf A$ and $X'\rightarrow\Spf B$ are admissible blow-ups (\ref{prop-cohrigidspace3}).
Hence, in this case, $B\rightarrow A'$ is also a strict weak isomorphism.
All the rest of the assertion follow from \ref{prop-affinoidblowupglobalsection}.
\end{proof}

\begin{prop}\label{prop-corlemtuavsrigidaff}
Let $A$ be an $I$-torsion free t.u.\ rigid-Noetherian\index{t.u. rigid-Noetherian ring@t.u.\ rigid-Noetherian ring}, where $I\subseteq A$ is an ideal of definition of $A$, such that $\Spec A\setminus V(I)$ is affine.
Set $\mathscr{X}=(\Spf A)^{\rig}$ and $\Spec B=\Spec A\setminus V(I)$.
Then the integral closure $A^{\int}$ of $A$ in $B$ is canonically isomorphic to $\Gamma(\mathscr{X},\O^{\int}_{\mathscr{X}})$.
\end{prop}

The proposition follows immediately from \ref{prop-lemmorbetweenaffnoid1} and \ref{prop-lemmorbetweenaffinoid2} if the ring $A$ is t.u.\ adhesive\index{t.u.a. ring@t.u.\ adhesive ring}.
The general case follows from GFGA comparison and the following observation: 
Let $A$ be a ring with a finitely generated ideal $I\subseteq A$.
We assume that $A$ is $I$-torsion free, and $U=X\setminus V(I)$ where $X=\Spec A$ is affine; we set $B=\Gamma(U,\O_U)$.
We need to compare the integral closure $A^{\int}$ of $A$ in $B$ with the ring 
$$
\til{A}=\varinjlim_{X'\rightarrow X}\Gamma(X',\O_{X'}),
$$
where $X'$ runs over all $U$-admissible blow-ups\index{admissible!U-admissible blow-up@$U$-{---} blow-up}\index{blow-up!U-admissible blow-up@$U$-admisible ---} (cf.\ \ref{dfn-Uadmissibleblowups} below) of $X$:
\begin{lem}\label{lem-corlemtuavsrigidaff}
We have $\til{A}=A^{\int}$.
\end{lem}

\begin{proof}
The inclusion $A^{\int}\subseteq\til{A}$ can be shown by an argument similar to that in the proof of \ref{prop-lemmorbetweenaffnoid1}.
For any $f\in\til{A}$, take a $U$-admissible blow-up $\pi\colon X'\rightarrow X$ such that $I\O_{X'}$ is invertible and $f\in\Gamma(X',\O_{X'})$.
Set $Z=\Spec A[f]$.
We want to show that $Z\rightarrow X$ is finite.
Since $A[f]$ is $I$-torsion free, $U$ is dence in $Z$.
We have a factorization $X'\stackrel{h}{\rightarrow}Z\stackrel{g}{\rightarrow}X$ of $\pi$.
Now since $h$ is proper and hence $h(X')$ is closed containing $U$, we have $h(X')=Z$, that is, $h$ is surjective.
Hence $g$ is proper.
Since $g$ is proper and affine, it is finite (cf.\ \ref{prop-appclassicalZR} below), as desired.
\end{proof}

\begin{thm}\label{thm-morbetweenaffinoid1}
Let $A$ and $B$ be adic rings with finitely generated ideals of definition $I\subseteq A$ and $J\subseteq B$, respectively.
Suppose that $A$ and $B$ are t.u.\ rigid-Noetherian\index{t.u. rigid-Noetherian ring@t.u.\ rigid-Noetherian ring} $($resp.\ t.u.\ adhesive\index{t.u.a. ring@t.u.\ adhesive ring} and that $A$ is $I$-torsion free and $B$ is $J$-torsion free$)$.
Then the canonical map
\begin{equation*}
\Bigg\{\begin{minipage}{5.7em}{$\xymatrix@C-4ex@R-3ex{&A'\\ A\ar[ur]^p&&B\ar[ul]_q}$}\end{minipage}
\,\Bigg|\,
\begin{minipage}{9.3em}
{\small {\rm $p$ is a strict (resp.\ finite) weak isomorphism\index{weak!weak isomorphism@--- isomorphism!strict weak isomorphism@strict --- ---}, and $q$ is an adic morphism.}}
\end{minipage}
\Bigg\}\big/\sim\ \longrightarrow\ \Hom_{\CRf}(\mathscr{X},\mathscr{Y}).
\end{equation*}
is a bijection.
Here the equivalence relation $\sim$ on the left-hand side is defined as follows$:$
$(A\rightarrow A'\leftarrow B)$ and $(A\rightarrow A''\leftarrow B)$ are equivalent if there exists third $(A\rightarrow A'''\leftarrow B)$ as above sitting in the commutative diagram
$$
\xymatrix@-2ex{
&A'\\
A\ar[ur]\ar[r]\ar[dr]&A'''\ar[u]\ar[d]&B\ar[ul]\ar[l]\ar[dl]\\
&A''\rlap{.}}
$$
\end{thm}

\begin{proof}
In view of \ref{prop-lemmorbetweenaffinoid2} we only need to show the injectivity.
Suppose we are given a diagram 
$$
\xymatrix@-3ex{
&A'\\
A\ar[ur]\ar[dr]&&B\ar[ul]\ar[dl]\\
&A''}
$$
where the left-hand maps are strict (resp.\ finite) weak isomorphisms.
Since $(\Spf A')^{\rig}=(\Spf A'')^{\rig}$, there exists admissible blow-ups $X'''\rightarrow\Spf A'$ and $X'''\rightarrow\Spf A''$.
Moreover, one can choose such an $X'''$ that admits an admissible blow-up $X'''\rightarrow\Spf A$ compatible with the above two admissible blow-ups.
Now, applying \ref{prop-affinoidblowupglobalsection}, we have the diagram as in the theorem.
\end{proof}

\subsection{Coherent sheaves on affinoids}\label{sub-affinoidscohsheaf}
Let $\mathscr{X}=(\Spf A)^{\rig}$ be an affinoid where $A$ is a t.u.\ rigid-Noetherian ring\index{t.u. rigid-Noetherian ring@t.u.\ rigid-Noetherian ring}, and $I\subseteq A$ a finitely generated ideal of definition.
Set $X=\Spf A$.
Recall that there exists a canonical exact categorical equivalence
$$
\Mod^{\pf}_A\stackrel{\sim}{\longrightarrow}\Mod^{\pf}_X, \qquad M\longmapsto M^{\Delta}
$$
({\bf \ref{ch-formal}}.\ref{thm-adicqcohpre1}).
On the other hand, if we set $Y=\Spec A$, then by \cite[(1.4.2), (1.4.3)]{EGAInew} we have an exact categorical equivalence
$$
\Mod^{\pf}_A\stackrel{\sim}{\longrightarrow}\Mod^{\pf}_Y, \qquad M\longmapsto \til{M}.
$$
Thus we have an exact equivalence of the categories
$$
\Mod^{\pf}_X\stackrel{\sim}{\longrightarrow}\Mod^{\pf}_Y, \qquad M^{\Delta}\longmapsto \til{M}.
$$

Set $U=Y\setminus V(I)$, and consider the diagram of categories
$$
\Coh_{\mathscr{X}}\stackrel{\cdot^{\rig}}{\longleftarrow}\Mod^{\pf}_X\stackrel{\sim}{\longrightarrow}\Mod^{\pf}_Y\longrightarrow\Coh_U,
$$
where the last functor is given by the restriction to $U$, $\mathscr{F}\mapsto\mathscr{F}|_U$.
By \ref{thm-tateacyclic20}, for any object $\mathscr{F}$ (resp.\ arrow $\varphi$) in $\Coh_{\mathscr{X}}$ there exists up to weak isomorphisms an object $\mathscr{F}_X$ (resp.\ an arrow $\varphi_X$) in $\Mod^{\pf}_X$ such that $\mathscr{F}^{\rig}_X=\mathscr{F}$ (resp.\ $\varphi^{\rig}_X=\varphi$), which gives rise to an object $\mathscr{F}_U$ (resp.\ arrow $\varphi_U$) in $\Coh_U$.
It is clear that this gives a well-defined functor 
$$
\Coh_{\mathscr{X}}\longrightarrow\Coh_U.\leqno{(\ast)}
$$
\begin{thm}\label{thm-affinoidscohsheaf}
The functor $(\ast)$ is an exact categorical equivalence. 
\end{thm}

\begin{proof}
By \ref{thm-tateacyclic20} it suffices to show the equivalence
$$
\Mod^{\pf}_Y/\{\textrm{weak isomorphisms}\}\stackrel{\sim}{\longrightarrow}\Coh_U, 
$$
where we mean by a weak isomorphism\index{weak isomorphism} a morphism $\varphi\colon \mathscr{F}\rightarrow\mathscr{G}$ in $\Mod^{\pf}_Y$ whose kernel and cokernel are annihilated by $I^s$ for some $s>0$.
But then, only essential surjectivity calls for the proof, which follows from \cite[(6.9.11)]{EGAInew}.
\end{proof}

\begin{ntn}\label{ntn-affinoidscohsheaf}{\rm 
The proof of \ref{thm-affinoidscohsheaf} shows that we have a quasi-inverse functor 
$$
\Coh_U\longrightarrow\Coh_{\mathscr{X}}
$$
defined similarly to $(\ast)$, denoted by $\mathscr{F}\mapsto\mathscr{F}^{\rig}$ by a slight abuse of notation.}
\end{ntn}

In general, for a ringed space $(X,\O_X)$ with the coherent structure sheaf, a coherent sheaf $\mathscr{F}$ is said to be {\em Noetherian}\index{Noetherian!Noetherian sheaf@--- sheaf} if the following condition is satisfied: for any point $x\in X$ there exists an open neighborhood $U$ such that any increasing sequence $\mathscr{G}_1\subseteq\mathscr{G}_2\subseteq\cdots$ of coherent subsheaves of $\mathscr{F}|_U$ terminates, that is, there exists a number $N$ such that $\mathscr{G}_N=\mathscr{G}_{N+1}=\cdots$\footnote{Note that the notion of Noetherness for sheaves as defined here is different from the one in \cite{KS}, Definition 11.1.1 (iii). Note also that, if the topological space $X$ is locally coherent, then the open set $U$ in the definition can be arbitrary coherent open neighborhood of $x$.}.
By \ref{thm-affinoidscohsheaf} we immediately have the following:
\begin{cor}\label{cor-structuresheafNoetherian}
Let $\mathscr{X}$ be a locally universally Noetherian rigid space\index{rigid space!universally Noetherian rigid space@universally Noetherian ---!locally universally Noetherian rigid space@locally --- ---}.
Then the sheaf $\O_{\mathscr{X}}$ is Noetherian. \hfill$\square$
\end{cor}

\subsection{Comparison theorem for affinoids}\label{sub-comparisonaffinoid}
\begin{thm}[Comparison theorem for affinoids]\label{thm-comparisonaffinoid}
Let $\mathscr{X}=(\Spf A)^{\rig}$ be an affinoid\index{affinoid} where $A$ is a t.u.\ rigid-Noetherian ring\index{t.u. rigid-Noetherian ring@t.u.\ rigid-Noetherian ring} {\rm ({\bf \ref{ch-formal}}.\ref{dfn-tuaringadmissible})}, and set $X=\Spf A$.
Consider the affine scheme $Y=\Spec A$, and set $U=Y\setminus V(I)$, where $I\subseteq A$ is a finitely generated ideal of definition.
Then for any coherent sheaf $\mathscr{F}$ on the Noetherian scheme $U$, we have a canonical isomorphism
$$
\H^q(\mathscr{X},\mathscr{F}^{\rig})\cong\H^q(U,\mathscr{F})
$$
for each $q\geq 0$.
\end{thm}

\begin{proof}
In view of \ref{prop-cohomologyrigidsp0} we may assume without loss of generality that $A$ is $I$-torsion free.
Let $\mathscr{I}_X=I^{\Delta}$, and $\mathscr{I}=(\sp^{-1}_X\mathscr{I}_X)\O^{\int}_{\mathscr{X}}$.
Let $\mathscr{G}_Y$ be a finitely presented sheaf on $Y$ such that $\mathscr{G}_Y|_U=\mathscr{F}$ (\cite[(6.9.11)]{EGAInew}).
The formal completion $\mathscr{F}_X=\widehat{\mathscr{G}}_Y$ gives a finitely presented formal model of $\mathscr{F}^{\rig}$ on $X$.
Set $\til{\mathscr{F}}=\sp^{\ast}_X\mathscr{F}_X$, which is an $\O^{\int}_{\mathscr{X}}$-module on $\ZR{\mathscr{X}}$.
By \ref{lem-cohsheavesZR01lem} we have
$$
\H^q(\mathscr{X},\mathscr{F}^{\rig})=\H^q(\mathscr{X},\varinjlim_{n\geq 0}\lHom(\mathscr{I}^n,\til{\mathscr{F}}))=\varinjlim_{n\geq 0}\H^q(\mathscr{X},\lHom(\mathscr{I}^n,\til{\mathscr{F}})),
$$
where the last equality is due to {\bf \ref{ch-pre}}.\ref{prop-directlimcohcoh11}.

Set $\mathscr{I}_Y=\til{I}$, which is a coherent ideal of $\O_Y$.
Let $\pi\colon X'\rightarrow X$ be the admissible blow-up along $\mathscr{I}_X$, and $p\colon Y'\rightarrow Y$ the blow-up of the scheme $Y$ along $\mathscr{I}_Y$.
Then $\pi$ is the formal completion of $p$.
Set $\mathscr{I}_{X'}=(\pi^{-1}\mathscr{I}_X)\O_{X'}$, which is an invertible ideal of definition.
Starting with $X_0=X'$, we construct the system $\{X_{\alpha}\}_{\alpha\in L}$ as in \ref{sit-cohsheavescoh}.
For any $\alpha$ one can take the corresponding blow-up $Y_{\alpha}\rightarrow Y'=Y_0$ whose formal completion is $X_{\alpha}\rightarrow X'=X_0$.
Moreover, for $\alpha\leq\beta$ there exists a morphism $p_{\alpha\beta}\colon Y_{\beta}\rightarrow Y_{\alpha}$ of schemes (in fact, a blow-up) whose formal completion coincides with $\pi_{\alpha\beta}$.
For any $\alpha\in L$, let $\mathscr{I}_{\alpha}=(\pi^{-1}_{0\alpha}\mathscr{I}_{X'})\O_{X_{\alpha}}$, and $\mathscr{F}_{\alpha}$ the pull-back of $\mathscr{F}_X$, which is the formal completion of the sheaf $\mathscr{G}_{\alpha}$, the pull-back of $\mathscr{G}_Y$ by the map $Y_{\alpha}\rightarrow Y$.
Then by \ref{cor-ZRstrsheaf2131} we have $\mathscr{I}=\sp_{\alpha}^{\ast}\mathscr{I}_{\alpha}$ for any $\alpha\in L$.
By \cite[$\mathbf{0}$, (5.4.9)]{EGAInew} and {\bf \ref{ch-pre}}.\ref{cor-ringedspacecohcoherentprojectivesystemindlimits01} we calculate
\begin{equation*}
\begin{split}
\varinjlim_{n\geq 0}\H^q(\mathscr{X},\lHom(\mathscr{I}^n,\til{\mathscr{F}}))&=\varinjlim_{n\geq 0}\H^q(\mathscr{X}, \sp^{\ast}_0\lHom_{\O_{X'}}(\mathscr{I}^n_0,\mathscr{F}_0))\\
&=\varinjlim_{n\geq 0}\varinjlim_{\alpha\in L}\H^q(X_{\alpha},\lHom_{\O_{X_{\alpha}}}(\mathscr{I}^n_{\alpha},\mathscr{F}_{\alpha})).
\end{split}
\end{equation*}
By comparison theorem ({\bf \ref{ch-formal}}.\ref{thm-GFGAcomweakcoherentspecial}) we have
\begin{equation*}
\begin{split}
\varinjlim_{n\geq 0}\varinjlim_{\alpha\in L}\H^q(X_{\alpha},\lHom_{\O_{X_{\alpha}}}(\mathscr{I}^n_{\alpha},\mathscr{F}_{\alpha}))&=\varinjlim_{\alpha\in L}\varinjlim_{n\geq 0}\H^q(Y_{\alpha},\lHom_{\O_{Y_{\alpha}}}(\mathscr{J}^n_{\alpha},\mathscr{G}_{\alpha}))\\
&=\varinjlim_{\alpha\in L}\H^q(Y_{\alpha},\varinjlim_{n\geq 0}\lHom_{\O_{Y_{\alpha}}}(\mathscr{J}^n_{\alpha},\mathscr{G}_{\alpha})),
\end{split}
\end{equation*}
where $\mathscr{J}_{\alpha}$ is the invertible ideal of $\O_{Y_{\alpha}}$ corresponding to $\mathscr{I}_{\alpha}$.

Set $U_{\alpha}=Y_{\alpha}\setminus V(\mathscr{J}_{\alpha})$.
Since $Y_{\alpha}$ is $\mathscr{J}_{\alpha}$-torsion free, $U_{\alpha}$ is a non-empty scheme.
By Deligne's formula (cf.\ Exercise \ref{exer-deligne}) we have
$$
\varinjlim_{n\geq 0}\lHom_{\O_{Y_{\alpha}}}(\mathscr{J}^n_{\alpha},\mathscr{G}_{\alpha})=j_{\alpha\ast}j^{\ast}_{\alpha}\mathscr{G}_{\alpha},
$$
where $j_{\alpha}\colon U_{\alpha}\hookrightarrow Y_{\alpha}$ is the open immersion.
Since the blow-up $Y_{\alpha}\rightarrow Y=\Spec A$ is isomorphic on $U=\Spec A\setminus V(I)$, and since $U\hookrightarrow Y_{\alpha}$ is an affine morphism, we have by {\bf \ref{ch-formal}}.\ref{thm-formalacyclicity} (2)
$$
\varinjlim_{\alpha\in L}\H^q(Y_{\alpha},\varinjlim_{n\geq 0}\lHom_{\O_{Y_{\alpha}}}(\mathscr{J}^n_{\alpha},\mathscr{G}_{\alpha}))=\H^q(U,\mathscr{G}_Y|_U).
$$
Consequently, we have $\H^q(\mathscr{X},\mathscr{F}^{\rig})=\H^q(U,\mathscr{F})$, as desired.
\end{proof}

\begin{cor}\label{cor-reducedrigidspaces}
Let $\mathscr{X}$ be a coherent universally Noetherian rigid space\index{rigid space!universally Noetherian rigid space@universally Noetherian ---}.
Then the following conditions are equivalent$:$
\begin{itemize}
\item[{\rm (a)}] $\mathscr{X}$ is reduced\index{rigid space!reduced rigid space@reduced ---} $(\ref{dfn-reducedrigidspaces});$
\item[{\rm (b)}] there exists a cofinal collection $\{X_{\lambda}\}_{\lambda\in\Lambda}$ of formal models of $\mathscr{X}$ consisting of reduced distinguished formal models\index{formal model!formal model of a coherent rigid space@--- (of a coherent rigid space)!distinguished formal model of a coherent rigid space@distinguished --- ---} {\rm (\ref{dfn-cohrigidspacedist})}$;$
\item[{\rm (c)}] any distinguished formal model of $\mathscr{X}$ is reduced.
\end{itemize}
\end{cor}

Here a formal scheme is said to be {\it reduced} if it is reduced as a ringed space; see {\bf \ref{ch-pre}}, \S\ref{subsub-ringedsplocalringedsp}.

\begin{proof}
Suppose $\mathscr{X}$ is reduced, and let $X$ be a distinguished formal model.
For any affine open subset $\Spf A$ of $X$, the ring $\H^0(\Spec A\setminus V(I),\O_{\Spec A})$ (where $I$ is an ideal of definition) is reduced due to \ref{thm-comparisonaffinoid}.
Hence the nilradical of $A$ is an $I$-torsion ideal, which is equal to zero, since $A$ is $I$-torsion free.
Thus we have (a) $\Rightarrow$ (c).
The implication (c) $\Rightarrow$ (b) is trivial.
If (b) holds, then by \ref{thm-comparisonaffinoid} one shows that for any point $x\in\ZR{\mathscr{X}}$ the stalk $\O_{\mathscr{X},x}$ is a reduced ring (since the filtered inductive limit of reduced rings is reduced), whence (b) $\Rightarrow$ (a).
\end{proof}

\subsection{Stein affinoids}\label{sub-steinaffinoid}
\index{affinoid!Stein affinoid@Stein ---|(}\index{Stein affinoid|(}
\subsubsection{Stein affinoids and Stein affinoid coverings}\label{subsub-steinaffinoids}
\begin{prop}\label{prop-serrecriterion}
Let $\mathscr{X}=(\Spf A)^{\rig}$ be an affinoid where $A$ is a t.u.\ rigid-Noetherian ring\index{t.u. rigid-Noetherian ring@t.u.\ rigid-Noetherian ring}.
Set $X=\Spf A$.
Then the following conditions are equivalent$:$
\begin{itemize}
\item[{\rm (a)}] $\H^1(\mathscr{X},\mathscr{F})=0$ for any coherent $\O_{\mathscr{X}}$-module $\mathscr{F}$ $($that is, $\mathscr{X}$ is `Stein'$);$
\item[{\rm (b)}] $\H^q(\mathscr{X},\mathscr{F})=0$ for $q\geq 1$ and for any coherent $\O_{\mathscr{X}}$-module $\mathscr{F};$
\item[{\rm (c)}] $\Spec A\setminus V(I)$ is an affine scheme, where $I$ is a finitely generated ideal of definition of $A;$
\item[{\rm (d)}] there exists an $I'$-torsion free t.u.\ rigid-Noetherian ring $A'$ $($where $I'\subseteq A'$ is a finitely generated ideal of definition$)$ such that $\mathscr{X}=(\Spf A')^{\rig}$ and that $\Spec A'\setminus V(I')$ is affine.
\end{itemize}
\end{prop}

\begin{proof}
First we show (c) $\Rightarrow$ (d).
Set $A'=A/A_{\Itor}$.
Then by \ref{prop-cohomologyrigidsp0} we have $\mathscr{X}=(\Spf A')^{\rig}$.
Since $\Spec A'\hookrightarrow\Spec A$ is affine, $\Spec A'\setminus V(IA')$ is affine, whence (d).

Next, let us show (d) $\Rightarrow$ (b).
We work with the notation as in the proof of \ref{thm-comparisonaffinoid}.
Take a coherent sheaf $\mathscr{H}$ on the Noetherian affine scheme $U'=\Spec A'\setminus V(I')$ such that $\mathscr{F}=\mathscr{H}^{\rig}$; this is possible due to \ref{thm-affinoidscohsheaf}.
Then by \ref{thm-comparisonaffinoid} we have $\H^q(\mathscr{X},\mathscr{F})=\H^q(U',\mathscr{H})$.
If $q>0$, then the right-hand cohomology is zero due to {\bf \ref{ch-pre}}.\ref{thm-vanishcohaff-0} (1), whence (b).

The implication (b) $\Rightarrow$ (a) is clear.
It remains to show (a) $\Rightarrow$ (c).
Set $U=\Spec A\setminus V(I)$.
We are to show that the Noetherian scheme $U$ is affine.
To this end, by Serre's criterion (\cite[$\mathbf{II}$, \S5.2]{EGA}) it suffices to show that for any coherent sheaf $\mathscr{H}$ on $U$ its first cohomology group on $U$ vanishes.
But this follows from \ref{thm-affinoidscohsheaf} and \ref{thm-comparisonaffinoid}.
\end{proof}

\danger{It is important to assume in \ref{prop-serrecriterion} that $\mathscr{X}$ is an affinoid; in fact, there exists an example of non-affinoid rigid space on which the higher cohomologies of any coherent sheaf vanish \cite{Liu}.}

\begin{dfn}\label{dfn-cohomologyrigidsp1str}{\rm 
A coherent universally Noetherian rigid space\index{rigid space!universally Noetherian rigid space@universally Noetherian ---} $\mathscr{X}$ is called a {\em Stein affinoid} if it is an affinoid of the form $(\Spf A)^{\rig}$ by a t.u.\ rigid-Noetherian ring $A$ that satisfies the conditions in \ref{prop-serrecriterion}.}
\end{dfn}

\begin{exa}\label{exa-nonStein}{\rm 
Let $R$ be a integrally closed Noetherian local domain of dimension $\geq 2$, and set $I=\m_R$.
Suppose $R$ is $I$-adically complete.
Then $R$ is a t.u.\ adhesive ring, and the rigid space $\mathscr{X}=(\Spf R)^{\rig}$ is an affinoid, but is not a Stein affinoid.}
\end{exa}

\begin{prop}\label{prop-cohomologyrigidsp1str}
Let $\mathscr{X}$ be a locally universally Noetherian\index{rigid space!universally Noetherian rigid space@universally Noetherian ---!locally universally Noetherian rigid space@locally --- ---} rigid space.
Then for any point $x\in\ZR{\mathscr{X}}$ there exists an affinoid open neighborhood $\mathscr{U}\hookrightarrow\mathscr{X}$ by a Stein affinoid $\mathscr{U}$.
Moreover, such affinoid open neighborhoods are cofinal in the set of all open neighborhoods of $x$.
\end{prop}

\begin{proof}
In the proof of \ref{prop-cohomologyrigidsp1}, replace $X$ by the admissible blow-up along an ideal of definition.
Since the ideal of definition of $X$ is invertible, the scheme $\Spec A\setminus V(I)$ for any affine open subset $U=\Spf A$ of $X$, where $I\subseteq A$ is the ideal of definition, is affine.
\end{proof}

Clearly, if the rigid space $\mathscr{X}$ in the situation in \ref{prop-cohomologyrigidsp1str} is locally universally adhesive\index{rigid space!universally adhesive rigid space@universally adhesive ---!locally universally adhesive rigid space@locally --- ---}, then one can take $\mathscr{U}$ as above to be of the form $(\Spf A)^{\rig}$ by a t.u.\ adhesive ring\index{t.u.a. ring@t.u.\ adhesive ring} $A$.

\begin{dfn}\label{dfn-cohomologyrigidsp12str}{\rm 
Let $\mathscr{X}$ be a locally universally Noetherian\index{rigid space!universally Noetherian rigid space@universally Noetherian ---!locally universally Noetherian rigid space@locally --- ---} rigid space.

(1) A {\em Stein affinoid open subspace}\index{Stein affinoid!Stein affinoid open subspace@--- open (subspace)} of $\mathscr{X}$ is an isomorphism class over $\mathscr{X}$ of objects $\mathscr{U}\hookrightarrow\mathscr{X}$ in the small site $\mathscr{X}_{\ad}$ (\ref{dfn-admissiblesite3gensmall}) such that $\mathscr{U}$ is a Stein affinoid.

(2) A {\em Stein affinoid covering}\index{Stein affinoid!Stein affinoid covering@--- covering} of $\mathscr{X}$ is a covering 
$$
\coprod_{\alpha\in L}\mathscr{U}_{\alpha}\longrightarrow\mathscr{X}
$$
of the site $\mathscr{X}_{\ad}$ such that each $\mathscr{U}_{\alpha}$ is a Stein affinoid.}
\end{dfn}

By the second assertion of \ref{prop-cohomologyrigidsp1str} we readily see:
\begin{prop}\label{prop-affinoidcov1str}
Let $\mathscr{X}$ be a locally universally Noetherian\index{rigid space!universally Noetherian rigid space@universally Noetherian ---!locally universally Noetherian rigid space@locally --- ---} rigid space.
Then any admissible covering of $\mathscr{X}$ is refined by a Stein affinoid covering. \hfill$\square$
\end{prop}

\subsubsection{Theorem A and Theorem B}\label{subsub-theoremAB}
\begin{thm}[Theorem A and Theorem B]\label{thm-cohomologyrigidsp1str}
Let $\mathscr{X}$ be a Stein affinoid, and $\mathscr{F}$ a coherent $\O_{\mathscr{X}}$-module.

{\rm (1)} The sheaf $\mathscr{F}$ is generated by global sections$;$ if $\mathscr{X}=(\Spf A)^{\rig}$, where $A$ is an $I$-torsion free t.u.\ rigid-Noetherian ring\index{t.u. rigid-Noetherian ring@t.u.\ rigid-Noetherian ring} $(I\subseteq A$ is an ideal of definition$)$ such that $U=\Spec A\setminus V(I)$ is affine, then there exists a finitely presented $A$-module $M$ such that $\til{M}^{\rig}\cong\mathscr{F}$ and 
$$
\H^0(\mathscr{X},\mathscr{F})=\varinjlim_{n\geq 0}\Hom_A(I^n,M).
$$
In particular, $\mathscr{F}\mapsto\Gamma(\mathscr{X},\mathscr{F})$ gives the quasi-inverse to the functor 
$$
\Coh_U\stackrel{\cdot^{\rig}}{\longrightarrow}\Coh_{\mathscr{X}}
$$
as in {\rm \ref{ntn-affinoidscohsheaf}}.

{\rm (2)} For $q>0$ we have $\H^q(\mathscr{X},\mathscr{F})=0$.
\end{thm}

\begin{proof}
(2) is already proved in \ref{prop-serrecriterion}.
As in the proof of \ref{prop-serrecriterion} we have $\H^q(\mathscr{X},\mathscr{F})=\H^q(U,\mathscr{H})$, where $\mathscr{H}$ is a coherent sheaf on $U$ such that $\mathscr{F}=\mathscr{H}^{\rig}$ (cf.\ \ref{thm-affinoidscohsheaf}).
If $q=0$, it is isomorphic to $\varinjlim_{n\geq 0}\Hom_A(I^n,M)$ by Deligne's formula\index{Deligne, P.}.
\end{proof}
\index{Stein affinoid|)}\index{affinoid!Stein affinoid@Stein ---|)}

\subsection{Associated schemes}\label{sub-associatedschemes}
\index{scheme!associated scheme to an affinoid@associated --- to an affinoid|(}
\index{affinoid!associated scheme to an affinoid@associated scheme to an ---|(}
\subsubsection{Definition and functoriality}\label{subsub-associatedscheme}
Let $\mathscr{X}=(\Spf A)^{\rig}$ be a universally Noetherian affinoid\index{affinoid!universally Noetherian affinoid@universally Noetherian ---}, where $A$ is a t.u.\ rigid-Noetherian\index{t.u. rigid-Noetherian ring@t.u.\ rigid-Noetherian ring} ring with a finitely generated ideal of definition $I\subseteq A$.
We set
$$
s(\mathscr{X})=\Spec A\setminus V(I), 
$$
which is a Noetherian scheme.
It follows from \ref{prop-lemmorbetweenaffinoid2} that the scheme $s(\mathscr{X})$ does not depend on the choice of $A$.
The Noetherian scheme $s(\mathscr{X})$ thus defined is said to be the {\em associated scheme}\index{scheme!associated scheme to an affinoid@associated --- to an affinoid}\index{affinoid!associated scheme to an affinoid@associated scheme to an ---} to the universally Noetherian affinoid $\mathscr{X}$.

Let $\mathscr{Y}=(\Spf B)^{\rig}$ be another universally Noetherian affinoid\index{affinoid!universally Noetherian affinoid@universally Noetherian ---} (where $B$ is a t.u.\ rigid-Noetherian\index{t.u. rigid-Noetherian ring@t.u.\ rigid-Noetherian ring} ring), and $\varphi\colon\mathscr{Y}\rightarrow\mathscr{X}$ a morphism of rigid spaces.
Then by \ref{prop-lemmorbetweenaffinoid2} we have the canonically induced morphism
$$
s(\varphi)\colon s(\mathscr{Y})\longrightarrow s(\mathscr{X})
$$
of schemes.
Thus the formation $s\colon\mathscr{X}\mapsto s(\mathscr{X})$ defines a functor $s$ from the category of all universally Noetherian affinoid\index{affinoid!universally Noetherian affinoid@universally Noetherian ---} to the category of Noetherian schemes.

\begin{prop}\label{prop-associatedschemesopenimm}
{\rm (1)} Let $\iota\colon\mathscr{Y}\hookrightarrow\mathscr{X}$ be an open immersion between universally Noetherian affinoids\index{affinoid!universally Noetherian affinoid@universally Noetherian ---}.
Then the morphism $s(\iota)\colon s(\mathscr{Y})\rightarrow s(\mathscr{X})$ is flat.

{\rm (2)} Let $\{\mathscr{U}_{\alpha}\}_{\alpha\in L}$ be a finite affinoid covering of a universally Noetherian affinoid\index{affinoid!universally Noetherian affinoid@universally Noetherian ---} $\mathscr{X}$.
Then the induced morphism $\coprod_{\alpha\in L}s(\mathscr{U}_{\alpha})\rightarrow s(\mathscr{X})$ is faithfully flat.
\end{prop}

Notice that the affinoids $\mathscr{U}_{\alpha}$ in (2) are also assumed to be of the form $\mathscr{U}_{\alpha}=(\Spf A_{\alpha})^{\rig}$ where $A_{\alpha}$ is a t.u.\ rigid-Noetherian ring (cf.\ our convention in the end of \S\ref{subsub-affinoidsintro}).
\begin{proof}
Let $\mathscr{Y}=(\Spf B)^{\rig}\hookrightarrow\mathscr{X}=(\Spf A)^{\rig}$ be an open immersion between universally Noetherian affinoids (where $A$ and $B$ are t.u.\ rigid-Noetherian rings).
By \ref{prop-zariskiriemanntoptop} and \ref{prop-blowups4111} we have a diagram of the following form:
$$
\xymatrix{U\,\ar[d]\ar@{^{(}->}[r]&X\ar[d]\\ \Spf B&\Spf A\rlap{,}}
$$
where the vertical arrows are admissible blow-ups, and the horizontal arrow is a quasi-compact open immersion.
Since admissible blow-ups over affine formal schemes are algebraizable and since such an admissible blow-up has an open basis consisting of algebraizable affine open subspaces, to show (1), we only need to show the following:
\begin{itemize}
\item[(i)] for a diagram of the form 
$$
\xymatrix{\Spec R\,\ar@{^{(}->}[r]&X\ar[d]\\ &\Spec A\rlap{,}}
$$
where the horizontal arrow is an open immersion and the vertical map is a blow-up along an admissible ideal of $A$, the map $\Spec\widehat{R}\setminus V(I\widehat{R})\rightarrow\Spec A\setminus V(I)$ is flat (where $I\subseteq A$ is a finitely generated ideal of definition);
\item[(ii)] for a set of finitely many diagrams $\{\Spec R_{\alpha}\hookrightarrow X\rightarrow\Spec A\}_{\alpha\in L}$ of the form as in (i) (with a fixed $X\rightarrow\Spec A$) such that $\{\Spf\widehat{R}_{\alpha}\}_{\alpha\in L}$ covers the formal completion of $X$, the morphism $\coprod_{\alpha\in L}\Spec\widehat{R}_{\alpha}\setminus V(I\widehat{R}_{\alpha})\rightarrow\Spec A\setminus V(I)$ is faithfully flat.
\end{itemize}

In (i), since $R$ is universally rigid-Noetherian, we know that the completion map $R\rightarrow\widehat{R}$ is flat, and hence the assertion is clear.
In (ii), the map in question is flat due to (i).
Hence we only need to show that the map is surjective, and hence only to show that all closed points are in the image (due to going-down theorem; cf.\ \cite[Theorem 9.5]{Matsu}).
For any closed point $x\in\Spec A\setminus V(I)$ one has a valuation ring $V$ and a map $\Spec V\rightarrow \Spec A$ such that the image of the generic point is $x$ and the image of the closed point lies in $V(I)$.
By valuative criterion we have a lift $\Spec V\rightarrow X$, and hence $\Spf \widehat{V}\rightarrow\widehat{X}$.
There exists $\alpha\in L$ such that this map factors though $\Spf\widehat{V}\rightarrow\Spf\widehat{R}_{\alpha}$.
Then we have $\Spec\widehat{V}\rightarrow\Spec\widehat{R}_{\alpha}$ by which the image of the generic point lies outside $V(I\widehat{R}_{\alpha})$ and is mapped to $x$.
This proves (1).
The proof of (2) is similar; one can reduce to the situation as in (ii).
\end{proof}

\subsubsection{The comparison map}\label{subsub-comparisonmapassaff}
Let $\mathscr{X}=(\Spf A)^{\rig}$ be a universally Noetherian affinoid\index{affinoid!universally Noetherian affinoid@universally Noetherian ---}.
The associated scheme $s(\mathscr{X})$ admits a canonical map, called the {\em comparison map}\index{comparison map},
$$
s\colon\ZR{\mathscr{X}}\longrightarrow s(\mathscr{X})\leqno{(\ast)}
$$
constructed as follows.

For $x\in\ZR{\mathscr{X}}$, take a rigid point\index{point!rigid point@rigid ---}\index{rigid point} $\alpha\colon\Spf V\rightarrow\ZR{\mathscr{X}}$, where $V$ is an $a$-adically complete valuation ring with $a\in\m_V\setminus\{0\}$, that maps the closed point to $x$ (\ref{dfn-ZRpoints2} (1)).
Consider the composition $\sp_X\circ\alpha\colon\Spf V\rightarrow X=\Spf A$, which is an adic morphism of affine formal schemes.
The last morphism induces\ the morphism $\Spec V\rightarrow\Spec A$ of affine schemes and hence the morphism
$$
{\textstyle \Spec V[\frac{1}{a}]}\longrightarrow\Spec A\setminus V(I)=s(\mathscr{X}).\leqno{(\ast\ast)}
$$
Since $V[\frac{1}{a}]$ is a field ({\bf \ref{ch-pre}}.\ref{prop-sep}), the map $(\ast)$ gives rise to a point of $s(\mathscr{X})$, which we denote by $s(x)$, and thus we obtain the desired map $(\ast)$ set-theoretically.
Notice that in the above construction the point $s(x)$ does not depend on the choice of the rigid point $\alpha$ due to {\bf \ref{ch-pre}}.\ref{prop-maxspe4} and {\bf \ref{ch-pre}}.\ref{prop-maxspe}.

\begin{rem}\label{rem-comparisonmapfactorsep}{\rm 
Notice that by the above construction we easily deduce the following: for $x,x'\in\ZR{\mathscr{X}}$ such that $x'\in G_x$ (that is, $x'$ is a generization\index{generization} of $x$) we have $s(x)=s(x')$.
Hence the map $s$ factors by the separation map (\S\ref{subsub-separation})
$$
\xymatrix{\ZR{\mathscr{X}}\ar[r]_{\sep_{\mathscr{X}}}\ar@/^1pc/[rr]^s&[\mathscr{X}]\ar[r]_(.45)s&s(\mathscr{X}),}
$$
where, by a slight abuse of notation, we denote the resulting continuous map $[\mathscr{X}]\rightarrow s(\mathscr{X})$ by the same notation $s$.}
\end{rem}

In case $\mathscr{X}$ is a Stein affinoid (\ref{dfn-cohomologyrigidsp1str}), the map $(\ast)$ admits another description as follows.
In this case, by \ref{thm-comparisonaffinoid} we have 
$$
s(\mathscr{X})=\Spec\Gamma(\mathscr{X},\O_{\mathscr{X}}),
$$
which is a Noetherian affine scheme.
For $x\in\ZR{\mathscr{X}}$ the associated rigid point $\alpha\colon\Spf V\rightarrow\ZR{\mathscr{X}}$ induces the map $\Gamma(\mathscr{X},\O_{\mathscr{X}})\rightarrow K_x=V[\frac{1}{a}]$, which is nothing but the one induced from the restriction map $\Gamma(\mathscr{X},\O_{\mathscr{X}})\rightarrow B_x=\O_{\mathscr{X},x}$.
Hence the point $s(x)$ coincides with the prime ideal that is the pull-back of the maximal ideal of $B_x$ by the last map.

\begin{prop}\label{prop-comparisonmapassaff}
The set-theoretic map $(\ast)$ extends canonically to a flat morphism of locally ringed spaces $($denoted by the same symbol$)$
$$
s\colon(\ZR{\mathscr{X}},\O_{\mathscr{X}})\longrightarrow(s(\mathscr{X}),\O_{s(\mathscr{X})}),
$$
where $\O_{s(\mathscr{X})}$ is the structure sheaf of the scheme $s(\mathscr{X})$.
\end{prop}

\begin{proof}
First notice that the issue is local on $\mathscr{X}$; more precisely, if $\mathscr{X}=\bigcup_{\alpha\in L}\mathscr{U}_{\alpha}$ is an affinoid covering, then it suffices to show that the morphism $\ZR{\mathscr{U}_{\alpha}}\rightarrow s(\mathscr{X})$ extends canonically to a flat morphism of locally ringed spaces for each $\alpha\in L$.
Since this morphism factors through $s(\mathscr{U}_{\alpha})$ and the morphism $s(\mathscr{U}_{\alpha})\rightarrow s(\mathscr{X})$ of schemes is flat (\ref{prop-associatedschemesopenimm} (1)), we only have to show that each $\ZR{\mathscr{U}_{\alpha}}\rightarrow s(\mathscr{U}_{\alpha})$ extends canonically to a flat morphism of locally ringed spaces.
Thus we may assume that $\mathscr{X}$ is a Stein affinoid.
But in this case, the morphism $s$ is nothing but the morphism corresponding to the identity map of $\Gamma(\mathscr{X},\O_{\mathscr{X}})$ by the correspondence in \cite[(1.6.3)]{EGAInew}, and hence the assertion is clear.
\end{proof}

The following corollary is easy to see by \ref{prop-descriptionlocalrings} (1); the last assertion follows from \ref{prop-associatedschemesopenimm} (1):
\begin{cor}\label{cor-comparisonmapassaffloc}
Let $\mathscr{X}=(\Spf A)^{\rig}$ be a universally Noetherian affinoid, $x\in\ZR{\mathscr{X}}$ a point, and $\{\mathscr{U}_{\alpha}=(\Spf A_{\alpha})^{\rig}\}_{\alpha\in L}$ a cofinal system of formal neighborhoods of $x$.
For each $\alpha\in L$, let $x_{\alpha}$ be the image of $x$ by the map $s\colon\ZR{\mathscr{U}_{\alpha}}\rightarrow s(\mathscr{U}_{\alpha})$.
Then we have
$$
\O_{\mathscr{X},x}=\varinjlim_{\alpha\in L}\O_{s(\mathscr{U}_{\alpha}),x_{\alpha}}.
$$
Moreover, for $\alpha\leq\beta$ the transition map $\O_{s(\mathscr{U}_{\alpha}),x_{\alpha}}\rightarrow\O_{s(\mathscr{U}_{\beta}),x_{\beta}}$ is flat. \hfill$\square$
\end{cor}

\begin{prop}\label{prop-comparisonmapassaffprop1}
Let $\mathscr{X}$ be a universally Noetherian affinoid\index{affinoid!universally Noetherian affinoid@universally Noetherian ---}.
Then the correspondence $\mathscr{F}\mapsto s^{\ast}\mathscr{F}$ establishes the categorical equivalence 
$$
\Coh_{s(\mathscr{X})}\stackrel{\sim}{\longrightarrow}\Coh_{\mathscr{X}}.
$$
Moreover, for any coherent sheaf $\mathscr{F}$ on $s(\mathscr{X})$ we have a canonical isomorphism between the cohomologies
$$
\H^q(s(\mathscr{X}),\mathscr{F})\stackrel{\sim}{\longrightarrow}\H^q(\mathscr{X},s^{\ast}\mathscr{F})
$$
for each $q\geq 0$.
\end{prop}

\begin{proof}
Since $s^{\ast}\mathscr{F}$ is nothing but $\mathscr{F}^{\rig}$ in the sense as in \ref{ntn-affinoidscohsheaf}, the assertions are rehashes of \ref{thm-affinoidscohsheaf} and \ref{thm-comparisonaffinoid}.
\end{proof}
\index{affinoid!associated scheme to an affinoid@associated scheme to an ---|)}
\index{scheme!associated scheme to an affinoid@associated --- to an affinoid|)}
\index{affinoid|)}

\addcontentsline{toc}{subsection}{Exercises}
\subsection*{Exercises}
\begin{exer}\label{exer-weierstrasslaurentrational}{\rm 
Show that the intersection of two Weierstrass (resp.\ Laurent, resp.\ rational) subdomains of an affinoid is again a Weierstrass (resp.\ Laurent, resp.\ rational) subdomain.}
\end{exer}

\begin{exer}\label{exer-cormorbetweenaffinoid1}{\rm 
Let $X$ be a locally universally rigid-Noetherian formal scheme\index{formal scheme!universally rigid-Noetherian formal scheme@universally rigid-Noetherian ---!locally universally rigid-Noetherian formal scheme@locally --- ---} {\rm ({\bf \ref{ch-formal}}.\ref{dfn-formalsch})}, and $\pi\colon X'\rightarrow X$ an admissible blow-up.
Then show that $\pi_{\ast}\O_{X'}$ is an adically quasi-coherent $\O_X$-algebra.}
\end{exer}

\begin{exer}\label{exer-associatedschememapsurjective}{\rm 
Let $\mathscr{X}$ be a universally Noetherian affinoid, and consider the comparison map $s\colon\ZR{\mathscr{X}}\rightarrow s(\mathscr{X})$.
Show that the image of $s$ contains all closed points of $s(\mathscr{X})$.}
\end{exer}



\section{Basic properties of morphisms of rigid spaces}\label{sub-basicmorproprigid}
In this section we discuss several properties of morphisms between rigid spaces. 
We give the definitions of those properties and establish some of the basic results, such as base change stability, interrelation with other properties, etc. 
In discussing separated morphisms and proper morphisms, we also give the valuative criterion. 
Notice that our properness in rigid geometry presented here is the so-called {\em Raynaud\index{Raynaud, M.} properness}, that is, the one defined by means of the properness of the formal models.
This differs, a priori, from the one introduced by Kiehl in \cite{Kieh1}, the so-called {\em Kiehl\index{Kiehl, R} properness}. 
It will be our objective in one of the later chapters to show that, at least in the adhesive situation, these notions of properness are actually equivalent to each other.
In \S\ref{subsub-finitudesrigid} we present the finiteness theorem for proper maps between universally adhesive rigid spaces.

\subsection{Quasi-compact and quasi-separated morphisms}\label{sub-qcompmorrigid}
\index{morphism of rigid spaces@morphism (of rigid spaces)!quasi-compact morphism of rigid spaces@quasi-compact ---|(}
\index{morphism of rigid spaces@morphism (of rigid spaces)!quasi-separated morphism of rigid spaces@quasi-separated ---|(}
\begin{prop}\label{prop-qcptmorphismrigid}
Let $\varphi\colon\mathscr{X}\rightarrow\mathscr{Y}$ be a morphism of rigid spaces.
Then the following conditions are equivalent$:$
\begin{itemize}
\item[{\rm (a)}] there exists a covering $\{\mathscr{V}_{\alpha}\rightarrow\mathscr{Y}\}_{\alpha\in L}$ in the small site $\mathscr{Y}_{\ad}$ with each $\mathscr{V}_{\alpha}$ a coherent rigid spaces such that $\mathscr{V}_{\alpha}\times_{\mathscr{Y}}\mathscr{X}$ is a quasi-compact rigid space\index{rigid space!quasi-compact rigid space@quasi-compact ---} {\rm (\ref{dfn-generalrigidspace2} (1))} for each $\alpha\in L;$
\item[{\rm (b)}] for any morphism $\mathscr{S}\rightarrow\mathscr{Y}$ of rigid spaces with $\mathscr{S}$ quasi-compact, $\mathscr{S}\times_{\mathscr{Y}}\mathscr{X}$ is quasi-compact$;$
\item[{\rm (c)}] the induced map $\ZR{\varphi}\colon\ZR{\mathscr{X}}\rightarrow\ZR{\mathscr{Y}}$ is quasi-compact\index{map@map (continuous)!quasi-compact map@quasi-compact ---}\index{quasi-compact!quasi-compact map@--- map} as a map of topological spaces $($cf.\ {\rm {\bf \ref{ch-pre}}.\ref{dfn-quasicompactness} (2))}$;$
\item[{\rm (d)}] the morphism of small admissible topoi $\varphi^{\sim}_{\ad}\colon\mathscr{X}^{\sim}_{\ad}\rightarrow\mathscr{Y}^{\sim}_{\ad}$ is quasi-compact $($cf.\ {\rm {\bf \ref{ch-pre}}.\ref{dfn-coherenttopos3})}.
\end{itemize}
\end{prop}

\begin{proof}
As (b) $\Rightarrow$ (a) is trivial, let us first show (a) $\Rightarrow$ (b).
One easily reduces to the situation where $\mathscr{Y}$ and $\mathscr{S}$ are coherent.
In this case, $\mathscr{X}$ is quasi-compact and hence is covered by finitely many coherent open rigid subspaces\index{rigid subspace!open rigid subspace@open ---} $\mathscr{U}_{\alpha}$.
Then $\mathscr{S}\times_{\mathscr{Y}}\mathscr{X}$ is covered by finitely many coherent open rigid spaces $\mathscr{U}_{\alpha}\times_{\mathscr{Y}}\mathscr{S}$ and hence is quasi-compact.
To show the equivalence of (a) and (c), again one can reduce to the case where $\mathscr{Y}$ is coherent.
Then the issue is to show that $\mathscr{X}$ is quasi-compact if and only if $\ZR{\mathscr{X}}$ is quasi-compact.
But this has been done in \ref{prop-thmgeneralrigidspace31-2}.
Finally, the equivalence of (c) and (d) follows from \ref{thm-generalrigidspace31-31}.
\end{proof}

\begin{prop}\label{prop-qcptmorphismrigidqs}
Let $\varphi\colon\mathscr{X}\rightarrow\mathscr{Y}$ be a morphism of rigid spaces.
Then the following conditions are equivalent$:$
\begin{itemize}
\item[{\rm (a)}] for any morphism $\mathscr{S}\rightarrow\mathscr{Y}$ from a quasi-separated rigid space\index{rigid space!quasi-separated rigid space@quasi-separated ---} {\rm (\ref{dfn-generalrigidspace2} (2))} $\mathscr{S}$, $\mathscr{S}\times_{\mathscr{Y}}\mathscr{X}$ is quasi-separated$;$
\item[{\rm (b)}] the induced map $\ZR{\varphi}\colon\ZR{\mathscr{X}}\rightarrow\ZR{\mathscr{Y}}$ is quasi-separated as a map of topological spaces$;$
\item[{\rm (c)}] the morphism of small admissible topoi $\varphi^{\sim}_{\ad}\colon\mathscr{X}^{\sim}_{\ad}\rightarrow\mathscr{Y}^{\sim}_{\ad}$ is quasi-separated $($cf.\ {\rm {\bf \ref{ch-pre}}.\ref{dfn-coherenttopos3})}$;$
\item[{\rm (d)}] the diagonal morphism $\mathscr{X}\rightarrow\mathscr{X}\times_{\mathscr{Y}}\mathscr{X}$ is quasi-compact.
\end{itemize}
\end{prop}

\begin{proof}
By \ref{thm-generalrigidspace31-31} we have the equivalence of (b) and (c).
The condition (c) is equivalent to (a) where the morphisms $\mathscr{S}\hookrightarrow\mathscr{Y}$ are assumed to be an open immersion with $\mathscr{S}$ coherent.
By \ref{prop-qcptmorphismrigid} (b), (c) is still equivalent to (a) without this assumption.

Suppose that (d) holds.
To deduce (c), it is enough to show that for any open immersion $\mathscr{S}\hookrightarrow\mathscr{Y}$ from a quasi-separated rigid space, the base change $\mathscr{X}_{\mathscr{S}}$ is quasi-separated.
To this end, we may assume that $\mathscr{Y}$ is quasi-separated and that $\mathscr{Y}=\mathscr{S}$.
For any quasi-compact open subspaces $\mathscr{U},\mathscr{V}\subseteq\mathscr{X}$, the intersection $\mathscr{U}\times_{\mathscr{X}}\mathscr{V}$ coincides with the pull-back of $\mathscr{U}\times_{\mathscr{Y}}\mathscr{V}$ by the diagonal and hence is quasi-compact.
Hence $\mathscr{X}$ is quasi-separated, thereby (c).

Conversely, since $\mathscr{Y}_{\ad}$ has a generating family consisting of coherent open subspaces, we deduce (c) $\Rightarrow$ (d), using \ref{prop-qcptmorphismrigid}.
\end{proof}

\begin{dfn}\label{dfn-qcptmorphismrigid1}{\rm 
A morphism $\varphi\colon\mathscr{X}\rightarrow\mathscr{Y}$ of rigid space is said to be {\em quasi-compact} (resp.\ {\em quasi-separated}) if it satisfies the conditions in \ref{prop-qcptmorphismrigid} (resp.\ \ref{prop-qcptmorphismrigidqs}).
If $\varphi$ is quasi-compact and quasi-separated, it is said to be {\em coherent}\index{morphism of rigid spaces@morphism (of rigid spaces)!coherent morphism of rigid spaces@coherent ---}.}
\end{dfn}

\begin{rem}\label{rem-consistencyopenimm3x}{\rm 
Notice that an open immersion $\iota\colon\mathscr{U}\hookrightarrow\mathscr{X}$ between coherent rigid spaces is always coherent in the above sense; cf.\ \ref{prop-consistencytermopenimm2}. This again explains the consistency of our terminology `coherent open immersion'\index{immersion!open immersion of rigid spaces@open --- (of rigid spaces)!coherent open immersion of rigid spaces@coherent --- ---} (defined in \ref{dfn-cohrigidspaceopenimm1}).}
\end{rem}

The following proposition is easy to see (cf.\ {\bf \ref{ch-pre}}.\ref{prop-basechangestable}):
\begin{prop}\label{prop-qcptmorphismrigid2}
{\rm (1)} A locally of finite type morphism\index{morphism of rigid spaces@morphism (of rigid spaces)!morphism of rigid spaces locally of finite type@--- locally of finite type} is of finite type if and only if it is quasi-compact.

{\rm (2)} The composition of two quasi-compact $($resp.\ quasi-separated, resp.\ coherent$)$ morphisms is quasi-compact $($resp.\ quasi-separated, resp.\ coherent$)$.

{\rm (3)} If $\varphi\colon\mathscr{X}\rightarrow\mathscr{X}'$ and $\psi\colon\mathscr{Y}\rightarrow\mathscr{Y}'$ are two quasi-compact $($resp.\ quasi-separated, resp.\ coherent$)$ morphisms over a rigid space $\mathscr{S}$, then the induced morphism $\varphi\times_{\mathscr{S}}\psi\colon\mathscr{X}\times_{\mathscr{S}}\mathscr{Y}\rightarrow\mathscr{X}'\times_{\mathscr{S}}\mathscr{Y}'$ is quasi-compact $($resp.\ quasi-separated, resp.\ coherent$)$.

{\rm (4)} If $\varphi\colon\mathscr{X}\rightarrow\mathscr{Y}$ is a quasi-compact $($resp.\ quasi-separated, resp.\ coherent$)$ morphism over a rigid space $\mathscr{S}$ and $\mathscr{S}'\rightarrow\mathscr{S}$ is a morphism of rigid spaces, then the induced morphism $\varphi_{\mathscr{S}'}\colon\mathscr{X}\times_{\mathscr{S}}\mathscr{S}'\rightarrow\mathscr{Y}\times_{\mathscr{S}}\mathscr{S}'$ is quasi-compact $($resp.\ quasi-separated, resp.\ coherent$)$.
\end{prop}
\index{morphism of rigid spaces@morphism (of rigid spaces)!quasi-separated morphism of rigid spaces@quasi-separated ---|)}
\index{morphism of rigid spaces@morphism (of rigid spaces)!quasi-compact morphism of rigid spaces@quasi-compact ---|)}

\subsection{Finite morphism}\label{finitemorrigid}
\index{morphism of rigid spaces@morphism (of rigid spaces)!finite morphism of rigid spaces@finite ---|(}
\begin{dfn}\label{dfn-finitemorcoherentrigid}{\rm 
(1) A morphism $\varphi\colon\mathscr{X}\rightarrow\mathscr{Y}$ of coherent rigid spaces is said to be {\em finite} if it has a finite\index{morphism of formal schemes@morphism (of formal schemes)!finite morphism of formal schemes@finite ---} {\rm ({\bf \ref{ch-formal}}.\ref{dfn-finitemorform1})} formal model $f\colon X\rightarrow Y$.

(2) A morphism $\varphi\colon\mathscr{X}\rightarrow\mathscr{Y}$ of rigid spaces is said to be {\em finite} if it is coherent\index{morphism of rigid spaces@morphism (of rigid spaces)!coherent morphism of rigid spaces@coherent ---} (\ref{dfn-qcptmorphismrigid1}) and for any open immersion $\mathscr{V}\hookrightarrow\mathscr{Y}$ from a coherent rigid space, the induced morphism $\mathscr{X}\times_{\mathscr{Y}}\mathscr{V}\rightarrow\mathscr{V}$ between coherent rigid spaces is finite in the sense of (1).}
\end{dfn}

It follows from {\bf \ref{ch-formal}}.\ref{prop-finitemorform2} (4) that for a morphism $\varphi\colon\mathscr{X}\rightarrow\mathscr{Y}$ of coherent rigid spaces the definition (2) of the finiteness is consistent with (1).
\begin{prop}\label{prop-finitemorrigid1}
Let $\varphi\colon\mathscr{X}\rightarrow\mathscr{Y}$ be a morphism of coherent universally Noetherian rigid spaces.
Then $\varphi$ is finite if and only if there exists a finite and distinguished formal model $f\colon X\rightarrow Y$ of $\varphi$.
\end{prop}

\begin{proof}
The `if' part is trivial.
Let $f\colon X\rightarrow Y$ be a finite formal model of $\varphi$, and $Y'\rightarrow Y$ an admissible blow-up such that $Y'$ is $\mathscr{I}$-torsion free, where $\mathscr{I}$ is an ideal of definition of $Y$.
Let $X'\hookrightarrow X$ be the strict transform (\ref{dfn-stricttransform2}).
Then by {\bf \ref{ch-formal}}.\ref{prop-finitemorform2} (2) (4) and {\bf \ref{ch-formal}}.\ref{prop-closedimmformal3} the morphism $X'\rightarrow Y'$ is finite.
\end{proof}

\begin{prop}\label{prop-finitemorrigid10}
Let $\mathscr{Y}$ be a coherent universally Noetherian rigid space, and $\mathscr{A}$ an $\O_{\mathscr{Y}}$-algebra that is coherent as an $\O_{\mathscr{Y}}$-module.
Then there exists uniquely up to $\mathscr{Y}$-isomorphisms a finite morphism $\varphi\colon\mathscr{X}\rightarrow\mathscr{Y}$ such that $\ZR{\varphi}_{\ast}\O_{\mathscr{X}}\cong\mathscr{A}$.
\end{prop}

\begin{proof}
Let $Y$ be a universally rigid-Noetherian formal model of $\mathscr{Y}$.
Then by \ref{thm-tateacyclic20} there exists an $\O_Y$-algebra $\mathscr{A}_Y$ that is finitely presented as an $\O_Y$-module.
Let $X=\Spf\mathscr{A}_Y$ ({\bf \ref{ch-formal}}.\ref{dfn-affinespectrumformalrel}).
Then the structural map $f\colon X\rightarrow Y$ is finite ({\bf \ref{ch-formal}}.\ref{prop-finitemorphismlocalconstruction}).
Let $\varphi=f^{\rig}\colon\mathscr{X}=X^{\rig}\rightarrow\mathscr{Y}$ be the associated morphism of coherent rigid spaces.
To show $\ZR{\varphi}_{\ast}\O_{\mathscr{X}}\cong\mathscr{A}$, we may assume that $Y$ (and hence $X$) is affine, and it suffices to show that $\Gamma(\mathscr{X},\O_{\mathscr{X}})=\Gamma(\mathscr{Y},\mathscr{A})$.
But this follows from \ref{thm-comparisonaffinoid}.

To show the uniqueness, suppose that there is another such finite morphism $\varphi'\colon\mathscr{X}'\rightarrow\mathscr{Y}$.
Let $X\rightarrow Y$ and $X'\rightarrow Y$ be universaly rigid-Noetherian formal models of $\mathscr{X}\rightarrow\mathscr{Y}$ and $\mathscr{X}'\rightarrow\mathscr{Y}$, respectively.
Replacing $Y$ by an affine open subset, we may assume that $\mathscr{Y}$ is affinoid of the form $(\Spf B)^{\rig}$, where $B$ is a t.u.\ rigid-Noetherian ring, and that $\mathscr{X}$ and $\mathscr{X}'$ are affinoids of the similar form.
Let $X=\Spf A$ (resp.\ $X'=\Spf A'$) be an affine formal model of $\mathscr{X}$ (resp.\ $\mathscr{X}'$) such that $A$ (resp.\ $A'$) is t.u.\ rigid-Noetherian ring.
By \ref{thm-comparisonaffinoid} $\Spec A'$ and $\Spec A$ are isomorphic outside the closed loci defined by an ideal of definition.
By \cite[Premi\`ere partie, (5.7.12)]{RG} (recorded below in \ref{prop-birationalgeom02}) we have admissible blow-ups $X''\rightarrow X$ and $X''\rightarrow X'$ and thus $\mathscr{X}\cong\mathscr{X}'$, as desired.
\end{proof}

In the situation as in \ref{prop-finitemorrigid10} the morphism $\varphi$ is called the {\em finite morphism associated to $\mathscr{A}$}.
By the uniqueness, one can also define the finite morphism associated to a coherent $\O_{\mathscr{X}}$-algebra on a general rigid space.

\begin{prop}\label{prop-finitemorrigid8}
Let $\varphi\colon\mathscr{X}\rightarrow\mathscr{Y}$ be a morphism of locally universally Noetherian rigid spaces.
Then $\varphi$ is finite if and only if it is the finite morphism associated to an $\O_{\mathscr{Y}}$-algebra that is coherent as an $\O_{\mathscr{X}}$-module.
\end{prop}

For the proof we need:
\begin{lem}\label{lem-finitemorrigid8}
Let $\varphi\colon\mathscr{X}\rightarrow\mathscr{Y}$ be a finite morphism of locally universally Noetherian rigid spaces, and $\mathscr{F}$ a coherent $\O_{\mathscr{X}}$-module.
Then $\ZR{\varphi}_{\ast}\mathscr{F}$ is a coherent $\O_{\mathscr{Y}}$-module.
\end{lem}

\begin{proof}
Take a finite formal model $f\colon X\rightarrow Y$ of $\varphi$.
Since finite morphisms are stable under base change ({\bf \ref{ch-formal}}.\ref{prop-finitemorform2} (4)), we may assume that $Y$ has an invertible ideal of definition $\mathscr{I}$ (replacing $Y$ by an admissible blow-up along an ideal of definition of finite type).
To show that $\ZR{\varphi}_{\ast}\mathscr{F}$ is coherent, we may work locally and hence may assume that $Y$ is affine; $Y=\Spf B$ with the invertible ideal of definition $I=(a)$.
Accordingly ({\bf \ref{ch-formal}}.\ref{prop-finitemorform2} (1)), $X$ is also affine; $X=\Spf A$.
By \ref{thm-affinoidscohsheaf} the sheaf $\mathscr{F}$ corresponds to a coherent $A[\frac{1}{a}]$-module $M$.
By \ref{thm-comparisonaffinoid} we see
$$
\Gamma(\ZR{\mathscr{X}},\O_{\mathscr{X}})=A[{\textstyle \frac{1}{a}}],\quad\Gamma(\ZR{\mathscr{X}},\mathscr{F})=M
$$
which are finite over the Noetherian ring $\Gamma(\ZR{\mathscr{Y}},\O_{\mathscr{Y}})=B[\frac{1}{a}]$, since $A$ is finite over $B$ ({\bf \ref{ch-formal}}.\ref{prop-finitemorform1}).
Since this holds for any sufficiently small affinoid open sets of $\mathscr{Y}$, we deduce that $\ZR{\varphi}_{\ast}\mathscr{F}$ is a coherent $\O_{\mathscr{Y}}$-module.
\end{proof}

\begin{proof}[Proof of Proposition {\rm \ref{prop-finitemorrigid8}}]
The `if' part is clear.
To show the converse, first notice that by \ref{lem-finitemorrigid8} the sheaf $\ZR{\varphi}_{\ast}\O_{\mathscr{X}}=\mathscr{A}$ is a coherent $\O_{\mathscr{Y}}$-module.
Then $\varphi$ is isomorphic to the one associated to $\mathscr{A}$, as one can verify by an argument similar to that in the proof of (the uniqueness of) \ref{prop-finitemorrigid10}.
\end{proof}

 \begin{prop}\label{prop-finitemorrigid2}
{\rm (1)} The composition of two finite morphisms between locally universally Noetherian rigid spaces is finite.

{\rm (2)} If $\varphi\colon\mathscr{X}\rightarrow\mathscr{X}'$ and $\psi\colon\mathscr{Y}\rightarrow\mathscr{Y}'$ are two finite morphisms over a rigid space $\mathscr{S}$, then the induced morphism $\varphi\times_{\mathscr{S}}\psi\colon\mathscr{X}\times_{\mathscr{S}}\mathscr{Y}\rightarrow\mathscr{X}'\times_{\mathscr{S}}\mathscr{Y}'$ is finite.

{\rm (3)} If $\varphi\colon\mathscr{X}\rightarrow\mathscr{Y}$ is a finite morphism over a rigid space $\mathscr{S}$ and $\mathscr{S}'\rightarrow\mathscr{S}$ is a morphism of rigid spaces, then the induced morphism $\varphi_{\mathscr{S}'}\colon\mathscr{X}\times_{\mathscr{S}}\mathscr{S}'\rightarrow\mathscr{Y}\times_{\mathscr{S}}\mathscr{S}'$ is finite.
\end{prop}

\begin{proof}
The statements (2) and (3) follow from {\bf \ref{ch-formal}}.\ref{prop-finitemorform2} (3), (4).
(Note that under the property (1) the properties (2) and (3) are equivalent due to {\bf \ref{ch-pre}}.\ref{prop-basechangestable}.)
To show (1), let $\varphi\colon\mathscr{X}\rightarrow\mathscr{Y}$ and $\psi\colon\mathscr{Y}\rightarrow\mathscr{Z}$ be finite morphisms between coherent rigid spaces.
One can take a diagram
$$
X\stackrel{f}{\longrightarrow}Y'\stackrel{\pi}{\longrightarrow}Y\stackrel{g}{\longrightarrow}Z
$$
consisting of coherent universally rigid-Noetherian formal schemes such that:
\begin{itemize}
\item $f$ (resp.\ $g$) is a finite formal model of $\varphi$ (resp.\ $\psi$);
\item $\pi$ is an admissible blow-up.
\end{itemize}
Since we may work locally on $\mathscr{Z}$, we may further assume that $Z$ is affine $Z=\Spf A$.
Then, by GFGA existence theorem ({\bf \ref{ch-formal}}.\ref{thm-GFGAexaweakcoherentspecial}), one can algebraize the diagram into the diagram of schemes 
$$
\til{X}\stackrel{\til{f}}{\longrightarrow}\til{Y}'\stackrel{\til{\pi}}{\longrightarrow}\til{Y}\stackrel{\til{g}}{\longrightarrow}\til{Z}=\Spec A,
$$
where $\til{f}$ and $\til{g}$ are finite morphisms and $\til{\pi}$ is a blow-up.
By algebraic flattening theorem (\cite{RG}) there exists a $U$-admissible blow-up ($U=\Spec A\setminus V(I)$, where $I\subseteq A$ is an ideal of definition) $\til{Z}'\rightarrow Z$ such that the strict transform $\til{X}'$ (resp.\ $\til{Y}''$) of $\til{X}$ (resp.\ $\til{Y}'$) is finite over $\til{Z}'$.
By passage to the $I$-adic completions of the resulting strict transforms, we obtain the diagram
$$
X'\stackrel{f'}{\longrightarrow}Y''\stackrel{g'}{\longrightarrow}Z'
$$
consisting of finite morphisms such that $\varphi=(f')^{\rig}$ and $\psi=(g')^{\rig}$.
Since $g'\circ f'$ is finite, $\psi\circ\varphi$ is finite, as desired.
\end{proof}

\begin{prop}\label{prop-finitemorrigid99}
Let $\varphi\colon\mathscr{X}\rightarrow\mathscr{Y}$ be a finite morphism of rigid spaces.
Then, for any point $y\in\ZR{\mathscr{Y}}$, the fiber $\ZR{\varphi}^{-1}(y)$ is a finite set.
\end{prop}

\begin{proof}
Let $\alpha_y\colon\Spf\widehat{V_y}\rightarrow\ZR{\mathscr{Y}}$ be the associated rigid point of $y$ (\ref{dfn-ZRpoints32}).
We need to show that $\ZR{\mathscr{X}\times_{\mathscr{Y}}(\Spf\widehat{V}_y)^{\rig}}$ is a finite set.
Hence, replacing $\mathscr{Y}$ by $(\Spf\widehat{V}_y)^{\rig}$, we may assume that $\mathscr{Y}$ is of the form $(\Spf V)^{\rig}$, where $V$ is an $a$-adically complete valuation ring.
Then $\mathscr{X}=(\Spf A)^{\rig}$, where $A$ is finite flat over $V$.
Since any morphism from $A$ to a valuation ring $V'$ is induced from a unique morphism from $A_{\red}$ to $V'$, we may further assume that $A=A_{\red}$; note that $A_{\red}$ for a topologically of finite type $V$-algebra $A$ is $a$-adically complete.
Since $A$ is finite over $V$, we have
$$
A\otimes_VK\cong L_1\times\cdots\times L_n,
$$
where $L_i/K$ is a finite extension field for $i=1,\ldots,n$.
Now, in view of \ref{prop-ZRpoints3}, points of $\ZR{\mathscr{X}}$ are in one to one correspondence with valuation subrings of $A\otimes_VK$ containing $A$ and dominating $V$.
To see there are only finitely many such objects, it suffices to invoke the classical fact that the number of valuation subrings of the finite extension $L_i$ ($i=1,\ldots,n$) of $K$ dominating $V$ is finite (e.g., \cite{Bourb1}, Chap.\ VI, \S8.3, Theorem 1).
\end{proof}
\index{morphism of rigid spaces@morphism (of rigid spaces)!finite morphism of rigid spaces@finite ---|)}

\subsection{Closed immersions}\label{sub-closedimmrigid}
\subsubsection{Definition and first properties}\label{subsub-closedimmrigid}
\index{immersion!closed immersion of rigid spaces@closed --- (of rigid spaces)|(}
\begin{prop}\label{prop-closedimmrigid1}
Let $\iota\colon\mathscr{Y}\rightarrow\mathscr{X}$ be a morphism of coherent universally Noetherian rigid spaces.
Then the following conditions are equivalent$:$
\begin{itemize}
\item[{\rm (a)}] there exists a formal model $i\colon Y\rightarrow X$ of $\varphi$ that is a closed immersion between coherent universally rigid-Noetherian formal schemes$;$
\item[{\rm (b)}] there exists a formal model $i\colon Y\rightarrow X$ of $\varphi$ that is a closed immersion of finite presentation between coherent universally rigid-Noetherian formal schemes$;$
\item[{\rm (c)}] there exists a cofinal family of formal models $\{i_{\lambda}\colon Y_{\lambda}\rightarrow X_{\lambda}\}$ of $\iota$ consisting of closed immersions between coherent universally rigid-Noetherian formal schemes.
\end{itemize}
\end{prop}

As a preparation for the proof, we first prove:
\begin{lem}\label{lem-closedimmrigid1fpmodel}
Let $i\colon Y\hookrightarrow X$ be a closed immersion of coherent universally rigid-Noetherian formal schemes.
Then $i$ factorizes as
$$
\xymatrix{Y\,\ar@{^{(}->}[r]^{i'}&Y'\,\ar@{^{(}->}[r]^{i''}&X,}
$$
where $i''$ is a closed immersion of finite presentation and $i'$ is a closed immersion defined by a bounded $\mathscr{I}_{Y'}$-torsion adically quasi-coherent ideal, where $\mathscr{I}_{Y'}$ is an ideal of definition of $Y'$ $($cf.\ {\rm {\bf \ref{ch-formal}}.\ref{prop-closedimmformal6}}$)$.
\end{lem}

\begin{proof}
Consider the surjective morphism $\O_X\rightarrow f_{\ast}\O_Y$ ({\bf \ref{ch-formal}}.\ref{cor-closedimmformal43}).
By {\bf \ref{ch-formal}}.\ref{prop-closedimmformal6} the kernel $\mathscr{K}$ of this morphism is an adically quasi-coherent ideal of $\O_X$.
Then using Exercise \ref{exer-extopenadiccoh1approx}, one finds an adically quasi-coherent subideal $\mathscr{K}'\subseteq\mathscr{K}$ of finite type such that $\mathscr{K}/\mathscr{K}'$ is bounded $\mathscr{I}_X$-torsion, where $\mathscr{I}_X$ is an ideal of definition of finite type on $X$.
Let $i''\colon Y'\hookrightarrow X$ be the closed immersion of finite presentation with the defining ideal $\mathscr{K}'$ ({\bf \ref{ch-formal}}.\ref{prop-closedimmformal6}).
Then we have the closed immersion $i'\colon Y\hookrightarrow Y'$ defined by $\mathscr{K}/\mathscr{K}'$.
\end{proof}

\begin{proof}[Proof of Proposition {\rm \ref{prop-closedimmrigid1}}]
First let us show (a) $\Rightarrow$ (b).
For a formal model $i\colon Y\rightarrow X$ of $\iota$ that is a closed immersion of coherent universally rigid-Noetherian formal schemes, we obtain a factorization as in \ref{lem-closedimmrigid1fpmodel}.
We need to show that the closed immersion $i''$ gives another formal model of $\varphi$.
To this end, we perform the admissible blow-up $Y''\rightarrow Y'$ of $Y'$ along an ideal of definition of finite type.
Such an admissible blow-up comes as the strict transform of the admissible blow-up of $X$ along the blow-up center of $Y''\rightarrow Y'$.
The strict transform of $Y$ is then isomorphic to $Y'$, thereby the claim.

The implication (b) $\Rightarrow$ (c) follows from \ref{prop-blowups15x}, \ref{lem-closedimmrigid1fpmodel}, and the fact that the strict transform of a closed immersion by an admissible blow-up is again a closed immersion.
The implication (c) $\Rightarrow$ (a) is trivial
\end{proof}

\begin{rem}\label{rem-closedimmrigid1coherent}{\rm 
Notice that, if the $X$ in \ref{prop-closedimmrigid1} (a) is coherent universally adhesive\index{formal scheme!universally adhesive formal scheme@universally adhesive ---}\index{adhesive!universally adhesive@universally ---!universally adhesive formal scheme@--- --- formal scheme} ({\bf \ref{ch-formal}}.\ref{dfn-formalsch}) (and hence so is $Y$, and $\mathscr{X}$ and $\mathscr{Y}$ are coherent universally adhesive rigid spaces), then one can always replace $X$ and $Y$ by distinguished formal models by taking admissible blow-up along an ideal of definition of finite type.
In this situation, the closed immersion $i\colon Y\hookrightarrow X$ is automatically finite presented.
In particular, one can further assume in (c) that each $i_{\lambda}$ is a closed immersion of finite presentation.}
\end{rem}

\begin{dfn}\label{dfn-closedimmrigid1coherent}{\rm 
A morphism $\iota\colon\mathscr{Y}\rightarrow\mathscr{X}$ of coherent universally Noetherian rigid spaces\index{rigid space!universally Noetherian rigid space@universally Noetherian ---!coherent universally Noetherian rigid space@coherent --- ---} (\ref{dfn-coherentuniversallyadhesiverigidspaces}) is said to be a {\em closed immersion} if it satisfies the conditions in \ref{prop-closedimmrigid1}.}
\end{dfn}

\begin{prop}\label{prop-closedimmrigid2}
Let $\mathscr{X}$ be a coherent universally Noetherian rigid space, and $\mathscr{K}\subseteq\O_{\mathscr{X}}$ a coherent ideal.
Then there exists uniquely up to canonical isomorphisms a closed immersion $\mathscr{Y}\hookrightarrow\mathscr{X}$ that induces an isomorphism of locally ringed spaces $(\ZR{\mathscr{Y}},\O_{\mathscr{Y}})\stackrel{\sim}{\rightarrow}(Z,\O_{\mathscr{X}}/\mathscr{K})$, where $Z$ is the support of the sheaf $\O_{\mathscr{X}}/\mathscr{K}$.
Moreover, any closed immersion $\mathscr{Y}\hookrightarrow\mathscr{X}$ is isomorphic to the one obtained in this way by a uniquely determined coherent ideal of $\O_{\mathscr{X}}$.
\end{prop}

The uniquely determined coherent ideal for a closed immersion $\mathscr{Y}\hookrightarrow\mathscr{X}$ stated in the proposition will be called the {\em defining ideal}.

\begin{proof}
Let $X$ be a coherent universally rigid-Noetherian formal model of $\mathscr{X}$, which we assume without loss of generality to have an invertible ideal of definition $\mathscr{I}_X$.
Then by \ref{thm-tateacyclic20} we have a finitely presented ideal $\mathscr{K}_X$ of $\O_X$ such that $\mathscr{K}^{\rig}_X=\mathscr{K}$; indeed, we have finitely presented sheaf $\mathscr{K}'_X$ with the map $\mathscr{K}'_X\rightarrow\O_X$ such that $\mathscr{K}^{\prime\rig}_X=\mathscr{K}$; then $\mathscr{K}_X=\mathscr{I}^n_X\mathscr{K'}_X$ for a sufficiently large $n>0$ is a finitely presented ideal of $\O_X$ having the same property.
Let $i\colon Y\hookrightarrow X$ be the closed immersion of finite presentation corresponding to $\mathscr{K}_X$ ({\bf \ref{ch-formal}}.\ref{prop-closedimmformal6}).
Then we have the morphism $\iota=i^{\rig}\colon\mathscr{Y}=Y^{\rig}\hookrightarrow\mathscr{X}$ of coherent rigid spaces.

As in the proof of \ref{prop-closedimmrigid1} one has a cofinal family of formal models $\{i_{\lambda}\colon Y_{\lambda}\rightarrow X_{\lambda}\}$ of $\iota$ consisting of closed immersions between coherent universally rigid-Noetherian formal schemes that dominates $i\colon Y\hookrightarrow X$.
We have $\ZR{\mathscr{Y}}=\varprojlim Y_{\lambda}$.
For each $\lambda$, let $\mathscr{K}_{X_{\lambda}}$ be the defining ideal of the closed immersion $i_{\lambda}$.
Set $\til{\mathscr{K}}=\varinjlim\sp^{-1}_{X_{\lambda}}\mathscr{K}_{X_{\lambda}}$, which is an ideal of $\O^{\int}_{\mathscr{X}}$ such that $\til{\mathscr{K}}\otimes\O_{\mathscr{X}}=\mathscr{K}$.
Since $\O^{\int}_{\mathscr{X}}/\til{\mathscr{K}}$ is easily seen to be $\mathscr{I}$-torsion free (where $\mathscr{I}$ is an ideal of definition (\ref{dfn-ZRstrsheaf2})), we deduce
$$
\Supp\O^{\int}_{\mathscr{X}}/\til{\mathscr{K}}=\Supp\O_{\mathscr{X}}/\mathscr{K}.
$$
Finally, by {\bf \ref{ch-pre}}.\ref{thm-projlimcohspacepres} (3) one sees that the map $\ZR{\mathscr{Y}}\rightarrow\ZR{\mathscr{X}}$ is closed, and then it is easy to see that it is actually an isomorphism onto the closed subset $\Supp \O^{\int}_{\mathscr{X}}/\til{\mathscr{K}}$.

To show the last assertion, let $\iota\colon\mathscr{Y}\hookrightarrow\mathscr{X}$ be a closed immersion, and take a formal model $i\colon Y\hookrightarrow X$ that is a closed immersion of finite presentation between coherent universally rigid-Noetherian formal schemes.
Since the question is local on $\mathscr{X}$, we may assume that $X$ is affine of the form $X=\Spf A$ where $A$ is a t.u.\ rigid-Noetherian ring\index{t.u. rigid-Noetherian ring@t.u.\ rigid-Noetherian ring} ({\bf \ref{ch-formal}}.\ref{dfn-tuaringadmissible}).
Then the closed immersion $i\colon Y\hookrightarrow X$ comes from a surjective map $A\rightarrow B$, where $Y=\Spf B$.
Let $K$ be the kernel of $A\rightarrow B$.
Then $K^{\Delta}$ is an adically quasi-coherent ideal of $\O_X$ that gives the kernel of $\O_X\rightarrow i_{\ast}\O_Y$.
It is then easy to check that the coherent ideal $\mathscr{K}=(K^{\Delta})^{\rig}$ of $\O_{\mathscr{X}}$ recovers, up to isomorphism, the closed immersion $\mathscr{Y}\hookrightarrow\mathscr{X}$.
The uniqueness is clear.
\end{proof}

\begin{cor}\label{cor-propclosedimmrigidx2x}
Let $\iota\colon\mathscr{Y}\rightarrow\mathscr{X}$ be a morphism of coherent universally Noetherian rigid spaces, and $\{\mathscr{U}_{\alpha}\}_{\alpha\in L}$ a covering in the site $\mathscr{X}_{\ad}$.
Then $\iota$ is a closed immersion if and only if for each $\alpha\in L$ the base change $\mathscr{Y}\times_{\mathscr{X}}\mathscr{U}_{\alpha}\rightarrow\mathscr{U}_{\alpha}$ is a closed immersion. \hfill$\square$
\end{cor}

The last corollary allows one to define closed immersions between locally universally Noetherian rigid spaces\index{rigid space!universally Noetherian rigid space@universally Noetherian ---!locally universally Noetherian rigid space@locally --- ---} (\ref{dfn-universallyadhesiverigidspaces}) consistently as follows:
\begin{dfn}\label{dfn-closedimmrigid1}{\rm 
(1) A morphism $\varphi\colon\mathscr{Y}\rightarrow\mathscr{X}$ of locally universally Noetherian rigid spaces\index{rigid space!universally Noetherian rigid space@universally Noetherian ---!locally universally Noetherian rigid space@locally --- ---} (\ref{dfn-universallyadhesiverigidspaces}) is said to be a {\em closed immersion} if it is coherent (\ref{dfn-qcptmorphismrigid1}) and $\mathscr{X}$ has a covering $\{\mathscr{X}_{\alpha}\}_{\alpha\in L}$ by coherent rigid spaces as in \ref{dfn-generalrigidspace1} such that for each $\alpha\in L$ the base change $\mathscr{Y}\times_{\mathscr{X}}\mathscr{X}_{\alpha}\rightarrow\mathscr{X}_{\alpha}$ is a closed immersion in the sense as in \ref{dfn-closedimmrigid1coherent}. 

(2) A {\em closed rigid subspace}\index{rigid subspace!closed rigid subspace@closed ---} of a rigid space $\mathscr{X}$ is an equivalence class (by $\mathscr{X}$-isomorphisms) of closed immersions $\mathscr{Y}\hookrightarrow\mathscr{X}$.}
\end{dfn}

By \ref{prop-closedimmrigid2} any closed immersion as in (1) comes from a uniquely determined coherent ideal, which we also call the {\em defining ideal}.

\begin{prop}\label{prop-corclosedimmrigid2}
A finite morphism $\iota\colon\mathscr{Y}\rightarrow\mathscr{X}$ of locally universally Noetherian rigid spaces\index{rigid space!universally Noetherian rigid space@universally Noetherian ---!locally universally Noetherian rigid space@locally --- ---} is a closed immersion if and only if the map $\O_{\mathscr{X}}\rightarrow\ZR{\iota}_{\ast}\O_{\mathscr{Y}}$ $($cf.\ {\rm \ref{lem-finitemorrigid8}}$)$ is surjective. \hfill$\square$
\end{prop}

\begin{prop}\label{prop-closedimmrigid2xx}
Let $\mathscr{X}$ be a Stein affinoid\index{affinoid!Stein affinoid@Stein ---} {\rm (\ref{dfn-cohomologyrigidsp1str})}, and $\mathscr{Y}\hookrightarrow\mathscr{X}$ a closed immersion. 
Then $\mathscr{Y}$ is a Stein affinoid.
\end{prop}

\begin{proof}
Take a formal model $X=\Spf A$ of $\mathscr{X}$ by a t.u.\ rigid-Noetherian ring $A$ such that $\Spec A\setminus V(I)$ is affine, where $I\subseteq A$ is an ideal of definition.
As in the proof of \ref{prop-closedimmrigid2} the defining ideal of $\mathscr{Y}\hookrightarrow\mathscr{X}$ comes from a finitely presented ideal of $\O_X$ and thus from a finitely presented ideal $K\subseteq A$.
Set $B=A/K$ and $Y=\Spf B$.
Then the closed immersion $Y\hookrightarrow X$ of finite presentation gives a formal model of $\mathscr{Y}\hookrightarrow\mathscr{X}$.
Clearly, $\Spec B\setminus V(IB)$ is affine.
\end{proof}

\begin{prop}\label{prop-opencomplement}
Let $\mathscr{Y}\hookrightarrow\mathscr{X}$ be a closed immersion of locally universally Noetherian rigid spaces.\index{rigid space!universally Noetherian rigid space@universally Noetherian ---!locally universally Noetherian rigid space@locally --- ---}
Then the induced map $\ZR{\mathscr{Y}}\rightarrow\ZR{\mathscr{X}}$ maps $\ZR{\mathscr{Y}}$ homeomorphically onto an overconvergent closed subset of $\ZR{\mathscr{X}}$.
\end{prop}

\begin{proof}
Let $x\in\ZR{\mathscr{Y}}$ and $y$ a generization of $x$ in $\ZR{\mathscr{X}}$.
We need to show that $y$ belongs to $\ZR{\mathscr{Y}}$.
Let $\alpha\colon\Spf\widehat{V}_x\rightarrow\ZR{\mathscr{X}}$ and $\beta\colon\Spf\widehat{W}_x\rightarrow\ZR{\mathscr{Y}}$ be the associated rigid points (\ref{dfn-ZRpoints32}).
The surjective map $\O^{\int}_{\ZR{\mathscr{X}},x}\rightarrow\O^{\int}_{\ZR{\mathscr{Y}}_x}$ (cf.\ {\bf \ref{ch-formal}}.\ref{cor-closedimmformal42}) gives rise to a local surjective map $h\colon V_x\rightarrow W_x$.
Since $h$ is, at the same time, $a$-adic, it is an injective map (cf.\ the proof of {\bf \ref{ch-pre}}.\ref{prop-maxspe4}).
Hence $h$ is an isomorphism, and thus $\widehat{V}_x\cong\widehat{W}_x$.
Now since $y$ belongs to the image of $\alpha$, we deduce that $y$ belongs to $\ZR{\mathscr{Y}}$, as desired.
\end{proof}

\begin{prop}\label{prop-closedimmrigid5}
{\rm (1)} A closed immersion is a finite morphism.

{\rm (2)} If $\varphi\colon\mathscr{Z}\rightarrow\mathscr{Y}$ and $\psi\colon\mathscr{Y}\rightarrow\mathscr{X}$ are closed immersions, then so is the composition $\psi\circ\varphi$.

{\rm (3)} If $\varphi\colon\mathscr{X}\rightarrow\mathscr{X}'$ and $\psi\colon\mathscr{Y}\rightarrow\mathscr{Y}'$ are two closed immersions over a rigid space $\mathscr{S}$ such that either $\mathscr{X}$ and $\mathscr{X}'$ or $\mathscr{X}$ and $\mathscr{Y}'$ are locally of finite type over $\mathscr{S}$, then the induced morphism $\varphi\times_{\mathscr{S}}\psi\colon\mathscr{X}\times_{\mathscr{S}}\mathscr{Y}\rightarrow\mathscr{X}'\times_{\mathscr{S}}\mathscr{Y}'$ is a closed immersion.

{\rm (4)} If $\varphi\colon\mathscr{X}\rightarrow\mathscr{Y}$ is a closed immersion over a rigid space $\mathscr{S}$ and $\mathscr{S}'\rightarrow\mathscr{S}$ is a morphism of rigid spaces such that either $\mathscr{X}$ and $\mathscr{Y}$ are locally of finite type over $\mathscr{S}$ or that $\mathscr{S}'$ is locally of finite type over $\mathscr{S}$, then the induced morphism $\varphi_{\mathscr{S}'}\colon\mathscr{X}\times_{\mathscr{S}}\mathscr{S}'\rightarrow\mathscr{Y}\times_{\mathscr{S}}\mathscr{S}'$ is a closed immersion.
\end{prop}

\begin{proof}
(1) is clear.
(2) can be shown by an argument similar to that in \ref{prop-finitemorrigid2} (1).
(3) follows easily from {\bf \ref{ch-formal}}.\ref{prop-closedimmformal5} (2).
Finally, (4) follows due to {\bf \ref{ch-pre}}.\ref{prop-basechangestable}.
\end{proof}
\index{immersion!closed immersion of rigid spaces@closed --- (of rigid spaces)|)}

\subsubsection{Irreducible rigid spaces}\label{subsub-irreduciblerigidspaces}
\begin{dfn}\label{dfn-irreduciblerigidspace}{\rm 
We say that a locally universally Noetherian rigid space $\mathscr{X}$ is {\it irreducible}\index{rigid space!irreducible rigid space@irreducible ---} (or, more adequately, {\it globally irreducible}) if the following condition is satisfied: if $\ZR{\mathscr{X}}=\ZR{\mathscr{Y}}\cup\ZR{\mathscr{Z}}$, where $\mathscr{Y}$ and $\mathscr{Z}$ are closed rigid subspaces\index{rigid subspace!closed rigid subspace@closed ---} (\ref{dfn-closedimmrigid1}) of $\mathscr{X}$, then either $\ZR{\mathscr{X}}=\ZR{\mathscr{Y}}$ or $\ZR{\mathscr{X}}=\ZR{\mathscr{Z}}$ holds.}
\end{dfn}

\subsubsection{Open complement}\label{subsub-closedimmrigidopencomp}
\begin{dfn}\label{dfn-opencomplement}{\rm 
Let $\mathscr{X}$ be a locally universally Noetherian rigid space\index{rigid space!universally Noetherian rigid space@universally Noetherian ---!locally universally Noetherian rigid space@locally --- ---}, and $\mathscr{Y}\hookrightarrow\mathscr{X}$ a closed subspace.
The {\em open complement}\index{open complement of a closed rigid subspace@open complement (of a closed rigid subspace)} of $\mathscr{Y}$ in $\mathscr{X}$ is the open subspace $\mathscr{U}$ of $\mathscr{X}$, denoted by $\mathscr{X}\setminus\mathscr{Y}$, such that $\ZR{\mathscr{U}}=\ZR{\mathscr{X}}\setminus\ZR{\mathscr{Y}}$.}
\end{dfn}

The Zariski-Riemann space $\ZR{\mathscr{U}}$ of the open complement $\mathscr{U}$ is, therefore, an overconvergent open subset\index{overconvergent!overconvergent subset@--- subset} of $\mathscr{X}$ due to \ref{prop-opencomplement}.
The actual construction of $\mathscr{U}$ is given as follows:
\begin{const}\label{const-opencomplement}{\rm 
It is enough to perform the construction in the case where $\mathscr{X}$ is coherent.
Take a cofinal family of formal models $\{i_{\lambda}\colon Y_{\lambda}\rightarrow X_{\lambda}\}$ of $\mathscr{Y}\hookrightarrow\mathscr{X}$ consisting of closed immersions between coherent universally rigid-Noetherian formal schemes (\ref{prop-closedimmrigid1}).
For any $\lambda$, let $U_{\lambda}$ be the open complement of $Y_{\lambda}$ in $X_{\lambda}$.
Then $\bigcup_{\lambda\in\Lambda}\sp^{-1}_{X_{\lambda}}(U_{\lambda})$ gives the complement $\ZR{\mathscr{X}}\setminus\ZR{\mathscr{Y}}$.

Now define for any $\alpha$ the coherent rigid space $\mathscr{U}_{\lambda}=U^{\rig}_{\lambda}$.
For $\lambda\leq\mu$ we have an obvious coherent open immersion $\mathscr{U}_{\lambda}\hookrightarrow\mathscr{U}_{\mu}$.
Hence $\{\mathscr{U}_{\lambda}\}$ is an increasing family if coherent rigid spaces and defines a quasi-separated rigid space (stretch of coherent rigid spaces; cf.\ \ref{dfn-admissiblesite31} (1)) $\mathscr{U}=\bigcup_{\lambda\in\Lambda}\mathscr{U}_{\lambda}$.
There exists an obvious open immersion $\mathscr{U}\hookrightarrow\mathscr{X}$.
By \ref{prop-zariskiriemanntop2gen} we have $\ZR{\mathscr{U}}=\ZR{\mathscr{X}}\setminus\ZR{\mathscr{Y}}$, as desired.}
\end{const}

\subsubsection{Closed subspaces of an affinoid}\label{subsub-closedsubspacesofaffinoids}
Let $\mathscr{X}=(\Spf A)^{\rig}$ be a universally Noetherian affinoid\index{affinoid!universally Noetherian affinoid@universally Noetherian ---}, where $A$ is a t.u.\ rigid-Noetherian\index{t.u. rigid-Noetherian ring@t.u.\ rigid-Noetherian ring} ring with a finitely generated ideal of definition $I\subseteq A$.
Consider the associated Noetherian scheme\index{scheme!associated scheme to an affinoid@associated --- to an affinoid}\index{affinoid!associated scheme to an affinoid@associated scheme to an ---} $s(\mathscr{X})=\Spec A\setminus V(I)$ (\S\ref{subsub-associatedscheme}) and the comparison map\index{comparison map} $s\colon\ZR{\mathscr{X}}\rightarrow s(\mathscr{X})$ (\S\ref{subsub-comparisonmapassaff}).

For any closed subscheme $Z\subseteq s(\mathscr{X})$, take a closed subscheme $\ovl{Z}\subseteq\Spec A$ of finite presentation such that $Z=\ovl{Z}\cap s(\mathscr{X})$.
Take the $I$-adic completion $\widehat{\ovl{Z}}$ of $\ovl{Z}$, which is a closed formal subscheme of $\Spf A$ of finite presentation.
Then we have the closed subspace $\mathscr{Z}=(\widehat{\ovl{Z}})^{\rig}$ of $\mathscr{X}$.
By the construction, if $\mathscr{J}$ is the defining coherent ideal of $Z$ on $s(\mathscr{X})$, then the coherent ideal $\mathscr{J}^{\rig}$ (\ref{ntn-affinoidscohsheaf}) gives the defining ideal of $\mathscr{Z}$ on $\mathscr{X}$.
Notice that, since the coherent ideal $\mathscr{J}^{\rig}$ is nothing but the pull-back sheaf $s^{\ast}\mathscr{J}=\mathscr{J}\O_{\mathscr{X}}$ (cf.\ \ref{prop-comparisonmapassaff} and \ref{prop-comparisonmapassaffprop1}), the locally ringed space $(\ZR{\mathscr{Z}},\O_{\mathscr{Z}})$ is canonically isomorphic to the fiber product 
$$
Z\times_{s(\mathscr{X})}\ZR{\mathscr{X}}
$$
in the category $\LRsp$.

\begin{prop}\label{prop-closedsubspaceaffinoid2}
{\rm (1)} For any closed subset $Z$ $($resp.\ open subset $U)$ of $s(\mathscr{X})$ the pull-back $s^{-1}(Z)$ $($resp.\ $s^{-1}(U))$ in $\ZR{\mathscr{X}}$ is an overconvergent closed $($resp.\ open$)$ subset\index{overconvergent!overconvergent subset@--- subset} {\rm (\ref{dfn-separation2})}.

{\rm (2)} For any closed subscheme $Z$ of $s(\mathscr{X})$ there exists uniquely a closed subspace $\mathscr{Z}$ of $\mathscr{X}$ such that $(\ZR{\mathscr{Z}},\O_{\mathscr{Z}})$ is isomorphic to the fiber product $Z\times_{s(\mathscr{X})}\ZR{\mathscr{X}}$ in the category of locally ringed spaces.
Moreover, the correspondence $Z\mapsto\mathscr{Z}$ establishes the bijection between the set of all closed subschemes $($resp.\ irreducible closed subschemes$)$ of $s(\mathscr{X})$ and the set of all closed subspaces $($resp.\ irreducible closed subspaces$)$ of $\mathscr{X}$.
\end{prop}

\begin{proof}
(1) follows from \ref{prop-separation21} and \ref{rem-comparisonmapfactorsep}.
Since closed subspaces of $\mathscr{X}$ are determined by its coherent defining ideal of $\O_{\mathscr{X}}$, (2) follows from \ref{prop-comparisonmapassaffprop1}.
\end{proof}

\begin{ntn}\label{ntn-pushpullclosed}{\rm 
For a closed subscheme $Z$ of $s(\mathscr{X})$ the corresponding closed subspace $\mathscr{Z}$ of $\mathscr{X}$ as in \ref{prop-closedsubspaceaffinoid2} (2) is denoted by 
$$
s^{\ast}Z.
$$
On the other hand, for a closed subspace $\mathscr{Z}$ of $\mathscr{X}$ the corresponding closed subscheme of $s(\mathscr{X})$ is (consistently) denoted by 
$$
s(\mathscr{Z}).
$$}
\end{ntn}

Notice that we have $\ZR{s^{\ast}Z}=s^{-1}(Z)$ as a topological space.

\subsection{Immersions}\label{sub-immersionrigid}
\subsubsection{Immersions and rigid subspaces}\label{subsub-immersionrigid}
\index{immersion!immersion of rigid spaces@--- (of rigid spaces)|(}
\begin{dfn}\label{dfn-immersionrigid}{\rm 
A morphism $\varphi\colon\mathscr{Y}\rightarrow\mathscr{X}$ of locally universally Noetherian rigid spaces\index{rigid space!universally Noetherian rigid space@universally Noetherian ---!locally universally Noetherian rigid space@locally --- ---} (\ref{dfn-universallyadhesiverigidspaces}) is said to be an {\em immersion} if it is a composition $\varphi=j\circ i$ where $i$ is a closed immersion and $j$ is an open immersion.}
\end{dfn}

\begin{prop}\label{prop-immersionrigid1}
Let $\varphi\colon\mathscr{Y}\rightarrow\mathscr{X}$ be a morphism of locally universally Noetherian rigid spaces\index{rigid space!universally Noetherian rigid space@universally Noetherian ---!locally universally Noetherian rigid space@locally --- ---}, and $\{\mathscr{V}_{\alpha}\hookrightarrow\mathscr{X}\}_{\alpha\in L}$ a family of open immersions such that $\varphi$ factors through the open immersion $\bigcup_{\alpha\in L}\mathscr{V}_{\alpha}\hookrightarrow\mathscr{X}$.
Then $\varphi$ is an immersion if and only if for any $\alpha\in L$ the base change $\mathscr{Y}\times_{\mathscr{X}}\mathscr{V}_{\alpha}\rightarrow \mathscr{V}_{\alpha}$ is an immersion.
\end{prop}

\begin{proof}
The `only if' part is clear.
Let us show the other part.
Take an open immersion $\mathscr{U}_{\alpha}\hookrightarrow\mathscr{V}_{\alpha}$ for each $\alpha\in L$ such that the immersion $\mathscr{Y}\times_{\mathscr{X}}\mathscr{V}_{\alpha}\rightarrow\mathscr{V}_{\alpha}$ factors by the closed immersion $\mathscr{Y}\times_{\mathscr{X}}\mathscr{V}_{\alpha}\hookrightarrow\mathscr{U}_{\alpha}$.
Set $\mathscr{U}=\bigcup_{\alpha\in L}\mathscr{U}_{\alpha}$, which is an open subspace of $\mathscr{X}$.
Then by \ref{cor-propclosedimmrigidx2x} the morphism $\mathscr{Y}\rightarrow\mathscr{U}$ is a closed immersion.
\end{proof}

\begin{prop}\label{prop-immersionrigid5}
{\rm (1)} If $\varphi\colon\mathscr{Z}\rightarrow\mathscr{Y}$ and $\psi\colon\mathscr{Y}\rightarrow \mathscr{X}$ are immersions, then so is the composition $\psi\circ\varphi$.

{\rm (2)} If $\varphi\colon\mathscr{X}\rightarrow\mathscr{X}'$ and $\psi\colon\mathscr{Y}\rightarrow\mathscr{Y}'$ are two immersions over a rigid space $\mathscr{S}$ such that either $\mathscr{X}$ and $\mathscr{X}'$ or $\mathscr{X}$ and $\mathscr{Y}'$ are locally of finite type over $\mathscr{S}$, then the induced morphism $\varphi\times_{\mathscr{S}}\psi\colon\mathscr{X}\times_{\mathscr{S}}\mathscr{Y}\rightarrow\mathscr{X}'\times_{\mathscr{S}}\mathscr{Y}'$ is an immersion.

{\rm (3)} If $\varphi\colon\mathscr{X}\rightarrow\mathscr{Y}$ is an immersion over a rigid space $\mathscr{S}$ and $\mathscr{S}'\rightarrow\mathscr{S}$ is a morphism of rigid spaces such that either $\mathscr{X}$ and $\mathscr{Y}$ are locally of finite type over $\mathscr{S}$ or that $\mathscr{S}'$ is locally of finite type over $\mathscr{S}$, then the induced morphism $\varphi_{\mathscr{S}'}\colon\mathscr{X}\times_{\mathscr{S}}\mathscr{S}'\rightarrow\mathscr{Y}\times_{\mathscr{S}}\mathscr{S}'$ is an immersion.
\end{prop}

To show (1), we need to prove the following lemma:
\begin{lem}\label{lem-immersionrigidcomposites}
Let $\mathscr{Z}\hookrightarrow\mathscr{Y}$ be an open immersion, and $\mathscr{Y}\hookrightarrow \mathscr{X}$ a closed immersion.
Then the composition $\mathscr{Z}\hookrightarrow\mathscr{X}$ is an immersion.
\end{lem}

\begin{proof}
We identify $\mathscr{Y}$ with a closed subspace of $\mathscr{X}$, and $\mathscr{Z}$ an open subspace of $\mathscr{Y}$.
There exists an open subset $\mathfrak{U}$ of $\ZR{\mathscr{X}}$ that contains $\ZR{\mathscr{Z}}$ as a closed subset.
Let $\mathscr{U}\subseteq\mathscr{X}$ be the open subspace supported on $\mathfrak{U}$.
Then $\mathscr{Z}\hookrightarrow\mathscr{U}$ is a closed subset given by the defining ideal of $\mathscr{Y}\hookrightarrow\mathscr{X}$ restricted on $\mathscr{U}$.
\end{proof}

\begin{proof}[Proof of Proposition {\rm \ref{prop-immersionrigid5}}]
(1) follows from \ref{lem-immersionrigidcomposites}.
The statement (2) follows from \ref{prop-closedimmrigid5} (3) and the corresponding (obvious) fact on open immersions.
Due to {\bf \ref{ch-pre}}.\ref{prop-basechangestable} the assertion (3) follows automatically.
\end{proof}

\begin{prop}\label{prop-immersionrigidsep}
Let $\varphi\colon\mathscr{X}\rightarrow\mathscr{Y}$ be a morphism of locally universally Noetherian rigid spaces\index{rigid space!universally Noetherian rigid space@universally Noetherian ---!locally universally Noetherian rigid space@locally --- ---}.
Then the diagonal morphism $\Delta_{\mathscr{X}}\colon\mathscr{X}\rightarrow\mathscr{X}\times_{\mathscr{Y}}\mathscr{X}$ is an immersion.
\end{prop}

\begin{proof}
As in the proof of {\bf \ref{ch-formal}}.\ref{prop-diagonalimmersionadic} we may assume in view of \ref{prop-immersionrigid1} that $\varphi$ has a formal model of the form $\Spf A\rightarrow\Spf B$.
The diagonal map $\Delta_{\mathscr{X}}$ in this case is a closed immersion due to {\bf \ref{ch-formal}}.\ref{cor-closedimmformal42}.
\end{proof}

\begin{dfn}\label{dfn-rigidsubspace}{\rm 
Let $\mathscr{X}$ be a locally universally Noetherian rigid space.
A {\em $($locally closed$)$ rigid subspace}\index{rigid subspace} of $\mathscr{X}$ is an $\mathscr{X}$-isomorphism class of immersions $\mathscr{Y}\hookrightarrow\mathscr{X}$.}
\end{dfn}
\index{immersion!immersion of rigid spaces@--- (of rigid spaces)|)}

\subsection{Separated morphisms and proper morphisms}\label{sub-seppropmorrigid}
\subsubsection{Closed morphisms}\label{subsub-seppropmorrigidclosed}
\index{morphism of rigid spaces@morphism (of rigid spaces)!closed morphism of rigid spaces@closed ---|(}
\begin{prop}\label{prop-rigidclosedmap}
Let $\varphi\colon\mathscr{X}\rightarrow\mathscr{Y}$ be a morphism between coherent rigid spaces.
Then the following conditions are equivalent$:$
\begin{itemize}
\item[{\rm (a)}] $\ZR{\varphi}\colon\ZR{\mathscr{X}}\rightarrow\ZR{\mathscr{Y}}$ is a closed map$;$
\item[{\rm (b)}] for any point $x\in\ZR{\mathscr{X}}$ we have
$$
\ZR{\varphi}(\ovl{\{x\}})=\ovl{\{\ZR{\varphi}(x)\}};
$$
\item[{\rm (c)}] any distinguished formal model $f\colon X\rightarrow Y$ of $\varphi$ is closed$;$
\item[{\rm (d)}] there exists a cofinal subset $\mathscr{C}$ of formal models of $\varphi$ consisting of closed maps.
\end{itemize}
\end{prop}

\begin{proof}
The implication (a) $\Rightarrow$ (b) is a consequence from general topology.
To show (b) $\Rightarrow$ (c), let $f\colon X\rightarrow Y$ be a distinguished formal model of $\varphi$, and take an ideal of definition $\mathscr{I}_Y$ of finite type of $Y$.
Set $\mathscr{I}_X=\mathscr{I}_Y\O_X$.
Denote by $X_0$ and $Y_0$ the closed subschemes of $X$ and $Y$ defined respectively by $\mathscr{I}_X$ and $\mathscr{I}_Y$, and let $f_0\colon X_0\rightarrow Y_0$ be the induced morphism.
We need to show that $f(Z)$ for any closed subset $Z\subset X$ is closed.
By \cite[$\mathbf{II}$, (7.2.2)]{EGA} it suffices to show that $f(Z)$ is closed under specialization.

Take $z=f(x)\in Z$.
By \ref{prop-ZRpoints4} we can find a point $u\in\ZR{\mathscr{X}}$ such that $\sp_X(u)=x$.
Since the specialization map $\sp_X$ is closed (\ref{thm-ZRcompact} (2)), we have $\sp_Y(\ovl{\{\ZR{\varphi}(u)\}})=\ovl{\{z\}}$.
It then follows from our hypothesis that
$$
\sp_Y(\ovl{\{\ZR{\varphi}(u)\}})=\sp_Y(\ZR{\varphi}(\ovl{\{u\}}))=f(\sp_X(\ovl{\{u\}})).
$$
Now since $Z$ is closed, we have $\sp_X(\ovl{\{u\}})\subseteq Z$.
Hence we have $\ovl{\{z\}}\subseteq F(Z)$, as desired.

The implication (c) $\Rightarrow$ (d) is obvious; (d) $\Rightarrow$ (a) follows from {\bf \ref{ch-pre}}.\ref{thm-projlimcohspacepres} (3).
\end{proof}

\begin{dfn}\label{dfn-rigidclosedmap}{\rm 
A morphism $\varphi\colon\mathscr{X}\rightarrow\mathscr{Y}$ of rigid spaces is said to be {\em closed} if the associated map $\ZR{\varphi}\colon\ZR{\mathscr{X}}\rightarrow\ZR{\mathscr{Y}}$ is a closed map.}
\end{dfn}

The following proposition follows immediately from \ref{prop-closedimmrigid2}:
\begin{prop}\label{prop-rigidclosedmapx1}
An immersion\index{immersion!immersion of rigid spaces@--- (of rigid spaces)} between locally universally Noetherian rigid spaces is a closed immersion\index{immersion!closed immersion of rigid spaces@closed --- (of rigid spaces)} if and only if it is closed.
\end{prop}

\begin{proof}
The ``only if'' part is clear.
Let $\varphi\colon\mathscr{Y}\hookrightarrow\mathscr{X}$ be an immersion, which we suppose to be closed.
By \ref{cor-propclosedimmrigidx2x}, we may assume that $\mathscr{X}$ is coherent.
Write $\varphi=j\circ\iota$, where $\iota\colon\mathscr{Y}\hookrightarrow\mathscr{U}$ is a closed immersion, and $j\colon\mathscr{U}\hookrightarrow\mathscr{X}$ is an open immersion.
Since the image of $\ZR{\mathscr{Y}}$ in $\ZR{\mathscr{X}}$ is closed, and hence is quasi-compact, we may assume that $\mathscr{U}$ is quasi-compact, and hence that it is a quasi-compact open subspace of $\mathscr{X}$.
Then one can take, using \ref{prop-blowups4111} if necessary, the coherent rigid-Noetherian distinguished formal models $Y\hookrightarrow U\hookrightarrow X$ consisting of a closed immersion followed by a coherent open immersion, so that the composition $Y\hookrightarrow X$ gives a formal model of $\phi$.
By \ref{prop-rigidclosedmap}, $Y\hookrightarrow X$ is closed, and is an immersion.
Then by {\bf \ref{ch-formal}}.\ref{cor-univclosedadicred0}, one deduces that $Y\hookrightarrow X$ is a closed immersion, and hence $\varphi$ is a closed immersion, as desired.
\end{proof}
\index{morphism of rigid spaces@morphism (of rigid spaces)!closed morphism of rigid spaces@closed ---|)}

\subsubsection{Separated morphisms and proper morphisms}\label{subsub-seppropmorrigidsep}
\index{morphism of rigid spaces@morphism (of rigid spaces)!separated morphism of rigid spaces@separated ---|(}
\index{morphism of rigid spaces@morphism (of rigid spaces)!proper morphism of rigid spaces@proper ---|(}
\begin{dfn}\label{dfn-seppropmorrigid1}{\rm 
Let $\varphi\colon\mathscr{X}\rightarrow\mathscr{Y}$ be a morphism of rigid spaces.

(1) The morphism $\varphi$ is said to be {\em universally closed}\index{morphism of rigid spaces@morphism (of rigid spaces)!universally closed morphism of rigid spaces@universally closed ---} if for any morphism $\mathscr{Z}\rightarrow\mathscr{Y}$ of rigid spaces the base change $\varphi_{\mathscr{Z}}\colon\mathscr{X}\times_{\mathscr{Y}}\mathscr{Z}\rightarrow\mathscr{Z}$ is closed.

(2) The morphism $\varphi$ is said to be {\em separated} if the diagonal morphism $\Delta_{\mathscr{X}}\colon \mathscr{X}\rightarrow\mathscr{X}\times_{\mathscr{Y}}\mathscr{X}$ is quasi-compact and closed.

(3) The morphism $\varphi$ is said to be {\em proper} if it is of finite type, separated, and universally closed.}
\end{dfn}

Notice that, due to \ref{prop-immersionrigidsep} and \ref{prop-rigidclosedmapx1}, a morphism $\varphi\colon\mathscr{X}\rightarrow\mathscr{Y}$ of locally universally Noetherian rigid spaces\index{rigid space!universally Noetherian rigid space@universally Noetherian ---!locally universally Noetherian rigid space@locally --- ---} (\ref{dfn-universallyadhesiverigidspaces}) is separated if and only if the diagonal morphism $\Delta_{\mathscr{X}}\colon \mathscr{X}\rightarrow\mathscr{X}\times_{\mathscr{Y}}\mathscr{X}$ is a closed immersion\index{immersion!closed immersion of rigid spaces@closed --- (of rigid spaces)}.
The following proposition follows immediately from the definition.

\begin{prop}\label{prop-seppropmorrigidxx1}
A separated morphism is quasi-separated\index{morphism of rigid spaces@morphism (of rigid spaces)!quasi-separated morphism of rigid spaces@quasi-separated ---}. \hfill$\square$
\end{prop}

\begin{thm}\label{thm-seppropmorrigid1}
Let $\varphi\colon\mathscr{X}\rightarrow\mathscr{Y}$ be a finite type morphism of coherent rigid spaces.
Then the following conditions are equivalent$:$
\begin{itemize}
\item[{\rm (a)}] $\varphi$ is separated$;$
\item[{\rm (b)}] any distinguished formal model $f\colon X\rightarrow Y$ of $\varphi$ is separated$;$
\item[{\rm (c)}] there exists a separated formal model $f\colon X\rightarrow Y;$
\item[{\rm (d)}] there exists a cofinal set of formal models of $\varphi$ consisting of separated morphisms.
\end{itemize}
\end{thm}

\begin{proof}
First let us prove the implication (a) $\Rightarrow$ (b).
Take a distinguished formal model $f\colon X\rightarrow Y$ of $\varphi$, and consider the commutative diagram
$$
\xymatrix{X'\ar[d]\ar[r]&Z\ar[d]\\ X\ar[r]&X\times_YX\rlap{,}}
$$
where the right vertical map is the admissible blow-up along an ideal of definition of finite type, and the left vertical map is the strict transform.
Since $X$ is a distinguished formal model of $\mathscr{X}$, the map $X'\rightarrow X$, which is again an admissible blow-up, is surjective.
Hence to show that $X\rightarrow X\times_YX$ is closed, it suffices to show that $X'\rightarrow X\times_YX$ is closed.
Since the admissible blow-up is proper, we only need to show that $X'\rightarrow Z$ is closed.
But since this is a distinguished formal model of $\mathscr{X}\rightarrow\mathscr{X}\times_{\mathscr{Y}}\mathscr{X}$, it is closed due to \ref{prop-rigidclosedmap}.

The implication (b) $\Rightarrow$ (c) is obvious.
To show (c) $\Rightarrow$ (d), take a separated formal model $f\colon X\rightarrow Y$ of $\varphi$.
For any formal model $g\colon X'\rightarrow Y'$ of $\varphi$ there exists a formal model $h$ dominating $g$ of the form 
$$
X''\stackrel{\pi}{\longrightarrow}X\times_YY''\stackrel{f\times_Y\id_{Y''}}{\longrightarrow}Y'',
$$
where $\pi$ and $Y''\rightarrow Y$ are admissible blow-ups.
The morphism $h$ is clearly separated ({\bf \ref{ch-formal}}.\ref{prop-sepmorformal2} (2)), whence the claim.

Finally, the implication (d) $\Rightarrow$ (a) follows from the following: by \ref{prop-rigidclosedmap} the diagonal morphism $\mathscr{X}\rightarrow\mathscr{X}\times_{\mathscr{Y}}\mathscr{X}$ is closed.
\end{proof}

\begin{prop}\label{prop-seppropmorrigidxxx11}
Let $\varphi\colon\mathscr{X}\rightarrow\mathscr{Y}$ be a separated morphism of locally universally adhesive rigid spaces\index{rigid space!universally adhesive rigid space@universally adhesive ---!locally universally adhesive rigid space@locally --- ---} {\rm (\ref{dfn-universallyadhesiverigidspaces})}, where $\mathscr{Y}=(\Spf B)^{\rig}$ is a Stein affinoid\index{affinoid!Stein affinoid@Stein ---} with $B$ a t.u.\ adhesive ring\index{t.u.a. ring@t.u.\ adhesive ring} {\rm ({\bf \ref{ch-formal}}.\ref{dfn-tuaringadmissible})}.
Then for any Stein affinoid open subspaces $\mathscr{U}=(\Spf P)^{\rig}$ and $\mathscr{V}=(\Spf Q)^{\rig}$ of $\mathscr{X}$ with $P$ and $Q$ t.u.\ adhesive, $\mathscr{U}\times_{\mathscr{X}}\mathscr{V}$ is a Stein affinoid open subspace of the form $(\Spf R)^{\rig}$ with $R$ t.u.\ adhesive.
\end{prop}

\begin{proof}
Consider the open subspace $\mathscr{U}\times_{\mathscr{Y}}\mathscr{V}$ of $\mathscr{X}\times_{\mathscr{Y}}\mathscr{X}$, which is a Stein affinoid, since we have $\mathscr{U}\times_{\mathscr{Y}}\mathscr{V}=(\Spf P'\times_{\Spf B}\Spf Q')^{\rig}$ where $P'$ and $Q'$ are as in \ref{prop-lemmorbetweenaffinoid2}.
Note that the schemes $\Spec B\setminus V(I)$ (where $I\subseteq B$ is an ideal of definition), $\Spec P'\setminus V(IP')$, and $\Spec Q'\setminus V(IQ')$ are affine.

First we show that $\mathscr{U}\times_{\mathscr{Y}}\mathscr{V}$ is a Stein affinoid.
It is an affinoid of the form $\Spf P'\widehat{\otimes}_BQ'$; the scheme $\Spec P'\widehat{\otimes}_BQ'\setminus V(I(P'\widehat{\otimes}_BQ'))$ is the pull-back of the affine scheme $\Spec P'\otimes_BQ'\setminus V(I(P'\otimes_BQ'))$ by the affine map $\Spec P'\widehat{\otimes}_BQ'\rightarrow\Spec P'\otimes_BQ'$ and hence is affine.
This shows that $\mathscr{U}\times_{\mathscr{Y}}\mathscr{V}$ is a Stein affinoid.

Now the space $\mathscr{U}\times_{\mathscr{X}}\mathscr{V}$ in question is the pull-back of $\mathscr{U}\times_{\mathscr{Y}}\mathscr{V}$ by the closed immersion $\mathscr{X}\rightarrow\mathscr{X}\times_{\mathscr{Y}}\mathscr{X}$ and hence is a closed subspace of the Stein affinoid $\mathscr{U}\times_{\mathscr{Y}}\mathscr{V}$.
Then the assertion follows from \ref{prop-closedimmrigid2xx}.
\end{proof}

\begin{prop}\label{prop-seppropgraphclosed}
Let $\varphi\colon\mathscr{X}\rightarrow\mathscr{Y}$ be a morphism of locally universally Noetherian rigid spaces\index{rigid space!universally adhesive rigid space@universally adhesive ---!locally universally adhesive rigid space@locally --- ---}, and $\psi\colon\mathscr{Y}\rightarrow\mathscr{Z}$ a separated and locally of finite type morphism between locally universally Noetherian rigid spaces.
Then the graph map
$$
\Gamma_{\varphi}\colon\mathscr{X}\longrightarrow\mathscr{X}\times_{\mathscr{Z}}\mathscr{Y}
$$
is a closed immersion\index{immersion!closed immersion of rigid spaces@closed --- (of rigid spaces)}.
\end{prop}

\begin{proof}
Since $\mathscr{Y}\rightarrow\mathscr{Z}$ is locally of finite type, the fiber product $\mathscr{X}\times_{\mathscr{Z}}\mathscr{Y}$ is locally universally Noetherian.
Then the proof goes similarly to {\bf \ref{ch-formal}}.\ref{prop-closedimmform111}.
\end{proof}

By an argument similar to that in {\bf \ref{ch-formal}}.\ref{prop-propertyP3} we have:
\begin{prop}\label{prop-propertyPrigid2}
Let $P$ be a property of morphisms in the category of locally universally Noetherian rigid spaces\index{rigid space!universally adhesive rigid space@universally adhesive ---!locally universally adhesive rigid space@locally --- ---} satisfying {\bf (I)}, {\bf (C)} in {\rm {\bf \ref{ch-pre}}, \S\ref{subsub-genpropertymorphism}} and the mutually equivalent conditions {\boldmath $(\mathbf{B}_i)$} for $i=1,2,3$ with $Q=$ `locally of finite type' in {\rm {\bf \ref{ch-pre}}, \S\ref{subsub-genpropertymorphismadic}}.
Suppose that any closed immersion satisfies $P$.
Then the following holds$:$ if $\varphi\colon\mathscr{X}\rightarrow\mathscr{Y}$ and $\psi\colon\mathscr{Y}\rightarrow\mathscr{Z}$ are morphisms of locally universally Noetherian rigid spaces such that $\psi\circ\varphi$ satisfies $P$ and $\psi$ is separated and locally of finite type, then $\varphi$ satisfies $P$. \hfill$\square$
\end{prop}

\begin{cor}\label{cor-propertyP3rigid}
Let $\varphi\colon\mathscr{X}\rightarrow\mathscr{Y}$ and $\psi\colon\mathscr{Y}\rightarrow\mathscr{Z}$ be morphisms of locally universally Noetherian rigid spaces\index{rigid space!universally adhesive rigid space@universally adhesive ---!locally universally adhesive rigid space@locally --- ---}, and suppose $\psi$ is separated and locally of finite type.
If $\psi\circ\varphi$ satisfies one of the following conditions, then so does $\varphi:$ 
\begin{itemize}
\item[{\rm (a)}] locally of finite type $($resp.\ of finite type$)$,
\item[{\rm (b)}] quasi-compact $($resp.\ quasi-separated, resp.\ coherent$)$,
\item[{\rm (c)}] finite,
\item[{\rm (d)}] closed immersion $($resp.\ immersion$)$,
\item[{\rm (e)}] closed $($universally closed$)$. \hfill$\square$
\end{itemize}
\end{cor}

\begin{prop}\label{prop-seppropmorrigid11x1}
{\rm (1)} The composition of two separated morphism of locally universally Noetherian rigid spaces\index{rigid space!universally adhesive rigid space@universally adhesive ---!locally universally adhesive rigid space@locally --- ---} is separated.

{\rm (2)} An open immersion of rigid spaces is separated. A closed immersion of rigid spaces is separated.

{\rm (3)} Let $\varphi\colon\mathscr{X}\rightarrow\mathscr{Y}$ be a morphism of rigid spaces, and $\{\mathscr{V}_{\alpha}\hookrightarrow\mathscr{Y}\}_{\alpha\in L}$ an admissible covering.
Then $\varphi$ is separated if and only if $\mathscr{X}\times_{\mathscr{Y}}\mathscr{V}_{\alpha}\rightarrow\mathscr{V}_{\alpha}$ is separated for each $\alpha\in L$.
\end{prop}

\begin{proof}
(1) Let $\varphi\colon\mathscr{X}\rightarrow\mathscr{Y}$ and $\psi\colon\mathscr{Y}\rightarrow\mathscr{Z}$ be separated morphisms of rigid spaces.
Since the diagonal $\mathscr{Y}\rightarrow\mathscr{Y}\times_{\mathscr{Z}}\mathscr{Y}$ is a closed immersion, so is the morphism $\mathscr{X}\times_{\mathscr{Y}}\mathscr{X}\rightarrow\mathscr{X}\times_{\mathscr{Z}}\mathscr{X}$ (as it is the base change of the former morphism by the canonciall morphism $\mathscr{X}\times_{\mathscr{Z}}\mathscr{X}\rightarrow\mathscr{Y}\times_{\mathscr{Z}}\mathscr{Y}$).
The diagonal $\mathscr{X}\rightarrow\mathscr{X}\times_{\mathscr{Z}}\mathscr{X}$ coincides with the composition of $\mathscr{X}\rightarrow\mathscr{X}\times_{\mathscr{Y}}\mathscr{X}$ followed by $\mathscr{X}\times_{\mathscr{Y}}\mathscr{X}\rightarrow\mathscr{X}\times_{\mathscr{Z}}\mathscr{X}$, whence the result.

(2) Clear.

(3) Set $\mathscr{U}_{\alpha}=\mathscr{X}\times_{\mathscr{Y}}\mathscr{V}_{\alpha}$ for each $\alpha\in L$.
Then the diagonal $\mathscr{U}_{\alpha}\rightarrow\mathscr{U}_{\alpha}\times_{\mathscr{V}_{\alpha}}\mathscr{U}_{\alpha}$ is the base change of $\mathscr{X}\rightarrow\mathscr{X}\times_{\mathscr{Y}}\mathscr{X}$ by the open immersion $\mathscr{U}_{\alpha}\times_{\mathscr{V}_{\alpha}}\mathscr{U}_{\alpha}\rightarrow\mathscr{X}\times_{\mathscr{Y}}\mathscr{X}$.
Then taking the associated maps between the Zariski-Riemann spaces, we deduce the desired equivalence of conditions by \ref{prop-rigidclosedmap}.
\end{proof}

\begin{cor}\label{cor-seppropmorrigid11}
All rigid spaces here are supposed to be locally universally Noetherian.

{\rm (1)} If $\varphi\colon\mathscr{X}\rightarrow\mathscr{X}'$ and $\psi\colon\mathscr{Y}\rightarrow\mathscr{Y}'$ are two separated morphisms over a rigid space $\mathscr{S}$, then the induced morphism $\varphi\times_{\mathscr{S}}\psi\colon\mathscr{X}\times_{\mathscr{S}}\mathscr{Y}\rightarrow\mathscr{X}'\times_{\mathscr{S}}\mathscr{Y}'$ is separated.

{\rm (2)} If $\varphi\colon\mathscr{X}\rightarrow\mathscr{Y}$ is a separated morphism over a rigid space $\mathscr{S}$ and $\mathscr{S}'\rightarrow\mathscr{S}$ is a morphism of rigid spaces, then the induced morphism $\varphi_{\mathscr{S}'}\colon\mathscr{X}\times_{\mathscr{S}}\mathscr{S}'\rightarrow\mathscr{Y}\times_{\mathscr{S}}\mathscr{S}'$ is separated.

{\rm (3)} If the composition $\psi\circ\varphi$ of two morphisms of rigid spaces is separated, then $\varphi$ is separated.
\end{cor}

\begin{proof}
By \ref{prop-seppropmorrigid11x1} (3), we may assume that the rigid spaces in question are all coherent.
Then the assertions follow easily from \ref{thm-seppropmorrigid1}, {\bf \ref{ch-formal}}.\ref{prop-sepmorformal2}, and {\bf \ref{ch-formal}}.\ref{prop-sepmorformal24}.
To show (3), we use the fact that an admissible blow-up is separated.
Note that the conditions (1) and (2) are equivalent due to {\bf \ref{ch-pre}}.\ref{prop-basechangestable}.
\end{proof}

\begin{thm}\label{thm-seppropmorrigid2}
Let $\varphi\colon\mathscr{X}\rightarrow\mathscr{Y}$ be a separated morphism of finite type of coherent rigid spaces.
Then the following conditions are equivalent$:$
\begin{itemize}
\item[{\rm (a)}] $\varphi$ is proper {\rm (\ref{dfn-seppropmorrigid1} (3));}
\item[{\rm (b)}] any distinguished formal model $f\colon X\rightarrow Y$ of $\varphi$ is proper$;$
\item[{\rm (c)}] there exists a proper formal model $f\colon X\rightarrow Y;$
\item[{\rm (d)}] there exists a cofinal set of formal models of $\varphi$ consisting of proper morphisms.
\end{itemize}
\end{thm}

\begin{proof}
First we show (a) $\Rightarrow$ (b).
Take a distinguished formal model $f\colon X\rightarrow Y$ of $\varphi$.
To see that $f$ is proper, it suffices by Exercise \ref{exer-propermorformalschemex1} to see that $f\times_Y\id_{\widehat{\A}^N_Y}\colon X\times_Y\widehat{\A}^N_Y\rightarrow\widehat{\A}^N_Y$ is closed for any $N\geq 0$.
But the morphism $f\times_Y\id_{\widehat{\A}^N_Y}$ is nothing but the formal model of $\varphi\times_{\mathscr{Y}}\id_{\D^N_{\mathscr{Y}}}$ (cf.\ \S\ref{subsub-unitdisk}), and hence the claim follows from \ref{prop-rigidclosedmap}.

The implication (b) $\Rightarrow$ (c) is obvious.
To show (c) $\Rightarrow$ (d), take a proper formal model $f\colon X\rightarrow Y$ of $\varphi$.
For any formal model $g\colon X'\rightarrow Y'$ of $\varphi$ there exists a model $h$ dominating $g$ of the form 
$$
X''\stackrel{\pi}{\longrightarrow}X\times_YY''\stackrel{f\times_Y\id_{Y''}}{\longrightarrow}Y'',
$$
where $\pi$ and $Y''\rightarrow Y$ are admissible blow-ups.
The morphism $h$ is clearly proper ({\bf \ref{ch-formal}}.\ref{prop-propermorformal2} (3)), whence the claim.

Finally, let us show (d) $\Rightarrow$ (a).
By \ref{prop-rigidclosedmap} the morphism $\varphi$ is closed.
It suffices then to see the condition (d) is preserved after base change.
But since (d) implies (c), and since (c) is preserved after base change, the claim follows.
\end{proof}

\begin{prop}\label{prop-seppropmorrigid11x2}
{\rm (1)} A finite morphism is proper.

{\rm (2)} The composition of two proper morphism between locally universally Noetherian rigid spaces\index{rigid space!universally adhesive rigid space@universally adhesive ---!locally universally adhesive rigid space@locally --- ---} is proper.

{\rm (3)} Let $\varphi\colon\mathscr{X}\rightarrow\mathscr{Y}$ be a morphism of rigid spaces, and $\{\mathscr{V}_{\alpha}\hookrightarrow\mathscr{Y}\}_{\alpha\in L}$ a covering in the small admissible site $\mathscr{Y}_{\ad}$.
Then $\varphi$ is proper if and only if $\mathscr{X}\times_{\mathscr{Y}}\mathscr{V}_{\alpha}\rightarrow\mathscr{V}_{\alpha}$ is proper for each $\alpha\in L$.
\end{prop}

\begin{proof}
(2) and (3) are easy to see.
We show (1).
Let $\varphi\colon\mathscr{X}\rightarrow\mathscr{Y}$ be a finite morphism of rigid spaces.
Since properness is a local condition due to (3), we may assume that $\mathscr{Y}$, and $\mathscr{X}$ also, is coherent.
Then it is proper by \ref{thm-seppropmorrigid2} and {\bf \ref{ch-formal}}.\ref{prop-propermorformal10}.
\end{proof}

\begin{cor}\label{cor-seppropmorrigid21}
{\rm (1)} If $\varphi\colon\mathscr{X}\rightarrow\mathscr{X}'$ and $\psi\colon\mathscr{Y}\rightarrow\mathscr{Y}'$ are two proper morphisms over a rigid space $\mathscr{S}$, then the induced morphism $\varphi\times_{\mathscr{S}}\psi\colon\mathscr{X}\times_{\mathscr{S}}\mathscr{Y}\rightarrow\mathscr{X}'\times_{\mathscr{S}}\mathscr{Y}'$ is proper.

{\rm (2)} If $\varphi\colon\mathscr{X}\rightarrow\mathscr{Y}$ is a proper morphism over a rigid space $\mathscr{S}$ and $\mathscr{S}'\rightarrow\mathscr{S}$ is a morphism of rigid spaces, then the induced morphism $\varphi_{\mathscr{S}'}\colon\mathscr{X}\times_{\mathscr{S}}\mathscr{S}'\rightarrow\mathscr{Y}\times_{\mathscr{S}}\mathscr{S}'$ is proper.

{\rm (3)} Suppose that the composition $\psi\circ\varphi$ of two morphisms between rigid spaces is proper.
If $\psi$ is separated, $\varphi$ is proper.
\end{cor}

\begin{proof}
We first claim the following: for a separated morphism $\mathscr{X}\rightarrow\mathscr{Y}$ where $\mathscr{Y}$ is quasi-separated, $\mathscr{X}$ is quasi-separated.
This indeed follows from \ref{prop-seppropmorrigidxx1}.
Hence, in particular, if it is proper and $\mathscr{Y}$ is coherent, then $\mathscr{X}$ is again coherent.
Combining this with \ref{prop-seppropmorrigid11x2} (3), we may assume that the rigid spaces in question are all coherent; then the claim follows easily from \ref{thm-seppropmorrigid2} and {\bf \ref{ch-formal}}.\ref{prop-propermorformal2}.
In view of {\bf \ref{ch-pre}}.\ref{prop-basechangestable} we also have (2).
Finally, (3) can be shown by an argument similar to that in the proof of \ref{cor-seppropmorrigid11} (3).
\end{proof}

\begin{cor}\label{cor-seppropmorrigid22}
Let $\varphi\colon\mathscr{X}\rightarrow\mathscr{Y}$ be a separated morphism of finite type of rigid spaces.
Then $\varphi$ is proper if and only if $\varphi\times_{\mathscr{Y}}\id_{\D^N_{\mathscr{Y}}}\colon\mathscr{X}\times_{\mathscr{Y}}\D^N_{\mathscr{Y}}\rightarrow\D^N_{\mathscr{Y}}$ is closed for any $N\geq 0$.
\end{cor}

\begin{proof}
One can assume that $\mathscr{Y}$ (and hence $\mathscr{X}$ also) is coherent.
Then the corollary follows from the argument in the proof of (a) $\Rightarrow$ (b) of \ref{thm-seppropmorrigid2}.
\end{proof}
\index{morphism of rigid spaces@morphism (of rigid spaces)!proper morphism of rigid spaces@proper ---|)}
\index{morphism of rigid spaces@morphism (of rigid spaces)!separated morphism of rigid spaces@separated ---|)}

\subsubsection{Valuative criterion}\label{subsub-seppropmorrigidcri}
\index{valuative!valuative criterion@--- criterion|(}
\begin{thm}\label{thm-valuativecriterionrigid}
Let $\varphi\colon\mathscr{X}\rightarrow\mathscr{Y}$ be a morphism of finite type between locally universally Noetherian rigid spaces\index{rigid space!universally adhesive rigid space@universally adhesive ---!locally universally adhesive rigid space@locally --- ---}.

{\rm (1)} The morphism $\varphi$ is separated if and only if the following condition holds$:$
let $\alpha\colon\mathscr{T}\rightarrow\mathscr{Y}$ be a rigid point\index{point!rigid point@rigid ---}\index{rigid point} $($cf.\ $\ref{dfn-ZRpoints2}\ (1))$ with a generization $\mathscr{T}'\rightarrow\mathscr{T}$, and suppose we are given a commutative diagram
$$
\xymatrix{
\mathscr{T}'\ar[r]^{\beta}\ar[d]&\mathscr{X}\ar[d]^{\varphi}\\
\mathscr{T}\ar[r]_{\alpha}&\mathscr{Y}\rlap{$;$}}
$$
then there exists at most one morphism $\mathscr{T}\rightarrow\mathscr{X}$ making the resulting diagram commute.

{\rm (2)} The morphism $\varphi$ is proper if and only if the following condition holds$:$
let $\alpha\colon\mathscr{T}\rightarrow\mathscr{Y}$ be a rigid point with a generization $\mathscr{T}'\rightarrow\mathscr{T}$, and suppose we are given a commutative diagram as above$;$
then there exists uniquely a morphism $\mathscr{T}\rightarrow\mathscr{X}$ making the resulting diagram commutative.
\end{thm}

Here by a generization of a rigid point $\alpha\colon\mathscr{T}=(\Spf V)^{\rig}\rightarrow\mathscr{Y}$ we mean the rigid point of the form $\mathscr{T}'=(\Spf V')^{\rig}\rightarrow\mathscr{T}\rightarrow\mathscr{Y}$ where $V'=\widehat{V}_{\mathfrak{p}}$ by a prime ideal $\mathfrak{p}\subseteq V$.

\begin{proof}
(1) Suppose that $\varphi$ is separated and that a rigid point $\alpha\colon\mathscr{T}=(\Spf V)^{\rig}\rightarrow\mathscr{Y}$ is given.
Consider the base change $\mathscr{X}\times_{\mathscr{Y}}\mathscr{T}\rightarrow\mathscr{T}$, which is again separated by \ref{cor-seppropmorrigid11} (2).
What to show is that, if we have two sections $\til{\beta}_1$ and $\til{\beta}_2$ having the common generization $\beta$ (as in the diagram below), we have $\til{\beta}_1=\til{\beta}_2$:
$$
\xymatrix@-1ex{
\mathscr{T}'\ar[r]^(.37){\beta}\ar[dr]&\mathscr{X}\times_{\mathscr{Y}}\mathscr{T}\ar[d]\\
&\mathscr{T}\rlap{.}\ar@<.5ex>@/_1pc/[u]\ar@<-.5ex>@/_1pc/[u]_{\til{\beta}_1,\til{\beta}_2}}
$$

To show this, we may assume that $\mathscr{Y}=\mathscr{T}$ and thus that $\mathscr{X}=\mathscr{X}\times_{\mathscr{Y}}\mathscr{T}$.
The two sections define the morphism $\mathscr{T}\rightarrow\mathscr{X}\times_{\mathscr{T}}\mathscr{X}$, which, restricted on $\mathscr{T}'$, factors through the diagonal $\mathscr{X}\rightarrow\mathscr{X}\times_{\mathscr{T}}\mathscr{X}$.
Since the diagonal morphism is a closed immersion, it gives rise to a closed immersion $\ZR{\mathscr{X}}\hookrightarrow\ZR{\mathscr{X}\times_{\mathscr{T}}\mathscr{X}}$.
Hence the map $\ZR{\mathscr{T}}\rightarrow\ZR{\mathscr{X}\times_{\mathscr{T}}\mathscr{X}}$ factors through $\ZR{\mathscr{X}}$.
By \ref{prop-ZRpointx1} we get a rigid point $\mathscr{T}\rightarrow\mathscr{X}$, which coincides with the section $\til{\beta}_1$ and $\til{\beta}_2$.
Hence we have $\til{\beta}_1=\til{\beta}_2$, as desired.

Suppose, conversely, that the condition in (1) is satisfied.
We want to show that the diagonal morphism $\Delta_{\mathscr{X}}\colon\mathscr{X}\rightarrow\mathscr{X}\times_{\mathscr{Y}}\mathscr{X}$ is closed.
By \ref{prop-rigidclosedmap} this is equivalent to that we have the equality
$$
\ZR{\Delta_{\mathscr{X}}}(\ovl{\{x\}})=\ovl{\{\ZR{\Delta_{\mathscr{X}}}(x)\}}
$$
for any $x\in\ZR{\mathscr{X}}$.
As the left-hand side is clearly contained in the right-hand side, we need to show the following: for any $x\in\ZR{\mathscr{X}}$ and a specialization $y'$ of $y=\ZR{\Delta_{\mathscr{X}}}(x)$ there exists a specialization $x'$ of $x$ that is mapped to $y'$.

Let $W=\widehat{V}_{x}$, and $\beta\colon\mathscr{T}'=(\Spf\widehat{V}_{x})^{\rig}\rightarrow\mathscr{X}$ the rigid point associated to $x$ (\ref{dfn-ZRpoints32}).
Consider the inclusion of valuation rings $\widehat{V}_{y'}\hookrightarrow\widehat{V}_y\hookrightarrow W$; the last injectivity is due to the fact that it is an adic homomorphism ({\bf \ref{ch-pre}}.\ref{prop-maxspe4}).
Let $V$ be the integral closure of $\widehat{V}_{y'}$ in $W$.
Then $V$ is again an $a$-adically complete valuation ring, and $W$ is a localization of $V$.
We have a rigid point $\til{\beta}\colon\mathscr{T}=(\Spf V)^{\rig}\rightarrow\mathscr{X}\times_{\mathscr{Y}}\mathscr{X}$.
The rigid point $\til{\beta}$ gives rise to two rigid points $\til{\beta}_1,\til{\beta}_2\colon\mathscr{T}\rightarrow\mathscr{X}$ by projections.
Let $\alpha\colon\mathscr{T}\rightarrow\mathscr{Y}$ be the composition $\varphi\circ\til{\beta}_1=\varphi\circ\til{\beta}_2$.
Then we have the commutative diagram 
$$
\xymatrix{
\mathscr{T}'\ar[r]^{\beta}\ar[d]&\mathscr{X}\ar[d]^{\varphi}\\
\mathscr{T}\ar[r]_{\alpha}&\mathscr{Y}\rlap{.}}
$$
together with two $\mathscr{T}$-sections $\til{\beta}_1$ and $\til{\beta}_2$.
Since $\mathscr{T}'$ is a generization of $\mathscr{T}$, we have $\til{\beta}_1=\til{\beta}_2$.
Let $x'$ be the image of the closed point by the map $\ZR{\til{\beta}_1}=\ZR{\til{\beta}_2}\colon\Spf V\rightarrow\ZR{\mathscr{X}}$.
Then by the construction the point $x'$ is a specialization of $x$ and is mapped to $y'$ by the map $\ZR{\Delta_{\mathscr{X}}}$.

(2) We first show the ``only if'' part.
Suppose that $\varphi$ is proper and that the diagram as above is given.
Since the properness is stable under base change (\ref{cor-seppropmorrigid21} (2)), we may assume $\mathscr{T}=\mathscr{Y}$.
Then we are given a diagram 
$$
\xymatrix@-1ex{
\mathscr{T}'\ar[r]^{\beta}\ar[dr]&\mathscr{X}\ar[d]_{\varphi}\\
&\mathscr{T}\rlap{,}}
$$
where $\beta$ is a $\mathscr{T}'$-section of $\varphi$, and what to prove is that there exists a section $\til{\beta}$ of $\varphi$ that extends $\beta$ (the uniqueness follows from separatedness and what we have already shown above).
Set $\mathscr{T}=(\Spf V)^{\rig}$ and $\mathscr{T}'=(\Spf W)^{\rig}$, where $W$ is a localization of $V$.
Let $x$ be the image of the closed point by the map $\ZR{\beta}\colon\Spf W\rightarrow\ZR{\mathscr{X}}$.
Since $\varphi$ is closed, we have $\ZR{\varphi}(\ovl{\{x\}})=\ovl{\{\ZR{\varphi}(x)\}}=\Spf V$.
Hence there exists a specialization $x'$ of $x$ that is mapped to the closed point of $\Spf V$ by $\ZR{\varphi}$.
By \ref{prop-ZRpoints52} and \ref{prop-ZRpoint32} (2) the map $\ZR{\beta}$ factors through the associated rigid point $\Spf\widehat{V}_{x'}\rightarrow\ZR{\mathscr{X}}$.
Hence we have a chain of maps
$$
\Spf W\longrightarrow\Spf\widehat{V}_{x'}\longrightarrow\Spf V;
$$
accordingly, we have 
$$
V\longrightarrow\widehat{V}_{x'}\longrightarrow W.
$$
All these maps are injective, since they are adic homomorphisms ({\bf \ref{ch-pre}}.\ref{prop-maxspe4}).
Moreover, the first map is local; that is, $\widehat{V}_{x'}$ dominates $V$.
Since $W$ is a localization of $V$, we have $\Frac(V)=\Frac(\widehat{V}_{x'})$.
Hence we have $V=\widehat{V}_{x'}$.
Thus we get the rigid point $\Spf V\rightarrow\ZR{\mathscr{X}}$ and a section $\mathscr{T}\rightarrow\mathscr{X}$, as desired.

Suppose, conversely, that the condition in (2) is satisfied.
Due to the uniqueness, $\varphi$ is separated; since the condition in question is stable under base change, it suffices to show that the morphism $\varphi$ is closed.
By \ref{prop-rigidclosedmap} this is equivalent to that for any $x\in\ZR{\mathscr{X}}$ we have the equality
$$
\ZR{\varphi}(\ovl{\{x\}})=\ovl{\{\ZR{\varphi}(x)\}}.
$$
As the left-hand side is clearly contained in the right-hand side, we need to show the following: for any $x\in\ZR{\mathscr{X}}$ and a specialization $y'$ of $y=\ZR{\varphi}(x)$ there exists a specialization $x'$ of $x$ that is mapped to $y'$.

Set $W=\widehat{V}_x$, and let $\beta\colon\mathscr{T}=(\Spf W)^{\rig}\rightarrow\mathscr{X}$ be the associated rigid point.
Similarly to the proof of (1) as above, we have inclusions of valuation rings $\widehat{V}_{y'}\hookrightarrow\widehat{V}_y\hookrightarrow W$.
Let $V$ be the integral closure of $\widehat{V}_{y'}$ in $W$.
Then we have a rigid point $\alpha\colon\mathscr{T}=(\Spf V)^{\rig}\rightarrow\mathscr{Y}$.
We have, by the hypothesis, a $\mathscr{T}$-section $\mathscr{T}\rightarrow\mathscr{X}$ of $\varphi$.
Let $x'$ be the image of the closed point by the corresponding map $\Spf V\rightarrow\ZR{\mathscr{X}}$.
Then $x'$ is a specialization of $x$ that is mapped to $y'$.
\end{proof}
\index{valuative!valuative criterion@--- criterion|)}

\subsubsection{Finiteness theorem}\label{subsub-finitudesrigid}
Let us include here a finiteness theorem of cohomologies of coherent sheaves for proper morphisms of universally adhesive rigid spaces.
First we note the following `boundedness' result for cohomologies of separated (\ref{dfn-seppropmorrigid1} (2)) and quasi-compact (\ref{dfn-qcptmorphismrigid1}) morphisms of rigid spaces:
\begin{prop}\label{prop-finitudesrigid1}
Let $\varphi\colon\mathscr{X}\rightarrow\mathscr{Y}$ be a separated and quasi-compact morphism between universally adhesive rigid spaces\index{rigid space!universally adhesive rigid space@universally adhesive ---} {\rm (\ref{dfn-universallyadhesiverigidspaces})}.
Then there exists an integer $r>0$ such that for any coherent $\O_{\mathscr{X}}$-module $\mathscr{F}$ and any $q\geq r$, we have $\RD^q\varphi_{\ast}\mathscr{F}=0$.
If, moreover, $\mathscr{Y}$ is a Stein affinoid, then one can take as $r$ the minimum number of Stein affinoids in a Stein affinoid covering of $\mathscr{X}$.
\end{prop}

\begin{proof}
This follows immediately from that fact that by \ref{prop-seppropmorrigidxxx11} and \ref{thm-cohomologyrigidsp1str} a Stein affinoid covering is a Leray covering for coherent sheaves.
\end{proof}

For a rigid space $\mathscr{X}$ we denote by $\DC^{\ast}(\mathscr{X})$ ($\ast=$``\ \ '', $+$, $-$, $\bd$) the derived category of the abelian category of $\O_{\mathscr{X}}$-modules.
The full subcategory of $\DC^{\ast}(\mathscr{X})$ consisting of objects with coherent cohomologies is denoted by $\DC^{\ast}_{\coh}(\mathscr{X})$.
The next theorem gives the finiteness of coherent cohomologies of proper (\ref{dfn-seppropmorrigid1} (3)) morphisms between quasi-compact universally adhesive rigid spaces:
\begin{thm}[Finiteness theorem for proper morphisms]\label{thm-finirigid}
Let $\varphi\colon\mathscr{X}\rightarrow\mathscr{Y}$ be a proper morphism between universally adhesive rigid spaces\index{rigid space!universally adhesive rigid space@universally adhesive ---} {\rm (\ref{dfn-universallyadhesiverigidspaces})}.
Then the functor $\RD\varphi_{\ast}$ maps $\DC^{\ast}_{\coh}(\mathscr{X})$ to $\DC^{\ast}_{\coh}(\mathscr{Y})$ for $\ast=$``\ \ '', $+$, $-$, $\bd$.
\end{thm}

To show the theorem, we need:
\begin{lem}\label{lem-exerpropprojlimcohsch52}
Let $\{X_i,p_{ij}\}_{i\in I}$ and $\{Y_i,q_{ij}\}_{i\in I}$ be projective system of ringed spaces indexed by a common directed set $I$, and $\{f_i\colon X_i\rightarrow Y_i\}_{i\in I}$ a morphism of these systems $($that is, $q_{ij}\circ f_j=f_i\circ p_{ij}$ for any $i\geq j)$ such that the following conditions are satisfied$:$
\begin{itemize}
\item for any $i\in I$ the underlying topological spaces of $X_i$ and $Y_i$ are coherent\index{space@space (topological)!coherent topological space@coherent ---}\index{coherent!coherent topological space@--- (topological) space} {\rm ({\bf \ref{ch-pre}}.\ref{dfn-quasicompact1})} and sober\index{space@space (topological)!sober topological space@sober ---} {\rm ({\bf \ref{ch-pre}}, \S\ref{subsub-sober})}$;$
\item for any $i\leq j$ the underlying continuous mapping of the transition maps $p_{ij}\colon X_j\rightarrow X_i$ and $q_{ij}\colon Y_j\rightarrow Y_i$ are quasi-compact {\rm ({\bf \ref{ch-pre}}.\ref{dfn-quasicompactness} (2))}$;$
\item for any $i\in I$ the underlying continuous mapping of $f_i\colon X_i\rightarrow Y_i$ is quasi-compact.
\end{itemize}
Let $0\in I$ be an index, and $\mathscr{F}_0$ an $\O_{X_0}$-module.
Then the canonical morphism 
$$
\varinjlim_{i\geq 0}q^{\ast}_i\RD^qf_{i\ast}p^{\ast}_{0i}\mathscr{F}_0\longrightarrow\RD^qf_{\ast}p^{\ast}_0\mathscr{F}_0
$$
is an isomorphism for any $q\geq 0$, where $q_i$ for $i\in I$ is the canonical projection $Y=\varprojlim_{j\in I}Y_j\rightarrow Y$.
\end{lem}

\begin{proof}
Since $\varinjlim_{i\in I}p^{-1}_ip^{\ast}_{0i}\mathscr{F}_0\cong p^{\ast}_0\mathscr{F}_0$ by {\bf \ref{ch-pre}}.\ref{lem-projlimcohsheaf1}, the desired result in $q=0$ follows from {\bf \ref{ch-pre}}.\ref{cor-coherentprojectivesystemindlimits3} and {\bf \ref{ch-pre}}.\ref{lem-corprojlimcohsheaf1}.
The general assertion can be shown by an argument similar to that in {\bf \ref{ch-pre}}, \S\ref{subsub-cohsheavestopologicalspacelimits}.
\end{proof}

\begin{proof}[Proof of Theorem {\rm \ref{thm-finirigid}}]
By a standard reduction process similarly to that in {\bf \ref{ch-formal}}, \S\ref{subsub-finitudesproof}, it suffices to show that for any coherent $\O_{\mathscr{X}}$-module $\mathscr{F}$ the sheaf $\RD^q\varphi_{\ast}\mathscr{F}$ is a coherent $\O_{\mathscr{Y}}$-module.
For this we may assume that $\mathscr{Y}$ is coherent (and then so is $\mathscr{X}$).
By \ref{thm-tateacyclic20} there exists a distinguished formal model $f\colon X\rightarrow Y$ of $\varphi$ and a coherent $\O_X$-module $\mathscr{F}_X$ such that $\mathscr{F}^{\rig}_X=\mathscr{F}$.
By \ref{thm-seppropmorrigid2} the morphism $f$ is proper.
Now applying {\bf \ref{ch-formal}}.\ref{thm-finiformal}, we know that $\RD^qf_{\ast}\mathscr{F}_X$ is a coherent $\O_Y$-module.
Since any distinguished model of $\varphi$ is proper and distinguished models are cofinal among all the formal models, we get the assertion by \ref{lem-exerpropprojlimcohsch52} and Deligne's formula (Exercise \ref{exer-deligne}).
\end{proof}

\subsection{Projective morphisms}\label{sub-projectiverigid}
\index{morphism of rigid spaces@morphism (of rigid spaces)!projective morphism of rigid spaces@projective ---|(}
Let $\mathscr{Y}$ be a coherent universally Noetherian rigid space\index{rigid space!universally Noetherian rigid space@universally Noetherian ---!coherent universally Noetherian rigid space@coherent --- ---}, and $\varphi\colon\mathscr{X}\rightarrow\mathscr{Y}$ a proper morphism (hence $\mathscr{X}$ is coherent).
Let $\mathscr{L}$ be an invertible $\O_{\mathscr{X}}$-module.
One can take by \ref{thm-corcohsheavesZR2101} an invertible $\O^{\int}_{\mathscr{X}}$-module $\mathscr{L}^+$ that gives an integral model\index{integral model} of $\mathscr{L}$.

\begin{dfn}\label{dfn-positivesheaf}{\rm 
The sheaf $\mathscr{L}^+$ is said to be {\em $\varphi$-positive} or {\em positive relative to $\varphi$} if there exist a formal model $f\colon X\rightarrow Y$ of $\varphi$ and an invertible sheaf $\mathscr{L}_X$ such that $\sp^{\ast}_X\mathscr{L}_X=\mathscr{L}^+$ and that $\mathscr{L}_X$ is ample, that is, $\mathscr{L}_0=\mathscr{L}/\mathscr{I}_X\mathscr{L}$ is $f_0$-ample, where $\mathscr{I}_X$ is an ideal of definition of finite type of $X$ and $f_0\colon X_0\rightarrow Y_0$ is the induced morphism of schemes obtained by mod $\mathscr{I}_X$-reduction.}
\end{dfn}

\begin{dfn}\label{dfn-projectiverigid1}{\rm 
A proper morphism $\varphi\colon\mathscr{X}\rightarrow\mathscr{Y}$ of coherent universally Noetherian rigid spaces\index{rigid space!universally Noetherian rigid space@universally Noetherian ---!coherent universally Noetherian rigid space@coherent --- ---} is said to be {\em projective} if there exists a pair $(\mathscr{L},\mathscr{L}^+)$ consisting of an invertible $\O_{\mathscr{X}}$-module $\mathscr{L}$ and an invertible integral model of it such that $\mathscr{L}^+$ is $\varphi$-positive.}
\end{dfn}

\begin{prop}\label{prop-projectiverigid1}
Let $\varphi\colon\mathscr{X}\rightarrow\mathscr{Y}$ be a morphism of coherent universally Noetherian rigid spaces\index{rigid space!universally Noetherian rigid space@universally Noetherian ---!coherent universally Noetherian rigid space@coherent --- ---}, where $\mathscr{Y}=(\Spf B)^{\rig}$ is an affinoid with $B$ t.u.\ rigid-Noetherian\index{t.u. rigid-Noetherian ring@t.u.\ rigid-Noetherian ring}.
Then $\varphi$ is projective if and only if there exists a finitely generated $B$-module $E$ and a $\mathscr{Y}$-closed immersion $\mathscr{X}\hookrightarrow\widehat{\P(E)}^{\rig}$.
\end{prop}

\begin{proof}
The `if' part is easy to see.
To show the `only if' part, we may assume that $B$ is $I$-torsion free, where $I$ is an ideal of definition of $B$.
Let $(\mathscr{L},\mathscr{L}^+)$ be consisting of an invertible $\O_{\mathscr{X}}$-module and an invertible integral model that is $\varphi$-positive.
Take a formal model $f\colon X\rightarrow Y'$ of $\varphi$ and an $f$-ample $\mathscr{L}_X$ such that $\sp^{\ast}_X\mathscr{L}_X=\mathscr{L}^+$.
By \cite[$\mathbf{II}$, (4.6.13) (iii)]{EGA} we may assume that between $Y'$ and $Y=\Spf B$ is an admissible blow-up $\pi\colon Y'\rightarrow Y$.
We may also assume that $\O_X$ is $I$-torsion free.
Let $\mathscr{J}$ be the admissible ideal of $Y$ such that $\pi$ is the admissible blow-up along $\mathscr{J}$.
Then $\mathscr{J}\O_{Y'}$ is $\pi$-ample, and hence $\mathscr{L}_X\otimes_{\O_Y}f^{\ast}(\mathscr{J}\O_{Y'})^{\otimes n}$ for some $n>0$ is $\pi\circ f$-ample (\cite[$\mathbf{II}$, (4.6.13) (ii)]{EGA}).
By an argument similar to that in the proof of \cite[$\mathbf{III}$, (5.4.3)]{EGA} applied to the proper map $\pi\circ f$, we get a $Y$-closed immersion $X\hookrightarrow\widehat{\P(E)}$ for some $B$-module $E$ of finite type, whence the desired result.
\end{proof}

\begin{prop}\label{prop-projectiverigid2}
{\rm (1)} A closed immersion is projective.

{\rm (2)} If $\varphi\colon\mathscr{Z}\rightarrow\mathscr{Y}$ and $\psi\colon\mathscr{Y}\rightarrow\mathscr{X}$ are projective, then so is the composition $\psi\circ\varphi$.
If $\psi\circ\varphi$ is projective and $\psi$ is separated, then $\varphi$ is projective.

{\rm (3)} If $\varphi\colon\mathscr{X}\rightarrow\mathscr{X}'$ and $\psi\colon\mathscr{Y}\rightarrow\mathscr{Y}'$ are two projective morphisms over a rigid space $\mathscr{S}$ such that either $\mathscr{X}$ and $\mathscr{X}'$ or $\mathscr{X}$ and $\mathscr{Y}'$ are locally of finite type over $\mathscr{S}$, then the induced morphism $\varphi\times_{\mathscr{S}}\psi\colon\mathscr{X}\times_{\mathscr{S}}\mathscr{Y}\rightarrow\mathscr{X}'\times_{\mathscr{S}}\mathscr{Y}'$ is projective.

{\rm (4)} If $\varphi\colon\mathscr{X}\rightarrow\mathscr{Y}$ is a projective morphism over a rigid space $\mathscr{S}$ and $\mathscr{S}'\rightarrow\mathscr{S}$ is a morphism of rigid spaces such that either $\mathscr{X}$ and $\mathscr{Y}$ are locally of finite type over $\mathscr{S}$ or that $\mathscr{S}'$ is locally of finite type over $\mathscr{S}$, then the induced morphism $\varphi_{\mathscr{S}'}\colon\mathscr{X}\times_{\mathscr{S}}\mathscr{S}'\rightarrow\mathscr{Y}\times_{\mathscr{S}}\mathscr{S}'$ is projective.
\end{prop}

\begin{proof}
(1) is clear.
To show (2), let $(\mathscr{L},\mathscr{L}^+)$ (resp.\ $(\mathscr{M},\mathscr{M}^+)$) be the pair as in \ref{dfn-projectiverigid1} for the morphism $\varphi$ (resp.\ $\psi$).
Let $f\colon X\rightarrow Y$ (resp.\ $g\colon Y'\rightarrow Z$) be a formal model of $\varphi$ (resp.\ $\psi$) with the relatively ample sheaf $\mathscr{L}_X$ (resp.\ $\mathscr{M}_{Y'}$) as in \ref{dfn-positivesheaf}.
By \cite[$\mathbf{II}$, (4.6.13) (iii)]{EGA} we may assume that $Y$ is an admissible blow-up of $Y'$.
Write the admissible blow-up as $\pi\colon Y\rightarrow Y'$, and let $\mathscr{J}$ be the blow-up center.
Then by \cite[$\mathbf{II}$, (4.6.13) (ii)]{EGA} the sheaf $\mathscr{L}_X\otimes f^{\ast}(\mathscr{J}\O_{Y'}\otimes\pi^{\ast}\mathscr{M}_{Y'}^{\otimes m})^{\otimes n}$ for sufficiently large $n$ and $m$ is ample relative to the composite map $g\circ\pi\circ f$.
Let this sheaf be $\mathscr{N}_X$.
Then the pair $((\sp^{\ast}_X\mathscr{N}_X)^{\rig},\sp^{\ast}_X\mathscr{N}_X)$ guarantees that the composition $\psi\circ\varphi$ is projective.
Thus the first half of (2) is proved.
The assertions (3) and (4) are proved by similar arguments with the aid of \cite[$\mathbf{II}$, (4.6.13) (iii) (iv)]{EGA}.
Finally, in view of \ref{prop-propertyPrigid2}, the second half of (2) follows.
\end{proof}
\index{morphism of rigid spaces@morphism (of rigid spaces)!projective morphism of rigid spaces@projective ---|)}

\addcontentsline{toc}{subsection}{Exercises}
\subsection*{Exercises}
\begin{exer}\label{exer-associatedschfibreprod}{\rm 
Let $\mathscr{X}\rightarrow\mathscr{Z}\hookleftarrow\mathscr{Y}$ be a diagram consisting of universally Noetherian affinoids, where the second morphism is a closed immersion.
Show that the canonical morphism of the associated schemes\index{scheme!associated scheme to an affinoid@associated --- to an affinoid}\index{affinoid!associated scheme to an affinoid@associated scheme to an ---} (cf.\ \ref{subsub-associatedscheme})
$$
s(\mathscr{X}\times_{\mathscr{Z}}\mathscr{Y})\longrightarrow s(\mathscr{X})\times_{s(\mathscr{Z})}s(\mathscr{Y})
$$
is an isomorphism.}
\end{exer}

\begin{exer}\label{exer-propvalutiverigid1}
{\rm Let $\mathscr{Y}$ be a rigid space, and $\alpha\colon\mathscr{T}\rightarrow\mathscr{Y}$ a rigid point.
Let us say that the rigid point $\alpha$ is {\em essentially of finite type}\index{rigid point!rigid point essentially of finite type@--- essentially of finite type} if $\mathscr{T}$ is $\mathscr{Y}$-isomorphic to the associated rigid point (\ref{dfn-ZRpoints32}) to a point $x$ of the disk $\D^N_{\mathscr{Y}}$ for some $N\geq 0$.
Show that in \ref{thm-valuativecriterionrigid} (2) it is enough to use rigid point of $\mathscr{Y}$ essentially of finite type.}
\end{exer}


\section{Classical points}\label{sec-localring}
The main aim of this section is to discuss the so-called {\em classical points}.
Classical points are, roughly speaking, the points of rigid spaces of type ($\mathrm{V_{\R}}$) (\ref{dfn-typeV})\index{rigid space!rigid space of typeV1@--- of type ($\mathrm{V_{\R}}$)} or of type (N) (\ref{dfn-typeN})\index{rigid space!rigid space of typeN@--- of type (N)} that are, so to speak, already dealt with in classical rigid geometry.
In \S\ref{sub-pointspoints} we first set up a useful concept, the {\em spectral functors}, which provides a general framework for dealing with several kinds of points on rigid spaces.
Spectral functors are useful not only in this section, but in later sections in the appendix when we discuss the classical rigid geometry (\S\ref{sec-berkovich}) and Berkovich analytic geometry (\S\ref{sub-vsbr-Berkovichanalyticspaces}).

In \S\ref{sub-classicalpoints} we define and discuss classical points.
The several basic facts discussed in this subsection will be at the basis of our later study of the comparison of our rigid geometry with the classical rigid geometry and with Berkovich analytic geometry.

In the final subsection \S\ref{sub-locringclassical}, we prove the Noetherness theorem (\ref{thm-noetherness}), which asserts that, if $\mathscr{X}$ is a rigid space of type (V)\index{rigid space!rigid space of typeV@--- of type (V)} or of type (N), then the local ring $B_x=\O_{\mathscr{X},x}$ at any point $x\in\ZR{\mathscr{X}}$ is a Noetherian ring.

\subsection{Spectral functor}\label{sub-pointspoints}
\index{functor!spectral functor@spectral ---|(}\index{spectral functor|(}
\subsubsection{Definitions}\label{subsub-pointspoints}
\begin{dfn}\label{dfn-spectral1}{\rm 
Let $\mathscr{C}$ be a subcategory of the category of rigid spaces $\Rf$. 
\begin{itemize}
\item[(1)] We say that $\mathscr{C}$ is {\em $\mathrm{O}$-stable}\index{stable!Ostable@O-{---}} if for any $\mathscr{X}\in\obj(\mathscr{C})$ any open immersion\index{immersion!open immersion of rigid spaces@open --- (of rigid spaces)} $\mathscr{U}\hookrightarrow\mathscr{X}$ belongs to $\mathscr{C}$. 
\item[(2)] We say that $\mathscr{C}$ is {\em $\mathrm{QCO}$-stable}\index{stable!QCOstable@QCO-{---}} if for any $\mathscr{X}\in\obj(\mathscr{C})$ any quasi-compact open immersion $\mathscr{U}\hookrightarrow\mathscr{X}$ belongs to $\mathscr{C}$. 
\end{itemize}}
\end{dfn}

\begin{dfn}\label{dfn-spectral2}{\rm 
Let $\mathscr{C}$ be a subcategory of $\Rf$, and
$$
S\colon\mathscr{C}\longrightarrow\Top
$$ 
a functor. 
Consider the following conditions:
\begin{itemize}
\item[(a)] there exists a natural transformation from $S$ to $\ZR{\,\cdot\,}$ that induces for any $\mathscr{X}\in\obj(\mathscr{C})$ an inclusion of sets $S (\mathscr{X})\subset\ZR{\mathscr{X}}$; 
\item[(b)] $S(\mathscr{X})=\emptyset$ if and only if $\mathscr{X}=\emptyset$; 
\item[(c)] for any quasi-compact open immersion $\mathscr{U}\hookrightarrow\mathscr{V}$ in $\mathscr{C}$ the induced map of topological spaces $S(\mathscr{U})\rightarrow S(\mathscr{V})$ maps $S(\mathscr{U})$ homeomorphically onto $S(\mathscr{V})\cap\ZR{\mathscr{U}}$ with the subspace topology induced from the topology on $S(\mathscr{V})$;
\item[(d)] for a quasi-compact open immersion $\mathscr{U}\hookrightarrow\mathscr{V}$ in $\mathscr{C}$ the equality $S(\mathscr{U})= S(\mathscr{V})$ holds if and only if $\mathscr{U}=\mathscr{V}$. 
\end{itemize}
The functor $S$ is called a {\em prespectral functor}\index{spectral functor!prespectral functor@pre{---}} if it satisfies the conditions (a), (b), and (c).
If it furthermore satisfies (d), then it is called a {\em spectral functor}.}
\end{dfn}

Note that we do not assume in (a) that the topology on $S(\mathscr{X})$ is the induced topology from $\ZR{\mathscr{X}}$. 
Note also that, if $S$ is prespectral (resp.\ spectral), then  $S|_{{\mathscr{C}}'} $ is prespectral (resp.\ spectral) for any subcategory ${\mathscr C}'$ of $\mathscr{C}$. 
\begin{prop}\label{prop-spectral4}
Let $\mathscr{C}$ be a $\mathrm{QCO}$-stable subcategory of $\Rf$, and $S\colon\mathscr{C}\rightarrow\Top$ a spectral functor. 
Then for any $\mathscr{X}$ in $\mathscr{C}$ and any quasi-compact open immersions $\mathscr{U}\hookrightarrow\mathscr{X}$ and $\mathscr{V}\hookrightarrow\mathscr{X}$, $S(\mathscr{U})=S(\mathscr{V})$ implies $\mathscr{U}=\mathscr{V}$.
\end{prop}

\begin{proof}
The open immersion $\mathscr{U}\cap\mathscr{V}\hookrightarrow\mathscr{U}$ is quasi-compact and hence belongs to $\mathscr{C}$ due to the $\mathrm{QCO}$-stability. 
By the condition \ref{dfn-spectral2} (c) we have $ S(\mathscr{U})=S(\mathscr{X})\cap\ZR{\mathscr{U}}$ and $S(\mathscr{V})=S(\mathscr{X})\cap\ZR{\mathscr{V}}$. 
We have $S(\mathscr{U})=S(\mathscr{U})\cap S(\mathscr{V})= S(\mathscr{X})\cap\ZR{\mathscr{U}}\cap\ZR{\mathscr{V}}=S(\mathscr{X})\cap\ZR{\mathscr{U}\cap\mathscr{V}}$ by our assumption.
Since the composition $\mathscr{U}\cap\mathscr{V}\hookrightarrow\mathscr{U}\hookrightarrow\mathscr{X}$ is a quasi-compact immersion, again by \ref{dfn-spectral2} (c) we have $S(\mathscr{U}\cap\mathscr{V})=S(\mathscr{X})\cap\ZR{\mathscr{U}\cap\mathscr{V}}$ and hence $S(\mathscr{U}\cap\mathscr{V})=S(\mathscr{U})$.  
Then $\mathscr{U}\cap\mathscr{V}=\mathscr{U}$ follows from the condition \ref{dfn-spectral2} (d). 
By exchanging the role of $\mathscr{U}$ and that of $\mathscr{V}$, we have $\mathscr{V}=\mathscr{U}\cap\mathscr{V}$.
Thus we get $\mathscr{U}=\mathscr{V}$, as desired. 
\end{proof}

\subsubsection{Continuity}\label{subsub-pointscontinuity}
Let $S\colon\mathscr{C}\rightarrow\Top$ be a prespectral functor, and consider an inductive system\index{system!inductive system@inductive ---} $\{\mathscr{U}_i\}_{i\in I}$ in $\mathscr{C}$ indexed by a directed set\index{set!directed set@directed ---}\index{directed set} $I$ consisting of quasi-compact open immersions $\mathscr{U}_i\hookrightarrow\mathscr{U}_j$ for $i\leq j$.
Then by \ref{dfn-spectral2} (c) one can define $S(\mathscr{U})$, where $\mathscr{U}=\varinjlim_{i\in I}\mathscr{U}_i$, by the formula
$$
S(\mathscr{U})=\varinjlim_{i\in I}S(\mathscr{U}_i).
$$
In other words, any prespectral functor behaves {\em continuously} under filtered limits by quasi-compact open immersions.

Let us discuss more general continuity properties of prespectral functors:
\begin{dfn}\label{dfn-spectral3}{\rm 
Let $\mathscr{C}$ be a subcategory of $\Rf$.  
A prespectral functor $S$ is said to be {\em continuous}\index{spectral functor!continuous spectral functor@continuous ---} if for any small category $\mathscr{D}$ and a functor $F\colon\mathscr{D}\rightarrow\mathscr{C}$ such that $F(f)$ is an open immersion for any morphism $f$ in $\mathscr{D}$, we have
$$
S(\varinjlim_{\mathscr{D}}F)=\varinjlim_{\mathscr{D}}S\circ F 
$$
in the category of topological spaces.}
\end{dfn}

The following proposition is easy to verify; we leave the checking to the reader:
\begin{prop}\label{prop-contnuityspectral}
Let $\mathscr{C}$ be an $\mathrm{O}$-stable\index{stable!Ostable@O-{---}} subcategory of $\Rf$. 
Suppose $\mathscr{C}$ is stable under disjoint sum.  
Let $S\colon\mathscr{C}\rightarrow\Top$ be a continuous prespectral functor.

{\rm (1)} The functor $S$ commutes with disjoint sum\index{disjoint sum}, that is, for $\mathscr{X}_{\alpha}\in\mathscr{C}$ indexed by a set $L$, we have
$$
S(\amalg_{\alpha\in L}\mathscr{X}_{\alpha})=\amalg_{\alpha\in L}S(\mathscr{X}_{\alpha}).
$$

{\rm (2)} The functor $S$ preserves the equivalence relation$;$ that is, for $\mathscr{X}\in\obj(\mathscr{C})$ and an open covering $\{\mathscr{U}_{\alpha}\}_{\alpha\in L}$ of $\mathscr X$, let $\mathscr{R}$ be the equivalence relation\index{equivalence relation} defining $\mathscr{X}$, that is, 
$$
\xymatrix@-1ex{\mathscr{R}\ar@<.5ex>[r]\ar@<-.5ex>[r]&\coprod_{\alpha\in L}\mathscr{U}_{\alpha}\ar[r]&\mathscr{X}}
$$
is exact$;$ then
$$
\xymatrix@-1ex{S(\mathscr{R})\ar@<.5ex>[r]\ar@<-.5ex>[r]&\coprod_{\alpha\in L}S(\mathscr{U}_{\alpha})\ar[r]&S(\mathscr{X})}
$$
is exact. \hfill$\square$
\end{prop}

\begin{rem}\label{rem-spectralfunctor2}{\rm 
If $S\colon\mathscr{C}\rightarrow\Top$ is a continuous spectral functor in the situation as in \ref{prop-contnuityspectral}, then the following holds: if $\mathscr{U}\hookrightarrow\mathscr{V}$ is an open immersion, then $S(\mathscr{U})=S(\mathscr{V})$ if and only if $\mathscr{U}=\mathscr{V}$.
The proof is given by an argument similar to that of \ref{prop-spectral4}.}
\end{rem}

\begin{prop}\label{prop-exampleprespectral1}
The functor 
$$ 
[\,\cdot\,]\colon\Rf\longrightarrow\Top
$$
$($the separated quotient$;$ cf.\ {\rm \S\ref{subsub-separation}}$)$ is a continuous prespectral functor. 
\end{prop}

\begin{proof}
It is easy to see that the functor $[\,\cdot\,]$ is prespectral.
The continuity follows from {\bf \ref{ch-pre}}.\ref{prop-structuresepquot2}.
\end{proof}

\begin{prop}\label{prop-pointspoints2}
Suppose $\mathscr{C}$ is $\mathrm{QCO}$-stable\index{stable!QCOstable@QCO-{---}}.
Let $S\colon\mathscr{C}\rightarrow\Top$ be a prespectral functor. 
Then the following conditions are equivalent$:$ 
\begin{itemize}
\item[{\rm (a)}] let $\mathscr{X}$ be a coherent space in $\mathscr{C}$, and $\mathscr{U}\subseteq\mathscr{X}$ a quasi-compact open subspace$;$ then $\ZR{\mathscr{U}}\cap S(\mathscr{X}) = S(\mathscr{X})$ implies $\mathscr{U}=\mathscr{X};$
\item[{\rm (b)}] let $\mathscr{X}$ be a non-empty coherent space in $\mathscr{C}$, and $X$ a distinguished formal model of $\mathscr{X};$ then for any non-empty closed formal subscheme $Y$ of $X$ of finite presentation, we have $\sp^{-1}_{X}(Y)\cap S(\mathscr{X})\neq\emptyset$. 
\end{itemize}
\end{prop}

\begin{proof}
Let us first show the implication (a) $\Rightarrow$ (b).
Let $\mathscr{X}$, $X$, and $Y$ be as in (b).
Set $U=X\setminus Y$.
Then $U$ is a quasi-compact open subset of $X$, and hence $\sp^{-1}_X(U)$ is of the form $\ZR{\mathscr{U}}$, where $\mathscr{U}=U^{\rig}$ is a quasi-compact open subspace of $\mathscr{X}$ (\ref{prop-zariskiriemanntoptop} (2)).
Suppose $S(\mathscr{X})\subseteq\ZR{\mathscr{U}}$.
Then we have $\mathscr{U}=\mathscr{X}$.
Since $X$ is distinguished, we have $U=X$ (\ref{prop-ZRpoints4}), which contradicts the assumption $Y\neq\emptyset$.

Next we show the converse.
Suppose $S(\mathscr{X})\subseteq\ZR{\mathscr{U}}$.
Take a distinguished formal model $X$ of $\mathscr{X}$ containing a quasi-compact open subset $U$ corresponding to $\mathscr{U}$.
If $\mathscr{U}\neq\mathscr{X}$, then by \ref{prop-ZRpoints4} the complement $Y=X\setminus U$ is non-empty.
Since $U$ is quasi-compact, $Y$ carries a structure of a formal scheme in such a way that we have a closed immersion $Y\hookrightarrow X$ of finite presentation.
But $\sp^{-1}_X(Y)\cap S(\mathscr{X})\neq\emptyset$ contradicts $S(\mathscr{X})\subseteq\ZR{\mathscr{U}}$.
\end{proof}

\subsubsection{Regularity}\label{subsub-spectralfunctreg}
Let $X$ be a topological space.
Then an open subset $U\subseteq X$ is said to be {\em regular} if $(\ovl{U})^{\circ}=U$.
Similarly, a closed subset $C$ of $X$ is said to be {\em regular} if $\ovl{C^{\circ}}=C$.
Note that, if $S$ is a regular open (resp.\ closed) subset of $X$, the complement $X\setminus S$ is a regular closed (resp.\ open) subset.

\begin{dfn}\label{dfn-realvalued}{\rm 
A prespectral functor $S$ is said to be {\em real valued}\index{spectral functor!real valued spectral functor@real valued ---} if $S(\mathscr X)$ is a subset of $[\mathscr{X}]$ for any $\mathscr{X}\in\obj(\mathscr{C})$.}
\end{dfn}

\begin{prop}\label{prop-tubes92}
Let $\mathscr{C}$ be a $\mathrm{QCO}$-stable\index{stable!QCOstable@QCO-{---}} subcategory of $\Rf$, and $S\colon\mathscr{C}\rightarrow\Top$ a real valued spectral functor. 

{\rm (1)}  For a quasi-compact open immersion $\mathscr{U}\hookrightarrow\mathscr{X}$ in $\mathscr{C}$, the open subset $\ZR{\mathscr{U}}$ and its closure $\ovl{\ZR{\mathscr{U}}}$ are regular in $\ZR{\mathscr{X}}$.

{\rm (2)} For a rigid space $\mathscr{X}$ in $\mathscr{C}$, any tube open subset\index{tube!tube open subset@--- open subset} {\rm (\ref{dfn-tubes1})} of $\ZR{\mathscr{X}}$ is regular.
\end{prop}

\begin{proof}
First we prove (1). 
By \ref{dfn-spectral2} we have $S(\mathscr{U}) = S(\mathscr{X})\cap\ZR{\mathscr{U}}$. 
Since $S$ is real valued, it is also equal to $S(\mathscr{X})\cap[\mathscr{U}]$. 
Hence by \ref{prop-spectral4}, for two quasi-compact open immersions $\mathscr{U}\hookrightarrow\mathscr{X}$ and $\mathscr{V}\hookrightarrow\mathscr{X}$, $[\mathscr{U}]=[\mathscr{V}]$ implies $\mathscr{U}=\mathscr{V}$.
Then the regularity of $\ZR{\mathscr{U}}$ follows by an argument similar to that in the proof of {\bf \ref{ch-pre}}.\ref{prop-lemlocallystronglycompactconverse}.
The regularity of $\ovl{\ZR{\mathscr{U}}}$ follows as a consequence. 

We prove (2). 
Since $\mathscr{C}$ is $\mathrm{QCO}$-stable, any retrocompact open set $\mathscr{U}$ of $\mathscr{X}$ belongs to $\mathscr{C}$. 
By (1), $\ZR{\mathscr{X}}\setminus \ovl{\ZR{\mathscr{U}}}$ is a regular open subset. 
Any tube open subset of $\mathscr{X}$ is of this form by \ref{prop-tubes91}, and the claim follows. 
\end{proof}

\subsubsection{Density argument}\label{subsub-densityargument}
\begin{thm}\label{thm-densityargument}
Let $\mathscr{X}$ be a coherent rigid space, and $S$ a spectral functor\index{functor!spectral functor@spectral ---}\index{spectral functor} {\rm (\ref{dfn-spectral2})} defined on the category of quasi-compact open subspaces of $\mathscr{X}$.
Then we have
$$
\Gamma(\ZR{\mathscr{X}},\O^{\int}_{\mathscr{X}})=\{f\in\Gamma(\ZR{\mathscr{X}},\O_{\mathscr{X}})\,|\,f_x\in\O^{\int}_{\mathscr{X},x}\ \textrm{for any}\ x\in S(\mathscr{X})\}.
$$
\end{thm}

We first prove the following lemma:
\begin{lem}\label{lem-densityargumentqcompact}
Let $\mathscr{X}$ be a coherent rigid space.
We define a subset $\mathfrak{W}_f\subseteq\ZR{\mathscr{X}}$ for $f\in\Gamma(\ZR{\mathscr{X}},\O_{\mathscr{X}})$ by
$$
\mathfrak{W}_f=\{x\in\ZR{\mathscr{X}}\,|\,f_x\in\O^{\int}_{\mathscr{X},x}\}.
$$
Then $\mathfrak{W}_f$ is a quasi-compact open subset of $\ZR{\mathscr{X}}$.
\end{lem}

\begin{proof}
Take $x\in\ZR{\mathscr{X}}$.
If $x\in\mathfrak{W}_f$, then it is clear that $\mathfrak{W}_f$ contains a quasi-compact open neighborhood of $x$.
Suppose $x\not\in\mathfrak{W}_f$.
With the notation as in \ref{ntn-ZRpoints}, since $(f\ \mathrm{mod}\ \m_{B_x})\in V_x$ is non-zero, $f$ is invertible in $B_x=\O_{\mathscr{X},x}$.
This implies that $g=\frac{1}{f}$ belongs to $A_x=\O^{\int}_{\mathscr{X},x}$.
Thus there exist a quasi-compact open neighborhood of $x$ of the form $\mathscr{U}=U^{\rig}$, where $U$ is a quasi-compact open subset of a formal model $X$ of $\mathscr{X}$, and an element $h\in\Gamma(U,\O_U)$ that gives $g$ via the inductive limit.
Replacing $U$ by a smaller one if necessary, we may assume that $U$ is affine $U=\Spf B$, that $hB$ is an open ideal of $B$, and that $fh=1$.
Then we have $\mathfrak{W}_f\cap\ZR{\mathscr{U}}=(\Spf B_{\{h\}})^{\rig}$.
This implies that the inclusion map $\mathfrak{W}_f\hookrightarrow\ZR{\mathscr{X}}$ is quasi-compact (cf.\ {\bf \ref{ch-pre}}.\ref{dfn-quasicompactness} (2)) and hence that $\mathfrak{W}_f$ is quasi-compact.
\end{proof}

\begin{proof}[Proof of Theorem {\rm \ref{thm-densityargument}}]
Suppose $f\in\Gamma(\ZR{\mathscr{X}},\O_{\mathscr{X}})$ belongs to the right-hand side of the equality.
Then we have $S(\mathscr{X})=S(\mathscr{W}_f)$, where $\mathscr{W}_f$ is the quasi-compact open subspace of $\mathscr{X}$ such that $\ZR{\mathscr{W}_f}=\mathfrak{W}_f$ as in \ref{lem-densityargumentqcompact}.
Since $S$ is a spectral functor, we have $\mathscr{W}_f=\mathscr{X}$, and hence the assertion follows.
\end{proof}
\index{spectral functor|)}\index{functor!spectral functor@spectral ---|)}

\subsection{Classical points}\label{sub-classicalpoints}
\subsubsection{Point-like rigid spaces}\label{subsub-pointlike}
\begin{dfn}\label{dfn-pointlike}{\rm 
A rigid space $\mathscr{Z}$ is said to be {\em point-like}\index{rigid space!pointlike rigid space@point-like ---} if it is coherent and reduced\index{rigid space!reduced rigid space@reduced ---} and there exists a unique minimal point in $\ZR{\mathscr{Z}}$.}
\end{dfn}

\begin{prop}\label{prop-pointlike1}
Let $\mathscr{Z}$ be a universally Noetherian point-like rigid space.

{\rm (1)} The set $\ZR{\mathscr{Z}}$ coincides with the set $G_x$ of all generizations\index{generization} of the minimal point $x$.
In particular, the set $\ZR{\mathscr{Z}}$ is totally ordered with respect to the ordering by generization.

{\rm (2)} Any distinguished formal model $Z$ of $\mathscr{Z}$ is reduced and has the underlying topological space totally ordered by the ordering by generization with the unique minimal point.
\end{prop}

\begin{proof}
(1) is clear by definition.
Let $Z$ be a distinguished formal model of $\mathscr{Z}$.
Then $Z$ is reduced by \ref{cor-reducedrigidspaces}.
By \ref{prop-ZRpoints4} the specialization map $\sp_Z\colon\ZR{\mathscr{Z}}\rightarrow Z$ is surjective.
From this the assertion (2) follows.
\end{proof}

\begin{exa}\label{exa-pointlike}{\rm 
If $V$ is an $a$-adically complete valuation ring, then $(\Spf V)^{\rig}$ is a point-like rigid space.}
\end{exa}

\subsubsection{Structure of point-like rigid spaces}\label{subsub-pointlikestr}
\begin{prop}\label{prop-pointlikestr1}
Any universally Noetherian point-like\index{rigid space!pointlike rigid space@point-like ---} rigid space is a Stein affinoid\index{affinoid!Stein affinoid@Stein ---}\index{Stein affinoid}.
\end{prop}

\begin{proof}
Let $\mathscr{Z}$ be a point-like rigid space, and $x\in\ZR{\mathscr{Z}}$ the minimal point.
By \ref{prop-cohomologyrigidsp1str} one can take a Stein affinoid neighborhood $\mathscr{V}$ of $x$.
Then by \ref{prop-pointlike1} (1) one sees that $\mathscr{Z}=\mathscr{V}$.
\end{proof}

\begin{prop}\label{prop-pointlikestr15}
Let $\mathscr{Z}=(\Spf A)^{\rig}$ be a universally Noetherian Stein affinoid, and suppose that $A$ is an $I$-torsion free t.u.\ rigid-Noetherian ring, where $I\subseteq A$ is an ideal of definition of $A$.
Set $\Spec B=s(\mathscr{Z})=\Spec A\setminus V(I)$ $($that is, $B=\Gamma(\mathscr{Z},\O_{\mathscr{Z}}))$, and suppose that $\Spec B$ is Jacobson.
Then $\mathscr{Z}$ is point-like if and only if the following conditions are satisfied$:$
\begin{itemize}
\item[{\rm (a)}] $A$ is a local integral domain, and $B=\Frac(A);$
\item[{\rm (b)}] the integral closure $A^{\int}$ of $A$ in $B$ is an $I$-adically separated henselian valuation ring.
\end{itemize}
Moreover, in this situation, we have $A^{\int}=\O^{\int}_{\mathscr{Z},z}$, where $z$ is the unique minimal point of $\ZR{\mathscr{Z}}$.
\end{prop}

\begin{proof}
Suppose that $\mathscr{Z}$ is point-like, and let $z$ be the unique minimal point of $\ZR{\mathscr{Z}}$.
Since $s(\mathscr{Z})=\Spec A\setminus V(I)=\Spec B$ is reduced and $A$ is $I$-torsion free, $A$ is a reduced ring.
Let $\alpha\colon (\Spf\widehat{V}_z)^{\rig}\rightarrow\mathscr{Z}$ be the associated rigid point at $z$.
Since $\ZR{\mathscr{Z}}$ coincides with the set of all generizations of $z$, and since $\Spf A$ is a distinguished formal model of $\mathscr{Z}$, the composite map $\Spf\widehat{V}_z\rightarrow\ZR{\mathscr{Z}}\rightarrow\Spf A$ is surjective (cf.\ \ref{prop-ZRpoints4}).

We claim that the underlying topological space of $s(\mathscr{Z})=\Spec B$ consists of one point.
Let $x\in s(\mathscr{Z})$ be a closed point, and suppose $s(\mathscr{Z})\neq\{x\}$.
Since $\Spec B$ is Jacobson, there exists another closed point $y\neq x$ of $s(\mathscr{Z})$.
But, as we have seen in \ref{prop-closedsubspaceaffinoid2} (2), non-empty closed subschemes of $s(\mathscr{Z})$ correspond bijectively to non-empty closed subspaces of $\mathscr{Z}$, which are, however, always supported on the whole topological space $\ZR{\mathscr{Z}}$ due to \ref{prop-closedsubspaceaffinoid2} (1), since $\mathscr{Z}$ is point-like.
But this is absurd, and hence we deduce by contradiction that $s(\mathscr{Z})=\Spec B$ consists of one point.
This shows that the map $\Spec\widehat{V}_z\rightarrow\Spec A$ is surjective, and hence $\Spec A$ is irreducible.
Hence the ring $A$ is a local integral domain such that $B=\Frac(A)$; in particular, we have a local injective homomorphism $A\hookrightarrow\widehat{V}_z$.

By \ref{prop-corlemtuavsrigidaff} we know that $A^{\int}=\Gamma(\mathscr{Z},\O^{\int}_{\mathscr{Z}})=\O^{\int}_{\mathscr{Z},z}$.
The composite map $A\hookrightarrow A^{\int}=\O^{\int}_{\mathscr{Z},z}\rightarrow V_z$ is, since $A\hookrightarrow\widehat{V}_z$ is injective, a local injective homomorphism.
Now since $A$ is a subring of $V_z$, and since $V_z$ is integrally closed, we have $A^{\int}\subseteq V_z$ and hence $A^{\int}=V_z$, which is an $I$-adically separated henselian valuation ring (cf.\ \ref{ntn-ZRpoints}).

Conversely, suppose that $\mathscr{Z}$ satisfies the conditions (a) and (b).
Set $Y=\Spec A$ and $U=\Spec B$, and consider the $U$-admissible blow-ups\index{admissible!U-admissible blow-up@$U$-{---} blow-up}\index{blow-up!U-admissible blow-up@$U$-admisible ---} (which are, by the passage to the formal completions, exactly corresponding to admissible blow-ups of $X=\Spf A$) of $Y$ (cf.\ \ref{dfn-Uadmissibleblowups}).
We denote by $\ZR{Y}$ the classical Zariski-Riemann space, that is, the projective limit of all $U$-admissible blow-ups (cf.\ \ref{dfn-classicalZRsp1}).
By the valuative criterion, one has a morphism $\Spec A^{\int}\rightarrow\ZR{Y}$; let $\til{z}$ be the image of the closed point.
Now we claim that any blow-up $Y'\rightarrow Y=\Spec A$ along an admissible ideal is affine.
Take an affine open subset $U'=\Spec A'\subseteq Y'$ that contains the image of $\til{z}$ by the specialization map $\ZR{Y}\rightarrow Y'$.
We have an inclusion $A'\hookrightarrow A^{\int}$ that factorizes $A\hookrightarrow A^{\int}$.
Suppose there exists a closed point $y\in Y'\setminus U'$.
There exists a valuation ring $V$ and a morphism $\Spec V\rightarrow Y'$ that dominates $y$.
Since $y$ is also mapped to $z$, which is the unique closed point of $\Spec A$, we have a map $A^{\int}\hookrightarrow V$ that factorizes $A\hookrightarrow V$.
Hence in $\ZR{Y}$ the the point $\til{z}$ is a specialization of the image of the closed point of $\Spec V\rightarrow\ZR{Y}$.
But this is absurd, since the specialization map $\ZR{Y}\rightarrow Y'$ is a closed map (cf.\ \ref{thm-classicalZRsp1}).

In particular, we deduce that any admissible blow-up $X'\rightarrow\Spf A$ is affine, say $X'=\Spf A'$, where $A'$ is contained in $A^{\int}$.
Hence we have the surjective map $\Spf\widehat{A^{\int}}\rightarrow X'$.
Since this holds for any admissible blow-up $X'$, by {\bf \ref{ch-pre}}.\ref{thm-projlimcohspacepres} (2) we have the surjective map $\Spf\widehat{A^{\int}}\rightarrow\ZR{\mathscr{Z}}$.
Since $\widehat{A^{\int}}$ is an $I$-adically complete valuation ring ({\bf \ref{ch-pre}}.\ref{thm-compval2006ver1}), we deduce that $\ZR{\mathscr{Z}}$ has a unique minimal point. By \ref{prop-corlemtuavsrigidaff}, we have $\O^{\int}_{\mathscr{Z},z}=A^{\int}=V_z$. Since the local ring at the maximal point is given by $B$, the fractional field of $\O^{\int}_{\mathscr{Z},z}=A^{\int}=V_z$, it follows by \ref{thm-fibersoverrigptsbehavior} that the local ring at each point of $\ZR{\mathscr{Z}}$ is a subring of $B$, and hence is an integral domain. Hence $\mathscr{Z}$ is reduced, and thus we have shown that $\mathscr{Z}$ is point-like.
\end{proof}

\begin{prop}\label{prop-pointlikestr2}
{\rm (1)} Let $V$ be an $a$-adically complete valuation ring\index{valuation!valuation ring@--- ring!a-adically complete valuation ring@$a$-adically complete --- ---}, and $\mathscr{Z}$ a point-like rigid space of finite type over $\mathscr{S}=(\Spf V)^{\rig}$.
Then $\mathscr{Z}$ is of the form $\mathscr{Z}=(\Spf W)^{\rig}$, where $W$ is quasi-finite, flat, and finitely presented over $V;$ if $V$ is of height one, then $W$ is finite over $V$.

{\rm (2)} Let $\mathscr{Z}$ be a point-like rigid space of type {\rm (N)} {\rm (\S\ref{subsub-examplesN})}.
Then $\mathscr{Z}$ is of the form $\mathscr{Z}=(\Spf W)^{\rig}$, where $W$ is an $a$-adically complete discrete valuation ring\index{valuation!valuation ring@--- ring!discrete valuation ring@discrete --- ---}.
\end{prop}

We show the following lemma before the proof:
\begin{lem}\label{lem-zariskianartinian}
Let $A$ be a Noetherian $I$-adically Zariskian ring.
Then the scheme $\Spec A\setminus V(I)$ is Jacobson {\rm (\cite[$\mathbf{IV}$, (10.4.1)]{EGA})}, and the closure $\ovl{\{\mathfrak{p}\}}$ in $\Spec A$ of any closed point $\mathfrak{p}\in\Spec A\setminus V(I)$ is of the form $\Spec B$ where $B$ is a $1$-dimensional semi-local ring$;$ moreover, we have $\ovl{\{\mathfrak{p}\}}\setminus V(I)=\{\mathfrak{p}\}$.
\end{lem}

\begin{proof}
Since $1+I\subseteq A^{\times}$, one sees easily that any non-empty closed subset of $\Spec A$ meets $V(I)$.
Hence it follows from \cite[$\mathbf{IV}$, (10.5.7)]{EGA} that $\Spec A\setminus V(I)$ is a Jacobson scheme.
Let $\mathfrak{p}$ be an ideal maximal in $\Spec A\setminus V(I)$, and set $Y=\ovl{\{\mathfrak{p}\}}$, where the closure is taken in $\Spec A$.
Since $\mathfrak{p}$ is not an open prime, it is not maximal in $A$ ({\bf \ref{ch-pre}}.\ref{prop-relpair32}).
Let $\mathfrak{q}$ be any prime ideal that strictly contains $\mathfrak{p}$.
Then $\mathfrak{q}$ is an open prime ideal, and hence $A_{\mathfrak{q}}$ is $IA_{\mathfrak{q}}$-adically Zariskian ({\bf \ref{ch-pre}}.\ref{prop-relpair32}).
By \cite[Chap.\ VIII, Theorem 10)]{ZSII} $A/\mathfrak{q}$ is an Artinian local ring.
Hence $\mathfrak{q}$ is actually a maximal ideal of $A$, and thus we deduce that the ring $B$ as above is of dimension $1$.
By \cite[$\mathbf{IV}$, (10.5.3)]{EGA} we deduce that $B$ is a semi-local ring.
\end{proof}

\begin{proof}[Proof of Proposition {\rm \ref{prop-pointlikestr2}}]
By \ref{prop-pointlikestr1} we have $\mathscr{Z}=(\Spf W)^{\rig}$, where $W$ is topologically finitely generated $V$-algebra in the case (1) or is a Noetherian adic ring in the case (2).
In both cases the scheme $\Spec W\setminus V(I)$ ($I\subseteq W$ is an ideal of definition) is Jacobson due to {\bf \ref{ch-pre}}.\ref{prop-classicalaffringjacobson} and \ref{lem-zariskianartinian}.
In view of \ref{prop-pointlikestr15} we may assume that $W$ is an $I$-adically complete local integral domain and that its integral closure $W^{\int}$ in $\Frac(W)$ is an $I$-adically separated henselian valuation ring.

(1) By \ref{thm-morbetweenaffinoid1} we may assume that the map $\mathscr{Z}\rightarrow\mathscr{S}$ has a distinguished formal model of the form $\Spf W\rightarrow\Spf V$.
Since we have the surjective map $\Spec W^{\int}\rightarrow\Spec W$, $W$ is a local ring.
By the assumption, $W$ is topologically of finite presentation over $V$ (\ref{prop-cohrigidspacefintype1}), and $a$-torsion free (hence $W$ is flat over $V$).

We claim that the map $\Spec W\rightarrow\Spec V$ is finite in case $V$ is of height one.
The map $\Spf W\rightarrow\Spf V$ is surjective, and the closed fiber of $\Spec W\rightarrow\Spec V$ is of dimension $0$, since any scheme of finite type over a field of positive dimension can never be totally ordered with respect to the order by generalization.
By Noether normalization ({\bf \ref{ch-pre}}.\ref{thm-noethernormalizationtype(V)}) we deduce that $W$ is finite over $V$.

For a general $V$ we use the technique alluded in {\bf \ref{ch-pre}}.\ref{rem-trick}.
Let $\mathfrak{q}=\sqrt{(a)}$ be the height one prime (cf.\ {\bf \ref{ch-pre}}.\ref{prop-maxspe}).
Since $W/\mathfrak{q}W$ is flat over $V/\mathfrak{q}$, and since the generic fiber of the map $\Spec W/\mathfrak{q}W\rightarrow\Spec V/\mathfrak{q}$ is of dimension $0$ (by the above-mentioned reasoning), we deduce that $W$ is of relative dimension $0$ over $V$.
Hence $W$ is quasi-finite over $V$, as desired.

(2) By \ref{lem-zariskianartinian} we know that $W$ is a local integral domain of dimension $1$.
Since $W$ is Noetherian, it follows that $I$ contains $\mathfrak{m}^n_W$ for some $n>0$; that is, the $I$-adic topology coincides with the topology defined by the maximal ideal.
By \cite[(32.2)]{Nagata1} we deduce that its normalization $W^{\int}$ is finite over $W$.
Hence in view of \ref{prop-lemmorbetweenaffnoid1} we may assume that $W$ is normal.
Then $W$ is a Noetherian $1$-dimensional integrally closed local domain, hence is a discrete valuation ring, as desired.
\end{proof}

\subsubsection{Classical points}\label{subsub-classicalpoints}
\begin{dfn}\label{dfn-classicalpoint}{\rm 
Let $\mathscr{X}$ be a rigid space of type (V)\index{rigid space!rigid space of typeV@--- of type (V)} (\ref{dfn-typeV}) or of type (N)\index{rigid space!rigid space of typeN@--- of type (N)} (\ref{dfn-typeN}).

(1) A {\em classical point}\index{point!classical point@classical ---}\index{classical point} of $\mathscr{X}$ is a point-like locally closed rigid subspace $\mathscr{Z}\subseteq\mathscr{X}$ that is retrocompact, that is, the immersion $\mathscr{Z}\hookrightarrow\mathscr{X}$ is quasi-compact.

(2) A classical point $\mathscr{Z}\subseteq\mathscr{X}$ is said to be {\em closed}\index{classical point!closed classical point@closed ---} if it is a closed subspace.}
\end{dfn}

Recall that a rigid space $\mathscr{X}$ is said to be of type {\rm ($\mathrm{V_{\R}}$)}\index{rigid space!rigid space of typeV1@--- of type ($\mathrm{V_{\R}}$)} if it is locally of finite type over a rigid space of the form $(\Spf V)^{\rig}$ where $V$ is an $a$-adically complete valuation ring ($a\in\m_V\setminus\{0\}$) {\em of height one} (\ref{dfn-typeV}).
\begin{prop}\label{prop-classicalpointsclosedheightone}
{\rm (1)} Any classical point of a rigid space of type {\rm ($\mathrm{V_{\R}}$)} is closed. 

{\rm (2)} Let $\mathscr{X}$ be a coherent rigid space of type {\rm (N)}.
Suppose $\mathscr{X}$ has a Noetherian distinguished formal model $X$ and an ideal of definition $\mathscr{I}$ such that the scheme $X_0=(X,\O_X/\mathscr{I})$ is Jacobson.
Then any classical point of $\mathscr{X}$ is closed.
\end{prop}

\begin{proof}
(1) follows immediately from \ref{prop-pointlikestr2} (1).
To show (2), let $\mathscr{Z}\hookrightarrow\mathscr{X}$ be a classical point; $\mathscr{Z}$ is a closed rigid subspace of a coherent open subspace $\mathscr{U}\subseteq\mathscr{X}$, and is of the form $\mathscr{Z}=(\Spf W)^{\rig}$ for an $a$-adically complete discrete valuation ring $W$; see \ref{prop-pointlikestr2} (2).
Passing to formal models, we have the sequence of morphisms 
$$
\Spf W\stackrel{i}{\longhookrightarrow}U\stackrel{j}{\longhookrightarrow}X'\stackrel{\pi}{\longrightarrow}X,
$$
where $U^{\rig}=\mathscr{U}$, $i$ is a closed immersion, $j$ is an open immersion, and $\pi$ is an admissible blow-up.
Dividing out by the ideal of definition $\mathscr{I}$, we get the sequence of morphisms of schemes
$$
\Spec W/aW\stackrel{i_0}{\longhookrightarrow}U_0\stackrel{j}{\longhookrightarrow}X'_0\stackrel{\pi}{\longrightarrow}X_0.
$$
Since $X'_0$ is Jacobson, the image of $\Spec W/aW$ in $X'_0$ is a closed point ({\rm \cite[$\mathbf{IV}$, (10.4.7)]{EGA}}).
Hence, by \ref{cor-univclosedadicred0}, $\Spec W\hookrightarrow X'$ is a closed immersion, and hence $\mathscr{Z}\hookrightarrow\mathscr{X}$ is a closed immersion.
\end{proof}

\begin{prop}\label{prop-classicalpointsexist}
Let $\mathscr{X}=(\Spf A)^{\rig}$ be an affinoid of type {\rm ($\mathrm{V_{\R}}$)}\index{rigid space!rigid space of typeV1@--- of type ($\mathrm{V_{\R}}$)} or of type {\rm (N)}\index{rigid space!rigid space of typeN@--- of type (N)}, where $A$ is $I$-torsion free for a finitely generated ideal of definition $I\subseteq A$.

{\rm (1)} For any closed point $x$ of the Noetherian scheme $s(\mathscr{X})$, there exists a unique closed classical point $\mathscr{Z}=Z^{\rig}\hookrightarrow\mathscr{X}$ of $\mathscr{X}$ such that the image of $s(\mathscr{Z})\rightarrow s(\mathscr{X})$ is $x$.

{\rm (2)} Suppose $($in type {\rm (N)} case$)$ that $\Spec A/I$ is Jacobson.
Then, for any classical point $\mathscr{Z}\hookrightarrow\mathscr{X}$ of $\mathscr{X}$, $s(\mathscr{Z})$ is a point, and the image of the map $s(\mathscr{Z})\rightarrow s(\mathscr{X})$ is a closed point of $s(\mathscr{X})$.
\end{prop}

\begin{proof}
(1) Let $x\in s(\mathscr{X})$ be a closed point.
The construction as in \S\ref{subsub-closedsubspacesofaffinoids} gives the closed subspace $\mathscr{Z}=\{x\}\times_{s(\mathscr{X})}\ZR{\mathscr{X}}\hookrightarrow\mathscr{X}$ corresponding to $x$.
We claim that this indeed gives a closed classical point of $\mathscr{X}$.
To describe $\mathscr{Z}$, let us take the reduced closure $Y=\Spec W$ of $\{x\}$ in $\Spec A$.
Then we have $\mathscr{Z}=(\Spf W)^{\rig}$.
In type {\rm ($\mathrm{V_{\R}}$)} case, by {\bf \ref{ch-pre}}.\ref{cor-noethernormalizationtype(V)}, $W$ is finite over $V$, and hence, by \ref{prop-corlemtuavsrigidaff} and \ref{prop-pointlikestr15}, $\mathscr{Z}=(\Spf W)^{\rig}$ is a closed classical point of $\mathscr{X}$.
In the type (N) case, $W$ is a Noetherian $1$-dimensional semi-local domain (\ref{lem-zariskianartinian}); by an argument similar to that in the proof of \ref{prop-pointlikestr2} (2), we deduce that $W^{\int}$, the integral closure of $W$, is finite over $W$, and is a complete discrete valuation ring.
Thus, similarly, $\mathscr{Z}=(\Spf W)^{\rig}$ is a closed classical point.
The uniqueness is easy to see, and is left to the reader.

(2) Let $\mathscr{Z}=(\Spf W)^{\rig}\hookrightarrow\mathscr{X}=(\Spf A)^{\rig}$ be a classical point.
One can replace, due to \ref{thm-morbetweenaffinoid1}, $W$ by a finite algebra isomorphic outside $I$, so that we have an immersion $\Spf W\hookrightarrow\Spf A$.
By \ref{prop-classicalpointsclosedheightone}, one can even do this in such a way that $\Spf W\hookrightarrow\Spf A$ is a closed immersion.
Then one obtains the closed immersion $\Spec W\rightarrow\Spec A$, and hence, $s(\mathscr{Z})=\Spec W[\frac{1}{a}]$ is a point mapped to a closed point of $s(\mathscr{X})=\Spec A[\frac{1}{a}]$.
\end{proof}

\begin{cor}\label{cor-classicalpointsexist0}
Let $\mathscr{X}=(\Spf A)^{\rig}$ be either an affinoid of type {\rm ($\mathrm{V_{\R}}$)}, or an affinoid of type {\rm (N)} having a distinguished Noetherian formal model $X$, together with an ideal of definition $\mathscr{I}$, such that the scheme $X_0=(X,\O_X/\mathscr{I})$ is Jacobson.
Then there exists a natural one to one correspondence between the set of all classical points of $\mathscr{X}$ and the set $s(\mathscr{X})^{\cl}$ of all closed points of the scheme $s(\mathscr{X})$.
\end{cor}

\begin{cor}\label{cor-classicalpointsexist3}
Let $\mathscr{X}$ be a rigid space of type {\rm ($\mathrm{V_{\R}}$)}\index{rigid space!rigid space of typeV1@--- of type ($\mathrm{V_{\R}}$)} or of type {\rm (N)}\index{rigid space!rigid space of typeN@--- of type (N)}, and $\mathscr{F}$ a coherent sheaf on $\mathscr{X}$.
Then $\mathscr{F}=0$ if and only if $\mathscr{F}|_{\mathscr{Z}}=0$ for any classical point $\mathscr{Z}\hookrightarrow\mathscr{X}$.
\end{cor}

\begin{proof}
We may assume that $\mathscr{X}$ is an affinoid, say $\mathscr{X}=(\Spf A)^{\rig}$.
By \ref{thm-affinoidscohsheaf} the sheaf $\mathscr{F}$ corresponds to a coherent sheaf $\mathscr{G}$ on the associated Noetherian scheme $s(\mathscr{X})=\Spec A\setminus V(I)$.
Clearly, we have $\mathscr{F}=0$ if and only if $\mathscr{G}=0$, which is further equivalent to that $\mathscr{G}_x=0$ for any {\em closed} point of $s(\mathscr{X})$ (due to \cite[$\mathbf{IV}$, (5.1.11)]{EGA}).
In view of the correspondence as in \ref{cor-classicalpointsexist0}, this is equivalent to $\mathscr{F}|_{\mathscr{Z}}=0$ for any classical point $\mathscr{Z}\hookrightarrow\mathscr{X}$ (here we have again used \ref{thm-affinoidscohsheaf}).
\end{proof}

\begin{cor}\label{cor-classicalpointsexist}
Let $\mathscr{X}$ be a non-empty rigid space of either type {\rm ($\mathrm{V_{\R}}$)}\index{rigid space!rigid space of typeV1@--- of type ($\mathrm{V_{\R}}$)} or type {\rm (N)}\index{rigid space!rigid space of typeN@--- of type (N)}.
Then $\mathscr{X}$ has a classical point.
Moreover, any non-empty rigid subspace $\mathscr{Y}$ of $\mathscr{X}$ has a classical point of $\mathscr{X}$.
\end{cor}

\begin{proof}
We may assume that $\mathscr{X}$ is affinoid.
Then the first assertion follows from \ref{cor-classicalpointsexist0}.
For the second, one can assume moreover that $\mathscr{Y}$ is an affinoid, and consider the immersion $\mathscr{Y}\hookrightarrow\mathscr{X}$.
If $\mathscr{Z}\hookrightarrow\mathscr{Y}$ is a classical point, then the composition $\mathscr{Z}\hookrightarrow\mathscr{X}$ defines a classical point of $\mathscr{X}$ (\ref{prop-immersionrigid5} (1)).
\end{proof}

In the sequel, we denote by 
$$
\ZR{\mathscr{X}}^{\cl}
$$
the set of all (isomorphism classes of) classical points of $\mathscr{X}$.
Notice that, if $\mathscr{X}$ is of type {\rm ($\mathrm{V_{\R}}$)} or {\rm (N)}, then due to \ref{prop-classicalpointsexist} and \ref{prop-closedsubspaceaffinoid2} (2) one sees that the set $\ZR{\mathscr{X}}^{\cl}$ can be naturally regarded as a subset of $\ZR{\mathscr{X}}$.
Moreover, if $\mathscr{Y}\subseteq\mathscr{X}$ is a rigid subspace, we have 
$$
\ZR{\mathscr{Y}}^{\cl}=\ZR{\mathscr{X}}^{\cl}\cap\ZR{\mathscr{Y}}.
$$

\subsubsection{Functoriality}\label{subsub-classicalpointsfunct}
\begin{prop}\label{prop-classicalpoint1}
Let $\mathscr{Y}$ be a rigid space of type {\rm ($\mathrm{V_{\R}}$)}\index{rigid space!rigid space of typeV1@--- of type ($\mathrm{V_{\R}}$)}, and $\varphi\colon\mathscr{X}\rightarrow\mathscr{Y}$ a morphism locally of finite type.
Let $\iota\colon\mathscr{Z}\hookrightarrow\mathscr{X}$ be a classical point of $\mathscr{X}$.
Then there exists uniquely a classical point $\iota'\colon\mathscr{Z}'\hookrightarrow\mathscr{Y}$ of $\mathscr{Y}$ that admits the commutative diagram
$$
\xymatrix{\mathscr{Z}\,\ar@{^{(}->}[r]^{\iota}\ar[d]_{\phi}&\mathscr{X}\ar[d]^{\varphi}\\ \mathscr{Z}'\,\ar@{^{(}->}[r]_{\iota'}&\mathscr{Y}\rlap{,}}
$$
where $\phi$ is a finite morphism.
\end{prop}

\begin{proof}
We may assume $\mathscr{Z}$, $\mathscr{X}$, and $\mathscr{Y}$ are affinoid, say $\mathscr{Z}=(\Spf W)^{\rig}$, $\mathscr{X}=(\Spf A)^{\rig}$, and $\mathscr{Y}=(\Spf B)^{\rig}$, where $W$, $A$, and $B$ are $a$-torsion free, and that $\iota$ is a closed classical point.
Moreover, we may assume that there exist a map $\Spf A\rightarrow\Spf B$ of finite presentation and a closed immersion $\Spf W\hookrightarrow\Spf A$ of finite presentation that give formal models of $\varphi$ and $\iota$, respectively.
By \ref{prop-classicalpointsexist} the closed immersion $\Spec W\hookrightarrow\Spec A$ defines a closed point $x$ of $s(\mathscr{X})=\Spec A[\frac{1}{a}]$.
Since $A[\frac{1}{a}]$ and $B[\frac{1}{a}]$ are classical affinoid algebras\index{algebra!affinoid algebra@affinoid ---!classical affinoid algebra@classical --- ---} over $K=\Frac(V)$ ({\bf \ref{ch-pre}}, \S\ref{subsub-classicalaffinoidalgebras}), the image $y$ of $x$ by the map $\Spec A\rightarrow\Spec B$ is a closed point of $s(\mathscr{Y})$ by {\bf \ref{ch-pre}}.\ref{cor-Tatealgebraprime} and {\bf \ref{ch-pre}}.\ref{thm-northernormaclassaff}.
Hence the point $y$ gives rise to a classical point $\iota'\colon\mathscr{Z}'\hookrightarrow\mathscr{Y}$, which comes from the formal model of the form $\Spf W'\rightarrow\Spf B$ as in \ref{prop-classicalpointsexist}.
By the construction we have the map $\Spf W\rightarrow\Spf W'$ that gives the map $\phi$ as above.
The morphism $\Spec W\rightarrow\Spec W'$ is finite due to {\bf \ref{ch-pre}}.\ref{cor-noethernormalizationtype(V)}.
The uniqueness follows immediately from the construction.
\end{proof}

It follows from the proposition that we have the functor
$$
\mathscr{X}\longmapsto\ZR{\mathscr{X}}^{\cl}
$$
defined on the category of rigid spaces of type {\rm ($\mathrm{V_{\R}}$)}\index{rigid space!rigid space of typeV1@--- of type ($\mathrm{V_{\R}}$)}.
By what we have seen in \S\ref{subsub-classicalpoints}, this is a prespectral functor (\ref{dfn-spectral1})\index{spectral functor!prespectral functor@pre{---}}.
We will see in the next paragraph that this is in fact a spectral functor\index{functor!spectral functor@spectral ---}\index{spectral functor}.

\subsubsection{Spectrality}\label{subsub-classicalpointsspectral}
\begin{prop}\label{prop-spectralfunctorseparatedquotients}
Let $X$ be either a coherent formal scheme of finite type over an $a$-adically complete valuation ring $V$ of height one $($where $a\in\m_V\setminus\{0\})$ or a coherent Noetherian formal scheme.
We assume that $X$ is $\mathscr{I}$-torsion free, where $\mathscr{I}$ is an ideal of definition of finite type on $X$.
Then the map
$$
\sp_X|_{[\mathscr{X}]}\colon[\mathscr{X}]\longrightarrow X
$$
by the specialization map {\rm (\S\ref{subsub-ZRdef})}\index{specialization map}, where $\mathscr{X}=X^{\rig}$, is surjective.
\end{prop}

Let us first prove the proposition in the {\rm ($\mathrm{V_{\R}}$)} case.
\begin{lem}\label{lem-spectralfunctorseparatedquotientsVR1}
Let $V$ be an $a$-adically complete valuation ring of height one $($where $a\in\m_V\setminus\{0\})$, $K=\Frac(V)$, and $A$ a topologically of finite type $V$-algebra.
Suppose $A$ is integral and $V$-flat and admits a finite injection $R=V\dl T_1,\ldots,T_n\dr\hookrightarrow A$.
Let $L=\Frac(R)$, and denote the localization of $R$ at $\sqrt{a}R$ by $W$ $($notice that $W$ is a valuation ring$)$.
Then the integral closure $A^{\int}$ of $A$ in $B=A[\frac{1}{a}]$ is described as
$$
{\textstyle 
A^{\int}=\{f\in A[\frac{1}{a}]\,|\,\textrm{$f$ is integral over $A\otimes_RW$}\}.}
$$
\end{lem}

\begin{proof}
Only the inclusion `$\supseteq$' calls for the proof.
Let $f\in A[\frac{1}{a}]$ be an element of the right-hand side, and $P=P(T)$ the characteristic polynomial of the multiplication map by $f$ in $\End_L(A\otimes_RL)$.
Since $f$ is integral over the Tate algebra\index{Tate, J.}\index{algebra!Tate algebra@Tate ---} $R[\frac{1}{a}]=R\otimes_VK$, all coefficients of $P(T)$ is integral over $R[\frac{1}{a}]$ (\cite[Chap.\ V, \S1.6, Prop.\ 17]{Bourb1}).
Since $R[\frac{1}{a}]$ is integrally closed (cf.\ {\bf \ref{ch-pre}}.\ref{prop-exerTateUFD}), we have $P(T)\in R[\frac{1}{a}][T]$.
Similarly, since $f$ is integral over $A\otimes_RW$, we have $P(T)\in W[T]$.
Hence we deduce that $P(T)\in R[T]$.
Now since $P(f)=0$ in $A\otimes_RL$, and since $A$ is integral, we deduce that $P(f)=0$ in $A$.
This implies that $f$ is integral over $R$ and hence integral over $A$.
\end{proof}

\begin{proof}[Proof of Proposition {\rm \ref{prop-spectralfunctorseparatedquotients}} in the type {\rm ($\mathrm{V_{\R}}$)} case]
First notice that we may assume that $\mathscr{X}$ is reduced and thus that $X$ is reduced (here we take the {\em reduced model} of $\mathscr{X}$, which will be discussed later in \S\ref{subsub-closedimmrigidred}).
Take $y\in X$, and let $Y$ be a closed subscheme of $X$ defined by an admissible ideal such that the underlying topological space of $Y$ coincides with $\ovl{\{y\}}$, the closure of the singleton set $\{y\}$ in $X$; such a $Y$ exists, for the scheme $X_{\red}$ is Noetherian.

First we assume that $y$ is a closed point of $X$.
We may replace $X$ by the admissible blow-up along the defining ideal of $Y$. Moreover,  we may replace $X$ by an affine neighborhood of $y$. In this way, we may assume that $X$ is affine, and $Y$ is defined by an element $f\in A$.
Since $A$ is an integral domain and $(f)\subseteq A$ is an admissible ideal, we have $1/f\in A[\frac{1}{a}]$.
By Noether normalization theorem ({\bf \ref{ch-pre}}.\ref{thm-noethernormalizationtype(V)}) we have a finite injection $R=V\dl T_1,\ldots,T_n\dr\hookrightarrow A$.
Let $L=\Frac(R)$, and $W$ the localization of $R$ at $\sqrt{a}R$.
Consider the finite $W$-algebra $D=A\otimes_RW$.
Since $D$ is $a$-torsion free, $D$ is finite flat over $W$, and $fD$ is an admissible ideal of $D$.
By \ref{lem-spectralfunctorseparatedquotientsVR1} the equality $fD=D$ would imply that $1/f$ is integral over $A$, which is absurd, since $Y$ is non-empty.
Hence we have $fD\neq D$, and we are reduced to the case $X=\Spf D\otimes_W\widehat{W}$ with $V$ replaced by $\widehat{W}$; the assertion in this case is obvious, for $D\otimes_W\widehat{W}$ is finite over $\widehat{W}$.

In the general case, we take an extension $V\rightarrow V'$ of $a$-adically complete valuation rings, and perform the base-change by $\Spf V'\rightarrow\Spf V$; let us denote $X_{V'}$ and $Y_{V'}$ the $V'$-formal schemes obtained by the base-change of $X$ and $Y$, respectively.
Let $k'$ be the residue field of $V'$.
If one can take such an extension $V\rightarrow V'$ so that there exists an $k'$-rational point of $Y_{k'}$ mapped by the canonical map to the point $y$, the situation is reduced to the above-treated case.
This is indeed possible as follows: 
The generic point $y=\Spec K'$ of $Y$ can be seen as a $K'$-rational point $y'$ of $Y\times_{\Spec k}\Spec K'$, where $k=V/\m_V$ is the residue field of $V$.
Take a local flat extension $V\rightarrow U$ such that the closed fiber of $\Spec U\rightarrow\Spec V$ is isomorphic to $\Spec K'\rightarrow\Spec k$. 
Taking a valuation ring dominating $U$ and the $a$-adic completion, we have the desired extension of valuation rings.
\end{proof}

The following lemma, which we will use in the proof in the type (N) case, is easy to see, and the proof is left to the reader:
\begin{lem}\label{lem-spectralfunctorseparatedquotients}
Let $\mathscr{X}$ be a coherent universally Noetherian rigid space, and $X$ a distinguished formal model with an ideal of definition $\mathscr{I}_X$ of finite type.
Suppose that for any point $x\in X$ there exist a valuation ring $V$ of height one and a local homomorphism $\O_{X,x}\rightarrow V$ such that $\mathscr{I}_{X,x}V\neq V,0$. 
Then $\sp_X|_{[\mathscr{X}]}\colon[\mathscr{X}]\rightarrow X$ is surjective. \hfill$\square$
\end{lem}

\begin{proof}[Proof of Proposition {\rm \ref{prop-spectralfunctorseparatedquotients}} in the type {\rm (N)} case]
By \ref{lem-spectralfunctorseparatedquotients} it suffices to prove the following: for a Noetherian local ring $A$ with the maximal ideal $\m_A$ and an ideal $I\subseteq A$ such that $A$ is $I$-torsion free, there is a height one valuation ring $V$ and a local homomorphism $A\rightarrow V$ such that $IV\neq V,0$. 
We may assume that $A$ is complete. 
Since $A$ is $I$-torsion free, $U=\Spec A\setminus V(\m_A)$ is non-empty.
Since $U$ is Jacobson (\cite[$\mathbf{IV}$, (10.5.9)]{EGA}), there is a closed point $x$ of $\Spec A\setminus V(I)$ that is also closed in $U$.
By replacing $\Spec A$ by the closure of $x$ in $\Spec A$, we may assume that $A$ is integral and one dimensional (\cite[$\mathbf{IV}$, (10.5.9)]{EGA}).
The claim is clear in this case. 
\end{proof}

Now let us state an immediate but important corollary of \ref{prop-spectralfunctorseparatedquotients}:
\begin{cor}\label{cor-spectralfunctorseparatedquotientscor}
Let $X$ be as in {\rm \ref{prop-spectralfunctorseparatedquotients}}, $Y$ a closed subscheme of $X$ defined by an admissible ideal of $X$, and $U=X\setminus Y$.
Set $\mathscr{X}=X^{\rig}$ and $\mathscr{U}=U^{\rig}$.
Then the following conditions are equivalent to each other$:$
\begin{itemize}
\item[{\rm (a)}] $Y$ is non-empty$;$
\item[{\rm (b)}] $[\mathscr{U}]\neq[\mathscr{X}];$ 
\item[{\rm (c)}] the tube open subset\index{tube!tube open subset@--- open subset} $C_{Y|X}$ {\rm (cf.\ \ref{prop-tubes91})} is non-empty. \hfill$\square$
\end{itemize}
\end{cor}

In particular, we deduce the following:
\begin{thm}\label{thm-spectralfunctorseparatedquotients}
The functor
$$
\mathscr{X}\longmapsto[\mathscr{X}]
$$
by the separated quotients {\rm (\S\ref{subsub-separation})} is a spectral functor on the category of rigid spaces of type {\rm ($\mathrm{V_{\R}}$)} or of type {\rm (N)}. \hfill$\square$
\end{thm}

Notice that the theorem implies, in particular, that the valuative space $\ZR{\mathscr{X}}$ is reflexive\index{valuative!valuative topological space@--- (topological) space!reflexive valuative topological space@reflexive --- ---} ({\bf \ref{ch-pre}}.\ref{dfn-reflexivevaluativespaces}).

\begin{prop}\label{prop-spectralfunctorseparatedquotients2}
Let $\mathscr{C}$ be a category of rigid spaces.
Suppose that $\mathscr{X}\mapsto [\mathscr{X}]$ is a spectral functor\index{functor!spectral functor@spectral ---}\index{spectral functor} on $\mathscr{C}$.
Then any prespectral functor\index{spectral functor!prespectral functor@pre{---}} $S$ on $\mathscr{C}$ is a spectral functor\index{functor!spectral functor@spectral ---}\index{spectral functor}. 
\end{prop}

\begin{proof}
Let $X$ be a distinguished formal model of a coherent rigid space $\mathscr{X}$ in $\mathscr{C}$, $Y$ a non-empty closed subscheme of $X$ defined by an admissible ideal of $X$, and $U=X\setminus Y$.
We need to show that $S(\mathscr{U})\neq S(\mathscr{X})$, where $\mathscr{U}=U^{\rig}$.
Note that the conditions (b) and (c) in \ref{cor-spectralfunctorseparatedquotientscor} are equivalent to each other also in our situation.
Since $[\mathscr{U}]\neq[\mathscr{X}]$ by the assumption, the tube open subset $\mathscr{T}=C_{Y|X}$ corresponding to $Y$ is non-empty.
This means that $S(\mathscr{T})$ is non-empty, since $S$ is a prespectral functor, and hence that
$S(\mathscr{U})\neq S(\mathscr{X})$ as $S(\mathscr{U})\cap S(\mathscr{T}) = S (\mathscr{U}\cap\mathscr{T})=\emptyset$.
\end{proof}

By \ref{thm-spectralfunctorseparatedquotients} and the proposition we have:
\begin{cor}\label{cor-spectralfunctorseparatedquotients2}
The functor $\mathscr{X}\mapsto\ZR{\mathscr{X}}^{\cl}$ defined on the category of rigid spaces of type {\rm ($\mathrm{V_{\R}}$)} is a spectral functor.
\end{cor}

The following statement, of which the type {\rm ($\mathrm{V_{\R}}$)} case is contained in \ref{thm-spectralfunctorseparatedquotients}, is an immediate corollary of \ref{cor-spectralfunctorseparatedquotientscor}:
\begin{cor}\label{cor-classicalpointsexist2}
Let $\mathscr{X}$ be a rigid space of type {\rm ($\mathrm{V_{\R}}$)}\index{rigid space!rigid space of typeV1@--- of type ($\mathrm{V_{\R}}$)} or of type {\rm (N)}\index{rigid space!rigid space of typeN@--- of type (N)}, and $\mathscr{U}$ a quasi-compact open subspace of $\mathscr{X}$.
If $\mathscr{U}\subsetneq\mathscr{X}$, then we have $\ZR{\mathscr{U}}^{\cl}\subsetneq\ZR{\mathscr{X}}^{\cl}$. \hfill$\square$
\end{cor}

\begin{cor}\label{cor-classicalpointsexist4}
Let $\mathscr{X}$ be a rigid space of type {\rm ($\mathrm{V_{\R}}$)}\index{rigid space!rigid space of typeV1@--- of type ($\mathrm{V_{\R}}$)} or of type {\rm (N)}\index{rigid space!rigid space of typeN@--- of type (N)}, and $\mathscr{U},\mathscr{V}$ quasi-compact open subspaces of $\mathscr{X}$.
Then $\ZR{\mathscr{U}}^{\cl}=\ZR{\mathscr{V}}^{\cl}$ implies $\mathscr{U}=\mathscr{V}$.
\end{cor}

\begin{proof}
Apply \ref{cor-classicalpointsexist2} to $\mathscr{U}\subseteq\mathscr{U}\cup\mathscr{V}$.
\end{proof}

\subsection{Noetherness theorem}\label{sub-locringclassical}
\subsubsection{Comparison of complete local rings}\label{subsub-locringclassicalcomp}
\begin{prop}\label{prop-complocringclassicalcomp}
Let $\mathscr{X}=(\Spf A)^{\rig}$ be an affinoid of type {\rm (V)} or {\rm (N)}, and $x\in\ZR{\mathscr{X}}$ a closed classical point.
Then for any $n\geq 1$ we have the canonical isomorphism
$$
\O_{\mathscr{X},x}/\m^n_{\mathscr{X},x}\cong\O_{s(\mathscr{X}),s(x)}/\m^n_{s(\mathscr{X}),s(x)}.
$$
In particular, we have
$$
\widehat{\O}_{\mathscr{X},x}\cong\widehat{\O}_{s(\mathscr{X}),s(x)},
$$
where $\widehat{\cdot}$ denotes the completions with respect to the maximal ideals.
\end{prop}

For the proof, we refer to \cite[(7.3.2/3)]{BGR}.
Here we offer a non-elementary proof, based on GAGA.
\begin{proof}
We may assume that $A$ is $a$-torsion free.
The closed classical point $x$ is a closed subspace $\mathscr{Z}\hookrightarrow\mathscr{X}$, and let $\mathscr{J}_{\mathscr{Z}}\subseteq\O_{\mathscr{X}}$ be the coherent ideal that defines $\mathscr{Z}$.
Take the corresponding coherent ideal $\mathscr{J}_{s(\mathscr{Z})}\subseteq\O_{s(\mathscr{X})}$ by the correspondence as in \ref{prop-comparisonmapassaffprop1}, which is the defining ideal of the closed point $s(\mathscr{Z})=s(x)$ on $s(\mathscr{X})$; see the proof of \ref{prop-classicalpointsexist} (1).
Then for any $n\geq 1$, we have
\begin{equation*}
\begin{split}
\O_{s(\mathscr{X}),s(x)}/\m^n_{s(\mathscr{X}),s(x)}&\cong\Gamma(s(\mathscr{X}),\O_{s(\mathscr{X})}/\mathscr{J}^n_{s(\mathscr{Z})})\\ &\cong\Gamma(\mathscr{X},\O_{\mathscr{X}}/\mathscr{J}^n_{\mathscr{Z}})\cong\O_{\mathscr{X},x}/\m^n_{\mathscr{X},x},
\end{split}
\end{equation*}
as desired.
\end{proof}

\begin{cor}\label{cor-complocringclassicalcomp}
Let $V$ be an $a$-adically complete valuation ring of height one, and $\mathscr{Y}\hookrightarrow\mathscr{X}$ an open immersion of affinoids of finite type over $(\Spf V)^{\rig}$.
Then the associated morphism $s(\mathscr{Y})\rightarrow s(\mathscr{X})$ of Noetherian schemes is flat, and maps the set of all closed points of $s(\mathscr{Y})$ injectively to the set of all closed points of $s(\mathscr{X})$.
Moreover, for any closed point $y\in s(\mathscr{Y})$ we have
$$
\widehat{\O}_{s(\mathscr{Y}),y}\cong\widehat{\O}_{s(\mathscr{X}),x},
$$
where $x$ is the image of $y$.
\end{cor}

\begin{proof}
The flatness is already shown in \ref{prop-associatedschemesopenimm} (1).
Closed points of $s(\mathscr{X})$ are in one to one correspondence with classical points of $\mathscr{X}$ due to \ref{cor-classicalpointsexist0}; the restrictions of them on $\mathscr{Y}$ are either empty or classical points of $\mathscr{Y}$, and hence the fiber of each closed point of $s(\mathscr{X})$ by the map $s(\mathscr{Y})\rightarrow s(\mathscr{X})$ consists at most one closed point (cf.\ \ref{prop-classicalpointsexist} (2)).
The other assertion follows immediately from \ref{prop-complocringclassicalcomp}.
\end{proof}

\subsubsection{Reducedness and irreducibility}\label{subsub-closedimmrigidred}
\begin{prop}\label{prop-reducedmodelrigid}
Let $\mathscr{X}$ be a rigid space of type {\rm (V)}, and $\mathscr{N}_{\mathscr{X}}$ the subsheaf of $\O_{\mathscr{X}}$ consisting of locally nilpotent sections.
Then $\mathscr{N}_{\mathscr{X}}$ is a coherent ideal of $\O_{\mathscr{X}}$.
\end{prop}

\begin{proof}
Considering affinoid coverings of $\mathscr{X}$, we may assume that $\mathscr{X}$ is an affinoid $\mathscr{X}=(\Spf A)^{\rig}$, where $A$ is $a$-torsion free and topologically of finite type algebra over an $a$-adically complete valuation ring $V$.
Considering the base change by the height one localization $V_{\mathfrak{p}}$ of $V$, where $\mathfrak{p}=\sqrt{(a)}$, one can assume without loss of generality that $V$ is of height one.
Let $N\subseteq A[\frac{1}{a}]$ be the nilpotent radical of the Noetherian ring $A[\frac{1}{a}]$, and consider the corresponding coherent ideal $\mathscr{N}\subseteq\O_{\mathscr{X}}$ on $\mathscr{X}=(\Spf A)^{\rig}$; see \S\ref{sub-affinoidscohsheaf}.
Clearly we have $\mathscr{N}\subseteq\mathscr{N}_{\mathscr{X}}$.
We want to show $\mathscr{N}=\mathscr{N}_{\mathscr{X}}$.
Considering $A/A\cap N$, which is again an $a$-torsion free topologically of finite type $V$-algebra, one deduces the other inclusion from the following statement: if $A[\frac{1}{a}]$ is reduced (that is, $N=0$), then $\mathscr{X}$ is reduced, that is, $\O_{\mathscr{X},x}$ is reduced for any $x\in\ZR{\mathscr{X}}$.
To show this, due to \ref{prop-descriptionlocalrings}, it suffices to show that, for any affinoid subdomain $\mathscr{U}=(\Spf B)^{\rig}\subseteq\mathscr{X}$, the ring $B[\frac{1}{a}]$ is reduced.
But for this, we only need to show that the local ring of $U=s(\mathscr{U})=\Spec B[\frac{1}{a}]$ at every closed point is reduced.
For any $y\in U^{\cl}$, the morphism $\O_{U,y}\rightarrow\widehat{\O}_{U,y}$ is faithfully flat, and hence it suffices to show that $\widehat{\O}_{U,y}$ is reduced.
Then, in view of \ref{prop-complocringclassicalcomp}, after all, what we need to show is: if $A[\frac{1}{a}]$ is reduced, then, for any classical point $x\in\ZR{\mathscr{X}}^{\cl}$ the complete local ring $\widehat{\O}_{\mathscr{X},x}$ is reduced.
But, again due to \ref{prop-complocringclassicalcomp}, this follows from the reducedness of $\widehat{\O}_{X,s(x)}$, where $X=s(\mathscr{X})=\Spec A[\frac{1}{a}]$, which is guaranteed by the fact that the classical affinoid algebra $A[\frac{1}{a}]$ is excellent ({\bf \ref{ch-pre}}.\ref{thm-propclassicalaffinoidjapanese}).
\end{proof}

By \ref{prop-closedimmrigid2} we have the closed immersion 
$$
\mathscr{X}_{\red}\longhookrightarrow\mathscr{X}
$$
with the defining ideal $\mathscr{N}_{\mathscr{X}}$.
The rigid space $\mathscr{X}_{\red}$ is determined up to canonical isomorphisms.
We call $\mathscr{X}_{\red}$ the {\it reduced model}\index{reduced model of a rigid space@reduced model (of a rigid space)}\index{rigid space!reduced model of a rigid space@reduced model of a ---} of $\mathscr{X}$.
By \ref{prop-ZRpoints3} (1) we deduce that the topological spaces $\ZR{\mathscr{X}}$ and $\ZR{\mathscr{X}_{\red}}$ coincide.

Let $\mathscr{X}$ be a rigid space of type {\rm (V)}, $\mathscr{J}$ a coherent ideal of $\O_{\mathscr{X}}$, and $\mathscr{Z}\subseteq\mathscr{X}$ the closed subspace defined by $\mathscr{J}$.
Applying \ref{prop-reducedmodelrigid} to $\mathscr{Y}$, one obtains a coherent ideal of $\O_{\mathscr{X}}$ containing $\mathscr{J}$ of which the stalk at any $x\in\ZR{\mathscr{X}}$ coincides with $\sqrt{\mathscr{J}_x}$.
We denote this sheaf by $\sqrt{\mathscr{J}}$.
This is a coherent ideal sheaf of $\O_X$, of which the corresponding closed subspace is $\mathscr{Z}_{\red}$.
In case $\mathscr{X}$ is an affinoid, say $\mathscr{X}=(\Spf A)^{\rig}$, then, if $\mathscr{I}$ corresponds to an ideal $I\subseteq A[\frac{1}{a}]$ as in \S\ref{sub-affinoidscohsheaf} (that is, $\mathscr{I}=I\O_{\mathscr{X}}$), then $\sqrt{\mathscr{I}}$ corresponds to $\sqrt{I}$.

\begin{cor}\label{cor-classicalpointsexist00}
Let $\mathscr{X}=(\Spf A)^{\rig}$ be an affinoid of type {\rm ($\mathrm{V_{\R}}$)}, and $\mathscr{Z}_1,\mathscr{Z}_2$ closed subspaces of $\mathscr{X}$.
Let $\mathscr{J}_1$ and $\mathscr{J}_2$ be the defining ideal of $\mathscr{Z}_1$ and $\mathscr{Z}_2$, respectively.
Then $\ZR{\mathscr{Z}_1}=\ZR{\mathscr{Z}_2}$ if and only if $\sqrt{\mathscr{J}_1}=\sqrt{\mathscr{J}_2}$.
\end{cor}

\begin{proof}
The `if' part is easy.
To show the `only if' part, let $Z_1$ and $Z_2$ be the closed subschemes of $s(\mathscr{X})=\Spec A[\frac{1}{a}]$ corresponding to $\mathscr{Z}_1$ and $\mathscr{Z}_2$, respectively, by the correspondence as in \ref{prop-closedsubspaceaffinoid2}.
If $\ZR{\mathscr{Z}_1}=\ZR{\mathscr{Z}_2}$, we have $\ZR{\mathscr{Z}_1}^{\cl}=\ZR{\mathscr{Z}_2}^{\cl}$, and hence the closed points in $Z_1$ coincides with those of $Z_2$ (\ref{cor-classicalpointsexist0}).
Since the classical affinoid algebra $A[\frac{1}{a}]$ is Jacobson ({\bf \ref{ch-pre}}.\ref{prop-classicalaffringjacobson}), this means, if $I_1$ and $I_2$ are the ideals of $A[\frac{1}{a}]$ defining $Z_1$ and $Z_2$, respectively, we have $\sqrt{I_1}=\sqrt{I_2}$.
Since $\mathscr{J}_i=I_i\O_{\mathscr{X}}$ for $i=1,2$, we have $\sqrt{\mathscr{J}_1}=\sqrt{\mathscr{J}_2}$.
\end{proof}

\begin{prop}\label{prop-reducedmodelrigid2}
Let $V$ be an $a$-adically complete valuation ring of height one, and $A$ a topologically finitely generated $V$-algebra.

{\rm (1)} The affinoid $\mathscr{X}=(\Spf A)^{\rig}$ is reduced\index{rigid space!reduced rigid space@reduced ---} if and only if so is the Noetherian scheme $s(\mathscr{X})=\Spec A[\frac{1}{a}]$.

{\rm (2)} The affinoid $\mathscr{X}=(\Spf A)^{\rig}$ is irreducible {\rm (\ref{dfn-irreduciblerigidspace})}\index{rigid space!irreducible rigid space@irreducible ---} if and only if so is the Noetherian scheme $s(\mathscr{X})=\Spec A[\frac{1}{a}]$.
\end{prop}

\begin{proof}
(1) Set $B=A[\frac{1}{a}]$ $(=\Gamma(\mathscr{X},\O_{\mathscr{X}}))$, and let $N$ be the nilpotent radical of the ring $B$.
As the proof of \ref{prop-reducedmodelrigid} indicates, $\mathscr{X}$ is reduced if and only if $N=0$, whence (1).

(2) We may assume $s(\mathscr{X})=\Spec A[\frac{1}{a}]$ is reduced, or what amounts to the same, the ring $A[\frac{1}{a}]$ is reduced.
Consider reduced closed rigid subspaces $\mathscr{Y},\mathscr{Z}\subseteq\mathscr{X}$, and the corresponding reduced closed subschemes $s(\mathscr{Y}),s(\mathscr{Z})$ of $s(\mathscr{X})$ (cf.\ \ref{prop-closedsubspaceaffinoid2}).
It follows from \ref{cor-classicalpointsexist00} that $\ZR{\mathscr{X}}=\ZR{\mathscr{Y}}\cup\ZR{\mathscr{Z}}$ if and only if $s(\mathscr{X})=s(\mathscr{Y})\cup s(\mathscr{Z})$, and similarly that, for example, $\ZR{\mathscr{Z}}=\ZR{\mathscr{X}}$ holds if and only if $s(\mathscr{X})=s(\mathscr{Z})$.
Hence we deduce that $\mathscr{X}$ is irreducible if and only if the scheme $s(\mathscr{X})$ is irreducible.
\end{proof}

Let $A$ and $\mathscr{X}$ be as in \ref{prop-reducedmodelrigid2}, and 
$$
s(\mathscr{X})=\bigcup^r_{i=1}X_i
$$
be the irreducible decomposition of the Noetherian scheme $s(\mathscr{X})=\Spec A[\frac{1}{a}]$.
For each $i=1,\ldots,r$, we have the uniquely determined closed rigid subspace $\mathscr{X}_i$ of $\mathscr{X}$ such that $s(\mathscr{X}_i)=X_i$ (cf.\ \ref{prop-comparisonmapassaffprop1}); in fact, if $\mathfrak{q}_i$ is the ideal of $B=A[\frac{1}{a}]$ corresponding to $X_i$, then $\mathscr{X}_i=(\Spf A/\mathfrak{q}_i\cap A)^{\rig}$.
By \ref{prop-reducedmodelrigid2}, each $Z_i$ ($i=1,\ldots,r$) is irreducible, and we have the irreducible decomposition
$$
\ZR{\mathscr{X}}=\bigcup^r_{i=1}\ZR{\mathscr{X}_i}.
$$

\subsubsection{Noetherness theorem}\label{subsub-noetherness}
\begin{thm}\label{thm-noetherness}
Let $\mathscr{X}$ be a rigid space of type {\rm (V)} or of type {\rm (N)}, and $x\in\ZR{\mathscr{X}}$ a point.
Then the local ring $\O_{\mathscr{X},x}$ is Noetherian.
\end{thm}

\begin{proof}[Proof of Theorem {\rm \ref{thm-noetherness}} in type {\rm (V)} case]
We may assume that $\mathscr{X}$ is an affinoid.
First we show the assertion in the following case: $\mathscr{X}$ is of finite type over $(\Spf V)^{\rig}$, where $V$ is an $a$-adically complete valuation ring of height one, and $x\in\ZR{\mathscr{X}}$ is a classical point.
We use the notation as in the proof of \ref{prop-complocringclassicalcomp} with $R_{\alpha}=\O_{s(\mathscr{U}_{\alpha}),x_{\alpha}}$ for $\alpha\in L$.
By \ref{cor-comparisonmapassaffloc} we only need to check the condition (c) in {\bf \ref{ch-pre}}.\ref{prop-nagataquintessence}.
As we have shown in the proof of \ref{prop-complocringclassicalcomp}, we have $R_{\alpha}/\m^n_{\alpha}\cong R_{\beta}/\m^n_{\beta}$ for any $\alpha\leq\beta$ and $n\geq 1$.
Hence we have $\m_{\alpha}R_{\beta}+\m^n_{\beta}=\m_{\beta}$ for any $n\geq 1$.
But since $R_{\beta}$ is Noetherian, every ideal is closed, and hence we have $\m_{\alpha}R_{\beta}=\m_{\beta}$, as desired.

To show the assertion in type (V) case in general, notice first that by \ref{thm-fibersoverrigptsbehavior} (3), {\bf \ref{ch-pre}}.\ref{cor-compval2006ver23}, and {\bf \ref{ch-pre}}.\ref{prop-convpreadhloc1}, we may assume that $V$ is of height one.
Consider for any point $x\in\ZR{\mathscr{X}}$ the associated rigid point\index{point!rigid point@rigid ---}\index{rigid point} $(\Spf V')^{\rig}\rightarrow\mathscr{X}$ (\ref{dfn-ZRpoints32}).
By \ref{thm-fibersoverrigptsbehavior} (3) we may assume that $x$ is of height one and hence that $V'$ is of height one.
Consider the base-change $\mathscr{X}'=\mathscr{X}\times_{(\Spf V)^{\rig}}(\Spf V')^{\rig}$ and the induced diagram
$$
\xymatrix@C-1ex{\mathscr{X}\ar[d]&\mathscr{X}'\ar[d]\ar[l]\\ (\Spf V)^{\rig}&(\Spf V')^{\rig}\rlap{.}\ar[l]\ar@/_1pc/[u]}
$$
The rigid point $(\Spf V')^{\rig}\rightarrow\mathscr{X}'$ determines a classical point $x'$.
By the special case treated above, we know that $\O_{\mathscr{X}',x'}$ is Noetherian.
Since $V'$ is faithfully flat over $V$, the map $\O_{\mathscr{X},x}\rightarrow\O_{\mathscr{X}',x'}$ is faithfully flat, and hence $\O_{\mathscr{X},x}$ is Noetherian, as desired.
\end{proof}

\begin{proof}[Proof of Theorem {\rm \ref{thm-noetherness}} in type {\rm (N)} case]
We may assume that $\mathscr{X}$ is an affinoid $\mathscr{X}=(\Spf A)^{\rig}$, where $A$ is an $I$-adically complete Noetherian ring by an ideal $I\subseteq A$.
Take a rigid point $(\Spf V)^{\rig}\rightarrow\mathscr{X}$, where $V$ is an $a$-adically complete valuation ring ($a\in\m_V\setminus\{0\}$), that maps the closed point to the point $x\in\ZR{\mathscr{X}}$.
Take a cofinal system of formal neighborhoods $\{\mathscr{U}_{\alpha}=(\Spf A_{\alpha})^{\rig}\}_{\alpha\in L}$ of $x$ in such a way that each $\Spf A_{\alpha}$ lies in an admissible blow-up of $\Spf A$ as an affine open subspace; notice that each $A_{\alpha}$ is Noetherian.
Let $x_{\alpha}$ for each $\alpha\in L$ be the image of $x$ by the specialization map $\ZR{\mathscr{X}}\rightarrow\Spf A_{\alpha}$, that is, the image of the closed point by $\Spec V\rightarrow\Spec A_{\alpha}$.
Furthermore, let $\xi_{\alpha}\in\Spec A_{\alpha}$ be the image of the generic point by the morphism $\Spec V\rightarrow\Spec A_{\alpha}$, that is, the image of $x_{\alpha}$ by the map $s\colon\ZR{(\Spf A_{\alpha})^{\rig}}\rightarrow\Spec A_{\alpha}\setminus V(IA_{\alpha})$.
We have $\O_{\mathscr{X},x}=\varinjlim_{\alpha\in L}\O_{\Spec A_{\alpha},\xi_{\alpha}}$, which we need to show to be Noetherian.
To this end, set $R_{\alpha}=\O_{\Spec A_{\alpha},x_{\alpha}}$; we denote also by $\xi_{\alpha}$ the image of $\xi_{\alpha}$ in $\Spec R_{\alpha}$.

Consider the completion $\widehat{R}_{\alpha}$ with respect to the maximal ideal; each $R_{\alpha}\rightarrow\widehat{R}_{\alpha}$ is faithfully flat, and for each $\alpha\leq\beta$ we have the induced local map $\widehat{R}_{\alpha}\rightarrow\widehat{R}_{\beta}$.
Notice that each $\widehat{R}_{\alpha}$ is quasi-excellent, since it is a Noetherian complete local ring.
Set $B_{\alpha}=\O_{\Spec\widehat{R}_{\alpha},\widehat{\xi}_{\alpha}}$ (where $\widehat{\xi}_{\alpha}\in\Spec\widehat{R}_{\alpha}$ is the point above $\xi_{\alpha}\in\Spec R_{\alpha}$); for $\alpha\leq\beta$ we have the induced local morphisms $B_{\alpha}\rightarrow B_{\beta}$.

Now for $\alpha\leq\beta$ the transition map $A_{\alpha}\rightarrow A_{\beta}$ comes from a Zariski open part of an admissible blow-up, hence so is $R_{\alpha}\rightarrow R_{\beta}$.
Since $R_{\alpha}\rightarrow\widehat{R}_{\alpha}$ is flat, $\widehat{R}_{\alpha}\rightarrow\widehat{R}_{\beta}$ is a composition of admissible blow-up (with respect to $I$-adic topology) followed by completion with respect to a maximal ideal.
Since $\widehat{R}_{\alpha}$ is quasi-excelent, the formal fiber $B_{\alpha}\rightarrow B_{\beta}$ is regular.
Now we fix an $\alpha_0\in L$ and replace $L$ by the cofinal subset $\{\alpha\in L\,|\,\alpha\geq\alpha_0\}$.
Then $\{\dim B_{\alpha}\}_{\alpha\in L}$ is upper-bounded by $\dim A_{\alpha_0}$.
Hence by {\bf \ref{ch-pre}}.\ref{cor-propnagataquintessencevar} we deduce that $B=\varinjlim_{\alpha\in L}B_{\alpha}$ is Noetherian.
Now since any $\O_{\Spec A_{\alpha},\xi_{\alpha}}\rightarrow B_{\alpha}$ is faithfully flat, we deduce that $\O_{\mathscr{X},x}\rightarrow B$ is faithfully flat; since $B$ is known to be Noetherian, so is $\O_{\mathscr{X},x}$, as desired.
\end{proof}



\section{GAGA}\label{sec-GAGA}
\index{GAGA|(}
In this section we discuss GAGA theorems in rigid geometry.
The first subsection (\S\ref{sub-GAGAfunctor}) is devoted to the definition of {\em GAGA functor}, which associates to any separated of finite type scheme $X$ over $U=\Spec A\setminus V(I)$, where $A$ is an adic ring with a finitely generated ideal of definition $I\subseteq A$, the `analytification' $X^{\an}$, a separated rigid space over the affinoid $(\Spf A)^{\rig}$.
Some of the basic properties of the GAGA functor will be discussed in \S\ref{subsub-GAGAfunctorprop}.
We also give a generalization of the GAGA functor to non-separated schemes (\S\ref{subsub-GAGAfunctornonsep}).

After discussing affinoid valued points (\S\ref{sub-affinoidvaluedpoints}), we introduce the so-called comparison maps and comparison functors in \S\ref{sub-GAGAcommap}, by which the GAGA theorems, GAGA comparison theorem and GAGA existence theorem, are formulated and proved in \S\ref{sub-GAGAcomprigid} and \S\ref{sub-GAGAexistrigid} from GFGA theorems.
Similarly to the GFGA theorems ({\bf \ref{ch-formal}}, \S\ref{sec-GFGAcom} and \S\ref{sec-GFGAexist}), our GAGA theorems are stated in terms of derived categorical language.

\subsection{Construction of GAGA functor}\label{sub-GAGAfunctor}
\subsubsection{The category $\Emb_{X|S}$}\label{subsub-GAGAfunctoremb}
Let $A$ be an adic ring of finite ideal type\index{adic!adic ring@--- ring!adic ring of finite ideal type@--- --- of finite ideal type}, and $I\subseteq A$ a finitely generated ideal of definition.
Set
$$
S=\Spec A\longhookleftarrow U=S\setminus D,
$$
where $D=V(I)$ is the closed subset corresponding to $I$, and 
$$
\mathscr{S}=(\Spf A)^{\rig}.
$$

Let $f\colon X\rightarrow U$ be a separated $U$-scheme of finite type.
We define the category $\Emb_{X|S}$ as follows:
\begin{itemize}
\item[$\bullet$] the objects are the commutative diagrams 
$$
\xymatrix{X\ \ar@{^{(}->}[r]\ar[d]_f&\ovl{X}\ar[d]^{\ovl{f}}\\ U\ \ar@{^{(}->}[r]&S\rlap{,}}
$$
where $\ovl{f}\colon\ovl{X}\rightarrow S$ is a proper $S$-scheme, and $X\hookrightarrow\ovl{X}$ is a birational open immersion\index{birational!birational open immersion@--- open immersion}, that is, an open immersion onto a dense open subspace of $X$ (cf.\ \ref{dfn-birationalgeom1} (3) below); 
\item[$\bullet$] a morphism $(X\hookrightarrow\ovl{X})\rightarrow(X\hookrightarrow\ovl{X}')$ is an $X$-admissible $S$-modification\index{admissible!U-admissible modification@$U$-{---} modification}\index{modification!U-admissible modification@$U$-admissible ---} $\ovl{X}\rightarrow\ovl{X}'$ (cf.\ \ref{dfn-birationalgeom1}).
\end{itemize}

Note that the category $\Emb_{X|S}$ is non-empty due to Nagata's\index{Nagata, M.} embedding theorem\index{Nagata's embedding theorem} (\ref{thm-nagataembedding}).

\begin{lem}\label{lem-GAGAfunctor1}
{\rm (1)} The category $\Emb_{X|S}$ is cofiltered.

{\rm (2)} Let $\ovl{X}=(X\hookrightarrow\ovl{X},\ovl{f}\colon\ovl{X}\rightarrow S)$ be an object of $\Emb_{X|S}$, and $\AId^{\ast}_{(\ovl{X},X)}$ be the set of all quasi-coherent ideals of $\O_{\ovl{X}}$ of finite type such that the corresponding closed subscheme is disjoint from $X$ $($cf.\ {\rm \S\ref{subsub-birationalgeomblowups}}$)$.
We introduce an ordering in the set $\AId^{\ast}_{(\ovl{X},X)}$ as follows$:$ $\mathscr{J}\preceq\mathscr{J}'$ if and only if there exists $\mathscr{J}''\in\AId^{\ast}_{(\ovl{X},X)}$ such that $\mathscr{J}=\mathscr{J}'\mathscr{J}''$.
Then $\AId^{\ast\opp}_{(\ovl{X},X)}$ is a directed set, and the functor 
$$
\AId^{\ast}_{(\ovl{X},X)}\longrightarrow\Emb_{X|S}
$$
that maps $\mathscr{J}$ to $\ovl{X}_{\mathscr{J}}=(X\hookrightarrow\ovl{X}_{\mathscr{J}}, \ovl{f}_{\mathscr{J}}\colon\ovl{X}_{\mathscr{J}}\rightarrow S)$, where $X_{\mathscr{J}}\rightarrow X$ is the blow-up along $\mathscr{J}$, is cofinal.
\end{lem}

\begin{proof}
(1) What to show are the following:
\begin{itemize}
\item[(a)] for $\ovl{X}_1=(X\hookrightarrow\ovl{X}_1, \ovl{f}_1\colon\ovl{X}_1\rightarrow S)$ and $\ovl{X}_2=(X\hookrightarrow\ovl{X}_2, \ovl{f}_2\colon\ovl{X}_2\rightarrow S)$, there exists an object $\ovl{X}_3=(X\hookrightarrow\ovl{X}_3, \ovl{f}_3\colon\ovl{X}_3\rightarrow S)$ and $X$-admissible $S$-modifications $\ovl{X}_3\rightarrow\ovl{X}_1$ and $\ovl{X}_3\rightarrow\ovl{X}_2$; 
\item[(b)] for $\ovl{X}_1=(X\hookrightarrow\ovl{X}_1, \ovl{f}_1\colon\ovl{X}_1\rightarrow S)$ and $\ovl{X}_2=(X\hookrightarrow\ovl{X}_2, \ovl{f}_2\colon\ovl{X}_2\rightarrow S)$ and two $X$-admissible $S$-modifications $q_0,q_1\colon \ovl{X}_2\rightarrow\ovl{X}_1$, there exists an object 
$\ovl{X}_3=(X\hookrightarrow\ovl{X}_3, \ovl{f}_3\colon\ovl{X}_3\rightarrow S)$ and an $X$-admissible $S$-modification $p\colon \ovl{X}_3\rightarrow\ovl{X}_2$ such that $q_0\circ p=q_1\circ p$.
\end{itemize}

To show (a), let $\ovl{X}_3=\ovl{X}_1\ast\ovl{X}_2$, the join\index{join} of $\ovl{X}_1$ and $\ovl{X}_2$ (cf.\ \ref{dfn-correspondencediagram}); that is, the closure of $X$ in the product $\ovl{X}_1\times_S\ovl{X}_2$.
Then $\ovl{X}_3\rightarrow S$ is proper. 
Clearly, the morphisms by projections $\ovl{X}_3\rightarrow\ovl{X}_1$ and $\ovl{X}_3\rightarrow\ovl{X}_2$ are $X$-admissible modifications.

To show (b), consider the Cartesian diagram of $S$-schemes:
$$
\xymatrix{\ovl{X}_2\ar[r]^(.34){(q_0,q_1)}&\ovl{X}_1\times_S\ovl{X}_1\\ Z\ar[u]\ar[r]&\ovl{X}_1\ar[u]_(.47){\Delta_{\ovl{X}_1}}\rlap{.}}
$$
Since $\ovl{X}_1$ is separated over $S$, the right-hand vertical map is a closed immersion, and hence so is the left-hand one.
The scheme $Z$ contains a copy of $X$.
Let $\ovl{X}_3$ be the scheme-theoretic closure of $X$ in $Z$.
Then $\ovl{X}_3$ defines an object of $\Emb_{X|S}$, and the map $p\colon\ovl{X}_3\rightarrow\ovl{X}_2$ is clearly an $X$-admissible $S$-modification such that $q_0\circ p=q_1\circ p$.

(2) follows immediately from \cite[Premi\`ere partie, (5.7.12)]{RG} (see \ref{prop-birationalgeom02}).
\end{proof}

\subsubsection{Construction of $X^{\an}$}\label{subsub-GAGAfunctorconst}
\begin{const}\label{const-GAGAfunctor}{\rm 
We continue with working in the situation as in \S\ref{subsub-GAGAfunctoremb}.
Let $f\colon X\rightarrow U$ be a separated $U$-scheme of finite type, and 
$$
\xymatrix{X\ \ar@{^{(}->}[r]\ar[d]_f&\ovl{X}\ar[d]^{\ovl{f}}\\ U\ \ar@{^{(}->}[r]&S}
$$
be an object of $\Emb_{X|S}$.
Set
\begin{equation*}
\begin{split}
Z&=(\ovl{X}\times_SU)\setminus X,\\
\ovl{Z}&=\textrm{the closure of}\ Z\ \textrm{in}\ \ovl{X},\\
\til{X}&=\ovl{X}\setminus\ovl{Z}.
\end{split}
\end{equation*}
Let $\widehat{\til{X}}\hookrightarrow\widehat{\ovl{X}}$ be the $I$-adic completion of the open immersion $\til{X}\hookrightarrow\ovl{X}$.
We have the open immersion 
$$
(\widehat{\til{X}})^{\rig}\longhookrightarrow(\widehat{\ovl{X}})^{\rig}
$$
of coherent rigid spaces.}
\end{const}

\begin{dfn}\label{dfn-GAGAfunctor2}{\rm 
Let the notation be as in \ref{const-GAGAfunctor}.
We define a sheaf $X^{\an}$ on the site $\Rf_{\mathscr{S},\ad}$ by
$$
X^{\an}=\varinjlim(\widehat{\til{X}})^{\rig},
$$
where the inductive limit is taken along the filtered category $\Emb^{\opp}_{X|S}$ or, equivalently, along the set $\AId^{\ast}_{(\ovl{X},X)}$.}
\end{dfn}

\begin{prop}\label{prop-GAGAfunctor1}
The sheaf $X^{\an}$ on $\Rf_{\mathscr{S},\ad}$ is a quasi-separated rigid space.
\end{prop}

\begin{proof}
Let $\ovl{X}=(X\hookrightarrow\ovl{X}, \ovl{f}\colon\ovl{X}\rightarrow S)$ be an object of $\Emb_{X|S}$, and consider an $X$-admissible blow-up $\ovl{X}_1\rightarrow\ovl{X}$, that is, a blow-up along an ideal in $\AId^{\ast}_{(\ovl{X},X)}$.
Then the induced morphism $\widehat{\ovl{X}}_1\rightarrow\widehat{\ovl{X}}$ is an admissible blow-up of coherent adic formal schemes of finite ideal type, and thus we have the canonical isomorphism $(\widehat{\ovl{X}}_1)^{\rig}\cong(\widehat{\ovl{X}})^{\rig}$.
We deduce by this that the corresponding $(\widehat{\til{X}}_1)^{\rig}$ sits in the following commutative diagram consisting of open immersions:
$$
\xymatrix@-3ex{
&(\widehat{\til{X}}_1)^{\rig}\ar@{^{(}->}[dr]\\
(\widehat{\til{X}})^{\rig}\ar@{^{(}->}[ur]\ar@{^{(}->}[rr]&&(\widehat{\ovl{X}})^{\rig}\rlap{.}}\\
$$
As the set $\AId^{\ast\opp}_{(\ovl{X},X)}$ is cofinal in $\Emb_{X|S}$, we deduce that $X^{\an}$ is a stretch of coherent rigid spaces and hence is quasi-separated (\ref{prop-generalrigidspace31-1}).
\end{proof}

Notice that by the construction we always have a canonical open immersion 
$$
(\widehat{\til{X}})^{\rig}\longhookrightarrow X^{\an}
$$
for any object $(X\hookrightarrow \ovl{X})$ of $\Emb_{X|S}$.
\begin{prop}\label{prop-GAGAfunctor4}
If $f\colon X\rightarrow U$ is proper, then we have $X^{\an}=(\widehat{\ovl{X}})^{\rig}$.
In particular, $X^{\an}$ is a coherent rigid space.
\end{prop}

\begin{proof}
Take arbitrary $(X\hookrightarrow\ovl{X})\in\obj(\Emb_{X|S})$, and let $\ovl{X}'$ be the closure of $X$ in $\ovl{X}$.
Then by \cite[Premi\`ere partie, (5.7.12)]{RG} (see \ref{prop-birationalgeom02} below) one sees that $(\widehat{\ovl{X}})^{\rig}\cong(\widehat{\ovl{X}'})^{\rig}$.
When we replace $\ovl{X}$ by $\ovl{X}'$, then we have $\ovl{Z}=\emptyset$, thereby the result.
\end{proof}

In particular, we have $U^{\an}=\mathscr{S}$ $(=(\Spf A)^{\rig})$.

\begin{const}\label{const-GAGAfunctor2}{\rm 
Next we construct, for a given $U$-morphism $h\colon Y\rightarrow X$ of $U$-schemes of finite type, the morphism of rigid spaces 
$$
h^{\an}\colon Y^{\an}\longrightarrow X^{\an}.
$$
Let us take arbitrary $(X\hookrightarrow\ovl{X})\in\obj(\Emb_{X|S})$ and $(Y\hookrightarrow\ovl{Y})\in\obj(\Emb_{Y|S})$.
Let $\ovl{Y}_1$ be the closure of the graph of $h$ in the product $\ovl{Y}\times_S\ovl{X}$; $\ovl{Y}_1$ is proper and admits an birational open immersion $Y\hookrightarrow\ovl{Y}_1$, that is, $(Y\hookrightarrow\ovl{Y}_1)$ defines an object of $\Emb_{Y|S}$ that sits in the commutative square
$$
\xymatrix{
Y\ar[r]\ar@{^{(}->}(-.5,-3);(-.5,-11)&X\ar@{^{(}->}(14.7,-3);(14.7,-11)\\ \ovl{Y}_1\ar[r]&\ovl{X}\rlap{.}}
$$
Thus we get $\til{Y}\rightarrow\til{X}$.
As $(X\hookrightarrow\ovl{X})\in\obj(\Emb_{X|S})$ and $(Y\hookrightarrow\ovl{Y})\in\obj(\Emb_{Y|S})$ are taken arbitrary, we get the desired morphism $h^{\an}$.}
\end{const}

Thus we get a functor 
$$
X\longmapsto X^{\an}
$$
from the category of separated $U$-schemes of finite type to the category of rigid spaces over $\mathscr{S}$, called the {\em GAGA functor}\index{GAGA!GAGA functor@--- functor}.
\begin{rem}\label{rem-algebraicspaceGAGA}{\rm 
Our construction of the GAGA functor relies on Nagata's\index{Nagata, M.} embedding theorem\index{Nagata's embedding theorem} \ref{thm-nagataembedding}.
The first author has proven a strong version of the Nagata's embedding theorem (\ref{thm-nagataembalgsp}) by which one can define the GAGA functor for algebraic space as follows.
In the situation as in \S\ref{subsub-GAGAfunctoremb}, we consider a separated $U$-algebraic space $f\colon X\rightarrow U$ of finite type.
We define the category $\Emb_{X|S}$ as before with the additional condition that the boundary $\ovl{X}\setminus X$ is a scheme (cf.\ \ref{thm-nagataembalgsp} (a)).
Now in \ref{const-GAGAfunctor}, since $\ovl{X}\times_SD$ and $\til{X}\times_SD$ are schemes, the formal completions $\widehat{\til{X}}$ and $\widehat{\ovl{X}}$ are formal schemes and hence define the rigid spaces $(\widehat{\til{X}})^{\rig}$ and $(\widehat{\ovl{X}})^{\rig}$.
Then one can just carry out the same construction to get the desired GAGA functor $X\mapsto X^{\an}$ from the category of separated $U$-algebraic spaces of finite type to the category of rigid spaces over $\mathscr{S}$.
See \cite{CT} for another approach to define GAGA functors for algebraic spaces.}
\end{rem}

Let us finally remark that, if the adic ring $A$ in \S\ref{subsub-GAGAfunctoremb} is t.u.\ rigid-Noetherian\index{t.u. rigid-Noetherian ring@t.u.\ rigid-Noetherian ring} (resp.\ t.u.\ adhesive\index{t.u.a. ring@t.u.\ adhesive ring}) ({\bf \ref{ch-formal}}.\ref{dfn-tuaringadmissible}), then the rigid space $X^{\an}$ is locally universally Noetherian\index{rigid space!universally Noetherian rigid space@universally Noetherian ---!locally universally Noetherian rigid space@locally --- ---} (resp.\ locally universally adhesive\index{rigid space!universally adhesive rigid space@universally adhesive ---!locally universally adhesive rigid space@locally --- ---}) (\ref{dfn-universallyadhesiverigidspaces}).

\subsubsection{Some basic properties of the GAGA functor}\label{subsub-GAGAfunctorprop}
\begin{prop}\label{prop-GAGAfunctor5x1}
The GAGA functor $X\mapsto X^{\an}$ maps a $U$-open immersion $V\hookrightarrow X$ to an $\mathscr{S}$-open immersion\index{immersion!open immersion of rigid spaces@open --- (of rigid spaces)} $V^{\an}\hookrightarrow X^{\an}$ and, in case the adic ring $A$ in {\rm \S\ref{subsub-GAGAfunctoremb}} is t.u.\ rigid-Noetherian\index{t.u. rigid-Noetherian ring@t.u.\ rigid-Noetherian ring}, a $U$-closed immersion $Y\hookrightarrow X$ to an $\mathscr{S}$-closed immersion\index{immersion!closed immersion of rigid spaces@closed --- (of rigid spaces)} $Y^{\an}\hookrightarrow X^{\an}$.
\end{prop}

\begin{proof}
First we prove that the functor $X\mapsto X^{\an}$ maps an open immersion to an open immersion.
As in \ref{const-GAGAfunctor2} one can consider morphisms of the form $\til{V}\rightarrow\til{X}$ starting from objects in $\Emb_{X|S}$ and $\Emb_{V|S}$.
By \cite[Premi\`ere partie, (5.7.11)]{RG} (included as \ref{prop-birationalgeom01} below) the morphism $\til{V}\rightarrow\til{X}$ can be replaced, by means of an $X$-admissible blow-up and the struct transform, by an open immersion.
Hence the $I$-adic completion $\widehat{\til{V}}\rightarrow\widehat{\til{X}}$ is an open immersion, and thus $(\widehat{\til{V}})^{\rig}\rightarrow(\widehat{\til{X}})^{\rig}$ is an open immersion of coherent rigid spaces.
As the rigid spaces $V^{\an}$ and $X^{\an}$ are unions of them, we deduce that the morphism $V^{\an}\hookrightarrow X^{\an}$ is an open immersion by the definition \ref{dfn-generalrigidspace4}.

Next we treat the assertion for a closed immersion $h\colon Y\hookrightarrow X$.
To show that $h^{\an}$ is a closed immersion, it suffices to show the following:
\begin{itemize}
\item[(a)] for any object $(X\hookrightarrow\ovl{X})$ of $\Emb_{X|S}$ the closure $\ovl{Y}$ of $Y$ in $\ovl{X}$ gives an object of $\Emb_{Y|S}$; 
\item[(b)] for any object $(Y\hookrightarrow\ovl{Y})$ of $\Emb_{Y|S}$ there exists an object $(X\hookrightarrow\ovl{X})$ of $\Emb_{X|S}$ such that the closure of $Y$ in $\ovl{X}$ dominates $\ovl{Y}$.
\end{itemize}
If these properties are satisfied, we get a cofinal family $\{\til{h}_{\lambda}\colon\til{Y}_{\lambda}\hookrightarrow\til{X}_{\lambda}\}$ of closed immersion over $S$ indexed by a directed set such that $\varinjlim_{\lambda}(\widehat{h}_{\lambda})^{\rig}=h^{\an}$, thereby the assertion.

The assertion (a) is clear; to see (b), take any $(X\hookrightarrow\ovl{X})$ and take the closure of the graph of the map $Y\hookrightarrow\ovl{X}$ in $\ovl{Y}\times_S\ovl{X}$, which clearly dominates $\ovl{Y}$.
\end{proof}

\begin{prop}\label{prop-GAGAfunctor5x2}
Suppose in the situation as in {\rm \S\ref{subsub-GAGAfunctoremb}} that $A$ is t.u.\ rigid-Noetherian.
Let $X$ be a separated $U$-scheme of finite type, and $Y$ a closed subscheme of $X$.
Then we have 
$$
(X\setminus Y)^{\an}=X^{\an}\setminus Y^{\an}
$$
$($cf.\ {\rm \ref{dfn-opencomplement}} for the notion of open complements of rigid spaces$)$.
\end{prop}

\begin{proof}
As in the proof of the second part of \ref{prop-GAGAfunctor5x1}, we have a cofinal family of $S$ closed immersions $\{\til{h}_{\lambda}\colon\til{Y}_{\lambda}\hookrightarrow\til{X}_{\lambda}\}$ such that $\varinjlim_{\lambda}(\widehat{h}_{\lambda})^{\rig}=h^{\an}$.
Then $\til{V}_{\lambda}=\til{X}_{\lambda}\setminus\til{Y}_{\lambda}$ gives the cofinal family of embeddings $V=X\setminus Y\hookrightarrow\til{V}$ as in \ref{const-GAGAfunctor}.
Hence we have $V^{\an}=\varinjlim_{\lambda}(\widehat{\til{V}}_{\lambda})^{\rig}=X^{\an}\setminus Y^{\an}$.
\end{proof}

\begin{prop}\label{prop-GAGAfunctor5x3}
The GAGA functor $X\mapsto X^{\an}$ is left-exact.
Moreover, it is compatible with base change$;$ that is, for any adic morphism $A\rightarrow A'$ with $S'=\Spec A'$ and $U'=S'\setminus V(IA')$ where $I\subseteq A$ is a finitely generated ideal of definition, and for any separated of finite type $U$-scheme $X$ $($resp.\ morphism $h$ between separated of finite type $U$-schemes$)$, we have $(X\times_UU')^{\an}\cong X^{\an}\times_{\mathscr{S}}\mathscr{S}'$ $($resp.\ $(h\times_UU')^{\an}\cong h^{\an}\times_{\mathscr{S}}\mathscr{S}')$, where $\mathscr{S}'=(\Spf A')^{\rig}$.
\end{prop}

\begin{proof}
To show that the GAGA functor is left-exact, it suffices to show that it preserves fiber products.
Consider a diagram $X\rightarrow Z\leftarrow Y$ of separated $U$-schemes of finite type.
We first suppose that $X,Y,Z$ are proper over $U$.
In this case, $(X\times_ZY)^{\an}$ is the associated rigid space of the formal completion of $\ovl{X}\times_{\ovl{Z}}\ovl{Y}$, where $\ovl{X}\rightarrow\ovl{Z}\leftarrow\ovl{Y}$ is a diagram consisting of Nagata compactifications.
Hence in this case, the assertion is clear.
In general, we first take Nagata compactifications $X$, $Y$, and $Z$ to embed them into a proper $U$-schemes and then apply \ref{thm-affinoidvaluedpoints} below to compare rigid points of $(X\times_ZY)^{\an}$ and those of $X^{\an}\times_{Z^{\an}}Y^{\an}$.
The compatibility with base change can be proved similarly.
\end{proof}

\begin{cor}\label{cor-GAGAfunctor51}
Suppose in the situation as in {\rm \S\ref{subsub-GAGAfunctoremb}} that $A$ is t.u.\ rigid-Noetherian.
Then the rigid space $X^{\an}$ for any separated $U$-scheme $X$ of finite type is separated.
\end{cor}

\begin{proof}
By \ref{prop-GAGAfunctor5x3} and \ref{prop-GAGAfunctor5x1} we know that the diagonal morphism $X^{\an}\rightarrow X^{\an}\times_{\mathscr{S}}X^{\an}=(X\times_UX)^{\an}$ is a closed immersion.
\end{proof}

\subsubsection{Some examples}\label{subsub-GAGAfunctorexa}
\begin{exas}\label{exas-gagm}{\rm 
(1) {\sl The additive group $\G_{\mathrm{a}}^{\an}$.} Consider the closed subscheme $D$ of $\P^1_U$ corresponding to the $\infty$-section, that is, $\P^1_U\setminus D=\G_{\mathrm{a},U}$.
Clearly, we have $D^{\an}\cong U^{\an}=\widehat{S}^{\rig}$.
Then by \ref{prop-GAGAfunctor5x2} we have
$$
\G_{\mathrm{a}}^{\an}=\P^{1,\an}_U\setminus D^{\an};
$$
notice that by \ref{prop-GAGAfunctor4} the projective space $\P^{1,\an}_U$ coincides with $\P^{1,\an}_{\mathscr{S}}=(\widehat{\P}^1_S)^{\rig}$ defined in \S\ref{subsub-examplesproj}.
Let $\mathscr{J}_D$ be the ideal defining $\widehat{D}$ (the $\infty$-section over $S$) in $\widehat{\P}^1_S$, and consider the admissible blow-up of $\widehat{\P}^1_S$ by the admissible ideal $\mathscr{J}_D+I^n\O_{\widehat{\P}^1_S}$ for each $n$.
Define $U_n$ to be the open part of the admissible blow-up where $\mathscr{J}_D$ generates the strict transform of $\mathscr{J}_D+I^n\O_{\widehat{\P}^1_S}$.
Then we have
$$
\G_{\mathrm{a}}^{\an}=\bigcup_{n\geq 0}U^{\rig}_n
$$
(cf.\ Exercise \ref{exer-GAGAfunctor1}).

(2) {\sl The multiplicative group $\G_{\mathrm{m}}^{\an}$.} 
This has a similar description as above where $D$ is replaced by the union of $\infty$-section and $0$-section.}
\end{exas}

\subsubsection{GAGA functor for non-separated schemes}\label{subsub-GAGAfunctornonsep}
Let $U$ and $S$ be as in \S\ref{subsub-GAGAfunctoremb}, and $X\rightarrow U$ a quasi-separated $U$-scheme locally of finite type.
(Notice that, in case the ring $A$ is t.u.\ rigid-Noetherian, any $U$-scheme of finite type $X$ is quasi-separated, since it is Noetherian.)
Let $X=\bigcup_{\alpha\in L}U_{\alpha}$ be a finite open covering by quasi-compact and separated schemes.
Then we define $X^{\an}$ to be the rigid space representing the sheaf on $\CRf_{S^{\an},\ad}$ that sits in the following exact sequence
$$\xymatrix@-1ex{
\coprod_{\alpha,\beta\in L}(U_{\alpha}\cap U_{\beta})^{\an}\ar@<.5ex>[r]\ar@<-.5ex>[r]&\coprod_{\alpha\in L}U^{\an}_{\alpha}\ar[r] &X^{\an}.}
$$
By \ref{prop-GAGAfunctor5x1} the rigid space $X^{\an}$ does not depend on the choice of the open covering.
By \ref{prop-GAGAfunctor5x1}, in case $X$ is separated, the above $X^{\an}$ coincides with the one defined in \ref{dfn-GAGAfunctor2}.

Similar construction can apply to define $h^{\an}$ for a morphism $h\colon X\rightarrow Y$ over $U$ between quasi-separated $U$-schemes of finite type.
Hence one get the extension of the GAGA functor\index{GAGA!GAGA functor@--- functor} defined on the category of all quasi-separated $U$-schemes of finite type.

\subsection{Affinoid valued points}\label{sub-affinoidvaluedpoints}
We continue with working in the situation as in \S\ref{subsub-GAGAfunctoremb}.
Let $\mathscr{T}$ be an affinoid\index{affinoid}, and suppose a morphism 
$$
\alpha\colon\mathscr{T}\longrightarrow X^{\an}
$$
of rigid spaces is given.
\begin{lem}\label{lem-affinoidvaluedpoints1}
There exists $\til{X}$ as in {\rm \ref{const-GAGAfunctor}} such that the morphism $\alpha$ factors through a morphism $\mathscr{T}\rightarrow(\widehat{\til{X}})^{\rig}$.
\end{lem}

\begin{proof}
Since $\mathscr{T}$ is quasi-compact, there exists a quasi-compact open subspace $\mathscr{U}$ of $X^{\an}$ such that $\alpha$ maps $\mathscr{T}$ to $\mathscr{U}$.
Since $X^{\an}$ is quasi-separated, $\mathscr{U}$ is coherent and hence is coherent in the small admissible site $(X^{\an})_{\ad}$.
By \cite[Expos\'e VI, 1.23]{SGA4-2} we find that there exists a $\til{X}$ as in \ref{const-GAGAfunctor} such that $\mathscr{U}\subseteq(\widehat{\til{X}})^{\rig}$, thereby the lemma.
\end{proof}

We henceforth denote the morphisms $\mathscr{T}\rightarrow(\widehat{\til{X}})^{\rig}$ and $\mathscr{T}\rightarrow(\widehat{\ovl{X}})^{\rig}$ also by $\alpha$.
We now assume that $\mathscr{T}$ is a universally adhesive affinoid of the form $\mathscr{T}=T^{\rig}$, where $T=\Spf B$; we may furthermore assume without loss of generality that $B$ is $J$-torsion free, where $J\subseteq B$ is an ideal of definition.

Consider the commutative diagram
$$
\xymatrix{
\mathscr{T}\ar[r]^{\alpha}\ar[dr]&(\widehat{\ovl{X}})^{\rig}\ar[d]^{(\widehat{\ovl{f}})^{\rig}}\\
&\widehat{S}^{\rig}\rlap{.}}
$$
There exists an admissible blow-up $T'\rightarrow T$ such that the diagram as above is induced from a commutative diagram of formal schemes of the following form:
$$
\xymatrix{
T'\ar[r]\ar[d]&\widehat{\ovl{X}}\ar[d]^{\widehat{\ovl{f}}}\\
T\ar[r]&\widehat{S}\rlap{,}}
$$
which induces
$$
\xymatrix{
T'\ar[r]\ar[dr]&\widehat{\ovl{X}}\times_{\widehat{S}}T\ar[d]\\
&T\rlap{.}}
$$
Since the morphisms in the last diagram are all proper (and of finite presentation, since $B$ is $J$-torsion free), by {\bf \ref{ch-formal}}.\ref{prop-GFGAexistapp1} we have
$$
\xymatrix{
Y'\ar[r]\ar[dr]&\ovl{X}\times_{S}Y\ar[d]\\
&Y\rlap{,}}
$$
where $Y=\Spec B$, and $Y'\rightarrow Y$ is the blow-up that gives the admissible blow-up $T'\rightarrow T$ by passage to the formal completions.

Since the morphism $T'\rightarrow\widehat{\ovl{X}}$ as above maps $T'$ to the open formal subscheme $\widehat{\til{X}}$, we see that the morphism $Y'\rightarrow \ovl{X}$ maps $Y'$ to $\til{X}$.
Hence we have 
$$
\til{\alpha}\colon Y_U=\Spec B\setminus V(J)=s(\mathscr{T})\longrightarrow X\leqno{(\ast)}
$$
by the base change; cf.\ \S\ref{subsub-associatedscheme} for the notion of the associated schemes\index{scheme!associated scheme to an affinoid@associated --- to an affinoid}\index{affinoid!associated scheme to an affinoid@associated scheme to an ---} $s(\mathscr{T})$.

\begin{thm}\label{thm-affinoidvaluedpoints}
Let $\mathscr{T}=(\Spf B)^{\rig}$ be an affinoid, where $B$ is a t.u.\ adhesive ring\index{t.u.a. ring@t.u.\ adhesive ring}.
Then the map 
$$
\left\{
\begin{minipage}{10em}
{\small morphism $\alpha\colon\mathscr{T}\rightarrow X^{\an}$ of rigid spaces}
\end{minipage}
\right\}
\longrightarrow
\left\{
\begin{minipage}{17em}
{\small pair $(\beta,h)$ consisting of $\beta\colon \mathscr{T}\rightarrow S^{\an}$ and $h\colon s(\mathscr{T})\rightarrow X$ such that the diagram ${\displaystyle \xymatrix@-3ex{s(\mathscr{T})\ar[r]\ar[dr]&X\ar[d]\\ &U}}$ commutes, where the arrow $s(\mathscr{T})\rightarrow s(S^{\an})=U$ is the one obtained from $\beta$ by {\rm \S\ref{subsub-associatedscheme}}}
\end{minipage}
\right\}
$$
given by
$$
\alpha\longmapsto (f^{\an}\circ\alpha,\til{\alpha})\leqno{(\ast\ast)}
$$
$($where $\til{\alpha}$ is the one as in $(\ast))$ is a bijection.
\end{thm}

\begin{proof}
We are going to construct the inverse mapping of the map $(\ast\ast)$.
Given a pair $(\beta,h)$ as above.
As before, we set $\mathscr{T}=T^{\rig}$ where $T=\Spf B$, and set $Y=\Spec B$; we suppose, moreover, that $B$ is $J$-torsion free, where $J\subseteq B$ is an ideal of definition.
For a Nagata compactification $\ovl{X}$ as in \ref{const-GAGAfunctor}, there exists a $Y_U$-admissible blow-up $Y'\rightarrow Y$ and a morphism $\ovl{h}\colon Y'\rightarrow\ovl{X}$ that gives rise to the map $h$ by base change.

So far, we get a morphism $\alpha\colon\mathscr{T}\rightarrow\ovl{X}^{\an}$ of rigid spaces.
We need to show that the image of this map lies in $X^{\an}\subseteq\ovl{X}^{\an}$.
Let $\ovl{Z}$ be as in \ref{const-GAGAfunctor}, and consider the base change
$$
\ovl{X}\times_SY'\longrightarrow Y'.
$$
The morphism $\ovl{h}$ gives a section of this morphism and hence gives a closed immerion $i_1\colon Y'\hookrightarrow\ovl{X}\times_SY'$.
On the other hand, we have the closed immersion $i_2\colon\ovl{Z}\times_SY'\hookrightarrow\ovl{X}\times_SY'$.
Let $\mathscr{J}_j$ be the defining ideal of the closed immerion $i_j$ for $j=1,2$.
These ideals are of finite type, since the schemes we are working with are all $J$-torsion free.
Since $i_1\times_SU$ and $i_2\times_SU$ has the disjoint images, the ideal $\mathscr{J}_1+\mathscr{J}_2$ is an open ideal with respect to the $J$-adic topology.
Let $\mathscr{J}$ be the push-out of the ideal $\mathscr{J}_1+\mathscr{J}_2$ by the projection map $\ovl{X}\times_SY'\rightarrow \ovl{X}$.
Then the blow-up along the ideal $\mathscr{J}$ is an $X$-admissible blow-up.
Replaing $\ovl{X}$ by the resulting one by the blow-up (and $Y'$ and $\ovl{Z}$ by the strict transforms), we see that the image of $\ovl{h}$ is disjoint from $\ovl{Z}$.
Hence we have the morphism $\ovl{h}\colon Y'\rightarrow \til{X}$, as desired.

By passage to the associated rigid spaces, we get 
$$
\alpha=(\widehat{\ovl{h}})^{\rig}\colon\mathscr{T}\longrightarrow(\widehat{\til{X}})^{\rig}\longhookrightarrow X^{\an}.
$$
Then one can show that this $\alpha$ depends only on the data $(\beta,h)$ and that the map $(\beta,h)\mapsto\alpha$ gives the inverse mapping of $(\ast\ast)$.
\end{proof}

\begin{exa}\label{exa-affinoidvaluedpoint1}{\rm 
Let $\mathscr{Z}$ be a locally universally adhesive rigid space over $S^{\an}$.
As usual, we denote by $\G_{\mathrm{a}}^{\an}(\mathscr{Z})$ the set of $S^{\an}$-morphisms $\mathscr{Z}\rightarrow\G_{\mathrm{a}}^{\an}$.
We have
$$
\G_{\mathrm{a}}^{\an}(\mathscr{Z})=\Gamma(\mathscr{Z},\O_{\mathscr{Z}}).
$$
Indeed, since $\G_{\mathrm{a}}^{\an}$ is a sheaf on the site $\CRf_{S^{\an},\ad}$, we may assume $\mathscr{Z}$ is a Stein affinoid (cf.\ \ref{prop-affinoidcov1str}), and then the equality follows from \ref{thm-affinoidvaluedpoints}.

Similarly, we have
$$
\G_{\mathrm{m}}^{\an}(\mathscr{Z})=\Gamma(\mathscr{Z},\O^{\times}_{\mathscr{Z}})
$$
for any rigid space $\mathscr{Z}$.}
\end{exa}

\subsection{Comparison map and comparison functor}\label{sub-GAGAcommap}
\subsubsection{Comparison map}\label{subsub-GAGAcommap}
We continue with working in the situation as in \S\ref{subsub-GAGAfunctoremb} where we assume that the adic ring $A$ is t.u.\ rigid-Noetherian\index{t.u. rigid-Noetherian ring@t.u.\ rigid-Noetherian ring}.
Let $f\colon X\rightarrow U$ be a separated $U$-scheme of finite type.
For any object $(X\hookrightarrow\widehat{X})$ of $\Emb_{X|S}$, consider the $S$-scheme $\til{X}$ as in \ref{const-GAGAfunctor}.
For any affine open subset $\Spec B\hookrightarrow\til{X}$, its open part $\Spec B\setminus V(IB)$ is an open subset of $X$.
In this situation, we have the comparison map\index{comparison map} $s\colon\ZR{(\Spf\widehat{B})^{\rig}}\rightarrow s((\Spf\widehat{B})^{\rig})=\Spec\widehat{B}\setminus V(I\widehat{B})$ (\S\ref{subsub-comparisonmapassaff}) and the composite morphism
$$
\ZR{(\Spf\widehat{B})^{\rig}}\longrightarrow\Spec\widehat{B}\setminus V(I\widehat{B})\longrightarrow\Spec B\setminus V(IB)\longhookrightarrow X\leqno{(\ast)}
$$
of locally ringed spaces.

One can glue the morphisms $(\ast)$ to a morphism
$$
\rho_X\colon (\ZR{X^{\an}},\O_{X^{\an}})\longrightarrow (X,\O_X)
$$
of locally ringed spaces.
Indeed, the set-theoretic map $\rho_X\colon\ZR{X^{\an}}\rightarrow X$ is actually obtained as follows. 
Let $x\in\ZR{X^{\an}}$.
Consider the associated rigid point 
$$
\alpha\colon\mathscr{T}\longrightarrow X^{\an},
$$
where $\mathscr{T}=(\Spf \widehat{V}_x)^{\rig}$ (\ref{dfn-ZRpoints32}).
Then by \ref{thm-affinoidvaluedpoints} one has the morphism of schemes $\Spec\widehat{V}_x[\frac{1}{a}]\rightarrow X$; since $\widehat{V}_x[\frac{1}{a}]$ is a field ({\bf \ref{ch-pre}}.\ref{prop-sep}), it defines a point $y$.
Thus we get the map $\ZR{X^{\an}}\rightarrow X$ by $x\mapsto y$.
It is clear that, compared with the construction of the comparison map in affinoid case (\S\ref{subsub-comparisonmapassaff}), this last map coincides with the one as in $(\ast)$ on each affinoid $(\Spf\widehat{B})^{\rig}$.
Hence we deduce that the maps $(\ast)$ for $\alpha\in L$ and $i\in J_{\alpha}$ glue together to a map $\rho_X$ of locally ringed spaces as above.
Notice that, simply by gluing, one can also define the morphism $\rho_X$ in the case where $X$ is not necessarily separated.

The morphism $\rho_X$ thus obtained for any $U$-scheme $X$ of finite type is called the {\em comparison map}\index{comparison map}.
\begin{prop}\label{prop-comparisonmap3}
The comparison map\index{comparison map} $\rho_X$ is flat.
\end{prop}

\begin{proof}
Due to \ref{prop-comparisonmapassaff} it suffices to show that the map 
$$
\Spec\widehat{B}\setminus V(I\widehat{B})\longrightarrow\Spec B\setminus V(IB)
$$
in $(\ast)$ is flat.
Since $B$ is $IB$-adically t.u.\ rigid-Noetherian ({\bf \ref{ch-formal}}.\ref{dfn-tuaringadmissible}), we know that the map $B\rightarrow\widehat{B}$ is flat ({\bf \ref{ch-pre}}.\ref{prop-btarf1}).
The assertion follows immediately from this.
\end{proof}

\subsubsection{Comparison functor}\label{subsub-GAGAcommapfunct}
\begin{prop}\label{prop-GAGAcommap1cor}
The comparison map $\rho_X$ gives rise to an exact functor
$$
\rho^{\ast}_X\colon\Mod_X\longrightarrow\Mod_{X^{\an}}.
$$
Moreover, it maps the full subcategory $\Coh_X$ to $\Coh_{X^{\an}}$.
\end{prop}

We call the functor $\rho^{\ast}_X$ the {\em comparison functor}\index{functor!comparison functor@comparison ---}.
\begin{proof}
The first half is clear by \ref{prop-comparisonmap3}.
Since $X$ is a Noetherian scheme, coherent $\O_X$-modules are exactly the $\O_X$-modules of finite presentation.
Similarly, coherent $\O_{X^{\an}}$-modules are nothing but finitely presented $\O_{\mathscr{X}}$-modules (\ref{cor-existlattice3}).
Hence the second assertion follows.
\end{proof}

\begin{rem}\label{rem-GAGAcommapfunct}{\rm 
Here is another description of the comparison functor for coherent sheaves; here we assume that the scheme $X$ is separated over $U$.
Let $\mathscr{F}$ be a coherent sheaf on the scheme $X$.
Consider the quasi-compact open immersion $X\hookrightarrow\til{X}$ as in \ref{const-GAGAfunctor}, and extend $\mathscr{F}$ to a finitely presented $\O_{\til{X}}$-module $\til{\mathscr{F}}$ (\cite[6.9.11]{EGAInew}). 
Then the formal completion $\widehat{\til{\mathscr{F}}}$ on $\widehat{\til{X}}$ is finitely presented, and thus we get $(\widehat{\til{\mathscr{F}}})^{\rig}$ (\S\ref{subsub-latticemodelsrig}), which is a coherent sheaf on $(\widehat{\til{X}}_{\alpha})^{\rig}$.
We get by gluing the desired coherent sheaf, which is nothing but $\rho^{\ast}_X\mathscr{F}$.}
\end{rem}

\subsection{GAGA comparison theorem}\label{sub-GAGAcomprigid}
\index{GAGA!GAGA comparison theorem@--- comparison theorem|(}
In this and next subsections, we work in the situation as in \S\ref{subsub-GAGAfunctoremb} with the additional assumption that
\begin{itemize}
\item[$(\ast)$] the adic ring $A$ is  t.u.\ rigid-Noetherian\index{t.u. rigid-Noetherian ring@t.u.\ rigid-Noetherian}.
\end{itemize}
The theorems in these subsections are announced in this situation, but the proofs will be done with the stronger assumption that $A$ is t.u.\ adhesive. 
The general case can be proven with the aid of what we discussed in {\bf \ref{ch-formal}}, \S\ref{sec-weakcoherency}; we wish to give the details in \cite{FK}.

Let 
$$
f\colon X\longrightarrow Y
$$
be a proper $U$-morphisms of separated and of finite type $U$-schemes.
We have the commutative diagram 
$$
\xymatrix{X\ar[d]_f&(\ZR{X^{\an}},\O_{X^{\an}})\ar[l]_(.66){\rho_X}\ar[d]^{f^{\an}}\\ Y&(\ZR{Y^{\an}},\O_{Y^{\an}})\ar[l]^(.66){\rho_Y}}
$$
of locally ringed spaces, where $\rho_X$ and $\rho_Y$ are the comparison maps (\S\ref{subsub-GAGAcommap}).

By {\bf \ref{ch-formal}}.\ref{thm-fini} we know that $\RD f_{\ast}$ maps $\DC^{\ast}_{\coh}(X)$ to $\DC^{\ast}_{\coh}(Y)$ for $\ast=$``\ \ '', $+$, $-$, $\bd$.
On the other hand, since the comparison functor $\rho^{\ast}_X$ is exact (\ref{prop-GAGAcommap1cor}), it induces an exact functor
$$
\DC^{\ast}(X)\longrightarrow\DC^{\ast}(X^{\an})
$$
(cf.\ {\bf \ref{ch-pre}}.\ref{prop-derivedcategory71}), where $\DC^{\ast}(X^{\an})=\DC^{\ast}(\ZR{X^{\an}},\O_{X^{\an}})$.
We write this functor as
$$
M\longmapsto M^{\rig}.
$$
Since $\rho_X$ is flat, one sees easily that this functor maps $\DC^{\ast}_{\coh}(X)$ to $\DC^{\ast}_{\coh}(X^{\an})$.
Similarly, we get a canonical functor from $\DC^{\ast}_{\coh}(Y)$ to $\DC^{\ast}_{\coh}(Y^{\an})$.

Thus we get the following diagram of triangulated categories:
$$
\xymatrix{\DC^{\ast}_{\coh}(X)\ar[r]^{\rig}\ar[d]_{\RD f_{\ast}}&\DC^{\ast}_{\coh}(X^{\an})\ar[d]^{\RD f^{\an}_{\ast}}\\ \DC^{\ast}_{\coh}(Y)\ar[r]_{\rig}&\DC^{\ast}(Y^{\an})\rlap{.}}
$$
Note that by \ref{prop-finitudesrigid1} the functor $\RD f^{\an}_{\ast}$ has finite cohomology dimension, whence the right-hand vertical arrow.

By an argument similar to that in {\bf \ref{ch-formal}}, \S\ref{sub-GFGAcomann}, one has the comparison map 
$$
\rho=\rho_f\colon\rig\circ\RD f_{\ast}\longrightarrow\RD f^{\an}_{\ast}\circ\rig,
$$
and thus we obtain the following diagram with a $2$-arrow:
$$
\xymatrix@-2ex{
\DC^{\ast}_{\coh}(X)\ar[rr]^{\rig}\ar[dd]_{\RD f_{\ast}}&&\DC^{\ast}_{\coh}(X^{\an})\ar[dd]^{\RD f^{\an}_{\ast}}\\
&&\ \\
\DC^{\ast}_{\coh}(Y)\ar[rr]_{\rig}&\ar@{=>}[ur]^{\rho}&\DC^{\ast}(Y^{\an})}\rlap{.}\leqno{(\ast\ast)}
$$

\begin{thm}[GAGA comparison theorem]\label{thm-GAGAcomprigid}
Suppose $f\colon X\rightarrow Y$ is proper.
Then the natural transformation $\rho$ gives a natural equivalence$;$ hence the diagram $(\ast\ast)$ is $2$-commutative.
\end{thm}

We give here the proof only in case the ring $A$ in \S\ref{subsub-GAGAfunctoremb} is t.u.\ adhesive\index{t.u.a. ring@t.u.\ adhesive ring}; for the proof of the general case, see \cite{FK}.
\begin{proof}
First note that, to show the theorem, we may restrict ourselves to considering only objects $M$ of $\DC^{\bd}_{\coh}(X)$ that is concentrated in degree $0$ by the similar reduction process as in {\bf \ref{ch-formal}}, \S\ref{subsub-GFGAcompf2} in the following way:
\begin{itemize}
\item[(1)] first we may assume $\ast=\bd$ (standard);
\item[(2)] by induction with respect to the amplitude of $M$, taking the distinguished triangle
$$
\tau^{\leq n}M\longrightarrow M\longrightarrow\tau^{\geq n+1}M\stackrel{+1}{\longrightarrow}
$$
into account, we may assume that $\amp(M)=0$;
\item[(3)] finally, by shift, we arrive at the hypothesis as above.
\end{itemize}
The theorem in this case will be shown in \ref{prop-GAGAcomprigid1} below.
\end{proof}

\begin{prop}\label{prop-GAGAcomprigid1}
Working in the situation as in {\rm \ref{thm-GAGAcomprigid}}, let $\mathscr{F}$ be a coherent $\O_X$-module.
Then the canonical morphism 
$$
(\RD^q f_{\ast}\mathscr{F})^{\rig}\longrightarrow\RD f^{\an}_{\ast}\mathscr{F}^{\rig}
$$
is an isomorphism for all $q\geq 0$.
\end{prop}

\begin{proof}
We first deal with the case where $Y$ is proper over $U$.
In this case, one can take suitable objects $X\hookrightarrow\ovl{X}$ in $\Emb_{X|S}$ and $Y\hookrightarrow\ovl{Y}$ in $\Emb_{Y|S}$ sitting in the Cartesian diagram
$$
\xymatrix{X\,\ar[d]_f\ar@{^{(}->}[r]&\ovl{X}\ar[d]^{\ovl{f}}\\ Y\,\ar@{^{(}->}[r]&\ovl{Y}\rlap{.}}
$$
Since $X^{\an}=(\widehat{\ovl{X}})^{\rig}$ and $Y^{\an}=(\widehat{\ovl{Y}})^{\rig}$, one can show the claim by the similar reasoning as in \ref{thm-comparisonaffinoid}.

In general, $Y^{\an}$ is the open complement (\ref{dfn-opencomplement}) of a closed subspace $\mathscr{Z}$ in $(\widehat{\ovl{Y}})^{\rig}$ (due to \ref{prop-GAGAfunctor5x2}).
Since $f$ is proper, $X^{\an}$ is, similarly, the open complement of $\mathscr{W}$, the pull-back of $\mathscr{Z}$, in $(\widehat{\ovl{X}})^{\rig}$.
Note that $\ovl{X}^{\an}_U=(\widehat{\ovl{X}})^{\rig}$ and $\ovl{Y}^{\an}_U=(\widehat{\ovl{Y}})^{\rig}$.
As the assertion is local on $Y$, one can reduce to the above case in the following way. 
First extend the coherent sheaf $\mathscr{F}$ on $X$ to a coherent sheaf $\mathscr{G}$ on the proper scheme $\ovl{X}_U$; this is possible, since $\ovl{X}_U$ is a Noetherian scheme.
Then due the first step, the assertion is true for the morphism $\ovl{X}_U\rightarrow\ovl{Y}_U$.

Now since $f$ is proper, the commutative diagram
$$
\xymatrix{X\,\ar@{^{(}->}[r]\ar[d]_f&\ovl{X}_U\ar[d]^{\ovl{f}_U}\\ Y\,\ar@{^{(}->}[r]&\ovl{Y}_U}
$$
is Cartesian, and by \ref{prop-GAGAfunctor5x3} we have the Cartesian diagram
$$
\xymatrix{X^{\an}\,\ar@{^{(}->}[r]\ar[d]_{f^{\an}}&\ovl{X}_U^{\an}\ar[d]^{\ovl{f}_U^{\an}}\\ Y^{\an}\,\ar@{^{(}->}[r]&\ovl{Y}_U^{\an}}
$$
of rigid spaces.
Hence we get the desired result by the base change.
\end{proof}

The following corollaries can be shown similarly, at least in the t.u.\ adhesive situation, to {\bf \ref{ch-formal}}, \S\ref{sub-GFGAcomext}:
\begin{prop}\label{prop-GAGArigidcom1}
In the situation as in {\rm \ref{thm-GAGAcomprigid}}, let $f\colon X\rightarrow U$ be a proper $U$-scheme.
Then for $M\in\obj(\DC^-_{\coh}(X))$ and $N\in\DC^+_{\coh}(X)$ we have the canonical isomorphism
$$
\RD\Hom_{\O_X}(M,N)\cong\RD\Hom_{\O_{X^{\an}}}(M^{\rig},N^{\rig})
$$
in $\DC^+(\Ab)$. \hfill$\square$
\end{prop}

\begin{cor}\label{cor-GAGArigidcom11}
In the situation as in {\rm \ref{thm-GAGAcomprigid}}, let $f\colon X\rightarrow U$ be a proper $U$-scheme.
Then the comparison functor 
$$
\DC^{\bd}_{\coh}(X)\stackrel{\rig}{\longrightarrow}\DC^{\bd}_{\coh}(X^{\an})
$$
is fully faithful. \hfill$\square$
\end{cor}
\index{GAGA!GAGA comparison theorem@--- comparison theorem|)}

\subsection{GAGA existence theorem}\label{sub-GAGAexistrigid}
\index{GAGA!GAGA existence theorem@--- existence theorem|(}
We continue to work in the situation described in the beginning of \S\ref{sub-GAGAcomprigid}.
\begin{thm}[GAGA existence theorem]\label{thm-GAGAexistrigid}
Let $f\colon X\rightarrow U$ be a proper $U$-scheme.
Then the comparison functor
$$
\DC^{\bd}_{\coh}(X)\stackrel{\rig}{\longrightarrow}\DC^{\bd}_{\coh}(X^{\an})
$$
is an exact equivalence of triangulated categories.
\end{thm}

Similarly to the comparison theorem, here we give the proof only in case where the ring $A$ is t.u.\ adhesive\index{t.u.a. ring@t.u.\ adhesive ring}; for the proof of the general case, see \cite{FK}.
\begin{proof}
As we have already seen that the functor in question is exact and fully faithful (\ref{cor-GAGArigidcom11}), we only need to show that it is essentially surjective.
But then by the reduction process carried out by induction with respect to the amplitudes of objects of $\DC^{\bd}_{\coh}(X^{\an})$ involving the distinguished triangle 
$$
\tau^{\leq n}M\longrightarrow M\longrightarrow\tau^{\geq n+1}M\stackrel{+1}{\longrightarrow}
$$
(similarly to that in {\bf \ref{ch-formal}}, \S\ref{subsub-GFGAexistproofproj}), it suffices to show that the functor $\rig\colon\Coh_X\rightarrow\Coh_{X^{\an}}$ is essentially surjective.

Take a diagram
$$
\xymatrix{X\ \ar@{^{(}->}[r]\ar[d]_f&\ovl{X}\ar[d]^{\ovl{f}}\\ U\ \ar@{^{(}->}[r]&S\rlap{,}}
$$
which is an object of $\Emb_{X|S}$; as we saw in the proof of \ref{prop-GAGAfunctor4}, we may, moreover, assume that this diagram is Cartesian.
We have $X^{\an}=(\widehat{\ovl{X}})^{\rig}$ (\ref{prop-GAGAfunctor4}).
Moreover, we may assume that $\ovl{X}$ is $I$-torsion free.

Let $\mathscr{F}$ be a coherent $\O_{X^{\an}}$-module.
Then by \ref{thm-tateacyclic20} there exists a coherent $\O_{\widehat{\ovl{X}}}$-module $\til{\mathscr{G}}$ such that $(\til{\mathscr{G}})^{\rig}=\mathscr{F}$.
By {\bf \ref{ch-formal}}.\ref{thm-GFGAexist} we get a coherent $\O_{\ovl{X}}$-module $\mathscr{G}$ of which the formal completion coincides with $\til{\mathscr{G}}$.
Now the sheaf $\mathscr{G}|_X$ is the desired one, for we have $(\mathscr{G}|_X)^{\rig}=\mathscr{F}$, as we have seen in \ref{rem-GAGAcommapfunct}.
\end{proof}

By an argument similar to that in {\bf \ref{ch-formal}}.\ref{prop-GFGAexistapp1} one can show:
\begin{cor}\label{cor-GAGAexistrigid1}
The GAGA-functor
$$
\cdot^{\an}\colon\PSch_U\longrightarrow\Rf_{\mathscr{S}}
$$
is fully faithful, where $\PSch_U$ denotes the category of proper $U$-schemes. \hfill$\square$
\end{cor}
\index{GAGA!GAGA existence theorem@--- existence theorem|)}

\subsection{Adic part for non-adic morphisms}\label{sub-adicpart}
\index{adic!adic part@--- part|(}
Finally, let us observe that one can apply Berthelot's\index{Berthelot, P.} construction of tubes \cite{Berth1} to construct canonical rigid spaces associated to non-adic morphisms of formal schemes.
The relation with GAGA will be shown in Exercise \ref{exer-GAGAfunctor2}.

\subsubsection{Adic part}\label{subsub-adicpart}
Let $f\colon X\rightarrow Y$ be a morphism between adic formal schemes of finite ideal type.
We do note assume that $f$ is adic.
Let $\mathscr{X}=X^{\rig}$ and $\mathscr{Y}=Y^{\rig}$ be the associated rigid spaces (\S\ref{subsub-extensionofrigfuntor}).
If $X$ is not adic over $Y$, $\mathscr{X}$ does not admit the structural map $\mathscr{X}\rightarrow\mathscr{Y}$. 
However, following \cite[\S0.2]{Berth1}, one can construct a canonical open rigid subspace $\mathscr{X}^{\adic}_{/\mathscr{Y}}\subseteq\mathscr{X}$, called the {\em adic part over $\mathscr{Y}$}, that admits the $\mathscr{Y}$-strcuture, as follows. 

Consider the subset $\mathfrak{U}$ of $\ZR{\mathscr{X}}$ defined as follows:
$$
\mathfrak{U}=\bigg\{x\in\ZR{\mathscr{X}}\,\bigg|\,\begin{minipage}{24em}{\small $f\circ\sp_X\circ\alpha_x\colon\Spf\widehat{V}_x\rightarrow Y$ is adic, where $\alpha_x\colon\Spf\widehat{V}_x\rightarrow\ZR{\mathscr{X}}$ is the associated rigid point of $x$ (\ref{dfn-ZRpoints32})}\end{minipage}\bigg\}. 
$$
\begin{prop}\label{prop-formallyfinitetype}
{\rm (1)} The subset $\mathfrak{U}$ is open in $\ZR{\mathscr{X}}$.

{\rm (2)} There exists a canonical open rigid subspace $\mathscr{X}^{\adic}_{/\mathscr{Y}}$ of $\mathscr{X}$ with the canonical morphism $\mathscr{X}^{\adic}_{/\mathscr{Y}}\rightarrow\mathscr{Y}$ such that $\ZR{\mathscr{X}^{\adic}_{/\mathscr{Y}}}=\mathfrak{U}$.
\end{prop}

\begin{proof}
To show (1), we can discuss in affine situation: Consider $f\colon X=\Spf A\rightarrow Y=\Spf B$, and let $I$ (resp.\ $J$) be a finitely generated ideal of definition of $A$ (resp.\ $B$) such that $JA\subseteq I$. 
We may moreover suppose, replacing $Y$ by the admissible blow-up along $J$ if necessary, $J$ is invertible and principal $J=aB$.
Notice that $A$ is $a$-adically complete due to {\bf \ref{ch-pre}}.\ref{prop-exertopologicallynilpotentadic} (b).
Suppose $x\in\mathfrak{U}$.
Since the composition $B\rightarrow A\rightarrow\widehat{V}_x$ is adic, $a\widehat{V}_x$ is an ideal of definition of $\widehat{V}_x$, and hence $I^n\widehat{V}_x\subseteq a\widehat{V}_x$ for some $n>0$.
Set $I^n=(f_1,\ldots,f_r)$ and consider the admissible blow-up $\til{X}$ of $X=\Spf A$ along the admissible ideal $JA+I^n=(a,f_1,\ldots,f_r)$; $\til{X}$ has the affine open part given by $\Spf A_{I^n}$ where 
$${\textstyle 
A_{I^n}=A\dl\frac{f_1}{a},\ldots,\frac{f_r}{a}\dr,}
$$
that is, the part where $(JA+I^n)\O_{\widetilde{X}}$ is generated by $J\O_{\widetilde{X}}$. 
The map $\sp_X\circ\alpha_x\colon\Spf\widehat{V}_x\rightarrow X=\Spf A$ factors through $\Spf A_{I^n}$.
Since $(\Spf A_{I^n})^{\rig}$ is an open subspace of $\mathscr{X}=X^{\rig}$ and since any point $y\in\ZR{(\Spf A_{I^n})^{\rig}}$ satisfies $I\widehat{V}_y\subseteq a\widehat{V}_y$, we have $\ZR{(\Spf A_{I^n})^{\rig}}\subseteq\mathfrak{U}$, and thus the openness of $\mathfrak{U}$ follows.

As the above argument shows, for any $x\in\mathfrak{U}$, there exists a positive integer $n>0$ such that $I^n\widehat{V}_x\subseteq a\widehat{V}_x$, and hence the open subset $\mathfrak{U}\subseteq\ZR{\mathscr{X}}$ is described as
$$
\mathfrak{U}=\bigcup_{n>0}\ZR{(\Spf A_{I^n})^{\rig}},
$$
where each $\ZR{(\Spf A_{I^n})}$ is naturally identified with a coherent open subset of $\ZR{\mathscr{X}}$.
Hence one can construct the desired rigid subspace in this case by
$$
\mathscr{X}^{\adic}_{/\mathscr{Y}}=\varinjlim_{n>0}(\Spf A_{I^n})^{\rig},
$$
which clearly does not depend on the choice of the ideal of definition $I$.
The construction in the general case is given by gluing, whence (2).
\end{proof}

\begin{rem}\label{rem-formallyfinitetypeopeninterior}{\rm 
In the situation as in the proof of \ref{prop-formallyfinitetype}, consider $Z=\Spf A$, where $A$ is now considered as an adic ring by the $J$-adic topology, and $W=\Spec A/I$.
Then, by the open interior formula (\ref{prop-tubes11}), we have $\mathfrak{U}=C_{W|Z}$; cf.\ Exercise \ref{exer-GAGAfunctor2}.}
\end{rem}

\begin{rem}[{\rm cf.\ \cite[(0.2.6)]{Berth1}}]\label{rem-formallyfinitetype}{\rm 
The rigid space $\mathscr{X}^{\adic}_{/\mathscr{Y}}$ in the affine case as in the proof of \ref{prop-formallyfinitetype} admits another open covering as follows.
Set $I=(f_1,\ldots,f_r)$.
Define, for any positive integer $n>0$, 
$${\textstyle 
A(n)=A\dl\frac{f^n_1}{a},\ldots,\frac{f^n_r}{a}\dr=A\dl T_1,\ldots,T_r\dr/(f^n_1-aT_1,\ldots,f^n_r-aT_r).}
$$
Then the rigid space $\mathscr{U}(n)=(\Spf A(n))^{\rig}$ is an open subspace of $\mathscr{X}$ such that 
$$
(\Spf A_{I^n})^{\rig}\subseteq\mathscr{U}(n)\subseteq(\Spf A_{I^{nr}})^{\rig}
$$
for any $n>0$.
Hence one has 
$$
\mathscr{X}^{\adic}_{/\mathscr{Y}}=\varinjlim_{n>0}\mathscr{U}(n).
$$
Note that, here, the open immersion $\mathscr{U}(n)\hookrightarrow\mathscr{U}(m)$ for $m>n$ is given by $A(m)\rightarrow A(n)$ that maps $f^m_i/a$ to $f^{m-n}(f^n_i/a)$ for $i=1,\ldots,r$.}
\end{rem}

\subsubsection{Functoriality}\label{subsub-adicpartfunctoriality}
In the situation as above, suppose $X$ is coherent. Let $Z$ be a coherent adic formal scheme of finite ideal type, and $g\colon Z\rightarrow Y$ an adic morphism.
Suppose we are given a morphism $\mathscr{Z}=Z^{\rig}\rightarrow\mathscr{X}^{\adic}_{/\mathscr{Y}}$ of rigid spaces over $\mathscr{Y}$.
Then the composition with the open immersion $\mathscr{X}^{\adic}_{/\mathscr{Y}}\hookrightarrow\mathscr{X}$ is induced by the rig-functor from a $Y$-morphism $Z'\rightarrow X$ of formal schemes, where $Z'$ is an admissible blow-up of $Z$.
Conversely, suppose we are given a $Y$-morphism of the form $Z'\rightarrow X$ from an admissible blow-up $Z'$ of $Z$.
For any admissible ideal $\mathscr{J}$ of $X$, since $Z'$ is adic over $Y$, $\mathscr{J}\O_{Z'}$ is an admissible ideal of $Z'$.
This means, in other words, for any admissible blow-up $X'\rightarrow X$, there exists an admissible blow-up $Z''\rightarrow Z'$ and a morphism $Z''\rightarrow X'$ such that the diagram
$$
\xymatrix{Z''\ar[r]\ar[d]&X'\ar[d]\\ Z'\ar[r]&X}
$$
commutes.
Hence we have a morphism $\mathscr{Z}\rightarrow\mathscr{X}$ of rigid spaces.
For any point $z\in\ZR{\mathscr{Z}}$, the map $\Spf\widehat{V}_z\rightarrow Y$ is adic, and hence the composition $\Spf\widehat{V}_z\rightarrow\ZR{\mathscr{Z}}\rightarrow\ZR{\mathscr{X}}$ has its image in the open part $\mathfrak{U}$ as in \S\ref{subsub-adicpart}.
We therefore obtain a morphism $\mathscr{Z}\rightarrow\mathscr{X}^{\adic}_{/\mathscr{Y}}$ of rigid spaces over $\mathscr{Y}$.

What we have seen is the existence of a canonical bijection 
$$
\varinjlim_{Z'}\Hom_Y(Z',X)\stackrel{\sim}{\longrightarrow}\Hom_{\Rf_{\mathscr{Y}}}(\mathscr{Z},\mathscr{X}^{\adic}_{/\mathscr{Y}}),
$$
where $Z'$ runs over all admissible blow-ups of $Z$.

\begin{prop}
Let $Y$ be an adic formal scheme of finite ideal type, and $\mathscr{Y}=Y^{\rig}$.
Then, by $X\mapsto(X^{\rig})^{\adic}_{/\mathscr{Y}}$, we have a natural functor $\Ac\Fs_Y\rightarrow\Rf_{\mathscr{Y}}$ from the category of adic formal schemes over $Y$ to the category of rigid spaces over $\mathscr{Y}$.
\end{prop}

\begin{proof}
We need to show that any morphism $W\rightarrow X$ in $\Ac\Fs_Y$ naturally induces a morphism $(W^{\rig})^{\adic}_{\mathscr{Y}}\rightarrow(X^{\rig})^{\adic}_{\mathscr{Y}}$.
For any open subspace $\mathscr{U}\subseteq\mathscr{X}=X^{\rig}$, we have $\mathscr{U}^{\adic}_{\mathscr{Y}}=\mathscr{U}\cap\mathscr{X}^{\adic}_{\mathscr{Y}}$.
Hence, by a standard patching argument, we may assume that $X$ is coherent.

The proof of \ref{prop-formallyfinitetype} shows that every point of $(W^{\rig})^{\adic}_{\mathscr{Y}}$ has a coherent open neighborhood of the form $Z^{\rig}$ by a formal subscheme $Z\subseteq W'$, adic over $Y$, of an admissible blow-up $W'$ of $W$.
By what we have seen above, we have canonically a morphism $Z^{\rig}\rightarrow\mathscr{X}^{\adic}_{\mathscr{Y}}$.
Then, by gluing, we have the desired morphism $(W^{\rig})^{\adic}_{\mathscr{Y}}\rightarrow\mathscr{X}^{\adic}_{\mathscr{Y}}$.
\end{proof}

\subsubsection{Adic part for formally locally of finite type morphisms}\label{subsub-adicpart-flft}
\begin{dfn}\label{dfn-formallyfinitetype}{\rm Let $f\colon X\rightarrow Y$ be a (not necessarily adic) morphism between adic formal schemes of finite ideal type. 
We say that $f$ is {\em formally locally of finite type}\index{morphism of formal schemes@morphism (of formal schemes)!morphism of formal schemes formally locally of finite type@--- formally locally of finite type} if, for any $x\in X$, there exists an open neighborhood $V\subseteq Y$ of $y=f(x)$ and an open neighborhood $U\subseteq X$ of $x$ such that $f(U)\subseteq V$, and that, if $\mathscr{I}$ (resp.\ $\mathscr{V}$) is an ideal of definition of $U$ (resp.\ $V$) of finite type such that $\mathscr{J}\O_X\subseteq\mathscr{I}$, the induced morphism $(X,\O_X/\mathscr{I})\rightarrow(Y,\O_Y/\mathscr{J})$ of schemes is locally of finite type. 
If, in addition, $f$ is quasi-compact, it is said to be {\em formally of finite type}\index{morphism of formal schemes@morphism (of formal schemes)!morphism of formal schemes formally of finite type@--- formally of finite type}.}
\end{dfn}

Notice that, in view of {\bf \ref{ch-formal}}.\ref{prop-cortopfintype11}, if $f\colon X\rightarrow Y$ is adic, then $f$ is formally locally of finite type (resp.\ formally of finite type) if and only if it is locally of finite type\index{morphism of formal schemes@morphism (of formal schemes)!morphism of formal schemes locally of finite type@--- locally of finite type} (resp.\ of finite type\index{morphism of formal schemes@morphism (of formal schemes)!morphism of formal schemes of finite type@--- of finite type}).

\begin{prop}\label{prop-formallyfinitetype3}
Let $f\colon X\rightarrow Y$ be a formally locally of finite type morphism between adic formal schemes of finite ideal type, and set $\mathscr{X}=X^{\rig}$ and $\mathscr{Y}=Y^{\rig}$. 
Then the rigid space $\mathscr{X}^{\adic}_{/\mathscr{Y}}$ is locally of finite type over $\mathscr{Y}$.
\end{prop}

\begin{proof}
We may discuss in affine situation.
In the notation as in the proof of \ref{prop-formallyfinitetype}, the assumption implies that $B/J\rightarrow A/I$ is of finite type.
We have $A_I/IA_I=A_I/aA_I$, which is, then, topologically of finite type over $B$.
Hence the assertion follows.
\end{proof}

\subsubsection{Examples}\label{subsub-adicpartexample}
\begin{exa}[{\rm Open unit disk}]\label{exa-openunitdisk}{\rm 
\index{unit disk!open unit disk@open ---}
Let $B$ be an adic ring of finite ideal type, with a finitely generated ideal of definition $J$, and consider the formal power series ring $A=B[\![X_1,\ldots,X_r]\!]$ with an ideal of definition $I=JA+(X_1,\ldots,X_r)$.
The rigid space $\mathscr{X}^{\adic}_{/\mathscr{Y}}$ (where $\mathscr{X}=(\Spf A)^{\rig}$ and $\mathscr{Y}=(\Spf B)^{\rig}$) in this situation is called the {\em open unit disk over $\mathscr{Y}$}.
To describe $\mathscr{X}^{\ast}_{/\mathscr{Y}}$ locally over $\mathscr{Y}$, we assume $J$ is ivertible $J=aB$.
As in \ref{rem-formallyfinitetype}, we have
$${\textstyle 
\mathscr{X}^{\adic}_{/\mathscr{Y}}=\bigcup_{n>0}(\Spf A\dl\frac{X^n_1}{a},\ldots,\frac{X^n_r}{a}\dr)^{\rig}}.
$$
Here, notice that we have 
$${\textstyle 
A\dl\frac{X^n_1}{a},\ldots,\frac{X^n_r}{a}\dr=B\dl X_1,\ldots,X_r,Y_1,\ldots,Y_r\dr/(X^n_1-aY_1,\ldots,X^n_r-aY_r),}
$$
and hence $(\Spf A\dl\frac{X^n_1}{a},\ldots,\frac{X^n_r}{a}\dr)^{\rig}$ is naturally viewed as a rational subdomain\index{affinoid!affinoid subdomain@--- subdomain!rational subdomain@rational subdomain}\index{rational subdomain} (\ref{exas-affinoidsubdomain}) 
$${\textstyle 
\D^r_{\mathscr{Y}}(\frac{X^n_1}{a},\ldots,\frac{X^n_r}{a})}
$$
of the unit disk\index{unit disk} $\D^r_{\mathscr{Y}}$ over $\mathscr{Y}$ (\S\ref{subsub-unitdisk}).
Thus, $\mathscr{X}^{\adic}_{/\mathscr{Y}}$ is isomorphic to the open subspace of $\D^r_{\mathscr{Y}}$ that is the union of all $\D^r_{\mathscr{Y}}(\frac{X^n_1}{a},\ldots,\frac{X^n_r}{a})$.}
\end{exa}

By the last description of the open unit disk $\mathscr{X}^{\adic}_{/\mathscr{Y}}$, if, for example, $B=V$ is an $a$-adically complete valuation ring, then it has the following set-theoretic description:
$$
\ZR{\mathscr{X}^{\adic}_{/\mathscr{Y}}}=\{x\in\ZR{\D^r_{\mathscr{Y}}}\,|\,\|X_i(x)\|_{a,c}<1\ \textrm{for}\ i=1,\ldots,r\},
$$
where $\|\cdot\|_{x,a,c}$ is the seminorm\index{seminorm} defined as in \S\ref{subsub-seminormsofpoints}.

\begin{exa}[{\rm Open annulus}]\label{exa-semistablecase}{\rm 
Let $V$ be an $a$-adically complete valuation ring, and set $A=V[\![X,Y]\!]/(XY-a)$.
Then one sees that the adic part $\mathscr{X}^{\adic}_{/\mathscr{Y}}$ (where $\mathscr{X}=(\Spf A)^{\rig}$ and $\mathscr{Y}=(\Spf V)^{\rig}$) is isomorphic to the open part of the unit disk\index{unit disk} $\D^1_{\mathscr{Y}}=(\Spf V\dl X\dr)^{\rig}$ characterized by
$$
\ZR{\mathscr{X}^{\adic}_{/\mathscr{Y}}}=\{x\in\ZR{\D^1_{\mathscr{Y}}}\,|\,c<\|X(x)\|_{a,c}<1\}.
$$
}
\end{exa}

\begin{exa}[{\rm Berthelot's\index{Berthelot, P.} tube (cf.\ \cite[\S1]{Berth1})}]{\rm 
Let $P$ be a locally of finite type scheme over $S=\Spec V$, where $V$ is an $a$-adically complete valuation ring of height one, and $X$ a closed subscheme of the scheme $P_0=(P,\O_P/a\O_P)$.
Consider the formal completion $\widehat{P}|_X$ of $P$ along $X$.
Then we have 
$$
((\widehat{P}|_X)^{\rig})^{\adic}_{\mathscr{S}}=]X[_P,
$$
where $\mathscr{S}=(\Spf V)^{\rig}$.}
\end{exa}
\index{adic!adic part@--- part|)}
\index{GAGA|)}

\addcontentsline{toc}{subsection}{Exercises}
\subsection*{Exercises}
\begin{exer}\label{exer-GAGAfunctor1}{\rm 
In the situation as in \S\ref{subsub-GAGAfunctoremb}, let $X$ be a separated $U$-scheme of finite type.
Show that $X^{\an}$ is the union of countably many affinoids.}
\end{exer}

\begin{exer}\label{exer-GAGAfunctor2}{\rm 
Let $A$ be an adic ring of finite ideal type, $I\subseteq A$ an ideal of definition, and $J\subseteq A$ a finitely generated ideal contained in $I$.
Consider $D=V(J)\subseteq X=\Spec A$ and the formal completions $\widehat{X}=\Spf A$ and $\widehat{X}|_D$ of $X$ along $I$ and $J$, respectively.
Set $\mathscr{X}=(\widehat{X})^{\rig}$ and $\mathscr{S}=(\widehat{X}|_D)^{\rig}$.
Then show that $\mathscr{X}^{\adic}_{/\mathscr{S}}$ is canonically isomorphic to $(X\setminus D)^{\an}$.}
\end{exer}


\section{Dimension of rigid spaces}\label{sec-dimension}
This section gives generalities on dimension and codimension in rigid geometry.
Like in complex analytic geometry, the dimension of rigid spaces is defined to be the supremum of local dimension at each point (\S\ref{subsub-dimensionatpoint}), which is defined as the Krull dimension of the local ring.
In \S\ref{subsub-dimensionfunction} we will briefly discuss the so-called {\it dimension function}, which will be used in the next section.
In case the rigid spaces in consideration are of type (N) or of type (V), the notion of dimension has several nice interpretations and properties, which will be discussed in \S\ref{subsub-comparisondim}.
We have, moreover, the GAGA comparison of dimensions, as discussed in \S\ref{subsub-dimaffinoidscor}.

In \S\ref{sub-codimention} we briefly discuss codimension of rigid subspaces.
The final section \S\ref{sub-relativedim} gives the notion of relative dimension for locally of finite type morphisms of rigid spaces.

\subsection{Dimension of rigid spaces}\label{sub-dimension}
\subsubsection{Dimension}\label{subsub-dimensionatpoint}
\begin{dfn}\label{dfn-dimensionatpoint}{\rm 
Let $\mathscr{X}$ be a rigid space.

(1) Let $x\in\ZR{\mathscr{X}}$ be a point.
The {\it dimension of $\mathscr{X}$ at $x$}\index{dimension of a rigid space@dimension (of a rigid space)}, denoted by 
$$
\dim_x(\mathscr{X}),
$$
is the Krull dimension $\dim(\O_{\mathscr{X},x})$ of the local ring $\O_{\mathscr{X},x}$ at $x$.

(2) The {\it dimension of $\mathscr{X}$}\index{dimension of a rigid space@dimension (of a rigid space)}, denoted by $\dim(\mathscr{X})$, is defined by 
$$
\dim(\mathscr{X})=\sup_{x\in\ZR{\mathscr{X}}}\dim_x(\mathscr{X}).
$$}
\end{dfn}

We set $\dim(\mathscr{X})=-\infty$ if $\mathscr{X}$ is empty.
Notice that, if $\mathscr{X}$ is locally universally Noetherian, and if $y$ is a generization of $x$, then we have $\dim_y(\mathscr{X})\geq\dim_x(\mathscr{X})$, since $\O_{\mathscr{X},y}$ is faithfully flat over $\O_{\mathscr{X},x}$ (\ref{thm-fibersoverrigptsbehavior} (3)).

\begin{exa}\label{exa-dimensionvaluationrings}{\rm 
Let $V$ be an $a$-adically complete valuation ring of arbitrary (positive) height, and set $\mathscr{X}=(\Spf V)^{\rig}$.
Then
$$
\dim(\mathscr{X})=0,
$$
since, in this case, $\ZR{\mathscr{X}}$ consists only of one point with the local ring being isomorphic to the fractional field of $V$.}
\end{exa}

\begin{prop}\label{prop-dimension2}
Let $\mathscr{X}$ be a universally Noetherian rigid space\index{rigid space!universally Noetherian rigid space@universally Noetherian ---!locally universally Noetherian rigid space@locally --- ---}, and $\mathscr{Y}$ a rigid subspace\index{rigid subspace} {\rm (\ref{dfn-rigidsubspace})} of $\mathscr{X}$.

{\rm (1)} For any point $x\in\ZR{\mathscr{Y}}$ we have $\dim_x(\mathscr{Y})\leq\dim_x(\mathscr{X})$.

{\rm (2)} We have $\dim(\mathscr{Y})\leq\dim(\mathscr{X})$.
\end{prop}

\begin{proof}
(1) follows immediately from the fact that $\O_{\mathscr{Y},x}$ is a quotient of $\O_{\mathscr{X},x}$; see \ref{prop-closedimmrigid2}.
(2) follows immediately from the definition of the dimensions.
\end{proof}

\begin{prop}\label{prop-dimensionatpoint1}
Let $\mathscr{X}$ be a rigid space, and $\{\mathscr{U}_{\alpha}\}_{\alpha\in L}$ a covering of $\mathscr{X}_{\ad}$ {\rm (\ref{dfn-admissiblesite3gensmall})}.
Then we have
$$
\dim(\mathscr{X})=\sup_{\alpha\in L}\dim(\mathscr{U}_{\alpha}).
$$
\end{prop}

\begin{proof}
By \ref{prop-dimension2} we have $\sup_{\alpha\in L}\dim(\mathscr{U}_{\alpha})\leq\dim(\mathscr{X})$.
On the other hand, for any $x\in\ZR{\mathscr{X}}$ there exists $\mathscr{U}_{\alpha}$ containing $x$.
Then $\dim_x(\mathscr{X})=\dim_x(\mathscr{U}_{\alpha})\leq\dim(\mathscr{U}_{\alpha})$.
Hence we have $\dim(\mathscr{X})\leq\sup_{\alpha\in L}\dim(\mathscr{U}_{\alpha})$.
\end{proof}

\begin{prop}\label{prop-dimension3}
Let $\mathscr{X}$ be a universally Noetherian rigid space, and $x\in\ZR{\mathscr{X}}$ a point.
Then we have
$$
\dim_x(\mathscr{X})\leq\sup_{x\in\mathscr{U}}\dim\Gamma(\mathscr{U},\O_{\mathscr{U}}),
$$
where $\mathscr{U}$ in the right-hand side runs through all Stein affinoid neighborhoods of $x$ in $\mathscr{X}$, and the dimension $\dim\Gamma(\mathscr{U},\O_{\mathscr{U}})$ is the Krull dimension of the ring $\Gamma(\mathscr{U},\O_{\mathscr{U}})$.
\end{prop}

\begin{proof}
Let $\mathfrak{p}_0\subsetneq\mathfrak{p}_1\subsetneq\cdots\subsetneq\mathfrak{p}_n$ be a strictly increasing chain of prime ideals of $\O_{\mathscr{X},x}$.
Let $\mathscr{U}$ be an arbitrary Stein affinoid open neighborhood of $x$.
Set $\mathscr{U}=(\Spf A)^{\rig}$, where $A$ is a t.u.\ rigid-Noetherian ring, and define $B$ by $\Spec A\setminus V(I)=\Spec B$, where $I\subseteq A$ is an ideal of definition (that is, $B=\Gamma(\mathscr{U},\O_{\mathscr{U}})$).
Consider $P_{i,\mathscr{U}}=\ker(B\rightarrow\O_{\mathscr{X},x}/\mathfrak{p}_i)$ for $i=0,\ldots,n$, which is a prime ideal of $B$.
Since
$$
\varinjlim_{x\in\mathscr{U}}P_{i,\mathscr{U}}=\mathfrak{p}_i,
$$
where the left-hand limit is taken along all Stein affinoid open neighborhoods of $x$, there exists a sufficiently small $\mathscr{U}$ such that
$$
P_{0,\mathscr{U}}\subsetneq P_{1,\mathscr{U}}\subsetneq\cdots\subsetneq P_{n,\mathscr{U}},
$$
from which the assertion follows.
\end{proof}

\subsubsection{Germs of rigid subspaces}\label{subsub-germsrigidsubspaces}
Let $\mathscr{X}$ be a locally universally Noetherian rigid space\index{rigid space!universally Noetherian rigid space@universally Noetherian ---!locally universally Noetherian rigid space@locally --- ---} (\ref{dfn-universallyadhesiverigidspaces}), and $x\in\ZR{\mathscr{X}}$ a point.
We define an equivalence relation on the set of all rigid subspaces\index{rigid subspace} of $\mathscr{X}$ as follows: two rigid subspaces $\mathscr{Z}_1$ and $\mathscr{Z}_2$ are equivalent if and only if there exists a quasi-compact open neighborhood $\mathscr{U}$ of $x$ such that $\mathscr{Z}_1\cap\mathscr{U}=\mathscr{Z}_2\cap\mathscr{U}$ as rigid subspaces of $\mathscr{U}$.
An equivalence class with respect to this equivalence relation is called a {\it germ of a rigid subspace at $x$}\index{germ of a rigid subspace@germ (of a rigid subspace)}.
Given a rigid subspace $\mathscr{Z}\subseteq\mathscr{X}$, the associated germ at $x$ will be denoted by 
$$
\mathscr{Z}_x.
$$
Clearly, for a rigid subspace $\mathscr{Z}$ and a quasi-compact open neighborhood $\mathscr{U}$ of $x$, we have $\mathscr{Z}_x=(\mathscr{Z}\cap\mathscr{U})_x$.
Hence, when considering a germ $\mathscr{Z}_x$, one can always replace $\mathscr{Z}$ by a rigid subspace of the form $\mathscr{Z}\cap\mathscr{U}$, which, moreover, can be assumed to be a {\it closed} rigid subspace of $\mathscr{U}$.

A germ $\mathscr{Z}_x$ is said to be {\it reduced}\index{germ of a rigid subspace@germ (of a rigid subspace)!reduced germ of a rigid subspace@reduced ---} if it comes from a reduced rigid subspace $\mathscr{Z}$ (cf.\ \ref{dfn-reducedrigidspaces}).

Let $\mathscr{Z}_{1,x}$ and $\mathscr{Z}_{2,x}$ be two germs of rigid subspaces at $x$.
We say that a germ $\mathscr{Z}_{1,x}$ is contained in $\mathscr{Z}_{2,x}$, written $\mathscr{Z}_{1,x}\subseteq\mathscr{Z}_{2,x}$, if there exists a quasi-compact open neighborhood $\mathscr{U}$ of $x$ such that $\mathscr{Z}_1\cap\mathscr{U}$ is a rigid subspace of $\mathscr{Z}_2\cap\mathscr{U}$.
In type (V) case, by \S\ref{subsub-closedimmrigidred}, any germ contains the uniquely determined {\it reduced model}\index{reduced model of a rigid space@reduced model (of a rigid space)}\index{rigid space!reduced model of a rigid space@reduced model of a ---}.

Let $\mathscr{Z}_x$ be a germ of a rigid subspace at $x$, and suppose $\mathscr{Z}$ is a closed subspace of a quasi-compact open neighborhood $\mathscr{U}$ of $x$.
Then we have the defining coherent ideal $\mathscr{J}$ of $\O_{\mathscr{U}}$.
Taking the stalk at $x$, we get a finitely generated ideal $\mathscr{J}_x$ of $\O_{\mathscr{X},x}$.
We denote this ideal by $I(\mathscr{Z}_x)$ and call the {\it ideal of $\mathscr{Z}_x$}.

\begin{prop}\label{prop-germsideal}
Let $\mathscr{X}$ be a locally universally Noetherian rigid space\index{rigid space!universally Noetherian rigid space@universally Noetherian ---!locally universally Noetherian rigid space@locally --- ---}.
Then the correspondence 
$$
\mathscr{Z}_x\longmapsto I(\mathscr{Z}_x)
$$
establishes a bijection between the set of all germs of rigid subspaces at $x$ and the set of all finitely generated ideals of $\O_{\mathscr{X},x}$.
Moreover, if $\mathscr{Z}_{1,x}$ and $\mathscr{Z}_{2,x}$ are two germs at $x$, then $\mathscr{Z}_{1,x}\subseteq\mathscr{Z}_{2,x}$ if and only if $I(\mathscr{Z}_{1,x})\supseteq I(\mathscr{Z}_{2,x})$.
In type {\rm (V)} case, reduced germs correspond to finitely generated radical ideals.
\end{prop}

\begin{proof}
The inverse to the above-mentioned correspondence is constructed as follows (cf.\ \S\ref{subsub-closedsubspacesofaffinoids}).
Let $J$ be a finitely generated ideal of $\O_{\mathscr{X},x}$.
Then there exist an affinoid open neighborhood $\mathscr{U}=(\Spf A)^{\rig}$ of $x$ (where $A$ is a t.u.\ rigid-Noetherian ring) and a finitely generated ideal $J'$ of $A[\frac{1}{a}]$ such that $J=J'\O_{\mathscr{X},x}$ (cf.\ \ref{prop-descriptionlocalrings} (1)). 
Take a finitely generated ideal $\til{J}$ of $A$ such that $\til{J}A[\frac{1}{a}]=J'$.
Let $Z\hookrightarrow\Spf A$ be the closed immersion of finite presentation corresponding to $\til{J}$, and consider the associated rigid space $\mathscr{Z}=Z^{\rig}$, which is the closed subspace of $\mathscr{U}$.
The germ $\mathscr{Z}_x$ of $\mathscr{Z}$ is easily seen to be independent of the choice of $\mathscr{U}$, $A$, and $\til{J}$, and thus we have the desired inverse to the correspondence $\mathscr{Z}_x\mapsto I(\mathscr{Z}_x)$.

In type (V) case, we know, thanks to \ref{thm-noetherness}, that the local ring $\O_{\mathscr{X},x}$ is Noetherian, and hence every ideal corresponds to a germ. 
Then it follows from \ref{prop-reducedmodelrigid} that a germ $\mathscr{Z}_x$ is reduced if and only if the corresponding ideal $I(\mathscr{Z}_x)$ is a radical ideal.
\end{proof}

Let $\mathscr{Z}_{1,x}$ and $\mathscr{Z}_{2,x}$ be two germs of rigid subspaces at $x$.
We may assume that these germs come from closed rigid subspaces $\mathscr{Z}_1$ and $\mathscr{Z}_2$ of a quasi-compact open neighborhood $\mathscr{U}$ of $x$.
Then the intersection $\mathscr{Z}_{1,x}\cap\mathscr{Z}_{2,x}$ is defined to be the germ of $\mathscr{Z}_1\times_{\mathscr{U}}\mathscr{Z}_2$ at $x$ (cf.\ \ref{prop-closedimmrigid5} (4)).
It is clear that we have the equality
$$
I(\mathscr{Z}_{1,x}\cap\mathscr{Z}_{2,x})=I(\mathscr{Z}_{1,x})+I(\mathscr{Z}_{2,x})
$$
of finitely generated ideals of $\O_{\mathscr{X},x}$.
One can define the union $\mathscr{Z}_{1,x}\cup\mathscr{Z}_{2,x}$ as the germ corresponding to the product $I(\mathscr{Z}_{1,x})\cdot I(\mathscr{Z}_{2,x})$.

A germ $\mathscr{Z}_x$ at $x$ is said to be {\it prime}\index{germ of a rigid subspace@germ (of a rigid subspace)!prime germ of a rigid subspace@prime ---} if the corresponding ideal $I(\mathscr{Z}_x)$ is a prime ideal of $\O_{\mathscr{X},x}$.
In type (V) case, this is equivalent to that $\mathscr{Z}_x$ is reduced and irreducible (that is, non-empty and whenever $\mathscr{Z}_x=\mathscr{Z}_{1,x}\cup\mathscr{Z}_{2,x}$ by reduced germs $\mathscr{Z}_{1,x},\mathscr{Z}_{2,x}$, we have either $\mathscr{Z}_x=\mathscr{Z}_{1,x}$ or $\mathscr{Z}_x=\mathscr{Z}_{2,x}$).

Let $\mathscr{X}$ be a locally universally Noetherian rigid space\index{rigid space!universally Noetherian rigid space@universally Noetherian ---!locally universally Noetherian rigid space@locally --- ---}, and $x\in\ZR{\mathscr{X}}$ a point such that the local ring $\O_{\mathscr{X},x}$ is Noetherian; the last assumption is always valid if $\mathscr{X}$ is of type (V) or of (N); see \ref{thm-noetherness}.
A {\it chain of prime germs of rigid subspaces at $x$}\index{chain of prime germs of rigid subspaces} is the diagram
$$
\mathscr{Z}_{0,x}\longhookrightarrow\mathscr{Z}_{1,x}\longhookrightarrow\cdots\longhookrightarrow\mathscr{Z}_{n,x}\leqno{(\ast)}
$$
consisting of prime germs at $x$ such that for $i=0,\ldots,n-1$ we have $\mathscr{Z}_{i,x}\neq\mathscr{Z}_{i+1,x}$.
The number $n$ in $(\ast)$ is called the {\it length}\index{chain of prime germs of rigid subspaces!length of a chain of prime germs of rigid subspaces@length of a ---} of the chain $(\ast)$.

\begin{prop}\label{prop-cordimensionfinite}
Let $\mathscr{X}$ be a rigid space, and $x\in\ZR{\mathscr{X}}$ a point such that he local ring $\O_{\mathscr{X},x}$ is Noetherian $($e.g.\ $\mathscr{X}$ is of type {\rm (V)} or of {\rm (N)}$)$.
Then the local ring $\O_{\mathscr{X},x}$ is Noetherian, and we have $\dim_x(\mathscr{X})<+\infty$.
Moreover, $\dim_x(\mathscr{X})$ coincides with the supremum of the lengths of all chains of prime germs of rigid subspaces at $x$.
\end{prop}

\begin{proof}
The local ring $\O_{\mathscr{X},x}$ is Noetherian due to \ref{thm-noetherness}, and hence is of finite Krull dimension.
The other assertion is clear by the definition of prime germs.
\end{proof}

\subsubsection{Dimension of rigid spaces of type (V) or of type (N)}\label{subsub-comparisondim}
\begin{thm}[Comparison of the dimensions for affinoids]\label{thm-dimrigidtypeVN1}
Let $\mathscr{X}=(\Spf A)^{\rig}$ be an affinoid of type {\rm (V)}\index{rigid space!rigid space of typeV@--- of type (V)} or of type {\rm (N)}\index{rigid space!rigid space of typeN@--- of type (N)}. 
Then for any closed classical point\index{point!classical point@classical ---}\index{classical point} $x\in\ZR{\mathscr{X}}^{\cl}$ {\rm (\ref{dfn-classicalpoint})} we have
$$
\dim_x(\mathscr{X})=\dim_{s(x)}(s(\mathscr{X})),
$$
where $s(\mathscr{X})$ is the Noetherian scheme associated to $\mathscr{X}$ {\rm (\S\ref{sub-associatedschemes})}, and $\dim_{s(x)}(s(\mathscr{X}))$ denotes the dimension of the Noetherian scheme $s(\mathscr{X})$ in the usual sense $($cf.\ {\rm \cite[$\mathbf{IV}$, \S5.1]{EGA}}$)$.
\end{thm}

\begin{proof}
By \ref{prop-complocringclassicalcomp} we have
$$
\dim(\O_{\mathscr{X},x})=\dim(\widehat{\O}_{\mathscr{X},x})=\dim(\widehat{\O}_{s(\mathscr{X}),s(x)})=\dim_{s(x)}(s(\mathscr{X})),
$$
which shows the assertion.
\end{proof}

\begin{prop}\label{prop-dimensionclassicalenough}
Let $\mathscr{X}$ be a rigid space of type {\rm (V)}\index{rigid space!rigid space of typeV@--- of type (V)} or of type {\rm (N)}\index{rigid space!rigid space of typeN@--- of type (N)}. 
Then we have 
$$
\dim(\mathscr{X})=\sup_{x\in\ZR{\mathscr{X}}^{\cl}}\dim_x(\mathscr{X})
$$
\end{prop}

\begin{proof}
The inequality $\sup_{x\in\ZR{\mathscr{X}}^{\cl}}\dim_x(\mathscr{X})\leq\sup_{x\in\ZR{\mathscr{X}}}\dim_x(\mathscr{X})$ is trivial.
We need to show the opposite inequality. 
For any $x\in\ZR{\mathscr{X}}$, there exists a Stein affinoid $\mathscr{U}=(\Spf A)^{\rig}$ neighborhood of $x$ such that 
$$
\dim_x(\mathscr{X})\leq\dim\Gamma(\mathscr{U},\O_{\mathscr{U}})=\sup_{z\in s(\mathscr{U})^{\cl}}\dim\O_{s(\mathscr{U}),z}
$$
due to \ref{prop-dimension3}.
Since, as in \ref{prop-classicalpointsexist} (1), any closed point of $s(\mathscr{U})$ gives rise to a closed classical point of $\mathscr{U}$, and hence a classical point of $\mathscr{X}$, we have, by \ref{thm-dimrigidtypeVN1},
$$
\sup_{z\in s(\mathscr{U})^{\cl}}\dim\O_{s(\mathscr{U}),z}\leq\sup_{y\in\ZR{\mathscr{X}}^{\cl}}\dim\O_{\mathscr{X},y},
$$
from which we deduce the desired inequality.
\end{proof}

\begin{cor}\label{cor-dimaffinoids}
Let $\mathscr{X}=(\Spf A)^{\rig}$ be an affinoid of type {\rm (V)}, or an affinoid of type {\rm (N)} having a distinguished Noetherian formal model $X$, together with an ideal of definition $\mathscr{I}$, such that the scheme $X_0=(X,\O_X/\mathscr{I})$ is Jacobson.
Then we have 
$$
\dim(\mathscr{X})=\dim(s(\mathscr{X})).
$$
\end{cor}

\begin{proof}
If $\mathscr{X}$ is a finite type affinoid over $(\Spf V)^{\rig}$, where $V$ is an $a$-adically complete valuation ring, we may replace $V$ by its height one localization $V_{\mathfrak{p}}$, where $\mathfrak{p}=\sqrt{(a)}$ is the associated height one prime, and thus may assume that $\mathscr{X}$ is of type {\rm ($\mathrm{V_{\R}}$)}.
Then, by \ref{cor-classicalpointsexist0}, \ref{prop-dimensionclassicalenough} and \ref{thm-dimrigidtypeVN1}, we have
$$
\dim(\mathscr{X})=\sup_{x\in\ZR{\mathscr{X}}^{\cl}}\dim_x(\mathscr{X})=\sup_{z\in s(\mathscr{X})^{\cl}}\dim\O_{s(\mathscr{X}),z}=\dim s(\mathscr{X}),
$$
as desired.
\end{proof}

\begin{cor}\label{cor-dimensionclosedfiber}
Let $\mathscr{X}$ be a coherent rigid space of type {\rm ($\mathrm{V_{\R}}$)}\index{rigid space!rigid space of typeV1@--- of type ($\mathrm{V_{\R}}$)}, and $X$ a distinguished formal model of $\mathscr{X}$, which is a finite type formal scheme over an $a$-adically complete valuation ring $V$ of height $1$.
Denote by $X_k=X\otimes_Vk$ the closed fiber of $X$, which is a finite type scheme over the residue field $k=V/\m_V$.
Then we have
$$
\dim(\mathscr{X})=\dim(X_k).
$$
\end{cor}

\begin{proof}
We may assume $\mathscr{X}=(\Spf A)^{\rig}$, where $A$ is an $a$-torsion free topologically of finite type $V$-algebra.
Take a finite injection $V\dl X_1,\ldots,X_d\dr\hookrightarrow A$ with the $V$-flat kernel ({\bf \ref{ch-pre}}.\ref{thm-noethernormalizationtype(V)}).
By \ref{cor-dimaffinoids}, we have $\dim(\mathscr{X})=\dim(\Spec A[\frac{1}{a}])=d$.
On the other hand, since we have a finite injection $k[X_1,\ldots,X_d]\hookrightarrow A_k=A\otimes_Vk$, we have $\dim(X_k)=d$, thereby the claim.
\end{proof}

\subsubsection{Calculation of the dimension}\label{subsub-dimensioncal}
\begin{prop}\label{prop-dimensioncal}
Let $V$ be an $a$-adically complete valuation ring, and consider the unit disk \index{unit disk} $\D^n_{\mathscr{S}}$ over $\mathscr{S}=(\Spf V)^{\rig}$ $($cf.\ {\rm \S\ref{subsub-unitdisk}}$)$.
We have
$$
\dim(\D^n_{\mathscr{S}})=n.
$$
\end{prop}

Once this proposition is established, one can then apply Noether normalization theorem ({\bf \ref{ch-pre}}.\ref{thm-noethernormalizationtype(V)}) to compute dimensions of rigid spaces of type (V).

To show the proposition, we consider the ring $V\dl X_1,\ldots,X_n\dr$ of restricted formal power series\index{restricted formal power series} ({\bf \ref{ch-pre}}, \S\ref{sub-powerseries}); we have $\D^n_{\mathscr{S}}=(\Spf V\dl X_1,\ldots,X_n\dr)^{\rig}$.
By \ref{cor-dimaffinoids} we have
$$
\dim(\D^n_{\mathscr{S}})=\dim(s(\D^n_{\mathscr{S}}))=\dim(\Spec V\dl X_1,\ldots,X_n\dr{\textstyle [\frac{1}{a}]}).
$$
Since $V\dl X_1,\ldots,X_n\dr{\textstyle [\frac{1}{a}]}=V_{\mathfrak{p}}\dl X_1,\ldots,X_n\dr{\textstyle [\frac{1}{a}]}$, where $\mathfrak{p}=\sqrt{(a)}$ (cf.\ {\bf \ref{ch-pre}}.\ref{cor-compval2006ver23}), we may assume that $V$ is of height one.
Now the ring $V\dl X_1,\ldots,X_n\dr{\textstyle [\frac{1}{a}]}$, the Tate algebra\index{Tate, J.}\index{algebra!Tate algebra@Tate ---} (cf.\ {\bf \ref{ch-pre}}, \S\ref{subsub-classicalaffinoidalgebras}), is Noetherian, and hence it suffices to show:
\begin{prop}\label{prop-dimensioncal2}
Let $V$ be an $a$-adically complete valuation ring of height one, and $K=V[\frac{1}{a}]$ $(=\Frac(V);$ cf.\ {\rm {\bf \ref{ch-pre}}.\ref{prop-sep}}$)$

{\rm (1)} Maximal ideals of $V\dl X_1,\ldots,X_n\dr{\textstyle [\frac{1}{a}]}$ are exactly the kernels of $K$-algebra homomorphisms of the form $V\dl X_1,\ldots,X_n\dr{\textstyle [\frac{1}{a}]}\rightarrow K'$, where $K'$ is a finite extension of $K$.

{\rm (2)} For any closed point $z$ of $\Spec V\dl X_1,\ldots,X_n\dr{\textstyle [\frac{1}{a}]}$ we have
$$
\dim_z(\Spec V\dl X_1,\ldots,X_n\dr{\textstyle [\frac{1}{a}]})=n.
$$
\end{prop}

\begin{proof}
The proof of (1) has been given in {\bf \ref{ch-pre}}.\ref{cor-Tatealgebraprime}, and that of (2) in {\bf \ref{ch-pre}}.\ref{prop-tatealgebradim}.
\end{proof}

\begin{cor}\label{cor-dimensioncal4}
Let $A$ be a topologically finitely generated algebra over an $a$-adically complete valuation ring $V$.
Let $V'$ be an $a'$-adically complete valuation ring, and $V\rightarrow V'$ is an adic homomorphism.
Set $A'=A\widehat{\otimes}_VV'$.
Then we have
$$
\dim((\Spf A)^{\rig})=\dim((\Spf A')^{\rig}).
$$
\end{cor}

\begin{proof}
We may assume that $A$ is $a$-torsion free (hence $V$-flat).
Let $\mathfrak{p}=\sqrt{(a)}$ be the associated height one prime of $V$, and consider the height one valuation ring $V_{\mathfrak{p}}$.
Since $V_{\mathfrak{p}}\otimes_VV'$ is nothing but the height one localization of $V'$, we may assume that $V$ and $V'$ are of height one.
By Noether normalization ({\bf \ref{ch-pre}}.\ref{thm-noethernormalizationtype(V)}) there exists an injective finite map $V\dl X_1,\ldots,X_n\dr\hookrightarrow A$.
Then we have $\dim((\Spf A)^{\rig})=\dim(\Spec A[\frac{1}{a}])=n$.
On the other hand, by applying $\widehat{\otimes}_VV'$, we have the injective finite homomorphism $V'\dl X_1,\ldots,X_n\dr\hookrightarrow A'$ ({\bf \ref{ch-formal}}.\ref{prop-finitemorform2} (4)); notice that the map $V\rightarrow V'$ is injective (since it is adic ({\bf \ref{ch-pre}}.\ref{prop-maxspe4})) and hence is flat.
Therefore, we have $\dim((\Spf A')^{\rig})=n$, as desired.
\end{proof}

The following corollary follows immediately from \ref{cor-dimensioncal4}:
\begin{cor}\label{cor-dimensioncal5}
Let $\mathscr{X}$ be a rigid space locally of finite type over $\mathscr{S}=(\Spf V)^{\rig}$, where $V$ is an $a$-adically complete valuation ring.
Let $V'$ be an $a'$-adically complete valuation ring, and $V\rightarrow V'$ an adic map$;$ set $\mathscr{S}'=(\Spf V')^{\rig}$.
Then we have
$$
\dim(\mathscr{X})=\dim(\mathscr{X}\times_{\mathscr{S}}\mathscr{S}').\eqno{\square}
$$
\end{cor}

\begin{cor}\label{cor-dimensioncal10}
Let $V$ be an $a$-adically complete valuation ring of height one, and $A$ an $a$-torsion free topologically finitely generated $V$-algebra, which is further assumed to be an integral domain.
Set $\mathscr{X}=(\Spf A)^{\rig}$, and $B=A[\frac{1}{a}]=\Gamma(\mathscr{X},\O_{\mathscr{X}})$.
Then $\dim_x(\mathscr{X})$ is constant for $x\in\ZR{\mathscr{X}}^{\cl}$, and is equal to $\dim(\mathscr{X})$, which is further equal to the Krull dimension $\dim B$.
\end{cor}

\begin{proof}
In view of \ref{thm-dimrigidtypeVN1}, we need to show that $\dim_z(\Spec B)$ is constant for $z\in(\Spec B)^{\cl}$; the other assertions follow from \ref{prop-dimensionclassicalenough} and \ref{cor-dimaffinoids}.
The claim is clear in the case $A=V\dl X_1,\ldots,X_n\dr$ by \ref{prop-dimensioncal2} (2).
In general, by Noether normalization\index{Noether normalization} ({\bf \ref{ch-pre}}.\ref{thm-northernormaclassaff}), we have a finite injective map $K\dl X_1,\ldots,X_n\dr\hookrightarrow B$.
Then, since $K\dl X_1,\ldots,X_n\dr$ is integrally closed ({\bf \ref{ch-pre}}.\ref{prop-exerTateUFD}), the ring extension $B/K\dl X_1,\ldots,X_n\dr$ satisfies going-up and going-down properties (cf.\ \cite[Theorem 9.4]{Matsu}), from which the desired constancy of the dimension follows.
\end{proof}

\begin{cor}\label{cor-dimensioncal101}
Let $V$ and $\mathscr{X}$ be as in {\rm \ref{cor-dimensioncal10}}, and $\mathscr{U}\subseteq\mathscr{X}$ a non-empty open subspace of $\mathscr{X}$.
Then we have
$$
\dim(\mathscr{X})=\dim(\mathscr{U}).\eqno{\square}
$$
\end{cor}

\subsubsection{GAGA comparison of the dimensions}\label{subsub-dimaffinoidscor}
\index{GAGA!GAGA comparison of dimensions@--- comparison of dimensions|(}
\begin{prop}\label{prop-compdimfintype}
Let $V$ be an $a$-adically complete valuation ring, and $A$ a $V$-algebra of finite type.
Consider the morphism of schemes
$${\textstyle 
\rho\colon\Spec\widehat{A}[\frac{1}{a}]\longrightarrow\Spec A[\frac{1}{a}],}
$$
where $\widehat{A}$ is the $a$-adic completion of $A$.

{\rm (1)} The morphism $\rho$ maps closed points of $\Spec\widehat{A}[\frac{1}{a}]$ injectively to closed points of $\Spec A[\frac{1}{a}]$.

{\rm (2)} For any closed point $x\in\Spec\widehat{A}[\frac{1}{a}]$ the induced map between the completed local rings is an isomorphism$:$
$$
\widehat{\O}_{\Spec\widehat{A}[\frac{1}{a}],x}\cong\widehat{\O}_{\Spec A[\frac{1}{a}],\rho(x)},
$$
where the completions are taken with respect to the adic topology defined by the maximal ideals.
In particular, we have 
$${\textstyle 
\dim_x(\Spec\widehat{A}[\frac{1}{a}])=\dim_{\rho(x)}(\Spec A[\frac{1}{a}]).}
$$
\end{prop}

\begin{proof}
Clearly, we may assume that $V$ is of height one and that $A$ is $a$-torsion free (that is, flat over $V$).
Write $A=V[X_1,\ldots,X_n]/J$ by an $a$-saturated (hence finitely generated) ideal $J$.
We have $\widehat{A}=V\dl X_1,\ldots,X_n\dr/JV\dl X_1,\ldots,X_n\dr$.
Hence the map in question is rewritten as
\begin{equation*}
\begin{split}
\Spec V\dl X_1,\ldots,X_n\dr{\textstyle [\frac{1}{a}]}/&JV\dl X_1,\ldots,X_n\dr{\textstyle [\frac{1}{a}]}\\
&{\textstyle \longrightarrow\Spec V[X_1,\ldots,X_n][\frac{1}{a}]/JV[X_1,\ldots,X_n][\frac{1}{a}].}
\end{split}
\end{equation*}
By \ref{prop-dimensioncal2} (1) maximal ideals of $\widehat{A}[\frac{1}{a}]$ are exactly the kernels of the map of the form $\widehat{A}[\frac{1}{a}]\rightarrow K'$, where $K'$ is a finite extension of $K=V[\frac{1}{a}]$.
Hence it defines a maximal ideal of $A[\frac{1}{a}]$, which is a $K$-algebra of finite type.
Since the map as above is uniquely determined by the images of $X_i$'s, this correspondence between maximal ideals is injective.
Thus we have shown (1).
On the other hand, we already know that (2) is true when $A=V[X_1,\ldots,X_n]$, since, in this case, the completed local rings in question are isomorphic to the ring of formal power series.
By taking the quotient modulo $J$, we get the desired assertion (since $J$ is finitely generated).
\end{proof}

Now let us consider the following situation.
Let $V$ be an $a$-adically complete valuation ring of height one, and set $K=\Frac(V)=V[\frac{1}{a}]$.
Let 
$$
f\colon X\longrightarrow\Spec K
$$
be a $K$-scheme of finite type.
Then by the GAGA functor (\S\ref{subsub-GAGAfunctornonsep}) one can consider the rigid space
$$
f^{\an}\colon X^{\an}\longrightarrow (\Spf V)^{\rig}
$$
of finite type over $V$.
\begin{prop}\label{prop-dimaffinoidsGAGA111}
For any classical point $x\in\ZR{X^{\an}}^{\cl}$ of $X^{\an}$ we have the canonical isomorphism
$$
\widehat{\O}_{X,\rho_X(x)}\cong\widehat{\O}_{X^{\an},x}
$$
of complete local rings.
$($Here $\rho_X\colon\ZR{X^{\an}}\rightarrow X$ is the comparison map $(${\rm \S\ref{subsub-GAGAcommap}}$).)$
In particular, we have
$$
\dim_x(X^{\an})=\dim_{\rho_X(x)}(X).
$$
\end{prop}

\begin{proof}
Since the question is local, we may assume that $X$ is affine; replacing by a compactification, we may furthermore assume that $X$ is proper.
One can choose an affinoid neighborhood $\mathscr{U}=(\Spf\widehat{A})^{\rig}$ of $x$, which comes from an affine open neighborhood $\Spec A\subseteq X$ of $\rho(x)$.
Then the question is to show that 
$$
\widehat{\O}_{\mathscr{U},x}\cong\widehat{\O}_{\Spec A[\frac{1}{a}],\rho(x)}.
$$
But this follows from \ref{prop-complocringclassicalcomp} and \ref{prop-compdimfintype} (2), since $s(\mathscr{U})=\Spec\widehat{A}[\frac{1}{a}]$.
\end{proof}

\begin{thm}[GAGA comparison of dimensions]\label{thm-dimaffinoidsGAGA}
Let $V$ be an $a$-adically complete valuation ring, and $X$ a finite type scheme over $K=\Frac(V)=V[\frac{1}{a}]$.
Then we have
$$
\dim(X^{\an})=\dim(X).
$$
\end{thm}

To show the theorem, we may assume that $V$ is of height one.
Then in view of \ref{prop-dimaffinoidsGAGA111} the theorem follows from the following lemma, which follows easily from \ref{thm-affinoidvaluedpoints}:
\begin{lem}\label{lem-dimaffinoidsGAGApt}
In the above situation, the comparison map $\rho_X\colon\ZR{X^{\an}}\rightarrow X$ maps the set of all classical points $\ZR{X^{\an}}^{\cl}$ bijectively onto the set of all closed points $X^{\cl}$. \hfill$\square$
\end{lem}
\index{GAGA!GAGA comparison of dimensions@--- comparison of dimensions|)}

\subsubsection{Dimension function}\label{subsub-dimensionfunction}
\begin{dfn}\label{dfn-dimensionfunction}{\rm 
Let $\mathscr{X}$ be a rigid space of type {\rm (V)} or of type {\rm (N)}.
We define the function 
$$
d=d_{\mathscr{X}}\colon\ZR{\mathscr{X}}\longrightarrow\Z_{\geq 0}\cup\{-\infty\},
$$
which we call the {\it dimension function} on $\mathscr{X}$, by
$$
d(x)=\inf_{x\in\ZR{\mathscr{U}}}\dim(\mathscr{U})
$$
for $x\in\ZR{\mathscr{X}}$, where $\mathscr{U}$ runs through all open neighborhood $\mathscr{U}$ of $x$.}
\end{dfn}

\begin{prop}\label{prop-dimensionfunction1}
Let $V$ be an $a$-adically complete valuation ring of height one, and $A$ an $a$-torsion free topologically finitely generated $V$-algebra.
Suppose that $\mathscr{X}=(\Spf A)^{\rig}$ is reduced $($that is, the nilpotent radical of $B=A[\frac{1}{a}]$ is zero; cf.\ {\rm \ref{prop-reducedmodelrigid2}}$)$.
Let 
$$
\mathscr{X}=\bigcup^r_{i=1}\mathscr{X}_i
$$
be the irreducible decomposition $($as in {\rm \S\ref{subsub-closedimmrigidred}}$)$.
Then, for any $x\in\ZR{\mathscr{X}}$, we have
$$
d(x)=\max\{\dim(\mathscr{X}_i):i=1,\ldots,r,\ x\in\ZR{\mathscr{X}_i}\}.
$$
\end{prop}

\begin{proof}
Since closed rigid subspaces are overconvergent closed subsets (\ref{prop-opencomplement}), we have $\{i:x\in\ZR{\mathscr{X}_i}\}=\{i:G_x\cap\ZR{\mathscr{X}_i}\neq\emptyset\}$ (where $G_x$ denotes, as before, the set of all generizations of $x$), which we denote by $I_x$.
Since $I_x$ is the set of indices $i=1,\ldots,r$ such that $\ZR{\mathscr{U}}\cap\ZR{\mathscr{X}_i}\neq\emptyset$ for any open neighborhood $\mathscr{U}$ of $x$, one can choose a sufficiently small affinoid open neighborhood $\mathscr{U}$ of $x$ such that (i) $d(x)=\dim(\mathscr{U})$ and (ii) $I_x=\{i\:\ZR{\mathscr{U}}\cap\ZR{\mathscr{X}_i}\neq\emptyset\}$ hold.
By \ref{cor-dimaffinoids}, we have $d(x)=\dim\Gamma(\mathscr{U},\O_{\mathscr{U}})$ (Krull dimension), and hence
\begin{equation*}
\begin{split}
d(x)&=\max_{i\in I_x}\dim\Gamma(\mathscr{X}_i\cap\mathscr{U},\O_{\mathscr{X}_i\cap\mathscr{U}})\\
&=\max_{i\in I_x}\dim\Gamma(\mathscr{X}_i,\O_{\mathscr{X}_i})=\max_{i\in I_x}\dim(\mathscr{X}_i),
\end{split}
\end{equation*}
where the second equality is due to \ref{cor-dimensioncal101}, and the third one again by \ref{cor-dimaffinoids}.
\end{proof}

\begin{cor}\label{cor-dimensionfunction1}
Let $\mathscr{X}$ be a rigid space of type {\rm ($\mathrm{V_{\R}}$)}.
The dimension function $d=d_{\mathscr{X}}$ is upper-semi continuous on $\ZR{\mathscr{X}}$; moreover, $\{x\in\ZR{\mathscr{X}}:d(x)\geq c\}$ for any $c\in\R$ is the underlying topological space of a closed rigid subspace of $\mathscr{X}$.
\end{cor}

\subsection{Codimension}\label{sub-codimention}
\begin{dfn}\label{dfn-codimensionatpoint}{\rm 
Let $\mathscr{X}$ be a locally universally Noetherian rigid space\index{rigid space!universally Noetherian rigid space@universally Noetherian ---!locally universally Noetherian rigid space@locally --- ---}.

(1) Let $x\in\ZR{\mathscr{X}}$ be a point, and $\mathscr{Y}_x$ a germ of a rigid subspace $\mathscr{Y}$ at $x$.
The {\em codimension of $\mathscr{Y}_x$ in $\mathscr{X}$ at $x$}\index{codimension of a rigid subspace@codimension (of a rigid subspace)}, denoted by 
$$
\codim_x(\mathscr{Y},\mathscr{X}),
$$
is the number defined as follows.
If $\mathscr{Y}_x$ is a prime germ, then it is defined to be the supremum of the lengths of all possible chains of prime germs of rigid subspaces at $x$
$$
\mathscr{Z}_{0,x}\longhookrightarrow\mathscr{Z}_{1,x}\longhookrightarrow\cdots\longhookrightarrow\mathscr{Z}_{n,x},
$$
where $\mathscr{Z}_{0,x}=\mathscr{Y}_x$; in general, it is the infimum of the codimensions in $\mathscr{X}$ of $\mathscr{Z}_x$ at $x$, where $\mathscr{Z}_x$ runs in the set of all prime germs that are contained in $\mathscr{Y}_x$.

(2) Let $\mathscr{Y}$ be a closed subspace of $\mathscr{X}$.
Then the {\em codimension of $\mathscr{Y}$ in $\mathscr{X}$}\index{codimension of a rigid subspace@codimension (of a rigid subspace)}, denoted by $\codim(\mathscr{Y},\mathscr{X})$, is defined by
$$
\codim(\mathscr{Y},\mathscr{X})=\inf_{x\in\ZR{\mathscr{X}}}\codim_x(\mathscr{Y},\mathscr{X}).
$$}
\end{dfn}

\subsection{Relative dimension}\label{sub-relativedim}
\begin{prop}\label{prop-relativedim1}
Let $\varphi\colon\mathscr{X}\rightarrow\mathscr{Y}$ be a locally of finite type morphism between locally universally Noetherian rigid spaces\index{rigid space!universally Noetherian rigid space@universally Noetherian ---!locally universally Noetherian rigid space@locally --- ---}, and $y\in\ZR{\mathscr{Y}}$.
Let $\alpha\colon\mathscr{S}=(\Spf V)^{\rig}\rightarrow\mathscr{Y}$ be a rigid point such that $\Spf V\rightarrow\ZR{\mathscr{Y}}$ maps the closed point to $y$.
Then the number $\dim(\mathscr{X}\times_{\mathscr{Y}}\mathscr{S})$ depends only on $\varphi$ and $y$ and does not depend on the choice of the rigid point $\alpha$.
\end{prop}

\begin{proof}
Since any rigid point factors through the associated rigid point $\alpha_y$ (\ref{dfn-ZRpoints32}), it suffices to show that the number in question is the same as the one with $\alpha$ replaced by $\alpha_y$.
However, this follows immediately from \ref{cor-dimensioncal5}.
\end{proof}

\begin{dfn}\label{dfn-relativedim}{\rm 
Let $\varphi\colon\mathscr{X}\rightarrow\mathscr{Y}$ be a morphism locally of finite type between rigid spaces, and $y\in\ZR{\mathscr{Y}}$.
Then the {\em relative dimension}\index{dimension of a rigid space@dimension (of a rigid space)!relative dimension of a rigid space@relative ---} of $\varphi$ at $y$, denoted by 
$$
\dim_y(\varphi)
$$
(or by $\dim_y\mathscr{X}$), is the number $\dim(\mathscr{X}\times_{\mathscr{Y}}\mathscr{S})$, where $\alpha\colon\mathscr{S}=(\Spf V)^{\rig}\rightarrow\mathscr{Y}$ is a rigid point such that $\Spf V\rightarrow\ZR{\mathscr{Y}}$ maps the closed point to $y$ (which is independent of the choice of $\alpha$ due to \ref{prop-relativedim1}).}
\end{dfn}

\begin{thm}[GAGA comparison of relative dimension]\label{thm-relativedimGAGA}\index{GAGA!GAGA comparison of dimensions@--- comparison of dimensions}
Let $S=\Spec A$ where $A$ is an adic ring with a finitely generated ideal of definition $I\subseteq A$, $D=V(I)$, and $U=S\setminus D$.
Let $X,Y$ be $U$-schemes of finite type, and $f\colon X\rightarrow Y$ a morphism over $U$.
Then for any $y\in\ZR{Y^{\an}}$ we have
$$
\dim_y(f^{\an})=\dim_{\rho_Y(y)}(f).
$$
\end{thm}

\begin{proof}
Let $\Spf\widehat{V}_y\rightarrow\ZR{Y^{\an}}$ be the rigid point associated to the point $y$.
Then $\rho_Y(y)$ is the image of $\Spec\widehat{V}_y[\frac{1}{a}]\rightarrow Y$, and hence what to prove is that the dimension of the rigid space $X^{\an}\times_{Y^{\an}}(\Spf\widehat{V}_y)$ of type (V) is equal to the dimension of the scheme $X\times_Y\Spec\widehat{V}_y[\frac{1}{a}]$ of finite type over the field $\widehat{V}_y[\frac{1}{a}]$.
Hence in view of \ref{prop-GAGAfunctor5x3} we may assume $B=V=\widehat{V}_y$ and $Y=\Spec V[\frac{1}{a}]$.
But the assertion in this case is nothing but \ref{thm-dimaffinoidsGAGA}.
\end{proof}



\section{Maximum modulus principle}\label{sec-pointsV1}
In this section we discuss the maximal modulus principle for coherent rigid spaces of type of type ($\mathrm{V_{\R}}$). 
The statement is given in \ref{thm-MMP} below.
To show the theorem, we first give a classification of points on rigid spaces of type ($\mathrm{V_{\R}}$) along the sprit of the classical classification of valuations as in {\bf \ref{ch-pre}}, \S\ref{sub-valexample}.
This part of the discussion may be interesting in its own right and can be read independently.
Especially, in the case of the unit disk (\S\ref{subsub-classifyingvaluationsdisk}), one finds a strong analogy between our classification and Berkovich's classification of points in \cite[1.4.4]{Berk1}.
The proof of the maximal modulus principle is based on the so-called spectral seminorm formula (\ref{thm-spectralseminormformula}), which says that the spectral seminorm takes the maximum value of the norms at at most finitely many divisorial points. 
\subsection{Classification of points}\label{sub-classifyingvaluations}
\subsubsection{Basic inequality}\label{subsub-classifyingvaluationsheights}
Let $V$ be an $a$-adically complete valuation ring of height one, $K=\Frac(V)$ the fractional field, and $k=V/\m_V$ the residue field.
We denote by $\Gamma$ the value group of $V$, that is, $\Gamma=K^{\times}/V^{\times}$.

Let $\mathscr{X}$ be a coherent rigid space of finite type over $\mathscr{S}=(\Spf V)^{\rig}$.
For any point $x\in\ZR{\mathscr{X}}$ we use the notation as in \ref{ntn-ZRpoints}.
We set:
\begin{itemize}
\item $d(x)$: the dimension function (\ref{dfn-dimensionfunction});
\item $t(x)=\mathrm{tr.deg}_kk_x$ ($k_x$ is the residue field of $\O^{\int}_{\mathscr{X},x}$, or equivalently, of $V_x$);
\item $\Gamma_x=K^{\times}_x/V^{\times}_x$, the value group of the valuation ring $V_x$;
\item $\mathrm{rat\textrm{-}rank}(V_x|V)=\dim_{\Q}(\Gamma_x/\Gamma)\otimes\Q$.
\end{itemize}
The valuation ring $V_x$ at $x$ is $a$-adically separated (\ref{cor-ZRstrsheaf212}) and has the associated height one prime $\mathfrak{p}_x=\sqrt{aV_x}$\index{associated height one prime} ({\bf \ref{ch-pre}}.\ref{dfn-maxspe2}).
We set 
\begin{itemize}
\item $\ovl{V}_x=V_x/\mathfrak{p}_x$,
\end{itemize}
which is a valuation ring such that $\mathrm{ht}(V_x)=\mathrm{ht}(\ovl{V}_x)+1$.
One sees easily that the inequality
$$
\mathrm{rat\textrm{-}rank}(\ovl{V}_x)\leq\mathrm{rat\textrm{-}rank}(V_x|V)
$$
holds.

\begin{thm}\label{thm-heighestimaterigidspace}
{\rm (1)} We have the inequality
$$
\mathrm{rat\textrm{-}rank}(V_x|V)+t(x)\leq d(x).
$$
In particular, $\mathrm{rat\textrm{-}rank}(\ovl{V}_x)$ and $\mathrm{ht}(V_x)=\mathrm{ht}(\ovl{V}_x)+1$ are finite, and we have
$$
\mathrm{ht}(V_x)+t(x)\leq d(x)+1.
$$

{\rm (2)} If $\mathrm{rat\textrm{-}rank}(V_x|V)+t(x)=d(x)$, then $\Gamma_x/\Gamma$ is a finitely generated $\Z$-module, and $k_x$ is a finitely generated extension of $k$.

{\rm (3)} If $\mathrm{rat\textrm{-}rank}(\ovl{V}_x)+t(x)=d(x)$, then the value group of $\ovl{V}_x$ is isomorphic as a group to $\Z^d$ for some $d\geq 0$.

{\rm (4)} If $\mathrm{ht}(V_x)+t(x)=d(x)+1$, then the value group of $\ovl{V}_x$ is isomorphic as an ordered group to $\Z^d$ $($with $d=\mathrm{ht}(\ovl{V}_x))$ equipped with the lexicographical order\index{ordering@order(ing)!lexicographical ordering@lexicographical ---} $($cf.\ {\rm {\bf \ref{ch-pre}}.\ref{exa-height})}.
\end{thm}

\begin{proof}
First let us show (1).
We take an affinoid open neighborhood $\mathscr{U}=(\Spf A)^{\rig}$ of $x$, where $A$ is an $a$-torsion free topologically finitely generated $V$-algebra, such that $d(x)=\dim(\mathscr{U})$.
Since we may work locally around $x$, we may set $\mathscr{X}=\mathscr{U}$ without loss of generality.
By Noether normalization {\bf \ref{ch-pre}}.\ref{thm-noethernormalizationtype(V)} we may assume that $A=V\dl X_1,\ldots,X_n\dr$.
We perform an inductive argument with respect to $n$ as follows.
Let $R=V\dl X_1,\ldots,X_{n-1}\dr$ so that $A=R\dl X\dr$ ($X=X_n$).
Set $\mathscr{Y}=(\Spf R)^{\rig}$, and let $y$ be the image of the point $x$ by the map $\ZR{\mathscr{X}}\rightarrow\ZR{\mathscr{Y}}$ induced from the canonical inclusion $R\hookrightarrow A=R\dl X\dr$; notice that the induced morphism $V_y\rightarrow V_x$ is $a$-adic and hence is injective.
By induction it suffices to show the following inequality:
$$
\dim_{\Q}(\Gamma_x/\Gamma_y)\otimes\Q+\mathrm{tr.deg}_{k_y}k_x\leq 1.\leqno{(\ast)}
$$

To show this, consider the valuation $v_x$ on $A$ associated to the valuation ring $V_x$, and its restriction $v'$ on $A'=R[X]$. 
Let $V'$ be the valuation ring of the valuation $v'$, and $\Gamma_{V'}$ and $k_{V'}$ the value group and the residue field of the valuation $v'$, respectively.

\medskip
{\sc Claim.} {\it We have $\Gamma_x=\Gamma_{V'}$ and $k_x=k_{V'}$.}

\medskip
Indeed, since $V'$ is $a$-adically dense in $V_x$ and $a\in\m_{V_x}$, we have $k_x=k_{V'}$.
Moreover, since $V_x$ is $a$-adically separated, for any $f\in V_x$ there exists $N\geq 1$ such that $v_x(f)<v_x(a^N)$.
If $f=g+a^Nh$ with $g\in V'$, then $v_x(f)=v_x(g)=v'(g)$.
Hence we have the other equality.

Now we apply {\bf \ref{ch-pre}}.\ref{prop-extensionvaluationestimate} to get the desired inequality $(\ast)$.
If we have the equality in $(\ast)$, then $k_x$ is a finitely generated field extension of $k_y$, and $\Gamma_x/\Gamma_y$ is finitely generated.
Hence by induction we deduce (2).

Let us show (3).
We may replace $\mathscr{X}$ by an affinoid neighborhood $\mathscr{X}=(\Spf A)^{\rig}$ by an $a$-torsion free $A$.
Let $K$ be the kernel of the morphism $A\rightarrow\ovl{V}_x$, and set $A_0=A/K$.
Since $\sqrt{aA}\in K$, $A_0$ is a finite type domain over $k$ with the valuation $A_0\hookrightarrow\ovl{V}_x$.
Let $R$ be the localization of $A_0$ at the prime ideal $\m_{\ovl{V}_x}\cap A_0$.
Then, as can be shown easily with the aid of Noether normalization {\bf \ref{ch-pre}}.\ref{thm-noethernormalizationtype(V)}, one has $\dim(R)\leq d(x)$.
Hence we have the desired result by {\bf \ref{ch-pre}}.\ref{thm-estimate} (2).
Similarly, (4) follows from {\bf \ref{ch-pre}}.\ref{thm-estimate} (3).
\end{proof}

\subsubsection{Divisorial points}\label{subsub-classifyingvaluationsdivisorial}
\index{point!divisorial point@divisorial ---|(}
Let $V$ be an $a$-adically complete valuation ring of height one, and $\mathscr{X}$ a coherent rigid space of finite type over $\mathscr{S}=(\Spf V)^{\rig}$.
We use the notation as in the previous paragraph.
\begin{dfn}\label{dfn-rigiddivisorialvaluations}{\rm
(1) A point $x\in\ZR{\mathscr{X}}$ is said to be {\em divisorial $($over $V)$} if $t(x)=d(x)$.

(2) Let $X$ be a distinguished formal model of $\mathscr{X}$, and $x\in\ZR{\mathscr{X}}$ a divisorial point.
We say that $x$ is {\em residually algebraic over $X$} if $k_x$ is an algebraic extension of the residule field $k_{\sp_X(x)}$ at $\sp_X(x)\in X$.}
\end{dfn}

By \ref{thm-heighestimaterigidspace} we have:
\begin{prop}\label{prop-rigiddivisorialvaluations}
If $x\in\ZR{\mathscr{X}}$ is divisorial, then $x$ is of height one {\rm (\ref{dfn-ZRpoints6} (1))} and $\Gamma_x/\Gamma$ is a finite group. 
If, moreover, $x$ is residually algebraic over a distinguished formal model $X$ of $\mathscr{X}$, then $k_x$ is a finite extension of $k_{\sp_X(x)}$. \hfill$\square$
\end{prop}
\index{point!divisorial point@divisorial ---|)}

\subsubsection{Example: Unit disk}\label{subsub-classifyingvaluationsdisk}
\index{unit disk|(}
Let $V$ be an $a$-adically complete valuation ring of height one, and $K=\Frac(V)$.
We consider the unit disk $\D^1_K=\D^1_{\mathscr{S}}=(\Spf V\dl T\dr)^{\rig}$ over $\mathscr{S}=(\Spf V)^{\rig}$ (cf.\ \S\ref{subsub-unitdisk}).

\begin{prop}\label{prop-classifyingvaluationsdisk}
{\rm (1)} For $x\in\ZR{\D^1_K}$ one and only one of the following cases occurs$:$
\begin{itemize}
\item[{\rm (a)}] $($divisorial case$)$ $x$ is divisorial\index{point!divisorial point@divisorial ---}$;$ in this case, $x$ is of height one, $t(x)=d(x)=1$, and $\Gamma_x/\Gamma$ is finite$;$
\item[{\rm (b)}] $($subject-to-divisorial case$)$ $x$ is of height two$;$ in this case, $t(x)=0$, and $\Gamma_x/\Gamma$ is finitely generated with $\mathrm{rat\textrm{-}rank}(V_x|V)=1;$ moreover, the maximal generization $\til{x}$ of $x$ {\rm (\ref{dfn-ZRpoints6} (2))} is divisorial$;$
\item[{\rm (c)}] $($irrational case$)$ $x$ is of height one, $t(x)=0$ and $\mathrm{rat\textrm{-}rank}(V_x|V)=1;$ in this case, $\Gamma_x/\Gamma$ is finitely generated$;$
\item[{\rm (d)}] $($limit case$)$ $x$ is of height one, $t(x)=0$, and $\mathrm{rat\textrm{-}rank}(V_x|V)=0$.
\end{itemize}

{\rm (2)} $($classical case$)$ If the kernel of $V\dl T\dr\rightarrow V_x$ is non-trivial, then $x$ falls in the case {\rm (d)}$;$ in this case $x$ is a classical point {\rm (\ref{dfn-classicalpoint})}, and $\Gamma_x/\Gamma$ is finite.
\end{prop}

\begin{rem}\label{rem-classifyingvaluationsdisk2}{\rm 
Notice the similarity between the classification as in (1) and the classification of valuations in {\bf \ref{ch-pre}}, \S\ref{subsub-casekrull2}.}
\end{rem}

\begin{proof}[Proof of Proposition {\rm \ref{prop-classifyingvaluationsdisk}}]
First suppose that the kernel $J$ of $A=V\dl T\dr\rightarrow V_x$ is non-trivial.
Since $J$ is clearly $a$-saturated, $J$ is a finitely generated prime ideal of $A$.
Since $\dim(\D^1_K)=\dim(\Spec K\dl T\dr)=1$ (\ref{prop-dimensioncal}), $J$ defines a closed point on the associated scheme $s(\D^1_K)=\Spec K\dl T\dr$.
It follows as in the proof of \ref{prop-classicalpointsexist} that $W=A/J$ is finite over $V$ and that $(\Spf W)^{\rig}\hookrightarrow\D^1_K$ defines a classical point, which is nothing but $x$.
By \ref{thm-heighestimaterigidspace} we have
$$
\mathrm{rat\textrm{-}rank}(V_x|V)+t(x)=\dim((\Spf W)^{\rig})=0,
$$
and hence $x$ is of height one, $t(x)=0$, and $\Gamma_x/\Gamma$ is finite, which proves (2).

Let us show (1).
Notice first that $\mathrm{ht}(V_x)\leq 2$ by \ref{thm-heighestimaterigidspace} (1) and that we always have $d(x)=1$.
If $\mathrm{ht}(V_x)=2$, then we have the equalities in \ref{thm-heighestimaterigidspace} with $t(x)=0$.
Since $\mathrm{rat\textrm{-}rank}(V_x|V)\geq 1$, we have $\mathrm{rat\textrm{-}rank}(V_x|V)=1$, and hence $\Gamma_x/\Gamma$ is finitely generated by \ref{thm-heighestimaterigidspace} (2).

Next, we assume $\mathrm{ht}(V_x)=1$.
Then $t(x)=0$ or $1$.
If $t(x)=1$, then $d(x)=t(x)$, and hence $x$ is divisorial.
We have the equality as in \ref{thm-heighestimaterigidspace} (2) with $\mathrm{rat\textrm{-}rank}(V_x|V)=0$, and hence $\Gamma_x/\Gamma$ is finite (case (a)).
Suppose $t(x)=0$.
Then $\mathrm{rat\textrm{-}rank}(V_x|V)=0$ or $1$.
If $\mathrm{rat\textrm{-}rank}(V_x|V)=1$, then $t(x)=0$, and hence $\Gamma_x/\Gamma$ is finitely generated (\ref{thm-heighestimaterigidspace} (2)), which is the case (c).
If $\mathrm{rat\textrm{-}rank}(V_x|V)=0$, we are in the case (d).

It remains to show that, if $\mathrm{ht}(V_x)=2$, then the maximal generization $\til{x}$ of $x$ is divisorial.
We have $\mathrm{ht}(V_{\til{x}})=1$.
Since $V_x$ is a composite of $V_{\til{x}}$ and a height one valuation ring for $k_{\til{x}}$ over $k$, $k_{\til{x}}$ is not algebraic over $k$.
Hence we have $t(\til{x})=1$, which means that $\til{x}$ is divisorial.
\end{proof}

\begin{rem}\label{rem-classifyingvaluationsdisk}{\rm 
The above classification is compared with the Berkovich's classification of points on a unit disk (\cite[1.4.4]{Berk1}) as follows.
Berkovich\index{Berkovich, V.G.} deals only with height one points, and hence the case (b) does not occur.
\begin{itemize}
\item Berkovich's type (1) and type (4) points fall in the limit case (d); moreover, type (1) points are precisely the points in the classical case (2);
\item type (2) points correspond to the divisorial case (a);
\item type (3) points correspond to the irrational case (c).
\end{itemize}}
\end{rem}

The comparison can be established as follows.
We assume for simplicity that $K$ is algebraically closed; notice that, in this situation, if $\Gamma_x/\Gamma$ is finite, then we have $\Gamma_x=\Gamma$.
We fix $a\in\m_V\setminus\{0\}$ and a real number $0<c<1$ and consider the valuation $|\cdot|\colon K\rightarrow\R_{\geq 0}$ associated to the valuation ring $V$ such that $|a|=c$ (cf.\ {\bf \ref{ch-pre}}.\ref{subsub-nonarchnorms}).

\medskip
$\bullet$ For $u,b\in V$ with $b\neq 0$ one has the `closed subdisk' $\D(u,|b|)=$`$\{\|T-u\|\leq|b|\}$' defined as follows: consider the admissible blow-up $X_{u,b}\rightarrow X=\Spf V\dl T\dr$ along the admissible ideal $J=(b,T-u)\in V\dl T\dr$; it has the affine part $U'=\Spf V\dl\frac{T-u}{b}\dr$, and we set $\D(u,|b|)=U^{\rig}$, which has the open immersion $\D(u,|b|)\hookrightarrow\D^1_K$; as an open subspace of $\ZR{\D^1_K}$, it coincides with $\sp^{-1}_{X'}(U)$.

\medskip
$\bullet$ It is easy to see the following: a (not necessarily height one) point $x\in\ZR{\D^1_K}$ lies in $\ovl{\ZR{\D(u,|b|)}}$, the closure of $\ZR{\D(u,|b|)}=\sp^{-1}_{X'}(U)$ in $\ZR{\D^1_K}$, if and only if 
$$
\|T-u\|_x\leq|b|,
$$
where $\|\cdot\|_x=\|\cdot\|_{x,\mathscr{I},c}$ (with $\mathscr{I}=a\O^{\int}_{\D^1_K}$) is the seminorm\index{seminorm} associated to the point $x$ (\S\ref{subsub-seminormsofpoints}).
This fact leads one to consider the overconvergent closed subsets of $\ZR{\D^1_K}$ of the form 
$$
D(u,r)=\{y\in\ZR{\D^1_K}\,|\,\|T-u\|_y\leq r\}
$$
for $u\in V$ and $r\in\R_{\geq 0}$.

\medskip
$\bullet$ For any height one point $x\in\ZR{\D^1_K}$, consider the family of subsets $\{E(f)\}$ indexed by $f\in A=V\dl T\dr$, where
$$
E(f)=\{y\in\ZR{\D^1_K}\,|\,\|f\|_y\leq\|f\|_x\}.
$$
By Weierstrass preparation theorem (Exercise \ref{exer-preparationthm}) this family is determined only by $E(f)$'s for polynomials $f\in V[T]$.
Moreover, if $f=g_1\cdots g_r$ in $V[T]$, then $E(f)$ contains the intersection $E(g_1)\cap\cdots\cap E(g_r)$.
Hence the family $\{E(f)\}$ is filtered, and the subfamily $\{E(T-u)\}_{u\in V}$ consisting only of those determined by linear polynomials is cofinal.
Since $\|T-u\|_x\leq\|T-v\|_x$ implies $E(T-u)\subseteq E(T-v)$, $\{E(T-u)\}_{u\in V}$ is a `family of embedded disks' (according to Berkovich's terminology \cite[1.4.4]{Berk1}).

Now recall the Berkovich's classification of points: {\it height one points of $\D^1_K$ are in one to one correspondence with equivalence classes of families of embedded disks} (see \cite[1.4.4]{Berk1} for more details).
Hence in our situation the collection $\{E(T-u)\}_{u\in V}$ defines a point in Berkovich's sense. 

To understand the types of the Berkovich's points thus obtained, we need to consider the set of radii $\{\|T-u\|_x\}_{u\in V}$; we set
$$
r_x=\inf\{\|T-u\|_x\,|\,u\in V\}.
$$
Notice that this value $r_x$ lies in $\|K\|_x$; indeed, if $\{T-u_n\}_{n\geq 0}$ is such that $\lim\|T-u_n\|_x=r_x$, then, since $|u_n-u_m|=\|u_n-u_m\|_x=\|(T-u_m)-(T-u_n)\|_x$ for $n,m\geq 0$, $\{u_n\}_{n\geq 0}$ is a Cauchy sequence in $V$, and we have $r_x=\|T-u\|_x$ for $u=\lim u_n$.
\begin{itemize}
\item[{\rm (i)}] Suppose that $r_x\not\in|K|$. Then since the value group of $V_x$ is strictly larger than $\Gamma$, the point $x$, which is of type (3) in Berkovich's classification, falls into the irrational case (c)
Notice that, in this case, for any $u\in V$ and for any $b\in V$ such that $\|T-u\|_x>|b|$, the rigid point $\Spf\widehat{V}_x\rightarrow X_{u,b}$ hits the closed fiber $X_{u,b}\otimes_Vk$ at the double point (the intersection of the two irreducible components).
\item[{\rm (ii)}] If $r_x\in|K^{\times}|$, then $x$ is of Berkovich's type (2).
For $u\in V$ and $b\in V\setminus\{0\}$ such that $r_x=|b|=\|T-u\|_x$, the rigid point $\Spf\widehat{V}_x\rightarrow X_{u,b}$ hits the generic point of $U'=\Spf V\dl\frac{T-u}{b}\dr$.
In particular, $x$ is divisorial\index{point!divisorial point@divisorial ---} (that is, $t(x)=1$).
\item[{\rm (iii)}] Suppose $r_x=0$. In this case, for any $b\in V\setminus\{0\}$ there exists $u\in V$ such that $\|T-u\|_x\leq|b|$; then the rigid point $\Spf\widehat{V}_x\rightarrow X_{u,b}$ hits the closed fiber $X_{u,b}\otimes_Vk$ at a closed point on the affine part $U'\otimes_Vk$.
The points in this case, of type (1) or (4) in Berkovich's classification, fall into (d) (limit case).
\end{itemize}
\index{unit disk|)}

\subsection{Maximum modulus principle}\label{sub-MMP}
\subsubsection{Spectral seminorm formula}\label{subsub-spectralseminormformula}
Let $V$ be an $a$-adically complete valuation ring of height one, and $\mathscr{X}$ a coherent rigid space of finite type over $\mathscr{S}=(\Spf V)^{\rig}$.
We fix once for all a real number $0<c<1$, which gives rise to a non-archimedean norm 
$$
|\cdot|\colon K={\textstyle V[\frac{1}{a}]}\longrightarrow\R_{\geq 0}
$$
with $|a|=c$, and the spectral seminorm\index{seminorm!spectral seminorm@spectral ---}\index{seminorm} $\|\cdot\|_{\Sp}=\|\cdot\|_{\Sp,\mathscr{I},c}$ (\S\ref{subsub-spectralseminorms}), where $\mathscr{I}=a\O^{\int}_{\mathscr{X}}$.

Let $X$ be a flat and of finite type formal scheme giving a distinguished formal model of $\mathscr{X}$, $X_k=X\otimes_Vk$ its closed fiber, where $k$ is the residue field $k=V/\m_V$ of $V$, and $(X_k)_{\red}$ the reduced model of $X_k$.
Notice that $X_k$ and $(X_k)_{\red}$ are finite type schemes over $k$.
We define
\begin{itemize}
\item $I(X)=$ the set of all points of $(X_k)_{\red}$ of which the closure is an irreducible component of $(X_k)_{\red}$;
\item $D(X)=$ the set of all divisorial points\index{point!divisorial point@divisorial ---} in $\ZR{\mathscr{X}}$ that are residually algebraic over $X$ (\ref{dfn-rigiddivisorialvaluations}).
\end{itemize}

Notice that, due to \ref{prop-rigiddivisorialvaluations}, one can equivalently define $D(X)$ as the set of all divisorial points of $\ZR{\mathscr{X}}$ that are {\it residually finite} over $X$.

\begin{thm}\label{thm-spectralseminormformula}
{\rm (1)} For any $z\in I(X)$, the set $\sp^{-1}_X(z)$, where $\sp_X\colon\ZR{\mathscr{X}}\rightarrow X$ is the specialization map, is non-empty, and is contained in $D(X)$.

{\rm (2)} We have 
$$
D(X)=\bigcup_{z\in I(X)}\sp^{-1}_X(z).
$$

{\rm (3)} {\rm (Spectral seminorm formula)}
The set $D(X)$ is finite. Moreover, for any $f\in\Gamma(\mathscr{X},\O_{\mathscr{X}})$, we have
$$
\|f\|_{\Sp}=\max_{x\in D(X)}\|f\|_x,
$$
where $\|\cdot\|_x=\|\cdot\|_{x,\mathscr{I},c}$ is the seminorm\index{seminorm} at the point $x$ {\rm (\S\ref{subsub-seminormsofpoints})}.
\end{thm}

To show the theorem, we need the following lemma.
\begin{lem}\label{lem-spectralseminormformula}
Let $A$ be an $a$-torsion free topologically of finite type integral domain over $V$, and set $\mathscr{X}=(\Spf A)^{\rig}$.
Then, for any proper closed subspace $\mathscr{Z}\subsetneq\mathscr{X}$, we have $\dim(\mathscr{Z})<\dim(\mathscr{X})$.
\end{lem}

\begin{proof}
Let $d=\dim(\mathscr{X})$, and take a finite injection $V\dl X_1,\ldots,X_d\dr\hookrightarrow A$ with $V$-flat cokernel (see {\bf \ref{ch-pre}}.\ref{thm-noethernormalizationtype(V)}).
Any proper closed subspace $\mathscr{Z}\subsetneq\mathscr{X}$ corresponds to a proper closed subscheme $Z\subseteq\Spec A[\frac{1}{a}]$ (\ref{prop-closedsubspaceaffinoid2}) such that $\dim(\mathscr{Z})=\dim(Z)$ (\ref{cor-dimaffinoids}), and hence to a $V$-flat proper closed subscheme $\ovl{Z}\subseteq\Spec A$.
Since a proper closed subscheme of $\Spec A$ does not contain the generic point of $\Spec A$, its image under the finite and dominant map $X\rightarrow\A=\Spf V\dl X_1,\ldots,X_d\dr$ does not contain the generic point, and hence is a proper closed subspace of $\A$.
Hence, to show the lemma, one can reduce to the case $X=\A$, and the assertion in this case is clear.
\end{proof}

\begin{proof}[Proof of Theorem {\rm \ref{thm-spectralseminormformula}}]
First notice that, for any open subset $U\subseteq X$, we have $I(U)=I(X)\cap (U_k)_{\red}$ and $D(U)=D(X)\cap\ZR{U^{\rig}}$.
In particular, if $X=\bigcup_{\alpha\in L}U_{\alpha}$ is a finite open covering, we have $I(X)=\bigcup_{\alpha\in L}I(U_{\alpha})$ and $D(X)=\bigcup_{\alpha\in L}D(U_{\alpha})$.
Hence the assertions (1), (2), and the finiteness (3) are reduced to those for each $U_{\alpha}$.
Moreover, since $[\mathscr{X}]=\bigcup_{\alpha\in L}[U^{\rig}_{\alpha}]$, the spectral seminorm $\|f\|_{\Sp}$ on $\mathscr{X}$ is equal to the maximum of the spectral seminorms of $f$ on $U_{\alpha}$'s.
Hence the last assertion of (3) can also be reduced to that on each $U_{\alpha}$.
We can therefore assume that $X$ is affine, say $X=\Spf A$, where $A$ is an $a$-torsion free topologically finitely generated $V$-algebra.

Let us first show the theorem under the additional assumption that $A$ is an integral domain.
If $d$ is the Krull dimension of $\mathcal{A}=A[\frac{1}{a}]$, then we have $\dim(\mathscr{X})=d$ by \ref{cor-dimaffinoids}.
Take a finite injection $V\dl X_1,\ldots,X_d\dr\hookrightarrow A$ with $V$-flat cokernel (see {\bf \ref{ch-pre}}.\ref{thm-noethernormalizationtype(V)}), and set $\A=\Spf V\dl X_1,\ldots,X_d\dr$.
Notice that, passing to the closed fibers, we still have finite injection $k[X_1,\ldots,X_d]\hookrightarrow A_k=A\otimes_Vk$.
Notice that, by \ref{prop-dimensionfunction1}, the dimension function $d(x)$ is constant, that is, $d(x)=d$ for all $x\in\ZR{\mathscr{X}}$.

Let us show (1). 
Since the specialization map $\sp_X\colon\ZR{\mathscr{X}}\rightarrow X$ is surjective (see \ref{prop-ZRpoints4}), $\sp^{-1}_X(z)$ for $z\in I(X)$ is non-empty.
Let us show the inclusion $\sp^{-1}_X(s)\subseteq D(X)$.
For any $z\in I(X)$, let $Z_z$ be the closure $\ovl{\{z\}}$ in $(X_k)_{\red}$, which is an irreducible component of $(X_k)_{\red}$.
Set 
$$
C_z=\bigcup_{z'\in I(X)\setminus\{z\}}Z_{z'},
$$
and let $U$ be the open complement of $C_z$ in $X$, which is an open formal subscheme of $X$.
Then $(U_k)_{\red}=Z_z\setminus C_z$ is irreducible and non-empty.
Since $\dim(U)^{\rig}=d$ holds by \ref{cor-dimensioncal101}, we have $\dim((U_k)_{\red})=d$ by \ref{cor-dimensionclosedfiber}.
Since $U$ satisfies $(U_k)_{\red}\cap I(X)=\{z\}$, we have $t(x)\geq d$ for any $x\in\sp^{-1}_X(z)$.
Such an $x$ is, since $d(x)=d$, always divisorial by the inequality in \ref{thm-heighestimaterigidspace}, hence belonging to $D(X)$.

Next we show (2). 
For $x\in D(X)$, we have, by constancy of the dimension function (\ref{prop-dimensionfunction1}), $d(x)=d$, hence $t(x)=d$.
This means that $x$ dominates the generic point of $\A_k$, the closed fiber of $\A$, hence belonging to $\sp^{-1}_X(z)$ for some $z\in I(X)$.
Thus we have $D(X)\subseteq\bigcup_{z\in I(X)}\sp^{-1}_X(z)$.
The other inclusion follows from (1).

Let us show (3).
We show that the set $D(X)$ is finite.
First, notice that $D(\A)$ consists of a single point, which we denote by $w$.
Since $\mathscr{X}\rightarrow\A^{\rig}$ is finite, we readily see that $D(X)$ is the preimage of $w$ by the map $\ZR{\mathscr{X}}\rightarrow\ZR{\A^{\rig}}$, which is finite due to \ref{prop-finitemorrigid99}

Next, let us show the spectral seminorm formula.
Since by definition $\|f\|_{\Sp}\geq\max_{x\in D(X)}\|f\|_x$, we need to show $\|f\|_{\Sp}\leq\max_{x\in D(X)}\|f\|_x$.
Suppose that for any $x\in D(X)$ we have $\|f^n\|_x\leq|a^m|$ for $n,m\geq 1$.
Since by \ref{lem-spectralfunctorseparatedquotientsVR1}
$$
A^{\int}=\{f\in\Gamma(\mathscr{X},\O_{\mathscr{X}}):\ \textrm{$f$ is integral at each $x\in D(X)$}\},
$$
this implies $f^n/a^m\in A^{\int}$.
Therefore, $\|f^n\|_{\Sp}\leq|a|^m$.
We have shown that, for any $n,m\geq 1$, $\max_{x\in D(X)}\|f\|_x\leq|a|^{\frac{m}{n}}$ implies $\|f\|_{\Sp}\leq|a|^{\frac{m}{n}}$.
By taking the limit $|a|^{\frac{m}{n}}\rightarrow\max_{x\in D(X)}\|f\|_x$, the desired inequality follows.

Finally, let us discuss the general case $X=\Spf A$, where $A$ is not necessarily an integral domain.
We may assume that $A$ is reduced.
Let 
$$
\Spec A=\bigcup_{i\in I}\til{Z}_i,\quad\til{Z}_i=\Spec A_i\ (i\in I),
$$
where $I$ is a finite set, be the irreducible decomposition; each $\til{Z}_i$ is the closure of an irreducible component of the Noetherian scheme $\Spec A[\frac{1}{a}]$.
For $i\in I$, set 
$$
\til{S}_i=\bigcup_{j\neq i}\big(\til{Z}_j\cap\til{Z}_i\big),
$$
and set $Z_i=\widehat{\til{Z}_i}$ and $S_i=\widehat{\til{S}_i}$.
Since $\til{Z}_i$ is irreducible and $\til{S}_i$ is a $V$-flat proper closed subscheme of $\til{Z}_i$, we have, by \ref{lem-spectralseminormformula}, $\dim(S^{\rig}_i)<\dim(Z^{\rig}_i)$.
In particular, we have $\dim((S_{i,k})_{\red})<\dim(Z^{\rig}_i)$ (\ref{cor-dimensionclosedfiber}).
By the above-discussed case of integral domains, any irreducible component of $(Z_{i,k})_{\red}$ is of dimension $d_i=\dim(Z^{\rig}_i)$, which implies that the complement $U_i$ of $(S_{i,k})_{\red}$ in $Z_i$ is dense in $Z_i$, and that $D(Z_i)=D(U_i)$.

We claim that the following equalities hold:
$$
D(X)=\bigcup_{i\in I}D(Z_i),\quad I(X)=\bigcup_{i\in I}I(Z_i).\eqno{(\ast)}
$$
Indeed, by \ref{prop-dimensionfunction1}, we have $D(X)\subseteq\bigcup_{i\in I}D(Z_i)$, and the other inclusion follows by
$$
\bigcup_{i\in I}D(Z_i)=\bigcup_{i\in I}D(U_i)=D\big({\textstyle \bigcup_{i\in I}U_i}\big)\subseteq D(X).
$$
SInce $U_i$ is dense in $Z_i$, we have $I(Z_{i,\red})=I(U_{i,\red})$ for $i\in I$.
Hence
$$
\bigcup_{i\in I}I(Z_i)=\bigcup_{i\in I}I(U_i)=I\big({\textstyle \bigcup_{i\in I}U_i}\big)=I(X),
$$
where the last equality follows from the fact that $\bigcup_{i\in I}U_i$ is dense in $X$.

By the equalities $(\ast)$ and the above argument on the case where $A$ is an integral domain, one deduces the assertions (1), (2), and the finitenes in (3).
To show the spectral seminorm formula, notice first that we have an irreducible decomposition
$$
\ZR{\mathscr{X}}=\bigcup_{i\in I}\ZR{Z^{\rig}_i}\eqno{(\ast\ast)}
$$
of the rigid space $\mathscr{X}=(\Spf A)^{\rig}$ (see \S\ref{subsub-closedimmrigidred}).
Then the decompositions $(\ast)$ and $(\ast\ast)$ and the spectral seminorm formula in the integral domain case imply the spectral seminorm formula in general.
\end{proof}

\begin{rem}\label{rem-spectralseminormformula}{\rm 
The theorem shows that there exists a canonical surjective map 
$$
D(X)\longrightarrow I(X),
$$
that is, any point in $D(X)$ dominates an irreducible component of $(X_k)_{\red}$, hence a point in $I(X)$.}
\end{rem}

\subsubsection{Maximum modulus principle}\label{subsub-MMP}
\index{Maximum modulus principle|(}
We continue to use the notation fixed in the beginning of \S\ref{subsub-spectralseminormformula}.
\begin{thm}[Maximum modulus principle]\label{thm-MMP}
Let $\mathscr{X}$ be a coherent rigid space over $\mathscr{S}=(\Spf V)^{\rig}$.

{\rm (1)} There exists a positive integer $e\geq 1$ such that for any $f\in\Gamma(\mathscr{X},\O_{\mathscr{X}})$ the value $\|f\|_{\Sp}$ belongs in $|K^{\times}|^{\frac{1}{e}}\cup\{0\}$.

{\rm (2)} For any $f\in\Gamma(\mathscr{X},\O_{\mathscr{X}})$ there exists a non-empty quasi-compact open subspace $\mathscr{U}\subset\mathscr{X}$ such that
$$
\ovl{\ZR{\mathscr{U}}}=\{x\in\ZR{\mathscr{X}}\,|\,\|f\|_x=\|f\|_{\Sp}\}.
$$

{\rm (3)} For any $f\in\Gamma(\mathscr{X},\O_{\mathscr{X}})$ the set $\ZR{\mathscr{X}}^{\cl}\cap\{x\in\ZR{\mathscr{X}}\,|\,\|f\|_x=\|f\|_{\Sp}\}$ is non-empty, that is, the maximum value of $\|f\|_x$ is attained at a classical point\index{point!classical point@classical ---}\index{classical point}.
\end{thm}

\begin{proof}
Since the group $\Gamma_x/\Gamma$ is finite at a divisorial point, (1) follows immediately from \ref{thm-spectralseminormformula}.
By \ref{cor-classicalpointsexist} the assertion (3) will follow from (2).
Thus it suffices to show (2).

By the reduction process similar to that in the beginning of the proof of \ref{thm-spectralseminormformula}, we may assume that $\mathscr{X}=X^{\rig}$ with $X=\Spf A$. 
Let $f\in\Gamma(\mathscr{X},\O_{\mathscr{X}})=A[\frac{1}{a}]$.
Since the case $\|f\|_{\Sp}=0$ is trivial, we assume $\|f\|_{\Sp}\neq 0$.
By (1) there exist $b\in K^{\times}$ and $e\geq 1$ such that $\|f^e\|_{\Sp}=|b|$.
We may replace $f$ by $f^e/b$ without changing the right-hand set of the desired equality and thus may assume that $\|f\|_{\Sp}=1$.
Then by \ref{lem-spectralfunctorseparatedquotientsVR1} we deduce that $f\in A^{\int}$.
By \ref{lem-corlemtuavsrigidaff} there exists an admissible blow-up $X'\rightarrow X$ such that $f\in\Gamma(X',\O_{X'})$.
We take the non-empty localization $X'_f$ by $f$, and consider the corresponding non-empty quasi-compact open subspace $\mathscr{U}=(X'_f)^{\rig}$ of $\mathscr{X}$.
It is then obvious that $\ovl{\ZR{\mathscr{U}}}=\{x\in\ZR{\mathscr{X}}\,|\,\|f\|_x=1\}$.
\end{proof}
\index{Maximum modulus principle|)}

\subsubsection{Reduction scheme}\label{subsub-reductionscheme}
\index{scheme!reduction scheme@reduction ---|(}
We use the same notation as in \S\ref{subsub-spectralseminormformula}.
The following theorem says that birational geometry (in our sense; see Introduction) is easier around residually finite divisorial points.
\begin{thm}\label{thm-reductionscheme}
Let $\mathscr{X}$ be a coherent rigid space of finite type over $\mathscr{S}=(\Spf V)^{\rig}$, $X$ a distinguished formal model of $\mathscr{X}$, and $x\in\ZR{\mathscr{X}}$ a divisorial point that is residually finite over $X$.
Then any admissible blow-up if $\pi\colon X'\rightarrow X$ is finite near $\sp_X(x)$.
\end{thm}

\begin{proof}
Let $U=\Spf A$ be an open neighborhood of $z=\sp_X(x)$ that does not intersect any other components $\ovl{\{z'\}}$ with $z'\in I(X)$ and $z'\neq z$.
Take a finite dominant morphism $U\rightarrow\A=\Spf V\dl X_1,\ldots,X_d\dr$ (by Noether normalization; see {\bf \ref{ch-pre}}.\ref{thm-noethernormalizationtype(V)}).
Note that $D(\A)$ consists of a single point, which we denote by $w$.
Set $W=\Spf \widehat{V_w}$, and consider the canonical morphism $W\rightarrow\A$.
Then $X'_W=X'\times_{\A}W$ is an admissible blow-up of $X_W=X\times_{\A}W$.
Since $X_W$ is finite flat over $W$, $X'_W$ is finite over $W$.
In particular, $X'_W\rightarrow X_W$ is finite.
Since $X'\rightarrow X$ is of finite type, there exists an open neighborhood $U'_k$ of $z$ in $U_k$ ($=$ the closed fiber of $U$) such that $\pi_k\colon X'_k\rightarrow X_k$ is finite over $U'_k$.
Hence $\pi$ is finite over the open subspace $U'$ of $U$ corresponding $U'_k$.
\end{proof}

\begin{dfn}\label{dfn-reductionscheme}{\rm 
Let $\mathscr{X}$ be a coherent rigid space of finite type over $\mathscr{S}=(\Spf V)^{\rig}$.

{\rm (1)} Two distinguished formal models $X_i$ ($i=1,2$) of $\mathscr{X}$ are said to be {\it strictly equivalent} if there exists a third distinguished formal model $X_3$ of $\mathscr{X}$ with finite morphisms $X_3\longrightarrow X_i$ ($i=1,2$).

{\rm (2)} For a distinguished formal model $X$ of $\mathscr{X}$, let $\mathscr{C}_X$ be the category of distinguished formal models $X'$ of $\mathscr{X}$ finite over $X$, and morphisms over $X$.
Set 
$$
R(\mathscr{X},X)=\varprojlim_{X'\in\mathscr{C}_X}(X'_k)_{\red}
$$
(where $(X'_k)_{\red}$ denotes the reduced model of the closed fiber $X'_k$ of $X'$), which is representable in the category of schemes.
We call $R(\mathscr{X},X)$ the {\it reduction scheme of $\mathscr{X}$ relative to $X$}.}
\end{dfn}

Notice that the strict equivalence is an equivalence relation, and that the category $\mathscr{C}_X$ is directed; indeed, for two finite $X'_1\rightarrow X$ and $X'_2\rightarrow X$ in $\mathscr{C}_X$, there exists an admissible blow-up $W\rightarrow X$ dominating both $X'_1$ and $X'_2$; then take the Stein factorization $W\rightarrow Z\rightarrow X$ (see {\bf \ref{ch-formal}}.\ref{prop-steinfactadqformalsch11}), and observe that $Z$ is finite over $X$, and that $Z$ dominates both $X'_i$ ($i=1,2$).

Notice also that the scheme $R(\mathscr{X},X)$ depends only on the strict equivalence class of $X$, and the construction is canonical; that is, for a finite type morphism $f\colon X\rightarrow Y$ between coherent formal schemes of finite type over $V$, we have the canonical morphism
$$
R(f)\colon R(\mathscr{X},X)\longrightarrow R(\mathscr{Y},Y),
$$
where $\mathscr{Y}=Y^{\rig}$.
Indeed, for any object $Y'\rightarrow Y$ in $\mathscr{C}_Y$, the objects in $\mathscr{C}_X$ dominating $X\times_YY'$ form a cofinal subcategory $\mathscr{C}'_X$, and hence one obtains
$$
R(\mathscr{X},X)=\varprojlim_{X'\in\mathscr{C}'_X}(X'_k)_{\red}\longrightarrow R(\mathscr{Y},Y)=\varprojlim_{Y'\in\mathscr{C}_Y}(Y'_k)_{\red}.
$$

\begin{rem}\label{rem-reductionscheme}{\rm 
If $\mathscr{X}=(\Spf A)^{\rig}$, then the reduction scheme $R(\mathscr{X},X)$ coincides with the spectrum of the $k$-algebra $\til{\mathcal{A}}$ (where $\mathcal{A}$ is the associated classical affinoid algebra $\mathcal{A}=A[\frac{1}{a}]$) discussed in \cite[\S6.3]{BGR}.
See \ref{cor-cordescpbelements1} in the appendix.}
\end{rem}

\begin{prop}\label{prop-reductionscheme1}
Let $\mathscr{X}$ be a coherent rigid space of finite type over $\mathscr{S}=(\Spf V)^{\rig}$, and $X$ a distinguished formal model of $\mathscr{X}$ over $V$.
Then there exists a distinguished formal model $X'$ with the following properties.
\begin{itemize}
\item[{\rm (a)}] The formal model $X'$ dominates $X$ by a finite morphism $X'\rightarrow X$.
\item[{\rm (b)}] The canonical map $D(X')\rightarrow I(X')$ induced by $\sp_{X'}$ is bijective.
\item[{\rm (c)}] For any $x\in D(X')$, the residue field $k_x$ is isomorphic to the residue field $k_{\sp_{X'}(x)}$ at $\sp_{X'}(x)\in X'$.
\end{itemize}
\end{prop}

\begin{proof}
Let us first try to find an admissible blow-up $Z\rightarrow X$ that satisfies (b) and (c).
Notice that, by the construction of the Zariski-Riemann space $\ZR{\mathscr{X}}$, for two distinct points $x,y\in\ZR{\mathscr{X}}$, there exists an admissible blow-up $Z\rightarrow X$ such that $\sp_Z(x)\neq\sp_Z(y)$.
Since $D(X)$ is a finite subset of $\ZR{\mathscr{X}}$, one can take an admissible blow-up $Z_0\rightarrow X$ such that $D(X)\rightarrow I(Z_0)$ is injective.
Notice that 
$$
\varinjlim_{Z}k_{\sp_Z(x)}=k_x
$$
holds, where $Z$ runs through distinguished formal model dominating $Z_0$.
Since, for any $x\in D(X)$, $k_x/k_{\sp_{Z_0}(x)}$ is a finite extension, we can take an admissible blow-up $Z$ dominating $Z_0$ such that $k_{\sp_Z(x)}=k_x$ for any $x\in D(X)$.
Thus we have the desired admissible blow-up $Z\rightarrow X$ with the properties (b) and (c).
Then, the formal model $X'$ in the Stein factorization $Z\rightarrow X'\rightarrow X$ (see {\bf \ref{ch-formal}}.\ref{prop-steinfactadqformalsch11}) satisfies all of the required properties, since $Z$ and $X'$ are isomorphic over an open subspace of $X$ that contains $I(X)=\sp_X(D(X))$; see \ref{thm-spectralseminormformula}.
\end{proof}

\begin{rem}\label{rem-reductionscheme3}{\rm 
The proposition and the proof of \ref{thm-spectralseminormformula} show that, if $X$ is affine, then the spectral seminorm formula (\ref{thm-spectralseminormformula} (3)) with $D(X)$ replaced by a proper subset of it, does not hold.}
\end{rem}

\begin{thm}\label{thm-reductionscheme2}
Let $\mathscr{X}$ be a coherent rigid space of finite type over $\mathscr{S}=(\Spf V)^{\rig}$, and $X$ a distinguished formal model of $\mathscr{X}$ over $V$.

{\rm (1)} $($Finiteness$)$ The reduction scheme $R(\mathscr{X},X)$ is a scheme of finite type over $k$.

{\rm (2)} $($Stability$)$ There exists a distinguished formal model $Z'$ finite over $X$ such that $(Z'_k)_{\red}=R(\mathscr{X},X)$.
Moreover, for any distinguished formal model $Z$ finite over $Z'$, we have $(Z_k)_{\red}=R(\mathscr{X},X)$.
$($We call such a formal model $Z'$ a {\rm stabilized model}.$)$
\end{thm}

\begin{proof}
Let $X'$ be a distinguished formal model with the properties (a), (b), and (c) in \ref{prop-reductionscheme1}.
For any object $Z$ of $\mathscr{C}_X$ dominating $X'$, the finite map $(Z_k)_{\red}\rightarrow (X'_k)_{\red}$ is birational, and thus $(Z_k)_{\red}$ is dominated by the normalization $W$ of $(X'_k)_{\red}$.
Hence $R(\mathscr{X},X)$ is dominated by $W$.
Since $W$ is finite over $(X'_k)_{\red}$, $R(\mathscr{X},X)$ is finite over $(X'_k)_{\red}$, which shows (1).
(2) is clear, since any increasing sequence of quasi-coherent algebras between $\O_{(X'_k)_{\red}}$ and $\O_W$ is stationary.
\end{proof}

\begin{thm}\label{thm-reductionscheme3}
For any morphism $f\colon X\rightarrow Y$ of coherent formal schemes of finite type over $V$, the following conditions are equivalent.
\begin{itemize}
\item[{\rm (a)}] $f$ is proper $($resp.\ separated, resp.\ affine, resp.\ finite$)$.
\item[{\rm (b)}] $R(f)$ is proper $($resp.\ separated, resp.\ affine, resp.\ finite$)$.
\end{itemize}
\end{thm}

\begin{proof}
Take a formal model $f'\colon X'\rightarrow Y'$ of $f^{\rig}\colon X^{\rig}\rightarrow Y^{\rig}$, where $X'$ and $Y'$ are stabilized model of $X^{\rig}$ and $Y^{\rig}$, respectively, so that $R(f)$ is given by $(f'_k)_{\red}\colon(X'_k)_{\red}\rightarrow(Y'_k)_{\red}$.
Then the assertion follows from the fact that a morphism $g\colon S\rightarrow T$ of schemes is proper (resp.\ separated, resp.\ affine, resp.\ finite) if and only if so is $g_{\red}\colon S_{\red}\rightarrow T_{\red}$.
For separatedness, see \cite{EGA}, {\bf I}, (5.5.1), and for properness, see \cite{EGA}, {\bf II}, (5.4.6).
As for affineness, one first uses the absolute Noetherian approximation (\cite[C.9]{TT}), \cite{EGA}, {\bf IV}, (8.10.5), and Serre's criterion \cite{EGA}, {\bf II}, (5.2.1).
Finally, the claim for finiteness follows from \ref{prop-appclassicalZR} in the appendix.
\end{proof}

\begin{rem}\label{rem-reductionscheme2}{\rm 
One can define the so-called {\it reduction map}
$$
\red(\mathscr{X},X)\colon\ZR{\mathscr{X}}\longrightarrow R(\mathscr{X},X)
$$
by $\sp_Z$, where $X$ is a stabilized model.
It is straightforward to see that the map $\red(\mathscr{X},X)$ is continuous; it is, moreover, surjective, even restricted on the subset $[\mathscr{X}]$ consisting of maximal points; see \ref{prop-spectralfunctorseparatedquotients}.}
\end{rem}
\index{scheme!reduction scheme@reduction ---|)}

\addcontentsline{toc}{subsection}{Exercises}
\subsection*{Exercises}
\begin{exer}\label{exer-spectralGauss}{\rm 
Let $V$ be an $a$-adically complete valuation ring of height one, $K=\Frac(V)=V[\frac{1}{a}]$ the fractional field, and $k=V/\m_V$ the residue field.
Set $A=V\dl T_1,\ldots,T_n\dr$ and $\mathcal{A}=A[\frac{1}{a}]=K\dl T_1,\ldots,T_n\dr$ (the Tate algebra), and consider the corresponding coherent rigid space $\D^1_K=(\Spf A)^{\rig}$ (the unit disk).
Let $|\cdot|\colon K\rightarrow\R_{\geq 0}$ be the non-archimedean norm such that $|a|=c$ for a fixed real number $0<c<1$, and consider the spectral seminorm $\|\cdot\|_{\Sp}=\|\cdot\|_{\Sp,\mathscr{I},c}$, where $\mathscr{I}=a\O^{\int}_{\D^1_K}$, and the Gauss norm\index{norm!Gauss norm@Gauss ---} $\|\cdot\|_{\mathrm{Gauss}}$ ({\bf \ref{ch-pre}}, \S\ref{subsub-classicalaffinoidalgebras}) on $\mathcal{A}=\Gamma(\D^1_K,\O_{\D^1_K})$ (cf.\ \ref{thm-comparisonaffinoid}).

(1) Show that 
$$
\|\cdot\|_{\Sp}=\|\cdot\|_{\mathrm{Gauss}}.
$$

(2) Let $x\in\ZR{\D^1_K}$ be such that the rigid point $\Spf\widehat{V}_x\rightarrow\Spf A$ hits the closed fiber $\Spec A\otimes_Vk=\Spec k[T_1,\ldots,T_n]$ at the generic point.
Then show that the seminorm $\|\cdot\|_x=\|\cdot\|_{x,\mathscr{I},c}$ coincides with the spectral seminorm $\|\cdot\|_{\Sp}$.}
\end{exer}

\begin{exer}\label{exer-spectralGauss2}{\rm 
Consider the situation as in Exercise \ref{exer-spectralGauss}.

{\rm (1)} For $f\in\mathcal{A}$, $\inf_{x\in\ZR{\mathscr{X}}}\|f\|_x$ is positive if and only if $f\in\mathcal{A}^{\times}$.

{\rm (2)} Consider the map $\|f(\cdot)\|\colon\ZR{\mathscr{X}}\rightarrow\R_{\geq 0}$ ($x\mapsto\|f\|_x$).
Then $f$ is unit in $\mathcal{A}$ if and only if the image of $\|f(\cdot)\|$ consists only of one non-zero value; otherwise, the image is the closed interval $[0,\|f\|_{\mathrm{Gauss}}]$.

{\rm (3)} Suppose $|\ovl{K}^{\times}|\neq\R_{>0}$, where $\ovl{K}$ is the algebraic closure of $K$.
Then for a non-zero non-unit $f\in\mathcal{A}$ and $r\in(0,\|f\|_{\mathrm{Gauss}}]\setminus|\ovl{K}^{\times}|$, the open subspace $\mathscr{U}$ such that $\ZR{\mathscr{U}}=\{x\in\ZR{\mathscr{X}}\,|\,\|f\|_x\neq r\}$ is a non-empty overconvergent open subspace that fails to enjoy the following property (the so-called {\em $\ZR{\cdot}^{\cl}$-admissibility}; cf.\ \S\ref{sub-admissibleopensubsets} below): for any coherent rigid space $\mathscr{V}$ of finite type over $\mathscr{S}=(\Spf V)^{\rig}$ and any $\mathscr{S}$-morphism $\varphi\colon\mathscr{V}\rightarrow\mathscr{X}$ that maps $\ZR{\varphi}(\ZR{\mathscr{V}}^{\cl})\subseteq\ZR{\mathscr{U}}^{\cl}$, we have $\varphi(\mathscr{V})\subseteq\mathscr{U}$.}
\end{exer}

\begin{exer}\label{exer-spectralGauss3}{\rm 
Let $V$ be an $a$-adically complete valuation ring of height one, and $K=V[\frac{1}{a}]$ the fractional field.
Let $\mathscr{X}$ be a rigid space of finite type over $\mathscr{S}=(\Spf V)^{\rig}$.
Consider an extension $V\subseteq V'$ of $a$-adically complete valuation rings, and let $\mathscr{X}'$ be the base change of $\mathscr{X}$ on $\mathscr{S}'=(\Spf V')^{\rig}$.

{\rm (1)} Show that, if the extension $V\subseteq V'$ is finite, then an open subspace $\mathscr{U}$ of $\mathscr{X}$ is a tube open subset if and only if its inverse image $\mathscr{U}'$ by the morphism $\mathscr{X}'\rightarrow\mathscr{X}$ is a tube open subset of $\mathscr{X}'$.

{\rm (2)} Find an example of an overconvergent open subspace $\mathscr{U}$ of $\mathscr{X}$ that is not a tube open subset but the inverse image $\mathscr{U}'$ is a tube open subset of $\mathscr{X}'$.}
\end{exer}

\begin{exer}\label{exer-reducedfibertheorem}{\rm 
Use \ref{thm-MMP} (1) and \ref{thm-comparisonaffinoid} to deduce the following elementary version of the reduced fiber theorem (cf.\ \cite{BLR2}). 
{\it Let $V$ be an $a$-adically complete valuation ring of height one, $X$ a coherent finite type flat formal scheme over $V$, and $\mathscr{X}=X^{\rig}$ the associated coherent rigid space over $K=\Frac(V)$.
Suppose that $\mathscr{X}$ is {\em geometrically reduced } over $K$, that is, for any finite extension $K'/K$ the base-change $\mathscr{X}_{K'}$ is reduced.
Then there exist a finite separable extention $K'/K$ and an admissible blow-up $X'\rightarrow X\otimes_VV'$ where $V'$ denotes the integral closure of $V$ in $K'$ such that $X'\rightarrow\Spf V'$ is flat and has the reduced geometric fiber.}}
\end{exer}


\setcounter{section}{0}
\renewcommand{\thesection}{\Alph{section}}
\renewcommand{\theexer}{{\bf \thechapter}.\Alph{section}.\arabic{exer}}
\section{Appendix: Adic spaces}\label{sec-adicspaces}
In this section we give the relationship between our rigid spaces and adic spaces \cite{Hube1}\cite{Hube2}\cite{Hube3}.
One of our goals is to show that the formation of Zariski-Riemann triple $(\ZR{\mathscr{X}},\O^{\int}_ {\mathscr{X}},\O_{\mathscr{X}})$ associated to a rigid space $\mathscr{X}$ gives a functor 
$$
\ZRT\colon\RNoe\Rf\longrightarrow\Acsp
$$
from the category of locally universally Noetherian rigid spaces\index{rigid space!universally Noetherian rigid space@universally Noetherian ---!locally universally Noetherian rigid space@locally --- ---} (\ref{dfn-universallyadhesiverigidspaces}) to the category of adic spaces.

\subsection{Triples}\label{sub-triples}
\index{triple|(}
As we saw in \S\ref{subsub-ZRstrsheaftriple}, to any rigid space $\mathscr{X}$ is canonically associated the Zariski-Riemann triple\index{Zariski-Riemann triple}\index{triple!Zariski Riemann triple@Zariski-Riemann ---} $\ZRT(\mathscr{X})=(\ZR{\mathscr{X}},\O^{\int}_{\mathscr{X}},\O_{\mathscr{X}})$ (\ref{dfn-ZRtriple}).
Huber's adic spaces are, on the other hand, defined as a similar kind of triples.
Therefore, it will be convenient, for the comparison between rigid spaces and adic spaces, to establish a general theory of triples.

\begin{dfn}\label{dfn-triplescategory}{\rm 
(1) A {\em triple} is a data $(X,\O^+_X,\O_X)$ consisting of a topological space $X$ and two sheaves $\O^+_X$ and $\O_X$ of topological rings on $X$ together with an injective morphism $\iota\colon\O^+_X\hookrightarrow\O_X$ of sheaves of topological rings such that 
\begin{itemize}
\item[{\rm (a)}] the injection $\iota$ maps $\O^+_X$ isomorphically onto an open subsheaf of $\O_X$;
\item[{\rm (b)}] $X^+=(X,\O^+_X)$ and $X=(X,\O_X)$ are topologically locally ringed spaces.
\end{itemize}

(2) A {\em morphism} of triples $(X,\O^+_X,\O_X)\rightarrow(Y,\O^+_Y,\O_Y)$ is a morphism of topologically locally ringed spaces $\varphi=(\varphi,h)\colon (X,\O_X)\rightarrow(Y,\O_Y)$ (where $h\colon \varphi^{-1}\O_Y\rightarrow\O_X$) such that the following conditions are satisfied:
\begin{itemize}
\item[{\rm (a)}] $h(\varphi^{-1}\O^+_Y)\subseteq\O^+_X$;
\item[{\rm (b)}] the induced morphism $\varphi^+=(\varphi,h)\colon(X,\O^+_X)\rightarrow(Y,\O^+_Y)$ is a morphism of topologically locally ringed spaces.
\end{itemize}}
\end{dfn}

We denote by $\Triples$ the category of triples.
The following notion of triples will be essential in discussing adic spaces:
\index{triple!valued triple@valued ---|(}
\begin{dfn}\label{dfn-Hubertriples}{\rm 
(1) A {\em valued triple} is a triple $(X,\O^+_X,\O_X)$ together with a set $\{v_x\}_{x\in X}$ consisting of, for each $x\in X$, a continuous (additive) valuation $v_x$ of $\O_{X,x}$ such that $\O^+_{X,x}=\{s\in\O_{X,x}\,|\,v_x(s)\geq 0\}$ and $\m^+_{X,x}$ ($=$ the maximal ideal of $\O^+_{X,x}$) $=\{s\in\O_{X,x}\,|\,v_x(s)>0\}$.
Here, a valuation $v$ on a topological ring $A$ is {\em continuous}\index{valuation!continuous valuation@continuous ---} if for any $\gamma\in\Gamma_v$, where $\Gamma_v$ is the value group\index{value group}, there exists an open neighborhood $U$ in $A$ of $0$ such that $v(x)>\gamma$ for every $x\in U$ (cf.\ \cite[\S3]{Hube1}).

(2) A morphism $\varphi=(\varphi,h)\colon ((X,\O^+_X,\O_X),\{v_x\}_{x\in X})\rightarrow((Y,\O^+_Y,\O_Y),\{v_y\}_{y\in Y})$ of valued triples is a morphism of triples satisfying the following condition$:$
\begin{itemize}
\item[$(\ast)$] for any $x\in X$, $v_{\varphi(x)}$ is equivalent\footnote{See, e.g., \cite[\S2]{Hube1} for the equivalence of valuations; in our case, the condition means that there exists an ordered isomorphism $\phi\colon\Gamma_{v_{\varphi(x)}}\cup\{\infty\}\stackrel{\sim}{\rightarrow}\Gamma_{v_x}\cup\{\infty\}$ (where $\Gamma_{v_x}$ etc.\ are taken to be the value groups\index{value group} ({\bf \ref{ch-pre}}, \S\ref{subsub-valuation}) of the valuations) such that $\phi\circ v_{\varphi(x)}=v_x\circ h_x$.} to $v_x\circ h_x$.
\end{itemize}}
\end{dfn}

We denote by $\VTriples$ the category of valued triples.

\begin{rem}\label{rem-Hubertriples}{\rm 
Notice that, if $((X,\O^+_X,\O_X),\{v_x\}_{x\in X})$ and $((Y,\O^+_Y,\O_Y),\{v_y\}_{y\in Y})$ are valued triples, any morphism of topologically locally ringed spaces $\varphi=(\varphi,h)\colon (X,\O_X)\rightarrow(Y,\O_Y)$ with the property $(\ast)$ as above gives automatically a morphism of valued triples.}
\end{rem}

\index{triple!analytic triple@analytic ---|(}
\begin{dfn}\label{dfn-analytictriples}{\rm 
A triple $(X,\O^+_X,\O_X)$ is said to be {\em analytic} if for any $x\in X$ there exist an open neighborhood $U\subseteq X$ of $x$ and an open ideal sheaf $\mathscr{I}$ of $\O^+_X|_U$ such that the following conditions are satisfied$:$
\begin{itemize}
\item[{\rm (a)}] for any $y\in U$, the ideal $\mathscr{I}_y\subseteq\O^+_{X,y}$ is finitely generated, and the topology of $\O^+_{X,y}$ is $\mathscr{I}_y$-adic$;$
\item[{\rm (b)}] for any $y\in U$, the ring $\O^+_{X,y}$ is $\mathscr{I}_y$-valuative\index{valuative!Ivaluative ring@$I$-{---} ring} {\rm ({\bf \ref{ch-pre}}.\ref{dfn-valuative1})}, and $\O_{X,y}=\varinjlim_{n\geq 1}\Hom(\mathscr{I}_y^n,\O^+_{X,y})$.
\end{itemize}}
\end{dfn}

The ideal sheaf $\mathscr{I}$ as above will be called an {\em ideal of definition}\index{ideal of definition} over $U$.
We denote by $\AnTriples$ the full subcategory of $\Triples$ consisting of analytic triples.

\begin{rem}\label{rem-analytictriples}{\rm 
If $(X,\O^+_X,\O_X)$ is an analytic triple, then for any $x\in X$ the topological ring $\O_{X,x}$ is an extremal f-adic ring\index{fadic ring@f-adic ring!extremal fadic ring@extremal ---}  ({\bf \ref{ch-pre}}, \S\ref{subsub-extremalfadicrings}) with $\O^+_{X,x}$ a ring of definition\index{ring of definition}.
In particular, by {\bf \ref{ch-pre}}.\ref{cor-typeRfadicrings}, for any morphism $\varphi\colon (X,\O^+_X,\O_X)\rightarrow(Y,\O^+_Y,\O_Y)$ of triples between analytic triples, the map $\O^+_{Y,\varphi(x)}\rightarrow\O^+_{X,x}$ for any $x\in X$ is adic.}
\end{rem}

Let $(X,\O^+_X,\O_X)$ be an analytic triple.
Then by {\bf \ref{ch-pre}}.\ref{thm-valuative}, the ring $B_x=\O_{X,x}$ for any $x\in X$ has the canonical valuation $v_x\colon\O_{X,x}\rightarrow\Gamma_x\cup\{\infty\}$ such that $\O^+_{X,x}=\{f\in B_x\,|\,v_x(f)\geq 0\}$ and $\m^+_{X,x}=\{f\in B_x\,|\,v_x(f)>0\}$.
Indeed, since $A_x=\O^+_{X,x}$ is $a$-valuative and $B_x=A_x[\frac{1}{a}]$ (where $\mathscr{I}_x=(a)$), the ring $V_x=A_x/J_x$, where $J_x=\bigcap_{n\geq 1}a^nA_x$, is an $a$-adically separated valuation ring of the fractional field $K_x=B_x/J_x$, whence inducing the valuation $v_x$ as above.
Moreover, it follows from \ref{rem-analytictriples} and Exercise \ref{exer-Ivaluativeringsmaps} that, for any morphism $(\varphi,h)\colon(X,\O^+_X,\O_X)\rightarrow(Y,\O^+_Y,\O_Y)$ of analytic triples and for any $x\in X$, the valuations $v_{\varphi(x)}$ and $v_x\circ h_x$ are equivalent.
Thus we have a canonical functor $\AnTriples\rightarrow\VTriples$, which is clearly fully faithful.

To sum up, we have the commutative diagram of categories
$$
\begin{xy}
(0,0)="VT"*{\VTriples},+<10em,0ex>="T"*{\Triples},+<0em,-10ex>="AT"*{\AnTriples\rlap{,}}
\ar"VT"+<1.6em,0ex>;"T"+<-1.2em,0ex>^{\mathrm{forgetful}}
\ar@{^{(}->}"AT"+<0em,2.6ex>;"T"+<0em,-2ex>
\ar@{_{(}->}"AT"+<-2em,1ex>;"VT"+<1.5em,-1.7ex>
\end{xy}
$$
where the arrows denoted by $\hookrightarrow$ are fully faithful.

\begin{thm}\label{thm-fullyfaithfultriples}
{\rm (1)} By $\mathscr{X}\mapsto\ZRT(\mathscr{X})=(\ZR{\mathscr{X}},\O^{\int}_{\mathscr{X}},\O_{\mathscr{X}})$, we have a functor, again defined by $\ZRT$, from the category $\Rf$ of rigid spaces to the category $\AnTriples$ of analytic triples.

{\rm (2)} The functor $\ZRT\colon\mathbf{RigNoeRf}\rightarrow\AnTriples$, restricted on the category of  locally universally Noetherian rigid spaces,  is faithful.
\end{thm}

\begin{proof}
(1) Clear.

(2) By the patching argument, it is enough to discuss morphisms between Stein affinoids of the form $\mathscr{X}=(\Spf A)^{\rig}$, where $A$ has an invertible ideal of definition $aA$.
By \ref{prop-corlemtuavsrigidaff} and \ref{thm-comparisonaffinoid}, we have $\Gamma(\mathscr{X},\O_X)=A[\frac{1}{a}]$ and $\Gamma(\mathscr{X},\O^{\int}_X)=A^{\int}$, where $A^{\int}$ is the integral closure of $A$ in $A[\frac{1}{a}]$.
Suppose we are given two morphisms $\varphi,\psi\colon\mathscr{X}=(\Spf A)^{\rig}\rightarrow\mathscr{Y}=(\Spf B)^{\rig}$ between Stein affinoids of the above form.
In view of \ref{prop-affinoidblowupglobalsection} and \ref{prop-lemmorbetweenaffinoid2}, we may assume that these morphisms are induced from two adic morphisms $f,g\colon B\rightarrow A$, respectively.
If $\ZRT(\varphi)$ and $\ZRT(\psi)$ gives the same morphisms of triples, then we have $f[\frac{1}{a}]=g[\frac{1}{a}]$, as a homomorphism $B[\frac{1}{a}]\rightarrow A[\frac{1}{a}]$.
Since $aA$ is an invertible ideal of $A$, $A\rightarrow A[\frac{1}{a}]$ is injective, and hence $f=g$.
\end{proof}
\index{triple!valued triple@valued ---|)}\index{triple!analytic triple@analytic ---|)}\index{triple|)}

\subsection{Rigid f-adic rings}\label{sec-adicsprev}
\index{adic!adic ring@--- ring!fadic ring@f-{---} ---|(}\index{fadic ring@f-adic ring|(}
\subsubsection{T.u.\ rigid-Noetherian f-adic rings}\label{subsub-turigidnoetherianfadicrings}
First recall that for a complete f-adic ring\index{fadic ring@f-adic ring!complete fadic ring@complete ---} $A$ ({\bf \ref{ch-pre}}, \S\ref{subsub-completeaffinoidrings}) we have defined in {\bf \ref{ch-pre}}, \S\ref{subsub-completeaffinoidrings} the restricted power series ring $A\dl X_1,\ldots,X_n\dr$.
If $A_0\subseteq A$ is a ring of definition\index{ring of definition} ({\bf \ref{ch-pre}}, \S\ref{subsub-fadicringsgeneralities}), then $A_0$ is an adic ring of finite ideal type\index{adic!adic ring@--- ring!adic ring of finite ideal type@--- --- of finite ideal type} ({\bf \ref{ch-formal}}.\ref{dfn-admissibleringoffiniteidealtype}), and we have 
$$
A\dl X_1,\ldots,X_n\dr=A_0\dl X_1,\ldots,X_n\dr\otimes_{A_0}A.
$$
Recall also that for an ideal of definition $I_0\subseteq A_0$ we have the equality
\begin{equation*}
\begin{split}
\Spec A\dl X_1,\ldots,X_n\dr&\setminus V(I_0A\dl X_1,\ldots,X_n\dr)\\ &=\Spec A_0\dl X_1,\ldots,X_n\dr\setminus V(I_0A_0\dl X_1,\ldots,X_n\dr)
\end{split}
\end{equation*}
due to {\bf \ref{ch-pre}}.\ref{prop-lemfadicbasicfact11}.
Since the topology on $A$ coincides with the one induced from the filtration $\{I^k_0\}_{k\geq 0}$, the scheme $\Spec A\dl X_1,\ldots,X_n\dr\setminus V(I_0A\dl X_1,\ldots,X_n\dr)$ does not depend on the choice of the ideal of definition $I_0$.
\begin{dfn}\label{dfn-turigidnoetherianfadicrings}{\rm 
A complete f-adic ring $A$ is said to be {\em t.u.\ rigid-Noetherian}\index{fadic ring@f-adic ring!topologically universally rigid Noetherian fadic ring@t.u.\ rigid-Noetherian ---} if $\Spec A\dl X_1,\ldots,X_n\dr\setminus V(I_0A\dl X_1,\ldots,X_n\dr)$ $($defined as above$)$ is a Noetherian scheme for any $n\geq 0$.}
\end{dfn}

Notice that, if $A$ is a bounded complete f-adic ring ({\bf \ref{ch-pre}}, \S\ref{subsub-extremalfadicrings}) or, equivalently, an adic ring of finite ideal type\index{adic!adic ring@--- ring!adic ring of finite ideal type@--- --- of finite ideal type}, then the above definition coincides with the one in {\bf \ref{ch-formal}}.\ref{dfn-tuaringadmissible}, since we have $A=A_0$ in the above notation.
Note also that, if $A$ is a Tate ring, this notion coincides with `strongly Noetherian' defined in \cite[\S2]{Hube2}.
The following proposition follows immediately from {\bf \ref{ch-pre}}.\ref{prop-lemfadicbasicfact11}:
\begin{prop}\label{prop-stronglynoetherianrigidnotherian}
A complete f-adic ring $A$ is t.u.\ rigid-Noetherian\index{fadic ring@f-adic ring!topologically universally rigid Noetherian fadic ring@t.u.\ rigid-Noetherian ---} if and only if any ring of definition $A_0\subseteq A$ is t.u.\ rigid-Noetherian\index{t.u. rigid-Noetherian ring@t.u.\ rigid-Noetherian ring} $($in the sense as in {\rm ({\bf \ref{ch-formal}}.\ref{dfn-tuaringadmissible})}$)$. \hfill$\square$
\end{prop}

\subsubsection{Finite type extensions}\label{subsub-finitetypeextensions}
In this paragraph we discuss finite type morphisms between complete f-adic rings in a special situation; for the general definition, see \cite[\S3]{Hube2}.

Let $A$ be a t.u.\ rigid Noetherian f-adic ring\index{fadic ring@f-adic ring!topologically universally rigid Noetherian fadic ring@t.u.\ rigid-Noetherian ---}, $A_0\subseteq A$ a ring of definition, and $I_0\subseteq A_0$ a finitely generated ideal of definition.
Consider an ideal $J\subseteq A$, and set $J_0=J\cap A_0$.
As we saw in {\bf \ref{ch-formal}}.\ref{rem-turigidnoetherianbasicproperties}, the ideal $J_0\subseteq A_0$ is closed in $A_0$, and the quotient $B_0=A_0/J_0$ is again t.u.\ rigid-Noetherian\index{t.u. rigid-Noetherian ring@t.u.\ rigid-Noetherian ring}.
The quotient $B=A/J$ is obviously an f-adic ring, and $B_0\subseteq B$ is a ring of definition. 
Then by {\bf \ref{ch-pre}}.\ref{prop-completefadicbasic1} we see that $B$ is again complete.
Since $B_0$ is t.u.\ rigid-Noetherian, $B$ is a t.u.\ rigid Noetherian f-adic ring\index{fadic ring@f-adic ring!topologically universally rigid Noetherian fadic ring@t.u.\ rigid-Noetherian ---} by \ref{prop-stronglynoetherianrigidnotherian}.

Let $\varphi\colon A\rightarrow B$ be a continuous homomorphism between t.u.\ rigid Noetherian\index{fadic ring@f-adic ring!topologically universally rigid Noetherian fadic ring@t.u.\ rigid-Noetherian ---} extremal f-adic rings; notice that by {\bf \ref{ch-pre}}.\ref{cor-typeRfadicrings} the map $\varphi$ is automatically adic.
We say that $\varphi$ is {\em of finite type} if the ring $B$ is, as an $A$-algebra, isomorphic to an $A$-algebra of the form
$$
A\dl X_1,\ldots,X_n\dr/\mathfrak{a}
$$
by an ideal $\mathfrak{a}\subseteq A\dl X_1,\ldots,X_n\dr$ (cf.\ \cite[3.3 (iii)]{Hube2}).
By the above observation and {\bf \ref{ch-pre}}.\ref{cor-openmappingtheorem}, the topological ring $B$ is topologically isomorphic, hence isomorphic as f-adic rings, to the above quotient of $A\dl X_1,\ldots,X_n\dr$.

\subsubsection{Rigidification of f-adic rings}\label{subsub-rigidificationfadicrings}
\begin{dfn}\label{dfn-strictlyequiv}{\rm 
Let $A$ be a t.u.\ rigid-Noetherian f-adic ring\index{fadic ring@f-adic ring!topologically universally rigid Noetherian fadic ring@t.u.\ rigid-Noetherian ---}, and $A_0,A_1\subseteq A$ two rings of definition of $A$.
We say that $A_0$ and $A_1$ are {\em strictly} (resp.\ {\em finitely}) {\em equivalent}\index{strictly equivalent}\index{finitely equivalent} if there exists a diagram $A_0\hookrightarrow A_2\hookleftarrow A_1$ consisting of strict (resp.\ finite) weak isomorphisms\index{weak!weak isomorphism@--- isomorphism!strict weak isomorphism@strict --- ---} (\ref{dfn-strictweakisomorphism}) between rings of definition of $A$.}
\end{dfn}

It follows from \ref{prop-affinoidblowupglobalsection} that, if $A_0$ and $A_1$ are t.u.\ adhesive\index{t.u.a. ring@t.u.\ adhesive ring} {\rm ({\bf \ref{ch-formal}}.\ref{dfn-tuaringadmissible})}, then strict equivalence is equivalent to finite equivalence.
For a t.u.\ rigid-Noetherian f-adic ring $A$ and a ring of definition $B\subseteq A$, we denote by $C(B)$ the strict equivalence class that contains $B$.
\begin{dfn}\label{dfn-rigidificationfadic}{\rm 
(1) A strict equivalence class of rings of definition of a t.u.\ rigid-Noetherian f-adic ring $A$ is called a {\em rigidification}\index{rigidification} of $A$.

(2) A pair $(A,C)$ consisting of t.u.\ rigid-Noetherian f-adic ring $A$ and a rigidification of $A$ is called an {\em f-r-pair}\index{pair!frpair@f-r-{---}}\index{frpair@f-r-pair}.}
\end{dfn}

\begin{exa}\label{exa-rigidificationvalution}{\rm 
Let $V$ be an $a$-adically complete valuation ring\index{valuation!valuation ring@--- ring!a-adically complete valuation ring@$a$-adically complete --- ---} ($a\in\m_V\setminus\{0\}$), and $K=\Frac(V)$.
Then $K$ with the topology defined by the filtration $\{a^nV\}_{n\geq 0}$ is a t.u.\ rigid-Noetherian f-adic ring having $V$ as a ring of definition.
Let $\mathfrak{p}=\sqrt{(a)}$ be the associated height one prime\index{associated height one prime} ({\bf \ref{ch-pre}}.\ref{dfn-maxspe2}), and consider the height one localization $\til{V}=V_{\mathfrak{p}}$.
Then $\til{V}$ coincides with the set $K^o$ of all power-bounded elements and is also a ring of definition; 
If $\mathrm{ht}(V)>1$, then $V$ and $\til{V}$ are not strictly equivalent, since $(\Spf V)^{\rig}\not\cong(\Spf\til{V})^{\rig}$.}
\end{exa}

If $(A,C)$ is an f-r-pair\index{pair!frpair@f-r-{---}}\index{frpair@f-r-pair}, then the affinoid $\mathscr{X}=(\Spf A_0)^{\rig}$ for $A_0\in C$ does not depend, up to isomorphisms, on the choice of the rings of definition $A_0$.
The affinoid $\mathscr{X}$ thus obtained will be called the {\em affinoid associated to the f-r-pair $(A,C)$}.

\begin{dfn}\label{dfn-rigidificationfadicmorphism}{\rm 
Let $(A,C)$ and $(A',C')$ be two f-r-pairs.\index{pair!frpair@f-r-{---}}\index{frpair@f-r-pair} 
A continuous homomorphism $\varphi\colon A\rightarrow A'$ is said to be {\em rigid} if there exist $A_0\in C$ and $A'_0\in C'$ such that $\varphi(A_0)\subseteq A'_0$ and the induced morphism $A_0\rightarrow A'_0$ is adic.}
\end{dfn}

Let $\varphi\colon A\rightarrow B$ be a finite type morphism between t.u.\ rigid-Noetherian extremal f-adic rings (\S\ref{subsub-finitetypeextensions}). 
If $A$ is equipped with a rigidification\index{rigidification} (\ref{dfn-rigidificationfadic}) $C$, then $B$ has a canonically induced rigidification defined as follows.
As we have seen in \S\ref{subsub-finitetypeextensions}, $B$ is topologically isomorphic to an f-adic ring of the form $A\dl X_1,\ldots,X_n\dr/\mathfrak{a}$.
Then for any $A_0\subseteq C$ the ring $A_0\dl X_1,\ldots,X_n\dr/\mathfrak{a}_0$, where $\mathfrak{a}_0=\mathfrak{a}\cap A_0\dl X_1,\ldots,X_n\dr$, is a ring of definition of $B$ and hence defines a rigidification of $B$.
Note that the rigidification of $B$ thus obtained does not depend on the choice of $A_0\in C$.
\index{adic!adic ring@--- ring!fadic ring@f-{---} ---|)}\index{fadic ring@f-adic ring|)}

\subsection{Adic spaces}\label{sub-affinoidringsadicspectrums}
\subsubsection{Affinoid rings}\label{subsub-affinoidringshuber}
\index{affinoid!affinoid ring@--- ring|(}
\begin{dfn}\label{dfn-affinoidringadicspace}{\rm 
An {\em affinoid ring}\index{affinoid!affinoid ring@--- ring} is a pair $\mathcal{A}=(\mathcal{A}^{\pm},\mathcal{A}^+)$ consisting of an f-adic ring $\mathcal{A}^{\pm}$ and a subring $\mathcal{A}^{+}\subset\mathcal{A}^{\pm}$ that is open, integrally closed in $\mathcal{A}^{\pm}$, and contained in $(\mathcal{A}^{\pm })^o$ (the subset of $\mathcal{A}^{\pm }$ consisting of power-bounded\index{power-bounded} elements ({\bf \ref{ch-pre}}, \S\ref{subsub-fadicringsgeneralities})).\footnote{In Huber's\index{Huber, R.} original notation (such as in \cite[\S 3]{Hube1}) affinoid rings are written as $\mathcal{A}=(\mathcal{A}^{\vartriangleright},\mathcal{A}^+)$.
Here in this book we prefer to denote the ring in the first entry by $\mathcal{A}^{\pm}$ instead of $\mathcal{A}^{\vartriangleright}$.}
The ring $\mathcal{A}^+$ is called the {\em subring of integral elements}\index{subring of integral elements} of $\mathcal{A}$.}
\end{dfn}

Let $\mathcal{A}=(\mathcal{A}^{\pm},\mathcal{A}^+)$ be an affinoid ring.
By a {\em ring of definition}\index{ring of definition} of $\mathcal{A}$ we mean a ring of definition of the f-adic ring $\mathcal{A}^{\pm}$ contained in $\mathcal{A}^+$. 
When $A_0$ is a ring of definition of $\mathcal{A}$, a finitely generated ideal of definition $I_0\subseteq A_0$ is called an {\em ideal of definition}\index{ideal of definition} of $\mathcal{A}$.

We say that an affinoid ring $\mathcal{A}=(\mathcal{A}^{\pm},\mathcal{A}^+)$ is {\em complete}\index{affinoid!affinoid ring@--- ring!complete affinoid ring@complete --- ---} (resp.\ {\em extremal}\index{affinoid!affinoid ring@--- ring!extremal affinoid ring@extremal --- ---}) if the f-adic ring $\mathcal{A}^{\pm}$ is complete\index{fadic ring@f-adic ring!complete fadic ring@complete ---} (resp.\ extremal f-adic)\index{fadic ring@f-adic ring!extremal fadic ring@extremal ---}.
For an affinoid ring $\mathcal{A}=(\mathcal{A}^{\pm},\mathcal{A}^+)$ its {\em completion} is given by the pair $\widehat{\mathcal{A}}=(\widehat{\mathcal{A}}^{\pm},\widehat{\mathcal{A}}^+)$ consisting of the completions of the topological rings.
It turns out that this is again an affinoid ring, hence a complete affinoid ring.

\begin{dfn}{\rm 
Let $\mathcal{A}=(\mathcal{A}^{\pm},\mathcal{A}^+) $, $\mathcal{B} = (\mathcal{B}^{\pm},\mathcal{B}^+)$ be affinoid rings. 
A {\em homomorphism} $\varphi\colon\mathcal{A}\rightarrow\mathcal{B}$ of affinoid rings is a continuous homomorphism $\varphi\colon\mathcal{A}^{\pm}\rightarrow\mathcal{B}^{\pm}$ such that $\varphi(\mathcal{A}^+)\subset\mathcal{B}^+$.
It is said to be {\em adic}\index{adic!adic morphism@--- morphism} if $\varphi\colon\mathcal{A}^{\pm}\rightarrow\mathcal{B}^{\pm}$ is adic ({\bf \ref{ch-pre}}, \S\ref{subsub-fadicringsgeneralities}).}
\end{dfn}
\index{affinoid!affinoid ring@--- ring|)}

\subsubsection{Adic spectrum}\label{subsub-adicspectrum}
Let $\mathcal{A}=(\mathcal{A}^{\pm},\mathcal{A}^+)$ be an affinoid ring.
The associated {\em adic spectrum}\index{adic spectrum} $\Spa\mathcal{A}$ is the topological space defined as follows:
\begin{itemize}
\item as a set, it is the set of all (equivalence classes of) valuations\footnote{Although in \cite{Hube1} valuations are written multiplicatively, here we prefer to write them additively.}\index{valuation} $v\colon\mathcal{A}^{\pm}\rightarrow\Gamma\cup\{\infty\}$ (cf.\ {\bf \ref{ch-pre}}.\ref{dfn-valuation1}) that satisfy $v(x)\geq 0$ for $x\in\mathcal{A}^+$ and are continuous\index{valuation!continuous valuation@continuous ---} (in the sense as in \ref{dfn-Hubertriples} (1)); 
\item the topology is the one generated by the subsets of the form $\{v\,|\,v(x)\geq v(y)\neq\infty\}$ for any $x,y\in\mathcal{A}^{\pm}$.
\end{itemize}

It is known that $\Spa\mathcal{A}$ is a coherent\index{space@space (topological)!coherent topological space@coherent ---}\index{coherent!coherent topological space@--- (topological) space} sober\index{space@space (topological)!sober topological space@sober ---} topological space (cf.\ \cite[Theorem 3.5 (i)]{Hube1}).

\begin{dfn}\label{dfn-rationalsubsetadic}{\rm 
A subset $U$ of $\Spa\mathcal{A}$ is said to be {\em rational}\index{rational subset} if there exist $f_0,\ldots,f_n\in\mathcal{A}^{\pm}$ such that the ideal $(f_0,\ldots,f_n)$ is open and 
$$
U=\{v\in\Spa\mathcal{A}\,|\,v(f_i)\geq v(f_0)\neq\infty,\ i=1,\ldots,n\}.
$$}
\end{dfn}

The rational subset as above is often denoted by $R(\frac{f_1,\ldots,f_n}{f_0})$.
Clearly, these subsets give a basis of the topology on $\Spa\mathcal{A}$.

Let $U=R(\frac{f_1,\ldots,f_n}{f_0})$ be a rational subset of $X=\Spa\mathcal{A}$.
Define the topological ring $\mathcal{A}^{\pm}(\frac{f_1,\ldots,f_n}{f_0})$ as follows:
\begin{itemize}
\item as a ring, it is $\mathcal{A}^{\pm}[\frac{1}{f_0}]$;
\item it has the ring of definition $B[\frac{f_1}{f_0},\ldots,\frac{f_n}{f_0}]$ with the ideal of definition $IB[\frac{f_1}{f_0},\ldots,\frac{f_n}{f_0}]$, where $B$ is a ring of definition of $\mathcal{A}^{\pm}$ with the ideal of definition $I\subseteq B$.
\end{itemize}

We denote the completion of $\mathcal{A}^{\pm}(\frac{f_1,\ldots,f_n}{f_0})$ by 
$$
\mathcal{A}^{\pm}\dl{\textstyle \frac{f_1,\ldots,f_n}{f_0}}\dr.
$$
It turns out that this ring is an f-adic ring determined only by $U$.
Define the presheaf $\O_X$ by 
$$
\O_X(U)=\mathcal{A}^{\pm}\dl{\textstyle \frac{f_1,\ldots,f_n}{f_0}}\dr
$$
for any rational subset $U\subseteq X=\Spf\mathcal{A}$ written as above, and $\O_X(V)$ for any open subset $V$ by the projective limit $\varprojlim_{U\subseteq V}\O_X(U)$ taken over all rational subsets contained in $V$ endowed with the projective limit topology.
We also define the presheaf $\O^+_X$ by
$$
\O^+_X(U)=\{f\in\O_X(U)\,|\,v(f)\geq 0\ \textrm{for any}\ v\in U\}
$$
for any open subset $U$.

If $\mathcal{A}$ is a complete affinoid ring\index{affinoid!affinoid ring@--- ring!complete affinoid ring@complete --- ---}, then we have $\Gamma(\Spa\mathcal{A},\O_X)=\mathcal{A}^{\pm}$ and $\Gamma(\Spa\mathcal{A},\O^+_X)=\mathcal{A}^+$ (cf.\ \cite[Prop.\ 1.6 (iv)]{Hube2}).

By definition, for any $x\in X=\Spa\mathcal{A}$, the stalk $\O_{X,x}$ is canonically equipped with the continuous valuation\index{valuation!continuous valuation@continuous ---} $v_x\colon\O_{X,x}\rightarrow\Gamma_x\cup\{\infty\}$, and we have $\O^+_{X,x}=\{f\in\O_{X,x}\,|\,v_x(f)\geq 0\}$.

\begin{prop}[{\rm \cite[Prop.\ 1.6 (i), (ii)]{Hube2}}]\label{prop-affinoidadicspace001}
{\rm (1)} For any $x\in X=\Spa\mathcal{A}$, the stalk $\O_{X,x}$ is a local ring with the maximal ideal $\m_{X,x}=\{f\in\O_{X,x}\,|\,v_x(f)=\infty\}$.

{\rm (2)} For any $x\in X=\Spa\mathcal{A}$, the stalk $\O^+_{X,x}$ is a local ring with the maximal ideal $\m^+_{X,x}=\{f\in\O_{X,x}\,|\,v_x(f)>0\}$. \hfill$\square$
\end{prop}

\subsubsection{Adic spaces}\label{subsub-adicspaces}
\index{adic space|(}
The presheaf $\O_X$ on $X=\Spa\mathcal{A}$ may not be a sheaf (cf.\ \cite[\S1]{Hube2}), but if it is a sheaf, then so is $\O^+_X$.
In this case, one obtains a valued triple\index{triple!valued triple@valued ---} (\ref{dfn-Hubertriples})
$$
((X=\Spa\mathcal{A},\O^+_X,\O_X),\{v_x\}_{x\in X}).\leqno{(\ast)}
$$

\begin{dfn}\label{dfn-affinoidadicspace}{\rm 
An {\em affinoid adic space}\index{adic space!affinoid adic space@affinoid ---} is a valued triple\index{triple!valued triple@valued ---} that is isomorphic in $\VTriples$ to a triple of the form $(\ast)$ as above, where $\O_X$ is assumed to be a sheaf.}
\end{dfn}

Notice that, if $X=\Spa\mathcal{A}$ is an affinoid adic space, then any rational subdomain $U\subseteq X$, endowed with the structures induced from $X$, is also an affinoid adic space.

\begin{rem}\label{rem-Hubercategoryconsistency}{\rm 
In \cite[\S2]{Hube2} and \cite[\S1.1]{Hube3}, Huber\index{Huber, R.} defined affinoid adic spaces as objects of a category $\mathscr{V}$, which is defined as follows: objects of $\mathscr{V}$ are triples of the form $(X,\O_X,\{v_x\}_{x\in X})$ consisting of a topological space $X$, a sheaf of complete topological rings $\O_X$, and a collection of valuations $\{v_x\}_{x\in X}$ of the stalks $\O_{X,x}$; a morphism $\varphi=(\varphi,h)\colon(X,\O_X,\{v_x\}_{x\in X})\rightarrow(Y,\O_Y,\{v_y\}_{y\in Y})$ in $\mathscr{Y}$ is a morphism of topologically ringed spaces satisfying the compatibility condition $(\ast)$ as in \ref{dfn-Hubertriples} (2).
However, by what we have remarked in \ref{rem-Hubertriples}, our definition is equivalent to Huber's.
The same remark also applies to the definitions of morphisms of affinoid adic spaces, of adic spaces, and of morphisms of adic spaces, which will come soon below.}
\end{rem}

Morphisms of affinoid adic spaces are defined to be morphisms in $\VTriples$, in view of the following:
\begin{prop}[{\rm \cite[Prop.\ 2.1 (i)]{Hube2}}]\label{prop-Huberfullyfaithfulness}
Let $X=\Spa\mathcal{A}$ and $Y=\Spa\mathcal{B}$ be affinoid adic spaces, where $\mathcal{A}$ is complete.

{\rm (1)} Every homomorphism $\mathcal{B}\rightarrow\mathcal{A}$ canonically induces a morphism $X\rightarrow Y$ of valued triples.

{\rm (2)} The mapping
$$
\{\mathcal{B}\rightarrow\mathcal{A}\,|\,\textrm{homomorphism of affinoid rings}\}\longrightarrow\Hom_{\VTriples}(X,Y)
$$
thus obtained is bijective. \hfill$\square$
\end{prop}

\begin{dfn}\label{dfn-adicspacesbytriples}{\rm 
(1) An {\em adic space} is a valued triple\index{triple!valued triple@valued ---} $((X,\O^+_X,\O_X),\{v_x\}_{x\in X})$ that is locally an affinoid adic space.

(2) A morphism $((X,\O^+_X,\O_X),\{v_x\}_{x\in X})\rightarrow((Y,\O^+_Y,\O_Y),\{v_y\}_{y\in Y})$ of adic spaces is a morphism of valued triples.}
\end{dfn}

We denote by $\Acsp$ the category of adic spaces.

\begin{prop}\label{prop-adicspacevaluativeHuber}
{\rm (1)} The underlying topological space of an adic space is a valuative space\index{valuative!valuative topological space@--- (topological) space}\index{space@space (topological)!valuative topological space@valuative ---} {\rm ({\bf \ref{ch-pre}}.\ref{dfn-valuativespace})}.

{\rm (2)} The underlying continuous map of a morphism of adic spaces is valuative\index{valuative!valuative map@--- map} {\rm ({\bf \ref{ch-pre}}.\ref{dfn-valuativemaps})} and locally quasi-compact\index{map@map (continuous)!quasi-compact map@quasi-compact ---!locally quasi-compact map@locally --- ---}\index{quasi-compact!quasi-compact map@--- map!locally quasi-compact map@locally --- ---} {\rm ({\bf \ref{ch-pre}}.\ref{dfn-locallyquasicompactmaps})}.
\end{prop}

\begin{proof}
(1) The underlying topological space of an adic space $X$ is clearly a locally coherent sober space.
It follows from \cite[Lemma 1.1.10 (i)]{Hube3} that, for any $x\in X$, the set $G_x$ of all generizations of $x$ is totally ordered.
Then the assertion follows from {\bf \ref{ch-pre}}.\ref{rem-valuativespace} (1).

(2) It follows from \cite[Lemma 1.1.10 (iv)]{Hube3} that the underlying continuous map of a morphism of adic spaces is valutive.
It is clearly locally quasi-compact, for morphisms between affinoid adic spaces are quasi-compact.
\end{proof}

\subsubsection{Analytic adic spaces}\label{subsub-analyticadicspaces}
\index{adic space!analytic adic space@analytic ---|(}
In the sequel, let us say that an affinoid ring $\mathcal{A}=(\mathcal{A}^{\pm},\mathcal{A}^+)$ is {\em extremal}\index{affinoid!affinoid ring@--- ring!extremal affinoid ring@extremal --- ---} (resp.\ {\em Tate})\index{affinoid!affinoid ring@--- ring!Tate affinoid ring@Tate --- ---} if the f-adic ring $\mathcal{A}^{\pm}$ is extremal\index{fadic ring@f-adic ring!extremal fadic ring@extremal ---} (resp.\ Tate)\index{Tate ring} ({\bf \ref{ch-pre}}, \S\ref{subsub-extremalfadicrings})

\begin{dfn}\label{dfn-analyticadicspaces}{\rm 
Let $X=((X,\O^+_X,\O_X),\{v_x\}_{x\in X})$ be an adic space.
A point $x\in X$ is said to be {\em analytic} if there exists an affinoid open neighborhood $U=\Spa\mathcal{A}$ of $x$ by an extremal affinoid ring $\mathcal{A}$.
If all points of $X$ are analytic, we say that the adic space $X$ is {\em analytic}.}
\end{dfn}

Notice that a point $x\in X$ is analytic if and only if there exists an affinoid open neighborhood $U=\Spa\mathcal{A}$ of $x$ by a Tate affinoid ring\index{affinoid!affinoid ring@--- ring!Tate affinoid ring@Tate --- ---} $\mathcal{A}$.
Indeed, $U=\Spa\mathcal{A}$ is an affinoid open neighborhood of $x$ by an extremal affinoid ring\index{affinoid!affinoid ring@--- ring!extremal affinoid ring@extremal --- ---} $\mathcal{A}$, and $I=(f_0,\ldots,f_n)$ is an ideal of definition of a ring of definition of $\mathcal{A}^{\pm}$, then the rational subsets $U(\frac{f_0,\ldots,f_n}{f_i})$ $(i=0,\ldots,n)$ cover $U$, and the corresponding f-adic rings $\mathcal{A}^{\pm}\dl \frac{f_0,\ldots,f_n}{f_i}\dr$ are Tate rings.

Notice also that any morphism $((X,\O^+_X,\O_X),\{v_x\}_{x\in X})\rightarrow((Y,\O^+_Y,\O_Y),\{v_y\}_{y\in Y})$ of analytic adic spaces is {\em adic}, that is, for any affinoid open $U\subseteq X$ and affinoid open $V\subseteq Y$ with $U\subseteq V$, $\O_Y(V)\rightarrow\O_X(U)$ is adic (due to {\bf \ref{ch-pre}}.\ref{cor-typeRfadicrings}; cf.\ \cite[Prop.\ 3.2 (i)]{Hube2}).

We denote by $\AnAcsp$ the full subcategory of $\Acsp$ consisting of analytic adic spaces.

\begin{prop}\label{prop-adicspacevaluativeHuber2}
Let $\mathcal{A}=(\mathcal{A}^{\pm},\mathcal{A}^+)$ be an extremal affinoid ring\index{affinoid!affinoid ring@--- ring!extremal affinoid ring@extremal --- ---}, and $I\subseteq\mathcal{A}^{\pm}$ be an ideal of definition of a ring of definition of $\mathcal{A}^{\pm}$.
Consider the associated affinoid adic space $((X=\Spa\mathcal{A},\O^+_X,\O_X),\{v_x\}_{x\in X})$, and $x\in X$.

{\rm (1)} The ring $\O^+_{X,x}$ is $I\O^+_{X,x}$-valuative\index{valuative!Ivaluative ring@$I$-{---} ring} {\rm ({\bf \ref{ch-pre}}.\ref{dfn-valuative1})}. Moreover, we have $\O_{X,x}=\varinjlim_{n\geq 1}\Hom(I^n\O^+_{X,x},\O^+_{X,x})$.

{\rm (2)} The valuation on $\O_{X,x}$ induced from the valuation ring $V_x=\O^+_{X,x}/J_x$ $($where $J_x=\bigcap_{n\geq 1}I^n\O^+_{X,x})$ is equivalent to $v_x$.
\end{prop}

\begin{proof}
Set $A_x=\O^+_{X,x}$ and $B_x=\O_{X,x}$.
Then $IA_x=(a)$ for a non-zero-divisor $a\in A_x$.
By {\bf \ref{ch-pre}}.\ref{prop-lemfadicbasicfact11} we have $B_x=\bigcup_{n\geq 0}[A_x:I^nA_x]$.
By this and $IB_x=B_x$, we have $B_x=\varinjlim_{n\geq 1}\Hom(I^nA_x,A_x)$.
Since the valuation $v_x$ is continuous, we have $J_x=\bigcap_{n\geq 1}I^nA_x=\{f\in B_x\,|\,v_x(f)=\infty\}$, which is a prime ideal both in $A_x$ and $B_x$.
By \cite[Prop.\ 1.6 (i)]{Hube2}, $J_x$ is a maximal ideal of $B_x$.
Hence $K_x=B_x/J_x$ is a fields, which contains $V_x=A_x/J_x$.
Since $V_x=\{s\in K_x\,|\,v_x(s)\geq 0\}$, $V_x$ is a valuation ring, which is clearly $a$-adically separated; in particular, we have $K_x=V_x[\frac{1}{a}]=\Frac(V_x)$.
Then by {\bf \ref{ch-pre}}.\ref{thm-valuative} (2), $A_x$ is $IA_x$-valuative, whence (1).
(2) is clear.
\end{proof}

It follows from the proposition that analytic adic spaces are analytic triples\index{triple!analytic triple@analytic ---} (\ref{dfn-analytictriples}).
By \ref{prop-adicspacevaluativeHuber2} and Exercise \ref{exer-Ivaluativeringsmaps} (cf.\ \ref{rem-analytictriples}), we have:
\begin{cor}\label{cor-adicspacevaluativeHuber2}
The forgetful functor by $((X,\O^+_X,\O_X),\{v_x\}_{x\in X})\mapsto(X,\O^+_X,\O_X)$ from the category $\AnAcsp$ of analytic adic spaces to the category $\AnTriples$ of triples is fully faithful. \hfill$\square$
\end{cor}

For the reader's convenience, we insert here the definition of locally of finite type morphisms between analytic adic spaces; see \cite[\S3]{Hube2} for the general definition.
For a complete affinoid ring $\mathcal{A}=(\mathcal{A}^{\pm},\mathcal{A}^+)$, we have the affinoid ring $\mathcal{A}\dl X_1,\ldots,X_n\dr=(\mathcal{A}^{\pm}\dl X_1,\ldots,X_n\dr,\mathcal{A}^+\dl X_1,\ldots,X_n\dr)$ defined as follows: $\mathcal{A}^{\pm}\dl X_1,\ldots,X_n\dr$ is the one already defined in {\bf \ref{ch-pre}}.\ref{subsub-completeaffinoidrings}, and $\mathcal{A}^+\dl X_1,\ldots,X_n\dr$ is the subring consisting of power series with coefficients in $\mathcal{A}^+$.
A homomorphism $\mathcal{A}\rightarrow\mathcal{B}$ of affinoid rings called a quotient mapping if $\mathcal{A}^{\pm}\rightarrow\mathcal{B}^{\pm}$ is surjective, continuous, and open, and the integral closure in $\mathcal{B}^{\pm}$ of the image of $\mathcal{A}^+$ coincides with $\mathcal{B}^+$.
A homomorphism $\mathcal{A}\rightarrow\mathcal{B}$ of complete extremal affinoid rings is said to be {\em $($topologically$)$ of finite type} if it factors through a quotient mapping $\mathcal{A}\dl X_1,\ldots,X_n\dr\rightarrow\mathcal{B}$ (cf.\ \cite[3.5 (iii)]{Hube2}; see also \S\ref{subsub-finitetypeextensions}); notice that, in this case, $\mathcal{B}^{\pm}$ is generated by a finitely many elements over $\mathcal{A}^{\pm}$ (cf.\ {\bf \ref{ch-pre}}, \S\ref{subsub-completeaffinoidrings}).
\begin{dfn}\label{dfn-adicspacefinitetype}{\rm 
A morphism $f\colon X\rightarrow Y$ between analytic adic spaces is said to be {\em locally of finite type} if for any $x\in X$ there exist an open affinoid neighborhood $U\cong\Spa\mathcal{A}$ of $x$ and an affinoid open subset $V\cong\Spa\mathcal{B}$ of $Y$ such that $f(U)\subseteq V$ and that the induced homomorphism $\mathcal{B}\rightarrow\mathcal{A}$ is topologically of finite type.
If, moreover, (the underlying continuous mapping of) $f$ is quasi-compact\index{map@map (continuous)!quasi-compact map@quasi-compact ---}\index{quasi-compact!quasi-compact map@--- map} ({\bf \ref{ch-pre}}.\ref{dfn-quasicompactness} (2)), we say that $f$ is {\em of finite type}.}
\end{dfn}
\index{adic space!analytic adic space@analytic ---|)}\index{adic space|)}

\subsection{Rigid geometry and affinoid rings}\label{sub-rgvsaffring}
\subsubsection{Affinoid rings associated to f-r-pairs}\label{subsub-affinoidringsassociatedtorigidfadicrings}
\begin{prop}\label{prop-affinoidrigidfadicrings}
Let $(A,C)$ be an f-r-pair {\rm (\ref{dfn-rigidificationfadic})}\index{pair!frpair@f-r-{---}}\index{frpair@f-r-pair}, and set
$$
\mathcal{A}^{\pm}=A,\qquad\mathcal{A}^+=\varinjlim_{B\in C}B.
$$
Then $\mathcal{A}=(\mathcal{A}^{\pm},\mathcal{A}^+)$ is a complete affinoid ring.
\end{prop}

\begin{proof}
It is clear that $\mathcal{A}^+$ is open in $\mathcal{A}^{\pm}$.
To show that $\mathcal{A}^+$ is integrally closed in $\mathcal{A}^{\pm}$, suppose $x\in\mathcal{A}^{\pm}$ is integral over $\mathcal{A}^+$.
Then $x$ is integral over some $B\in C$.
The homomorphism $B\rightarrow B[x]$ is finite, and hence $B[x]$ is open and bounded in $\mathcal{A}^{\pm}$.
For a finitely generated ideal of definition $I\subseteq B$, we have by {\bf \ref{ch-pre}}.\ref{prop-lemfadicbasicfact11} the equality $\Spec B\setminus V(I)=\Spec B[x]\setminus V(IB[x])$, and hence $B\rightarrow B[x]$ is isomorphic outside $I$\index{isomorphic outside I@isomorphic outside $I$} (\ref{dfn-morbetweenaffinoid1}).
Hence by \ref{prop-lemmorbetweenaffnoid1} we deduce $B[x]\in C$ and thus $x\in\mathcal{A}^+$.
To verify that $\mathcal{A}^+$ is contained in $(\mathcal{A}^{\pm})^o$, we only have to remark that any element of $\mathcal{A}^+$ is contained in a ring of definition and hence is power-bounded.
\end{proof}

\begin{dfn}\label{dfn-rigidification}{\rm 
(1) An {\em a-r-pair}\index{arpair@a-r-pair}\index{pair!arpair@a-r-pair} is a pair $\mathcal{R}= (\mathcal{A}_{\mathcal{R}},\mathcal{C}_{\mathcal{R}})$ of a complete affinoid ring $\mathcal{A}_{\mathcal{R}}=(\mathcal{A}^{\pm}_{\mathcal{R}},\mathcal{A}^+_{\mathcal{R}})$ with $\mathcal{A}^{\pm}_{\mathcal{R}}$ t.u.\ rigid-Noetherian\index{fadic ring@f-adic ring!topologically universally rigid Noetherian fadic ring@t.u.\ rigid-Noetherian ---}, and a rigidification $\mathcal{C}_{\mathcal{R}}$ of $\mathcal{A}^{\pm}_{\mathcal{R}}$ that satisfies the following compatibility: $\mathcal{A}^+_{\mathcal{R}}=\varinjlim_{B\in\mathcal{C}_{\mathcal{R}}}B$.

(2) For two a-r-pairs $\mathcal{R}=(\mathcal{A}_{\mathcal{R}},\mathcal{C}_{\mathcal{R}})$ and $\mathcal{R}'=(\mathcal{A}_{\mathcal{R}'},\mathcal{C}_{\mathcal{R}'})$ a {\em homomorphism of a-r-pairs} is a homomorphism $\varphi\colon\mathcal{A}_{\mathcal{R}}\rightarrow\mathcal{A}_{\mathcal{R}'}$ enjoying the following property: there exist $B\in\mathcal{C}_{\mathcal{R}}$ and $B'\in\mathcal{C}_{\mathcal{R}'}$ such that $\varphi(B)\subset B'$ and the induced homomorphism $B\rightarrow B'$ is adic (cf.\ \ref{dfn-rigidificationfadicmorphism}).}
\end{dfn}

Note that any homomorphism of a-r-pairs is adic.
The proposition \ref{prop-affinoidrigidfadicrings} shows that any f-r-pair\index{pair!frpair@f-r-{---}}\index{frpair@f-r-pair} $(A,C)$ canonically induces an a-r-pair $\mathcal{R}= (\mathcal{A}_{\mathcal{R}},\mathcal{C}_{\mathcal{R}})$ with $\mathcal{A}^{\pm}_{\mathcal{R}}=A$ and $\mathcal{C}_{\mathcal{R}}=C$.
We call this a-r-pair the {\em a-r-pair associated to the f-r-pair $(A,C)$}.

\subsubsection{Stein affinoids and analytic affinoid pairs}\label{subsub-steinaffrigidaffinoid}
Let $\mathscr{X}$ be a Stein affinoid\index{affinoid!Stein affinoid@Stein ---}\index{Stein affinoid} (\ref{dfn-cohomologyrigidsp1str}) $\mathscr{X}=(\Spf A)^{\rig}$ with $\Spec A\setminus V(I)$ affine, where $I\subseteq A$ a finitely generated ideal of definition.
We may assume that $A$ is $I$-torsion free (\ref{prop-cohomologyrigidsp0}).
As we saw in {\bf \ref{ch-pre}}, \S\ref{subsub-fadicringsgeneralities}, the ring $B$ with $\Spec B=\Spec A\setminus V(I)$ is a complete extremal f-adic ring ({\bf \ref{ch-pre}}, \S\ref{subsub-extremalfadicrings}) having $A$ (resp.\ $I$) as a ring (resp.\ an ideal) of definition.
By \ref{thm-comparisonaffinoid} we have
$$
\Gamma(\mathscr{X},\O_{\mathscr{X}})=B.
$$
Moreover, we have by \ref{prop-affinoidblowupglobalsection},
$$
\Gamma(\mathscr{X},\O^{\int}_{\mathscr{X}})=\varinjlim_{X'\rightarrow\Spf A}\Gamma(X',\O_{X'})=\varinjlim_{A'\in C(A)}A',
$$
where $X'\rightarrow\Spf A$ runs through all admissible blow-ups of $\Spf A$ (see \S\ref{subsub-rigidificationfadicrings} for the definition of $C(A)$).
By \ref{prop-corlemtuavsrigidaff} we have the following results:
\begin{prop}\label{prop-corsteinaffvsadic1corint}
In the situation as above, $\Gamma(\mathscr{X},\O^{\int}_{\mathscr{X}})$ coincides with the integral closure of $A$ in $\Gamma(\mathscr{X},\O_{\mathscr{X}})$. \hfill$\square$
\end{prop}

\begin{thm}\label{thm-steinaffvsadic1}
Let $\mathscr{X}$ be a Stein affinoid\index{affinoid!Stein affinoid@Stein ---}\index{Stein affinoid}. 

{\rm (1)} The pair $\mathcal{A}_{\mathscr{X}}=(\mathcal{A}^{\pm}_{\mathscr{X}},\mathcal{A}^+_{\mathscr{X}})$ defined by 
$$
\mathcal{A}^{\pm}_{\mathscr{X}}=\Gamma(\mathscr{X},\O_{\mathscr{X}}),\qquad\mathcal{A}^+_{\mathscr{X}}=\Gamma(\mathscr{X},\O^{\int}_{\mathscr{X}})
$$
is a complete extremal affinoid ring\index{affinoid!affinoid ring@--- ring!complete affinoid ring@complete --- ---} such that $\mathcal{A}^{\pm}_{\mathscr{X}}$ is t.u.\ rigid-Noetherian.

{\rm (2)} Any distinguished affine formal model\index{formal model!formal model of a coherent rigid space@--- (of a coherent rigid space)!distinguished formal model of a coherent rigid space@distinguished --- ---} $X=\Spf A$ of $\mathscr{X}$ induces a strict equivalence\index{strictly equivalent} class $C(A)$ of $\mathcal{A}^{\pm}_{\mathscr{X}}$ defined by $A$, which only depends on $\mathscr{X};$ moreover, if we denote by $\mathcal{C}_{\mathscr{X}}=C(A)$ the rigidification thus obtained, then 
$$
\Aff(\mathscr{X})=(\mathcal{A}_{\mathscr{X}},\mathcal{C}_{\mathscr{X}})
$$
is an a-r-pair\index{arpair@a-r-pair}\index{pair!arpair@a-r-pair}. \hfill$\square$
\end{thm}

\begin{dfn}\label{dfn-rigidanalyticring}{\rm 
(1) An affinoid ring isomorphic to $\mathcal{A}_{\mathscr{X}}$ as in \ref{thm-steinaffvsadic1} (1) is called an {\em analytic affinoid ring}\index{affinoid!affinoid ring@--- ring!analytic affinoid ring@analytic --- ---} associated to $\mathscr{X}$.

(2) An a-r-pair isomorphic to $\Aff(\mathscr{X})$ as defined in \ref{thm-steinaffvsadic1} (2) is called a {\em analytic affinoid pair}\index{affinoid!analytic affinoid pair@analytic --- pair} associated to $\mathscr{X}$.}
\end{dfn}

It is clear that, for any morphism $\varphi\colon\mathscr{X}\rightarrow\mathscr{Y}$ of Stein affinoids, we have canonically an induced homomorphism $\Aff(\varphi)\colon\Aff(\mathscr{Y})\rightarrow\Aff(\mathscr{X})$ of a-r-pairs.
Thus $\Aff$ gives rise to a functor from the opposite category of the category of Stein affinoids to the category of a-r-pairs.
\begin{thm}\label{thm-steinaffvsadic3}
The functor $\Aff$ thus obtained is fully faithful with the essential image being the totality of analytic affinoid pairs\index{affinoid!analytic affinoid pair@analytic --- pair}.
\end{thm}

\begin{proof}
Let $\mathscr{X}$ and $\mathscr{Y}$ be Stein affinoids, and consider the map
$$
\Hom_{\CRf}(\mathscr{X},\mathscr{Y})\longrightarrow\Hom(\Aff(\mathscr{Y}),\Aff(\mathscr{X})),
$$
where the last set is the set of homomorphisms of a-r-pairs.
It is not difficult to see that this map is injective (left to the reader).
To show the surjectivity, let us be given a homomorphism $\psi\colon\Aff(\mathscr{Y})\rightarrow\Aff(\mathscr{X})$ of a-r-pairs.
Let $\Spf A$ and $\Spf B$ be distinguished affine formal models of $\mathscr{X}$ and $\mathscr{Y}$, respectively.
On the other hand, let $A_0$ and $B_0$ be rings of definition in $\mathcal{C}_{\mathscr{X}}$ and $\mathcal{C}_{\mathscr{Y}}$, respectively, such that $\psi(B_0)\subseteq A_0$ and that the map $B_0\rightarrow A_0$ is adic.
One can replace $B_0$ and $A_0$ by strict weak equivalence and thus may assume $B=B_0$ and $A=A_0$ (cf.\ \ref{prop-lemmorbetweenaffinoid2} and \ref{prop-stronglynoetherianrigidnotherian}).
Hence we get the map $\varphi\colon\mathscr{X}\rightarrow\mathscr{Y}$ induced by the adic map $B\rightarrow A$.
It is straightforward to check that $\Aff(\varphi)=\psi$.
\end{proof}

\subsubsection{Visualization and adic spectrum}\label{subsub-visualizationadicspace}
\begin{thm}\label{thm-steinaffvsadic2}
Let $\mathscr{X}$ be a Stein affinoid\index{Stein affinoid}\index{affinoid!Stein affinoid@Stein ---}, and consider the analytic affinoid pair\index{affinoid!analytic affinoid pair@analytic --- pair} $\Aff(\mathscr{X})=(\mathcal{A}_{\mathscr{X}},\mathcal{C}_{\mathscr{X}})$ associated to $\mathscr{X}$.
Then the adic spectrum\index{adic spectrum} $X=\Spa\mathcal{A}_{\mathscr{X}}$ is canonically homeomorphic to $\ZR{\mathscr{X}};$ moreover, by this identification, the sheaf $\O^{\int}_{\mathscr{X}}$ $($resp.\ $\O_{\mathscr{X}})$ coincides with the presheaf $\O^+_X$ $($resp.\ $\O_X)$; see {\rm \S\ref{subsub-adicspectrum}}.
\end{thm}

\begin{proof}
Set $B=\Gamma(\mathscr{X},\O_{\mathscr{X}})$.
To construct the map $\ZR{\mathscr{X}}\rightarrow X=\Spa\mathcal{A}_{\mathscr{X}}$, take $x\in\ZR{\mathscr{X}}$.
Then one has the associated rigid point $\alpha_x\colon\Spf\widehat{V}_x\rightarrow\ZR{\mathscr{X}}$ (\ref{dfn-ZRpoints32}).
Let $\Spf A$ be a distinguished formal affine model of $\mathscr{X}$, and $I$ a finitely generated ideal of definition of $A$.
Then the rigid point induces an adic map $A\rightarrow\widehat{V}_x$.
Since $\widehat{V}_x$ is $a$-adically separated where $I\widehat{V}_x=(a)$, and since $B=\varinjlim_{n\geq 1}\Hom(I^n,A)$, this adic map gives rise to a continuous valuation $v_x\colon B\rightarrow\Gamma_x\cup\{\infty\}$, where $\Gamma_x$ is the value group of the valuation ring $\widehat{V}_x$, whence the desired map $x\mapsto v_x$.
It is easy to see that this map is continuous.

To construct the inverse map $X=\Spa\mathcal{A}_{\mathscr{X}}\rightarrow\ZR{\mathscr{X}}$, take a continuous valuation $v\colon B\rightarrow\Gamma\cup\{\infty\}$, and consider the composition $v\colon A\rightarrow\Gamma\cup\{\infty\}$.
The kernel $v^{-1}(0)$ is a prime ideal, and thus the map factors through a field; taking the inverse image of the subset consisting of positive values, one gets a valuation ring $V$ that factors through $v$.
Set $IV=(a)$; we have $a\neq 0$ and $a\not\in V^{\times}$.
Then $V$ is $a$-adically separated due to the continuity of $v$.
By {\bf \ref{ch-pre}}.\ref{thm-compval2006ver1} the $a$-adic completion $\widehat{V}$ is an $a$-adically complete valuation ring having the same value group as $V$, and thus we get the rigid point $\alpha\colon\Spf\widehat{V}\rightarrow\Spf A$, whence $\alpha\colon\Spf\widehat{V}\rightarrow\ZR{\mathscr{X}}$ (cf.\ \ref{prop-ZRpoints3}).
It is straightforward to see that the maps thus obtained are continuous and the inverse maps to each other, and thus the first assertion of \ref{thm-steinaffvsadic2} is shown.

Next we compare the structural presheaves.
To this end, since the issue is local, we may replace $A$ by an affine part of the admissible blow-up along the ideal of definition $I$, and thus we may assume that $I=(a)$ is invertible.
Take $f_0,\ldots,f_n\in B$ such that the ideal $(f_0,\ldots,f_n)$ in $B$ is the unit ideal (hence is open).
In this situation we consider the rational subdomains
$$
{\textstyle \mathscr{U}_k=\mathscr{X}(\frac{f_0}{f_k},\ldots,\frac{f_{k-1}}{f_k},\frac{f_{k+1}}{f_k},\ldots,\frac{f_n}{f_k})}
$$
for $k=0,\ldots,n$ (\ref{exas-affinoidsubdomain} (3)), which gives a covering $\ZR{\mathscr{X}}=\bigcup^n_{k=0}\ZR{\mathscr{U}_k}$ by quasi-compact open subsets.
Notice that the rational subdomains $\ZR{\mathscr{U}_k}$ correspond to the rational subsets $U_k=\{v\in\Spa\mathcal{A}_{\mathscr{X}}\,|\,v(f_i)\geq v(f_k)\neq\infty,\ i\neq k\}$ of $X=\Spa\mathcal{A}_{\mathscr{X}}$ by the map constructed above.
Since this covering comes from the admissible blow-up along the admissible ideal $(a^lf_0,\ldots,a^lf_n)$ (with $l\geq 0$ sufficiently large), the open coverings of this form are cofinal, and hence it is enough to show that the presheaves take the same values on the subsets $\ZR{\mathscr{U}_k}$.

For $k=0,\ldots,n$ we have
$$
\Gamma(U_k,\O_X)=B\dl{\textstyle \frac{f_0,\ldots,f_{k-1},f_{k+1},\ldots,f_n}{f_k}}\dr\leqno{(\ast)}
$$
by definition (\S\ref{subsub-adicspectrum}).
On the other hand, we have
$$
\Gamma(\ZR{\mathscr{U}_k},\O_{\mathscr{X}})=A\dl{\textstyle\frac{a^lf_0}{a^lf_k},\ldots,\frac{a^lf_{k-1}}{a^lf_k},\frac{a^lf_{k+1}}{a^lf_k},\ldots,\frac{a^lf_n}{a^lf_k}}\dr\otimes_AB\leqno{(\ast\ast)}
$$
by \ref{thm-cohomologyrigidsp1str} (1), which is nothing but the right-hand side of $(\ast)$.
Hence we have $\O_X=\O_{\mathscr{X}}$.
Since $\Gamma(U_k,\O^+_X)$ consists of elements in $\Gamma(U_k,\O_X)$ that have non-negative values at any $v\in U_k$, it coincides with $\Gamma(\ZR{\mathscr{U}_k},\O^{\int}_{\mathscr{X}})$, since in view of {\bf \ref{ch-pre}}.\ref{thm-valuative} the inequality $v_x(f_x)\geq 0$ for $f\in\Gamma(\ZR{\mathscr{U}_k},\O_{\mathscr{X}})$ holds if and only if $f_x\in\O^{\int}_{\mathscr{X},x}$.
\end{proof}

\begin{cor}\label{cor-steinaffvsadic2}
Let $\mathcal{A}_{\mathscr{X}}$ be the analytic affinoid ring\index{affinoid!affinoid ring@--- ring!analytic affinoid ring@analytic --- ---} associated to a Stein affinoid $\mathscr{X}$, and consider $X=\Spa\mathcal{A}_{\mathscr{X}}$. Then the presheaves $\O_X$ and $\O^+_X$ are sheaves, thereby defining the affinoid adic space\index{adic space!affinoid adic space@affinoid ---} {\rm (\ref{dfn-affinoidadicspace})} $((X,\O^+_X,\O_X),\{v_x\}_{x\in X})$, which is, as a valued triple\index{triple!valued triple@valued ---}, canonically isomorphic to the Zariski-Riemann triple\index{triple!Zariski Riemann triple@Zariski-Riemann ---} $\ZRT(\mathscr{X})$ associated to $\mathscr{X}$.
\end{cor}

Note that the Zariski-Riemann triple $\ZRT(\mathscr{X})$ is an analytic triple (\ref{thm-fullyfaithfultriples}), and hence is naturally regarded as a valued triple.

\begin{proof}
The first assertion follows immediately from \ref{thm-steinaffvsadic2}.
One can check the second assertion by the calculation as in the proof of \ref{thm-steinaffvsadic2}.
\end{proof}

\subsubsection{Description of power-bounded elements}\label{subsub-descpbelements}
Let $\mathscr{X}=(\Spf A)^{\rig}$ be a Stein affinoid\index{affinoid!Stein affinoid@Stein ---}, where $A$ is an $I$-torsion free (where $I\subseteq A$ is a finitely generated ideal of definition) t.u.\ rigid-Noetherian ring\index{t.u. rigid-Noetherian ring@t.u.\ rigid-Noetherian} such that $\Spec A\setminus V(I)$ is affine.
In this paragraph, we are interested in the subring $\Gamma(\mathscr{X},\O_{\mathscr{X}})^o$ of power-bounded\index{power-bounded} elements ({\bf \ref{ch-pre}}, \S\ref{subsub-fadicringsgeneralities}) in the f-adic ring $\Gamma(\mathscr{X},\O_{\mathscr{X}})$.
We first remark that the following inclusion holds:
$$
\Gamma(\mathscr{X},\O_{\mathscr{X}})^o\subseteq\{f\in\Gamma(\mathscr{X},\O_{\mathscr{X}})\,|\,f_x\in\O^{\int}_{\mathscr{X},x}\ \textrm{for any}\ x\in[\mathscr{X}]\}.
$$
Indeed, if $x\in[\mathscr{X}]$, then the valuation at $x$ (cf.\ \ref{ntn-ZRpoints}) is of height one, and may take values in $\R$.
If $f\in\Gamma(\mathscr{X},\O_{\mathscr{X}})$ is power-bounded, then the value of $f_x$ is non-negative, which implies that $f_x\in\O^{\int}_{\mathscr{X},x}$.

\begin{prop}\label{prop-cordescpbelements1}
Suppose that there exists a real valued spectral functor\index{spectral functor!real valued spectral functor@real valued ---} {\rm (\ref{dfn-realvalued})} defined on the category of quasi-compact open subspaces of $\mathscr{X}$.
Then we have the equality 
$$
\Gamma(\mathscr{X},\O_{\mathscr{X}})^o=\Gamma(\mathscr{X},\O^{\int}_{\mathscr{X}}).
$$
Moreover, $\Gamma(\mathscr{X},\O_{\mathscr{X}})^o$ coincides with the integral closure of $A$ in $\Gamma(\mathscr{X},\O_{\mathscr{X}})$.
\end{prop}

\begin{proof}
The first part of the assertion follows from the above remark and \ref{thm-densityargument}.
The second part follows from \ref{prop-corsteinaffvsadic1corint}.
\end{proof}

By \ref{thm-spectralfunctorseparatedquotients} we have:
\begin{cor}\label{cor-cordescpbelements1}
Let $\mathscr{X}$ be a Stein affinoid of type {\rm ($\mathrm{V_{\R}}$)}\index{rigid space!rigid space of typeV1@--- of type ($\mathrm{V_{\R}}$)} or of type {\rm (N)}\index{rigid space!rigid space of typeN@--- of type (N)}.
Then
$$
\Gamma(\mathscr{X},\O_{\mathscr{X}})^o=\Gamma(\mathscr{X},\O^{\int}_{\mathscr{X}}).
$$
Moreover, if $X=\Spf A$ is a distinguished affine formal model of $\mathscr{X}$, where $A$ is a topologically of finite type algebra over an $a$-adically complete height one valuation ring in type {\rm ($\mathrm{V_{\R}}$)} case, or, in type {\rm (N)} case, an $I$-adically complete Noetherian ring $($where $I\subseteq A$ is a finitely generated ideal of definition$)$ such that $X\setminus V(I)$ is affine, then $\Gamma(\mathscr{X},\O_{\mathscr{X}})^o$ is the integral closure of $A=\Gamma(X,\O_X)$ in $\Gamma(\mathscr{X},\O_{\mathscr{X}})$.
\end{cor}

\subsubsection{Rigidification and finite type extensions}\label{subsub-rigidificationfinitetypeext}
\begin{thm}\label{thm-rigidificationfinitetypeext}
Let $A$ be a t.u.\ rigid-Noetherian extremal f-adic ring, and $C$ a rigidification\index{rigidification} {\rm (\ref{dfn-rigidificationfadic})} of $A$.
Let $\mathcal{A}=(\mathcal{A}^{\pm},\mathcal{A}^+)$ be the associated affinoid ring {\rm (\ref{prop-affinoidrigidfadicrings})}.
Then the following categories are canonically equivalent to each other$:$
\begin{itemize}
\item[{\rm (a)}] the category of $A$-algebras of finite type $($cf.\ {\rm \S\ref{subsub-finitetypeextensions}}$)$ with rigid {\rm (\ref{dfn-rigidificationfadicmorphism})} $A$-algebra homomorphisms $($notice that, as we have seen in {\rm \S\ref{subsub-rigidificationfadicrings}}, any finite type extension of $A$ has the canonical rigidification induced from $C);$
\item[{\rm (b)}] the category of finite type affinoid rings {\rm (\S\ref{subsub-analyticadicspaces})} over $\mathcal{A};$
\item[{\rm (c)}] the category of finite type affinoid rings over $\mathcal{A}$ with {\em discrete} $A$-algebra homomorphisms that preserve the positive parts.
\end{itemize}
\end{thm}

We will use the following lemma for the proof of the theorem:
\begin{lem}\label{lem-rigidificationfinitetypeext}
Let $A$ be as in {\rm \ref{thm-rigidificationfinitetypeext}}, and $B$ an $A$-algebra of finite type.
Let $A_0\subseteq A$ and $B_0\subseteq B$ be respective rings of definition, and $I_0\subseteq A_0$ a finitely generated ideal of definition.
Consider elements $f_1,\ldots,f_n\in B^{\int}_0$ in the integral closure of $B_0$ in $B$.
Then the subring $B_1$ of $B$ weakly generated over $B_0$ by $f_1,\ldots,f_n$ {\rm ({\bf \ref{ch-pre}}, \S\ref{subsub-completeaffinoidrings})} is a ring of definition of $B$.
Moreover, $B_1$ is generated over $B_0$ by $f_1,\ldots,f_n$, and defines the same rigidification as $B_0$.
\end{lem}

\begin{proof}
Since $f_1,\ldots,f_n\in B^{\int}_0$, the ring $B_0[f_1,\ldots,f_n]$ is finite over $B_0$ and hence is $I_0$-adically complete.
Hence we have $B_1=B_0[f_1,\ldots,f_n]$, which is a ring of definition.
It is clear that $B_1$ is generated by $f_1,\ldots,f_n$ and $B_0$.
Since $B_1$ is finite over $B_0$, they define the same rigidification by \ref{prop-lemmorbetweenaffnoid1}.
\end{proof}

\begin{proof}[Proof of Theorem {\rm \ref{thm-rigidificationfinitetypeext}}]
It is clear that there exist the canonical functors, one from the category (a) to the category (b) and another from the category (b) to the category (c).
We are going to show that there exists the canonical functor from (c) to (a), which gives a quasi-inverse to the composite functor from category (a) to category (c).
Let $\mathcal{B}=(\mathcal{B}^{\pm},\mathcal{B}^+)$ and $\mathcal{B}'=(\mathcal{B}^{\prime\pm},\mathcal{B}^{\prime +})$ be affinoid rings of finite type over $\mathcal{A}$, and $\varphi\colon\mathcal{B}^{\pm}\rightarrow\mathcal{B}^{\prime\pm}$ an $A$-algebra homomorphism such that $\varphi(\mathcal{B}^+)\subseteq\mathcal{B}^{\prime +}$.
Replacing $A_0$ in the same rigidification class if necessary, one can take topologically finitely generated rings of definition $B_0\subseteq\mathcal{B}^+$ and $B'_0\subseteq\mathcal{B}^{\prime +}$ over $A_0$.
Then replace $B'_0$ by the one generated by $B'_0$ and $\varphi(B_0)$, which is, by \ref{lem-rigidificationfinitetypeext}, again a ring of definition that gives rise to the same rigidification.
Since $\varphi(B_0)\subseteq B'_0$, it defines a morphism in the category (a).
\end{proof}

\begin{cor}\label{cor-rigidificationfinitetypeext}
Let $A$ be a t.u.\ rigid-Noetherian extremal f-adic ring, and $C$ a rigidification\index{rigidification} of $A$.
Let $\mathcal{A}=(\mathcal{A}^{\pm},\mathcal{A}^+)$ be the associated affinoid ring.
Suppose that there exists a real valued spectral functor\index{spectral functor!real valued spectral functor@real valued ---} {\rm (\ref{dfn-realvalued})} defined on the category of all rigid spaces of finite type over the affinoid associated to $(A,C)$ {\rm (\S\ref{subsub-rigidificationfadicrings})}.
Then the following categories are equivalent to each other$:$
\begin{itemize}
\item[{\rm (a)}] the category of affinoid rings of finite type over $\mathcal{A};$
\item[{\rm (b)}] the category of $A$-algebras of finite type with continuous $A$-algebra homomorphisms.
\end{itemize}
\end{cor}

\begin{proof}
By \ref{thm-rigidificationfinitetypeext} it suffices to check that morphisms in the category (b) preserve the positive part.
This follows from \ref{prop-cordescpbelements1}.
\end{proof}

The significance of the corollary lies in that the latter category does not refer to the rigidifications; in other words, the rigidifications in this situation are {\em canonical}.
Notice that by \ref{thm-spectralfunctorseparatedquotients} the hypothesis in \ref{cor-rigidificationfinitetypeext} is automatic in type {\rm ($\mathrm{V_{\R}}$)}\index{rigid space!rigid space of typeV1@--- of type ($\mathrm{V_{\R}}$)} and type {\rm (N)}\index{rigid space!rigid space of typeN@--- of type (N)} cases.

\subsubsection{Analytic rings of type (N)}\label{subsub-noetheriananalyticrings}
\begin{dfn}\label{dfn-noetheriananalyticring}{\rm 
An f-adic ring $A$ is called an {\em analytic ring of type {\rm (N)}}\index{analytic ring!analytic ring of typeN@--- of type (N)} if it is complete, extremal ({\bf \ref{ch-pre}}, \S\ref{subsub-extremalfadicrings}), and admits a Noetherian ring of definition.}
\end{dfn}

\begin{prop}\label{prop-noetheriananalyticring2}
Let $A$ be an analytic ring of type {\rm (N)}\index{analytic ring!analytic ring of typeN@--- of type (N)}, and $A_0$ a Noetherian ring of definition of $A$. 
For any ring of definition $R$ of $A$ contained in $A^o$, the subring generated by $A_0$ and $R$ in $A^o$ is finite over $A_0$. 
\end{prop}

\begin{proof} 
Let $J$ be an ideal of definition of $R$. 
Replacing $J$ by a power of $J$ if necessary, we may assume that $J\subset A_0$ and that $JA_0$ is an ideal of definition of $A_0$.
The subring $R$ is a subset of $M_0=\{x\in A^{\int}_0\,|\, Jx\subseteq A_0\}$, which is an $A_0$-submodule of $A^{\int}_0$.
Let $X$ be the blow-up of $\Spec A$ along the ideal $JA$. 
We have $J\O_X=\O_X(1)$. 
Then $M_0$ is an $A_0$-submodule of $\Gamma(X,\O(-1))$, and the last module is a finitely generated $A_0$-module. 
Since $A_0$ is Noetherian, $M_0$ is a finitely generated $A_0$-module, and hence the ring generated by $A_0$ and $R$ is finite over $A_0$. 
\end{proof}

\begin{cor}\label{cor-noetheriananalyticring2}
Let $A$ be an analytic ring of type {\rm (N)}.
Then the strict equivalence\index{strictly equivalent} class containing a Noetherian ring of definition is unique. \hfill$\square$
\end{cor}

For an analytic ring $A$ of type (N) we denote by $C_0(A)$ the finite equivalence class (\ref{dfn-strictlyequiv}) containing a Noetherian ring of definition.

\begin{prop}\label{prop-noetheriananalyticring3}
Let $A$ be an analytic ring $A$ of type {\em (N)}.
Then every ring of definition of $A$ belonging to $C_0(A)$ is Noetherian.
\end{prop}

\begin{proof}
The assertion follows immediately from Eakin-Nagata theorem, which guarantees that for a finite ring extension $A\subseteq B$ of rings with $B$ Noetherian, $A$ is Noetherian.
\end{proof}

\begin{cor}\label{cor-noetheriananalyticring3}
Let $\mathscr{X}$ be a universally adhesive affinoid.  
If $\mathscr{X}$ has at least one Noetherian affine formal model, any distinguished affine formal model $X$ of $\mathscr{X}$ is Noetherian. 
\end{cor}

\begin{proof}
This follows from \ref{prop-corlemtuavsrigidaff}, \ref{thm-morbetweenaffinoid1}, and \ref{prop-noetheriananalyticring3}.
\end{proof}

\begin{prop}\label{prop-noetheriananalyticring4}
Let $A$ and $A'$ be analytic rings of type {\rm (N)}. 
Suppose we are given rigidifications\index{rigidification} induced from $C_0(A)$ and $C_0(A')$ of $A$ and $A'$, respectively.
Then any homomorphism of affinoid rings $\varphi\colon(A,A^o)\rightarrow(A',{A'}^o)$ is a homomorphism of a-r-pairs {\rm (\ref{dfn-rigidification})}. 
\end{prop}

\begin{proof}
By {\bf \ref{ch-pre}}.\ref{prop-lempowerboundedfadic} (2) and {\bf \ref{ch-pre}}.\ref{cor-typeRfadicrings} we know that $\varphi$ is adic.
Let $A_0$ and $A'_0$ be Noetherian rings of definition of $A$ and $A'$, respectively.
Let $A''_0$ be the subring of $A'$ generated by $A'_0$ and $\varphi(A_0)$.
We claim that $A''_0$ is finite over $A'_0$.
Indeed, since $\varphi(A_0)$ is bounded, there exists $n\geq 0$ such that $I^{\prime n}_0\varphi(A_0)\subseteq A'_0$, where $I'_0\subseteq A'_0$ is an ideal of definition.
This means that $A''_0\subseteq[A'_0:I^{\prime n}_0]$.
But since $[A'_0:I^{\prime n}_0]$ is obviously finite over $A'_0$ (since $\Spec A'=\Spec A'_0\setminus V(I'_0)$), we see that $A''_0$ is finite.
\end{proof}

By what we have seen in this paragraph we have (cf.\ \ref{cor-rigidificationfinitetypeext} and \ref{thm-spectralfunctorseparatedquotients}):
\begin{thm}\label{thm-rigidificationnoetherian}
The following categories are canonically equivalent to each other$:$
\begin{itemize}
\item[{\rm (a)}] the category of a-r-pairs by analytic rings of type {\rm (N)} with homomorphisms of a-r-pairs$;$
\item[{\rm (b)}] the category of affinoid rings given by analytic rings of type {\rm (N)} $($without rigidifications$);$
\item[{\rm (c)}] the category of analytic rings of type {\rm (N)} with continuous ring homomorphisms. \hfill$\square$
\end{itemize}
\end{thm}

\subsubsection{Canonical rigidifications of classical affinoid algebras}\label{subsub-canonicalrigidificationclassicalaffinoidalg}
Let $V$ be an $a$-adically complete valuation ring of height one $(a\in\m_V\setminus\{0\})$, and $K=\Frac(V)$.
Recall that any classical affinoid algebra\index{algebra!affinoid algebra@affinoid ---!classical affinoid algebra@classical --- ---} $\mathcal{A}$ over $K$ ({\bf \ref{ch-pre}}.\ref{dfn-classicalaffinoidalgebras}) is a complete Tate ring\index{Tate ring} ({\bf \ref{ch-pre}}, \S\ref{subsub-extremalfadicrings}).
For any subring $A\subseteq\mathcal{A}$ topologically of finite type over $V$ such that $A[\frac{1}{a}]=\mathcal{A}$, the rigid space $\mathscr{X}=(\Spf A)^{\rig}$ is a Stein affinoid, and hence by \ref{thm-steinaffvsadic1} (1) and \ref{cor-cordescpbelements1}, the pair $(\mathcal{A},\mathcal{A}^o)$ is an affinoid ring\index{affinoid!affinoid ring@--- ring} (\ref{dfn-affinoidringadicspace}).

\begin{dfn}[Canonical subring]\label{dfn-canonicalsubring}{\rm 
For a classical affinoid algebra $\mathcal{A}$ over $K$ we define
$$
\mathcal{A}^{\can}=\{f\in\mathcal{A}\,|\,\textrm{$(f$ mod $\m)$ is integral over $V$ for any maximal ideal $\m\subseteq\mathcal{A}$}\},
$$
and call it the {\em canonical subring} of $\mathcal{A}$.}
\end{dfn}

Notice that the subring $\mathcal{A}^{\can}$ only depends on the $K$-algebra structure and does not depend on the topological (f-adic) structure of $\mathcal{A}$.
\begin{prop}[Uniqueness of affinoid ring structure]\label{prop-canonicalsubring}
Let $\mathcal{A}$ be a classical affinoid algebra over $K$, and choose a ring of definition $A\subseteq\mathcal{A}$ that is topologically finitely generated over $V$.
Let $\mathscr{X}=(\Spf A)^{\rig}$ be the associated Stein affinoid.
Then we have the equalities
$$
\mathcal{A}^{\can}=\Gamma(\mathscr{X},\O_{\mathscr{X}})^o\ (=\mathcal{A}^o)=\Gamma(\mathscr{X},\O^{\int}_{\mathscr{X}}).
$$
\end{prop}

The significance of the proposition lies in that the affinoid ring structure $(\mathcal{A},\mathcal{A}^o)$ actually does not depend on the topology on $\mathcal{A}$; combined with the previously obtained results in \S\ref{subsub-rigidificationfinitetypeext}, such an $\mathcal{A}$ also possesses the canonical rigidification, that is to say, the assignment $\mathcal{A}\leadsto\mathscr{X}$ is canonical, only depending on the $K$-algebra structure of $\mathcal{A}$.

\begin{proof}
For any $x\in\ZR{\mathscr{X}}^{\cl}$, denote by $\m_x\subseteq\mathcal{A}$ the corresponding maximal ideal (\ref{cor-classicalpointsexist0}).
Let $f\in\mathcal{A}=\Gamma(\mathscr{X},\O_{\mathscr{X}})$.
If $f\in\mathcal{A}^{\can}$, then the image of $f$ in $\O_{\mathscr{X},x}/\m_{\mathscr{X},x}$ is finite over $V_x=\O^{\int}/\m_{\mathscr{X},x}$, the valuation ring at $x$, and hence belongs to $V_x$ (here we partially used the notation as in \ref{ntn-ZRpoints}).
Hence the germ $f_x\in\O_{\mathscr{X},x}$ lies in $\O^{\int}_{\mathscr{X},x}$.
Now by \ref{thm-densityargument} and \ref{cor-spectralfunctorseparatedquotients2} we have $f\in\Gamma(\mathscr{X},\O^{\int}_{\mathscr{X}})$.
This shows the inclusion $\mathcal{A}^{\can}\subseteq\Gamma(\mathscr{X},\O^{\int}_{\mathscr{X}})$, while the opposite inclusion is easy to see.
The other equality follows from \ref{cor-cordescpbelements1}.
\end{proof}

\begin{cor}\label{cor-canonicalsubring}
For any $K$-algebra homomorphism $\varphi\colon\mathcal{A}\rightarrow\mathcal{A}'$ between classical affinoid algebras over $K$, we have $\varphi(\mathcal{A}^o)\subseteq\mathcal{A}^{\prime o}$.
In particular, $\varphi$ induces a morphism of affinoid rings $(\mathcal{A},\mathcal{A}^o)\rightarrow (\mathcal{A}',\mathcal{A}^{\prime o})$.
\end{cor}

Thus one deduces, in particular, that any $K$-algebra homomorphism between classical affinoid algebras is continuous, the result that is considered to be one of the most fundamental ones in classical rigid analytic geometry.
\begin{proof}
In view of \ref{cor-classicalpointsexist0} and \ref{prop-classicalpoint1}, the inclusion $\varphi(\mathcal{A}^o)\subseteq\mathcal{A}^{\prime o}$ follows immediately from \ref{prop-canonicalsubring}.
For the second, what to do is to check that the topologies on $\mathcal{A}$ and $\mathcal{A}'$ are uniquely determined and that a ring of definition of finite type over $V$ is mapped inside a ring of definition of finite type over $V$.
This follows from \ref{thm-rigidificationfinitetypeext}.
\end{proof}

\begin{cor}\label{cor-banachfadicringstate3}
In the situation as in {\rm \ref{cor-canonicalsubring}}, the morphism $\varphi\colon(\mathcal{A},\mathcal{A}^o)\rightarrow (\mathcal{A}',\mathcal{A}^{\prime o})$ is a homomorphism of a-r-pairs with respect to the canonical rigidifications $($cf.\ {\rm \S\ref{subsub-finitetypeextensions}}$)$. \hfill$\square$
\end{cor}

Compiling all results so far obtained, we have:
\begin{thm}\label{thm-canonicityclassicalaffinoid}
Let $V$ be an $a$-adically complete valuation ring of height one, and $K=\Frac(V)$ its fractional field.
Set $\mathscr{S}=(\Spf V)^{\rig}$.
Then the following categories are canonically equivalent to each other$:$
\begin{itemize}
\item[{\rm (a)}] the category of affinoids of finite type over $\mathscr{S}$ and morphisms over $\mathscr{S};$
\item[{\rm (b)}] the opposite category of the category of classical affinoid algebras over $K$ and $K$-algebra homomorphisms.\hfill$\square$
\end{itemize}
\end{thm}

\subsection{Rigid geometry and adic spaces}\label{sub-rigidgeomadicspaces}
\index{adic space|(}
We already know that the Zariski-Riemann triple\index{Zariski-Riemann triple}\index{triple!Zariski Riemann triple@Zariski-Riemann ---} $\ZRT(\mathscr{X})$ for a Stein affinoid\index{affinoid!Stein affinoid@Stein ---} $\mathscr{X}$ is canonically identified with the affinoid adic space\index{adic space!affinoid adic space@affinoid ---} 
$$
((X=\Spa\mathcal{A}_{\mathscr{X}},\O^+_X,\O_X),\{v_x\}_{x\in X})
$$
associated to the analytic affinoid pair\index{affinoid!analytic affinoid pair@analytic --- pair} $\Aff(\mathscr{X})= (\mathcal{A}_{\mathscr{X}},\mathcal{C}_{\mathscr{X}})$ (\ref{cor-steinaffvsadic2}).
In particular, we know that in this situation the presheaves $\O^+_X$ and $\O_X$ are sheaves.
Therefore, by the definition of adic spaces\index{adic space} (\ref{dfn-adicspacesbytriples}) we have:
\begin{thm}\label{thm-rigidgeomadicspaces}
{\rm (1)} Let $\mathscr{X}$ be a locally universally Noetherian rigid space\index{rigid space!universally Noetherian rigid space@universally Noetherian ---!locally universally Noetherian rigid space@locally --- ---}. Then the associated valued triple\index{triple!valued triple@valued ---} $\ZRT(\mathscr{X})=(\ZRT(\mathscr{X}),\{v_x\}_{x\in X})$ is an analytic adic space\index{adic space!analytic adic space@analytic ---} {\rm (\ref{dfn-analyticadicspaces})}.

{\rm (2)} Let $\varphi\colon\mathscr{X}\rightarrow\mathscr{Y}$ be a morphism between locally universally Noetherian rigid spaces. Then the morphism $\ZRT(\varphi)\colon\ZRT(\mathscr{X})\rightarrow\ZRT(\mathscr{Y})$ of valued triples is a morphism of adic spaces.
\end{thm}

\begin{proof}
(1) follows immediately from \ref{cor-steinaffvsadic2}.
Since morphisms of adic spaces are, by definition, morphisms of valued triples, (2) follows.
\end{proof}

Thus we have the functor
$$
\ZRT\colon\mathbf{RigNoeRf}\longrightarrow\AnAcsp
$$
from the category of locally universally Noetherian rigid spaces\index{rigid space!universally Noetherian rigid space@universally Noetherian ---!locally universally Noetherian rigid space@locally --- ---} to the category of analytic adic spaces\index{adic space!analytic adic space@analytic ---}, which is faithful due to \ref{thm-fullyfaithfultriples} (2).

\begin{thm}\label{thm-rigidgeomadicspacesff1}
Let $\mathscr{S}$ be a locally universally Noetherian rigid space. Then the visualization functor $\ZRT$ gives a categorical equivalence from the category of locally of finite type rigid spaces over $\mathscr{S}$ to the category of adic spaces locally of finite type {\rm (\ref{dfn-adicspacefinitetype})} over $\ZRT(\mathscr{S})$.
\end{thm}

\begin{proof}
In view of \ref{cor-rigidificationfinitetypeext}, it suffices to show that the functor in question is fully faithful, for, then, the essential surjectivity follows by a standard patching argument.
It is obvious that, for any locally of finite type $\mathscr{X}$ over $\mathscr{S}$, the adic space $\ZRT(\mathscr{X})$ is locally of finite type over $\ZRT(\mathscr{S})$.
We need to show that, for any $\mathscr{X}$ and $\mathscr{Y}$ locally of finite type over $\mathscr{S}$, the map 
$$
\Hom_{\mathscr{S}}(\mathscr{X},\mathscr{Y})\longrightarrow\Hom_{\ZRT(\mathscr{S})}(\ZRT(\mathscr{X}),\ZRT(\mathscr{Y})),
$$
already known to be injective, is bijective.
We may assume that $\mathscr{S}$, $\mathscr{X}$, and $\mathscr{Y}$ are Stein affinoids.
Since the functor $\mathcal{A}\mathcal{R}$ is fully faithful (\ref{thm-steinaffvsadic3}), we only need to check that any morphism $\varphi\colon\ZRT(\mathscr{X})\rightarrow\ZRT(\mathscr{Y})$ canonically gives a map $\mathcal{A}\mathcal{R}(\mathscr{Y})\rightarrow\mathcal{A}\mathcal{R}(\mathscr{X})$ of a-r-pairs (which is automatically adic).
Since $\mathcal{A}_{\mathscr{X}}=(\Gamma(\mathscr{X},\O_{\mathscr{X}}),\Gamma(\mathscr{X},\O^{\int}_{\mathscr{X}}))$ etc., the morphism $\varphi$ induces a homomorphism $\mathcal{A}_{\mathscr{Y}}\rightarrow\mathcal{A}_{\mathscr{X}}$ of affinoid rings.
As we have seen in \S\ref{subsub-rigidificationfadicrings}, the rigidifications $\mathcal{C}_{\mathscr{X}}$ and $\mathcal{C}_{\mathscr{Y}}$ are canonically determined by the rigidification $\mathcal{C}_{\mathscr{S}}$ of $\mathcal{A}\mathcal{R}(\mathscr{S})$, and it is clear that, in view of the definition of topologically finite type affinoid rings (cf.\ \cite[\S3, pp.\ 534ff]{Hube2}), these rigidifications are preserved by $\ZRT(\mathscr{S})$-morphisms of adic spaces.
\end{proof}

\begin{thm}\label{thm-rigidgeomadicspacesff2}
The category of analytic spaces of type {\rm (N)} and the category of rigid spaces of type {\rm (N)}\index{rigid space!rigid space of typeN@--- of type (N)} $(${\rm \ref{dfn-typeN}}$)$ are equivalent. 
\end{thm}

Here by an {\em analytic space of type} {\rm (N)} we mean an adic space that is locally isomorphic to the adic spectrum of a complete affinoid ring $\mathcal{A}=(\mathcal{A}^{\pm},\mathcal{A}^+)$, where $\mathcal{A}^{\pm}$ is an analytic ring of type {\rm (N)}\index{analytic ring!analytic ring of typeN@--- of type (N)} (\ref{dfn-noetheriananalyticring}) and $\mathcal{A}^+=A^o$.
\begin{proof}
Due to \ref{prop-noetheriananalyticring4} one can show by an argument similar to that in the proof of \ref{thm-rigidgeomadicspacesff1} that the functor $\ZRT$ is fully faithful from the category of rigid spaces of type (N) to the category of analytic spaces of type {\rm (N)}.
It is clear that the functor in question is essentially surjective.
\end{proof}

\begin{rem}\label{rem-rigidgeomadicspacesff}{\rm 
Let $X$ be a coherent universally rigid-Noetherian formal scheme $X$ such that $\O_X$ is $\mathscr{I}$-torsion free where $\mathscr{I}$ is an ideal of definition.
We do not know the following: if the associated rigid space $X^{\rig}$ is of type (N), then is $X$ Noetherian as a formal scheme?
Notice that this is true if $X$ is universally adhesive.}
\end{rem}
\index{adic space|)}



\section{Appendix: Tate's rigid analytic geometry}\label{sec-berkovich}
\subsection{Admissibility}\label{sec-admissibility}
\subsubsection{Admissibility with respect to a spectral functor}\label{sub-admissibleopensubsets}
Throughout this paragraph we fix:
\begin{itemize}
\item an $\mathrm{O}$-stable\index{stable!Ostable@O-{---}} subcategory $\mathscr{C}$ of the category of rigid spaces (\ref{dfn-spectral1} (1)) consisting of quasi-separated rigid spaces,
\item a continuous spectral functor\index{functor!spectral functor@spectral ---}\index{spectral functor}\index{spectral functor!continuous spectral functor@continuous ---} $S$ (\ref{dfn-spectral2}, \ref{dfn-spectral3}) defined on $\mathscr{C}$.
\end{itemize}
We assume that the data $(\mathscr{C},S)$ satisfy the following weak left-exactness:
\begin{itemize}
\item[{\rm (a)}] for any arrow $\varphi\colon\mathscr{X}\rightarrow\mathscr{Y}$ in $\mathscr{C}$ and any open subspace $\mathscr{V}\subseteq\mathscr{Y}$ we have $S(\varphi^{-1}(\mathscr{V}))=S(\varphi)^{-1}(S(\mathscr{V}))$;
\item[{\rm (b)}] for any family of open subspaces $\{\mathscr{U}_{\alpha}\}_{\alpha\in L}$ of $\mathscr{X}$ belonging to $\mathscr{C}$ we have $S(\bigcup_{\alpha\in L}\mathscr{U}_{\alpha})=\bigcup_{\alpha\in L}S(\mathscr{U}_{\alpha})$.
\end{itemize}

\begin{dfn}\label{dfn-CSadmissibility}{\rm 
Let $\mathscr{X}\in\obj(\mathscr{C})$.
An open subspace $\mathscr{U}\subseteq\mathscr{X}$ is said to be {\em $(\mathscr{C},S)$-admissible} if the following condition is satisfied:
\begin{itemize}
\item for any arrow $\varphi\colon\mathscr{V}\rightarrow\mathscr{X}$ in $\mathscr{C}$ where $\mathscr{V}$ is coherent and $S(\varphi)(S(\mathscr{V}))\subseteq S(\mathscr{U})$, we have $\varphi(\mathscr{V})\subseteq\mathscr{U}$.
\end{itemize}}
\end{dfn}

It is clear that, if $\psi\colon\mathscr{Y}\rightarrow\mathscr{X}$ is a morphism in $\mathscr{C}$ and $\mathscr{U}\subseteq\mathscr{X}$ is a $(\mathscr{C},S)$-admissible open subspace, then $\psi^{-1}(\mathscr{U})$ is a $(\mathscr{C},S)$-admissible open subspace of $\mathscr{Y}$.
Notice that, for any $\mathscr{X}\in\obj(\mathscr{C})$, any quasi-compact open subspace $\mathscr{U}\subseteq\mathscr{X}$ is $(\mathscr{C},S)$-admissible; indeed, for $\varphi\colon\mathscr{V}\rightarrow\mathscr{X}$ as above, since $S(\varphi^{-1}(\mathscr{U}))=S(\varphi)^{-1}(S(\mathscr{U}))=S(\mathscr{V})$, we have $\varphi^{-1}(\mathscr{U})=\mathscr{V}$ due to \ref{rem-spectralfunctor2}.

\begin{prop}\label{prop-csadmissiblemaximal}
Let $\mathscr{X}$ be a rigid space in $\mathscr{C}$, and $\mathscr{U}$ a $(\mathscr{C},S)$-admissible open subspace of $\mathscr{X}$.
Then $\mathscr{U}$ is the maximal open subspace of $\mathscr{X}$ among open subspaces $\mathscr{V}\subseteq\mathscr{X}$ such that $S(\mathscr{V})\subseteq S(\mathscr{U})$.
\end{prop}

\begin{proof}
For any quasi-compact open subspace $\mathscr{V}\subseteq\mathscr{X}$ such that $S(\mathscr{V})\subseteq S(\mathscr{U})$, we have $\mathscr{V}\subseteq\mathscr{U}$.
By continuity of $S$ we deduce that $\mathscr{U}$ is the union of all such quasi-compact subspaces.
\end{proof}

\begin{rem}\label{rem-csadmissiblemaximal}{\rm 
One can similarly show the following: for $\mathscr{X}\in\obj(\mathscr{C})$, an open subspace $\mathscr{U}\subseteq\mathscr{X}$ is $(\mathscr{C},S)$-admissible if and only if for any $\varphi\colon\mathscr{V}\rightarrow\mathscr{X}$ in $\mathscr{C}$ where $\mathscr{V}$ is coherent, $\varphi^{-1}(\mathscr{U})$ is the maximal open subspace of $\mathscr{V}$ among open subspaces $\mathscr{W}\subseteq\mathscr{V}$ such that $S(\mathscr{W})=S(\varphi^{-1}(\mathscr{U}))$.}
\end{rem}

For $\mathscr{X}\in\obj(\mathscr{C})$ we denote by $\mathfrak{A}(\mathscr{X})$ the category of $(\mathscr{C},S)$-admissible open subspaces in $\mathscr{X}$ and open immersions.
\begin{cor}\label{cor-csadmissiblemaximal1}
The spectral functor $S$ is conservative on $\mathfrak{A}(\mathscr{X})$, that is, for $\mathscr{U},\mathscr{U}'\in\obj(\mathfrak{A}(\mathscr{X}))$, $S(\mathscr{U})=S(\mathscr{U}')$ implies $\mathscr{U}=\mathscr{U}'$.\hfill$\square$
\end{cor}

\begin{dfn}\label{dfn-CSadmissibilitycov}{\rm 
Let $\mathscr{X}$ be a rigid space in $\mathscr{C}$, and $\mathscr{U}\subseteq\mathscr{X}$ a $(\mathscr{C},S)$-admissible open subspace.
An open covering $\mathscr{U}=\bigcup_{\alpha\in L}\mathscr{U}_{\alpha}$ is said to be {\em $(\mathscr{C},S)$-admissible} if every $\mathscr{U}_{\alpha}$ $(\alpha\in L)$ is $(\mathscr{C},S)$-admissible.}
\end{dfn}

By \ref{cor-csadmissiblemaximal1} and the weak left-exactness of $S$ we have:
\begin{prop}\label{prop-scadmissiblegtopology}
By attaching to each $\mathscr{X}\in\obj(\mathscr{C})$ the category $\mathfrak{A}(\mathscr{X})$ together with the notion of $(\mathscr{C},S)$-admissible coverings, we have a Grothendieck topology on the category $\mathscr{C}$.
The associated topos is equivalent to the large admissible topos {\rm (\ref{dfn-admissiblesite3genlarge})} restricted to $\mathscr{C}$. \hfill$\square$
\end{prop}

As we remarked before, any quasi-compact open subspace is $(\mathscr{C},S)$-admissible.
Here is another example of $(\mathscr{C},S)$-admissible subspaces: 
\begin{prop}\label{prop-csadmissibletube}
Any tube open subset\index{tube!tube open subset@--- open subset} {\rm (\ref{dfn-tubes1})}  of $\mathscr{X}$ is $(\mathscr{C},S)$-admissible.
\end{prop}

\begin{proof}
Let $T=\ZR{\mathscr{X}}\setminus\ovl{\ZR{\mathscr{U}}}$ be a tube open subset\index{tube!tube open subset@--- open subset} (\ref{dfn-tubes1}), where $\mathscr{U}\subseteq\mathscr{X}$ is a retrocompact open subspace.
Since the inverse image of a tube open subset by a morphism of coherent rigid spaces is again a tube open subset (cf.\ Exercise \ref{exer-tubesubsetscoherentmaps}), it suffices to show, by \ref{rem-csadmissiblemaximal}, that $T$ is the maximal among all open subsets $\mathscr{W}\subseteq\mathscr{X}$ such that $S(\mathscr{W})=S(T)$.
For any quasi-compact open subspace $\mathscr{V}\subseteq\mathscr{X}$ such that $S(\mathscr{V})\subseteq S(T)$, we have $S(\mathscr{U}\cap\mathscr{V})=S(\mathscr{U})\cap S(\mathscr{V})=\emptyset$.
This shows that $\mathscr{U}\cap\mathscr{V}=\emptyset$ and hence that $\ZR{\mathscr{V}}\subseteq T$.
Then the assertion follows from continuity of $S$.
\end{proof}

\subsubsection{G-topology on a topological space}\label{subsub-gtopology}
Let us recall the definition of {\em G-topologies} on a topological space $X$ (cf.\ \cite[(9.1.1)]{BGR}).
Let $\Open(X)$ be the category of open subsets of $X$; the objects of $\Open(X)$ are the open subsets of $X$, and for two open subsets $U,V\subseteq X$ 
\begin{equation*}
\Hom_{\Open(X)}(U,V)=\begin{cases}\{\ast\}&\textrm{if $U\subseteq V$,}\\ \emptyset&\textrm{otherwise.}\end{cases}
\end{equation*}
A G-topology on $X$ is a Grothendieck topology\index{topology!Grothendieck topology@Grothendieck ---}\index{Grothendieck!Grothendieck topology@--- topology} $\tau$ on a full subcategory $\mathfrak{A}$ of $\Open(X)$ such that $U,V\in\obj(\mathfrak{A})$ implies $U\cap V\in\obj(\mathfrak{A})$; open subsets in $\mathfrak{A}$ (resp.\ coverings of $\tau$) are called {\em $\tau$-admissible open subsets} (resp.\ {\em $\tau$-admissible open covering}).

For example, the topology on $X$ itself gives rise to a G-topology by an obvious manner. 
We call this G-topology the {\em canonical G-topology}.

If $\tau$ is a G-topology on a topological space $X$ and $Y\subseteq X$ is a subset, then one can consider the {\em induced G-topology} on $Y$, denoted by $\tau|_Y$, as follows:
\begin{itemize}
\item $\tau|_Y$-admissible open subsets are the open subsets of $Y$ of the form $U\cap Y$ by a $\tau$- admissible open subset $U$;
\item $\tau|_Y$-admissible coverings of an $\tau|_Y$-admissible open subset $U\cap Y$ is the open covering of the form $\{Y\cap U_{\alpha}\}_{\alpha\in L}$ by a $\tau$-admissible covering $\{U_{\alpha}\}_{\alpha\in L}$ of $U$.
\end{itemize}

\danger{Let $X$ be a topological space, and $\tau$ a G-topology on $X$ such that the associated topos $\top(\tau)$ is equivalent to the canonical topos $\top(X)$.
Then it may be the case, for a subset $Y\subseteq X$, that the topoi $\top(\tau|_Y)$ and $\top(Y)$ are not equivalent.
For example, consider $X=\R$ with the $\tau$ by the open intervals with rational extremities; although $\tau$ gives the topos equivalent to the one by $\R$ itself, the topoi $\top(\tau|_{\Q})$ and $\top(\Q)$ are different from each other, for the former is equivalent to $\top(\tau)$ (and hence the rational line $\Q$ with this G-topology is connected (cf.\ \cite[(9.1.1)]{BGR})), whereas the topological space $\Q$ is, as well-known, totally disconnected.}

In general, for a topological space $X$, a subset $Y\subseteq X$, and a G-topology $\tau$ on $X$ of which the associated topos is equivalent to the canonical one, the topos $\top(\tau|_Y)$ is equivalent to $\top(X)$ if $\tau$ satisfies the following condition: if $U,V$ are $\tau$-admissible open subsets, then $Y\cap U=Y\cap V$ implies $U=V$.

\subsubsection{G-topology associated to a spectral functor}\label{sub-gtopologyspectralfunctor}
Now we return to the situation as in the beginning of \S\ref{sub-admissibleopensubsets}.
We introduce the notion of $(\mathscr{C},S)$-admissibility (and thus a G-topology) on each topological space $S(\mathscr{X})$ for $\mathscr{X}\in\obj(\mathscr{C})$:
\begin{dfn}\label{dfn-csadmissibilityS}{\rm 
For $\mathscr{X}\in\obj(\mathscr{C})$ a $($not necessarily open$)$ subset $U\subseteq S(\mathscr{X})$ is said to be {\em $(\mathscr{C},S)$-admissible} if the following conditions are satisfied$:$
\begin{itemize}
\item[{\rm (a)}] there exists an open subspace $\mathscr{U}\subseteq\mathscr{X}$ such that $S(\mathscr{U})=U$;
\item[{\rm (b)}] for any morphism $\varphi\colon\mathscr{V}\rightarrow\mathscr{X}$ in $\mathscr{C}$ where $\mathscr{V}$ is coherent and $S(\varphi)(S(\mathscr{V}))\subseteq U$, there exists a collection of quasi-compact open subspaces $\{\mathscr{U}_{\alpha}\}_{\alpha\in L}$ such that $U=\bigcup_{\alpha\in L}S(\mathscr{U}_{\alpha})$ and that the covering $\{S(\varphi^{-1}(\mathscr{U}_{\alpha}))\}_{\alpha\in L}$ of $S(\mathscr{V})$ $($which is a covering due to the weak left-exactness assumption in \S\ref{sub-admissibleopensubsets}$)$ is refined by a finite covering of the form $\{S(\mathscr{V}_j)\}_{j\in J}$ by quasi-compact open subspaces $\mathscr{V}_j\subseteq\mathscr{V}$.
\end{itemize}}
\end{dfn}

Notice here that the collection $\{\mathscr{U}_{\alpha}\}_{\alpha\in L}$ of quasi-compact open subspaces of $\mathscr{X}$ may depend on the map $\varphi$; we will show in \ref{lem-scadmissibilitycriterion} that such a collection can be taken independently on $\varphi$.

\begin{prop}\label{prop-scadmissibilitycriterion}
For $\mathscr{X}\in\obj(\mathscr{C})$ a subset $U\subseteq S(\mathscr{X})$ is $(\mathscr{C},S)$-admissible if and only if it is of the form $U=S(\mathscr{U})$ by a $(\mathscr{C},S)$-admissible open subspace $\mathscr{U}\subseteq\mathscr{X}$.
\end{prop}

Notice that due to \ref{cor-csadmissiblemaximal1} the $(\mathscr{C},S)$-admissible open subspace $\mathscr{U}$ is uniquely determined by $U$.
To prove the proposition, we need the following lemma:
\begin{lem}\label{lem-scadmissibilitycriterion}
Let $\mathscr{X}$ be a coherent rigid space in $\mathscr{C}$, and $\{\mathscr{U}_{\alpha}\}_{\alpha\in L}$ a collection of quasi-compact open subspaces of $\mathscr{X}$.
Then $\{\mathscr{U}_{\alpha}\}_{\alpha\in L}$ covers $\mathscr{X}$ if and only if the following conditions are satisfied$:$
\begin{itemize}
\item[{\rm (a)}] $S(\mathscr{X})=\bigcup_{\alpha\in L}S(\mathscr{U}_{\alpha});$
\item[{\rm (b)}] $\{S(\mathscr{U}_{\alpha})\}_{\alpha\in L}$ admits a finite refinement of the form $\{S(\mathscr{V}_j)\}_{j\in J}$, where $\mathscr{V}_j$ are quasi-compact open subspaces of $\mathscr{X}$.
\end{itemize}
\end{lem}

\begin{proof}
The `only if' part is clear.
Let us show the `if' part.
Take for each $j\in J$ an index $\alpha(j)\in L$ such that $S(\mathscr{V}_j)\subseteq S(\mathscr{U}_{\alpha(j)})$.
Since $S$ is spectral and since $\mathscr{V}_j$ and $\mathscr{U}_{\alpha(j)}$ are quasi-compact, we have $\mathscr{V}_j\subseteq\mathscr{U}_{\alpha(j)}$.
From $S(\mathscr{X})=\bigcup_{j\in J}S(\mathscr{V}_j)=S(\bigcup_{j\in J}\mathscr{V}_j)$, we deduce $\mathscr{X}=\bigcup_{j\in J}\mathscr{V}_j$, since $S$ is spectral, and thus we have $\mathscr{X}=\bigcup_{\alpha\in L}\mathscr{U}_{\alpha}$, as desired.
\end{proof}

\begin{proof}[Proof of Proposition {\rm \ref{prop-scadmissibilitycriterion}}]
Suppose there exists a $(\mathscr{C},S)$-admissible open subspace $\mathscr{U}\subseteq\mathscr{X}$ such that $U=S(\mathscr{U})$.
Then any open covering $\{\mathscr{U}_{\alpha}\}_{\alpha\in L}$ of $\mathscr{U}$ consisting of quasi-compact open subspaces fulfills the condition \ref{dfn-csadmissibilityS} (b), and thus the `if' part is clear.
To show the `only if' part, consider the set $\{\mathscr{U}\,|\,S(\mathscr{U})=U\}$ (non-empty by the definition \ref{dfn-csadmissibilityS} (a)).
Due to the continuity of $S$, this set has the maximal element $\mathscr{U}$.
We need to show that this $\mathscr{U}$ is $(\mathscr{C},S)$-admissible.
Let $\mathscr{V}$ be a coherent rigid space in $\mathscr{C}$, and $\varphi\colon\mathscr{V}\rightarrow\mathscr{X}$ a morphism in $\mathscr{C}$ such that $S(\varphi)(S(\mathscr{V}))\subseteq U$.
Let $\{\mathscr{U}_{\alpha}\}_{\alpha\in L}$ be as in \ref{dfn-csadmissibilityS} (b), and set $\mathscr{U}'=\bigcup_{\alpha\in L}\mathscr{U}_{\alpha}$.
Since $S(\mathscr{U})=S(\mathscr{U}')=U$, we have $\mathscr{U}'\subseteq\mathscr{U}$.
Since the covering $\{S(\varphi^{-1}(\mathscr{U}_{\alpha}))\}_{\alpha\in L}$ of $S(\mathscr{V})$ admits a finite refinement $\{S(\mathscr{V}_j)\}_{j\in J}$ as in \ref{dfn-csadmissibilityS} (b), by \ref{lem-scadmissibilitycriterion} we have $\mathscr{V}=\bigcup_{\alpha\in L}\varphi^{-1}(\mathscr{U}_{\alpha})=\varphi^{-1}(\mathscr{U}')$ and thus $\varphi(\mathscr{V})\subseteq\mathscr{U}'\subseteq\mathscr{U}$.
\end{proof}

\begin{dfn}\label{dfn-csadmissibilitycovtop}{\rm 
Let $\mathscr{X}$ be a rigid space in $\mathscr{C}$, and $U\subseteq S(\mathscr{X})$ a $(\mathscr{C},S)$-admissible subset.
A collection $\{U_{\alpha}\}_{\alpha\in L}$ of subsets of $S(\mathscr{X})$ is called a {\em $(\mathscr{C},S)$-admissible covering} of $U$ if the following conditions are satisfied:
\begin{itemize}
\item[{\rm (a)}] each $U_{\alpha}$ is a $(\mathscr{C},S)$-admissible subset of $S(\mathscr{X})$;
\item[{\rm (b)}] for any morphism $\varphi\colon\mathscr{V}\rightarrow\mathscr{X}$ in $\mathscr{C}$ from a coherent rigid space such that $S(\varphi)(S(\mathscr{V}))\subseteq U$, $\{S(\varphi)^{-1}(U_{\alpha})\}_{\alpha\in L}$ admits a finite refinement of the form $\{S(\mathscr{V}_j)\}_{j\in J}$, where each $\mathscr{V}_j\subseteq\mathscr{V}$ is a quasi-compact open subspace, that covers $S(\mathscr{V})$.
\end{itemize}}
\end{dfn}

One has the following proposition, which can be shown by an argument similar to that in \ref{lem-scadmissibilitycriterion}:
\begin{prop}\label{prop-admissibilitycovcompatibility}
Let $\mathscr{X}$ be a rigid space in $\mathscr{C}$, and $\{\mathscr{U}_{\alpha}\}_{\alpha\in L}$ a family of $(\mathscr{C},S)$-admissible open subspaces of $\mathscr{X}$.
Then $\{\mathscr{U}_{\alpha}\}_{\alpha\in L}$ is a $(\mathscr{C},S)$-admissible covering of $\mathscr{X}$ if and only if $\{S(\mathscr{U}_{\alpha})\}_{\alpha\in L}$ is a $(\mathscr{C},S)$-admissible covering of $S(\mathscr{X})$. \hfill$\square$
\end{prop}

\subsection{Rigid analytic geometry}\label{sub-rigidanalyticgeometry}
\index{rigid analytic space|(}
\subsubsection{Classical affinoids}\label{subsub-classicalaffinoidvarieties}
Let $V$ be an $a$-adically complete valuation ring of height one ($a\in\m_V\setminus\{0\}$), and $K=\Frac(V)$.
As in {\bf \ref{ch-pre}}, \S\ref{subsub-nonarchnorms}, the field $K$ has a non-archimedean valuation $\|\cdot\|\colon K\rightarrow\R_{\geq 0}$ and is complete with respect to the associated metric topology.

Let $\mathbf{Aff}_K$ be the category of classical affinoid algebras\index{algebra!affinoid algebra@affinoid ---!classical affinoid algebra@classical --- ---} over $K$ ({\bf \ref{ch-pre}}.\ref{dfn-classicalaffinoidalgebras}) and $K$-algebra homomorphisms.
We consider the dual category $\mathbf{Aff}^{\opp}_K$, called the category of (classical) {\em affinoids} over $K$.
For a classical affinoid algebra $\mathcal{A}$, we consider the set $\Spm\mathcal{A}$ of all closed points in $\Spec\mathcal{A}$; recall that any classical affinoid algebra is a Jacobson ring\index{Jacobson ring} ({\bf \ref{ch-pre}}.\ref{prop-classicalaffringjacobson}).
Then for any morphism $\varphi\colon\mathcal{A}\rightarrow\mathcal{B}$ in $\mathbf{Aff}_K$, the corresponding morphism in $\mathbf{Aff}^{\opp}_K$ is interpreted into the map $\Spm\mathcal{B}\rightarrow\Spm\mathcal{A}$ ($\m\mapsto\varphi^{-1}(\m)$), of which the existence is guaranteed by {\bf \ref{ch-pre}}.\ref{cor-Tatealgebraprime}.

Subsets of an affinoid $\Spm\mathcal{A}$ of the form
$$
\{x\in\Spm\mathcal{A}\,|\,\|f_1(x)\|\leq 1,\ldots,\|f_r(x)\|\leq 1\}
$$
for $f_1,\ldots,f_r\in\mathcal{A}$ form a basis of a topology on the set $\Spm\mathcal{A}$, called the {\em canonical topology} (cf.\ \cite[\S9]{Tate1}\cite[7.2.1]{BGR}).
Here, for an element $f\in\mathcal{A}$ and a point $x=\m\in\Spm\mathcal{A}$, we write $f(x)=(f\ \mathrm{mod}\ \m)$, which lies in the finite extension $K'=\mathcal{A}/\m$ of $K$ ({\bf \ref{ch-pre}}.\ref{cor-Tatealgebraprime}), and denote, by a slight abuse of notation, the unique extension of the norm $\|\cdot\|$ on $K$ to $K'$ by the same symbol. 
Notice that any morphism $\Spm\mathcal{B}\rightarrow\Spm\mathcal{A}$ in $\mathbf{Aff}^{\opp}_K$ is continuous with respect to the canonical topologies.

\subsubsection{Affinoid subdomains}\label{subsub-classicalaffinoidsubvarieties}
\begin{dfn}\label{dfn-classicalaffinoidsubvarieties}{\rm 
A subset $U\subseteq\Spm\mathcal{A}$ is said to be an {\em affinoid subdomain}\index{affinoid!affinoid subdomain@--- subdomain} if the functor $F_U\colon\mathbf{Aff}_K\rightarrow\Sets$ defined by
$$
F_U(\mathcal{B})=\{\varphi\in\Hom_K(\mathcal{A},\mathcal{B})\,|\,\varphi^{\ast}(\Spm\mathcal{B})\subseteq U\}
$$
for $\mathcal{B}\in\obj(\mathbf{Aff}_K)$ is representable by a classical affinoid algebra over $K$; in other words, these exists a $K$-algebra homomorphism $\mathcal{A}\rightarrow\mathcal{A}_U$ of classical affinoid algebras over $K$ that satisfies the following universal mapping property: for any $K$-algebra homomorphism $\mathcal{A}\rightarrow\mathcal{B}$ of classical affinoid algebras over $K$ such that the image of the induced map $\Spm\mathcal{B}\rightarrow\Spm\mathcal{A}$ lies in $U$, there exists a unique $K$-algebra homomorphism $\mathcal{A}_U\rightarrow\mathcal{B}$ such that the diagram
$$
\xymatrix@C-3ex@R-1.5ex{\mathcal{A}\ar[rr]\ar[dr]&&\mathcal{A}_U\ar@{-->}[dl]\\ &\mathcal{B}}
$$
commutes.}
\end{dfn}

The classical affinoid algebra $\mathcal{A}_U$ as above is uniquely determined up to unique isomorphisms.
It is, moreover, known that the map $\Spm\mathcal{A}_U\rightarrow\Spm\mathcal{A}$ gives a homeomorphism onto $U$ (cf.\ \cite[7.2.2, 7.2.5]{BGR}).

Let $\mathcal{A}$ be a classical affinoid algebra over $K$, and $A$ a topologically finitely geberated $V$-algebra such that $\mathcal{A}=A[\frac{1}{a}]$ (cf.\ {\bf \ref{ch-pre}}.\ref{dfn-classicalaffinoidalgebras}).
Consider the affinoid $\mathscr{X}=(\Spf A)^{\rig}$ in the sense of \ref{dfn-cohomologyrigidsp1}.
Then by \ref{cor-classicalpointsexist0} the set $\Spm\mathcal{A}$ coincides with the set of all classical points $\mathscr{X}_0$ of $\mathscr{X}$.
It is straightforward to see, by the aid of \ref{thm-canonicityclassicalaffinoid}, that for any affinoind subdomain (in the sense of \ref{dfn-affinoidsubdomain1}) $\mathscr{U}=(\Spf B)^{\rig}\subseteq\mathscr{X}$ the subset $U=\mathscr{U}_0\subseteq\mathscr{X}_0$ of classical points gives an affinoid subdomain of $\Spm\mathcal{A}$, with the classical affinoid algebra $\mathcal{A}_U=B[\frac{1}{a}]$.
The famous Gerritzen-Grauert theorem \cite{GG} states the converse\begin{thm}[{\rm cf.\ Gerritzen-Grauert \cite{GG}}]\label{thm-GGref}
In the situation as above, a subset $U\subseteq\Spm\mathcal{A}$ is an affinoid subdomain if and only if there exists an affinoid subdomain $\mathscr{U}=(\Spf B)^{\rig}\subseteq\mathscr{X}=(\Spf A)^{\rig}$ $($in the sense of {\rm \ref{dfn-affinoidsubdomain1}}$)$ such that $U$ coincides with the set of classical points $\mathscr{U}_0$ of $\mathscr{U};$ in this case, the related classical affinoid algebra is given by $\mathcal{A}_U=B[\frac{1}{a}]$.
\end{thm}

We will prove this theorem independently from \cite{GG} in the next volume (cf.\ Introduction).
But, for the reader's convenience, let us include here the argument to deduce this theorem from the classical Gerritzen-Grauert theorem (e.g., \cite{BGR}, 7.3.5).

Suppose $U=\Spm\mathcal{B}\hookrightarrow X=\Spm\mathcal{A}$ is an affinoid subdomain, and take $\mathscr{U}=(\Spf B)^{\rig}\rightarrow\mathscr{X}=(\Spf A)^{\rig}$ such that $\mathscr{U}_0\rightarrow\mathscr{X}_0$ is the given $U\hookrightarrow X$.
Take a finite covering by rational subdomains $X=\bigcup_{\alpha\in L}X_{\alpha}$, where $X_{\alpha}=\Spm\mathcal{A}_{\alpha}$, such that $U\cap X_{\alpha}=\Spm\mathcal{A}_{\alpha}\widehat{\otimes}_{\mathcal{A}}\mathcal{B}$ is a Weierstrass domain of $X_{\alpha}=\Spm\mathcal{A}_{\alpha}$ (\cite{BGR}, 7.3.5/2).
Since $X_{\alpha}$ is a rational subdomain of $X$, there exists $A_{\alpha}$ as in \ref{exas-affinoidsubdomain} (3) such that $\mathscr{X}_{\alpha}=(\Spf A_{\alpha})^{\rig}$ gives a rational subdomain of $\mathscr{X}$ and $X_{\alpha}=(\mathscr{X}_{\alpha})_0$ for $\alpha\in L$.
Since $U\cap X_{\alpha}$ is a Weierstrass domain, there exists $B_{\alpha}$ as in  \ref{exas-affinoidsubdomain} (1) such that $\mathscr{U}_{\alpha}=(\Spf B_{\alpha})^{\rig}$ gives a rational subdomain of $\mathscr{X}_{\alpha}$ and $U\cap X_{\alpha}=(\mathscr{U}_{\alpha})_0$ for $\alpha\in L$.
Since $\mathscr{U}=\bigcup_{\alpha\in L}\mathscr{U}_{\alpha}$, and $\mathscr{U}_{\alpha}\hookrightarrow\mathscr{X}_{\alpha}$ and $\mathscr{X}_{\alpha}\hookrightarrow\mathscr{X}$ are all open immersions for $\alpha\in L$, we deduce that $\mathscr{U}\hookrightarrow\mathscr{X}$ is an open immersion, as desired.

\subsubsection{Rigid analytic spaces}\label{subsub-rigidanalyticvarieties}
We introduce the so-called weak G-topology and strong G-topology on affinoids as follows:
\begin{itemize}
\item (\cite[9.1.4]{BGR}) the {\em weak G-topology} on an affinoid $X=\Spm\mathcal{A}$ is defined as follows:
\begin{itemize}
\item admissible open subsets are affinoid subdomains;
\item admissible coverings are finite coverings by affinoid subdomains;
\end{itemize}
\item (\cite[9.1.4/2]{BGR}) the so-called {\em strong G-topology} on an affinoid $X=\Spm\mathcal{A}$ is defined as follows:
\begin{itemize}
\item a subset $U\subseteq X$ is an admissible open subset if and only if it has a covering $U=\bigcup_{\alpha\in L}U_{\alpha}$ by affinoid subdomains of $X$ such that for any morphism $f\colon Y\rightarrow X$ of affinoids with $f(Y)\subseteq U$ the induced covering $\{f^{-1}(U_{\alpha})\}_{\alpha\in L}$ is refined by an admissible covering of weak G-topology on $Y$;
\item a covering $U=\bigcup_{\alpha\in L}V_{\alpha}$ of an admissible open subset $V$ by admissible open subsets is an admissible covering if and only if for any morphism $f\colon Y\rightarrow X$ of affinoids with $f(Y)\subseteq U$ the induced covering $\{f^{-1}(V_{\alpha})\}_{\alpha\in L}$ is refined by an admissible covering of weak G-topology on $Y$.
\end{itemize}
\end{itemize}

The structure presheaf $\O_X$ on the affinoid $X=\Spm\mathcal{A}$ is defined as follows: For any affinoid subdomain $U\subseteq X$, we set $\O_X(U)=\mathcal{A}_U$. It follows from the universal mapping property of affinoid subdomains that this indeed gives the definition of a presheaf $\O_X$.
\begin{thm}[{Tate's acyclicity theorem \cite[\S8]{Tate1}}]\label{thm-tateacyclicitytheoremclassical} The presheaf $\O_X$ is a sheaf with respect to the weak G-topology, and hence uniquely extends to a sheaf with respect to the strong G-topology $($cf.\ {\rm \cite[9.2.3]{BGR})}. \hfill$\square$
\end{thm}

\begin{dfn}[{cf.\ \cite[9.3.1/4]{BGR}}]\label{dfn-tatesrigidanalyticspaceclassical}{\rm 
{\rm (1)} An {\em affinoid space} over $K$ is a locally G-ringed space $(X,\O_X)$ over $K$, that is, a pair consisting of a G-topological space $X$ and a sheaf of local rings over $K$ on $X$ with respect to the G-topology, that is $K$-isomorphic to the locally G-ringed space of the form $(\Spm\mathcal{A},\O_{\Spm\mathcal{A}})$ considered with the strong G-topology.

{\rm (2)} A {\em rigid analytic space} over $K$ is locally G-ringed space $(X,\O_X)$ over $K$ such that the following conditions are satisfied$:$
\begin{itemize}
\item[{\rm (a)}] the G-topology of $X$ enjoys the following properties$:$
\begin{itemize}
\item[{\rm ($\mathrm{G_0}$)}] $\emptyset$ and $X$ are admissible open subsets of $X;$
\item[{\rm ($\mathrm{G_1}$)}] for an admissible open subset $U\subseteq X$ and a subset $V\subseteq U$, $V$ is admissible open in $X$ if there exists an admissible covering $\{U_{\alpha}\}_{\alpha\in L}$ of $U$ such that each $V\cap U_{\alpha}$ is admissible open in $X;$
\item[{\rm ($\mathrm{G_2}$)}] if an open covering of an admissible open subset $U$ consisting of admissible open subsets admits a refinement by admissible covering, it is an admissible covering of $U$.
\end{itemize}
\item[{\rm (b)}] there exists an admissible covering $\{X_{\alpha}\}_{\alpha\in L}$ of $X$ such that, for each $\alpha\in L$, $(X_{\alpha},\O_X|_{X_{\alpha}})$ is an affinoid space over $K$.
\end{itemize}

{\rm (3)} A {\em morphism} of rigid analytic spaces $(X,\O_X)\rightarrow(Y,\O_Y)$ over $K$ is a $K$-morphism of locally G-ringed spaces over $K$.}
\end{dfn}

\subsubsection{Comparison with rigid spaces}\label{subsub-comparisonrigidanalyticvarieties}
Let $V$ and $K$ be as in the beginning of \S\ref{subsub-rigidanalyticvarieties}.
Set $\mathscr{S}=(\Spf V)^{\rig}$, and let $\mathscr{C}$ be the category of rigid spaces locally of finite type over $\mathscr{S}$.
By \ref{cor-spectralfunctorseparatedquotients2} the functor $S=\ZR{\,\cdot\,}^{\cl}$ is a continuous spectral functor defined on $\mathscr{C}$ satisfying the properties (a) and (b) in the beginning of \S\ref{sub-admissibleopensubsets}.

For a rigid space $\mathscr{X}$ in $\mathscr{C}$, let 
$$
\mathscr{X}_0=(\mathscr{X}_0,\tau_0,\O_{\mathscr{X}_0})
$$
be a G-ringed space over $K$ defined as follows:
\begin{itemize}
\item as a set, $\mathscr{X}_0$ is the set of all classical points\index{point!classical point@classical ---}\index{classical point} of $\mathscr{X}$:
$$
\mathscr{X}_0=\ZR{\mathscr{X}}^{\cl};
$$
\item the G-topology $\tau_0$ is the one on $\mathscr{X}_0=S(\mathscr{X})$ by admissible open subsets as in \ref{dfn-csadmissibilityS} and admissible coverings as in \ref{dfn-csadmissibilitycovtop}. (It is clear that $(\mathscr{C},S)$-admissible subsets of $\mathscr{X}_0=S(\mathscr{X})$ are open subsets with respect to the subspace topology on $\mathscr{X}_0$ induced from the topology on the Zariski-Riemann space $\ZR{\mathscr{X}}$);
\item $\O_{\mathscr{X}_0}$ is the sheaf of local rings defined as follows: for any quasi-compact open subspace $\mathscr{U}\subseteq\mathscr{X}$
$$
\O_{\mathscr{X}_0}(\ZR{\mathscr{U}}^{\cl})=\O_{\mathscr{X}}(\ZR{\mathscr{U}})
$$
(notice that quasi-compact open subspaces are $(\mathscr{C},S)$-admissible).
\end{itemize}

It is clear by the construction that for an $\mathscr{S}$-morphism $\varphi\colon\mathscr{X}\rightarrow\mathscr{Y}$ of locally of finite type rigid spaces, one has the uniquely induced morphism $\varphi_0\colon\mathscr{X}_0\rightarrow\mathscr{Y}_0$ of G-ringed spaces over $K$.
Moreover, if $\varphi$ is an open immersion, then so is $\varphi_0$.

\begin{thm}\label{thm-rigidanalyticcomparison}
For a locally of finite type $\mathscr{S}$-rigid space $\mathscr{X}$, the G-ringed space $\mathscr{X}_0=(\mathscr{X}_0,\tau_0,\O_{\mathscr{X}_0})$ thus obtained is a rigid analytic space over $K$ in the sense of {\rm \cite[(9.3.1/4)]{BGR}}.
Thus we have a functor 
$$
\mathscr{X}\longmapsto\mathscr{X}_0
$$
from the category of locally of finite type $\mathscr{S}$-rigid spaces to the category of rigid analytic spaces over $K$.
Moreover, this functor gives a categorical equivalence from the category of quasi-separated locally of finite type $\mathscr{S}$-rigid spaces to the category of quasi-separated rigid analytic spaces over $K$.
\end{thm}

Here, a rigid analytic space is said to be {\it quasi-separated} (in the sense of Tate) if the intersection of two affinoid subdomains is a finite union of affinoid subdomains.
Note that the quasi-separatedness in the last assertion is necessary due to the following example. Let $\mathscr{X}$ be any non-empty finite type affinoid over $\mathscr{S}$, and $\mathscr{U}=\mathscr{X}\setminus\{x\}$, where $x$ is a non-classical closed point of $\mathscr{X}$. Then the gluing of two copies of $\mathscr{X}$ along $\mathscr{U}$ is a non-quasi-separated $\mathscr{S}$-rigid space, of which the associated rigid analytic space is the same as that of $\mathscr{X}$.

\begin{proof}[Proof of Theorem {\rm \ref{thm-rigidanalyticcomparison}}]
Let us first show step by step that $\mathscr{X}_0=(\mathscr{X}_0,\tau_0,\O_{\mathscr{X}_0})$ is a rigid analytic space over $K$ and, at the same time, that the functor $\mathscr{X}\mapsto\mathscr{X}_0$ is fully faithful.

\medskip
{\sc Step 1.} Suppose $\mathscr{X}$ is an affinoid $\mathscr{X}=(\Spf A)^{\rig}$, where $A$ is a topologically finitely generated\index{finitely generated!topologically finitely generated@topologically ---} algebra over $V$.
In this situation the ring $\mathcal{A}=A[\frac{1}{a}]$ is a classical affinoid algebra\index{algebra!affinoid algebra@affinoid ---!classical affinoid algebra@classical --- ---} ({\bf \ref{ch-pre}}, \S\ref{subsub-classicalaffinoidalgebras}).
By \ref{cor-classicalpointsexist0} the set $\mathscr{X}_0$ is canonically identified with the set $\Spm\mathcal{A}$ of all closed points in $\Spec\mathcal{A}$.

\medskip
{\sc Claim.} {\it The G-topology $\tau_0$ coincides with the strong G-topology on the affinoid $\mathscr{X}_0=\Spm\mathcal{A}$.}

\medskip
To show the claim, first notice the following facts:
\begin{itemize}
\item for a quasi-compact open subspace $\mathscr{U}\subseteq\mathscr{X}$, $\mathscr{U}_0\subseteq\mathscr{X}_0$ is an affinoid subdomain in the sense of \ref{dfn-classicalaffinoidsubvarieties} if and only if $\mathscr{U}$ is an affinoid subdomain of $\mathscr{X}$ in the sense of \ref{dfn-affinoidsubdomain1} (Gerritzen-Grauert theorem \cite{GG}; cf.\ \ref{thm-GGref});
\item for a finite collection $\{\mathscr{U}_{\alpha}\}_{\alpha\in L}$ of quasi-compact open subspaces of $\mathscr{X}$, $\mathscr{X}_0=\bigcup_{\alpha\in L}(\mathscr{U}_{\alpha})_0$ if and only if $\mathscr{X}=\bigcup_{\alpha\in L}\mathscr{U}_{\alpha}$; this follows from that $S=\ZR{\,\cdot\,}^{\cl}$ is a spectral functor.
\end{itemize}
Then the claim is evident due to the definitions of admissible open subsets (\ref{dfn-csadmissibilityS}) and admissible coverings (\ref{dfn-csadmissibilitycovtop}).

Now it is immediate from \ref{thm-comparisonaffinoid} and the definition of the sheaf $\O_{\mathscr{X}_0}$ as above that the resulting G-ringed space $(\mathscr{X}_0,\tau_0,\O_{\mathscr{X}_0})$ is an affinoid space in the sense as in \ref{dfn-tatesrigidanalyticspaceclassical} (1).

Finally we show that the functor $\mathscr{X}\mapsto\mathscr{X}_0$ restricted on the full subcategory consisting of finite type affinoids over $\mathscr{S}$ is fully faithful (with the essential image being the full subcategory consisting of affinoid spaces over $K$).
Consider affinoids $\mathscr{X}=(\Spf A)^{\rig}$ and $\mathscr{Y}=(\Spf B)^{\rig}$, where $A,B$ are topologically finitely generated $V$-algebras.
It is well-known in classical rigid analytic geometry that morphisms $\mathscr{X}_0=\Spm\mathcal{A}\rightarrow\mathscr{Y}_0=\Spm\mathcal{B}$ (where $\mathcal{A}=A[\frac{1}{a}]$ and $\mathcal{B}=B[\frac{1}{a}]$) of affinoid spaces over $K$ are in canonical one-to-one correspondence with $K$-algebra homomorphisms $\mathcal{B}\rightarrow\mathcal{A}$. 
Based on this fact, the desired fully faithfulness follows from \ref{thm-canonicityclassicalaffinoid}.

\medskip
{\sc Step 2.} Suppose $\mathscr{X}$ is a quasi-compact open subspace of an affinoid.
In this case, since $\mathscr{X}$ is the union of finitely many affinoids $\mathscr{U}_{\alpha}$ in such a way that any intersection $\mathscr{U}_{\alpha}\cap\mathscr{U}_{\beta}$ is again an affinoid, it is easy to see that $(\mathscr{X}_0,\tau_0,\O_{\mathscr{X}_0})$ is actually a rigid analytic space over $K$.
One can also show that the functor $\mathscr{X}\mapsto\mathscr{X}_0$ restricted on the full subcategory consisting of rigid spaces of this type is fully faithful; indeed, this follows from the affinoid case by standard patching argument.

\medskip
{\sc Step 3.} Suppose $\mathscr{X}$ is coherent.
In this case $\mathscr{X}$ is the union of finitely many affinoids $\mathscr{U}_{\alpha}$ in such a way that any intersection $\mathscr{U}_{\alpha}\cap\mathscr{U}_{\beta}$ is quasi-compact.
Hence the assertion follows from {\sc Step 1} and {\sc Step 2}.

\medskip
{\sc Step 4.} If $\mathscr{X}$ is quasi-separated, then by \ref{prop-generalrigidspace31-1} it is a stretch of coherent rigid spaces (\ref{dfn-admissiblesite31} (1)), and thus the assertion follows easily.
Then the assertion in the general case also follows in view of the definition of general rigid spaces as in \ref{dfn-generalrigidspace1}.

\medskip
Next we show the other assertion, that is, the essential surjectivity.
By \cite[(9.3.2), (9.3.3)]{BGR} (by patching of analytic spaces and analytic mappings) we may restrict to the affinoid case.
But then the assertion is clear.
\end{proof}

The following is a corollary of the proof:
\begin{cor}\label{cor-rigidanalyticcomparison}
Under the categorical equivalence in {\rm \ref{thm-rigidanalyticcomparison}}, affinoids $($resp.\ separated spaces, resp.\ coherent spaces$)$ corresponds to affinoid spaces $($resp.\ separated analytic spaces, resp.\ coherent analytic spaces$)$. \hfill$\square$
\end{cor}
\index{rigid analytic space|)}

\subsubsection{Coherent sheaves}\label{subsub-rigidanalyticcoherentsheaves}
Let $\mathscr{X}$ be a rigid space locally of finite type over $\mathscr{S}=(\Spf V)^{\rig}$, and $\mathscr{X}_0=(\mathscr{X}_0,\tau_0,\O_{\mathscr{X}_0})$ the associated rigid analytic space.
By \ref{cor-classicalpointsexist4} the subset $\mathscr{X}_0=\ZR{\mathscr{X}}^{\cl}$ of $\ZR{\mathscr{X}}$ is very dense.
Hence it follows from the definitions of the G-topology $\tau_0$ and the structure sheaf $\mathscr{O}_{\mathscr{X}_0}$ that we have the canonical equivalence between the locally ringed topoi
$$
(\mathscr{X}^{\sim}_{\ad},\O_{\mathscr{X}})\cong(\mathscr{X}^{\sim}_{0,\ad},\O_{\mathscr{X}_0})
$$
(induced from the obvious inclusion $\mathscr{X}_0\hookrightarrow\mathscr{X}$ of sites), where $\mathscr{X}^{\sim}_{\ad}$ denotes the small admissible topos\index{topos!admissible topos@admissible ---} associated to $\mathscr{X}$ (\ref{dfn-admissiblesite3gensmall}) and $\mathscr{X}^{\sim}_{0,\ad}$ the topos associated to $\tau_0$.
In particular, the functor $\mathscr{F}\mapsto\mathscr{F}_0=\mathscr{F}|_{\mathscr{X}_0}$ gives rise to categorical equivalences 
$$
\Mod_{\mathscr{X}}\stackrel{\sim}{\rightarrow}\Mod_{\mathscr{X}_0},\quad\Coh_{\mathscr{X}}\stackrel{\sim}{\rightarrow}\Coh_{\mathscr{X}_0}.
$$

\section{Appendix: Non-archimedean analytic space of Banach type}\label{sec-valspecbanach}
In this appendix, we will introduce and discuss metrized analytic spaces, and the relationship with Berkovich's analytic geometry. 
Our construction is based on a spectral theory of filtered rings (\S\S\ref{sub-vsbr-banachrings}--\ref{sub-vsbr-valuativespectrum}).
It will turn out that the new notion of spectra, the so-called {\it valuative spectra},  is, in certain situation, equivalent to that of the so-called reified adic spectra by Kedlaya \cite{Kedlaya} and can be seen as a globalization of Temkin's work \cite{Temkin1}.
The interpretation of Berkovich spaces is given in \S\ref{sub-vsbr-Berkovichanalyticspaces}.

\subsection{Seminorms and norms}\label{sub-vsbr-banachrings}
Let us first recall some basic notions and terminologies related to seminormed and normed rings (cf.\ \cite{BGR}, \cite{Berk1}).

Let $M$ be an abelian group, written additively.
A {\it seminorm}\index{seminorm} on $M$ is a function of the form 
$$
|\cdot|\colon M\longrightarrow\R_{\geq 0},\qquad x\longmapsto |x|
$$
satisfying the following conditions:
\begin{itemize}
\item[(a)] $|0|=0$
\item[(b)] for $x,y\in M$, $|x-y|\leq |x|+|y|$.
\end{itemize}
A seminorm $|\cdot|$ is called a {\it norm} if $\{x\in M:|x|=0\}$ is coincides with $\{0\}$.
A seminorm or a norm $|\cdot|$ is said to be {\it non-archimedean} if it satisfies the following condition:
\begin{itemize}
\item[(b$)'$] for $x,y\in M$, $|x-y|\leq\max\{|x|,|y|\}$.
\end{itemize}

Two seminorms $|\cdot|$ and $|\cdot|'$ on $M$ are said to be {\it equivalent} if there exist real numbers $C,C'>0$ such that, for any $x\in A$,
$$
|x|'\leq C|x|\leq C'|x|'.
$$

A {\it zero seminorm} is a seminorm $|\cdot|$ such that $|x|=0$ for all $x\in M$.
A {\it trivial norm} is a norm $|\cdot|$ idefined by, for $x\in M$, 
$$
|x|=\begin{cases}0&(\textrm{if $x=0$}),\\ 1&(\textrm{otherwise}).\end{cases}
$$

A {\it seminormed group} (resp.\ {\it normed group}) is a pair $(M,|\cdot|)$ consisting of an abelian group $M$ and a seminorm (resp.\ norm) $|\cdot|$ on $M$.

Let $A$ be ring.
A {\it $($ring$)$ seminorm}\index{seminorm} on $A$ is a seminorm on the underlying additive group of $A$ such that 
\begin{itemize}
\item[(c)] $|1|\leq1$;
\item[(d)] for $x,y\in A$, $|xy|\leq |x|\cdot |y|$.
\end{itemize}

By the first condition, we see that any ring seminorm on a non-zero ring is not the zero seminorm.

A seminorm $|\cdot|$ on a ring $A$ is said to be {\it power-multiplicative} if $|1|=1$ and $|x^n|=|x|^n$ for any $x\in A$ and $n\geq 1$.
It is {\it multiplicative} if $|1|=1$ and $|xy|=|x|\cdot |y|$ for any $x,y\in A$.
A multiplicative norm is called a {\it valuation}.
If $|\cdot|$ and $|\cdot|'$ are seminorms on $A$ with $|\cdot|$ power-bounded, and if there exists a real number $C>0$ such that $|x|\leq C|x|'$ for any $x\in A$, then $|x|\leq |x|'$ for any $x\in A$ (that is, $C$ can be chosen to be $1$).
Indeed, for any positive integer $n$, we have $|x|^n\leq C|x|^{\prime n}$, hence $|x|\leq\sqrt[n]{C}|x|'$, thereby the claim.
In particular, two power-bounded seminorms on a ring are equivalent if and only if they are equal.

A {\it seminormed ring} (resp.\ {\it normed ring}) is a pair $(A,|\cdot|)$ consisting of a ring $A$ and a seminorm (resp.\ norm) $|\cdot|$ on $A$.

For a seminormed group $(M,|\cdot|)$ and a surjective homomorphism $\pi\colon M\rightarrow N$ of abelian groups, the {\it residue seminorm} $|\cdot|_{\mathrm{res}}$ on $N$ induced from $|\cdot|$ is the seminorm on $N$ defined by, for $y\in N$, 
$$
|y|_{\mathrm{res}}=\inf\{|x|: \pi(x)=y\}.
$$
When $\pi$ is a ring homomorphism, then the residue seminorm induced from a ring seminorm is again a ring seminorm.

Let $(A,|\cdot|_A)$ and $(B,|\cdot|_B)$ be seminormed groups (resp.\ seminormed rings), and $f\colon A\rightarrow B$ a group (resp.\ ring) homomorphism.
\begin{itemize}
\item $f$ is said to be {\it bounded} if there exists a real number $C>0$ such that 
$$
|f(x)|_B\leq C|x|_A
$$
for any $x\in A$. (By an argument similar to the one as above, the number $C$ can be $1$, if $|\cdot|_B$ is power-bounded.)
\item $f$ is said to be {\it admissible} if the residue seminorm on the image $f(A)$ induced from $|\cdot|_A$ is equivalent to the restriction of $|\cdot|_B$.
Admissible injective (resp.\ surjective) homomorphisms are often called {\it admissible monomorphisms} (resp.\ {\it admissible epimorphisms}).
\end{itemize}

A seminormed module over a seminormed ring $(A,|\cdot|_A)$ is a seminormed group $(M,|\cdot|_M)$ such that 
\begin{itemize}
\item[(e)] there exists $C>0$ such that $|ax|_M\leq C|a|_A|x|_M$ for any $a\in A$ and any $x\in M$.
\end{itemize}

Let $(M,|\cdot|_M)$ and $(N,|\cdot|_N)$ be non-archimedean seminormed modules over a non-archimedean seminormed ring $(A,|\cdot|_A)$.
Then, on tensor product $M\otimes_AN$, we have the following seminorm, called the {\it tensor-product seminorm}: for $x\in M\otimes_AN$, 
$$
|x|=\inf\Big\{\max_{1\leq i\leq n}|m_i|_M|n_i|_N\,\Big|\, x=\sum^n_{i=1}m_i\otimes n_i,\ m_i\in M,\ n_i\in N\Big\}.
$$
The completion of $M\otimes_AN$ with respect to the tensor-product seminorm will be denoted by 
$$
M\,\widehat{\otimes}_A\,N,
$$
and called the {\it complete tensor product}.

A {\it Banach ring} is a normed ring $(A,|\cdot|)$ such that $A$ is complete with respect to the metric by the norm.
Note that the zero ring $\{0\}$ with the zero map as its norm is a Banach ring, which can be characterized as a Banach ring $(A,|\cdot|_A)$ such that $|1|_A<1$.
Moreover, any ring $A$ with the trivial norm $|\cdot|_{\mathrm{triv}}$ is a Banach ring.

\subsection{Graded valuations}\label{sub-vsbr-gradedvaluations}
\index{valuation!valuation ring@--- ring!graded valuation ring@graded --- ---|(}
\subsubsection{Graded rings and modules}\label{subsub-vsbr-gradedrings}
Many notions on rings and modules admit `graded version,' which generalizes the classical theory of rings and modules. 
Here we recall some of the basics on graded rings and modules, especially on graded fields and graded valuation rings (cf.\ \cite{NasOys1}, \cite{NasOys2}, \cite{Temkin1}).

Let $\Delta$ be an abelian group, written multiplicatively, with the unit element $1\in\Delta$.
A {\it $\Delta$-graded ring} is a commutative ring with $1=1_G\in G$ of the form
$$
G=\bigoplus_{d\in\Delta}G_d,
$$
where $G_d$ for $d\in\Delta$ are additive subgroups of $G$, such that $G_d\cdot G_{d'}\subseteq G_{dd'}$ for any $d,d'\in\Delta$.
Similarly, a {\it $\Delta$-graded $G$-module} is a $G$-module of the form $M=\bigoplus_{d\in\Delta}M_d$ satisfying $G_d\cdot M_{d'}\subseteq M_{dd'}$.
Note that, in this situation, the unit-element part $G_1$ is a commutative ring, and all pieces $G_d$ and $M_d$ for $d\in\Delta$ are $G_1$-modules; in fact, one can easily show that $1_G\in G_1$ (cf.\ \cite{NasOys2}, 1.1.1). 
Note also that, for a graded ideal $I=\bigoplus_{d\in\Delta}I_d$ of $G$, $I=G$ if and only if $I_1=G_1$.

Notice that, if $\Delta=\{1\}$ (trivial group), then all what follows throughout this subsection boils down to the classical theory of rings and modules.

For $\Delta$-graded rings $G,G'$, a homomorphism $f\colon G\rightarrow G'$ of $\Delta$-graded rings is defined to be a ring homomorphism that preserves the grading.
One can similarly define morphisms of $\Delta$-graded $G$-modules.
One has thus categories of $\Delta$-graded rings and modules.
If $\Delta'$ is a subgroup of $\Delta$, then, by restriction of the grading, one has the functor from the category of $\Delta$-graded rings to the category of $\Delta'$-graded rings.
A left adjoint to this functor is given by ``$0$-extension'' of the gradings, that is, to any $\Delta'$-graded ring $G'$, we associate $G=\bigoplus_{d\in\Delta}G_d$ with 
$$
G_d=\begin{cases}G_d&(d\in\Delta'),\\ \{0\}&(d\not\in\Delta')\end{cases}
$$
(cf.\ \cite{NasOys2}, 1.1.2).

A {\it $\Delta$-graded integral domain} is a non-zero $\Delta$-graded ring $G=\bigoplus G_d$ such that the product $ab$ of any non-zero homogeneous elements $a,b\in G$ is non-zero.
A {\it $\Delta$-graded field} is a non-zero $\Delta$-graded ring $G$ such that any non-zero homogeneous element is invertible.
Note that a $\Delta$-graded integral domain (resp.\ $\Delta$-graded field) $G$ itself may not be an integral domain (resp.\ a field), but the unit-element part $G_1$ is always an integral domain (resp.\ a field) in the usual sense.

A homogeneous ideal $\mathfrak{p}\subseteq G$ of a $\Delta$-graded ring $G$ is {\it prime} if $\mathfrak{p}\neq G$ and the set $S_{\mathfrak{p}}$ of all homogeneous elements of $G\setminus\mathfrak{p}$ is closed under multiplication, or equivalently, the quotient $\Delta$-graded ring $G/\mathfrak{p}=\bigoplus_{d\in\Delta}G_d/(\mathfrak{p}\cap G_d)$ is a $\Delta$-graded integral domain.
A maximal element in the set of all homogeneous, and not equal to $G$ itself, ideals in $G$ is called a {\it maximal} homogeneous ideal.
It follows that a homogeneous ideal $\m\subseteq G$ is maximal if and only if $G/\m$ is a $\Delta$-graded field.

By a standard argument using Zorn's lemma, one can show that any non-zero $\Delta$-graded ring has a maximal homogeneous ideal.

Let $G$ be a $\Delta$-graded ring, and $\m\subseteq G$ a homogeneous ideal.
Then $\m$ is the unique maximal homogeneous ideal of $G$ if and only if $1\not\in\m$ and any non-invertible homogeneous element of $G$ belongs to $\m$.
In this situation, $G$ is said to be {\it local}.
In other words, a local $\Delta$-graded ring is a $\Delta$-graded ring that has a unique maximal homogeneous ideal.

Let $G$ be a $\Delta$-graded ring, and $S\subseteq G$ a multiplicative subset consisting of homogeneous elements.
Then one has the natural grading by $\Delta$ on the localization $G_S=S^{-1}G$ by $(G_S)_d=\sum_{d'=dd''}(S\cap G_{d''})^{-1}G_{d'}$ for $d\in\Delta$ (cf.\ \cite{NasOys2}, 8.1).
For example, if $\mathfrak{p}\subseteq G$ is a prime homogeneous ideal, then one has the set $S_{\mathfrak{p}}$ as above, and form the localization $S^{-1}_{\mathfrak{p}}G$, which we denote by $G_{\mathfrak{p}}$.
Notice that $G_{\mathfrak{p}}$ is a local $\Delta$-graded ring, where $\mathfrak{p}G_{\mathfrak{p}}$ is the unique maximal homogeneous ideal.
In particular, if $G$ is a $\Delta$-graded integral domain, then $G_{(0)}$, the graded localization by the prime homogeneous ideal $(0)$, is a $\Delta$-graded field, called the {\it graded fractional field}, and denoted by $\Frac_{\Delta}(G)$.

\subsubsection{Graded valuation rings}\label{subsub-vsbr-gradedvaluationrings}
Let $\Delta$ be an abelian group. 
We only consider gradings by $\Delta$, and drop `$\Delta$' from the notation.

\begin{dfn}[\cite{Temkin1}, \S1]\label{dfn-vsbr-gradedvaluationrings}{\rm 
(1) Let $K$ be a graded field. A graded subring $V$ of $K$ is said to be a {\it graded valuation ring for $K$} if for any homogeneous non-zero $x\in K$, either $x$ or $x^{-1}$ is contained in $V$.

(2) A {\it graded valuation ring} is a graded integral domain $V$ that is a graded valuation ring for its graded fractional field $K=\Frac_{\Delta}(V)$.}
\end{dfn}

Notice that, in the situation as in (1), any graded subring $W$ of $K$ that contains $V$ is again a graded valuation ring.
Notice also that, if $V$ is a graded valuation ring for $K$, then the unit-element part $V_1$ is a valuation ring for $K_1$ (in the usual sense; see {\bf \ref{ch-pre}}, \S\ref{subsub-valdfn}).

Let $(G_1,\m_{G_1})$ and $(G_2,\m_{G_2})$ be local graded subrings of a graded field $K$.
We say that $G_1$ {\it dominates} $G_2$ if $G_2\subseteq G_1$ and $\m_{G_2}=\m_{G_1}\cap G_2$.

\begin{prop}\label{prop-vsbr-gradedvaluationrings}
Let $K$ be a graded field, and $V$ a graded subring of $K$.
Then the following conditions are equivalent to each other:
\begin{itemize}
\item[{\rm (a)}] $V$ is a graded valuation ring for $K$;
\item[{\rm (b)}] $\Frac_{\Delta}(V)=K$, and the set of all homogeneous principal ideals of $V$ is totally ordered by inclusion order;
\item[{\rm (c)}] $\Frac_{\Delta}(V)=K$, and the set of all homogeneous ideals of $V$ is totally ordered by inclusion order;
\item[{\rm (d)}] $V$ is local, and $(V,\m_V)$ is the maximal element, with respect to the order by domination, in the set of all local graded subrings of $K$.
\end{itemize}
\end{prop}

We provide the proof, which mimics the proof of the classical one (cf.\ \cite{Bourb1}, Chapter VI, Theorem 1), for the reader's convenience.
\begin{proof}
(a) $\Rightarrow$ (b): For two non-zero principal homogeneous ideals $aV,bV\subseteq V$, we have $aV\subseteq bV$ or $bV\subseteq aV$, according respectively as $a/b\in V$ or $b/a\in V$.

(b) $\Rightarrow$ (c): Let $I,J\subseteq V$ be two homogeneous ideals such that $I\not\subseteq J$.
Take a homogeneous $a\in I\setminus J$. For any non-zero $b\in J$, either $a/b$ or $b/a$ belongs to $V$. But $b/a\in V$ has to hold, since $a\not\in J$.
Hence $b=a\cdot(b/a)\in I$, which shows $J\subseteq I$.

(c) $\Rightarrow$ (a): For any non-zero homogeneous $a/b\in K$, where $a,b$ are homogeneous elements in $V$, we have $a/b\in V$ or $b/a\in V$, according respectively as $aV\subseteq bV$ or $bV\subseteq aV$.

(c) $\Rightarrow$ (d): Take a maximal homogeneous ideal $\m$ of $V$, which is unique due to the assumption. Hence $V$ is a local graded ring.
Let $W$ be a local graded subring of $K$ that dominates $V$.
For any non-zero homogeneous $x\in W$, either $x\in V$ or $x^{-1}\in V$ holds.
Since $\m_VW\subseteq\m_W\neq W$ and $x^{-1}\not\in\m_W$, we have $x\in V$, which implies $W=V$.

(d) $\Rightarrow$ (a): Let $x\in K\setminus V$ be a homogeneous element. 
We want to show that $x^{-1}\in V$.
We first show that $x$ is not integral over $V$.
Indeed, if $x$ is integral over $V$, then the graded subring $W=V[x]$ of $K$ is finite over $V$.
One can take a maximal ideal $\m$ (in the usual sense) of $V$ containing $\m_V$; since $\m W\neq W$, we have $\m_VW\neq W$.
(The last inequality also follows from the graded version of Nakayama's lemma; see \cite{NasOys1}, 8.4.)
Thus there exists maximal homogeneous ideal $\m'$ of $W$ that contains $\m_VW$.
Since $W_{\m'}$ dominates $V$, we have $x\in W_{\m'}=V$ by the assumption, hence contradiction.

Thus $x$ is not integral over $V$, and hence by \cite{Bourb1}, Chapter VI, \S1.2, Lemma 1, the ideal of $V[x^{-1}]$ generated by $\m_V$ and $x^{-1}$ is not equal to $V[x^{-1}]$ itself, and hence there exists a maximal homogeneous ideal $\m''$ of $V[x^{-1}]$ that dominates $V$.
Therefore, we have $x^{-1}\in V[x^{-1}]_{\m''}=V$, as desired.
\end{proof}

The last part of the above proof shows that the following statement holds.
\begin{prop}\label{prop-vsbr-gradedvaluationrings1}
Let $V$ be a graded valuation ring for $K=\Frac_{\Delta}(V)$.
Then $V$ is graded integrally closed, viz., any homogeneous element $x\in K$ integral over $V$ lies in $V$.
\end{prop}

Finally, let us state the existence of graded valuation rings that dominate given local graded rings.
\begin{prop}\label{prop-vsbr-gradedvaluationrings2}
Let $K$ be a graded field, and $A\subseteq K$ a local graded subring.
Then there exists a graded valuation ring $V$ for $K$ that dominates $A$.
\end{prop}

\begin{proof}
By a standard argument using Zorn's lemma, one can show that there exists a maximal local graded subring $V\subseteq K$ that dominates $A$.
By \ref{prop-vsbr-gradedvaluationrings}, $V$ is a graded valuation ring.
\end{proof}

\subsubsection{Graded valuations}\label{subsub-vsbr-valuationgraded}
For a graded ring $G$, let us denote by $h(G)$ the set of all homogeneous elements
$$
h(G)=\bigcup_{d\in\Delta}G_d.
$$
\begin{dfn}\label{dfn-vsbr-valuationgraded}{\rm 
Let $K$ be a graded field, and $\Gamma$ a totally ordered commutative group\index{ordered!totally@totally ---!totally ordered commutative group@--- --- commutative group} ({\bf \ref{ch-pre}}, \S\ref{subsub-ordab}).
A graded valuation\index{valuation!graded valuation@graded ---} $v$ of $K$ with the {\it value target group} $\Gamma$ is a mapping
$$
v\colon h(K)\longrightarrow\Gamma\cup\{\infty\}
$$
that satisfies the following conditions:
\begin{itemize}
\item[{\rm (a)}] $v(xy)=v(x)+v(y)$ for $x,y\in h(K)$; 
\item[{\rm (b)}] $v(x+y)\geq\inf\{v(x),v(y)\}$ for $x,y\in h(K)$ of the same degree; 
\item[{\rm (c)}] $v(1)=0$ and $v(0)=\infty$.
\end{itemize}}
\end{dfn}

If $v$ is a graded valuation of $K$, then, for $d\in\Delta$, 
$$
V_d=\{x\in K_d: v(x)\geq 0\}
$$
is a subgroup of $K_d$ (due to (b) and (c)), and $V=\bigoplus_{d\in\Delta}V_d$ is a graded subring of $K$ (due to (a)).
One deduces, using the conditions (a), (b), and (c), that, for any non-zero homogeneous $x\in K$, we have $v(x)\geq 0$ or $v(x^{-1})\geq 0$.
Hence $V$ is a graded valuation ring for $K$.
The maximal homogeneous ideal of $V$ is the homogeneous ideal generated by $\{x\in h(K): v(x)>0\}$.

Conversely, for a graded valuation ring $V$, one has a graded valuation $v$ of the graded fractional field $K=\Frac_{\Delta}(V)$, which gives back, in the manner as above, the given $V$ up to isomorphism; the construction is as follows.
Let $U(K)$ be the multiplicative group of the non-zero homogeneous elements of $K$, and $U(V)$ the multiplicative group of the invertible homogeneous elements of $V$.
By \ref{prop-vsbr-gradedvaluationrings} (b), one can show that 
$$
\Gamma_V=U(K)/U(V)
$$
is a totally ordered commutative group in such a way that the image of $h(V)\setminus\{0\}$ is the set of all positive elements.
Then the mapping
$$
{\textstyle v\colon h(K)\longrightarrow\Gamma_V\cup\{\infty\}},\quad v(x)=
\begin{cases}
[x]\ (=x\ \mathrm{mod}\ U(V))&\textrm{if}\ x\neq 0,\\
\infty&\textrm{if}\ x=0,
\end{cases}
$$
is a graded valuation of $K$.
It is clear that $V$ is the `non-negative' part; that is, $V_d=\{x\in K_d: v(x)\geq 0\}$ for all $d\in\Delta$.

As in the classical case, one can define the notion of equivalence of graded valuations in such a way that, by the above correspondence, the set of all equivalence classes of graded valuations of $K$ is mapped bijectively onto the set of all graded valuation rings for $K$.

\subsubsection{Generization and specialization of graded valuations}\label{subsub-vsbr-genspe}
Similarly to the classical case, we have the following result. 
\begin{prop}\label{prop-vsbr-genspe}
Let $V$ be a graded valuation ring for a graded field $K$. 
Then there exist canonical order-preserving bijections among the following sets:
\begin{itemize}
\item[{\rm (a)}] the set of all homogeneous prime ideals of $V$ with the inclusion order;
\item[{\rm (b)}] the set of all graded subrings $(\neq V)$ lying between $V$ and $K$ $($which are automatically graded valuation rings$)$ with the reversed inclusion order;
\item[{\rm (c)}] the set of all proper isolated subgroups\index{isolated subgroup} of $\Gamma_V$ with the reversed inclusion order.
\end{itemize}
\end{prop}

The bijections are described similarly to the ones in {\bf \ref{ch-pre}}.\ref{prop-height2}.
\begin{cor}\label{cor-vsbr-genspe}
Let $V$ be a graded valuation ring for a graded field $K$. 
Then any graded subring $W$ of $K$ containing $V$ is a graded valuation ring for $K$, and the set of all such subrings is totally ordered with respect to the inclusion order.
\end{cor}

Based on these results, we can define the following notions similarly to the classical case.
\begin{dfn}\label{dfn-vsbr-genspe}{\rm Let $K$ be a graded field.

(1) For two graded valuation rings $V,W$ for $K$, we say $V$ is a {\it specialization} of $W$, or $W$ is a {\it generization} of $V$, if $V\subseteq W$.

(2) The {\it height}\index{height!height of a graded valuation ring@--- (of a graded valuation (ring))} of a graded valuation ring $V$, denoted by $\mathrm{ht}_{\Delta}(V)$, is the height of the value group $\Gamma_V$.}
\end{dfn}

Note that $\mathrm{ht}_{\Delta}(V)=0$ if and only if $V=K$.

\subsubsection{Unit-element part}\label{subsub-vsbr-unitelementpart}
As mentioned before, the unit-element part $V_1$ of a graded valuation ring $V$ for the graded field $K$ is a valuation ring (in the usual sense) for the unit-part $K_1$ of $K$.
The maximal ideal $\m_{V_1}$ of $V_1$ is equal to $\m_V\cap V_1$.

In general, for a $\Delta$-graded field $F$, we set
$$
\Delta_F=\{d\in\Delta: F_d\neq\{0\}\}.
$$
Then $\Delta_F$ is a subgroup of $\Delta$ (since $F$ is a graded field).
Note that, if $d\in\Delta_F$, we have $F_d\cdot F_{d^{-1}}=F_1$, and $F_d$ is a one dimensional vector space over $F_1$.

\begin{prop}\label{prop-vsbr-unitelementpart1}
Let $V$ be a graded valuation ring for a graded field $K$, and $k=V/\m_k$ the graded residue field of $V$.

{\rm (1)} We have 
$$
\Delta_k=\{d\in\Delta: V_d\cdot V_{d^{-1}}=V_1\}.
$$
In particular, $\Delta_k$ is a subgroup of $\Delta_K$ $($which may not coincide with $\Delta_K)$.

{\rm (2)} For any $d\in\Delta_k$, $V_d$ is a free $V_1$-module of rank one.

{\rm (3)} For any $d\in\Delta$, 
$$
(\m_V)_d=\begin{cases} \m_{V_1}V_d&(d\in\Delta_k),\\ V_d&(d\not\in\Delta_k).\end{cases}
$$
\end{prop}

\begin{proof}
(1) If $d\in\Delta_k$, we have $(\m_V)_d=\m_V\cap V_d\neq V_d$, and hence there exists $u\in V_d$ that is invertible in $V$.
For any $x\in V_1$, we have $x=(xu)\cdot u^{-1}\in V_d\cdot V_{d^{-1}}$, hence $V_1=V_d\cdot V_{d^{-1}}$.
Conversely, if $V_d\cdot V_{d^{-1}}=V_1$, there exist $u\in V_d$ invertible in $V$, and hence $k_d=V_d/(\m_V\cap V_d)\neq 0$.

(2) Since $\m_{V_1}V_d\subseteq \m_V\cap V_d$, we have $\m_{V_1}V_d\neq V_d$ for $d\in\Delta_k$.
Take $x\in V_d\setminus\m_{V_1}V_d$.
Let us show that $x$ freely generates $V_d$.
For $y\in V_d$, we have $y=ux$ or $x=uy$ for some $u\in V_1$ according respectively as $y/x\in V$ or $x/y\in V$.
If $y\in V_d\setminus\m_{V_1}V_d$, then $u$ must be invertible in $V$, and thus $y\in xV_1$.
If $y\in\m_{V_1}V_d$, then $y=ux$ should hold, and hence $y\in xV_1$.
Now, since $V$ is a graded integral domain, $V_d$ is $V_1$-torsion free.
Hence $V_d$ is freely generated by $x$.

(3) If suffices to show that $\m_V\cap V_d=\m_{V_1}V_d$ for $d\in\Delta_k$.
The inclusion $\m_{V_1}V_d\subseteq m_V\cap V_d$ is clear.
In the proof of (2), we have shown that $V_d/\m_{V_1}V_d$ is, as a vector space over $k_1=V_1/\m_{V_1}$, generated by the residue class of $x$; as it is non-zero (since $d\in\Delta_k$), it has to be of dimension one. 
Hence the surjection $V_d/\m_{V_1}V_d\rightarrow V_d/(\m_V\cap V_d)$ is an isomorphism, which shows that $\m_V\cap V_d=\m_{V_1}V_d$, as desired.
\end{proof}

\begin{dfn}\label{dfn-vsbr-unitelementpart1}{\rm 
Let $V$ be a graded valuation ring for a graded field $K$, and $k=V/\m_k$ the graded residue field of $V$.
We say $V$ is {\it non-degenerate} if the equality $\Delta_k=\Delta_K$ holds.}
\end{dfn}

As a corollary of \ref{prop-vsbr-unitelementpart1}, we have the following result.
\begin{cor}\label{cor-vsbr-unitelementpart1}
In the situation as in {\rm \ref{prop-vsbr-unitelementpart1}}, $V$ is non-degenerate if and only if $\m_{V_1}V=\m_V$.
\end{cor}

\begin{prop}\label{prop-vsbr-unitelementpart2}
Let $V$ be a graded valuation ring for a graded field $K$.
Suppose $V$ is non-degenerate.

{\rm (1)} Any generization $W$ of $V$ is again non-degenerate.

{\rm (2)} The canonical map $V_d\otimes_{V_1}V_{d'}\rightarrow V_{dd'}$ is an isomorphism, whenever $V_d$ and $V_{d'}$ are non-zero.
\end{prop}

\begin{proof}
(1) Any generization $W$ of $V$ is of the form $W=V_{\mathfrak{p}}$ for a homogeneous prime ideal $\mathfrak{p}\subseteq V$.
The graded residue field $k_W$ of $W$ is then the graded fractional field of $V/\mathfrak{p}$.
For $d\in\Delta$, $(k_V)_d\neq 0$ implies $(V/\mathfrak{p})_d\neq 0$, and hence $(k_W)_d\neq 0$.

(2) By \ref{prop-vsbr-unitelementpart1} (2), it suffices to show that 
$$
(V_d/\m_{V_1}V_d)\otimes_{V_1/\m_{V_1}}(V_{d'}/\m_{V_1}V_{d'})\longrightarrow V_{dd'}/\m_{V_1}V_{dd'}
$$
is surjective.
This follows from the equality $\m_{V_1}V=\m_V$ (\ref{cor-vsbr-unitelementpart1}) and that $V/\m_V=k_V$ is a graded field.
\end{proof}

\begin{rem}\label{rem-vsbr-unitelementpart1}{\rm 
Note that a graded valuation ring $V$ for a graded field $K$ is non-degenerate if and only if it enjoys the following properties:
\begin{itemize}
\item[(a)] each $V_d$ for $d\in\Delta_K$ is a free $V_1$-module;
\item[(b)] for any $d,d'\in\Delta_K$, $V_d\otimes_{V_1}V_{d'}\cong V_{dd'}$.
\end{itemize}
Indeed, the `only if' part is already shown in \ref{prop-vsbr-unitelementpart1} (2) and \ref{prop-vsbr-unitelementpart2}(2).
Conversely, if the conditions (a) and (b) hold, then for any $d\in\Delta_K$, one can take the free generator $x$ of $V_d$ as a $V_1$-module.
The condition (b) with $d=d^{-1}$ implies that $x$ is invertible in $V$.
Then, for any $ax\in(\m_V)_d$, since $v(ax)=v(a)>0$, we have $a\in\m_{V_1}$, showing that $(\m_V)_d=\m_{V_1}V_d$.
Hence $V$ is non-degenerate by \ref{cor-vsbr-unitelementpart1}.}
\end{rem}

\begin{prop}\label{prop-vsbr-unitelementpart3}
Let $V$ be a graded valuation ring for a graded field $K$.
Suppose $V$ is non-degenerate.
Then the inclusion $K^{\times}_1\hookrightarrow U(K)$ induces an isomorphism
$$
K^{\times}_1/V^{\times}_1\stackrel{\sim}{\longrightarrow}\Gamma_V.
$$
$($See {\rm \S\ref{subsub-vsbr-valuationgraded}} for the notation.$)$
\end{prop}

Notice that, since $U(V)\cap K^{\times}_1=V^{\times}_1$, the map in question is always injective.
But, as the following proof indicates, non-degeneracy of $V$ is essential for $i$ to be surjective.
\begin{proof}
Denote the map in question $K^{\times}_1/V^{\times}_1\rightarrow\Gamma_V$ by $i$.
To prove that $i$ is surjective, it suffices to show that the mod $U(V)$ class of any non-zero $x\in V_d$ with $d\in\Delta_K$ is contained in the image of $i$.
Since $d\in\Delta_K=\Delta_k$, one can take a free generator $u$ of $V_d$ as a $V_1$-module (see \ref{prop-vsbr-unitelementpart1} (2)).
Note that $u$ is invertible in $V$, that is, $u\in U(V)$.
Then $x/u\in K^{\times}_1$, and its modulo $V^{\times}_1$ class does not depend on the choice of $u$.
Thus $x$ mod $U(V)$ is contained in the image of $i$.
\end{proof}

By {\bf \ref{ch-pre}}.\ref{prop-height2}, \ref{prop-vsbr-genspe}, and \ref{prop-vsbr-unitelementpart3}, we have the following result.
\begin{cor}\label{cor-vsbr-unitelementpart3}
Let $V$ be a graded valuation ring for a graded field $K$.
Suppose $V$ is non-degenerate.

{\rm (1)} The mapping $\mathfrak{p}\mapsto\mathfrak{p}_1=\mathfrak{p}\cap V_1$ gives a bijection from the set of all homogeneous prime ideals of $V$ to the set of prime ideals of $V_1$.

{\rm (2)} The mapping $W_1\mapsto W_1\otimes_{V_1}V=\bigoplus_{d\in\Delta}(W_1\otimes_{V_1}V_d)$ gives a bijection from the set of all valuation rings for $K_1$ containing $V_1$ to the set of all graded valuation rings for $K$ containing $V$.
\end{cor}

\subsubsection{The space of graded valuations}\label{subsub-vsbr-spacegradedvals}
We continue the discussion of graded objects with the grading by a fixed abelian group $\Delta$.
Let $A$ be a graded ring, and $K$ a graded $A$-algebra that is a graded field; notice that we do not assume that the homomorphism $A\rightarrow K$ is injective.
Define
$$
\mathrm{ZR}(K,A)=\bigg\{V\subseteq K\,\bigg|\ \begin{minipage}{15em}{\small 
$V$ is a graded valuation ring for $K$, and is a graded $A$-subalgebra of $K$}
\end{minipage}\bigg\}.
$$
For any graded $A$-subalgebra $B$ of $K$, set
$$
U(B)=\mathrm{ZR}(K,B),
$$
which is a subset of $\mathrm{ZR}(K,A)$.
The subsets of the form $U(B)$ are closed under finite intersection. 
Indeed, if $B_1,B_2$ are graded $A$-subalgebras of $K$, then 
$$
U(B_1)\cap U(B_2)=U(B),\eqno{(\ast)}
$$
where $B$ is the graded $A$-subalgebra of $K$ generated by (homogeneous generators of) $B_1$ and $B_2$.

Note that, in case $\Delta=\{1\}$ and $A$ is a subring of $K$, then $\mathrm{ZR}(K,A)$ is the classical Zariski-Riemann space\index{Zariski-Riemann space!classical Zariski-Riemann space@classical ---} associated to $X=\Spec A$ (see \ref{dfn-classicalclassicalZR} below).
In general, similarly to the classical case, we consider the topology on $\mathrm{ZR}(K,A)$ generated by subsets of the form $U(B)$ for $B$ finite type (viz., generated by finitely many homogenous elements) over $A$.

\begin{thm}[cf.\ \cite{Temkin1}, \S2]\label{thm-vsbr-spacegradedvals}
The topological space $\mathrm{ZR}(K,A)$ is a coherent valuative space\index{valuative!valuative topological space@--- (topological) space}\index{space@space (topological)!valuative topological space@valuative ---} $($see {\rm {\bf \ref{ch-pre}}.\ref{dfn-valuativespace}} for the definition of valuative spaces$)$.
\end{thm}

\begin{proof}
Set $X=\mathrm{ZR}(K,A)$, and consider the set $h(K)$ of homogeneous elements of $K$.
For any non-zero $f\in h(K)$, the two open sets $U(A[f])$ and $U(A[f^{-1}])$ cover $X$.
More precisely, one has
$$
X\setminus U(A[f])=\{V\in U(A[f^{-1}]): f^{-1}\in\m_V\}.
$$

For any $f\in h(K)$, set $C(f)=X\setminus U(A[f])$, and for any family of elements $\{f_{\alpha}\}_{\alpha\in L}$ of $h(K)$, set
$$
C(\{f_{\alpha}\}_{\alpha\in L})=\bigcap_{\alpha\in L}C(f_{\alpha}).
$$
If none of $f_{\alpha}$'s is zero, then let $J$ be the ideal of $B=A[f_{\alpha}:\alpha\in L]$ generated by $f^{-1}_{\alpha}$ for all $\alpha\in L$.
It follows from the above description that 
$$
C(\{f_{\alpha}\}_{\alpha\in L})=\begin{cases}\{V\in X: J\subseteq\m_V\}&\textrm{(if $f_{\alpha}\neq 0$ for all $\alpha\in L$),}\\ \emptyset&\textrm{(otherwise).}\end{cases}
$$
In particular, in view of \ref{prop-vsbr-gradedvaluationrings2}, we deduce that the family of open subsets $\{U(A[f_{\alpha}])\}_{\alpha\in L}$ covers $X$ if and only if $J=B$.

Now let us show that $X$ is quasi-compact.
Take a family $\{C_{\alpha}\}_{\alpha\in L}$ of closed subsets of $X$ with the finite intersection property, that is, $\bigcap_{\alpha\in L'}C_{\alpha}\neq\emptyset$ for any finite subset $L'\subseteq L$.
We want to show that $\bigcap_{\alpha\in L}C_{\alpha}$ is non-empty.
Assume contrary.
As in \cite{Bourb4}, Chap.\ I, \S9, Exer.\ 1, we may reduce to the case where $C_{\alpha}=C(f_{\alpha})$ by $f_{\alpha}\in h(K)$ for $\alpha\in L$.
Then, by what we have seen above, we have $f_{\alpha}\neq 0$ for all $\alpha\in L$, and there is an expression 
$$
1=\sum_{\alpha\in L'}c_{\alpha}f^{-1}_{\alpha}
$$
by a finite subset $L'\subset L$ and $c_{\alpha}\in A[f^{-1}_{\alpha}:\alpha\in L']$.
But this implies $\bigcap_{\alpha\in L'}C_{\alpha}=\emptyset$, contradicting the finite intersection property, which we have assumed. 
Hence it follows that $X$ is quasi-compact.
Moreover, since this implies that each $U(B)$, where $B$ is a finite type graded $A$-subalgebra of $K$, is quasi-compact, we have also shown that $X$ has an open basis consisting of quasi-compact open subsets.

It is then easy to see that $X$ is quasi-separated.
Indeed, to see this, it suffices to show that the intersection of open subsets of the form $U(B_1)\cap U(B_2)$, where $B_1,B_2$ are finite type graded $A$-subalgebras of $K$, is quasi-compact, which follows from what we have seen in $(\ast)$ above.

So far, we have shown that $X$ is a coherent topological space.
Next, let us show that $X$ is sober.
It is easy to see that $X$ is a $\mathrm{T}_0$-space ({\bf \ref{ch-pre}}, \S\ref{subsub-sober}).
Let $Z$ be an irreducible closed subset of $X$.
Let $F_Z$ be the set of all open subsets $U$ of $X$ such that $U\cap Z\neq \emptyset$.
Then, since $Z$ is irreducible, $F_Z$ is a prime filter\index{filter (of lattices)!prime filter of lattices@prime ---} (see {\bf \ref{ch-pre}}, \S\ref{subsub-strcohtopsp}).
Now, referring to the fact that $F_Z$ has to be generated by open subsets of the form $U(B)$ by a finite type graded $A$-subalgebra $B$ of $K$, we define $V$ to be the union of all such $B$'s such that $U(B)\in F_Z$.
Since $U(A[f])\cap Z\neq\emptyset$ for any $f\in B$ with $U(B)\in F_Z$, it follows that $V$ is a graded $A$-subalgebra of $K$.
We claim that $V$ is a graded valuation ring for $K$. 
Indeed, as we have seen above, any non-zero $f\in h(K)$ give an open covering $\{U(A[f]),U(A[f^{-1}])\}$ of $X$, and hence at least one of $f$ and $f^{-1}$ belongs to $V$.
Since the point in $X$ corresponding to $V$ belongs to all open subsets $U$ such that $U\cap Z\neq\emptyset$, it is the unique generic point of $Z$.
Hence we have shown that $X$ is a sober space.

Finally, let us remark that, to show that $X$ is valuative, we only have to invoke \ref{cor-vsbr-genspe}, for it asserts that the set of generizations of a point in $X$ is totally ordered; see {\bf \ref{ch-pre}}.\ref{rem-valuativespace} (1).
\end{proof}
\index{valuation!valuation ring@--- ring!graded valuation ring@graded --- ---|)}

\subsection{Filtered valuations}\label{sub-vsbr-filteredvaluations}
\index{valuation!filtered valuation@filtered ---|(}
Throughout this subsection, we set $\Delta=\R_+$, the multiplicative group of positive real numbers, and all gradings and filtrations are indexed by $\R_+$, often without mentioning it.

\subsubsection{Filtered rings}\label{subsub-vsbr-filteredrings}
In the sequel, we will consider filtered objects indexed by $\R_+$.
For each filtered object $(M,F)$, where $F=\{F_r\}_{r> 0}$ is an increasing (viz., $F_r\subseteq F_{r'}$ for $r\leq r'$) filtration by subobjects of $M$ indexed by $\R_+$, the associated graded object $\mathrm{Gr}_FM$ is defined by
$$
\mathrm{Gr}_FM=\bigoplus_{r\in\R_+}F_r/F_{<r},\quad\textrm{where}\quad F_{<r}=\bigcup_{r'<r}F_{r'}.
$$
For $f\in F_r$, we set
$$
[f]_r=(f\ \mathrm{mod}\ F_{<r}),
$$
which is an element of $\mathrm{Gr}_{F,r}M=F_r/F_{<r}$.

\begin{ntn}\label{ntn-vsbr-filteredrings}{\rm 
For two filtration $F=\{F_r\}_{r>0}$ and $F'=\{F'_r\}_{r>0}$, we will often write
$$
F\subseteq F',
$$
or say ``$F$ is contained in $F'$'', to mean $F_r\subseteq F'_r$ for any $r>0$.}
\end{ntn}

\begin{dfn}[{\rm Cf.\ \cite{NasOys1}, Chap.\ I, 1.1}]\label{dfn-vsbr-filteredfield}{\rm 
(1) An {\it $\R_+$-filtered ring}, or {\it filtered ring} for short, is a pair $(A,F)$ consisting of a (commutative) ring $A$ with $1=1_A$ and a ascending filtration $F$ by additive subgroups satisfying the following conditions:
\begin{itemize}
\item[{\rm (a)}] $F$ is $\R_+$-multiplicative, that is, $1\in F_1$, and $F_r\cdot F_{r'}\subseteq F_{rr'}$ for $r,r'>0$;
\item[{\rm (b)}] $F$ is exhaustive (cf.\ {\bf \ref{ch-pre}}, \S\ref{subsub-topfromfil}), that is, $\bigcup_{r>0}F_r=A$.
\end{itemize}

(2) A {\it morphism} $f\colon (A,F_A)\rightarrow(B,F_B)$ between filtered rings is a ring homomorphism $f\colon A\rightarrow B$ such that $f((F_A)_r)\subseteq (F_B)_r$ for any $r>0$.
In this situation, as usual, we often say that $(B,F_B)$ is an $(A,F_A)$-algebra.}
\end{dfn}

Given a filtered ring $(A,F)$, one can endow $A$ with the linear topology\index{filtration by submodules@filtration (by submodules)!topology defined by a filtration by submodules@topology defined by a ---}\index{topology!topology defined by a filtration@--- defined by a filtration} defined by the filtration $F=\{F_r\}_{r\in\R_+}$ (see {\bf \ref{ch-pre}}, \S\ref{subsub-topologicalringsmodules}).
We call such a topological ring an {\it $\R_+$-linearly topologized ring}.
For an $\R_+$-linearly topologized ring $A$, a {\it filtration of definition} is a multiplicative and exhaustive filtration $F$ on $A$ from which the induced linear topology coincides with the topology of $A$.

It can be shown that:
\begin{itemize}
\item if $F_1,F_2$ are filtration of definition, and $F_1\subseteq F\subseteq F_2$, then $F$ is also a filtration of definition;
\item if $F_1,F_2$ are filtration of definition, then $F_1\cap F_2=\{(F_1)_r\cap(F_2)_r\}_{r\in\R_+}$ is also a filtration of definition.
\end{itemize}

We say that two filtration of definition $F,F'$ are {\it equivalent} if there are positive numbers $c,c'\geq 0$ such that $F_r\subseteq F'_{c'r}$ and $F'_r\subseteq F_{cr}$ hold for any $r>0$.

Let $(A,F)$ be a filtered ring.
We set
$$
F_0=\bigcap_{r>0}F_r.
$$
We say that $(A,F)$ is {\it separated} if $F_0=\{0\}$.

Let $(A,F)$ be a filtered ring, and $I\subseteq A$ an ideal.
Then $A/I$ has the induced filtration $\ovl{F}$ (cf.\ {\bf \ref{ch-pre}}, \S\ref{subsub-topfromfil}) defined by $\ovl{F}_r=(F_r+I)/I$ for $r>0$, which makes the pair $(A/I,\ovl{F})$ an $(A,F)$-algebra.

A {\it filtered multiplicative subset} of $(A,F)$ is a multiplicative subset $S$ of $A$ with the decomposition $S=\bigcup_{r>0}S_r$ such that the following conditions are satisfied:
\begin{itemize}
\item $S_r\subseteq F_r$ for $r>0$;
\item $1\in S_1$, and $S_rS_{r'}\subseteq S_{rr'}$ for $r,r'>0$.
\end{itemize}
In this situation, one has the induced filtration $S^{-1}F$ on $S^{-1}A$ given by
$$
(S^{-1}F)_r=\sum_{r'>0}S^{-1}_{r'}F_{rr'}
$$
for $r>0$, and thus we obtain an $(A,F)$-algebra $(S^{-1}A,S^{-1}F)$.
Note that the canonical map
$$
(\ovl{S})^{-1}\mathrm{Gr}_FA\longrightarrow\mathrm{Gr}_{S^{-1}F}S^{-1}A,
$$
where $\ovl{S}$ is the multiplicative subset of $\mathrm{Gr}_FA$ induced from $S$, is an isomorphism of graded rings.

The {\it completion} of a filtered ring $(A,F)$ is a filtered ring $(\widehat{A},\widehat{F})$ given by
$$
\widehat{A}=\varprojlim_{r>0}A/F_r,\qquad \widehat{F}_r=\varprojlim_{r>r'}F_r/F_{r'}
$$
for $r>0$.
There exists a canonical homomorphism $(A,F)\rightarrow(\widehat{A},\widehat{F})$ of filtered rings.
If this morphism is an isomorphism, we say that $(A,F)$ is {\it complete}.

More generally, if $F'$ is another filtration on $A$ with $F\subseteq F'$, one has the induced filtration $F^{\prime\wedge}_F$ on $\widehat{A}$ given by $(F^{\prime\wedge}_F)_r=\varprojlim_{r>r'}F'_r/F_{r'}$ for $r>0$.

If $(B_1,F_1)$ and $(B_2,F_2)$ are $(A,F)$-algebras, then $B_1\otimes_AB_2$ has the filtration $\til{F}$ defined by 
$$
\til{F}_r=\sum_{st=r}\image((F_1)_s\otimes(F_2)_t\rightarrow B_1\otimes_AB_2)
$$
for any $r>0$, which gives rise to a filtered ring $(B_1\otimes_AB_2,\til{F})$ sitting in the canonical commutative square
$$
\xymatrix{(B_1,F_1)\ar[r]&(B_1\otimes_AB_2,\til{F})\\ (A,F)\ar[u]\ar[r]&(B_2,F_2)\rlap{.}\ar[u]}
$$
The filtration $\til{F}$ is called the {\it tensor-product filtration}.
The completion of $(B_1\otimes_AB_2,\til{F})$, denoted by 
$$
(B_1\widehat{\otimes}_AB_2,F_1\widehat{\otimes}_FF_2),
$$
is called the {\it complete tensor product} of $(B_1,F_1)$ and $(B_2,F_2)$ over $(A,F)$.

\subsubsection{Filtrations and seminorms}\label{subsub-vsbr-seminorms}
Let $(A,F)$ be a filtered ring.
Define a new filtration $F^+$ by 
$$
F^+_r=\bigcap_{r'>r}F_{r'}
$$
for any $r>0$.
The resulting object $(A,F^+)$ is a filtered ring, which is a filtered $(A,F)$-algebra.
Moreover, we have
$$
F^+_{<r}=F_{<r}
$$
for any $r>0$, and hence the induced morphism $\mathrm{Gr}_FA\rightarrow\mathrm{Gr}_{F^+}A$ is always injective, viz., we may regard $\mathrm{Gr}_FA$ as a subring of $\mathrm{Gr}_{F^+}A$.

Let $(A,F)$ be a filtered ring.
Since the filtration $F$ is exhaustive, one can define a function $\nu\colon A\rightarrow\R_{\geq 0}$ by
$$
\nu(x)=\inf\{r: x\in F_r\}
$$
for $x\in A$.
Due to the fact that $F$ is multiplicative, one can easily show that the mapping $\nu$ is a non-archimedean seminorm\index{seminorm} (\S\ref{sub-vsbr-banachrings}) on the ring $A$.
It is a norm if and only if $(A,F)$ is separated, and is the trivial norm if and only if $F$ is trivial; here, we say that the filtration $F$ is {\it trivial} if $F_r=\{0\}$ for $0<r<1$ and $F_r=A$ for $r\geq 1$.

Conversely, if we are given a ring $A$ and a non-archimedean seminorm $\nu$ on $A$, one can define the filtration $F_{\nu}$ by
$$
(F_{\nu})_r=\nu^{-1}([0,r])=\{f\in A: \nu(f)\leq r\}
$$
for $r>0$, which makes the pair $(A,F_{\nu})$ a filtered ring.
It is easy to see that, starting from a filtered ring $(A,F)$, the resulting filtration $F_{\nu}$, where $\nu$ is the associated seminorm of $F$, coincides with $F^+$, and that this gives rise to the bijection between the following sets:
\begin{itemize}
\item the set of all filtrations $F$ on $A$ such that $F=F^+$;
\item the set of all non-archimedean seminorms $\nu$ on $A$.
\end{itemize}
In this way, one sees that the notion of filtered rings gives a refinement of the notion of non-archimedean seminormed rings.
Let us say that a filtration $F$ in a filtered ring $(A,F)$ is {\it of seminorm type} if $F=F^+$.

If $F$ is a filtration of definition of an $\R_+$-linearly topologized ring $A$, then so is $F^+$, and $F$ and $F^+$ are equivalent.
Indeed, we have $F\subseteq F^+$ and $(F^+)_r\subseteq F_{cr}$ for any $c>1$ and $r>0$.
Thus, it follows that equivalence classes of filtrations of definition are in one-to-one correspondence with the equivalence classes of seminorms on $A$.

Note that the condition $f((F_A)_r)\subseteq (F_B)_r$ ($r>0$) for filtered morphisms as in \ref{dfn-vsbr-filteredfield} (2) implies $\nu_B(f(x))\leq\nu_A(x)$ for any $x\in A$, where $\nu_A,\nu_B$ is the associated seminorms of $F_A, F_B$, respectively.
In particular, filtered homomorphisms induce bounded homomorphisms between the corresponding seminormed rings.

Note also that the completion $(\widehat{A},\widehat{F})$ of $(A,F)$ is isomorphic to the completion of $(A,F^+)$, and hence $\widehat{A}$ is a Banach ring, isomorphic to the completion of $A$ with respect to the seminorm $\nu$ associated to $F$. (Caution: the induced filtration $\widehat{F}$ may not be of seminorm type.)

For a filtered ring $(A,F_A)$ and a surjective ring homomorphism $A\rightarrow B$, consider the induced filtration $F_B$ on $B$ as in \S\ref{subsub-vsbr-filteredrings}.
If $\nu_A$ is the seminorm corresponding to $F^+_A$, then the seminorm on $B$ corresponding to $F^+_B$ is the residue seminorm induced from $\nu_A$ (\S\ref{sub-vsbr-banachrings}).

Let $(B_1,F_1)$ and $(B_2,F_2)$ be $(A,F)$-algebras, and consider the tensor-product filtration $\til{F}$.
Then the seminorm on $B_1\otimes_AB_2$ corresponding to $\til{F}^+$ is nothing but the tensor-product seminorm (\S\ref{sub-vsbr-banachrings}) from the seminorms $\nu_1$, $\nu_2$, and $\nu$ corresponding to $F^+_1$, $F^+_2$, and $F^+$, respectively.
Hence, in the situation as in the end of \S\ref{subsub-vsbr-filteredrings}, the complete tensor product $B_1\widehat{\otimes}_AB_2$ coincides with the complete tensor product of the corresponding seminormed rings (as in \S\ref{sub-vsbr-banachrings}).

\subsubsection{Filtered polynomial and power series algebras}\label{subsub-vsbr-filteredpolynomials}
Let $X$ be a set, and fix a function $w\colon X\rightarrow\R_+$, which we call a {\it weight function}.
Let $\Z[X]$ be the $\Z$-polynomial ring with the indeterminacies in $X$, where the indeterminacy corresponding to $x\in X$ is denoted by $X_x$, and $M_X$ the free (additive) commutative monoid with the basis $X$, which is identified with the set of all monomials by the correspondence
$$
M_X\ni\alpha=\sum_{x\in X}m_x[x]\quad\longleftrightarrow\quad X^{\alpha}=\prod_{x\in X}(X_x)^{m_x}\in\Z[X].
$$
We extend the function $w$ to $M_X$, hence on the set of all monomials in $\Z[X]$, in such a way that $w(X^{\alpha+\beta})=w(X^{\alpha})w(X^{\beta})$ holds for any $\alpha,\beta\in M_X$.
For any $f=\sum_{\alpha}a_{\alpha}X^{\alpha}\in\Z[X]$, we define the {\it weight} of $f$ by
$$
w(f)=\begin{cases}\sup\{w(X^{\alpha}): a_{\alpha}\neq 0\}&(f\neq 0),\\ 0&(f=0),\end{cases}
$$
which induces a filtration $F_w$ on $\Z[X]$ defined by
$$
(F_w)_r=\{f\in\Z[X]: w(f)\leq r\}
$$
for any $r>0$.

For a filtered ring $(A,F)$, consider the tensor product $A[X]=A\otimes_{\Z}\Z[X]$ with the tensor product filtration (\S\ref{subsub-vsbr-filteredrings}), also denoted by $F_w$.
Explicitly, $f=\sum_{\alpha}a_{\alpha}X^{\alpha}\in A[X]$ belong to $(F_w)_r$ if and only if $a_{\alpha}\in F_{r/w(X^{\alpha})}$ for any $\alpha$.
The associated seminorm, which corresponds to $F^+_w$ by the correspondence in \S\ref{subsub-vsbr-seminorms}, is nothing but the Gauss seminorm $\|\cdot\|_{\mathrm{Gauss}}$, given by
$$
\|f\|_{\mathrm{Gauss}}=\max\{w(X^{\alpha})\|a_{\alpha}\|: \alpha\in M_X\}
$$
for $f=\sum_{\alpha}a_{\alpha}X^{\alpha}\in A[X]$, where $\|\cdot\|$ denotes the seminorm on $A$ associated to the filtration $F$.

Note that, by construction, we have
$$
(A[X],F_w)\cong\varinjlim_{X'\subseteq X}(A[X'],F_{w|_{X'}}),
$$
where $X'$ runs through all finite subsets of $X$.
Note also that we have an isomorphism of graded rings
$$
\mathrm{Gr}_{F_w}A[X]\stackrel{\sim}{\longrightarrow}(\mathrm{Gr}_FA)[X],
$$
where the grading on $(\mathrm{Gr}_FA)[X]$ is defined by the grading of $\mathrm{Gr}_FA$ and the weight, by $aX^{\alpha}$ ($a\in F_r$) mapped to $[a]_rX^{\alpha}$.

If $(A,F)$ is complete, then the completion of $(A[X], F_w)$ is called the {\it weighted power series algebra} over $(A,F)$, and is denoted by 
$$
A\dl (X,w)\dr.
$$
As an $A$-algebra, $A\dl (X,w)\dr$ is isomorphic to the completion of $A[X]$ with respect to the Gauss seminorm $\|\cdot\|_{\mathrm{Gauss}}$.

In the sequel, these materials will be considered often in the situation as $X=\{1,2,\ldots,n\}$ and $w(i)=r_i>0$ for $i=1,2,\ldots,n$.
In this situation, as usual, the weighted power series algebra $A\dl (X,w)\dr$ will be denoted by 
\begin{equation*}
\begin{split}
A\dl r^{-1}X\dr&=A\dl r^{-1}_1X_1,\ldots,r^{-1}_nX_n\dr\\
&=\bigg\{\sum_{m\in\N^n}a_mX^m\in A[\![X]\!]\,\bigg|\, r^m\|a_m\|_A\rightarrow 0\ \textrm{as}\ |m|\rightarrow\infty\bigg\},
\end{split}
\end{equation*}
where $r^m=r^{m_1}_1\cdots r^{m_n}_n$, $X^m=X^{m_1}_1\cdots X^{m_n}_n$ and $|m|=m_1+\cdots+m_n$ for $m=(m_1,\ldots,m_n)\in\N^n$.
This algebra will play an important role from \S\ref{subsub-vsbr-Rfinitetype} below.

\subsubsection{Filtered valuation fields}\label{subsub-vsbr-filteredvaluationfields}
\index{valuation!filtered valuation@filtered ---!filtered valuation field@--- --- field|(}
\begin{dfn}\label{dfn-vsbr-filteredfields}{\rm 
A {\it filtered field} is a separated filtered ring $(K,F)$ with $K$ being a field (in the usual sense).}
\end{dfn}

\begin{dfn}\label{dfn-vsbr-filteredvaluationfield11}{\rm 
A {\it filtered valuation field} is a filtered field $(K,F)$ that satisfies the following conditions:
\begin{itemize}
\item[(a)] $\mathrm{Gr}_{F^+}K$ is a graded field (\S\ref{subsub-vsbr-gradedrings});
\item[(b)] $\mathrm{Gr}_FK$ is a graded valuation ring for $\mathrm{Gr}_{F^+}K$ (\ref{dfn-vsbr-gradedvaluationrings} (1)).
\end{itemize}}
\end{dfn}

Note that, as we have seen in \S\ref{subsub-vsbr-seminorms}, $\mathrm{Gr}_FK$ is, in general, regarded as a graded subring of $\mathrm{Gr}_{F^+}K$.
Note also that, in (b) above, we allow the case $\mathrm{Gr}_FK=\mathrm{Gr}_{F^+}K$, that is, $F=F^+$. 
In this case, we say that the filtered valuation field $(K,F)$ is {\it of maximal type}.
If $(K,F)$ is a filtered valuation field of maximal type, then $F=F^+$ is of seminorm type, that is, $F$ corresponds to a seminorm $\nu\colon K\rightarrow\R_{\geq 0}$ on $K$ (\S\ref{subsub-vsbr-seminorms}).
Since the filtration $F$ is separated, and $K$ is a field, one sees that the seminorm $\nu$ is a valuation (absolute value).
Note that $\nu$ may be a trivial one, that is, $\nu(x)=1$ for all non-zero $x\in K$.

Conversely, as we saw in \S\ref{subsub-vsbr-seminorms}, any non-archimedean valued field $(K,\nu)$ gives rise to a filtered valuation field $(K,F_{\nu})$ of maximal type.
Thus we see that filtered valuation fields of maximal type are virtually the same objects as non-archimedean valued fields.
The notion of filtered valuation fields in general gives, therefore, a refinement of the notion of non-archimedean valued fields.
In general, for a filtered valuation field $(K,F)$, $(K,F^+)$ is a filtered valuation field of maximal type, called the {\it associated maximal type} of $(K,F)$.
The valuation $\nu$ corresponding to the maximal type will be called the {\it associated absolute value} of the filtered valuation field $(K,F)$.

\begin{prop}\label{prop-vsbr-filteredvaluationfield1}
If $(K,F)$ is a filtered valuation field, then the seminorm $\nu$ corresponding to the filtration $F^+$ is a non-archimedean absolute value on $K$.
Conversely, given a non-archimedean absolute value $\nu$ on $K$ and a graded valuation subring $G$ for $\mathrm{Gr}_{F_{\nu}}K$, there exists a unique filtration $F$ on $K$ that makes the pair $(K,F)$ a filtered valuation field such that $G=\mathrm{Gr}_FK$ and that $F_{\nu}=F^+$.
\end{prop}

\begin{proof}
The first assertion is easy to see, as we have already discussed above.
Let $\nu$ be a non-archimedean absolute value on $K$, and $G$ a graded valuation ring for $\mathrm{Gr}_{F_{\nu}}K$.
Define the filtration $F=\{F_r\}_{r\geq 0}$ so that, for each $r>0$, $F_r$ is the preimage of $G_r$ under the map 
$$
(F_{\nu})_{r}\longrightarrow \mathrm{Gr}_{F_\nu,r}K=(F_{\nu})_{r}/(F_{\nu})_{<r}.
$$
Then, clearly, $(K,F)$ is a filtered valuation field with $\mathrm{Gr}_FK=G$.
One sees easily that $F^+=F_{\nu}$.
The uniqueness is clear.
\end{proof}

\begin{cor}\label{cor-vsbr-filteredvaluationfield1}
Let $(K,\nu)$ be a non-archimedean valued field with the absolute value $\nu\colon K\rightarrow\R_{\geq 0}$, and $\widehat{K}$ the completion of $K$ with respect to $\nu$.
Then there exists a natural bijection between the set of all filtrations $F$ on $K$ that make the resulting pair $(K,F)$ a filtered valuation field with the associated absolute value $\nu$, and the set of all filtrations $\til{F}$ on $\widehat{K}$ that make the pair $(\widehat{K},\til{F})$ a filtered valuation field with the associated absolute value equal to the extension $\widehat{\nu}$ of $\nu$ on $\widehat{K}$.
\end{cor}

\begin{proof}
This follows from \ref{prop-vsbr-filteredvaluationfield1} and the isomorphism $\mathrm{Gr}_{F_v}K\cong\mathrm{Gr}_{F_{\widehat{v}}}\widehat{K}$.
\end{proof}

\begin{rem}\label{rem-vsbr-filteredvaluationfield1}{\rm 
Note that, if $(K,F)$ is a filtered valuation field, then $F_1$ is a valuation ring for $K$ in the usual sense.
Indeed, $F_1/F_{<1}$ is a valuation ring for the field $(F_{\nu})_1/(F_{\nu})_{<1}$, where $\nu$ is the associated absolute value (see \S\ref{subsub-vsbr-gradedvaluationrings}).
On the other hand, $(F_{\nu})_1=\{x\in K: \nu(x)\leq 1\}$ is clearly a valuation ring for $K$, with residue field $(F_{\nu})_1/(F_{\nu})_{<1}$.
Hence, by the patching argument ({\bf \ref{ch-pre}}.\ref{prop-composition2}), $F_1$ is a valuation ring.}
\end{rem}

Based on this, we call $F_1$ the {\it associated valuation ring} of the filtered valuation field $(K,F)$.
It is a valuation ring for the field $K$ (in the usual sense).

Proposition \ref{prop-vsbr-filteredvaluationfield1} suggests that a filtered valuation field $(K',F')$ should be said to dominate $(K,F)$ by a morphism $(K,F)\rightarrow(K',F')$ when the absolute value $\nu_{F'}$ extends $\nu_F$ and $\mathrm{Gr}_{F'}K$ dominates $\mathrm{Gr}_FK$ in the sense as in \S\ref{subsub-vsbr-gradedvaluationrings}.
This last condition can be boiled down to the following one.
\begin{dfn}\label{dfn-vsbr-filteredvaluationfielddominate}{\rm 
Let $(K,F)$ and $(K',F')$ be filtered valuation field.
We say that $(K',F')$ {\it dominates} $(K,F)$ by a morphism $(K,F)\rightarrow(K',F')$ if $F'_r\cap K=F_r$ for any $r>0$.}
\end{dfn}

\begin{dfn}\label{dfn-vsbr-filteredvaluationfield3}{\rm 
Let $(K,F)$ is a filtered valuation field, and $\nu=\nu_F$ the associated absolute value on $K$.
The {\it height} of $(K,F)$, denoted by $\mathrm{ht}(K,F)$, is the pair $(\mathrm{ht}(\nu),\mathrm{ht}(\mathrm{Gr}_FK))$.}
\end{dfn}

Here, $\mathrm{ht}(\nu)$ denotes the height of the corresponding valuation ring $(F_{\nu})_1=\{x\in K:\nu(x)\leq 1\}$, which is $0$ or $1$, according as $\nu$ is trivial or non-trivial.

\begin{dfn}\label{dfn-vsbr-filteredvaluationfield4}{\rm 
Let $K$ be a field, and suppose we have two filtrations $F_1$ and $F_2$ on $K$ such that $(K,F_i)$ ($i=1,2$) are filtered valuation fields with the same maximal type.
We say $F_1$ is a {\it specialization} of $F_2$, or $F_2$ is a {\it generization} of $F_1$, if $F_1\subseteq F_2$.}
\end{dfn}

Notice that the `same maximal type' condition is equivalent to that these filtrations induce the same absolute value $\nu$ on $K$.
Thus $F_1$ is a specialization of $F_2$ if an only if the valuation ring $\mathrm{Gr}_{F_1}K$ is a specialization of the valuation ring $\mathrm{Gr}_{F_2}K$.
In particular, the set of all generization of a fixed $(K,F)$ is a totally ordered with respect to the inclusion order (see \ref{cor-vsbr-genspe}).
\index{valuation!filtered valuation@filtered ---!filtered valuation field@--- --- field|)}

\subsubsection{Filtered valuation via valuation}\label{subsub-vsbr-filvalviavaluation}
Let us first remark that, when we consider absolute values $\nu\colon K\rightarrow\R_{\geq 0}$ as a valuation (of height $0$ or $1$), the ordering of $\R_+$ is the reverse of the standard order, hence is reverse to the index-ordering for the filtrations $F=\{F_r\}_{r>0}$; in fact, $\R_+$ as a value target group is the one isomorphic to the standard $\R$ via
$$
(\R_+,\leq^{\mathrm{rev}})\stackrel{\sim}{\longrightarrow}(\R,\leq),\qquad x\longmapsto\log_c x
$$
for some $0<c<1$.

\begin{dfn}\label{dfn-vsbr-filvalviavaluation1}{\rm 
Let $K$ be a field.

(1) An {\it $\R_+$-valuation}\index{valuation!Rplus valuation@$\R_+$ ---} $v$ of $K$ is a valuation of the form 
$$
v\colon K\longrightarrow (\R_+\times\Gamma)\cup\{(0,+\infty)\},
$$
where $\Gamma$ is a totally ordered commutative group, and $\R_+\times\Gamma$ is endowed with the lexicographic order.
The induced valuation $v^{\R_+}\colon K\rightarrow\R_{\geq 0}$ (resp.\ $v^{\mathrm{Gr}}\colon K\rightarrow\Gamma\cup\{+\infty\}$) given by the first (resp.\ second) projection is called the {\it absolute value} (resp.\ {\it graded part}) associated to $v$.

(2) Let $(K,v)$ and $(K',v')$ be two fields with $\R_+$-valuations with the value target groups $\Gamma_v=\R_+\times\Gamma$ and $\Gamma_{v'}=\R_+\times\Gamma'$, respectively.
A {\it homomorphism} $(\varphi,\phi)\colon (K,v)\rightarrow(K',v')$ of $\R_+$-valuation fields consists of a field homomorphism $\varphi\colon K\rightarrow K'$ and a homomorphism of ordered groups $\phi\colon\Gamma_v\rightarrow\Gamma_{v'}$ such that the following conditions are satisfied:
\begin{itemize}
\item[(a)] $\phi$ is of the form $\phi=\id_{\R_+}\times\psi$, where $\psi\colon\Gamma\rightarrow\Gamma'$ is a homomorphism of ordered groups;
\item[(b)] $v'\circ\varphi=\phi\circ v$.
\end{itemize}
We say that $v'$ {\it dominates} $v$ if, for $f\in K$, $v^{\mathrm{Gr}}(f)>0$ implies $v^{\prime\mathrm{Gr}}(\varphi(f))>0$.

(3) Two $\R_+$-valuations $v_1$ and $v_2$ on $K$ are said to be {\it equivalent} if there exists an $\R_+$-valuation $v$ on $K$ and homomorphisms $(K,v)\rightarrow (K_i,v_i)$ by which $(K_i,v_i)$ dominates $(K,v)$ for $i=1,2$.}
\end{dfn}

Let $(K,F)$ be a filtered valuation field.
One has the associated absolute value $\nu=\nu_F\colon K\rightarrow\R_{\geq 0}$, and the graded valuation ring $\mathrm{Gr}_FK$ for $\mathrm{Gr}_{F^+}K$, hence the value group $\Gamma$ of the associated graded valuation $v_F$ (see \S\ref{subsub-vsbr-valuationgraded}).
Then one can define an $\R_+$-valuation on $K$ with the value target group $\R_+\times\Gamma$ by
$$
K\ni f\longmapsto v(f)=(\nu(f),v_F([f]_{\nu(f)}))\in(\R_+\times\Gamma)\cup\{(0,+\infty)\}.
$$
Notice that, from this $\R_+$-valuation $v$, the filtration $F$ is recovered by 
$$
F_r=\{0\}\cup v^{-1}((0,r)\times\Gamma)\cup v^{-1}(\{r\}\times\Gamma_{\geq 0})
$$
for $r>0$.
By this, one can show the following result.
\begin{prop}\label{prop-vsbr-filvalviavaluation1}
The above construction gives rise to a bijection from the set of all filtrations $F$ on $K$ that make the pair $(K,F)$ a filtered valuation field, to the set of all equivalence classes of $\R_+$-valuations on $K$.
\end{prop}

\begin{rem}\label{rem-vsbr-filvalviavaluation1}{\rm 
$\R_+$-valuations can be regarded as so-called {\it reified valuations} by K.\ Kedlaya\index{Kedlaya, K.} \cite{Kedlaya} with an extra structure.
According to \cite{Kedlaya}, 5.1, a reified valuation on a field $K$ is a valuation $v\colon K\rightarrow\Gamma\cup\{0\}$, where $\Gamma$ is written multiplicatively, with a {\it reification}, viz., an order preserving homomorphism $r\colon \R_+\rightarrow\Gamma$.
In terms of this, an $\R_+$-valuation on $K$ is interpreted as a reified valuation $(v,r)$ on $K$ with an order-preserving homomorphism $p\colon\Gamma\rightarrow\R_+$ such that $p\circ r=\id_{\R_+}$.}
\end{rem}

\subsubsection{Non-degenerate filtered valuations}\label{subsub-vsbr-filteredvaluationnondeg}
\begin{prop}\label{prop-vsbr-filteredvaluationnondeg0}
Let $(K,F)$ be a filtered valuation field, $\nu=\nu_F$ the associated absolute value on $K$, and $v\colon K\rightarrow(\R_+\times\Gamma)\cup\{(0,+\infty)\}$ the associated $\R_+$-valuation.
Then the following conditions are equivalent:
\begin{itemize}
\item[{\rm (a)}] the graded valuation ring $\mathrm{Gr}_FK$ is non-degenerate $($see {\rm \ref{dfn-vsbr-unitelementpart1}}$)$;
\item[{\rm (b)}] for any $r\in\nu(K^{\times})$, there exists $f_r\in K$ such that $v(f_r)=(r,0)$;
\item[{\rm (c)}] the value group $v(K^{\times})$ of $v$ is isomorphic to the product $\nu(K^{\times})\times v^{\mathrm{Gr}}(K^{\times})$ $($with the lexicographic order$)$ as a totally ordered group.
\end{itemize}
\end{prop}

\begin{proof}
(a) $\Rightarrow$ (b): Let $r\in\nu(K^{\times})$.
Then $\mathrm{Gr}_{F^+,r}K\neq\{0\}$.
Since $\mathrm{Gr}_FK$ is non-degenerate, there exists an element in $\mathrm{Gr}_{F,r}K$, which is invertible in $\mathrm{Gr}_FK$ (see \ref{prop-vsbr-unitelementpart1} (1)).
This means that there exists $f_r\in K$ such that $v(f_r)=(r,0)$.

(b) $\Rightarrow$ (c): By assumption, we have a splitting $\nu(K^{\times})\hookrightarrow v(K^{\times})$ by $r\mapsto (r,0)$, from which the assertion follows by standard argument.

(c) $\Rightarrow$ (b) is clear.

(b) $\Rightarrow$ (a): For any $r\in\nu(K^{\times})$, there exists $f_r\in K$ such that $\nu(f)=r$ and $[f]_r$ is invertible in $\mathrm{Gr}_FK$, from which the non-degeneracy of $\mathrm{Gr}_FK$ follows by \ref{prop-vsbr-unitelementpart1} (1).
\end{proof}

\begin{dfn}\label{dfn-vsbr-filteredvaluationnondeg}{\rm 
A filtered valuation field $(K,F)$ is said to be {\it non-degenerate} if it satisfies one, hence all, of the conditions in {\rm \ref{prop-vsbr-filteredvaluationnondeg0}}.}
\end{dfn}

Note that there exist filtered valuation fields that are not non-degenerate.
Indeed, they can be obtained by $\R_+$-valuation $v\colon K\rightarrow(\R_+\times\Gamma)\cup\{(0,+\infty)\}$ with the value group not being a product of totally ordered abelian groups; cf.\ {\bf \ref{ch-pre}}.\ref{thm-estimate}.

\begin{prop}\label{prop-vsbr-filteredvaluationnondeg1}
Let $(K,F)$ be a non-degenerate filtered valuation field.
Then any generization of $(K,F)$ is again non-degenerate.
Moreover, there exists a natural bijection between the set of all generizations of $(K,F)$ and the set of all generizations $($in the usual sense$)$ of $F_1$.
\end{prop}

\begin{proof}
As we have remarked at the end of \S\ref{subsub-vsbr-filteredvaluationfields}, $F'$ is a generization of $F$ if and only if $\mathrm{Gr}_{F'}K$ is a generization of $\mathrm{Gr}_FK$.
Hence the last assertion of the proposition follows from \ref{cor-vsbr-unitelementpart3} (2) and what we have seen in \ref{rem-vsbr-filteredvaluationfield1}.
The first assertion follows from \ref{prop-vsbr-unitelementpart2} (1).
\end{proof}

\begin{prop}\label{prop-vsbr-filteredvaluationnondeg2}
Let $(K,F)$ be a non-degenerate filtered valuation field, and $\nu$ the associated absolute value.

{\rm (1)} For any $r\in\Gamma_{\nu}=\nu(K^{\times})$, $F_r$ is a free $F_1$-module of rank one.

{\rm (2)} For any $r,r'\in\Gamma_{\nu}=\nu(K^{\times})$, the map by multiplication $F_r\otimes_{F_1}F_{r'}\rightarrow F_{rr'}$ is an isomorphism.
\end{prop}

\begin{proof}
Take $f_r$ as in \ref{prop-vsbr-filteredvaluationnondeg0} (b).
For any $x\in F_r$, we have $f^{-1}_rx\in F_1$, and hence $F_r=F_1f_r$.
Thus we have shown the first assertion, and also the second assertion follows.
\end{proof}

\begin{rem}\label{rem-vsbr-filteredvaluationnondeg}{\rm 
In summary, a non-degenerate filtered valuation field amounts to the data as follows (cf.\ \ref{rem-vsbr-unitelementpart1}):
\begin{itemize}
\item a field $K$ equipped with a non-archimedean absolute value $\nu\colon K\rightarrow\R_{\geq 0}$,
\item for any $r\geq 0$, a submodule $F_r$ of $(F_{\nu})_r=\{x\in K: \nu(x)\leq r\}$, such that 
\begin{itemize}
\item[(a)] for any $r\geq 0$, $F_r=\bigcup_{r'\in\Gamma_{\nu}\cap [0,r]}F_{r'}$;
\item[(b)] $F_1$ is a valuation ring for $K$, and the valuation ring $(F_{\nu})_1$ of $\nu$ is the generization of $F_1$;
\item[(c)] for $r\in\Gamma_{\nu}=\nu(K^{\times})$, $F_r$ is a free $F_1$-submodule of $(F_{\nu})_r$ of rank one, and $F_r\otimes_{F_1}(F_{\nu})_1\cong(F_{\nu})_r$;
\item[(d)] for $r,r'\in\Gamma_{\nu}$, $F_r\otimes_{F_1}F_{r'}\cong F_{rr'}$ compatibly with the multiplication map $(F_{\nu})_r\otimes_{(F_{\nu})_1}(F_{\nu})_{r'}\cong(F_{\nu})_{rr'}$.
\end{itemize}
\end{itemize}}
\end{rem}

Finally, let us show that, in a certain situation, the formation of filtered valuations reduces to that of usual valuations.
Let $(K,F)$ be a non-degenerate filtered valuation field with the associated absolute value $\nu_K$, and $L/K$ a field extension with an absolute value $\nu_L$ that extends $\nu_K$.
We assume $\nu_K(K^{\times})=\nu_L(L^{\times})$.
Take a valuation ring $V$ for $L$ such that
\begin{itemize}
\item[(a)] $V$ contains $F_1$;
\item[(b)] $\nu_L$ is a generization of $V$.
\end{itemize}
Then the filtration $F_L$ of $L$ defined by
$$
(F_L)_r=F_r\cdot V\cong F_r\otimes_{F_1}V
$$
for $r\in\R_+$ gives a filtered valuation field $(L,F_L)$ such that $(F_L)_1=V$.
Since each $(F_L)_r$ for $r\in\nu_L(L^{\times})$ is free of rank one over $V$, and $(F_L)_{rr'}=(F_L)_r\cdot(F_L)_{r'}$, $(L,F_L)$ is non-degenerate.

\begin{prop}\label{prop-vsbr-examplefilteredvaluations}
Let $(K,F)$ be a non-degenerate filtered valuation field with the associated absolute value $\nu_K$, and $L/K$ a field extension with an absolute value $\nu_L$ that extends $\nu_K$.
Suppose $\nu_K(K^{\times})=\nu_L(L^{\times})$.
Then any filtration $\til{F}$ of $L$ that makes the pair $(L,\til{F})$ a filtered valuation field over $(K,F)$ is obtained as above.
\end{prop}

\begin{proof}
By assumption, for any $r\in\nu_L(L^{\times})$, $\til{F}_r$ contains $F_r$, and $F_r\cdot \til{F}_1\subseteq\til{F}_r$.
Hence $(L,\til{F})$ is a generization of $(L,F_L)$ obtained as above from $V=\til{F}_1$.
By \ref{prop-vsbr-filteredvaluationnondeg1}, both of them are non-degenerate.
Moreover, they have the same `$F_1$'-part.
Hence, again by \ref{prop-vsbr-filteredvaluationnondeg1}, they coincide with each other, which is what we wanted to show.
\end{proof}

\subsubsection{Examples of filtered valuations}\label{subsub-vsbr-examplefilteredvaluations}
Let $(K,F)$ be a filtered valuation field, 
$$
v\colon K\longrightarrow(\R_+\times\Gamma)\cup\{(0,+\infty)\}
$$
the associated $\R_+$-valuation\index{valuation!Rplus valuation@$\R_+$ ---}, and $\nu=v^{\R_+}$ the associated absolute value.
Let $X$ be a set with a weight function $w\colon X\rightarrow\R_+$, and consider the polynomial algebra $(K[X],F_w)$ as in \S\ref{subsub-vsbr-filteredpolynomials}.
The seminorm corresponding to $F^+_w$ is the Gauss norm $\|\cdot\|_{\mathrm{Gauss}}$, and we have an isomorphism
$$
\mathrm{Gr}_{F^+_w}(K[X])\stackrel{\sim}{\longrightarrow}(\mathrm{Gr}_{F_{\nu}}K)[X], 
$$
mapping $aX^{\alpha}$ ($\nu(a)=r$) to $[a]_rX^{\alpha}$.
For any non-zero $f=\sum_{\alpha\in M_X}a_{\alpha}X^{\alpha}\in K[X]$, we have
$$
[f]=\sum_{\alpha\in M_X}[a_{\alpha}]_{w(\alpha)^{-1}r}\cdot X^{\alpha},
$$
where $[f]=[f]_r$ with $r=\|f\|_{\mathrm{Gauss}}$.
We have
\begin{itemize}
\item[(a)] $[fg]=[f][g]$ for non-zero $f,g\in K[X]$;
\item[(b)] if $\|f\|_{\mathrm{Gauss}}=\|g\|_{\mathrm{Gauss}}$ and $[f]+[g]\neq 0$, then $[f+g]=[f]+[g]$ and the function 
$$
v_{X,w}(f)=\begin{cases}(\|f\|_{\mathrm{Gauss}},v^{\mathrm{Gr}}(\Cont^{\mathrm{Gr}}([f])))&(f\neq 0),\\ (0,+\infty)&(f=0),\end{cases}
$$
defined on $K[X]$ extends to an $\R_+$-valuation of the quotient field $K(X)$ with the value target group $\R_+\times\Gamma$.
\end{itemize}
Here $\Cont^{\mathrm{Gr}}(H)$ for homogeneous $H\in(\mathrm{Gr}_{F_{\nu}}K)[X]$ is the graded $\mathrm{Gr}_FK$-submodule of $\mathrm{Gr}_{F_{\nu}}K$ generated by the classes of coefficients of $H$, and $v^{\mathrm{Gr}}(\Cont^{\mathrm{Gr}}(H))$ is defined to be $v^{\mathrm{Gr}}(a)$, where $a$ is a generator of $\Cont^{\mathrm{Gr}}(H)$.
(One checks easily, as in the classical case, that $\Cont^{\mathrm{Gr}}(HH')=\Cont^{\mathrm{Gr}}(H)\cdot\Cont^{\mathrm{Gr}}(H')$ for homogeneous $H,H'$.)
The $\R_+$-valuation $v_{X,w}$ will be called the {\it Gauss valuation} of $K(X)$, and often denoted by $v_{\mathrm{Gauss}}$.
The absolute value part of $v_{\mathrm{Gauss}}$ is the Gauss norm $\|\cdot\|_{\mathrm{Gauss}}$.

Note that the equality (a), which simply follows from that $(\mathrm{Gr}_{F_{\nu}}K)[X]$ is a graded integral domain, implies the Gauss lemma for the Gauss seminorm.

\begin{exa}\label{exa-vsbr-Kstandard}{\rm
We apply the above construction to an extreme case $X=\R_+$ and $w=\id_{\R_+}$.
In this situation, the Gauss norm $\|\cdot\|_{\mathrm{Gauss}}$ on $K(X)$ has the value group equal to the whole $\R_+$.
Moreover, if we denote by $X_r$ the indeterminacy corresponding to $r\in\R_+$, then we have $v_{\mathrm{Gauss}}(X_r)=(r,0)$, which shows that the filtered valuation field given by $(K(X),v_{\mathrm{Gauss}})$ is non-degenerate (due to \ref{prop-vsbr-filteredvaluationnondeg0}).
We denote the filtered valuation field over $(K,F)$ thus obtained by
$$
K^{\mathrm{st}}=(K(\R_+),v_{\mathrm{Gauss}}),
$$
and call it the {\it standard extension}.
The standard extension is a non-degenerate filtered valuation field over any given $(K,F)$, having $\R_+$ as its value group of the absolute value part.
Moreover, the formation $K\mapsto K^{\mathrm{st}}$ is functorial.}
\end{exa}

Note that the existence of standard extension also shows that any filtered valuation field can be embedded into a non-degenerate isometric extension.
\index{valuation!filtered valuation@filtered ---|)}

\subsection{Valuative spectrum of non-archimedean Banach rings}\label{sub-vsbr-valuativespectrum}
In this subsection, we define and discuss valuative spectra of non-archinedean commutative Banach rings.

\subsubsection{Gelfand-Berkovich spectrum}\label{subsub-berkovichanalyticspacespectrum}
\begin{dfn}[{\cite{Berk1}, \S1.2}]\label{dfn-berkovichspectrum}{\rm 
Let $(A,\|\cdot\|_A)$ be a commutative Banach ring.
The {\it Gelfand-Berkovich spectrum}\index{Berkovich spectrum} $\mathscr{M}(A)$ of $A$ is
\begin{itemize}
\item the set of all bounded multiplicative seminorms\index{seminorm!multiplicative seminorm@multiplicative ---} on $A$, 
\item with the weakest topology such that all functions of the form
$$
|f(\cdot)|\colon\mathscr{M}(A)\longrightarrow\R_{\geq 0},\quad z=|\cdot|_z\longmapsto |f|_z
$$
for $f\in A$ are continuous.
\end{itemize}}
\end{dfn}

Here, a seminorm $|\cdot|\colon A\rightarrow\R_{\geq 0}$ is {\it bounded} if there exists $C>0$ such that $|f|\leq C\|f\|_A$ for any $f\in A$.

As indicated above, for any point $z\in\mathscr{M}(A)$, the corresponding seminorm on $A$ is often denoted by $|\cdot|_z$.
The kernel $\ker(|\cdot|_z)=\{f\in A:|f|_z=0\}$ is a closed prime ideal of $A$, and the seminorm $|\cdot|_z$ induces a valuation on the integral domain $A/\ker(|\cdot|_z)$, and hence on its fractional field.
The Hausdorff completion of this field with respect to the valuation is denoted by $\mathcal{H}(z)$, called the {\it complete residue field at $z$}.
As usual, the image of an element $f\in A$ in $\mathcal{H}(z)$ is denoted by $f(z)$; thus, $|f(z)|$ denotes the value $|f|_z$.

The topological space $\mathscr{M}(A)$ is compact, and is non-empty if $A\neq\{0\}$ (\cite{Berk1}, 1.2.1).
Note that $f\in A$ is invertible in $A$ if and only if the function $z\mapsto|f(z)|$ is everywhere non-zero on $\mathscr{M}(A)$.
If $A$ is non-archimedean, then $\mathcal{H}(z)$ for any $z\in\mathscr{M}(A)$ is non-archimedean. 

In the sequel, all seminorms and norms are assumed to be non-archimdean, unless otherwise clearly stated.
Note that the norm $|\cdot|_z$ on the complete residue field $\mathcal{H}(z)$ may be a trivial one.

For $f\in A$, we will denote by $\|f\|_{\Sp}$ the {\it spectral seminorm}
$$
\|f\|_{\Sp}=\sup_{z\in\mathscr{M}(A)}|f(z)|,
$$
which is equal to the {\it spectral redius} 
$$
\rho(f)=\inf_{n\geq 1}\|f^n\|_A^{\frac{1}{n}}=\lim_{n\rightarrow\infty}\|f^n\|_A^{\frac{1}{n}};
$$
see \cite{Berk1}, \S1.3.
By this we can define an $\R_+$-filtration $F^{\Sp}_A$ on $A$, called the {\it spectral filtration}, as
$$
(F^{\Sp}_A)_r=\{f\in A:\|f\|_{\Sp}\leq r\},
$$
by which we obtain an $\R_+$-filtered ring (\ref{dfn-vsbr-filteredfield}) $(A,F^{\Sp}_A)$.
It follows that, for any $z\in\mathscr{M}(A)$, the map $f\mapsto f(z)$ maps $(F^{\Sp}_A)_r$ to $(F_{|\cdot|_z})_r$ for any $r>0$, that is, the residue map $A\rightarrow\mathcal{H}(z)$ is a morphism of filtered rings (\ref{dfn-vsbr-filteredfield}).
Note that the filtered field $\mathcal{H}(z)$ with the spectral filtration coincides with $(\mathcal{H}(z),F_{|\cdot|_z})$, hence is a filtered valuation field of maximal type (\ref{dfn-vsbr-filteredvaluationfield11}).

In general, any bounded homomorphism $f\colon A\rightarrow B$ induces the morphism $f\colon (A,F^{\Sp}_A)\rightarrow(B,F^{\Sp}_B)$ of filtered rings, and hence the morphism $\mathrm{Gr}_{F^{\Sp}_A}A\rightarrow\mathrm{Gr}_{F^{\Sp}_B}B$ of graded rings.

For a commutative Banach ring $A$, denote by $A^{\Sp}$ the completion of $A$ with respect to the spectral filtration $F^{\Sp}_A$ (the so-called {\it uniform completion}).
The canonical morphism $A\rightarrow A^{\Sp}$ induces an isomorphism
$$
\mathscr{M}(A^{\Sp})\stackrel{\sim}{\longrightarrow}\mathscr{M}(A), 
$$
and we have $(A^{\Sp})^{\Sp}=A^{\Sp}$.
A Banach ring $A$ is said to be {\it uniform}, or a {\it Banach function ring}, if $A\rightarrow A^{\Sp}$ is a bounded isomorphism.

Let $(A,F)$ be a filtered ring, and consider filtered homomorphisms of the form $(A,F)\rightarrow(K,F_{\nu})$, where $K$ is a filtered valuation field of maximal type. 
For two such homomorphisms $(A,F)\rightarrow(K,F_{\nu})$ and $(A,F)\rightarrow(K',F_{\nu'})$, consider the relation by domination, that is, the relation given by a commutative diagram
$$
\xymatrix@R-3ex{&(K,F_{\nu})\ar[dd]\\ (A,F)\ar[ur]\ar[dr]\\ &(K',F_{\nu'})\rlap{,}}
$$
where by the down-arrow $(K',F_{\nu'})$ dominates $(K,F)$ (see \ref{dfn-vsbr-filteredvaluationfielddominate}).
Consider the equivalence relation generated by this relation, and set
$$
\mathscr{M}(A,F)=\bigg\{\begin{minipage}{25em}{\small equivalence class of filtered homomorphisms of the form $(A,F)\rightarrow(K,F_{\nu})$ to a filtered valuation field of maximal type}\end{minipage}\bigg\}.
$$
Consider the completion $\widehat{A}=\widehat{A}_F$ of $A$ with respect to the filtration $F$.
Then $\widehat{A}$ is a Banach ring by the norm induced from the seminorm $\nu_F$ associated to $F$.
For any $(A,F)\rightarrow (K,F_{\nu})$, the norm $\nu$ on $K$ induces a seminorm on $\widehat{A}$, which is clearly multiplicative and bounded.
Hence we have the natural map
$$
\mathscr{M}(A,F)\longrightarrow\mathscr{M}(\widehat{A}).\eqno{(\ast)}
$$
\begin{thm}\label{thm-vsbr-GromovBerkovichspectrum}
The map $(\ast)$ is a bijection.
In particular, if $\mathrm{Gr}_FA\neq\{0\}$, then $\mathscr{M}(A,F)$ is non-empty.
\end{thm}

\begin{proof}
Let $z\in\mathscr{M}(\widehat{A})$ be given by $\alpha\colon\widehat{A}\rightarrow\mathcal{H}(z)$. 
We want to show that the composition 
$$
A\longrightarrow\widehat{A}\stackrel{\alpha}{\longrightarrow}\mathcal{H}(z)
$$
extends to a filtered map $(A,F)\rightarrow(\mathcal{H}(z),F_{|\cdot|_z})$.
Since $A\rightarrow\widehat{A}$ obviously respects the filtrations, we may assume $A=\widehat{A}$.
Consider the uniform completion $A\rightarrow A^{\Sp}$, which factorize $A\rightarrow\mathcal{H}(z)$.
The morphism $A\rightarrow A^{\Sp}$ respects filtration, for we have $\|f\|_A\geq\|f\|_{A^{\Sp}}$.
The morphism $A^{\Sp}\rightarrow\mathcal{H}(z)$ also respects the filtration, since the norms on both sides are power-multiplicative, and hence $\|f\|_{A^{\Sp}}\geq |f(z)|$.
Hence we have shown the first assertion.
Since $\mathrm{Gr}_FA\neq\{0\}$ implies $\widehat{A}\neq\{0\}$, the second assertion follows from \cite{Berk1}, 1.2.1.
\end{proof}

Here let us give a few technical corollaries for our later purpose.
\begin{cor}\label{cor-vsbr-GromovBerkovichspectrum1}
Let $(A,F)$ be a filtered ring.
For a graded prime ideal $\mathfrak{p}$ of $\mathrm{Gr}_FA$, there exists a filtered valuation field $(K,F_K)$ over $(A,F)$ such that the graded valuation ring $\mathrm{Gr}_{F_K}K$ dominates $(\mathrm{Gr}_FA)_{\mathfrak{p}}$.
\end{cor}

\begin{proof}
Let $\ovl{S}_{\mathfrak{p}}$ be the graded multiplicative system corresponding to $\mathfrak{p}$ (cf.\ \S\ref{subsub-vsbr-gradedrings}), and $S_{\mathfrak{p}}$ the filtered multiplicative subset of $A$ defined by 
$$
S_{\mathfrak{p}}=\bigcup_{r>0}(S_{\mathfrak{p}})_r, 
$$
where $(S_{\mathfrak{p}})_r$ is the preimage of $(\ovl{S}_{\mathfrak{p}})_r$ under the map $F_r\rightarrow F_r/F_{<r}=\mathrm{Gr}_{F,r}A$.
The filtered localization $A'=S^{-1}_{\mathfrak{p}}A$ has $(\mathrm{Gr}_FA)_{\mathfrak{p}}$ as its associated graded ring, which is non-zero.
By \ref{thm-vsbr-GromovBerkovichspectrum}, there exists a filtered homomorphism $(A',F')\rightarrow K$, where $K=(K,F_{\nu})$ is a filtered valuation field of maximal type.
The image of $\mathrm{Gr}_{F'}A'$ in $\mathrm{Gr}_{F_{\nu}}K$ is a local graded subring.
Now by \ref{prop-vsbr-gradedvaluationrings2} and \ref{prop-vsbr-filteredvaluationfield1}, there exists a filtered valuation $(K,F_K)$ over $(A',F')$ such that $\mathrm{Gr}_{F_K}K$ dominates $(\mathrm{Gr}_FA)_{\mathfrak{p}}$.
\end{proof}

\begin{cor}\label{cor-vsbr-GromovBerkovichspectrum2}
Let $(K,F_K)$ be a filtered valuation field, and $(M,F_M)$ and $(L,F_L)$ filtered valuation fields that dominate $(K,F_K)$.
Then there exists a filtered valuation field $(N,F_N)$ sitting in the commutative diagram
$$
\xymatrix{(M,F_M)\ar[r]&(N,F_N)\\ (K,F_K)\ar[u]\ar[r]&(L,F_L)\ar[u]}
$$
consisting of dominating filtered homomorphisms.
\end{cor}

\begin{proof}
Let us denote by $\nu_K,\nu_L,\nu_M$ the associated absolute values of $K,L,M$, respectively, and set $V_K=(F_K)_1$, $V_L=(F_L)_1$, and $V_M=(F_M)_1$.
Replacing $K,L,M$ by $K^{\mathrm{st}},L^{\mathrm{st}},M^{\mathrm{st}}$ constructed as in \ref{exa-vsbr-Kstandard}, respectively, we may assume $K,L,M$ are all non-degenerate, and having the whole $\R_+$ as their value groups.
By \ref{prop-vsbr-examplefilteredvaluations}, we have
$$
(F_L)_r=(F_K)_r\otimes_{V_K}V_L,\quad (F_M)_r=(F_K)_r\otimes_{V_K}V_M
$$
for all $r>0$.
Set $A=L\otimes_KM$ with the tensor product filtration $F_A$.
Since $(F_A)_r=(F_K)_r\otimes_{V_K}(V_L\otimes_{V_K}V_M)$ for $r>0$, the canonical surjection 
$$
\mu\colon\mathrm{Gr}_{F_L}L\otimes_{\mathrm{Gr}_{F_M}M}\mathrm{Gr}_{F_M}M\longrightarrow\mathrm{Gr}_{F_A}A
$$
is an isomorphism.

Now, take a graded prime ideal $\mathfrak{p}$ of $\mathrm{Gr}_{F_A}A$ lying above the graded maximal ideals of $\mathrm{Gr}_{F_L}L$ and of $\mathrm{Gr}_{F_M}M$.
By \ref{cor-vsbr-GromovBerkovichspectrum1}, one has a filtered homomorphism $(A,F_A)\rightarrow (N,F_N)$ to a filtered valuation field such that $\mathrm{Gr}_{F_N}N$ dominates $(\mathrm{Gr}_{F_A}A)_{\mathfrak{p}}$.
This $(N,F_N)$ gives the commutative square as asserted.
\end{proof}

\subsubsection{$\R_+$-finite type algebras}\label{subsub-vsbr-Rfinitetype}
Let $A$ be a Banach ring.
Recall that, for an $n$-tuple $r=(r_1,\ldots,r_n)$ of positive real numbers, the $\R_+$-power series ring of variables $T=(T_1,\ldots,T_n)$ with coefficients in $A$ is
\begin{equation*}
\begin{split}
A\dl r^{-1}T\dr&=A\dl r^{-1}_1T_1,\ldots,r^{-1}_nT_n\dr\\
&=\bigg\{\sum_{m\in\N^n}a_mT^m\in A[\![T]\!]\,\bigg|\, r^m\|a_m\|_A\rightarrow 0\ \textrm{as}\ |m|\rightarrow\infty\bigg\},
\end{split}
\end{equation*}
where $r^m=r^{m_1}_1\cdots r^{m_n}_n$, $T^m=T^{m_1}_1\cdots T^{m_n}_n$ and $|m|=m_1+\cdots+m_n$ for $m=(m_1,\ldots,m_n)\in\N^n$.
It is a Banach ring with the Gauss norm
$$
\|f\|_{A\dl r^{-1}T\dr}=\sup_{m\in\N^n}r^m\|a_m\|_A
$$
for $f=\sum_{m\in\N_n}a_mT^m$; see \S\ref{subsub-vsbr-filteredpolynomials}.
By spectral seminorm formula (\ref{thm-spectralseminormformula} (3)), we have
$$
\|f\|_{A\dl r^{-1}T\dr,\Sp}=\sup_{n\in\N^n}r^m\|a_m\|_{A,\Sp}
$$
for $f=\sum_{m\in\N_n}a_mT^m$ (Exercise \ref{exer-vsbr-spectralseminormformula}), and the filtered ring $(A\dl r^{-1}T\dr,F^{\Sp}_{A\dl r^{-1}T\dr})$ is isomorphic to the completion of $(A[X],F^{\Sp}_w)$, where $F^{\Sp}_w$ is the tensor product filtration of $F^{\Sp}_A$ on $A$ and $F_w$ on $\Z[X]$; see \S\ref{subsub-vsbr-filteredpolynomials}.
Passing to the associated graded algebras,  we have an isomorphism
$$
\mathrm{Gr}_{F^{\Sp}_{A\dl r^{-1}T\dr}}A\dl r^{-1}T\dr\stackrel{\sim}{\longrightarrow}(\mathrm{Gr}_{F^{\Sp}_A}A)[X],\quad [T_i]_{r_i}\longmapsto X_i,
$$
where each $X_i$ on the right-hand side is of degree $r_i$ ($i=1,\ldots,n$).

\begin{dfn}\label{dfn-sbr-Rfinitetype}{\rm 
A Banach $A$-algebra $B$ is said to be {\it of $\R_+$-finite type} if it is isomorphic to a Banach $A$-algebra of the form
$$
A\dl r^{-1}T\dr/\mathfrak{a}
$$
by a bounded isomorphism, where $\mathfrak{a}\subseteq A\dl r^{-1}T\dr$ is a closed ideal.
If one can take such an isomorphism with $r=(1,\ldots,1)$, then we say $B$ is {\it of finite type}.}
\end{dfn}

An important class of $\R_+$-finite type $A$-algebra is as follows.
Let $f=(f_0,f_1,\ldots,f_n)$ be an $(n+1)$-tuple of elements of $A$ such that $(f_0,f_1,\ldots,f_n)=A$, and $r=(r_1,\ldots,r_n)$ an $n$-tuple of positive real numbers.
The {\it $\R_+$-rational algebra over $A$ corresponding $(f,r)$} is 
$$
A\Big\langle\!\!\Big\langle r^{-1}_1\frac{f_1}{f_0},\ldots,r^{-1}_n\frac{f_n}{f_0}\Big\rangle\!\!\Big\rangle=A\dl r^{-1}T\dr/\ovl{(f_0T_1-f_1,\ldots,f_0T_n-f_n)},
$$
which is a Banach $A$-algebra by the residue norm.
An $\R_+$-finite type $A$-algebra isomorphic to the one of this form will be called an {\it $\R_+$-rational localization} of $A$.

Note that 
$$
\mathbf{U}_0(f,r)=\mathscr{M}\Big(A\Big\langle\!\!\Big\langle r^{-1}_1\frac{f_1}{f_0},\ldots,r^{-1}_n\frac{f_n}{f_0}\Big\rangle\!\!\Big\rangle\Big)
$$
is identified with the closed subset 
$$
\{z\in\mathscr{M}(A): |f_i(z)|\leq r_i\cdot |f_0(z)|\}
$$
of $\mathscr{M}(A)$, the so-called {\it $\R_+$-rational subdomain} of $\mathscr{M}(A)$.

\subsubsection{Integrally closed filtrations}\label{subsub-vsbr-integrallyclosedfiltrations}
\begin{dfn}\label{dfn-vsbr-integrallyclosedfiltrationsr1}{\rm 
Let $(A,F)$ be a filtered ring, and $r\in\R_+$.
We say $x\in A$ is {\it $r$-integral} over $F$ if there exists an equality
$$
x^n+a_1x^{n-1}+\cdots+a_{n-1}x+a_n=0\eqno{(\ast)}
$$
for some $n\geq 1$ and $a_i\in F_{r^i}$ for $i=1,\ldots,n$.}
\end{dfn}

\begin{lem}\label{lem-vsbr-integrallyclosedfiltrations1}
Let $(A,F)$ be a filtered ring, and $x\in A$.
Then $x$ is $r$-integral $(r>0)$ if and only if $x\in F^+_r$ and $[x]_r\in\mathrm{Gr}_{F^+}A$ is integral over the graded subring $\mathrm{Gr}_FA$.
\end{lem}

\begin{proof}
To show the `only if' part, suppose $x$ is $r$-integral, and the an equality $(\ast)$ above holds by $a_i\in F_{r^i}$ ($i=1,\ldots,n$).
Let $\|\cdot\|$ be the seminorm corresponding to $F^+$ (see \S\ref{subsub-vsbr-seminorms}).
Since $\|a_i\|\leq r^i$ for $i=1,\ldots,n$, we have
$$
\|x^n\|\leq\sup_{i=1,\ldots,n}r^i\|x^{n-i}\|,
$$
from which $\|x\|\leq r$ follows.
Hence $x\in F^+_r$.
It is then clear that $[x]_r$ is integral over $\mathrm{Gr}_FA$.

To show the converse, choose $a_i\in F_{r^i}$ for $i=1,\ldots,n$ such that 
$$
[x]^n_r+[a_1]_r[x]^{n-1}_r+\cdots+[a_{n-1}]_{r^{n-1}}[x]_r+[a_n]_{r^n}=0.
$$
Then $\alpha=x^n+a_1x^{n-1}+\cdots+a_{n-1}x+a_n$ belongs to $F_{<r^n}$, and the equality
$$
x^n+a_1x^{n-1}+\cdots+a_{n-1}x+a_n-\alpha=0
$$
implies that $x$ is $r$-integral.
\end{proof}

\begin{dfn}\label{dfn-vsbr-integrallyclosedfiltrations1}{\rm 
Let $(A,F)$ be a filtered ring.

(1) For $r>0$, define 
$$
F^{\int}_r=\{x\in A: \textrm{$x$ is $r$-integral over $F$}\}.
$$
Then $F^{\int}=\{F^{\int}_r\}_{r>0}$ is, by \ref{lem-vsbr-integrallyclosedfiltrations1}, an $\R_+$-multiplicative filtration on $A$ containing $F$.
The filtration $F$ is said to be {\it integrally closed} if $F=F^{\int}$.

(2) Let $F'$ be an $\R_+$-multiplicative (see \ref{dfn-vsbr-filteredfield} (a)) filtration of $A$ containing $F$.
The {\it integral closure of $F$ in $F'$}, denoted by $F^{\int, F'}$, is defined by
$$
F^{\int,F'}_r=F'_r\cap F^{\int}_r
$$
for any $r>0$.}
\end{dfn}

\begin{prop}\label{prop-vsbr-integrallyclosedfiltrations1}
{\rm (1)} If $(A,F)$ is a filtered ring with $F$ integrally closed, then for any filtered multiplicative subset $S=\bigcup_{r>0}S_r$ of $A$ $($see {\rm \S\ref{subsub-vsbr-filteredrings})}, the filtration $S^{-1}F$ on $S^{-1}A$ is integrally closed.

{\rm (2)} Let $(A,F_0)$ be a filtered ring, and $F$ an integrally closed $\R_+$-multiplicative filtration containing $F_0$.
Then the completion $F^{\wedge}_{F_0}$ of $F$ with respect to $F_0$ $($see {\rm \S\ref{subsub-vsbr-filteredrings})} is integrally closed.
\end{prop}

\begin{proof}
(1) can be shown similarly to the classical case.
To show (2), let $x$ be an element of $(F^{\wedge}_{F_0})^+_r$ satisfying $x^n+a_1x^{n-1}+\cdots+a_{n-1}x+a_n=0$ for $a_i\in(F^{\wedge}_{F_0})_{r^i}$ ($i=1,\ldots,n$).
Take $z\in A$ such that $x-z\in(F^{\wedge}_{F_0})_{<r}$, and $b_i\in F_{r^i}$ such that $a_i-b_i\in(F^{\wedge}_{F_0})_{<r^i}$ ($i=1,\ldots,n$).
Then we have 
$$
[z]^n_r+[b_1]_r[z]^{n-1}_r+\cdots+[b_{n-1}]_{r^{n-1}}[z]_r+[b_n]_{r^n}=0
$$
in $\mathrm{Gr}_{F^{\wedge}_{F_0}}A^{\wedge}_{F_0}=\mathrm{Gr}_FA$.
Since $F$ is integrally closed, $z\in F_r$ by \ref{lem-vsbr-integrallyclosedfiltrations1}, whence $x\in (F^{\wedge}_{F_0})_r$, as desired.
\end{proof}

\begin{prop}\label{prop-vsbr-integrallyclosedfiltrations2}
Let $(A,F_0)$ be a filtered ring, and $F$ an $\R_+$-multiplicative filtration containing $F_0$.
Consider the set $\mathscr{M}(A,F_0)$ as in {\rm \ref{subsub-berkovichanalyticspacespectrum}}, and the associated spectral filtration $F^{\Sp}_{(A,F_0)}$.

{\rm (1)} If $F$ is contained in $F^{\Sp}_{(A,F_0)}$, then $F^{\int}\subseteq F^{\Sp}_{(A,F_0)}$.

{\rm (2)} For any $r>0$, any element of $(F^{\Sp}_{(A,F_0)})_{<r}$ is $r$-integral over $F$.
In particular, $F^+=F^{\Sp}_{(A,F_0)}$ if $F$ is integrally closed and $F\subseteq F^{\Sp}_{(A,F_0)}$.

{\rm (3)} If $\til{F}$ is an $\R_+$-multiplicative filtration containing $F$ such that $\til{F}^+=F^{\Sp}_{(A,F_0)}$, then 
$$
(F^{\int,\til{F}})_r=\{x\in\til{F}_r: \textrm{$[x]_r$ is integral over $\mathrm{Gr}_{F,r}A$}\}\eqno{(\ast\ast)}
$$
for $r>0$.
\end{prop}

\begin{proof}
(1) Take a point $z\colon (A,F_0)\rightarrow K_z$ of $\mathscr{M}(A,F_0)$, where $K_z$ is a filtered valuation field of maximal type.
By assumption, the image of $F$ is contained in $F_{K_z}$.
Observe that the homomorphism $(A,F)\rightarrow(K_z,F_{K_z})$ factors through $(A,F^+)$.
Since $r$-integral elements belong to $F^+_r$ for any $r>0$ (\ref{lem-vsbr-integrallyclosedfiltrations1}), it follows that the image of $F^{\int}$ is contained in $F_{K_z}$, from which the assertion follows.

(2) Let $\nu$ be the seminorm on $A$ corresponding to $F^+_0$. 
If $x\in(F^{\Sp}_{(A,F_0)})_{<r}$, then $\nu(x^N)<r^N$ for some $N\geq 1$ (see \S\ref{subsub-berkovichanalyticspacespectrum}).
Then $\alpha=x^n$ is contained in $(F_0)_{r^n}$, and $x^n-\alpha=0$ shows that $x$ is $r$-integral over $F_0$, hence over $F$.
If $F$ is integrally closed and $F\subseteq F^{\Sp}_{(A,F_0)}$, then we have $F_{<r}=(F^{\Sp}_{(A,F_0)})_{<r}$ for any $r>0$, which implies $F^+=F^{\Sp}_{(A,F_0)}$.

(3) Let $x$ be an element in the right-hand set of $(\ast\ast)$.
As in the proof of \ref{lem-vsbr-integrallyclosedfiltrations1}, one has $a_i\in F_{r^i}$ for $i=1,\ldots,n$ such that $\alpha=x^n+a_1x^{n-1}+\cdots+a_{n-1}x+a_n$ belongs to $F_{<r^n}$.
Since $F_{<r^n}=(F^{\Sp}_{(A,F_0)})_{<r^n}$, there exists $N\geq 1$ such that $\beta=\alpha^N\in F_{<r^{nN}}$ by an argument similar to that in the proof of (2).
Then $(x^n+a_1x^{n-1}+\cdots+a_{n-1}x+a_n)^N-\beta=0$ shows that $x\in (F^{\int,\til{F}})_r$.
The other inclusion is easy to see.
\end{proof}

\begin{thm}\label{thm-vsbr-filtrationintegral}
Let $(A,F_0)$ be a filtered ring, and $F$ and $\til{F}$ $\R_+$-multiplicative filtrations on $A$ such that $F_0\subseteq F\subseteq\til{F}$ and $\til{F}^+=F^{\Sp}_{(A,F_0)}$.
Then 
$$
(F^{\int,\til{F}})_r=\bigg\{x\in\til{F}_r\,\bigg|\,\begin{minipage}{20em}{\small for any morphism $(\mathcal{A},F)\rightarrow (K,F_K)$ to a filtered valuation field, the image of $x$ is contained in $(F_K)_r$}\end{minipage}\bigg\}
$$
for $r>0$.
\end{thm}

\begin{proof}
First note that the filtration $F_K$ in a filtered valuation field $(K,F_K)$ is integrally closed, due to \ref{lem-vsbr-integrallyclosedfiltrations1} and \ref{prop-vsbr-gradedvaluationrings1}.
By this, one sees that $(F^{\int,\til{F}})_r$ is contained in the right-hand set.
To show the other inclusion, take $x\in\til{F}_r$ such that $[x]_r$ is not integral over $\mathrm{Gr}_F\mathcal{A}$ in $\mathrm{Gr}_{\til{F}}\mathcal{A}$.
By \cite{Bourb1}, Chap.\ VI, \S1.2, Lemma 1, $(\mathrm{Gr}_{\til{F}}\mathcal{A})_{[x]_r}$ is non-zero, and $[x]^{-1}_r$ generates a proper graded ideal of $(\mathrm{Gr}_F\mathcal{A})[[x]^{-1}_r]$.

Let $S=\bigcup_{r>0}S_r$ be the filtered multiplicative subset (\S\ref{subsub-vsbr-filteredrings}) generated by $x$, such that $x\in S_r$ and $A_x=S^{-1}A$ with the filtration $F_x$ induced from $F$.
We have $\mathrm{Gr}_{F_x}A_x\cong(\mathrm{Gr}_FA)_{[x]_r}$.
Since the graded ideal of $\mathrm{Gr}_{F_x}A_x$ generated by $[x]^{-1}_r$ is a proper ideal, by \ref{cor-vsbr-GromovBerkovichspectrum1}, there exists a filtered valuation field $(K,F_K)$ over $(A_x,F_x)$ such that the graded maximal ideal of $\mathrm{Gr}_{F_K}K$ contains the image of $[x]^{-1}_r$, hence the image of $x$ is not contained in $(F_K)_r$, since, if it is, then the image of $x$ in $\mathrm{Gr}_{F_K}K$ is invertible.
\end{proof}

\begin{cor}\label{cor-vsbr-integrallyclosedfiltrations1}
Let $(A,F_0)$ be a filtered ring, $F$ an $\R_+$-multiplicative filtration on $A$ such that $F_0\subseteq F\subseteq F^{\Sp}_{(A,F_0)}$.
Suppose $(A,F)$ is integrally closed.
Then the filtration $F_w$ in the weighted polynomial algebra $(A[X],F_w)$ as in {\rm \S\ref{subsub-vsbr-filteredpolynomials}} is integrally closed.
Similarly, the filtration $\widehat{F}_w$ in the weighted power series algebra $(A\dl (X,w)\dr,\widehat{F}_w)$, the completion of $(A[X],F_w)$ with respect to the filtration $F_{0,w}$ induced from $F_0$, is integrally closed.
\end{cor}

\begin{proof}
Let us show the assertion under the assumption that $(A,F)$ is a filtered valuation field; the general case can be reduced to this case by \ref{thm-vsbr-filtrationintegral}.
Let $v$ be the corresponding $\R_+$-valuation, and $v_{\mathrm{Gauss}}$ the Gauss valuation of $A[X]$ (see \S\ref{subsub-vsbr-examplefilteredvaluations}).
Take $f\in A[X]$ such that $f^n+a_1f^{n-1}+\cdots+a_0=0$ for $a_i\in (F_w)_{r^i}$.
From this and $v_{\mathrm{Gauss}}(a_i)\geq (r^i,0)$, we have $v_{\mathrm{Gauss}}(f)\geq (r,0)$, hence $f\in (F_w)_r$.
This shows that $F_w$ is integrally closed.

Finally, $\widehat{F}_w$ in $(A\dl (X,w)\dr,\widehat{F}_w)$ is integrally closed due to \ref{prop-vsbr-integrallyclosedfiltrations1} (2).
\end{proof}

\subsubsection{Power bounded filtration}\label{subsub-vsbr-powerboundedfiltration}
\begin{dfn}\label{dfn-vsbr-powerboundedfiltration1}{\rm
(1) Let $(A,\|\cdot\|)$ be a seminormed ring, and $r>0$ a positive real number.
An element $f\in A$ is said to be {\it $r$-power bounded} if
$$
\sup_{n\geq 1}r^{-n}\|f^n\|<+\infty.
$$

(2) Let $(A,F)$ be a filtered ring, and $\|\cdot\|$ the seminorm corresponding to the filtration $F^+$ (see \S\ref{subsub-vsbr-seminorms}).
Then we say $f\in A$ is {\it $r$-power bounded} if it is $r$-power bounded as an element of the seminormed ring $(A,\|\cdot\|)$.}
\end{dfn}

Similarly to the classical notion of power-boundedness, the $r$-power-boundedness is closely related to the universal property of power series algebras.
\begin{prop}\label{prop-vsbr-powerboundedfiltration1}
Let $A$ be a Banach ring.
For $X=(X_1,\ldots,X_n)$, $r=(r_1,\ldots,r_n)$, and a Banach ring $B$, there exists a natural bijection between the following sets:
\begin{itemize}
\item[{\rm (a)}] the set of all bounded homomorphisms of the form $A\dl r^{-1}X\dr\rightarrow B$;
\item[{\rm (b)}] the set of all $n$-tuples $f=(f_1,\ldots,f_n)$ of elements of $B$, where $f_i$ is $r_i$-power bounded for $i=1,\ldots,n$.
\end{itemize}
\end{prop}

This follows from the following lemma.
\begin{lem}\label{lem-vsbr-powerboundedfiltration1}
For $f\in A$ and $r>0$, the following conditions are equivalent to each other:
\begin{itemize}
\item[{\rm (a)}] $f$ is $r$-power bounded;
\item[{\rm (b)}] the canonical homomorphism
$$
A\longrightarrow A\dl r^{-1}T\dr/\overline{(T-f)}
$$
is a bounded isomorphism.
\end{itemize}
\end{lem}

\begin{proof}
Suppose (b) holds, and set $B=A\dl r^{-1}T\dr/\overline{(T-f)}$.
The residue norm on $B$ induced from the Gauss norm on $A\dl r^{-1}T\dr$ is equivalent to the norm on $A$.
Since $\|T^n\|_B\leq r^n$ for $n\geq 1$, $f$ is $r$-power bounded.
Conversely, suppose $f$ is $r$-power bounded, and consider the homomorphism
$$
A\dl r^{-1}S\dr\longrightarrow A\dl r^{-1}T\dr,\quad S\longmapsto T-f,
$$
and observe that it is a bounded isomorphism with the inverse arrow by $T\mapsto S+f$.
Passing to the quotient by $\ovl{(S)}$, we get the bounded isomorphism
$$
A\stackrel{\sim}{\longrightarrow}A\dl r^{-1}S\dr/\ovl{(S)}\stackrel{\sim}{\longrightarrow}A\dl r^{-1}T\dr/\overline{(T-f)},
$$
as desired.
\end{proof}

The following is a corollary of the proposition.
\begin{cor}\label{cor-vsbr-powerboundedfiltration1}
Let $A$ be a Banach ring, $f=(f_0,\ldots,f_n)$ an $(n+1)$-tuple of elements of $A$, which generate the unit ideal $A$, and $r=(r_1,\ldots,r_n)$ an $n$-tuple of positive real numbers.
Then, for a Banach $A$-algebra $B$, there exists a natural bijection between the following sets:
\begin{itemize}
\item[{\rm (a)}] the set of all bounded homomorphisms of the form
$$
A\Big\langle\!\!\Big\langle r^{-1}_1\frac{f_1}{f_0},\ldots,r^{-1}_n\frac{f_n}{f_0}\Big\rangle\!\!\Big\rangle\longrightarrow B;
$$
\item[{\rm (b)}] the set of all $n$-tuples $z=(z_1,\ldots,z_n)$ of elements of $B$ such that, for $i=1,\ldots,n$,  $z_i$ is $r_i$-power bounded and $f_0z_i=f_i$.
\end{itemize}
In particular, any bounded $A$-algebra homomorphism between $\R_+$-rational localization over $A$ $($see {\rm \S\ref{subsub-vsbr-Rfinitetype})} is an epimorphism in the category of Banach rings.
\end{cor}

\begin{dfn}[{\rm Power bounded filtration}]\label{dfn-vsbr-powerboundedfiltration2}{\rm 
For a seminormed ring $(A,\|\cdot\|)$, the {\it power bounded filtration} on $A$, denote by $F^o_A$, is a filtration $\{F^o_A\}_{r\in\R_+}$ defined by
$$
(F^o_A)_r=\{x\in A: \textrm{$x$ is $r$-power bounded}\}.
$$
For a filtered ring $(A,F=F_A)$, the power bounded filtration on $A$ is defined similarly, in reference to the $r$-power-boundedness as in \ref{dfn-vsbr-powerboundedfiltration1} (2), and is denoted by $F^o=F^o_A$.}
\end{dfn}

Note that the power bounded filtration $F^o$ is an exhaustive and $\R_+$-multiplicative filtration, hence defining a filtered ring $(A,F^o)$ (see \ref{dfn-vsbr-filteredfield}).
Moreover, for a filtered ring $(A,F)$, we have
$$
F^o\subseteq F^{\Sp},\quad F^{\Sp}_{<r}\subseteq F^o_r
$$
for $r>0$, since the spectral seminorm $\|\cdot\|_{\Sp}$ upper-bounds the seminorm $\|\cdot\|$, and is power-multiplicative.
It follows from the following proposition that $F^o$ is integrally closed (\ref{dfn-vsbr-integrallyclosedfiltrations1}), and, for any filtration of definition $F_0$, the inclusions
$$
(F_0)^{\int}\subseteq F^o\subseteq F^{\Sp}
$$
hold.

\begin{prop}\label{prop-vsbr-powerboundedfiltration2}
Let $(A,\|\cdot\|)$ be a seminormed ring, and $F_0$ a filtration of definition of $A$.
Let $F$ be an $\R_+$-multiplicative filtration on $A$ such that $F_0\subseteq F\subseteq F^o$.
Then we have $F^{\int}\subseteq F^o$. 
In particular, $F^o$ is integrally closed.
\end{prop}

\begin{proof}
Let $x\in F^+_r$ satisfy an equality of the form 
$$
x^n+a_1x^{n-1}+\cdots+a_{n-1}x+a_n=0,
$$
where $a_i\in F_{r^i}$ for $i=1,\ldots,n$.
For any $k\geq n$, we have an expression
$$
x^k=\sum^n_{i=1}\Big(\sum_{\lambda}C_{\lambda}a^{\lambda}\Big)x^{n-i}
$$
where $C_{\lambda}\in\Z$ and $\lambda$ runs over all partitions of $k-(n-i)$ such as $\lambda=(\lambda_1,\ldots,\lambda_n)\in\Z^n_{\geq 0}$ and $\sum^n_{j=1}\lambda_j\cdot j=k-(n-i)$, and $a^{\lambda}=\prod^n_{j=1}a^{\lambda_j}_j$.
By the assumption, there exists $C\geq 1$ such that $\|a^l_i\|\leq Cr^{il}$.
Then
$$
\|x^k\|\leq Dr^k
$$
with $D=(r^{-1}C)^n\sup_{i=1,\ldots,n}r^i\|x^{n-i}\|$, which shows that $x\in F^o_r$.
\end{proof}

\begin{cor}\label{cor-vsbr-powerboundedfiltration2}
Let $(A,F_0)$ be a filtered ring.
Let $F$ be an $\R_+$-multiplicative filtration on $A$ such that $F_0\subseteq F\subseteq F^o$.
Then $F^{\int}\subseteq F^o$.
In particular, $F^o$ is integrally closed.
\end{cor}

\begin{prop}\label{prop-vsbr-powerboundedfiltration3}
Let $(A,F)$ be a filtered ring.

{\rm (1)} For a filtered multiplicative subset $S=\bigcup_{r>0}S_r$ $($see {\rm \S\ref{subsub-vsbr-filteredrings})}, the filtration $S^{-1}F^o$ on $S^{-1}A$ is contained in $(S^{-1}F)^o$.

{\rm (2)} Passing to the completion with respect to $F$, we have $(\widehat{F})^o=(F^o)^{\wedge}_F$.

{\rm (3)} The induced filtration $(F^o)_w$ on the weighted polynomial algebra $($resp.\ power series algebra$)$ $A[X]$ $($resp.\ $A\dl (X,w)\dr)$ coincides with the power bounded filtration of $A[X]$ $($resp.\ $A\dl (X,w)\dr)$.
\end{prop}

\begin{proof}
(1) and (2) are easy to see.
For (3), it suffices, due to (2), to discuss the case of the polynomial algebra.
For a monomial $cX^{\alpha}$ in $(F_w)_r$, we have $c\in F_s$, where $s=w(X^{\alpha})^{-1}r$, and $X^{\alpha}$ is $w(X^{\alpha})$-power bounded.
By assumption, $c$ is $s$-power bounded, and hence the claim follows.
\end{proof}

\subsubsection{$\R_+$-affinoid rings}\label{subsub-vsbr-Rplusaffinoidrings}
\begin{dfn}\label{dfn-vsbr-Rplusaffinoidrings}{\rm 
An {\it $\R_+$-affinoid ring}\index{affinoid!affinoid ring@--- ring!Rplus affinoid ring@$\R_+$- --- ---}, or an {\it affinoid Banach ring}\index{affinoid!affinoid ring@--- ring!Banach affinoid ring@Banach --- ---}, is a triple $\mathcal{A}=(A,F_0,F)$ consisting of the following data:
\begin{itemize}
\item a Banach ring $A$ with a filtration $F_0$ of definition;
\item an integrally closed filtration $F$ such that $F^{\int}_0\subseteq F\subseteq F^{\Sp}$.
\end{itemize}}
\end{dfn}

Note that, by \ref{thm-vsbr-filtrationintegral}, $F_r$ contains $(F^{\Sp})_{<r}$, and any $\R_+$-affinoid ring structure of $A$ is obtained from that of the uniform completion $A^{\Sp}$ by taking the preimage under $A\rightarrow A^{\Sp}$.

In the sequel, we will often drop the filtration of definition $F_0$ from the notation of $\R_+$-affinoid rings, and write them as a pair $\mathcal{A}=(A,F_{\mathcal{A}})$.
The underlying Banach ring $A$ of an $\R_+$-affinoid ring $\mathcal{A}$ will be referred to as
$$
\mathcal{A}^{\mathrm{B}}.
$$
The graded ring $\mathrm{Gr}_{F_{\mathcal{A}}}\mathcal{A}^B$ will be simply denoted by
$$
\mathrm{Gr}\,\mathcal{A}.
$$

\begin{dfn}\label{dfn-vsbr-Rplusaffinoidringshom}{\rm 
Let $\mathcal{A}=(A,F_{\mathcal{A}})$ and $\mathcal{B}=(B,F_{\mathcal{B}})$ be two $\R_+$-affinoid rings.
A {\it homomorphism} $\varphi\colon \mathcal{A}\rightarrow \mathcal{B}$ of $\R_+$-affinoid rings is a bounded homomorphism $\varphi^{\mathrm{B}}\colon A\rightarrow B$ that respects the filtrations, that is, $\varphi^{\mathrm{B}}((F_{\mathcal{A}})_r)\subseteq(F_{\mathcal{B}})_r$ for all $r>0$.
In this situation, the Banach affinoid ring $\mathcal{B}$ is also called an {\it $\mathcal{A}$-affinoid algebra}.}
\end{dfn}

We denote by
$$
\R_+\textrm{-}\mathbf{Aff}
$$
the category of $\R_+$-affinoid rings.
We have a forgetful functor 
$$
\mathcal{A}\longmapsto\mathcal{A}^{\mathrm{B}}
$$
from $\R_+\textrm{-}\mathbf{Aff}$ to the category of Banach rings with bounded homomorphisms.
Note that, if $\varphi\colon\mathcal{A}\rightarrow\mathcal{B}$ is a homomorphism of $\R_+$-affinoid rings, then, replacing the filtration of definition $F_{A,0}$ of $A$ by $F_{A,0}\cap(\varphi^{\mathrm{B}})^{-1}(F_{B,0})$, which is also a filtration of definition of $A$, one may always assume that $\varphi^{\mathrm{B}}$ respects filtrations of definition.

Note also that any complete filtered valuation field $(K=\widehat{K},F)$ can be uniquely regarded as an $\R_+$-affinoid ring, since $F_0=F^{\int}_0=F\subseteq F^{\Sp}_K=F^+$, where the first equality is due to \ref{prop-vsbr-gradedvaluationrings1} and \ref{lem-vsbr-integrallyclosedfiltrations1}.

The category of $\R_+$-affinoid rings has tensor products.
For two $\R_+$-affinoid rings $\mathcal{A}_i=(A_i,F_{\mathcal{A}_i})$ ($i=1,2$) over $\mathcal{A}=(A,F_{\mathcal{A}})$, we first take, as above, filtrations of definition of $A$, $A_1$, and $A_2$ such that $(A_i,F_{A_i,0})$ are filtered algebra over $(A,F_{A_0})$, and then take the complete tensor product $B=A_1\widehat{\otimes}_AA_2$ with respect to these filtrations of definition; then the completion of $F_{\mathcal{A}_1}\otimes_{F_{\mathcal{A}}}F_{\mathcal{A}_2}$ with respect to the filtration of definition of $B$ gives a filtration, of which the integral closure $F_{\mathcal{B}}$ gives an $\R_+$-affinoid ring $\mathcal{B}$. 
Then pair $\mathcal{B}=(B,F_{\mathcal{B}})$ gives the desired tensor product of the $\R_+$-affinoid rings $\mathcal{A}_1$ and $\mathcal{A}_2$ over $\mathcal{A}$.

For an $\R_+$-affinoid ring $\mathcal{A}=(A,F_{\mathcal{A}})$, the power series algebra $A\dl r^{-1}X\dr$ as in \S\ref{subsub-vsbr-Rfinitetype} has the canonical $\mathcal{A}$-affinoid ring structure by $F_{\mathcal{A},w}$; see \S\ref{subsub-vsbr-filteredpolynomials} for the definition of the filtration $F_{\mathcal{A},w}$.

For any $\R_+$-finite type $A$-algebra with a presentation (that is, an admissible epimorphism (see \S\ref{sub-vsbr-banachrings}) of $A$-algebras)
$$
\phi\colon A\dl r^{-1}X\dr\longrightarrow B
$$
has an $\mathcal{A}$-affinoid ring structure by the integral closure of the filtration induced from $F_{\mathcal{A},w}$; note that the $\mathcal{A}$-affinoid algebra structure on $B$ may depend on the choice of the presentation $\phi$.
In particular, any $\R_+$-rational localization $B=A\dl r^{-1}_1\frac{f_1}{f_0},\ldots,r^{-1}_n\frac{f_n}{f_0}\dr$ is equipped with a structure of $\mathcal{A}$-affinoid algebra, which we henceforth denote by 
$$
\mathcal{A}\Big(\frac{r^{-1}_1f_1,\ldots,r^{-1}_nf_n}{f_0}\Big),
$$
and call it the {\it $\R_+$-rational localization} of $\mathcal{A}$.

\begin{prop}\label{prop-cor-vsbr-powerboundedfiltration3}
Let $\mathcal{A}=(A,F_{\mathcal{A}})$ be an $\R_+$-affinoid ring\index{affinoid!affinoid ring@--- ring!Rplus affinoid ring@$\R_+$- --- ---} {\rm (\ref{dfn-vsbr-Rplusaffinoidrings})} with $F_{\mathcal{A}}\subseteq(F_{\mathcal{A}})^o$.
Then, for any $\R_+$-finite type $A$-algebra $B$ and any $\mathcal{A}$-affinoid algebra structure $\mathcal{B}=(B,F_{\mathcal{B}})$ on $B$, we have $F_{\mathcal{B}}\subseteq(F_{\mathcal{B}})^o$.
\end{prop}

\begin{proof}
Let $\mathcal{A}\dl r^{-1}X\dr\rightarrow\mathcal{B}$ be a presentation.
By \ref{prop-vsbr-powerboundedfiltration3} (3), we have $(F_{\mathcal{A}})_w\subseteq (F^o_{\mathcal{A}})_w$.
Passing to the quotient, the image of $(F_{\mathcal{A}})_w$, and hence also the integral closure $F_{\mathcal{B}}$, is contained in $(F_{\mathcal{B}})^o$.
\end{proof}

\subsubsection{Valuative spectrum}\label{subsub-vsbr-valuativespectrum}
\index{valuative!valuative spectrum@--- spectrum|(}
In the sequel, for an $\R_+$-affinoid ring $\mathcal{A}$, we write $\mathscr{M}(\mathcal{A})=\mathscr{M}(\mathcal{A}^{\mathrm{B}})$.

\begin{dfn}\label{dfn-vsbr-valuativespectrum}{\rm 
Let $\mathcal{A}=(A,F_{\mathcal{A}})$ be an $\R_+$-affinoid ring.

(1) A {\it valuation} $v$ of $\mathcal{A}$ is a pair $v=(z,F)$, where $z\in\mathscr{M}(\mathcal{A})$ and a filtration $F$ on the complete residue field $\mathcal{H}(z)$ that makes the pair $(\mathcal{H}(z),F)$ a filtered valuation field over $(A,F_{\mathcal{A}})$.

(2) We set
$$
\Spec^{\mathrm{val}}\mathcal{A}=\textrm{the set of all valuation of $\mathcal{A}$,}
$$
and call it the {\it valuative spectrum} of the $\R_+$-affinoid ring $\mathcal{A}=(A,F)$.

(3) For a filtered ring $(A,F_0)$, consider an integrally closed multiplicative filtration $F$ on $A$ such that $F^{\int}_0\subseteq F\subseteq F^{\Sp}_0$.
We define
$$
\Spec^{\mathrm{val}}(A,F)=\Spec^{\mathrm{val}}(\widehat{A},\widehat{F}).
$$
}
\end{dfn}

We are going to define a topology on the set $\Spec^{\mathrm{val}}\mathcal{A}$.
For an $(n+1)$-tuple $f=(f_0,f_1,\ldots,f_n)$ of elements in $A=\mathcal{A}^{\mathrm{B}}$ such that $(f_0,f_1,\ldots,f_n)=A$ and $n$-tuple $r=(r_1,\ldots,r_n)$ of positive real numbers, define a subset 
$$
U_0(f,r)=\bigg\{x=(z,F)\in\Spec^{\mathrm{val}}\mathcal{A}\,\bigg|\,\frac{f_j(z)}{f_0(z)}\in F_{r_j}\ \textrm{for}\ j=1,\ldots,n\bigg\}
$$
of $\Spec^{\mathrm{val}}\mathcal{A}$; note that $f_0(z)\neq 0$ on $U_0(f,r)$.
We call a subset of this form a {\it rational subdomain} (or more precisely, {\it $\R_+$-rational subdomain}).
For example, for $f,g\in A$ and $r,s>0$, 
$$
U_0((g,f,1),(r,s))=\bigg\{(z,F)\,\bigg|\,g(z)\neq 0,\ \frac{f(z)}{g(z)}\in F_r,\ \frac{1}{g(z)}\in F_s\bigg\}
$$
is a rational subdomain.

Note that the intersection of finitely many rational subdomains is again a rational subdomain; e.g., for $f=(f_0,f_1,\ldots,f_n)$, $r=(r_1,\ldots,r_n)$ and $f'=(f'_0,f'_1,\ldots,f'_m)$, $r'=(r'_1,\ldots,r'_m)$, we have
$$
U_0(f,r)\cap U_0(f',r')=U_0(h,s),
$$
where 
\begin{equation*}
\begin{split}
h&=(f_0f'_0,f_if'_j:0\leq i\leq n,0\leq j\leq m,(i,j)\neq (0,0)),\\
s&=(r_ir'_j:0\leq i\leq n,0\leq j\leq m,(i,j)\neq (0,0)).
\end{split}
\end{equation*}

\begin{dfn}\label{dfn-vsbr-valuativespectrum2}{\rm 
We endow a topology on $\Spec^{\mathrm{val}}\mathcal{A}$ in such a way that the set of all $\R_+$-rational subdomains forms an open basis.}
\end{dfn}

Note also that the homomorphism $\mathcal{A}\rightarrow\mathcal{A}(\frac{r^{-1}_1f_1,\ldots,r^{-1}_nf_n}{f_0})$ of $\mathcal{A}$-affinoid algebras induces a bijection
$$
\Spec^{\mathrm{val}}\mathcal{A}\Big(\frac{r^{-1}_1f_1,\ldots,r^{-1}_nf_n}{f_0}\Big)\stackrel{\sim}{\longrightarrow}U_0(f,r).
$$
One can show, similarly to the classical case as in \cite{BGR}, 7.2.4, that rational subdomain of a rational subdomain is, via the bijection as above, again a rational subdomain.
From this, one sees that the above bijection is in fact a homeomorphism.

If $\mathcal{B}$ and $\mathcal{B}'$ are $\R_+$-rational localizations of $\mathcal{A}$, the complete tensor product $\mathcal{B}\widehat{\otimes}_{\mathcal{A}}\mathcal{B}'$ gives the $\R_+$-rational localization that gives the intersection $U\cap U'$, where $U=\Spec^{\mathrm{val}}\mathcal{B}$ and $U'=\Spec^{\mathrm{val}}\mathcal{B}'$ are the corresponding $\R_+$-rational subdomains.

Subsets of the following form are useful for our later argument; for $f,g\in\mathcal{A}$ and $r>0$, define
$$
B(g,f,r)=\bigg\{(z,F)\in\Spec^{\mathrm{val}}\mathcal{A}\,\bigg|\,g(z)\neq 0,\ \frac{f(z)}{g(z)}\in F_r\bigg\},
$$
and call it a {\it basic subset}.
A basic subset is an open subset, since
$$
B(g,f,r)=\bigcup_{s>0}U_0((g,f,1),(r,s)).
$$
Moreover, since
$$
U_0(f,r)=\bigcap^n_{i=1}B(f_0,f_i,r_i)
$$
holds for $f=(f_0,f_1,\ldots,f_n)$ and $r=(r_1,\ldots,r_n)$, the basic subsets form a subbase of the topology on $\Spec^{\mathrm{val}}\mathcal{A}$.

There is a canonical map
$$
\sep_{\mathcal{A}}\colon\Spec^{\mathrm{val}}\mathcal{A}\longrightarrow\mathscr{M}(\mathcal{A}),
$$
which maps $x=(z,F)$ to $z$.
This map is surjective; indeed, there exists a (non-continuous, in general) section of $\sep_{\mathcal{A}}$, which maps $z$ to $(z,F_z)$, where $F_z$ is the filtration on $\mathcal{H}(z)$ corresponding to the absolute value of $\mathcal{H}(z)$, that is, the filtration such that $(\mathcal{H}(z),F_z)$ is of maximal type.
We can view $\mathscr{M}(\mathcal{A})$ as a subset of $\Spec^{\mathrm{val}}\mathcal{A}$ in this manner.
For a rational subdomain $U=U_0(f,r)$ as above, the image $\sep_{\mathcal{A}}(U)$ coincides with $\mathscr{M}(\mathcal{A}(\frac{r^{-1}_1f_1,\ldots,r^{-1}_nf_n}{f_0}))$; in particular, $\sep_{\mathcal{A}}(U)$ is a closed subset in $\mathscr{M}(\mathcal{A})$.

\begin{prop}\label{prop-vsbr-valuativespectrum}
The map $\sep_{\mathcal{A}}$ is a continuos map.
\end{prop}

\begin{proof}
By the definition of the topology of $\mathscr{M}(\mathcal{A})$, subsets of the form
$$
\mathbf{X}(f^{\epsilon},r^{\epsilon})=\{z\in\mathscr{M}(\mathcal{A}): |f(z)|^{\epsilon}<r^{\epsilon}\}
$$
for $f\in\mathcal{A}$, $\epsilon=\pm 1$ (where, if $\epsilon=-1$, we tacitly suppose $f(z)\neq 0$ in the left-hand set), and $r>0$, form an open basis of the topology of $\mathscr{M}(\mathcal{A})$.
In accordance with this, define a subset $X(f^{\epsilon},r^{\epsilon})$ of $\Spec^{\mathrm{val}}\mathcal{A}$ by
$$
X(f^{\epsilon},r^{\epsilon})=\begin{cases}
\bigcup_{s<r}U_0((1,f),s)&(\epsilon=1),\\ \bigcup_{s^{-1}<r^{-1}}U_0((f,1),s^{-1})&(\epsilon=-1).\end{cases}
$$
Since for any $(z,F)\in\Spec^{\mathrm{val}}\mathcal{A}$ and $r>0$ the equality $F_{<r}=(F_z)_{<r}$ holds, we have
$$
\sep^{-1}_{\mathcal{A}}(\mathbf{X}(f^{\epsilon},r^{\epsilon}))=X(f^{\epsilon},r^{\epsilon}),
$$
from which the assertion follows.
\end{proof}
\index{valuative!valuative spectrum@--- spectrum|)}

\subsubsection{Basic properties of valuative spectrum}\label{subsub-vsbr-basicpropertiesvalspec}
We first note the following topological property of valuative spectrums.
\begin{lem}\label{lem-vsbr-toppropertyvalspec1}
Let $\mathcal{A}$ be an $\R_+$-affinoid ring, and $\mathcal{A}_1$ and $\mathcal{A}_2$ $\mathcal{A}$-affinoid algebras.
Then the canonical mapping of topological spaces
$$
\Spec^{\mathrm{val}}\mathcal{A}_1\widehat{\otimes}_{\mathcal{A}}\mathcal{A}_2\longrightarrow\Spec^{\mathrm{val}}\mathcal{A}_1\times_{\Spec^{\mathrm{val}}\mathcal{A}}\Spec^{\mathrm{val}}\mathcal{A}_2
$$
is surjective.
\end{lem}

\begin{proof}
This follows from \ref{cor-vsbr-GromovBerkovichspectrum2}.
\end{proof}

For $z\in\mathscr{M}(\mathcal{A})$, we set
$$
C_z=\sep^{-1}_{\mathcal{A}}(z)=\{x=(z,F)\in\Spec^{\mathrm{val}}\mathcal{A}\}
$$
with the subspace topology as a subspace of $\Spec^{\mathrm{val}}\mathcal{A}$.
\begin{lem}\label{lem-vsbr-toppropertyvalspec3}
There exists a canonical homeomorphism 
$$
\psi_z\colon C_z\stackrel{\sim}{\longrightarrow}\mathrm{ZR}(\mathrm{Gr}_{F_z}\mathcal{H}(z),\mathrm{Gr}\,\mathcal{A}),
$$
where $F_z$ is the filtration on $\mathcal{H}(z)$ induced from the absolute value; see {\rm \S\ref{subsub-vsbr-spacegradedvals}}.
In particular, $C_z$ is a valuative space.
\end{lem}

\begin{proof}
The last assertion is a consequence of \ref{thm-vsbr-spacegradedvals}.
We need to construct a homeomorphism as in the first assertion.
Set $\til{\mathcal{H}}(z)=\mathrm{Gr}_{F_z}\mathcal{H}(z)$, which is a graded field, and $\til{\mathcal{A}}=\mathrm{Gr}\,\mathcal{A}$.
Define a map
$$
\psi_z\colon C_z\longrightarrow\mathrm{ZR}(\til{\mathcal{H}}(z),\til{\mathcal{A}}),\quad (z,F)\mapsto(\til{\mathcal{A}}\rightarrow\mathrm{Gr}_F\mathcal{H}(z)).
$$
We first claim that $\psi_z$ is bijective.
First, notice that, since $f_z\colon \mathcal{A}\rightarrow(\mathcal{H}(z),F_z)$ is a filtered homomorphism, it induces a graded homomorphism $\til{f}_z\colon\til{\mathcal{A}}\rightarrow\til{\mathcal{H}}(z)$.
For a filtered valuation $F$ on $\mathcal{H}(z)$, $f_z((F_{\mathcal{A}})_r)\subseteq F_r$ if and only if $(\til{f}_z)_r(\til{\mathcal{A}}_r)\subseteq(\mathrm{Gr}_{F_z}(\mathcal{H}(z)))_r$.
Hence the claim follows from \ref{prop-vsbr-filteredvaluationfield1}.

Next, we show that $\psi_z$ is an open map.
It suffices to show that $\psi_z(U)$ is open for $U=C_z\cap U_0(f,r)$, where $f=(f_0,f_1,\ldots,f_n)$ and $r=(r_1,\ldots,r_n)$.
For $x=(z,F)\in C_z$, $x\in U$ if and only if
\begin{itemize}
\item[{\rm (a)}] $f_0(z)\neq 0$ and $|f_i(z)/f_0(z)|\leq r_i$ for $i=1,\ldots,n$,
\item[{\rm (b)}] $f_i(z)/f_0(z)\in F_{r_i}$.
\end{itemize}
Notice that, under (a), the condition (b) is equivalent to the following one:
\begin{itemize}
\item[{\rm (b$)^{\prime}$}] the residue class $[f_i(z)/f_0(z)]_{r_i}\in(\til{\mathcal{H}}(z))_{r_i}$ belongs to $(\mathrm{Gr}_F\mathcal{H}(z))_{r_i}$ for $i=1,\ldots,n$.
\end{itemize}
Hence, if we set a graded $\til{\mathcal{A}}$-subalgebra $B$ of $\til{\mathcal{H}}(z)$ by
$$
B=\begin{cases}\til{\mathcal{A}}[\frac{f_1(z)}{f_0(z)},\ldots,\frac{f_n(z)}{f_0(z)}]&(\textrm{if (a) holds}),\\ \til{\mathcal{H}}(z)&(\textrm{otherwise}),\end{cases}
$$
we have $\psi_z(U)=U(B)$, which shows the claim.

Finally, to conclude the proof, we need to show that any open subset of $\mathrm{ZR}(\mathrm{Gr}_{F_z}\mathcal{H}(z),\mathrm{Gr}\,\mathcal{A})$ can be obtained in this way.
Since subsets of the form $U(\til{\mathcal{A}}[\alpha])$ for $\alpha\in\til{\mathcal{H}}(z)$ give an open basis of the topology of $\mathrm{ZR}(\mathrm{Gr}_{F_z}\mathcal{H}(z),\mathrm{Gr}\,\mathcal{A})$, we need to find an open subset $U$ of $C_z$ such that $\psi_z(U)=U(\til{\mathcal{A}}[\alpha])$.
Let $A$ be the image of $\mathcal{A}$ in $\mathcal{H}(z)$.
Since the fractional field of $A$ is dense in $\mathcal{H}(z)$, there exists $f,g\in\mathcal{A}$ such that $g(z)\neq 0$ and $f(z)/g(z)$ represents $\alpha$.
Then we have $\psi_z(B(g,f,r))=U(\til{\mathcal{A}}[\alpha])$.
\end{proof}

In the next paragraph, we will show the following theorem.
\begin{thm}\label{thm-vsbr-toppropertyvalspec}
Let $\mathcal{A}$ be an $\R_+$-affinoid ring.
Then $\Spec^{\mathrm{val}}\mathcal{A}$ is a coherent valuative space\index{valuative!valuative topological space@--- (topological) space}\index{space@space (topological)!valuative topological space@valuative ---}.
\end{thm}

By this theorem and \ref{lem-vsbr-toppropertyvalspec3}, we have the following corollary.
\begin{cor}\label{cor-vsbr-toppropertyvalspec1}
Let $\mathcal{A}$ be an $\R_+$-affinoid ring, and $x=(z,F)\in\Spec^{\mathrm{val}}\mathcal{A}$.
Then the homomorphism $\mathcal{A}\rightarrow(\mathcal{H}(z),F)$ induces a homeomorphism from $\Spec^{\mathrm{val}}(\mathcal{H}(z),F)$ to the set $G_x$ of all generizations of $x$ with the subspace topology.
In particular, $\sep_{\mathcal{A}}\colon\Spec^{\mathrm{val}}\mathcal{A}\rightarrow\mathscr{M}(\mathcal{A})$ gives the $\mathrm{T}_1$-quotient {\rm ({\bf \ref{ch-pre}}, \S\ref{subsub-separationgen})}.
\end{cor}

Finally, by {\bf \ref{ch-pre}}.\ref{thm-reflexivization}, we have the following statement.
\begin{cor}\label{cor-vsbr-toppropertyvalspec2}
$(\Spec^{\mathrm{val}}\mathcal{A})^{\mathrm{ref}}$ is a coherent valuative space.
\end{cor}

\subsubsection{Proof of Theorem \ref{thm-vsbr-toppropertyvalspec}}\label{subsub-vsbr-toppropertyvalspec}
\begin{lem}\label{lem-vsbr-toppropertyvalspec4}
Let $z\in\mathscr{M}(\mathcal{A})$.
If a finite family $\{U_j\}_{j\in J}$ of rational subdomains of $\Spec^{\mathrm{val}}\mathcal{A}$ covers $C_z$, then there exists an open neighborhood $\mathbf{V}$ of $z$ in $\mathscr{M}(\mathcal{A})$ such that $\{U_j\}_{j\in J}$ covers $\sep^{-1}_{\mathcal{A}}(\mathbf{V})$.
\end{lem}

\begin{proof}
Set $X=\Spec^{\mathrm{val}}\mathcal{A}$.
Write each $U_j$ as a finite intersection of basic subsets
$$
U_j=\bigcap_{i\in I_j}B_i
$$
for $j\in J$.
For $k=(i_j)_{j\in J}\in\prod_{j\in J}I_j$, we set $B_k=\bigcup_{j\in J}B_{i_j}$.
We have
$$
\bigcup_{j\in J}U_j=\bigcup_{j\in J}\bigcap_{i\in I_j}B_i=\bigcap_{k\in K}B_k.
$$
Hence, to show the lemma, we only need to show the following claim.

\medskip
{\sc Claim.} {\it Let $B$ be a finite union of basic subsets of $X$.
If $B$ contains $\sep^{-1}_{\mathcal{A}}(z)$, then there exists an open neighborhood $\mathbf{V}$ of $z$ in $\mathscr{M}(\mathcal{A})$ such that $\sep^{-1}_{\mathcal{A}}(\mathbf{V})$ is contained in $B$.}

\medskip
To show the claim, write $B$ as a finite union of basic subsets
$$
B=\bigcup_{i\in I}B_i,
$$
where $B_i=(B,g_i,f_i,r_i)$ for $i\in I$.
If there is $i_0\in I$ such that $|f_{i_0}(z)/g_{i_0}(z)|<r_{i_0}$, then 
$$
\mathbf{V}=\{w\in\mathscr{M}(\mathcal{A}): |f_{i_0}(w)/g_{i_0}(w)|<r_{i_0}\}
$$
enjoys the desired property.

We may thus assume $|f_i(z)/g_i(z)|\geq r_i$ for any $i\in I$.
Let $T'=(T'_i)_{i\in I}$ and $r'=(r_i)_{i\in I}$.
For $z'\in\sep_{\mathcal{A}}(X)$ with $f_i(z')\neq 0$ and $|g_i(z')/f_i(z')|\leq r^{-1}_i$ for any $i\in I$, define
$$
\varphi_{z'}\colon\mathcal{A}\dl r'T'\dr\longrightarrow\mathcal{H}(z'),\quad T'_i\longmapsto\frac{g_i(z')}{f_i(z')},
$$
and set 
$$
\til{\mathcal{A}\dl r'T'\dr}=\mathrm{Gr}_{(F_{\mathcal{A}})_w}\mathcal{A}\dl r'T'\dr.
$$
Let $G_{z'}$ be the image of $\til{\mathcal{A}\dl r'T'\dr}$ in $\til{\mathcal{H}}(z')=\mathrm{Gr}_{F_{z'}}\mathcal{H}(z')$ by the homomorphism induced from $\varphi_{z'}$.
Set $h_i(z')=\varphi_{z'}(T'_i)$ for $i\in I$.
Then $h_i(z')$ belongs to $(F_{z'})_{r^{-1}_i}$.

Now, the following two conditions are equivalent to each other:
\begin{itemize}
\item[{\rm (a)}] $\{B_i\}_{i\in I}$ covers $C_{z'}$;
\item[{\rm (b)}] the homogeneous ideal of $G_{z'}$ generated by $[h_i(z')]_{r^{-1}_i}$ for $i\in I$ is $G_{z'}$ itself.
\end{itemize}
Indeed, by the map $\psi_{z'}$ as in \ref{lem-vsbr-toppropertyvalspec3}, we have $\psi_{z'}(B_i)=B(\til{\mathcal{A}\dl r'T'\dr}[\alpha_i])$, where $\alpha_i$ is the class of $f_i(z')/g_i(z')$ in $\til{\mathcal{H}}(z')$.
If (a) holds, then for the subset $I'=\{i\in I:|f_i(z')/g_i(z')|=r_i\}$, $\{B_i\}_{i\in I'}$ already covers $C_{z'}$, as we have assumed $|h_i(z')|=|g_i(z')/f_i(z')|\leq r^{-1}_i$.
Since $\alpha_i\neq 0$ for $i\in I'$, we obtain (b) as in the proof of \ref{thm-vsbr-spacegradedvals}.
Conversely, if (b) holds, then similarly, for $I''=\{i\in I:[h_i(z')]_{r^{-1}_i}\neq 0\}$, $\{B_i\}_{i\in I''}$ covers $C_{z'}$.

Now, by the assumption that $\{B_i\}_{i\in I}$ covers $C_z$, we know that there exists an expression
$$
1=\sum_{i\in I}\ovl{a}_i[h_i(z)]_{r^{-1}_i}
$$
for some degree $r_i$ element $\ovl{a}_i$ of $\til{\mathcal{A}\dl r'T'\dr}$ for $i\in I$.
One can take the preimage $a_i\in(F_{\til{\mathcal{A}\dl r'T'\dr}})_{r_i}$ of $\ovl{a}_i$ by polynomials in $T'$.
Then the polynomial
$$
P(T')=\Big(\sum_{i\in I}a_iT'_i\Big)-1
$$
in $T'$ belongs to $(F_{\til{\mathcal{A}\dl r'T'\dr}})_1$ and satisfies
$$
|\varphi_z(P)(z)|<1.
$$

Define
$$
\mathbf{V}=\bigg\{z'\in\mathscr{M}(\mathcal{A})\,\bigg|\,f_i(z')\neq 0\ (i\in I),\ \Big|\Big(P\big|_{T'_i=\frac{g_i}{f_i}}\Big)(z')\Big|<1\bigg\}
$$
which is an open neighborhood of $z$ in $\mathscr{M}(\mathcal{A})$.
To show that this $\mathbf{V}$ has the desired property, suppose contrary.
That is, there exists $x'=(z',F)\in\sep^{-1}_{\mathcal{A}}(\mathbf{V})$ not covered by $\{B_i\}_{i\in I}$.
Since $x'\in \bigcap_{i\in I}(X\setminus B_i)$, we have $|g_i(z')/f_i(z')|\leq r^{-1}_i$ for any $i\in I$, and hence we have $\varphi_{z'}\colon\til{\mathcal{A}\dl r'T'\dr}\rightarrow \mathcal{H}(z')$ as above.
Moreover, $|\varphi_{z'}(P)(z')|<1$ holds by the definition of $\mathbf{V}$.
For $a_i(z')=\varphi_{z'}(a_i)\in (F_{z'})_{r_i}$, the classes $[a_i(z')]_{r_i}$ and $[h_i(z')]_{r^{-1}_i}$ satisfy
$$
\sum_{i\in I}[a_i(z')]_{r_i}[h_i(z')]_{r^{-1}_i}-1=\big(\varphi_{z'}(P)(z')\ \mathrm{mod}\ (F_{z'})_{<1}\big)=0.
$$
Thus we see that the homogeneous ideal in $G_{z'}$ generated by $[h_i(z')]_{r^{-1}_i}$ for $i\in I$ is $G_{z'}$ itself.
But this is absurd in view of the equivalence of the conditions (a) and (b) above. 
We have therefore proved the claim, and hence the lemma.
\end{proof}

\begin{proof}[Proof of Theorem {\rm \ref{thm-vsbr-toppropertyvalspec}}]
Set $X=\Spec^{\mathrm{val}}\mathcal{A}$.
We first show that any rational subdomain $U\subseteq X$ is quasi-compact.
To show this, we may assume $U=X$ without loss of generality.
Since $\mathscr{M}(\mathcal{A})$ is compact, it suffices to show that $\pi=\sep_{\mathcal{A}}$ is proper.
Note that, for any point $z\in\mathscr{M}(\mathcal{A})$, the fiber $\pi^{-1}(z)$ is quasi-compact due to \ref{lem-vsbr-toppropertyvalspec3} and \ref{thm-vsbr-spacegradedvals}.
Hence, by \cite{Bourb4}, Chap.\ I, \S10.2, Theorem 1, it suffices to show that $\pi$ is a closed map.
Assume contrary. 
Let $C\subseteq X$ be a closed subset such that $\pi(C)$ is not closed in $\mathscr{M}(\mathcal{A})$, and take $z\in\ovl{\pi(C)}\setminus\pi(C)$, where $\ovl{\pi(C)}$ denotes the closure of $\pi(C)$.
Let $\{W_j\}_{j\in J}$ be an open covering of $W=X\setminus C$ by rational subdomains. 
By our assumption, $\pi^{-1}(z)\subseteq C$, and there is a finite subset $J_z\subseteq J$ such that $\bigcup_{j\in J_z}W_j$ contains $\pi^{-1}(z)$.
By \ref{lem-vsbr-toppropertyvalspec4}, there exists an open neighborhood $\mathbf{V}_z$ of $z$ in $\mathscr{M}(\mathcal{A})$ such that $\bigcup_{j\in J_z}W_j$ contains $\pi^{-1}(\mathbf{V}_z)$.
Hence $C\cap\pi^{-1}(\mathbf{V}_z)=\emptyset$, that is, $z\in\ovl{\pi(C)}$, which is absurd.
Hence we have shown that $\pi$ is closed, and hence that any rational subdomain of $X$ is quasi-compact.
Notice that rational subdomains are closed under finite intersection

Next, we show that $X$ is $\mathrm{T}_0$.
Take $x=(z,F)$ and $z'=(z',F')$. If $z\neq z'$, then we can find open neighborhoods $\mathbf{V}_z$ and $\mathbf{V}_{z'}$ of $z$ and $z'$, respectively, such that $\mathbf{V}_z\cap\mathbf{V}_{z'}=\emptyset$.
then $\pi^{-1}(\mathbf{V}_z)$ and $\pi^{-1}(\mathbf{V}_{z'})$ separate the points $x$ and $x'$.
If $z=z'$, then we only need to invoke the $\mathrm{T}_0$-ness of $C_z\cong\mathrm{ZR}(\mathrm{Gr}_{F_z}\mathcal{H}(z),\mathrm{Gr}_F\mathcal{A})$.

Hence, we have shown that $X$ is a coherent topological space.
Next, we show that $X$ is sober.
Let $Z$ be an irreducible closed subset of $X$.
Since $\pi=\sep_{\mathcal{A}}$ is, as we have seen above, a closed map, $\pi(Z)$ is an irreducible closed subset of $\mathscr{M}(\mathcal{A})$.
Since $\mathscr{M}(\mathcal{A})$ is Hausdorff, $\pi(Z)=\{z\}$ for some $z\in\mathscr{M}(\mathcal{A})$.
Since $C_z=\pi^{-1}(z)$ is sober (\ref{thm-vsbr-spacegradedvals}), $Z$ has the unique generic point, as desired.

Finally, by \ref{lem-vsbr-toppropertyvalspec3} and \ref{thm-vsbr-spacegradedvals}, we deduce that $X$ is valuative.
\end{proof}

\subsubsection{Relation with adic spectrum}\label{subsub-vsbr-adicspectrum}
Let $A$ be a complete f-adic ring\index{adic!adic ring@--- ring!fadic ring@f-{---} ---}\index{fadic ring@f-adic ring}\index{fadic ring@f-adic ring!complete fadic ring@complete ---}.
Then $A$ is a Banach ring with respect to the following norm.
Let $A_0$ be a ring of definition of $A$ with an ideal of definition $I_0$ (see {\bf \ref{ch-pre}}.\ref{dfn-fadic}), and fix a real number $0<c<1$.
Define an decreasing filtration $\{F^m_{(A_0,I_0)}\}_{m\in\Z}$ indexed by integers by
$$
F^m_{(A_0,I_0)}=\begin{cases}
I^m_0&(m\geq 0),\\
[A_0:I^{-m}_0]&(m<0).
\end{cases}
$$
Then the norm $\|\cdot\|_A$ on $A$ is defined by, for $f\in A$, 
$$
\|f\|_A=\begin{cases}
c^n,\ \textrm{where}\ n=\inf\{m\in\Z:f\in F^m_{(A_0,I_0)}\}&(f\neq 0),\\ 
0&(f=0).
\end{cases}
$$

Note that the norm $\|\cdot\|_A$ actually depends on the choice of $(A_0,I_0)$.
Note also that a homomorphism $\varphi\colon A\rightarrow B$ inducing an adic morphism $(A_0,I_0)\rightarrow (B_0,I_0B_0)$ between the rings of definition, respects the filtrations $F_{(A_0,I_0)}$ and $F_{(B_0,I_0B_0)}$, and gives rise to a bounded homomorphism 
$$
\varphi\colon (A,\|\cdot\|_A)\longrightarrow (B,\|\cdot\|_B).
$$
For a Huber's affinoid ring $(\mathcal{A}^{\pm},\mathcal{A}^+)$  (in the sense as in {\rm \ref{dfn-affinoidringadicspace}}), an ideal $I$ of $\mathcal{A}^{+}$ is called an {\it ideal of definition}, if there exists a pair $(A_0, I_0)$ consisting of a ring of definition and an ideal of definition, such that $I = I_0\mathcal{A}^{+}$. 
Then the equivalence class of the norm defined as above from $(A_0, I_0)$ depends only on the data $((\mathcal{A}^{\pm},\mathcal{A}^+), I )$.

\begin{dfn}\label{dfn-vsbr-Banachringfadictype}{\rm
(1) A Banach ring is said to be {\it of f-adic type} if it is isomorphic to an f-adic ring as a topological ring.

(2) An $\R_+$-affinoid ring $\mathcal{A}=(\mathcal{A}^{\mathrm{B}},F_{\mathcal{A}})$ is said to be {\it of adic type} if there exist an affinoid ring\index{affinoid!affinoid ring@--- ring} $(\mathcal{A}^{\pm},\mathcal{A}^+)$ in the sense as in {\rm \ref{dfn-affinoidringadicspace}} and an topological isomorphism $\mathcal{A}^{\mathrm{B}}\cong\mathcal{A}^{\pm}$ that induces $(F_{\mathcal{A}})_1\cong\mathcal{A}^+$ and $\mathscr{M}(\mathcal{A})\cong[\Spa (\mathcal{A}^{\pm},\mathcal{A}^+)]$.}
\end{dfn}

For an affinoid ring $(\mathcal{A}^{\pm},\mathcal{A}^{+})$ and an ideal of definition $I\subset\mathcal{A}^{+}$, one can construct an $\R_+$-affinoid ring $\mathcal{A}=(\mathcal{A}^{\mathrm{B}},F_{\mathcal{A}})$ of adic type such that $\mathcal{A}^{\mathrm{B}}\cong\mathcal{A}^{\pm}$, $(F_{\mathcal{A}})_1\cong\mathcal{A}^+$, and $\mathscr{M}(\mathcal{A})\cong[\Spa (\mathcal{A}^{\pm},\mathcal{A}^+)]$ as follows.

Set $A=\mathcal{A}^{\pm}$, and choose $(A_0, I_0)$ such that $I=I_0\mathcal{A}^{+}$. 
We consider the norm $\|\cdot\|$ on $A=\mathcal{A}^{\pm}$ defined as above.
Consider the filtration $F_0=F^+_0$ on $A$ corresponding to this norm.
Then we define the multiplicative filtration $F_{\mathcal{A}}$ on $A$ as the integral closure of the one generated by $\mathcal{A}^+$ and $F_0$, which gives an $\R_+$-affinoid ring $\mathcal{A}=(A,F_{\mathcal{A}})$ such that $(F_{\mathcal{A}})_1\cong\mathcal{A}^+$.
Note that the filtration $F_{\mathcal{A}}$ does not depend on a choice of $(A_0, I_0)$, and contains $F^o_A$.

We need to check $\mathscr{M}(\mathcal{A})\cong[\Spa (\mathcal{A}^{\pm},\mathcal{A}^+)]$.
Observe first that, for any $x\in\Spa(\mathcal{A}^{\pm},\mathcal{A}^+)$, the maximal generization $\til{x}$ of $x$ corresponds to a continuous valuation of height zero or one (since the topology is adic).
Let $\varphi_x\colon A\rightarrow K_x$ be the morphism to the corresponding valuation field.
Since it is continuous, then there exists a ring of definition $A'_0$ of $A$ such that $\varphi_x(A'_0)\subseteq V_x$, where $V_x$ is the valuation ring of $K_x$.
Since $\varphi_x$ maps $A^o$ to $V_x$, and $A_0\subseteq A^o$, it follows from $\varphi_x$ induces a continuous homomorphism $A_0\rightarrow V_x$, which is, moreover, adic, since the height of $V$ is one or zero.

Since $A_0\rightarrow V$ is adic, the valuation $\|\cdot\|_x$ on $K_x$ has the following description.
Consider the norm $|\cdot|_x$ on $K_x$ determined as above by the ring of definition $V_x$ and the ideal of definition $I_0V_x$.
Then, for any $f\in K_x$,
$$
\|f\|_x=\lim_{n\rightarrow\infty}|f^n|^{\frac{1}{n}}_x=\inf_{n\geq 1}|f^n|^{\frac{1}{n}}_x,
$$
that is, $\|\cdot\|_x$ is the associated power-multiplcative norm of $|\cdot|$.

Hence $\til{x}$ gives a multiplicative and bounded seminorm\index{seminorm!multiplicative seminorm@multiplicative ---}
$$
\|\cdot\|_{\til{x}}\colon A\longrightarrow\R_{\geq 0}.
$$
Thus we have the map
$$
\Spa(\mathcal{A}^{\pm},\mathcal{A}^+)\longrightarrow\mathscr{M}(A),\qquad x\longmapsto\|\cdot\|_{\til{x}}.
$$
It is clear by definition that this map factors through the separated quotient $[\Spa(A,A^o)]$.

\begin{prop}\label{prop-vsbr-berkovichspectrum}
The resulting map
$$
[\Spa(\mathcal{A}^{\pm},\mathcal{A}^+)]\longrightarrow\mathscr{M}(A)\eqno{(\ast)}
$$
is a homeomorphism.
\end{prop}

\begin{proof}
One checks that the map $(\ast)$ is injective by the boundedness of $\|\cdot\|_{\til{x}}$.
For $y\in\mathscr{M}(A)$, the corresponding bounded multiplicative seminorm\index{seminorm!multiplicative seminorm@multiplicative ---} $y=|\cdot|_y\colon A\rightarrow\R_{\geq 0}$ defines a (maximal) point of $\Spa(\mathcal{A}^{\pm},\mathcal{A}^+)$; note that, since boundedness implies $|\cdot|_y\leq\|\cdot\|_A$, we have $|f|_y\leq 1$ for any $f\in A^o$.
Thus the map $(\ast)$ is bijective.
To compare the topologies of $[\Spa(\mathcal{A}^{\pm},\mathcal{A}^+)]$ and $\mathscr{M}(A)$, it suffices to show that the subsets of the form $\{x\in[\Spa(\mathcal{A}^{\pm},\mathcal{A}^+)]:|f(x)|<1\}$ with $f\in A^o$ form an open basis of the topology of $[\Spa(\mathcal{A}^{\pm},\mathcal{A}^+)]$, which is easy to see.
\end{proof}

If $\mathcal{A}=(\mathcal{A}^{\mathrm{B}},F_{\mathcal{A}})$ is an $\R_+$-affinoid ring of adic type, we write, by a slight abuse of notation, 
$$
F_1({\mathcal{A}})=(\mathcal{A}^{\mathrm{B}},(F_{\mathcal{A}})_1),
$$
which is an affinoid ring in the sense of Huber. The notation indicates that Huber's affinoid structure is obtained by taking the $F_1$-part of $\R_+$-affinoid structures. 
We have a continuous mapping
$$
\Spec^{\mathrm{val}}\mathcal{A}\longrightarrow\Spa F_1(\mathcal{A}),\eqno{(\ast)}
$$
called the {\it $1$-restriction map}, which maps $x=(z,F)$ to the induced valuation $(F_{\mathcal{A}})_1\rightarrow F_1$ on $\mathcal{A}^+=(F_{\mathcal{A}})_1$.

\begin{prop}\label{prop-vsbr-adicspectrum}
The $1$-restriction map $(\ast)$ is surjective.
\end{prop}

\begin{proof}
For $x\in\Spa F_1(\mathcal{A})$, let $V$ be the corresponding valuation ring, and $K=\Frac(V)$.
Let $V_1$ be the maximal generization, which corresponds to the point $\sep(x)$ in $\mathscr{M}(\mathcal{A})\cong[\Spa F_1(\mathcal{A})]$, and $\nu\colon K\rightarrow\R_+$ the corresponding to absolute value.
Consider the subring $G_1=V/(F_{\nu})_{<0}$ of $\ovl{K}_{\nu}=\mathrm{Gr}_{F_{\nu}}K$ in degree $1$, and let $G$ be the graded subring over $\mathrm{Gr}\,\mathcal{A}$ generated by $G_1$ in $\ovl{K}_{\nu}$.
Note that $G_1$ is a valuation ring, and is the degree $1$-part of the graded ring $G$.
Take a graded valuation ring $\til{V}$ of $K_{\nu}$ such that $\m_{G_1}G\subseteq \m_{\til{V}}$.
Then $\til{V}$ is a graded valuation ring with $\til{V}_1=G_1$, and the associated filtered valuation $F$ of $K$ gives a lift of $x$ to $\Spec^{\mathrm{val}}\mathcal{A}$.
\end{proof}

Finally, let us describe the relation between valuative spectra and reified adic spectra in \cite{Kedlaya}.
In \cite{Kedlaya}, 6.1, Kedlaya\index{Kedlaya, K.} introduced the notion of {\it affinoid seminormed ring}, which can be interpreted in our language as a filtered ring $(A,F)$ with a seminormed type integrally closed filtration $F$ contained in $F^o_A$.
The Banach ring case, the so-called {\it affinoid Banach rings}, are pairs of the form 
$$
(A,A^{\mathrm{Gr}})
$$
consisting of a non-archimedean Banach ring $A$ and a graded integrally closed subring $A^{\mathrm{Gr}}$ of the graded ring $\mathrm{Gr}_{F^o}A$ by the power-bounded filtration $F^o$.\footnote{It was Temkin\index{Temkin, M.} \cite{Temkin1} who first considered filtrations on Berkovich's $K$-affinoid algebras of this kind and, in particular, the associated graded rings.}
A morphism $(A,A^{\mathrm{Gr}})\rightarrow(B,B^{\mathrm{Gr}})$ between affinoid Banach rings of this sense is a bounded homomorphism $A\rightarrow B$ that induces $A^{\mathrm{Gr}}\rightarrow B^{\mathrm{Gr}}$.

If $(A,A^{\mathrm{Gr}})$ is an affinoid Banach ring in the sense of Kedlaya, then one can define a multiplicative filtration $F$ on $A$ by
$$
F_r=\textrm{the preimage of $(A^{\mathrm{Gr}})_r$ under $F^o_r\rightarrow(\mathrm{Gr}_{F^o}A)_r=F^o_r/F^o_{<r}$.}
$$
By \ref{prop-vsbr-powerboundedfiltration2}, the resulting filtration $F$ is integrally closed, whence giving an $\R_+$-affinoid ring $(A,F)$.
In this way, one can establish a categorical equivalence between the category of affinoid Banach rings in the sense of Kedlaya and the category of $\R_+$-affinoid rings $\mathcal{A}=(A,F_{\mathcal{A}})$ such that $F_{\mathcal{A}}\subseteq F^o_A$.

Given an affinoid Banach ring $(A,A^{\mathrm{Gr}})$ in the sense of Kedlaya, the {\it reifed adic spectra} $\mathrm{Spra}(A,A^{\mathrm{Gr}})$ is the set of reified valuations, that is, a pair $(v,r)$ consisting of a valuation
$$
v\colon A\longrightarrow\Gamma\cup\{0\}
$$
with a value target group $\Gamma$ (written multiplicatively), and a order-preserving homomorphism (called a {\it reification})
$$
r\colon\R_+\longrightarrow\Gamma,
$$
such that the following conditions are satisfied: 
\begin{itemize}
\item[{\rm (a)}] for any $r>0$, and $f\in F_r$ (where $F$ is the filtration on $A$ constructed from $A^{\mathrm{Gr}}$ as above), $v(f)\leq r$;
\item[{\rm (b)}] for any $\gamma\in\Gamma_v$ ($=$ the ordered subgroup generated by $r(\R_+)$ and the image of $v$), there exists $r\in R_+$ such that $r\leq\gamma$.
\end{itemize}
By (b), one has for any $\gamma\in\Gamma_v$ a positive real number
$$
p(\gamma)=\inf\{r\in\R_+:\gamma\leq r\},
$$
which gives a splitting $p\colon\Gamma_v\rightarrow\R_+$ of the reification $r$.
Hence, as we have seen in \ref{rem-vsbr-filvalviavaluation1}, the reified valuation $v$ gives an $\R_+$-valuation (\ref{dfn-vsbr-filvalviavaluation1}).
Taking the condition (a) into account, one can thus establish a canonical bijection
$$
\mathrm{Spra}(A,A^{\mathrm{Gr}})\stackrel{\sim}{\longrightarrow}\Spec^{\mathrm{val}}(A,F).
$$
Moreover, comparing the notions of `rational subdomains' on both sides, one can show that this bijection gives a homeomorphism.

\subsubsection{$\R_+$-affinoid algebras of $\R_+$-finite type over $K$}\label{subsub-vsbr-berkovichaffinoidalgebras}
Let $K=(K,\nu=|\cdot|)$ be a non-archimedean Banach field.
We allow the case where the valuation $|\cdot|$ is a trivial one.

We call a Banach $K$-algebra of the form 
$$
K\dl r^{-1}X\dr=K\dl r^{-1}_1X_1,\ldots,r^{-1}_nX_n\dr
$$
(defined in \S\ref{subsub-vsbr-filteredpolynomials}) a {\it Berkovich algebra}\index{algebra!Berkovich algebra@Berkovich ---}, which is a Tate algebra\index{algebra!Tate algebra@Tate ---} ({\bf \ref{ch-pre}}, \S\ref{subsub-classicalaffinoidalgebras}), if $r_1=\cdots=r_n=1$. 

\begin{dfn}[{\cite{Berk1}, \S2.1}]\label{dfn-berkovichaffinoicalgebras}{\rm (1) A {\it Berkovich's $K$-affinoid algebra}\index{algebra!affinoid algebra@affinoid ---!Berkovichs affinoid algebra@Berkovich's --- ---} (or simply, {\it $K$-affinoid algebra}, if there is no danger of confusion) is a commutative Banach $K$-algebra $A$ that admits an admissible epimorphism (called {\it presentation}) 
$$
K\dl r^{-1}_1X_1,\ldots,r^{-1}_nX_n\dr\longrightarrow A.\eqno{(\ast)}
$$

(2) If the presentation $(\ast)$ can be found with 
$$
r_1=\cdots=r_n=1, 
$$
or equivalently, with 
$$
r_1,\ldots,r_n\in\sqrt{|K^{\times}|}=\{r\in\R_+:r^n\in|K^{\times}|\ \textrm{for some $n\geq 1$}\}
$$
(cf.\ \cite{BGR}, 6.1.5/4), then $A$ is called a {\it strictly $K$-affinoid algebra}.}
\end{dfn}

We denote by
$$
\mathbf{Aff}^{\mathbf{B}}_K
$$
the category of Berkovich's $K$-affinoid algebra\index{algebra!affinoid algebra@affinoid ---!Berkovichs affinoid algebra@Berkovich's --- ---} with bounded $K$-algebra homomorphisms. 
The category $\mathbf{Aff}^{\mathbf{B}}_K$ contains $\mathbf{Aff}_K$, the category of classical affinoid algebras\index{algebra!affinoid algebra@affinoid ---!classical affinoid algebra@classical --- ---} over $K$ ({\bf \ref{ch-pre}}.\ref{dfn-classicalaffinoidalgebras}) and $K$-algebra homomorphisms, as a full subcategory (cf.\ \cite{BGR}, 6.1.3/1).

Note that, if the norm $|\cdot|$ of $K$ is non-trivial, then strictly $K$-affinoid algebras are nothing but classical affinoid algebras\index{algebra!affinoid algebra@affinoid ---!classical affinoid algebra@classical --- ---} discussed already in {\bf \ref{ch-pre}}, \S\ref{subsub-classicalaffinoidalgebras}.

Berkovich algebra $K\dl r^{-1}X\dr$ has the Gauss norm $\|\cdot\|_{\mathrm{Gauss}}$, and hence the filtration $F_w=F^{\Sp}$ (see \S\ref{subsub-vsbr-filteredpolynomials}), which makes the pair 
$$
(K\dl r^{-1}X\dr, F_{K\dl r^{-1}X\dr})
$$
(where $F_{K\dl r^{-1}X\dr}=F^{\Sp}_{K\dl r^{-1}X\dr}$) an $\R_+$-affinoid algebra over $K=(K,F_{\nu})$.

For a Berkovich's $K$-affinoid algebra $A$ with a presentation as in $(\ast)$, one has a filtration on $A$ induced from the filtration $F_{K\dl r^{-1}X\dr}$ (see \S\ref{subsub-vsbr-filteredrings}), and thus can be regarded as an $\R_+$-affinoid algebra of $\R_+$-finite type over $K=(K,F_{\nu})$.
The next theorem shows that this $\R_+$-affinoid algebra structure does not depend on the choice of the presentation.
\begin{thm}\label{thm-vsbr-berkovichaffinoidalgebras}
Let $\mathcal{A}=(A,F_{\mathcal{A}})$ be an $\R_+$-affinoid algebra of $\R_+$-finite type over $K=(K,F_{\nu})$.
Then $F_{\mathcal{A}}=F^{\Sp}_A$.
\end{thm}

\begin{proof}
Choose a presentation $K\dl r^{-1}X\dr\rightarrow A$ as in $(\ast)$.
We first assume that $|K^{\times}|=\R_+$.
Take elements $a_i\in K$ such that $|a_i|=r_i$ for $i=1,\ldots,n$, which give rise to an isometry 
$$
K\dl r^{-1}X\dr\stackrel{\sim}{\longrightarrow} K\dl Y\dr,
$$
where $Y=(Y_1,\ldots,Y_n)$, mapping $X_i$ to $a_iY_i$ ($i=1,\ldots,n$).
In this way, we may assume assume that the affinoid structure comes from a presentation as in $(\ast)$ with $r_1=\cdots=r_n=1$.
In this case, $(F_{\mathcal{A}})_1$ is the integral closure of $A_0=$ the image of $V\dl X\dr$ in $A$, where $V=(F_K)_1$ is the associated valuation ring of $K$.
Moreover, it is easy to see that $(F^{\Sp}_A)_r=(F_K)_r\otimes_{(F_K)_1}(F^{\Sp}_A)_1$ for any $r>0$.
By \ref{prop-canonicalsubring}, we have $(F_{\mathcal{A}})_1=(F^{\Sp}_A)_1$, and hence we have $F_{\mathcal{A}}=F^{\Sp}_A$.

In general, take an extension of Banach fields $L/K$ such that $|L^{\times}|=\R_+$; e.g., $K^{\mathrm{st}}$ as in \ref{exa-vsbr-Kstandard}.
Suppose $F_{\mathcal{A}}\neq F^{\Sp}_A$.
Then, by \ref{thm-vsbr-filtrationintegral}, there exists a filtered valuation field $(M,F_M)$ over $\mathcal{A}$ such that the image of $F^{\Sp}_A$ is not contained in $F_M$.
Consider the base change $\mathcal{A}_L=(A\widehat{\otimes}_KL,F_L)$ with the tensor product filtration (see \S\ref{subsub-vsbr-filteredrings}).
By \ref{cor-vsbr-GromovBerkovichspectrum2}, there exists a filtered valuation field $(N,F_N)$ dominating $(M,F_M)$, which sits in the commutative diagram
$$
\xymatrix{(M,F_M)\ar[r]&(N,F_N)\\ \mathcal{A}\ar[u]\ar[r]&\mathcal{A}_L\ar[u].}
$$
Since, as we have already seen, $F_L=F^{\Sp}_{A_L}$, $\mathcal{A}\rightarrow (N,F_N)$ factors through $(\mathcal{A},F^{\Sp}_A)$. 
This means that the image of $F^{\Sp}_A$ in $N$ is contained in $F_M$, which is absurd.
\end{proof}

When considering Berkovich's $K$-affinoid algebras, we henceforth regard them as an $\R_+$-affinoid algebra of $\R_+$-finite type over $K$ in the canonical way as in the theorem.
Note that one can also regard Berkovich's $K$-affinoid algebras as $\R_+$-affinoid algebra of adic type (\ref{dfn-vsbr-Banachringfadictype} (2)), if the valuation of $K$ is non-trivial.
\begin{cor}\label{cor-vsbr-berkovichaffinoidalgebras1}
Let $\mathcal{A}$ be a Berkovich's $K$-affinoid algebra, and $L/K$ an extension of Banach fields. 
Set $\mathcal{A}_L=\mathcal{A}\widehat{\otimes}_KL$.
Then the canonical map 
$$
\Spec^{\mathrm{val}}\mathcal{A}_L\longrightarrow\Spec^{\mathrm{val}}\mathcal{A}
$$
 is surjective.
 \end{cor}

\begin{proof}
The $\R_+$-affinoid ring structure structure on $\mathcal{A}_L$ is the integral closure of the complete tensor product filtration of $F_{\mathcal{A}}$ and $F_L$.
Hence the claim follows from \ref{lem-vsbr-toppropertyvalspec1}.
\end{proof}

\begin{prop}\label{prop-vsbr-berkovichaffinoidalgebras1}
Let $\mathcal{A}$ be a Berkovich's $K$-affinoid algebra.
Suppose $|K^{\times}|=\R_+$.
Then the $1$-restriction map {\rm (\S\ref{subsub-vsbr-adicspectrum})}
$$
\Spec^{\mathrm{val}}\mathcal{A}\longrightarrow\Spa F_1(\mathcal{A})
$$
is a homeomorphism.
\end{prop}

\begin{proof}
By \ref{prop-vsbr-examplefilteredvaluations}, the map in question is bijective.
For $f=(f_0,f_1,\ldots,f_n)$ and $r=(r_1,\ldots,r_n)$, one can take $a_i\in K$ such that $|a_i|=r_i$ for $i=1,\ldots,n$, hence
\begin{equation*}
\begin{split}
\Spec^{\mathrm{val}}\mathcal{A}\Big(\frac{r^{-1}_1f_1,\ldots,r^{-1}_nf_n}{f_0}\Big)&=\Spec^{\mathrm{val}}\mathcal{A}\Big(\frac{1^{-1}g_1,\ldots,1^{-1}g_n}{f_0}\Big)\\
&\cong\Spa F_1(\mathcal{A})\Big(\frac{g_1,\ldots,g_n}{f_0}\Big),
\end{split}
\end{equation*}
where $g_i=a_if_i$ for $i=1,\ldots,n$.
This means that rational subsets on both sides are the same; see \S\ref{subsub-adicspectrum}.
\end{proof}

\begin{cor}\label{cor-vsbr-berkovichaffinoidalgebras2}
For a Berkovich's $K$-affinoid algebra $\mathcal{A}$, the valuative space $\Spec^{\mathrm{val}}\mathcal{A}$ is reflexive\index{valuative!valuative topological space@--- (topological) space!reflexive valuative topological space@reflexive --- ---} {\rm ({\bf \ref{ch-pre}}.\ref{dfn-reflexivevaluativespaces})}.
\end{cor}

\begin{proof}
Take an extension $L/K$ of Banach fields with $|L^{\times}|=\R_+$.
By \ref{cor-vsbr-berkovichaffinoidalgebras1}, $\Spec^{\mathrm{val}}\mathcal{A}_L\rightarrow\Spec^{\mathrm{val}}\mathcal{A}$ is surjective, and $\Spec^{\mathrm{val}}\mathcal{A}_L$ is reflexive by \ref{prop-vsbr-berkovichaffinoidalgebras1}, \ref{thm-rigidgeomadicspacesff1}, and \ref{thm-spectralfunctorseparatedquotients}.
\end{proof}

\begin{cor}\label{cor-vsbr-berkovichaffinoidalgebras3}
Suppose $K$ has a non-trivial valuation.
Then, for a Berkovich's $K$-affinoid algebra $\mathcal{A}$, $\Spa F_1(\mathcal{A})$ is reflexive.
\end{cor}

\subsubsection{Reflexivity of valuative spectrum}\label{subsub-vsbr-reflexivity}
\begin{thm}[{\rm Local affineness}]\label{thm-vsbr-reflexivity}
Let $\mathcal{A}$ be an $\R_+$-affinoid ring, $X=\Spec^{\mathrm{val}}\mathcal{A}$, and $U\hookrightarrow V$ an open immersion of quasi-compact open subsets of $X$.
Then, for any $x\in V$, there exists a quasi-compact open neighborhood $W_x$ of $x$ in $V$ such that $W_x$ and $U\cap W_x$ are $\R_+$-rational subdomains of $X$.
\end{thm}

To show the theorem, we need the following lemma.
\begin{lem}\label{lem-vsbr-reflexivity}
Let $X$ be coherent sober space, $U,V\subseteq X$ quasi-compact open subsets, and $x\in X$.
Suppose $U\cap G_x=V\cap G_x$, where $G_x$ denotes the set of all generizations of $x$ in $X$.
Then there exists a quasi-compact open neighborhood $W_x$ of $x$ such that $U\cap W_x=V\cap W_x$.
\end{lem}

\begin{proof}
By considering $U\cap V$ instead of $U$, we may assume $U\subseteq V$.
Set $C=X\setminus U$.
Suppose the assertion is false. 
Then $C\cap V\cap W\neq\emptyset$ for any quasi-compact open neighborhood $W$ of $x$.
Since $C$ is coherent and sober, this implies $G_x\neq\emptyset$ ({\bf \ref{ch-pre}}.\ref{thm-projlimcohsch1} (2)), which is absurd.
\end{proof}

\begin{proof}[Proof of Theorem {\rm \ref{thm-vsbr-reflexivity}}]
We may assume $V=X$.
Set $x=(z,F_x)$.
The set $G_x$ of generizations of $x$ in $X$ is canonically identified with $\Spec^{\mathrm{val}}(K_z,F_x)$ regarded as a subset of $X$.
Since $U\cap G_x$ is a quasi-compact open subset of $G_x$, there exists $g,f\in\mathcal{A}^{\mathrm{B}}$ and $r>0$ such that the basic subset(see \S\ref{subsub-vsbr-valuativespectrum}) $B(g,f,r)$ contains $x$ and $U\cap G_x=B(f,g,r^{-1})\cap G_x$.
By \ref{lem-vsbr-reflexivity}, there exists an $\R_+$-rational neighborhood $W_x$ of $x$ contained in $B(g,f,r)$ such that $U\cap W_x=B(f,g,r^{-1})\cap W_x$.
Notice that, since $W_x\subseteq B(g,f,r)$, $B(f,g,r^{-1})\cap W_x$ is an $\R_+$-rational subdomain of $W_x$, and hence of $X$.
\end{proof}

\begin{rem}\label{rem-vsbr-reflexivity0}{\rm 
The proof indicates that, in the situation as in \ref{thm-vsbr-reflexivity}, one can take $W_x=\Spec^{\mathrm{val}}\mathcal{B}$ such that $U\cap W_x$ is an $\R_+$-rational subdomain of the form $\Spec^{\mathrm{val}}\mathcal{B}(\frac{r\cdot 1}{f})$ for some $f\in(F_{\mathcal{B}})_r$.}
\end{rem}

\begin{prop}\label{prop-vsbr-reflexivity}
Let $\mathcal{A}$ be an $\R_+$-affinoid ring.
Then the following conditions are equivalent.
\begin{itemize}
\item[{\rm (a)}] $X=\Spec^{\mathrm{val}}\mathcal{A}$ is reflexive.
\item[{\rm (b)}] For any $\R_+$-rational localization $\mathcal{B}=(B,F_{\mathcal{B}})$ of $\mathcal{A}$, $\Spec^{\mathrm{val}}(B,F^{\Sp}_B)\rightarrow\Spec^{\mathrm{val}}\mathcal{B}$ is bijective.
\end{itemize}
\end{prop}

\begin{proof}
Suppose (1) holds, and let $\mathcal{B}=(B,F_{\mathcal{B}})$ be an $\R_+$-rational localization of $\mathcal{A}$.
By assumption, the quasi-compact open subspace $Y=\Spec^{\mathrm{val}}\mathcal{B}$ of $X$ is also reflective.
Since $[Y]=\mathscr{M}(B)=[\Spec^{\mathrm{val}}(B,F^{\Sp}_B)]$, $Y^{\mathrm{ref}}=(\Spec^{\mathrm{val}}(B,F^{\Sp}_B))^{\mathrm{ref}}\hookrightarrow\Spec^{\mathrm{val}}(B,F^{\Sp}_B)$.
This means that $Y^{\mathrm{ref}}\hookrightarrow Y$ factors as $Y^{\mathrm{ref}}\hookrightarrow\Spec^{\mathrm{val}}(B,F^{\Sp}_B)\hookrightarrow Y$, by which we have $Y\cong\Spec^{\mathrm{val}}(B,F^{\Sp}_B)$.

Conversely, suppose (2) holds, and take an open immersion $U\hookrightarrow V$ of quasi-compact open subsets of $X$ such that $[U]=[V]$.
By \ref{thm-vsbr-reflexivity}, one has a finite open covering $\{V_{\alpha}\}_{\alpha\in L}$ of $V$ such that $V_{\alpha}$ and $U\cap V_{\alpha}$ are $\R_+$-rational subdomains of $X$ for any $\alpha\in L$.
We may assume that $U_\cap V_{\alpha}\hookrightarrow V_{\alpha}$ is induced from a homomorphism $\mathcal{B}_{\alpha}\rightarrow\mathcal{B}'_{\alpha}$ of $\R_+$-rational localizations of $\mathcal{A}$.
Our assumption $[U]=[V]$ implies $[\Spec^{\mathrm{val}}\mathcal{B}_{\alpha}]=[\Spec^{\mathrm{val}}\mathcal{B}'_{\alpha}]$.
It follows from the assumption (1) that $\mathcal{B}_{\alpha}$ and $\mathcal{B}'_{\alpha}$ are of spectral type, viz., $F_{\mathcal{B}_{\alpha}}=F^{\Sp}_{B_{\alpha}}$ and similarly for $F_{\mathcal{B}'_{\alpha}}$.
From this and $[\Spec^{\mathrm{val}}\mathcal{B}_{\alpha}]=[\Spec^{\mathrm{val}}\mathcal{B}'_{\alpha}]$ it follows that $\Spec^{\mathrm{val}}\mathcal{B}_{\alpha}=\Spec^{\mathrm{val}}\mathcal{B}'_{\alpha}$ for any $\alpha\in L$ (cf \ref{rem-vsbr-reflexivity} below), thereby $U=V$.
\end{proof}

Note that the proposition, together with \ref{thm-vsbr-berkovichaffinoidalgebras}, gives another proof of \ref{cor-vsbr-berkovichaffinoidalgebras2}, the reflexiveness of $\Spec^{\mathrm{val}}\mathcal{A}$ for a Berkovich's $K$-affinoid algebra $\mathcal{A}$.

\begin{rem}\label{rem-vsbr-reflexivity}{\rm 
For an $\R_+$-rational localization $\mathcal{A}=(A,F_{\mathcal{A}})\rightarrow\mathcal{A}'=(A',F_{\mathcal{A}'})$, the following conditions are equivalent.
\begin{itemize}
\item[{\rm (a)}] $A^{\Sp}\cong A^{\prime\Sp}$;
\item[{\rm (b)}] $\mathscr{M}(A)$ is homeomorphic to $\mathscr{M}(A')$;
\item[{\rm (c)}] $\Spec^{\mathrm{val}}(A',F^{\Sp}_{A'})\rightarrow\Spec^{\mathrm{val}}(A,F^{\Sp}_A)$ is bijective.
\end{itemize}}
\end{rem}

\subsection{Non-archimedean analytic space of Banach type}\label{sub-vsbr-NAASBT}
\subsubsection{Admissible site of $\R_+$-affinoid rings}\label{subsub-vsbr-admsiteRaffinoidrings}
Let $\mathcal{A}=(A,F_{\mathcal{A}})$ be an $\R_+$-affinoid ring.
We denote by $R_{\mathcal{A}}$ the category of $\R_+$-rational localizations of $\mathcal{A}$ (see \S\ref{subsub-vsbr-Rplusaffinoidrings}).
Morphisms in $R_{\mathcal{A}}$ (resp.\ $R^{\opp}_{\mathcal{A}}$) are all epimorphisms (resp.\ monomorphisms) (\ref{cor-vsbr-powerboundedfiltration1}).
When we regard an object $\mathcal{B}$ of $R_{\mathcal{A}}$ as an object of $R^{\opp}_{\mathcal{A}}$, we denote it by $S(\mathcal{B})$; the same convention is applied also to arrows in $R_{\mathcal{A}}$.

We define a notion of coverings in $R^{\opp}_{\mathcal{A}}$ as follows: a collection of arrow $\{S(\mathcal{B}_{\alpha})\rightarrow S(\mathcal{B})\}_{\alpha\in L}$ is a covering of $S(\mathcal{B})$ if 
$$
\Spec^{\mathrm{val}}\mathcal{B}=\bigcup_{\alpha\in L}\Spec^{\mathrm{val}}\mathcal{B}_{\alpha}
$$
holds.
It is straightforward to check that this notion of covering defines a Grothendieck topology $J_{\mathcal{A}}$ on $R^{\opp}_{\mathcal{A}}$.
We set 
$$
D_{\mathcal{A}}=(R^{\opp}_{\mathcal{A}},J_{\mathcal{A}}),
$$
and call it the {\it admissible site} of $\mathcal{A}$.

\begin{prop}\label{prop-vsbr-admsiteRaffinoidrings1}
There exists an equivalence of the associated topoi 
$$
D^{\sim}_{\mathcal{A}}\cong\top(\Spec^{\mathrm{val}}\mathcal{A})
$$
$($see {\rm {\bf \ref{ch-pre}}, \S\ref{subsub-toposasstopsp}} for the notation$)$.
\end{prop}

\begin{proof}
Let $R$ be the set of all $\R_+$-rational subdomains of $\Spec^{\mathrm{val}}\mathcal{A}$ (\S\ref{subsub-vsbr-valuativespectrum}).
The set $R$ is partially ordered by the inclusion order, and is regarded as a category. 
We further regard the poset $R$ as a site in a standard manner; a covering of $a\in R$ is a finite collection $\{b_1,\ldots,b_r\}$ of elements in $R$ such that $a=\sup\{b_1,\ldots,b_r\}$; note that $R$ is closed under finite intersections. 
Since $R$ generates the topology of $\Spec^{\mathrm{val}}\mathcal{A}$ (\ref{dfn-vsbr-valuativespectrum2}) and any $U\in R$ is coherent, the morphism of site $\Spec^{\mathrm{val}}\mathcal{A}\rightarrow R$ induces an equivalence of the associated topoi.
Clearly, there exists a morphism of sites $j\colon R\rightarrow D_{\mathcal{A}}$ associated to the functor $D_{\mathcal{A}}\rightarrow R$ that maps $S(\mathcal{B})$ to $\Spec^{\mathrm{val}}\mathcal{B}$.
Suppose $S(\mathcal{B})\rightarrow S(\mathcal{B}')$ induces $\Spec^{\mathrm{val}}\mathcal{B}\cong\Spec^{\mathrm{val}}\mathcal{B}'$.
Due to this and the fact that $S(\mathcal{B})\rightarrow S(\mathcal{B}')$ is a monomorphism (whence having $S(\mathcal{B})\times_{S(\mathcal{B}')}S(\mathcal{B})\cong S(\mathcal{B})$), $S(\mathcal{B})\rightarrow S(\mathcal{B}')$ is a covering arrow, and hence, for any sheaf $\mathscr{F}$ on $D_{\mathcal{A}}$, we have $\mathscr{F}(S(\mathcal{B}))\cong\mathscr{F}(S(\mathcal{B}'))$ by the definition of the Grothendieck topology $J_{\mathcal{A}}$.
Thus $j^{-1}$ induces an equivalence of the associated topoi.
Then we obtain the desired equivalence of topoi by composition of the equivalences obtained above.
\end{proof}

We define a presheaf $\til{\O}_{\mathcal{A}}$ on the site $D_{\mathcal{A}}$ by
$$
\til{\O}_{\mathcal{A}}(S(\mathcal{B}))=\mathcal{B}^{\mathrm{B}}.
$$
This is a presheaf of Banach rings on $D_{\mathcal{A}}$, and we denote by $\O_{\mathcal{A}}$ the sheaf of rings given by sheafification of $\til{\O}_{\mathcal{A}}$ seen as a presheaf of rings.
By \ref{prop-vsbr-admsiteRaffinoidrings1}, one can regard $\O_{\mathcal{A}}$ as a sheaf on $\Spec^{\mathrm{val}}\mathcal{A}$.

For $x\in\Spec^{\mathrm{val}}\mathcal{A}$, the stalk $\mathcal{A}_x=\O_{\mathcal{A},x}$ at $x$ is a ring that allows the following description.
Let $S_x$ be the category of $\R_+$-rational localization $\mathcal{B}$ such that $x$ lies in the image of $\Spec^{\mathrm{val}}\mathcal{B}$.
Then, by what we have seen in \S\ref{subsub-vsbr-valuativespectrum}, $S_x$ is directed, and 
$$
\mathcal{A}_x=\varinjlim_{\mathcal{B}\in\obj(S_x)}\mathcal{B}^{\mathrm{B}},
$$
where inductive limit is taken in the category of rings.
Notice that the ring $\mathcal{A}_x$ comes with a canonical ring homomorphism 
$$
\mathcal{A}_x\longrightarrow\mathcal{H}(z),
$$
where $x=(z,F)$.
Let us denote by $\m_x$ the kernel of this map, and set
$$
k(x)=\mathcal{A}_x/\m_x.
$$

\begin{prop}\label{prop-vsbr-admsiteRaffinoidrings2}
Let $\mathcal{A}$ be an $\R_+$-affinoid ring, and $x=(z,F)\in\Spec^{\mathrm{val}}\mathcal{A}$.

{\rm (1)} An element $h\in\mathcal{A}_x$ is invertible if and only if its image in $\mathcal{H}(z)$ is non-zero. In particular, $(\mathcal{A}_x,\m_x)$ is a local ring, and thus $(\Spec^{\mathrm{val}}\mathcal{A},\O_{\mathcal{A}})$ is a locally ringed space.

{\rm (2)} The field $k(x)$, view as a subfield of $\mathcal{H}(z)$, is dense in $\mathcal{H}(z)$.
\end{prop}

\begin{proof}
(1) Let $\mathcal{B}$ be an object of $S_x$, and let $H\in\mathcal{B}$ represent $h\in\mathcal{A}_x$.
Replacing $\mathcal{B}$ further by a rational localization around $x$ if necessary, we may assume that $|H(w)|$, viewed as a real-valued function of $w\in\mathscr{M}(\mathcal{B})$, is non-zero everywhere on $\mathscr{M}(\mathcal{B})$.
Then $H$ is invertible in $\mathcal{B}$, and hence $h$ is invertible in $\mathcal{A}_x$.

(2) Let $A'$ be the image of $\mathcal{A}\rightarrow\mathcal{H}(z)$.
Then, by definition of the complete residue field $\mathcal{H}(z)$, the fractional field $K'$ of $A'$ is dense in $\mathcal{H}(z)$.
Clearly, $K'$ is a subfield of $k(x)$, whence the claim.
\end{proof}

Let $x=(z,F)\in\Spec^{\mathrm{val}}\mathcal{A}$.
By \ref{prop-vsbr-admsiteRaffinoidrings2} (2), there exists a unique filtration $F_x$ on $k(x)$ of which the completion with respect to the seminorm $\|\cdot\|_z$ coincides with $F$, viz., there exists a unique filtered valuation field of the form $(k(x),F_x)$ of which the completion is the given $(\mathcal{H}(z),F)$.
Define a filtration $F_{\mathcal{A}_x}$ on $\mathcal{A}_x$ by
$$
(F_{\mathcal{A}_x})_r=\textrm{the preimage of $(F_x)_r$ under $\mathcal{A}_x\rightarrow k(x)$}.
$$
for $r>0$, which obviously gives rise to a filtered ring $(\mathcal{A}_x,F_{\mathcal{A}_x})$.

\begin{prop}\label{prop-vsbr-admsiteRaffinoidrings3}
{\rm (1)} The filtered ring $(\mathcal{A}_x,F_{\mathcal{A}_x})$ coincides with 
$$
\varinjlim_{\mathcal{B}\in\obj(S_x)}\mathcal{B},
$$
where the inductive limit is taken in the category of filtered rings.

{\rm (2)} For $r>0$,
$$
(F_{\mathcal{A}})_r=\{f\in\mathcal{A}: \textrm{$f\in (F_{\mathcal{A}_x})_r$ for any $x\in\Spec^{\mathrm{val}}\mathcal{A}$}\}.
$$
\end{prop}

\begin{proof}
First notice that the category of filtered rings has small inductive limits, which commutes with the inductive limits in the category of rings by the forgetful functor.
Let $(\mathcal{A}_x,F')$ be the filtered inductive limit as in (1).
By definition, we have $F'\subseteq F_{\mathcal{A}_x}$.
We first check $F'_{<r}=(F_{\mathcal{A}_x})_{<r}$ for $r>0$.
Let $h\in(F_{\mathcal{A}_x})_{<r}$, and take $\mathcal{B}$ from $S_x$ and $H\in\mathcal{B}$ such that $h$ is represented by $H$.
Since $|H(z)|<r$, there is an overconvergent open neighborhood $V$ of $z$ in $\Spec^{\mathrm{val}}\mathcal{A}$ such that $|H(w)|<r$ for all $w\in\sep_{\mathcal{A}}(V)\subseteq\mathscr{M}(\mathcal{A})$.
Replacing $\mathcal{B}$ by another one in $S_x$ with smaller valuative spectrum, we may assume that $\Spec^{\mathrm{val}}\mathcal{B}$ is contained in $V$, and hence that $H\in(F_{\mathcal{B}})_{<r}=(F^{\Sp}_{\mathcal{B}})_{<r}$.
Thus we have $h\in F'_{<r}$, as desired.

To show $F'_r=(F_{\mathcal{A}_x})_r$ for $r>0$, take $\alpha\in(F_{\mathcal{A}_x})_r$.
Then, by what we have seen in the proof of \ref{prop-vsbr-admsiteRaffinoidrings2} (2), $[\alpha]_r\in\mathrm{Gr}_{F_{\mathcal{A}_x}}\mathcal{A}_x$ is represented by $\ovl{a}/\ovl{b}$, where $a,b\in A$, $\ovl{a}$ and $\ovl{b}$ denote the image of $a$ and $b$ in $\mathrm{Gr}\,\mathcal{A}$, and $\ovl{b}\neq 0$.
Let $s$ be a positive real number such that $s>1/|b(z)|$.
Then $a/b$ in 
$$
\mathcal{B}=\mathcal{A}\Big(\frac{r^{-1}a,s^{-1}1}{b}\Big)
$$
represents $\ovl{a}/\ovl{b}$, and $\alpha-a/b\in F'_{<r}$.
Hence we have $\alpha\in F'_r$, as desired, and thus (1) is proved.
(2) follows immediately from \ref{thm-vsbr-filtrationintegral}.
\end{proof}

The following proposition determines $\R_+$-rational localizations $\mathcal{A}\rightarrow\mathcal{B}$ with the same valuative spectrum, and thus clarify the dependence of $\R_+$-rational subdomains on their presentations as the valuative spectrums of $\R_+$-rational localizations. 

\begin{prop}\label{prop-vsbr-admsiteRaffinoidrings4}
Let $\mathcal{A}\rightarrow\mathcal{B}$ be an $\R_+$-rational localization.
Then $\Spec^{\mathrm{val}}\mathcal{B}=\Spec^{\mathrm{val}}\mathcal{A}$ if and only if $\mathcal{B}$ is isomorphic to an $\R_+$-rational localization of the form 
$$
\mathcal{A}\Big(\frac{r^{-1}_1g_1,\ldots,r^{-1}_ng_n}{1}\Big)
$$
for $g_i\in(F_{\mathcal{A}})_{r_i}$ $(i=1,\ldots,n)$.
\end{prop}

\begin{proof}
Only the ``if'' part calls for a proof.
Set 
$$
\mathcal{B}=\mathcal{A}\Big(\frac{r^{-1}_1f_1,\ldots,r^{-1}_nf_n}{f_0}\Big),
$$
as an $\mathcal{A}$-affinoid algebra.
The element $f_0$ takes no zero on $\mathscr{M}(\mathcal{B})$, and hence on $\mathscr{M}(\mathcal{A})$, due to the assumption.
This implies that $f_0$ is invertible in $\mathcal{A}$.
Set $g_i=f^{-1}_0f_i$ for $i=1,\ldots,n$.
The image of $g_i$ in $\mathcal{B}$ is in $(F_{\mathcal{B}})_{r_i}$ for $i=1,\ldots,n$, and thus $\mathcal{B}$ is isomorphic to the $\R_+$-rational localization as above.
By the equality $\Spec^{\mathrm{val}}\mathcal{B}=\Spec^{\mathrm{val}}\mathcal{A}$ and \ref{prop-vsbr-admsiteRaffinoidrings3} (2), we have $g_i\in(F_{\mathcal{A}})_{r_i}$ for $i=1,\ldots,n$.
\end{proof}

\subsubsection{Sheaf condition of Banach type}\label{subsub-vsbr-sheafconditionBanachtype}
Let $\mathcal{A}=(\mathcal{A},F_{\mathcal{A}})$ be an $\R_+$-affinoid ring, and consider the presheaf $\til{\O}_{\mathcal{A}}$ on $D_{\mathcal{A}}$ as in \S\ref{subsub-vsbr-admsiteRaffinoidrings}, which is a presheaf of Banach rings with all restriction maps being bounded homomorphisms.

\begin{dfn}\label{dfn-vsbr-sheafconditionBanachtype1}{\rm 
We say that $\til{\O}_{\mathcal{A}}$ satisfies the {\it sheaf condition of Banach type}, or that the $\R_+$-affinoid ring $\mathcal{A}$ is {\it $\R_+$-sheafy}, if, for any finite covering $\{S(\mathcal{B}_{\alpha})\}_{\alpha\in L}$ of $S(\mathcal{B})$, 
$$
\til{\O}_{\mathcal{A}}(S(\mathcal{B}))\longrightarrow\ker\Big(\xymatrix{{\displaystyle \prod_{\alpha\in L}\til{\O}_{\mathcal{A}}(S(\mathcal{B}_{\alpha}))}\ar@<.5ex>[r]\ar@<-.5ex>[r]&{\displaystyle \prod_{\alpha,\beta\in L}\til{\O}_{\mathcal{A}}(S(\mathcal{B}_{\alpha}\widehat{\otimes}_{\mathcal{B}}\mathcal{B}_{\beta}))}}\Big)
$$
is a bounded isomorphism of Banach rings.
(Notice that the right-hand side is a closed set of the Banach ring $\prod_{\alpha\in L}\til{\O}_{\mathcal{A}}(S(\mathcal{B}_{\alpha}))$, and hence is a Banach ring.)}
\end{dfn}

\begin{rem}\label{rem-vsbr-sheafconditionBanachtype1}{\rm 
The `$\R_+$-sheafy' condition is, in our situation, equivalent to `sheafy' in \cite{Kedlaya}, 3.19.
Note that, as we have seen in the end of \S\ref{subsub-vsbr-adicspectrum}, affinoid Banach rings in the sense of Kedlaya are nothing but $\R_+$-affinoid rings of the form $\mathcal{A}=(A,F_{\mathcal{A}})$ such that $F_{\mathcal{A}}\subseteq F^o_{\mathcal{A}}$.
As the next proposition shows, the last condition in our situation is rather a consequence of $\R_+$-sheafiness.}
\end{rem}

We denote by 
$$
\R_+\textrm{-}\mathbf{Aff}^{\mathrm{Sh}}
$$
the full subcategory of $\R_+\textrm{-}\mathbf{Aff}$, the category of $\R_+$-affinoid rings, consisting of $\R_+$-sheafy $\R_+$-affinoid rings.

\begin{prop}\label{prop-vsbr-sheafconditionBanachtype1}
Let $\mathcal{A}=(\mathcal{A},F_{\mathcal{A}})$ be an $\R_+$-affinoid ring, and suppose $\mathcal{A}$ is $\R_+$-sheafy.

{\rm (1)} The inclusion $F_{\mathcal{A}}\subseteq F^o_{\mathcal{A}}$ holds.

{\rm (2)} The functor $\mathcal{B}\mapsto\Spec^{\mathrm{val}}\mathcal{B}$ gives an equivalence from the site of $\R_+$-rational subdomains on $\Spec^{\mathrm{val}}\mathcal{A}$ to the site $D_{\mathcal{A}}$.
\end{prop}

\begin{proof}
(1) Let $f\in(F_{\mathcal{A}})_r$.
Then $\Spec^{\mathrm{val}}\mathcal{A}=\Spec^{\mathrm{val}}\mathcal{A}\dl r^{-1}f\dr$.
Since, by the sheaf condition of Banach type, $\mathcal{A}\rightarrow\mathcal{A}\dl r^{-1}f\dr$ is a bounded homomorphism.
Hence $f$ is $r$-power bounded.

(2) By (1) and \ref{prop-cor-vsbr-powerboundedfiltration3}, we have $F_{\mathcal{B}}\subseteq F^o_{\mathcal{B}}$ for any $\R_+$-rational localization $\mathcal{B}$.
Thus for any $r>0$ and $f\in(F_{\mathcal{B}})_r$, the morphisms $\mathcal{B}\rightarrow\mathcal{B}(\frac{r^{-1}f}{1})$ is a bounded homomorphism.
Hence for any rational localization $\mathcal{B}'$ of $\mathcal{B}$, $\Spec^{\mathrm{val}}\mathcal{B}'=\Spec^{\mathrm{val}}\mathcal{B}$ implies the existence of bounded isomorphism $\mathcal{B}\cong\mathcal{B}'$ by \ref{prop-vsbr-admsiteRaffinoidrings4}.
\end{proof}

By \ref{prop-vsbr-sheafconditionBanachtype1} (2), we see that an $\R_+$-rational subdomain $U$ of $\Spec^{\mathrm{val}}\mathcal{A}$ determines an $\R_+$-rational localization $\mathcal{B}$ up to isomorphism.
In particular, for an $\R_+$-sheafy $\R_+$-affinoid ring $\mathcal{A}$, one can define a presheaf $\O_{\mathcal{A}}$ on $\Spec^{\mathrm{val}}\mathcal{A}$ such that, for any $\R_+$-rational subdomain $U=\Spec^{\mathrm{val}}\mathcal{B}$, 
$$
\O_{\mathcal{A}}(U)=\mathcal{B}^{\mathrm{B}},
$$
and for any open subset $V$, 
$$
\O_{\mathcal{A}}(V)=\varprojlim_U\O_{\mathcal{A}}(U),
$$
where $U$ runs through all $\R_+$-rational subdomains contained in $V$.
Note that this presheaf $\O_{\mathcal{A}}$ is a sheaf, and is nothing but the one corresponding to the sheaf $\til{\O}_{\mathcal{A}}$ on the site $D_{\mathcal{A}}$.

An important example of $\R_+$-sheafy $\R_+$-affinoid rings is Berkovich's $K$-affinoid algebras (uniquely regarded as $\R_+$-affinoid rings as in \S\ref{subsub-vsbr-berkovichaffinoidalgebras}).
\begin{prop}\label{prop-vsbr-BerkovichRmetrizedanalyticspace1}
Any Berkovich's $K$-affinoid algebra $\mathcal{A}$ is $\R_+$-sheafy {\rm (\ref{dfn-vsbr-sheafconditionBanachtype1})}.
\end{prop}

\begin{proof}
Note first that $X=\Spec^{\mathrm{val}}\mathcal{A}$ is reflexive (\ref{cor-vsbr-berkovichaffinoidalgebras2}), and that, for any Berkovich's $K$-affinoid algebra $\mathcal{B}$,  $F_{\mathcal{B}}=F^o_{\mathcal{B}}=F^{\Sp}_{\mathcal{B}}$ due to \ref{thm-vsbr-berkovichaffinoidalgebras}.
Hence we have a well-defined presheaf $\til{\O}$ on the site of rational subdomains and finite coverings by rational subdomains of $\mathscr{M}(\mathcal{A})$ defined by, for any $\R_+$-rational subdomain $U=\Spec^{\mathrm{val}}\mathcal{B}$ of $X$, $\til{\O}([U])=\til{\O}_{\mathcal{A}}(U)=\mathcal{B}^{\mathrm{B}}$.
The sheaf condition for $\til{\O}$ as in \cite{Berk1}, 2.2.5 is nothing but the sheaf condition of Banach type for $\til{\O}_{\mathcal{A}}$.
\end{proof}

For more information on the sheaf condition, see \cite{Kedlaya}, \S8.

\subsubsection{Metrized Banach ringed spaces}\label{subsub-vsbr-MBRS}
\begin{dfn}\label{dfn-vsbr-MBRS}{\rm 
(1) A locally ringed space $(X,\O_X)$ is a {\it Banach ringed space} if it enjoys the following properties:
\begin{itemize}
\item[{\rm (a)}] $X$ is a sober locally coherent space;
\item[{\rm (b)}] for a coherent open subset $U\subseteq X$, $\O_X(U)$ is a Banach ring;
\item[{\rm (c)}] the restriction map $\O_X(V)\rightarrow\O_X(U)$ for $U\subseteq V$, where $U$ and $V$ are coherent open subsets, is a bounded homomorphism of Banach rings;
\item[{\rm (d)}] for any finite open covering $U=\bigcup_{\alpha\in L}U_{\alpha}$ of a coherent open subset $U$ by coherent open subsets, the isomorphism
$$
\O(U)\stackrel{\sim}{\longrightarrow}\ker\Big(\xymatrix{{\displaystyle \prod_{\alpha\in L}\O(U_{\alpha})}\ar@<.5ex>[r]\ar@<-.5ex>[r]&{\displaystyle \prod_{\alpha,\beta\in L}\O(U_{\alpha}\cap U_{\beta})}}\Big)
$$
is a bounded isomorphism. (Note that the right-hand ring is a closed subring of $\prod_{\alpha\in L}\O(U_{\alpha})$.)
\end{itemize}

(2) A {\it morphism} $f\colon (X,\O_X)\rightarrow (Y,\O_Y)$ of Banach ringed spaces is a morphism of locally ringed spaces, with locally quasi-compact underlying continuous mapping, having the following property: 
\begin{itemize}
\item[{\rm (e)}] for a coherent open subset $U$ of $X$ and a coherent open subset $V$ of $Y$ such that $U\subseteq f^{-1}(V)$, the induced homomorphism $\O_Y(V)\rightarrow\O_X(U)$ is a bounded homomorphism.
\end{itemize}}
\end{dfn}

By definition, any open subset of a Banach ringed space is naturally a Banach ringed space, which we call an open subspace.

\begin{dfn}\label{dfn-vsbr-MBRSval}{\rm 
Let $X=(X,\O_X)$ be a Banach ringed space satisfying the following conditions:
\begin{itemize}
\item[{\rm (a)}] the underlying topological space $X$ is a valuative space;
\item[{\rm (b)}] if $x\in X$ is a generization of $y\in X$, then the natural ring homomorphism $\O_{X,y}\rightarrow\O_{X,x}$ is local.
\end{itemize}
An {\it $\R_+$-valuation $v$ of $X=(X,\O_X)$} is a family of $\R_+$-filtered valuation fields $v=\{v_x=(k(x),F_x)\}_{x\in X}$, where $k(x)$ is the residue field at $x$ for $x\in X$, enjoying the following properties:
\begin{itemize}
\item[{\rm (c)}] if $x$ is a maximal point, then $v_x$ is a filtered valuation field of maximal type;
\item[{\rm (d)}] if $x$ is a generization of $y\in X$, then the local homomorphism $\O_{X,y}\rightarrow\O_{X,x}$ induces a filtered homomorphism $v_y=(k(y),F_y)\rightarrow v_x=(k(x),F_x)$.
\end{itemize}
A Banach ringed space satisying (a) and (b) equipped with an $\R_+$-filtered valuation $\mathscr{X}=((X,\O_X),v)$ is called an {\it $\R_+$-metrized Banach ringed space}.}
\end{dfn}

In this situation, the structure sheaf $\O_X$ has the multiplicative filtration $F_{\mathscr{X}}$ defined as follows.
For any open subset $U\subseteq X$ and $r>0$, 
$$
F_{\mathscr{X}}(U)_r=\{f\in\O_X(U): (f_x\ \mathrm{mod}\ \m_{X,x})\in(F_x)_r\ \textrm{for any $x\in U$}\}.
$$

By definition, any open subset of an $\R_+$-metrized Banach ringed space is naturally an $\R_+$-metrized Banach ringed space, which we call an open subspace.

\begin{dfn}\label{dfn-vsbr-MBRSval2}{\rm 
A {\it morphism} 
$$
f\colon\mathscr{X}=((X,\O_X),v)\longrightarrow\mathscr{X}'=((X',\O_{X'}),v')
$$
of $\R_+$-metrized Banach ringed spaces is a morphism of Banach ringed spaces $f\colon (X,\O_X)\rightarrow(Y,\O_Y)$ satisfying the following conditions:
\begin{itemize}
\item[{\rm (a)}] the underlying continuous mapping of $f$ is valuative\index{valuative!valuative map@--- map} ({\bf \ref{ch-pre}}.\ref{dfn-valuativemaps});
\item[{\rm (b)}] for $x\in X$, the local homomorphism $\O_{X',f(x)}\rightarrow\O_{X,x}$ induces a dominating filtered homomorphism $(k(f(x)),F_{f(x)})\rightarrow(k(x),F_x)$. 
\end{itemize}}
\end{dfn}

Note that a morphism $f\colon\mathscr{X}\rightarrow\mathscr{X}'$ of $\R_+$-metrized Banach ringed spaces respects the filtrations $F_{\mathscr{X}}$ and $F_{\mathscr{X}'}$.

\begin{thm}\label{thm-vsbr-MBRSval}
Let $\mathcal{A}$ be an $\R_+$-affinoid ring, and suppose it is $\R_+$-sheafy.
Let $X=(\Spec^{\mathrm{val}}\mathcal{A},\O_{\mathcal{A}})$ be the resulting local ringed space.

{\rm (1)} The locally ringed space $X$ is a Banach ringed space.
Moreover, there exists a natural $\R_+$-valuation $v_{\mathcal{A}}$ of $X$, which make the pair $((X=\Spec^{\mathrm{val}}\mathcal{A},\O_{\mathcal{A}}),v_{\mathcal{A}})$ an $\R_+$-metrized Banach ringed space.

{\rm (2)} The functor 
$$
\mathcal{A}\longmapsto((\Spec^{\mathrm{val}}\mathcal{A},\O_{\mathcal{A}}),v_{\mathcal{A}}),
$$
from the category of $\R_+$-sheafy $\R_+$-affinoid rings to the category of $\R_+$-metrized Banach ringed spaces, is fully faithful.
\end{thm}

\begin{proof}
(1) First, note that $X$ is a locally ringed space due to \ref{prop-vsbr-admsiteRaffinoidrings2}.
To show (1), we first check that $X$ is a Banach ringed space.
The condition (a) of \ref{dfn-vsbr-MBRS} is clear, since $X$ is coherent.
For any coherent open subset $U\subseteq X$, choose a finite open covering $\{U_{\alpha}\}_{\alpha\in L}$ by $\R_+$-rational subdomains, and define
$$
\O(U)=\ker\Big(\xymatrix{{\displaystyle \prod_{\alpha\in L}\O_{\mathcal{A}}(U_{\alpha})}\ar@<.5ex>[r]\ar@<-.5ex>[r]&{\displaystyle \prod_{\alpha,\beta\in L}\O_{\mathcal{A}}(U_{\alpha}\cap U_{\beta})}}\Big).
$$
Then $\O(U)$ is a Banach ring, as a closed subring of $\prod_{\alpha\in L}\O_{\mathcal{A}}(U_{\alpha})$, and coincides with $\O_{\mathcal{A}}(U)$ if $U$ is an $\R_+$-rational subdomain.
It is straightforward to see that, in general, $\O(U)$ is independent of the choice of the finite covering $\{U_{\alpha}\}_{\alpha\in L}$, and that the restriction map $\O(V)\rightarrow\O(U)$ for $U\subseteq V$ is bounded, due to the sheaf condition of Banach type.
Thus, so far, we have shown that the conditions (b), (c), (d) of \ref{dfn-vsbr-MBRS} are satisfied, and thus that $X$ is a Banach ringed space.

Let us show the rest of (1).
The condition (a) of \ref{dfn-vsbr-MBRSval} follows from \ref{thm-vsbr-toppropertyvalspec}.
The condition (b) is satisfied due to \ref{prop-vsbr-admsiteRaffinoidrings2} (1).
Now, for each $x\in\Spec^{\mathrm{val}}\mathcal{A}$, we attach the filtered valuation field $v_x=(k(x),F_x)$ given in \S\ref{subsub-vsbr-admsiteRaffinoidrings}.
Then $v_{\mathcal{A}}=\{v_x\}_{x\in X}$ satisfies the conditions (c) and (d) in \ref{dfn-vsbr-MBRSval}.
Hence $((X=\Spec^{\mathrm{val}}\mathcal{A},\O_{\mathcal{A}}),v_{\mathcal{A}})$ an $\R_+$-metrized Banach ringed space.

(2) To check that any morphism $\varphi\colon\mathcal{B}\rightarrow\mathcal{A}$ of $\R_+$-sheafy $\R_+$-affinoid rings induces a morphism of $\R_+$-metrized Banach ringed spaces 
$$
f\colon X=((\Spec^{\mathrm{val}}\mathcal{A},\O_{\mathcal{A}}),v_{\mathcal{A}})\longrightarrow Y=((\Spec^{\mathrm{val}}\mathcal{B},\O_{\mathcal{B}}),v_{\mathcal{B}}),
$$
we only need to verify the following: if $f(U)\subseteq V$, where $U=\Spec^{\mathrm{val}}\mathcal{A}'$ (resp.\ $V=\Spec^{\mathrm{val}}\mathcal{B}'$) is an $\R_+$-rational subdomain of $\Spec^{\mathrm{val}}\mathcal{A}$ (resp.\ $\Spec^{\mathrm{val}}\mathcal{B}$), then there exists uniquely a morphism $\varphi'\colon\mathcal{B}'\rightarrow\mathcal{A}'$ of $\R_+$-affinoid rings such that 
$$
\xymatrix{\mathcal{B}\ar[r]^{\varphi}\ar[d]_{\mathrm{res}}&\mathcal{A}\ar[d]^{\mathrm{res}}\\ \mathcal{B}'\ar[r]_{\varphi'}&\mathcal{A}'}
$$
commutes.
Set $\mathcal{B}'=\mathcal{B}(\frac{r^{-1}_1g_1,\ldots,r^{-1}_ng_n}{g_0})$.
Let us show that $\varphi(g_i)/\varphi(g_0)\in(F_{\mathcal{A}'})_{r_i}$ for $i=1,\ldots,n$.
Indeed, for any $x\in V=\Spec^{\mathrm{val}}\mathcal{A}'$, $k(x)$ contains the image of $\varphi(g_i)/\varphi(g_0)$ by the assumption, and, since $\varphi$ respects filtrations, it lies in $(F_x)_{r_i}$.
Hence $\varphi(g_i)/\varphi(g_0)\in(F_{\mathcal{A}'})_{r_i}$ for $i=1,\ldots,n$.
Now, set $\mathcal{A}''=\mathcal{A}'(\frac{r^{-1}_1\varphi(g_1),\ldots,r^{-1}_n\varphi(g_n)}{\varphi(g_0)})$, which allows a morphism $\mathcal{B}'\rightarrow\mathcal{A''}$.
Since $\varphi(g_i)/\varphi(g_0)\in(F_{\mathcal{A}'})_{r_i}$ for $i=1,\ldots,n$, we have $U=\Spec^{\mathrm{val}}\mathcal{A}'=\Spec^{\mathrm{val}}\mathcal{A}''$ due to \ref{prop-vsbr-admsiteRaffinoidrings4}.
Then, by \ref{prop-vsbr-sheafconditionBanachtype1} (2), we have $\mathcal{A}''\cong\mathcal{A}'$, whence the claim.

Conversely, suppose a morphism $f$ as above is given.
By definition, $f$ induces a bounded homomorphism of Banach rings $\varphi\colon\mathcal{B}^{\mathrm{B}}\rightarrow\mathcal{A}^{\mathrm{B}}$.
Let $x=(z,F)\in\Spec^{\mathrm{val}}\mathcal{A}$.
By the condition (b) of \ref{dfn-vsbr-MBRSval2}, the composition
$$
\mathcal{B}^{\mathrm{B}}\longrightarrow\mathcal{A}^{\mathrm{B}}\longrightarrow\mathcal{H}(z)
$$
and $F$ on $\mathcal{H}(z)$ define the point $f(x)$ of $\Spec^{\mathrm{val}}\mathcal{B}$.
Based on this, let us show that $\varphi$ preserves the filtrations $F_{\mathcal{B}}$ and $F_{\mathcal{A}}$, hence giving a morphism of $\R_+$-affinoid rings $\mathcal{B}\rightarrow\mathcal{A}$.
For $h\in(F_{\mathcal{B}})_r$ and $y=f(x)\in\Spec^{\mathrm{val}}\mathcal{B}$, let $\ovl{h}_y$ be the image of $h$ in the residue field $k(y)$.
We know that $\ovl{h}_y\in(F_y)_r$.
Then the image of $\varphi(h)$ in $k(x)$ coincides with the image of $\ovl{h}_y$, which belongs to $(F_x)_r$ by the condition (b) of \ref{dfn-vsbr-MBRSval2}.
Hence $\varphi(h)\in(F_{\mathcal{A}})_r$.

Hence $\varphi$ induces a morphism $g\colon X\rightarrow Y$ of $\R_+$-metrized Banach ringed spaces.
We need to show that $f=g$.
For $x=(z,F)\in X$, consider $f(x)=(f(z),F')$ and $g(x)=(g(z),F'')$.
Observe that the diagrams
$$
\xymatrix{\mathcal{B}\ar[r]^{\varphi}\ar[d]&\mathcal{A}\ar[d]\\ (k(f(z)),F')\ar[r]&(k(z),F),}\qquad\xymatrix{\mathcal{B}\ar[r]^{\varphi}\ar[d]&\mathcal{A}\ar[d]\\ (k(g(z)),F'')\ar[r]&(k(z),F)}
$$
are commutative with dominating lower horizontal arrows (see \ref{dfn-vsbr-MBRSval2} (b)).
Then, both $f(x)$ and $g(x)$ are the filtered valuations on $\mathcal{B}$ induced from the composition
$$
\mathcal{B}\longrightarrow\mathcal{A}\longrightarrow(k(z),F)\longrightarrow(\mathcal{H}(z),F),
$$
from which $f(x)=g(x)$ follows.

Finally, we need to show that the morphism of the structure sheaves is the same as the one induced by $\varphi$.
For two $\R_+$-rational subdomains $U=\Spec^{\mathrm{val}}\mathcal{A}'$ and $V=\Spec^{\mathrm{val}}\mathcal{B}$ of $X$ and $Y$, respectively, such that $f(U)\subseteq V$, the square of bounded homomorphisms of Banach rings
$$
\xymatrix{\mathcal{B}\ar[r]^{\varphi}\ar[d]_{\mathrm{res}}&\mathcal{A}\ar[d]^{\mathrm{res}}\\ \mathcal{B}'\ar[r]_{\psi}&\mathcal{A}'}
$$
is commutative, where $\psi$ is the bounded homomorphism induced by $f$.
Since $\mathcal{B}\rightarrow\mathcal{B}'$ is an epimorphism (\ref{cor-vsbr-powerboundedfiltration1}), $\psi$ is equal to the homomorphism induced from $\varphi$.
Since this is valid for all $\R_+$-rational subdomains $U,V$ with $f(U)\subseteq V$, one sees, by considering filters of $\R_+$-rational subdomains, that $f=g$, as desired.
\end{proof}

\begin{ntn}\label{ntn-vsbr-MBRSval}{\rm 
In the sequel, for an $\R_+$-affinoid ring $\mathcal{A}=(\mathcal{A},F_{\mathcal{A}})$, we denote simply by
$$
\Spec^{\mathrm{val}}\mathcal{A}=((\Spec^{\mathrm{val}}\mathcal{A},\O_{\mathcal{A}}),v_{\mathcal{A}}),
$$
unless any confusion is expected to occur.}
\end{ntn}

\begin{dfn}\label{dfn-vsbr-MBRSval3}{\rm 
An $\R_+$-metrized Banach ringed space isomorphic to $\Spec^{\mathrm{val}}\mathcal{A}$ by an $\R_+$-sheafy $\R_+$-affinoid ring $\mathcal{A}$ is called an {\it $\R_+$-metrized affinoid space}.
Especially, $\Spec^{\mathrm{val}}\mathcal{A}$ is called the {\it $\R_+$-metrized affinoid space attached to $\mathcal{A}$}.}
\end{dfn}

Any $\R_+$-rational subdomain $\Spec^{\mathrm{val}}\mathcal{B}$ of $\Spec^{\mathrm{val}}\mathcal{A}$ is an $\R_+$-metrized affinoid space.
Note that, in this situation, even if $F_{\mathcal{A}}=F^o_{\mathcal{A}}$, $F_{\mathcal{B}}$ and $F^o_{\mathcal{B}}$ may not be equal to each other.

\begin{dfn}\label{dfn-vsbr-MBRSval32}{\rm 
An $\R_+$-metrized Banach ringed space $X$ is an {\it $\R_+$-metrized analytic space} if it is covered by open subspaces that are $\R_+$-metrized affinoid spaces.}
\end{dfn}

By what we have mentioned above, any open subset of an $\R_+$-metrized analytic space is again an $\R_+$-metrized analytic space, called an {\it open subspace}.

We denote by
$$
\mathscr{M}\mathbf{Ansp}^{\R_+}
$$
the full subcategory of the category of $\R_+$-metrized Banach ringed spaces consisting of $\R_+$-metrized analytic spaces.
By what we have seen above, there exists a fully faithful functor
$$
(\R_+\textrm{-}\mathbf{Aff}^{\mathrm{Sh}})^{\opp}\longrightarrow\mathscr{M}\mathbf{Ansp}^{\R_+},\qquad \mathcal{A}\longmapsto\Spec^{\mathrm{val}}\mathcal{A},
$$
with the essential image being $\R_+$-metrized affinoid spaces.

For an $\R_+$-metrized analytic space $S$, one can consider the comma category
$$
\mathscr{M}\mathbf{Ansp}^{\R_+}_S
$$
of $\R_+$-metrized analytic spaces over $S$ with morphisms over $S$.

\begin{dfn}\label{dfn-vsbr-MBRSvalfintype}{\rm 
A morphism $f\colon X\rightarrow Y$ of $\R_+$-metrized analytic space is said to be {\it locally of $\R_+$-finite type} if there exists an open covering $Y=\bigcup_{i\in I}V_i$ and, for each $i\in I$, an open covering $f^{-1}(V_i)=\bigcup_{j\in J_i}U_{ij}$ such that the following conditions are satisfied:
\begin{itemize}
\item[{\rm (a)}] for $i\in I$ and $j\in J_i$, $V_i$ and $U_{ij}$ are $\R_+$-metrized affinoid space; say $V_i\cong\Spec^{\mathrm{val}}\mathcal{B}_i$ and $U_{ij}=\Spec^{\mathrm{val}}\mathcal{A}_{ij}$;
\item[{\rm (b)}] by the homomorphism $\mathcal{B}_i\rightarrow\mathcal{A}_{ij}$ corresponding to $U_{ij}\rightarrow V_i$, $\mathcal{A}_{ij}$ is of $\R_+$-finite type over $\mathcal{B}_i$.
\end{itemize}
A morphism {\it of $\R_+$-finite type} is a locally of finite type morphism such that the underlying continuous map is quasi-compact.}
\end{dfn}

\subsubsection{Relation with adic spaces}\label{subsub-vsbr-relationadicspaces}
\begin{dfn}\label{dfn-vsbr-relationadicspaces1}{\rm 
An $\R_+$-metrized analytic space $X$ is said to be {\it of adic type} if there exists an open covering $X=\bigcup_{\alpha\in L}U_{\alpha}$ by open $\R_+$-affinoid subspaces with the following properties:
\begin{itemize}
\item[{\rm (a)}] for each $\alpha\in L$, $U_{\alpha}\cong\Spec^{\mathrm{val}}\mathcal{A}_{\alpha}$, where $\mathcal{A}_{\alpha}$ is an $\R_+$-affinoid ring of adic type (see \ref{dfn-vsbr-Banachringfadictype} (2));
\item[{\rm (b)}] for $\alpha,\beta\in L$, $U_{\alpha}\cap U_{\beta}$ has an open covering $U_{\alpha}\cap U_{\beta}=\bigcup_{\lambda}V_{\lambda}$ with $V_{\lambda}=\Spec^{\mathrm{val}}\mathcal{B}_{\lambda}$, where $\mathcal{B}_{\lambda}$ is a {\it finite type} $\R_+$-rational localization of both $\mathcal{A}_{\alpha}$ and $\mathcal{A}_{\beta}$, viz., $\mathcal{B}_{\lambda}=\mathcal{A}_{\alpha}(\frac{1\cdot f_1,\ldots,1\cdot f_n}{f_0})$ and similarly over $\mathcal{A}_{\beta}$.
\end{itemize}}
\end{dfn}

If $X$ is an $\R_+$-metrized analytic space of adic type with an open affinoid covering $X=\bigcup_{\alpha\in L}U_{\alpha}$ as above, each affine piece $U_{\alpha}=\Spec^{\mathrm{val}}\mathcal{A}_{\alpha}$ corresponds, by $1$-restriction (\S\ref{subsub-vsbr-adicspectrum}), an adic space $\Spa F_1(\mathcal{A}_{\alpha})$; notice that the affinoid ring $F_1(\mathcal{A}_{\alpha})$ is sheafy due to the $\R_+$-sheafiness of $\mathcal{A}_{\alpha}$.
This construction can be globalized by patching, and we get an adic space, denoted by 
$$
F_1(X)
$$
as the gluing of $\{\Spf F_1(\mathcal{A}_{\alpha})\}_{\alpha\in L}$.

We denote by
$$
\mathscr{M}\mathbf{Ansp}^{\R_+,\mathrm{adic}}
$$
the full subcategory of $\mathscr{M}\mathbf{Ansp}^{\R_+}$ consisting of $\R_+$-metrized analytic space of adic type. 
The above construction gives a functor
$$
(\cdot)_1\colon\mathscr{M}\mathbf{Ansp}^{\R_+,\mathrm{adic}}\longrightarrow\mathbf{Adsp},\qquad X\longmapsto F_1(X),
$$
called the {\it $1$-restriction functor} to the category of adic spaces.
Note that there exists a continuous surjective map
$$
X\longrightarrow F_1(X)
$$
by \ref{prop-vsbr-adicspectrum}.

As indicated in \S\ref{subsub-vsbr-adicspectrum}, from any affinoid ring $(\mathcal{A}^{\pm},\mathcal{A}^+)$ in the sense of Huber equipped with an ideal of definition $I$ (in the sense therein), one can construct in a standard way an $\R_+$-affinoid ring of adic type that gives back the original affinoid ring by $1$-restriction.
This can also be globalized in the following way.
Consider an adic space $X$ {\it with an ideal of definition} $\mathscr{I}$; here, by an ideal of definition of $X$, we mean an ideal sheaf of $\O^+_X$ such that there exists an open affinoid covering $X=\bigcup_{\alpha\in L}U_{\alpha}$ with the property that each $\mathscr{I}|_{U_{\alpha}}$ comes from an ideal of definition of $\mathcal{A}_{\alpha}$ (where $U_{\alpha}=\Spf\mathcal{A}_{\alpha}$).
For example, if the adic space comes from a rigid space by the functor $\ZRT$ (\S\ref{sub-rigidgeomadicspaces}), it is an ideal of definition in the sense of \ref{dfn-ZRstrsheaf2}.
We also have to fix a real number $0<c<1$.
Then one can construct an $\R_+$-metrized analytic space $\til{X}$ such that $F_1(\til{X})=X$ by globalizing the local construction as above, and thus we have a functor 
$$
\bigg\{\begin{minipage}{17em}{\small adic spaces with ideal of definition and morphisms respecting ideals of definition}\end{minipage}\bigg\}\longrightarrow \mathscr{M}\mathbf{Ansp}^{\R_+,\mathrm{adic}}.
$$
In particular, adic spaces that admit ideal of definition lie in the essential image of the functor $(\cdot)_1$.

In this way, $\R_+$-metrized analytic spaces are regarded as `adic spaces (rigid spaces) with extra higher structure,' where the extra higher structure is an analogue of metrics.

An important example of the last-mentioned functor is the following.
Let $K=(K,\nu=|\cdot|)$ be a non-archimedean Banach field with non-trivial valuation $|\cdot|\colon K\rightarrow\R_{\geq 0}$, and $a\in K$ with $0<|a|<1$, and set $c=|a|$ (the following construction actually does not depend on the choice of $a$, since the norm $|\cdot|$ on $K$ is fixed).
Let $\mathscr{X}=\Spa\mathcal{A}$ be an affinoid adic space of finite type (\ref{dfn-adicspacefinitetype}) over $K$.
Then $\mathcal{A}^{\pm}$ is a classical affinoid algebra over $K$ (or equivalently, strictly $K$-affinoid algebra), and $\mathcal{A}^+=\mathcal{A}^o$ (\ref{prop-canonicalsubring}).
Moreover, by \ref{thm-vsbr-berkovichaffinoidalgebras}, $\mathcal{A}$ can be viewed uniquely as an $\R_+$-affinoid algebra with the filtration $F_{\mathcal{A}}=F^o_{\mathcal{A}}=F^{\Sp}_{\mathcal{A}}$.
Hence we have the functor (called the {\it metrization functor})
$$
\mathscr{X}=\Spa\mathcal{A}\longmapsto\mathscr{X}^{\mathrm{met}}=\Spec^{\mathrm{val}}(\mathcal{A},F^o_{\mathcal{A}})
$$
from the category of affinoid adic spaces of finite type over $K$ to the category of $\R_+$-metrized affinoid spaces of finite type over $K$.
Note that the sheaf condition on the adic space side (equivalent to Tate's acyclicity) implies the sheaf condition of Banach type, since continuity implies boundedness for $K$-linear maps (cf.\ \cite{BGR}, 2.1.8).
By \ref{thm-vsbr-MBRSval} (2), we have by gluing a canonical fully faithful functor
$$
(\cdot)^{\mathrm{met}}\colon\mathbf{Adsp}^{\mathrm{lft}}_K\longrightarrow\mathscr{M}\mathbf{Ansp}^{\R_+}_K
$$
from the category of locally of finite type adic spaces over $K$ to the category of $\R_+$-metrized analytic spaces over $K$.

\begin{cor}\label{cor-Banalytichuberadicspaces11}
The functor $\mathscr{X}\mapsto\mathscr{X}^{\mathrm{met}}$ thus obtained is a fully faithful functor from the category of locally of finite type adic spaces over $K$ to the category of $\R_+$-metrized analytic spaces $\R_+$-finite type over $K$.
The essential image consists of $\R_+$-finite type $\R_+$-metrized analytic spaces $X$ over $K$ having an open $\R_+$-affinoid covering 
$$
X=\bigcup_{\alpha\in L}U_{\alpha},\quad U_{\alpha}=\Spec^{\mathrm{val}}\mathcal{A}_{\alpha}
$$
such that 
\begin{itemize}
\item[{\rm (a)}] for each $\alpha\in L$, $U_{\alpha}$ is an $\R_+$-metrized affinoid space of finite type over $K$, viz., $\mathcal{A}_{\alpha}$ is of finite type over $K$ $($see {\rm \ref{dfn-sbr-Rfinitetype})}, and
\item[{\rm (b)}] for each $\alpha,\beta\in L$, $U_{\alpha}\cap U_{\beta}$ is also covered by $\R_+$-affinoid subdomains of finite type over $K$.
\end{itemize}
Moreover, we have a canonical continuous surjection $\mathscr{X}^{\mathrm{met}}\rightarrow\mathscr{X}$ $($by the $1$-restriction map {\rm (\S\ref{subsub-vsbr-adicspectrum})}. 
\end{cor}

\subsection{Berkovich's analytic geometry}\label{sub-vsbr-Berkovichanalyticspaces}
\subsubsection{Gerritzen-Grauert theorem}\label{subsub-vsbr-GG}
Let $A$ be a Berkovich's $K$-affinoid algebra\index{algebra!affinoid algebra@affinoid ---!Berkovichs affinoid algebra@Berkovich's --- ---} (\ref{dfn-berkovichaffinoicalgebras}), and set $X=\mathscr{M}(A)$.
\begin{dfn}[{\rm \cite{Berk1}, 2.2.1}]\label{dfn-Berkovichaffinoidsubdomains}{\rm 
A closed subset $U\subseteq\mathscr{M}(A)$ is said to be an {\it affinoid subdomain} if these exists a bounded $K$-algebra homomorphism $A\rightarrow A_U$ of Berkovich's $K$-affinoid algebras such that the following universal mapping property holds: for any Berkovich's $L$-affinoid algebra $B$, where $L=(L,|\cdot|)$ is a complete isometric extension of $K$, and any bounded $K$-algebra homomorphism $A\rightarrow B$ such that the image of the induced map $\mathscr{M}(B)\rightarrow\mathscr{M}(A)$ lies in $U$, there exists a unique bounded $K$-algebra homomorphism $A_U\rightarrow B$ such that the diagram
$$
\xymatrix@C-3ex@R-1.5ex{A\ar[rr]\ar[dr]&&A_U\ar@{-->}[dl]\\ &B}
$$
commutes.}
\end{dfn}

One sees, as in the case of classical rigid geometry, that the Berkovich's $K$-affinoid algebra $A_U$ is uniquely determined by $U$ up to unique isomorphisms, and that the map $\mathscr{M}(A_U)\rightarrow\mathscr{M}(A)$ gives a homeomorphism onto $U$ (\cite{Berk1}, 2.2.4).

In the rest of this paragraph, we give an argument for reducing Gerritzen-Grauert type theorem for Berkovich's $K$-affinoid algebras to that for Tate's affinoid algebras, based on generalities of locally coherent spaces. 
\begin{lem}\label{lem-vsbr-GG1}
Let $X$ be a locally coherent space, and $W\subseteq X$ a subset of $X$ satisfying the following conditions:
\begin{itemize}
\item[{\rm (a)}] $X\setminus W$ with the subspace topology is sober and retrocompact {\rm ({\bf \ref{ch-pre}}.\ref{dfn-retrocompact})} in $X$;
\item[{\rm (b)}] $W$ is stable under generizations.
\end{itemize}
Then $W$ is open in $X$.
\end{lem}

\begin{proof}
We may assume $X$ is coherent.
Set $C=X\setminus W$.
For any quasi-compact open subset $V\subseteq X$, $C\cap V$ is quasi-compact by (a).
Since open sets of $C$ of this type form a basis of topology, which is stable under finite intersection, $C$ is coherent and sober.
In particular, $C\cap V$ is coherent, for any quasi-compact open subset $V\subseteq X$.

Let $x\in W$.
Since $C\cap G_x=\emptyset$ by (b) (where $G_x$ denotes, as before, the set of all generizations of $x$ in $X$), there exists a quasi-compact open neighborhood $V_x$ of $x$ such that $C\cap V_x=\emptyset$; indeed, otherwise $C\cap G_x=\varprojlim_VC\cap V$, the inverse limit of non-empty coherent subsets with quasi-compact transition maps, has to be non-empty due to {\bf \ref{ch-pre}}.\ref{thm-projlimcohsch1}.
Thus we know, for any $x\in W$, the existence of open neighoborhood $V_x$ contained in $W$, which shows the assertion.
\end{proof}

\begin{cor}\label{cor-vsbr-GG1}
Let $f\colon X\rightarrow Y$ be a quasi-compact surjective map between locally coherent sober spaces.
Assume $W\subseteq Y$ is stable under generizations, and $f^{-1}(W)$ is open in $X$.
Then $W$ is open in $Y$.
\end{cor}

\begin{proof}
We may assume $Y$ is coherent.
Set $C=Y\setminus W$.
We are going to check the conditions (a) and (b) in \ref{lem-vsbr-GG1}.
It is easy to see that $C\hookrightarrow Y$ is quasi-compact, since $f^{-1}(C)=X\setminus f^{-1}(W)$ is coherent and $f$ is surjective.
The subspace $C$ is a $\mathrm{T}_0$-space, as $Y$ is $\mathrm{T}_0$.
Let $Z\subseteq C$ be an irreducible closed subset.
Consider the set $E$ of all quasi-compact open subsets of $C$, and set $E'=\{U\in E: U\cap Z\neq\emptyset\}$.
Then $E'$ is a filter, and $\bigcap_{U\in E'}(U\cap Z)\neq\emptyset$ holds; indeed, $f^{-1}(\bigcap_{U\in E'}(U\cap Z))=\bigcap_{U\in E'}f^{-1}(U\cap Z)\neq\emptyset$ ({\bf \ref{ch-pre}}.\ref{thm-projlimcohsch1}), since $f$ is surjective and quasi-compact.
It is then easy to see that $Z=\ovl{\{z\}}$ for $z\in\bigcap_{U\in E'}(U\cap Z)$.
Hence $C$ is sober, whence the claim.
\end{proof}

\begin{thm}\label{thm-vsbr-GG}
Let $K$ be a non-archimedean Banach field, $\mathcal{A}$ a Berkovich's $K$-affinoid algebra, and $\mathcal{B}$ another Berkovich's $K$-affinoid algebra that gives an affinoid subdomain in the sense of {\rm \ref{dfn-Berkovichaffinoidsubdomains}}.
Then the induced mapping
$$
\Spec^{\mathrm{val}}\mathcal{B}\longrightarrow\Spec^{\mathrm{val}}\mathcal{A}
$$
is an open immersion of topological spaces, viz., a homeomorphism onto an open subset of $\Spec^{\mathrm{val}}\mathcal{A}$.
\end{thm}

\begin{proof}
Take an extension $L/K$ of Banach fields such that $|L^{\times}|=\R_+$.
We first check:
\begin{itemize}
\item[{\rm (a)}] $f\colon\Spec^{\mathrm{val}}\mathcal{B}\rightarrow\Spec^{\mathrm{val}}\mathcal{A}$ is injective;
\item[{\rm (b)}] $f_L\colon\Spec^{\mathrm{val}}\mathcal{B}_L\rightarrow\Spec^{\mathrm{val}}\mathcal{A}_L$ is an open immersion,
\item[{\rm (c)}] $f_L(\Spec^{\mathrm{val}}\mathcal{B}_L)$ is the inverse image of $f(\Spec^{\mathrm{val}}\mathcal{B})$ under the canonical map $\Spec^{\mathrm{val}}\mathcal{A}_L\rightarrow\Spec^{\mathrm{val}}\mathcal{A}$.
\end{itemize}

For (a), first observe that $\mathcal{H}(f(z))\cong\mathcal{H}(z)$ for $z\in\mathscr{M}(\mathcal{B})$.
Then filtered valuation field structures on $\mathcal{H}(z)$ are the same as those of $\mathcal{H}(f(z))$, whence the injectivity in question.

For (b), first note that $\mathcal{B}_L$ gives an affinoid subdomain of $\mathscr{M}(\mathcal{A}_L)$ in the sense of Berkovich, and both $\mathcal{A}_L$ and $\mathcal{B}_L$ are Tate's affinoid algebras over $L$.
Thus, (b) is the consequence of the Gerritzen-Grauert theorem \ref{thm-GGref} for Tate's affinoid spaces, which we assume in this book, due to \ref{prop-vsbr-berkovichaffinoidalgebras1} and \ref{thm-rigidgeomadicspacesff1}.

To check (c), take $x=(z,F)\in\Spec^{\mathrm{val}}\mathcal{A}_L$, which maps to a point in $f(\Spec^{\mathrm{val}}\mathcal{B})$.
By the argument in checking (a) above, the filtered homomorphism $(\mathcal{A},F_{\mathcal{A}})\rightarrow(\mathcal{H}(z),F)$ factors through $\mathcal{B}$, and $x$ lifts to a unique point of $\Spec^{\mathrm{val}}\mathcal{B}_L$.
Then (c) follows by \ref{thm-vsbr-berkovichaffinoidalgebras}.

To conclude, we need to check that the subset $W=f(\Spec^{\mathrm{val}}\mathcal{B})$ is stable under generizations.
This also follows from the argument in checking (a) above.
Hence $W=f(\Spec^{\mathrm{val}}\mathcal{B})$ is an open subset of $\Spec^{\mathrm{val}}\mathcal{A}$ by \ref{cor-vsbr-berkovichaffinoidalgebras1} and \ref{cor-vsbr-GG1}.
Replacing $\mathcal{B}$ by Berkovich's $K$-affinoid algebras giving affinoid subdomains of $\mathscr{M}(\mathcal{B})$, one deduces that the map $f$ is an open map.
\end{proof}

\subsubsection{Berkovich analytic spaces}\label{subsub-vsbr-Berkovichanalyticspaces}
\index{Berkovich analytic space|(}
\begin{dfn}\label{dfn-berkovichanalyticspacesdefinition}{\rm 
(1) Let $X$ be a locally Hausdorff topological space, and $\tau$ a net\index{net} ({\bf \ref{ch-pre}}.\ref{dfn-netsfortopologicalspaces} (2)) on $X$ consisting of compact subsets.
A {\it $K$-affinoid atlas $\mathscr{A}$ on $X$ with the net $\tau$} is given by the following data:
\begin{itemize}
\item[{\rm (a)}] a Berkovich's $K$-affinoid algebra $A_U$ and a homeomorphism $U\stackrel{\sim}{\rightarrow}\mathscr{M}(A_U)$ for each $U\in\tau$;
\item[{\rm (b)}] for each pair $U,U'\in\tau$ with $U\subseteq U'$, a bounded $K$-algebra homomorphism $\rho^V_U\colon A_{U'}\rightarrow A_U$ that identifies $U$ with an affinoid subdomain in $U'\cong\mathscr{M}(A_{U'})$ with the corresponding $K$-affinoid algebra $A_U$;
\end{itemize}
these data are assumed to satisfy the following cocycle condition: 
\begin{itemize}
\item[{\rm (c)}] for any triple $U,U',U''\in\tau$ with $U\subseteq U'\subseteq U''$, we have $\rho^{U''}_U=\rho^{U'}_U\circ\rho^{U''}_{U'}$.
\end{itemize}

(2) A {\it $($Berkovich's$)$ $K$-analytic space} is a triple $X=(X,\tau,\mathscr{A})$ consisting of a locally Hausdorff space $X$, a net $\tau$ on $X$, and a $K$-affinoid atlas $\mathscr{A}$ on $X$ with the net $\tau$.
For example, for any Berkovich's $K$-affinoid algebra $A$, the set $X=\mathscr{M}(A)$ with the $K$-affinoid atlas given by affinoid subdomains gives a Berkovich's $K$-analytic space, called {\it $($Berkovich's$)$ $K$-affinoid space}.

(3) A {\it strong morphism} of Berkovich's $K$-analytic spaces
$$
X=(X,\tau,\mathscr{A})\rightarrow X'=(X',\tau',\mathscr{A}')
$$
consists of:
\begin{itemize}
\item[{\rm (a)}] a continuous map $\varphi\colon X\rightarrow X'$ such that for each $U\in\tau$ there exists $U'\in\tau'$ such that $\varphi(U)\subseteq U'$ and
\item[{\rm (b)}] for each pair $(U,U')$ with $U\in\tau$, $U'\in\tau'$, and $\varphi(U)\subseteq U'$, a bounded $K$-algebra homomorphism $\phi_{U/U'}\colon A'_{U'}\rightarrow A_U$ such that the induced map $\mathscr{M}(A_U)\rightarrow\mathscr{M}(A'_{U'})$ is identified with $\varphi|_U\colon U\rightarrow U'$ via $\mathscr{M}(A_U)\cong U$ and $\mathscr{M}(A'_{U'})\cong U'$.
\end{itemize}

(4) A strong morphism $\varphi\colon X=(X,\tau,\mathscr{A})\rightarrow X'=(X',\tau',\mathscr{A}')$ is said to be a {\it quasi-isomorphism} if $\varphi\colon X\rightarrow X'$ is a homeomorphism and, for any pair $(U,U')$ as in (3) (b), the map $\mathscr{M}(A_U)\rightarrow\mathscr{M}(A'_{U'})$ identifies $U\cong\mathscr{M}(A_U)$ with an affinoid subdomain in $U'$.

(5) The {\it category of $($Berkovich's$)$ $K$-analytic spaces}, denoted by $\mathbf{Bsp}_K$, is the quotient category of the category of $K$-analytic spaces and strong morphisms by quasi-isomorphisms. 

(6) The {\it category of $($Berkovich's$)$ strictly $K$-analytic spaces}, denoted by $\mathbf{Bsp}^{\ast}_K$, is defined in the similar way, where all $K$-affinoid algebras in the $K$-affinoid atlas are strictly $K$-affinoid.}
\end{dfn}

Let $X=(X,\tau,\mathcal{A})$ be a $K$-analytic space. 
A subset $W\subseteq X$ is said to be {\it $\tau$-special} if it is compact and there exists a finite covering $W=\bigcup^n_{i=1}W_i$ such that $W_i\cap W_j$ ($i,j=1,\ldots,n$) belongs to $\tau$ and that $A_{W_i}\widehat{\otimes}_AA_{W_j}\rightarrow A_{W_i\cap W_j}$ is an admissible epimorphism.
In this case, one has a commutative Banach $K$-algebra $A_W$ defined by
$$
A_W=\ker\Big({\displaystyle \prod^n_{i=1}A_{W_i}\rightarrow\prod_{i,j}A_{W_i\cap W_j}}\Big).
$$
By generalized Tate acyclicity theorem (\cite{Berk1}, 2.2.5), the Banach $K$-algebra $A_W$ does not depend, up to isomorphism, on the choice of the covering $\{W_i\}^n_{i=1}$.
We define
\begin{itemize}
\item the collection $\ovl{\tau}$ to be the set of all subsets $W\subseteq X$ satisfying the following condition: there exists $U\in\tau$ with $W\subseteq U$ such that $W$ is identified with an affinoid subdomain of $\mathscr{M}(A_U)$;
\item the collection $\widehat{\tau}$ to be the set of all $\ovl{\tau}$-special subsets $W$ such that the algebra $A_W$ is $K$-affinoid and, for any finite covering $W=\bigcup^n_{i=1}W_i$ as above, each $W_i$ is an affinoid subdomain of $\mathscr{M}(A_W)$ with the corresponding $K$-affinoid algebra $A_{W_i}$.
\end{itemize}
Then these collections are nets on $X$, and the $K$-affinoid atlas $\mathscr{A}$ extends uniquely to $K$-affinoid atlases $\ovl{\mathscr{A}}$ and $\widehat{\mathscr{A}}$ with respect to the nets $\ovl{\tau}$ and $\widehat{\tau}$, respectively (\cite{Berk2}, 1.2.6, 1.2.13).
\begin{itemize}
\item Elements of $\widehat{\tau}$ are called {\it $K$-affinoid domains} in $X$.
\item The $\widehat{\tau}$-special subsets are called {\it $K$-special domains} in $X$.
\end{itemize}

Similarly, one defines the notions of {\it strictly $K$-affinoid domains} and {\it strictly $K$-special domains}.
\index{Berkovich analytic space|)}

\subsubsection{G-topology on Berkovich analytic spaces}\label{subsub-Gtopologyberkovichanalyticspace}
\begin{dfn}\label{dfn-berkovichspacesanalyticdomains}{\rm 
A subset $Y\subseteq X$ of a $K$-analytic space $X=(X,\tau,\mathscr{A})$ is called a {\it $K$-analytic domain} if $\widehat{\tau}|_Y=\{V\in\widehat{\tau}:V\subseteq Y\}$ is a net on $Y$.}
\end{dfn}

In this case, $Y$ endowed with the subspace topology carries the induced structure of a $K$-analytic space $Y=(Y,\widehat{\tau}|_Y,\mathscr{A}|_Y)$, where $\mathscr{A}|_Y$ is the induced atlas with the net $\widehat{\tau}|_Y$, together with the canonical morphism $Y\rightarrow X$ of $K$-analytic spaces.
Moreover,
\begin{itemize}
\item any open subset $U\subseteq X$ is a $K$-analytic domain (but not conversely in general);
\item the intersection of finitely many $K$-analytic domains is again a $K$-analytic domain.
\end{itemize}

We now define a Grothendieck topology\index{topology!Grothendieck topology@Grothendieck ---}\index{Grothendieck!Grothendieck topology@--- topology} on $X$ generated by the pretopology defined as follows:
\begin{itemize}
\item admissible open subsets are $K$-analytic domains;
\item an admissible covering of a $K$-analytic domain $Y\subseteq X$ is a family $\{Y_{\alpha}\}_{\alpha\in L}$ of analytic domains in $Y$ that is a quasi-net on $Y$.
\end{itemize}
We call this Grothendieck topology\index{topology!Grothendieck topology@Grothendieck ---}\index{Grothendieck!Grothendieck topology@--- topology} the G-topology on $X$; note that its admissible opens are not necessarily open subsets of the topological space $X$.
One has then a natural sheaf $\O_{X_{\mathrm{G}}}$, called the structure sheaf, on the resulting site $X_{\mathrm{G}}$ in such a way that, for any $K$-affinoid domain $W\subseteq X$, we have
$$
\O_{X_{\mathrm{G}}}(W)=A_W.
$$

\begin{dfn}[{\rm \cite[p.\ 22]{Berk2}}]\label{dfn-Berkovichanalyticgoodspaces}{\rm 
A $K$-analytic space $X=(X,\tau,\mathscr{A})$ is said to be {\it good} if every point $x\in X$ has a neighborhood $W$ from $\widehat{\tau}$.}
\end{dfn}

By definition, $K$-affinoid spaces are good. 
It is known that any $K$-analytic space coming from a $K$-scheme of locally of finite type through the GAGA functor (cf.\ \cite{Berk1}, \S3.4) is good.
But there are many interesting non-good $K$-analytic spaces.

If $X$ is good, then one can restrict the G-topology as above to open subsets and open coverings to recover the usual topology of $X$. The resulting sheaf, obtained from $\O_{X_{\mathrm{G}}}$ by the restriction, is denoted by $\O_X$.
It is known that the category of coherent $\O_{X_{\mathrm{G}}}$-modules and that of coherent $\O_X$-modules are equivalent to each other (\cite{Berk2}, 1.3.4).

\subsubsection{Berkovich analytic spaces and $\R_+$-metrized analytic spaces}\label{subsub-vsbr-BerkovichRmetrizedanalyticspace}
Let $\mathcal{A}$ be a Berkovich's $K$-affinoid algebra.
We regard $\mathcal{A}$ as an $\R_+$-affinoid ring of $\R_+$-finite type over $K$ in the unique manner as indicated in \S\ref{subsub-vsbr-berkovichaffinoidalgebras}.

Let $X=(X,\tau,\mathscr{A})$ be a Berkovich's $K$-analytic space\index{Berkovich analytic space} (\ref{dfn-berkovichanalyticspacesdefinition} (2)).
For any $U\in\tau$, one has the $K$-affinoid algebra $\mathcal{A}_U$, and the corresponding $\R_+$-metrized affinoid space $\Spec^{\mathrm{val}}\mathcal{A}_U$.
Let $D_U$ be the distributive lattice of quasi-compact open subsets of $\Spec^{\mathrm{val}}\mathcal{A}_U$.
Then, $D_U$ gives a valuation\index{valuation!valuation of a compact space@--- of a compact space} ({\bf \ref{ch-pre}}.\ref{dfn-valuationscompactHausdorff}) of $U\cong\mathscr{M}(\mathcal{A}_U)$.
Then $\Spec^{\mathrm{val}}\mathcal{A}_U$ is homeomorphic to $\Spec D_U$, since the underlying topological space of $\Spec^{\mathrm{val}}\mathcal{A}_U$ is a reflexive valuative space (\ref{cor-vsbr-berkovichaffinoidalgebras2}).
These data constitute a pre-valuation $v_X=(\tau,\{D_U\}_{U\in\tau})$ ({\bf \ref{ch-pre}}.\ref{dfn-valuationlocallyHausdorff} (1)) of the underlying topological space $X$.
By {\bf \ref{ch-pre}}.\ref{thm-associatedvaluationtopologicalspaces2}, one has the reflexive valuative space\index{valuative!valuative topological space@--- (topological) space!reflexive valuative topological space@reflexive --- ---} $\Spec v_X$ by 
$$
\Spec v_X=\varinjlim_{U\in\tau}\Spec D_U
$$
so that $[\Spec v_X]\cong X$.
By {\bf \ref{ch-pre}}.\ref{prop-RLSCVsp1}, $\Spec v_X$ is locally strongly compact\index{valuative!valuative topological space@--- (topological) space!locally strongly compact valuative topological space@locally strongly compact --- ---} ({\bf \ref{ch-pre}}.\ref{dfn-locallycompactspace}).
By \ref{thm-vsbr-GG} and \ref{prop-vsbr-BerkovichRmetrizedanalyticspace1}, the structures of $\R_+$-metrized analytic spaces on $\Spec^{\mathrm{val}}\mathcal{A}_U$ glue by open immersions, and hence we obtain an $\R_+$-metrized analytic space locally of $\R_+$-finite type over $K$ supported on $\Spec v_X$, which we denote by
$$
X^{\mathrm{met}}.
$$

It is obvious that any strong morphism (\ref{dfn-berkovichanalyticspacesdefinition} (3)) $\varphi\colon X=(X,\tau,\mathcal{A})\rightarrow X'=(X',\tau',\mathcal{A}')$ induces a morphism
$$
\varphi^{\mathrm{met}}\colon X^{\mathrm{met}}\longrightarrow X^{\prime\mathrm{met}}
$$
of $\R_+$-metrized analytic spaces over $K$.\footnote{Notice that the underlying continuous map $\Spec v_X\rightarrow\Spec v_{X'}$ of $\varphi^{\mathrm{met}}$ is the one induced, by the functor considered in {\bf \ref{ch-pre}}.\ref{thm-RLSCVsp1}, from the morphism $(X,v^{\mathrm{val}}_X)\rightarrow(X',v^{\mathrm{val}}_{X'})$ of valued locally Hausdorff spaces ({\bf \ref{ch-pre}}.\ref{dfn-valuationlocallyHausdorffmorph}), which can be constructed from the strong morphism $\varphi$.}
Moreover, if $\varphi$ is a quasi-isomorphism (\ref{dfn-berkovichanalyticspacesdefinition} (4)), then $\varphi^{\mathrm{met}}$ is an isomorphism, since the underlying topological spaces are reflexive.
Thus we obtain the functor 
$$
(\cdot)^{\mathrm{met}}\colon X\longmapsto X^{\mathrm{met}}
$$
from the category $\Bsp_K$ of Berkovich's $K$-analytic spaces to the category $\mathscr{M}\mathbf{Ansp}^{\R_+\textrm{-}\mathrm{lft}}_K$ of $\R_+$-metrized analytic spaces over $K$.

\begin{thm}\label{thm-comparisonberkovichanalyticspacesbis1}
{\rm (1)} The functor $(\cdot)^{\mathrm{met}}$ gives a fully faithful functor from the category $\Bsp_K$ to the category $\mathscr{M}\mathbf{Ansp}^{\R_+\textrm{-}\mathrm{lft}}_K$.
Moreover, the essential image of $(\cdot)^{\mathrm{met}}$ consists of $\R_+$-metrized analytic spaces locally of $\R_+$-finite type with locally strongly compact underlying topological spaces.

{\rm (2)} For any Berkovich's $K$-analytic space $X$, the admissible topos $X^{\sim}_G$ is equivalent to the admissible topos of $X^{\mathrm{met}}$ $($cf.\ {\rm \S\ref{subsub-vsbr-admsiteRaffinoidrings})} as a ringed topos. In particular, it is spacial, and locally ringed.

{\rm (3)} For any Berkovich's $K$-analytic space $X$, the underlying topological space of $X^{\mathrm{met}}$ is a reflexive valuative space, and $[X^{\mathrm{met}}]$ is canonically homeomorphic to $X$.
Moreover, the underlying topological space of $X^{\mathrm{met}}$ is quasi-separated $($resp.\ paracompact, resp.\ coherent$)$ if and only if $X$ is Hausdorff $($resp.\ paracompact, resp.\ compact$)$.
\end{thm}

\begin{proof}
(1) First consider the functor $(\cdot)^{\mathrm{met}}$ restricted on Berkovich's $K$-affinoid spaces $\mathscr{M}(\mathcal{A})$.
In this case, the functor is given by $\mathscr{M}(\mathcal{A})\mapsto\Spec^{\mathrm{val}}\mathcal{A}$, and the fully-faithfulness in this case follows immediately from \ref{thm-vsbr-MBRSval} (2). 

In general, we need to check the following conditions.
\begin{itemize}
\item[(a)] If $\varphi,\psi\colon X=(X,\tau,\mathscr{A})\rightarrow X'=(X',\tau',\mathscr{A}')$ are two strong morphisms of Berkovich $K$-analytic spaces such that the resulting morphisms $f=\varphi^{\mathrm{met}}$ and $g=\psi^{\mathrm{met}}$ coincide, then $\varphi=\psi$.
\item[(b)] For any valuative and locally quasi-compact morphism $f\colon X^{\mathrm{met}}\rightarrow X^{\prime\mathrm{met}}$ of $\R_+$-metrized analytic spaces locally of $\R_+$-finite type over $K$, there exists a diagram 
$$
X\longleftarrow X''\stackrel{\varphi}{\longrightarrow}X'
$$
of strong morphisms of Berkovich's $K$-analytic spaces such that the first morphism is a quasi-isomorphism supported on the identity map $\id_X$ and that $\varphi^{\mathrm{met}}=f$.
\end{itemize}

(a) Since $[X^{\mathrm{met}}]\cong X$ etc.\ as a topological space, $\varphi$ and $\psi$ coincide as continuous mappings of the underlying topological spaces.
For any $U\in\tau$ and $U'\in\tau'$ such that $\varphi(U)\subseteq U'$ (and hence $\psi(U)\subseteq U'$), the maps $f$ and $g$ give the same morphism $\Spec^{\mathrm{val}}\mathcal{A}_U\rightarrow\Spec^{\mathrm{val}}\mathcal{A}_{U'}$ as a morphism of $\R_+$-metrized analytic spaces.
Hence, by what we have already shown, $\varphi=\psi$ as a strong morphism.

(b) Since $f$ is valuative, it induces a continuous mapping $[f]\colon X\rightarrow X'$.
Since $f$ is locally quasi-compact, one can find a refinement $\tau''$ of $\tau$ and the corresponding $K$-affinoid atlas $\mathscr{A}''$ with the property that, for any $U\in\tau''$ there exists $U'\in\tau'$ such that $f(U)\subseteq U'$.
We then have $\Spec^{\mathrm{val}}\mathcal{A}_U\rightarrow\Spec^{\mathrm{val}}\mathcal{A}_{U'}$, which induces $\mathcal{A}'_{U'}\rightarrow\mathcal{A}_U$ (\ref{thm-vsbr-MBRSval} (2)).
Hence we have a strong morphism $\varphi\colon X''=(X,\tau'',\mathscr{A}'')\rightarrow X'$.
Since $[\varphi^{\mathrm{met}}]=[f]$ as a continuous mapping $X\rightarrow X'$, by {\bf \ref{ch-pre}}.\ref{cor-proplemlocallyquasicompactmapsgen}, we have $\varphi^{\mathrm{met}}=f$ as a continuous mapping $X^{\mathrm{met}}\rightarrow X^{\prime\mathrm{met}}$.
It is then obvious, by the affinoid case discussed above, that $\varphi^{\mathrm{met}}=f$ as a morphism of $\R_+$-metrized analytic spaces.
Since
$$
X''=(X,\tau'',\mathscr{A}'')\rightarrow X=(X,\tau,\mathscr{A})
$$
is clearly a quasi-isomorphism supported on the identity map $\id_X$, we have verified the claim.

Finally, by {\bf \ref{ch-pre}}.\ref{prop-RLSCVsp1} the $\R_+$-metrized analytic space $X^{\mathrm{met}}$ is locally strongly compact.
Conversely, if an $\R_+$-metrized analytic space $W$ locally of $\R_+$-finite type over $K$ is locally strongly compact, by {\bf \ref{ch-pre}}.\ref{thm-RLSCVsp1}, $W$ is homeomorphic to $\Spec v$ by a valued locally Hausdorff space $(X=[W],v=(\tau(v),\{v_S\}))$, where $\tau(v)=\{[U]:U\in\QCOuv(W)\}$ and $v_{[U]}=\{[V]\subseteq[U]:V\in\QCOuv(U)\}$.
Hence, for any $\R_+$-affinoid subdomain $U\subseteq W$, one can set $\mathcal{A}_{[U]}$ to be the corresponding $K$-affinoid algebra, and thus we get a Berkovich's $K$-analytic space $(X,\tau,\mathcal{A})$ such that $X^{\mathrm{met}}\cong W$.

(2) By \ref{thm-vsbr-GG}, the site $X_G$ is generated by rational subdomains, and hence the associated topos is equivalent to the admissible topos of $X^{\mathrm{met}}$.

(3) That $X^{\mathrm{met}}$ is reflexive and that $[X^{\mathrm{met}}]\cong X$ are the consequences from the construction of $X^{\mathrm{met}}$.
The other assertions follow from {\bf \ref{ch-pre}}.\ref{thm-RLSCVsp1} and {\bf \ref{ch-pre}}.\ref{prop-paracompactseparatedquotient}.
\end{proof}

Thus, in case the valuation on $K$ is non-trivial, we have, combined with the functor $\mathscr{X}\mapsto\mathscr{X}^{\mathrm{met}}$ considered in \ref{cor-Banalytichuberadicspaces11}, the 2-cartesian diagram of categories
$$
\begin{xy}
(0,0)="BS"*{\Bsp_K},+<10em,0ex>="BA"*{\mathscr{M}\mathbf{Ansp}^{\R_+}_K},+<0em,-10ex>="A"*{\Acsp^{\mathrm{lft}}_K},+<-10em,0ex>="LSA"*{\Acsp^{\mathrm{lsc,lft}}_K}
\ar@{^{(}->}"BS"+<1.8em,-.2ex>;"BA"+<-2.5em,-.2ex>^{(\cdot)^{\mathrm{met}}}
\ar@{^{(}->}"LSA"+<2.7em,-.2ex>;"A"+<-2em,-.2ex>
\ar@{^{(}->}"LSA"+<0em,2.4ex>;"BS"+<0em,-2ex>
\ar@{^{(}->}"A"+<0em,2.4ex>;"BA"+<0em,-2ex>_{(\cdot)^{\mathrm{met}}}
\end{xy}
$$
consisting of fully faithful functors, where $\Acsp^{\mathrm{lft}}_K$ the category of locally of finite type adic spaces over $K$, and $\Acsp^{\mathrm{lsc,lft}}_K$ the full subcategory of $\Acsp^{\mathrm{lft}}_K$ consisting of adic spaces with locally strongly compact underlying topological spaces.
The following corollary is easy to see.
\begin{cor}\label{cor-comparisonberkovichanalyticspacesbis1}
The category $\Acsp^{\mathrm{lsc,lft}}_K$ of locally strongly compact and locally of finite type adic spaces over $K$ is naturally categorically equivalent to the category $\Bsp^{\ast}_K$ of Berkovich's strictly $K$-analytic spaces. 
\end{cor}

Note that the corollary is consistent with \cite{Hube3}, 8.3.1, for the combination of ``locally strongly compact'' and ``quasi-separated'' is equivalent to ``{\it taut}'' (cf.\ {\bf \ref{ch-pre}}.\ref{rem-taut}).
Note also that, since the category $\Acsp^{\mathrm{lft}}_K$ of locally of finite type adic spaces over $K$ is equivalent to the category of locally of finite type rigid spaces over $(\Spf V)^{\rig}$ (where $V$ is the valuation ring of $K$) by the functor $\ZRT$ (\ref{thm-rigidgeomadicspacesff1}), the above comparison also applies to the situation where $\Acsp^{\mathrm{lft}}_K$ (resp.\ $\Acsp^{\mathrm{lsc,lft}}_K$) is replaced by the category of locally of finite type (resp.\ and locally quasi-compact\index{rigid space!locally quasi-compact@locally quasi-compact ---} (\ref{dfn-locallycompactspacerigid})) rigid spaces over $(\Spf V)^{\rig}$; see \ref{thm-berkovichcomparison} below.

\begin{rem}[{\rm cf.\ \cite{Hell}, 4.4}]\label{rem-comparisonberkovichanalyticspacesbis1}{\rm 
Note that the fully faithful functor
$$
\mathscr{B}\colon\Bsp_K\longhookrightarrow\mathscr{M}\mathbf{Ansp}^{\R_+\textrm{-}\mathrm{lft}}_K
$$
is not a categorical equivalence.
Indeed, there exists an $\R_+$-metrized analytic space locally of $\R_+$-finite type (e.g., a locally of finite type adic space) over $K$ whose underlying topological space is not locally strongly compact.
For example, let $\mathscr{X}$ be a locally strongly compact and locally of finite type adic space over $K$, and consider $Y=X\setminus\{x\}$, where $x\in X$ is a closed point of height larger than one.
Then $Y$ is still a locally of finite type adic space, but is not locally strongly compact.
Indeed, $X$ and $Y$ give rise to the same valuation on the locally Hausdorff space $[X]=[Y]$, and hence, if $Y$ is locally strongly compact, then by {\bf \ref{ch-pre}}.\ref{thm-RLSCVsp1}, $X$ and $Y$ have to be homeomorphic by a valuative map, which is absurd.}
\end{rem}

\subsubsection{Comparison with rigid spaces}\label{subsub-comparisonberkovichanalyticspacesbis}
Let $V$ be an $a$-adically complete valuation ring of height one ($a\in\m_V\setminus\{0\}$), and $K=\Frac(V)$.
Let $\mathscr{X}$ be a locally quasi-compact\index{rigid space!locally quasi-compact@locally quasi-compact ---} {\rm (\ref{dfn-locallycompactspacerigid})} locally of finite type rigid space over $\mathscr{S}=(\Spf V)^{\rig}$; recall that the locally-quasi-compactness assumption is satisfied if, for example, $\mathscr{X}$ is coherent ({\bf \ref{ch-pre}}.\ref{prop-locallycompactspacecoherent}).
By \ref{thm-proplocstringcompsp1} we know that the separated quotient $X=[\mathscr{X}]$ is locally compact (and hence is locally Hausdorff).
If, moreover, $\mathscr{X}$ is quasi-separated\index{rigid space!quasi-separated rigid space@quasi-separated ---}, then $[\mathscr{X}]$ is Hausdorff.

Let $\mathscr{U}=\{\mathscr{U}_{\alpha}=(\Spf A_{\alpha})^{\rig}\}_{\alpha\in L}$ be an open covering of $\mathscr{X}$ consisting of affinoids such that $\tau_{\mathscr{U}}=\{[\mathscr{U}_{\alpha}]\}_{\alpha\in L}$ is a net on $[\mathscr{X}]$ (cf.\ {\bf \ref{ch-pre}}.\ref{prop-coveringandnets1topspace}).
Then one has the strictly $K$-analytic space 
$$
\Berk{\mathscr{X}}=([\mathscr{X}],\tau_{\mathscr{U}},\mathcal{A}_{\mathscr{U}})
$$
with the strictly $K$-affinoid atlas $\mathcal{A}_{\mathscr{U}}$ given by the affinoid algebra $\mathcal{A}_{\alpha}=A_{\alpha}[\frac{1}{a}]$ and a homeomorphism $[\mathscr{U}_{\alpha}]\stackrel{\sim}{\rightarrow}\mathscr{M}(\mathcal{A}_{\alpha})$ (\ref{prop-vsbr-berkovichspectrum} and \ref{prop-canonicalsubring}) for each $\alpha\in L$.

It is clear by the construction that, if $\mathscr{V}=\{\mathscr{V}_{\lambda}\}_{\lambda\in\Lambda}$ is another affinoid open covering of $\mathscr{X}$ that gives a refinement of $\{\mathscr{U}_{\alpha}\}_{\alpha\in L}$, then the similarly defined strictly $K$-analytic space $\mathscr{X}'_{\mathrm{B}}=([\mathscr{X}],\tau_{\mathscr{V}},\mathcal{A}_{\mathscr{V}})$ admits a canonical quasi-isomorphism $\mathscr{X}'_{\mathrm{B}}\rightarrow\mathscr{X}_{\mathrm{B}}$.
Hence the strictly $K$-analytic space $\mathscr{X}_{\mathrm{B}}$ is uniquely determined up to quasi-isomorphism only by the rigid space $\mathscr{X}$.
It is then clear that, for an $\mathscr{S}$-morphism $\varphi\colon\mathscr{X}\rightarrow\mathscr{Y}$ of locally quasi-compact locally of finite type rigid spaces, we have the induced morphism $\Berk{\varphi}\colon\Berk{\mathscr{X}}\rightarrow\Berk{\mathscr{Y}}$ in the category of strictly $K$-analytic spaces.

\begin{rem}\label{rem-berkovich1}{\rm 
Let $\mathscr{X}$ and $\tau=\tau_{\mathscr{U}}$ be as above.
We define collection of subsets $\ovl{\tau}$ and $\widehat{\tau}$ of $X=[\mathscr{X}]$ as follows: 
\begin{itemize}
\item $\ovl{\tau}=\{[\mathscr{V}]:\textrm{$\mathscr{V}$ is an affinoid subdomain of some $\mathscr{U}_{\alpha}\in\tau$}\}$;
\item $\widehat{\tau}=\{[\mathscr{V}]:\textrm{$\mathscr{V}$ is a quasi-compact separated open subspace of $\mathscr{X}$}\}$.
\end{itemize}
Then these are nets on $X=[\mathscr{X}]$; cf.\ \cite{Berk2}, \S1.2.}
\end{rem}

Note here that, if $j\colon\mathscr{U}\hookrightarrow\mathscr{X}$ is an open immersion of locally quasi-compact locally of finite type rigid spaces over $\mathscr{S}$, then $\Berk{j}\colon\Berk{\mathscr{U}}\hookrightarrow\Berk{\mathscr{X}}$ is an open immersion in the sense as in \cite{Berk2}, \S1.3.
Indeed, the net $\ovl{\tau}$ on $X=[\mathscr{X}]$ as above induces by restriction a net $\sigma$ of $U=[\mathscr{U}]$; more precisely, $\sigma=\{[\mathscr{V}]\in\ovl{\tau}\,|\,\mathscr{V}\subseteq\mathscr{U}\}$.

\begin{thm}[cf.\ {\rm \cite{Berk2}, 1.6.1}]\label{thm-berkovichcomparison}
The functor $\mathscr{X}\mapsto\Berk{\mathscr{X}}$ establishes a categorical equivalence from the category of all locally quasi-compact locallly of finite type rigid spaces over $\mathscr{S}=(\Spf V)^{\rig}$ to the category $\Bsp^{\ast}_K$ of all strictly $K$-analytic spaces.
Moreover, $\Berk{\mathscr{X}}$ is Hausdorff $($resp.\ paracompact Hausdorff, resp.\ compact$)$ if and only if $X$ is quasi-separated\index{rigid space!quasi-separated rigid space@quasi-separated ---} $($resp.\ paracompact\index{rigid space!paracompact rigid space@paracompact ---} and quasi-separated, resp.\ coherent$)$.
\end{thm}

Note that by \ref{prop-paracompact1r} paracompact quasi-separated rigid spaces are locally quasi-compact, and hence one has the functor $\mathscr{X}\mapsto\Berk{\mathscr{X}}$ for paracompact rigid spaces as in the theorem.

\begin{proof}
Observe that the functor $\mathscr{X}\mapsto\Berk{\mathscr{X}}$ followed by $(\cdot)^{\mathrm{met}}\colon\Bsp_K\rightarrow\mathscr{M}\mathbf{Ansp}^{\R_+}_K$ in \S\ref{subsub-vsbr-BerkovichRmetrizedanalyticspace} coincides (up to natural equivalent) to the composition of $\ZRT$ (cf.\ \ref{thm-rigidgeomadicspacesff1}) followed by $\mathscr{X}\mapsto\mathscr{X}^{\mathrm{met}}$ (cf.\ \ref{cor-Banalytichuberadicspaces11}).
Hence the first assertion follows immediately by \ref{cor-comparisonberkovichanalyticspacesbis1}.
The other assertions follow from {\bf \ref{ch-pre}}.\ref{thm-RLSCVsp1} and {\bf \ref{ch-pre}}.\ref{prop-paracompactseparatedquotient}.
\end{proof}

\begin{prop}\label{prop-goodstrictlyanalyticspace}
Let $\mathscr{X}$ be a quasi-separated locally quasi-compact locally of finite type rigid space over $\mathscr{S}$.
Then the following conditions are equivalent.
\begin{itemize}
\item[{\rm (a)}] The associated strictly $K$-analytic space $\Berk{\mathscr{X}}$ is good {\rm (\ref{dfn-Berkovichanalyticgoodspaces})}.
\item[{\rm (b)}] For any point $x\in\ZR{\mathscr{X}}$, there exists a pair $(\mathscr{U},\mathscr{V})$ of affinoid open neighborhoods of $x$ such that $\ovl{\ZR{\mathscr{U}}}\subseteq\ZR{\mathscr{V}}$, where $\ovl{\,\cdot\,}$ denotes the closure in $\ZR{\mathscr{X}}$.
\end{itemize}
\end{prop}

\begin{proof}
Let us first show (a) $\Rightarrow$ (b).
We may assume that $x$ is a height-one point.
There exists an affinoid open neighborhood $\mathscr{V}$ of $x$ such that $x$ considered as a point of $[\mathscr{X}]$ lies in the interior of $[\mathscr{V}]$.
This implies that $x\in\int_{\mathscr{X}}(\mathscr{V})$.
Take an affinoid open neighborhood $\mathscr{U}$ of $x$ contained in the open subset $\int_{\mathscr{X}}(\mathscr{V})$.
Then by \ref{thm-ZRpoints61} and \ref{cor-interior4} we have $\ovl{\ZR{\mathscr{U}}}\subseteq\ZR{\mathscr{V}}$, as desired.
To show the converse, take any height-one point $x\in[\mathscr{X}]$ and a pair $(\mathscr{U},\mathscr{V})$ of affinoid open neighborhoods as above.
Since $\ovl{\ZR{\mathscr{U}}}\subseteq\ZR{\mathscr{V}}$, $x$ lies in the overconvergent interior $\int_{\mathscr{X}}(\mathscr{V})$ (\ref{cor-interior4}).
Hence $[\mathscr{V}]$ gives a strictly affinoid neighborhood of $x$.
\end{proof}

\addcontentsline{toc}{subsection}{Exercises}
\subsection*{Exercises}
\begin{exer}\label{exer-vsbr-gradedtrivial}{\rm 
Show that, for a filtered ring $(A,F)$, the following conditions are equivalent to each other:
\begin{itemize}
\item[{\rm (a)}] the associated graded ring $\mathrm{Gr}_FA$ is zero;
\item[{\rm (b)}] $F_0=A$;
\item[{\rm (c)}] $1\in F_r$ for some $r<1$.
\end{itemize}
(In this situation, we say that $(A,F)$ is {\it graded trivial}.)}
\end{exer}

\begin{exer}\label{exer-vsbr-R+valuations}{\rm 
Let $K$ be a field, and $V$ a valuation ring for $K$.
Suppose $K$ has a non-archimedean absolute value $\nu\colon K\rightarrow\R_{\geq 0}$ that is a localization of the valuation associated to $V$ (of height $0$ or $1$).
Show that there exists an $\R_+$-valuation $v$ on $K$ with $V$ equal to the associated valuation ring.}
\end{exer}

\begin{exer}\label{exer-vsbr-spectralseminormformula}{\rm 
Let $\mathcal{A}$ be a commutative Banach ring.
Prove the spectral seminorm formula
$$
\|f\|_{\mathcal{A}\dl r^{-1}T\dr,\Sp}=\sup_{n\in\N^n}r^m\|a_m\|_{\mathcal{A},\Sp}
$$
for $f=\sum_{m\in\N_n}a_mT^m$.}
\end{exer}

\begin{exer}\label{exer-vsbr-basicsubsets}{\rm 
(1) Let $\mathcal{A}$ be an $\R_+$-affinoid ring. For $f=(f_0,f_1,\ldots,f_n)$ that generates $\mathcal{A}$, and for an $n$-tuple $r=(r_1,\ldots,r_n)$ of positive real numbers, show that there exist $m_1,\ldots,m_n>0$ such that 
$$
U_0(f,r)=\bigcap^n_{i=1}B(f_0,f_i,r_i)=\bigcap^n_{i=1}U_0((f_0,f_i,1),(r_i,s_i))
$$
for $s_i>m_i$ $(i=1,\ldots,n)$.

(2) Show that a finite intersection of basic subsets is a rational subdomain if and only if its image under the separation map $\sep_{\mathcal{A}}$ is closed in $\mathscr{M}(\mathcal{A})$.}
\end{exer}

\begin{exer}\label{exer-vsbr-valuativespectrum}{\rm 
For a $\Delta$-graded ring $G$, let $\Spec_{\Delta}G$ be the graded spectrum, that is, the set of all homogeneous prime ideals of $G$.

(1) For a homogeneous element $f\in h(G)$, define $D(f)=\{\mathfrak{p}\in\Spec_{\Delta}G:f\not\in\mathfrak{p}\}$. The Zariski topology of $\Spec_{\Delta}G$ is the topology generated by $D(f)$ for $f\in h(G)$.
Show that $\Spec_{\Delta}G$ is coherent and sober.

(2) (Reduction map) Let $\mathcal{A}$ be an $\R_+$-affinoid ring. 
For $x=(z,F)\in\Spec^{\mathrm{val}}\mathcal{A}$, let $P_x$ the kernel of the induced map $\mathrm{Gr}\,mathcal{A}\rightarrow\mathrm{Gr}_F\mathcal{H}(z)\rightarrow k_x$, where $k_x$ is the graded residue field of the graded valuation ring $\mathrm{Gr}_F\mathcal{H}(z)$.
Consider the map 
$$
\sp_{\mathcal{A}}\colon\Spec^{\mathrm{val}}\mathcal{A}\longrightarrow\Spec_{\R_+}\mathrm{Gr}\,\mathcal{A},\quad x\longmapsto P_x
$$
called the {\it $($graded$)$ reduction map}.
Show that $\sp_{\mathcal{A}}$ is continuous, quasi-compact, and surjective.

(3) If $\mathcal{A}$ is a Berkovich's $K$-affinoid algebra, then
$$
\sp_{\mathcal{A}}|_{\mathscr{M}(\mathcal{A})}\colon\mathscr{M}(\mathcal{A})\rightarrow\Spec_{\R_+}\mathrm{Gr}\,\mathcal{A}
$$
is surjective.}
\end{exer}

\section{Appendix: Rigid Zariskian spaces}\label{sec-rigidzariskiansp}
\index{rigid Zariskian space|(}
\subsection{Admissible blow-ups}\label{sub-cohrigidzariskiansp}
In this section we will only deal with Zariskian schemes\index{scheme!Zariskian scheme@Zariskian ---}\index{Zariskian!Zariskian scheme@--- scheme} ({\bf \ref{ch-formal}}, \S\ref{sec-zariskianschemes}) {\em of finite ideal type}\index{scheme!Zariskian scheme@Zariskian ---!Zariskian scheme of finite ideal type@--- --- of finite ideal type}\index{Zariskian!Zariskian scheme@--- scheme!Zariskian scheme of finite ideal type@--- --- of finite ideal type}; a Zariskian scheme $X$ is said to be of finite ideal type if for any affine open set $U=\Spz A$ of $X$ the ring $A$ has a finitely generated ideal of definition.
In the sequel, with a slight abuse of notation, we denote by $\CZs^{\ast}$ the category of coherent Zariskian schemes of finite ideal type and adic morphisms.
Similarly to the formal schemes case ({\bf \ref{ch-formal}}.\ref{cor-extension2}) one can show the following:
\begin{prop}\label{prop-Zariskianscheme2}
Let $X$ be an object of $\CZs^{\ast}$.
Then $X$ has an ideal of definition $\mathscr{I}$ of finite type. \hfill$\square$
\end{prop}

For an affine Zariskian scheme $X=\Spz A$ associated to a pair $(A,I)$ and an admissible ideal $J$, one has the so-called {\em admissible blow-up}\index{blow-up!admissible blow-up@admissible ---} $\pi\colon X'\rightarrow X$; the Zariskian scheme $X'$ is the one given by $X'=\zat{Y^{\prime}}|_{V(I)}$ (cf.\ {\bf \ref{ch-formal}}, \S\ref{subsub-zariskianschemes}), where $Y'\rightarrow Y=\Spec A$ is the blow-up along the ideal $J$.
This construction gives, by gluing, the general concept of {\em admissible blow-ups}\index{blow-up!admissible blow-up@admissible ---} $\pi\colon X'\rightarrow X$ of arbitrary Zariskian schemes and admissible ideals; here, an ideal $\mathscr{J}$ of $\O_X$ is said to be {\em admissible} if it is quasi-coherent of finite type and open or, equivalently, if it is locally written in the form $J^{\lozenge}$ (in the notation as in {\bf \ref{ch-formal}}, \S\ref{subsub-zariskianrings}), where $J$ is a finitely generated $I$-adically open ideal of $A$.

Let $X$ be an object of $\CZs^{\ast}$.
We define the category $\BL_X$ to be the category of all admissible blow-ups of $X$ as follows:
\begin{itemize}
\item[$\bullet$] objects of $\BL_X$ are admissible blow-ups $\pi\colon X'\rightarrow X$; 
\item[$\bullet$] a morphism $\pi'\rightarrow\pi$ between two objects $\pi\colon X'\rightarrow X$ and $\pi'\colon X''\rightarrow X$ means the morphism that makes the resulting traingle 
$$
\xymatrix@C-4ex{X''\ar[rr]\ar[dr]_{\pi'}&&X'\ar[dl]^(.56){\pi}\\ &X}
$$
commutes.
\end{itemize}

By an argument similar to that in \ref{prop-catblowups1} one has the following:
\begin{prop}\label{prop-admissibleblowupzariski31}
{\rm (1)} The category $\BL_X$ is cofiltered $($cf.\ {\rm {\bf \ref{ch-pre}},\ \S\ref{subsub-finalcofinal}}$)$, and $\id_X$ gives the final object.

{\rm (2)} Let us define the ordering on the set $\AId_X$ of all admissible ideals of $\O_X$ as follows$:$ $\mathscr{J}\preceq\mathscr{J}'$ if and only if there exists an admissible ideal $\mathscr{J}''$ such that $\mathscr{J}=\mathscr{J}'\mathscr{J}''$.
Then $\AId^{\opp}_X$ is a directed set, and the functor 
$$
\AId_X\longrightarrow\BL_X
$$
that maps $\mathscr{J}$ to the admissble blow-up along $\mathscr{J}$ is cofinal. \hfill$\square$
\end{prop}

\subsection{Coherent rigid Zariskian spaces}\label{sub-coherentRZS}
\subsubsection{The category of coherent rigid Zariskian spaces}\label{subsub-coherentRZScat}
\begin{dfn}\label{dfn-coherentRZS2}{\rm 
(1) We define the category $\CRz$ as follows:
\begin{itemize}
\item[$\bullet$] objects of $\CRz$ are the same as those of $\CZs^{\ast}$; that is, $\obj(\CRz)=\obj(\CZs^{\ast})$.
For an object $X$ of $\CZs^{\ast}$ we denote by $X^{\rig}$ the same object regarded as an object of $\CRz$; 
\item[$\bullet$] for $X,X'\in\obj(\CZs^{\ast})$ we set
$$
\Hom_{\CRz}(X^{\rig},X^{\prime\rig})=\varinjlim\Hom_{\CZs^{\ast}}(\textrm{--},X'),
$$
where $\Hom_{\CZs^{\ast}}(\textrm{--},X')$ is the functor 
$$
\Hom_{\CZs^{\ast}}(\textrm{--},X')\colon\BL^{\opp}_X\longrightarrow\mathbf{Sets}
$$
that maps $\pi\colon X''\rightarrow X$ to the set $\Hom_{\CZs^{\ast}}(X'',X')$.
\end{itemize}}
\end{dfn}

By \ref{prop-admissibleblowupzariski31} the inductive limit in the above definition can be replaced by a filtered inductive limit along the directed set $\AId^{\opp}_X$.
The composition law for morphsms in the category $\CRz$ 
$$
\Hom_{\CRz}(X^{\rig},X^{\prime\rig})\times\Hom_{\CRz}(X^{\prime\rig},X^{\prime\prime\rig})\longrightarrow\Hom_{\CRz}(X^{\rig},X^{\prime\prime\rig})
$$
is described similarly to the case of rigid (formal) spaces (\S\ref{subsub-cohrigidspace}).

\begin{dfn}\label{dfn-coherentRZS3}{\rm 
(1) An object of the category $\CRz$ is called a {\em coherent rigid Zariskian space}\index{rigid Zariskian space!coherent rigid Zariskian space@coherent ---}.

(2) For an object $X$ of $\CZs^{\ast}$ the coherent rigid Zariskian space $X^{\rig}$ is called the {\em associated $($coherent$)$ rigid space}.
Similarly, for a morphism $f\colon X\rightarrow X'$ in $\CZs^{\ast}$ the {\em associated morphism} of rigid spaces is denoted by $f^{\rig}\colon X^{\rig}\rightarrow X^{\prime\rig}$.}
\end{dfn}

We often denote by 
$$
Q\colon\CZs^{\ast}\longrightarrow\CRz, \qquad X\longmapsto Q(X)=X^{\rig}
$$
the canonical quotient functor.

Similarly to the case of rigid (formal) spaces, one also defines consistently the comma category $\CRz_{\mathscr{S}}$ of coherent rigid Zariskian spaces over a fixed coherent rigid Zariskian space $\mathscr{S}$.

The following assertions are shown similarly to \ref{prop-cohrigidspace3}, \ref{cor-cohrigidspace31}, and \ref{cor-cohrigidspace32}:
\begin{prop}\label{prop-coherentRZS31}
Let $X$ and $X'$ be objects of $\CZs^{\ast}$, and consider the corresponding rigid space $X^{\rig}$ and $Y^{\rig}$.
Then there exists an isomorphism $X^{\rig}\rightarrow Y^{\rig}$ in $\CRz$ if and only if it is represented by a diagram in $\CZs^{\ast}$
$$
\xymatrix@-2ex{
&X''\ar[dr]\ar[dl]\\
X&&X'\rlap{,}}
$$
where both arrows are admissible blow-ups.\hfill$\square$
\end{prop}

\begin{cor}\label{cor-coherentRZS311}
Let $f\colon X\rightarrow X'$ be a morphism in $\CZs^{\ast}$.
Then $f^{\rig}\colon X^{\rig}\rightarrow X^{\prime\rig}$ is an isomorphism if and only there exists a commutative diagram 
$$
\xymatrix@-2ex{
&X''\ar[dr]\ar[dl]\\
X\ar[rr]_{f}&&X'\rlap{,}}
$$
where both $X''\rightarrow X$ and $X''\rightarrow X'$ are admissible blow-ups.\hfill$\square$
\end{cor}

\begin{cor}\label{cor-coherentRZS312}
Consider the diagram 
$$
\xymatrix{
X\ar[d]_f&X\times_{X''}X'\ar[d]^h\ar[l]\\
X''&X'\rlap{.}\ar[l]}
$$
If $f^{\rig}$ is an isomorphism in $\CRz$, then so is $h^{\rig}$.\hfill$\square$
\end{cor}

\begin{dfn}\label{dfn-coherentRZS4}{\rm 
(1) Let $\mathscr{X}$ be a coherent rigid Zariskian space.
A {\em model}\index{model (of a coherent rigid Zariskian space)} of $\mathscr{X}$ is a couple $(X,\phi)$ consisting of $X\in\obj(\CZs^{\ast})$ and an isomorphism $\phi\colon X^{\rig}\stackrel{\sim}{\rightarrow}\mathscr{X}$.
A model $(X,\phi)$ is said to be {\em distinguished} if $X$ is $\mathscr{I}$-torsion free, where $\mathscr{I}$ denotes an ideal of definition of $X$.

(2) Let $\varphi\colon\mathscr{X}\rightarrow\mathscr{X}'$ be a morphism of coherent rigid Zariskian spaces.
A {\em model} of $\varphi$ is a triple $(f,\phi,\psi)$ consisting of a morphism $f\colon X\rightarrow X'$ of $\CZs^{\ast}$ and isomorphisms $\phi\colon X^{\rig}\stackrel{\sim}{\rightarrow}\mathscr{X}$ and $\psi\colon X^{\prime\rig}\stackrel{\sim}{\rightarrow}\mathscr{X}'$ such that the resulting square
$$
\xymatrix{X^{\rig}\ar[d]_{\phi}\ar[r]^{f^{\rig}}&X^{\prime\rig}\ar[d]^{\psi}\\ \mathscr{X}\ar[r]_{\varphi}&\mathscr{X}'}
$$
commutes.
A model $(f,\phi,\psi)$ is said to be {\em distinguished}\index{model (of a coherent rigid Zariskian space)!distinguished model (of a rigid Zariskian space)@distinguished ---} if $X$ and $X'$ are distinguished models of $\mathscr{X}$ and $\mathscr{X}'$, respectively.}
\end{dfn}

Similarly to the case of coherent rigid (formal) spaces, one can define the category $\FM_{\mathscr{X}}$ of models of $\mathscr{X}$ and the category $\FM_{\varphi}$ of models of $\varphi$; the details are left to the reader.
It is easy to see that these categories are cofiltered.
We denote by $\FM^{\dist}_{\mathscr{X}}$ (resp.\ $\FM^{\dist}_{\varphi}$) the full subcategory of $\FM_{\mathscr{X}}$ (resp.\ $\FM_{\varphi}$) consisting of distinguished models.
The following propositions are the counterparts of \ref{prop-cohrigidspacedist100} and \ref{prop-cohrigidspacedist1}, which are shown by a similar argument:
\begin{prop}\label{prop-coherentRZS41}
{\rm (1)} If $(X,\phi)$ is a distinguished model of $\mathscr{X}$ and $\pi\colon X'\rightarrow X$ is an admissible blow-up, then $(X',\phi\circ\pi^{\rig})$ is a distinguished model of $\mathscr{X}$.

{\rm (2)} If, moreover, $\mathscr{I}$ $($an ideal of definition of $X)$ is an invertible ideal, then $\mathscr{I}\O_{X'}$ is invertible. \hfill$\square$
\end{prop}

\begin{prop}\label{prop-coherentRZS42}
The categories $\FM^{\dist}_{\mathscr{X}}$ and $\FM^{\dist}_{\varphi}$ are cofiltered, and the inclusions $\FM^{\dist}_{\mathscr{X}}\hookrightarrow\FM_{\mathscr{X}}$ and $\FM^{\dist}_{\varphi}\hookrightarrow\FM_{\varphi}$ are cofinal. \hfill$\square$
\end{prop}

\subsubsection{Visualization}\label{subsub-coherentRZS}
Let $\mathscr{X}$ be a coherent rigid Zariskian space.
One can define the {\em associated Zariski-Riemann space}\index{Zariski-Riemann space} $\ZR{\mathscr{X}}$ in entirely analogous way as in \S\ref{sub-ZRdef}.
Let 
$$
\sp_X\colon\ZR{\mathscr{X}}\longrightarrow X
$$
be the specialization map\index{specialization map}, where $X$ is a model of $\mathscr{X}$.
Similarly to \ref{thm-ZRcompact} and \ref{prop-zariskiriemanntoptop} we have:
\begin{thm}\label{thm-ZRzariskiancompact}
Let $\mathscr{X}$ be a coherent rigid Zariskian space, and $\ZR{\mathscr{X}}$ the associated Zariski-Riemann space.

{\rm (1)} The topological space $\ZR{\mathscr{X}}$ is coherent and sober $({\bf \ref{ch-pre}}.\ref{dfn-quasicompact1})$.

{\rm (2)} For any model $X$ of $\mathscr{X}$ the specialization map $\sp_X$ is quasi-compact $({\bf \ref{ch-pre}}.\ref{dfn-quasicompactness}\ (2))$ and closed. \hfill$\square$
\end{thm}

\begin{prop}\label{prop-Zariskianzariskiriemanntoptop}
{\rm (1)} For any quasi-compact open subset $\mathfrak{U}$ of $\ZR{\mathscr{X}}$ there exists a model $X$ and a quasi-compact open subset $U$ of $X$ such that $\mathfrak{U}=\sp^{-1}_X(U)$.

{\rm (2)} Let $X$ be an object of $\CZs^{\ast}$, and $U$ a quasi-compact open subset of $X$.
Set $\mathscr{X}=X^{\rig}$ and $\mathscr{U}=U^{\rig}$.
Then the induced map $\ZR{\mathscr{U}}\rightarrow \ZR{\mathscr{X}}$ maps $\ZR{\mathscr{U}}$ homeomorphically onto the quasi-compact open subset $\sp^{-1}_X(U)$. \hfill$\square$
\end{prop}

One can define general rigid Zariskian spaces and their assocaited Zariski-Riemann triples in a similar way as in the case of rigid (formal) spaces.
The details are left to the reader.
\index{rigid Zariskian space|)}

\section{Appendix: Classical Zariski-Riemann spaces}\label{sec-classicalZR}
\subsection{Birational geometry}\label{sub-birationalgeom}
\subsubsection{Basic terminologies}\label{subsub-birationalgeomterm}
Throughout this section we fix once for all a coherent ($=$ quasi-compact and quasi-separated) scheme $S$ and work entirely in the category $\CAs_S$ of coherent $S$-algebraic spaces\index{coherent!coherent algebraic space@--- algebraic space}\index{algebraic space!coherent algebraic space@coherent ---} with $S$-morphisms.
Notice that all arrows in $\CAs_S$ are automatically coherent.

Our main objects in this subsection are the pairs of the form $(X,U)$ consisting of an object $X\in\obj(\CAs_S)$ and a quasi-compact open subspace $U$ of $X$; notice that $U$ is again coherent, hence belonging to $\CAs_S$.
With the objects of this kind, we are to consider the theory of the so-called {\em $U$-admissible birational geometry}\index{admissible!U-admissible birational geometry@$U$-{---} birational geometry}.
Here, in order to include the usual birational geometry (without reference to any $U$), we admit the case where $U=\emptyset$.

\begin{dfn}\label{dfn-Uadmissiblemap}{\rm 
Let $X$ be a coherent $S$-algebraic space, and $U\subseteq X$ a quasi-compact open subspace.
A morphism $f\colon Y\rightarrow X$ in $\CAs_S$ is said to be {\em $U$-admissible} if the induced morphism $f^{-1}(U)\rightarrow U$ is an isomorphism.}
\end{dfn}

\begin{dfn}[cf.\ {\cite[$\mathbf{I}$, \S7]{EGA}}]\label{dfn-birationalgeom1}{\rm 
Let $X$ and $Y$ be coherent $S$-algebraic spaces, and $U\subseteq X$ a quasi-compact open subspace.

(1) A {\em $(U$-admissible$)$ rational map}\index{admissible!U-admissible rational map@$U$-{---} rational map} $f\colon Y\dashrightarrow X$ is an equivalence class of a pair $(V,f_V)$ consisting of a quasi-compact dense open subspace $V$ of $Y$ and a $U$-admissible $S$-morphism $f_V\colon V\rightarrow X$, where the equivalence relation is defined as follows:
$(V,f_V)$ and $(V',f_{V'})$ are equivalent if the morphisms $f_V$ and $f_{V'}$ coincide on a quasi-compact dense open subspace $W$ of $Y$ such that $W\subseteq V\cap V'$; since the intersection of two quasi-compact dense open subspaces of $Y$ is again a quasi-compact open dense subspace, this defines an equivalence relation.

(2) A $U$-admissible rational map $f\colon Y\dashrightarrow X$ is said to be {\em birational}\index{birational!U-admissible birational map@$U$-admissible --- map} if it is the equivalence class of a pair $(V,f_V)$ such that $f_V\colon V\rightarrow X$ is a $U$-admissible open immersion onto a quasi-compact dense open subspace of $X$.

(3) A {\em $(U$-admissible$)$ birational morphism}\index{birational!birational morphism@--- morphism} $f\colon Y\rightarrow X$ is a $U$-admissible $S$-morphism that is birational, that is, there exists a quasi-compact open dense subspace $V$ of $Y$ such that $f|_V\colon V\rightarrow X$ is an open immersion onto a quasi-compact dense open subspace of $X$.

(4) A $U$-admissible $S$-morphism $f\colon Y\rightarrow X$ is called a {\em $(U$-admissible$)$ $S$-modification}\index{admissible!U-admissible modification@$U$-{---} modification}\index{modification!U-admissible modification@$U$-admissible ---} (or simply {\em modification}) if it is proper and birational.}
\end{dfn}

The composition $g\circ f$ of a $U$-admissible rational map $f\colon Y\dashrightarrow X$ and an $f^{-1}(U)$-admissible rational map $g\colon Z\dashrightarrow Y$ is defined in an obvious way and is a $U$-admissible rational map.
A $U$-admissible rational map $f\colon Y\dashrightarrow X$ is birational if and only if there exists an $f^{-1}(U)$-admissible rational map $g\colon X\dashrightarrow Y$ such that $g\circ f$ (resp.\ $f\circ g$) is represented by the identity map on a quasi-compact open dense subspace of $Y$ (resp.\ $X$).
Note that, if $X$ and $Y$ are schemes having finitely many irreducible components, our definition of birational morphisms coincides with that in \cite[(2.3.4)]{EGAInew}.

\subsubsection{$U$-admissible blow-ups}\label{subsub-birationalgeomblowups}
\index{admissible!U-admissible blow-up@$U$-{---} blow-up|(}
\index{blow-up!U-admissible blow-up@$U$-admisible ---|(}
We continue with working in the situation as above; let $X$ be a coherent $S$-algebraic space, and $U$ a quasi-compact open subspace of $X$.
\begin{dfn}\label{dfn-Uadmissibleideal}{\rm 
A {\em $U$-admissible ideal}\index{admissible!U-admissible ideal@$U$-{---} ideal} is a quasi-coherent ideal $\mathscr{J}\subset\O_X$ of finite type such that the closed subspace $V(\mathscr{J})$ of $X$ corresponding to $\mathscr{J}$ is disjoint from $U$.}
\end{dfn}

We denote by $\AId_{(X,U)}$ the set of all $U$-admissible ideals.
One can introduce the ordering $\preceq$ to the set $\AId_X$ by a similar manner as in \ref{prop-catblowups1}.
Moreover, similarly to {\bf \ref{ch-formal}}.\ref{prop-admissibleideal1x}, for an $S$-morphism $Y\rightarrow X$ of $S$-algebraic spaces one has the induced map 
$$
\AId_{(X,U)}\longrightarrow\AId_{(Y,f^{-1}(U))},\qquad\mathscr{J}\longmapsto\mathscr{J}\O_Y.
$$

\begin{dfn}\label{dfn-Uadmissibleblowups}{\rm 
Let $X$ be a coherent $S$-algebraic space, and $U$ a quasi-compact open subspace of $X$.
A {\em $U$-admissible blow-up} of $X$ is a blow-up of $X$ along a $U$-admissible ideal $\mathscr{J}$, that is, the morphism
$$
X_{\mathscr{J}}=\Proj\bigoplus_{n\geq 0}\mathscr{J}^n\longrightarrow X.
$$}
\end{dfn}

Clearly, the $U$-admissible blow-up is $U$-admissible as a morphism between $S$-algebraic spaces.
Notice that the $U$-admissible blow-up is proper, since $\mathscr{J}$ is of finite type.
It is, in particular, birational if the open complement of the closed subscheme $V(\mathscr{J})$ in $X$ is dense.

\begin{dfn}\label{dfn-Uadmissiblestricttrans}{\rm 
Let $X$ be a coherent $S$-algebraic space, $U\subseteq X$ a quasi-compact open subspace, and $\mathscr{J}$ a $U$-admissible ideal.
Let $\pi\colon X'\rightarrow X$ the $U$-admissible blow-up of $X$ along $\mathscr{J}$.

(1) For an $\O_X$-module $\mathscr{F}$ the {\em strict transform}\index{strict transform} of $\mathscr{F}$ by $\pi$ is the $\O_{X'}$-module given by 
$$
\pi'\mathscr{F}=\pi^{\ast}\mathscr{F}/(\pi^{\ast}\mathscr{F})_{\textrm{$\mathscr{J}$-}\mathrm{tor}}.
$$
Notice that $\pi'\mathscr{F}$ is quasi-coherent if $\mathscr{F}$ is quasi-coherent.

(2) For a morphism $f\colon Y\rightarrow X$ of coherent $S$-algebraic spaces the {\em strict transform}\index{strict transform} of $f$ by $\pi$ is the composite morphism $Y'\hookrightarrow Y\times_XX'\rightarrow X'$, where the first arrow is the closed immersion given by the ideal $(\O_{Y\times_XX'})_{\textrm{$\mathscr{J}$-}\mathrm{tor}}$.}
\end{dfn}

Notice that in the situation as in (2) the morphism $Y'\rightarrow Y$ is the $f^{-1}(U)$-admissible blow-up along the ideal $\mathscr{J}\O_Y$.

The following fact will be frequently used (often tacitly) in the sequel:
\begin{prop}\label{prop-birationalgeom023}
Let $X$ be a coherent algebraic space, and $U$ a quasi-compact open subspace of $X$.
Let $Y\hookrightarrow X$ be a quasi-compact open immersion, and set $V=Y\cap U$.

{\rm (1)} For any $U$-admissible blow-up $X'\rightarrow X$ the strict transform $Y'\rightarrow Y$ coincides with the base change $X'\times_XY\rightarrow Y$ and is the $V$-admissible blow-up of $Y$ along the restriction of the blow-up center of $X'\rightarrow X$.

{\rm (2)} {\rm (Extension of $U$-admissible blow-up)} Conversely, for any $V$-admissible blow-up $Y'\rightarrow Y$ there exists a $U$-admissible blow-up $X'\rightarrow X$ that admits the Cartesian diagram
$$
\xymatrix{Y'\,\ar@{^{(}->}[r]\ar[d]&X'\ar[d]\\ Y\,\ar@{^{(}->}[r]&X.}
$$
\end{prop}

\begin{proof}
(1) is clear.
To show (2), let $\mathscr{J}$ be the blow-up center of $Y'\rightarrow Y$.
The direct image $j_{\ast}\mathscr{J}$ by the open immersion $j\colon Y\rightarrow X$ is a quasi-coherent ideal of $\O_X$.
By {\bf \ref{ch-pre}}.\ref{thm-grusonraynaudlimit}, we have a quasi-coherent ideal $\mathscr{I}$ of $\O_X$ of finite type such that the support of the corresponding closed subspace $Z=V(\mathscr{I})$ is $X\setminus U$ and that $\mathscr{I}|_V\subseteq\mathscr{J}$.
The closed subspace $Z$ is a coherent algebraic space, and $Z\cap Y$ is a quasi-compact open subspace of $Z$.
Consider the quasi-coherent ideal $\mathscr{J}/(\mathscr{I}|_V)$ of finite type on $Z\cap Y$, and extend it to a quasi-coherent ideal of finite type on $Z$ ({\bf \ref{ch-pre}}.\ref{thm-grusonraynaudlimit} and {\bf \ref{ch-pre}}.\ref{cor-grusonraynaudlimit}).
Pulling it back to $\O_X$, we get a $U$-admissible ideal $\til{\mathscr{J}}$ that extends $\mathscr{J}$.
Then the $U$-admissible blow-up $X'\rightarrow X$ along $\til{\mathscr{J}}$ satisfies the properties as in (2).
\end{proof}

Here we include some facts on $U$-admissible blow-ups quoted from \cite[Premi\`ere partie, \S 5]{RG}, which are fundamental in the $U$-admissible birational geometry:
\begin{prop}[{\cite[Premi\`ere partie, (5.1.4)]{RG}}]\label{prop-blowup3a}
Let $X$ be a coherent scheme, and $U\subseteq X$ a quasi-compact open subset.
Let $\pi\colon X'\rightarrow X$ be an $U$-admissible blow-up along a quasi-coherent ideal $\mathscr{J}\subseteq\O_X$ of finite type, and $X''\rightarrow X'$ an $\pi^{-1}(U)$-admissible blow-up along a quasi-coherent ideal $\mathscr{J}'\subseteq\O_{X'}$ of finite type.
Then there exists a quasi-coherent ideal $\mathscr{J}''$ of $\O_X$ of finite type such that $\mathscr{J}''\O_{X'}$ coincides with $\mathscr{J}^m\cdot\mathscr{J}^{\prime n}\O_{X'}$ for some positive integers $m$ and $n$ and that the composition $X''\rightarrow X$ coincides up to canonical isomorphisms with the $U$-admissible blow-up along the ideal $\mathscr{J}\cdot\mathscr{J}''$. \hfill$\square$
\end{prop}

The proposition says, in particular, that the composition of two $U$-admissible blow-ups is again a $U$-admissible blow-up.

\begin{prop}\label{prop-birationalgeom01}
Let $X$ be a coherent algebraic space, $U$ a quasi-compact open subspace of $X$, and $f\colon Y\rightarrow X$ a separated morphism of algebraic spaces of finite type.
Suppose that the induced morphism $f^{-1}(U)\rightarrow U$ is an open immersion.
Then there exists an $U$-admissible blow-up $\pi\colon X'\rightarrow X$ such that, if $Y'$ denotes the strict transform of $Y$, the resulting morphism $Y'\rightarrow X'$ is an open immersion$:$
$$
\xymatrix{Y'\,\ar@{^{(}->}[r]^{f'}\ar[d]_{\pi'}&X'\ar[d]^{\pi}\\ Y\,\ar[r]_f&X.}
$$
Moreover, if $U$ is dense in $X$, then the map $\pi'\colon Y'\rightarrow Y$ is an $f^{-1}(U)$-admissible modification. \hfill$\square$
\end{prop}

The first assertion is a special case of \cite[Premi\`ere partie, (5.7.11)]{RG}.
The second is clear, since $Y'\rightarrow Y$ is proper.

\begin{prop}\label{prop-birationalgeom02}
Let $X$ be a coherent algebraic space, $U$ a quasi-compact open subspace of $X$, and $f\colon Y\rightarrow X$ a $U$-admissible proper morphism.
Then there exists a $U$-admissible blow-up $\pi\colon X'\rightarrow X$ such that, if $Y'$ denotes the strict transform of $Y$, the resulting morphism $Y'\rightarrow X'$ is an isomorphism$:$
$$
\xymatrix@-2ex{&X'\ar[dl]_{\pi}\ar@/^1pc/[ddl]^{\pi'}\\ Y\ar[d]_f\\ X.}
$$
In other words, there exists an $f^{-1}(U)$-admissible blow-up $\pi\colon X'\rightarrow Y$ such that the composition $X'\rightarrow Y\rightarrow X$ is a $U$-admissible blow-up. \hfill$\square$
\end{prop}

This is a special case of \cite[Premi\`ere partie, (5.7.12)]{RG}.
\index{blow-up!U-admissible blow-up@$U$-admisible ---|)}
\index{admissible!U-admissible blow-up@$U$-{---} blow-up|)}

\subsubsection{The correspondence diagram}\label{subsub-Uadmissiblejoin}
\index{correspondence diagram|(}
Let $S$ be a coherent scheme, and consider the diagram 
$$
\begin{xy}
(0,0)="X"*{X},
"X"+<-4em,-6ex>="U"*{U},
"X"+<0em,-12ex>="Y"*{Y}
\ar@{^{(}->}^{j}"U"+<1em,1.5ex>;"X"+<-.7em,-1.2ex>
\ar"U"+<.6em,-1.5ex>;"Y"+<-.7em,1.2ex>
\end{xy}\leqno{(\ast)}
$$
of coherent $S$-algebraic spaces, where $j$ is a quasi-compact open immersion, and $Y$ is separated of finite type over $S$.
Consider the induced map
$$
i\colon U\longrightarrow X\times_SY.
$$
Since the composition $U\stackrel{i}{\rightarrow}X\times_SY\stackrel{p_1}{\rightarrow}X$, where $p_1$ is the first projection, is the open immersion $j$, the morphism $i$ is a quasi-compact immersion (cf.\ \cite[(4.3.6) (iv)]{EGAInew}).
\begin{dfn}\label{dfn-correspondencediagram}{\rm 
The scheme-theoretic closure (cf.\ \cite{Knu}, II.4.6) of the image of $U$ in $X\times_SY$ is called the {\em join}\index{join} of $X$ and $Y$ along $U$ and is denoted by $X\ast^UY$.}
\end{dfn}

The join $X\ast^UY$ comes with the morphisms $p\colon X\ast^UY\rightarrow X$ and $q\colon X\ast^UY\rightarrow Y$ induced, respectively, from the first and the second projections.
The following proposition is clear, and the proof is left to the reader:
\begin{prop}\label{prop-correspondencediagram}
The diagram $(\ast)$ extends to the commutative diagram
$$
\begin{xy}
(0,0)="X"*{X},
"X"+<-4em,-8ex>="U"*{U},
"X"+<4em,-8ex>="XUY"*{X\rlap{$\ast^UY$}},
"X"+<0em,-16ex>="Y"*{Y}
\ar@{^{(}->}"U"+<1em,1.5ex>;"X"+<-.7em,-1.2ex>
\ar@{^{(}->}"U"+<1em,-.5ex>;"XUY"+<-.7em,-.5ex>
\ar"U"+<.6em,-1.5ex>;"Y"+<-.7em,1.2ex>
\ar_{p}"XUY"+<-.6em,1.2ex>;"X"+<.7em,-1.2ex>
\ar^{q}"XUY"+<-.6em,-1.5ex>;"Y"+<.7em,1.2ex>
\end{xy}
$$
of $S$-algebraic spaces, where the horizontal arrow is an open immersion.
Moreover, we have$:$
\begin{itemize}
\item[{\rm (1)}] the morphism $p$ is $U$-admissible$;$ if $U\rightarrow Y$ is also an open immersion, then $q$ is $U$-admissible$;$
\item[{\rm (2)}] if $X$ $($resp.\ $Y)$ is proper over $S$, then $q$ $($resp.\ $p)$ is proper. \hfill$\square$
\end{itemize}
\end{prop}

Notice that, if $X$ and $Y$ are schemes, then all what we have done (and what we will do in this paragraph) can be done within the category of schemes.

\begin{prop}\label{prop-correspondencediagramfurther}
Consider the diagram $(\ast)$ of coherent $S$-algebraic spaces.

{\rm (1)} There exists a $U$-admissible blow-up $X'\rightarrow X$ sitting in a commutative diagram of the form
$$
\begin{xy}
(0,0)="X"*{X'},
"X"+<-4em,-8ex>="U"*{U},
"X"+<4em,-8ex>="XUY"*{Z\rlap{,}},
"X"+<0em,-16ex>="Y"*{Y}
\ar@{^{(}->}"U"+<1em,1.5ex>;"X"+<-.8em,-1.2ex>
\ar@{^{(}->}"U"+<1em,-.5ex>;"XUY"+<-.7em,-.5ex>
\ar"U"+<.6em,-1.5ex>;"Y"+<-.7em,1.2ex>
\ar@{_{(}->}_{p'}"XUY"+<-1em,1.5ex>;"X"+<.6em,-1.2ex>
\ar^{q'}"XUY"+<-.6em,-1.5ex>;"Y"+<.7em,1.2ex>
\end{xy}
$$
where $p'$ is the strict transform\index{strict transform} of $p\colon X\ast^UY\rightarrow X$, which is an open immersion.
Moreover$:$
\begin{itemize}
\item[{\rm (a)}] if $X$ is proper over $S$, then $q'$ is proper$;$
\item[{\rm (b)}] if $Y$ is proper over $S$, then one can take $X'$ as above such that the map $p'$ is an isomorphism.
\end{itemize}

{\rm (2)} Suppose that the map $U\rightarrow Y$ is an open immersion.
Then there exist a $U$-admissible blow-up\index{admissible!U-admissible blow-up@$U$-{---} blow-up}\index{blow-up!U-admissible blow-up@$U$-admisible ---} $X'$ of $X$ and a $U$-admissible blow-up $Y'$ of $Y$ sitting in a commutative diagram of the form
$$
\begin{xy}
(0,0)="X"*{X'},
"X"+<-4em,-8ex>="U"*{U},
"X"+<4em,-8ex>="XUY"*{Z\rlap{,}},
"X"+<0em,-16ex>="Y"*{Y'}
\ar@{^{(}->}"U"+<1em,1.5ex>;"X"+<-.8em,-1.2ex>
\ar@{^{(}->}"U"+<1em,-.5ex>;"XUY"+<-.7em,-.5ex>
\ar@{^{(}->}"U"+<.7em,-2.1ex>;"Y"+<-.7em,1.2ex>
\ar@{_{(}->}_{p'}"XUY"+<-1em,1.5ex>;"X"+<.6em,-1.2ex>
\ar@{_{(}->}^{q'}"XUY"+<-.7em,-2.5ex>;"Y"+<.7em,1.2ex>
\end{xy}
$$
consisting of open immersions, where $Z$ is a $U$-admissible blow-up of the join $X\ast^UY$.
Moreover, if $X$ $($resp.\ $Y)$ is proper, then one can take $X'$ and $Y'$ as above such that $q'$ $($resp.\ $p')$ is an isomorphism.
\end{prop}

\begin{proof}
(1) By \ref{prop-birationalgeom01} there exists a $U$-admissible blow-up $X'$ of $X$ such that the strict transform $Z\rightarrow X'$ of the map $X\ast^UY\rightarrow X$ is an open immersion.
The assertion (a) is clear.
If $Y$ is proper, then the map $X\ast^UY\rightarrow X$ is proper and $U$-admissible.
Then by \ref{prop-birationalgeom02}, one can take a $U$-admissible blow-up $X'$ of $X$ such that the strict transform $p'\colon Z\rightarrow X'$ of $X\ast^UY\rightarrow X$ is an isomorphism, whence (b).

(2) First, as in (1), we take $X'$ and $Z$ such that the resulting map $Z\rightarrow X$ is an open immersion.
Since $Z\rightarrow Y$ is also $U$-admissible, one can do the same for this morphism to obtain the strict transform $Z'\rightarrow Y'$ that is an open immersion.
Then by \ref{prop-birationalgeom023} there exists a $U$-admissible blow-up $X''\rightarrow X'$ that induces the $U$-admissible blow-up $Z'\rightarrow Z$.
Hence we obtain the diagram of the form as in (2).
The other assertion is shown similarly to the assertion (b) of (1); here we use \ref{prop-birationalgeom023} to keep the resulting diagram consisting only of open immersions.
\end{proof}
\index{correspondence diagram|)}

\subsubsection{Birational category}\label{subsub-birationalcategory}
Let $S$ be a coherent scheme with finitely many irreducible components.
Consider the category $\mathscr{C}$ defined as follows:
\begin{itemize}
\item objects are separated $S$-algebraic spaces $X$ of finite type;
\item arrows $Y\rightarrow X$ between objects $Y$ and $X$ are given by ($\emptyset$-admissible) rational maps $f\colon Y\dashrightarrow X$.
\end{itemize}
We define the category $\Bir_S$ as the localized category of $\mathscr{C}$ by the set of all ($\emptyset$-admissible) birational maps.

For any object $X$ of the category $\mathscr{C}$, that is, a separated $S$-algebraic space $X$ of finite type, we denote by 
$$
k(X)
$$
the same object considered as an object in $\Bir_S$, following the customary used notation in classical birational geometry.

\begin{dfn}\label{dfn-birationalgeom2models}{\rm 
Let $\mathscr{X}$ be an object of $\Bir_S$.
A {\em model}\index{model (in birational geometry)} of $\mathscr{X}$ is a pair $(X,\phi)$ consisting of a separated $S$-algebraic space $X$ of finite type and an isomorphism 
$$
\phi\colon k(X)\stackrel{\sim}{\longrightarrow}\mathscr{X}
$$
in the category $\Bir_S$.}
\end{dfn}

Notice that, if $X_1$ and $X_2$ are models of $\mathscr{X}$, then there exists a separated $S$-algebraic space $U$ of finite type and birational open immersions
$$
X_1\stackrel{j_1}{\longhookleftarrow}U\stackrel{j_2}{\longhookrightarrow}X_2,
$$
that is, the images in $j_1(U)$ and $j_2(U)$ are dense in $X_1$ and $X_2$, respectively.
Hence one can consider the join $X_1\ast X_2=X_1\ast^UX_2$, which is again a model of $\mathscr{X}$.
Notice that, as we discussed in \S\ref{subsub-Uadmissiblejoin}, the join\index{join} $X_1\ast X_2$ depends only on the closures of $U$ (that is, $X_1$ and $X_2$), and hence is independent on the choice of $U$.

\subsection{Classical Zariski-Riemann spaces}\label{sub-classicalZRsp}
\subsubsection{The cofiltered category of modifications}\label{subsub-filtcatmodifications}
Let $S$ be a coherent scheme, $X$ a coherent $S$-algebraic space, and $U\subseteq X$ a quasi-compact open subspace.
We define the category $\MD_{(X,U)}$ as follows: the objects are $U$-admissible modifications\index{modification!U-admissible modification@$U$-admissible ---} (\ref{dfn-birationalgeom1} (4)) $X'\rightarrow X$ and arrows between two such morphisms $X'\rightarrow X$ and $X''\rightarrow X$ are $X$-morphisms $X'\rightarrow X''$.
\begin{prop}\label{prop-classicalZRsp01}
{\rm (1)} The category $\MD_{(X,U)}$ is cofiltered, and $\id_X$ gives the final object.

{\rm (2)} The opposite ordered set $\AId^{\opp}_{(X,U)}$ of the set of all $U$-admissible ideals with the above-mentioned ordering is a directed set, and the functor 
$$
\AId_{(X,U)}\longrightarrow\MD_{(X,U)}
$$
that maps $\mathscr{J}$ to the blow-up $X_{\mathscr{J}}\rightarrow X$ along $\mathscr{J}$ is cofinal.
\end{prop}

\begin{proof}
(1) What to show are the following:
\begin{itemize}
\item[(a)] for two $U$-admissible modifications $X'\rightarrow X$ and $X''\rightarrow X$ there exists a $U$-admissible modification $X'''\rightarrow X$ and $X$-morphisms $X'''\rightarrow X'$ and $X'''\rightarrow X''$;
\item[(b)] for two $U$-admissible modifications $X'\rightarrow X$ and $X''\rightarrow X$ and two $X$-morphisms $f_0,f_1\colon X''\rightarrow X'$ there exist a $U$-admissible modification $X'''\rightarrow X$ and an $X$-morphism $g\colon X'''\rightarrow X''$ such that $f_0\circ g=f_1\circ g$.
\end{itemize}

Let us prove (a). 
By \ref{prop-birationalgeom02} there exists $U\times_XX'$-admissible (resp.\ $U\times_XX''$-admissible) blow-up $Y'\rightarrow X'$ (resp.\ $Y''\rightarrow X''$) such that the composition $Y'\rightarrow X'\rightarrow X$ (resp.\ $Y''\rightarrow X''\rightarrow X$) is an $U$-admissible blow-up.
Let $\mathscr{J}'$ (resp.\ $\mathscr{J}''$) be the blow-up center of $Y'\rightarrow X$ (resp.\ $Y''\rightarrow X$), and set $\mathscr{J}'''=\mathscr{J}'\mathscr{J}''$.
Let $X'''\rightarrow X$ be the $U$-admissible blow-up along $\mathscr{J}'''$.
Then there are $X$-morphisms $X'''\rightarrow Y'$ and $X'''\rightarrow Y''$ by the universality of blow-ups (cf.\ \cite[Premi\`ere partie, (5.1.2) (i)]{RG}).
The compositions $X'''\rightarrow Y'\rightarrow X'$ and $X'''\rightarrow Y''\rightarrow X''$ are the desired morphisms.

Next let us show (b). 
Take $Y'$, $Y''$, $\mathscr{J}'$, and $\mathscr{J}''$ as above.
The blow-up $Y'\rightarrow X'$ (resp.\ $Y''\rightarrow X''$) is along $\mathscr{J}'\O_{X'}$ (resp.\ $\mathscr{J}''\O_{X''}$).
Hence there exists two maps $h_0,h_1\colon Y''\rightarrow Y'$ induced respectively from $f_0$ and $f_1$.
Take the $U$-admissible blow-up $X'''\rightarrow X$ along $\mathscr{J}'''=\mathscr{J}'\mathscr{J}''$.
Then by the uniqueness in the universality of blow-ups, the two compositions $X'''\rightarrow Y''\rightarrow Y'$ where the second maps are given by $h_0$ and $h_1$ are the same.
Hence if one define $g$ to be the composition $X'''\rightarrow Y''\rightarrow X''$, we have $f_0\circ g=f_1\circ g$, as desired.

(2) is clear by the argument above.
\end{proof}

\subsubsection{The classical Zariski-Riemann spaces}\label{subsub-classicalZRspace}
\index{Zariski-Riemann space!classical Zariski-Riemann space@classical ---|(}
\begin{dfn}\label{dfn-classicalZRsp1}{\rm 
Let $X$ be a coherent $S$-{\em scheme}, and $U$ a quasi-compact open subset of $X$.
Consider the functor $S_{(X,U)}\colon\MD_{(X,U)}\rightarrow\LRsp$ that maps $X'\rightarrow X$ to the underlying locally ringed space of $X'$.
The (classical) {\em Zariski-Riemann space} associated to the pair $(X,U)$, denoted by $\ZR{X}_{U}$, is the underlying topological space of the limit 
$$
\varprojlim S_{(X,U)}.
$$
The limit is canonically identified with the limit taken along a directed set
$$
\varprojlim_{\mathscr{J}\in\AId_{(X,U)}}X_{\mathscr{J}},
$$
which guarantees the existence of the limit (cf.\ {\bf \ref{ch-pre}}.\ref{prop-limLRS}).
We have for any modification $X'\rightarrow X$ of $X$ the canonical projection ({\em specialization map}\index{specialization map}) $\ZR{X}_{U}\rightarrow X'$, which we denote by $\sp_{X'}$.
The inductive limit sheaf 
$$
\O_{\ZR{X}_{U}}=\varinjlim_{\mathscr{J}\in\AId_{(X,U)}}\sp^{-1}_{X_{\mathscr{J}}}\O_{X_{\mathscr{J}}}
$$
is a sheaf of local rings, called the {\em structure sheaf} of $\ZR{X}_{U}$.}
\end{dfn}

Any inclusion $U_1\hookrightarrow U_2$ of quasi-compact dense open subsets of $X$ induces a map
$$
\ZR{X}_{U_1}\longrightarrow\ZR{X}_{U_2}, 
$$
which is quasi-compact due to {\bf \ref{ch-pre}}.\ref{thm-projlimcohspacepres} (1) and (3).

Notice that, in the above construction, one can actually replace $X$ by the closure $\ovl{U}$ of $U$ in $X$ without changing the formation of the locally ringed space $\ZR{X}_U$; indeed, any irreducible component that does not touch $U$ can be effaced by a $U$-admissible blow-up.
In particular, if $U=\emptyset$, then the resulting space $\ZR{X}_U$ is empty, while if $U\neq\emptyset$, we have $\ZR{X}_U\neq\emptyset$, since it contains $U$.
Because of this, we should have a special treatment for what we really want to mean by $\ZR{X}_{\emptyset}$, for which we offer the following definition.

\begin{dfn}\label{dfn-classicalclassicalZR}{\rm 
If $X$ is irreducible, we define
$$
\ZR{k(X)}=\varprojlim_U\ZR{X}_{U},
$$
where $U$ runs over all non-empty quasi-compact open subsets, and call it the (classical) Zariski-Riemann space associated to $X$.}
\end{dfn}

\begin{rem}\label{rem-classicalZRspace}{\rm 
Our classical Zariski-Riemann spaces defined as above generalizes the so-called {\em abstract Riemann manifolds}\index{abstract Riemann manifold} introduced by O.\ Zariski\index{Zariski, O.} (cf.\ \cite{Zar4}) in the case where $S$ is the spectrum of an algebraically closed field $k$ and $X$ is an algebraic variety over $k$.
Note that the first embryonic idea of abstract Riemann manifolds came out from the notion of abstract Riemann surfaces in the works of Dedekind-Weber in the 19th century.}
\end{rem}

\begin{thm}\label{thm-classicalZRsp1}
Let $X$ be a non-empty coherent $S$-scheme, and $U$ a quasi-compact open subset of $X$. $($We allow $U=\emptyset$, if $X$ is irreducible, with the Zariski-Riemann space $\ZR{X}_{\emptyset}$ replaced by $\ZR{k(X)}$ defined as in {\rm \ref{dfn-classicalclassicalZR}}.$)$
\begin{itemize}
\item[{\rm (1)}] {\rm (Zariski \cite{Zar4})} The topological space $\ZR{X}_{U}$ is coherent and sober {\rm ({\bf \ref{ch-pre}}.\ref{dfn-quasicompact1})}
\item[{\rm (2)}] For any object $X'\rightarrow X$ of $\MD_{(X,U)}$ the specialization map $\sp_{X'}$ is quasi-compact and closed.
\end{itemize}
\end{thm}

\begin{proof}
The proof is given by an argument similar to that of \ref{thm-ZRcompact} when $U\neq\emptyset$.
The other case follows from {\bf \ref{ch-pre}}.\ref{thm-projlimcohsch1} (1), {\bf \ref{ch-pre}}.\ref{thm-projlimcohspacepres} (1), (3), and the fact that $\ZR{X}_{U_1}\rightarrow\ZR{X}_{U_2}$ induced from an inclusion $U_1\hookrightarrow U_2$ of quasi-compact open subsets of $X$ is quasi-compact and closed.
\end{proof}

By an argument similar to that in \ref{prop-zariskiriemanntoptop} we have:
\begin{prop}\label{prop-classicalZRsp1}
A subset $\mathfrak{V}$ of $\ZR{X}_{U}$ is a quasi-compact open subset if and only if there exists an object $X'$ of $\MD_{(X,U)}$ and a quasi-compact open subset $V'$ of $X'$ such that $\mathfrak{V}=\sp^{-1}_{X'}(V')$. \hfill$\square$
\end{prop}

\begin{prop}\label{prop-classicalZRsp31}
Let $X$ be a coherent $S$-scheme, and $U\subseteq X$ a quasi-compact open subset.

{\rm (1)} Let $Y\subset X$ be a quasi-compact open subset.
Then we have $\ZR{Y}_{U\cap Y}=\sp^{-1}_X(Y)$, where $\sp_X\colon\ZR{X}_U\rightarrow X$ is the specialization map, and a canonical open immersion 
$$
(\ZR{Y}_{U\cap Y},\O_{\ZR{Y}_{U\cap Y}})\longhookrightarrow(\ZR{X}_U,\O_{\ZR{X}_U})
$$
of locally ringed spaces.

{\rm (2)} If $X'\rightarrow X$ is a proper $U$-admissible morphism, then the induced morphism
$$
(\ZR{X'}_U,\O_{\ZR{X'}_U})\longrightarrow(\ZR{X}_U,\O_{\ZR{X}_U})
$$
is an isomorphism.
\end{prop}

\begin{proof}
(1) follows from \ref{prop-birationalgeom023}, and (2) from \ref{prop-birationalgeom02}.
\end{proof}
\index{Zariski-Riemann space!classical Zariski-Riemann space@classical ---|)}

\subsubsection{Comparison maps}\label{subsub-classicalZRspcompmaps}
Let $X_1$ and $X_2$ be coherent $S$-schemes, and consider the diagram 
$$
X_1\longhookleftarrow U\longhookrightarrow X_2
$$
of quasi-compact $S$-open immersions onto dense open subsets.
By \ref{prop-correspondencediagramfurther} (2) and \ref{prop-birationalgeom023} we have:
\begin{prop}\label{lem-classicalZRsp4}
Suppose that $X_1$ is separated of finite type over $S$, and that $X_2$ is proper over $S$.
Then there exists a canonical open immersion 
$$
(\ZR{X_1}_{U},\O_{\ZR{X_1}_U})\longhookrightarrow(\ZR{X_2}_{U},\O_{\ZR{X_2}_{U}})
$$
of locally ringed spaces. \hfill$\square$
\end{prop}

We call the open immersion thus obtained the {\em comparison map}\index{comparison map} of the models $X_1$ and $X_2$.

\begin{cor}\label{cor-classicalZRsp4}
In the above situation, suppose that both $X_1$ and $X_2$ are proper over $S$.
Then the associated Zariski-Riemann space $\ZR{X_1}_{U}$ and $\ZR{X_2}_{U}$ are canonically isomorphic. \hfill$\square$
\end{cor}

\subsubsection{Relation with rigid Zariskian spaces}\label{subsub-relationRZS}
Let $X$ be a coherent $S$-scheme, and $U\subseteq X$ a quasi-compact dense open subset.
Set $Z=X\setminus U$, which we regard as a closed subscheme of $X$ defined by a quasi-coherent ideal $\mathscr{I}\subset\O_X$ of finite type.
Let $Y=\zat{X}|_Z$ be the Zariskian completion ({\bf \ref{ch-formal}}, \S\ref{subsub-zariskianschemes}) with the ideal of definition $\mathscr{I}\O_Y$ ({\bf \ref{ch-formal}}.\ref{dfn-zariskianschemesideal}), and $\mathscr{X}=Y^{\rig}$ the associated rigid Zariskian space\index{rigid Zariskian space} (\S\ref{subsub-coherentRZScat}).
The open immersion $j\colon U\hookrightarrow X$ lifts to an open immersion $\til{j}\colon U\hookrightarrow\ZR{X}_{U}$ of locally ringed spaces such that the diagram
$$
\xymatrix{&\ZR{X}_{U}\ar[d]^{\sp_X}\\ U\,\ar@{^{(}->}[ur]^{\til{j}}\ar@{^{(}->}@<-.3ex>[r]_j&X}
$$
commutes.
Consider the closed subset $\ZR{X}_{U}\setminus U$ and the inclusion map $\til{i}\colon\ZR{X}_{U}\setminus U\hookrightarrow\ZR{X}_{U}$.

\begin{prop}\label{prop-relclassicalZRZRS}
There exists a canonical isomorphism
$$
(\ZR{X}_{U}\setminus U,\til{i}^{-1}\O_{\ZR{X}_{U}})\stackrel{\sim}{\longrightarrow}(\ZR{\mathscr{X}},\O^{\int}_{\mathscr{X}})
$$
of locally ringed spaces.
\end{prop}

\begin{proof}
Let $\pi\colon X'\rightarrow X$ be a $U$-admissible blow-up, and $\mathscr{J}$ the blow-up center.
The morphism $\pi$ induces a morphism $\zat{\pi}\colon Y'=\zat{X'}|_{\pi^{-1}(Z)}\rightarrow Y=\zat{X}|_Z$ of Zariskian schemes, which is in fact an admissible blow-up along the admissible ideal $\mathscr{J}\O_Y$.
The topological space $\ZR{X}_{U}\setminus U$ is nothing but the projective limit of the underlying topological spaces of the Zariskian schemes of the form $Y'=\zat{X'}|_{\pi^{-1}(Z)}$, where $\pi$ runs through the set of all $U$-admissible blow-ups of $X$.
Hence in order to show that there exists a canonical homeomorphism between the topological spaces $\ZR{X}_{U}\setminus U$ and $\ZR{\mathscr{X}}$, it suffices to show that any admissible blow-up of the Zariskian scheme $Y$ is naturally extended to a $U$-admissible blow-up of $X$.

Let $Y'\rightarrow Y$ be an admissible blow-up of the Zariskian scheme $Y$ along the admissible ideal $\mathscr{J}$.
Take $k\geq 0$ such that $\mathscr{J}\supseteq\mathscr{I}^{k+1}\O_Y$.
Consider the scheme $Y_k=(Y,\O_Y/\mathscr{I}^{k+1}\O_Y)$, which is canonically a closed subscheme of $X$.
Let $i_k\colon Y_k\hookrightarrow X$ be the closed immersion.
Consider the surjective map $\O_X\rightarrow i_{k\ast}\O_{Y_k}$ of sheaves on $X$, and take the pull-back by this map of the ideal $\mathscr{J}_k=\mathscr{J}/\mathscr{I}^{k+1}\O_Y$.
Denote the pull-back ideal by $\til{\mathscr{J}}$.
Then clearly we have $i^{-1}\til{\mathscr{J}}=\mathscr{J}$, where $i\colon Z\hookrightarrow X$ is the closed immersion; moreover, $\til{\mathscr{J}}$ is a $U$-admissible ideal.
Hence if we denote by $X'\rightarrow X$ the $U$-admissible blow-up along $\til{\mathscr{J}}$, then it induces the admissible blow-up $Y'\rightarrow Y$ that we have begun with.
Hence the claim is proved.

Since the sheaf pull-back $i^{-1}$ commutes with filtered inductive limit, we have $\til{i}^{-1}\O_{\ZR{X}_{U}}\cong\O^{\int}_{\mathscr{X}}$, and thus the proposition is shown.
\end{proof}

\subsubsection{Points of the Zariski-Riemann space}\label{subsub-pointsZRsp}
Let $X$ be a coherent $S$-scheme, and $U$ a quasi-compact dense open subset of $X$.
Let $\mathscr{I}\subset\O_X$ be a quasi-coherent ideal of finite type such that $X\setminus U=V(\mathscr{I})$.
Consider the set 
$$
\left\{
\begin{minipage}{19em}
{\small $X$-isomorphism classes of morphisms of the form $\alpha\colon\Spec V\rightarrow X$, where $V$ is an $\mathscr{I}V$-adically separated valuation ring\index{valuation!valuation ring@--- ring}, which maps the generic point to a point of $U$}
\end{minipage}
\right\}.
$$
There exists a map from this set to the Zariski-Riemann space $\ZR{X}_{U}$ constructed as follows.
Let $\alpha\colon\Spec V\rightarrow X$ be given, and set $K=\Frac(V)$.
For any $U$-admissible modification $X'\rightarrow X$, we have a morphism $\Spec K\rightarrow X'$ such that the diagram
$$
\xymatrix{\Spec K\ar[d]\ar[r]&X'\ar[d]\\ \Spec V\ar[r]&X}
$$
commutes, since $X'\rightarrow X$ is isomorphic on $U$.
Since $X'\rightarrow X$ is proper, we have a unique lift $\alpha_{X'}\colon\Spec V\rightarrow X'$.
Thus we have a projective system of morphisms $\{\alpha_{X'}\colon\Spec V\rightarrow X'\}_{X'}$ and hence the morphism of locally ringed spaces $\Spec V\rightarrow\ZR{X}_{U}$.
The desired point of $\ZR{X}_{U}$ is the image of the closed point of $\Spec V$.

We introduce an equivalence relation $\approx$ on the above set, which is generated by the relation $\sim$ defined as follows: for two elements $\alpha\colon\Spec V\rightarrow X$ and $\beta\colon\Spec W\rightarrow X$, $\alpha\sim\beta$ if there exists a local injection $f\colon V\hookrightarrow W$ such that $\alpha\circ\Spec f=\beta$ (cf.\ {\bf \ref{ch-pre}}.\ref{prop-schematic1} (2)).
There exists a map
$$
\left\{
\begin{minipage}{19em}
{\small $X$-isomorphism classes of morphisms of the form $\alpha\colon\Spec V\rightarrow X$, where $V$ is an $\mathscr{I}V$-adically separated valuation ring, which maps the generic point to a point of $U$}
\end{minipage}
\right\}\big/_{{\textstyle \approx}}\longrightarrow\ZR{X}_{U}\leqno{(\ast)}
$$
induced from the mapping defined above.
We topologize the left-hand set by the weakest topology such that the map $(\ast)$ is continuous.
By \ref{thm-classicalZRsp1} (3) and \ref{prop-classicalZRsp1} this topology coincides with the one generated by subsets of the following form: 
$$
\left\{
\begin{minipage}{16em}
{\small $\approx$-equivalence class of $\alpha\colon\Spec V\rightarrow X$ that extends to $\Spec V\rightarrow X'$ such that the image is contained in $Y'$}
\end{minipage}
\right\},
$$
where $X'\rightarrow X$ is a $U$-admissible modification and $Y'$ is an affine open subset of $X'$.

\begin{thm}\label{thm-pointsZRsp1}
The map $(\ast)$ is a homeomorphism.
\end{thm}

The theorem can be verified by an argument similar to that in the proof of \ref{prop-ZRpoints3}, with the aid of the following lemma, which can be shown similarly to \ref{prop-ZRstrsheaf21}:
\begin{lem}\label{lem-pointsZRsp11}
For any point $\xi\in\ZR{X}_{U}$, the ring $\O_{\ZR{X}_{U},\xi}$ is $\mathscr{I}\O_{\ZR{X}_{U},\xi}$-valuative $(${\rm {\bf \ref{ch-pre}}.\ref{dfn-valuative1}}$)$. \hfill$\square$
\end{lem}

Finally, let us include here, as a corollary, the classical version (where $U=\emptyset$) of the theorem, which is well-known in classical birational geometry:
\begin{cor}\label{cor-pointsZRsp1}
Let $X$ be a coherent and integral $S$-scheme of finite type, and $k(X)$ its function field.
Then the topological space $\ZR{X}$ is identified with the following space: as a set, it is the set of all valulation rings for $k(X)$ dominating the local ring $\O_{X,x}$ of a point $x$ of $X;$ the topology is generated by subsets of the following form$:$
$$
\left\{
\begin{minipage}{10em}
{\small valuation ring for $k(X)$ that contains $A$}
\end{minipage}
\right\},
$$
where $A$ varies among subrings of $k(X)$ arising from a dominant and finite type morphism $\Spec A\rightarrow X$ over $S$.
\end{cor}

\begin{proof}
The Zariski-Riemann space $\ZR{X}$ is the projective limit of all Zariski-Riemann spaces of the form $\ZR{X}_U$, where $U$ are non-empty quasi-compact open subsets of $X$.
Hence it is the projective limit of the spaces in the left-hand side of $(\ast)$ as above.
Suppose $\alpha\colon\Spec V\rightarrow X$ belongs to this set.
Then the image of the generic point of $V$ lies in any quasi-compact open subset of $X$ and hence is the generic point of $X$.
Hence the fractional field $K=\Frac(V)$ of $V$ contains the function field $k(X)$.
Set $W=V\cap k(X)$.
Then $W$ is a valuation ring for $k(X)$, and the map $W\hookrightarrow V$ is local.
Moreover, we have the morphism $\beta\colon\Spec W\rightarrow X$ that is $\approx$-equivalent to $\alpha$.
Conversely, if we are given a valuation ring $W$ for $k(X)$ dominating a point of $X$, then the induced morphism $\beta\colon\Spec W\rightarrow X$ lies in $\ZR{X}_U$ for all $U$.
\end{proof}

\section{Appendix: Nagata's embedding theorem}\label{sec-nagataembedding}
\index{Nagata's embedding theorem|(}
\subsection{Announcement of the theorem}\label{sub-nagataembedding}
In this section we give a complete proof of the following theorem, which generalizes the famous theorem due to M.\ Nagata\index{Nagata, M.} (cf.\ \cite{Nagata2}):
\begin{thm}[Nagata's embedding theorem]\label{thm-nagataembedding}
Let $Y$ be a coherent scheme, and $f\colon X\rightarrow Y$ a separated morphism of finite type between schemes. 
Then there exists a proper $Y$-scheme $\ovl{f}\colon\ovl{X}\rightarrow Y$ that admits a dense open immersion $X\hookrightarrow\ovl{X}$ over $Y$.
\end{thm}

\subsection{Preparation for the proof}\label{sub-nagataembeddingpre}
\subsubsection{Canonical compactification}\label{subsub-canonicalcompnagata}
Let $f\colon X\rightarrow Y$ a separated morphism of finite type between algebraic spaces, where $Y$ is a coherent scheme.
Let us say that $f$ (or $X$) is {\em compactifiable}\index{compactifiable} if there exists a commutative diagram
$$
\xymatrix@C-2ex@R-1ex{X\,\ar@{^{(}->}[rr]^j\ar[dr]_f&&\til{X}\ar[dl]^{\til{f}}\\ &Y}
$$
such that $\til{f}$ is proper and $j$ is an open immersion.
If this is the case, since the scheme-theoretic closure of $X$ in $\til{X}$ is proper over $Y$, one can take $\til{f}$ such that $j$ is a dense open immersion.
In this situation, $\til{X}$ is called a {\em compactification} of $X$ over $Y$.
It is clear that, if $f$ is affine, then $f$ is compactifiable by a proper $Y$-scheme.
If $X$ is a scheme, $X$ is covered by open subschemes affine over $Y$, and hence is {\em locally compactifiable}.

Suppose that $f\colon X\rightarrow Y$, where $X$ is a scheme, is compactifiable by a proper $Y$-scheme $\ovl{f}\colon\ovl{X}\rightarrow Y$; here, we assume without loss of generality that $\ovl{X}$ contains $X$ as a quasi-compact dense open subset.
Then consider the Zariski-Riemann space\index{Zariski-Riemann space!classical Zariski-Riemann space@classical ---} $\ZR{\ovl{X}}_{X}$, which is the projective limit of all $X$-admissible blow-ups of the scheme $\ovl{X}$.
Let $f'\colon\ovl{X}'\rightarrow Y$ be another compactification of $f$ by a proper $Y$-scheme.
Then by \ref{cor-classicalZRsp4} the topological spaces $\ZR{\ovl{X}}_{X}$ and $\ZR{\ovl{X}'}_{X}$ are canonically isomorphic to each other, that is, the space $\ZR{\ovl{X}}_{X}$ is independent up to canonical isomorphisms on the choice of the compactification.

\begin{dfn}[Canonical compactification]\label{dfn-canonicalcompnagata}{\rm 
Let $f\colon X\rightarrow Y$ be compactifiable by a proper $Y$-scheme $\ovl{f}\colon\ovl{X}\rightarrow Y$.
We write
$$
\ZR{X}_{\cpt}=\ZR{\ovl{X}}_{X},
$$
and call it the {\em canonical compactification}\index{compactification!canonical compactification@canonical ---} of $X$ over $Y$.}
\end{dfn}

Notice that the topological space $\ZR{X}_{\cpt}$ is coherent and sober (\ref{thm-classicalZRsp1} (1)) and has $X$ as a quasi-compact open subspace.
By \ref{prop-relclassicalZRZRS} we have
$$
\ZR{X}_{\cpt}\setminus X\stackrel{\sim}{\longrightarrow}\ZR{(\zat{\ovl{X}}|_Z)^{\rig}},
$$
where $Z$ is a finitely presented closed subscheme of $\ovl{X}$ with the underlying topological space $\ovl{X}\setminus X$.
The topological space $\ZR{X}_{\cpt}$ has the structure sheaf 
$$
\O_{\ZR{X}_{\cpt}}=\O_{\ZR{\ovl{X}}_{X}}.
$$

Note also that, if $Z\hookrightarrow X$ is a closed immersion of separated and finite type $Y$-schemes, and if $X$ is compactifiable by a proper $Y$-scheme, then $Z$ is compactifiable by a proper $Y$-scheme, and $\ZR{Z}_{\cpt}$ coincides with the closure of $Z$ in $\ZR{X}_{\cpt}$ with respect to the canonical inclusion $Z\hookrightarrow\ZR{X}_{\cpt}$.

\begin{prop}\label{prop-canonicalcompnagata}
{\rm (1)} Let $X\hookrightarrow X'$ be a $Y$-open immersion of separated $Y$-schemes of finite type.
Suppose that $X'$ is compactifiable by a proper $Y$-scheme.
Then $X$ is compactifiable by a proper $Y$-scheme, and we have the canonical closed map
$$
(\ZR{X}_{\cpt},\O_{\ZR{X}_{\cpt}})\longrightarrow(\ZR{X'}_{\cpt},\O_{\ZR{X'}_{\cpt}})
$$
of locally ringed spaces.

{\rm (2)} Let $f\colon U\rightarrow Y$ be a separated $Y$-scheme of finite type that is compactifiable by a proper $Y$-scheme, and $j\colon U\hookrightarrow X$ an open immersion into a separated $Y$-scheme of finite type.
Then we have the canonical open immersion
$$
(\ZR{X}_U,\O_{\ZR{X}_U})\longhookrightarrow(\ZR{U}_{\cpt},\O_{\ZR{U}_{\cpt}})
$$
of locally ringed spaces.
\end{prop}

\begin{proof}
Let $\ovl{X}'$ be a compactification of $X'$ over $Y$.
Since the closure $\ovl{X}$ of $X$ in $\ovl{X}'$ gives a compactification of $X$ over $Y$, $X$ is compactifiable by a proper $Y$-scheme.
Since any $X'$-admissible blow-up $\til{X}'\rightarrow\ovl{X}'$ induces an $X$-admissible blow-up $\til{X}\rightarrow\ovl{X}$ by the strict transform, we have the desired map by passage to the projective limits.
By {\bf \ref{ch-pre}}.\ref{thm-projlimcohspacepres} (3) this map is closed, whence (1).
(2) follows from \ref{lem-classicalZRsp4}.
\end{proof}

Notice that the quasi-compact open subset $\ZR{X}_U$ in $\ZR{U}_{\cpt}$ depends only on the closure of $U$ in $X$.

\subsubsection{General construction}\label{subsub-canonicalcompgen}
Let $f\colon X\rightarrow Y$ be a separated morphism of finite type between coherent schemes.
We are to construct the space $\ZR{X}_{\cpt}$, again called the {\it canonical compactification}\index{compactification!canonical compactification@canonical ---}, without assuming that $X$ is compactifiable\index{compactifiable}.

Let $U$ be a quasi-compact open subset of $X$ that is compactifiable by a proper $Y$-scheme (e.g., affine over $Y$).
By \ref{prop-canonicalcompnagata} (2) the space $\ZR{X}_U$ is regarded as a quasi-compact open subset of $\ZR{U}_{\cpt}$.

\begin{dfn}[Partial compactification]\label{dfn-partialcompnagata}{\rm 
Set
$$
\ZR{U}^X_{\pc}=\ZR{U}_{\cpt}\setminus\ovl{\ZR{X}_U\setminus U},
$$
where the closure in the right-hand side is taken in $\ZR{U}_{\cpt}$, and call it the {\it partial compactification of $U$ relative to $X$}\index{compactification!partial compactification@partial ---}.}
\end{dfn}

\begin{prop}\label{prop-partialcompnagata0}
Let $X$ be a separated of finite type $Y$-scheme, and $U\subseteq X$ a quasi-compact open subset, which is assumed to be compactifiable by a proper $Y$-scheme.

{\rm (1)} For any proper $U$-admissible morphism $X'\rightarrow X$, we have $\ZR{U}^X_{\pc}=\ZR{U}^{X'}_{\pc}$.

{\rm (2)} There exists a proper $U$-admissible morphism $X'\rightarrow X$ such that $X'$ is compactifiable by a proper $Y$-scheme.
\end{prop}

The proposition says that, whenever discussing the partial compactification $\ZR{U}^X_{\pc}$, one may assume $X$ is compactifiable without loss of generality. 

\begin{proof}
(1) The assertion is clear if $X'\rightarrow X$ is a $U$-admissible blow-up, since, then, the images of $\ZR{X}_U$ and $\ZR{X'}_U$ in $\ZR{U}_{\cpt}$ coincide with each other.
The general case reduces to this situation due to \ref{prop-birationalgeom02}.

(2) Take a proper $Y$-scheme $\ovl{U}$ that contains $U$ as an open subset.
Consider the join $X\ast^U\ovl{U}$, and take a $U$-admissible blow-up $\ovl{U}'\rightarrow\ovl{U}$ such that the strict transform, denoted by $X'\rightarrow\ovl{U}'$, of $X\ast^U\ovl{U}\rightarrow\ovl{U}$ is an open immersion, as in \ref{prop-correspondencediagramfurther} (1).
Then $X'$ is compactifiable by a proper $Y$-scheme (\ref{prop-canonicalcompnagata} (1)), and $X'\rightarrow X$ is proper $U$-admissible (\ref{prop-correspondencediagramfurther} (1) (a)).
\end{proof}

\begin{lem}\label{lem-partialcompnagata00}
Let $X$ be a separated of finite type $Y$-scheme, and $U\subseteq X$ a quasi-compact open subset.
We assume that $X$ is compactifiable by a proper $Y$-scheme.

{\rm (1)} The partial compactification $\ZR{U}^X_{\pc}$ coincides with the maximal open subset of $\ZR{U}_{\cpt}$ on which the canonical map $\ZR{U}_{\cpt}\rightarrow\ZR{X}_{\cpt}$ is an isomorphism.

{\rm (2)} Set $Z=X\setminus U$. Then we have
$$
\ZR{U}^X_{\pc}=\ZR{X}_{\cpt}\setminus\ovl{Z},
$$
where $\ovl{Z}$ denotes the closure of $Z$ in $\ZR{X}_{\cpt}$ with respect to the canonical inclusion $Z\hookrightarrow\ZR{X}_{\cpt}$.
\end{lem}

\begin{proof}
(1) Take a compactification $\ovl{X}$ of $X$ and the closure $\ovl{U}$ of $U$ in $\ovl{X}$.
We have $\ZR{X}_{\cpt}=\ZR{\ovl{X}}_X$ and $\ZR{U}_{\cpt}=\ZR{\ovl{U}}_U=\ZR{\ovl{X}}_U$.
For any $U$-admissible blow-up $\ovl{X}_1\rightarrow\ovl{X}$, we have $\sp_{\ovl{X}_1}(\ZR{X}_U\setminus U)=X_1\setminus U$, where $X_1=X\times_{\ovl{X}}\ovl{X}_1$, and 
$$
\ovl{\ZR{X}_U\setminus U}=\bigcap_{\ovl{X}_1\rightarrow\ovl{X}}\sp^{-1}_{\ovl{X}_1}(\ovl{X_1\setminus U}), \eqno{(\ast)}
$$
where $\ovl{X}_1$ of the left-hand side runs through all $U$-admissible blow-ups of $\ovl{X}$; see {\bf \ref{ch-pre}}.\ref{lem-projlimclosedmapslem} (2).

Let $W$ be a quasi-compact open subset of $\ZR{U}_{\cpt}$, and take a $U$-admissible blow-up $\ovl{X}_1\rightarrow\ovl{X}$ and a quasi-compact open subset $W_1\subseteq\ovl{X}_1$ such that $W=\sp^{-1}_{\ovl{X}_1}(W_1)$; see \ref{prop-classicalZRsp1}.
In view of \ref{prop-classicalZRsp31} (1), the map $\ZR{U}_{\cpt}\rightarrow\ZR{X}_{\cpt}$ is an isomorphism on $W$ if and only if 
$$
\ZR{W_1}_{U\cap W_1}=\ZR{W_1}_{X_1\cap W_1},
$$
which is, furthermore, equivalent to that any $X_1\cap W_1$-admissible blow-up of $W_1$ is dominated by a $U\cap W_1$-admissible blow-up.
This is possible if and only if $U\cap W_1=X_1\cap W_1$, or equivalently, $W_1\cap(X_1\setminus U)=\emptyset$.

To sum up, we have shown the following: $\ZR{U}_{\cpt}\rightarrow\ZR{X}_{\cpt}$ is an isomorphism on $W$ if and only if $W_1\cap(X_1\setminus U)=\emptyset$ for any $U$-admissible blow-up $\ovl{X}_1\rightarrow\ovl{X}$ and any quasi-compact open subset $W_1\subseteq\ovl{X}_1$ such that $W=\sp^{-1}_{\ovl{X}_1}(W_1)$.
By the equality $(\ast)$, the last condition is equivalent to that $W\subseteq\ZR{U}^X_{\pc}$.

(2) Restricting to $X$-admissible blow-ups, one can similarly show the following: a quasi-compact open $W\subseteq\ZR{X}_{\cpt}=\ZR{\ovl{X}}_X$ lies in the image of $\ZR{U}^X_{\pc}$ if and only if, for any $X$-admissible blow-up $\ovl{X}_1\rightarrow\ovl{X}$ and a quasi-compact open subset $W_1\subseteq\ovl{X}_1$ such that $W=\sp_{\ovl{X}_1}(W_1)$, we have $W_1\cap (X\setminus U)=\emptyset$, which is equivalent to that $W\cap Z=\emptyset$.
\end{proof}

Let $X$ be a separated of finite type $Y$-scheme, which we assume to be compactifiable by a proper $Y$-scheme $\ovl{X}$.
Since the open immersion $X\hookrightarrow\ovl{X}$ is quasi-compact, one can take a quasi-coherent ideal $\mathscr{I}\subseteq\O_{\ovl{X}}$ of finite type that gives a closed subscheme with the underlying topological space $\ovl{X}\setminus X$.
Then, for any $x\in\ZR{X}_{\cpt}=\ZR{\ovl{X}}_X$, the local ring $\O_{\ZR{X}_{\cpt},x}$ is $\mathscr{I}_x\O_{\ZR{X}_{\cpt},x}$-valuative; see \ref{lem-pointsZRsp11}.
Set 
$$
J_x=\bigcap_{n\geq 1}\mathscr{I}^n_x,
$$
which is a prime ideal of $\O_{\ZR{X}_{\cpt},x}$, giving the associated valuation ring $V_x=\O_{\ZR{X}_{\cpt},x}/J_x$; see {\bf \ref{ch-pre}}, \S\ref{subsub-valuativecomp}.
By the morphism as in \S\ref{subsub-pointsZRsp}, the generic point of $\Spec V_x$ defines a generization $\til{x}\in\ZR{X}_{\cpt}$ of $x$.

\begin{lem}\label{lem-partialcompnagata000}
In the above situation, let $U\subseteq X$ be a quasi-compact open subset.
Then we have 
$$
\ZR{U}^X_{\pc}=\{x\in\ZR{X}_{\cpt}:\til{x}\in U\}.
$$
\end{lem}

\begin{proof}
Let $x\in\ZR{X}_{\cpt}$.
For any $X$-admissible blow-up $\ovl{X}_1\rightarrow\ovl{X}$, we have $\sp_{\ovl{X}_1}(\til{x})\in X$.
Due to \ref{lem-partialcompnagata00} (2), $x\in\ZR{U}^X_{\pc}$ if and only if $\sp_{\ovl{X}_1}(\til{x})\not\in Z=X\setminus U$, that is, $\sp_{\ovl{X}_1}(\til{x})\in U$, whence the assertion.
\end{proof}

\begin{prop}\label{prop-partialcompnagata}
Let $X$ be a separated of finite type $Y$-scheme, and $U_1,U_2$ quasi-compact open subsets of $X$.

{\rm (1)} Suppose that $U_2$ is compactifiable by a proper $Y$-scheme, and that $U_1\subseteq U_2$.
Then there exists a canonical open immersion $\ZR{U_1}^X_{\pc}\hookrightarrow\ZR{U_2}^X_{\pc}$ that extends the inclusion map $U_1\hookrightarrow U_2$.

{\rm (2)} Suppose $X$ is compactifiable by a proper $Y$-scheme. 
Then we have the following equalities in $\ZR{X}_{\cpt}$:
$$
\ZR{U_1}^X_{\pc}\cap\ZR{U_2}^X_{\pc}=\ZR{U_1\cap U_2}^X_{\pc},\quad\ZR{U_1}^X_{\pc}\cup\ZR{U_2}^X_{\pc}=\ZR{U_1\cup U_2}^X_{\pc}.
$$
\end{prop}

\begin{proof}
(1) By \ref{prop-partialcompnagata0}, we may assume without loss of generality that $X$ is compactifiable by a proper $Y$-scheme.
Then, by \ref{lem-partialcompnagata000}, we have the inclusion $\ZR{U_1}^X_{\pc}\hookrightarrow\ZR{U_2}^X_{\pc}$ in $\ZR{X}_{\cpt}$.

(2) follows immediately from \ref{lem-partialcompnagata000}.
\end{proof}

Let $f\colon X\rightarrow Y$ be a separated and of finite type $Y$-scheme, and $\{U_{\alpha}\}_{\alpha\in L}$ a finite open covering of $X$ such that each $U_{\alpha}$ is compactifiable by a proper $Y$-scheme.
Define the topological space $\ZR{X}_{\cpt}$ by the following cokernel diagram in the category $\Top$ of topological spaces
$$
\xymatrix{\ZR{X}_{\cpt}&\coprod_{\alpha\in L}\ZR{U_{\alpha}}^X_{\pc}\ar[l]&\coprod_{\alpha,\beta\in L}\ZR{U_{\alpha}\cap U_{\beta}}^X_{\pc},\ar@<-.5ex>[l]\ar@<.5ex>[l]}
$$
that is, by patching the partial compactifications $\ZR{U_{\alpha}}^X_{\pc}$ along the open subspaces $\ZR{U_{\alpha}\cap U_{\beta}}^X_{\pc}$.
The space $\ZR{X}_{\cpt}$ has the structure sheaf $\O_{\ZR{X}_{\cpt}}$ by patching the structure sheaves $\O_{\ZR{U_{\alpha}}_{\cpt}}|_{\ZR{U_{\alpha}}^X_{\pc}}$.

\begin{prop}\label{prop-lempartialcompnagata1}
{\rm (1)} For any quasi-compact open subset $U$ of $X$ that is compactifiable by a proper $Y$-scheme, $\ZR{U}^X_{\pc}$ is canonically an open subspace of $\ZR{X}_{\cpt}$.
In particular, the formation of $\ZR{X}_{\cpt}$ does not depend on the choice of the finite open covering $\{U_{\alpha}\}_{\alpha\in L}$.

{\rm (2)} If $X$ is compactifiable by a proper $Y$-scheme, then the above-defined $\ZR{X}_{\cpt}$ coincides with the one defined in {\rm \ref{dfn-canonicalcompnagata}}.
\end{prop}

\begin{proof}
First notice that (2) follows immediately from \ref{prop-partialcompnagata}.
As for (1), since $U$ and all $U_{\alpha}$ are compactifiable, we have $\ZR{U}^X_{\pc}=\ZR{\bigcup_{\alpha\in L}U\cap U_{\alpha}}^X_{\pc}=\bigcup_{\alpha\in L}\ZR{U\cap U_{\alpha}}^X_{\pc}$ by \ref{prop-partialcompnagata} (2), and since each $\ZR{U\cap U_{\alpha}}^X_{\pc}$ is canonically an open subspace of $\ZR{U_{\alpha}}^X_{\pc}$ by \ref{prop-partialcompnagata} (1), $\ZR{U}^X_{\pc}$ is an open subspace of $\ZR{X}_{\cpt}$.
To show that the formation of $\ZR{X}_{\cpt}$ does not depend on the choice of the finite open covering $\{U_{\alpha}\}_{\alpha\in L}$, one employs the inductive argument based on the following statement, which follows from what we have proven just now: for any quasi-compact open subset $U\subset X$, the space $\ZR{X}_{\cpt}$ constructed from the covering $\{U_{\alpha}\}_{\alpha\in L}$ is canonically isomorphic to the one from the covering $\{U_{\alpha}\}_{\alpha\in L}\cup\{U\}$.
\end{proof}

\begin{dfn}\label{dfn-partialcompnagata1}{\rm 
The topological space 
$$
\ZR{X}_{\cpt}
$$
thus defined is called the {\em canonical compactification}\index{compactification!canonical compactification@canonical ---} of $X$ over $Y$.}
\end{dfn}

Notice that, by the construction, $\ZR{X}_{\cpt}$ contains $X$ as an open subset and that there exists the canonical morphism $\ZR{X}_{\cpt}\rightarrow Y$ of locally ringed spaces such that the diagram
$$
\xymatrix@C-2ex@R-1ex{X\,\ar@{^{(}->}[rr]^(.44)j\ar[dr]_f&&\ZR{X}_{\cpt}\ar[dl]\\ &Y}
$$
commutes.

Let $U$ be a quasi-compact open subset of $X$.
For a finite open covering $\{U_{\alpha}\}_{\alpha\in L}$ of $X$ as above, the collection of spaces $\{\ZR{U\cap U_{\alpha}}^X_{\pc}\}_{\alpha\in L}$ defines an open subset of $\ZR{X}_{\cpt}$, which we consistently denote by
$$
\ZR{U}^X_{\pc}.
$$
Notice that the space $\ZR{U}^X_{\pc}$ is, at the same time, an open subspace of $\ZR{U}_{\cpt}$.

If $X\hookrightarrow X'$ is an open immersion of separated $Y$-schemes of finite type, then one sees by an obvious patching argument that $\ZR{X'}_X$ is regarded as an open subset of $\ZR{X}_{\cpt}$.

\subsubsection{Properties of canonical compactification}\label{subsub-propcancompnagata}
\begin{prop}\label{prop-propcancompnagata1}
Any $Y$-open immersion $X\hookrightarrow X'$ of separated $Y$-schemes of finite type canonically induces a closed map
$$
(\ZR{X}_{\cpt},\O_{\ZR{X}_{\cpt}})\longrightarrow(\ZR{X'}_{\cpt},\O_{\ZR{X'}_{\cpt}})
$$
of locally ringed spaces.
Moreover, the open subspace $\ZR{X}^{X'}_{\pc}$ is the maximal one among the open subspaces of $\ZR{X}_{\cpt}$ on which the restrictions of the map $\ZR{X}_{\cpt}\rightarrow\ZR{X'}_{\cpt}$ are open immersions.
\end{prop}

\begin{proof}
The first assertion follows from \ref{prop-canonicalcompnagata} (1) by patching.
The last assertion follows from the proof of \ref{prop-partialcompnagata} (2).
\end{proof}

\begin{cor}\label{cor-propcancompnagata1}
The topological space $\ZR{X}_{\cpt}$ is quasi-compact.
\end{cor}

\begin{proof}
Take a finite open covering $\{U_{\alpha}\}_{\alpha\in L}$ of $X$ such that each member $U_i$ is compactifiable by a proper $Y$-scheme.
Since $\ZR{X}_{\cpt}$ is covered by the open subsets $\ZR{U_{\alpha}}^X_{\pc}$, the map
$$
{\textstyle \coprod_{\alpha\in L}\ZR{U_{\alpha}}_{\cpt}\longrightarrow\ZR{X}_{\cpt}}
$$
by \ref{prop-propcancompnagata1} is surjective.
Since each $\ZR{U_{\alpha}}_{\cpt}$ is quasi-compact, the result follows.
\end{proof}

\begin{prop}[Valuative criterion]\label{prop-propcancompval}
Let $V$ be a valuation ring, and set $K=\Frac(V)$.
Suppose we are given a commutative diagram
$$
\xymatrix{\Spec K\ar[r]^(.6){\beta}\ar[d]&X\ar[d]^f\\ \Spec V\ar[r]_(.6){\alpha}&Y.}
$$
Then there exists uniquely a morphism $\Spec V\rightarrow\ZR{X}_{\cpt}$ of locally ringed spaces such that the resulting diagram
$$
\xymatrix@-1ex{&\ZR{X}_{\cpt}\ar[d]\\ \Spec V\ar[r]_(.6){\alpha}\ar[ur]&Y}
$$
commutes.
\end{prop}

\begin{proof}
The assertion is clear if $X$ is compactifiable.
In general, take a quasi-compact open subset $U$ of $X$ that contains the image of $\beta$.
We may assume that $U$ is compactifiable by a proper $Y$-scheme.
Then we have $\gamma\colon\Spec V\rightarrow\ZR{U}_{\cpt}$ that lifts $\alpha$.
On the other hand, since the open subset $\ZR{X}_U$ of $\ZR{U}_{\cpt}$ is quasi-compact, $\gamma^{-1}(\ZR{X}_U)$ has the minimal point $\mathfrak{p}$ ({\bf \ref{ch-pre}}.\ref{prop-lemvalaval}).
We may then replace $U$ by a compactifiable quasi-compact open subset of $X$ that contains the image $x$ of $\gamma(\mathfrak{p})$ by the specialization map $\sp_X\colon\ZR{X}_U\rightarrow X$.
Then $\ovl{\{\mathfrak{p}\}}=\Spec V/\mathfrak{p}\subseteq\Spec V$ is mapped to $\ZR{U}_{\cpt}$ by the map $\gamma$.
It is easy to see by the construction that $\ovl{\{\mathfrak{p}\}}$ is actually mapped into $\ZR{U}^X_{\pc}$.
Then by patching the valuation rings $\Spec V_{\mathfrak{p}}=\gamma^{-1}(\ZR{X}_U)\rightarrow X$ and $\Spec V/\mathfrak{p}=\ovl{\{\mathfrak{p}\}}\rightarrow\ZR{U}^X_{\pc}\subseteq\ZR{X}_{\cpt}$ (cf.\ {\bf \ref{ch-pre}}, \S\ref{sub-composition}), one has the desired lifting.
\end{proof}

By this and the usual valuative criterian of properness we have:
\begin{cor}\label{cor-propcancompval}
Let $X\hookrightarrow X'$ be a dense $Y$-open immersion of separated $Y$-schemes of finite type.
If $\ZR{X'}_X=\ZR{X}_{\cpt}$, then $X'$ is proper over $Y$.\hfill$\square$
\end{cor}

\subsection{Proof of Theorem \ref{thm-nagataembedding}}\label{sub-nagataembeddingproof}
\subsubsection{Lemmas}\label{subsub-nagataembeddinglemmas}
\begin{lem}[Intersection lemma]\label{lem-nagataembedding2}
Let $X\hookrightarrow X_1$ and $X\hookrightarrow X_2$ be dense $Y$-open immersions between separated $Y$-schemes of finite type, and consider the join $X_1\ast^XX_2$ $($cf.\ {\rm \ref{dfn-correspondencediagram}}$)$.
Then we have 
$$
\ZR{X_1\ast^XX_2}_X=\ZR{X_1}_X\cap\ZR{X_2}_X
$$
in $\ZR{X}_{\cpt}$.
\end{lem}

\begin{proof}
It follows from \ref{prop-correspondencediagramfurther} (2) and \ref{prop-birationalgeom023} that there exists projective systems $\{X_{i,\lambda}\}_{\lambda\in\Lambda}$ ($i=1,2$) and $\{Z_{\lambda}\}_{\lambda\in\Lambda}$ indexed by a directed set $\Lambda$ and dense open immersions 
$$
X_{1,\lambda}\longhookleftarrow Z_{\lambda}\longhookrightarrow X_{2,\lambda}
$$
that are compatible with the projection maps such that $\ZR{X_i}_X=\varprojlim_{\lambda\in\Lambda}X_{i,\lambda}$ for $i=1,2$ and $\ZR{X_1\ast^XX_2}_X=\varprojlim_{\lambda\in\Lambda}Z_{\lambda}$.
Then the desired equality follows from the left-exactness of projective limits.
\end{proof}

\begin{lem}[Patching lemma]\label{lem-nagataembedding4}
Let $X\hookrightarrow X_1$ and $X\hookrightarrow X_2$ be dense $Y$-open immersions between separated $Y$-schemes of finite type.
Then there exists another dense $Y$-open immersion $X\hookrightarrow Z$ of separated $Y$-schemes of finite type such that 
$$
\ZR{Z}_X=\ZR{X_1}_X\cup\ZR{X_2}_X
$$
in $\ZR{X}_{\cpt}$.
\end{lem}

\begin{proof}
By \ref{prop-correspondencediagramfurther} (2) we have a diagram 
$$
\til{X}_1\stackrel{\til{p}_1}{\longhookleftarrow}\til{W}\stackrel{\til{p}_2}{\longhookrightarrow}\til{X}_2
$$
consisting of $X$-admissible quasi-compact $Y$-open immersions, where $\til{X}_1$ and $\til{X}_2$ are admissible blow-ups of $X_1$ and $X_2$, respectively.
The desired $Z$ should be the gluing of $\til{X}_1$ and $\til{X}_2$ along $\til{W}$ in the usual sense.
Note that the scheme $\til{W}$ is proper over the join $W=X_1\ast^XX_2$.

To conclude, we need to show that the model $Z$ thus obtained 
is a separated $S$-scheme; the key-point is that the intersection $\langle X_1\rangle_X\cap\langle X_2\rangle_X$ is represented by the join $X_1\ast^XX_2$, as we have seen in \ref{lem-nagataembedding2}.
To show the claim, we use the valuative criterion. Let $V$ be a 
valuation ring with $K$ the field of fractions. Suppose we have two 
morphisms $\alpha_1,\alpha_2\colon\Spec V\rightarrow Z$ 
dominating the same maxmial
point such that $\alpha_1\otimes_VK=\alpha_2\otimes_VK$.
We need to show that $\alpha_1=\alpha_2$.
When the closed point of $\Spec V$ is mapped by both $\alpha_1$ and 
$\alpha_2$ into either one of $\til{X}_1,\til{X}_2$, the claim 
immediately follows, since $X_1,X_2$ are assumed to be separated.
If, to the contrary, say, $\alpha_1(\Spec V)\subset\til{X}_1$ and 
$\alpha_2(\Spec V)\subset\til{X}_2$, since $\alpha_1$ and $\alpha_2$ 
coincides on $\Spec K$, we have a morphism
$$
(\alpha_1,\alpha_2)\colon\Spec V\longrightarrow X_1\ast^XX_2,
$$
and hence $\beta\colon\Spec V\rightarrow\til{W}$, since $\til{W}$ is 
proper over $W=X_1\ast^XX_2$.
Now by the construction the compositions
$$
\Spec 
V\stackrel{\beta}{\longrightarrow}\til{W}\longhookrightarrow\til{X}_i\longhookrightarrow\til{X}
$$
($i=1,2$) are nothing but $\alpha_i$, thereby the claim.
\end{proof}

\begin{lem}[Local extension lemma]\label{lem-nagataembedding5}
For any point $y\in\ZR{X}_{\cpt}$ there exists a dense $Y$-open immersion $X\hookrightarrow X_y$ of separated $Y$-schemes of finite type such that $\ZR{X_y}_X$ contains the point $y$.
\end{lem}

\begin{proof}
Take a quasi-compact open subset $U$ of $X$ that is compactifiable by a proper $Y$-scheme $\ovl{U}$ such that $y$ lies in $\ZR{U}^X_{\pc}$.
By the construction of $\ZR{U}^X_{\pc}$ we have $y\not\in\ovl{\ZR{X}_U\setminus U}$.
It follows therefore that there exists a quasi-compact open neighborhood $\mathfrak{V}$ of $y$ in $\ZR{U}_{\cpt}$ that is disjoint from $\ZR{X}_U\setminus U$; note that $\mathfrak{V}\subseteq\ZR{X}_{\cpt}$.
Replacing the compactification $\ovl{U}$ of $U$ by a $U$-admissible blow-up if necessary, we may assume that there exists a quasi-compact open subset $V$ of $\ovl{U}$ such that $\sp^{-1}_{\ovl{U}}(V)=\mathfrak{V}$ ({\bf \ref{ch-pre}}.\ref{prop-projlimcohtopspqcptopen}).
Replacing $V$ by $V\cup U$ if necessary, we may assume that $V$ contains $U$.
By the construction we have $\ZR{X}_U\cap\ZR{V}_U=U$.
Patching $X$ and $V$ along $U$ birationally, we get a separated $Y$-scheme of finite type $X_y$ as in \ref{lem-nagataembedding4} such that $\ZR{X_y}_U=\ZR{X}_U\cup\mathfrak{V}$.
(Here we used the fact that the specialization map $\sp_{\ovl{U}}$ is surjective ({\bf \ref{ch-pre}}.\ref{thm-projlimcohspacepres} (2)) and hence $\ZR{V}_U=\mathfrak{V}$.)
Since $\ZR{X_y}_X\cap\ZR{U}^X_{\pc}=\ZR{X_y}_U\cap\ZR{U}^X_{\pc}$, $\ZR{X_y}_X$ contains $y$, as desired.
\end{proof}

\subsubsection{Proof of the theorem}\label{subsub-nagataembthmproof}
Now we proceed to the proof of \ref{thm-nagataembedding}.
Let $f\colon X\rightarrow Y$ be a separated $Y$-scheme of finite type, where $Y$ is a coherent scheme.
First, for any $y\in\ZR{X}_{\cpt}$ we take $X_y$ as in \ref{lem-nagataembedding5} containing $X$ such that $y\in\ZR{X_y}_X$.
Now the quasi-compactness of $\ZR{X}_{\cpt}$ (\ref{cor-propcancompnagata1}) implies that there 
exists a {\em finite} set $I$ of $\ZR{X}_{\cpt}$ such that $\{\ZR{X_y}_X\}_{y\in I}$ gives a covering of $\ZR{X}_{\cpt}$.
Then, applying \ref{lem-nagataembedding4} successively, we get a separated $Y$-scheme $\ovl{X}$ of finite type containing $X$ such that
$$
{\textstyle \ZR{\ovl{X}}_X=\bigcup_{y\in I}\ZR{X_y}_X=\ZR{X}_{\cpt}}.
$$
By \ref{cor-propcancompval} we conclude that $\ovl{X}$ is proper over $Y$, which therefore gives the desired compactification. \hfill$\square$
\index{Nagata's embedding theorem|)}

\subsection{Application: Removing Noetherian hypothesis}\label{sub-appclassicalZR}
\begin{prop}\label{prop-appclassicalZR}
Let $Y$ be a quasi-compact scheme, and $f\colon X\rightarrow Y$ a morphism.
If $f$ is proper and affine, then it is finite.
\end{prop}

To show the proposition, we need several preparatory results:
\begin{lem}\label{lem-appclassicalZR}
Let $X$ be a coherent scheme, and $U$ a quasi-compact open subset of $X$.
Let $\{Z_{\lambda}\}_{\lambda\in\Lambda}$ be a filtered projective system consisting of closed subscheme of $X$ indexed by a directed set $\Lambda$ such that for each $\lambda\geq\mu$ the transition maps $i_{\lambda\mu}\colon Z_{\lambda}\rightarrow Z_{\mu}$ are the closed immersions over $X$.
Then we have the equalities
$$
\bigcap_{\lambda\in\Lambda}\ovl{U\cap Z_{\lambda}}=\ovl{\bigcap_{\lambda\in\Lambda}(U\cap Z_{\lambda})}=\ovl{U}\cap(\bigcap_{\lambda\in\Lambda}Z_{\lambda}).
$$
\end{lem}

\begin{proof}
The inclusion $\bigcap_{\lambda\in\Lambda}\ovl{U\cap Z_{\lambda}}\supseteq\ovl{\bigcap_{\lambda\in\Lambda}(U\cap Z_{\lambda})}$ is obvious.
Take any $x\in\bigcap_{\lambda\in\Lambda}\ovl{U\cap Z_{\lambda}}$ and any quasi-compact open neighborhood $V$ of $x$.
Then $V\cap (U\cap Z_{\lambda})\neq\emptyset$ for any $\lambda\in\Lambda$.
By {\bf \ref{ch-pre}}.\ref{thm-projlimcohsch1} (2) we deduce that $V\cap (\bigcap_{\lambda\in\Lambda}U\cap Z_{\lambda})\neq\emptyset$.
Since this holds for any quasi-compact open neighborhood $V$ of $x$, we have $x\in \ovl{\bigcap_{\lambda\in\Lambda}(U\cap Z_{\lambda})}$, whence the first equality.
The other equality $\bigcap_{\lambda\in\Lambda}\ovl{U\cap Z_{\lambda}}=\ovl{U}\cap(\bigcap_{\lambda\in\Lambda}Z_{\lambda})$ is similarly shown.
\end{proof}

\begin{prop}\label{prop-appclassicalZR1}
Let $Y$ be a coherent scheme, and $f\colon X\rightarrow Y$ a separated morphism of finite type. 
Let $\{Z_{\lambda}\}_{\lambda\in\Lambda}$ be a filtered projective system consisting of closed subscheme of $X$ indexed by a directed set $\Lambda$ such that for each $\lambda\geq\mu$ the transition maps $i_{\lambda\mu}\colon Z_{\lambda}\rightarrow Z_{\mu}$ are the closed immersions over $X$.
If the projective limit $\varprojlim_{\lambda\in\Lambda}Z_{\lambda}$ is proper over $Y$, then there exists $\lambda_0\in\Lambda$ such that $Z_{\lambda}$ is proper over $Y$ for any $\lambda\geq\lambda_0$.
\end{prop}

\begin{proof}
Take a proper $Y$-scheme $\ovl{f}\colon\ovl{X}\rightarrow Y$ together with an dense $Y$-open immersion $X\hookrightarrow\ovl{X}$ (\ref{thm-nagataembedding}).
Let $\ovl{Z}_{\lambda}$ be the scheme-theoretic closure of $Z_{\lambda}$ in $\ovl{X}$ for each $\lambda\in\Lambda$.
Suppose that $\ovl{Z}_{\lambda}\cap(\ovl{X}\setminus X)\neq\emptyset$ for any $\lambda\in\Lambda$.
Then by quasi-compactness of $X$ we have $\bigcap_{\lambda\in\Lambda}\ovl{Z}_{\lambda}\cap(\ovl{X}\setminus X)\neq\emptyset$.
On the other hand, we have by \ref{lem-appclassicalZR} the equality $\bigcap_{\lambda\in\Lambda}\ovl{Z}_{\lambda}=\ovl{\bigcap_{\lambda\in\Lambda}Z_{\lambda}}$.
Hence, we have $\bigcap_{\lambda\in\Lambda}\ovl{Z}_{\lambda}\cap(\ovl{X}\setminus X)=\ovl{(\bigcap_{\lambda\in\Lambda}Z_{\lambda})}\cap(\ovl{X}\setminus X)\neq\emptyset$.
But since $\bigcap_{\lambda\in\Lambda}Z_{\lambda}$ is proper over $Y$, it is closed in $\ovl{X}$, and hence $\ovl{(\bigcap_{\lambda\in\Lambda}Z_{\lambda})}\cap(\ovl{X}\setminus X)=\emptyset$, which is absurd.
Hence there exists $\alpha_0$ such that $\ovl{Z_{\lambda}}\subseteq X$ for $\lambda\geq\lambda_0$, that is, $Z_{\lambda}$ is proper.
\end{proof}

\begin{cor}\label{cor-appclassicalZR1}
Let $Y$ be a coherent scheme, and $f\colon X\rightarrow Y$ a separated morphism of finite type.
Let $Z$ be a closed subscheme of $X$ that is proper over $Y$.
Then there exists a closed subscheme $Z'$ of $X$ of finite presentation that is proper over $Y$ and contains $Z$ scheme-theoretically.
\end{cor}

\begin{proof}
This follows immediately from \cite[(6.9.15)]{EGAInew} and \ref{prop-appclassicalZR1}.
\end{proof}

\begin{proof}[Proof of Proposition {\rm \ref{prop-appclassicalZR}}]
We may assume that $Y$ is affine.
If $f$ is of finite presentation, then there exist a proper and affine morphism $f'\colon X'\rightarrow Y'$ between Noetherian schemes (\cite[$\mathbf{IV}$, (8.10.5)]{EGA}) and a morphism $Y\rightarrow Y'$ such that $f=f'_Y$.
Since the assertion is well-known in the Noetherian case, we have the desired result in case $f$ is finitely presented.

In general, write $X$ as a closed subscheme of $\A^n_Y$ for some $n\geq 0$.
By \ref{cor-appclassicalZR1} there exists a closed subscheme $X'$, proper over $Y$, of $\A^n_Y$ of finite presentation that contains $X$ scheme-theoretically.
Since $X'$ is finite over $Y$, so is $X$.
\end{proof}

\subsection{Nagata embedding for algebraic spaces}\label{sub-nagataembalgsp}
\index{Nagata's embedding theorem!Nagata's embedding theorem for algebraic spaces@--- for algebraic spaces|(}
Finally, we include here the statement of the compactification theorem for algebraic spaces (\cite{Fujiw3}, \cite{CLO}):
\begin{thm}[Nagata's embedding theorem for algebraic spaces]\label{thm-nagataembalgsp}
Let $Y$ be a coherent algebraic space, and $f\colon X\rightarrow Y$ a separated $Y$-algebraic space\index{algebraic space} of finite type.
Then there exists a proper $Y$-algebraic space $\ovl{f}\colon\ovl{X}\rightarrow Y$ that admits a dense open immersion $X\hookrightarrow\ovl{X}$ over $Y;$ moreover$:$
\begin{itemize}
\item[{\rm (a)}] there exists a quasi-coherent ideal $\mathscr{I}$ of $\O_{\ovl{X}}$ of finite type such that the underlying topological subspace of the corresponding closed subspace coincides with the boundary $\partial X=\ovl{X}\setminus X$ and that $(\partial X,\O_{\ovl{X}}/\mathscr{I})$ is a scheme$;$
\item[{\rm (b)}] if $X$ is a scheme, then there exists a compactification $\ovl{f}\colon\ovl{X}\rightarrow Y$ as above such that $\ovl{X}$ is a scheme. \hfill$\square$
\end{itemize}
\end{thm}

Notice that, in view of {\bf \ref{ch-pre}}.\ref{cor-KnudsonSerrecriterionabsolute}, the contents of (a) can be rephrased without ambiguity as that the boundary $\partial X=\ovl{X}\setminus X$ `is' a scheme.
\index{Nagata's embedding theorem!Nagata's embedding theorem for algebraic spaces@--- for algebraic spaces|)}

\addcontentsline{toc}{subsection}{Exercises}
\subsection*{Exercises}
\begin{exer}\label{exer-partialcomp}{\rm 
Let $Y$ be a coherent scheme, $f\colon X\rightarrow Y$ a separated $Y$-scheme of finite type, and $U\subseteq X$ a quasi-compact open subset.
Then the open immersion $j\colon U\hookrightarrow X$ induces a closed map $\ZR{j}_{\cpt}\colon\ZR{U}_{\cpt}\rightarrow\ZR{X}_{\cpt}$.
Moreover, $\ZR{j}_{\cpt}$ maps the open subset $\ZR{U}^X_{\pc}$ isomorphically onto an open subset of $\ZR{X}_{\cpt}$.}
\end{exer}

\begin{exer}\label{exer-nagatachow}{\rm 
Let $Y$ be a coherent scheme, and $f\colon X\rightarrow Y$ a separated $Y$-algebraic space of finite type.
Let $U$ be a quasi-compact open subspace of $X$ that is a scheme.
Then there exists a $U$-admissible blow-up $X'$ of $X$ that is a scheme.}
\end{exer}




%
\addcontentsline{toc}{chapter}{\rm Solutions and Hints for Exercises}
\markboth{Solutions and Hints}{Solutions and Hints}
\chapter*{Solutions and Hints for Exercises}
\section*{Chapter \ref{ch-pre}}
{\bf Exercise \ref{exer-finalcountable}.} We may assume that $I$ is countable.
To skip the trivial cases, we may assume, moreover, that $I$ does not have maximal elements.
For any fixed $\alpha_0\in I$ we may replace $I$ by the final subset $\{\alpha\in I\,|\,\alpha\geq \alpha_0\}$ and thus may assume that $I$ has the minimum element $\alpha_0$.
Write $I=\{\alpha_0,\alpha_1,\alpha_2,\ldots\}$.
We set $L(0)=\alpha_0$ and $L(1)=\alpha_1$ and define $L(k+1)=\alpha_{l(k+1)}$ for $k\geq 1$ inductively as follows: $l(k+1)$ is the smallest number such that $\alpha_i<\alpha_{l(k+1)}$ for $l(k-1)+1\leq i\leq l(k)$.
Then the resulting map $L\colon\N\rightarrow I$ is final.

\medskip\noindent
{\bf Exercise \ref{exer-projectivelimitaffine}.} The projective limit $X$ is the underlying topological space of a scheme affine over a coherent scheme; see \cite[$\mathbf{II}$, \S8.2]{EGA}.

\medskip\noindent
{\bf Exercise \ref{exer-generizationmaximal}.} 
(1) Consider for each $i\in I$ the set $G_{x_i}$ of all generizations of $x_i$, which is coherent and sober due to {\bf \ref{ch-pre}}.\ref{cor-projlimcohspacelattice1}.
For $i\leq j$ there exists the canonical inclusion map $G_{x_j}\hookrightarrow G_{x_i}$, which is quasi-compact due to {\bf \ref{ch-pre}}.\ref{cor-quasicompactness}.
Then the intersection $\bigcap_{i\in I}G_{x_i}$ is non-empty due to {\bf \ref{ch-pre}}.\ref{thm-projlimcohsch1}.

(2) Use (1) and apply Zorn's lemma.

\medskip\noindent
{\bf Exercise \ref{exer-coherentprojlimopencovering}.} 
Set $U=\bigcup^n_{k=1}U_k$, and apply {\bf \ref{ch-pre}}.\ref{cor-projlimcohsch11} to the situation $V=X_i$.

\medskip\noindent
{\bf Exercise \ref{exer-limitcoherentsoberspaceconnected}.}
Suppose that $X$ is not connected, and take non-empty open subsets $U_0,U_1\subseteq X$ such that $X=U_0\cup U_1$ and $U_0\cap U_1=\emptyset$.
Since $U_0,U_1$ are closed and $X$ is quasi-compact ({\bf \ref{ch-pre}}.\ref{thm-projlimcohsch1} (1)), $U_0,U_1$ are quasi-compact.
By {\bf \ref{ch-pre}}.\ref{prop-projlimcohtopspqcptopen} there exist quasi-compact open subsets $U_{0i},U_{1i}\subseteq X_i$ for some $i\in I$ such that $U_0=p^{-1}_i(U_{0i})$ and $U_1=p^{-1}_i(U_{1i})$, where $p_i\colon X\rightarrow X_i$ is the projection map.
Since $p_i$ is surjective ({\bf \ref{ch-pre}}.\ref{cor-projlimcohspacelattice2}), one deduces easily that $X_i=U_{0i}\cup U_{1i}$ and $U_{0i}\cap U_{1i}=\emptyset$.

\medskip\noindent
{\bf Exercise \ref{exer-valuativespaceT1}.} 
Let $x,y\in[X]$ with $x\neq y$.
Since $\sep^{-1}_X(x)=\ovl{\{x\}}$ is an overconvergent closed subset of $X$, there exists by {\bf \ref{ch-pre}}.\ref{prop-separation21a} a unique open subset $U\subseteq[X]$ such that $\sep^{-1}_X(U)=X\setminus\ovl{\{x\}}$.
In particular, we have $y\in U$ and $x\not\in U$.

\medskip\noindent
{\bf Exercise \ref{exer-locallystronglycompact1}.} 
(1) follows from {\bf \ref{ch-pre}}.\ref{thm-locallystronglycompacttheorem}.
To show (2), we may assume in view of {\bf \ref{ch-pre}}.\ref{thm-locallystronglycompacttheorem} that $X$ is quasi-separated, since the question is local on $[X]$.
Take for any $x\in U\cap V$, in view of {\bf \ref{ch-pre}}.\ref{prop-locallycompactnessformerdef}, a pair $(U_x,V_x)$ of coherent open neighborhoods in $V$ of the closure of $\{x\}$ in $V$ such that $V_x$ contains the closure $\ovl{U_x}$ of $U_x$ in $V$.
Then $(U\cap U_x,U\cap V_x)$ gives a pair of open neighborhoods of $x$ in $U\cap V$.
Since $U\cap U_x\hookrightarrow U_x$ and $U\cap V_x\hookrightarrow V_x$ are quasi-compact, $U\cap U_x$ and $U\cap V_x$ are coherent.
As we have $U\cap\ovl{U_x}\subseteq U\cap V_x$, $(U\cap U_x,U\cap V_x)$ satisfies the condition as in {\bf \ref{ch-pre}}.\ref{prop-locallycompactnessformerdef}.

\medskip\noindent
{\bf Exercise \ref{exer-locallystronglycompact2}.} 
(1) Since the question is local on $Y$, we may assume that $Y$ is Hausdorff. 
Replacing $Y$ by the image of $f$ endowed with the subspace topology from $Y$, we may also assume that $f$ is bijective.
For any $x\in X$, take a relatively compact open neighborhood $U_x$ of $x$ in $X$.
Then, since $\ovl{U_x}$ is compact, $f\vert_{\ovl{U_x}}$ is a homeomorphism onto its image.
Hence $f$ is a local homeomorphism.
Since $f$ is bijective, we deduce that $f$ is a homeomorphism.

(2) It suffices to show that a locally compact subspace of $Y$ is open in its closure.
One can reduce by localization to the case where $Y$ is Hausdorff, and the claim in this situation is well-known.

(3) By {\bf \ref{ch-pre}}.\ref{thm-locallystronglycompacttheorem} the spaces $[U]$ and $[X]$ are locally compact locally Hausdorff spaces.

\medskip\noindent
{\bf Exercise \ref{exer-locallystronglycompact3}.} 
(1) To show the `if' part, we may assume in view of {\bf \ref{ch-pre}}.\ref{thm-locallystronglycompacttheorem} that $X$ is quasi-separated, since the question is local on $[X]$.
Let $V\subseteq X$ be a coherent open subset.
We need to show that $U\cap V$ is quasi-compact.
Since $[V]$ is compact and $[U]$ is closed, $[Z]\cap[U]=[Z\cap U]$ is compact.
On the other hand, since $U\cap V$ is locally strongly compact (Exercise \ref{exer-locallystronglycompact1}), the separation map $\sep_{U\cap V}$ is proper ({\bf \ref{ch-pre}}.\ref{thm-locallystronglycompacttheorem}), and hence $U\cap V$ is quasi-compact, as desired.

To show the `only if' part, consider an open covering $X=\bigcup_{\alpha\in L}U_{\alpha}$ by coherent open subsets.
Since $U\cap U_{\alpha}$ is quasi-compact, $[U]\cap[U_{\alpha}]=[U\cap U_{\alpha}]$ is compact, and hence is closed in the Hausdorff space $[U_{\alpha}]$ ({\bf \ref{ch-pre}}.\ref{cor-tubes2acor}) for each $\alpha\in L$.
Then it follows from {\bf \ref{ch-pre}}.\ref{cor-propstructuresepquot} that $[U]$ is closed in $[X]$.

(2) Take an open subset $\mathscr{Z}\subseteq[X]$ of $[X]$ that contains $[U]$ as a closed subset.
Then $Z=\sep^{-1}_X(\mathscr{Z})$ satisfies the desired property.

\medskip\noindent
{\bf Exercise \ref{exer-locallystronglycompact4}.} 
Let $U,V$ be locally strongly compact open subsets of a locally strongly compact valuative space $X$.
We want to show that $U\cap V$ is locally strongly compact. (Here, recall that $X$ has an open basis consisting of coherent open subsets.)
Since the problem is local on $[X]$, we may assume in view of Exercise \ref{exer-locallystronglycompact2} (3) that $[U]$ is closed in $[X]$.
By Exercise \ref{exer-locallystronglycompact3} (1) the inclusion map $U\hookrightarrow X$ is quasi-compact.
Then apply Exercise \ref{exer-locallystronglycompact1} (2).

\medskip\noindent
{\bf Exercise \ref{exer-locallystronglycompact5}.} 
We want to show that, for any coherent open subset $U\subseteq X$, the inclusion map $U\hookrightarrow X$ is quasi-compact.
Since $U$ is locally strongly compact ({\bf \ref{ch-pre}}.\ref{prop-locallycompactspacecoherent}), $[U]$ is identified with a subspace of $[X]$ (Exercise \ref{exer-locallystronglycompact2} (3)).
Since $[U]$ is compact and $[X]$ is Hausdorff, $[U]$ is closed in $[X]$.
Hence $U\hookrightarrow X$ is quasi-compact by Exercise \ref{exer-locallystronglycompact3} (1).

\medskip\noindent
{\bf Exercise \ref{exer-injlimcohcoh}.} 
See \cite{Kempf}, \S2.

\medskip\noindent
{\bf Exercise \ref{exer-coherencedirectlimit}.}
Consider the homomorphism of the form $A^{\oplus p}\rightarrow A$. 
We want to show that its kernel is finitely generated.
There exist $i\in I$ and a map $A^{\oplus p}_i\rightarrow A_i$ that induces the above map by the tensor product with $A$.
Since $A_i$ is coherent, its kernel $K$ is finitely generated.
Since $A$ is flat over $A_i$, the kernel of $A^{\oplus p}\rightarrow A$ is given by $K\otimes_{A_i}A$.

\medskip\noindent
{\bf Exercise \ref{exer-algintegernonNoe}.} 
The ring in question is coherent due to Exercise \ref{exer-coherencedirectlimit}.
To see it is not Noetherian, consider a prime number $p$ and the sequence of ideals
$$
(p)\subseteq (\sqrt{p})\subseteq ({}^4\!\!\!\sqrt{p})\subseteq\cdots\subseteq ({}^{2^n}\!\!\!\!\sqrt{p})\subseteq\cdots,
$$
which one can show is strictly increasing using the $p$-adic valuation. 

\medskip\noindent
{\bf Exercise \ref{exer-generizationmapsch}.} We may assume that $X$ is affine $X=\Spf A$ where $A$ is a Noetherian adic ring with an ideal of definition $I\subseteq A$.
Let $\mathfrak{p}$ (resp.\ $\mathfrak{q}$) be the open prime ideal corresponding to $x$ (resp.\ $y$); we have $\mathfrak{p}\supseteq\mathfrak{q}$.
Then $\O_{X,x}=A_{\{S\}}$ and $\O_{X,y}=A_{\{T\}}$ (in the notation as in \cite[$\mathbf{0}_{\mathbf{I}}$, \S7.6]{EGA}; notice that these are Noetherian rings (\cite{EGA}, $\mathbf{0}_{\mathbf{I}}$, (7.6.18))), where $S=A\setminus\mathfrak{p}$ and $T=A\setminus\mathfrak{q}$.
We have a commutative square of rings
$$
\xymatrix@-2ex{(\widehat{A}_{\mathfrak{p}})_{\mathfrak{q}\widehat{A}_{\mathfrak{p}}}\ar[r]&\widehat{A}_{\mathfrak{q}}\\ (A_{\{S\}})_{\mathfrak{q}A_{\{S\}}}\ar[u]\ar[r]& A_{\{T\}}\ar[u]\rlap{,}}
$$
where $\widehat{A}_{\mathfrak{p}}$ (resp.\ $\widehat{A}_{\mathfrak{q}}$) denotes the $I$-adic completion of the localization $A_{\mathfrak{p}}$ (resp.\ $A_{\mathfrak{q}}$).
What to show is that the lower horizontal arrow is faithfully flat.
By \cite[$\mathbf{0}_{\mathbf{I}}$, (7.6.2), (7.6.18)]{EGA} the vertical arrows are faithfully flat.
Moreover, since the $I$-adic completion of $(\widehat{A}_{\mathfrak{p}})_{\mathfrak{q}\widehat{A}_{\mathfrak{p}}}$ is nothing but $\widehat{A}_{\mathfrak{q}}$, the upper horizontal arrow is faithfully flat (\cite[Chap.\ III, \S3.5, Prop.\ 9]{Bourb1}).
Hence $(A_{\{\mathfrak{p}\}})_{\mathfrak{q}A_{\{\mathfrak{p}\}}}\rightarrow A_{\{\mathfrak{q}\}}$ is faithfully flat.

\medskip\noindent
{\bf Exercise \ref{exer-limitflatnesspreserve}.} 
Using {\bf \ref{ch-pre}}.\ref{prop-limLRS}, reduce the questions into those of inductive limits of rings.
Then use the results in {\bf \ref{ch-pre}}, \S\ref{subsub-directlimitspre}.

\medskip\noindent
{\bf Exercise \ref{exer-limitcohesivepreserve}.} Let $p_i\colon X\rightarrow X_i$ be the projection for each $i\in I$.
We need to show that, for any quasi-compact open subset $U\subseteq X$, the kernel of a morphism of $\O_U$-modules of the form
$$
\varphi\colon\O^{\oplus p}_U\longrightarrow\O_U
$$
is of finite type. 
To show this, in view of {\bf \ref{ch-pre}}.\ref{prop-projlimcohtopspqcptopen} and {\bf \ref{ch-pre}}.\ref{cor-lemlimLRS}, we may assume $U=X$.
By {\bf \ref{ch-pre}}.\ref{thm-injlimmodpf} there exists $i\in I$ and a morphism $\varphi_i\colon\mathscr{F}_i\rightarrow\mathscr{G}_i$ of finitely presented $\O_{X_i}$-modules such that $p^{\ast}_i\varphi_i\cong\varphi$.
Since $X_i$ is cohesive, these sheaves are coherent $\O_{X_i}$-modules, and hence $\mathscr{K}_i=\ker\varphi_i$ is again coherent.
Now by Exercise \ref{exer-limitflatnesspreserve} (1) it follows that the induced sequence $0\rightarrow p^{\ast}_i\mathscr{K}_i\rightarrow\O^{\oplus p}_X\rightarrow\O_X$ is exact, where $p^{\ast}_i\mathscr{K}_i$ is clearly of finite type.

\medskip\noindent
{\bf Exercise \ref{exer-lemqcptcomplementsupp}.} The `if' part is clear.
Suppose $U$ is quasi-compact, and consider a closed subscheme $Y$ with the underlying topological space $X\setminus U$.
Let $\mathscr{J}$ be the defining ideal of $Y$ in $X$.
By \cite[$\mathbf{I}$, (9.4.9) \& $\mathbf{IV}$, (1.7.7)]{EGA} one can write $\mathscr{J}$ as a filtered inductive limit of quasi-coherent subideals $\mathscr{J}_i$ of $\O_X$ of finite type.
Set $U_i=X\setminus V(\mathscr{J}_i)$.
Then $U=\bigcup U_i$.
Since $U$ is quasi-compact, we have $U=U_i$ for some $i$.

\medskip\noindent
{\bf Exercise \ref{exer-closedimmersionnilpotentthickening}.}
By \cite[(2.3.5)]{EGAInew} we know that $f$ is an affine morphism.
Thus we can assume that $X$ and $Y$ are affine $X=\Spec A$ and $Y=\Spec B$.
Let $I\subseteq A$ be the nilpotent finitely generated ideal such that $\til{I}=\mathscr{I}$.
By induction with respect to $n\geq 1$ such that $I^n=0$, we may assume that $I^2=0$.
The assumption is that the induced map $A/I\rightarrow B/IB$ is surjective.
By \cite[Theorem 8.4]{Matsu} (see {\bf \ref{ch-pre}}.\ref{prop-complpair1}) $B$ is finitely generated as an $A$-module.
To show that $A\rightarrow B$ is surjective, it suffices to show that the map $I\rightarrow IB$ is surjective (by snake lemma).
Let $C$ be the cokernel of $I\rightarrow IB$.
Then $C$ is a finitely generated $A$-module such that $IC=0$.
By Nakayama's lemma we have $C=0$.

\medskip\noindent
{\bf Exercise \ref{exer-derivedformalism3}.} We mimic the proof of \cite[Expos\'e II, Cor.\ 2.2.2.1]{SGA6}, which refers to \cite[Expos\'e II, Prop.\ 2.2.2]{SGA6}, in which the Noetherian hypothesis is used only in the end of the proof, where one refers to \cite[$\mathbf{I}$, (9.4.9)]{EGA}.
It uses the fact that quasi-coherent sheaves of finite type on locally Noetherian schemes are coherent.
We replace this part of the proof by the following reasoning:
for any quasi-coherent sheaf $\mathscr{F}$ of finite type, there exists a coherent sheaf $\mathscr{G}$ and a surjective morphism $\mathscr{G}\rightarrow\mathscr{F}$.
This follows from \cite[$\mathbf{II}$, (2.7.9)]{EGA}.

\medskip\noindent
{\bf Exercise \ref{exer-valtorsionflat}.}  (3) Embed $\Spec A$ into $(\P^1_V)^n$, and show that the closure is finite using the projections $(\P^1_V)^n\to\P^1_V$.

\medskip\noindent
{\bf Exercise \ref{exer-valflatatorfree}.} By Exercise \ref{exer-valtorsionflat} (1) it suffices to show that $M$ is torsion free if and only if $M$ is $a$-torsion free.
This follows from the following observation: for any $x\in V\setminus\{0\}$, since $\bigcap_{n\geq 0}(a^n)=0$, there exists $n\geq 0$ such that $x\not\in(a^n)$; then we have $a^n=xy$ for some $y\in V$.

\medskip\noindent
{\bf Exercise \ref{exer-valuationdrill2}.} (1) By {\bf \ref{ch-pre}}.\ref{prop-isolated3} the totally ordered commutative group $\Gamma_V$ is ordered isomorphic to an ordered subgroup of $\R^d$ (where $d=\mathrm{ht}(V)$) with the lexicographical order (cf.\ {\bf \ref{ch-pre}}.\ref{exa-height}).
Take $a\in V$ such that $v(a)=(i_1,\ldots,i_d)$ with $i_1>0$.
Then for any element $b\in V\setminus\{0\}$ there exists $n\geq 0$ such that $v(b)<v(a^n)$, that is, $b\not\in a^nV$.
Hence we have $\bigcap_{n\geq 0}a^nV=\{0\}$.

(2) Consider the quotient ring $V/\mathfrak{p}$, which is again a valuation ring such that $0<\mathrm{ht}(V/\mathfrak{p})<+\infty$ ({\bf \ref{ch-pre}}.\ref{prop-composition1} (1) and {\bf \ref{ch-pre}}.\ref{prop-composition3}).
Then apply (1).

\medskip\noindent
{\bf Exercise \ref{exer-filtrationtopclosure}.} An element $x\in M$ lies in $\ovl{N}$ if and only if $(x+F^{\lambda})\cap N\neq\emptyset$ for any $\lambda\in\Lambda$, where the last condition is equivalent to $x\in N+F^{\lambda}$ for any $\lambda\in\Lambda$.

\medskip\noindent
{\bf Exercise \ref{exer-principalclosed}.}
Consider the closure $J=\bigcap_{\lambda\in\Lambda}(gA+I^{(\lambda)})$ of the ideal $gA$.
Let $x\in J$.
For any $\lambda\in\Lambda$, we have $a_{\lambda}\in A$ and $b_{\lambda}\in I^{(\lambda)}$ such that $x=ga_{\lambda}+b_{\lambda}$.
For $\lambda\leq\mu$, we have $g(a_{\lambda}-a_{\mu})=b_{\mu}-b_{\lambda}\in I^{(\lambda)}$.
Since $(g$ mod $I^{(\lambda)})$  is a non-zero-divisor in $A/I^{(\lambda)}$, we have $a_{\lambda}-a_{\mu}\in I^{(\lambda)}$ for any $\mu\geq\lambda$.
Hence $\{a_{\lambda}\}_{\lambda\in\Lambda}$ is a Cauchy sequence, converging to an element $a\in A$.
Since $x-ga\in\bigcap_{\lambda\in\Lambda}I^{(\lambda)}=\{0\}$, we have $x\in gA$.

\medskip\noindent
{\bf Exercise \ref{exer-filtrationcompletion}.}
(1) By {\bf \ref{ch-pre}}.\ref{cor-ML4} the exact sequences
$$
0\longrightarrow N/F^{(n)}\longrightarrow M/F^{(n)}\longrightarrow M/N\longrightarrow 0
$$
induces, by passage to the projective limits, the exact sequence
$$
0\longrightarrow N^{\wedge}_{F^{\bullet}}\longrightarrow M^{\wedge}_{F^{\bullet}}\longrightarrow M/N\longrightarrow 0
$$
(cf.\ {\bf \ref{ch-pre}}, \S\ref{subsub-completionexactsequence}).

(2) If $F^{\bullet}$ is separated, then the canonical maps $M\rightarrow M^{\wedge}_{F^{\bullet}}$ and $N\rightarrow N^{\wedge}_{F^{\bullet}}$ are injective.
Since $M/N\cong M^{\wedge}_{F^{\bullet}}/N^{\wedge}_{F^{\bullet}}$, we have $N^{\wedge}_{F^{\bullet}}\cap M=N$.

\medskip\noindent
{\bf Exercise \ref{exer-topologicaltensorproducts}.}
See \cite[$\mathbf{0}_{\mathbf{I}}$, \S7.7]{EGA}.

\medskip\noindent
{\bf Exercise \ref{exer-counterexaflatness1}.}
The $a$-adic completion $\widehat{V}$ coincides with the $a$-adic completion of $V/\mathfrak{p}$, where $\mathfrak{p}=\bigcap_{n\geq 0}(a^n)$.

\medskip\noindent
{\bf Exercise \ref{exer-APverifyinglemma}.}
We show that the condition implies {\bf (AP)} (resp.\ {\bf (APf)}).
Let $M$ be a finitely generated $A$-module, and $N\subseteq M$ an $A$-module (resp.\ a finitely generated $A$-submodule).
For any fixed $n\geq 0$, consider $\ovl{N}=N/I^nN$, which is an $A$-submodule of $\ovl{M}=M/I^nN$.
As $\ovl{N}\cap I^m\ovl{M}=0$ for some $m\geq 0$, we have $N\cap I^mM\subseteq I^nN$.

\medskip\noindent
{\bf Exercise \ref{exer-auxiliaryexeroncompletion}.}
Consideran exact sequence of the form
$$
0\longrightarrow K\longrightarrow A^{\oplus n}\longrightarrow M\longrightarrow 0
$$
and the induced commutative diagram with exact rows
$$
\xymatrix@-1ex{0\ar[r]&K^{\wedge}_{I^{\bullet}A^{\oplus n}\cap K}\ar[r]&\widehat{A}^{\oplus n}\ar[r]&\widehat{M}\ar[r]&0\\
&K\otimes_A\widehat{A}\ar[u]\ar[r]_(.55){(\ast)}&\widehat{A}^{\oplus n}\ar@{=}[u]\ar[r]&M\otimes_A\widehat{A}\ar[u]\ar[r]&0\rlap{,}}
$$
where $K^{\wedge}_{I^{\bullet}A^{\oplus n}\cap K}$ denotes the Hausdorff completion of $K$ with respect to the induced filtration $I^{\bullet}A^{\oplus n}\cap K$.
By {\bf \ref{ch-pre}}.\ref{cor-propARconseq1-2} the image of $(\ast)$ is closed in $\widehat{A}^{\oplus n}$ with respect to the $I$-adic topology.
This implies that the left-hand vertical map is surjective.
The right-hand vertical map is an isomorphism by snake lemma.

\medskip\noindent
{\bf Exercise \ref{exer-zriskiancovering}.}
We need to show that the map $\prod^r_{i=1}\Spec\zat{(A_{f_i})}\rightarrow\Spec A$ is surjective.
It suffices to show the following: for any prime ideal $\mathfrak{p}\subset A$ there exists $i$ $(1\leq i\leq r)$ such that $f_i^n+ax\not\in\mathfrak{p}$ for any $n\geq 0$, $a\in I$, and $x\in A$.
Suppose contrarily that for each $1\leq i\leq r$ we have $f^{n_i}_i+a_ix_i\in\mathfrak{p}$ for some $n_i$, $a_i$, and $x_i$.
One can assume $n_1=\cdots=n_r=n$.
Since $(f^n_1,\ldots,f^n_r)=A$, there exists $b_1,\ldots,b_n\in A$ such that $\sum^r_{i=1}b_if^n_i=1$.
Then we have $1+\sum^r_{i=1}a_ib_ix_i\in\mathfrak{p}\cap(1+I)$, which is absurd, for $A$ is $I$-adically Zariskian.

\medskip\noindent
{\bf Exercise \ref{exer-aadictopcompleteseparatedclosed}.}
Set $L=\image(f)$.
By {\bf \ref{ch-pre}}.\ref{lem-ARprincipalgeneralized} the topology on $L$ defined by the filtration $\{L\cap a^nN\}_{n\geq 0}$ coincides with the $a$-adic topology.
Since $N$ is $a$-adically separated, $L$ is $a$-adically separated, and hence the kernel $\ker(f)$ is closed in $M$ with respect to the $a$-adic topology.
By {\bf \ref{ch-pre}}.\ref{cor-qconsistency1111} we deduce that $L$ is $a$-adically complete and hence is Hausdorff complete with respect to the filtration $\{L\cap a^nN\}_{n\geq 0}$.
This implies that $L$ is closed in $N$ (cf.\ Exercise \ref{exer-filtrationtopclosure}).

\medskip\noindent
{\bf Exercise \ref{exer-closuresaturated}.} 
We may assume $I\subseteq (f)$.
Suppose $f^mx\in\ovl{\mathfrak{a}}$ for $x\in A$ and $m\geq 1$.
By Exercise \ref{exer-filtrationtopclosure}, for any $n\geq 1$ we have $f^mx=a_n+y_n$, where $a_n\in\mathfrak{a}$ and $y_n\in I^n$.
Since $I\subseteq (f)$, one can write $y_n=f^mz_n$ for any $n>m$, where $z_n\in I^{n-m}$.
Then for $n>m$ we have $f^m(x-z_n)=y_n\in\mathfrak{a}$.
Since $\mathfrak{a}$ is $f$-saturated, we have $x-z_n\in\mathfrak{a}$ for any $n>m$, whence $x\in\ovl{\mathfrak{a}}$.

\medskip\noindent
{\bf Exercise \ref{exer-adhesivecounterexamples}.} Due to Exercise \ref{exer-valuationdrill2} (2) there exists $a\in\m_V$ such that $\bigcap_{n\geq 0}a^n\til{V}$ is the height one prime ideal.
Show that $\til{V}[\frac{1}{a}]=W$ and apply {\bf \ref{ch-pre}}.\ref{prop-exaadhesiveval}.

\medskip\noindent
{\bf Exercise \ref{exer-intersecfg}.} Consider the exact sequence $0\rightarrow N\cap P\rightarrow N\oplus P\rightarrow N+P\rightarrow 0$.
Since $N+P\subseteq M$ is $I$-torsion free, it is finitely presented.

\medskip\noindent
{\bf Exercise \ref{exer-Ivaluativeringsmaps}.} 
(1) Since $h$ is adic, we have $h(J)\subseteq J'$, hence $h$ induces a local homomorphisms $V\rightarrow V'$.
Moreover, $h$ induces $g=h[\frac{1}{a}]\colon B\rightarrow B'$, which is local, for $h(J)\subseteq J'$.
Since $K=B/J\rightarrow B'/J'=K'$, being a homomorphism between fields, is injective, $V\rightarrow V'$ is also injective.
Moreover, since $V\rightarrow V'$ is local between valuation rings, one can easily show $V=K\cap V'$ in $K'$.

(2) Since $A=\{f\in B\,|\,(f\ \mathrm{mod}\ J)\in V\}$ and $V=K\cap V'$, we have $g^{-1}(A')=A$.

\medskip\noindent
{\bf Exercise \ref{exer-valuapproxfiniteheight0}.}
We first show the claim in the case where both $V$ and $V'$ are of height one. 
In this case, the real valued valuation on $V$ restricts to that of $V'$, and hence Cauchy sequences in $V$ are Cauchy sequences in $V'$. 
Hence we have $\widehat{V}\subseteq\widehat{V'}$ in this case.
In general, take the associated height one primes $\mathfrak{p}\subseteq V$ and $\mathfrak{p}'\subseteq V'$, respectively.
By {\bf \ref{ch-pre}}.\ref{cor-compval2006ver23}, we have $\widehat{V}\subseteq\widehat{V_{\mathfrak{p}}}$ and similarly for $V'$.
By this and the above-discussed height one case, we have the desired result.

\medskip\noindent
{\bf Exercise \ref{exer-valuapproxfiniteheight}.} (1) Write $V$ as the filtered inductive limit of subrings $V=\varinjlim_{\lambda\in\Lambda}A_{\lambda}$, where each $A_{\lambda}$ is finitely generated over $\Z$ (and hence is Noetherian).
Localizing at the prime ideal $\m_V\cap A_{\lambda}$, we may assume that $A_{\lambda}$ are local rings and that the maps $A_{\lambda}\rightarrow A_{\mu}$ and $A_{\lambda}\rightarrow V$ are local.
Let $K_{\lambda}=\Frac(A_{\lambda})$ for $\lambda\in\Lambda$, and consider the composite valuation $K_{\lambda}\hookrightarrow K\rightarrow\Gamma_V\cup\{\infty\}$, which defines, as the subset of all elements in $K_{\lambda}$ with non-negative values, a valuation ring $V_{\lambda}$ for $K_{\lambda}$ that dominates the local ring $A_{\lambda}$.
Since $A_{\lambda}$ is Notherian local, we know in view of {\bf \ref{ch-pre}}.\ref{thm-estimate} and {\bf \ref{ch-pre}}.\ref{prop-rationalrank1} that each $V_{\lambda}$ is of finite height.
Clearly, we have $V=\varinjlim_{\lambda\in\Lambda}V_{\lambda}$, since $A_{\lambda}\subseteq V_{\lambda}\subseteq V$ for each $\lambda\in\Lambda$.

(2) Suppose $V$ is $a$-adically separated, and write $V=\varinjlim_{\lambda\in\Lambda}V_{\lambda}$ as in (1).
We may assume that $a\in V_{\lambda}$ for any $\lambda\in\Lambda$.
Then each $V_{\lambda}$ is $a$-adically separated, since $\bigcap_{n\geq 0}a^nV_{\lambda}\subseteq\bigcap_{n\geq 0}a^nV=\{0\}$.
Suppose $V$ is $a$-adically complete.
Then by Exercise \ref{exer-valuapproxfiniteheight0}, we have $\widehat{V_{\lambda}}\subset V$, each of which is of finite height (see {\bf \ref{ch-pre}}.\ref{thm-compval2006ver1} (5)), and $V=\varinjlim_{\lambda\in\Lambda}\widehat{V_{\lambda}}$.

\medskip\noindent
{\bf Exercise \ref{exer-contentidealfg}.} First note that for any non-zero element $b\in V$ there exists $m\geq 0$ such that $a^m$ divides $b$ but $a^{m+1}$ does not; in fact, since $\bigcap_{m\geq 0}a^mV=0$, we have $b\in a^mV\setminus a^{m+1}V$ for some $m\geq 0$.
Let $f=\sum_{\nu_1,\ldots,\nu_n}b_{\nu_1,\ldots,\nu_n}X^{\nu_1}_1\cdots X^{\nu_n}_n$.
Then the numbers $m$ such that $b_{\nu_1,\ldots,\nu_n}\in a^mV\setminus a^{m+1}V$ increase as $|\nu_1+\cdots+\nu_n|\rightarrow\infty$.
Then the ideal $\Cont(f)$ is generated by $b_{\nu_1,\ldots,\nu_n}$'s that have the minimal $m$.

\medskip\noindent
{\bf Exercise \ref{exer-BGRlemma2pp261}.} 
(1) The case in $k=0$ follows from \cite[(7.1.1/2)]{BGR}; then the equality $\m^{k+1}_0\mathcal{A}=\m^{k+1}$ for any $k\geq 0$ follows immediately.
For $k>0$ we look at the following commutative diagram with exact rows:
$$
\xymatrix@C-1ex{0\ar[r]&\m^k/\m^{k+1}\ar[r]&\mathcal{A}/\m^{k+1}\ar[r]&\mathcal{A}/\m^k\ar[r]&0\\
0\ar[r]&\m^k_0/\m^{k+1}_0\ar[r]\ar[u]&A_0/\m^{k+1}_0\ar[r]\ar[u]&A_0/\m^k_0\ar[r]\ar[u]&0\rlap{.}}
$$
The rightmost vertical arrow is an isomorphism by induction with respect to $k$.
It suffices to show that the map $\m^k_0/\m^{k+1}_0\rightarrow\m^k/\m^{k+1}$ between finite-dimensional vector spaces over the residue field $\mathcal{A}/\m=A_0/\m_0$ is an isomorphism.
Since $\m^{k+1}_0\mathcal{A}=\m^{k+1}$, it is surjective.
Since one sees easily that $\m^{k+1}_0=\m^{k+1}\cap A_0$, it is injective, too.

(2) For any maximal ideal $\m\subseteq\mathcal{A}$, the $\m$-adic completion $\widehat{\mathcal{A}}_{\m}$ of $\mathcal{A}_{\m}$ is isomorphic to the $\m_0$-adic completion of the regular local ring $(A_0)_{\m_0}$, where $\m_0=\m\cap A_0$, and hence $\mathcal{A}_{\m}$ is regular. 
Then apply Serre's theorem (\cite[Theorem 19.3]{Matsu}).

\medskip\noindent
{\bf Exercise \ref{exer-exercisefurthertech2}.} 
Let us set $A'=A\otimes_VV'$ and $\mathfrak{p}=\sqrt{aV}$, the associated height one prime of $V$. Note that $V'=V_{\mathfrak{p}}$.

(1) Let us first assume that $A$ is $V$-flat.
If $V$ is of height one, then the desired results follow from Noether normalization ({\bf \ref{ch-pre}}.\ref{thm-noethernormalizationtype(V)}).
In general, we know that $A'$ is finite over $V'$, and that $A/\mathfrak{p}A$ is of finite type over $V/\mathfrak{p}$.
Since $A$ is $V$-flat, $A$ is of finite type over $V$ due to {\bf \ref{ch-pre}}.\ref{prop-patching2} (2).

In case $A$ is not necessarily $V$-flat, let $J=A_{\ator}$ and $B=A/J$.
Notice that $J$ is finitely generated.
Employing the inductive argument with respect to $n$ for which $a^nJ=0$, one may assume $aJ=0$.
Then one has the exact sequence:
$$
\Tor_1(B,V/a)\longrightarrow J\longrightarrow A/aA\longrightarrow B/aB\longrightarrow 0,\eqno{(\ast)}
$$
with $\Tor_1(B,V/a)=0$.
We know $B$ is of finite type over $V$, and take a generating set $\{x_1,\ldots,x_n\}$.
Take $y_i\in A$ ($i=1,\ldots,n$) that is mapped to $x_i$ in $B$.
Let $z_1,\ldots,z_m$ be the topological generator of $A$ over $V$, and consider the subring $C=V[z_1,\ldots,z_m,y_1,\ldots,y_n]\subseteq A$.
In view of $(\ast)$, we deduce that $J\subseteq C$ and $C/J\cong B$.

(3) By height-one localization and patching argument as above, one can reduce to the situation where $V$ is of height one.
Since $B=A/J$ is $V$-flat and finite outside $aV$, it is quasi-finite over $V$.
According to \cite{EGA}, {\bf IV}, (18.12.3), one has a decomposition $B=B'\times B''$, where $B'$ is finite over $V$, and $B''\otimes_Vk=0$ ($k=V/\m_V$ is the residue field).
We have $\widehat{B'}=\widehat{B}$ and $B''=\ker(B\rightarrow\widehat{B})$.
Now, set $A''=\ker(A\rightarrow\widehat{A})$; $A''$ is an ideal of $A$ consisting of $a$-divisible elements, that is, for any $f\in A''$ and any $n\geq 1$, there exists $g\in A''$ such that $f=a^ng$.
Since the $a$-torsion part $J$ of $A$ is bounded, we have $A''\cap J=\{0\}$.
By this and $J=\ker(\widehat{A}\rightarrow\widehat{B})$, one has $A''\cong B''$, which gives a section to the surjective morphism $A\rightarrow B''$ by the composition $A\rightarrow B=B'\times B''\rightarrow B''$.
Hence we have a decomposition $A\cong A'\times A''$ with $A'\hookrightarrow\widehat{A}$.
Since $A''$ is $a$-divisible, we have $A''\otimes_V(V/aV)=0$.
Moreover, one has the exact sequence
$$
0\longrightarrow J\longrightarrow A'\longrightarrow B'\longrightarrow 0.
$$
In particular, $A'$ is $a$-adically complete, since so are $J$ and $B'$.

\medskip\noindent
{\bf Exercise \ref{exer-preparationthm}.} Apply {\bf \ref{ch-pre}}.\ref{thm-division} to get $X^{\nu_1(f)}=uf+r$, where $r$ has no exponent in $(\nu_1(f),0,\ldots,0)+\N^n$, that is, $r$ is a polynomial in $X_1$ of degree $<\nu_1(f)$.
Set $g=X^{\nu_1(f)}_1-r$.
We have $uf=g$.
Dividing out by $\m_V$, we have the equality of polynomials $\ovl{u}\ovl{f}=\ovl{g}$ with coefficients in $k=V/\m_V$; since the leading degrees in $X_1$ of $\ovl{f}$ and $\ovl{g}$ coincides and their leading coefficients are unit, we deduce that $\ovl{u}\in k^{\times}$ and hence that $u$ is a unit.
For the uniqueness, observe that for a given $f$ the equality $X^{\nu_1(f)}=uf+r$ with $r$ being a polynomial in $X_1$ of degree $<\nu_1(f)$ determines $u$ and $r$.

\section*{Chapter \ref{ch-formal}}
{\bf Exercise \ref{exer-idealofdefinitionintersection}.} We may assume that $X=\Spf A$ by an admissible ring $A$ and that $\mathscr{I}=I^{\Delta}$ and $\mathscr{I}'=I^{\prime\Delta}$ by ideals of definition $I,I'\subseteq A$.
It is then easy to see that $I\cap I'$ is an ideal of definition.
If $F^{\bullet}=\{F^{\lambda}\}_{\lambda\in\Lambda}$ is a descending filtration of ideals that defines the topology on $A$, we have
\begin{equation*}
\begin{split}
(I\cap I')^{\Delta}&=\varprojlim_{F^{\lambda}\subseteq I\cap I'}I\cap I'/F^{\lambda}=\varprojlim_{F^{\lambda}\subseteq I\cap I'}I/F^{\lambda}\cap I'/F^{\lambda}\\ &=(\varprojlim I/F^{\lambda})\cap(\varprojlim I'/F^{\lambda})=I^{\Delta}\cap I^{\prime\Delta},
\end{split}
\end{equation*}
which shows that $\mathscr{I}\cap\mathscr{I}'$ is an ideal of definition.

\medskip\noindent
{\bf Exercise \ref{exer-adicformalschemesidealsofdefinitionfinitude}.} We may assume that $X=\Spf A$ and $\mathscr{I}=I^{\Delta}$ where $I$ is an ideal of definition of the adic ring $A$.
Since $A$ is $I$-adically complete, $I$ is finitely generated if so is $I/I^2$ ({\bf \ref{ch-pre}}.\ref{prop-complpair1}).
Then apply {\bf \ref{ch-formal}}.\ref{cor-adicformalschemeidealofdefinitionfinitetype11}.

\medskip\noindent
{\bf Exercise \ref{exer-disjointsumproduct}.} Confer \cite[(3.2.4)]{EGAInew}.

\medskip\noindent
{\bf Exercise \ref{exer-adqformalprop1}.} Take $g\in A$ such that $f-j(g)$ belongs to $I\widehat{A}$.
Then for an open prime ideal $\mathfrak{p}$ of $\widehat{A}$ the element $f$ does not belong to $\mathfrak{p}$ if and only if $j(g)$ does not belong to $\mathfrak{p}$, for $\mathfrak{p}$ contains the ideal of definition $I\widehat{A}$.

\medskip\noindent
{\bf Exercise \ref{exer-adicformalschemesinductivelimitadicmaps}.} By {\bf \ref{ch-formal}}.\ref{prop-formalindlimschadic} the space $X$ is an adic formal schemes of finite ideal type.
To show that $X\rightarrow Y$ is adic, we may work in the affine situation, and the claim follows from {\bf \ref{ch-pre}}.\ref{prop-exeramazingfactonalgebra}.

\medskip\noindent
{\bf Exercise \ref{exer-adicallyqcohalgebrafinitetypecriterion}.} We may assume that $X$ is affine $X=\Spf A$ where $A$ is an adic ring with the finitely generated ideal of definition $I$ such that $\mathscr{I}=I^{\Delta}$.
Then by {\bf \ref{ch-formal}}.\ref{prop-topqcoh1} we have $\mathscr{B}=B^{\Delta}$ for an $I$-adically complete $A$-algebra $B$.
It follows from {\bf \ref{ch-formal}}.\ref{thm-adicqcoh1} (2) that $\mathscr{B}$ is adically quasi-coherent sheaf of finite type if and only if $B$ is finitely generated as an $A$-module, which is further equivalent to that $B/IB$ is finitely generated ({\bf \ref{ch-pre}}.\ref{prop-complpair1}).
Since $\mathscr{B}/\mathscr{I}\mathscr{B}=\til{B/IB}$ ({\bf \ref{ch-formal}}.\ref{prop-adicqcohaff1f}), the claim follows from \cite[(1.4.3)]{EGAInew}.

\medskip\noindent
{\bf Exercise \ref{exer-adicquotientbynonadic1}.}
We may assume that $X$ is affine $X=\Spf A$ where $A$ is a t.u.\ rigid-Noetherian ring and that $\mathscr{I}=I^{\Delta}$ where $I\subseteq A$ is a finitely generated ideal of definition.
By {\bf \ref{ch-formal}}.\ref{thm-adicqcohpre1} the morphism $\varphi$ comes from a morphism $f\colon M\rightarrow N$ between finitely generated $A$-modules, and we have $\coker(\varphi)=\coker(f)^{\Delta}$, which is an adically quasi-coherent sheaf of finite type; in particular, its $\mathscr{I}$-torsion part is bounded, since $\coker(f)_{\Itor}$ is bounded $I$-torsion.
Let $K=\ker(f)$.
By {\bf \ref{ch-formal}}.\ref{prop-lemdeltasheafadicallyquasicoherent1-3} we deduce that $K^{\Delta}=\ker(\varphi)$.
Since $M_{\Itor}$ is bounded $I$-torsion, $K_{\Itor}$ is also bounded $I$-torsion.

\medskip\noindent
{\bf Exercise \ref{exer-extopenadiccoh1}.}
Set $X_0=(X,\O_X/\mathscr{I})$, and define $U_0$ similarly.
Then $\mathscr{F}_0=\mathscr{F}/\mathscr{I}\mathscr{F}$ is a quasi-coherent sheaf on $X_0$, and $\mathscr{G}_0=\mathscr{G}/\mathscr{I}\mathscr{F}|_U$ is a quasi-coherent subsheaf of $\mathscr{F}_0$ of finite type.
Since $X_0$ is a coherent scheme ({\bf \ref{ch-formal}}.\ref{prop-qsepformal1}) and $U_0$ is a quasi-compact open subset of $X_0$, we may apply \cite[$\mathbf{I}$, (9.4.7) \& $\mathbf{IV}$, (1.7.7)]{EGA} to deduce that there exists a quasi-coherent subsheaf $\mathscr{G}'_0$ of $\mathscr{F}_0$ of finite type which extends $\mathscr{G}_0$.
Let $\mathscr{G}'$ be the inverse image sheaf of $\mathscr{G}'_0$ by the canonical map $\mathscr{F}\rightarrow\mathscr{F}_0$.
By {\bf \ref{ch-formal}}.\ref{prop-exeradicallyqcohsheavesbypullback} it is an adically quasi-coherent subsheaf of $\mathscr{F}$ of finite type.

\medskip\noindent
{\bf Exercise \ref{exer-extopenadiccoh1approx}.} 
For any open subset $U\subseteq X$, consider the set $S_U$ of all finite type adically quasi-coherent subsheaves $\mathscr{H}$ of $\mathscr{G}|_U$ such that $\mathscr{I}^s|_U(\mathscr{G}|_U/\mathscr{H})=0$ for some $s>0$.
If $U$ is affine, one has $S_U\neq\emptyset$.
To show that $S_X\neq\emptyset$ by induction with respect to the number of quasi-compact open subsets in a finite covering of $X$, we may assume that $X$ is covered by two quasi-compact open subsets $U_1$ and $U_2$ such that $S_{U_1}$ and $S_{U_2}$ are non-empty.
Take $\mathscr{H}_i\in S_{X_i}$ for $i=1,2$.
Let $\mathscr{H}_{12}$ be the adically quasi-coherent subsheaf of finite type of $\mathscr{G}|_{U_1\cap U_2}$ on $U_1\cap U_2$ generated by $\mathscr{H}_1|_{U_1\cap U_2}$ and $\mathscr{H}_2|_{U_1\cap U_2}$.
Since $\mathscr{H}_{12}/\mathscr{H}_1|_{U_1\cap U_2}$ is bounded $\mathscr{I}$-torsion, by Exercise \ref{exer-extopenadiccoh1} one may extend the quasi-coherent sheaf $\mathscr{H}_{12}$ $(i=1,2)$ onto $U_i$ to get an adically quasi-coherent subsheaf of finite type $\til{\mathscr{H}}_i$ of $\mathscr{G}|_{U_i}$ that contains $\mathscr{H}_i$.
The sheaves $\til{\mathscr{H}}_1$ and $\til{\mathscr{H}}_2$ patch together to an element of $S_X$.
For $\mathscr{H}\in S_X$ the sheaf $\mathscr{G}/\mathscr{H}$ is annihilated by $\mathscr{I}^n$ for some $n>0$.
Hence it defined a quasi-coherent sheaf on the scheme $(X,\O_X/\mathscr{I}^n)$.
As it is an inductive limit of quasi-coherent subsheaves of finite type, we get the desired result by the pull-back (cf.\ {\bf \ref{ch-formal}}.\ref{prop-exeradicallyqcohsheavesbypullback}).

\medskip\noindent
{\bf Exercise \ref{exer-lemdeltasheafadicallyquasicoherent1-4general}.} 
Let $M\subseteq L$, where $L$ is a finitely generated $A$-module.
Consider the commutative diagram with exact rows (due to {\bf \ref{ch-formal}}.\ref{prop-lemdeltasheafadicallyquasicoherent1-3}):
$$
\xymatrix@-1ex{0\ar[r]&(M/N)^{\Delta}\ar[r]&(L/N)^{\Delta}\ar[r]&(L/M)^{\Delta}\ar[r]&0\\ 0\ar[r]&M^{\Delta}\ar[u]\ar[r]&L^{\Delta}\ar[u]\ar[r]&(L/M)^{\Delta}\ar@{=}[u]\ar[r]&0\rlap{.}}
$$
The desired exact sequence follows by snake lemma.

\medskip\noindent
{\bf Exercise \ref{exer-lemdeltasheafadicallyquasicoherent1-application}.}
We may assume that $X$ is affine $X=\Spf A$ where $A$ is a t.u.\ rigid-Noetherian ring with a finitely generated ideal of definition $I\subseteq A$.

(1) Consider for each $k\geq 0$ the induced exact sequence
$$
0\longrightarrow\mathscr{F}/\mathscr{F}\cap\mathscr{I}^{k+1}\mathscr{G}\longrightarrow\mathscr{G}/\mathscr{I}^{k+1}\mathscr{G}\longrightarrow\mathscr{H}/\mathscr{I}^{k+1}\mathscr{H}\longrightarrow 0,\eqno{(\ast)}
$$
where $\mathscr{I}=I^{\Delta}$.
Since each $\mathscr{F}/\mathscr{F}\cap\mathscr{I}^{k+1}\mathscr{G}$ is quasi-coherent, we deduce from {\bf \ref{ch-formal}}.\ref{lem-vanishingcohomologyadicallyuseful} that the sequence
$$
0\longrightarrow\Gamma(X,\mathscr{F})\longrightarrow\Gamma(X,\mathscr{G})\longrightarrow\Gamma(X,\mathscr{H})\longrightarrow 0
$$
is exact.
By {\bf \ref{ch-formal}}.\ref{prop-lemdeltasheafadicallyquasicoherent1-3} we have $\mathscr{F}=\Gamma(X,\mathscr{F})^{\Delta}$, which is adically quasi-coherent.

(2) Let $N=\Gamma(X,\mathscr{F})\hookrightarrow M=\Gamma(X,\mathscr{G})$; $M$ is finitely generated.
We have $N^{\Delta}=\mathscr{F}$ and $M^{\Delta}=\mathscr{G}$.
Hence the claim follows from {\bf \ref{ch-formal}}.\ref{prop-lemdeltasheafadicallyquasicoherent1-3}.

\medskip\noindent
{\bf Exercise \ref{exer-adicallyqcohlocallnoetherian}.} (1) See \cite[(10.10.2.7)]{EGAInew}.

(2) Use \cite[(10.10.2.9)]{EGAInew} and {\bf \ref{ch-formal}}.\ref{thm-adicqcohpre1}.

(3) It suffices to show that $\mathscr{J}$ is of finite type.
We may assume that $X$ is affine $X=\Spf A$ where $A$ is a Noetherian adic ring.
Set $J=\Gamma(X,\mathscr{J})$, which is an ideal of $A$, and hence is finitely generated.
We have the canonical morphism $J^{\Delta}\rightarrow\mathscr{J}$.
Then show that this is an isomorphism.

\medskip\noindent
{\bf Exercise \ref{exer-propadmissibleideal2}.} 
We may assume that $X$ is affine $X=\Spf A$, where $A$ is a t.u.\ adhesive ring with a finitely generated ideal of definition $I\subseteq A$ and $A$ is $I$-torsion free.
We look at the exact sequence
$$
0\longrightarrow\mathscr{J}\cap\mathscr{J}'\longrightarrow\mathscr{J}\oplus\mathscr{J}'\longrightarrow\mathscr{J}+\mathscr{J}'\longrightarrow 0,
$$
and compare it with the exact seqence
$$
0\longrightarrow J\cap J'\longrightarrow J\oplus J'\longrightarrow J+J'\longrightarrow 0,
$$
where $J=\Gamma(X,\mathscr{J})$ and $J'=\Gamma(X,\mathscr{J}')$, which are $I$-admissible ideals of $A$.
By Exercise \ref{exer-intersecfg} the intersection $J\cap J'$ is finitely generated.
Hence, applying the exact functor $(\ast)$ in {\bf \ref{ch-formal}}.\ref{thm-adicqcohpre1}, we get the first exact sequence from the second one.
This implies that $\mathscr{J}\cap\mathscr{J}'=(J\cap J')^{\Delta}$, which is an admissible ideal.

\medskip\noindent
{\bf Exercise \ref{exer-affinemorphism2idealofdefinition}.} We may assume that $X$ and $Y$ are affine $X=\Spf A$ and $Y=\Spf B$ and that the morphism $f\colon X\rightarrow Y$ comes from an adic map $B\rightarrow A$ between adic rings of finite ideal type.
Let $I\subseteq B$ be a finitely generated ideal of definition of $B$ such that $I^{\Delta}=\mathscr{I}$, and $M$ be the $A$-module such that $\mathscr{F}=M^{\Delta}$.
Then we have $f_{\ast}M^{\Delta}=M_{[B]}^{\Delta}$, where $M_{[B]}$ is the module $M$ regarded as a $B$-module ({\bf \ref{ch-formal}}.\ref{lem-exeradicallydeltasheavesunited}).
By {\bf \ref{ch-formal}}.\ref{lem-deltasheafadicallyquasicoherent1-2} we have $f_{\ast}\mathscr{I}^{k+1}M^{\Delta}=(I^{k+1}M)_{[B]}^{\Delta}=I^{k+1}M_{[B]}^{\Delta}=I^{k+1}f_{\ast}M^{\Delta}$, which shows the assertion for the module sheaf $\mathscr{F}$. 
The other case is similar.

\medskip\noindent
{\bf Exercise \ref{exer-closedimmformal61}.} All the assertion follows from {\bf \ref{ch-formal}}.\ref{prop-lemdeltasheafadicallyquasicoherent1-3}.

\medskip\noindent
{\bf Exercise \ref{exer-adicnessadicallyflatdescent}.} We consider the conormal cones 
$$
{\textstyle \gr^{\bullet}_I(A)=\bigoplus_{n\geq 0}I^n/I^{n+1}\quad\textrm{and}\quad\gr^{\bullet}_{IB}(B)=\bigoplus_{n\geq 0}I^nB/I^{n+1}B}
$$
(cf.\ {\bf \ref{ch-pre}}, \S\ref{sub-Reescone}).
Since $\gr^{\bullet}_{IB}(B)=\gr^{\bullet}_I(A)\otimes_AB=\gr^{\bullet}_I(A)\otimes_{A_0}B_0$, where $A_0=A/I$ and $B_0=B/IB$, and since $\Spf B\rightarrow\Spf A$ is adically faithfully flat, we deduce that the map $\gr^{\bullet}_I(A)\rightarrow\gr^{\bullet}_{IB}(B)$ is faithfully flat.
Now consider 
$$
{\textstyle \gr^{\bullet}_F(A)=\bigoplus_{n\geq 0}F^n/F^{n+1}\quad\textrm{and}\quad\gr^{\bullet}_{FB}(B)=\bigoplus_{n\geq 0}F^nB/F^{n+1}B,}
$$
which are a graded $\gr^{\bullet}_I(A)$-module and a graded $\gr^{\bullet}_{IB}(B)$-module, respectively.
Since we have $\gr^{\bullet}_{FB}(B)=\gr^{\bullet}_F(A)\otimes_{\gr^{\bullet}_I(A)}\gr^{\bullet}_{IB}(B)$, $\gr^{\bullet}_F(A)$ is finitely generated as a $\gr^{\bullet}_I(A)$-module if and only if $\gr^{\bullet}_{FB}(B)$ is finitely generated as a $\gr^{\bullet}_{IB}(B)$-module.
Hence the assertion follows from {\bf \ref{ch-pre}}.\ref{prop-complpair2}.

\medskip\noindent
{\bf Exercise \ref{exer-propermorformalschemex1}.} Set $X_k=(X,\O_X/I^{k+1}\O_X)$ and $Y_k=(Y,\O_Y/I^{k+1}\O_Y)$ for $k\geq 0$ (where $I\subseteq A$ is a finitely generated ideal of definition), and consider the induced map $f_k\colon X_k\rightarrow Y_k$.
Then $f$ is proper if and only if $f_0$ is proper. 
Then the claim follows from \cite[$\mathbf{II}$, (5.6.3)]{EGA} (where we first reduce to the case where $Y$ is affine, and use \cite[$\mathbf{IV}$, (8.10.5.1)]{EGA} instead of \cite[$\mathbf{II}$, (5.6.1)]{EGA}).

\medskip\noindent
{\bf Exercise \ref{exer-propetaleequivquote2}.} Use the technique as in the proof of {\bf \ref{ch-formal}}.\ref{thm-etaleequivquote}; cf.\ \cite[II.1.7]{Knu}.

\medskip\noindent
{\bf Exercise \ref{exer-idealofdefinitioneffectivedescent}.} We may assume that $Y$ has an ideal of definition of finite type $\mathscr{I}$. 
By {\bf \ref{ch-formal}}.\ref{prop-exeridealofdefinitioncompletepullback} we have $\widehat{f^{\ast}}\mathscr{I}^n=\mathscr{I}^n\O_X$ for $n\geq 1$. 
Since $\widehat{f^{\ast}}\mathscr{J}$ is an ideal of definition of $X$, there exists $n\geq 1$ such that $\widehat{f^{\ast}}\mathscr{I}^n\subseteq\widehat{f^{\ast}}\mathscr{J}\subseteq\O_X$.
By this, one has an a.q.c.\ ideal sheaf $\mathscr{H}$ containing $\mathscr{I}^n$ such that $\widehat{f^{\ast}}\mathscr{J}=\widehat{f^{\ast}}\mathscr{H}$, whence $\mathscr{H}=\mathscr{J}$.
Notice that $\mathscr{J}$ is an open ideal of $\O_Y$.
Hence, in particular, $\mathscr{J}^m$ for any $m\geq 1$ is again an open adically quasi-coherent ideal.
For a sufficiently large $m\geq 1$ we have $(\widehat{f^{\ast}}\mathscr{J})^m\subseteq\widehat{f^{\ast}}\mathscr{I}$.
Combined with the canonical morphism $\widehat{f^{\ast}}\mathscr{J}^m\rightarrow(\widehat{f^{\ast}}\mathscr{J})^m$, we conclude again by {\bf \ref{ch-formal}}.\ref{prop-etaledescent3} that $\mathscr{J}^m\subseteq\mathscr{I}$.

\medskip\noindent
{\bf Exercise \ref{exer-formalalgebraicspaceszariskicoveringqcpt}.} Confer \cite[II.3.13]{Knu}.

\medskip\noindent
{\bf Exercise \ref{exer-pointsformalalgebraicspaces2}.} Reduce to the situation where $X$ has an ideal of definition, and apply \cite[II.6.7]{Knu} and {\bf \ref{ch-formal}}.\ref{prop-localcriterionrepresentability}; cf.\ \cite[Premi\`ere partie, (5.7.7)]{RG}.

\medskip\noindent
{\bf Exercise \ref{exer-fini1}.} For any coherent sheaf $\mathscr{F}$ on $X$ there exists by {\bf \ref{ch-formal}}.\ref{prop-fini1} an integer $n$ such that $(f^{\ast}f_{\ast}\mathscr{F}(n))(-n)\rightarrow\mathscr{F}$ is surjective.
Let $M$ be the $B$-module such that $\til{M}=f_{\ast}\mathscr{F}(n)$.
Find a finitely generated $B$-submodule $N$ of $M$ such that the map $(f^{\ast}\til{N})(-n)\rightarrow\mathscr{F}$ is surjective.
Then replace $N$ by a free $B$-module of finite rank.
Prove in this way that any coherent sheaf on $X$ admits a resolution by coherent sheaves of the form $\O(n)^{\bigoplus m}$, which are projective objects in the category $\QCoh_X$.
Then use {\bf \ref{ch-pre}}.\ref{cor-subcategoryderived31} (2).


\medskip\noindent
{\bf Exercise \ref{exer-desalgsp}.} Take a quasi-compact representable \'etale covering $Y\rightarrow X$, and set $R=Y\times_XY$.
The projections $p_1,p_2\colon R\rightarrow Y$ give an \'etale equivalence relation (defining $X$) such that $(p_1,p_2)\colon R\hookrightarrow Y\times_AY$ is an immersion (resp.\ closed immersion).
Since $Y$ and $R$ are finitely presented over $A$, there exist a filtered family $\{A_{\lambda}\}_{\lambda\in\Lambda}$ of subrings of $A$ of finite type over $\Z$ and a projective system $\{p_{1,\lambda},p_{2,\lambda}\colon R_{\lambda}\rightarrow Y_{\lambda}\}_{\lambda\in\Lambda}$ of diagrams that converges to $p_1,p_2\colon R\rightarrow Y$ such that for each $\lambda\in\Lambda$ $R_{\lambda}$ and $Y_{\lambda}$ are finite type over $A_{\lambda}$.
Then by Exercise \ref{exer-equivalencerelationinvolution} $p_{1,\lambda},p_{2,\lambda}\colon R_{\lambda}\rightarrow Y_{\lambda}$ defines an \'etale equivalence relation for sufficiently large $\lambda\in\Lambda$.

\medskip\noindent
{\bf Exercise \ref{exer-formalcompletionfaithflness}.} The `if' part is clear. To show the converse, by a reduction argument similar to that in \cite[$\mathbf{III}$, (4.6.8)]{EGA} we may assume that $Y=Z$ and that it is affine $Y=Z=\Spec A$.
Let $I\subseteq A$ be the finitely generated ideal that defines $W\subseteq Z$.
Since $f^{-1}(y)$ is a singleton set for any $y\in W$, there exists an open neighborhood $U$ of $W$ such that $f^{-1}(U)\rightarrow U$ is quasi-finite (this follows from \cite[$\mathbf{III}$, (4.4.11)]{EGA} applied with the standard limit argument (cf.\ \cite[$\mathbf{IV}$, \S8.10]{EGA})).
Hence we may assume $f$ is quasi-finite; since $f$ is proper of finite presentation, one deduces by Zariski's Main Theorem\index{Zariski, O.}\index{Zariski's Main Theorem} (\cite[$\mathbf{IV}$, (8.12.6)]{EGA}) that $f$ is finite.
Hence $X=\Spec B$, where $B$ is a finite $A$-algebra.
Now since $A$ is $I$-adically universally adhesive (and so is $B$), the map $\zat{A}\rightarrow\widehat{A}$ is faithfully flat (and so is $\zat{B}\rightarrow\widehat{B}$), where $\zat{A}$ is the associated Zariskian ring\index{Zariskian!associated Zariskian@associated ---} ({\bf \ref{ch-pre}}.\ref{prop-relpair21} (2) and {\bf \ref{ch-pre}}.\ref{prop-btarf1} (2)).
Moreover, by {\bf \ref{ch-pre}}.\ref{prop-btarf1} (1) we know that $\widehat{B}=\zat{B}\otimes_{\zat{A}}\widehat{A}$.
Hence the map $\widehat{A}\rightarrow\widehat{B}$ is an isomorphism (resp.\ surjective) if and only if $\zat{A}\rightarrow\zat{B}$ is an isomorphism (resp.\ surjective).
Since $\zat{A}$ is the inductive limit of the rings of the form $A_{(1+a)}$, where $a\in I$, we have the desired result.

\medskip\noindent
{\bf Exercise \ref{exer-weakisomcomplete}.} It suffices to show the assertion in the following two cases:
\begin{itemize}
\item[(a)] $f$ is injective;
\item[(b)] $f$ is surjective.
\end{itemize}

Case (a). Since $M/N$ is bounded $I$-torsion, we have $I^nM\subseteq N$ for a sufficiently large $n\gg 0$.
This implies that the subspace topology on $N$ induced from the $I$-adic topology on $M$ is the $I$-adic topology, since $I^{n+k}N\subseteq I^{n+k}M=I^{n+k}M\cap N\subseteq I^kN$ for any $k\geq 0$.
Hence by {\bf \ref{ch-pre}}.\ref{cor-ML4} and {\bf \ref{ch-pre}}.\ref{lem-bt1vis} the sequence
$$
0\longrightarrow\widehat{N}\longrightarrow\widehat{M}\longrightarrow M/N\longrightarrow 0
$$
is exact, where $\widehat{\,\cdot\,}$ denotes the $I$-adic completion.
The assertion follows from this.

Case (b). In this case, $K=\ker(f)$ is bounded $I$-torsion.
Since the subspace topology on $K$ induced from the $I$-adic topology on $N$ is the $I$-adic topology, we have the exact sequence 
$$
0\longrightarrow K\longrightarrow\widehat{N}\longrightarrow\widehat{M}\longrightarrow 0
$$
by the similar reasoning as before.

\medskip\noindent
{\bf Exercise \ref{exer-FPapproximationfunctoriality}.} Take any FP-approximations $\alpha\colon\mathscr{F}'\rightarrow\mathscr{F}$ and $\beta'\colon\mathscr{G}''\rightarrow\mathscr{G}$, and consider the fiber product $\mathscr{K}$ of the maps $\varphi\circ\alpha$ and $\beta'$, which is an adically quasi-coherent subsheaf of the direct sum $\mathscr{F}'\oplus\mathscr{G}''$.
Then mimic the argument as in the proof of {\bf \ref{ch-rigid}}.\ref{lem-categoryFPthickeningfiltered} (use Exercise \ref{exer-extopenadiccoh1approx} in the formal scheme case).

\section*{Chapter \ref{ch-rigid}}

\medskip\noindent
{\bf Exercise \ref{exer-propblowups4a}.} Let $\mathscr{J}$ be an admissible ideal that gives the admissible blow-up $\pi\colon X'\rightarrow X$.
Consider the admissible ideals $\mathscr{J}\O_Y$ on $Y$ and $\mathscr{J}\O_{Y'}$ on $Y'$, and let $Z\rightarrow Y$ and $Z'\rightarrow Y'$ the respective admissible blow-ups.
We want to show that $Z$ and $Z'$ are isomorphic.
By {\bf \ref{ch-rigid}}.\ref{prop-blowups1} (3) there exists $Z'\rightarrow Z$ that makes the resulting square commutative.
Similarly, by {\bf \ref{ch-rigid}}.\ref{prop-blowups1} (3) we have $Z\rightarrow X'$, which induces $Z\rightarrow Y'$.
Then one sees again by {\bf \ref{ch-rigid}}.\ref{prop-blowups1} (3) that there exists an arrow $Z\rightarrow Z'$.
It is then easy to see that these arrow $Z\rightarrow Z'$ and $Z'\rightarrow Z$ are inverse to each other.

\medskip\noindent
{\bf Exercise \ref{exer-Jtorsionpartquoteaqc}.} Consider $U=\Spf A_{\{f\}}\subseteq X$ for any $f\in A$.
The module $\Gamma(U,\mathscr{F}_{\textrm{$\mathscr{J}$-}\mathrm{tor}})$ is the $J$-torsion part of $\Gamma(U,\mathscr{F})=M\widehat{\otimes}_AA_{\{f\}}=M\otimes_AA_{\{f\}}$ ({\bf \ref{ch-pre}}.\ref{prop-btarf1} (1)). Since $A_{\{f\}}$ is flat over $A$ ({\bf \ref{ch-pre}}.\ref{prop-btarf1} (2)), we have $\Gamma(U,\mathscr{F}_{\textrm{$\mathscr{J}$-}\mathrm{tor}})=M_{\Jtor}\otimes_AA_{\{f\}}$.
On the other hand, since $J$ is $I$-admissible, we have $M_{\Jtor}\subseteq M_{\Itor}$, and $M_{\Jtor}\otimes_AA_{\{f\}}$ is bounded $I$-torsion.
Then by {\bf \ref{ch-pre}}.\ref{lem-bt1} and {\bf \ref{ch-formal}}.\ref{prop-adicqcohaff1g} we have $\Gamma(U,(M_{\Jtor})^{\Delta})=M_{\Jtor}\otimes_AA_{\{f\}}$.
Hence we have $\mathscr{F}_{\textrm{$\mathscr{J}$-}\mathrm{tor}}=(M_{\Jtor})^{\Delta}$, as desired.

\medskip\noindent
{\bf Exercise \ref{exer-Jtorsionpartquoteaqc2}.} We may work in the affine situation $X=\Spf A$ where $A$ is a t.u.\ rigid-Noetherian ring.
Then the claim follows from Exercise \ref{exer-Jtorsionpartquoteaqc} and {\bf \ref{ch-formal}}.\ref{prop-lemdeltasheafadicallyquasicoherent1-3}.

\medskip\noindent
{\bf Exercise \ref{exer-secondproofpropstricttransformaqctoaqcrigidNoetherian}.} By {\bf \ref{ch-formal}}.\ref{prop-completepullbackaqcsheaves0} we know that $\pi^{\ast}\mathscr{F}$ is an adically quasi-coherent sheaf of finite type on $X'$.
By {\bf \ref{ch-formal}}.\ref{prop-admissibleideal1x} the sheaf $\mathscr{J}\O_{X'}$ is an admissible ideal of $\O_{X'}$.
Hence the assertion follows from Exercise \ref{exer-Jtorsionpartquoteaqc2}.

\medskip\noindent
{\bf Exercise \ref{exer-admissibleidealstricttransform}.} Since $B$ is $J$-torsion free, we have $K=(J\otimes_AB)_{\Jtor}$.
In particular, we can deduce similarly as in the hint for Exercise \ref{exer-Jtorsionpartquoteaqc} that $K$ is a bounded $I$-torsion module and hence is $I$-adically complete.
By this, by an argument similar to that in the hint for Exercise \ref{exer-Jtorsionpartquoteaqc}, one can show that $\mathscr{K}=K^{\Delta}$, whence (1).
Then by Exercise \ref{exer-Jtorsionpartquoteaqc} one has (2).

\medskip\noindent
{\bf Exercise \ref{exer-adequateformaltorigid}.} For a quasi-separated $X$, cover $X$ by coherent open subsets, and then construct $X^{\rig}$ as a stretch of coherent rigid spaces.
In general, cover $X$ by coherent open subsets (e.g.\ affine), and mimic the construction as in {\bf \ref{ch-rigid}}.\ref{dfn-generalrigidspace1}.

\medskip\noindent
{\bf Exercise \ref{exer-deligne}.} Mimic the proof of \cite[(6.9.17)]{EGAInew} with the following modifications: 
\begin{itemize}
\item consider the modules $R$ and $S$ as in \cite[p.\ 323]{EGAInew}; it follows that $R/S$ is $\mathfrak{I}$-torsion; since $R/S$ is the submodule in the finitely generated $A$-module $A^{\oplus mn}/S$, the $\mathfrak{I}$-torsion is bounded;
\item use {\bf (AP)} instead of Artin-Rees lemma in \cite[p.\ 324]{EGAInew}.
\end{itemize}

\medskip\noindent
{\bf Exercise \ref{exer-propZRstrsheaf2}.} We may assume that $\mathscr{X}$ is a coherent rigid space.
Then there exists a formal model $X$ of $\mathscr{X}$ on which there exist ideals of definition $\mathscr{I}_X$ and $\mathscr{I}'_X$ such that $\mathscr{I}=(\sp^{-1}_X\mathscr{I}_X)\O^{\int}_{\mathscr{X}}$ and $\mathscr{I}'=(\sp^{-1}_X\mathscr{I}'_X)\O^{\int}_{\mathscr{X}}$.
Then the assertion follows from the fact that there exist $n,m>0$ such that $\mathscr{I}^m_X\subseteq\mathscr{I}^{\prime n}_X\subseteq\mathscr{I}_X$.

\medskip\noindent
{\bf Exercise \ref{exer-proptypeN}.} If $\mathscr{X}$ has a Noetherian formal model $X$, then any admissible blow-up $X'$ of $X$ is again Noetherian, since $X'$ is of finite type over $X$ (cf.\ \cite[$\mathbf{I}$, (10.13.2)]{EGA}), whence (a) $\Rightarrow$ (b).
Suppose (b) holds, and let $\mathscr{U}\subseteq\mathscr{X}$ be a quasi-compact open subspace.
Then there exists a Noetherian formal model $X$ of $\mathscr{X}$ having a quasi-compact open subset $U$ that corresponds to $\mathscr{U}$ ({\bf \ref{ch-rigid}}.\ref{prop-zariskiriemanntoptop}), whence (b) $\Rightarrow$ (c).
The implication (c) $\Rightarrow$ (d) is clear.
To show (d) $\Rightarrow$ (a), take a finite open covering $\{\mathscr{U}_{\alpha}\}_{\alpha\in L}$ as in (d).
By {\bf \ref{ch-rigid}}.\ref{prop-zariskiriemanntoptop} there exists a formal model $X$ and an open covering $\{U_{\alpha}\}_{\alpha\in L}$ of $X$ such that each $U_{\alpha}$ gives a formal model of $\mathscr{U}_{\alpha}$ (cf.\ {\bf \ref{ch-rigid}}, \S\ref{sub-comparisontopoi}).
Applying (a) $\Rightarrow$ (b) to each $U_{\alpha}$ and replacing $X$ by an admissible blow-up if necessary, we may assume that each $U_{\alpha}$ is Noetherian.
Then $X$ is Noetherian, as desired.

\medskip\noindent
{\bf Exercise \ref{exer-tubesubsetscoherentmaps}.} Cf.\ {\bf \ref{ch-pre}}.\ref{cor-valuativemapstubes2}.

\medskip\noindent
{\bf Exercise \ref{exer-corlatticemodels41}.} Replacing $X$ by admissible blow-ups and the sheaves by their strict transforms, one can assume that $\mathscr{I}^n\mathscr{F}_X|_U\subseteq\mathscr{I}^m\mathscr{G}_X|_U\subseteq\mathscr{F}_X|_U$ holds for positive integers $n,m$.
Then since $\mathscr{F}_X+\mathscr{G}_X$ is obviously adically quasi-coherent of finite type that is $\mathscr{I}_X$-torsion free, it gives a lattice model.
By the exact sequence 
$$
0\longrightarrow\mathscr{F}_X\cap\mathscr{G}_X\longrightarrow\mathscr{F}_X\oplus\mathscr{G}_X\longrightarrow\mathscr{F}_X+\mathscr{G}_X\longrightarrow 0,
$$
$\mathscr{F}_X\cap\mathscr{G}_X$ is of finite type, since $\mathscr{F}_X+\mathscr{G}_X$ if finitely presented.
Hence the assertion for $\mathscr{F}_X\cap\mathscr{G}_X$ follows.

\medskip\noindent
{\bf Exercise \ref{exer-corexistlatticelocfree}.} We can start from the following situation (similarly to that in {\sc Step 3} of the proof of {\bf \ref{ch-rigid}}.\ref{lem-existlattice1}): 
$X$ is a distinguished formal model of $\mathscr{X}$ having an invertible ideal of definition $\mathscr{I}_X$, $\ZR{\mathscr{X}}=\bigcup^n_{\alpha=1}\mathfrak{U}_{\alpha}$ is a finite open covering, and for each $\alpha$ we have
\begin{itemize}
\item an admissible blow-up $X_{\alpha}\rightarrow X$ and a quasi-compact open subset $U_{\alpha}$ of $X_{\alpha}$ such that $\mathfrak{U}_{\alpha}=\sp^{-1}_{X_{\alpha}}(U_{\alpha})$; 
\item a positive integer $s$ (independent on $\alpha$); 
\item a formal model $\mathscr{F}_{X_{\alpha}}$ of $\mathscr{F}$ on $X_{\alpha}$ and a weak isomorphism $\til{\varphi}_{\alpha}\colon\mathscr{I}^s_X\O^{\oplus r}_X|_{U_{\alpha}}\rightarrow\mathscr{F}_{X_{\alpha}}|_{U_{\alpha}}$ such that $\til{\varphi}_{\alpha}^{\rig}=\varphi|_{\mathfrak{U}_{\alpha}}$.
\end{itemize}
Notice that, since $\mathscr{I}_X$ is invertible, each $\mathscr{I}^s_X\O^{\oplus r}_X|_{U_{\alpha}}$ is locally free.
By an argument similar to that in {\sc Step 3} of the proof of {\bf \ref{ch-rigid}}.\ref{lem-existlattice1}, one finds a further admissible blow-up on which the strict transforms of $\mathscr{I}^s_X\O^{\oplus r}_X|_{U_{\alpha}}$ and $\mathscr{F}_{X_{\alpha}}|_{U_{\alpha}}$ glue together and gives rise to a weak isomorphism between the resulting lattice models.
Taking the admissible blow-up to be a distinguished formal model of $\mathscr{X}$, one sees that the resulting gluing of the strict transforms of $\mathscr{I}^s_X\O^{\oplus r}_X|_{U_{\alpha}}$ is locally free.

\medskip\noindent
{\bf Exercise \ref{exer-weierstrasslaurentrational}.} See \cite[(7.2.3/7)]{BGR}.

\medskip\noindent
{\bf Exercise \ref{exer-cormorbetweenaffinoid1}.} Let $A'=\Gamma(X',\O_{X'})$, and consider the morphism $q\colon X''=\Spf A'\rightarrow\Spf A$ induced from the strict weak isomorphism $A\rightarrow A'$ ({\bf \ref{ch-rigid}}.\ref{prop-affinoidblowupglobalsection}).
Then we have $\pi_{\ast}\O_{X'}=q_{\ast}\O_{X''}$. 
Since $q$ is affine, $q_{\ast}\O_{X''}$ is adically quasi-coherent by {\bf \ref{ch-formal}}.\ref{cor-effectivedescaffineadic}.

\medskip\noindent
{\bf Exercise \ref{exer-associatedschememapsurjective}.} Set $\mathscr{X}=(\Spf A)^{\rig}$.
For any closed point $x\in s(\mathscr{X})=\Spec A\setminus V(I)$, take a map $\Spec V\rightarrow\Spec A$ from the spectrum of a valuation ring such that the generic point is mapped to $x$ and that the closed points is mapped in $V(I)$.
Then we have a rigid point $\Spf\widehat{V}\rightarrow\ZR{\mathscr{X}}$.

\medskip\noindent
{\bf Exercise \ref{exer-propvalutiverigid1}.} Easy by {\bf \ref{ch-rigid}}.\ref{cor-seppropmorrigid22}.

\medskip\noindent
{\bf Exercise \ref{exer-GAGAfunctor1}.} Take an object $(X\hookrightarrow\ovl{X})$ of $\Emb_{X|S}$, and let $Z$ and $\ovl{Z}$ be as in {\bf \ref{ch-rigid}}.\ref{const-GAGAfunctor}.
Then by {\bf \ref{ch-rigid}}.\ref{prop-GAGAfunctor5x2} we have $X^{\an}=(\ovl{X}_U)^{\an}\setminus Z^{\an}$.
Let $\mathscr{J}$ be the ideal defining $\ovl{Z}$ in $\ovl{X}$, which we may assume to be invertible; we may moreover assume $I\O_{\ovl{X}}$ is invertible.
Let $\ovl{X}_n$ for each $n\geq 1$ be the blow-up of $\ovl{X}$ along the ideal $\mathscr{J}^n+I\O_{\ovl{X}}$, and $\til{X}_n$ the maximal open set of $\ovl{X}_n$ where $\mathscr{J}^n\O_{\ovl{X}_n}$ generates the ideal $\mathscr{J}^n\O_{\ovl{X}_n}+I\O_{\ovl{X}_n}$.
Then $\til{X}_n\rightarrow\ovl{X}$ is affine for any $n\geq 1$ and we have $\varinjlim_n(\widehat{\til{X}}_n)^{\rig}=X^{\an}$.

\medskip\noindent
{\bf Exercise \ref{exer-vsbr-gradedtrivial}.}
Only (a) $\Rightarrow$ (c) calls for a hint. 
Suppose $\mathrm{Gr}_FA=\{0\}$.
Then one can show $F=F^+$.
Consider the seminorm $\nu$ asssociated to $F$.
Since $1\in F_{<1}=F^+_{<1}$, $\nu(1)<1$, which implies (c).

\medskip\noindent
{\bf Exercise \ref{exer-vsbr-R+valuations}.}
Let $F_{\nu}$ be the filtration induced from $\nu$.
Then $(K,F_{\nu})$ is a filtered valuation field of maximal type, and we have $V\subseteq(F_{\nu})_1$.
Consider $\ovl{V}=V/(F_{\nu})_{<1}$, which is a subring of $\mathrm{Gr}_{F_{\nu},1}K$.
We regard $\ovl{V}$ as a graded local subring of $\mathrm{Gr}_{F_{\nu}}K$ by $0$-extension (cf.\ {\bf \ref{ch-rigid}}, \S\ref{subsub-vsbr-gradedrings}).
By {\bf \ref{ch-rigid}}.\ref{prop-vsbr-gradedvaluationrings2}, there exists a graded valuation subring of $\mathrm{Gr}_{F_{\nu}}K$, whose unit-element part coincides with $\ovl{V}$.
Hence by {\bf \ref{ch-rigid}}.\ref{prop-vsbr-filteredvaluationfield1}, we have a filtration $F$ on $K$ such that $(K,F)$ is a filtered valuation field with $F_1=V$.

\medskip\noindent
{\bf Exercise \ref{exer-vsbr-spectralseminormformula}.}
Reduce to the case where $\mathcal{A}$ is a Banach field by the proper continuous mapping $\mathscr{M}(\mathcal{A}\dl r^{-1}T\dr)\rightarrow\mathscr{M}(\mathcal{A})$.

\medskip\noindent
{\bf Exercise \ref{exer-vsbr-valuativespectrum}.}
(3) Consider $L=K^{\mathrm{st}}$ as in \ref{exa-vsbr-Kstandard}.
We first show that $\Spec_{\R_+}\mathrm{Gr}\,\mathcal{A}_L\rightarrow\Spec_{\R_+}\mathrm{Gr}\,\mathcal{A}$ is surjective.
This follows from the surjectivity of $\sp_{\mathcal{A}}$ and of $\Spec^{\mathrm{val}}\mathcal{A}_L\rightarrow\Spec^{\mathrm{val}}\mathcal{A}$.
Then the claim reduces to the strict case \ref{rem-reductionscheme2}.

\medskip\noindent
{\bf Exercise \ref{exer-nagatachow}.} Take a proper $Y$-scheme $\ovl{U}$ that contains $U$ as an open subset.
Then by {\bf \ref{ch-rigid}}.\ref{prop-correspondencediagramfurther} (2), replacing $\ovl{U}$ by a $U$-admissible blow-up if necessary, we can find a $U$-admissible blow-up $X'$ of $X$ that admits an open immersion $X'\hookrightarrow\ovl{U}$.
Since $\ovl{U}$ is a scheme, $X'$ is a scheme.

\addcontentsline{toc}{chapter}{\rm List of Notations}
\markboth{List of Notations}{List of Notations}
\chapter*{List of Notations}
\section*{Categories}
\newcommand{\categorylist}[4]{\noindent
\begin{minipage}{6.5em}#1\end{minipage}
\begin{minipage}{12.3em}#2\end{minipage}\hspace{.4em}
\begin{minipage}{14em}#3\end{minipage}\hspace{.4em}
\begin{minipage}{6em}#4\end{minipage}}

\subsection*{--- Sets and Spaces}\ 

{\small 
\categorylist{{\sl Name}}{{\sl Objects}}{{\sl Morphisms}}{{\sl See}}

\categorylist{$\Sets$}{Sets}{Maps}{{\bf \ref{ch-pre}}, \S\ref{subsub-categorynot}}

\categorylist{$\Top$}{Topological spaces}{Continuous maps}{{\bf \ref{ch-pre}}, \S\ref{subsub-categorynot}}

\categorylist{$\STop$}{Sober top.\ spaces}{Continuous maps}{{\bf \ref{ch-pre}}, \S\ref{subsub-sober}}

\categorylist{$\CSTop$}{Coherent sober top.\ spaces}{Quasi-compact maps}{{\bf \ref{ch-pre}}.\ref{thm-stonerepresentationthm}}

\categorylist{$\Vsp$}{Valuative spaces}{Valuative loc.\ q-cpt maps}{{\bf \ref{ch-pre}}.\ref{subsub-reflexivization}}

\categorylist{$\RVsp$}{Reflexive ---}{Valuative loc.\ q-cpt maps}{{\bf \ref{ch-pre}}.\ref{subsub-reflexivization}}

\categorylist{$\Rsp$}{Ringed spaces}{Morphisms of ringed spaces}{{\bf \ref{ch-pre}}, \S\ref{subsub-ringedsplocalringedsp}}

\categorylist{$\LRsp$}{Locally ringed spaces}{Local morphisms}{{\bf \ref{ch-pre}}, \S\ref{subsub-ringedsplocalringedsp}}

\bigskip

\categorylist{$\Triples$}{Triples}{Morphism of triples}{{\bf \ref{ch-rigid}}.\ref{dfn-triplescategory}}

\categorylist{$\VTriples$}{Valued triples}{Morphism of valued triples}{{\bf \ref{ch-rigid}}.\ref{dfn-Hubertriples}}

\categorylist{$\AnTriples$}{Analytic triples}{Morphism of triples}{{\bf \ref{ch-rigid}}.\ref{dfn-analytictriples}}

\bigskip

\categorylist{$\Sch$}{Schemes}{Morphisms of schemes}{{\bf \ref{ch-pre}}, \S\ref{subsub-schemesbasics}}

\categorylist{$\Sch_S$}{$S$-schemes}{$S$-morphisms}{{\bf \ref{ch-pre}}, \S\ref{subsub-schemesbasics}}

\categorylist{$\As$}{Algebraic spaces}{Morphisms of algebraic spaces}{{\bf \ref{ch-pre}}, \S\ref{subsub-algebraicspacesconv}}

\categorylist{$\CAs$}{Coherent ---}{Morphisms of algebraic spaces}{{\bf \ref{ch-rigid}}, \S\ref{subsub-birationalgeomterm}}

\categorylist{$\As_S$}{$S$-algebraic spaces}{$S$-morphisms}{{\bf \ref{ch-pre}}, \S\ref{subsub-algebraicspacesconv}}

\bigskip

\categorylist{$\Zs$}{Zariskian schemes}{Morphisms of Zariskian schemes}{{\bf \ref{ch-formal}}, \S\ref{subsub-zariskianschemes}}

\categorylist{$\AZs$}{Affine ---}{Morphisms of Zariskian schemes}{{\bf \ref{ch-formal}}, \S\ref{subsub-zariskianschemes}}

\categorylist{$\CZs$}{Coherent ---}{Morphisms of Zariskian schemes}{{\bf \ref{ch-formal}}, \S\ref{subsub-zariskianschemes}}

\bigskip

\categorylist{$\Fs$}{Formal schemes}{Morphisms of formal schemes}{{\bf \ref{ch-formal}}, \S\ref{subsub-categoryformalschnotation}}

\categorylist{$\Ac\Fs$}{Adic --- of finite ideal type}{Morphisms of formal schemes}{{\bf \ref{ch-formal}}, \S\ref{subsub-categoryformalschnotation}}

\categorylist{$\Ac\Fs^{\ast}$}{Adic --- of finite ideal type}{Adic morphisms}{{\bf \ref{ch-formal}}, \S\ref{subsub-categoryformalschnotation}}

\categorylist{$\Ac\Fs^{\ast}_{/S}$}{--- adic over $S$}{Adic morphisms}{{\bf \ref{ch-formal}}, \S\ref{subsub-categoryformalschnotation}}

\categorylist{$\RNoe\Fs$}{Univ.\ rigid-Noetherian ---}{Morphisms of formal schemes}{{\bf \ref{ch-formal}}, \S\ref{subsub-ntnadequatecategory}}

\categorylist{$\RNoe\Fs^{\ast}$}{Univ.\ rigid-Noetherian ---}{Adic morphisms}{{\bf \ref{ch-formal}}, \S\ref{subsub-ntnadequatecategory}}

\categorylist{$\RNoe\Fs^{\ast}_{/S}$}{--- adic over $S$}{Adic morphisms}{{\bf \ref{ch-formal}}, \S\ref{subsub-ntnadequatecategory}}

\categorylist{$\Adh\Fs$}{Universally adhesive ---}{Morphisms of formal schemes}{{\bf \ref{ch-formal}}, \S\ref{subsub-ntnadequatecategory}}

\categorylist{$\Adh\Fs^{\ast}$}{Universally adhesive ---}{Adic morphisms}{{\bf \ref{ch-formal}}, \S\ref{subsub-ntnadequatecategory}}

\categorylist{$\Adh\Fs^{\ast}_{/S}$}{--- adic over $S$}{Adic morphisms}{{\bf \ref{ch-formal}}, \S\ref{subsub-ntnadequatecategory}}

\categorylist{$\CFs$}{Coherent formal schemes}{Morphisms of formal schemes}{{\bf \ref{ch-formal}}, \S\ref{subsub-notationcategorycoherentformaschemes}}

\categorylist{$\Ac\CFs$}{Adic --- of finite ideal type}{Morphisms of formal schemes}{{\bf \ref{ch-formal}}, \S\ref{subsub-notationcategorycoherentformaschemes}}

\categorylist{$\Ac\CFs^{\ast}$}{Adic --- of finite ideal type}{Adic morphisms}{{\bf \ref{ch-formal}}, \S\ref{subsub-notationcategorycoherentformaschemes}}

\categorylist{$\Ac\CFs^{\ast}_{/S}$}{--- adic over $S$}{Adic morphisms}{{\bf \ref{ch-formal}}, \S\ref{subsub-notationcategorycoherentformaschemes}}

\categorylist{$\RNoe\CFs$}{Univ.\ rigid-Noetherian ---}{Morphisms of formal schemes}{{\bf \ref{ch-formal}}, \S\ref{subsub-ntnadequatecategory}}

\categorylist{$\Adh\CFs$}{Universally adhesive ---}{Morphisms of formal schemes}{{\bf \ref{ch-formal}}, \S\ref{subsub-ntnadequatecategory}}

\categorylist{$\Af\Fs$}{Affine formal schemes}{Morphisms of formal schemes}{{\bf \ref{ch-formal}}, \S\ref{subsub-categoryformalschnotation}}

\categorylist{$\Ac\FAs_S$}{Formal algebraic spaces}{Morphisms of sheaves}{{\bf \ref{ch-formal}}, \S\ref{subsub-formalalgebraicspacesdef}}

\categorylist{$\Ac\CFAs_S$}{Coherent ---}{Morphisms of sheaves}{{\bf \ref{ch-formal}}, \S\ref{subsub-formalalgebraicspacesdef}}

\bigskip

\categorylist{$\Rf$}{Rigid spaces}{Morphisms of rigid spaces}{{\bf \ref{ch-rigid}}, \S\ref{subsub-generalrigidspace}}

\categorylist{$\ARf$}{Affinoids}{Morphisms of rigid spaces}{{\bf \ref{ch-rigid}}, \S\ref{subsub-affinoidsintro}}

\categorylist{$\CRf$}{Coherent rigid spaces}{Morphisms of rigid spaces}{{\bf \ref{ch-rigid}}.\ref{dfn-cohrigidspace1}}

\bigskip

\categorylist{$\Acsp$}{Adic spaces}{Morphism of adic spaces}{{\bf \ref{ch-rigid}}.\ref{dfn-adicspacesbytriples}}

\categorylist{$\AnAcsp$}{Analytic adic spaces}{Morphism of adic spaces}{{\bf \ref{ch-rigid}}.\ref{dfn-analyticadicspaces}}

\bigskip

\categorylist{$\mathbf{Bsp}_K$}{Berkovich $K$-analytic spaces}{Morphisms of $K$-analytic spaces}{{\bf \ref{ch-rigid}}.\ref{dfn-berkovichanalyticspacesdefinition} (5)}

\categorylist{$\mathbf{Bsp}^{\ast}_K$}{Strictly ---}{Morphisms of $K$-analytic spaces}{{\bf \ref{ch-rigid}}.\ref{dfn-berkovichanalyticspacesdefinition} (6)}

}

\subsection*{--- Sheaves}

{\small 
\categorylist{{\sl Name}}{{\sl Objects}}{{\sl Morphisms}}{{\sl See}}

\categorylist{$\ASh_X$}{Abelian sheaves on $X$}{Sheaf homomorphisms}{{\bf \ref{ch-pre}}, \S\ref{subsub-canflasqueres}}

\categorylist{$\Coh_X$}{Coherent sheaves}{$\O_X$-linear morphisms}{{\bf \ref{ch-pre}}, \S\ref{subsub-ringedspsheavesmodules}}

\categorylist{$\Mod_X$}{$\O_X$-modules}{$\O_X$-linear morphisms}{{\bf \ref{ch-pre}}, \S\ref{subsub-ringedspsheavesmodules}}

\categorylist{$\QCoh_X$}{Quasi-coherent sheaves}{$\O_X$-linear morphisms}{{\bf \ref{ch-pre}}, \S\ref{subsub-ringedspsheavesmodules}}

\categorylist{$\AQCoh_X$}{Adically quasi-coh.\ sheaves}{$\O_X$-linear morphisms}{{\bf \ref{ch-formal}}, \S\ref{subsub-adicallyqcohdefdef}}

}

\subsection*{--- Algebras, modules, etc.}

{\small 
\categorylist{{\sl Name}}{{\sl Objects}}{{\sl Morphisms}}{{\sl See}}

\categorylist{$\Ab$}{Abelian groups}{Group homomorphisms}{{\bf \ref{ch-pre}}, \S\ref{subsub-categorynot}}

\categorylist{$\Alg_A$}{$A$-algebras}{$A$-algebra homomorphisms}{{\bf \ref{ch-pre}}, \S\ref{subsub-categorynot}}

\categorylist{$\Coh_A$}{$A$-modules}{$A$-module homomorphisms}{{\bf \ref{ch-pre}}, \S\ref{subsub-categorynot}}

\categorylist{$\DLat$}{Distributive lattices}{Lattice homomorphisms}{{\bf \ref{ch-pre}}, \S\ref{subsub-strcohtopsp}}

\categorylist{$\Mod_A$}{$A$-modules}{$A$-module homomorphisms}{{\bf \ref{ch-pre}}, \S\ref{subsub-categorynot}}

\categorylist{$\Coh_A$}{Coherent ---}{$A$-module homomorphisms}{{\bf \ref{ch-pre}}. \ref{dfn-cohringsmodules1}}

}

\subsection*{--- Others}

{\small 
\categorylist{{\sl Name}}{{\sl Objects}}{{\sl Morphisms}}{{\sl See}}

\categorylist{$\BL_X$}{admissible blow-ups of $X$}{$X$-morphisms}{{\bf \ref{ch-rigid}}, \S\ref{sub-categoryadmblow-up}}

\categorylist{$\Emb_{X|S}$}{Nagata embeddings of $X/S$}{$X$-admissible $S$-modifications}{{\bf \ref{ch-rigid}}, \S\ref{subsub-GAGAfunctoremb}}

\categorylist{$\MD_{(X,U)}$}{$U$-admissible modifications}{$X$-morphisms}{{\bf \ref{ch-rigid}}, \S\ref{subsub-filtcatmodifications}}

}

\section*{Other symbols}
\newcommand{\symbollist}[3]{\noindent
\begin{minipage}{7em}#1\end{minipage}
\begin{minipage}{26.2em}#2\end{minipage}
\begin{minipage}{6em}#3\end{minipage}}

\subsection*{--- Operator symbols}

{\small 
\symbollist{{\sl Symbol}}{{\sl Meaning}}{{\sl See}}

\symbollist{$\Ac^{\ast}(\mathscr{A})$}{Full subcategory of acyclic objects}{{\bf \ref{ch-pre}}, \S\ref{sub-derivedcategory}}

\symbollist{$\Aff(\mathscr{X})$}{Rigid analytic ring associated to $\mathscr{X}$}{{\bf \ref{ch-rigid}}.\ref{dfn-rigidanalyticring}}

\symbollist{$\AId_X$}{Set of admissible ideals}{{\bf \ref{ch-formal}}.\ref{dfn-admissibleideal}}

\symbollist{$\AId_{(X,U)}$}{Set of $U$-admissible ideals}{{\bf \ref{ch-rigid}}, \S\ref{subsub-birationalgeomblowups}}

\symbollist{$\amp(F)$}{Amplitude of $F$}{{\bf \ref{ch-pre}}.\ref{dfn-derivedcategory6}}

\symbollist{$\Aut_{\mathscr{C}}(x)$}{Automorphism group}{{\bf \ref{ch-pre}}, \S\ref{subsub-categoryconv}}

\symbollist{$\mathscr{C}^{\bullet}(X,\,\cdot\,)$}{Canonical flasque resolution}{{\bf \ref{ch-pre}}, \S\ref{subsub-canflasqueres}}

\symbollist{$\CC^{\ast}(\mathscr{A})$}{Category of complexes of objects in $\mathscr{A}$}{{\bf \ref{ch-pre}}, \S\ref{subsub-complexcategorydef}}

\symbollist{$\card(x)$}{Cardinality of $x$}{{\bf \ref{ch-pre}}, \S\ref{subsub-existuniverse}}

\symbollist{$\Cart_{/I}(E,F)$}{Category of cartesian functors from $E$ to $F$}{{\bf \ref{ch-pre}}, \S\ref{subsub-fiberedtopos}}

\symbollist{$\codim(\mathscr{Y},\mathscr{X})$}{Codimension}{{\bf \ref{ch-rigid}}.\ref{dfn-codimensionatpoint} (2)}

\symbollist{$\codim_x(\mathscr{Y},\mathscr{X})$}{Codimension at $x$}{{\bf \ref{ch-rigid}}.\ref{dfn-dimensionatpoint} (1)}

\symbollist{$\cone(f)^{\bullet}$}{Mapping cone}{{\bf \ref{ch-pre}}.\ref{dfn-homotopypre2}}

\symbollist{$\Cont(f)$}{Content ideal}{{\bf \ref{ch-pre}}, \S\ref{subsub-convpreadhnormalization}}


\symbollist{$\DC^{\ast}(\mathscr{A})$}{Derived category of $\mathscr{A}$}{{\bf \ref{ch-pre}}.\ref{dfn-derivedcategory3}}

\symbollist{$\DC^{\ast}(X)$}{Derived category of $\O_X$-modules}{{\bf \ref{ch-pre}}, \S\ref{subsub-ringedspacederivedcatformalism}}

\symbollist{$\DC^{\ast}_{\qcoh}(X)$}{--- with quasi-coherent cohomologies}{{\bf \ref{ch-pre}}, \S\ref{subsub-ringedspacederivedcatformalism}}

\symbollist{$\DC^{\ast}_{\coh}(X)$}{--- with coherent cohomologies}{{\bf \ref{ch-pre}}, \S\ref{subsub-ringedspacederivedcatformalism}}

\symbollist{$\D^n_{\mathscr{S}}$}{Unit disk over $\mathscr{S}$}{{\bf \ref{ch-rigid}}, \S\ref{subsub-unitdisk}}

\symbollist{$\Dercont_A(\cdot,\cdot)$}{Continuous derivations}{{\bf \ref{ch-formal}}, \S\ref{subsub-continuousderivations}}

\symbollist{$\dim(\mathscr{X})$}{Dimension}{{\bf \ref{ch-rigid}}.\ref{dfn-dimensionatpoint} (2)}

\symbollist{$\dim_x(\mathscr{X})$}{Dimension at $x$}{{\bf \ref{ch-rigid}}.\ref{dfn-dimensionatpoint} (1)}

\symbollist{$\End_{\mathscr{C}}(x)$}{Set of endomorphisms}{{\bf \ref{ch-pre}}, \S\ref{subsub-categoryconv}}

\symbollist{$G_x$}{Set of all generizations of $x$}{{\bf \ref{ch-pre}}, \S\ref{subsub-genspetopsp}}

\symbollist{$\Gamma_V$}{Value group}{{\bf \ref{ch-pre}}, \S\ref{subsub-valuation}}

\symbollist{$\Gamma_X$}{Global section functor}{{\bf \ref{ch-pre}}.\ref{prop-projlimsheafleftexact0}}

\symbollist{$\FN_{\mathscr{X},x}$}{Category of formal neighborhoods}{{\bf \ref{ch-rigid}}, \S\ref{subsub-fibersoverrigptslocalrings}}

\symbollist{$\Hom_{\mathscr{C}}(x,y)$}{Class of arrows from $x$ to $y$ in the category $\mathscr{C}$}{{\bf \ref{ch-pre}}, \S\ref{subsub-categoryconv}}

\symbollist{$\Homcont_A(\cdot,\cdot)$}{Continuous $A$-linear maps}{{\bf \ref{ch-formal}}.\ref{prop-genpairdifferential1}}

\symbollist{$\mathrm{ht}(\Gamma)$}{Height}{{\bf \ref{ch-pre}}, \S\ref{subsub-ordab}}

\symbollist{$\mathrm{ht}(V)$}{Height}{{\bf \ref{ch-pre}}.\ref{dfn-height1}}

\symbollist{$\Ht(F^{\bullet},G^{\bullet})$}{Set of arrows homotopic to zero}{{\bf \ref{ch-pre}}, \S\ref{subsub-homopotypre}}

\symbollist{$\Id(A)$}{Set of all ideals of $A$}{{\bf \ref{ch-pre}}, \S\ref{subsub-strcohtopsp}}

\symbollist{$\int_{X}F$}{Overconvergent interior}{{\bf \ref{ch-pre}}, \S\ref{subsub-gentopoverconvint}}

\symbollist{$\Isol(\Gamma)$}{Set of all isolated subgroups}{{\bf \ref{ch-pre}}, \S\ref{subsub-ordab}}

\symbollist{$\Isom_{\mathscr{C}}(x,y)$}{Set of isomorphisms}{{\bf \ref{ch-pre}}, \S\ref{subsub-categoryconv}}

\symbollist{$\KC^{\ast}(\mathscr{A})$}{Categories of complexes up to homotopy}{{\bf \ref{ch-pre}}.\ref{dfn-homotopypre2}}


\symbollist{$\varprojlim$}{Limit, projective limit}{{\bf \ref{ch-pre}}, \S\ref{subsub-limdefuniv}}

\symbollist{$\varinjlim$}{Colimit, inductive limit}{{\bf \ref{ch-pre}}, \S\ref{subsub-limdefuniv}}

\symbollist{$\varprojLim$}{$2$-categorical limit}{{\bf \ref{ch-pre}}, \S\ref{subsub-fiberedtopos}}

\symbollist{$\varinjLim$}{$2$-categorical colimit}{{\bf \ref{ch-pre}}, \S\ref{subsub-fiberedtopos}}

\symbollist{$\LT(f)$}{Leading term}{{\bf \ref{ch-pre}}, \S\ref{subsub-stbasis}}

\symbollist{$\FM_{\mathscr{X}}$}{Category of formal models}{{\bf \ref{ch-rigid}}.\ref{dfn-formalmodelcat}}

\symbollist{$\FM^{\dist}_{\mathscr{X}}$}{--- of distinguished formal models}{{\bf \ref{ch-rigid}}.\ref{dfn-cohrigidspacedist}}

\symbollist{$N(A)$}{Topologically nilpotent elements}{{\bf \ref{ch-pre}}, \S\ref{subsub-fadicringsgeneralities}}

\symbollist{$\mathscr{N}_{Y/X}$}{Conormal sheaf}{Ex.\ \ref{exer-conormalsheaf}}

\symbollist{$\nu(g)$}{Leading degree}{{\bf \ref{ch-pre}}, \S\ref{subsub-stbasis}}

\symbollist{$\O_{\mathscr{X}}$}{Rigid structure sheaf}{{\bf \ref{ch-rigid}}.\ref{dfn-ZRstrsheaf3}}

\symbollist{$\O^{\int}_{\mathscr{X}}$}{Integral structure sheaf}{{\bf \ref{ch-rigid}}.\ref{dfn-ZRstrsheaf1}}

\symbollist{$\obj(\mathscr{C})$}{Class of objects in the category $\mathscr{C}$}{{\bf \ref{ch-pre}}, \S\ref{subsub-categoryconv}}

\symbollist{$\Omega^1_{A/B}$}{Differential module}{{\bf \ref{ch-formal}}, \S\ref{subsub-differentialscantop}}

\symbollist{$\widehat{\Omega}^1_{X/Y}$}{Complete differential module}{{\bf \ref{ch-formal}}, \S\ref{subsub-diffcomplete}}

\symbollist{$\Omega^1_{X/Y}$}{Sheaf of differentials}{{\bf \ref{ch-formal}}, \S\ref{subsub-differentialinvformalsch}}

\symbollist{$\open(E)$}{Set of all isom.\ classes of subobj's of a final object}{{\bf \ref{ch-pre}}, \S\ref{subsub-topoipopints}}

\symbollist{$\Ouv(X)$}{Set of all open subsets of $X$}{{\bf \ref{ch-pre}}.\ref{exa-latticeexafoherent}}


\symbollist{$\P^{n,\an}_{\mathscr{S}}$}{Rigid analytic projective $n$-space over $\mathscr{S}$}{{\bf \ref{ch-rigid}}, \S\ref{subsub-examplesproj}}

\symbollist{$\Proj S$}{Homogeneous prime spectrum of $S$}{{\bf \ref{ch-pre}}, \S\ref{subsub-cohprerequi-0-projective}}

\symbollist{$\pts(E)$}{Set of all isomorphism classes of points of $E$}{{\bf \ref{ch-pre}}, \S\ref{subsub-topoipopints}}

\symbollist{$\QCOuv(X)$}{Set of all quasi-compact open subsets of $X$}{{\bf \ref{ch-pre}}.\ref{exa-latticeexafoherent}}

\symbollist{$\Qis^{\ast}(\mathscr{A})$}{Set of all quasi-isomorphisms}{{\bf \ref{ch-pre}}, \S\ref{subsub-derivedcategorydef}}

\symbollist{$s(\mathscr{X})$}{The associated scheme}{{\bf \ref{ch-rigid}}, \S\ref{sub-associatedschemes}}

\symbollist{$\sp(E)$}{Associated topological space of $E$}{{\bf \ref{ch-pre}}, \S\ref{subsub-topoipopints}}

\symbollist{$\Spa\mathscr{A}$}{Adic spectrum of $\mathscr{A}$}{{\bf \ref{ch-rigid}}, \S\ref{subsub-adicspectrum}}

\symbollist{$\Spec A$}{Prime spectrum of $A$}{{\bf \ref{ch-pre}}, \S\ref{subsub-strcohtopsp}}

\symbollist{$\Spf A$}{Formal spectrum of $A$}{{\bf \ref{ch-formal}}, \S\ref{subsub-formalnotformalspec}}

\symbollist{$\Spz A$}{Zariskian spectrum of $A$}{{\bf \ref{ch-formal}}.\ref{dfn-zariskianspectrum}}

\symbollist{$R(A,I)$}{Rees algebra}{{\bf \ref{ch-pre}}, \S\ref{sub-Reescone}}

\symbollist{$\mathrm{rat\textrm{-}rank}(\Gamma)$}{Rational rank}{{\bf \ref{ch-pre}}, \S\ref{subsub-ordab}}

\symbollist{$\mathrm{rat\textrm{-}rank}(V)$}{Rational rank}{{\bf \ref{ch-pre}}.\ref{dfn-height1}}

\symbollist{$\sigma^{\leq n}$, $\sigma^{\geq n}$}{Stupid truncations}{{\bf \ref{ch-pre}}, \S\ref{subsub-complexcategorytruncations}}

\symbollist{$\tau^{\leq n}$, $\tau^{\geq n}$}{Truncations}{{\bf \ref{ch-pre}}, \S\ref{subsub-complexcategorytruncations}}

\symbollist{$\top(X)$}{Topos associated to $X$}{{\bf \ref{ch-pre}}, \S\ref{subsub-toposasstopsp}}

\symbollist{$\mathrm{tr.deg}$}{Transcendental degree}{{\bf \ref{ch-pre}}, \S\ref{sub-estimate}}

\symbollist{$\widehat{\mathbf{V}}(\mathscr{E})$}{Vector bundle}{Ex.\ \ref{exer-vectorbundleformal}}

\symbollist{$\ZRT(\mathscr{X})$}{Zariski-Riemann triple associated to $\mathscr{X}$}{{\bf \ref{ch-rigid}}.\ref{dfn-ZRtriple}}

}

\subsection*{--- Spaces}

{\small 
\symbollist{$\A^n_S$}{Affine $n$-space over $S$}{Ex.\ \ref{exer-affinespaceformal} (1)}

\symbollist{$\widehat{\A}^n_S$}{Affine $n$-space over $S$}{Ex.\ \ref{exer-affinespaceformal} (1)}

\symbollist{$\G_{\mathrm{a}}^{\an}$}{Additive group}{{\bf \ref{ch-rigid}}.\ref{exas-gagm} (1)}

\symbollist{$\G_{\mathrm{m}}^{\an}$}{Multiplicative group}{{\bf \ref{ch-rigid}}.\ref{exas-gagm} (2)}

\symbollist{$\P^n_S$}{Projective $n$-space over $S$}{{\bf \ref{ch-pre}}, \S\ref{subsub-cohprerequi-0-projective}}

\symbollist{$\widehat{\P}^n_S$}{Projective $n$-space over $S$}{Ex.\ \ref{exer-affinespaceformal} (2)}
}

\subsection*{--- Maps, morphisms, arrows}

{\small 
\symbollist{$\dashrightarrow$}{Rational map}{{\bf \ref{ch-rigid}}.\ref{dfn-birationalgeom1} (1)}

\symbollist{$\sep_X$}{Separation map}{{\bf \ref{ch-pre}}, \S\ref{subsub-separationgen}}

\symbollist{$\sep_{\mathscr{X}}$}{Separation map}{{\bf \ref{ch-rigid}}, \S\ref{subsub-separation}}

\symbollist{$\sp_X$}{Specialization map}{{\bf \ref{ch-rigid}}, \S\ref{subsub-ZRdef}}

}

\subsection*{--- Super- and subscripts, brackets, etc.}

{\small 
\symbollist{{\sl Symbol}}{{\sl Meaning}}{{\sl See}}

\symbollist{$\widehat{A}$}{Adic completion of $A$}{{\bf \ref{ch-pre}}, \S\ref{subsub-Iadiccompletioncomplete}}

\symbollist{$\het{A}$}{Henselization of $A$}{{\bf \ref{ch-pre}}, \S\ref{subsub-henselianpairs}}

\symbollist{$A^o$}{Power-bounded elements}{{\bf \ref{ch-pre}}, \S\ref{subsub-fadicringsgeneralities}}

\symbollist{$\zat{A}$}{Associated Zariskian of $A$}{{\bf \ref{ch-pre}}, \S\ref{subsub-zariskianpairs}}

\symbollist{$\mathscr{C}^{\opp}$}{Opposite category of $\mathscr{C}$}{{\bf \ref{ch-pre}}, \S\ref{subsub-categoryconv}}

\symbollist{$(\mathscr{D}^{\leq 0},\mathscr{D}^{\geq 0})$}{$t$-structure}{{\bf \ref{ch-pre}}, \S\ref{sub-triangulatedcategory}}

\symbollist{$[k]$}{Shift functor}{{\bf \ref{ch-pre}}, \S\ref{subsub-complexcategoryshifts}}

\symbollist{$\widehat{M}$}{Adic completion}{{\bf \ref{ch-pre}}.\ref{dfn-Iadiccompletiondelicate1}}

\symbollist{$M^{\wedge}_{F^{\bullet}}$}{Hausdorff completion}{{\bf \ref{ch-pre}}, \S\ref{subsub-completionfiltration}}

\symbollist{$M^{\Delta}$}{$\Delta$-sheaf}{{\bf \ref{ch-formal}}, \S\ref{subsub-thedeltasheaves}}

\symbollist{$M^{\for}$}{Formal completion}{{\bf \ref{ch-formal}}, \S\ref{subsub-GFGAcomannconst}}

\symbollist{$M_{\Itor}$}{$I$-torsion part}{{\bf \ref{ch-pre}}, \S\ref{subsub-pairstorsionsandsaturation}}

\symbollist{$[X]$}{Separated quotient}{{\bf \ref{ch-pre}}, \S\ref{subsub-separationgen}}

\symbollist{$\widehat{X}|_Y$}{Formal completion along $Y$}{{\bf \ref{ch-pre}}, \S\ref{sub-schpair}}

\symbollist{$\het{X}|_Y$}{Henselization along $Y$}{{\bf \ref{ch-pre}}, \S\ref{sub-schpair}}

\symbollist{$X^{\an}$}{Analytification}{{\bf \ref{ch-rigid}}, \S\ref{subsub-GAGAfunctorconst}}

\symbollist{$X^{\opp}$}{Ordered set with the reversed ordering}{{\bf \ref{ch-pre}}, \S\ref{subsub-orderings}}

\symbollist{$X^{\rig}$}{Associated rigid space}{{\bf \ref{ch-rigid}}.\ref{dfn-cohrigidspace1}}

\symbollist{$X^{\mathrm{ref}}$}{Reflexivization}{{\bf \ref{ch-pre}}.\ref{thm-reflexivization}}

\symbollist{$X^{\sob}$}{Associated sober space}{{\bf \ref{ch-pre}}, \S\ref{subsub-sober}}

\symbollist{$\zat{X}|_Y$}{Associated Zariskian aalong $Y$}{{\bf \ref{ch-pre}}, \S\ref{sub-schpair}}

\symbollist{$\mathscr{X}_{\ad}$}{Small admissible site}{{\bf \ref{ch-rigid}}.\ref{dfn-admissiblecovering12}}

\symbollist{$\mathscr{X}^{\adic}_{/\mathscr{Y}}$}{Adic part}{{\bf \ref{ch-rigid}}.\ref{subsub-adicpart}}

\symbollist{$\ZR{\mathscr{X}}$}{Associated Zariski-Riemann space}{{\bf \ref{ch-rigid}}, \S\ref{subsub-ZRdef}}

\symbollist{$\ZR{\mathscr{X}}^{\cl}$}{Classical points}{{\bf \ref{ch-rigid}}, \S\ref{subsub-classicalpoints}}

\symbollist{$\ZR{X}_U$}{Classical Zariski-Riemann space}{{\bf \ref{ch-rigid}}.\ref{dfn-classicalZRsp1}}

\symbollist{$\ZR{U}^X_{\pc}$}{Partial compatification}{{\bf \ref{ch-rigid}}.\ref{dfn-partialcompnagata}}

\symbollist{$\ZR{X}_{\cpt}$}{Canonical compatification}{{\bf \ref{ch-rigid}}, \S\ref{subsub-canonicalcompnagata}}

\symbollist{$\|\cdot\|_{x,\mathscr{I},c}$}{Seminorm at $x$}{{\bf \ref{ch-rigid}}, \S\ref{subsub-seminormsofpoints}}

\symbollist{$\|\cdot\|_{\Sp,\mathscr{I},c}$}{Spectral seminorm}{{\bf \ref{ch-rigid}}, \S\ref{subsub-spectralseminorms}}

}

\subsection*{--- Special symbols}

{\small 
\symbollist{{\sl Symbol}}{{\sl Meaning}}{{\sl See}}

\symbollist{$\mathbf{2}$}{Boolean lattice formed by two elements $0$ and $1$}{{\bf \ref{ch-pre}}, \S\ref{subsub-strcohtopsp}}

\symbollist{$\mathsf{U}$}{Grothendieck universe}{{\bf \ref{ch-pre}}, \S\ref{subsub-existuniverse}}

}

\addcontentsline{toc}{chapter}{\rm Index}
\printindex
\end{document}